\documentclass{amsbook}
\usepackage{amsmath,a4wide,amssymb,xcolor,pict2e}
\usepackage{graphics,bm,array}

\usepackage{imakeidx}
\usepackage{pgf}
\usepackage{hyperref}

\usepackage[all]{xy}

\newcommand*{\parallelogram}{%
  \rlap{\rotatebox{-30}{\rule[.05ex]{.4pt}{.77em}}}%
  \kern.04em%
  \rlap{\kern.36em\raisebox{0.649519052835em}{\rule{.7em}{.4pt}}}%
  \rule{.7em}{.4pt}\kern-.04em%
  \rotatebox{-30}{\rule[.05ex]{.4pt}{.77em}}}

\newcommand{\subparallel}{\,|\mbox{\hskip-1pt}\big|\,}

\newcommand{\ssq}{\subseteq}

\def\defterm#1{{\color{purple}\em #1}}

\newcommand{\defeq}{:=} 

\newcommand{\Ra}{\Rightarrow}

\newcommand{\IP}{\mathbb P}
\newcommand{\IN}{\mathbb N}
\newcommand{\IR}{\mathbb R}
\newcommand{\IZ}{\mathbb Z}
\newcommand{\IQ}{\mathbb Q}
\newcommand{\IC}{\mathbb C}

\newcommand{\IA}{\mathbb A}
\newcommand{\IF}{\mathbb F}
\newcommand{\IH}{\mathbb H}
\newcommand{\A}{\mathcal A}
\newcommand{\C}{\mathcal C}
\newcommand{\F}{\mathcal F}
\newcommand{\I}{\mathcal I}
\newcommand{\M}{\mathcal M}

\newcommand{\e}{\varepsilon}
\newcommand{\Sym}{\mathsf{Sym}}
\newcommand{\Alt}{\mathsf{Alt}}
\newcommand{\w}{\omega}
\newcommand{\two}{2}
\newcommand{\Fix}{\mathsf{Fix}}
\newcommand{\LT}{\dddot X}

\newcommand{\Hull}{\mathsf H}

\newcommand{\prop}{\Join}

\newcommand{\projeq}{\;{=}\hskip-6pt{=}\hskip-8pt\raisebox{-1.8pt}{$\wedge$}\;}
\newcommand{\projeqind}{\!{}_{{=}\hskip-5.5pt{}_\wedge}}
\newcommand{\projupind}{\!{}^{{=}\hskip-5.5pt{}_\wedge}}
\newcommand{\pr}{\mathrm{pr}}
\newcommand{\supp}{\mathrm{supp}}

\newcommand{\Tau}{\mathcal T}
\newcommand{\vecv}{\vec{\boldsymbol v}}
\newcommand{\har}{\mathsf{char}}
\newcommand{\Ker}{\mathrm{Ker}}
\newcommand{\End}{\mathrm{End}}
\newcommand{\LB}{\mathsf{LB}}
\newcommand{\Lenz}{\mathsf{L}}

\newcommand{\vertiii}[1]{{\left\vert\kern-0.25ex\left\vert\kern-0.25ex\left\vert #1 
    \right\vert\kern-0.25ex\right\vert\kern-0.25ex\right\vert}}

\newcommand{\Liners}{\mathsf{Liners}}
\makeatletter
\DeclareRobustCommand{\cev}[1]{%
  {\mathpalette\do@cev{#1}}%
}
\newcommand{\do@cev}[2]{%
  \vbox{\offinterlineskip
    \sbox\z@{$\m@th#1 x$}%
    \ialign{##\cr
      \hidewidth\reflectbox{$\m@th#1\vec{}\mkern4mu$}\hidewidth\cr
      \noalign{\kern-\ht\z@}
      $\m@th#1#2$\cr
    }%
  }%
}
\makeatother

\newcommand{\dom}{\mathsf{dom}}
\newcommand{\rng}{\mathsf{rng}}

\def\lines#1{\overline{\,{#1}\,}}

\def\Aline#1#2{\overline{#1\phantom{\mbox{\footnotesize$|$}}#2}}
\def\Line#1#2{\overline{#1\phantom{\mbox{\footnotesize$|$}}#2}}

\def\overvector#1{\overrightarrow{#1}}
\def\vecbold#1{\vec{\boldsymbol {#1}}}
\newcommand\Af{\mathsf L}

\newcommand{\Aut}{\mathsf{Aut}}
\newcommand{\SAut}{\mathsf{SAut}}
\newcommand{\Aff}{\mathsf{Aff}}
\newcommand{\EAff}{\mathsf{EAff}}
\newcommand{\SAff}{\mathsf{SAff}}

\newcommand{\Dil}{\mathsf{Dil}}
\newcommand{\Trans}{\mathsf{Tra}}
\newcommand{\Hol}{\mathsf{Hol}}
\newcommand{\AGL}{\mathsf{AGL}}

\newcommand{\Proj}{\mathsf{Proj}}
\newcommand{\TT}{\mathsf{I}}
\newcommand{\NT}{\mathsf{N}}
\newcommand{\CT}{\mathsf{X}}
\newcommand{\ST}{\mathsf{T}}
\newcommand{\TSH}{\mathsf{INX}}
\newcommand{\TS}{\mathsf{IN}}

\newcommand{\SC}{\mathsf{T}}
\newcommand{\INXT}{\mathsf{INXT}}

\newcommand{\Syl}{\mathsf{Syl}}
\newcommand{\Det}{\mathrm{Det}}
\newcommand{\Stab}{\mathsf{Stab}}

\newcommand{\fdom}{\mbox{\rm\textborn}}
\newcommand{\fran}{\mbox{\rm\textdied}}
\newcommand{\Mor}{\mathsf{Mor}}
\newcommand{\Ob}{\mathsf{Ob}}
\newcommand{\LG}{\mathsf{LG}}

\newcommand{\PSC}{\mathsf F}
\newcommand{\invol}{
\begin{picture}(9,10)(0,2)
\put(5,0){\circle*{1}}
\put(0,5){\circle*{1}}
\put(6,1){\circle*{1}}
\put(1,6){\circle*{1}}
\put(7,2){\circle*{1}}
\put(2,7){\circle*{1}}
\put(8,3){\circle*{1}}
\put(3,8){\circle*{1}}
\put(9,4){\circle*{1}}
\put(4,9){\circle*{1}}
\end{picture}
}
\newcommand{\Invol}{
\begin{picture}(10,10)(0,2)
\put(5,0){\line(1,1){5}}
\put(0,5){\line(1,1){5}}
\end{picture}
}

\newcommand{\Dscale}{
\begin{picture}(9,9)(0,2)
    \linethickness{0.5mm}
  \polygon*(0,1)(0,9)(8,9)
  \polygon*(1,0)(9,0)(9,8)
\end{picture}
}

\newcommand{\discale}
{\begin{picture}(8,8)(0,2)\put(1,8){\line(1,-1){6}}\end{picture}}

\newcommand{\reshitka}
{\begin{picture}(8,8)(0,1.3)\put(2,0){\line(0,1){8}}\put(6,0){\line(0,1){8}}
\put(0,2){\line(1,0){8}}\put(0,6){\line(1,0){8}}
\end{picture}}

\newcommand{\didots}
{\begin{picture}(8,8)(0,1.3)\put(1.5,6.5){\circle*{2}}\put(6.5,1.5){\circle*{2}}\end{picture}}

\newcommand{\cdidots}
{\begin{picture}(11,8)(0,1.3)\put(2,4){\circle{3.5}}\put(6,7){\circle*{2}}\put(10,2){\circle*{2}}\end{picture}}

\newcommand{\Boolean}{
\begin{picture}(9,8)(0,1)
\linethickness{=0.6pt}
\put(0,0){\color{teal}\line(1,0){8}}
\put(0,8){\color{teal}\line(1,0){8}}
\put(0,0){\color{cyan}\line(0,1){8}}
\put(8,0){\color{cyan}\line(0,1){8}}
\put(0,0){\color{red}\line(1,1){8}}
\put(0,8){\color{red}\line(1,-1){8}}
\end{picture}
}
\newcommand{\Pentagon}{\begin{picture}(13,12)(-7,-3.5)
\linethickness{0.7pt}
\put(-5.7,1.9){\color{teal}\line(1,0){11.4}}
\put(-3.5,-4.9){\color{teal}\line(1,0){7}}
\put(-5.7,1.9){\color{orange}\line(57,41){5.7}}
\put(-3.5,-4.9){\color{orange}\line(57,41){9.2}}
\put(5.7,1.9){\color{olive}\line(-57,41){5.7}}
\put(3.5,-4.9){\color{olive}\line(-57,41){9.2}}
\put(-3.5,-4.9){\color{blue}\line(35,109){3.5}}
\put(3.5,-4.9){\color{blue}\line(11,34){2.2}}
\put(3.5,-4.9){\color{violet}\line(-35,109){3.5}}
\put(-3.5,-4.9){\color{violet}\line(-11,34){2.2}}

\put(0,6){\circle*{1.5}}
\put(-5.7,1.9){\circle*{1.5}}
\put(-3.5,-4.9){\circle*{1.5}}
\put(3.5,-4.9){\circle*{1.5}}
\put(5.7,1.9){\circle*{1.5}}
\end{picture}
}
\newcommand{\puls}{\,\raisebox{0pt}{\small$\heartsuit$}\,}
\newcommand{\plus}{+}

\title{Linear Geometry and Algebra}
\author{Taras Banakh}

\newtheorem{theorem}{Theorem}[section]
\numberwithin{theorem}{section}
\newtheorem{theorems}[theorem]{Theorem$^\dag$}
\newtheorem{corollary}[theorem]{Corollary}
\newtheorem{corollarys}[theorem]{Corollary$^\dag$}
\newtheorem{proposition}[theorem]{Proposition}
\newtheorem{propositions}[theorem]{Proposition$^\dag$}
\newtheorem{lemma}[theorem]{Lemma}
\newtheorem{lemmas}[theorem]{Lemma$^\dag$}
\newtheorem{claim}[theorem]{Claim}

\newtheorem{problem}[theorem]{Problem}
\newtheorem{conjecture}[theorem]{Conjecture}

\theoremstyle{definition}

\newtheorem{exercise}[theorem]{Exercise}
\newtheorem{Exercise}[theorem]{Exercise$^\star$}
\newtheorem{remark}[theorem]{Remark}
\newtheorem{remarks}[theorem]{Remark$^\dag$}
\newtheorem{definition}[theorem]{Definition}
\newtheorem{example}[theorem]{Example}
\newtheorem{examples}[theorem]{Example$^\dag$}

\newtheorem{question}[theorem]{Question}

\makeindex[name=person,title={Person Index}]%
\makeindex[name=note,title={Notation Index}]
\makeindex

\begin{document}

\maketitle

\setcounter{tocdepth}{1}
\tableofcontents

\mainmatter

\chapter*{Preface}

\rightline{\em Make things as simple as possible, but not simpler.}

\rightline{\index[person]{Einstein}Albert Einstein\footnote{{\bf Albert Einstein} (1879 -- 1955) was a German-born theoretical physicist who is widely held to be one of the greatest and most influential scientists of all time. Best known for developing the theory of relativity, Einstein also made important contributions to quantum mechanics, and was thus a central figure in the revolutionary reshaping of the scientific understanding of nature that modern physics accomplished in the first decades of the twentieth century.}}

\vskip15pt

%

\rightline{\em Entia non sunt multiplicanda praeter necessitatem\footnote{Entities must not be multiplied beyond necessity.}.}

\rightline{\index[person]{Occam}William Occam\footnote{{\bf William of Occam} (1287 -- 1347) was an English Franciscan friar, scholastic philosopher, apologist, and Catholic theologian, who is believed to have been born in Ockham, a small village in Surrey. He is considered to be one of the major figures of medieval thought and was at the centre of the major intellectual and political controversies of the 14th century. He is commonly known for Occam's razor, the methodological principle that bears his name, and also produced significant works on logic, physics and theology. William is remembered in the Church of England with a commemoration on the 10th of April.}}

\vskip15pt
\rightline{\em Let no one ignorant of geometry enter here.}

\rightline{--- Inscription over the entrance to Plato’s Academy}
\vskip15pt








\vskip40pt

Linear Geometry studies geometric properties which can be expressed via the notion of a line. All information about lines is encoded in a ternary relation called a line relation. A set endowed with a line relation is called a liner. So, Linear Geometry studies liners. Imposing some additional axioms on a liner, we obtain some special classes of liners:  regular, projective, affine, proaffine, etc. Linear Geometry includes Affine and Projective Geometries and is a part of Incidence Geometry.

The aim of this book is to present a self-contained logical development of Linear Geometry, starting with some intuitively acceptable geometric axioms and ending with algebraic structures that necessarily arise from studying the structure of geometric objects that satisfy those simple and intuitive geometric axioms. We shall meet many quite exotic algebraic structures that arise this way: magmas, loops, ternars, quasi-fields, alternative rings, procorps, profields, etc. We strongly prefer (synthetic) geometric proofs and use tools of analytic geometry only when no purely geometric proof is available. 

Also our approach is more affine oriented: we first clarify the affine situation as more intutive and afterwards extend the obtained affine results to projective case. Such an affine-oriented approach allowed us to prove that the double degeneration of the Desargues Axiom is equivalent to the disjunction of the Fano and Moufang Axioms and also to the commutativity of plus and puls operations in all ternars of a projective plane. All the truth about the interplay between those properties was not clearly elaborated since times of Ruth Moufang who believed that the single and double degenerations of the Desargues Axiom are equivalent (by the way, this is still an open problem, whether Fano projective planes are Moufang, a general belief is that an example of non-Moufang Fano projective plane should exist but nobody has succeeded to construct such a pathological plane, yet).
Many results in this book are proved for proaffine liners which include affine and projective liners as special cases. In particular, we tried to find definitions of many well-known geometric properties (like those of Pappus, Desargues, Moufang, or Fano) that behave equally well is both affine and projective contexts. Linear Geometry also is fascinating as a rather old-fashioned subject that still has many very difficult and simply-formulated unsolved problems. For example, nobody knows if a projectve plane of order 12 exists, or whether every homogeneous finite projective plane is Desarguesian or whether a finite rigid projective plane exists (by the way, a finite rigid affine plane does exist, of order 25). The list of such selected open problems is included in the last Chapter~\ref{ch:open-prob}. 

Theorems whose proof involves some external results (unproven in this book)  are marked with the $\dag$ symbol. Such theorems$^\dag$ appear in Chapter~\ref{ch:elementary-Moufang} analyzing the structure of elementary Moufang loops and also at the very end of Chapter~\ref{ch:permutation} that discusses the classification of $2$-transitive finite groups (whose proof involves the Classification of Finite Simple Groups).  
\smallskip

Linear Geometry has been developed by many great mathematicians since times of Antiquity (Thales, Euclides, Proclus, Pappus), through Renaissance (Descartes, Desargues), Early Modernity (Playfair, Gauss, M\"obius, Lobachevski,  Steiner, von Staudt, Poncelet, Bolyai), Late Modernity Times (Beltrami, Jordan, Klein, Galucci, Hilbert, Steinitz, Fano, Hessenberg, Veblen, Wedderburn, Moufang, Lenz, Barlotti) till our contemporaries (Hartshorne, Hall, Buekenhout, Gleason, Kantor, Doyen, Hubault, Dembowski, Klingenberg, Grundh\"ofer, M\"uller, Nagy). In footnotes in relevant places we add a short biographic information on geometers that contributed to development of this fascinating field of mathematics, which is very helpful and instructive for understanding the internal and historical logic of development of Linear Geometry.

I would like to express my thanks for the help in writing this text to Ivan Hetman, Alex Ravsky, Vlad Pshyk, Oksana Skygar, Alex Ilchuk, the active participants of the seminar in Foundations of Geometry at Ivan Franko National University of Lviv, and its online ZOOM-continuation (Yurii Yarosh, Zoriana Novosad, Iryna Yegorchenko, Jauhien Piatlicki).

Finally, I thank the KSE Program ``Talents for Ukraine''\footnote{\tt https://150.foundation.kse.ua/en/about} for a grant support of this interesting but extremely time-consuming project.

\bigskip

\rightline{\em Lviv, March $2023$ -- March $2026$}

\chapter{Liners}

We recal that an \index{$n$-ary relation}\defterm{$n$-ary relation} on a set $X$ is a subset of the $n$-power $X^n$ of $X$. Elements of $X^n$ are functions $x:n\to X$ from the natural number $n\defeq\{0,\dots,n-1\}$, which can be identified with $n$-tuples $(x_0,x_1,\dots,x_{n-1})$ that will be denoted by $x_0x_1\cdots x_{n-1}$ (without commas and parentheses). For a relation $R\subseteq X^n$ and an $n$-tuple $x_0\dots x_{n-1}$ we write $Rx_0\cdots x_{n-1}$ iff $(x_0,\dots,x_{n-1})\in R$. So, we can think of an $n$-ary relation $R$ as the $n$-ary predicate $Rx_0\cdots x_{n-1}$, i.e., a  function $R:X^n\to\{{\sf false},{\sf true}\}$ with logical values.

\section{Line relations}\label{s:line-relations}

\begin{definition}\label{d:liner} A ternary relation $\Af\subseteq X^3 $ on a set $X$ is called \index[note]{$\mathsf Lxyz$}\index{line relation}\defterm{a line relation} or else a \index{liner structure}\defterm{liner structure} on $X$ if it satisfies the following axioms
\begin{itemize}
\item[{\sf (IL)}] \index{Axiom!Identity}{\sf Identity:} $\forall x,y\in X\;\;\big(\Af xyx\;\to\; x=y\big)$;
\item[{\sf (RL)}] \index{Axiom!Reflexivity}{\sf Reflexivity:} $\forall x,y\in X\;\;(\Af xxy\;\wedge\;\Af xyy)$;
\item[{\sf (EL)}] \index{Axiom!Exchange}{\sf Exchange:} $\forall a,b,x,y\in X\;\big(( \Af axb\wedge \Af ayb\wedge x\ne y)\to (\Af xay\wedge\Af xby)\big)$.
\end{itemize}
A set $X$ endowed with a line relation $\Af\subseteq X^3$ is called a \index{liner}\defterm{liner}.
\end{definition}

\begin{remark} Intuitively, the relation $\Af xyz$ between points $x,y,z$ of a liner $(X,\Af)$ means that the point $y$ belongs to the line determined by the points $x$ and $z$. If $x=z$, this line degenerates to the singleton $\{x,z\}=\{x\}=\{z\}$.
\end{remark}

\begin{proposition}\label{p:permutation} Let $x,y,z\in X$ be three pairwise distinct points in a liner $(X,\Af)$. If $\Af xyz$, then $\Af xzy\wedge \Af yxz\wedge\Af yzx\wedge \Af zxy\wedge\Af zyx$. 
\end{proposition}

\begin{proof} Assume that $\Af xyz$. The axiom {\sf (RL)} ensures that $\Af xzz$. By the Exchange Axiom {\sf(EL)}, $\Af xyz\wedge \Af xzz$ implies $\Af yxz\wedge\Af zxy$. By the same reason, $\Af xyz\wedge \Af xxz$ implies $\Af xzy\wedge \Af yzx$. Finally, $\Af yzx\wedge \Af zzx$ implies $\Af zyx$, by the axiom {\sf(EL)}.
\end{proof}

\begin{exercise} Prove that a ternary relation $\Af\subseteq X^3$ on a set $X$ is a line relation if and only if  it satisfies the axioms {\sf(IL)}, {\sf(RL)} and 
\begin{itemize}
\item[{\sf (EL')}] $\forall a,b,x,y\in X\;\big(( \Af axb\wedge \Af ayb\wedge x\ne y)\to \Af xay\big)$.
\end{itemize}
\end{exercise}

\begin{corollary}\label{c:SA} For any points $x,y,z\in X$  in a liner $(X,\Af)$, we have $\Af xyz\leftrightarrow\Af zyx$.
\end{corollary}

\begin{proof} It suffices to prove that $\Af xyz\Ra \Af zyx$. So, assume that $\Af xyz$. If $x=z$, then the axiom {\sf (IL)} ensures that $x=y=z$ and the axiom {\sf (RL)} ensures that $\Af zyx$. So, we assume that $x\ne z$. If $y\in\{x,z\}$, then $\Af zyx$ by the axiom {\sf (RL)}. If $y\notin\{x,z\}$, then $x,y,z$ are pairwise distinct points and $\Af zyx$ follows from Proposition~\ref{p:permutation}.
\end{proof}





\begin{exercise}\label{ex:subspace} Let $(X,\Af)$ be a liner. Prove that for every subset $Y\subseteq X$,  the relation $\Af{\restriction}_Y\defeq\Af\cap Y^3$ is a line relation on $Y$, and the pair $(Y,\Af\cap Y^3)$ is a liner.
\end{exercise} 

\begin{exercise} Prove that the axioms {\sf (IL)}, {\sf (RL)} and {\sf (EL)} are independent.
\end{exercise}

Exercise~\ref{ex:subspace} allows us to introduce the definition of a lineal subspace of a liner.

\begin{definition} Let $(X,\Af)$ be a liner and $Y\subseteq X$ be a subset of $X$. The liner $(Y,\Af\cap Y^3)$ is called a \index{subliner}\defterm{subliner} of the liner $(X,\Af)$.
\end{definition}

\section{Lines in liners}

\begin{definition} For two points $a,b\in X$ of a liner $(X,\Af)$, the set
$$\overline{a\Af b}\defeq\{x\in X:\Af axb\}$$is called  \defterm{the $\Af$-line} passing through the points $a,b$. If the line relation $\Af$ is clear from the context, then we shall omit the symbol $\Af$ in the notation of an $\Af$-line, a shall write $\Aline ab$ instead of $\overline{a\Af b}$.\index[note]{$\Aline xy$}
\end{definition}

\begin{theorem}\label{t:Alines} If $X$ is a liner, then 
\begin{enumerate}
\item $\forall x\in X\;\;\Aline xx=\{x\}$;
\item $\forall x,y\in X\;\;\{x,y\}\subseteq\Aline xy=\Aline yx$;
\item $\forall x,y\in X\;\;\forall u,v\in\Aline xy\;\;(u\ne v\Ra \Aline uv=\Aline xy)$.
\end{enumerate}
\end{theorem}

\begin{proof} 1. For every $x\in X$, the equality $\Aline xx=\{x\}$ follows from the axiom {\sf(IL)}.
\smallskip

2. For every $x,y\in X$, the equality  $\Aline xy=\Aline yx$ follows from Corollary~\ref{c:SA}. The inclusion $\{x,y\}\subseteq\Aline xy$ follows from the axiom {\sf (RL)}. 
\smallskip

3. Given any points $x,y\in X$ and $u,v\in\Aline xy$, assume that $u\ne v$. In this case, the axiom {\sf(IL)} ensures that $x\ne y$. To show that $\Aline xy\subseteq\Aline uv$, take any point $z\in\Aline xy$. If $z\in\{u,v\}$, then $z\in\Aline uv$ by the axiom {\sf (RL)}. So, assume that $z\notin\{u,v\}$. In this case, $\Af xuy\wedge \Af xzy$ implies $\Af uyz$, by the axiom {\sf (EL)}. By {\sf (RL)}, $\Af uuz$. By {\sf(EL)}, $\Af uuz\wedge \Af uyz$ implies $\Af uzy$. 

Next, we show that $\Af uvy$. Since $u\ne z$, the axiom {\sf(IL)} ensures that $u\ne y$. By the axiom {\sf (EL)}, $\Af xuy$ and $\Af xvy$ imply $\Af uyv$. If $y=v$, then $\Af uyv$ implies $\Af uvy$. If $y\ne v$, then the points $u,v,y$ are pairwise distinct and $\Af uvy$ follows from $\Af uyv$, by Proposition~\ref{p:permutation}. 

Therefore, we have $\Af uzy\wedge \Af uvy$. Applying the axiom {\sf (EL)}, we obtain $\Af vuz$. By Proposition~\ref{p:permutation}, $\Af vuz$ implies $\Af uzv$ and hence $z\in\Aline uv$. This completes the proof of the inclusion $\Aline xy\subseteq\Aline uv$. By analogy we can prove that $\Aline uv\subseteq\Aline xy$.
\end{proof}

\begin{definition} A subset $L$ of a liner $X$ is called a \index{line}\defterm{line} if $L=\Aline xy$ for some distinct points $x,y\in X$.
\end{definition} 

\begin{definition} A liner $X$ is called \index{line-finite liner}\index{liner!line-finite}\defterm{line-finite} if every line in $X$ is finite.
\end{definition}

\begin{exercise} Construct an example of an infinite line-finite liner.
\smallskip

{\em Hint:} Take any infinite-dimensional linear space over a finite field.
\end{exercise}


\begin{definition} Let $\kappa$ be a cardinal. A liner $(X,\Af)$ is defined to be \index{$\kappa$-long liner}\index{liner!$\kappa$-long}\defterm{$\kappa$-long} if every line $L$ in $X$ has cardinality $|L|\ge\kappa$.
\end{definition}

\begin{remark} The definition of a line implies that every liner is $2$-long.
\end{remark}

Theorem~\ref{t:Alines}(3) implies that every line $L$ in a liner $(X,\Af)$ is equal to the $\Af$-line $\Aline xy$ passing through any distinct points $x,y\in L$. Therefore, any line in a liner is uniquely determined by any distinct points of the line. 

\begin{definition} Two lines $L,\Lambda$ in a liner are defined to be \index{concurrent lines}\index{lines!concurrent}\defterm{concurrent} if $L\cap\Lambda$ is a singleton.
\end{definition}

Observe that two lines are concurrent if and only if they are distinct and have a common point.

\begin{definition} Points $x_1,x_2,\dots,x_n$ in a liner are defined to be \index{colinear points}\index{points!colinear}\defterm{colinear} if there exists a line $L$ in $X$ such that $\{x_1,x_2,\dots,x_n\}\subseteq L$.
\end{definition}

The following theorem shows that the family of all lines carries all information about a liner.

\begin{theorem}\label{t:L1+L2} For every liner $(X,\Af)$, the family $\mathcal L$ of all lines in $X$ has the following three properties:
\begin{enumerate}
\item[{\sf(L1)}] any distinct points $x,y\in X$ belong to a unique line $L\in\mathcal L$;
\item[{\sf(L2)}] every line $L\in\mathcal L$ contains at least two distinct points of $X$;
\item[{\sf(L3)}] $\Af=\big\{(x,y,z)\in X^3:(x=z\to y=x)\wedge\big(x\ne z\to (\exists L\in\mathcal L\;\;\{x,y,z\}\subseteq L)\big)\big\}$.
\end{enumerate}
Conversely, for every family of subsets $\mathcal L$ of a set $X$ satisfying the axiom {\sf(L1)}, {\sf(L2)}, the ternary relation $\Af$ defined by the formula {\sf(L3)} is a line relation on $X$ such that $\mathcal L$ coincides with the family of all lines of the liner $(X,\Af)$.
\end{theorem}

\begin{proof} The properties {\sf (L1)} and {\sf(L2)} follow from Theorem~\ref{t:Alines}(3,2), and ${\sf(L3)}$ follows from the definition of a line in the liner $(X,\Af)$.

Now assume that a family $\mathcal L$ of subsets of a set $X$ satisfies the axioms {\sf(L1)} and {\sf(L2)}, and define a ternary relation $\Af$ by the formula {\sf(L3)}, which implies that the relation $\Af$ satisfies the axiom {\sf(IL)}. The axiom {\sf(L1)} implies that the relation $\Af$ satisfied the axiom {\sf (RL)}. To see that the relation $\Af$ satisfies the Exchange Axiom {\sf(EL)}, take any points $a,b\in X$ and distinct points $x,y\in\Aline ab$. The definition of $\Af$ ensures that $a\ne b$ and there exist sets $L_x,L_y\in\mathcal L$ such that $\{a,x,b\}\subseteq L_x$ and $\{a,y,b\}\subseteq L_y$. Since $\{a,b\}\subseteq L_x$ and $\{a,b\}\subseteq L_y$, the axiom {\sf(L1)} ensures that $L_x=L_y$ and hence the set $L\defeq L_x=L_y\in\mathcal L$ witnesses that $\Af xay\wedge \Af xby$. Therefore, $\Af$ is a line relation. 

The equality $\mathcal L=\{\Aline ab:a,b\in X\;\wedge\; a\ne b\}$ can be easily derived from the axioms {\sf(L1), (L2)} and  the following lemma.

\begin{lemma} For every set $L\in\mathcal L$ and any distinct points $a,b\in L$ we have $L=\Aline ab$.
\end{lemma}

\begin{proof} The inclusion $L\subseteq\Aline ab$ follows from the definition of the relation $\Af$. To see that $\Aline ab\subseteq L$, take any point $x\in \Aline ab$. By the definition of the relation $\Af$, there exists a set $L_x\in\mathcal L$ such that $\{a,x,b\}\subseteq L$. Since $\{a,b\}\subseteq L\cap L_x$, the axiom {\sf(L1)} guarantees that $x\in L_x=L$ and hence $\Aline ab\subseteq L$.
\end{proof}
\end{proof} 

\begin{remark}Theorem~\ref{t:L1+L2} provides an alternative way of defining a liner as a set $X$ endowed with a family of lines $\mathcal L$  satisfying the axioms (L1), (L2) of Theorem~\ref{t:L1+L2}. This alternative way is widely used in literature and also will be used later in this book for defining various exaples of liners.
\end{remark}

\section{Closed sets in liners}

\begin{definition} A subset $A$ of a liner $X$ is called \defterm{closed} if for any distinct points $a,b,c,d\in A$ with $\Aline ab\ne \Aline cd$, we have $\Aline ab\cap\Aline cd\subseteq A$.
\end{definition}

\begin{exercise} Let $X$ be a liner. Show that any subset $A\subseteq X$ of cardinality $\le 3$ is closed in $X$. 
\end{exercise}

\begin{exercise} Show that the intersection of any family of closed subsets of a liner is closed in the liner.
\end{exercise}

\begin{definition} Let $A$ be a subset of a liner $X$. The smallest closed subset $\langle A\rangle$ of $X$ that contains $A$ is called the \defterm{closure} of $A$ in $X$. The closure $\langle A\rangle$ is endowed with the structure of a liner inherited from $X$.
\end{definition}

\begin{exercise} Analyse the structure of the closures of $4$-element sets in the real plane.
\end{exercise}

\begin{definition} A subset $D$ of a liner $X$ is called \defterm{dense} in $X$ if $\langle D\rangle=X$.
\end{definition}

\begin{Exercise} Prove that the $4$-element set $\{(0,0),(6,0),(0,6),(2,2)\}$ is dense in the rational plane $\IQ\times\IQ$.
\end{Exercise}

\begin{definition} The \defterm{density} $d(X)$ of a liner $X$ is the smallest cardinality of a dense subset of $X$, i.e., 
$$d(X)\defeq\min\{|A|:A\subseteq X\;\wedge\;\langle A\rangle=X\}.$$
\end{definition}

\begin{Exercise} Calculate the density of the rational plane $\IQ\times\IQ$.
\end{Exercise}
 
\section{Flat sets in liners}

\begin{definition}\label{d:affine} A set $A\subseteq X$ in a liner $(X,\Af)$ is called \index{flat}\index{set!flat}\defterm{flat} if $$\forall x,y\in A\;\;(\Aline xy\subseteq A).$$Flat subsets in liners are called \defterm{flats}. A flat $A$ in $X$ is \index{proper flat}\index{flat!proper}\defterm{proper} if $A\ne X$.
\end{definition} 

Theorem~\ref{t:Alines} implies the following proposition.

\begin{proposition}\label{p:lines-are-affine} Every line in a liner is a flat.
\end{proposition}


\begin{proposition}\label{p:cov-aff} Let $A$ be a flat in a liner $X$ and $\F$ be a finite nonempty family of flats in $X$ such that $A\subseteq\bigcup\F$. If $X$ is $(|\F|+1)$-long, then $A\subseteq F$ for some $F\in\F$.
\end{proposition}

\begin{proof} Since the family $\F$ is finite, we can replace $\F$ by a smaller family and assume that $A\subseteq\bigcup\F$ but $A\not\subseteq\bigcup\mathcal E$ for any proper subfamily $\mathcal E\subset\F$. If $|\F|=1$, then $A\subseteq\bigcup\F$ implies that $A\subseteq F$ for the unique set $F\in\F$. So, we assume that $|\F|>1$. Choose any set $E\in\F$ and consider the family $\mathcal E=\F\setminus\{E\}$. Since $A\not\subseteq\bigcup\mathcal E$, there exists a point $p\in A\setminus\bigcup\mathcal E$. It follows from $p\in \bigcup F\setminus \bigcup\mathcal E$ that $p\in E$. Since the family $\F$ is minimal and has cardinality $|\F|>1$, there exists a point $q\in A\setminus E$. The choice of the points $p,q$ ensures that $\{p,q\}\not\subseteq F$ for every $F\in\F$. Applying Theorem~\ref{t:Alines}, we conclude that $|\Aline pq\cap F|\le 1$ for every flat $F\in\F$ and hence $|\Aline pq|\le|\F|$, which implies that $X$ is not $(|\F|+1)$-long.
\end{proof}

\begin{definition} Let $A,B$ be two flats in a liner $X$. The flat $A$ is called a \index{hyperplane}\defterm{hyperplane} in $B$ if $A\subsetneq B$ and every flat $C$ with $A\subseteq C\subseteq B$ is equal to $A$ or $B$.
\end{definition}

\begin{exercise} Construct an example of a liner containing no hyperplanes.
\end{exercise}



\section{Flat hulls of sets in liners}

It is clear that the intersection of flats in a liner $X$ is a flat. Consequently, every subset $A\subseteq X$ of a liner $X$ is a  subset of the smallest flat that contains $A$. This smallest flat is called the \index{flat hull}\defterm{flat hull} of $A$ in $X$ and is denoted by \index[note]{$\overline A$}\defterm{$\overline A$}. It is equal to the intersection of all flats in $X$ that contain $A$. 

\begin{exercise} Show that for any sets $A\subseteq B$ in a liner $X$, we have $\langle A\rangle\subseteq\overline A\subseteq \overline B$.
\end{exercise}

\begin{proposition}\label{p:aff-finitary} For every set $A\subseteq X$ in a liner $X$ and every point $x\in\overline A$, there exists a finite set $F_x\subseteq A$ such that $x\in \overline {F_x}$.
\end{proposition}

\begin{proof} Let $[A]^{<\w}$ be the family of finite subsets of $A$. Observe that the set $\bigcup_{F\in[A]^{<\w}}\overline F$ contains $A$ and is flat in $X$. Indeed, given any distinct points $x,y\in \bigcup_{F\in[A]^{<\w}}\overline F$, we can find two finite sets $F_x,F_y\in[A]^{<\w}$ such that $x\in \overline{F_x}$ and $y\in\overline{F_y}$. Then for the finite set $F_x\cup F_y$ we have $\{x,y\}\subseteq \overline{F_x}\cup\overline{F_y}\subseteq \overline{F_x\cup F_y}$ and hence $\Aline xy\subseteq\overline{F_x\cup F_y}\subseteq\bigcup_{F\in[A]^{<\w}}\overline F$, witnessing that the set $\bigcup_{F\in[A]^{<\w}}\overline F$ is flat. Taking into account that $\overline A$ is the smallest flat containing $A$, we conclude that $\overline A\subseteq\bigcup_{F\in[A]^{<\w}}\overline F$. On the other hand, for every $F\in[A]^{<\w}$, the inclusion $F\subseteq A$ implies $\overline F\subseteq \overline A$ and hence $\bigcup_{F\in[A]^{<\w}}\overline F\subseteq \overline A$ and finally, $\overline A=\bigcup_{F\in[A]^{<\w}}\overline F$.
\end{proof}

The flat hull of any set in a liner can be described using the operation $\Aline ST$ assigning to every subsets $S,T\subseteq X$ of a liner $(X,\Af)$ the set
$$\Aline ST\defeq\bigcup_{s\in S}\bigcup_{t\in T}\Aline st,$$which is the union of lines connecting points of the sets $S$ and $T$.   

\begin{exercise} Prove that for every set $A\subseteq X$ in a liner $X$, its flat hull $\overline A$ is equal to the union $\bigcup_{n\in\w}A_n$ of an the increasing sequence $(A_n)_{n\in\w}$ of sets, defined by the recursive formula: $A_0=A$ and $A_{n+1}=\Aline {A_n}{A_n}$ for all $n\in\w$. 
\end{exercise}

\begin{exercise} Let $V$ be an $R$-module over a corps $R$. Show that for every set $A\subseteq V$ its flat hull coincides with the set of affine combinations 
$$\textstyle\bigcup_{n\in \w}\big\{\sum_{i\in n}r_i\cdot a_i:(r_i)_{i\in\w}\in R^n,\;(a_i)_{i\in n}\in A^n,\;\sum_{i\in n}r_i=1\big\}.$$
\end{exercise}

\begin{exercise} Let $\mathbb V$ be the projective space of an $R$-module over a corps $R$. Show that for every set $\A\subseteq \mathbb PV$ its flat hull in $\mathbb PV$ coincides with the set  
$$\textstyle\bigcup_{n\in \w}\big\{L\in\mathbb LV:L\subseteq \sum_{i\in n}L_i:(L_i)_{i\in n}\in \A^n\big\}.$$
\end{exercise}

\section{Flat relations and functions}

\begin{definition}\label{d:flat-relation} Let $X$ and $Y$ be two liners. A relation $R\subseteq X\times Y$ is defined to be \index{flat relation}\index{relation!flat}\defterm{flat} if for any flats $A\subseteq X$ and $B\subseteq Y$, the sets\index[note]{$R[A]$}
$$R[A]\defeq\{y:\exists a\in A\;\;(a,y)\in R\}\subseteq Y\quad\mbox{and}\quad R^{-1}[B]\defeq\{x:\exists b\in B\;\;(x,b)\in R\}\subseteq X$$are flat.
\end{definition}

Definition~\ref{d:flat-relation} implies that for any flat relation $R\subseteq X\times Y$ between liners, its domain $\dom[R]\defeq R^{-1}[Y]$ and range $\rng[R]\defeq R[X]$ are flat sets in the liners $X$ and $Y$, respectively.

Since functions are relations, the definition of a flat relation is also applicable to functions. We recall that a \index{function}\defterm{function} is a relation $F$ such that for every $x\in\dom[F]$ there exists a unique $y\in\rng[F]$ such that $(x,y)\in F$. This unique element $y$ is called the \defterm{value} of $F$ at $x$ and is denoted by $F(x)$. 

A function $F$ is called \index{function!injective}\defterm{injective} if the relation $F^{-1}\defeq\{(y,x):(x,y)\in F\}$ is a function.
Injective functions are called \index{injection}\defterm{injections}. Given a function $F$ and sets $X,Y$, we write $F:X\to Y$ and say that $F$ is a function from $X$ to $Y$ if $\dom[F]=X$ and $\rng[F]\subseteq Y$.  A function $F:X\to Y$ is \index{bijective function}\index{function!bijection}\defterm{bijective} if $F$ is injective and $\rng[F]=Y$. Bijective functions are called \defterm{bijections}.

Flat bijections between subsets of liners admit the following simple characterization. 

\begin{proposition}\label{p:flat-injection} Let $X$ and $Y$ two liners. An injective function $F\subseteq X\times Y$ is flat if and only if its domain $\dom[F]$ is a flat in $X$ and for every line $L\subseteq \dom[F]$, its image $F[L]$ is a line in $Y$.
\end{proposition}

\begin{proof} If $F$ is flat, then the set $\dom[F]=F^{-1}[X]$ is flat by Definition~\ref{d:flat-relation}. Moreover, for every line $L\subseteq\dom[F]$, the set $F[L]$ is a flat. The injectivity of $F$ ensures that $|F[L]|=|L|\ge 2$ and hence $F[L]$ contains at least two distinct points $a,b$. Since $F[L]$ is flat, $\Aline ab\subseteq F[L]$.  To show that $F[L]=\Aline ab$, take any point $c\in F[L]$ and find a point $z\in L$ such that $F(z)=c$. Since $F$ is flat, the preimage $F^{-1}[\Aline ab]$ is a flat in $X$, containing the distinct points $x,y\in L$. Then $z\in L=\Aline xy\subseteq F^{-1}[\Aline ab]$ and $c=F(z)\in F[F^{-1}[\Aline ab]]\subseteq\Aline ab$, witnessing that $F[L]$ coincides with the line $\Aline ab$. 
\smallskip

To prove the ``if'' part, assume that the set $\dom[F]$ is a flat in $X$ and for every line $L\subseteq \dom[F]$, the set $F[L]$ is a line in $Y$. Given flats $A\subseteq X$ and $B\subseteq Y$, we should prove that the sets $F[A]$ and $F^{-1}[B]$ are flat in $Y$ and $X$, respectively. To see that the set $F[A]$ is flat in $Y$, take any distinct points $a,b\in F[A]$ and choose any (necessarily distinct) points $x=A\cap F^{-1}(a)$ and $y=A\cap F^{-1}(b)$. By our assumption, $F[\Aline xy]$ is a line containing the points $a,b$ and hence $F[\Aline xy]=\Aline ab$. Then $\Aline ab=F[\Aline xy]\subseteq F[A]$, witnessing that the set $F[A]$ is flat.

To shows that the set $F^{-1}[B]$ is flat, take any distinct points $x,y\in F^{-1}[B]\subseteq\dom[F]$. By the injectivity of $F$, the points $a=F(x)$ and $b\in F(y)$ are distinct. Then $\Aline ab$ is a line in the flat $B$. Taking into account that $\dom[F]$ is flat in $X$, we conclude that $\Aline xy\subseteq\dom[F]$. By our assumption, $F[\Aline xy]$ is a line containing the points $a,b\in B$ and hence $F[\Aline xy]=\Aline ab\subseteq B$ and $\Aline xy\subseteq F^{-1}[B]$, witnessing that the set $F^{-1}[B]$ is flat.
\end{proof}

\begin{corollary}\label{c:line-bijections} Let $X$ and $Y$ be two liners. Any bijection $F:L_1\to L_2$ between two lines $L_1\subseteq X$ and $L_2\subseteq Y$ is a flat function.
\end{corollary}

A function $F:X\to Y$ between the underlying sets of two liners $(X,\Af_X)$ and $(Y,\Af_Y)$ is called 
\begin{itemize}
\item a \index{liner morphism}\index{liner!morphism of}\defterm{liner morphism} if $\{Fxyz:xyz\in \Af_X\}\subseteq\Af_Y$;
\item a \index{liner isomorphism}\index{liner!isomorphism}\defterm{liner isomorphism} if $F$ is a bijective map such that $F$ and $F^{-1}$ are liner morphisms;
\item a \index{liner embedding}\index{liner!embedding}\defterm{liner embedding} if $F$ is an liner isomorphism of $X$ into the subliner $F[X]$ of $Y$.
\end{itemize}

%

Proposition~\ref{p:flat-injection} implies the following characterization of liner isomorphisms.

\begin{theorem}\label{t:liner-isomorphism<=>} For a bijective function $F:X\to Y$ between two liners $X,Y$, the following conditions are equivalent:
\begin{enumerate}
\item $F$ is flat;
\item $F$ is a liner isomorphism;
\item for every line $L$ in $X$, the set $F[L]$ is a line in $Y$;
\item for every line $L$ in $Y$, the set $F^{-1}[L]$ is a line in $X$.
\end{enumerate}
\end{theorem}

\begin{proposition} Let $F:X\to Y$ be an isomorphism between two liners. For every set $A\subseteq X$, we have $F[\overline{A}]=\overline{F[A]}$.
\end{proposition}

\begin{proof} Since the function $F$ is flat and the set $\overline A$ is flat in $X$, the set $F[\overline A]$ is flat in $Y$ and hence $\overline{F[A]}\subseteq  F[\overline A]$ (because $F[A]\subseteq F[\overline A]$ and $\overline{F[A]}$ is the smallest flat set containing $F[A]$). On the other hand, the flatness of the function $F$ and the flatness of the set $\overline{F[A]}$ in $X$ imply that the flatness of the set $F^{-1}[\overline{F[A]}]$ in $X$. Taking into account that $A=F^{-1}[F[A]]\subseteq F^{-1}[\overline{F[A]}]$, we conclude that $\overline A\subseteq F^{-1}[\overline{F[A]}]$ and hence $F[\overline A]\subseteq F[F^{-1}[\overline{F[A]}]]=\overline{F[A]}$. Therefore, $F[\overline A]=\overline{F[A]}$.
\end{proof}

Now we introduce some algebraic terminology related to semigroups of transformations. We recall that a \index{semigroup}\defterm{semigroup} is a set  $X$ endowed with an associative binary operation $X\times X\to X$, $(x,y)\mapsto xy$. A semigroup $X$ is a \index{monoid}\defterm{monoid} if there exists a (necessarily unique) element $e\in X$ such that $\forall x\;(xe=x=ex)$. This unique element $e$ is called the \index{identity}\defterm{identity element} of the monoid $X$. A subset $A\subseteq X$ of a monoid $X$ is called a \index{submonoid}\defterm{submonoid} if $e\in A$ and $AA=\{xy:x,y\in A\}\subseteq X$. A semigroup $X$ is called \index{inverse semigroup}\index{semigroup!inverse}\defterm{inverse} if for every element $x\in X$ there exists a unique element $x^{-1}\in X$ such that $xx^{-1}x=x$ and $x^{-1}xx^{-1}=x^{-1}$.

For every set $X$, the set \index[note]{$\mathcal R_X$}$\mathcal R_X$ of all relations $R\subseteq X\times X$ is a monoid with respect to the operation of composition of relations: for two relations $F,G$, their composition is the relation 
$$GF\defeq\big\{(x,z):\exists y\;\big((x,y)\in F\;\wedge\;(y,z)\in G\big)\big\}.$$
The identity element of the monoid $\mathcal R_X$ is the identity function $1_X:X\to X$.

Since the composition of two functions is a function, the set \index[note]{$\mathcal F_X$}$\mathcal F_X$ of functions $F$ with $\dom[F]\cup\rng[F]\subseteq X$ is a submonoid of $\mathcal R_X$. The monoid $\mathcal F_X$ contains the submonoid $\mathcal I_X$ consisting of injective functions $F$ with $\dom[F]\cup\rng[F]\subseteq X$. We recall that a function $F$ is \index{function!injective}\defterm{injective} if the relation $F^{-1}\defeq\{(y,x):(x,y)\in F\}$ is a function. The monoid \index[note]{$\mathcal I_X$}$\mathcal I_X$ is called \index{symmetric inverse monoid}\defterm{the symmetric inverse monoid} of the set $X$. 
The set \index[note]{$\Sym(X)$}$\Sym(X)\defeq\{F\in\mathcal I_X:\dom[F]=X=\rng[F]\}$ is a subgroup of $\mathcal I_X$ called the \index{symmetric group}\defterm{symmetric group} on the set $X$. 

The symmetric group $\Sym(X)$  contains an important subgroup $\Alt(X)$ that consists of all permutations of $X$ that are products of even number of \index{transposition}\defterm{transpositions} (i.e., bijections $F:X\to X$ with two-element \defterm{support} $\supp(F)\defeq\{x\in X:F(x)\ne x\})$. The group $\Alt(X)$ is called the \index{alternating group}\defterm{alternating group} on the set $X$ and its elements are called \defterm{even permutations} of the set $X$. 

The symmetric and alternating groups $\Sym(n)$ and $\Alt(n)$ on a number $n\defeq\{0,1,\dots,n-1\}$ will be denoted by $S_n$ and $A_n$, respectively. 

Since the composition of two flat relations $F,G\subseteq X\times X$ on a liner $X$ is a flat relation, the set \index[note]{$\mathcal{FR}_X$}$\mathcal{FR}_X$ of flat relations on $X$ is a submonoid of the monoid $\mathcal R_X$. The set \index[note]{$\mathcal {FF}_X$}$\mathcal{FF}_X\defeq\mathcal{FR}_X\cap\mathcal F_X$ of flat functions between flat subsets of the liner  $X$ is a submonoid of the monoid $\mathcal F_X$, and the set $\mathcal{FI}_X\defeq \mathcal{FF}_X\cap\mathcal I_X$ of flat injections is an inverse submonoid of the symmetric inverse monoid $\mathcal I_X$. Elements of the monoid \index[note]{$\mathcal{FI}_X$}$\mathcal{FI}_X$ are flat bijections between flat subsets of the liner $X$. By Corollary~\ref{c:line-bijections}, every bijection between lines in a liner $X$ is an element of the monoid $\mathcal{FI}_X$. Flat bijections between liners are called \index{collineation}\defterm{collineation}. By Proposition~\ref{p:flat-injection}, a bijective function $F:X\to Y$ between liners $X,Y$ is a collineation if and only if for every lines $L\subseteq X$ and $L'\subseteq Y$, the sets $F[L]$ and $F^{-1}[L']$ are lines in the liners $Y$ and $X$, respectively.

\begin{definition} A \index{line bijection}\defterm{line bijection} in a liner $X$ is an injective function $F\subseteq X\times X$ whose domain and range $\dom[F]$ and $\rng[F]$ are lines in $X$.
\end{definition} 

So, line bijections in a liner $X$ are just bijective functions between lines in $X$. By Proposition~\ref{c:line-bijections}, every line bijection is a flat function.

For two line bijections $F, G$ in a liner $X$, the composition $GF$ is a line bijection if and only if $\rng[F]=\dom[G]$ if and only if $|\rng[F]\cap\dom[G]|\ge 2$. In the opposite case, $|\rng[F]\cap\dom[G]|\le 1$, which means that $GF$ is a trivial bijection in $X$.

\begin{definition} A \index{trivial bijection}\index{bijection!trivial}\defterm{trivial bijection} in a set $X$ is a bijective function $F:A\to B$ between subsets $A,B\subseteq X$ of cardinality $|A|=|B|\le 1$.
\end{definition}

The set of trivial bijections in a set $X$ is an ideal in the inverse symmetric monoid $\mathcal{I}_X$.

Let \index[note]{$\mathcal I_X^\ell$}$\mathcal I_X^\ell$ be the smallest submonoid of $\I_X$, containing all line bijections. Every element of $\I_X^\ell$ is either a line bijection or a trivial bijection or the identity map of $X$. 

\begin{exercise}  Show that $\mathcal{I}_X^{\ell}$ is an inverse submonoid of the inverse monoid $\mathcal{I}^\ell_X$.
\smallskip

\noindent{\em Hint:} Use the equality $(FG)^{-1}=G^{-1} F^{-1}$ holding for any relations  $F,G$.\end{exercise}

\begin{exercise} Show that the line relation $\Af$ of a liner $X$ can be recovered from the algebraic structure of the semigroup $\I_X^\ell$.
\smallskip

\noindent{\em Hint:} Look at the semilattice of idempotents $\{E\in\I_X^\ell:EE=E\}$ of $\I_X^\ell$ and derive the information about points and lines from the algebraic structure of this semilattice.
\end{exercise}

An \index{automorphism}\index{liner!automorphism of}\defterm{automorphism} of a liner $X$ is any isomorphism $F:X\to X$. Automorphisms of a liner form a subgroup $\Aut(X)$ in the inverse monoid $\mathcal{FI}_X$ of flat injections of $X$. Moreover, \index[note]{$\Aut(X)$}$\Aut(X)$ is the maximal subgroup of $\mathcal{FI}_X$ containing the identity map $1_X:X\to X$ of $X$.

\section{The category of liners}


Let us recall \cite{BanakhCST} that a {\em category}\index{category} is a $6$-tuple $\C=(\Ob,\Mor,\fdom,\fran,\mathsf 1,\circ)$ consisting of
\begin{itemize}
\item a class $\Ob$ whose elements are called {\em objects}\index{object of category}\index{category!object of} of the category $\C$ (briefly, {\em $\C$-objects});
\item a class $\Mor$ whose elements are called {\em morphisms}\index{morphism of category}\index{category!morphism of} of the category $\C$ (briefly, {\em $\C$-morphisms});
\item two functions $\fdom\colon\Mor\to\Ob$ and $\fran\colon\Mor\to\Ob$ assigning to each morphism $f\in\Mor$ its \index{source $\fdom$}{\em source} $\fdom(f)\in\Ob$ and \index{target $\fran$}{\em target} $\fran(f)\in\Ob$;
\item a function $\mathsf 1:\Ob\to\Mor$ assigning to each object $X\in\Ob$ a morphism $\mathsf 1_X\in\Mor$, called the \index{identity morphism} {\em identity morphism} of $X$, and satisfying the equality $\fdom(\mathsf 1_X)=X=\fran(\mathsf 1_X)$;
\item a function $\circ$ with domain $\dom[\circ]=\{\langle f,g\rangle\in \Mor\times\Mor:\fran(g)=\fdom(f)\}$ and range $\rng[\circ]\subseteq\Mor$ assigning to any $\langle f,g\rangle\in\dom[\circ]$ a morphism $f\circ g\in\Mor$ such that $\fdom(f\circ g)=\fdom(g)$ and $\fran(f\circ g)=\fran(f)$ and the following axioms are satisfied:
\begin{itemize}
\item [(A)] for any morphisms $f,g,h\in\Mor$ with $\fdom(f)=\fran(g)$ and $\fdom(g)=\fran(h)$ we have\\ $(f\circ g)\circ h=f\circ(g\circ h)$;
\item [(U)] for any morphism $f\in\Mor$ we have $\mathsf 1_{\fran(f)}\circ f=f=f\circ\mathsf 1_{\fdom(f)}$.
\end{itemize}
\end{itemize}

The function $\circ$ is called the \index{composition of morphisms}{\em operation of composition} of morphisms and the axiom $(A)$ is called the {\em associativity} of the composition. For two objects $X,Y\in\Ob$ of the category $\C$ the set $\Mor(X,Y)\defeq\{f\in\Mor:\fdom(f)=X\;\wedge\;\rng(f)=Y\}$ is called the set of morphisms from $X$ to $Y$. The set $\Mor(X,X)$ carries the structure of a monoid with identity element $1_X$.

More information on Category Theory can be found in the books \cite{MacLane} and \cite{BanakhCST}.

Liners and their isomorphisms form a category denoted by $\Liners$. Objects of this category are liners and morphisms are liner isomorphisms $\Phi:X\to Y$  (which are triples $(X,\Phi,Y)$ consisting of liners $X,Y$ and an isomorphism $\Phi:X\to Y$). The functions $\fdom$ and $\fran$ assign to each liner isomorphism $\Phi:X\to Y$ its domain $X$ and range $Y$, respectively. The map $\mathsf 1$ assigns to each liner $X$ its identity isomorphism $\mathsf 1_X:X\to X$. The composition map $\circ$ is the standard composition of liner isomorphisms $\Phi:X\to Y$ and $\Psi:Y\to Z$. 

For two liners $X,Y$ the set $\Mor(X,Y)$ of morphisms from  $X$ to $Y$ in the category $\Liners$ coincides with the set of isomorphisms between $X$ and $Y$. The set $\Mor(X,X)$ carries the structure of a group and is denoted by $\Aut(X)$, the group of automorphisms of the liner $X$.

\begin{definition} A liner $X$ is defined to be \index{rigid liner}\index{liner!rigid}\defterm{rigid} if its automorphism group $\Aut(X)$ is trivial.
\end{definition}
\chapter{Exchange and ranks}

\section{Exchange Properties in liners}

\begin{definition} A liner $(X,\Af)$ is defined to have the \index{Exchange Property}\index{property!Exchange}\defterm{Exchange Property} if for every flat $A\subseteq X$ and points $x\in X\setminus A$, $y\in \overline{A\cup\{x\}}\setminus  A$ we have $x\in\overline{A\cup\{y\}}$. 
\end{definition}

\begin{remark} The Exchange Property first appeared in the Steinitz Exchange Lemma, used by  \index[person]{Steinitz}Steinitz\footnote{{\bf Ernst Steinitz} (1871 -- 1928) was a German mathematician. In 1910 Steinitz published the very influential paper ``{\em Algebraische Theorie der K\"orper}'' laying fundamentals of the theory of fields. Bourbaki called this article ``a basic paper which may be considered as having given rise to the current conception of Algebra''. Also he is the author of the famous Steinitz Exchange Lemma, which allows to define dimension of vector spaces.} for introducing the dimension of vector spaces.
\end{remark}

\begin{exercise} Let $X$ be a liner with the Exchange Property. Prove that for every set $A\subseteq X$ and points $x\in X$ and $y\in\overline{A\cup\{x\}}\setminus\overline{A}$ we have $\overline{A\cup\{x\}}=\overline{A\cup\{y\}}$.
\end{exercise}


The Exchange Property can be quantified as follows.

\begin{definition}\label{d:k-EP} A liner $X$ has a \index{$\kappa$-Exchange Property}\index{property!$\kappa$-Exchange}\defterm{$\kappa$-Exchange Property} for a cardinal $\kappa$ if for any set $A\subseteq X$ of cardinality $|A|<\kappa$ and any points $x\in X$ and $y\in\overline{A\cup\{x\}}\setminus\overline A$, we have $x\in\overline{A\cup \{y\}}$. 
\end{definition}

It is clear that a liner has the Exchange Property if and only if it has the $\kappa$-Exchange Property for every cardinal $\kappa$. 

\begin{exercise} Show that every liner has the $2$-Exchange Property.
\end{exercise}

\begin{proposition}\label{p:fh-EP<=>wEP} A liner has the Exchange Property if and only if it has the $n$-Exchange Property for every $n\in\w$.
\end{proposition}

\begin{proof} The ``only if'' part is trivial. To prove the ``if'' part, assume that a liner $X$ has the $n$-Exchange Property for every $n\in\w$. To prove the Exchange Property for $X$, fix any set $A\subseteq X$ and points $x\in X\setminus\overline A$ and $y\in\overline{A\cup\{x\}}\setminus\overline A$. By Proposition~\ref{p:aff-finitary}, there exists a finite set $F\subseteq A$ such that $y\in\overline{F\cup\{x\}}$ and hence $y\in\overline{F\cup\{x\}}\setminus\overline{A}\subseteq\overline{F\cup\{x\}}\setminus\overline F$. Since $X$ has the $(|F|+1)$-Exchange Property, $y\in\overline{F\cup\{x\}}\setminus\overline F$ implies $x\in\overline{F\cup\{y\}}\subseteq\overline{A\cup\{y\}}$.
\end{proof}

\begin{exercise} Let $F:X\to Y$ be a liner isomorpism between two liners $X,Y$, and $\kappa$ be a cardinal. Prove that the liner $X$ has the $\kappa$-Exchange Property if and only if $Y$ has the $\kappa$-Exchange Property.
\end{exercise}

\begin{example}[Terence Tao\footnote{See ({\tt mathoverflow.net/a/450472/61536}).}]\label{ex:Tao}\index[person]{Tao} There exists a liner $(X,\Af)$ that does not have the $3$-Exchange Property.
\end{example}

\begin{proof} Let $\IF_3=\{-1,0,1\}$ be a 3-element field, $n,m$ be positive integers and $f:\IF_3^n\to \IF_3^m$ be an even function with $f(\{0\}^n)=\{0\}^m$.  On the set $X=\IF_3^n\times \IF_3^m$, consider the ternary relation
$$
\begin{aligned}
\Af&=\{((x,x'),(y,y'),(z,z'))\in X^3:(y,y')\in\{(x,x'),(z,z')\}\}\\
&\cup\{((x,x'),(y,y'),(z,z'))\in [X]^3:x+y+z=0\;\wedge\;x'+y'+z'=f(x-z)\}.
\end{aligned}
$$
Here $[X]^3=\{(x,y,z)\in X^3:|\{x,y,z\}|=3\}$.
We claim $\Af$ is a line relation on $X$. The definition of $\Af$ ensures that it satisfies the axioms {\sf(IL)} and {\sf(RL)}. It remains to check the Exchange Axiom {\sf(EL)}. Take any points $(a,a'),(b,b'),(x,x'),(y,y')\in X=\IF_3^n\times \IF_3^m$ with $((a,a'),(x,x'),(b,b'))\in\Af$, $((a,a'),(y,y'),(b,b'))\in\Af$ and $(x,x')\ne (y,y')$.
Then $(a,a')\ne(b,b')$. We should prove that $\{((x,x'),(a,a'),(y,y')),(x,x'),(b,b'),(y,y'))\}\subseteq\Af$. This follows from the axiom {\sf(RL)} if $\{(x,x'),(y,y')\}=\{(a,a'),(b,b')\}$. So, we assume that $\{(x,x'),(y,y')\}\ne\{(a,a'),(b,b')\}$.

We claim that $\{(x,x'),(y,y')\}\cap\{(a,a'),(b,b')\}\ne\varnothing$. In the opposite case, the definition of the relation $\Af$ ensures that $a+x+b=0=a+y+b$ and $a'+x'+b'=f(a-b)=a'+y'+b'$, which implies $(x,x')=(y,y')$ and contradicts the choice of the pairs $(x,x')\ne(y,y')$.
This contradiction shows that $\{(x,x'),(y,y')\}\cap\{(a,a'),(b,b')\}\ne\varnothing$. Without loss of generality we can assume that $(x,x')=(a,a')$ and $(y,y')\ne(b,b')$. Then $((x,x'),(a,a'),(y,y'))\in\Af$ by the axiom {\sf(RL)}. Since $((x,x'),(y,y'),(b,b'))=((a,a'),(y,y'),(b,b'))\in[X]^3$, we have $x+y+b=0$ and $x'+y'+b'=f(x-b)$. Observe that $y-x=y+x+x=-b+x$ and hence $x'+b'+y'=x'+y'+b'=f(x-b)=f(y-x)=f(x-y)$, witnessing that $((x,x'),(b,b'),(y,y'))\in\Af$.

Therefore, $(X,\Af)$ is a liner. Now we select a special even function $f$ for which this liner does not have the $3$-Exchange Property. For $n=2$ and $m=1$, consider the even function $f:\IF_3^2\to \{0,1\}\subseteq \IF_3$ such that $f^{-1}(1)=\{(0,1),(0,-1)\}$. For such a function $f$, the set $L=\{(x,0,0):x\in \IF_3\}$ is a line in $(X,\Af)$ such that for the point $\vecbold{u}=(0,1,0)$ the flat hull $\overline{L\cup\{\vecbold{u}\}}$ coincides with $X$. On the other hand, for the point $\vecbold{v}=(0,0,1)\in\overline{L\cup\{\vecbold{u}\}}\setminus L$ the flat hull $\overline{L\cup\{\vecbold{v}\}}$ coincides with the set $\IF_3\times\{0\}\times \F_3$, which is strictly smaller than $\overline{L\cup\{\vecbold{u}\}}=X$. This means that the liner $(X,\Af)$ does not have the $3$-Exchange Property.

Observe that $L_1\defeq\{(-1,1,0),(0,1,0),(1,1,0)\}$ and $L_2\defeq\{(0,-1,-1),(0,0,-1),(0,1,0)\}$ are two distinct lines in the plane $X=\overline{L\cup\{\vecbold{u}\}}$ that contain the point $\vecbold{u}$ and are disjoint with the line $L$. Also for the lines $L_2$ and $L_3\defeq\{(-1,-1,-1),(-1,-1,0),(-1,-1,1)\}$ differ by the following properties: $\forall x\in L_3 \;(L_3\subseteq \overline{L_2\cup\{x\}})\quad\mbox{but}\quad \forall x\in L_2\;(L_2\not\subseteq \overline{L_3\cup\{x\}})$.



\end{proof}

\section{The rank and dimension in liners}

\begin{definition}\label{d:aff-rank} The \index[note]{$\|X\|$}\index{rank}\defterm{rank} $\|A\|$ of a set $A$ in a liner $X$ is the smallest cardinality of a set $B\subseteq X$ such that $A\subseteq\overline B$.
\end{definition}

For a cardinal $\kappa$ let $\kappa-1$ be a unique cardinal or integer number such that $\kappa=1+(\kappa-1)$. The definition of $\kappa-1$ implies that $\kappa-1=\kappa$ for every infinite cardinal $\kappa$. 

\begin{definition}\label{d:aff-dimension} For a subset $A$ of a liner $X$, the  number
$$\dim(A)\defeq\|A\|-1$$
is called the \index[note]{$\dim(X)$}\defterm{dimension}\index{dimension} of $A$. 
\end{definition}

According to Definition~\ref{d:aff-dimension},  singletons in liners have dimension zero and lines have dimension 1. The empty set has dimension $\dim(\varnothing)=-1$.

\begin{proposition}\label{p:rank-EP} If a liner $X$ has the $\kappa$-Exchange Property for some cardinal $\kappa$, then for every subset $A\subseteq X$ of rank $\|A\|\le \kappa$, there exists a set $A'\subseteq A$ of cardinality $|A'|=\|A\|$ such that $A\subseteq\overline {A'}$.
\end{proposition}

\begin{proof} We divide the proof into two claims. 

\begin{claim}\label{cl:rank-EP} For every subset $A\subseteq X$ of finite rank $\|A\|\le\kappa$, there exists a set $B\subseteq A$ of cardinality $|B|=\|A\|$ such that $A\subseteq\overline{B}$.
\end{claim}

\begin{proof} By the definion of the rank $\|A\|$, there exists a set $B\subseteq X$ of cardinality $|B|=\|A\|$ such that $A\subseteq\overline B$. We can assume that the set $B$ is chosen so that the cardinal $|B\setminus A|$ is the smallest possible.  Assuming that $B\not\subseteq A$, we can fix a point $b\in B\setminus A$. Since the cardinal $\|A\|=|B|$ is finite, the set $B_{\mbox{-}b}\defeq B\setminus\{b\}$ has cardinality $|B_{\mbox{-}b}|<|B|=\|A\|$ and hence $A\not\subseteq\overline{B_{\mbox{-}b}}$. Fix any point $a\in A\setminus\overline{B_{\mbox{-}b}}$. 
 Since $a\in A\subseteq\overline{B}=\overline{B_{\mbox{-}b}\cup\{b\}}$ and $|B_{\mbox{-}b}|<|B|=\|A\|\le\kappa$, the $\kappa$-Exchange Property of $X$ ensures that $b\in \overline{B_{\mbox{-}b}\cup\{a\}}$. Then the set $B'\defeq B_{\mbox{-}b}\cup\{a\}$ has the following properties: $|B'|=|B|=\|A\|$, $A\subseteq \overline{B_{\mbox{-}b}\cup\{b\}}\subseteq\overline{B_{\mbox{-}b}\cup\overline{B_{\mbox{-}b}\cup \{a\}}}=\overline{B_{\mbox{-}b}\cup\{a\}}=\overline{B'}$ and $|B'\setminus A|=|(B\setminus A)\setminus\{b\})|<|B\setminus A|$, which contradicts the  minimality of the cardinal $|B\setminus A|$. This contradiction shows that $B\subseteq A$. 
\end{proof}  

\begin{claim} For every subset $A\subseteq X$ of infinite rank $\|A\|\le\kappa$, there exists a set $B\subseteq A$ of cardinality $|B|=\|A\|$ such that $A\subseteq\overline{B}$.
\end{claim}

\begin{proof} By the definition of the rank $\|A\|$, there exists a set $C\subseteq X$ of (infinite) cardinality $|C|=\|A\|$ such that $A\subseteq \overline{C}$. Consider the family $[C]^{<\w}$ of finite subsets of $C$ and observe that $|[C]^{<\w}|=|C|=\|A\|$. For every finite set $F\in[C]^{<\w}$, consider the set $A_F\defeq A\cap\overline{F}$ and observe that the set $A_F$ has finite rank $|A_F|_\Hull\le |F|$. By Claim~\ref{cl:rank-EP}, there exists a set $B_F\subseteq A_F$ of cardinality $|B_F|=\|A_F\|\le|F|$ such that $A_F\subseteq\overline{B_F}$. We claim that the set $A'\defeq\bigcup_{F\in[C]^{<\w}}B_F\subseteq A$ has the required properties.
Indeed,
$$
A=A\cap\overline C=A\cap\bigcup_{F\in[C]^{<\w}}\overline F=\bigcup_{F\in[C]^{<\w}}A_F\subseteq\bigcup_{F\in[C]^{<\w}}\overline{B_F}\subseteq\overline{A'}$$
and
$$|A'|\le\sum_{F\in[C]^{<\w}}|B_F|\le\sum_{F\in[C]^{<\w}}|F|\le |[C]^{<\w}|\cdot \w=|C|=\|A\|.$$
Since $A'\subseteq A$ and $A\subseteq\overline{A'}$, the definition of the rank $\|A\|$ ensures that $\|A\|\le|A'|$ and hence $|A'|=\|A\|$. 
\end{proof}
\end{proof}

\begin{exercise}[Yurii Yarosh] Construct an $\w$-long liner $X$ containing a subset $A\subseteq X$ of rank $\|A\|=3$ such that $A\not\subseteq B$ for every finite subset $B\subseteq A$.
\smallskip

\noindent{\em Hint:} Consider a suitable subliner of the Euclidean plane.
\end{exercise}

Proposition~\ref{p:rank-EP} has the following self-improvement.

\begin{proposition}\label{p:rank-EP2} If a liner $X$ has the $\kappa$-Exchange Property for some cardinal $\kappa$, then for every subsets $A,B\subseteq X$ with $\|A\cup B\|\le \kappa$, there exist sets $A'\subseteq A$ and $B'\subseteq B\setminus\overline A$ such that $A\subseteq \overline {A'}$,
 $B\subseteq \overline{A'\cup B'}$ and $|A'\cup B'|=\|A\cup B\|$.
\end{proposition}

\begin{proof}  We divide the proof into two cases, considered in the following two claims.

\begin{claim}\label{cl:rank-EP2} For every subsets $A,B\subseteq X$ with finite rank $\|A\cup B\|\le\kappa$, there exist sets $A'\subseteq A$ and $B'\subseteq B\setminus\overline A$  such that $A\subseteq \overline {A'}$,
 $B\subseteq \overline{A'\cup B'}$ and $|A'\cup B'|=\|A\cup B\|$.
\end{claim}

\begin{proof} By Proposition~\ref{p:rank-EP}, there exists a set $C\subseteq A\cup B$ of finite cardinality $|C|=\|A\cup B\|$ such that $A\cup B\subseteq\overline C$. We can assume that the cardinality $A\cap C$ is the largest possible. In this case we will show that $A\subseteq \overline{A\cap C}$. Assuming that $A\not\subseteq \overline{A\cap C}$, choose a point $a\in A \setminus\overline{A\cap C}$. Since $a\in A\subseteq \overline C$, we can choose a subset $M\subseteq C$ of the smallest possible cardinality such that $a\in \overline M$ and $A\cap C\subseteq M$.
The choice of the point $a\notin\overline{A\cap C}$ ensures that $A\cap C\ne M$. So, we can find a point $x\in M\setminus (A\cap C)=M\setminus A$. Assuming that $x\in \overline{A\cap C}$, we conclude that $a\in \overline{M}= \overline{(M\setminus \{x\})\cup \{x\}\cup (A\cap C)}=\overline{(M\setminus\{x\})\cup (A\cap C)}=\overline{M\setminus\{x\}}$, which contradicts the minimality of $M$. Therefore, $x\notin\overline{A\cap C}$. The minimality of $M$ ensures that $a\notin\overline{M\setminus\{x\}}$. Since $a\in\overline{(M\setminus\{x\})\cup\{x\}}\setminus\overline{M\setminus\{x\}}$  and $|M\setminus\{x\}|<|M|\le C\le\|A\cup B\|\le\kappa$, the $\kappa$-Exchange Property of $X$ ensures that $x\in\overline{(M\setminus\{x\}) \cup\{a\}}$. Now concider the set $C'=(C\setminus\{x\})\cup\{a\}$ and observe that $|C'|=|C|=\|A\cup B\|$ and 
$$
\begin{aligned}
A\cup B&\subseteq \overline C=\overline{(C\setminus M)\cup(M\setminus\{x\})\cup\{x\}}\subseteq\overline{(C\setminus M)\cup(M\setminus\{x\})\cup\overline{(M\setminus\{x\})\cup\{a\})}}\\
&=\overline{(C\setminus M)\cup(M\setminus\{x\})\cup\{a\}}=\overline{(C\setminus\{x\})\cup\{a\}}=\overline{C'}.
\end{aligned}
$$
Since $|A\cap C'|=|(A\cap C)\cup\{a\}|=|A\cap C|+1$, we obtain a contradiction with the maximality of $|A\cap C|$. This contradiction shows that $A\subseteq\overline{A\cap C}$. Then the sets $A'\defeq A\cap C$ and $B'=C\setminus\overline A\subseteq B\setminus\overline A$ have the required properties.
\end{proof}

\begin{claim}\label{cl:rank-EP3} For every sets $A,B\subseteq X$ with infinite rank $\|A\cup B\|\le\kappa$, there exist sets $A'\subseteq A$ and $B'\subseteq B\setminus\overline A$ such that $A\subseteq \overline {A'}$,
 $B\subseteq \overline{A'\cup B'}$ and $|A'\cup B'|=\|A\cup B\|$.
\end{claim}

\begin{proof}  By Proposition~\ref{p:rank-EP}, there exists a set $D\subseteq A\cup B$ of (infinite) cardinality $|D|=\|A\cup B\|$ such that $A\cup B\subseteq \overline{D}$. Consider the family $[D]^{<\w}$ of finite subsets of $D$ and observe that $|[D]^{<\w}|=|D|=\|A\|$. For every finite set $F\in[D]^{<\w}$, consider the sets  $A_F\defeq A\cap\overline{F}$ and $B_F\defeq B\cap\overline F$ and observe that  $\|A_F\cup B_F\|\le |F|<\omega$. By Claim~\ref{cl:rank-EP2}, there exist sets $A_F'\subseteq A_F$ and $B'_F\subseteq B_F$ such that $A_F\subseteq\overline{A'_F}$, $B_F\subseteq \overline{A'_F\cup B'_F}$, and $|A'_F\cup B'_F|=\|A_F\cup B_F\|$. Then the sets  $A'\defeq\bigcup_{F\in[D]^{<\w}}A'_F$ and $B'\defeq \bigcup_{F\in[D]^{<\w}}B_F'$ have the required  properties:
$$
\begin{aligned}
&A=A\cap\overline D=A\cap\bigcup_{F\in[D]^{<\w}}\overline F=\bigcup_{F\in[D]^{<\w}}A_F\subseteq\bigcup_{F\in[D]^{<\w}}\overline{A'_F}\subseteq \overline{A'},\\
&B=B\cap\overline D=B\cap\bigcup_{F\in[D]^{<\w}}\overline F=\bigcup_{F\in[D]^{<\w}}B_F\subseteq\bigcup_{F\in[D]^{<\w}}\overline{A'_F\cup B'_F}\subseteq\overline{A'\cup B'},\\
&|A'\cup B'|\le\sum_{F\in[D]^{<\w}}|A'_F\cup B'_F|\le\sum_{F\in[D]^{<\w}}|F|\le |[D]^{<\w}|\cdot \w=|D|=\|A\cup B\|.
\end{aligned}
$$
Since $A\cup B\subseteq\overline{A'\cup B'}$, the definition of the rank $\|A\cup B\|$ ensures that $\|A\cup B\|\le|A'\cup B'|$ and hence $|A'\cup B'|=\|A\cup B\|$.
\end{proof}
\end{proof}

\begin{proposition}\label{p:line=dim1} A subset $L$ of a liner $(X,\Af)$ is a line if and only if $L$ is a flat of rank $\|L\|=2$.
\end{proposition}

\begin{proof} If $L$ is a line, then $L=\Aline xy$ for two distinct points $x,y\in L$, witnessing that $\|L\|\le 2$. Assuming that $\|L\|<2$, we can find a point $z\in X$ such that $L\subseteq\overline{\{z\}}=\{z\}$ and conclude that $x=z=y$, which contradicts the choice of the points $x,y$. This contradiction shows that $\|L\|=2$. By Proposition~\ref{p:lines-are-affine}, the line $L$ is a flat.

Now assume that $L$ is a flat of rank $\|L\|=2$. Then $L\subseteq\overline{\{x,y\}}$ for some distinct points $x,y\in X$. By Proposition~\ref{p:lines-are-affine}, $\overline{\{x,y\}}=\Aline xy$. Since $\|L\|=2>1$, the set $L$ contain two distinct points $a,b$. Since $L$ is a flat, $\Aline ab\subseteq L\subseteq\overline{\{x,y\}}=\Aline xy$. By Theorem~\ref{t:Alines}, the inclusion $\{a,b\}\subseteq \Aline ab\subseteq L\ssq\Aline xy$ implies $\Aline xy=L=\Aline ab$, which means that $L$ is a line.
\end{proof}

\begin{definition} A subset $P$ of a liner $X$ is called a \index{plane}\defterm{plane} in $X$ if $P$ is a flat of  rank $\|P\|=3$ (and dimension $\dim(L)=\|P\|-1=2$). A liner $X$ is called a \defterm{plane} if it has rank $\|X\|=3$.
\end{definition} 

\begin{proposition}\label{p:L+p=plane} For every line $L$ in a liner $X$ and every point $x\in X\setminus L$, the flat $\overline{L\cup\{x\}}$ is a plane in $X$.
\end{proposition}

\begin{proof} It is clear that $P\defeq\overline{L\cup\{x\}}$ is a flat of rank $2=\|L\|\le \|P\|\le \|L\|+1=3$. Assuming that $\|P\|=2$ and applying Proposition~\ref{p:line=dim1}, we conclude that $P$ is a line and hence $x\in P=L$ by Theorem~\ref{t:Alines}. But this contradicts the choice of the point $x\notin L$. This contradiction shows that $\|P\|=3$ and hence $P$ is a plane.
\end{proof}  

\begin{definition} Two lines $L,\Lambda$ in a liner $X$ are called 
\begin{itemize}
\item \defterm{coplanar}\index{lines!coplanar}\index{coplanar lines} if $\|L\cup \Lambda\|\le 3$;
\item \index{skew lines}\index{lines!skew}\defterm{skew} if $\|L\cup \Lambda\|=4$.
\end{itemize} 
\end{definition}

\begin{exercise} Show that  any skew lines are disjoint. 
\end{exercise}

\section{Ranked liners}

\begin{definition}\label{d:ranked} A liner $X$ is  called \index{ranked liner}\index{liner!ranked}\defterm{ranked} if any two flats $A\subseteq B$ in $X$ of the same finite rank $\|A\|=\|B\|<\w$ are equal.
\end{definition}

\begin{remark} A liner $X$ is ranked if and only if its rank function $\|\cdot\|$ is strictly monotone in the sense that distinct flats $A\subset B$ of finite rank have distinct ranks $\|A\|<\|B\|$.
\end{remark}

Definition~\ref{d:ranked} can be quantified as follows.

\begin{definition}\label{d:k-ranked} A liner $X$ is called \index{$\kappa$-ranked liner}\index{liner!$\kappa$-ranked}\defterm{$\kappa$-ranked} for a cardinal $\kappa$ if any two flats $A\subseteq B$ in $X$ of finite rank $\|A\|=\|B\|\le \kappa$ are equal.
\end{definition}

It is clear that a liner is ranked if and only if it is $\kappa$-ranked for every cardinal $\kappa$ if and only if it is $n$-ranked for every $n\in\w$.

\begin{exercise}\label{ex:hull-0ranked} Prove that every liner is $2$-ranked.
\end{exercise}

The following theorem shows that the rankedness of liners in equivalent to the Exchange Property.

\begin{theorem}\label{t:ranked<=>EP} Let $\kappa$ be a cardinal. A liner is $\kappa$-ranked if and only if it has the $\kappa$-Exchange Property.
\end{theorem}

\begin{proof} For finite cardinals $\kappa$, the equivalence of the $\kappa$-rankedness and the $\kappa$-Exchange Property will be proved by induction on $\kappa$. Definition~\ref{d:k-EP} ensures that every liner has the $0$-Exchange Property. On the other hand, every liner $X$ is $0$-ranked. Indeed, for every flats $A\subseteq B$ in $X$ with $\|A\|=\|B\|=0$, we have $A=\emptyset=B$.\smallskip

Assuming that for some $\kappa\in\w$ the $\kappa$-rankedness is equivalent to the $\kappa$-Exchange Property, we shall prove that the $(\kappa+1)$-rankedness  is equivalent to the $(\kappa+1)$-Exchange Property. 

Take any $(\kappa+1)$-ranked liner $X$. It is $\kappa$-ranked and hence has the $\kappa$-Exchange property by the inductive assumption. To prove that $X$ has the $(\kappa+1)$-Exchange Property, take any set $A\subseteq X$ of cardinality $|A|<\kappa+1$ and any points $x\in X$ and $y\in\overline{A\cup\{x\}}\setminus \overline{A}$.  Observe that $\|A\|\le|A|\le\kappa$ and $\|A\|\le \|A\cup\{x\}\|\le\|A\|+1$. Since $\overline{A}\ne\overline{A\cup\{x\}}$, the $\kappa$-rankedness of the liner $X$ implies that $\|\overline{A\cup\{x\}}\|=\kappa+1$. By analogy we can show that $\|\overline{A\cup\{y\}}\|=\kappa+1$. Since $\overline{A\cup\{y\}}\subseteq\overline{A\cup\{x\}}$ and $\|\overline{A\cup\{y\}}\|=\|\overline{A\cup\{x\}}\|$, the $(\kappa+1)$-rankedness of $X$ implies that $x\in \overline{A\cup\{x\}}=\overline{A\cup\{y\}}$.
\smallskip

Now assume that a liner $X$ has the $(\kappa+1)$-Exchange Property. To prove that $X$ is $(\kappa+1)$-ranked, take any flats $A\subseteq B\subseteq X$ of rank $\|A\|=\|B\|=\kappa+1$. Assuming that $A\ne B$, choose a point $b\in B\setminus A$ and consider the set $A\cup\{b\}$, which has rank $$\kappa+1=\|A\|\le\|A\cup\{b\}\|\le\|B\|=\kappa+1,$$by the monotonicity of the rank.
By Proposition~\ref{p:rank-EP}, there exists a set $C\subseteq A\cup\{b\}$ of cardinality $|C|=\|A\cup\{b\}\|=\kappa+1$ such that $A\cup\{b\}\subseteq \overline{C}$. Assuming that $b\notin C$, we conclude that $b\in\overline{C}\subseteq A$, which contradicts the choice of $b\notin A$. This contradiction shows that $b\in C$ and hence the set $C_{-b}\defeq C\setminus\{b\}$ has cardinality $|C_{-b}|<|C|=\kappa+1=\|A\|$. The definition of the rank $\|A\|$ ensures that there exists a point $a\in A\setminus\overline{C_{-b}}$. Since $a\in \overline{C_{-b}\cup\{b\}}\setminus\overline{C_{-b}}$, the $(\kappa+1)$-Exchange Property of $X$ ensures that $b\in\overline{C_{-b}\cup \{a\}}\subseteq \overline{A}=A$, which contradicts the choice of the point $b\notin A$. This contradiction shows that $A=B$, witnessing that the liner $X$ is $(\kappa+1)$-ranked.
\smallskip

By the Principle of Mathematical Induction, the $\kappa$-rankedness is equivalent to the $\kappa$-Exchange Property for all $\kappa\in\w$. Now let $\kappa$ be an infinite cardinal. 

Assume that a liner $X$ is $\kappa$-ranked. By Proposition~\ref{p:fh-EP<=>wEP}, the $\kappa$-Exchange property of $X$ will follow as soon as we prove that $X$ has the $n$-Exchange Property for every $n\in\w$. Given any number $n\in\w$, take a set $A\subseteq X$ of cardinality $|A|<n$ and points $x\in X$ and $y\in\overline{A\cup\{x\}}\setminus\overline{A}$. The  $\kappa$-rankedness of the liner $X$ implies the $n$-rankedness of $X$, which is equivalent to the $n$-Exchange Property, by the already proved finite case of the theorem. Then $y\in\overline{A\cup\{x\}}\setminus\overline{A}$ implies $x\in\overline{A\cup\{y\}}$. Therefore, the liner $X$ has the $n$-Exchange Property. By Proposition~\ref{p:fh-EP<=>wEP}, $X$ has the Exchange Property and hence $X$ has the $\kappa$-Exchange Property.

Now assume that a liner $X$ has the $\kappa$-Exchange Property. To prove that $X$ is $\kappa$-ranked, take any flats $A\subseteq B$ in $X$ of finite rank $\|A\|=\|B\|\le\kappa$. 
Consider the finite cardinal $n\defeq\|A\|=\|B\|$. The $\kappa$-Exchange Property of $X$ implies the $n$-Exchange Property, which is equivalent to the $n$-rankedness of $X$, by the (already proved) finite case of the theorem. Since $\|A\|=\|B\|\le n$, the $n$-rankedness of $X$ ensures that $A=B$, witnessing that the liner $X$ is $\kappa$-ranked. 
\end{proof}

Theorem~\ref{t:ranked<=>EP} implies the following characterization of  ranked liners.

\begin{corollary}\label{c:ranked<=>EP} A liner is ranked if and only if it has the Exchange Property.
\end{corollary}

The notion of a ranked liner allows us to prove the following characterization of planes in $3$-ranked liners.

\begin{proposition}\label{p:SEP-plane}
Let $X$ be a $3$-ranked liner.
\begin{enumerate}
\item A set $A\subseteq X$ is a plane if and only if $A=\overline{\{x,y,z\}}$ for some points $x\in X$, $y\in X\setminus\{x\}$, and $z\in X\setminus\Aline xy$.
\item If a set $A\subseteq X$ is a plane, then $A=\overline{\{x,y,z\}}$ for every points $x\in A$, $y\in A\setminus\{x\}$, and $z\in X\setminus\Aline xy$.
\end{enumerate}
\end{proposition}

\begin{proof} By Theorem~\ref{t:ranked<=>EP}, the $3$-ranked liner $X$ has the $3$-Exchange Property.
\smallskip

 1. If $A$ is a plane in $X$, then $\|A\|=3$ and by Proposition~\ref{p:rank-EP},  there exist points $a,b,c\in A$ such that $A\subseteq\overline{\{a,b,c\}}$. Since $A$ is flat, $\overline{\{a,b,c\}}\subseteq\overline{A}=A\subseteq\overline{\{a,b,c\}}$ and hence $A=\overline{\{a,b,c\}}$. It follows from $\|\{a,b,c\}\|=\|A\|=3$ that  $b\ne a$ and $c\in\Aline ab$.

If $A=\overline{\{x,y,z\}}$ for some points $x\in X$, $y\in X\setminus\{x\}$, and $z\in X\setminus\Aline xy$, then $A$ is a plane by Proposition~\ref{p:L+p=plane}. 
\smallskip

2. Now assume that $A$ is a plane and take any points $x\in A$, $y\in A\setminus\{x\}$, and $z\notin\Aline xy$.  By Proposition~\ref{p:L+p=plane}, $\overline{\{x,y,z\}}=\overline{\Aline xy\cup\{z\}}$ is a plane and hence $\|\overline{\{x,y,z\}}\|=3=\|A\|$. Since $\overline{\{x,y,z\}}\subseteq A$, the $3$-rankedness of the liner $X$ implies $\overline{\{x,y,z\}}=A$.
\end{proof}

\begin{exercise} Prove that a liner $X$ is $3$-ranked if and only if every plane $A$ is equal to the flat hull $\overline{\{x,y,z\}}$ of any points $x\in A$, $y\in A\setminus\{x\}$ and $z\in A\setminus\Aline xy$.
\end{exercise}

\begin{exercise} Find a liner $X$ containing a plane $A\ne X$ such that $\|A\|=3=\|X\|$.
\smallskip

\noindent{\em Hint:} Look at the liner from Example~\ref{ex:Tao}.
\end{exercise}

\section{The relative rank and codimension of sets in liners}

\begin{definition} For two sets $A,B$ in a liner $X$, let $\|A\|_B$ be the smallest cardinality of a set $C\subseteq X$ such that $A\subseteq\overline{B\cup C}$. The \index[note]{$\|A\|_B$}cardinal $\|A\|_B$ is called the \index{relative rank}\defterm{$B$-relative rank} (or just the \defterm{$B$-rank}) of the set $A$ in $X$. If $B\subseteq A$, then the $B$-rank $\|A\|_B$ is denoted by \index[note]{$\dim_B(A)$}$\dim_B(A)$ and called the \index{codimension}\defterm{codimension} of $B$ in $A$.
\end{definition}

Observe that every set $A$ in a liner has rank $\|A\|=\|A\|_\emptyset$. So, the rank is a special case of the $B$-rank for $B=\emptyset$. 

\begin{exercise} Show that the relative rank is bimonotone in the sense that $\|A\|_C\le\|B\|_D$ for any sets $A\subseteq B$ and $D\subseteq C$ in a liner.
\end{exercise}

\begin{proposition}\label{p:rank-EP3} If a liner $X$ has the $\kappa$-Exchange Property for some cardinal $\kappa$, then for every subsets $A,B\subseteq X$ with $\|A\cup B\|\le \kappa$, there exist sets $A'\subseteq A$ and $B'\subseteq B\setminus\overline A$ such that $A\subseteq \overline {A'}$,
 $B\subseteq \overline{A'\cup B'}$, $|A'\cup B'|=\|A\cup B\|$, and $|B'|=\|B\|_A$.
\end{proposition}

\begin{proof}  We divide the proof into two cases, considered in the following two claims.

\begin{claim}\label{cl:rank-EP4} For every sets $A,B\subseteq X$ with finite relative rank $\|B\|_A\le\kappa$, there exist sets $A'\subseteq A$ and $C\subseteq B\setminus\overline A$  such that $A\subseteq \overline {A'}$,
 $B\subseteq \overline{A'\cup C}$, $|A'\cup C|=\|A\cup B\|$, and $|C|=\|B\|_A$.
\end{claim}

\begin{proof} By Proposition~\ref{p:rank-EP2}, there exist  sets $A'\subseteq A$ and $B'\subseteq B\setminus\overline A$ such that $A\subseteq\overline{A'}$, $B\subseteq\overline{A'\cup B'}$ and $|A'\cup B'|=\|A\cup B\|$. Then $B\subseteq\overline{A'\cup B'}\subseteq\overline{A\cup B'}$ and hence $\|B\|_A\le|B'|$.  By the definition of the $A$-rank $\|B\|_A$, there exists a set $C\subseteq X\setminus\overline A$ of cardinality $|C|=\|B\|_A$ such that $B\subseteq\overline{A\cup C}$. We can assume that the intersection $|C\cap B'|$ has the largest possible cardinality. We claim that in this case $C\subseteq B'$. To derive a contradiction, assume that $C\not\subseteq B'$ and choose a point $c\in C\setminus B'$. Assuming that $B'\subseteq \overline{A\cup(C\setminus\{c\})}$, we conclude that $B\subseteq\overline{A'\cup B'}\subseteq\overline{A\cup(C\setminus\{c\})}$ and hence $\|B\|_A\le|C\setminus\{c\}|<|C|=\|B\|_A$, which is acontradiction showing that the set  $B'\setminus\overline{A\cup(C\setminus\{c\})}$ contains some point $b$. Since $b\in B'\subseteq B\subseteq\overline{A\cup C}=\overline{A'\cup C}=\overline{A'\cup(C\setminus\{c\})\cup\{c\}}$, by Proposition~\ref{p:aff-finitary}, there exists a finite set $F\subseteq A'\cup(C\setminus\{c\})$ such that $b\in\overline{F\cup\{c\}}$. It follows from $F\subseteq A'\cup(C\setminus\{c\})$ and $b\notin\overline{A\cup(C\setminus\{c\})}$ that $b\notin\overline F$. Since $F\cup\{c\}\subseteq A'\cup C$, $|C|=\|A\|_B\le|B'|$ and $A'\cap C\subseteq A\cap(X\setminus\overline A)=\emptyset=A'\cap B'$, we have the inequality
$$|F|<|F\cup\{c\}|\le |A'\cup C|=|A'|+|C|\le|A'|+|B'|=|A'\cup B'|=\|A\cup B\|\le\kappa.$$ Since $b\in\overline{F\cup\{c\}}\setminus\overline F$, the $\kappa$-Exchange Property of $X$ implies $c\in \overline{F\cup \{b\}}$.
 Now concider the set $C'=(C\setminus\{c\})\cup\{b\}$ and observe that $|C'|=|C|=\|B\|_A$ and 
$$
\begin{aligned}
B&\subseteq \overline {A\cup C}=\overline{A\cup(C\setminus\{c\})\cup\{c\}}\subseteq\overline{A\cup(C\setminus\{c\})\cup\overline{F\cup\{b\}}}\\
&\subseteq\overline{A\cup(C\setminus\{c\})\cup\overline{A'\cup(C\setminus\{c\})\cup\{b\}}}=\overline{A\cup(C\setminus\{c\})\cup\{b\}}=\overline{A\cup C'}.
\end{aligned}
$$
Since $|C'\cap B'|=|((C\setminus\{c\})\cup\{b\})\cap B'|=|C\cap B'|+1$, we obtain a contradiction with the maximality of $|C\cap B'|$. This contradiction shows that $C\subseteq B'$. It is clear that $|A'\cup C|=|A'|+|C|\le|A'|+|B'|=|A'\cup B'|=\|A\cup B\|$. Since $A\cup B\subseteq\overline{A\cup C}=\overline{A'\cup C}$, the definition of the rank $\|A\cup B\|$ ensures that $\|A\cup B\|\le|A'\cup C|$ and hence $|A'\cup C|=\|A\cup B\|$.
\end{proof}

\begin{claim}\label{cl:rank-EP5} For any sets $A,B\subseteq X$ with infinite relative rank $\|B\|_A\le\kappa$, there exist sets $A'\subseteq A$ and $B'\subseteq B\setminus\overline A$ such that $A\subseteq \overline {A'}$,
 $B\subseteq \overline{A'\cup B'}$, $|A'\cup B'|=\|A\cup B\|$, and $|B'|=\|B\|_A$.
\end{claim}

\begin{proof} By Proposition~\ref{p:rank-EP2}, there exist sets $A''\subseteq A$ and $B''\subseteq B\setminus\overline A$ such that $A\subseteq \overline{A''}$, $B\subseteq\overline{A''\cup B''}$ and $|A''\cup B''|=\|A\cup B\|$. By Proposition~\ref{p:rank-EP}, there exists a set $D\subseteq X\setminus A$ of (infinite) cardinality $|D|=\|B\|_A$ such that $B\subseteq \overline{A\cup D}$. Consider the family $[D]^{<\w}$ of finite subsets of $D$ and observe that $|[D]^{<\w}|=|D|=\|B\|_A$. For every finite set $F\in[D]^{<\w}$, consider the set $B_F\defeq B''\cap\overline{A\cup  F}$ and observe that  $\|B_F\|_A\le |F|<\omega$. By Claim~\ref{cl:rank-EP4}, there exist sets $A_F'\subseteq A''$ and $B'_F\subseteq B_F$ such that $A''\subseteq\overline{A'_F}$, $B_F\subseteq \overline{A'_F\cup B'_F}$, $|A'_F\cup B'_F|=\|A''\cup B_F\|$ and $|B'_F|=\|B_F\|_{A''}=\|B_F\|_A$. Then the sets $A'\defeq\bigcup_{F\in[D]^{<\w}}A'_F\subseteq A''$ and $B'\defeq \bigcup_{F\in[D]^{<\w}}B_F'\subseteq B''$ have the required  properties:
$$
\begin{aligned}
&A\subseteq \overline{A''}\subseteq \bigcup_{F\in[D]^{<\w}}\overline{A'_F}\subseteq \overline{A'},\\
&B''=B''\cap\overline{A\cup D}=\bigcup_{F\in[D]^{<\w}}B''\cap\overline{A\cup F}=\bigcup_{F\in[D]^{<\w}}B_F\subseteq\bigcup_{F\in[D]^{<\w}}\overline{A'_F\cup B'_F}\subseteq\overline{A'\cup B'},\\
&B\subseteq\overline{A''\cup B''}\subseteq\overline{\overline{A'}\cup\overline{A'\cup B'}}=\overline{A'\cup B'},\\
&|A'\cup B'|\le|A''\cup B''|=\|A\cup B\|,\\
&|B'|\le\sum_{F\in[D]^{<\w}}|B'_F|\le\sum_{F\in[D]^{<\w}}|F|\le|[D]^{<\w}|\cdot\w=|D|=\|B\|_A
\end{aligned}
$$
Since $A\cup B\subseteq\overline{A'\cup B'}$, the definition of the cardinal $\|A\cup B\|$ ensures that $\|A\cup B\|\le|A'\cup B'|$ and hence $|A'\cup B'|=\|A\cup B\|$.
\end{proof}
\end{proof}

\section{Independence in liners}

\begin{definition} A subset $A$ of a liner $X$ is \index{independent subset}\index{subset!independent}\defterm{independent} if $a\notin\overline{A\setminus\{a\}}$ for every $a\in A$. 
\end{definition}

\begin{exercise} Show that any subset of an independent set in a liner is independent.
\end{exercise}

\begin{exercise} Find an example of a liner containing two independent sets whose union is not independent.
\end{exercise}

\begin{exercise} Let $X$ be a liner and $Y$ be a linear subspace of $X$. Show that a set $A\subseteq Y$ is independent in $Y$ if $A$ is independent in $X$. Is the converse true?
\end{exercise}

\begin{proposition}\label{p:inter-indep-flats} If $A$ is an independent set in a ranked liner $X$, then for every nonempty finite set $F\subseteq A$, we have $\overline{A\setminus F}=\bigcap_{a\in F}\overline{A\setminus\{a\}}$.
\end{proposition}

\begin{proof} It is clear that  $\overline{A\setminus F}\subseteq\bigcap_{a\in F}\overline{A\setminus\{a\}}$. Assuming that  $\overline{A\setminus F}\ne\bigcap_{a\in F}\overline{A\setminus\{a\}}$, we can take a nonempty finite set $F\subseteq A$ of smallest possible cardinality such that $\overline{A\setminus F}\ne\bigcap_{a\in F}\overline{A\setminus\{a\}}$. This inequality ensures that $|F|\ge 2$. Choose any point $e\in F$ and consider the nonempty finite set $E\defeq F\setminus\{e\}$. The minimality of $F$ ensures that $\overline{A\setminus E}=\bigcap_{a\in E}\overline{A\setminus\{a\}}$. 
Since $$\overline{A\setminus F}\ne\bigcap_{a\in F}\overline{A\setminus\{a\}}=\overline{A\setminus\{e\}}\cap \bigcap_{a\in E}\overline{A\setminus\{a\}}=\overline{A\setminus\{e\}}\cap 
\overline{A\setminus E},$$ there exists a point $x\in (\overline{A\setminus\{e\}}\cap \overline{A\setminus E})\setminus\overline{A\setminus F}$. Since $x\in \overline{A\setminus E}\setminus \overline{A\setminus F}=\overline{(A\setminus F)\cup\{e\}}\setminus\overline{A\setminus F}$, the Exchange Property of the ranked liner $X$ ensures that $e\in \overline{(A\setminus F)\cup\{x\}}\subseteq \overline{A\setminus\{e\}}$, which contradicts the independence of the set $A$.
\end{proof}



The notion of an independent set is a partial case of a more general notion of a $B$-independent set.

\begin{definition} Let $B$ be a set in a liner $X$. A set $A\subseteq X$ is called \index{subset!independent}\defterm{$B$-independent} if $a\notin \overline{B\cup(A\setminus\{a\})}$ for every $a\in A$.
\end{definition}

Therefore, a set $A$ in a liner $X$ is independent if and only if it is $\emptyset$-independent.

\begin{exercise} Let $B$ be a set in a liner $X$. Show that for every set $A\subseteq X$ of finite $B$-rank $\|A\|_B$, there exists a $B$-independent set $C\subseteq X\setminus\overline B$ of cardinality $|C|=\|A\|_B$ such that $A\subseteq\overline{C}$.
\end{exercise}


%

%

Now we prove two propositions on unions of indepenedent sets.

\begin{proposition}\label{p:union-of-independent} Let $A$ be a set in a liner $X$, $I$ be an $A$-independent set in $X$ and $J$ is an $(A\cup I)$-independent set in $X$. If the liner $X$ has the $\|A\cup I\cup J\|$-Exchange Property, then the set $I\cup J$ is $A$-independent.
\end{proposition}

\begin{proof} To derive a contradiction, assume that $I\cup J$ is not $A$-independent and find a point $x\in I\cup J$ such that $x\in \overline{A\cup((I\cup J)\setminus\{x\})}$. If $x\in J$, then $x\in\overline{A\cup I\cup(J\setminus\{x\})}$ contradicts the $(A\cup I)$-independence of the set $J$. So, $x\in I$. By Proposition~\ref{p:rank-EP}, there exists a set $B\subseteq A\cup(I\setminus\{x\})\cup J$ of cardinality $|B|=\|A\cup (I\setminus\{x\})\cup J\|\le\|A\cup I\cup J\|$ such that $A\cup (I\setminus\{x\})\cup J\subseteq\overline{B}$. 
By Proposition~\ref{p:aff-finitary}, there exists a finite set $F\subseteq B$ such that 
$x\in \overline{F}$. We can assume that $F$ has the smallest possible cardinality. Since $x\in F\subseteq B\subseteq A\cup(I\setminus\{x\})\cup J$, the $A$-independence of the set $I$ ensures that $J\cap F$ is not empty and hence contains some point $y$. The minimality of the set $F$ guarantees  that $x\notin\overline{F\setminus\{y\}}$. Then $x\in\overline{(F\setminus\{y\})\cup\{y\})}\setminus\overline{F\setminus\{y\}}$. Since $|F\setminus\{y\}|<|F|\le |B|\le\|A\cup I\cup J\|$, the $\|A\cup I\cup J\|$-Exchange Property of $X$ ensures that $$y\in\overline{(F\setminus\{y\})\cup\{x\}}\subseteq \overline{A\cup(I\setminus\{x\})\cup(J\setminus\{y\})\cup\{x\}}=\overline{(A\cup I)\cup(J\setminus\{y\})},$$which contradicts the $(A\cup I)$-independence of the set $J$. This contradiction shows that the set $I\cup J$ is $A$-independent.
\end{proof}

\begin{proposition}\label{p:add-point-to-independent} Let $A,I$ be two sets in a liner $X$. If the set $I$ is $A$-independent and the liner $X$ is $\|A\cup I\|$-ranked, then for every point $x\in X\setminus\overline{A\cup I}$, the set $I\cup\{x\}$ is $A$-independent.
\end{proposition}

\begin{proof}  By Theorem~\ref{t:ranked<=>EP}, the $\|A\cup I\|$-ranked  liner $X$ has the $\|A\cup I\|$-Exchange Property.

Assuming that the set $I\cup\{x\}$ is not $A$-independent, we can find a point $y\in I\cup\{x\}$ such that $y\in \overline{A\cup ((I\cup\{x\})\setminus\{y\})}$. 
If $y=x$, then $x=y\in \overline{A\cup ((I\cup\{x\})\setminus\{y\})}=\overline{A\cup I}$, which contradicts the choice of the point $x\notin \overline{A\cup I}$. This contradiction shows that $x\ne y$ and hence $y\in (I\cup\{x\})\setminus\{x\}=I$. Then $A$-independence of the set $I$ ensures that $y\notin \overline{A\cup(I\setminus\{y\})}$ and hence $\overline{A\cup(I\setminus\{y\})}\ne\overline{A\cup I}$.

If the cardinal $\|A\cup I\|$ is infinite, then the liner $X$ has the Exchange Property, by Proposition~\ref{p:fh-EP<=>wEP}. In this case $y\in \overline{A\cup(I\setminus \{y\})\cup\{x\}}\setminus\overline{A\cup(I\setminus\{y\})}$ implies $x\in\overline{A\cup(I\setminus\{y\})\cup\{y\}}=\overline{A\cup I}$, which contradicts the choice of the point $x$.

This contradiction shows that the cardinal $\|A\cup I\|$ is finite. Since $\overline{A\cup(I\setminus\{y\})}\ne\overline{A\cup I}$, the $\|A\cup I\|$-rankedness of the liner $X$ ensures that $\| A\cup(I\setminus\{y\})\|<\|A\cup I\|$. By Proposition~\ref{p:rank-EP}, there exists a set $F\subseteq A\cup(I\setminus\{y\})$ of cardinality $|F|=\|A\cup(I\setminus\{y\})\|<\|A\cup I\|$ such that $\overline{A\cup(I\setminus\{y\})}=\overline F$.

 Since $y\in\overline{A\cup(I\setminus\{y\})\cup\{x\}}\setminus\overline{A\cup(I\setminus\{y\})}=\overline{F\cup\{x\}}\setminus\overline F$ and $|F|<\|A\cup I\|$, the $\|A\cup I\|$-Exchange Property of $X$ ensures that $x\in\overline{F\cup\{y\}}\subseteq\overline{A\cup(I\setminus\{y\})\cup\{y\}}=\overline{A\cup I}$, which contradicts the choice of $x$. This contradiction shows that the set $I\cup \{x\}$ in $A$-independent.
\end{proof}

\begin{remark} Observe that Proposition~\ref{p:add-point-to-independent} does not follow from Proposition~\ref{p:union-of-independent}.
\end{remark}

\begin{definition} Let $B,A$ be two sets in a liner $X$. A $B$-independent set $I\subseteq A$ is called a \index{maximal independent set}\index{subset!maximal independent} \defterm{maximal $B$-independent in} $A$ if $I=J$ for any $B$-independent set $J$ with $I\subseteq J\subseteq A$.
\end{definition} 

\begin{exercise} Let $A$ be a set in a liner and $I$ be an $A$-independent set in $X$. Show that $I\cap\overline A=\emptyset$.
\end{exercise}

\begin{lemma}\label{l:Max-indep} For every sets $A,B$ in a liner $X$, every $B$-independent set $I\subseteq A$ can be enlarged to a maximal $B$-independent set $J\subseteq A$.
\end{lemma}

\begin{proof} Consider the family $\mathcal I$ if all $B$-independent sets in $A$ that contain the $B$-independent set $I$. The family $\mathcal I$ is endowed with the inclusion partial order. Using the Kuratowski--Zorn Lemma, choose a maximal chain $\C$ in the partially ordered set $\mathcal I$. We claim that the union $J\defeq \bigcup\C\subseteq A$  is an $B$-independent set in $X$. Assuming that $J$ is not $B$-independent, we can find a point $x\in J=J\setminus\overline{B}$ such that $x\in\overline{B\cup(J\setminus \{x\})}$. By Proposition~\ref{p:aff-finitary},  there exists a finite set $F\subseteq J\setminus\{x\}$ such that $x\in\overline{B\cup F}$. Since $\C$ is chain of sets with $\{x\}\cup F\subseteq J=\bigcup\C$, there exists a $B$-independent set $C\in\C$ such that $\{x\}\cup F\subseteq C$ and hence $F\subseteq C\setminus\{x\}$. Then $x\in \overline{B\cup F}\subseteq\overline{B\cup(C\setminus\{x\})}$, which contradicts the $B$-independence of the set $C$. This contradiction shows that the set $J$ is $B$-independent.
\end{proof}

\begin{theorem}\label{t:Max=codim} Let $A,B$ be any sets in a liner $X$. If the liner $X$ has the $\|A\cup B\|$-Exchange Property, then the $A$-rank $\|B\|_A$ of the set $B$ is equal to the cardinality of any maximal $A$-independent set in $B$. 
\end{theorem}

\begin{proof} Let $M$ be any maximal $A$-independent set in $B$. The $A$-independence of $M$ ensures that $M\cap\overline A=\emptyset$. We claim that $B\subseteq \overline{A\cup M}$. In the opposite case, we can choose a point $x\in B\setminus\overline{A\cup M}$. Since $\|A\cup M\|\le\|A\cup B\|$ and the liner $X$ is $\|A\cup B\|$-ranked, by Proposition~\ref{p:add-point-to-independent}, the set $M\cup\{x\}$ is $A$-independent, which contradicts the maximality of the $B$-independent set $M$. This contradiction shows that $B\subseteq\overline{A\cup M}$ and hence $\|B\|_A=\|M\|_A$.

By Proposition~\ref{p:rank-EP3}, there exist sets $A'\subseteq A$ and $B'\subseteq M$ such that $A\subseteq \overline{ A'}$, $M\subseteq \overline{A'\cup B'}$, $|A'\cup B'|=\|A\cup M\|$ and $|B'|=\|M\|_A=\|B\|_A$. Assuming that $B'\ne M$, we can choose a point $x\in M\setminus B'$ and observe that $x\in M\subseteq\overline{A\cup B'}\subseteq\overline{A\cup(M\setminus\{x\})}$, which contradicts the $A$-independence of the set $M$. This contradiction shows that $B'=M$ and hence $|M|=|B'|=\|M\|_A=\|B\|_A$.
\end{proof}

Theorem~\ref{t:Max=codim} implies the following corollary.

\begin{corollary}\label{c:Max=dim} Let $A$ be a set in a liner $X$. If the liner $X$ has the $\|A\|$-Exchange Property, then the rank $\|A\|$ of the set $A$ is equal to the cardinality of any maximal independent set in $A$. 
\end{corollary}

\begin{proposition}\label{p:rank+} For every sets $A,B,C$ in a liner $X$, 
$$\|C\|_A\le\|B\|_{A}+\|C\|_{B}.$$If $\overline A\subseteq \overline B\subseteq \overline C$ and the liner $X$ has the $\|C\|$-Exchange Property, then $\|C\|_A=\|B\|_A+\|C\|_B$.
\end{proposition}

\begin{proof} By the definition of the relative ranks $\|B\|_{A}$ and $\|C\|_{B}$, there exist sets $B',C'\subseteq X$ of cardinality $|B'|=\|B\|_{A}$ and $|C'|=\|C\|_{B}$ such that $B\subseteq\overline{A\cup B'}$ and $C\subseteq \overline{B\cup C'}$. Since
$$C\subseteq\overline{B\cup C'}\subseteq\overline{\overline{A\cup B'}\cup C'}\subseteq\overline{A\cup B'\cup C'},$$we obtain the desired inequality $\|C\|_{A}\le|B'\cup C'|\le|B'|+|C'|=\|B\|_{A}+\|C\|_{B}$.
\smallskip

Now assume that $\overline A\subseteq \overline B\subseteq \overline C$ and the liner $X$ has the $\|C\|$-Exchange Property.  Using Lemma~\ref{l:Max-indep}, choose a maximal $A$-independent set $I\subseteq B\setminus\overline A$ and  a maximal $B$-independent set $J\subseteq C\setminus\overline B$. 
Since $\|A\cup I\|\le \|C\|$, Proposition~\ref{p:add-point-to-independent} implies that $\overline B=\overline{A\cup I}$ and hence the $B$-independent set $J$ is $(A\cup I)$-independent. Since the liner $X$ has the $\|C\|$-Exchange Property and $\|A\cup I\cup J\|\le\|C\|$, Proposition~\ref{p:union-of-independent} implies that the set $I\cup J$ is $A$-independent. Since $\|B\cup J\|\le\|C\|$, Proposition~\ref{p:add-point-to-independent} implies that $\overline{A\cup I\cup J}=\overline{B\cup J}=\overline C$ and hence the set $I\cup J$ is maximal $A$-independent in $C$. Theorem~\ref{t:Max=codim} ensures that $\|B\|_A=|I|$, $\|C\|_B=|J|$ and 
$$\|C\|_A=|I\cup J|=|I|+|J|=\|B\|_A+\|C\|_B.$$
\end{proof}

\begin{corollary}\label{c:rank+} For every sets $A,B$ in a liner $X$, 
$\|B\|\le\|A\|+\|B\|_{A}$. If $\overline A\subseteq \overline B$ and the liner $X$ has the $\|B\|$-Exchange Property, then $\|B\|=\|A\|+\|B\|_A$ and $\dim(B)=\dim(A)+\dim_A(B)$.
\end{corollary}

\begin{exercise}
 Find an example of two flats $A\subseteq B$ in a liner with $\|B\|<\|A\|+\|B\|_A$.

{\em Hint:} Look at the liner  from Example~\ref{ex:Tao}.
\end{exercise}

We recall that a flat $A$ is a {\em hyperplane} in a flat $B$ if $A\subset B$ and every flat $C$ with $A\subseteq C\subseteq B$ is equal either to $A$ or to $B$.

\begin{proposition}\label{p:hyperplane} A flat $B$ in a flat $A$ in a ($\|A\|$-ranked) liner $X$ is a hyperplane in $A$ if (and only if) $\|A\|_B=1$.
\end{proposition}

\begin{proof} If $B$ is a hyperplane in $A$, then $B\ne A=\overline{B\cup\{x\}}$ for any point $x\in A\setminus B$ and hence $\|A\|_B=1$. This proves the ``only if'' part of the characterization.
\smallskip

To prove the ``if'' part, assume that $\|A\|_B=1$ and the liner $\|A\|$-ranked. By Theorem~\ref{t:ranked<=>EP}, $X$ has the $\|A\|$-Exchange Property. It follows from $\|A\|_B=1>0$ that $B\ne A$. To show that $B$ is a hyperplane in $A$, take any flat $C$ with $B\subset C\subseteq A$. Fix any point $x\in C\setminus B$ and observe that $\{x\}$ is an $B$-independent set. By Theorem~\ref{t:Max=codim}, the singleton $\{x\}$ is a maximal $B$-independent in $A$. By Proposition~\ref{p:add-point-to-independent}, $A=\overline{B\cup\{x\}}\subseteq C\subseteq A$ and hence $C=A$, witnessing that $C$ is a hyperplane in $A$.
\end{proof}

\begin{exercise} Find an example of a liner $X$ and a flat $B$ in $X$ such that $\|X\|_B=1$ but $B$ is not a hyperplane in $X$.
\smallskip

\noindent{\em Hint:} Look at the liner from Example~\ref{ex:Tao}.
\end{exercise}



\section{The ranks of flat functions in liners} 

For a relation $R\subseteq X\times X$ in a liner $X$, the cardinal\index[note]{$\|R\|$} 
$$\|R\|\defeq \max\{\|\dom[R]\|,\|\rng[R]\|\}$$ is called the \index{rank}\index{relation!rank of}\defterm{rank} of $R$.
It is clear that $\|R\|=\|R^{-1}\|$. Also for any  relations $F,G\subseteq X\times X$ we have $\dom[FG]\subseteq\dom[G]$ and $\rng[FG]\subseteq\rng[F]$, which implies
$$\|FG\|=\max\{\|\dom[FG]\|,\|\rng[FG]\|\}\le\max\{\|\rng[F]\|,\|\dom[G]\|\}\le\max\{\|F\|,\|G\|\}.$$ Consequently, for every cardinal $\kappa$, the sets
$$\{F\in\mathcal{R}_{X}:\|F\|\le\kappa\}\quad\mbox{and}\quad\{F\in\mathcal{FI}_{X}:\|F\|<\kappa\} $$are  subsemigroups of the monoid $\mathcal{R}_{X}$ of relations on $X$, and so are their intersections with the semigroups $\mathcal{FR}_X$, $\mathcal{FF}_X$, $\mathcal{FI}_X$.

We recall that a relation $R\subseteq X\times Y$ between liners $X,Y$ is {\em flat} if for every flats $A\subseteq X$ and $B\subseteq Y$ the sets $$R[A]\defeq\{y\in Y:\exists x\in A\;(x,y)\in R\}\quad\mbox{and}\quad R^{-1}[B]\defeq\{x\in X:\exists y\in B\;(x,y)\in R\}$$are flat in the liners $Y$ and $X$, respectively.

\begin{proposition}\label{p:rank-function-monotone} Let $X$ be a ranked liner and $F\subseteq X\times X$ be a flat function. For every flat $A$ in $X$ we have $\|F[A]\|\le\|A\|$. 
\end{proposition}

\begin{proof}  Let $M$ be a maximal independent set in $F[A]$. Choose a set $M'\subseteq A$ such that the function $F{\restriction}_{M'}:M'\to M$ is bijective. We claim that the set $M'$ is independent in $X$. Indeed, take any element $a\in M'$ and let $b=F(a)\in M$ be its image under the function $F$. The independence of the set $M$ implies that $b\notin \overline{M\setminus\{b\}}$. Then $F^{-1}(b)\cap F^{-1}[\overline{M\setminus\{b\}}]=\varnothing$ and hence $a\notin F^{-1}[\overline{M\setminus\{b\}}]$. Since the function $F$ is flat, the set $F^{-1}[\overline{M\setminus\{b\}}]$ is a flat in $X$ containing the set $M'\setminus\{a\}$. Then $$\overline{M'\setminus\{a\}}\subseteq F^{-1}[\overline{M\setminus\{b\}}]\subseteq F^{-1}[X\setminus\{b\}]\subseteq X\setminus\{a\},$$ witnessing that the set $M'\subseteq A$ is independent and hence $\|A\|\ge|M'|=|M|=\|F[A]\|$, see Theorem~\ref{t:Max=codim}.
\end{proof}

\begin{corollary}\label{c:flat=injection=>rank=} If $X$ is a ranked liner and $F\subseteq X\times X$ is a flat injective function, then $\|F[A]\|=\|A\|$ for every flat $A\subseteq\dom[F]$. 
\end{corollary}

\begin{proof} By Proposition~\ref{p:rank-function-monotone}, for every flat $A\subseteq\dom[F]$, we have $\|B\|=\|F[A]\|\le \|A\|$.  Since the function $F$ is flat and injective, its inverse $F^{-1}$ is a flat function in $X$ and the set $B\defeq F[A]$ is a flat in $X$. Moreover, $A=F^{-1}[B]$ by the injectivity of $F$ and the inclusion $A\subset\dom[F]$. By Proposition~\ref{p:rank-function-monotone}, $$\|A\|=|F^{-1}[B]\le\|B\|=\|F[A]\|\le\|A\|$$and hence $\|F[A]\|=\|A\|$.
\end{proof}

\begin{corollary}\label{c:injection=dim} Every flat injection $F\in\mathcal{FI}_{X}$ in a ranked liner $X$ has rank $$\|F\|=\|\dom[F]\|=\|\rng[F]\|.$$
\end{corollary}

\chapter{Regularity and Parallelity Axioms}

\section{Regularity Axioms}

In this section we introduce four properties of liners, describing the structure of the hull $\overline{A\cup\{a\}}$ of a flat $A$ with an attached singleton $\{a\}$.

For two subsets $A,B$ of a liner $X$ and a point $x\in X$, let\index[note]{$\Aline Ax$}\index[note]{$\Aline xA$}\index[note]{$\Aline AB$}
$$\Aline Ax=\Aline xA\defeq\bigcup_{a\in A}\Aline ax\quad\mbox{and}\quad\Aline AB\defeq\bigcup_{x\in A}\bigcup_{y\in B}\Aline xy=\bigcup_{a\in A}\Aline aB=\bigcup_{b\in B}\Aline Ab.$$

\begin{definition}\label{d:regular} A liner $X$ is called
\begin{itemize}
\item \index{strongly regular liner}\index{liner!strongly regular}\defterm{strongly regular} if for every  nonempty flat $A\subseteq X$ and point $b\in X\setminus A$, we have\newline $\overline{A\cup\{b\}}=\Aline Ab$;
\item \index{regular liner}\index{liner!regular}\defterm{regular} if for every flat $A\subseteq X$ and points $a\in A $, $b\in X\setminus A$ we have\newline $\overline{A\cup\{b\}}=\bigcup_{y\in\Aline ab}\Aline Ay$;
\item \index{proregular liner}\index{liner!proregular}\defterm{proregular} if for every flat $A\subseteq X$ and points $a\in A $, $b\in X\setminus A$ with $\Aline ab\ne\{a,b\}$, we have\newline $\overline{A\cup\{b\}}=\bigcup_{y\in\Aline ab}\Aline Ay$;
\item \index{weakly regular liner}\index{liner!weakly regular}\defterm{weakly regular} if for every flat $A\subseteq X$ and points $a\in A $, $b\in X\setminus A$ we have\newline $\overline{A\cup\{b\}}=\bigcup_{x\in A}\overline{\{a,b,x\}}$.
\end{itemize}
\end{definition}

Definition~\ref{d:regular} implies the following simple proposition.

\begin{proposition}\label{p:sr=>r=>wr} Let $X$ be a liner.
\begin{enumerate}
\item If $X$ is strongly regular, then $X$ is regular.
\item If $X$ is regular, then $X$ is proregular and weakly regular.
\end{enumerate}
\end{proposition}



The notion of a (pro)regular liner can be quantified as follows.

\begin{definition}\label{d:k-regular} Let $\kappa$ be a cardinal. A liner $X$ is called \index{$\kappa$-regular liner}\index{liner!$\kappa$-regular}\defterm{$\kappa$-regular} (resp. \index{proregular liner}\index{liner!proregular}\defterm{$\kappa$-proregular}\/) if for every set $A\subseteq X$ of cardinality $|A|<\kappa$ and every points $o\in\overline A $, $p\in X\setminus\overline A$ (resp. with $\Aline op\ne\{o,p\}$), we have $\overline{\{p\}\cup A}=\bigcup_{u\in\overline{o\,p}}\bigcup_{a\in \overline A}\Aline ua$.
\end{definition}

Observe that a liner is (pro)regular if and only if it is $\kappa$-(pro)regular for every cardinal $\kappa$. 

\begin{exercise} Show that a $3$-long liner $X$ is $\kappa$-regular for some cardinal $\kappa$ if and only if $X$ is $\kappa$-proregular.
\end{exercise}

\begin{exercise} Show that every liner is $2$-regular.
\end{exercise}

\begin{proposition}\label{p:k-regular<=>2ex} Let $\kappa$ be a cardinal. Every $\kappa$-proregular liner is $\kappa$-ranked and has the $\kappa$-Exchange Property.
\end{proposition}

\begin{proof} To prove the $\kappa$-Exchange Property of the $\kappa$-proregular liner $X$, take any set $A\subseteq X$ of cardinality $|A|<\kappa$, and any points $x\in X$ and $y\in\overline{A\cup\{x\}}\setminus\overline A$. If for every $o\in\overline A$ the line $\Aline ox$ coincides with $\{o,x\}$, then the set $\overline A\cup\{x\}$ is flat and $\overline{A\cup\{x\}}=\overline A\cup\{x\}$ and hence $y\in \overline{A\cup\{x\}}\setminus \overline A\subseteq\{x\}$. Then $x=y\in\overline{A\cup\{y\}}$ and we are done. So, assume that $\Aline ox\ne\{o,x\}$ for some point $o\in\overline A$. In this case, the $\kappa$-proregularity of $X$ ensures that $y\in \Aline\alpha\beta$ for some $\alpha\in\overline A$ and $\beta\in\Aline ox$. It follows from $y\in\Aline \alpha\beta\overline A$ that $\beta\notin\overline A$ and $\alpha\ne y$. Then $\beta\in\Aline\alpha\beta=\Aline \alpha y$ and $x\in\Aline ox=\Aline o\beta\subseteq\overline{\{o,\alpha,y\}}\subseteq \overline{\overline A\cup\{y\}}=\overline{A\cup\{y\}}$, witnessing that $X$ has the $\kappa$-Exchange Property.
\end{proof}

\begin{example}\label{ex:Steiner13} The $13$-element group $\IZ_{13}$ endowed with the line relation
$$\Af\defeq\big\{(x,y,z)\in \IZ_{13}^3:\{x,y,z\}\in\big\{a+\{0,3,4\},a+\{0,5,7\}:a\in\IZ_{13}\big\}\big\}$$ is an example of a finite liner, which is ranked but not $3$-proregular.
\end{example}

Proposition~\ref{p:k-regular<=>2ex} implies the following corollary.

\begin{corollary}\label{c:proregular=>ranked} Every proregular liner is ranked and has the Exchange Property.
\end{corollary}

The Exchange Property in weakly regular liners is equivalent to the $3$-Exchange Property.

\begin{proposition}\label{p:wr-ex<=>3ex} A weakly regular liner $X$ has the Exchange Property if and only if it has the $3$-Exchange Property.
\end{proposition}

\begin{proof} The ``only if'' part is trivial. To prove the ``if'' part, assume that a weakly regular liner $X$ has the $3$-Exchange Property. To show that $X$ has the Exchange Property, take any flat $A\subseteq X$ and points $x\in X\setminus A$ and $y\in\overline{A\cup\{x\}}\setminus A$. 
If $A=\emptyset$, then $y\in\overline{A\cup\{x\}}\setminus A=\{x\}$ and hence $x=y\in\overline{A\cup\{x\}}$. So, assume that $A\setminus \emptyset$ and fix a point $o\in A$. If $y\in \Aline ox$, then $x=\Aline ox=\Aline oy\subseteq\overline{A\cup\{y\}}$.

So, we assume that $y\notin\Aline ox$. By the weak regularity of $X$, there exists a point $a\in A$ such that $y\in \overline{\{a,o,x\}}$. Since $y\in\overline{\{o,a,x\}}\setminus A\subseteq\overline{\Aline oa\cup\{x\}}\setminus\Aline oa$, the $3$-Exchange Property of $X$ ensures that $x\in\overline{\Aline oa\cup\{y\}}\subseteq \overline{A\cup\{y\}}$.
\end{proof}

\begin{proposition}\label{p:reg<=>wreg+3-reg} A liner $X$ is regular if and only if it is weakly regular and $3$-regular.
\end{proposition}

\begin{proof} The ``only if'' part is trivial. To prove the ``if'' part, assume that a liner $X$ is weakly regular and $3$-regular. The regularity of $X$ will be established as soon as we show that for every flat $A\subset X$ and points $o\in A$, $b\in X\setminus A$, the equality $\overline{A\cup\{b\}}=\bigcup_{y\in\Aline ob}\Aline Ay$ holds. Since $X$ is weakly regular, for every $x\in \overline{A\cup \{b\}}$, there exists a point $a\in A$ such that $x\in \overline{\{a,o,b\}}$. Since $X$ is $3$-regular and $\|\Aline oa\|<3$, there exist points $y\in\Aline ob$ and $z\in \Aline oa\subseteq A$ such that $x\in \Aline zy\subseteq \Aline Ay$, witnessing that  $\overline{A\cup\{b\}}=\bigcup_{y\in\Aline ob}\Aline Ay$.
\end{proof}

The following characterization shows that regularity of liners can be be equivalently defined using a suitable first-order axiom.

\begin{theorem}\label{t:HA} A liner $(X,\Af)$ is regular if and only if it is $4$-regular if and only if it satisfies the {\sf Axiom of Regularity}
\begin{itemize}
\item[{\sf (AR)}] $\forall o,a,b,u,v,x,y,z\in X\;$\newline
$( \Af ovu\wedge\Af axv\wedge\Af byu\wedge \Af xzy)\;\Rightarrow\;\exists s,t,c,w\;(\Af aso\wedge\Af bto\wedge \Af sct\wedge\Af owu\wedge \Af czw)$.
\end{itemize}
\end{theorem}

\begin{proof} It is clear that every regular liner is $4$-regular. To prove that a $4$-regular liner $(X,\Af)$ satisfies the axiom {\sf (AR)}, take any points $o,a,b,u,v,x,z,y\in X$  with $ \Af ovu\wedge \Af axv\wedge\Af byu\wedge \Af xzy$. It follows that 
$z\in \Aline xy\subseteq \overline{\Aline av\cup\Aline bu}=\overline{\{a,b,u,v\}}\subseteq\overline{\{a,b,o,u\}}$. Consider the flat  $A\defeq\overline{\{a,o,b\}}$.

If $b\in\Aline oa$, then $A=\overline{\{a,o,b\}}=\Aline oa=\Aline{\Aline oa}{\Aline ob}$. If $b\notin\Aline oa$, then the 3-regularity of $(X,\Af)$ ensures that $A=\overline{\{a,o,b\}}=\Aline{\Aline oa}{\Aline ob}$. In both cases, we obtain that $A=\Aline{\Aline oa}{\Aline ob}$.

\begin{picture}(200,125)(-150,-50)

\put(0,0){\line(1,0){90}}
\put(0,0){\line(0,1){60}}
\put(0,0){\line(-1,-1){45}}
\put(0,60){\line(1,-1){60}}
\put(0,50){\color{red}\line(1,-2){34}}
\put(0,30){\line(-1,-2){30}}
\put(30,30){\line(-2,-1){40}}
\put(-45,-45){\color{red}\line(3,1){135}}

\put(0,0){\circle*{3}}
\put(-1,-8){$o$}
\put(0,30){\circle*{3}}
\put(-7,29){$v$}
\put(0,50){\color{red}\circle*{3}}
\put(-10,48){$w$}
\put(0,60){\circle*{3}}
\put(-6,63){$u$}
\put(-30,-30){\circle*{3}}
\put(-37,-30){$a$}
\put(30,30){\circle*{3}}
\put(32,32){$y$}
\put(90,0){\color{red}\circle*{3}}
\put(93,-2){$t$}
\put(60,0){\circle*{3}}
\put(61,2){$b$}
\put(-45,-45){\color{red}\circle*{3}}
\put(-52,-49){$s$}
\put(14,22){\color{blue}\circle*{3}}
\put(14,27){$z$}
\put(-10,10){\circle*{3}}
\put(-18,8){$x$}
\put(34.2,-18.6){\color{red}\circle*{3}}
\put(33,-28){$c$}
\end{picture}

If $u\in A$, then $x\in\Aline ua\subseteq A$, $y\in\Aline vb\subseteq \overline{\{o,u,b\}}\subseteq\overline{A\cup\{u\}}=A$ and $z\in\Aline xy\subseteq A$. Since $A=\Aline{\Aline oa}{\Aline ob}$ there exist points $s\in \Aline oa$ and $t\in\Aline ob$ such that $z\in \Aline st$. Then the points $s,t,c\defeq z$ and $w\defeq u$ have the required properties.

If $u\notin A$, then by the $4$-regularity, for the point $z\in\Aline xy\subseteq\overline{\{o,a,b,u\}}=\overline{A\cup\{u\}}$ there exist points $w\in\Aline ou$ and $c\in A$ such that $z\in\Aline wc$. Since $c\in A=\Aline{\Aline oa}{\Aline ob}$, there exist points $s\in\Aline oa$ and $t\in\Aline ob$ such that $c\in\Aline st$. Then the points $s,t,c,w$ witness that the axiom {\sf(AR)} is satisfied. 
\smallskip

Now assume that a liner $(X,\Af)$ satisfies the axiom {\sf (AR)}. To show that $(X,\Af)$ is regular, take any flat $A\subseteq X$ and points $o\in A$, $p\in X\setminus A$. Consider the line $L=\Aline op$. First we prove that the set $\Aline AL$ is flat, which means that $\Aline xy\subseteq \Aline AL$ for any points $x,y\in \Aline AL$. Fix any point $z\in\Aline xy$. Since $x,y\in\Aline AL$, there exist points $a,b\in A$  and $u,v\in L$ such that $\Af axv$ and $\Af byu$. 

If the set $\{o,v,u\}$ is a singleton, then  $\Af ovu$. 
If $\{o,v,u\}$ is not a singleton, then $o\ne u$ or $o\ne v$. Replacing the points $x,y$ by their places, if necessary, we can assume that $o\ne u$ and hence $v\in L=\Aline ou$. By the axiom {\sf (AR)}, there exist points $s,t,c,w\in X$ such that $\Af aso\wedge\Af bto\wedge \Af sct\wedge\Af owu \wedge \Af czw$. Then $c\in\Aline st\subseteq \Aline{\Aline ao}{\Aline ob}\subseteq \overline{\{a,o,b\}}=A$ and $w\in \Aline ou\subseteq  L$, witnessing the that  $z\in\Aline cw\in \Aline AL$ and the set $\Aline AL$ is flat. Then $\overline{A\cup \{p\}}=\Aline AL$ and for every point $q\in \overline{A\cup\{p\}}$, there exist points $\alpha\in A$ and $\lambda\in L=\Aline op$ such that $q\in\Aline \alpha\lambda$, witnessing that the liner $(X,\Af)$ is regular.
\end{proof}

\begin{corollary}\label{c:reg<=>||X||-reg} A liner $X$ is regular if and only if $X$ is $\|X\|$-regular.
\end{corollary}

\begin{proof} The ``only if'' part is trivial. To prove the ``if'' part, assume that a liner $X$ is $\|X\|$-regular. If $\|X\|\ge 4$, then $X$ is $4$-regular and regular, by Theorem~\ref{t:HA}. So, assume that $\|X\|<4$. To prove that $X$ is regular, fix any flat $A\subset X$ and points $o\in A$ and $b\in X\setminus A$. By Proposition~\ref{p:k-regular<=>2ex}, the $\|X\|$-regular liner $X$ is $\|X\|$-ranked. Since $A\ne X$, the $\|X\|$-rankedness of $X$ implies $\|A\|<\|X\|$. The $\|X\|$-regularity of $X$ implies the desired equality $\overline{A\cup \{b\}}=\bigcup_{y\in\Aline ob}\overline{Ay}$, witnessing that the liner $X$ is regular. 
\end{proof}  

\begin{exercise} Prove that a liner $X$ is regular if and only if  for every lines $A,B,C\subseteq X$ with $A\cap B\cap C\ne\varnothing$ the set $\Aline{\Aline AB}C$ is flat.
\end{exercise}

By analogy with Theorem~\ref{t:HA}, we can prove the following first-order characterization of the proregularity.

\begin{proposition}\label{p:proregular<=>} A liner $(X,\mathsf L)$ is proregular if and only if it is $4$-proregular if and only if it satisfies the 
{\sf Axiom of Proregularity:}  $\forall o,a,b,u,v,x,y,z,p\in X$\newline
$( \Af ovu\wedge\Af axv\wedge\Af byu\wedge \Af xzy\wedge \Af opu\wedge o\ne p\ne u)\Rightarrow\exists s,t,c,w\;(\Af aso\wedge\Af bto\wedge \Af sct\wedge\Af owu\wedge \Af czw)$.
%
\end{proposition}

Also the $3$-proregularity can be written by the following first-order formula.

\begin{proposition}\label{p:3-proregular<=>} A liner $(X,\mathsf L)$ is $3$-proregular if and only if its satisfies the\newline {\sf Axiom of 3-Proregularity:}  $\forall o,a,b,u,v,x,y,z,p\in X$\newline
$( \Af ovu\wedge\Af axu\wedge\Af byv\wedge \Af oba\wedge\Af xzy\wedge \Af opu\wedge o\ne p\ne u)\Rightarrow\exists s,w\;(\Af osa\wedge\Af owu\wedge \Af szw)$.
\end{proposition}

\begin{exercise} Write down the proof of Proposition~\ref{p:3-proregular<=>}.
\end{exercise}

\begin{exercise} Prove that a liner $X$ is weakly regular if and only if 
$$\forall o,x,y,z\in X\;\forall p\in \overline{\{o,x,z\}}\;\forall q\in\overline{\{o,y,z\}}\;\forall r\in\Aline pq\;\exists s\in\overline{\{x,o,y\}}\;\;(r\in\overline{\{o,s,z\}}).$$
\end{exercise}

\begin{exercise} Let $F:X\to Y$ be a liner isomorphism between two liners $X,Y$. Prove that the liner $X$ is regular (resp. strongly regular, weakly regular, proregular, $\kappa$-regular, $\kappa$-proregular) if and only if so is the liner $Y$.
\end{exercise}

\section{Parallelity Postulates and Axioms}

In this section we introduce many parallelity postulates and parallelity axioms that will be studied in the next sections. We start with four classical Parallelity Postulates.

\begin{definition}\label{d:PPBL} A liner $X$ is defined to be
\begin{itemize}
\item \index[person]{Proclus}\index{Proclus liner}\index{liner!Proclus}\defterm{Proclus}\footnote{{\bf Proclus Lycius} (412 -- 485), called Proclus the Successor ($\Pi\rho\acute o\kappa\lambda o\zeta\;\acute o\;\Delta\i\acute\alpha\delta o\chi o\zeta$), was a Greek Neoplatonist philosopher, one of the last major classical philosophers of late antiquity. He set forth one of the most elaborate and fully developed systems of Neoplatonism and, through later interpreters and translators, exerted an influence on Byzantine philosophy, Early Islamic philosophy, Scholastic philosophy, and German Idealism, especially Hegel, who called Proclus's Platonic Theology ``the true turning point or transition from ancient to modern times, from ancient philosophy to Christianity''. In ``{\em Commentary on Euclid's Elements}'' Proclus discusses the Euclides Parallel Postulate and suggests a simpler version of the postulate, known as the Proclus Parallelity Postulate.} 
if for every plane $P\subseteq X$, line $L\subseteq P$ and point $x\in P\setminus L$ there exists at most one line $\Lambda$ in $X$ such that $x\in \Lambda\subseteq P\setminus L$;
\item \index[person]{Playfair}\index{Playfair liner}\index{liner!Playfair}\defterm{Playfair}\footnote{{\bf John Playfair} (1748 -- 1819) was a professor of natural philosophy at the University of Edinburgh. Playfair   In 1783 he was a co-founder of the Royal Society of Edinburgh. He served as General Secretary to the society 1798–1819. Playfair's textbook {\em Elements of Geometry} made a brief expression of Euclid's parallel postulate known now as Playfair's axiom.} if for every plane $P\subseteq X$, line $L\subseteq P$ and point $x\in P\setminus L$ there exists a unique line $\Lambda$ in $X$ such that $x\in \Lambda\subseteq P\setminus L$;
\item \index[person]{Bolyai}\index{Bolyai liner}\index{liner!Bolyai}\defterm{Bolyai}\footnote{{\bf J\'anos Bolyai} (1802 -- 1860) was a Hungarian mathematician who developed absolute geometry -- a geometry that includes both Euclidean geometry and hyperbolic geometry. In the absolute geometry, the Bolyai Parallel Postulate is satisfied. The discovery of a consistent alternative geometry that might correspond to the structure of the universe helped to free mathematicians to study abstract concepts irrespective of any possible connection with the physical world.} if for every plane $P\subseteq X$, line $L\subseteq P$ and point $x\in P\setminus L$ there exists at least one line $\Lambda$ in $X$ such that $x\in \Lambda\subseteq P\setminus L$;
\item \index[person]{Lobachevsky}\index{liner!Lobachevsky}\index{Lobachevsky liner}\defterm{Lobachevsky}\footnote{{\bf Nikolai Ivanovich Lobachevsky} (1792 -- 1856) was a Russian mathematician and geometer, known primarily for his work on hyperbolic geometry, otherwise known as Lobachevskian geometry, and also for his fundamental study on Dirichlet integrals, known as the Lobachevsky integral formula. William Kingdon Clifford called Lobachevsky the ``Copernicus of Geometry'' due to the revolutionary character of his work.
}  if for every plane $P\subseteq X$, line $L\subseteq P$ and point $x\in P\setminus L$ there exist at least two distinct lines $\Lambda$ in $X$ such that $x\in \Lambda\subseteq P\setminus L$.
\end{itemize}
\end{definition}

These four definitions are partial cases of the following Parallel Postulates involving a cardinal parameter $\kappa$.

\begin{definition}\label{d:k-PPB} Let $\kappa$ be a cardinal. A liner $X$ is defined to be 
\begin{itemize}
\item \index{$\kappa$-hypoparallel liner}\index{liner!hypoparallel}\defterm{$\kappa$-hypoparallel} if for every plane $P\subseteq X$, line $L\subseteq P$ and point $x\in P\setminus L$ there exist at most $\kappa$ lines $\Lambda$ such that $x\in \Lambda\subseteq P\setminus L$;
\item \index{$\kappa$-parallel liner}\index{liner!$\kappa$-parallel}\defterm{$\kappa$-parallel} if for every plane $P\subseteq X$, line $L\subseteq P$ and point $x\in P\setminus L$ there exist exactly $\kappa$ lines $\Lambda$ such that $x\in \Lambda\subseteq P\setminus L$;
\item \index{$\kappa$-hyperparallel liner}\index{liner!$\kappa$-hyperparallel}\defterm{$\kappa$-hyperparallel} if for every plane $P\subseteq X$, line $L\subseteq X$ and point $x\in P\setminus L$ there exist at least $\kappa$  lines $\Lambda$ such that $x\in \Lambda\subseteq P\setminus L$.
\end{itemize}
\end{definition}

It is clear that a liner is $\kappa$-parallel if and only if it is both $\kappa$-hypoparallel and $\kappa$-hyperparallel.

Observe that a line $X$ is Proclus, Playfair, Bolyai, Lobachevsky if and only if it is $1$-hypoparallel, $1$-parallel, $1$-hyperparallel, $2$-hyperparallel, respectively.

 %
%

In Theorems~\ref{t:Proclus<=>} and \ref{t:Playfair<=>}, we shall prove that  Proclus and Playfair Parallel Postulates are first-order properties of liners, whose first-order characterizations involve the following parallelity axioms.

\begin{definition}\label{d:PP} A liner $X$ is defined to be
\begin{itemize}
\item \index{projective liner}\index{liner!projective}\defterm{projective} if $\forall o,x,y\in X\;\;\forall p\in\Aline xy\;\;\forall v\in\Aline oy\setminus\{p\}\;\;(\Aline vp\cap\Aline ox\ne \varnothing)$;
\item \index{proaffine liner}\index{liner!proaffine}\defterm{proaffine} if $\forall o,x,y\in X\;\forall p\in\Aline xy\setminus\Aline ox\;\exists u\in\Aline oy\;\forall v\in\Aline oy\setminus \{u\}\;\;(\Aline vp\cap\Aline ox\ne \varnothing)$;
\item \index{affine liner}\index{liner!affine}\defterm{affine} if $\forall o,x,y\in X\;\forall p\in\Aline xy\setminus\Aline ox\;\exists u\in\Aline oy\;\forall v\in\Aline oy\;\;(u=v\;\Leftrightarrow\;\Aline vp\cap\Aline ox=\overline \varnothing)$;
\item \index{hyperaffine liner}\index{liner!hyperaffine}\defterm{hyperaffine} if $\forall o,x,y\in X\;\forall p\in\Aline xy\setminus\Aline ox\;\exists u\in\Aline oy\;\;\;(\Aline up\cap\Aline ox= \varnothing)$;
\item \index{hyperbolic liner}\index{liner!hyperbolic}\defterm{hyperbolic} if $\forall o,x,y\in X\,\forall p\in\Aline xy\setminus(\Aline ox\,{\cup}\,\Aline oy)\,\exists u,\!v{\in}\Aline oy\;(u\ne v\,\wedge\,\Aline up\cap\Aline ox{=}\varnothing{=}\Aline up\cap\Aline ox)$;
\item \index{injective liner}\index{liner!injective}\defterm{injective} if $\forall o,x,y\in X\;\forall p\in\Aline xy\setminus(\Aline ox\,{\cup}\,\Aline oy)\;\forall u\in\Aline oy\setminus\{o,y\}\;\;(\Aline up\cap\Aline ox=\varnothing)$.
\end{itemize} The first-order formulas defining projective, proaffine, affine, hyperaffine, and hyperbolic liners will be called the \index{Axiom!projectivity}\index{projectivity axiom}\defterm{projectivity}, \index{Axiom!proaffinity}\index{proaffinity axiom}\defterm{proaffinity}, \index{Axiom!affinity}\index{affinity axiom}\defterm{affinity}, \index{Axiom!hyperaffinity}\index{hyperaffinity axiom}\defterm{hyperaffinity}, \index{Axiom!hyperbolocity}\index{hyperbolicity axiom}\defterm{hyperbolicity}, and \index{Axiom!injectivity}\index{injectivity axiom}\defterm{injectivity axioms}, respectively. 
\end{definition}

\begin{remark} The proaffinity, affinity, hyperaffinity, and hyperbolicity axioms are first order counterparts of the Parallel Postulates of Proclus, Playfair, Bolyai, and Lobachevsky, respectively. The projectivity axiom has been explicitely formulated  \index[person]{Veblen}Veblen \footnote{{\bf Oswald Veblen} (1880 -- 1960) was an American mathematician, geometer and topologist, whose work found application in atomic physics and the theory of relativity. He proved the Jordan curve theorem in 1905. He suggested an elegant axiomatization of projective geometry in the paper [O.Veblen, J.Young, {\em A Set of Assumptions for Projective Geometry}, American J. Math. 30{:4} (1908), 347-380].} in his axiom system for Projective Geometry. The injectivity axiom (as an antonym to the projectivity) is well-known in the theory of unitals, as the absence of O'Nan configurations (known also as Pasch configurations), see Section~\ref{s:hyperbolic} for more information on the classical unitals.
\end{remark}

\begin{exercise} Draw illustrations to the parallelity postulates introduced in Definitions~\ref{d:PPBL} and \ref{d:k-PPB}, \ref{d:PP}.
\end{exercise}

\begin{exercise} Prove that every $3$-long  injective liner is hyperaffine, and every $4$-long injective liner is hyperbolic.
\end{exercise}

\begin{exercise} Prove the following implications between various parallelity postulates and  axioms.
$$
\xymatrix{
\mbox{$0$-hypoparallel}\ar@{=>}[r]\ar@{<=>}[d]&\mbox{Proclus}\ar@{<=>}[d]&\mbox{Playfair}\ar@{=>}[l]\ar@{<=>}[d]\ar@{=>}[r]&
\mbox{Bolyai}\ar@{<=>}[d]&\mbox{Lobachevsky}\ar@{=>}[l]\ar@{<=>}[d]\\
\mbox{$0$-parallel}\ar@{=>}[r]\ar@{<=>}[d]&\mbox{$1$-hypoparallel}\ar@{=>}[d]&\mbox{$1$-parallel}\ar@{=>}[l]\ar@{=>}[d]\ar@{=>}[r]&
\mbox{$1$-hyperparallel}&\mbox{$2$-hyperparallel}\ar@{=>}[l]\\
\mbox{projective}\ar@{=>}[r]&\mbox{proaffine}&\mbox{affine}\ar@{=>}[l]\ar@{=>}[r]&\mbox{hyperaffine}\ar@{=>}[u]&\mbox{hyperbolic}\ar@{=>}[l]\ar@{=>}[u]
}
$$

{\em Hint:} For non-trivial implications, see Theorems~\ref{t:projective<=>}, \ref{t:Proclus<=>}.
\end{exercise}

\begin{proposition}\label{p:k-parallel=>3-ranked} If a liner $X$ is $\kappa$-parallel for some finite cardinal $\kappa$, then $X$ is $3$-ranked.
\end{proposition}

\begin{proof} Assuming that $X$ is not $3$-ranked, we can find two planes $P\subset \Pi$ in $X$ such that $P\ne \Pi$. Fix any line $\Lambda\subset P$ and choose a point $x\in P\setminus \Lambda$. Let $\mathcal L$ be the family of lines in the liner $X$. Since $X$ is $\kappa$-parallel, for any plane $F\in\{P,\Pi\}$, the family $\mathcal L_F\defeq\{L\in\mathcal L:x\in L\subseteq F\setminus \Lambda\}$ has cardinality $\kappa$. The inclusion $P\subseteq\Pi$ implies $\mathcal L_P\subseteq\mathcal L_\Pi$ and hence $\mathcal L_P=\mathcal L_\Pi$ as both families $\mathcal L_P$ and $\mathcal L_\Pi$ have the same finite cardinality $\kappa$. On the other hand, for every point $y\in \Pi\setminus P$, we have $P\cap\Aline xy=\{x\}$ and hence $\Aline xy\in\mathcal L_\Pi\setminus \mathcal L_P$, which is a contradiction showing that the $\kappa$-parallel liner $X$ is $3$-ranked.
\end{proof}
\begin{proposition} A liner $X$ is injective if and only if among any four distinct lines in $X$ two are disjoint.
\end{proposition}

\begin{proof} To prove the ``if'' part, assume among any four lines in $X$, two are disjoint. To prove that $X$ is injective, take any points $o,x,y\in X$, $p\in\Aline xy\setminus(\Aline ox\cup\Aline oy)$ and $u\in\Aline oy\setminus\{o,y\}$. Then the lines $\Aline ox,\Aline oy,\Aline xy$ are pairwise concurrent. Since $u\in \Aline up\cap\Aline oy$ and $p\in\Aline up\cap\Aline xy$, the line $\Aline up$ is disjoint with the line $\Aline xy$, witnessing that the liner $X$ is injective.
\smallskip

 To prove the ``only if'' part, assume that $X$ contains four pairwise concurrent lines  $L_1,L_2,L_3,L_4$. Then there exist unique points $o,x,y,p,u,v\in X$ such that $\{o\}=L_1\cap L_2$, $\{x\}=L_1\cap L_3$, $\{y\}=L_1\cap L_2$, $\{p\}=L_3\cap L_4$, $\{u\}=L_2\cap L_4$ and $\{v\}=L_1\cap L_4$. Then $p\in L_3\setminus(L_1\cup L_2)=\Aline xy\setminus(\Aline ox\cup\Aline oy)$,$u\in L_2\setminus(L_1\cup L_3)=\Aline oy\setminus\{o,y\}$ and $v\in L_4\cap L_1=\Aline up\cap \Aline ox$, witnessing that $X$ is not injective. 
\end{proof}



\begin{exercise}\label{ex:Beltrami-Klein} Show that for every $n\ge 2$, the subliner $B=\{x\in \IR^n:\sum_{i\in n}x_i^2<1\}$ of the Euclidean space $\IR^n$ is a hyperbolic liner. This hyperbolic liner is called the \index[person]{Beltrami}\index[person]{Klein}\defterm{Beltrami\footnote{{\bf Eugenio Beltrami} (1835 -- 1900) was an Italian mathematician notable for his work concerning differential geometry and mathematical physics. His work was noted especially for clarity of exposition. He was the first to prove consistency of non-Euclidean geometry by modeling it on a surface of constant curvature, the pseudosphere, and in the interior of an $n$-dimensional unit sphere, the so-called Beltrami–Klein model.}--Klein\footnote{{\bf Christian Felix Klein}  (1849 -- 1925) was a German mathematician and mathematics educator, known for his work with group theory, complex analysis, non-Euclidean geometry, and on the associations between geometry and group theory. His 1872 Erlangen program, classifying geometries by their basic symmetry groups, was an influential synthesis of much of the mathematics of the time.} model}\index{Beltrami--Klein model} of hyperbolic geometry.
\end{exercise}
 
\begin{exercise} Find an example of a proaffine liner which is neither projective nor affine.

{\em Hint:} Remove a point from a suitable projective liner.
\end{exercise}

\begin{exercise} Find an example of a proaffine  liner which is not $3$-regular.
\smallskip

{\em Hint:} Consider the following 6-element subliner of the Euclidean plane:

\begin{picture}(100,80)(-130,-7)

\put(0,0){\line(1,0){120}}
\put(60,60){\line(1,-1){60}}
\put(60,60){\line(-1,-1){60}}

\put(0,0){\circle*{3}}
\put(60,60){\circle*{3}}
\put(30,30){\circle*{3}}
\put(90,30){\circle*{3}}
\put(60,0){\circle*{3}}
\put(120,0){\circle*{3}}
\end{picture}
\end{exercise}

\begin{exercise} Find an example of an affine liner which is not $3$-regular.
\smallskip

{\em Hint:} Look at Example~\ref{ex:Z15}.
\end{exercise}

\begin{exercise} Find an example of an affine liner which is $3$-regular but not $4$-regular.

{\em Hint:} Look at Example~\ref{ex:HTS}.
\end{exercise}


\begin{exercise} Find an example of an affine ranked plane, which is not $3$-regular.

{\em Hint:} Look at the $13$-element liner in Example~\ref{ex:Steiner13}.
\end{exercise}

\begin{exercise} Let $F:X\to Y$ be a liner isomorphism between liners $X,Y$. Prove that the liner $X$ is  projective (resp. proaffine, affine, hyperaffine, hyperbolic, injective, Proclus, Playfair, Bolyai, Lobachevsky, $\kappa$-hypoparallel, $\kappa$-parallel, $\kappa$-hyperparallel) if and only if so is the liner $Y$.
\end{exercise}

\section{Balanced liners}

\begin{definition}\label{d:balanced} Given two cardinals $\kappa,\mu$, we say that a liner $X$ is \index{$\kappa$-balanced liner}\index{liner!$\kappa$-balanced}\defterm{$\kappa$-balanced} with 
\index[note]{$\vert X\vert_{\kappa}$}
$|X|_\kappa=\mu$ if every flat $A\subseteq X$ of rank $\|A\|=\kappa$ has cardinality $|A|=\mu$. A liner $X$ is \index{balanced liner}\index{liner!balanced}\defterm{balanced} if it is $\kappa$-balanced for every cardinal $\kappa$.
\end{definition}

\begin{remark} For a $\kappa$-balanced liner $X$, the cardinal $|X|_\kappa$ is well-defined if $X$ contains a flat $F\subseteq X$ of rank $\|F\|=\kappa$. In this case $|X|_\kappa=|F|$ for every flat $F\subseteq X$ of rank $\|F\|=\kappa$. 
\end{remark} 

It is clear that every liner $X$ is $1$-balanced with $|X|_1=1$. A liner $X$ is $2$-balanced if and only if all lines have the same cardinality $|X|_2$.  

\begin{theorem}\label{t:2-balanced-dependence} Let $X$ be a $2$-balanced liner and $\mathcal L$ be the family of all lines in $X$. If $|X|_2<|X|$, then for every point $x\in X$, the cardinality of the family $\mathcal L_x\defeq\{L\in\mathcal L:x\in L\}$ satisfies two equations
$$|X|-1=(|X|_2-1)\cdot |\mathcal L_x|\quad\mbox{and}\quad |X|_2\cdot|\mathcal L|=|X|\cdot|\mathcal L_x|$$and hence $|\mathcal L_x|$ does not depend on $x$. If $X$ is finite, then $|X|_2$ divides the number $(|\mathcal L_x|-1)\cdot|\mathcal L_x|$.
\end{theorem} 

\begin{proof} Fix any point $x\in X$. Since $$X=\{x\}\cup\bigcup_{L\in\mathcal L_x}(L\setminus\{x\})$$ and $L\cap\Lambda=\{x\}$ for any distinct lines $L,\Lambda\in\mathcal L_x$, we have
$$|X|=1+\sum_{L\in\mathcal L_x}|L\setminus\{x\}|=1+\sum_{L\in\mathcal L_x}(|X_2|-1)=1+|\mathcal L_x|\cdot(|X_2|-1).$$
Since $|X|_2<|X|$, the cardinal $|\mathcal L_x|$ is uniquely determined by the equation $|X|-1=|\mathcal L_x|\cdot(|X|_2-1)$ and hence it does not depend on $x$.

Next, consider the set $\mathcal P\defeq\{(x,L)\in X\times\mathcal L:x\in L\}$ and its projections 
$$\pr_X:\mathcal P\to X,\;\;\pr_X:(x,L)\mapsto x,\quad\mbox{and}\quad 
\pr_{\mathcal L}:\mathcal P\to \mathcal L,\;\;\pr_{\mathcal L}:(x,L)\mapsto L.$$ Observe that for every $x\in X$ and $L\in\mathcal L$ we have $|\pr_X^{-1}(x)|=|\mathcal L_x|$ and $|\pr_{\mathcal L}^{-1}(L)|=|L|=|X|_2-1$. Then
$$|\mathcal P|=\sum_{x\in X}|\pr_X^{-1}(x)|=\sum_{x\in X}|\mathcal L_x|=|X|\cdot|\mathcal L_x|.$$
On the other hand, 
$$|\mathcal P|=\sum_{L\in\mathcal L}|\pr_{\mathcal L}^{-1}(L)|=\sum_{L\in\mathcal L}|L|=|\mathcal L|\cdot|X|_2.$$
Therefore, $|X|_2\cdot|\mathcal L|=|\mathcal P|=|X|\cdot|\mathcal L_x|$.
\smallskip

Now assume that $X$ is finite. Write the number $|X|_2$ as $|X|_2=b\cdot d$ where $d$ is the largest common divisor of the numbers $|X|_2$ and $|\mathcal L_x|$ where $x$ is any point of $X$. Since $b$ is coprime with $|\mathcal L_x|$ and $b$ divides $|X|_2\cdot|\mathcal L|=|X|\cdot|\mathcal L_x|$, it divides $|X|=1+(|X|_2-1)\cdot|\mathcal L_x|$ and hence $b$ divides $|\mathcal L_x|-1$. Then $|X|_2=b\cdot d$ divides $(|\mathcal L_x|-1)\cdot|\mathcal L_x|$.
\end{proof}

\begin{remark} Theorem~\ref{t:2-balanced-dependence} implies that the cardinality $v\defeq|X|$ of any finite $2$-balanced liner $X$ with $k\defeq |X|_2<|X|$ satisfies the equality $v=1+(k-1)\cdot r$ for some number $r\ge k$ such that $k$ divides $(r-1){\cdot}r$ (denoted by $k\,|\,(r-1)r)$. The latter divisibility condition implies (and in fact, is equivalent to) the divisibilities $(k-1)\,|\,(v-1)$ and $(k-1)k\,|\,(v-1)v$. Pairs of positive integer numbers $(k,v)$ satisfying the divisibility conditions $$(k-1)\,|\,(v-1)\quad\mbox{and}\quad(k-1)k\,|\,(v-1)v$$ are called \defterm{admissible}. Therefore, for every finite $2$-balanced liner $X$, the pair $(|X|_2,|X|)$ is admissible. The admissibility of a pair $(k,v)$ does not imply the existence of a $2$-balanced liner $X$ with $(|X|_2,|X|)=(k,v)$. 
For example, the pairs $(6,36)$ and $(7,43)$ are admissible, but liners $X$ with  $(|X|_2,|X|)\in\{(6,36),(7,43)\}$ do not exist, by Corollary~\ref{c:no6order} of the Bruck--Ryser's Theorem~\ref{t:Bruck-Ryser}. The pair $(6,46)$ is admissible but no liner $X$ with $(|X|_2,|X|)=(6,46)$ exists, by a result of Houghten, Thiel, Janssen,  and \index[person]{Lam}Lam\footnote{{\bf Clement Wing Hong Lam} is a Canadian mathematician  specializing in combinatorics.  He is famous for the computer proof, with Larry Thiel and S. Swiercz, of the nonexistence of a finite projective plane of order 10. Lam earned his PhD in 1974 under Herbert Ryser at Caltech with thesis ``{\em Rational G-Circulants Satisfying the Matrix Equation $A^2=dI+\lambda J$}''. He is a Professor Emeritus in the Department of Computer Science and Software Engineering at Concordia University in Montreal, Quebec. In 1992 he received the Lester Randolph Ford Award for the article ``{\em The search for a finite projective plane of order $10$}''. In 2006 he received the Euler medal.}  \cite{HTJL2001}, proved by a computer search. 
 On the other hand, the admissibility condition is asymptotically sufficient for the existence of finite $2$-balanced liners, according to the following fundamental theorem, proved by \index[person]{Wilson}Wilson
\footnote{{\bf Richard Michael Wilson} (born in 1945) is a mathematician and a professor emeritus at the California Institute of Technology. 
Wilson was educated at Indiana University where he was awarded a Bachelor of Science degree in 1966 followed by a Master of Science degree from Ohio State University in 1968. His PhD, also from Ohio State University was awarded in 1969 for research ``An Existence Theory For Pairwise Balanced Designs'' supervised by Dijen K. Ray-Chaudhuri. His breakthrough in pairwise balanced designs, and orthogonal Latin squares  built upon the groundwork set before him, by R. C. Bose, E. T. Parker, S. S. Shrikhande, and Haim Hanani are widely referenced in combinatorial design theory and coding theory. 
Wilson is also one of the world's top experts on historical flutes.} 
 in \cite{Wilson1}, \cite{Wilson2}, \cite{Wilson3}.
\end{remark}

\begin{theorems}[Wilson, 1975]\label{t:Wilson} For every integer numbers $k$ there exists a number $v_k$ such that for every admissible pair $(k,v)$ with $v>v_k$, there exists a $2$-balanced liner $X$ such that $|X|_2=k$ and $|X|=v$. 
\end{theorems}

\begin{remark} We can assume that the number $v_k$ in Wilson's Theorem~\ref{t:Wilson} is the smallest possible. In this case it is determined uniquely. The Wilson proof of Theorem~\ref{t:Wilson} does not provide an explicit formula for finding the number $v_k$, it just claims that it does exist. However, for $k\le 9$ the following upper and lower bounds for the numbers $v_k$ are known (see \cite[\S II.3.1]{HCD}):
$$
\begin{array}{|c|cccccccc|}
\hline
k&2&3&4&5&6&7&8&9\\
\hline
v_k\ge&2&3&4&5&46&43&8&9\\
v_k\le &2&3&4&5&801&2605&3753&16497\\
\hline
\end{array}
$$As we have already mentioned, liners $X$ with $|X|_2=6$ and $X\in\{36,46\}$ do not exist. The existence of liners $X$ with $|X|_2=6$ is undecided for the following values $v$ (see Table 3.4 in \cite[II.3.1]{HCD} and \cite{BHR2025}): 51, 61, 81, 166, 231, 256, 261, 286, 316, 321, 346, 351, 376, 406, 411, 436, 471, 501, 561, 591,
616, 646, 651, 676, 771, 796, 801.
\end{remark}
 
By analogy with Theorem~\ref{t:2-balanced-dependence} we can prove the following theorem describing the structure of planes in $2$-balanced $\kappa$-parallel liners.

\begin{theorem}\label{t:2-balance+k-parallel=>3-balance} If a liner $X$ is $2$-balanced and $\kappa$-parallel for some cardinal $\kappa$, then $X$ is $3$-balanced with $|X|_3=1+(\kappa+|X|_2)(|X|_2-1)$. 
\end{theorem}

\begin{proof} Given any plane $P\subseteq X$, we should show that $|P|=1+(\kappa+|X|_2)(|X|_2-1)$. Fix any  line $\Lambda\subseteq P$ and point $x\in P\setminus\Lambda$. Consider the family $\mathcal L_{x}$ of all lines in the plane $P$ that contain the point $x$, and observe that 
$$P={\textstyle\bigcup}\mathcal L_{x}=\{x\}\cup\bigcup_{L\in\mathcal L_{x}}(L\setminus\{x\})\quad\mbox{and hence}\quad|P|=1+\sum_{L\in\mathcal L_{x}}|L\setminus\{x\}|.$$
 The family $\mathcal L_{x}$ is the union of two subfamilies
$$\mathcal L_{x,\Lambda}^\circ\defeq\{L\in\mathcal L_{x}:L\cap\Lambda=\varnothing\}\quad\mbox{and}\quad\mathcal L_{x,\Lambda}^\bullet\defeq\{L\in\mathcal L_{x}:L\cap\Lambda\ne\varnothing\}=\{\Aline xy:y\in\Lambda\}.$$
Since $X$ is $\kappa$-parallel and $2$-balanced, $|\mathcal L_{x}|=|\mathcal L_{x,\Lambda}^\circ|+|\mathcal L_{x,\Lambda}^\bullet|=\kappa+|\Lambda|=\kappa+|X|_2$.
Then
$$|P|=1+\sum_{L\in\mathcal L_{x}}|L\setminus\{x\}|=1+|\mathcal L_{x}|(|X|_2-1)=1+(\kappa+|X|_2)(|X|_2-1),$$witnessing that the liner $X$ is $3$-balanced with $|X|_3=1+(\kappa+|X|_2)(|X|_2-1)$.
%
\end{proof}

\begin{proposition}\label{p:23-balance=>k-parallel} If a liner $X$ of rank $\|X\|\ge 3$ is $2$-balanced and $3$-balanced, then there exists a unique cardinal $\lambda\ge|X|_2$ such that $|X|_3-1=\lambda\cdot(|X|_2-1)$. If $|X|_2<|X|_3$, then $X$ is $\kappa$-parallel for the unique cardinal $\kappa$ such that $\lambda=\kappa+|X|_2$. If $|X|_3$ is finite, then $|X|_2$ divides the numbers $(\lambda-1)\lambda$ and $(\kappa-1)\kappa$.
\end{proposition} 

\begin{proof} Assume that a liner $X$ is $2$-balanced and $3$-balanced. Since $\|X\|\ge 3$, there exists a plane $P\subseteq X$. Fix any line $\Lambda\subseteq P$ and point $x\in P\setminus\Lambda$.  Let  $\mathcal L_x$ be the family of lines in the plane $P$ that contain the point $x$. The family $\mathcal L_x$ is the union of two subfamilies $$\mathcal L_{x,\Lambda}^\circ\defeq\{L\in\mathcal L_{x}:L\cap\Lambda=\varnothing\}\quad\mbox{and}\quad\mathcal L_{x,\Lambda}^\bullet\defeq\{L\in\mathcal L_{x}:L\cap\Lambda\ne\varnothing\}=\{\Aline xy:y\in\Lambda\}.$$ Consider the cardinal $\lambda\defeq|\mathcal L_x|$ and observe that $$\lambda=|\mathcal L_x|=|\mathcal L_{x,\Lambda}^\circ|+|\mathcal L_{x,\Lambda}^\bullet|=|\mathcal L_{x,\Lambda}^\circ|+|\Lambda|=|\mathcal L_{x,\Lambda}^\circ|+|X|_2\ge|X|_2$$and $|X|_3=|P|=|\textstyle{\bigcup}\mathcal L_x|=\big|\{x\}\cup\bigcup_{L\in\mathcal L_{x}}L\setminus\{x\}\big|=1+|\mathcal L_x|(|X|_2-1)=1+\lambda(|X|_2-1)$. 
If the cardinal $|X|_3$ is finite, then the equation $|X|_3 =1+\lambda(|X|_2-1)$ uniquely determines the (finite) cardinal $\lambda$. If $|X|_3$ is infinite and $|X|_2=|X_3|$, then $|X|_2\le\lambda\le 1+\lambda(|X|_2-1)=|X|_3$ implies $\lambda=|X_2|=|X|_3$. If $|X|_3$ is infinite and $|X|_2<|X|_3$, then $|X|_3=1+\lambda(|X|_2-1)=\max\{\lambda,|X|_2-1\}>|X|_2$ implies $\lambda=|X|_3$. Therefore, if the cardinal $|X|_3$ is infinite, then the cardinal $\lambda$ equals $|X|_3$ and hence is uniquely determined.

Now assume that $|X|_2<|X|_3$. If $|X|_3$ is finite, then so is the cardinal $\lambda\ge|X|_2$ and we can find a unique cardinal $\kappa$ such that $\lambda=\kappa+|X|_2$. If $|X|_3$ is infinite, then $\lambda=|X|_3$ and $\kappa=\lambda=|X|_3$ is a unique cardinal such that $\lambda=\kappa+|X|_2$. Since $\kappa+|X|_2=\lambda=|\mathcal L_x|=|\mathcal L_{x,\Lambda}^\circ|+|X|_2$, 
the uniqueness of the cardinal $\kappa$ implies $|\mathcal L_{x,\Lambda}^\circ|=\kappa$, which means that the liner $X$ is $\kappa$-parallel.
\smallskip

If $|X|_3$ is finite, then we can apply Theorem~\ref{t:2-balanced-dependence} and conclude that $|X|_2$ divides the numbers $(\lambda-1)\lambda$ and $(\kappa-1)\kappa$. 
\end{proof}

\begin{theorem}\label{t:k-parallel=>2-balance} Assume that a liner $X$ is $\kappa$-parallel for some finite cardinal $\kappa$. If $\kappa>0$ or $X$ is $3$-long, then $X$ is $2$-balanced, $3$-balanced and $|X|_3=1+(|X|_2-1)(\kappa+|X|_2)$. If $|X|_2$ is finite, then $|X|_2$ divides $(\kappa-1)\kappa$.
\end{theorem}

\begin{proof} To see that $X$ is $2$-balanced, it suffices to check that $|\Lambda|=|\Lambda'|$ for any concurrent lines $\Lambda,\Lambda'$ in $X$. In this case, the flat hull $P\defeq\overline{\Lambda\cup\Lambda'}$ is a plane. Since the lines $\Lambda,\Lambda'$ are distinct, there exist points $x\in\Lambda\setminus\Lambda'$ and $y\in\Lambda'\setminus\Lambda$.

\begin{claim} There exists a point $z\in P\setminus(\Lambda\cup\Lambda')$. 
\end{claim}

\begin{proof} If $\kappa>0$, then the plane $P$ contains a line $L_x$ such that $x\in L_x$ and $L_x\cap \Lambda'=\varnothing$ and hence $L_x\ne \Lambda$ and $L_x\cap\Lambda=\{x\}$. Choose any point $z\in L_x\setminus\{x\}$ and observe that $z\in P\setminus(\Lambda\cup\Lambda')$. 

If $X$ is $3$-long, then the line $\Aline xy$ contains some point $z\in\Aline xy\setminus\{x,y\}$. 
Assuming that $z\in\Lambda$, we conclude that $y\in\Aline xy=\Aline xz\subseteq\Lambda$, which contradicts the choice of $y$. This contradiction shows that $z\notin\Lambda$. By analogy we can prove that $z\notin\Lambda'$.
\end{proof}

Let $\mathcal L_z$ be the family of all lines in the plane $P=\overline{L\cup\{x\}}$ that contain the point $z$.  Observe that $$\mathcal L_z=\mathcal L_{z,\Lambda}^\circ\cup\mathcal L_{z,\Lambda}^\bullet=\mathcal L_{z,\Lambda'}^\circ\cup\mathcal L_{z,\Lambda'}^\bullet,$$where
 $$
\begin{aligned}
\mathcal L_{z,\Lambda}^\circ&\defeq\{L\in\mathcal L_{z}:L\cap\Lambda=\varnothing\},&\mathcal L_{z,\Lambda}^\bullet&\defeq\{L\in\mathcal L_{z}:L\cap\Lambda\ne\varnothing\}=\{\Aline zp:p\in\Lambda\},\\
\mathcal L_{z,\Lambda'}^\circ\!&\defeq\{L\in\mathcal L_{z}:L\cap\Lambda'=\varnothing\},&\mathcal L_{z,\Lambda'}^\bullet\!&\defeq\{L\in\mathcal L_{z}:L\cap\Lambda'\ne\varnothing\}=\{\Aline zp:p\in\Lambda'\}.
\end{aligned}
$$
Since $X$ is $\kappa$-parallel,
$|\mathcal L_z|=|\mathcal L_{z,\Lambda}^\circ|+|\mathcal L_{z,\Lambda}^\bullet|=\kappa+|\Lambda|$ and $|\mathcal L_z|=\kappa+|\Lambda'|$. Therefore,
$$\kappa+|\Lambda|=|\mathcal L_z|=\kappa+|\Lambda'|.$$
Since the cardinal $\kappa$ is finite, the equality $\kappa+|\Lambda|=\kappa+|\Lambda'|$ implies $|\Lambda|=|\Lambda'|$, witnessing that the liner $X$ is $2$-balanced. By Theorem~\ref{t:2-balance+k-parallel=>3-balance}, the $2$-balanced $\kappa$-parallel liner is $3$-balanced and $|X|_3=1+(|X|_2-1)(\kappa+|X|_2)$.

If $|X|_2$ is finite, then $|X|_3=1+(|X|_2-1)(\kappa+|X|_2)$ is finite and $|X|_2$ divides $(\kappa-1)\kappa$, by Proposition~\ref{p:23-balance=>k-parallel}.
\end{proof}

Theorems~\ref{t:2-balance+k-parallel=>3-balance} and Proposition~\ref{p:23-balance=>k-parallel} imply the following characterization.


\begin{corollary} For a $3$-long finite liner $X$ of rank $\|X\|\ge 3$, the following conditions are equivalent:
\begin{enumerate}
\item $X$ is $\kappa$-parallel for some cardinal $\kappa$;
\item $X$ is $2$-balanced and $3$-balanced;
\item $X$ is $2$-balanced, $3$-balanced and $\kappa$-parallel for a unique cardinal $\kappa$  satisfying the equation $|X|_3-1=(\kappa+|X|_2)(|X|_2-1)$.
\end{enumerate}
\end{corollary}

\begin{exercise} Give an example of a $0$-parallel liner which is $3$-balanced but not $2$-balanced.
\end{exercise}

\begin{proposition}\label{p:2-balanced=>w-balanced} If a liner $X$ is $2$-balanced, then for every infinite cardinal $\kappa$, the liner $X$ is $\kappa$-balanced with $|X|_\kappa=\max\{\kappa,|X|_2\}$.
\end{proposition}

\begin{proof} First we show that for every finite set $A\subseteq X$, its flat hull $\overline{A}$ has cardinality $|A|\le\max\{\w,|X|_2\}$. Consider the sequence $(A_n)_{n\in\w}$ defined by the recursive formula: $$A_0=A\quad\mbox{and}\quad A_{n+1}=\Aline {A_n}{A_n}\quad\mbox{for $n\in\w$}.$$  It is easy to see that $\overline A=\bigcup_{n\in\w}A_n$. Observe that $|A_0|=|A|<\w\le\max\{\w,|X|_2\}$. Assume that for some $n\in\w$ the upper bound $|A_n|\le\max\{\w,|X|_2\}$ has been proved. Then $$|A_{n+1}|\le\sum_{x,y\in A_n}|\Aline xy|\le|A_n|\cdot|A_n|\cdot |X|_2\le\max\{\w,|X|_2\}.$$
Then $|A_n|\le\max\{\w,|X|_2\}$ for all $n\in\w$ and hence $|\overline A|\le\big|\bigcup_{n\in\w}A_n\big|=\sup_{n\in\w}|A_n|\le\max\{\w,|X|_2\}$. 

Now take any flat $I\subseteq X$ of infinite rank $\kappa$ and find a set $A\subseteq X$ of cardinality $|A|=\kappa$ such that $I\subseteq\overline A$. The family $[A]^{<\w}$ of finite subsets of the infinite set $A$ has cardinality $|[A]^{<\w}|=|A|=\kappa$. By Proposition~\ref{p:aff-finitary}, $\overline A=\bigcup_{F\in[A]^{<\w}}\overline F$ and hence
$$|I|\le|\overline A|\le\sum_{F\in[A]^{<\w}}|\overline F|\le|[A]^{<\w}|\cdot\max\{\w,|X|_2\}=\kappa\cdot\max\{\w,|X|_2\}=\max\{\kappa,|X|_2\}.$$
Since $I$ if a flat of infinite rank, it contains a line and hence $|I|\ge|X_2|$. Also $\kappa=\|I\|\le|I|$ and hence $\max\{\kappa,|X|_2\}\le|I|\le\max\{n,|X|_2\}$, witnessing that $|X|_\kappa=\max\{\kappa,|X|_2\}$.
\end{proof}

\begin{theorem}\label{t:wr+k-parallel=>n-balanced} If a weakly regular $3$-ranked liner $X$ is $2$-balanced and $\kappa$-parallel for some cardinal $\kappa$, then $X$ is balanced with$$|X|_{n}=1+(|X|_2-1)\sum_{r=0}^{n-2}(\kappa+|X|_2-1)^r$$
for every finite cardinal $n\ge 2$.
\end{theorem}

\begin{proof} By Proposition~\ref{p:2-balanced=>w-balanced}, the $2$-balanced liner $X$ is $\kappa$-balanced for every infinite cardinal $\kappa$. It remains to prove that every flat $F\subseteq X$ of finite rank $\|F\|=n\ge 2$ has cardinality $|F|=1+(|X|_2-1)\sum_{r=0}^{n-2}(\kappa+|X|_2-1)^r$.  This equality will be proved by induction on $n$. If $n=2$, then $F$ is a line and $$|F|=|X|_2=1+(|X|_2-1)\sum_{r=0}^{n-2}(\kappa+|X|_2-1)^r.$$ If $n=3$, then $F$ is a plane and 
$$|F|=|X|_3=1+(\kappa+|X|_2)(|X|_2-1)=1+(|X|_2-1)\sum_{r=0}^{n-2}(\kappa+|X|_2-1)^r,$$
by Theorem~\ref{t:2-balance+k-parallel=>3-balance}.

Assume that for some finite cardinal $n\ge 3$, the equality  $$|F|=1+(|X|_2-1)\sum_{r=0}^{\|F\|-2}(\kappa+|X|_2-1)^r$$ has been proved for all flats $F\subseteq X$ of rank $\|F\|=n$. Take any flat $F\subseteq X$ of rank $\|F\|=n+1$. By Theorem~\ref{t:ranked<=>EP} and Proposition~\ref{p:wr-ex<=>3ex}, the weakly regular $3$-ranked liner $X$ is ranked. By Corollary~\ref{c:Max=dim}, $F=\overline I$ for some independent set $I\subseteq F$ of cardinality $|I|=\|F\|=n+1$. Choose any distinct points $o,b\in I$ and consider the hyperplane $A\defeq \overline{I\setminus\{b\}}$ in the flat $F$. By the inductive assumption,
 $$|A|=1+(|X|_2-1)\sum_{r=0}^{n-2}(\kappa+|X|_2-1)^r.$$ Let $\mathcal L_o\defeq\{L\in\mathcal L:o\in L\subseteq A\}$ be the family of all lines in $A$ that contain the point $o$. 
Observe that
\begin{equation}\label{eq:|A|=}
|A|=|\{o\}\cup\bigcup_{L\in\mathcal L_o}(L\setminus\{o\})|=1+|\mathcal L_o|(\lambda-1).
\end{equation}

For any distinct lines $L,L'\in\mathcal L_o$, the flats $\overline{L\cup \{b\}}$ and $\overline{L'\cup\{b\}}$ are planes such that $\Aline ob\subseteq \overline{L\cup\{b\}}\cap\overline{L'\cup\{b\}}$. Assuming that $\Aline ob\ne \overline{L\cup\{b\}}\cap\overline{L'\cup\{b\}}$, we obtain $$2=\|\Aline ob\|<\|\overline{L\cup\{b\}}\cap\overline{L'\cup\{b\}}\|\le 3\mbox{ \ and hence \ }\overline{L\cup\{b\}}\cap\overline{L'\cup\{b\}}=\overline{L\cup\{b\}}=\overline{L'\cup\{b\}},$$ by the rankedness of $X$. Then $L\cup L'\subseteq A\cap(\overline{L\cup\{b\}}\cup\overline{L'\cup\{b\}})=A\cap\overline{L\cup\{b\}}$. Since $3=\|L\cup L'\|\le\|A\cap\overline{L\cup\{b\}}\|\le\|L\cup\{b\}\|=3$, the rankedness of the liner $X$ ensures that $b\in\overline{L\cup\{b\}}=\overline{L\cup L'}\subseteq A=\overline{I\setminus\{b\}}$, which contradicts the independence of the set $I$. This contradiction shows that $\overline{L\cup\{b\}}\cap\overline{L'\cup\{b\}}=\Aline ob$ for any distinct lines $L,L'\in\mathcal L_o$. Since $X$ is  weakly regular,
$$
F=\overline{A\cup\{b\}}=\bigcup_{a\in A}\overline{\{a,o,b\}}=\bigcup_{L\in\mathcal L_o}\overline{L\cup\{b\}}=\Aline ob\cup\bigcup_{L\in\mathcal L_o}(\overline{L\cup\{b\}}\setminus\Aline ob)
$$and hence
$$\begin{aligned}
|F|&=|\Aline ob|+\sum_{L\in\mathcal L_o}|\overline{L\cup\{b\}}\setminus\Aline ob|=|X|_2+|\mathcal L_o|(1+(\kappa+|X|_2)(|X|_2-1)-|X|_2)\\
&=|X|_2+|\mathcal L_o|(|X|_2-1)(\kappa+|X|_2-1)=|X|_2+(|A|-1)(\kappa+|X|_2-1)\\
&=|X|_2+\big((|X|_2-1)\sum_{r=0}^{n-2}(\kappa+|X|_2-1)^n\big)(\kappa+|X|_2-1)\\
&=|X|_2+(|X|_2-1)\sum_{r=1}^{n-1}(\kappa+|X|_2-1)^r=1+(|X|_2-1)\sum_{r=0}^{(n+1)-2}(\kappa+|X|_2-1)^r,
\end{aligned}
$$
by Theorem~\ref{t:2-balance+k-parallel=>3-balance}, the equality (\ref{eq:|A|=}) and the inductive hypothesis.
\end{proof}

\begin{question} For which cardinals $\kappa$ there exist weakly regular $3$-long $\kappa$-balanced ranked liners of rank $\|X\|>3$?
\end{question} 

\begin{proposition} Every $2$-balanced liner $X$ with $|X|_2<|X|<|X|_2^3-2{\cdot}|X|_2^2+2{\cdot}|X|_2$ is a ranked plane. 
\end{proposition}

\begin{proof} Assuming that $X$ is not a ranked plane, we can find a plane $P$ such that $P\ne X$. Choose any line $L\subset P$ and points $a\in P\setminus L$ and $b\in X\setminus P$. Observe that 
$$|P|\ge \Big|\{a\}\cup\bigcup_{x\in L}(\Aline xa\setminus\{a\})\Big|=1+|L|{\cdot}(|X|_2-1)=|X|_2^2-|X|_2+1$$and 
$$|X|\ge \Big|\{b\}\cup\bigcup_{y\in P}(\Aline yb\setminus\{b\})\Big|\ge 1+(|X|_2^2-|X|_2+1){\cdot}(|X|_2-1)=|X|_2^3-2{\cdot}|X|_2^2+2{\cdot}|X|_2,$$which contradicts our assumption.
\end{proof}

\begin{theorem}\label{t:affine=>Avogadro} Every affine liner $X$ is $2$-balanced.
\end{theorem}

\begin{proof} First we prove the equality $|L|=|\Lambda|$ for two intersecting lines $L,\Lambda$ in $X$. This equality trivially holds if $L=\Lambda$. So, we assume that $L\ne \Lambda$. Let $o$ be the unique point of the intersection $L\cap\Lambda$. Choose any points $a\in L\setminus \{o\}=L\setminus\Lambda$ and $b\in \Lambda\setminus \{o\}=\Lambda\setminus L$. Consider the sets $L^\circ\defeq L\setminus\{o,a\}$ and $\Lambda^\circ\defeq \Lambda\setminus\{o,b\}$ and the relation 
$$F\defeq\{(x,y)\in L^\circ\times \Lambda^\circ:\Aline xy\cap \Aline ab=\varnothing\}.$$
The affinity axiom implies that the relation $F$ is a bijective function between the sets $L^\circ$ and $\Lambda^\circ$. Then $|L|=|L^\circ|+2=|\Lambda^\circ|+2=|\Lambda|$.

Now take any disjoint lines $L,\Lambda$ in $X$. Choose any points $a\in L$ and $b\in \Lambda$. By the preceding case (of intersecting lines), we have  $|L|=|\Aline ab|=|\Lambda|$. Therefore, all lines in $X$ have the same cardinality, which means that $X$ is $2$-balanced.
\end{proof}

\begin{remark}\label{r:order-affine} By Theorem~\ref{t:affine=>Avogadro}, all lines in an affine liner have the same cardinality, which is called the \index{order}\index{affine liner!order of}\defterm{order} of the affine liner. By the Veblen--Young Theorem (1908), the order of any line-finite affine regular liner of dimension $\ge 3$ is a power of a prime number. There is a (still unproved) conjecture that the same result is true for line-finite affine planes: their orders are powers of primes. By a famous Bruck-Ryser Theorem~\ref{t:Bruck-Ryser}, if a number $n\in (4\IZ+1)\cup(4\IZ+2)$ is an order of an affine plane, then $n=x^2+y^2$ for some integer numbers $x,y$. This theorem implies that the order of an affine plane cannot be equal to $6$, $14$, or $22$. On the other hand, the number $10=4\cdot 2+2=1^2+9^2$ is the sum of two squares but no affine plane of order 10 exists, as was shown by Lam, Thiel and Swiercz~\cite{LTS1989} in 1989, using heavy computer calculations. The problem of the (non-)existence of an affine plane of order $12$ remains open, out of reach of modern (super)computers.
\end{remark}

\begin{exercise} Show that every affine regular liner of order $n\ge 3$ and dimension $d$ contains exactly $n^d$ points.
\end{exercise}

\begin{exercise} Find an example of an affine liner $X$ which is $2$-balanced with $|X|_2=3$ but not $3$-balanced.
\smallskip

{\em Hint:} Look at Example~\ref{ex:Tao}.
\end{exercise}

\begin{exercise} Find an example of an affine $3$-regular liner $X$ which is $3$-balanced but not $4$-balanced.
\smallskip

{\em Hint:} Look at Example~\ref{ex:HTS}.
\end{exercise}


Next, we evaluate  the cardinality of lines in proaffine and projective liners.

\begin{proposition}\label{p:Avogadro-proaffine} Let $L,\Lambda$ be two intersecting lines in a (projective) proaffine liner $X$. If $\overline{L\cup\Lambda}\ne L\cup \Lambda$, then $|L|\le|\Lambda|+1$ (and $|L|=|\Lambda|$).
\end{proposition}

\begin{proof} Assume that $\overline{L\cup\Lambda}\ne L\cup \Lambda$ and choose a point $c\in \Aline L\Lambda\setminus(L\cup\Lambda)$. The inequality $\overline{L\cup\Lambda}\ne L\cup\Lambda$ implies that $L\ne \Lambda$ and hence $L\cap\Lambda=\{o\}$ for a unique point $o$. Assuming that $\Aline L\Lambda=L\cup\Lambda$, we conclude that the set $L\cup\Lambda$ is flat and hence $\overline{L\cup\Lambda}=L\cup\Lambda$, which contradicts our assumption. This contradiction shows that there exists a point $c\in\Aline L\Lambda\setminus(L\cup\Lambda)$. Consider the relation 
$$F\defeq\{(x,y)\in L\times \Lambda:c\in \Aline xy\}.$$
For a set $S$, let $$F[S]=\{y\in \Lambda:\exists x\in L\cap S\;\;(x,y)\in F\}\;\;\mbox{and}\;\;F^{-1}[S]=\{x\in L:\exists y\in \Lambda\cap S\;\;(x,y)\in F\}$$ be the image and the preimage of the set $S$ under the relation $F$, respectively.  It is easy to see that $F$ is a bijective function between the sets $F^{-1}[\Lambda]$ and $F[L]$. It follows from $c\notin L\cup\Lambda$ that $o\notin F^{-1}[\Lambda]\cup F[L]$. 

If the space $X$ is proaffine, then there exist a point $u\in L\setminus\{o\}$ such that for every point  $x\in L\setminus\{u\}$, the intersection $\Aline xc\cap\Lambda$ contains some point $y$. If $x\ne o$, then $x\ne y$ and $y\in\Aline xc$ implies $c\in\Aline xy$ and hence $(x,y)\in F$ and $x\in F^{-1}[\Lambda]$.
Therefore, $L\setminus\{o,u\}\subseteq F^{-1}[\Lambda]$. By analogy we can find a point $v\in \Lambda\setminus\{o\}$ such that $\Lambda\setminus\{o,v\}\subseteq F[L]$.  
Then $$|L|=|L\setminus\{o,u\}|+2\le|F^{-1}[\Lambda]|+2=|F[L]|+2\le|\Lambda\setminus\{o\}|+2=|\Lambda|+1.$$ 

If the liner $X$ is projective, then $L\setminus\{o\}=F^{-1}[\Lambda]$ and $\Lambda\setminus\{o\}=F[L]$, which implies
$$|L|=|L\setminus\{o\}|+1=|F^{-1}[\Lambda]|+1=|F[L]|+1=|\Lambda\setminus\{o\}|+1=|\Lambda|.$$
\end{proof}


\begin{proposition}\label{p:3-long=>not=line+line} If a liner $X$ is $3$-long, then $\overline{L\Lambda}\ne L\cup\Lambda$ for any distinct lines $L,\Lambda$ in $X$.
\end{proposition}

\begin{proof} If two lines $L,\Lambda\subseteq X$ are distinct, then $|L\cap\Lambda|\le 1$ and we can choose points $x\in L\setminus\Lambda$ and $y\in\Lambda\setminus L$. Since $X$ is $3$-long, there exists a point $z\in\Aline xy\setminus\{x,y\}\subseteq\overline{L\cup\Lambda}$. Assuming that $z\in L$, we conclude that $y\in\Aline xy=\Aline xz\subseteq L$, which contradicts the choice of $y$. This contradiction shows that $z\notin L$. By analogy we can prove that $z\notin\Lambda$.
\end{proof}

Propositions~\ref{p:Avogadro-proaffine} and \ref{p:3-long=>not=line+line} imply the following corollaries.

\begin{corollary}\label{c:Avogadro-projective} Every $3$-long projective liner $X$ is $2$-balanced.
\end{corollary}

\begin{proof} Take any lines $L,\Lambda$ and choose any points $x\in L$ and $y\in\Lambda$. By Propositions~\ref{p:Avogadro-proaffine} and \ref{p:3-long=>not=line+line}, $|L|=|\Aline xy|=|\Lambda|$.
\end{proof}

If a liner $X$ of rank $\|X\|\ge 2$ is not $3$-long, then it contains a line of length $2$ and there is no further restrictions on the length of lines  in such liners.

\begin{example} For every set $K$ of non-zero cardinals with $1,2\in K$, there exists a projective liner $X$ of cardinality $|X|=\sum_{k\in K}k$ and rank $2\cdot|K|-1$ such that $K=\{|\Aline xy|:x,y\in X\}$.
\end{example}

\begin{proof} For every cardinal $k\in K$, consider the set $L_k\defeq \{k\}\times k$ and observe that $|L_k|=k$ and $L_\alpha\cap L_\beta=\varnothing$ for every distinct cardinals $\alpha,\beta\in K$. Let $X\defeq \bigcup_{k\in K}L_k$ and $$\Af\defeq\{xyz\in X:y\in \{x,z\}\vee (|\{x,y,z\}|=3\wedge \exists k\in K\setminus\{1\}\;\{x,y,z\}\subseteq L_k\}.$$
It is easy to see that $(X,\Af)$ is a projective liner of cardinality $|X|=\sum_{k\in K}k$ such that $\|X\|=2\cdot|K|-1$ and $\{|\Aline xy|:x,y\in X\}=\{|L_k|:k\in K\}=K$. 
\end{proof}

\begin{corollary}\label{c:L<L+2} Any lines $L,\Lambda$ in a $3$-long proaffine liner $X$ have $|L|\le|\Lambda|+2$.
\end{corollary}

\begin{proof} Choose any points $x\in L$ and $y\in\Lambda$. By Propositions~\ref{p:Avogadro-proaffine} and \ref{p:3-long=>not=line+line}, $$|L|\le|\Aline xy|+1=(|\Lambda|+1)+1=|\Lambda|+2.$$
\end{proof}

\begin{corollary} Every $\w$-long proaffine liner is $2$-balanced.
\end{corollary}

\begin{exercise} Prove a liner $X$ is proaffine if every line in $X$ contains at most $3$ points.
\end{exercise} 

\begin{question}\label{q:L+1} Let $X$ be a $3$-long proaffine liner. Is $|L|\le|\Lambda|+1$ for any lines $L,\Lambda$ in $X$?
\end{question}

\begin{remark} Question~\ref{q:L+1} has an affirmative answer for $3$-long Proclus liners, see Theorem~\ref{t:card-lines-in-Proclus}.
\end{remark}

\begin{remark} If all lines in a projective liner have the same finite cardinality $\lambda$, then the number $\lambda-1$ is called the \index{order}\index{projective liner!order of}\defterm{order} of the projective liner. By the Veblen--Young Theorem, the order of a finite projective liner of dimension $\ge 3$ is a power of a prime number. There is a (still unproved) conjecture that the same result holds for projective planes: their order is a power of a prime number. Using projective completions (discussed in Chapter~\ref{s:completions}), it can be shown that a  projective plane of order $n$ exists if and only if an affine plane of order $n$ exists. For information on possible orders of finite affine planes, see Remark~\ref{r:order-affine}.
\end{remark}

\section{Projective liners}

Let us recall that two lines $A,B$ in a liner $X$ are {\em skew} if $\|A\cup B\|=4$. It is clear that any skew lines are disjoint. The converse is true in projective liners.

\begin{theorem}\label{t:projective<=>} For every liner $X$, the following conditions are equivalent:
\begin{enumerate}
\item $X$ is $0$-hypoparallel;
\item $X$ is $0$-parallel;
\item any disjoint lines in $X$ are skew;
\item $X$ projective;
\item $X$ is strongly regular.
\end{enumerate}
\end{theorem}

\begin{proof} The implication $(1)\Ra(2)$ is trivial.
\smallskip

$(2)\Ra(3)$ Assume that a liner $X$ is $0$-parallel and take any disjoint lines $L,\Lambda$ in $X$. Assuming that $L,\Lambda$ are not skew, we conclude that $\|L\cup\Lambda\|\ne 4$ and hence $P=\overline{L\cup\Lambda}$ is a plane in $X$ containing two disjoint lines, which contradicts the $0$-parallelity of $X$. 
\smallskip 

$(3)\Ra(4)$ Assume that any disjoint lines in a liner $X$ are skew. To prove that $X$ is projective, take any points $o,x,y\in X$, $p\in\Aline xy$ and $v\in\Aline oy\setminus\{p\}$. We have to prove that $\Aline vp\cap\Aline ox\ne\varnothing$. To derive a contradiction, assume that $\Aline vp\cap\Aline ox=\varnothing$. Assuming that $o=x$, we conclude that $\Aline xy=\Aline oy=\Aline vp$ and hence $o\in\Aline vp\cap\Aline ox$, which contradicts our assumption. Therefore, $o\ne x$ and $\Aline ox$ is a line. By our assumption, the disjoint lines $\Aline vp$ and $\Aline ox$ are skew and hence $4=\|\Aline vp\cup\Aline ox\|\le\|\overline{\{o,y,x\}}\|\le 3$, which is a contradiction showing that $\Aline vp\cap\Aline ox\ne\varnothing$ and the liner $X$ is projective.
\smallskip

$(4)\Ra(5)$  Assume that a liner $X$ is projective. To show that $X$ is strongly regular, take any flat $A\subseteq X$ and point $b\in X\setminus A$. We need to show that $\overline{A\cup\{b\}}=\Aline Ab$. This equality will follow as soon as we check that the set $\Aline Ab$ is flat.
 Given two points $x,y\in \Aline Ab$ and $z\in\Aline xy$, we need to show that $z\in \Aline Ab$. If $\{x,y\}\subseteq A$, then $z\in\Aline xy\subseteq A\subseteq\Aline Ab$ and we are done. So, assume that $\{x,y\}\not\subseteq A$.
Since $x,y\in\Aline Ab$, there exist points $p,q\in A$ such that $x\in\Aline pb$ and $y\in\Aline qb$.

If $x=p$, then $\{x,y\}\not\subseteq A$ implies that $y\ne q$ and hence $b\in\Aline qb=\Aline yq$. Since $z\in \Aline xy$, by the projectivity axiom, there exists a point $a\in \Aline xq=\Aline pq\subseteq A$ such that 
$z\in\Aline ab\subseteq\Aline Ab$.

\begin{picture}(200,95)(-200,-15)
\put(0,0){\line(0,1){60}}
\put(0,60){\line(-1,-2){30}}
\put(0,60){\line(1,-2){30}}
\put(-40,0){\line(1,0){80}}
\put(-30,0){\line(3,2){45}}

\put(-30,0){\circle*{3}}
\put(-34,-8){$x$}
\put(0,0){\color{red}\circle*{3}}
\put(-2,-8){\color{red}$a$}
\put(30,0){\circle*{3}}
\put(28,-8){$q$}
\put(0,20){\circle*{3}}
\put(2,15){$z$}
\put(15,30){\circle*{3}}
\put(19,28){$y$}
\put(0,60){\circle*{3}}
\put(-2,63){$b$}
\put(45,-3){$A$}
\end{picture}

By analogy we can prove that $y=q$ implies $z\in\Aline Ab$.

So, we assume that $x\ne p$ and $y\ne q$. Since $b\in\Aline pb=\Aline px$ and $z\in\Aline xy$, by the projectivity axiom, there exists a point $w\in\Aline py$ such that $z\in\Aline wb$. Since  $b\in\Aline qb=\Aline qy$ and  $w\in \Aline qy$, by the projectivity axiom, there exists a point $a\in\Aline pq$ such that $w\in\Aline qb$. Then $z\in \Aline wb\subseteq \Aline ab\subseteq\Aline Ab$.
\smallskip

\begin{picture}(200,95)(-200,-15)
\put(0,0){\line(0,1){60}}
\put(0,60){\line(-1,-2){30}}
\put(0,60){\line(1,-2){30}}
\put(-40,0){\line(1,0){80}}
\put(-15,30){\line(1,0){30}}
\put(-30,0){\line(3,2){45}}

\put(-30,0){\circle*{3}}
\put(-34,-8){$p$}
\put(0,0){\color{red}\circle*{3}}
\put(-2,-8){\color{red}$a$}
\put(30,0){\circle*{3}}
\put(28,-8){$q$}
\put(-15,30){\circle*{3}}
\put(-23,27){$x$}
\put(0,30){\circle*{3}}
\put(1,32){$z$}
\put(15,30){\circle*{3}}
\put(19,28){$y$}
\put(0,60){\circle*{3}}
\put(-2,63){$b$}
\put(45,-3){$A$}
\put(0,20){\color{blue}\circle*{3}}
\put(2,15){\color{blue}$w$}
\end{picture}

$(5)\Ra(1)$ Assume that a liner $X$ is strongly regular.  Assuming that $X$ is not $0$-hypoparallel,  we can find a plane $P\subseteq X$, lines $L,\Lambda\subseteq P$  and a point $x\in X\setminus L$ such that $x\in\Lambda\subseteq P\setminus L$. Choose any point $y\in \Lambda\setminus\{x\}$. By Corollary~\ref{c:proregular=>ranked} and Theorem~\ref{t:ranked<=>EP}, the strongly regular liner $X$ has the Exchange Property and is ranked. Then $P=\overline{L\cup\{x\}}$.  The strong regularity of $X$ implies that $y\in \Lambda\subseteq P=\overline{L\cup\{x\}}=\bigcup_{a\in L}\Aline ax$ and hence $y\in\Aline ax$ for some $a\in L$. It follows from $y\in\Lambda\setminus L$ that $y\ne a$ and hence $\Aline ax=\Aline yx=\Lambda$ and $a\in\Lambda\cap L=\varnothing$, which is a contradiction showing that the strongly regular liner $X$ is $0$-proaffine.
\end{proof}

\begin{corollary}\label{c:proj=>ranked} Every projective liner is  ranked and has the Exchange Property.
\end{corollary}

Theorem~\ref{t:projective<=>} and Corollary~\ref{c:proregular=>ranked} imply the following corollary.

\begin{corollary}\label{c:Steiner-projective<=>} A $2$-balanced liner $X$ with $|X|_2<\w$ is projective if and only if $X$ is $3$-balanced and $|X|_3=|X|_2^2-|X|_2+1$.
\end{corollary}

\begin{proof} If the balanced liner $X$ is projective, then by Theorem~\ref{t:projective<=>}, $X$ is $0$-parallel and by Theorem~\ref{t:2-balance+k-parallel=>3-balance}, $X$ is $3$-balanced with $|X|_3=1+|X|_2(|X|_2-1)=|X|_2^2-|X|_2+1$.

Next, assume that $X$ is $3$-balanced and $|X|_3=|X|_2^2-|X|_2+1$. It follows from $|X|_2<\w$ that $|X|_3=1+|X|_2(|X|_2-1)>|X|_2(|X|_2-1)\ge |X|_2$. Since $|X|_3-1=(0+|X|_2)(|X|_2-1)$, by Proposition~\ref{p:23-balance=>k-parallel}, $X$ is $0$-parallel and by Theorem~\ref{t:projective<=>}, $X$ is projective.
\end{proof}

Theorems~\ref{t:projective<=>}, \ref{t:wr+k-parallel=>n-balanced}, and Corollary~\ref{c:proj=>ranked} imply the following corollary.

\begin{corollary}\label{c:projective-order-n} Every finite projective liner $X$ of order $\ell$ contains exactly $$\sum_{k=0}^{\|X\|-1}\ell^k=\frac{\ell^{\|X\|}-1}{\ell-1}$$ points.
\end{corollary}

By Corollary~\ref{c:Avogadro-projective}, every $3$-long projective liner is $2$-balanced.  Now we are going to describe the structure of projective liners which are not $3$-long. This descriprion involves the notion of a maximal $n$-long flat in a liner. 

\begin{definition}\label{d:max-long} Let $\kappa$ be a cardinal number. A flat $B$ in a liner $X$ is called \index{flat!$\kappa$-long}\index{$\kappa$-long flat}\defterm{$\kappa$-long} if $|\Aline xy|\ge \kappa$ for every distinct points $x,y\in B$. A $\kappa$-long flat $B$ in $X$ is called \index{flat!maximal $\kappa$-long}\index{maximal $\kappa$-long flat}\defterm{maximal $\kappa$-long} if every $\kappa$-long flat $A\subseteq X$ containing $B$ coincides with $B$.
\end{definition}

\begin{remark} By Definition~\ref{d:max-long}, every flat $F\subseteq X$ of cardinality $|F|\le 1$ in a liner $X$ is $\kappa$-long for every cardinal $\kappa$. The Kuratowski-Zorn Lemma implies that every $\kappa$-long flat in a liner $X$ is a subset of a maximal $\kappa$-long flat in $X$. In particular, every point $x$ of a liner $X$ belongs to some maximal $\kappa$-long flat (which can be equal to the singleton $\{x\}$, if the cardinal $\kappa$ is too large).
\end{remark}

\begin{lemma}\label{l:ox=2} If $M$ is a maximal $3$-long flat in a projective liner $X$, then $\Aline ox=\{o,x\}$ for every points $o\in M$ and $x\in X\setminus M$.
\end{lemma}

\begin{proof} To derive a contradiction, assume that $\Aline ox\ne\{o,x\}$ and hence $|\Aline ox|\ge 3$. Since $\Aline ox$ is a $3$-long flat containing the singleton $\{o\}$, the maximal $3$-long flat $M$ is not equal to $\{o\}$ and hence $|M|\ge 3$.  It is easy to see that the flat $M$ in the projective liner $X$ is a projective liner (as a subliner of $X$). By Corollary~\ref{c:Avogadro-projective}, the $3$-long projective liner $M$ is $2$-balanced and hence $|M|_2\ge 3$ is a well-defined cardinal.

\begin{claim}\label{cl:oxM2} $|\Aline ox|=|M|_2$.
\end{claim}

\begin{proof} Since $\Aline ox\ne \{o,x\}$, there exists a point $x'\in\Aline ox\setminus\{o,x\}$. Since the flat $M$ is $3$-long, there exist points $y\in M\setminus\{o\}$ and $y'\in\Aline oy\setminus\{o,y\}$. Since $X$ is projective, there exists a point $z\in \Aline xy\cap\Aline {x'}{y'}$. It is easy to see that $z\notin\Aline ox\cup\Aline oy$, witnessing that $\overline{\{o,x,y\}}\ne\Aline ox\cup\Aline oy$.
Applying Proposition~\ref{p:Avogadro-proaffine}, we conclude that $|\Aline ox|=|\Aline oy|=|M|_2$. 
\end{proof}

\begin{claim}\label{cl:axM2} For every $a\in M$, the line $\Aline ax$ has cardinality $|\Aline ax|=|M|_2$.
\end{claim}

\begin{proof} If $a=o$, then $|\Aline ax|=|\Aline ox|=|M|_2$, by Claim~\ref{cl:oxM2}. So, assume that $a\ne o$. It follows from $x\notin M$ that $\Aline ox\cap M=\{o\}$ and $\Aline ax\cap M=\{a\}$. Since the flat $M$ is $3$-long, $\Aline oa\not\subseteq \{o,a\}=M\cap(\Aline ox\cup\Aline oa)$ and hence $\overline{\Aline ax\cup\Aline ox}=\overline{\{a,o,x\}}\not\subseteq\Aline ox\cup\Aline ax$.   Applying Proposition~\ref{p:Avogadro-proaffine} and Claim~\ref{cl:oxM2}, we conclude that $|\Aline ax|=|\Aline ox|=|M|_2$.
\end{proof}

\begin{claim}\label{cl:Mx=M2} Every line $\Lambda$ in $\overline{M\cup\{x\}}$ has cardinality $|\Lambda|=|M|_2\ge 3$. 
\end{claim}

\begin{proof} Given any line $\Lambda$ in $\overline{M\cup\{x\}}$, choose any distinct points $u,v\in\Lambda$. By Theorem~\ref{t:projective<=>}, the projective liner $X$ is strongly regular and hence $\overline{M\cup\{x\}}=\bigcup_{a\in M}\Aline ax$. Then there exist points $a,b\in M$ such that $u\in\Aline ax$ and $v\in \Aline bx$.  

If $v\in\Aline ax$, then $\Aline uv=\Aline ax$ and $|\Aline uv|=|\Aline ax|=|M|_2$, by Claim~\ref{cl:axM2}. If $u\in\Aline bx$, then $\Aline uv=\Aline bx$ and $|\Aline uv|=|\Aline bx|=|M|_2$, by Claim~\ref{cl:axM2}.

So, we assume that $u\notin\Aline bx$ and $v\notin\Aline ax$. In this case $a\ne b$, $u\ne v$, and $u\ne x\ne v$. If $u=a$ and $v=b$, then $\Lambda=\Aline ab\subseteq M$ and $|\Lambda|=|M|_2$. 

If $u=a$ and $v\ne b$, then $x\in \Aline vb\setminus(\Aline ab\cup\Aline av)$ and hence $\overline{\Aline ab\cup\Aline av}\ne\Aline ab\cup\Aline av$. By Proposition~\ref{p:Avogadro-proaffine}, 
$|\Lambda|=|\Aline uv|=|\Aline av|=|\Aline ab|=|M|_2$. By analogy we can show that $|\Lambda|=|M|_2$ if $u\ne a$ and $v=b$. 

So, we assume that $u\ne a$ and $v\ne b$. Since the liner $X$ is projective, the lines $\Aline uv$ and $\Aline ab$ in the plane $\overline{\{a,x,b\}}$ have a common point. Observe that $\Aline xa\cap\Aline uv=\{u\}$ and $\Aline xa\cap\Aline ab=\{a\}$ and hence $x\notin \Aline uv\cup\Aline ab$. Since $x\in \Aline ua\subseteq\overline{\Aline uv\cup\Aline ab}$ and $\Aline uv\cap\Aline ab\ne\varnothing$, Proposition~\ref{p:Avogadro-proaffine} ensures that $|\Lambda|=|\Aline uv|=|\Aline ab|=|M|_2$.\end{proof}

Claim~\ref{cl:Mx=M2} implies that the flat $\overline{M\cup\{x\}}$ is $3$-long, which contradicts the maximality of the $3$-long flat $M$. This contradiction shows that $\Aline ox=\{o,x\}$, which completes the proof of the lemma.
\end{proof}

A family of sets $\F$ is called \index{disjoint family}\index{family!disjoint}\defterm{disjoint} if $A\cap B=\varnothing$ for every distinct sets $A,B\in\F$. The following theorem shows that the line structure of a projective liner is uniquely deteremined by the line structures of its maximal $3$-long flats.

\begin{theorem}\label{t:projective=+3long} The family $\mathcal M$ of maximal $3$-long flats in a projective liner $X$ is disjoint, and the line relation $\Af$ of $X$ is equal to
$$\big\{xyz\in X^3:y\in\{x,z\}\big\}\cup\bigcup_{M\in\mathcal M}(\Af\cap M^3)$$
\end{theorem}

\begin{proof} Choose any distinct maximal $3$-long flats $M,M'$ in $X$. The maximality of $M$ ensures that $M\not\subseteq M'$ and hence there exists a point $x\in M\setminus M'$. Assuming that the flats $M,M'$ have a common point $o\in M\cap M'$, we conclude that the line $\Aline ox$ in the $3$-long liner $M$ has cardinality $|\Aline ox|\ge 3$. On the other hand, Lemma~\ref{l:ox=2} ensures that $\Aline xo=\{x,o\}$. This is a contradiction showing that $M\cap M'=\varnothing$. Therefore, the family $\mathcal M$ of maximal $3$-long flats in $X$ is disjoint.

It remains to show that the line relation $\Af$ of $X$ is equal to the relation $$\Af'\defeq\{xyz\in X^3:y\in\{x,z\}\}\cup\bigcup_{M\in\mathcal M}(\Af\cap M^3).$$ It is clear that $\Af'\subseteq\Af$. To show that $\Af\subseteq\Af'$, take any triple $xyz\in\Af$. If $|\{x,y,z\}|\le 2$, then $y\in\{x,z\}$ and hence $xyz\in\Af'$. If $|\{x,y,z\}|=3$, then $\Aline xz$ is a $3$-long flat in $X$. By the Kuratowski--Zorn Lemma, there exists a maximal $3$-long flat $M\in\mathcal M$ in $X$ that contains the $3$-long flat $\Aline xz$. Then $xyz\in\Af\cap M^3$. 
\end{proof}

Next, we describe the structure of projective planes which are not $3$-long.

\begin{definition} A liner $X$ is called a \index{near-pencil}\index{liner!near-pencil}\defterm{near-pencil} if $X$ contains a line $L\subseteq X$ such that $X\setminus L$ is a singleton.
\end{definition}

In the following theorem by a projective plane we understand a projective liner of rank $3$. 

\begin{theorem} A liner $X$ is a near-pencil if and only if $X$ is projective plane which is not $3$-long.
\end{theorem}

\begin{proof} To prove the ``if'' part, assume that $X$ is a near-pencil. Then $X$ contains a line $L$ such that $X\setminus L=\{x\}$ for some point $x\in X$. It follows that the liner $X$ has rank 3 and hence is a plane. Every line in $X$ coincides with $L$ or with the line $\Aline xy$ for some $y\in L$. Consequently, any two lines in $X$ have a common point and hence the liner $X$ is projective, by Theorem~\ref{t:projective<=>}. To see that the near-pencil $X$ is not $3$-long, take any point $y\in L$ and observe that $$\Aline xy\setminus\{y\}=\Aline xy\setminus (\Aline xy\cap L) =(\Aline xy\cap X)\setminus L=\Aline xy\cap(X\setminus L)=\Aline xy\cap\{x\}=\{x\}$$ and hence  $\Aline xy=\{x,y\}$. 
\smallskip

To prove the ``only if'' part, assume that $X$ is a projective plane and $X$ is not $3$-long. Let $\mathcal M$ be the family of maximal $3$-long flats in $X$. If every flat $M\in\mathcal M$ is a singleton, then every line in $X$ consists of two points and every set in $X$ is flat. Then $|X|=\|X\|=3$ and $X=\{x,y,z\}$ for some distinct points $x,y,z$. Since $L\defeq \Aline xy=\{x,y\}$ and $X\setminus L=\{z\}$, the liner $X$ is a near-pencil. If some flat $M\in\mathcal M$ is not a singleton, then $M$ contains a $3$-long line $L$. Since the flat $M$ is $3$-long and the liner $X$ is not $3$-long, $M\ne X$ and hence $2\le\|L\|\le \|M\|<\|X\|=3$ and $L=M$, by the rankedness of the projective liner $X$ (which follows from  Theorem~\ref{t:projective<=>} and Corollary~\ref{c:proregular=>ranked}). Now take any point $x\in X\setminus M$. Lemma~\ref{l:ox=2} implies that $M\cup\{x\}=L\cup\{x\}$ is a flat in $X$. Since $\|L\cup\{x\}\|=3=\|X\|$, the rankedness of the projective liner $X$ ensures that $X=L\cup\{x\}$ and hence $X$ is a near-pencil.    
\end{proof}



%
%

\begin{proposition}\label{p:projective<} Let $A,B$ be two flats in a projective liner $X$. If $\|A\|_B<\|A\|$, then $A\cap B\ne\emptyset$.
\end{proposition}

\begin{proof} By Theorem~\ref{t:projective<=>}, the projective liner $X$ is strongly regular and by  Corollary~\ref{c:proregular=>ranked}, the space $X$ is ranked and has the Exchange Property. By Proposition~\ref{p:rank-EP3}, there exists a set $C\subseteq A\setminus\overline B$ of cardinality $|C|=\|A\|_B$ such that $A\subseteq \overline{B\cup C}$. Assuming that $\|A\|_B<\|A\|$, we conclude that $A\not\subseteq \overline{C}$. Since $A\cap\overline{B\cup C}=A\ne\overline{C}$, by Proposition~\ref{p:aff-finitary}, there exists a finite set $F\subseteq B\cup C$ such that $A\cap\overline F\not\subseteq \overline C$. We can assume that $F$ has the smallest possible cardinality. 

We claim that $F\cap C=\emptyset$. To derive a contradiction, assume that $F\cap C\ne\emptyset$ and fix a point $c\in F\cap C$. By the choice of $F$, there exists a point $a\in (A\cap \overline F)\setminus\overline C$. The minimality of the set $F$ ensures that $A\cap\overline{F\setminus\{c\}}\subseteq \overline C$ and hence $a\in\overline{(F\setminus\{c\})\cup\{c\}}\setminus\overline{F\setminus\{c\}}$. By the strong regularity of $X$, there exists a point $x\in \overline{F\setminus\{c\}}$ such that $a\in\Aline xc$. The choice of $a\in A\setminus\overline C\subseteq A\setminus\{c\}$ ensures that $a\ne c$ and hence $x\in\Aline xc=\Aline ac\subseteq A$. Then $x\in A\cap\overline{F\setminus\{c\}}\subseteq \overline C$ and $a\in\Aline xc\subseteq\overline C$, which contradicts the choice of the point $a\notin \overline C$. This contradiction show that $F\cap C=\emptyset$. Then $F\subseteq (B\cup C)\setminus C\subseteq B$ and hence $a\in A\cap\overline F\setminus\overline C\subseteq A\cap B\ne\varnothing$.
\end{proof}

\begin{corollary}\label{c:projective<} If $A,B$ are two disjoint flats in a projective liner $X$, then $\|A\|_B=\|A\|$.
\end{corollary}

\begin{corollary}\label{c:line-meets-hyperplane} Every line $L$ in a projective liner $X$ has nonempty intersection with every hyperplane $H$ in $X$.
\end{corollary}

\begin{proof} Since $\|X\|_H=1<2=\|L\|$, the intersection $L\cap H$ is not empty, by Proposition~\ref{p:projective<}.
\end{proof}


Theorem~\ref{t:projective<=>} implies that projective liners can be equivalently defined by the following four first-order axioms.

\begin{theorem}\label{t:proj-1-axioms} A set $X$ endowed with a ternary relation $\Af\subseteq X^3$ is a projective liner if and only if the following four axioms are satisfied:
\begin{itemize}
\item[{\sf (I)}] {\sf Identity:} $\forall x,y\in X\;\;\big(\Af xyx\;\Ra\;x=y)$;
\item[{\sf (R)}] {\sf Reflexivity:} $\forall x,y\in X\;\;(\Af xxy\;\wedge\;\Af xyy)$;
\item[{\sf (E)}] {\sf Exchange:} $\forall a,b,x,y\in X\;\big(( \Af axb\wedge \Af ayb\wedge x\ne y)\Rightarrow (\Af xay\wedge\Af xby)\big)$;
\item[{\sf (V)}] {\sf Veblen's Axiom:} $\forall o,x,y,p,u\in X\;\big(\Af ypx\;\wedge\;p\ne y\;\wedge\; \Af ouy)\;\Ra\;\exists v\;(\Af ovx\;\wedge \Af upv)$.
\end{itemize}
\end{theorem}

\section{Projective subplanes in projective liners}

In this section we prove a nice result of Bruck\footnote{{\bf Richard Hubert Bruck} (1914 -- 1991) was an American mathematician best known for his work in the field of algebra, especially in its relation to projective geometry and combinatorics. Bruck studied at the University of Toronto, where he received his doctorate in 1940 under the supervision of Richard Brauer. He spent most his career as a professor at University of Wisconsin--Madison, advising at least 31 doctoral students. He is best known for his 1949 paper coauthored with H. J. Ryser, the results of which became known as the Bruck--Ryser theorem concerning the possible orders of finite projective planes. In 1946, he was awarded a Guggenheim Fellowship. In 1956, he was awarded the Chauvenet Prize for his article ``Recent Advances in the Foundations of Euclidean Plane Geometry''. In 1962, he was an invited speaker at the International Congress of Mathematicians in Stockholm. In 1963, he was a Fulbright Lecturer at the University of Canberra. In 1965 a Groups and Geometry conference was held at the University of Wisconsin in honor of Bruck's retirement. Dick Bruck and his wife Helen were supporters of the fine arts. They were patrons of the regional American Players Theatre in Wisconsin.} on the interplay between the orders of projective planes and their projective subplanes.

\begin{theorem}[Bruck, 1955]\label{t:Bruck55} Let $X$ be a projective liner of finite order $n$ and $P$ be a projective subliner of order $p<n$ in $X$. If $\|P\|\ge 3$, then either $n=p^2$ or else $n\ge p^2+p$.
\end{theorem}  

\begin{proof} Let $\mathcal L$ be the family of lines in $X$. Since $X$ is a projective liner of order $n$, $|L|=n+1$ for every $L\in\mathcal L$. Fix any line $\Lambda\in\mathcal L$ such that $|\Lambda\cap P|\ge 2$, and observe that $\Lambda\cap P$ is a line in the liner $P$. Since the projective liner $P$ has order $p<n$, $|\Lambda\cap P|=p+1<n+1=|\Lambda|$ and hence there exists a point $\lambda\in\Lambda\setminus P$. Since $\|P\|\ge 3$, there exists a point $a\in P\setminus\Lambda$. Consider the plane $P'=\overline{(\Lambda\cap P)\cup\{a\}}$ in $P$.
By Corollary~\ref{c:projective-order-n}, $|P'|=1+p+p^2$ and hence $|P'\setminus\Lambda|=p^2$. We claim that $\Aline\lambda x\cap P'=\{x\}$ for every $x\in P'\setminus \Lambda$. Assuming that $|\Aline \lambda x\cap P'|\ge 2$, we conclude that $\Aline \lambda x\cap P'$ is a line in the projective plane $P'$. By the projectivity of $P'$, the lines $\Lambda\cap P'=\Lambda\cap P$ and $\Aline \lambda x\cap P'$ have a common point $y\in P'\cap\Lambda\cap\Aline \lambda x\subseteq \Lambda\cap\Aline \lambda x=\{\lambda\}$, which contradicts the choice of $\lambda\in\Lambda\setminus P$. This contradiction shows that $|\Aline\lambda x\cap P'|=1$ and hence the family $\mathcal L'_\lambda\defeq\{L\in\mathcal L:(\lambda\in L)\;\wedge\;(L\cap P'\ne\varnothing)\}$ has cardinality $|\mathcal L'_\lambda|=1+|P'\setminus\Lambda|=1+p^2$. Now consider the plane $\Pi\defeq\overline{\Lambda\cup\{a\}}$ in $X$ and observe that $P'\subseteq \Pi$. Since the projective liner $X$ has order $n$, $1+n\ge|\mathcal L'_\lambda|=1+p^2$ and hence $n\ge p^2$. 

It remains to prove that $n\ne p^2$ implies $n\ge p^2+p$. Assuming that $n>p^2$, we conclude that the family $\mathcal L_\lambda\defeq\{L\in\mathcal L:\lambda\in L\subseteq \Pi\}$ has cardinality $|\mathcal L_\lambda|=1+n>1+p^2=|\mathcal L'_\lambda|$ and hence there exists a line $\Lambda'\in\mathcal L_\lambda\setminus\mathcal L'_\lambda$. Then $\Lambda'\cap P'=\varnothing$. 
Since $P'$ is a projective plane of order $p$, it contains $p^2+p+1$ points and $p^2+p+1$ lines. Consequently, the family $\mathcal L'\defeq\{L\in\mathcal L:|L\cap P'|\ge 2\}$ has cardinality $p^2+p+1$. By the projectivity of the plane $\Pi$, for every line $L\in\mathcal L'$, the intersection $\Lambda'\cap L$ is a singleton. Moreover, for any distinct lines $L,L'\in\mathcal L'$, the intersections $L\cap P'$ and $L'\cap P'$ are distinct lines in the projective plane $P'$, which implies that $L\cap L'=(L\cap P')\cap (L'\cap P')$ is a singleton in $P'$ and hence the singletons $L\cap\Lambda'$ and $L'\cap\Lambda'$ are distinct. Then $1+n=|\Lambda'|\ge|\mathcal L'|=1+p+p^2$ and hence $n\ge p^2+p$. Let $k\defeq n-p^2-p\ge 0$.

\end{proof}  

Theorem~\ref{t:Bruck55} motivates the following long-standing open problem.

\begin{problem} Is there a projective plane $X$ of composite order $p^2+p$ for some $p$? In particular, is there a projective plane of order $12=3^2+3$?
\end{problem}

\section{Proaffine and Proclus liners} 

In this section we study the interplay between proaffine and Proclus liners, and will prove that Proclus liners are exactly proaffine $3$-proreglar liners. 

\begin{theorem}\label{t:Proclus<=>}For a liner $X$, the following conditions are equivalent:
\begin{enumerate}
\item $X$ is Proclus;
\item $X$ is proaffine and $3$-proregular;
\item  for every line $L\subseteq X$ and points $o\in L$, $x\in X\setminus L$,  $y\in\overline{L\cup\{x\}}$,\newline the set $\{u\in\Aline ox:\Aline uy\cap L=\varnothing\}$ contains at most one point.
\end{enumerate}
\end{theorem}

\begin{proof} We shall prove the implications $(1)\Ra(2)\Ra(3)\Ra(1)$. 
\smallskip

$(1)\Ra(2)$ Assume that a liner $X$ is Proclus. We divide the proof of the proaffinity and $3$-proregularity into two claims.

\begin{claim} The Proclus liner $X$ is proaffine.
\end{claim}

\begin{proof} Assuming that $X$ is not proaffine, we can find points $o,x,y\in X$, $p\in\Aline yx\setminus\Aline ox$ such that the set $I\defeq\{u\in\Aline oy:\Aline up\cap\Aline ox=\emptyset\}$ contains two distinct points $v,w$. In this case $p\notin\Aline oy$ (otherwise $I=\{p\}$ cannot contain two distinct points). 

 It follows from $p\notin\Aline yx\setminus\Aline oy$ that $x\notin \Aline oy$ and hence $L\defeq\Aline ox$ is a line. It follows from $p\in\Aline yx\setminus\Aline ox$ that  $y\in \Aline yp=\Aline px$. 
Since $p\notin\Aline oy$ and $v,w\in\Aline oy$, the flats $\Lambda\defeq\Aline vp$ and  $\Lambda'\defeq\Aline wp$ are lines  such that 
$$\Lambda\cup\Lambda'\subseteq\overline{\{v,p,w\}}\subseteq\overline{\{o,y,p\}}\subseteq\overline{\{o,p,x\}}=\overline{\Aline ox\cup\{p\}}=\overline{L\cup\{p\}}.$$Since $X$ is Proclus, $\Lambda=\Lambda'$ and hence $\{v,w\}\subseteq\Aline vp\cap\Aline wp\cap\Aline oy$ and $o\in \Aline op=\Aline wv=\Aline vp=\Aline wp$, which contradicts $v,w\in I$. This contradiction shows that the Proclus liner $X$ is proaffine.
\end{proof}

\begin{claim}\label{cl:Proclus=>3-regular} The Proclus liner $X$ is $3$-proregular.
\end{claim}

\begin{proof} Take any set $A\subseteq X$ of cardinality $|A|<3$ and  points $o\in\overline A$ and $p\in X\setminus\overline A$ with $\Aline op\ne\{o,p\}$. We have to show that the flat hull $\overline{A\cup\{p\}}$ coincides with the set $\Lambda\defeq\bigcup_{x\in \overline A}\bigcup_{y\in \overline{o\,p}}\Aline xy$. If $|A|\le1$, then $o\in \overline{A}=A=\{o\}$ and $\overline{A\cup\{p\}}=\Aline op=\Lambda$. So, we assume that $|A|=2$. In this case $\overline A$ is a line in $X$. 
The equality $\overline{A\cup\{p\}}=\Lambda$ will follow as soon as we check that the set $\Lambda$ is flat. Given any points $a,b\in\Lambda$ and $c\in\Aline ab$, we have to show that $c\in\Lambda$. To derive a contradiction, assume that $c\notin \Lambda$.  In this case $\{a,b\}\not\subseteq \overline A\cup\Aline op$. We lose no generality assuming that $b\notin\overline A\cup\Aline op$. 
Two cases are possible.
\smallskip

First assume that the line $\Aline ab$ has a common point $\alpha$ with the flat $\overline A$. Since $\alpha,c\in\Aline ab$ and $\alpha\ne c$ (because $\alpha\in\overline A$ and $c\notin\overline A$), we have $b\in\Aline ab=\Aline \alpha c$.

Since $b\in\Lambda$, there exist points $x\in\overline A$ and $y\in\Aline op$ such that $b\in\Aline  xy$. It follows from $x\in\overline A$ and $b\in \Aline xy\setminus\overline A$ that $y\notin \overline A$, $b\ne x$ and hence $y\in\Aline xy=\Aline xb\subseteq \overline{\{x,\alpha,c\}}$. If $\alpha=x$, then $x\ne y\in\overline{\{x,\alpha,c\}}=\Aline xc$ implies $c\in\Aline xy\subseteq\Lambda$ are we are done. So, we assume that $\alpha\ne x$. In this case $L\defeq\Aline \alpha x$ is a line in $\overline A$ and  $o\in \overline A=L$. Then  $y\in\overline{\{\alpha,x,c\}}=\overline{L\cup\{c\}}$ and hence $\Aline op=\Aline oy\subseteq \overline{L\cup\{c\}}$. Since $c\notin\Aline op$, for every $v\in\Aline op$, the flat $\Aline vc$ is a line in the plane $\overline{L\cup\{c\}}$. Since $X$ is Proclus, the set $I\defeq\{v\in\Aline op:\Aline vc\cap L=\varnothing\}$ contains at most one point. By our assumption,   $\Aline op\ne\{o,p\}$, so there exists a point $v\in\Aline op\setminus(\{o\}\cup I)$. Since $v\notin I$, there exists a point $\gamma\in\Aline vc\cap L$. It follows from $\gamma\in L=\overline A$ and $c\notin\overline A$ that $c\in\Aline vc=\Aline v\gamma\subseteq \Lambda$, which contradicts our assumption.
\smallskip

This contradiction shows that $\Aline ab\cap \overline A=\varnothing$. Assuming that $\Aline op\cap\Aline ab=\varnothing$, we conclude that $\overline A$ and $\Aline op$ are two distinct lines in the plane $\overline{A\cup \{p\}}$ that contain the point $o$ and are disjoint with the line $\Aline ab$, which contradicts the Proclus Parallel Postulate for $X$. This contradiction shows that $\Aline op\cap\Aline ab$ is not empty and hence contains some point $v$. Since the line $\Aline op$ contains more than two points,  there exists a point $y\in\Aline op\setminus\{o,v\}$. Assuming that $y\in\overline A$, we conclude that $p\in\Aline oy\subseteq\overline A$ and hence $b\in\Lambda\subseteq\overline A$, which contradicts our assumption. It follows from $c\notin\Aline op$ that $\Aline vc$ is a line such that $\Aline vc\subseteq \Aline ab\subseteq X\setminus\overline A$. Since $X$ is Proclus,  $\Aline vc$ is a unique line in the plane $\overline{A\cup \{p\}}$ that is disjoint with the line $\overline A$. Then the line $\Aline yc$ has a common point $x$ with the line $\overline A$. Since $y\notin \overline A$, the points $x,y$ are distinct and hence $c\in\Aline yc=\Aline xy\subseteq \Lambda$, witnessing that the set $\Lambda$ is flat. Then $\overline{A\cup\{u\}}=\Lambda$ and the liner $X$ is $3$-proregular.
\end{proof}



$(2)\Ra(3)$ Assume that the liner $X$ is proaffine and $3$-proregular. To check the condition (3), fix any  line $L\subseteq X$ and points $o\in L$, $y\in X\setminus L$, and  $z\in\overline{L\cup\{y\}}$. We should prove that the set $I\defeq\{u\in\Aline oy:\Aline uz\cap L=\emptyset\}$ contains at most one point.
If $z\in L$, then $I=\varnothing$ and we are done. So, assume that $z\notin L$. If $\Aline oy=\{o,y\}$, then $I\subseteq\{y\}$ and we are done. So, assume that $\Aline oy\ne\{o,y\}$. 

Since the liner $X$ is $3$-proregular, $\overline{L\cup\{y\}}=\bigcup_{s\in\overline{o\,y}}\Aline sL$ and hence $z\in \Aline vx$ for some $v\in\Aline oy$ and $x\in L$. Assuming that $v=o$, we conclude that $z\in\Aline vx=\Aline ox\subseteq L$, which is a contradiction showing that $v\ne o$ and hence $\Aline ov=\Aline oy$.

Since the liner $X$ is proaffine, the set $J\defeq\{u\in\Aline ov:\Aline uz\cap\Aline ox=\emptyset\}$ contains at most one point. It follows from $\Aline ov=\Aline oy$ and $\Aline ox\subseteq L$ that $I\subseteq J$. So, the set $I$ contains at most one point.
\smallskip

$(3)\Ra(1)$ Assume that the liner $X$ satisfies the condition $(3)$.

\begin{claim}\label{cl:3=>3-ranked} The liner $X$ has the $3$-Exchange Property and hence is $3$-ranked.
\end{claim}

\begin{proof} Take any set $A\subseteq X$ of cardinality $|A|<3$ and any points $x\in X$ and $y\in\overline{A\cup\{x\}}\setminus\overline A$.

We have to prove that $x\in\overline{A\cup\{y\}}$. If for every $a\in \overline A$ the line $\Aline ax$ coincides with $\{a,x\}$, then the set $\overline A\cup\{x\}$ is flat and $\overline{A\cup\{x\}}=\overline A\cup\{x\}$. Then $y\in\overline{A\cup\{x\}}\setminus\overline A=(\overline A\cup\{x\})\setminus\overline A\subseteq\{x\}$ and hence $x=y\in\overline{A\cup\{y\}}$.

So, we assume that $\Aline ax\ne\{a,x\}$ for some $a\in \overline A$. If $A=\{a\}$, then $y\in\overline{A\cup\{x\}}\setminus\overline A=\Aline ax\setminus\{a\}$ implies $x\in\Aline ax=\Aline ay=\overline{A\cup\{y\}}$, by Theorem\ref{t:Alines}. It remains to consider the case of $|A|=2$. In this case, $\overline A$ is a line in $X$. If the line $\Aline xy$ has a common point $z$ with the line $\overline A$, then $z,y\in\Aline xy$ implies $x\in \Aline xy=\Aline zy\subseteq\overline{A\cup\{y\}}$ and we are done. So, assume that $\Aline xy\cap\overline A=\varnothing$. The condition (3) ensures that the set $I=\{u\in\Aline ax:\Aline uy\cap \overline A=\varnothing\}$ coincides with the singleton $\{x\}$. By the assumption, the line $\Aline ax$ contains a point $u\in\Aline ax\setminus\{a,x\}\subseteq \Aline ax\setminus I$. Assuming that $u\in\overline A$, we conclude that $\Aline ax=\Aline au\subseteq\overline A$ and $y\in\overline {A\cup\{x\}}=\overline A$, which contradicts the choice of $y$. This contradiction shows that $u\notin\overline A$. Since $u\notin I$,  the line $\Aline uy$ has a common point $\alpha$ with the line $\overline A$. Taking into account that $y\notin\overline A$, we conclude that $u\in\Aline \alpha y$ and $x\in\Aline ax=\Aline au\subseteq\overline{\{a,\alpha,y\}}\subseteq\overline{A\cup\{y\}}$, witnessing that $X$ has the $3$-Exchange Property. By Theorem~\ref{t:ranked<=>EP}, the liner $X$ is $3$-ranked. 
\end{proof}
 
To derive a contradiction, assume that the liner $X$ is not Proclus. Then there exists a plane $P\subseteq X$, a line $L\subseteq P$, a point $p\notin X\setminus L$ and two distinct lines $\Lambda,\Lambda'$ such that $p\in\Lambda\cap\Lambda'$, $\Lambda\cup\Lambda'\subseteq P\setminus L$.

Assuming that $\Aline op=\{o,p\}$ for every point $o\in L$, we conclude that the set $L\cup\{p\}$ is a plane that coincides with the plane $P$ because $X$ is $3$-ranked. Then $\Lambda\subseteq P\setminus L=\{p\}$ is not a line. This contradiction shows that for some point $o\in L$, the line $\Aline op$ contains more than two points.

\begin{claim}\label{cl:proaffine-Lambda} For every $\lambda\in(\Lambda\cup\Lambda')\setminus\{p\}$ and $y\in\Aline op\setminus\{o,p\}$, we have $\Aline \lambda y\cap L\not\subseteq\{o\}$.
\end{claim}

\begin{proof} We lose no generality assuming that $\lambda\in\Lambda$. By Claim~\ref{cl:3=>3-ranked}, the liner $L$ is $3$-ranked and hence $P=\overline{L\cup\{p\}}$. Since $o\in L$, $\lambda\in\Lambda\subseteq P=\overline{L\cup \{p\}}$ and $\Aline p\lambda\cap L\subseteq\Lambda\cap L=\varnothing$, the condition (3) ensures that the set $I\defeq \{v\in\Aline op:\Aline v\lambda\cap L=\varnothing\}$ coincides with the singleton $\{p\}$. Then for every $y\in\Aline op\setminus\{o,p\}\subseteq \Aline op\setminus I$, there exists a point $x\in\Aline y\lambda\cap L$. Assuming that $x=o$, we obtain $\Aline y\lambda=\Aline yo=\Aline \lambda o$. Then $\lambda\in \Aline\lambda y=\Aline yo=\Aline yp=\Aline op$ and $o\in\Aline op\cap\Aline ox\subseteq \Aline p\lambda\cap L=\Lambda\cap L=\varnothing$, which is a contradiction showing that $x\ne o$ and hence $x\in\Aline \lambda y\cap L\not\subseteq\{o\}$.
\end{proof}

\begin{claim}\label{cl:proaffine-Lambda'} For every $x\in L\setminus\{o\}$ and $y\in\Aline op\setminus\{o,p\}$, we have $\Aline xy\cap\Lambda\not\subseteq\{p\}$.
\end{claim}

\begin{proof} Choose any point $\lambda\in\Lambda\setminus\{p\}$. By Claim~\ref{cl:proaffine-Lambda}, there exists a point $z\in \Aline \lambda y\cap L\setminus\{o\}$. Then $L=\Aline oz$ and $\Aline op=\Aline yp$. Since $x\in L=\Aline zo\subseteq\Aline{\Aline \lambda y}{\Aline py}=\overline{\{\lambda,p,y\}}\subseteq\overline{\{\lambda,p,o\}}=\overline{\Lambda\cup\{o\}}$ and $\Aline ox\cap\Lambda\subseteq L\cap\Lambda=\varnothing$, the condition (3) ensures that
the set $I\defeq\{v\in \Aline po:\Aline vx\cap\Lambda=\varnothing\}$ is equal to the singleton $\{o\}$. Since $y\ne o$, we obtain that that $\Aline yx\cap\Lambda\ne\varnothing$. Fix any point $q\in\Lambda\cap\Aline yx$.  Assuming that $q=p$ and taking into account that $y\ne p= q$, we conclude that
$\Aline yx=\Aline yq=\Aline yp=\Aline ox$ and hence $q\in \Lambda\cap\Aline ox\subseteq \Lambda\cap L=\varnothing$, which is a contradiction showing that $q\ne p$ and hence $\Aline xy\cap\Lambda\not\subseteq\{p\}$.
\end{proof}

\begin{picture}(100,100)(-100,0)

\put(10,20){\line(1,0){100}}
\put(115,16){$L$}
\put(60,60){\line(2,-1){50}}
\put(115,30){$\Lambda'$}
\put(60,60){\line(-2,1){50}}
\put(60,60){\line(2,1){50}}
\put(115,84){$\Lambda$}
\put(60,60){\line(-2,-1){50}}
\put(60,20){\line(0,1){40}}
\put(40,20){\line(1,1){60}}

\put(60,20){\circle*{3}}
\put(58,12){$o$}
\put(40,20){\circle*{3}}
\put(35,23){$x$}
\put(60,40){\circle*{3}}
\put(51,39){$y$}
\put(60,60){\circle*{3}}
\put(58,63){$p$}
\put(73.5,53.5){\circle*{3}}
\put(69,42){$\lambda'$}
\put(100,80){\circle*{3}}
\put(97,70){$\lambda$}

\end{picture}

Now we can finish the proof of the implication $(3)\Ra(1)$. Fix any point $\lambda'\in\Lambda'\setminus\{p\}$. The choice of the point $o$ ensures that there exists a point $y\in\Aline op\setminus\{o,p\}$. Then $\Aline op=\Aline yp=\Aline oy$. By Claim~\ref{cl:proaffine-Lambda}, there exists a point $x\in\Aline {\lambda'}y\cap L\setminus\{o\}$ and by Claim~\ref{cl:proaffine-Lambda'}, there exists a point $\lambda\in \Aline xy\cap\Lambda\setminus\{p\}$. It follows from $\lambda\in \Lambda\setminus L\subseteq \Lambda\setminus\{x\}$ that $y\in \Aline xy=\Aline x\lambda$. Since the space $X$ is proaffine, the set $I\defeq\{v\in\Aline x\lambda:\Aline vp\cap \Aline xo=\varnothing\}$ contains at most one point. Taking into account that  $\lambda,\lambda'\in I$, we conclude that $\lambda=\lambda'$ and hence $\Lambda=\Aline p\lambda=\Aline p{\lambda'}=\Lambda'$, which contradicts the choice of the lines $\Lambda$ and $\Lambda'$.
\end{proof}

Theorem~\ref{t:Proclus<=>} and Proposition~\ref{p:3-proregular<=>} imply the following first-order characterization of Proclus liners.

\begin{proposition} A liner $(X,\Af)$ is Proclus if and only if it satisfies two first-order axioms:
\begin{enumerate}
\item {\sf proaffine:} $\forall o,x,y,p\;\big((\Af xpy\wedge\neg\Af opx)\Ra (\exists u \;\forall v\;\; (\Af ovy\wedge v\ne u)\Ra \exists t \;(\Af vtp\wedge \Af otx))\big)$;
\item {\sf $3$-proregular:} $\forall o,a,b,u,v,x,y,z,p\in X$\newline
$( \Af ovu\wedge\Af axu\wedge\Af byv\wedge \Af oba\wedge\Af xzy\wedge \Af opu\wedge o\ne p\ne u)\Rightarrow\exists s,w\;(\Af osa\wedge\Af owu\wedge \Af szw)$.
\end{enumerate}
\end{proposition}

\begin{exercise} Find an example of an affine liner, which is not $3$-regular.
\smallskip

{\em Hint:} Look at the liner in Example~\ref{ex:Tao}.
\end{exercise}

\begin{exercise} Find an example of a proaffine ranked liner, which is not $3$-regular.
\smallskip

{\em Hint:} Look at the line in Example~\ref{ex:Steiner13}.
\end{exercise}

The following proposition is a version of the classical Proclus Parallelity Postulate.

\begin{proposition}[Proclus Postulate]\label{p:Proclus-Postulate} Let $P$ be a plane in a Proclus liner $X$ and $L,L',\Lambda$ be three lines in the plane $P$. If $L\cap L'=\varnothing$ and $|\Lambda\cap L|=1$, then $|\Lambda\cap L'|=1$.
\end{proposition}

\begin{proof}  Let $o$ be the unique point of the intersection $\Lambda \cap L$. Assuming that $\Lambda\cap L'=\varnothing$, we conclude that $L,\Lambda$ are two distinct lines in $P\setminus L'$ such that $o\in L\cap \Lambda$, which contradicts the Proclus Parallelity Postulate. This contradiction shows that $|\Lambda\cap L'|>0$. Assuming that $|\Lambda\cap L'|>1$ and applying Theorem~\ref{t:Alines}, we conclude that $\Lambda=L'$ and hence $1=|\Lambda\cap L|=|L'\cap L|=0$, which is a contradiction showing that $|\Lambda\cap L'|=1$.
\end{proof}

\begin{corollary}\label{c:Proclus-par=} Let $L,\Lambda$ be two disjoint lines in a  Proclus liner $X$. If $\|L\cup\Lambda\|\le 3$, then $|L|=|\Lambda|$.
\end{corollary}

\begin{proof} By Theorem~\ref{t:Proclus<=>}, the Proclus liner $X$ is $3$-proregular and by Proposition~\ref{p:k-regular<=>2ex}, the $3$-proregular liner $X$ is $3$-ranked. 

Since $\|L\cup\Lambda\|\le 3$, the flat hull $P\defeq\overline{L\cup\Lambda}$ is a plane. 
If $P=L\cup\Lambda$, then for every $x\in L$ and $y\in\Lambda$, the line $\Aline xy$ coincides with the doubleton $\{x,y\}$, which implies that the set $L\cup\{y\}$ is a plane. Since $X$ is $3$-ranked, the planes $L\cup\{y\}$ and $P=L\cup\Lambda$ coincide and hence $\Lambda=\{y\}$, which is a contradiction showing that $P\ne L\cup\Lambda$.

Then there exists a point $o\in P\setminus(L\cup\Lambda)$.  Proposition~\ref{p:Proclus-Postulate} implies that the relation $$F\defeq\{(x,y)\in L\times\Lambda:\Aline ox=\Aline oy\}$$is a bijective function between the lines $L,\Lambda$, witnessing that $|L|=|\Lambda|$.
\end{proof}  



%

\section{Affine and Playfair liners}

In this section we study the interplay between affine and Playfair liners, and show that Playfair liners are exactly affine $3$-regular $3$-long liners.


\begin{theorem}\label{t:affine-char1} For a liner $X$, the following conditions are equivalent:
\begin{enumerate}
\item $X$ is affine and $3$-proregular;
\item $X$ is affine and $3$-regular;
\item for every line $L\subseteq X$ and points $o\in L$, $y\in X\setminus L$,  $z\in\overline{L\cup\{y\}}\setminus L$, the set $\{u\in\Aline oy:\Aline uz\cap L=\emptyset\}$ is a singleton;
\item for every line $L\subseteq X$ and point $p\in X\setminus L$ with $\overline{L\cup\{p\}}\ne L\cup\{p\}$, there exists a unique line $\Lambda$ such that $p\in \Lambda\subseteq\overline{\{p\}\cup L}\setminus L$.
\end{enumerate}
\end{theorem}

\begin{proof} We shall prove the implications $(1)\Ra(2)\Ra(3)\Ra(4)\Ra(1)$. 
\smallskip

$(1)\Ra(2)$ Assume that $X$ is affine and $3$-proregular. By Theorem~\ref{t:affine=>Avogadro}, the affine liner $X$ is $2$-balanced. If $X$ is not $3$-long, then all lines in $X$ have length $2$, which implies that all sets in $X$ are flat. Then for every flat $A\subseteq X$ and point $x\in X\setminus A$ we have $\overline{A\cup\{x\}}=A\cup\{x\}=\bigcup_{a\in A}\Aline ax$, witnessing that the liner $X$ is strongly regular and hence $3$-proregular. If $X$ is $3$-long, then the $3$-proregularity of $X$ implies the $3$-regularity of $X$.
\smallskip 

$(2)\Ra(3)$ Assume that $X$ is affine and $3$-regular. To prove the condition (3), take any  line $L\subseteq X$ and points $o\in L$, $y\in X\setminus L$,  $z\in\overline{L\cup\{y\}}\setminus L$. Being affine, the liner $X$ is proaffine. So, we can apply Theorem~\ref{t:Proclus<=>}and conclude that the set $I\defeq\{u\in \Aline oy:\Aline uz\cap L=\varnothing\}$ contains at most one point. So, it remains to prove that the set $I$ is not empty. If $z\in\Aline oy$, then $z\in I$. So, assume that $z\notin\Aline oy$. 
 By the $3$-regularity of $X$, there exist points $v\in\Aline oy$ and $x\in L$ such that $z\in\Aline vx$. It follows from $z\in \Aline vx\setminus\Aline oy$ that $x\ne o$ and hence $\Aline ox=L$. Since $X$ is affine, there exists a point $u\in\Aline oy$ such that $\Aline uz\cap \Aline ox=\varnothing$. Since $\Aline uz\cap L= \Aline uz\cap \Aline ox=\varnothing$, the point $u$ belongs to the set $I$. 
\smallskip

$(3)\Ra(4)$ Assume that $X$ satisfies the condition (3).  To check the condition (4), take any line $L$ and point $p\in X\setminus L$ such that $\overline{L\cup \{p\}}\ne L\cup\{p\}$. At first we prove that there exists a line $\Lambda$ such that $p\in\Lambda\subseteq\overline{L\cup\{p\}}\setminus L$. By our assumption, there exists a point $z\in \overline{L\cup\{p\}}\setminus L\cup\{p\}$. If $\Aline zp\cap L=\varnothing$, then $\Lambda\defeq\Aline zp$ is a required line with $p\in\Lambda\subseteq\overline{L\cup\{p\}}\setminus L$.
So, we assume that $\Aline zp\cap L$ contains some point $x$. Choose any point $o\in L\setminus\{x\}$.  It follows from $z\in X\setminus L\subseteq X\setminus\{x\}$ that $p\in \Aline pz=\Aline zx\subseteq\overline{L\cup\{z\}}$. By the condition (3), there exists a point $u\in \Aline oz$ such that $\Aline up\cap L=\varnothing$.  Assuming that $u=p$, we conclude that $p\in\Aline oz\cap\Aline xz=z$, which contradicts the choice of $z$. This contradiction shows that $u\ne p$ and hence $\Lambda\defeq\Aline up$ is a line such that $p\in\Lambda\subseteq\overline{L\cup\{p\}}\setminus L$. 
Applying Theorem~\ref{t:Proclus<=>}, we conclude that $\Lambda$ is a unique line with $p\in\Lambda\subseteq\overline{L\cup\{p\}}\setminus L$.
\smallskip

$(4)\Ra(1)$ Assume that the condition (4) is satisfied. This condition implies the condition (3) in Theorem~\ref{t:Proclus<=>}, which implies that the liner $X$ is proaffine and $3$-proregular.  To prove that the liner $X$ is affine, take any points $o,x,y\in X$ and $p\in\Aline yx\setminus\Aline ox$. Since $X$ is proaffine, it suffices to find a point $u\in\Aline oy$ such that $\Aline up\cap\Aline ox= \varnothing$. If $p\in\Aline oy$, then the point $u\defeq p$ has the required property. So, assume that $p\notin\Aline oy$. In this case $p\in\Aline yx\setminus\Aline ox$ implies $x\ne o$ and hence  the flat $L\defeq \Aline ox$ is a line. Since $p\in\Aline yx\setminus (L\cup\{y\})\subseteq \overline{L\cup\{y\}}\setminus(L\cup\{y\})$, we can apply the condition (4) and find a line $\Lambda$ such that $p\in\Lambda\subseteq\overline{L\cup\{p\}}\setminus L$. By Proposition~\ref{p:k-regular<=>2ex}, the $3$-proregular liner $X$ is $3$-ranked, which implies that $\overline{\Lambda\cup\{o\}}=\overline{L\cup\{p\}}=\overline{L\cup\{y\}}$. We claim that $\Lambda\cap\Aline oy\ne\varnothing$. In the opposite case, $\Aline ox$ and $\Aline oy$ are two lines in the plane $\overline{\Lambda\cup\{o\}}$ such that $o\in\Aline ox\cap\Aline oy$ and $\Aline ox\cap\Lambda=\varnothing=\Aline oy\cap\Lambda$, which contradicts the condition (4) applied to the line $\Lambda$ and point $o$.
\end{proof}

\begin{theorem}\label{t:Playfair<=>} For a liner $X$ of cardinality $|X|>2$, the following conditions are equivalent:
\begin{enumerate}
\item $X$ is Playfair;
\item $X$ is affine, $3$-regular, and $3$-long;
\item for every line $L\subseteq X$ and point $x\in X\setminus L$ there exists a unique line $\Lambda$ such that $x\in \Lambda\subseteq\overline{L\cup\{x\}}\setminus L$.
\end{enumerate}
\end{theorem}

\begin{proof} We shall prove the implications $(1)\Ra(2)\Ra(3)\Ra(1)$.
\smallskip

$(1)\Ra(2)$ Assume that $X$ is Playfair. Then $X$ is Proclus. By Theorem~\ref{t:affine-char1}, $X$ is affine and $3$-regular. By Theorem~\ref{t:affine=>Avogadro}, all lines in $X$ have the same cardinality. Assuming that $X$ is not $3$-long, we conclude that all lines in $X$ have cardinality $2$ and hence all subsets of $X$ are flat. Since $|X|>2$, there exists a set $P\subseteq X$ of cardinality $|X|=3$. The set $X$ is a plane containing a $2$-element line $L\subseteq P$ and a unique point $x\in P\setminus L$. For this point $x$ there exists no lines $\Lambda$ such that $x\in\Lambda\subseteq P\setminus L$, which contradicts the Parallelity Postulate of Playfair. This contradiction shows that the liner $X$ is $3$-long.
\smallskip

$(2)\Ra(3)$ Assume that $X$ is affine, $3$-regular and $3$-long. To prove that the condition (3) holds, take any line $L\subseteq X$ and point $x\in X\setminus L$. Choose any point $o\in L$. Since $X$ is $3$-long, there exists a point $y\in\Aline ox\setminus\{o,x\}$. Then $y\in \overline{L\cup\{x\}}\setminus(L\cup\{x\})$. By Theorem~\ref{t:affine-char1}, there exists a unique line $\Lambda$ such that $x\in\Lambda\subseteq\overline{L\cup\{x\}}\setminus L$.
\smallskip

$(3)\Ra(1)$ Assume that the liner $X$ satisfies the condition (3). By Theorem~\ref{t:affine-char1}, the liner $X$ is Proclus. To prove that $X$ is Playfair, choose any plane $P\subseteq X$, line $L\subseteq P$ and point $x\in P\setminus L$. By the condition (3), there exists a line $\Lambda$ such that $x\in \Lambda\subseteq\overline{L\cup\{x\}}\setminus L\subseteq P\setminus L$. Since $X$ is Proclus, the line $\Lambda$ is unique, witnessing that $X$ is Playfair.
\end{proof}

\begin{corollary}\label{c:affine-cardinality} Every $3$-long affine regular liner $X$ of finite dimension $\dim(X)\ge 1$ is $2$-balanced and has cardinality $|X|=\big(|X|_2\big)^{\dim(X)}$.
\end{corollary}  

\begin{proof} By Theorem~\ref{t:Playfair<=>}, the liner $X$ is Playfair and hence $1$-parallel. By Theorem~\ref{t:wr+k-parallel=>n-balanced}, 
$$|X|=1+(|X|_2-1)\sum_{r=0}^{\|X\|-2}(1+|X|_2-1)^r=1+(|X|_2-1)\sum_{r=0}^{\|X\|-2}|X|_2^{r}=|X|_2^{\|X\|-1}=|X|_2^{\dim(X)}.$$
\end{proof} 

Let us recall that a liner $X$ is called \defterm{line-finite} if every line in $X$ contains finitely many points.

\begin{corollary}\label{c:Playfair<=>23balanced} A (line-finite) liner $X$ is Playfair if (and only if) $X$ is $2$-balanced, $3$-balanced and $|X|_3=(|X|_2)^2$.
\end{corollary}

\begin{proof} If $X$ is Playfair, then it is $3$-long and $1$-parallel. If $\|X\|\le1$, then $X$ contains no lines and no planes and hence $X$ is (vacuously) $2$-balanced and $3$-balanced and satisfies the equality $|X|_3=|X|_2^2$ (by the same vacuous reason).  If $\|X\|=2$, then $X$ is a line and hence $X$ is $2$-balanced. The linear $X$ contains no planes, so is vacuously $3$-balanced and $|X|_3=|X|_2^2$ (by the same vacuous reason). So, assume that $\|X\|\ge 3$. By Theorem~\ref{t:Playfair<=>}, the Playfair liner $X$ is affine and $3$-long. By Theorem~\ref{t:affine=>Avogadro}, the affine linear $X$ is $2$-balanced. Being Playfair, the liner $X$ is $1$-parallel. By Theorem~\ref{t:2-balance+k-parallel=>3-balance}, the $2$-balanced $1$-parallel liner is $3$-balanced with $|X|_3= 1+(1+|X|_2)(|X|_2-1)=1+(|X|_2^2-1)=|X|_2^2$.
\smallskip

Now assume that $X$ is line-finite, $2$-balanced, $3$-balanced and $|X|_3=|X|^2_2$. If $\|X\|<3$, then $X$ contains no planes and hence is Playfair vacuously. So, assume that $\|X\|\ge 3$. Since $X$ is line-finite, the cardinal $|X|_2$ is finite. Then every plane $P$ in $X$ has finite cardinality $|P|=|X|_3=|X|_2^2$. Assuming that $|X|_2=2$, we conclude that all lines contain exactly two points, which implies that every subset of $X$ is flat. Consequently, every plane in $X$ contains exactly three points, which contradicts the equality $|X|_3=|X|_2^2=4$. This contradiction shows that $|X|_2\ge 3$, which means that $X$ is $3$-long. By Proposition~\ref{p:23-balance=>k-parallel}, the liner $X$ is $\kappa$-parallel for a unique cardinal $\kappa$ that satisfies the equation $|X|_2^2-1=|X|_3=1+(\kappa+|X|_2)(|X|_2-1)$. The unique solution of this equation is $\kappa=1$, which means that $X$ is $1$-parallel and hence Playfair.
\end{proof}

\begin{proposition}\label{p:Playfair-plane<=>} A liner $X$ is a Playfair plane if and only if $X$ is $3$-long, $\|X\|>2$, and for every line $L\subseteq X$ and point $x\in X\setminus L$ there exists a unique line $L_x$ in $X$ such that $x\in L_x\subseteq X\setminus L$.
\end{proposition}

\begin{proof}  The ``only if'' part follows from Theorem~\ref{t:Playfair<=>} and the definition of a Playfair plane.

To prove the ``if'' part, assume that $X$ is $3$-long, $\|X\|>2$, and for every line $L\subseteq X$ and point $x\in X\setminus L$ there exists a unique line $L_x$ in $X$ such that $x\in L_x\setminus X\setminus L$.  To see that $X$ is Playfair, take any plane $P\subseteq X$, line $L\subseteq P$ and point $x\in P\setminus X$. By our assumption, there exists a unique line $L_x$ such that $x\in L_x\subseteq X\setminus L$. Assuming that $X\ne P$, we can apply Proposition~\ref{p:cov-aff} and find a point $y\in X\setminus(P\cup L_x)$. Then $\Aline xy\ne L_x$ is a line in $X$ such that $\Aline xy\cap P=\{x\}$ and hence $\Aline xy\cap L=\varnothing$, which contradicts the uniqueness of the line $L_x$. This contradiction shows that $X=P$ and hence $L_x$ is a unique line in $P$ such that $x\in L_x\subseteq P\setminus L$, witnessing that the liner $X$ is a Playfair plane.
\end{proof}

Theorem~\ref{t:Playfair<=>} implies the following first-order characterization of Playfair liners. 

\begin{proposition}\label{p:Playfair<=>3+a+3} A liner $(X,\mathsf L)$ is Playfair if and only if it satisfies the following three first-order axioms:
\begin{enumerate}
\item {\sf $3$-Long:} $\forall x,y\;\big((x\ne y\;\wedge\; \exists z\;(x\ne z\ne y))\Ra \exists z\;(\Af xzy\wedge x\ne z\ne y)\big)$;
\item {\sf Affine:} $\forall o,x,y,p\,\big((\Af xpy\,{\wedge}\, \neg\Af opx)\,{\Ra}\,\exists u\,(\Af ouy\,{\wedge}\, \forall v \,(\Af ovy\,{\Ra}\,(u\ne v\,{\Leftrightarrow}\,\exists t\,(\Af vtp\,{\wedge}\,\Af otx))))\big)$;
\item {\sf $3$-Regular:}  $\forall o,\!a,\!b,\!u,\!v,\!x,\!y,\!z\,
\big((\Af ovu{\wedge}\Af axu{\wedge}\Af oba{\wedge} \Af byv{\wedge}\Af xzy)\,{\Rightarrow}\,\exists s,\!w\,(\Af osa{\wedge}\Af owu{\wedge} \Af szw)\big)\!\,.$
\end{enumerate}
\end{proposition}


\begin{exercise} Find an example of an affine liner $X$ such that for every line $L\subseteq X$ and point $p\in X\setminus L$ there exist five distinct lines in the plane $\overline{L\cup\{p\}}$ that contain the point $p$ and are disjoint with the line $L$.
\smallskip

\noindent{\em Hint:} Look at the liner from Example~\ref{ex:Z15}.
\end{exercise}

\section{Hyperaffine, hyperbolic, and injective liners}\label{s:hyperbolic}

In this section we present some (non-trivial) examples of finite hyperaffine, hyperbolic, and injective liners, found by Ivan Hetman among known balanced incomplete block designs. 


  
\begin{example}\label{ex:hyperaffine-nonPlayfair} There exists a hyperaffine non-affine $2$-balanced $4$-parallel plane with $|X|_2=4$ and $|X|=|X|_3=25$. All its 50 lines are given in the table taken from \cite[1.34.1]{HCD}.
\begin{center}
{\tt 00000000111111122222223333344445555666778899aabbil\\
134567ce34578cd34568de468bh679f78ag79b9aabcddecejm\\
298dfbhkea6g9kf7c9afkg5cgfihdgifchi8ejjcjdfhgfghkn\\
iaolgmjnmbohnljonblhmjjdlknmeklnekmkinlimimonooloo
}
\end{center}
\end{example}

\begin{Exercise} Show that every hyperaffine liner $X$ of cardinality $|X|<25$ is affine.
\end{Exercise}

\begin{example}[Hetman \cite{HetmanS4}]\label{ex:non-hyperaffine} There exists a $4$-parallel balanced ranked liner $X$ with $|X|_2=4$ and $|X|=|X|_3=25$, which is not hyperaffine. It is the group $\IZ_5\times\IZ_5$ endowed with the family of lines $\mathcal L\defeq\{B+z:B\in\big\{\{00,10,01,22\},\{13,31,33,44\},\;z\in \IZ_5\times \IZ_5\}$. The non-hyperaffinity of this liner is witnessed by the points $o\defeq 00$, $x\defeq 41$, $y\defeq 01$, and $p\defeq 13$.
\end{example}

\begin{example}[Hetman\footnote{See {\tt mathoverflow.net/a/459391/61536}}]\label{ex:Hetman-hyperbolic91} \index[person]{Hetman}There exists a hyperbolic $8$-parallel plane $X$ with $|X|_2=7$ and $|X|=|X|_3=91$. It is the ring $\IZ_{91}$ endowed with the family of lines $\mathcal L\defeq\{x+L:x\in \IZ_{91},\;L\in\mathcal B\}$, where
$\mathcal B\defeq\big\{\{0, 8, 29, 51, 54, 61, 63\},\{0, 11, 16, 17, 31, 35, 58\},\{0, 13, 26, 39, 52, 65, 78\}\big\}.$
\end{example}

\begin{example}[Hetman\footnote{See, {\tt mathoverflow.net/a/458712/61536}}]\index[person]{Hetman}   There exists a hyperbolic $22$-parallel plane $X$ with $|X|_2=7$ and $|X|=|X|_3=175$. It is the ring $X\defeq \IZ_{7}\times\IZ_5\times\IZ_5$ endowed with the family of lines $\mathcal L\defeq\{x+L:x\in X,\;L\in\mathcal B\}$, where
$\mathcal B\defeq\big\{\{000,100,200,300,400,500,600\},\{000,113,142,222,233,420,430\}$, $\{000,\hskip-1pt 134,\hskip-1pt 121,\hskip-1pt 223,\hskip-1pt 232,\hskip-1pt 402,\hskip-1pt 403\},\{000,\hskip-1pt 112,\hskip-1pt 143,\hskip-1pt 211,\hskip-1pt 244,\hskip-1pt 401,\hskip-1pt 404\},\{000,\hskip-1pt 231,\hskip-1pt 124,\hskip-1pt 241,\hskip-1pt 214,\hskip-1pt 410,\hskip-1pt 440\}\big\}.
$
\end{example}

\begin{example}[Hetman\footnote{See {\tt mathoverflow.net/a/457561/61536}}]\label{ex:Hetman-hyperbolic217}\index[person]{Hetman} There exists a hyperbolic $29$-parallel plane $X$ with $|X|_2=7$ and $|X|=|X|_3=217$. It is the ring $\IZ_{217}$ endowed with the family of lines $\mathcal L\defeq\{x+L:x\in \IZ_{217},\;L\in\mathcal B\}$, where
$\mathcal B\defeq\big\{\{0,1,37,67,88,92,149\},\{0,15,18,65,78,121,137\},\{0,8,53,79,85,102,107\}$,\break $\{0,11,86,100,120,144,190\},\{0,29,64,165,198,205,207\},\{0,31,62,93,124,155,186\}\big\}.$
\end{example}

\begin{example}[Hetman\footnote{See \tt https://mathoverflow.net/a/457894/61536}\label{ex:Hetman-hyperbolic}]\index[person]{Hetman} For every $n\in\{3,4,5,7,8,9\}$, there exists an injective $2$-balanced $(n^2-n-1)$-parallel plane $X$ with $|X|_2=n+1$ and $|X|=|X|_3=n^3+1$. The liner $X$ contains $n^2(n^2-n+1)$ lines. For $n=3$, the liner $X$ has $3^2(3^2-3+1)=63$ lines, which are encoded in the columns of following table:
\begin{center}\tt
0000000001111111122222222333333334444455556666777788899aabbcgko\\
14567ghij4567cdef456789ab456789ab59adf8bce9bcf8ade9decfdfcedhlp\\
289abklmnba89lknmefdchgjijighfecd6klhilkgjnmhjmngiajgihigjheimq\\
3cdefopqrghijrqopqrponmklporqklmn7romnqpnmqoklrplkbopporqqrfjnr
\end{center}
\end{example}

The injective (and hence hyperbolic) liners found by Ivan Hetman in Example~\ref{ex:Hetman-hyperbolic} are classical unitals.

\begin{definition} A \index{unital}\defterm{unital} is a $2$-balanced liner $X$ with $|X|=n^3+1$ and $|X|_2=n+1$ for some finite cardinal $n\ge 2$. 
\end{definition}


\begin{definition} A \index{classical unital}\index{unital!classical}\defterm{classical unital} is the subliner $$U_q\defeq \{\IF_{q^2}(x,y,z):(x,y,z)\in\IF_{q^2}^3\setminus\{0\}^3\;\;(x^{q+1}+y^{q+1}+z^{q+1}=0)\}$$ 
of the projective plane $\mathbb P\IF_{q^2}^3$ over a field $\IF_{q^2}$ whose order $|\IF_{q^2}|=q^2$ is a square of some number $q$ (which is necessarily a prime power).
\end{definition}

\begin{Exercise}
Show that every classical unital $U_q$ is an injective $2$-balanced $(q^2-q-1)$-parallel liner with $|U_q|_2=q+1$ and $|U_q|=|U_q|_3=q^3+1$. If $q=2$, then the liner $U_q$ is affine. If $q\ge 3$, then the liner $U_q$ is hyperbolic.
\smallskip

\noindent{\em Hint:} For the proof, see Lemma 7.42 in \cite{UPP}.
\end{Exercise} 

Let us illustrate the geometry of the classical unital $U_3$, which is a subliner of the projective plane $\mathbb P\IF_9^3$ over the $9$-element field 
$$\IF_9=\{-1-i,-1,-1+i,-i,0,i,1-i,1,1+i\}$$where $i$ is an element of the field with $i^2=-1$.  The addition in the filed $\IF_9$ is by modulo $3$. In the field $\IF_9$ the elements $1$ and $-1$ have four roots of the $4$-th order:
$$
\begin{aligned}
\sqrt[4]{1}&\defeq\{z\in\IF_9:z^4=1\}=\{1,-1,i,-i\}\quad\mbox{and}\\
\sqrt[4]{-1}&\defeq\{z\in\IF:z^4=-1\}=\{1+i,1-i,-1+i,-1-i\}.
\end{aligned}
$$

The projective plane $\mathbb P\IF^3_9$ can be identified with the set $$(\{1\}\times\IF_9^2)\cup(\{(0,1)\}\times\IF_9\})\cup\{(0,0,1)\}$$ in $\IF_9^3$. The set $\{1\}\times\IF_9^2$ can be identified with the affine plane $\IF_9^2$, and $(\{(0,1)\}\times\IF_9)\cup\{(0,0,1)\}$ with the projective line $\IF_9\cup\{\infty\}$. Then the classical unital $U_3$ can be identified with the $28$-element set
$$
\begin{aligned}
&\{(x,y)\in \IF_9^2:x^4+y^4=-1\}\cup\{x\in \IF_9:x^4=-1\}\\
&=(\{0\}\times\sqrt[4]{-1})\cup(\sqrt[4]{-1}\times\{0\})\cup(\sqrt[4]{1}\times\sqrt[4]{1})\cup\sqrt[4]{-1}
\end{aligned}
$$in $\IF_9^2\cup\IF_9\cup\{\infty\}$.
\smallskip

In the following picture we draw these 28 points and also check the injectivity axiom on an example of the points $\color{blue}o=(-1,-1)$, $\color{blue}x=(1,-1)$, $y=(-1,1)$ and $p=(i,-i)\in\Aline xy$. For the points $u=(-1,-i)$ and $v=(-1,i)$ on the line $\Aline oy$ we can see that the lines $\Aline up$ and $\Aline vp$ indeed are disjoint with the line $\Aline ox$.



In fact, injective liners are not exotic at all, and as subliners are present in every $\w$-long plane. We recall that a {\em plane} is a liner of rank $3$.

\begin{theorem}\label{t:injective-subliner} Let $X$ be an $\w$-long plane. For every cardinal $\lambda\in [3,\w]$, there exists an injective subliner $Y\subseteq X$ such that 
\begin{enumerate}
\item $|Y|=\|Y\|=\w$;
\item $Y$ is $n$-balanced for every positive cardinal $n$;
\item $|Y|_2=\lambda$ and $|Y|_n=\w$ for every cardinal $n\ge 3$.
\end{enumerate}
\end{theorem}

\begin{proof} Given any cardinal $\lambda\in[3,\w]$, consider the sequence of finite cardinals $(\lambda_n)_{n\in\w}$   defined by the formula
$$
\lambda_n=\begin{cases}
\lambda&\mbox{if $\lambda<\w$};\\
\lambda+2&\mbox{if $\lambda=\w$}.
\end{cases}
$$


We shall inductively construct an increasing sequence of finite sets $(R_n)_{n\in\IN}$ and $(Y_n)_{n\in\IN}$ such that $Y_0=R_0$ and for every $n\in\IN$, the following conditions are satisfied:
\begin{enumerate}
\item[$(1_n)$] $|R_n|=2+n$, $R_{n-1}\subseteq R_n$ and $R_n\setminus R_{n-1}\subseteq X\setminus\bigcup_{x,y\in Y_{n-1}}\Aline xy$;
\item[$(2_n)$] $Y_{n-1}\subseteq Y_n\subseteq \bigcup_{x,y\in Y_{n-1}\cup R_{n}}\Aline xy$;
\item[$(3_n)$] $|Y_{n}\cap\Aline xy|=\lambda_n$ for every distinct points $x,y\in Y_{n-1}\cup R_{n}$;
\item[$(4_n)$] $|Y_n\cap\Aline xy|\le\lambda_n$ for every $x,y\in Y_n$.
\item[$(5_n)$] for all points $x,y\in Y_{n-1}\cup R_n$, $p\in \Aline xy\setminus (Y_{n-1}\cup R_n)$, $ q\in Y_n\setminus \Aline xy$, we have $\Aline pq\cap Y_n=\{p,q\}$.
\end{enumerate}

Assume that for some $n\in\IN$, the finite sets $Y_{n-1}$, $R_{n-1}$ such that satisfying the inductive conditions $(1_{n-1})$--$(5_{n-1})$ have been constructed (if $n=1$, then we require only that $Y_0=R_0$ and $\|\{x,y,z\}\|=3$ for any distinct points $x,y,z\in Y_0=R_0$). By Proposition~\ref{p:cov-aff}, the $\w$-long liner $X$ is not the union of finitely many lines. Consequently, there exists a set $R_n\subseteq X$ of cardinality $|R_n|=n+2$ such that $R_{n-1}\subseteq R_n$, $|R_n\setminus R_{n-1}|=1$ and $(R_n\setminus R_{n-1})\cap \bigcup_{x,y\in Y_{n-1}}\Aline xy=\varnothing$. 

The inductive condition $(4_{n-1})$ and the choice of the set $R_n\setminus R_{n-1}\subseteq X\setminus\bigcup_{x,y\in Y_{n-1}}\Aline xy$ ensure that for every $x,y\in Y_{n-1}$ and $u,v\in R_n\setminus R_{n-1}$  we have the following upper bounds:
$$
\begin{aligned}
&|(Y_{n-1}\cup R_n)\cap\Aline xy|=|Y_{n-1}\cap\Aline xy|\le\lambda_{n-1};\\
&|(Y_{n-1}\cup R_n)\cap\Aline xu|=|\{x,u\}|=2\le\lambda_{n-1};\\
&|(Y_{n-1}\cup R_n)\cap\Aline uv|=|\{u,v\}|=1\le\lambda_{n-1}.
\end{aligned}
$$Then the set $Z_0\defeq Y_{n-1}\cup R_n$ has the following property:
$\forall x,y\in Z_0\;\;|Z_0\cap\Aline xy|\le\lambda_{n-1}\le\lambda_n.$

Write the set $[Z_0]^2\defeq\{D\subseteq Z_0:|D|=2\}$ of $2$-elements subsets of $Z_0$ as $[Z_0]^2=\{D_k\}_{k=1}^m$ where $m=|[Z_0]^2|$. 

By finite induction we shall construct an increasing sequence of finite sets $(Z_k)_{k=1}^m$ such that  for every $k\in \{1,\dots,m\}$, the following conditions are satisfied:
\begin{enumerate}
\item[(i$_k$)] $Z_{k-1}\subseteq Z_k$;
\item[(ii$_k$)] $|Z_k\cap\overline {D_k}|=\lambda_n$;
\item[(iii$_k$)] $Z_k\setminus Z_{k-1}\subseteq \overline{D_k}\setminus \bigcup\{\Aline xy:x,y\in Z_{k-1},\;\Aline xy\not\subseteq \overline {D_k}\}$;
\item[(iv$_k$)] $|Z_k\cap\Aline xy|\le\lambda_n$ for every $x,y\in Z_{k}$.
\end{enumerate}

Assume that for some $k\in\{1,\dots,m\}$, the set $Z_{k-1}$ has been constructed.  The inductive condition (iv$_{k-1})$ ensures that $|Z_{k-1}\cap\overline {D_k}|\le\lambda_{n}$. Since the set $Z_{k-1}$ is finite and the line $\overline{D_n}$ is infinite, there exists a finite set $Z_k$ satisfying the inductive conditions (i$_k$)--(iii$_k$). 

We claim that $Z_k$ satisfies the inductive condition (iv$_k$). Given any points $x,y\in Z_k$, we need to check that $|Z_k\cap\Aline xy|\le\lambda_n$. If $\{x,y\}\subseteq\overline{D_k}$, then $|Z_k\cap\Aline xy|\le|Z_k\cap\overline{D_k}|=\lambda_n$ and we are done. So, we assume that $\{x,y\}\not\subseteq \overline{D_k}$ and hence $\{x,y\}\not\subseteq Z_k\setminus Z_{k-1}$. If $\{x,y\}\subseteq Z_{k-1}$, then $|Z_k\cap\Aline xy|=|Z_{k-1}\cap\Aline xy|\le\lambda_n$, by the inductive conditions (iii$_k$) and (iv$_k$). It remains to consider the case $\{x,y\}\cap Z_{k-1}\ne\varnothing\ne\{x,y\}\cap(Z_k\setminus Z_{k-1})$. We lose no generality assuming that $x\in Z_{k-1}$ and $y\in Z_{k}\setminus Z_{k-1}$. In this case $Z_k\cap \Aline xy=\{x,y\}$. Indeed, assuming that $Z_k\cap\Aline xy$ contains some point $z\notin\{x,y\}$, we conclude that $z\in Z_{k-1}$ or $z\in Z_k\setminus Z_{k-1}$. If $z\in Z_{k-1}$, then $y\in (Z_k\setminus Z_{k-1})\cap \Aline xz$ implies $\Aline xz\subseteq\overline{D_k}$, by the inductive condition (iii$_k$). Then 
$|Z_k\cap\Aline xy|=|Z_k\cap\Aline xz|=|Z_k\cap \overline{D_k}|=\lambda_n$. If $z\in Z_{k}\setminus Z_{k-1}$, then $\{y,z\}\subseteq Z_k\setminus Z_{k-1}\subseteq\overline{D_n}$ and hence $|Z_k\cap\Aline xy|=|Z_k\cap\Aline zy|=|Z_k\cap\overline{D_k}|=\lambda_n$. This completes the proof of the inductive condition (iv$_k$). 


After completing the inductive construction of the sequence $(Z_k)_{k=1}^m$, we obtain the finite set $Z_m$ satisfying the inductive conditions (i$_m$)--(iv$_m$), which imply that  the set $Y_n\defeq Z_m$ satisfies the inductive conditions $(2_n)$--$(4_n)$. The choice of the set $R_n$ ensures that the inductive condition $(1_n)$ is satisfied, too. It remains to check the inductive condition $(5_n)$.
Take any points $x,y\in Y_{n-1}\cup R_n$, $p\in\Aline xy\cap Y_n\setminus(Y_{n-1}\cup R_n)$, and $q\in Y_n\setminus\Aline pq$. We need to check that $\Aline pq\cap Y_n=\{p,q\}$. To derive a contradiction, assume that $\Aline pq\cap Y_n$ contains some point $z\notin \{p,q\}$. Let $k,i,j\in\{0,\dots,m\}$ be the  smallest numbers such that $p\in Z_k$, $q\in Z_i$, $z\in Z_j$. It follows from $p\notin Y_{n-1}\cup R_n=Z_0$ that $k>0$ and hence $p\in Z_k\setminus Z_{k-1}\subseteq\overline{D_k}$.  Since $x,y\in Y_{n-1}\cup R_n=Z_0$ and $p\in \Aline xy$, the inductive condition (iii$_k$) implies that $\Aline xy=\overline{D_k}$.

 If $\max\{i,j\}<k$, then 
$z,q\in Z_{k-1}$. Since $p\in (Z_k\setminus Z_{k-1})\cap\Aline qz$, the inductive condition (iii$_k$) ensures that $q\in\Aline qz=\overline{D_k}=\Aline xy$, which contradicts the choice of $q$. Therefore, $\max\{i,j\}\ge k$. 

If $i=k$, then $q\in Z_i\setminus Z_{i-1}=Z_k\setminus Z_{k-1}\subseteq\overline{D_k}=\Aline xy$, which contradicts the choice of $q$.  If $j=k$, then $z\in Z_j\setminus Z_{j-1}=Z_k\setminus Z_{k-1}\subseteq\overline{D_k}$ and $q\in\Aline pz\subseteq\overline{D_k}=\Aline xy$, which contradicts the choice of $q$. Those two contradictions show that $\max\{i,j\}>k$.

If $i<j$, then $p,q\in Z_{j-1}$ and $z\in (Z_j\setminus Z_{j-1})\cap \Aline pq$ imply $\Aline pq\subseteq\overline{D_j}$, by the inductive condition (iii$_j$). Since $p\in (Z_k\setminus Z_{k-1})\cap\overline{D_j}$, the inductive condition (iii$_k$) implies that $q\in \overline{D_j}=\Aline xy$, which contradicts the choice of the point $q$.

If $j<i$, then $p,z\in Z_{i-1}$ and $q\in (Z_i\setminus Z_{i-1})\cap\Aline pz$  imply $q\in\Aline pz\subseteq\overline{D_i}$, by the inductive condition (iii$_i$). Since $p\in (Z_k\setminus Z_{k-1})\cap\overline{D_i}$, the inductive condition (iii$_k$) implies that $q\in\overline{D_i}=\Aline xy$, which contradicts the choice of the point $q$. 

Therefore, $i=j>k$. Then $q,z\in Z_{i}\setminus Z_{i-1}\subseteq \overline{D_i}$ and hence $p\in (Z_k\setminus Z_{k-1})\cap\Aline qz=(Z_k\setminus Z_{k-1})\cap\overline{D_i}$. The inductive condition (iii$_k$) ensures that $q\in \overline{D_i}=\Aline xy$, which contradicts the choice of $q$.
This is a final contradiction showing that the inductive condition $(5_n)$ holds.

This completes the inductive step of the construction of the sets $R_n$ and $Y_n$. After completing the inductive construction, we obtains the sequences of finite sets $(R_n)_{n\in\w}$ and $(Y_n)_{n\in\w}$ satisfying the inductive conditions $(1_n)$--$(5_n)$.

It remains to prove that the subliner $Y\defeq\bigcup_{n\in\w}Y_n$ is injective and has the properties (1)--(3) of Theorem~\ref{t:injective-subliner}. This will be done in the following claims.

\begin{claim} The liner $Y$ is injective.
\end{claim}

\begin{proof} Assuming that $Y$ is not injective, we can find four lines $L_1,L_2,L_3,L_4\subseteq Y$ such that for every distinct numbers $i,j\in\{1,2,3,4\}$ there exists a unique point $y_{ij}\in Y$ such that $L_i\cap L_j=\{y_{ij}\}$ and the points $y_{12},y_{13},y_{14},y_{23},y_{24},y_{34}$ are pairwise distinct. For every distinct numbers $i,j\in\{1,2,3,4\}$, let $n_{ij}$ be the smallest number such that $y_{ij}\in Y_{n_{ij}}\cup R_{n_{ij}+1}$. Let $n=\max\{n_{12},n_{13},n_{14},n_{23},n_{24},n_{34}\}$. Since $|Y_0|=2<6$, the number $n$ is positive. We lose no generality assuming that $n=n_{12}$.  Since $\|\{y_{12},y_{13},y_{14}\}\|=3=\|\{y_{21},y_{23},y_{24}\}\|$, the inductive condition  $(1_{n+1})$ implies that $(\{y_{12},y_{13},y_{14}\}\cup\{y_{12},y_{23},y_{24}\})\cap (R_{n+1}\setminus Y_{n})=\varnothing$. Since $y_{12}\in Y_n$ and $n>0$, by the inductive condition $(2_n)$, there exist points $u,v\in Y_{n-1}\cup R_n$ such that $y_{12}\in\Aline uv$. Assuming that $y_{13}\notin\Aline uv$, we can apply the inductive condition $(5_n)$ and conclude that $y_{14}\in\Aline{y_{12}}{y_{13}}\cap Y_n=\{y_{12},y_{13}\}$, which contradicts the choice of the lines $L_1,L_2,L_3,L_4$. This contradiction shows that $y_{13}\in\Aline uv$ and hence $L_1=\Aline {y_{12}}{y_{13}}=\Aline uv$. By analogy we can prove that $L_2=\Aline uv=L_1$, which contradicts the choice of the lines $L_1,L_2$. This contradiction completes the proof of the injectivity of the liner $Y$.
\end{proof}




\begin{claim}\label{cl:rank=w} $\|Y\|=\w$.
\end{claim}

\begin{proof} It is clear that $\|Y\|\le|Y|=\w$.  Assuming that the rank of $Y$ is finite, we can find a finite set $F\subseteq Y$ whose flat hull $\overline F$ in $Y$ coincides with $Y$. 
Since $Y=\bigcup_{n\in\w}Y_n$, there exists a finite cardinal $n\ge 2$ such that $F\subseteq Y_{n-1}$.
Then the flat hull $\overline{Y_{n-1}}$ of the set $Y_{n-1}$ in $Y$ coincides with $Y$. Let $H_{n-1}\defeq Y_{n-1}$ and for every $m\ge n$, let $H_m\defeq\bigcup_{x,y\in H_{m-1}}(Y_m\cap\Aline xy)$. We claim that $\overline{Y_{n-1}}=\bigcup_{m=n}^\infty H_m$. Indeed, for every points $x,y\in \bigcup_{m=n}^\infty H_m$ and $z\in Y\cap\Aline xy$, we can find $m\ge n$ such that $x,y\in H_m$ and $z\in Y_m$. Then $z\in H_{m+1}$, witnessing that the set $\bigcup_{m=n}^\infty H_m$ is flat in the liner $Y$ and hence coincides with the flat hull $\overline{Y_{n-1}}$ of the set $Y_{n-1}$ in $Y$.

By induction we shall prove that for every $m\ge n-1$ we have $H_{m+1}\cap Y_m=H_m$. The equality $H_n\cap Y_{n-1}=Y_{n-1}=H_{n-1}$ follows from the definition of the sets $H_{n-1}\defeq Y_{n-1}$ and $H_n\defeq\bigcup_{x,y\in H_{n-1}}Y_n\cap\Aline xy\supseteq H_{n-1}=Y_{n-1}$.  

Assume that for some $m>n$, we have proved that $H_{m}\cap Y_{m-1}=H_{m-1}$. We shall prove the equality  $H_{m+1}\cap Y_m=H_m$. Assuming that this equality does not hold,  we can find a point $z\in H_{m+1}\cap Y_m\setminus H_m$. Since $z\in H_{m+1}$, there exist points $x,y\in H_m$ such that $z\in Y_{m+1}\cap\Aline xy$. Then $|\{x,y,z\}|=3$. Assuming that $x,y\in H_{m-1}$, we conclude that $z\in Y_m\cap\Aline xy\in H_m$, which contradicts the choice of $z$. This contradiction shows that $\{x,y\}\not\subseteq H_{m-1}=H_m\cap Y_{m-1}$. Then $x\notin Y_{m-1}$ or $y\notin Y_{m-1}$. We lose no generality assuming that $x\notin Y_{m-1}$. Since $x\in H_m$, there exist points $u,v\in H_{m-1}$ such that $x\in Y_m\cap \Aline uv$. Since $z\in Y_m$, there exist points $s,t\in Y_{m-1}\cup R_m$ such that $z\in\Aline st$. Assuming that $\Aline st\not\subseteq \Aline uv$, we can apply the inductive condition $(5_m)$ and conclude that $y\in Y_m\cap\Aline xz=\{x,z\}$, which contradicts $|\{x,y,z\}|=3$. This contradiction shows that $z\in \Aline st\subseteq \Aline uv$ and hence $z\in Y_m\cap\Aline uv\subseteq H_m$, which contradicts the choice of $z$. This contradiction shows that $H_{m+1}\cap Y_m=H_m$.

The inductive condition $(1_n)$ implies $|R_n\setminus R_{n-1}|=r_n-r_{n-1}=1$ and hence $R_n\setminus R_{n-1}=\{z\}$ for some point $z$. Since $z\in Y=\overline{Y_{n-1}}=\bigcup_{m=n-1}H_m$, there exists a number $m\ge n-1$ such that $z\in H_m$. We can assume that $m$ is the smallest number with $m\ge n-1$ and $z\in H_m$. The inductive condition $(1_n)$ implies that $z\notin \bigcup_{x,y\in Y_{n-1}}\Aline xy$ and hence $z\notin H_n$. Therefore, $m>n$. Since $z\in H_m$, there exist points $x,y\in H_{m-1}$ such that $z\in\Aline xy$. It follows from $z\in \Aline xy\setminus H_{m-1}$ and $x,y\in H_{m-1}$ that $|\{x,y,z\}|=3$. Assuming that $x,y\in H_{m-2}$, we conclude that $z\in Y_n\cap\Aline xy\subseteq Y_{m-1}\cap\Aline xy\subseteq H_{m-1}$, which contradicts the minimality of $m$. This contradiction shows that $x\notin H_{m-2}$ or $y\notin H_{m-2}$. We lose no generality assuming that $x\notin H_{m-2}$. It follows from $H_{m-1}\cap Y_{m-2}=H_{m-2}$ that $x\notin Y_{m-2}$. Since $x\in H_{m-1}\setminus H_{m-2}$, there exist points $u,v\in H_{m-2}\subseteq Y_{m-2}$ such that $x\in \Aline uv$. Assuming that $z\in \Aline uv$, we conclude that $z\in Y_{m-1}\cap\Aline uv\subseteq H_{m-1}$, which contradicts the minimality of $m$. This contradiction shows that $z\notin \Aline uv$ and hence $\Aline zz\not\subseteq\Aline uv$. The inductive condition $(5_{m-1})$ ensures that $y\in Y_{m-1}\cap\Aline xz=\{x,z\}$, which contradicts $|\{x,y,z\}|=3$. This contradiction completes the proof of Claim~\ref{cl:rank=w}.
\end{proof}

\begin{claim} The liner $Y$ is $2$-balanced and $|Y|_2=\lambda$.
\end{claim}

\begin{proof} For every distinct points $x,y\in Y$ and every $n\in\w$ with $x,y\in Y_{n-1}$, the inductive condition $(3_n)$ ensures that $|Y_n\cap\Aline xy|=\lambda$ and hence $|Y\cap\Aline xy|=\sup_{n\in\w}\lambda_n=\lambda$, witnessing that the liner $Y$ is $2$-balanced with $|Y|_2=\lambda$.
\end{proof}

\begin{claim} Every flat $F\subseteq Y$ of rank $\|F\|\ge 3$ is infinite and hence $|Y|_n=\w$ for every cardinal $n\ge 3$.
\end{claim} 

\begin{proof} Since $F$ has rank $\|F\|\ge 3$, there exist points $x,y,z\in F$ such that $\|\{x,y,z\}\|=3$. Find the smallest number $m\in\w$ such that $\{x,y,z\}\subseteq Y_{m}$. Since $|Y_0|=2$, the number $n$ is positive. The minimality of $n$ ensures that $\{x,y,z\}\not\subseteq Y_{m-1}$.  We lose no generality assuming that $z\notin Y_{m-1}$. Moreover, if $\{x,y,z\}\cap (R_{n}\setminus Y_{n-1})\ne\varnothing$, then we will assume that $z\in R_{m}\setminus Y_{m-1}$. Since $R_{m}\setminus Y_{m-1}$ is a singleton, $x,y\notin R_{m}\setminus Y_{n-1}$. 

 By induction we shall construct a sequence of points $(z_n)_{n=m}^\infty$ such that $z_n\in (F\cap Y_n)\setminus Y_{n-1}$ and $\|\{x,y,z_n\}\|=3$ for every $n\ge m$. We start the inductive construction letting $z_m=z$. Assume that for some $n\ge m$, a point $z_{n}\in F\cap Y_{n}\setminus Y_{n-1}$ with $\|\{x,y,z_n\}\|=3$ has been constructed. Since $z_n\in Y_n$, the inductive condition $(2_n)$ ensures that $z_n\in\Aline uv$ for some points $u,v\in Y_{n-1}\cup R_n$.

 It follows from $\|\{x,y,z_n\}\|=3$ that  $x\notin\Aline uv$ or $y\notin\Aline uv$. We lose no generality assuming that $x\notin \Aline uv$. Since $x\in Y_m\subseteq Y_n$, the inductive condition $(5_n)$ ensures that  $\Aline x{z_n}\cap Y_n=\{x,z_n\}$. By the inductive condition  $(3_{n+1})$, the intersection $Y_{n+1}\cap\Aline x{z_n}$ has cardinality $\lambda_n\ge 3$ and hence there  exists a point $z_{n+1}\in Y_{n+1}\cap\Aline x{z_n}\setminus\{x,z_n\}$. Taking into account that $F$ is a flat in $Y$ and $Y_n\cap\Aline x{z_n}=\{x,z_n\}$, we conclude that $z_{n+1}\in F\cap Y_{n+1}\setminus Y_n$.
Assuming that $\|\{x,y,z_{n+1}\}\|<3$, we conclude that $z_{n+1}\in\Aline xy$ and $z_n\in\Aline x{z_n}\subseteq \Aline xy$, which contradicts $\|\{x,y,z_n\}\|=3$. This completes the inductive step and also the inductive construction of the sequence $(z_n)_{n=m}^\infty$ in $F$. Assuming that the set $F$ is finite, we can find two numbers $n,k\in \w$ such that $m\le n<k$ and $z_n=z_k$. Then $z_n=z_k\in Y_n\cap (Y_{k}\setminus Y_{k-1})\subseteq Y_{k-1}\cap(Y_k\setminus Y_{k-1})=\varnothing$, which is a contradiction showing that the flat $F$ is infinite and $|F|=\w$. 

Since all flats of rank $\ge 3$ in $Y$ are countable, $|Y|_n=\w$ for any cardinal $n\ge 3$.
\end{proof}
\vskip-5pt
\end{proof}

\begin{question} Does the injective liner $Y$ in Theorem~\ref{t:injective-subliner} contain flats of arbitrary finite rank $r$?
\end{question}



\begin{question}\label{q:hyperaffine} Is a regular ranked liner $X$ hyperaffine if for every line $L\subseteq X$ and point $p\in X\setminus L$, there exists a line $\Lambda$ with $p\in \Lambda\subseteq\overline{\{p\}\cup L}\setminus L$?
\end{question} 

\begin{question}\label{q:hyperbolic} Is a regular ranked liner $X$ hyperbolic if for every line $L\subseteq X$ and point $p\in X\setminus L$, there exist two distinct lines $\Lambda_1,\Lambda_2$ such that $p\in \Lambda_i\subseteq\overline{\{p\}\cup L}\setminus L$ for $i\in\{1,2\}$?
\end{question}

\begin{remark} The regularity is essential in Questions~\ref{q:hyperaffine}, \ref{q:hyperbolic} because the proaffine $13$-element ranked liner from Example~\ref{ex:Steiner13} is not hyperbolic but for every line $L\subseteq X$ and point $x\in X\setminus L$ there exist at least three\footnote{Calculated on computer by Ivan Hetman.} distinct lines in $\overline{L\cup\{x\}}\setminus L=X\setminus L$ that contain $x$.
\end{remark}

\chapter{Regularity of affine and Playfair liners} 

In this chapter we shall prove two non-trivial theorems on the regularity of affine  liners: Theorem~\ref{t:4-long-affine} saying that $4$-long affine liners are regular and Theorem~\ref{t:Playfair<=>regular} saying that Playfair balanced liners are regular. To deduce Theorem~\ref{t:Playfair<=>regular} from Theorem~\ref{t:4-long-affine}, we have to consider the exeptional case of (Playfair) liners whose lines have cardinality 3. Such liners are well-known in Finite Geometry as Steiner (Hall) liners. Steiner liners carry an algebraic structure of an idempotent commutative involutory magma, which generates a structure of a loop of order 3. For Hall liners this loop is Moufang, which allows us to apply the powerful theory of commutative Moufang loops for analyzing the structure of Hall liners and establishing that the rankedness and regularity of Hall liners is equivalent to the associativity of the corresponding Moufang loop. 

\section{Regularity of $4$-long affine liners}

The following theorem explains why the length of the lines in Examples~\ref{ex:Tao} and \ref{ex:Steiner13} is $3$.

\begin{theorem}\label{t:4-long-affine} Every $4$-long affine liner is Playfair and regular.
\end{theorem}

\begin{proof} Let $X$ be a $4$-long affine liner. First we prove that the liner $X$ is $3$-regular. We separate the proof of the $3$-regularity of $X$ into a series of claims.

\begin{claim}\label{cl:4-long1} Let $o,a,b$ be points of $X$ such that $o\notin \Aline ab$. Let $\alpha\in\Aline oa$, $\beta\in\Aline ob$, $c\in\Aline ab$ be points such that $\Aline ab\cap\Aline\alpha\beta=\varnothing=\Aline oc\cap\Aline\alpha\beta$. Then $|X|_2=4$ and there exists a point $\gamma\in\Aline oc\cap\Aline \alpha b$ such that the intersections $\Aline \gamma \beta\cap\Aline cb$ and $\Aline \gamma\beta\cap\Aline oa$ are empty.
\end{claim}

\begin{picture}(200,150)(-200,-15)

\put(0,40){\line(0,1){80}}
\put(-60,0){\line(1,2){60}}
\put(60,0){\line(-1,2){60}}
\put(-100,60){\line(1,0){200}}
\put(-100,0){\line(1,0){200}}
\put(-60,0){\line(3,2){90}}
\put(60,0){\line(-3,2){60}}
\qbezier(0,40)(-28,57)(-95,57)
\put(-34,52){\circle{4}}
\put(-34,52){\color{white}\circle*{3}}

\put(0,120){\circle*{3}}
\put(-2,123){$o$}
\put(60,0){\circle*{3}}
\put(58,-10){$\beta$}
\put(-60,0){\circle*{3}}
\put(-64,-10){$\alpha$}
\put(-30,60){\circle*{3}}
\put(-36,63){$a$}
\put(30,60){\circle*{3}}
\put(31,63){$b$}
\put(0,60){\circle*{3}}
\put(3,63){$c$}
\put(0,40){\circle*{3}}
\put(-3,31){$\gamma$}
\end{picture}

\begin{proof} It follows from $c\in\Aline ab$, $\{\alpha,\beta\}\subseteq\Aline oa\cup\Aline ob$ and $\Aline oc\cap\Aline\alpha\beta=\varnothing$ that $a\ne c\ne b$. Since 
the liner $X$ is affine, there exists a unique point $\gamma\in\Aline oc$ such that $\Aline \gamma\beta\cap\Aline cb=\varnothing$ and hence $\gamma\in \Aline oc\setminus(\Aline oa\cup\Aline ob\cup\Aline ab)$. If $\gamma\in \Aline \alpha\beta$, then $\gamma\in\Aline oc\cap\Aline\alpha\beta=\varnothing$, which is a contradiction showing that $\gamma\notin\Aline \alpha\beta$. Assuming that $\Aline \gamma\beta\cap\Aline oa$ contains some point $\alpha'$, we observe that 
$$\Aline {\alpha'}\beta\cap\Aline ab=\Aline \gamma\beta\cap\Aline cb=\varnothing=\Aline\alpha\beta\cap\Aline ab$$ and hence $\alpha'=\alpha$ by the proaffinity of $X$. Then $\gamma\in\Aline {\alpha'}\beta=\Aline\alpha\beta$, which is a contradiction showing that $\Aline \gamma\beta\cap\Aline oa=\varnothing$.

It remains to prove that $\gamma\in\Aline\alpha b$ and $|X|_2=4$. Theorem~\ref{t:affine=>Avogadro} implies that the affine liner $X$ is $2$-balanded and hence the cardinal number $|X|_2$ is well-defined. Since $X$ is $4$-long, $|X|_2\ge 4$.

By the affinity of $X$, the set $I\defeq\{v\in \Aline oc:\Aline vb\cap\Aline oa=\varnothing\}$ is a singleton. Since $X$ is $4$-long, there exists a point $p\in \Aline oc\setminus(\{o,c\}\cup I)$. If $\{u,c,\gamma\}\cup I\ne \Aline oc$, then we can additionally assume that $p\ne \gamma$. Since $p\notin I$, there exists a point $q\in \Aline pb\cap\Aline oa$. Assuming that $q=a$, we conclude that $p\in\Aline bq\cap\Aline oc=\Aline ba\cap\Aline oc=\{c\}$, which contradicts the choice of $p$. This contradiction shows that $q\ne a$. Assuming that $\beta\in\Aline pb$, we conclude that $p\in \Aline b\beta\cap\Aline oc=\{o\}$, which contradicts the choice of $p$. This contradiction shows that $\beta\notin\Aline pb$.

Taking into account that $q\ne a$, $X$ is affine and $\Aline ab\cap\Aline \alpha\beta=\varnothing$, we conclude that the intersection $\Aline qb\cap\Aline \alpha\beta$ contains some point $s$. 
 The proaffinity of $X$ ensures that $\{v\in\Aline oc:\Aline vb\cap\Aline \gamma\beta=\varnothing\}=\{c\}\ne\{p\}$ and hence $\Aline pb\cap\Aline \gamma\beta$ contains some point $t$. 

\begin{picture}(200,190)(-160,-50)

\put(-60,0){\line(1,0){180}}
\put(0,120){\line(-1,-2){60}}
\put(0,120){\line(1,-2){60}}
\put(-40,40){\line(1,0){80}}
\put(40,40){\line(-2,1){64}}
\put(40,40){\line(2,-1){80}}
\put(0,120){\line(0,-1){150}}
\put(0,-30){\line(2,1){90}}
\put(60,0){\line(-1,1){60}}

\put(0,120){\circle*{3}}
\put(-2,123){$o$}
\put(0,40){\circle*{3}}
\put(-7,42){$c$}
\put(-60,0){\circle*{3}}
\put(-64,-10){$\alpha$}
\put(60,0){\circle*{3}}
\put(58,-10){$\beta$}
\put(-40,40){\circle*{3}}
\put(-47,41){$a$}
\put(40,40){\circle*{3}}
\put(42,41){$b$}
\put(0,60){\circle*{3}}
\put(3,63){$p$}
\put(0,-30){\circle*{3}}
\put(-2,-38){$\gamma$}
\put(-24,72){\circle*{3}}
\put(-32,73){$q$}
\put(120,0){\circle*{3}}
\put(118,-8){$s$}
\put(90,15){\circle*{3}}
\put(89,18){$t$}
\put(20,40){\circle*{3}}
\put(16,32){$u$}
\end{picture}

We claim that $p=\gamma$. To derive a contradiction, assume that $p\ne \gamma$. Since $\{v\in\Aline oc:\Aline v\beta\cap\Aline cb=\varnothing\}=\{\gamma\}\ne\{p\}$, the intersection $\Aline p\beta\cap\Aline cb$ is not empty and contains some point $u$. Assuming that $u=b$, we conclude that $\beta\in\Aline pu=\Aline qb=\Aline pb$, which is not true. Therefore, $u\ne b$ and $\Aline ub=\Aline cb=\Aline ab$. Since $s,t\in \Aline pb$, $u\in\Aline p\beta$ and  $$\Aline \beta s \cap\Aline ub=\Aline \alpha\beta\cap\Aline ab=\varnothing=\Aline \gamma\beta\cap\Aline cb=\Aline\beta t\cap\Aline ub,$$ the proaffinity of $X$ ensures that $s=t$ and hence $\gamma\in \Aline \beta t=\Aline\beta s=\Aline \alpha\beta$, which is a desired contradiction showing that $p=\gamma$. In this case, the choice of $p$ ensures that $\Aline oc=\{o,c,\gamma\}\cup I$ and hence $|X|_2=|\Aline oc|=4$.

Next, we show that $q=\alpha$. Indeed, assuming that $q\notin\alpha$, we conclude that $\Aline \beta q\cap\Aline ba$ contains some point $w$. It follows from $\beta\notin\Aline pb=\Aline qb$ that $w\ne b$ and hence $\Aline cb=\Aline bw$. Since $$\Aline \beta s\cap\Aline bw=\Aline\beta\alpha\cap\Aline ba=\varnothing=\Aline \beta\gamma\cap\Aline ca=\Aline\beta\gamma\cap\Aline bw,$$ the proaffinity of $X$ implies that $\gamma=s\in\Aline \alpha\beta$, which is a contradiction. This contradiction shows that $q=\alpha$ and hence $\gamma=p\in \Aline qp=\Aline\alpha b$.
\end{proof}

\begin{claim}\label{cl:4-long2} Let $o,a,b,c,\alpha,\beta$ be points of $X$ such that $o\notin \Aline ab$, $\alpha\in\Aline oa$, $\beta\in\Aline ob$, and $c\in\Aline ab$. If $\Aline ab\cap\Aline \alpha\beta=\varnothing=\Aline oc\cap\Aline\alpha\beta$, then there exist points $a',b',\gamma,\delta\in X$ such that $\Aline oa=\{o,a,\alpha,a'\}$, $\Aline ob=\{o,b,\beta,b'\}$, $\Aline oc=\{o,c,\gamma,\delta\}$, $\gamma\in \Aline \alpha b\cap\Aline a{b'}$, $\delta\in \Aline a\beta\cap\Aline {a'}b$ and the intersection $\Aline \gamma\beta\cap\Aline cb$ is empty.
\end{claim}

\begin{proof} By Claim~\ref{cl:4-long1}, $|X|_2=4$ and there exist points $\gamma\in\Aline oc\cap\Aline \alpha b$ and $\delta\in\Aline oc\cap\Aline \beta a$ such that the sets $\Aline \gamma\beta\cap\Aline cb$, $\Aline \gamma\beta\cap\Aline oa$, $\Aline \delta\alpha\cap\Aline ca$, $\Aline \delta\alpha\cap\Aline ob$ are empty. Assuming that $\gamma=\delta$, we conclude that $b\in \Aline\alpha\gamma\cap\Aline ob=\Aline\alpha\delta\cap\Aline ob=\varnothing$, which is a contradiction showing that $\gamma\ne\delta$ and hence $\Aline oc=\{o,c,\gamma,\delta\}$. Since $|X|_2=4$, there exists points $a',b'\in X$ such that $\Aline oa=\{o,a,\alpha,a'\}$ and $\Aline ob\setminus\{o,b,\beta,b'\}$. By the affinity, there exists a unique point $b''\in\Aline ob$ such that $\Aline {a'}{b''}\cap\Aline ab=\varnothing$.
The latter condition implies $b''\notin\{o,a\}$. Assuming that $b''=\beta$, we obtain that $a',\alpha$ are two distinct points on the line $\Aline oa$ such that $\Aline {a'}\beta\cap\Aline ab=\varnothing=\Aline\alpha\beta\cap\Aline ab$, which contradicts the proaffinity of $X$. This contradiction shows that $b''\ne\beta$ and hence $b''\in\Aline ob\setminus\{o,b,\beta\}=\{b'\}$. Therefore, $\Aline{a'}{b'}\cap\Aline ab=\varnothing$. By analogy we can prove that $\Aline{a'}{b'}\cap\Aline\alpha\beta=\varnothing$. 

Assuming that $\Aline {a'}{b'}\cap\Aline oc\ne\varnothing$, we conclude that either $\gamma\in\Aline {a'}{b'}$ or $\delta\in\Aline{a'}{b'}$. If $\gamma\in\Aline {a'}{b'}$, then  
$$\{b',\beta\}\subseteq\{v\in\Aline ob:\Aline v\gamma\cap\Aline cb=\varnothing\}$$ and by the proaffinity of $X$, $b'=\beta$, which contradicts the choice of $b'$. If $\delta\in\Aline {a'}{b'}$, then $\{a',\alpha\}\subseteq\{v\in\Aline oa:\Aline v\delta\cap\Aline ac=\varnothing\}$ and by the proaffinity of $X$, $a'=\alpha$, which contradicts the choice of $a'$. Those contradictions show that $\Aline {a'}{b'}\cap\Aline oc=\varnothing$. 

By Claim~\ref{cl:4-long1}, there exist points $\gamma'\in\Aline oc\cap\Aline{a'}{b}$ and $\delta'\in\Aline oc\cap\Aline{b'}{a}$ such that the sets $\Aline {\gamma'}{b'}\cap\Aline cb$, $\Aline {\gamma'}{b'}\cap\Aline oa$, $\Aline{\delta'}{a'}\cap\Aline ca$, $\Aline{\delta'}{a'}\cap\Aline ob$ are empty. It follows that $\Aline oc=\{o,c,\gamma,\delta\}=\{o,c,\gamma',\delta'\}$ and hence $(\gamma,\delta)=(\gamma',\delta')$ or $(\gamma,\delta)=(\delta',\gamma')$. 

If $(\gamma,\delta)=(\gamma',\delta')$, then  $\Aline oc\cap\Aline{a'}{b}=\{\gamma'\}=\{\gamma\}=\Aline oc\cap\Aline\alpha b$ implies $\{a'\}=\Aline oa\cap\Aline b{\gamma'}=\Aline oa\cap\Aline b{\gamma}=\{\alpha\}$, which contradicts the choice of $a'\ne\alpha$.

Therefore, $(\gamma,\delta)=(\delta',\gamma')$. In this case, the points $a',b'$ and $\gamma=\delta'$, $\delta=\gamma'$ have the required properties.
\end{proof} 

\begin{claim}\label{cl:4-long3} Let $o,a,b,c,\alpha,\beta$ be points of $X$ such that $o\notin \Aline ab$, $\alpha\in\Aline oa$, $\beta\in\Aline ob$, and $c\in\Aline ab$. If $\Aline ab\cap\Aline \alpha\beta=\varnothing$, then $\Aline oc\cap\Aline\alpha\beta\ne\varnothing$.
\end{claim}

\begin{proof} To derive a contradiction, assume that $\Aline ab\cap\Aline \alpha\beta=\varnothing=\Aline oc\cap\Aline\alpha\beta$. By Claim~\ref{cl:4-long2}, there exist points $a',b',\gamma,\delta\in X$ such that 
$\Aline oa=\{o,a,\alpha,a'\}$, $\Aline ob=\{o,b,\beta,b'\}$, $\Aline oc=\{o,c,\gamma,\delta\}$, $\gamma\in \Aline \alpha b\cap\Aline a{b'}$, $\delta\in \Aline a\beta\cap\Aline {a'}b$ and $\Aline \gamma \beta\cap\Aline ab=\varnothing$.

Since $X$ is $4$-long, there exists a point $e\in \Aline ab\setminus\{a,b,c\}$. If $\Aline oe\cap\Aline \alpha\beta=\varnothing$, then by Claim~\ref{cl:4-long2}, $\Aline oe=\{o,e,\gamma',\delta'\}$ for some points $\gamma'\in \Aline \alpha b\cap\Aline a{b'}=\{\gamma\}$, $\delta'\in\Aline a\beta\cap\Aline{a'}b=\{\delta\}$, which implies that $\Aline oe=\Aline o{\gamma'}=\Aline o\gamma=\Aline oc$ and $e=c$. But this contradicts the choice of $e$. This contradiction shows that the intersection $\Aline oe\cap\Aline\alpha\beta$ contains some point $\varepsilon$. Since $\Aline oc\cap\Aline \varepsilon\beta=\Aline oc\cap\Aline\alpha\beta=\varnothing$, Claim~\ref{cl:4-long2} yields points $e',\gamma',\delta'\in X$ such that $\Aline oe=\{o,e,\e,e'\}$, $\Aline oc=\{o,c,\gamma',\delta'\}$, $\gamma'\in\Aline \e b\cap\Aline e{b'}$ and $\delta'\in \Aline e\beta\cap\Aline {e'}b$.   It follows from $\{o,c,\gamma,\delta\}=\Aline oc=\{o,c,\gamma',\delta'\}$ that $(\gamma,\delta)=(\gamma',\delta')$ or $(\gamma,\delta)=(\delta',\gamma')$. 

If $(\gamma,\delta)=(\gamma',\delta')$, then $\{\e\}=\Aline \alpha\beta\cap\Aline {\gamma'}b=\Aline\alpha\beta\cap \Aline\gamma b=\{\alpha\}$ and hence $e\in \Aline ab\cap\Aline o\e=\Aline ab\cap\Aline o\alpha=\{a\}$, which contradicts the choice of $e$. This contradiction shows that $(\gamma,\delta)=(\delta',\gamma')$ and then $\gamma=\delta'\in\Aline e\beta\cap\Aline{e'}b$ and hence $e\in\Aline \gamma\beta\cap\Aline ab=\varnothing$, which is a final contradiction completing the proof of the claim.
\end{proof}

\begin{claim}\label{cl:4-long4} The liner $X$ is $3$-regular.
\end{claim}

\begin{proof} Given two concurrent lines $A,\Lambda$ in $X$, we need to show that the set $\Aline A\Lambda$ is flat. In the opposite case, we can find points $a,b\in\Lambda$ and $c\in\Aline ab\setminus\Aline A\Lambda$. Then $\{a,b\}\not\subseteq A\cup\Lambda$ and we lose no generality assuming that $b\notin A\cup \Lambda$. Let $o$ be the unique point of the intersection $A\cap\Lambda$. Since $X$ is proaffine, the sets $I_a\defeq \{v\in\Lambda:\Aline va\cap A=\varnothing\}$ and $I_b\defeq \{v\in\Lambda:\Aline vb\cap A=\varnothing\}$ have cardinality $\max\{|I_a|,|I_b|\}\le 1$. Since $X$ is $4$-long, there exists a point $\lambda\in \Lambda\setminus(I_a\cup I_b\cup L)$. Since $\lambda\notin I_a\cup I_b$, there exist points $\alpha\in\Aline \lambda a\cap A$ and $\beta\in\Aline \lambda b\cap A$.

If there exists a point $\gamma\in\Aline \lambda c\cap\Aline\alpha\beta\subseteq\Aline \alpha\beta\subseteq A$, then $c\in\Aline \gamma\lambda\subseteq \Aline A\Lambda$, which contradicts the choice of $c$. This contradiction shows that $\Aline \lambda c \cap\Aline \alpha\beta=\varnothing$. By  Claim~\ref{cl:4-long3}, there exists a point $p\in\Aline ab\cap\Aline\alpha\beta\subseteq \Aline ab\cap A$. Assuming that $p=\beta$, we conclude that $a\in \Aline pb=\Aline\beta b=\Aline \lambda\beta$ and $c\in\Aline ab\subseteq \Aline \beta\lambda\subseteq\Aline A\Lambda$, which contradicts the choice of $c$. This contradiction shows that $p\ne\beta$ and hence $\Aline p\beta=A$. 

If $\Lambda\cap\Aline ab$ contains some point $q$, then $\{v\in \Aline \lambda q:\Aline vc\cap\Aline po=\varnothing\}=\{\lambda\}$, by the proaffinity of $X$. Since $X$ is $3$-long, there exists a point $v\in \Aline \lambda q\setminus\{\lambda,o\}$. For this point, the intersection $\Aline vc\cap\Aline po$ contains some point $s$. Since $v\in\Aline \lambda q\subseteq \Aline \lambda o\subseteq\Lambda$ and $s\in\Aline po\subseteq\Aline \alpha\beta\subseteq A$, the point $c\in \Aline sv$ belongs to the set $\Aline A\Lambda$, which contradicts the choice of $c$. 
This contradiction shows that $$\Lambda\cap\Aline ab=\varnothing.$$

\begin{picture}(200,150)(-180,-15)

\put(-60,0){\line(1,0){120}}
\put(0,120){\line(-1,-2){60}}
\put(0,120){\line(1,-2){60}}
\put(-15,0){\line(1,2){37.5}}
\put(0,120){\line(0,-1){120}}
\put(60,0){\line(-1,1){80}}

\put(0,120){\circle*{3}}
\put(-2,123){$\lambda$}
\put(22.5,75){\circle*{3}}
\put(26,78){$b$}
\put(-60,0){\circle*{3}}
\put(-64,-10){$o$}
\put(60,0){\circle*{3}}
\put(58,-10){$\beta$}
\put(0,30){\circle*{3}}
\put(-8,28){$a$}
\put(10,50){\circle*{3}}
\put(13,48){$c$}
\put(-20,80){\circle*{3}}
\put(-25,83){$v$}
\put(-15,0){\circle*{3}}
\put(-18,-8){$p$}
\put(0,0){\circle*{3}}
\put(-3,-8){$\alpha$}
\end{picture}

Then $o\in A=\Aline\beta p$ and $\Aline pb\cap\Aline \lambda o=\Aline ab\cap\Lambda=\varnothing$ imply $\Aline \beta c\cap\Lambda\ne\varnothing$, by Claim~\ref{cl:4-long3}. Let $v$ be the unique point of the intersection $\Aline \beta c\cap\Lambda$. Assuming that $v=\beta$, we conclude that $b\in\Aline \lambda\beta\cap\Aline ab=\Aline \lambda v\cap\Aline ab\subseteq\Lambda\cap\Aline ab=\varnothing$, which is a contradiction showing that $v\ne\beta$ and hence $c\in \Aline \beta v\subseteq \Aline A\Lambda$, which contradicts the choice of the point $c$. This is the final contradiction showing that the set $\Aline A\Lambda$ is flat and the liner $X$ is $3$-regular. 
\end{proof}

By Claim~\ref{cl:4-long4} and Theorem~\ref{t:Playfair<=>}, the $3$-regular $4$-long affine liner $X$ is Playfair and by Proposition~\ref{p:k-regular<=>2ex}, the $3$-regular liner $X$ is $3$-ranked. To prove that the liner $X$ is regular, fix any flat $A\subseteq X$, point $o\in A$ and a line $\Lambda\subseteq X$ such that $A\cap\Lambda=\{o\}$. We have to show that $\overline{A\cup\Lambda}=\Aline A\Lambda$. This equality will follow as soon as we prove that the set $\Aline A\Lambda$ is flat.

\begin{claim}\label{cl:4-long5} For every points $x,y\in \Aline A\Lambda$, if $\Aline xy\cap A=\varnothing$, then $\Aline xy\subseteq\Aline A\Lambda$.
\end{claim}

\begin{proof} Assume that $\Aline xy\cap A=\varnothing$. Given any point $z\in\Aline xy$, we should prove that $z\in\Aline A\Lambda$. Since $x,y\in\Aline A\Lambda$, there exist points  $a,b\in A$ and $u,v\in\Lambda$ such that $x\in\Aline au$ and $y\in\Aline bv$. Since the liner $X$ is proaffine, the sets $I\defeq\{\lambda\in \Aline ou:\Aline\lambda x\cap \Aline oa=\varnothing\}\subseteq\Lambda$ and $J\defeq\{\lambda\in\Aline ov:\Aline\lambda y\cap\Aline ob=\varnothing\}\subseteq \Lambda$ have cardinality $\max\{|I|,|J|\}\le 1$. Since the liner $X$ is $4$-long, there exists a point $\lambda\in \Lambda\setminus(\{o\}\cup I\cup J)$. By the choice of $\lambda\notin I\cup J$, there exist points $\alpha\in\Aline \lambda x\cap \Aline oa\subseteq A$ and $\beta\in\Aline\lambda y\cap\Aline ob\subseteq A$. If $\lambda\in \Aline xy$, then $\alpha\in\Aline \lambda x\subseteq \Aline xy$, which contradicts $\Aline xy\cap A=\varnothing$. This contradiction shows that $\lambda\notin\Aline xy$. Since $\Aline xy\cap\Aline \alpha\beta\subseteq \Aline xy\cap A=\varnothing$, Claim~\ref{cl:4-long3} ensures that $\Aline \lambda z\cap\Aline\alpha\beta\ne\varnothing$ and hence $z\in\Aline A\lambda\subseteq \Aline A\Lambda$.
\end{proof} 

\begin{claim}\label{cl:4-long6} For every points $x\in A$ and $y\in \Aline A\Lambda$, we have $\Aline xy\subseteq \Aline A\Lambda$.
\end{claim}

\begin{proof} To derive a contradiction, assume that $\Aline xy\not\subseteq\Aline A\Lambda$ and find a point $z\in\Aline xy\setminus\Aline A\Lambda$. If $y\in A$, then $z\in\Aline xy\subseteq A\subseteq\Aline A\Lambda$, which contradicts the choice of $z$. This contradiction shows that $y\notin A$. If $y\in \Lambda$, then $z\in \Aline xy\subseteq\overline{\{x\}\cup\Lambda}$ and hence $z\in\Aline {\Aline ox}\Lambda\subseteq \Aline A\Lambda$, by the $3$-regularity of the liner $X$. By this contradicts the choice of $z$. This contradiction shows that $y\notin\Lambda$.

Since $y\in\Aline A\Lambda$, there exist points $b\in A$ and $\lambda\in \Lambda$ such that $y\in\Aline  b\lambda$. It follows from $y\notin A\cup \Lambda$ that $b\ne o\ne\lambda$ and hence $\Lambda=\Aline o\lambda$. The $3$-regularity of the liner $X$ ensures that $\overline{\Lambda\cup\{b\}}=\Aline {\Aline ob}\Lambda\subseteq\Aline A\Lambda$. If $x\in\overline{\Lambda\cup\{b\}}$, then $z\in\Aline xy\subseteq \overline{\Lambda\cup\{b\}}\subseteq\Aline A\Lambda$, which contradicts the choice of $z$. This contradiction shows that  $x\notin\overline{\Lambda\cup\{b\}}$.

Consider the plane $P\defeq\overline{\{o,x,y\}}$ and observe that $z\in\Aline xy\subseteq P$. Since $X$ is Playfair, there exists a unique line $L$ such that $z\in L\subseteq P\setminus \Aline ox$. Since $X$ is $4$-long, there exists point $p\in\Aline xy\setminus\{x,y,z\}$. Since $L$ and $\Aline ox$ are disjoint lines in the plane $P$ in the Playfair space $X$ and  $\Aline op\ne\Aline ox\ne\Aline oy$,  there exist points $u\in\Aline oy\cap L$ and $q\in\Aline op\cap L$. Assuming that $u=q$, we conclude that $\{y\}=\Aline ou\cap\Aline xy=\Aline oq\cap\Aline xy=\{p\}$, which contradicts the choice of $p$. This contradiction shows that $u\ne q$ and hence $L=\Aline uq$.

\begin{picture}(120,140)(-180,-45)

\put(0,0){\color{red}\line(-1,0){60}}
\put(0,0){\line(0,1){80}}
\put(-3,83){$\Lambda$}
\put(0,0){\line(1,-1){30}}
\put(-60,0){\color{blue}\line(3,-1){90}}
\put(30,-30){\color{blue}\line(-1,3){30}}
{\linethickness{1pt}
\put(-60,0){\color{magenta}\line(7,3){105}}
}
\put(0,60){\color{blue}\line(3,-1){45}}
\put(45,45){\color{red}\line(-1,0){95}}
\put(-60,42){\color{red}$L$}
\put(0,0){\color{red}\line(-1,3){15}}
\put(0,0){\color{red}\line(1,3){15}}
\put(0,60){\color{blue}\line(-1,-5){15}}
\put(0,0){\line(-1,-1){15}}

\put(0,0){\color{red}\circle*{3}}
\put(3,-1){$o$}
\put(-60,0){\color{magenta}\circle*{3}}
\put(-69,-2){\color{magenta}$x$}
\put(30,-30){\color{blue}\circle*{3}}
\put(33,-35){\color{blue}$b$}
\put(0,60){\color{blue}\circle*{3}}
\put(2,62){\color{blue}$\lambda$}
\put(45,45){\color{magenta}\circle*{3}}
\put(48,43){\color{magenta}$z$}
\put(15,45){\color{red}\circle*{3}}
\put(11,48){\color{red}$u$}
\put(-15,45){\color{red}\circle*{3}}
\put(-17,50){\color{red}$q$}
\put(-15,-15){\color{blue}\circle*{3}}
\put(-19,-24){\color{blue}$a$}
\put(-7.5,22.5){\color{magenta}\circle*{3}}
\put(-16,23){\color{magenta}$p$}
\put(10,30){\color{magenta}\circle*{3}}
\put(14,25){\color{magenta}$y$}
\end{picture}

  Observe that $\Aline ox\subseteq P\cap A$. Assuming that $\Aline ox\ne P\cap A$ and using the $3$-rankedness of $X$, we conclude that $P\cap A=P$ and hence $z\in P\subseteq A$, which contradicts the choice of $z$. This contradiction shows that $P\cap A=\Aline ox$ and hence $L\cap A=(L\cap P)\cap A=L\cap(P\cap A)=L\cap\Aline ox=\varnothing$.

Next, consider the plane $\overline{\{x,b,\lambda\}}$ and observe that $z\in \Aline xy\subseteq\overline{\{x,b,\lambda\}}$. If there exists a point $c\in\Aline \lambda z\cap\Aline xb$, then $z\in \Aline c\lambda\subseteq\Aline A\Lambda$, which contradicts the choice of $z$. This contradiction shows that $\Aline \lambda z\cap \Aline xb=\varnothing$. It follows from $z\notin\Lambda$ that $\Aline \lambda z$ is a line in the plane $\overline{\{x,b,\lambda\}}$, parallel to the line $\Aline xb\subseteq A$. Assuming that $\lambda\in \Aline xz$, we conclude that $\{\lambda\}\in\Aline \lambda z\cap\Aline xz=\{z\}$ and hence $z=\lambda\in\Lambda\subseteq\Aline A\Lambda$, which contradicts the choice of $z$. This contradiction shows that $\lambda\notin\Aline xz$. Since $\Aline \lambda z$ and $\Aline xb$ are parallel lines in the Playfair liner $X$ and $\Aline \lambda p\ne \Aline \lambda z$, there exists a point $a\in \Aline \lambda p\cap \Aline xb$. Then $p\in\Aline a\lambda \subseteq\overline{\{o,a,\lambda\}}$. The $3$-regularity of $X$ ensures that $q\in \Aline op\subseteq\overline{\{o,a,\lambda\}}=\Aline {\Aline oa}{\Aline o\lambda}=\Aline A\Lambda$. By the same reason, $u\in\Aline oy\subseteq\overline{\{o,b,\lambda\}}=\Aline{\Aline ob}{\Aline o\lambda}\subseteq\Aline A\Lambda$.  Since the line $\Aline uq=L$ is disjoint with the flat $A$ and $\{u,q\}\subseteq\Aline A\Lambda$, Claim~\ref{cl:4-long5} ensures that $z\in L=\Aline uq\subseteq\Aline A\Lambda$, which contradicts the choice of $z$. This is a final contradiction completing the proof of Claim~\ref{cl:4-long6}.
\end{proof}

Now we are ready to prove that the set $\Aline A\Lambda$ is flat. Given any points $x,y\in\Aline A\Lambda$, we need to show that $\Aline xy\subseteq\Aline A\Lambda$. If $x=y$, then $\Aline xy=\{x\}\subseteq\Aline A\Lambda$ and we are done. So, assume that $x\ne y$. If $\Aline xy\cap A=\varnothing$, then the inclusion $\Aline xy\subseteq\Aline A\Lambda$ follows from Claim~\ref{cl:4-long5}. So, assume that $\Aline xy\cap A$ contains some point $\alpha$. Since $x\ne y$, either $x\ne\alpha$ or $y\ne\alpha$. We lose no generality assuming that $y\ne\alpha$. Then $\Aline xy=\Aline \alpha y\subseteq \Aline A\Lambda$, by Claim~\ref{cl:4-long6}.
\end{proof} 

\begin{exercise} For every cardinal $\kappa$ construct a non-regular $\kappa$-long Proclus liner $X$ of rank $\|X\|=4$.
\smallskip

{\em Hint:} Take any $(\kappa+1)$-long projective liner $Y$ of rank $\|Y\|=4$, choose a line $L\subseteq Y$ and a point $x\in L$. Then the subliner $X\defeq Y\setminus(L\setminus\{x\})$ is nonregular, $\kappa$-long, Proclus, and has rank $\|X\|=\|Y\|=4$, according to Theorem~\ref{t:regular<=>flat4}, Proposition~\ref{p:projective-minus-proflat} and Corollary~\ref{c:procompletion-rank}.
\end{exercise}

\begin{question} Is every $4$-long proaffine liner Proclus?
\end{question}

\begin{question} Assume that a liner $X$ is $\kappa$-parallel and $(3+\kappa)$-long for some finite cardinal $\kappa$. Is $X$ $3$-regular?
\end{question}

The following famous example of \index[person]{Hall}Hall\footnote{{\bf Marshall Hall Jr.} (1910 -- 1990) was an American mathematician who made significant contributions to group theory and combinatorics. He wrote a number of papers of fundamental importance in group theory, including his solution of Burnside's problem for groups of exponent 6, showing that a finitely generated group in which the order of every element divides 6 must be finite. His work in combinatorics includes an important paper of 1943 on projective planes, which for many years was one of the most cited mathematics research papers. In this paper he constructed a family of non-Desarguesian planes which are known today as Hall planes. He also worked on block designs and coding theory.} shows that Theorem~\ref{t:4-long-affine} cannot be generalized to $3$-long affine liners.

\begin{example}[Hall, 1960]\label{ex:HTS} There exists a Playfair liner, which is not $4$-ranked.
\end{example}

\begin{proof} Let $\IF_3$ be the $3$-element field, $X\defeq\IF_3^4$ be the vector space over $\IF_3$ with origin $\boldsymbol o$ and the standard basis $\boldsymbol e_0,\boldsymbol e_1,\boldsymbol e_2,\boldsymbol e_3$. Let  $\mathcal L$ be the family of all 3-elements sets $\{\boldsymbol x,\boldsymbol y,\boldsymbol z\}\subseteq \IF_3^4$ such that $\boldsymbol x+\boldsymbol y+\boldsymbol z=(x_1-y_1)(x_2y_3-x_3y_2)\boldsymbol e_0$, where $(x_0,x_1,x_2,x_3)$ and $(y_0,y_1,y_2,y_3)$ are the coordinates of the vectors $\boldsymbol x,\boldsymbol y$ in the basis $\boldsymbol e_0,\boldsymbol e_1,\boldsymbol e_2,\boldsymbol e_3$. The family $\mathcal L$ determines the line relation
$$\Af=\{(x,y,z)\in X^3:x=y=z\}\cup\{(x,y,z)\in X^3:\{x,y,z\}\in\mathcal L\}.$$
It can be shown that the liner $(X,\Af)$ is affine, $3$-regular and $3$-long (so Playfair), but not $4$-ranked. This liner contains pairwise disjoint lines ${\color{teal}L_1}=\boldsymbol e_1+\boldsymbol e_2+\IF_3\boldsymbol e_3$, ${\color{red}L_2}=2\boldsymbol e_1+2\boldsymbol e_2+\IF_3\boldsymbol e_3$,
${\color{blue}L_3}=2\boldsymbol e_1+\IF_3\boldsymbol e_3$ such that the set $\Aline{\color{teal}L_1}{\color{red}L_2}$ is a plane, equal to the union of three disjoint lines $\color{teal}L_1$, $\color{red}L_2$ and ${\color{purple}L_4}\defeq\IF_3(\boldsymbol e_0+\boldsymbol e_3)$. On the other hand, the set $\Aline {\color{teal}L_1}{\color{blue}L_3}$ is the union of the lines $\color{teal}L_1$, $\color{blue}L_2$ and the plane ${\color{cyan}P}\defeq\IF_3\boldsymbol e_0+2\boldsymbol e_2+\IF_3\boldsymbol e_3$. Consequently, $\overline{{\color{teal}L_1}\cup{\color{blue}L_3}}=X$ and $\dim(X)=3$. On the other hand, $\{0\}\times\IF_3^3$ is a proper $3$-dimensional flat in $X$, witnessing that the liner $X$ is not $4$-ranked and hence is not regular. The plane $\overline{{\color{teal}L_1}\cup{\color{red}L_2}}$ has one-point intersection with the plane $\IF_3\boldsymbol e_2+\IF_3\boldsymbol e_4$, which is impossible in 3-dimensional regular spaces, see Theorem~\ref{t:w-modular<=>}.


\end{proof}


\section{Steiner liners}

By Theorem~\ref{t:4-long-affine}, every $4$-long affine liner is regular. If an affine liner $X$ is not $4$-long, then all lines in $X$ have the same length $|X|_2\in\{2,3\}$. 
Finite $2$-balanced liners $X$ with $|X|_2=3$ are well-known in Combinatorics as Steiner Triple Systems. This motivates the following definition.

\begin{definition} A liner $X$ is called \index[person]{Steiner}\defterm{Steiner}\footnote{{\bf Jakob Steiner} (1796  -- 1863) was a Swiss mathematician who worked primarily in geometry.  This he treated synthetically, to the total exclusion of analysis, which he hated, and he is said to have considered it a disgrace to synthetic geometry if equal or higher results were obtained by analytical geometry methods. In his own field he surpassed all his contemporaries. His investigations are distinguished by their great generality, by the fertility of his resources, and by the rigor in his proofs. He has been considered the greatest pure geometer since Apollonius of Perga.} if $X$ is $2$-balanced with $|X|_2=3$.
\end{definition}

Therefore, all lines in Steiner liners contain exactly three points. 

\begin{exercise} Prove that every Steiner liner is proaffine.
\end{exercise}

\begin{exercise} Show that a Steiner liner is affine if and only if it is injective.
\end{exercise}

Theorem~\ref{t:injective-subliner} implies the following corollary supplying us with many examples of infinite affine Steiner liners.

\begin{corollary} Let $r\in\{3,\w\}$ Every $\w$-long plane $P$ contains an infinite affine Steiner subliner $X$ of rank $\|X\|=r$. 
\end{corollary}

On the other hand, Theorem~\ref{t:2-balanced-dependence} implies that the following corollary imposing restrictions on possible cardinalities of finite Steiner liners.\index[person]{Kirkman}

\begin{corollary}[Kirkman\footnote{{\bf Thomas Penyngton Kirkman} (1806 -- 1895) was a British mathematician and ordained minister of the Church of England. Despite being primarily a churchman, he maintained an active interest in research-level mathematics, and was listed by Alexander Macfarlane as one of ten leading 19th-century British mathematicians. In the 1840s, he obtained an existence theorem for Steiner triple systems that founded the field of combinatorial design theory, while the related Kirkman's schoolgirl problem is named after him.}, 1846]\label{c:Kirkman} Every finite Steiner liner $X$ has cardinality $$|X|\in (6\IN+1)\cup(6\IN-3).$$
\end{corollary} 

\begin{proof} By Theorem~\ref{t:2-balanced-dependence}, there exists a finite cardinal $k$ such that $$|X|-1=(|X|_2-1)k=2k\quad\mbox{and}\quad 3|\mathcal L|=|X|_2\cdot k=|X|\cdot k=(2k+1)k,$$ where $\mathcal L$ is the family of lines in $X$. Then either $k$ or $2k+1$ is divisible by $3$. If $k\in 3\IN$, then $|X|=2k+1\in 6\IN+1$. If $k\notin3\IN$, then $2k+1\in 3m$ for some $m\in\IN$ and then $m$ is odd and hence $2k+1\in 3m\in 3(2\IN-1)=6\IN-3$.
\end{proof}

\begin{remark} The necessary condition given by Corollary~\ref{c:Kirkman} is also sufficient: for every number $n\in (6\IN+1)\cup(6\IN-3)$ there exists a Steiner liner of cardinality $n$, see \cite{Bose1939}, \cite{Skolem1958}, \cite{RCW1971}. Moreover, if $n\in (6\IN+13)\cup(6\IN-3)$, then there exists an affine Steiner liner of cardinality $n$, see \cite{GGW}.
\end{remark}

\begin{theorem}\label{t:Steiner-3reg} Every Steiner $3$-regular liner $X$ is $p$-parallel for some $p\in\{0,1\}$. Consequently, $X$ is projective or affine.
\end{theorem}

\begin{proof} By Proposition~\ref{p:k-regular<=>2ex}, the $3$-regular liner $X$ is $3$-ranked.
We claim that every plane $\Pi\subseteq X$ has cardinality $|\Pi|\le 9$. Choose any points $o\in\Pi$, $a\in \Pi\setminus\{o\}$ and $b\in\Pi\setminus\Aline oa$.  Since $X$ is Steiner, $\Aline oa=\{o,a,a'\}$ and  $\Aline ob=\{o,b,b'\}$ for some points $a',b'\in \Pi$.

The $3$-regularity of the liner $X$ ensures that $$\Pi=\overline{\{o,a,b\}}=\bigcup_{x\in\overline{oa}}\bigcup_{y\in\overline{ob}}\Aline xy=\Aline oa\cup\Aline ob\cup\bigcup_{x\in\{a,a'\}}\bigcup_{y\in\{b,b'\}}(\Aline xy\setminus\{x,y\})$$ and hence $|\Pi|\le 5+4=9$. 

The $3$-rankedness of $X$ implies that the liner $\Pi$ is $3$-balanced. By Proposition~\ref{p:23-balance=>k-parallel}, the $3$-balanced Steiner liner $\Pi$ is $p$-parallel for some $p\in\w$ such that $|X|_3=1+(|X|_2-1)(|X_2|+p)=1+2(3+p)=2p+7$. It follows from $2p+7=|X|_3=|\Pi|\le 9$ that $p\in\{0,1\}$. If $p=0$, then the liner $X$ is projective, by Theorem~\ref{t:projective<=>}. If $p=1$, then the liner $X$ is affine, by Theorem~\ref{t:Playfair<=>}.
\end{proof}

 Every Steiner liner $X$ carries the commutative binary operation $\circ:X\times X\to X$ assigning to every pair $(x,z)\in X\times X$ the unique point $y\in X$ such that $\{x,y,z\}=\Aline xz$. This operation will be called the \index{midpoint operation}\index{operation!midpoint}\defterm{midpoint operation} on $X$. The midpoint operation satisfies the identities 
$$ x\circ x=x,\quad x\circ y=y\circ x,\quad (x\circ y)\circ y=x,$$
turning $X$ into an involutory idempotent commutative magma. 
\smallskip

Now we recall the definitions of some algebraic structures, starting with the most basic notion of a magma. A \index{magma}\defterm{magma} is a set $X$ endowed with a binary operation $\cdot:X\times X\to X$, $\cdot:(x,y)\mapsto xy$.

A magma $X$ is 
\begin{itemize}
\item \index{commutative magma}\index{magma!commutative}\defterm{commutative} if $xy=yx$ for every $x,y\in X$;
\item \index{associative magma}\index{magma!associative}\defterm{associative} if $(xy)z=x(yz)$ for every $x,y,z\in X$;
\item \index{idempotent magma}\index{magma!idempotent}\defterm{idempotent} if $xx=x$ for every $x\in X$;
\item \index{involutory magma}\index{magma!involutory}\defterm{involutory} if $x(xy)=y$ for every $x,y\in X$;
\item \index{self-distributive magma}\index{magma!self-distributive}\defterm{self-distributive} if $x(yz)=(xy)(xz)$;
\item \index{unital magma}\index{magma!unital}\defterm{unital} if there exists an element $e\in X$ such that $xe=x=ex$ for every $x\in X$;
\item a \index{quasigroup}\index{magma!quasigroup}\defterm{quasigroup} if for every $a,b\in X$ there exist unique $x,y\in X$ such that $ax=b=ya$;
\item a \index{loop}\index{magma!loop}\defterm{loop} if $X$ is a unital quasigroup;
\item a \index{group}\index{magma!group}\defterm{group} of $X$ is an associative loop;
\item \index{Moufang}\index{magma!}\defterm{Moufang} if $(xy)(zx)=(x(yz))x=x((yz)x)$ for every $x,y\in X$.
\end{itemize}

A subset $S$ of a magma $X$ is called a \index{submagma}\defterm{submagma} of the magma $X$ if $\{xy:x,y\in S\}\subseteq S$. A subset $S$ of a loop $X$ is called a \index{subloop}\defterm{subloop} of $X$ if for every $a,b\in S$, $ab\in S$ and there exist points $x,y\in S$ such that $ax=b=ya$.

\begin{theorem}\label{t:Steiner<=>midpoint} Any Steiner liner endowed with the midpoint operation is an involutory idempotent commutative magma. Conversely, every involutory idempotent commutative magma $(M,\circ)$ endowed with the line relation
$$\Af\defeq\{(x,y,z)\in M^3:y\in\{x,z,x\circ z\}\}$$ is a Steiner liner. Moreover, a subset $F\subseteq X$ is flat in the liner $X$ if and only if $F$ is a submagma of the magma $(X,\circ)$.
\end{theorem}

\begin{proof} Let $X$ be a Steiner liner and $\circ:X\times X\to X$ be its midpoint operation, whose  definition ensures that $\{x,y,x\circ y\}=\Aline xy$ for any $x,y\in X$. In particular, $\Aline xx=\{x\}=\{x,x,x\circ x\}$ and hence $x\circ x=x$, which means that the midpoint operation is idempotent.

Since $\Aline xy=\Aline yx$, the midpoint operation is commutative. To see that it is involutory, observe that $\{x,y,x\circ y\}=\{x,x\circ y,y\}=\{x,x\circ y,x\circ(x\circ y)\}$ and hence $y=x\circ(x\circ y)$ if $|\{x,y,x\circ y\}|=3$. If $|\{x,y,x\circ y\}|<3$, then $\{x,y,x\circ y\}=\{x,x\circ y,x\circ(x\circ y)\}$ is a singleton and again $x\circ(x\circ y)=y$.
\smallskip

Now assume that $(M,\circ)$ is an involutory idempotent commutative magma $(X,\circ)$ and consider the relation $\Af\defeq\{(x,y,z)\in M^3:y\in\{x,z,x\circ z\}\}$. We claim that $\Af$ is a line relation on $M$. We need to check the axioms {\sf(IL)}, {\sf(RL)}, {\sf(EL)} from Definition~\ref{d:liner}. Indeed, if $(x,y,z)\in\Af$ and $x=z$, then $y\in\{x,z,x\circ z\}=\{x\}$, witnessing that the axiom {\sf(IL)} is satisfied. For every $x,z\in X$, we have $\{x,z\}\subseteq\{x,z,x\circ z\}$ and hence $(x,x,z),(x,z,z)\in \Af$, witnessing that the axiom {\sf(RL)} is satisfied. To see that $\Af$ satisfies the Exchange Axiom {\sf(EL)}, take any points  $a,b,x,y\in X$ such that $(a,x,b),(a,y,b)\in\Af$ and $x\ne y$. Then $\{x,y\}\subseteq\{a,b,a\circ b\}$.
If $\{x,y\}=\{a,b\}$, then $\{a,b\}\subseteq\{x,y,x\circ y\}$ and hence $(x,a,y),(x,b,y)\in\Af$, by the definition of $\Af$. If $\{x,y\}=\{a,a\circ b\}$, then the involutarity of the midpoint operation ensures that  $b=a\circ (a\circ b)=x\circ y$ and hence $\{a,b\}\subseteq\{x,y,x\circ y\}$ and $(x,a,y),(x,b,y)\in\Af$. If $\{x,y\}=\{b,a\circ b\}$, then $a=b\circ (b\circ a)=x\circ y$ and hence $\{b,a\}\subseteq\{x,y,x\circ y\}$ and $(x,a,y),(x,b,y)\in\Af$. Therefore, $\Af$ is a line relation. To see that $(M,\Af)$ is a Steiner liner, take any distinct points $x,y\in M$ and consider the line $\Aline xy=\{x,y,x\circ y\}$. Assuming that $x=x\circ y$, we conclude that $x=x\circ x=x\circ(x\circ y)=y$ by the idempotence and involutarity of the operation $\circ$. But the equality $x=y$ contradicts the choice of $x,y$. This contradiction shows that $x\ne x\circ y$. By analogy we can show that $y\ne x\circ y$ and hence $|\Aline xy|=|\{x,y,x\circ y\}|=3$, witnessing that $(M,\Af)$ is a Steiner liner. 
\smallskip

Now take any subset $F\subseteq X$. If $F$ is flat in the liner $X$, then for every $x,y\in F$, $x\circ y\in\Aline xy\subseteq F$, which means that $F$ is a submagma of the magma $(X,\circ)$. If $X$ is a submagma of the magma $(X,\circ)$, then for every $x,y\in F$, the line $\Aline xy=\{x,y,x\circ y\}$ is a subset of $F$, witnessing that $F$ is flat in the liner $X$.  
\end{proof}

A less evident is the relation of Steiner liners to commutative loops. 
By definition, every loop $X$ contains an element $e\in X$ such that $ex=x=xe$ for every $x\in X$. This element is unique and is called the \index{identity}\index{identity element}\defterm{identity element} of the loop. A loop $X$ is defined to be of \index{loop!of exponent 3}\defterm{exponent $3$} if 
$$(xx)x=x(xx)=e$$
 for every element $x\in X$.

\begin{proposition}\label{p:Steiner=>loop3} For every element $e$ of a Steiner liner $X$, the binary operation
$$\cdot:X\times X\to X,\quad \cdot:(x,y)\mapsto xy\defeq e\circ(x\circ y),$$
turns $X$ into a commutative loop of exponent $3$ with identity element $e$.
\end{proposition}

\begin{proof} Let $X$ be a Steiner liner and $e$ be any point of $X$. The commutativity of the midpoint operation on $X$ implies the commutativity of the binary operation $$\cdot:X\times X\to X,\quad \cdot:(x,y)\mapsto xy\defeq e\circ(x\circ y).$$ To see that this binary operation turns $X$ into a quasigroup, observe that for every $a,b\in X$ the point $x=a\circ(e\circ b)$ is the unique solution of the equation $ax=b$ because
$$ax=e\circ(a\circ  x)=e\circ (a\circ (a\circ (e\circ b)))=e\circ (e\circ b)=b.$$
If $y\in X$ is any other point with $ay=b$,
then 
$$x=a\circ(e\circ b)=a\circ (e\circ ay)=a\circ(e\circ (e\circ (a\circ y))=a\circ (a\circ y)=y,$$
by the involutarily of the midpoint operation. Therefore, the equation $ax=b$ has a unique solution $x=a\circ(e\circ b)$ in the commutative magma $(X,\cdot)$, which means that $(X,\cdot)$ is a quasigroup.

The point $e$ is the identity element of the quasigroup $(X,\cdot)$ because for every $x\in X$ the involutarity of the midpoint operation ensures that $ex=e\circ(e\circ x)=x$. 

The magma $(X,\cdot)$ has exponent $3$ because the properties of the midpoint operation ensure that
$$(xx)x=e\circ (xx\circ x)=e\circ ((e\circ (x\circ x))\circ x)=e\circ ((e\circ x)\circ x)=e\circ e=e$$and
$$(xx)(xx)=e\circ(xx\circ xx)=e\circ xx=e\circ(e\circ (x\circ x))=e\circ(e\circ x)=x$$for every $x\in X$.
\end{proof}

By Proposition~\ref{p:Steiner=>loop3}, for every element $e$ of a Steiner liner $X$, the binary operation 
$$\cdot:X\times X\to X,\quad \cdot:(x,y)\mapsto e\circ(x\circ y),$$
turns $X$ into a commutative loop with identity element $e$. We shall denote this loop by $X_e$. 

\begin{proposition}\label{p:flat<=>subloop} For any set $F$ in a Steiner liner $X$ and any element $e\in F$, the following conditions are equivalent:
\begin{enumerate}
\item $F$ is a flat in the liner $X$;
\item $F$ is a subloop of the loop $X_e$.
\item $F$ is a submagma of the magma $X_e$;
\item $F$ is a submagma of the magma $(X,\circ)$.
\end{enumerate}
\end{proposition}

\begin{proof} $(1)\Ra(2)$ If $F$ is a flat in the liner $X$, then by Theorem~\ref{t:Steiner<=>midpoint}, $F$ is submagma of the magma $(X,\circ)$ and hence for every $x,y\in F$ the point $xy=e\circ (x\circ y)$ belongs to $F$, witnessing that $F$ is a submagma of the magma $X_e$. 
To see that $F$ is a subloop of the loop $X_e$, observe that for every $a,b\in F$, the unique solution $x=a\circ (e\circ b)$ of the equation $ax=b$ belong to $F$ (because $F$ is a submagma of the magma $(X,\circ)$).
\smallskip

The implication $(2)\Ra(3)$ is trivial.
\smallskip

$(3)\Ra(4)$ If $F$ is a submagma of the magma $X_e$, then for every $a,b\in F$, we have
$$a\circ b=e\circ(e\circ (a\circ b))=e\circ ab=e\circ(ab\circ ab)=(ab)^2\in F,$$
which means that $F$ is a submagma of the magma $(X,\circ)$.

The implication $(4)\Ra(1)$ follows from Theorem~\ref{t:Steiner<=>midpoint}.
\end{proof}

\begin{question} For which magmas $X$, the relation $$\Af=\{(x,y,z)\in X:y\in \{x,z,(xz)^2\}\}$$ is a line relation on $X$?
\end{question}

Now we shall characterize projective Steiner liners via properties of the midpoint operation.

\begin{theorem}\label{t:Steiner<=>projective} For a Steiner liner $X$ the following conditions are equivalent:
\begin{enumerate}
\item $X$ is projective;
\item $X$ is $0$-parallel;
\item $X$ is $3$-balanced with $|X|_3=7$;
\item for any distinct points $x,y,z\in X$ we have $(x\circ y)\circ (x\circ z)=y\circ z$.
\end{enumerate}
\end{theorem}

\begin{proof} The equivalence $(1)\Leftrightarrow(2)\Leftrightarrow(3)$ follows from Theorems~\ref{t:projective<=>} and Corollary~\ref{c:Steiner-projective<=>}.
\smallskip

$(2)\Ra(4)$ Assume that the Steiner liner $X$ is $0$-parallel. Given three distinct points $x,y,z\in X$, we should prove the equality $(x\circ y)\circ(x\circ z)=y\circ z$. If $\{x,y,z\}$ is a line in $X$, then $$(x\circ y)\circ(x\circ z)=z\circ y=y\circ z.$$ If $\{x,y,z\}$ is not a line, then  $\{x,y,x\circ y\}$ and $\{x,z,x\circ z\}$ are two lines in the plane $P=\overline{\{x,y,z\}}$ such that $\{x,y,x\circ y\}\cap\{x,z,x\circ z\}=\{x\}$. Since the liner $X$ is $0$-parallel, the intersection of the lines  $\{y,z,y\circ z\}$ and $\{x\circ y,x\circ z,(x\circ y)\circ(x\circ z)\}$ contains a unique point $p$. Since $\{y,z\}\cap\{x\circ z,x\circ z\}=\varnothing$, the point $p$ is equal to the points $y\circ z$ and $(x\circ y)\circ(x\circ z)$ and hence $(x\circ y)\circ(x\circ z)=p=y\circ z$.
\smallskip

$(4)\Ra(1)$. Assume that the condition (4) is satisfied. To prove that the liner $X$ is projective, take any points $o,x,y\in X$ and $p\in\Aline xy\setminus\Aline ox$.  Given any point $u\in \Aline oy\setminus\{p\}$, we need to prove that $\Aline up\cap\Aline ox\ne\varnothing$. If $p\in\Aline oy$, then $\Aline up=\Aline oy$ and $o\in\Aline up\cap\Aline ox=\varnothing$. So, we assume that $p\notin\Aline oy$. In this case $p\in\Aline yx\setminus(\Aline ox\cup\Aline oy)$ implies $x\ne o$ and $y\notin\Aline ox$. If $u=y$, then $\Aline up=\Aline yp=\Aline yx$ and hence $x\in \Aline up\cap\Aline ox\ne\varnothing$. If $u=o$, then $o\in \Aline up\cap\Aline ox\ne\varnothing$. So, we assume that $o\ne u\ne y$. In this case $\Aline yu=\{y,u,o\}=\{y,u,y\circ u\}$ and $y\circ u=o$.  It follows from $p\in\Aline yx\setminus(\Aline ox\cup\Aline oy)$ that $\Aline xy=\{x,p,y\}$ and $y\circ p=x$. The condition (4) implies that $o\circ x=(y\circ u)\circ (y\circ p)=u\circ p$ and hence $u\circ p=o\circ x\in\Aline up\cap\Aline ox$, witnessing that the liner $X$ is projective.
\end{proof}








In Theorem~\ref{t:Hall<=>Moufang} we shall prove that a Steiner liner $X$ is Playfair if and only if for every $e\in X$, the commutative loop $X_e$ is Moufang. This characterization motivates studying commutative Moufang loops in more details.

\section{Commutative Moufang loops}$\color{white}.$\label{s:comMoufang}
\medskip

\rightline{\em Unfortunately, for the wide audience}

\rightline{\em whose interest in geometry}

\rightline{\em was awakened in high school or university,}

\rightline{\em answers to the deeper questions in geometry,}

\rightline{\em if they have been given at all,}

\rightline{\em require a long excursion into abstract algebra.}

\rightline{\em I know no remedy for this situation.}
\smallskip

\rightline{\index[person]{Bruck}R.H. Bruck\footnote{The citation is taken from Bruck's paper ``Recent Advances in the Foundations of Euclidean Geometry'' avarded by Chauvenet Prize, an annual award given by the Mathematical Association of America in recognition of an outstanding expository article on a mathematical topic. It consists of a prize of \$1,000 and a certificate.}, 1955}
\bigskip

In this section we establish some properties of commutative Moufang loops, which will be used in Section~\ref{s:Hall} devoted to analyzing the structure of Hall liners.

Let us recall that a loop $X$ is \index[person]{Moufang}{\em Moufang}\footnote{{\bf Ruth Moufang} (1905--1977), a German mathematician, Professor at University of Frankfurt, a student of Max Dehn who was a student of David Hilbert. In 1933 Moufang showed that Desargues' Theorem does not hold in the projective plane over octonions, thus initiating a new branch of geometry studying  Moufang planes. Many mathematical concepts are named after Ruth Moufang: Moufang loops, Moufang plane, Moufang-Lie algebras, Moufang polygons, Moufang sets.} if its satisfies the \index{Moufang identity}\defterm{Moufang identity}
$$(xy)(zx)=(x(yz))x$$ for every elements $x,y,z\in X$. 

For every element $a$ of a magma $X$, by $L_a$ and $R_a$ we shall denote the left and right shifts
$$L_a:X\to X,\;L_a:x\mapsto ax,\quad\mbox{and}\quad R_a:X\to X,\;R_a:x\mapsto xa,$$respectively. 

If the magma $X$ is a quasigroup, then the shifts $L_a,R_a$ are permutations (i.e., bijective maps) of the set $X$.

\begin{proposition}\label{p:loop-inverse} For every elements $x$ of a Moufang loop $X$, there exists a unique element $x^{-1}\in X$ such that $xx^{-1}=e=x^{-1}x$ and $(x^{-1})^{-1}=x$. Moreover, for every $x,y\in X$ we have
$$x^{-1}(xy)=y=(yx)x^{-1}.$$
\end{proposition}

\begin{proof} For every element $x\in X$ the bijectivity of the right shift $R_x:X\to X$ yields a unique element $x^{-1}\in X$ such that $x^{-1}x=R_x(x^{-1})=e$.
 
For every $y\in X$, the Moufang identity implies
$yx^{-1}=e(yx^{-1})=(x^{-1}x)(yx^{-1})=(x^{-1}(xy))x^{-1}$ and hence 
\begin{equation}\label{eq:loop1}
y=x^{-1}(xy),
\end{equation}
by the injectivity of the right shift  $R_{x^{-1}}:X\to X$.

Applying the identity $y=x^{-1}(xy)$ to the element $y=x^{-1}$, we obtain $x^{-1}e=x^{-1}=x^{-1}(xx^{-1})$ and hence $xx^{-1}=e$, by the injectivity of the left shift $L_{x^{-1}}$.

The equality $xx^{-1}=e=(x^{-1})^{-1}x^{-1}$ and the injectivity of the right shift $R_{x^{-1}}$ 
 imply 
\begin{equation}\label{eq:loop2}
x=(x^{-1})^{-1}.
\end{equation}

For every $x,y\in X$, the equalities (\ref{eq:loop1}), (\ref{eq:loop2}) and the Moufang identity imply
$$ey^{-1}=y^{-1}=x(x^{-1}y^{-1})=(y^{-1}(yx))(x^{-1}y^{-1})=(y^{-1}((yx)x^{-1}))y^{-1}$$and hence 
$e=y^{-1}((yx)x^{-1})$, by the injectivity of the right shift $R_{y^{-1}}$. Since $y^{-1}y=e=y^{-1}((yx)x^{-1})$, the injectivity of the left shift $L_{y^{-1}}$ implies $y=(yx)x^{-1}$.
\end{proof}

\begin{lemma}\label{l:(xy)y=x(yy)} Every elements $x,y$ of a commutative Moufang loop $X$ satisfy the identities $$(xx)y=x(xy)\quad\mbox{and}\quad y(xx)=(yx)x.$$
\end{lemma}

\begin{proof} Proposition~\ref{p:loop-inverse}, the commutativity of $X$ and the Moufang identity ensure that
$$(xx)y=(xx)((yx^{-1})x)=(x(x(yx^{-1})))x=(x(x(x^{-1}y)))x=(xy)x=x(xy)$$
and hence  
$$(yx)x=x(xy)=(xx)y=y(xx),$$
by the commutativity of $X$.
\end{proof}

\begin{lemma}\label{l:Moufang-xyy} Let $x,y,z$ be elements of a commutative Moufang loop $X$. If $|\{x,y,z\}|\le 2$, then $(xy)z=x(yz)$.
\end{lemma}

\begin{proof} Since $|\{x,y,z\}|\le 2$, three cases are possible.

1. If $x=z$, then $(xy)z=z(yx)=x(yz)$, by the commutativty of $X$.

2. If $y=z$, then $(xy)z=(xz)z=x(zz)=x(yz)$, by Lemma~\ref{l:(xy)y=x(yy)}.

3. If $x=y$, then $(xy)z=(xx)z=x(xz)=x(yz)$, by Lemma~\ref{l:(xy)y=x(yy)}.
\end{proof}

\begin{proposition}\label{p:com-Moufang<=>} A commutative loop $X$ is Moufang if and only if $$(xx)(yz)=(xy)(xz)$$ for every $x,y,z\in X$.
\end{proposition}

\begin{proof} If $X$ is Moufang, then we can apply the Moufang identity and Lemma~\ref{l:(xy)y=x(yy)} to obtain the equality
$$(xy)(xz)=(xy)(zx)=(x(yz))x=((yz)x)x=(yz)(xx)=(xx)(yz).$$

Now assume that a commutative loop $X$ satisfies the identity $(xx)(yz)=(xy)(xz)$ for every $x,y,z\in X$. Then for every $x,y,z,a\in X$ we have $$(xa)x=(xa)(xe)=(xx)(ae)=(xx)a.$$ Applying the equality $(xx)a=(xa)x$ to the element $a=yz$, we obtain the equality $$(xy)(zx)=(xy)(xz)=(xx)(yz)=(x(yz))x,$$witnessing that the loop $X$ is Moufang.
\end{proof}

 A self-map $h:X\to X$ of a magma $X$ is called an \index{endomorphism}\index{magma!endomorphism of}\defterm{endomorphism} of the magma $X$ if $$h(xy)=h(x)h(y)$$ for every $x,y\in X$. If $h$ is bijective, then $h$ is an \index{automorphism}\index{magma!automorphism of}\defterm{automorphism} of the magma $X$.


\begin{lemma}\label{l:Fix-subloop} For every endomorphism $h:X\to X$ of a loop $X$, the set $$S\defeq\{x\in X:h(x)=x\}$$is a subloop of the loop $X$.
\end{lemma}

\begin{proof} Given any elements $x,y\in S$, observe that $h(xy)=h(x)h(y)=xy$ and hence $xy\in S$. Since $X$ is a quasigroup, for every $a,b \in S$ there exist unique elements $x,y\in X$ such that $ax=b$ and $ya=b$. Then $ax=b=h(b)=h(ax)=h(a)h(x)=ah(x)$ and hence $x=h(x)$ by the uniqueness of $x$. By analogy we can show that $h(y)=y$ and hence $x,y\in S$, witnessing that $S$ is a subloop of the loop $X$.
\end{proof} 

\begin{lemma}\label{l:associator-subloop} For every elements $a,b$ of a commutative Moufang loop $X$, the set $$\{z\in X:(ab)z=a(bz)\}$$ is a subloop of the loop $X$.
\end{lemma}

\begin{proof} For every $a\in X$, the left shift $L_a:X\to X$, $L_a:x\mapsto ax$, is a permutation of $X$ because the loop $X$ is a quasigroup. Then $L_a$ is an element of the permutation group $S_X$ of the set $X$. By Proposition~\ref{p:loop-inverse}, for every $a,x\in X$ we have
$$L_{a^{-1}}(L_a(x))=a^{-1}(ax)=x,$$which means that $$L_a^{-1}=L_{a^{-1}}.$$ Lemma~\ref{l:Moufang-xyy} implies $a(ax)=(aa)x=a^2x$ and hence $L_aL_a=L_{a^2}$. The identity $(ax)(ay)=a^2(xy)$ proved in Proposition~\ref{p:com-Moufang<=>} can be rewritten as the equality 
\begin{equation}\label{eq:Bruck2}
L_{a^2}(xy)=L_a(x)L_a(y).
\end{equation} 

Let $G$ be the smallest subgroup of the permutation group $S_X$, containing all left shifts $L_a$, $a\in X$. Since $L_a^{-1}=L_{a^{-1}}$, every element of the group $G$ can be written as the composition of finitely many left shifts. Given any elements $a_1,a_2,\dots,a_n\in G$, consider the permutations
$$h\defeq L_{a_1}\cdots L_{a_n}\quad\mbox{and}\quad H\defeq L_{a_1^2}\cdots L_{a_n^2}.$$
By repeated application of (\ref{eq:Bruck2}), we obtain that $H (xy)=h(x)h(y)$ for every $x,y\in X$. If $h(e)=e$, then $H(x)=H(xe)=h(x)h(e)=h(x)e=h(x)$ and hence $H=h$ and $h(xy)=H(xy)=h(x)h(y)$, which means that $h$ is an automorphism of the loop $X$. Therefore, the subgroup $$G_e\defeq\{h\in G:h(e)=e\}$$ of the group $G$ consists of automorphisms of the loop $X$.

Now fix any elements $a,b$ of the commutative Moufang loop $X$. The equation $$L_{ab}(e)=(ab)(e)=ab=a(be)=L_aL_b(e)$$ implies $L_{b^{-1}}L_{a^{-1}}L_{ab}(e)=e$ and $L_{b^{-1}}L_{a^{-1}}L_{ab}\in G_e$.  By Lemma~\ref{l:Fix-subloop}, the set 
$$\{z\in X:(ab)z=a(bz)\}=\{z\in X:L_{ab}(z)=L_aL_b(z)\}=\{z\in X:L_b^{-1}L^{-1}_aL_{ab}(z)=z\}$$is a subloop of $X$.
\end{proof}

\begin{lemma}\label{l:loop-bac} Let $a,b,c$ be elements of a commutative Moufang loop $X$. If $(ab)c=a(bc)$, then $(ba)c=b(ac)$.
\end{lemma}

\begin{proof}
Let $S$ be the smallest subloop of $X$  containing the elements $a,b,c$. By Lemma~\ref{l:Fix-subloop}, the set $S'\defeq\{z\in S:(ab)z=a(bz)\}$ is a subloop of $S$, containing the element $c$. Then $c^{-1}\in S'$ and hence  $(ab)c^{-1}=a(bc^{-1})$. By Propositions~\ref{p:loop-inverse} and \ref{p:com-Moufang<=>}, $$(ba)c=((ab)c)(c^{-1}c)=((ab)c^{-1})(cc)=(a(bc^{-1}))(cc)=(ac)((bc^{-1})c)=(ac)b=b(ac).$$
\end{proof}

\begin{lemma}\label{l:Moufang-abc} Let $a,b,c$ be elements of a commutative Moufang loop. If $(ab)c=a(bc)$, then $(xy)z=x(yz)$ for any points $x,y,z\in \{a,b,c\}$.
\end{lemma}

\begin{proof}  By Lemma~\ref{l:loop-bac}, the equality $(ab)c=a(bc)$ implies the equality $(ba)c=b(ac)$. By the commutativity of $X$, the equalities $(ab)c=a(bc)$ and $(ba)c=b(ac)$ imply the equalities $(cb)a=c(ba)$ and $(ca)b=c(ab)$. By Lemma~\ref{l:loop-bac}, the equalities $(cb)a=c(ba)$ and $(ca)b=c(ab)$ imply $(bc)a=b(ca)$ and $(ac)b=a(cb)$. Therefore, $(xy)z=x(yz)$ for any distinct points $x,y,z\in\{a,b,c\}$. If $|\{x,y,z\}|\le 2$, then the equality $(xy)z=x(yz)$ follows from Lemma~\ref{l:Moufang-xyy}. 
\end{proof}

\begin{theorem}[Moufang]\label{t:Moufang-loop} Let $a,b,c\in X$ be elements of a commutative Moufang loop $X$. If $(ab)c=a(bc)$, then the set $\{a,b,c\}$ is contained in a subgroup of the loop $X$.
\end{theorem}

\begin{proof} 
Let $a,b,c\in X$ be any elements of a commutative Moufang loop $X$ such that $(ab)c=a(bc)$. Let $S$ be the smallest subloop containing the elements $a,b,c$. By Lemma~\ref{l:associator-subloop}, the set $S'\defeq\{z\in S:(ab)z=a(bz)\}$ is a subloop of $S$, containing the element $c$. Lemma~\ref{l:Moufang-xyy} ensures that $a,b\in S'$ and hence $S'=S$ because $S$ is the smallest subloop of $X$ that contains the elements $a,b,c$. Therefore, $(ab)z=a(bz)$ for all elements $z\in S$.

By Lemma~\ref{l:Moufang-abc}, the equality $(ab)c=a(bc)$ implies $(cb)a=c(ba)$. Repeating the above argument, we can show that $x(bc)=(cb)x=c(bx)=(xb)c$ for every $x\in S$. By Lemma~\ref{l:associator-subloop}, for every $x\in S$ the set $S_x=\{z\in S:x(bz)=(xb)z\}$ is a subloop of $S$. This subloop contains the elements $a$ and $c$. By Lemma~\ref{l:Moufang-xyy}, the subloop $S_x$ contains the element $b$.  Then $S_x=S$ because $S$ is the smallest subloop containing the elements $a,b,c$.  Therefore, $x(bz)=(xb)z$ for all $x,z\in S$. 

By Lemma~\ref{l:Moufang-abc}, the equality $(ab)c=a(bc)$ implies the equalities $(ba)c=b(ac)$ and $(ac)b=a(cb)$. Repeating the above argument, we can prove that $(xa)z=x(az)$ and $(xc)z=x(cz)$ for all $x,z\in S$.  

For every $x,z\in S$, consider the set $S_{x,z}=\{y\in S:(xy)z=x(yz)\}$ and observe that $\{a,b,c\}\subseteq S_{x,z}$. Lemma~\ref{l:Moufang-abc} implies that $S_{x,z}\defeq\{y\in S:(xy)z=x(yz)\}=\{y\in S:(xz)y=x(zy)\}$. By Lemma~\ref{l:associator-subloop}, the set $\{y\in S:(xz)y=x(zy)\}=S_{x,z}$ is a subloop of $S$. Since $\{a,b,c\}\subseteq S_{x,z}$, the subloop $S_{x,z}$ coincides with the smallest subloop $S$ of $X$, containing the points $a,b,c$. Then for every $y\in S=S_{x,z}$ we have $(xy)z=x(yz)$, which means that the loop $S$ is associative and hence $S$ is a subgroup of the loop $X$.
\end{proof}


For two elements $x,y$ of a loop $X$, their \index{commutator}\defterm{commutator} is the element $$[xy]=((xy)x^{-1})y^{-1}$$of the loop $X$. Proposition~\ref{p:loop-inverse} implies that a loop $X$ is commutative if and only if $[x,y]=e$ for every elements $x,y\in X$.

The \index{associator}\defterm{associator} $[x,y,z]$ of elements $x,y,z\in X$ of a loop $X$ is defined as 
$$[x,y,z]\defeq ((xy)z)(x(yz))^{-1}.$$

 Proposition~\ref{p:loop-inverse} implies that a loop $X$ is associative if and only if $[x,y,z]=e$ for every elements $x,y,z\in X$.

\begin{lemma}\label{l:associator} Let $X$ be a commutative Moufang loop of exponent $3$. For any elements $x,y,z\in X$ the following identities hold:
\begin{enumerate}
\item $[x,y,z]=[z^{-1},y^{-1},x^{-1}]$;
\item $[x,y,z]=[x,z^{-1},y]$;
\item $[x,y,[x,y,z]]=e$.
\end{enumerate}
\end{lemma}

\begin{proof} In the subsequent calculations, for an element $x$ of the loop $X$, it will be convenient to denote its inverse $x^{-1}$ by $\bar x$. So, $x\bar x=\bar xx=e$. Since the loop $X$ has order $3$, $\bar x=x^2$ and hence $x(yz)=x^4(yz)=(x^2y)(x^2z)=(\bar xy)(\bar xz)$ for every $y,z\in X$. We shall use the identities 
$$x(yz)=(\bar xy)(\bar xz)\quad\mbox{and}\quad (xy)z=(x\bar z)(y\bar z)$$ in the subsequent calculations without special reference. By Lemma~\ref{l:Moufang-xyy} and the Moufang Theorem~\ref{t:Moufang-loop}, any two elements $x,y$ of $X$ belong to some (commutative) subgroup of $X$. Therefore, $(xy)^{-1}=y^{-1}x^{-1}=\bar x\bar y$ for every elements $x,y\in X$. Now, given any elements $x,y,z\in X$, we prove the statements (1)--(3) of the lemma.
\smallskip

1. By the commutativity of $X$,
$$[x,y,z]=((xy)z)(x(yz))^{-1}=((xy)z)(\bar x(\bar y\bar z))=((\bar z\bar y)\bar x)(z(yx))=[\bar z,\bar y,\bar x]=[z^{-1},y^{-1},x^{-1}].$$

2. The second statement of the lemma follows from the following series of equalities:
$$
\begin{aligned}
[x,y,z]&=\big((xy)z\big)\big(\bar x(\bar y\bar z)\big)=\big((xy)(x(yz))\big)\big(z(x(yz))\big)=\big(((\bar x\bar y)x)((\bar x\bar y)(yz))\big)\big((\bar zx)(\bar z(zy))\big)\\
&=\big(\bar y((\bar x\bar y)(yz))\big)\big((\bar zx)y\big)=\big((y(\bar x\bar y))(y(yz))\big)\big((\bar zx)y)\big)=\big(\bar x(\bar yz)\big)\big((\bar zx)y\big)\\
&=\big((x\bar z)y\big)\big(\bar x(z\bar y)\big)=[x,\bar z,y]=[x,z^{-1},y].
\end{aligned}
$$

3. Next, we show that $[x,y,[x,y,z]]=e$, which is equivalent to the equality $(xy)[x,y,z]=x(y[x,y,z])$.
Observe that
$$
\begin{aligned}
(xy)[x,y,z]&=(xy)\big(((xy)z)(\bar x(\bar y\bar z))\big)=\big((\bar x\bar y)((xy)z)\big)\big((\bar x\bar y)(\bar x(\bar y\bar z))\big)=z\big(\big((xy)\bar x\big)\big((xy)(\bar y\bar z)\big)\big)\\
&=z\big(y\big((xy)(\bar y\bar z)\big)\big)=z\big(\big(\bar y(xy)\big)\big(\bar y(\bar y\bar z)\big)\big)=z(x(y\bar z))=(\bar zx)(\bar zy\bar z)=(x\bar z)(yz).
\end{aligned}
$$
On the other hand, 
$$
\begin{aligned}
x(y[x,y,z])&=x\big(y\big(((xy)z)(\bar x(\bar y\bar z))\big)\big)=x\big(\big(\bar y((xy)z)\big)\big(\bar y(\bar x(\bar y\bar z))\big)\big)=x\big(\big((y(xy))(yz)\big)\big((y\bar x)(y(\bar y\bar z))\big)\big)\\
&=x\big(\big((x\bar y)(yz)\big)\big((y\bar x)\bar z\big)\big)=x\big(\big((y\bar x)\bar z\big)\big((yz)(x\bar y)\big)\big)=x\big(\big((\bar y x)z)(yz)\big)\big((\bar yx)z)(x\bar y)\big)\big)\\
&=x\big(\big(((\bar yx)(\bar y\bar z))(z(\bar y\bar z))\big)\big((\bar xy)z\big)\big)=x\big(\big(((\bar yx)(\bar y\bar z))\bar y\big)\big((\bar xy)z\big)\big)\\
&=x\big(\big((y(x\bar z))\bar y\big)\big((\bar xy)z\big)\big)=x\big(\big(x\bar z\big)\big((\bar xy)z\big)\big)=\big(\bar x(x\bar z)\big)\big(\bar x((\bar xy)z)\big)\\
&=\bar z\big((x(\bar xy))(xz)\big)=\bar z(y(xz))=(zy)(z(xz))=(yz)(x\bar z),
\end{aligned}
$$
and hence $(xy)[x,y,z]=(x\bar z)(yz)=x(y[x,y,z])$ and $[x,y,[x,y,z]]=e$.
\end{proof}

The \index{centre}\index{Moufang loop!centre of}\defterm{\em centre} of a Moufang loop $X$ is the set
$$\big\{z\in X:\forall x,y\in X\;\big([x,y,z]=e=[x,z]\big)\big\}.$$
By Lemma~\ref{l:associator-subloop}, the centre of a commutative Moufang loop $X$ is a commutative subgroup of $X$. A subloop $H$ of a Moufang loop $X$ is called \index{central subloop}\index{subloop!central}\defterm{central} if it is contained in the centre of the loop $X$. It is easy to see that every central subloop of a Moufang loop is associative and hence is a subgroup of the loop $X$. Every  central subgroup $H$ of a Moufang loop $X$ determines the quotient Moufang loop $X/H$ whose elements are cosets $xH=Hx$. The binary operation on cosets can be defined using the pointwise product of two sets.  

For any sets $A,B\subseteq X$ in a magma $X$ and any point $x\in X$, consider the sets
$$xA\defeq\{xa:a\in A\},\quad Ax\defeq\{ax:a\in A\}\quad\mbox{and}\quad AB\defeq\{ab:a\in A,\;b\in B\}.$$ 

\begin{lemma}\label{l:loop-cosets} Let $H$ be a central subloop of a Moufang loop $X$. For every $x,y\in X$,
\begin{enumerate}
\item $xH=Hx$;
\item $xH=zH$ for every element $z\in xH$;
\item $xH=yH$ if and only if $xH\cap yH\ne\varnothing$;
\item $(xH)(yH)=(xy)H$.
\end{enumerate}
\end{lemma}

\begin{proof} Fix any elements $x,y\in X$.
\smallskip

1. To see that $xH=Hx$, it suffices to show that $xh=hx$ for every $h\in X$. Since $H$ is a subloop of the centre of $X$, we have $e=[xh]=((xh)x^{-1})h^{-1}$. Multiplying this equality by $h$ from the right and applying Proposition~\ref{p:loop-inverse}, we obtain $h=eh=(((xh)x^{-1})h^{-1})h=(xh)x^{-1}$. Multiplying the equality $h=(xh)x^{-1}$ by $x$ from the right and applying Proposition~\ref{p:loop-inverse}, we obtain $hx=((xh)x^{-1})x=xh$, witnessing that the elemements $x,h$ commute and hence $xH=Hx$.
\vskip3pt

2. For every element $z\in xH$, we have $zH\subseteq (xH)H\subseteq x(HH)=xH$ because $[x,h,z]=e$ for every $h\in H$ and $x,z\in X$. On the other hand, $z=xh$ for some $h\in H$ and hence $x=(xh)h^{-1}=zh^{-1}\subseteq zH$. Then $xH\subseteq zH$ and hence $xH=zH$.
\smallskip

3. If $xH\cap yH$ contains some point $z\in X$, then $xH=zH=yH$ by the preceding statement.
\smallskip

4. Since $H$ is a subset of the centre of $X$, we have $(xy)H=x(yH)=(xe)(yH)\subseteq (xH)(yH)$.
On the other hand, for every $g,h\in H$, we have $(xg)(yh)=x(g(yh))=x((yh)g)=x(y(hg))=(xy)(hg)\in (xy)H$, because $g\in H$ belongs to the centre of the loop $X$. Therefore, $(xH)(yH)=(xy)H$.
\end{proof}

Given a subgroup $H$ of  the centre of a Moufang loop $X$, consider the set
$$X/H\defeq\{xH:x\in X\}=\{Hx:x\in X\},$$endowed with the operation of pointwise multiplication $(xH)(yH)=(xy)H$. By Lemma~\ref{l:loop-cosets}, this operation is well-defined. 
So, $X/H$ is a magma. In Proposition~\ref{p:loop-quotient} we shall prove that the magma $X/H$ is a Moufang loop. This Moufang loop is called the \index{quotient loop}\index{loop!quotient}\defterm{quotient loop} of $X$ by the central subgroup $H$.

\begin{proposition}\label{p:loop-quotient} For every central subgroup $H$ of a (commutative) Moufang loop $X$, the magma $X/H$ is a (commutative) Moufang loop.
\end{proposition}

\begin{proof} Since $(eH)(xH)=(ex)H=xH=(xe)H=(xH)(eH)$ for every $x\in X$, the magma $X/H$ is unital. To see that $X$ is a quasigroup, take any cosets $aH,bH\in X/H$. Since $X$ is a quasigroup, there exist unique elements $x,y\in X$ such that $ax=b=ya$. Then $(aH)(xH)=bH=(yH)(aH)$. It remains to show that $xH$ and $yH$ are unique elements of the magma $X/H$ with $(aH)(xH)=bH=(yH)(aH)$. Assume that $x'H,y'H\in X/H$ are elements of $X/H$ such that $(aH)(x'H)=bH=(y'H)(aH)$. Then $(ax')H=(aH)(x'H)=bH$ and hence $ax'\in bH=(ax)H$ and $ax'=(ax)h$ for some $h\in H$. Since $h$ belongs to the centre of the loop $X$, $ax'=(ax)h=a(xh)$. Multiplying the equality $ax'=(ax)h$ by $a$ from the left, and applying Proposition~\ref{p:loop-inverse}, we obtain $x'=a^{-1}(ax')=a^{-1}(a(xh))=xh$ and hence $x'H=(xh)H=x(hH)=xH$. By analogy we can show that $y'H=yH$, witnssing that the unital magma $X/H$ is quasigroup and hence a loop.

For every $x,y,z\in X$, the Moufang identity implies $$((xH)(yH))((zH)(xH))=((xy)(zx))H=
((x(yz))x)H=\big(xH((yH)(zH))\big)xH,$$which means that the loop $X/H$ is Moufang.

If the loop $X$ is commutative, then for every $x,y\in X$, we have
$$(xH)(yH)=(xy)H=(yx)H=(yH)(xH),$$
which means that the quotient loop $X/H$ is commutative.
\end{proof}

\begin{theorem}\label{t:associator} Let $X$ be a commutative Moufang loop of exponent $3$. For every elements $x,y,z\in X$, their associator $[x,y,z]$ belongs to the centre of the  subloop $S\subseteq X$, generated by the elements  $x,y,z$.
\end{theorem} 

\begin{proof} Given any elements $x,y,z\in X$, denote their associator $[x,y,z]$ by $a$.
By Lemma~\ref{l:associator}(3), $(xy)a=x(ya)$. By Lemma~\ref{l:Moufang-abc}, the equality $(xy)a=x(ya)$ implies $(ax)y=a(xy)$. By Lemma~\ref{l:associator}(2), $a=[x,y,z]=[x,z^{-1},y]$ and by Lemma~\ref{l:associator}(3), $(xz^{-1})a=x(z^{-1}a)$. By Theorem~\ref{t:Moufang-loop}, the elements $x,z^{-1},a$ belong to an associative subloop of $X$ and hence $[a,x,z]=e=[a,x,x]$. Then the set $S_x\defeq\{s\in S:[a,x,s]=e\}$ contains the elements $x,y,z$.  By Lemma~\ref{l:associator-subloop}, the set $S_x=\{s\in S:[a,x,s]=e\}$ is a subloop of $X$. Since the subloop $S_x$ contains the elements $x,y,z$, it coincides with the subloop $S$, generated by $x,y,z$. Therefore, $[a,x,s]=e$ for every element $s\in S$. By analogy we can prove that $[a,y,s]=[a,z,s]=e$ for every $s\in S$. Applying Lemma~\ref{l:Moufang-abc}, for every $s\in S$, the equalities $[a,x,s]=[a,y,s]=[a,z,s]=e$ imply $[a,s,x]=[a,s,y]=[a,s,z]=e$. Therefore, the elements $x,y,z$ belong to the set $S_{a,s}\defeq\{t\in S:[a,s,t]=e\}$. By  Lemma~\ref{l:associator-subloop}, the set $S_{a,s}$ is a subloop of $S$ and hence $S_{a,s}=S$, which means that $[a,s,t]=e$ and hence the element $a=[x,y,z]$ belongs to the centre of the loop $S$.
\end{proof}

\begin{theorem}\label{t:Moufang81} Let $X$ be a commutative Moufang loop of exponent $3$, let $x,y,z$ be any elements in $X$, let $S$ be the subloop of $X$, generated by the elements $x,y,z$, and $A\subseteq S$ be the subgroup of $S$, generated by the associator $a\defeq [x,y,z]$. If $a\ne e$, then 
\begin{enumerate}
\item the quotient loop $S/A$ is a commutative group of exponent $3$ and cardinality $|S/A|=27$;
\item $|S|=81$;
\item for any distinct points $u,v\in\{x,y,z\}$, the subloop of $X$, generated by the elements $u,v,a$, is a group of cardinality $27$.
\end{enumerate}
\end{theorem}

\begin{proof} By Theorem~\ref{t:associator}, the associator $a$ belongs to the centre of the commutative loop $S$. Lemma~\ref{l:associator-subloop} implies that the center of $S$ is a subloop of $S$ and hence $A$ is a central subgroup in $S$. So, we can consider the quotient loop $S/A$, which is a commutative Moufang loop, by Proposition~\ref{p:loop-quotient}. It follows from $[x,y,z]=a\in A$ that $[xA,yA,zA]=eA$. By Moufang Theorem~\ref{t:Moufang-loop}, the quotient loop $S/A$ is associative and hence $S/A$ is a commutative group. Since the loop $X$ has exponent $3$, the group $S/A$ has exponent $3$, too.  Then $S/A$ is a commutative group of exponent $3$ with three generators $xA,yA,zA$, and hence $|S/A|=3^n$ for some $n\le 3$. Assuming that $|S/A|<27$, we can find two elements $u,v\in\{x,y,z\}$ such that $\{xA,yA,zA\}\subseteq\{(u^iA)(v^jA):i,j\in\{-1,0,1\}\}$. We lose no generality assuming that $u=x$ and $v=y$. Then $zA=(x^iA)(y^jA)$ for some $i,j\in\{-1,0,1\}$. By Lemma~\ref{l:associator}(3), $[x,y,a]=e$ and by Moufang Theorem~\ref{t:Moufang-loop}, the elements $x,y,a$ are contained in an associative subloop $G$ of the loop $X$. Then $z\in x^iy^jA\subseteq G$ and hence $a=[x,y,z]=e$, which contradicts our assumption. This contradiction shows that $|S/A|=27$ and then $S=|S/A|\cdot |A|=27\cdot 3=81$, by Lemma~\ref{l:loop-cosets}(3).

Now take any distinct elements $u,v\in\{x,y,z\}$ and consider the subloop $H$ generated by the elements $u,v,a$. Let $w$ be the unique element of the set $\{x,y,z\}\setminus \{u,v\}$. By Theorem~\ref{t:associator}, the element $a=[x,y,z]$ belongs to the centre of the loop $S$. Therefore $[u,v,a]=e$ and by Moufang Theorem~\ref{t:Moufang-loop}, the loop $H$ is associative. Therefore, $H$ is a commutative group of exponent $3$, generated by three elements $u,v,a$. Then $|H|=3^n$ for some $n\le 3$. 
Assuming that $|H|<27$ and taking into account that $a\ne e$, we can find an element $g\in\{u,v\}$ such that $\{u,v\}\subseteq \{g^ia^j:i,j\in\{-1,0,1\}\}$. We lose no generality assuming that $g=v$. Then $u\in v^ia^j$ for some $i,j\in\{-1,0,1\}$ and $uA\in v^ia^jA=v^iA$. Then the group $S/A$ is generated by two elements $vA$ and $wA$, which contradicts the equality $|S/A|=27$. This contradictsion shows that $|H|=27$.  
\end{proof}

A magma $X$ is defined to be \index{magma!$n$-generated}\index{$n$-generated magma}\defterm{$n$-generated} for a number $n\in\IN$ if $X$ is generated by $n$ elements but is not generated by $(n-1)$-element.

\begin{corollary}\label{c:Moufang81} Any two non-associative $3$-generated commutative Moufang loops of exponent $3$ are isomorphic and contain exactly $81$ element.
\end{corollary}

\begin{proof} Let $X,X'$ be two  non-associative $3$-generated commutative Moufang loops of exponent 3. Let $x,y,z$ be generators of the loop $X$ and $x',y',z'$ be generators of the loop $X'$. Let $e,e'$ be the identity elements of the loops $X,X'$, respectively. Since $X$ is not associative, Moufang Theorem~\ref{t:Moufang-loop} ensures that $[x,y,z]\ne e$. Then $|X|=81$, by Theorem~\ref{t:Moufang81}. By analogy we can prove that $|X'|=81$. 
In the product $X\times X'$ consider the subloop $Y$, generated by the elements $(x,x'),(y,y'),(z,z')$. It is clear that the projections $p:Y\to X$, $p:(u,u')\mapsto u$, and $p':Y\to X'$, $p':(u,u')\mapsto u'$, are surjective homomorphisms from the loop $Y$ onto the loops $X$ and $X'$, respectively. Then $[(x,x'),(y,y'),(z,z')]=([x,y,z],[x',y',z'])\ne(e,e')$. By Theorem~\ref{t:Moufang81}, the loop $Y$ has cardinality $|Y|=81$, which implies that the surjective homomorphism $p:Y\to X$ and $p':Y\to X'$ are bijective and hence are isomorphisms. Then the loops $X,X'$ are isomorphic to $Y$ and hence are isomorphic.
\end{proof} 

\begin{remark} The Moufang Theorem~\ref{t:Moufang-loop} was first proved by Ruth Moufang in 1935, by complicated induction. The proof of Theorem~\ref{t:Moufang-loop} presented here is due to Bruck in \cite{Bruck1951}. In fact, this theorem remains true for non-commutative Moufang loops as well. Also Theorem~\ref{t:associator} has its far generalizations to non-commutative Moufang loops, see the monography of Bruck \cite{Bruck1958}. We presented here only the small piece of the theory of Moufang loops, providing necessary tools for studying the structure of Hall liners and presenting a self-contained proof of Theorem~\ref{t:Playfair<=>regular}.
\end{remark}

\section{Hall liners}\label{s:Hall}

\begin{definition} A liner $X$ is called \index[person]{Hall}\index{Hall line}\index{liner!Hall}\defterm{Hall}  if $X$ is $2$-balanced and $3$-balanced with $$|X|_3=(|X|_2)^2=9.$$
\end{definition}

The liner from Example~\ref{ex:HTS} is an example of a Hall liner. Hall liners are well-known in Theory of Combinatorial Designs as \index{Hall Triple System}\defterm{Hall Triple Systems}, see \cite[\S VI.28]{HCD}. Corollary~\ref{c:Playfair<=>23balanced} implies the following characterization of Hall liners.

\begin{theorem}\label{t:Hall<=>Playfair+Steiner} A liner $X$ is Hall if and only if $X$ is Steiner and Playfair.
\end{theorem}

Hall liners also have an interesting algebraic characterization. Let us recall that every Steiner liner $X$ carries the midpoint operation $\circ:X\times X\to X$ assigning to every pair $(x,y)\in X\times X$ the unique point $x\circ y$ of $X$ such that $\Aline xy=\{x,y,x\circ y\}$. By Theorem~\ref{p:Steiner=>loop3}, for every point $e\in X$, the midpoint operation determines the binary operation $\cdot:X\times X\to X$, $\cdot:(x,y)\mapsto e\circ(x\circ y)$ turning $X$ into a commutative loop of exponent $3$ with identity element $e$. We denote this loop by $X_e$. 

Let us recall that 
a magma $X$ is {\em self-distributive} if $x(yz)=(xy)(xz)$ for all $x,y,z\in X$.

\begin{theorem}\label{t:Hall<=>Moufang} For a Steiner liner $X$, the following conditions are equivalent:
\begin{enumerate}
\item $X$ is Hall;
\item the mid-point operation on $X$ is self-distributive;
\item for every point $e\in X$, the commutative loop $X_e$ is Moufang.
\end{enumerate}
\end{theorem}

\begin{proof} We shall prove the implications $(1)\Ra(2)\Ra(3)\Ra(1)$.
\smallskip

$(1)\Ra(2)$ Assume that the Steiner liner $X$ is Hall. To prove that the mid-point operation $\circ$ on $X$ is self-distibutive, take any points $x,y,z\in X$. We have to prove that $x\circ(y\circ z)=(x\circ y)\circ(x\circ z)$. If $y=z$, then $$x\circ (y\circ z)=x\circ (y\circ y)=x\circ y=(x\circ y)\circ(x\circ y)=(x\circ y)\circ(x\circ z),$$by the idempotence of the midpoint operation.
So, assume that $y\ne z$. In this case $y\circ z$ is the unique point $p$ of the set $\Aline yz\setminus\{y,z\}$. Depending on the location of the point $x$, two cases are possible. 

First assume that $x\in \Aline yz=\{y,z,p\}$. This case has three subcases.

If $x=y$, then $x\circ (y\circ z)=y\circ p=z$ and $(x\circ y)\circ(x\circ z)=y\circ p=z=x\circ(y\circ z)$.

If $x=z$, then $x\circ (y\circ z)=z\circ p=y$ and $(x\circ y)\circ (x\circ z)=(z\circ y)\circ z=p\circ z=y=x\circ(y\circ z)$.

If $x=p$, then $x\circ (y\circ z)=p\circ p=p$ and $(x\circ y)\circ (x\circ z)=(p\circ y)\circ(p\circ z)=z\circ y=p=x\circ(y\circ z)$.

Next, assume that $x\notin\Aline yz=\{y,z,p\}$ and consider the plane $P=\overline{\{x,y,z\}}$ and the line $L_y\defeq\Aline xy$ in the plane $P$.  By Theorem~\ref{t:Hall<=>Playfair+Steiner}, the Hall liner is Playfair. Then there exist unique lines $L_p,L_z\subseteq P$ such that $p\in L_p\subseteq P\setminus L_y$ and $z\in L_z\subseteq P\setminus L_y$.
Assuming that $p\in L_z$, we conclude that $y=p\circ z\in L_z$, which contradicts $L_z\cap L=\varnothing$. This contradiction shows that $p\notin L_z$, which implies that the lines $L_p$ and $L_z$ are distinct. Since these two lines are disjoint with the line $L_y$ in the plane $P$, the Playfair Parallel Postulate ensures that $L_z\cap L_p=\varnothing$.  It follows from $|P|=|X|_3=9$ that $P=L_y\cup L_p\cup L_z$. 

Next, consider the line $\Lambda_z\defeq\Aline xz$. By the Playfair Postulate, there exist unique lines $\Lambda_y,\Lambda_p\subseteq P$ such that $y\in \Lambda_y\subseteq P\setminus \Lambda_z$ and $p\in \Lambda_p\subseteq P\setminus \Lambda_z$. Repeating the preceding argument, we can show that the lines $\Lambda_y,\Lambda_z,\Lambda_p$ are disjoint that $P=\Lambda_y\cup\Lambda_z\cup\Lambda_p$. 

\begin{picture}(150,160)(-150,-35)

\multiput(-20,0)(0,40){3}{\line(1,0){120}}
\multiput(0,-20)(40,0){3}{\line(0,1){120}}
\multiput(0,0)(40,0){3}{\multiput(0,0)(0,40){3}{\circle*{3}}}
\put(0,80){\line(1,-1){80}}

\put(-9,-9){$x$}
\put(83,-9){$y$}
\put(-9,83){$z$}
\put(32,32){$p$}
\put(105,-4){$L_y$}
\put(105,36){$L_p$}
\put(105,76){$L_z$}
\put(-3,105){$\Lambda_z$}
\put(37,105){$\Lambda_p$}
\put(77,105){$\Lambda_y$}
\put(80,80){\color{red}\circle*{4}}
\end{picture}

Assuming that the point $x\circ p=x\circ(y\circ z)$ belongs to the line $L_y=\{x,y,x\circ y\}$, we conclude that $p=x\circ (x\circ p)\subseteq L_y$ and hence $p\in L_y\cap L_p=\varnothing$, which is a contradiction showing that $x\circ p\notin L_y$. Assuming that $x\circ p\in L_p$, we conclude that $x=p\circ(p\circ x)\in L_p$ and hence $x\in L_y\cap L_p=\varnothing$, which is a contradiction showing that $x\circ p\notin L_p$. Therefore, $x\circ p\in P\setminus(L_y\cup L_p)=L_z$. By analogy we can prove that $x\circ p\in P\setminus(\Lambda_z\cup\Lambda_p)=\Lambda_y$.

Taking into account that $x\circ y\in L_y\setminus L_z$ and $x\circ z\in \Lambda_z\setminus\Lambda_y$, we can show that $(x\circ y)\circ(x\circ z)\in P\setminus (L_y\cup L_p\cup\Lambda_z\cup\Lambda_p)=L_z\cap \Lambda_y=\{x\circ p\}$ and hence $$x\circ(y\circ z)=x\circ p=(x\circ y)\circ(x\circ z),$$ witnessning that the midpoint operation $\circ$ on $X$ is self-distibutive.
\smallskip

$(2)\Ra(3)$ Assume that the midpoint operation on $X$ is self-distributive, fix any point $e\in X$, and consider the binary operation $\cdot:X\times X\to X$, $\cdot:(x,y)\mapsto e\circ(x\circ y)$, turning $X$ into a commutative loop of exponent $3$ with identity element $e$. We denote this loop by $X_e$. To prove that the loop $X_e$ is Moufang, fix any point $x,y,z\in X$. We need to prove that $(xy)(zx)=(x(yz))x$.
The definition of the binary operation on the loop $X_e$ and the commutativity, self-distributivity and involutarity of the midpoint operation $\circ$ ensure that
\begin{multline*}
(xy)(zx)=(xy)(xz)=e\circ\big((xy)\circ (xz)\big)=e\circ\big((e\circ(x\circ y))\circ (e\circ(x\circ z)\big)\\
=e\circ\big(e\circ((x\circ y)\circ(x\circ z))\big)=(x\circ y)\circ (x\circ z)=x\circ(y\circ z).
\end{multline*}
On the other hand, 
\begin{multline*}
$$(x(yz))x=e\circ(x(yz)\circ x)=(e\circ x(yz))\circ(e\circ x)=(e\circ(e\circ (x\circ yz)))\circ (e\circ x)\\
=(x\circ yz)\circ (x\circ e)=x\circ(yz\circ e)=x\circ (e\circ yz)=x\circ (e\circ (e\circ(y\circ z)))=x\circ(y\circ z).
\end{multline*}
Therefore,
$$(xy)(zx)=x\circ(y\circ z)=(xy)(zx),$$witnessing that the loop $X_e$ is Moufang.
\smallskip

$(3)\Ra(1)$ Assume that for every point $e\in X$, the commutative loop $X_e$ is Moufang. To prove that the Steiner liner $X$ is Hall, it suffices to check that every plane  in $X$ contains exactly $9$ points. First we prove this fact for planes which are flat hulls of $3$-element sets.

\begin{claim}\label{cl:plane=9} For every points $e,x,y\in X$ with $\|\{e,x,y\}\|=3$, the plane $P\defeq\overline{\{e,x,y\}}$ contains exactly $9$ points.
\end{claim}

\begin{proof} The equality $\|\{e,x,y\}\|=3$ implies $e\ne x$ and $y\notin\Aline ex$. By our assumption, the commutative loop $X_e$ is Moufang. By Lemma~\ref{l:Moufang-xyy} and Moufang Theorem~\ref{t:Moufang-loop}, the subloop $G$ of $X_e$ generated by the elements $x,y$ is associative. Since the loop $X_e$ has exponent $3$, $G$ is a commutative group of exponent $3$, generated by two elements $x,y$ such that $x\ne e$ and $y\notin\Aline ex=\{e,x,e\circ x\}=\{e,x,e\circ(x\circ x)\}=\{e,x,x^2\}$. Such group has cardinality $|G|=9$. By Proposition~\ref{p:flat<=>subloop}, the set $G$ is flat in the liner $X$ and hence $P=\overline{\{e,x,y\}}\subseteq G$. By Proposition~\ref{p:flat<=>subloop}, the plane $P$ is a subloop of the loop $X_e$ and hence $P$ is a subgroup of the group $G$. Taking into account that the subgroup $G$ is generated by the elements $x,y\in P\subseteq G$, we conclude that $G=P$ and hence $|P|=|G|=9$.
\end{proof}

 Now take any plane $P$ in $X$. Then $P\subseteq\overline{\{a,b,c\}}$ for some points $a,b,c\in X$ such that $3=\|P\|\le\|\{a,b,c\}\|\le 3$. Since the plane $P$ is not contained a line, we can choose points $p,q,r\in P$ such that $\|\{p,q,r\}\|=3$. Then $\overline{\{p,q,r\}}\subseteq P\subseteq\overline{\{a,b,c\}}$. By Claim~\ref{cl:plane=9}, $|\overline{\{p,q,r\}}|=9=|\overline{\{a,b,c\}}|$ and hence
$\overline{\{p,q,r\}}=P=\overline{\{a,b,c\}}$ and $|P|=9$.
\end{proof}

\begin{theorem}\label{t:Hall<=>regular} For a Hall liner $X$, the following conditions are equivalent:
\begin{enumerate}
\item $X$ is regular;
\item $X$ is weakly regular;
\item $X$ is $4$-regular;
\item $X$ is $4$-ranked;
\item $X$ is $4$-balanced;
\item $X$ is $4$-balanced with $|X|_4=27$;
\item $X$ contains no flat $F\subseteq X$ of rank $\|F\|=4$ and cardinality $|F|=81$;
\item for every point $e\in X$, the loop $X_e$ is associative.
\end{enumerate}
\end{theorem}

\begin{proof} By Theorem~\ref{t:Hall<=>Playfair+Steiner}, the Hall liner $X$ is Playfair, by Theorem~\ref{t:Playfair<=>}, the Playfair liner $X$ is $3$-regular, and by Proposition~\ref{p:k-regular<=>2ex}, the $3$-regular liner $X$ is $3$-ranked. We shall prove  the implications $(1)\Ra(2)\Ra(6)\Ra(4\wedge 5)\Ra(4\vee 5)\Ra(6)\Ra(7)\Ra(8)\Ra(3)\Ra(1)$. 
\smallskip

The  the implication $(1)\Ra(2)$ is trivial, see Proposition~\ref{p:sr=>r=>wr}.
\smallskip


$(2)\Ra(6)$ Assume that the Hall liner $X$ is weakly regular. By definition, the Hall liner $X$ is $2$-balanced with $|X|_2=3$. By definition, the Playfair liner $X$ is $1$-parallel. By Theorem~\ref{t:wr+k-parallel=>n-balanced}, the weakly regular $2$-balanced $3$-ranked $1$-parallel liner $X$ is $4$-balanced with
$$|X|_4=1+(|X|_2-1)\sum_{n=0}^{4-2}(1+|X|_2-1)^n=1+(3-1)\sum_{n=0}^23^n=1+2(1+3+9)=27.$$

The implication $(6)\Ra(5)$ is trivial. To see that $(6)\Ra(4)$, assume that $X$ is $4$-balanced with $|X|_4=27$. To see that $X$ is $4$-ranked, take any flats $A\subseteq B$ in $X$ of the same rank $\|A\|=\|B\|\le4$. We have to prove that $A=B$. This equality will follow from $A\subseteq B$ as soon as we check that $|A|=|B|<\w$. By definition, the Hall liner $X$ is $2$-balanced and $3$-balanced with $|X|_2\le|X|_3\le|X|_4=27$. Then for the cardinal $\kappa\defeq\|A\|=\|B\|\le 4$, we have $|A|=|B|=|X|_\kappa\le 27<\w$ and hence $A=B$.
\smallskip


$(4\vee 5)\Ra(6)$ Assume that the Hall liner $X$ is $4$-ranked or $4$-balanced. 

\begin{claim}\label{cl:Hall-27} For every points $e,a,b,c\in X$ with $\|\{e,a,b,c\}\|=4$ the flat $H\defeq \overline{\{e,a,b,c\}}$ contains exactly $27$ points. 
\end{claim}

\begin{proof} By Theorem~\ref{t:Hall<=>Moufang}, the loop $X_e$ is Moufang. Let $S$ be the subloop of the Moufang loop $X_e$, generated by the points $a,b,c$. By Proposition~\ref{p:flat<=>subloop}, the subloop $S$ of the loop $X_e$ is flat in the liner $X$ and the flat $H$ is a subloop of the loop $X_e$. Then $H=\overline{\{e,a,b,c\}}\subseteq S$ (because $S$ is flat) and $S\subseteq H$ (because $H$ is a subloop of $X$ containing the points $a,b,c$). Therefore, $H=S$. If the associator $[a,b,c]$ is not equal to $e$, then by Theorem~\ref{t:Moufang81}, $S$ has cardinality $|S|=81$ and contains a subgroup $G\subseteq S$ of cardinality $|G|=27$. By Proposition~\ref{p:flat<=>subloop}, the group $G$ is flat in the liner $X$. Assuming that $\|G\|<4$, we conclude that $27=|G|\le|X|_3=9$, which is a contradiction showing that $\|G\|=4$. Then $G\subseteq H$ are two flats in $X$ of the same rank $\|G\|=\|H\|=4$ and different cardinality, witnessing that the liner $X$ is not $4$-balanced and not $4$-ranked. But this contradicts our assumption. This contradiction shows that $[a,b,c]=e$. By Theorem~\ref{t:Moufang-loop}, the loop $S$ is associative and hence $S$ is a subgroup of the loop $X_e$. Taking into account that $S$ is a commutative group of exponent $3$, generated by three elements $a,b,c$, we conclude that $|H|=|S|=3^n$ for some $n\le 3$. Since $\|\{e,a,b,c\}\|=4$, the flat $P\defeq\overline{\{e,a,b\}}$ is a plane and hence $|P|=|X|_3=9$. Since $\|H\|=4\ne 3=\|P\|$, $3^2=|P|<|H|=3^n$ and hence $n=3$ and $|H|=3^3=27$.
\end{proof}

Now take any flat $F\subseteq X$ of rank $\|F\|=4$. Then $F\subseteq \overline{\{e,a,b,c\}}$ for some points $e,a,b,c\in X$ with $\|\{e,a,b,c\}\|=4$. Since $\|F\|=4$, there exist points $x,y,z,u\in F$ such that $y\ne x$, $z\notin\Aline xy$ and $u\notin\overline{\{x,y,z\}}$. The choice of the point $z\notin\Aline xy$ ensures that $\|\{x,y,z,u\}\|\ge\|\{x,y,z\}\|= 3$. Assuming that $\|\{x,y,z,u\}\|\ne 4$, we conclude that $\overline{\{x,y,z,u\}}=\overline{\{x,y,z\}}$, by the $3$-rankedness of $X$. But  the inclusion $u\in\overline{\{x,y,z,u\}}=\overline{\{x,y,z\}}$ contradicts the choice of the point $u$. This contradiction shows that $\|\{x,y,z,u\}\|=4$. By Claim~\ref{cl:Hall-27}, the flats $\overline{\{x,y,z,u\}}$ and $\overline{\{e,a,b,c\}}$ have cardinality $27$. Taking into account that $\overline{\{x,y,z,u\}}\subseteq F\subseteq\overline{\{e,a,b,c\}}$, we conclude that  $|\overline{\{x,y,z,u\}}|= |F|=|\overline{\{e,a,b,c\}}|=27$. Therefore, every flat $F\subseteq X$ of rank $\|F\|=4$ has cardinality $|F|=27$, witnessing that $X$ is balanced with $|X|_4=27$.
\smallskip

The implication $(6)\Ra(7)$ is trivial. 
\smallskip

$(7)\Ra(8)$ Assuming that for some point $e\in X$ the loop $X_e$ is non-associative, we can find three points $a,b,c\in X$ whose associator $[a,b,c]$ is not equal to $e$. By Theorem~\ref{t:Moufang81}, the subloop $S$ generated by the elements $a,b,c$ has cardinality $|S|=81$. By Proposition~\ref{p:flat<=>subloop}, the set $S$ is flat in the liner and hence $\overline{\{e,a,b,c\}}\subseteq S$. By Proposition~\ref{p:flat<=>subloop}, the flat $\overline{\{e,a,b,c\}}$ is a subloop of the loop $X_e$ and hence $S\subseteq \overline{\{e,a,b,c\}}$ (because $S$ is the smallest subloop containing the points $e,a,b,c$). Therefore, $S=\overline{\{e,a,b,c\}}$ and $\|S\|\le|\{e,a,b,c\}|\le 4$. Assuming that $\|S\|<4$, conclude that $|S|\le|X|_3=9$, which contradicts $|S|=81$. This contradiction shows that $S$ is a flat of rank $\|S\|=4$ and cardinality $|S|=81$, which contradicts (7).
\smallskip

$(8)\Ra(3)$ Assume that for every $e\in X$, the loop $X_e$ is associative. To prove that the liner $X$ is $4$-regular, take any flat $A\subseteq X$ of rank $\|A\|<4$ and any points $e\in A$ and $z\in X\setminus A$. We have to show that $\overline{A\cup\{z\}}=\Aline A{\Aline ez}$. If $\|A\|<3$, then  $\overline{A\cup\{z\}}=\Aline A{\Aline ez}$ by the $3$-regularity of the Hall liner $X$. So, we assume that $\|A\|=3$. In this case $A$ is a plane of cardinality $|A|=|X|_3=9$. Let $\mathcal L_e$ be the family of lines in the plane $A$ that contain the point $e$. Taking into account that $L\cap \Lambda=\{e\}$ for any distinct lines, we conclude that $$9=|A|=\big|\{e\}\cup\bigcup_{L\in\mathcal L_e}(L\setminus\{x\})\big|=1+|\mathcal L_e|(|X_2|-1)=1+|\mathcal L_e|\cdot 2$$ and hence $|\mathcal L_e|=4$.

For every line $L\in\mathcal L_e$, consider the plane $\overline{L\cup\{z\}}$ and observe that $|\overline{L\cup\{z\}}|=|X|_3=9$. The $3$-rankedness of the Hall liner $X$ ensures that for any distinct lines $L,\Lambda\in\mathcal L_e$, the intersection of the planes $\overline{L\cup\{z\}}$ and $\overline{\Lambda\cup\{z\}}$ coincides with the line $\Aline ez$. Then the set 
$$\Pi\defeq\bigcup_{L\in\mathcal L_e}\overline{L\cup\{z\}}=\Aline ez\cup\bigcup_{L\in\mathcal L_e}(\overline{L\cup\{z\}}\setminus\Aline ez)$$ has cardinality $$|\Pi|=|\Aline ez|+\sum_{L\in\mathcal L_o}|\overline{L\cup\{z\}}\setminus\Aline oz|=|X|_2+|\mathcal L_o|\cdot (|X|_3-|X|_2)=3+4\cdot 6=27.$$
Let $S$ be the subloop of the Moufang loop $X_e$, generated by the elements $e,x,y,z$. Since the loop $X_e$ is associative, $S$ is a commutative group of rank $3$, generated by $3$ elements. Then $S$ has cardinality $|S|\le 27$. By Proposition~\ref{p:flat<=>subloop}, the set $S$ is flat in the liner and hence $\Pi\subseteq \overline{\{e,x,y,z\}}\subseteq S$ and $27=|\Pi|\le |S|\le 27$, which implies 
$\Pi=\overline{\{e,x,y,z\}}=\overline{A\cup\{z\}}=S$. Then for every $p\in\overline{A\cup\{z\}}=\Pi$, there exists a line $L\in\mathcal L_o$ such that $p\in\overline{L\cup\{z\}}$. The $3$-regularity of the Hall liner $X$ ensures that $p\in\overline{L\cup\{z\}}=\Aline L{\Aline oz}\subseteq \Aline A{\Aline oz}$, witnessing that $\overline{A\cup\{z\}}=\Aline A{\Aline oz}$. This completes the proof of the $4$-regularity of the Hall liner $X$.
\smallskip

The final implication $(3)\Ra(1)$ follows from Theorem~\ref{t:HA}.
\end{proof}

A group $X$ is defined to be \index{nilpotent group}\index{group!nilpotent of class $n$}\defterm{nilpotent of class $n$} for $n\in\IN$ if the quotient group $X/Z$ of $X$ by its centre $Z\defeq\{z\in X:\forall x\in X\;(xz=zx)\}$ is nilpotent of class  $n-1$. The trivial group is nilpotent of class $0$. So, non-trivial commutative groups are nilpotent of class $1$. A group $X$ is nilpotent of class at most $2$ if and only if $[[x,y],z]=1$. Groups satisfying the identity $[[x,y],y]=1$ are called $2$-Engel.  
It is known \cite{Burnside1901}, \cite{Levi1942} that every group $X$ of exponent $3$ is $2$-Engel and every $2$-Engel group is nilpotent of class at most $3$. In particular, for every $m\in \IN$, the free Burnside group $B(m,3)=\langle x_1,\dots,x_m:x^3=1\rangle$ of exponent $3$ with $m$ generators is nilpotent of class at most $3$. The Burnside group $B(1,3)$ is nilpotent of class $1$ and has cardinality $3$, the Burnside group $B(2,3)$ is nilpotent of class 2 with $|B(2,3)|=3^3=27$, and for every $m\ge 3$, the Burnside group $B(m,3)$ is nilpotent of class $3$ with $|B(m,3)|=3^{C^1_m+C^2_m+C^3_m}$.

\begin{example} Every group $X$ of exponent $3$ endowed with the binary operation $$\circ:X\times X\to X,\quad\circ:(x,y)\mapsto x\circ y\defeq xy^{-1}x,$$ is a self-distributive idempotent involutory commutative magma. Consequently, $X$ endowed with the line relation $$\Af\defeq\big\{(x,y,z)\in X^3:y\in\{x,z,xz^{-1}x\}\big\}$$is a Hall liner. This Hall liner is regular if and only if the group $X$ is nilpotent of class at most $2$. 
\end{example}

\begin{proof} Given any elements $x,y,z\in X$, observe that
$$
\begin{aligned}
&x\circ x=xx^{-1}x=x,\\
&x\circ (x\circ y)=x\circ xy^{-1}x=x(x^{-1}yx^{-1})x=y,\\
&x\circ y=xy^{-1}x=xy^{-1}x(x^{-1}y)^3=xy^{-1}xx^{-1}yx^{-1}yx^{-1}xy=yx^{-1}y=y\circ x,\\
&(x\circ y)\circ(x\circ z)=(xy^{-1}x)\circ(xz^{-1}x)=xy^{-1}x(x^{-1}zx^{-1})xy^{-1}x=xy^{-1}zy^{-1}x=x\circ(y\circ z),
\end{aligned}
$$witnessing that the binary operation $\circ$ turns $X$ into an idempotent involutory commutative self-distributive magma. By Theorems~\ref{t:Steiner<=>midpoint} and \ref{t:Hall<=>Moufang}, the set $X$ endowed with the line relation
$$\Af\defeq\big\{(x,y,z)\in X^3:y\in\{x,z,xz^{-1}x\}\big\}$$is a Hall liner.

By Proposition~\ref{p:Steiner=>loop3} and Theorem~\ref{t:Hall<=>Moufang}, for every point $e\in X$, the binary operation
$$*:X\times X\to X,\;*:(x,y)\mapsto xy\defeq x*y\defeq e\circ (x\circ y)=e\circ (xy^{-1}x)=ex^{-1}yx^{-1}e,$$
turns $X$ into a commutative Moufang loop of exponent 3. 

If the group $X$ is nilpotent of class at most $2$, then all commutators in $X$ belong to the center of $X$. In particular,  for every $x,y,z\in X$, the commutator $(ze^{-1})(xe^{-1})(ze^{-1})^{-1}(xe^{-1})^{-1}=ze^{-1}xz^{-1}ex^{-1}$ commutes with all elements of $X$ and hence 
$$
\begin{aligned}
(x*y)*z&=z*(x*y)=ez^{-1}(x*y)z^{-1}e=ez^{-1}(ex^{-1}yx^{-1}e)z^{-1}e\\
&=\big(ex^{-1}ez^{-1}ze^{-1}xe^{-1}\big)ez^{-1}ex^{-1}yx^{-1}ez^{-1}e\\
&=ex^{-1}ez^{-1}\big(ze^{-1}xz^{-1}ex^{-1}\big)(yx^{-1}ez^{-1})e=ex^{-1}ez^{-1}(yx^{-1}ez^{-1})\big(ze^{-1}xz^{-1}ex^{-1}\big)e\\
&=ex^{-1}ez^{-1}yz^{-1}ex^{-1}e=ex^{-1}ez^{-1}yz^{-1}ex^{-1}e=x*(z*y)=x*(y*z),
\end{aligned}
$$
witnessing that the loop $X_e$ is associative. By Theorem~\ref{t:Hall<=>regular}, the Hall liner $X$ is regular.

Now assume conversely, that the Hall liner $X$ is regular. By Theorem~\ref{t:Hall<=>regular}, for every $e\in X$ the loop $X_e$ is associative. In particular, $X_e$ is associative for the identity element $e$ of the group $X$. To prove that the group $X$ is nilpotent of class at most 2, we need to show that for every elements $x,y\in X$ their commutator $xyx^{-1}y^{-1}$ belongs to the center of the group $X$ and hence commutes with every element $z\in X$. This follows from the chain of equalities:
$$
\begin{aligned}
(xyx^{-1}y^{-1})z&=
(xyx^{-1}y^{-1})(zxy)y^{-1}x^{-1}=xy\big(x^{-1}y^{-1}(zxy)y^{-1}x^{-1}\big)=xy\big(x*(y*(zxy))\big)\\
&=xy\big((y*(zxy))*x\big)=xy\big(y*((zxy)*x)\big)=xy\big(y*(x*(zxy))\big)\\
&=xy\big(y^{-1}x^{-1}(zxy)x^{-1}y^{-1}\big)=z(xyx^{-1}y^{-1}).
\end{aligned}
$$
\end{proof}

\begin{remark}\label{r:Moufang81} By a deep result of Bruck and Slaby \cite[Theorem 10.1]{Bruck1958}, every finitely generated commutative Moufang loop $X$ of exponent $3$ is centrally nilpotent, which implies that $|X|=3^n$ for some $n\in\IN$. This algebraic result implies that every Hall liner $X$ of finite rank has cardinality $|X|=3^n$ for some $n\in\IN$. This fact was also proved by Hall \cite{Hall1980} by group-theoretic methods, without involving the machinery of Moufang loops. Corollary~\ref{c:Moufang81} implies that every non-regular Hall liners of rank 4 are isomorphic and have cardinality $81$. The following table taken from \cite[VI.28.10]{HCD} shows the number $\hbar$ of non-isomorphic non-regular Halls liners of cardinality $3^n$.
\begin{center}
\begin{tabular}{c||c|c|c|c|c|c|c|c|c}
$n$&1&2&3&4&5&6&7&8&9\\
\hline
$\hbar$&0&0&0&1&1&3&12&$\ge 45$&?
\end{tabular}
\end{center}
The liner in Example~\ref{ex:HTS} is the unique non-regular Hall liner of cardinality $81$.
\end{remark}

\section{The regularity of Playfair liners}

Now we are able to prove the promised characterization of regular Playfair liners.

\begin{theorem}\label{t:Playfair<=>regular} For a Playfair liner $X$ the following conditions are equivalent:
\begin{enumerate}
\item $X$ is regular;
\item $X$ is weakly regular;
\item $X$ is ranked;
\item $X$ is balanced;
\item $X$ is $4$-ranked or $4$-balanced;
\item $X$ contains no flat $F\subseteq X$ of rank $\|F\|=4$ and cardinality $|F|=81$.
\end{enumerate}
\end{theorem}

\begin{proof} By Theorem~\ref{t:Playfair<=>}, the Playfair liner $X$ is affine, $3$-long and $3$-regular, by Proposition~\ref{p:k-regular<=>2ex}, the $3$-regular liner $X$ is $3$-ranked, and by Theorem~\ref{t:affine=>Avogadro}, the affine line $X$ is $2$-balanced. By definition, the Playfair liner $X$ is $1$-parallel.

We shall prove the implications $(1)\Ra(2)\Ra(3\wedge 4)\Ra(5)\Ra(6)\Ra(1)$. 
\smallskip

The implication $(1)\Ra(2)$ is trivial, see Proposition~\ref{p:sr=>r=>wr}.
\smallskip

$(2)\Ra(3\wedge 4)$. Assume that the liner $X$ is weakly regular. By Proposition~\ref{p:wr-ex<=>3ex} and Theorem~\ref{t:ranked<=>EP}, the weakly regular $3$-ranked liner $X$ is ranked. By Theorem~\ref{t:wr+k-parallel=>n-balanced}, the weakly regular $3$-ranked $2$-balanced $1$-parallel liner $X$ is balanced.
\smallskip

The implication $(3\wedge 4)\Ra(5)$ is trivial.
\smallskip

$(5)\Ra(6)$ Assume that the Playfair liner $X$ is $4$-ranked or $4$-balanced. Given any flat $F\subseteq X$ of rank $\|F\|=4$, we have to check that $|F|\ne 81$. If $X$ is $4$-long, then $X$ is regular, by Theorem~\ref{t:4-long-affine}. By Corollary~\ref{c:affine-cardinality}, $|F|=|F|_2^{\dim(F)}=|X|_2^3\ne 81$. If $X$ is not $4$-long, then $|X|_2=3$ and hence $X$ is a Steiner liner. By Theorem~\ref{t:Hall<=>Playfair+Steiner}, the Playfar Steiner liner $X$ is Hall. By Theorem~\ref{t:Hall<=>regular}, $|F|\ne 81$.
\vskip3pt

$(6)\Ra(1)$. Assume that the Playfair liner $X$ contains no flat $F$ of rank $\|F\|=4$ and cardinality $|F|=81$. If $X$ is $4$-long, then $X$ is regular, by Theorem~\ref{t:4-long-affine}. If $X$ is not $4$-long, then $X$ is Steiner and Hall, by Theorem~\ref{t:Hall<=>Playfair+Steiner}. By Theorem~\ref{t:Hall<=>regular}, the Hall liner $X$ is regular. 
\end{proof}

In Proposition~\ref{p:Playfair<=>3+a+3} we presented a first-order characterization of Playfair liners. In fact two first-order axioms (the affinity and $3$-regularity) in that characterization can be replaced by a single axiom, called the biaffinity.

\begin{definition} A liner $X$ is called \index{biaffine liner}\index{liner!biaffine}\defterm{biaffine} if for every points $o,a,c\in X$, $b\in\Aline ac$ and $y\in\Aline ob\setminus\Aline ac$, there exist unique points $x\in\Aline oa$ and $z\in\Aline oc$ such that $y\in\Aline xz$ and $\Aline xz\cap\Aline ac=\varnothing$.
\end{definition}

\begin{picture}(200,80)(-200,-10)
\put(-30,0){\line(1,0){60}}
\put(-30,0){\line(1,2){30}}
\put(30,0){\line(-1,2){30}}
\put(-20,20){\line(1,0){40}}
\put(0,0){\line(0,1){60}}

\put(0,0){\circle*{3}}
\put(-3,-11){$b$}
\put(-30,0){\circle*{3}}
\put(-34,-9){$a$}
\put(30,0){\circle*{3}}
\put(29,-9){$c$}
\put(0,60){\circle*{3}}
\put(-2,63){$o$}
\put(-20,20){\color{red}\circle*{3}}
\put(-29,18){$x$}
\put(0,20){\circle*{3}}
\put(2,13){$y$}
\put(20,20){\color{red}\circle*{3}}
\put(23,18){$z$}
\end{picture}

\begin{theorem}\label{t:Playfair2<=>} A liner $X$ is Playfair if and only if $X$ is biaffine and $3$-long.
\end{theorem}

\begin{proof} First assume that the liner $X$ is Playfair. By Theorem~\ref{t:Playfair<=>}, $X$ is $3$-long and affine. To prove that $X$ is biaffine, fix any points $o,a,c\in X$, $b\in\Aline ac$ and $y\in\Aline ob\setminus\Aline ac$.  Since $X$ is affine, there exist unique points $x\in\Aline oa$ and $z\in\Aline oc$ such that $\Aline xy\cap\Aline ab=\varnothing$ and $\Aline yz\cap\Aline bc=\varnothing$.  If $a=c$, then $b=a=c$ and $\Aline xy\cap \Aline ab=\varnothing=\Aline yz\cap\Aline bc$ implies $x=y=z$. So, $x=y$ and $z=y$ are unique points with $y\in\Aline xz$ and $\Aline xz\cap\Aline ac=\varnothing$. 

So, assume that $a\ne c$. If $b=a$, then $\Aline xy\cap\Aline ab=\varnothing$ implies $x=y$ and hence $y=x\in\Aline xz$ and $\Aline xz\cap\Aline ac=\Aline yz\cap\Aline bc=\varnothing$. The uniqueness of the point $z$ implies that $x=y$ and $z$ are unique points with $y\in\Aline xz$ and $\Aline xz\cap \Aline ac=\varnothing$.

If $b=c$, then $\Aline yz\cap\Aline bc=\varnothing$ implies $y=z$ and hence $y=z\in\Aline xz$ and $\Aline xz\cap\Aline ac=\Aline xy\cap\Aline ab=\varnothing$. The uniqueness of the point $x$ implies that $x$ and $z=y$ are unique points with $y\in\Aline xz$ and $\Aline xz\cap \Aline ac=\varnothing$.

So, we assume that $a\ne b\ne c$, which implies $a\ne c$.  If $y=o$, then $x\in\Aline oa$, $z\in\Aline ac$ and $\Aline xy\cap\Aline ab=\varnothing=\Aline yz\cap\Aline bc$ imply $x=y=z$ and hence $x=y$ and $z=y$ are unique points with $y\in\Aline xz$ and $\Aline xz\cap\Aline ac=\{y\}\cap\Aline ac=\varnothing$. 

So, we assume that $y\ne o$ and hence $y\notin\Aline oa\cup\Aline oc$. In this case $x\ne y\ne z$ and hence $\Aline xy$ and $\Aline yz$ are two lines in the plane $\overline{\{o,a,c\}}$ that contain the point $y$ and are disjoint with the line $\Aline ac$. Since $X$ is Playfair, $\Aline xy=\Aline yz$ and hence $y\in\Aline xy=\Aline yz=\Aline xz$ and $\Aline xz\cap\Aline ac=\Aline xy\cap\Aline ab=\varnothing$. The uniqueness of the points $x,z$ implies from the affinity of $X$.
\smallskip

Now assume that the liner $X$ is biaffine and $3$-long. First we show that $X$ is affine. Given any points $o,b,c\in X$ and $y\in \Aline ob\setminus\Aline bc$, we should prove that there exists a unique point $z\in\Aline oc$ such that $\Aline yz\cap\Aline bc=\varnothing$. By the biaffinity of $X$, for the point $a\defeq b$, there exist unique points $x\in\Aline oa=\Aline ob$ and $z\in\Aline oc$ such that $y\in\Aline xz$ and $\Aline xz\cap\Aline ac=\varnothing$. Assuming that $x\ne y$, we conclude that $\Aline xz=\Aline xy=\Aline ob=\Aline oa$ and hence $a\in\Aline xz\cap\Aline ac=\varnothing$, which is a contradiction showing that $x=y$. Then $\Aline yz\cap\Aline bc=\Aline xz\cap\Aline ac=\varnothing$. 
Assuming that $z'\in\Aline oc$ is another point such that $\Aline y{z'}\cap\Aline bc=\varnothing$, we conclude that $y=x\in\Aline x{z'}$ and $\Aline x{z'}\cap\Aline ac=\Aline y{z'}\cap\Aline bc=\varnothing$. The uniqueness of the points $x,z$ implies $z'=z$. Therefore, $z$ is a unique point on the line $\Aline oc$ such that $\Aline yz\cap\Aline bc=\varnothing$, witnessing that the liner $X$ is affine.

By Theorem~\ref{t:affine=>Avogadro}, the affine line $X$ is $2$-balanced. If $|X|_2\ge 4$, then $X$ is Playfair, by Theorem~\ref{t:4-long-affine}. So, assume that $|X|_2<4$, which means that $X$ is Steiner. We shall prove that $X$ is Hall. 

\begin{claim}\label{cl:hull3=9} For every points $o,a,c\in X$ with $\|\{o,a,c\}\|=3$, the plane $\overline{\{o,a,c\}}$ has cardinality $|\overline{\{o,a,c\}}|=9$.
\end{claim}

\begin{proof} It follows from $\|\{o,a,c\}\|=3$ that $a\ne c$ and $o\notin\Aline ac$. Since $|X|_2=3$, there exist unique points $b\in\Aline ac\setminus\{a,c\}$ and $y\in\Aline ob\setminus\{o,b\}$. Since $\Aline oa\cap\Aline ob=\{o\}\ne\{y\}$, the point $y$ does not belong to the lines $\Aline oa$ and $\Aline oc$. By the biaffinity of $X$, there exist unique points $x\in\Aline oa$ and $y\in\Aline oc$ such that $y\in\Aline xz$ and $\Aline xz\cap\Aline ac=\varnothing$. The latter condition implies that $x\ne a$ and $z\ne c$. Assuming that $x=o$, we conclude that $z\in\Aline xy\setminus\{x,y\}=\Aline oy\setminus\{o,y\}=\{b\}\subseteq\Aline ac$, which contradicts $\Aline xz\in\Aline ac=\varnothing$. This contradiction shows that $x\ne o$ and hence $\Aline oa=\{o,x,a\}$. By analogy we can prove that $\Aline oc=\{o,y,c\}$. By the biaffinity of $X$, there exist unique points $p\in\Aline bx$ and $q\in\Aline by$ such that $o\in\Aline pq$ and $\Aline pq\cap\Aline xz=\varnothing$. Repeating the above argument, we can prove that $\Aline bx=\{b,x,p\}$ and $\Aline bz=\{b,z,q\}$.

\begin{picture}(200,160)(-180,-20)
\put(-60,0){\line(1,0){120}}
\put(-60,0){\line(1,2){60}}
\put(60,0){\line(-1,2){60}}
\put(-30,60){\line(1,0){60}}
\put(0,0){\line(0,1){120}}
\put(0,0){\line(1,2){60}}
\put(0,0){\line(-1,2){60}}
\put(-60,120){\line(1,0){120}}
\put(-60,0){\line(1,1){120}}
\put(60,0){\line(-1,1){120}}

\put(0,0){\circle*{3}}
\put(-3,-10){$b$}
\put(-60,0){\color{orange}\circle*{3}}
\put(-64,-10){$a$}
\put(60,0){\color{red}\circle*{3}}
\put(59,-10){$c$}
\put(0,120){\circle*{3}}
\put(-2,125){$o$}
\put(-30,60){\color{red}\circle*{3}}
\put(-39,58){$x$}
\put(0,60){\circle*{3}}
\put(2,68){$y$}
\put(30,60){\color{orange}\circle*{3}}
\put(33,58){$z$}
\put(-60,120){\color{orange}\circle*{3}}
\put(-65,125){$p$}
\put(60,120){\color{red}\circle*{3}}
\put(60,125){$q$}
\end{picture}

We claim that $p,q\notin \Aline ac=\{a,b,c\}$. Indeed, assuming that $p=a$, we conclude that $q\in\Aline po=\Aline ao$ and hence $q\in \Aline ao\setminus\{o,x\}=\{a\}=\{p\}$ and $o=q=p\in\Aline bx\cap\Aline by=\{b\}$, which contradicts the choice of $o$. This contradiction shows that $p\ne a$. Assuming that $p=b$, we conclude that $q\in\Aline po=\Aline bo$ and $z\in\Aline bq\subseteq \Aline bo$. Since $|X|_2=3$, $z\in\Aline bo\setminus\{b,o\}=y$ and then $c\in\Aline oz\cap\Aline ac=\Aline oy\cap\Aline ac=\{b\}$, which contradicts the choice of $b$. This contradiction shows that $p\ne b$. Assuming that $b=c$, we conclude that $x\in\Aline pb=\Aline cb=\Aline ac$, which contradicts the choice of $x$. This contradiction shows that $p\ne c$ and hence $p\notin\Aline ac$. By analogy we can prove that $q\notin\Aline ac$ and hence $\Aline pq\cap\Aline ac=\{p,o,q\}\cap\Aline ac=\varnothing$. Applying the biaffinity of $X$ to the points $q,p,x$, $b\in\Aline px$ and $z\in\Aline qb\setminus\Aline px$ and taking into account that $o\in\Aline qp$ and $|X|_2=3$, we conclude $\{q,c,x\}$ is a line in $X$.
Applying the biaffinity of $X$ to the points $p,q,z$, $b\in\Aline qz$ and $x\in\Aline pb\setminus\Aline qz$ and taking into account that $o\in\Aline pq$ and $|X|_2=3$, we conclude $\{p,a,z\}$ is a line in $X$.

Applying the biaffinity of $X$ to the points $p,a,c$, $b\in\Aline ac$ and $x\in\Aline pb\setminus\Aline ac$ and taking into account that $\{p,z,a\}$ is a line, we conclude $\{p,y,c\}$ is a line in $X$.
By analogy we can prove that $\{q,y,a\}$ is a line. Now we see that the 9-point subset $\{a,b,c,x,y,z,p,o,q\}$ of the liner $X$ is flat and hence $|\overline{\{o,a,c\}}|=|\{a,b,c,x,y,z,p,o,q\}|=9$.
\end{proof}

Now we are able to show that the liner $X$ is $3$-balanced with $|X|_3=9$. Given any plane $P\subseteq X$, find points $a,b,c\in X$ such that $P\subseteq\overline{\{a,b,c\}}$. It follows from $3=\|P\|\le\|\{a,b,c\}\|\le 3$ that $\|\{a,b,c\}\|=3$ and hence $|\overline{\{a,b,c\}}|=9$, by Claim~\ref{cl:hull3=9}. Since $P$ is a plane, there exist points $x,y,z\in P$ such that $x\ne y$ and $z\notin\Aline xy$ and hence $\|\{x,y,z\}\|=3$. By Claim~\ref{cl:hull3=9}, $|\overline{\{x,y,z\}}|=9$. Taking into account that $\overline{\{x,y,z\}}\subseteq P\subseteq\overline{\{a,b,c\}}$ and $|\overline{\{x,y,z\}}|=9=|\overline{\{a,b,c\}}|$, we conclude that $\overline{\{x,y,z\}}=P=\overline{\{a,b,c\}}$ and $|P|=|\overline{\{a,b,c\}}|=9$. This shows that the biaffine Steiner liner $X$ is Hall. By Theorem~\ref{t:Hall<=>Playfair+Steiner}, the Hall liner $X$ is Playfair.
\end{proof}

\begin{exercise} Prove that a liner $(X,\Af)$ is  affine and regular if and only if for every points $o,x,y,z\in X$ and $p\in\overline{\{o,x,y,z\}}\setminus\overline{\{o,x,y\}}$ there exists a unique point $q\in\Aline oz$ such that $\Aline pq\cap\overline{\{o,x,y\}}=\varnothing$.
\end{exercise}

\begin{question} Is every finite balanced affine liner $X$ of rank $4$ regular?
\end{question}

\begin{example}[Ivan Hetman]\label{ex:Z15} There exists a non-regular ranked balanced affine Steiner liner $X$ with $\|X\|=3$ and $|X|=15$. It is the ring $\IZ_{15}$ endowed with the family of lines $$\mathcal L=\{x+L:x\in\IZ_{15},\;L\in\mathcal B\}\mbox{ \  where \ }\mathcal B=\big\{\{0,1,4\},\{0,6,8\},\{0,5,10\}\big\}.$$
\end{example}

\chapter{Modularity in liners}

\section{The submodularity of the rank in liners}

The inequality established in the following theorem is called the \defterm{submodularity} of the rank in ranked liners.

\begin{theorem}\label{t:submodular} For any flats $A,B$ in an ranked liner $X$, we have the inequality
$$\|A\cap B\|+\|A\cup B\|\le\|A\|+\|B\|.$$
If $\|A\cup B\|$ is infinite, then $\|A\cap B\|+\|A\cup B\|=\|A\|+\|B\|$.
\end{theorem} 

\begin{proof}By Corollary~\ref{c:ranked<=>EP}, the ranked liner $X$ has the Exchange Property.  By Lemma~\ref{l:Max-indep}, there exists a maximal independent set $M$ in the set $A\cap B$. By Theorem~\ref{t:Max=codim}, $|M|=\|A\cap B\|$. By Proposition~\ref{p:add-point-to-independent}, $\overline{M}=\overline{A\cap B}$, and hence $$\|M\|=\|A\cap B\|=|M|.$$  By Lemma~\ref{l:Max-indep}, there exist maximal $M$-independent sets $I\subseteq A$ and $J\subseteq B$. By Theorem~\ref{t:Max=codim}, $|I|=\|A\|_M$ and $|J|=\|B\|_M$. The $M$-independence of the sets $I$ and $J$ ensures that $I\cap \overline{M}=\varnothing=\overline{M}\cap J$. 
Also 
$I\cap J=(I\cap A)\cap J=I\cap (A\cap J)\subseteq I\cap (A\cap B)\subseteq I\cap \overline{M}=\varnothing$. Therefore, the sets $I,J,M$ are pairwise disjoint and hence $|I\cup M\cup J|=|I|+|M|+|J|$. 

By the maximality of the $M$-independent sets $I,J$ and Proposition~\ref{p:add-point-to-independent}, $A\subseteq \overline{M\cup I}$ and $B\subseteq \overline{M\cup J}$. Then
$$A\cup B\subseteq \overline{M\cup I}\cup\overline{M\cup J}\subseteq\overline{M\cup I\cup J},$$ which implies $\|A\cup B\|\le |M\cup I\cup J|=|I|+|M|+|J|$ and 
$$\|A\cap B\|+\|A\cup B\|=|M|+\|A\cup B\|\le|M|+|I|+|M|+|J|=\|M\|+\|A\|_M+\|M\|+\|B\|_M=\|A\|+\|B\|,$$
according to Corollary~\ref{c:rank+}.
\smallskip

If $\|A\cup B\|$ is infinite, then  $$\|A\cup B\|\le\|A\cup B\|+\|A\cap B\|\le\|A\|+\|B\|\le\|A\cup B\|+\|A\cup B\|=\|A\cup B\|$$ and hence $\|A\cup B\|+\|A\cap B\|=\|A\|+\|B\|$ by the Cantor-Bernstein Theorem.
\end{proof}

\section{Modular pairs of flats in liners}

\begin{definition} Let $A,B$ be two flats in a liner $X$. The pair $(A,B)$ is called \index{modular pair of flats}\defterm{modular} if $\|A\cup B\|+\|A\cap B\|=\|A\|+\|B\|$.
\end{definition}



\begin{lemma}\label{l:modular} A pair $(A,B)$ of flats $A,B$  in a ranked liner $X$ is modular if
$$A\cap \overline{B\cup C}=\overline{(A\cap B)\cup C}\quad\mbox{and}\quad\overline{D\cup A}\cap B=\overline{D\cup (A\cap B)}$$for every sets $C\subseteq A$ and $D\subseteq B$. 
\end{lemma} 

\begin{proof} By Lemma~\ref{l:Max-indep}, there exist maximal $(A\cap B)$-independent sets $I\subseteq A$ and $J\subseteq B$. By Theorem~\ref{t:Max=codim}, $|I|=\|A\|_{A\cap B}$ and $|J|=\|B\|_{A\cap B}$. By Proposition~\ref{p:add-point-to-independent}, $\overline{(A\cap B)\cup I}=A$ and $\overline{(A\cap B)\cup J}=B$. The $(A\cap B)$-independence of the set $I$ ensures that $I\cap(A\cap B)=\varnothing$ and hence $I\cap J\subseteq (I\cap A)\cap B=I\cap(A\cap B)=\varnothing$.

We claim that the set $I\cup J$ is $(A\cap B)$-independent. In the opposite case, we can find a point $x\in I\cup J$ such that $x\in\overline{(A\cap B)\cup((I\cup J)\setminus \{x\})}$. Since $I\cap J=\varnothing$, either $x\in I\setminus J$ or $x\in J\setminus I$. 

If $x\in I\setminus J$, then for the set $C\defeq I\setminus\{x\}$, the assumption of the lemma guarantees that  
$$
x\in A\cap\overline{(A\cap B)\cup C\cup J}\subseteq A\cap \overline{B\cup C}
 =\overline{(A\cap B)\cup C},
$$
which contradicts the $(A\cap B)$-independence of the set $I$.

If $x\in J\setminus I$, then for the set $D\defeq J\setminus\{x\}$, the assumption of the lemma guarantees that  
$$
x\in \overline{(A\cap B)\cup I\cup D}\cap B\subseteq  \overline{D\cup A}\cap B
 =\overline{D\cup (A\cap B)},
$$
which contradicts the $(A\cap B)$-independence of the set $J$.

 Those contradictions show that the set $I\cup J\subseteq A\cup B$ is $(A\cap B)$-independent and hence $|I\cup J|\le \|I\cup J\|_{A\cap B}\le \|A\cup B\|_{A\cap B}$, by 
Theorem~\ref{t:Max=codim}.  By Corollary~\ref{c:rank+} and Theorem~\ref{t:ranked<=>EP}, 
$$
\begin{aligned}
\|A\|+\|B\|&=\|A\cap B\|+\|A\|_{A\cap B}+\|A\cap B\|+\|B\|_{A\cap B}=\|A\cap B\|+|I|+\|A\cap B\|+|J|\\
&=\|A\cap B\|+|I\cup J|+\|A\cap B\|\le \|A\cap B\|+\|A\cup B\|_{A\cap B}+\|A\cap B\|\\
&=\|A\cup B\|+\|A\cap B\|\le\|A\|+\|B\|,
\end{aligned}$$
and hence $\|A\|+\|B\|=\|A\cap B\|+\|A\cup B\|.$
 \end{proof}

\begin{lemma}\label{l:aff-wmodular} Let $A,B$ be two  flats in a weakly regular liner $X$ possessing the Exchange Property. If $A\cap B\ne\varnothing$, then for every $n\in\w$ and every set $F\subseteq B$ of cardinality $|F|\le n$, we have 
$$
\overline{F\cup A}\cap B\subseteq \overline{F\cup(A\cap B)}.$$
\end{lemma}

\begin{proof} If $n=0$, then $F=\varnothing$ and $\overline{F\cup A}\cap B=A\cap B\subseteq\overline{F\cup(A\cap B)}$. Assume that  for some $n\in\w$, the inclusion $\overline{F\cup A}\cap B\subseteq \overline{F\cup(A\cap B)}$ has been proved for all sets $F\subseteq B$ of cardinality $|F|\le n$. 

Take any set $F\subseteq B$ of cardinality $|F|=n+1$. To prove that $\overline{F\cup A}\cap B\subseteq \overline{F\cup (A\cap B)}$, choose any element $x\in \overline{F\cup A}\cap B$. If $x\in \overline{E\cup A}$ for some subset $E\subset F$ of cardinality $|E|<|F|$, then by the inductive assumption, $x\in \overline{E\cup A}\cap B\subseteq \overline{E\cup (A\cap B)}\subseteq \overline{F\cup(A\cap B)}$ and we are done. So, we assume that $x\notin \overline{E\cup A}$ for every proper subset $E\subset F$. This implies $x\notin A$ and $F\cap A=\varnothing$.

Choose any point $y\in F$ and consider the set $E\defeq F\setminus\{y\}$ and the flat $\Lambda\defeq\overline{E\cup A}$. It follows from $|E|<|F|$ that  $x\notin\Lambda$ but $x\in\overline{F\cup A}=\overline{E\cup  A\cup\{y\}}\subseteq\overline{\Lambda\cup\{y\}}$ and hence $y\in\overline{\Lambda\cup\{x\}}\setminus\Lambda$, by the Exchange Property. Fix any point $o\in A\cap B$. The weak regularity of the liner $X$ ensures that $y\in\overline{\{\lambda, x,o\}}$ for some $\lambda\in\Lambda$. Since $y\in\overline{\{\lambda,x,o\}}\setminus\Lambda\subseteq\overline{\{\lambda,o,x\}}\setminus\Aline\lambda o$, the Exchange Property implies  $x\in\overline{\{ \lambda, o,y\}}\setminus\Lambda\subseteq\overline{\{\lambda,o,y\}}\setminus
\Aline \lambda o$. If $x\in\Aline oy$, then $x\in \Aline yo\subseteq\overline{F\cup(A\cap B)}$ are we are done. So, assume that  $x\notin\Aline oy$. Then $x\in \overline{\{o,y,\lambda\}}\setminus\Aline oy$ and the Exchange Property ensures that $\lambda\in\overline{\{o,y,x\}}\subseteq B$ and hence $\lambda\in \Lambda\cap B=\overline{E\cup A}\cap B$. Since $|E|=n$, the inductive assumption ensures that $\lambda\in \overline{E\cup A}\cap B\subseteq \overline{E\cup(A\cap B)}$. Then $x\in \overline{\{y,\lambda,o\}}\subseteq
\overline{F\cup\overline{E\cup(A\cap B)}}= \overline{F\cup(A\cap B)}$.
\end{proof}

\begin{proposition}\label{p:normal=>modular} For every interecting flats $A,B$ in a weakly regular ranked liner $X$, the pair $(A,B)$ is modular.
\end{proposition}

\begin{proof} By Corollary~\ref{c:ranked<=>EP}, the ranked liner $X$ has the Exchange Property. Fix any two intersecting flats $A,B$ in $X$.  The modularity of the pair $(A,B)$ will follow from Lemma~\ref{l:modular} as soon as we check that $\overline{E\cup A}\cap B=\overline{E\cup(A\cap B)}$ for every set $E\subseteq B$.  Take any point $x\in\overline{E\cup A}\cap B$. By Proposition~\ref{p:aff-finitary}, there exists a finite set $F\subseteq E$ such that $x\in\overline{F\cup A}$. By Lemma~\ref{l:aff-wmodular}, $x\in \overline{F\cup A}\cap B\subseteq\overline{F\cup(A\cap B)}\subseteq\overline{E\cup(A\cap B)}$ and hence $\overline{E\cup A}\cap B\subseteq\overline{E\cup(A\cap B)}$. The converse inclusion $\overline{E\cup(A\cap B)}\subseteq \overline{E\cup A}\cap B$ is trivial.
\end{proof}

\section{Modular and weakly modular liners}

\begin{definition} A liner $X$ is called (\index{weakly modular liner}\index{liner!weakly modular}\defterm{weakly}) \index{modular liner}\index{liner!modular}\defterm{modular} if its rank function is (\defterm{weakly}) \defterm{modular} in the sense that
$$\|A\cup B\|+\|A\cap B\|=\|A\|+\|B\|$$
 for any flats $A,B$ in $X$ (with $A\cap B\ne\varnothing$).
\end{definition}

Observe that a liner $X$ is (weakly) modular if and only if every pair of flats $A,B\subseteq X$ (with $A\cap B\ne \varnothing$) is modular.

\begin{theorem}\label{t:modular} For a ranked liner $X$ the following conditions are equivalent:
\begin{enumerate}
\item $X$ is (weakly) modular;
\item $A\cap\overline{B\cup C}=\overline{(A\cap B)\cup C}$ for any flats $A,B\subseteq X$ and set $C\subseteq A$ (with $A\cap B\ne\varnothing$);
\item $A\cap\overline{B\cup C}=\overline{(A\cap B)\cup C}$ for any flats $A,B\subseteq X$ and set $C\subseteq A$ with $\|A\cup B\|<\w$ (and $A\cap B\ne\varnothing$).
\end{enumerate}
\end{theorem}

\begin{proof} We shall prove the implications $(1)\Ra(3)\Ra(2)\Ra(1)$.
\smallskip

$(1)\Ra(3)$ Assume that the liner $X$ is (weakly) modular and take any flats $A,B\subseteq X$ and set $C\subseteq A$ such that $\|A\cup B\|<\w$ (and $A\cap B\ne\varnothing$). We should prove the equality $A\cap\overline{B\cup C}=\overline{(A\cap B)\cup C}$. We lose no generality assuming that $A\cap B\subseteq C$. In this case $\overline{(A\cap B)\cup C}=\overline{C}$. The equality $\overline{(A\cap B)\cup C}=A\cap\overline{B\cup C}$ will follow from the obvious inclusion $\overline{C}\subseteq A\cap\overline{B\cup C}$ and the rankedness of $X$ as soon as we check that 
$\|C\|\ge \|A\cap\overline{B\cup C}\|$.
Take a maximal independent set $M\subseteq A\cap B$, a maximal $(A\cap B)$-independent set $I\subseteq C$ and a maximal $(A\cap B)$-independent set $J\subseteq B$. The $(A\cap B)$-independence of $I,J$ ensures that $I\subseteq C\setminus(A\cap B)\subseteq A\setminus B$ and $J\subseteq B\setminus A$.
By the (weak) modularity,
$$
\begin{aligned}
\|A\cap\overline{B\cup C}\|&=\|A\|+\|\overline{B\cup C}\|-\|A\cup B\cup C\|=
\|A\|+\|B\cup \overline{C}\|-\|A\cup B\|\\
&=\|A\|+(\|B\|+\|\overline{C}\|-\|B\cap\overline{C}\|)-(\|A\|+\|B\|-\|A\cap B\|)\\
&=\|C\|-\|B\cap \overline{C}\|+\|A\cap B\|\le\|C\|
\end{aligned}
$$
because $A\cap B\subseteq B\cap C\subseteq B\cap\overline{C}$ and hence $\|A\cap B\|\le\|B\cap\overline{C}\|$.
\smallskip

$(3)\Ra(2)$ Given any flats $A,B\subseteq X$ and set $C\subseteq A$ (with $A\cap B\ne\varnothing$), we should prove that  $A\cap\overline{B\cup C}=\overline{(A\cap B)\cup C}$. The inclusion $\overline{(A\cap B)\cup C}\subseteq A\cap\overline{B\cup C}$ is trivial. To prove the inclusion $A\cap\overline{B\cup C}\subseteq \overline{(A\cap B)\cup C}$, take any point $a\in A\cap\overline{B\cup C}$. By Proposition~\ref{p:aff-finitary}, there exist finite sets $F\subseteq B$ and $C'\subseteq C$ such that $a\in\overline{F\cup C'}$. If $A\cap B\ne\varnothing$, then we can additionally assume that $F\cap A\ne\varnothing$. Consider the flats $A'=\overline{\{a\}\cup C'\cup(A\cap F)}$ and $B'=\overline{F}$, and observe that $a\in A'\cap\overline{B'\cup C'}$ and $\|A'\cup B'\|=\|\{a\}\cup C'\cup F\|<\w$. By (3),
$a\in A'\cap\overline{B'\cup C'}=\overline{(A'\cap B')\cup C'}\subseteq \overline{(A\cap B)\cup C}$.
\smallskip

The implication $(2)\Ra(1)$ follows from Lemma~\ref{l:modular}.
\end{proof}

\begin{exercise} Let $F:X\to Y$ be a liner isomorphism between liners $X,Y$. Prove that the liner $X$ is (weakly) modular if and only if so is the liner $Y$.
\smallskip

\noindent{\em Hint:} Apply Corollary~\ref{c:flat=injection=>rank=}.
\end{exercise}

\section{Characterizing weakly modular liners}

\begin{theorem}\label{t:w-modular<=>}For a liner $X$ the following conditions are equivalent:
\begin{enumerate}
\item $X$ is weakly modular;
\item $X$ is ranked and for every planes $P\cup\Pi$ with $\|P\cup\Pi\|=4$, the intersection $P\cap\Pi$ is not a singleton;
\item $X$ is ranked and weakly regular.
\end{enumerate}
\end{theorem}

\begin{proof} $(1)\Ra(2)$ Assume that $X$ is weakly modular. To prove that $X$ is ranked, take any flats $A\subseteq B$ of finite rank $\|A\|=\|B\|$ in $X$. If $\|A\|=\|B\|=0$, then $A=\varnothing=B$ and we are done. So, assume that $\|A\|>0$ and choose a point $a\in A$. Assuming that $A\ne B$, choose any points $b\in B\setminus A$. By the weak modularity of $X$, $$1=\|\Aline ab\cap A\|=\|\Aline ab\|+\|B\|-\|A\cup\Aline ab\|\ge 2+\|B\|-\|B\|=2,$$which is a desired contradiction showing that $A=B$. Therefore, the liner $X$ is ranked. By Corollary~\ref{c:ranked<=>EP}, the ranked liner $X$ has the Exchange Property.

Now we check that for every planes $P,\Pi$ in $X$ with $\|P\cup\Pi\|=4$, the intersection $P\cap\Pi$ is not a singleton. Asuming that $P\cap\Pi$ is not empty, we can apply the weak modularity and conclude that
$$\|P\cap\Pi\|=\|P\|+\|\Pi\|-\|P\cup\Pi\|=3+3-4=2,$$
witnessing that $P\cap\Pi$ is not a singleton.
\smallskip

$(2)\Ra(3)$. Assume that the liner $X$ satisfies the condition (2). To prove that $X$ is weakly regular, we need to check that for every flat $A\subseteq X$ and points $o\in A$ and $b\in X\setminus A$, we have $\overline{A\cup\{b\}}=\bigcup_{a\in A}\overline{\{a,o,b\}}$. This equality will follow as soon as we check that the set $\Lambda\defeq\bigcup_{a\in A}\overline{\{a,o,b\}}$ is flat. Given any points $p,q\in \Lambda$ and $r\in\Aline pq$, we should find a point $a\in A$ such that $r\in\overline{\{a,o,b\}}$. If $r\in\Aline ob$, then the point $a\defeq o$ has the required property. So, we assume that $r\notin\Aline ob$ and hence $\overline{\{r,o,b\}}$ is a plane.

 For the points $p,q$, find points $x,y\in A$ such that $p\in\overline{\{x,o,b\}}$ and $q\in\overline{\{y,o,b\}}$. If $y\in \overline{\{x,o,b\}}$, then $r\in\overline{\{x,o,b,y\}}=\overline{\{x,o,b\}}$ and the point $a\defeq x$ has the required property. If $x\in \overline{\{y,o,b\}}$, then $r\in\Aline pq\subseteq \overline{\{x,o,b,y\}}=\overline{\{y,o,b\}}$ and the point $a\defeq y$ has the required property.
So, we assume that $x\notin\overline{\{y,o,b\}}$ and $y\notin\overline{\{x,o,b\}}$.  Then $\overline{\{x,o,y\}}$ and $\overline{\{o,b,r\}}$ are two distinct planes in $X$ with $o\in\overline{\{x,o,y\}}\cap\overline{\{o,b,r\}}$. Since $\overline{\{x,o,y\}}\cup\overline{\{r,o,b\}}\subseteq\overline{\{x,o,y,p,q,b\}}\subseteq\overline{\{x,o,y,b\}}$, the rankedness of $X$ ensures that $\|\overline{\{x,o,y\}}\cup\overline{\{r,o,b\}}\|=4$. By the condition (2), the intersection $\overline{\{x,o,y\}}\cap\overline{\{r,o,b\}}$ contains some point $a\ne o$. Then $\Aline oa$ is a line in the plane $\overline{\{x,o,y\}}\subseteq A$. Since $b\notin A$, the flat $\overline{\{o,a,b\}}$ is a plane in the plane $\overline{\{r,o,b\}}$. The rankedness of $X$ guarantees that $r\in\overline{\{r,o,b\}}=\overline{\{a,o,b\}}$, witnessing that that the set $\Lambda$ is flat and the liner $X$ is weakly regular.
\smallskip

The implication $(3)\Ra(1)$ follows from  Proposition~\ref{p:normal=>modular}.  
\end{proof}


Theorems~\ref{t:w-modular<=>} and \ref{t:Playfair<=>regular} implies the following characterization of weakly modular Playfair liners.

\begin{corollary} For a Playfair liner $X$ the following conditions are equivalent:
\begin{enumerate}
\item $X$ is weakly modular;
\item $X$ is regular;
\item $X$ is weakly regular;
\item $X$ is ranked;
\item $X$ is balanced;
\item $X$ is $4$-ranked or $4$-balanced;
\item $X$ contains no flat $F$ of rank $\|F\|=4$ and cardinality $|F|=81$.
\end{enumerate}
\end{corollary}

\begin{exercise} Find an example of an affine liner which is not weakly modular.
\smallskip

\noindent{\em Hint:} Look at the liner from Example~\ref{ex:HTS}.
\end{exercise}

\begin{example} There exists a proaffine ranked liner which is not weakly modular and not balanced.
\end{example}

\begin{proof} Let $H$ be a Hall liner of cardinality $|H|=9$ and let $\Af_H$ be the line relation of $H$. Take any set $X$ containing $H$ so that $X\setminus H=\{a\}$ for some $a$, and consider the line relation $$\Af=\Af_H\cup\{(x,y,z):a\in \{x,y,z\}\;\wedge\;|\{x,y,z\}|\le 2\}.$$It is easy to see that $X$ endowed with the line relation $\Af$ is a liner. Since every line in $X$ has cardinality at most $3$, the liner $X$ is proaffine. Since $X$ contains lines of cardinality $2$ and $3$, $X$ is not $2$-balanced. The liner $X$ contains planes of cardinality $9$ and $4$ and hence $X$ is not $3$-balanced. It is easy to see that the liner $X$ is ranked and has rank $\|X\|=4$. Choose any disjoint lines $L,\Lambda$ in the Hull liner $H$ and observe that $P\defeq L\cup\{a\}$ and $\Pi\defeq\Lambda\cup\{a\}$ are two planes in $X$ such that $P\cap \Pi=\{a\}$ and $\overline{L\cup\Pi}=X$. By Theorem~\ref{t:w-modular<=>}, the liner $X$ is not weakly regular.
\end{proof}

\begin{question} Is every affine ranked liner weakly modular?
\end{question}

\section{Characterizing modular liners}

\begin{theorem}\label{t:modular<=>}
For a liner $X$, the following conditions are equivalent:
\begin{enumerate}
\item $X$ is modular;
\item $X$ is strongly regular;
\item $X$ is projective.
\end{enumerate}
\end{theorem}

\begin{proof} The equivalence $(2)\Leftrightarrow(3)$ has been proved in Theorem~\ref{t:projective<=>}. So, it remains to prove that $(1)\Leftrightarrow(2)$.
\smallskip

$(1)\Ra(2)$. Assume that the liner $X$ is modular. By Theorem~\ref{t:w-modular<=>} and \ref{t:ranked<=>EP}, the modular liner $X$ is ranked and has the Exchange Property. To prove that $X$ is strongly regular, take any nonempty flat $A$ and points $b\in X\setminus A$ and $x\in\overline{A\cup\{b\}}$. We have to find a point $a\in A$ such that $x\in\Aline ab$. If $x=b$, then any point $a\in A$ has the required property. So, we assume that $x\ne b$.  By Proposition~\ref{p:aff-finitary}, there exists a finite set $F\subseteq A$ such that $x\in\overline{F\cup \{b\}}$. We can assume that the cardinality of the set $F$ is the smallest possible. Then $F$ is independent. It follows from $x\in\overline{\{b\}\cup F}$ that $\|\Aline xb\cup F\|=\|\{b\}\cup F\|\le 1+\|F\|$. The modularity of $X$ implies the equality
$$\|\Aline xb\cap \overline F\|=\|\Aline xb\|+\|F\|-\|\Aline xb\cup F\|=2+\|F\|-\|\{b\}\cup F\|\ge2+\|F\|-(1+\|F\|)=1>0,$$witnessing that $\Aline xb\cap\overline F$ is not empty and contains some point $a\in A\subseteq X\setminus\{b\}$. By Theorem~\ref{t:Alines}, $a\in\Aline xb$ implies $x\in\Aline xb=\Aline ab$, witnessing that the liner $X$ is strongly regular.
\smallskip

$(2)\Ra(1)$ Assume that the liner $X$ is strongly regular. Then $X$ is weakly regular and ranked, by Corollary~\ref{c:proregular=>ranked}. By Theorem~\ref{t:w-modular<=>}, $X$ is weakly modular. To show that $X$ is modular, take any flats $A,B\subseteq X$. If $A\cap B\ne\emptyset$, then $\|A\cap B\|+\|A\cup B\|=\|A\|+\|B\|$, by the weak modularity of $X$. So, assume that $A\cap B=\emptyset$. By Corollary~\ref{c:projective<}, $\|B\|_A=\|B\|$, and by Corollary~\ref{c:rank+}, 
$$\|A\cup B\|+\|A\cap B\|=\|A\|+\|B\|_A+0=\|A\|+\|B\|,$$
witnessing that the liner $X$ is modular.
\end{proof}

\section{The regularity of balanced liners}

In this section we shall apply the weak modularity to prove that that a finite balanced liner of rank $\ge 4$ is regular if and only if it is Proclus. We start with several lemmas having an independent value. 

\begin{lemma}\label{l:P=union-of-lines} Let $X$ be a weakly modular liner of rank $\|X\|\ge 4$. If some plane $P\subseteq X$ contains two disjoint lines, then $P$ is the union of a family of pairwise disjoint lines.
\end{lemma}

\begin{proof} Let $L,L'$ be two disjoint lines in the plane $P$. Since $\|X\|\ge 4>3=\|P\|$, there exists a point $o\in X\setminus P$. By Theorem~\ref{t:w-modular<=>}, the weakly modular liner $X$ is ranked.

\begin{claim}\label{cl:intersectP} For every line $\Lambda\subseteq P$ we have $\overline{\Lambda\cup\{o\}}\cap P=\Lambda$.
\end{claim}

\begin{proof} It is clear that $\Lambda\subseteq\overline{\Lambda\cup\{o\}}\cap P$ and $2=\|\Lambda\|\le\|\overline{\Lambda\cup\{o\}}\|\le 3$. Assuming that $\Lambda\ne\overline{\Lambda\cup\{o\}}\cap P$, we conclude that $2=\|\Lambda\|<\|\overline{\Lambda\cup\{o\}}\cap P\|\le\|\overline{\Lambda\cup\{o\}}\|\le 3$, by the rankedness of $X$. Since $\|\overline{\Lambda\cup\{o\}}\cap P\|=\|\overline{\Lambda\cup\{o\}}\|=3$, the rankedness of $X$ ensures that $o\in \overline{\Lambda\cup\{o\}}=\overline{\Lambda\cup\{o\}}\cap P\subseteq P$, which contradicts the choice of the point $o$. This contradiction shows that $\overline{\Lambda\cup\{o\}}\cap P=\Lambda$.
\end{proof}

By Theorem~\ref{t:w-modular<=>}, the flat $\Lambda=\overline{L\cup\{o\}}\cap\overline{L'\cup\{o\}}$ has rank $\|\Lambda\|\ge 2$. Assuming that $\|\Lambda\|>2$, we conclude that $3\le\|\Lambda\|\le\max\{\|\overline{L\cup\{o\}}\|,\|\overline{L'\cup\{o\}}\|\}\le 3$ and hence $\Lambda=\overline{L\cup\{o\}}=\overline{L'\cup\{o\}}$, by the rankedness of the liner $X$.

By Claim~\ref{cl:intersectP}, $L=\overline{L\cup\{o\}}\cap P=\overline{L'\cup\{o\}}\cap P=L'$, 
which contradicts the choice of the disjoint lines $L,L'$. This contradiction shows that $\|\Lambda\|=2$ and hence $\Lambda$ is a line. Claim~\ref{cl:intersectP} implies
$$\Lambda\cap P=(\overline{L\cup\{o\}}\cap\overline{L'\cup\{o\}})\cap P=(\overline{L\cup\{o\}}\cap P)\cap(\overline{L'\cup\{o\}}\cap P)=L\cap L'=\varnothing.$$
The rankedness of $X$ and Theorem~\ref{t:w-modular<=>} imply that for every point $p\in P$, the flat
$$L_p\defeq \overline{\Lambda\cup\{p\}}\cap P$$is a line.
Observe that
$$L_p=\overline{\Lambda\cup\{p\}}\cap P\subseteq\overline{\Lambda\cup L_p}\cap P\subseteq\overline{\Lambda\cup \overline{\Lambda\cup\{p\}}}\cap P=\overline{\Lambda\cup\{p\}}\cap P=L_p$$and hence $L_p=\overline{\Lambda\cup L_p}\cap P$. 
Then for every point $x\in L_p$, we have 
$$L_x=\overline{\Lambda\cup\{x\}}\cap P\subseteq \overline{\Lambda\cup L_p}\cap P=L_p.$$
Then for any points $p,q\in X$ with $L_p\cap L_q\ne\varnothing$, we can choose a point $x\in L_p\cap L_q$ and conclude that $L_p=L_x=L_q$. This shows that $\{L_p:p\in P\}$ is the required family of parvise disjoint lines covering the plane $P$.
\end{proof}

\begin{corollary}\label{c:x2|x3} Let $X$ be a balanced weakly modular liner of rank $\|X\|\ge 4$. If $X$ is not projective, then $|X|_3=|X|_2\cdot\kappa$ for some cardinal $\kappa\ge |X|_2$.
\end{corollary}

\begin{proof} By Theorem~\ref{t:w-modular<=>}, the weakly modular liner $X$ is ranked. Since $\|X\|\ge 4$, the balanced liner $X$ contains flats of every rank $r\le 4$, so the cardinals $|X|_2$ and $|X|_3$ are well-defined. It is clear that $2\le|X|_2\le|X|_3$. If the cardinal $|X|_3$ is infinite, then $|X|_3=|X|_2\cdot\kappa$ for a unique cardinal $\kappa\in\{|X|_2,|X|_3\}$. So, assume that $|X|_3$ is finite and so is the cardinal $|X|_2\le |X|_3$. If $X$ is not projective, then $X$ is not $0$-parallel, by Theorem~\ref{t:projective<=>}. And hence there exists a plane $P\subseteq X$ containing two disjoint lines. By Lemma~\ref{l:P=union-of-lines}, the plane $P$ is the union of a disjoint family $\F$ of lines, which implies that $|X|_3=|P|=|X|_2\cdot|\F|$. Proposition~\ref{p:cov-aff} implies that $|\F|\ge|X|_2$.
\end{proof}

\begin{lemma}\label{l:disjoint-planes} Let $X$ be a $2$-balanced and $3$-balanced weakly modular liner. If $|X|_3<\w$ and $X$ is not projective, then for every flat $Y\subseteq X$ of rank $\|Y\|=4$, every plane $P\subseteq Y$ and every point $y\in Y\setminus P$, there exists a plane $\Pi\subseteq Y$ such that $y\in\Pi\subseteq Y\setminus P$.
\end{lemma}

\begin{proof} Let $Y$ be any flat of rank $\|Y\|=4$ in $X$, $P$ be any plane in $Y$ and $y\in Y\setminus P$ be any point. Assume that $|X|_3<\w$ and the liner $X$ is not projective. By Proposition~\ref{p:23-balance=>k-parallel}, the liner $X$ is $p$-parallel for a unique (finite) cardinal $p$ such that $|X|_3=1+(|X|_2+p)(|X|_2-1)$. 
By Theorem~\ref{t:w-modular<=>}, the weakly modular liner $X$ is ranked and weakly regular. By Theorem~\ref{t:wr+k-parallel=>n-balanced}, for every finite cardinal $n\le\|X\|$, the weakly regular liner $X$ is $n$-balanced and
$$|X|_n=1+(|X|_2-1)\sum_{r=0}^{n-2}(|X|_2+p-1)^r.$$
It will be convenient to denote the cardinal $|X|_n$ by $x_n$ and the cardinal $|X|_2+p-1$ by $\ell$.
With those notations, we have $$x_n=1+(x_2-1)\sum_{r=0}^{n-2}\ell^r$$for every $n\le\|X\|$.

Let $\mathcal P_y$ be the family of all planes in $Y$ that contain the point $y$. Since the liner $X$ is ranked, every plane $\Pi\in\mathcal P_y$ is equal to $\overline{\{x,y,z\}}$ for some points $x\in Y\setminus\{y\}$ and $z\in Y\setminus\Aline xy$. This implies the equality 
$$|\mathcal P_y|=\frac{(x_4-1)(x_4-x_2)}{(x_3-1)(x_3-x_2)}.$$

Assuming that every plane $\Pi\in\mathcal P_y$ intersects the plane $P$, we can apply Theorem~\ref{t:w-modular<=>} and conclude that $\Pi=\overline{L\cup\{y\}}$ for some line $L$ in $P$, which is equal to $\Pi\cap P$. This implies that the cardinality of $\mathcal P_y$ is equal to the cardinality of the family $\mathcal L_P$ of all lines in the plane $P$ and hence
$$|\mathcal P_y|=|\mathcal L_P|=\frac{x_3(x_3-1)}{x_2(x_2-1)}.$$
Therefore, we have the equation
$$\frac{(x_4-1)(x_4-x_2)}{(x_3-1)(x_3-x_2)}=|\mathcal P_y|=|\mathcal L_P|=\frac{x_3(x_3-1)}{x_2(x_2-1)}.$$
Taking into account that $$x_3=1+(x_2-1)(1+\ell)=x_2(1+\ell)-\ell\quad\mbox{and}\quad
x_4=1+(x_2-1)(1+\ell+\ell^2),$$
we conclude that
$$
\begin{aligned}
x_3-x_2&=(x_2-1)(1+\ell)-(x_2-1)=(x_2-1)\ell\quad\mbox{and}\\
x_4-x_2&=(x_2-1)(1+\ell+\ell^2)-(x_2-1)=(x_2-1)\ell(1+\ell)
\end{aligned}
$$
and hence
$$
\begin{aligned}
1+\ell+\ell^2&=\frac{(x_2-1)(1+\ell+\ell^2)(x_2-1)\ell(1+\ell)}{(x_2-1)(1+\ell)(x_2-1)\ell}=
\frac{(x_4-1)(x_4-x_2)}{(x_3-1)(x_3-x_2)}\\
&=\frac{x_3(x_3-1)}{x_2(x_2-1)}=\frac{(x_2(1+\ell)-\ell)(x_2-1)(1+\ell)}{x_2(x_2-1)}=\frac{(x_2(1+\ell)-\ell)(1+\ell)}{x_2}
\end{aligned}
$$
and $$x_2(1+\ell+\ell^2)=x_2(1+2\ell+\ell^2)-\ell(1+\ell).$$
After cancellations, we obtain
$$x_2=1+\ell=1+(x_2+p-1)=x_2+p$$
and hence $p=0$. Therefore, the liner $X$ is $0$-parallel. By Theorem~\ref{t:projective<=>}, the liner $X$ is projective, which contradicts our assumption. This contradiction shows that some plane $\Pi\in\mathcal P_y$ is disjoint with the plane $P$.
\end{proof} 

The following lemma is a higher-dimensional counterpart of Lemma~\ref{l:P=union-of-lines}.

\begin{lemma}\label{l:Y=union-of-planes} Let $X$ be a weakly modular liner of rank $\|X\|\ge 5$. If some flat $Y\subseteq X$ of rank $\|Y\|=4$  contains two disjoint planes, then $Y$ is the union of a family of pairwise disjoint planes.
\end{lemma}

\begin{proof} Let $Y\subseteq X$ be a flat of rank $\|Y\|=4$ and $P,P'$ be two disjoint planes in $Y$. Since $\|X\|\ge 5>4=\|Y\|$, there exists a point $o\in X\setminus Y$. By Theorem~\ref{t:w-modular<=>}, the weakly modular liner $X$ is ranked.

\begin{claim}\label{cl:intersectY} For every plane $\Pi\subseteq Y$ we have $\overline{\Pi\cup\{o\}}\cap Y=\Pi$.
\end{claim}

\begin{proof} It is clear that $\Pi\subseteq\overline{\Pi\cup\{o\}}\cap Y$ and $3=\|\Pi\|\le\|\overline{\Pi\cup\{o\}}\|\le 4$. Assuming that $\Pi\ne\overline{\Pi\cup\{o\}}\cap Y$, we conclude that $3=\|\Pi\|<\|\overline{\Pi\cup\{o\}}\cap Y\|\le\|\overline{\Pi\cup\{o\}}\|\le 4$, by the rankedness of $X$. Since $\|\overline{\Pi\cup\{o\}}\cap Y\|=\|\overline{\Pi\cup\{o\}}\|=4$, the rankedness of $X$ ensures that $o\in \overline{\Pi\cup\{o\}}=\overline{\Pi\cup\{o\}}\cap Y\subseteq Y$, which contradicts the choice of the point $o$. This contradiction shows that $\overline{\Pi\cup\{o\}}\cap Y=\Pi$.
\end{proof}

Since $P\subseteq\overline{P\cup P'}\subseteq Y$ and $P\ne\overline{P\cup P'}$, the rankedness of $X$ implies that $3=\|P\|<\|\overline{P\cup P'}\|\le\|Y\|=4$ and hence $\|\overline{P\cup P'}\|=\|Y\|=4$. The rankedness of $X$ ensures that $\overline{P\cup P'}=Y\ne\overline{P\cup P'\cup\{o\}}$ and hence $\|P\cup P'\cup\{o\}\|=5$. Since $P\ne\overline{P\cup\{o\}}$, the rankedness of $X$ ensures that $\|P\cup\{o\}\|=\|P\|+1=4$.  
Since the liner $Y$ is weakly modular, the flat $\Pi\defeq\overline{P\cup\{o\}}\cap\overline{P'\cup\{o\}}$ has rank 
$$\|\Pi\|=\|P\cup\{o\}\|+\|P'\cup\{o\}\|-\|P\cup P'\cup\{o\}\|=4+4-5=3,$$
which means that the flat $\Pi$ is a plane in $X$. 

Claim~\ref{cl:intersectY} implies
$$\Pi\cap Y=(\overline{P\cup\{o\}}\cap\overline{P'\cup\{o\}})\cap Y=(\overline{P\cup\{o\}}\cap Y)\cap(\overline{P'\cup\{o\}}\cap Y)=P\cap P'=\varnothing.$$
Since $$Y\ne \Pi\cup Y\subseteq\overline{P\cup\{o\}}\cup Y\subseteq\overline{Y\cup\{o\}}\ne Y,$$the rankedness of $X$ ensures that $$4=\|Y\|<\|\Pi\cup Y\|\le\|Y\cup\{o\}\|\le \|Y\|+1=5$$ and hence $$\|\Pi\cup Y\|=\|Y\cup\{o\}\|=5.$$ 
By the rankedness of $X$, for every point $y\in Y$, the flat $\overline{\Pi\cup\{y\}}\ne\Pi$ has rank $\|\Pi\cup\{y\}\|=\|\Pi\|+1=4$ and by the weak modularity of $X$,  the flat
$$P_y\defeq \overline{\Pi\cup\{y\}}\cap Y$$has rank
$$\|P_y\|=\|\overline{\Pi\cup\{y\}}\|+\|Y\|-\|\overline{\Pi\cup\{y\}}\cup Y\|=4+4-\|\Pi\cup Y\|=8-5=3,$$
which means that $P_y$ is a plane in $Y$.

Observe that $P_y=\overline{\Pi\cup\{y\}}\cap Y\subseteq\overline{\Pi\cup P_y}\cap Y\subseteq\overline{\Pi\cup \overline{\Pi\cup\{y\}}}\cap Y=\overline{\Pi\cup\{y\}}\cap Y=P_y$ and hence $P_y=\overline{\Pi\cup P_y}\cap Y$. 
Then for every point $x\in P_y$, we have 
$$P_x=\overline{\Pi\cup\{x\}}\cap Y\subseteq \overline{\Pi\cup P_y}\cap Y=P_y.$$
Given any points $y,z\in Y$ with $P_y\cap P_z\ne\varnothing$, we can choose a point $z\in P_y\cap P_z$ and conclude that $P_y=P_x=P_z$. This shows that $\{P_y:y\in Y\}$ is the required family of parvise disjoint planes covering the flat $Y$.
\end{proof}

Next, we prove a higher-dimensional counterpart of Corollary~\ref{c:x2|x3}

\begin{corollary}\label{c:x3|x4} Let $X$ be a balanced weakly modular liner of rank $\|X\|\ge 5$. If $X$ is not projective, then $|X|_4=|X|_3\cdot\kappa$ for some cardinal $\kappa\ge |X|_2$.
\end{corollary}

\begin{proof} By Theorem~\ref{t:w-modular<=>}, the weakly modular liner $X$ is ranked. Since $\|X\|\ge 5$, the balanced liner $X$ contains flats of every rank $r\le 5$, so the cardinals $|X|_2$, $|X|_3$ and $|X|_4$ are well-defined. It is clear that $2\le|X|_2\le|X|_3\le|X|_4$. If the cardinal $|X|_4$ is infinite, then $|X|_4=|X|_3\cdot|X|_4=|X|_3\cdot\kappa$ for some cardinal $\kappa\ge |X|_2$. So, assume that $|X|_4$ is finite and so are the cardinals $|X|_2$ and $|X|_3$. Fix any flat $Y\subseteq X$ of rank $\|Y\|=4$. If $X$ is not projective, then the flat $Y$ contains two disjoint planes, by Lemma~\ref{l:disjoint-planes}. By Lemma~\ref{l:Y=union-of-planes}, the flat $Y$ is the union of a disjoint family $\F$ of planes, which implies $|X|_4=|Y|=|X|_3\cdot|\F|$. Proposition~\ref{p:cov-aff} ensures that $|\F|\ge|X|_2$.
\end{proof}

Now we are able to prove mains results of this section, which are due to \index[person]{Doyen}Doyen\footnote{{\bf Jean Doyen}, a professor of Mathematics at Universit\'e Libre de Bruxelles, in 1970 defended his Ph.D. Thesis ``Sur les syst\`mes de Steiner'' advised by Paul Libois.}  and \index[person]{Hubaut}Hubaut\footnote{{\bf Xavier Hubaut} (1936--2020), worked in Universit\'e Libre de Bruxelles. In 1964 he defened the Ph.D. Dissertation ``Trois questions de math\'ematique conforme'', advised by Paul Libois.} \cite{DH1971}.

\begin{theorem}[Doyen--Hubaut, 1971]\label{t:DH1} Let $X$ be a balanced weakly modular liner of rank $\|X\|\ge 5$. If $|X|_3<\w$, then $X$ is $p$-parallel for some $p\in\{0,1\}$.
\end{theorem}

\begin{proof} By Proposition~\ref{p:23-balance=>k-parallel}, the balanced liner $X$ is $p$-parallel for a unique cardinal $p$ such that $|X|_3=1+(|X|_2-1)(|X_2|+p)$. We have to prove that $p\le 1$. To derive a contradiction, assume that $p>1$. This implies that $X$ is not projective. For every $n\le\|X\|$, let $x_n\defeq |X|_n$, and consider the cardinal $\ell\defeq x_2+p-1$. Then 
$$x_3=1+(x_2-1)(x_2+p)=1+(x_2-1)(1+\ell)=x_2(1+\ell)-\ell.$$
Corollary~\ref{c:x2|x3} ensures that the number $x_2$ divides the number $x_3=x_2(1+\ell)-\ell$, which implies that $\ell=x_2a$ for some positive integer number $a$.

By Theorem~\ref{t:wr+k-parallel=>n-balanced}, $$x_4=1+(x_2-1)(1+\ell+\ell^2)=x_2(1+\ell+\ell^2)-\ell(1+\ell)=x_2(1+\ell+\ell^2-a(1+\ell)).$$By Corollary~\ref{c:x3|x4}, the number $x_3=x_2(1+\ell-a)$ divides the number $x_4=x_2(1+\ell+\ell^2-a(1+\ell))$ and hence $1+\ell-a$ divides $1+\ell+\ell^2-a(1+\ell)$. Since
$$1+\ell+\ell^2-a(1+\ell)=(1+\ell-a)(1+\ell)-\ell,$$the number $1+\ell-a=1+a(x_2-1)$ divides $\ell=ax_2$ and hence $ax_2=(1+a(x_2-1))b$ for some positive integer number $b$.
The latter equation is equivalent to $ax_2(b-1)=(a-1)b$ and to $x_2\frac{b-1}b=\frac{a-1}a$.
If $a>1$, then $b>1$ and 
$$\frac32\le x_2 \frac12\le x_2\frac{b-1}b=\frac{a-1}a<1,$$which is a contradiction showing that $a=1$ and hence $x_2+p-1=\ell=ax_2=x_2$ and $p=1$.
\end{proof}

\begin{corollary}\label{c:wm5=>regular} Every balanced weakly regular liner $X$ with $\|X\|\ge 5$ and $|X|_3<\w$ is regular and $p$-parallel for some $p\in\{0,1\}$.
\end{corollary}

\begin{proof} It follows from $|X|_3<\w$ that the liner $X$ is $3$-ranked. By Theorem~\ref{t:ranked<=>EP} and Proposition~\ref{p:wr-ex<=>3ex}, the $3$-ranked weakly regular liner $X$ is ranked. By Theorem~\ref{t:w-modular<=>}, the ranked weakly regular liner $X$ is weakly modular. By Theorem~\ref{t:DH1}, the liner $X$ is $p$-parallel for some $p\in\{0,1\}$. If $p=0$, then $X$ is projective and strongly regular, by Theorem~\ref{t:projective<=>}. If $p=1$, then the $1$-parallel liner $X$ is Playfair (by Definitions~\ref{d:PPBL} and \ref{d:k-PPB}) and by Theorem~\ref{t:Playfair<=>}, the Playfair liner $X$ is $3$-regular, affine, and $3$-long. If $|X|_2\ge 4$, then the affine liner $X$ is regular, by Theorem~\ref{t:4-long-affine}. If $|X|_2=3$, then the weakly regular Steiner Playfair liner $X$ is regular, by Theorems~\ref{t:Hall<=>Playfair+Steiner} and \ref{t:Hall<=>regular}. 
\end{proof}

\begin{theorem}[Doyen--Hubaut, 1971]\label{t:DH2} Let $X$ be a balanced weakly modular liner of rank $\|X\|= 4$. If $|X|_3<\w$, then $X$ is $p$-parallel for some $p\in\{0,1,|X|_2^2-|X|_2+1,|X|^3_2+1\}$.
\end{theorem}

\begin{proof} By Proposition~\ref{p:23-balance=>k-parallel}, the balanced liner $X$ with $|X|_2<\w$ is $p$-parallel for some finite cardinal $p$. Assuming that $p$ is positive, we shall prove that $p\in\{1,|X|_2^2-|X|_2+1,|X|_2^3+1\}$. It will be convenient to denote the cardinals $|X|_n$ by $x_n$, and put $\ell\defeq x_2+p-1$. Theorem~\ref{t:wr+k-parallel=>n-balanced} ensures that 
$$x_3=1+(x_2-1)(1+\ell)=x_2(1+\ell)-\ell\quad\mbox{and}\quad x_4=1+(x_2-1)(1+\ell+\ell^2).$$
By Corollary~\ref{c:x2|x3}, the number $x_2$ divides the number $x_3=x_2(1+\ell)-\ell$ and hence divides the number $\ell$. Then $\ell=ax_2$ for some positive integer number $a$. 

Observe that the number of planes in $X$ that contain some fixed point of $X$ equals
$$\frac{(x_4-1)(x_4-x_2)}{(x_3-1)(x_3-x_2)}=\frac{(x_2-1)(1+\ell+\ell^2)(x_2-1)(\ell+\ell^2)}{(x_2-1)(1+\ell)(x_2-1)\ell}=1+\ell+\ell^2.$$
Then the total number of planes in $X$ equals 
$\frac{x_4(1+\ell+\ell^2)}{x_3}$, which implies that $x_3$ divides the number $x_4(1+\ell+\ell^2)$ and also the number $(x_4-x_3)(1+\ell+\ell^2)=(x_2-1)\ell^2(1+\ell+\ell^2)$.

Since $\ell=\frac{x_3-x_2}{x_2-1}$, the number $x_3$ divides the number
$$
(x_2-1)^4\ell^2(1+\ell+\ell^2)=(x_3-x_2)^2((x_2-1)^2+(x_3-x_2)(x_2-1)+(x_3-x_2)^2)$$ and hence $x_3=x_2(1+\ell)-\ell=x_2(1+\ell-a)$ divides the number
$$x_2^2((x_2-1)^2-x_2(x_2-1)+x_2^2)=x_2^2(x_2^2-x_2+1).$$
Then the number $1+\ell-a=1+a(x_2-1)$ divides the number 
$x_2(x_2^2-x_2+1)=(x_2-1)(x_2^2+1)+1$ and hence
$$(1+a(x_2-1))(b(x_2-1)+1)=(x_2-1)(x_2^2+1)+1$$
for some non-negative integer number $b$. After cancellations, we obtain the equality
$$(a+b)+ab(x_2-1)=x_2^2+1=2+(x_2+1)(x_2-1).$$
It follows that $a+b=2+(x_2-1)c$ for some positive integer number $c$. 
Then $(x_2-1)(c+ab)=-2+a+b+ab(x_2-1)=(x_2+1)(x_2-1)$ and hence $c+ab=x_2+1$.
If $b=0$, then $c=c+ab=x_2+1$ and $a=a+b=2+(x_2-1)(x_2+1)=x_2^2+1$.
If $b$ is positive, then the equality $a+b=2+(x_2-1)c$ implies $x_2\ge x_2+1-c=ab\ge 1+(x_2-1)c$ and $c=1$. In this case $a+b=x_2+1$ and $ab=x_2$ is possible only if $(a,b)\in\{(1,x_2),(x_2,1)\}$. Therefore, $(a,b)\in\{(1,x_2),(x_2,1),(x_2^2+1,0)\}$ and hence
$x_2+p-1=\ell=ax_2\in\{x_2,x_2^2,x_2^3+x_2\}$, and finally,
$$p\in\{1,x_2^2-x_2+1,x_2^3+1\}.$$
\end{proof}

\begin{problem} Is every finite balanced weakly regular liner of rank $\ge 4$ proaffine?
\end{problem}

\begin{theorem}\label{t:balanced<=>ranked} For a balanced liner $X$ with $\|X\|\ge 4$ and $|X|_3<\w$, the following conditions  are equivalent:
\begin{enumerate}
\item $X$ is regular;
\item $X$ is $p$-parallel for some $p\in\{0,1\}$;
\item $X$ is projective or Playfair;
\item $X$ is Proclus;
\item $X$ is $3$-proregular and proaffine.
\end{enumerate}
If $\|X\|\ge5$, then the conditions \textup{(1)--(5)} are equivalent to:
\begin{enumerate}
\item[\textup{(6)}] $X$ is weakly regular;
\item[\textup{(7)}] $X$ is weakly modular.
\end{enumerate}
\end{theorem} 

\begin{proof} First we prove the implications $(1)\Ra(2)\Ra(3)\Ra(4)\Ra(5)\Ra(1)$.
\smallskip

$(1)\Ra(2)$ Assume that the liner $X$ is regular. By Corollary~\ref{c:proregular=>ranked} and Theorem~\ref{t:w-modular<=>}, the regular liner $X$ is weakly modular.  By Theorems~\ref{t:DH1} and \ref{t:DH2}, the balanced weakly modular liner $X$ of rank $\|X\|\ge4$ is $p$-parallel for some $p\in\{0,1,|X|_2^2-|X|_2+1,|X|_2^3+1\}$. By Theorem~\ref{t:2-balance+k-parallel=>3-balance}, $|X|_3=1+(|X|_2-1)(|X|_2+p)$.  Fix any plane $P\subseteq X$ and two concurrent lines $A,B$ in $P$. By the regularity of $X$, for every point $x\in P\setminus(A\cup B)$ there exist points $a\in A$ and $b\in B$ such that $x\in\Aline ab$. It follows from $x\notin A\cup B$ that $a\notin B$ and $b\notin A$. Then $$P=A\cup B\cup\bigcup_{a\in A\setminus B}\bigcup_{b\in B\setminus A}(\Aline ab\setminus\{a,b\})$$ and hence
$$
\begin{aligned}
1+(|X|_2-1)(|X|_2+p)&=|X|_3=|P|\le |A|+|B|-|A\cap B|+|A\setminus B|\cdot|B\setminus A|\cdot(|X|_2-2)\\
&=2|X|_2-1+(|X|_2-1)^2(|X|_2-2)=1+(|X|_2-1)(2-3|X|_2+|X|_2^2),
\end{aligned}
$$
which implies $$p\le 2-4|X|_2+|X|_2^2<1-|X|_2+|X|_2^2<|X|_2^3+1.$$
Therefore, $p\in\{0,1\}$.
\smallskip

$(2)\Ra(3)$ If $X$ is $p$-parallel for some $p\in\{0,1\}$, then $X$ is projective or Playfair by 
Definitions~\ref{d:PPBL} and \ref{d:k-PPB}.
\smallskip

$(3)\Ra(4)$ If $X$ is projective or Playfair, then $X$ is Proclus, by Definition~\ref{d:PPBL}.
\smallskip

$(4)\Ra(5)$ If $X$ is Proclus, then $X$ is $3$-proregular and proaffine, by Theorem~\ref{t:Proclus<=>}.
\smallskip

$(5)\Ra(1)$ If  $X$ is $3$-proregular and proaffine, then $X$ is Proclus, by Theorem~\ref{t:Proclus<=>}. Since $|X|_3<\w$ and $\|X\|\ge 3$, by Proposition~\ref{p:23-balance=>k-parallel}, the balanced liner $X$ is $p$-parallel for some finite cardinal $p$. Since $X$ is Proclus, $p\le 1$. If $p=0$, then the $0$-parallel liner $X$ is (strongly) regular, by 
Theorem~\ref{t:projective<=>}. If $p=1$, then the $1$-parallel liner $X$ is Playfair, by Definition~\ref{d:PPBL}. By Theorem~\ref{t:Playfair<=>}, the Playfair liner $X$ is $3$-long, $3$-regular and affine. If $|X|_3\ge 4$, then the $4$-long affine liner $X$ is regular, by 
Theorem~\ref{t:4-long-affine}. If $|X|_3=3$, then $X$ is a Hall liner, by Theorem~\ref{t:Hall<=>Playfair+Steiner}.
By Theorem~\ref{t:Hall<=>regular}, the balanced Hall liner $X$ is regular.
\smallskip

$(1)\Ra(6)$ The implication $(1)\Ra(6)$ follows from Proposition~\ref{p:sr=>r=>wr}.
\smallskip

$(6)\Ra(7)$ Being $3$-balanced with $|X|_3<\w$, the liner $X$ is $3$-ranked. If $X$ is weakly regular, then $X$ is ranked by Proposition~\ref{p:wr-ex<=>3ex} and Theorem~\ref{t:ranked<=>EP}. By Theorem~\ref{t:w-modular<=>}, the ranked weakly regular liner $X$ is weakly modular.
\smallskip

$(7)\Ra(1)$ If $X$ is weakly modular and $\|X\|\ge 5$, then $X$ is regular, by Corollary~\ref{c:wm5=>regular}.
\end{proof}

\begin{corollary} For a finite liner $X$ of rank $\|X\|\ge 5$, the following conditions are equivalent:
\begin{enumerate}
\item $X$ is regular and Playfair;
\item $X$ is balanced, weakly regular and not strongly regular.
\item $X$ is balanced, weakly modular and not modular.
\end{enumerate}
\end{corollary}

\begin{proof} $(1)\Ra(2)$ Assume that $X$ is regular and Playfair. By Proposition~\ref{p:sr=>r=>wr} and Corollary~\ref{c:proregular=>ranked}, the regular liner $X$ is weakly regular and ranked. By Theorem~\ref{t:Playfair<=>}, the Playfair liner $X$ is affine. By Theorem~\ref{t:affine=>Avogadro}, the affine liner $X$ is $2$-balanced. By Definition~\ref{d:PPBL}, the Playfair liner $X$ is $1$-parallel. Since $\|X\|\ge 3$, the $1$-parallel liner $X$ is not $0$-parallel and hence $X$ is not strongly regular, by Theorem~\ref{t:projective<=>}. By Theorem~\ref{t:2-balance+k-parallel=>3-balance}, the $2$-balanced $1$-parallel liner $X$ is $3$-balanced. By Theorem~\ref{t:wr+k-parallel=>n-balanced}, the $2$-balanced $3$-balanced ranked weakly regular finite liner $X$ is balanced.
\smallskip

$(2)\Ra(3)$ Assume that $X$ is balanced, weakly regular and not strongly regular. By Theorem~\ref{t:modular<=>}, the liner $X$ is not modular. Being finite and balanced, the liner $X$ is ranked. By Theorem~\ref{t:w-modular<=>}, the ranked weakly regular liner $X$ is weakly modular.
\smallskip

$(3)\Ra(1)$ Assume that the liner $X$ is balanced, weakly modular and not modular. By Theorems~\ref{t:modular<=>} and \ref{t:projective<=>}, the liner $X$ is not $0$-parallel. By Theorem~\ref{t:DH1}, the balanced finite liner $X$ of rank $\|X\|\ge 5$ is $p$-parallel for some $p\in\{0,1\}$. Since $X$ is not $0$-parallel, $X$ is $1$-parallel and hence $X$ is Playfair, see Definitions~\ref{d:PPBL} and \ref{d:k-PPB}. By Theorem~\ref{t:balanced<=>ranked}, the balanced Playfair liner $X$ is regular.
\end{proof}

\begin{problem} Is every finite balanced finite affine liner $X$ of rank $4$ regular?
\end{problem}
\index[person]{Hetman}

\begin{Exercise}[Ivan Hetman] Show that every balanced regular liner $X$ with $|X|_2\le 4$ is proaffine.
\end{Exercise}

\begin{examples}[Ivan Hetman\footnote{{\tt mathoverflow.net/a/459557/61536}}]\label{ex:Denniston} There exists a hyperbolic (and hence not proaffine) regular balanced liner $X$ with $|X|_2=8$ and $|X|=|X|_3=120$.
\end{examples}

\begin{proof} Let $F$ be a (unique) 16-element field and $r\in F$ be an element such that $r^4+r+1=0$. Consider the $8$-element subgroup $$H\defeq\{0,r^3,r^2,r^2+r^3,r,r+r^3,r+r^2,r+r^2+r^3\}$$of the additive group of the field $F$. It can be shown that the subliner $$X\defeq\{(x,y)\in F\times F:x^2+rxy+y^2\in H\}$$of the affine plane $F\times F$ is balanced, regular and has $|X|=|X|_3=120$ and $|X|_2=8$. 
Computer calculations show that for every points $o,x,y\in X$ and $p\in\Aline oy\setminus (\Aline ox\cup\Aline oy)$ the set $\{u\in\Aline oy:\Aline up\cap\Aline ox=\varnothing\}$ has cardinality between $2$ and $6$, which implies that the liner $X$ is hyperbolic and hence not proaffine. 
\end{proof}

\begin{remark} The liner mentioned in Example~\ref{ex:Denniston} is one of maximal arcs, discovered by Denniston \cite{Denniston1969}.
\end{remark}

\chapter{Parallelity in liners}\label{s:parallelity-in-matroids}

\section{Subparallelity}\label{ss:subparallelity-in-matroids}


\begin{definition}\label{d:subparallel}  Given two flats $A$ and $B$ in a liner, we write\index[note]{$\subparallel$} $A\subparallel B$ and say that $A$ is \index{subparallel flats}\index{flat!subparallel to}\defterm{subparallel to} $B$ if $\forall a\in A\;\;(A\subseteq \overline{\{a\}\cup B})$.
\end{definition}

\begin{lemma}\label{l:subparallel+intersect=>subset} For two intersecting flats $A,B$ in a liner, we have $A\subparallel B$ if and only if $A\subseteq B$.
\end{lemma}

\begin{proof} Take any point $o\in A\cap B$. If $A\subparallel B$, then $A\subseteq \overline{B\cup\{o\}}=B$. If $A\subseteq B$, then for every $a\in A$ we have $A\subseteq B=\overline{\{a\}\cup B}$, which means that $A\subparallel B$.
\end{proof}


\begin{exercise} Show that a flat $A$ in a liner $X$ is subparallel to a flat $B\subseteq X$ if and only if for every $a,b\in A$ the flat $\Aline ab$ is subparallel to $B$.
\end{exercise}

\begin{theorem}\label{t:subparallel-char} Let $\kappa$ be a cardinal and $X$ be a $\kappa$-ranked liner. For two flats $A,B\subseteq X$ with $\|B\|<\kappa$, the following conditions are equivalent:
\begin{enumerate}
\item $A\subparallel B$;
\item either $A\subseteq B$ or  $A\subseteq \overline{\{a\}\cup B}\setminus B$ for some $a\in A$;
\item either $A\subseteq B$ or  $A\cap B=\varnothing$ and $\|A\|_B=1$.
\end{enumerate}
If $\|B\|<\w$, then the conditions \textup{(1)--(3)} are equivalent to the condition
\begin{itemize}
\item[(4)] either $A\subseteq B$ or $A\cap B=\varnothing$ and $\|A\cup B\|=\|B\|+1$.
\end{itemize}
\end{theorem}

\begin{proof} $(1)\Ra(2)$ Assume that $A\subparallel B$. If $A=\varnothing$, then $A=\varnothing\subseteq\overline B=B$ and we are done. So, assume that $A\ne\varnothing$. If $A\cap B$ contains some point $x$, then $A\subseteq \overline{\{x\}\cup B}\subseteq \overline{B}=B$ and we are done. So, assume  $A\ne\varnothing=A\cap B$ and choose any point $a\in A$. It follows from $A\subparallel B$ and $A\cap B=\varnothing$ that $A\subseteq\overline{\{a\}\cup B}\setminus B$.
\smallskip

The implication $(2)\Ra(3)$ follows from the definition of the rank $\|A\|_B$.
\smallskip

$(3)\Ra(1)$ Assume that the condition (3) is satisfied. If $A\subseteq B$, then for every $x\in A$ we have $A\subseteq B\subseteq\overline{\{x\}\cup B}$, witnessing that $A\subparallel B$. So, assume that $A\not\subseteq B$. The condition (3) ensures that $A\cap B=\varnothing\ne A$ and $\|A\|_B=1$. Then $A\subseteq\overline{B\cup\{x\}}$ for some $x\in X\setminus B$. By Theorem~\ref{t:ranked<=>EP}, the $\kappa$-ranked liner $X$ has the $\kappa$-Exchange Property. Since $\|B\|<\kappa$, for every $a\in A$ we have $a\in \overline{B\cup\{x\}}\setminus B$. The $\kappa$-Exchange Property of $X$ ensures that $\overline{B\cup\{x\}}=\overline{B\cup\{a\}}$ and hence $A\subseteq \overline{\{x\}\cup B}=\overline{\{a\}\cup B}$, witnessing that $A\subparallel B$.
\smallskip

$(3)\Rightarrow(4)$ Assume that either $A\subseteq B$ or $A\cap B=\emptyset$ and $\|A\|_B=1$. Then $\|A\cup B\|\le\|B\|+\|A\|_B=\|B\|+1\le\kappa$. Since $X$ has the $\kappa$-Exchange Property, $\|A\cup B\|=\|B\|+\|A\|_B=\|B\|+1$ by Corollary~\ref{c:rank+}.

$(4)\Ra(3)$ Assume that $\|B\|<\w$ and the condition (4) is satisfied. We have to prove that $\|A\|_B=1$ if  $A\cap B=\emptyset$ and $\|A\cup B\|=\|B\|+1$. Since $\|B\|<\kappa$, we obtain $\|A\cup B\|=\|B\|+1\le\kappa$ and hence $\|A\cup B\|=\|B\|+\|A\|_B$, by Corollary~\ref{c:rank+}. Then $\|B\|+\|A\|_B=\|A\cup B\|=\|B\|+1$ and since the cardinal $\|B\|$ is finite, $\|A\|_B=1$. 
\end{proof}

\begin{corollary}\label{c:subparallel} Let $\kappa$ be a finite cardinal and $A,B$ are two flats in a $\kappa$-ranked liner $X$. If $\|A\|=\|B\|<\kappa$, then $A\subparallel B\;\Leftrightarrow\; B\subparallel A$.
\end{corollary}

\begin{proof} Assume that $\|A\|=\|B\|$ and $A\subparallel B$. By Theorem~\ref{t:subparallel-char}, either $A\subseteq B$ or else $A\cap B=\varnothing$ and $\|A\cup B\|=1+\|B\|$. If $A\subseteq B$, then $A=B$ (as $X$ is $\kappa$-ranked) and $B\subparallel A$. If $A\not\subseteq B$, then $A\cap B=\varnothing$ and  $\|A\cup B\|=1+\|B\|=1+\|A\|$. Applying Theorem~\ref{t:subparallel-char}, we conclude that $B\subparallel A$. By analogy we can prove that $B\subparallel A$ implies $A\subparallel B$.
\end{proof}

\begin{exercise} Find an example of a liner $X$ containing two lines $L,\Lambda$ such that $L\subparallel \Lambda$ but not $\Lambda\subparallel L$.
\smallskip

\noindent{\em Hint:} Look at the lines $L_2$ and $L_3$ in Example~\ref{ex:Tao}.
\end{exercise} 

\begin{corollary}\label{c:para+intersect=>coincide} Let $A,B$ be two flats in a $\|B\|$-ranked liner $X$ such that $A\subparallel B$ and $\|A\|=\|B\|<\w$. If $A\cap B\ne\varnothing$, then $A=B$.
\end{corollary}

\begin{proof} By Lemma~\ref{l:subparallel+intersect=>subset}, $A\subseteq B$. Since $\|A\|=\|B\|<\w$ and the liner $X$ is $\|B\|$-ranked, the inclusion $A\subseteq B$ implies the equality $A=B$.
\end{proof} 
 
\begin{proposition}\label{p:subparallel=>dim} Let $\kappa$ be a cardinal and $A,B$ be two flats in a $\kappa$-ranked liner $X$. If $\|B\|<\kappa$, $A\subparallel B$ and $B\ne\varnothing$, then $\|A\|\le \|B\|$.
\end{proposition}

\begin{proof} By Theorem~\ref{t:subparallel-char}, either $A\subseteq B$ and hence $\|A\|\le\|B\|$ or else $A\subseteq \overline{\{a\}\cup B}\setminus B$ for some $a\in A$, which implies $\|A\|_B\le 1$. In the second case, by Corollary~\ref{c:rank+}, $\|A\|\le\|A\cup B\|\le\|A\|_B+\|B\|=1+\|B\|$. If the cardinal $\|B\|$ is infinite, then $\|A\|\le\|A\cup B\|\le 1+\|B\|=\|B\|$ and we are done. 

So, assume that the rank $\|B\|$ is finite. Assuming that $\|A\|>\|B\|$ and taking into account that $\|A\|\le\|\overline{\{a\}\cup B}\|\le 1+\|B\|$, we conclude that $\|A\|=\|\{a\}\cup B\|=1+\|B\|$. Since the liner $X$ is ranked, $A= \overline{\{a\}\cup B}$, which contradicts the assumption $A\cap B=\varnothing\ne B$.
\end{proof}


\begin{proposition}\label{p:subparallel-function} Let $F:X\to Y$ be a flat function between liners and $A,B$ be two flats in $X$. If  $A\subparallel B$, then $F[A]\subparallel F[B]$.
\end{proposition}

\begin{proof} Assume that $A\subparallel B$. To show that $F[A]\subparallel F[B]$, take any point  $y\in F[A]$ and find a point $x\in A$ such that $y=F(x)$. It follows from $A\subparallel B$ that $A\subseteq\overline{B\cup\{x\}}$. Since the function $F$ is flat, the set $F^{-1}[\overline{F[B]\cup\{y\}}]$ is flat and contains the set $B\cup\{x\}$. Taking into account that $\overline{B\cup\{x\}}$ is the smallest flat containing $B\cup\{x\}$, we conclude that $A\subseteq \overline{B\cup\{x\}}\subseteq F^{-1}[\overline{F[B]\cup\{y\}}]$ and hence $F[A]\subseteq F[F^{-1}[\overline{F[B]\cup\{y\}}]]$. Since the function $F$ is flat, the image $F[X]$ is a flat in $Y$ and hence $\overline{F[B]\cup\{y\}}\subseteq F[X]$. Then $F[A]\subseteq F^{-1}[F[\overline{F[B]\cup\{y\}}]=\overline{F[B]\cup\{y\}}$, witnessing that $F[A]\subparallel F[B]$.
\end{proof}

\section{Parallelity}\label{ss:parallelity-in-matroids}

\begin{definition}\label{d:parallel} Given two flats $A$ and $B$ in a liner $X$, we write\index[note]{$\parallel$} $A{\parallel}B$ and say that $A$ and $B$ are \index{parallel flats}\defterm{parallel} if $A$ is subparallel to $B$  and $B$ is subparallel to $A$. Therefore,
$$A{\parallel}B\;\Leftrightarrow\;(A\subparallel B\;\wedge\; B\subparallel A).$$For two flats $A,B$, the negation of $A\parallel B$ is denoted by $A\nparallel B$.
\end{definition}

Theorem~\ref{t:subparallel-char} implies the following characterization.

\begin{theorem}\label{t:parallel-char} Let $\kappa$ be a cardinal, $X$ be a $\kappa$-ranked liner and $A,B$ be two nonempty flats in $X$. If $\max\{\|A\|,\|B\|\}<\kappa$, then the following conditions are equivalent:
\begin{enumerate}
\item $A\parallel B$;
\item either $A=B$ or $A\subseteq \overline{\{a\}\cup B}\setminus B$ and $B\subseteq\overline{\{b\}\cup A}\setminus A$ for some points $a\in A$ and $b\in B$;
\item either $A=B$ or else $A\cap B=\varnothing$ and $\|A\|_B=1=\|B\|_A$.
\end{enumerate}
If $\min\{\|A\|,\|B\|\}<\w$, then the conditions \textup{(1)--(3)} are equivalent to the condition
\begin{enumerate}
\item[(4)] either $A=B$ or else $A\cap B=\varnothing$ and $1+\|A\|=\|A\cup B\|=1+\|B\|$.
\end{enumerate}
\end{theorem}

\begin{corollary}\label{c:singletons||} Any two singletons in a liner are parallel.
\end{corollary}

Proposition~\ref{p:subparallel=>dim} implies

\begin{corollary}\label{c:parallel} Let $\kappa$ be a cardinal, $X$ be a $\kappa$-ranked liner and $A,B$ be two flats in $X$. If $0<\min\{\|A\|,\|B\|\}<\kappa$ and $A\parallel B$, then $\|A\|=\|B\|$.
\end{corollary}

Lemma~\ref{l:subparallel+intersect=>subset} implies the following simple (but helpful) proposition.

\begin{proposition}\label{p:para+intersect=>coincide} Two intersecting flats in a liner are parallel if and only if they are equal.
\end{proposition}

Theorem~\ref{t:parallel-char} implies the following characterization of the parallelity of lines in $3$-ranked liner, which often is taken as the definition of the parallelity of lines.

\begin{corollary}\label{c:parallel-lines<=>} Two lines $L,\Lambda$ is a $3$-ranked liner are parallel if and only if $\|L\cup\Lambda\|\le 3$ and either $L=\Lambda$ or $L\cap\Lambda=\varnothing$.
\end{corollary}

\begin{corollary}\label{c:parallel-in-Proclus} Let $A,B,C$ be three lines in a Proclus liner $X$. If $A\parallel B$, $B\parallel C$ and $\|A\cup C\|\le 3$, then $A\parallel C$.
\end{corollary}

\begin{proof} By Theorem~\ref{t:Proclus<=>}, the Proclus liner $X$ is $3$-proregular and by Proposition~\ref{p:k-regular<=>2ex}, the $3$-proregular liner $X$ is $3$-ranked. 

If $A\cap B\ne\varnothing$, then $A\parallel B$ implies $A=B$ and hence $B\parallel C$ implies $A\parallel C$. By analogy we can show that $A\parallel C$ if $B\cap C\ne\varnothing$.

So, we assume that $A\cap B=\varnothing=B\cap C$.  If the lines $A$ and $C$ have a common point $x$, then $A,C$ are two lines in the plane $\overline{B\cup\{x\}}$ that contain the point $x\in A\cap C\setminus B$ and are disjoint with the line $B$. Since $X$ is Proclus, $A=C$ and hence $A\parallel C$.

So, we assume that the lines $A,B$ are disjoint. Since $\|A\cup C\|\le 3$, we can apply Corollary~\ref{c:parallel-lines<=>} and conclude that $A\parallel C$.
\end{proof}

\begin{proposition}\label{p:parallel=>3-long}  If a $2$-balanced liner $X$ contains two disjoint parallel lines, then $X$ is $3$-long.
\end{proposition}

\begin{proof} If a $2$-balanced liner $X$ is not $3$-long, then all lines in $X$ contain exactly 2 points and hence any subset of $X$ flat. Then for any disjoint lines $L,\Lambda$ and any point $x\in\Lambda$, we obtain $\overline{L\cup\{x\}}=L\cup\{x\}$ and hence $\Lambda\not\subseteq \overline{L\cup\{x\}}$, witnessing that $L\nparallel\Lambda$.
\end{proof} 

Proposition~\ref{p:subparallel-function} implies the following corollary on preservation of parallel flats by flat functions.

\begin{corollary}\label{c:parallel-function} Let $F:X\to Y$ be a flat function between liners and $A,B$ be two flats in $X$. If  $A\parallel B$, then $F[A]\parallel F[B]$.
\end{corollary}

\section{Parallelity in modular liners}

\begin{proposition} For two flats $A,B$ in a modular liner $X$, the following conditions are equivalent:
\begin{enumerate}
\item $A\subparallel B$;
\item $A\subseteq B$ or $\|A\|=1$.
\item $A\subseteq B$ or $\|A\|\le 1$.
\end{enumerate}
\end{proposition}

\begin{proof} By Theorem~\ref{t:modular<=>}, the modular liner $X$ is projective and strongly regular. Since every regular liner is proregular, the (strongly) regular liner $X$ is proregular. By Corollary~\ref{c:proregular=>ranked}, the proregular liner $X$ is ranked.
\smallskip

$(1)\Ra(2)$ Assume that $A\subparallel B$. By Theorem~\ref{t:subparallel-char}, either $A\subseteq B$ or $A\cap B=\varnothing$ and $\|A\|_B=1$.  If $A\not\subseteq B$, then $A\ne\varnothing$ and hence $\|A\|_B=1\le\|A\|$. Assuming that $1<\|A\|$, we conclude that $\|A\|_B<\|A\|$. Since the liner $X$ is projective, we can apply Proposition~\ref{p:projective<} and conclude that the intersection $A\cap B$ contains some point $a\in A\cap B$. It follows from $A\subparallel B$ that $A\subseteq\overline{\{a\}\cup B\}}=\overline B=B$, which contradicts our assumption. This contradiction shows that $A\not\subseteq B$ implies $\|A\|=1$.
\vskip5pt

The implication $(2)\Ra(3)$ is trivial.
\smallskip

$(3)\Ra(1)$ Assume that $A\subseteq B$ or $\|A\|\le 1$. If $A\subseteq B$, then $A\subparallel B$ by Theorem~\ref{t:subparallel-char}. So, assume that $A\not\subseteq B$. In this case  there exists a point $a\in A\setminus B$ and the condition (2) implies that $\|A\|\le 1$. It follows from $\|A\|\le 1$ that $A=\{a\}$ and hence $A=\{a\}\subseteq\overline{\{a\}\cup B}\setminus B$. Applying Theorem~\ref{t:subparallel-char}, we conclude that $A\subparallel B$.
\end{proof}

\begin{corollary}\label{c:parallel-in-modular} For two flats $A,B$ in a modular liner $X$, the following conditions are equivalent:
\begin{enumerate}
\item $A\parallel B$;
\item $A=B$ or $\max\{\|A\|,\|B\|\}=1$;
\item $A=B$ or $\max\{\|A\|,\|B\|\}\le 1$.
\end{enumerate}
\end{corollary}

\section{Parallelity in weakly modular liners}

\begin{proposition} Let $A,B,C$ be flats in a weakly modular liner $X$. 
 If $A\subparallel B$ and $B\cap C\ne\varnothing$, then $A\cap C\subparallel B\cap C$.
 \end{proposition}
 
\begin{proof}  If $A\subparallel B$, then for every $a\in A\cap C $, by Theorem~\ref{t:modular}, $$A\cap C\subseteq\overline{\{a\}\cup B}\cap C=\overline{\{a\}\cup(B\cap C)},$$
witnessing that $A\cap C\subparallel B\cap C$.
\end{proof}

\begin{corollary}\label{c:paraintersect} Let $A,B,C$ be flats in a weakly modular liner $X$. 
 If $A\parallel B$ and $A\cap C\ne\varnothing\ne B\cap C$, then $A\cap C\parallel B\cap C$.
\end{corollary}
 
\begin{proposition}\label{p:subparallel-char4} For two disjoint flats $A,B$ in a weakly modular liner $X$, the following conditions are equivalent:
\begin{enumerate}
\item $A\subparallel B$;
\item for every $b\in B$, the flats $A$ and $B\cap\overline{\{b\}\cup A}$ are parallel;
\item there exists a flat $C\subseteq B$ such that $A\parallel C$;
\item $\forall a,x\in A\;\forall b\in B\;\exists y\in B\;\;(\Aline ax\parallel \Aline by)$.
\end{enumerate}
\end{proposition}

\begin{proof} We shall prove the implications $(1)\Ra (4)\Ra(2)\Ra(3)\Ra(1)$. By Theorems~\ref{t:w-modular<=>} and \ref{t:ranked<=>EP}, the weakly modular liner $X$ is ranked and has the Exchange Property.
\smallskip

$(1)\Ra(4)$ Assume that $A\subparallel B$. Given any point $a,x\in A$ and $b\in B$, we should find a point $y\in B$ such that $\Aline ax\parallel \Aline by$. If $x=a$, then for $y\defeq b$ we have $\Aline ax=\{a\}\parallel \{b\}=\Aline by$.
So, assume that $x\ne a$. Then $x\notin \Aline ab$ and hence $\|\{a,b,x\}\|=3$. 

Since $x\in A\subseteq\overline{\{a\}\cup B}$, by Proposition~\ref{p:aff-finitary}, there exists a finite set $F\subseteq B$ such that $b\in F$ and $x\in\overline{\{a\}\cup F}$. It follows from $x\in A=A\setminus B$ that $a\notin \overline F$ and hence $\|\{a\}\cup F\|=1+\|F\|$, by Proposition~\ref{p:add-point-to-independent} and Theorem~\ref{t:Max=codim}. The weak modularity of the liner $X$ guarantees that $$
\begin{aligned}
3+\|F\|&=\|\overline{\{a,b,x\}}\|+\|\overline F\|=\|\overline{\{a,b,x\}}\cup \overline F\|+\|\overline{\{a,b,x\}}\cap\overline{F}\|\\
&=\|\{a\}\cup F\|+\|\overline{\{a,b,x\}}\cap\overline{F}\|=1+\|F\|+\|\overline{\{a,b,x\}}\cap\overline F\|,
\end{aligned}
$$
and hence $\|\overline{\{a,b,x\}}\cap \overline{F}\|=2$. So, there exists a point $y\in \overline{\{a,b,x\}}\cap\overline{F}\setminus\{b\}$. 
Since 
$$\|\Aline ax\cup\Aline by\|=\|\overline{\{a,x,b,y\}}\|=\|\overline{\{a,b,x\}}\|=3=1+\|\Aline ax\|=1+\|\Aline by\|,$$ we can apply Corollary~\ref{c:parallel} and conclude that $\Aline ax\parallel \Aline by$. 
\smallskip

$(4)\Ra(2)$ Assume that (4) holds. To prove (2), fix any point $b\in B$ and consider the flat $C\defeq B\cap \overline{\{b\}\cup A}$. We need to check that $A\parallel C$, which is equivalent to $A\subparallel C$ and $C\subparallel A$.

To prove $A\subparallel C$, it suffices to check that $A\subseteq\overline{\{a\}\cup C}$ for every $a\in A$. Fix any point $a\in A$. Given any $x\in A$, apply (4) and find $y\in B$ such that $\Aline ax\parallel \Aline by$ and hence $x\in\Aline ax\subseteq\overline{\{a\}\cup\Aline by}=\overline{\{a,b,y\}}$. If $x\in\Aline ab$, then $x\in\Aline ab\subseteq\overline{\{a\}\cup C}$. If $x\notin\Aline ab$, then by the Exchange Property, the inclusion $x\in \overline{\{a,b,y\}}\setminus\Aline ab$ implies $y\in\overline{\{a,b,x\}}\cap B\subseteq\overline{\{b\}\cup A}\cap B=C$ and hence $x\in \overline{\{a,b,y\}}\subseteq\overline{\{a\}\cup C}$. This completes the proof of the  relation $A\subparallel C$.

To show that $C\subparallel A$, we should check that $C\subseteq\overline{\{c\}\cup A}$  for every $c\in C$. Since $c\in C=B\cap\overline{\{b\}\cup A}\subseteq\overline{\{b\}\cup A}\setminus A$, the Exchange Property ensures that $b\in \overline{\{c\}\cup A}$ and hence $C\subseteq\overline{\{b\}\cup A}\subseteq\overline{\overline{\{c\}\cup A}\cup A}=\overline{\{c\}\cup A}$, witnessing that $C\subparallel A$.
\smallskip

The implication $(2)\Ra(3)$ is trivial.  
\smallskip

$(3)\Ra(1)$ Assume that $A\parallel C$ for some flat $C\subseteq B$. Then for every $a\in A$ we have $A\subseteq\overline{\{a\}\cup C}\subseteq\overline{\{a\}\cup B}$, witnessing that $A\subparallel B$. 
\end{proof}

\section{Parallelity in proaffine regular liners}

Let us recall that a liner $X$ is {\em proaffine} if $\forall o,x,y\in X\;\forall p\in\Aline yx\setminus\Aline ox\;\exists u\in\Aline oy\;\forall v\in \Aline oy\setminus\{u\}\;\;(\Aline vp\cap\Aline ox\ne\varnothing)$.

\begin{theorem}\label{t:proaffine-char} For a regular liner $X$, the following conditions are equivalent:
\begin{enumerate}
\item the liner $X$ is proaffine;
\item for every flat $A\subseteq X$ and points $a\in A$, $b\in X\setminus A$, $p\in \overline{A\cup\{b\}}\setminus A$, there exists a point $u\in\Aline ab$ such that $\forall v\in\Aline ab\setminus\{u\}\;\;(\Aline vp\cap A\ne\varnothing)$;
\item for every flat $A\subseteq X$ and point $b\in X\setminus A$ there exists at most one flat $B\subseteq X$ such that $b\in B$ and $B\parallel A$.
\end{enumerate}
\end{theorem}

\begin{proof} By Proposition~\ref{p:sr=>r=>wr} and Corollary~\ref{c:proregular=>ranked}, the regular liner $X$ is weakly regular, ranked and has the Exchange Property, and by Theorem~\ref{t:w-modular<=>}, the weakly regular ranked liner $X$ is weakly modular.
\smallskip 

$(1)\Ra(2)$ Given a flat $A\subseteq X$ and points  $a\in A$, $b\in X\setminus A$, and $p\in \overline{A\cup\{b\}}\setminus A$, we need to find a point $u\in\Aline ab$ such that $\Aline vp\cap A\ne\varnothing$ for every $v\in\Aline ab\setminus\{u\}$. By the regularity of $X$, there exist points $x\in A$ and $y\in\Aline ab$ such that $p\in\Aline xy$. Since the liner $X$ is proaffine, there exists a point $u\in\Aline ab$ such that $\Aline vp\cap\Aline ax\ne\varnothing$ for every $v\in \Aline ab\setminus\{u\}$. Since $A$ is a flat, $\Aline ax\subseteq A$ and hence for every $v\in\Aline ab\setminus\{u\}$, we have $\varnothing\ne \Aline vp\cap \Aline ax\subseteq \Aline vp\cap A$.  
\smallskip

$(1)\Ra(3)$  To prove the condition (3), take any flat $A\subseteq X$ and points $a\in A$, $b\in X\setminus A$,  and assume that $B,B'$ are two flats in $X$ such that $b\in B\cap B'$ and $B\parallel A\parallel B'$. We need to show that $B=B'$. By symmetry, it suffices to check that $B\subseteq B'$. Given any point $y\in B$, we should prove that  $y\in B'$. If $y=b$, then $y=b\in B'$ and we are done. So, assume that  $y\ne b$. 

Since $B\subparallel A$ and $b\in B\setminus A$, by Theorem~\ref{t:subparallel-char} and Proposition~\ref{p:subparallel-char4},  there exists a point $x\in A$ such that $\Aline ax\parallel \Aline by$. By Corollary~\ref{c:parallel}, $\|\Aline ax\|=\|\Aline by\|=2$, which means that $L\defeq \Aline ax$ is a line. It follows from $\Aline by\parallel L$ that $\Aline by\subseteq\overline{L\cup\{b\}}$.

Since $A\subparallel B'$, by Proposition~\ref{p:subparallel-char4}, there exists a point $y'\in B'$ such that $\Aline b{y'}\parallel \Aline ax$ and hence $\Aline b{y'}\subseteq \overline{\Aline ax\cup\{b\}}=\overline{L\cup\{b\}}$. Since $\Aline by\parallel\Aline ax\parallel\Aline b{y'}$ and  $b\in X\setminus A\subseteq X\setminus\Aline ax$, we can apply Theorem~\ref{t:parallel-char} and conclude that $\Aline by\cap\Aline ax=\varnothing=\Aline b{y'}\cap\Aline ax$. Then $\Aline by$ and $\Aline b{y'}$ are two lines in the plane $\overline{L\cup\{b\}}$ with $\Aline by\cap L=\varnothing=\Aline b{y'}\cap L$. By Theorem~\ref{t:Proclus<=>}, $\Aline by=\Aline b{y'}$ and hence $y\in \Aline by=\Aline b{y'}\subseteq B'$. 
\smallskip

The implications $(3)\Ra(1)\Leftrightarrow(2)$ follow immediately from Theorem~\ref{t:Proclus<=>}and Corollary~\ref{c:parallel}.
\end{proof}

In proaffine regular liners, the parrallelism is an equivalence relation (which is a generalization of  Proposition I.30 on parallel lines in Euclid's ``Elements'').
 
\begin{theorem}\label{t:Proclus-lines} For any lines $A,B,C$ in a proaffine regular liner $X$, we have the implication $(A{\parallel}B\wedge B{\parallel}C)\Ra (A{\parallel}C)$.
\end{theorem}

\begin{proof}  Assuming that $A{\parallel}B\wedge B{\parallel}C$, we should prove that $A{\parallel}C$. This is trivially true if $A=B$ or $B=C$. So, we assume that $A\ne B\ne C$. Applying Proposition~\ref{p:para+intersect=>coincide}, we conclude that $A\cap B=\varnothing=B\cap C$.

If $A\cap C\ne\varnothing$, then $A=C$ by Theorem~\ref{t:proaffine-char}. So, assume that $A\cap C=\varnothing$. In this case $2=\|A\|<\|A\cup C\|$. If $A\cap\overline{B\cup C}\ne \varnothing$, then fix any point $a\in A\cap\overline{B\cup C}$ and observe that the relation $A\parallel B$ implies $A\subseteq\overline{\{a\}\cup B}\subseteq\overline{\overline{B\cup C}\cup B}=\overline{B\cup C}$. Then $A\cup C\subseteq\overline{B\cup C}$ and hence $\|A\cup C\|\le \|B\cup C\|$. By Theorem~\ref{t:parallel-char}, $B\parallel C$ and $B\cap C=\varnothing$ imply $\|B\cup C\|=3$. Taking into account that $2<\|A\cup C\|\le\|B\cup C\|=3$, we conclude that $\|A\cup C\|=3$ and hence $A\parallel C$ by Theorem~\ref{t:parallel-char}.

It remains to consider the case $A\cap\overline{B\cup C}=\varnothing$.  Fix any point $a\in A$ and observe that $A\parallel B$ implies $A\subseteq\overline{\{a\}\cup B}$. Then  $A\cup B\cup C\subseteq\overline{\{a\}\cup B\cup C}$ and hence $\|A\cup B\cup C\|\le 1+\|B\cup C\|=4$. On the other hand, by the rankedness of $X$, the inequality $\overline{B\cup C}\ne\overline{A\cup B\cup C}$ ensure that $3=\|B\cup C\|<\|A\cup B\cup C\|\le 4$, which implies 
$$\|A\cup B\cup C\|=\|\{a\}\cup B\cup C\|=4.$$
It follows from $a\notin \overline{B\cup C}$ that $B\ne\overline{\{a\}\cup B}$ and $C\ne\overline{\{a\}\cup C}$ and hence $\|\overline{\{a\}\cup B}\|=3=\|\overline{\{a\}\cup C}\|$, by the rankedness of $X$.  Consider the flat $\Lambda\defeq\overline{\{a\}\cup B}\cap\overline{\{a\}\cup C}$. By  Proposition~\ref{p:sr=>r=>wr}, Corollary~\ref{c:proregular=>ranked} and Theorem~\ref{t:w-modular<=>}, the regular liner $X$ is weakly modular. Then$$\|\Lambda\|=\|\overline{\{a\}\cup B}\|+\|\overline{\{a\}\cup C}\|-\|\overline{\{a\}\cup B}\cup \overline{\{a\}\cup C}\|=3+3-4=2,$$
which means that $\Lambda$ is a line.

\begin{claim}\label{cl:B||Lambda||C} $B\parallel\Lambda$ and $\Lambda\parallel C$.
\end{claim}

\begin{proof} First we show that $\Lambda\cap B=\varnothing$. In the opposite case, we could find a point $\lambda\in\Lambda\cap B$ and conclude that $\lambda\in \Lambda\cap B\subseteq \overline{\{a\}\cup C}\setminus C$. By the Exchange Property, $a\in \overline{\{\lambda\}\cup C}\subseteq \overline{B\cup C}$, which contradicts our assumption $A\cap \overline{B\cup C}=\varnothing$. This contradiction shows that $\Lambda\cap B=\varnothing$. Since $B\ne B\cup\Lambda\subseteq\overline{\{a\}\cup B}$, the rankedness of $X$ implies  $2=\|B\|<\|B\cup\Lambda\|\le\|\{a\}\cup B\|\le 3$ and hence $\|B\cup\Lambda\|=3$. By Theorem~\ref{t:parallel-char}, $B\parallel \Lambda$. By analogy we can prove that $C\parallel \Lambda$. 
\end{proof}

Since $a\in A\cap\Lambda$ and both lines $A$ and $\Lambda$ are parallel to $B$, Theorem~\ref{t:proaffine-char} ensures that $A=\Lambda$. Applying Claim~\ref{cl:B||Lambda||C}, we finally obtain $A\parallel C$.
\end{proof}
 
\begin{corollary}\label{c:subparallel-transitive} For any flats $A,B,C$ in a proaffine regular liner $X$ we have the implication $(A{\subparallel}B\wedge B{\subparallel}C)\Ra (A\subparallel C)$. 
\end{corollary}

\begin{proof} Assume that $A{\subparallel}B\wedge B{\subparallel}C$. If $A\cap B\ne\varnothing$, then $A\subseteq B$, by Theorem~\ref{t:subparallel-char}. Then for every $a\in A\subseteq B$, the relation $B\subparallel C$ ensures that $A\subseteq B\subseteq\overline{\{a\}\cup C}$, witnessing that $A\subparallel C$. So, we assume that $A\cap B=\varnothing$.  

If $B\cap C\ne\varnothing$, then $B\subseteq C$, by Theorem~\ref{t:subparallel-char}. Then for every $a\in A$, the relation $A\subparallel B$ ensures that $A\subseteq \overline{\{a\}\cup B}\subseteq\overline{\{a\}\cup C}$, witnessing that $A\subparallel C$. So, we assume that $B\cap C=\varnothing$.  

If $B=\varnothing$, then the relation $A\subparallel B$ implies that for every $a\in A$,
$$A\subseteq \overline{\{a\}\cup B}=\{a\}\subseteq\overline{\{a\}\cup C},$$witnessing that $A\subparallel C$. So, we can assume that the flat $B$ contains some point $b$.

\begin{claim}\label{cl:ax||cz} For every points $a,x\in A$ and $c\in C$, there exists a point $z\in C$ such that $\Aline ax\parallel \Aline cz$.
\end{claim}

\begin{proof} Since $A\subparallel B$ and $A\cap B=\varnothing$, we can apply Proposition~\ref{p:subparallel-char4} and find a point $y\in B$ such that $\Aline ax\parallel \Aline by$. Since $B\subparallel C$ and $B\cap C=\varnothing$, we can apply Proposition~\ref{p:subparallel-char4} and find a point $z\in C$ such that $\Aline by\parallel \Aline cz$. 

By Corollary~\ref{c:parallel}, $\|\Aline ax\|=\|\Aline by\|=\|\Aline cz\|\le 2$. 
If $\|\Aline ax\|=1$, then $\|\Aline cz\|=\|\Aline ax\|=1$ and $\Aline ax=\{a\}\parallel \{c\}=\Aline cz$ are parallel singletons.
If $\|\Aline ax\|=2$, then $\|\Aline cz\|=\|\Aline by\|=\|\Aline ax\|=2$
and hence the flats $\Aline cz$, $\Aline by$, $\Aline ax$ are lines.  By Theorem~\ref{t:Proclus-lines},  $\Aline ax\parallel \Aline cz$.
\end{proof}
If $A\cap C\ne\overline \varnothing$, then we can fix a point $a\in A\cap C$. By Claim~\ref{cl:ax||cz}, for every $x\in A $, there exists $z\in C $ such that  $\Aline ax\parallel\Aline az$. Since $a\in\Aline ax\cap\Aline az$, Theorem~\ref{t:parallel-char} implies that $x\in\Aline ax=\Aline az\subseteq C$ and hence $A\subseteq C$ and $A\subparallel C$

If $A\cap C=\varnothing$, then Claim~\ref{cl:ax||cz} and Proposition~\ref{p:subparallel-char4} imply $A\subparallel C$.
\end{proof}

Corollary~\ref{c:subparallel-transitive} and Definition~\ref{d:parallel} imply the following corollary saying that for proaffine regular liners, the parallelism relation $\parallel$ is an equivalence relation. 

\begin{corollary}\label{c:parallel-transitive} For every flats $A,B,C$ in a proaffine regular liner $X$, we have the implication  $(A{\parallel}B\wedge B{\parallel}C)\Ra (A\parallel C).$ 
\end{corollary}

\begin{exercise} Find an example of an affine $3$-regular space containing lines $L_1,L_2,L_3$ such that $L_1\parallel L_2\parallel L_3$ but $L_1\nparallel L_3$.
\smallskip

\noindent{\em Hint:} Look at Example~\ref{ex:HTS}.
\end{exercise}



The following theorem generalizes the Proclus Parallelity Postulate, proved in Proposition~\ref{p:Proclus-Postulate}.

\begin{theorem}\label{t:Proclus2} Let $A,B$ be parallel flats in a proaffine regular liner $X$ and $L$ be a line in the flat $\overline{A\cup B}$. If $|L\cap A|=1$, then $|L\cap B|=1$.
\end{theorem}

\begin{proof}  To derive a contradiction, assume that $|L\cap A|=1$ and $|L\cap B|\ne 1$. 
If $A\cap B\ne\varnothing$, then $A\parallel B$ implies $A=B$, by Theorem~\ref{t:parallel-char}. In this case $L\subseteq\overline{A\cup B}=\overline A=A$, which contradicts $|L\cap A|=1$ and shows that $A\cap B=\varnothing$. Assuming that $L\cap B\ne\emptyset$ and taking into account that $|L\cap B|>1$, we conclude that $L\subseteq B$, which contradicts $|A\cap L|=1$. This contradiction shows that  $L\cap B=\varnothing$. 

Since $|L\cap A|=1$, there exists a unique point $a\in L\cap A $. Choose any point $c\in L\setminus\{a\}$. Assuming set $B$, we conclude that $L\subseteq \overline{A\cup B}=\overline A=A$ and hence $|L\cap A|=|L|>1$, which contradicts our assumption. This contradiction shows that $B\ne\varnothing$. Then we can fix a point $b\in B $. It follows from $L\cap B=\varnothing$  that $b\notin L$ and hence $P\defeq\overline{L\cup\{b\}}$ is a plane, by Proposition~\ref{p:L+p=plane}. Since $A\parallel B$ and $L\subseteq\overline{A\cup B}\subseteq \overline{\{a\}\cup B}$, there exists a finite set $F\subseteq B$ such that $b\in F$ and  $L\subseteq \overline{\{a\}\cup F}$. We can assume that $F$ has the smallest possible cardinality. In this case, the set $\{a\}\cup F$ is independent and hence $\|P\cup F\|=\|\{a\}\cup F\|=\|\{a\}\cup F\|=1+|F|$. By Proposition~\ref{p:sr=>r=>wr} and Theorem~\ref{t:w-modular<=>}, the regular liner $X$ is locally modular. Then 
$$\|P\cap \overline F\|=\|P\|+\|F\|-\|P\cup F\|=3+|F|-\|\{a\}\cup F\|=3+|F|-(1+|F|)=2$$and hence $\Lambda\defeq P\cap\overline F$ is a line in $B$ such that $b\in\Lambda\subseteq P=\overline{L\cup\{b\}}$ and $\Lambda\cap L\subseteq B\cap L=\varnothing$. By Theorem~\ref{t:subparallel-char}, $\Lambda\subparallel L$.

Since $A\parallel B$ and $A\cap B=\varnothing$, we can apply Proposition~\ref{p:subparallel-char4} and find a line $L'\subseteq A$ such that $a\in L'$ and $L'\parallel \Lambda$. Then $L$ and $L'$ are two lines in $X$ that are parallel to the line $\Lambda$ and contain the point $a$.  By Theorem~\ref{t:proaffine-char}, $L=L'\subseteq A$, which contradicts $|L\cap A|=1$. This contradiction shows that $|L\cap B|=1$.  
\end{proof}

\begin{theorem}\label{t:subparallel-via-base} A nonempty flat $A$ in a proaffine regular liner $X$ is subparallel to a flat $B$ in $X$ if and only if there exists a point $a\in A $ and a set $\Lambda\subseteq A$ such that $A=\overline\Lambda$ and $\forall x\in \Lambda\;(\Aline ax\subparallel B)$.
\end{theorem}

The proof of Theorem~\ref{t:subparallel-via-base} is based on the following two lemmas (that will be used also in the proof of Theorem~\ref{t:affine-char2}).

\begin{lemma}\label{l:subparallel+1} Let $X$ be a proaffine regular liner and $A,B$ be two flats in $X$ such that $A\subparallel B$ and $\|A\|<\w$. If $L$ is a line in $X$ such that $L\cap A\ne\varnothing$ and $L\subparallel B$, then $\overline{A\cup L}\subparallel B$.
\end{lemma}

\begin{proof} Since $A\cap L\ne\varnothing$, there exists a point $a\in A\cap L $. If $a\in B$, then $a\in A\subparallel B$ and $a\in L\subparallel B$ imply $A\cup L\subseteq B$, $\overline{A\cup L}\subseteq B$, and $\overline{A\cup L}\subparallel B$. So, we assume that $a\notin B$ and hence $A\not\subseteq B$. Since $A\subparallel B$, $L\subparallel B$ and $a\in A\cap L\setminus B$, we can apply Theorem~\ref{t:subparallel-char} and conclude that $A\cap B=\varnothing=L\cap B$.

If $L\subseteq A$, then $\overline{A\cup L}=\overline{A}=A\subparallel B$ and we are done. So, assume that $L\not\subseteq A$. It follows from  $L\subparallel B$ that $B$ contains some point $b\in B $. By Proposition~\ref{p:subparallel-char4}, the flats $A'\defeq B\cap\overline{A\cup\{b\}}$ and $L'\defeq B\cap\overline{L\cup\{b\}}$ are parallel to the flats $A$ and $L$, respectively. By Theorem~\ref{t:subparallel-char}, $A'\cap A=\varnothing=L'\cap L$ and $A'\cup L'\subseteq \overline{A\cup L\cup\{b\}}$. We claim that the flat $I\defeq \overline{A\cup L}\cap\overline{A'\cup L'}$ is empty. To derive a contradiction, assume that $I\ne\varnothing$. 

Since $L\not\subseteq A$, there exists a point $c\in L\setminus A\subseteq L\setminus\{a\}$. 
Proposition~\ref{p:add-point-to-independent} and Theorem~\ref{t:Max=codim} ensure that $\|A\cup L\|=\|A\cup\{c\}\|=\|A\|+1$. Choose any point $c'\in L'\setminus\{b\}$. Assuming that $c'\in A'$, we can apply Proposition~\ref{p:subparallel-char4}(4) and find a point $c''\in A$ such that $\Aline a{c''}\parallel \Aline b{c'}$. Since $\Aline a{c''}\parallel \Aline b{c'}=L'\parallel L=\Aline ac$, Theorem~\ref{t:Proclus-lines} ensures that $\Aline a{c''}\parallel \Aline ac$ and hence $c\in \Aline ac=\Aline a{c''}\subseteq A$, which contradicts the choice of the point $c\notin A$. This contradiction shows that $c'\notin A'$. Then $\|A'\cup L'\|=\|A'\cup\{c'\}\|=\|A'\|+1=\|A\|+1$, see Propositions~\ref{p:add-point-to-independent} and Corollary~\ref{c:parallel}. It follows from $A\parallel A'$, $L\parallel L'$ and $b\in A'\cap L'$ that $A'\cup L'\subseteq\overline{A\cup L\cup\{b\}}$ and hence 
$$\|\overline{A\cup L}\cup\overline{A'\cup L'}\|=\|A\cup L\cup A'\cup L'\|=\|A\cup L\cup\{b\}\|\le\|A\cup L\|+1=(\|A\|+1)+1=\|A\|+2.$$

By Proposition~\ref{p:sr=>r=>wr} and Theorem~\ref{t:w-modular<=>}, the regular liner $X$ is weakly modular. Since  $I=\overline{A\cup L}\cap\overline{A'\cup L'}\ne\varnothing$, we can apply the weak modularity of the liner $X$ and conclude that  
$$
\|I\|=\|A\cup L\|+\|A'\cup L'\|-\|A\cup L\cup A'\cup L'\|\ge
\|A\|+1+\|A'\|+1-(\|A\|+2)=\|A\|.
$$
Therefore, $A\cap I\subseteq A\cap B=\varnothing$, and  $A$ and $I$ are disjoint flats of rank $\ge \|A\|$ in the flat $\overline{A\cup L}$ that has rank $\|A\|+1$. The rankedness of the regular liner $X$ ensures that $\|I\|=\|A\|$. By Theorem~\ref{t:parallel-char}, $A\parallel I$, and by Theorem~\ref{t:Proclus2}, $L\cap I\ne\varnothing$, which is impossible as $L\cap I\subseteq L\cap B=\varnothing$. This contradiction shows that $\overline{A\cup L}\cap  \overline{A'\cup L'}=I=\varnothing$. 
Then $\overline{A\cup L}$ and $ \overline{A'\cup L'}$ are two flats of rank $\|A\|+1$  in the flat $\overline{A\cup L\cup\{b\}}$ that has rank $\|A\|+2$. Applying Theorem~\ref{t:parallel-char}, we conclude that $\overline{A\cup L}\parallel \overline{A'\cup L'}$ and hence $\overline{A\cup L}\subparallel B$, by Proposition~\ref{p:subparallel-char4}.
\end{proof}

\begin{lemma}\label{l:subparallel+2} Let $A,B$ be flats in a proaffine regular liner $X$. Let $L$ be a line in $X$ such that $A\cap L\ne\varnothing$. If $A\subparallel B$ and $L\subparallel B$, then $\overline{A\cup L}\subparallel B$.
\end{lemma}

\begin{proof} If $A\cap B\ne\varnothing$, then $A\subparallel B$ implies $A\subseteq B$, by Theorem~\ref{t:subparallel-char}. Then $\varnothing\ne L\cap A\subseteq A\subseteq B$ and $L\subparallel B$ implies $L\subseteq B$, by Theorem~\ref{t:subparallel-char}. It follows from $A\cup L\subseteq B$ that $\overline{A\cup L}\subseteq B$ and hence $\overline{A\cup L}\subparallel B$. 
By analogy we can show that $L\cap B\ne\varnothing$ implies $\overline{L\cup A}\subparallel B$.

So, assume that $A\cap B=\varnothing=L\cap B$. Since $A\cap L\ne\varnothing$, there exists a point $a\in A\cap L $. The subparallelity $A\subparallel B$ and $L\subparallel B$ implies $A\cup L\subseteq\overline{\{a\}\cup B}$ and hence $\overline{A\cup L}\subseteq\overline{\{a\}\cup B}$. The subparallelity $\overline{A\cup L}\subparallel B$ will follow from Theorem~\ref{t:subparallel-char} as soon as we prove that $\overline{A\cup L}\cap B=\varnothing$. To derive a contradiction, assume that the set $\overline{A\cup L}\cap B$ contains some point $x$. By Proposition~\ref{p:aff-finitary}, there exists a flat $F\subseteq A$ of finite rank such that $a\in F$ and $x\in \overline{F\cup L}$. It follows from $ F\subseteq A\subparallel B$ that $F\subparallel B$. By Lemma~\ref{l:subparallel+1}, $\overline{F\cup L}\subparallel B$. Since $a\in\overline{F\cup L}\setminus B$, Theorem~\ref{t:subparallel-char} ensures that $x\in \overline{F\cup L}\cap B=\varnothing$, which contradicts the choice of the point $x$. This contradiction shows that $\overline{A\cup L}\cap B=\varnothing$ and hence  $\overline{A\cup L}\subparallel B$.
\end{proof} 

\begin{proof}[Proof of Theorem~\ref{t:subparallel-via-base}.] The ``only if'' part of Theorem~\ref{t:subparallel-via-base} is trivial. To prove the ``if'' part, take any flats $A,B$ in a proaffine regular liner $X$, any point $a\in A $ and any set $\Lambda\subseteq A$ such that $A=\overline{\Lambda}$ and $\Aline xa\subparallel B$ for every $x\in\Lambda$. It follows from $a\in A=\overline\Lambda$ that $\Lambda$ is not empty.

 If $a\in B$, then for every $x\in\Lambda$, the subparallelity $\Aline xa\subparallel B$ implies $x\in\Aline xa\subseteq\overline{B\cup\{a\}}=B$ and hence $\Lambda\subseteq B$ and $A=\overline{\Lambda}\subseteq \overline{B}=B$. Then $A\subparallel B$ and we are done.

So, assume that $a\notin B$.  Write the set $\Lambda\setminus\{a\}$ as $\{x_\alpha\}_{\alpha\in\kappa}$ for some cardinal $\kappa$. For every ordinal $\alpha\le\kappa$, consider the flat $A_{\alpha}\defeq\overline{\{a\}\cup\{x_\beta:\beta<\alpha\}}$ and observe that $(A_{\alpha})_{\alpha\le\kappa}$ is an increasing sequence of flats in $X$ such that $A_0=\{a\}$ and $A_\kappa=A$.  Proposition~\ref{p:aff-finitary} implies the equality $A_\alpha=\bigcup_{\beta<\alpha}A_\beta$ for every limit ordinal $\alpha>0$.

By transfinite induction, we shall prove that $A_\alpha\subparallel B$ for every ordinal $\alpha\le\kappa$. For $\alpha=0$, we have $A_0=\{a\}\subseteq \overline{B\cup\{a\}}\setminus B$ and hence $A_0\subparallel B$ by Theorem~\ref{t:subparallel-char}. Assume that for some nonzero ordinal $\alpha\le\kappa$ and all ordinals $\beta<\alpha$ we know that $A_\beta\subparallel B$. Since $a\in A_\beta\setminus B$, Theorem~\ref{t:subparallel-char} ensures that $A_\beta\subseteq\overline{B\cup\{a\}}\setminus B$. If the ordinal $\alpha$ is limit, then Proposition~\ref{p:aff-finitary} implies  $A_\alpha=\bigcup_{\beta<\alpha}A_\beta\subseteq\overline{B\cup\{a\}}\setminus B$. Applying Theorem~\ref{t:subparallel-char}, we conclude that $A_\alpha\subparallel B$. 
If the ordinal $\alpha$ is not limit, then $\alpha=\beta+1$ for some ordinal $\beta<\alpha$. In this case $A_\alpha=\overline{A_\beta\cup \{x_\beta\}}=\overline{A_\beta\cup\Aline a{x_\beta}}$. Since $x_\beta\in\Lambda\setminus\{a\}$, $\Aline a{x_\beta}$ is a line subparallel to $B$. By the inductive assumption, $A_\beta\subparallel B$. Applying Lemma~\ref{l:subparallel+1}, we conclude that $A_\alpha=\overline{A_\beta\cup\Aline a{x_\beta}}\subparallel B$. This completes the inductive step. By the Principle of Transifinite Induction, $A_\alpha\subparallel B$ for all ordinals $\alpha\le\kappa$. In particular, $A=A_\kappa\subparallel B$.
\end{proof}



\section{Parallelity in affine and Playfair liners}

Let us recall that a liner $X$ is {\em affine} if  
$\forall o,x,y\in X\;\;\forall u\in\Aline oy\;\;\forall p\in\Aline ux\setminus\Aline ox\;\;\exists v\in\Aline oy\;\;\forall w\in\Aline oy\;\;\;(u=v\;\Leftrightarrow\;\Aline vp\cap\Aline ox=\varnothing)$.

\begin{theorem}\label{t:affine-char2} For a regular liner $X$, the following conditions are equivalent:
\begin{enumerate}
\item the liner $X$ is affine;
\item for every flat $A\subseteq X$ and points $a\in A$, $b\in X\setminus A$, $p\in \overline{A\cup\{b\}}\setminus A$, there exists a unique point $u\in\Aline ab$ such that $\Aline up\cap A=\varnothing$;
\item for every flat $A\subseteq X$ and points $a\in A$, $b\in X\setminus A$ with $\Aline ab\ne\{a,b\}$, there exists a unique flat $B\subseteq X$ such that $b\in B$ and $B\parallel A$.
\end{enumerate}
\end{theorem}

\begin{proof} By Proposition~\ref{p:sr=>r=>wr}, Corollary~\ref{c:proregular=>ranked} and Theorem~\ref{t:w-modular<=>}, the regular liner $X$ is weakly modular, ranked and has the Exchange Property.
\smallskip

 $(1)\Ra(2)$ Assume that the liner $X$ is affine. To prove the condition (2), take any 
flat $A\subseteq X$ and points $a\in A $, $b\in X\setminus A$, $p\in \overline{A\cup\{b\}}\setminus A$.
By Theorem~\ref{t:proaffine-char}, the set $I\defeq\{u\in\Aline ab:\Aline up\cap A=\emptyset\}$ contains at most one point. So, it remains to prove that the set $I$ is not empty.

If $p\in\Aline ab$, then $p\in I$. So, assume that $p\notin\Aline ab$. In this case the flat $P\defeq\overline{\{a,b,p\}}$ is a plane. Since $p\in\overline{A\cup\{b\}}=\bigcup_{F\in[A]^{<\w}}\overline{F\cup\{b\}}$, there exists a finite set $F\subseteq A$ such that $a\in F$ and $p\in\overline{F\cup\{b\}}$. Proposition~\ref{p:add-point-to-independent} and Theorem~\ref{t:Max=codim} imply that $\|F\cup P\|=\|F\cup\{b\}\|=\|F\|+1$. The weak modularity of $X$ ensures that
$$\|P\cap \overline F\|=\|P\|+\|\overline F\|-\|P\cup\overline F\|=3+\|F\|-(\|F\|+1)=2.$$
Therefore, $\Lambda\defeq P\cap\overline F\subseteq P\cap A$ is a line in the plane $P=\overline{\{a,b,p\}}$. Since $p\in X\setminus A\subseteq X\setminus\Lambda$, the flat $\overline{\Lambda\cup\{p\}}$ is a plane in the plane $P$. The rankedness of the liner $X$ guarantees that $\overline{\Lambda\cup\{p\}}=P$ and hence $b\in P\setminus A=\overline{\{p\}\cup \Lambda}\setminus\Lambda$. The Exchange Property ensures that $p\in\overline{\Lambda\cup\{b\}}\setminus\Lambda$. By Theorem~\ref{t:affine-char1}, there exists a point $u\in\Aline ab$ such that $\Aline up\cap \Lambda=\varnothing$. 
It follows from $b\in P\setminus A$ that $\Lambda=P\cap\overline F\subseteq P\cap A\subsetneq P$. The rankedness of the liner $X$ ensures that $2=\|\Lambda\|=\|P\cap A\|<\|P\|=3$ and  $\Lambda=P\cap A$. 
Then $\varnothing=\Aline up\cap \Lambda=\Aline up\cap(P\cap A)=(\Aline up\cap P)\cap A=\Aline up\cap A$ and hence $u\in I$.
\smallskip

$(2)\Ra(3)$ Assume that the condition (2) holds. To prove that the condition (3) is satisfied, take any flat $A\subseteq X$ and points $a\in A $ and $b\in X\setminus A$ with $\Aline ab\ne\{a,b\}$. 
 Consider the family $\mathcal F$ of all flats $B$ in $X$ such that $b\in B\subseteq \overline{A\cup\{b\}}\setminus A$. Since $b\notin A$, the singleton $\{b\}$ belongs to the family $\F$. Endow the family $\mathcal F$ with the inclusion partial order and take any maximal linearly ordered chain $\mathcal M$ in $\mathcal F$. Proposition~\ref{p:aff-finitary} implies that the union $M\defeq\bigcup\mathcal M$ is a flat in $X$. It is clear that $b\in M\subseteq \overline{A\cup\{b\}}\setminus A$ and hence $M\subparallel A$, by Theorem~\ref{t:subparallel-char}. 

By Proposition~\ref{p:subparallel-char4}, the flat $A'\defeq A\cap\overline{M\cup\{a\}}$ is parallel to the flat $M$.
Assuming that $A'\ne A$, we can choose a point $c\in A\setminus A'\subseteq A\setminus\{a\}$ and conclude that  $L\defeq \Aline ac$ is a line in the flat $A$. By our assumption, there exists a point $p\in \Aline ab\setminus\{a,b\}$. Assuming that $p\in \Aline ca$, we conclude that $b\in \Aline ab=\Aline pa=\Aline ca\subseteq A$, which contradicts the choice of the point $b$. Therefore, $p\notin\Aline ca$ and hence $p\in\overline{\{a,c,b\}}\ne\Aline ac\cup\{b\}$.  
 By Theorem~\ref{t:affine-char1}, there exists a line $\Lambda$ such that $b\in \Lambda\subseteq \overline{L\cup\{b\}}\setminus L$. By Theorem~\ref{t:subparallel-char}, 
$\Lambda\subparallel L\subseteq A$ and by Corollary~\ref{c:subparallel}, $\Lambda\parallel L$. Since  
$\Lambda\subparallel A$ and $b\in M\cup\Lambda$, we can apply Lemma~\ref{l:subparallel+2} and conclude that $M\subseteq \overline{M\cup \Lambda}\subparallel A$. The maximality of the chain $\mathcal M$ guarantees that $L\subseteq \overline{M\cup L}=M$. The parallelity $L\parallel \Lambda\subseteq M$ ensures that $c\in\Aline ca=L \subseteq A\cap \overline{\Lambda\cup\{a\}}\subseteq A\cap\overline{M\cup\{a\}}=A'$, which contradicts the choice of the point $c$. This contradiction shows that $A=A'$ and hence $A\parallel M$. Theorem~\ref{t:proaffine-char}, ensures that $M$ is the unique flat such that $b\in M$ and $M\parallel A$.
\smallskip

$(3)\Ra(1)$ Assume that the liner $X$ satisfies the condition (3). We claim that $X$ satisfies the condition (4) of Theorem~\ref{t:affine-char1}. Fix a line $L\subseteq X$ and a point $p\in X$ with $\overline{L\cup\{p\}}\ne L\cup\{p\}$. Then there exists a point $z\in \overline{L\cup\{p\}}\setminus(L\cup\{p\})$. 

If $\Aline pz\cap L=\varnothing$, then $\Lambda\defeq\Aline pz\subseteq\overline{L\cup\{p\}}$ is a line such that $p\in \Lambda$ and $\Lambda\parallel L$. Fix any point $o\in L $. By the regularity of $X$, there exist points $a\in L$ and $b\in\Aline op$ such that $z\in \Aline ab$.  Assuming that $b=p$, we conclude that $z\in\Aline ab\subseteq\Aline ap$ and hence $a\in \Aline pz\cap L=\varnothing$, which is a contradiction showing that $b\ne p$. Assuming that $b=a$, we conclude that $z\in\Aline ba\subseteq L$, which contradicts the choice of $z$. Therefore $b\in\Aline ap\setminus\{a,p\}$. The condition (3) ensures that $\Lambda$ is a unique line such that $p\in\Lambda\subseteq \overline{L\cup\{p\}}\setminus L$.

Next, assume that $\Aline pz\cap L\ne\varnothing$. Fix any point $a\in\Aline pz\cap L $. Then $z\in\Aline pa\setminus\{p,a\}$ and by the condition (3), there exists a unique flat $\Lambda$ such that $p\in \Lambda$ and $\Lambda\parallel L$. By Theorem~\ref{t:subparallel-char} and Corollary~\ref{c:subparallel}, $\Lambda$ is a unique line such that $p\in\Lambda\subseteq\overline{L\cup\{p\}}\setminus L$. 

This shows that the condition (4) of Theorem~\ref{t:affine-char1} is satisfied and the liner $X$ is affine. 
\end{proof} 


\begin{theorem}\label{t:Playfair} For every line $L$ in a Playfair liner $X$ and every point $x\in X\setminus L$, there exists a unique line $\Lambda$ in $X$ such that  $x\in\Lambda$ and $\Lambda\parallel L$.
\end{theorem}

\begin{proof} By Theorem~\ref{t:Playfair<=>}, there exists a unique line $\Lambda$ in $X$ such that  $x\in\Lambda\subseteq \overline{L\cup\{x\}}\setminus L$. By Theorem~\ref{t:Playfair}, the Playfair liner $X$ is $3$-regular and by Proposition~\ref{p:k-regular<=>2ex}, the $3$-regular liner $X$ is $3$-ranked. By Theorem~\ref{t:subparallel-char}, $\Lambda\subparallel L$, and by Corollary~\ref{c:subparallel}, $\Lambda\parallel L$.
\end{proof}

For regular Playfair liners, we have a higher dimensional counterpart of Theorem~\ref{t:Playfair}, which follows from Theorem~\ref{t:affine-char2}.

\begin{corollary}\label{c:Playfair} For any flat $A$ in a regular Playfair liner $X$ and every point $x\in X$ there exists a unique flat $B\subseteq X$ such that $x\in B$ and $B\parallel A$.
\end{corollary}

Let us recall that a liner $X$ is {\em Bolyai} if for every plane $P\subseteq X$, line $L\subseteq P$ and point $x\in P\setminus L$, there exists a line $\Lambda$ in $X$ such that $x\in \Lambda\subseteq P\setminus L$.

\begin{theorem}\label{t:Play-reg<=>Bolyai} For every liner $X$ the following conditions are equivalent:
\begin{enumerate}
\item $X$ is Playfair and regular;
\item $X$ is Bolyai, $3$-ranked, and for every lines $A,B,C\subseteq X$, if $A\parallel B$ and $B\parallel C$, then $A\parallel C$.
\end{enumerate}
\end{theorem}

\begin{proof} The implication $(1)\Ra(2)$ follows from Theorem~\ref{t:Proclus-lines} and Corollary~\ref{c:proregular=>ranked}.
\smallskip

$(2)\Ra(1)$. Assume that $X$ is Bolyai, $3$-ranked and the parallelity relation on lines in $X$ is transitive.  To show that the liner $X$ is Proclus, take any plane $P\subseteq X$, line $L\subseteq P$ and point $x\in P\setminus L$. Let $\Lambda_1,\Lambda_2$ be two lines in $X$ such that $x\in\Lambda_i\subseteq P\setminus L$ for every $i\in\{1,2\}$. Since the liner $X$ is $3$-ranked, Corollary~\ref{c:parallel-lines<=>} implies that $\Lambda_1\parallel L\parallel \Lambda_2$ and hence $\Lambda_1\parallel \Lambda_2$. Since $x\in\Lambda_1\cap\Lambda_2\ne\varnothing$, Corollary~\ref{c:parallel-lines<=>} ensures that $\Lambda_1=\Lambda_2$, witnessing that the Bolyai liner $X$ is Proclus and hence Playfair. By Theorem~\ref{t:Playfair<=>}, the Playfair liner $X$ is affine, $3$-regular and $3$-long.
\smallskip

To prove that the liner $X$ is regular, take any flat $A\subseteq X$ and points $o\in A$ and $u\in X\setminus A$. We have to prove that the set $\Lambda\defeq\bigcup_{x\in A}\bigcup_{y\in\overline{ou}}\Aline xy$ is flat. Assuming that $\Lambda$ is not flat, we can find points $x,y\in\Lambda$ and $z\in\Aline xy\setminus\Lambda$. In this case $\{x,y\}\not\subseteq A\cup\Aline ou$. Without loss of generality, we can assume that $y\notin A\cup\Aline ou$. Since $y\in\Lambda$, there exist points $u'\in\Aline ou$ and $b\in A$ such that $y\in \Aline b{u'}$. Replacing $u$ by $u'$, we can assume that $y\in\Aline bu$. Since $X$ is affine, there exists a unique point $v\in\Aline ou$ such that $\Aline vy\cap\Aline ob=\varnothing$. By Corollary~\ref{c:parallel-lines<=>}, $\Aline vy\parallel \Aline ob$. It follows from $y\notin A\cup\Aline ou$ that $b\notin\Aline ou$ and $u\notin A$. 

If $x\in\overline{\{o,u,b\}}$, then $z\in\Aline xy\subseteq\overline{\{o,u,b\}}$ and by the $3$-regularity of $X$, $z\in \Aline {\Aline ob}{\Aline ou}\subseteq\Aline A{\Aline ou}=\Lambda$, which contradicts the choice of $z$. This contradiction shows that $x\notin\overline{\{o,u,b\}}$. 

 Assuming that $z\in \Aline vy$, we conclude that $x\in\Aline zy\subseteq\Aline vy\subseteq \overline{\{o,u,b\}}$, which is not true. Therefore, $z\notin \Aline vy$.

If $\Aline xy\cap A\ne\varnothing$, then we lose no generality assuming that $x\in A$. Since $X$ is Playfair, there exists a unique line $L$ in $X$ such that $x\in L$ and $L\parallel \Aline vy$. It follows from $L\parallel \Aline vy\parallel \Aline ob$ that $L\parallel \Aline ob$ and hence $L\subseteq\overline{\{o,b,x\}}\subseteq A$. Since $\Aline vz\ne \Aline vy\parallel L$, the Proclus property of $X$ ensures that $\Aline vz\cap L$ contains some point $\alpha$. Then $z\in\Aline v\alpha\subseteq \Lambda$, which contradicts the choice of $z$. This contradiction shows that $\Aline xy\cap A=\varnothing$.  

Since $x\in \Lambda$, there exist points $w\in \Aline ou$ and $a\in A$ such that $x\in \Aline wa$. Since $X$ is Playfair, there exists a unique line $L$ such that $a\in L$ and $L\parallel \Aline vy$. It follows from $L\parallel \Aline vy\parallel \Aline ob$ that $L\parallel \Aline ob$ and hence $L\subseteq\overline{\{o,b,a\}}\subseteq A$.  We claim that $w\ne v$. To derive a contradiction, assume that $w=v$. By the Proclus Parallel Postulate, $\Aline wz=\Aline vz\ne\Aline wy=\Aline vy\parallel L$ implies $\Aline wz\cap L\ne\varnothing$. Then $z\in \Aline Lw\subseteq \Aline Aw\subseteq\Lambda$, which contradicts the choice of $z$. This contradiction shows that $w\ne v$. Since $\Aline vy\parallel \Aline ob$, the proaffinity of $X$ ensures that $\Aline wy\cap\Aline ob\ne\varnothing$. Replacing $u$ by $w$, we can assume that $u=w$. Then $z\in\Aline xy\subseteq\overline{\{w,a,b\}}$. It follows from $\Aline xy\cap\Aline ab\subseteq \Aline xy\cap A=\varnothing$ that $\Aline xy\parallel \Aline ab$. Since $X$ is Proclus, $\Aline wz\cap \Aline ab\ne\varnothing$ and hence $z\in\Aline {\Aline ab}w\subseteq \Aline Aw\subseteq\Lambda$, which contradicts the choice of $z$. This is the final contradiction completing the proof of the theorem.
\end{proof}

\begin{proposition}\label{p:ABC-regular} For any parallel lines $A,B,C$ in a Playfair liner $X$, the flat $\overline{A\cup B\cup C}$ is a regular subliner of $X$.
\end{proposition}

\begin{proof} By Theorem~\ref{t:Playfair<=>}, the Playfair liner $X$ is $3$-long, $3$-regular and affine. If $X$ is $4$-long, then by Theorem~\ref{t:4-long-affine}, the $4$-long affine liner $X$ is regular, and so is the flat $\overline{A\cup B\cup C}$ in $X$. So, assume that $X$ is not $4$-long. In this case, $X$ is a Hall liner, by Theorem~\ref{t:Hall<=>Playfair+Steiner}. Fix any points $e\in B$, $b\in B\setminus \{e\}$, $a\in A$ and $c\in C$. It follows from $A\parallel B\parallel C$ that $\overline{A\cup B\cup C}=\overline{\{e,a,b,c\}}$. 

Given any points $x,y\in X$, let $x\circ y$ be a unique point such that $\{x,y,x\circ y\}=\Aline xy$. Let $X_e$ be the set $X$ endowed with the binary operation $xy\defeq e\circ(x\circ y)$. By Proposition~\ref{p:Steiner=>loop3} and Theorem~\ref{t:Hall<=>Moufang}, the magma $X_e$ is a commutative Moufang loop of exponent 3. By Proposition~\ref{p:flat<=>subloop}, the flat $S\defeq\overline{\{e,a,b,c\}}$ is a subloop of the loop $X_e$, generated by the elements $a,b,c$. We claim that the subloop $S$ is associative. To derive a contradiction, assume that $S$ is not associative. In this case, Theorem~\ref{t:Moufang-loop} implies that the associator $[a,b,c]$ is not equal to $e$. Let $H=\{e,[a,b,c],[a,b,c]^2\}$ be the subgroup of the loop $S$, generated by the associator $[a,b,c]$. By Theorem~\ref{t:associator}, the associator $[a,b,c]$ belongs to the centre of the loop $S$. By Theorem~\ref{t:Moufang81}, the quotient loop $S/H$ is a commutative group of cardinality $|S/H|=27$. Let $q:S\to S/H$, $q:x\mapsto xH$, be the quotient homomorphism.
Taking into account that the group $S/H$ has cardinality $27$ and is generated by the elements $aH,bH,cH$, we conclude that $bH\ne H$.

By Lemma~\ref{l:Moufang-xyy} and Theorem~\ref{t:Moufang-loop}, the subloops $\overline{\{e,a,b\}}$ and $\overline{\{e,b,c\}}$ of the loop $X_e$ are associative. 
It follows from $A\parallel B\parallel C$ that $A=\{a,ab,ab^2\}$ and $C=\{c,bc,b^2c\}$. Since $A\parallel C$, either $(ab)\circ c=a\circ (bc)$ or $(ab)\circ c=a\circ (b^2c)$. We claim that $(ab)\circ c=a\circ(bc)$. In the opposite case, $(ab)\circ c=a\circ(b^2c)$. Observe that for every $x,y\in X$ the midpoint $x\circ y$ is equal to $(xy)^2$. Then $(ab)\circ c=a\circ(b^2c)$ implies $((ab)c)^2=(ab)\circ c=a\circ (b^2c)=(a(b^2c))^2$ and hence $(ab)c=((ab)c)^4=(a(b^2c))^4=a(b^2c)$ because the Moufang loop $X_e$ has exponent $3$. Applying to the equality $(ab)c=a(b^2c)$ the quotient homomorphism $q$ and taking into account that $S/H$ is a commutative group, we conclude that $aHbHcH=aHbHbHcH$ and after cancellation, $H=bH$, which contradicts the inequality $bH\ne H$ established in the preceding paragraph. This contradiction shows that $(ab)\circ c=a\circ(bc)$ and hence $((ab)c)^2=(ab)\circ c=a\circ(bc)=(a(bc))^2$. Since the loop $X_e$ has exponent 3, $(ab)c=((ab)c)^4=(a(bc))^4=a(bc)$ and hence $[a,b,c]=e$, which contradicts our assumption. This contradiction shows that the loop $S$ is associative. Then $|S|\le 27$ and the Hall liner $S=\overline{A\cup B\cup C}$ is regular, by Theorem~\ref{t:Hall<=>regular}.
\end{proof} 

\begin{exercise} Explore the subparallelity and parallelity relations for the liners in Examples~\ref{ex:Tao} and \ref{ex:HTS}.
\end{exercise}

\section{Triangles and parallelograms in liners}

\begin{definition} A triple of points $abc\in X^3$ in a liner $X$ is called a \index{triangle}\defterm{triangle} in $X$ if $\|\{a,b,c\}\|=3$. The points $a,b,c$ are called the \index{triangle!vertices}\defterm{vertices} of the triangle $abc$, and the lines $\Aline ab,\Aline bc,\Aline ca$ are called the \index{triangle!sides}\defterm{sides} of the triangle $abc$.
\end{definition}

\begin{picture}(100,50)(-170,-5)
\put(0,0){\line(1,0){50}}
\put(0,0){\line(1,1){40}}
\put(40,40){\line(1,-4){10}}
\put(0,0){\circle*{3}}
\put(-3,-8){$a$}
\put(40,40){\circle*{3}}
\put(40,43){$b$}
\put(50,0){\circle*{3}}
\put(48,-8){$c$}
\end{picture}

\begin{exercise} Show that every triangle $abc$ in a liner has $|\{a,b,c\}|=\|\{a,b,c\}\|=3$.
\end{exercise}

\begin{definition} A quadruple of points $abcd\in X^4$  in a liner $X$ is called a \index{parallelogram}\defterm{parallelogram} in $X$ if $\Aline ab\parallel\Aline cd\ne \Aline ab$ and $\Aline bc\parallel\Aline ad\ne \Aline bc$. The points $a,b,c,d$ are called the \index{parallelogram!vertices}\defterm{vertices}, the lines $\Aline ab,\Aline bc,\Aline cd,\Aline da$ are called the \index{parallelogram!sides}\defterm{sides}, and the lines $\Aline ac,\Aline bd$ are called the \index{parallelogram!diagonal}\defterm{diagonals} of the parallelogram $abcd$.
\end{definition}

\begin{picture}(300,60)(-140,10)
{\linethickness{0.8pt}
\put(15,20){\color{teal}\line(1,0){50}}
\put(25,60){\color{teal}\line(1,0){50}}
\put(15,20){\color{cyan}\line(1,4){10}}
\put(65,20){\color{cyan}\line(1,4){10}}
}
\put(15,20){\line(3,2){60}}
\put(65,20){\line(-1,1){40}}
\put(15,20){\circle*{3}}
\put(8,12){$a$}
\put(65,20){\circle*{3}}
\put(67,12){$d$}
\put(25,60){\circle*{3}}
\put(21,63){$b$}
\put(75,60){\circle*{3}}
\put(78,63){$c$}
\end{picture}

Corollary~\ref{c:parallel-lines<=>} implies the following useful characterization of parallelograms in $3$-ranked liners.

\begin{proposition}\label{p:tobe-parallelogram} A quadruple of points $abcd$ in a $3$-ranked liner $X$ is a parallelogram if and only if  $\Aline ab\cap\Aline cd=\varnothing=\Aline bc\cap\Aline ad$ and $\|\{a,b,c,d\}\|\le 3$.
\end{proposition}

\begin{exercise} Let $abcd$ be a parallelogram in a liner $X$. Show that $$adcb,badc,bcda,cbad,cdab,dabc,dcba$$ are parallelograms.
\end{exercise}

\begin{exercise} Let $abcd$ be a parallelogram in a liner $X$. Show that for any distinct points $x,y,z\in \{a,b,c,d\}$, the triple $xyz$ is a triangle in $X$.
\end{exercise}

\begin{proposition}\label{p:k-parallel=>parallelogram} If a liner $X$ is  $\kappa$-parallel for some finite nonzero cardinal $\kappa$, then for every triangle $abc$ there exists a point $d\in X$ such that $abcd$ is a parallelogram.
\end{proposition}

\begin{proof} Let $abc$ be any triangle in $X$ and $P\defeq\overline{\{a,b,c\}}$ be the plane generated by this triangle in $X$. Let $\mathcal L$ be the family of all lines in $P$. For every point $x\in P$ and line $\Lambda\subset P$, consider the families $$\mathcal L_x\defeq\{L\in\mathcal L:x\in L\},\quad\mathcal L^\bullet_{x,\Lambda}\defeq\{L\in\mathcal L_x:L\cap\Lambda\ne\varnothing\},\quad\mbox{and}\quad\mathcal L^\circ_{x,\Lambda}\defeq\{L\in\mathcal L_x:L\cap\Lambda=\varnothing\}.$$
Since $X$ is $\kappa$-parallel, $|\mathcal L^\circ_{x,\Lambda}|=\kappa$ whenever $x\notin\Lambda$. 
 
Consider the lines $H\defeq \Aline ab$ and $V\defeq\Aline bc$. Since $|\mathcal L^\circ_{c,H}|=\kappa>0$, there exists a line $H'\in\mathcal L^\circ_{c,H}$. Since $H\in \mathcal L^\circ_{b,H'}\setminus \mathcal L^\circ_{b,V}$ and $|\mathcal L^\circ_{b,H'}|=|\mathcal L^\circ_{b,V}|=\kappa<\w$, there exists a line $V'\in \mathcal L^\circ_{b,V}\setminus\mathcal L^\circ_{b,H'}$, which has a common point $d$ with the line $H'$. Then $\Aline ab\cap\Aline cd=H\cap H'=\varnothing=V\cap V'=\Aline ac\cap \Aline bd$. By Proposition~\ref{p:k-parallel=>3-ranked}, the $\kappa$-parallel liner $X$ is $3$-ranked. By Proposition~\ref{p:tobe-parallelogram}, $abcd$ is a parallelogram in $X$.  
\end{proof}

\begin{theorem}\label{t:parallelogram3+1}  A liner $X$ is Playfair if and only if $X$ is affine, $3$-ranked, and for any triangle $abc$ in $X$ there exists a unique point $d\in X$ such that $abcd$ is parallelogram.
\end{theorem}

\begin{proof} To prove the ``only if'' part, assume that a liner $X$ is Playfair.
By Theorem~\ref{t:Playfair<=>}, $X$ is affine and $3$-regular. By Proposition~\ref{p:k-regular<=>2ex}, the $3$-regular liner $X$ is $3$-ranked.  Take any triangle $abc$ in $X$ and consider the plane $P\defeq\overline{\{a,b,c\}}$. Since Playfair liners are $1$-parallel, we can apply Proposition~\ref{p:k-parallel=>parallelogram} and conclude that there exists a point $d\in P$ such that $abcd$ is a parallelogram. To prove that the parallelogram $abcd$ is unique, assume that $d'$ is another point such that $abcd'$ is a parallelogram. The parallelity relation $\Aline c{d'}\parallel \Aline ab$ implies $d'\in \overline{\{a,b,c\}}=P$. Since $abcd$ and $abcd'$ are parallelograms in the Playfair plane $P$, the parallelity relations $\Aline cd\parallel \Aline ab\parallel \Aline c{d'}$ impliy $\Aline cd=\Aline c{d'}$. By the reason, $\Aline b{d'}\parallel\Aline ac\parallel \Aline bd$ imply $\Aline bd=\Aline b{d'}$. Then $\{d\}=\Aline cd\cap\Aline bd=\Aline c{d'}\cap \Aline b{d'}=\{d'\}$, witnessing that $abcd=abcd'$ and the parallelogram $abcd$ is unique.

\smallskip

To prove the ``only if'' part, assume that $X$ is an affine, $3$-ranked liner such that for any triangle $abc$ in $X$, there exists a unique point $d\in X$ such that $abcd$ is a parallelogram. Then the triangle $abc$ is contained in the plane $\overline{\{a,b,c,d\}}$ of cardinality $\ge 4$, which implies $|X|_2\ge 3$ (by  Theorem~\ref{t:affine=>Avogadro}, the affine liner $X$ is $2$-balanced, so the cardinal $|X|_2$ is well-defined). If $|X|_2\ge 4$, then by Theorem~\ref{t:4-long-affine}, the $4$-long affine liner $X$ is regular and by Theorem~\ref{t:Playfair<=>}, the $4$-long regular affine liner is Playfair. So, we assume that $|X|_2=3$, which means that $X$ is a Steiner liner. 

We claim that the midpoint operation on the Steiner liner $X$ is self-distributive. Given any points $o,x,y\in X$, we need to show that the point $u\defeq o\circ(x\circ y)$ is equal to the point $v\defeq (o\circ x)\circ(o\circ y)$. By our assumption, there exists a unique point $p\in X$ such that $xoyp$ is a parallelogram. 

Let us consider 5 possible cases.
\smallskip

1. If $x=y$, then by the idempotence of the midpoint operation, we obtain
$$o\circ (x\circ y)=o\circ (x\circ x)=o\circ x=(o\circ x)\circ (o\circ x)=(o\circ x)\circ(o\circ y).$$
\smallskip

2. If $o=x$, then by the idempotence of the midpoint operation, we obtain
$$o\circ(x\circ y)=(o\circ o)\circ(x\circ y)=(o\circ x)\circ (o\circ y).$$
\smallskip

3. If $o=y$, then by the idempotence and the commutativity of the midpoint operation we obtain 
$o\circ(x\circ y)=(x\circ y)\circ o=(y\circ x)\circ(o\circ o)=(o\circ x)\circ (o\circ y).$
\smallskip

4. If $|\{o,x,y\}|=3$ and $\|\{o,x,y\}\|=2$, then the set $\{x,o,y\}$ is a line in the Steiner liner $X$ and hence $o\circ (x\circ y)=o\circ o=o=y\circ x=(o\circ x)\circ (o\circ y)$. 
\smallskip

5. So, assume that $\|\{o,x,y\}\|=3$, which means that $xoy$ is a triangle. We recall that $u=o\circ(x\circ y)$ and $v=(o\circ x)\circ(o\circ y)$.

\begin{claim}\label{cl:u=p} The quadruple $xoyu$ is a paralelogram and hence $u=p$.
\end{claim}

\begin{proof} First observe that $x\circ y=o\circ u\in\Aline xy\cap\Aline ou$. Since $u=o\circ(x\circ y)\in \overline{\{o,x\circ y\}}\subseteq\overline{\{o,\overline{\{x,y\}}}=\overline{\{0,x,y\}}$, the set $\{o,x,y,u\}$ has rank $3$.
Consider the point $z\defeq x\circ y=o\circ u\in\Aline xy\cap\Aline ou$ and observe that $z\ne o\circ z=u$ and the line $\Aline xy$ coincides with the set $\{x,y,z\}$. 
Assuming that $u\in\Aline ox$, we conclude that $z\in \Aline ou\subseteq \Aline ox$ and $y\in\Aline  xz\subseteq\Aline ox$, which contradicts $\|\{x,o,y\}\|=3$. This contradiction shows that $u\notin \Aline ox$. 

Since the liner $X$ is affine, for the point $u\in \Aline zo\setminus \Aline xo$, there exists a point $w\in\Aline zx$ such that $\Aline wu\cap\Aline xo=\varnothing$. It follows that $w\in\Aline zx\setminus\{x,z\}=\{y\}$ and hence $\Aline yu\cap\Aline xo=\Aline wu\cap\Aline xo=\varnothing$.

By analogy we can prove that $\Aline xu\cap\Aline yo=\varnothing$. Since $z\in\Aline xy\cap\Aline ou$, the quadruple $xoyu$ is a parallelogram in $X$, by Proposition~\ref{p:tobe-parallelogram}. The uniqueness of the point $p$ ensures that $u=p$. 
\end{proof}

\begin{claim}\label{cl:v=p} The quadruple $xoyv$ is a paralellogram and hence $v=p$.
\end{claim}

\begin{proof} 
Consider the points $x_o\defeq o\circ x=x\circ o$ and $y_o\defeq o\circ y=y\circ o$, and observe that  $y\in \Aline {y_o}o\setminus \Aline xo=\Aline{y_o}o\setminus \Aline {x_o}o$. Since the liner $X$ is affine, there exists a point $w\in\Aline {x_o}{y_o}$ such that $\Aline wy\cap \Aline {x_o}o=\varnothing$. Taking into account that $w\in\Aline{x_o}{y_o}\setminus\{x_o,y_o\}=\{v\}$, we conclude that $w=v$ and hence $\Aline vy\cap\Aline xo=\Aline wy\cap\Aline{x_o}o=\varnothing$. By analogy we can prove that $\Aline vx\cap\Aline yo=\varnothing$. Since $v=(o\circ x)\circ(o\circ y)$, the set $\{x,o,y,v\}$ has rank at most $3$. By Proposition~\ref{p:tobe-parallelogram},  $xoyv$ is a parallelogram in $X$. The uniqueness of the point $p$ ensures that $v=p$.
\end{proof}

By Claims~\ref{cl:u=p} and \ref{cl:v=p}, 
$o\circ(x\circ y)=u=p=v=(o\circ x)\circ(o\circ y),$ witnessing that the midpoint operation on the Steiner liner $X$ is self-distributive. By Theorem~\ref{t:Hall<=>Moufang}, the Steiner liner $X$ is Hall, and by Theorem~\ref{t:Hall<=>Playfair+Steiner}, the Hall liner $X$ is Playfair.
\end{proof} 

\section{Boolean parallelograms and Boolean liners}

\begin{definition} A parallelogram $abcd$ in a liner $X$ is called \index{Boolean parallelogram} \index{parallelogram!Boolean}\defterm{Boolean} if its diagonals $\Aline ac$ and $\Aline bd$ are parallel. 
\end{definition}

For $3$-ranked Steiner affine liners, Proposition~\ref{p:k-parallel=>parallelogram} and Theorem~\ref{t:parallelogram3+1} are completed by the following proposition.

\begin{proposition}\label{p:Steiner-paragram} For any triangle $abc$ in a $3$-ranked Steiner affine liner $X$, there exists a unique point $d\in X$ such that $abcd$ is a non-Boolean parallelogram.
\end{proposition}

\begin{proof}  Since $X$ is Steiner, there exist unique points $o\in \Aline ac\setminus\{a,c\}$ and $d\in\Aline ob\setminus\{o,b\}$. Since $X$ is affine, there exists a point $x\in \Aline ob$ such that $\Aline xc\cap\Aline ab=\varnothing$. Taking into account that the lines $\Aline oc$ and $\Aline bc$ intersect $\Aline ab$, we conclude that $x=d$ and hence $\Aline ab\cap\Aline cd=\Aline ab\cap\Aline cx=\varnothing$.

\begin{picture}(300,65)(-140,10)
{\linethickness{0.8pt}
\put(15,20){\color{red}\line(1,0){50}}
\put(25,60){\color{red}\line(1,0){50}}
\put(15,20){\color{cyan}\line(1,4){10}}
\put(65,20){\color{cyan}\line(1,4){10}}
}
\put(15,20){\line(3,2){60}}
\put(65,20){\line(-1,1){40}}
\put(15,20){\circle*{3}}
\put(8,12){$b$}
\put(65,20){\circle*{3}}
\put(66,12){$a$}
\put(25,60){\circle*{3}}
\put(21,63){$c$}
\put(75,60){\circle*{3}}
\put(71,63){$x$}
\put(78,58){$d$}
\put(45,40){\circle*{3}}
\put(43,43){$o$}

\end{picture}

By analogy we can show that $\Aline ad\cap \Aline bc=\varnothing$.   By Proposition~\ref{p:tobe-parallelogram}, $abcd$ is a parallelogram. Since $\Aline ac\cap\Aline bd=\{o\}$, the parallelogram $abcd$ is not Boolean. 

Now assume that $\delta$ is another point of $X$ such that $abc\delta$ is a non-Boolean parallelogram. Then the unique common point of its diagonals $\Aline ac$ and $\Aline b\delta$ coincides with the unique point $o$ of the line $\Aline ac$, which is distinct from $a$ and $c$. Then $\delta=\Aline bo\setminus\{b,o\}=\{b,o,d\}\setminus\{b,o\}=\{d\}$, witnessing that the point $d$ is unique.
\end{proof}

\begin{corollary} For every triangle $abc$ in a $3$-long $3$-ranked affine liner $X$, there exists a point $d\in X$ such that $abcd$ is a parallelogram.
\end{corollary} 
 
\begin{proof} If $|X|_2=3$, then the point $d$ exists, by Proposition~\ref{p:Steiner-paragram}. So, assume that $|X|_2\ge 4$. By Theorem~\ref{t:4-long-affine}, the $4$-long affine liner $X$ is Playfair. By Theorem~\ref{t:parallelogram3+1}, there exists a point $d$ such that $abcd$ is a parallelogram.
\end{proof}

\begin{corollary}\label{c:Steiner-Playfair-Boolean} A Steiner liner is Playfair if and only if it is affine, $3$-ranked, and contains no Boolean parallelograms.
\end{corollary}

\begin{proof} Assume that a Steiner liner $X$ is Playfair. By Theorem~\ref{t:Playfair<=>} and Proposition~\ref{p:k-regular<=>2ex}, the Playfair liner $X$ is affine and $3$-ranked. To derive a contradiction, assume that $X$ contains a Boolean parallegram $abcd$. By Proposition~\ref{p:Steiner-paragram}, there exists a point $\delta\in X$ such that $abc\delta$ is a non-Boolean parallelogram. By Theorem~\ref{t:parallelogram3+1}, the triangle $abc$ can be completed to a unique parallelogram, which implies $d=\delta$. Then the parallelogram $abcd=abc\delta$ is simultaneously Boolean and non-Boolean, which is a contradiction showing that the Steiner Playfair plane contains no Boolean parallelograms.
\smallskip

Now assume that the Steiner liner $X$ is affine, $3$-ranked, and contains no Boolean parallelograms. Assuming that $X$ is not Playfair, we can apply Theorem~\ref{t:parallelogram3+1} and find two distinct points $d,d'$ such that $abcd$ and $abcd'$ are parallelograms. Since $X$ contains no Boolean parallelograms, those two parallelograms are not Boolean. Then $d,d'$ are two distinct points such that $abcd$ and $abcd'$ are two distinct non-Boolean parallelograms, which contradicts Proposition~\ref{p:Steiner-paragram}.
\end{proof}

Let us recall that a liner is \index{line-finite liner}\index{liner!line-finite}\defterm{line-finite} if all its lines are finite sets. Observe that a $2$-balanced liner $X$ is line-finite if and only if the cardinal $|X|_2$ is finite.

\begin{proposition}\label{p:large=>Boolean-parallelogram} If a line-finite $2$-balanced $3$-balanced $3$-ranked liner $X$ has $|X|_3\ge |X|_2^2$ \textup{(}and $|X|_3>3{\cdot}|X|_2^2-9{\cdot}|X|_2+9$\textup{)}, then for every triangle $abc$ there exists a point $d\in X$ such that $abcd$ is a \textup{(}Boolean\textup{)} parallelogram.
\end{proposition}

\begin{proof} Given any triangle $abc$ in $X$, consider the plane $P\defeq\overline{\{a,b,c\}}$ and observe that $|P|=|X|_3$.
\smallskip

 First assume that $|X|_3>3{\cdot}|X|_2^2-9{\cdot}|X|_2+9$.
Observe that the set $$D\defeq \bigcup_{x\in \Aline ac}\Aline bx\cup\bigcup_{y\in \Aline ab}\Aline cy\cup\bigcup_{z\in \Aline cb}\Aline az$$is contained in the plane $P$ and has cardinality $$
\begin{aligned}
|D|&=|\{a,b,c\}|+|(\Aline ab\cup\Aline bc\cup\Aline ac)\setminus\{a,b,c\}|\\
&+\sum_{x\in \Aline ac\setminus\{a,c\}}|\Aline bx\setminus\{b,x\}|+\sum_{y\in \Aline ab\setminus\{a,b\}}|\Aline cy\setminus\{c,y\}|+\sum_{z\in \Aline bc\setminus\{b,c\}}|\Aline az\setminus\{a,z\}|\\
&=3+3{\cdot}(|X|_2-2)+3{\cdot}(|X|_2-2)^2=3{\cdot}|X|_2^2-9{\cdot}|X|_2+9<|X|_3=|P|.
\end{aligned}
$$ So we can choose a point $d\in P\setminus D$ and observe that the intersections $\Aline da\cap\Aline bc$, $\Aline dc\cap\Aline ab$ and $\Aline db\cap\Aline ac$ are empty. The $3$-rankedness of $X$ ensures that $abcd$ is a Boolean parallelogram in the plane $P$, by Proposition~\ref{p:tobe-parallelogram}.
\smallskip

Next, assume that $|X|_2^2\le |X|_3\le 3{\cdot}|X|_2^2-9{\cdot}|X|_2+9$. In this case the cardinal $|X|_3=|P|$ is finite. Let $\mathcal L$ be the family of all lines in $P$. For any point $x\in P$ and line $\Lambda$ in $P$, consider the families $$\mathcal L_x\defeq\{L\in\mathcal L:x\in L\},\quad\mathcal L^\bullet_{x,\Lambda}\defeq\{L\in\mathcal L_x:L\cap\Lambda\ne\varnothing\}\quad\mbox{and}\quad\mathcal L^\circ_{x,\Lambda}\defeq\{L\in\mathcal L_x:L\cap\Lambda=\varnothing\}$$ 
It is easy to see that for any line $\Lambda\in\mathcal L$ and point $x\in X\setminus\Lambda$, the set $\mathcal L_{x,\Lambda}^\bullet$ has cardinality $|\mathcal L_{x,\Lambda}^\bullet|=|X|_2$ and hence $|\mathcal L_{x,\Lambda}^\circ|=|\mathcal L_x|-|X|_2$.

Theorem~\ref{t:2-balanced-dependence} ensures $|X|_3=|P|=1+(|X|_2-1)\cdot\lambda$, where $\lambda=|\mathcal L_x|$ for any point $x\in X$. The inequality $|X|\ge|X|_2^2$ ensures that $\lambda>|X|_2$, and hence $|\mathcal L_{x,\Lambda}^\circ|=|\mathcal L_x|-|\mathcal L_{x,\Lambda}^\bullet|=\lambda-|X|_2>0$ for any point $x\in X$ and line $\Lambda\in \mathcal L\setminus\mathcal L_x$. 
 
Consider the lines $H\defeq \Aline ab$ and $V\defeq\Aline bc$ in the plane $P$. Since $|\mathcal L^\circ_{c,H}|=\lambda-|X|_2>0$, there exists a line $H_c\in\mathcal L^\circ_{c,H}$. Since $H\in \mathcal L^\circ_{b,H'}\setminus \mathcal L^\circ_{b,V}$ and $|\mathcal L^\circ_{b,H'}|=|\mathcal L^\circ_{b,V}|=\lambda-|X|_2$, there exists a line $V'\in \mathcal L^\circ_{b,V}\setminus\mathcal L^\circ_{b,H'}$, which has a common point $d$ with the line $H'$. Then $\Aline ab\cap\Aline cd=H\cap H'=\varnothing=V\cap V'=\Aline ac\cap \Aline bd$ and hence $abcd$ is a parallelogram in the $3$-ranked liner $X$, according to Proposition~\ref{p:tobe-parallelogram}.  
\end{proof}

\begin{corollary} For every triangle $abc$ in a Steiner ranked plane $X$ of cardinality $|X|>9$, there exists a point $d\in X$ such that $abcd$ is a Boolean parallelogram.
\end{corollary}

\begin{definition}
A liner $X$ is defined to be \index{Boolean liner}\index{liner!Boolean}\defterm{Boolean} if for any quadruple $abcd\in X^4$, $\Aline ab\cap\Aline cd=\varnothing=\Aline bc\cap\Aline ad$ implies $\Aline ac\cap\Aline bd=\varnothing$.
\end{definition}

\begin{theorem}\label{t:Boolean<=>} A $3$-ranked liner $X$ is Boolean if and only if every parallelogram in $X$ is Boolean.
\end{theorem}

\begin{proof} Assume that $X$ is Boolean and take any parallelogram $abcd$ in $X$. Then $\Aline ab\cap\Aline cd=\varnothing=\Aline bc\cap\Aline ad$. Since $X$ is Boolean, the latter condition implies $\Aline ac\cap\Aline bd=\varnothing$. Taking into account that $\Aline ac$ and $\Aline bd$ are disjoint lines with $\|\Aline ac\cup\Aline bd\|=\{a,b,c,d\}\|\le 3$, we can apply Corolary~\ref{c:parallel-lines<=>} and conclude that $\Aline ac\parallel \Aline bd$, witnessing that the parallelogram $abcd$ is Boolean.

Now assume that every parallelogram in $X$ is Boolean. Given any points $a,b,c,d\in X$ with $\Aline ab\cap\Aline cd=\varnothing=\Aline bc\cap\Aline ad$, we shall prove that $\Aline ac\cap\Aline bd=\varnothing$. To derive a contradiction, assume that $\Aline ac\cap\Aline bd$ contains some point $o$. Then $|\{a,b,c,d\}\|=\|\{a,b,c\}\|\le 3$. Applying Proposition~\ref{p:tobe-parallelogram}, we conclude that $abcd$ is a parallelogram. Since $X$ is Boolean, this parallelogram has parallel diagonals $\Aline ac\parallel\Aline bd$. Assuming that $\Aline ac\cap\Aline bd\ne\varnothing$, we conclude that $\Aline ac=\Aline bd$ and hence the points $a,b,c,d$ are collinear. In this case the equality $\Aline ab\cap\Aline cd=\varnothing$ implies $a=b$ and $c=d$ and hence $\Aline bc=\Aline ad$, which contradicts $\Aline bc\cap\Aline ad=\varnothing$. This contradiction shows that $\Aline ac\cap\Aline bd=\varnothing$, witnessing that the liner $X$ is Boolean.
\end{proof}

\begin{exercise} Let $V$ be a vector space over a field $F$ of characteristic two. Show that $V$ endowed with its canonical line relation is a Boolean liner.
\end{exercise}




\section{Hyper-Bolyai liners}

Motivated by the problem of extending triangles to parallelograms, we introduce and study in this section the class of hyper-Bolyai liners, intermediate between the Playfair and Bolyai liners. 

Let us recall that a liner $X$ is \index{Bolyai liner}\index{liner!Bolyai}\defterm{Bolyai} if for any plane $P\subseteq X$, line $L\subseteq P$ and point $p\in P\setminus L$ there exists a line $\Lambda$ such that $p\in \Lambda\subseteq P\setminus L$.

\begin{definition}\label{d:hyper-Bolyai} A liner $X$ is defined to be \index{hyper-Bolyai liner}\index{liner!hyper-Bolyai}\defterm{hyper-Bolyai} if 
for any plane $P\subseteq X$, concurrent lines $A,B\subseteq P$ and point $p\in P\setminus B$, there exists a line $\Lambda\subseteq X$ such that $p\in \Lambda$ and $|\Lambda\cap A|=1$ and $\Lambda\parallel B$.
\end{definition}

\begin{picture}(100,80)(-150,-25)
\linethickness{0.7pt}
\put(0,0){\color{teal}\line(1,0){70}}
\put(75,-3){\color{teal}$B$}
\put(0,30){\color{teal}\line(1,0){70}}
\put(75,27){\color{teal}$\Lambda$}
\put(20,-20){\color{cyan}\line(0,1){70}}
\put(10,50){\color{cyan}$A$}
\put(20,30){\color{red}\circle*{3}}
\put(50,30){\circle*{3}}
\put(48,22){$p$}
\put(20,0){\circle*{3}}
\end{picture}

Hyper-Bolyai liners relate to some other types of liners as follows.
$$\xymatrix{
\mbox{Proclus}\ar@{=>}[d]&\mbox{Playfair}\ar@{=>}[l]\ar@{=>}[r]\ar@{=>}[d]&\mbox{hyper-Bolyai}\ar@{=>}[r]\ar@{=>}[d]&\mbox{Bolyai}\\
\mbox{proaffine}&\mbox{affine}\ar@{=>}[l]\ar@{=>}[r]&\mbox{hyperaffine}\ar_{\tiny\mbox{$3$-long}}[ru]
}
$$
Non-trivial implications in this diagram are proved in Theorem~\ref{t:Proclus<=>}, \ref{t:Playfair<=>}, and the following proposition.

\begin{proposition}\label{p:hyper-Bolyai-interplay}
\begin{enumerate}
\item Every hyper-Bolyai liner is Bolyai.
\item Every hyper-Bolyai liner is hyperaffine.
\item Every hyperaffine $3$-long is Bolyai.
\item A liner is Playfair if and only if it is Proclus and hyper-Bolyai.
\end{enumerate}
\end{proposition}

\begin{proof} 1. Definitions~\ref{d:PPBL} and \ref{d:hyper-Bolyai} imply that every hyper-Bolyai liner is Bolyai.
\smallskip

2. Assume that a liner $X$ is hyper-Bolyai.  To prove that $X$ is hyperaffine, take any points $o,x,y\in X$ and $p\in \Aline xy\setminus\Aline ox$. We have to find a point $u\in \Aline oy$ such that $\Aline up\cap \Aline ox=\varnothing$. If $p\in \Aline oy$, then the point $u\defeq p\in \Aline oy$ has the required property. So, assume that $p\notin\Aline oy$, which implies $\Aline oy\cap\Aline xy=\{y\}$ and $o\ne x$.  Consider the concurrent lines $A\defeq\Aline oy$, $B\defeq\Aline ox$, and the plane $P\defeq\overline{A\cup B}$. Observe that  $p\in \Aline xy\setminus\{x\}\subseteq P\setminus B$. Since $X$ is hyper-Bolyai, there exists a line $L$ in $X$ such that $p\in L$, $L\parallel B$, and $L\cap A=\{u\}$ for some point $u\in A=\Aline oy$. It follows from $p\in L\setminus B$ and $L\parallel B$ that $L\cap B=\varnothing$ and hence $\Aline up\cap B=L\cap B=\varnothing$, witnessing that the liner $X$ is hyperaffine.
\smallskip

3. Assume that the liner $X$ is $3$-long and hyperaffine. To prove that $X$ is Bolyai, take any plane $P\subseteq X$, line $L\subset P$, and point $p\in P\setminus L$. Choose any distinct points $o,x\in L$. Since $X$ is $3$-long, there exists a point $y\in \Aline xp\setminus\{x,p\}$. Since $X$ is hyperaffine, there exists a point $u\in \Aline oy$ such that $\Aline up\cap \Aline ox=\varnothing$. Then the line $\Lambda\defeq \Aline up$ contains the point $p$ and does not intersect the line $L=\Aline ox$, witnessing that the liner $X$ is Bolyai.
\smallskip

4. If $X$ is Playfair, then $X$ is Proclus, by Definition~\ref{d:PPBL}. To show that $X$ is hyper-Bolyai, take any plane $P\subseteq X$, concurrent lines $A,B\subset P$ and point $p\in P\setminus B$. Since $X$ is Playfair, there exists a line $L$ such that $p\in L\subseteq P\setminus B$. Since $|A\cap B|=1$ and $L\cap B=\varnothing$, the Proclus Postulate~\ref{p:Proclus-Postulate} ensures that $|L\cap A|=1$, witnessing that the liner $X$ is hyper-Bolyai.

Now assume conversely that the liner $X$ is Proclus and hyper-Bolyai. Then $X$ is Bolyai, by Proposition~\ref{p:hyper-Bolyai-interplay}(1). By Definition~\ref{d:PPBL}, the Proclus Bolyai liner $X$ is Playfair.
\end{proof}

\begin{remark} By Example~\ref{ex:non-hyperaffine}, there exists a $4$-long Bolyai plane, which is not hyperaffine. By Example~\ref{ex:hyperaffine-nonPlayfair}, there exists a  non-Playfair hyperaffine $4$-long ranked plane, which is hyper-Bolyai, by Theorem~\ref{t:hyper-Bolyai<=>}.
\end{remark}

In the proof of Theorem~\ref{t:hyper-Bolyai<=>}, we shall apply the following characterization of hyperaffine liners.

\begin{proposition}\label{p:hyperaffine<=>} A liner $X$ is hyperaffine if and only if for any concurrent lines $A,B$ in $X$ and point $p\in \overline{A\cup B}\setminus (A\cup B)$, there exists a point $a\in A$ such that $\Aline ap\cap B=\varnothing$.
\end{proposition}

\begin{picture}(100,80)(-170,-20)
\linethickness{0.7pt}
\put(0,0){\color{teal}\line(1,0){70}}
\put(75,-3){\color{teal}$B$}
\put(0,30){\color{teal}\line(1,0){70}}
\put(20,-20){\color{cyan}\line(0,1){75}}
\put(10,50){\color{cyan}$A$}
\put(20,30){\color{red}\circle*{3}}
\put(22,32){$a$}
\put(50,30){\circle*{3}}
\put(48,22){$p$}
\put(20,0){\circle*{3}}
\put(13,-7){$o$}
\end{picture}

\begin{proof} To prove the ``only if'' part, assume that $X$ is hyperaffine. Given any concurrent lines $A,B\subset X$ and point $p\in \overline{A\cup B}\setminus(A\cup B)$, we should find a point $a\in A$ such that $\Aline ap\cap B=\varnothing$. Let $o$ be the unique common point of the concurrent lines $A$ and $B$. Choose any point $y\in A\setminus\{o\}$. If $\Aline yp\cap B=\varnothing$, then the point $a\defeq y$ has the required property. So, assume that $\Aline yp\cap B\ne\varnothing$ and choose a point $b\in \Aline yp\cap B$. Then $p\in \Aline yb\setminus\Aline ob$ and by the hyperaffinity of $X$, there exists a point $a\in \Aline oy=A$ such that $\varnothing =\Aline ap\cap\Aline ob=\Aline ap\cap B$.

Now we prove the ``only if'' part. Assume that for any concurrent lines $A,B$ in $X$ and point $p\in \overline{A\cup B}\setminus (A\cup B)$, there exists a point $a\in A$ such that $\Aline ap\cap B=\varnothing$. Given any points $o,x,y\in X$ and $p\in \Aline yx\setminus\Aline ox$, we should find a point $u\in \Aline oy$ such that $\Aline up\cap\Aline ox=\varnothing$.

If $p\in \Aline oy$, then the point $u\defeq p\in \Aline oy$ has the required property $\Aline up\cap \Aline ox=\{p\}\cap\Aline ox=\varnothing$. So, assume that $\notin \Aline oy$. Then $A\defeq \Aline oy$ and $B\defeq\Aline ox$ are concurrent lines and $p\in \Aline xy\setminus(\Aline ox\cup\Aline oy)\subseteq\overline{A\cup B}\setminus(A\cup B)$. By the assumption, there exists a point $u\in A=\Aline oy$ such that $\varnothing =\Aline up\cap B=\Aline up\cap\Aline ox$, witnessing that the liner $X$ is hyperaffine.
\end{proof}

Now we are in a position to present a characterization of hyper-Bolyai liners, in the spirit of the characterizations of Proclus and Playfair liners given in Theorems \ref{t:Proclus<=>} and \ref{t:Playfair<=>}.

\begin{theorem}\label{t:hyper-Bolyai<=>} A liner $X$ of cardinality $|X|>2$ is hyper-Bolyai if and only if $X$ is $3$-long, hyperaffine, and $3$-ranked.
\end{theorem}

\begin{proof} To prove the ``if'' part, assume that $X$ is a $3$-long hyperaffine $3$-ranked liner. Given any plane $P\subseteq X$, concurrent lines $A,B\subseteq P$ and point $p\in P\setminus B$, we should find a line $L$ in $X$ such that $p\in L$, $|L\cap A|=1$ and $L\parallel B$. By the $3$-rankedness of $X$, the plane $P$ coincides with the plane $\overline{A\cup B}$. Let $o$ be the unique common point of the lines $A,B$. If $p\notin A$, then 
$p\in P\setminus(A\cup B)=\overline{A\cup B}\setminus(A\cup B)$. Applying Proposition~\ref{p:hyperaffine<=>} to the hyperaffine liner $X$, find a point $a\in A$ such that $\Aline ap\cap B=\varnothing$. Consider the line $L\defeq\Aline ap$ and observe that $L\cap A=\{a\}$. Since $X$ is $3$-ranked and $\|L\cup B\|=\|\{a,p\}\cup B\|\le\|A\cup B\|=3$, Corollary~\ref{c:parallel-lines<=>} ensures that the disjoint lines $L$ and $B$ in the plane $\overline{A\cup B}$ are parallel.  
If $p\in A$, then take any points $x\in B\setminus\{o\}$ and $y\in \Aline px\setminus\{p,x\}$. The point $y$ exists because the liner $X$ is $3$-long. Since $p\in \Aline xy\setminus B=\Aline xy\setminus \Aline ox$, by the hyperaffinity of $X$, there exists a point $a\in \Aline oy$ such that $\Aline ap\cap B=\varnothing$. Since $\Aline oy\cap\Aline xy=\{y\}\ne\{p\}$, the flat $\Lambda\defeq\Aline ap$ is a line such that $\Lambda\cap A=\{a\}$ and $\Lambda\parallel B$, by the $3$-rankeness of $X$, see Corollary~\ref{c:parallel-lines<=>}. This completes the proof of the ``if'' part.
\smallskip

To prove the ``only if'' part, assume that $X$ is a hyper-Bolyai liner with $|X|>2$. By Proposition~\ref{p:hyper-Bolyai-interplay}, the hyper-Bolyai liner $X$ is hyperaffine. It remains to prove that $X$ is $3$-long and $3$-ranked.  If $\|X\|\le 2$, then $|X|>2$ implies that $X$ is $3$-long. In this case, $X$ is hyperaffine and $3$-ranked, because $X$ contains no triangles and no planes. So, assume that $\|X\|\ge 3$.  

To prove that $X$ is $3$-long, we have to show that every line $A$ in $X$ contains at least three points. Choose any points $o,a\in A$. Since $\|X\|\ge 3$, there exists a point $b\in X\setminus A$. Then the line $B\defeq\Aline ob$ is concurrent to the line $A$. Moreover, $P\defeq\overline{A\cup B}=\overline{\{a,o,b\}}$ is a plane containing the concurrent lines $A,B$. Since $X$ is hyper-Bolyai, for the point $p\defeq a\in P$, there exists a line $L\subseteq X$ such that $p\in L$ and $L\parallel B$. Then $L\subseteq\overline{B\cup \{a\}}\setminus B$, which implies $\overline{B\cup\{a\}}\ne B\cup\{a\}$ and hence $\Aline ax\ne\{a,x\}$ for some $x\in B$. Then we can choose a point $p\in \Aline xa\setminus\{x,a\}$. If $x=o$, then $o,p,a$ are three distinct points on the line $A=\Aline oa=\Aline xa$ and we are done. So assume that $x\ne o$. By the  hyperaffinity of $X$, there exists a point $u\in \Aline oy=A$ such that $\Aline up\cap B=\varnothing$. Since $\Aline py\cap B=\{x\}$ and $\Aline op\cap B=\{o\}$, the point $u$ is distinct from the points $o$ and $y$, which implies that $|A|\ge|\{o,u,y\}|=3$ and witnesses that the liner $X$ is $3$-long.

Assuming that the liner $X$ is not $3$-ranked, we can find two distinct planes $P\subset \Pi$. Choose any point $p\in\Pi\setminus P$ and fix any concurrent lines $A,B$ in the plane $P$. By the assumption, there exists a line $L$ in $X$ such that $p\in L$, $L\parallel B$ and $L\cap A=\{a\}$ for some point $a$. The parallelity relation $L\parallel B$ ensures that $p\in L\subseteq \overline{B\cup\{a\}}\subseteq P$, which contradicts the choice of the point $p$. This contradiction shows that the liner $X$ is $3$-ranked.
\end{proof}

The following proposition (which motivated our study of hyper-Bolyai liners) show that triangles in hyper-Bolyai liners can be extended to parallelograms.

\begin{proposition}\label{p:hyper-Bolyai=>parallelogram} For every triangle $abc$ in a hyper-Bolyai liner $X$, there exists a point $d\in X$ such that $abcd$ is a parallelogram.
\end{proposition}

\begin{proof} Given any triangle $abc$ in $X$, consider the plane $P\defeq\overline{\{a,b,c\}}$. Applying Definition~\ref{d:hyper-Bolyai} to the concurrent lines $\Aline ba$ and $\Aline bc$ and the point $a$ in the plane $P$, find a line $A$ in $X$ such that $a\in A$ and $A\parallel \Aline bc$. Then $A\subseteq \overline{\{a,b,c\}}=P$.  Applying Definition\ref{d:hyper-Bolyai} to the concurrent lines $A,\Aline ab$ and the point $c$ in the plane $P$, find a line $C$ in $X$ such that $c\in C$, $C\parallel \Aline ba$ and $C\cap A=\{d\}$ for some point $d$.

\begin{picture}(60,70)(-150,-15)
\linethickness{0.7pt}

\put(0,0){\color{teal}\line(1,0){50}}
\put(0,0){\color{blue}\line(0,1){40}}
\put(50,0){\color{blue}\line(0,1){40}}
\put(53,16){\color{blue}$C$}
\put(0,40){\color{teal}\line(1,0){50}}
\put(22,43){\color{teal}$A$}

\put(0,0){\circle*{3}}
\put(-7,-7){$b$}
\put(50,0){\circle*{3}}
\put(50,-7){$c$}
\put(0,40){\circle*{3}}
\put(-8,39){$a$}
\put(50,40){\color{red}\circle*{3}}
\put(53,39){$d$}

\end{picture}

It follows from $a\in A\setminus\Aline bc$ and $c\in C\setminus \Aline ba$ that $A\parallel \Aline bc\ne A$ and $C\parallel\Aline ab\ne C$, witnessing that $abcd$ is a parallelogram. 
\end{proof}

\begin{remark} By Example~\ref{ex:non-hyperaffine}, there exists a $4$-parallel plane $X$ with $|X|_2=4$ and $|X|=|X|_3=25$, which is not hyperaffine and hence not hyper-Bolyai.  By Proposition~\ref{p:k-parallel=>parallelogram}, every triangle in $X$ can be completed to a parallelogram. This example shows that Proposition~\ref{p:hyper-Bolyai=>parallelogram} cannot be turned into a characterization of hyper-Bolyai liners.
\end{remark}

\chapter{Spread completions of liners}

In this chapter we define the construction of the spread completion of a $3$-ranked liners and will characterize $3$-ranked liners whose spread completions are $4$-long and projective. This characterization includes notions of a para-Playfair liner and bi-Bolyai liner, studied in Sections~\ref{s:para-Playfair} and \ref{s:bi-Bolyai}. 

\section{Para-Playfair liners}\label{s:para-Playfair}

\begin{definition}\label{d:para-Playfair} A liner $X$ is defined to be \index{para-Playfair liner}\index{liner!para-Playfair}\defterm{para-Playfair} if for every plane $P\subseteq X$, disjoint lines $L,L'\subseteq P$ and a point $x\in P\setminus L$, there exists a unique line $L_x$ such that $x\in L_x\subseteq P\setminus L$.
\end{definition} 

\begin{proposition}\label{p:Playfair=>para-Playfair=>Proclus} Every Playfair liner is para-Playfair and every para-Playfair liner is Proclus.
\end{proposition}

\begin{proof} Definitions~\ref{d:PPBL} and \ref{d:para-Playfair} imply that every Playfair liner is para-Playfair. Assuming that some para-Playfair liner $X$ is not Proclus, we can find a plane $\Pi\subseteq X$, a line $L\subseteq \Pi$, a point $x\in \Pi\setminus L$ and two distinct lines $L',L''\subseteq \Pi\setminus\Lambda$ containing the point $x$. Since $L\cap L'=\varnothing$, the para-Playfair property of $X$ ensures that $L'$ is a unique line in the plane that contains $x$ and is disjoint with $L$. Therefore, $L''=L'$, which contradicts the choice of the (distinct) lines $L',L''$ and shows that the para-Playfair liner $X$ is Proclus.
\end{proof}

Definitions~\ref{d:PPBL}, \ref{d:PP}, \ref{d:para-Playfair}, Theorems~\ref{t:projective<=>}, \ref{t:Proclus<=>}, \ref{t:Playfair<=>}, and Proposition~\ref{p:Playfair=>para-Playfair=>Proclus} imply that for any liner the  following implications hold:
$$\xymatrix{
\mbox{strongly regular}\ar@{=>}[r]\ar@{<=>}[d]&\mbox{\color{red}para-Playfair}\ar@{=>}[d]&&\mbox{Playfair}\ar@{=>}[ll]\ar@{=>}[d]\\
\mbox{projective}\ar@{=>}[r]&\mbox{Proclus}\ar@{=>}[r]&\mbox{proaffine}&\mbox{affine}\ar@{=>}[l].
}
$$

\begin{theorem}\label{t:4-para-Playfair} Every proaffine regular liner $X$ of rank $\|X\|\ne 3$ is para-Playfair. 
\end{theorem}

\begin{proof} By Theorems~\ref{t:Proclus<=>} and \ref{t:w-modular<=>}, the proaffine regular liner $X$ is Proclus and weakly modular. If $\|X\|\le 2$, then $X$ contains no planes and hence is para-Playfair vacuously. So, assume that $\|X\|\ge 4$.

\begin{claim}\label{cl:parallel-escorted} For any disjoint coplanar lines $L,L'$ in $X$ and  point $x\in X\setminus \overline{L\cup L'}$, there exists a line $L_x$ in $X$ such that $x\in L_x$ and $L_x\parallel L$.
\end{claim}

\begin{proof} Since the lines $L,L'$ are disjoint and coplanar,  the flat hull $\Pi\defeq\overline{L\cup L'}$ is a plane in $X$. Since $x\notin \Pi$, the $3$-rankedness of $X$ ensures that $\overline{L'\cup\{x\}}\cap \Pi=L'$. By Theorem~\ref{t:w-modular<=>}, the intersection $L_x\defeq\overline{L\cup\{x\}}\cap \overline{L'\cup\{x\}}$ is a line containing $x$. Observe that 
$$
L_x\cap L=(\overline{L\cup\{x\}}\cap\overline{L'\cup\{x\}})\cap L\subseteq \overline{L'\cup\{x\}}\cap(\Pi\cap L)=(\overline{L'\cup\{x\}})\cap \Pi)\cap L=L'\cap L=\varnothing.$$ By Corollary~\ref{c:parallel-lines<=>}, the disjoint coplanar lines $L$ and $L_x$ are parallel.
\end{proof} 

Now we are able to prove that the liner $X$ is para-Playfair. Given a plane $P\subseteq X$, two disjoint lines $L,L'\subseteq P$ and a point $x\in P\setminus L$, we have to find a unique line $L_x\subseteq X$ such that $x\in L_x\subseteq P\setminus L$. 
Since $\|X\|\ge 4$, there exists a point $y\in X\setminus P$. By Corollary~\ref{c:parallel-lines<=>}, the disjoint coplanar lines $L$ and $L'$ are parallel. By Claim~\ref{cl:parallel-escorted}, there exists a line $L_y$ in $X$ such that $y\in L_y$ and $L_y\parallel L$. By Theorem~\ref{t:Proclus-lines}, $L_y\parallel L\parallel L'$ implies $L_y\parallel L'$. Therefore, the flat hulls $\overline{L\cup L_y}$ and $\overline{L'\cup L_y}$ are planes in $X$. The $3$-rankedness of the regular liner $X$ implies that $P\cap\overline{L\cup L_y}=L$ and  $P\cap\overline{L'\cup L_y}=L'$. 

It follows from $x\in P\setminus L=P\setminus\overline{L\cup L_y}$ that $x\notin \overline{L\cup L_y}$. By Claim~\ref{cl:parallel-escorted}, there exists a line $L_x$ in $X$ such that $x\in L_x$ and $L_x\parallel L$. Then $L_x\subseteq \overline{L\cup\{x\}}\subseteq P$. The uniqueness of the line $L_x$ follows from the Proclus Axiom.
\end{proof}

The following example shows that Theorem~\ref{t:4-para-Playfair} does not hold for proaffine regular liners of rank $\|X\|=3$.

\begin{example}\label{ex:not-para-Playfair} {\em  There exists a countable $\w$-long Proclus plane, which is not para-Playfair.}
\end{example}

\begin{proof} 
Fix a family $(I_k)_{k\in 5}$ consisting of five pairwise disjoint countable infinite sets, and consider the liner $X_0\defeq\bigcup_{k\in 5}I_k$ endowed with the line relation
$$\Af_0\defeq \{xyz\in X_0^3:y\in \{x,z\}\}\cup\bigcup_{k\in 5}\{xyz\in I_k^3:x\ne z\}.$$ The sets $I_0,\dots,I_4$ are unique lines in the liner $X_0$ that contain more than two points. Let $L_0\defeq I_0$ and $L_0'\defeq I_1$.

For every $n\in\IN$ we shall inductively construct a liner $X_n$ containing two disjoint lines $L_n,L_n'$ so that the following conditions are satisfied:
\begin{enumerate}
\item[$(1_n)$] $|X_n|=\w$;
\item[$(2_n)$] $X_{n-1}$ is a subliner of the liner $X_n$;
\item[$(3_n)$] $L_{n-1}\subseteq L_n$ and $L_{n-1}'\subseteq L_n'$;
\item[$(4_n)$] for any points $x,y\in X_n$ with $\{x,y\}\not\subseteq X_{n-1}$, there exists a line $L$ in $X_n$ such that $L\cap\{x,y\}=\varnothing$;
\item[$(5_n)$] for every disjoint lines $L,\Lambda$ in the liner $X_{n}$ with $\{L,\Lambda\}\ne \{L_n,L_n'\}$,  the flat hulls $\overline L,\overline\Lambda$ of the sets $L,\Lambda$ in the liner $X_{n+1}$ have a common point.
\end{enumerate}

Assume that for some $n\in\w$ we have defined an increasing sequences of liners $(X_k)_{k\le n}$ and two sequences of lines $(L_k)_{k\le n}$ and $(L_k')_{k\le n}$ satisfying the conditions $(1_k)$--$(4_k)$ for all $k\le n$ and also the conditions $(5_k)$ for all $k<n$. Let $\Af_n$ be the line relation of the liner $X_n$, and $\mathcal L_n$ be the family of all lines in $X_n$. Let $$\mathcal D_n\defeq\big\{\{L,\Lambda\}:(L,\Lambda\in\mathcal L_n)\;\wedge\;(L\cap\Lambda=\varnothing)\;\wedge\;(\{L,\Lambda\}\ne\{L_n,L_n'\})\big\}$$ be the family of all unordered pairs of disjoint lines in the liner $X_n$, which are distinct from the selected pair of disjoint lines $\{L_n,L_n'\}$.
Observe that $$|\mathcal L_n|\le|X_n\times X_n|=\w\quad\mbox{and}\quad|\mathcal D_n|\le |\mathcal L_n\times\mathcal L_n|=\w.$$ 

 Choose any set $X_{n+1}$ containing $X_n$ so that the sets $|X_{n+1}\setminus X_n|$ and $\mathcal D_n$ have the same cardinality. Then $|X_{n+1}|=|X_n|+|\mathcal D_n|=\w$. Fix any bijection $p_n:\mathcal D_n\to X_{n+1}\setminus X_n$ and for any pair $\{A,B\}\in\mathcal D_n$, consider the point $p_{A,B}\defeq p_n(\{L,\Lambda\})\in X_{n+1}\setminus X_n$. Endow the set $X_{n+1}$ with the line relation
$$
\begin{aligned}
\Af_{n+1}&\defeq \Af_n\cup\big\{xyz\in X^3_{n+1}:y\in\{x,z\}\big\}\\
&\quad\;\;\cup \bigcup_{\{A,B\}\in\mathcal D_n}\bigcup_{C\in \{A,B\}}\{\big((C^2\times\{p_{A,B}\})\cup(\{p_{A,B}\}\times C^2)\cup\{xyz\in C\times\{p_{A,B}\}\times C:x\ne z\}\big),
\end{aligned}
$$ 
where $\Af_n$ is the line relation of the liner $X_n$. It is easy to see that $X_{n+1}$ is a liner containing $X_n$ as a subline and for every line $L$ in the liner $X_n$, the set 
$$L\cup\{p_{A,B}:L\in\{A,B\}\in\mathcal D_n\}$$is the unique line in the liner $X_{n+1}$ extending the line $L$. The definition of the family $\mathcal D_n$ and the bijectivity of the function $p_n$ ensure that  $$L_{n+1}\defeq L_n\cup\{p_{L,\Lambda}:L_n\in\{L,\Lambda\}\in\mathcal D_n\}\quad\mbox{and}\quad L'_{n+1}\defeq L_n\cup\{p_{L,\Lambda}:L'_n\in\{L,\Lambda\}\in\mathcal D_n\}$$are two disjoint lines in the liner $X_{n+1}$ satisfying the inductive condition $(3_{n+1})$.
 
It is easy to see that the the inductive conditions $(1_{n+1})$--$(3_{n+1})$ and $(5_n)$ are satisfied. To check the condition $(4_{n+1})$, fix any points $x,x'\in X_{n+1}$ such that $\{x,x'\}\not\subseteq X_n$. We lose no generality assuming that $x\notin X_n$ and hence $x=p_{L,\Lambda}$ for some pair $\{L,\Lambda\}\in\mathcal D_n$. Find the smallest number $i\le n+1$ such that $x'\in X_i$.

First assume that $i=0$ and hence $x'\in X_0=\bigcup_{k\in 5}I_k$. Since $I_0,\dots,I_4$ are disjoint lines in the liner $X_0\subseteq X_{n+1}$, there exists a number $k\in 5$ such that $I_k\not\subseteq L$, $I_k\not\subseteq \Lambda$ and $x'\notin I_k$. Then the flat hull $\overline {I_k}$ of the set $I_k$ in the liner $X_{n+1}$ is a line, disjoint with the doubleton $\{x,x'\}$.

Next, assume that $i>0$. Then $x'\in X_i\setminus X_{i-1}$ and hence  $x'=p_{L',\Lambda'}$ for some pair $\{L',\Lambda'\}\in\mathcal D_i$. In this case, we can find a number $k\in 5$ such that $I_k\not\subseteq J$ for any line $J\in\{L,L',\Lambda,\Lambda'\}$ and conclude that the flat hull $\overline{I_k}$ of the set $I_k$ in the liner $X_{n+1}$ is a line in $X_{n+1}$, disjoint with the doubleton $\{x,y\}$.
\smallskip

After completing the inductive construction, consider the liner $X\defeq\bigcup_{n\in\w}X_n$ endowed with the line relation $\Af\defeq\bigcup_{n\in\w}\Af_n$. The inductive conditions $(3_n)_{n\in\w}$ ensure that $L\defeq\bigcup_{n\in\w}L_n$ and $L'\defeq\bigcup_{n\in\w}L_n'$ are disjoint lines in the liner $X$. The required properties of the liner $X$ will be proved in separate claims.  

\begin{claim}\label{cl:all-intersects} Any lines $\Lambda,\Lambda'$ in $X$ with $\{\Lambda,\Lambda'\}\ne\{L,L'\}$ have a common point $p\in \Lambda\cap\Lambda'$.
\end{claim}

\begin{proof} To derive a contradiction, assume that $\Lambda\cap\Lambda'=\varnothing$. Find $n\in\w$ such that $|\Lambda\cap X_n|\ge 2$ and $|\Lambda'\cap X_n|\ge 2$. Then $\Lambda\cap X_n$ and $\Lambda'\cap X_n$ are disjoint lines in the liner $X_n$. If $\{\Lambda\cap X_n,\Lambda'\cap X_n\}=\{L_n,L_n'\}$, then $$\{\Lambda,\Lambda'\}=\{\overline{\Lambda\cap X_n},\overline{\Lambda'\cap X_n}\}=\{\overline{L_n},\overline{L_n'}\}=\{L,L'\},$$ which contradicts the choice of the lines $\Lambda,\Lambda'$. This contradiction shows that  $\{\Lambda\cap X_n,\Lambda'\cap X_n\}\ne\{L_n,L_n'\}$ and hence $\{\Lambda\cap X_n,\Lambda'\cap X_n\}\in\mathcal D_n$. The inductive condition $(5_n)$ ensures that the lines $\Lambda\cap X_{n+1}$ and $\Lambda'\cap X_{n+1}$ has a common point in the liner $X_{n+1}$, which is impossible as $\Lambda\cap\Lambda'=\varnothing$.
\end{proof}

\begin{claim}\label{cl:Y-Proclus} The liner $X$ is a Proclus plane.
\end{claim}

\begin{proof} Since the infinite liner $X_0$ is not a line, the liner $X$ is  not a line and hence $\|X\|\ge 3$. To see that $\|X\|=3$, fix any distinct points $a,b\in L$ and $c\in L'$. For every $x\in L'\setminus\{c\}$, the pair of lines $\{\Aline ac,\Aline bx\}$ is distinct from the pair $\{L,L'\}$. By Claim~\ref{cl:all-intersects}, there exists a point $p\in\Aline ac\cap\Aline bx$. It is easy to see that $p\ne b$ and hence $x\in \Aline pb\subseteq\overline{\{a,c,b\}}$. Therefore, $L\cup L'\subseteq\overline{\{a,b,c\}}$. For every $y\in Y\setminus (L\cup L')$, the pair of lines $\{\Aline ac,\Aline by\}$ is distinct from the pair $\{L,L'\}$. By Claim~\ref{cl:all-intersects}, there exists a point $q\in \Aline ac\cap\Aline by$. It is easy to see that $q\ne b$ and hence $y\in\Aline bq\subseteq\overline{\{b,a,c\}}$. Therefore, $X=\overline{\{a,b,c\}}$ is a plane. Claim~\ref{cl:all-intersects} implies that the plane $X$ is Proclus.
\end{proof}

\begin{claim}\label{cl:Y-3long} The Proclus plane $X$ is $\w$-long.
\end{claim}

\begin{proof} First we prove that the liner $X$ is $3$-long. In the oppsite case,  we could find two distinct points $x,y\in X$ such that the line $\Aline xy$ in $X$ coincides with the  doubleton $\{x,y\}$. Find the smallest number $n\in\w$ such that $\{x,y\}\subseteq  X_n$ and observe that or every $m\ge n$, the doubleton $\{x,y\}=\Aline xy\cap X_m$ is a line in the liner $X_m$. If $n>0$, then the minimality of $n$ ensures that $\{x,y\}\not\subseteq X_{n-1}$. By the inductive condition $(4_n)$, there exists a line $\Lambda\subseteq X_n$, disjoint with the doubleton $\{x,y\}$. Let $\overline \Lambda$ be the flat hull of the line $\Lambda$ in the liner $X_{n+1}$. The inductive condition $(5_n)$ ensures that the lines $\Aline xy\cap X_{n+1}=\{x,y\}$ and $\overline \Lambda$ have a common point $$p\in\{x,y\}\cap\overline \Lambda\subseteq \{x,y\}\cap X_n\cap\overline \Lambda=\{x,y\}\cap \Lambda=\varnothing,$$ which is a contradiction showing that $n=0$. Then $\{x,y\}\subseteq X_0$. Since the lines $I_2,I_3,I_4$ are disjoint, there exists a number $k\in\{2,3,4\}$ such that $\{x,y\}\cap I_k=\varnothing$. Then $\{\{x,y\},I_k\}\ne\{I_0,I_1\}=\{L_0,L_0'\}$ and hence the pair $\{\{x,y\},I_k\}$ belongs to the family $\mathcal D_0$. The inductive condition $(5_1)$ ensures that the flat hull $\overline{I_k}$ of the line $I_k$ in the liner $X_1$ has a common point $p$ with the line $\Aline xy\cap X_1=\{x,y\}$. Then $$p\in \{x,y\}\cap\overline {I_k}=\{x,y\}\cap X_n\cap\overline{I_k}=\{x,y\}\cap I_k=\varnothing,$$which is afinal contradiction showing that the liner $X$ is $3$-long.

By Claims~\ref{cl:Y-Proclus}, $X$ is a Proclus liner, and  Theorem~\ref{t:Proclus<=>}, the Proclus liner $X$ is proaffine. Observe that the line $L$ in $X$ contains the infinite set $I_0$ and hence is infinite. By  Corollary~\ref{c:L<L+2}, for every line $\Lambda$ in the $3$-long proaffine liner $X$, we have $\w=|L|\le|\Lambda|+2$, which implies that $|\Lambda|\ge\w$. Therefore, the Proclus plane $X$ is $\w$-long.
\end{proof}

Now we are able to prove that the Proclus plane $X$ has is not para-Playfair. Since $L,L'$ are two proper flats in the $\w$-long liner $X$, there exists a point $x\in Y\setminus(L\cup L')$. Assuming that $X$ is para-Playfair, we can find a line $\Lambda$ such that $x\in \Lambda\subseteq X\setminus L$. Then $L,\Lambda$ is a pair of disjoint lines such that $\{L,\Lambda\}\ne\{L,L'\}$, which contradicts Claim~\ref{cl:Y-3long}.
\end{proof}

\section{Bolyai lines in liners}

Let us recall that a liner $X$ is \index{Bolyai liner}\index{liner!Bolyai}\defterm{Bolyai} if for every plane $P\subseteq X$, line $L\subseteq P$ and point $x\in P\setminus L$ there exists at least one line $\Lambda$ such that $x\in \Lambda\subseteq P\setminus L$.

\begin{definition} A line $L$ in a liner $X$ is defined to be  \index{Bolyai line}\index{line!Bolyai}\defterm{Bolyai} if for every $x\in X\setminus L$ there exists a line $\Lambda$ in $X$ such that $x\in\Lambda\subseteq\overline{L\cup\{x\}}\setminus L$.
\end{definition}

\begin{exercise} Show that a line $L$ in a liner $X$ is Bolyai if and only if for every plane $P\subseteq X$ that contains the line $L$ and every point $x\in P\setminus L$, there exists a line $\Lambda\subset P$ such that $x\in\Lambda$ and $\Lambda\cap L=\varnothing$.
\end{exercise}

Observe that a liner $X$ is Bolyai if and only if every line in $X$ is Bolyai. This fact implies the following characterization of Playfair liners.

\begin{proposition} For a liner $X$ the following conditions are equivalent:
\begin{enumerate}
\item $X$ is Playfair;
\item $X$ is Proclus and Bolyai;
\item $X$ is Proclus and every line in $X$ is Bolyai.
\end{enumerate}
\end{proposition}

\begin{proposition}\label{p:lines-in-para-Playfair} A line $L$ in a para-Playfair regular liner $X$ of rank $\|X\|\ge 3$ is Bolyai if and only if there exists a line $\Lambda$ such that $L\cap\Lambda=\varnothing$ and the lines $L,\Lambda$ are coplanar.
\end{proposition}

\begin{proof} To prove the ``only if'' part, assume that the line $L$ is Bolyai in $X$. Since $\|X||\ge 3$, there exists a point $x\in X\setminus L$. Since $L$ is Bolyai, there exists a line $\Lambda$ in $X$ such that $x\in \Lambda\subseteq\overline{L\cup\{x\}}\setminus L$. Since $\overline{L\cup\{x\}}$ is a plane containing the lines $L,\Lambda$, the disjoint lines $L$ and $\Lambda$ are coplanar.
\smallskip

To prove the ``if'' part, assume that $\Lambda$ is a line in $X$ such that $L\cap\Lambda=\varnothing$ and the lines $L$ and $\Lambda$ are coplanar. Let $\Pi$ be a plane, containing both lines $L$ and $\Lambda$. To prove that the line $L$ is Bolyai in $X$, take any point $x\in X\setminus L$. We have to find a line $L_x\subseteq\overline{L\cup\{x\}}$ such that $x\in L_x$ and $L_x\cap L=\varnothing$. By Corollary~\ref{c:proregular=>ranked}, the regular liner $X$ is ranked. If $x\in \Pi$, then by the para-Playfair property of $X$, there exists a line $L_x$ in $X$ such that $x\in L_x\subseteq \Pi\setminus L$. Since the liner $X$ is $3$-ranked, $\Pi=\overline{L\cup\{x\}}$ and hence $L_x\subseteq\overline{L\cup\{x\}}$.

Next, assume that $x\notin \Pi$. The $3$-rankedness of the regular liner $X$ implies that $\Pi\cap\overline{L\cup\{x\}}=L$ and $\Pi\cap\overline{\Lambda\cup\{x\}}=\Lambda$. By Theorem~\ref{t:w-modular<=>}, the regular liner $X$ is weakly modular.  By Theorem~\ref{t:w-modular<=>}, the intersection $L_x\defeq\overline{L\cup\{x\}}\cap \overline{\Lambda\cup\{x\}}$ is a line containing $x$. Observe that 
$$
L_x\cap L=(\overline{L\cup\{x\}}\cap\overline{\Lambda\cup\{x\}})\cap (\Pi\cap L)\subseteq (\overline{\Lambda\cup\{x\}}\cap \Pi)\cap L=\Lambda\cap L=\varnothing,
$$
witnessing that the line $L$ is Bolyai.
\end{proof}

In fact, the regularity of the para-Playfair liner in Proposition~\ref{p:lines-in-para-Playfair} can be replaced by the $3$-long property of the liner.

\begin{proposition}\label{p:Bolyai-in-para-Playfair} A line $L$ in a $3$-long para-Playfair liner $X$ of rank $\|X\|\ge 3$ is Bolyai if and only if there exists a line $L'$ such that $L\cap L'=\varnothing$ and the lines $L,L'$ are coplanar.
\end{proposition}

\begin{proof} The ``only if'' part is proved in Proposition~\ref{p:lines-in-para-Playfair}. To prove the ``if'' part, assume that  there exists a line $L'$ such that $L\cap L'=\varnothing$ and the lines $L,L'$ lie in some plane $\Pi$. Given any point $x\in X\setminus L$, we need to find a line $L_x$ in $X$ such that $x\in L_x\subseteq\overline{L\cup\{x\}}\setminus L$.

By Proposition~\ref{p:Playfair=>para-Playfair=>Proclus}, the para-Playfair liner $X$ is Proclus, by Theorem~\ref{t:Proclus<=>}, the Proclus liner $X$ is proaffine and $3$-proregular, and by Proposition~\ref{p:k-regular<=>2ex}, the $3$-proregular line $X$ is $3$-ranked. 

Consider the plane $P\defeq\overline{L\cup\{x\}}$ and observe that $L\subseteq\Pi\cap P$. If $x\in \Pi$, then $\Pi=P$, by the $3$-rankedness of $X$. The  para-Playfair property of $X$ ensures the existence of a line $L_x$ such that $x\in L_x\subseteq\Pi\setminus L=\overline{L\cup\{x\}}\setminus L$. So, assume that $x\notin \Pi$, which implies $\Pi\cap P=L$, by the $3$-rankedness of the liner $X$.

Fix any point $y\in L'$. Since $X$ is $3$-long, there exists a point $c\in \Aline xy\setminus\{x,y\}$. Consider the set $I\defeq\{p\in \Pi:\Aline pc\cap P=\varnothing\}$. The choice of the points $y$ and $c\in \Aline xy$ ensure that $y\notin I$. By the proaffinity of $X$, for every $z\in L$, the set $I_z\defeq\{p\in \Aline yz:\Aline pc\cap\Aline zx=\varnothing\}$ contains at most one point. Then for every point  $p\in\Aline zc\setminus I_z$, we have  $\varnothing\ne\Aline pc\cap\Aline zx\subseteq \Aline pc\cap P$ and hence the intersection $I\cap\Aline yz\subseteq I_z$ contains at most one point. 
\smallskip

Two cases are possible.
\smallskip

1. First assume that $L'\setminus\{y\}\not\subseteq I$ and hence there exists a point $y'\in L'\setminus (I\cup\{y\})$. By the definition of the set $I\not\ni y'$, there exists a point $x'\in \Aline c{y'}\cap P$. Assuming that $x=x'$, we conclude that $y'\in \Aline c{x'}\cap\Pi=\Aline xy\cap\Pi=\{y\}$, which contradicts the choice of $y'\ne y$. Therefore, $L_x\defeq\Aline x{x'}$ is a line in the plane $P$. Assuming that $L_x\cap L\ne\varnothing$, we can find a point $z\in L\cap L_x$ and consider the plane $\overline{\{x,y,z\}}$. Since $x\notin \Pi$, the intersection $\overline{\{x,y,z\}}\cap\Pi$ coincides with the line $\Aline yz$, by the $3$-rankedness of the liner $X$. On the other hand, the plane $\overline{\{x,y,z\}}$ contains the points $c\in\Aline xy$, $x'\in\Aline xz$ and $y'\in \Aline c{x'}$. Then $\Aline y{y'}\subseteq \overline{\{x,y,z\}}\cap \Pi=\Aline yz$ and hence $\Aline yz=\Aline y{y'}$ and $z\in\Aline y{y'}\cap L=L'\cap L=\varnothing$, which is a contradiction witnessing that $L_x\cap L=\varnothing$. 
\smallskip

2. Next, assume that $L'\setminus\{y\}\subseteq I$. Fix any point $z\in L$. Taking into account that $X$ is $3$-long, we can choose two distinct point $y',y''\in L'\setminus\{y\}$ and find a point $u'\in \Aline z{y'}\setminus\{z,y'\}$. Since $u'\notin \{y'\}=\Aline z{y'}\cap I$, there exists a point $v'\in \Aline c{u'}\cap P$. Since $X$ is para-Playfair, there exists a line $\Lambda$ in $X$ such that $u'\in \Lambda\subseteq \Pi\setminus L'$. Since $X$ is Proclus and $\Lambda\cap L'=\varnothing$, there exist a point $u''\in\Aline z{y''}\cap\Lambda$.   Since $u''\notin \{y''\}=\Aline z{y'}\cap I$, there exists a point $v''\in \Aline c{u''}\cap P$. By analogy with the preceding case, we can show that $L''\defeq\Aline {v'}{v''}$ is a line in the plane $P$, disjoint with the line $L$. Then the plane $P$ in the para-Playfair liner $X$ contains two disjoint lines $L$ and $L''$ and hence there exists a line $L_x\subseteq P$ such that $x\in L_x$ and $L_x\cap L=\varnothing$. 
\end{proof}

\section{Bi-Bolyai liners}\label{s:bi-Bolyai}

\begin{definition}\label{d:bi-Bolyai} A liner $X$ is called \index{bi-Bolyai liner}\index{liner!bi-Bolyai}\defterm{bi-Bolyai} if for every concurrent Bolyai lines $L,\Lambda$ in $X$, any line in the plane $\overline{L\cup\Lambda}$ is Bolyai in $X$.
\end{definition}

Definitions~\ref{d:PPBL} and \ref{d:bi-Bolyai} imply that for every liner the following implications hold:
$$\mbox{Playfair}\Ra\mbox{Bolyai}\Ra\mbox{bi-Bolyai}.$$

The main result of this section is the following regularity theorem.

\begin{theorem}\label{t:4long+pP+bB=>regular} Every $4$-long para-Playfair bi-Bolyai liner $X$ is regular.
\end{theorem}

\begin{proof} By Proposition~\ref{p:Playfair=>para-Playfair=>Proclus}, the para-Playfair liner $X$ is Proclus. By Theorem~\ref{t:Proclus<=>}, the Proclus liner $X$ is proaffine and $3$-proregular. Being $3$-proregular and $4$-long, the liner $X$ is $3$-regular and $3$-ranked. 

To prove that the liner $X$ is $4$-regular, fix any plane $A\subseteq X$ and line $\Lambda\subseteq X$ such that $A\cap\Lambda=\{o\}$ for some point $o$. We have to show that the set $\Aline A\Lambda\defeq\bigcup_{x\in A}\bigcup_{y\in\Lambda}\Aline xy$ is flat. First we prove seven claims describing some properties of the set $\Aline A\Lambda$.

\begin{claim}\label{cl:Bolyai-parallel} For any distinct points $a,b\in\Aline A\Lambda$ with $\Aline ab\cap A=\varnothing$, we have $\Aline ab\subseteq\Aline A\Lambda$.
\end{claim}

\begin{proof}  
Since $a,b\in \Aline A\Lambda\setminus A$, there exist points $\alpha,\beta\in A$ and $u,v\in\Lambda\setminus\{o\}$ such that $a\in\Aline\alpha u$ and $b\in\Aline \beta v$. Since $X$ is proaffine, the subsets $U\defeq\{\lambda\in\Aline ou:\Aline\lambda a\cap\Aline o\alpha=\varnothing\}$ and  $V\defeq\{\lambda\in\Aline ov:\Aline\lambda b\cap\Aline o\beta=\varnothing\}$ of the line $\Lambda$ have cardinality $|U\cup V|\le|U|+|V|\le 2$. Since the liner $X$ is $4$-long, there exists a point $\lambda\in\Lambda\setminus(\{o\}\cup U\cup V)$. The definitions of the sets $U,V$ ensure that there exist points $\alpha'\in \Aline \lambda a\cap\Aline o\alpha\subseteq A$ and $\beta'\in \Aline \lambda b\cap\Aline o\beta\subseteq A$.
Consider the plane $\Pi=\overline{\{\lambda,a,b\}}$ and observe that $\Aline ab$ and $\Aline {\alpha'}{\beta'}$ are two disjoint lines in $\Pi$. By the Proclus Postulate~\ref{p:Proclus-Postulate}, for every point $c\in \Aline ab$, there exists a point $\gamma\in\Aline\lambda c\cap\Aline{\alpha'}{\beta'}\subseteq A$. Then $c\in \Aline \gamma\lambda\subseteq \Aline A\Lambda$ and hence $\Aline ab\subseteq \Aline A\Lambda$.
\end{proof}



\begin{claim}\label{cl:Bolyai-coplanar} If for distinct points $a,b\in\Aline A\Lambda$, the lines $\Aline ab$ and $\Lambda$ are coplanar, then $\Aline ab\subseteq\Aline A\Lambda$.
\end{claim}

\begin{proof} Let $\Pi$ be a plane containing the coplanar lines $\Aline ab$ and $\Lambda$. If $\Aline ab\cap A=\varnothing$, then $\Aline ab\subseteq \Aline A\Lambda$, by Claim~\ref{cl:Bolyai-parallel}. So, assume that $\Aline ab\cap A\ne\varnothing$. In this case, we lose no generality assuming that $a\in A$ and $b\notin A\cup\Lambda$. Since $b\in\Aline A\Lambda\setminus(A\cup\Lambda)$, there exist points $\beta\in A\setminus\Lambda$ and $\lambda\in\Lambda\setminus A$ such that $b\in\Aline\beta\lambda$. Then $\beta\in \Aline \lambda b\cap A\subseteq \Pi\cap A$ and $b\in\Aline L\Lambda$ for the line $L\defeq\Aline o\beta\subseteq \Pi\cap A$. By the $3$-regularity of $X$, $\Aline ab\subseteq\Pi=\Aline L\Lambda\subseteq\Aline A\Lambda$. 
\end{proof}

Assuming that the set $\Aline A\Lambda$ is not flat in $X$, we can find a point $c\in \Aline{\Aline A\Lambda}{\Aline A\Lambda}\setminus\Aline A\Lambda$.

\begin{claim}\label{cl:Bolyai-tangent} $\overline{\Lambda\cup\{c\}}\cap A=\{o\}$.
\end{claim}

\begin{proof} To derive a contradiction, assume that  $\overline{\Lambda\cup\{c\}}\cap A\ne\{o\}$ and hence $\overline{\Lambda\cup\{c\}}\cap A$ contains some point $\alpha\notin\Lambda$. By the $3$-rankedness of the $3$-regular liner $X$, $\overline{\Lambda\cup\{\alpha\}}=\overline{\Lambda\cup\{c\}}$. By the $3$-regularity of the plane $\overline{\Lambda\cup\{\alpha\}}$, there exist points $\gamma\in\Aline o\alpha$ and $\lambda\in \Lambda$ such that $c\in\Aline \alpha\lambda\subseteq\Aline A\Lambda$, which contradicts the choice of the point $c$. This contradiction shows that  $\overline{\Lambda\cup\{c\}}\cap A=\{o\}$.
\end{proof}

\begin{claim}\label{cl:Bolyai-disjoint} For any points $a,\beta\in A$, $\lambda\in \Lambda\setminus\{o\}$, and $b\in\Aline \beta\lambda$ with $c\in\Aline ab$, the lines $\Aline \lambda c$ and $\Aline a\beta$ in the plane $\Pi\defeq\overline{\{a,b,\lambda\}}$ are disjoint.
\end{claim}

\begin{proof} By Claim~\ref{cl:Bolyai-coplanar}, the lines $\Aline ab$ and $\Lambda$ are not coplanar and hence $o\notin \Pi$ and $o\notin \Aline\lambda c\subseteq\Pi$. Claim~\ref{cl:Bolyai-tangent} ensures that $\Aline \lambda c\cap\Aline a\beta\subseteq (\overline{\Lambda\cup\{c\}}\cap A)\setminus\{o\}=\varnothing$. Therefore, $\Aline \lambda c$ and $\Aline a\beta$ are two disjoint lines in the plane $\Pi$.
\end{proof}

\begin{claim}\label{cl:Bolyai-absorb} For any point $a\in A$ with $\Aline ac\cap \Aline A\Lambda\ne\{a\}$, we have $\Aline ac\setminus\{c\}\subseteq \Aline A\Lambda$.
\end{claim}

\begin{proof} Assume that $a\in A$ is a point such that the set $\Aline ac\cap\Aline A\Lambda$ contains some point $b\ne a$. Since $b\in \Aline A\Lambda$, there exist  points $\lambda\in\Lambda$ and $\beta\in A$ such that $b\in\Aline\beta\lambda$.  It follows from $c\in\Aline ab\setminus \Aline A\Lambda$ that $b\notin A\cup\Lambda$ and hence $\beta\ne o\ne\lambda$. Given any point $d\in \Aline ac\setminus\{c\}$, we need to show that $d\in\Aline A\Lambda$. Claim~\ref{cl:Bolyai-disjoint} ensures that the lines $\Aline a\beta$ and $\Aline \lambda c$ in the plane $\overline{\{\lambda,a,c\}}$ are disjoint. Taking into account that the liner $X$ is Proclus and $\Aline \lambda d\ne\Aline \lambda c$, we conclude that there exists a point $\delta\in \Aline \lambda d\cap\Aline a\beta\subseteq \Aline a\beta\subseteq A$. Then $d\in \Aline \delta\lambda\subseteq\Aline A\Lambda$.
\end{proof}

\begin{claim}\label{cl:Bolyai-concurrent} For any point $a\in A$ with $\Aline ac\cap \Aline A\Lambda\ne\{a\}$, every line $L\subset\overline{\{a,o,c\}}$ containing the point $c$ is concurrent with the line $\Aline oa$.
\end{claim}

\begin{proof} To derive a contradiction, assume that the plane $\overline{\{a,o,c\}}$ contains a line $L $ such that $c\in L$ and $L\cap\Aline ao=\varnothing$. By Claim~\ref{cl:Bolyai-coplanar}, $o\notin \Aline ac$. Since the liner $X$ is $4$-long, there exist two distinct points $b,d\in\Aline ac\setminus\{a,c\}$. By Claim~\ref{cl:Bolyai-absorb}, $\{b,d\}\subseteq \Aline ac\setminus\{c\}\subseteq\Aline A\Lambda$. 
Since $\Aline ac\cap\Aline ao=\{a\}$, the lines $\Aline ob$ and $\Aline od$ are distinct from the line $\Aline oa$, which is disjoint with the line $L$ in the plane $\overline{\{a,o,c\}}$. By the Proclus Postulate~\ref{p:Proclus-Postulate}, there exist points $\beta\in \Aline ob\cap L$ and $\delta\in\Aline od\cap L$. It follows from $b\ne d$ that $\beta\ne \delta$ and hence $L=\Aline \beta\delta$. By Claim~\ref{cl:Bolyai-coplanar}, $\{\beta,\delta\}\subseteq\Aline ob\cup\Aline od\subseteq \Aline A\Lambda$. By Claim~\ref{cl:Bolyai-parallel}, $c\in \Lambda=\Aline\beta\delta\subseteq\Aline A\Lambda$, which contradicts the choice of the point $c$. 
\end{proof}

\begin{claim}\label{cl:abeta-Bolyai} For any points $a,\beta\in A$, $\lambda\in \Lambda$, and $b\in\Aline \beta\lambda$ with $c\in\Aline ab$, the line $\Aline a\beta$ is Bolyai.
\end{claim}

\begin{proof} It follows from $a\in A$ and $c\in \Aline ab\setminus \Aline A\Lambda$ that $b\notin A\cup\Lambda$ and hence $\beta\ne o\ne\lambda$. Since $X$ is Proclus, the set $I\defeq\{x\in \Aline ob:\Aline \lambda x\cap\Aline o\beta=\varnothing\}$ contains at most one point. Since the liner $X$ is $4$-long, there exists a point $b'\in\Aline ob\setminus(\{o,b\}\cup I)$. Since $b'\notin I$, there exists a point $\beta'\in\Aline\lambda{b'}\cap\Aline o\beta\subseteq A$. Then $b'\in \Aline {\beta'}\lambda\subseteq\Aline A\Lambda$.
By Claim~\ref{cl:Bolyai-concurrent}, the line $\Aline c{b'}$ in the plane $\overline{\{a,o,c\}}$ has a common point $a'$ with the line $\Aline oa$.
Claim~\ref{cl:Bolyai-disjoint} ensures that $\Aline a\beta\cap\Aline \lambda c=\varnothing=\Aline{a'}{\beta'}\cap\Aline \lambda c$. Assuming that the line  $L\defeq\Aline{a'}{\beta'}$ has a common point $\alpha$ with the line $\Aline a\beta$, we conclude that $c,\lambda,\alpha\in \overline{\{a,\beta,\lambda\}}\cap\overline{\{a',\beta',\lambda\}}$ and hence $\overline{\{a,\beta,\lambda\}}=\overline{\{c,\lambda,\alpha\}}=\overline{\{a',\beta',\lambda\}}=\overline{\{a,b,a',b',\lambda\}}$ and $o\in \Aline b{b'}\subseteq \overline{\{a,b,\lambda\}}$. Then the lines $\Lambda=\Aline o\lambda$ and $\Aline ab$ are coplanar, which contradicts Claim~\ref{cl:Bolyai-coplanar}. This contradiction shows that $L\cap\Aline a\beta=\varnothing$. By Proposition~\ref{p:Bolyai-in-para-Playfair}, the line $\Aline a\beta$ is Bolyai in $X$.
\end{proof}

By the choice of the point $c\in\Aline{\Aline A\Lambda}{\Aline A\Lambda}\setminus\Aline A\Lambda$, there exist distinct points $a,b\in\Aline A\Lambda$ such that $c\in\Aline ab$. Claim~\ref{cl:Bolyai-parallel} ensures that $\Aline ab\cap A\ne\varnothing$, so we can assume that $a\in A$ and $b\notin A\cup\Lambda$. By Claim~\ref{cl:Bolyai-coplanar}, the lines $\Aline ab$ and $\Lambda$ are not coplanar, which implies that $o\notin\Aline ab$ and hence $\Pi\defeq\overline{\{o,a,b\}}$ is a plane. Since $b\in\Aline A\Lambda$, there exist points $\beta\in A$ and $\lambda\in \Lambda$ such that $b\in\Aline \beta\lambda$. By Claim~\ref{cl:abeta-Bolyai}, $\Aline a\beta$ is a Bolyai line in $X$.

Since the Proclus liner $X$ is proaffine, the set $I\defeq\{y\in\Lambda:\Aline by\cap\Aline o\beta=\varnothing\}$ contains at most one point. Since the liner $X$ is $4$-long, there exists a point $\mu\in\Lambda\setminus(\{o,\lambda\}\cup I)$. Then there exists a point $\gamma\in \Aline \mu b\cap\Aline o\beta\subseteq\Aline o\beta\subseteq A$. The inequality $\lambda\ne \mu$ implies the inequality $\beta\ne\gamma$ and hence $\Aline o\beta\ne\Aline o\gamma$. Since $b\in\Aline \mu\gamma$, we can apply Claim~\ref{cl:abeta-Bolyai} and conclude that $\Aline a\delta$ is a Bolyai line in $X$. Since $X$ is bi-Bolyai and the Bolyai lines $\Aline a\beta$ and $\Aline a\gamma$ are concurrent, any line in theplane $\overline{\Aline a\beta\cup\Aline a\gamma}=A$ is Bolyai. In particular, the line $\Aline ao\subseteq A$ is Bolyai in $X$ and hence  there exists a line $L_c\subseteq \overline{\{a,o,c\}}$ such that $L_c\cap \Aline oa=\varnothing$, which contradicts Claim~\ref{cl:Bolyai-parallel}. This contradiction completes the proof of the $4$-regularity of the liner $X$. By Theorem~\ref{t:HA}, the $4$-regular liner $X$ is regular.
\end{proof}

\begin{remark} In Theorem~\ref{t:spread=projective2} we shall prove that every {\em finite} $4$-long para-Playfair liner is bi-Bolyai and regular.
\end{remark}

\begin{Exercise} Construct an example of an $\w$-long para-Playfair liner which is not bi-Bolyai.
\smallskip

{\em Hint:} See Example~\ref{ex:para-Playfair}.
\end{Exercise} 

\section{Spreading lines in liners}

\begin{definition}
Let $X$ be a liner and $\mathcal L$ be the family of lines in $X$. 
A subfamily $\mathcal S\subseteq \mathcal L$ is called a \index{spread of lines}\defterm{spread of lines} in $X$ if every point $x\in X$ belongs to a unique line $L\in\mathcal S$. 

A line $L\in\mathcal L$ is called \index{spreading line}\index{line!spreading}\defterm{spreading} if the family\index[note]{$L_\parallel$} 
$$L_\parallel \defeq\{\Lambda\in\mathcal L:\Lambda\parallel L\}$$ is a spread of lines in $X$ such that $L_\parallel =\Lambda_\parallel$ for every line $\Lambda\in L_\parallel$.
\end{definition}

\begin{exercise} Show that a family of lines $\mathcal S$ in a liner $X$ is a spread of lines if and only if $\bigcup\mathcal S=X$ and any two distinct lines in $\mathcal S$ are disjoint.
\end{exercise}

\begin{theorem}\label{t:spreading<=>} For a line $L$ in a proaffine regular liner $X$, the following conditions are equivalent:
\begin{enumerate}
\item $L$ is spreading;
\item $L$ is Bolyai;
\item $\bigcup L_\parallel=X$.
\end{enumerate}
\end{theorem}

\begin{proof} $(1)\Ra(2)$ If the line $L$ is spreading, then for every $x\in X\setminus L$, there exists a line $L_x\in L_\parallel$ such that $x\in L_x$. Assuming that $L\cap L_x\ne\varnothing$, we conclude that $x\in L_x=L$, which contradicts the choice of the point $x$. This contradiction shows that $L_x\cap L=\varnothing$. It follows from $L_x\parallel L$ and $L_x\cap L\ne\varnothing$ that $L_x\subseteq \overline{L\cup\{x\}}\setminus L$, witnessing that the line $L$ is Bolyai.
\smallskip

$(2)\Ra(3)$ If the line $L$ is Bolyai, then for every $x\in X\setminus L$, there exists a line $L_x$ such that $x\in L_x\subseteq\overline{L\cup\{x\}}\setminus L$.  Corollary~\ref{c:parallel-lines<=>} ensures that $L_x\in L_\parallel$, which implies $\bigcup L_\parallel=X$.
\smallskip 

$(3)\Ra(1)$ Assume that $\bigcup L_\parallel =X$. To prove that the line $L$ is spreading, we need to check that $L_\parallel $ is a spread of lines in $X$ and $L_\parallel=\Lambda_\parallel $ for every $\Lambda\in L_\parallel$.

Theorem~\ref{t:Proclus<=>} ensures that the proaffine regular liner $X$ is Proclus. 
To see that $L_\parallel$ is a spread of lines in $X$, we have to show that every point $x\in X$ is contained in a unique line $L_x\in L_\parallel$. Since $\bigcup L_\parallel =X$, there exists a line $L_x\in L_\parallel $ such that $x\in L_x$. If $L_x=X$, then $X$ is a line and $L_x$ is a unique line in the family $L_\parallel=\{X\}$ that contains $x$. So, assume that $L_x\ne X$ and choose any point $y\in X\setminus L_x$. Since $y\in X=\bigcup L_\parallel$, there exists a line $L_y\in L_\parallel$ such that $y\in L_y$. Since $L_x\parallel L\parallel L_y$, Theorem~\ref{t:Proclus-lines} ensures that $L_x\parallel L_y$ and hence $L_x\cap L_y=\varnothing$, see Proposition~\ref{p:para+intersect=>coincide}.
Assuming that the family $L_\parallel$ contains another line $L_x'\in L_\parallel$ with $x\in L_x'$, we can apply Theorem~\ref{t:Proclus-lines} and conclude that $L_x'\parallel L_y$ and hence $L_x'\cap L_y=\varnothing$. Then $L_x,L_x'$ are two lines in the plane $\overline{\{x\}\cup L_y}$ that contain the point $x$ and are disjoint with the line $L_y$. The Proclus Axiom ensures that $L_x=L_x'$, witnessing that the family $L_\parallel$ is a spread of lines.
Theorem~\ref{t:Proclus-lines} implies that for every line $\Lambda\in L_\parallel$, the families $L_\parallel$ and $\Lambda_\parallel$ are equal, witnessing that the line $L$ is spreading.
\end{proof}

\begin{theorem}\label{t:Playfair<=>spreading} A liner $X$ is regular and Playfair if and only if $X$ is $3$-ranked and every line in $X$ is speading.
\end{theorem}

\begin{proof} To prove the ``only if'' part, assume that the liner $X$ is regular and Playfair. By Theorem~\ref{t:Playfair<=>}, the Playfair liner $X$ is $3$-regular and by Proposition~\ref{p:k-regular<=>2ex}, the $3$-regular liner $X$ is $3$-ranked. By Theorems~\ref{t:Playfair}, for every line $L$, the family $L_\parallel $ is a  spread of lines, and by Theorem~\ref{t:Play-reg<=>Bolyai}, $L_\parallel=\Lambda_\parallel$ for every line $\Lambda\in L_\parallel$.

To prove the ``if'' part, assume that the liner $X$ is $3$-ranked and every line in $X$ is spreading. To show that $X$ is Playfair, take any plane $P\subseteq X$, line $L\subseteq P$ and point $x\in P\setminus L$. Since the line $L$ is spreading, there exists a line $L_x\in L_\parallel$ such that $x\in L_x$. Then $L_x\subseteq\overline{\{x\}\cup L}\setminus L\subseteq P\setminus L$. Assuming that $L_x'$ is another line such that $x\in L_x'\subseteq P\setminus L$, we can apply Corollary~\ref{c:parallel-lines<=>} and conclude that $L_x'\parallel L$. Therefore, $L_x,L_x'\in L_\parallel$ are two lines containing the point $x$. Since $L_\parallel$ is a spread of lines, $L_x'=L_x$. Therefore, $L_x$ is a unique line in $X$ with $x\in L_x\subseteq P\subseteq L$, witnessing that the liner $X$ is Playfair and hence Bolyai. By Theorem~\ref{t:Playfair<=>} and Proposition~\ref{p:k-regular<=>2ex}, the Playfair liner $X$ is $3$-regular and $3$-ranked. By Theorem~\ref{t:Play-reg<=>Bolyai}, the regularity of $X$ will follow as soon as we check that for every lines $A,B,C\subseteq X$, $A\parallel B\parallel C$ implies $A\parallel C$. It follows from $A\parallel B\parallel C$ that $B\in A_\parallel$ and $C\in B_\parallel$. Since the line $A$ is spreading, $C\in B_\parallel=A_\parallel$ and hence $A\parallel C$. 
\end{proof}

\begin{remark} Let $A,B$ be two spreading lines $A,B$ in a liner $X$. If the spreads $A_\parallel $ and $B_\parallel$ contain some common line $L$, then $A_\parallel=L_\parallel=B_\parallel$. Therefore, for two spreading lines $A,B$, the spreads $A_\parallel$ and $B_\parallel$ either coincide or else are disjoint.
\end{remark}

\begin{proposition}\label{p:spread-3-long} If $L$ is a spreading line in a $3$-ranked liner $X$, then every line $\Lambda\notin L_\parallel$ in $X$ is $3$-long.
\end{proposition}

\begin{proof} To derive a contradiction, assume that some line $\Lambda\notin L_\parallel$ in $X$ has cardinality $|\Lambda|=2$. Find two points $x,y\in X$ such that $\Lambda=\{x,y\}$ and choose unique lines $L_x,L_y\in L_\parallel$ such that $x\in L_x$ and $y\in L_y$. Assuming that $L_x\cap L_y\ne\varnothing$, we can apply Proposition~\ref{p:para+intersect=>coincide} and conclude that $L_x=L_y$ and hence $\Lambda=L_x=L_y\in L_\parallel$, which contradicts the choice of the line $\Lambda\notin L_\parallel$.  
This contradiction shows that $L_x\cap L_y=\varnothing$. Assuming that for every $u\in L_x$ we have $\Aline yu=\{y,u\}$, we conclude that the set $L_x\cup\{y\}$ is flat and hence $L_y\subseteq \overline{L_x\cup\{y\}}=L_x\cup\{y\}$, which contradicts $L_x\cap L_y=\varnothing$. This contradiction shows that for some point $u\in L_x$, the line $\Aline uy$ contains a point $v\in \Aline uy\setminus\{u,y\}$. Let $L_v\in L_\parallel$ be the unique line in the spread $L_\parallel$ containing the point $v$. Then $L_x,L_y,L_v$ are three disjoint parallel lines in the plane $P=\overline{L\cup\{y\}}$, and $\Lambda,L_v$ are two disjoint lines in the plane $P$. Corollary~\ref{c:parallel-lines<=>} ensures that $\Lambda\parallel L_v$ and hence $\Lambda\in L_\parallel$, which contradicts the choice of $\Lambda$. 
\end{proof}

\begin{exercise} Find an example of a Proclus plane $X$ of cardinality $|X|=6$ containing a spreading line $L$ of cardinality $|L|=2$.
\end{exercise}  

\section{The spread completion of a liner}

Let $X$ be a liner, $\mathcal L$ be the family of all lines in $X$, and $\mathcal S$ be the family of all spreading lines in $X$. The set\index{$\partial X$}  
$$\partial X\defeq\{L_\parallel:L\in\mathcal S\}$$of spreads of spreading lines in $X$ is called \index{boundary}
\index{liner!boundary of}\defterm{the boundary} of $X$. Since $\bigcup(\partial X)=\mathcal S$, the family $\mathcal S$ can be recovered from the boundary of $X$. Elements $L_\parallel\in\partial X$ of the boundary will be called \index{direction}\defterm{directions} in the liner $X$.

We claim that $X\cap\partial X=\varnothing$. Indeed, assuming that some point $x\in X$ is equal to the spread $L_{\parallel}$ for some speading line $L\in\mathcal S$, we can find a line $L_x$   such that $x\in L_x\in L_{\parallel}=x$, which contradicts the Axiom of Foundation in Set Theory\footnote{The Axiom of Foundation says that every nonempty set $x$ contains an element $y\in x$ such that $y\cap x=\varnothing$.}

Attaching to the liner $X$ its boundary $\partial X$, we obtain the set\index[note]{$\overline X$} $$\overline X\defeq X\cup\partial X,$$which is the underlying set of the spread completion of $X$. Endow the set $\overline X$ with the ternary relation\index[note]{$\overline{\mathsf L}$} $$
\begin{aligned}
\overline\Af&\defeq\Af\cup\{(x,y,\boldsymbol z)\in X\times X\times \partial X:x=y\;\vee\;\Aline xy\in \boldsymbol z\}\\
&\cup \{(x,\boldsymbol y,z)\in X\times \partial X\times X:\Aline xz\in\boldsymbol y\}\cup\{({\boldsymbol x},y,z)\in \partial X\times X\times X:y=z\;\vee\;\Aline yz\in {\boldsymbol x}\}\\
&\cup\{(x,{\boldsymbol y},{\boldsymbol z})\in X\times\partial X\times\partial X:{\boldsymbol y}={\boldsymbol z}\}\cup\{({\boldsymbol x},{\boldsymbol y},z)\in X\times\mathcal S_\parallel\times\mathcal S_\parallel:{\boldsymbol x}={\boldsymbol y}\}\\
&\cup\{({\boldsymbol x},{\boldsymbol y},{\boldsymbol z})\in \partial X\times\partial X\times\partial X:\forall A\in{\boldsymbol x}\;\forall B\in {\boldsymbol y}\;\forall C\in {\boldsymbol z}\;(A\cap B\cap C\ne\varnothing\;\Ra\; B\subseteq\overline{A\cup C})\},
\end{aligned}
$$
where $\Af$ is the line relation of the liner $X$.

\begin{lemma}\label{l:overlineLisLineal} If the liner $X$ is $3$-ranked, then $\overline\Af$ is a line relation on the set $\overline X$.
\end{lemma}

\begin{proof}  The axioms {\sf(IL)} and {\sf(RL)} follow immediately from the definition of $\overline L$. To prove that  $\overline L$ satisfies the Exchange Axiom, fix any points  $a,b,x,y\in\overline X$ with  $\overline\Af axb\;\wedge\;\overline \Af ayb\;\wedge\;x\ne y$. We should prove that $\overline \Af xay\wedge\overline\Af xby$.
Depending on the position of the points $x,y,a,b$ in the set $\overline X=X\cup\partial X$, we should consider $2^4=16$ cases.

1. If $a,b,x,y\in X$, then $\{xay,aby\}\subseteq\Af\subseteq\overline\Af$ by the Exchange Axiom for the line relation $\Af$.

2. If $a,b,x\in X$ and $y\in\partial X$, then $\overline\Af axb\;\wedge\;\overline \Af ayb$ implies $x\in\Aline ab\in y$. Then $\{xay,xby\}\subseteq\overline\Af$, by the definition of $\Af$.

3. If $a,b,y\in X$ and $x\in\partial X$, then $\overline\Af axb\;\wedge\;\overline \Af ayb$ implies $y\in\Aline ab\in x$. Then $\{xay,xby\}\subseteq\overline\Af$, by the definition of $\Af$.

4. If $a,b\in X$ and $x,y\in\partial X$, then $\overline\Af axb\;\wedge\;\overline \Af ayb$ implies $\Aline ab\in x\cap y$ and hence $x=y$, which contradicts our assumption $x\ne y$.

5. If $a,x,y\in X$ and $b\in\partial X$, then $\overline\Af axb\;\wedge\;\overline \Af ayb$ implies $a=x$ or $\Aline ax\in b$ and $a=y$ or $\Aline ay\in b$.  If $a=x$, then $xay\in\Af\subseteq\overline\Af$ by the axiom {\sf (RL)}. Taking into account that $x\ne y$, we conclude that $a\ne y$ and hence $\Aline xy=\Aline ay\in b$ and $xby\in \overline\Af$. Therefore, $a=x$ implies $\{xay,xby\}\subseteq\overline\Af$. By analogy we can prove that $a=y$ implies   $\{xay,xby\}\subseteq\overline\Af$. So, assume that $x\ne a\ne y$. Then $\Aline ax\in b$ and $\Aline ay\in b$ and hence $\Aline ax=\Aline ay$ and $a\in a\Aline xy=\Aline ax=\Aline ay\in b$ and $\{xay,xby\}\subseteq \overline\Af$.

6. If $a,x\in X$ and $b,y\in\partial X$, then $\overline\Af axb\;\wedge\;\overline \Af ayb$ implies $a=x$ or $\Aline ax\in b$ and $y=b$. The last equality implies $xby\in\overline\Af$. If $a=x$, then $xay\in\overline \Af$. If $\Aline xa\in b$, then $\Aline xa\in b=y$ and hence $xay\in\overline\Af$.

7. The case $a,y\in X$ and $b,x\in\partial X$ can be considered by analogy with the case 6.

8. If $a\in X$ and $b,x,y\in\partial X$, then $\overline\Af axb\;\wedge\;\overline \Af ayb$ implies $x=b=y$, which contradicts the choice of $x\ne y$.

9. The case $a\in\partial X$ and $b,x,y\in X$ can be considered by analogy with case 5.

10. The case $a,y\in\partial X$ and $b,x\in X$ can be considered by analogy with case 6.

11. The case $a,x\in\partial X$ and $b,y\in X$ can be considered by analogy with case 6.

12.  The case $a,x,y\in\partial X$ and $b\in X$ can be considered by analogy with case 8.

13,14,15. The definition of $\overline\Af$ implies that the case $a,b\in\partial X$ and $(x\in X\vee y\in X)$ is impossible.

16. If $a,b,c,x\in\partial X$, then $\overline\Af axb\;\wedge\;\overline \Af ayb$
imply that for every point $p\in X$ there exist lines $A,A'\in a$, $B,B'\in b$, $L\in x$ and $\Lambda\in y$ such that $p\in A\cap A'\cap B\cap B'\cap L\cap \Lambda$, $L\subseteq\overline{A\cup B}$ and $\Lambda\subseteq\overline{A'\cup B'}$.  Since $a,b,x,y$ are spreads of lines, $A=A'$ and $B=B'$. Since $x\ne y$, the lines $L\in x$ and $\Lambda\in y$ are not parallel and hence $\overline{L\cup\Lambda}$ is a plane in the plane $\overline{A\cup B}$. The $3$-rankedness of $X$ ensures that $A\cup B\subseteq \overline{A\cup B}=\overline{L\cup\Lambda}$ witnessing that $\overline \Af xay\wedge \overline \Af xby$.
\end{proof}

\begin{definition} For a $3$-ranked liner $X$, the liner $\overline X=X\cup\partial X$ endowed with the line relation $\overline\Af$ is called  \index{spread completion}\defterm{the spread completion} of the liner $X$.
\end{definition}

The following lemma describes the structure of lines of the spread completions of proaffine regular liners.

\begin{lemma}\label{l:LinesInProjCompletion} Let $X$ be a proaffine regular liner and $\overline X$ be its spread completion.
\begin{enumerate}
\item For every spreading line $L$ in $X$, the set $L\cup\{L_\parallel\}$ is a line in the liner $\overline X$. 
\item If $\Lambda$ is a line in $\overline X$ and the set $L\defeq \Lambda\cap X$ is not empty, then $L$ is a spreading line in $X$ and $\Lambda=L\cup\{L_\parallel\}$.
\item If $\Lambda$ is a line in $\overline X$ and $\Lambda\cap X=\varnothing$, then for every  distinct directions ${\boldsymbol a},{\boldsymbol b}\in \Lambda$ and concurrent lines $A\in {\boldsymbol a}$ and $B\in {\boldsymbol b}$, we have $\Lambda=\{L_{\parallel}:L\in\mathcal S\;\wedge\;L\subset\overline{A\cup B}\}$, where $\mathcal S=\bigcup(\partial X)$ is the family of spreading lines in $X$.
\end{enumerate}
\end{lemma}

\begin{proof} 1, 2. The first two statements of the lemma follow immediately from the definition of the line relation $\overline\Af$. 
\smallskip

3. Assume that $\Lambda$ is a line in $\overline X$ such that $L\cap X=\varnothing$. Choose any two distinct directions ${\boldsymbol a},{\boldsymbol b}\in \Lambda\subseteq \overline X\setminus X=\partial X$ and  lines $A\in {\boldsymbol a}$ and $B\in {\boldsymbol b}$ such that $A\cap B$ contains some point $x$. Since ${\boldsymbol a}\ne {\boldsymbol b}$, the lines $A,B$ are not parallel and hence $\overline{A\cup B}$ is a plane in $X$. We have to prove that $\Lambda=\{L_\parallel:L\in\mathcal S\;\wedge\; L\subseteq \overline{A\cup B}\}$. By definition of the line relation $\overline\Af$, for every direction ${\boldsymbol c}\in \Aline {\boldsymbol a}{\boldsymbol b}$ there exists a line $C\in {\boldsymbol c}$ such that $A\cap B\cap C\ne\varnothing$ and $C\subseteq\overline{A\cup B}$, witnessing that $\Lambda\subseteq\{L_\parallel:L\in\mathcal S;\wedge\;L\subseteq \overline{A\cup B}\}$. To prove that   $\{L_\parallel:L\in\mathcal S\;\wedge\;L\subseteq \overline{A\cup B}\}\subseteq \Lambda$, take any direction $c\in \{L_\parallel:L\in\mathcal S\;\wedge\;L\subseteq P\}$. To prove that $\overline L{\boldsymbol a}{\boldsymbol c}{\boldsymbol b}$, take any lines $A'\in {\boldsymbol a}$, $B'\in {\boldsymbol b}$, and $C'\in {\boldsymbol c}$ such that $A'\cap B'\cap C'$ contains some point $x'$. By the choice of ${\boldsymbol c}$, there exists a line $L\in {\boldsymbol c}$ such that $L\subseteq \overline{A\cup B}$. Since ${\boldsymbol c}$ a spread of lines, there exists a unique line $C\in {\boldsymbol c}$ such that $x\in C$. It follows from $C\parallel L\subseteq P$ that $C\subseteq\overline{\{x\}\cup L}\subseteq \overline{A\cup B}$. By Theorem~\ref{t:subparallel-via-base}, $A'\parallel A$ and $B'\parallel B$ and $C'\parallel C\subseteq \overline{A\cup B}$ imply $\overline{A'\cup B'\cup C'}\subparallel \overline{A\cup B}\parallel \overline{A'\cup B'}$ and hence $C'\subseteq\overline{\{x'\}\cup A'\cup B'}=\overline{A'\cup B'}$, witnessing that $\overline L {\boldsymbol a}{\boldsymbol c}{\boldsymbol b}$ and hence ${\boldsymbol c}\in \Aline {\boldsymbol a}{\boldsymbol b}=\Lambda$. 

\end{proof}

\begin{proposition}\label{p:iso-spread-completion} For every isomorphism $A:X\to Y$ between $3$-ranked liners there exists a unique isomorphism $\bar A:\overline X\to\overline Y$ of the spead completions of the liners $X,Y$ such that $A=\bar A{\restriction}_X$.
\end{proposition}

\begin{proof} Define the function $\bar A:\overline X\to\overline Y$ assigning to every point $x\in X$ the point $y\defeq A(x)$ and to every direction ${\boldsymbol \delta}\in\partial X$, the direction $\bar A({\boldsymbol \delta})\defeq\{A[L]:L\in \{\boldsymbol \delta\}\in \partial Y$. The definition of the liner structure in the liners $\overline X$ and $\overline X$ implies that the function $\bar  A:\overline Y\to \overline Y$ is an isomorphism of the liners $\overline X$ and $\overline Y$. The definition of the function $B$ ensures that $A=\bar  A{\restriction}_X$. To prove the uniqueness of $\bar A$, assume that $B:\overline X\to \overline Y$ is another automorphism of the liners $\overline X$ and $\overline Y$ such that $A=B{\restriction}_X$. Then for every direction ${\boldsymbol \delta}\in\partial X$ and every line $L\in{\boldsymbol \delta}$, the set $\overline L\defeq L\cup\{{\boldsymbol \delta}\}$ is a unique line in $\overline X$ such that $L=\overline L\cap X$. Then $B[\overline L]=A[L]\cup \{B({\boldsymbol \delta})\}$ is a unique line in $\overline Y$ such that $B[\overline L]\cap Y=A[L]$. By the definition of the liner structure in the spread completion $\overline Y$, this unique line is equal to $A[L]\cup\{A[L]_\parallel\}=A[L]\cup \{\bar  A({\boldsymbol \delta})\}$ and hence $B({\boldsymbol \delta})=\bar A({\boldsymbol \delta})$, witnessing that $\bar  A=B$.
\end{proof}

Proposition~\ref{p:iso-spread-completion} justifies the following definition.

\begin{definition} For an isomorphism $A:X\to Y$ between $3$-ranked liners, its \index{isomorphism!spread completion of}\index{spread completion}\defterm{spread completion} is a unique isomorphism $\bar A:\overline X\to\overline Y$ such that $A=\bar A{\restriction}_X$. 
\end{definition}

\section{Completely regular liners}

\begin{definition}\label{d:comp-regular} A liner $X$ is called \index{completely regular liner}\index{liner!completely regular}\defterm{completely regular} if $X$ is $3$-ranked and its spread completion $\overline X$ is strongly regular.
\end{definition}

The following theorem is one of two main theorems characterizing completely regular liners (the other one is Theorem \ref{t:spread=projective2}).

\begin{theorem}\label{t:spread=projective1} A liner $X$ is completely regular if and only if $X$ is regular, para-Playfair, and bi-Bolyai.
\end{theorem}

\begin{proof} To prove the ``only if'' part, assume that the liner $X$ is completely regular.  Then $X$ is $3$-ranked and its spread completion $\overline X$ is strongly regular and hence projective. 

\begin{claim}\label{cl:comp-reg=>reg} The liner $X$ is regular.
\end{claim} 

\begin{proof} To show that $X$ is regular, fix any flat $A\subseteq X$,  and points $o\in A$ and $b\in X\setminus A$. Given any point $z\in \overline{A\cup\{b\}}\setminus(A\cup\Aline ob)$, we have to find points $x\in A$ and $y\in\Aline ob$ such that $z\in\Aline xy$. 
The definition of the line relation on the liner $\overline X$ ensures that the set $\overline A=A\cup\{L_\parallel:L\in\bigcup(\partial X)\;\wedge\;L\subseteq A\}$ is the flat hull of $A$ in $\overline X$. By the strong regularity of the projective liner $\overline X$, there exists a point $a\in\overline A$ such that $z\in \Aline ab$. If $a\in A$, then $z\in \Aline ab\cap X$ and we are done. So, assume that $a\in\overline A\setminus A$. Then $a$ is a spread of parallel lines in $X$ and $a$ contains a spreading line $L\in a$ such that $o\in L$.  It follows from $L\in a\in \overline A$  that $L\subseteq A$ and hence $\Aline ob\notin L_\parallel$. Proposition~\ref{p:spread-3-long} ensures that the line $\Aline ob$ contains a point $y\in \Aline ob\setminus\{o,b\}$.  By the projectivity of $\overline X$, the lines $\Aline yz$ and $\Aline oa$ in the plane $\overline{\{o,a,b\}}\subseteq\overline X$ have a common point $x$. Assuming that $x=a$, we conclude that $\{y\}=\Aline ob\cap\Aline xz=\Aline ob\cap\Aline az=\{b\}$, which contradicts the choice of $y$. This contradiction shows that $x\in \Aline oa\setminus \{a\}\subseteq X$.  Then $x,y$ are required points of $X$ with $z\in \Aline xy$, witnessing that the completely liner $X$ is regular.
\end{proof}

\begin{claim}\label{cl:comp-reg=>para-Playfair} The liner $X$ is para-Playfair.
\end{claim}

\begin{proof}  To prove that $X$ is para-Playfair, fix a plane $\Pi\subseteq X$ and two disjoint lines $L,\Lambda\subseteq \Pi$. We need to prove that for every $x\in \Pi\setminus L$ there exists a unique line $L_x$ in $X$ such that $x\in L_x\subseteq X\setminus L$. Let $\overline L,\overline \Lambda$ be the flat hulls of the lines $L,\Lambda$ in the projective liner $\overline X$. Since the flat hull $\overline \Pi$ of the plane $\Pi$ in the projective liner $\overline X$ is a plane and $\overline L\cup\overline\Lambda\subseteq\overline \Pi$, the lines $\overline L$ and $\overline\Lambda$ have a common point ${\boldsymbol a}$. Since $\{{\boldsymbol a}\}\cap X\subseteq \overline L\cap\overline \Lambda\cap X=L\cap\Lambda=\varnothing$, the point ${\boldsymbol a}$ belongs to the horizon $\partial X$ and hence ${\boldsymbol a}=L_\parallel=\Lambda_\parallel$, which means that the lines $L$ and $\Lambda$ are spreading. The spread $L_\parallel$ contains a unique line $L_x$ with $x\in L_x$. Assuming that $L'$ is another line in $X$ with $x\in L'\subseteq \Pi\setminus L$, we can apply Claim~\ref{cl:comp-reg=>reg} and Corollary~\ref{c:parallel-lines<=>} and conclude that $L'\in L_\parallel$ and hence $L'=L_x$.
\end{proof}

\begin{claim} The liner $X$ is bi-Bolyai.
\end{claim}

\begin{proof} We need to show that for every concurrent Bolyai lines $A,B$ in $X$, every line $L$ in the plane $\overline{A\cup B}$ is Bolyai. By Claims~\ref{cl:comp-reg=>reg} and \ref{cl:comp-reg=>para-Playfair}, the completely regular liner $X$ is regular and para-Playfair. By Proposition~\ref{p:Playfair=>para-Playfair=>Proclus}, the para-Playfar liner $X$ is Proclus, and by Theorem~\ref{t:Proclus<=>}, the Proclus liner is proaffine. By Theorem~\ref{t:spreading<=>}, the Bolyai lines $A,B$ in $X$ are spreading. Since the spreadin glines $A,B$ are concurrent, their spreads $a\defeq A_\parallel$ and $b\defeq B_\parallel$ are two distinct points of the spead completion $\overline X$ of $X$. Then the flat hull $\Aline ab$ of the doubleton $\{a,b\}$ is a line in the projective liner $\overline X$. Taking into account that the  the set $\overline X\setminus X=\partial X$ is flat in $\overline X$, we conclude that $\Aline ab\subseteq\partial X$. Since $X$ is a subliner of the liner $\overline X$, the flat hull $\overline L$ of the line $L$ in the projective liner $\overline X$ is a line in $\overline X$ such that $\overline L\cap X=L$. Since $\overline {A\cup B}$ is a plane in $X$, its flat hull in $\overline X$ is a plane containing the lines $\Aline ab$ and $\overline L$. By the projectivity of $\overline X$, the lines $\Aline ab$ and $\overline L$ have a common point $c\in \Aline ab\subseteq\partial X$. It follows from $c\in\overline L=\{c\}\cup L$ that $L$ is a spreading line in $X$. By Theorem~\ref{t:spreading<=>}, the spreading line $L$ is Bolyai.
\end{proof}

To prove the ``if'' part, assume that the liner $X$ is regular, para-Playfair, and bi-Bolyai. To prove that the liner $X$ is completely regular, it suffices to check that the spread completion $\overline X$ of $X$ is $0$-parallel, see Theorem~\ref{t:projective<=>}. 
Assuming that $\overline X$ is not $0$-parallel, find a plane $\Pi$ in $\overline X$ containing two disjoint lines.  Since $\Pi$ is a plane in $\overline X$, there exist three points $a,b,c\in\overline X$ such that $\Pi\subseteq\overline{\{a,b,c\}}$ in $\overline X$. Since $\Pi$ is a plane, the set $\{a,b,c\}$ has rank 3 in the liner $\overline X$.

\begin{claim}\label{cl:ProjectiveCompletion1} If $\{a,b,c\}\cap X\ne\varnothing$, then there exists a plane $P$ in $X$ such that $\Pi\subseteq\overline{\{a,b,c\}}\subseteq P\cup\{L_{\parallel}:L\in\mathcal S\;\wedge\;L\subset P\}$, where $\mathcal S$ is the family of spreading lines in $X$.
\end{claim}

\begin{proof} If $a,b,c\in X$, then the plane $P\defeq\overline{\{a,b,c\}}$ has the desired property because the set $P\cup\{L_\parallel:L\in\mathcal S\;\wedge\;L\subset P\}$ is flat in the liner $\overline X$ and $a,b,c\in P$.
\smallskip

If $|\{a,b,c\}|\cap X=2$, then we lose no generality assuming that $a,b\in X$ and $c\notin X$. Then $c$ is a spread of lines containing a line $C\in c$ such that $a\in C$. Since the set $\{a,b,c\}$ has rank 3 in $\overline X$, the point $b$ does not belong to the line $C$ and hence  the flat $P\defeq\overline{C\cup\{b\}}$ is a plane in $X$ having the desired property. 
\smallskip

 If $|\{a,b,c\}\cap X|=1$,  then we  lose no generality assuming that $a\in X$ and $b,c\notin X$. Since $b,c$ are spreads of lines, there exist unique lines $B\in b$ and $C\in c$ such that $a\in B\cap C$. Since the set $\{a,b,c\}$ has rank 3  in $(\overline X,\overline\Af)$, the spreads $b,c$ are distinct, the lines $B,C$ are distinct and hence the flat $P\defeq\overline{B\cup C}$ is a plane in $X$ having the desired property. 
\end{proof}

\begin{claim}\label{cl:abc} $\{a,b,c\}\subseteq \overline X\setminus X$.
\end{claim}

\begin{proof} In the opposite case, $\{a,b,c\}\cap X\ne\varnothing$. By Claim~\ref{cl:ProjectiveCompletion1},  there exists a plane $P$ in the liner $X$ such that $\Pi\subseteq\overline{\{a,b,c\}}\subseteq P\cup\{L_\parallel:L\in\mathcal S\;\wedge\;L\subset P\}$. By the choice of the plane $\Pi$, there exist two disjoint lines $A$ and $B$ in $\Pi\subseteq  P\cup\{L_\parallel:L\in\mathcal S\;\wedge\;L\subset P\}$. 

If $A\cap X\ne\varnothing \ne B\cap X$, then by the definition of the line relation $\overline \Af$, the intersections $A\cap X$ and $B\cap X$ are disjoint spreading lines in the plane $P$. By Corollary~\ref{c:parallel-lines<=>}, the disjoint lines $A\cap X,B\cap X$ in the plane $P$ are parallel  and hence $(A\cap X)_{\parallel}=(B\cap X)_\parallel\in A\cap B=\varnothing$, which is a desired contradiction.

If $A\cap X\ne\varnothing=B\cap X$, then $A\cap X\subseteq P$ and $B\subseteq\overline X\setminus X=\partial X$. Choose any distinct spreads $u,v\in B$, any point $x\in P$, and find unique spreading lines $U\in u$ and $v\in V$ such that $x\in U\cap V$. The choice of the plane $P$ ensures that $U\cup V\subseteq P$ and hence $\overline{U\cup V}=P$, by the $3$-rankedness of $X$. By our assumption, the line $A\cap X$ in the plane $\overline{U\cup V}=P$ is spreading. Then its spread $(A\cap X)_\parallel\in A$ belongs to the line $B=\Aline uv=\{L_\parallel:L\in\mathcal S\;\wedge\; L\subseteq\overline{U\cup V}\}$, according to Lemma~\ref{l:LinesInProjCompletion}. Therefore, $(A\cap X)_\parallel\in A\cap B=\varnothing$,  which is a desired contradiction.

Finally, assume that $A\cap X=\varnothing=B\cap X$. Applying Lemma~\ref{l:LinesInProjCompletion}, we can show that $A=\{L_\parallel:L\in\mathcal S\;\wedge\;L\subseteq P\}=B$, which contradicts $A\cap B=\varnothing$.

Those contradictions complete the proof of the inclusion $\{a,b,c\}\subseteq\overline X\setminus X$.
\end{proof}

By Claim~\ref{cl:abc}, $\{a,b,c\}\subseteq\overline X\setminus X$ and hence $\Pi\cap X\subseteq \overline{\{a,b,c\}}\cap X=\varnothing$. Fix any point $x\in X$ and find unique lines $A\in a$, $B\in b$, $C\in c$ such that $x\in A\cap B\cap C$. Then $\|A\cup B\cup C\|\le 4$. Assuming that $\|A\cup B\cup C\|<4$, we conclude that the set $\{a,b,c\}\subseteq \overline{A\cup B\cup C}\setminus X$ has rank $<3$ in $\overline X$, which contradicts the choice of the set $\{a,b,c\}$. This contradiction shows that $\|A\cup B\cup C\|=4$. It follows that $\Pi\subseteq\overline{\{a,b,c\}}\subseteq\overline{A\cup B\cup C}\setminus X\subseteq\{L_\parallel:L\in\mathcal S\;\wedge\;L\subseteq \overline{A\cup B\cup C}\}$. By our assumption, the plane $\Pi$ contains two disjoint lines $L$ and $\Lambda$. Fix any distinct points $p,q\in L$ and $s,t\in \Lambda$. The points $p,q,s,t$ are spreads of parallel lines in $X$. Then there exist uniques lines $P\in p$, $Q\in q$, $S\in s$, $T\in t$ such that $x\in P\cap Q\cap S\cap T$. By Corollary~\ref{c:proregular=>ranked} and Theorem~\ref{t:w-modular<=>}, the regular liner $X$ is locally modular and hence
\begin{multline*}
\|\overline{P\cup Q}\cap\overline{S\cup T}\|=\|\overline{P\cup Q}\|+\|\overline{S\cup T}\|-\|\overline{P\cup Q}\cup\overline{S\cup T}\|\\
=3+3-\|P\cup Q\cup S\cup T\|\ge 6-\|A\cup B\cup C\|=6-4=2,
\end{multline*}
which means that the intersection $\overline{P\cup Q}\cap\overline{S\cup T}\ni x$ contains some line $\Gamma\ni x$. By our assumption, the line $\Gamma$ in the plane $\overline{P\cup Q}$ is spreading.  Then $\Gamma_\parallel\in\Aline pq\cap\Aline st=L\cap\Lambda=\varnothing$, which is a desired contradiction completing the proof of the $0$-parallelity of the liner $\overline X$. 
\end{proof}   

Theorems~\ref{t:spread=projective1} and \ref{t:w-modular<=>} imply that for every liner we have the implications:
$$\mbox{strongly regular}\Ra\mbox{completely regular}\Ra\mbox{regular}\Ra\mbox{weakly modular}\Ra\mbox{weakly regular}.$$

Theorems~\ref{t:spread=projective1}, \ref{t:4long+pP+bB=>regular}, and Corollary~\ref{c:proregular=>ranked} imply the following characterization of $4$-long completely regular liners.

\begin{corollary} A $4$-long liner is completely regular if and only if it is para-Playfair and bi-Bolyai.
\end{corollary}

\begin{theorem}\label{t:proaffine3=>compregular} Every proaffine regular liner $X$ of rank $\|X\|\ne 3$ is completely regular.
\end{theorem}

\begin{proof} By Theorems~\ref{t:4-para-Playfair} and \ref{t:w-modular<=>}, the proaffine regular liner $X$ is para-Playfair and weakly modular. Let $\overline X=X\cup\partial X$ be the spread completion of $X$ and $\mathcal S=\bigcup(\partial X)$ be the family of spreading lines in $X$. If $\|X\|\le 1$, then $X$ contains no lines and hence $\overline X=X$ is $0$-parallel and strongly regular. If $\|X\|=2$, then $X$ is a projective line. In this case $\mathcal S=\{X\}$ and $\overline X=X\cup\{\{X\}\}$ is a projective line.
So, assume that $\|X\|\ge 4$. 

We claim that the liner $X$ is bi-Bolyai. 
Given any concurrent Bolyai lines $L,L'$ in $X$, we need to show that every line in the plane $P\defeq\overline{L\cup L'}$ is Bolyai. Since $\|X\|\ge 4>3=\|P\|$, there exists a point $x\in X\setminus P$. By Theorem~\ref{t:spreading<=>}, the Bolyai lines $L,L'$ are spreading. Then there exist unique lines $L_x\in L_\parallel$ and $L_x'\in L'_\parallel$ such that $x\in L_x\cap L_x'$ and hence $\|L\cup L'\cup L_x\cup L_x'\|=\|L\cup L'\cup\{x\}\|=4$. 

Given any line  $\Lambda$ on the plane $P$, consider the plane $\overline{\Lambda\cup\{x\}}$. Since $\|\overline{L_x\cup L_x'}\cup\overline{\Lambda\cup\{x\}}\|= \|L\cup L'\cup\{x\}\|\le 4$, the intersection $\Lambda_x\defeq \overline{\Lambda\cup\{x\}}\cap\overline{L_x\cup L_x'}$ is a line, by the weak modularity of the liner $X$. Then $\Lambda$ and $\Lambda_x$ are two disjoint lines in the plane $\overline{\Lambda\cup\{x\}}$. By Proposition~\ref{p:lines-in-para-Playfair}, the line $\Lambda$ is Bolyai, witnesing the para-Playfair regular liner $X$ is bi-Bolyai. By Theorem~\ref{t:spread=projective1}, the bi-Bolyai para-Playfair regular liner $X$ is completely regular.
\end{proof}

\begin{corollary}\label{c:affine-spread-completion} Every  affine regular liner $X$ is completely regular.
\end{corollary}

\begin{proof} If $\|X\|\ne 3$, then $X$ is completely regular, by Theorem~\ref{t:proaffine3=>compregular}. So, assume that $\|X\|=3$. By Theorem~\ref{t:affine=>Avogadro}, the affine liner $X$ is $2$-balanced. If $|X|_2=2$, then $X$ contains no parallel lines and hence $\overline X=X$ is projective. If  $|X|_2=3$, then the affine regular liner  $X$ is $3$-long. By Theorem~\ref{t:Playfair<=>}, the $3$-long affine regular liner $X$ is Playfair and hence para-P ayfairand bi-Bolyai. By Theorems~\ref{t:Playfair<=>spreading} and \ref{t:spread=projective1}, the liner $X$ is completely regular.
\end{proof}

\begin{proposition}\label{p:spread-3long} If a liner $X$ is completely regular and not projective, then the spread completion $\overline X$ of $X$ is projective and $3$-long.
\end{proposition}

\begin{proof} Assume that a liner $X$ is completely regular. Then its spread completion $\overline X$ is strongly regular and hence projective, by Theorem~\ref{t:projective<=>}. If $X$ is not projective, then $X\ne\overline X$ and hence $X$ contains a spreading line $S$.  To see that $\overline X$ is $3$-long, take any line $L$ in $\overline X$. If $L\cap X\ne\varnothing$, then the definition of the line relation $\overline \Af$ on the liner $\overline X$ ensures that $L\cap X$ is a line in $X$.
If $L\cap X\notin S_\parallel$, then $|L|\ge|L\cap X|\ge 3$, by Proposition~\ref{p:spread-3-long}. If $L\cap X\in S_\parallel$, then $L=(L\cap X)\cup S_\parallel$ and hence $|L|=|L\cap X|+1\ge 3$.

 Next, assume that $L\cap X=\varnothing$ and choose two distinct directions ${\boldsymbol a},{\boldsymbol b}\in L\subseteq \overline X\setminus X$. Then ${\boldsymbol a},{\boldsymbol b}$ are two spreads of parallel spreading lines in $X$. Choose any point $x\in X$ and find two spreading lines $A\in {\boldsymbol a}$ and $B\in {\boldsymbol b}$ such that $x\in A\cap B$. Since ${\boldsymbol a}\ne {\boldsymbol b}$, the lines $A,B$ are not parallel. Then $A\cap B=\{x\}$ and there exist points $y\in A\setminus B$ and $z\in B\setminus A$. Since the liner $\overline X$ is projective, the lines $\Aline {\boldsymbol a}{\boldsymbol b}$ and $\Aline yz$ in the plane  $P=\overline{\{x,y,z\}}\subseteq \overline X$ have a common point ${\boldsymbol c}\in\Aline yz\cap \Aline {\boldsymbol a}{\boldsymbol b}\subseteq \overline X\setminus X$. Since the line $\Aline yz\cap X$ is concurrent with the lines $A,B$, the direction ${\boldsymbol c}$ is distinct from the directions ${\boldsymbol a},{\boldsymbol b}$. Then $|L|=|\Aline {\boldsymbol a}{\boldsymbol b}|\ge|\{{\boldsymbol a},{\boldsymbol b},{\boldsymbol c}\}|=3$, witnessing that the liner $\overline X$ is $3$-long. 
\end{proof} 

Proposition~\ref{p:spread-3long} implies the following corollary.

\begin{corollary}\label{c:spread-3long} The spread completion of any $3$-long completely regular liner $X$ is projective and $3$-long.
\end{corollary}

\section{Complete properties of liners}

In this section we develop some terminology related to properties of liners.
In the sequel we shall strive to define properties of liners which are complete in the following sense.

\begin{definition} A property $\mathcal P$ of liners is called \index{property!complete}\index{complete property}\defterm{complete} if a completely regular liner $X$ has property $\mathcal P$ if and only if the spread completion $\overline X$ of $X$ has property $\mathcal P$.
\end{definition}

\begin{exercise} Show that the properties of a liner to be completely regular,  regular, proregular, weakly regular, weakly modular, ranked, Proclus, para-Playfair, proaffine, hyperbolic are complete.
\end{exercise}

\begin{exercise} Show that the properties of a liner to be strongly regular, modular, projective, affine, hyperaffine, injective, Playfair, Bolyai, Lobachevsky, Boolean are not complete. 
\end{exercise}


There is a simple method of turning a property of projective liners to a complete property of arbitrary (not necessarily projective) liners.

\begin{definition} Let $\mathcal P$ be a property of projective liners. A liner $X$ is called \index{completely $\mathcal P$ liner}\index{liner!completely $\mathcal P$}\defterm{completely $\mathcal P$} if $X$ is completely regular and its spread completion $\overline X$ has property $\mathcal P$.
\end{definition}

\begin{proposition} For every property $\mathcal P$ of projective liners, the property of liners to be completely $\mathcal P$ is complete.
\end{proposition}

\begin{proof} Let $X$ be a completely regular liner and $\overline X$ be its spread completion, which is strongly regular, by Definition~\ref{d:comp-regular}. By Theorem~\ref{t:modular<=>} and Corollary~\ref{c:proregular=>ranked}, the strongly regular liner $\overline X$ is projective and ranked. Since the projective liner $\overline X$ contains no spreading lines, the spread completion of $\overline X$ coincides with $\overline X$ and hence the strongly regular liner $\overline X$ is completely regular.

  We have to show that the liner $X$ is completely $\mathcal P$ if and only if its spread completion $\overline X$ is completely $\mathcal P$.

If $X$ is completely $\mathcal P$, then its spread completion $\overline X$ has the property $\mathcal P$. Since the liner $\overline X$ is ranked and projective, it is completely regular and its spead completion coincides with $\overline X$, which implies that the liner $\overline X$ is completely $\mathcal P$.

If the liner $\overline X$ is completely $\mathcal P$, then its spread completion (which is $\overline X$) has the property $\mathcal P$ and then, the liner $X$ is completely $\mathcal P$.
\end{proof}

In the sequel, we shall often consider instances of the following general problem.

\begin{problem}\label{prob:inner-completely-P} Given a property $\mathcal P$ of projective liners, find an inner characterization of completely $\mathcal P$ liners.
\end{problem}

\begin{remark} Corollaries~\ref{c:completelyD<=>}, \ref{c:completelyP<=>}, \ref{c:completelyF<=>}, \ref{c:completelyM<=>} answer Problem~\ref{prob:inner-completely-P} for the properties of projective liners to be Desarguesian, Pappian, Fano, Moufang.
\end{remark}

Next, we consider two methods of turning a property of affine liners into a property of projective liners. By Proposition~\ref{p:projective-minus-hyperplane}, for every hyperplane $H$ in a projective liner $Y$, the subliner $X\defeq Y\setminus H$ is affine and regular.

\begin{definition} Let $\mathcal P$ be a property of affine liners. A projective liner $X$ is defined to be 
\begin{enumerate}
\item \index{everywhere $\mathcal P$ projective liner}\index{projective liner!everywhere $\mathcal P$}\defterm{everywhere $\mathcal P$} if for every hyperplane $H\subset X$, the affine liner $X\setminus H$ has property $\mathcal P$;
\item \index{somewhere $\mathcal P$ projective liner}\index{projective liner!somewhere $\mathcal P$}\defterm{somewhere $\mathcal P$} if for some hyperplane $H\subset X$, the affine liner $X\setminus H$ has property $\mathcal P$.
\end{enumerate}
\end{definition}

\begin{definition} A property $\mathcal P$ of liners is called a \defterm{liner property} if for any isomorphic liners $X,Y$, the liner $X$ has property $\mathcal P$ if and only if $Y$ has property $\mathcal P$.
\end{definition}

\begin{theorem}\label{t:everywhereP<=>somewhereP} Let $\mathcal P$ be a complete liner property. For any  $4$-long projective liner $X$, the following conditions are equivalent:
\begin{enumerate}
\item $X$ has property $\mathcal P$;
\item $X$ is everywhere $\mathcal P$;
\item $X$ is somewhere $\mathcal P$.
\end{enumerate}
\end{theorem}

\begin{proof} Let $X$ be a $4$-long projective liner. 
\smallskip

$(1)\Ra(2)$ Assume that the projective liner $X$ has property $\mathcal P$. Given any hyperplane $H\subset X$, we need to show that the affine liner $A\defeq X\setminus H$ has property $\mathcal P$.
 By Proposition~\ref{p:projective-minus-hyperplane}, the subliner $A=X\setminus H$ is affine and regular. Since the projective liner $X$ is $4$-long, the affine liner $A$ is $3$-long. By Corollary~\ref{c:affine-spread-completion}, the affine regular liner $A$ is completely regular. By Corollary~\ref{c:pcompletion=scompletion}, the identity inclusion $A\to X$ extends to an isomorphism $\Phi:\overline A\to X$. Since $\mathcal P$ is a liner property, the projective liner $\overline A$ has property $\mathcal P$, being an isomorphic copy of the projective liner $X$ with property $\mathcal P$. Since the spread completion $\overline A$ of the affine liner $A$ has property $\mathcal P$ and the property $\mathcal P$ is complete, the affine liner $A$ has property $\mathcal P$.
 \smallskip
 
The implication $(2)\Ra(3)$ is obvious.
 \smallskip
 
$(3)\Ra(1)$ Assume that the projective liner $X$ is somewhere $\mathcal P$. Then there exists a hyperplane $H\subset X$ such that the affine liner $A\defeq X\setminus H$ has property $\mathcal P$. 
 By Proposition~\ref{p:projective-minus-hyperplane}, the subliner $A=X\setminus H$ is affine and regular. Since the projective liner $X$ is $4$-long, the affine liner $A$ is $3$-long. By Corollary~\ref{c:affine-spread-completion}, the affine regular liner $A$ is completely regular. Since the completely regular liner $A$ has property $\mathcal P$, its spread completion $\overline A$ has property $\mathcal P$, by the completeness of $\mathcal P$. By Corollary~\ref{c:pcompletion=scompletion}, the identity inclusion $A\to X$ extends to an isomorphism $\Phi:\overline A\to X$. Since $\mathcal P$ is a  linear property, the projective liner $X$ has property $\mathcal P$, being an isomorphic copy of the projective liner $\overline A$ possessing the property $\mathcal P$. 
\end{proof} 

\begin{remark} Everywhere Boolean projective liners are called Fano and will be studied in Chapter~\ref{ch:Fano}.
\end{remark}

\chapter{Projective completions of  liners}\label{s:completions}

In this chapter we shall study completions of liners.

\section{Projective and normal completions of liners}\label{s:proj-completions}


\begin{definition}\label{d:procompletion} A \index{completion}\index{liner!completion of}\defterm{completion} of a liner $X$ is a $3$-long liner $Y$ that contains $X$ as a subliner such that $\overline{Y\setminus X}\ne Y$. If the liner $Y$ is projective, then $Y$ is called a \index{projective completion}
\index{liner!projective completion of}\defterm{projective completion} of the liner $X$. If $X$ is a line, then we shall additionally require that the remainder $Y\setminus X$ is not empty, in order to make the projective completion of a line unique. 
\end{definition}

\begin{example} By Proposition~\ref{p:spread-3long}, the spread completion $\overline X$ of any (completely regular) nonempty $3$-long $3$-ranked  liner $X$ is a (projective) completion of $X$.
\end{example} 

\begin{exercise} Let $Y$ be a projective completion of a liner $X$. Prove that $Y=X$ if $\|X\|\le 1$. 
\end{exercise}

\begin{exercise} Prove that a projective completion $Y$ of a projective liner $X$ of rank $\|X\|\ne 2$ coincides with $X$.
\end{exercise}

\begin{proposition}\label{p:trace-of-line} Let $Y$ be a completion of a liner $X$. For every line $L\subseteq Y$ with $L\not\subseteq\overline{Y\setminus X}$, the intersection $L\cap X$ is a line in the liner $X$.
\end{proposition}

\begin{proof} Since $L\not\subseteq \overline{Y\setminus X}$, the intersection $L\cap\overline{Y\setminus X}$ contains at most one point. Since the liner $Y$ is $3$-long, the set $L\setminus\overline{Y\setminus X}\subseteq L\setminus(Y\setminus X)=L\cap X$ contains two distinct points $x,y$ and hence $L\cap X=\Aline xy\cap X$ is a line in $X$, by the definition of the induced line structure on the subliner $X$.
\end{proof} 

\begin{exercise}\label{ex:Fano-trace} Find an example of a liner $X$ and a projective completion $Y$ of $X$ such that the liner $Y$ contains a plane $P$ such that $P\not\subseteq\overline{Y\setminus X}$ but $P\cap X$ is not a plane in $X$.
\smallskip

{\em Hint:} Let $Y$ be a Steiner projective plane, $L$ be a line in $Y$ and $X\defeq Y\setminus L$. Then $Y$ is a projective completion of $X$, and $P\defeq Y$ is a plane such that $P\cap X=X$ is not a plane, because every line in $X$ contains exactly two elements and so $X$ has rank $\|X\|=|X|=4>3$.
\end{exercise}

Proposition~\ref{p:trace-of-line} and Exercise~\ref{ex:Fano-trace} suggest the following natural definition.

\begin{definition} A completion $Y$ of a liner $X$ is called \index{normal completion}\index{completion!normal}\defterm{normal} if for every plane $P\subseteq Y$ with $P\not\subseteq\overline{Y\setminus X}$, the intersection $P\cap X$ is a plane in $X$.
\end{definition}

\begin{proposition}\label{p:pcomp3=>normal} Any projective completion $Y$ of any liner $X$ of rank $\|X\|\le 3$ is normal.
\end{proposition}

\begin{proof} Given any plane $P\subseteq Y$, observe that $3=\|P\|\le \|Y\|\le \|X\|\le 3$ and hence $P=Y$, by the rankendness of the projective linr $Y$. Then $P\cap X=X$ is a plane in $X$. 
\end{proof}

Next, we shall show that any projective completion of a $3$-long liner is normal.
This important fact will be deduced from the following lemma, which will be used also in the proof of Proposition~\ref{p:projective-minus-flat}.

\begin{lemma}\label{l:trace-flat} Let $Y$ be a projective completion of a $3$-long liner $X$, $A$ be a flat in $X$ and  $\overline A$ be the flat hull of $A$ in the liner $Y$. If $A\not\subseteq \overline{Y\setminus X}$, then $A=\overline A\cap X$.
\end{lemma}

\begin{proof} The inclusion $A\subseteq \overline A\cap X$ is obvious. So, it remains to prove that $\overline A\cap X\subseteq A$. Since $A\not\subseteq \overline{Y\setminus X}$, there exists a point $o\in A\setminus\overline{Y\setminus X}$. Consider the set $$\Aline oA\defeq\bigcup_{a\in A}\Aline oa$$ in the projective liner $Y$.  We claim that the set $\Aline oA$ is flat in $Y$. Given any distinct points $u,v\in\Aline oA$, we should prove that $\Aline uv\subseteq \Aline oA$. Since $u,v\in\Aline oA$, there exist points $x,y\in A$ such that $u\in\Aline ox$ and $v\in\Aline oy$. If $u\in\Aline oy$, then $\Aline uv\subseteq\Aline oy\subseteq \Aline oA$ and we are done. By analogy we can show that $\Aline uv\subseteq \Aline oA$ if $v\in\Aline ox$. So, we assume that $u\notin\Aline oy$ and $v\notin\Aline ox$ and hence $x\ne o\ne y$. To show that $\Aline uv\subseteq\Aline oA$, fix any point $w\in \Aline uv\setminus\{u,v\}$. Since $o\notin\overline{Y\setminus X}$ and $\overline{Y\setminus X}$ is a flat in $Y$, the intersection $\Aline ow\cap \overline{Y\setminus X}$ contains at most one point. Since the liner $Y$ is 3-long, there exists a point $z\in\Aline ow\setminus (\{o\}\cup \overline{Y\setminus X})\subseteq (\Aline ow\cap X)\setminus\{o\}$. 
Assuming that $z\in\Aline ox$, we conclude that $v\in \Aline uw\subseteq \overline{\{u,o,z\}}\subseteq \overline{\{u,o,x\}}=\Aline ox$, which contradicts our assumption. This contradiction shows that $z\notin\Aline ox$. By analogy we can prove that $z\notin\Aline oy$. 

\begin{picture}(200,80)(-150,-10)

\put(0,0){\line(1,0){90}}
\put(0,0){\line(0,1){60}}
\put(30,0){\line(-1,1){30}}
\put(60,0){\line(-1,1){60}}
\put(0,45){\line(2,-1){90}}
\put(0,0){\line(1,1){30}}

\put(0,0){\circle*{3}}
\put(-2,-8){$o$}
\put(30,0){\circle*{3}}
\put(28,-8){$u$}
\put(60,0){\circle*{3}}
\put(58,-8){$x$}
\put(90,0){\circle*{3}}
\put(88,-8){$p$}
\put(0,30){\circle*{3}}
\put(-8,28){$v$}
\put(0,45){\circle*{3}}
\put(-11,43){$y_p$}
\put(0,60){\circle*{3}}
\put(-8,58){$y$}
\put(30,30){\circle*{3}}
\put(30,33){$z$}
\put(15,15){\circle*{3}}
\put(19,12){$w$}

\end{picture}

Since $z\notin\Aline ox$,  for every $p\in \Aline  ox$, the flat $\Aline pz$ is a line in the plane $\overline{\{o,x,y\}}$. Since the space $Y$ is projective, the line $\Aline pz\subseteq\overline{\{o,x,y\}}$ intersects the line $\Aline oy=\Aline ov$ at some point $y_p$. It follows from $z\notin\Aline ox\cup\Aline oy$ that $\{p\}=\Aline ox\cap\Aline z{y_p}$. So, the map $\Aline ox\to\Aline oy$, $p\mapsto y_p$, is injective. Taking into account that the set $\overline{Y\setminus X}$ is flat in $X$ and $o\notin\overline{Y\setminus X}$ we conclude that $|\overline{Y\setminus X}\cap\Aline oy|\le 1$. Then at most one point $p\in\Aline ox$ has $y_p\in\overline{Y\setminus X}$. Since the liner $Y$ is  $3$-long and $|\overline{Y\setminus X}\cap\Aline ox|\le 1$, the intersection $\Aline ox\cap X$ is a line in the liner $X$. Since the liner $X$ is $3$-long,  there exists a point $p\in(\Aline ox\cap X)\setminus\{o\}$ such that $y_p\notin \overline{Y\setminus X}$. Then $\{p,y_p\}\subseteq(\Aline ox\cup\Aline oy)\setminus \overline{Y\setminus X}\subseteq A$ and $z\in\Aline p{y_p}\setminus \overline{Y\setminus X}\subseteq A$ by the flatness of the set $A$ in the liner $X$. Then $w\in\Aline oz\subseteq\Aline oA$, witnessing that the set $\Aline oA$ is flat in $Y$.  Then $\overline A\subseteq \Aline aA$ and $\overline A\cap X\subseteq \Aline oA\cap X=A$, by the flatness of the set $A$ in the liner $X$. 
\end{proof}

\begin{theorem}\label{t:procompletion=>normal} Any projective completion $Y$ of a $3$-long liner $X$ is normal.
\end{theorem}

\begin{proof} Given any plane $P\subseteq Y$ with $P\not\subseteq \overline{Y\setminus X}$, we should prove that $P\cap X$ is a plane in the liner $X$. Since $P\not\subseteq \overline{Y\setminus X}$, there exists a point $o\in P\setminus\overline{Y\setminus X}$. Since $P$ is a plane in $Y$, there exist points $a\in P\setminus\{o\}$ and $b\in P\setminus\Aline oa$. Since $o\notin\overline{Y\setminus X}$, the intersection $\Aline oa\cap\overline{Y\setminus X}$ contains at most one point. Since the liner $Y$ is $3$-long, there exists a point $x\in \Aline oa\setminus(\{o\}\cup \overline{Y\setminus X})\subseteq \Aline oa\cap X$. By analogy we can find a point $y\subseteq(\Aline ob\cap X)\setminus\{o\}$. It follows from $b\notin\Aline oa$ that $\Aline ox\cap\Aline oy=\Aline oa\cap\Aline ob=\{o\}$ and hence $y\notin\Aline ox$. Then the set $\{x,o,y\}$ has rank $3$ in the liner $X$ and hence its flat hull $A$ in the liner $X$ is a plane in $X$. Lemma~\ref{l:trace-flat} ensures that $A=\overline{A}\cap X$. 

By Theorem~\ref{t:projective<=>}, the projective liner $Y$ is (strongly) regular, and by Corollary~\ref{c:proregular=>ranked}, the regular liner $Y$ is ranked. Since the set $\{x,o,y\}\subseteq P$ has rank $3$ in the liner $Y$, its flat hull $\overline{\{x,y,z\}}$ in the liner $Y$ coincides with the plane $P$. Then $P=\overline{\{x,o,y\}}\subseteq\overline{A}\subseteq P$ and hence $\overline A=P$ and $P\cap X=\overline{A}\cap X=A$ is a plane in $X$, witnessing that the liner $Y$ is a normal completion of the liner $X$.
\end{proof}

\begin{proposition}\label{p:spread-comp-normal} The spread completion $\overline X$ of any  $3$-ranked liner $X$ is normal.
\end{proposition}

\begin{proof} Given any plane $P\subseteq \overline X$ with $P\cap X\ne\varnothing$, we should prove that $P\cap X$ is a plane in $X$. Since $P$ is a plane, there exist non-collinear points $a,b,c\in \overline X$ such that $P\subseteq\overline{\{a,b,c\}}$. It follows from $P\cap X\ne\varnothing$ that $\{a,b,c\}\cap X\ne\varnothing$. We lose no generality assuming that $a\in X$. If $b\in X$, then let $B\defeq\{b\}$. If $b\notin X$, then $b$ is a spread of lines in $X$, containing a unique line $B$ such that $a\in B$. By analogy define the set $C$. If $c\in X$, then put $C\defeq\{c\}$. If $c\notin X$, then let $C\in c$ be a unique line such that $a\in C$. Observe that the flat hull $\Pi\defeq\overline{\{a\}\cup B\cup C}$ of the set $\{a\}\cup B\cup C$ in $X$ is a plane.  The definition of the line relation on the liner $\overline X$ ensures that the set $\tilde \Pi\defeq P\cup\{L_\parallel:L\in\mathcal S,\;L\subseteq \Pi\}$ is flat in $\overline X$. Taking into account that this flat contains the set $\{a,b,c\}$, we conclude that $P\subseteq\tilde \Pi$ and hence $P\cap X\subseteq\tilde \Pi\cap X=\Pi$. Then $3=\|\{c\}\cup A\cup B\|\le \|P\cap X\|\le\|\Pi\|\le 3$, which means that $P\cap X$ is a plane in $X$.
\end{proof} 

Next, we prove that every projective completion of a non-projective Proclus liner is normal. For this we need the following theorem describing the structure of non-projective Proclus planes which are not $3$-long.

Let us recall that a projective plane $X$ is called \defterm{Steiner} if $|X|_2=3$ (i.e., each line in $X$ contains exactly three points). A Steiner projective plane contains $7$ points, $7$ lines and is unique up to an isomorphism.

\begin{theorem}\label{t:Proclus-not-3long} If a Proclus plane $X$ is not projective and not $3$-long, then $X$ has a projective completion and $|X|\le 6$. More precisely, $X=P\setminus H$ for some Steiner projective plane $P$  and some  set $H\subseteq P$ of cardinality $|H|\in\{1,2\}$. If $X$ is regular or para-Playfair, then $|X|=6$ and $|H|=1$.
\end{theorem}

\begin{proof} Assume that a Proclus plane $X$ is not projective and not $3$-long.   By Theorem~\ref{t:Proclus<=>}, the Proclus liner $X$ is proaffine and $3$-proregular. By Proposition~\ref{p:k-regular<=>2ex}, the $3$-proregular liner $X$ is $3$-ranked. 
Since the liner $X$ is not $3$-long, there exists a line $A$ of length $|A|=2$. Since $X$ is a plane, there exists a point $b\in X\setminus L$. Since $X$ is a $3$-ranked plane, $X=\overline{L\cup\{b\}}$. Assuming that $\Aline xb=\{x,b\}$ for every $x\in A$, we conclude that the plane $X=\overline{A\cup\{b\}}=A\cup\{b\}$ is projective, which contradicts our assumption. Therefore, there exists a point $o\in A$ such that the line $L\defeq\Aline ob$ has cardinality $|L|>2$. Let $a$ be a unique point of the singleton $A\setminus\{o\}$. If $X=L\cup A=L\cup\{a\}$, then the liner $X$ is projective, which contradicts our assumption. This contradiction shows that $X=\overline{L\cup A}\ne L\cup A$. 
 By Proposition~\ref{p:Avogadro-proaffine}, $|L|\le|A|+1=3$ and hence $L=\Aline ob=\{o,b,c\}$ for some point $c\in X$. Since $|L|=3$, the $3$-proregularity of $X$ guarantees that $X=\bigcup_{x\in A}\bigcup_{y\in L}\Aline xy=\{o\}\cup\Aline ab\cup\Aline ac$. Proposition~\ref{p:Avogadro-proaffine} ensures that $\max\{|\Aline ab|,|\Aline ac|\}\le|A|+1=3$. Since $X\ne L\cup A$, one of the lines $\Aline ab$ or $\Aline ac$ contains three points. We lose no generality assuming that $|\Aline ab|=3$ and hence $\Aline ab=\{a,b,d\}$ for some point $d$. 
 
If $\Aline ac=\{a,c\}$, then we can take any distinct points $p,q\notin X=\{o,a,b,c,d\}$ and consider the liner $Y\defeq X\cup\{p,q\}$ endowed with the family of lines
$$\mathcal L\defeq\big\{\{o,b,c\},\{a,b,d\},\{a,c,q\},\{o,a,p\},\{o,d,q\},\{c,d,p\},\{b,p,q\}\big\}.$$

\begin{picture}(80,90)(-150,-15)
\linethickness{0.6pt}
\put(0,0){\line(1,0){60}}
\put(0,0){\line(0,1){60}}
\put(0,60){\line(1,-1){60}}
\put(0,0){\line(1,1){30}}
\put(0,60){\line(1,-2){20}}
\put(60,0){\line(-2,1){40}}
\multiput(0,30)(4,-2){5}{\line(2,-1){3}}
\multiput(30,0)(-2,4){5}{\line(-1,2){1.5}}

\put(0,0){\circle*{3}}
\put(-2,-8){$o$}
\put(60,0){\circle*{3}}
\put(58,-8){$a$}
\put(0,60){\circle*{3}}
\put(-3,63){$d$}
\put(20,20){\circle*{3}}
\put(17,11){$c$}
\put(30,30){\circle*{3}}
\put(30,33){$b$}
\put(30,0){\circle*{3.5}}
\put(30,0){\color{white}\circle*{3}}
\put(28,-9){$p$}
\put(0,30){\circle*{3.5}}
\put(0,30){\color{white}\circle*{3}}
\put(-9,28){$q$}

\end{picture}

It is easy to see that $Y$ is a Steiner projective  plane containing $X$ as a subliner such that $|Y\setminus X|=|\{p,q\}|=2$. Consider the line $\Lambda\defeq\{o,d\}$ in $X$. Since $c\in\overline{\{A\cup\Lambda\}}\setminus\Aline A\Lambda$, the liner $X$ is not regular. 
Observe also that $\{a,c\}$ and $\{o,d\}$ are two disjoint lines in $X$, but every line containing the point $b$ intersects the line $\{a,c\}$, witnessing that the liner $X$ is not para-Playfair.
\smallskip

Next, assume that $|\Aline ac|=3$ and find a point $e\in\Aline ac\setminus\{a,c\}$. We claim that $e\in\Aline od$. In the opposite case, the distinct lines $\{o,b,c\}$ and $\{o,a\}$ contain the point $o$ and are disjoint with the line $\Aline de$, which contradicts the Proclus property of the liner $X$. This contradiction shows that $\{o,d,e\}$ is a line in $X$. 
We claim that $\Aline cd=\{c,d\}$ and $\Aline be=\{b,e\}$. Assuming that the line $\Aline cd$ contains a point $x\notin \{c,d\}$, we conclude that two concurrent lines $\{o,b,c\}$ and $\{o,d,e\}$ are disjoint with the line $\Aline ax$, which contradicts the Proclus property of $X$. This contradiction shows that $\Aline cd=\{c,d\}$. By analogy we can show that $\Aline be=\{b,e\}$.

\begin{picture}(80,92)(-150,-15)
\linethickness{0.6pt}
\put(0,0){\line(1,0){60}}
\put(0,0){\line(0,1){60}}
\put(0,60){\line(1,-1){60}}
\put(0,0){\line(1,1){30}}
\put(0,60){\line(1,-2){20}}
\put(60,0){\line(-2,1){60}}
\put(0,30){\line(1,0){30}}
\multiput(30,0)(-2,4){5}{\line(-1,2){1.5}}

\put(15,30){\circle*{3.5}}
\put(15,30){\color{white}\circle*{3}}

\put(0,0){\circle*{3}}
\put(-2,-8){$o$}
\put(60,0){\circle*{3}}
\put(58,-8){$a$}
\put(0,60){\circle*{3}}
\put(-3,63){$d$}
\put(20,20){\circle*{3}}
\put(17,11){$c$}
\put(30,30){\circle*{3}}
\put(30,33){$b$}
\put(30,0){\circle*{3.5}}
\put(30,0){\color{white}\circle*{3}}
\put(28,-9){$p$}
\put(0,30){\circle*{3}}
\put(-9,28){$e$}

\end{picture}

Choose any point $p\notin X=\{o,a,b,c,d,e\}$ and consider the liner $Y=X\cup\{p\}$ endowed with the family of lines 
$$\mathcal L\defeq\big\{\{o,b,c\},\{a,b,d\},\{a,c,e\},\{o,d,e\},\{o,a,p\},\{b,e,p\},\{c,d,p\}\big\}.$$It is easy to see that $Y$ is a Steiner projective  plane containing $X$ as a subliner such that $|Y\setminus X|=|\{p\}|=1$.
\end{proof}

\begin{corollary}\label{c:proj-comp-Proclus} Any projective completion $Y$ of a non-projective Proclus plane $X$ is normal.
\end{corollary}

\begin{proof} If the Proclus liner $X$ is $3$-long, then the normality of the projective completion $Y$ follows from Theorem~\ref{t:procompletion=>normal}. So, assume that $X$ is not $3$-long. By Theorem~\ref{t:Proclus-not-3long}, $X=P\setminus D$ for some Steiner projective plane $P$ and some set $D\subseteq P$ of cardinality $0<|D|\le 2$. It follows that $\|X\|=3$. By Proposition~\ref{p:pcomp3=>normal}, the projective completion $Y$ of the plane $X$ is normal.
\end{proof}

\section{The uniqueness of projective completions}

In this section we shall prove that a projective completion of a $3$-long liner is unique up to an isomorphism.

We recall that a function $F:X\to Y$ between the underlying sets of two liners $(X,\Af_X)$ and $(Y,\Af_Y)$ is called 
\begin{itemize}
\item a \index{liner morphism}\index{liner!morphism of}\defterm{liner morphism} if $\{Fxyz:xyz\in \Af_X\}\subseteq\Af_Y$;
\item a \index{line isomorphism}\index{liner!isomorphism}\defterm{liner isomorphism} if $F$ is a bijective map such that $F$ and $F^{-1}$ are liner morphisms;
\item a \index{liner embedding}\index{liner!embedding}\defterm{liner embedding} if $F$ is an liner isomorphism of $X$ into the subliner $F[X]$ of $Y$.
\end{itemize}

\begin{theorem}\label{t:projective=>extension} Let $Y$ be a normal completion of a $3$-long liner $X$. For any projective completion $Z$ of the liner $X$, there exists a unique injective liner morphism $F:Y\to Z$ such that $F(x)=x$ for all $x\in X$. If the liner $Y$ is $3$-ranked, then the map $F:Y\to Z$ is a liner embedding. 
\end{theorem}

\begin{proof} If $\|X\|\le 1$, then $Y=X=Z$ and the identity map $F:Y\to Z$ is a unique liner embedding such hat $F(x)=x$ for all $x\in X$. If $\|X\|=2$, then $\overline{Y\setminus X}\ne Y$ implies $|Y\setminus X|=1$. On the other hand, the definition of a projective completion ensures that $|Z\setminus X|=1$. Then there exists a unique injective map $F:Y\to Z$ such that $F(x)=x$ for all $x\in X$. Since $\|Y\|=2=\|Z\|$, the injective map $F:Y\to Z$ is a liner embedding. So, we assume that $\|X\|\ge 3$. To distinguish between the flat hulls of sets in the liners $X,Y,Z$, for sets $A\subseteq X$, $B\subseteq Y$, $C\subseteq Z$ we denote by ${\overline A}^X$, $\overline B^Y$ and $\overline C^Z$ the flat hulls of the sets $A,B,C$ in the liners $X,Y,Z$, respectively. If $A$ (resp. $B$, $C$) equals to set $\{a,b\}$ for some points $a,b$, then the flat $\overline{A}^X$ (resp $\overline{B}^Y$, $\overline{C}^Z$) will be denoted by $X_{ab}$ (resp. $Y_{ab}$, $Z_{ab}$). Since $X$ is a subliner of the liners $Y$ and $Z$, for any points $a,b\in X$, we have $Y_{ab}\cap X=X_{ab}=Z_{ab}\cap X$.

Given any point $y\in Y\setminus X$ and $o\in X$, consider the line $Y_{oy}\defeq\overline{\{o,y\}}^Y$ in the liner $Y$. Since the projective completion $Y$ of $X$ is $3$-long, $|Y_{oy}|\ge 3$. If $o\notin\overline{Y\setminus X}$, then $Y_{oy}\cap\overline{Y\cap X}=\{y\}$ and hence $X_{oy}\defeq Y_{oy}\cap X=Y_{oy}\setminus\{y\}$ is a line in the liner $X$. Then its flat hull $Z_{oy}\defeq \overline{X_{oy}}^Z$ is a line in the projective liner $Z$. 

\begin{claim}\label{cl:Zxy-Zzy} For any $y\in Y\setminus X$ and $o,x\in X\setminus\overline{Y\setminus X}$ the intersection $Z_{oy}\cap Z_{xy}$ is not empty.
\end{claim}

\begin{proof} Since the liner $Y$ is a normal completion of the liner $X$, the set $X_{oy}\cup X_{xy}$ has rank $\le 3$ in the liner $X$. Then the flat hull $P\defeq\overline{X_{oy}\cup X_{xy}}^Z$ has rank $\|P\|\le 3$ in the liner $Z$. By the projectivity of $Z$, the lines $Z_{oy}$ and $Z_{xy}$ in the flat $P$ have nonempty intersection.
\end{proof}

\begin{claim}\label{cl:Zxy-singleton} For every points $y\in Y\setminus X$ and $o\in X\setminus(\overline{Y\setminus X}\cup\overline{Z\setminus X})$, the set $Z_{oy}\cap \overline{Z\setminus X}=Z_{oy}\setminus X$ is a singleton.
\end{claim}

\begin{proof} It follows from $o\in Z_{oy}\setminus \overline{Z\setminus X}$ that the intersection $Z_{oy}\cap \overline{Z\setminus X}$ contains at most one point. Since the liner $X$ has rank $\|X\|\ge 3$, the line $X_{oy}\defeq Y_{oy}\cap X$ is a proper flat in $X$. Then $\overline{Y\setminus X}\cap X$ and $X_{oy}$ are two proper flats in the liner $X$. Since $X$ is $3$-long, there exists a point $x\in X\setminus (\overline{Y\setminus X}\cup X_{oy})$, see Proposition~\ref{p:cov-aff}. Then $Y_{xy}\cap Y_{oy}=\{y\}$ and hence the lines $X_{xy}\defeq Y_{xy}\cap X$ and $X_{oy}\defeq Y_{oy}\cap X$ are disjoint. By Claim~\ref{cl:Zxy-Zzy}, the intersection $Z_{xy}\cap Z_{oy}$ is not empty. Since $Z_{xy}\cap Z_{oy}\cap X=X_{xy}\cap X_{oy}=\varnothing$, the non-empty set $Z_{xy}\cap Z_{oy}$ is a subset of the set $Z_{oy}\setminus X$. Then $0<|Z_{oy}\setminus X|\le|Z_{oy}\cap\overline{Z\setminus X}|\le 1$ and hence the set $Z_{py}\setminus X=Z_{oy}\cap\overline{Z\setminus X}$ is a singeton.
\end{proof}

\begin{claim}\label{cl:Zy-singleton} For any point $y\in Y\setminus X$, the intersection 
$$Z_y\defeq \bigcap\{Z_{xy}:x\in X\setminus\overline{Y\setminus X}\}$$is a singleton that belongs to the set $Z\setminus X$.
\end{claim}

\begin{proof} Since $Y$ and $Z$ are completions of the liner $X$, the sets $\overline{Y\setminus X}\cap X$ and $\overline{Z\setminus X}\cap X$ are proper flats in the liner $X$. Since the liner $X$ is $3$-long, there exists a point $o\in X\setminus(\overline{Y\setminus X}\cup\overline{Z\setminus X})$, see Proposition~\ref{p:cov-aff}. Claim~\ref{cl:Zxy-singleton} ensures that the set $Z_{oy}\cap\overline{Z\setminus X}=Z_{oy}\setminus X$ contains a unique point $z_y$. Since the sets $X_{oy}\defeq Y_{oy}\cap X$ and $X\cap\overline{Y\setminus X}$ are proper flats in the $3$-long liner $X$, there exists a point $p\in X\setminus (Y_{oy}\cup\overline{Y\setminus X})$. 
By Claim~\ref{cl:Zxy-Zzy}, the intersection $Z_{oy}\cap Z_{py}$ is not empty. Since $$Z_{oy}\cap Z_{py}\cap X=X_{oy}\cap X_{py}=(X\cap Y_{oy})\cap(Y_{py}\cap X)=X\cap \{y\}=\varnothing,$$ the non-empty set $Z_{oy}\cap Z_{py}$ is a subset of the set $$Z_{oy}\setminus X=\{z_y\}.$$ Then $Z_{oy}\cap Z_{py}=\{z_y\}$ and hence $Z_y\subseteq Z_{oy}\cap Z_{py}\subseteq\{z_y\}$. It remains to show that $z_y\in Z_{xy}$ for every point $x\in X\setminus\overline{Y\setminus X}$. If $x\in Y_{oy}$, then $Y_{xy}=Y_{oy}$, $X_{xy}=Y_{xy}\cap X=Y_{oy}\cap X=X_{oy}$ and $z_y\in Z_{oy}=Z_{xy}$. If $x\notin Y_{oy}$, then $Y_{xy}\cap Y_{oy}=\{y\}$ and $X_{xy}\cap X_{oy}=\varnothing$ and hence $Z_{xy}\cap Z_{oy}\cap X=X_{xy}\cap X_{oy}=(X\cap Y_{xy})\cap(Y_{oy}\cap X)=X\cap\{y\}=\varnothing$ and $Z_{xy}\cap Z_{oy}\subseteq Z_{oy}\setminus X=\{z_y\}$. By Claim~\ref{cl:Zxy-Zzy}, the set $Z_{xy}\cap Z_{oy}$ is not empty, which implies $Z_{xy}\cap Z_{oy}=\{z_y\}$ and $z_y\in Z_{xy}$.
\end{proof}

By Claim~\ref{cl:Zy-singleton}, for every point $y\in Y\setminus X$, the set
$$Z_y\defeq\bigcap\{Z_{xy}:x\in X\setminus\overline{Y\setminus X}\}$$contains a unique point $z_y$ and this point belongs to the set $Z\setminus X$. Consider the function $F:Y\to Z$ defined by the formula
$$F(x)=\begin{cases} x&\mbox{if $x\in X$};\\
z_y&\mbox{if $y\in Y\setminus X$}.
\end{cases}
$$

\begin{claim} The function $F$ is injective.
\end{claim}

\begin{proof} Given any distinct points $y,u\in Y$, we should prove that $F(y)\ne F(u)$. If $y,u\in X$, then $F(y)=y\ne u=F(u)$ and we are done. If $y\in X$ and $y\notin Y\setminus X$, then $F(y)=y\in X$ and $F(u)=z_u\in Z\setminus X$ and hence $F(y)\ne F(u)$. By analogy we can show that $F(y)\ne F(u)$ if $y\in Y\setminus X$ and $u\in X$.

So, assume that $u,y\in Y\setminus X$. Since $X\cap\overline{Y\setminus X}$ and $Y_{uy}\cap X$ are proper flats in the $3$-long liner $X$, there exists a point $o\in X\setminus(Y_{uy}\cup\overline{Y\setminus X})$, by Proposition~\ref{p:cov-aff}. Then $Y_{oy}\cap \overline{Y\setminus X}=\{y\}=Y_{oy}\cap Y_{uy}$. Since the liner $X$ is $3$-long, the line $X_{oy}\defeq Y_{oy}\setminus\{y\}$ in $X$ contains a point $x\in X_{oy}\setminus\{z_y\}\subseteq X_{oy}=Y_{oy}\setminus\{y\}\subseteq Y_{oy}\setminus Y_{uy}$.   It follows from $x\notin Y_{uy}$ and $u\ne y$ that $Y_{xy}\cap Y_{xu}=\{x\}$ and then $X_{xy}\cap X_{xu}=\{x\}=Z_{xy}\cap Z_{xu}$ and finally, 
$\{z_y\}\cap\{z_u\}\subseteq  Z_{xy}\cap Z_{xu}=\{x\}$, which implies $F(y)=z_y\ne z_u=F(u)$ because $x\ne z_y$.
\end{proof}

\begin{claim} The injective function $F:Y\to Z$ is a liner morphism.
\end{claim}

\begin{proof} Let $\Af_X,\Af_Y,\Af_Z$ be the line relations of the liners $X,Y,Z$, respectively. Given any collinear triple $abc\in \Af_Y$, we should prove that the triple $a'b'c'\defeq Fabc$ belongs to the line relation $\Af_Z$. If $b\in\{a,c\}$, then $b'\in\{a',c'\}$ and hence $a'b'c'\in \Af_Z$, by the Reflexivity Axiom {\sf(RL)}, see Definition~\ref{d:liner}. So, we assume that $a\ne b\ne c$. In this case $a\ne c$ and the injectivity of the function $F$ ensures that $a',b',c'$ are pairwise distinct points of the liner $Z$.  

If $\{a,b,c\}\subseteq X$, then $a'b'c'=abc\in \Af_X\subseteq\Af_Z$ and we are done. So, assume that $\{a,b,c\}\not\subseteq X$. We lose no generality assuming that $c\in Y\setminus X$. Then $c'=F(c)=z_c\in Z\setminus X$.  If $a\notin \overline{Y\setminus X}$, then $b\in\Aline ac\setminus\{c\}\subseteq Y\setminus\overline{Y\setminus X}\subseteq X$ and hence $a'=F(a)=a$, $b'=F(b)=b$ and $c'\in Z_{ac}=\overline{X_{ac}}^Z=\overline{X_{ab}}^Z=Z_{ab}=Z_{a'b'}$ and finally, $a'b'c'\in \Af_Z$. So, we assume that $a\in\overline{Y\setminus X}$ and hence $b\in \Aline ac\subseteq\overline{Y\setminus X}$. 

Since $\overline{Y\setminus X}\cap X$ and $\overline{Z\setminus X}\cap X$ are proper flats in the $3$-long liner $X$, there exists a point $o\in X\setminus(\overline{Y\setminus X}\cup\overline{Z\setminus X})$, see Proposition~\ref{p:cov-aff}. 

Observe that $P\defeq\overline{\{o,a,b,c,\}}^Y=\overline{\{o,a,c\}}^Y$ is a plane in the liner $Y$. Since $Y$ is a normal completion of $X$, the intersection $P\cap X$ is a plane in the liner $X$.  Then its flat hull $\Pi\defeq \overline{P\cap X}^Z$ in $Z$ is a plane in $Z$. Since the plane $P\cap X$ contain the lines $X_{oa},X_{ob}, X_{oc}$, its flat hull $\Pi$ contains the lines $Z_{oa},Z_{ob},Z_{oc}$. 
The definition of the function $F$ ensures that $a'=F(a)\in Z_{oa}=\overline{X_{oa}}^Z\subseteq \Pi$, $b'=F(b)\in Z_{ob}\subseteq\Pi$ and $c'=F(c)\in Z_{oc}\subseteq \Pi$. Since $o\in \Pi\setminus\overline{Z\setminus X}$, the $3$-rankedness of the projective liner $Z$ ensures that $\|\Pi\cap\overline{Z\setminus X}\|<\|\Pi\|=3$. If $\{a',b',c'\}\subseteq \overline{Z\setminus X}$, then 
$\|\{a',b',c'\}\|\le\|\Pi\cap\overline{Z\setminus X}\|\le 2$ and hence $a'b'c'\in\Af_Z$.

So, we assume that $\{a',b',c'\}\not\subseteq \overline{Z\cap X}$. Since $c'\in Z\setminus X$, we lose no generality assuming that $a'\in Z\setminus \overline{Z\cap X}\subseteq X$. Since $F[X]=X$, the injectivity of the function $F$ ensures that $a\in X$ and $a'=F(a)$.  By the projectivity of the liner $Z$, the lines $Z_{a'b'}$ and $Z_{oc'}$ in the plane $\Pi$ have a common point $c''$.  If $c''=c'$, then $b'\in Z_{ac''}=Z_{a'c'}$ and hence $a'b'c'\in \Af_Z$.
It remains to show that the case $c''\ne c'$ is impossible. To derive a contradiction, assume that $c''\ne c'$. Then $c''\in  Z_{oc}\setminus\{c'\}=X_{oc}\subseteq X$. If $b\in X$, then $$c''\in Z_{a'b'}\cap Z_{oc'}\cap X=Z_{ab}\cap Z_{oc}\cap X=X_{ab}\cap X_{oc}=X\cap Y_{ab}\cap Y_{oc}=X\cap \{c\}=\varnothing,$$which is a desired contradiction. Next, assume that $b\notin X$ and hence $b'=F(b)=z_b$.  Assuming that $c''\notin\overline{Y\setminus X}$, we conclude that $a=a'\in Z_{c''b'}\setminus\{b'\}=Z_{c''b}\setminus\{b'\}=X_{c''b}=Y_{c''b}\setminus\{b\}$ and hence $a\in \overline{Y\setminus X}\cap (Y_{c''b}\setminus\{b\})=\varnothing$, which is a desired contradiction showing that $c''\in\overline{Y\setminus X}$. Then $c''\in \overline{Y\setminus X}\cap(Z_{oc}\setminus\{c'\})=\overline{Y\setminus X}\cap X_{oc}=\emptyset$, which is a desired contradiction. 
\end{proof}

The uniqueness of the liner morphism $F$ is proved in the following claim.

\begin{claim} Let $G:X\to Y$ be any injective liner morphism such that $G(x)=x$ for all $x\in X$. Then $G=F$.
\end{claim}

\begin{proof}  Given any point $y\in Y\setminus X$, we should show that $F(y)=G(y)$. Taking into account that the  injective functions $F,G$ are identity on $X$, we conclude that $F(y),G(y)\in Z\setminus X$. Since $\overline{Y\setminus X}\cap X$ and $\overline{Z\setminus X}\cap X$ are two proper flats in the $3$-long liner $X$, there exists a point $x\in X\setminus(\overline{Y\setminus X}\cup\overline{Z\setminus X})$, see Proposition~\ref{p:cov-aff}. Since the liner $Y$ is $3$-long and $|Y_{xy}\cap\overline{Y\setminus X}|\le 1$, the intersection $X_{xy}\defeq Y_{xy}\cap X$ is a line in $X$ and its flat hull $Z_{xy}\defeq\overline{X_{xy}}^Z$ in $Z$ is a line in the liner $Z$. Choose any point $x'\in X_{xy}$. Since the functions $F,G$ are liner morphisms, the inclusion $xyx'\in\Af_Y$ implies $\{F(y),G(y)\}\subseteq Z_{xx'}\setminus X=Z_{xy}\setminus X=\{z_y\}$ and hence $F(y)=z_y=G(y)$.
\end{proof}

\begin{claim} If the liner $Y$ is $3$-ranked, then the function $F:Y\to Z$ is a liner embedding.
\end{claim}

\begin{proof} Since $F$ is a liner morphism, it suffices to show that for every triple $abc\in Y^3$ and its image $a'b'c'\defeq Fabc$ under the map $F$, if $a'b'c'\in \Af_Z$, then $abc\in \Af_Y$. So, take any triple $abc\in Y^3$ and assume that the triple $a'b'c'\defeq Fabc$ belongs to the line relation $\Af_Z$.
We have to show that $abc\in\Af_Y$. This follows from the Reflexivity Axiom {\sf (RL)} if $b\in\{a,c\}$. So, assume that $b\notin\{a,c\}$. In this case $a\ne c$. Indeed, assuming that $a=c$, we conclude that $a'=c'$ and then $a'b'c'\in\Af_Z$ imply $a'=b'=c'$ and hence $a=b=c$, which contradicts our assumption.  Therefore, $a\ne c$ and hence $a,b,c$ are three pairwise distinct points in the liner $Y$.

Since $\overline{Y\setminus X}\cap X$ and $\overline{\{a',b',c'\}}\cap X$ are proper flats in the $3$-long liner $X$, there exists a point $o\in X\setminus(\overline{Y\setminus X}\cup\overline{\{a',b',c'\}})$.  Since the points $a',b',c'$ are collinear in the liner $Z$, the flat hull $\Pi\defeq\overline{\{o,a',b',c'\}}^Z$ is a plane in $Z$. By Theorem~\ref{t:procompletion=>normal}, the projective completion $Z$ of $X$ is normal and hence the intersection $X\cap \Pi$ is a plane in $X$. Then the set $P\defeq \overline{X\cap\Pi}^Y$ is a plane in the liner $Y$. If $a\in X$, then $a=a'\in\Pi\cap X\subseteq P$. If $a\in Y\setminus X$, then $X_{oa}\subseteq Z_{oa}\cap X=Z_{oa'}\cap X\subseteq \Pi\cap X$ and again $a\in \overline{X_{oa}}^Y\subseteq \overline{X\cap\Pi}^Y=P$. By analogy we can show that $\{b,c\}\subseteq P$.  Since the liner $Y$ is $3$-ranked and $o\in P\setminus\overline{Y\setminus X}$, the flat $P\cap\overline{Y\setminus X}$ has rank $\|P\cap\overline{Y\setminus X}\|<\|P\|=3$. If $\{a,b,c\}\subseteq \overline{Y\setminus X}$, then $\|\{a,b,c\}\|\le\|P\cap\overline{Y\setminus X}\|\le 2$ and hence $abc\in\Af_Y$. So, assume that $\{a,b,c\}\not\subseteq\overline{Y\setminus X}$. If $\{a,b,c\}\subseteq X$, then $abc=a'b'c'\in X^3\cap\Af_Z=\Af_X\subseteq \Af_Y$ and we are done.
\smallskip

So, assume that $\{a,b,c\}\not\subseteq X$ and hence $\{a,b,c\}\cap\overline{Y\cap X}=\{a,b,c\}\setminus X$ is a singleton. In this case we lose no generality assuming that $a,b\in Y\setminus \overline{Y\cap X}\subseteq X$ and $c\in Y\setminus X$. If $b\in X$, then $c'\in Z_{ac'}=Z_{a'c'}$, $b=b'\in Z_{a'c'}\cap X=Z_{ac'}\cap X=X_{ac}=Y_{ac}\cap X\subseteq Y_{ac}$ and hence $abc\in \Af_Y$. If $b\notin X$, then $X_{ab}\cup X_{ac}\subseteq X\cap\overline{\{a',b',c'\}}^Z$ and hence $\|X_{ab}\cup X_{ac}\|\le 2$ and $\|Y_{ab}\cup Y_{ac}\|\le 2$, which implies $abc\in\Af_Y$.
\end{proof}
\end{proof}

\begin{theorem}\label{t:extend-to-procompletions} For any projective completions $Y,Z$ of a $3$-long liner $X$, there exists a unique liner isomorphism $F:Y\to Z$ such that $F(x)=x$ for all $x\in X$.
\end{theorem}

\begin{proof} By Theorem~\ref{t:procompletion=>normal}, the projective completions $Y,Z$ of $X$ are normal, and by Theorem~\ref{t:projective=>extension}, there exist unique injective liner homorphisms $F:Y\to Z$ and $G:Z\to Y$ such that $F(x)=x=G(x)$ for all $x\in X$. Then $GF:Y\to Y$ and $FG:Z\to Z$ are injective liner homomorphisms such that $GF(x)=x=FG(x)$ for every $x\in X$. The uniqueness part of Theorem~\ref{t:projective=>extension} ensures that $GF$ and $FG$ are the identity maps of the liners $Y$ and $Z$, respctively. Therefore, $F:Y\to Z$ is a liner isomorphism with $F^{-1}=G$.
\end{proof}

Theorems~\ref{t:spread=projective1}, \ref{t:extend-to-procompletions} and Proposition~\ref{p:spread-3long} imply the following corollary.

\begin{corollary}\label{c:pcompletion=scompletion} Let $X$ be a $3$-long completely regular liner and $\overline X$ be the spread completion of $X$. For every projective completion $Y$ of $X$, there exists a unique liner isomorphism $F:Y\to\overline X$ such that $F(x)=x$ for every $x\in X$.
\end{corollary}

\begin{theorem}\label{t:extend-isomorphism-to-completions} Let $X,Y$ be two $3$-long liners and $\widetilde X$, $\widetilde Y$ be projective completions of the liners $X,Y$, respectively. Every isomorphism $F:X\to Y$ uniquely extends to an isomorphism $\widetilde F:\widetilde  X\to\widetilde Y$.
\end{theorem}

\begin{proof} Choose any set $Z$ containing the set $Y$ so that $|Z\setminus Y|=|\widetilde X\setminus X|$. Choose any bijection $\varphi:\widetilde X\setminus X\to Z\setminus Y$.  The maps $F$ and $\varphi$ compose a bijective function $\Phi:\widetilde X\to Z$ defined by 
$$\Phi(x)\defeq
\begin{cases}
F(x)&\mbox{if $x\in X$};\\
\varphi(x)&\mbox{if $x\in \widetilde X\setminus X$}.
\end{cases}
$$
Endow the set $Z$ with the line relation $$\Af_Z\defeq \{\Phi xyz:xyz\in \Af_{\widetilde X}\},$$where $\Af_{\widetilde X}$ is the line relation on the liner $\widetilde X$. Then $\Phi:X\to Z$ is an isomorphism of the liners $\widetilde X$ and $Z$. Since $F$ is an isomorphism of the liners $X$ and $Y$, the line relation $\Af_Y$ of the liner $Y$ is equal to $\{Fxyz:xyz\in\Af_X\}=\{\Phi xyz:xyz\in\Af_{\widetilde X}\cap X^3\}=\Af_Z\cap Y^3$, which means that $Y$ is a subliner of the liner $Z$. Taking into account that $\widetilde X$ is a projective completion of $X$ and $\Phi:\widetilde X\to Z$ is a liner isomorphism with $\Phi[X]=Y$, we conclude that $Z$ is a projective completion of the liner $Y$. By Theorem~\ref{t:extend-to-procompletions}, there exists a liner isomorphism $\Psi:Z\to \widetilde Y$ such that $\Psi(y)=y$ for all $y\in Y$. Then $\widetilde F\defeq\Psi\circ\Phi:\widetilde X\to\widetilde Y$ is a liner isomorhism extending the isomorphism $F$. 

To see that $\widetilde F$ is unique, assume that $\overline F:\widetilde X\to\widetilde Y$ is another liner isomorphism extending the isomorphism $F$. Then $I\defeq(\overline F)^{-1}\widetilde F:\widetilde X \to\widetilde X$ is a liner isomorphism such that $I(x)=x$ for every $x\in X$. The uniqueness part of  Theorem~\ref{t:extend-to-procompletions} ensures that $I=(\overline F)^{-1}\widetilde F$ is the identity map of the projective liner $\widetilde X$. Then $\overline F=\overline FI=\overline F(\overline F)^{-1}\widetilde F=\widetilde F$, witnessing that the isomorphism $\widetilde F:\widetilde X\to\widetilde Y$ extending the isomorphism $F$ is unique. 
\end{proof}


\section{Interplay between a liner and its horizon}

For a $3$-long liner $X$ possessing a projective completion $Y$, the complement $Y\setminus X$ is called the \index{horizon}\index{completion!horizon of}\defterm{horizon} of $X$. Since a projective completion of $X$ is unique (up to an isomorphism), the horizon of $X$ is also uniquely determined. In this section we discuss the interplay betwen properties of a liner and its horizon.

We start with a duality between the proaffinity a liner and the proflat property of its horizon.

\begin{definition} A subset $A$ of a liner $X$ is called \index{proflat set}\index{subset!proflat}\defterm{proflat} in $X$ if $|\Aline xy\setminus A|\le 1$ for every points $x,y\in A$. It is clear that every flat in a liner $X$ is a proflat in $X$.
\end{definition}

\begin{exercise} Show that for every flats $A\subseteq B$ in a liner $X$, the set $B\setminus A$ is  proflat in $X$.
\end{exercise}

\begin{exercise} Let $(A_n)_{n\in\IZ}$ be a sequence of flats in a liner $X$ such that $A_n\subseteq A_{n+1}$ for all $n\in\IZ$. Show that the set $\bigcup_{n\in\w}A_{2n+1}\setminus A_{2n}$ is  proflat in $X$.
\end{exercise}

\begin{exercise} Let $H $ be a proflat set in a liner $X$. Show that $|\Aline xy\cap H|\le 1$ for any points $x,y\in X\setminus H$.
\end{exercise}

Let us recall that a liner $X$ is called {\em Proclus} if for every plane $P\subseteq X$, line $L\subseteq P$ and point $x\in P\setminus L$ there exists at most one line $\Lambda$ in $X$ such that $x\in\Lambda\subseteq P\setminus L$.  

\begin{proposition}\label{p:projective-minus-proflat} For every proflat set $H$ in a  projective liner $Y$, the subliner $X\defeq Y\setminus H$ of $Y$ is Proclus.
\end{proposition}

\begin{proof} To distinguish between flat hulls of sets in the liners $X$ and $Y$, for subsets $A\subseteq X$ and $B\subseteq Y$, we denote by $\overline{A}'$ and $\overline{B}$ the flat hulls of the sets in the liners $X$ and $Y$, respectively. If $A$ (resp. $B$) is equal to $\{a,b\}$ for some points $a,b$, then the flat hull $\overline{\{a,b\}}'$ (resp. $\overline{\{a,b\}}$) will be denoted by $\Aline ab'$ (resp. $\Aline ab$). If $a,b\in X$, then $\Aline ab'=\Aline ab\cap X$. 

By Theorem~\ref{t:Proclus<=>}, a liner is Proclus if and only if it is proaffine and $3$-proregular. 
To prove that the liner $X$ is  proaffine, it suffices to check that for any points $o,x,y\in X$ and $p\in \Aline xy'\setminus(\Aline ox'\cup \Aline oy')$, the set $I\defeq\{v\in \Aline oy':\Aline vp'\cap \Aline ox'=\varnothing\}=\{v\in \Aline oy':\Aline vp\cap \Aline ox\subseteq H\}$ contains at most one point. Assuming that $|I|>1$, choose two distinct points $u,v\in I$. Since $p\notin \Aline xy\setminus \Aline oy$, the flats $\Aline up$, $\Aline vp$, $\Aline ox$ are lines in the liner $Y$.

The projectivity of the liner $Y$ ensures that the intersections $\Aline up\cap \Aline ox$ and $\Aline vp\cap \Aline ox$ are not empty and hence they contain some points $u'\in \Aline up\cap \Aline ox$ and $v'\in \Aline vp\cap \Aline ox$. It follows from $u,v\in I$ and $u\ne v$ that $u',v'\in H$ and $u'\ne v'$. Since the set  $H$ is  proflat, $|\Aline{u'}{v'}\setminus H|\le 1$. Since $\{o,x\}\subseteq \Aline ox'\subseteq \Aline{u'}{v'}\setminus H$, the points $o$ and $x$ coincide  and hence $p\in \Aline yx'=\Aline yo'$, which contradicts the choice of the point $p$. This contradiction shows that $|I|\le1$ and hence the liner $X$ is proaffine.
\smallskip

To show that the proaffine liner $X$ is $3$-proregular, fix any flat $A\subseteq X$ of rank $\|A\|<3$ and any points $o\in A$ and $b\in X\setminus A$ such that $|X\cap \Aline ob|\ge 3$. Given any point $z\in {\overline{A\cup\{b\}}}'$, we need to find points $x\in A$ and $y\in X\cap \Aline ob$ such that $z\in \Aline xy$. If $z\in A$, then the points $x\defeq z$ and $y\defeq b$ have the required properties. So, we assume that $z\notin A$. 

If $z\in \Aline ob$, then the points $x\defeq o$ and $y\defeq b$ have the required properties. So, we assume that $z\notin \Aline ob$. In this case $A\ne\{o\}$ and hence $A=\Aline oa$ for some point $a\in A$ (let us recall that the set $A$ has rank $<3$ in the liner $X$). Since the set $X\cap\overline{A\cup \{b\}}$ is flat in $X$, $x\in\overline{A\cup\{b\}}'\subseteq X\cap\overline{A\cup\{b\}}=X\cap\overline{\{o,a,b\}}$. Since $|\Aline oa|\ge|\{o,a\}|= 2$ and the set $H$ is proflat, the intersection $\Aline oa\cap H$ contains at most one point.  Since the liner $X$ is projective, for every $y\in \Aline ob$, the lines $\Aline yz$ and $\Aline oa$ have a common point. Since $|\Aline oa\cap H|\le 1$, the set $I=\{y\in \Aline ob:\Aline yx\cap \Aline oa\subseteq H\}$ contains at most one point.  Since the liner $X$ is $3$-long,  there exists a point $y\in (\Aline ob\cap X)\setminus(\{o\}\cup  I)$. Since $y\notin I$, the unique point $x$ of the intersection $\Aline yz\cap \Aline oa$ does not belong to the set $H$ and hence $x,y\in X$ are required points with $z\in \Aline xy$ witnessing that the liner $X$ is $3$-proregular.
\end{proof}

\begin{theorem}\label{t:proaffine<=>proflat} For a projective completion $Y$ of a liner $X$, the following conditions are equivalent:
\begin{enumerate}
\item the liner $X$ is Proclus;
\item the liner $X$ is proaffine;
\item the horizon $Y\setminus X$ of $X$ is proflat in $Y$.
\end{enumerate}
\end{theorem}

\begin{proof} The implications $(3)\Ra(1)$ and $(1)\Ra(2)$ follow from Proposition~\ref{p:projective-minus-proflat} and Theorem~\ref{t:Proclus<=>}, respectively. So, it remains to prove that $(2)\Ra(3)$, which is equivalent to $\neg(3)\Ra\neg(2)$.  

If the horizon $H\defeq Y\setminus X$ is not proflat in $X$, then  $H$ contains two points $a,b$ such that the set $\Aline ab\setminus H$ contains two distinct points $o$ and $x$. Observe that $o,x\in\Aline ab\subseteq\overline H$. Since $Y$ is a completion of $X$, $\overline{H}\ne Y$ and hence there exists a point $p\in Y\setminus\overline H$. Since the liner $Y$ is $3$-long, there exists a point $y\in \Aline xp\setminus\{x,p\}$. Assuming that $y\in \overline H$, we conclude that $p\in\Aline xy\subseteq \overline{H}$, which contradicts the choice of the point $p$. This contradiction shows that $y\notin\overline H$ and hence $o,x,y\in Y\setminus H=X$. It follows from $y\notin\Aline ab=\Aline ox$ that $\Aline yo\cap\Aline ab=\{o\}$ and hence $a,b\notin\Aline oy$.

\begin{picture}(200,85)(-150,-10)

\put(0,0){\line(1,0){90}}
\put(30,0){\line(0,1){60}}
\put(60,0){\line(-1,2){30}}
\put(90,0){\line(-7,4){60}}
\put(0,0){\line(2,1){48}}

\put(0,0){\color{red}\circle*{3}}
\put(-9,-3){$a$}
\put(30,0){\circle*{3}}
\put(28,-8){$o$}
\put(60,0){\circle*{3}}
\put(58,-8){$x$}
\put(90,0){\color{red}\circle*{3}}
\put(93,-3){$b$}
\put(30,15){\circle*{3}}
\put(22,14){$v$}
\put(30,34.3){\circle*{3}}
\put(22,31){$u$}
\put(30,60){\circle*{3}}
\put(22,58){$y$}
\put(48,24){\circle*{3}}
\put(51,26){$p$}
\end{picture}

Since $Y$ is projective, there exist points $u\in\Aline oy\cap \Aline ap$ and $v\in\Aline oy\cap\Aline bp$. It follows from $a,b\in H$ and $p\notin\overline H$ that $u,v\notin\overline H$. Then $u,v$ are two distinct  points in the line $\Aline oy\cap X$ such that $\Aline up\cap\Aline ox=\{a\}$ and $\Aline vp\cap\Aline ox=\{b\}$, which implies $(\Aline up\cap X)\cap(\Aline xo\cap X)=\varnothing=(\Aline vp\cap X)\cap(\Aline ox\cap X)$ and witnesses that the liner $X$ is not proaffine.
\end{proof}

\begin{theorem}\label{t:flat-horizon} Let $Y$ be a projective completion of a liner $X$. The liner $X$ is para-Playfair if and only if the set $H\defeq Y\setminus X$ is flat in $Y$.
\end{theorem}

\begin{proof} Assume that the liner $X$ is para-Playfair. By Proposition~\ref{p:Playfair=>para-Playfair=>Proclus}, the para-Playfair liner $X$ is Proclus. By Theorem~\ref{t:proaffine<=>proflat}, the set $H=Y\setminus X$ is proflat in $X$. Assuming that $H$ is not flat, we can find distinct points $x,y\in H$ and $z\in \Aline xy\setminus H$. Since $\overline H\ne Y$, there exists a point $o\in Y\setminus\overline H$. Since the liner $Y$ is $3$-long, there exists a points $a\in \Aline ox\setminus\{o,x\}$ and $b\in\Aline oy\setminus\{o,y\}$. Assuming that $o\in\Aline ab$, we conclude that $o\in\Aline ab=\Aline xy\subseteq\overline H$, which contradicts the chocie of the point $o$. This contradiction shows that $o\notin\Aline ab$. Since the liner $Y$ is projective, the lines $\Aline ab$ and $\Aline oz$ in the plane $\overline{\{o,x,y\}}$ have a common point $c$, which is distinct from the point $o\notin\Aline ab$. By the projectivity of $Y$, the lines $\Aline xb$ and $\Aline oz$ have a common point $d$. It follows from $x,y,z\in\overline H$ and $o\notin\overline H$ that $a,b,c,d\in Y\setminus\overline H\subseteq X$. 

\begin{picture}(200,95)(-180,-15)
\linethickness{0.6pt}
\put(0,0){\line(1,1){60}}
\put(0,0){\line(-1,1){60}}
\put(0,0){\line(0,1){60}}
\put(-60,60){\line(3,-1){90}}
\put(-60,60){\color{blue}\line(1,0){120}}
\put(60,60){\line(-3,-1){90}}
\put(-30,30){\line(1,0){60}}

\put(0,0){\circle*{3}}
\put(-3,-8){$o$}
\put(30,30){\circle*{3}}
\put(28,33){$b$}
\put(60,60){\circle*{3}}
\put(63,58){$y$}
\put(-30,30){\circle*{3}}
\put(-32,33){$a$}
\put(0,30){\circle*{3}}
\put(1,23){$c$}
\put(0,40){\circle*{3}}
\put(1,43){$d$}
\put(0,60){\color{red}\circle*{3}}
\put(-2,63){$z$}
\put(-60,60){\circle*{3}}
\put(-68,58){$x$}
\end{picture}

Consider the plane $\Pi\defeq\overline{\{o,a,b\}}$ in the liner $X$ and observe that $c\in \Aline ab\cap \Aline oz\subseteq \Aline ab\cap X\subseteq\Pi$ and $z,d\in \Aline oc\subseteq \Pi$. Taking into account that the liner $X$ is para-Playfair, and $\Aline oa\cap X=\Aline ox\cap X$, $\Aline bd\cap X=\Aline bx\cap X$ are two disjoint lines in the plane $\Pi\subseteq X$, we can find a line $L\subseteq\Pi$ such that $z\in L$ and $L\cap(\Aline ox\cap X)=\varnothing$. By the projectivity of $Y$, the flat hull $\overline L$ of the line $L$ in the projective liner $Y$ has a common point $p$ with the line $\Aline ox$. It follows from $\{p\}\cap X\subseteq \Aline ox\cap\overline L\cap X=(\Aline ox\cap X)\cap L=\varnothing$ that $p\in \Aline ox\setminus X=\{x\}$ and hence $\overline L=\Aline xz=\Aline xy$. Since the set $H$ is proflat in $Y$, $L=\overline L\cap X=\Aline xy\setminus H=\{z\}$, which is impossible because lines in liners contain at least two points. This contradiction shows that the horizon $H$ of $X$ is flat.
\smallskip

Now assuming that the horizon $H$ of $X$ in $Y$ is flat, we shall prove that the liner $X$ is para-Playfair. By Theorems~\ref{t:proaffine<=>proflat}, the liner $X$ is Proclus. 
To show that $X$ is para-Playfair, take any plane $P\subseteq X$, disjoint lines $L,\Lambda\subseteq P$ and  point $x\in P\setminus L$.  Let $\overline L,\overline \Lambda,\overline P$ be the flat hulls of the sets $L,\Lambda,P$ in the projective liner $Y$. Then $\overline P$ is a plane in $Y$ and $\overline L,\overline \Lambda$ are lines in the plane $\overline P$. By the projectivity of the liner $Y$, there exists a point $y\in \overline L\cap\overline\Lambda$. Assuming that $y\in X$, we conclude that $$y\in X\cap\overline L\cap\overline\Lambda=(X\cap\overline L)\cap(X\cap\overline\Lambda)=L\cap\Lambda=\varnothing,$$which is a contradiction showing that $y\in Y\setminus X=H$. Since $x\in P\setminus L\subseteq X=Y\setminus H$ and $y\in H$, the intersection $\Aline xy\cap H$ coincides with the singleton $\{y\}$.
Since the projective liner $Y$ is $3$-long, the intersection $L_x\defeq \Aline xy\cap X=\Aline xy\setminus\{y\}$ is a line in the liner $X$.  By Corollary~\ref{c:proj-comp-Proclus}, the projective completion $Y$ of the Proclus liner $X$ is normal and hence $\overline P\cap X=P$ and $L_x\subseteq X\cap\overline P=P$. It follows from $x\notin L$ that $\Aline x y\cap \overline L=\{y\}$ and hence $L_x\cap L=\varnothing$.  The uniqueness of the line $L_x$ follows from the Proclus property of the liner $X$. Therefore, the liner $X$ is para-Playfair. 
\end{proof}

\begin{proposition}\label{p:projective-minus-flat} Let $Y$ be a projective completion of a $3$-long liner $X$. If the horizon $Y\setminus X$ of $X$ is flat in $Y$, then
\begin{enumerate}
\item for every set $A\subseteq X$ and its flat hull $\overline{A}$ in the liner $Y$, the intersection $\overline{A}\cap X$ coincides with the flat hull of $A$ in the liner $X$;
\item the liner $X$ is regular;
\item the liners $X$ and $Y$ have the same rank.
\end{enumerate}
\end{proposition}

\begin{proof} Assume that the horizon $Y\setminus X$ of $X$ in $Y$ is flat in $Y$. Then the liner $X$ is proaffine, according to Theorem~\ref{t:proaffine<=>proflat}.  To distinguish between flat hulls of sets in the liners $X$ and $Y$, for subsets $A\subseteq X$ and $B\subseteq Y$, we denote by $\overline{A}'$ and $\overline{B}$ the flat hulls of the sets in the liners $X$ and $Y$, respectively.
\smallskip

1. Given a set $A\subseteq X$, we should prove that $\overline{A}'=X\cap\overline A$. By Lemma~\ref{l:trace-flat}, the flat $B\defeq\overline A'$ is equal to $X\cap\overline{B}$, where $\overline B$ is the flat hull of the set $B$ in the liner $Y$. Then $$\overline{A}'\subseteq X\cap \overline{A}\subseteq X\cap\overline{B}=B=\overline A'$$ and hence $\overline{A}'=X\cap\overline{A}$.
\smallskip

2. To show that the liner $X$ is regular, fix any flat $A$ in $X$ and points $o\in A$, $b\in X\setminus A$, and $z\in \overline{A\cup\{b\}}'$. We have to find  points $x\in A$ and $y\in X\cap \Aline ob$ such that $z\in \Aline xy$. If $z\in\Aline ob$, then the points $x\defeq o$ and $y\defeq b$ have the required property. If $z\in A$, then the points $x\defeq z$ and $y\defeq b$ have the required property. So, we assume that $z\notin A\cup\Aline ob$. By the preceding statement $z\in\overline{A\cup\{b\}}'=\overline{A\cup\{b\}}\cap X$ and $A=\overline{A}\cap X$. By Theorem~\ref{t:projective<=>}, the projective liner $Y$ is strongly regular. Then for the point $z\in\overline{A\cup\{b\}}$,  there exists a point $a\in\overline{ A}$ such that $z\in\Aline ab$.  Since $z\in\Aline ab\setminus \Aline ob$, the point $a$ does not belong to 
the line $\Aline ob$ and hence $\overline{\{o,a,b\}}$ is a plane in the projective liner $Y$. Assuming that $z\in\Aline oa$, we conclude that $z\in \Aline oa\cap X\subseteq \overline{A}\cap X=A$, which contradicts our assumption. Therefore, the point $z$ does not belong to the lines $\Aline oa$ and $\Aline ob$ in the plane $\overline{\{o,a,b\}}$.

Since the horizon $Y\setminus X$ of $X$ is flat in $Y$ and $o\in X$, the set $\Aline ob\cap (Y\setminus X)$ contains at most one point and hence the set $I\defeq\{x\in\Aline oa:\Aline xz\cap\Aline ob\subseteq Y\setminus X\}$ also contains at most one point (here we use the fact that $z\notin\Aline ob$). Since the liner $Y$ is $3$-long, the intersection $\Aline oa\cap X$ contains at least two distinct points and hence it is a line in $X$. Since the liner $X$ is $3$-long, there exists a point $x\in(\Aline oa\cap X)\setminus(I\cup\{o\})$. The choice of $x\ne o$ ensures that $a\in\Aline oa=\Aline ox$ and hence $z\in\Aline ab\subseteq\overline{\{o,b,x\}}=\overline{\Aline ob\cup\{x\}}$. By the strong regularity of the projective liner $Y$, there exists a point $y\in\Aline ob$ such that $z\in\Aline  xy$. Since $x\notin I$, the point $y$ does not belong to the flat $Y\setminus X$. Then $x\in\Aline oa\cap X\subseteq{\overline A}\cap X=A$ and $y\in\Aline ob\cap X$ are two points with $z\in \Aline xy$, witnessing that the liner $X$ is regular. 
\smallskip

3. Since $Y\setminus X=\overline{Y\setminus X}\ne Y$, the liner $X$ contains some point $o\in X$. Since the liner $Y$ is $3$-long, the choice of the point $o\notin Y\setminus X$ implies the equality $Y=\bigcup_{x\in X}\Aline ox=\overline{X}$. By (already proved) Proposition~\ref{p:projective-minus-flat}(2), the liner $X$ is regular, and by Corollary~\ref{c:proregular=>ranked}, the regular liner $X$ is ranked. Using Lemma~\ref{l:Max-indep}, choose any maximal independent set $M$ in the ranked liner $X$. Proposition~\ref{p:add-point-to-independent} implies that ${\overline M}'=X$. Proposition~\ref{p:projective-minus-flat}(1) implies that the set $M$ remains independent in the projective liner $Y$. Since $Y=\overline{X}={\overline M}$, the set $M$ is maximal independent in the liner $Y$. By Theorem~\ref{t:Max=codim}, $\|Y\|=|M|=\|X\|$.
\end{proof}

\begin{corollary}\label{c:procompletion-rank}The rank $\|X\|$ of a $3$-long liner $X$ coincides with the  rank $\|Y\|$ of any projective completion $Y$ of $X$.
\end{corollary}

\begin{proof} Let $X$ be a $3$-long liner and $Y$ be a projective completion of $X$. If the horizon $H\defeq Y\setminus X$ of $X$ is flat in $Y$, then the ranks of the liners $X$ and $Y$ coincide, by Proposition~\ref{p:projective-minus-flat}(3). 
So, assume that $H$ is not flat and consider the subliner $Z\defeq Y\setminus\overline{H}$ of the liners $X$ and $Y$. By Corollary~\ref{c:Avogadro-projective}, the $3$-long projective liner $Y$ is $2$-balanced. Since the set $H$ is not flat in $Y$, it is not empty and hence contains some point $u\in H$. Then for any point  $o\in Z=Y\setminus\overline{H}$, the line $\Aline ou\subseteq Y$ has $\Aline ou\cap\overline{H}=\{u\}$ and hence $\Aline ou\setminus\{u\}$ is a line in $X$. Since the liner $X$ is $3$-long, $|Y|_2=|\Aline oy|=|\Aline oy\setminus\{y\}|+1\ge 3+1=4$. This shows that the projective liner $Y$ is $4$-long and hence the liner $Z=Y\setminus\overline{H}$ is $3$-long. By Proposition~\ref{p:projective-minus-flat}, the $3$-long liner $Z$ is regular and $\|Z\|=\|Y\|$.

Since $o\in Z=Y\setminus\overline{Y\setminus X}$, for every $y\in Y$ the set $\Aline oy\cap Z$ is a line in $Z$, which implies $Y=\bigcup_{z\in Z}\Aline oz\subseteq \overline{Z}$ and hence $Y=\overline{Z}=\overline{X}$. By the definition of the rank $\|X\|$, there exists a set $B\subseteq X$ of cardinality $|B|=\|X\|$ whose flat hull $\overline{B}'$ in $X$ coincides with $X$. Then $Y=\overline{X}={\overline B}$ and hence $\|Y\|\le|B|=\|X\|$. By analogy we can show that $\|X\|\le\|Z\|$. Then $\|Y\|\le\|X\|\le\|Z\|=\|Y\|$ and hence $\|X\|=\|Z\|=\|Y\|$.
\end{proof}

By Proposition~\ref{p:projective-minus-proflat}, for any proflat set $H$ in a projective liner $Y$, the subliner $X\defeq Y\setminus H$ is Proclus and hence $3$-proregular. The following example shows that in general the $3$-proregularity of $Y\setminus H$ cannot be improved to the $3$-regularity.

\begin{example} Every $2$-element set $H$ is a Steiner projective plane $P$ is proflat in $P$ but the liner $P\setminus H$ is not $3$-regular.

\begin{picture}(80,80)(-160,-10)
\linethickness{0.6pt}
\put(0,0){\line(1,0){60}}
\put(26,-10){$A$}
\put(0,0){\line(0,1){60}}
\put(-9,28){$\Lambda$}
\put(0,60){\line(1,-1){60}}
\put(0,0){\line(1,1){30}}
\put(0,60){\line(1,-2){20}}
\put(60,0){\line(-2,1){40}}

\put(0,0){\circle*{3}}
\put(60,0){\circle*{3}}
\put(0,60){\circle*{3}}
\put(20,20){\color{red}\circle*{3}}
\put(30,30){\circle*{3}}
\end{picture}

\end{example} 

Nonetheless we have the following proposition.

\begin{proposition} If a  $3$-long liner $X$ of rank $\|X\|\le3$ has a projective completion, then  $X$ is regular.
\end{proposition}

\begin{proof} Let $Y$ be a projective completion of $X$ and let $H\defeq Y\setminus X$ be a horizon in $X$. By Corollary~\ref{c:procompletion-rank}, $\|Y\|=\|X\|\le 3$. Since the projective liner $Y$ is ranked, $\|H\|<\|Y\|=3$.

 To prove that $X$ is regular, fix any flat $A\subseteq X$ and line $\Lambda$ in $X$ such that $A\cap\Lambda=\{o\}$ for some point $o\in X$. We have to prove that the set $\Aline A\Lambda\defeq\bigcup_{a\in A}\bigcup_{b\in\Lambda}(\Aline ab\cap X)$ is flat in $X$. If $A=\{o\}$, then $\Aline A\Lambda=\Lambda$ is flat. So, assume that $A\ne\{o\}$ and hence $\overline A$ is a line in the projective plane $Y$ and $A=\overline A\cap X$ is a line in $X$. 
 
 We shall prove that $\Aline A\Lambda=X$. Given any point $x\in X$, we need to check that $x\in\Aline A\Lambda$. This inclusion is trivial if $x\in A\cup\Lambda$. So, we assume that $x\notin A\cup\Lambda$.
 
 Depending on the location of the point $o$, we consider two cases.
 \smallskip
 
1. First we assume that $o\notin\overline H$.  In this case $|\overline \Lambda\cap\overline H|\le 1$ and $|\overline A\cap\overline H|\le 1$. Next, consider three subcases.
\smallskip

1.1. If  $|\overline \Lambda\cap\overline H|=1$, then let $\lambda$ be the unique point of the intersection $\overline\Lambda\cap\overline H$. Since $Y$ is a projective plane, there exists a unique point $\alpha\in \overline{\{\lambda,x\}}\cap\overline A\subseteq Y$. Since $X$ is $3$-long, there exists a point $a\in A\setminus\{o,\alpha\}$. By the projectivity of $Y$, there exists a point $b\in \overline\Lambda\cap\overline{\{a,x\}}$. It follows from $a\ne\alpha$ and $x\notin A\cup\Lambda$ that $b\in \overline\Lambda\setminus\{\lambda\}=\Lambda$. Then $x\in\Aline ab\subseteq\Aline A\Lambda$ and we are done.
\smallskip

1.2. By analogy we can prove that $x\in \Aline A\Lambda$ if $|\overline A\cap\overline H|=1$.
\smallskip

1.3. If $\overline \Lambda\cap\overline H=\varnothing=\overline A\cap\overline H$, then $\overline \Lambda=\Lambda$ and $\overline A=A$. Choose any point $a\in A\setminus\{o\}$. By the projectivity of $Y$, there exists a point $b\in \Lambda\cap\overline{\{a,x\}}$. Then $x\in\Aline ab\subseteq\Aline A\Lambda$ and we are done.
\smallskip

2. Next, assume that $o\in\overline H$. If $A\subseteq \overline H$, then choose any point $a\in A\setminus\{o\}$. By the projectivity of $Y$, the lines $\overline{\{a,x\}}$ and $\overline\Lambda$ have a common point $b\in Y$. Assuming that $b\in \overline H$, we conclude that $x\in\overline{\{a,b\}}\subseteq\overline H=\overline A$, which contradicts the choice of $x\notin A\cup\Lambda$. Therefore, $b\in\overline\Lambda\setminus\overline H\subseteq \Lambda$. Then $x\in\Aline ab\subseteq\Aline A\Lambda$. By analogy we can prove that $x\in\Aline A\Lambda$ if $\overline\Lambda\subseteq\overline H$.  So, assume that $A\not\subseteq\overline H$ and $\Lambda\not\subseteq\overline H$. Choose any point $a\in A\setminus\{o\}$. By the projectivity of $Y$, the lines $\overline{\{a,x\}}$ and $\overline\Lambda$ have a common point $b$. Assuming that $b\notin X$, we conclude that $b\in \overline\Lambda\cap H\subseteq\{o\}\subseteq X$, which is a contradiction showing that $b\in \overline\Lambda\cap X=\Lambda$. Then $x\in\Aline ab\subseteq\Aline A\Lambda$.
\end{proof}

\begin{theorem}\label{t:regular<=>flat4} Let $Y$ be a projective completion of a $3$-long proaffine liner $X$ of rank $\|X\|\ge 4$. The liner $X$ is regular if and only if $X$ is weakly regular if and only if the horizon $Y\setminus X$ is flat in $Y$.
\end{theorem}

\begin{proof} To distinguish between flat hulls of sets in the liners $X$ and $Y$, for sets $A\subseteq X$ and $B\subseteq Y$, we denote by $\overline{A}^X$ and $\overline{B}$ the flat hulls of the sets $A$ and $B$ in the liners $X$ and $Y$, respectively. Observe that for any points $a,b\in X$, we have $\overline{\{a,b\}}^X=\overline{\{a,b\}}\cap X=\Aline ab\cap X$.
\smallskip

If $X$ is regular, then $X$ is weakly regular. If the horizon $Y\setminus X$ is flat in $Y$, then the liner $X$ is regular, by Proposition~\ref{p:projective-minus-flat}. It remains to prove that the weak regularity of $X$ implies that the horizon $H\defeq Y\setminus X$ is flat in $Y$. To derive a contradiction, assume that $H$ is not flat in $Y$.  Then there exist two distinct points $p,q\in H$ such that $\Aline pq\not\subseteq H$. Since $X$ is proaffine, the horizon $H$ is proflat (by Theorem~\ref{t:proaffine<=>proflat}) and hence $\Aline pq\setminus H$ is a singleton containing a unique point $o$. Since $\overline H\ne Y$, there exists a point $a\in Y\setminus\overline{H}$. Since $Y$ is $3$-long, there exists a point $b\in \Aline ap\setminus\{a,p\}$. It follows from $p\in\overline H$ and $a\notin\overline H$ that $b\notin\overline H$. 

Corollary~\ref{c:procompletion-rank} ensures that $\|Y\|=\|X\|\ge 4$, which implies that the plane $A\defeq\overline{\{a,o,b\}}$ does not coincide with $Y$. Taking into account that $Y$ is $3$-long, we can apply Proposition~\ref{p:cov-aff} and find a point $u\in Y\setminus(\overline{H}\cup A)$. It follows from $p\in\overline{H}$ and $u\notin\overline{H}$ that $\Aline up\cap\overline{H}=\{p\}$. Since the liner $Y$ is $3$-long, there exists a point $z\in \Aline up\setminus \{u,p\}$. It follows from $u\notin\overline H$ that $|\Aline ua\cap\overline H|\le 1$ and $|\Aline ub\cap\overline H|\le 1$. Observe that $\{o,p\}\subseteq A\cap\overline{\{z,o,u\}}$ and hence $\Aline op\subseteq A\cap\overline{\{z,o,u\}}$. We claim that $\Aline op=A\cap\overline{\{z,o,u\}}$. In the opposite case, by the rankedness of the projective liner $Y$, we obtain $A=\overline{\{z,o,u\}}$, which contradicts the choice of $u$. This contradiction shows that $A\cap\overline{\{z,o,u\}}=\Aline op$.

\begin{picture}(200,180)(-130,-75)

\put(0,0){\line(1,0){120}}
\put(0,0){\line(1,-1){60}}
\put(60,-60){\line(1,1){60}}
\put(0,0){\color{red}\line(3,-1){90}}
\put(0,0){\color{red}\line(-3,1){30}}
\put(0,0){\line(0,1){90}}
\put(60,-60){\line(-2,5){60}}
\put(90,-30){\line(-3,4){90}}
\put(120,0){\line(-4,3){120}}
\put(30,15){\line(1,1){30}}

\put(0,0){\circle*{3}}
\put(-3,-8){$o$}
\put(-30,10){\color{red}\circle*{3}}
\put(-38,8){$q$}
\put(30,15){\circle*{3}}
\put(21,11){$x$}
\put(45,30){\circle*{3}}
\put(49,28){$z$}
\put(60,45){\circle*{3}}
\put(62,48){$y$}
\put(60,-60){\circle*{3}}
\put(58,-68){$a$}
\put(90,-30){\color{red}\circle*{3}}
\put(93,-33){$p$}
\put(120,0){\circle*{3}}
\put(123,-3){$b$}
\put(0,90){\circle*{3}}
\put(-2,93){$u$}
\end{picture}

Since the liner $Y$ is projective and $z\in \Aline up\setminus(\Aline ua\cup\Aline ub)\subseteq\overline{\{a,u,b\}}\setminus(\Aline ua\cup\Aline ub)$, for every $x\in\Aline ua$, the line $\Aline xz$ intersects the line $\Aline ub$. Since $|\Aline ub\cap\overline H|\le 1$, the set $I\defeq\{x\in\Aline ua:\Aline xz\cap\Aline ub\subseteq\overline H\}$ has cardinality $|I|\le 1$. Since $u\notin\overline H$, the intersection $\Aline ua\cap\overline{H}$ contains at most one points. Since the liner $Y$ is $3$-long, the intersection $\Aline ua\cap X$ is a line in $X$. Since the liner $X$ is $3$-long, there exists a point $x\in (\Aline ua\cap X)\setminus(\{u\}\cup I)$. Let $y$ be the unique point of the intersection $\Aline xz\cap\Aline ub$ of the lines $\Aline xz$ and $\Aline ub$ in the projective plane $\overline{\{u,a,b\}}$. The choice of $x\in \Aline ua\setminus I$ ensures that $y\notin \overline H$. Then the points $o,a,b,u,x,y,z$ belong to the liner $X$ and hence $z\in X\cap\Aline xy\subseteq \overline{(A\cap X)\cup\{u\}}^X$. Since the liner $X$ is weakly regular, there exists  a point $c\in A\cap X$ such that $z\in\overline{\{c,o,u\}}^X\subseteq\overline{\{c,o,u\}}$. Since the projective liner $Y$ is ranked,  $$c\in  X\cap A\cap\overline{\{z,o,u\}}\subseteq X\cap(A\cap\overline{\{z,o,u\}})=X\cap \Aline op=\{o\}$$ and hence $z\in\Aline up\cap\overline{\{c,o,u\}}=\Aline up\cap \overline{\{o,u\}}=\{u\}$, which contradicts the choice of $z$. This contradiction shows that the set $H$ is flat in $X$.
\end{proof}

\begin{proposition}\label{p:projective-minus-hyperplane} For every hyperplane $H$ in a projective liner $Y$, the subliner $X\defeq Y\setminus H$ of $X$ is affine and regular.
\end{proposition}

\begin{proof} To see that the liner $X$ is affine, we need to show that for any points $o,x,y\in X$ and $p\in X\cap\Aline yx\setminus\Aline ox$, the set $I\defeq\{v\in X\cap\Aline oy:(X\cap\Aline vp)\cap(X\cap \Aline ox)=\varnothing\}$ is a singleton. By Proposition~\ref{p:projective-minus-proflat}, the liner $X$ is proaffine and hence $|I|\le 1$. So, it remains to show that $I$ is not empty. If $p\in\Aline oy$, then $p\in I$ and we are done. So, assume that $p\notin \Aline oy$. It follows from $p\in \Aline yx\setminus\Aline oy$ that $o\ne x$ and hence $\Aline ox$ is a line in the projective liner $X$. Since $H$ is a hyperplane in $X$, $\|X\|_H=1$ and hence $\|\Aline ox\|_H\le\|X\|_H=1<2=\|\Aline ox\|$. 
By  Proposition~\ref{p:projective<}, the intersection $\Aline ox\cap H$ contains some point $z$. It follows from $p\in {\Aline yx}\setminus{\Aline ox}$ that $y\notin {\Aline ox}$ and hence $\Aline oy\cap\Aline ox=\{o\}$. Since $\Aline zp$ and $\Aline oy$ are two lines in the plane $\overline{\{o,x,y\}}$, the $0$-parallelity of the projective liner $X$ ensures that the intersection $\Aline zp\cap\Aline oy$ contains some point $v$. Assuming that $v=z$, we conclude that $v=z\in\Aline ox\cap\Aline oy\cap H=\{o\}\cap H=\varnothing$, which is a contradiction showing that $v\ne z$ and hence $p\in\Aline zp=\Aline vz$. Assuming that $v\in H$, we conclude that $p\in\Aline vz\subseteq H$, which contradicts the choice of $p\in X$. This contradiction shows that $v\in Y\setminus H=X$. Assuming that $v\in\Aline ox$, we conclude that $p\in\Aline vz\cap X\subseteq \Aline ox\cap X$, which contradicts the choice of the point $p$. This contradiction shows that $v\notin \Aline ox$ and hence $\Aline vp\cap\Aline ox=\Aline vz\cap\Aline ox=\{z\}$. Then $$(\Aline vp\cap X)\cap(\Aline ox\cap X)=(\Aline vp\cap\Aline ox)\cap X=\{z\}\cap X=\emptyset,$$ witnessing that $v\in I\ne\varnothing$.

Next, we show that the affine liner $X$ is regular. By Theorem~\ref{t:affine=>Avogadro}, the affine liner $X$ is $2$-balanced. If $|X|_2=2$, then the liner $X$ is projective and hence (strongly) regular. So, assume that $|X|_2\ge 3$. Fix any point $o\in X$ and using the Kuratowski--Zorn Lemma, find a maximal $3$-long flat $M$ containing the point $o$ in the  projective liner $Y$. Assuming that $X\not\subseteq M$, we could find a point $x\in X\setminus M\subseteq X\setminus\{o\}$. Lemma~\ref{l:ox=2} ensures that  $\Aline o x=\{o,x\}$, which contradicts $|X|_2\ge 3$. This contradiction shows that $X\subseteq M$. Then the flat hull $\overline X$ in $M$ is a $3$-long projective liner and $\overline X\setminus X=\overline X\cap(Y\setminus X)=\overline X\cap H$ is a flat in the $3$-long projective liner $\overline X$. Then $\overline X$ is a projective completion of the $3$-long affine liner $X$ and the horizon $\overline X\setminus X$ of $X$ in $\overline X$ is flat. By Proposition~\ref{p:projective-minus-flat}(2), the affine liner $X$ is regular.
\end{proof}

\begin{theorem}\label{t:regular-horizon3} Let $Y$ be a projective completion of a liner $X$ such that $|Y|_2=3$. The liner $X$ is regular if and only if $|Y\setminus X|\le 1$ or $Y\setminus X$ is a hyperplane in $Y$.
\end{theorem}

\begin{proof} Let $H\defeq Y\setminus X$ be the horizon of $X$ in $Y$. To prove the ``if'' part, assume that $|H|\le 1$ or $H$ is a hyperplane in $Y$. In the latter case, $X$ is regular, by Proposition~\ref{p:projective-minus-hyperplane}.  

\begin{claim}\label{cl:Pi-minus-H=regular} If $|H|\le 1$, then for every plane $\Pi\subseteq T$, the subliner $Z\defeq \Pi\setminus H$ of $Y$ is regular.
\end{claim}

\begin{proof}  To prove that $Z$ is regular, take any flat $A$ in $Z$ and any line $\Lambda\subseteq Z$ such that $\Lambda\cap A=\{o\}$ for some point $o\in Z$. We have to prove that the set $\Aline A\Lambda\defeq\bigcup_{x\in A}\bigcup_{y\in\Lambda}\Aline xy\cap Z$ is flat in $Z$.
In the opposite case, we can find points $a,b\in\Aline A\Lambda$ such that $\Aline ab\cap X\not\subseteq \Aline A\Lambda$ and hence $c\notin\Aline A\Lambda$ for some point $c\in\Aline ab\cap X$. It follows that $\{a,b\}\not\subseteq A\cup\Lambda$ and we lose no generality assuming that $b\notin A\cup\Lambda$.  Since $b\in\Aline A\Lambda$, there exist points $\beta\in A$ and $\lambda\in\Lambda$ such that $b\in\Aline \beta \lambda$. It follows from $b\notin A\cup\Lambda$ that $\beta\ne o\ne\lambda$.  

Since the liner $Y$ is Steiner, there exists a unique point $\mu\in\Aline o\lambda\setminus\{o,\lambda\}\subseteq \Pi$. Consider the plane $P\defeq\overline{\{o,\beta,\lambda\}}\subseteq\Pi$ and observe that $b\in\Aline\beta\lambda\subseteq P$. Since the projective liner $Y$ is ranked, $\|P\|=3=\|\Pi\|$ implies $P=\Pi$ and hence $c\in Z\subseteq \Pi=P$. By projectivity of the plane $P=\Pi$, there exist points $\gamma\in \Aline \lambda c\cap\Aline o\beta$ and $\delta\in\Aline\mu c\cap\Aline o\beta$. It follows from $c\notin\Lambda\cup A$ that $c\notin \Aline o\lambda\cup \Aline o\beta$ and hence $\gamma\ne\delta$. Assuming that $\gamma\in Z$, we conclude that $c\in\Aline\gamma\lambda\cap Z\subseteq\Aline A\Lambda$, which contradicts the choice of $c$. This contradiction shows that $\gamma\in \Pi\setminus Z\subseteq H$ and hence $\gamma\in H$ and $H=\{\gamma\}$ (because $|H|\le 1$). Then $\delta\in \Pi\setminus H=Z$,  $\mu \in\Aline o\lambda\setminus H\subseteq Z$, and $c\in\Aline \mu\delta\cap Z\subseteq\Aline A\Lambda$, which contradicts the choice of $c$. This contradiction shows that the set $\Aline A\Lambda$ is flat and the liner $Z$ is regular.
\end{proof}

\begin{claim} If $|H|\le 1$, then the subliner $X=Y\setminus H$ of $Y$ is regular.
\end{claim}

\begin{proof}  To prove that $X$ is regular, take any flat $A$ in $X$ and any line $\Lambda\subseteq X$ such that $\Lambda\cap A=\{o\}$ for some point $o\in X$. We have to prove that the set $\Aline A\Lambda\defeq\bigcup_{x\in A}\bigcup_{y\in\Lambda}\Aline xy\cap X$ is flat in $X$.
In the opposite case, we can find points $a,b\in\Aline A\Lambda$ such that $\Aline ab\cap X\not\subseteq \Aline A\Lambda$ and hence $c\notin\Aline A\Lambda$ for some point $c\in\Aline ab$. It follows that $\{a,b\}\not\subseteq A\cup\Lambda$ and we lose no generality assuming that $b\notin A\cup\Lambda$.  Since $a,b\in\Aline A\Lambda$, there exist points $\alpha,\beta\in A$ and $u,v\in\Lambda$ such that $a\in\Aline \alpha u$ and $b\in\Aline \beta v$. It follows from $b\notin A\cup\Lambda$ that $\beta\ne o\ne v$. If $a\in\Lambda$, then we can assume  that $u=a$ and $\alpha=o$.  

Consider the plane $\Pi\defeq\overline{\{\beta,o,v\}}$ in $Y$. Since $\Pi\setminus X\subseteq Y\setminus X=H$ and $|H|\le 1$, the subliner $X\cap\Pi$ of the plane $\Pi$ is regular, by Claim~\ref{cl:Pi-minus-H=regular}. If $\alpha\in \Pi$, then by the regularity of the liner $X\cap\Pi$, there exist points $\lambda\in\Lambda$ and $\gamma\in A\cap\Pi$ such that $c\in\Aline \gamma\lambda\subseteq\Aline A\Lambda$, which contradicts the choice of $c$. This contradiction shows that $\alpha\notin \Pi$. In this case also $a\notin\Lambda$. It follows from $\alpha\notin\Pi=\overline{\{\beta,o,v\}}$ that $P\defeq\overline{\{\alpha,o,\beta\}}$ is a plane in $Y$. Since $|P\setminus X|\le|H|\le 1$, we can apply Claim~\ref{cl:Pi-minus-H=regular} and conclude that the subliner $P\cap X$ of $P$ is regular and hence $X\cap P$ coincides with the flat hull of the set $(\Aline o\alpha\cup\Aline o\beta)\cap X\subseteq A\cap X$ in $X\cap P$, which implies $X\cap P\subseteq A$ and $c\notin P$. Since the flat $S\defeq\overline{\{\alpha,o,\beta,v\}}$ has rank $4$, the plane $P$ is a hyperplane in $S$. By Corollary~\ref{c:line-meets-hyperplane}, the line $\{a,b,c\}\subseteq S$ intersects $P$. It follows from $P\cap X\subseteq A$ and $b,c\notin A$ that $b,c\notin P$ and hence $a\in P\cap X\subseteq A$. By the modularity of the projective liner $S$, the intersection of the planes $P$ and $\Gamma\defeq\overline{\{o,v,c\}}$ is a line $L$. It follows from $|H|\le 1$ that $L\cap X$ is a line in $P\cap X\subseteq A$. Then there exists a point $\gamma\in (L\cap X)\setminus\{o\}\subseteq A$. It follows from $|H|\le 1$ that the set $\Gamma\cap X$ coincides with the flat hull of the set $(\Aline ov\cup\Aline o\gamma)\cap X$ in $\Gamma\cap X$. By Claim~\ref{cl:Pi-minus-H=regular}, the subliner $X\cap\overline{\{o,v,c\}}$ of the projective plane $\overline{\{o,v,c\}}$ is regular. So, there exist points $\lambda\in \Aline ov\cap X=\Lambda$ and $\delta\in\Aline o\gamma\cap X\subseteq A$ such that $c\in \Aline \delta\lambda\subseteq\Aline A\Lambda$, which contradicts the choice of $c$.
This is the final contradiction completing the proof of the regularity of the liner $X$.
\end{proof}
    
To prove the ``only if'' part, assume that the liner $X$ is regular. We need to prove that $|Y\setminus X|\le 1$ or $Y\setminus X$ is a hyperplane in $Y$. To derive a contradiction, assume that $|Y\setminus X|\ge 1$ and $Y\setminus X$ is not a hyperplane in $Y$.

\begin{claim}\label{cl:horizon-flat} The horizon $H\defeq Y\setminus X$ is flat in $Y$. 
\end{claim}

\begin{proof} Assuming that $H$ is not flat in $Y$, we can find points $\alpha,\beta\in H$ and $\gamma\in\Aline ab\setminus H\subseteq X$. Since $\overline H\ne Y$, there exists     
a point $o\in Y\setminus\overline H\subseteq X$. Since $|Y|_2=3$, there exist points $a,b,c\in Y$ such that $\Aline o\alpha=\{o,\alpha,a\}$, $\Aline o\beta=\{o,\beta,b\}$ and $\Aline o\gamma=\{o,\gamma,c\}$. It follows from $\alpha,\beta,\gamma\in\overline H$ and $o\notin\overline H$ that $a,b,c\in Y\setminus\overline H\subseteq X$. Since $|Y|_2=3$, the plane $\overline{\{o,\alpha,\beta\}}$ contains exactly 7 points $\{o,a,b,c,\alpha,\beta,\gamma\}$ and $\{a,b,\gamma\}$ is a line in $\overline{\{o,\alpha,\beta\}}$. Consider the lines $A\defeq \{o,a\}$, $\Lambda\defeq\{o,b\}$ in $X$ and observe that $c\in\Aline o\gamma\subseteq \overline{A\cup\Lambda}$ but $c\notin\Aline A\Lambda$ in $X$, witnessing that the liner $X$ is not regular. This contradiction shows that the horizon $H$ is flat in $Y$.
\end{proof}

By Claim~\ref{cl:horizon-flat}, the horizon $H\defeq Y\setminus X$ is flat in $Y$. 
Since $|H|>1$, we can fix two distinct points $\{u,a_1\}\subseteq H$. Since $H=\overline H\ne Y$, there exists a point $o\in Y\setminus H=X$. Since $|Y|_2=3$, there exists a unique point $v\in Y$ such that $\Aline ou=\{o,u,v\}$. Since the flat $H$ is not a hyperplane in $Y$, the flat $\overline{H\cup\{o\}}$ does not coincide with $Y$ and hence there exists a point $a_2\in Y\setminus\overline{H\cup\{o\}}\subseteq Y\setminus H=X$. Find a unique point $a_3\in Y$ such that $\Aline {a_1}{a_2}=\{a_1,a_2,a_3\}$ and observe that $a_3\in Y\setminus\overline{H\cup\{o\}}\subseteq Y\setminus H=X$. For every $i\in\{1,2,3\}$, find points $b_i,c_i,d_i\in Y$ such that $\Aline u{a_i}=\{u,a_i,b_i\}$, $\Aline v{a_i}=\{v,a_i,c_i\}$ and $d_i\in \Aline u{c_i}\cap \Aline o{a_i}$. Since $H$ is flat, $b_1\in\Aline u{a_1}\subseteq H$. On the other hand, for every $i\in\{2,3\}$, it follows from $a_i\notin\overline{H\cup\{o\}}$ and $u\in H$ that $b_i,c_i,d_i\notin\overline{H\cup\{o\}}$.

\begin{picture}(80,85)(-150,-10)

\put(0,0){\line(1,0){60}}
\put(0,0){\line(0,1){60}}
\put(0,60){\line(1,-1){60}}
\put(0,0){\line(1,1){30}}
\put(0,60){\line(1,-2){30}}
\put(60,0){\line(-2,1){60}}

\put(0,0){\circle*{3}}
\put(-2,-9){$o$}
\put(60,0){\circle*{3}}
\put(63,-9){$a_i$}
\put(30,0){\circle*{3}}
\put(27,-10){$d_i$}
\put(0,30){\circle*{3}}
\put(-9,28){$v$}
\put(0,60){\circle*{3}}
\put(-2,63){$u$}
\put(20,20){\circle*{3}}
\put(15,11){$c_i$}
\put(30,30){\circle*{3}}
\put(33,30){$b_i$}
\end{picture}

\begin{claim}\label{cl:c1c2c3} The point $c_1,c_2,a_3$ are collinear.
\end{claim}

\begin{proof} Consider the plane $\Pi\defeq\overline{\{o,a_2,a_3\}}$ in $Y$ and observe that the points $a_1\in\Aline{a_2}{a_3}$, $d_2\in \Aline o{a_2}$, and $d_3\in \Aline o{a_3}$ belong to $\Pi$. By the projectivity of the plane $\Pi$, the lines $\Aline{d_1}{d_2}$ and $\Aline {a_1}{a_2}=\{a_1,a_2,a_3\}$ have a common point, equal to the point $a_3$. 
Now consider the plane $P\defeq\overline{\{u,d_1,a_3\}}$ in $Y$ and observe that $d_2\in\Aline {d_1}{a_3}\subseteq P$ and $\{c_1,c_2\}\subseteq\Aline u{d_1}\cup\Aline u{d_2}\subseteq P$. It follows from $u\notin \overline{\{o,a_1,a_2\}}$ that $c_1,c_2\notin \overline{\{o,a_1,a_2\}}$ and hence $\{c_1,c_2\}\cap\{d_1,d_2\}=\varnothing$. By the projectivity of the plane $P$, the lines $\Aline {c_1}{c_2}$ and $\Aline {d_1}{d_2}=\{d_1,d_2,a_3\}$ have a common point, equal to the point $a_3$. Therefore, the points $c_1,c_2,a_3$ are collinear.
\end{proof}

Consider the flat $A\defeq\overline{\{o,a_2,a_3\}}\cap X$ and the line $\Lambda\defeq\Aline ou\cap X=\{o,v\}$ in $X$. The rankedness of the projective liner $Y$ ensures that $\Lambda\cap A=\{o\}$. We claim that the set $\Aline A\Lambda\defeq\bigcup_{x\in A}\bigcup_{y\in\Lambda}\Aline xy\cap X$ is not flat in $X$, witnessing that the liner $X$ is not regular. To derive a contradiction, assume that $\Aline A\Lambda$ is flat in $X$. Then $v\in \Lambda\setminus\{o\}$ and $a_2\in A$ implies $c_2\in \Aline {a_2}v\cap X\subseteq \Aline A\Lambda$ and $c_2\notin \overline{\{o,a_2,a_3\}}$. Claim~\ref{cl:c1c2c3} ensures that $c_1\in \Aline {a_3}{c_2}$. Since the set $\Aline A\Lambda$ is flat in $X$, $c_1\in \Aline a\lambda$ for some $a\in A\subseteq X$ and $\lambda\in\Lambda=\{o,v\}$. It follows from  $c_1\notin A$ that $\lambda\in\Lambda\setminus\{o\}=\{v\}$. Then $a\in\Aline v {c_1}\cap A=\{a_1\}$ and hence $a=a_1\in X\cap H=\varnothing$, which is a contradiction completing the proof of the theorem.
\end{proof}

\begin{theorem}\label{t:affine<=>hyperplane} For a projective completion $Y$ of a liner $X$, the following conditions are equivalent:
\begin{enumerate}
\item the liner $X$ is affine and regular;
\item the liner $X$ is affine;
\item $Y\setminus X$ is a hyperplane in $Y$.
\end{enumerate}
\end{theorem}

\begin{proof}  The implication $(1)\Ra(2)$ is trivial. 
\smallskip

$(2)\Ra(3)$ Assume that the horizon $H\defeq Y\setminus X$ of $X$ is not a hyperplane in $Y$. If $H\ne\overline H$, then we can choose a point $x\in\overline{H}\setminus H$ and conclude that $\overline{H\cup\{x\}}\subseteq\overline H\ne Y$. If $H=\overline H$, then $\overline H$ is not a hyperplane in $Y$ and hence there exists a point $x\in Y\setminus \overline H$ such that $\overline{H\cup\{x\}}\ne Y$. In both cases, there exists  a point $x\in Y\setminus H$ such that $\overline{H\cup\{x\}}\ne Y$. Choose any point $o\in Y\setminus\overline{H\cup\{x\}}$ and observe that $\Aline xo\cap \overline{H\cup\{x\}}=\{x\}$ and hence $\Aline xo\cap H=\varnothing$.  Since the liner $Y$ is $3$-long, we can apply Proposition~\ref{p:cov-aff} and find a point $y\in Y\setminus(\overline{H\cup\{x\}}\cup\Aline ox)$. It follows that $\Aline yx\cap\overline{H \cup\{x\}}=\{x\}$ and hence $\Aline yx\cap H=\varnothing$. Since the liner $Y$ is $3$-long, there exists a point $p\in\Aline xy\setminus\{x,y\}$. The choice of the point $y\notin\Aline ox$ ensures that $p\notin \Aline ox$. For every point $v\in\Aline oy\cap X$, the line $\Aline vp$ has non-empty intersection with the line $\Aline ox$, by the projectivity of the liner $Y$. Since $\Aline ox\cap H=\varnothing$, the nonempty intersection $\Aline ox\cap\Aline vp$ coincides with the intersection of the lines $\Aline vp\cap X$ and $\Aline ox\cap X=\Aline ox$, witnessing that the liner $X$ is not affine.
\smallskip

The implication $(3)\Ra(1)$ follows from Proposition~\ref{p:projective-minus-hyperplane}.
\end{proof}

\section{The (non)existence of projective completions}

In this section we detect liners that have (or do not have) projective completions.

\begin{theorem}\label{t:procompletion<=>normcompletionA} A non-projective $3$-long $3$-ranked liner $X$ has a projective completion if and only if $X$ is a normal completion of some $3$-long affine regular liner $A\subseteq X$.
\end{theorem}

\begin{proof} To prove the ``if'' part, assume that a $3$-long $3$-ranked liner $X$ is a normal completion of some $3$-long affine regular liner $A\subseteq X$. By Corollary~\ref{c:affine-spread-completion}, the spread completion $\overline A$ of $A$ is a projective completion of $A$. By Theorem~\ref{t:projective=>extension}, there exists a liner embedding $F:X\to \overline A$ such that $F(a)=a$ for all $a\in A$. Choose any set $Y$ containing the set $X$ so that the remainder $Y\setminus X$ admits a bijective map $f:Y\setminus X\to \overline A\setminus F[X]$. Consider the function $\overline F:Y\to \overline A$ defined by $$\overline F(y)=\begin{cases}F(y)&\mbox{if $y\in X$};\\
f(y)&\mbox{if $y\in Y\setminus X$}.
\end{cases}
$$
Endow the set $Y$ with the line relation
$$\Af_Y\defeq\{xyz\in Y^3:\overline Fxyz\in \overline \Af\},$$where $\overline \Af$ stands for the line relation on the liner $\overline A$. Such definition of the line relation $\Af_Y$ ensures that $\overline F:Y\to \overline A$ is a liner isomorphism between the liners $Y$ and $\overline A$. Since the liner $\overline A$ is projective, so is its isomorphic copy $Y$. Moreover, $\overline{\overline A\setminus A}=\overline A\setminus A\ne \overline A$ implies $\overline{Y\setminus X}\subseteq\overline{Y\setminus A}\ne Y$, witnessing that the projective liner $Y$ is a projective completion of the liner $X$.
\smallskip

To prove the ``only if'' part, assume that a non-projective $3$-long $3$-ranked liner $X$ has a projective completion $Y$. By Corollary~\ref{c:Avogadro-projective},  the $3$-long projective liner $Y$ is $2$-balanced. Since $X$ is not projective, the horizon $H\defeq Y\setminus X$ of $X$ in $Y$ is not empty and hence contains some point $y\in H=Y\setminus X$. Choose any point $x\in X$ and observe that $|Y|_2\ge |\Aline xy|=|\Aline xy\cap X|+1\ge 4$. By the Kuratowski-Zorn Lemma, the proper flat $\overline H$ in $Y$ can be enlarged to a maximal proper flat $M$ in the liner $Y$.  The maximality of $M$ ensures that $M$ is a hyperplane in $Y$. By Proposition~\ref{p:projective-minus-hyperplane}, the liner $A\defeq Y\setminus M$ is affine. Since the projective liner $Y$ is $4$-long and $M$ is a flat in $Y$, the affine liner $A=Y\setminus M$ is $3$-long. By Proposition~\ref{p:projective-minus-flat}, the liner $A=Y\setminus M$ is regular. Since $\varnothing\ne Y\setminus M\subseteq Y\setminus H=X$, the completement $X\setminus M=Y\setminus M$ is not empty. Observe that the set $X\setminus A=X\cap M$ is flat in $X$ and hence $X$ is a completion of the affine liner $A$. 

To see that $X$ is a normal completion of the liner $A$, take any plane $P$ in $X$ such that $P\not\subseteq  \overline{X\setminus A}=M\cap X$. We have to prove that $P\cap A$ is a plane in the liner $A$. Choose any points $a\in A\cap P$, $b\in P\setminus\{a\}$ and $c\in P\setminus\Aline ab$.
Since $a\notin M$, the intersection $\Aline ab\cap M$ contains at most one point. Since the liner $X$ is $3$-long, we can replace the point $a$ by any point in the (nonempty) sets $\Aline ab\setminus(\{a\}\cup M)$ and assume that $a\in P\setminus M$. By analogy we can replabe the point $b$ by a suitable point in the line $\Aline ac\setminus(\{a\}\cup M)$ and assume that $c\in P\cap A$.
 
Since the liner $X$ is $3$-ranked, the plane $P$ coincides with the flat hull of the set $\{a,b,c\}$ in $X$.
Let $\Pi$ be the flat hull of the set $\{a,b,c\}$ in the projective liner $Y$. It is easy to see that $P\subseteq\Pi$.  Since the points $a,b,c$ are non-colinear, the flat $\Pi$ is a plane in $Y$. Since $\{a,b,c\}\subseteq P\cap A=P\setminus M\subseteq \Pi\setminus M$, the plane $\Pi$ is not contained in the flat $M$. By Theorem~\ref{t:procompletion=>normal}, the projective completion $Y$ of the $3$-long liner $A$ is normal and hence $\Pi\cap A$ is a plane in $A$. Since the points $a,b,c$ are non-collinear and $\{a,b,c\}\subseteq P\cap A\subseteq \Pi\cap A$, the set $P\cap A$ is a plane in the liner $A$, witnessing that the completion $X$ of the liner $A$ is normal.
\end{proof}





\begin{example}\label{ex:Proclus-without-completion} {\em There exists an $\w$-long Proclus plane having no projective completions.}
\end{example}

\begin{proof} By Example~\ref{ex:not-para-Playfair}, there exists an $\w$-long Proclus plane $X$ containing two disjoint lines $L,L'$ such that any two lines $A,B\subseteq X$ with $\{A,B\}\ne\{L,L'\}$ have nonempty intersection. We claim that the liner $X$ has no projective completions. To derive a contradiction, assume that the liner $X$ has a projective completion $Y$. Since the $3$-long liner $X$ is a plane, so is its projective completion $Y$. Then the flat hulls $\overline L$ and $\overline{L'}$ of the disjoint lines $L,L'$ in the liner $X$ have a common point $y\in\overline L\cap\overline{L'}\subseteq Y\setminus X$ in the projective liner $Y$. Since $L$, $L'$ and $X\cap\overline{Y\setminus X}$ are proper flats in the $\w$-long liner $X$, there exists a point $x\in X\setminus (L\cup L'\cup\overline{Y\setminus X})$, according to Proposition~\ref{p:cov-aff}. Then the line $\Aline xy\cap X$ in the liner $X$ is disjoint with the lines $L,L'$, which contradicts the choice of the lines $L,L'$.
\end{proof}


\begin{example}\label{ex:para-Playfair}{\em There exists an $\w$-long para-Playfair plane having no projective completions.}
\end{example}

\begin{proof} First we fix a suitable terminology.  Every function $F:X\to Y$ is identified with its graph $\{(x,y)\in X\times Y:y=F(x)\}$. A function $F$ is {\em injective} if the relation $F^{-1}\defeq\{(y,x):(x,y)\in F\}$ is a function. Two functions $F,G$ are called {\em transversal} if $|F\cap G|=1$. A family of functions $\F$ is called {\em transversal} if any two distinct functions $F,G\in\F$ are transversal. A bijective function $F:\w\to\w$ will be called an {\em $\w$-bijection}. 

\begin{lemma}\label{l:transversal} Let $\F$ be a finite transversal family of $\w$-bijections. For any numbers $x,x'\in\w$ and $y,y'\in\w$ with $x\ne x'$ and $y\ne y'$ there exists an $\w$-bijection $\Phi$ such that
$\{(x,y),(x',y')\}\subseteq\Phi$ and the family $\F\cup\{\Phi\}$ is transversal.
\end{lemma}

\begin{proof} If $\{(x,y),(x',y')\}\subseteq F$ for some $F\in\F$, then put $\Phi\defeq F$ and finish the proof. So, assume that $\{(x,y),(x',y')\}\not\subseteq F$ for every $F\in\F$. Let $m\defeq |\F|$ and $\{F_i\}_{i\in m}=\F$ be an enumeration of the set $\F$. 

Let $x_0\defeq x$, $x_1\defeq x'$, $y_0\defeq y$, $y_1\defeq y'$. For every $i\in m$ choose a pair $(x_{i+2},y_{i+2})\in F_i\setminus(\bigcup_{j<i}F_j)$ such that $x_i\notin\{x_j:j<i\}$ and $y_i\notin\{y_j:j<i\}$. The choice of the pair $(x_{i+2},y_{i+2})$ is always possible because the family $\F$ is transversal and hence the set $\bigcup_{j<i}(F_i\cap F_j)$ is finite and the set $F_i\setminus\bigcup_{j<i}F_j$ is infinite. After completing the inductive construction, we obtain two sequences $(x_i)_{i\in m+2}$ and $(y_i)_{i\le m+2}$ such that $\varphi\defeq \{(x_i,y_i)\}_{i\in m+2}$ is an injective function such that $\{(x,y),(x',y')\}\subseteq \varphi$ and $|\varphi\cap F|=1$ for every $F\in\F$.  

Next, choose a sequence $(x_i)_{i=m+2}^\infty$ of pairwise distinct elements of $\w$ such that $\{x_i:i\ge m+2\}=\w\setminus\{x_i\}_{i<m+2}$. For every number $i\ge m+2$ let $y_i$ be the smallest number in the infinite set $\w\setminus(\{y_j:j<i\}\cup\{F(x_i):F\in\F\})$. Then $\Phi\defeq\{(x_i,y_i)\}_{i\in\w}$ is an injective function from $\w$ to $\w$ such that $\varphi\subseteq\Phi$ and $|\Phi\cap F|=|\varphi\cap F|=1$ for every $F\in\F$. It remains to show that the function $\Phi$ is bijective. In the opposite case, the set $\w\setminus\Phi[\w]$ is not empty and we can consider its smallest element $s\defeq\min(\w\setminus\Phi[\w])$. Since the finite set $\{0,\dots,s-1\}$ is contained in $\Phi[\w]$, there exists a number $k\ge m+2$ such that $\{0,\dots,s-1\}\subseteq \{\Phi(x_i):i<k\}$. For every $i\ge k$, the choice of $s\ne y_i=\min(\w\setminus(\{y_j:j<i\}\cup\{F(x_i):F\in\F\})$ ensures that $s=\{F(x_i):F\in\F\}$ and hence $s\in \Psi_i(x_i)$ for some $\Phi_i\in\F$. Since $|\F|<\w=|\w\setminus k|$, there exist two distinct numbers $i,j\in \w\setminus k=\w\setminus\{0,1,\dots,k-1\}$ such that $\Psi_i=\Psi_j$ and hence $\Psi_i(x_i)=s=\Psi_j(x_j)=\Psi_i(x_j)$, which implies that the function $\Psi_i=\Psi_j\in\F$ is not injective, which contradicts the choice of the family $\F$. This contradiction shows that the function $\Phi:\w\to \w$ is bijective. 
\end{proof}

Let $\mathcal B$ be a set of injective functions $B$ such that $|B|=2$ and $\dom[B]\cup\rng[B]\subseteq \w$.
Every function $B\in\mathcal B$ is a set of the form $\{(x,y),(x',y')\}$ for some numbers $x,x',y,y'\in\w$ with $x\ne x'$ and $y\ne y'$. Let $\{B_n\}_{n\in\w}$ be an enumeration of the countable set $\mathcal B$ such that $B_n\ne B_m$ for any distinct numbers $n,m\in\w$. Using Lemma\ref{l:transversal}, choose a sequence of $\w$-bijections $(F_n)_{n\in\w}$ such that for every $n\in\w$, the family $\{F_i\}_{i\le n}$ is transversal and $B_n\subseteq F_n$.

Now consider the liner $X\defeq \w\times\w$ endowed with the family of lines $$\mathcal L\defeq\{\w\times\{y\}:y\in\w\}\cup\{\{x\}\times \w:x\in\w\}\cup\{F_n:n\in\w\}.$$ The transversality of the family $\{F_n\}_{n\in\w}$ ensures that any distinct lines $L,\Lambda\in\{F_n:n\in\w\}$ are concurent.  The liner $X$ has only two spreads of parallel lines: $\{\w\times\{y\}:y\in\w\}$ and $\{\{x\}\times\w:x\in\w\}$, witnessing that the liner $X$ is para-Playfair. To see that the liner $X$ is a plane, consider the points $o\defeq(0,0)$, $a\defeq(1,0)$ and $b\defeq(0,1)$. We claim that $X=\overline{\{a,o,b\}}$. Indeed, given any point $p\in X\setminus(\Aline oa\cup\Aline ob)$, find a number $y\in\w$ such that $p\in\w\times\{y\}$ and choose any point $u\in\Aline ob\setminus(\w\times\{y\})$. By the choice of the family $(F_n)_{n\in\w}$ there exists $n\in\w$ such that $\{p,u\}\subseteq F_n\in\mathcal L$. Since the function $F_n:\w\to \w$ is surjective, there exists a point $v\in F_n\cap(\w\times\{0\})=F_n\cap\Aline oa$. Then $p\in F_n=\Aline uv\subseteq\overline{\{a,o,b\}}$, witnessing that the liner $X$ is a plane.

We claim that the liner $X$ has no projective completions. To derive a contradiction, assume that the liner $X$ has a projective completion $Y$. Since the liner $X$ is not projective, the horizon $H\defeq Y\setminus X$ of $X$ in $Y$ is not empty. Since $X$ is Proclus, the horizon $H\defeq Y\setminus X$ of the completion $Y$ is proflat, by Theorem~\ref{t:proaffine<=>proflat}. By Corollary~\ref{c:procompletion-rank}, $\|Y\|=\|X\|=2$, which means that $Y$ is a projective plane and hence any two lines in $Y$ have a common point.  For a subset $A\subseteq Y$ we denote by $\overline A$ the flat hull of the set $A$ in the liner $Y$.  
If $A$ is a line in $X$, then $\overline{A}$ is a line in $Y$ such that $\overline{A}\cap X=A$.
By the projectivity of $Y$, the lines $\overline{\w\times\{0\}}$ and $\overline{\w\times\{1\}}$ in $Y$ have a common point $h\in H$. Since $H$ is proflat in $Y$, $\overline{\w\times\{0\}}=(\w\times\{0\})\cup\{h\}$. By the projectivity of $Y$, for every $y\in\w\setminus\{0\}$, the lines $\overline{\w\times\{0\}}$ and $\overline{\w\times\{y\}}$ have a common point and this common point must be equal to $h$ (because $\overline{\w\times\{0\}}\setminus(\w\times\{0\})=\{h\})$. Then $\overline{\w\times\{y\}}=(\w\times\{y\})\cup\{h\}$. By analogy we can show that there exists a point $v\in Y$ such that $\overline{\{x\}\times\w}=(\{x\}\times\w)\cup\{v\}$ for every $x\in\w$. Assuming that $v=h$, we conclude that for the point $o\defeq(0,0)$, we have $\w\times\{0\}=\Aline oh\cap X=\Aline ov\cap X=\{0\}\times\w$, which is not true. Therefore, $h\ne v$.  

By Corollary~\ref{c:Avogadro-projective}, the $3$-long projective liner $Y$ is $2$-balanced and hence $\w$-long. Then $\Aline vh$ is an infinite line in $Y$. Since $H$ is proflat in $Y$, the set $\Aline hv\cap X$ contains at most one point and hence  there exists a point $d\in\Aline hv\setminus(\{h,v\}\cup X)$.
Choose any point $p\in X\setminus\overline{H}$ and consider the line $\Aline pd\cap X$ in $X$.
Assuming that $\Aline pd\cap X=\{x\}\times\w$ for some $x\in\w$, we conclude that $d\in\Aline pd\setminus X=\overline{\{x\}\times\w}\setminus X=\{v\}$, which contradicts the choice of the point $d$. This contradiction shows that $\Aline pd\cap X\ne\{x\}\times\w$ for every $x\in \w$. By analogy we can prove that $\Aline pd\cap X\ne\w\times\{y\}$ for every point $y\in \w$. Then $\Aline pd\cap X=F_n$ for some $n\in\w$. By Proposition~\ref{p:cov-aff}, there exists a point $q\in X\setminus (\Aline pd\cup\overline{H})$. Repeating the above argument, we can find a number $m\in\w$ such that $\Aline qd\cap X=F_m$. The choice of the point $q\notin F_n$ ensures that $F_n\ne F_m$. The transversality of the family $\{F_k\}_{k\in\w}$ ensures that $F_n\cap F_m=\{z\}$ for some point $z\in X$. Then $\Aline pd=\Aline zd=\Aline qd$, which contradicts the choice of the point $q$. This is a final contradiction showing that the para-Playfair liner $X$ has not projective completions. 
\end{proof}

\begin{exercise} For every $n\in\IN$ construct an $\w$-long para-Playfair plane $X$ with $|\overline X\setminus X|=n$.

{\em Hint:} Modify the construction from Example~\ref{ex:para-Playfair}.
\end{exercise}

In contrast to Examples~\ref{ex:Proclus-without-completion} and \ref{ex:para-Playfair}, finite proaffine regular liners do have projective completions.
 
\begin{theorem}\label{t:procompletion-finite} Every non-projective finite proaffine regular liner $X$ has a projective completion.
\end{theorem}

\begin{proof} By Theorem~\ref{t:Proclus<=>}, the proaffine regular liner $X$ is Proclus. If $X$ is not $3$-long, then Theorem~\ref{t:Proclus-not-3long} ensures that $X$ has a projective completion.  So, assume that $X$ is $3$-long. If $\|X\|\ne 3$, then the spread completion of $X$ is a projective completion of $X$, by Theorem~\ref{t:proaffine3=>compregular} and Proposition~\ref{p:spread-3long}. So, assume that $\|X\|=3$ and hence $X$ is a finite Proclus plane, by Theorem~\ref{t:Proclus<=>}. By Proposition~\ref{p:k-regular<=>2ex}, the regular liner $X$ is $3$-ranked. Taking into account that $X$ is finite and has rank $\|X\|=3$, we conclude that the $3$-ranked liner $X$ is $3$-balanced with $|X|_3=|X|$. 

If the liner $X$ is $2$-balanced, then $X$ is $p$-parallel for some $p\in\w$, by Proposition~\ref{p:23-balance=>k-parallel}. Since $X$ is Proclus, $p\le 1$. If $p=0$, then $X$ is a projective completion of $X$. If $p=1$, then $X$ is Playfair and affine.  By Corollary~\ref{c:affine-spread-completion} and Proposition~\ref{p:spread-3long}, the spread completion $\overline X$ of the  $3$-long affine regular liner $X$ is a projective completion of $X$.

So, we assume that the $3$-long Proclus plane $X$ is not $2$-balanced. Let $\mathcal L$ be the family of lines of the liner $X$ and let $$\ell\defeq\min\{|L|:L\in \mathcal L\}\ge 3.$$

\begin{claim}\label{cl:lines-short-or-long} Every line $L\in\mathcal L$ has cardinality $|L|\in\{\ell,\ell+1\}$.
\end{claim}

\begin{proof} Let $\Lambda$ be a line of cardinality $|\Lambda|=\ell$. If $L\cap\Lambda\ne \varnothing$, then $|L|\le|\Lambda|+1=\ell+1$, by  Propositions~\ref{p:Avogadro-proaffine} and \ref{p:3-long=>not=line+line}. If $L\cap\Lambda=\varnothing$, then $|L|=|\Lambda|=\ell\le \ell+1$, by Corollary~\ref{c:Proclus-par=}.
\end{proof}

By Claim~\ref{cl:lines-short-or-long}, every line in $X$ has cardinality $\ell$ or $\ell+1$. Lines of cardinality $\ell$ will be called {\em short} and lines of cardinality $\ell+1$ will be called {\em long}.

For every point $x\in X$, consider the family of lines $$\mathcal L_x\defeq\{L\in\mathcal L:x\in L\}$$ containing the point $x$.

\begin{claim}\label{cl:Lx<L+1} For every line $L\subseteq X$ and point $x\in X\setminus L$, the family $\mathcal L_x$ has $|\mathcal L_x|\in\{|L|,|L|+1\}$.
\end{claim}

\begin{proof} Since the map $L\to\mathcal L_x$, $y\mapsto\Aline xy$, is injective, the set $\mathcal L_x$ has cardinality $|\mathcal L_x|\ge|\{\Aline xy:y\in L\}|=|L|$. Since $X$ is a Proclus plane, the family $\mathcal L_x$ contains at most one line that is disjoint with $L$ and hence does not belong to the subfamily $\{\Aline xy:y\in L\}$. Then $$|L|\le|\mathcal L_x|\le|\{\Aline xy:y\in L\}|+1=|L|+1$$ and hence $|\mathcal L_x|\in\{|L|,|L|+1\}$.
\end{proof}


\begin{claim}\label{cl:Lx<ell+1} Every point $x\in X$ has  $|\mathcal L_x|\in\{\ell,\ell+1\}$.
\end{claim}

\begin{proof} If some line $L\in\mathcal L\setminus\mathcal L_x$ has cardinality $|L|=\ell$, then $|\mathcal L_x|\in\{|L|,|L|+1\}=\{\ell,\ell+1\}$, by Claim~\ref{cl:Lx<L+1}. 

So, assume that $|L|=\ell+1$ for every line $L\in\mathcal L\setminus \mathcal L_x$. Since $X$ is a plane, there exists a line that does not contain the point $x$. This line has cardinality $\ell+1$ and Claim~\ref{cl:Lx<L+1} implies that $|\mathcal L_x|\in\{\ell+1,\ell+2\}$. 

Fix any line $\Lambda\in\mathcal L$ of cardinality $|\Lambda|=\ell$. Our assumption ensures that $\Lambda\in\mathcal L_x$.  By Claim~\ref{cl:Lx<L+1}, every point $y\in X\setminus \Lambda$ has $|\mathcal L_y|\le |\Lambda|+1=\ell+1$. If some point $y\in X\setminus \Lambda$ has $|\mathcal L_y|=\ell+1$, then there exists a line $L\in \mathcal L_y$ disjoint with the line $\Lambda$. Since $x\in \Lambda$, the line $L$ does not belong to the family $\mathcal L_x$ and hence $|L|=\ell+1>|\Lambda|$, which contradicts Corollary~\ref{c:Proclus-par=}. This contradiction shows that $|\mathcal L_y|=\ell$ for every $y\in X\setminus \Lambda$. 

 Since $X$ is a plane, there exists a point $y\in X\setminus \Lambda$. Taking into account that $X\setminus\{y\}=\bigcup_{L\in\mathcal L_y}|L\setminus\{y\}|$ and $|L|=\ell+1$ for every line $L\in\mathcal L_y\setminus\{\Aline xy\}$, we conclude that $$|X|=1+|\mathcal L_y\setminus\{\Aline xy\}|\cdot \ell+|\Aline xy\setminus\{y\}|=1+(\ell-1)\ell+(|\Aline xy|-1)=\ell(\ell-1)+|\Aline xy|.$$

On the other hand, $$X\setminus\{x\}=\bigcup_{L\in\mathcal L_x}(L\setminus \{x\})=(\Aline xy\setminus\{x\})\cup\bigcup_{L\in\mathcal L_x\setminus\{\Aline xy\}}(L\setminus\{x\})$$ and hence 
$$\ell(\ell-1)+|\Aline xy|= |X|\ge 1+|\Aline xy\setminus\{x\}|+(|\mathcal L_x|-1)(\ell-1)=|\Aline xy|+(\mathcal L_x-1)(\ell-1),$$which implies $|\mathcal L_x|\le \ell+1$.
\end{proof}

Two cases are possible.
\smallskip

I. First assume that some point $x\in X$ has $|\mathcal L_x|=\ell$. 

\begin{claim}\label{cl:L=ell} Every line $L\in\mathcal L\setminus \mathcal L_x$ has cardinality $|L|=\ell$.
\end{claim}

\begin{proof} By Claim~\ref{cl:Lx<L+1}, every line $L\in\mathcal L\setminus\mathcal L_x$ has cardinality $\ell\le|L|\le |\mathcal L_x|=\ell$, which implies $|L|=\ell$.
\end{proof}

\begin{claim}\label{cl:2L=ell+1} There exists two points $y,z\in X\setminus\{x\}$ such that $\Aline xy\ne\Aline xz$ and $|\Aline xy|=\ell+1=|\Aline xz|$.
\end{claim}

\begin{proof}  Since $X$ is not $2$-balanced, there exists a line $\Lambda\in\mathcal L$ of cardinality $|\Lambda|=\ell+1$. Claim~\ref{cl:L=ell} ensures that $x\in\Lambda$.  We claim that $|\Aline xy|=\ell+1$ for some $y\in X\setminus\Lambda$. To derive a contradiction, assume that $|\Aline xy|=\ell$ for every $y\in X\setminus\Lambda$. Then $X\setminus\{x\}=\bigcup_{L\in\mathcal L_x}(L\setminus\{x\})$ implies
$$|X|=1+|\Lambda\setminus\{x\}|+(|\mathcal L_x\setminus\{\Lambda\})(\ell-1)=\ell+1+(\ell-1)(\ell-1)=\ell^2-\ell+2.$$ On the other hand, for every $y\in X\setminus\Lambda$, Claims~\ref{cl:Lx<L+1} and \ref{cl:Lx<ell+1} imply $\ell+1=|\Lambda|\le|\mathcal L_y|\le \ell+1$ and hence
$$\ell^2-\ell+2=|X|=1+\sum_{L\in\mathcal L_y}|L\setminus\{y\}|=1+(\ell+1)(\ell-1)=\ell^2,$$
which implies $\ell=2$ and contradicts our assumption. This contradiction shows that $|\Aline xy|=\ell+1$ for some $y\in X\setminus\Lambda$. Choose any point $z\in\Lambda\setminus \{x\}$ and observe that $|\Aline xz|=|\Lambda|=\ell+1$ and $\Aline xz\ne \Aline xy$.
\end{proof}

\begin{claim}\label{cl:Ly=ell+1} Every point $y\in X\setminus \{x\}$ has $|\mathcal L_y|=\ell+1$.
\end{claim}

\begin{proof} By Claim~\ref{cl:2L=ell+1}, for every point $y\in X\setminus  \{x\}$, there exists a line $\Lambda\in\mathcal L_x$ such that $|\Lambda|=\ell+1$ and $y\notin\Lambda$. Applying Claims~\ref{cl:Lx<L+1} and \ref{cl:Lx<ell+1}, we conclude that $\ell+1=|\Lambda|\le|\mathcal L_y|\le \ell+1$ and hence $|\mathcal L_y|=\ell+1$.
\end{proof}

\begin{claim}\label{cl:Lx=ell+1} Every line $L\in\mathcal L_x$ has cardinality $\ell+1$.
\end{claim}

\begin{proof} By Claim~\ref{cl:2L=ell+1}, the family $\mathcal L_x$ contains  a  line $\Lambda$ of cardinality $|\Lambda|=\ell+1$. To derive a contradiction, assume that some line $\Lambda'\in\mathcal L_x$ has cardinality $|\Lambda'|=\ell$. Choose any points $y\in \Lambda\setminus\{x\}$ and $z\in\Lambda'\setminus\{x\}$. Claim~\ref{cl:Ly=ell+1} ensures that $|\mathcal L_y|=\ell+1=|\mathcal L_z|$. Claim~\ref{cl:L=ell} implies that every line $L\in\mathcal L_z\cup(\mathcal L_y\setminus\{\Lambda\})$ has cardinality $|L|=\ell$. Then $X\setminus\{z\}=\bigcup_{L\in\mathcal L_z}(L\setminus\{z\})$ implies $$|X|=1+|\mathcal L_z|(\ell-1)=1+(\ell+1)(\ell-1)=\ell^2.$$
On the other hand, $X\setminus\{y\}=\bigcup_{L\in\mathcal L_y}(L\setminus\{y\})$ implies
$$\ell^2=|X|= 1+|\Lambda\setminus\{y\}|+(|\mathcal L_y\setminus\{\Lambda\}|)(\ell-1)=1+\ell+\ell(\ell-1)=\ell^2+1,$$
which is a desired contradiction showing that every line $\Lambda\in\mathcal L_x$ has cardinality $|\Lambda|=\ell+1$.
\end{proof}

By Claim~\ref{cl:Lx=ell+1}, every line $L\in\mathcal L_x$ has cardinality $|L|=\ell+1$. Then the equality $X\setminus\{x\}=\bigcup_{L\in\mathcal L_x}(L\setminus\{x\})$ implies $$|X|=1+|\mathcal L_x|\cdot \ell=1+\ell^2.$$
Then the subliner $A\defeq X\setminus \{x\}$ of the liner $X$ has cardinality $|A|=|X|-1=\ell^2$. Claims~\ref{cl:L=ell}, \ref{cl:Lx=ell+1}, \ref{cl:Ly=ell+1} ensure that the liner $A$ is $2$-balanced with $|A|_2=\ell$ and $1$-parallel, so Playfair. By Theorem~\ref{t:2-balance+k-parallel=>3-balance}, the Playfair liner $X$ is $3$-balanced with $|A|_3=|A|_2^2=\ell^2=|A|$. The equality $|A|_3=|A|$ implies $\|A\|=3$. By Theorem~\ref{t:Playfair<=>}, the Playfair liner $A$ is $3$-long, $3$-regular and affine. Since $\|A\|=3$, the $3$-regular liner $A$ is regular. Then $X$ is a completion of the affine regular liner $A$. Let us show that $X$ is a normal completion of $A$. Given any plane $P$ in $X$, we need to show that $P\cap A$ is a plane in $X$. The $3$-rankedness of the Proclus plane $X$ ensures that $X=P$ and hence $P\cap A=A$ is a plane in $A$. Therefore, the non-projective $3$-long $3$-ranked liner $X$ is a normal completion of the $3$-long affine regular liner $A$. By Theorem~\ref{t:procompletion<=>normcompletionA}, the liner $X$ has a projective completion.
\smallskip

II. Next, consider the second case, when $|\mathcal L_x|=\ell+1$ for all points $x\in X$.

\begin{claim}\label{cl:Lx||} For every line $L\subseteq X$ of cardinality $|L|=\ell$ and every $x\in X\setminus L$, there exists a line $L'\in\mathcal L_x$ such that $L'\cap L=\varnothing$ and $|L'|=\ell$;
\end{claim} 

\begin{proof} Since $|L|=\ell<\ell+1=|\mathcal L_x|$, there exists a line $L'\in\mathcal L_x$ such that $L'\cap L=\varnothing$. Corollary~\ref{c:Proclus-par=} ensures that $|L'|=|L|=\ell$.
\end{proof}

\begin{claim}\label{cl:|X|=ell2+ell} The liner $X$ has cardinality  $|X|=\ell^2+\ell$.
\end{claim}

\begin{proof} Since the liner $X$ is not $2$-balanced, there exist lines $L,\Lambda$ in $X$ of cardinality $|L|=\ell$ and $|\Lambda|=\ell+1$. By Claim~\ref{cl:Lx||}, for every point $x\in \Lambda\setminus L$, there exists a line $L_x\in\mathcal L_x$ such that $L_x\cap L=\varnothing$ and $|L_x|=|L|$. If $L\cap\Lambda$ contains some point $x$, then put $L_x\defeq L$. Taking into account that $X$ is a Proclus plane, we conclude that  the family $(L_x)_{x\in\Lambda}$ consists of pairwise disjoint lines. We claim that $X=\bigcup_{x\in\Lambda}L_x$. Indeed, for every point $y\in X\setminus L$, by the Claim~\ref{cl:Lx||}, there exists a line $L_y\in\mathcal L_y$ such that $L_y\cap L=\varnothing$. Assuming that $\Lambda\cap L_y=\varnothing$, we conclude that $|\mathcal L_y|\ge|\Lambda|+1=\ell+2$, which contradicts Claim~\ref{cl:Lx<ell+1}. This contradiction shows that the intersection $\Lambda\cap L_y$ contains some point $z$. Then $L_z$ and $L_y$ are two lines in the  plane $X$ that contain the point $z$ and are  disjoint with the line $L$. Since $X$ is Proclus, $L_z=L_y$ and hence $y\in L_z\subseteq\bigcup_{x\in\Lambda} L_x$. This completes the proof of the equality $X=\bigcup_{x\in\Lambda}L_x$. Therefore $\{L_x:x\in\Lambda\}$ is a disjoint cover of $X$ by lines, each of cardinality $\ell$. This implies that $|X|=|\Lambda|\cdot|L|=(\ell+1)\ell=\ell^2+\ell$. 
\end{proof}

\begin{claim}\label{cl:!Lx=ell} For every $x\in X$, there exists a unique line $L\in\mathcal L_x$ with $|L|=\ell$.
\end{claim}

\begin{proof} Claim~\ref{cl:|X|=ell2+ell} ensures that $|X|=\ell^2+\ell$ and Claim~\ref{cl:Lx||} implies that for every $x\in X$ the family $\mathcal L_x$ contains a line of cardinality $\ell$. Assuming that $\mathcal L_x$ contains two distinct lines of  cardinality $\ell$, we conclude that
$$\ell^2+\ell=|X|=1+\sum_{L\in\mathcal L_x}|L\setminus\{x\}|\le 1+2(\ell-1)+(|\mathcal L_x|-2)\ell=2\ell-1+(\ell-1)\ell=\ell^2+\ell-1,$$
which is a contradiction showing that the family $\mathcal L_x$ contains a unique line of cardinality $\ell$.
\end{proof}

\begin{claim}\label{cl:short-lines-are-disjoint} Let $L,\Lambda\in\mathcal L$ be any distint lines in the liner $X$. The lines $L,\Lambda$ are disjoint if and only if $|L|=\ell=|\Lambda|$.
\end{claim}

\begin{proof} To prove the ``if'' part, assume that $|L|=\ell=|\Lambda|$. Assuming that the lines $L,\Lambda$ have have a common point $x$, we conclude that the family $\mathcal L_x$ contains two distinct lines of cardinality $\ell$, which contradicts Claim~\ref{cl:!Lx=ell}. Therefore, $L\cap\Lambda=\varnothing$. 

To prove the ``only if'' part, assume that $L\cap\Lambda=\varnothing$. Assuming that $|L|=\ell+1$, choose any point $p\in\Lambda$ and observe that $|\mathcal L_p|\ge|\{\Lambda\}\cup\{\Aline xp:x\in L\}|=|L|+1=\ell+2$, which contradicts Claim~\ref{cl:Lx<ell+1}. This contradiction shows that $|L|=\ell$. By analogy we can prove that $|\Lambda|=\ell$.
\end{proof}

By Claims~\ref{cl:!Lx=ell} and \ref{cl:short-lines-are-disjoint}, the family $\mathcal S\defeq\{L\in\mathcal L:|L|=\ell\}$ is a spread of lines in $X$. Moreover, every spreading line in $X$ belongs to the spread $\mathcal S$. So, $X$ is para-Playfair and $X$ contains no concurrent spreading lines. By Theorem~\ref{t:spread=projective1} and Proposition~\ref{p:spread-3long}, the spread completion $\overline X$ of the Proclus plane $X$ is a  projective completion of $X$.
\end{proof}

\section{The Kuiper--Dembowski classification of finite Proclus planes}

In this section we shall prove the following \index[person]{Kuiper}Kuiper\footnote{{\bf Nicolaas Hendrik Kuiper} (1920 --1994) was a Dutch mathematician, known for Kuiper's test (in statistics) and Kuiper's theorem (on the weak contractibility of the group of bounded linear transformations of an infinite-dimensional complex Hilbert space). He also contributed to the Nash embedding theorem. Kuiper studied at University of Leiden in 1937-41, and worked as a secondary school teacher of mathematics in Dordrecht in 1942-47. He completed his Ph.D. in differential geometry from the University of Leiden in 1946. In 1947 he came to the United States at the invitation of Oscar Veblen, where he stayed at the Institute for Advanced Study for one year as Veblen's assistant. He became professor of pure mathematics at the University of Amsterdam in 1962. In 1969-70 he made a second visit at the Institute for Advanced Study. At his return from Princeton, he gave a talk at the International Congress of Mathematicians organised in Nice, during which he was appointed in the executive committee of the International Mathematical Union for 1971–1975. He finally served as director of the Institut des Hautes \'Etudes Scientifiques from 1971 until his retirement in 1985, then stayed there as a long-term visitor for six years. In 1990, he was appointed chairman of the program committee of the International Congress of Mathematicians held at Kyoto.}--Dembowski\footnote{
{\bf Heinz Peter Dembowski} (1928 -- 1971) was a German mathematician, specializing in combinatorics. He is known for the Dembowski-Wagner theorem and for Dembowski-Ostrom polynomials. The primary focus of Dembowski's research was finite geometries and their interrelations with group theory, about which he wrote an authoritative textbook. He proved the theorem, famous in finite geometry, that every inversive plane of even order $n$ is isomorphic to the system of points and plane sections of an ovoid in a three-dimensional projective space over the $n$-element field. 
In 1962 he was an approved speaker (but not an invited speaker) with half-hour talk ``Partial planes with parallelism'' at the International Congress of Mathematicians in Stockholm. His doctoral students include William Kantor (who proved an important classification theorem for $2$-transitive finite balanced liners).}\index[person]{Dembowski} classification of finite Proclus planes.

\begin{theorem}[Kuiper--Dembowski, 1962]\label{t:Kuiper-Dembowski} Every finite Proclus plane $X$ is equal to the complement $P\setminus H$ of a proflat set $H$ of rank $\|H\|\le 2$ in some projective plane $P$. The proflat set $H$ is one of the following:
\begin{itemize}
\item the empty set, in which case $X$ is a projective plane;
\item a singleton, in which case $X$ is a punctured projective plane;
\item a line, in which case $X$ is an affine plane;
\item a line with a removed point, in which case $X$ is an affine liner with attached point at infinity.
\end{itemize}
\end{theorem}

\begin{proof} Let $X$ be a finite Proclus plane. If $X$ is projective, then $X\defeq P\setminus H$ for the projective plane $P=X$ and the proflat set $H\defeq \varnothing$ of rank $\|H\|=0\le 2$.

So, let us assume that the liner $X$ is not projective. If the liner $X$ is not $3$-long, then by Theorem~\ref{t:Proclus-not-3long}, $X=P\setminus H$ for a Steiner projective plane $P$ and a set $H\subseteq P$ of cardinality $|H|\in\{1,2\}$. If $|H|=1$, then $X=P\setminus H$ is a punctured projective plane.
If $|H|=2$, then $\overline H$ is a line in $P$ and $\overline H\setminus H$ contains a unique point $o$. Observe that the subliner $A\defeq P\setminus\overline H$ is affine and $X=P\setminus A=A\cup\{o\}$ is an affine liner with attached point at infinity.

Next, consider the case of a $3$-long liner $X$. By Theorem~\ref{t:procompletion-finite}, the $3$-long finite Proclus plane $X$ has a projective completion $Y$. By Corollary~\ref{c:procompletion-rank}, $\|Y\|=\|X\|=3$ and hence $Y$ is a $3$-long projective plane and the horizon $H\defeq Y\setminus X$ has rank $\|H\|<\|Y\|=3$ (because $\overline H\ne Y$).

Since $X$ is not projective, the horizon $H$ is not empty. Then we can fix points $y\in H$ and $x\in Y\setminus\overline H\subseteq X$ and observe that the intersection $\Aline xy\cap\overline H$ contains at most one point. Since the projective plane $Y$ is $3$-long, the intersection $\Aline xy\cap X$ is a line in the liner $X$. Since the liner $X$ is $3$-long, $|\Aline xy|=|\Aline xy\cap X|+|\{y\}|\ge 3+1=4$. 
By Corollary~\ref{c:Avogadro-projective}, the $3$-long projective liner $Y$ is $2$-balanced. Then $|Y|_2=|\Aline xy|\ge 4$, which means that $Y$ is $4$-long. 

Now let us return to the analysis of the structure of the horizon $H$. If $\|H\|=1$, then $X=Y\setminus H$ is a punctured projective plane.

So, assume that $\|H\|=2$. In this case, the $3$-rankedness of the projective liner $Y$ ensures that $\overline H$ is a hyperplane in $Y$. By Theorem~\ref{t:affine<=>hyperplane} and Corollary~\ref{c:procompletion-rank}, the subliner $A=Y\setminus \overline H$ of $Y$ is an affine plane.
If $H=\overline H$, then $H$ is a line and $X=Y\setminus H=Y\setminus\overline H=A$ is an affine plane.

Finally, assume that $H$ is not flat. Then $H$ contains two distinct points $p,q\in H$ such that $\Aline pq\not\subseteq H$. It follows from $\|H\|=2$ that $\Aline pq=\overline{H}$. 
 By Theorem~\ref{t:proaffine<=>proflat}, the horizon $H$ is proflat in $Y$, which implies that $\Aline pq\setminus H=\overline H\setminus H$ contains a unique point $o$. Then $X=Y\setminus H=A\cup\{o\}$ is the affine plane $A$ with attached point $o$ at infinity. 
\end{proof}

Applying the classification of finite Proclus planes given in Theorem~\ref{t:Kuiper-Dembowski} we can now answer Question~\ref{q:L+1} about the cardinality of lines in $3$-long Proclus liners.

\begin{theorem}\label{t:card-lines-in-Proclus} For any  lines $L,\Lambda$ in a $3$-long Proclus liner $X$, we have $|L|\le1+|\Lambda|$.
\end{theorem}

\begin{proof} Let $\mathcal L$ be the family of lines in $X$ and $\ell\defeq\min\{|L|:L\in\mathcal L\}$.
Fix any line $L_0\in\mathcal L$ of cardinality $|L_0|=\ell$. We have to prove that every line in $X$ has cardinality $\ell$ or $\ell+1$. To derive a contradiction, assume that some line $L_2\in\mathcal L$ has cardinality $|L_2|>\ell+1$. By Corollary~\ref{c:L<L+2}, $|L_2|\le |L_0|+2=\ell+2$. Then $\ell+1<|L_2|\le\ell+2$ implies that the cardinal $\ell$ is finite and hence $|L_2|=\ell+2$.

Fix any points $a\in L_0$ and $b\in L_2\setminus\{a\}$ and consider the line $L_1\defeq\Aline ab$. Proposition~\ref{p:Avogadro-proaffine} ensures that $\ell+2=|L_2|\le|L_1|+1\le (|L_0|+1)+1=\ell+2$ and hence $|L_1|=\ell+1$.

By Theorem~\ref{t:Kuiper-Dembowski}, the Proclus plane $\overline{L_1\cup L_2}$ is either a punctured projective plane or an affine plane with an attached point at infinity.

If $\overline{L_1\cup L_2}$ is a punctured projective plane, then for every point $c\in L_2\setminus\{b\}$, the line $\Aline ac$ is concurrent with the ``short'' line $L_1=\Aline ab$ of the plane $\overline{L_1\cup L_2}$ and hence $|\Aline ac|=|L_1|+1=\ell+2$. On the other hand, since the lines $L_0$ and $\Aline ac$ are concurrent, we can apply Proposition~\ref{p:Avogadro-proaffine} and conclude that $\ell+2=|\Aline ac|\le |L_0|+1=\ell+1$, which is a contradiction showing that $\overline{L_1\cup L_2}$ is not a punctured projective plane and hence it is an affine plane with an attached point $p$ at infinity. This point $p$ has the property that every line in the plane $\overline{L_1\cup L_2}$ that contain $p$ has length $\ell+2=|L_2|$. In particular, the line $\Aline ap$ has lenth $|\Aline ap|=\ell+2=|L_0|+2$, which contradicts Proposition~\ref{p:Avogadro-proaffine}.

In both cases we obtain a contradiction showing that every line in $X$ has cardinality $\ell$ or $\ell+1$, and hence $|L|\le|\Lambda|+1$ for any two lines $L,\Lambda\in \mathcal L$.
\end{proof}

\section{Characterizing completely regular liners}

We recall that a liner $X$ is defined to be \index{completely regular line}\index{liner!completely regular}\defterm{completely regular} if $X$ is $3$-ranked and the spread completion $\overline X$ of $X$ is strongly regular. 
In this section we present a characterization of completely regular (finite) liners, extending the characterization given in Theorem~\ref{t:spread=projective1}. 

\begin{theorem}\label{t:spread=projective2} For a non-projective liner $X$, the following conditions are equivalent:
\begin{enumerate}
\item $X$ is completely regular;
\item $X$ is a regular para-Playfair liner possessing a projective completion;
\item $X$ is regular and has a projective completion $Y$ with flat horizon $Y\setminus X$;
\end{enumerate}
If the liner $X$ is finite, then the conditions \textup{(1)--(3)} are equivalent to
\begin{itemize}
\item[\textup{(4)}] $X$ is regular and para-Playfair.
\end{itemize}
If the liner $X$ is finite and $4$-long, then the conditions \textup{(1)--(4)} are equivalent to
\begin{itemize}
\item[\textup{(5)}] $X$ is para-Playfair.
\end{itemize}
\end{theorem}

\begin{proof} $(1)\Ra(2)$ Assume that a liner $X$ is completely regular, which means that $X$ is $3$-ranked and its spread completion $\overline X$ is strongly regular. By Theorem~\ref{t:spread=projective1}, the completely regular liner $X$ is regular and para-Playfair. By Proposition~\ref{p:spread-3long}, the strongly regular liner $\overline X$ is projective and $3$-long. So, $\overline X$ is a projective completion of $X$.
\smallskip

The implication $(2)\Ra(3)$ follows from Theorem~\ref{t:flat-horizon}.
\smallskip

$(3)\Ra(1)$ Assume that the liner $X$ regular and has a projective completion $Y$ with flat horizon $H\defeq Y\setminus X$. By Theorem~\ref{t:flat-horizon}, the liner $X$ is para-Playfair, and by Proposition~\ref{p:Playfair=>para-Playfair=>Proclus} and Theorem~\ref{t:Proclus<=>}, the para-Playfair liner $X$ is Proclus and proaffine. 

Next, we show that the liner $X$ is bi-Bolyai. Given any concurrent Bolyai lines $A,B\subseteq X$, we need to check that every line $L$ in the plane $\Pi\defeq\overline{A\cup B}\subseteq X$ is Bolyai. Observe that the flat hull $\overline \Pi$ of the plane $\Pi$ in the projective liner $Y$ is a plane in $Y$. Since the lines $A,B$ are Bolyai, there exist two lines $A',B'\subseteq \Pi$ such that $A\cap A'=\varnothing= B\cap B'$. Let $\overline L,\overline A,\overline{A'},\overline{B},\overline{B'}$ be the flat hulls of the lines $L,A,A',B,B'$ in the projective liner $Y$. By the projectivity of the liner $Y$, there exist points $a\in \overline{A}\cap\overline{A'}$ and $b\in\overline{B}\cap\overline{B'}$. It follows from $A\cap A'=\varnothing=B\cap B'$ that $a,b\in \overline\Pi\setminus X=\overline\Pi\cap H$. Since the lines $A,B$ are concurrent, so are the lines $\overline A$ and $\overline B$. Then $a\ne b$ and by the projectivity of the liner $Y$, there exists a point $c\in\overline L\cap\Aline ab\in H\cap\overline \Pi\supseteq\{a,b\}$. Since the liner $Y$ is $3$-long, for every $x\in X\setminus L$, the flat $L_x\defeq \Aline xc\cap X=\Aline xc\setminus\{c\}$ is a line in $X$ disjoint with the line $L=\overline L\setminus \{c\}$. By Corollary~\ref{c:proj-comp-Proclus}, $Y$ is a normal completion of $X$ and hence the lines $L_x$ and $L$ in $X$ are coplanar. By Corollary \ref{c:parallel-lines<=>}, the disjoint coplanar lines $L_x$ and $L$ are parallel in $X$. Therefore, $X=\bigcup L_\parallel$. By Theorem~\ref{t:spreading<=>}, the line $L$ is Bolyai. By  Theorem~\ref{t:spread=projective1}, the regular para-Playfair bi-Bolyai liner $X$ is completely regular.
\smallskip

The implication $(2)\Ra(4)$ is trivial, and $(4)\Ra(2)$ follows from Theorem~\ref{t:procompletion-finite} (to see that $X$ is proaffine, apply  Proposition~\ref{p:Playfair=>para-Playfair=>Proclus} and Theorems~\ref{t:Proclus<=>}).
\smallskip

The implication $(4)\Ra(5)$ is trivial.
\smallskip

$(5)\Ra(1)$ Assume that the liner $X$ is finite, $4$-long, and para-Playfair. We claim that the liner $X$ is bi-Bolyai. Given two concurrent Bolyai lines $A,B\subseteq X$, we need to check that every line $L$ in the plane $P\defeq\overline{A\cup B}$ is Bolyai. By Proposition~\ref{p:Playfair=>para-Playfair=>Proclus}, the para-Playfair plane $P$ is Proclus. Since $P$ contains disjoint lines, it is not projective. By Theorem~\ref{t:Proclus<=>}, the Proclus liner $P$ is proaffine and $3$-proregular. Being $3$-long, the $3$-proregular liner $P$ is $3$-regular. Since $\|P\|=3$, the $3$-regular liner $P$ is regular. By the (already proved) implication $(4)\Ra(1)$, the finite regular para-Playfair liner $P$ is completely regular, and by Theorem~\ref{t:spread=projective1}, the completely regular liner $P$ is bi-Bolyai. Then for the line $L$ in the plane $P$ there exists a line $L'\subseteq P\setminus L$. By Proposition~\ref{p:Bolyai-in-para-Playfair}, the line $L$ in the $3$-long para-Playfair liner $X$ is Bolyai, witnessing that the liner $X$ is Bolyai. By Theorem~\ref{t:4long+pP+bB=>regular}, the $4$-long para-Playfair bi-Bolyai liner $X$ is regular and by Theorem~\ref{t:spread=projective1}, the regular para-Playfair bi-Bolyai liner $X$ is completely regular.
\end{proof}

\begin{remark} Example~\ref{ex:HTS} of a non-regular Hall liner shows that in Theorem~\ref{t:spread=projective2} the condition (5) is not equivalent to the conditions (1)--(4) for $3$-long Playfair liners.
\end{remark}

\chapter{The free closures and free projectivizations of liners}

In this chapter we introduce and study the (canonical) constructions of the free closure and free projectivization of a liner. Those constructions were first introduced by Marshall Hall in his seminal paper \cite{Hall1943}, published in 1943, during the Second World War\footnote{In this paper Hall writes: ``The author has made an attempt to investigate projective planes as systematically
as possible. A considerable portion of these investigations is not
included here either because of the incomplete and unsatisfactory nature of
the results or because of lack of generality. War duties have forced the postponement
of the completion of this work.'' According to Wikipedia, Marshall Hall worked in Naval Intelligence during World War II, including six months in 1944 at Bletchley Park, the center of British wartime code breaking.}. In our presentation, the ``canonization'' of the constructions of the free closure and free projectivization are achieved with the help of the Axiom of Foundation accepted in many popular systems of axioms (like ZFC or NBG) in Foundations of Mathematics. The Axiom of Foundation says that every nonempty set $S$ contains an element $s\in S$ such that $s\cap S=\varnothing$. This axiom forbids the existence of sequences of sets $(x_n)_{n\in\w}$ with $x_{n+1}\in x_n$ for all $n\in\w$.

\section{The free closure of a liner}

Let $(X,\mathsf L)$ be a liner, $\mathcal L$ be the family of lines in $X$, and $\mathcal D\defeq\big\{\{A,B\}:A,B\in \mathcal L\;\wedge\;A\cap B=\varnothing\big\}$ be the family of pairs of disjoint lines in $X$. For any pair $\{A,B\}\in\mathcal D$, consider the ``point'' $p_{A,B}\defeq\{X,A,B\}$ and observe that $p_{A,B}\notin X$. Indeed, assuming that $p_{A,B}\in X$, we obtain the contradiction $X\in\{X,A,B\}=p_{A,B}\in X$ with the Axiom of Foundation.

The \index{free closure}\index{liner!free closure}\defterm{free closure} of the liner $(X,\mathsf L)$ is the set of points
$$\hat X\defeq X\cup\{p_{A,B}:\{A,B\}\in\mathcal D\},$$endowed with the line relation
$$
\begin{aligned}
\mathsf L_{\hat X}&\defeq \mathsf L\cup\big\{xyz\in \hat X:y\in\{x,z\}\big\}\\
&\quad\;\;\cup \bigcup_{\{A,B\}\in\mathcal D}\bigcup_{C\in \{A,B\}}\{\big((C^2\times\{p_{A,B}\})\cup(\{p_{A,B}\}\times C^2)\cup\{xyz\in C\times\{p_{A,B}\}\times C:x\ne z\}\big).
\end{aligned}
$$ 
It follows from $X\cap\{p_{A,B}:\{A,B\}\in\mathcal D\}=\varnothing$ that ${\mathsf L}_{\hat X}$ is indeed a line relation on $\hat X$ such that $\mathsf L={\mathsf L}_{\hat X}\cap X^3$ and hence $(X,\mathsf L)$ is a subliner of its free closure $(\hat X,\mathsf L_{\hat X})$. For every line $A\in\mathcal L$ its flat hull $\overline A$ in the free closure $\hat X$ coincides with the set of points $A\cup \{p_{A,B}:B\in \mathcal L\;\wedge\;A\cap B=\varnothing\}$.

It is easy to see that the free closure of a liner $X$ coincides with $X$ if and only if $X$ contains no disjoint lines if and only if $X$ is either empty, a singleton, a line or a projective plane.

\begin{definition} A subliner $X$ of a liner $Y$ is defined to be \index{preclosed set}\index{set!preclosed}\defterm{preclosed in} $Y$ if for any lines $A,B\subseteq X$ their flat hulls in $Y$ have a common point $y\in \overline A\cap\overline B\subseteq Y$.
\end{definition}

The definition of the free closure implies that every liner is preclosed in its free closure. Moreover, the free closures have the following universality property.

\begin{proposition}\label{p:line-in-2closure} Let $X$ be a liner, $Y\defeq \hat X$ be its free closure and $Z\defeq\hat Y$ be the free closure of $Y$. If $X\ne Y$, then
\begin{enumerate}
\item $Y\ne Z$;
\item no line in $X$ remains a line in $Z$.
\end{enumerate}
\end{proposition}

\begin{proof} For two distinct points $x,y\in Z$, let $\Aline xy$ be the line in the liner $Z$, containing the points $x,y$. Assume that $X\ne Y$.
\smallskip

1. Since $X\ne Y$, the liner $X$ contains disjoint lines $A,B$.  Fix any distinct points $a,a'\in A$ and $b,b'\in B$. Since $X$ is preclosed in $Y$, there exist points $p\in (\Aline a{a'}\cap\Aline b{b'}\cap Y)\setminus X$ and $y\in\Aline a{b'}\cap\Aline {a'}b\subseteq Y$. It is easy to see that $p,y\notin \Aline ab$. The definition of the liner relation in the liner $Y$ ensures that the doubleton $\{p,y\}$ is a line in $Y$, disjoint with the line $\Aline ab\cap Y$, which implies $Y\ne Z$.
\smallskip

2. To derive a contradiction, assume that some line $L\subseteq X$ remains a line in the liner $Z$. Then $L$ is a line in $Y$. The definition of the free closures $Y=\hat X$ and $Z=\hat Y$ ensures that the line $L$ intersects every line in the liners $X$ or $Y$.  Since $X\ne Y$, the liner $X$ contains two disjoint lines $A,B$. Since $L$ intersects every line in $X$, there exist unique points $a\in A\cap L$ and $b\in B\cap L$.
Choose any distinct points $a,a'\in A$ and $b,b'\in B$. Since $X$ is preclosed in $Y$, there exist points $p\in (\Aline a{a'}\cap\Aline b{b'}\cap Y)\setminus X$ and $y\in \Aline {a'}b\cap \Aline a{b'}\cap Y$. By the definition of the line relation in the free closure $Y$ of $X$, the doubleton $\{p,y\}$ is a line in $Y$ which is disjoint with the line $L$. Since the liner $Y$ is preclosed in the liner $Z$, there exists a point $q\in \Aline ap\cap \overline L=\Aline ap\cap L\subseteq L\subseteq X$. Then the line $\{x,p\}=\Aline ap\cap Y$ contains the point $q\notin\{x,p\}$, which is impossible.
\end{proof}

\begin{proposition}\label{p:extend-to-freeclosures} Let $\Phi:X\to Y$ be a liner embedding. If the set $\Phi[X]$ is preclosed in $Y$, then there exists a unique liner morphism $\hat \Phi:\hat X\to Y$ such that $\hat \Phi(x)=\Phi(x)$ for all $x\in X$.
\end{proposition} 

\begin{proof} Observe that for any pair $\{A,B\}\in \mathcal D$ of disjoint lines in the liner $X$, the images $\Phi[A]$ and $\Phi[B]$ are disjoint lines in the subliner $\Phi[X]$ of $Y$. Since $\Phi[X]$ is preclosed in $Y$, the lines $\overline{\Phi[A]}$ and $\overline {\Phi[B]}$ have a unique common point $y_{A,B}\in \overline{\Phi[A]}\cap\overline {\Phi[B]}\subseteq Y$. Extend the liner embedding $\Phi:X\to Y$ to the map $\hat\Phi:\hat X\to Y$ assigning to every point $p_{A,B}\in \hat X\setminus X$ the point $y_{A,B}$. The definition of the line relation on the free closure $\hat X$ of $X$ ensures that $\hat \Phi:\hat X\to Y$ is a liner morphism.
\end{proof}

Proposition~\ref{p:extend-to-freeclosures} implies the following corollary.

\begin{corollary}\label{c:extend-to-freeclosures} For every isomorphism $\Phi:X\to Y$ between liners, there exists a unique isomorphism $\hat \Phi:\hat X\to\hat Y$ such that $\hat\Phi{\restriction}_X=\Phi$.
\end{corollary} 

\begin{definition} A liner $X$ is defined to be \defterm{$3$-wide} if for every point $x\in X$ there exist three distinct lines $A,B,C$ in $X$ such that $x\in A\cap B\cap C$ and $\min\{|A|,|B|,|C|\}\ge 3$.
\end{definition}

\begin{exercise} Show that every $3$-long liner $X$ of rank $\|X\|\ge 3$ is $3$-wide.
\end{exercise}

\begin{proposition}\label{p:3long-in-freeclosure} Let $X$ be a liner and $\hat X$ be its free closure. Any $3$-wide subliner $S$ of $\hat X$ is contained in $X$.
\end{proposition}

\begin{proof} Assuming that $S\not\subseteq X$, fix a point $s\in S\setminus X$. Then $s=p_{A,B}$ for a unique pair $\{A,B\}\in\mathcal D$ and every line in $\hat X$ containing $s$ coincides with $\overline A$ or $\overline B$, witnessing that the subliner $S$ is not $3$-wide.
\end{proof}

\begin{corollary}\label{c:free-restriction} For every isomorphism $\Phi:\hat X\to\hat Y$ between the  free closures of $3$-wide liners $X,Y$, the  restriction $\Phi{\restriction}_X:X\to Y$ is a well-defined isomorphism between the liners $X,Y$.
\end{corollary}

\begin{proof} Since $\Phi$ is a liner isomorphism, the image $\Phi[X]$ of the $3$-wide liner $X$ is a $3$-wide subliner of $Y$. Applying Propositin~\ref{p:3long-in-freeclosure}, we conclude that $\Phi[X]\subseteq Y$. By analogy we can show that $\Phi^{-1}[Y]\subseteq X$ and hence $\Phi[X]=Y$.  Therefore, the  restriction $\Phi{\restriction}_X:X\to Y$ is a well-defined isomorphism between the liners $X,Y$.
\end{proof}

Corollary~\ref{c:extend-to-freeclosures} implies that the construction of free closure is a functor in the category $\Liners$ of liners and their isomorphisms. This functor assigns to each liner $X$ its free closure $\hat X$, and to each isomorphism $\Phi:X\to Y$ between liners the unique liner isomorphism $\hat \Phi:\hat X\to\hat Y$ such that $\hat \Phi{\restriction}_X=\Phi$. For a liner $X$, the monoid of morphisms from $X$ to $X$ in the category $\Liners$ coinsides with the automorphism group $\Aut(X)$ of $X$.

 Corollaries~\ref{c:extend-to-freeclosures} and \ref{c:free-restriction} imply the following theorem.

\begin{theorem}\label{t:Aut-freeclosure} For every $3$-wide liner $X$, the map $\Aut(X)\to\Aut(\hat X)$, $\Phi\mapsto\hat\Phi$, assigning to each automorphism $\Phi:X\to X$ its unique exension $\hat\Phi:\hat X\to\hat Y$ is a well-defined isomorphism between the automorphism groups $\Aut(X)$ and $\Aut(\hat X)$. Its inverse is the restriction operator $\Aut(\hat X)\to\Aut(X)$, $\Phi\mapsto\Phi{\restriction}_X.$
\end{theorem}

A liner $X$ is \index{rigid liner}\index{liner!rigid}\defterm{rigid} if its automorphism group $\Aut(X)$ is trivial. Theorem~\ref{t:Aut-freeclosure} implies the following corollary. 

\begin{corollary} A $3$-wide liner is rigid if and only if its free closure is a rigid liner.
\end{corollary}


\section{The free projectivization of a liner}

Given any liner $X$, consider the increasing sequence of liners $(X_n)_{n\in\w}$ where $X_0\defeq X$ and for every $n\in\w$, the liner $X_{n+1}$ is the free closure of the liner $X_n$. For every $n\in\w$, let $\mathsf L_n$ be the line relation of the liner $X_n$. The \index{free projectivization}\index{liner!free projectivization}\defterm{free projectivization} of the liner $X$ is the liner $$\widehat X\defeq\bigcup_{n\in\w}X_n$$ endowed with the line relation $\widehat{\mathsf L}\defeq\bigcup_{n\in\w}\mathsf L_n$.

\begin{theorem}\label{t:free-projectivization} Let $X$ be a liner and let $\widehat X$ be its free projectivization.  
\begin{enumerate}
\item The liner $\widehat X$ has cardinality $|\widehat X|=\max\{|X|,\w\}$.
\item If the liner $X$ contains no disjoint lines, then $\widehat X=X=\hat X$.
\item If $X$ contains disjoint lines, then its free projectivization $\widehat X$ is an $\w$-long projective plane.
\item The liner $\hat X$ has rank $\|\widehat X\|=\min\{\|X\|,3\}$.
\item Every finite $3$-wide subliner of $\widehat X$ is contained in $X$.
\end{enumerate}
\end{theorem}

\begin{proof} Let $(X_n)_{n\in\w}$ be the increasing sequence of liners such that $X_0=X$ and for every $n\in\w$, the liner $X_{n+1}$ is the free closure of the liner $X_n$. For every $n\in\w$, let $\mathcal L_n$ be the family of lines in the liner $X_n$, and $\mathcal D_n\defeq\big\{\{A,B\}:A,B\in\mathcal L_n\;\wedge\;A\cap B=\varnothing\big\}$ be the family of  pairs of disjoint lines in $X_n$. Since every line is uniquely determined by two distinct points on the line, $|\mathcal L_n|\le |X_n\times X_n|$ and hence $|\mathcal D_n|\le|\mathcal L_n\times \mathcal L_n|\le|X_n|^4$. 
\smallskip

1. By induction, we shall prove that $|X_n|\le\max\{|X|,\w\}$ for every $n\in\w$. For $n=0$, this follows from the equality $X_0=X$. Assume that for some $n\in\w$ we have proved that $|X_n|\le\max\{|X|,\w\}$. The definition of the free closure ensures that $$|X_{n+1}|=|X_n\cup\{p_{A,B}:\{A,B\}\in\mathcal D_n\}|=|X_n|+|\mathcal D_n|\le |X_n|+|X_n|^4\le\max\{|X_n|,\w\}=\max\{|X|,\w\}.$$Then $|\widehat X|=\big|\bigcup_{n\in\w}X_n\big|\le\sum_{n\in\w}|X_n|\le\max\{|X|,\w\}
$.
\smallskip

2. If $X$ contains no disjoint lines, then $\mathcal D_0=\varnothing$ and the free closure $X_1\defeq X\cup\{\{X,A,B\}:\{A,B\}\in \mathcal D_0\}=X$ of $X_0=X$ coincides with $X$ and also contains no disjoint lines. Proceeding by induction, we conclude that $X_n=X$ for all $n\in\w$ and hence $\widehat X=X$.
\smallskip

3. Assume that the liner $X$ contains two disjoint lines and hence $\|X\|\ge 3$.  Then $X=X_0\ne X_1$. Applying Proposition~\ref{p:line-in-2closure}(1), we conclude that $X_n\ne X_{n+1}$ for every $n\in\w$. Since $X$ is a subliner of its free projectivization $\widehat X$, the flat hulls of the disjoint lines in $X$ are distinct lines in $\widehat X$, which implies $\|\widehat X\|\ge 3$. To see that $\widehat X$ is a projective plane, it suffices to show that any two lines in $X$ have a common point. To derive a contradiction, assume that the liner $\widehat X$ contains two disjoint lines $L$ and $\Lambda$. Since $L\cup\Lambda\subseteq \widehat X=\bigcup_{n\in\w}X_n$, there exists $n\in\w$ such that $|L\cap X_n|\ge 2$ and $|\Lambda\cap X_n|\ge 2$ and hence $A\defeq L\cap X_n$ and $B\defeq \Lambda\cap X_n$ are lines in the liner $X_n$. Since $X_{n+1}$ is the free closure of the liner $X_n$, the point $p_{A,B}\defeq\{X,A,B\}$ is the common point of the flat hulls $\overline A=L\cap X_{n+1}$ and $\overline B=\Lambda\cap X_{n+1}$ of the lines $A,B$ in the liner $X_{n+1}$. Then $$p_{A,B}\in \overline A\cap\overline B=(L\cap X_{n+1})\cap(\Lambda\cap X_{n+1})\subseteq L\cap \Lambda=\varnothing,$$
which is a contradiction showing that the liner $\widehat X$ contains no disjoint lines and hence $\widehat X$ is projective and has rank $\|\widehat X\|=3$, so $\widehat X$ is a projective plane.

It remains to show that the projective plane $\widehat X$ is $\w$-long. In the opposite case, $\widehat X$ contains a finite line $L$. Let $\ell\in\w$ be the smallest number such that $L\subseteq X_\ell$. Since $X_{\ell}\ne X_{\ell+1}$, we can apply Proposition~\ref{p:line-in-2closure}(2) and conclude that no line in the liner $X_\ell$ remains a line in the liner $X_{\ell+2}$. In particular, the line $L=L\cap X_{\ell}=L\cap X_{\ell+2}$ is not a line in $X_{\ell+2}$, which is a desired contradiction showing that the projecive plane $\widehat X$ is $\w$-long.
\smallskip

4. If $X$ contains no disjoint lines, then $\|X\|\le 3$ and  $\widehat X=X$, by a preceding statement. In this case, $\|\widehat X\|=\|X\|=\min\{\|X\|,3\}$. If $X$ contains disjoint lines, then $\|X\|\ge 3$ and $\widehat X$ is a projective plane, by the preceding item. In this case $\|\widehat X\|=3=\min\{\|X\|,3\}$.
\smallskip

5. Assume that $S$ is a finite $3$-wide subliner of $\widehat X$. Let $n\in\w$ be the smallest number such that $S\subseteq X_n$. Proposition~\ref{p:3long-in-freeclosure} implies that $n=0$ and hence $S\subseteq X_0=X$.
\end{proof}

\begin{theorem}\label{t:freeprojectivization-on-morphisms} For any isomorphism $\Phi:X\to Y$ between liners $X,Y$, there exists a unique isomorphism $\widehat\Phi:\widehat X\to\widehat Y$ such that $\widehat \Phi{\restriction}_X=\Phi$. 
\end{theorem}

\begin{proof} Let $X_0\defeq X$, $Y_0\defeq Y$  and for every $n\in\w$, let $X_{n+1}$ and $Y_{n+1}$ be the free closures of the liners $X_n$ and $Y_n$, respectively. Let $\Phi_0=\Phi$ and for every $n\in\IN$, let $\Phi_{n+1}:X_{n+1}\to Y_{n+1}$ be the unique isomorphism such that $\Phi_{n+1}{\restriction}{X_n}=\Phi_n$. The isomorphism $\Phi_{n+1}$ is well-defined, according to Corollary~\ref{c:extend-to-freeclosures}.



Then $\widehat\Phi\defeq\bigcup_{n\in\w}\Phi_n$ is an isomorphism of the liners $\widehat X=\bigcup_{n\in\w}X_n$ and $\widehat Y=\bigcup_{n\in\w}Y_n$ extending the isomorphism $\Phi$. Let us show that the isomorphism $\widehat \Phi$ is unique. Assume that $\Psi:\widehat X\to\widehat Y$ is another liner isomorphism such that $\Psi{\restriction}_X=\Phi$. Assuming that $\widehat\Psi\ne\Phi$, find the smallest number $n\in\w$ such that the liner $X_{n+1}$ contains a point $x\in X_{n+1}$ such that $\widehat\Phi(x)\ne\Psi(x)$. Then $\Psi{\restriction}_{X_n}=\widehat\Phi{\restriction}_{X_n}=\Phi_n$. Since the set $\Psi[X_n]=\Phi_n[X_n]=Y_n$ is preclosed in $Y_{n+1}\subseteq\widehat Y$, we can apply Proposition~\ref{p:extend-to-freeclosures} and conclude that $\Psi{\restriction}_{X_{n+1}}=\Phi_{n+1}$ and hence $\Psi(x)=\Phi_{n+1}(x)=\widehat\Phi(x)$, which contradicts the choice of $x$. This contradiction shows that the isomorphism $\widehat \Phi$ equals $\Psi$ and hence is unique.   
\end{proof}

\begin{remark} The construction of free projectivization is a functor in the category of liners and their isomorphisms. This functor assigns to each liner $X$ is free projectivization $\widehat X$ and to each liner isomorpism $\Phi:X\to Y$  the unique isomorphism $\widehat \Phi:\widehat X\to\widehat Y$ such that $\widehat \Phi{\restriction}_X=\Phi$. The isomorphism $\widehat \Phi$ is well-defined, by Theorem~\ref{t:freeprojectivization-on-morphisms}.
\end{remark}

\begin{theorem}\label{t:3wide-freeprojectivization} Let $X,Y$ be liners and $\Phi:\widehat X\to\widehat Y$ be an isomorphism between their free projectivizations. If the liners $X,Y$ are finite and $3$-wide, then $\Phi[X]=Y$ and hence the restriction $\Psi=\Phi{\restriction}_X$ is an isomorphism of the liners $X,Y$. Moreover, $\Phi=\widehat\Psi$.
\end{theorem}

\begin{proof}  Since the liner $X$ is finite and $3$-long, its isomorphic image $\Phi[X]$ is a finite $3$-wide subliner of the liner $\widehat Y$. Theorem~\ref{t:free-projectivization}(5) ensures that $\Phi[X]\subseteq Y$. By analogy we can show that $\Phi^{-1}[Y]\subseteq X$. Then the restriction $\Psi\defeq \Phi{\restriction}_X$ is an isomorphism of the liners $X,Y$. By Theorem~\ref{t:freeprojectivization-on-morphisms}, $\Phi=\widehat\Psi$.
\end{proof}

Theorems~\ref{t:freeprojectivization-on-morphisms} and \ref{t:3wide-freeprojectivization} imply the following corollary.

\begin{corollary} Two finite $3$-wide liners are isomorphic if and only if their free projectivizations are isomorphic.
\end{corollary}

Corollaries~\ref{c:extend-to-freeclosures} and \ref{c:free-restriction} imply the following theorem.

\begin{theorem}\label{t:Aut-freeprojectivization}For every finite $3$-wide liner $X$, the map $\Aut(X)\to\Aut(\widehat X)$, $\Phi\mapsto\widehat\Phi$, assigning to each automorphism $\Phi:X\to X$ its unique exension $\widehat\Phi:\widehat X\to\widehat Y$ is a well-defined isomorphism between the automorphism groups $\Aut(X)$ and $\Aut(\widehat X)$. Its inverse is the restriction operator $\Aut(\widehat X)\to\Aut(X)$, $\Phi\mapsto\Phi{\restriction}_X.$
\end{theorem}

Theorems~\ref{t:Aut-freeclosure} and \ref{t:Aut-freeprojectivization} imply the following corollary.

\begin{theorem}\label{t:rigid-3wide} For a finite $3$-wide liner $X$, the following conditions are equivalent:
\begin{enumerate}
\item the liner $X$ is rigid;
\item the free closure $\hat X$ of $X$ is rigid;
\item the free projectivization $\widehat X$ of $X$ is rigid.
\end{enumerate}
\end{theorem} 

\begin{examples}[Ivan Hetman]\label{ex:rigid25} There exists at least 28 rigid affine planes of order $25$.
\end{examples}

\begin{proof} The list of known projective planes of order 25 can be found at the site of Eric \index[person]{Moorhouse}Moorhouse ({\tt https://ericmoorhouse.org/pub/planes25/}). This list contains  28 projective planes $\Pi$ of order 25 that contain a line $L\subseteq \Pi$ whose orbit under the action of the automorphism group has length 500, which is equal to the cardinality of the automorphism group of the plane $\Pi$. Then the affine plane $\Pi\setminus L$ is rigid (because each non-trivial automorphism of the projective plane $\Pi$ moves the line $L$).
\end{proof}

By Example~\ref{ex:rigid25}, there exists a $3$-wide rigid liner of cardinality $625$.

\begin{examples}[Ivan Hetman]\label{ex:rigid15} By \cite[II.1.28]{HCD}, there exist exactly 80 nonisomorphic Steiner liners $X$ of cardinality $|X|=15$. All of them are $3$-wide (because every point of a Steiner liner of cardinality $15$ is contained in $7$ lines). By \cite[II.1.29]{HCD}, exactly 36 of those 80 Steiner liners are rigid.
\end{examples}

\begin{problem} What is the smallest cardinality of a $3$-wide rigid liner? 
\end{problem}

\begin{theorems}\label{t:rigid-projplane} There exists a rigid countable $\w$-long projective plane.
\end{theorems}

\begin{proof} By Example~\ref{ex:rigid25} or \ref{ex:rigid15}, there exists a rigid finite $3$-wide liner $X$. By Theorems~\ref{t:rigid-3wide} and \ref{t:free-projectivization}, its free projectivization $\widehat X$ is a rigid countable $\w$-long projective plane.
\end{proof}
  
Theorem~\ref{t:rigid-projplane} motivates the following well-known open problem. 
  
\begin{problem} Is there a rigid finite projective plane?
\end{problem}  

\begin{corollary}\label{c:proj-embed} Every liner $X$ is a subliner of an $\w$-long projective plane $Y$ such that every finite $3$-wide subliner $S\subseteq Y$  is contained in $X$.
\end{corollary}

\begin{proof} If the liner $X$ contains two disjoint lines, then its free projectivization $Y\defeq \widehat X$ has the required properties, by Theorem~\ref{t:free-projectivization}.

Now assume that the liner $X$ does no contain disjoint lines. Choose any set $Z$ containing the set $X$ so that $|Z\setminus X|=4$ and endow $Z$ with the family of lines $$\mathcal L_Z\defeq \mathcal L_X\cup\big\{\{x,y\}:(x,y)\in (Z\times Z)\setminus(X\times X)\big\},$$
where $\mathcal L_X$ is the family of lines of the liner $X$. It is easy to see that the liner $Z$ contains two disjoint lines $A,B\subseteq Z\setminus X$ (each of cardinality $2$). By Theorem~\ref{t:free-projectivization}, the free projectivization $ \widehat Z$ of the liner $Z$ is an $\w$-long projective plane such that every finite $3$-wide subliner $S\subseteq \widehat Z$ is contained in the liner $Z$. Assuming that $S\not\subseteq X$, we can find a point $s\in S\setminus X$. Since $S$ is $3$-wide, there exists a point $z\in S\setminus\{s\}$ such that the line $\Aline sz\cap S$ contains at least $3$-point. On the other hand, the definition of the liner structure on the liner $Z$ ensures that every line in $Z$ that contains the point $s\in Z\setminus X$ is a doubleton. This contradiction shows that $S\subseteq X$ and hence the liner $Y\defeq \widehat Z$ has the required property.
\end{proof}  

Corollary~\ref{c:proj-embed} suggests the following known open problem (probably due to Marshall Hall).

\begin{problem} Is every finite liner a subliner of a finite projective plane?
\end{problem}

\chapter{Desarguesian liners}

In this section we introduce and study Desarguesian
liners, whose definition is suggested by the \index[person]{Desargues}Desargues\footnote{{\bf Girard Desargues} (1591--1661), a French mathematician and engineer, who is considered one of the founders of projective geometry.}
Theorems for affine, proaffine and projective liners.

\section{Desargues Theorem for affine liners}

In this section we prove the Desargues Theorem for affine  liners. 

First we estabish some terminology.

Lines $L_1,\dots,L_n$ in a liner are called
\begin{itemize}
\item \index{concurrent lines}\index{lines!concurrent}\defterm{concurrent} if $\bigcap_{i=1}^nL_i$ is a singleton;
\item \index{parallel lines}\index{lines!parallel}\defterm{parallel} if for any $i,j\in\{1,\dots,n\}$ the lines $L_i$ and $L_j$ are parallel;
\item \index{paraconcurrent lines}\index{lines!paraconcurrent}\defterm{paraconcurrect} if the lines $L_1,\dots,L_n$ are either parallel or concurrent.
\end{itemize}

\begin{theorem}[Desargues]\label{t:Desargues-affine} Let $A,B,C$ be paraconcurrent lines in an affine regular liner $X$ of rank $\|X\|\ge 4$, and $a,a'\in A\setminus (B\cup C)$, $b,b'\in B\setminus(A\cup C)$, $c,c'\in C\setminus(A\cup B)$ be any points.
If $\Aline ab\parallel \Aline { a'}{b'}$ and $\Aline bc\parallel \Aline {b'}{c'}$, then $\Aline ac\parallel\Aline {a'}{c'}$.
\end{theorem}

\begin{proof} The choice of the points $a,b,c$ ensures that the lines $A,B,C$ are pairwise distinct. By Theorem~\ref{t:affine=>Avogadro}, the affine liner $X$ is $2$-balanced, so all lines in $X$ have the same cardinality. If $X$ is not $3$-long, then every line in $X$ contains exactly two points. By Proposition~\ref{p:parallel=>3-long}, $X$ contains no disjoint parallel line. In this case the lines $A,B,C$ are concurrent and $A\setminus (B\cup C)$, $B\setminus(A\cup C)$ and $C\setminus(A\cup B)$ are singletons. Then $a=a'$ and $c=c'$ and hence $\Aline ac=\Aline {a'}{c'}$ and $\Aline ac\parallel \Aline{a'}{c'}$. 

So, we assume that the liner $X$ is $3$-long. 
Since the lines $A,B,C$ are paraconcurrent, $\|A\cup B\|=3=\|B\cup C\|$. 

\begin{claim}\label{cl:Desargues:Q} There exists a line $D\subseteq X$ such that the lines $A,B,C,D$ are paraconcurrent and $D\cap \overline{A\cup C}=A\cap B\cap C$.
\end{claim}

\begin{proof} Since the lines $A,C$ are paraconcurrent, $\|A\cup C\|\le 3<4\le\|X\|$ and hence there exists a point $d\in X\setminus \overline{A\cup C}$. If the lines $A,B,C$ are concurrent, then the intersection $A\cap B\cap C$ contains some point $o$. In this case we put $D\defeq\Aline od$. It is clear that the lines $A,B,C,D$ have common point $o$ and hence are concurrent. Assuming that $D\cap\overline{A\cup C}\ne A\cap B\cap C$, we can find a point $x\in D\cap\overline{A\cup B}\setminus (A\cap B\cap C)\subseteq Q\cap\overline{A\cup C}\setminus\{o\}$ and conclude that $d\in D=\Aline ox\subseteq\overline{A\cup C}$, which contradicts the choice of the point $d$.

Next, assume that the lines $A,B,C$ are parallel. By Theorem~\ref{t:Playfair}, there exists a line $D$ such that $d\in D$ and $D\parallel B$. By Theorem~\ref{t:Proclus-lines}, $D\parallel A$ and $D\parallel C$. So, the lines $A,B,C,D$ are parallel. Assuming that $D\cap\overline{A\cup C}\ne\varnothing$, we can find a point $x\in D\cap\overline{A\cup C}$. It follows from $D\parallel A$ that $d\in D\subseteq\overline{\{x\}\cup A}\subseteq\overline{A\cup C}$, which contradicts the choice of the point $d\notin\overline{A\cup C}$. This contradiction shows that $D\cap\overline{A\cup C}=\varnothing$. 

In both case, we have found a line $D$ in $X$ such that the lines $A,B,C,D$ are paraconcurrent and $D\cap\overline{A\cup C}=A\cap B\cap C$.
\end{proof}

\begin{claim}\label{cl:Desargues:q} There exist a line $D\subseteq X$ and points $d,d'\in Q\setminus \overline{A\cup C}$ such that $\Aline ad\parallel \Aline {a'}{d'}$, $\Aline cd\parallel \Aline {c'}{d'}$ and the lines $A,B,C,D$ are paraconcurrent.
\end{claim}

\begin{proof} Two cases are possible.
\smallskip

1. If $B\not\subseteq \overline{A\cup C}$, then put $D\defeq B$, $d\defeq b$ and $d'\defeq b'$. It is clear that the lines $A,B,C,D$ are paraconcurrent, $\Aline ad=\Aline ab\parallel\Aline{a'}{b'}=\Aline {a'}{d'}$ and   $\Aline cd=\Aline cb\parallel\Aline{c'}{b'}=\Aline {c'}{d'}$. It remains to prove that $d,d'\notin \overline{A\cup C}$.
\smallskip

Two subcase are possible. 
\smallskip

1.1. If the lines $A,B,C$ are concurrent, then the set $A\cap B\cap C$ contains some point $o$. The choice of the points $b,b'\in B\setminus(A\cap C)$ ensure that $b\ne o\ne b'$. Assuming that $d\in\overline{A\cup C}$, we conclude that $B=\Aline ob=\Aline od\subseteq\overline{A\cup C}$, which contradicts our assumption. By analogy we can show that $d'\notin\overline{A\cup C}$.
\smallskip

1.2. If the lines $A,B,C$ are parallel, then $B\cap\overline{A\cup C}=\varnothing$. Indeed, otherwise, we can find a point $x\in B\cap\overline{A\cup C}$ and use the parallelity $B\parallel A$ to conclude that $B\subseteq\overline{\{x\}\cup A}\subseteq\overline{A\cup C}$, which contradicts our assumption. This contradiction shows that $B\cap\overline{A\cup C}=\varnothing$ and hence $d,d'\notin \overline{A\cup C}$.
\smallskip 

2. Next, assume that $B\subseteq\overline{A\cup C}$. By Claim~\ref{cl:Desargues:Q}, there exists a line $D\subseteq X$ such that the lines $A,B,C,D$ are paraconcurrent and $D\cap\overline{A\cup C}=A\cap B\cap C$. Choose any point $d\in D\setminus (A\cap B\cap C)$. Such point exists because $|D|\ge 2>1\ge |A\cap B\cap C|$. 

By Corollary~\ref{c:Playfair}, there exists a line $L$ such that $b'\in L$ and $L\parallel \Aline bd$. Then $L\subseteq\overline{\{b'\}\cup \Aline bd}=\overline{\{b',b,d\}}\subseteq\overline{B\cup\{d\}}$.  Assuming that $L\cap D=\varnothing$, and taking into account that $\|L\cup D\|\le\|\overline{\{d\}\cup B}\|\le 3$, we can apply Theorem~\ref{t:parallel-char} and conclude that $L\parallel D$. By Theorem~\ref{t:Proclus-lines}, $\Aline bq\parallel L$ and $L\parallel D$ imply $\Aline bd\parallel D$ and hence $D=\Aline bd$, which contradicts the paraconcurrence of the lines $A,B,C,D$. This contradiction shows that $L\cap D$ contains some point $d'$. Assuming that $d'\in\overline{A\cup C}$, we conclude that $d'\in Dd\cap\overline{A\cup C}= A\cap B\cap C$ and hence $\Aline{b'}{d'}=B$, which is not true as the line $B$ is not parallel to the line $\Aline bd$.
\smallskip

It remains to prove that $\Aline ad\parallel \Aline {a'}{d'}$ and $\Aline cd\parallel \Aline{c'}{d'}$. First we check that $\Aline ad\parallel \Aline {a'}{d'}$.

Since $\Aline ab\parallel \Aline{a'}{b'}$ and $\Aline bd\parallel \Aline{b'}{d'}$, Theorem~\ref{t:subparallel-via-base} ensures that the flats $\overline{\{a,b,d\}}$ and $\overline{\{a',b',d'\}}$ are parallel. Consider the flats
 $I\defeq\overline{A\cup D}\cap \overline{\{a,b,d\}}$ and $I'\defeq\overline{A\cup D}\cap \overline{\{a',b',d'\}}$ and observe that $\overline{\{a,d\}}\subseteq I\subseteq \overline{\{a,b,d\}})\setminus\{b\}$ and $\overline{\{a',d'\}}\subseteq I'\subseteq \overline{\{a',b',d'\}}\setminus\{d'\}$. The $2$-rankedness of the regular liner $X$ implies $2=\|\{a,d\}\|\le \|I\|<\|\{a,b,d\}\|\le 3$ and $2=\|\{a',d'\}\|\le \|I'\|<\|\{a',b',d'\}\|\le 3$. By the $2$-rankedness,  $I=\overline{\{a,d\}}=\Aline ad$ and $I'=\overline{\{a',d'\}}=\Aline {a'}{d'}$. By Corollary~\ref{c:paraintersect}, the lines $\Aline ad=I= \overline{A\cup D}\cap\overline{\{a,b,d\}}$ and  $\Aline {a'}{d'}=I'= \overline{A\cup D}\cap\overline{\{a',b',d'\}}$ are parallel. 
 
By analogy we can show that the lines  $\Aline cd$ and $\Aline{c'}{d'}$ are parallel.
\end{proof}

\begin{picture}(300,180)(-100,-55)

\put(0,0){\line(1,0){180}}
\put(185,-4){$B$}
\put(0,0){\line(4,1){200}}
\put(205,48){$C$}
\put(0,0){\line(4,-1){160}}
\put(165,-45){$A$}
\put(0,0){\line(1,1){110}}
\put(113,113){$D$}

\put(30,0){\line(2,1){30}}
{\linethickness{1pt}
\put(30,0){\color{violet}\line(1,-1){10}}
\put(30,0){\color{green}\line(0,1){30}}
\put(40,-10){\color{blue}\line(-1,4){10}}
\put(40,-10){\color{red}\line(4,5){20}}
\put(30,30){\color{cyan}\line(2,-1){30}}

\put(120,-30){\color{blue}\line(-1,4){30}}
\put(120,-30){\color{red}\line(4,5){60}}
\put(90,90){\color{cyan}\line(2,-1){90}}
\put(90,0){\color{green}\line(0,1){90}}
\put(90,0){\color{violet}\line(1,-1){30}}
}
\put(80,50){$L$}
\put(90,0){\line(2,1){90}}

\put(0,0){\circle*{3}}
\put(-8,-3){$o$}
\put(30,0){\circle*{3}}
\put(23,-3){$b$}
\put(30,30){\circle*{3}}
\put(25,33){$d$}
\put(40,-10){\circle*{3}}
\put(36,-18){$a$}
\put(60,15){\circle*{3}}
\put(59,18){$c$}

\put(90,0){\circle*{3}}
\put(85,-10){$b'$}
\put(90,90){\circle*{3}}
\put(85,93){$d'$}
\put(120,-30){\circle*{3}}
\put(116,-40){$a'$}
\put(180,45){\circle*{3}}
\put(179,48){$c'$}
\end{picture}

Let $D$ be the line and $d,d'\in D$ the points given by Claim~\ref{cl:Desargues:q}.
Theorem~\ref{t:subparallel-via-base} implies that the flats $\overline{\{a,c,d\}}$ and $\overline{\{a',c',d'\}}$ are parallel. Consider the flats
 $J\defeq\overline{A\cup C}\cap \overline{\{a,c,d\}}$ and $J'\defeq\overline{A\cup C}\cap \overline{\{a',c',d'\}}$ and observe that $\Aline ac\subseteq J\subseteq \overline{\{a,c,d\}}\setminus\{d\}$ and $\overline{\{a',c'\}}\subseteq J'\subseteq \overline{\{a',c',d'\}}\setminus\{d'\}$. The rankedness of the liner $X$ ensures that $2=\|\{a,c\}\|\le \|J\|<\|\{a,c,d\}\|\le 3$ and $2=\|\{a',c'\}\|\le \|J'\|<\|\{a',c',d'\}\|\le 3$ and hence $J=\overline{\{a,c\}}=\Aline ac$ and $J'=\overline{\{a',c'\}}=\Aline {a'}{c'}$. By Corollary~\ref{c:paraintersect}, the lines $\Aline ac=J= \overline{A\cup C}\cap\overline{\{a,c,d\}}$ and  $\Aline {a'}{c'}=J'= \overline{A\cup C}\cap\overline{\{a',c',d'\}}$ are parallel. 
\end{proof}

\section{Desargues Theorem for Proclus liners}

\begin{theorem}\label{t:Des-Proclus} Let $A,B,C$ be three paraconcurrent lines in a Proclus weakly modular liner $X$ and $a,a'\in A\setminus(B\cup C)$, $b,b'\in B\setminus(A\cup C)$, $c,c'\in C\setminus(A\cup B)$ be distinct points. If $\overline{\{a,b,c\}}$ and $\overline{\{a',b',c'\}}$ are two distinct planes, then the set $T\defeq(\Aline ab\cap\Aline{a'}{b'})\cup(\Aline bc\cap\Aline {b'}{c'})\cup(\Aline ac\cap\Aline{a'}{c'})$ has rank $\|T\|=\{0,2\}$.
\end{theorem}

\begin{proof} By Theorems~\ref{t:w-modular<=>}, the weakly modular liner $X$ is ranked. By the assumption, $\overline{\{a,b,c\}}$ and $\overline{\{a',b',c'\}}$ are two distinct planes in $X$ and hence their intersection $L\defeq \overline{\{a,b,c\}}\cap\overline{\{a',b',c'\}}$ has rank $\|L\|<3$, by the rankedness of $X$. Since $T\subseteq L$, the set $T$ has rank $\|T\|\le\|L\|\le 2$. It remains to check that $\|T\|\ne 1$.
This inequality will follow as soon as we show that $T\ne\varnothing$ implies $|T|\ge 2$.

So, assume that $T\ne\varnothing$. Then the flat $L$ is not empty.

\begin{claim}\label{cl:Des-Proclus} $\{a,b,c\}\cap\overline{\{a',b',c'\}}=\varnothing=\{a',b',c'\}\cap\overline{\{a,b,c\}}$.
\end{claim}

\begin{proof} Assuming that $a\in \overline{\{a',b',c'\}}$, we conclude that $A=\Aline a{a'}\subseteq\overline{\{a',b',c'\}}$. Since the lines $A,B,C$ are paraconcurrent, $B\cup C\subseteq \overline{A\cup\{b',c'\}}\subseteq \overline{\{a',b',c'\}}$ and hence $\overline{\{a,b,c\}}\subseteq\overline{A\cup B\cup C}\subseteq\overline{\{a',b',c'\}}$, which contradicts $\overline{\{a,b,c\}}\ne\overline{\{a',b',c'\}}$ and the rankedness of $X$. This contradiction shows that $a\notin\overline{\{a',b',c'\}}$. By analogy we can prove that $b,c\notin \overline{\{a',b',c'\}}$ and $a',b',c'\notin \overline{\{a,b,c\}}$.
\end{proof}

The paraconcurrence of the lines $A,B,C$ and Claim~\ref{cl:Des-Proclus} imply $\|\overline{\{a,b,c,a',b',c'\}}\|=4$. Applying Theorem~\ref{t:w-modular<=>}, we conclude that the flat $L$ is a line (being the non-empty intersection of two planes whose union has rank $4$ in the weakly modular liner $X$). Since $X$ is Proclus, at most one side of the triangle $abc$ is parallel to the line $L$ in the plane $\overline{\{a,b,c\}}$. So, we lose no generality assuming that the lines $\Aline ab$ and $\Aline bc$ are not parallel to the line $L$. Then there exist unique points $x\in L\cap \Aline ab$ and $y\in L\cap \Aline bc$. 

Assuming that $x=y$, we conclude that $\{x\}=\{y\}=(L\cap\Aline ab)\cap (L\cap\Aline bc)=L\cap \{b\}$ and hence $b\in L\subseteq \overline{\{a',b',c'\}}$, which contradicts Claim~\ref{cl:Des-Proclus}. This contradiction shows that $x\ne y$. 

The paraconcurence of the lines $A,B$, the $3$-rankedness of $X$, and Claim~\ref{cl:Des-Proclus} imply that $\overline{\{a,a',b,b'\}}\cap\overline{\{a,b,c\}}=\Aline ab$ and $\overline{\{a,a',b,b'\}}\cap\overline{\{a',b',c'\}}=\Aline {a'}{b'}$. Then  
$$x\in\Aline ab\cap L=(\overline{\{a,a',b,b'\}}\cap\overline{\{a,b,c\}})\cap(\overline{\{a,b,c\}}\cap\overline{\{a',b',c'\}})\subseteq\overline{\{a,a',b,b'\}}\cap\overline{\{a',b',c'\}}=\Aline {a'}{b'}$$ and hence $x\in\Aline ab\cap\Aline {a'}{b'}\subseteq T$. By analogy, we can show that $y\in T$ and hence $|T|\ge|\{x,y\}|=2$.
\end{proof}

\begin{corollary}\label{c:Des-proaffine} Let $A,B,C$ be three paraconcurrent lines in a proaffine regular liner $X$ and $a,a'\in A\setminus(B\cup C)$, $b,b'\in B\setminus(A\cup C)$, $c,c'\in C\setminus(A\cup B)$ be distinct points. If $\overline{\{a,b,c\}}$ and $\overline{\{a',b',c'\}}$ are two distinct planes, then the set $T\defeq(\Aline ab\cap\Aline{a'}{b'})\cup(\Aline bc\cap\Aline {b'}{c'})\cup(\Aline ac\cap\Aline{a'}{c'})$ has rank $\|T\|=\{0,2\}$.
\end{corollary}

\begin{proof} By Theorems~\ref{t:Proclus<=>}, the proaffine regular liner $X$ is Proclus, and by Corollary~\ref{c:proregular=>ranked} and Theorem~\ref{t:w-modular<=>}, the regular liner $X$ is ranked and weakly modular. So, we can appy Theorem~\ref{t:Des-Proclus}.
\end{proof}

\begin{exercise} Find an example of a non-Desarguesian projective liner $X$ of rank $\|X\|=4$.
\smallskip

{\em Hint:} Consider any non-Desarguesian projective plane with an attachd point.
\end{exercise}

\section{Desargues Theorem for projective liners}

In this section we prove the Desargues Theorem for projective liners, which is a projective counterpart of the Desargues Theorem~\ref{t:Desargues-affine}  for affine liners. It will be convenient to formulate this theorem in terms of centrally perspective triangles.

\begin{definition} A triple $abc\in X^3$ in a liner is called a \index{triangle}\defterm{triangle} if $\|\{a,b,c\}\|=3$.
\end{definition}

\begin{definition}\label{d:perspective-trianges} Two triangles $abc$ and $a'b'c'$ in a liner $X$ are called
\begin{itemize}
\item \index{disjoint triangles}\index{triangles!disjoint}\defterm{disjoint} if $\{a,b,c\}\cap\{a',b',c'\}=\varnothing$.
\item  \index{perspective triangles}\index{triangles!perspective}\defterm{perspective from a point} $o\in X$ if $\Aline oa=\Aline o{a'}$, $\Aline ob=\Aline o{b'}$ and $\Aline oc=\Aline o{c'}$ are three pairwise distinct lines in $X$;
\item \index{centrally perspective triangles}\index{triangles!centrally perspective}\defterm{centrally perspective} if the triangles $abc,a'b'c'$ are perspective from some point $o\in X$ (called the \index{perspector}\index{triangles!perspector of}\defterm{perspector} of the trianges $abc,a'b'c'$).
\end{itemize}
\end{definition}

\begin{theorem}[Desargues]\label{t:Desargues-projective} If $abc$ and $a'b'c'$ are two disjoint centrally perspective triangles in a $3$-long projective liner $X$ of rank $\|X\|\ne 3$, then the set 
$$T\defeq (\Aline ab\cap\Aline{a'}{b'})\cup(\Aline bc\cap\Aline {b'}{c'})\cup(\Aline ac\cap\Aline{a'}{c'})$$ has cardinality $|T|=3$ and rank $\|T\|=2$.
\end{theorem}

\begin{proof} By Theorem~\ref{t:projective<=>} and Corollary~\ref{c:proregular=>ranked}, the projective liner $X$ is strongly regular, ranked, and has the Exchange Property. Take any centrally perspective disjoint triangles $abc$ and $a'b'c'\in X^3$. Then there exists a point $o\in X$ such that $$A\defeq \Aline oa=\Aline o{a'}=\Aline a{a'},\quad B\defeq\Aline ob=\Aline o{b'}=\Aline{b}{b'},\quad\mbox{and}\quad C\defeq\Aline oc=\Aline o{c'}=\Aline c{c'}$$ are distinct lines in $X$. This implies that 
$$o\notin\Aline ab\cup\Aline bc\cup \Aline ac\cup\Aline {a'}{b'}\cup \Aline {b'}{c'}\cup\Aline {a'}{c'}.$$ Since the liner $X$ is projective and the centrally perspective triangles $abc,a'b'c'$ are disjoint, the sets $$T_a\defeq \Aline bc\cap\Aline{b'}{c'},\quad T_b\defeq \Aline ac\cap\Aline{a'}{c'},\quad T_c\defeq\Aline ac\cap\Aline{a'}{c'}$$are singletons. We have to prove that the set $T\defeq T_a\cup T_b\cup T_c$ has cardinality $|T|=3$ and rank $\|T\|=2$. 

\begin{picture}(140,150)(-130,-15)

\put(0,0){\color{blue}\line(1,0){120}}

\put(120,0){\line(0,1){120}}
\put(0,60){\color{red}\line(1,0){120}}
\put(60,0){\color{blue}\line(1,2){60}}

{\linethickness{1pt}
\put(0,0){\color{green}\line(0,1){120}}
\put(0,120){\color{green}\line(1,-1){80}}
}

\put(0,0){\color{blue}\circle*{3}}
\put(-4,-10){\color{blue}$a'$}

\put(80,40){\color{blue}\circle*{3}}
\put(82,34){\color{blue}$c'$}

\put(60,60){\color{red}\circle*{3}}
\put(57,51){\color{red}$c$}

\put(0,120){\color{green}\circle*{3}}
\put(-3,123){\color{green}$o$}

{\linethickness{1.5pt}
\put(0,0){\color{blue}\line(1,0){60}}
\put(0,0){\color{blue}\line(2,1){80}}
\put(60,0){\color{blue}\line(1,2){20}}

\put(0,60){\color{red}\line(1,0){60}}
\put(0,60){\color{red}\line(2,-1){40}}
\put(40,40){\color{red}\line(1,1){20}}
}

{\linethickness{1pt}
\put(60,0){\color{green}\line(-1,2){60}}
}

\put(0,60){\color{red}\circle*{3}}
\put(-8,57){\color{red}$a$}
\put(40,40){\color{red}\circle*{3}}
\put(33,32){\color{red}$b$}

\put(60,0){\color{blue}\circle*{3}}
\put(57,-10){\color{blue}$b'$}

\put(0,0){\color{blue}\line(2,1){120}}
\put(120,0){\color{red}\line(-2,1){120}}
\put(120,120){\color{red}\line(-1,-1){80}}

\put(120,0){\circle*{3}}
\put(123,-4){$T_c$}
\put(120,60){\circle*{3}}
\put(123,57){$T_b$}
\put(120,120){\circle*{3}}
\put(123,118){$T_a$}
\end{picture}

First we prove that $|T|=3$. Assuming that $T_a=T_b$, we conclude that $$T_a=T_b=\Aline bc\cap\Aline {b'}{c'}\cap\Aline ac\cap\Aline {a'}{c'}\subseteq (\Aline bc\cap\Aline ac)\cap\Aline{a'}{c'}=\{c\}\cap\Aline{a'}{c'}=\varnothing,$$
which is a contradiction showing that $T_a\ne T_b$. By analogy we can show that $T_a\ne T_c$ and $T_b\ne T_c$ and hence $|T|=|T_a\cup T_b\cup T_c|=3$ and $\|T\|\ge 2$.

It remains to prove that $\|T\|\le 2$. To derive a contradiction, assume that $\|T\|>2$.  Observe that $T\subseteq\overline{\{a,b,c\}}\cap\overline{\{a',b',c'\}}$ and hence $\|T\|=3$. Moreover, the $3$-rankedness of $X$ ensures that $\overline{\{a,b,c\}}=\overline{T}=\overline{\{a',b',c'\}}$ and hence $o\in \Aline{a}{a'}\subseteq \overline T$. Taking into account that the liner $X$ contains distinct lines $A,B,C$ and has rank $\|X\|\ne 3$, we conclude that $\|X\|\ge 4>3=\|\overline T\|$. Then there exists a point $o'\in X\setminus\overline T$. Since the liner $X$ is $3$-long, there exists a point $d\in \Aline{o'}{c}\setminus\{o',c\}$. The Exchange Axiom ensures that $o'\in \Aline {o'}c=\Aline dc$. Assuming that $d\in\overline T$, we conclude $o'\in\Aline dc\subseteq\overline T$, which contradicts the choice of the point $o'$. This contradiction shows that $d\notin\overline T $ and hence $d\ne o$. Then $\Aline od$ and $\Aline {o'}{c'}$ are two lines in the plane $\overline{\{o,c,o'\}}$. By the $0$-parallelity of projective liners, there exists a point $d'\in\Aline od\cap\Aline{o'}{c'}$. Assuming that $d'=d$, we conclude that $c'\ne d'=d\ne o'$ and hence $c'\in\Aline {o'}{d'}=\Aline {o'}d$. Since $c,c'\in\Aline {o'}d$ and $c\ne c'$, $o'\in \Aline {o'}{c}=\Aline {c'}{c}\subseteq\overline T$, which contradicts the choice of the point $o'$. This contradiction shows that $d'\ne d$.

\begin{claim}\label{cl:abd} $\|\overline{\{a,b,d\}}\cap\overline{\{a',b',d'\}}\|\le 2$.
\end{claim}

\begin{proof} Assuming that $\|\overline{\{a,b,d\}}\cap\overline{\{a',b',d'\}}\|>2$, we obtain $\overline{\{a,b,d\}}=\overline{\{a,b,d\}}\cap\overline{\{a',b',d'\}}=\overline{\{a',b',d'\}}$, by the $3$-rankedness of the projective liner $X$. Then $\Aline ab\cup\Aline{a'}{b'}\subseteq \overline T\cap\overline{\{a,b,d\}}=\overline T\cap\overline{\{a',b',d'\}}$. Since $\|\Aline ab\cup\Aline{a'}{b'}\|>2$, the $3$-rankedness of $X$ ensures that $\overline{\Aline ab\cup\Aline{a'}{b'}}=\overline T\cap\overline{\{a,b,d\}}=\overline T=\overline{\{a,b,d\}}$, which contradicts  $d\notin\overline T$.
\end{proof}

Consider the set $T'\defeq (\Aline ab\cap\Aline{a'}{b'})\cup(\Aline ad\cap\Aline{a'}{d'})\cup(\Aline bd\cap\Aline{b'}{d'})$ and observe that $T'\subseteq\overline{\{a,b,d\}}\cap\overline{\{a',b',d'\}}$. Claim~\ref{cl:abd}, implies $\|T'\|\le \|\overline{\{a,b,d\}}\cap\overline{\{a',b',d'\}}\|\le 2.$

Observe that $\Aline ad$ and $\Aline{a'}{d'}$ are two lines in the plane $\overline{\{o,a,d\}}$, and  $\Aline bd$, $\Aline{b'}{d'}$ are two lines in the plane $\overline{\{o,b,d\}}$. By the $0$-parallelity of projective liners, there exist points $s_a\in \Aline ad\cap\Aline{a'}{d'}$ and $s_b\in\Aline bd\cap\Aline {b'}{d'}$.

Observe that $s_a\in\Aline ad\cap\Aline {a'}{d'}\subseteq\overline{\{a,o',c\}}\cap\overline{\{a',o',c'\}}$ and hence $\Aline {o'}{s_a}\subseteq\overline{\{a,c,o'\}}\cap\overline{\{a',c',o'\}}$. Assuming that  $o'=s_a$, we conclude that $o'=s_a\in\overline{\{o,a,d\}}$ and hence $c\in\Aline{o'}{d}\subseteq\overline{\{o,a,d\}}$. Since $\|\{o,a,c\}\|=3=\|\{o,a,d\}\|$, the $3$-rankedness of $X$ ensures that $d\in\overline{\{o,a,d\}}=\overline{\{o,a,c\}}\subseteq\overline T$, which contradicts $d\notin\overline T$. This contradiction shows that $o'\ne s_a$ and hence $$2=\|\Aline {o'}{s_a}\|\le\|\overline{\{a,c,o'\}}\cap\overline{\{a',c',o'\}}\|\le 3=\|\overline{\{a,c,o'\}}\|=\|\overline{\{a',c',o'\}}\|.$$
Assuming that $\|\overline{\{a,c,o'\}}\cap\overline{\{a',c',o'\}}\|=3$, we conclude that $\overline{\{a,c,o'\}}=\overline{\{a,c,o'\}}\cap\overline{\{a',c',o'\}}=\overline{\{a',c',o'\}}$ and hence $\Aline ac\cup\Aline{a'}{c'}\subseteq \overline T\cap\overline{\{a,c,o'\}}$. Since $\Aline ac\ne \Aline {a'}{c'}$, the $3$-rankedness of the liner $X$ implies  $\overline T=\overline{\Aline ac\cup\Aline {a'}{c'}}=\overline{\{a,c,o'\}}$, which contradicts the choice of $o'\notin\overline T$. This contradiction implies that $\|\overline{\{a,c,o'\}}\cap\overline{\{a',c',o'\}}\|=2$ and hence $\Aline{o'}{s_a}=\overline{\{a,c,o'\}}\cap\overline{\{a',c',o'\}}$, by the $2$-rankedness of $X$. Then $T_b=\Aline ac\cap\Aline {a'}{c'}\subseteq \overline{\{a,c,o'\}}\cap\overline{\{a',c',o'\}}=\Aline{o'}{s_a}$. By analogy we can prove that $T_a\subseteq\Aline{o'}{s_b}$. Then 
$$T=T_a\cup T_b\cup T_c\subseteq \Aline {o'}{s_b}\cup\Aline{o'}{s_a}\cup T_c\subseteq \overline{\{o'\}\cup T'}.$$ Since 
$3=\|T\|\le\|\overline{\{o'\}\cup T'}\|\le \|T'\|+1\le 3$, the $3$-rankedness of $X$ implies the equality $\overline T=\overline{\{o'\}\cup T'}$, which contradicts the choice of the point $o'\notin \overline T$. This is a final contradiction completing the proof of Theorem~\ref{t:Desargues-projective}.
\end{proof}

\section{Desarguesian liners}

 Theorems~\ref{t:Desargues-affine}, \ref{t:Des-Proclus} and 
 \ref{t:Desargues-projective} motivate the following definition.

\begin{definition}\label{d:Desarguesian} A liner $X$ is called \index{Desarguesian liner}\index{liner!Desarguesian}\defterm{Desarguesian} if for every plane $\Pi\subseteq X$ and centrally perspective disjoint triangles $abc,a'b'c'$ in $\Pi$, the set    
$T\defeq (\Aline ab\cap\Aline{a'}{b'})\cup(\Aline bc\cap\Aline {b'}{c'})\cup(\Aline ac\cap\Aline{a'}{c'})$ has rank $\|T\|\in\{0,2\}$.
\end{definition}

\begin{picture}(140,150)(-130,-15)

\put(0,0){\color{blue}\line(1,0){120}}

\put(120,0){\line(0,1){120}}
\put(0,60){\color{red}\line(1,0){120}}
\put(60,0){\color{blue}\line(1,2){60}}

{\linethickness{1pt}
\put(0,0){\color{green}\line(0,1){120}}
\put(0,120){\color{green}\line(1,-1){80}}
}

\put(0,0){\color{blue}\circle*{3}}
\put(-4,-10){\color{blue}$a'$}

\put(80,40){\color{blue}\circle*{3}}
\put(82,34){\color{blue}$c'$}

\put(60,60){\color{red}\circle*{3}}
\put(57,51){\color{red}$c$}

\put(0,120){\color{green}\circle*{3}}

\put(0,0){\color{blue}\line(2,1){120}}

{\linethickness{1.5pt}
\put(0,60){\color{red}\line(1,0){60}}
\put(0,60){\color{red}\line(2,-1){40}}
\put(40,40){\color{red}\line(1,1){20}}

\put(0,0){\color{blue}\line(2,1){80}}

{\linethickness{1pt}
\put(60,0){\color{green}\line(-1,2){60}}
}

\put(0,0){\color{blue}\line(1,0){60}}
\put(60,0){\color{blue}\line(1,2){20}}
}

\put(0,60){\color{red}\circle*{3}}
\put(-8,57){\color{red}$a$}
\put(40,40){\color{red}\circle*{3}}
\put(33,32){\color{red}$b$}

\put(60,0){\color{blue}\circle*{3}}
\put(57,-10){\color{blue}$b'$}

\put(120,0){\color{red}\line(-2,1){120}}
\put(120,120){\color{red}\line(-1,-1){80}}

\put(120,0){\circle*{3}}
\put(125,-4){$\Aline ab\cap\Aline{a'}{b'}$}
\put(120,60){\circle*{3}}
\put(125,57){$\Aline ac\cap\Aline{a'}{c'}$}
\put(120,120){\circle*{3}}
\put(125,118){$\Aline bc\cap\Aline{b'}{c'}$}
\end{picture}



Theorem~\ref{t:Des-Proclus} and Definition~\ref{d:Desarguesian} imply the following useful characterization of Desarguesian Proclus weakly regular liners.

\begin{theorem}\label{t:D-Proclus-wr<=>} A Proclus weakly regular liner $X$ is Desarguesian if and only if for every centrally perspective disjoint triangles $abc$ and $a'b'c'$ in $X$, the set $$T\defeq(\Aline ab\cap\Aline{a'}{b'})\cup(\Aline bc\cap\Aline{b'}{c'})\cup(\Aline ac\cup\Aline{a'}{c'})$$ has rank $\|T\|\in\{0,2\}$.
\end{theorem} 

\begin{theorem}\label{t:Desargues<=>planeD} A Proclus liner $X$ is Desarguesian if and only if every $3$-long plane in $X$ is Desarguesian.
\end{theorem}

\begin{proof} The ``only if'' part is trivial. To prove the ``if'' part, assume that every $3$-long plane in a Prolcus liner $X$ is Desarguesian. To prove that $X$ is Desarguesian, take any plane $P\subseteq X$ and centrally perspective disjoint triangles $abc,a'b'c'$ in $P$. We have to prove that the set $T=(\Aline ab\cap\Aline {a'}{b'})\cup(\Aline bc\cap\Aline{b'}{c'})\cup(\Aline ac\cap\Aline{a'}{c'})$ has rank $\|T\|\in\{0,2\}$.  

The triangles $abc$ and $a'b'c'$ are perspective from the unique point $o\in\Pi$ of the intersection $\Aline a{a'}\cap\Aline b{b'}\cap\Aline c{c'}$. It follows that $|P|\ge|\{o,a,b,c,a',b',c'\}|=7$. If the Proclus plane $P$ is not projective, then $P$ is $3$-long, by Theorem~\ref{t:Proclus-not-3long}.
If  the plane $P$ is projective, then we can consider the maximal $3$-long flat $M\subseteq P$ that contains the point $o$. By Lemma~\ref{l:ox=2}, $\Aline ox=\{o,x\}$ for every $x\in P\setminus M$, which implies that $\{o,a,b,c,a',b',c'\}\subseteq M$ and hence $3\le\|\{a,b,c\}\|\le\|M\|\le\|P\|=3$. 
By Theorem~\ref{t:Proclus<=>} and Proposition~\ref{p:k-regular<=>2ex}, the Proclus plane is $3$-proregular and $3$-ranked. The $3$-rankedness of the plane $\Pi$ ensures that $M=P$ and hence the plane $P=M$ is $3$-long. 

In both cases, the plane $P$ is $3$-long. Since $3$-long planes in $X$ are Desarguesian, the set $T$ has rank $\|T\|\in\{0,2\}$ in $P$ and also in $X$. 
\end{proof}

\begin{theorem}\label{t:pD-minus-flat} For every flat $H$ in a Desarguesian projective liner $Y$, the subliner $X\defeq Y\setminus H$ of $Y$ is Desarguesian.
\end{theorem}

\begin{proof} To prove that the liner $X$ is Desarguesian, take any plane $P\subseteq X$ and disjoint centrally perspective traingles $abc$ and $a'b'c'$ in $P$. Observe that the flat hull $\overline P$ of the plane $P$ is a plane in the projective liner $Y$. Since the projective liner $Y$ is Desarguesian, the set $T\defeq(\Aline ab\cap \Aline {a'}{b'})\cup(\Aline bc\cap\Aline{b'}{c'})\cup(\Aline ac\cap\Aline{a'}{c'})$ has rank $\|T\|\in\{0,2\}$ in $\overline P$ and also in $Y$. Then the set $T\cap X$ has rank $\|T\cap X\|\le 2$ in the liner $X$. It remains to prove that $\|T\cap X\|\ne 1$.

 Repeating the argument from the proof of Theorem~\ref{t:Desargues-projective}, we can show that the set $T$ has cardinality $|T|=3$. 
Then $\overline{T}$ is a line in the projective liner $Y$. If $|T\cap H|\ge 2$, then $T\subseteq\overline{T}\subseteq H$ because $F$ is a flat. In this case $|T\cap X|=0$. If $|T\cap H|\le 1$, then $|T\cap X|\ge|T|-1=2$ and we are done. 
\end{proof}

Now we consider some important examples of Desarguesian liners. 

\begin{theorem}\label{t:proaffine-Desarguesian} Every $3$-long proaffine regular liner $X$ of rank $\|X\|\ne 3$ is Desarguesian.
\end{theorem}

\begin{proof} If $\|X\|\le 2$, then $X$ is Desarguesian because $X$ contains no (centrally perspective) triangles. So, assume that $\|X\|\ge 3$ and hence $\|X\|\ge 4$ (because $\|X\|\ne 3$). By Theorem~\ref{t:proaffine3=>compregular} the proaffine regular space $X$ is completely regular and by Corollary~\ref{c:spread-3long}, the spread completion $\overline X$ of $X$ is projective and $3$-long.  By Corollary~\ref{c:procompletion-rank}, $\|\overline X\|=\|X\|\ne 3$, and by Theorem~\ref{t:Desargues-projective}, the projective liner $\overline  X$ is Desarguesian. Since $\overline X\setminus X$ is flat in $\overline X$, the liner $X$ is Desarguesian, by Theorem~\ref{t:pD-minus-flat}.
\end{proof}

\begin{theorem} Every  balanced regular liner $X$ with $\|X\|\ne 3$ and $|X|_3<\w$ is Desarguesian.
\end{theorem}

\begin{proof} If $\|X\|\le 2$, then $X$ is Desarguesian because it contains no triangles. So, assume that $\|X\|>2$. It follows from $\|X\|\ne 3$ that $\|X\|\ge 4$. If $|X|_2=2$, then $X$ is Desarguesian because $X$ contains no disjoint centrally perspective triangles. So, assume that $|X|_2\ge 3$. By Theorem~\ref{t:balanced<=>ranked}, the $3$-long balanced liner $X$ of rank $\|X\|\ge 4$ is $p$-parallel for some $p\in\{0,1\}$ and hence $X$ is proaffine.  By Theorem~\ref{t:proaffine-Desarguesian}, the $3$-long proaffine liner $X$ of rank $\|X\|\ne 3$ is Desarguesian.
\end{proof} 

\begin{proposition}\label{p:Steiner+affine=>Desargues} Every affine liner $X$ with $|X|_2\le 3$ is Desarguesian.
\end{proposition}

\begin{proof} To prove that the liner $X$ is Desarguesian, take any plane $P\subseteq X$ and centrally perspective disjoint triangles $abc$ and $a'b'c'$ in $P$. We shall prove the set $T\defeq(\Aline ab\cap\Aline{a'}{b'})\cup(\Aline ac\cap\Aline{a'}{c'})\cup(\Aline bc\cap\Aline{b'}{c'})$ is empty and hence  $\|T\|=0\in\{0,2\}$. Let $o\in\Aline a{a'}\cap\Aline b{b'}\cap\Aline c{c'}\subseteq X$ be a unique perspector of the triangles $abc$ and $a'b'c'$. Since $|X|_2\le 3$, $\Aline oa=\{o,a,a'\}$, $\Aline ob=\{o,b,b'\}$, $\Aline oc=\{o,c,c'\}$.

Since $X$ is affine, for the point $b'\in\Aline ob\setminus \Aline ab$ there exists a unique point $u\in\Aline o{a}$ such that $\Aline u{b'}\cap\Aline ab=\varnothing$. Taking into account that $u\in\Aline o{a}\setminus \{o,a\}=\{a'\}$, we conclude that $\Aline {a'}{b'}\cap\Aline ab=\Aline u{b'}\cap\Aline ab=\varnothing$. By analogy we can show that $\Aline bc\cap\Aline{b'}{c'}=\varnothing=\Aline ac\cap\Aline {a'}{c'}$ and hence $T=\varnothing$, witnessing that the affine  liner $X$ is Desarguesian.
\end{proof}

\begin{corollary}\label{c:affine-Desarguesian} Every affine liner $X$ of rank $\|X\|\ne 3$ is Desarguesian.
\end{corollary}

\begin{proof} By Theorem~\ref{t:affine=>Avogadro}, the affine liner $X$ is $2$-balanced. If $|X|_2\le3$, then $X$ is Desarguesian, by Proposition~\ref{p:Steiner+affine=>Desargues}. So, assume that $|X|_2\ge 4$. By Theorem~\ref{t:4-long-affine}, the $4$-long affine liner $X$ is regular. By Theorem~\ref{t:proaffine-Desarguesian}, $3$-long affine regular liner $X$ of rank $\|X\|\ne 3$ is Desarguesian.
\end{proof}

\begin{proposition}\label{p:Steiner+projective=>Desargues} Every Steiner $3$-regular liner $X$ is  Desarguesian.
\end{proposition}

\begin{proof} By Theorem~\ref{t:Steiner-3reg}, the Steiner $3$-regular liner is $p$-parallel for some $p\in\{0,1\}$.

If $p=1$, then $X$ is affine. By Proposition~\ref{p:Steiner+affine=>Desargues}, the affine Steiner liner $X$ is Desarguesian.

So, assume that $p=0$ and hence the $0$-parallel liner $X$ is projective and (strongly) regular, by Theorem~\ref{t:projective<=>}. To show that $X$ is Desarguesian, take any plane $\Pi\subseteq X$ and two centrally perspective disjoint triangles $abc$ and $a'b'c'$ in $\Pi$.  Let $o\in\Aline a{a'}\cap\Aline b{b'}\cap\Aline c{c'}\subseteq \Pi$ be a unique perspector of the triangles $abc$ and $a'b'c'$. Since $X$ is Steiner, $\Aline oa=\{o,a,a'\}$, $\Aline ob=\{o,b,b'\}$, $\Aline oc=\{o,c,c'\}$.
Theorem~\ref{t:k-parallel=>2-balance} ensures that $|\Pi|=|X|_3=1+(|X|_2-1)(|X|_2+p)=1+2(3+p)=7$ and hence $\Pi=\{o,a,a',b,b',c,c'\}$ is a Stwiner projective plane. Since the triangles $abc$ and $a'b'c'$ are disjoint and centrally perspective, $\varnothing\ne \Aline ab\cap\Aline {a'}{b'}\subseteq X\setminus (\Aline a{a'}\cup\Aline b{b'})=\{c,c'\}$. 

If $\Aline ab\cap\Aline {a'}{b'}=\{c\}$, then $\{a,b,c\}$ is a line in $X$, which contradicts the choice of the triangle $abc$. If $\Aline ab\cap\Aline {a'}{b'}=\{c'\}$, then $\{a',b',c'\}=\Aline{a'}{b'}$ is a line in $X$, which contradicts the choice of the triangle $a'b'c'$. Those contradictions show that the Steiner projective plane $\Pi$ does not contain disjoint centrally perspective triangles and hence is vacuously Desarguesian.
\end{proof}

\begin{example}[Ivan Hetman] There exists a non-Desarguesian Steiner ranked plane.
\end{example}

\begin{proof} It can be shown that the ring $\IZ_{31}$ endowed with the family of lines $\mathcal L\defeq\{L+x:x\in\IZ_{31},\;L\in\mathcal B\}$ where
$\mathcal B\defeq\big\{\{0, 14, 18\}, \{0, 15, 21\}, \{0, 19, 20\}, \{0, 22, 24\}, \{0, 23, 26\}\big\}$ is a non-Desarguesian Steiner ranked plane. The triangles $abc\defeq(14,27,11)$ and $a'b'c'\defeq(18,13,30)$ are perspective from the point $o\defeq 0$, but the set $T\defeq(\Aline ab\cap\Aline{a'}{b'})\cup(\Aline bc\cap\Aline {b'}{c'})\cup(\Aline ac\cap\Aline{a'}{c'})=\{10,17,19\}$ has rank $\|T\|=3$.
\end{proof}

\begin{Exercise}\label{Ex:Des8} Prove that every $2$-balanced projective liner with $|X|_2\le 9$ is Desarguesian.
\smallskip

{\em Hint:} This follows from the uniqueness of a projective plane of order $\le 8$. The authors of \cite{HSW1956} write that the uniqueness of projective planes of order $\le 5$ is an easy exercise (see, however Corollary \ref{c:4-Pappian}). The uniqueness of projective planes of order $7$ and $8$ was proved with the help of computers in 1956, see \cite{HSW1956}. By Bruck-Ryser Theorem~\ref{t:Bruck-Ryser}, no projective plane of order $6$ exists.
\end{Exercise}  

\begin{exercise} Find an example of a $2$-balanced non-Desarguesian projective plane $X$ with $|X|_2=10$.
\smallskip

{\em Hint:} Consider the projective completion of the non-Desarguesian affine plane $\mathbb J_9\times\mathbb J_9$, considered in Section~\ref{s:J9xJ9}.
\end{exercise}


\begin{exercise} Find an example of a non-Desarguesian projective plane $X$ with $|X|_2=\mathfrak c$.
\smallskip

{\em Hint:} Consider the projective completion of the Moulton plane, considered in Section~\ref{s:Moulton}.
\end{exercise}






\section{Desarguesian affine liners}

The main result of this section is the following characterization of Desarguesian affine liners.\index[person]{Pra\.zmowska}
\vskip20pt

\begin{theorem}[Pra\.zmowska\footnote{{\bf Ma\l gorzata Pra\.zmowska} is an adjunkt in Bia\l ystok University in Poland. The proof of Theorem~\ref{t:ADA<=>} is taken from the paper \cite{Prazmowska2004}.}, 2004]\label{t:ADA<=>} An affine liner $X$ is Desarguesian\newline if and only if $X$ satisfies  the \index{Affine Desargues Axiom}\index{Axiom!Affine Desargues}\defterm{Affine Desargues Axiom}:
\begin{itemize}
\item[{\sf(ADA)}] for every plane $P\subseteq X$ and centrally perspective triangles\\ $abc,a'b'c'$ in $P$, 
 if $\Aline ab\cap\Aline{a'}{b'}=\varnothing=\Aline bc\cap\Aline {b'}{c'}$, then $\Aline ac\cap\Aline{a'}{c'}=\varnothing$.
\end{itemize}

\begin{picture}(100,0)(-370,-35)

\put(0,0){\line(0,1){33}}
\put(0,0){\line(0,-1){33}}
\put(0,0){\line(7,4){28}}
\put(0,0){\line(-7,-4){28}}
\put(0,0){\line(7,-4){28}}
\put(0,0){\line(-7,4){28}}

{\linethickness{1pt}
\put(28,-16){\color{blue}\line(-4,7){28}}
\put(-28,-16){\color{cyan}\line(4,7){28}}
\put(28,16){\color{cyan}\line(-4,-7){28}}
\put(-28,16){\color{blue}\line(4,-7){28}}
\put(28,-16){\color{red}\line(-1,0){56}}
\put(28,16){\color{red}\line(-1,0){56}}
}

\put(0,0){\circle*{3}}
\put(28,16){\circle*{3}}
\put(28,19){$a'$}
\put(28,-16){\circle*{3}}
\put(28,-24){$c$}
\put(-28,16){\circle*{3}}
\put(-32,19){$c'$}
\put(-28,-16){\circle*{3}}
\put(-32,-24){$a$}
\put(0,33){\circle*{3}}
\put(-2,36){$b$}
\put(0,-33){\circle*{3}}
\put(-2,-43){$b'$}

\end{picture}
\end{theorem}
\smallskip

We divide the proof of Theorem~\ref{t:ADA<=>} into four lemmas.

\begin{lemma}\label{l:D=>ADA} Every Desarguesian liner $X$ satisfies the Affine Desargues Axiom.
\end{lemma}

\begin{proof} To prove that $X$ satisfies the Affine Desargues Axiom, take any plane $P\subseteq X$ and two centrally perspective  triangles $abc$ and $a'b'c'$ in $X$ such that $\Aline ab\cap\Aline {a'}{b'}=\varnothing=\Aline bc\cap\Aline{b'}{c'}$. The latter equalities imply that the centrally perspective triangle $abc$ and $a'b'c'$ are disjoint. Since $X$ is Desarguesian, the set $T\defeq (\Aline ab\cap\Aline {a'}{b'})\cup(\Aline bc\cap\Aline{b'}{c'})\cup(\Aline ac\cap\Aline{a'}{c'})=\Aline ac\cap\Aline {a'}{c'}$ has cardinality $|T|\ne 1$ and rank $\|T\|\le 2$. The perspectivity of the disjoint triangles implies that $|T|=|\Aline ac\cap\Aline {a'}{c'}|\le 1$. Since $|T|\ne 1$, the set $T=\Aline ac\cap\Aline{a'}{c'}$ is empty, witnessing that $X$ satisfies the Affine Desargues Axiom.
\end{proof}

\begin{lemma}\label{l:ADA<=>} For any $3$-ranked affine liner $X$, the following conditions are equivalent:
\begin{enumerate}
\item $X$  satisfies the Affine Desargues Axiom;
\item  for any coplanar concurrent lines $A,B,C$ in $X$ and points $a,a'\in A\setminus(B\cup C)$, $b,b'\in B\setminus(A\cup C)$, $c,c'\in C\setminus(A\cup B)$, if $\Aline ab\parallel \Aline{a'}{b'}$ and $\Aline bc\parallel \Aline{b'}{c'}$, then $\Aline ac\parallel \Aline{a'}{c'}$.
\end{enumerate}
\end{lemma}

\begin{proof} Let $X$ be a $3$-ranked affine liner. By Theorem~\ref{t:parallel-char}, two disjoint lines $L,L'$ in the $3$-ranked liner $X$ are parallel if and only if they are coplanar.
\smallskip

$(1)\Ra(2)$. Assume that the liner $X$ satisfies the Affine Desargues Axiom. Given any coplanar concurrent lines $A,B,C\subseteq X$ and points $a,a'\in A\setminus(B\cup C)$, $b,b'\in B\setminus(A\cup C)$, $c,c'\in C\setminus(A\cup B)$ with $\Aline ab\parallel \Aline {a'}{b'}$ and $\Aline bc\parallel \Aline {b'}{c'}$, we have to prove that $\Aline ac\parallel\Aline {a'}{c'}$. Since the lines $A,B,C$ are concurrent, there exists a point $o\in A\cap B\cap C$. The choice of the points $a,b,c$ ensures that the concurrent lines $A,B,C$ are distinct. 



If $b=b'$, then $\Aline ab\parallel\Aline{a'}{b'}$ and $\Aline bc\parallel \Aline {b'}{c'}$ imply $\Aline ab=\Aline{a'}{b'}$ and $\Aline bc=\Aline {b'}{c'}$. Then $\{a\}=A\cap\Aline ab=A\cap\Aline{a'}{b'}=\{a'\}$ and $\{c\}=C\cap\Aline bc=C\cap\Aline{b'}{c'}=\{c'\}$ and hence $\Aline ac=\Aline{a'}{c'}$ and $\Aline ac\parallel\Aline{a'}{c'}$.

So, assume that $b\ne b'$. In this case $\Aline ab\cap\Aline {a'}{b'}=\varnothing=\Aline bc\cap\Aline{b'}{c'}$. Indeed, assuming that $\Aline ab\cap\Aline {a'}{b'}\ne\varnothing$ and taking into account that $\Aline ab\parallel \Aline {a'}{b'}$, we conclude that $\Aline ab=\Aline {a'}{b'}$ and hence $\{b\}=B\cap\Aline ab=B\cap\Aline{a'}{b'}=\{b'\}$, which contradicts the assumption $b\ne b'$. By analogy we can show that $b\ne b'$ implies $\Aline bc\cap\Aline{b'}{c'}=\varnothing$. Now the Affine Desargues Axiom ensures that $\Aline ac\cap\Aline{a'}{c'}=\varnothing$ and hence $\Aline ac\parallel \Aline{a'}{c'}$, by Theorem~\ref{t:parallel-char} (because $\|\Aline ac\cup\Aline{a'}{c'}\|=\|A\cup C\|=3$). 
\smallskip

$(2)\Ra(1)$ Assume that the liner $X$ satisfies the condition (2). To prove that $X$ is satisfies the Affine Desargues Axiom, take any plane $P$ and centrally perspective triangles $abc,a'b'c'$ in $P$, and assume that $\Aline ab\cap\Aline{a'}{b'}=\varnothing=\Aline bc\cap\Aline {b'}{c'}$. We have to prove that $\Aline ac\cap\Aline {a'}{c'}=\varnothing$. By Corollary~\ref{c:parallel-lines<=>}, $\Aline ab\cap\Aline{a'}{b'}=\varnothing=\Aline bc\cap\Aline {b'}{c'}$ imply $\Aline ab\parallel\Aline {a'}{b'}$ and $\Aline bc\parallel\Aline {b'}{c'}$. The condition (2) ensures that $\Aline ac\parallel \Aline {a'}{c'}$. Assuming that $\Aline ac\cap\Aline {a'}{c'}\ne\varnothing$, we conclude that the parallel lines $\Aline ac$ and $\Aline {a'}{c'}$ are equal and hence $\{a\}=A\cap \Aline ac=A\cap\Aline{a'}{c'}=\{a'\}$, which contradicts $\Aline ab\cap\Aline{a'}{b'}=\varnothing$. This contradiction shows that $\Aline ac\cap\Aline{a'}{c'}=\varnothing$, witnesing that the affine $3$-ranked liner $X$ satisfies the Affine Desargies Axiom.
\end{proof}

\begin{lemma}\label{l:ADA+P} Let $X$ be an affine plane satisfying the Affine Desargues Axiom, and let $abc$, $a'b'c'$ be two centrally perspective triangles in $X$ such that $\Aline ac\cap \Aline {a'}{c'}=\varnothing$ and $\Aline ab\cap\Aline{a'}{b'}=\{x\}$  for some point $x\in X$. Then $\Aline bc\cap\Aline {b'}{c'}=\{y\}$ for some point $y\in X$ such that $\Aline xy\subparallel \Aline ac$.
\end{lemma}

\begin{proof} Since the triangles $abc$ and $a'b'c'$ are centrally perspective, there exists a point $o\in X$ such that $A\defeq\Aline oa=\Aline o{a'}$, $B\defeq\Aline ob=\Aline o{b'}$, $C\defeq\Aline oc=\Aline o{c'}$ are parwise distinct lines in $X$. Then $o\notin\{a,a',b,b',c,c'\}$ and hence $a,a'\in A\setminus(B\cup C)$, $b,b'\in B\setminus(A\cup C)$, $c,c'\in C\setminus(A\cup B)$.

By Theorem~\ref{t:affine=>Avogadro}, the affine liner $X$ is $2$-balanced. 
It follows from $\Aline ac\cap\Aline{a'}{c'}=\varnothing$ that $|X|_2\ge|\Aline oa|\ge|\{o,a,a'\}|\ge 3$.
Since $X$ is affine and $a'\in \Aline oa\setminus(\{o\}\cup\Aline ab)$, there exists a unique point  $u\in\Aline ob\setminus\{o\}$ such that $\Aline {a'}{u}\cap\Aline ab=\varnothing$. If $|X|_3=3$, then $u\in\Aline ob\setminus\{o,b\}=\{b'\}$ and hence $\{x\}=\Aline ab\cap\Aline{a'}{b'}=\Aline ab\cap\Aline {a'}u=\varnothing$, which is a contradiction showing that $|X|_2\ge 4$. By Theorem~\ref{t:4-long-affine}, the $4$-long affine liner $X$ is regular. 

If $\Aline bc\cap\Aline {b'}{c'}=\varnothing$, then the Affine Desargues Axiom implies $\Aline ab\cap\Aline {a'}{b'}=\varnothing$, which contradicts $\Aline ab\cap\Aline{a'}{b'}=\{x\}$. This contradiction shows that there exists a point $y\in\Aline bc\cap\Aline{b'}{c'}$. Assuming that $\Aline bc\cap\Aline{b'}{c'}\ne\{y\}$, we conclude that $\Aline bc=\Aline {b'}{c'}$ and then $\{c\}=C\cap\Aline bc=C\cap\Aline {b'}{c'}=\{c'\}$, which contradicts $\Aline ac\cap\Aline {a'}{c'}=\varnothing$. Therefore, $\Aline bc\cap\Aline{b'}{c'}=\{y\}$. It remains to show that $\Aline xy\subparallel \Aline ac$.

 If $b=b'$, then $b=b'\in \Aline ab\cap\Aline{a'}{b'}\cap\Aline bc\cap\Aline{b'}{c'}=\{x\}\cap\{y\}$ and hence the singleton $\Aline xy=\{x\}=\{y\}$ is subparallel to the line $\Aline ac$.
So, we assume that $b\ne b'$.

\begin{picture}(140,155)(-100,-15)

{
\put(60,0){\line(1,2){60}}
\put(0,0){\line(0,1){120}}
\put(0,120){\line(1,-1){90}}
\put(0,40){\color{cyan}\line(3,1){240}}
\put(0,0){\line(2,1){240}}
\put(120,120){\line(-1,-1){80}}
\put(60,0){\line(-1,2){60}}

}

{\linethickness{1pt}
\put(0,40){\color{red}\line(1,0){40}}
\put(0,40){\color{cyan}\line(3,1){60}}
\put(40,40){\color{blue}\line(1,1){20}}
\put(0,0){\line(2,1){80}}
\put(60,0){\line(1,2){20}}

}

{\linethickness{1.5pt}
\put(0,0){\color{cyan}\line(3,1){90}}
\put(0,0){\color{red}\line(1,0){60}}
\put(60,60){\color{cyan}\line(3,1){180}}
\put(120,120){\color{red}\line(1,0){120}}
\put(60,60){\color{blue}\line(1,1){60}}
\put(60,0){\color{blue}\line(1,1){30}}
}

{\linethickness{1pt}
\put(60,0){\line(1,2){20}}
}

{
\put(60,0){\line(-1,2){60}}
\put(0,0){\line(2,1){240}}

}

\put(0,120){\circle*{3}}
\put(-3,123){$o$}

\put(0,0){\circle*{3}}
\put(-4,-8){$a$}
\put(80,40){\circle*{3}}
\put(77,44){$b$}
\put(60,60){\circle*{3}}
\put(58,64){$b'$}
\put(0,40){\circle*{3}}
\put(-10,38){$a'$}
\put(40,40){\circle*{3}}
\put(45,39){$c'$}
\put(60,0){\circle*{3}}
\put(57,-8){$c$}

\put(240,120){\circle*{3}}
\put(244,118){$x$}

\put(90,30){\circle*{3}}
\put(94,28){$b''$}
\put(120,120){\circle*{3}}
\put(118,124){$y$}
\end{picture}

Since $X$ is affine, there exists a point $b''\in \Aline ob$ such that $\Aline {a}{b''}\cap \Aline {a'}{b'}=\varnothing$ and hence $\Aline a{b''}\parallel \Aline{a'}{b'}$, by Theorem~\ref{t:parallel-char}.  Since $\Aline {a}{b}$ intersects the line $\Aline {a'}{b'}$, the point $b''$ is not equal to the points $b$ and $b'$. By Lemma~\ref{l:ADA<=>}, $\Aline {b''}{c}\parallel \Aline {b'}{c'}$. Observe that the triangles $acb''$ and $xyb'$ are perspective from the point $b$. Since $\Aline a{b''}\parallel \Aline {a'}{b'}=\Aline{b'}x$ and $\Aline {b''}c\parallel \Aline{b'}{c'}$, Lemma~\ref{l:ADA<=>} ensures that $\Aline xy\parallel \Aline ac$.
\end{proof}

\begin{lemma}\label{t:ADA=>Desargues} If an affine  liner $X$ satisfies the Affine Desargues Axiom, then $X$ is Desarguesian.
\end{lemma}

\begin{proof} Assume that an affine liner $X$ satisfies the Affine Desargues Axiom. To prove that $X$ is Desarguesian, take any plane $P$ and centrally perspective disjoint triangles $abc$ and $a'b'c'$ in $P$. We have to prove that the set $T\defeq (\Aline ab\cap\Aline {a'}{b'})\cup(\Aline bc\cap\Aline{b'}{c'})\cup(\Aline ac\cap\Aline{a'}{c'})$ has cardinality $|T|\ne 1$ and rank $\|T\|\le 2$. Let $o\in P$ be a perspector of the centrally perspective triangles $abc$ and $a'b'c'$. Then  $$A\defeq\Aline oa=\Aline o{a'}=\Aline a{a'},\quad B\defeq\Aline ob=\Aline o{b'}=\Aline b{b'}\quad\mbox{and}\quad C\defeq\Aline oc=\Aline o{c'}=\Aline c{c'}$$ are distinct lines with $a,a'\in A\setminus(B\cup C)$, $b,b'\in B\setminus(A\cup C)$ and $c,c'\in C\setminus(A\cup B)$.
 
Since the triangles $abc$ and $a'b'c'$ are disjoint, $\Aline ab\ne\Aline{a'}{b'}$, $\Aline bc\ne\Aline{b'}{c'}$ and $\Aline ac\ne \Aline {a'}{c'}$. In this case $$\max\{|\Aline ab\cap\Aline {a'}{b'}|,|\Aline bc\cap\Aline {b'}{c'}|,|\Aline ac\cap\Aline{a'}{c'}|\}\le 1.$$ If  two of the sets  $\Aline ab\cap\Aline{a'}{b'}$, $\Aline bc\cap\Aline {b'}{c'}$, $\Aline ac\cap\Aline{a'}{c'}$ are empty, then the Affine Desargues Axiom guarantees that the third set is empty and hence $|T|=0\ne 1$ and $\|T\|=0\le 2$.

If just one of the sets  $\Aline ab\cap\Aline{a'}{b'}$, $\Aline bc\cap\Aline {b'}{c'}$, $\Aline ac\cap\Aline{a'}{c'}$ is empty and two other are not empty, then $|T|=2$ and $\|T\|=2$.

 So, assume that the sets $\Aline ab\cap\Aline{a'}{b'}$, $\Aline bc\cap\Aline {b'}{c'}$, $\Aline ac\cap\Aline{a'}{c'}$ are nonempty and find unique points $x,y,z\in X$ such that $\Aline ab\cap\Aline{a'}{b'}=\{x\}$, $\Aline bc\cap\Aline{b'}{c'}=\{y\}$, and $\Aline ac\cap\Aline{a'}{c'}=\{z\}$. Then $T=\{x,y,z\}$ and we have to prove that  $|\{x,y,z\}|\le 1$ and $\|\{x,y,z\}\|\le 2$. Repeating the argument from the proof of Theorem~\ref{t:Desargues-projective}, we can show that $|\{x,y,z\}|=3$. It remains to show that $\|\{x,y,z\}|\le 2$.  To derive a contradiction, assume that $\|T\|=3$.


 By Theorem~\ref{t:affine=>Avogadro}, the affine liner $X$ is $2$-balanced. Observe that $|X|_2\ge|A|\ge|\{o,a,a'\}|=3$. Since $X$ is affine and $a'\in \Aline oa\setminus(\{o\}\cup \Aline ab)$, by the affinity of $X$, there exists a unique point $b''\in \Aline ob$ such that $\Aline {a'}{b''}\cap\Aline ab=\varnothing$. If $|X|_2=3$, then $b''\in\Aline ob\setminus\{o,b\}=\{b'\}$ and hence $\{x\}=\Aline ab\cap\Aline {a'}{b'}=\Aline ab\cap\Aline a'{b''}=\varnothing$, which is a contradiction showing that $|X|_2\ge 4$. By Theorem~\ref{t:4-long-affine}, the $4$-long affine liner $X$ is regular and hence Playfair.

Since $X$ is affine and regular, there exist points $a''\in \Aline oa\setminus \{o\}$ and $c''\in\Aline oc\setminus\{o\}$ such that $\Aline {a''}{c'}\cap \Aline ac=\varnothing=\Aline ac\cap \Aline{a'}{c''}$ and hence $\Aline{a''}{c'}\parallel \Aline ac\parallel \Aline{a'}{c''}$ and $\Aline {a''}{c'}\parallel \Aline {a'}{c''}$, by Theorem~\ref{t:Proclus-lines}. Assuming that $b'\in\Aline {a''}{c'}\cap\Aline{a'}{c''}$, we conclude that $\Aline {a''}{c'}=\Aline{a'}{c''}$ and hence $\|\{a',b',c'\}\|\le 2$, which contradicts the  definition of the triangle $a'b'c'$. This contradiction shows that $b'\notin\Aline {a''}{c'}$ or $b'\notin\Aline{a'}{c''}$. 
We lose no generality assuming that $b'\notin\Aline{a'}{c''}$ (in the opposite case, rename the points $a,c,a',c'$ by $c,a,c',a'$). Then $a'b'c''$ is a triangle, which is perspective to the triangle $zyc$ from the point $c'$.

\begin{picture}(140,150)(-100,-15)

{\linethickness{0.75pt}
\multiput(0,5)(0,10){6}{\color{cyan}\line(0,1){5}}
\multiput(0,0)(0,10){6}{\color{blue}\line(0,1){5}}
\put(0,60){\color{cyan}\line(2,-1){120}}

\put(120,0){\line(0,1){120}}
\put(0,60){\color{cyan}\line(1,0){120}}
\put(60,0){\line(1,2){60}}
\put(0,60){\color{blue}\line(0,1){60}}
\put(0,120){\color{blue}\line(1,-1){80}}
\put(0,40){\line(3,1){240}}
\put(0,0){\line(2,1){240}}
\put(120,120){\line(-1,-1){80}}
\put(120,120){\color{red}\line(1,0){120}}
\put(60,0){\color{blue}\line(-1,2){60}}

}

{\linethickness{1pt}
\put(60,0){\color{red}\line(1,0){60}}
\put(120,0){\line(0,1){60}}
\put(0,0){\line(2,1){120}}
\put(0,40){\color{red}\line(1,0){40}}
\put(0,40){\line(3,1){60}}

}

{\linethickness{1.5pt}
\put(0,0){\color{red}\line(1,0){60}}
\put(0,0){\line(2,1){80}}
\put(60,0){\line(1,2){20}}

\put(0,60){\line(1,0){60}}
\put(0,60){\line(2,-1){40}}
\put(40,40){\line(1,1){20}}
}

\put(0,120){\circle*{3}}
\put(-3,123){$o$}

\put(0,0){\circle*{3}}
\put(-4,-8){$c$}
\put(80,40){\circle*{3}}
\put(77,44){$b$}
\put(60,60){\circle*{3}}
\put(58,64){$b'$}
\put(0,40){\circle*{3}}
\put(-10,38){$c''$}
\put(0,60){\circle*{3}}
\put(-9,58){$c'$}
\put(40,40){\circle*{3}}
\put(45,39){$a'$}
\put(60,0){\circle*{3}}
\put(57,-8){$a$}

\put(240,120){\circle*{3}}
\put(244,118){$y'$}

\put(120,0){\circle*{3}}
\put(123,-4){$z$}
\put(120,60){\circle*{3}}
\put(123,54){$y$}
\put(120,120){\circle*{3}}
\put(123,112){$x$}
\put(123,123){$x'$}
\end{picture}

Now consider the centrally perspective triangles $abc$ and $a'b'c''$. Since $\Aline ab\cap\Aline{a'}{b'}=\{x\}$ and $\Aline ac\cap\Aline{a'}{c''}=\varnothing$, Lemma~\ref{l:ADA+P}, ensures that $\Aline bc\cap\Aline {b'}{c''}=\{y'\}$ for some point $y'\in X$ such that  $\Aline x{y'}\subparallel \Aline ac$.  Assuming that $x=y'$, we obtain a contradiction
$$\{x\}=\{y'x\}=\Aline ab\cap \Aline {a'}{b'}\cap \Aline bc\cap\Aline {b'}{c''}\subseteq (\Aline ab\cap\Aline bc)\cap\Aline{a'}{b'}=\{b\}\cap\Aline {a'}{b'}=\varnothing,$$showing that $x\ne y'$. By Corollary~\ref{c:subparallel}, $\Aline x{y'}\subparallel \Aline ac$ implies $\Aline x{y'}\parallel \Aline ac$.

Finally, consider two triangles $a'b'c''$ and $zyc$, which are perspective from the point $c'$, and observe that $\Aline {b'}{c''}\cap\Aline yc=\Aline {b'}{c''}\cap \Aline bc=\{y'\}$. Since $\Aline {a'}{c''}\cap \Aline ac=\varnothing$,  Lemma~\ref{l:ADA+P} ensures that $\Aline {a'}{b'}\cap\Aline zy=\{x'\}$ for some point $x'\in X$ such that $\Aline {x'}{y'}\subparallel \Aline ac$. Assuming that $x'=y'$ and taking into account that $a'b'c''$ is a triangle, we conclude that
$$\{x'\}=\{y'\}=\Aline {a'}{b'}\cap\Aline zy\cap \Aline bc\cap\Aline{b'}{c''}\subseteq (\Aline {b'}{a'}\cap\Aline {b'}{c''})\cap\Aline bc=\{b'\}\cap\Aline bc=\varnothing,$$
which is a contradiction showing that $x'\ne y'$. Then $\Aline {x'}{y'}\subparallel \Aline ac$ implies $\Aline {x'}{y'}\parallel \Aline ac$, according to Corollary~\ref{c:subparallel}. It follows from $\Aline x{y'}\parallel \Aline ac\parallel \Aline {x'}{y'}$ that $\Aline xy=\Aline {x'}{y'}$.  It follows from $y'\in \Aline {b'}{c''}\setminus\{b'\}=\Aline{b'}{c''}\setminus \Aline {a'}{b'}$ that 
$$\{x\}=\Aline{a'}{b'}\cap\Aline x{y'}=\Aline {a'}{b'}\cap\Aline{x'}{y'}=\{x'\}\subseteq \Aline zy$$ and hence $\|T\|=\|\{x,y,z\}\|\le\|\{y,z\}\|\le 2$.
\end{proof}



\section{Moufang liners}

In this section we introduce Moufang liners, prove that every Desarguesian liner is Moufang and every Moufang proaffine regular liner is completely regular.  

\begin{definition}\label{d:Moufang} A liner $X$ is called \defterm{Moufang} (and also \defterm{uno-Desarguesian}) if it satisfies the \defterm{uno-Desargues Axiom}
\begin{itemize}
\item[$(D_1$)] for every disjoint centrally perspective triangles $abc$ and $a'b'c'$ with $b'\in \Aline ac$, the set $$T\defeq(\Aline ab\cap\Aline {a'}{b'})\cup(\Aline bc\cap\Aline {b'}{c'})\cup(\Aline ac\cap\Aline {a'}{a'})$$ has rank $\|T\|\in\{0,2\}.$
\end{itemize}
\end{definition}

\begin{remark} Liners satisfying the the uno-Desargues Axiom $(D_1)$ are called Moufang in honour of Ruth Moufang who deeply studied such liners and proved that affine (or projective) Moufang planes are isomorphic to (projective completions of) coordinate planes of alternative rings. We shall consider those results of Ruth Moufang in Section~\ref{s:Moufang-proj1} and Chapter~\ref{ch:Moufang}.
\end{remark}

The following two characterizations of Moufang liners are counterparts of Theorem~\ref{t:Desargues<=>planeD} characterizing Desarguesian liners.

\begin{proposition}\label{p:Moufang<=>planeM} A liner $X$ is Moufang if and only if every plane in $X$ is Moufang.
\end{proposition}

\begin{proof} The ``only if'' part is trivial. To prove the ``if'' part, assume that every  plane in $X$ is Moufang. To prove that $X$ is Moufang, take any two disjoint centrally perspective triangles $abc$ and $a'b'c'$ in $X$ such that $b'\in \Aline ac$. Since the triangles $abc$, $a'b'c'$ are cemtrally perspective, there exists a unique point $o\in \Aline a{a'}\cap\Aline{b}{b'}\cap\Aline c{c'}$. Consider the plane $\Pi\defeq\overline{\{a,o,c\}}$ and observe that $b'\in \Aline ac\subseteq\Pi$ and $\{a',b,c'\}\subseteq\Aline oa\cup\Aline o{b'}\cup\Aline oc\subseteq\Pi$.
By our assumption, the plabe $\Pi$ is Moufang and hence the set
$$T\defeq(\Aline ab\cap \Aline{a'}{b'})\cup(\Aline bc\cap\Aline{b'}{c'})\cup(\Aline ac\cap\Aline{a'}{c'})$$has rank $\|T\|\in\{0,2\}$, witnessing that the liner $X$ is Moufang.
\end{proof}

For Proclus liners, Proposition~\ref{p:Moufang<=>planeM} can be improved as follows.

\begin{proposition}\label{p:Moufang<=>3plane} A Proclus liner $X$ is Moufang if and only if every $3$-long plane in $X$ is Moufang.
\end{proposition}

\begin{proof} The ``only if'' part is trivial. To prove the ``only if'' part, assume that every $3$-long plane in a Proclus liner $X$ is Moufang. By Theorem~\ref{t:Proclus<=>}, the Proclus liner $X$ is $3$-proregular, and by Proposition~\ref{p:k-regular<=>2ex}, the $3$-proregular liner $X$ is $3$-ranked. By Proposition~\ref{p:Moufang<=>planeM}, the Moufang property of $X$ will follow as soon as we show that every plane $P$ in $X$ is Moufang. So, take any disjoint centrally perspective triangles $abc,a'b'c'$ in $P$ and find the perspector $o\in \Aline a{a'}\cap\Aline b{b'}\cap\Aline c{c'}\subseteq P$ of the triangles $abc$ and $a'b'c'$. Observe that $|P|\ge|\{o,a,b,c,a',b',c'\}|\ge 7$.

If the Proclus plane $P$ is not projective, then we can apply Theorem~\ref{t:Proclus-not-3long} and conclude that the plane $P$ is $3$-long. 
If  the plane $P$ is projective, then we can consider the maximal $3$-long flat $M\subseteq P$ that contains the point $o$. By Lemma~\ref{l:ox=2}, $\Aline op=\{o,p\}$ for all points $p\in P\setminus M$, which implies that $\{a,a',b,b',c,c'\}\subseteq M$ and hence $3=\|A\cup B\cup C\cup D\|\le\|M\|\le\|P\|=3$. The $3$-rankedness of the Proclus liner $X$ ensures that $M=P$ and hence the plane $P=M$ is $3$-long.

Therefore, in both cases the Proclus plane $P$ is $3$-long. By our assumption, the $3$-long plane $P$ is Moufang and hence the set
$$T\defeq(\Aline ab\cap \Aline{a'}{b'})\cup(\Aline bc\cap\Aline{b'}{c'})\cup(\Aline ac\cap\Aline{a'}{c'})$$has rank $\|T\|\in\{0,2\}$, witnessing that the liner $X$ is Moufang.
\end{proof}

\begin{proposition}\label{p:Desarg=>Moufang} Every Desarguesian liner is Moufang.
\end{proposition}

\begin{proof} Let $X$ be a Desarguesian liner. To prove that $X$ is Moufang,  take any disjoint centrally perspective triangles $abc$ and $a'b'c'$ with $b'\in \Aline ac$. We have to prove that the set $$T\defeq(\Aline ab\cap\Aline {a'}{b'})\cup(\Aline bc\cap\Aline {b'}{c'})\cup(\Aline ac\cap\Aline {a'}{a'})$$has rank $\|T\|\in\{0,2\}$.

 Let $o\in \Aline a{a'}\cap \Aline b{b'}\cap\Aline c{c'}$ be the perspector of the prespective triangles $abc$ and $a'b'c'$. Observe that the plane $P\defeq\overline{\{a,o,c\}}$ contains points $a'\in \Aline oa$, $c'\in \Aline oc$, $b'\in \Aline ac$ and $b'\in \Aline ob$. Therefore, $abc$ and $a'b'c'$ are two centrally perspective triangles in the plane $P$. Since the liner $X$ is Desarguesian, the set $T$ has rank $\|T\|\in\{0,2\}$. 
\end{proof}

\begin{proposition}\label{p:Moufang-minus-flat} For every flat $H$ in a Moufang projective liner $Y$, the subliner $X\defeq Y\setminus H$ of $Y$ is Moufang.
\end{proposition}

\begin{proof} To prove that the liner $X$ is Moufang, take any disjoint centrally perspective traingles $abc$ and $a'b'c'$ in $X$ with $b'\in \Aline ac$. Observe that the flat hull $\overline P$ of the plane $P$ is a plane in the projective liner $Y$. Since the projective liner $Y$ is Moufang, the set $T\defeq(\Aline ab\cap \Aline {a'}{b'})\cup(\Aline bc\cap\Aline{b'}{c'})\cup(\Aline ac\cap\Aline{a'}{c'})$ has rank $\|T\|\in\{0,2\}$ in the liner $Y$. Then the set $T\cap X$ has rank $\|T\cap X\|\le 2$ in the liner $X$. It remains to prove that $\|T\cap X\|\ne 1$.

 Repeating the argument from the proof of Theorem~\ref{t:Desargues-projective}, we can show that the set $T$ has cardinality $|T|=3$. 
Then $\overline{T}$ is a line in the projective liner $Y$. If $|T\cap H|\ge 2$, then $T\subseteq\overline{T}\subseteq H$ because $F$ is a flat. In this case $|T\cap X|=0$. If $|T\cap H|\le 1$, then $|T\cap X|\ge|T|-1=2$ and we are done. 
\end{proof}

\begin{proposition}\label{p:nD=>pP} Let $Y$ be a projective completion of a $3$-long liner $X$. If $X$ is Moufang, then the horizon $Y\setminus X$ of $X$ in $Y$ is flat and hence the liner $X$ is completely regular.
\end{proposition}

\begin{proof}  To prove that the horizon $H\defeq Y\setminus X$ of $X$ in $Y$  is flat, take any points $x,y\in H$ and a point $z\in \Aline xy$. If $z\in\{x,y\}$, then $z\in H$ and we are done. So, assume that $z\notin\{x,y\}$, which implies $x\ne z\ne y$. Since $Y$ is a completion of $X$, there exists a point $a\in Y\setminus\overline H\subseteq X$. It follows from $x\in H$ and $a\notin\overline H$ that $\Aline ax\cap\overline H=\{x\}$. Since $Y$ is $3$-long, the intersection $\Aline ax\cap X$ contains more than one point and hence is a line in $X$. Since the liner $X$ is $3$-long, $|\Aline ax|=|\Aline ax\cap X|+1\ge 3+1=4$. By Corollary~\ref{c:Avogadro-projective}, the $3$-long projective liner $Y$ is $2$-balanced and hence $|Y|_2=|\Aline ax|\ge 4$. Therefore, the projective liner $Y$ is $4$-long.

\begin{picture}(120,150)(-120,-15)
\put(0,0){\color{cyan}\line(1,0){120}}
\put(120,0){\line(0,1){120}}
\put(0,120){\line(1,0){120}}
\put(60,0){\line(1,6){20}}
\put(60,0){\line(1,2){60}}
\put(120,0){\line(-1,3){40}}
\put(0,0){\line(2,1){120}}
\put(0,0){\line(4,3){96}}
\put(0,120){\line(2,-1){120}}
\put(0,120){\line(1,-1){80}}

\put(0,0){\color{cyan}\circle*{3}}
\put(-2,-8){$y$}
\put(60,0){\color{red}\circle*{3}}
\put(58,-8){$z$}
\put(120,0){\color{cyan}\circle*{3}}
\put(118,-8){$x$}
\put(0,120){\circle*{3}}
\put(-2,122){$o$}
\put(120,120){\circle*{3}}
\put(120,123){$a$}
\put(120,60){\circle*{3}}
\put(123,57){$b$}
\put(68.6,51.4){\circle*{3}}
\put(72,48){$c'$}
\put(80,40){\circle*{3}}
\put(80,33){$c$}
\put(96,72){\circle*{3}}
\put(100,71){$b'$}
\put(80,120){\circle*{3}}
\put(78,123){$a'$}

\end{picture}

Consider the plane $P\defeq\overline{\{a,x,y\}}$ in the projective liner $Y$. By Theorem~\ref{t:procompletion=>normal}, $P\cap X$ is a plane in the $3$-long liner $X$. The $3$-rankedness of the (strongly) regular projective liner $Y$ ensures that $P\cap \overline H=\Aline xy$. By Proposition~\ref{p:cov-aff}, there exists a point $o\in P\setminus(\Aline xy\cup\Aline ax\cup\Aline az)$. Choose any point $b\in\Aline ax\setminus\{a,x\}$.  Since the plane $P$ is projective, there exist unique points $c\in \Aline az\cap\Aline yb$, $b'\in \Aline az\cap\Aline ob$, $a'\in \Aline oa\cap\Aline x{b'}$, and $c'\in \Aline oc\cap\Aline y{b'}$. The choice of the points $a,b,c,a',b',c'$ guarantees that $abc$ and $a'b'c'$ are two disjoint centrally perspective triangles in the plane $P\cap X$ such that $b'\in \Aline ac$, $\Aline ab\cap\Aline {a'}{b'}=\{x\}$, $\Aline bc\cap\Aline {b'}{c'}=\{y\}$ and $z\in \Aline ac$. Since the liner $X$ is Moufang, the set 
$$T\defeq (\Aline ab\cap\Aline {a'}{b'}\cap X)\cup(\Aline bc\cap\Aline {b'}{c'}\cap X)\cup(\Aline ac\cap\Aline{a'}{c'}\cap X)$$has rank $\|T\|\in\{0,2\}$ in the liner $X$. Since $x,y\notin X$, the set $T$ is equal to the set $\Aline ac\cap\Aline{a'}{c'}\cap X$ which contains at most one point. Then $\|T\|=0$ and hence $\Aline ac\cap\Aline{a'}{c'}\cap X=\varnothing$. Since the liner $Y$ is projective, the set $\Aline ac\cap\Aline {a'}{c'}$ is not empty. Then $\varnothing\ne\Aline ac\cap\Aline {a'}{c'}\subseteq \Aline ac\setminus X\subseteq \{z\}$ and hence $z\in Y\setminus X=H$, witnessing that the horizon $H$ of $X$ is flat in $Y$. By Theorem~\ref{t:spread=projective2}, the liner $X=Y\setminus H$ is completely regular.
\end{proof}

\begin{theorem}\label{t:Moufang=>compreg} Every Moufang proaffine regular liner $X$ is completely regular.
\end{theorem}

\begin{proof} Let $X$ be a Desarguesian proaffine regular liner. If $X$ is projective or has rank $\|X\|\ne 3$, then $X$ is completely regular, by Corollary~\ref{c:proregular=>ranked} and Theorem~\ref{t:proaffine3=>compregular}. So, assume that the liner $X$ is not projective and $\|X\|=3$. By Theorem~\ref{t:Proclus<=>}, the proaffine regular liner $X$ is Proclus. If $X$ is not $3$-long, then by Theorem~\ref{t:Proclus-not-3long}, $X=Y\setminus \{p\}$ for some Steiner projective plane $Y$ and some point $p\in Y$. Then $Y$ is a projective completion of $X$ with with flat horizon $H=\{p\}$. It is easy to see that the liner $X$ is $3$-ranked and its spread completion is isomorphic to $Y$ and hence is projective. So, we assume that the liner $X$ is $3$-long. If $X$ is finite, then the proaffine regular liner $X$ has a projective completion $Y$, by Theorem~\ref{t:procompletion-finite}. By Corollary~\ref{c:Avogadro-projective}, the $3$-long projective liner $Y$ is $2$-balanced. Since $X\ne Y$ and the liner $X$ is $3$-long, the projective liner $Y$ is $4$-long. By Proposition~\ref{p:nD=>pP}, the horizon $Y\setminus X$ is flat in $Y$. Applying Theorem~\ref{t:spread=projective2}, we conclude that $X$ is completely regular.

So, assume that the Moufang Proclus plane $X$ is infinite.  We divide the  further proof in a series of claims and lemmas.

\begin{claim}\label{cl:w-long} The liner $X$ is $\w$-long.
\end{claim}

\begin{proof} If $X$ is not $\w$-long, then $X$ contains a finite line $L$. By Proposition~\ref{p:Avogadro-proaffine}, every line in $X$ has cardinality $\le|L|+1$ and hence is finite. Since $X$ is a plane, there exists a line $\Lambda\subseteq X$, concurrent with $L$. The regularity and rankedness of the plane $X$ ensures that $X=\overline{L\cup\Lambda}=\bigcup_{x\in L}\bigcup_{y\in\Lambda}\Aline xy$ and hence
$$|X|\le\sum_{x\in L}\sum_{y\in\Lambda}|\Aline xy|\le |L|\cdot|\Lambda|\cdot (|L|+1)<\w,$$which contradicts our assumption. This contradiction shows that the liner $X$ is $\w$-long.
\end{proof}

\begin{claim}
For any disjoint lines $L,L'\subseteq X$, the set $$I\defeq \{x\in X\setminus  L':\forall y\in X\setminus\{x\}\;\;(\Aline xy\cap L'\ne\varnothing)\}$$ contains at most one point.
\end{claim}

\begin{proof} To derive a contradiction, assume that the set $I$ contains two distinct points $a$ and $\alpha$. Since $L\cap L'=\varnothing$, the definition of the set $I$ ensures that the points $a$ and $\alpha$ do not belong to the line $L$. Since $a,\alpha\in I$, the line $\Aline a\alpha$ intesects the line $L'$ and also the line $L$, which is parallel to the line $L$, by Corollary~\ref{c:parallel-lines<=>}. Then there exist unique points $x\in \Aline a\alpha\cap L$ and $b\in \Aline a\alpha\cap L'$. Choose any point $a'\in L'\setminus\{b\}$. Let $\mathcal L_b$ be the family of lines containing the point $b$. Since the liner $X$ is Proclus, the set $\mathcal L_b'\defeq\{L\in\mathcal L_b:L\cap \Aline a{a'}=\varnothing\}\cup\{L\in\mathcal L_b:L\cap\Aline \alpha{a'}=\varnothing\}$ contains a most two elements. Choose any point $b'\in \Aline x{a'}\setminus\{x,a'\}$ such that $\Aline o{b'}\notin \mathcal L_b'$. Such a choice of the point $b'$ guarantees the existence of points $o\in \Aline b{b'}\cap \Aline a{a'}$ and $o'\in \Aline b{b'}\cap\Aline\alpha{a'}$. It follows from $a\ne\alpha$ that $o\ne o'$. Choose any point $y\in L\setminus(\{x\}\cup \Aline b{b'})$. By Proclus Postulate~\ref{p:Proclus-Postulate}, there exists a point $c'\in L'\cap\Aline y{b'}$. Since $c'\notin \Aline b{b'}$, the lines $\Aline o{c'}$ and $\Aline {o'}{c'}$  
are distinct. Since $X$ is Proclus, one of these two lines intersects the line $by$. We lose no generality assuming that $\Aline o{c'}\cap \Aline by$ contains some point $c$. 

\begin{picture}(240,145)(-200,-15)
\put(-70,25){\line(6,7){70}}
\linethickness{=0.6pt}
\put(0,0){\line(0,1){120}}
\put(0,0){\line(2,3){40}}
\put(0,0){\line(-2,3){40}}
\put(-70,45){\color{teal}\line(1,0){140}}
\put(-70,60){\color{teal}\line(1,0){140}}
\put(-60,30){\color{teal}\line(1,0){120}}
\put(-70,25){\line(2,1){70}}
\put(60,30){\line(-2,1){60}}
\put(60,30){\line(-2,3){60}}
\put(-60,30){\line(2,3){60}}

\put(0,106.7){\circle*{2.5}}
\put(-8,104){$o'$}
\put(-70,25){\circle*{3}}
\put(-73,17){$\alpha$}
\put(0,0){\circle*{3}}
\put(-2,-9){$b'$}
\put(-30,45){\circle*{3}}
\put(-34,37){$x$}
\put(30,45){\circle*{3}}
\put(28,37){$y$}
\put(-60,30){\circle*{3}}
\put(-62,22){$a$}
\put(60,30){\circle*{3}}
\put(60,22){$c$}
\put(-40,60){\circle*{3}}
\put(-45,62){$a'$}
\put(40,60){\circle*{3}}
\put(40,63){$c'$}
\put(0,60){\circle*{3}}
\put(1,62){$b$}
\put(0,120){\circle*{3}}
\put(-2,123){$o$}
\end{picture}

Then $abc$ and $a'b'c'$ are two disjoint centrally perspective triangles such that $b\in \Aline {a'}{c'}$, $\Aline ab\cap\Aline{a'}{b'}=\{x\}$ and $\Aline bc\cap\Aline {b'}{c'}=\{y\}$. Since $X$ is Moufang, the nonempty set 
$$T\defeq(\Aline ab\cap\Aline{a'}{b'})\cup(\Aline bc\cap\Aline{b'}{c'})\cup(\Aline ac\cap\Aline{a'}{c'})=\{x,y\}\cup(\Aline ac\cap L')$$has rank $\|T\|=2$ and hence 
$\Aline ac\cap L'\subseteq \Aline xy\cap L'=L\cap L'=\varnothing$, which contradicts the inclusion $a\in I$.
\end{proof}

\begin{lemma}\label{l:Moufang=>para-Playfair} The Moufang Proclus plane $X$ is para-Playfair.
\end{lemma}

\begin{proof} To derive a contradiction, assume that $X$ is not para-Playfair. Then there exist two disjoint lines $L,L'$ in $X$ and a point $a\in X\setminus (L\cup L')$ such that for every point $x\in X\setminus\{a\}$, the line $\Aline ax$ intersects the parallel lines $L,L'$. 

\begin{claim}\label{cl:near-Des1} For every points $b\in L'$ and $y\in L\setminus\Aline ab$, the set
$$I\defeq\{z\in L'\setminus \{b\}:\forall x\in X\setminus\{z\}\;\;(\Aline xz\cap \Aline by\ne\varnothing)\}$$contains at most one point. 
\end{claim}

\begin{proof} To derive a contradiction, assume that the set $I$ contains two distinct points $c'$ and $c''$. Since $\Aline ab\cap L=\{b\}$ and $L\parallel L'$, we can apply  the Proclus Postulate~\ref{p:Proclus-Postulate} and find a point $x\in \Aline ab\cap L$.  The choice of the point $y\notin \Aline ab$ ensures that $x\ne y$. 
Since the liner $X$ is Proclus and  $\w$-long, there exists a line $\Lambda$ in $X$ such that $b\in \Lambda\ne L'$, $\Aline xb\ne \Lambda\ne \Aline yb$ and $\Aline y{c'}\nparallel \Lambda\nparallel \Aline y{c''}$. Then there exist points $b'\in \Aline y{c'}\cap\Lambda$ and $b''\in \Aline y{c''}\cap\Lambda$ and those points are distinct because $c'\ne c''$. Since $\Lambda\cap L\ne\{y\}$, the points $b',b''$ do not belong to the line $L$. By Proclus Postulate~\ref{p:Proclus-Postulate}, there exist points $a'\in L'\cap\Aline x{b'}$ and $a''\in L'\cap\Aline x{b''}$. Those points are distinct because $b'\ne b''$ and $x\notin\Lambda$. Also $a'\ne b\ne a''$ because $x\notin\Lambda$.  Since the plane $X$ is Proclus, one of the lines $\Aline a{a'}$ or $\Aline a{a''}$ intersects the line $\Lambda$. We lose no generality assuming that $\Aline a{a'}\cap\Lambda$ is not empty and hence contains some point $o$.
The inequality $a'\ne b$ implies that $o\notin L'$. Since $c'\in I$, there exists a point $c\in \Aline o{c'}\cap\Aline by$.

\begin{picture}(240,150)(-200,-15)
\put(55,60){\line(-5,-3){55}}
\put(-55,60){\line(5,-3){55}}
\put(-60,30){\line(1,6){15}}

\linethickness{=0.6pt}
\put(0,0){\color{cyan}\line(0,1){120}}
\put(2,85){\color{cyan}$\Lambda$}
\put(0,0){\line(2,3){40}}
\put(0,0){\line(-2,3){40}}
\put(-70,45){\color{teal}\line(1,0){140}}
\put(75,42){\color{teal}$L$}
\put(-70,60){\color{teal}\line(1,0){140}}
\put(75,58){\color{teal}$L'$}
\put(-60,30){\color{teal}\line(1,0){120}}
\put(-60,30){\line(2,1){60}}
\put(60,30){\line(-2,1){60}}
\put(60,30){\line(-2,3){60}}
\put(-60,30){\line(2,3){60}}

\put(0,27){\circle*{2.5}}
\put(2,20){$b''$}
\put(0,0){\circle*{3}}
\put(-3,-9){$b'$}
\put(-30,45){\circle*{3}}
\put(-34,37){$x$}
\put(30,45){\circle*{3}}
\put(28,37){$y$}
\put(-60,30){\circle*{3}}
\put(-62,22){$a$}
\put(60,30){\circle*{3}}
\put(60,22){$c$}
\put(-40,60){\circle*{3}}
\put(-45,62){$a'$}
\put(40,60){\circle*{3}}
\put(40,63){$c'$}
\put(55,60){\circle*{2.5}}
\put(53,63){$c''$}
\put(-55,60){\circle*{2.5}}
\put(-65,63){$a''$}

\put(0,60){\circle*{3}}
\put(1,62){$b$}
\put(0,120){\circle*{3}}
\put(-2,123){$o$}
\end{picture}

Then $abc$ and $a'b'c'$ are two centrally perspective triangles in $X$ such that $b\in L'=\Aline {a'}{c'}$, $\Aline ab\cap\Aline {a'}{b'}=\{x\}$ and $\Aline bc\cap\Aline{b'}{c'}=\{y\}$. Since $X$ is Moufang, the nonempty set $$T\defeq(\Aline ab\cap\Aline {a'}{b'})\cup(\Aline bc\cap\Aline{b'}{c'})\cup(\Aline ac\cap\Aline {a'}{c'})=\{x,y\}\cup(\Aline ac\cap L')$$has rank $\|T\|=2$, which implies $\Aline ac\cap L'\subseteq \Aline xy\cap L'=L\cap L'=\varnothing$. But the equality $\Aline ac\cap L'=\varnothing$ contradicts the choice of the point $a$.
\end{proof}

\begin{claim}\label{cl:near-Des2} For every point $b\in L'$, the set
$$J\defeq\{z\in L'\setminus \{b\}:\forall x\in X\setminus\{z\}\;\;(\Aline xz\cap \Aline ab\ne\varnothing)\}$$contains at most one point. 
\end{claim}

\begin{proof} To derive a contradiction, assume that the set $J$ contains two distinct points $a'$ and $a''$. Choose any point $\alpha\in \Aline ab\setminus(\{a,b\}\cup L)$. By Claim~\ref{cl:near-Des1}, the liner $X$ contains a line $L''$ such that $\alpha\in L''$ and $L''\parallel L'$. Since $X$ is Proclus and $\w$-long, there exists a line $\Lambda$ in $X$ such that $b\in\Lambda$, $L'\ne\Lambda\ne\Aline ab$ and $\Aline \alpha {a'}\nparallel \Lambda\nparallel\Aline \alpha{a''}$. By the choice of $\Lambda$, there exist points $o'\in \Aline \alpha{a'}\cap\Lambda$ and $o''\in \Aline \alpha{a''}\cap\Lambda$. It follows from $a'\ne b\ne a''$ and $\Lambda\ne\Aline ab$ that the points $o',o''$ are distinct and do not belong to $L'\cup L''$. Choose any point $\gamma\in L''\setminus(\{\alpha\}\cup\Lambda)$ and find a unique point $y\in\Aline b\gamma\cap L$. Since $o',o''\notin L''$, there exist unique points $c'\in \Aline{o'}\gamma\cap L'$ and $c''\in \Aline{o''}\gamma\cap L''$. It follows from $o'\ne o''$ and $\gamma\notin\Lambda$ that $c'\ne c''$. By Claim~\ref{cl:near-Des1}, the set $I\defeq\{c\in L':\forall z\in X\setminus\{c\}\;\;(\Aline cz\cap \Aline by\ne\varnothing)\}$ contains at most one point. Then either $c'\notin I$ or $c''\notin I$. We lose no generality assuming that $c'\notin I$. Then there exist $z\in X\setminus \{c'\}$ such that $\Aline z{c'}\cap\Aline by=\varnothing$. Since $\Aline z{c'}\parallel \Aline by$ and $\Aline by\cap\Lambda=\{b\}$, there exists a unique point $b'\in \Aline z{c'}\cap\Lambda$. Then $\alpha b\gamma$ and $a'b'c'$ are two centrally perspective triangles such that $b\in L'=\Aline {a'}{c'}$. Since $X$ is Moufang, the set
$$T\defeq(\Aline \alpha b\cap\Aline{a'}{b'})\cup(\Aline b\gamma\cap\Aline {b'}{c'})\cup(\Aline \alpha\gamma\cap\Aline {a'}{c'})=(\Aline \alpha b\cap\Aline{a'}{b'})\cup \varnothing\cup \varnothing=\Aline \alpha b\cap\Aline{a'}{b'}$$
has rank $\|T\|\in\{0,2\}$ and hence $\Aline ab\cap\Aline {a'}{b'}=\Aline \alpha b\cap\Aline{a'}{b'}=\varnothing$, which contradicts the inclusion $a'\in J$.

\begin{picture}(240,130)(-200,-15)
\linethickness{=0.6pt}
\put(0,15){\color{cyan}\line(0,1){87.3}}
\put(-75,0){\line(11,15){75}}
\put(75,0){\line(-11,15){75}}
\put(-75,0){\line(1,1){75}}
\put(75,0){\line(-1,1){75}}
\put(-75,0){\color{teal}\line(1,0){150}}
\put(-75,60){\color{teal}\line(1,0){150}}
\put(80,57){\color{teal}$L$}
\put(-75,75){\color{teal}\line(1,0){150}}
\put(80,72){\color{teal}$L'$}
\put(0,15){\line(-1,3){20}}
\put(0,15){\line(1,3){20}}

\put(-45,30){\circle*{3}}
\put(-44,24){$a$}
\put(0,75){\circle*{3}}
\put(2,77){$b$}
\put(-75,0){\circle*{3}}
\put(-78,-8){$\alpha$}
\put(75,0){\circle*{3}}
\put(73,-8){$\gamma$}
\put(0,15){\circle*{3}}
\put(-3,5){$b'$}
\put(-15,60){\circle*{3}}
\put(-12,54){$x$}
\put(15,60){\circle*{3}}
\put(7,54){$y$}
\put(-20,75){\circle*{3}}
\put(-25,78){$a'$}
\put(20,75){\circle*{3}}
\put(20,78){$c'$}
\put(0,102.3){\circle*{3}}
\put(-2,105){$o'$}
\end{picture}

\end{proof}

Now we are able to complete the proof of Lemma~\ref{l:Moufang=>para-Playfair}. Choose any point $b\in L'$. By Claim~\ref{cl:near-Des2}, the set $J\defeq\{x\in L'\setminus\{b\}:\forall y\in X\setminus\{x\}\;\;(\Aline xy\cap \Aline ab\ne\varnothing\}$ contains at most one point. Since $X$ is $\w$-long, there exists a point $a'\in L'\setminus(\{b\}\cup I)$. By the definition of the set $I$, there exists a line $L_{a'}$ such that $a'\in L_{a'}\subseteq X\setminus\Aline ab$. Take any point $y\in L\setminus \Aline ab$. By Claim~\ref{cl:near-Des1}, the set  $I\defeq\{y\in L'\setminus\{b\}:\forall z\in X\setminus\{y\}\;\;(\Aline yz\cap \Aline by\ne\varnothing\}$ contains at most one point. Since $X$ is $\w$-long, there exist two distinct points $c',c''\in L'\setminus(\{b,a'\}\cup J)$. By the definition of the set $J$, there exist lines $L_{c'},L_{c''}$ such that $c'\in L_{c'}$, $c''\in L_{c''}$ and $L_{c'}\cap\Aline by=\varnothing=L_{c''}\cap\Aline by$. Since $\Aline yb\cap\Aline ab=\{b\}$, the Proclus Postulate~\ref{p:Proclus-Postulate} ensures that there exist points $b'\in L_{a'}\cap L_{c'}$ and $b''\in L_{a'}\cap L_{c''}$. Since the points $c',c''$ are distinct, so are the points $b',b''$. Then the lines $\Aline b{b'}$ and $\Aline b{b''}$ are distinct and at least one of them intersects the line $\Aline a{a'}$ (because the plane $X$ is Proclus). We lose no generality assuming that $\Aline a{a'}\cap\Aline b{b'}$ is not empty and hence contains some point $o$. Since $\Aline o{c'}\cap \Aline {c'}{b'}=\{c'\}$ and $\Aline {c'}{b'}\parallel \Aline by$, there exists a unique point $c\in \Aline by\cap\Aline o{c'}$. 

\begin{picture}(100,125)(-140,-15)
\put(100,40){\color{violet}\line(-1,1){50}}
\put(-10,30){\color{blue}\line(1,1){60}}
\put(-22,23){\color{blue}$L_{a'}$}
\put(40,40){\color{violet}\line(-1,1){30}}
\put(0,40){\color{teal}\line(1,0){120}}
\put(0,70){\color{teal}\line(1,0){120}}
\put(40,0){\color{violet}\line(-1,1){40}}
\put(40,0){\color{blue}\line(1,1){40}}
\put(40,0){\line(0,1){80}}
\put(40,80){\color{violet}\line(1,-1){60}}
\put(100,10){\color{violet}$L_{c'}$}
\put(50,90){\color{violet}\line(1,-1){65}}
\put(117,17){\color{violet}$L_{c''}$}

{\linethickness{=1pt}
\put(0,40){\color{teal}\line(1,0){80}}
\put(0,40){\color{blue}\line(1,1){40}}
\put(20,20){\color{teal}\line(1,0){40}}
\put(20,20){\color{blue}\line(1,1){20}}
\put(40,40){\color{violet}\line(1,-1){20}}
\put(80,40){\color{violet}\line(-1,1){40}}
}

\put(125,37){\color{teal}$L'$}
\put(125,67){\color{teal}$L$}

\put(40,0){\circle*{3}}
\put(38,-8){$o$}
\put(0,40){\circle*{3}}
\put(-8,40){$a'$}
\put(40,40){\circle*{3}}
\put(41,42){$b$}
\put(80,40){\circle*{3}}
\put(79,43){$c'$}
\put(100,40){\circle*{2.5}}
\put(99,43){$c''$}
\put(40,0){\circle*{3}}
\put(38,-8){$o$}
\put(20,20){\circle*{3}}
\put(13,13){$a$}
\put(60,20){\circle*{3}}
\put(61,13){$c$}
\put(10,70){\circle*{2.5}}
\put(8,74){$y$}
\put(40,80){\circle*{3}}
\put(37,83){$b'$}
\put(50,90){\circle*{2.5}}
\put(48,93){$b''$}
\end{picture}

Then $abc$ and $a'b'c'$ are two centrally perspective triangles in $X$ with $b\in L'=\Aline {a'}{c'}$. Since the liner $X$ is Moufang, the set $$T\defeq(\Aline ab\cap\Aline{a'}{b'})\cup(\Aline bc\cap\Aline{b'}{c'})\cup(\Aline ac\cap \Aline {a'}{c'})=\varnothing\cup \varnothing\cup(\Aline ac\cap L')$$has rank  
$\|T\|\in\{0,2\}$ and hence is empty. On the other hand the equality $\Aline ac\cap L'=\varnothing$ contradicts the choice of the point $a$.
\end{proof}

\begin{lemma}\label{l:near-Des=>bi-Bolyai} The Moufang Proclus plane $X$ is bi-Bolyai.
\end{lemma}

\begin{proof} Given two concurrent Bolyai lines $A,B$, we need to show that every line $L$ in the plane $\overline{A\cup B}=X$ is Bolyai. This is clear if $L\in A_\parallel\cup B_\parallel$. So, assume that $L\notin A_\parallel\cup B_\parallel$. Given any point $a\in X\setminus L$, we need to find a line $\Lambda$ such that $a\in \Lambda\subseteq X\setminus L$. Replacing the Bolyai line $A$ by a suitable parallel line, we lose no generality assuming that $a\in A$. Since the lines $L$ and $A$ are consurrent, there exist a unique point $b\in A\cap L$. Replacing the Bolyai line $B$ by a suitable parallel line, we lose no generality assuming that $b\in B$. Fix any point $a'\in L\setminus\{b\}$ and any distinct points $c',c''\in L\setminus\{a',b\}$. Since the lines $A,B$ are Bolyai, there exist lines $A'\in A_\parallel$ and $B',B''\in B_\parallel$ such that $a'\in A'$, $c'\in B'$, $c''\in B''$. By the Proclus Postulate~\ref{p:Proclus-Postulate}, there exist points $b'\in A'\cap B'$ and $b''\in A'\cap B''$. Since the points $c'$ and $c''$ are distinct, the parallel lines $B',B''$ are disjoint and hence $b'\ne b''$.   
 Then the lines $\Aline b{b'}$ and $\Aline b{b''}$ are distinct and at least one of them intersects the line $\Aline a{a'}$ (because the plane $X$ is Proclus). We lose no generality assuming that $\Aline a{a'}\cap\Aline b{b'}$ is not empty and hence contains some point $o$. Since $\Aline o{c'}\cap \Aline {c'}{b'}=\{c'\}$ and $\Aline {c'}{b'}=B'\parallel B$, there exists a unique point $c\in B\cap\Aline o{c'}$. 

\begin{picture}(100,125)(-140,-15)
\put(100,40){\color{violet}\line(-1,1){50}}
\put(-10,30){\color{blue}\line(1,1){60}}
\put(-22,23){\color{blue}$A'$}
\put(40,40){\color{blue}\line(-1,-1){32}}
\put(-2,1){\color{blue}$A$}
\put(40,40){\color{violet}\line(1,-1){32}}
\put(72,-1){\color{violet}$B$}
\put(0,40){\color{red}\line(1,0){120}}
\put(40,0){\color{violet}\line(-1,1){40}}
\put(40,0){\color{blue}\line(1,1){40}}
\put(40,0){\line(0,1){80}}
\put(40,80){\color{violet}\line(1,-1){60}}
\put(100,10){\color{violet}$B'$}
\put(50,90){\color{violet}\line(1,-1){62}}
\put(114,19){\color{violet}$B''$}

{\linethickness{=1pt}
\put(0,40){\color{red}\line(1,0){80}}
\put(0,40){\color{blue}\line(1,1){40}}
\put(20,20){\color{red}\line(1,0){40}}
\put(20,20){\color{blue}\line(1,1){20}}
\put(40,40){\color{violet}\line(1,-1){20}}
\put(80,40){\color{violet}\line(-1,1){40}}
}

\put(125,37){\color{red}$L$}

\put(40,0){\circle*{3}}
\put(38,-8){$o$}
\put(0,40){\circle*{3}}
\put(-8,40){$a'$}
\put(40,40){\circle*{3}}
\put(41,42){$b$}
\put(80,40){\circle*{3}}
\put(79,43){$c'$}
\put(100,40){\circle*{2.5}}
\put(99,43){$c''$}
\put(40,0){\circle*{3}}
\put(38,-8){$o$}
\put(20,20){\circle*{3}}
\put(11,18){$a$}
\put(60,20){\circle*{3}}
\put(63,18){$c$}
\put(40,80){\circle*{3}}
\put(37,83){$b'$}
\put(50,90){\circle*{2.5}}
\put(48,93){$b''$}
\end{picture}

Then $abc$ and $a'b'c'$ are two centrally perspective triangles in $X$ with $b\in L'=\Aline {a'}{c'}$. Since the liner $X$ is Moufang, the set $$T\defeq(\Aline ab\cap\Aline{a'}{b'})\cup(\Aline bc\cap\Aline{b'}{c'})\cup(\Aline ac\cap \Aline {a'}{c'})=\varnothing\cup \varnothing\cup(\Aline ac\cap L)$$has rank  
$\|T\|\in\{0,2\}$ and hence is empty. Then $\Aline ac$ is a required line that contains the poit $a$ and does not intersect the line $L$, witnessing that the line $L$ is Bolyai and the liner $X$ is bi-Bolyai.
\end{proof}

By Lemmas~\ref{l:Moufang=>para-Playfair} and \ref{l:near-Des=>bi-Bolyai}, the $\w$-long regular liner $X$ is para-Playfair and bi-Bolyai. By Theorem~\ref{t:spread=projective1}, the liner $X$ is completely regular.
\end{proof}

Proposition~\ref{p:Desarg=>Moufang} and Theorem~\ref{t:Moufang=>compreg} imply the following important corollary.

\begin{corollary}\label{c:Desarg=>compreg} Every Desarguesian regular liner in completely regular.
\end{corollary}




\section{Inverse Desargues Theorems}

In this section we shall prove the inverse Desargues Theorems~\ref{t:ID1}, \ref{t:ID2}, \ref{t:ID3}. First, we prove a counterpart of Lemma~\ref{l:ADA+P} for Desarguesian Proclus planes.

\begin{theorem}\label{t:pD=>PP} Let $abc$ and $a'b'c'$ be two centrally perspective disjoint triangles in a Desarguesian Proclus plane $X$. If $\Aline ac\cap\Aline{a'}{c'}=\varnothing$ and $\Aline ab\cap\Aline {a'}{b'}=\{x\}$ for some point $x\in X$, then there exists a point $y\in X\setminus\{x\}$ such that $\Aline bc\cap\Aline{b'}{c'}=\{y\}$ and $\Aline xy\cap \Aline ac=\varnothing=\Aline xy\cap\Aline{a'}{c'}$.
\end{theorem}

\begin{proof} By Theorem~\ref{t:Proclus<=>}, the Proclus plane $X$ is $3$-proregular, and by Proposition~\ref{p:k-regular<=>2ex}, the $3$-proregular liner $X$ is $3$-ranked.  Since the liner $X$ is Desarguesian, the nonempty set
$$(\Aline ab\cap \Aline {a'}{b'})\cup(\Aline bc\cap\Aline{b'}{c'})\cup(\Aline ac\cap\Aline{a'}{c'})=\{x\}\cup(\Aline bc\cap\Aline {b'}{c'})\cup\varnothing$$ has rank $2$ and hence $\Aline bc\cap\Aline {b'}{c'}$ contains some point $y$. Assuming that $\Aline bc\cap\Aline {b'}{c'}\ne\{y\}$, we conclude that $\Aline bc=\Aline {b'}{c'}$ and hence $\{c\}=\Aline c{c'}\cap \Aline bc=\Aline c{c'}\cap\Aline {b'}{c'}=\{c'\}$, which contradicts the choice of the disjoint triangles $abc,a'b'c'$. This contradiction shows that $\Aline bc\cap\Aline {b'}{c'}=\{y\}$. Assuming that $x=y$, we conclude that $\Aline ab=\Aline xb=\Aline yb=\Aline bc$ and then the points $a,b,c$ are collinear and cannot form a triangle. This contradiction shows that $x\ne y$ and hence $\Aline xy$ is a line. 

It remains to prove that $\Aline xy\cap \Aline ac=\varnothing=\Aline xy\cap\Aline{a'}{c'}$. 
Let $o\in \Aline a{a'}\cap\Aline b{b'}\cap\Aline c{c'}$ be the perspector of the centrally perspective triangles $abc$ and $a'b'c'$. Observe that the plane $X$ contains two disjoint lines $\Aline ac$ and $\Aline{a'}{c'}$ and hence $X$ is not projective. If $X$ is not $3$-long, then by Theorem~\ref{t:Proclus-not-3long}, $X$ is a proper subliner of some Steiner projective plane and hence $|X|<7$, which is not possible because $|X|\ge |\{o,a,b,c,a',b',c'\}|=7$. This contradiction shows that the Proclus plane $X$ is $3$-long and hence the $3$-proregular plane $X$ is $3$-regular and regular.

By Corollary~\ref{c:Desarg=>compreg}, the Desarguesian $3$-long Proclus plane $X$ is completely regular, and by  Theorem~\ref{t:spread=projective1}, the completely regular liner $X$ is para-Playfair and bi-Bolyai. Two cases are possible.
\smallskip

1. If the plane $X$ is Playfair, then $\Aline xy\parallel \Aline ac$, by Lemma~\ref{l:ADA+P}. By Corollary~\ref{c:parallel-lines<=>}, the disjoint lines $\Aline ac$ and $\Aline{a'}{c'}$ in the plane $X$ are parallel and hence $\Aline xy\parallel \Aline{a'}{c'}$. If $\Aline xy\cap\Aline ac\ne\varnothing$, then $\Aline xy\parallel \Aline ac$ implies $\Aline xy=\Aline ac$ and hence $b\in \Aline xa\subseteq \Aline ac$, which is impossible becuase $abc$ is a triangle. This contradiction shows that $\Aline xy\cap\Aline ac=\varnothing$. By analogy we can prove that $\Aline xy\cap\Aline{a'}{c'}=\varnothing$. 
\smallskip

2. Next, assume that the plane $X$ is not Playfair. Since the liner $X$ is para-Playfair, the disjoint lines $\Aline ac$ and $\Aline {a'}{c'}$ are Bolyai in $X$.  Since $X$ is bi-Bolyai and not Playfair,  any Bolyai line in the plane $X$ is parallel to the lines $\Aline ac,\Aline{a'}{c'}$. Since the liner $X$ is para-Playfair and regular, any two disjoint lines in $X$ are Bolyai (by Proposition~\ref{p:lines-in-para-Playfair}) and hence any disjoint lines in $X$ are parallel to the line $\Aline ac$. We use this fact to derive a contradiction from the assumption $\Aline ac\cap\Aline xy\ne\varnothing$. Assume that the lines $\Aline xy$ and $\Aline ac$ have a common point $z$.  It is easy to show that $z\notin\{a,c\}$.

\begin{picture}(300,145)(-230,-50)

\put(-105,75){\color{red}\line(1,0){210}}
\put(0,0){\line(1,0){30}}
\put(0,0){\line(-1,0){30}}
\put(0,30){\line(1,0){15}}
\put(0,30){\line(-1,0){200}}
\put(-30,0){\line(1,2){30}}
\put(30,0){\line(-1,2){30}}
\put(0,0){\line(0,1){60}}
\put(0,0){\line(0,-1){30}}
\put(0,22.5){\line(2,1){105}}
\put(0,-30){\line(1,1){105}}
\put(0,22.5){\line(-2,1){105}}
\put(0,-30){\line(-1,1){105}}
{\color{red}\qbezier(-105,75)(-150,75)(-200,30)}
\put(0,-30){\line(-33,47){57}}

\put(-200,30){\line(23,-3){230}}

\put(0,-30){\circle*{3}}
\put(-2,-40){$b'$}
\put(-26.4,7.2){\circle*{3}}
\put(-23,8){$a''$}
\put(0,22.5){\circle*{3}}
\put(2,15){$b$}
\put(0,60){\circle*{3}}
\put(-2,63){$o$}
\put(30,0){\circle*{3}}
\put(31,-9){$c'$}
\put(15,30){\circle*{3}}
\put(15,34){$c$}
\put(-30,0){\circle*{3}}
\put(-39,-10){$a'$}
\put(-15,30){\circle*{3}}
\put(-20,34){$a$}
\put(-200,30){\circle*{3}}
\put(-210,28){$z$}
\put(-57,51){\circle*{3}}
\put(-59,54){$x'$}

\put(105,75){\circle*{3}}
\put(103,79){$y$}
\put(-105,75){\circle*{3}}
\put(-108,78){$x$}
\end{picture}

Since the lines $\Aline oa$ and $\Aline z{c'}$ are not parallel to the line $\Aline ac$, they are concurrent and hence there exists a unique point $a''\in\Aline oa\cap\Aline z{c'}$. It follows from $z\in \Aline ac\setminus\{a\}=\Aline ac\setminus \Aline oa$ and $c'\notin \Aline ac$ that $a''\ne a$. Then $abc$ and $a''b'c'$ are two perspective  triangles from the point $o$. Since the plane  $X$ is Desarguesian, the nonempty set
$$T\defeq(\Aline ab\cap\Aline{a''}{b'})\cup(\Aline bc\cap\Aline{b'}{c'})\cup(\Aline ac\cap\Aline{a''}{c'})=(\Aline ab\cap\Aline{a''}{b'})\cup\{y\}\cup\{z\}$$has rank $\|T\|=2$ and hence $\Aline ab\cap\Aline{a''}{b'}\subseteq\Aline yz=\Aline xz$. Since the line $\Aline ab=\Aline xa$ is not parallel to the line $\Aline ac$, the lines $\Aline ab$ and $\Aline {a''}{b'}$ have a common point $x'\in \Aline ab\cap\Aline {a''}{b'}\subseteq \Aline xb\cap\Aline xz$. It follows from $\Aline ac\cap \Aline {a''}{c'}=\{z\}\ne\varnothing=\Aline ac\cap\Aline{a'}{c'}$ that $a'\ne a''$. Then $\{b'\}=\Aline {a'}{b'}\cap\Aline{a''}{b'}=\Aline x{b'}\cap\Aline{x'}{b'}$ implies $x\ne x'$. In its turn, $\{x,x'\}\subseteq \Aline xb\cap \Aline xz$ implies $\Aline ab=\Aline xb=\Aline x{x'}=\Aline xz=\Aline xy$ and hence $y\in \Aline ab$ and $c\in \Aline ac\cap \Aline by\subseteq \Aline ac\cap\Aline ab=\{a\}$, which is impossible as $abc$ is a triangle. This contradiction shows that $\Aline xy\cap \Aline ac=\varnothing$. By analogy we can show that $\Aline xy\cap\Aline{a'}{c'}=\varnothing$. 
\end{proof}

Now we are able to prove three Inverse Desargues Theorems for proaffine regular liners.
Let us recall that three distinct lines $A,B,C$ are called \index{paraconcurrent lines}\index{lines!paraconcurrent}\defterm{paraconcurrent} if they are parallel or $A\cap B\cap C$ is a singleton. 

\begin{definition} Two triangles $abc$ and $a'b'c'$ in a liner $X$ are called \index{paraperspective triangles}\index{triangles!paraperspective}\defterm{paraperspective} if there exist distinct paraconcurrent lines $A,B,C$ in $X$ such that $a,a'\in A$, $b,b'\in B$, $c,c'\in C$.
\end{definition}

\begin{exercise} Show that two disjoint triangles $abc$ and $a'b'c'$ in a liner $X$ are paraconcurrent if and only if the lines $\Aline a{a'}$, $\Aline b{b'}$, $\Aline c{c'}$ are distinct and paraconcurrent.
\end{exercise} 

\begin{theorem}\label{t:ID1} Let $X$ be a proaffine regular liner satisfying the Affine Desargues Axiom and let $abc$ and $a'b'c'$ be two disjoint triangles such that the lines $\Aline a{a'},\Aline b{b'},\Aline c{c'}$ are distinct. If $\Aline ab\parallel \Aline {a'}{b'}$, $\Aline bc\parallel \Aline{b'}{c'}$ and $\Aline ac\parallel \Aline{a'}{c'}$, then the triangles $abc$ and $a'b'c'$  are paraperspective.
\end{theorem} 

\begin{proof} Assume that  $\Aline ab\parallel \Aline {a'}{b'}$, $\Aline bc\parallel \Aline{b'}{c'}$, $\Aline ac\parallel \Aline{a'}{c'}$. We have to prove that the lines $A\defeq\Aline a{a'}$, $B\defeq\Aline b{b'}$, $C\defeq\Aline c{c'}$ are parallel or concurrent. It follows from $\Aline ab\parallel \Aline{a'}{b'}$ that $\overline{A\cup B}=\overline{\{a,a',b,b'\}}$ is a plane. If the lines $A,B$ are disjoint, then they are parallel, by Corollary~\ref{c:parallel-lines<=>}. By analogy we can prove that $B\cap C=\varnothing=A\cap C$ implies $B\parallel C$ and $A\parallel C$. Therefore, the lines $A,B,C$ are parallel if the set $O\defeq (A\cap B)\cup (B\cap C)\cup(A\cap C)$ is empty. 

Next, assume that the set $O$ contains some point $o$. We lose no generality assuming that $o\in A\cap B$ and hence $A\cap B=\{o\}$ (because the lines $A,B,C$ are distinct). Let us check that $o\notin\{a,a',b,b',c,c'\}$. Assuming that $o=a$, we conclude that $o\notin\{a',b,b'\}$. Taking into account that $b'\in B=\Aline ob=\Aline ab\parallel \Aline{a'}{b'}$, we conclude that $a'\in\Aline {a'}{b'}=B$ and hence $A=\Aline a{a'}=\Aline o{a'}=B$, which contradicts our assumption. This contradiction shows that $o\ne a$. By analogy we can prove that $o\notin\{a',b,b'\}$. Assuming that $o=c$, we conclude that $\Aline{a'}{a}=\Aline ao=\Aline ac\parallel\Aline{a'}{c'}$ and hence $c'\in \Aline{a'}{c'}=\Aline {a'}{a}$ and $\Aline c{c'}=\Aline a{a'}$, which contradicts the choice of the triangles $abc$ and $a'b'c'$. This contradiction shows that $o\ne c$. By analogy we can prove that $o\ne c'$. Therefore, $o\notin\{a,b,c,a',b',c'\}$.

Assuming that $o\in\Aline bc$, we conclude that $b'\in \Aline ob=\Aline bc$. Taking into account that $\Aline bc\parallel\Aline{b'}{c'}$, we conclude that $\Aline {b'}{c'}=\Aline bc$ and hence $\Aline b{b'}=\Aline bc=\Aline c{c'}$, which contradicts the choice of the triangles $abc$ and $a'b'c'$. This contradiction shows that $o\notin \Aline bc$. By analogy we can prove that $o\notin\Aline{b'}{c'}$.

We claim that $\Aline ab\cap\Aline {a'}{b'}=\varnothing=\Aline bc\cap\Aline{b'}{c'}$. Assuming that the parallel lines  $\Aline ab$ and $\Aline{a'}{b'}$ have a common point,  we conclude that $\Aline ab=\Aline {a'}{b'}$ and hence $\{a\}=A\cap\Aline ab=A\cap\Aline{a'}{b'}=\{a'\}$, which contradicts the choice of the distinct points $a,a'$. This contradiction shows that $\Aline ab\cap\Aline{a'}{b'}=\varnothing$. By analogy we can prove that $\Aline bc\cap\Aline{b'}{c'}=\varnothing$ and $\Aline ac\cap\Aline{a'}{c'}=\varnothing$. 

It follows from $\Aline bc\parallel \Aline {b'}{c'}$ that $P=\overline{\{b,b',c,c'\}}$ is a plane containing the lines $B=\Aline b{b'}$ and $\Aline oc$. Taking into account that $\Aline oc\cap\Aline bc=\{c\}$ and $\Aline bc\cap\Aline{b'}{c'}=\varnothing$, we can apply Proposition~\ref{p:Proclus-Postulate} and conclude that  there exists a unique point $c''\in \Aline o{c}\cap \Aline{b'}{c'}$. Then $abc$ and $a'b'c''$ are two disjoint perspective triangles from the point $o$ such that $\Aline ab\cap\Aline {a'}{b'}=\varnothing=\Aline bc\cap\Aline{b'}{c'}=\Aline bc\cap\Aline {b'}{c''}$. By Theorem~\ref{t:Proclus<=>}, the proaffine regular liner $X$ is Proclus and weakly regular. 
If the planes $\overline{\{a,b,c\}}$ and $\overline{\{a',b',c''\}}$ are distinct, then by Theorem~\ref{t:Des-Proclus}, the set
$$T\defeq (\Aline ab\cap\Aline{a'}{b'})\cup(\Aline bc\cap\Aline{b'}{c''})\cup(\Aline ac\cap\Aline{a'}{c''})=\varnothing \cup\varnothing\cup(\Aline ac\cap\Aline{a'}{c''}$$has rank $\|T\|\in\{0,2\}$ and hence $\Aline ac\cap\Aline{a'}{c''}=\varnothing$. So, in both cases we obtain that the lines $\Aline ac$ and $\Aline{a'}{c''}$ are disjoint.

If $\overline{\{a,b,c\}}=\overline{\{a',b',c''\}}$, then we can apply the the Affine Desargues Axiom and conclude  that $\Aline ac\cap \Aline {a'}{c''}=\varnothing$.  It follows from $\Aline ac\parallel \Aline {a'}{c'}$ and $o\in\Aline a{a'}$ that the disjoint lines $\Aline ac$ and $\Aline{a'}{c''}$ belong to the plane $\overline{\{o,a,c,a',c'\}}$. By Corollary~\ref{c:parallel-lines<=>}, the lines $\Aline ac$ and $\Aline {a'}{c''}$ are parallel. By Theorem~\ref{t:Proclus-lines}, $\Aline {a'}{c'}\parallel \Aline ac\parallel \Aline {a'}{c''}$ implies $\Aline{a'}{c'}\parallel \Aline{a'}{c''}$ and hence $\Aline {a'}{c'}=\Aline{a'}{c''}$. Then $c''\in \Aline {a'}{c''}\cap \Aline{b'}{c'}=\Aline {a'}{c'}\cap\Aline{b'}{c'}=\{c'\}$ and hence $o\in \Aline c{c''}=\Aline c{c'}=C$ and finally, $O=\{o\}\in A\cap B\cap C$, which means that the lines $A,B,C$ are concurrent.  
\end{proof}

Lemma~\ref{l:D=>ADA} and Theorem~\ref{t:ID1} imply the following corollary.

\begin{corollary}Let $abc$ and $a'b'c'$ be two disjoint triangles in a Desarguesian proaffine regular liner $X$ such that the lines $\Aline a{a'},\Aline b{b'},\Aline c{c'}$ are distinct. If $\Aline ab\parallel \Aline {a'}{b'}$, $\Aline bc\parallel \Aline{b'}{c'}$ and $\Aline ac\parallel \Aline{a'}{c'}$, then the triangles $abc$ and $a'b'c'$ are paraperspective.
\end{corollary}

\begin{theorem}\label{t:ID2} Let  $abc$ and $a'b'c'$ be two disjoint triangles in a Desarguesian proaffine regular liner $X$ such that $\Aline ac\parallel \Aline {a'}{c'}$ and $\Aline ab\cap\Aline {a'}{b'}=\{x\}$, $\Aline bc\cap \Aline{b'}{c'}=\{y\}$ for some points $x,y$. If $\Aline xy\parallel \Aline ac$, then the triangles $abc$ and $a'b'c'$ are paraperspective. 
\end{theorem}

\begin{proof}  Consider the lines $A\defeq\Aline a{a'}$, $B\defeq\Aline b{b'}$, $C\defeq\Aline c{c'}$ and observe that they are pairwise coplanar. Assuming that $(A\cap B)\cup(B\cap C)=\varnothing$, we conclude that $A\parallel B\parallel C$ and hence the lines $A,B,C$ are parallel, by Corollary~\ref{c:parallel-lines<=>} and  Theorem~\ref{t:Proclus-lines}. So, assume that the set $(A\cap B)\cup (B\cap C)$ contains some point $o$. Without loss of generality, we can assume that $o\in A\cap B$. We have to prove that $o\in C$. To derive a contradiction, assume that $o\notin C$. It follows from $\Aline xy\parallel \Aline ac$ and $a\ne c$ that $x\ne y$.

Let us check that $o\notin\{a,b,a',b'\}$. Assuming that $o=a$, we conclude that $a=o\in B$ and  $x\in\Aline ab=B$. It follows from $a\in B\ne A=\Aline a{a'}$ that $a'\notin B$. Then $x\in B\cap\Aline {a'}{b'}=\{b'\}$, $y\in\Aline {b'}{c'}=\Aline x{c'}$ and $\Aline{a'}{c'}\parallel \Aline ac\parallel \Aline xy=\Aline{b'}{c'}$, which impies the equality $\Aline{a'}{c'}=\Aline{b'}{c'}$ contradicting the choice of the triangle $a'b'c'$. By analogy we can prove that $o\ne a'$. Assuming that $o=b$, we conclude that $b=o\in A\ne B$ and hence $b'\notin A$. Then $x\in\Aline ab\cap \Aline{a'}{b'}=A\cap\Aline{a'}{b'}=\{a'\}$ and $\Aline{a'}y=\Aline xy\parallel \Aline ac\parallel \Aline{a'}{c'}$ implies $\Aline {a'}y=\Aline{a'}{c'}$. Since $a'b'c'$ is a triangle, $y\in\Aline{a'}{c'}\cap\Aline{b'}{c'}=\{c'\}$. Then $c'=y\in \Aline bc=\Aline oc$ and hence $o\in \Aline c{c'}=C$, which contradicts our assumption. This contradiction shows that $o\ne b$. By analogy we can prove that $o\ne b'$.
Therefore, $o\notin\{a,b,a',b'\}$.

Assuming that $o\in\Aline ac$, we conclude that $A=\Aline ao=\Aline ac\parallel \Aline{a'}{c'}$ and hence $\Aline {a'}{c'}=A=\Aline ac$ and $\Aline c{c'}=\Aline a{a'}$, which contradicts our assumption. This contradiction shows that $o\notin\Aline ac$. By analogy we can prove that $o\notin\Aline{a'}{c'}$.
It follows from $B\ne C$ that $c\notin B$ or $c'\notin B$. We lose no generality assuming that $c\notin B$.

 It follows from $o\in \Aline a{a'}\setminus(\Aline ac\cup\Aline{a'}{c'})$ that the parallel lines $\Aline ac$ and $ \Aline{a'}{c'}$ are disjoint. Applying Proposition~\ref{p:Proclus-Postulate}, we can find a unique point $c''\in \Aline oc\cap\Aline {a'}{c'}$. Then the triangles $abc$ and $a'b'c''$ are perspective from the point $o$. By Theorem~\ref{t:pD=>PP}, there exists a point $y'\in \Aline bc\cap\Aline {b'}{c''}$ such that $\Aline x{y'}\parallel \Aline ac$. 
Since $\Aline xy\parallel \Aline ac$, the lines $\Aline xy$ and $\Aline x{y'}$ coincide, by the Proclus Axiom. Since $abc$ is a triangle, the lines $\Aline bc$ and $\Aline ac$ are not parallel and hence the lines $\Aline bc$ and $\Aline xy=\Aline x{y'}$ are not parallel as well. Then $\{y\}=\Aline bc\cap \Aline xy=\Aline bc\cap \Aline x{y'}=\{y'\}$. We claim that the point $y=y'$ is not equal to $b'$. In the opposite case, $b\ne b'=y\in B\cap \Aline bc\cap\Aline{b'}{c'}$ implies $c\in B$, which contradicts our assumption. This contradiction shows that $y=y'\ne b'$. Then $c''\in\Aline{a'}{c'}\cap \Aline {b'}{y'}=\Aline {a'}{c'}\cap \Aline {b'}y=\Aline {a'}{c'}\cap\Aline{b'}{c'}=\{c'\}$, and finally, $o\in \Aline c{c''}=\Aline c{c'}=C$, which contradicts our assumption.

The contradiction shows that $o\in A\cap B\cap C$ and hence the distinct lines $A,B,C$ are concurrent.
\end{proof}

\begin{theorem}\label{t:ID3} Let $abc$ and $a'b'c'$ be two disjoint triangles in a Desarguesian proaffine regular liner $X$ such that the lines $\Aline a{a'},\Aline b{b'},\Aline c{c'}$ are distinct, and $\Aline ab\cap\Aline {a'}{b'}=\{x\}$, $\Aline bc\cap \Aline{b'}{c'}=\{y\}$ and $\Aline ac\cap \Aline{a'}{c'}=\{z\}$ for some points $x,y,z\in X\setminus\{a,a',b,b',c,c'\}$. If the points $x,y,z$ are collinear, then the triangles $abc$ and $a'b'c'$ are paraperspective.
\end{theorem}

\begin{proof} Assume that the points $x,y,z$ are collinear. Let us show that $|\{x,y,z\}|=3$. Assuming that $x=y$, we conclude that $x=y\in\Aline ab\cap \Aline bc\cap\Aline{a'}{b'}\cap\Aline {b'}{c'}=\{b\}\cap\{b'\}=\varnothing$, which is a contradiction showing that $x\ne y$. By analogy we can prove  that $y\ne z$ and $x\ne z$. Consider the lines $A\defeq\Aline a{a'}$, $B\defeq\Aline b{b'}$ and $C\defeq \Aline c{c'}$.

\begin{claim}\label{cl:abc-notin-ABC} $a,a'\notin B\cup C$, $b,b'\notin A\cup C$, $c,c'\notin A\cup B$.
\end{claim}

\begin{proof} Assuming that $a\in B$ and taking into account that $a\in B\ne A=\Aline a{a'}$, we conclude that $a'\notin B$ and hence $x\in \Aline ab\cap\Aline {a'}{b'}=B\cap\Aline {a'}{b'}=\{b'\}$, which contradicts our assumption. By analogy we can prove that $a\notin C$ and also that $a'\notin B\cup C$, $b,b'\notin A\cup C$ and $c,c'\notin A\cup B$.
\end{proof}

Since the lines $\Aline ab$ and $\Aline{a'}{b'}$ are concurent, they are coplanar, which implies that the lines $A=\Aline a{a'}$ and $B=\Aline b{b'}$ are coplanar. If $A\cap B=\varnothing$, then $A\parallel B$, by Corollary~\ref{c:parallel-lines<=>}. The same argument applies to the pairs $B,C$ and $A,C$. Therefore, if the set $O\defeq(A\cap B)\cup(B\cap C)\cup (A\cap C)$ is empty, then  the lines $A,B,C$ are parallel and the triangles $abc, a'b'c'$ are paraperspective.

So, assume that the set $O$ contains some point $o$. We lose no generality assuming that $o\in A\cap B$. We have to prove that $o\in C$.  Since $a,a'\in A\setminus  B$ and $b,b'\in B\setminus A$, the point $o\in A\cap B$ does not belong to the set $\{a,a',b,b'\}$.

Since the lines $\Aline ac$ and $\Aline{a'}{c'}$ are concurrent, the flat $\overline{\{a,c,a',c'\}}$ is a plane containing the lines $\Aline {a'}{c'}$ and $\Aline oc$. By the same argument, the plane $\overline{\{b,c,b',c'\}}$ contains the lines $\Aline{b'}{c'}$ and $\Aline oc$. Assuming that $\Aline oc\cap (\Aline {a'}{c'}\cup\Aline{b'}{c'})=\varnothing$, we conclude that $\Aline{a'}{c'}\parallel \Aline oc\parallel \Aline{b'}{c'}$ and hence $\Aline{a'}{c'}\parallel \Aline{b'}{c'}$ and $\Aline{a'}{c'}=\Aline{b'}{c'}$, by Theorem~\ref{t:Proclus-lines} and Proposition~\ref{p:para+intersect=>coincide}. The equality $\Aline{a'}{c'}=\Aline{b'}{c'}$ implies that $a'b'c'$ is not a triangle. This contradiction shows that the line $\Aline oc$ has a common point $c''$ with the set $\Aline {a'}{c'}\cup\Aline{b'}{c'}$. We lose no generality assuming that $c''\in \Aline{a'}{c'}$.

Let us show that $c''\notin\{o,c\}$. Assuming that $c''=o$, we conclude that $o=c''\in\Aline{a'}{c'}$ and $c'\in \Aline o{a'}=A$, which contradicts 
 Claim~\ref{cl:abc-notin-ABC}. Assuming that $c''=c$, we conclude that $c=c''\in \Aline{a'}{c'}$ and hence $a'\in \Aline c{c'}=C$, which contradicts  Claim~\ref{cl:abc-notin-ABC}. Therefore, $c''\notin\{o,c\}$.

\begin{picture}(140,150)(-130,-15)

\put(0,0){\line(1,0){120}}

\put(120,0){\line(0,1){120}}
\put(0,60){\line(1,0){120}}
\put(60,0){\line(1,2){60}}

{
\put(0,0){\line(0,1){120}}
\put(0,120){\line(1,-1){80}}
}

\put(0,0){\circle*{3}}
\put(-4,-10){$a'$}

\put(80,40){\circle*{3}}
\put(82,34){$c''$}
\put(68,37){$c'$}

\put(60,60){\circle*{3}}
\put(57,51){$c$}

\put(0,120){\circle*{3}}
\put(-3,123){$o$}

{
\put(0,0){\line(1,0){60}}
\put(0,0){\line(2,1){80}}
\put(60,0){\line(1,2){20}}

\put(0,60){\line(1,0){60}}
\put(0,60){\line(2,-1){40}}
\put(40,40){\line(1,1){20}}
}

{
\put(60,0){\line(-1,2){60}}
}

\put(0,60){\circle*{3}}
\put(-8,57){$a$}
\put(40,40){\circle*{3}}
\put(33,32){$b$}

\put(60,0){\circle*{3}}
\put(57,-10){$b'$}

\put(0,0){\line(2,1){120}}
\put(120,0){\line(-2,1){120}}
\put(120,120){\line(-1,-1){80}}

\put(120,0){\circle*{3}}
\put(123,-4){$x$}
\put(120,60){\circle*{3}}
\put(123,57){$z$}
\put(120,120){\circle*{3}}
\put(123,118){$y$}
\end{picture}

Therefore, $c''\notin \{o,c\}$. Then $abc$ and $a'b'c''$ are two triangles, perspective from the point $o$. Theorem~\ref{t:D-Proclus-wr<=>} implies that the nonempty set
$$(\Aline ab\cap\Aline{a'}{b'})\cup(\Aline bc\cap\Aline{b'}{c''})\cap(\Aline ac\cap\Aline{a'}{c'})=\{x\}\cup(\Aline bc\cap\Aline{b'}{c''})\cup\{z\}$$has rank 2 and hence $\Aline bc\cap\Aline{b'}{c''}\subseteq\Aline xz$. The collinearity of the points $x,y,z$ implies that $\Aline bc\cap\Aline {b'}{c'}=\{y\}\subseteq \Aline xz$. Assuming that $b\in\Aline xz$, we conclude that $a\in \Aline xb\subseteq \Aline xz$ and $c\in \Aline bc=\Aline by=\Aline xz$, which implies that the points $a,b,c$ are collinear.
But this contradicts the choice of the triangle $abc$. Therefore, $b\notin \Aline xz$.

We claim that  $\Aline bc\cap\Aline{b'}{c''}\ne\varnothing$. To derive a contradiction, assume that  $\Aline bc\cap\Aline {b'}{c''}=\varnothing$. Observe that $\Aline ab\cap\Aline {a'}{b'}=\{x\}$ and $\Aline ac\cap\Aline{a'}{c''}=\Aline ac\cap\Aline{a'}{c'}=\{z\}$. By Theorem~\ref{t:pD=>PP}, $\Aline xz\cap \Aline bc=\varnothing$, which contradicts $y\in \Aline xz\cap\Aline bc$. This contradiction shows that $\Aline bc\cap\Aline{b'}{c''}\ne\varnothing$ and hence $$\varnothing\ne \Aline bc\cap\Aline{b'}{c''}=\Aline {b'}{c''}\cap\Aline {b}{c}\cap\Aline xz=\Aline{b'}{c''}\cap\{y\}\subseteq \{y\}$$ and  $y\in\Aline{b'}{c''}$. Then $c''\in\Aline {b'}{y}\cap\Aline{a'}{c'}=\Aline{b'}{c'}\cap\Aline{a'}{c'}=\{c'\}$ and finally, $o\in \Aline c{c''}=\Aline c{c'}=C$. Therefore, $o\in A\cap B\cap C$ and the lines $A,B,C$ are concurrent and the triangles $abc$ and $a'b'c'$ are paraperspective.
\end{proof}

\section{Parallel Desargues Theorems}

In this section we prove three theorems on properties of triangles whose vertices lie on three parallel lines.

\begin{theorem}\label{t:PD1} Let $A,B,C$ be three distinct parallel lines in a Desarguesian proaffine regular liner $X$ and $a,a'\in A$, $b,b'\in B$, $c,c'\in C$ are any points. If $\Aline ab\parallel \Aline {a'}{b'}$ and $\Aline bc\parallel \Aline{b'}{c'}$, then $\Aline ac\parallel \Aline {a'}{c'}$.
\end{theorem}

\begin{proof} If $b\in\Aline ac$, then $\Aline {a'}{b'}\parallel \Aline ab=\Aline bc\parallel \Aline{b'}{c'}$ implies $\Aline {a'}{b'}=\Aline{b'}{c'}=\Aline {a'}{c'}$ and hence $\Aline{a'}{c'}=\Aline {a'}{b'}\parallel \Aline ab=\Aline ac$. By analogy we can show that $b'\in\Aline {a'}{c'}$ implies $\Aline ac\parallel \Aline{a'}{c'}$. So, assume that $b\notin\Aline ac$ and $b'\notin\Aline{a'}{c'}$, which means that $abc$ and $a'b'c'$ are two triangles in $X$. Since $A,C$ are distinct parallel lines, their flat hull $P\defeq\overline{A\cup C}$ is a plane in $X$. 

Since $\Aline ab\parallel \Aline{a'}{b'}$ and $\Aline bc\parallel \Aline{b'}{c'}$, we can apply  Theorem~\ref{t:subparallel-via-base} and conclude that $\overline{\{a,b,c\}}$ and $\overline{\{a',b',c'\}}$ are two parallel planes in $X$. Observe that $\Aline ac\subseteq P\cap\overline{\{a,b,c\}}$ and $\Aline{a'}{c'}\subseteq P\cap\overline{\{a',b',c'\}}$. If $b\notin P$, then $P\cap\overline{\{a,b,c\}}=\Aline ac$, by the $3$-rankedness of the regular liner $X$. By Theorem~\ref{t:w-modular<=>}, the regular liner $X$ is weakly modular, and by Corollary~\ref{c:paraintersect}, the flat $I\defeq P\cap\overline{\{a',b',c'\}}$ is parallel to the line $\Aline ac=P\cap\overline{\{a,b,c\}}$.  Then $\Aline{a'}{c'}=I\parallel \Aline ac$, by Corollary~\ref{c:parallel}. By analogy we can prove that $b'\notin P$ implies $\Aline ac\parallel\Aline{a'}{c'}$.

So, assume that $b,b'\in P$. 
If $b=b'$, then $\Aline ab\parallel\Aline {a'}{b'}$ and $\Aline bc\parallel \Aline{b'}{c'}$ imply $\Aline ab=\Aline{a'}{b'}$ and $\Aline bc=\Aline{b'}{c'}$. In this case $\{a\}=A\cap\Aline ab=A\cap\Aline{a'}{b'}=\{a'\}$, $\{c\}=C\cap\Aline bc=C\cap\Aline{b'}{c'}=\{c'\}$, $\Aline ac=\Aline{a'}{c'}$ and finally, $\Aline ac\parallel \Aline {a'}{c'}$.

So, we assume that $b\ne b'$. In this case the parallelity relations $\Aline ab\parallel\Aline {a'}{b'}$ and $\Aline bc\parallel \Aline{b'}{c'}$ imply $a\ne a'$ and $c\ne c'$. Then $|\{a,b,c,a',b',c'\}|=6$. 
By Theorem~\ref{t:Proclus<=>}, the proaffine regular liner $P$ is Proclus. Assuming that $P$ is not $3$-long, we can apply Theorem~\ref{t:Proclus-not-3long} and conclude that $P$ is a subliner of some Steiner projective plane. Since the plane $P$ contains disjoint lines, it is not projective and hence $|P|<7$ and $P=\{a,b,c,a',b',c'\}$ is a punctured Steiner projective plane. Then $\{A,B,C\}$ is a unique spread of parallel lines in $P$ and hence the disjoint lines $\Aline ab,\Aline {a'}{b'}$ belong to the spread $\{A,B,C\}$, which is a contradiction showing that the Proclus plane $P$ is $3$-long.   

By Corollary~\ref{c:Desarg=>compreg} and Theorem~\ref{t:spread=projective1}, the Desarguesian proaffine regular liner $P$ is completely regular and hence para-Playfair and bi-Bolyai.
Then the disjoint lines $\Aline ab$ and $\Aline {a'}{b'}$ in $P$ are Bolyai. By the same reason, the disjoint lines $\Aline bc$ and $\Aline{b'}{c'}$ are Bolyai in $P$. Since the Bolyai lines $\Aline ab$ and $\Aline bc$ are concurrent,  the line $\Aline ac\subseteq P=\overline{\Aline ab\cup\Aline bc}$ is Bolyis in $P$, by the bi-Bolyai property of $P$. Therefore,  there exists a unique line $L\subseteq P$ such that $a'\in L$ and $L\parallel \Aline ac$.

\begin{picture}(300,100)(-120,-15)

\put(0,0){\line(1,0){150}}
\put(0,32){\line(1,0){150}}
\put(0,64){\line(1,0){150}}
\put(16,32){\line(1,-2){16}}
\put(16,32){\line(1,2){16}}
\put(100,0){\line(0,1){64}}
\put(32,0){\line(0,1){64}}
\put(84,32){\line(1,-2){16}}
\put(84,32){\line(1,2){16}}

\put(153,-3){$C$}
\put(153,29){$B$}
\put(153,61){$A$}
\put(103,45){$L$}

\put(102,-9){$c''$}
\put(100,0){\circle*{3}}
\put(91,-9){$c'$}
\put(32,0){\circle*{3}}
\put(29,-8){$c$}
\put(100,64){\circle*{3}}
\put(98,68){$a'$}
\put(84,32){\circle*{3}}
\put(76,34){$b'$}
\put(16,32){\circle*{3}}
\put(10,34){$b$}
\put(32,64){\circle*{3}}
\put(29,68){$a$}

\end{picture}

It follows from $a'\notin \Aline ac$ that $\|L\cup\Aline ac\|=3=\|P\|$ and hence $P=\overline{\Aline ac\cup L}$, by the $3$-rankedness of the regular liner $X$.  Since $|\Aline bc\cap\Aline ac|=1$ and $\Aline bc\subseteq P=\overline{\Aline ac\cup L}$, we can apply  Theorem~\ref{t:Proclus-lines} and conclude that the line $\Aline bc$ intersects the line $L$ at a unique point. Taking into account that $b'\notin\Aline bc$, we conclude that $\|\Aline bc\cup\Aline{b'}{c'}\|=3$ and hence $L\subseteq P=\overline{\Aline bc\cup\Aline{b'}{c'}}$.  Since $|L\cap\Aline bc|=1$, we can apply Theorem~\ref{t:Proclus-lines} and conclude that the line $L$ intersects the line $\Aline {b'}{c'}$ at a single point $c''$. Assuming that $c''=a'$, we conclude that $a'\in\Aline {b'}{c'}$ and hence $b'\in\Aline {a'}{c'}$, which contradicts our assumption. This contradiction shows that $c''\ne a'$ and hence $L=\Aline {a'}{c''}$. Assuming that $c''=b'$, we conclude that $L=\Aline {a'}{c''}=\Aline {a'}{b'}\parallel \Aline ab$. Then $\Aline ac\parallel L\parallel\Aline ab$ implies $\Aline ac=\Aline ab$, which contradicts $b\notin\Aline ac$.
This contradiction shows that $c''\ne b'$ and hence $\Aline {b'}{c''}=\Aline{b'}{c'}\parallel \Aline bc$. 
Then $abc$ and $\Aline {a'}{b'}{c''}$ are two triangles such that the lines $\Aline a{a'}=A$, $\Aline b{b'}=B$ and $\Aline{c}{c''}$ are pairwise distinct and $\Aline ab\parallel \Aline {a'}{b'}$, $\Aline bc\parallel \Aline {b'}{c''}$ and $\Aline ac\parallel \Aline{a'}{c''}$. The Inverse Desargues Theorem~\ref{t:ID1} ensures that the lines $A,B,\Aline c{c''}$ are paraconcurrent. Since the lines $A,B$ are parallel, the line $\Aline c{c''}$ is parallel to the lines $A,B$ and hence $\Aline c{c''}=C$ because the liner $X$ is Proclus. Then $\{c'\}=\Aline {b'}{c'}\cap C=\Aline{b'}{c''}\cap C=\{c''\}$ and hence $\Aline {a'}{c'}=\Aline {a'}{c''}=L\parallel \Aline ac$.
\end{proof}

\begin{theorem}\label{t:PD2} Let $A,B,C$ be three distinct parallel lines in a Desarguesian proaffine regular liner $X$ and $a,a'\in A$, $b,b'\in B$, $c,c'\in C$ be any points. If $\Aline ac\parallel \Aline {a'}{c'}$ and $\Aline ab\cap \Aline{a'}{b'}=\{x\}$ for some point $x\in X$, then there exists a unique point $y\in X$ such that $\Aline bc\cap\Aline {b'}{c'}=\{y\}$ and $\Aline xy\subparallel \Aline ac$.
\end{theorem}

\begin{proof} Assuming that the lines $\Aline {b}{c}$ and $\Aline {b'}{c'}$ in the plane $B\cup C$ are disjoint, we conclude that $\Aline bc\parallel\Aline{b'}{c'}$. Since $\Aline ac\parallel\Aline{a'}{c'}$ and $\Aline bc\parallel \Aline{b'}{c'}$, we can apply Theorem~\ref{t:PD1} and conclude that $\Aline ab\parallel \Aline {a'}{b'}$, which contradicts $\Aline ab\cap\Aline{a'}{b'}=\{x\}$. Therefore, there exists a point $y\in \Aline bc\cap\Aline {b'}{c'}$. Assuming that $\Aline bc\cap\Aline{b'}{c'}\ne\{y\}$, we conclude that $\Aline bc=\Aline{b'}{c'}$ and hence $\{b\}=B\cap\Aline bc=B\cap \Aline{b'}{c'}=\{b'\}$ and $\{c\}=C\cap\Aline bc=C\cap\Aline{b'}{c'}=\{c'\}$. Then $\Aline ac\parallel \Aline{a'}{c'}$ implies $\Aline ac=\Aline{a'}{c'}$ and hence $a=A\cap\Aline ac=A\cap\Aline{a'}{c'}=\{a'\}$ and $\Aline ab=\Aline{a'}{b'}=\Aline ab\cap\Aline {a'}{b'}=\{x\}$, which is a contradiction showing that $\Aline bc\cap\Aline {b'}{c'}=\{y\}$.

It remains to prove that $\Aline xy\subparallel \Aline ac$. This suparallel relation trivially holds if $x=y$. So, assume that $x\ne y$. In this case $b\ne b'$. Indeed, assuming that $b=b'$, we conclude that $b=b'\in \Aline ab\cap\Aline{a'}{b'}\cap\Aline bc\cap\Aline {b'}{c'}=\{x\}\cap\{y\}$ and hence $x=b=y$.

If the line $B$ is disjoint with the plane $P\defeq \overline{A\cup C}$, then $\overline{\{a,b,c\}}$ and $\overline{\{a',b',c'\}}$ are two planes in $X$ that have at least two distinct common points $x$ and $y$. The $3$-rankedness of $X$ and $b,b'\notin P$ imply that $P\cap\overline{\{a,b,c\}}=\Aline ac$ and $P\cap\overline{\{a',b',c'\}}=\Aline{a'}{c'}$. Assuming that $\Aline xy\ne\overline{\{a,b,c\}}\cap\overline{\{a',b',c'\}}$, we conclude that $\overline{\{a,b,c\}}=\overline{\{a',b',c'\}}=\overline{\{a,a',c,c'\}}=P$, which contradicts our assumption. This contradiction shows that $\Aline xy=\overline{\{a,b,c\}}\cap\overline{\{a',b',c'\}}$ and hence
$$\Aline xy\cap\Aline ac=(\overline{\{a,b,c\}}\cap\overline{\{a',b',c'\}})\cap (P\cap\overline{\{a,b,c\}})=(P\cap \overline{\{a,b,c\}})\cap(P\cap\overline{\{a',b',c'\}})=\Aline ac\cap\Aline{a'}{c'}=\varnothing.$$By Corollary~\ref{c:parallel-lines<=>}, the disjoint lines  $\Aline xy$ and $\Aline ac$ in the plane $\overline{\{a,b,c\}}$ are parallel.   

Next, assume that $B\cap P$ contains some point $\beta$. Then $A\parallel B$ implies $B\subseteq\overline{\{\beta\}\cup A}\subseteq P$. So, $A,B,C$ are three parallel lines in the plane $P$. If the parallel lines $\Aline ac$ and $\Aline{a'}{c'}$ have a common point, then they coincide and hence $\{a\}=A\cap\Aline ac=A\cap\Aline{a'}{c'}=\{a'\}$ and $\{c\}=C\cap\Aline ac=C\cap\Aline{a'}{c'}=\{c'\}$. Then $a=a'\in\Aline ab\cap\Aline{a'}{b'}=\{x\}$, $c=c'\in\Aline bc\cap\Aline{b'}{c'}=\{y\}$ and finally, $\Aline xy=\Aline ac\parallel \Aline ac$. 

So, assume that the parallel lines $\Aline ac$ and $\Aline {a'}{c'}$ are disjoint. If $b\in\Aline ac$, then $x\in \Aline ab=\Aline ac$ and $y\in \Aline bc=\Aline ac$ and hence $\Aline xy=\Aline ac\parallel\Aline ac$. By analogy we can prove that $b'\in\Aline {a'}{c'}$ implies $\Aline xy\parallel \Aline {a'}{c'}\parallel \Aline ac$. So, assume that $b\notin\Aline ac$ and $b'\notin\Aline {a'}{c'}$. It follows from $a\ne a'$ and $b\ne b'$ that $x\notin A\cup B$. By the same reason, $b\ne b'$ and $c\ne c'$ imply  $y\notin B\cup C$.

 Then $|P|\ge |\{a,b,c,a',b',c'\}|=6$. By Theorem~\ref{t:Proclus<=>}, the proaffine regular plane $P$ is Proclus. Since $P$ contains disjoint lines, the plane $P$ is not projective.  Assuming that $P$ is not $3$-long, we can apply Theorem~\ref{t:Proclus-not-3long} and conclude that $P$ is a  proper subliner of some Steiner projective plane. 
Then $|P|<7$ and hence $P=\{a,b,c,a',b',c'\}$ is a punctured Steiner projective plane containing a single spread of parallel lines $\{A,B,C\}$. Then the disjoint lines $\Aline ac$ and $\Aline{a'}{b'}$ must be parallel to the lines $A,B,C$, which is not true. This contradiction shows that the Proclus plane $P$ is $3$-long. By Corollary~\ref{c:Desarg=>compreg} and Theorem~\ref{t:spread=projective1}, the $3$-long Desarguesian Proclus plane $P$ is para-Playfair. Then the plane $P$ contains a line $L_x$ such that $x\in L_x$ and $L_x\parallel \Aline ac$. Taking into account that $\Aline bc\cap\Aline ac=\{c\}$, we can apply Proposition~\ref{p:Proclus-Postulate} and find a unique point $y'\in L_x\cap\Aline bc$. 

\begin{claim}\label{cl:y'b'} $y'\notin\{x\}\cup B$ and $b'\notin\Aline x{y'}$.
\end{claim}

\begin{proof} Assuming that $y'=x$, we conclude that $x\in \Aline ab\cap \Aline bc=\{b\}$
 and hence $b=x\in\Aline {a'}{b'}\cap B=\{b'\}$, which contradicts our assumption. This contradiction shows that $x\ne y'$ and hence $\Aline x{y'}=L_x\parallel \Aline ac$. Assuming that $y'\in B$, we conclude that $y'\in B\cap \Aline bc=\{b\}$ and hence $a\in \Aline xb=\Aline x{y'}=L_x\parallel \Aline ac$ implies $b=y'\in L_x=\Aline ac$, which contradicts our assumption. This contradiction shows that $y'\notin B$. 

If $b'\in \Aline x{y'}=L_x$, then 
$\Aline{a'}{b'}=\Aline x{y'}=L_x\parallel \Aline {a'}{c'}$ and hence $b'\in\Aline {a'}{b'}=\Aline{a'}{c'}$, which contradicts our assumption. This contradiction shows that $b'\notin\Aline{x}{y'}$.
\end{proof}

By Claim~\ref{cl:y'b'}, $b'\notin\Aline x{y'}$ and hence $\Aline{y'}{b'}\cap\Aline x{y'}=\{y'\}$. Since $\Aline x{y'}\parallel \Aline {a'}{c'}$, we can apply Proposition~\ref{p:Proclus-Postulate} and find a unique point $c''\in \Aline{a'}{c'}\cap\Aline{y'}{b'}$. 

\begin{picture}(200,130)(-200,-20)

\put(-5,95){$L_x$}
\put(0,-8){\color{cyan}\line(0,1){53}}
\put(-30,-8){\color{cyan}\line(0,1){68}}
\put(30,-8){\color{cyan}\line(0,1){68}}
\put(0,0){\line(1,1){90}}
\put(0,45){\line(2,1){90}}
\put(0,0){\line(-1,1){90}}
\put(0,45){\line(-2,1){90}}
\put(-30,60){\color{red}\line(1,0){60}}
\put(-30,30){\color{red}\line(1,0){60}}
\put(-90,90){\color{red}\line(1,0){180}}

\put(0,0){\circle*{3}}
\put(3,-6){$b'$}
\put(0,45){\circle*{3}}
\put(-2,48){$b$}
\put(30,30){\circle*{3}}
\put(22,32){$c'$}
\put(33,24){$c''$}
\put(-30,30){\circle*{3}}
\put(-27,32){$a'$}
\put(30,60){\circle*{3}}
\put(28,63){$c$}
\put(-30,60){\circle*{3}}
\put(-33,63){$a$}
\put(90,90){\color{red}\circle*{3}}
\put(85,94){$y'$}
\put(93,84){$y$}
\put(-90,90){\circle*{3}}
\put(-98,88){$x$}

\put(-34,-20){$A$}
\put(-4,-20){$B$}
\put(26,-20){$C$}
\end{picture}

\begin{claim} $c''\notin\{a,b,c,a',b'\}$.
\end{claim}

\begin{proof} It follows from $c''\in\Aline{a'}{c'}$ and $\{a,c,b'\}\cap\Aline{a'}{c'}=\varnothing$ that $c''\notin\{a,c,b'\}$. Assuming that $c''=b$, we conclude that $b=c''\in B\cap \Aline {y'}{b'}=\{b'\}$, which contradicts our assumption. Assuming that $c''= a'$, we conclude that $a'=c''\in\Aline {y'}{b'}$, $x\in \Aline {a'}{b'}= \Aline {y'}{b'}$ and $b'\in\Aline x{y'}$, which contradicts Claim~\ref{cl:y'b'}.
\end{proof}

The preceding claim ensures that $c''\notin\{a,b,c,a',b'\}$ and hence $abc$ and $a'b'c''$ are two disjoint triangles such that $\Aline ab\cap\Aline{a'}{b'}=\{x\}$, $\Aline bc\cap\Aline {b'}{c''}=\{y'\}$ and $\Aline x{y'}\parallel\Aline ac\parallel \Aline{a'}{c''}$. By Theorem~\ref{t:ID2}, the lines $\Aline a{a'},\Aline b{b'},\Aline c{c''}$ are paraconcurrent. Since the lines $A=\Aline a{a'}$ and $B=\Aline b{b'}$ are parallel, the line $\Aline c{c''}$ is parallel to the lines $A,B,C$ and hence $\Aline c{c''}=C$, by the Proclus Axiom. Then $c''\in C\cap\Aline{a'}{c'}=\{c'\}$, $y'\in \Aline bc\cap \Aline {b'}{c''}=\Aline bc\cap\Aline{b'}{c'}=\{y\}$ and $\Aline xy=\Aline x{y'}=L_x\parallel \Aline ac$.  
\end{proof} 

\begin{theorem}\label{t:PD3} Let $A,B,C$ be three distinct parallel lines in a Desarguesian proaffine regular liner $X$, and let $a,a'\in A$, $b,b'\in B$, $c,c'\in C$, $x,y,z\in X$ be points such that $\Aline ab\cap \Aline {a'}{b'}=\{x\}$, $\Aline bc\cap \Aline{b'}{c'}=\{y\}$, $\Aline ac\cap \Aline{a'}{c'}=\{z\}$. Then the points $x,y,z$ are collinear.
\end{theorem}

\begin{proof} To derive a contradiction, assume that the points $x,y,z$ are not collinear. Then $|\{x,y,z\}|=\|\{x,y,z\}\|=3$. 

\begin{claim}\label{cl:PD3-6}
$|\{a,a',b,b',c,c'\}|=6$.
\end{claim}

\begin{proof}  If $a=a'$, then $a=a'\in \Aline ab\cap\Aline{a'}{b'}=\{x\}$ implies $a=a'=x$. Also $a=a'\in \Aline ac\cap\Aline {a'}{c'}=\{z\}$ implies $z=a=a'=x$, which contradicts $|\{x,y,z\}|=3$.
By analogy we can prove that $b\ne b'$ and $c\ne c'$. Taking into account that the parallel lines $A,B,C$ are disjoint, we conclude that $|\{a,a',b,b',c,c'\}|=6$.
\end{proof}

\begin{claim}\label{cl:PD3-xyz}
$x\notin A\cup B$, $y\notin B\cup C$, $z\notin A\cup C$.
\end{claim}

\begin{proof} Assuming that $x\in A$, we conclude that $$\{x\}=\{x\}\cap A=\Aline ab\cap\Aline {a'}{b'}\cap A=(\Aline ab\cap A)\cap(\Aline{a'}{b'}\cap A)=\{a\}\cap\{a'\}$$ and hence $a=x=a'$, which contradicts Claim~\ref{cl:PD3-6}. This contradiction shows that $x\notin A$. By analogy we can prove that $x\notin B$ and also that $y\notin B\cup C$ and $z\notin A\cup C$.
\end{proof}

\begin{claim}\label{cl:PD3-bb'} $b\notin\Aline ac$ and $b'\notin\Aline{a'}{c'}$.
\end{claim}

\begin{proof} If $b\in\Aline ac$, then  $\{x,y,z\}\subseteq \Aline ab\cup\Aline bc\cup \Aline ac=\Aline ac$ and hence $\|\{x,y,z\}\|\le\|\{a,c\}\|=2$, which contradicts our assumption. This contradiction shows that $b\notin\Aline ac$. By analogy we can prove that $b'\notin\Aline{a'}{c'}$.
\end{proof}

Claims~\ref{cl:PD3-6} and \ref{cl:PD3-bb'} imply that $abc$ and $a'b'c'$ are disjoint triangles in $X$.

\begin{claim}\label{cl:PD3-6not3} $\{a,a',b,b',c,c'\}\cap(\Aline xy\cup\Aline yz\cup\Aline xz)=\varnothing$.
\end{claim}

\begin{proof} Assuming that $a\in \Aline xy$ and taking into account that $x\notin A$ and $y\notin B$, we conclude that $b\in\Aline xa\subseteq \Aline xy$ and $c\in \Aline by=\Aline xy$. Therefore, $\{a,b,c\}\subseteq \Aline xy$, which contradicts Claim~\ref{cl:PD3-bb'}.

 Assuming that $a\in\Aline yz$ and taking into account that $z\notin A$ and $y\notin C$, we conclude that $c\in \Aline za\subseteq \Aline yz$ and $b\in \Aline yc\subseteq \Aline yz$. Therefore, $\{a,b,c\}\subseteq \Aline yz$, which contradicts Claim~\ref{cl:PD3-bb'}.

Assuming that $a\notin \Aline xz$ and taking into account that $x\notin A$ and $z\notin C$, we conclude that $b\in \Aline xa\subseteq \Aline xz$ and $c\in \Aline az\subseteq \Aline xz$.  Therefore, $\{a,b,c\}\subseteq \Aline xz$, which contradicts Claim~\ref{cl:PD3-bb'}. Therefore, $a\notin \Aline xy\cup\Aline yz\cup\Aline xz$. By analogy we can prove that $a',b,b',c,c'\notin\Aline xy\cup\Aline yz\cup\Aline xz$.
\end{proof}

Since $A,C$ are two disjoint parallel lines, the flat $P\defeq\overline{A\cup C}$ is a plane in $X$.

\begin{claim}\label{cl:PD3-BP} $B\subseteq P$.
\end{claim}

\begin{proof} Assuming that $B\not\subseteq P$ and taking into account that $A\parallel B$, we conlude that $B\cap P=\varnothing$. Then $\overline{\{a,b,c\}}$ and $\overline{\{a',b',c'\}}$ are two planes in $X$ such that $\overline{\{a,b,c\}}\cap P=\Aline ac$ and $\overline{\{a',b',c'\}}\cap P=\Aline{a'}{c'}$, by the $3$-rankedness of the regular liner $X$. Since the lines $\Aline ac$ and $\Aline {a'}{c'}$ are distinct, the planes $\overline{\{a,b,c\}}$ and $\overline{\{a',b',c'\}}$ are distinct, too. Then the intersection $\overline{\{a,b,c\}}\cap\overline{\{a',b',c'\}}$ is a line containing the points $x,y,z$, which implies that the points $x,y,z$ are collinear. But this contradicts our assumption. This contradiction shows that $B\subseteq P$.
\end{proof}

Claim~\ref{cl:PD3-BP} ensures that $\{x,y,z\}\subseteq \overline{\{a,b,c,a',b',c'\}}=\overline{A\cup B\cup C}=P$. The lines $\Aline ac$ and $\Aline{a'}{c'}$ are concurrent and hence they cannot be both parallel to the line $\Aline xy$. Therefore, $\Aline ac\nparallel \Aline xy$ or $\Aline {a'}{c'}\nparallel\Aline xy$. We lose no generality assuming that $\Aline {a'}{c'}\nparallel \Aline xy$. Then there exists a unique point $z'\in \Aline {a'}{c'}\cap\Aline xy$.  

\begin{claim}\label{cl:PD3-z'} $z'\notin A\cup \{b,b'\}\cup C$.
\end{claim} 

\begin{proof} If $z'\in A$, then $z'\in A\cap\Aline {a'}{c'}\cap \Aline xy=\{a'\}\cap \Aline xy$ implies $a'=z'\in \Aline xy$, which contradicts Claim~\ref{cl:PD3-6not3}. By analogy, we can derive a contradiction assuming that $z'\in C$. By Claim~\ref{cl:PD3-6not3}, $z'\in \Aline xy\subseteq P\setminus\{b,b'\}$.
\end{proof}

Two cases are possible. 
\smallskip

1. First we assume that the lines $\Aline {a}{z'}$ and $\Aline y{b}=\Aline{b}{c}$ have a common point $c''$.
\begin{claim}\label{cl:PD3-c''} $c''\notin A\cup B\cup\{c'\}\cup \Aline {a}{b}\cup\{x,y,z'\}$.
\end{claim}

\begin{proof} By Claim~\ref{cl:PD3-z'}, $z'\notin A\cup C\cup\{b,b'\}$. Assuming that $c''\in A$, we conclude that $c''\in A\cap\Aline {a}{z'}\cap\Aline{b}{c}=\{a\}\cap\Aline {b}{c}$, which contradicts Claim~\ref{cl:PD3-bb'}. This contradiction shows that $c''\notin A$. 

Assuming that $c''\in B$, we conclude that $c''\in B\cap \Aline {b}{c}\cap \Aline {a}{z'}=\{b'\}\cap\Aline{a}{z'}$ and hence $c''=b\in\Aline {a}{z'}$ and 
$z'\in \Aline{a}{b}\cap \Aline xy=\Aline x{a}\cap \Aline xy=\{x\}$, according to Claim~\ref{cl:PD3-6not3}.
Then $x=z'\in \Aline {a'}{c'}$ and hence $b'\in\Aline {a'}x=\Aline {a'}{c'}$, which contradicts Claim~\ref{cl:PD3-bb'}. This contradiction shows that  $a''\notin B$.

Assuming that $c''=c'$, we conclude that $c'=c''\in C\cap\Aline {b}{c}=\{c\}$, which contradicts Claim~\ref{cl:PD3-6}.

Since $abc$ is a triangle, $a\notin\Aline bc=\Aline {b}{c''}$, which implies $c''\notin\Aline{a}{b}$. 

Assuming that $c''=x$ and taking into account that $x\notin A$, we conclude that $x=c''\in \Aline {a}{z'}$ and hence $\Aline {a}{z'}=\Aline {a}x=\Aline{a}{b}$ and $x=c''\in \Aline{a}{z'}\cap\Aline{b}{c}=\Aline {a}{b}\cap\Aline{b}{c}=\{b\}$, which contradicts Claim~\ref{cl:PD3-6not3}.

Assuming that $c''=y$ and taking into account that $c''\in \Aline {a}{z'}\setminus A$ and $a\notin\Aline xy$, we conclude that $z'\in \Aline {a}{c''}\cap \Aline xy=\Aline {a}y\cap\Aline xy=\{y\}$. Then $y=z'\in \Aline {b'}{c'}\cap\Aline {a'}{c'}=\{c'\}$, which contradicts Claim~\ref{cl:PD3-6not3}. This contradiction shows that $c''\ne y$. 

Assuming that $c''=z'$, we conclude that $c''=z'\in \Aline y{b}\cap\Aline xy\cap\Aline {a'}{c'}=\{y\}\cap\Aline {a'}{c'}$ and hence $y\in \Aline {a'}{c'}$ and $b'\in \Aline y{c'}\subseteq\Aline {a'}{c'}$, which contradicts Claim~\ref{cl:PD3-bb'}.
\end{proof}

\begin{picture}(200,210)(-200,-50)

\put(0,-40){\color{red}\line(0,1){180}}
\put(0,0){\line(1,0){120}}
\put(0,40){\line(1,0){120}}
\put(0,80){\line(1,0){120}}
\put(0,-40){\line(1,2){60}}
\put(0,-10){\line(2,1){100}}
\put(100,40){\line(-1,1){100}}

\put(0,0){\line(-1,0){120}}
\put(0,40){\line(-1,0){120}}
\put(0,80){\line(-1,0){120}}
\put(0,-40){\line(-1,2){60}}
\put(0,-10){\line(-2,1){100}}
\put(-100,40){\line(1,1){100}}

\put(0,140){\circle*{3}}
\put(-3,143){$x$}
\put(0,-40){\color{red}\circle*{3}}
\put(-8,-45){$z$}
\put(2,-45){$z'$}
\put(0,-10){\circle*{3}}
\put(2,-16){$y$}
\put(20,0){\circle*{3}}
\put(15,3){$c'$}
\put(100,40){\circle*{3}}
\put(100,43){$b'$}
\put(60,80){\circle*{3}}
\put(60,83){$a'$}
\put(-20,0){\circle*{3}}
\put(-25,-8){$c$}
\put(-19,3){$c''$}
\put(-100,40){\circle*{3}}
\put(-105,43){$b$}
\put(-60,80){\circle*{3}}
\put(-65,83){$a$}
\put(125,-4){$C$}
\put(125,36){$B$}
\put(125,76){$A$}
\end{picture}

Claims~\ref{cl:PD3-6} and \ref{cl:PD3-bb'} ensure that $a'b'c'$ and $abc''$ are two disjoint triangles such that
$\Aline {a'}{b'}\cap\Aline {a}{b}=\{x\}$, $\Aline {b'}{c'}\cap\Aline {b}{c''}=\Aline {b'}{c'}\cap \Aline{b}{c}=\{y\}$ and $\Aline {a'}{c'}\cap\Aline {a}{c''}=\{z'\}$.  Claims~\ref{cl:PD3-6not3}, \ref{cl:PD3-z'}, \ref{cl:PD3-c''} imply that $\{x,y,z'\}\cap\{a',b',c',a,b,c''\}=\varnothing$. By the Inverse Desargues Theorem~\ref{t:ID3}, the lines $\Aline {a'}a,\Aline {b'}b,\Aline {c'}{c''}$ are paraconcurrent. Since the lines $A=\Aline {a'}a$ and $B=\Aline {b'}b$ are parallel, the line $\Aline {c'}{c''}$ is parallel to the lines $A,B,C$ and hence $\Aline {c'}{c''}=C=\Aline {c'}c$. 
Then $c''\in C\cap \Aline {b}{c}=\{c\}$ and $z'\in \Aline {a'}{c'}\cap\Aline {a}{c''}=\Aline {a'}{c'}\cap\Aline{a}{c}=\{z\}$. Then $\|\{x,y,z\}\|=\|\{x,y,z'\}\|=2<3$, which contradicts our assumption. 
\smallskip

2. Next, assume the second case when $\Aline {a}{z'}\cap\Aline{b}{c}=\varnothing$. Assuming that the Proclus plane $P$ is not $3$-long, we can apply Theorem~\ref{t:Proclus-not-3long} and conclude that $P$ is a proper subliner of a Steiner projective plane and hence $|P|<7$ and $P=\{a,b,c,a',b',c'\}$ is a punctured Steiner projective plane. Then $P$ has a unique spread $\{A,B,C\}$ of parallel lines and any disjoint lines in $P$ should belong to this spread, which is not true for the disjoint lines $\Aline{a}{z'}$ and $\Aline{b}{c}$. This contradiction shows that the Proclus plane $P$ is $3$-long. By Corollary~\ref{c:Desarg=>compreg} and Theorem~\ref{t:spread=projective1}, the $3$-long Desarguesian Proclus plane $P$ is comletely regular and hence para-Playfair and bi-Bolyai. Then the disjoint lines $\Aline{a}{z'}$ and $\Aline{b}{c}$ are Bolyai. Since $z'\notin A$, the Bolyai lines $A$ and $\Aline{a}{z'}$ are concurrent. Since $P$ is bi-Bolyai,  every line in the plane $\overline{A\cup\Aline{a}{z'}}=P$ is Bolyai, which implies that $P$ is a Playfair plane. In this case we can proceed as in Lemma~\ref{t:ADA=>Desargues}.   

Since $P$ is Playfair, there exist points $a''\in A$ and $c''\in C$ such that $\Aline{a''}{c}\parallel \Aline {a'}{c'}\parallel \Aline{a}{c''}$ and hence $\Aline {a''}{c}\parallel \Aline {a}{c''}$. Assuming that $b\in\Aline {a''}{c}\cap\Aline{a}{c''}$, we conclude that $b\in \Aline {a''}{c}=\Aline{a}{c''}$ and hence $\|\{a,b,c\}\|\le 2$, which contradicts Claim~\ref{cl:PD3-bb'}. This contradiction shows that $b\notin\Aline {a''}{c}$ or $b\notin\Aline{a}{c''}$. 
We lose no generality assuming that $b\notin\Aline{a}{c''}$. Assuming that $c''=c$, we conclude that $\Aline ac=\Aline a{a''}\parallel {a'}{c'}$, which contradicts $\Aline ac\cap\Aline{a'}{c'}=\{z\}$. This contradiction shows that $c''\ne c$. Then $abc''$ is a triangle, which is perspective to the triangle $zy{c'}$ from the point $c$.

Now consider the  triangles ${a'}{b'}{c'}$ and $abc''$. Since $\Aline {a'}{b'}\cap\Aline{a}{b}=\{x\}$ and $\Aline {a'}{c'}\cap\Aline{a}{c''}=\varnothing$, Theorem~\ref{t:PD2}, ensures that $\Aline {b'}{c'}\cap\Aline {b}{c''}=\{y'\}$ for some point $y'\in X$ such that  $\Aline x{y'}\subparallel \Aline {a'}{c'}$.  Assuming that $x=y'$, we obtain a contradiction
$$\{x\}=\{y'\}=\Aline {a}{b}\cap \Aline {a'}{b'}\cap \Aline {b'}{c'}\cap\Aline {b}{c''}\subseteq \Aline ab\cap(\Aline {a'}{b'}\cap\Aline {b'}{c'})=\Aline ab\cap\{b'\}=\varnothing,$$showing that $x\ne y'$. By Corollary~\ref{c:subparallel}, $\Aline x{y'}\subparallel \Aline {a'}{c'}\parallel \Aline ac$ implies $\Aline x{y'}\parallel \Aline ac$.

Finally, consider two triangles $abc''$ and $zy{c'}$, which are perspective from the point $c$, and observe that $\Aline {b}{c''}\cap\Aline y{c'}=\Aline {b}{c''}\cap \Aline {b'}{c'}=\{y'\}$. Since $\Aline {a}{c''}\cap \Aline {a'}{c'}=\varnothing$, Theorem~\ref{t:PD2} ensures that $\Aline {a}{b}\cap\Aline zy=\{x'\}$ for some point $x'\in X$ such that $\Aline {x'}{y'}\subparallel \Aline {a'}{c'}$. Assuming that $x'=y'$ and taking into account that $abc''$ is a triangle, we conclude that
$$\{x'\}=\{y'\}=\Aline {a}{b}\cap\Aline zy\cap \Aline {b'}{c'}\cap\Aline{b}{c''}\subseteq (\Aline ba\cap\Aline {b}{c''})\cap\Aline {b'}{c'}=\{b\}\cap\Aline {b'}{c'}=\varnothing,$$
which is a contradiction showing that $x'\ne y'$. Then $\Aline {x'}{y'}\subparallel \Aline {a'}{c'}$ implies $\Aline {x'}{y'}\parallel \Aline {a'}{c'}$, according to Corollary~\ref{c:subparallel}. It follows from $\Aline x{y'}\parallel \Aline {a'}{c'}\parallel \Aline {x'}{y'}$ that $\Aline xy=\Aline {x'}{y'}$.  It follows from $y'\in \Aline {b}{c''}\setminus\{b\}=\Aline{b}{c''}\setminus \Aline {a}{b}$ that 
$$\{x\}=\Aline{a}{b}\cap\Aline x{y'}=\Aline {a}{b}\cap\Aline{x'}{y'}=\{x'\}\subseteq \Aline zy$$ and hence $\|\{x,y,z\}\|\le\|\{y,z\}\|\le 2$, which contradicts our assumption. 
\end{proof} 

Theorems~\ref{t:Des-Proclus}, \ref{t:PD3} and \ref{t:ID3} imply the following characterization.

\begin{corollary}  Let $abc$ and $a'b'c'$ be two disjoint triangles in a Desarguesian proaffine regular liner $X$ such that the lines $\Aline a{a'},\Aline b{b'},\Aline c{c'}$ are distinct, and $\Aline ab\cap\Aline {a'}{b'}=\{x\}$, $\Aline bc\cap \Aline{b'}{c'}=\{y\}$,  $\Aline ac\cap \Aline{a'}{c'}=\{z\}$ for some points $x,y,z\in X\setminus\{a,a',b,b',c,c'\}$. The points $x,y,z$ are collinear if and only if the triangles $abc$ and $a'b'c'$ are paraperspective.
\end{corollary}

\section{Completions of Desarguesian proaffine liners}

By Corollary~\ref{c:Desarg=>compreg}, every Desarguesian proaffine regular liner is completely regular. 

\begin{theorem}\label{t:Desargues-completion} The spread completion of any  Desarguesian completely regular liner is a Desarguesian projective liner.
\end{theorem}

\begin{proof} 
Let $\overline X$ be the spread completion of a Desarguesian completely regular liner $X$. We have to prove that the projective liner $\overline X$ is Desarguesian. If $\|\overline X\|\le 2$, then $\overline X$ is Desarguesian because it contains no triangles. So, assume that $\|\overline X\|\ge 3$. If $X$ is projective, then $\overline X=X$ is a Desarguesian projective liner. So, assume that the liner $X$ is not projective. By Proposition~\ref{p:spread-3long}, the spread completion $\overline X$ of $X$ is a projective $3$-long liner. If $\|\overline X\|\ne 4$, then the $3$-long projective line $\overline X$ is Desarguesian, by Theorem~\ref{t:Desargues-projective}. So, assume that $\|\overline X\|=3$. If $|\overline X|_2=3$, then the Steiner projective liner $\overline X$ is Desarguesian, by Proposition~\ref{p:Steiner+projective=>Desargues}. So, we assume that $|\overline X|_2\ge 4$, which implies that the liner $X$ is $3$-long. By Theorem~\ref{t:spread=projective1}, the completely regular liner $X$ is para-Playfair and hence Proclus. By Corollary~\ref{c:procompletion-rank}, $\|X\|=\|\overline X\|=3$ and hence $X$ is a Proclus plane.
 
To prove that the projective plane $\overline X$ is Desarguesian, take any two centrally perspective triangles $abc$ and $a'b'c'$ in $\overline X$. Let $o\in\Aline a{a'}\cap\Aline b{b'}\cap \Aline c{c'}$ be the unique perspector of the triangles $abc$ and $a'b'c'$ in $\overline X$. By the projectivity of the plane $\overline X$, there exist unique points $x,y,z\in \overline X$ such that $\Aline ab\cap\Aline{a'}{b'}=\{x\}$, $\Aline bc\cap\Aline{b'}{c'}=\{y\}$ and $\Aline ac\cap\Aline{a'}{c'}=\{z\}$. We have to prove that the set $\{x,y,z\}$ has rank $\|\{x,y,z\}\|=2$. 

\begin{claim}\label{cl:xyz-D} $x\notin \Aline oa\cup\Aline ob$,  $y\notin \Aline ob\cup\Aline oc$ and $z\notin \Aline oa\cup\Aline oc$.
\end{claim}

\begin{proof} Assuming that $x\in\Aline oa$, we conclude that $\{x\}=\Aline oa\cap\{x\}=\Aline oa\cap(\Aline ab\cap\Aline{a'}{b'})=(\Aline oa\cap\Aline{ab})\cap(\Aline o{a'}\cap\Aline {a'}{b'})=\{a\}\cap\{a'\}=\varnothing$. By analogy we can prove that $x\notin\Aline ob$. Therefore, $x\notin\Aline oa\cup\Aline ob$. By analogy we can prove that $y\notin \Aline ob\cup\Aline oc$ and $z\notin \Aline oa\cup\Aline oc$.
\end{proof}

\begin{claim}\label{cl:xyz3} $|\{x,y,z\}|=3$.
\end{claim}

\begin{proof} Assuming that $x=y$ and taking into account that $x,y\notin \Aline ob$, we conclude that $a\in \Aline bx=\Aline by=\Aline bc$, which contradicts the choice of the triangle $abc$. By analogy we can prove that $x\ne z\ne y$.
\end{proof}

\begin{claim}\label{cl:abc-not-xyz} $\{a,b,c,a',b',c'\}\cap(\Aline xy\cup \Aline yz\cup\Aline xz)=\varnothing$.
\end{claim}

\begin{proof} Assuming that $a\in \Aline xy$, we conclude that $b\in \Aline ab=\Aline xy$ and $c\in \Aline by=\Aline xy$. So, $\{a,b,c\}\subseteq \Aline xy$, which is impossible as $abc$ is a triangle. By analogy we can prove that $a\notin \Aline yz\cup \Aline xz$ and also that $b,c,a',b',c'\notin \Aline xy\cup\Aline yz\cup\Aline xz$.
\end{proof} 

If $\{x,y,z\}\subseteq \overline X\setminus X$, then 
$$\|\{x,y,z\}\|\le\|\overline X\setminus X\|<\|\overline X\|=3$$and hence the points $x,y,z$ are collinear. So, we assume that $\{x,y,z\}\not\subseteq\overline X\setminus X$.

Claims~\ref{cl:xyz-D} and \ref{cl:xyz3} imply that 
$$|\{o,a,b,c,a',b',c',x,y,z\}|=10.$$

Depending on the cardinality of the intersection $$I\defeq \{o,a,b,c,a',b',c',x,y,z\}\cap(\overline X\setminus X),$$ we shall separately consider five cases.
\smallskip

0. If $|I|=0$, then $\{o,a,b,c,a',b',c',x,y,z\}\subseteq X$ and the Desargues Axiom ensures that the nonempty set $\{x,y,z\}$ has rank $2$ in $X$ and hence in $\overline X$.
\smallskip

1. Next, assume that $|I|=1$. This case has three subcases.
\smallskip

1.1. If $I=\{o\}$, then $A\defeq \Aline oa\cap X$, $B=\Aline ob\cap X$, $C=\Aline oc\cap X$ are three distinct parallel lines in $X$. Applying Theorem~\ref{t:PD3} the paraperspective triangles $abc,a'b'c'$ in the Desarguesian liner $X$, we conclude that the points $x,y,z$ are collinear in $X$ and hence in $\overline X$.
\smallskip

1.2. If $I\cap \{x,y,z\}\ne\varnothing$, then we lose no generality assuming that $I=\{z\}$. In this case, $abc$ and $a'b'c'$ are two centrally perspective triangles in the Desarguesian liner $X$ such that $(\Aline ab\cap X)\cap(\Aline {a'}{b'}\cap X)=\{x\}$, $(\Aline bc\cap X)\cap(\Aline{b'}{c'}\cap X)=\{y\}$ and the lines $(\Aline ac\cap X)\cap(\Aline {a'}{c'}\cap X)$ are parallel in $X$. By Theorem~\ref{t:pD=>PP}, $(\Aline xy\cap X)\parallel (\Aline ac\cap X)$ and hence the intersection $\Aline xy\cap\Aline ac$ contains some point $z'\in \overline X\setminus X$. Since $\Aline ac=(\Aline ac\cap X)\cup\{z\}$, the point $z'$ is equal to $z$. Then $z=z'\in\Aline xy$, which means that the points $x,y,z$ are collinear.
\smallskip

1.3. If $I\cap \{a,b,c,a',b,c'\}\ne\varnothing$, then we lose no generality assuming that $I=\{c'\}$. Then $A\defeq \Aline {a'}{c'}\cap X$, $B\defeq\Aline {b'}{c'}\cap X$, $C\defeq \Aline o{c'}\cap X$ are three distinct parallel lines in the Desarguesian liner $X$ such that $z,a'\in A$, $y,b'\in B$ and $o,c\in C$. 

\begin{picture}(100,110)(-200,-15)
\put(0,0){\line(0,1){80}}
\put(-70,0){\line(1,0){140}}
\put(-70,40){\line(1,0){140}}
\put(-70,80){\line(1,0){140}}
\put(0,0){\line(1,1){40}}
\put(20,0){\line(-1,1){80}}
\put(20,0){\line(0,1){40}}
\put(0,80){\line(1,-1){60}}
\put(60,20){\line(-2,1){120}}
\put(0,20){\line(1,0){60}}

\put(0,0){\circle*{3}}
\put(-3,-8){$o$}
\put(20,0){\circle*{3}}
\put(17,-8){$c$}
\put(0,20){\circle*{3}}
\put(-8,16){$a$}
\put(20,20){\circle*{3}}
\put(15,22){$b$}
\put(60,20){\color{red}\circle*{3}}
\put(60,13){$x$}
\put(63,21){$x'$}
\put(40,40){\circle*{3}}
\put(40,43){$b'$}
\put(20,40){\circle*{3}}
\put(18,45){$y$}
\put(0,80){\circle*{3}}
\put(-2,83){$a'$}
\put(-60,80){\circle*{3}}
\put(-63,83){$z$}
\put(-80,-4){$C$}
\put(-80,36){$B$}
\put(-80,76){$A$}

\end{picture}

Consider the paraperspective triangles $a'ob',zcy$ and observe that  $\Aline {a'}o\cap \Aline zc=\Aline o{a'}\cap \Aline ac=\{a\}$ and $\Aline o{b'}\cap\Aline cy=\Aline o{b}\cap \Aline bc=\{b\}$. Assuming that $(\Aline {a'}{b'}\cap X)\cap(\Aline zy\cap X)=\varnothing$, we can apply Theorem~\ref{t:PD2} and conclude that $(\Aline ab\cap X)\cap(\Aline {a'}{b'}\cap X)=\varnothing$, which contradicts $\Aline ab\cap\Aline{a'}{b'}=\{x\}\subseteq X$. This contradiction shows that the lines $\Aline {a'}{b'}\cap X$ and $\Aline zy\cap X$ have a common point $x'\in X$. Applying Theorem~\ref{t:PD3}, we conclude that $x'\in \Aline ab\cap X$. Then $x'\in \Aline ab\cap \Aline{a'}{b'}=\{x\}$ and hence $x=x'\in \Aline zy$, so the points $x,y,z$ are collinear in $\overline X$.     
\smallskip

2. Next, assume that $|I|=2$. This case has three subcases.
\smallskip

2.1. If $o\in I$, then $|I|=2$ implies $I\cap\{a,b,c,a',b',c'\}=\varnothing$ and hence $I\cap\{x,y,z\}\ne\varnothing$. In this case we lose no generality assuming that $I\cap\{x,y,z\}=\{z\}$. Observe that the lines $A\defeq\Aline oa\cap X$, $B\defeq\Aline ob\cap X$, $C\defeq\Aline oc\cap X$ in the liner $X$ are distinct and parallel, and $abc,a'b'c'$ are paraconcurrent triangles in $X$ such that $(\Aline ab\cap X)\cap(\Aline{a'}{b'}\cap X)=\{x\}$, $(\Aline bc\cap X)\cap(\Aline {b'}{c'}\cap X)=\{y\}$ and $(\Aline ac\cap X)\cap(\Aline{a'}{c'}\cap X)=\{z\}\cap X=\varnothing$. By Theorem~\ref{t:PD2}, the line $\Aline xy\cap X$ is parallel to the line $\Aline ac\cap X$, which implies that $\varnothing \ne\Aline xy\cap \Aline ac\setminus X\subseteq\Aline ac\setminus X=\{z\}$ and hence the points $x,y,z$ are collinear in $\overline X$.
\smallskip

2.2. If $o\notin I$ and $I\subseteq \{a,b,c,a',b',c,c'\}$, then we lose no generality assuming that $I=\{b',c'\}$ or $I=\{b,c'\}$.  If $I=\{b',c'\}$, then $y\in \Aline {b'}{c'}\subseteq \overline X\setminus X$ and hence $y\in I$, which contradicts $I=\{b',c'\}$. Therefore, the case $I=\{b',c'\}$ is impossible and thus $I=\{b,c'\}$.

\begin{picture}(100,110)(-170,-15)
\put(0,0){\line(0,1){80}}
\put(-20,0){\line(1,0){90}}
\put(20,40){\line(1,0){50}}
\put(0,80){\line(1,0){70}}
\put(0,0){\line(1,1){40}}
\put(-20,0){\line(1,1){40}}
\put(40,40){\line(-1,1){40}}
\put(-20,0){\line(1,2){40}}
\put(0,40){\line(1,1){20}}
\put(20,40){\line(0,1){40}}

\put(0,0){\circle*{3}}
\put(-3,-10){$b'$}
\put(-20,0){\circle*{3}}
\put(-23,-8){$y$}
\put(0,40){\circle*{3}}
\put(-9,43){$x'$}
\put(-10,33){$x$}
\put(40,40){\circle*{3}}
\put(40,43){$o$}
\put(20,40){\circle*{3}}
\put(22,42){$c$}
\put(20,60){\circle*{3}}
\put(22,60){$a$}
\put(20,80){\circle*{3}}
\put(18,83){$z$}
\put(0,80){\circle*{3}}
\put(-2,83){$a'$}
\put(80,-4){$B$}
\put(80,36){$C$}
\put(80,76){$A$}
\end{picture}

Since $c'\in \overline X\setminus X$, the intersections $A\defeq \Aline {a'}{c'}\cap X$, $B\defeq \Aline {b'}{c'}\cap X$, $C\defeq\Aline o{c'}\cap X$ are distinct parallel lines in the Desarguesian liner $X$ such that $y\in\Aline{b'}{c'}=B$ and $z\in \Aline{a'}{c'}=A$. Then the triangles $oa'b'$ and $czy$ are paraperspective and $(\Aline o{a'}\cap X)\cap (\Aline cz\cap X)=\Aline oa\cap \Aline ca\cap X=\{a\}\cap X=\{a\}$ and $(\Aline o{b'}\cap X)\cap(\Aline cy\cap X)=(\Aline o{b}\cap \Aline bc)\cap X=\{b\}\cap X=\varnothing$, which means that the lines $\Aline o{b'}\cap X$ and $\Aline cy\cap X$ are parallel. By Theorem~\ref{t:PD2}, there exists a point $x'\in (\Aline {a'}{b'}\cap X)\cap(\Aline zy\cap X)$ such that $\Aline a{x'}\cap X\parallel \Aline o{b'}\cap X$. On the other hand, $(\Aline ax\cap X)\cap(\Aline o{b'}\cap X)=(\Aline ab\cap\Aline ob)\cap X=\{b\}\cap X=\varnothing$, which means that  $\Aline ax\cap X\parallel \Aline o{b'}\cap X\parallel \Aline a{x'}\cap X$ and hence $\Aline ax\cap X\parallel \Aline a{x'}\cap X$, by Proposition~\ref{p:Proclus-Postulate}. The Proclus Axiom implies that $\Aline ax\cap X=\Aline a{x'}\cap X$ and hence $x'\in \Aline {a'}{b'}\cap \Aline a{x'}=\Aline {a'}{b'}\cap\Aline ax=\Aline {a'}{b'}\cap\Aline ab=\{x\}$ and hence $x=x'\in \Aline zy\cap X$. So, the points $x,y,z$ are collinear in $X$ and also in $\overline X$. 
\smallskip

2.3. If $o\notin I$ and $I\cap\{a,b,c,a',b',c'\}\ne\varnothing\ne I\cap\{x,y,z\}$, then we lose no generality assuming that $I\cap\{x,y,z\}=\{z\}$. Assuming that $a\in I$, we conclude that $c\in \Aline az\subseteq\overline X\setminus X$ and hence $c\in I$, which contradicts $I=\{a,z\}$. This contradiction shows that $a\notin I$. By analogy we can prove that $a',c,c'\notin I$. Then $I=\{b,z\}$ or $I=\{b',z\}$. We lose no generality assuming that $I=\{b',z\}$. 

\begin{picture}(200,120)(-180,-60)

\put(0,0){\line(1,0){80}}
\put(0,0){\line(1,1){40}}
\put(0,0){\line(1,-1){40}}
\put(60,0){\line(-2,1){80}}
\put(60,0){\line(-2,-1){80}}
\put(-20,40){\line(1,0){100}}
\put(-20,-40){\line(1,0){100}}
\put(-20,-40){\line(0,10){80}}
\put(20,-20){\line(0,1){40}}
\put(40,-40){\line(0,1){80}}

\put(0,0){\circle*{3}}
\put(-8,-2){$o$}
\put(60,0){\circle*{3}}
\put(61,2){$b$}
\put(20,20){\circle*{3}}
\put(17,23){$a$}
\put(20,-20){\circle*{3}}
\put(17,-28){$c$}
\put(40,40){\circle*{3}}
\put(38,43){$a'$}
\put(40,-40){\circle*{3}}
\put(38,-49){$c'$}
\put(-20,-40){\circle*{3}}
\put(-28,38){$x$}
\put(-20,40){\circle*{3}}
\put(-28,-42){$y$}
\put(90,-4){$B$}
\put(90,-44){$C$}
\put(90,36){$A$}
\end{picture}

Consider the parallel lines $A\defeq \Aline{a'}{b'}\cap X$, $B\defeq \Aline o{b'}\cap X$ and $C\defeq\Aline {c'}{b'}\cap X$ in $X$, and observe that $x,a'\in A$, $o,b\in B$, $y,c'\in C$. Then the triangles $a'oc'$ and $xby$ are paraperspective in $X$ and have $(\Aline {a'}o\cap X)\cap(\Aline xb\cap X)=(\Aline{a}o\cap\Aline ab)\cap X=\{a\}$ and $(\Aline o{c'}\cap X)\cap(\Aline by\cap X)=(\Aline o{c}\cap\Aline bc)\cap X=\{c\}$. By Theorem~\ref{t:PD3}, the set 
$$
T\defeq (\Aline {a'}o\cap\Aline xb\cap X)\cup(\Aline o{c'}\cap\Aline by\cap X)\cup(\Aline {a'}{c'}\cap \Aline xy\cap X)=\{a,c\}\cup(\Aline{a'}{c'}\cap\Aline xy\cap X)
$$
has rank $\|T\|=2$ and hence $\Aline{a'}{c'}\cap\Aline xy\cap X\subseteq \Aline {a}{c}\cap \Aline {a'}{c'}\cap X=\{z\}\cap X=\varnothing$. By the projectivity of $\overline X$, we have $$\varnothing\ne \Aline {a'}{c'}\cap \Aline xy\subseteq (\Aline{a'}{c'}\cap\Aline xy\cap X)\cup(\Aline{a'}{c'}\setminus X)=\varnothing\cup\{z\}.$$ Then  $z\in \Aline xy$ and the points $x,y,z$ are collinear.
\smallskip

3. Next, assume that $|I|=3$. Since $\{x,y,z\}\not\subseteq I$, this case has five subcases.
\smallskip

3.1. If $o\in I$ and $I\cap\{a,b,c,a',b',c'\}\ne\varnothing$, then we lose no generality that $c\in I$. In this case $c'\in\Aline oc\subseteq \overline X\setminus X$ and hence $c'\in I$ and $I=\{o,c,c'\}$ and $\{a,b,a',b',x,y,z\}\subseteq X$. Consider the triangles $aza'$ and $byb'$ in $X$ and observe that $(\Aline az\cap X)\cap(\Aline by\cap X)=(\Aline ac\cap\Aline bc)\cap X=\{c\}\cap X=\varnothing$, $(\Aline z{a'}\cap X)\cap(\Aline y{b'}\cap X)=(\Aline {c'}{a'}\cap \Aline{c'}{b'})\cap X=\{c'\}\cap X=\varnothing$ and $(\Aline a{a'}\cap X)\cap (\Aline b{b'}\cap X)=(\Aline a{a'}\cap\Aline b{b'})\cap X=\{o\}\cap X=\varnothing$.

\begin{picture}(100,110)(-200,-55)

\put(0,0){\line(1,0){60}}
\put(0,0){\line(1,1){40}}
\put(0,0){\line(1,-1){40}}
\put(60,0){\line(-1,2){20}}
\put(60,0){\line(-1,-2){20}}
\put(40,-40){\line(0,1){80}}

\put(0,0){\line(-1,0){60}}
\put(0,0){\line(-1,-1){40}}
\put(0,0){\line(-1,1){40}}
\put(-60,0){\line(1,-2){20}}
\put(-60,0){\line(1,2){20}}
\put(-40,-40){\line(0,1){80}}

\put(0,0){\circle*{3}}
\put(-3,-10){$x$}
\put(60,0){\circle*{3}}
\put(63,-3){$y$}
\put(40,-40){\circle*{3}}
\put(38,-50){$b'$}
\put(40,40){\circle*{3}}
\put(38,43){$b$}
\put(-60,0){\circle*{3}}
\put(-69,-3){$z$}
\put(-40,40){\circle*{3}}
\put(-42,-48){$a$}
\put(-40,-40){\circle*{3}}
\put(-42,43){$a'$}
\end{picture}

Claim~\ref{cl:abc-not-xyz} implies that the lines $\Aline ab\cap X, \Aline{a'}{b'}\cap X$ and $\Aline zy\cap X$ are distinct. Applying Theorem~\ref{t:ID1}, we conclude that the triangles $aza'$ and $byb'$ are paraconcurrent and hence $\{x\}=(\Aline ab\cap X)\cap(\Aline{a'}{b'}\cap X)=\Aline ab\cap\Aline {a'}{b'}\cap\Aline zy$, which implies that the points $x,y,z$ are collinear.
\smallskip

3.2. If $o\in I$ and $I\cap\{a,b,c,a',b',c'\}=\varnothing$, then $|I\cap\{x,y,z\}|=2$ and we lose no generality assuming that $I\cap\{x,y,z\}=\{x,y\}$. Then $A\defeq \Aline oa\cap X$, $B\defeq\Aline ob\cap X$, $C\defeq\Aline oc\cap X$ are three distinct parallel lines in $X$ and $abc$, $a'b'c'$ are two paraconcurrent triangles in $X$ such that $(\Aline ab\cap X)\cap(\Aline{a'}{b'}\cap X)=(\Aline ab\cap \Aline{a'}{b'})\cap X=\{x\}\cap X=\varnothing$ and $(\Aline bc\cap X)\cap(\Aline{b'}{c'}\cap X)=(\Aline bc\cap \Aline{b'}{c'})\cap X=\{y\}\cap X=\varnothing$. Applying Theorem~\ref{t:PD1}, we conclude that $\varnothing =(\Aline ac\cap X)\cap(\Aline{a'}{c'}\cap X)=(\Aline ac\cap\Aline{a'}{c'})\cap X=\{z\}\cap X=\{z\}$, which is a contradiction showing that this subcase is impossible.
\smallskip

3.3. If $o\notin I$ and $\{x,y,z\}\cap I=\varnothing$, then $I$ has at most one-point intersection with each of the lines $\Aline oa$, $\Aline ob$, $\Aline oc$. Assuming that $I=\{a,b,c\}$, we conclude that the points $a,b,c\subseteq \overline X\setminus X$ collinear, which contradicts the choice of the triangle $abc$. By analogy we obtain a contradiction assuming that $I=\{a',b',c'\}$. Then $I\cap\{a,b,c\}\ne\varnothing\ne I\cap\{a',b',c'\}$ and we lose no generality assuming that $I=\{a',b,c'\}$. Then $z\in \Aline {a'}{c'}\subseteq \overline X\setminus X$ and hence $z\in I$, which contradicts the equality $I=\{a',b,c'\}$. Therefore, the case $o\notin I$ and $\{x,y,z\}\cap I=\varnothing$ is impossible.
\smallskip

3.4. If $o\notin I$ and $|\{x,y,z\}\cap I|=1$, then we lose no generality assuming that $I\cap\{x,y,z\}=\{z\}$. Then $|I\cap \{a,b,c,a',b',c'\}|=2$. Since $o\notin I$, the set $I$ has at most one-point intersection with every line $\Aline oa$, $\Aline ob$ or $\Aline oc$. Then $I\cap(\Aline oa\cap \Aline oc)\ne\varnothing$. We lose no generality assuming that $a'\in I$. Then $c'\in \Aline {a'}z\subseteq \overline X\setminus X$ and hence $I=\{a',c'\}$.
 Claim~\ref{cl:abc-not-xyz} ensures that $b'\notin\Aline xy$ and hence $xb'y$ is a triangle in $X$. 

\begin{picture}(200,100)(-120,-50)

\put(0,0){\line(1,0){105}}
\put(0,0){\line(1,1){40}}
\put(0,0){\line(1,-1){40}}
\put(40,-40){\line(0,10){80}}
\put(40,-40){\line(1,2){35}}
\put(40,40){\line(1,-2){35}}
\put(75,-30){\line(0,1){60}}
\put(105,0){\line(-1,1){30}}
\put(105,0){\line(-1,-1){30}}

\put(0,0){\circle*{3}}
\put(-8,-2){$o$}
\put(40,40){\circle*{3}}
\put(38,43){$c$}
\put(40,-40){\circle*{3}}
\put(38,-48){$a$}
\put(60,0){\circle*{3}}
\put(58,-11){$b$}
\put(105,0){\circle*{3}}
\put(108,-3){$b'$}
\put(75,30){\circle*{3}}
\put(73,33){$x$}
\put(75,-30){\circle*{3}}
\put(73,-38){$y$}
\end{picture}

Then $aoc$ and $xb'y$ are two triangles in $X$ such that $(\Aline ao\cap X)\cap(\Aline x{b'}\cap X)=(\Aline a{a'}\cap \Aline{a'}{b'}\cap X=\{a'\}\cap X=\varnothing$ and $(\Aline oc\cap X)\cap(\Aline{b'}y\cap X)=(\Aline c{c'}\cap \Aline {b'}{c'})\cap X=\{c'\}\cap X=\varnothing$. Observe also that
$\Aline ax\cap \Aline o{b'}\cap\Aline cy=\Aline ab\cap\Aline ob\cap \Aline bc=\{b\}$, which means that the triangles $aoc$ and $xb'y$ are perspective from the point $b$. The Affine Desargues Axiom ensures that $(\Aline ac\cap X)\parallel (\Aline xy\cap X)$ and hence $\varnothing\ne \Aline ac\cap \Aline xy\subseteq (\Aline ac\cap \Aline xy\cap X)\cup(\Aline ac\setminus X)=\varnothing\cup\{z\}$ and $z\in \Aline xy$.
\smallskip

3.5. If $o\notin I$ and $|I\cap\{x,y,z\}|=2$, then we lose no generality assuming that $I\cap \{x,y,z\}=\{x,y\}$. Taking into account that $\overline X\setminus X$ is a line in $\overline X$, we conclude that $\overline X\setminus X=\Aline xy$ and hence $\varnothing\ne \{a,b,c,a',b,c'\}\cap I=\{a,b,c,a',b',c'\}\cap\Aline xy$, which contradicts Claim~\ref{cl:abc-not-xyz}. This contradiction shows that the case $o\notin I$ and $|I\cap\{x,y,z\}|=2$ is impossible.
\smallskip

4. Finally assume that $|I|\ge 4$. By our assumption, $|I\cap \{x,y,z\}|\le 2$ and hence $I\cap\{a,b,c,a',b',c'\}\ne \varnothing$. Repeating the argument from the case 3.5, we can show that $|I\cap \{x,y,z\}|\ne 2$ and hence $|I\cap\{x,y,z\}|\le 1$. If  $I\cap\{x,y,z\}=\varnothing$, then $|I\cap\{o,a,b,c,a',b',c'\}|\ge 4$ and $I$ contains one of the sets $\{o,a,a'\}$, $\{o,b,b'\}$, $\{o,c,c'\}$. We lose no generality assuming that $I=\{o,c,c',b'\}$ and then $z\in\Aline {b'}{c'}\subseteq\overline X\setminus X$ and hence $z\in I$, which contradicts $I\cap\{x,y,z\}=\varnothing$. Therefore, $I\cap\{x,y,z\}$ is a singleton and we lose no generality assuming that $I\cap\{x,y,z\}=\{z\}$. Then $|I\cap\{o,a,b,c,a',b',c'\}|\ge 3$. If $I\cap\{a,c,a',c'\}\ne\varnothing$, then we lose no generality assuming that $a'\in I$. Then $c'\in \Aline {a'}z\subseteq\overline X\setminus X$ and hence $\{a',c'\}\subseteq I$. 
Taking into account that  $\|I\|\le\|\overline X\setminus X\|\le 2$ and $\Aline oa\cap\Aline oc=\{o\}$, we conclude that $\{o,a,c\}\cap I=\varnothing$. Assuming that $b'\in I$, we conclude that $\{x,y\}\subseteq\Aline {a'}{b'}\cup\Aline {b'}{c'}\subseteq\overline X\setminus X$ and hence $\{x,y,z\}\subseteq I$, which contradicts our assumption. Therefore, $\{o,a,c,b'\}\cap I=\varnothing$ and hence $I=\{z,a',c',b\}$ and $\{o,a,b',c,x,y\}\subseteq X$. Consider the triangles $aoc$ and $xb'y$ in $X$ and observe that $(\Aline ax\cap X)\cap(\Aline o{b'}\cap X)=(\Aline ab\cap\Aline o{b})\cap X=\{b\}\cap X=\varnothing$ and $(\Aline o{b'}\cap X)\cap (\Aline cy\cap X)=(\Aline ob\cap\Aline cb)\cap X=\{b\}\cap X=\varnothing$.

\begin{picture}(100,90)(-150,-45)

\put(0,0){\line(1,0){50}}
\put(0,0){\line(1,1){30}}
\put(0,0){\line(1,-1){30}}
\put(30,30){\line(1,0){50}}
\put(30,-30){\line(1,0){50}}
\put(50,0){\line(1,1){30}}
\put(50,0){\line(1,-1){30}}
\put(30,-30){\line(0,1){60}}
\put(80,-30){\line(0,1){60}}

\put(0,0){\circle*{3}}
\put(-8,-2){$o$}
\put(50,0){\circle*{3}}
\put(47,3){$b'$}
\put(30,-30){\circle*{3}}
\put(28,-38){$c$}
\put(30,30){\circle*{3}}
\put(28,33){$a$}
\put(80,30){\circle*{3}}
\put(78,-38){$y$}
\put(80,-30){\circle*{3}}
\put(78,33){$x$}

\end{picture}

Therefore, the triangles $aoc$ and $xb'y$ are paraperspective. Taking into account that  $(\Aline ao\cap X)\cap(\Aline x{b'}\cap X) =(\Aline o{a'}\cap\Aline {a'}{b'})\cap X=\{a'\}\cap X=\varnothing$ and $(\Aline oc\cap X)\cap(\Aline {b'}y\cap X)=(\Aline o{c'}\cap\Aline {b'}{c'})\cap X=\{c'\}\cap X=\varnothing$ and applying the Affine Desargues Axiom, we conclude that $(\Aline ac\cap X)\cap(\Aline xy\cap X)=\varnothing$. The projectivity of the liner $\overline X$ ensures that $\varnothing\ne \Aline ac\cap\Aline xy\subseteq (\Aline xc\cap\Aline xy\cap X)\cup(\Aline ac\setminus X)=\varnothing\cup\{z\}$ and hence $z\in\Aline xy$ and the points $x,y,z$ are collinear.
 \end{proof}

\begin{corollary}\label{c:completion-Desarguesian} Any projective completion of a Desarguesian liner is a Desarguesian projective liner.
\end{corollary}

\begin{proof} Let $Y$ be a projective completion of a Desarguesian liner $X\subseteq Y$. Then $Y$ is a $3$-long projective liner and $\overline{Y\setminus X}\ne Y$. By Corollary~\ref{c:Avogadro-projective}, the $3$-long projective liner $Y$ is $2$-balanced and hence the cardinal number $|Y|_2\ge 3$ is well-defined. If $|Y|_2=3$, then the Steiner projective liner $Y$ is Desarguesian, by Proposition~\ref{p:Steiner+projective=>Desargues}. So, assume that the liner $Y$ is $4$-long. Since the liner $X$ is Desarguesian, its horizon $H=Y\setminus X$ is flat in the $4$-long liner $Y$, by Proposition~\ref{p:nD=>pP}. By Theorem~\ref{t:proaffine<=>proflat} and Proposition~\ref{p:projective-minus-flat}, the liner $X=Y\setminus H$ is proaffine and regular. Since $H$ is flat in $Y$, $|Y|_2\ge 4$ implies that the liner $X$ is $3$-long. 
By Corollary~\ref{c:Desarg=>compreg}, the Desarguesian proaffine regular liner $X$ is completely regular, and by Theorem~\ref{t:Desargues-completion}, the spread completion $\overline X$ of $X$ is a Desarguesian projective liner. By Corollary~\ref{c:pcompletion=scompletion}, the projective completion $Y$ of $X$ is isomorphic to the Desarguesian projective liner $\overline X$ and hence $Y$ is Desarguesian.
\end{proof}



Corollary~\ref{c:Desarg=>compreg} and Theorem~\ref{t:Desargues-completion} imply the following characterization of  Desarguesian liners.

\begin{theorem}\label{t:Desarguesian-liner<=>} For any liner $X$, the following conditions are equivalent:
\begin{enumerate}
\item $X$ is Desarguesian, proaffine and regular;
\item $X$ is completely regular and its spread completion is Desarguesian;
\item $X$ is $3$-ranked and the spread completion of $X$ is a Desarguesian projective liner.
\end{enumerate}
\end{theorem}

Let us recall that a liner $X$ is \defterm{completely $\mathcal P$} for some property $\mathcal P$ of projective liners if $X$ is completely regular and its spread completion $\overline X$ has property $\mathcal P$. Theorem~\ref{t:Desarguesian-liner<=>} implies the following corollary answering Problem~\ref{prob:inner-completely-P} for the Desargues property.

\begin{corollary}\label{c:completelyD<=>} A liner is completely Desarguesian if and only if it is  Desarguesian, proaffine, and regular.
\end{corollary}

Let us recall that a property $\mathcal P$ of liners is called \defterm{complete} if a completely regular liner has property $\mathcal P$ if and only if its spread completion has property $\mathcal P$. Theorem~\ref{t:Desarguesian-liner<=>} implies the following important fact.

\begin{corollary}\label{c:Desargues=>complete} The property of a liner to be Desarguesian is complete.
\end{corollary}

Next, we prove that a projective liner is Desarguesian if and only if it is everywhere Desarguesian.

\begin{theorem}\label{t:projDes<=>} For a projective liner $X$, the following conditions are equivalent:
\begin{enumerate}
\item $X$ is Desarguesian;
\item for every flat $H$ in $X$, the liner $X\setminus H$ is Desarguesian;
\item for every hyperplane $H\subset X$, the affine liner $X\setminus H$ is Desarguesian.
\end{enumerate}
If the projective liner $X$ is $3$-long, then the conditions \textup{(1)--(3)} are equivalent to the conditions:
\begin{enumerate}
\item[(4)] for some hyperplane $H$ in $X$, the affine liner $X\setminus H$ is Desarguesian;
\item[(5)] for some flat $H\ne X$, the liner $X\setminus H$ is Desarguesian.
\end{enumerate}
\end{theorem}

\begin{proof} The implication $(1)\Ra(2)$ follows from Theorem~\ref{t:pD-minus-flat}, and the implication $(2)\Ra(3)$ is trivial.
\vskip5pt

$(3)\Ra(1)$ Assume that for every hyperplane $H\subseteq X$, the liner $X\setminus H$ is Desarguesian. By Theorem~\ref{t:Desargues<=>planeD}, to prove that the liner $X$ is Desarguesian, it suffices to check that every $3$-long plane $P$ in $X$ is Desarguesian. By Corollary~\ref{c:Avogadro-projective}, the $3$-long projective liner $P$ is $2$-balanced and hence the cardinal $|P|_2$ is well-defined. If $|P|_2=3$, then the Steiner projective plane $P$ is Desarguesian, by Proposition~\ref{p:Steiner+projective=>Desargues}. So, assume that $|P|_2\ge 4$. Fix any independent set $a,b,c\in P$. By the Kuratowski--Zorn Lemma, there exists a maximal independent set $I\subseteq X$ such that $\{a,b,c\}\subseteq I$. The maximal independence of $I$ ensures that the flat $H\defeq\overline{I\setminus\{a\}}$ is a hyperplane in $X$. By the assumption, the liner $X\setminus H$ is Desarguesian and so is the liner $(X\setminus H)\cap P=P\setminus H$. Since $a\notin H$ and $\{b,c\}\subseteq P\cap H$, the flat $H\cap P$ is a hyperplane in $P$. Then $P$ is a projective completion of the Desarguesian liner  $P\setminus H$. By Corollary~\ref{c:completion-Desarguesian}, the projective liner $P$ is Desarguesian.
\smallskip

The implications $(3)\Ra(4)\Ra(5)$ are trivial, and the implication $(5)\Ra(1)$ follows from Corollary~\ref{c:completion-Desarguesian}.
\end{proof}

Theorem~\ref{t:projDes<=>} implies the following corollary (which can be also deduced from Corollary~\ref{c:Desargues=>complete}, Theorem~\ref{t:everywhereP<=>somewhereP} and Proposition~\ref{p:Steiner+projective=>Desargues}).

\begin{corollary} A ($3$-long) projective liner $X$ is Desarguesian if and only if $X$ is everywhere Desarguesian (if and) only if $X$ is somewhere Desarguesian.
\end{corollary}

\chapter{Thalesian liners}

In this section we study Thalesian liners, named in honour of \index[person]{Thales}Thales of Miletus\footnote{{\bf Thales of Miletus} ($Theta\alpha\lambda\tilde\eta\zeta\;o\;M\iota\lambda\acute\eta\sigma \iota o\zeta$, 
624/620 -- 548/545 BC) was an Ancient Greek pre-Socratic philosopher from Miletus in Ionia, Asia Minor. Thales is a founder of the theory of parallels in geometry.
Thales was one of the Seven Sages, founding figures of Ancient Greece, and credited with the saying ``know thyself" which was inscribed on the Temple of Apollo at Delphi.\newline \indent
Many regard him as the first philosopher in the Greek tradition, breaking from the prior use of mythology to explain the world and instead using natural philosophy. He is thus otherwise credited as the first to have engaged in mathematics, science, and deductive reasoning. The first philosophers followed him in explaining all of nature as based on the existence of a single ultimate substance. Thales theorized that this single substance was water. Thales thought the Earth floated on water.\newline 
\indent Thales has been credited with the discovery of five geometric theorems: (1) that a circle is bisected by its diameter, (2) that angles in a triangle opposite two sides of equal length are equal, (3) that opposite angles formed by intersecting straight lines are equal, (4) that the angle inscribed inside a semicircle is a right angle, and (5) that a triangle is determined if its base and the two angles at the base are given. His mathematical achievements are difficult to assess, however, because of the ancient practice of crediting particular discoveries to men with a general reputation for wisdom.
\newline
\indent Thales was said to have calculated the heights of the pyramids and the distance of ships from the shore. In science, Thales was an astronomer who reportedly predicted the weather and a solar eclipse. He was also credited with discovering the position of the constellation Ursa Major as well as the timings of the solstices and equinoxes. Thales was also an engineer; credited with diverting the Halys River.}, who is considered to be a founding father of the deductive geometry.
, and also study the interplay between Moufang and Thalesian liners.

\section{Thalesian liners}

\begin{definition}\label{d:para-Desargues} A liner $X$ is called \index{Thalesian liner}\index{liner!Thalesian}\defterm{Thalesian}
 if $X$ satisfies the \index{Thales Axiom}\index{Axiom!Thales}\defterm{\sf Thales Axiom}:
\begin{itemize}
\item[{\sf (TA)}] {\em for any plane $P\subseteq X$, pairwise disjoint  lines $A,B,C$ in $P$ and any points $a,a'\in A$, $b,b'\in B$, $c,c'\in C$, if $\Aline ab\cap\Aline {a'}{b'}=\varnothing=\Aline bc\cap\Aline {b'}{c'}$, then $\Aline ac\cap\Aline {a'}{c'}=\varnothing$.}
\end{itemize}
\end{definition}

\begin{picture}(200,80)(-120,-10)
\put(0,0){\line(1,0){140}}
\put(145,-3){$A$}
\put(0,30){\line(1,0){140}}
\put(145,27){$B$}
\put(0,60){\line(1,0){140}}
\put(145,57){$C$}

{\linethickness{1pt}
\put(40,0){\color{red}\line(0,1){60}}
\put(40,0){\color{blue}\line(-1,1){30}}
\put(40,60){\color{cyan}\line(-1,-1){30}}

\put(120,0){\color{red}\line(0,1){60}}
\put(120,0){\color{blue}\line(-1,1){30}}
\put(120,60){\color{cyan}\line(-1,-1){30}}
}

\put(40,0){\circle*{3}}
\put(38,-8){$a$}
\put(40,60){\circle*{3}}
\put(38,63){$c$}
\put(10,30){\circle*{3}}
\put(7,33){$b$}
\put(120,0){\circle*{3}}
\put(118,-9){$a'$}
\put(120,60){\circle*{3}}
\put(118,63){$c'$}
\put(90,30){\circle*{3}}
\put(87,33){$b'$}
\end{picture}
 
Thalesian $3$-ranked  liners can be characterized as follows.

\begin{proposition}\label{p:para-Desarguesian<=>} A $3$-ranked liner $X$ is Thalesian if and only if for any plane $P\subseteq X$, any distinct parallel lines $A,B,C$ in $P$ and points $a,a'\in A$, $b,b'\in B$, $c,c'\in C$, if $\Aline ab\parallel \Aline{a'}{b'}$ and $\Aline bc\parallel \Aline{b'}{c'}$, then $\Aline ac\parallel \Aline{a'}{c'}$.
\end{proposition}

\begin{proof} Let $X$ be a $3$-ranked liner.
\smallskip

To prove the ``only if'' part, assume that  $X$ is Thalesian and take any plane $P\subseteq X$, distinct parallel lines $A,B,C\subseteq P$ and  points $a,a'\in A$, $b,b'\in B$, $c,c'\in C$ such that $\Aline ab\parallel \Aline{a'}{b'}$ and $\Aline bc\parallel \Aline{b'}{c'}$. We have to prove that $\Aline ac\parallel \Aline{a'}{c'}$.

If $b=b'$, then  $\Aline ab\parallel \Aline{a'}{b'}$ and $\Aline bc\parallel \Aline{b'}{c'}$ imply $\Aline ab=\Aline{a'}{b'}$ and $\Aline bc=\Aline{b'}{c'}$, according to  Proposition~\ref{p:para+intersect=>coincide}. Then $\{a\}=A\cap\Aline ab=A\cap\Aline{a'}{b'}=\{a'\}$, $\{c\}=C\cap\Aline bc=C\cap\Aline{b'}{c'}=\{c'\}$ and hence $\Aline ac=\Aline {a'}{c'}$ and $\Aline ac\parallel \Aline{a'}{c'}$.

So, assume that $b\ne b'$. In this case $\Aline ab\parallel \Aline {a'}{b'}$ implies $\Aline ab\cap\Aline {a'}{b'}=\varnothing$. Indeed, assuming that $\Aline ab\cap\Aline {a'}{b'}\ne\varnothing$ and applying Proposition~\ref{p:para+intersect=>coincide}, we conclude that $\Aline ab=\Aline{a'}{b'}$ and then $\{b\}=B\cap\Aline ab=B\cap\Aline{a'}{b'}=\{b'\}$, which contradicts the assumption $b\ne b'$. This contradiction shows that $\Aline ab\cap\Aline{a'}{b'}=\varnothing$. By analogy we can show that $\Aline bc\cap\Aline{b'}{c'}=\varnothing$. The Thales Axiom implies $\Aline ac\cap\Aline{a'}{c'}=\varnothing$. By Corollary~\ref{c:parallel-lines<=>}, the disjoint lines $\Aline ac$ and $\Aline{a'}{c'}$ in the plane $P$ are parallel.
\smallskip

To prove the ``if'' part, take any plane $P\subseteq X$, pairwise disjoint  lines $A,B,C\subseteq P$ and points $a,a'\in A$, $b,b'\in B$ and $c,c'\in C$ such that $\Aline ab\cap\Aline {a'}{b'}=\varnothing=\Aline bc\cap\Aline{b'}{c'}$. Since $X$ is $3$-ranked, the coplanar disjoint lines $A,B,C$ are parallel, by Corollary~\ref{c:parallel-lines<=>}. By the same Corollary~\ref{c:parallel-lines<=>}, $\Aline ab\cap\Aline {a'}{b'}=\varnothing=\Aline bc\cap\Aline{b'}{c'}$ implies $\Aline ab\parallel\Aline {a'}{b'}$ and $\Aline bc\parallel\Aline{b'}{c'}$. Applying the condition of the ``if'' part, we conclude that $\Aline ac\parallel \Aline {a'}{c'}$. Assuming that $\Aline ac\cap\Aline{a'}{c'}\ne\varnothing$, we can apply Proposition~\ref{p:para+intersect=>coincide} and conclude that $\Aline ac=\Aline {a'}{c'}$ and hence $\{a\}=A\cap\Aline ac=A\cap\Aline{a'}{c'}=\{a'\}$, which contradicts $\Aline ab\cap\Aline {a'}{b'}=\varnothing$. This contradiction shows that the parallel lines $\Aline ac$ and $\Aline{a'}{c'}$ are disjoint. Therefore, the  liner $X$ is Thalesian. 
\end{proof}

Next we characterize Thalesian Proclus liners.

\begin{proposition}\label{p:Thales<=>planeT} A Proclus liner $X$ is Thalesian if and only if every $3$-long plane in $X$ is a Thalesian liner.
\end{proposition}

\begin{proof} The ``only if'' part is trivial. To prove the ``if'' part, assume that every $3$-long plane in a Proclus liner $X$ is Thalesian. To prove that $X$ is Thalesian, take any plane $P\subseteq X$, disjoint lines $A,B,C$ in $P$ and any points $a,a'\in A$, $b,b'\in B$, $c,c'\in C$ such that $\Aline ab\cap\Aline{a'}{b'}=\varnothing$ and $\Aline bc\cap\Aline{b'}{c'}=\varnothing$. We have to prove that $\Aline ac\cap\Aline{a'}{c'}=\varnothing$. 

Since the Prolcus line $P$ is plane contains disjoint lines, it is not projective. Assuming that the plane $\Pi$ is not $3$-long, we can apply Theorem~\ref{t:Proclus-not-3long} and conclude that $\Pi=Y\setminus \{p\}$ for some Steiner projective plane $Y$ and some point $p\in Y$. It follows from the parallelity of the lines $A,B,C$ that $p\in\overline A\cap\overline B\cap\overline C$. On the other hand, $\Aline ab\cap\Aline {a'}{b'}=\varnothing$ implies $p\in\overline{\{a,b\}}\cap\overline{\{a',b'\}}$. Then $b\in\overline{\{a,p\}}=A$, which is a contradiction showing that the plane $\Pi$ is $3$-long. By our assumption, the $3$-long plane $\Pi$ is Pappian and hence $\Aline ac\cap\Aline{a'}{c'}=\varnothing$. 
\end{proof}

\begin{proposition}\label{p:Thales-regular<=>} A proaffine regular liner $X$ is Thalesian if and only if every distinct parallel lines $A,B,C$ and points $a,a'\in A$, $b,b\in B$, $c,c'\in C$, $\Aline ab\cap\Aline {a'}{b'}=\varnothing=\Aline bc\cap\Aline {b'}{c'}$ implies $\Aline ac\cap\Aline{a'}{c'}=\varnothing$.
\end{proposition}

\begin{proof} The ``only if'' part is trivial. To prove the ``if'' part, assume that every $3$-long plane in a Proclus liner $X$ is Thalesian. To prove that $X$ is Thalesian, take any plane $P\subseteq X$, disjoint lines $A,B,C$ in $P$ and any points $a,a'\in A$, $b,b'\in B$, $c,c'\in C$ such that $\Aline ab\cap\Aline{a'}{b'}=\varnothing$ and $\Aline bc\cap\Aline{b'}{c'}=\varnothing$. We have to prove that $\Aline ac\cap\Aline{a'}{c'}=\varnothing$. 

Since the lines $A,B$ are parallel, the flat hull $\overline{A\cup B}$ is a plane containing two disjoint lines $\Aline ab$ and $\Aline {a'}{b'}$. By Corollary~\ref{c:parallel-lines<=>}, the lines $\Aline ab$ and $\Aline {a'}{b'}$ are parallel. By analogy we can show that $\Aline bc\parallel \Aline {b'}{c'}$. 
If $\Aline ab=\Aline bc$, then $\Aline {a'}{b'}\parallel \Aline ab=\Aline bc\parallel\Aline{b'}{c'}$ and hence $\Aline {a'}{b'}=\Aline{b'}{a'}=\Aline {a'}{c'}$. In this case $\Aline ac\cap\Aline{a'}{c'}=\Aline ab\cap\Aline {a'}{b'}=\varnothing$. By analogy we can show that $\Aline {a'}{b'}=\Aline{b'}{c'}$ implies $\Aline ac\cap\Aline{a'}{c'}\ne\varnothing$. 

So, assume that $\Aline ab\ne\Aline bc$ and $\Aline {a'}{b'}\ne\Aline {b'}{c'}$. Then  $\overline{\{a,b,c\}}$ and $\overline{\{a',b',c'\}}$ are planes in the liner $X$. By Theorem~\ref{t:subparallel-via-base}, these two planes are parallel. If they are disjojnt, then  $\Aline ac\cap\Aline {a'}{c'}\subseteq \overline{\{a,b,c\}}\cap\overline{\{a',b',c'\}}=\varnothing$ and we are done.

So, assume that the parallel planes $\overline{\{a,b,c\}}$ and $\overline{\{a',b',c'\}}$ have a common point and hence they coincide. Applying the Thales Axiom to the plane $\Pi=\overline{A\cup B}=\overline{\{a,b,c\}}=\overline{\{a',b',c'\}}$, we conclude that $\Aline ac\cap\Aline{a'}{c'}=\varnothing$. 
\end{proof}

\begin{theorem}\label{t:ADA=>AMA} Every Desarguesian Proclus liner is Thalesian.
\end{theorem}

\begin{proof} To prove that $X$ is Thalesian, take any plane $P\subseteq X$ and pairwise disjoint lines $A,B,C$ in $P$ and any points $a,a'\in A$, $b,b'\in B$, $c,c'\in C$ such that $\Aline ab\cap \Aline {a'}{b'}=\varnothing=\Aline bc\cap \Aline{b'}{c'}$, which implies $\Aline ab\parallel \Aline{a'}{b'}$ and $\Aline bc\parallel \Aline {b'}{c'}$, by Corollary~\ref{c:parallel-lines<=>}. The same Corollary~\ref{c:parallel-lines<=>}, the disjoint coplanar lines $A,B,C$ are parallel. By Theorem~\ref{t:PD1}, the lines $\Aline ac$ and $\Aline {a'}{c'}$ are parallel.  Assuming that $\Aline ac\cap\Aline {a'}{c'}\ne\varnothing$, we conclude that $\Aline ac=\Aline{a'}{c'}$ and hence $\{a\}=A\cap\Aline ac=A\cap\Aline{a'}{c'}=\{a'\}$, which contradicts $\Aline ab\cap\Aline {a'}{b'}=\varnothing$. This contradiction shows that $\Aline ac\cap\Aline {a'}{c'}=\varnothing$, witnessing that the liner $X$ is Thalesian. 
\end{proof}

The following example found by Ivan Hetman shows that Theorem~\ref{t:ADA=>AMA} cannot be extended to (pro)affine lines.

\begin{example}[Ivan Hetman, 2024]\label{ex:Z15D} The ring $X:=\IZ_{15}$ endowed with the family of lines $$\mathcal L=\big\{x+L:x\in\IZ_{15},\;L\in\{\{0,1,4\},\{0,2,9\},\{0,5,10\}\big\}$$ is a $4$-parallel affine balanced Steiner liner, which is Desarguesian, but not Thalesian.
\end{example}

\begin{proof} It is well-known that the ring $\IZ_{15}$ is isomorphic to the product $\IZ_5\times\IZ_3$, which can be visualized as a $5\times 3$ rectangle. Elements of $\IZ_5\times \IZ_3$ can be labeled by ordered pairs $xy$ of numbers $x\in 5=\{0,1,2,3,4\}$ and $y\in 3=\{0,1,2\}$. The family of lines $\mathcal L$ coincides with all possible shifts of the (disjoint and hence parallel) lines 
$$L_1\defeq\{02,11,22\},\quad L_2\defeq\{01,20,41\},\quad L_3\defeq\{30,31,32\}.$$

\begin{picture}(120,90)(-150,-15)
\multiput(0,0)(0,30){3}{\multiput(0,0)(30,0){5}{\color{gray}\circle*{3}}}

\put(0,60){\color{cyan}\line(1,-1){30}}
\put(30,30){\color{cyan}\line(1,1){30}}

\put(0,30){\color{blue}\line(2,-1){60}}
\put(60,0){\color{blue}\line(2,1){60}}
\put(90,0){\color{magenta}\line(0,1){60}}
\put(87,64){\color{magenta}$L_3$}
\put(125,26){\color{blue}$L_2$}
\put(-15,56){\color{cyan}$L_1$}

\put(0,60){\color{cyan}\circle*{3}}
\put(30,30){\color{cyan}\circle*{3}}
\put(60,60){\color{cyan}\circle*{3}}

\put(0,30){\color{blue}\circle*{3}}
\put(60,0){\color{blue}\circle*{3}}
\put(120,30){\color{blue}\circle*{3}}

\put(90,0){\color{magenta}\circle*{3}}
\put(90,30){\color{magenta}\circle*{3}}
\put(90,60){\color{magenta}\circle*{3}}

\end{picture}

It can be shown that the Steiner liner $X$ is affine, ranked, balanced, and $4$-parallel. To see that the liner $X$ is not Thalesian, consider the points
$$a\defeq 12,\quad a'\defeq 20,\quad b\defeq 31,\quad b'\defeq 42,\quad c\defeq 41,\quad c'\defeq 21$$ and observe that the lines $$\Aline a{a'}=\{00,12,20\},\quad \Aline b{b'}=\{22,32,42\},\quad \Aline c{c'}=\{21,30,41\}$$ are disjoint and hence parallel in the ranked liner $X$. Also $$\Aline ab\cap\Aline{a'}{b'}=\{12,31,02\}\cap \{20,42,10\}=\varnothing,\quad \Aline bc\cap\Aline{b'}{c'}=\{31,41,10\}\cap\{42,21,02\}=\varnothing,$$
but $$\Aline ac\cap\Aline {a'}{c'}=\{12,41,22\}\cap\{20,21,22\}=\{22\}\ne\varnothing,$$witnessing that the liner $X$ is not Thalesian.

\begin{picture}(120,80)(-150,-5)
\multiput(0,0)(0,30){3}{\multiput(0,0)(30,0){5}{\color{gray}\circle*{2}}}

\put(0,0){\color{brown}\line(1,2){30}}
\put(30,60){\color{brown}\line(1,-2){30}}
\put(60,60){\color{olive}\line(1,-1){30}}
\put(90,30){\color{olive}\line(1,1){30}}

\put(60,30){\color{orange}\line(1,-1){30}}
\put(90,0){\color{orange}\line(1,1){30}}

{\linethickness{1pt}
\put(0,60){\color{cyan}\line(1,0){30}}
\put(30,60){\color{cyan}\line(2,-1){60}}

\put(30,0){\color{cyan}\line(1,0){30}}
{\color{cyan}
\qbezier(60,0)(105,15)(120,60)
}

\put(30,0){\color{blue}\line(2,1){60}}
\put(90,30){\color{blue}\line(1,0){30}}

\put(60,30){\color{blue}\line(2,1){60}}
\put(60,30){\color{blue}\line(-2,1){60}}

\put(30,60){\color{red}\line(1,0){30}}
\put(60,60){\color{red}\line(2,-1){60}}
\put(60,0){\color{magenta}\line(0,1){60}}
}
\put(0,60){\circle*{2}}
\put(30,0){\circle*{2}}

\put(30,60){\color{brown}\circle*{3}}
\put(27,64){\color{brown}$a$}
\put(60,0){\color{brown}\circle*{3}}
\put(57,-10){\color{brown}$a'$}
\put(0,0){\color{brown}\circle*{3}}

\put(90,30){\color{olive}\circle*{3}}
\put(88,20){\color{olive}$b$}
\put(120,60){\color{olive}\circle*{3}}
\put(123,57){\color{olive}$b'$}
\put(60,60){\color{magenta}\circle*{4}}
\put(60,60){\color{red}\circle*{2}}

\put(90,0){\color{orange}\circle*{3}}
\put(120,30){\color{orange}\circle*{3}}
\put(123,27){\color{orange}$c$}
\put(60,30){\color{orange}\circle*{3}}
\put(65,26){\color{orange}$c'$}

\end{picture}

\end{proof}

%

The following corollary of the Thales Axiom is called the Inverse Thales Theorem.

\begin{theorem}\label{t:Inverse-Thales} Let $abc,a'b'c'$ be two disjoint triangles in a Thalesian Para-Playfair liner $X$ such that $\Aline ab\parallel \Aline{a'}{b'}$, $\Aline bc\parallel \Aline{b'}{c'}$ and $\Aline ac\parallel \Aline{a'}{c'}$. If $\Aline a{a'}\cap\Aline b{b'}=\varnothing$, then $\Aline {a}{a'}\parallel \Aline c{c'}\parallel \Aline b{b'}$.
\end{theorem}

\begin{proof}   By Proposition~\ref{p:Playfair=>para-Playfair=>Proclus}, Theorem~\ref{t:Proclus<=>} and Proposition~\ref{p:k-regular<=>2ex}, the para-Playfair liner $X$ is Proclus, $3$-proregular and $3$-ranked. Since $\Aline ab\parallel \Aline {a'}{b'}$, the flat $\overline{\{a,b,a',b'\}}$ is a plane containing the disjoint lines $\Aline a{a'}$ and $\Aline b{b'}$, which are parallel, by Corollary~\ref{c:parallel-lines<=>}. 

If $c$ or $c'$ belongs to the line $\Aline a{a'}$, then $\Aline ac\parallel \Aline {a'}{c'}$ implies $\Aline c{c'}=\Aline a{a'}\parallel \Aline b{b'}$.
If $c$ or $c'$ belongs to the line $\Aline b{b'}$, then $\Aline ac\parallel \Aline {a'}{c'}$ implies $\Aline c{c'}=\Aline b{b'}\parallel \Aline a{a'}$.

So, assume that $c,c'\notin\Aline a{a'}\cup\Aline b{b'}$. First we prove that $\Aline a{a'}\cap\Aline c{c'}=\varnothing$. To derive a contradiction, assume that $\Aline a{a'}\cap\Aline c{c'}$ contains some point $o$. The parallelity relation $\Aline bc\parallel \Aline {b'}{c'}$ implies that $\Pi\defeq \overline{\{b,c,b',c'\}}$ is a plane containing the point $o\in \Aline c{c'}\cap\Aline a{a'}$. It follows from $\Aline a{a'}\parallel \Aline b{b'}$ that $\Aline a{a'}\subseteq \overline{\{b,b',o\}}\subseteq\Pi$. Therefore, the plane $\Pi$ contains the points $a,a',b,b',c,c',o$ and hence $|\Pi|\ge 7$. Applying Theorem~\ref{t:Proclus-not-3long}, we conclude that the Proclus plane $\Pi$ is $3$-long and being $3$-proregular, $\Pi$ is regular.  Since $X$ is para-Playfair, there exists a line $L\subseteq \Pi$ such that $c\in L\subseteq\Pi\setminus \Aline a{a'}$ and hence $L\parallel \Aline a{a'}$. By Proposition~\ref{p:Proclus-Postulate}, there exists a point $c''\in L\cap\Aline {a'}{c'}$. By Theorem~\ref{t:Proclus-lines}, $\Aline {c}{c''}=L\parallel \Aline a{a'}\parallel \Aline b{b'}$ implies $\Aline {c}{c''}\parallel \Aline{b}{b'}$.  Therefore, $\Aline a{a'},\Aline b{b'}$ and $\Aline{c}{c''}$ are three parallel lines in the plane $\Pi$.  Since $\Aline ab\parallel \Aline {a'}{b'}$, $\Aline ac\parallel \Aline {a'}{c'}=\Aline{a'}{c''}$, we can apply  Proposition~\ref{p:para-Desarguesian<=>} and conclude that $\Aline bc\parallel\Aline {b'}{c''}$. Since $X$ is Proclus,  $\Aline {b'}{c''}\parallel \Aline bc\parallel \Aline {b'}{c'}$ implies $\Aline {b'}{c''}=\Aline {b'}{c'}$. Since $a'b'c'$ is a triangle, $c''\in \Aline {b'}{c''}\cap\Aline {a'}{c''}=\Aline {b'}{c'}\cap\Aline{a'}{c'}=\{c'\}$ and hence $\Aline c{c'}=\Aline c{c''}=L\subseteq \Pi\setminus \Aline a{a'}$, which contradicts $o\in\Aline a{a'}\cap \Aline c{c'}$. This contradiction shows that $\Aline a{a'}\cap\Aline c{c'}=\varnothing$. It follows from $\Aline ac\parallel \Aline{a'}{c'}$ that $\overline{\{a,c,a',c'\}}$ is a plane containing the disjoint lines $\Aline a{a'}$ and $\Aline c{c'}$. By Corollary~\ref{c:parallel-lines<=>}, these disjoint lines are parallel.   By analogy we can prove that $\Aline b{b'}\parallel \Aline c{c'}$.
\end{proof}




\begin{remark} Theorems~\ref{t:Desargues-affine} and \ref{t:ADA=>AMA} imply that every proaffine regular liner of dimension $\ge 3$ is Desarguesian and hence Thalesian. On the other hand, there are many examples of (finite) Thalesian affine regular planes, which are not Desarguesian, see \cite{HTP}. Thalesian affine regular planes are known in the literature as translation affine  planes. By Veblen-Wedderburn Theorem \ref{t:VW-Thalesian<=>quasifield} such planes are isomorphic to affine planes over quasi-fields.
\end{remark}   

\section{Moufang projective liners are everywhere Thalesian}

In this section we reveal the interplay between Moufang projective liners and Thalesian affine liners showing that Moufang projective liners are exactly everywhere Thalesian projective liners. Let us recall that for a property $\mathcal P$ of affine regular liners, a projective liner $Y$ is defined to be \defterm{everywhere $\mathcal P$} if for every hyperplane $H\subseteq Y$, the affine regular liner $Y\setminus H$ has property $\mathcal P$.

\begin{theorem}\label{t:Moufang<=>everywhere-Thalesian} A projective liner $Y$ is Moufang if and only if for every hyperplane $H$ in $Y$, the affine liner $Y\setminus H$ is Thalesian.
\end{theorem}

\begin{proof} Assume that a projective liner $Y$ is Moufang and take any hyperplane $H$ in $Y$. By Theorem~\ref{t:affine<=>hyperplane}, the subliner $X\defeq Y\setminus H$ of $Y$ is affine and regular. We have to prove that the liner $X$ is Thalesian. Take any plane $\Pi\subseteq X$, disjoint lines $A,B,C\subset \Pi$ and distinct points $a,a'\in A$, $b,b'\in B$, $c,c'\in C$ such that $\Aline ab\cap\Aline{a'}{b'}=\varnothing=\Aline bc\cap\Aline{b'}{c'}$. We have to prove that $\Aline ac\cap\Aline{a'}{c'}=\varnothing$. Since the affine plane $\Pi$ contains three disjoint lines $A,B,C$, it is $3$-long, and by Theorem~\ref{t:Playfair<=>}, the $3$-long affine regular liner $X$ is Playfair. Since $\Pi$ is an affine plane in $X$, its flat hull $\overline \Pi$ is a plane in the projective liner $Y$. By Theorem~\ref{t:affine<=>hyperplane}, the intersection $\overline\Pi\cap H$ is a line $\overline \Pi$. For a line $L$ in the plane $\Pi$, let $\overline L$ denote the flat hull of $L$ in $Y$. By the projectivity of the plane $\overline \Pi$, the liners $\overline L$ and $H$ have a unique commopn point. To distinguish lines in the liners $X$ and $Y$, for two distinct points $x,y\in Y$ we shall denote by $\overline{\{x,y\}}$ the line in $Y$ that contains the points $x,y$. 

If $b\in\Aline ac$, then $\Aline ab\cap\Aline{a'}{b'}=\varnothing=\Aline bc\cap\Aline{b'}{c'}$ implies $\Aline {a'}{b'}\parallel \Aline ab=\Aline bc\parallel \Aline {b'}{c'}$ and hence $\Aline {a'}{b'}=\Aline {b'}{c'}$, by the Proclus property of the Playfair plane $\Pi$. In this case $\Aline ac\cap\Aline{a'}{c'}=\Aline ab\cap\Aline{a'}{b'}=\varnothing$. By analogy we can show that $b'\in\Aline {a'}{c'}$ implies $\Aline ac\cap\Aline{a'}{c'}=\varnothing$.

So, assume that $b\notin\Aline ac$ and $b'\notin\Aline{a'}{c'}$, which means that $abc$ and $a'b'c'$ are two triangles in the plane $\Pi$.  By the projectivity of the plane $\overline X$, there exists a point $\beta'\in\overline A\cap\overline C$. Since the lines $A,C$ are disjoint, the point $\beta'$ belongs to the hyperplane $H$. By the projectivity of the plane, $\overline \Pi$, we also have $\varnothing\ne \overline A\cap\overline B\subseteq \overline A\setminus X=\{\beta'\}$, which implies that $\overline A\cap\overline B\cap\overline C=\{\beta'\}$. By the projectivity of $\overline \Pi$, there exist points $\alpha\in \overline{\{b,c\}}\cap\overline{\{b',c'\}}$, $\beta\in \overline{\{a,c\}}\cap\overline{\{a',c'\}}$ and $\gamma\in  \overline{\{a,b\}}\cap\overline{\{a',b'\}}$. It follows from $\varnothing=\Aline bc\cap \Aline{b'}{c'}=\overline{\{b,c\}}\cap X\cap\overline{\{b',c'\}}=X\cap\{\alpha\}$ that $\alpha\in \overline \Pi\setminus X=\overline \Pi\cap H$. By analogy we can show that $\gamma\in \overline\Pi\cap H$. Since $abc$ and $a'b'c'$ are disjoint triangles, the points $\alpha,\gamma\in\overline \Pi\cap H$ are distinct, which implies that $\overline \Pi\cap H$ is a line in the projective plane $\overline\Pi$ and hence $\beta'\in \overline\Pi\cap H=\overline{\{\alpha,\gamma\}}$. 



Now consider the triangles $\alpha b\gamma$ and $a'\beta'c'$ in the plane $\overline \Pi$ and observe that these two triangles are perspective from the point $b'$. Since the liner $Y$ is Moufang, the set
$$T\defeq (\overline{\{\alpha,b\}}\cap\overline{\{a',\beta'\}})\cup  (\overline{\{b,\gamma\}}\cap\overline{\{\beta',c'\}})\cup(\overline{\{\alpha,\gamma\}}\cap\overline{\{a',c'\}})=\{a,c\}\cup(\overline{\{\alpha,\gamma\}}\cap\overline{\{a',c'\}})$$ has rank $\|T\|=2$ and hence
$\overline{\{\alpha,\gamma\}}\cap\overline{\{a',c'\}}\subseteq\overline{\{a,c\}}$ and 
$$\varnothing\ne \overline{\{\alpha,\gamma\}}\cap\overline{\{a',c'\}}= \overline{\{\alpha,\gamma\}}\cap\overline{\{a',c'\}}\cap\overline{\{a,c\}}=(\overline\Pi\cap H)\cap\{\beta\},$$ which implies $\beta\in H$ and finally 
$$\Aline ac\cap\Aline{a'}{c'}=\overline{\{a',c'\}}\cap\overline{\{a,c\}}\cap X=\{\beta\}\cap X\subseteq H\cap X=\varnothing.$$
\smallskip

Next, assume that for every hyperplane $H$ in $Y$, the liner $X\setminus H$ is Thalesian. To prove that the projective liner $Y$ is Moufang, take any centrally perspective disjoint triangles $abc$ and $a'b'c'$ in $Y$ with $b'\in\Aline ac$. We have to prove that the set $$T\defeq(\Aline ab\cap\Aline{a'}{b'})\cup (\Aline bc\cap\Aline{b'}{c'})\cup(\Aline ac\cap\Aline{a'}{c'})$$ has rank $\|T\|\in\{0,2\}$ in the liner $Y$. 
Since the triangles $abc$ and $a'b'c'$ are centrally perspective, there exists a unique point $o\in\Aline a{a'}\cap\Aline b{b'}\cap\Aline c{c'}$. Consider the plane $\Pi\defeq\overline{\{a,o,c\}}$ and observe that $b'\in \Aline ac\subseteq \Pi$ and $\{a',b,c'\}\subseteq \Aline oa\cup\Aline o{b'}\cup \Aline oc\subseteq\Pi$.  
Let $F$ be a maximal $3$-long flat in $Y$ that contains the point $o$. By Lemma~\ref{l:ox=2}, for every points $x\in F$ and $y\in Y\setminus F$, the line $\Aline xy$ coincides with the doubleton $\{x,y\}$. Taking into account that $|\Aline oa|\ge|\{o,a,a'\}|\ge 3$ and $|\Aline oc|\ge|\{o,c,c'\}|\ge 3$, we conclude that $a,c\in F$ and hence the plane $\Pi\subseteq F$ is $3$-long. By Corollary~\ref{c:Avogadro-projective}, the $3$-long projective plane  $\Pi$ is $2$-balanced. If $|\Pi|_2=3$, then the Steiner projective plane $\Pi$ is Desarguesian (by Proposition~\ref{p:Steiner+projective=>Desargues}) and hence $\|T\|\in\{0,2\}$. So, assume that $|\Pi|_2\ge 4$. 

By the Kuratowski--Zorn Lemma, the independent set $\{o,a,c\}$ can be enlarged to a maximal independent set $M$ in the projective liner $Y$. The Exchange Property of the projective liner $Y$ giaranees that the flat $H\defeq\overline{M\setminus\{o\}}$ is a hyperplane in $Y$. By our assumption, the subliner $Y\setminus H$ of $Y$ is Thalesian and then the flat $P\defeq \Pi\setminus H$ in $Y\setminus H$ also is Thalesian. Since $|\Pi|_2\ge 4$, the Thalesian liner $P=\Pi\setminus H$ is $3$-long. By Theorems~\ref{t:affine<=>hyperplane}, \ref{t:Playfair<=>}, and Corollary~\ref{c:procompletion-rank}, the $3$-long liner $P$ is a Playfair plane. The independence of the set $M$ implies $o\in \Pi\setminus H=P$ and hence $$\{a',b,c'\}\in (\Aline oa\setminus\{a\})\cup(\Aline o{b'}\setminus\{b'\})\cup(\Aline oc\setminus\{c\})\subseteq \Pi\setminus H=P.$$ By the projectivity of the plane $\Pi$, there exist unique points $x\in \Aline ab\cap\Aline {a'}{b'}$, $y\in\Aline bc\cap\Aline {b'}{c'}$ and $z\in \Aline ac\cap\Aline{a'}{c'}$. It is easy to see that $x\ne b'\ne y$ and hence $\{x,y\}\subseteq(\Aline {a'}{b'}\setminus\{b'\})\cup(\Aline{b'}{c'}\setminus\{b'\})\subseteq\Pi\setminus H=P$.

\begin{picture}(300,150)(-100,-15)
\linethickness{=0.6pt}
\put(0,0){\color{red}\line(1,0){240}}
\put(60,120){\color{violet}\line(1,-2){60}}
\put(60,120){\color{blue}\line(-1,-2){60}}
\put(60,120){\color{cyan}\line(0,-1){120}}
\put(30,60){\color{red}\line(7,-2){210}}
\put(30,60){\color{cyan}\line(1,-2){30}}
\put(100,40){\color{cyan}\line(-1,-1){40}}
\put(60,40){\color{violet}\line(3,-2){60}}
\put(60,40){\color{blue}\line(-3,-2){60}}
\put(45,30){\color{red}\line(13,-2){195}}

\put(84,24){\circle*{3}}
\put(81,30){$y$}
\put(45,30){\circle*{3}}
\put(44,34){$x$}
\put(0,0){\color{blue}\circle*{3}}
\put(-4,-8){\color{blue}$a$}
\put(60,0){\color{cyan}\circle*{3}}
\put(57,-10){\color{cyan}$b'$}
\put(120,0){\color{violet}\circle*{3}}
\put(118,-8){\color{violet}$c$}
\put(240,0){\color{red}\circle*{3}}
\put(240,-8){\color{red}$z$}
\put(30,60){\circle*{3}}
\put(22,60){$a'$}
\put(100,40){\circle*{3}}
\put(101,42){$c'$}
\put(60,40){\circle*{3}}
\put(62,40){$b$}
\put(60,120){\circle*{3}}
\put(58,123){$o$}

\end{picture}

Consider the triangles $xby$ and $a'oc'$ in the Playfair plane $P$ and observe that $b'\in \Aline x{a'}\cap\Aline bo\cap\Aline y{c'}$ implies that $\Aline x{a'}\cap P$, $\Aline bo\cap P$, $\Aline y{c'}\cap P$ are pairwise disjoint lines in the plane $P$. Since the plane $P$ is Thalesian, $\Aline xb\cap\Aline{a'}o\cap P=\{a\}\cap P=\varnothing=\{c\}\cap P=\Aline by\cap\Aline o{c'}\cap P$, implies $\Aline xy\cap \Aline{a'}{c'}\cap P=\varnothing$. The projectivity of the plane $\Pi$ ensures that $\varnothing\ne \Aline xy\cap \Aline{a'}{c'}\subseteq\Pi\setminus P=\Pi\cap H=\Aline ac$ and hence $\varnothing \ne\Aline xy\cap\Aline{a'}{c'}=\Aline xy\cap\Aline{a'}{c'}\cap\Aline ac=\Aline xy\cap\{z\}$ and finally, $z\in \Aline xy$. Then the set $$T=(\Aline ab\cap\Aline {a'}{b'})\cup(\Aline bc\cap\Aline {b'}{c'})\cup(\Aline ac\cap\Aline{a'}{c'})=\{x,y,z\}$$has rank $\|T\|=\|\{x,y\}\|=2$, witnessing that the projective liner $Y$ is Moufang.
\end{proof}

\begin{remark} By Theorem~\ref{t:Moufang<=>everywhere-Thalesian}, a projective liner is Moufang if and only if it is everywhere Thalesian. Somewhere Thalesian projective planes are called \defterm{translation projective planes}, see Theorem~\ref{t:paraD<=>translation}. We shall meet translation projective planes in the Lenz  classification of projective planes, presented in Section~\ref{s:Lenz}.
\end{remark}

\section{The Moulton plane}\label{s:Moulton} 

The Moulton plane is an example of an affine plane which is not Thalesian (even not uno-Thalesian). The Moulton plane is named after Forest Ray Moulton, american astronomer, who suggested this example in \cite{Moulton1902}.

The Moulton plane is the set $X\defeq\IR\times\IR$ endowed with the family of lines
$$\mathcal L\defeq\{L_{a,b}:a,b\in\IR\}\cup\{L_c:c\in\IR\}$$where 
$$L_{a,b}\defeq\begin{cases}\{(x,ax+b):x\in\IR\}&\mbox{if $a\ge 0$};\\
\{(x,\tfrac12ax+b):x\le 0\}\cup\{(x,ax+b):x\ge 0\}&\mbox{if $a<0$};
\end{cases}$$
and $L_c\defeq\{(c,y):y\in\IR\}$ for real numbers $a,b,c\in\IR$. 

In the following picture we draw the lines $L_{1,6}$, $L_{-1,8}$ and $L_4$ on the Moulton plane.

\begin{picture}(300,180)(-100,-10)

\put(100,0){\color{red}\vector(0,1){160}}
\put(0,20){\color{red}\vector(1,0){200}}
{\linethickness{0.5pt}
\put(100,100){\line(-2,1){90}}
\put(40,112){$L_{-1,8}$}
\put(100,100){\line(1,-1){100}}
\put(40,41){$L_{1,6}$}
\put(140,0){\line(0,1){150}}
\put(20,0){\line(1,1){145}}
}
\put(143,85){$L_4$}

\put(92,10){$\mathbf 0$}
\put(93,145){$y$}
\put(195,12){$x$}
\end{picture}

It is easy to check that the family $\mathcal L$ satisfies the axioms {\sf (L1)} and {\sf(L2)} of Theorem~\ref{t:L1+L2} and hence determines a unique line relation $\Af$ on the set $X$, whose family of lines coincides with $\mathcal L$. 

\begin{exercise}\label{ex:Moulton-affine} Show that the Moulton plane is an affine regular liner of rank $\|X\|=3$.
\end{exercise}

\begin{exercise}\label{ex:Moulton} Show that the Moulton plane is not Thalesian.
\smallskip

\noindent{\em Hint:} To see that the Moulton plane is not Thalesian, look at the following picture, at which the lines $A,B,C$ are parallel,  $\Aline ab\parallel \Aline{a'}{b'}$, $\Aline bc\parallel\Aline {b'}{c'}$ but $\Aline ac\nparallel \Aline{a'}{c'}$.
\end{exercise}

\begin{picture}(300,180)(-100,-40)

\put(100,-20){\color{green}\vector(0,1){150}}

\put(0,20){\line(1,0){200}}
\put(0,40){\line(1,0){200}}
\put(0,60){\line(1,0){200}}

{\linethickness{1pt}
\put(100,20){\color{green}\line(0,1){40}}
\put(100,20){\color{blue}\line(1,1){20}}
\put(120,40){\color{red}\line(-1,1){20}}

\put(80,20){\color{green}\line(0,1){40}}
\put(80,20){\color{blue}\line(1,1){20}}
\put(100,40){\color{red}\line(-1,1){20}}
}

\put(100,40){\color{red}\line(1,-2){30}}
\put(100,40){\color{red}\line(-1,1){70}}
\put(60,-20){\color{blue}\line(1,1){120}}
\put(40,-20){\color{blue}\line(1,1){120}}
\put(80,-20){\color{green}\line(0,1){140}}

\put(100,60){\color{red}\line(-2,1){70}}
\put(120,40){\color{red}\line(1,-1){60}}

\put(-10,17){$B$}
\put(-10,37){$C$}
\put(-10,57){$A$}

\put(80,20){\circle*{3}}
\put(73,22){$b$}
\put(100,20){\circle*{3}}
\put(102,11){$b'$}
\put(80,60){\circle*{3}}
\put(73,53){$a$}
\put(100,60){\circle*{3}}
\put(102,62){$a'$}
\put(100,40){\circle*{3}}
\put(103,34){$c$}
\put(120,40){\circle*{3}}
\put(118,44){$c'$}
\end{picture}

The construction of the Moulton plane can be generalized as follows. Take any strictly increasing bijections $f,g:\IR\to\IR$ of the real line. The \index{$(f,g)$-Moulton plane}\defterm{$(f,g)$-Moulton plane} is the plane $\IR\times\IR$ endowed with the family of lines $\mathcal L_{f,g}\defeq\{L_{a,b}:a,b\in \IR\}\cup\{L_c:c\in\IR\}$ where 
$$L_{a,b}\defeq\begin{cases}\{(x,ax+b):x\in\IR\}&\mbox{if $a\ge 0$};\\
\{(f(x),g(ax+b)):x\in\IR\}&\mbox{if $a<0$};
\end{cases}
$$
and $L_c\defeq\{(c,y):y\in\IR\}$ for real numbers $a,b,c\in\IR$.

\begin{exercise} Check that $(\IR\times\IR,\mathcal L_{f,g})$ is a Playfair plane.
\end{exercise}

\begin{exercise} Explore the interplay between the properties of the $(f,g)$-Moulton plane $(\IR\times\IR,\mathcal L_{f,g})$ and the properties of the increasing bijections $f,g$.
\end{exercise}

\section{Two non-Desarguesian planes of order 9}\label{s:J9xJ9}

In this section we present two examples of non-Desarguesian  affine planes of the smallest possible order $9$. One of them is Thalesian and the other is not. Both planes are constructed with the help of  the near-field $\mathbb J_9$, defined as follows.

Let $\mathbb F_3=\{-1,0,1\}$ be the $3$-element field and $\mathbb J_9\defeq\mathbb F_3^2$ be the two-dimensional vector space over the field $\IF_3$. In the vector space \index[note]{$\mathbb J_9$}$\mathbb J_9$, consider the vectors
$$\mathbf 0\defeq (0,0),\quad \mathbf 1\defeq(1,0),\quad\boldsymbol i\defeq (0,1),\quad \boldsymbol j\defeq(1,1),\quad\boldsymbol k\defeq(1,-1).$$
Then 
$$
\begin{aligned}
\mathbb J_9&=\{(-1,-1),(-1,0),(-1,1),(0,-1),(0,0),(0,1),(1,-1),(1,0),(1,1)\}\\
&=\{-\mathbf 1-\boldsymbol i,-\mathbf 1,-\mathbf 1+\boldsymbol i,-\boldsymbol i,\mathbf 0,\boldsymbol i,\mathbf 1-\boldsymbol i,\mathbf 1,\mathbf 1+\boldsymbol i\}\\
&=\{-\boldsymbol j,-\mathbf 1,-\boldsymbol k,-\boldsymbol i,\mathbf 0,\boldsymbol i,\boldsymbol k,\mathbf 1,\boldsymbol j\}
\end{aligned}
$$
Endow the set $\mathbb J_9$ with the operation of multiplication of quaternion units, which is uniquely determined by the famous Hamilton formula $$\boldsymbol i^2=\boldsymbol j^2=\boldsymbol k^2=\boldsymbol i\boldsymbol j\boldsymbol k=-\mathbf 1.$$

\begin{exercise} Show that $(x+y)\cdot z=x\cdot z+y\cdot z$ for every $x,y,z\in\mathbb J_9$.
\end{exercise}  

Consider the liner $H\defeq\mathbb J_9\times\mathbb J_9$ endowed with the family of lines 
$$\mathcal L\defeq\{L_{a,b}:a,b\in\mathbb J_9\}\cup\{L_c:c\in\mathbb J_9\}$$
where $L_{a,b}=\{(x,xa+b):x\in\mathbb J_9\}$ and $L_c\defeq\{(c,y):y\in\mathbb J_9\}$.
 
\begin{exercise} Show that the liner $H$ is a Thalesian affine regular liner.
\end{exercise}

To see that the liner $H$ is not Desarguesian, consider the points  $o\defeq(\mathbf 0,\mathbf 0)$ and 
$$a\defeq (\mathbf 1,\mathbf 0),\;\;b\defeq(\mathbf 1,\boldsymbol i),\;\;c\defeq(\boldsymbol i,\boldsymbol i),\;\;a'\defeq(\boldsymbol j,\mathbf 0),\;\;b'\defeq(\boldsymbol j,-\boldsymbol k),\;\;c'\defeq(-\boldsymbol{k},-\boldsymbol k).
$$
\begin{exercise} Show that 
$$\Aline a{a'}=\Aline oa=\Aline o{a'}=L_{\mathbf 0,\mathbf 0},\quad \Aline b{b'}=\Aline ob=\Aline o{b'}=L_{\boldsymbol i,\mathbf 0},\quad \Aline c{c'}=\Aline oc=\Aline o{c'}= L_{\mathbf 1,\mathbf 0},$$which implies that  the triangles $abc$ and $a'b'c'$ are perspective from the point $o\defeq(\mathbf 0,\mathbf 0)$.
\end{exercise}

\begin{exercise} Show that $\Aline ab=L_{\mathbf 1}\parallel L_{\boldsymbol j}= \Aline{a'}{b'}$ and $\Aline bc=L_{\mathbf 0,\boldsymbol i}\parallel L_{\mathbf 0,-\boldsymbol k}=\Aline{b'}{c'}$.
\end{exercise} 

\begin{exercise} Show that $\Aline ac=L_{\boldsymbol j,-\boldsymbol j}$, $\Aline {a'}{c'}=L_{-\boldsymbol k,\boldsymbol i}$, and hence $\Aline ac\nparallel \Aline {a'}{c'}$.
\end{exercise}

The lines $\Aline oa,\Aline ob,\Aline oc,\Aline ab,\Aline bc,\Aline ac,\Aline{a'}{b'},\Aline{b'}{c'},\Aline{a'}{c'}$ are drawn in the following picture showing that the affine liner $H$ is not Desarguesian.



This Thalesian affine plane $H$ has a dual plane $H'$, which is not Thalesian. This dual plane $H'$ is the set $\mathbb J_9\times\mathbb J_9$ endowed with the family of lines $$\mathcal L'\defeq\{L'_{a,b}:a,b\in \mathbb J_9\}\cup\{L'_c:c\in\mathbb J_9\},$$ where $L'_{a,b}\defeq\{ax+b:x\in\mathbb J_9\}$ and $L'_c\defeq\{(c,y):y\in\mathbb J_9\}$.

\begin{exercise} Show that the liner $H'\defeq (\mathbb J_9\times\mathbb J_9,\mathcal L')$ is affine and regular.
\end{exercise}

To see that the affine liner $H'$ is not Thalesian, consider the points 
$$
\begin{aligned}
& a\defeq (\mathbf 0,\mathbf 1),\quad b\defeq(\mathbf 0,\boldsymbol 0),\quad 
c\defeq(\boldsymbol i,\boldsymbol i),\quad\mbox{and}\\
&a'\defeq(\boldsymbol k,\mathbf 1),\quad b'\defeq(\boldsymbol k,\mathbf 0),\quad c'\defeq(1,\boldsymbol i).
\end{aligned}
$$

\begin{exercise} Show that
$$\Aline a{a'}=L'_{\mathbf 0,\mathbf 1},\quad \Aline b{b'}=L'_{\mathbf 0,\mathbf 0},\quad \Aline c{c'}= L'_{\mathbf 0,\boldsymbol i},$$which implies that the lines $\Aline a{a'},\Aline b{b'},\Aline c{c'}$ are parallel.
\end{exercise}

\begin{exercise} Show that $\Aline ab=L'_{\mathbf 0}\parallel L'_{\mathbf 1}= \Aline{a'}{b'}$ and $\Aline bc=L'_{\mathbf 1,\mathbf 0}\parallel L'_{\mathbf 1,-\boldsymbol k}=\Aline{b'}{c'}$.
\end{exercise} 

\begin{exercise} Show that $\Aline ac=L'_{\boldsymbol j,\mathbf 1}$, $\Aline {a'}{c'}=L'_{-\boldsymbol i,-\boldsymbol i}$, and hence $\Aline ac\nparallel \Aline {a'}{c'}$.
\end{exercise}

The lines $\Aline aa',\Aline bb',\Aline cc',\Aline ab,\Aline bc,\Aline ac,\Aline{a'}{b'},\Aline{b'}{c'},\Aline{a'}{c'}$ are drawn in the following picture showing that the affine liner $H'$ is not Thalesian.



The non-Desarguesian affine planes $H'$ and $H$ are known in Geometry as  \index{Hall plane}\index{plane!Hall}\index{dual Hall plane}\index{plane!dual Hall}\index[person]{Hall}\defterm{the Hall planes} and \defterm{the dual Hall plane} of order 9. We shall refer to the planes $H$ and $H'$ as the \defterm{right-Hall} and \defterm{left-Hall affine planes of order $9$}. So, the right-Hall affine plane of order 9 is Thalesian but not Desarguesian, and the left-Hall affine plane of order 9 is not Thalesian. The projective completions of the right-Hall and left-Hall affine planes of order 9 are called the \defterm{right-Hall} and \defterm{left-Hall projective planes of order $9$}. Those two projective planes are not Desarguesian, by Theorem~\ref{t:Desarguesian-liner<=>}. 

In the tables throughout the book, the right-Hall and left-Hall affine planes of order 9 are called {\tt Thales} and {\tt Hall}, respectively. The unique Desarguesian affine plane of order 9 is called {\tt Desarg}. The unique non-Thalesian affine plane whose projective completion is the right-Hall projective plane of order 9 is called {\tt dhall}. By {\tt hall} we denote the affine plane which is not isomorphic to the left-Hall affine plane {\tt Hall} but has projective completion, which is the left-Hall projective plane of order $9$. The planes {\tt Hall} and {\tt hall} have automorphism groups of orders  31040 and 3840, respectively (so the plane {\tt Hall} has larger automorphism group comparing to the plane {\tt hall}). Finally, two affine planes that have projective completions isomorphic to the the Hughes projective plane are called {\tt Huges} and {\tt hughes}. Those affine planes have automorphism groups of order 2592 and 432, respectively.
The following table lists all seven affine planes of order 9 together with the cardinalities of their automorphism groups (calculated by Ivan Hetman).
$$
\begin{array}{|c|c|c|c|c|c|c|c|}
\hline
\mbox{Affine plane $\Pi$:}&\mbox{\tt Desarg}&\mbox{\tt Thales}&\mbox{\tt Hall}&\mbox{\tt hall}&\mbox{\tt dhall}&\mbox{\tt Hughes}&{\tt hughes}\\
\hline
|\Aut(\Pi)|:&933120&311040&311040&3840&3456&2592&432\\
\hline
\end{array}
$$

\chapter{Special automorphisms of liners} 

In this chapter we study some special automorphism of liners: dilation, translation, central, paracentral, and hyperfixed. Let us recall that an \defterm{automorphism} of a liner $X$ is any isomorphism $A:X\to X$ of $X$ to itself. The set $\Aut(X)$ of all automorphisms of a liner $X$ is a group with respect to the operation of composition. A subgroup $H$ is a group $G$ is called \index{normal subgroup}\defterm{normal} in $G$ if $xHx^{-1}=H$ for all $x\in G$.

\section{Dilations}

\begin{definition} An automorphism $A:X\to X$ of a liner $X$ is called a \index{dilation}\defterm{dilation} if for every line $L$ in $X$, the line $A[L]$ is parallel to $L$.
\end{definition}

\begin{proposition}\label{p:Dil(X)-subgroup} Let $X$ be a proaffine regular liner. The set\index[note]{$\Dil(X)$} $\Dil(X)$ of dilations of $X$ is a normal subgroup in the automorphism group $\Aut(X)$ of $X$.
\end{proposition} 

\begin{proof} Given any dilations $A,B\in\Dil(X)$ we first show that the automorphism $AB^{-1}:X\to X$ is a dilation. Given any line $L$ in $X$, we need to show that the line $AB^{-1}[L]$ is parallel to $L$. Since $B$ is a dilation of $X$, so is the automorphism $B^{-1}:X\to X$. Then  the line $\Lambda\defeq B^{-1}[L]$ is parallel to $L$. Since $A$ is a dilation of $X$, the line $A[\Lambda]$ is parallel to the line $\Lambda$. Applying Theorem~\ref{t:Proclus-lines}, we conclude that the line $AB^{-1}[L]=A[\Lambda]$ is parallel to the line $L$, witnessing that the automorphism $AB^{-1}$ is a dilation. Therefore, the set of dilation $\Dil(X)$ is a subgroup of the automorphism group $\Aut(X)$.
\smallskip

Next, we check that the subgroup $\Dil(X)$ is normal in the group $\Aut(X)$. Given a dilation $D:X\to X$ and an automorphism $A:X\to X$, we need to check that the automorphism $ADA^{-1}$ is a dilation of $X$. Let $L$ be any line in $X$ and $\Lambda\defeq A^{-1}[L]$ be its image under the automorphism $A^{-1}$. Since $D$ is a dilation, the line $\Lambda'\defeq D[\Lambda]$ is parallel to $\Lambda$. Since $A$ is an automorphism of the liner $X$, the lines $A[\Lambda']=ADA^{-1}[L]$ and $A[\Lambda]=L$ are parallel, witnessing that the automorphism $ADA^{-1}$ is a dilation of $X$.
\end{proof}

Let $F:X\to X$ be a self-map of a set $X$. A point $x\in X$ is called a \index{fixed point}\defterm{fixed point} of the map $F$ if $F(x)=x$. We denote by $$\Fix(F)\defeq\{x\in X:F(x)=x\}$$the set of fixed points of the map $F$.

\begin{lemma}\label{l:dilation-fixed2} If a dilation $A:X\to X$ of a liner $X$ has two distinct fixed points $a,b$, then $A(x)=x$ for every $x\in X\setminus\Aline ab$.
\end{lemma}

\begin{proof} Given any point $x\in X\setminus\Aline ab$, consider the lines $L\defeq\Aline ax$ and $\Lambda\defeq\Aline bx$. It follows from $x\notin\Aline ab$ that $L\cap\Lambda=\{x\}$. Taking into account that $A$ is a dilation of $X$, we conclude that $A[L]\parallel L$ and $A[\Lambda]\parallel \Lambda$. Taking into account that $a=A(a)\in L\cap A[L]$  and $b=A(b)\in \Lambda\cap A[\Lambda]$, we can apply Proposition~\ref{p:para+intersect=>coincide} and conclude that $A[L]=L$ and $A[\Lambda]=\Lambda$. Then $A(x)\in A[L]\cap A[\Lambda]=L\cap\Lambda=\{x\}$ and hence $A(x)=x$. 
\end{proof}

\begin{theorem}\label{t:dilation-atmost1} Every non-identity dilation $A:X\to X$ of a liner $X$ of rank $\|X\|\ge 3$ has at most one fixed point.
\end{theorem}

\begin{proof} Assume that $A$ has two distinct fixed points $a,b\in X$. Since $\|X\|\ge 3$, there exists a point $c\in X\setminus \Aline ab$. Lemma~\ref{l:dilation-fixed2} ensures that $c$ is a fixed point of $X$. We claim that $A(x)=x$ for every point $x\in X$. If $x\in X\setminus\Aline ab$, then $A(x)=x$, by Lemma~\ref{l:dilation-fixed2}. So, assume that $x\in \Aline ab$. If $x\in \{a,b\}$, then $A(x)=x$ because $a,b$ are fixed points of $A$. So assume that $x\in \Aline ab\setminus\{a,b\}$. Then $x\in X\setminus\Aline ac$ and $A(x)=x$, by Lemma~\ref{l:dilation-fixed2} because $a$ and $c$ are fixed points of the dilation $A$.
\end{proof}

\begin{proposition} If a completely regular liner $X$ admits a non-identity dilation, then $X$ is Playfair.
\end{proposition}

\begin{proof}  By Theorem~\ref{t:spread=projective1}, the completely regular liner $X$ is regular, para-Playfair and hence Proclus. Let $D:X\to X$ be a non-identity dilation of the liner $X$. To prove that $X$ is Playfair, it suffices to show that every line $L$ in $X$ is spreading. If $X=L$, then $L_\parallel=\{L\}=\{X\}$ is a spread of parallel lines in $X$ and the line $L$ is spreading. So, assume that $X\ne L$ and hence $\|X\|\ge 3$.

By Theorem~\ref{t:dilation-atmost1}, the non-identity dilation $D:X\to X$ has at most one fixed point. So, we can choose a point $a\in L$ such that $a'\defeq D(a)\ne a$.

\begin{claim}\label{cl:Lab-spreading} For every point $b\in X\setminus\Aline a{a'}$, the line $\Aline ab$ is spreading in $X$.
\end{claim}

\begin{proof} Consider the point $b'\defeq D(b)$ and observe that $b'=D(b)\ne D(a)=a'$. Since $D$ is a dilation, the line $\Aline ab$ is parallel to the line $\Aline{a'}{b'}$. Assuming that $\Aline ab=\Aline{a'}{b'}$, we conclude that $a'\in\Aline a{a'}\cap\Aline {a'}{b'}=\Aline a{a'}\cap\Aline ab=\{a\}$, which contradicts the choice of the point $a$. Therefore, the parallel lines $\Aline ab$ and $\Aline{a'}{b'}$ are disjoint. Since $X$ is para-Playfair, the line $\Aline ab$ is spreading, by Proposition~\ref{p:lines-in-para-Playfair}.
\end{proof} 

If $a'\notin L$, then choose any point $b\in L\setminus\{a\}=L\setminus\Aline a{a'}$ and applying Claim~\ref{cl:Lab-spreading} conclude that the line $L$ is spreading. So, assume that $a'\in L$, which implies that $L=\Aline a{a'}$. 
Since $\|X\|\ge 3$, there exists a point $b\in X\setminus\Aline a{a'}$. By Claim~\ref{cl:Lab-spreading}, the line $\Aline ab\subseteq \Pi$ is spreading. For the point $b'=D(b)\ne D(a)=a'$, the line $\Aline {a'}{b'}=D[\Aline ab]$ is parallel to the line $\Aline ab$. Since $D$ is a dilation, the line $D[\Aline a{a'}]$ is parallel to the line $\Aline a{a'}$. Since the parallel lines $\Aline a{a'}$ and $D[\Aline a{a'}]$ have a common point $a'$, they coincide. Then $b\notin \Aline a{a'}$ implies $b'=D(b)\notin D[\Aline a{a'}]=\Aline a{a'}$. By Claim~\ref{cl:Lab-spreading}, the line $\Aline a{b'}$ is spreading. Since $\Aline ab\cap\Aline a{a'}=\{a\}\ne\{a'\}\subseteq\Aline a{a'}\cap\Aline{a'}{b'}$, the parallel lines $\Aline ab$ and $\Aline {a'}{b'}$ are distinct and hence disjoint. Then $b'\notin\Aline ab$ and the spreading lines $\Aline ab$ and $\Aline a{b'}$ in the plane $\Pi=\overline{\{a,b,a',b'\}}$ are concurrent.
Since $X$ is completely regular, the line $L\subseteq\Pi$ is spreading, by Theorem~\ref{t:spread=projective1}. 
\end{proof}

\begin{exercise} Find an example of proaffine regular liner $X$, which is not Playfair but admits a non-identity dilation $D:X\to X$.
\end{exercise}

Theorem~\ref{t:dilation-atmost1} motivates the following classification of dilations.

\begin{definition} A non-identity automorphism $A:X\to X$ of a liner $X$ is called 
\begin{enumerate}
\item a \index{translation}\defterm{translation} of $X$ if $A$ is a dilation of $X$ having no fixed points;
\item a \index{homothety}\defterm{homothety} of $X$ if $A$ is a dilation of $X$ having a unique fixed point $o$ called the \index{centre}\defterm{centre} of the homothety.
\end{enumerate}
By definition, the identity automorphism of $X$ is both a translation and a homothety of $X$.
\end{definition}

\section{Translations}

This section is devoted to studying translations of liners. Given any automorphism $A:X\to X$ of a liner $X$, consider the family of flats
$$\lines{A}\defeq\{\Aline xy:xy\in A\}.$$
It is easy to see that  an automorphism $A:X\to X$ is the identity automorphism of $X$ if and only if $\lines{A}$ is the family of all singletons in $X$.

\begin{proposition}\label{p:Trans-spread} For every nonidentity translation $T:X\to X$ of a liner $X$, the family $\lines{T}$ is a spread of lines in $X$ such that $T[L]=L$ for every $L\in\lines{T}$. If the liner $X$ is $3$-ranked, then $\lines{T}$ is a spread of parallel lines in $X$.
\end{proposition}

\begin{proof} Since the translation $T$ has no fixed points, for every $xy\in T$, the flat $\Aline x y$ is a line in $X$ containing the point $x$. Next, we prove that $T[L]=L$ for every $L\in\lines{T}$. Given a line $L\in\lines{T}$, find a pair $xy\in T$ such that $L=\Aline xy$, and consider the line $T[L]=\Aline{Tx}{Ty}$. Since $Tx=y\in L\cap T[L]\ne\varnothing$, the parallel lines $L$ and $T[L]$ coincide. 

If two lines $L,\Lambda\in\lines T$ have a common point $x$, then $L=T[L]=\Aline x{Tx}=T[\Lambda]=\Lambda$, witnessing that $\lines{T}$ is a spread of lines in $X$.

Now assume that the liner $X$ is $3$-ranked. Given two pairs $xy,uv\in T$, we need to check that the lines $\Aline xy$ and $\Aline uv$ are parallel. If $x\in \Aline uv$, then $\Aline xy=\Aline uv$ because $\lines{T}$ is spread of lines. So, assume that $x\notin\Aline uv$. In this case $\Aline xu$ and $\Aline yv$ are parallel lines (because $T$ is a dilation) and $\overline{\{x,u,y,v\}}$ is a plane in $X$. By Corollary~\ref{c:parallel-lines<=>}, the disjoint lines  $\Aline xy$ and $\Aline uv$ in the plane $\overline{\{x,u,y,v\}}$ are parallel.
\end{proof}

\begin{proposition} Let $T$ be a non-identity translation of a liner $X$. For every $n\in\IN$ with $T^n\ne 1_X$, the spreads of lines $\lines{T^n}$ and $\lines{T}$ coincide.
\end{proposition}

\begin{proof} The  implication $(T^n\ne1_X)\;\Ra\;(\lines{T^n}=\lines{T})$ will be proved by induction on $n$. For $n=1$, this implication is trivially true. Assume that for some $n\in\IN$ we know that $\lines{T^k}=\lines{T}$ for all $k<n$ with $T^k\ne 1_X$. We have to prove that $\lines{T^n}=\lines{T}$ if $T^n\ne 1_X$. If $T^{n-1}=1_X$, then $T^n=T$ and $\lines{T^n}=\lines T$. So, assume that $T^{n-1}\ne1_X$. The induction hypothesis ensures that $\lines{T^{n-1}}=\lines T$. Given any line $L\in\lines{T^n}$, find a pair $oy\in T^n$ such that $L=\Aline oy$. Consider the point $x=T^{n-1}(o)$ and observe that $y=T^n(o)=T(T^{n-1}(o))=T(x)$. It follows from $\lines{T^{n-1}}=\lines T$ that $\Aline ox=\Aline xy$ and hence $L=\Aline oy=\Aline xy\in \lines T$. Therefore, $\lines{T^n}\subseteq\lines T$. To prove that $\lines  T\subseteq\lines{T^n}$, take any line $\Lambda\in\lines T$ and find a pair $uv\in T$ such that $\Aline uv=\Lambda$. Consider the point $w=T^n(u)$ and the line $\Aline uw\in\lines{T^n}$. Since $\lines{T^n}$ and $\lines{T}$ are two spreads of lines with $\lines{T^n}\subseteq\lines{T}$, the lines $\Aline uw$ and $\Aline uv$ coincide and hence $\Lambda=\Aline uv=\Aline uw\in\lines{T^n}$.   
\end{proof}

\begin{exercise} Prove that every para-Playfair regular liner admitting a non-identity translation is  Playfair.
\end{exercise}

\begin{proposition}\label{p:Ax=Bx=>A=B} Let $X$ be a Proclus liner of rank $\|X\|\ge 3$. Two translations $A,B$ of $X$ are equal if and only if $A(x)=B(x)$ for some point $x\in X$.
\end{proposition}

\begin{proof} The ``only if'' part of this characterization is trivial. The ``if'' part will be deduced from the following claim.

\begin{claim}\label{cl:ABx=>A=B} If $x'=A(x)=B(x)$ for some points $x,x'\in X$, then $A(y)=B(y)$ for every $y\in X\setminus \Aline x {x'}$.
\end{claim}

\begin{proof} If $x'=A(x)=B(x)=x$, then the translations $A,B$ are identity maps of $X$ and hence $A(y)=y=B(y)$ for every $y\in X$. So, assume that $x'\ne x$ and consider the line $\Aline x{x'}$. Given any point $y\in X\setminus\Aline x{x'}$, we need to prove that $A(y)=B(y)$. Consider the points $y'\defeq A(y)$ and $y''\defeq B(y)$. Since $A$ is a dilation, $\Aline xy\parallel \Aline {x'}{y'}$ and hence $\Aline{x'}{y'}\subseteq \overline{\Aline xy\cup\{x'\}}=\overline{\{x,y,x'\}}$. By Theorem~\ref{t:Proclus<=>}, the Proclus liner $X$ is $3$-proregular and by Proposition~\ref{p:k-regular<=>2ex}, the $3$-proregular liner $X$ is $3$-ranked. By Proposition~\ref{p:Trans-spread}, the lines $\Aline x{x'}$ and $\Aline y{y'}$ in the plane $\overline{\{x,y,x'\}}$ are parallel. By analogy we can prove that $\Aline x{x'}$ and $\Aline y{y''}$ are two disjoint parallel lines in the plane $\overline{\{x,y,x'\}}$. Since the liner $X$ is Proclus,  the lines $\Aline y{y'}$ and $\Aline y{y''}$ coincide. Taking into account that $A,B$ are dilations with $A(x)=x'=B(x)$,   $y'=A(y)$ and $y''=B(y)$, we conclude that $\Aline {x'}{y'}\parallel\Aline xy\parallel\Aline{x'}{y''}$  and hence $\Aline{x'}{y'}=\Aline {x'}{y''}$, by the Proclus property of the liner $X$. Then $$\{A(y)\}=\{y'\}=\Aline {y}{y'}\cap\Aline {x'}{y'}=\Aline y{y''}\cap\Aline{x'}{y''}=\{y''\}=\{B(y)\}.$$
\end{proof}

Now we are ready two prove that two translations $A,B$ of $X$ coincide if $A(x)=B(x)$ for some point $x\in X$. Consider the point $x'\defeq A(x)=B(x)$. Given any point $y\in X$, we have to prove that $A(y)=B(y)$. If $y\notin \Aline x{x'}$, then the equality $A(y)=B(y)$ follows from Claim~\ref{cl:ABx=>A=B}. So, assume that $y\in \Aline x{x'}$. If $y=x$, then $A(y)=A(x)=B(x)=B(y)$ by the choice of $x$. So, assume that $y\ne x$. Since $y\in\Aline x{x'}\setminus\{x\}$, the points $x'=A(x)=B(x)$ is not equal to $x$, which implies that $A,B$ are non-identity translations of $X$.  

Since $\|X\|\ge 3$, there exists a point $z\in X\setminus \Aline x{x'}=\Aline xy$. By Claim~\ref{cl:ABx=>A=B}, $A(z)=B(z)$. Consider the point $z'\defeq A(z)=B(z)$. Since $A,B$ are non-identity translations of $X$, the point $z'\defeq A(z)=B(z)$ is not equal to the point $z$. By Proposition~\ref{p:Trans-spread}, the lines $\Aline x{x'}$ and $\Aline z{z'}$ are disjoint and hence $y\notin \Aline z{z'}$. Since $z'=A(z)=B(z)$, we can apply Claim~\ref{cl:ABx=>A=B} and conclude that $A(y)=B(y)$. 
\end{proof}

\begin{proposition}\label{p:Trans(X)isnormal} Let $X$ be a proaffine regular liner of rank $\|X\|\ge 3$. The set\index[note]{$\Trans(X)$} $\Trans(X)$ of translations is a normal subgroup of the automorphism group $\Aut(X)$ of $X$.
\end{proposition}

\begin{proof} To prove that $\Trans(X)$ is a subgroup of the automorphism group $\Aut(X)$, we have to show that for any translations $A,B\in\Trans(X)$, the automorphism $A^{-1}B$ is a translation of $X$.  
Assuming that $A^{-1}B$ is not a translation, we can apply Theorem~\ref{t:dilation-atmost1} and conclude that the dilation $A^{-1}B$ has a unique fixed point $o$. Consider the point $b\defeq B(o)$ and observe that $o=A^{-1}B(o)=A^{-1}(b)$ implies $A(o)=b=B(o)$. Applying Proposition~\ref{p:Ax=Bx=>A=B}, we conclude that $A=B$ and hence $A^{-1}B$ is the identity map of $X$, which contradicts our assumption. This contradiction shows that $A^{-1}B$ is a translation and $\Trans(X)$ is a subgroup of the group $\Aut(X)$.
\smallskip

Next, we show that the subgroup $\Trans(X)$ is normal in the group $\Aut(X)$.
Given any translation $T\in\Trans(X)$ and automorphim $A\in\Aut(X)$, we have to prove that the automorphism $B\defeq A^{-1}TA$ is a translation of $X$. By Proposition~\ref{p:Dil(X)-subgroup}, the automorphism $B=A^{-1}TA$ is a dilation of $X$. Assuming that $B$ is not a translation, we conclude that $B$ has a unique fixed point $o$. Then the point $a\defeq A(o)$ is a fixed point of the translation $T$:
$$T(a)=TA(o)=AA^{-1}TA(o)=A(o)=a,$$
which implies that $T$ is the identity automorphism $X$ and so $T$ is the automorphism, which contradicts our assumption. This contradiction shows that the automorphism $A^{-1}TA$ is a translation of $X$ and the group $\Trans(X)$ is normal in $\Aut(X)$.
\end{proof}

%
A group $G$ is called \index{elementary group}\index{group!elementary}\defterm{elementary} if there exists a prime number $p$ such that $g^p=1$ for every $g\in G$.\index[person]{Baer}

\begin{theorem}[Baer\footnote{{\bf  Reinhold Baer} (1902 -- 1979) was a German mathematician, known for his work in algebra. He introduced injective modules in 1940. He is the eponym of Baer rings and Baer groups.}]\label{t:Trans-commutative}  Let $X$ be a proaffine regular liner. If $\|\{T(o):T\in\Trans(X)\}\|\ne 2$ for some point $o\in X$ (and the group $\Trans(X)$ contains a non-identity element of finite order), then the group $\Trans(X)$ is commutative (and elementary).
\end{theorem}

\begin{proof} Assume that $\|\{T(o):T\in\Trans(X)\}\|\ne 2$ for some point $o\in X$. If $\|\{T(o):T\in\Trans(X)\}\|=1$, then $\{T(o):T\in\Trans(X)\}=\{o\}$ and  $\Trans(X)$ is a trivial group. So, assume that $\|\{T(o):T\in\Trans(X)\}\|\ne 1$ and hence $\|X\|\ge\|\{T(o):T\in\Trans(X)\}\|\ge 3$.

\begin{claim}\label{cl:AB=BA} Let $A,B$ be two nonidentity translations of the proaffine liner $X$. If  $\lines A\ne\lines B$, then $AB=BA$.
\end{claim}

\begin{proof} Fix any point $o\in X$ and consider the points $a\defeq A(o)$ and $b\defeq B(o)$. By Proposition~\ref{p:Trans-spread}, $\lines A,\lines B$ are two spreads of parallel lines in $X$. By our assumption, those two spread are distinct, which implies that $b\notin\Aline oa$. 

Consider the points $b_a\defeq B(a)=BA(o)$ and $a_b\defeq A(b)=AB(o)$. 
 Since $A,B$ are dilations, $\Aline oa\parallel \Aline b{b_a}$ and $\Aline ob\parallel \Aline a{a_b}$. Then $\Aline oa,\Aline ob,\Aline b{b_a},\Aline a{a_b}$ are lines in the plane $\overline{\{o,a,b\}}$. Since $\lines A$ and $\lines B$ are spreads of parallel lines, $\Aline oa\parallel \Aline b{a_b}$ and $\Aline ob\parallel \Aline a{b_a}$. By Theorem~\ref{t:Proclus<=>}, the proaffine regular space $X$ is Proclus. Then $\Aline b{b_a}\parallel \Aline oa\parallel \Aline b{a_b}$ and $\Aline a{a_b}\parallel \Aline ob\parallel \Aline a{b_a}$ imply $\Aline b{b_a}=\Aline b{a_b}$ and $\Aline a{a_b}=\Aline  a{b_a}$.

It follows from $b\notin \Aline oa$ that $\{a_b\}=\Aline {b}{a_b}\cap\Aline {a}{a_b}=\Aline b{b_a}\cap\Aline a{b_a}=\{b_a\}$ and hence $AB(o)=A(b)=a_b=b_a=B(a)=BA(o)$ and $AB=BA$, by Proposition~\ref{p:Ax=Bx=>A=B}.
\end{proof}

\begin{claim} The group $\Trans(X)$ is commutative.
\end{claim}

\begin{proof} Given any translations $A,B\in \Trans(X)$, we should prove that $AB=BA$. The equality $AB=BA$ trivially holds if $A$ or $B$ is the identity map. So, assume that $A$ and $B$ are non-identity translations. Since $\Trans(X)$ is a subgroup of the group $\Aut(X)$, the compositions $AB$ and $BA$ are translations of $X$. Since $\|\{T(o):T\in\Trans(X)\}\|\ge 3$,  there exist translations $C,D\in\Trans(X)$ such that $\|\{o,C(o),D(o)\}\|=3$.  If $\|\{o,A(o),B(o)\}\|=3$, then $AB=BA$, by Claim~\ref{cl:AB=BA}. So, assume that $\|\{o,A(o),B(o)\}\|<3$, which implies that $L\defeq\overline{\{o,A(o),B(o)\}}$ is a line in $X$. Since $\|\{o,C(o),D(o)\}\|=3$, there exists a translation $T\in\{C,D\}$ such that $T(o)\notin L$. Then $AT$ is a translation such that $AT(o)\notin L=\overline{\{o,B(o)\}}=\overline{\{o,A(o)\}}=\overline{\{0,AB(o)\}}$ and hence $B(AT)=(AT)B$, $AT=TA$ and $T(AB)=(AB)T$, by Claim~\ref{cl:AB=BA}. Then $$BA=BA(TT^{-1})=B(AT)T^{-1}=(AT)BT^{-1}=(TA)BT^{-1}=T(AB)T^{-1}=(AB)TT^{-1}=AB.$$
 \end{proof}



\begin{claim}\label{cl:An=1=>Bn=1} Let $A,B\in\Trans(X)$ be two non-identity translations with $\lines A\ne\lines B$. Then $$\forall n\in\IN\;\;(A^n=1_X\;\Leftrightarrow\;B^n=1_X).$$
\end{claim}

\begin{proof} Given any $n\in\IN$, it suffices to prove that $A^n=1_X$ implies $B^n=1_X$. To derive a contradiction, assume that $A^n=1_X\ne B^n$. By Claim~\ref{cl:AB=BA}, $AB=BA$ and hence $(BA^{-1})^n=B^nA^{-n}=B^n\ne 1_X$. Consider the points $a\defeq A(o)$, $b\defeq B(o)=BA^{-1}(a)$. It follows from $\lines A\ne\lines B$ that $b\notin\Aline oa$ and hence the lines $\Aline ob$ and $\Aline ab$ are concurrent. By Propositions~\ref{p:Trans-spread}, the families $\lines A,\lines B,\lines{BA^{-1}}$ are spreads of lines in $X$ such that $\Aline ob\in \lines B=\lines{B^n}=\lines{(BA^{-1})^n}=\lines {BA^{-1}}$. Then the lines $\Aline ob,\Aline ab\in\lines{BA^{-1}}$ are parallel, which is not true. This is a contradiction showing that $A^n=1_X$ implies $B^n=1_X$.
\end{proof}

\begin{claim} If the group $\Trans(X)$ contains a non-identity element of finite order, then $X$ is elementary.
\end{claim}

\begin{proof} Since $\Trans(X)$ contains a non-identity element of finite order, we can consider the smallest number $p$ for which there exists a non-identity translation $T\in\Trans(X)$ with $T^p=1_X$. The minimality of $p$ ensures that $p$ is a prime number. Applying Claim~\ref{cl:An=1=>Bn=1}, we conclude that $A^p=1_X$ for all $A\in\Trans(X)$.
\end{proof}
\end{proof}

A group $G$ is called \index{group!divisible}\index{divisible group}\defterm{divisible} if the equality $x^n=g$ has a solution in the group $G$ for every $g\in G$ and $n\in\IN$.

\begin{problem} Let $X$ be a proaffine regular liner such that $\|\{T(o):T\in\Trans(X)\}\|\ge 3$ for some point $o\in X$, and every non-identity translation of $X$ has infinite order in the group $\Trans(X)$. Is the group $\Trans(X)$ divisible?
\end{problem}

The following theorem shows that the condition $\|\{T(o):T\in\Trans(X)\}\|\ne2$ is essential in Theorem~\ref{t:Trans-commutative}. The idea of the proof of this theorem was suggested by \index[person]{Krauz}Yoav Krauz\footnote{{\bf Yoav Krauz}, a mathematician from Israel, the absolute winner of the International Mathematical Olympiads for undergraduate students in 2014, see {\tt https://www.imc-math.org.uk/imc2014/imc2014-scores.html}.} in his answer to the problem {\tt https://mathoverflow.net/q/467100/61536}, posed by the author at {\tt MathOverflow}. 
\begin{theorem}\label{t:Krauz} For every infinite group $G$, there exists a Playfair plane $X$ of cardinality $|X|=|G|$ such that $\Dil(X)=\Trans(X)$ and the group $\Dil(X)=\Trans(X)$  is isomorphic to the group $G$.
\end{theorem}

\begin{proof} Let  $\kappa\defeq |G|$ be the cardinality of the infinite group $G$.  Every function $F:\kappa\to G$ will be identified with its graph $\{(x,y)\in \kappa\times G:y=F(x)\}$ in $\kappa\times G$. Given any function $F:\kappa\to G$ and an element $g\in G$ we denote by $gF$ the function $\{(x,gy):(x,y)\in F\}\in G^\kappa$. A family of functions $\F\subseteq G^\kappa$ is called \index{family!transversal}\index{transversal family}\defterm{$G$-transversal} if for every functions $F,F'\in\F$ and every $g\in G$, either $F=gF'$ or $F\cap gF'$ is a singleton.

\begin{lemma}\label{l:Krauz1} Let $\F\subseteq G^\kappa$ be a nonempty $G$-transversal family of bijections of cardinality $|\F|<\kappa$, and let $\varphi\subseteq \kappa\times G$ be an injective function of cardinality $|\varphi|<\kappa$ such that $|\varphi\cap gF|\le 1$ for every $F\in \F$ and every $g\in G$. Then there exists a bijective function $\Phi:\kappa\to G$ such that $\varphi\subseteq \Phi$ and the family $\F\cup\{\Phi\}$ is $G$-transversal.
\end{lemma}

\begin{proof} Write the set  of ordinals $[|\varphi|,\kappa)\defeq\{\alpha\in\kappa:|\varphi|\le\alpha\}$ as the union $\Omega\cup\Omega'\cup\Omega''$ of three pairwise disjoint sets $\Omega,\Omega',\Omega''$ of cardinality $|\Omega|=|\Omega'|=|\Omega''|=\kappa$. Since the set $G\times\F$ has cardinality $\kappa$, there exists a function $\xi:\kappa\to G\times \F$ such that $\xi[\Omega]=G\times \F$. For every ordinal $\alpha\in\kappa$, consider the pair $(g_\alpha,F_\alpha)\defeq\xi(\alpha)$. 

Fix a strict well-order $\prec$ on the group $G$ such that for every $g\in G$, the set ${\downarrow}g\defeq\{x\in G:x\prec g\}$ has cardinality $<\kappa$. For a nonempty subset $A\subseteq G$ we denote by $\min A$ the smallest element of the set $A$ with respect to the well-order $\prec$ on $G$.

\begin{claim}\label{cl:Krauz1} There exist transifinite sequences $(x_\alpha)_{\alpha\in\kappa}\in \kappa^\kappa$ and $(y_\alpha)_{\alpha\in\kappa}\in G^\kappa$ such that $\{(x_\alpha,y_\alpha):\alpha\in|\varphi|\}=\varphi$ and for every $\alpha\in[|\varphi|,\kappa)$ the following conditions are satisfied:
\begin{enumerate}
\item[$(1_\alpha)$] if $\alpha\in\Omega'$, then $x_\alpha\defeq\min(\kappa\setminus\{x_\beta:\beta<\alpha\})$ and\newline $y_\alpha\notin \{y_\beta:\beta<\alpha\}\cup\{y_\beta F(x_\beta)^{-1}F(x_\alpha):\beta<\alpha,\;F\in\F\}$;
\item[$(2_\alpha)$] if $\alpha\in\Omega''$, then $y_\alpha\defeq\min (G\setminus\{y_\beta:\beta<\alpha\})$ and\newline $x_\alpha\notin \{x_\beta:\beta<\alpha\}\cup\bigcup_{\beta\in\alpha}\bigcup_{F\in\F}\{x\in\kappa:y_\alpha= y_\beta F(x_\beta)^{-1} F(x)\}$;
\item[$(3_\alpha)$] if $\alpha\in\Omega$ and the set $C_\alpha\defeq \{\beta<\alpha:y_\beta=g_\alpha F_\alpha(x_\beta)\}$ is not empty,\newline then $(x_\alpha,y_\alpha)\defeq(x_\gamma,y_\gamma)$ where $\gamma=\min C_\alpha$; 
\item[$(4_\alpha)$] if $\alpha\in\Omega$ and the set $C_\alpha$ is empty, then $x_\alpha\notin\{x_\beta:\beta<\alpha\}$ and\newline $y_\alpha=g_\alpha F_\alpha(x_\alpha)\notin\{y_\beta:\beta<\kappa\}\cup\{y_\beta F(x_\beta)^{-1}F(x_\alpha):\beta<\alpha,\;F\in\F\}$.
\end{enumerate}
\end{claim}

\begin{proof} Write the injective function $\varphi\subseteq \kappa\times G$ as $B=\{(x_\alpha,y_\alpha):\alpha\in|\varphi|\}$ for some transfinite sequences $(x_\alpha)_{\alpha\in|\varphi|}$ and $(y_\alpha)_{\alpha\in|\varphi|}$. Assume that for some ordinal $\alpha\in[|\varphi|,\kappa)$ we have constructed sequences $(x_\beta)_{\beta<\alpha}$ and $(y_\beta)_{\beta<\alpha}$ satisfying the conditions $(1_\beta)$--$(4_\beta)$ for every ordinal $\beta\in[|\varphi|,\alpha)$.

If $\alpha\in\Omega'$, then consider the point $x_\alpha\defeq\min(\kappa\setminus\{x_\beta:\beta\in\alpha\})$ and choose any point $y_\alpha\in G$ that does not belong to the set
$$\{y_\beta:\beta\in\alpha\}\cup\{y_\beta\cdot F(x_\beta)^{-1}F(x_\alpha):\beta<\alpha,\;F\in\F\},$$which has cardinality $|\alpha|+|\alpha\times\F|<\kappa$, so the point $y_\alpha$ indeed exists.

If $\alpha\in\Omega''$, then consider the point $y_\alpha\defeq\min (G\setminus\{y_\beta:\beta\in\alpha\})$ and choose any point $x_\alpha\in\kappa$ that does not belong to the set $\{x_\beta:\beta\in\alpha\}\cup\bigcup_{\beta\in\alpha}\bigcup_{F\in\F}\{x\in\kappa:y_\alpha= y_\beta\cdot F(x_\beta)^{-1}\cdot F(x)\}$, which has cardinality $|\alpha|+|\alpha\times\F|<\kappa$, so the point $x_\alpha$ indeed exists.

If $\alpha\in\Omega$ and the $C_\alpha\defeq \{\beta<\alpha:y_\beta=g_\alpha F_\alpha(x_\beta)\}$ is not empty, then put $(x_\alpha,y_\alpha)\defeq(x_\gamma,y_\gamma)$ where $\gamma=\min C_\alpha$.

Finally, assume that $\alpha\in\Omega$ and $C_\alpha=\varnothing$. For every $F\in\F$ and every ordinal $\beta<\alpha$, consider the set
$$X_{F,\beta}\defeq\{x\in \kappa:g_\alpha F_\alpha(x)=y_\beta  F(x_\beta)^{-1} F(x)\}.$$ The $G$-transversality of the family $\F$ ensures that for every $F\in\F\setminus\{F_\alpha\}$ the set $X_{F,\beta}$ is a singleton. On the other hand, $C_\alpha=\varnothing$ implies that for every $\beta<\alpha$, the set $X_{F_\alpha,\beta}$ is empty. Then the set $$X_\alpha\defeq\bigcup_{F\in\F}\bigcup_{\beta<\alpha}X_{F,\alpha}$$has cardinality $|X_\alpha|\le|\F|\cdot|\alpha|<\kappa$. Since $F_\alpha:\kappa\to G$ is a bijection, the set $$X_\alpha'\defeq\bigcup_{\beta<\alpha}\{x\in X:y_\alpha=g_\alpha F_\alpha(x)\}$$ has cardinality $\le|\alpha|<\kappa$. Choose any point $$x_\alpha\in \kappa\setminus(X_\alpha\cup X_\alpha'\cup\{x_\beta:\beta<\alpha\})$$and observe that the points $x_\alpha$ and $y_\alpha\defeq g_\alpha F_\alpha(x_\alpha)$ satisfy the inductive condition $(4_\alpha)$. This completes the inductive step of the construction. 
\end{proof}

Let $(x_\alpha)_{\alpha\in\kappa}$ and $(y_\alpha)_{\alpha\in\kappa}$ be the transfinite sequences provided by Claim~\ref{cl:Krauz1}. The inductive conditions $(1_\alpha)_{\alpha\in\Omega'}$ and $(2_\alpha)_{\alpha\in\Omega''}$ imply that $\{x_\alpha:\alpha\in\kappa\}=\kappa$ and $\{y_\alpha:\alpha\in \kappa\}=G$.

We claim that the relation 
$$\Phi\defeq\{(x_\alpha,y_\alpha):\alpha\in\kappa\}$$ is an injective function  such that $\varphi\subseteq \Phi$ and the family $\F\cup\{\Phi\}$ is $G$-transversal.

Assuming that $\Phi$ is not a function, we can find two ordinals $\beta<\alpha$ in $\kappa$ such that $x_\beta=x_\alpha$ and $y_\beta\ne y_\alpha$. We can assume that $\alpha$ is the smallest ordinal for which there exists an ordinal $\beta<\alpha$ such that $x_\beta=x_\alpha$ and $y_\beta\ne y_\alpha$. Since $\varphi=\{(x_\alpha,y_\gamma):\gamma\in|\varphi|\}$ is a function, the ordinal $\alpha$ belongs to the set $\Omega\cup\Omega'\cup\Omega''$. Since $x_\alpha=x_\beta$, the inductive conditions $(1_\alpha)$, $(2_\alpha)$ and $(4_\alpha)$ of Claim~\ref{cl:Krauz1} imply that $\alpha\in\Omega$ and the set $C_\alpha=\{\gamma\in\alpha:y_\gamma=g_\alpha\cdot F_\alpha(x_\alpha)\}$ is not empty. Then $(x_\alpha,y_\alpha)=(x_\gamma,y_\gamma)$ where $\gamma=\min C_\alpha<\alpha$. Since $x_\gamma=x_\alpha=x_\beta$ and $y_\gamma=y_\alpha\ne y_\beta$, we obtain a contradiction with the minimality of the ordinal $\alpha$. This contradiction shows that $\Phi$ is a function. 

Assuming that the function $\Phi$ is not injective,  we can find two ordinals $\beta<\alpha$ in $\kappa$ such that $x_\beta\ne x_\alpha$ and $y_\beta=y_\alpha$. We can assume that $\alpha$ is the smallest ordinal for which there exists an ordinal $\beta<\alpha$ such that $x_\beta\ne x_\alpha$ and $y_\beta=y_\alpha$. Since $\varphi=\{(x_\alpha,y_\gamma):\gamma\in|\varphi|\}$ is an injective function, the ordinal $\alpha$ belongs to the set $\Omega\cup\Omega'\cup\Omega''$. Since $y_\beta=y_\alpha$,   the inductive conditions $(1_\alpha)$, $(2_\alpha)$ and $(4_\alpha)$ of Claim~\ref{cl:Krauz1} imply that $\alpha\in\Omega$ and the set $C_\alpha=\{\gamma\in\alpha:y_\gamma=g_\alpha\cdot F_\alpha(x_\alpha)\}$ is not empty. Then $(x_\alpha,y_\alpha)=(x_\gamma,y_\gamma)$ where $\gamma=\min C_\alpha<\alpha$. Since $x_\gamma=x_\alpha\ne x_\beta$ and $y_\gamma=y_\alpha=y_\beta$, we obtain a contradiction with the minimality of the ordinal $\alpha$. This contradiction shows that the function $\Phi$ is injective.

Since $\{x_\alpha:\alpha\in\kappa\}=\kappa$ and $\{y_\alpha:\alpha\in\kappa\}=G$, the injective function $\Phi:\kappa\to G$ is bijective.
The choice of the transfinite sequences $(x_\alpha)_{\alpha\in|\varphi|}$ and $(y_\alpha)_{\alpha\in|\varphi|}$ ensures that 
$\varphi=\{(x_\alpha,y_\alpha):\alpha\in |\varphi|\}\subseteq \Phi$.

It remains to check that the family $\F\cup\{\Phi\}$ is $G$-transversal. Since $\F$ is $G$-transversal, it suffices to show that for every $g\in G$ and $F\in\F$ the intersection $\Phi\cap gF$ is a singleton. 

First we show that the set $\Phi\cap gF$ is not empty.
Since $(g,F)\in G\times\F=\xi[\Omega]$, there exists an ordinal $\alpha\in\Omega$ such that $(g,F)=\xi(\alpha)=(g_\alpha,F_\alpha)$. Consider th set $$C_\alpha\defeq\{\beta\in \alpha:y_\beta=g_\alpha F_\alpha(x_\beta)\}=\{\beta\in\alpha:y_\beta=g F(x_\beta)\}.$$ If $C_\alpha\ne\varnothing$, then for the ordinal $\gamma\defeq \min C_\alpha$, the inductive condition $(3_\alpha)$ ensures that
$$(x_\alpha,y_\alpha)=(x_\gamma,y_\gamma)=(x_\gamma,g F(x_\gamma))\in \Phi\cap gF\ne\varnothing.$$ If $C_\alpha=\varnothing$, then the inductive condition $(4_\alpha)$ ensures that $y_\alpha=g_\alpha F_\alpha(x_\alpha)=g F(x_\alpha)$ and hence
$(x_\alpha,y_\alpha)\in\Phi\cap g F\ne\varnothing$.

Assuming that the  intersection $\Phi\cap gF$ is not a singleton, we can find two ordinals $\beta<\alpha$ in $\kappa$ such that $x_\beta\ne x_\alpha$ and $\{(x_\beta,y_\beta),(x_\alpha,y_\alpha)\}\subseteq \Phi\cap gF$.
We can assume that $\alpha$ is the smallest ordinal with this property. 
It follows that $y_\beta=g F(x_\beta)$ and
$$y_\alpha=g F(x_\alpha)=y_\beta F(x_\beta)^{-1} F(x_\alpha).$$
If $\alpha<|\varphi|$, then $(x_\beta,y_\beta)$ and $(x_\alpha,y_\alpha)$ are two distinct elements of the set $\varphi\cap g F$, which has cardinality $|\varphi\cap g F|\le 1$, by our asumption. This is a contradiction showing that $\alpha\in [|\varphi|,\kappa)=\Omega\cup\Omega'\cup\Omega'$. Since $y_\alpha=y_\beta F(x_\beta)^{-1} F(x_\alpha)$, the inductive conditions $(1_\alpha)$--$(4_\alpha)$ imply that $\alpha\in \Omega$ and the set $C_\alpha=\{\gamma\in\alpha:y_\alpha=g_\alpha F_\alpha(x_\gamma)\}$ is not empty. In this case, for the ordinal $\gamma=\min C_\alpha$ we have $(x_\gamma,y_\gamma)=(x_\alpha,y_\alpha)\ne(x_\beta,y_\beta)$ and hence
$$\{(x_\gamma,y_\gamma),(x_\beta,y_\beta)\}\subseteq \Phi\cap g F,$$ which contradicts the minimality of the ordinal $\alpha$. This contradiction show that the family $\F\cup\{\Phi\}$ is $G$-transversal.
\end{proof}

Let $\F_2$ be the set of injective functions $\varphi\subseteq \kappa\times G$ such that $|\varphi|=2$. Let $\F_3$ be the set of all sequences $\big((x_1,y_1),(x_2,y_2),(x_3,y_3)\big)\in (\kappa\times G)^3$ such that $|\{0,x_1,x_2,x_3\}|=4$. Let $\F_4$ be the set of sequences $(x_0,y_0,y_1,y_2,y_3)\in \kappa\times G^4$ such that $y_0\notin\{y_1,y_2,y_3\}$ and $y_1\ne y_2\ne y_3$.


Write the set $\kappa$ as the union $\kappa=\Omega_2\cup\Omega_3\cup\Omega_4\cup(\Omega_4+1)$ of four pairwise disjoint sets $\Omega_2,\Omega_3,\Omega_4$ and $\Omega_4+1\defeq\{\alpha+1:\alpha\in\Omega_4\}$ of cardinality $\kappa$.

Since $|\F_2\cup\F_3\cup \F_4|=\kappa$, there exists a transfinite sequence  $(F_{\alpha})_{\alpha\in\kappa}$ such that $\F_i=\{F_\alpha\}_{\alpha\in \Omega_i}$ for $i\in\{2,3,4\}$ and $F_{\alpha+1}=F_\alpha$ for all $\alpha\in\Omega_4$. Moreover, we can assume that for every $F\in\F_2\cup\F_3\cup\F_4$, the set $\{\alpha\in\kappa:F_\alpha=F\}$ has cardinality $\kappa$.

\begin{lemma}\label{l:Krauz2} There exists a transfinite sequence $(\Phi_\alpha)_{\alpha\in\kappa}$ of bijective functions $\Phi_\alpha:\kappa\to G$ such that for every $\alpha\in\kappa$ the following conditions are satisfied:
\begin{enumerate}
\item[$(1_\alpha)$] the family $\{\Phi_\beta\}_{\beta\le\alpha}$ is $G$-transversal;
\item[$(2_\alpha)$] if $\alpha\in\Omega_2$, then $F_\alpha\subseteq g\Phi_\alpha$ for some $g\in G$;
\item[$(3_\alpha)$] if $\alpha\in\Omega_3$, then for the sequence $\big((x_1,y_1),(x_2,y_2),(x_3,y_3)\big)\defeq F_\alpha$ we have $\Phi_\alpha(0)=1_G$ and 
$\Phi_\alpha(x_2)\Phi_\alpha(x_1)\ne \Phi_\alpha(x_3)\ne\Phi_\alpha(x_2)y_2^{-1}\Phi_\alpha(x_1)y_1^{-1}y_3$;
\item[$(4_\alpha)$] if $\alpha\in\Omega_4$, then for the sequence $(x_0,y_0,y_1,y_2,y_3)\defeq F_\alpha$,  there exist distinct points $x_1,x_2\in \kappa$ such that $\{(x_0,y_0),(x_1,y_1),(x_2,y_2)\}\subseteq \Phi_\alpha$, $\{(x_0,y_0),(x_1,y_2)\}\in \Phi_{\alpha+1}$ and $(x_2,y_3)\notin \Phi_{\alpha+1}$.
\end{enumerate}
\end{lemma}

\begin{proof} Assume that for some ordinal $\alpha\in\kappa$ a $G$-transversal family $\{\Phi_\beta\}_{\beta<\alpha}$ has been constructed. Consider the element $F_\alpha$ of the transfinite sequence $(F_\beta)_{\beta\in\kappa}$. If $\alpha\in\Omega_2$ and $F_\alpha\subseteq g\Phi_\gamma$ for some $g\in G$ and $\gamma<\alpha$, then we put $\Phi_\alpha\defeq g\Phi_\gamma$ and conclude that the family $\{\Phi_\beta\}_{\beta\le\alpha}=\{\Phi_\beta\}_{\beta<\alpha}$ is $G$-transversal.  If $\alpha\in\Omega_2$ and $F_\alpha\not\subseteq g\Phi_\gamma$ for all $g\in G$ and $\gamma<\alpha$, then $|F_\alpha\cap g\Phi_\gamma|<|F_\alpha|=2$ and we can apply Lemma~\ref{l:Krauz1} to find a bijective function $\Phi_\alpha:\kappa\to G$ such that $F_\alpha\subseteq \Phi_\alpha$ and the family $\{\Phi_\beta\}_{\beta\le\alpha}$ is $G$-transversal. 

Next, assume that $F_\alpha\in\F_3$ and consider the sequence $\big((x_1,y_1),(x_2,y_2),(x_3,y_3)\big)\defeq F_\alpha\in( \kappa\times G)^3$. The definition of the set $\F_3$ ensures that $|\{0,x_1,x_2,x_3\}|=4$. Let $(x_0,y'_0)\defeq(0,1_G)$ and choose points $y'_1,y'_2,y'_3\in G$ such that
$$
y'_i\notin\{y'_j:j<i\}\cup\{y'_j\Phi_\beta(x_j)^{-1}\Phi_\beta(x_i):j<i,\;\beta<\alpha\}$$
for every $i\in\{1,2,3\}$, and also
$$
y_3'\notin\{y_2'y_1'\}\cup\{y_2'y_2^{-1}y_1'y_1^{-1}y_3\}.$$
The choice of the points $y_0',y_1',y_2',y_3'$ ensures that $F'\defeq\{(x_0,y_0'),(x_1,y'_1),(x_2,y_2'),(x_3,y_3')\}$ is an injective function with $F'(0)=F(x_0)=y_0'=1_G$. We claim that $|F'\cap g\Phi_\beta|\le 1$ for every $g\in G$ and $\beta<\alpha$. To derive a contradiction, assume that $|F'\cap g\Phi_\beta|\ge 2$ for some $g\in G$ and $\beta<\alpha$. Then there exist two numbers $j<i\le 3$ such that $\{(x_i,y_i'),(x_j,y_j')\}\subseteq g\Phi_\beta$. It follows that $y_j'=g\Phi_\beta(x_j)$ and $y_i'=g\Phi_\beta(x_i)=y_j'\Phi_\beta(x_j)^{-1}\Phi_\beta(x_i)$, which contradicts the choice of the point $y_i'$. This contradiction shows that $|F'\cap g\Phi_\beta|\le 1$ for every $g\in G$ and $\beta<\alpha$. By Lemma~\ref{l:Krauz1}, there exists a bijective function $\Phi_\alpha:\kappa\to G$ such that $F'\subseteq \Phi_\alpha$ and the family $\{\Phi_\beta\}_{\beta\le\alpha}$ is $G$-transversal. Observe that $\Phi_\alpha(0)=F'(x_0)=y_0'=1_G$ and $$
\Phi_\alpha(x_3)=y_3'\notin\{y_2'y_1',y_2'y_2^{-1}y_1'y_1^{-1}y_3\}=\{\Phi_\alpha(x_2)\Phi_\alpha(x_1),\Phi_\alpha(x_2)y_2^{-1}\Phi_\alpha(x_1)y_1^{-1}y_3\}.
$$

Finally, assume that $\alpha\in\Omega_4$. Consider the sequence $(x_0,y_0,y_1,y_2,y_3)\defeq F_\alpha$. The definition of the set $\F_4$ ensures that $y_0\notin\{y_1,y_2,y_3\}$ and $y_1\ne y_2\ne y_3$.
Choose any point $x_1\in\kappa\setminus\{x_0\}$ such that  $$x_1\notin \bigcup_{\beta<\alpha}\bigcup_{i=1}^2\{x\in\kappa:y_0\Phi_\beta(x_0)^{-1}\Phi_\beta(x)=y_i\},$$and then choose a  point $x_2\in\kappa\setminus\{x_0,x_1\}$ such that $$x_2\notin\bigcup_{\beta<\alpha}\bigcup_{i=0}^1\{x\in \kappa:y_i\Phi_\beta(x_i)^{-1}\Phi_\beta(x)=y_2\}.$$
Consider the injective function $\varphi\defeq\{(x_0,y_0),(x_1,y_1),(x_2,y_2)\}$. The choice of the points $x_1,x_2$ ensures that $|\varphi\cap g\Phi_\beta|\le 1$ for every $g\in G$ and $\beta<\alpha$. By Lemma~\ref{l:Krauz1}, there exist a bijective function $\Phi_\alpha:\kappa\to G$ such that $\varphi\subseteq\Phi_\alpha$ and the family of bijections $\{\Phi_\beta\}_{\beta\le\alpha}$ is $G$-transversal.

Choose any point $y_3'\in G\setminus\{y_0,y_2,y_3\}$ such that
$$y_3'\notin\bigcup_{\beta\le\alpha}\{y_{0}\Phi_\beta(x_0)^{-1}\Phi_\beta(x_2),y_2\Phi_\beta(x_1)^{-1}\Phi_\beta(x_2)\},$$and consider the injective function $\varphi'\defeq\{(x_0,y_0),(x_1,y_2),(x_2,y_3')\}$.  The choice of the points $x_1,y_3'$ ensures that $|\varphi'\cap g\Phi_\beta|\le 1$ for every $g\in G$ and $\beta<\alpha$. By Lemma~\ref{l:Krauz1}, there exist a bijective function $\Phi_{\alpha+1}:\kappa\to G$ such that $\varphi'\subseteq\Phi_{\alpha+1}$ and the family of bijections $\{\Phi_\beta\}_{\beta\le\alpha+1}$ is $G$-transversal. It is easy to see that the functions $\Phi_\alpha$ and $\Phi_{\alpha+1}$ satisfy the inductive condition $(4_\alpha)$.

 \end{proof} 

Let $(\Phi_\alpha)_{\alpha\in\kappa}$ be the $G$-transversal family of bijections satisfying the conditions of Lemma~\ref{l:Krauz2}. Consider the liner $X=\kappa\times G$ endowed with the family of lines $$\mathcal L\defeq\big\{\{\alpha\}\times G:\alpha\in\kappa\big\}\cup\big\{\kappa\times\{g\}:g\in G\big\}\cup\big\{g\Phi_\alpha:g\in G,\;\alpha\in\kappa\big\}.$$The $G$-transversality of the family $\{\Phi_\alpha\}_{\alpha\in\kappa}$ and the conditions $(2_\alpha)_{\alpha\in\kappa}$ of Lemma~\ref{l:Krauz2} ensure that every distinct points of the set $X$ are contained in a unique line, so $X$ is indeed a liner. The $G$-transversality of the family $\{\Phi_\alpha\}_{\alpha\in\kappa}$ implies that $X$ is a Playfair plane such that for every $g\in G$ the map $T_g:X\to X$, $T_g:(x,y)\mapsto (x,gy)$, is a translation. Then $\{T_g:g\in G\}$ is a subgroup of $\Trans(X)$, isomorphic to $G$. We claim that $\Trans(X)=\{T_g:g\in G\}$. To derive a contradiction, assume that there exists a translation $T\in\Trans(X)\setminus\{T_g:g\in G\}$.
Consider the point $z_0\defeq(0,1_G)\in \kappa\times G=X$. Assuming that $T(z_0)\in (0,g)$ for some $g\in G$, we conclude that $T(z_0)=T_g(z_0)$ and hence $T=T_g$, which contradicts the choice of the translation $T$. Therefore, $z_1\defeq T(z_0)\in X\setminus(\{0\}\times G)$. 

Consider the line $L\defeq\Aline {z_0}{z_1}$ in the liner $X$ and observe that $L\cap(\{0\}\times G)=\{z_0\}$. Since $T$ is a translation with $T(z_0)=z_1$, $T[L]=L$. Choose a point $z_2\in L\setminus\{z_0,z_1,T^{-1}(z_0),T^{-1}(z_1)\}$ and put $z_3\defeq T(z_2)$. For every $i\in 4$, write the point $z_i\in X=\kappa\times G$ as $z_i=(x_i,y_i)$ for some points $x_i\in\kappa$ and $y_i\in G$. Then the sequence $$F\defeq \big((x_1,y_1),(x_2,y_2),(x_3,y_3)\big)$$ is an element of the family $\F_3$. So, there exists an ordinal $\alpha\in\Omega_3$ such that $F_\alpha=F$. 
The inductive condition $(3_\alpha)$ ensures that $z_0\defeq(0,1_G)\in\Phi_\alpha$.

\begin{picture}(100,160)(-130,-15)

\put(0,0){\line(2,1){80}}
\put(0,0){\line(0,1){100}}
\put(-3,105){$G$}
\put(0,20){\line(2,1){80}}
\put(85,40){$L$}
\put(0,0){\line(2,3){50}}
\put(50,25){\line(0,1){50}}
\put(20,10){\line(0,1){20}}
\put(70,35){\line(0,1){78}}

\qbezier(50,75)(70,105)(75,130)
\put(80,130){$\Phi_\alpha$}

\put(70,113){\color{red}\circle*{3}}
\put(0,0){\vector(0,1){20}}
\put(0,0){\vector(2,1){20}}

\put(50,25){\vector(0,1){20}}
\put(50,25){\vector(2,1){20}}
\put(0,0){\vector(2,3){20}}
\put(50,25){\vector(2,3){20}}
\put(50,75){\vector(2,3){20}}

\put(0,0){\circle*{3}}
\put(-3,-8){$z_0$}
\put(0,20){\circle*{3}}
\put(-8,18){$g$}
\put(17,2){$z_1$}
\put(20,10){\circle*{3}}
\put(50,25){\circle*{3}}
\put(47,17){$z_2$}
\put(70,35){\circle*{3}}
\put(67,27){$z_3$}
\end{picture}

Consider the element $g\defeq \Phi_\alpha(x_1)y_1^{-1}$ of the group $G$, and the translation $T_g:X\to X$, $T_g:(x,y)\mapsto (x,gy)$.  Theorem~\ref{t:Trans-commutative} implies that $TT_g=T_gT$ is a translation. Observe that $T_gT(z_0)=T_g(z_1)=T_g(x_1,y_1)=(x_1,gy_1)=(x_1,\Phi_\alpha(x_1)y_1^{-1}y_1)=(x_1,\Phi_\alpha(x_1))\in\Phi_\alpha$ and hence $\overline{\{z_0,T_gT(z_0)\}}=\Phi_\alpha$. By Proposition~\ref{p:Trans-spread}, $\overline{\{z_2,T_gT(z_2)\}}\parallel\Phi_\alpha$ and hence $\{z_2,T_g(z_3)\}=\{z_2,T_gT(z_2)\}\subseteq y_2\Phi_\alpha(x_2)^{-1}\Phi_\alpha$. Then $gy_3=y_2\Phi_\alpha(x_2)^{-1}\Phi_\alpha(x_3)$ and
$$\Phi_\alpha(x_3)=\Phi_\alpha(x_2)y_2^{-1}gy_3=\Phi_\alpha(x_2)y_2^{-1}\Phi_\alpha(x_1)y_1^{-1}y_3,$$
which contradicts the condition $(3_\alpha)$ of Lemma~\ref{l:Krauz2}.
This is a final contradiction showing that $\Trans(X)=\{T_g:g\in G\}$.
\smallskip

Assuming that $\Dil(X)\ne\Trans(X)$, we can find a dilation $D\in\Dil(X)\setminus \Trans(X)$, which has a unique fixed point $(x_0,y_0)=D(x_0,y_0)\in \kappa\times G$. Choose any point $y_1\in G\setminus\{y_0\}$. 

\begin{picture}(100,100)(-150,-15)

\put(0,0){\line(1,0){100}}
\put(105,-3){$\kappa$}
\put(0,0){\line(0,1){80}}
\put(-3,85){$G$}
\put(0,0){\line(1,1){70}}
\put(0,36){\line(1,0){70}}
\put(0,48){\line(1,0){70}}
\put(0,64){\line(1,0){70}}
\put(64,0){\line(0,1){70}}
\put(48,0){\line(0,1){70}}
\qbezier(0,0)(40,15)(64,48)
\put(64,48){\line(3,4){10}}

\put(0,0){\circle*{3}}
\put(-12,0){$y_0$}
\put(-1,-8){$x_0$}
\put(48,0){\circle*{3}}
\put(45,-8){$x_2$}
\put(64,0){\circle*{3}}
\put(61,-8){$x_1$}
\put(0,36){\circle*{3}}
\put(-12,33){$y_3$}
\put(0,48){\circle*{3}}
\put(-12,45){$y_2$}
\put(0,64){\circle*{3}}
\put(-12,61){$y_1$}

\put(36,36){\circle*{3}}
\put(48,48){\circle*{3}}
\put(64,64){\circle*{3}}
\put(64,48){\circle*{3}}
\put(48,30){\color{red}\circle*{3}}

\put(72,73){$\Phi_\alpha$}
\put(78,60){$\Phi_{\alpha+1}$}

\end{picture}

Taking into account that $D$ is a dilation with fixed point $(x_0,y_0)\in \{x_0\}\times G\in\mathcal L$, we conclude that $D[\{x_0\}\times G]=\{x_0\}\times G$ and hence $D(x_0,y_1)=(x_0,y_2)$ for some point $y_2\in G\setminus\{y_0,y_1\}$ and $D(x_0,y_2)=(x_0,y_3)$ for some point $y_3\in G\setminus\{y_0,y_2\}$. Then $(x_0,y_0,y_1,y_2,y_3)\in\F_4$ and hence $(x_0,y_0,y_1,y_2,y_3)=F_\alpha$ for some $\alpha\in\Omega_4$. By the inductive condition $(4_\alpha)$, there exist  distinct points $x_1,x_2\in \kappa$ such that $\{(x_0,y_0),(x_1,y_1),(x_2,y_2)\}\subseteq\Phi_\alpha$, $\{(x_0,y_0),(x_1,y_2)\}\subseteq\Phi_{\alpha+1}$ and $(x_2,y_3)\notin\Phi_{\alpha+1}$. Taking into account that $D$ is a dilation with $(x_0,y_0)\in \Phi_\alpha\cap\Phi_{\alpha+1}$, we conclude that $D[\Phi_\alpha]=\Phi_\alpha$ and $D[\Phi_{\alpha+1}]=\Phi_{\alpha+1}$. On the other hand, $D(x_0,y_1)=(x_0,y_2)$ implies $D[\kappa\times\{y_1\}]=\kappa\times\{y_2\}$ and hence $D(x_1,y_1)\in D[(\kappa\times\{y_1\})\cap\Phi_\alpha]=(\kappa\times\{y_2\})\cap\Phi_\alpha=\{(x_2,y_2)\}$.
Then $D[\{x_1\}\times G]=\{x_2\}\times G$ and hence $$D(x_1,y_2)\in D[(\{x_1\}\times G)\cap(\kappa\times\{y_2\})\cap\Phi_{\alpha+1}]=(\{x_2\}\times G)\cap(\kappa\times \{y_3\})\cap\Phi_{\alpha+1}=\{(x_2,y_3)\}\cap\Phi_{\alpha+1}=\varnothing,$$which is a desired contradiction showing that $\Dil(X)=\Trans(X)$. Now we see that the group $\Dil(X)=\Trans(X)=\{T_g:g\in G\}$ is isomorphic to the group $G$.
\end{proof}

\begin{remark} The first exampe of a Playfair plane with non-commutative group of translations was constructed by \index[person]{Pickert}Pickert\footnote{{\bf G\"unter Pickert} (1917--2015), a German mathematician. Pickert went to school in Eisenach and from 1933 studied mathematics and physics at the University of Göttingen (among other things, he heard David Hilbert's lecture on the basics of geometry in 1933/34) and the TH Danzig. In 1939 he received his doctorate in G\"ottingen under Helmut Hasse (New Methods in the Structural Theory of Commutative-Associative Algebras, Mathematical Annals Vol. 116, p. 217). During the Second World War he was a soldier in Poland, Russia (Stalingrad) and Tunisia, most recently as a first lieutenant. From 1943 to 1946, Pickert was a prisoner of war in the USA, where he gave lectures at the``camp university”. From 1946 he was an assistant at the University of T\"ubingen, where he completed his habilitation in 1948, became a lecturer and in 1953 an adjunct professor (he also gave substitute lectures in G\"ottingen in 1950 and in Heidelberg in 1951). From 1962 he was a full professor and director of the Mathematical Institute at the University of Giessen, where he retired in 1985. Among other things, Pickert dealt with finite geometries, which he also brought into mathematics didactics. In 1955 his influential textbook on projective planes was published. In 1991 he was awarded an honorary doctorate from the Julius Maximilian University of Würzburg. In 1988 he received the Federal Cross of Merit, 1st Class.
One of his doctoral students was Johannes Andr\'e who known for the description of the structure of Thalesian affine planes and also for constructing Andr\'e planes, examples of non-Desarguesian Thalesian planes of order $p^n$. The first example of a Playfair plane with non-commutative translation group was constructed by Pickert in [G.~Pickert, {\em  Nichtkommutative cartesische Gruppen}, Arch. Math. {\bf 3} (1952), 335--342].} in 1952. 
\end{remark}

Theorem~\ref{t:Krauz} suggests two open problems.

\begin{problem} Let $G$ be a finite group. Is there a (finite) Playfair plane $X$ whose translation group $\Trans(X)$ is isomorphic to the group $G$?
\end{problem}

\begin{problem} Let $G$ be a group. Is there a Playfair plane $X$ such that $\Aut(X)=\Trans(X)$ and the group $\Aut(X)=\Trans(X)$ is isomorphic to the group $G$?
\end{problem}

\section{Translation and dilation liners}

\begin{definition} A liner $X$ is called a \index{translation liner}\index{liner!translation}\defterm{translation} if for every points $x,y\in X$ there exists a translation $T:X\to Y$ such that $T(x)=y$.
\end{definition}

\begin{proposition}\label{p:2=>translation} Every liner $X$ of rank $\|X\|\le 2$ is translation.
\end{proposition}

\begin{proof} Given two points $e,a\in X$, we need to find a translation $T:X\to X$ such that $T(e)=a$. Since there exist abelian groups of any cardinality, there exists a binary operation $\cdot:X\times X\to X$ such that $(X,\cdot)$ is a group with identity $e$. Then the function $T:X\to X$, $T:x\mapsto x\cdot a$, is bijective and $T(e)=e\cdot a=a$. It remains to prove that $T$ is a translation of the liner $X$. If $e=a$, then $T$ is the identity map of $X$ and hence $T$ is the identity translation of $X$. If $e\ne a$, then $2=|\{e,a\}|\le|X|\le\|X\|\le 2$ and $X$ is a unique line of the liner $X$. Since $T[X]=X\parallel X$, the bijection is a dilation of $X$. Since $T$ has no fixed points, the dilation $T$ is a translation of $X$.
\end{proof}

\begin{proposition}\label{p:translation=>Bolyai} Every translation (Proclus) liner $X$ is Bolyai (and Playfair).
\end{proposition}

\begin{proof} Assume that $X$ is a translation liner. To prove that $X$ is Bolyai, we need to show that for every plane $\Pi$, line $L\subseteq \Pi$ and point $x\in\Pi\setminus L$ there exists a line $\Lambda$ such that $x\in \Lambda\subseteq \Pi\setminus L$.  Choose any point $o\in L$. Since  $X$ is a translation liner, there exists a translation $T:X\to X$ such that $T(o)=x$. Since $T$ is a dilation, $\Lambda\defeq T[L]$ is a line in $X$ such that $\Lambda\parallel L$. Observe that $x=T(o)\in T[L]=\Lambda$ and $\Lambda\subseteq \overline{L\cup\{x\}}\subseteq \Pi$. Assuming that $\Lambda\cap L\ne\varnothing$, we can apply Proposition~\ref{p:para+intersect=>coincide} and conclude that $x\in \Lambda=L$, which contradicts the choice of the point $x$. This contradiction shows that $x\in\Lambda\subseteq \Pi\setminus L$, witnessing that the liner $X$ is Bolyai. If $X$ is Proclus, then $X$ is Playfair because a liner is Playfair if and only if it is Proclus and Bolyai, according to Definition~\ref{d:PPBL}. 
\end{proof}

\begin{theorem}\label{t:translation=>Thalesian} Every translation Proclus liner is Thalesian.
\end{theorem}

\begin{proof} Let $X$ be a translation Proclus liner. To prove that $X$ is Thalesian, take any distinct parallel coplanar lines $A,B,C$ and any points $a,a'\in A$, $b,b'\in B$, $c,c'\in C$ such that $\Aline ab\cap\Aline{a'}{b'}=\varnothing=\Aline bc\cap\Aline{b'}{c'}$ and hence $a\ne a'$, $b\ne b'$, and $c\ne c'$. Since the liner $X$ is translation, there exists a translation $T:X\to X$ such that $T(b)=b'$. Let $a''\defeq T(a)$. By Proposition~\ref{p:Trans-spread}, the lines $\Aline a{a''}$, $\Aline b{b''}$ belong to the spread of parallel lines $\overline T\defeq\{\Aline xy:xy\in T\}$ and hence $\Aline a{a''}\parallel \Aline b{b'}=B$. Since $a\in A\cap\Aline a{a''}$ and $A\parallel B\parallel \Aline a{a''}$, the lines $A$ and $\Aline a{a''}$ coincide by the Proclus Axiom. Since $T$ is a translation, the lines $\Aline ab$ and $\Aline {a''}{b'}=T[\Aline ab]$ are parallel. Then $\Aline {a'}{b'}\parallel \Aline ab\parallel \Aline{a''}{b'}$ and hence $\Aline{a'}{b'}=\Aline{a''}{b'}$, by the Proclus Axiom. Now we see that $a''=\Aline a{a''}\cap\Aline {a''}{b'}=A\cap\Aline{a'}{b'}=\{a'\}$ and hence $a''=a'$. By analogy we can prove that $c'=c''\defeq T(c)$. Since $T$ is a translation, $\Aline {a'}{c'}=\Aline{a''}{c''}=T[\Aline ac]\parallel \Aline ac$. Assuming that $\Aline {a'}{c'}\cap\Aline ac\ne\varnothing$, we can apply Proposition~\ref{p:para+intersect=>coincide} and conclude that $\Aline {a'}{c'}=\Aline ac$ and then $\{a\}=A\cap\Aline ac=A\cap\Aline{a'}{c'}=\{a'\}$, which is a contradiction showing that the parallel lines $\Aline ac$ and $\Aline{a'}{c'}$ are disjoint. 
\end{proof}

\begin{remark} In Theorem~\ref{t:paraD<=>translation}, we shall prove that a $3$-long affine regular liner is translation if and only if it is Thalesian.
\end{remark}

Theorem~\ref{t:Trans-commutative} implies the following corollary.

\begin{corollary}\label{c:Trans-commutative} The translation group $\Trans(X)$ of any translation proaffine regular liner $X$ of rank $\|X\|\ge 3$ is commutative.
\end{corollary}

For every point $o$ in a proaffine regular liner $X$, the set\index[note]{$\Dil_o(X)$}
$$
\Dil_o(X)\defeq\{D\in \Dil(X):D(o)=o\}$$
is a subgroup of the dilation group $\Dil(X)$. The definition of a translation implies the equality $\Dil_o(X)\cap\Trans(X)=\{1_X\}$.


\begin{proposition}\label{p:Dilo=Dil/Trans} Let $X$ be a translation proaffine regular liner of rank $\|X\|\ge 3$. For every point $o\in X$, $$\Dil(X)=\Trans(X)\cdot \Dil_o(X)=\Dil_o(X)\cdot\Trans(X)$$ and the group $\Dil_o(X)$ is isomorphic to the quotient group $\Dil(X)/\Trans(X)$.
\end{proposition}

\begin{proof} The equality $\Dil(X)=\Trans(X)\cdot\Dil_o(X)=\Dil_o(X)\cdot\Trans(X)$ will follow as soon as we find for every dilation $D\in\Dil(X)$ a homothety $H\in\Dil_o(X)$ and translations $T,T'\in\Trans(X)$ such that $D=TH=HT'$. So, fix any dilation $D$ of $X$.  If $D$ is a translation of $X$, then the translations $T=T'=D$ and the homothety $H\defeq 1_X$ have the required property: $D=TH=HT'$. So, assume that $D$ is not a translation and hence $D(p)=p$ for a unique point $p\in X$, according to Theorem~\ref{t:dilation-atmost1}. Since $X$ is a translation liner, there exists a translation $P:\Aut(X)$ such that $P(o)=p$. Then the dilation $H\defeq P^{-1}DP$ has $H(o)=P^{-1}DP(o)=P^{-1}D(p)=P^{-1}(p)=o$ and hence $H\in\Dil_o(X)$. By Proposition~\ref{p:Trans(X)isnormal}, the group $\Trans(X)$ is normal in the group $\Aut(X)$ and hence the automorphisms $T\defeq P(HP^{-1}H^{-1})$ and $T'\defeq (H^{-1}PH)P^{-1}$ are translations of $X$ such that  $D=PHP^{-1}=PHP^{-1}H^{-1}H=TH$ and $D=PHP^{-1}=HH^{-1}PHP^{-1}=HT'$.
Therefore, $\Dil(X)=\Trans(X)\cdot\Dil_o(X)=\Dil_o(X)\cdot \Trans(X)$.

Since $\Trans(X)$ is a normal subgroup of the group $\Dil(X)$, we can consider the quotient group $\Dil(X)/\Trans(X)$, the  quotient homomorphism $Q:\Dil(X)\to\Dil(X)/\Trans(X)$, and its restriction $Q_o\defeq Q{\restriction}_{\Dil_o(X)}$. Since $\Trans(X)\cap\Dil_o(X)=\{1_X\}$ and $\Dil(X)=\Trans(X)\cdot\Dil_o(X)$, the homomorphism $Q_o:\Dil_o(X)\to\Dil(X)/\Trans(X)$ is bijective and hence $Q_o$ is an isomorphism of the groups $\Dil_o(X)$ and $\Dil(X)/\Trans(X)$.
\end{proof} 

\begin{example}[Ivan Hetman] There exists a Playfair plane of order 9 whose translation group is trivial and the dilation group $\Dil(X)$ has only two elements. This space contains a unique point $o\in X$ such that $\Dil(X)=\Dil_o(X)$ and for every point $p\in X\setminus\{o\}$ we have $\Dil_p(X)=\Trans(X)\ne \Dil(X)$.
\end{example}

\begin{problem} Find an example of an affine plane $X$ such that $\Dil(X)\ne \Trans(X)\cdot\Dil_o(X)$ for every point $o\in X$.
\end{problem}

\begin{definition} A liner $X$ is called a \index{dilation liner}\index{liner!dilation}\defterm{dilation liner} if for every distrinct collinear points $x,y,o\in X$ there exists a dilation $D:X\to X$ such that $Dox=oy$.
\end{definition}

\begin{remark} In Theorem~\ref{t:Des<=>dilation} we shall prove that an affine liner $X$ is dilation if and only if $X$ is Desarguesian.
\end{remark}

\section{Central and paracentral automorphisms of liners}

\begin{definition} Let $A:X\to X$ be an automorphism of a liner $X$. A point $c\in X$ is called a \index{centre}\defterm{centre} of the automorphism $A$ if $A[\Aline xc]=\Aline xc$ for every point $x\in X$. An automorphism $A$ of a liner $X$ is called \index{central automorphism}\index{automorphism!central}\defterm{central} if it has a centre.
\end{definition}

\begin{example} Every point $c\in X$ of a liner $X$ is a centre of the identity automorphism $1_X:X\to X$. 
\end{example}

\begin{proposition}\label{p:cenralauto-has1-center} Let $X$ be a liner of rank $\|X\|\ge 3$. A non-identity automorphism $A:X\to X$ can have at most one centre.
\end{proposition}

\begin{proof} Assuming that an automorphism $A:X\to X$ has two distinct centres 
$a,b\in X$, we shall prove that $A(x)=x$ for all $x\in X$. 

\begin{claim}\label{cl:Ax=x} For every $x\in X\setminus \Aline ab$, we have $A(x)=x$.
\end{claim}

\begin{proof} It follows from $x\notin\Aline ab$ that $\{x\}=\Aline xa\cap \Aline xb$ and hence $$\{A(x)\}=A[\Aline xa\cap \Aline xb]=A[\Aline xa]\cap A[\Aline xb]=\Aline xa\cap\Aline xb=\{x\}.$$ 
\end{proof}

\begin{claim}\label{cl:A-centre} Every point $c\in X\setminus\Aline ab$ is a centre of the automorphism $A$.
\end{claim}

\begin{proof} Given any point $x\in X$, we need to show that $A[\Aline xc]=\Aline xc$. If $x=c$, then $A[\Aline xc]=A[\{c\}]=\{A(c)\}=\{c\}=\Aline xc$, by Claim~\ref{cl:Ax=x}. So, assume that $x\ne c$. Since the liner $X$ is $3$-long, the line $\Aline xc$ contains some point $y \notin \Aline ab\cup\{c\}$. Claim~\ref{cl:Ax=x} implies that $A[\Aline xc]=A[\Aline yc]=\Aline yc=\Aline xc$, witnessing that $c$ is a centre of the automorphism $A$.
\end{proof}

Now we are ready to prove that $A(x)=x$ for every $x\in X$. If $x\notin\Aline ab$, then the equality $A(x)=x$ follows from Claim~\ref{cl:Ax=x}. So, assume that $x\in \Aline ab$. Since $\|X\|\ge 3$, there exist points $c\in X\setminus \Aline ab$ and $d\in \{a,b\}\setminus\{x\}$. By Claim~\ref{cl:A-centre}, the point $c$ is a centre of the automorphism $A$. Since $x\notin \Aline cd$ and the points $c,d$ are centres of the automorphism $A$, we can apply Claim~\ref{cl:Ax=x} and conclude that $A(x)=x$.
\end{proof}

Let us recall that a liner $X$ is called \defterm{hyperaffine} (resp. \defterm{affine}) if for every points $o,x,y\in X$ and $p\in \Aline xy\setminus \Aline ox$ there exists a (unique) point $u\in \Aline oy$ such that $\Aline up\cap\Aline ox=\varnothing$. It is clear that every affine liner is hyperaffine. 

\begin{theorem}\label{t:central<=>homothety} An automorphism $A:X\to X$ of a \textup{(}$3$-ranked hyperaffine\textup{)} liner $X$ is central if \textup{(}and only if\textup{)} $A$ is a homothety.
\end{theorem}

\begin{proof} To prove the ``if'' part, assume that $A$ be a homothety and let $c$ be its centre. Given any point $x\in X$, we should prove that $A[\Aline xc]=\Aline xc$. If $x=c$, then $A[\Aline xc]=A[\{c\}]=\{A(c)\}=\{c\}=\Aline xc$. If $x\ne c$, then $\Aline xc$ is a line. Since the automorphism $A$ is a dilation, $A[\Aline xc]\parallel \Aline xc$. Since $c\in \Aline ac\cap A[\Aline xc]$, the parallelity relation $A[\Aline xc]\parallel \Aline xc$ implies $A[\Aline xc]=\Aline xc$, witnessing that the automorphism $A$ is central.
\smallskip

To prove the ``only if'' part, assume that the automorphism $A$ is central, and the liner $X$ is  $3$-ranked and hyperaffine. Let $c$ be a centre of $A$. Then $A[L]=L$ for every line $L\subseteq X$ containing the point $c$. To prove that $A$ is a dilation of $X$, we need to check that $A[L]\parallel L$ for every line $L\subseteq X$. If $A[L]=L$, then $A[L]\parallel L$ and we are done. So, assume that $A[L]\ne L$. In this case $c\notin L\cup A[L]$.  We claim that $A[L]\cap L=\varnothing$. In the opposite case, there exists a unique point $y\in A[L]\cap L$. If $L=\{y,A^{-1}(y)\}$, then $A[L]=\{y,A(y)\}\subseteq A[\Aline yc]=\Aline yc$ and $c\in \Aline yc=A[L]$, which contradicts our assumption. So, $L\ne \{y,A^{-1}(y)\}$ and hence there exists a point $x\in L\setminus\{y,A^{-1}(y)\}$ whose image $x'\defeq A(x)$ is distinct form the points $x$ and $y$. Assuming that $y\in \Aline xc$, we conclude that $c=\Aline xy=L$, which contradicts our assumption. This contradiction shows that $y\notin \Aline xc$. Let $y^-\defeq A^{-1}(y)\in L=\Aline xy$. Since $A$ is central, $y^-=A^{-1}(y)\in A^{-1}[\Aline yc]=\Aline yc$ and hence $y^-\in \Aline xy\cap \Aline yc=\{y\}$ and $y=A(y^-)=A(y)$. Since the liner $X$ is hyperaffine, there exists a point $z\in \Aline {x'}{y}$ such that $\Aline {c}{z}\cap\Aline xy=\varnothing$. Then 
$$z\in \Aline {x'}{y}\cap\Aline cz=A[\Aline xy]\cap A[\Aline cz]=A[\Aline xy\cap\Aline cz]=A[\varnothing]=\varnothing,$$which is a contradiction showing that $A[L]\cap L=\varnothing$. Choose any point $u\in L\setminus \{x\}$ and consider its image $u'\defeq A(u)\in A[\Aline uc]=\Aline uc$. Observe that $A[L]=A[\Aline xu]=\Aline {x'}{u'}\in\overline{\{x,c,u\}}$ and hence $\|L\cup A[L]\|\le\|\{x,c,u\}\|=3$. Applying Corollary~\ref{c:parallel-lines<=>}, we conclude that the disjoint coplanar lines $L$ and $A[L]$ in the $3$-ranked liner $L$ are parallel, witnessing that the automorphism $A$ is a dilation. Since $A(c)=c$, the dilation $A$ is a homothety. 
\end{proof}

A family of lines $\F$ in a liner $X$ is called
\begin{itemize}
\item \index{lines!concurrent}\index{concurrent lines}\defterm{concurrent} if $\exists c\in X\;\;\forall L\in\F\;\;(c\in L)$;
\item \index{lines!parallel}\index{parallel lines}\defterm{parallel} if $\forall L,\Lambda\in\F\;\;(L\parallel \Lambda)$;
\item \index{lines!paraconcurrent}\index{paraconcurrent lines}\defterm{paraconcurrent} $\F$ is parallel or concurrent.
\end{itemize}

\begin{proposition} An automorphism $A:X\to X$ of a $3$-long liner $X$ of rank $\|X\|\ge 3$ is central if and only if the family of lines $\overline A\defeq \{\Aline xy:xy\in A\setminus 1_X\}$ is concurrent.
\end{proposition}

\begin{proof} If an automorphism $A:X\to X$ is central, then there exists a point $c\in X$ such that $A[\Aline xc]=\Aline xc$ for all $x\in X$. Then for every pair $xy\in A\setminus 1_X$, we have $y=A(x)\in A[\Aline xc]=\Aline xc$ and hence $c\in \Aline xc=\Aline xy$, witnessing that the family of lines $\overline A$ is concurrent.   

Now assume conversely that the family of lines $\overline A$ is concurrent and find a point $c\in X$ such that $c\in\Aline xy$ for every pair $xy\in A\setminus 1_X$.

\begin{claim} $A(c)=c$.
\end{claim}

\begin{proof} To derive a contradiction, assume that $c'\defeq A(c)\ne c$. Then $c''\defeq A(c')\ne c'=A(c)$.  We claim that there exists a point $b\in X\setminus \Aline c{c'}$ such that $A(b)\ne b$. Since $\|X\|\ge 3>2=\|\Aline c{c'}\|$, there exists a point $a\in X\setminus\Aline c{c'}$. If $A(a)\ne a$, then put $b\defeq a$. If $F(a)=a$, then take any point $b\in \Aline ac\setminus\{a,c\}$. Such a point $b$ exists because the liner $X$ is $3$-long. Then $b'\defeq A(b)\in A[\Aline ac\setminus\{a\}]=\Aline a{c'}\setminus\{a\}=\Aline a{c'}\setminus \Aline ac\subseteq X\setminus\{b\}$ and $b''\defeq A(b')\ne A(b)=b'$.  The choice of the point $c$ ensures that $$c\in\Aline c{c'}\cap\Aline {c'}{c''}\cap\Aline b{b'}\cap\Aline {b'}{b''}\subseteq \Aline {c'}{c''}\cap\Aline{b'}{b''}=A[\Aline c{c'}]\cap A[\Aline b{b'}]=A[\Aline c{c'}\cap\Aline b{b'}]=A[\{c\}]=\{A(c)\},$$which contradicts our assumption.
\end{proof}

To see that the automorphism $A$ is central, it suffices to show that $A[\Aline xc]=\Aline xc$. If $x\in\Fix(A)$, then $Axc=xc$ and hence $A[\Aline xc]=\Aline xc$. If $x\notin\Fix(A)$, then for the point $x'\defeq A(x)$, we have $xx'\in A\setminus 1_X$. The choice of the point $c$ ensures $c\in \Aline x{x'}$ and hence $A[\Aline xc]=\Aline {x'}c=\Aline xc$.
\end{proof}


\begin{definition} An automorphism $A:X\to X$ of a liner $X$ is called \index{paracentral automorphism}\index{automorphism!paracentral}\defterm{paracentral} if the family $\overline A\defeq\{\Aline xy:xy\in A\setminus 1_X\}$ is parallel. A spread of parallel lines containing the family $\overline A$ is called a \defterm{central direction} of the paracenral automorphism. A non-identity paracentral automorphism $A$ of a liner $X$ can have at most one central direction, which will be called \index{central direction}\index{paracentral automorphism!central direction}\defterm{the central direction} of $X$.
\end{definition}

\begin{theorem} An automorphism $A:X\to X$ of a $3$-ranked liner $X$ is a translation if and only if $A$ is paracentral and $\Fix(A)\in\{\varnothing,X\}$.
\end{theorem}

\begin{proof} The ``only if'' part follows from Proposition~\ref{p:Trans-spread}. To prove the ``if'' part, assume that the automorphism $A$ is paracentral and $\Fix(A)\in \{\varnothing, X\}$. If $\Fix(A)=X$, then $A$ is the identity automorphism and hence $A$ is a translation of $X$, by definition. So, assume that $\Fix(A)=\varnothing$. 

Given any line $L\subseteq X$, we should prove that $A[L]\parallel L$. Fix any distinct points $x,y\in L$ and consider their images $x'\defeq A(x)$ and $y'\defeq A(y)$. Since $A$ is paracentral and $x,y\in X=X\setminus\Fix(A)$, the lines $\Aline x{x'}$ and $\Aline y{y'}$ are parallel and hence $L=\Aline xy$ and $A[L]=\Aline{x'}{y'}$ are two lines in the flat $\Pi=\overline{\{x,x',y,y'\}}=\overline{\{x,x',y\}}=\overline{\{y,y',x\}}$ of rank $\|P\|\le 3$. If $\|P\|=2$, then $A[L]=\Pi=L$ and hence $A[L]\parallel L$. So, assume that $\|\Pi\|=3$, which implies $\Aline x{x'}\cap\Aline y{y'}=\varnothing$. Assuming that $A[L]\cap L\ne\varnothing$, we can find a point $z'\in A[L]\cap L=\Aline xy\cap \Aline {x'}{y'}$. Consider the point $z\defeq A^{-1}(z')\in L=\Aline xy$ and observe that $L=\Aline z{z'}\parallel \Aline x{x'}$ and hence $\Aline xy=L=\Aline x{x'}$, which contradicts $\Aline y{y'}\cap\Aline x{x'}=\varnothing$. This contradition shows that $A[L]\cap L=\varnothing$ and hence the disjoint lines $A[L],L$ in the plane $\Pi$ are parallel, by Corollary~\ref{c:parallel-lines<=>}. Therefore, $A$ is a dilation. Since $\Fix(A)=\varnothing$, the dilation $A$ is a translation.
\end{proof}

\begin{proposition}\label{p:central|=(para)central} Let $Y$ be a projective completion of a $3$-long $3$-ranked liner $X$ of rank $\|X\|\ge 3$, and $A:Y\to Y$ be an automorphism of the projective space $Y$ such that $A[X]=X$. The automorphism $A:Y\to Y$ of $Y$ is central if and only if the automorphism $A{\restriction}_X:X\to X$ of $X$ is central or paracentral.
\end{proposition}

\begin{proof} Assume that the automorphism $A:Y\to Y$ is central and let $c$ be a centre of $A$. If $c\in X$, then for every $x\in X$ we have $A[\Aline xc\cap X]=A[\Aline xc]\cap A[X]=\Aline xc\cap X$, witnessing that $c$ is a centre of the automorphism $A{\restriction}_X$. If $c\in Y\setminus X$, then for every point $x\in X\setminus \Fix(A)$, the point $x'\defeq A(x)$ belongs to the set $(\Aline xc\cap X)\setminus\{x'\}$ and hence $\Aline xc\cap X=\Aline x{x'}\cap X$ is a line in $X$. Given any points $x,y\in X\setminus\Fix(A)$ and their images $x'\defeq A(x)$ and $y'\defeq A(y)$, we need to check that the lines $\Aline x{x'}\cap X=\Aline xc\cap X$ and $\Aline y{y'}\cap X=\Aline yc\cap X$ are parallel. If $\Aline xc=\Aline yc$, then the lines $\Aline x{x'}\cap X$ and $\Aline y{y'}\cap X$ coincide and hence are parallel. If $\Aline xc\ne\Aline yc$, then $$(\Aline x{x'}\cap X)\cap(\Aline y{y'}\cap X)=(\Aline xc\cap\Aline yc)\cap X=\{c\}\cap X=\varnothing$$and hence the lines $\Aline x{x'}\cap X$ and $\Aline y{y'}\cap X$ in $X$ are disjoint. On the other hand, by Theorem~\ref{t:procompletion=>normal}, the set $\Pi\defeq\overline{\{x,o,y\}}\cap X$ is a plane in the $3$-ranked liner $X$ containing the disjoint lines $\Aline x{x'}\cap X=\Aline xc\cap X$ and $\Aline y{y'}\cap X=\Aline yc\cap X$. By Corollary~\ref{c:parallel-lines<=>}, the disjoint coplanar lines $\Aline x{x'}\cap X$ and $\Aline y{y'}\cap X$ are parallel, witnessing that the automorphism $A{\restriction}_X:X\to X$ is paracentral.
\smallskip

Now assume that the automorphism $A{\restriction}_X:X\to X$ is central and let $c\in X$ be a centre of the automorphism $A{\restriction}_X$. Given any point $y\in Y$, we need to check that $A[\Aline yc]=\Aline yc$. This equality is trivially true if $y=c$. So, assume that $y\ne c$.  If $|\Aline yc\cap X|\ge 2$, then $\Aline yc\cap X$ is a line in $X$. Taking into account that $c$ is a centre of the automorphism $A{\restriction}_X$, we conclude that $A[\Aline yc\cap X]=\Aline yc\cap X$. Then  the lines $A[\Aline yc]$ and $\Aline yc$ coincide because $|A[\Aline yc]\cap\Aline yc|\ge|\Aline yc\cap X|\ge 2$. So, assume that $\Aline yc\cap X=\{c\}$.
Since $\overline{Y\setminus X}\ne Y$, there exists a point $o\in Y\setminus\overline{Y\setminus X}\subseteq X$. For every $z\in \Aline yc\setminus\{c\}\subseteq Y\setminus X$ we have $\Aline oz\setminus\{z\}=\Aline oz\cap X$ and hence $A[\Aline oz\setminus \{z\}]=A[\Aline oz\cap X]=\Aline oz\cap X$, $A[\Aline oz]=\Aline oz$ and finally, $A[\{z\}]=A[\Aline oz\setminus X]=\Aline oz\setminus X=\{z\}$. Therefore, $A(z)=z$ for every point $z\in \Aline yc$ and hence $A[\Aline yc]=\Aline yc$, witnessing that $c$ is a centre of the automorphism $A$.

Finally, assume that the automorphism $A{\restriction}_X:X\to X$ is paracentral but not central. Then $\Fix(A{\restriction}_X)\ne X$ and for every $xx',yy'\in A{\restriction}_X\setminus 1_X$, the lines $\Aline x{x'}\cap X$ and $\Aline y{y'}\cap X$ in $X$ are parallel. Since $\Fix(A{\restriction}_X)\ne X$, there exists a point $a\in X\setminus\Fix(A)$. Two cases are possible. 

First assume that there exists a point $b\in X\setminus (\Aline a{a'}\cup\Fix(A))$. In this case the parallel lines $\Aline a{a'}\cap X$ and $\Aline b{b'}\cap X$ are disjoint.  Since the liner $Y$ is projective, there exists a point $c\in (\Aline a{a'}\cap\Aline b{b'})\setminus X$. We claim that $c$ is a centre of the automorphism $A$. Given any point $y\in Y$, we need to show that $A[\Aline yc]=\Aline yc$. First we show that $c$ is a fixed point of the automorphism $A$. It follows from $a'=A(a)\ne a$ that $a''\defeq A(a')\ne a'$. It follows from $\Aline a{a'}\parallel \Aline{a'}{a''}$ that $\Aline a{a'}=\Aline {a'}{a''}$ and hence $A[\Aline a{a'}]=\Aline {a'}{a''}=\Aline a{a'}$ and $A[\Aline b{b'}]=\Aline b{b'}$. Then $A(c)\in A[\Aline a{a'}\cap \Aline b{b'}]=A[\Aline a{a'}]\cap A[\Aline b{b'}]=\Aline a{a'}\cap\Aline b{b'}=\{c\}$ and hence $c\in\Fix(A)$. 

If $y=c$, then $A[\Aline yc]=A[\{c\}]=\{c\}=\Aline yc$. So, assume that $y\ne c$. If $y\in \Fix(A)$, then $y,c\in\Fix(A)$ implies $A[\Aline yc]=\Aline yc$ and we are done. If $y\in X\setminus\Fix(A)$, then for the point $y'\defeq A(y)$ the parallelity relation  $\Aline y{y'}\parallel\Aline a{a'}$ and the projectivity of the liner $Y$ implies  $\varnothing\ne (\Aline y{y'}\cap\Aline a{a'})\setminus X\subseteq \Aline a{a'}\setminus X=\{c\}$ and hence $c\in \Aline y{y'}$ and $A[\Aline yc]=\Aline {y'}c=\Aline yc$. So, assume that $y\in Y\setminus(X\cup\Fix(A))\subseteq Y\setminus X$. Since $\overline{Y\setminus X}\ne Y$, there exists a point $o\in Y\setminus\overline{Y\setminus X}\subseteq X$. Choose any point $z\in \Aline oy\setminus\{o,y\}$. By the strong regularity of the projective liner $Y$, the plane $\Pi\defeq\overline{\{o,z,c\}}\subseteq Y$ is equal to $\bigcup_{x\in \Aline oz}\Aline xc$. As we have already proved, $A[\Aline xc]=\Aline xc$ for every $x\in X$. In particular, the points $o'\defeq A(o)$ and $z'\defeq A(z)$ belong to the set $A[\Aline oc\cup\Aline zc]=\Aline oc\cup\Aline zc\subseteq \Pi$. Then $A[\Pi]=A[\overline{\{o,z,c\}}]=\overline{\{o',z',c\}}\subseteq\Pi$ and hence $A[\Pi]=\Pi$, by the $3$-rankedness of the projective liner $Y$. Observe that $\Aline yc=\Pi\setminus\bigcup_{x\in \overline{o\,y}\setminus\{y\}}(\Aline xc\setminus\{c\})$, which implies 
$$A[\Aline yc]=A\big[\Pi\setminus\bigcup_{x\in\overline{o\,y}\setminus\{y\}}(\Aline xc\setminus\{c\})\big]=A[\Pi]\setminus\bigcup_{x\in\overline{o\,y}\setminus\{y\}}A[\Aline xc\setminus\{c\}]=\Pi\setminus\bigcup_{x\in\overline{o\,y}\setminus\{y\}}(\Aline xc\setminus\{c\})=\Aline yc,$$
witnessing that the point $c$ is the centre of the automorphism $A$.
\end{proof}

\section{Hyperfixed automorphisms of liners}

\begin{definition} An automorphism $A:X\to X$ of a liner $X$ is called a \index{hyperfixed automorphism}\index{automorphism!hyperfixed}\defterm{hyperfixed} if its set of fixed points $$\Fix(A)\defeq\{x\in X:A(x)=x\}$$ contains a hyperplane of $X$.
\end{definition}

\begin{example} The identity automorphism of any liner is  hyperfixed.
\end{example}

\begin{example}\label{ex:hyperfixed} An automorphism $A:X\to X$ of a liner $X$ of rank $\|X\|\le 2$ is hyperfixed if and only if $|\Fix(A)|\ge\|X\|-1$.
\end{example}

In the following proposition, by a \defterm{space} we understand any $3$-long regular liner of rank $\ge 3$.

\begin{proposition}\label{p:hyperfixed=>!H} If a non-identity automorphism $A:X\to X$ of a proaffine space $X$ is hyperfixed, then its set of fixed point $\Fix(A)$ contains a unique hyperplane $H$.\newline Moreover, $|\Fix(A)\setminus H|\le 1$.
\end{proposition}

\begin{proof} Since $A$ is hyperfixed, there exists a hyperplane $H\subseteq\Fix(A)$. 

\begin{claim}\label{cl:oaFix} For any distinct points $a,b\in \Fix(A)\setminus H$, we have $\bigcup_{o\in H\setminus\overline{a\,b}}\Aline oa\subseteq \Fix(A)$.
\end{claim}

\begin{proof}  Given any point $o\in H\setminus\Aline ab$, we need to prove that $\Aline oa\subseteq \Fix(A)$. It follows from $\{o,a\}\subseteq \Fix(A)$ that $A[\Aline oa]=\Aline oa$. Since $b\in X=\overline{H\cup \{a\}}$, we can apply the regularity of $X$ and find points $u\in \Aline oa$ and $v\in H$ such that $b\in \Aline uv$.
Since $X$ is proaffine, the set $I=\{x\in \Aline ou:\Aline xb\cap \Aline ov=\varnothing\}$ has cardinality $|I|\le 1$. By the definition of the set $I$, for every point $x\in \Aline ou\setminus I=\Aline oa\setminus I$, there exists a point $y\in \Aline ov\subseteq H$ such that $b\in \Aline xy$. Then $x\in \Aline by$ and $A(x)\in A[\Aline oa\cap\Aline by]=A[\Aline oa]\cap A[\Aline by]=\Aline oa\cap\Aline by=\{x\}$. Therefore, $\Aline oa\setminus I\subseteq \Fix(A)$. Assuming that $\Aline oa\not\subseteq \Fix(A)$, we can find a point $x\in \Aline oa\setminus\Fix(A)$ and conclude that $\{x\}=I$ and $A(x)\in \Aline oa\setminus\{x\}=\Aline oa\setminus I\subseteq \Fix(A)$ and $A(A(x))=A(x)$. Since the function $A$ is bijective, the equality $A(A(x))=A(x)$ implies $A(x)=x$, which contradicts the choice of $x$. This contradiction shows that $\Aline oa\subseteq \Fix(A)$.
\end{proof}

Now we are able to prove that $|\Fix(A)\setminus H|\le 1$. In the opposite case, we can find two distinct points $a,b\in \Fix(A)\setminus H$. Since the automorphism $A$ is not identity, $A(x)\ne x$ for some point $x\in X$. Claim~\ref{cl:oaFix} ensures that $x\notin \bigcup_{o\in H\setminus\overline{a\,b}}(\Aline oa\cup\Aline ob)$. Since $\|X\|=\|H\|-1\ge 3$, there exists a point $o\in H\setminus\Aline ab$. Then $x\notin \Aline oa\cup\Aline ob$. Since $a\in X=\overline{H\cup\{x\}}$, by the regularity of $X$, there exist points $u\in \Aline ox$ and  $v\in H$ such that $a\in \Aline uv$.  It follows from $a\notin \Aline ox=\Aline ou$ that $v\ne o$. Since $u,a\in\Fix(A)$ are two distinct points with $o\notin ua$, Claim~\ref{cl:oaFix} ensures that $x\in \Aline ou\subseteq \Fix(A)$, which contradicts the choice of $x$. This contradiction shows that $|\Fix(A)\setminus H|\le 1$ for every hyperplane $H\subseteq \Fix$. 

Assuming that $H'$ is another hyperplane in $\Fix(A)$, we conclude that $|H'\setminus H|\le|\Fix(A)\setminus H|\le 1$ and hence $H'\subseteq H$ because the liner $X$ is $3$-long. Assuming that $H'\ne H$, we can take any point $x\in H\setminus H'$ and conclude that $X=\overline{H'\cup\{x\}}\subseteq H\ne X$, which is a contradiction showing that $H=H'$. Therefore, the set $\Fix(H)$ contains a unique hyperplane of $X$.
\end{proof}

\begin{proposition}\label{p:hyperfixed=>central} Let $A:X\to X$ be a hyperfixed automorphism of a proaffine regular liner $X$. If $\Fix(A)$ is not a hyperplane in $X$, then the automorphism $A$ is central.
\end{proposition}

\begin{proof} Since $A$ is hyperfixed, there exists a hyperplane $H\subseteq X$ such that $H\subseteq \Fix(A)$. If $\Fix(A)$ is not a hyperplane in $X$, then $H\ne\Fix(A)$ and hence $\Fix(A)\setminus H$ contains some point $c$. It follows that $A[\Aline cc]=A[\{c\}]=\{c\}=\Aline cc$. Next, we show that $A[\Aline xc]=\Aline xc$ for every $x\in X\setminus \{c\}$. If $\|X\|\le 2$, then $A[\Aline xc]=A[X]=X=\Aline xc$. So, assume that $\|X\|\ge 3$. If the intersection $H\cap\Aline xc$ contains some point $o$, then $\{o,c\}\subseteq \Fix(A)$ and $A[\Aline xc]=A[\Aline oc]=\Aline oc=\Aline xc$.

So, assume that $\Aline xc\cap H=\varnothing$. Choose any point $o\in H$. Since $x\in X=\overline{H\cup\{c\}}$, by the regularity of $X$, there exist points $u\in \Aline oa$ and $v\in H$ such that $x\in \Aline uv$. By Theorem~\ref{t:Proclus<=>}, the proaffine regular liner $X$ is Proclus. Then $\Aline xc$ is a unique line in the plane $\Pi\defeq\overline{\{v,o,c\}}$, which contains the point $c$ and is disjoint with the line $\Aline ov$. Then for every $z\in \Pi\setminus\Aline xc$, there exists a point $y\in \Aline zc\cap\Aline ov\subseteq H$. Observe that $A[\Aline yc]=\Aline yc$ and 
$$A[\Aline xc]=A\big[\Pi\setminus\bigcup_{y\in\Pi\setminus\overline{x\,c}}(\Aline yc\setminus\{c\})\big]=\Pi\setminus\bigcup_{y\in\Pi\setminus\overline{x\,c}}\Pi[\Aline yc\setminus\{c\}]=\Pi\setminus\bigcup_{y\in\Pi\setminus\overline{x\,c}}(\Aline yc\setminus\{c\})=\Aline xc,$$
witnessing that $c$ is a center of the automorphism $A$.
\end{proof} 

\begin{theorem}\label{t:hyperfixed=>(para)central} Every hyperfixed automorphism of a proaffine space is central or paracentral.
\end{theorem}

\begin{proof} Let $A:X\to X$ be a hypefixed automorphism of a proaffine space $X$. Assuming that $A$ is not central, we shall prove that $A$ is paracentral.  Proposition~\ref{p:hyperfixed=>central} ensures that the set $H\defeq \Fix(A)$ is a hyperplane in $X$.

\begin{claim}\label{cl:Axx'=xx'} For every point $x\in X$ and its image $x'\defeq A(x)$ we have $A[\Aline x{x'}]=\Aline x{x'}$.
\end{claim}

\begin{proof} If $x\in \Fix(A)$, then $x'=A(x)=x$ and $A[\Aline x{x'}]=A[\{x\}]=\{A(x)\}=\{x\}=\Aline x{x'}$. So, assume that $x\notin\Fix(A)$ and hence $\Aline x{x'}$ is a line in $X$. If the line $\Aline x{x'}$ has a common point $o$ with the hyperplane  $H=\Fix(A)$, then $\Aline x{x'}=\Aline xo=\Aline {x'}o$ and $A[\Aline x{x'}]=A[\Aline xo]=\Aline {x'}o=\Aline x{x'}$. So, assume that $\Aline x{x'}\cap H=\varnothing$. Choose any point $o\in H$ and consider the plane $\Pi\defeq\overline{\{x,x',o\}}$. Since $x'\in X=\overline{H\cup\{x\}}$, we can apply  the regularity of $X$ and find points $u\in \Aline ox$ and $v\in H$ such that $x'\in \Aline uv$. Assuming that $v=o$, we conclude that $x'\in \Aline uv=\Aline uo\subseteq \Aline xo$ and hence $o\in \Aline x{x'}\cap H=\varnothing$, which is a contradiction showing that $v\ne o$. Then $v\in (H\cap \Aline u{x'})\setminus\{o\}\subseteq (H\cap\Pi)\setminus\{o\}$. By the rankedness of the regular liner $X$, $\overline{\{o,v,x\}}=\Pi=\overline{\{o,v,x'\}}$ and $A[\Pi]=A[\overline{\{o,v,x\}}]=\overline{\{o,v,x'\}}=\Pi$. By the Proclus property of the proaffine space $X$, $\Aline x{x'}$ is a unique line in the plane $\Pi$ that contains the point $x'$ and is disjoint with the line $\Aline ov$. Asuming that $A[\Aline x{x'}]\ne \Aline x{x'}$, we conclude that the intersection $A[\Aline x{x'}]\cap \Aline ov$ contains some point $w$. Taking into account that $w\in\Aline ov\subseteq \Fix(A)$, we conclude that $A(w)=w\in A[\Aline x{x'}]$ and hence $w\in\Aline x{x'}\cap H$, which contradicts our assumption. This contradiction shows that $A[\Aline x{x'}]=\Aline x{x'}$.
\end{proof}

\begin{claim}\label{cl:xx'yy'-colinear} Let $x,y\in X\setminus H$ be any points and $x'\defeq A(x)$, $y'\defeq A(y)$ be their images. If $\Aline x{y}\cap H\ne\varnothing$, then the lines $\Aline x{x'}$ and $\Aline {y}{y'}$ are coplanar.
\end{claim}

\begin{proof} Since $\Aline xy\cap H\ne\varnothing$, there exists a unique point $o\in \Aline xy\cap H$. If $x'\in \Aline xo$, then $A[\Aline x{x'}]=A[\Aline xo]=\Aline {x'}{o}=\Aline xo$ and $y'=A(y)\in A[\Aline xo]=\Aline xo$ and finally $\Aline y{y'}=\Aline xo=\Aline x{x'}$. Therefore, the lines $\Aline x{x'}=\Aline y{y'}$ are coplanar.

So, assume that $x'\notin \Aline xo$. In this case $\Pi\defeq\overline{\{o,x,x'\}}$ is a plane in $X$. By the regularity of the liner $X=\overline{H\cup\{x\}}$, there exist points $u\in \Aline ox$ and $v\in H$ such that $x'\in \Aline uv$ and hence $v\in \Aline u{x'}\subseteq \Pi$. Assuming that $v=o$, we conclude that $x'\in \Aline uv=\Aline uo\subseteq \Aline xo$, which contradicts our assumption. Therefore, $v\ne o$ and hence $\overline{\{o,v,x\}}=\Pi=\overline{\{o,v,x'\}}$, by the rankedness of the regular liner $X$. Then $A[\Pi]=A[\overline{\{o,v,x\}}]=\overline{\{o,v,x'\}}=\Pi$ and $\Aline x{x'}\cup \Aline y{y'}\subseteq \Pi$, which means that the lines $\Aline x{x'}$ and $\Aline y{y'}$ are coplanar.
\end{proof}

\begin{claim}\label{cl:xx'|yy'} Let $x,y\in X\setminus H$ be any points and $x'\defeq A(x)$, $y'\defeq A(y)$ be their images. If $\Aline x{x'}\cap H=\varnothing \ne\Aline x{y}\cap H$, then the lines $\Aline x{x'}$ and $\Aline {y}{y'}$ are parallel.
\end{claim}

\begin{proof} By Claim~\ref{cl:xx'yy'-colinear}, the lines $\Aline x{x'}$ and $\Aline y{y'}$ are coplanar. Assuming that these lines are not parallel, we can find a unique point $c\in \Aline x{x'}\cap\Aline y{y'}$. Claim~\ref{cl:Axx'=xx'} ensures that $\{A(c)\}=A[\Aline x{x'}\cap\Aline y{y'}]=A[\Aline x{x'}]\cap A[\Aline y{y'}]=\Aline x{x'}\cap\Aline y{y'}=\{c\}$ and hence $c\in \Fix(A)=H$ and $c\in \Aline x{x'}\cap H=\varnothing$, which is a contradiction completing the proof of the claim.
\end{proof}

\begin{claim}\label{cl:xx'||yy'} Let $x,y\in X\setminus H$ be any points and $x'\defeq A(x)$, $y'\defeq A(y)$ be their images. If $\Aline x{x'}\cap H=\varnothing$, then the lines $\Aline x{x'}$ and $\Aline {y}{y'}$ are parallel.
\end{claim}

\begin{proof} If $x=y$, then the lines $\Aline x{x'}$ and $\Aline y{y'}$ coincide and hence are parallel. So, assume that $x\ne y$. If $\Aline xy\cap H\ne\varnothing$, then Claim~\ref{cl:xx'|yy'} implies $\Aline x{x'}\parallel \Aline y{y'}$. So, assume that $\Aline xy\cap H=\varnothing$.  Choose any point $o\in H$ and use the regularity of the liner $X=\overline{H\cup\{x\}}$ to find points $u\in\Aline ox$ and $v\in H$ such that $y\in\Aline uv$.  Assuming that $u\in H$, we conclude that $y\in \Aline uv\subseteq H$, which contradicts the choice of $y$. Therefore, $u\in X\setminus H$ and $u'\defeq A(u)\ne u$.  Since $x\ne y$ and $\Aline xy\cap H=\varnothing$, the point $u$ is not equal to the point $x$ and hence $o\in \Aline xu\cap H\ne\varnothing$. Since $\Aline xu\cap H\ne \varnothing\ne\Aline uy\cap H$, we can apply Clacm~\ref{cl:xx'|yy'} and conclude that $\Aline x{x'}\parallel \Aline u{u'}\parallel \Aline y{y'}$, and hence $\Aline x{x'}\parallel \Aline y{y'}$, by Theorem~\ref{t:Proclus-lines}.
\end{proof}

\begin{claim}\label{cl:xx'yy'-coplanar2} Let $x,y\in X\setminus H$ be any points and $x'\defeq A(x)$, $y'\defeq A(y)$ be their images. If $\Aline x{x'}\cap H\ne\varnothing$, then the lines $\Aline x{x'}$ and $\Aline {y}{y'}$ are coplanar.
\end{claim}

\begin{proof} By our assumption, there exists a point $o\in \Aline x{x'}\cap H$. If $y\in \Aline xo$, then $y'=A(y)\in A[\Aline xo]=\Aline {x'}{o}=\Aline xo$ and $\Aline y{y'}=\Aline xo=\Aline x{x'}$. Therefore, the coinciding lines $\Aline x{x'}=\Aline y{y'}$ are coplanar.

So, assume that $y\notin \Aline xo$. In this case $\Pi\defeq\overline{\{o,x,y\}}$ is a plane in $X$. By the regularity of $X=\overline{H\cup\{x\}}$, there exist points $u\in \Aline ox$ and $v\in H$ such that $y\in \Aline uv$ and hence $v\in \Aline uy\subseteq \Pi$. Assuming that $v=o$, we conclude that $y\in \Aline uv=\Aline uo\subseteq \Aline xo$, which contradicts our assumption. Therefore, $v\ne o$ and hence $\overline{\{o,v,x\}}=\Pi=\overline{\{o,v,y\}}$, by the rankedness of the regular liner $X$. Then $y'=A(y)\in A[\Pi]=A[\overline{\{o,v,x\}}]=\overline{\{o,v,x'\}}=\Pi$ and $\Aline x{x'}\cup \Aline y{y'}\subseteq \Pi$, which means that the lines $\Aline x{x'}$ and $\Aline y{y'}$ are coplanar.
\end{proof}

Now we are able to complete the proof of the theorem. Given any points $x,y\in X\setminus H$ and their images $x'\defeq A(x)$ and $y'\defeq A(y)$, we should prove that the lines $\Aline x{x'}$ and $\Aline y{y'}$ are parallel. If for some point $z\in X\setminus H$ and its image $z'\defeq A(z)$ the line $\Aline z{z'}$ is disjoint with the hyperplane $H$, then by Claim~\ref{cl:xx'||yy'}, the lines $\Aline x{x'},\Aline z{z'},\Aline y{y'}$ are parallel and we are done. So, assume that $\Aline z{z'}\cap H\ne\varnothing$ for every $z\in X\setminus H$. In this case Claim~\ref{cl:xx'yy'-coplanar2} ensures that the lines $\Aline x{x'}$ and $\Aline y{y'}$ are coplanar. Assuming that these two coplanar lines are not parallel, we conclude that $\Aline x{x'}\cap\Aline y{y'}=\{c\}$ for some point $c\in X$. Claim~\ref{cl:Axx'=xx'} ensures that $\{A(c)\}=A[\Aline x{x'}\cap\Aline y{y'}]=A[\Aline x{x'}]\cap A[\Aline y{y'}]=\Aline x{x'}\cap \Aline y{y'}=\{c\}$ and hence $c\in \Fix(A)=H$. Since the automorphism $A$ is not central, there exists a point $z\in X$ such that $A[\Aline zc]\ne\Aline zc$. By Claim~\ref{cl:xx'yy'-coplanar2}, the lines $\Aline x{x'}$ and $\Aline z{z'}$ are coplanar. Assuming that the coplanar lines $\Aline x{x'}$ and $\Aline z{z'}$ are not parallel, we can find a unique point $o\in \Aline x{x'}\cap \Aline z{z'}$ and conclude that $\{A(o)\}=A[\Aline x{x'}]\cap A[\Aline y{y'}]=\Aline x{x'}\cap \Aline y{y'}=\{o\}$ and hence $o\in \Fix(A)=H$. Then $o\in \Aline x{x'}\cap H=\{c\}$ and $\Aline z{z'}=\Aline zo=\Aline zc$. Applying Claim~\ref{cl:Axx'=xx'}, we conclude $A[\Aline zc]=A[\Aline z{z'}]=\Aline z{z'}=\Aline zo=\Aline zc$, which contradicts the choice of the point $z$. This contradiction shows that the lines $\Aline z{z'}$ and $\Aline x{x'}$ are parallel. By analogy we can prove that the lines $\Aline z{z'}$ and $\Aline y{y'}$ are parallel. Then $\Aline x{x'}\parallel \Aline z{z'}\parallel \Aline y{y'}$ and $\Aline x{x'}\parallel \Aline y{y'}$, by Theorem~\ref{t:Proclus-lines}.
\end{proof}

\section{Central and hyperfixed automorphisms of projective spaces}

\begin{theorem}\label{t:central<=>hyperfixed} An automorphism $A:X\to X$ of a projective space $X$ is central if and only if $A$ is hyperfixed.
\end{theorem}

\begin{proof} To prove the ``if'' part, assume that the automorphism $A:X\to X$ is hyperfixed and hence the set $\Fix(A)\defeq\{x\in X:A(x)=x\}$ contains some hyperplane $H$ in $X$. If $H\ne\Fix(A)$, then the automorphism $A$ is central, by Proposition~\ref{p:hyperfixed=>central}. So, assume that $H=\Fix(A)$. Choose any point $x\in X\setminus \Fix(A)$ and consider its image $x'\defeq A(x)$. By Proposition~\ref{p:cov-aff},  the $3$-long liner $X$ contains a point $y\notin H\cup\Aline x{x'}$. Assuming that the hyperfixed automorphism $A$ is not central, we can apply Theorem~\ref{t:hyperfixed=>(para)central} and conclude that $\Aline x{x'}$ and $\Aline y{y'}$ are distinct parallel lines in the projective liner $X$, which is impossible (by Corollary~\ref{c:parallel-in-modular} and Theorem~\ref{t:modular<=>}). This contradiction shows that the hyperfixed automorphism $A$ of the projective space $X$ is central.
\smallskip

To prove the ``only if' part, assume that the automorphism $A$ is central. If $A$ is the identity automorphism of $X$, then $\Fix(A)=X$ and $A$ is hyperfixed. So, assume that $A$ is not identity. By Proposition~\ref{p:cenralauto-has1-center}, the central automorphism $A$ has a unique center $c$. 

\begin{claim}\label{cl:Fix-meets-every-line} The set $\Fix(A)$ has nonempty intersection with every line $L\subseteq X$.
\end{claim}

\begin{proof} To derive a contradiction, assume that some line $L$ is disjoint with the set $\Fix(A)\defeq \{x\in X:A(x)=x\}$. Given any distinct points $x,y\in L$, consider the points $x'\defeq A(x)$ and $y'\defeq A(y)$. Since $c$ is a center of $A$, $x'=A(x)\in A[\Aline xc]=\Aline xc$ and $y'=A(y)\in A[\Aline yc]=\Aline yc$. Since $x,y\in L=L\setminus \Fix(A)$ and $c\in \Fix(A)$, the points $x',y'$ do not belong to the line $L$. Observe that $A[L]=A[\Aline xy]=\Aline {x'}{y'}\subseteq \overline{\{x,c,y\}}$ and hence $L$ and $A[L]$ are two distinct lines in the plane $\overline{\{x,y,c\}}$. Since the liner $X$ is projective, the coplanar lines $L$ and $A[L]$ have a unique common point $o\in L\cap A[L]$. Consider its image $o'\defeq A(o)\in A[L]$ and observe that $L\cap\Fix(A)=\varnothing$ implies $o\ne o'$ and hence $c\in \Aline o{o'}=A[L]$ and $c=A(c)\in A[L]$ implies $c\in L\cap\Fix(A)$, which contradicts our assumption.
\end{proof}

\begin{claim}\label{cl:xy-subset-FixA} For every points $x,y\in \Fix(A)$, if $c\notin \Aline xy$, then $\Aline xy\subseteq\Fix(A)$.
\end{claim}

\begin{proof} Take any point $z\in \Aline xy$ and observe that $\{z\}=\Aline xy\cap\Aline zc$, which implies $\{A(z)\}=A[\Aline xy]\cap A[\Aline zc]=\Aline xy\cap\Aline zc=\{z\}$ and hence $z\in \Fix(A)$.
\end{proof}

\begin{claim}\label{cl:hyperplane1} If $\Aline xc\subseteq \Fix(A)$ for every $x\in \Fix(A)$, then the set $\Fix(A)$ is a hyperplane in $X$.
\end{claim}

\begin{proof} First we prove that the set $\Fix(A)$ is flat in $X$. Given any distinct points $x,y\in\Fix(A)$, we should prove that $\Aline xy\subseteq \Fix(A)$. If $c\notin \Aline xy$, then  $\Aline xy\subseteq\Fix(A)$, by Claim~\ref{cl:xy-subset-FixA}. If $c\in \Aline xy$, then $\Aline xy=\Aline xc\subseteq\Fix(A)$, by the assumption. In both cases, $\Aline xy\subseteq \Fix(A)$, witnesing that the set $\Fix(A)$ is flat. 
 
Assuming that the flat $\Fix(A)$ is not a hyperplane, we can find a point $x\in X\setminus \Fix(A)$ such that $\overline{\Fix(A)\cup\{x\}}\ne X$. Then for any point $y\in X\setminus\overline{\Fix(A)\cup\{x\}}$, we have $\Aline xy\cap \Fix(A)=\Aline xy\cap\Fix(A)\cap\overline{\Fix(A)\cup\{x\}}=\Fix(A)\cap \{x\}=\varnothing$, which contradicts Claim~\ref{cl:Fix-meets-every-line}.
\end{proof}

\begin{claim}\label{cl:hyperplane2} If $\Aline ac\not\subseteq \Fix(A)$ for some $a\in\Fix(A)$, then the set $\Fix(A)\setminus\{c\}$ is a hyperplane in $X$.
\end{claim}

\begin{proof} Using the Kuratowski--Zorn Lemma, choose a maximal flat subset $M\subseteq \Fix(A)$ of $X$ such that $a\in M$. It follows from $\Aline ac\not\subseteq\Fix(A)$ that $c\notin M$. Choose any point $b\in\Aline ac\setminus\Fix(A)$. Assuming that $\overline{M\cup\{c\}}\ne X$, we can find a point $d\in X\setminus \overline{M\cup\{c\}}$. By Claim~\ref{cl:Fix-meets-every-line}, there exists a point $e\in \Aline bd\cap \Fix(A)$. We claim that $c\notin \Aline xe$ for all $x\in M$. In the opposite case, $e\in \Aline xc\subseteq \overline{M\cup\{c\}}$ and hence $e\in \Aline bd\cap\overline{M\cup\{c\}}=\{b\}\not\subseteq \Fix(A)$, which contradicts the choice of the point $e$. This contradiction shows that $c\notin \Aline xe$ for all $x\in M$. Applying Claim~\ref{cl:xy-subset-FixA}, we conclude that $\bigcup_{x\in M}\Aline xe\subseteq\Fix(A)$. The strong regularity of the projective space $X$ ensures that the set $\bigcup_{x\in M}\Aline xe\subseteq \Fix(A)$ is flat, which contradicts the maximality of the flat $M$.  This contradiction shows that $\overline {M\cup\{c\}}=X$ and hence $\|X\|_M=1$ and $M$ is a hyperplane in $X$, by Proposition~\ref{p:hyperplane}. Proposition~\ref{p:hyperfixed=>!H} ensures that $\Fix(A)\setminus M=\{c\}$ and hence $\Fix(A)\setminus\{c\}=M$ is a hyperplane in $X$.
\end{proof}  

Claims~\ref{cl:hyperplane1} and \ref{cl:hyperplane2} imply that the central automorphism $A$ is hyperfixed.
\end{proof}

\begin{corollary} An automorphism $A:X\to X$ of a (non-affine) completely regular space $X$ is hyperfixed (if and) only if $A$ is central or paracentral.
\end{corollary}

\begin{proof} By Theorem~\ref{t:spread=projective1}, the completely regular space $X$ is para-Playfair and hence proaffine. If an automorphism $A:X\to X$ is hyperfixed, then $A$ is central or paracentral, by Theorem~\ref{t:hyperfixed=>(para)central}. Now assume that an automorphism $A:X\to X$ is central or paracentral. 

Since the space $X$ is completely regular, its spread completion $\overline X$ is a projective space. By Theorem~\ref{t:extend-isomorphism-to-completions}, the automorphism $A:X\to X$ entends to an automorphism $\bar A:\overline X\to\overline X$ of the spread completion $\overline X$ of $X$. By Proposition~\ref{p:central|=(para)central}, the automorphism $\bar A$ of the projective space $\overline X$ is central, and by Theorem~\ref{t:central<=>hyperfixed}, the central automorphism $\bar A:\overline X\to\overline X$ of the projective space $\overline X$ is hyperfixed. Then there exists a hyperplane $H$ in $\overline X$ such that $H\subseteq \Fix(\bar A)$. Taking into account that the space $X$ is not affine, we can apply Theorem~\ref{t:affine<=>hyperplane} and conclude that the flat $\partial X$ is not a hyperplane in $\overline X$. Then $H\not\subseteq\partial X$ and hence $H\cap X\ne\varnothing$. Lemma~\ref{l:trace-flat} implies that $H\cap X\subseteq\Fix(A)$ is a hyperplane in $X$, witnessing that the automorphism $A$ of $X$ is hyperfixed.
\end{proof}

\section{Hyperfixed and paracentral automorphisms of affine spaces}

\begin{proposition}\label{p:hyperfixed=>paracentral} Every hyperfixed automorphism $A:X\to X$ of an affine space $X$ is  paracentral and hence $\Fix(A)$ is a flat in $X$.
\end{proposition}

\begin{proof} Assuming that the hyperfixed automorphism $A:X\to X$ is not paracentral, we can apply Theorem~\ref{t:hyperfixed=>(para)central} and conclude that $A$ is central and not identity. By Theorem~\ref{t:central<=>homothety}, $A$ is a homothety and by Theorem~\ref{t:dilation-atmost1}, the set $\Fix(A)$ is a singleton and hence $A$ cannot be hyperfixed. This contradiction shows the hyperfixed automorphism $A$ of the affine space is paracentral.

Next, we show that the set $\Fix(A)$ is flat in $X$. If the hyperfixed automorphism $A$ is identity, then $\Fix(A)=X$ is flat. If $A$ is not identity, then the paracentral property of $A$ implies that $A$ is not central and hence $\Fix(A)$ is a hyperplane in $X$, by Proposition~\ref{p:hyperfixed=>central}. 
\end{proof}

\begin{definition}\label{d:hypershear-hyperscale} A hyperfixed automorphism $A:X\to X$ of an affine space $X$ is called a 
\begin{itemize}
\item a \index{hypershear}\defterm{hypershear} if for every point $x\in X\setminus\Fix(A)$ and its image $y\defeq A(x)$, the line $\Aline xy$ is subparallel to the flat $\Fix(A)$;
\item a \index{hyperscale}\defterm{hyperscale} if for every point $x\in X\setminus\Fix(A)$ and its image $y\defeq A(x)$, the line $\Aline xy$ is not subparallel to the flat $\Fix(A)$.
\end{itemize}
\end{definition}

Definition~\ref{d:hypershear-hyperscale} implies that the identity map of an affine space is both a hypershear and a hyperscale. Also Definition~\ref{d:hypershear-hyperscale} and Proposition~\ref{p:hyperfixed=>paracentral} imply the following (trivial) classification of hyperfixed automorphisms of affine spaces. 

\begin{proposition}\label{p:hyperfixed<=>shear-or-scale} An automorphism $A:X\to X$ of an affine space $X$ is hyperfixed if and only if $A$ is a hypershear or a hyperscale.
\end{proposition}

Finally, we present a classification of paracentral automorphisms of affine spaces.

\begin{proposition}\label{p:aff-paracentral<=>} An automorphism $A:X\to X$ of an affine space  $X$ is paracentral if and only if $A$ is a translation, a hypershear or a hyperscale.
\end{proposition}

\begin{proof} The ``if'' part follows from Propositions~\ref{p:Trans-spread} and  \ref{p:hyperfixed=>paracentral}. To prove the ``only if'' part, assume that an automorphism $A:X\to X$ of an affine space $X$ is paracentral. By Corollary~\ref{c:affine-spread-completion}, the affine space $X$ is completely regular and hence its spread completion $\overline X$ is a projective space, which is a projective completion of $X$. By Theorem~\ref{t:extend-isomorphism-to-completions}, the automorphism $A$ extends to an automorphism $\bar A:\overline X\to\overline X$ of the projective space $\overline X$. By Proposition~\ref{p:central|=(para)central}, the automorphism $\bar A$ of the projective space $\overline X$ is central. By Theorem~\ref{t:central<=>hyperfixed}, the central automorphism $\bar A$ is hyperfixed.  If $\bar A$ is the identity automorphism of $\overline X$, then $A$ is the identity automorphism of $X$ and hence $A$ is simultaneously a translation and a hypershear and a hyperscale. So, assume that the central automorphism $\bar A$ is not idenity. Then $\Fix(A)=H\cup\{c\}$ for a unique hyperplane $H\subseteq\overline X$ and a unique center $c$ of $\bar A$. By Theorem~\ref{t:affine<=>hyperplane}, the boundary $\partial X\defeq \overline X\setminus X$ is a hyperplane in $\overline X$. If $H=\partial X$, then the paracentral property of $A$ and the central property of $\bar A$ imply that $c\in\partial A$ and $A$ is a translation. So, assume that $H\ne\partial X$. In this case Lemma~\ref{l:trace-flat} implies that $H\cap X$ is a hyperplane in $X$ and hence the automorphism $A:X\to X$ is hyperfixed. By Proposition~\ref{p:hyperfixed<=>shear-or-scale}, the hyperfixed automorphism $A$ of $X$ is a hypershear or a hyperscale.
\end{proof}

\section{Shear liners}

\begin{definition} A liner $X$ is called a \index{shear liner}\index{liner!shear}\defterm{shear liner} if for every hyperplane $H\subseteq X$ and distinct points $x,y\in X\setminus H$ with $\Aline xy\cap H=\varnothing$, there exists an automorphism $A:X\to X$ of $X$ such that $A(x)=y$ and $\Fix(A)=H$.
\end{definition}

\begin{theorem}\label{t:shear=>translation} Every shear Playfair plane is translation.
\end{theorem}

\begin{proof} We shall deduce this theorem from the following lemma, which will be also applied in the proof of Theorem~\ref{t:vtrans+vshear=vshears}.

\begin{lemma}\label{l:hypershear2=translation} Let $A,B:X\to X$ be two hypershears of a Playfair plane $X$ such that $\Fix(A)\cap \Fix(B)=\varnothing$. If there exist points $a\in \Fix(B)$ and $b\in\Fix(A)$ such that $\Aline ab\parallel \Aline{A(a)}{B(b)}$, then the automorphism $AB=BA$ is a translation of $X$.
\end{lemma}

\begin{proof}   Proposition~\ref{p:hyperfixed=>paracentral} ensures that the hypershears $A,B$ are paracentral and moreover the sets $\Fix(A)$ and $\Fix(B)$ are flats in $X$. Since $\Fix(A)\cap\Fix(B)=\varnothing$, the flats $\Fix(A)$ and $\Fix(B)$ are parallel lines in the Playfair plane $X$. Since $a\in\Fix(B)\subseteq X\setminus\Fix(A)$, the point $a'=A(a)$ is not equal to the point $a\in\Fix(B)$ and hence $\Aline a{a'}$ is a line parallel to the line $\Fix(A)\parallel \Fix(B)$. Then $a'\in\Aline a{a'}=\Fix(B)$. By analogy we can show that the point $b'\defeq B(b)$ is not equal to $b$ and $\Aline b{b'}=\Fix(A)$. By the assumption, $\Aline ab\parallel \Aline {a'}{b'}$.

Since the automorphisms $A,B$ are paracentral with $\Fix(A)\parallel\Fix(B)$,  for every point $x\in X$ and its images $y\defeq A(x)$ and $z\defeq B(x)$ we have $\Aline xy\subparallel \Aline a{a'}=\Fix(B)\parallel \Fix(A)$ and $\Aline xz\subparallel \Aline b{b'}=\Fix(A)\parallel \Fix(B)$. This implies that the automorphisms $AB$ and $BA$ are paracentral.

Assuming that $AB$ is not a translation, we can apply Proposition~\ref{p:aff-paracentral<=>} and conclude that the paracentral automorphism $AB$ is hyperfixed and not identity. By Proposition~\ref{p:hyperfixed=>paracentral},  the set $\Fix(AB)$ is a line in the Playfair plane $X$. Assuming that the line $\Aline a{a'}=\Fix(B)$ has a common point $o$ with the line $\Fix(AB)$, we conclude that $B(o)=o=AB(o)$ and hence $A(o)=AB(o)=o$ and $o\in\Fix(A)\cap\Fix(B)=\varnothing$. This contradiction shows that $\Aline a{a'}\cap \Fix(AB)=\varnothing$. Assuming that the line $\Aline b{b'}=\Fix(A)$ has a common point $o$ with the line $\Fix(AB)$, we conclude that $A(o)=o=AB(o)$ and hence $B(o)=o$ and $o\in \Fix(A)\cap\Fix(B)=\varnothing$. This contradiction shows that $\Aline b{b'}\cap\Fix(AB)=\varnothing$. Therefore, the lines $\Aline a{a'},\Aline b{b'}$ and $\Fix(AB)$ are parallel. Since the automorphisms $A,B$ are paracentral, $A[\Fix(AB)]=\Fix(AB)=B[\Fix(AB)]$.

By the Proclus Axiom, the lines $\Aline a{b'}$ and $\Fix(AB)$ have a unique common point $o$. Let $o'\defeq A(o)$. Since the automorphism $A$ is paracentral, $\Aline o{o'}\parallel \Aline a{a'}\parallel\Fix(AB)$ and hence $o'\in\Fix(AB)$ and $AB(o')=o'=A(o)$ and $B(o')=o$. 
Observe that $o'=A(o)\in A[\Aline a{b'}]=\Aline {a'}{b'}$. On the other hand, $o'=B^{-1}(o)\in B^{-1}[\Aline a{b'}]=\Aline ab$ and hence $o'\in\Aline ab\cap\Aline {a'}{b'}$ and the parallel lines $\Aline ab$ and $\Aline {a'}{b'}$ coincide. Then $\Fix(B)=\Aline a{a'}=\Aline {b}{b'}=\Fix(A)$, which contradicts the assumption. This contradiction shows that the automorphism $AB$ is a translation of $X$. By analogy we can prove that the automorphism $BA$ is a translation. Since $\{a,a'\}\subseteq \Fix(B)$, $AB(a)=A(a)=a'=B(a')=BA(a)$, the translations $AB$ and $BA$ coincide by Proposition~\ref{p:Ax=Bx=>A=B}.
\end{proof}

Now we can present a proof of Theorem~\ref{t:shear=>translation}. Given two distinct points $a,b$ of a shear Playfair plane $X$, we need to construct a translation $T:X\to X$ such that $T(a)=a'$. Choose any point $b\in X\setminus\Aline a{a'}$. Since the plane $X$ is Playfair,  there exists a unique point $b'\in X$ such that $\Aline b{b'}\parallel \Aline a{a'}$ and $\Aline a{a'}\parallel \Aline b{b'}$. 
Since $X$ is a shear plane, there exist automorphisms $A,B:X\to X$ such that $A(a)=a'$, $B(b)=b'$, $\Aline a{a'}=\Fix(A)$ and $\Aline b{b'}=\Fix(B)$. 
Then $\Fix(A)\cap\Fix(B)=\Aline b{b'}\cap\Aline a{a'}=\varnothing$. By Lemma~\ref{l:hypershear2=translation}, the automorphism $T\defeq BA=BA$ is a translation of the plane $X$ with $T(a)=AB(a)=A(a)=a'$. 
\end{proof}

\begin{remark}\label{rem:translation=2shears} The proof of Theorem~\ref{t:shear=>translation} actually shows that every translation of a shear Playfaor plane is a composition of two shears.
\end{remark} 

\begin{exercise} Prove that every shear affine space is translation.
\smallskip

{\em Hint:} Apply Theorems~\ref{t:shear=>translation}, \ref{t:Desargues-affine}, \ref{t:ADA=>AMA}, and \ref{t:paraD<=>translation}.
\end{exercise}

\begin{proposition}\label{p:shear=>2-homogen} Let $X$ be a shear Playfair plane. For any pairs $xy,x'y'\in X^2\setminus 1_X$, there exists an automorphism $A:X\to X$ such that $Axy=x'y'$.
\end{proposition}

\begin{proof} By Theorem~\ref{t:shear=>translation}, the shear Playfair plane $X$ is translation and hence there exists a translation $T:X\to X$ such that $T(x)=x'$. Consider the point $z\defeq T(y)$ and choose any point
 $z'\in X\setminus(\Aline {x'}{z}\cup\Aline{x'}{y'})$. Since $X$ is Playfair, there exist lines $L,L'\subseteq X$ such that $x'\in L\cap L'$, $L\parallel {z}{z'}$ and $L'\parallel \Aline {z'}{y'}$. Since $X$ is shear, there exist hypershears $S,S'$ of $X$ such that $L\subseteq \Fix(S)$, $L'\subseteq\Fix(S')$ and $S(z)=z'$ and $S(z')=y'$. Then the automorphism $A\defeq S'ST$ of $X$ has the required properties: $A(x)=S'ST(x)=S'S(x')=S'(x')=x'$ and $A(y)=S'ST(y)=S'S(z)=S'(z)=y'$.
 \end{proof} 
 
\begin{Exercise} Show that for any triangles $xyz$ and $x'y'z'$ in a shear Playfair plane $X$, there exists an automorphism $A$ of $X$ such that $Axyz=x'y'z'$.
\smallskip

{\em Hint:} This is difficult and the known proofs exploit some deep results on the algebraic structure of alternative rings. 
\end{Exercise} 

\section{Hypersymmetric liners}

A bijection $F$ of a set $X$ is called \index{involutive bijection}\defterm{involutive} if $F(F(x))=x$ for every $x\in X$.

\begin{definition} A liner $X$ is called \index{hypersymmetric liner}\index{liner!hypersymmetric}\defterm{hypersymmetric} if for every triangle $abc$ in $X$ there exists a hyperfixed involutive automorphism $A:X\to X$ such that $Aabc=cba$.
\end{definition} 

In this definition by $Aabc$ we understand the ordered triple $(A(a),A(b),A(c))$.

The main result of this section is the following theorem whose proof exploits the ideas of \index[person]{Burn}R.P.~Burn\footnote{{\bf Robert Pemberton Burn} (1934 -- 2024), British mathematics educator, priest. Ordained deacon Church South India, 1963; ordained priest Church of England, 1981. Member Mathematics Association, Association Teachers of Mathematics, British Society History of Math.}

\begin{theorem}\label{t:hypersymmetric=>shear} Every hypersymmetric Playfair plane is shear.
\end{theorem}

\begin{proof} By Theorem~\ref{t:parallelogram3+1}, for every triangle $xyz$ in a Playfair  liner $X$ there exists a unique point $p\in X$ such that $xyzp$ is a parallelogram. We shall denote this point $p$ by $\diamond xyz$. It is clear that $\diamond xyz=\diamond zyx$ for every triangle $xyz$ in $X$.

\begin{lemma}\label{l:diamond-preserve} For every triangle $xyz$ in a Playfair liner $X$ and every automorphism $A:X\to X$ with $Axyz=zyx$, we have $A(\diamond xyz)=\diamond xyz$.
\end{lemma}

\begin{proof} Let $p\defeq \diamond xyz$ and $p'\defeq A(p)$. Since $xyzp$ is a parallelogram, $\Aline xy\parallel \Aline zp$ and $\Aline xp\parallel \Aline yz$. Since $Axyz=zyx$ and the automorphism $A$ preserves the parallelity relation, 
$$\Aline x{p'}=A[\Aline zp]\parallel A[\Aline xy]=\Aline zy\parallel \Aline xp$$ and hence $\Aline x{p'}=\Aline xp$. By analogy we can show that $\Aline z{p'}=\Aline zp$. Then $p'\in \Aline x{p'}\cap\Aline z{p'}=\Aline xp\cap\Aline zp=\{p\}$ and hence $A(\diamond xyz)=A(p)=p'=p=\diamond xyz$. 
\end{proof}

\begin{lemma}\label{l:atmost1hypersymmetry} For every triangle $abc$ in a Playfair plane $X$, there exists at most one hyperfixed automorphism $A$ of $X$ such that $Aabc=cba$.
\end{lemma}

\begin{proof} Assume that $A,B$ are two hyperfixed automorphisms of $X$ such that $Aabc=cba=Babc$. By Proposition~\ref{p:hyperfixed=>paracentral}, the sets $\Fix(A)$ and $\Fix(B)$ are lines in $X$ with $b\in\Fix(A)\cap\Fix(B)$.  Consider the point $d\defeq \diamond abc$. Lemma~\ref{l:diamond-preserve} ensures that $d\in \Fix(A)\cap\Fix(B)$ and hence $\Aline bd=\Fix(A)=\Fix(B)$. Assuming that $A\ne B$, we can find a point $x\in X$ such that $A(x)\ne B(x)$. Since $\Fix(A)=\Aline bd=\Fix(B)$, the point $x$ does not belong to the line $\Aline bd$. 

Since $X$ is a Playfair plane, there exist unique points $b',d'\in\Aline bd$ such that $\Aline {x}{b'}\parallel \Aline ab$ and $\Aline{x}{d'}\parallel \Aline ad$. By Theorem~\ref{t:parallelogram3+1}, there exists a unique point $y\in X$ such that $\Aline {b'}y\parallel \Aline bc$ and $\Aline {d'}y\parallel \Aline dc$. Consider the point $z\defeq A(x)$ and observe that $\{b,b'\}\subseteq \Aline bd=\Fix(A)$ and $\Aline {x}{b'}\parallel \Aline ab$ imply $\Aline {z}{b'}=A[\Aline x{b'}]\parallel A[\Aline ab]=\Aline cb\parallel \Aline {y}{b'}$ and hence $\Aline{z}{b'}=\Aline{y}{b'}$. By analogy we can show that $\Aline{z}{d'}=\Aline{y}{d'}$. Then $\{z\}=\Aline {z}{b'}\cap\Aline{z}{d'}=\Aline{y}{b'}\cap\Aline{y}{d'}=\{y\}$ and hence $A(x)=z=y$. By analogy we can show that $B(x)=y=A(x)$, which contradicts the choice of the point $x$ and shows that $A=B$.    
\end{proof}

\begin{lemma}\label{l:central-parallelograms} Let $abcd$ and $a'b'c'd'$ be two parallelograms in a hypersymmetric Playfair plane $X$. If $ab=a'b'$ and $\Aline ac\cap\Aline bd=\Aline {a'}{c'}\cap\Aline{b'}{d'}\ne\varnothing$, then $cd=c'd'$.
\end{lemma}

\begin{proof} Assume that $a'b'=ab$. The condition  $\Aline ac\cap\Aline bd=\Aline {a'}{c'}\cap\Aline{b'}{d'}\ne\varnothing$ implies that $\Aline ac\cap\Aline bd=\Aline {a'}{c'}\cap\Aline{b'}{d'}=\{o\}$ for some point $o\in X$. It follows from $\Aline ab\cap\Aline cd=\varnothing=\Aline ad\cap\Aline bc$ that $o\notin(\Aline ab\cup\Aline bc\cup\Aline cd\cup\Aline ad)$. By analogy we can show that $o\notin (\Aline {a'}{b'}\cup\Aline {b'}{c'}\cup\Aline {c'}{d'}\cup\Aline {a'}{d'})$.

 Consider the point $p\defeq\diamond aod$. We claim that $\Aline po\parallel \Aline ab$. 
 
\begin{picture}(100,90)(-150,-15)
\put(0,0){\line(1,0){30}}
\put(0,0){\line(1,1){30}}
\put(0,0){\line(0,1){30}}
\put(30,0){\line(0,1){60}}
\put(30,0){\line(1,1){30}}
\put(0,30){\line(1,-1){30}}
\put(0,30){\line(1,0){60}}
\put(0,30){\line(1,1){30}}
\put(60,30){\line(-1,1){30}} 

\put(0,0){\circle*{3}}
\put(-2,-8){$p$}
\put(30,0){\circle*{3}}
\put(28,-8){$a$}
\put(0,30){\circle*{3}}
\put(-8,28){$d$}
\put(30,30){\circle*{3}}
\put(32,23){$o$}
\put(60,30){\circle*{3}}
\put(63,27){$b$}
\put(30,60){\circle*{3}}
\put(28,63){$c$}
\end{picture}

 Since the affine space $X$ is hypersymmetric, there exist hyperfixed automorphism $A$ of $X$ such that $Aapd=dpa$. By Proposition~\ref{p:hyperfixed=>paracentral}. the hyperfixed automorphism $A$ is paracentral and $\Fix(A)$ is a line in $X$. Lemma~\ref{l:diamond-preserve} ensures that $A(o)=A(\diamond apd)=\diamond apd=o$ and hence $\Fix(A)=\Aline op$.  Since $b\in \Aline do\setminus \Aline op$, the point $b'\defeq A(b)$ belongs to the set $A[\Aline do\setminus \Aline op]=\Aline ao\setminus \Aline op$ and hence $b\notin \Aline op=\Fix(A)$ and $b'\ne b$.
 Since the automorphism $A$ is paracentral, $\Aline b{b'}\parallel \Aline ad\parallel \Aline bc$ and hence $\Aline b{b'}=\Aline bc$ and $b'\in \Aline bc\cap\Aline ao=\{c\}$. Therefore, $A(b)=b'=c$. Assuming that $\Aline op\nparallel\Aline ab$, we can find a common point $q\in \Aline op\cap\Aline ab$ and conclude that $q=A(q)\in A[\Aline ab]=\Aline dc$ and hence $q\in \Aline ab\cap\Aline cd$, which contradicts the definition of a parallelogram (in our case $abcd$). This contradiction shows that $\Aline op\parallel \Aline ab$. 
 
 By analogy we can prove that for the point $p'\defeq\diamond aod'$, the line $\Aline o{p'}$ is parallel to the line $\Aline {a'}{b'}=\Aline ab$ and hence $\Aline o{p'}=\Aline op\parallel\Aline ab$. Then $p'=\diamond abo=p$ and $d'=\diamond aop'=\diamond oap=d$, and finally, $c'=\diamond bad'=\diamond abd=c$.
 \end{proof}
 
Now we can present a proof of Theorem~\ref{t:hypersymmetric=>shear}. Assume that a Playfair plane $X$ is hypersymmetric. To prove that the plane $X$ is shear, take any line $L\subseteq X$ and any distinct points $a,c\in X\setminus L$ such that $\Aline ac\cap L=\varnothing$. We have to construct an automorphism $A:X\to X$ such that $L=\Fix(A)$ and $A(a)=c$. Choose any points $b\in \Aline ac\setminus\{a,c\}$ and $o\in L$. Since $X$ is hypersymmetric, there exist hyperfixed involutive automorphisms $A',A''$ of the liner $X$ such that $A'aob=boa$ and $A''boc=coa$.
By Proposition~\ref{p:hyperfixed=>paracentral}, the hyperfixed automorphisms $A',A''$ are paracentral.  Taking into account that the line $\Aline ac=\Aline ab=\Aline bc$ is parallel to the line $L$ and the authormorphisms $A',A''$ are paracentral, we conclude that for every pair $xy\in A'\cup A''$ the flat $\Aline xy$ is subparallel to the line $L$. Then the automorphism $A\defeq A''A'$ has the same property: for every pair $xy\in A$, the flat $\Aline xy$ is subparallel to the line $L$. This means that the automorphism $A$ is paracentral. Since $A(o)=A''A'(o)=o$, the set $\Fix(A)$ is not empty. By Proposition~\ref{p:aff-paracentral<=>}, the paracentral automorphism $A$ of the affine space $X$ is hyperfixed, and by Proposition~\ref{p:hyperfixed=>paracentral}, the set $\Fix(A)$ is a line in $X$. Assuming that the line $\Fix(A)$ is not equal to the line $L$, we conclude that the lines $\Fix(A)$ and $\Aline ac$ are concurrent and hence contain a common point $x=A(x)=A''A'(x)\in\Aline ac$. Applying to the equality $x=A''A'(x)$ the involutive bijection $A''$, we conclude that $A''(x)=A'(x)$. Consider the point $y\defeq A''(x)=A'(x)$ and observe that $A''(y)=x=A'(y)$, by the involutivity of the automorphisms $A'$ and $A''$. If $x\ne y$, then $xoy$ is a triangle such that $A'xoy=yox=A''xoy$ and hence $A'=A''$, by Lemma~\ref{l:atmost1hypersymmetry}. Then $a=A'(b)=A''(b)=c$, which contradicts the choice of the points $a\ne c$. This contradiction shows that $x=y=A'(x)=A''(x)$ and hence $x\in\Fix(A')\cap\Fix(A'')$. Then $\Fix(A')=\Aline ox=\Fix(A'')$. Consider the points $o'\defeq \diamond aob$ and $o''\defeq \diamond boc$. Lemma~\ref{l:diamond-preserve} ensures that $\Aline o{o'}=\Fix(A')=\Fix(A'')=\Aline o{o''}$. Then $aobo'$ and $cobo''$ are two parallelograms such that $\{x\}=\Aline ab\cap\Aline o{o'}=\Aline bc\cap\Aline o{o''}$. Applying Lemma~\ref{l:central-parallelograms}, we conclude that $a=c$, which contradicts the choice of the (distinct) points $a,c$. This contradiction shows that $\Fix(A)=L$ and hence $A$ is a desired automorphism of the plane $X$ such that $A(a)=c$ and $\Fix(A)=L$.
\end{proof}

\begin{remark} In Theorem~\ref{t:inversive-dot} we shall prove that a Playfair plane is hypersymmetric if and only if it is shear.
\end{remark}

\chapter{Parallel projections in affine spaces}\label{ch:LTrans}

\section{Affine spaces}

\begin{definition} An \index{affine space}\index{space!affine}\defterm{affine space} is any $3$-long affine regular liner of rank $\|X\|\ge 3$.
\end{definition}

Theorems~\ref{t:HA} and \ref{t:affine=>Avogadro} imply that affine spaces admit the following first-order axiomatization.

\begin{theorem}\label{t:affine-first-order} A set $X$ endowed with a ternary relation $\Af\subseteq X^3$ is an affine space if and only if it satisfies the following six first-order axioms:
\begin{itemize}
\item[{\sf (IL)}] {\sf Identity:} $\forall x,y\in X\;\;\big(\Af xyx\;\to\; x=y\big)$;
\item[{\sf (RL)}] {\sf Reflexivity:} $\forall x,y\in X\;\;(\Af xxy\;\wedge\;\Af xyy)$;
\item[{\sf (EL)}] {\sf Exchange:} $\forall a,b,x,y\in X\;\big(( \Af axb\wedge \Af ayb\wedge x\ne y)\to (\Af xay\;\wedge\;\Af xby)\big)$;
\item[{\sf (AR)}] {\sf Regularity:}
$\forall o,a,b,u,v,x,y,z\in X\;$\newline
$( \Af ovu\wedge\Af axv\wedge\Af byu\wedge \Af xzy)\;\Rightarrow\;\exists s,t,c,w\;(\Af aso\wedge\Af bto\wedge \Af sct\wedge\Af owu\wedge \Af czw)$;
\item[{\sf (AF)}] {\sf Affinity:} $\forall o,x,y,p\in X$\newline
$(\Af xpy\wedge \neg \Af opx)\to\exists u\, \big(\Af ouy\wedge \forall v \,(\Af ovy\to (u\ne v\leftrightarrow \exists a\,(\Af vap\wedge\Af oax)))\big)$;
\item[{\sf (ND)}] {\sf Non-Degeneracy:} $\exists x,y,u,v\in X\;(x\ne u\ne y\wedge \Af xuy\wedge \neg\Af xvy)$.
\end{itemize}
\end{theorem}

\begin{remark} By Theorem~\ref{t:4-long-affine}, for a $4$-long liner $(X,\Af)$, the Axiom of Regularity {\sf(AR)} in Theorem~\ref{t:affine-first-order} can be removed. 
\end{remark}



\begin{proposition}\label{p:lines3+} For any parallel lines $L,L'$ in an affine space $X$, there exists a line $L''$ in $X$ such that $L''\parallel L$ and $L''\cap(L\cup L')=\varnothing$.
\end{proposition}

\begin{proof} By Proposition~\ref{p:cov-aff}, there exists a point $x\in X\setminus (L\cup L')$. By Theorem~\ref{t:Playfair}, there exists a line $L''$ that contains the point $x$ and is parallel to the line $L$. By Theorem~\ref{t:Proclus-lines}, $L''\parallel L'$. 
Assuming that $L''\cap L\ne\varnothing$, we can apply Proposition~\ref{p:para+intersect=>coincide} and conclude that $L''=L$. Then $x\in L''=L$, which contradicts the choice of the point $x$. This contradiction shows that $L''\cap L=\varnothing$. By analogy we can prove that $L''\cap L'=\varnothing$.
\end{proof}

By Theorem~\ref{t:Playfair<=>}, every affine space $X$ is Playfair and hence every line $L$ in $X$ is spreading and determines the direction $L_\parallel\in\partial X$. For a direction $\delta\in\partial X$ in an affine space $X$ and a point $x\in X$ we denote by $\overline{x\delta}$ a unique line in the spread $\delta$ that contains the point $x$.

 We say that a direction $\delta\in\partial X$ is \index{direction!subparallel to a flat}\defterm{subparallel} to a flat $A\subseteq X$ if every line $L\in\delta$ is subparallel to $A$. By Corollary~\ref{c:parallel-transitive}, a direction $\delta\in\partial X$ is subparallel to a flat $A$ if and only if some line $L\in\delta$ is subparallel to $A$ if and only if some line $L\in\delta$ is contained in the flat $A$.

Let us recall that a relation $R\subseteq X\times X$ on a liner $X$ is \index{flat relation}\index{relation!flat}\defterm{flat} if for every flat $A,B\subseteq X$, the sets $$R[A]\defeq\{y:\exists a\in A\;(a,y)\in R\}\quad\mbox{and}\quad R^{-1}[B]\defeq\{x:\exists b\in B\;(a,b)\in R\}$$are flat in $X$.

In the following proposition we present an important example of a flat relation, which will be applied to  in the next section to studying parallel projections between flats in affine spaces.

\begin{proposition}\label{p:flat-relation} Let $X$ be an affine space. For two flats $A,B\subseteq X$ and a direction $\delta\in\partial X$, the relation
$$F\defeq\{(x,y)\in A\times B:\Aline x\delta=\Aline y\delta\}$$ is flat and has the following properties.
\begin{enumerate}
\item If the direction $\delta$ is not subparallel to the flat $B$, then $F$ is a function.
\item If $F$ is a nonempty function, then the direction $\delta$ is not subparallel to $B$.
\item  If the direction $\delta$ is not subparallel to the flat $A$, then $F^{-1}$ is a function.
\item If $F^{-1}$ is a nonempty function, then the direction $\delta$ is not subparallel to $A$.
\item If $A\cap B=\varnothing$, then $F$ is an injective function and $\dom[F],\rng[F]$ are disjoint parallel flats in $X$.
\end{enumerate}
\end{proposition}

\begin{proof} 0. To prove that the relation $F$ is flat, fix any flat $C\subseteq X$. We should prove that the sets $F[C]$ and $F^{-1}[C]$ are flat in $X$. To show that $F[C]$ is flat, fix any distinct points $x,y\in F[C]$. We have to check that $\Aline xy\subseteq F[C]$. Fix any point $z\in\Aline xy$. Since $x,y\in F[C]$, there exist points $a,b\in C$ such that $(a,x),(b,y)\in F$. The definition of $F$ ensures that $a,b\in A$, $x,y\in B$ and $\Aline a\delta=\Aline x\delta\in\delta$ and $\Aline b\delta=\Aline y\delta\in\delta $. Since the sets $A,B,C$ are flat, $\Aline ab\subseteq A\cap C$ and $\Aline xy\subseteq B$. Since $\delta$ is a spread of parallel lines, there exists a  unique line $L\in\delta$ such that $z\in L$.

\begin{claim}\label{cl:Lcapxy} $L\cap\Aline ab\ne\varnothing$.
\end{claim}

\begin{proof} If $\Aline a\delta=\Aline b\delta$, then $\Aline x\delta=\Aline a\delta=\Aline b\delta=\Aline y\delta$ and $z\in \Aline xy\subseteq \Aline a\delta$. Since $L,\Aline a\delta\in\delta$ are two parallel lines containing the point $z$, $L=\Aline a\delta$ and hence $a\in L\cap\Aline ab\ne\varnothing$.

So, assume that the parallel lines $\Aline a\delta$ and $\Aline b\delta$ are distinct and hence $\Pi\defeq\overline{\{\Aline a\delta\cup\Aline b\delta\}}$ is a plane containing the points $x,y,z$ and also the line $L$. Since the lines $\Aline ab$ and $\Aline a\delta$ are concurrent and $\Aline a\delta\parallel L$, the line $\Aline ab$ is not parallel to the line $L$ and hence $\Aline ab\cap L\ne\varnothing$, by Corollary~\ref{c:parallel-lines<=>}.
\end{proof} 

By Claim~\ref{cl:Lcapxy}, there exists a point $c\in L\cap\Aline ab\subseteq A\cap C$. Then $\Aline c\delta=\Aline z\delta=L\in\delta$ and hence  $(c,z)\in F$, witnessing that $z\in F[C]$ and $\Aline ab\subseteq F[C]$. So, the set $F[C]$ is flat. By analogy we can prove that the set $F^{-1}[C]$ is flat in $X$.
\smallskip

1. Now we shall prove that the relation $F$ is a function if the direction $\delta$ is not subparallel to the flat $B$.  Given two pairs $(x,y),(x,z)\in F$, we should check that $y=z$. To derive a contradiction, assume that $y\ne z$. It follows from $(x,y),(x,z)\in F$ that $\Aline y\delta=\Aline x\delta=\Aline z\delta\in\delta$ and hence $\Aline yz=\Aline x\delta\in\delta$. Since $\Aline yz\subseteq B$, the direction $\delta$ is subparallel to the flat $B$, which contradicts our assumption. This contradiction shows that $y=z$ and hence $F$ is a function.
\smallskip

2. Assume that $F$ is a nonempty function and fix a pair $(x,y)\in F$. Then $(x,y)\in A\times B$ and $\Aline x\delta=\Aline y\delta\in\delta$. Assuming that $\delta$ is subparallel to the flat $B$, we conclude that $\Aline x\delta=\Aline y\delta\subseteq B$ and hence $\{x\}\times\Aline y\delta\subseteq F$, which means that $F$ is not a function. This contradiction shows that the line $\Lambda$ is not subparallel to $B$.
\smallskip

3,4. The statements 3,4 follow from the statements 1,2 applied to the sets $B,A$ and the function $F^{-1}$ instead of $A,B$ and the function $F$.
\smallskip

5. Assume that $A\cap B=\varnothing$. If $F=\varnothing$, then $F$ is an injective function and the flats $\dom[F]=\varnothing$ and $\rng[F]=\varnothing$ are parallel. So, assume that $F\ne \varnothing$ and fix any pair $(a,b)\in F$. Observe that $a\ne b$ and $\Aline ab=\Aline a\delta=\Aline b\delta\in\delta$. Assuming that the direction $\delta$ is subparallel to the flat $A$, we conclude that $a\in \Aline ab=\Aline b\delta\subseteq A$, which contradicts $A\cap B=\varnothing$. This contradiction shows that the direction $\delta$ is not subparallel to $A$. By analogy we can show that $\delta$ is not subparallel to $B$. By the statements 2 and 4, $F$ and $F^{-1}$ are functions and hence $F$ is an injective function. 

Since the relation $F$ is flat, its domain $\dom[F]=F^{-1}[X]$ and range $\rng[F]=F[X]$ are flats. It remains to check that these flats are parallel. 

To show that $\dom[F]\subparallel \rng[F]$, fix any point $a\in\dom[F]$. Given any point $x\in \dom[F]$, we should prove that $x\in\overline{\{a\}\cup\rng[F]}$. This is clear if $x=a$. So, we assume that $x\ne a$. Let $y\defeq F(x)$ and $b\defeq F(a)$. Since $\dom[F]\cap\rng[F]=\varnothing$, the points $a,b$ are distinct and hence  $\Aline ab\in\delta$. By analogy we can show that $\Aline xy\in\delta$ and hence $\Aline xy\parallel\Aline ab$. Then $\overline{\{x,y,a,b\}}$ is a plane containing two disjoint lines $\Aline ax\subseteq\dom[F]$ and $\Aline by\subseteq\rng[F]$. By Theorem~\ref{t:parallel-char}, $\Aline ax\parallel \Aline by$ and hence $x\in\overline{\{a\}\cup\Aline by}\subseteq\overline{\{a\}\cup \rng[F]}$, witnessing that $\dom[F]\subparallel \rng[F]$. By analogy we can show that $\rng[F]\subparallel \dom[F]$. Therefore, $\dom[F]\parallel\rng[F]$.
\end{proof}

\section{Parallel projections in affine spaces}

We recall that for a point $x$ in an affine space $X$ and a direction $\delta\in\partial X$, we denote by $\Aline x\delta$ a unique line $L\in\delta$ that contains the point $x$. 

\begin{definition}\label{d:parproj} Let $X$ be an affine space.  A function $F\subseteq X\times X$ is called a \index{parallel projection}\index{projection!central}\defterm{parallel projection in a direction $\delta\in\partial X$} if $$F=\{(x,y)\in \dom[F]\times \rng[F]:\Aline x\delta=\Aline y\delta\}.$$ A function $F\subseteq X\times X$ is called \defterm{parallel projection} if $F$ is a parallel projection in some direction $\delta\in\partial X$. If $\dom[F]$ and $\rng[F]$ are lines in $X$, then $F$ is called a \index{line projection}\index{line!projection}\defterm{line projection in $X$}. 
\end{definition}

\begin{picture}(200,110)(-200,-60)

{\linethickness{1pt}
\put(0,0){\color{blue}\line(2,1){70}}
\put(0,0){\color{cyan}\line(2,-1){70}}
\put(0,0){\color{cyan}\line(-2,1){70}}
\put(0,0){\color{blue}\line(-2,-1){70}}
}

\put(-2,40){\color{red}$\delta$}
\put(77,-38){\color{cyan}$\rng[F]$}
\put(75,34){\color{blue}$\dom[F]$}
\put(60,30){\color{red}\vector(0,-1){60}}
\put(-60,-30){\color{red}\vector(0,1){60}}

\put(30,15){\color{red}\vector(0,-1){30}}
\put(-30,-15){\color{red}\vector(0,1){30}}

\put(60,30){\color{blue}\circle*{3}}
\put(60,-30){\color{cyan}\circle*{3}}
\put(-60,30){\color{cyan}\circle*{3}}
\put(-60,-30){\color{blue}\circle*{3}}

\put(30,15){\color{blue}\circle*{3}}
\put(30,-15){\color{cyan}\circle*{3}}
\put(-30,15){\color{cyan}\circle*{3}}
\put(-30,-15){\color{blue}\circle*{3}}
\end{picture}

Applying Proposition~\ref{p:flat-relation}, we obtain the following general example of a flat parallel projection.

\begin{example} Let $A,B$ be two flats in an affine regular liner $X$ and a direction $\delta\in\partial X$ is not subparallel to $B$. Then $$F\defeq\{(x,y)\in A\times B:\Aline x\delta=\Aline y\delta\}$$ is a flat parallel projection in the direction $\delta$. The function $F$ is injective if  the direction $\delta$ is not subparallel to the flat $A$. If $A\cap B=\varnothing$, then the function $F$ is injective and flats $\dom[F],\rng[F]$ are disjoint and parallel.
\end{example}

\begin{example}\label{ex:para-identity} For every flat $A\ne X$ in a $3$-ranked space $(X,\Af)$, the identity map $1_A$ is a parallel projection.
\end{example}

\begin{proof} Since $A\ne X$, there exist points $a\in A$ and $b\in X\setminus A$, which determine a direction $\delta\defeq(\Aline ab)_\parallel$, which is not subparallel to the flat $A$. Then $1_A=\{(x,y)\in A\times A:\Aline x\delta=\Aline y\delta\}$ is a parallel projection in the direction $\delta$. 
\end{proof}

Sometimes it will be convenient to use the following alternative desription of parallel projections (that does not involve directions).

\begin{proposition} Let $X$ be an affine space. A function $F\subseteq X\times X$ is a parallel projection (in a direction $\delta\in\partial X$) if and only if there exists a line $L$ (with $L\in\delta$) such that $$F=\{(x,y)\in\dom[F]\times\rng[F]:\Aline xy\subparallel L\}.$$
\end{proposition}

\begin{proposition}\label{p:para-projection} Let $X$ be an affine space and $F\subseteq X\times X$ be a parallel projection in a direction $\delta\in\partial X$. 
\begin{enumerate}
\item For every $x\in\dom[F]\cap\rng[F]$ we have $F(x)=x$.
\item If $\dom[F]=\rng[F]$, then $F$ is the identity map of the set $\dom[F]=\rng[F]$.
\item The function $F$ is flat if and only if $\dom[F]$ and $\rng[F]$ are flats.
\item If $\dom[F]$ and $\rng[F]$ are disjoint flats, then the function $F$ is injective and the flats $\dom[F]$ and $\rng[F]$ are parallel. 
\item If $\dom[F]$ and $\rng[F]$ are two lines, then the function $F$ is injective.
\end{enumerate}
\end{proposition}

\begin{proof} By Definition~\ref{d:parproj}, $F=\{(x,y)\in\dom[F]\times\rng[F]:\Aline x\delta=\Aline y\delta\}$. 
\smallskip

1. For every $x\in\dom[F]\cap\rng[F]$ we have $\Aline x\delta=\Aline x\delta$ and hence $(x,x)\in F$ and $F(x)=x$.
\smallskip

2. The second statement follows immediately from the first one.
\smallskip

3. If the function $F$ is flat, then its domain $\dom[F]=F^{-1}[X]$ and range $\rng[F]=F[X]$ are flat. On the other hand, if $\dom[F]$ and $\rng[F]$ are flat, then the function $F=\{(x,y)\in\dom[F]\times\rng[F]:\Aline xy\subparallel \Lambda\}$ is flat by Proposition~\ref{p:flat-relation}.
\smallskip

4. If $\dom[F]$ and $\rng[F]$ are disjoint flats, then the function $F$ is injective and the flats $\dom[F],\rng[F]$ are parallel, by Proposition~\ref{p:flat-relation}.
\smallskip

5. Assume that $\dom[F]$ and $\rng[F]$ are two lines. By Proposition~\ref{p:flat-relation}, the injectivity of the function $F$ will follow as soon as we check that $\dom[F]\notin\delta$. To derive a contradiction, assume that $\dom[F]\in\delta$. Fix any pair $(a,b)\in F$ and observe that $\Aline a\delta=\Aline b\delta\in \delta$ and $a\in \dom[F]\in\delta$ implies $\dom[F]=\Aline a\delta=\Aline b\delta$. Then for every $x\in\dom[F]$ we have $\Aline x\delta=\dom[L]=\Aline b\delta$ and hence $(x,b)\in F$ and $\dom[F]\times\{b\}=F$ and $\rng[F]=\{b\}$ is not a line, which contradicts our assumption. This contradiction shows that the direction $\delta$ is not subparallel $\dom[F]$. By Proposition~\ref{p:flat-relation},  the function $F$ is injective and hence $F$ is a bijection between the lines $\dom[F]$ and $\rng[F]$. 
\end{proof}

\begin{definition} A parallel projection $P\subseteq X\times X$ in a liner $X$ is called a \index{parallel shift}\defterm{parallel shift} if $\dom[P]$ and $\rng[P]$ are parallel flats in $X$. A function $T\subseteq X\times X$ is called a \index{parallel translation}\defterm{parallel translation} if $T$ is the composition of finitely many parallel shifts.
\end{definition}

By Proposition~\ref{p:para-projection}(2,4), every parallel shift $S$ in an affine space $X$ is a bijective flat function between the parallel flats $\dom[S]$ and $\rng[S]$. Consequently, every parallel translation $T$ in $X$ is a bijective flat function between the flats $\dom[T]$ and $\rng[T]$.

\begin{proposition}\label{p:parshift=}  For every parallel translation $T\subseteq X\times X$ in an affine space $X$, the flats $\dom[T]$ and $\rng[T]$ and parallel.
\end{proposition}

\begin{proof} Write $T$ as the composition $P_n\cdots P_1$ of parallel shifts $P_1,\dots,P_n$. For every $k\in\{1,\dots,n\}$ consider the parallel translation $F_k=P_k\cdots P_1$. 
Since every parallel shift is flat, the functions $F_1,\dots,F_n$ are flat. Consequently, $\dom[T]=\dom[F_n]\subseteq\dom[P_1]$ is a flat in $X$. For every $k\in\{1,\dots,n\}$, consider the function $$P'_k\defeq P_k\cap(X\times F_k[\dom[T]]),$$ and observe that $\dom[P'_1]=\dom[T]$ and $\dom[P'_k]=\rng[P'_{k-1}]$ for every $k\in\{2,\dots,n\}$. 

If $\dom[P_k]\cap\rng[P_k]\ne\varnothing$ for some $k\in\{1,\dots,n\}$, then $\dom[P_k]\parallel \rng[P_k]$ implies $\dom[P_k]=\rng[P_k]$ and Proposition~\ref{p:para-projection}(2) ensures that $P_k$ is the identity map of the flat $\dom[P_k]=\rng[P_k]$. Then $P'_k$ is the identity map of the flat $\dom[P_k]\cap F_k[\dom[T]]$, and hence $\dom[P'_k]=\rng[P'_k]$ and $\dom[P'_k]\parallel\rng[P'_k]$.

Next, assume that $\dom[P_k]\cap\rng[P_k]=\varnothing$. Since $P_k$ is a parallel projection, there exists a direction $\delta_k$ in $X$ such that $$P_k=\{(x,y)\in\dom[P_k]\times\rng[P_k]:\Aline x{\delta_k}=\Aline y{\delta_k}\}.$$Then  $$P_k'=\{(x,y)\in \dom[P_k]\times(\rng[P_k]\cap F_k[\dom[T]]):\Aline x{\delta_k}=\Aline y{\delta_k}\}.$$ Applying Proposition~\ref{p:flat-relation}, we conclude that $P_k'$ is a flat function. Since $\dom[P'_k]\cap\rng[P_k']\subseteq\dom[P_k]\cap\rng[P_k]=\varnothing$, the flats $\dom[P_k']$ and $\rng[P_k']$ are parallel, by Proposition~\ref{p:flat-relation}(5).

Then $$\dom[T]=\dom[P'_1]\parallel \rng[P_1']=\dom[P'_2]\parallel \rng[P_2']=\dots=\dom[P'_n]\parallel\rng[P_n']=\rng[T]$$and hence $\dom[T]\parallel \rng[T]$, by Corollary~\ref{c:parallel-transitive}.
\end{proof}

\begin{exercise} Let $F:X\to Y$ be an isomorphism between liners $X,Y$. Show that for every (flat) parallel projection $P\subseteq X\times X$, the function $FPF^{-1}\subseteq Y\times Y$ is a (flat) parallel projection.
\end{exercise}

\section{Line projections in affine spaces}

\begin{definition} A \index{line projection}\defterm{line projection} in a liner $X$ is a parallel projection $P\subseteq X\times X$ such that $\dom[P]$ and $\rng[P]$ are lines in $X$. 
\end{definition} 

\begin{exercise}\label{ex:conj-lineprojection} Let $F:X\to Y$ be an isomorphism between two liners $X,Y$. For every line projection $P$ in $X$, the function $FPF^{-1}$ is a line projection in the liner $Y$.
\end{exercise}

Let us recall that a {\em line bijection} is a bijective function between two lines in a liner. By Proposition~\ref{p:para-projection}, every line projection $P$ in a liner $X$ is a line bijection and hence $P$ is an element of the inverse monoid $\I_X^\ell$. Let us recall that $\I_X^\ell$ is the smallest submonoid of the inverse symmetric monoid $\I_X$, containing all line bijections in $X$. The monoid $\I_X^\ell$ consists of the identity map of $X$, all line bijections in $X$, and all trivial bijections, i.e., bijections between flats of cardinality $\le 1$. 

The smallest submonoid of $\I_X^\ell$ that contains all line projections in $X$ is denoted by \index[note]{$\I_X^{\prop}$}$\I_X^{\prop}$. Line bijections that belong to the monoid $\I_X^{\prop}$ are called \defterm{line affinities}. Alternatively, line affinities can be defined as follows. 

\begin{definition}\label{d:line-affinity} A line bijection $F$ in a liner $X$ is called a \index{line affinity}\defterm{line affinity} if $F$ is the composition of finitely many line projections in $X$.
\end{definition}

\begin{exercise} Let $F:X\to Y$ be an isomorphism of two liners $X,Y$. Show that for every line affinity $A$ between two lines in $X$, the function $FAF^{-1}$ is a line affinity between two lines in the liner $Y$.
\end{exercise} 


The aim of the next three sections will be to understand the structure of line affinities in Desargesian affine spaces.

\begin{exercise} Show that $\mathcal{I}_X^{\prop}$ is an inverse submonoid of the inverse monoid $\mathcal{I}^\ell_X$.
\smallskip

\noindent{\em Hint:} Use the equality $(F G)^{-1}=G^{-1} F^{-1}$ holding for any functions $F,G$.
\end{exercise}

Proposition~\ref{p:flat-relation}(5) implies the following important fact.

\begin{proposition}\label{p:paraconcurrent} For every line projection $P$ in an affine liner, the lines $\dom[P]$ and $\rng[P]$ are either concurrent or parallel.
\end{proposition}

Proposition~\ref{p:paraconcurrent} motivates the following definition.

\begin{definition}\label{d:shift-shear} A line projection $P$ in a liner is called
\begin{itemize}
\item a \index{line shift}\defterm{line shift} if $\dom[P]\parallel \rng[P]$;
\item a \index{line shear}\defterm{line shear} if $\dom[P]\cap\rng[P]$ is a singleton, called the \defterm{shear center} of $P$;
\end{itemize}
\end{definition}

Propositions~\ref{p:paraconcurrent}  imply that every line projection in an affine regular space $X$ is either a line shift or a line shear. 

\begin{exercise} Let $F:X\to Y$ be an isomorphism of two liners $X,Y$. Show that for every line shift (resp. line shear) $P$ between two lines in $X$, the function $FPF^{-1}$ is a line shift (resp. line shear) between two lines in $Y$.
\end{exercise} 

For a liner $X$ we denote by \index[note]{$\I_X^{\#}$}$\I_X^\#$ the smallest submonoid of $\I_X^{\prop}$, containing all line shifts in $X$. Line bijections that belong to the monoid $\I_X^\#$ are called {\em line translations}. Alternatively, line translations can be defined as follows. 

\begin{definition} A line bijection $F$ in a liner $X$ is called a \index{line translation}\defterm{line translation} if $F$ is the finite composition of line shifts. 
\end{definition}

\begin{exercise} Let $F:X\to Y$ be an isomorphism of two liners $X,Y$. Show that for every line translation $T$ between two lines in $X$, the function $FTF^{-1}$ is a line translation between two lines in $Y$.
\end{exercise} 

\begin{exercise} Show that $\mathcal{I}_X^{\#}$ is an inverse submonoid of the inverse monoid $\mathcal{I}^{\prop}_X$.
\end{exercise}

Therefore, for every liner $X$, we obtain the following increasing chain of inverse monoids:
$$\I_X^\#\subseteq \I_X^{\prop}\subseteq \I_X^{\ell}\subseteq \I_X.$$

For every line $L$ in a liner $X$, the subset $$\Sym^\#_X(L)\defeq\{T\in\I_X^\#:\dom[T]=L=\rng[T]\}\subseteq\Sym(L)$$of all line translations $T:L\to L$ is a subgroup of the semigroup $\I_X^\#$, called the \index{group of line translations}\index{$\Sym^\#_X(L)$}\defterm{group of line translations} of the line $L$. The study of the structure of the groups $\Sym_X^\#(L)$ will be continued in Section~\ref{s:IX[L;u,v]}.
\smallskip

Recall that two lines in a liner are \index{paraconcurrent lines}\index{lines!paraconcurrent}\defterm{paraconcurrent} if they are parallel or concurrent.

\begin{theorem}\label{t:paraproj-exists} For every paraconcurrent lines $L,L'$ in an affine space $X$ and every points $a\in L\setminus L'$ and $a'\in L'\setminus L$, there exists a unique line projection $P:L\to L'$ such that $P(a)=a'$.
\end{theorem}

\begin{proof} The choice of the (distinct) points $a,a'$ ensures that the lines $L,L'$ are distinct and $\delta\defeq (\Aline a{a'})_\parallel$ is a direction in $X$.  We claim that the relation $P\defeq\{(x,y)\in L\times L':\Aline x\delta=\Aline y\delta\}$ is a required line projection with $\dom[P]=L$, $\rng[P]=L'$ and $P(a)=a'$. First, we show that $P$ is a function. Since the line $\Aline a{a'}$ is not subparallel to the lines $L,L'$, Proposition~\ref{p:flat-relation} guarantees that $P$ is an injective function. 



Next, we show that $\dom[P]=L$. The definition of $P$ ensures that $\dom[P]\subseteq L$. To prove that $L\subseteq \dom[P]$, choose any point $x\in L$. If $x\in L'$, then $\Aline xx\subparallel \Aline a{a'}$ ensures that $(x,x)\in P$ and hence $x\in \dom[P]$. If $x=a$, then $\Aline a\delta=\Aline a{a'}=\Aline {a'}\delta$ ensures that $(a,{a'})\in P$ and hence  $x=a\in\dom[P]$. It remains to consider the case of $x\in L\setminus(\{a\}\cup L')$. Assuming that $x\in\Aline a{a'}$, we conclude that $a'\in\Aline a{a'}=\Aline ax=L$, which contradicts the choice of the point $a'$. This contradiction shows that $x\notin\Aline a{a'}$. Consider the line $\Aline x\delta$ and observe that  $\Aline x\delta\subseteq\overline{\{x,a,a'\}}\subseteq\overline{L\cup L'}$.  Since the lines $L,L'$ are distinct and paraconcurrent, their flat hull $\overline{L\cup L'}$ is a plane in the affine space $X$. By Corollary~\ref{c:proregular=>ranked}, the affine space is ranked and hence $L'\subseteq \overline{L\cup L'}=\overline{L_x\cup\Aline a{a'}}$. By Theorem~\ref{t:Proclus2}, there exists a unique  point $y\in \Aline x\delta$. It follows from $\Aline x\delta=\Aline y\delta$ that $(x,y)\in P$ and hence $x\in\dom[P]$ and $\dom[P]=L$. By analogy we can prove that $\rng[P]=L'$.

The definition of the function $P$ ensure that $(a,a')\in P$ and hence $P(a)=a'$. To see that $P$ is a unique parallel projection with $\dom[P]=L$, $\rng[P]=L'$ and $P(a)=a'$, assume that $F$ is another parallel projection with  $\dom[F]=L$, $\rng[F]=L'$ and $F(a)=a'$. Then $F=\{(x,y)\in L\times L':\Aline x{\delta'}=\Aline y{\delta'}\}$ for some direction $\delta'\in\partial X$. It follows from $F(a)=a'\ne a$ that $\Aline a{a'}\in\delta'\cap\delta$ and hence $\delta'=\delta$ and $F=P$.
\end{proof}

For concurrent lines in affine spaces, Theorem~\ref{t:paraproj-exists} implies the following corollary on the existence of line shears with prescribed values.

\begin{corollary}\label{c:shear-exists} For every lines $L,L'$ in an affine space and points $o\in L\cap L'$, $a\in L\setminus L'$ and $a'\in L'\setminus L$, there exists a unique line shear $R:L\to L'$ such that $Roa=oa'$.
\end{corollary}

Next, we derive from Theorem~\ref{t:paraproj-exists} two corollaries on the existense of line shifts and line translations in affine spaces.

\begin{corollary}\label{c:parashift-exists} For every line $L$ in an affine space $X$ and every points $x\in L$ and $x'\in X\setminus L$ there exists a unique line shift $S$ in $X$ such that $\dom[S]=L$ and $S(x)=x'$.
\end{corollary}

\begin{proof} By Theorem~\ref{t:Playfair}, there exists a unique line $L'$ in $X$ such that $x'\in L'\parallel L$. It follows from $x'\notin L$ and $L\parallel L'$ that $L'\cap L=\varnothing$. By Theorem~\ref{t:paraproj-exists}, there exists a line projection $S:L\to L'$ such that $S(x)=x'$. Since $\dom[S]=L\parallel L'=\rng[S]$, the line projection $S$ is a line shift.

Now assume that $S'$ is another line shift such that $\dom[S']=L$ and $S'(x)=x'$. By Definition~\ref{d:shift-shear}, $\rng[S']\parallel \dom[S']=L\parallel L'$ and hence $\rng[S']\parallel L'$ and $\rng[S']=L'$ as $x'\in\rng[S']\cap L'\ne\varnothing$. The uniqueness part of Theorem~\ref{t:paraproj-exists} guarantess that $S=S'$.
\end{proof}

\begin{corollary}\label{c:par-trans} For every parallel lines $L,L'$ in an  affine space $X$ and every points $x\in L$ and $x'\in L'$ there exists a line translation $T:L\to L'$ such that $T(x)=x'$. Moreover, the line translation $T$ is the composition of two line shifts  in $X$.
\end{corollary}

\begin{proof}  By Proposition~\ref{p:lines3+}, there exists a line $\Lambda$ in $X$ such that $\Lambda\parallel L$ and $\Lambda\cap (L\cup L')=\varnothing$. Choose any point $y\in \Lambda$. By Theorem~\ref{t:paraproj-exists}, there exist line shifts $S:L\to\Lambda$ and $S':\Lambda\to L'$ such that  $P(x)=y$ and $P'(y)=x'$. Then the composition $T\defeq P' P:L\to L'$ is a desired line translation such that   $T(x)=P'(P(x))=P'(y)=x'$.
\end{proof}

The following proposition implies that for an affine space $X$, the monoid $\mathcal{I}^{\prop}_X$ is generated by line shears.

\begin{proposition}\label{p:shift=2shears} Let $S$ be a line shift in an affine space $X$. Then for every points $a\in\dom[S]$ and $b'\in\rng[S]\setminus\{S(a)\}$ there exist unique line shears $R_1:\dom[S]\to\Aline a{b'}$ and $R_2:\Aline a{b'}\to\rng[S]$ such that $S=R_2 R_1$. 
\end{proposition}

\begin{proof} Fix any points $a\in\dom[S]$ and $b'\in\rng[S]\setminus\{S(a')\}$. Consider the points $a'=S(a)\ne b'$ and $b=S^{-1}(b')\ne a$. By Theorem~\ref{t:paraproj-exists}, there exist  line shears $R_1:\dom[S]\to\Aline a{b'}$ and $R_2:\Aline a{b'}\to\rng[S]$ such that $R_1(b)=b'$ and $R_2(a)=a'$.
We claim that $S=R_2 R_1$. 

\begin{picture}(200,80)(-100,-20)

{\linethickness{1pt}
\put(0,0){\line(1,0){140}}
\put(0,40){\line(1,0){140}}
\put(10,50){\line(2,-1){120}}
}
\put(30,40){\vector(0,-1){40}}
\put(110,40){\vector(0,-1){40}}
\put(60,40){\vector(0,-1){40}}

\put(145,38){$\dom[S]$}
\put(147,-2){$\rng[S]$}

\put(30,40){\circle*{4}}
\put(28,43){$a$}
\put(60,40){\circle*{3}}
\put(58,43){$x$}
\put(110,40){\circle*{3}}
\put(108,43){$b$}
\put(60,25){\circle*{3}}
\put(63,28){$z$}
\put(30,0){\circle*{3}}
\put(27,-10){$a'$}
\put(60,0){\circle*{3}}
\put(52,-10){$y{=}y'$}
\put(110,0){\circle*{4}}
\put(106,-10){$b'$}

\end{picture}

Fix any point $x\in S$ and consider the points $y\defeq S(x)$, $z=R_1(x)$ and $y'=R_2(z)$. Taking into account that $S$, $R_1$ and $R_2$ are line projections, we can apply the transitivity of the subparallelity relation in affine spaces (see Corollary~\ref{c:subparallel-transitive}) and conclude that $\Aline a{a'}\parallel \Aline xy\parallel \Aline b{b'}$, $\Aline xz\subparallel \Aline b{b'}$ and $\Aline z{y'}\subparallel \Aline a{a'}$. Then also $\Aline x z\subparallel \Aline xy$ and $\Aline z{y'}\subparallel \Aline xy$. If $z=x$, then $\Aline x{y'}=\Aline z{y'}\subparallel \Aline xy$. If $z=y'$, then $\Aline x{y'}=\Aline xz\subparallel \Aline xy$. If $x\ne z\ne y'$, then $\Aline xz$ and $\Aline z{y'}$ are two intersecting lines that are subparallel to the line $\Aline xy$. By Corollary~\ref{c:subparallel}, $\Aline xz\parallel \Aline xy$ and $\Aline z{y'}\parallel\Aline xy$. The proaffinity of the affine space $X$ ensures that $\Aline xz=\Aline z{y'}$ and hence $\Aline x{y'}=\Aline xz\subparallel \Aline xy$. Therefore, in all cases we obtain that $\Aline x{y'}\subparallel \Aline xy$ and hence $\Aline x{y'}=\Aline xy$. Then $\{y\}=\rng[S]\cap\Aline xy=\rng[S]\cap\Aline x{y'}=\{y'\}$ and hence $S(x)=y=y'=R_2(R_1(x))$, witnessing that $S=R_2 R_1$.

To see that the line shears $R_2$ and $R_1$ are unique, assume that $R_1':\dom[S]\to\Aline a{b'}$ and $R'_2:\Aline a{b'}\to\rng[S]$ are two line shears such that $S=R'_2 R'_1$. Proposition~\ref{p:para-projection}(2) ensures that $R_1'(a)=a$ and $R_2'(b')=b'$. Then $R_2(a)=a'=S(a)=R_2'(R_1'(a))=R_2'(a)$ and $R_2=R_2'$, by the uniqueness part of Theorem~\ref{t:paraproj-exists}. On the other hand, by the injectivity of the line shear $R_2'$, the equations $R'_2(b')=b'=S(b)=R_2'(R_1'(b))$ imply $R_1'(b)=b'=R_1(b)$ and then $R_2=R_2'$, by the uniqueness part of Theorem~\ref{t:paraproj-exists}.
\end{proof}
 
Now we show that a trivial bijection in affine spaces can be represented as compositions of two line shifts or two line shears.

\begin{theorem}\label{t:singleton-bijection} Let $F:A\to B$ be a bijection between sets of cardinality $|A|=|B|\le 1$ in an affine space $X$. Then 
\begin{enumerate}
\item $F=S_2 S_1$ for some line shifts $S_1,S_2$ in $X$;
\item $F=R_2 R_1$ for some line shears $R_1,R_2$ in $X$;
\item $F=R_3 S_3$ for some line shift $S_3$ and line shear $R_3$ in $X$;
\item $F=S_4 R_4$ for some line shear $R_4$ and line shift $S_4$ in $X$.
\end{enumerate}
\end{theorem}

\begin{proof} Three cases are possible.

1. $A=B=\varnothing$. The space $X$, being affine of dimension $>1$, contains two disjoint parallel lines $L$ and $\Lambda$. By Theorem~\ref{t:paraproj-exists}, there exists a line shift $S:L\to \Lambda$. Then for the line shift $S_2=S_1=S$, we obtain $S_2 S_1=\varnothing=F$. Since the liner $X$ has dimension $>1$, there exist lines $L'$ and $\Lambda'$ such that the intersections $L\cap L'$ and $\Lambda\cap\Lambda'$ are singletons.  By Theorem~\ref{t:paraproj-exists}, there exist  line shears $R_1:L'\to L$ and $R_2:\Lambda\to \Lambda'$. Since the lines $L$ and $\Lambda$ are disjoint, $R_2 R_1=\varnothing =F$. Also for the line shift $S_3\defeq S$ and the line shear $R_3\defeq R_1^{-1}$ we obtain $R_3 S_3=\varnothing=F$. Finally, for the line shear $R_4\defeq R_1$ and the line shift $S_4\defeq S^{-1}$ we obtain $S_4 R_4=\varnothing =F$.  
\smallskip

2. $A=B\ne\varnothing$. Let $a$ be the unique point of the singleton $A=B$. Fix any line $L$ is $X$ containing the point $a$ and choose a point $b\in X\setminus L$. By Theorem~\ref{t:Playfair}, there exists a line $\Lambda$ such that $b\in\Lambda$ and $\Lambda\parallel L$. By  Theorem~\ref{t:paraproj-exists}, there exists a line shift $S_1:L\to\Lambda$ such that $S_1(a)=b$. By Proposition~\ref{p:lines3+}, there exists a point $c\in X\setminus(L\cup\Lambda)$. By Theorem~\ref{t:Playfair}, there exists a line $\Lambda'$ such that $b\in\Lambda'$ and $\Lambda'\parallel\Aline ac$. Then $\Lambda'\cap\Lambda=\{b\}$.  By  Theorem~\ref{t:paraproj-exists}, there exists a line shift $S_2:\Lambda'\to\Aline ac$ such that $S_2(b)=a$. Then $S_2 S_1=1_A=F$. 

By Theorem~\ref{t:paraproj-exists}, there exists a shear $R:L\to \Aline ab$. Then for the shears $R_2= R_1=R$, we obtain $R_2 R_1=R R=1_A=F$.   By Theorem~\ref{t:paraproj-exists}, there exists a shear $R_3:\Aline ab\to \Lambda$. Then for the line shift $S_3\defeq S_1:L\to\Lambda$, we obtain $R_3 S_3=1_A=F$. Finally, for the line shear $R_4\defeq R$ and the line shift $S_4\defeq S_1$, we obtain $R_4 S_4=1_A=F$.
\smallskip

3. $A\ne B$. Let $a\in A$ and $b\in B$ be unique points of the singletons $A$ and $B$, respectively. Since $\dim(X)>1$, there exists a point $c\in X\setminus\Aline ab$. By Theorem~\ref{t:Playfair}, there exist unique lines $L_a,L_b,L_c$  such that $a\in L_a\parallel \Aline bc$, $b\in L_b\parallel \Aline ac$, and $c\in L_c\parallel \Aline ab$. By Theorem~\ref{t:paraproj-exists}, there exist line shifts $S:L_a\to\Aline bc$, $S_1:L_a\to\Aline bc$, and $S_2:L_c\to\Aline ab$ such that $S(a)=b$, $S_1(a)=c$ and $S_2(c)=b$. Then $F=S_2 S_1$.

 By Theorem~\ref{t:paraproj-exists}, there exist a line shear $R_1:\Aline ac\to\Aline bc$ such that $R_1(a)=b$, and a line shears $R_2:\Aline ab\to\Aline bc$and $R_4:\Aline ab\to\Aline ac$. It is easy to see that $R_2 R_1=F$. Also for the line shifts $S_3=S_4\defeq S:L_a\to\Aline cb$ and the line shear $R_3 \defeq R_2:\Aline ab\to\Aline bc$ we have $R_3 S_3=F=S_4 R_4$.
\end{proof}

Theorem~\ref{t:singleton-bijection}(1) implies

\begin{corollary} For any affine space $X$, the semigroups $\mathcal{I}^\#_X\subseteq\mathcal{I}^{\prop}_X$ contain all bijections between subsets of cardinality $\le 1$ in $X$.
\end{corollary}

\begin{theorem}\label{t:aff-trans} For every affine space $X$ and points $a,b,a',b'\in X$ with $a\ne b$ and $a'\ne b'$, there exists a line affinity $A:\Aline ab\to\Aline{a'}{b'}$ such that $F(a)=a'$ and $F(b)=b'$. Moreover, the line affinity $A$ is the composition of two line shears or two line shifts.
\end{theorem}

\begin{proof} Using  Proposition~\ref{p:cov-aff}, choose a point $c\in X\setminus(\Aline ab\cup \Aline{a'}{b'})$. Five cases are possible. 

1. If $a=a'$, then using Theorem~\ref{t:paraproj-exists}, find line shears $R:\Aline ab\to\Aline ac$ and $R': \Aline {a'}c\to \Aline{a'}{b'}$ such that $R(b)=c$ and $R'(c)=b'$, and observe that the affine function $F\defeq R' R$ has the required property: $(F(a),F(b))=(a',b')$.

2.  If $b=b'$, then using Theorem~\ref{t:paraproj-exists}, find line shears $R:\Aline ba\to\Aline bc$ and $R':\Aline {b'}c\to \Aline{b'}{a'}$ such that $R(a)=c$ and $R'(c)=a'$, and observe that the affine function $F\defeq R' R$ has the required property: $(F(a),F(b))=(a',b')$.

3. If $a'\notin \Aline ab$ and $b\notin \Aline{a'}{b'}$, then using Theorem~\ref{t:paraproj-exists},  find line shears $R:\Aline ba\to\Aline b{a'}$ and $R': \Aline {a'}b\to \Aline{a'}{b'}$ such that $R(a)=a'$ and $R'(b)=b'$. Then the affine function $F\defeq R' R$ has the required property: $(F(a),F(b))=(a',b')$.

4. If $b'\notin \Aline ab$ and $a\notin \Aline{a'}{b'}$, then using Theorem~\ref{t:paraproj-exists}, find line shears $R:\Aline ab\to\Aline a{b'}$ and $R':\Aline {b'}a\to \Aline{b'}{a'}$ such that $R(b)=b'$ and $R'(a)=a'$. Then the affine function $F\defeq R' R$ has the required property: $(F(a),F(b))=(a',b')$.

5. If none of the above four cases holds, then $\Aline ab=\Aline{a'}{b'}$. By Theorem~\ref{t:Playfair}, there exist unique lines $L$ and $L'$ in $X$ such that $b\in L\parallel \Aline ac$ and $b'\in L'\parallel \Aline c{a'}$. It follows that the lines $L$ and $L'$ are subsets of the plane $\overline{\{a,b,c\}}=\overline{L\cup\Aline ac}=\overline{L'\cup\Aline c{a'}}$. Since $\Aline ac\cap\Aline c{a'}=\{c\}$, we can apply Theorem~\ref{t:Proclus2} (two times) and conclude that the intersections $L\cap \Aline c{a'}$ and $L\cap L'$ are singletons. Let $c'$ be the unique point of the intersection $L\cap L'$. 
Consider the parallel projections $$P\defeq\{(x,y)\in \Aline ab\times\Aline c{c'}:\Aline xy\subparallel L\}\quad\mbox{and}\quad P'\defeq\{(x,y)\in\Aline c{c'}\times \Aline{a'}{b'}:\Aline xy\subparallel L'\}$$and observe that the affine transformation $T\defeq P_2 P_1$ have the required property:
$$(T(a),T(b))=(P'(P(a)),P'(P(b)))=(P'(c),P(c'))=(a',b')).$$ 
If $\Aline c{c'}\parallel \Aline ab$, then $P,P'$ are line shifts. In the other case, $P,P'$ are line shears.
\end{proof}

\begin{exercise} Show that for any line $L$ in an affine space,  any non-identity line translation $T:L\to L$ cannot be writen as the composition of two line shears.
\end{exercise}

\begin{Exercise} Show that every line bijection in an affine space $X$ is the composition of three parallel projections between suitable subsets of $X$.
\end{Exercise}

\section{Line translations in Thalesian affine spaces}

We recall that an affine space $X$ is \defterm{Thalesian} if and only if for every parallel lines $A,B,C$ in $X$ and points $a,a'\in A\setminus (B\cup C)$,  $a,a'\in B\setminus (A\cup C)$,  $a,a'\in C\setminus (A\cup B)$, if $\Aline ab\parallel \Aline {a'}{b'}$ and $\Aline bc\parallel \Aline {b'}{c'}$, then $\Aline ac\parallel \Aline{a'}{c'}$.

\begin{theorem}\label{t:shift=AMA} An affine space $X$ is Thalesian if and only if for any line shifts $F:L_1\to L_2$ and $G:L_2\to L_3$ between lines $L_1,L_2,L_3$ in $X$ with $L_1\ne L_3$, the composition $G F:L_1\to L_3$ is a line shift.
\end{theorem}

\begin{proof} To prove the ``only if'' part, assume that an affine space $X$ is Thalesian. Given any
line shifts $F:L_1\to L_2$ and $G:L_2\to L_3$ in $X$ with $L_1\ne L_3$, we should check that the composition $G F:L_1\to L_3$ is a line shift. By Proposition~\ref{p:parshift=},
$L_1\parallel L_2\parallel L_3$ and by Corollary~\ref{c:parallel-transitive}, $L_1\parallel L_3$.
The definition of a line shift ensures that $L_1=\dom[F]\ne\rng[F]=L_2$ and $L_2=\dom[G]\ne\rng[F]=L_3$. Therefore, $L_1,L_2,L_3$ are three distinct parallel lines in $X$.

 Fix any point $a_1\in L_1$ and consider the points $a_2=F(a)\in L_2$ and $a_3=F(a_2)\in L_3$. By Theorem~\ref{t:paraproj-exists}, there exists a line shift $S:L_1\to L_3$ such that $S(a_1)=a_3$. We claim that that $G F=S$. Given any point $x_1\in L_1$, we should prove that $GF(x_1)=S(x_1)$.
Consider the points $x_2=F(x_1)\in L_2$, $x_3=G(x_2)\in L_3$, and $x_3'=S(x_3)$. Since $F$, $G$ and $S$ are line shifts, $\Aline{x_1}{x_2}\parallel \Aline{a_1}{a_2}$, $\Aline{x_2}{x_3}\parallel\Aline {a_2}{a_3}$, and $\Aline{x_1}{x_3'}\parallel\Aline{a_1}{a_3}$, by  Corollary~\ref{c:parallel-transitive}. Applying the Thales Axiom, we conclude that $\Aline{x_1}{x_3}\parallel \Aline{a_1}{a_3}\parallel \Aline{x_1}{x_3'}$ and hence $\Aline{x_1}{x_3}=\Aline{x_1}{x_3'}$, by the proaffinity of the liner $X$. Then $\{x_3\}=L_3\cap\Aline{x_1}{x_3}=L_3\cap\Aline{x_1}{x_3'}=\{x_3'\}$ and hence $P(x_1)=x_3'=x_3=GF(x_1)$, witnessing that $G F=S$ is a line shift. 
\smallskip 

To prove the ``if'' part, assume that  for any line shifts $F:L_1\to L_2$ and $G:L_2\to L_3$ in $X$ with $L_1\ne L_3$, the composition $G F:L_1\to L_3$ is a line shift. To check that $X$ is Thalesian, take any parallel lines $A,B,C\subseteq X$ and points $a,a'\in A\setminus(B\cup C)$, $b,b'\in B\setminus(A\cup C)$, $c,c'\in C\setminus(A\cup B)$ such that $\Aline ab\parallel \Aline{a'}{b'}$ and $\Aline bc\parallel \Aline{b'}{c'}$. The choice of the points $a,b,c$ ensures that the lines $A,B,C$ are distinct. By Theorem~\ref{t:paraproj-exists}, there exist line shifts $F:A\to B$ and $G:B\to C$ such that  $F(a)=b$ and $G(b)=c$. By our assumption, the composition $G F$ is a line shift. Taking into account that $\Aline {a'}{b'}\parallel \Aline ab$ and $\Aline{b'}{c'}\parallel \Aline bc$, we conclude that $F(a')=b'$ and $G(b')=c'$. Then $G F(a')=c'$. Since $G F$ is a parallel projection with $\{(a,c),(a',b')\}\in G F$, the definition of a parallel projection and Corollary~\ref{c:parallel-transitive} imply that $\Aline {a'}{c'}\parallel \Aline ac$, witnessing that the Thales Axiom holds and the affine space $X$ is Thalesian.
\end{proof}

Now we prove several lemmas on composition of line shifts in Thalesian affine spaces.

\begin{lemma}\label{l:shift=shift+shift} Let $S$ be a line shift in a Thalesian affine space $X$. For every line shift $S_1$ in $X$ with $\dom[S_1]=\dom[S]$ and $\rng[S_1]\ne\rng[S]$, there exists a line shift $S_2$ such that $S=S_2 S_1$.
\end{lemma}

\begin{proof} By Proposition~\ref{p:paraconcurrent}, we have $\rng[S]\parallel \dom[S]=\dom[S_1]\parallel \rng[S_1]$ and hence $\rng[S]\parallel \rng[S_1]$, according to Corollary~\ref{c:parallel-transitive}.
Fix any point $x\in \dom[S]=\dom[S_1]$ and consider the points $y\defeq S(x)$ and $y_1\defeq S_1(x)$. By Theorem~\ref{t:paraproj-exists}, there exists a line shift $S_2:\rng[S_1]\to\rng[S]$ such that $S_2(y_1)=y$. By Theorem~\ref{t:shift=AMA}, the composition $S_2 S_1$ is a line shift with $S_2 S_1(x)=S_2(y_1)=y=S(x)$, and by Corollary~\ref{c:parashift-exists}, $S=S_2 S_1$.
\end{proof}



\begin{lemma}\label{l:3shifts=2shifts} The composition of three line shifts in a Thalesian affine space $X$ is equal to the composition of two lines shifts in $X$.
\end{lemma}

\begin{proof} Let $P_1,P_2,P_3$ be line shifts in a Thalesian affine space $X$, and $P=P_3 P_2 P_1$. If $\|\dom[P]\|\le 1$, then by Theorem~\ref{t:singleton-bijection}, $P$ is the composition of two line shifts in $X$.

So, we assume that $\|\dom[P]\|=2$. In this case $\dom[P]=\dom[P_1]$, $\rng[P_1]=\dom[P_2]$, $\dom[P_2]=\rng[P_3]$ and $\rng[P_3]=\rng[P]$.  

If $\dom[P_1]\ne\rng[P_2]$, then by Theorem~\ref{t:shift=AMA}, the composition $S=P_2 P_1$ is a line shift and hence $P=P_3  S$ is the composition of two line shifts.

If $\dom[P_2]\ne\rng[P_3]$, then by Theorem~\ref{t:shift=AMA}, the composition $S'=P_3 P_2$ is a line shift and hence $P=S' P_1$ is the composition of two line shifts. 

So, we assume that $\dom[P_1]=\rng[P_2]$ and $\dom[P_2]=\rng[P_3]$.  Then $$L\defeq\dom[P]= \dom[P_1]=\rng[P_2]=\dom[P_3]\quad\mbox{and}\quad \Lambda\defeq\rng[P]=\rng[P_3]=\dom[P_2]=\rng[P_1]$$ are two distinct parallel lines in $X$.  By Proposition~\ref{p:lines3+}, there exists a line $L'$ in $X\setminus(L\cup \Lambda)$, which is parallel to the lines $L$ and $\Lambda$. Fix any points $x\in L$ and $x'\in L'$ and consider the point $y=P_1(x)\in\Lambda$. By Theorem~\ref{t:paraproj-exists}, there exist line shifts $S$ and $S'$ such that $\dom[S]=\dom[P_1]$, $\rng[S]=L'=\dom[S']$, $\rng[S']=\rng[P_1]$ and $S(x)=x'$, $S'(x')=y$. By Theorem~\ref{t:shift=AMA}, the composition $S' S:\dom[P_1]\to\rng[P_1]$ is a line shift such that $S'S(x)=P_1(x)$. By Theorem~\ref{t:paraproj-exists}, $S' S=P_1$ and hence $P=P_3 P_2 P_1=P_3 P_2 S' S$. By Theorem~\ref{t:shift=AMA}, the composition $P_2 S'$ is a line shift and so is the compositions $P_3(P_2 S')$. Therefore $P$ is equal to the composition $(P_3(P_2 S')) S$ of two line shifts $P_3(P_3 S')$ and $S$.
\end{proof} 

\begin{corollary}\label{c:many-shifts=2shifts} Let $X$ be a Thalesian affine space. Every line  translation $T$ in $X$ is the composition of two line shifts  in $X$.
\end{corollary}

\begin{proof} The definition of a line translation ensures that $T=P_n\cdots P_1$ for some line shifts in $X$. We can assume that the number $n$ is the smallest possible. Lemma~\ref{l:3shifts=2shifts} ensures that $n\le 2$. If $n=2$, then $T=P_2 P_1$ is the composition of two line shifts.
If $n=1$, then $T=P_1$ is the composition of two line shifts by Lemma~\ref{l:shift=shift+shift}.
\end{proof} 

\begin{theorem}\label{t:unique-translation} An affine space $X$ is Thalesian if and only if for every line $L$ in $X$ and every points $x\in L$ and $x'\in X$ there exists a unique line translation $T$ in $X$ such that $\dom[T]=L$ and $T(x)=x'$.
\end{theorem}

\begin{proof} To prove the ``if'' part, assume that the affine space $X$ is not Thalesian. By Theorem~\ref{t:shift=AMA}, there exist three distinct parallel lines $L_1,L_2,L_3$ and two line shifts $P_1:L_1\to L_2$ and $P_2:L_2\to L_3$ such that the composition $T\defeq P_2 P_1$ is not a line shift. Fix any point $x\in L_1$ and consider the point $x'=T(x)\in L_3\subseteq X\setminus L_1$.
By Theorem~\ref{t:paraproj-exists}, there exists a line shift $P:L_1\to L_3$ such that $P(x)=x'$. Then $P$ and $T$ are two distinct line translations with $\dom[P]=L_1=\dom[T]$ and $P(x)=x'=T(x)$.
\smallskip

To prove the ``only if'' part, assume that the affine space $X$ is Thalesian. Fix a line $L$ and points $x\in L$ and $x'\in X$. Since $X$ is an affine space, there exists a unique line $L'$ in $X$ such that $x'\in L'\parallel L$. By Corollary~\ref{c:par-trans}, there exists a parallel translation $T:L\to L'$ such that $T(x)=x'$. It remains to prove that $T$ is a unique parallel translation with $\dom[T]=L$ and $T(x)=x'$. Assume that $T'$ is another parallel translation with $\dom[T']=L$ and $T'(x)=x'$.  By Proposition~\ref{p:parshift=}, $\rng[T']\parallel \dom[T]=L$ and hence $\rng[T']\parallel L'$. Since $x'\in L'\cap\rng[T']$, the lines $\rng[T']$ and $L'=\rng[T]$ coincide.

If $L\ne L'$, then by Corollary~\ref{c:many-shifts=2shifts} and Theorem~\ref{t:shift=AMA}, the parallel translations $T$ and $P$ are line shifts and by Corollary~\ref{c:parashift-exists}, $T=P$.

Next, assume that $L=L'$. In this case, Corollary~\ref{c:many-shifts=2shifts} ensures that $T=T_2 T_1$ and $T'=T'_2 T'_1$ for some line shifts $T_1,T_2,T'_1,T'_2$ in $X$. Consider the lines $\Lambda\defeq\rng[T_1]=\dom[T_2]$ and $\Lambda'=\rng[T_1]=\dom[T_2']$.
Proposition~\ref{p:parshift=} ensures that $\Lambda\parallel \dom[T_1]=L=\dom[T_1']\parallel\Lambda'$ and hence the lines $\Lambda$ and $\Lambda'$ are parallel.
Two cases are possible.
\smallskip

1. First assume that $\Lambda\ne \Lambda'$. Let $y=T_1(x)$ and $z=T_1'(x)$. By Theorem~\ref{t:paraproj-exists}, there exists a line shift $S:\Lambda\to\Lambda'$ such that $S(y)=z$.
Theorem~\ref{t:shift=AMA} ensures that $S T_1$ is a line shift. Since $\dom[S T_1]=L=\dom[T_1']$ and $S T_1(x)=z=T_1'(x)$, Corollary~\ref{c:parashift-exists} guarantees that $S T_1=T_1'$ and hence $T'=T_2' T_1'=T_2' S T_1$. By Theorem~\ref{t:shift=AMA}, the composition $T_2' S:\Lambda\to L'=L$ is a line shift. Since 
$\dom[T_2' S]=\Lambda=\dom[T_2]$ and $T_2' S(y)=T_2'(z)=T_2'(T_1'(x))=x'=T_2(T_1(x))=T_2(y)$, Corollary~\ref{c:parashift-exists} guarantees that $T_2' S=T_2$ and hence $$T'=T'_2 T'_1=T_2' (S T_1')=(T_2' S) T_1=T_2 T_1=T.$$
\smallskip

2. Next, assume that $\Lambda=\Lambda'$. By Proposition~\ref{p:lines3+}, there exists a line $I\subseteq X\setminus(L\cup \Lambda)$, which is parallel to the lines $L$ and $\Lambda$.  By Theorem~\ref{t:paraproj-exists}, there exists a line shift $P:L\to I$. Consider the point $p=P(x)\in I$.
 By Theorem~\ref{t:paraproj-exists}, there exist line shifts $S,S'$ in $X$ such that $\dom[S]=I=\dom[S']$, $\rng[S]=\Lambda=\Lambda'=\rng[S']$, $S(p)=T_1(x)$ and $S'(p)=T_1'(x)$. By Theorem~\ref{t:shift=AMA}, the compositions $S P$ and $S' P$ are line shifts and Corollary~\ref{c:parashift-exists}, $S P=T_1$ and $S' P=T_1'$. Then $T=T_2 T_1=T_2 S P$ and $T'=T_2' T'_1=T_2' S' P$. By Theorem~\ref{t:shift=AMA}, the composition $T_2 S$ and $T'_2 S'$ are line shifts. Since $T_2(S(p))=T_2(T_1(x))=T(x)=x'=T_2'(T_1'(x))=T_2'(S'(p))$ and Corollary~\ref{c:parashift-exists}, $T_2 S=T_2' S'$. Then
$$T=T_2 T_1=T_2 (S P)=(T_2 S) P=(T_2' S') P=T_2'\circ(S' P)=T_2' T_1'=T'.$$
\end{proof}

\section{Translations versus line translations}

In this section we study the interplay between translations and line shifts in affine spaces.

\begin{proposition}\label{p:translation=>shift} Let $T:X\to X$ be a  translation of an affine space $X$. For every line $L\subseteq X$ with $T[L]\ne L$, the restriction $T{\restriction}_L$ is a line shift of $X$.
\end{proposition}

\begin{proof}  Since $T$ is a dilation, $L\ne T[L]\parallel L$ and hence $T[L]\cap L=\varnothing$. Fix any point $a\in L$ and consider the point $b\defeq T(a)\in T[L]$. By Proposition~\ref{p:Trans-spread}, the family $\lines T\defeq\{\Aline xy:xy\in T\}$ is a spread of parallel lines. Then $T{\restriction}_L=\{xy\in L\times T[L]:\Aline xy\parallel \Aline ab\}$ is a line shift.
\end{proof}

\begin{proposition}\label{p:translation|=>line-translation} Let $X$ be an affine space such that $\|\{T(o):T\in\Trans(X)\}\|\ne 2$ for some point $o\in X$. Then for every translation $T\in\Trans(X)$ and every line $L\subseteq X$, the restriction $T{\restriction}_L:L\to T[L]$ is a line translation.
\end{proposition}

\begin{proof} If $T(o)=o$, then $T$ is an identity map of $X$, by Proposition~\ref{p:Ax=Bx=>A=B}. In this case $T{\restriction}_L$ is the identity map of $L$ and hence a line translation of $L$. So, assume that $T(o)\ne o$. If $T[L]\ne L$, then $T{\restriction}_L$ is a line shift and hence a line translation, by Proposition~\ref{p:translation=>shift}. So, assume that $T[L]=L$. Since $\|\{A(o):A\in\Trans(X)\}\|\ne 2$, and $\{o,T(o)\}\subseteq \{A(o):A\in\Trans(X)\}$, there exists a nonidentity translation $A\in\Trans(X)$ such that $\lines A\ne\lines T$ and hence $A[L]\ne L\ne A^{-1}[L]$. Consider the line $\Lambda\defeq A[L]$ and observe that $\Lambda\parallel L$ and $\Lambda\cap L=\varnothing$. Since $\Trans(X)$ is a subgroup of the automorphism group $\Aut(X)$, the automorphism  $AT$ is a translation with $AT[L]=A[L]=\Lambda\ne L$. By Proposition~\ref{p:translation=>shift}, the restrictions $AT{\restriction}_L:L\to\Lambda$ and $A^{-1}{\restriction}_\Lambda:\Lambda\to L$ are line shifts. Then their composition $A^{-1}{\restriction}_\Lambda\circ AT{\restriction}_L:L\to L$ is a line translation, equal to the restriction $T{\restriction}_L$.
\end{proof} 

\begin{question} Is there a translation $T:X\to X$ of an affine space $X$ such that for some line $L\subseteq X$, the restriction $T{\restriction}_L$ is not a line translation?
\end{question}

\begin{lemma}\label{l:TS=ST} For every line shift $S$ in an affine space $X$, and every translation $T:X\to X$ with $T[\dom[S]]=\dom[S]$, we have $TS=ST$.
\end{lemma}

\begin{proof} The equality $T[\dom[S]]=\dom[S]$ implies $\dom[TS]=\dom[S]=T^{-1}[\dom[S]]=\dom[ST]$. Therefore, the equality $TS=ST$ will follow as soon as we check that $TS(o)=ST(o)$ for every point $o\in \dom[S]$. Given any point $o\in \dom[S]$, consider the points $t\defeq T(o)\in T[\dom[S]]=\dom[S]$, $s\defeq S(o)\in\Lambda$, $s_t\defeq S(t)=ST(o)$ and $t_s\defeq T(s)=TS(o)$. Since $T$ is a dilation, $\Aline os\parallel \Aline t{t_s}$ and hence $\Pi\defeq\overline{\{o,s,t,t_s\}}$ is a plane containing the lines $\Aline ot=\dom[S]$ and $\Aline s{t_s}$. By Proposition~\ref{p:Trans-spread}, the lines $\Aline ot$ and $\Aline s{t_s}$ in the plane $\Pi$ are disjoint and parallel. Since $S$ is a shift between the lines $\dom[S]$ and $\rng[S]$, those lines are parallel. Since $s\in \rng[S]\cap\Aline s{t_s}$, the lines $\rng[S]$ and $\Aline s{t_s}$ coincide by the Playfair property of the affine regular liner $X$. 

Since $T$ is a dilation and $S$ is a line shift, $\Aline t{t_s}\parallel \Aline os\parallel \Aline t{s_t}$ and hence $\Aline t{t_s}=\Aline t{s_t}$, by the Playfair property of the affine regular liner $X$. Then $\{t_s\}=\Lambda\cap\Aline t{t_s}=\Lambda\cap\Aline t{s_t}=\{s_t\}$ and hence $TS(o)=t_s=s_t=ST(o)$.
\end{proof}

\begin{theorem}\label{t:TS=ST} For every line translation $S$ in an affine space $X$ and every translation $T:X\to X$ with $T[\dom[S]]=\dom[S]$, we have $TS=ST$.
\end{theorem}

\begin{proof} To derive a contradiction, assume that there exist a line translation $S$ and a translation $T:X\to X$ such that $T[\dom[S]]=\dom[S]$  and $ST\ne TS$. The inequality $ST\ne TS$ implies that $T$ is a nonidentity translation of $X$ and $S$ is a non-identity translation of the line $\dom[S]$.
By Proposition~\ref{p:Trans-spread}, the family $\overline T=\{\Aline xy:xy\in T\}$ is a spread of parallel lines in $X$. Write the line translation $S$ as the composition $S=S_1S_2\cdots S_n$ of line shifts. We can assume that the number $n$ is the smallest possible. Lemma~\ref{l:TS=ST} implies that $n\ge 2$. The minimality of $n$ ensures that $S_2\cdots S_nT=TS_2\cdots S_n$. Since $S_2\cdots S_n$ is a line translation, $\dom[S_1]=\rng[S_2\cdots S_n]\parallel \dom[S_n]$, by Proposition~\ref{p:parshift=}. Since $T[\dom[S]]=\dom[S]=\dom[S_n]\parallel \dom[S_1]$, the line $\dom[S_1]$ belongs to the spread $\overline T=\{\Aline xy:xy\in T\}$ and hence $T[\dom[S_1]]=\dom[S_1]$. In this case Lemma~\ref{l:TS=ST} ensures that $TS_1=S_1T$ and hence
$$ST=S_1S_2\cdots S_nT=S_1TS_2\cdots S_n=TS_1S_2\cdots S_n=TS,$$
which contradicts the choice of $T$ and $S$. 
\end{proof}


 
The following important theorem was proved as Lemma 2.1 in a paper \cite{Gleason1956} of \index[person]{Gleason}Andrew Gleason\footnote{{\bf Andrew Mattei Gleason} (1921 -- 2008) was an American mathematician who made fundamental contributions to widely varied areas of mathematics, including the solution of Hilbert's fifth problem, and was a leader in reform and innovation in math­e­mat­ics teaching at all levels. Gleason's theorem in quantum logic and the Greenwood--Gleason graph, an important example in Ramsey theory, are named for him.\newline \indent
As a young World War II naval officer, Gleason broke German and Japanese military codes. After the war he spent his entire academic career at Harvard University, from which he retired in 1992. His numerous academic and scholarly leadership posts included chairmanship of the Harvard Mathematics Department and the Harvard Society of Fellows, and presidency of the American Mathematical Society. He continued to advise the United States government on cryptographic security, and the Commonwealth of Massachusetts on math­e­mat­ics education for children, almost until the end of his life.
\newline \indent Gleason won the Newcomb Cleveland Prize in 1952 and the Gung--Hu Distinguished Service Award of the American Mathematical Society in 1996. He was a member of the National Academy of Sciences and of the American Philosophical Society, and held the Hollis Chair of Mathematics and Natural Philosophy at Harvard.
\newline \indent He was fond of saying that math­e­mat­ic­al proofs 
``really aren't there to convince you that something is true --- 
they're there to show you why it is true.'' 
The Notices of the American Mathematical Society called him ``one of the quiet giants of twentieth-century mathematics, the consummate professor dedicated to scholarship, teaching, and service in equal measure.''
}, a very interesting and important person in the history of American mathematics.

\begin{theorem}[Gleason, 1956]\label{t:Gleason1} Let $X$ be an affine space, $L$ be a line in $X$. For a bijection $P:L\to L$ the following conditions are equivalent:
\begin{enumerate}
\item there exists a translation $T:X\to X$ such that $T{\restriction}_L=P$;
\item for every line translation $S:L\to L$ we have $SP=PS$;
\item for every line $\Lambda\subseteq X$ and line shifts $A:L\to\Lambda$ and $B:\Lambda\to L$ we have $(BA)P=P(BA)$;
\item for every line $\Lambda\subseteq X$ and line shifts $A,B:\Lambda\to L$, we have $A^{-1}PA=B^{-1}PB$;
\end{enumerate}   
\end{theorem}

\begin{proof} The implication $(1)\Ra(2)$ follows from Theorem~\ref{t:TS=ST}, and $(2)\Ra(3)$ is trivial.
\smallskip

To prove that $(3)\Ra(4)$, assume that  for every line $\Lambda\subseteq X$ and line shifts $A:L\to\Lambda$ and $B:\Lambda\to L$ we have $(BA)P=P(BA)$. Then for every  line shifts $A,B:\Lambda\to L$, we have $BA^{-1}P=PBA^{-1}$ and hence $A^{-1}PA=B^{-1}BA^{-1}PA=B^{-1}PBA^{-1}A=B^{-1}PB$.
\smallskip

$(4)\Ra(1)$ Assume that the condition (4) is satisfied. Given any point $x\in X$,  find a unique line $L_x\in L_\parallel$ containing the point $x$, and fix any line shift  $A_x:L_x\to L$ between the parallel lines $L_x$ and $L$. If $x\in L$, then $L_x=L$ and $A_x$ is the identity map of the line $L=L_x$. 
Consider the function $T:X\to X$ defined by $T(x)\defeq A_x^{-1}PA_x(x)$ for every $x\in X$. The condition (4) ensures that the point $T(x)$ does not depend on the choice of the line shift $A_x:L_x\to L$. It is easy to see that $T:X\to X$ is a bijective function with $T^{-1}(y)=A_y^{-1}P^{-1}A_y(y)$ and $T{\restriction}_L=P$. 

\begin{claim}\label{cl:x'y'||xy} For every distinct points $x,y\in X$ and the points $x'\defeq T(x)$ and $y'\defeq T(y)$, the lines $\Aline{x'}{y'}$ and $\Aline xy$ are parallel.
\end{claim}

\begin{proof} We recall that $L_x,L_y$ are unique lines in the spread of parallel lines $L_\parallel$ such that $x\in L_x$ and $y\in L_y$, and $A_x:L_x\to L$, $A_y:L_y\to L$ are line shifts. Six cases are possible. 
\smallskip

1. If $L_x=L=L_y$, then $\Aline xy=L$ and $\Aline{x'}{y'}=L=\Aline xy$.
\smallskip

2. If $L_x=L\ne L_y$, then we lose no generality assuming that $A_y(y)=x$. In this case $x'=T(x)=P(x)$ and $y'=T(y)=A_y^{-1}PA_y(y)=A_y^{-1}(x')$. Since $A_y$ is a line shift with $A_y(y)=x$ and $A_y(y')=x'$, the line $\Aline xy$ is parallel to the line $\Aline {x'}{y'}$. \smallskip

3. If $L_x\ne L=L_y$, then by analogy with the preceding case we can show that the line  $\Aline{x'}{y'}$ is parallel to the line $\Aline xy$.
\smallskip

4. If $L_x=L_y\ne L$, then $\Aline {x'}{y'}=L_x=L_y=\Aline xy$ and hence $\Aline{x'}{y'}\parallel \Aline xy$.
\smallskip

5. If the parallel lines $L,L_x,L_y$ are distinct and coplanar in $X$, then the line $\Aline xy$ has a unique common point $z$ with the line $L$. In this case we lose no generality assuming that $A_x(x)=z=A_y(y)$. Consider the point $z'\defeq T(z)=P(z)$ and observe that $x'=T(x)=A_x^{-1}PA_x(x)=A_x^{-1}(z')$ and $y'=T(y)=A_y^{-1}PA_y(y)=A_y^{-1}(z')$. Since $A_x$ and $A_y$ are line shifts, $\Aline {x'}{z'}\parallel\Aline xz$ and $\Aline {y'}{z'}\parallel \Aline yz$. Since $\Aline xz=\Aline yz$, the parallel lines $\Aline {x'}{z'}$ and $\Aline {y'}{z'}$ coincide and hence $\Aline {x'}{y'}=\Aline{x'}{z'}\parallel \Aline xz=\Aline xy$.
\smallskip

6. Finally assume that the parallel lines $L,L_x,L_y$ are non-coplanar. In this case the affine space has rank $\|X\|\ge 4$ and is Desarguesian, by Corollary~\ref{c:affine-Desarguesian}. We lose no generality assuming that $A_x(x)=z=A_y(x)$ for some point $z\in L$. Consider the point $z'\defeq P(z)=T(z)$ and observe that $x'=T(x)=A_x^{-1}PA_x(x)=A_x^{-1}(z')$ and $y'=T(y)=A_y^{-1}PA_y(y)=A_y^{-1}(z')$. Since $A_x,A_y$ are line shifts, $\Aline xz\parallel \Aline{x'}{z'}$ and $\Aline yz\parallel \Aline{y'}{z'}$, and hence $\Aline{x'}{y'}\parallel \Aline {x}{y}$ by Theorem~\ref{t:Desargues-affine}.
\end{proof}

Now we shall prove that the bijective function $T:X\to X$ is a dilation of $X$. Given any line $\Lambda\subseteq X$, we should prove that its image $T[\Lambda]$ is a line, parallel to the line $\Lambda$. Fix any distinct points $x,y\in\Lambda$ and consider their images $x'\defeq T(x)$ and $y'\defeq T(y)$, which are distinct because $T$ is a bijection of $X$. By Claim~\ref{cl:x'y'||xy}, the line $\Aline{x'}{y'}$ is parallel to the line $\Aline xy$. So, it suffices to show that $T[\Lambda]=\Aline {x'}{y'}$. Given any point $z\in \Aline xy\setminus\{x\}$, consider its image $z'\defeq T(z)$. Claim~\ref{cl:x'y'||xy} ensures that $\Aline {x'}{z'}\parallel \Aline xz=\Aline xy\parallel \Aline {x'}{y'}$ and hence $z'\in\Aline{x'}{z'}=\Aline {x'}{y'}$. This proves that $T[\Lambda]\subseteq\Aline {x'}{y'}$. Applying the same argument to the bijection $T^{-1}$, we can prove that $T^{-1}[\Aline{x'}{y'}]\subseteq \Aline xy=\Lambda$ and hence $T[\Lambda]=\Aline {x'}{y'}$ is a line, parallel to the line $\Lambda=\Aline xy$. Therefore, $T$ is a dilation of the affine space $X$. 

\begin{claim}\label{cl:T(o)=o=>T=Id} If $T(o)=o$ for some $o\in X\setminus L$, then $T$ is the identity map of $X$.
\end{claim}

\begin{proof} For every point $z\in L$, we can choose a line shift $A:L_o\to L$ such that $A(o)=z$. The condition (4) ensures that $o=T(o)=A_o^{-1}PA_o(o)=A^{-1}PA(o)=A^{-1}P(z)$ and hence $z=A(o)=AA^{-1}P(z)=P(z)$, which means that $P$ is the identity map of the line $L$ and then $T$ is the identity map of $X$.
\end{proof}

\begin{claim}\label{cl:T(o)=o=>T=I} If $T(o)=o$ for some $o\in X$, then $T$ is the idenity map of $X$.
\end{claim}

\begin{proof} If $o\notin L$, then $T$ is the identity map of $X$, by Claim~\ref{cl:T(o)=o=>T=Id}. So, assume that $T(o)=o\in L$. Since the affine space $X$ has rank $\|X\|\ge 3$, there exists a point $x\in X\setminus L$. We lose no generality assuming that $A_x(x)=o$. Then $T(x)=A_x^{-1}PA_x(x)=A_x^{-1}P(o)=A_x^{-1}T(o)=A_x^{-1}(o)=x$. Applying Claim~\ref{cl:T(o)=o=>T=Id}, we conclude that $T$ is the identity map of $X$.
\end{proof}

Claim~\ref{cl:T(o)=o=>T=I} ensures that the dilation $T$ is a translation of $X$. 
\end{proof}

Let us recall that a liner $X$ is called {\em translation} if for every points $x,y\in X$ there exists a translation $T:X\to X$ such that $T(x)=y$.

\begin{theorem}\label{t:paraD<=>translation|} For an affine space $X$, the following conditions are equivalent:
\begin{enumerate}
\item the liner $X$ is  translation;
\item for every line shift $S:L\to\Lambda$ between lines $L$ and $\Lambda$ in $X$, there exists a translation $T:X\to X$ such that $S=T{\restriction}_L$;
\item for every line translation $S:L\to\Lambda$ between lines $L$ and $\Lambda$ in $X$, there exists a translation $T:X\to X$ such that $S=T{\restriction}_L$;
\item for every line translation $S:L\to L$ of a line $L$  in $X$, there exists a translation $T:X\to X$ such that $S=T{\restriction}_L$.
\end{enumerate}
\end{theorem}

\begin{proof} $(1)\Ra(2)$. Assume that the affine space $X$ is translational and fix any line shift $S:L\to\Lambda$ between lines $L$ and $\Lambda$ in $X$. By Proposition~\ref{p:paraconcurrent}, the lines $L,\Lambda$ are parallel. If $L=\Lambda$, then the line shift $S$ is the identity map of $L$, according to Proposition~\ref{p:para-projection}(1). Then $S=T{\restriction}_L$ for the identity translation $T$ of $X$. So, assume that $L\ne\Lambda$. In this case $L\cap\Lambda=\varnothing$. Take any point $a\in L$ and consider the point $b\defeq S(a)$. Since $X$ is a translation affine space, there exists a translation $T:X\to X$ such that $T(a)=b$. Since $T$ is a translation of $X$, the line $T[L]$ is parallel to the line $L$ and hence is parallel to the line $\Lambda$, by Theorem~\ref{t:Proclus-lines}. Since $b\in T[L]\cap\Lambda$, the parallel lines $T[L]$ and $\Lambda$ coincide. By Proposition~\ref{p:Trans-spread}, for every pair $xy\in T$ the line $\Aline xy$ is parallel to the line $\Aline ab$. Then the restriction $T{\restriction}_L=\{(x,y)\in L\times\Lambda:\Aline xy\parallel \Aline ab\}$ coincides with the line shift $S$.
\smallskip

$(2)\Ra(3)$ Let $S:L\to\Lambda$ be a line translation between lines $L$ and $\Lambda$ in $X$. By definition of a line translation, there exist lines shifts $S_1,\dots,S_n$ in $X$ such that $S=S_1\cdots S_n$ and hence $\dom[S_n]=\dom[S]=L$ and $\dom[S_i]=\rng[S_{i+1}]$ for every $i\in\{1,\dots,n-1\}$.
 By the condition (2), for every $i\in\{1,\dots,n\}$ there exists a translation $T_i\in\Trans(X)$ such that $S_i=T_i{\restriction}_{\dom[S_i]}=T_iS_i^{-1}S_i$. Since $\Trans(X)$ is a group, the composition $T\defeq T_1\cdots T_n$ is a translation of $X$. We claim that the translation $T$ has the required property: $T{\restriction}_L=S$. 

By downward induction we shall prove that for every $k\in\{2,\dots,n\}$ we have the equality 
$$T{\restriction}_L=T_1\cdots T_{k-1}S_k\cdots S_n.
$$
 For $k=n$ the equality $$T{\restriction}_L=T_1\cdots T_n{\restriction}_L=T_1\cdots T_nS_n^{-1}S_n=T_1\cdots T_{n-1}S_n$$holds because $T_n{\restriction}_L=S_n=T_nS_n^{-1}S_n$. Assume that for some $k\in\{3,\dots,n\}$ we know that $T{\restriction}_L=T_1\cdots T_{k-1}S_{k}\cdots S_n$. It follows from $\rng[S_k]=\dom[S_{k-1}]$ that $S_k=S_{k-1}^{-1}S_{k-1}S_k$ and hence  $$T{\restriction}_L=T_1\cdots T_{k-1}S_k\cdots S_n=T_1\cdots T_{k-1}S_{k-1}^{-1}S_{k-1}S_k\cdots S_n=T_1\cdots T_{k-2}\cdot S_{k-1}S_k\cdots S_n.$$
This completes the inductive step. For $k=2$ we have the equality 
$$T{\restriction}_L=T_1S_2S_3\cdots S_n=T_1(S_1^{-1}S_1S_2)S_3\cdots S_n=(T_1S_1^{-1}S_1)S_2\cdots S_n=S_1S_2\cdots S_n=S.$$ 

The implication $(3)\Ra(4)$ is trivial.
\smallskip

$(4)\Ra(1)$ Assume that the affine space $X$ satisfies the condition (4). Given any points $x,y\in X$, fix any line $L\subseteq X$ containing  the points $x,y$.  By Corollary~\ref{c:par-trans}, there exists a line translation $S:L\to L$ such that $S(x)=y$. By the condition (4), there exists a translation $T\in\Trans(X)$ such that $S=T{\restriction}_L$ and hence $T(x)=S(x)=y$, and the liner $X$ is translational. 
\end{proof}

We recall that for an affine space $X$, we denote  by $\I^{\#}_X$ the smallest submonoid of the symmetric inverse monoid $\I_X$, containing all line shifts in $X$.

\begin{theorem} For an affine space $X$, the following conditions are equivalent:
\begin{enumerate}
\item $X$ is  translation;
\item $\I_X^{\#}=\{1_X\}\cup\{T{\restriction}_L:T\in\Trans(X),\;\overline L=L\subseteq X,\;\|L\|\le 2\}$;
\item every subgroup in the semigroup $\I_X^{\#}$ is commutative;
\item for every line shifts $A,B,C,D$ in $X$ with $\dom[A]=\rng[B]=\dom[C]=\rng[D]$ we have $(DC)(BA)=(BA)(DC)$.
\end{enumerate}
\end{theorem}

\begin{proof} Denote the set $\{1_X\}\cup \{T{\restriction}_L:T\in\Trans(X),\;\overline L=L\subseteq X,\;\|L\|\le 2\}$ by $\Tau$.
\smallskip 

$(1)\Ra(2)$ Assuming that the affine space $X$ is translation, we shall show that  $\I_X^{\#}=\Tau$.  The definition of the monoid $\I_X^{\#}$ implies that every element $A$ of $\I_X^{\#}$ is either the identity map $1_X$ of $X$ or $A=S_1\cdots S_n$ for some line shifts $S_1,\dots,S_n$ of $X$. In the latter case, $\dom[A]$ is a flat of rank $\le 2$. If $\|\dom[A]\|=0$, then $A=\varnothing=1_X{\restriction}_{\varnothing}\in\Tau$. If $\|\dom[A]\|=1$, then $A=\{ab\}$ for some pair $ab\in X^2$. Since $X$ is a translation affine space, there exists a translation $T\in\Trans(X)$ such that $T(a)=b$. Then $A=T{\restriction}_{\{a\}}\in\Tau$. If $\|\dom[A]\|=2$, then $A$ is a line translation. By Theorem~\ref{t:paraD<=>translation|}, there exists a translation $T\in\Trans(X)$ such that $A=T{\restriction}_{\dom[A]}\in\Tau$. Therefore, $\I_X^{\#}\subseteq \Tau$. To see that $\Tau\subseteq\I^{\#}_X$, take any element $A\in\Tau$. If $A=1_X$, then $A\in\I_X^{\#}$ because $\I_X^{\#}$ is a submonoid of the symmetric inverse monoid $\I_X$. If $A\ne 1_X$, then $A=T{\restriction}_L$ for some translation $T\in\Trans(X)$ and some flat $L\subseteq X$ of rank $\|L\|\le 2$. If $\|L\|\le 1$, the the injection $T{\restriction}_L$ belongs to $\I_X^{\#}$, by Theorem~\ref{t:singleton-bijection}(1). So, assume that $\|L\|=2$. Since $X$ is a translation affine space, $\|\{B(o):B\in\Trans(X)\}\|=\|X\|\ge 3$ for every point $o\in X$. By Proposition~\ref{p:translation|=>line-translation}, the restriction $A=T{\restriction}_L$ is a line translation and hence $A\in\I^{\#}_X$.
\smallskip

$(2)\Ra(3)$ Assuming that $\I^{\#}_X=\Tau$, we shall prove that every subgroup $G$ of the monoid $\I^{\#}_X$ is commutative. Let $e$ be the idempotent of the group $G$. Since $e\in G\subseteq \I_X^{\#}\subseteq \mathcal{FI}_X$, the idempotent $e$ is the identity transformation  $1_L$ of some flat $L$ in $X$. If $L=X$ or $\|L\|\le 1$, then the group $G=\{e\}=\{1_L\}$ is trivial and hence commutative. So, assume that $\|L\|\le 2$ and $L\ne X$. Since the monoid $\I_X^{\#}$ is generated by line shifts, $\|L\|=2$. Since $\|X\|\ge 3$, there exist three points $o,a,b\in X$ such that $\|\{o,a,b\}\|=3$. By Corollary~\ref{c:parashift-exists}, there exist line shifts $S_a,S_b\in\I^{\#}_X$ such that $S_a(o)=a$ and $S_b(o)=b$. Since $S_a,S_b\in\I_X^{\#}=\Tau$, there exist translations $T_a,T_b\in\Trans(X)$ such that $S_a=T_a{\restriction}_{\dom[S_a]}$ and $S_b=T_b{\restriction}_{\dom[S_b]}$. Then $\|\{T(o):T\in\Trans(X)\|\ge\|\{1_X(o),T_a(o),T_b(o)\}\|=\|\{o,a,b\}\|=3$. By Theorem~\ref{t:Trans-commutative}, the group $\Trans(X)$ is commutative. For every element $g\in G$, we have $g=ege=1_Lg1_L$, which implies that $\dom[g]=\rng[g]=L$. Since $g\in G\subseteq\I_X^{\#}\subseteq \Tau$, there exists a translation $T_g\in\Trans(X)$ such that $g=T_g{\restriction}_L=T_g1_L$. It follows that $1_LT_g1_L=1_Lg=g=T_g1_L$. By the commutativity of the group $\Trans(X)$, for every elements $g,h\in G$ we have $$
\begin{aligned}
gh&=(T_g1_L)(T_h1_L)=T_g(1_LT_h1_L)=T_g(T_h1_L)=(T_gT_h)1_L\\
&=(T_hT_g)1_L=T_h(T_g1_L)=T_h(1_LT_g1_L)=(T_h1_L)(T_g1_L)=hg,
\end{aligned}
$$
witnessing that the group $G$ is commutative.
\smallskip

$(3)\Ra(4)$ Assume that every subgroup of the monoid $\I^{\#}_X$ is commutative. 
Let $A,B,C,D$ be line shifts in $X$ such that $\dom[A]=\rng[B]=\dom[C]=\rng[D]$.  Since $A,B,C,D$ are line shifts of the (regular) affine space $X$, the lines $\rng[A]$, $\dom[B]$, $\rng[C]$, $\dom[D]$ are parallel to the line $L\defeq \dom[A]=\rng[B]=\dom[C]=\rng[D]$. If $\rng[A]\cap\dom[B]=\varnothing$, then $BA=\varnothing$ and $$(DC)(BA)=(DC)\varnothing=\varnothing=\varnothing(DC)=(BA)(DC).$$ By analogy, $\rng[C]\cap\dom[D]=\varnothing$ implies $(DC)(BA)=(BA)(DC)$. So, assume that $\rng[A]\cap\dom[B]\ne\varnothing\ne\rng[C]\cap\dom[D]$. In this case $\rng[A]=\dom[B]$, $\rng[C]=\dom[D]$, and the line translations $BA$ and $DC$ are elements of the group $G\defeq\{g\in\I^{\#}_X:g1_L=g=1_Lg\}$. By our assumption, the group $G$ is commutative and hence $(DC)(BA)=(BA)(DC)$.  
\smallskip

$(4)\Ra(1)$ Assume that $(DC)(BA)=(BA)(DC)$ for every line shifts $A,B,C,D$ in $X$ with $\dom[A]=\rng[B]=\dom[C]=\rng[D]$. 
To show that $X$ is a translation affine space, fix any points $x,y\in X$. Choose any line $L\subseteq X$ containing the points $x,y$ and any line $L'\subseteq X$ such that $L'\parallel L$ and $L'\cap L=\varnothing$. Choose any point $x'\in L'$.  By Corollary~\ref{c:parashift-exists}, there exist line shifts $A:L\to L'$ and $B:L'\to L$ such that $A(x)=x'$ and $B(x')=y$. Our assumption ensures that for every line $\Lambda$ in $X$ and line shifts $C:L\to\Lambda$, $D:\Lambda\to L$ we have the equality $(BA)(DC)=(DC)(BA)$. By Theorem~\ref{t:Gleason1}, there exists a translation $T\in\Trans(X)$ such that $BA=T{\restriction}_L$. Then $T(x)=BA(x)=y$, witnessing that $X$ is a translation affine space. 
\end{proof}

\begin{theorem} For an affine space $X$, the following conditions are equivalent:
\begin{enumerate}
\item $X$ is Thalesian;
\item for every line $L$ in $X$, the group $\Sym^\#_X(L)$ is Abelian;
\item for some concurrent lines $L,\Lambda$ in $X$ the groups $\Sym^\#_X(L)$ and $\Sym^\#_X(\Lambda)$ are Abelian.
\end{enumerate}
\end{theorem}

\begin{proof} $(1)\Ra(2)$ Assume that the affine space $X$ is Thalesian. By Theorem~\ref{t:grp-line-trans}(2), for every line $L$ in $X$, the group $\Sym^\#_X(L)$ is Abelian.
\smallskip

The implication $(2)\Ra(3)$ is trivial.
\smallskip

$(3)\Ra(1)$ assume that $L$ and $\Lambda$ are two concurrent lines in $X$ such that the groups  $\Sym^\#_X(L)$ and $\Sym^\#_X(\Lambda)$ are Abelian. We have to prove that the afffine space $X$ is Thalesian. If $\|X\|\ne 3$, then the affine space $X$ is Desarguesian and Thalesian, by Corollary~\ref{c:affine-Desarguesian} and Theorem~\ref{t:ADA=>AMA}. So, assume that $\|X\|=3$, which means that $X$ is a Playfair plane. Let $o$ be the unique common point of the concurrent lines $L,\Lambda$. By Theorem~\ref{t:paraD<=>translation}, it suffices to show that for every point $p\in X$ there exists a translation $T:X\to X$ such that $T(o)=p$.

\begin{claim}\label{cl:T(o)=p} For every point $p\in L\cup\Lambda$ there exists a translation $T:X\to X$ such that $T(o)=p$.
\end{claim}

\begin{proof} First assume that $p\in L$. Choose any line $A$ in $X$ such that $A\cap L=\varnothing$. Choose any point $a\in A$ and consider the directions $\boldsymbol u\defeq(\Aline oa)_\parallel$ and $\boldsymbol v\defeq(\Aline ap)_\parallel$ in the Playfair plane $X$. By definition of the group $\Sym^\#_X(L)$, the permutation $P\defeq \boldsymbol v_{L,A}\boldsymbol u_{A,L}$ of the line $L$ belongs to the group $\Sym^\#_X(L)$. Observe that $P(o)=\boldsymbol v_{L,A}\boldsymbol u_{A,L}(o)=\boldsymbol v_{L,A}(a)=p$. Since the group $\Sym^\#_X(L)$ is Abelian, we can apply Gleason's Theorem~\ref{t:Gleason1} and conclude that $P=T{\restriction}_L$ for some translation $T:X\to X$. Then $T(o)=P(o)=p$. By analogy we can find a translation $T:X\to X$ with $T(o)=p$ if $p\in\Lambda$.
\end{proof}

Now we are ready to prove that for every point $p\in X$ there exists a translation $T:X\to X$ such that $T(o)=p$. If $p\in L\cup \Lambda$, then the translation $T$ exists, by Claim~\ref{cl:T(o)=p}. So, assume that $p\notin L\cup\Lambda$. Since $X$ is a Playfair plane, there exist points $x\in L\setminus\{o\}$ and $y\in\Lambda\setminus\{o\}$ such that $oxpy$ is a parallelogram. By Claim~\ref{cl:T(o)=p}, there exist translations $T_x,T_y$ of the plane $X$ such that $T_x(o)=x$ and $T_y(o)=y$. Then $T=T_xT_y$ is a translation of $X$ such that $T(o)=p$, witnessing that the Playfair plane $X$ is translation. By Theorem~\ref{t:paraD<=>translation}, the translation Playfair plane $X$ is Thalesian.
\end{proof}

\section{$\partial$-Translation affine spaces}

\begin{definition} An affine space $X$ is defined to be \defterm{$\partial$-translation} if for every spread of parallel lines $\delta\in\partial X$, there exists a translation $T:X\to X$ such that $\{\Aline xy:xy\in T\}=\delta$.
\end{definition}

\begin{theorem}\label{t:partial-translation+prime=>translation} An affine space $X$ of prime order is translation if and only if it is $\partial$-translation.
\end{theorem}

\begin{proof} The ``only if'' part is trivial. To prove the ``if'' part, assume that the affine space $X$ is $\partial$-translation. To prove that $X$ is translation, we should  check that the group of translations $\Trans(X)$ acts transitively on $X$. Given any distinct points $a,b\in X$, we should find a translation $T\in\Trans(X)$ such that $T(a)=b$. Consider the line $L\defeq \Aline ab$ and the subgroup $G\defeq\{T\in\Trans(X):T[L]=L\}$ of the group $\Trans(X)$. Since the affine space $X$ is $\partial$-translation, the group $G$ is not trivial and hence $|G|>1$. 

Given any point $x\in L$, consider the map $\alpha_x:G\to L$, $\alpha_x:T\mapsto T(x)$. Proposition~\ref{p:Ax=Bx=>A=B} ensures that the map $\alpha_x$ is injective and hence $|\alpha_x[G]|=|G|$. Observe that the family $\{\alpha_x[G]:x\in L\}$ is a partition of the line $L$ into subsets of cardinality $|G|$. This implies that $|G|$ divides $|L|$. Taking into account that the order $|X|_2=|L|$ of the affine space $X$ is prime and $|G|>1$, we conclude that $|G|=|L|$ and hence $b\in L=\alpha_a[G]$ and $b=T(a)$ for some translation $T\in G\subseteq\Trans(X)$, witnessing that $X$ is a translation affine plane.
\end{proof}

\begin{problem} Is every $\partial$-translation finite Playfair plane translation?
\end{problem}

\begin{proposition} If an affine space $X$ is $\partial$-translation, then
\begin{enumerate}
\item the translation group $\Trans(X)$ of $X$ is commutative and has cardinality $\Trans(X)\ge 2+|X|_2$.
\item If $\Trans(X)$ contains a non-identity element of finite order, then $\Trans(X)$ is elementary Abelian.
\end{enumerate}
\end{proposition} 

\begin{proof} Assume that an affine space $X$ is $\partial$-translation. 
Fix any point $o\in X$ and let $\mathcal L_o$ be the family of lines in $X$ that contain the point $o$. Since the affine space $X$ is Playfair, $|\mathcal L_o|=1+|X|_2$. Since $X$ is $\partial$-translation, for every line $L\in\mathcal L_o$ the intersection $L\cap\{T(o):T\in\Trans(X)\}$ contains at least one point, distinct from the point $o$. This implies $|\Trans(X)|\ge 1+|\mathcal L_o|=2+|X|_2$ and $\|\{T(o):T\in\Trans(X)\}\|=\|X\|\ge 3$.  By Baer's Theorem~\ref{t:Trans-commutative}, the group $\Trans(X)$ is commutative. Moreover, if $\Trans(X)$ contains a non-identity element of finite order, then $\Trans(X)$ is elementary Abelian.
\end{proof}


Let us recall that for a line $L$ in a liner $X$, we denote by $\Sym^\#_X(L)$ the group of all line translation $T:L\to L$. For a magma $(M,\cdot)$ its \defterm{center} is the set $$Z(M)\defeq\bigcap_{x\in M}\{z\in M:x\cdot z=z\cdot x\}$$of all elements of $M$ that commute with all other elements of $M$. Observe that a magma $M$ is comutative if and only if $M=Z(M)$. If $G$ is a group, the its centre $Z(G)$ is a normal commutative subgroup of $G$. 
 
Gleason's Theorem~\ref{t:Gleason1} implies the following group characterizations of translation and $\partial$-translation affine spaces.

\begin{corollary}\label{c:partialT<=>} An affine space $X$ is
\begin{enumerate}
\item translation if and only if for every line $L\subseteq X$, the group $\Sym^\#_X(L)$ is commutative;
\item  $\partial$-translation if and only if for every line $L\subseteq X$ the group $\Sym_X^\#(L)$ has non-trivial center.
\end{enumerate}
\end{corollary}

Corollary~\ref{c:partialT<=>}(2) motivates the problem of detecting groups with non-trivial center. A classical result in this direction is due to \index[person]{Burnside}Burnside\footnote{{\bf William Burnside} (1852 -- 1927) was an English mathematician,  known mostly as an early researcher in the theory of finite groups. 
\newline
Burnside was born in London in 1852. He went to school at Christ's Hospital until 1871 and attended St. John's and Pembroke Colleges at the University of Cambridge, where he was the Second Wrangler (bracketed with George Chrystal) in 1875. He lectured at Cambridge for the following ten years, before being appointed professor of mathematics at the Royal Naval College in Greenwich. While this was a little outside the main centres of British mathematical research, Burnside remained a very active researcher, publishing more than 150 papers in his career.
\newline
Burnside's early research was in applied mathematics. This work was of sufficient distinction to merit his election as a fellow of the Royal Society in 1893, though it is little remembered today. Around the same time as his election his interests turned to the study of finite groups. This was not a widely studied subject in Britain in the late 19th century, and it took some years for his research in this area to gain widespread recognition.
\newline 
The central part of Burnside's group theory work was in the area of group representations, where he helped to develop some of the foundational theory, complementing, and sometimes competing with, the work of Ferdinand Georg Frobenius, who began his research in the subject during the 1890s. One of Burnside's best known contributions to group theory is his $p^aq^b$ theorem, which shows that every finite group whose order is divisible by fewer than three distinct primes is solvable.
\newline
In 1897 Burnside's classic work ``Theory of Groups of Finite Order'' was published. The second edition (pub. 1911) was for many decades the standard work in the field. A major difference between the editions was the inclusion of character theory in the second.
\newline
Burnside is also remembered for the formulation of Burnside's problem that concerns the question of bounding the size of a group if there are fixed bounds both on the order of all of its elements and the number of elements needed to generate it, and also for Burnside's lemma (a formula relating the number of orbits of a permutation group acting on a set with the number of fixed points of each of its elements) though the latter had been discovered earlier and independently by Frobenius and Augustin Cauchy.
\newline
He received an honorary doctorate (D.Sc.) from the University of Dublin in June 1901. 
\newline
In addition to his mathematical work, Burnside was a noted rower. While he was a lecturer at Cambridge, he also coached the rowing crew team. In fact, his obituary in ``The Times'' took more interest in his athletic career, calling him ``one of the best known Cambridge athletes of his day''.
\newline
He is buried at the West Wickham Parish Church in South London.}.

A number $n\in\IN$ is called a \defterm{prime power}\index{prime power} if $n=p^k$ for some prime number $p$ and some number $k\in\IN$. 

\begin{theorem}[Burnside, 1897]\label{t:Burnside-center} A finite group $G$ has non-trivial center if the order of $G$ is a prime power.
\end{theorem}

\begin{proof} Assume that $|G|=p^n$ for some prime number $p$ and some $n$. For every $g\in G$, consider the function $\gamma_g:G\to G$, $\gamma_g:x\mapsto xgx^{-1}$, and observe that it is a group homomorphism and hence its kernel $Z_G(g)\defeq\{x\in G:g=xgx^{-1}\}$ is a subgroup of $G$  and hence the cardinal $|Z_G(g)|$ divides the cardinality $|G|=p^n$ of the group $G$ and the cardinal $|\gamma_g[G]|=\frac{|G|}{|Z_G(g)|}$ also divides $|G|=p^n$. The image $\gamma_g[G]$ is the orbit of the element $g$ under the conjugating action of the group $G$ on $G$. Observe that an element $c\in G$ belongs to the center $Z(G)$ of $G$ if and only if its orbit $\gamma_c[G]$ is a singleton if and only if $|\gamma_c[G]|$ is not divisible by $p$. 

It is easy to see that for two elements $a,b\in G$, the sets $\gamma_a[G]$ and $\gamma_b[G]$ are either disjoint or coincide. Then $\{\gamma_g[G]:g\in G\}$ is a partition of the group $G$ into pairwise disjoint sets. Choose a set $A\subseteq G$ such that $A\cap \gamma_g[G]$ is a singleton for every element $g\in G$. Then 
$$p^n=|G|=\sum_{a\in A}|\gamma_a[G]|=|Z(G)|+\sum_{a\in A\setminus Z(G)}|\gamma_a[G]|$$
and hence the cardinal $|Z(G)|$ is divisible by $p$ (because for every $a\in A\setminus Z(G)|$ the cardinal $|\gamma_g[G]|$ is divisible by $p$. Taking into account that $Z(G)$ is a subgroup if $G$, we conclude that $|Z(G)|=p^k$ for some $k\ge 1$ (because $|Z(G)|$ is divisible by $p$ and $|Z(G)|$ divides $|G|=p^n$).
\end{proof}


 
\section{The group of line translations}\label{s:IX[L;u,v]}

In this section we prove some basic results on the structure of the group $\Sym^\#_X(L)$ of line translations of a line $L$ in a Playfair plane $X$.
First, we describe generators of the group $\Sym^\#_X(L)$. We say that a group $G$ is \defterm{generated by a set $S\subseteq G$} if $G$ coincides with the smallest subgroup of $G$ that contains the set $S$.

Given any line $\Lambda\in L_\parallel$ and any direction ${\boldsymbol \delta}\in \partial X\setminus\{L_\parallel\}$, consider the line shift ${\boldsymbol \delta}_{\Lambda,L}:L\to\Lambda$ assigning to every point $x\in L$ the unique point $y\in \Lambda\cap \Aline x{\boldsymbol \delta}$ (the point $y$ exists because $X$ is a Playfair plane and $L\notin{\boldsymbol \delta}$).

For two directions ${\boldsymbol u},{\boldsymbol v}\in\partial X\setminus\{L_\parallel\}$, consider the subset
$$\Sym_X^\#[L;\boldsymbol u,\boldsymbol v]=\{{\boldsymbol u}_{L,\Lambda} {\boldsymbol v}_{\Lambda,L}:\Lambda\in L_\parallel\}$$
of the group $\Sym_X^\#(L)$.

\begin{exercise} Show that $\Sym_X^\#[L;\boldsymbol u,\boldsymbol v]^{-1}=\Sym_X^\#[L;{\boldsymbol v},{\boldsymbol u}]$, where $A^{-1}\defeq\{a^{-1}:a\in A\}$ for a subset $A$ of a group $G$.
\end{exercise}

\begin{proposition}\label{p:|IX|=|L|} For every line $L$ in a Playfair plane $X$, any distinct  directions ${\boldsymbol u},{\boldsymbol v}\in\partial X\setminus\partial L_\parallel$ and any  point $o\in L$, the function $\alpha:\Sym_X^\#[L;\boldsymbol u,\boldsymbol v]\to L$, $\alpha:T\mapsto T(o)$, is bijective and hence $|\Sym_X^\#[L;\boldsymbol u,\boldsymbol v]|=|L|=|X|_2$.
\end{proposition}

\begin{proof} Given any point $y\in L$, find a unique point $z\in \Aline o{\boldsymbol v}\cap\Aline y{\boldsymbol u}$ (the point $z$ exists because $\boldsymbol u$ and $\boldsymbol v$ are two distinct directions in the Playfair plane $X$). Next find a unique line $\Lambda\in L_\parallel$ containing the points $z$ and observe that ${\boldsymbol u}_{L,\Lambda}{\boldsymbol v}_{\Lambda,L}$ is a unique line translation $T\in \Sym_X^\#[L;\boldsymbol u,\boldsymbol v]$ such that $T(o)=y$. This shows that the function $\alpha$ is bijective and hence $|\Sym_X^\#[L;\boldsymbol u,\boldsymbol v]|=|L|=|X|_2$.
\end{proof}

\begin{lemma}\label{l:IX-generators} For every line $L$ in an affine space $X$, the group $\Sym_X^\#(L)$ is generated by the set $$\bigcup_{{\boldsymbol u},{\boldsymbol v}\in\partial X\setminus\{L_\parallel\}}\Sym_X^\#[L;\boldsymbol u,\boldsymbol v].$$
\end{lemma}

\begin{proof} By definition, line translations $T\in \Sym^\#_X(L)$ are compositions of line shifts. Therefore, it suffices to prove that for every $n\in\IN$ and line shifts $S_0,\dots,S_n$ in $X$ with $S_n\cdots S_0\in \Sym_X^\#(L)$, the line translation $S_n\cdots S_0$ can be written as the product $T_n\cdots T_1$ of line translations $T_1,\dots, T_n\in \bigcup_{{\boldsymbol u},{\boldsymbol v}\in\partial X\setminus\{L_\parallel\}}\Sym_X^\#[L;\boldsymbol u,\boldsymbol v]$. 

For $n=1$, the inclusion $S_n\cdots S_0\in\Sym^\#_X(L)$ implies $\dom[S_0]=L$, $\dom[S_0]=\rng[S_1]$ and $\rng[S_1]=L$. Consider the line $\Lambda\defeq\rng[S_0]=\dom[S_1]$. Since $S_0,S_1$ are line shifts, $\Lambda\in L_\parallel$ and $S_0={\boldsymbol u}_{\Lambda,L}$, $S_1={\boldsymbol v}_{L,\Lambda}$ for some directions ${\boldsymbol u},{\boldsymbol v}\in\partial X\setminus\{L_\parallel\}$. Then $S_1\cdot S_0=T_1\defeq {\boldsymbol v}_{L,\Lambda}{\boldsymbol u}_{\Lambda,L}\in\Sym_X^\#[L;{\boldsymbol u},{\boldsymbol v}]$ and we are done.

Now assume that for some $n\ge 2$ we know that for every $k<n$ and every line shifts $S_0,\dots,S_k$ with $S_k\cdots S_0\in \Sym_X^\#(L)$, the line translation $S_k\cdots S_0$ can be written as the product $T_k\cdots T_1$ of line translations $T_1,\dots, T_k\in \bigcup_{{\boldsymbol u},{\boldsymbol v}\in\partial X\setminus\{L_\parallel\}}\Sym_X^\#[L;\boldsymbol u,\boldsymbol v]$. Take any line shifts $S_0,\dots, S_n$ with $S_n\cdots S_0\in\Sym^\#_X(L)$. Since $S_0, S_{1}$ are line shifts, the lines $L_0\defeq\rng[S_0]=\dom[S_1]$ and $L_1\defeq\rng[S_1]$ are parallel to the line $\dom[S_0]=L$. Moreover,  $S_0={\boldsymbol v}_{L_0,L}$ and $S_1={\boldsymbol u}_{L_1,L_0}$ for some directions ${\boldsymbol u},{\boldsymbol v}\in\partial X\setminus\{L_\parallel\}$. Observe that $S_1'\defeq S_1{\boldsymbol u}_{L_0,L}={\boldsymbol u}_{L_1,L_0}{\boldsymbol u}_{L_0,L}={\boldsymbol u}_{L_1,L}$ is a line shift such that $S_n\cdots S_2S_1'\in \Sym^\#_X(L)$. By the inductive assumption, there exist line translations $T_2,\dots,T_n\in \bigcup_{{\boldsymbol u},{\boldsymbol v}\in\partial X\setminus\{L_\parallel\}}\Sym_X^\#[L;\boldsymbol u,\boldsymbol v]$ such that $S_n\cdots S_2S_1'=T_n\cdots T_2$. Then for the line translation $T_1\defeq {\boldsymbol u}_{L,L_0}{\boldsymbol v}_{L_0,L}\in \Sym_X^\#[L;{\boldsymbol u},{\boldsymbol v}]$, we have
\begin{multline*}
S_n\cdots S_2S_1S_0=S_n\cdots S_2{\boldsymbol u}_{L_0,L_1}{\boldsymbol v}_{L,L_0}=S_n\cdots S_2{\boldsymbol u}_{L_1,L_0}1_{L_0}{\boldsymbol v}_{L_0,L}\\
=S_n\cdots S_2{\boldsymbol u}_{L_1,L_0}{\boldsymbol u}_{L_0,L}{\boldsymbol u}_{L,L_0}{\boldsymbol v}_{L_0,L}=S_n\cdots S_2S_1'T_1=T_n\cdots T_2T_1.
\end{multline*}
\end{proof}

For two subsets $A,B$ of a group $G$ we denote by 
$$A\cdot B\defeq\{a\cdot b:(a,b)\in A\times B\}$$their pointwise product in $G$.

\begin{lemma}\label{IXuv=IXud+IXdv}  For every line $L$ in a Playfair plane $X$ and directions ${\boldsymbol u},{\boldsymbol v},{\boldsymbol \delta}\in\partial X\setminus\{L_\parallel\}$, 
$$\Sym_X^\#[L;\boldsymbol u,\boldsymbol v]\subseteq \Sym_X^\#[L;{\boldsymbol u},{\boldsymbol \delta}]\cdot \Sym_X^\#[L;{\boldsymbol \delta},{\boldsymbol v}].$$
\end{lemma}

\begin{proof} Given any element $T\in \Sym_X^\#[L;\boldsymbol u,\boldsymbol v]$, find a line $\Lambda\in L_\parallel\setminus\{L\}$ such that $T={\boldsymbol u}_{L,\Lambda}{\boldsymbol v}_{\Lambda,L}$. Taking into account that the composition ${\boldsymbol \delta}_{\Lambda,L}{\boldsymbol \delta}_{L,\Lambda}:\Lambda\to \Lambda$ is the identity map of the line $\Lambda$, we conclude that
$$
T={\boldsymbol u}_{L,\Lambda}{\boldsymbol v}_{\Lambda,L}=
{\boldsymbol u}_{L,\Lambda}({\boldsymbol \delta}_{\Lambda,L}{\boldsymbol \delta}_{L,\Lambda}) {\boldsymbol v}_{\Lambda,L}
=({\boldsymbol u}_{L,\Lambda}{\boldsymbol \delta}_{\Lambda,L})\cdot({\boldsymbol \delta}_{L,\Lambda} {\boldsymbol v}_{\Lambda,L})\in
\Sym^\#_X[L;{\boldsymbol u},{\boldsymbol \delta}]\cdot\Sym_X^\#[L;{\boldsymbol \delta},{\boldsymbol v}].
$$
\end{proof}

Lemmas~\ref{l:IX-generators} and \ref{IXuv=IXud+IXdv} imply the following corollary.

\begin{corollary}\label{c:IX-generators} For every line $L$ in a Playfair plane $X$ and every direction ${\boldsymbol \delta}\in\partial X\setminus\{L_\parallel\}$, the group $\Sym_X^\#(L)$ is generated by the set $$\bigcup_{\boldsymbol u\in\partial X\setminus\{L_\parallel\}}\big(\Sym_X^\#[L;{\boldsymbol u},{\boldsymbol \delta}]\cup\Sym_X^\#[L;{\boldsymbol \delta},{\boldsymbol u}]\big).$$
\end{corollary}

\begin{lemma}\label{l:IX-commute} Let $L$ be a line in a Playfair plane $X$ and ${\boldsymbol u},{\boldsymbol v},{\boldsymbol \delta}\in\partial X\setminus \{L_\parallel\}$ be directions such that the sets $\Sym^\#_X[L;{\boldsymbol v},{\boldsymbol \delta}],\Sym^\#_X[L;{\boldsymbol \delta},{\boldsymbol u}],\Sym^\#_X[L;{\boldsymbol v},{\boldsymbol u}]$ are subgroups of the group $\Sym_X^\#(L)$. Then 
$$\Sym_X^\#[L;{\boldsymbol u},{\boldsymbol \delta}]\cdot\Sym_X^\#[L;{\boldsymbol \delta},{\boldsymbol v}]= \Sym_X^\#[L;{\boldsymbol \delta},{\boldsymbol v}]\cdot\Sym_X^\#[L;{\boldsymbol u},{\boldsymbol \delta}].$$ 
\end{lemma}

\begin{proof} Given any elements $\alpha\in \Sym_X^\#[L;\boldsymbol u,\boldsymbol \delta]$ and $\beta\in\Sym_X^\#[L;\boldsymbol \delta,\boldsymbol v]$, we shall prove that $\beta\alpha \in \Sym_X^\#[L;\boldsymbol u,\boldsymbol \delta]\cdot\Sym_X^\#[L;\boldsymbol \delta,\boldsymbol v]$. Indeed, the definition of the sets $\Sym_X^\#[L;\boldsymbol u,\boldsymbol \delta]$ and $\Sym_X^\#[L;\boldsymbol \delta,\boldsymbol v]$ ensures that $\alpha={\boldsymbol u}_{L,A}{\boldsymbol \delta}_{A,L}$ and  $\beta={\boldsymbol \delta}_{L,B}{\boldsymbol v}_{B,L}$ for some lines $A,B\in L_\parallel$. Since $\Sym_X^\#[L;{\boldsymbol u},{\boldsymbol v}]$ is a group, there exists a line $C\in L_\parallel$ such that $({\boldsymbol u}_{L,B}{\boldsymbol v}_{B,L})({\boldsymbol u}_{L,A}{\boldsymbol v}_{A,L})={\boldsymbol u}_{L,C}{\boldsymbol v}_{C,L}$. 
Then $${\boldsymbol u}_{L,B}{\boldsymbol \delta}_{B,L}({\boldsymbol \delta}_{L,B}{\boldsymbol v}_{B,L}{\boldsymbol u}_{L,A}{\boldsymbol \delta}_{A,L}){\boldsymbol \delta}_{L,A}{\boldsymbol v}_{A,L}={\boldsymbol u}_{L,B}{\boldsymbol v}_{B,L}{\boldsymbol u}_{L,A}{\boldsymbol v}_{A,L}={\boldsymbol u}_{L,C}{\boldsymbol v}_{C,L}={\boldsymbol u}_{L,C}{\boldsymbol \delta}_{C,L}{\boldsymbol \delta}_{L,C}{\boldsymbol v}_{C,L}$$ and hence
$$
\begin{aligned}
\beta\alpha&={\boldsymbol \delta}_{L,B}{\boldsymbol v}_{B,L}{\boldsymbol u}_{L,A}{\boldsymbol \delta}_{A,L}=({\boldsymbol u}_{L,B}{\boldsymbol \delta}_{B,L})^{-1}{\boldsymbol u}_{L,C}{\boldsymbol \delta}_{C,L}{\boldsymbol \delta}_{L,C}{\boldsymbol v}_{C,L}({\boldsymbol \delta}_{L,A}{\boldsymbol v}_{A,L})^{-1}\\
&\in \Sym_X^\#[L;\boldsymbol u,\boldsymbol \delta]\cdot\Sym_X^\#[L;\boldsymbol \delta,\boldsymbol v]
\end{aligned}
$$ because $\Sym_X^\#[L;\boldsymbol u,\boldsymbol \delta]$ and $\Sym_X^\#[L;\boldsymbol \delta,\boldsymbol v]$ are groups. Therefore,
$$ \Sym_X^\#[L;\boldsymbol \delta,\boldsymbol v]\cdot\Sym_X^\#[L;\boldsymbol u,\boldsymbol \delta]\subseteq \Sym_X^\#[L;\boldsymbol u,\boldsymbol \delta]\cdot\Sym_X^\#[L;\boldsymbol \delta,\boldsymbol v]$$ and after inversion,
\begin{multline*}
\Sym_X^\#[L;\boldsymbol u,\boldsymbol \delta]\cdot \Sym_X^\#[L;\boldsymbol \delta,\boldsymbol v]
=\Sym_X^\#[L;\boldsymbol u,\boldsymbol \delta]^{-1}\cdot \Sym_X^\#[L;\boldsymbol \delta,\boldsymbol v]^{-1}\\
=
(\Sym_X^\#[L;\boldsymbol \delta,\boldsymbol v]\cdot\Sym_X^\#[L;\boldsymbol u,\boldsymbol \delta])^{-1}
\subseteq \big(\Sym_X^\#[L;\boldsymbol u,\boldsymbol \delta]\cdot\Sym_X^\#[L;\boldsymbol \delta,\boldsymbol v]\big)^{-1}\\ =\Sym_X^\#[L;\boldsymbol \delta,\boldsymbol v]^{-1}\cdot\Sym_X^\#[L;\boldsymbol u,\boldsymbol \delta]^{-1}=\Sym_X^\#[L;\boldsymbol \delta,\boldsymbol v]\cdot\Sym_X^\#[L;\boldsymbol u,\boldsymbol \delta]
\end{multline*}
because $\Sym_X^\#[L;\boldsymbol u,\boldsymbol \delta]$ and $\Sym_X^\#[L;{\boldsymbol \delta},{\boldsymbol v}]$ are groups and hence $\Sym_X^\#[L;\boldsymbol u,\boldsymbol \delta]=\Sym_X^\#[L;\boldsymbol u,\boldsymbol \delta]^{-1}$ and $\Sym_X^\#[L;{\boldsymbol \delta},{\boldsymbol v}]=\Sym_X^\#[L;{\boldsymbol \delta},{\boldsymbol v}]^{-1}$.
Thus we obtaine the desired equality 
$$ \Sym_X^\#[L;\boldsymbol \delta,\boldsymbol v]\cdot\Sym_X^\#[L;\boldsymbol u,\boldsymbol \delta]=\Sym_X^\#[L;\boldsymbol u,\boldsymbol \delta]\cdot\Sym_X^\#[L;\boldsymbol \delta,\boldsymbol v].$$
\end{proof}

\begin{proposition}\label{p:IX=n**(n-1)} Let $L$ be a line in a Playfair plane $X$ such that for every directions ${\boldsymbol u},{\boldsymbol v}\in\partial X\setminus\{L_\parallel\}$, the set $\Sym^\#_X[L;\boldsymbol u,\boldsymbol v]$ is a subgroup of the group $\Sym_X^\#(L)$. If the cardinal $n\defeq |X|_2$ is finite, then the cardinality of the group $\Sym_X^\#(L)$ divides the the number $n^{n-1}$.
\end{proposition}

\begin{proof} Fix any direction ${\boldsymbol \delta}\in\partial X\setminus\{L_\parallel\}$. Since the cardinal $n\defeq|X|_2$ is finite, the horizon $\partial X$ is finite of cardinality $|X|_2+1=n+1$. Then $\partial X=\{{\boldsymbol \delta}_0,\dots,{\boldsymbol \delta}_{n}\}$ for some directions ${\boldsymbol \delta}_0,\dots,{\boldsymbol \delta}_n$ such that ${\boldsymbol \delta}_0=L_\parallel$ and ${\boldsymbol \delta}_1=\delta$. By our assumption, for every $i\in\{1,\dots, n\}$, the set $\Sym_X^\#[L;{\boldsymbol \delta},{\boldsymbol \delta}_i]$ is a subgroup of the group $\Sym_X^\#(L)$. Then $\Sym_X^\#[L,{\boldsymbol \delta}_i,{\boldsymbol \delta}]=\Sym_X^\#[L;{\boldsymbol \delta}_i,{\boldsymbol \delta}]^{-1}=\Sym_X^\#[L;{\boldsymbol \delta},{\boldsymbol \delta}_i]$. 

If $i=1$, then the set $\Sym_X^\#[L;{\boldsymbol \delta},{\boldsymbol \delta}_1]$ is a singleton containing the identity map of the line $L$. For $i\in\{2,\dots,n\}$, the group $\Sym_X^\#[L;{\boldsymbol \delta},{\boldsymbol \delta}_i]$ has cardinality $|L|=|X|_2=n$, by Proposition~\ref{p:|IX|=|L|}. By Corollary~\ref{c:IX-generators}, the group $\Sym_X^\#(L)$ is generated by the set $\bigcup_{i=2}^n\Sym_X^\#[L;{\boldsymbol \delta},{\boldsymbol \delta}_i]$. By Lemma~\ref{l:IX-commute}, for any numbers $i,j\in\{1,\dots,n\}$ the subgroups $\Sym_X^\#[L;{\boldsymbol \delta},{\boldsymbol \delta}_i]$ and $\Sym_X^\#[L;{\boldsymbol \delta},{\boldsymbol \delta}_j]$ mutually commute in the sense that $$\Sym_X^\#[L;{\boldsymbol \delta},{\boldsymbol \delta}_i]\cdot\Sym_X^\#[L;{\boldsymbol \delta},{\boldsymbol \delta}_j]=
\Sym_X^\#[L;{\boldsymbol \delta},{\boldsymbol \delta}_j]\cdot\Sym_X^\#[L;{\boldsymbol \delta},{\boldsymbol \delta}_i].$$

Consider the sequence of the sets $(H_k)_{k=1}^n$ defined by the recursive formula $H_1=\Sym_X^\#[L;{\boldsymbol \delta},{\boldsymbol \delta}_1]=\{1_L\}$ and $H_{k+1}=H_k\cdot \Sym_X^\#[L;{\boldsymbol \delta},{\boldsymbol \delta}_k]$ for $k\in\{1,\dots,n-1\}$. The mutual commutativity of the subgroups $\Sym_X^\#[L;{\boldsymbol \delta},{\boldsymbol \delta}_i]$ for $i\in\{1,\dots,n\}$ implies that every set $H_i$ is a subgroup of the group $\Sym_X^\#(L)$ and $H_n=\Sym^\#_X(L)$. It is easy  to see that for every $i\in\{2,\dots,n\}$ the index\footnote{The {\em index} of a subgroup $H$ is a group $G$ is the cardinality of the set $\{Hx:x\in G\}$ of all cosets $Hx\defeq \{hx:h\in H\}$.} of the subgroup $H_{i-1}$ in $H_i$ is equal to the index of the subgroup $H_{i-1}\cap\Sym_X^\#[L;{\boldsymbol \delta},{\boldsymbol \delta}_i]$ in the group $\Sym_X^\#[L;{\boldsymbol \delta},{\boldsymbol \delta}_i]$ and hence the index of the group $H_{i-1}$ in $H_i$ divides the number $n=|\Sym^\#_X[L;{\boldsymbol \delta},{\boldsymbol \delta}_i]|$. By induction this implies that for every $i\in\{1,\dots,n\}$, the cardinality of the group $H_i$ divides the number $n^{i-1}$. In particular, the cardinality of the group $H_n=\Sym_X^\#(L)$ divides the number $n^{n-1}$.
\end{proof}

Proposition~\ref{p:IX=n**(n-1)} and Theorem~\ref{t:Burnside-center}  imply the following corollary.

\begin{corollary}\label{c:IX-center} Let $L$ be a line in a Playfair plane $X$ such that for every directions ${\boldsymbol u},{\boldsymbol v}\in\partial X\setminus\{L_\parallel\}$ the set $\Sym_X^\#[L;\boldsymbol u,\boldsymbol v]$ is a subgroup of the group $\Sym_X^\#(L)$. If $|X|_2$ is a prime power, then the group $\Sym_X^\#(L)$ has non-trivial center.
\end{corollary}

Corollaries~\ref{c:IX-center} and \ref{c:partialT<=>}  imply the following main result of this section.

\begin{theorem}[Gleason, 1956]\label{t:partial-translation<=} A Playfair plane $X$ is $\partial$-translation whenever $|X|_2$ is a prime power and for every line $L\subset X$ and distinct directions ${\boldsymbol u},{\boldsymbol v}\in\partial X\setminus\{L_\parallel\}$, the set $\Sym_X^\#[L;\boldsymbol u,\boldsymbol v]$ is a subgroup of the group $\Sym_X^\#(L)$. 
\end{theorem}

\begin{problem} Is every $\partial$-translation finite Playfair plane translation?
\end{problem}

We shall also need a modification of Gleason's Theorem~\ref{t:partial-translation<=} in which the groups $\Sym_X^\#[L;\boldsymbol u,\boldsymbol v]$ are replaced by groups 
$\Sym_X^\#[L;{\boldsymbol h},\Lambda]$, defined as follows.

For two parallel lines $L,\Lambda$ in a Playfair plane $X$ and a direction ${\boldsymbol h}\in\partial X\setminus \{L_\parallel\}$, consider the set
$$\Sym^\#_X[L;{\boldsymbol h},\Lambda]\defeq\big\{{\boldsymbol u}_{L,\Lambda}{\boldsymbol h}_{\Lambda,L}:{\boldsymbol u}\in \partial X\setminus\{L_\parallel\}\big\}.$$

\begin{proposition}\label{p:|IX|=|L|+} For any parallel lines $L,\Lambda$ in a Playfair plane $X$, any direction ${\boldsymbol h}\in\partial X\setminus L_\parallel$ and point $o\in L$, the function $\alpha:\Sym_X^\#[L;{\boldsymbol h},\Lambda]\to L$, $\alpha:T\mapsto T(o)$, is bijective and hence $|\Sym_X^\#[L;{\boldsymbol h},\Lambda]|=|L|=|X|_2$.
\end{proposition}

\begin{proof} Since $\Lambda\notin{\boldsymbol h}$, there exists a unique point $o'\in\Lambda\cap\Aline o{\boldsymbol h}$. Then for the direction $\boldsymbol u\defeq(\Aline {o'}y)_\parallel$, the line translation 
$T\defeq {\boldsymbol u}_{L,\Lambda}{\boldsymbol h}_{L,\Lambda}$ is a unique element of the set $\Sym_X^\#[L;{\boldsymbol h},\Lambda]$ such that $T(o)=y$. This shows that the function $\alpha$ is bijective and hence $|\Sym_X^\#[L;{\boldsymbol h},\Lambda]|=|L|=|X|_2$.
\end{proof}

\begin{lemma}\label{l:IXuv=IXd-IXd}  For every line $L$ in a Playfair plane $X$ and directions ${\boldsymbol u},{\boldsymbol v},{\boldsymbol h}\in\partial X\setminus\{L_\parallel\}$, 
$$\Sym_X^\#[L;\boldsymbol u,\boldsymbol v]\subseteq \bigcup_{\Lambda\in L_\parallel}\Sym_X^\#[L;{\boldsymbol h},\Lambda]\cdot\Sym_X^\#[L;{\boldsymbol h},\Lambda]^{-1}.$$
\end{lemma}

\begin{proof} Given any element $T\in \Sym_X^\#[L;\boldsymbol u,\boldsymbol v]$, find a line $\Lambda\in L_\parallel\setminus\{L\}$ such that $T={\boldsymbol u}_{L,\Lambda}{\boldsymbol v}_{\Lambda,L}$. Taking into account that the composition ${\boldsymbol h}_{\Lambda,L}{\boldsymbol h}_{L,\Lambda}:\Lambda\to \Lambda$ is the identity map of the line $\Lambda$, we conclude that
$$
\begin{aligned}
T&={\boldsymbol u}_{L,\Lambda}{\boldsymbol v}_{\Lambda,L}=
{\boldsymbol u}_{L,\Lambda}({\boldsymbol h}_{\Lambda,L}{\boldsymbol h}_{L,\Lambda}) {\boldsymbol v}_{\Lambda,L}
=({\boldsymbol u}_{L,\Lambda}{\boldsymbol h}_{\Lambda,L})\cdot({\boldsymbol h}_{L,\Lambda} {\boldsymbol v}_{\Lambda,L})\\
&=({\boldsymbol u}_{L,\Lambda}{\boldsymbol h}_{\Lambda,L})\cdot({\boldsymbol v}_{L,\Lambda} {{\boldsymbol h}}_{\Lambda,L})^{-1}
\in
\Sym^\#_X[L;{\boldsymbol h},\Lambda]\cdot\Sym_X^\#[L;{\boldsymbol h},\Lambda]^{-1}.
\end{aligned}
$$
\end{proof}

Lemmas~\ref{l:IX-generators} and \ref{l:IXuv=IXd-IXd} imply the following corollary.

\begin{corollary}\label{c:IX-generators+} For every line $L$ in a Playfair plane $X$ and every direction ${\boldsymbol h}\in\partial X\setminus\{L_\parallel\}$, the group $\Sym_X^\#(L)$ is generated by the set $$\bigcup_{\Lambda\in L_\parallel}\Sym_X^\#[L;{\boldsymbol h},\Lambda].$$
\end{corollary}

\begin{lemma}\label{l:IX-commute+} Let $L$ be a line in a Playfair plane $X$, $A,B\in L_\parallel$ be two lines and ${\boldsymbol h}\in\partial X\setminus \{L_\parallel\}$ be a direction such that the sets $\Sym^\#_X[L;{\boldsymbol h},A],\Sym^\#_X[L;{\boldsymbol h},B],\Sym^\#_X[B;{\boldsymbol h},A]$ are subgroups of the group $\Sym_X^\#(L)$. Then 
$$\Sym_X^\#[L;{\boldsymbol h},A]\cdot\Sym_X^\#[L;{\boldsymbol h},B]= \Sym_X^\#[L;{\boldsymbol h},B]\cdot\Sym_X^\#[L;{\boldsymbol h},A].$$ 
\end{lemma}

\begin{proof} Given any elements $\alpha\in \Sym_X^\#[L;{\boldsymbol h},A]$ and $\beta\in\Sym_X^\#[L;{\boldsymbol h},B]$, we first prove that $\alpha\beta \in \Sym_X^\#[L;{\boldsymbol h},B]\cdot\Sym_X^\#[L;{\boldsymbol h},A]$. The definition of the set $\Sym_X^\#[L;\boldsymbol u,\boldsymbol h]$ ensures that $\alpha={\boldsymbol a}_{L,A}{\boldsymbol h}_{A,L}$ for some direction $\boldsymbol a\in\partial X\setminus\{L_\parallel\}$. Since $\Sym_X^\#[L,{\boldsymbol h},B]$ is a group, $\beta^{-1}\in \Sym_X^\#[L,{\boldsymbol h},B]$ and hence $\beta^{-1}={\boldsymbol b}_{L,B}{{\boldsymbol h}}_{B,L}$ for some direction $\boldsymbol b\in\partial X\setminus\{L_\parallel\}$.
Then $\beta=(\boldsymbol b_{L,B}{\boldsymbol h}_{B,L})^{-1}={\boldsymbol h}_{B,L}^{-1}\boldsymbol b^{-1}_{L,B}={\boldsymbol h}_{L,B}\boldsymbol b_{B,L}$ and 
$$\alpha\beta={\boldsymbol a}_{L,A}{\boldsymbol h}_{A,L}{\boldsymbol h}_{L,B}\boldsymbol b_{B,L}=\boldsymbol a_{L,B}\boldsymbol a_{B,A}{\boldsymbol h}_{A,B}\boldsymbol b_{B,A}\boldsymbol b_{A,L}=\boldsymbol a_{L,B}\boldsymbol a_{B,A}{\boldsymbol h}_{A,B}\boldsymbol b_{B,A}{\boldsymbol h}_{A,B}{\boldsymbol h}_{B,A}\boldsymbol b_{A,L}$$
Since $\Sym^\#_X[B;{\boldsymbol h},A]$ is a group, there exists a direction $\boldsymbol c\in\partial X\setminus L_\parallel$ such that $$(\boldsymbol a_{B,A}{\boldsymbol h}_{A,B})(\boldsymbol b_{B,A}{\boldsymbol h}_{A,B})=\boldsymbol c_{B,A}{\boldsymbol h}_{A,B}.$$ Then
$$
\begin{aligned}
\alpha\beta&=\boldsymbol a_{L,B}(\boldsymbol c_{B,A}{\boldsymbol h}_{A,B}{\boldsymbol h}_{B,A})\boldsymbol b_{A,L}=
\boldsymbol a_{L,B}\boldsymbol c_{B,A}\boldsymbol b_{A,L}\\
&=\boldsymbol a_{L,B}\boldsymbol c_{B,L}\boldsymbol c_{L,A}\boldsymbol b_{A,L}=
(\boldsymbol a_{L,B}{\boldsymbol h}_{B,L}{\boldsymbol h}_{L,B}\boldsymbol c_{B,L})(\boldsymbol c_{L,A}{\boldsymbol h}_{A,L}{\boldsymbol h}_{L,A}\boldsymbol b_{A,L})\\
&=(\boldsymbol a_{L,B}{\boldsymbol h}_{B,L}(\boldsymbol c_{L,B}{\boldsymbol h}_{B,L})^{-1})(\boldsymbol c_{L,A}\delta_{A,L}(\boldsymbol b_{L,A}{\boldsymbol h}_{A,L})^{-1})\in \Sym^\#_X[L;{\boldsymbol h},B]\cdot\Sym^\#_X[L,{\boldsymbol h},A]
\end{aligned}
$$
because $\Sym_X^\#[L;{\boldsymbol h},B]$ and $\Sym_X^\#[L;{\boldsymbol h},A]$ are groups. Therefore,
$$ \Sym_X^\#[L;{\boldsymbol h},A]\cdot\Sym_X^\#[L;{\boldsymbol h},B]\subseteq \Sym_X^\#[L;{\boldsymbol h},B]\cdot\Sym_X^\#[L;{\boldsymbol h},A]$$ and after inversion,
\begin{multline*}
\Sym_X^\#[L;{\boldsymbol h},B]\cdot \Sym_X^\#[L;{\boldsymbol h},A]=\Sym_X^\#[L;{\boldsymbol h},B]^{-1}\cdot \Sym_X^\#[L;{\boldsymbol h},A]^{-1}\\=
(\Sym_X^\#[L;{\boldsymbol h},A]\cdot\Sym_X^\#[L;{\boldsymbol h},B])^{-1}
\subseteq (\Sym_X^\#[L;{\boldsymbol h},B]\cdot\Sym_X^\#[L;{\boldsymbol h},A])^{-1}\\=\Sym_X^\#[L;{\boldsymbol h},A]^{-1}\cdot\Sym_X^\#[L;{\boldsymbol h},B]^{-1}=\Sym_X^\#[L;{\boldsymbol h},A]\cdot\Sym_X^\#[L;{\boldsymbol h},B]
\end{multline*}
because $\Sym_X^\#[L;{\boldsymbol h},A]=\Sym_X^\#[L;{\boldsymbol h},A]^{-1}$ and $\Sym_X^\#[L;{\boldsymbol h},B]=\Sym_X^\#[L;{\boldsymbol h},B]^{-1}$ are groups.
Thus we obtain the desired equality 
$ \Sym_X^\#[L;{\boldsymbol h},A]\cdot\Sym_X^\#[L;{\boldsymbol h},B]=\Sym_X^\#[L;{\boldsymbol h},B]\cdot\Sym_X^\#[L;{\boldsymbol h},A].$
\end{proof}

\begin{proposition}\label{p:IX=n**(n-1)+} Let $\boldsymbol h,\boldsymbol v$ be two distinct directions in a Playfair plane $X$ such that for every lines $L,\Lambda\in\boldsymbol v$, the set $\Sym^\#_X[L;\boldsymbol h,\Lambda]$ is a subgroup of the group $\Sym_X^\#(L)$. If the cardinal $n\defeq |X|_2$ is finite, then for every line $L\in\boldsymbol v$, the cardinality of the group $\Sym_X^\#(L)$ divides the the number $n^{n-1}$.
\end{proposition}

\begin{proof} Fix any line $L\in\boldsymbol v$. Observe that  $n=|X|_2=|L_\parallel|$ and hence $L_\parallel =\{L_1,L_2,\dots, L_{n}\}$ for some lines $L_0,\dots,L_{n}$ such that $L_1=L$. By our assumption, for every $i,j\in\{1,\dots, n\}$, the set $\Sym_X^\#[L_i;{\boldsymbol h},L_j]=(\Sym_X^\#[L_i;{\boldsymbol h},L_j])^{-1}$ is a subgroup of the group  $\Sym_X^\#(L)$. 

If $i=1$, then the set $\Sym_X^\#[L;{\boldsymbol h},L_1]$ is a singleton containing the identity map of the line $L$. For $i\in\{2,\dots,n\}$, the group $\Sym_X^\#[L;{\boldsymbol h},L_i]$ has cardinality $|X|_2=n$, by Proposition~\ref{p:|IX|=|L|+}. By Corollary~\ref{c:IX-generators+}, the group $\Sym_X^\#(L)$ is generated by the set $\bigcup_{i=2}^n\Sym_X^\#[L;{\boldsymbol h},L_i]$. By Lemma~\ref{l:IX-commute+}, for any numbers $i,j\in\{1,\dots,n\}$ the subgroups $\Sym_X^\#[L;{\boldsymbol h},L_i]$ and $\Sym_X^\#[L;{\boldsymbol h},L_j]$ mutually commute in the sense that $ \Sym_X^\#[L;{\boldsymbol h},L_i]\cdot\Sym_X^\#[L;{\boldsymbol h},L_j]=
\Sym_X^\#[L;{\boldsymbol h},L_j]\cdot\Sym_X^\#[L;{\boldsymbol h},L_i].$

Consider the sequence of the sets $(H_k)_{k=1}^n$ defined by the recursive formula $$H_1=\Sym_X^\#[L;{\boldsymbol h},L_1]=\{1_L\}\quad\mbox{and}\quad H_{k+1}=H_k\cdot \Sym_X^\#[L;{\boldsymbol h},L_k]$$ for $k\in\{1,\dots,n-1\}$. The mutual commutativity of the subgroups $\Sym_X^\#[L;{\boldsymbol h},L_i]$ for $i\in\{1,\dots,n\}$ implies that every set $H_i$ is a subgroup of the group $\Sym_X^\#(L)$ and $H_n=\Sym^\#_X(L)$. It is easy  to see that for every $i\in\{2,\dots,n\}$ the index of the subgroup $H_{i-1}$ in $H_i$ is equal to the index of the subgroup $H_{i-1}\cap\Sym_X^\#[L;{\boldsymbol h},L_i]$ in the group $\Sym_X^\#[L;{\boldsymbol h},L_i]$ and hence the index of the group $H_{i-1}$ in $H_i$ divides the number $n=|\Sym^\#_X[L;{\boldsymbol h},L_i]|$. By induction this implies that for every $i\in\{1,\dots,n\}$, the cardinality of the group $H_i$ divides the number $n^{i-1}$. In particular, the cardinality of the group $H_n=\Sym_X^\#(L)$ divides the number $n^{n-1}$.
\end{proof}

Proposition~\ref{p:IX=n**(n-1)+} and Theorem~\ref{t:Burnside-center}  imply the following corollary.

\begin{corollary}\label{c:IX-center+} Let $\boldsymbol h,\boldsymbol v$ be two distinct direction in a Playfair plane $X$ such that for every lines $L,\Lambda\in \boldsymbol v$, the set $\Sym_X^\#[L;\boldsymbol h,\Lambda]$ is a subgroup of the group $\Sym_X^\#(L)$. If $|X|_2$ is a prime power, then for every $L\in\boldsymbol v$, the group $\Sym_X^\#(L)$ has non-trivial center.
\end{corollary}

Corollaries~\ref{c:IX-center+} and \ref{c:partialT<=>}  imply the following modification of Gleason's Theorem~\ref{t:partial-translation<=}.

\begin{theorem}\label{t:partial-translation<=+} A Playfair plane $X$ is $\partial$-translation whenever $|X|_2$ is a prime power and for distinct directions ${\boldsymbol h},{\boldsymbol v}\in\partial X$, and every lines $L,\Lambda\in \boldsymbol v$, the set $\Sym_X^\#[L;\boldsymbol h,\Lambda]$ is a subgroup of the group $\Sym_X^\#(L)$. 
\end{theorem}

\chapter{Vectors in affine spaces}\label{ch:vectors}

This chapter is devoted to (functional) vectors in affine spaces and their relation to translations. (Functional) vectors are equivalence classes of (Thalesian) pairs under the translation equivalence.  

\section{Translation equivalence of pairs}\label{s:trans-eq}

In this section we discuss the notion of translation equivalence of pairs in an affine space $X$. This equivalence is induced by the action of the inverse monoid $\I_X^\#$ on $X^2$. We recall that $\I_X^\#$ is the smallest submonoid of the symmetric inverse monoid $\I_X$ that contains all line shifts.
Line bijections that are finite compositions of line shifts are called line translations.

Every pair $x_0x_1\in X^2$ is (by definition) the function $\{(0,x_0),(1,x_1)\}:2\to X$ on the number $2\defeq\{0,1\}$. So, for any function $F$, we can consider the composition 
$$Fx_0x_1\defeq\{(i,F(x_i)):i\in 2,\;x_i\in\dom[F]\}$$ of the functions $x_0x_1$ and $F$. The function $Fx_0x_1$ has domain $\dom[Fx_0x_1]=\{i\in 2:x_i\in\dom[F]\}$. This domain equals $2$ if and only if $\{x_0,x_1\}\subseteq\dom[F]$.



\begin{definition} Let $X$ be an affine space. Given two pairs $ab,xy\in X^2$, we write \index[note]{$ab\#xy$}$ab\# xy$ and say that $ab$ and $xy$ are \index{translation equivalent pairs}\index{pairs!translation equivalent}\defterm{translation equivalent} if there exists a line translation $T$ in $X$ such that $Tab=xy$. 
\end{definition}

\begin{proposition}\label{p:transeq} For any points $a,b,u,v,x,y\in X$ of an affine space $X$, the translation equivalence has the following properties:
\begin{enumerate}
\item $ab\#ab$;
\item $ab\#uv\;\Rightarrow\;uv\#ab$;
\item $(ab\#uv\;\wedge\;uv\#xy)\;\Ra\;ab\#xy$;
\item $ab\#uv\;\Rightarrow ba\#vu$;
\item $aa\#xx$;
\item If $ab\#uv$ and $a=b$, then $u=v$.
\end{enumerate}
\end{proposition}

\begin{exercise} Prove Proposition~\ref{p:transeq}.
\end{exercise}

In the following theorem we characterize the translation equivalence without mentioning line translations.

\begin{theorem}\label{t:trans-eq<=>} Two pairs $ab,uv$ in a (Thalesian) affine space $X$ are translation equivalent if and only if  there exist a number $n\in\w$ (with $n\le 2$) and pairs $x_0y_0,x_1y_1,\dots,x_ny_n\in X^2$ such that $x_0y_0=ab$, $x_ny_n=uv$, and  for every positive $i\le n$, $$\Aline {x_{i-1}}{x_i}\parallel\Aline {y_{i-1}}{y_i}\quad \mbox{and}\quad  \Aline {x_{i-1}}{y_{i-1}}\cap \Aline{x_i}{y_i}=\varnothing.$$
\end{theorem}

\begin{proof} To prove the ``only if'' part, assume that the pairs $ab$ and $uv$ are translation equivalent. If $ab=uv$, then for $n=0$ and the pair $x_0y_0\defeq ab$, we obtain $ab=x_0y_0=x_ny_n=uv$ (and since no positive $i\le n=0$ exists the conditions $\Aline {x_{i-1}}{x_i}\parallel\Aline {y_{i-1}}{y_i}$ and $\Aline {x_{i-1}}{y_{i-1}}\cap \Aline{x_i}{y_i}=\varnothing$ holds vacuously).

Now assume that $ab\ne uv$. Since $ab\#uv$, there exists a line translation $T$ in $X$ such that $Tab=uv$. Write $T$ is the composition $S_n\cdots S_1$ of non-identity line shifts in $X$. We can assume that the number $n$ is the smallest posssible. If the affine space $X$ is Thalesian, then $n\le 2$, by  Corollary~\ref{c:many-shifts=2shifts}. By the definition of a line shift, for every $i\in\{1,\dots,n\}$, there exists a direction $\delta_i$ such that $S_i=\{(x,y)\in\dom[S_i]\times\rng[S_i]:\Aline x{\delta_i}=\Aline y{\delta_i}\}$. Let $x_0y_0\defeq ab$ and for every $i\in\{1,\dots,n\}$, let $x_iy_i\defeq S_ix_{i-1}y_{i-1}$. For every $i\in\{1,\dots,k\}$, we have 
$$\Aline {x_{i-1}}{y_{i-1}}\cap\Aline {x_i}{y_i}\subseteq\dom[S_i]\cap\rng[S_i]=\varnothing,$$
which implies that $x_{i-1}\ne x_i$ and $y_{i-1}\ne y_i$. The choice of the direction $\delta_i$ ensures that the lines $\Aline {x_{i-1}}{x_i}$ and $\Aline {y_{i-1}}{y_i}$ belong to the direction $\delta_i$ and hence $\Aline{x_{i-1}}{x_i}\parallel \Aline{y_{i-1}}{y_i}$.
This completes the proof of the ``only if'' part.
\smallskip

\begin{picture}(200,85)(-80,-15)

\put(0,0){\line(0,1){40}}
\put(0,0){\line(2,1){30}}
\put(0,40){\line(2,1){30}}
\put(30,15){\line(0,1){40}}
\put(30,15){\line(2,-1){30}}
\put(30,55){\line(2,-1){30}}
\put(60,0){\line(0,1){40}}
\put(60,0){\line(2,1){20}}
\put(60,40){\line(2,1){20}}
\put(100,30){\circle*{2}}
\put(110,30){\circle*{2}}
\put(120,30){\circle*{2}}
\put(140,10){\line(2,-1){20}}
\put(140,50){\line(2,-1){20}}
\put(160,0){\line(0,1){40}}
\put(160,0){\line(2,1){30}}
\put(160,40){\line(2,1){30}}
\put(190,15){\line(0,1){40}}
\put(190,15){\line(2,-1){30}}
\put(190,55){\line(2,-1){30}}
\put(220,0){\line(0,1){40}}

\put(0,0){\circle*{3}}
\put(-8,-3){$a$}
\put(2,-6){$x_0$}
\put(0,40){\circle*{3}}
\put(-8,38){$b$}
\put(2,33){$y_0$}
\put(30,15){\circle*{3}}
\put(27,6){$x_1$}
\put(30,55){\circle*{3}}
\put(27,60){$y_1$}
\put(60,0){\circle*{3}}
\put(56,-8){$x_2$}
\put(60,40){\circle*{3}}
\put(56,45){$y_2$}
\put(160,0){\circle*{3}}
\put(150,-8){$x_{n-2}$}
\put(160,40){\circle*{3}}
\put(150,48){$y_{n-2}$}
\put(190,15){\circle*{3}}
\put(180,5){$x_{n-1}$}
\put(190,55){\circle*{3}}
\put(180,60){$y_{n-1}$}
\put(220,0){\circle*{3}}
\put(210,-6){$x_n$}
\put(223,-3){$u$}
\put(220,40){\circle*{3}}
\put(208,35){$y_n$}
\put(223,38){$v$}
\end{picture}

To prove the ``if'' part, assume that for the pairs $ab$ and $uv$ there exist a number $n\in\w$ and pairs $x_0y_0,\dots,x_ny_n\in X^2$ such that $x_0y_0=ab$, $x_ny_n=uv$ and  $\Aline {x_{i-1}}{x_i}\parallel\Aline {y_{i-1}}{y_i}$, $\Aline {x_{i-1}}{y_{i-1}}\cap \Aline{x_i}{y_i}=\varnothing$ for every positive $i\le n$.

\begin{claim}\label{cl:xi=yi} For every $i\in\{1,\dots,n\}$, $x_{i-1}=y_{i-1}\;\Leftrightarrow\;x_i=y_i.$
\end{claim}

\begin{proof}  Assume that $x_{i-1}=y_{i-1}$. Then the parallel lines  $\Aline{x_{i-1}}{x_i}$ and $\Aline{y_{i-1}}{y_i}$ coincide, by Corollary~\ref{c:para+intersect=>coincide}. 

Assuming that $x_i\ne y_i$, we conclude that $\Aline{x_i}{y_i}=\Aline{x_{i-1}}{x_i}=\Aline{y_{i-1}}{y_i}$ and hence  $\Aline{x_{i-1}}{y_{i-1}}\subseteq \Aline{x_i}{y_i}$, which contradicts the condition $\Aline{x_{i-1}}{y_{i-1}}\cap\Aline {x_i}{y_i}=\varnothing$. This contradiction shows that $x_i=y_i$, and completes the proof of the implication $x_{i-1}=y_{i-1}\;\Ra\;x_i=y_i$. 

By analogy we can prove that $x_i=y_i\;\Ra\; x_{i-1}=y_{i-1}$.
\end{proof}

If $a\ne b$, then Claim~\ref{cl:xi=yi} implies that $x_i\ne y_i$ for all $i\in\{0,\dots,n\}$. In this case, for every $i\in\{1,\dots,k\}$ the line $\Lambda_i\defeq\Aline {y_{i-1}}{y_i}$ is not subparallel to the disjoint parallel lines $\Aline {x_{i-1}}{y_{i-1}}$ and $\Aline{x_i}{y_i}$. By Proposition~\ref{p:flat-relation}, the relation $$S_i\defeq\{(x,y)\in\Aline {x_{i-1}}{y_{i-1}}\times\Aline{x_i}{y_i}:\Aline xy\subparallel \Lambda_i\}$$is a line projection, which is a line shift because $\dom[S_i]=\Aline{x_{i-1}}{y_{i-1}}$ and $\rng[S_i]=\Aline {x_i}{y_i}$ are disjoint. Then the composition $T=S_n\dots S_1$ is a line translation witnessing that $Tab=uv$ and hence $ab\#uv$.

If $a=b$, then Claim~\ref{cl:xi=yi} implies that $x_i=y_i$ for all $i\in\{0,\dots,n\}$. In particular, $u=x_n=y_n=v$. Take any line $\Lambda\subseteq X$ containing the points $a=b$ and $u=v$. Next, choose any line $L\subseteq X$ such that $L\cap\Lambda=\{a\}$. Since the space $X$ is affine, there exists a line $L'$ in $X$ such that $u\in L'$ and $L'\parallel L$. By Proposition~\ref{p:flat-relation}, the relation $T\defeq\{(x,y)\in L\times L':\Aline xy\subparallel \Lambda\}$ is a line translation with $Tab=uv$, witnessing that $ab\#uv$.
\end{proof}

\begin{exercise}\label{ex:treq} Let $F:X\to Y$ be an isomorphism of affine spaces $X,Y$, and $xy,uv\in X^2$, $x'y',u'v'\in Y^2$ be pairs such that $uv=Fxy$ and $u'v'=Fx'y'$. Show that the pairs $xy$ and $uv$ are translation equivalent in the liner $X$ if and only if the pairs $x'y'$ and $u'v'$ are translation equivalent in $Y$.
\end{exercise}

\section{Thalesian pairs in affine spaces}\label{s:Thales-pairs} 

Corollary~\ref{c:par-trans} and Theorem~\ref{t:unique-translation} imply the following result on the existence of translation equivalent pairs in affine spaces.

\begin{corollary}\label{c:transeq-exists} For every pair of points $xy\in X^2$ in a (Thalesian) affine space $X$ and every point $a\in X$, there exists a (unique) point $b\in X$ such that $ab\#xy$.
\end{corollary}

This corollary motivates the following definition.

\begin{definition}\label{d:Thalesian-pair} A pair of points $xy\in X^2$ in an affine space $X$ is called \index{Thalesian pair}\index{pair!Thalesian}\defterm{Thalesian} if for
every point $a \in X$ there exists a {\em unique} point $b\in X$ such that $ab\#xy$. 
\end{definition}

\begin{example}  For every point $x\in X$ the pair $xx$ is Thalesian.
\end{example}

Theorem~\ref{t:unique-translation} implies the following characterization of Thalesian affine spaces in terms of Thalesian pairs.

\begin{corollary}\label{c:Thalesian<=>pair-Thalesian} An affine space $X$ is Thalesian if and only if every pair $xy\in X^2$ is
Thalesian.
\end{corollary}

\begin{exercise} Let $ab, xy\in X^2$ be two translation equivalent pairs in an affine space $X$. Assuming that the pair $ab$ is Thalesian, show that the pair $xy$ is Thalesian.
\smallskip

{\em Hint:} Use the transitivity of the translation equivalence.
\end{exercise}

\begin{proposition}\label{p:xy-Thales=>yx-Thales} If a pair $xy\in X^2$ in an affine space $X$ is Thalesian, then the pair $yx$
is Thalesian, too.
\end{proposition}

\begin{proof} To prove that the pair $yx$ is Thalesian, we have two show that every point $x'\in X$ with $yx'\#yx$ coincides with $x$. If $x=y$, then we can apply Proposition~\ref{p:transeq}(6) and conclude that
$x'=y=x$. So, assume that $x\ne y$. Since $yx\#yx'$, there exists a line translation $T:\Aline yx\to\Aline y{x'}$ such that $Tyx = yx'$. Proposition~\ref{p:parshift=} implies that $\Aline y x\parallel \Aline y{x'}$ and hence $\Aline yx=\Aline y{x'}$. Consider the line $L\defeq \Aline x y = \Aline x{y'}$. Since $X$ is an affine space, there exists a line $\Lambda\subset X$ such that $\Lambda\parallel L$ and $\Lambda\cap L = \varnothing$. Fix any line shift $S :L\to\Lambda$ and consider the pair $ab\defeq Sx'y$. By Corollary~\ref{c:parashift-exists},
there exist line shifts $S':\Lambda\to L$ such that $x=S'(a)= S'S(x') = S'ST(x)$. Since the pair $xy$ is Thalesian, $y = S'ST(y) = S'S(y) = S'(b)$. It follows from $Sx'y = ab$ and $S'ab = xy$ that
$\Aline {x'}a\parallel \Aline yb\parallel \Aline xa$ and hence $\Aline {x'}a = \Aline xa$ and finally, $x'\in L \cap\Aline{x'}a = L\cap\Aline xa = \{x\}$.
\end{proof}

\begin{theorem}\label{t:aa'-Thalesian<=>} A pair $aa'\in X^2$ in an affine space $X$ is Thalesian if and only if
$$\forall b,b',c,c'\in X\;\;\big((\Aline ab\parallel \Aline {a'}{b'})\;\wedge\;(\Aline bc\parallel \Aline {b'}{c'})\;\wedge\;(\Aline a{a'}\cap\Aline b{b'}=\varnothing=\Aline b{b'}\cap\Aline c{c'})\big)\;\Ra\;(\Aline a{a'}\parallel \Aline c{c'}).$$
\end{theorem}

\begin{proof} First assume that a pair $aa'$ is Thalesian and take any pairs 
$bb',cc'\in X^2$ with $\Aline ab\parallel \Aline {a'}{b'}$, $\Aline bc\parallel \Aline {b'}{c'}$, and $\Aline a{a'}\cap\Aline b{b'}=\varnothing=\Aline b{b'}\cap\Aline c{c'}$. We have to prove that $\Aline a{a'}\parallel \Aline c{c'}$. Applying Theorem~\ref{t:trans-eq<=>}, we conclude that $a{a'}\#cc'$. 
If $a=a'$, then $c=c'$, by Proposition~\ref{p:transeq}(6). In this case, $\Aline a{a'}=\{a\}\parallel \{c\}=\Aline c{c'}$.

So, assume that $a\ne a'$.  Proposition~\ref{p:transeq}(6) implies that $b\ne b'$ and $c\ne c'$. It follows from $\Aline ab\parallel \Aline {a'}{b'}$, $\Aline bc\parallel \Aline {b'}{c'}$, and $\Aline a{a'}\cap\Aline b{b'}=\varnothing=\Aline b{b'}\cap\Aline c{c'}$ that $\Aline a{a'}\parallel \Aline b{b'}\parallel \Aline c{c'}$ and hence $\Aline a{a'}\parallel \Aline c{c'}$. If $c\in \Aline a{a'}$, then $\Aline c{c'}=\Aline a{a'}$ and hence $\Aline a{a'}\parallel \Aline c{c'}$. If $c\notin \Aline a{a'}$, then the parallel lines $\Aline a{a'}$ and $\Aline c{c'}$ are disjoint. Since the affine space $X$ is Playfair, there exists a point $c''\in \Aline c{c'}$ such that $\Aline ac\parallel \Aline {a'}{c''}$. Then $cc''\#aa'\#cc'$ and hence $c''=c'$ because the pair $aa'$ is Thalesian. Then $\Aline ac\parallel \Aline {a'}{c''}=\Aline {a'}{c'}$.
\smallskip

Now assume that
\begin{equation}\label{eq:Thalesian=>}
\forall bb',cc'\in X^2\;\;\big(\Aline ab\parallel \Aline {a'}{b'})\;\wedge\;(\Aline bc\parallel \Aline {b'}{c'})\;\wedge\;(\Aline a{a'}\cap\Aline b{b'}=\varnothing=\Aline b{b'}\cap\Aline c{c'})\big)\;\Ra\;(\Aline a{a'}\parallel \Aline c{c'}).
\end{equation}
 We have to prove that the pair $aa'$ is Thalesian. In the opposite case, we can find a point $a''\ne a'$ such that $aa'\#aa''$. By Theorem~\ref{t:trans-eq<=>}, there exists a number $n\in\IN$ and a sequence of pairs $x_0y_0,x_1y_1,\dots,x_ny_n\in X^2$ such that $x_0y_0=aa'$, $x_ny_n=aa''$ and for every positive $i\le n$ we have 
$\Aline {x_{i-1}}{x_i}\parallel \Aline {y_{i-1}}{y_i}$ and $\Aline {x_{i-1}}{y_{i-1}}\cap\Aline {x_i}{y_i}=\varnothing$. We can assume that the number $n$ is the smallest possible. It follows from $a'\ne a''$ that $n\ge 2$. If $\Aline {x_0}{y_0}\cap \Aline {x_2}{y_2}=\varnothing$, then the assumption (\ref{eq:Thalesian=>}) ensures that $\Aline {x_0}{x_2}\parallel \Aline {y_0}{y_2}$ and then the pair $x_1y_1$ can be removed from the sequence $x_0y_0,x_1y_1,\dots,x_ny_n$, which contradicts the minimality of $n$. This contradiction shows that $\Aline {x_0}{y_0}=\Aline {x_2}{y_2}$. Assuming that $n=2$, we conclude that $x_2y_2=aa''$ and hence $\Aline {a''}{y_1}=\Aline {y_2}{y_1}\parallel \Aline  {x_2}{x_1}=\Aline {x_0}{x_1}\parallel \Aline {y_0}{y_1}=\Aline {a'}{y_1}$ and $a''\in \Aline {a''}{y_1}\cap\Aline {x_0}{y_0}=\Aline{a'}{y_1}\cap\Aline{x_0}{y_0}=\{a'\}$, which contradicts the choice of $a''$. This contradiction shows that $n\ge 3$.

\begin{picture}(200,110)(-160,-15)

\put(0,0){\line(1,0){40}}
\put(0,20){\line(1,0){40}}
\put(40,0){\line(2,1){40}}
\put(40,20){\line(2,1){40}}
\put(0,0){\line(0,1){20}}
\put(40,0){\line(0,1){20}}
\put(80,20){\line(0,1){20}}
\put(80,20){\line(-2,1){40}}
\put(80,40){\line(-2,1){40}}
\put(0,0){\line(1,1){70}}
\put(0,20){\line(1,1){70}}
\put(70,70){\line(0,1){20}}
\put(40,40){\line(0,1){20}}

\put(0,0){\circle*{3}}
\put(-3,-9){$x_1'$}
\put(0,20){\circle*{3}}
\put(-6,26){$y_1'$}
\put(40,0){\circle*{3}}
\put(36,-9){$x_0$}
\put(40,20){\circle*{3}}
\put(35,24){$y_0$}
\put(80,20){\circle*{3}}
\put(83,17){$x_1$}
\put(80,40){\circle*{3}}
\put(83,40){$y_1$}
\put(40,40){\circle*{3}}
\put(37,32){$x_2$}
\put(40,60){\circle*{3}}
\put(35,65){$y_2$}
\put(70,70){\circle*{3}}
\put(73,67){$x_3$}
\put(70,90){\circle*{3}}
\put(73,90){$y_3$}

\end{picture}

The minimality of $n$ ensures that $\Aline {x_1}{x_2}\ne\Aline{x_2}{x_3}$. Choose any point $x_1'\in \Aline {x_2}{x_3}\setminus\{x_2,x_3\}$ and using Theorem~\ref{t:parallelogram3+1}, find a unique point $y_1'\in X$ such that $\Aline {x_1'}{y'_1}\parallel \Aline {x_0}{y_0}$ and $\Aline {y_0}{y'_1}\parallel \Aline {x_0}{x_1'}$. It follows from $\Aline {x_2}{y_2}\cap\Aline {x_3}{y_3}=\varnothing$ that $\Aline {x_1'}{y_1'}\cap \Aline {x_0}{y_0}=\Aline{x_1'}{y_1'}\cap\Aline{x_2}{y_2}=\varnothing$. We claim that $y_2\in \Aline {y_1'}{y_3}$ and hence $\Aline {y_1'}{y_3}=\Aline {y_2}{y_3}\parallel \Aline {x_2}{x_3}=\Aline {x_1'}{x_3}$. In the opposite case, we can find a point $y_2'\in \Aline {y_1}{y_2}\setminus\{y_2\}$ such that $\Aline {y_1'}{y_2'}\parallel \Aline {x_1'}{x_2}$. Since $X$ is Playfair, there exist unique points $x_2'\in \Aline {x_1}{x_2}$ and $x_2''\in \Aline{x_1'}{x_2}$ such that $\Aline {x_2'}{y_2'}\parallel \Aline {x_2}{y_2}\parallel \Aline{x_2''}{y_2}$. The assumption (\ref{eq:Thalesian=>}) ensures that $\Aline {x_0}{x_2'}\parallel\Aline {y_0}{y_2'}\parallel \Aline{x_0}{x_2''}$ and hence $\Aline {x_0}{x_2'}=\Aline{x_0}{x_2''}=\Aline{x_2'}{x_2''}$, which is not possible because the line $\Aline {x_2'}{x_2''}=\Aline{x_2'}{y_2'}=\Aline{x_2''}{y_2'}$ is disjoint with the line $\Aline {x_2}{y_2}=\Aline{x_0}{y_0}$. This contradiction shows that $\Aline{x'_1}{x_3}\parallel\Aline {y_1'}{y_3}$ and then we can replace the sequence $x_0y_0,x_1y_1,x_2y_2,\dots,x_ny_n$ by the shorter sequence $x_0y_0,x_1'y_1',x_3y_3,\dots,x_ny_n$, which contradicts the minimality of $n$. This contradiction shows that the pair $aa'$ is Thalesian.
\end{proof}

\section{Vectors and functional vectors}

\begin{definition}\label{d:vector} For points $x,y$ of an affine space $X$, the equivalence class \index[note]{$\overvector{xy}$}$$\overvector{xy}\defeq\{uv\in X^2:uv\#xy\}$$ is called the \index{vector}\defterm{vector} determined by the pair $xy$.
\end{definition}

The set of vectors in an affine space $X$ is denoted by \index[note]{$X^2_{\#}$}\defterm{$X^2_{\#}$}. Therefore, $$X^2_\#\defeq\{\overvector{xy}:xy\in X\}.$$

\begin{definition} A vector $\vecv\in X^2_{\#}$
in an affine space $X$ is called a \index{functional vector}\index{vector!functional}\defterm{functional vector} (briefly,
a \index{funvector}\defterm{funvector}) if $\vecv$ is a function, i.e., for every point $x\in X$ there exists a unique point $y\in X$ such
that $xy\in\vecv$.
\end{definition}

The set of all functional vectors in an affine space $X$ is denoted by \index[note]{$\overvector X$}\defterm{$\overvector X$}.

Definition~\ref{d:Thalesian-pair} implies the following characterization of funvectors.

\begin{proposition} A vector $\vecv\in X^2_{\#}$ in an affine space $X$ is functional if and only if every pair $xy\in\vecv$ is Thalesian if and only if some pair $xy\in \vecv$ is Thalesian.
\end{proposition}

Corollary~\ref{c:Thalesian<=>pair-Thalesian} implies the following characterization of Thalesian affine spaces.

\begin{corollary}\label{c:Thalesian<=>vector=funvector}For an affine space $X$, the following conditions are equivalent:
\begin{enumerate}
\item the affine space $X$ is Thalesian;
\item every vector in $X$ is functional;
\item $\overvector X=X^2_{\#}$;
\item for every point $x\in X$ and vector $\vecv\in X^2_{\#}$, there exists a unique point $y\in X$ such that $\vecv=\overvector{xy}$.
\end{enumerate}
\end{corollary}

The set of functional vectors $\overvector X$ contains the distinguished element \index{zero vector}\index{vector zero}\index[note]{$\vec{\mathbf 0}$}$$\mbox{\defterm{$\vec{\mathbf 0}$}}\defeq\{xx:x\in X\},$$ called the \defterm{zero vector}. 

The set $X^2_{\#}$ of vectors admits a natural unary operation
$$-:X^2_{\#}\to X^2_{\#},\quad -:\vecv\mapsto -\vecv\defeq\{yx:xy\in\vecv\}.$$
The vector $-\vecv = \{\overvector{yx}: xy\in\vecv\}$ is called the \index{opposite vector}\index{vector!opposite}\defterm{opposite vector} to the vector $\vecv$.
Proposition~\ref{p:transeq}(4) implies that this operation is well-defined. It is clear that $-\vec{\mathbf 0}=\vec{\mathbf 0}$, and the operation of taking the opposite vector is involutive in the sense that $-(-\vecv)=\vecv$ for every vector $\vecv\in X^2_{\#}$.  Proposition~\ref{p:xy-Thales=>yx-Thales} implies that for every functional vector $\vecv\in\overvector X$, its opposite vector $-\vecv$ is a functional vector. So, the operation of taking
the opposite vector is also an unary operation on the set of functional vectors $\overvector X$.

%

\begin{exercise} Is every vector $\vec {\boldsymbol v}$ with $-\vec{\boldsymbol v}=\vec{\boldsymbol v}$ equal to $\vec{\mathbf 0}$?
\smallskip

\noindent{\em Hint:} No, because in vector spaces over a 2-element field, the equality $-\vec{\boldsymbol v}=\vec{\boldsymbol v}$ holds for every vector $\vec{\boldsymbol v}$.
\end{exercise}


\begin{proposition}\label{p:=vectors=>parallel-sides} Let $x,y,u,v$ be points of an  affine space $X$. If $\overvector{xy}=\overvector{uv}$ (and $\overvector{xy}=\overvector{uv}\in\overvector X$), then $\Aline xy\parallel \Aline uv$ (and $\Aline xu\parallel \Aline yv$).
\end{proposition} 

\begin{proof} Assuming that $\overvector{xy}=\overvector{uv}$, find a line translation $T$ in $X$ such that $Txy=uv$. 

It $x=y$, then $u=v$, by Proposition~\ref{p:transeq}(6). Then $\Aline xy=\{x\}\parallel \{u\}=\Aline uv$ and $\Aline xu=\Aline yv$.

So, assume that $x\ne y$ and hence $u\ne v$. Then $\Aline xy=\dom[T]\parallel\rng[T]=\Aline uv$, by Proposition~\ref{p:parshift=}. 

Finally, assuming that $\overvector{xy}=\overvector{uv}\in\overvector X$, we shall prove that $\Aline xu\parallel \Aline yv$. If $x=u$, then the equality $\overvector{xy}=\overvector{uv}\in\overvector X$ implies $y=v$ and then $\Aline xu=\{x\}\parallel \{y\}=\Aline yv$. If $x\ne u$ and $\Aline xy=\Aline uv$, then $y\ne v$ and then  
$\Aline xu = \Aline xy =\Aline uv = \Aline yv$ and we are done. So, assume that $\Aline xy\ne\Aline uv$. By Corollary~\ref{c:parashift-exists}, there
exists a line shift $S : \Aline xy\to \Aline uv$ such that $S(x) = u$. Let $w\defeq S(y)$ and observe that $uv\#xy\#Sxy = uw$.
Since $\overvector{uv}\in \overvector X$, the pair $uv$ is Thalesian and hence $v = w$. Since $Sxy=uw$, the definition of the line shift $S$ implies  $\Aline yv = \Aline yw\parallel \Aline x u$.
\end{proof}

\begin{proposition}\label{p:parallelogram=>vectors=} Let $p,q,u,v$ be points in an affine space $X$. If $\Aline pu\parallel \Aline qv$ and $\Aline pq\cap\Aline uv=\varnothing$, then $\overvector{pq}=\overvector{uv}$.
\end{proposition}

\begin{proof} If $p=q$, then $\Aline pu\parallel \Aline qv$ implies $\Aline pu=\Aline qv$. Assuming that $u\ne v$, we conclude that $\Aline uv=\Aline pu=\Aline qv$ and hence $\Aline pq\cap\Aline uv=\Aline pq\ne\varnothing$, which contradicts our assumption. This contradiction shows that $u=v$. Then $\overvector{pq}=\vec{\mathbf 0}=\overvector{uv}$. 
By analogy we can show that $u=v$ implies $\overvector{uv}=\vec{\mathbf 0}=\overvector{pq}$. 

So, we assume that $p\ne q$ and $u\ne v$.  It follows from $\Aline pq\cap\Aline uv=\varnothing$ that $\Lambda\defeq\Aline qv$ is a line, which is not subparallel to the lines $\Aline pq$ and $\Aline uv$. By Proposition~\ref{p:flat-relation}, the relation 
$$P\defeq\{(x,y)\in\Aline pq\times\Aline uv:\Aline xy\subparallel \Lambda\}$$is a parallel projection. It follows from $\{p,q\}\subseteq\dom[P]\subseteq\Aline pq$ and $\{u,v\}\subseteq\rng[P]\subseteq\Aline uv$ that $\dom[P]=\Aline pq$ and $\rng[P]=\Aline uv$, which means that $P$ is a line projection. Since  $\dom[P]\cap\rng[P]=\Aline pq\cap\Aline uv=\varnothing$, the line projection $P$ is a line shift such that $Ppq=uv$, witnessing that $pq\#uv$ and hence $\overvector{pq}=\overvector{uv}$.
\end{proof}




\begin{remark} For every point $o$ in an affine space,
 the function
$$\overvector X\to X,\quad\overvector{ox}\mapsto x,$$
is well-defined and injective and the function
$$X\to X^2_{\#},\quad x\mapsto \overvector{ox},$$
is surjective. For a point $x\in X$, the vector $\overvector{ox}$ is called \index{radius-vector}\defterm{the radius-vector} of the point $x$.
\end{remark}


\begin{exercise} Explore (functional) vectors on the Moulton plane (which is not Thalesian).
\end{exercise}

\section{Functional vectors versus translations}\label{s:vectors=translations}

Observe that every vector $\vecv$ in an affine space $X$ is the set of pairs $xy\in X^2$. Every pair $xy = \{(0, x), (1, y)\}\in X^2$ can be canonically identified with the ordered pair $(x, y)\in X\times X$,
which allows us to think of vectors as relations (which are sets of ordered pairs). The definition
of a functional vector implies that for every $x\in X$ there exists a unique point $y\in X$ such that
$xy\in \vecv$, which means that $\vecv$ is a well-defined function, assingning to every point $x\in X$ the unique
point $y\in X$ such that $xy\in \vecv$. Therefore, it is legal to think of functional vectors as functions from $X$ to $X$.
It turns out that such functions are exactly translations of $X$.

\begin{theorem}\label{t:Trans(X)=vecX} The equality $\overvector X=\Trans(X)$ holds for every affine space $X$.
\end{theorem}

\begin{proof} First we prove that $\overvector X\subseteq\Trans(X)$. Fix any functional vector $\vecv\in\overvector X$. If $\vecv=\vec{\mathbf 0}\defeq \{xx:x \in X\} = 1_X$, then $\vecv$ is the identity translation of $X$ and we are done. So, assume that $\vecv\ne \vec{\mathbf 0}$.

\begin{claim} The function $\vecv:X\to X$ is bijective.
\end{claim}

\begin{proof} By Proposition~\ref{p:xy-Thales=>yx-Thales}, the set $\vecv^{-1}\defeq\{yx : xy\in\vecv\}$ is a functional vector in $X$. Taking into account that $\vecv$ and $\vecv^{-1}$ are functions with $\dom[\vecv] = X = \dom[\vecv^{-1}]$, we conclude that $\vecv\circ\vecv^{-1}=1_X=\vecv^{-1}\circ\vecv$, which means that $\vecv:X\to X$ is a bijective function.
\end{proof}

\begin{claim} The functional vector $\vecv$ is a dilation of $X$.
\end{claim}

\begin{proof} Given any line $L$ in $X$, we have to show that the set $\vecv[L]$ is a line, parallel to the line $L$. Fix any point $o\in L$ and consider the point $v\defeq \vecv(o)$. 

If $v\in L$, then for every $x\in L$, the points $y\defeq\vecv(x)$ and $z\defeq\vecv^{-1}(x)$ have $\overvector{xy}=\vecv=\overvector{ov}$ and $\overvector{xz}=\vecv^{-1}=\overvector{vo}$. By Proposition~\ref{p:=vectors=>parallel-sides}, $\Aline xy\parallel \Aline ov=L$  and $\Aline xz\parallel \Aline vo= L$. Hence $\vecv(x)=y\in\Aline xy=L$
and $\vecv^{-1}(x)=z\in\Aline xz=L$. Therefore, $\vecv[L]\subseteq L$ and $\vecv^{-1}[L]\subseteq L$, which implies that $\vecv[L]=L$ is a line parallel to the line $L$. 

Next,  assume that $v\notin L$. Since $X$ is an affine space, there exists a unique line $\Lambda$ such that $v\in\Lambda$ and $\Lambda\parallel L$. For every point $x\in L\setminus\{o\}$, the point $y\defeq\vecv(x)$ has the property $\overvector{xy}=\vecv=\overvector{ov}$.
By Proposition~\ref{p:=vectors=>parallel-sides}, $\overvector{xy}= \overvector{ov}\in\overvector X$ implies $\Aline vy\parallel\Aline ox=L$. Since the affine space $X$ is Proclus,
$\Aline vy\parallel L\parallel\Lambda$ implies $y\in\Aline vy=\Lambda$. Therefore, $\vecv[L]\subseteq\Lambda$. By analogy we can prove that $\vecv^{-1}[\Lambda]\subseteq L$.
Therefore, $\vecv [L]=\Lambda\parallel L$, witnessing that $\vecv$ is a dilation of $X$.
\end{proof}

It remains to show that the dilation $\vecv$ of $X$ is a translation of $X$. In the opposite case, there exists  a point $x\in X$ whose image $y\defeq\vecv(x)$ coincides with $x$. Then $\vecv=\overvector{xy}=\overvector {xx}=\vec{\mathbf 0}$, which contradicts our assumption. Therefore, $\vecv$ is a translation of $X$ and $\overvector X\subseteq\Trans(X)$.

To see that $\Trans(X)\subseteq \overvector X$, choose any translation $T\in\Trans(X)$.

\begin{lemma}\label{l:Tra-equivalent} Every two pairs $xy,x'y'\in T$ are translation equivalent.
\end{lemma}

\begin{proof} By Proposition~\ref{p:Trans-spread}, the family $\overline T \defeq \{\Aline x y : xy\in T\}$ is a spread of parallel lines in $X$.

\begin{claim}\label{cl:Trans-eq} If $xy, x'y'\in T$ are two pairs with $\Aline xy\ne \Aline{x'}{y'}$, then $xy\#x'y'$.
\end{claim}

\begin{proof} Since $T$ is a dilation, the lines $\Aline x{x'}$ and $\Aline {y}{y'} = T[\Aline x{x'}]$ are parallel, which implies that $\Pi\defeq\{x, x', y, y'\}$ is a plane. Since $\Aline x y$ and $\Aline{x'}{y'}$ are distinct lines in the spread of parallel lines $\overline T$, they are disjoint and parallel. Taking into account that $\Aline x y \parallel\Aline {x'}{y'}$ and $\Aline x{x'}\parallel\Aline y{y'}$, we conclude that $xy\#x'y'$.
\end{proof}

Now we are ready to prove that any pairs $xy, x'y'\in T$ are translation equivalent. If $\Aline xy\ne\Aline{x'}{y'}$, then $xy\#x'y'$, by Claim~\ref{cl:Trans-eq}. If $\Aline x y = \Aline {x'}{y'}$, then chose any point $x''\in X \setminus\Aline xy$ and consider
the point $y''\defeq T(x'')$. Since $\Aline {x'}{y'}=\Aline x y\ne\Aline {x''}{y''}$, we can apply Claim~\ref{cl:Trans-eq} and
conclude that $xy\#x''y''\#x'y'$, which implies that $xy\#x'y'$, by the transitivity of the translation equivalence.
\end{proof}

\begin{claim}\label{cl:T=vec(xy)} $T = \overvector{xy}$ for every pair $xy\in T$.
\end{claim}

\begin{proof} Fix any pair $xy\in T$. If $x = y$, then $T(x) = y = x = 1_X(x)$ and hence $T = 1_X =\overvector{xx}=\overvector{xy}$, by Proposition~\ref{p:Ax=Bx=>A=B}. So, assume that $x\ne y$. Lemma\ref{l:Tra-equivalent} ensures that $T\subseteq\overvector{xy}$. To see that $\overvector{xy}\subseteq T$, take any pair $uv\in\overvector{xy}$. By Proposition~\ref{p:transeq}(6), the inequality $x\ne y$ implies
$u\ne v$. Since $xy\#uv$, there exists a line translation $S :\Aline x y\to\Aline u v$ such that $Sxy = uv$. By Theorem~\ref{t:TS=ST}, $TS = ST$ and hence $T(u) = TS(x) = ST(x) = S(y) = v$ and finally $uv \in T$.
Therefore, $\overvector{xy}=T$.
\end{proof}

Claim~\ref{cl:T=vec(xy)} implies that for every $xy\in T$, the set $T = \overvector{xy}$ is a vector in $X$. Since $T$ is a function with $\dom[T] = X$, for every $x\in X$ there exists a unique point $y\in X$ with $xy\in T$, which
implies that the vector $T$ is functional. Therefore, $\Trans(X) \subseteq \overvector X$ and finally, $\overvector X=\Trans(X)$.
\end{proof}

\begin{remark} Theorem~\ref{t:Trans(X)=vecX} is related to an old discussion about the nature of vectors. Are vectors equivalence classes of directed segments or translations? Theorem~\ref{t:Trans(X)=vecX} shows that
functional vectors are simultaneously translations and equivalence classes of directed segments.
This dual nature of functional vectors can be considered as a geometric counterart of the famous wave--particle duality in quantum mechanics
\footnote{Wave--particle duality is the concept in quantum mechanics that quantum entities exhibit particle or wave
properties according to the experimental circumstances. It expresses the inability of the classical concepts such as
particle or wave to fully describe the behavior of quantum objects. During the 19th and early 20th centuries, light
was found to behave as a wave, and then later discovered to have a particulate character, whereas electrons were
found to act as particles, and then later discovered to have wavelike aspects. The concept of duality arose to name these contradictions.}.
\end{remark}

Recall that a liner $X$ is called \index{translation liner}\index{liner!translation}\defterm{translation} if for every points $a, b\in X$, there exists a translation $T:X\to X$  such that $T(a) = b$. Therefore, an affine space $X$ is translation if and only if its translation group $\Trans(X)$ acts transitively on $X$.

\begin{theorem}\label{t:paraD<=>translation} A $3$-long affine regular liner $X$ is Thalesian if and only if $X$ is translation.
\end{theorem}

\begin{proof} If $\|X\|\le 2$, then $X$ is both Thalesian and translation, see Proposition~\ref{p:2=>translation}. So, assume that $\|X\|\ge 3$. In this case $X$ is an a affine space.

If $X$ is Thalesian, then every pair $xy\in X^2$ is Thalesian and hence $\overvector{xy}\in \overvector X$ is a
translation of $X$ such that $\overvector{xy}(x) = y$, witnessing that $X$ is a translation affine space.

If $X$ is a translation affine space, then for every points $x, y\in X$, there exists a translation $T\in \Trans(X) =\overvector X$ with $xy\in T$, which means that the pair $xy$ is Thalesian and the affine space $X$
 is Thalesian, by Corollary~\ref{c:Thalesian<=>pair-Thalesian}.
\end{proof}

Theorems~\ref{t:aa'-Thalesian<=>} and \ref{t:Trans(X)=vecX} imply the following characterization.
 
\begin{theorem}\label{t:aa'-translation<=>} For two points $a,a'\in X$ of an affine space $X$, there exists a translation $T:X\to X$ with $T(a)=a'$ if and only if
 $$\forall b,b',c,c'\in X\;\big((\Aline ab\parallel \Aline {a'}{b'})\;\wedge\;(\Aline bc\parallel \Aline {b'}{c'})\;\wedge\;(\Aline a{a'}\cap\Aline b{b'}=\varnothing=\Aline b{b'}\cap\Aline c{c'})\big)\;\Ra\;(\Aline a{a'}\parallel \Aline c{c'}).$$
\end{theorem}

\section{The addition of functional vectors}\label{s:addition-vectors}

Since $\overvector X=\Trans(X)$, the algebraic structure on the group of translations $\Trans(X)$ can be used to define an operation of addition of functional vectors.

\begin{theorem}\label{t:vector-addition} Let $X$ be an affine space. There exists a unique binary operation $$+:\overvector X\times\overvector X\to \overvector X,\quad  + : (\vec{\boldsymbol x},\vec{\boldsymbol y})\mapsto \vec{\boldsymbol x}+\vec{\boldsymbol y},$$
such that $\overvector{xy}+\overvector{yz}=\overvector{xz}$ for every points $x, y, z\in X$ with $\overvector{xy},\overvector{yz}\in \overvector X$. The binary operation $+$ has the following properties for every functional vectors $\vec{\boldsymbol x},\vec{\boldsymbol y},\vec{\boldsymbol z}\in \overvector X$:
\begin{enumerate}
\item $\vec{\boldsymbol x}+\vec{\boldsymbol y}=\vec{\boldsymbol y}\circ\vec{\boldsymbol x}$;
\item $(\vec{\boldsymbol x}+\vec{\boldsymbol y})+\vec{\boldsymbol z}=\vec{\boldsymbol x}+(\vec{\boldsymbol y}+\vec{\boldsymbol z})$;
\item $\vec{\boldsymbol x}+\vec{\mathbf 0}=\vec{\boldsymbol x}=\vec{\mathbf 0}+\vec{\boldsymbol x}$;
\item $\vec{\boldsymbol x}+(-\vec{\boldsymbol x})=(-\vec{\boldsymbol x})+\vec{\boldsymbol x}$.
\end{enumerate}
Therefore, $(\overvector X,+)$ is a group with neutral element $\vec{\mathbf 0}$. If for some point $o \in X$, the set $\overvector X(o)\defeq\{\vecv(o):\vecv\in \overvector X\}$ has rank $\|\overvector X(o)\|\ne 2$, then the group $(\overvector X,+)$ is commutative.
\end{theorem}

\begin{proof} By Proposition~\ref{p:Trans(X)isnormal}, $\Trans(X)$ is a subgroup of the automorphism group $\Aut(X)$ of the affine space $X$. Define the binary operation $+:\overvector X\times\overvector X\to\overvector X$ letting $\vec{\boldsymbol u}+\vec{\boldsymbol v}\defeq \vec{\boldsymbol v}\circ \vec{\boldsymbol u}$. Observe that for every points $x, y, z\in X$ with $\overvector{xy},\overvector{yz}\in \overvector X$ we have
$$(\overvector{xy}+\overvector{yz})(x) = (\overvector{yz}\circ\overvector{xy})(x) = \overvector{yz}(y) = z,$$
which implies $\overvector{xy}+\overvector{yz}=\overvector{xz}$. 
To see that the latter condition uniquely determines the binary operation $+$, assume that $\oplus:\overvector X\times\overvector X\to\overvector X$ is another binary operation on $\overvector X$ such that
\begin{equation}\label{eq:xy+yz=xz}
\overvector{xy}\oplus\overvector{yz}=\overvector{xz}
\end{equation}
for every points $x, y, z\in X$ with $\overvector{xy},\overvector{yz}\in \overvector X$. To see that $+ = \oplus$, take any vectors $\vec{\boldsymbol u}, \vec{\boldsymbol v}\in\overvector X$ and fix any point $x\in X$. Consider the points $y\defeq\vec{\boldsymbol u}(x)$ and $z := \vec{\boldsymbol v}(y)$. Then $\overvector{xy}= \vec{\boldsymbol u}\in\overvector X$ and
$\overvector {yz}=\vec{\boldsymbol v}\in\overvector X$. The equality (\ref{eq:xy+yz=xz}) ensures that
$$(\vec{\boldsymbol u}\oplus \vec{\boldsymbol v})(x) = (\overvector{xy}\oplus \overvector{yz})(x) = \overvector{xz}(x) = z = \overvector{yz}(\overvector{xy}(x)) = (\overvector{yz}\circ\overvector{xy})(x)=  (\overvector{xy} + \overvector{yz})(x) = (\vec{\boldsymbol u} + \vec{\boldsymbol v})(x)
$$
and hence the vectors $\vec{\boldsymbol u}\oplus \vec{\boldsymbol v}$ and $\vec{\boldsymbol u}+ \vec{\boldsymbol v}$  coincide.

Since $\Trans(X)$ is a group, the definition of the addition operation $+$ on $\overvector X$ ensures that $(\overvector X,+)$ is a group, which is anti-isomorphic to the group $\Trans(X)$. Taking into account that for every
functional vector $\overvector{xy}\in \overvector X$ its opposite vector $-\overvector{xy}=\overvector{yx}=\overvector{xy}^{-1}$ coincides with the inverse element $\overvector{xy}^{-1}$ to $\overvector{xy}$ in the group $\Trans(X)$, we conclude that the conditions (1)–(4) are satisfied.
If for some point $o\in  X$, the set 
$\overvector X(o)=\{\vecv(o):\vecv\in \overvector X\}=\{T(o) : T\in\Trans(X)\}$ has rank $\|\overvector X(o)\|\ne 2$, then by Baer’s Theorem~\ref{t:Trans-commutative}, the group $\Trans(X)$ is commutative and so is the group
$(\overvector X,+)$.
\end{proof}

\begin{corollary}\label{c:Thalesian-vectors-commutative} For every Thalesian affine space $X$, the vector addition has the following
properties for every vectors $\vec{\boldsymbol x},\vec{\boldsymbol y}, \vec{\boldsymbol z}\in \overvector X$:
\begin{enumerate}
\item $\vec{\boldsymbol x}+\vec{\boldsymbol y}=\vec{\boldsymbol y}+ \vec{\boldsymbol x}=\vec{\boldsymbol y}\circ \vec{\boldsymbol x}$;
\item $(\vec{\boldsymbol x}+\vec{\boldsymbol y})+\vec{\boldsymbol z}=\vec{\boldsymbol x}+(\vec{\boldsymbol y}+\vec{\boldsymbol z})$;
\item $\vec{\boldsymbol x}+\vec{\mathbf 0}=\vec{\boldsymbol x}=\vec{\mathbf 0}+\vec{\boldsymbol x}$;
\item $\vec{\boldsymbol x}+(-\vec{\boldsymbol x})=\vec{\mathbf 0}=(-\vec{\boldsymbol x})+\vec{\boldsymbol x}$.
\end{enumerate}
Therefore, $(\overvector X,+)=(X^2_{\#},+)$ is a commutative group.
\end{corollary}

\begin{corollary}\label{c:vectors-on-parallelogram} Let $a, b, x, y$ be points in a Thalesian affine space $X$. If $\overvector{ab}=\overvector{xy}\in\overvector X$, then $\overvector{ax} =\overvector{by}$.
\end{corollary}

\begin{proof} By Theorem~\ref{t:vector-addition}, the addition of functional vectors in $X$ is commutative. It follows from
$\overvector{ab}=\overvector{xy}$ that $\overvector{ba}=\overvector{yx}$. The commutativity of the funvector addition ensures that $$\overvector{ax}=
\overvector{ay}+\overvector{yx}=\overvector{ay}+\overvector{ba}=\overvector{ba}+\overvector{ay} =\overvector{by}.$$
\end{proof}

Theorem~\ref{t:vector-addition} also implies the following properties of the natural action of the group $\overvector X$ of functional vectors on the affine space $X$.

\begin{corollary}\label{c:vector-action} For every affine space $X$, the map
$$+: X \times \overvector X\to X,\quad +: (x,\vecv)\mapsto x+\vecv\defeq\vecv(x),$$
has the following properties:
\begin{enumerate}
\item $\forall x,y\in X\;\;(x+\overvector{xy}=y)$;
\item $\forall x\in X\;\;(x+\vec{\mathbf 0}=x)$;
\item $\forall x\in X\;\;\forall \vec{\boldsymbol y},\vec{\boldsymbol z}\in\overvector X\;\;\;\big(x+(\vec{\boldsymbol y}+\vec{\boldsymbol z})=(x+\vec{\boldsymbol y})+\vec{\boldsymbol z}\big)$.
\end{enumerate}
\end{corollary}

The following question actually asks whether the operation of addition of funvectors can be extended to an action of the group $\overvector X$ on the set of vectors $X^2_{\#}$.

\begin{question} Let $\vecv$ be a funvector in an affine space $X$ and $ab, xy$ be two translation
equivalent pairs in $X$. Let $b'\defeq\vecv (b)$ and $y'\defeq\vecv(y)$. Are the pairs $ab'$ and $xy'$ translation
equivalent?
\end{question}

\section{The subparallelity of vectors to flats} 

\begin{definition} Let $X$ be an affine space. Given a vector $\vec{\boldsymbol v}\in X^2_{\#}$ and a flat $A\subseteq X$,  we write\index[note]{$\subparallel$} $\vec{\boldsymbol v}\subparallel A$ and say that $\vec{\boldsymbol v}$ is \index{vector!subparallel}\defterm{subparallel} to $A$ if $\forall xy\in\vec{\boldsymbol v}\;\;\Aline xy\subparallel A$. 
\end{definition}

Proposition~\ref{p:=vectors=>parallel-sides} and Corollary~\ref{c:subparallel-transitive} imply the following characterization of the subparallelity of vectors to flats.

\begin{proposition}\label{p:vect-subparallel<=>} A vector $\vec{\boldsymbol v}$ in a proaffine space $X$ is subparallel to a flat $A\subseteq X$ if and only if there exist a pair $xy\in\vec{\boldsymbol v}$ such that $\Aline xy\subparallel A$.
\end{proposition}

\begin{proposition}\label{p:vector-subparallel} Let $A$ be a flat in a Thalesian affine space $X$, and  $\vec{\boldsymbol u},\vec{\boldsymbol v}\in\overvector{X}$ be two functional vectors. If $\vec{\boldsymbol u}\subparallel A$ and $\vec{\boldsymbol v}\subparallel A$, then $(\vec{\boldsymbol v}+\vec{\boldsymbol u})\subparallel A$.
\end{proposition}

\begin{proof} If $A=\varnothing$, then $\vec{\boldsymbol u}\subparallel A=\varnothing$ and $\vec{\boldsymbol v}\subparallel A=\varnothing$ imply $\vec{\boldsymbol u}=\vec{\mathbf 0}=\vec{\boldsymbol v}$ and $\vec{\boldsymbol u}+\vec{\boldsymbol v}=\vec{\mathbf 0}\subparallel A$.

So, assume that $A\ne\varnothing$ and fix any point $x\in A$. Let $y\defeq\vecv(x)$ and $z\defeq\vec{\boldsymbol u}(y)$. Then $\vecv+\vec{\boldsymbol u}=\overvector{xy}+\overvector{yz}=\overvector{xz}$.

The subparallelity $\vecv \subparallel A$ implies $\Aline xy\subparallel A$ and hence $y\in\overline{\{x\}\cup A}=A$. In its turn, the subparallelity $\vec{\boldsymbol u}\subparallel A$ implies $\Aline yz\subparallel A$ and hence $z\in\overline{\{y\}\cup A}=A$. Then $\Aline xz\subseteq A$ and $\Aline xz\subparallel A$, which implies $\vecv+\vec{\boldsymbol u}=\overvector{xz}\subparallel A$, by Proposition~\ref{p:vect-subparallel<=>}.
\end{proof}

\begin{corollary}\label{c:vector-subparallel} Let $A$ be a flat in a Thalesian affine space $X$ and let  $o,x,y,z\in X$ be points such that  $\overvector{ox}+\overvector{oy}=\overvector{oz}$. If $o,x,y\in A$, then $z\in A$.
\end{corollary}

\begin{proof} It follows from $o,x,y\in A$ that and Proposition~\ref{p:vect-subparallel<=>} that $\overvector{ox}\subparallel A$ and $\overvector{oy}\subparallel A$. By Proposition~\ref{p:vector-subparallel}, $\overvector{oz}\subparallel A$ and hence $z\in\overline{A\cup\{o\}}=\overline A=A$.
\end{proof}




Now we shall apply Proposition~\ref{p:vector-subparallel} to prove that line affinities in Thalesian affine spaces preserve the operation of addition of vectors.

\begin{theorem}\label{t:laf+vectors}  Let $A$ be a line affinity in a Thalesian affine space, $o,x,y,z$ be points in $\dom[A]$ and $o',x',y',z'\in\rng[A]$ be their images under the function $A$, respectively. If $\overvector{ox}+\overvector{oy}=\overvector{oz}$, then $\overvector{o'x'}+\overvector{o'y'}=\overvector{o'z'}$.
\end{theorem}

\begin{proof} Since every line affinity is the composition of finitely many line projections, it suffices to prove this theorem assuming that $A$ is a line projection. In this case there exists a line $\Lambda$ in $X$ such that $A=\{(a,b)\in \dom[A]\times\rng[A]:\Aline ab\subparallel \Lambda\}$. Let $L$ be a unique line such that $z\in L$ and $L\parallel\Lambda$. Corollary~\ref{c:parallel-transitive} implies that  $$A=\{(u,v)\in \dom[A]\times\rng[A]:\Aline uv\subparallel L\}.$$ By Proposition~\ref{p:flat-relation}, the line $L$ is not parallel to the lines $\dom[A]$ and $\rng[A]$.

It follows from $Aoxyz=o'x'y'z'$ that the flats $\Aline o{o'}$, $\Aline x{x'}$, $\Aline y{y'}$, $\Aline z{z'}$ are subparallel to the line $L$. Then the vectors $\overvector{o' o}$, $\overvector{x'x}$, $\overvector{y'y}$ and $\overvector{z'z}$ also are subparallel to $L$. Let $z''\in X$ be a unique point such that such that $\overvector{o'z''}=\overvector{o'x'}+\overvector{o'y'}$. Corollary~\ref{c:vector-subparallel}, ensures that the point $z''$ belongs to the flat $\rng[A]$. 
Observe that 
$$
\begin{aligned}
\overvector{z''z}&=\overvector{z''o'}+\overvector{o'o}+\overvector{oz}=-\overvector{o'z''}+\overvector{o'o}+\overvector{oz}\\
&=-(\overvector{o'x'}+\overvector{o'y'})+\overvector{o'o}+(\overvector{ox}+\overvector{oy})=\overvector{o'o}+(\overvector{ox}-\overvector{o'x'})+(\overvector{oy}-\overvector{o'y'})\\
&=\overvector{o'o}+(\overvector{oo'}+\overvector{o'x'}+\overvector{x'x}-\overvector{o'x'})+(\overvector{oo'}+\overvector{o'y'}+\overvector{y'y}-\overvector{o'y'})\\
&=(\overvector{o'o}+\overvector{oo'}+\overvector{oo'})+\overvector{x'x}+\overvector{y'y}=\overvector{oo'}+\overvector{x'x}+\overvector{y'y}\subparallel L,
\end{aligned}
$$
by Proposition~\ref{p:vector-subparallel}.   It follows from $\overvector{z'z}\subparallel L$, $\overvector{z''z}\subparallel L$ and $z\in L$ that  $\{z',z''\}\subseteq L$, see Corollary~\ref{c:vector-subparallel}. Taking into account that the lines $\rng[A]$ and $L$ are not parallel, we conclude that $L\cap\rng[A]$ is a singleton and so is the set $\{z',z''\}\subseteq L\cap\rng[A]$. Then $z'=z''$ and hence $\overvector{o'z'}=\overvector{o'z''}=\overvector{o'x'}+\overvector{o'y'}$.
\end{proof} 

\begin{definition}\label{d:vectors-parallel} Given two vectors $\vecbold{v},\vecbold{u}\in\overvector{X}$ in an affine space $X$ we write\index[note]{$\subparallel$}\index[note]{$\parallel$} $\vecbold{v}\subparallel\vecbold u$ (resp. $\vecbold{v}\parallel\vecbold u$) and say that the vector $\vecbold{v}$ is \index{subparallel vectors}\index{vectors!subparallel}\index{parallel vectors}\index{vectors!parallel}\defterm{subparallel} (resp. \defterm{parallel\/}) to the vector $\vecbold{u}$ if $\Aline xv\subparallel \Aline yu$ (resp. $\Aline xv\parallel \Aline yu$)  for all pairs $xv\in\vecbold{v}$ and $yu\in\vecbold{u}$.
\end{definition}

Proposition~\ref{p:=vectors=>parallel-sides} and Corollary~\ref{c:parallel-transitive} imply the following characterization.

\begin{corollary} Let $X$ be an affine space. A vectors $\vecbold{v}\in X^2_{\#}$ is (sub)parallel to a vector $\vecbold{u}\in X^2_{\#}$  if and only if the flat $\Aline xv$ is (sub)parallel  to the flat $\Aline yu$ for some pairs $xv\in\vecbold{v}$ and $yu\in\vecbold{u}$.
\end{corollary}

\section{Transformation groups preserving vectors}

The notion of subparallelity of vectors allows us to define a natural structure of a liner on the group of functional vectors $\overvector X$ in an affine space $X$.
Namely, for every affine space $X$, endow the set of funvector
$\overvector X=\Trans(X)$ with the line relation
$$\vec \Af\defeq\{\vec{\boldsymbol x},\vec{\boldsymbol y}, \vec{\boldsymbol z}\in \overvector X^3: \vec{\boldsymbol x}\ne \vec{\boldsymbol z}\;\;\wedge\;\;(\vec{\boldsymbol y}-\vec{\boldsymbol x})\subparallel (\vec{\boldsymbol z}-\vec{\boldsymbol x})\}.$$

\begin{exercise} Check that $\vec\Af$ is indeed a line relation on
$\overvector X$.
\end{exercise}

\begin{exercise}\label{ex:overvector-liner} Show that for every point $o$ in an affine space $X$, the function $\overvector X\to X$, $\vecv\to\vecv(o)$, is an isomorphism of the liner $(\overvector X,\vec\Af)$ onto the subliner $\overvector X(o)\defeq\{\vecv(o):\vecv\in\overvector X\}$ of the affine space $X$. Show that $\|\overvector X\|\ge 3$ if and only if $\|\overvector X(o)\|\ge 3$.
\end{exercise}

\begin{remark} The liner $\overvector X$ can have quite intricated geometry. For example, there exists an affine plane $X$ of order $16$ for which the liner  $\overvector X$ has 32 points and contains four parallel lines of length $8$, $32$ lines of length $4$ and $256$ lines of length $2$.
\end{remark}

Next, we consider three transformation groups on an affine space $X$, which coincide with the
group of translations $\Trans(X)$, if the affine space $X$ is Thalesian. Those transformation groups
consist of linevections, vections and funvections, respectively.

\begin{definition} An automorphism $A : X \to X$ of an affine space $X$ is defined to be
\begin{itemize}
\item a \index{linevection}\index{automorphism!linevection}\defterm{linevection} if for every line $L\subset X$, the restriction $A{\restriction}_L$ is a line translation;
\item a  \index{vection}\index{automorphism!vection}\defterm{vection} if $Axy\#xy$ for every pair $xy\in X^2$;
\item a  \index{funvection}\index{automorphism!funvection}\defterm{funvection} if $Axy\#xy$ for every Thalesian pair $xy\in  X^2$.
\end{itemize}
\end{definition}

\begin{proposition}\label{p:TVFT} Let $A:X\to X$ be an automorphism of an affine space $X$.
\begin{enumerate}
\item If $A$ is a linevection, then $A$ is a vection;
\item If $A$ is a vection, then $A$ is a funvection;
\item If $A$ is a vection, then $A$ is a dilation;
\item $A$ is a funvection if and only if $AT = TA$ for every translation $T\in\Trans(X)$.
\item If $\|\overvector X\|\ne 0$ and $A$ is a funvection and a dilation, then $A$ is a translation.
\item If $A$ is a translation and $\|\overvector X\|\ne 2$, then $A$ is a linevection.
\item If $\|\overvector X\|\ge 3$, then $A$ is a translation if and only if $A$ is a linevection if and only if it is
a vection if and only if $A$ is a funvection and a dilation.
\item If $X = \bigcup_{T\in\Trans(X)}\overline{\{o, T(o)\}}$ for some $o\in X$ and $A$ is a funvection, then $A$ is a dilation.
\end{enumerate}
\end{proposition}

\begin{proof} Let $A:X\to X$ be an automorphism of an affine space $X$. 
\smallskip

1. Assume that $A$ is a linevection. To prove that $A$ is a vection, we need to show that $Axy\in\overvector{xy}$ for every pair $xy\in A$. Choose any line $L\subset  X$ containing the points $x,y$. Since $A$ is linevection, the restriction $T\defeq A{\restriction}_L$ is a line translation with $Axy = Txy\#xy$, witnessing that
$A$ is a vection.
\smallskip

2. If $A$ is a vection, then $A$ is a funvection because every funvector is a vector.
\smallskip

3. Assume that $A$ is a vection of $X$. To prove that $A$ is a dilation, take any line $L\subset X$. Choose
any distinct points $x, y\in L$. Consider the pair $uv = Axy$ and observe that $A[L] = A[\Aline x y] = \Aline u v$.
Since $A$ is a vection of $X$, $uv = Axy\#xy$. By Proposition~\ref{p:=vectors=>parallel-sides}, $uv\#xy$ implies $L = \Aline x y \parallel\Aline u v =
A[L]$, witnessing that the automorphism $A$ is a dilation of $X$.
\smallskip

4. Assume that $A$ is a funvection of $X$. Given any translation $T$ of $X$, we should prove that $AT = TA$. This equality is trivially true if $T$ is the identity translation of $X$. So, assume that
$T\ne 1_X=\vec{\mathbf 0}$. Take any point $x\in A$ and concider the point $y\defeq T(x)$. Since $\overvector{xy} = T\ne 1_X=\vec{\mathbf 0}$, the point $y$ is not equal to $x$ and hence $\Aline xy$ is a line. By Theorem~\ref{t:Trans(X)=vecX}, $\overvector{xy}=T\in\Trans(X) =\overvector X$ and hence the pair $xy$ is Thalesian. Since $A$ is a funvection, the pair $uv\defeq Axy$ is translation
equivalent to the Thalesian pair $xy$. Then there exists a line translation $S:\Aline xy \to\Aline uv$ such that
$uv = Sxy$. Theorem~\ref{t:TS=ST} ensures that $TS = ST$ and hence
$$AT(x) = A(y) = v = S(y) = ST(x) = TS(x) = T(u) = TA(x),$$
witnessing that $AT = TA$.

Now assuming that $AT = TA$ for every translation $T\in\Trans(X)$,  we shall prove that the automorphism $A$ is a funvection. Given any Thalesian pair $xy\in X^2$, we should prove that
$Axy\in\overvector{xy}$. By Theorem~\ref{t:Trans(X)=vecX}, the vector $T\defeq\overvector{xy}$ is a translation of $X$. Consider the pair
$uv\defeq Axy$. By our assumption, $AT=TA$ and hence $v=A(y) = AT(x) = TA(x) = T(u)$ and
hence $Axy = uv \in T = \overvector{xy}$.
\smallskip

5. Assume that $\|\overvector X\|\ne0$ and $A$ is a funvection and a dilation. Assuming that the dilation $A$ is not a translation, we can consider its unique fixed point $x = A(x)$. Since $\|\overvector X\|\ne0$, there exists a nonzero vector $\vecv\in \overvector X\setminus\{\mathbf 0\}$. Let $y\defeq\vecv(x)\ne x$ and $a\defeq A(y)$. Since $A$ is a funvection, $xa=Axy\#xy$. Then $A(y)=a=\overvector{xa}(x)=\overvector{xy}(x)=y$ and hence the dilation $A$ has two fixed points $x$ and $y$. By Proposition~\ref{p:Ax=Bx=>A=B},
$A$ is the identity translation, which contradicts our assumption. This contradiction shows that the dilation $A$ is a translation.
\smallskip

6. If $A$ is a translation and $\|\overvector X\|\ne 2$, then $A$ is a linevection, by Theorem~\ref{t:Trans(X)=vecX} and Proposition~\ref{p:translation|=>line-translation}.
\smallskip

7. Assume that $\|\overvector X\|\ge 3$ for some point $o\in X$. If $A$ is a translation, then $A$ is a linevection,
by Proposition~\ref{p:TVFT}(6). If $A$ is a linevection, then $A$ is a vection, by Proposition~\ref{p:TVFT}(1). If $A$ is a vection, then $A$ is a funvection and dilation, by Proposition~\ref{p:TVFT}(2,3). If $A$ is a funvection and dilation, then $A$ is a translation, by Proposition~\ref{p:TVFT}(5).
\smallskip

8. Assume that $X =\bigcup_{T\in\Trans(X)}\overline{\{o, T(o)\}}$ for some point $o\in X$. Since $\| X\|\ge 3$, this equality implies 
$\|\overvector X(o)\|=\|\{T(o):T\in\Trans(X)\}\|=\|X\|\ge 3$. Assume that $A$ is a funvection. To prove that $A$ is a dilation,
take any line $L\subset X$. Since $X$ is an affine space, there exists a line $L_o$ such that $o\in L_o$ and $L_o\parallel L$.
Since $X =\bigcup_{T\in\Trans(X)}\overline{\{o, T(o)\}}$, there exists a translation $T\in\Trans(X)$ such that $T(o)\in L_o \setminus\{o\}$. Then for the point $o'\defeq T(o)$, the pair $oo'\in  T\in\Trans(X) =\overvector X$ is Thalesian. By Corollary~\ref{c:transeq-exists},
there exists a point $y\in L$ such that $xy\#oo'$ and hence $\overvector{xy}=\overvector{oo'}\in \overvector X$. Since $A$ is a funvection, $Axy\#xy$ and hence $A[L] = A[\Aline x y]\parallel \Aline x y = L$, witnessing that the automorphism $A$ is a dilation.
\end{proof}

For an affine space $X$, let $\Trans_L(X), \Trans_V(X)$, and $\Trans_{F}(X)$ denote the subsets of the group $\Aut(X)$, consisting of linevections, vections, and funvections, respectively. It is easy to see
that $\Trans_L(X), \Trans_V(X)$ and $\Trans_{F}(X)$ are normal subgroups of the automorphism group $\Aut(X)$.
Proposition~\ref{p:TVFT} implies that
$$
\begin{aligned}
&\Trans_L(X)\subseteq \Trans_V(X)\subseteq \Trans_{F}(X)\cap \Dil(X)\quad\mbox{and}\\
&\Trans_{F}(X)=Z(\Trans(X))\defeq\{A\in \Aut(X):\forall T\in\Trans(X),\;\; TA = AT\}.
\end{aligned}
$$

Proposition~\ref{p:TVFT} also implies the following corollaries.

\begin{corollary} Let $X$ be an affine space such that $\|\overvector X\|\ge 3$. Then
$$\overvector X= \Trans(X) =\Trans_L(X)=\Trans_V(X)=\Trans_{F}(X)\cap\Dil(X)\subseteq \Trans_{F}(X) = Z(\Trans(X)).
$$
\end{corollary}

\begin{corollary} Let $X$ be an affine space such that $X=\bigcup_{T\in\Trans(X)}\overline{\{o, T(o)\}}$ for some
$o\in X$. Then
$$\overvector X=\Trans(X)=\Trans_L(X)=\Trans_V(X)=\Trans_{F}(X)=Z(\Trans(X)).
$$
\end{corollary}

\begin{corollary} If an affine space X is Thalesian, then
$$\overvector X=\Trans(X)=\Trans_L(X)=\Trans_V(X)=\Trans_{F}(X)=Z(\Trans(X)).
$$
\end{corollary}

\begin{remark}For non-Thalesian affine planes the groups $\Trans(X),\Trans_L(X),\Trans_V(X),\Trans_{F}(X)$ can be distinct, as witnessed by following table that shows the cardinalities of those groups, calculated by Ivan \index[person]{Hetman}Hetman\footnote{{\tt https://docs.google.com/spreadsheets/d/1a0uWy1gdCPF4Z7xzkYivSwxG62oNUbZ1vyPCmZlVQJU/edit\#gid=0}} for all seven (pairwise nonisomorphic) affine planes of order 9.

\begin{center}
\begin{tabular}{|c|c|c|c|c|c|c|c|c|}
\hline
$X$&\phantom{\Large$|^{|^|}_{|}$}$\!\!\!\!\!\|\overvector X\|$&$\Trans(X)$&$\Trans_L(X)$&$\Trans_V(X)$&$\Trans_F(X)$&$\Dil(X)$&$\Aut(X)$&$\Aut(\overline X)$\\
\hline
{\tt Desarg}&3&81&81&81&81&648&933120&84913920\\
{\tt Thales}&3&81&81&81&81&162&311040&311040\\
{\tt Hall}&2&9&9&9&9&72&311040&311040\\
{\tt Hughes}&3&9&9&9&9&18&2592&33696\\
{\tt dhall}&1&1&?&8&3456&8&3456&311040\\
{\tt hall}&1&1&?&2&3840&2&3840&311040\\
{\tt hughes}&1&1&1&1&432&1&432&33696\\
\hline
\end{tabular}
\end{center}
\end{remark}

\chapter{Line affinities and homotheties in affine spaces}

In this chapter we continue studying line affinities, started in Chapter~\ref{ch:LTrans}. Let us recall that a bijection $A:L\to\Lambda$ between two lines in an affine space $X$ is called a {\em line affinity} if $A=P_1\cdots P_n$ for some line projections $P_1,\dots,P_n$. A bijection $P:L\to L'$ between two lines $L,L'$ is called a {\em line projection} if there exists a line $\Lambda$ such that $P=\{(x,y)\in L\times L':\Aline xy\subparallel \Lambda\}$. A line projection $P:L\to L'$ is called a \index{line shift}\index{line!shift}\defterm{line shift} if $L\parallel L'$, and a \index{line shear}\index{line!shear}\defterm{line shear} if $L\cap L'=\{o\}$ for a unique point $o$, called the \index{shear center}\defterm{shear center}.

\section{Line affinities in affine spaces}

In this section we reveal the structure of line affinities in affine spaces.

\begin{proposition}\label{p:shear=shift+shear+shift} For every line shear $R$ in an affine space $X$ and every point $o'\in X$, there exist line shifts $S,S'$ and a line shear $R'$ such that $R=S'R'S$ and $R'(o')=o'$. 
\end{proposition}

\begin{proof} Since the line shear $R$ is a line projection, there exists a line $\Lambda$ in $X$ such that $$R=\{(x,y)\in\dom[R]\times\rng[R]:\Aline xy\subparallel\Lambda\}.$$ By Proposition~\ref{p:flat-relation}, the line $\Lambda$ is not parallel to the lines $\dom[R]$ and $\rng[R]$. Since the affine space $X$ is Playfair, there exist unique lines $L,L'$ in $X$ such that $o'\in L\cap L'$, $L\parallel\dom[R]$ and $L'\parallel \rng[R]$. By Theorem~\ref{t:subparallel-via-base}, the plane $\Pi\defeq\overline{\dom[R]\cup\rng[R]}$ is parallel to the plane $\overline{L\cup L'}$ and hence the line $\Lambda$ is subparallel to the plane $\overline{L\cup L'}$ but is not parallel to the lines $L,L'$. Applying Proposition~\ref{p:flat-relation}(1,3) and Corollary~\ref{c:parallel-lines<=>}, we conclude that the relation
$$R'\defeq\{(x,y)\in L\times L':\Aline xy\subparallel \Lambda\}$$is a line shear with $\dom[R']=L$, $\rng[R']=L'$, and $R'(o')=o'$. 

If the point $o'$ belongs to the plane $\Pi$, then the line shifts 
$$S\defeq\{(x,y)\in \dom[R]\times L:\Aline xy\subparallel\Lambda\}\quad\mbox{and}\quad S'\defeq\{(x,y)\in L'\times \rng[R]:\Aline xy\subparallel \Lambda\}$$have the desired property $R=S'R'S.$

So, assume that $o\notin \Pi$ and hence $\|X\|>\|\Pi\|=3$.
By Corollaries~\ref{c:affine-Desarguesian} and Theorem~\ref{t:ADA=>AMA}, the affine space $X$ is Desarguesian and Thalesian. Let $o$ be the center of the line shear $R$. By Corollary~\ref{c:par-trans}, there exist line translations $S:\dom[R]\to L$ and $S':L'\to\rng[R]$ such that $S(o)=o'$ and $S'(o')=o$.

We claim that $R=S'R'S$. Given any $x\in \dom[R]$, we should check that $R(x)=S'R'S(x)$. If $x=o$, then $S'R'S(x)=S'R'S(o)=S'R'(o')=S'(o')=o=R(o)=R(x)$ and we are done.
 So, assume that $x\notin o$. Consider the points $y\defeq R(x)\in\rng[R]\setminus\{o\}$, $x'\defeq S(x)$, and $y'\defeq (S')^{-1}(y)$. Taking into account that $S,S'$ are line shifts, we conclude that $\Aline x{x'}\parallel \Aline o{o'}\parallel \Aline y{y'}$. Since $\Aline ox=\dom[R]\parallel L=\Aline{o'}{x'}$, $\Aline oy=\rng[R]\parallel L'=\Aline{o'}{y'}$, the Thalesian property of $X$ ensures that $\Aline {x'}{y'}\parallel \Aline xy\subparallel\Lambda$ and hence $R'(x')=y'$ and $R(x)=y=S'(y')=S'R'(x')=S'R'S(x)$. 
 
\begin{picture}(150,120)(-150,-15)

{\linethickness{0.75pt}
\put(0,0){\color{blue}\vector(1,0){100}}
\put(0,0){\color{cyan}\line(0,1){90}}
\put(0,0){\color{teal}\line(-1,2){30}}
\put(100,0){\color{cyan}\line(0,1){90}}
\put(100,0){\color{teal}\line(-1,2){30}}
\put(0,90){\color{blue}\vector(1,0){100}}
\put(70,60){\color{red}\vector(1,1){30}}
\put(-30,60){\color{red}\vector(1,1){30}}
\put(-30,60){\color{blue}\vector(1,0){100}}
}

\put(0,0){\circle*{3}}
\put(-2,-8){$o$}
\put(0,90){\circle*{3}}
\put(-2,95){$y$}
\put(100,0){\circle*{3}}
\put(98,-8){$o'$}
\put(-30,60){\circle*{3}}
\put(-38,57){$x$}
\put(70,60){\circle*{3}}
\put(73,57){$x'$}
\put(100,90){\circle*{3}}
\put(98,95){$y'$}

\put(-27,73){\color{red}$R$}
\put(27,63){\color{blue}$S'$}
\put(47,93){\color{blue}$S$}
\put(73,73){\color{red}$R'$}
\end{picture}
\end{proof}

\begin{corollary}\label{c:affinity=TRT} Let $A$ be a line affinity in an affine space $X$. For every point $o\in X$, there exist line translations $T_0,T_1,\dots,T_n$ and line shears $R_1,\dots,R_n$ in $X$ such that 
$$A=T_nR_nT_{n-1}\cdots T_1R_1T_0$$and $R_i(o)=o$ for all $i\in \{1,\dots,n\}$.
\end{corollary}

\begin{proof} Write the line affinity $A$ as the composition $P_n\dots P_1$ of line projections. Since every line shift is a composition of two line shears, we lose no generality assuming that every line projection $P_i$ is a line shear. By Proposition~\ref{p:shear=shift+shear+shift}, for every $i\in\{1,2,\dots,n\}$ there exist line shifts $S_i,S_i'$ and a line shear $R_i$ such that $P_i=S'_iR_iS_i$. Put $T_0\defeq S_1$, $T_n\defeq S_n'$ and $T_i\defeq S_{i+1}S'_{i}$ for every $i\in\{1,\dots,n-1\}$. Then $T_0,\dots,T_n$ are line translations such that
$$
\begin{aligned}
A&=P_n\cdots P_1=(S'_nR_nS_n)(S'_{n-1}R_{n-1}S_{n-1})\cdots (S_1'R_1S_1)\\
&=S'_nR_n(S_nS'_{n-1})R_{n-1}(S_{n-1}S'_{n-1})\cdots R_2(S_2S_1')R_1S_1=T_nR_nT_{n-1}R_{n-1}\cdots T_1R_1T_0.
\end{aligned}
$$
\end{proof}

\section{Line affinities in Thalesian affine spaces}

\begin{proposition}\label{p:RTR=T} Let $R$ be a line shear in a Thalesian affine space $X$. For every line translation $T:\dom[R]\to\dom[R]$ the function $RTR^{-1}:\rng[R]\to\rng[R]$ is a line translation.
\end{proposition}

\begin{proof} Since $R$ is a line shear, there exists a line $\Lambda$ in $X$ such that $$R=\{(x,y)\in\dom[R]\times\rng[R]:\Aline xy\subparallel \Lambda\}.$$ By Proposition~\ref{p:flat-relation}(2,4), the line $\Lambda$ is not parallel to the lines $\dom[R]$ and $\rng[R]$. Let $o\in \dom[R]\cap\rng[R]$ be the center of the shear $R$, and $\Lambda_os$ be a unique line such that $o\in\Lambda_o$ and $\Lambda_o\parallel \Lambda$. By Proposition~\ref{p:subparallel-char4}, the line $\Lambda$ is subparallel to the plane $\Pi\defeq\overline{\dom[R]\cup\rng[R]}$ and hence $\Lambda_o\subset \Pi$.

If the point $o'\defeq T(o)$ is equal to the point $o$, then the line translation $T$ is the identity map of the line $\dom[R]$, by Theorem~\ref{t:unique-translation}. In this case, $RTR^{-1}=RR^{-1}$ is the identity translation of the line $\rng[R]$ and we are done. 

So, assume that $o'\ne o$ and hence $o'\in\dom[R]\setminus\rng[R]$. Let $o''\in \Lambda_o$ be a unique point such that $\Aline{o'}{o''}\parallel \rng[R]$. By Theorem~\ref{t:paraD<=>translation}, the Thalesian affine space $X$ is a translation space. So, there exist translations $T',T'':X\to X$ such that $T'(o)=o'$ and $T''(o)=o''$. Then $\Lambda_o=\Aline o{o''}\in\overline{T''}\defeq\{\Aline xy:xy\in T''\}$ and hence $\Aline xy\parallel \Lambda$ for every pair $xy\in T''$. By Proposition~\ref{p:translation|=>line-translation}, the restriction $T'{\restriction}_{\dom[R]}$ 
is a line translation, equal to the line translation $T$.

 By Proposition~\ref{p:translation=>shift}, the restrictions $S'\defeq T'{\restriction}_{\rng[R]}$ and $S\defeq T''{\restriction}_{\rng[R]}$ are line shifts. 
It follows from $S[\rng[R]]\parallel \rng[R]\parallel S'[\rng[R]]$ and $\Aline {o'}{o''}\parallel \rng[R]$ that $\Aline {o'}{o''}=S[\rng[R]]=S'[\rng[R]]$. Then $S^{-1}S'=(T'')^{-1}T'{\restriction}_{\rng[R]}$ is a line translation. It remains to prove that $RTR^{-1}=S^{-1}S'$. Given any point $x\in\rng[R]$, we should check that $RTR^{-1}(x)=S^{-1}S'(x)$. Consider the points $y\defeq R^{-1}(x)$, $z\defeq T(y)$, $a\defeq S'(x)=T'(x)$ and $b\defeq S^{-1}(a)$.

If $x=o$, then $y=R^{-1}(x)=R^{-1}(o)=o$,  $$z=T(y)=T(o)=o'=T'(o)=S'(o)=S'(x)=a=S(b),$$ which implies $\Aline bz\parallel \Aline o{o''}\subseteq \Lambda_o$ and hence $b=R(z)$ and 
$$RTR^{-1}(x)=RT(y)=R(z)=b=S^{-1}(a)=S^{-1}S'(x).$$

\begin{picture}(200,140)(-180,-20)

\put(-80,0){\color{red}\line(1,0){160}}
\put(85,-3){\color{red}$\Lambda_o$}
\linethickness{0.75pt}

\put(-20,0){\color{red}\line(1,0){40}}
\put(20,0){\color{cyan}\line(1,1){100}}
\put(-20,0){\color{cyan}\line(1,1){100}}
\put(20,0){\color{blue}\line(-1,1){100}}
\put(80,60){\color{red}\vector(-1,0){120}}
{\linethickness{1pt}
\put(20,0){\color{red}\vector(-1,0){40}}
\put(60,80){\color{red}\vector(1,0){40}}
\put(20,0){\color{blue}\vector(-1,1){20}}
\put(-40,60){\color{blue}\vector(-1,1){20}}
\put(-60,80){\color{red}\vector(1,0){160}}
\put(80,60){\color{blue}\vector(-1,1){20}}
}

\put(110,105){\color{cyan}$\rng[R]$}
\put(-96,105){\color{blue}$\dom[R]$}

\put(20,0){\circle*{3}}
\put(20,-7){$o$}
\put(11,11){\color{blue}$T$}
\put(-20,0){\circle*{3}}
\put(-25,3){$o''$}
\put(-4,-10){\color{red}$T''$}
\put(0,20){\circle*{3}}
\put(-2,24){$o'$}
\put(80,60){\circle*{3}}
\put(81,53){$x$}
\put(74,68){\color{blue}$S'$}
\put(100,80){\circle*{3}}
\put(103,75){$b$}
\put(76,83){\color{red}$S^{-1}$}
\put(60,80){\circle*{3}}
\put(55,83){$a$}
\put(-40,60){\circle*{3}}
\put(-48,56){$y$}
\put(-46,68){\color{blue}$T$}
\put(-60,80){\circle*{3}}
\put(-68,76){$z$}
\put(10,83){\color{red}$R$}
\put(5,50){\color{red}$R^{-1}$}

\end{picture}

Next, assume that $x\ne o$. In this case $y=R^{-1}(x)\ne x$ and $\Aline xy\parallel \Lambda$.
 Since $X$ is a translation liner, there exists a translation $F:X\to X$ such that $F(x)=y$. By Corollary~\ref{c:Trans-commutative}, the group $\Trans(X)$ is commutative and hence $FT'=T'F$. Then $$z=T(y)=T'(y)=T'F(x)=FT'(x)=FS'(x)=F(a)$$ and hence $\Aline az\parallel \Aline xy\parallel \Lambda\parallel \Aline ab$ and $\Aline zb\parallel \Lambda$. Since $(z,b)\in\dom[R]\times\rng[R]$, the definition of the line shear $R$ ensures that $b=R(z)$ and hence
 $$S^{-1}S'(x)=S^{-1}(a)=b=R(z)=RT(y)=RTR^{-1}(x).$$
\end{proof}

\begin{theorem}\label{t:A=T'RRT} Let $A$ be a line affinity in a Thalesian affine space $X$. For every point $o\in X$, there exist line translations $T,T'$ and line shears $R_1,\dots,R_n$ in $X$ such that 
$$A=T'R_nR_{n-1}\cdots R_2R_1T,$$ and $R_i(o)=o$ for all $i\in \{1,\dots,n\}$. Moreover, if $o\in\dom[A]$, then $T$ is the identity translation of $\dom[A]$.  
\end{theorem}

\begin{proof} By Corollary~\ref{c:affinity=TRT},  $A=T_nR_nT_{n-1}\cdots T_1R_1T_0$ for some line translations $T_0,T_1,\dots,T_n$ and line shears $R_1,\dots,R_n$ in $X$ such that $R_i(o)=o$ for all $i\in\{1,\dots,n\}$. If $o\in\dom[A]$, then we can additionally assume that $T_0$ is the identity translation of the line $\dom[A]$. For every $i\in\{1,\dots,n-1\}$ we have $ \rng[R_i]=\dom[T_i]\parallel\rng[T_i]=\dom[R_{i+1}]$. Taking into account that $o\in\rng[R_i]\cap\dom[R_{i+1}]$, we conclude that $\rng[R_i]=\dom[T_i]=\rng[T_i]=\dom[R_{i+1}]$. 

If $o\in \dom[A]=\dom[T_0]$, then also $\dom[T_0]=\rng[T_0]$.  

 Let $T\defeq T_0$, $T'_1\defeq T_1$ and for every $i\in\{2,\dots,n\}$, let $T_i'\defeq T_iR_iT_{i-1}'R_i^{-1}$. By Proposition~\ref{p:RTR=T}, the function $T_i':\rng[R_i]\to\rng[R_i]$ is a line translation. By induction we shall prove that for every $k\in\{1,\dots,n-1\}$, the following equality holds.
\begin{itemize} 
\item[$(*_k)$] $A=T_nR_nT_{n-1}\cdots T_{k+1}R_{k+1}T'_{k}R_kR_{k-1}\cdots R_1T.$
\end{itemize}

For $k=1$ this equality follows from the equalities $T=T_0$ and $T_1'=T_1$. Assume that for some $k\in\{2,\dots,n-1\}$, the equality $(*_{k-1})$ holds.
Then
$$
\begin{aligned}
A&=T_nR_nT_{n-1}\cdots T_{k+1}R_{k+1}T_{k}R_{k}T'_{k-1}R_{k-1}\cdots R_1T\\
&=T_nR_nT_{n-1}\cdots T_{k+1}R_{k+1}T_kR_kT'_{k-1}R_k^{-1}R_kR_{k-1}\cdots R_1T\\
&=T_nR_nT_{n-1}\cdots T_{k+1}R_{k+1}T_k'R_kR_{k-1}\cdots R_1T,
\end{aligned}
$$
witnessing that the equality $(*_k)$ holds. 

Then for the line translation $T'\defeq T_n'$, the equality $(*_{n-1})$ implies the desired equality
$$
\begin{aligned}
A&=T_nR_nT'_{n-1}R_{n-1}\cdots R_2R_1T\\
&=T_nR_nT'_{n-1}R_n^{-1}R_nR_{n-1}\cdots R_2R_1T\\
&=T_n'R_nR_{n-1}\dots R_1T=T'R_nR_{n-1}\cdots R_1T.
\end{aligned}
$$
If $o\in\dom[A]$, then $T=T_0$ is the identity translation of the line $\dom[A]$.
\end{proof}

\section{Line shears in Desarguesian affine spaces}

In this section we study the structure of compositions of line shears around the same point in a Desarguesian affine spaces. By Theorem~\ref{t:Desargues-affine} and Lemma~\ref{l:ADA<=>}, an affine space $X$ is  Desarguesian if and only if for every concurrent lines $A,B,C$ in $X$ and points $a,a'\in A\setminus (B\cup C)$,  $a,a'\in B\setminus (A\cup C)$,  $a,a'\in C\setminus (A\cup B)$, if $\Aline ab\parallel \Aline {a'}{b'}$ and $\Aline bc\parallel \Aline {b'}{c'}$, then $\Aline ac\parallel \Aline{a'}{c'}$. We recall that lines $A,B,C$ in a liner are {\em concurrent} if $L_1\cap L_2\cap L_3$ is a singleton. 

The following theorem is a ``shear'' counterpart of Theorem~\ref{t:shift=AMA}.

\begin{theorem}\label{t:shear=ADA} An affine space $X$ is Desarguesian if and only if for any concurrent lines $L_1,L_2,L_3$ with $L_1\ne L_3$ and every line shears $P:L_1\to L_2$ and $R:L_2\to L_3$, the composition $R P:L_1\to L_3$ is a line shear.
\end{theorem} 

\begin{proof}  To prove the ``only if'' part, assume that the affine space $X$ is Desarguesian. Given any
line shears $P:L_1\to L_2$ and $R:L_2\to L_3$ between concurrent lines $L_1,L_2,L_3$ with $L_1\ne L_3$, we should check that the composition $R P:L_1\to L_3$ is a line shear. 
Let $o$ be the unique point of the intersection $L_1\cap L_2\cap L_3$. 

The definition of a line shear ensures that $L_1\cap L_2$ and $L_2\cap L_3$ are singletons and hence $L_1\cap L_2=\{o\}=L_2\cap L_3$. Therefore, $L_1,L_2,L_3$ are three distinct concurrent lines in $X$.
By Proposition~\ref{p:para-projection}(2), $P(o)=o=R(o)$. 

 Fix any point $a_1\in L_1\setminus\{o\}$, and consider the points $a_2=P(a_1)\in L_2\setminus\{o\}$ and $a_3=R(a_2)\in L_3\setminus\{o\}$. By Theorem~\ref{t:paraproj-exists}, there exists a line shear $T:L_1\to L_3$ such that  $T(a_1)=a_3$. By Proposition~\ref{p:para-projection}(2), $T(o)=o$.  We claim that that $R P=T$. Given any point $x_1\in L_1$, we should prove that $RP(x_1)=T(x_1)$. If $x_1=o$, then $RP(x_1)=RP(o)=R(o)=o=T(o)=T(x_1)$. So, we assume that $x_1\ne o$.
Consider the points $x_2=P(x_1)\in L_2\setminus\{o\}$, $x_3=R(x_2)\in L_3\setminus\{o\}$, and $x_3'=T(x_3)\in L_3\setminus\{o\}$. Applying Theorem~\ref{t:Desargues-affine} and Lemmas~\ref{l:D=>ADA},  \ref{l:ADA<=>}, we conclude that $\Aline{x_1}{x_3}\parallel \Aline{a_1}{a_3}\parallel \Aline{x_1}{x_3'}$ and hence $\Aline{x_1}{x_3}=\Aline{x_1}{x_3'}$, by the Proclus Axiom (holding for affine spaces). Then $\{x_3\}=L_3\cap\Aline{x_1}{x_3}=L_3\cap\Aline{x_1}{x_3'}=\{x_3'\}$ and hence $T(x_1)=x_3'=x_3=RP(x_1)$, witnessing that $R P=T$ is a line shear. 
\smallskip 

To prove the ``if'' part, assume that  for any line shears $P:L_1\to L_2$ and $R:L_2\to L_3$ between concurrent lines $L_1,L_2,L_3$ with $L_1\ne L_3$, the composition $R P:L_1\to L_3$ is a line shear. To check that the affine space $X$ is Desarguesian, take any concurrent lines $A,B,C$ in $X$ and points $a,a'\in A\setminus(B\cup C)$, $b,b'\in B\setminus(A\cup C)$, $c,c'\in C\setminus(A\cup B)$ such that $\Aline ab\parallel \Aline{a'}{b'}$ and $\Aline bc\parallel \Aline{b'}{c'}$. The choice of the points $a,b,c$ ensures that the lines $A,B,C$ are distinct. By Theorem~\ref{t:paraproj-exists}, there exist line shears $P:A\to B$ and $R:B\to C$ such that $P(a)=b$, and $R(b)=c$. By our assumption, the composition $R P$ is a line shear. Taking into account that $\Aline {a'}{b'}\parallel \Aline ab$ and $\Aline{b'}{c'}\parallel \Aline bc$, we conclude that $R(a')=b'$ and $P(b')=c'$. Then $R P(a')=c'$. Since $R P$ is a line projection with $\{(a,c),(a',c')\}\in R P$, the definition of a parallel projection and Corollary~\ref{c:parallel-transitive} imply that $\Aline {a'}{c'}\parallel \Aline ac$, witnessing that the affine space $X$ is Desarguesian.
\end{proof}

\begin{lemma}\label{l:shear=shear+shear} Let $R,R'$ be line shears around a point $o$ in a Desarguesian affine space $X$. If $\dom[R_1]=\dom[R]$ and $\rng[R_1]\ne\rng[R]$, then there exists a line shear $R_2$ such that $R=R_2 R_1$.
\end{lemma}

\begin{proof} Fix any point $x\in \dom[R]=\dom[R_1]$ with $x\ne o$, and consider the points $y\defeq R(x)$ and $y_1\defeq R_1(x)$. By Theorem~\ref{t:paraproj-exists}, there exists a line shear $R_2:\rng[R_1]\to\rng[R]$ such that $R_2oy_1=oy$. By Theorem~\ref{t:shear=ADA}, the composition $R_2 R_1$ is a line shear with $R_2 R_1 ox=R_2oy_1=oy=Rox$, and by Theorem~\ref{t:paraproj-exists}, $R=R_2 R_1$.
\end{proof}

\begin{lemma}\label{l:shear=2shears} Every line shear around a point $o$ of a Desarguesian affine space $X$ is the composition of two line shears around the point $o$.
\end{lemma}

\begin{proof} Let $R$ be a line shear around a point $o\in X$. Then $\dom[R]\cap\rng[R]=\{o\}$. Since the affine space $X$ is $3$-long, there exists a line $\Lambda\subseteq X$ such that $\dom[R]\cap\Lambda=\Lambda\cap\rng[R]=\{o\}$. By Theorem~\ref{t:paraproj-exists}, there exists a line shear $R_1:\dom[R]\to\Lambda$ and by Lemma~\ref{l:shear=2shears}, there exists a line shear $R_2:\Lambda\to\rng[R]$ such that $R=R_2 R_1$.
\end{proof}

\begin{lemma}\label{l:3shears=2shears} The composition of three line shears around a point $o$ of a Desarguesian affine space $X$ is equal to the composition of two line shears around the point $o$.
\end{lemma}

\begin{proof} Let $P_1,P_2,P_3$ be line shears around a point $o$ in a Desarguesian affine space $X$, and $P=P_3 P_2 P_1$. If $\|\dom[P]\|\le 1$, then by Theorem~\ref{t:singleton-bijection}, $P$ is the composition of two line shears around $o$ in $X$.

So, we assume that $\|\dom[P]\|=2$. In this case $\dom[P]=\dom[P_1]$, $\rng[P_1]=\dom[P_2]$, $\dom[P_2]=\rng[P_3]$ and $\rng[P_3]=\rng[P]$.  

If $\dom[P_1]\ne\rng[P_2]$, then by Theorem~\ref{t:shear=ADA}, the composition $R=P_2 P_1$ is a line shear and hence $P=P_3  R$ is the composition of two line shears.

If $\dom[P_2]\ne\rng[P_3]$, then by Theorem~\ref{t:shear=ADA}, the composition $R'=P_3 P_2$ is a line shear and hence $P=R' P_1$ is the composition of two line shears. 

So, we assume that $\dom[P_1]=\rng[P_2]$ and $\dom[P_2]=\rng[P_3]$.  Then $$L\defeq\dom[P]= \dom[P_1]=\rng[P_2]=\dom[P_3]\quad\mbox{and}\quad \Lambda\defeq\rng[P]=\rng[P_3]=\dom[P_2]=\rng[P_1]$$ are two distinct lines with common point $o$. Since $X$ is $3$-long,  there exists a line $L'$ in $X$ such that $L\cap L'=L'\cap \Lambda$. Fix any points $x\in L\setminus\{o\}$ and $x'\in L'\setminus\{o\}$ and consider the point $y=P_1(x)\in\Lambda$. By Theorem~\ref{t:paraproj-exists}, there exist line shears $R$ and $R'$ such that $\dom[R]=\dom[P_1]$, $\rng[R]=L'=\dom[R']$, $\rng[R']=\rng[P_1]$ and $Rox=ox'$, $R'ox'=oy$. By Theorem~\ref{t:shear=ADA}, the composition $R' R:\dom[P_1]\to\rng[P_1]$ is a line shear such that $R'R(x)=P_1(x)$. By Theorem~\ref{t:paraproj-exists}, $R' R=P_1$ and hence $P=P_3 P_2 P_1=P_3 P_2 R' R$. By Theorem~\ref{t:shear=ADA}, the composition $P_2 R'$ is a line shear and so is the compositions $P_3(P_2 R')$. Therefore $P$ is equal to the composition $(P_3(P_2 R')) R$ of two line shears $P_3(P_3 R')$ and $R$.
\end{proof} 

Lemmas~\ref{l:shear=shear+shear} and \ref{l:3shears=2shears} imply the following corollary.

\begin{corollary}\label{c:n-shears=2shears} The composition of finitely many line shears around a point $o$ in a Desarguesian affine space $X$ is equal to the composition of two lines shears around the point $o$.
\end{corollary}



\section{Line affinities in Desarguesian affine spaces}

\begin{lemma}\label{l:affine+o=translation+shear} Let $A$ be a line affinity in a Desarguesian affine space $X$. If $\dom[A]\nparallel \rng[A]$, then for every point $o\in \dom[A]$, there exist a line translation $T$ and a line shear $R$ in $X$ such that $A=TR$ and $R(o)=o$.
\end{lemma}

\begin{proof} By Theorem~\ref{t:A=T'RRT}, $A=TR_n\cdots R_1$ for some  line translation $T$ and some  line shears $R_1,\dots,R_n$ around the point $o$. By Corollary~\ref{c:n-shears=2shears} and Theorem~\ref{t:shear=ADA}, the composition $R\defeq R_n\cdots R_1$ is a line shear around the point $o$. Then $A=TR_n\cdots R_1=TR$.
\end{proof}

\begin{lemma}\label{l:2transitive} Let $L,\Lambda$ be two non-parallel lines in a Desarguesian affine space $X$. Two line affinities $F,G$ with $\dom[F]=L=\dom[G]$ and $\rng[F]=\Lambda=\rng[G]$ are equal if and only if $(F(x),F(y))=(G(x),G(y))$ for some distinct points $x,y\in L$.
\end{lemma}

\begin{proof} The ``only if'' part is trivial. To prove the ``if'' part, assume that $(F(x),F(y))=(G(x),G(y))$ for some distinct points $x,y\in L$. By Lemma~\ref{l:affine+o=translation+shear}, there exist line translations $T,T'$ and line shears $R,R'$ such that $F=TR$, $G=T'R'$ and $R(x)=x=R'(x)$. Consider the point $$x'\defeq T'(x)=T'R'(x)=G(x)=F(x)=TR(x)=T(x)$$ and observe that $(T')^{-1}(x')=x=T^{-1}(x')$ and $$
\dom[T^{-1}]=\rng[T]=\rng[F]=\Lambda=\rng[G]=\rng[T']=\dom[(T')^{-1}].$$ By Theorem~\ref{t:unique-translation}, $T^{-1}=(T')^{-1}$ and hence $T=T'$.  Then $\rng[R]=\dom[T]=\dom[T']=\rng[R']$ and  
$$R(y)=T^{-1}F(y)=T^{-1}G(y)=R'(y).$$
Now the uniqueness part of Theorem~\ref{t:paraproj-exists} ensures that $R=R'$ and hence $F=R T=R' T'=G$.
\end{proof}

The following theorem is the main result of this section.

\begin{theorem}\label{t:affine-2transitive} Two line affinities in a Desarguesian affine space $X$ coincide if and only if $\dom[F]=\dom[G]$, $\rng[F]=\rng[G]$ and $(F(x),F(y))=(G(x),G(y))$ for some distinct points $x,y\in \dom[F]=\dom[G]$.
\end{theorem}

\begin{proof} The ``only if'' part is trivial. To prove the ``if'' part, assume that   $\dom[F]=\dom[G]$, $\rng[F]=\rng[G]$ and $(F(x),F(y))=(G(x),G(y))$ for some distinct points $x,y\in X$. If the lines   $\dom[F]=\dom[G]$ and $\rng[F]=\rng[G]$ are non-parallel, then the equality $F=G$ follows from Lemma~\ref{l:2transitive}. So, assume that $\dom[F]=\dom[G]\parallel \rng[F]=\rng[G]$. Choose any shear $R$ in $X$ with $\dom[R]=\rng[F]=\rng[G]$. Then $F'\defeq R F$ and $G'\defeq R G$ are line affinities such that the lines $\dom[F']=\dom[F]=\dom[G]=\dom[G']$ and $\rng[F']=\rng[R]=\rng[G']$ are not parallel and $(F'(x),F'(y))=(RF(x),RF(y)=(RG(x),RG(y)=(G'(x),G'(y))$. By Lemma~\ref{l:2transitive}, $F'=G'$ and hence $$F=1_{\rng[F]} F=R^{-1} R F=R^{-1} F'=R^{-1} G'=R^{-1} R G=1_{\rng[G]} G=G.$$
\end{proof}

\begin{corollary}\label{c:Desarg<=>unique-aff} An affine space $X$ is Desarguesian if and only if for any points $x,x',y,y'\in X$ with $x\ne y$ and $x'\ne y'$ there exists a unique line affinity $A:\Aline xy\to\Aline{x'}{y'}$ such that $A(x)=x'$ and $A(y)=y'$.
\end{corollary}

\begin{proof} The ``only if'' part follows from Theorems~\ref{t:aff-trans} and Theorem~\ref{t:affine-2transitive}, and the ``if'' part follows from Theorem~\ref{t:shear=ADA}.
\end{proof}

\begin{corollary}\label{c:affinity=2projections} Every line affinity $A$ in a Desarguesian affine space $X$ is the composition of two line shears or two line shifts.
\end{corollary}

\begin{proof} Fix any distinct points $x,y\in\dom[A]$ and consider the distinct points $x'\defeq A(x)$ and $y'\defeq A(y)$. By Theorem~\ref{t:aff-trans}, there exist two line projections $P,R$ such that $P R(x)=x'$ and $P R(y)=y'$, and moreover either $P,R$ are two line shears or $P,R$ are two line shifts. By Theorem~\ref{t:affine-2transitive}, $A=P R$.
\end{proof} 
 
\begin{proposition}  Every line affinity $A$ in a Desarguesian affine space $X$ is the composition of three line shears.
\end{proposition}

\begin{proof} Fix any distinct points $x,y\in\dom[A]$ and consider the distinct points $x'\defeq A(x)$ and $y'\defeq A(y)$. By Proposition~\ref{p:cov-aff}, there exists a point $c\in X\setminus(\dom[A]\cup\rng[A])$. By Theorem~\ref{t:paraproj-exists}, there exist line shears $R_1,R_2,R_3$ such that $R_1xy=xc$, $R_2cx=cx'$ and $R_3cx'=y'x'$. By Theorem~\ref{t:affine-2transitive}, $A=R_3R_2R_1$.
\end{proof}

\section{Homotheties versus line shears}

\begin{theorem}\label{t:RH=HR} Let $X$ be an affine space and $H:X\to X$ be a homothety of $X$ with center $o$. Then  $HR=RH$ for every line shear $R$ in $X$ around the point $o$.
\end{theorem}

\begin{proof} Taking into account that $H$ is a homothety with center $o\in \dom[R]$, we conclude that $H[\dom[R]]=\dom[R]$ and hence $\dom[HR]=\dom[R]=\dom[RH]$. Given any point $x\in \dom[HR]=\dom[RH]$, we need to show that $HR(x)=RH(x)$. If $x=o$, then $HR(x)=HR(o)=H(o)=o=R(o)=RH(o)=RH(x)$ and we are done. So, assume that $x\ne o$.

Consider the point $y\defeq R(x)$ and the points $x'\defeq H(x)\in\Aline ox=\dom[R]$ and $y'\defeq H(y)\in \Aline oy=\rng[R]$. Since $H$ is a translation with center $o$, $\Aline {x'}{y'}\parallel \Aline xy$ and $y'\in\Aline o{y}=\rng[R]$, witnessing that $y'=R(x')$ and hence $HR(x)=H(y)=y'=R(x')=RH(x)$.   
\end{proof}

The following theorem is a ``shear'' counterpart of the Gleason's Theorem~\ref{t:Gleason1}.

\begin{theorem} Let $X$ be an affine space, $L$ be a line in $X$ and $o\in L$ be a point. For a bijection $A:L\to L$ with $A(o)=o$, the following conditions are equivalent:
\begin{enumerate}
\item $A=H{\restriction}_L$ for some homothety $H:X\to X$;
\item for every line $\Lambda\subseteq X$ with $L\cap \Lambda=\{o\}$ and any line shears $P:L\to\Lambda$ and $R:\Lambda\to L$ we have $A(RP)=(RP)A$;
\item for every line $\Lambda\subseteq X$ with $L\cap \Lambda=\{o\}$ and any line shears $P,R:\Lambda\to L$, we have $P^{-1}AP=R^{-1}AR$.
\end{enumerate}
\end{theorem}

\begin{proof} $(1)\Ra(2)$. Assume that $A=H{\restriction}_L$ for some homothety $H:X\to X$. Then $H(o)=A(o)=o$, so $o$ is the center of the homothety. Given any line $\Lambda\subset X$ with $L\cap\Lambda=\{o\}$ and line shears $P:L\to\Lambda$ and $R:\Lambda\to L$, we need to check that $A(RP)=(RP)A$. Let $1_L$ be the identity map of the line $L$. Then $$A(RP)=H1_LRP=HRP=RHP=RPH=R(P1_L)H=(RP)H1_L=(RP)A,$$
by Theorem~\ref{t:RH=HR}.
\smallskip

$(2)\Ra(3)$ Assume that the condition (2) is satisfied. To check the condition (3), take any line $\Lambda\subseteq X$ with $L\cap \Lambda=\{o\}$ and any line shears $P,R:\Lambda\to L$. The condition (2) ensures that $A(PR^{-1})=(PR^{-1})A$. Multiplying this equality by $P^{-1}$ from the left and $R$ from the right, we obtain the desired equality
$$P^{-1}AP=P^{-1}AP1_\Lambda=P^{-1}APR^{-1}R=P^{-1}PR^{-1}AR=1_\Lambda R^{-1}AR=R^{-1}AR.$$
\smallskip

$(3)\Ra(1)$ Assume that the line bijection $A:L\to L$ satisfies the condition (3). For every point $x\in X\setminus L$, fix any line shear $R_x:\Aline ox\to L$ with $R_x(o)=o$. For every $x\in L$, let $R_x:L\to L$ be the identity map of the line $L$. Consider the function $H:X\to X$ assigning to every point $x\in X$ the point  $H(x)\defeq R_x^{-1}AR_x(x)$. It is clear that $H{\restriction}_L=A$. The condition (3) ensures that  for any $x\in X\setminus L$, the value $H(x)=R_x^{-1}AR_x(x)$ does not depend on the choice of the line shear $R_x:\Aline ox\to L$.

\begin{claim}\label{cl:H-=RAR} The function $H:X\to X$ is bijective and $H^{-1}(y)=R_y^{-1}A^{-1}R_y(y)$ for every $y\in X$.
\end{claim}

\begin{proof} To see that $H$ is injective, take any distinct point $x,y\in X$. We have to show that the points $x'\defeq H(x)$ and $y'\defeq H(y)$ are distinct. 

If $x=o$, then $x'=H(x)=H(o)=o=R_y^{-1}AR_y(o)\ne R_y^{-1}AR_y(y)=y'$. 

If $y=o$, then $y'=H(y)=H(o)=o=R_x^{-1}AR_x(o)\ne R_x^{-1}AR_x(x)=x'$.

So, assume that $o\notin\{x,y\}$, which implies $x'\ne o\ne y'$. 

If $\Aline ox\ne\Aline oy$, then $\{x'\}\cap\{y'\}\subseteq \Aline ox\cap\Aline oy=\{o\}$ and hence $x'\ne y'$. 

So, assume that $\Aline ox=\Aline oy$. In this case we can take the line shears $R_x:\Aline ox\to L$ and $R_y:\Aline oy\to L$ to be equal. Since the line shear $R_x=R_y$ is injective, the inequality $x\ne y$ and the condition (3) imply
$x'=R_x^{-1}AR_x(x)\ne R_x^{-1}AR_x(y)=R_y^{-1}AR_y(y)=y'$.

Next, we show that the map $H:X\to X$ is surjective. Given any point $y\in X\setminus \{x\}$, consider the point $x\defeq R_y^{-1}A^{-1}R_y(y)$, and observe that $\Aline ox=\Aline oy$ and hence $\dom[R_x]=\dom[R_y]$. The condition (3) ensures that 
$$
\begin{aligned}
H(x)&=R_x^{-1}AR_x(x)=R_y^{-1}AR_y(x)
=R_y^{-1}AR_yR_y^{-1}A^{-1}R_y(y)\\
&=R_y^{-1}A1_LA^{-1}R_y(y)=R_y1_LR_y^{-1}(y)=y,
\end{aligned}
$$
witnessing that the injective function $H:X\to X$ is bijective and hence $H^{-1}(y)=R_y^{-1}AR_y(y)$ for every $y\in H$.
\end{proof}

\begin{claim}\label{cl:Hxy||xy} For every distinct points $x,y\in X$ and the points $x'\defeq H(x)$ and $y'\defeq H(y)$, the lines $\Aline{x'}{y'}$ and $\Aline xy$ are parallel.
\end{claim}

\begin{proof} If $x=o$, then $x'=H(x)=H(o)=o$, $y'=H(y)\in \Aline oy$ and $\Aline {x'}{y'}=\Aline o{y'}=\Aline oy=\Aline xy$.

If $y=o$, then $y'=H(y)=H(o)=o$, $x'=H(x)\in \Aline ox$ and $\Aline {x'}{y'}=\Aline {x'}{o}=\Aline xo=\Aline xy$.
So, assume that $o\notin\{x,y\}$. If $\Aline ox=\Aline oy$, then $\Aline {x'}{y'}=\Aline ox\cup\Aline oy=\Aline xy$. In all these cases, $\Aline{x'}{y'}=\Aline xy$ and hence $\Aline{x'}{y'}\parallel \Aline xy$.

So, assume that $\Aline ox\ne \Aline oy$ and $o\notin\{x,y\}$. Five cases are possible. 
\smallskip

1. If $\Aline ox=L\ne \Aline oy$, then we lose no generality assuming that $R_y(y)=x$. In this case $x'=H(x)=A(x)$ and $y'=H(y)=R_y^{-1}AR_y(y)=R_y^{-1}(x')$. Since $R_y$ is a line shear with $R_y(y)=x\ne o\ne x'=R_y(y')$, the line $\Aline xy$ is parallel to the line $\Aline {x'}{y'}$. \smallskip

2. If $\Aline ox\ne L=\Aline oy$, then by analogy with the preceding case we can show that the line  $\Aline{x'}{y'}$ is parallel to the line $\Aline xy$.
\smallskip

3. If $\Aline ox\ne L\ne\Aline oy$ and $\Aline xy\cap L=\{z\}$ for some point $z\ne o$, then we lose no generality assuming that $R_x(x)=z=R_y(y)$. Consider the point $z'\defeq H(z)=A(z)$ and observe that $x'=H(x)=R_x^{-1}AR_x(x)=R_x^{-1}(z')$ and $y'=H(y)=R_y^{-1}AR_y(y)=R_y^{-1}(z')$. Since $R_x$ and $R_y$ are line shears, $\Aline {x'}{z'}\parallel\Aline xz$ and $\Aline {y'}{z'}\parallel \Aline yz$. Since $\Aline xz=\Aline yz$, the parallel lines $\Aline {x'}{z'}$ and $\Aline {y'}{z'}$ coincide and hence $\Aline {x'}{y'}=\Aline{x'}{z'}\parallel \Aline xz=\Aline xy$.
\smallskip

4. If $\Aline ox\ne L\ne\Aline oy$ and the lines $L$ and $\Aline xy$ are disjoint and coplanar, then the line $\Aline {x'}{y'}$ in contained in the plane $\overline{L\cup\{x,y\}}$. Assuming that $\Aline {x'}{y'}\nparallel L$, we can find a point $z'\in L\cap\Aline{x'}{y'}$. We lose no generality assuming that $R_{x'}(x')=z'=R_{y'}(y')$.   
Consider the point $z\defeq H^{-1}(z')=A^{-1}(z')$. Applying Claim~\ref{cl:H-=RAR}, we conclude that $x=H^{-1}(x')=R_{x'}^{-1}A^{-1}R_{x'}(x')=R_x^{-1}A^{-1}(z')=R^{-1}_{x'}(z)$ and $y=H^{-1}(y')=R_{y'}^{-1}A^{-1}R_{y'}(y')=R_{y'}^{-1}A^{-1}(z')=R_{y'}^{-1}(z)$. Since $R_{x'}$ and $R_{y'}$ are line shears with $R_{x'}xx'=zz'=R_{y'}yy'$, we have the parallelity relations $\Aline {x}{z}\parallel\Aline {x'}{z'}=\Aline{x'}{y'}=\Aline {y'}{z'}\parallel \Aline yz$ implying $\Aline xz=\Aline yz=\Aline xy$, which contradict $L\cap \Aline xy\ne\varnothing$. This contradiction shows that $\Aline {x'}{y'}\cap L=\varnothing$. By Corollary~\ref{c:parallel-lines<=>}, $\Aline{x'}{y'}\parallel L\parallel \Aline xy$.
\smallskip

5. Finally assume that the lines $L$ and $\Aline xy$ are skew. In this case the affine space $X$ has rank $\|X\|\ge 4$ and is Desarguesian, by Corollary~\ref{c:affine-Desarguesian}. We lose no generality assuming that $R_x(x)=z=R_y(x)$ for some point $z\in L$. Consider the point $z'\defeq H(z)=A(z)$ and observe that $x'=H(x)=R_x^{-1}AR_x(x)=R_x^{-1}(z')$ and $y'=H(y)=R_y^{-1}AR_y(y)=R_y^{-1}(z')$. Taking into account that $R_x,R_y$ are line shears around the point $o\in\Aline x{x'}\cap\Aline y{y'}\cap\Aline z{z'}$ and $\Aline xz\parallel \Aline{x'}{z'}$ and $\Aline yz\parallel \Aline{y'}{z'}$, we can apply Desargues' Theorem~\ref{t:Desargues-affine} and conclude that $\Aline{x'}{y'}\parallel \Aline {x}{y}$.
\end{proof}

Now we shall prove that the bijective function $H:X\to X$ is a homothety of $X$. Given any line $\Lambda\subseteq X$, we should prove that its image $T[\Lambda]$ is a line, parallel to the line $\Lambda$. Fix any distinct points $x,y\in\Lambda$ and consider their images $x'\defeq H(x)$ and $y'\defeq H(y)$, which are distinct because $H$ is a bijection of $X$. By Claim~\ref{cl:Hxy||xy}, the line $\Aline{x'}{y'}$ is parallel to the line $\Aline xy$. So, it suffices to show that $H[\Lambda]=\Aline {x'}{y'}$. Given any point $z\in \Aline xy\setminus\{x\}$, consider its image $z'\defeq H(z)$. Claim~\ref{cl:Hxy||xy} ensures that $\Aline {x'}{z'}\parallel \Aline xz=\Aline xy\parallel \Aline {x'}{y'}$ and hence $z'\in\Aline{x'}{z'}=\Aline {x'}{y'}$. This proves that $H[\Lambda]\subseteq\Aline {x'}{y'}$. Applying the same argument to the bijection $H^{-1}$, we can prove that $H^{-1}[\Aline{x'}{y'}]\subseteq \Aline xy=\Lambda$ and hence $H[\Lambda]=\Aline {x'}{y'}$ is a line, parallel to the line $\Lambda=\Aline xy$. Therefore, $H$ is a homothety of the affine space $X$. 
\end{proof}

\chapter{Portions and scalars in affine spaces}

\section{Affine equivalence of line triples in affine spaces}

In this section we discuss the notion of affine equivalence of line triples in an affine space $X$. This equivalence is induced by the action of the inverse monoid $\I_X^{\prop}$ on the set\index[note]{$\dddot X$} 
$$\LT\defeq\{xyz\in X^3:y\in\Aline xz\;\wedge\;x\ne z\}$$ whose elements will be called \index{line triple}\defterm{line triples} in $X$. The inverse monoid $\I_X^{\prop}$ is the smallest submonoid of the symmetric inverse monoid $\I_X$ that contains all line projections (i.e., parallel projections between lines in $X$). Line bijections that belong to the monoid $\I_X^{\prop}$ are called line affinities. They are finite compositions of line projections.
  
Every triple $x_0x_1x_2\in X^3$ is (by definition) the function $\{(0,x_0),(1,x_1),(2,x_2)\}:3\to X$ on the number $3\defeq\{0,1,2\}$. So, for any function $F$, we can consider the composition 
$$Fx_0x_1x_2\defeq\{(i,F(x_i)):i\in 3,\;x_i\in\dom[F]\}$$ of the functions $x_0x_1x_3$ and $F$. The function $Fx_0x_1x_2$ has domain $\dom[Fx_0x_1x_2]=\{i\in 3:x_i\in\dom[F]\}$. This domain equals $3$ if and only if $\{x_0,x_1,x_2\}\subseteq\dom[F]$. 

\begin{definition}\label{d:affine-equivalence} Let $X$ be an affine space. Given two line triples $abc,xyz$ in $X$, we write\index[note]{$abc\prop xyz$} $abc\prop xyz$ and say that $abc$ and $xyz$ are \index{line triples!affinely equivalent}\index{affinely equivalent line triples}\defterm{affinely equivalent} if there exists a line affinity $A$ in $X$ such that $Aabc=xyz$. 
\end{definition}

\begin{proposition}\label{p:aff-eq} Let $X$ be an affine space. For any line triples $abc,uvw,xyz$ in $X$, the following properties hold:
\begin{enumerate}
\item $abc\prop abc$;
\item $abc\prop uvw\;\Leftrightarrow\;uvw\prop abc$;
\item $(abc\prop uvw\;\wedge\;uvw\prop xyz)\;\Ra\;abc\prop xyz$;
\item $abc\prop uvw\;\Leftrightarrow\; cba\prop wvu$;
\item $aac\prop uuw\;\wedge\;acc\prop uww$;
\item If $abc\prop uvw$, then $(a=b\;\Leftrightarrow\;u=v)\;\wedge  (b=c\;\Leftrightarrow\;v=w)$.
\end{enumerate}
\end{proposition}

\begin{exercise} Prove Proposition~\ref{p:aff-eq}.
\smallskip

\noindent{\em Hint:} In the proof of Proposition~\ref{p:aff-eq} apply Theorem~\ref{t:aff-trans}.
\end{exercise}

In the following theorem we characterize the affine equivalence without mentioning line affinities.

\begin{theorem}\label{t:aff-eq<=>} Two line triples $abc,uvw$ in a (Desarguesian) affine space $X$ are affinely equivalent if and only if  there exist a number $n\in\w$ (with $n\le 2$), line triples $x_0y_0z_0,x_1y_1z_1,\dots$, $x_ny_nz_n$, and lines $\Lambda_1,\dots,\Lambda_n$ such that $x_0y_0z_0=abc$, $x_ny_nz_n=uvw$, and  for every positive $i\le n$, $$\Aline {x_{i-1}}{x_i}\subparallel\Lambda_i,\quad \Aline {y_{i-1}}{y_i}\subparallel \Lambda_i,\quad \Aline {z_{i-1}}{z_i}\subparallel \Lambda_i,  \quad \mbox{and}\quad  \Aline {x_{i-1}}{z_{i-1}}\ne \Aline{x_i}{z_i}.$$
\end{theorem}

\begin{proof} To prove the ``only if'' part, assume that the  line triples $abc$ and $uvw$ are affinely equivalent. If $abc=xyz$, then for $n=0$ and the triple $x_0y_0z_0\defeq abc$, we obtain $abc=x_0y_0z_0=x_ny_nz_n=uvw$ (and since no positive $i\le n=0$ exists the other conditions hold vacuously).

Now assume that $abc\ne uvw$. Since $abc\prop uvw$, there exists a line affinity $A$ in $X$ such that $Aabc=uvw$. Write $A$ is the composition $P_n\cdots P_1$ of line projections $P_1,\dots,P_n$ in $X$. We can assume that the number $n$ is the smallest possible. Since $abc\ne uvw$, the line affinity $A$ is not an identity map. The minimality of $n$ ensures that all line projections $P_i$ are not identity maps and hence $\dom[P_i]\ne\rng[P_i]$, by Proposition~\ref{p:para-projection}(2).   If the affine space $X$ is Desarguesian, then $n\le 2$, by  Corollary~\ref{c:affinity=2projections}.

By the definition of a line projection, for every $i\in\{1,\dots,n\}$, there exists a line $\Lambda_i$ such that $P_i=\{(x,y)\in\dom[P_i]\times\rng[P_i]:\Aline xy\subparallel \Lambda_i\}$. Let $x_0y_0z_0\defeq abc$ and for every $i\in\{1,\dots,n\}$, let $x_iy_iz_i\defeq P_ix_{i-1}y_{i-1}z_{i-1}$. 
Then $$\Aline {x_{i-1}}{x_i}\subparallel\Lambda_i,\quad \Aline {y_{i-1}}{y_i}\subparallel \Lambda_i,\quad\mbox{and}\quad \Aline {z_{i-1}}{z_i}\subparallel \Lambda_i,$$by the definition of the line projection $P_i$ and the choice of the line triples $x_{i-1}y_{i-1}z_{i-1}$ and $x_iy_iz_i=P_ix_{i-1}y_{i-1}z_{i-1}$. 
It follows from $a\ne c$, $u\ne w$ and $Aabc=uvw$ that $\Aline {x_{i-1}}{z_{i-1}}=\dom[P_i]\ne\rng[P_i]=\Aline {x_i}{z_i}$. This completes the proof of the ``only if'' part.
\smallskip

To prove the ``if'' part, assume that for the pairs $ab$ and $uv$ there exist a number $n\in\w$,  line triples $x_0y_0z_0,\dots,x_ny_nz_n$ in $X$ and lines $\Lambda_1,\dots,\Lambda_n$ such that $x_0y_0z_0=abc$, $x_ny_nz_n=uvw$ and 
 $$\Aline {x_{i-1}}{x_i}\subparallel\Lambda_i,\quad \Aline {y_{i-1}}{y_i}\subparallel \Lambda_i,\quad \Aline {z_{i-1}}{z_i}\subparallel \Lambda_i,  \quad \mbox{and}\quad  \Aline {x_{i-1}}{z_{i-1}}\ne \Aline{x_i}{z_i}$$
for every positive number $i\le n$.

\begin{claim}\label{cl:xi=zi} For every $i\in\{0,\dots,n\}$,  $\Aline {x_{i-1}}{z_{i-1}}\nparallel \Lambda_i\nparallel\Aline{x_i}{z_i}$.
\end{claim}

\begin{proof} To derive a contradiction, assume that $\Lambda_i\parallel \Aline {x_i}{z_i}$. Since $\Aline {x_{i-1}}{z_{i-1}}\ne\Aline {x_i}{z_i}$, either $x_{i-1}\ne x_i$ or $z_{i-1}\ne z_i$. If $x_{i-1}\ne x_i$, then  by Corollary~\ref{c:subparallel}, $\Aline {x_{i-1}}{x_i}\subparallel \Lambda_i$ implies $ \Aline {x_{i-1}}{x_i}\parallel \Lambda_i\parallel \Aline{x_i}{z_i}$ and hence $x_{i-1}\in \Aline{x_{i-1}}{x_i}= \Aline{x_i}{z_i}$. It follows from $\Aline{z_{i-1}}{z_i}\subparallel \Lambda_i\parallel \Aline {x_i}{z_i}$ that $z_{i-1}\in\Aline {x_i}{z_i}$ and hence $\Aline {x_{i-1}}{z_{i-1}}=\Aline {x_i}{z_i}$, which contradicts the choice of the line triples $x_{i-1}y_{i-1}z_{i-1}$ and $x_iy_iz_i$. 

This contradiction shows that $x_{i-1}=x_i$ and hence $z_{i-1}\ne z_i$. By Corollary~\ref{c:subparallel}, $\Aline {z_{i-1}}{z_i}\subparallel \Lambda_i$ implies $ \Aline {z_{i-1}}{z_i}\parallel \Lambda_i\parallel \Aline{x_i}{z_i}$ and hence $z_{i-1}\in \Aline{z_{i-1}}{z_i}= \Aline{x_i}{z_i}$. It follows from $\Aline{x_{i-1}}{x_i}\subparallel \Lambda_i\parallel \Aline {x_i}{z_i}$ that $x_{i-1}\in\Aline {x_i}{z_i}$ and hence $\Aline {x_{i-1}}{z_{i-1}}=\Aline {x_i}{z_i}$, which contradicts the choice of the line triples $x_{i-1}y_{i-1}z_{i-1}$ and $x_iy_iz_i$. 
This contradiction shows that $\Lambda_i\nparallel \Aline {x_i}{z_i}$.

By analogy one can prove $\Lambda_i\nparallel \Aline{x_{i-1}}{z_{i-1}}$. 
\end{proof}

 By Proposition~\ref{p:flat-relation}, the relation $$P_i\defeq\{(x,y)\in\Aline {x_{i-1}}{z_{i-1}}\times\Aline{x_i}{z_i}:\Aline xy\subparallel \Lambda_i\}$$is a line projection. Then the composition $A=P_n\cdots P_1$ is a line affinity witnessing that $Aabc=uvw$ and hence $abc\prop uvw$.
\end{proof}

\section{Desarguesian triples in affine spaces}

\begin{proposition}\label{p:triple-exists} For every line triple $abc$ in an affine space $X$ and every distinct points $x,z$ in $X$ there exists a point $y\in\Aline xz$ such that $abc\prop xyz$.
\end{proposition}
 
\begin{proof} By Theorem~\ref{t:aff-trans}, there exists a line affinity $A:\Aline ac\to\Aline xz$ such that $Aac=xz$. Then the point $y=A(b)$ has the required property: $Aabc=xyz$ and hence $abc\prop xyz$.
\end{proof}

\begin{theorem}\label{t:triple-unique} An affine space $X$ is Desarguesian if and only if for every line triple $abc$ in $X$ and every distinct points $x,z\in X$ there exists a unique point $y\in \Aline xz$ such that $xyz\prop abc$.
\end{theorem}

\begin{proof} Assume that the affine space $X$ is Desarguesian and fix a line triple $abc$ in $X$ and two distinct points $x,z\in X$. By Proposition~\ref{p:triple-exists}, there exists a point $y\in\Aline xz$ such that $xyz\prop abc$. To show that the point $y$ is unique, assume that $y'\in\Aline xz$ is another point such that $xy'z\prop abc$. Since $xyz\prop abc\prop xy'z$, there exists a line affinity $A:\Aline xz\to\Aline xz$ such that $Axyz=xy'z$. By Theorem~\ref{t:affine-2transitive}, $A$ is the identity map of the line $\Aline xz$ and hence $y'=A(y)=y$.
\smallskip

If the affine space $X$ is not Desarguesian, then by Corollary~\ref{c:Desarg<=>unique-aff},  there exist points $a\ne c$ and $x\ne z$ in $X$ and two distinct line affinities $A,A'$ in $X$ such that $Aac=xz=A'ac$. It follows from $\{a,c\}\in\dom[A]\cap\dom[A']$ and $a\ne c$ that $\dom[A]=\Aline ac=\dom[A']$. Since $A\ne A'$, there exists a point $b\in\dom[A]=\dom[A']=\Aline ac$ such that $A(b)\ne A'(b)$. Then $y\defeq A(b)$ and $y'\defeq A'(b)$ are two distinct points in the line $\Aline xz=\rng[A]=\rng[A']$ such that $xy'z\prop abc\prop xyz$.
\end{proof}

Proposition~\ref{t:triple-unique} motivates the following definition.

\begin{definition}\label{d:Desarg-triple} A line triple $abc$ in an affine space $X$ is \index{line triple!Desarguesian}\index{Desarguesian line triple}\defterm{Desarguesian} if for every distinct points $x,z\in X$, there exists a unique point $y\in\Aline xz$ such that $xyz\prop abc$. 
\end{definition}

Proposition~\ref{p:aff-eq}(5) implies the following important example of Desarguesian line triples.

\begin{example} For every distinct points $a,c$ in an affine space $X$, the line triples $aac$ and $acc$ are Desarguesian.
\end{example}

In terms of Desarguesian triples, Theorem~\ref{t:triple-unique} can be rewritten as follows.
 
\begin{theorem}\label{t:Desargues<=>3Desargues} An affine space $X$ is Desarguesian if and only if every line triple in $X$ is Desarguesian.
\end{theorem}

\begin{proposition}\label{p:Desarg-triple<=>} A line triple $abc$ in an affine space $X$  is Desarguesian if and only if $A(b)=b$ for any line affinity $A$ with $Aac=ac$.
\end{proposition}

\begin{proof} The ``only if'' part follows immediately from the definition of a Desarguesian line triple.

To prove the ``if'' part, assume that the line triple $abc$ is not Desarguesian. Then there exists two two distinct points $x,z\in X$ and two distinct points $y,y'\in\Aline xz$ such that $xyz\prop abc\prop xy'z$. Find two line affinities $F,F'$ such that $Fabc=xyz$ and $F'xy'z=abc$. Then the line affinity $A=F'F$ has the required property: $Aac=F'Fab=F'xz=ac$ and $A(b)=F'(F(b))=F'(y)\ne F'(y')=b$.
\end{proof}

\begin{proposition}\label{p:Desarg-triples} Let $a,b,c,d$ be four points in an affine space $X$.
\begin{enumerate}
\item If the triple $abc$ is Desarguesian, then the triple $cba$ is Desarguesian;
\item If the triple $abc$ is Desarguesian and $a\ne b$, then the triple $acb$ is Desarguesian.
\item If the triples $abc$ and $acd$ are Desarguesian, then the triple $abd$  is Desarguesian.
\end{enumerate}
\end{proposition}

\begin{proof} 1. Assume that $abc$ is a Desarguesian line triple. By Proposition~\ref{p:Desarg-triple<=>}, to prove that the triple $cba$ is Desarguesian, it suffices to check that $A(b)=b$ for every line affinity $A$ with $Aca=ca$. It follows from $Aca=ca$ that $Aac=ac$. Since the triple $abc$ is Desarguesian, $A(b)=b$, witnessing that the triple $cba$ is Desarguesian.
\smallskip

2. Assume that the triple $abc$ is Desarguesian and $a\ne b$. By Proposition~\ref{p:Desarg-triple<=>}, to prove that the triple $acb$ is Desarguesian, it suffices to check that $A(c)=c$ for every line affinity $A$ with $Aab=ab$. 
So, fix a line affinity $A$ with $Aab=ab$ and put $c'\defeq A(c)$. Choose any point $x\in X\setminus\Aline ac$. By Theorem~\ref{t:paraproj-exists}, there exist line shears $P,R$ such that $Pac'=ax$ and $Rax=ac$. Then $RPA$ is a line affinity such that $RPAac=ac$. Since the triple $abc$ is Desarguesian, $RPA(b)=b=A(b)$ and hence $RP(b)=b$. Then for the point $b'\defeq P(b)=R^{-1}(b)$ we have $\Aline {c'}x\parallel \Aline b{b'}\parallel\Aline {c}{x}$, which implies $\Aline {c'}{x}=\Aline c{x}$ and $\{c'\}=\Aline o{c'}\cap \Aline {c'}{x}=\Aline oc\cap \Aline cx=\{c\}$. Then  $A(c)=c'=c$.
\smallskip

3. Assume that $abc$ and $acd$ are Desarguesian line triples. By Proposition~\ref{p:Desarg-triple<=>}, to prove that the triple $abd$ is Desarguesian, it suffices to check that  $A(b)=b$ for every line affinity $A:\Aline ad\to\Aline ad$ with $Aad=ad$.
Since $abc$ and $acd$ are line triples, $\Aline ac=\Aline ad=\dom[A]=\rng[A]$. Since the line triple $acd$ is Desarguesian, the equality $Aad=ad$ implies the equality $Aac=ac$. Since the line triple $abc$ is Desarguesian, the latter equality implies $A(b)=b$. 
\end{proof}

A quarduple $abcd$ in a liner $X$ is called a \index{line quadruple}\defterm{line quadruple} if the points $a,b,c,d$ are colinear and $|\{a,c,d\}|=3$, i.e., the points $a,c,d$ are pairwise distinct.

The following lemma will help us to define the operation of multiplication of portions and scalars in Lemma~\ref{t:scalar-by-portion}.

\begin{lemma}\label{l:invariance-scalar-multiplication} Let $abcd$ and $xyzw$ be two line quadruples in an affine space $X$ such that  $abc\prop xyz$ and $acd\prop xzw$. If at least one of the line triples $abc$ or $acd$ is Desarguesian, then $abd\prop xyw$.
\end{lemma}

\begin{proof} Assume that $abc\prop xyz$ and $acd\prop xzw$. Then there exist line affinities $A,A'$ such that $Aabc=xyz$ and $A'acd=xzw$. Since $abcd$ and $xyzw$ are line quadruples, $\dom[A]=\Aline ac=\Aline ad=\dom[A']$, $\rng[A]=\Aline xz=\Aline xw=\rng[A']$, and $|\{a,c,d\}|=3=|\{x,z,w\}|$.

If the line triple $abc$ is  Desarguesian, then $A'ac=xz=Aac$ implies $A'(b)=A(b)=y$.
If the line triple $acd$ is Desarguesian, then the line triple $adc$ is Desarguesian, by Proposition~\ref{p:Desarg-triples}(2). Then the equality $Aac=xz=A'ac$ implies $A(d)=A'(d)=w$. Then $Aabd=xyw$, witnessing that $abd\prop xyw$.
\end{proof}

\section{Portions and scalars in affine spaces}

We recall that $\dddot X$ stands for the set $\{xyz\in X^3:y\in\Aline xz\;\wedge\;x\ne z\}$ of all line triples in a liner $X$.

\begin{definition} For a line triple $zvu$ in an affine space $X$, its equivalence class\index[note]{$\overvector{zvu}$} $$\overvector{zvu}\defeq\{abc\in\LT:abc\prop zvu\}$$is called the \index{portion}\defterm{portion} determined by the line triple $zvu$. The portion $\overvector{zvu}$ will be also denoted by $zv/zu$ and called the \index{ratio}\defterm{ratio} of the pairs $zv$ and $zu$. 

 The portion $\overvector{zvu}$ is called a \index{scalar}\defterm{scalar} if the line triple $zvu$ is Desarguesian. In this case every line triple $abc\in \overvector{zvu}$ is Desarguesian.
\end{definition}

In a line triple $zvu$, the points $z,v,u$ play non-symmetric roles. The notations $z,v,u$ are abbreviations of {\em zero}, {\em  variable}, {\em unit}, respectively.

\begin{picture}(100,35)(-100,-15)

\put(0,0){\line(1,0){150}}

\put(30,0){\circle*{3}}
\put(28,3){$0$}
\put(27,-8){$z$}
\put(60,0){\color{red}\circle*{3}}
\put(57,-8){$v$}
\put(120,0){\circle*{3}}
\put(117,-8){$u$}
\put(118,3){$1$}
\end{picture}

For an affine space $X$, we denote by\index[note]{${\dddot X}_{\Join}$}
$${\color{magenta}\dddot X_{\!\Join}}\defeq\{\overvector{zvu}:zvu\in\dddot X\}$$the set of portions in the affine space $X$, and by \index[note]{$\IR_X$}\defterm{$\IR_X$} the set of scalars in $X$, i.e., the set of portions $\overvector{zvu}\in\dddot X_{\!\Join}$ (or ratios $zv/zu$) of Desarguesian line triples  $zvu$ in $X$. Therefore, $\IR_X\subseteq\dddot X_{\!\Join}$. By Theorem~\ref{t:Desargues<=>3Desargues}, an affine space $X$ is Desarguesian if and only if $\IR_X=\dddot X_{\!\Join}$.  

For every affine space $X$, the set $\IR_X$ contains two important scalars\index[note]{$0$}\index[note]{$1$} 
$$\mbox{\defterm{$0$}}\defeq\{xyz\in X^3:x=y\ne z\}\quad\mbox{and}\quad \mbox{\defterm{$1$}}\defeq\{xyz\in X^3:x\ne y=z\}.$$

Proposition~\ref{p:aff-eq}(8) implies that $0$ and $1$ are indeed scalars, which can be characterized as follows.

\begin{proposition}\label{p:sczalars01} For a line triple $zvu$ in an affine space $X$ 
\begin{enumerate}
\item $\overvector{zvu}=0$ if and only if $v=z$;
\item $\overvector{zvu}=1$ if and only if $v=u$.
\end{enumerate}
\end{proposition}

In the sets $\dddot X_{\!\Join}$ and $\IR_X$, consider two important subsets
$${\color{magenta}\dddot X_{\!\Join}^*}\defeq\dddot X_{\!\Join}\setminus\{0\}\quad\mbox{and}\quad{\color{magenta}\IR_X^*}\defeq\IR_X\setminus\{0\}$$of \index[note]{${\dddot X}^*_{\Join}$}
\index[note]{$\IR_X^*$}\defterm{nonzero portions} and \defterm{nonzero scalars}, respectively.

\begin{remarks} In Corollary~\ref{c:portions=3} we shall prove that for every non-Thalesian finite affine space $X$, the set $\dddot X_{\Join}$ contains exactly three elements: $0$, $1$, and the portion $\overvector{zvu}$ for any line triple $zvu$ with $z\ne v\ne u$.
\end{remarks} 

\section{Unary operations on portions}\label{s:unary-RX}

In this section we consider two natural unary operations on portions, which are induced by permutations of line triples. The permutatin $zvu\mapsto zuv$ induces the operation of inversion of nonzero portions. 

\begin{proposition}\label{p:portion-inversion} For every affine space $X$, there exists a unique unary operation 
$$(\cdot)^{-1}:\dddot X_{\Join}\to\dddot X_{\!\Join},\quad (\cdot)^{-1}:\alpha\mapsto \alpha^{-1},$$ such that $(\overvector{zvu})^{-1}=\overvector{zuv}$ for every line triple $zvu\in\dddot X_{\!\Join}$ with $v\ne z$. 
\end{proposition}

\begin{proof} Define the unitary operation $(\cdot)^{-1}:\dddot X_{\!\Join}^*\to\dddot X_{\!\Join}^*$ assigning the every nonzero portion $\alpha\in \dddot X_{\!\Join}^*$ the portion $\overvector{zuv}$ where $zvu$ is any line triple in $\alpha$. Let us show that the operation $(\cdot)^{-1}$ is well-defined. Since $zvu\in\alpha\ne 0$, the point $v$ is not equal to the point $z$ and hence $zuv$ is a line triple determining the portion $\overvector{zuv}=\alpha^{-1}$. Since $z\ne u$, the portion $\alpha^{-1}=\overvector{zuv}$ is nonzero.

Let us show that the portion $\alpha^{-1}$ does not depend on the choice of the line triple $zvu\in\alpha$. Any other line triple $oae\in\alpha$ is affinely equivalent to the line triple $zvu\in\alpha$ and hence $Aoae=zvu$ for some line affinity $A$ in $X$, witnessing that $Azvu=oae$ and hence $\overvector{oea}=\overvector{zuv}$.
\end{proof}

\begin{definition} The unary operation $(\cdot)^{-1}$ on $\dddot X^*_{\Join}$ (resp. $\IR_X$) introduced in Proposition~\ref{p:portion-inversion} will be called the \index{inversion}\index{portion!inversion of}\index{scalar!inversion of}\defterm{operation of inversion} of portions (resp.  scalars) in the affine space $X$. 
\end{definition} 

Propositions~\ref{p:portion-inversion} and \ref{p:Desarg-triples}(2) imply the following corollary.

\begin{corollary}\label{c:inversionR=>involution} Let $X$ be an affine space. The operation of inversion on the set $\dddot X_{\!\Join}^*$ has the following properties:
\begin{enumerate}
\item $(p^{-1})^{-1}=p$ for every $p\in \dddot X_{\!\Join}^*$;
\item for every nonzero scalar $p\in \IR^*_X$ its inverse $p^{-1}$ is a nonzero scalar.
\end{enumerate}
\end{corollary}

By Corollary~\ref{c:inversionR=>involution}(2), the operation of inversion is an involution on the sets $\dddot X^*_{\!\Join}$ and $\IR^*_X$. 

\begin{definition} A self-map $F:X\to X$ of a set $X$ is called an \index{involution}\defterm{involution} if $F(F(x))=x$ for every $x\in X$.
\end{definition}

Proposition~\ref{p:aff-eq}(4) implies that the following proposition.  
 
\begin{proposition}\label{p:1-portion} For every affine space $X$, there exists a unique unary operation $$1\mbox{-}:\dddot X_{\!\Join}\to\dddot X_{\!\Join},\quad 1\mbox{-}:p\mapsto 1\mbox{-}\,p,$$such that
$1\mbox{-}\,\overvector{zvu}=\overvector{uvz}$ for every line triple $zvu\in \dddot X$. 
\end{proposition}

The unary operation $1\mbox{-}:\dddot X_{\!\Join}\to\dddot X_{\!\Join}$,  introduced in Proposition~\ref{p:1-portion}  is called the \index{$01$-involution}\index[note]{$1$-}\defterm{$01$-involution} of the set of portions in the affine space $X$. 

\begin{proposition} The $01$-involution on an affine space $X$ has following properties:
\begin{enumerate}
\item $1\mbox{-}(1\mbox{-}\,p)=p$ for every portion $p\in\dddot X_{\!\Join}$;
\item $1\mbox{-}\,0=1$ and $1\mbox{-}\,1=0$;
\item for every scalar $p\in\bar\IR_X$, the portion $1\mbox{-}\,p$ is a scalar.
\end{enumerate}
\end{proposition} 

\begin{proof} The properties (1) and (2) follow from Proposition~\ref{p:1-portion} and the definition of the scalars $0$ and $1$;  (3) follows from Proposition~\ref{p:Desarg-triples}(1).
\end{proof}

\section{Scaling vectors by scalars}

The importance of scalars for affine geometry is explained by the fact that the set of scalars $\IR_X$ admits a natural action on the set $X^2$ or pairs, which allows to ``scale distances". This action is the function $\IR_X\times X^2\to X^2$ assigning to every scalar $\sigma\in\IR_X$ and pair $xz\in X^2$ the pair $xy$, where $y\in\Aline xz$ is a unique point such that $xyz\in\{xxx\}\cup \sigma$.

What is even more important is that this action of $\IR_X$ on $X^2$ induces an action of $\IR_X$ on the set of vectors $X^2_{\#}$, which allows one to scale vectors in an affine space $X$.

\begin{theorem}\label{t:scalar-by-vector} For every affine space $X$, there exists a unique function  
$$\mbox{$\cdot:\IR_X\times X^2_{\#}\to X^2_{\#},\quad(\sigma,\vec{\boldsymbol v})\mapsto\sigma\cdot\vec{\boldsymbol v}$}$$such that for every scalar $\sigma\in\IR_X$ and triple $xyz\in X^3$ with $xyz\in\{xxx\}\cup\sigma$ we have $\sigma\cdot\overvector{xz}=\overvector{xy}$.
\end{theorem}

\begin{proof} Define the function $$\cdot:\IR_X\times X^2_{\#}\to X^2_{\#},\quad(\sigma,\vec{\boldsymbol v})\mapsto\sigma\cdot\vec{\boldsymbol v},$$assigning to every scalar $\sigma\in\IR_X$ and vector $\vec{\boldsymbol v}\in X^2_{\#}$ the vector $\sigma\cdot\vec{\boldsymbol v}\defeq\overvector{xy}$ where $xyz\in X^3$ is a triple with $xyz\in\{xxx\}\cup \sigma$ and $xz\in\vec{\boldsymbol v}$.

Let us first check that such a triple $xyz$ exists. Fix any pair $xz\in\vec{\boldsymbol v}$. If $x=z$, then observe that the triple $xyz=xxx$ has the required properties: $xyz\in\{xxx\}\cup\sigma$ and $\overvector{xz}=\overvector{xx}=\vec{\mathbf 0}=\vec{\boldsymbol v}$. If $x\ne z$, then by Proposition~\ref{p:triple-exists}, there exists a point $y\in\Aline xz$ such that $xyz\in\sigma$.

Next, we show that $\sigma\cdot\vec{\boldsymbol v}=\overvector{xz}$ does not depend on the choice of the triple $xyz$. Let $uvw\in X^3$ be another triple such that $uvw\in\{uuu\}\cup\sigma$ and $uw\in\vec{\boldsymbol v}$. If $\vec{\boldsymbol v}=\vec{\mathbf 0}$, then $x=z$ and $u=w$, by Proposition~\ref{p:transeq}(6). In this case $xyz=xxx$ and $uvw=uuu$ and hence $\overvector{xy}=\vec{\mathbf 0}=\overvector{uv}$. 

So, we assume that $\vec{\boldsymbol v}\ne\vec{\mathbf 0}$. In this case, Proposition~\ref{p:transeq}(6) ensures that $x\ne z$ and $u\ne w$. Since $\overvector{xz}=\vec{\boldsymbol v}=\overvector{uw}$, there exists a line translation $T$ in $X$ such that $Txz=uw$. Since $\sigma=\overvector{xyz}=\overvector{uvw}$ is a scalar, the triples $xyz,uvw$ are Desarguesian, and hence $Txyz=uvw$, which implies $Txy=uv$ and $\overvector{xy}=\overvector{uv}$.
\end{proof}

\begin{proposition}\label{p:scalar-by-vector} Let $X$ be an affine space. For every scalar $\sigma\in \IR_X$ and vector $\vecv\in X^2_{\#}$, the following properties hold;
\begin{enumerate}
\item $\sigma\cdot\vec{\mathbf 0}=\vec{\mathbf 0}$;
\item $0\cdot\vecv=\vec{\mathbf 0}$;
\item $1\cdot\vecv=\vecv$;
\item If $\vecv\in\overvector X$, then $\sigma\cdot\vecv\in\overvector X$.
\end{enumerate}
\end{proposition}

\begin{proof} 1. Choose any point $x\in X$ and observe that $xxx\in\{xxx\}\cup\sigma$. Applying Theorem~\ref{t:scalar-by-vector}, we obtain $\sigma\cdot \vec{\mathbf 0}=\sigma\cdot\overvector{xx}=\overvector{xx}=\vec{\mathbf 0}$.
\smallskip

2. Choose any pair $xz\in\vecv$. If $x=z$, then $\vecv=\vec{\mathbf 0}$ and $0\cdot \vecv=0\cdot\vec{\mathbf 0}=\vec{\mathbf 0}=\vecv$, by the preceding item. If $x\ne z$, then $xxz\in 0$ and by Theorem~\ref{t:scalar-by-vector}, $0\cdot\vecv=\overvector{xxz}\cdot\overvector{xz}=\overvector{xx}=\vec{\mathbf 0}$.
\smallskip

3. Choose any pair $xz\in\vecv$. If $x=z$, then $\vecv=\vec{\mathbf 0}$ and $0\cdot \vecv=0\cdot\vec{\mathbf 0}=\vec{\mathbf 0}=\vecv$, by the first  item. If $x\ne z$, then $xzz\in 1$ and by Theorem~\ref{t:scalar-by-vector}, $1\cdot\vecv=\overvector{xzz}\cdot\overvector{xz}=\overvector{xz}=\vecv$.
\smallskip

4. Assuming that $\vecv\in\overvector X$, we shall prove that $\sigma\cdot\vecv\in\overvector X$. If $\sigma=0$ or $\vecv=\vec{\mathbf 0}$, then $\sigma\cdot\vecv=\vec{\mathbf 0}\in\overvector X$ and we are done. So, assume that $\sigma\ne 0$ and $\vecv\ne\vec{\mathbf 0}$. Choose any pair $xz\in\vecv$. Since the vector $\vecv$ is functional, the pair $xz\in\vecv$ is Thalesian. Since $\sigma$ is a scalar, there exists a unique point $y\in\Aline xz$ such that $xyz\in\sigma$.  By Theorem~\ref{t:scalar-by-vector}, $\sigma\cdot\vecv=\overvector{xyz}\cdot\overvector{xz}=\overvector{xy}$. To prove that the vector $\sigma\cdot\vecv=\overvector{xy}$ is functional, take any point $u\in X$ such that $xu\# xy$. We should prove that $u=y$. Since $xu\#xy$, there exists a line translation $T$ in $X$ such that $Txy=xu$. Consider the point $w\defeq T(z)$ and observe that $Txz=xw$ and hence $xz\# xw$. Since the pair $xz$ is Thalesian, the translation equivalence $xz\# xw$ implies $z=w$. Then $Txyz=xuw=xuz$ implies $xyz\Join xuz$. Since $\sigma$ is a scalar, the triple $xyz\in\sigma$ is Desarguesian. In this case, the relation $xyz\Join xuz$ implies $y=u$, witnessing that the vector $\overvector{xy}$ is functional.
\end{proof}

\section{Multiplication of scalars and portions}

Scalars also can scale portions, which leads to the operation of multiplication of scalars and portions. 

\begin{theorem}\label{t:scalar-by-portion} For every affine space $X$, there exists a unique binary operation $$\mbox{$\cdot:\dom[\cdot]\to\IR_X,\quad \cdot:(\alpha,\beta)\mapsto\alpha\cdot\beta$, \ on the set \ $\dom[\cdot]\defeq(\IR_X\times\dddot X_{\!\Join})\cup(\dddot X_{\!\Join}\times\IR_X)$}$$ such that for every pair of nonzero portions $(\alpha,\beta)\in\dom[\cdot]$ the following two conditions are satisfied:
\begin{enumerate}
\item if $\beta=0$, then $\alpha\cdot\beta=0$;
\item $\forall o,x,y,u\in X\;\;(oxy\in\alpha\;\wedge\;oyu\in\beta)\;\Ra\;(\overvector{oxy}\cdot\overvector{oyu}=\overvector{oxu})$.
\end{enumerate}
\end{theorem}

\begin{proof} Given any pair of portions $(\alpha,\beta)\in\dom[\cdot]\defeq(\IR_X\times\dddot X_{\Join})\cup(\dddot X_{\Join}\times\IR_X)$, define the product $\alpha\cdot\beta$ as follows. If $\beta=0$, then put $\alpha\cdot\beta\defeq 0$. If $\beta\ne0$, then choose any line triple $oye\in \beta$. It follows from $\beta\ne 0$ that $y\ne o$. Applying Proposition~\ref{p:triple-exists}, find a point $x\in\Aline oy$ such that $oxy\in\alpha$. Then put $\alpha\cdot\beta\defeq\overvector{oxe}$. Lemma~\ref{l:invariance-scalar-multiplication} ensures that the  scalar $\alpha\cdot\beta$ does not depend on the choice of the Desarguesian triples $oye\in \beta$ and $oxy\in\alpha$. It is clear that the conditions (1), (2) uniquely determine the binary operation $\cdot$.
\end{proof}

\begin{definition} The operation $\cdot:\dom[\cdot]\to \dddot X_{\!\Join}$ defined in Theorem~\ref{t:scalar-by-portion} is called the \index{multiplication}\index{portion!multiplication}\index{scalar!multiplication}\defterm{operation of multiplication of scalars and portions} in $X$.
\end{definition} 

\begin{theorem}\label{t:multiplication-scalar} For every affine space $X$, the multiplication of scalars by portions has the following properties:
\begin{enumerate}
\item $\forall \alpha\in\dddot X_{\!\Join}\;\;(\alpha\cdot 0=0=0\cdot \alpha)$;
\item $\forall \alpha\in\dddot X_{\!\Join}\;\;(\alpha\cdot 1=\alpha=1\cdot \alpha)$;
\item $\forall \alpha\in\IR^*_X\;\;(\alpha\cdot \alpha^{-1}=1=\alpha^{-1}\cdot\alpha)$;
\item $\forall \alpha,\beta\in\IR_X\;\;(\alpha\cdot\beta\in\IR_X)$;
\item $\forall \alpha,\beta\in\IR_X\;\forall\gamma\in \dddot X_{\!\Join}\;\;(\alpha\cdot \beta)\cdot \gamma=\alpha\cdot (\beta\cdot \gamma)$;
\item $\forall \alpha\in\dddot X_{\!\Join}\;\forall \beta,\gamma\in\IR_X\;\;(\alpha\cdot \beta)\cdot \gamma=\alpha\cdot (\beta\cdot \gamma)$.
\end{enumerate}
\end{theorem}

\begin{proof} 1. Fix any portion $\alpha\in\dddot X_{\!\Join}$. The equality $\alpha\cdot 0=0$ follows Theorem~\ref{t:scalar-by-portion}(1). To prove that $0\cdot\alpha=0$, choose any line triple $zae\in\alpha$. If $z=a$, then $\alpha=\overvector{zae}=\overvector{zze}=0$ and then $0\cdot\alpha=0\cdot 0=0$ by Theorem~\ref{t:scalar-by-portion}(1). If $z\ne a$, then $zza\in 0$ and  $0\cdot\alpha=\overvector{zza}\cdot\overvector{zae}=\overvector{zze}=0$, by Theorem~\ref{t:scalar-by-portion}.
\smallskip 

2. Fix any portion $\alpha\in \dddot X_{\!\Join}$ and choose any line triple $zau\in \alpha$. By Theorem~\ref{t:scalar-by-portion}(2), $\alpha\cdot 1=\overvector{zau}\cdot\overvector{zuu}=\overvector{zau}=\alpha$. If $\alpha=0$, then $1\cdot \alpha=1\cdot 0=0=\alpha$, by Theorem~\ref{t:scalar-by-portion}(1). If $\alpha\ne 0$, then $z\ne a$ and hence $zaa\in 1$, by definition of the scalar $1$. Then $1\cdot\alpha=\overvector{zaa}\cdot\overvector{zau}=\overvector{zau}=\alpha$, by Theorem~\ref{t:scalar-by-portion}.
\smallskip

3. Given any nonzero scalar $\alpha\in\IR^*_X$, choose any line triple $oau\in\alpha$. It follows from $\alpha\ne 0$ that $o\ne a$. Then $\alpha^{-1}=\overvector{oua}$, by Proposition~\ref{p:portion-inversion}. Applying Theorem~\ref{t:scalar-by-portion}, we conclude that
 $\alpha\cdot\alpha^{-1}=\overvector{oau}\cdot\overvector{oua}=\overvector{oaa}=1$ and $\alpha^{-1}\cdot\alpha=\overvector{oua}\cdot\overvector{oau}=\overvector{ouu}=1$. 
\smallskip

4. Given any scalars $\alpha,\beta\in\IR_X$, choose any line triple $obu\in\beta$. If $o=b$, then 
$$\alpha\cdot\beta=\alpha\cdot\overvector{obu}=\alpha\cdot\overvector{oou}=\alpha\cdot 0=0\in\IR_X,$$
by the first item. So, assume that $o\ne b$. In this case we can find a unique point $a\in\Aline ob$ such that $oab\in\alpha$. Since the line triples $oab$ and $obu$ are Desarguesian, the triple $oau$ is Desarguesian, by Proposition~\ref{p:Desarg-triples}(3). By Theorem~\ref{t:scalar-by-portion},
$$\alpha\cdot\beta=\overvector{oab}\cdot\overvector{obu}=\overvector{oau}$$and hence the portion $\alpha\cdot\beta$ is a scalar.
\smallskip

5. Fix any scalars $\alpha,\beta\in\IR_X$ and portion $\gamma\in X_{\!\Join}$. By the preceding item $\alpha\cdot\beta\in\IR_X$. If one of the portions $\alpha,\beta,\gamma$ is zero, then $(\alpha\cdot\beta)\cdot\gamma=0=\alpha\cdot(\alpha\cdot\gamma)$, by the first item. So, assume that the portions $\alpha,\beta,\gamma$ are nonzero. Fix any line triple $ocu\in \gamma$. The inequality $\gamma\ne0$ implies $o\ne c$. By Proposition~\ref{p:triple-exists}, there exists a point $b\in\Aline oc$ such that $obc\in\beta$. The inequality $\beta\ne0$ implies $o\ne b$.  By Proposition~\ref{p:triple-exists}, there exists a point $a\in\Aline ob$ such that $oab\in\alpha$. By Theorem~\ref{t:scalar-by-portion},
$$
(\alpha\cdot\beta)\cdot\gamma=(\overvector{oab}\cdot\overvector{obc})\cdot\overvector{ocu}=\overvector{oac}\cdot\overvector{ocu}
=\overvector{oau}=\overvector{oab}\cdot\overvector{obu}=\overvector{oab}\cdot(\overvector{obc}\cdot\overvector{ocu})=\alpha\cdot(\beta\cdot\gamma).$$
\smallskip

6. The last statement can be proved by analogy with the preceding statement.
\end{proof}
 


\begin{remark} Theorem~\ref{t:multiplication-scalar} implies that for every affine space $X$, 
the set of scalars $\IR_X$ endowed with the operation of multiplication of scalars is a monoid with zero, and $\IR_X^*$ is a subgroup of this monoid. 
\end{remark}

Finally, we show that the actions of the monoid $\IR_X$ on the sets $X^2$ and $X^2_{\#}$ agree with the algebraic structure of the monoid $\IR_X$.

\begin{theorem}\label{t:multiplicative-action} For every points $x,y\in X$ in an affine space $X$ and every scalars $\alpha,\beta\in\IR_X$, the following equalities hold:
\begin{enumerate}
\item $\alpha\cdot xx=xx$ and $\alpha\cdot\vec{\mathbf 0}=\vec{\mathbf 0}$;
\item $0\cdot xy=xx$ and $0\cdot\overvector{xy}=\vec{\mathbf 0}$;
\item $1\cdot xy=xy$ and $1\cdot\overvector{xy}=\overvector{xy}$;
\item $\alpha\cdot(\beta\cdot xy)=(\alpha\cdot\beta)\cdot xy$ and $\alpha\cdot(\beta\cdot \overvector{xy})=(\alpha\cdot\beta)\cdot \overvector{xy}$.
\end{enumerate}
\end{theorem}

\begin{proof} It should be mentioned that in the above equalities we use the symbol $\cdot$ for denoting three distinct operations: the multiplication of scalars, the action of scalars on pairs, and the action of scalars on vectors. 
\smallskip

1. Observe that $a\defeq x$ is the unique point such that $xax\in\{xxx\}\cup\alpha$. Then 
$$\alpha\cdot xx=xa=xx\quad\mbox{and}\quad\alpha\cdot\vec{\mathbf 0}=\alpha\cdot \overvector{xx}=\overvector{xa}=\vec{\mathbf 0},$$
by the definition of the actions of $\IR_X$ on $X^2$ and $X^2_{\#}$.
\smallskip

2. If $x=y$, then by the preceding item, $0\cdot xy=xx$ and $0\cdot\overvector{xy}=0\cdot\vec{\mathbf 0}=\vec{\mathbf 0}$. If $x\ne y$, then  $$0\cdot xy=\overvector{xxy}\cdot xy=xx\quad\mbox{and}\quad 0\cdot\overvector{xy}=\overvector{xxy}\cdot\overvector{xy}=\overvector{xx}=\vec{\mathbf 0},$$
by the definition of the actions of $\IR_X$ on $X^2$ and $X^2_{\#}$.
\smallskip

3. If $x=y$, then the first statement ensures that
$$1\cdot xy=1\cdot xx=xx=xy\quad\mbox{and}\quad 1\cdot\overvector{xy}=1\cdot\vec{\mathbf 0}=\vec{\mathbf 0}=\overvector{xy}.$$ 
If $x\ne y$, then  $$1\cdot xy=\overvector{xyy}\cdot xy=xy\quad\mbox{and}\quad 1\cdot\overvector{xy}=\overvector{xyy}\cdot\overvector{xy}=\overvector{xy},$$ by the definition of the actions of $\IR_X$ on $X^2$ and $X^2_{\#}$.
\smallskip

4. If $\beta=0$, then $\alpha\cdot \beta=\alpha\cdot 0=0$, by Theorem~\ref{t:scalar-by-portion}(1). Applying the first and the second items, we conclude that
$$\alpha\cdot(\beta\cdot xy)=\alpha\cdot(0\cdot xy)=\alpha\cdot xx=xx=0\cdot xy=(\alpha\cdot\beta)\cdot 
xy$$ and
$$\alpha\cdot(\beta\cdot \overvector{xy})=\alpha\cdot(0\cdot \overvector{xy})=\alpha\cdot \vec{\mathbf 0}=\vec{\mathbf 0}=0\cdot \overvector{xy}=(\alpha\cdot\beta)\cdot 
\overvector{xy}.$$

So, assume that $\beta\ne 0$. If $x=y$, then by the first item,
$$\alpha\cdot(\beta\cdot xy)=\alpha\cdot(\beta\cdot xx)=\alpha\cdot xx=xx=(\alpha\cdot\beta)\cdot xx=(\alpha\cdot\beta)\cdot xy$$and
$$\alpha\cdot(\beta\cdot \overvector{xy})=\alpha\cdot(\beta\cdot \vec{\mathbf 0})=\alpha\cdot\vec{\mathbf 0}=\vec{\mathbf 0}=(\alpha\cdot\beta)\cdot \vec{\mathbf 0}=(\alpha\cdot\beta)\cdot \overvector{xy}.$$
So, we assume that $x\ne y$. In this case, by Theorem~\ref{t:triple-unique}, there exists a unique point $b\in\Aline xy$ such that $xby\in\beta$. Since $\beta\ne0$, Proposition~\ref{p:aff-eq}(6) ensures that $x\ne b$. Then Theorem~\ref{t:triple-unique} yields a unique point $a\in\Aline xb$ such that $xab\in\alpha$. 
Theorem~\ref{t:scalar-by-portion}(2) ensures that $\alpha\cdot\beta=\overvector{xab}\cdot\overvector{xby}=\overvector{xay}$. Then $$\alpha\cdot(\beta\cdot xy)=\overvector{xab}\cdot(\overvector{xby}\cdot xy) =\overvector{xab}\cdot xb= xa=\overvector{xay}\cdot xy=(\alpha\cdot\beta)\cdot xy$$
and 
$$\alpha\cdot(\beta\cdot \overvector{xy})=\overvector{xab}\cdot(\overvector{xby}\cdot \overvector{xy}) =\overvector{xab}\cdot \overvector{xb}= \overvector{xa}=\overvector{xay}\cdot \overvector{xy}=(\alpha\cdot\beta)\cdot \overvector{xy}$$
by the definition of the action of $\IR_X$ on $X^2$ and $X^2_{\#}$.
\end{proof}



\section{The commutativity of multiplication of scalars}


\begin{lemma}\label{l:RX-commutative} Let $X$ be an affine space. For any scalars $\alpha,\beta\in\IR_X$, the following conditions are equivalent:
\begin{enumerate}
\item $\alpha\cdot\beta=\beta\cdot\alpha$;
\item for any concurrent lines $L,L'\subseteq X$ and points $a,b,c\in L\setminus L'$, $a',b',c'\in L'\setminus L$ and $o\in L\cap L'$ with $oab\in\alpha$ and $ocb\in\beta$, the parallelity relations  $\Aline a{b'}\parallel \Aline {a'}{b}$ and $\Aline b{c'}\parallel \Aline {b'}c$ imply $\Aline a{c'}\parallel \Aline {a'}c$.
\end{enumerate}
\end{lemma}

\begin{proof} To prove that $(2)\Ra(1)$, assume that the condition (2) is satisfied. We need to show that $\alpha\cdot\beta=\beta\cdot\alpha$. If $\alpha$ or $\beta$ equals zero, then $\alpha\cdot\beta=0=\beta\cdot\alpha$, by Theorem~\ref{t:multiplication-scalar}(1). So, we assume that $\alpha\ne0\ne\beta$. Choose any distinct points $o,a\in X$. By Proposition~\ref{p:triple-exists}, there exist points $b\in\Aline oa\setminus\{o\}$ and $c\in\Aline o{b}\setminus\{o\}$ such that $oba\in\alpha\ne0$ and $ocb\in\beta\ne 0$. Since the affine space $X$ has dimension $>1$, there exists a point $a'\in X\setminus \Aline oa$. Since the liner $X$ is affine,  there exist unique points $b',c'\in\Aline o{a'}$ such that $\Aline {b'}a\parallel \Aline b{a'}$ and $\Aline {c'}b\parallel \Aline {b'}c$. The latter parallelity relations imply that $oa'b'\Join oba$ and $ob'c'\Join ocb$, and hence $\alpha=\overvector{oba}=\overvector{oa'b'}$ and $\beta=\overline{ocb}=\overvector{ob'c'}$. The condition (2) ensures that $\Aline a{c'}\parallel \Aline {a'}c$ and hence $oca\Join o{a'}{c'}$ and $\overvector{oca}=\overvector{oa'c'}$. Then $$\alpha\cdot\beta=\overvector{oa'b'}\cdot \overvector{ob'c'}=\overvector{oa'c'}=\overvector{oca}=\overvector{ocb}\cdot\overvector{oba}=\beta\cdot\alpha.$$

\begin{picture}(140,135)(-160,-15)

{\linethickness{1pt}
\put(30,0){\color{red}\line(-1,1){30}}
\put(90,0){\color{red}\line(-1,1){90}}
}

\put(0,0){\line(1,0){100}}
\put(105,-4){$L$}
\put(0,0){\line(0,1){100}}
\put(-4,105){$L'$}
\put(0,90){\color{blue}\line(1,-2){45}}
\put(30,0){\color{blue}\line(-1,2){30}}
\put(0,30){\color{cyan}\line(3,-2){45}}
\put(90,0){\color{cyan}\line(-3,2){90}}

\put(0,0){\circle*{3}}
\put(-5,-10){$o$}
\put(30,0){\circle*{3}}
\put(27,-10){$a$}
\put(45,0){\circle*{3}}
\put(42,-10){$b$}
\put(90,0){\circle*{3}}
\put(87,-10){$c$}
\put(0,30){\circle*{3}}
\put(-10,28){$c'$}
\put(0,60){\circle*{3}}
\put(-10,58){$b'$}
\put(0,90){\circle*{3}}
\put(-10,88){$a'$}
\end{picture}

To prove that $(1)\Ra(2)$, assume that $\alpha\cdot\beta=\beta\cdot\alpha$.
Fix any concurrent lines $L,L'$ in $X$ and points $a,b,c\in L\setminus L'$,  $a',b',c'\in L'
\setminus L$, and $o\in L\cap L'$ such that $oba\in\alpha$, $ocb\in\beta$, $\Aline a{b'}\parallel \Aline {a'}b$ and $\Aline {b'}c\parallel \Aline {b}{c'}$. The latter parallelity relations imply $\overline{oa'b'}=\overline{oba}=\alpha$ and $\overline{ob'c'}=\overline{ocb}=\beta$. By Proposition~\ref{p:Desarg-triples}(3) and Theorem~\ref{t:scalar-by-portion}, the portions $\overvector{oca}$ and $\overvector{oa'c'}$ are scalars such that
$$\overvector{oca}=\overvector{ocb}\cdot\overvector{oba}=\beta\cdot\alpha=\alpha\cdot\beta=\overvector{oa'b'}\cdot\overvector{ob'c'}=\overvector{oa'c'}$$
and hence $oca\Join o{a'}{c'}$, which means that $Foca=o{a'}{c'}$ for some line affinity $F:\Aline oa\to\Aline o{c'}$. By Theorem~\ref{t:paraproj-exists}, there exists a line shear $R:\Aline oa\to\Aline o{c'}$ such that $R(a)=c'$. Since the line triple $oca$ is Desarguesian, $o{a'}{c'}=Foca=Roca$. The definition of a line shear ensures that $\Aline{a'}c\parallel \Aline {c'}a$.
\end{proof}

Lemma~\ref{l:RX-commutative} implies the following characterization of affine spaces with commutative monoid of scalars.

\begin{theorem}\label{t:RX-commutative}  For an affine space $X$, the multiplicative monoid of scalars $\IR_X$ is commutative if and only if for every concurrent lines $L,L'\subseteq X$ and points $a,b,c\in L\setminus L'$, $a',b,c'\in L'\setminus L$ and $o\in L\cap L'$, if $\Aline {a'}b\parallel \Aline a{b'}$, $\Aline b{c'}\parallel \Aline{b'}c$, and the line triples $oba$, $ocb$ are Desarguesian, then $\Aline a{c'}\parallel \Aline{a'}c$.
\end{theorem}




\section{Scalars and homotheties in affine spaces}

In this section we study the interplay between scalars and homotheties of (Thalesian) affine spaces. Let us recall that a \index{homothety}\defterm{homothety} is a dilation possessing a fixed point. For every point $o$ in an affine space $X$, the set $
\Dil_o(X)\defeq\{D\in \Dil(X):D(o)=o\}$ is a subgroup of the dilation group $\Dil(X)$ of the affine space $X$. If the affine space $X$ is Thalesian, then group $\Dil_o(X)$ is isomorphic to the quotient group $\Dil(X)/\Trans(X)$, by Proposition~\ref{p:Dilo=Dil/Trans}.

Let $\IR^*_X=\IR_X\setminus\{0\}$ be the multiplicative group of non-zero scalars in an affine space $X$. 

Given a nonzero scalar $s\in\IR^*_X$ and a point $o\in X$, consider the function 
$s_o:X\to X$ assigning to every point $x\in X$ the unique point $y\in X$ such that $oyx\in s\cup\{ooo\}$.

\begin{proposition} For every point $o$ in an affine space $X$ and every scalar $s\in\IR_X^*$, the function $s_o:X\to X$ is a homothety of $X$ with $s_o(o)=o$.
\end{proposition}

\begin{proof} The definition of the function $s_o$ ensures that $s_o(o)=o$. To see that the function $s_o:X\to X$ is bijective, take any triple $oab\in s$ and consider the portion $t\defeq \overvector{oba}$, which is a scalar, by Proposition~\ref{p:Desarg-triples}. Then for the function $t_o:X\to X$, we obtain $t_os_o=1_X=s_ot_o$. Indeed, since $s$ is a scalar, for every point $x\in X\setminus\{o\}$ there exist a unique point $y\in\Aline ox$ such that $oyx\in s=\overvector{oab}$ and hence $\overvector{oxy}=\overvector{oba}=t$. It follows that  $t_os_o(x)=t_o(y)=x$ and hence $t_os_o=1_X$. By analogy we can prove that $s_ot_o=1_X$. Therefore, $s_o$ is a bijective function of $X$.

To see that $s_o$ is a dilation of $X$, it suffices to show that for every distinct points $x,x'\in X$ and the points $y\defeq s_o(x)$ and $y'\defeq s_o(x')$, the lines  $\Aline x{x'}$ and $\Aline y{y'}$ are parallel. If $x=o$, then $y=s_o(x)=s_o(o)=o$ and $\Aline y{y'}=\Aline o{y'}=\Aline o{x'}=\Aline x{x'}$, so the lines $\Aline x{x'}$ and $\Aline y{y'}$ are parallel. By analogy we can prove that $\Aline x{x'}\parallel \Aline y{y'}$ if $x'=o$.
If $x\ne o\ne x'$ and $o\in\Aline x{x'}$, then $\{y,y'\}\subseteq \Aline ox\cup\Aline o{x'}=\Aline x{x'}$ and hence $\Aline y{y'}=\Aline x{x'}$ and $\Aline y{y'}\parallel \Aline x{x'}$. So assume that $o\notin \Aline x{x'}$. By Theorem~\ref{t:paraproj-exists}, there exists a line shear $R$ such that $Rox=ox'$. Let $y''\defeq R(y)\in\Aline o{x'}$. Then $o{y''}{x'}\Join oyx\in s$ and $o{y'}{x'}\in s$ imply $y''=y'$, because $s^{-1}$ is a scalar. Taking into account that $R$ is a parallel projection with $yy'=yy''=Rxx'$, we conclude that $\Aline y{y'}=\Aline y{y''}\parallel \Aline x{x'}$, witnessing that $s_o$ is a dilation of $X$.  Since $s_o(o)=o$, the dilation $s_o$ is a homothety of $X$.
\end{proof}

\begin{theorem}\label{t:RX=Dilo} For every point $o$ in an affine space $X$, the function $*_o:\IR_X^*\to\Dil_o(X)$, $*_o:s\mapsto s_o$, is an injective group homomorphism. If the affine space $X$ is Thalesian, then the function $*_o:\IR_X^*\to \Dil_o(X)$ is a group isomorphism and hence the groups $\IR^*_X$ and $\Dil(X)/\Trans(X)$ are isomorphic.
\end{theorem}

\begin{proof} To see that the function $*_o:\IR^*_X\to\Dil_o(X)$ is injective, take two distinct nonzero scalars $s,s'\in \IR^*_X$. Fix any point $x\in X\setminus \{o\}$ and find two points $y,y'\in\Aline ox$ such that $oyx\in s$ and $o{y'}x\in s'$. Since $s\ne s'$, the points $y,y'$ are distinct. Then $s_o(x)=y\ne y'=s'_o(x)$ and hence $*_o(s)=s_o\ne s'_o=*_o(s')$, witnessing that the function $*_o$ is injective.
\smallskip

By Theorem~\ref{t:multiplication-scalar}, the set $\IR_X^*$ endowed with the operation of multiplication of scalars is a group. To see that $*_o:\IR_X^*\to\Dil(X)$ is a group homomorphism, take any nonzero scalars $s,t\in \IR^*_X$ and consider their product $p\defeq s\cdot t$. We have to prove that $p_o=s_o\circ t_o$. Given any point $x\in X$ we have to check that $p_o(x)=s_o(t_o(x))$. If $x=o$, then $p_o(x)=o=s_o(t_o(o))=s_o(t_o(x))$ and we are done. So, assume that $x\ne o$. Find a point $y\in\Aline ox\setminus\{o\}$ such that $oyx\in t$, and a point $z\in\Aline oy=\Aline ox$ such that $ozy\in s$. The equality $p=s\cdot t=\overvector{ozy}\cdot \overvector{oyx}=\overvector{ozx}$ implies  $p_o(x)=z=s_o(y)=s_o(t_o(x))$, witnessing that $p_o=s_o\cdot t_o$ and $*_o:\IR^*_X\to\Dil_o(X)$ is a group homomorphism.
\smallskip

Next, assume that the affine space $X$ is Thalesian. In this case we shall prove that the group homomorphism $*_o:\IR^*_X\to\Dil_o(X)$ is surjective. Given any homothety $H\in\Dil_o(X)$, we should find a scalar $s\in\IR^*_X$ such that $s_o=H$. Fix any point $x\in X\setminus \{o\}$ and consider the point $y\defeq H(x)\ne H(o)=o$. Taking into account that $H$ is a homothety  with $H(o)=o$, we conclude that $\Aline oy=\Aline ox$ and hence $s\defeq \overvector{oyx}$ is a well-defined portion. To show that $s$ is a scalar, take any line affinity $A$ in $X$ with $Aox=ox$. By Theorem~\ref{t:A=T'RRT}, $A=TR_n\cdots R_1$ for some line translation $T$ and some line shears $R_1,\dots,R_n$ with $R_i(o)=o$ for all $i\in\{1,\dots,n\}$. Then $o=A(o)=TR_n\cdots R_1(o)=T(o)$ and hence the line translation $T$ is the identity map of the line $\Aline ox$. By Theorem~\ref{t:RH=HR}, $HR_i=R_iH$ for every $i\in\{1,\dots,n\}$ and hence
$$A(y)=R_n\cdots R_1(y)=R_n\cdots R_1H(x)=HR_n\cdots R_1(x)=HA(x)=H(x)=y,$$witnessing that the line triple $oyx$ is Desarguesian and $s=\overvector{oyx}$ is a nonzero scalar. Observe that $s_o(x)=y=H(x)$ and hence $H^{-1}s_o$ is a dilation such that $H^{-1}s_o(o)=o$ and $H^{-1}s_o(x)=x$. By Theorem~\ref{t:dilation-atmost1}, the dilation $H^{-1}s_o$ is the identity map of $X$ and hence $s_o=H$, witnessing that the injective homomorphism $*_o:\IR^*_X\to\Dil_o(X)$ is surjective and hence $*_o$ is a group isomorphism. By Proposition~\ref{p:Dilo=Dil/Trans}, the group $\Dil_o(X)$ is isomorphic to the quotient group $\Dil(X)/\Trans(X)$ and hence the multiplicative group $\IR^*_X$ is isomorphic to the group $\Dil_o(X)/\Trans(X)$. 
\end{proof}


Proposition~\ref{p:Dilo=Dil/Trans} and Theorems~\ref{t:Trans(X)=vecX}, \ref{t:RX=Dilo} imply the following description of the algebraic structure of the group of dilations of a Thalesian affine space.

\begin{corollary}\label{c:Dil=vecXxRX}  For every Thalesian affine space $X$, the group of dilations $\Dil(X)$ is a semidirect product $\overvector X\rtimes\IR_X^*$ of the groups $(\overvector X,+)$ and $(\IR_X^*,\cdot)$.
\end{corollary}


\begin{theorem}\label{t:RXne01=>paraD} Let $X$ be an affine space. If $\IR_X\ne\{0,1\}$, then $X$ is Thalesian.
\end{theorem}

\begin{proof} Assume that the monoid of scalars $\IR_X$ contains some scalar $s\in\IR_X\setminus\{0,1\}$. We shall prove $X$ is a translation affine space. Given any distinct points $x,y\in X$, we should find a translation $T:X\to X$ such that $T(x)=y$.
Choose any point $z\in X\setminus\Aline xy$. Let us recall that the $01$-involution $1\mbox{-}:\IR_X\to\IR_X$ assings to every scalar $\overvector{uvw}$ the scalar $\overvector{wvu}$. Since $s\notin\{0,1\}$, the scalars $s^{-1}$, $(1\mbox{-}s)^{-1}$ and $(1\mbox{-}s^{-1})^{-1}$ are well-defined. Then there exist unique points $a,b\in X$ such that $xaz\in (1\mbox{-}s)^{-1}$ and $zby\in (1\mbox{-}s^{-1})^{-1}$.  The equality $(1\mbox{-}s)^{-1}=\overvector{xaz}$ implies $1\mbox{-}s=\overvector{xza}$, $s=\overvector{azx}$, and $s_a(x)=z$. On the other hand, the equality $(1\mbox{-}s^{-1})^{-1}=\overvector{zby}$ implies $1\mbox{-}s^{-1}=\overvector{zyb}$, $s^{-1}=\overvector{byz}$, and $s^{-1}_b(z)=y$.  Then  $T\defeq s^{-1}_b\circ s_a:X\to X$ is a dilation such that $T(x)=s^{-1}_b(s_a(x))=s_b^{-1}(z)=y$. 

It remains to prove that the dilation $T$ is a translation. In the opposite case, the dilation $T$ has a fixed point $o=T(o)$. Then $o=T(o)=s_b^{-1}(s_a(o))$ and hence $s_b(o)=s_a(o)$. Observe that for the point $o'\defeq s_b(o)=s_a(o)$, we have $\overvector{ao'o}=s=\overvector{bo'o}$, $\overvector{oo'a}=1\mbox{-}s=\overvector{oo'b}$, $\overvector{oao'}=(1\mbox{-}s)^{-1}=\overvector{obo'}$, according to Proposition~\ref{p:Desarg-triples}. Since the triples $oao',obo'$ are Desarguesian, the equality $\overvector{oao'}=\overvector{obo'}$ implies $a=b$ and hence $a=b\in \Aline xz\cap \Aline yz=\{z\}$. Then $1\mbox{-}s=\overvector{xza}=\overvector{xzz}=1$ and $s=0$, which contradicts the choice of the scalar $s\notin\{0,1\}$. This contradiction shows that the dilation $T$ has no fixed points and hence $T$ is a translation. Therefore, $X$ is a translation affine space. By Theorem~\ref{t:paraD<=>translation}, the translation affine space $X$ is Thalesian.
\end{proof}

The following characterization of Desarguesian affine spaces is a counterpart of Theorem~\ref{t:paraD<=>translation} characterizing Thalesian affine spaces as translation spaces. Let us recall that an affine space $X$ is called \index{liner!dilation}\defterm{dilation} if for every distinct collinear points $o,x,y\in X$, there exists a homothety $H\in\Dil_o(X)$ such that $H(x)=y$.

\begin{theorem}\label{t:Des<=>dilation} An affine space $X$ is Desarguesian if and only if $X$ is a dilation space.
\end{theorem}

\begin{proof} If $X$ is a Desarguesian, then for every distinct collinear points $o,x,y$ and the scalar $s\defeq\overvector{oyx}$, the homothety $H\defeq s_o\in\Dil_o(X)$ maps $x$ to $y$, witnessing that the affine space $X$ is dilation.

Next, assume that $X$ is dilation. We claim that the affine space $X$  satisfies the Affine Desargues Axiom. Indeed, take any centrally perspective triangles $abc$ and $a'b'c'$ in $X$ such that $\Aline ab\cap\Aline{a'}{b'}=\varnothing=\Aline bc\cap\Aline{b'}{c'}$. Let $o$ be a unique perspector of the centrally perspective triangles $abc$ and $a'b'c'$. By our assumption, there exists a homothety $H\in\Dil_o(X)$ such that $H(b)=b'$. Since $H$ is a dilation, $H[\Aline ab]\parallel \Aline {a'}{b'}$ and hence $H[\Aline ab]=\Aline {a'}{b'}$. Then $a'\in H[\Aline ab]\cap H[\Aline oa]=H[\Aline ab\cap \Aline oa]=\{H(a)\}$ and hence $a'=H(a)$. By analogy we can prove that $c'=H(c)$. Since $D$ is a dilation of $X$, $\Aline {a'}{c'}=H[\Aline ac]\parallel \Aline ac$. Assuming that $\Aline {a'}{c'}\cap\Aline ac=\varnothing$, we conclude that $\Aline {a'}{c'}=\Aline ac$ and hence $\{c\}=\Aline oc\cap\Aline ac=\Aline o{c'}\cap\Aline{a'}{c'}=\{c'\}$, which contradicts $\Aline bc\cap\Aline{b'}{c'}=\varnothing$. This contradictsion shows that $\Aline {a'}{c'}\cap\Aline ac=\varnothing$, witnessing that the affine space $X$ satisfies the Affine Desargues Axiom. By Theorem~\ref{t:ADA<=>}, the affine liner $X$ is Desarguesian.
\end{proof}

By analogy with Theorem~\ref{t:aa'-translation<=>} one can prove the following ``local'' version of Theorem~\ref{t:Des<=>dilation}.

\begin{theorem}\label{t:local-homothety} For any distinct collinear points $o,a,a'$ in an affine space $X$, there exists a homothety $H:X\to X$ with $Hoa=oa'$ if and only if for any points $b,c\in X\setminus\{o\}$ and $b'\in \Aline ob$ and $c'\in \Aline oc$ with $\Aline ab\parallel \Aline{a'}{b'}$, $\Aline bc\parallel \Aline {b'}{c'}$ and $\Aline oa\cap\Aline ob=\{o\}=\Aline ob\cap\Aline oc$, we have $\Aline ac=\Aline{a'}{c'}$.
\end{theorem}

\begin{problem} Is a finite affine space $X$ Thalesian if for every $x,y\in X$ there exists an automorphism $A:X\to X$ such that $A(x)=y$ (and $\Dil(X)\ne\{1_X\})$?
\end{problem}

\section{Scalars in Thalesian affine spaces}\label{s:scalar-addition}

By Theorem~\ref{t:RXne01=>paraD}, every affine space $X$ with $\IR_X\ne\{0,1\}$ is Thalesian. In this section we shall prove that for any Thalesian affine space $X$ its set of scalars $\IR_X$ has a natural addition operation turning the monoid $\IR_X$ into a corps (= skew-field). 

By Corollary~\ref{c:Thalesian<=>vector=funvector}, every vector in a Thalesian affine space $X$ is functional, and by Corollary~\ref{c:Thalesian-vectors-commutative}, the set $\overvector X=X^2_{\#}$ of (functional) vectors in $X$ is a commutative group with respect to the operation of addition. We shall use this operation of vector addition for defining the operation of addition of scalars in Thalesian affine spaces.

\begin{theorem}\label{t:scalar-addition} For every Thalesian affine space $X$, there exists a unique binary operation $$+:\IR_X\times \IR_X\to\IR_X,\quad +:(\alpha,\beta)\mapsto\alpha+\beta,$$ such that for every distinct points $o,e\in X$ and points $x,y,z\in\Aline oe$ with $\overvector{ox}+\overvector{oy}=\overvector{oz}$, we have $\overvector{oxe}+\overvector{oye}=\overvector{oze}$.
\end{theorem}

\begin{proof} Given any scalars $\alpha,\beta\in\IR_X$, choose any distinct points $o,e\in X$, and find  unique points $x,y\in\Aline ae$ such that $oxe\in\alpha$ and $oye\in\beta$. Since the vector $\overvector{ox}+\overvector{oy}$ is functional, there exists a unique point $z\in X$  such that $\overvector{oz}=\overvector{ox}+\overvector{oy}$. Corollary~\ref{c:vector-subparallel} ensures that $z\in\Aline oe$, so we can put  $\alpha+\beta\defeq \overvector{oze}$. 

Let us show that the portion $\alpha+\beta$ does not depend on the choice of the points $o,e$. 
Choose any distinct points $o',e'\in X$ and find unique points $x',y',z'\in\Aline {o'}{e'}$ such that $o'x'e'\in\alpha$, $o'y'e'\in\beta$, and $\overvector{o'z'}=\overvector{o'x'}+\overvector{o'y'}$.

By Theorem~\ref{t:aff-trans}, there exists a line affinity $A$ in $X$ such that $Aoe=o'e'$. Since the triples $oxe$ and $oye$ are Desarguesian, $A(x)=x'$ and $A(y)=y'$. Theorem~\ref{t:laf+vectors} ensures that $A(z)=z'$ and hence $\overvector{oze}=\overvector{o'z'e'}$.

Finally, we show that the line triple $oze$ is Desarguesian and hence the portion $\alpha+\beta=\overvector{oze}$ is a scalar. By Proposition~\ref{p:Desarg-triple<=>}, it suffices to check that $A(z)=z$ for every line affinity $A$ with $Aoe=oe$. Let $o'x'y'z'\defeq Aoxyz$.  Since $\overvector{oz}=\overvector{ox}+\overvector{oy}$, Theorem~\ref{t:laf+vectors} implies that $\overvector{o'z'}=\overvector{o'x'}+\overvector{o'y'}$. Since $\alpha$ and $\beta$ are scalars, the triples $oxe$ and $oye$ are Desarguesian and hence $o'e'=Aoe=oe$ implies $x'y'=xy$. Then $\overvector{oz'}=\overvector{o'z'}= \overvector{o'x'}+\overvector{o'y'}=\overvector{ox}+\overvector{oy}=\overvector{oz}$ and $A(z)=z'=z$ by the functionality of the vector $\overvector{oz'}=\overvector{oz}$.
\end{proof}

\begin{theorem}\label{t:addition-scalars}
For every Thalesian affine space $X$, the addition of scalars has the following properties:
\begin{enumerate}
\item {\sf Associativity:} $\forall \alpha,\beta,\gamma\in\IR_X\;\;(\alpha+\beta)+\gamma=\alpha+(\beta+\gamma)$;
\item {\sf Commutativity:} $\forall \alpha,\beta\in\IR_X\;\;\alpha+\beta=\beta+\alpha$;
\item {\sf Zero:} $\forall \alpha\in\IR_X\;\alpha+0=\alpha$;
\item {\sf Invertibilty:} $\forall \alpha\in\IR_X\;\exists \beta\in\IR_X\;\;\alpha+\beta=0$.
\end{enumerate}
\end{theorem}

\begin{proof} By Corollaries~\ref{c:Thalesian<=>vector=funvector} and \ref{c:Thalesian-vectors-commutative}, the set $\overvector X=X^2_{\#}$ is a commutative group with respect the operation of addition of (functional) vectors. Fix any distinct points $o,d$ in $X$.
\smallskip

1. Given any scalars $\alpha,\beta,\gamma\in \IR_X$, find unique points $a,b,c\in\Aline od$ such that $oad\in \alpha$, $obd\in\beta$, $ocd\in\gamma$. By the associativity of the addition of (functional) vectors, we have $(\overvector{oa}+\overvector{ob})+\overvector{oc}=\overvector{oa}+(\overvector{ob}+\overvector{oc})$.  Since every vector in the Thalesian affine space $X$ is functional, there exist unique points $x,y,z\in\Aline od$ such that $\overvector{ox}=\overvector{oa}+\overvector{ob}$ and $\overvector{oy}=\overvector{ob}+\overvector{oc}$. Then $$\overvector{oz}=\overvector{oa}+\overvector{ob}+\overvector{oc}=\overvector{ox}+\overvector{oc}=\overvector{oa}+\overvector{oy}.$$ 
By Theorem~\ref{t:scalar-addition}, 
\begin{multline*}(\alpha+\beta)+\gamma=(\overvector{oad}+\overvector{obd})+\overvector{ocd}=\overvector{oxd}+\overvector{ocd}=\overvector{ozd}
=\overvector{oad}+\overvector{oyd}=\overvector{oad}+(\overvector{obd}+\overvector{ocd})=\alpha+(\beta+\gamma).
\end{multline*}
\smallskip

2. Given any scalars $\alpha,\beta\in \IR_X$, find unique points $a,b\in\Aline od$ such that $oad\in \alpha$ and $obd\in\beta$. By the commutativity of the addition of vectors, we have $\overvector{oa}+\overvector{ob}=\overvector{ob}+\overvector{oa}$. Since every vector in the Thalesian affine space $X$ is functional, there exists a unique point $c\in\Aline od$ such that $\overvector{oc}=\overvector{oa}+\overvector{ob}=\overvector{ob}+\overvector{oa}$.
By Theorem~\ref{t:scalar-addition}, $\alpha+\beta=\overvector{oad}+\overvector{obd}=\overvector{ocd}=\overvector{obd}+\overvector{oad}=\beta+\alpha$.  
\smallskip

3. Given any scalar $\alpha\in\IR_X$, find a unique point $a\in\Aline od$ such that $oad\in\alpha$. Since $0=\overvector{ood}$ and $\overvector{oa}+\overvector{oo}=\overvector{oa}$, we can apply Theorem~\ref{t:scalar-addition} and conclude that $\alpha+0=\overvector{oad}+\overvector{ood}=\overvector{oad}=\alpha$. 
\smallskip

4. Given any scalar $\alpha\in\IR_X$, find a unique point $a\in\Aline od$ such that $oad\in\alpha$. Since $\overvector{ao}$ is a functional vector, there exists a unique point $b\in X$ such that $\overvector{ob}=\overvector{ao}$. The definition of the addition of functional vectors in Theorem~\ref{t:vector-addition} ensures that $$\overvector{oa}+\overvector{ob}=\overvector{ob}\circ \overvector{oa}=\overvector{ao}\circ\overvector{oa}=(\overvector{oa})^{-1}\circ \overvector{oa}=1_X=\vec{\mathbf 0}.$$ Then the scalar $\beta\defeq \overvector{obd}$ has the required property: $\alpha+\beta=\overvector{oad}+\overvector{obd}=\overvector{ood}=0.$
\end{proof}

Next, we show that the multiplication of scalars in a Thalesian affine space is distributive over the addition of scalars.

\begin{theorem}\label{t:scalar-distributive} Every scalars $\sigma,\alpha,\beta\in \IR_X$ is a Thalesian affine space $X$ satisfy the following distributivity laws:
$$\sigma\cdot(\alpha+\beta)=\sigma\cdot\alpha+\sigma\cdot\beta\quad\mbox{and}\quad(\alpha+\beta)\cdot\sigma=\alpha\cdot\sigma+\beta\cdot\sigma.$$
\end{theorem}

\begin{proof} 
If $\sigma=0$, then by Theorems~\ref{t:multiplication-scalar}(1) and \ref{t:addition-scalars}(3),
$$\sigma\cdot(\alpha+\beta)=0\cdot(\alpha+\beta)=0=0+0=0\cdot\alpha+0\cdot\beta=\sigma\cdot\alpha+\sigma\cdot\beta$$and
 $$(\alpha+\beta)\cdot\sigma=(\alpha+\beta)\cdot 0=0=0+0=\alpha\cdot 0+\beta\cdot 0=\alpha\cdot\sigma+\beta\cdot\sigma.$$
So, we assume that $\sigma\ne 0$. 

\begin{lemma} $(\alpha+\beta)\cdot\sigma=\alpha\cdot\sigma+\beta\cdot\sigma$.
\end{lemma}

\begin{proof} Choose any line triple $ose\in\sigma$. It follows from $\sigma\ne0$ that $s\ne o$. Since $\alpha$ and $\beta$ are scalars,  there exist unique points $a,b\in\Aline os$ such that $oas\in \alpha$ and $obs\in\beta$. Since every vector in the Thalesian affine space is functional,  there exists a unique point $c\in X$ such that such that $\overvector{oc}=\overvector{oa}+\overvector{ob}$. Corollary~\ref{c:vector-subparallel} ensures that $c\in\Aline os=\Aline oe$. By Theorem~\ref{t:scalar-addition}, $\alpha+\beta=\overvector{oas}+\overvector{obs}=\overvector{ocs}$ and by Theorem~\ref{t:scalar-by-portion},
$$
(\alpha+\beta)\cdot\sigma=(\overvector{oas}+\overvector{obs})\cdot \overvector{ose}=\overvector{ocs}\cdot\overvector{ose}=\overvector{oce}
=\overvector{oae}+\overvector{obe}=\overvector{oas}\cdot\overvector{ose}+\overvector{obs}\cdot\overvector{ose}=\alpha\cdot\sigma+\beta\cdot\sigma.
$$
\end{proof}

\begin{lemma}\label{l:alpha+beta=0} If $\alpha\ne 0\ne\beta$ and $\alpha+\beta=0$, then $\sigma\cdot\alpha+\sigma\cdot\beta=0$.
\end{lemma}

\begin{proof} Choose any distinct points $o,e$ in $X$, and find unique points $a,b\in\Aline oe$ such that $oae\in\alpha$ and $obe\in\beta$. The assumption $\alpha\ne0\ne\beta$ ensures $a\ne o\ne b$. It follows from $\alpha+\beta=0$ that $\overvector{oa}+\overvector{ob}=\vec{\mathbf 0}$. Choose any  point $t\in X\setminus\Aline oe$ and find a unique point $s\in\Aline ot$ such that $ost\in \sigma$. The point $s$ is not equal to $o$ because $\sigma\ne0$. By Theorem~\ref{t:paraproj-exists}, there exist line shears $R_a,R_b$ in $X$ such that $R_aot=oa$ and $R_bot=ob$. Consider the points $a'\defeq R_a(s)$ and $b'\defeq R_b(s)$ on the line $\Aline oe=\Aline oa=\Aline ob$.

Since the vector $\overvector{ob}$ is functional, there exists a unique point $u\in X$ such that  $\overvector{tu}=\overvector{ob}$. By Proposition~\ref{p:=vectors=>parallel-sides}, $\Aline tu\parallel \Aline ob$ and $\Aline  ot\parallel \Aline bu$. By Proposition~\ref{p:=vectors=>parallel-sides}, the equality $\overvector{tu}=\overvector{ob}=-\overvector{oa}=\overvector{ao}$ implies $\Aline at\parallel \Aline ou$. By Theorem~\ref{t:paraproj-exists}, there exist line shears $P_t,P_b$ in $X$ such that $P_tot=ou=P_bob$. Since the line triple $ost\in \sigma$ is Desarguesian and $P_tot=ou=P_bR_bot$, the points $P_t(s)$ and $P_bR_b(s)$ coincide. So, we can consider the point $v\defeq P_t(s)=P_bR_b(s)=P_b(b')\in\Aline ou$ and conclude that $\Aline sv\parallel \Aline tu\parallel \Aline ob'$ and $\Aline {b'}v\parallel \Aline bu\parallel \Aline ot\parallel \Aline {o}s$. Applying Proposition~\ref{p:parallelogram=>vectors=}, we can show that $\overvector{ob'}=\overvector{sv}$.

\begin{picture}(160,100)(-150,-20)

\put(0,0){\line(1,0){80}}
\put(0,0){\line(-1,0){60}}
\put(0,0){\line(0,1){60}}

\put(-60,0){\line(1,1){60}}
\put(60,0){\line(-1,1){60}}

\put(-20,0){\line(1,1){20}}
\put(20,0){\line(-1,1){20}}

\put(20,0){\line(0,1){20}}
\put(0,0){\line(1,1){60}}
\put(0,20){\line(1,0){20}}
\put(60,0){\line(0,1){60}}
\put(0,60){\line(1,0){60}}

\put(-60,0){\circle*{3}}
\put(-63,-8){$a$}
\put(-20,0){\circle*{3}}
\put(-23,-9){$a'$}
\put(0,0){\circle*{3}}
\put(-2,-9){$o$}
\put(20,0){\circle*{3}}
\put(18,-10){$b'$}
\put(60,0){\circle*{3}}
\put(58,-10){$b$}
\put(80,0){\circle*{3}}
\put(77,-10){$e$}

\put(0,60){\circle*{3}}
\put(-1,63){$t$}
\put(0,20){\circle*{3}}
\put(1,22){$s$}

\put(20,20){\circle*{3}}
\put(18,23){$v$}
\put(60,60){\circle*{3}}
\put(58,63){$u$}



\end{picture}

 Since $\Aline ov=\Aline ou\parallel \Aline at\parallel \Aline {a'}s$ and
 $\Aline {a'}o\parallel \Aline sv$, Proposition~\ref{p:parallelogram=>vectors=} implies that $\overvector{sv}=\overvector{a'o}=-\overvector{oa'}$ and hence $\overvector{ob'}=\overvector{sv}=-\overvector{oa'}$. Then  $\overvector{oa'}+\overvector{ob'}=\vec{\mathbf 0}=\overvector{oo}$. It follows from $R_aost=oa'a$ and $R_bost=ob'b$ that $\overvector{oa'a}=\overvector{ob'b}=\overvector{ost}=\sigma$ and hence
$$
\sigma\cdot \alpha+\sigma\cdot \beta=\overvector{oa'a}\cdot\overvector{oae}+\overvector{ob'b}\cdot\overvector{obe}\\
=\overvector{oa'e}+\overvector{ab'e}=\overvector{ooe}=0.
$$
\end{proof}

\begin{lemma} $\sigma\cdot(\alpha+\beta)=\sigma\cdot\alpha+\sigma\cdot\beta$.
\end{lemma}

\begin{proof} If $\alpha=0$, then $$\sigma\cdot(\alpha+\beta)=\sigma\cdot(0+\beta)=\sigma\cdot\beta=0+\sigma\cdot\beta=\sigma\cdot 0+\sigma\cdot\beta=\sigma\cdot\alpha+\sigma\cdot\beta.$$
By analogy we can prove that $\beta=0$ implies $\sigma\cdot(\alpha+\beta)=\sigma\cdot\alpha+\sigma\cdot\beta$. 

So, assume that $\alpha\ne 0\ne\beta$. If $\alpha+\beta=0$, then by Lemma~\ref{l:alpha+beta=0}, 
$$\sigma\cdot \alpha+\sigma\cdot\beta=0=\sigma\cdot 0=\sigma\cdot(\alpha+\beta).$$
So, we assume that $\alpha+\beta\ne 0$. Choose any distinct points $o,e\in X$, and find unique points $a,b\in\Aline oe$ such that $oae\in\alpha$ and $obe\in\beta$. By  Corollary~\ref{c:vector-subparallel}, there exists a point $c\in\Aline oe$ such that $\overvector{oc}=\overvector{oa}+\overvector{ob}$. Choose any point $t\in X\setminus\Aline oe$ and find a unique point $s\in\Aline ot\setminus\{o\}$ such that $ost\in\sigma\ne 0$. Since the vector $\overvector{oa}$ is functional, there exists a unique point $u\in X$ such that $\overvector{tu}=\overvector{oa}$. By Proposition~\ref{p:=vectors=>parallel-sides}, $\overvector{tu}=\overvector{oa}$ implies that $\Aline ot\parallel \Aline au$. It follows from $\overvector{oa}+\overvector{ob}=\overvector{oc}$ that $\overvector{bc}=\overvector{oa}=\overvector{tu}$.  By Proposition~\ref{p:=vectors=>parallel-sides}, $\Aline bc\parallel \Aline tu$ and $\Aline bt\parallel \Aline cu$. By Theorem~\ref{t:paraproj-exists}, there exist unique line rotatations $R_a,R_b,R_c,P_t,P_a,P_c$ in $X$ such that $R_aot=oa$, $R_bot=ob$, and $R_cot=oc$, and $ou=P_tot=P_aoa=P_coc$. Consider the points $a'\defeq R_a(s)$, $b'\defeq R_b(s)$, $c'\defeq R_c(s)$ on the line $\Aline oe$, and the point $v\defeq P_t(s)$ on the line $\Aline ou$. Then $\Aline sv\parallel\Aline tu\parallel\Aline oe$.

\begin{picture}(200,100)(-110,-20)
\put(0,0){\line(1,0){180}}

\put(0,0){\circle*{3}}
\put(-3,-10){$o$}
\put(20,0){\circle*{3}}
\put(17,-10){$a'$}
\put(30,0){\circle*{3}}
\put(28,-10){$b'$}
\put(50,0){\circle*{3}}
\put(47,-10){$c'$}
\put(60,0){\circle*{3}}
\put(58,-10){$a$}
\put(90,0){\circle*{3}}
\put(88,-10){$b$}
\put(150,0){\circle*{3}}
\put(147,-10){$c$}
\put(180,0){\circle*{3}}
\put(178,-10){$e$}

\put(0,0){\line(1,2){30}}
\put(60,0){\line(1,2){30}}
\put(30,60){\line(1,0){60}}
\put(30,60){\line(1,-1){60}}
\put(90,60){\line(1,-1){60}}
\put(30,60){\line(2,-1){120}}
\put(30,60){\line(1,-2){30}}

\put(0,0){\line(1,2){10}}
\put(20,0){\line(1,2){10}}
\put(10,20){\line(1,0){20}}
\put(10,20){\line(1,-1){20}}
\put(30,20){\line(1,-1){20}}

\put(10,20){\line(2,-1){40}}
\put(10,20){\line(1,-2){10}}

\put(0,0){\line(3,2){90}}

\put(10,20){\circle*{3}}
\put(3,19){$s$}
\put(30,20){\circle*{3}}
\put(28,22){$v$}
\put(30,60){\circle*{3}}
\put(24,59){$t$}
\put(90,60){\circle*{3}}
\put(88,63){$u$}

\end{picture}

Observe that $P_cR_cot=ou=P_tot$. Since the line triple $ost\in\sigma\in\IR_X$ is Desarguesian, $P_c(c')=P_cR_c(s)=P_t(s)=v$ and hence $\Aline {c'}v\parallel \Aline cu\parallel \Aline bt$. Since $R_bost=ob'b$, the line $\Aline sb'$ is parallel to $\Aline tb$ and hence $\Aline s{b'}\parallel \Aline {c'}v$.
Applying Proposition~\ref{p:parallelogram=>vectors=}, we can show that $\overvector{sv}=\overvector{b'c'}$.

Next, we show that $\Aline {a'}v\parallel \Aline os$. Observe that $P_aR_aot=P_aoa=ou=P_tot$. Since the line triple $ost\in\sigma\in\IR_X$ is Desarguesian, $P_a(a')=P_aR_a(s)=P_t(s)=v$ and hence $\Aline {a'}v\parallel \Aline au \parallel ot\parallel os$. Since $\Aline o{a'}\parallel \Aline sv$, $\Aline os\parallel \Aline{a'}v$ and $\Aline sv\cap\Aline o{a'}=\varnothing$, we can apply Proposition~\ref{p:parallelogram=>vectors=} and conclude that $\overvector{oa'}=\overvector{sv}=\overvector{b'c'}$. Then $$\overvector{oc'}=\overvector{ob'}+\overvector{b'c'}=\overvector{ob'}+\overvector{oa'}=\overvector{oa'}+\overvector{ob'}.$$
 It follows from $R_aost=oa'a$, $R_bost=ob'b$, $R_cost=oc'c$ that $\overvector{oa'a}=\overvector{ob'b}=\overvector{oc'c}=\overvector{ost}=\sigma$.
Then
$$
\sigma\cdot\alpha+\sigma\cdot\beta=\overvector{oa'a}\cdot\overvector{oae}+\overvector{ob'b}\cdot\overvector{obe}=
\overvector{oa'e}+\overvector{ob'e}=\overvector{oc'e}=\overvector{oc'c}\cdot\overvector{oce}=\sigma\cdot(\overvector{oae}+\overvector{obe})=\sigma\cdot(\alpha+\beta).
$$
\end{proof}
\end{proof}

Theorems~\ref{t:addition-scalars} and \ref{t:scalar-distributive} motivate the following definition.

\begin{definition} A \index{corps}\defterm{corps} is a set $F$ endowed with two binary operations $+,\cdot:F\times F\to F$ and two distinct elements $0,1\in F$ satisfying the following axioms:
\begin{enumerate}
\item $\forall x,y,z\in F\;\;x+(y+z)=(x+y)+z$;
\item $\forall x,y\in F\;\;x+y=y+x$;
\item $\forall x\in F\;\;x+0=x=0+x$;
\item $\forall x\in F\;\exists y\in F\;\;x+y=0=y+x$;
\item $\forall x,y,z\in F\;\;x\cdot(y\cdot z)=(x\cdot y)\cdot z$;
\item $\forall x\in F\;\;x\cdot 1=x=1\cdot x$;
\item $\forall x\in F\setminus\{0\}\;\exists y\in F\;\;x\cdot y=1=y\cdot x$;
\item $\forall a,x,y\in F\;\;a\cdot(x+y)=a\cdot x+a\cdot y$;
\item $\forall a,x,y\in F\;\;(x+y)\cdot a=x\cdot a+y\cdot a$.
\end{enumerate} 
A corps $F$ is called a \index{field}\defterm{field} if $x\cdot y=y\cdot x$ for all $x,y\in F$.
\end{definition}


\begin{exercise} Find an example of a corps which is not a field.
\smallskip

\noindent{\em Hint:} The corps $\mathbb H$ of quaternions.
\end{exercise}

\begin{theorem}\label{t:01*=>+} The addition operation of a corps $R$ is uniquely determined by the operations of multiplication and $01$-involution $1\mbox{-}:R\to R$, $1\mbox{-}:x\mapsto 1-x$.
\end{theorem}

\begin{proof} Observe that the neutral element $0$ of the additive operation is a unique element $z\in R$ such that $z\cdot x=z=x\cdot z$ for every element $x\in R$. Also the neutral element $1$ of the multiplicative operation is a unique element $u\in R$ such that $u\cdot x=x=x\cdot u$ for every $x\in R$. So, $0$ and $1$ can be uniquely recovered from the multiplicative structure of the ring. Since $R^*\defeq R\setminus\{0\}$ is a multiplicative group, for every $x\in R^*\defeq R\setminus\{0\}$ there exists a unique element $x^{-1}\in R^*$ such that $x^{-1}\cdot x=1=x^{-1}\cdot x$.

Observe that for every $y\in R\setminus\{1\}$ the additive inverse $-y$ of $y$ is equal to $(1-y)\cdot(1-(1-y)^{-1})$. The additive inverse $-1$ to $1$ is the unique element $v\in R$ of the set $R\setminus\{ x\cdot(1-x^{-1}):x\in R^*\}$. Therefore, the operations of multiplication and $01$-involution uniquely determine the unary operation of taking the additive inverse in $R$. Finally, observe that for every $x,y\in R$ we have the equality 
$$x+y=\begin{cases}y&\mbox{if $x=0$};\\
x\cdot(1-x^{-1}\cdot(-y))&\mbox{if $x\ne 0$}
\end{cases}
$$witnessing that the operations of multiplication of $01$-involution uniquely determine the operation of addition on $R$.
\end{proof}

%

Theorems~\ref{t:multiplication-scalar}, \ref{t:addition-scalars}, and \ref{t:scalar-distributive} imply

\begin{theorem}\label{t:paraD=>hasRX} For every Thalesian affine space $X$, its set of scalars $\IR_X$ is a corps with respect to the operations of addition and multiplication of scalars. 
\end{theorem} 

By Theorem~\ref{t:01*=>+}, the operation of addition of scalars in a Thalesian affine space is uniquely determined by the operations of scalar multiplication and the involution $1{-}:\IR_X\to\IR_X$, $1{-}:s\mapsto 1-s$. Let us show that this involution agrees with the operation of $01$-involution $1\mbox{-}:X^2_{\#}\to X^2_{\#}$, $1\mbox{-}:\overvector{zvu}\mapsto\overvector{uvz}$, introduced in Section~\ref{s:unary-RX}.

\begin{theorem}\label{t:1-=1-} For every Desarguesian triple $xyz$ in an  affine space $X$, we have the equality 
$$\overvector{xyz}+\overvector{zyx}=\overvector{xzz}=1,$$ which implies  $1\mbox{-}\,s=1-s$ for every scalar $s\in \IR_X$.
\end{theorem}

\begin{proof} If $y=x$, then $$\overvector{xyz}+\overvector{zyx}=\overvector{xxz}+\overvector{zxx}=0+1=1=\overvector{xzz}.$$
If $y=z$, then 
$$\overvector{xyz}+\overvector{zyx}=\overvector{xzz}+\overvector{zzx}=1+0=1=\overvector{xzz}.$$
So, assume that $y\notin \{x,z\}$. In this case $\overvector{xyz}\in \IR_X\setminus\{0,1\}$ and the affine space $X$ is Thalesian, by Theorem~\ref{t:RXne01=>paraD}.

Choose any point $o\in X\setminus\Aline xz$. By Corollary~\ref{c:shear-exists}, there exist unique line shears $R,R',R''$ in $X$ such that $Rzx=zo$, $R'oz=ox$, and $R''xo=xz$. Consider the points $u\defeq R(y)$, $v\defeq R'(u)$, and $w\defeq R''(v)$. The line shears $R,R',R''$ witness that $zyx\Join zuo\Join ovx\Join xwz$ and hence $\overvector{zyx}=\overvector{xwz}$.

\begin{picture}(120,90)(-135,-15)

{\linethickness{0.75pt}
\put(0,0){\color{blue}\vector(1,0){90}}
\put(0,0.8){\color{red}\vector(1,0){30}}
\put(90,0){\color{red}\vector(1,0){30}}
\put(45,45){\color{red}\vector(1,0){30}}
}

\put(0,0){\line(1,1){60}}
\put(120,0){\line(-1,1){60}}
\put(30,0){\line(1,1){45}}
\put(90,0){\line(-1,1){45}}

\put(0,0){\circle*{3}}
\put(-2,-8){$x$}
\put(30,0){\circle*{3}}
\put(28,-8){$y$}
\put(90,0){\circle*{3}}
\put(87,-8){$w$}
\put(120,0){\circle*{3}}
\put(118,-8){$z$}
\put(45,45){\circle*{3}}
\put(36,44){$v$}
\put(75,45){\circle*{3}}
\put(78,44){$u$}
\put(60,60){\circle*{3}}
\put(58,63){$o$}
\end{picture}

 The definitions of the line shears $R,R',R''$ ensure that $xy\# vu\#wz$ and hence $\overvector{xy}=\overvector{wz}$. Then  $\overvector{xy}+\overvector{xw}=\overvector{wz}+\overvector{xw}=\overvector{xw}+\overvector{wz}=\overvector{xz}.$ By Theorem~\ref{t:scalar-addition},
$$\overvector{xyz}+\overvector{zyx}=\overvector{xyz}+\overvector{xwz}=\overvector{xzz}=1.$$

Now we can prove that $1\mbox{-}s=1-s$ for every scalar $s\in\IR_X$. Choose any line triple $xyz\in s$. Since $s$ is a scalar, the line triple $xyz$ is Desarguesian and satisfies the equality $\overvector{xyz}+\overvector{zyx}=1$, which implies the desired equality
$$1\mbox{-}\,s=1\mbox{-}\,\overvector{xyz}=\overvector{zyx}=1-\overvector{xyz}=1-s,$$
by the associativity and commutativity of the group $(\overvector X,+)$.  
\end{proof}

\section{The scalar corps of an affine space}\label{s:RXis-corps}

Theorems~\ref{t:01*=>+} and \ref{t:1-=1-} allow us to introduce a canonical structure of a corps on the set of scalars $\IR_X$ of any (not necessarily Thalesian) affine space $X$. This can be done as follows.

Using the operations of multiplication and $01$-involution, for every $y\in \IR_X\setminus\{1\}$ we can define the opposite scalar $\mbox{-}y$ to $y$ by the formula $$\mbox{-}y\defeq (1\mbox{-}y)\cdot(1\mbox{-}\,(1\mbox{-}y)^{-1}),$$where $(1\mbox{-}y)^{-1}$ is the inverse of the nonzero element $1\mbox{-}y$ in the  group $\IR^*_X$. The opposite scalar $\mbox{-}1$ to $1$ is defined as the unique element of the set $\IR_X\setminus\{y\cdot (1\mbox{-}\,y^{-1}):y\in\IR_X^*\}$. The following lemma shows that the element $\mbox{-}1$ is well-defined.

\begin{lemma} For every affine space $X$, the set $\IR_X\setminus\{y\cdot(1\mbox{-}\,y^{-1}):y\in\IR_X^*\}$ is a singleton.
\end{lemma}
 
\begin{proof} If $\IR_X=\{0,1\}$, then $$\{y\cdot(1\mbox{-}\,y^{-1}):y\in\IR^*_X\}=\{1\cdot (1\mbox{-}1^{-1})\}=\{1\cdot (1\mbox{-}1)\}=\{1\cdot 0\}=\{0\}$$ and hence 
$$\IR_X\setminus\{y\cdot(1\mbox{-}\,y^{-1}):y\in\IR_X^*\}=\IR_X\setminus\{0\}=\{1\}$$is a singleton.

If $\IR_X\ne\{0,1\}$, then by Theorem~\ref{t:RXne01=>paraD}, the affine space $X$ is Thalesian and by Theorem~\ref{t:paraD=>hasRX}, the set of scalars $\IR_X$ endowed with the operations of addition and multiplication is a corp. Moreover, by Theorem~\ref{t:1-=1-}, $1\mbox{-}\,s=1-s$ for every scalar $s\in\IR_X$. In this case, the set
$$\IR_X\setminus\{y\cdot(1\mbox{-}\,y^{-1}):y\in\IR_X^*\}=\IR_X\setminus\{y-1:y\in\IR_X^*\}=\{-1\}$$ 
is the singleton containing the opposite element $-1$ to $1$ in the additive group $(\IR_X,+)$ of the corps $\IR_X$.
\end{proof}
 
Endow $\IR_X$ with the binary operation $+:\IR_X\times\IR_X\to\IR_X$, assigning to every scalars $x,y\in \IR_X$ the scalar 
\begin{equation}\label{eq:dotaddition}
x+y\defeq \begin{cases} y&\mbox{if $x=0$};\\
x{\cdot}(1\mbox{-}(x^{-1}{\cdot}(\mbox{-}y)))&\mbox{if $x\ne 0$}.
\end{cases}
\end{equation}

Theorem~\ref{t:1-=1-} and (the proof of) Theorem~\ref{t:01*=>+} imply that  for every Thalesian affine space $X$, the addition operation defined by the formula~(\ref{eq:dotaddition})  coincides with the operation of addition of scalars introduced in Section~\ref{s:scalar-addition}. Because of that we use the same notation for both operations of addition of scalars on Thalesian affine spaces.

\begin{theorem}\label{t:RX-corps} For every affine space $X$, the set of scalars $\IR_X$ endowed with the operations of addition and multiplication is a corp.
\end{theorem}

\begin{proof} If $\IR_X\ne\{0,1\}$, then the affine space $X$ is Thalesian, by Theorem~\ref{t:RXne01=>paraD}. In this case $\IR_X$ is a corp by Theorem~\ref{t:paraD=>hasRX}.

It remains to prove that $\IR_X$ is a corps if $\IR_X=\{0,1\}$.
In this case,
$\mbox{-}0=(1\mbox{-}\,0){\cdot}(1\mbox{-}(1\mbox{-}\,0)^{-1})=1\cdot (1\mbox{-}\,1)=0$ and $\IR_X\setminus\{y\cdot(1\mbox{-}\,y^{-1}):y\in\IR^*_X\}=\IR_X\setminus\{1{\cdot}(1\mbox{-}(1\mbox{-}\,1)^{-1}\}=\IR_X\setminus\{0\}=\{1\}$, which implies that $\mbox{-}1=1$. By the definition of the addition, $0+x=x$ for every $x\in\IR_X=\{0,1\}$. On the other hand, $1+0=1{\cdot}(1\mbox{-}(1^{-1}{\cdot}(\mbox{-}0)))=(1\mbox{-}(1{\cdot} 0))=1$ and $1+ 1=1{\cdot} (1\mbox{-}(1^{-1}{\cdot} (\mbox{-}1)))=1\mbox{-}(1{\cdot} 1)=0$. Therefore, the addition of scalars in $\IR_X$ coincides with the addition modulo 2. Since $0\cdot 0=0$, $1\cdot 1=1$ and $0\cdot 1=0=1\cdot 0$, the multiplication of scalars in $\IR_X$ coincides with the multiplication modulo 2.
Therefore, $\IR_X$ is the field of residues modulo 2.
\end{proof}

\begin{definition}\label{d:scalar-corps-for-affine-space} For every affine space $X$, the corp of scalars $\IR_X$ is called the \defterm{scalar corp} of the affine space $X$. 
\end{definition}

\begin{exercise} Find an example of a Desarguesian affine space $X$ whose scalar corps $\IR_X$ is not a field.
\smallskip

\noindent{\em Hint:} The quaternion plane $\mathbb H\times\mathbb H$.
\end{exercise}

\section{The scalar corps of an affine liner}

It will be convenient to extend the definition of the scalar corps from affine spaces to all affine liners of rank $\ge 3$. A proper way of such an extension is suggested by the following characterization.

\begin{proposition}\label{p:Steiner-linetriple} Two line triples $abc,xyz\in\dddot X$ in a Steiner affine space $X$ are affinely equivalent if and only if there exists a bijection $F:\Aline xz\to\Aline ac$ such that $Fxyz=abc$.
\end{proposition}

\begin{proof} The ``only if'' part is trivial and follows from Definition~\ref{d:affine-equivalence} and bijectivity of line affinities. To prove the ``if'' part, assume that for two line triples $abc,xyz\in\dddot X$ in a Steiner affine space $X$, there exists a bijection $F:\Aline ac\to\Aline xz$ such that $Fxyz=abc$. By Theorem~\ref{t:aff-trans}, there exists a line affinity $A:\Aline xz\to \Aline ac$ such that $Axz=ac$. We claim that $A(y)=F(y)$. If $y=x$, then $A(y)=A(x)=a=F(x)=F(y)$. If $y=z$, then $A(y)=A(z)=c=F(z)=F(y)$. So, assume that $x\ne y\ne z$. The injectivity of the functions $F$ and $G$ ensures that $a=F(x)\ne F(y)\ne F(z)=c$ and $a=A(x)\ne A(y)\ne A(z)=c$. Then $\{A(y),F(y)\}\subseteq \Aline ac\setminus\{a,c\}$. Since the liner $X$ is Steiner, the set $\Aline ac\setminus\{a,c\}$ is a singleton, which implies $A(y)=F(y)$. Therefore, $Axyz=abc$ and the triples $xyz$ and $abc$ are affinely equivalent.
\end{proof}

Proposition~\ref{p:Steiner-linetriple} suggests the following definition of the scalar corps for any affine liner of rank $\ge 3$. 

Let $X$ be an affine liner of rank $\ge 3$. By Proposition~\ref{t:affine=>Avogadro}, $X$ is $2$-balanced and hence the cardinal $|X|_2\ge 2$ is well-defined. If $X$ is an affine space (i.e., a $3$-long affine regular liner of rank $\ge 3$), then its scalar corps $\IR_X$ has been defined in Definition~\ref{d:scalar-corps-for-affine-space}. So, assume that the affine liner $X$ is not an affine space (which means that either $|X|_2<3$ or $X$ is not regular). Theorem~\ref{t:4-long-affine} ensures that $|X|_2\le 3$. Consider the set $$\dddot X\defeq\{xyz\in X^3:x\ne z\;\wedge\;y\in \Aline xz\}$$whose elements are called \defterm{line triples} in $X$. 

Given two line triples $abc,xyz$ in $X$, we write\index[note]{$abc\prop xyz$} $abc\prop xyz$ and say that $abc$ and $xyz$ are \index{line triples!affinely equivalent}\index{affinely equivalent line triples}\defterm{affinely equivalent} if there exists a bijection $B:\Aline xz\to\Aline ac$ such that $Bxyz=abc$. A triple $zvu\in \dddot X$ is called \defterm{Desarguesian} if $Bzvu=zvu$ for any bijection $B:\Aline zu\to\Aline zu$ with $Bzu=zu$. 

\begin{proposition}\label{p:Desarg-triple} A line triple $zvu\in\dddot X$ is Desarguesian if and only if $v\in\{z,u\}$ or $|X|_2=3$.
\end{proposition}

\begin{exercise} Write down the proof of Proposition~\ref{p:Desarg-triple}.
\end{exercise}

For a line triple $zvu\in \dddot X$, the set $$\overvector{zvu}\defeq\{abc\in\dddot X:abc\prop zvu\}$$is called the \defterm{portion} determined by the line triple $zvu$. The portion $\overvector{zvu}$ is called a \defterm{scalar} if the line triple $zvu$ is Desarguesian. In this case all line triples in $\overvector{zuv}$ are Desarguesian, by Proposition~\ref{p:Desarg-triple}. 

Let \index[note]{${\dddot X}_{\Join}$}
$${\color{magenta}\dddot X_{\!\Join}}\defeq\{\overvector{zvu}:zvu\in\dddot X\}$$be the set of portions, and \index[note]{$\IR_X$}\defterm{$\IR_X$} be the set of scalars in $X$, i.e., the set of portions $\overvector{zvu}\in\dddot X_{\!\Join}$ of Desarguesian line triples  $zvu$ in $X$. Therefore, $\IR_X\subseteq\dddot X_{\!\Join}$.

The set $\IR_X$ contains two important scalars\index[note]{$0$}\index[note]{$1$} 
$$\mbox{\defterm{$0$}}\defeq\{xyz\in X^3:x=y\ne z\}\quad\mbox{and}\quad \mbox{\defterm{$1$}}\defeq\{xyz\in X^3:x\ne y=z\}.$$
Therefore, $\{0,1\}\subseteq\IR_X\subseteq \dddot X_{\!\Join}$.
If $|X|_2>2$, then $\dddot X_{\!\Join}=\{0,1,2\}$, where  
$$2\defeq\{xyz\in\dddot X:x\ne y\ne z\}.$$

The structure of the sets $\IR_X$ and $\dddot X_{\!\Join}$ is described by the following proposition whose proof is left to the reader as an exercise.

\begin{proposition} \begin{enumerate}
\item[\textup{(1)}] If $|X|_2=2$, then $\{0,1\}=\IR_X=\dddot X_{\!\Join}$;
\item[\textup{(2)}] If $|X|_2=3$, then $\IR_X=\dddot X_{\!\Join}=\{0,1,2\}$;
\end{enumerate}
\end{proposition}

The set $\IR_X=\big\{\{0,1\},\{0,1,2\}\big\}$ admits unique addition and multiplication operations that turn $\IR_X$ into a field with zero $0$ and unit $1$. Those unique operations are the operation of addition and multiplication modulo $|\IR_X|=|X|_2$. The obtained field $\IR_X$ is called the \defterm{scalar corps} of the liner $X$.

\begin{exercise} Show that for any Desarguesian line triples $abc$ and $acd$ in $X$, the following identities hold:
\begin{enumerate}
\item $1-\overvector{abc}=\overvector{cba}$;
\item If $a\ne b$, then $\overvector{abc}^{-1}=\overvector{acb}$;
\item $\overvector{abc}\cdot\overvector{acd}=\overvector{abd}$,
\end{enumerate}
where the operations of subtraction, inversion and multiplication are taken in the field $\IR_X$.
\end{exercise}

Thefore, if an affine liner $X$ of rank $\|X\|\ge 3$ is not an affine space, then its scalar corps $\IR_X$ has cardinality $|\IR_X|=3$ if and only if $X$ is Steiner and $|\IR_X|=2$, otherwise. Consequently, the scalar corps of any Steiner affine liner of rank $\ge 3$ is a $3$-element field. 

\chapter{Vector algebra in affine spaces}

\section{The vector space of a Thalesian affine space} 

In this section we shall prove that the set of vectors $\overvector X$ in a Thalesian affine space $X$ carries a natural structure of a module over the scalar corps $\IR_X$ of $X$.

\begin{definition} An \index{$R$-module}\defterm{$R$-module} over a corps $R$ is a nonempty set $\mathbf X$ endowed with binary operations $$+:\mathbf X\times \mathbf X\to \mathbf X,\quad +:(\boldsymbol{x},\boldsymbol{y})\mapsto \boldsymbol{x}+\boldsymbol{y},\quad\mbox{and}\quad \cdot:R\times \mathbf X\to \mathbf X,\quad\cdot:(\alpha,\boldsymbol{x})\mapsto\alpha{\cdot}\boldsymbol{x},$$ that satisfy the axioms:
\begin{enumerate}
\item $\forall \boldsymbol{x},\boldsymbol{y},\boldsymbol{z}\in \mathbf X\;\;(\boldsymbol{x}+\boldsymbol{y})+\boldsymbol{z}=\boldsymbol{x}+(\boldsymbol{y}+\boldsymbol{z})$;
\item $\forall  \boldsymbol{x},\boldsymbol{y}\in \mathbf X\;\;0{\cdot}\boldsymbol x=0{\cdot}\boldsymbol y$;
\item $\forall \alpha,\beta\in R\;\forall \boldsymbol{x}\in \mathbf X\;\;\alpha{\cdot}(\beta{\cdot}\boldsymbol{x})=(\alpha{\cdot}\beta){\cdot}\boldsymbol{x}$;
\item $\forall \boldsymbol{x}\in \mathbf X\;\;1{\cdot} \boldsymbol{x}=\boldsymbol{x}$;
\item $\forall \alpha,\beta\in R\;\forall \boldsymbol{x}\in \mathbf X\;\;(\alpha+\beta){\cdot}\boldsymbol{x}=\alpha{\cdot}\boldsymbol{x}+\beta{\cdot}\boldsymbol{x}$;
\item $\forall \alpha\in R\;\forall \boldsymbol{x},\boldsymbol{y}\in \mathbf X\;\;\alpha{\cdot}(\boldsymbol{x}+\boldsymbol{y})=\alpha{\cdot}\boldsymbol{x}+\alpha{\cdot}\boldsymbol{y}$.
\end{enumerate}
Modules over fields are also called \index{vector space}\defterm{vector spaces}.
\end{definition}

\begin{exercise}\label{ex:R-module+Abelian-group} Let  $\mathbf X$ be an $R$-module over a corps $R$. Prove that $\mathbf X$ endowed with the operation of addition is a commutative group and for every $\boldsymbol x\in \mathbf X$ the element $0{\cdot}\boldsymbol x$ is the neutral element of the group $(\mathbf X,+)$ and the element $(-1){\cdot}\boldsymbol x$ is the inverse element to $\boldsymbol x$ in the group $(\mathbf X,+)$. 
\end{exercise}

\begin{theorem}\label{t:paraD=>RX-module} For every Thalesian affine space $X$, its set of vectors  $\overvector X$ endowed with the operations of addition  $+:\overvector X\times\overvector X\to\overvector X$ and multiplication $\cdot:\IR_X\times \overvector X\to\overvector X$ is a vector space over the scalar corps $\IR_X$.
\end{theorem}

\begin{proof} For the set $\overvector X$ endowed with the operations of addition of vectors and multiplication of scalars by vectors, the axioms (1)--(4) of a vector space have been proved in Theorems~\ref{t:addition-scalars} and \ref{t:multiplicative-action}. It remains to prove the axioms (5) and (6) of distributivity.
\smallskip

5. Fix any scalars $\alpha,\beta\in \IR_X$ and a vector $\vec{\boldsymbol v}\in\overvector{X}$. Choose any pair $oe\in\vec{\boldsymbol v}$. If $o=e$, then $\vec{\boldsymbol v}=\vec{\mathbf 0}$ and hence
 $$(\alpha+\beta)\cdot\vec{\boldsymbol v}=(\alpha+\beta)\cdot \vec{\mathbf 0}=\vec{\mathbf 0}=\vec{\mathbf 0}+\vec{\mathbf 0}=\alpha\cdot \vec{\mathbf 0}+\beta\cdot \vec{\mathbf 0}=\alpha\cdot\vec{\boldsymbol v}+\beta\cdot\vec{\boldsymbol v},$$
by Theorems~\ref{t:multiplicative-action}.

So,  assume that $o\ne e$. Since $\alpha$ and $\beta$ are scalars,  there exist unique points $a,b\in\Aline oc$ such that $oae\in \alpha$ and $obe\in\beta$. Since $\overvector{oa}+\overvector{ob}$ is a functional vector, there exists a unique point $s\in X$ such that  $\overvector{os}=\overvector{oa}+\overvector{ob}$.  Corollary~\ref{c:vector-subparallel} ensures that $s\in\Aline oe$.    By Theorem~\ref{t:scalar-addition}, $\alpha+\beta=\overvector{oae}+\overvector{obe}=\overvector{ose}$ and by Theorem~\ref{t:scalar-by-vector},
$$
(\alpha+\beta)\cdot\vec{\boldsymbol v}=(\overvector{oae}+\overvector{obe})\cdot \overvector{oe}=\overvector{ose}\cdot\overvector{oe}=\overvector{os}
=\overvector{oa}+\overvector{ob}=\overvector{oae}\cdot\overvector{oe}+\overvector{obe}\cdot\overvector{oe}=\alpha\cdot\vec{\boldsymbol v}+\beta\cdot\vec{\boldsymbol v}.
$$

6. The left distributivity law will be derived from its restricted versions for opposite and non-parallel vectors.

We recall that two vectors $\vecbold{v},\vecbold{u}\in\overvector X$ are {\em parallel} if $\Aline pv\parallel \Aline  qu$ for any pairs $pv\in\vecbold{v}$ and $qu\in\vecbold{u}$.

\begin{claim}\label{cl:sx(u+v)} For every non-zero scalar $\sigma\in\IR_X$ and any non-parallel vectors $\vecbold{u},\vecbold{v}\in\overvector X\setminus\{\vec{\mathbf 0}\}$, 
$$\sigma\cdot(\vecbold{u}+\vecbold{v})=\sigma\cdot\vecbold{u}+\sigma\cdot\vecbold{v}.$$
\end{claim}

\begin{proof} Fix any point $o\in X$ and choose points $u,v\in X$ such that $\overvector{ou}=\vecbold{u}$ and $\overvector{ov}=\vecbold{v}$. The non-parallelity of the vectors $\vecbold{u},\vecbold{v}$ ensures that $\Aline ou\nparallel\Aline ov$ and hence $\Aline ou\cap\Aline ov=\{o\}$. It follows from $\vecbold{u}\ne\vec{\mathbf 0}\ne\vecbold{v}$ that $u\ne o\ne v$.

Let $w\in X$ be a unique point such that $\overvector{ow}=\vecbold{v}+\vecbold{u}$. Assuming that $w=o$, we conclude that $\overvector{ow}=\vec{\mathbf 0}$ and hence $\overvector{ov}=\vecbold{v}=-\vecbold{u}=-\overvector{ou}=\overvector{uo}$. By Proposition~\ref{p:=vectors=>parallel-sides}, the equality $\overvector{ov}=\overvector{uo}$ implies $\Aline ov\parallel \Aline uo=\Aline ou$, which means that the vectors $\vecbold{v}$ and $\vecbold{u}$ are parallel. But this  contradicts the choice of the vectors $\vecbold{u},\vecbold{v}$. This contradiction shows that $w\ne o$. Since $\sigma$ is a scalar,  there exists a unique point $s\in\Aline ow$ such that $osw\in\sigma$. The choice of $\sigma\ne0$ ensures that $s\ne o$.

\begin{picture}(160,100)(-150,-20)

\put(0,0){\vector(1,0){60}}
\put(0,0){\vector(0,1){60}}



\put(20,0){\line(0,1){20}}
\put(0,0){\line(1,1){60}}
\put(0,20){\line(1,0){20}}
\put(60,0){\line(0,1){60}}
\put(0,60){\line(1,0){60}}

\put(0,0){\circle*{3}}
\put(-6,-8){$o$}
\put(20,0){\circle*{3}}
\put(18,-10){$u'$}
\put(60,0){\circle*{3}}
\put(60,-8){$u$}

\put(0,60){\circle*{3}}
\put(-7,63){$v$}
\put(0,20){\circle*{3}}
\put(-9,16){$v'$}

\put(20,20){\circle*{3}}
\put(22,16){$s$}
\put(60,60){\circle*{3}}
\put(62,62){$w$}



\end{picture}

By Theorem~\ref{t:paraproj-exists}, there exist line rotations $R_u,R_v$ such that $R_uow=ou$ and $R_vow=ov$. Consider the points $u'\defeq R_u(s)$ and $v'\defeq R_v(s)$. It follows from $o\ne s$ that $u'\ne o\ne v'$. Taking into account that  $R_uosw=ou'u$, $R_vosw=ov'v$ and $\Aline ou\cap\Aline ov=\{o\}$, we conclude that $\Aline {u'}s\parallel \Aline uw\parallel \Aline ov=\Aline o{v'}$ and $\Aline{v'}s\parallel \Aline vw\parallel \Aline ou= \Aline o{u'}$. It follows from $\Aline o{u'}=\Aline ou\nparallel \Aline ov=\Aline o{v'}$ that $u'\notin \Aline o{v'}$ and hence the parallel lines $\Aline {u'}s$ and $\Aline o{v'}$ are disjoint. By Proposition~\ref{p:parallelogram=>vectors=}, 
$\overvector{u's}=\overvector{ov'}$. Then 
$$\overvector{os}=\overvector{ou'}+\overvector{u's}=\overvector{ou'}+\overvector{ov'}.$$
It follows from $R_uosw=ou'u$ and $R_vosw=ov'v$ that $\overvector{ou'u}=\overvector{ov'v}=\overvector{osw}=\sigma$. Applying Theorem~\ref{t:scalar-by-vector}, we conclude that
$$\sigma\cdot\vecbold{u}+\sigma\cdot\vecbold{v}=\overvector{ou'u}\cdot\overvector{ou}+\overvector{ov'v}\cdot\overvector{ov}=\overvector{ou'}+\overvector{ov'}=\overvector{os}=\overvector{osw}\cdot\overvector{ow}=\sigma\cdot(\overvector{ou}+\overvector{ov})=\sigma\cdot(\vecbold{u}+\vecbold{v}).$$
\end{proof}

Two vectors $\vecbold{v},\vecbold{u}\in\overvector X$ are called {\em opposite} if $\vecbold{v}=-\vecbold{u}$.

\begin{claim}\label{cl:u+v=0} For any non-zero scalar $\sigma\in\IR_X$ and any opposite nonzero  vectors $\vecbold{u},\vecbold{v}\in\vec X$, we have $\sigma\cdot\vecbold{u}+\sigma\cdot\vecbold v=\vec{\mathbf 0}$.
\end{claim}

\begin{proof} Choose any pair $ou\in \vecbold{u}$. Find a unique point $v\in X$ such that $\overvector{ov}=\vecbold{v}=-\vecbold{u}=\overvector{uo}$. 
By Proposition~\ref{p:transeq}(6), $\overvector{ou}=\vecbold{u}\ne\vec{\mathbf 0}\ne\vecbold{v}=\overvector{ov}$ implies $u\ne o\ne v$. It follows from $\overvector{ov}=\overvector{uo}$ that $\Aline ov\parallel \Aline uo$ and hence $\Aline ov=\Aline ou$.

Choose any  point $t\in X\setminus\Aline ou$ and find a unique point $s\in\Aline ot$ such that $ost\in \sigma$. The point $s$ is not equal to $o$ because $\sigma\ne0$. By Theorem~\ref{t:paraproj-exists}, there exist line rotations $R_u,R_v$ in $X$ such that $R_uot=ou$ and $R_vot=ov$. Consider the points $u'\defeq R_u(s)$ and $v'\defeq R_v(s)$ on the line $\Aline ou=\Aline ov$. The inequality $s\ne o$ implies $u'\ne o\ne v'$. Taking into account that $R_uost=ou'u$ and $R_vost=ov'v$, we conclude that $\Aline s{u'}\parallel \Aline tu$ and $\Aline s{v'}\parallel \Aline tv$.

Since the vector $\overvector{ov}$ is functional, there exists a unique point $b\in X$ such that  $\overvector{tb}=\overvector{ov}$. By Proposition~\ref{p:=vectors=>parallel-sides}, $\Aline tb\parallel \Aline ov$ and $\Aline  ot\parallel \Aline vb$. By Proposition~\ref{p:=vectors=>parallel-sides}, the equality $\overvector{tb}=\overvector{ov}=-\overvector{ou}=\overvector{uo}$ implies $\Aline ut\parallel \Aline ob$. By Theorem~\ref{t:paraproj-exists}, there exist line rotations $P_t,P_v$ in $X$ such that $P_tot=ob=P_vov$. Since the line triple $ost\in \sigma$ is Desarguesian and $P_tot=ob=P_vR_vot$, the points $P_t(s)$ and $P_vR_v(s)$ coincide. So, we can consider the point $a\defeq P_t(s)=P_vR_v(s)=P_v(v')\in\Aline ob$ and conclude that $\Aline sa\parallel \Aline tb\parallel \Aline ov'$ and $\Aline {v'}a\parallel \Aline vb\parallel \Aline ot\parallel \Aline {o}s$. Applying Proposition~\ref{p:parallelogram=>vectors=}, we can show that $\overvector{ov'}=\overvector{sa}$.

\begin{picture}(160,100)(-150,-20)

\put(0,0){\vector(1,0){60}}
\put(0,0){\vector(-1,0){60}}
\put(0,0){\line(0,1){60}}

\put(-60,0){\line(1,1){60}}
\put(60,0){\line(-1,1){60}}

\put(-20,0){\line(1,1){20}}
\put(20,0){\line(-1,1){20}}

\put(20,0){\line(0,1){20}}
\put(0,0){\line(1,1){60}}
\put(0,20){\line(1,0){20}}
\put(60,0){\line(0,1){60}}
\put(0,60){\vector(1,0){60}}

\put(-60,0){\circle*{3}}
\put(-63,-8){$u$}
\put(-20,0){\circle*{3}}
\put(-23,-9){$u'$}
\put(0,0){\circle*{3}}
\put(-2,-9){$o$}
\put(20,0){\circle*{3}}
\put(18,-10){$v'$}
\put(60,0){\circle*{3}}
\put(58,-9){$v$}

\put(0,60){\circle*{3}}
\put(-1,63){$t$}
\put(0,20){\circle*{3}}
\put(1,22){$s$}

\put(20,20){\circle*{3}}
\put(17,23){$a$}
\put(60,60){\circle*{3}}
\put(58,63){$b$}



\end{picture}

 Since $\Aline oa=\Aline ob\parallel \Aline ut\parallel \Aline {u'}s$ and $\Aline {u'}o\parallel \Aline sa$, Proposition~\ref{p:parallelogram=>vectors=} implies that $\overvector{sa}=\overvector{u'o}=-\overvector{ou'}$ and hence $\overvector{ov'}=\overvector{sa}=-\overvector{ou'}$. Then  $\overvector{ou'}+\overvector{ov'}=\vec{\mathbf 0}$. It follows from $R_aost=ou'u$ and $R_bost=ov'v$ that $\overvector{ou'u}=\overvector{ov'v}=\overvector{ost}=\sigma$ and hence
$$
\sigma\cdot \vecbold{u}+\sigma\cdot \vecbold{v}=\overvector{ou'u}\cdot\overvector{ou}+\overvector{ov'v}\cdot\overvector{ov}\\
=\overvector{ou'}+\overvector{ov'}=\vec{\mathbf 0}.
$$
\end{proof}

With Claims~\ref{cl:sx(u+v)} and \ref{cl:u+v=0} at our disposition, we can present a proof of the equality
$$\sigma\cdot(\vecbold{u}+\vecbold v)=\sigma\cdot\vecbold u+\sigma\cdot\vecbold v$$ for every scalar $\sigma\in\IR_X$ and vectors $\vecbold u,\vecbold v\in\overvector X$.

If $\sigma=0$, then by Theorems~\ref{t:multiplicative-action}(2) and \ref{t:vector-addition}(2),
$$\sigma\cdot(\vecbold{a}+\vecbold{b})=0\cdot(\vecbold{a}+\vecbold{b})=\vec{\mathbf 0}=\vec{\mathbf 0}+\vec{\mathbf 0}=0\cdot\vecbold{a}+0\cdot\vecbold{b}=\sigma\cdot\vecbold{a}+\sigma\cdot\vecbold{b}.$$
So, we assume that $\sigma\ne 0$. 

 If $\vecbold u=\vec{\mathbf 0}$, then $$\sigma\cdot(\vecbold{u}+\vecbold v)=\sigma\cdot(\vec{\mathbf 0}+\vecbold v)=\sigma\cdot\vecbold v=\vecbold {\mathbf 0}+\sigma\cdot\vecbold v=\sigma\cdot \vec{\mathbf 0}+\sigma\cdot\vecbold v=\sigma\cdot\vecbold u+\sigma\cdot\vecbold v,$$
by Theorems~\ref{t:vector-addition}(2) and \ref{t:multiplicative-action}(1).  
By analogy we can prove that $\vecbold v=\vec{\mathbf 0}$ implies $\sigma\cdot(\vecbold u+\vecbold v)=\sigma\cdot\vecbold u+\sigma\cdot\vecbold v$. 

So, assume that $\vecbold{u}\ne \vec{\mathbf 0}\ne\vecbold{v}$. If the vectors $\vecbold u$ and $\vecbold v$ are not parallel, then the equality $\sigma\cdot(\vecbold{u}+\vecbold v)=\sigma\cdot\vecbold u+\sigma\cdot\vecbold v$ follows from Claim~\ref{cl:sx(u+v)}.

So, assume that the vectors $\vecbold{u}$ and $\vecbold{v}$ are parallel. 
If $\vecbold{u}+\vecbold{v}=\vec{\mathbf 0}$, then the equality $$\sigma\cdot(\vecbold{u}+\vecbold v)=\sigma\cdot\vec{\mathbf 0}=\vec{\mathbf 0}=\sigma\cdot\vecbold u+\sigma\cdot\vecbold v$$ follows from  Theorem~\ref{t:multiplicative-action}(1) and  Claim~\ref{cl:u+v=0}.

So, assume that $\vecbold{u},\vecbold{v}$ are parallel non-zero vectors with $\vecbold{u}+\vecbold{v}\ne\vec{\mathbf 0}$. Fix any point $o\in X$ and find unique points $u,v\in X$ such that $\overvector{ou}=\vecbold{u}$ and $\overvector{ov}=\vecbold v$. Since the nonzero vectors $\overvector{ou}=\vecbold u$ and $\overvector{ov}=\vecbold v$ are parallel, the lines $\Aline ou$ and $\Aline ov$ are equal. Choose any point $x\in X\setminus\Aline ou$ and consider the vector $ \overvector{ox}$. It is clear that this vector is not parallel to the vectors $\vecbold u$ and $\vecbold v$. By Theorem~\ref{t:vector-addition}(1,2), $$\vecbold u+\vecbold v=\vecbold u+\vec{\mathbf 0}+\vecbold u=\overvector{ou}+(\overvector{ox}+\overvector{xo})+\overvector{ov}=(\overvector{ou}+\overvector{ox})+(\overvector{xo}+\overvector{ov}).$$
We claim that the vectors $\overvector{ou}+\overvector{ox}$ and $\overvector{xo}+\overvector{ov}=\overvector{xv}$ are not parallel.  In the opposite case, their sum is parallel to the vector $\overline{xv}$. On the other hand, this sum $(\overvector{ou}+\overvector{ox})+(\overvector{xo}+\overvector{ov})=\overvector{ou}+\overvector{ov}=\vecbold u+\vecbold v$ is parallel to the line $\Aline ou$, which is not parallel to the vector $\overvector{xv}$. This contradiction shows that the vectors $\overvector{ou}+\overvector{ox}$ and $\overvector{xo}+\overvector{ov}=\overvector{xv}$ are not parallel. Applying Claims~\ref{cl:sx(u+v)}, \ref{cl:u+v=0}, and Theorem~\ref{t:vector-addition}(1,2), we conclude that
$$
\begin{aligned}
\sigma\cdot(\vecbold u+\vecbold v)&=\sigma\cdot\big((\overvector{ou}+\overvector{ox})+(\overvector{xo}+\overvector{ov})\big)=\sigma\cdot(\overvector{ou}+\overvector{ox})+\sigma\cdot(\overvector{xo}+\overvector{ov})\\
&=(\sigma\cdot\overvector{ou}+\sigma\cdot\overvector{ox})+(\sigma\cdot\overvector{xo}+\sigma\cdot\overvector{ov})=\sigma\cdot\vecbold u+(\sigma\cdot\overvector{xo}+\sigma\cdot\overvector{xo})+\sigma\cdot\vecbold u\\
&=\sigma\cdot\vecbold u+\vec{\mathbf 0}+\sigma\cdot\vecbold v=\sigma\cdot\vecbold u+\sigma\cdot\vecbold v.
\end{aligned}
$$
\end{proof}



\section{The algebraic structure of lines and flats in affine spaces}

In this section we prove that the line relation $\Af$ of a Desarguesian affine space $X$ can be recovered from the vector space structure of $\overvector X$ and the map $\overvector{**}:X\to\overvector X$, $\overvector{**}:(x,y)\mapsto\overvector{xy}$. First we prove the following general fact holding for any affine space.

\begin{proposition}\label{p:lines'=>lines} For every affine space $X$,
$$\{xyz\in X^3:\exists \alpha\in \IR_X\;\;\overvector{xy}=\alpha\cdot\overvector{xz}\}\subseteq\Af.$$
\end{proposition}

\begin{proof} Take any triple $xyz\in X^3$ such that $\overvector{xy}=\alpha{\cdot}\overvector{xz}$ for some scalar $\alpha\in\IR_X$. If $x=z$ or $\alpha=0$, then $\overvector{xy}=\alpha{\cdot}\overvector{xz}=\vec{\mathbf 0}$, by Theorem~\ref{t:multiplicative-action}(1,2). By Proposition~\ref{p:transeq}(6), the equality $\overvector{xy}=\vec{\mathbf 0}$ implies $x=y$. Then the triple $xyz$ belongs to the line relation $\Af$ by the Reflexivity Axiom {\sf(RL)}, see Definition~\ref{d:liner}.

Next, assume that $x\ne z$ and $\alpha\ne 0$. By Propositions~\ref{p:triple-exists}, there exists a point $u\in \Aline xz\setminus\{x\}$ such that $xuz\in\alpha$. Theorem~\ref{t:scalar-by-vector} ensures that $\overvector{xy}=\alpha{\cdot}\overvector{xz}=\overvector{xuz}\cdot\overvector{xz}=\overvector{xu}$. Applying Proposition~\ref{p:=vectors=>parallel-sides}, we conclude that $y\in\Aline xy\parallel \Aline xu=\Aline xz$ and hence $xyz\in\Af$.
\end{proof}

\begin{theorem}\label{t:paraD<=>Desargues} A Thalesian affine space $X$ is Desarguesian if and only if 
$$\Af=\{xyz\in X^3:\exists \alpha\in \IR_X\;\;\overvector{xy}=\alpha{\cdot}\overvector{xz}\}.$$
\end{theorem}

\begin{proof} Consider the set $\Af'\defeq\{xyz\in X^3:\exists \alpha\in \IR_X\;\;\overvector{xy}=\alpha{\cdot}\overvector{xz}\}$. Proposition~\ref{p:lines'=>lines} ensures that $\Af'\subseteq \Af$.
\smallskip

If $\Af'=\Af$, then for every line triple $xyz\in\Af=\Af'$, there exists a scalar $\alpha\in\IR_X$ such that $\overvector{xy}=\alpha{\cdot}\overline{xz}$. By Proposition~\ref{p:triple-exists}, there exists a point $u\in\Aline xz$ such that $xuz\in\alpha$.  Theorem~\ref{t:scalar-by-vector} ensures that $\overvector{xy}=\alpha{\cdot}\overvector{xz}=\overvector{xuz}\cdot\overvector{xz}=\overvector{xu}$. Since the affine space $X$ is Thalesian, the vector $\overvector{xy}=\overvector{xu}$ is functional and hence $u=y$, witnessing that the triple $xyz=xuz\in\alpha$ is Desargesian. Since every line triple in $X$ is Desarguesian, the affine space $X$ is Desarguesian, by Theorem~\ref{t:Desargues<=>3Desargues}.
\smallskip

Next, assuming that the affine space $X$ is Desarguesian, we shall prove that $\Af'=\Af$. Given any triple $xyz\in\Af$, we have to show that $xyz\in \Af'$. If $x=z$, then $x=y=z$ and hence $\overvector{xy}=0\cdot\overvector{xz}$, which implies $xyz\in\Af'$. So, assume that $x\ne z$. Since $X$ is Desarguesian, $xyz$ is a Desarguesian triple and hence $\alpha\defeq\overline{xyz}\in\IR_X$. 
Theorem~\ref{t:scalar-by-vector} ensures that $\overvector{xy}=\overvector{xyz}{\cdot}\overvector{xz}=\alpha{\cdot}\overvector{xz}$, witnessing that $xyz\in\Af'$.
\end{proof}

\begin{exercise} Show that for any distinct points $o,e$ in a Desarguesian affine space $X$, we have the equality 
$$\Aline oe=\{x\in X:\exists \alpha\in\IR_{X}\;\;\overvector{ox}=\alpha{\cdot}\overvector{oe}\}.$$
\end{exercise}

\begin{exercise} Show that for any points $o,u,v$ in a Desarguesian affine space $X$, we have the equality
$$\overline{\{o,u,v\}}=\{x\in X:\exists \alpha,\beta\in\IR_X\;\;\overvector{ox}=\alpha{\cdot}\overvector{ou}+\beta\cdot\overvector{ov}\}.$$
\end{exercise}

\begin{definition}\label{d:R-submodule} Let $\mathbf Y$ be an $R$-module over a corps $R$. A subset $\mathbf X\subseteq \mathbf Y$ is called an \index{$R$-submodule}\defterm{$R$-submodule} of $\mathbf Y$ if it has two properties:
\begin{enumerate}
\item $\forall \boldsymbol x,\boldsymbol y\in\mathbf X\;\;(\boldsymbol x+\boldsymbol y\in\mathbf X)$;
\item $\forall \alpha\in R\;\forall \boldsymbol x\in\mathbf X\;\;(\alpha\cdot\boldsymbol x\in \mathbf X)$.
\end{enumerate}
If $R$ is a field, then an $R$-submodule $\mathbf X$ of an $R$-module $\mathbf Y$ is also called a \index{vector subspace}\defterm{vector subspace} of the vector space $\mathbf Y$.
\end{definition}

\begin{exercise} Let $R$ be a corps and $\mathbf Y$ be an $R$-module.  Show that every $R$-submodule $\mathbf X$ of $\mathbf Y$ is an $R$-module with respect to the operation of addition and multiplication, inherited from the $R$-module $\mathbf Y$.
\end{exercise}

\begin{proposition}\label{p:flat=>R-submodule} For every nonempty flat $A$ in a Thalesian affine space $X$, the set $$\vec A\defeq\{\vec{\boldsymbol x}\in \overvector X:\vec{\boldsymbol x}\subparallel A\}$$is an $\IR_X$-submodule of the $\IR_X$-module $\overvector X$.
\end{proposition} 

\begin{proof} To show that the set $\vec A$ is an $\IR_X$-submodule of the $\IR_X$-module $\overvector X$, we have to check two conditions of Definition~\ref{d:R-submodule}. By Proposition~\ref{p:vector-subparallel}, for every vectors $\vec{\boldsymbol x},\vec{\boldsymbol y}\in\vec A$, their sum $\vec{\boldsymbol x}+\vec{\boldsymbol y}$ belongs to the set $\vec A$. It remains to show that for every vector $\vec{\boldsymbol a}\in\vec A$ and every scalar $s\in\IR_X$, the vector $s\cdot\vec{\boldsymbol a}$ belongs to the set $\vec A$. Fix any pair $xz\in\vec{\boldsymbol a}$. It follows from $\vec{\boldsymbol a}\in\vec A$ that $\vec{\boldsymbol a}\subparallel A$ and hence $\Aline xz\subparallel A$. If $x=z$, then $\vec{\boldsymbol a}=\overvector{xz}=\vec{\mathbf 0}$, $s\cdot\vec{\boldsymbol a}=\vec{\mathbf 0}\subparallel A$ and hence $s\cdot\vec{\boldsymbol a}\in\vec A$.

If $x\ne z$, then by Proposition~\ref{p:triple-exists}, there exists a point $y\in \Aline xz$ such that $xyz\in s$ and then $s\cdot\vec{\boldsymbol a}=\overvector{xyz}\cdot\overvector{xz}=\overvector{xy}$. It follows from $\Aline xy\subseteq \Aline xz\subparallel A$ that $\Aline xy\subparallel A$ and hence $s\cdot\vec{\boldsymbol a}=\overvector{xy}\in\vec A$, by Proposition~\ref{p:vect-subparallel<=>}.
\end{proof}

\begin{theorem}\label{t:flat<=>R-submodule} A subset $A$ of a Desarguesian affine space $X$ is flat in $X$ if and only if the set $\vec A\defeq\{\vecv\in\overvector X:\vecv\subparallel A\}$ is an $\IR_X$-submodule of the $\IR_X$-module $\overvector X$.
\end{theorem}

\begin{proof} By Theorem~\ref{t:ADA=>AMA}, the Desarguesian affine space $X$ is Thalesian. The ``only if'' part  follows from Proposition~\ref{p:flat=>R-submodule}.
\smallskip

To prove the ``if'' part, ssume that the set $\vec A$ is a submodule of $\overvector X$. To prove that the set $A$  is flat in $X$, choose any distinct points $x,z\in A$ and a point $y\in \Aline xz$. We have to show that $y\in A$. By Theorem~\ref{t:Desargues<=>3Desargues}, the line triple $xyz$ is Desaguesian and hence $\overvector{xyz}\in\IR_X$. Since $\vec A$ is an $\IR_X$-submodule of the $\IR_X$-module,
$$\overvector{oy}=\overvector{ox}+\overvector{xy}=\overvector{ox}+\overvector{xyz}\cdot\overvector{xz}=1\cdot\overvector{ox}+\overvector{xyz}\cdot(\overvector{oz}-\overvector{ox})=(1-\overvector{xyz})\cdot\overvector{ox}+\overvector{xyz}\cdot\overvector{oz}\in \vec A$$and hence $y\in A$, by the functionality of the vector $\overvector{oy}$.
\end{proof}

\section{The minimal field of an affine space}

A field $F$ is called \index{minimal field}\index{field!minimal}\defterm{minimal} if $F$ contains no proper subfields. Every minimal field is isomorphic either to the field $\IQ$ of rational numbers or to the field $\IZ/p\IZ$ of residues modulo a prime number $p$. 

Every corps $R$ contains a unique minimal subfield $\underline{R}$, which is the smallest subfield of $R$ that contains $0$ and $1$. The cardinality $|\underline{R}|$ of the minimal subfield of $R$ is denoted by $\har(R)$ and called the \index{characteristic of a corps}\defterm{characteristic} of the corps $R$. The characteristic $\har(R)$ of every corps is an element of the set $\IR\cup\{\w\}$ where $\IP$ is the set of all prime numbers.   If a corps $R$ has finite characteristic $p$, then its minimal field $\underline{R}$ coincides with the set $\{0,1,2,\dots,p-1\}\subseteq R$ and is isomorphic to the $p$-element field $\IZ/p\IZ$. If $\har(R)=\w$, then the minimal field $\underline{R}$ of $R$ coincides with the set of fractions $\{\frac{m\cdot 1}{n\cdot 1}:m\in\IZ,\;n\in\IN\}$ and is isomorphic to the field of rationals $\IQ$. 

For every Thalesian space $X$, its vector space $\overvector X$ is an $\IR_X$-module and also an $\underline{\IR}_X$-module, so $\overvector X$ is a vector space over the minimal field $\underline{\IR}_X$. This fact implies the following important theorem.

\begin{theorem}\label{t:cardinality-pD} Every finite Thalesian affine space $X$ has order $|X|_2=p^n$ and cardinality $|X|=p^{n(\|X\|-1)}$ for some prime number $p$ and some number $n\in\IN$. If $n=1$, then the affine space $X$ is Desarguesian.
\end{theorem}

\begin{proof} Since $X$ is Thalesian, its vector space $\overvector X$ is a vector space over the minimal field $\underline{\IR}_X$, which is finite and its cardinality $p=\har(\IR_X)$ is a prime number. Moreover, by Theorem~\ref{t:flat<=>R-submodule}, for any line $L\subseteq X$, the set $$\overvector L\defeq\{\overvector x\in\overvector X:\overvector x\subparallel L\}$$ is an $\IR_X$-submodule of the $\IR_X$-module $\overvector X$ and hence an $\overvector L$ is a $\underline{\IR}_X$-submodule of $\overvector X$. Then $|\overvector L|=p^n$, where $n$ is the dimension of the vector space $\overvector L$ over the field $\underline{\IR}_X$. 

Fix any distinct points $o,e\in L$. Since the affine space $X$ is Thalesian, the map $L\to\overvector L$, $x\mapsto \overvector {ox}$, is bijective and hence $|X|_2=|L|=|\overvector L|=p^n$. By Corollary~\ref{c:affine-cardinality}, $|X|=|X|_2^{\|X\|-1}=p^{n(\|X\|-1)}$.

If $n=1$, then $|X|_2=|L|=|\overvector L|=p=|\underline{\IR}_X|\le|\IR_X|\le|X|_2$ and hence $|\IR_X|=|X|_2$, which implies that the injective map $\IR_X\to L$, $r\mapsto x+r\cdot\overvector {oe}$, is bijective and hence the Thalesian affine space $X$ is Desarguesian, by Theorem~\ref{t:paraD<=>Desargues}.
\end{proof}

\begin{definition} A liner $X$ is defined to be \index{liner!prime}\defterm{prime} if $X$ is completely regular and its projective completion $\overline X$ has order $|\overline X|_2-1=p$ for some prime number $p$.
\end{definition}

Corollary~\ref{c:affine-spread-completion} implies that an affine liner $X$ is prime if and only if $X$ is regular and has prime order.

Theorems~\ref{t:cardinality-pD}, \ref{t:paraD<=>translation} and \ref{t:partial-translation+prime=>translation} imply the following characterization of Desarguesian prime affine spaces.

\begin{theorem}\label{t:partial-translation<=>} For a prime affine space $X$, the following conditions are equivalent:
\begin{enumerate}
\item $X$ is Desarguesian;
\item $X$ is Thalesian;
\item $X$ is translation;
\item $X$ is $\partial$-translation.
\end{enumerate}
\end{theorem}

Theorem~\ref{t:partial-translation<=>} motivates the following old open problem,

\begin{problem} Is every prime affine space Desarguesian?
\end{problem}

\begin{definition} A liner $X$ is defined to be \index{liner!prime-power}\defterm{prime-power} if $X$ is completely regular and its projective completion $\overline X$ has order $|\overline X|_2-1=p^n$ for some prime number $p$ and some $n\in\IN$.
\end{definition}

Theorem~\ref{t:cardinality-pD} implies that every finite Thalesian affine space is prime-power. This fact suggests the following old open problem.

\begin{problem} Is every finite affine space prime-power?
\end{problem}

Let us observe that some information on the cardinality of a finite (not necessarily Thalesian) affine space $X$ can be deduced from the structure of its translation group $\Trans(X)$. By Proposition~\ref{p:Ax=Bx=>A=B}, the action of
the group $\Trans(X)$ on $X$ is free, which implies that the cardinal $|\Trans(X)| = |\overvector X|$ divides the cardinal $|X|$ of  $X$.
If $\|\overvector X\|\ge 3$, then the group $\Trans(X)$ is elementary Abelian, by Theorem~\ref{t:Trans-commutative}. In this case there exists a unique prime number $p$ such that $T^p = 1_X$ for every
translation $T\in \Trans(X)$. Then $\Trans(X)$ is a $(\IZ/p\IZ)$-module over the prime field $\IZ/p\IZ$ and hence
$|\overvector X|=|\Trans(X)| = p^n$ for some $n\ge 2$.





\begin{problem}  Is there a (necessarily infinite) affine space $X$ such that $\|\overvector X\|\ge 3$ and the group $(\overvector X,+)$ is torsion free and non-divisible?
\end{problem}

\section{The Andr\'e description of Thalesian affine spaces}

In this section we present an algebraic description of Thalesian affine spaces, due to \index[person]{Andr\'e}Andr\'e\footnote{{\bf Johaness Andr\'e} (1925 -- 2011) was a German mathematician, born in Hamburg. He studied in Universities of Hamburg and T\"ubingen. In 1954 he defended his Ph.D. Thesis ``\"Uber nicht-Desarguessche Ebenen mit transitiver Translationsgruppe'' under supervision of G\"unter Pickert (who was a student of Helmut Hasse). Since 1963 Andr\'e worked in Universit\"at des Saarlandes.}, who invented it in 1954. An equivalent geometric description was also suggested by \index[person]{Bruck}Bruck and \index[person]{Bose}Bose\footnote{{\bf Raj Chandra Bose} (1901 -- 1987) was an Indian American mathematician and statistician best known for his work in design theory, finite geometry and the theory of error-correcting codes in which the class of BCH codes is partly named after him. He also invented the notions of partial geometry, association scheme, and strongly regular graph and started a systematic study of difference sets to construct symmetric block designs. He was notable for his work along with S.S.~Shrikhande and E.T.~Parker in their disproof of the famous conjecture made by Leonhard Euler dated 1782 that for no $n$ do there exist two mutually orthogonal Latin squares of order $4n + 2$.} in 1964.

The Andr\'e description of Thalesian affine spaces is based on the notion of a spread of $R$-submodules of an $R$-module.

\begin{definition} Let $X$ be an $R$-module over a corps $R$. A family $\mathcal S$ of $R$-submodules is called a \index{spread of $R$-modules}\index{spread}\defterm{spread of $R$-submodules} in $X$ if $X=\bigcup\mathcal S$ and $A\cap B=\{0\}$ for every distinct $R$-submodules $A,B\in\mathcal S$. A spread of $R$-submodules $\mathcal S$ in $X$ is called \index{spread!direct}\index{direct spread}\defterm{direct} if $|\mathcal S|\ge 3$ and $X=A+B\defeq\{a+ b:a\in A,\;b\in B\}$ for any distinct $R$-submodules $A,B\in\mathcal S$.
\end{definition}

Two $R$-modules $X,Y$ are \index{isomorphic $R$-modules}\defterm{isomorphic} if there exists a bijective function $F:X\to Y$ such that $F(x+y)=F(x)+F(y)$ and $F(r\cdot x)=r\cdot F(x)$ for every $x,y\in X$ and $r\in R$. These two properties of $F$ are called the \defterm{additivity} and the \defterm{$R$-homogeneity} of the isomorphism $F$.

\begin{lemma}\label{l:spread-Rmodules} Let $R$ be a corps and $\mathcal S$ be a direct spread of $R$-submodules in an $R$-module $X$. Then any two $R$-submodules $A,B\in\mathcal S$ are isomorphic. Consequently, $\{0\}\ne S\ne X$ for every $S\in\mathcal S$.
\end{lemma}

\begin{proof} Choose any distinct $R$-submodules $A,B\in\mathcal S$. Since $|\mathcal S|\ge 3$, there exists an $R$-submodule $C\in\mathcal S\setminus\{A,B\}$.
Consider the relation $F\defeq\{(a,b)\in A\times B:a-b\in C\}$. We claim that $F$ is an isomorphism of the $R$-modules.

 To see that $\dom[F]=A$, take any point $a\in A\subseteq X=B+C$ and find a point $b\in B$ such that $a\in b+C$ and conclude that $(a,b)\in F$ and $a\in\dom[F]$.   To see that $\rng[F]=B$, take any  point $b\in B\subseteq X=A+C$ and find a point $a\in A$ such that $b\in a+C$ and conclude that $a\in b-C=b+C$ and hence $(a,b)\in F$ and $b\in\rng[F]$. To see that $F$ is a function, take any pairs $(a,b),(a,b')\in F$ and conclude that $a\in (b+C)\cap (b'+C)$ and hence $b-b'\in (a-C)-(a'-C)=C-C=C$ and $b-b'\in B\cap C=\{0\}$, which implies $b=b'$.
By analogy we can show that $F^{-1}$ is a function. Therefore, $F$ is a  bijective function between the $R$-submodules $A$ and $B$ of $X$.

To see that the biejction $F$ is additive, take any pairs $(x,y),(x',y')\in F$ and consider the pair $(x+x',y+y')\in A\times B$. It follows from $(x,y),(x',y')\in F$ that $x-y\in C$ and $x'-y'\in C$ and hence $(x+x')-(y+y')=(x-y)+(x'-y')\in C+C=C$ and $(x+x',y+y')\in F$, which means that $F(x+x')=y+y'=F(x)+F(x')$.

To see that $F$ is $R$-homogeneous, take any elements $r\in R$ and pair $(x,y)\in F$. Then $y-x\in C$. Since $C$ is an $R$-module, $r\cdot y-r\cdot x=r\cdot (y-x)\in r\cdot C\subseteq C$ and hence $(r\cdot x,r\cdot y)\in F$, which means that $F(r\cdot x)=r\cdot F(x)$.  Therefore, the bijective function $F$ is an isomorphism of the $R$-modules $A$ and $B$.
\smallskip

Assuming that $A=\{0\}$, we conclude that $B=\{0\}$ and hence $X=A+B=\{0\}$ and $\mathcal S=\big\{\{0\}\big\}$, which contradicts $|\mathcal S|\ge 3$. Assuming that $A=X$, we conclude that $B=X\cap B=A\cap B=\{0\}$, which is a contradiction showing that $\{0\}\ne A\ne X$.
\end{proof}

\begin{theorem}[Andre\'e, 1954]\label{t:Andre1} Let $X$ be a Thalesian affine space and $\mathcal L$ be the family of lines in $X$. For every line $L\subseteq X$, the set $\vec L\defeq\{\boldsymbol x\in\overvector X:\boldsymbol x\subparallel L\}$ is an $\IR_X$-submodule of the $\IR_X$-module $\overvector X$. The family $\vec{\mathcal L}\defeq\{\vec L:L\in\mathcal L\}$ is a spread of $\IR_X$-submodules of the $\IR_X$-module $\overvector X$. If $X$ is a plane, then the spread of $\IR_X$-submodules $\vec{\mathcal L}$ is direct and $|\vec{\mathcal L}|\ge 3$.
\end{theorem}

\begin{proof} By Proposition~\ref{p:flat=>R-submodule}, for every line $L\in\mathcal L$, the set $\vec L\defeq\{\boldsymbol x\in\overvector X:\boldsymbol x\subparallel L\}$ is an $\IR_X$-submodule of the $\IR_X$-module $\overvector X$. To see that the family $\vec{\mathcal L}\defeq\{\vec L:L\in\mathcal L\}$ is a spread of $\IR_X$-submodules of $\overvector X$, we need to check that $\bigcup\vec{\mathcal L}=\overvector X$ and either $\vec L\cap\vec\Lambda=\{\vec{\mathbf 0\}}$ or $\vec L=\vec\Lambda$ for every lines $L,\Lambda\in\mathcal L$.

To see that $\overvector X=\bigcup\vec{\mathcal L}$, take any vector $\vec{\boldsymbol v}\in \overvector X$ and choose any pair $xy\in\vec{\boldsymbol v}$. If $x=y$, then $\vec{\boldsymbol v}=\overvector{xy}=\vec{\mathbf 0}\in\vec L$ for every line $L\in\mathcal L$. If $x\ne y$, then $\vec{\boldsymbol v}=\overvector{xy}\in \vec L$ for the line $L\defeq\Aline xy$.

Next, take any lines $L,\Lambda$. If $\vec L\cap\vec\Lambda\ne\{\vec{\mathbf 0}\}$, then there exists a non-zero vector $\vec{\boldsymbol v}\in\vec L\cap\vec\Lambda$. Choose any pair $xy\in\vec{\boldsymbol v}$ and conclude that $\Aline xy\subparallel L$ and $\Aline xy\subparallel \Lambda$. By Corollary~\ref{c:parallel} and Theorem~\ref{t:Proclus-lines}, the latter subparallelity relations imply  $L\parallel \Aline xy\parallel \Lambda$, $L\parallel \Lambda$ and $\vec L=\vec\Lambda$.

Therefore, $\vec{\mathcal L}$ is a spread of $R$-submodules in $\overvector X$. Now assumming that $X$ is a plane, we shall prove that the spread $\vec{\mathcal L}$ is direct. Fix any lines $L,\Lambda\in\mathcal L$ such that $\vec L\ne\vec\Lambda$. In this case, the lines $L,\Lambda$ are not parallel and have a unique common point $o$. Given any vector $\vec{\boldsymbol z}\in\overvector X$, find a unique point $z\in X$ such that $\vec{\boldsymbol z}=\overvector{oz}$. If $z\in L\cup \Lambda$, then $\vec{\boldsymbol v}\in \vec L\cup\vec \Lambda\subseteq \vec L+\vec\Lambda$. So, assume that $z\notin L\cup\Lambda$.

Since the affine space $X$ is regular, it is Playfair. Then there exist unique lines $L_z,\Lambda_z$ such that $z\in L_z\subseteq X\setminus L$ and $z\in\Lambda_z\subseteq X\setminus \Lambda$. The uniqueness of the lines $L_z,\Lambda_z$ implies that $L_z\cap \Lambda\ne\varnothing\ne \Lambda_z\cap L$. So, we can find unique points $x\in L\cap \Lambda_z$ and $y\in \Lambda\cap L_z$. Corollary~\ref{c:parallel-lines<=>} ensures that $\Aline xz=\Lambda_z\parallel \Lambda$. Then $\vec{\boldsymbol v}=\overvector{oz}=\overvector{ox}+\overvector{xz}\in \vec L+\vec\Lambda$.
\end{proof}

Next, we prove the Inverse Andr\'e Theorem, which gives a method of construction of non-Desarguesian Thalesian affine spaces. The proof of the Inverse Andr\'e Theorem uses the following simple lemma on parallelograms in $R$-modules endowed with a spread of $R$-submodules.

\begin{lemma}\label{l:parallelogram-Andre} Let $R$ be a corps and $\mathcal S$ be a  spread of submodules in an $R$-module $X$. For any distinct $R$-modules $A,B\in \mathcal S$ and any points $a,b,c,d\in X$, if $\{b-a,c-d\}\subseteq A$ and $\{c-b,d-a\}\subseteq B$, then  $b-a=c-d$ and $d-a=c-b$.
\end{lemma}

\begin{proof} Observe that $b-a=(d-a)+(c-d)+(b-c)$ and hence
$$(b-a)-(c-d)=(d-a)-(c-d)\in A\cap B=\{0\},$$which implies the desired equalities $b-a=c-d$ and $d-a=c-b$.
\end{proof}

\begin{theorem}[Andr\'e, 1954]\label{t:Andre2} Let $R$ be a corps  and $\mathcal S$ be a direct spread of $R$-submodules of an $R$-module $X$ such that $|S|\ge 3$ for every $R$-submodule $S\in\mathcal S$. Then $X$ endowed with with the family of lines $\mathcal L=\{S+x:S\in\mathcal S,\;x\in X\}$ is a Thalesian Playfair plane such that $\mathcal S=\{L\in\mathcal L:o\in L\}$, where $o$ is the neutral element of the additive group of the $R$-module $X$. 
\end{theorem}

\begin{proof} First we show that the family of lines $\mathcal L$ satisfies the properties {\sf(L1)} and {\sf(L2)} from Theorem~\ref{t:L1+L2} and hence uniquely determines a liner structure on $X$.

Given any distinct points $x,y\in X=\bigcup\mathcal S$, find an $R$-submodule $S\in \mathcal S$ such that $y-x\in S$. Then the line $S+x\in\mathcal L$ contains both points $x$ and $y$. Next, assume that $L\in\mathcal L$ is another line containing the points $x,y$. Find an $R$-submodule $S'\in\mathcal S$ and a point $z\in X$ such that $x,y\in L=S'+z$. Then $y-x\in (S'+z)-(S'+z)=S'-S'=S'$ and hence $0\ne y-x\in S\cap S'$. Since $\mathcal S$ is a spread of $R$-submodules, $S=S'$ and $x\in S'+z=S+z$ implies $x+S=z+S=L$, witnessing that the property {\sf(L1)} is satisfied.

To show that the property {\sf(L2)} is satisfied, take any line $L\in\mathcal L$ and find an $R$-submodule $S\in\mathcal S$ and point $x\in X$ such that $L=S+x$. Since $(X,+)$ is a group, $|L|=|S+x|=|S|\ge 2$. Therefore, the family $\mathcal L$ has the properties {\sf (L1),(L2)} and determines a unique liner structure on $X$. 

Our assumption ensures that every $R$-submodule $S\in\mathcal S$ has cardinality $|S|\ge 3$, which implies that every line $L\in\mathcal L$ has length $|L|\ge 3$ and the liner $X$ is $3$-long. Lemma~\ref{l:spread-Rmodules} implies that the liner $X\ne L$ for every  $L\in\mathcal L$ and hence $\|X\|\ge 3$. To prove that $X$ is a Playfair plane, it suffices to show that for every line $L\in\mathcal L$ and every $x\in X\setminus L$ there exists a unique line $L_x\in\mathcal L$ such that $x\in L_x\subseteq X\setminus L$. Since $L\in\mathcal L$, there exists an $R$-submodule $S\in\mathcal S$  and a point $y\in X$ such that $L=S+y$. Then $L_x\defeq S+x$ is a line containing $x$ and disjoint with the line $L$. Assuming that $L_x'\in\mathcal L\setminus\{L_x\}$ is another line with $x\in L_x'\subseteq X\setminus L$, we can find an $R$-submodule $S'\in\mathcal S$ and a point $x'\in X$ such that $L_x'=S'+x'$. Then $x\in S'+z$ implies $x+S'\subseteq S'+z+S'=z+S'\subseteq x-S'+S'=x+S'$ and $x\in S'+z=S'+x$. It follows from $S'+x=S'+z=L_x'\ne L_x=S+x$ that $S'\ne S$ and hence $X=S'+S=S'-S$. Then there exist elements $s\in S$ and $s'\in S'$ such that $y-x=s'-s$ and hence $y+s=x+s'\in (y+S)\cap (x+S')=L\cap L_x'=\varnothing$, which contradicts the choice of the line $L_x'$. This contradiction shows that $L_x=S+x$ is a unique line that contain the point $x$ and is disjoint with the line $L$, witnessing that the liner $X$ is a Playfair plane.   

To see that the liner $X$ is Thalesian, take any distinct parallel lines $A,B,C\in\mathcal L$ and points $a,a'\in A$, $b,b'\in B$, $c,c'\in C$ such that $\Aline ab\cap\Aline{a'}{b'}=\varnothing=\Aline bc\cap\Aline{b'}{c'}$. 
 Applying Lemma~\ref{l:parallelogram-Andre}, we conclude that $$a-b=a'-b',\quad\mbox{and}\quad c-b=c'-b'.$$
Then $a-c=(a-b)+(b-c)=(a'-b')+(b'-c')=a'-c'$. Find a unique $R$-submodule $S\in\mathcal S$ containing the nonzero element $a-c=a'-c'$ of $R$ and conclude that $\Aline ac=a+S$ and $\Aline {a'}{c'}=a'+S$, witnessing that $\Aline ac\parallel\Aline{a'}{c'}$ and the liner $(X,\mathcal L)$ is Thalesian.

The equality $\mathcal S=\{L\in\mathcal L:0\in L\}$ follows from the definition of the family $\mathcal L$.
\end{proof}

\section{Andr\'e planes}

In this section we apply Theorem~\ref{t:Andre2} to produces examples of non-Desarguesian Thalesian planes, called Andr\'e planes. 

Let $R$ be a corps containing at least three elements and $F$ be a subfield of $R$. Then $R$ can be considered as an $F$-module. 
 Observe that for every $a\in R^*\defeq\{0\}$, the set $S_a\defeq \{(x,x{\cdot}a):x\in R\}$ is an $F$-submodule of the $F$-module  $R\times R$, and $\mathcal S=\{S_a:a\in R^*\}\cup\{\{0\}\times R\}$ is a direct spread of $R$-submodules that determines the standard liner structure of the affine (Desarguesian) plane $R\times R$.  

Now we modify the spread $\mathcal S$ to produce a Thalesian Playfair plane which is not necessarily Desarguesian.

Fix any normal subgroup $H$ of the multiplicative group $R^*\defeq R\setminus\{0\}$ and let  $\Aut(R)$ be the group of $F$-linear automorphisms of the $F$-module $R$. In the group $\Aut(R)$ consider the subgroup $\Aut_H(R)$ consisting of automorphisms $\sigma\in\Aut(R)$ such that $\sigma(x)\in H{\cdot}x=x{\cdot}H$ for all $x\in R^*$.  Let $\sigma_*:R^*\to \Aut_H(R)$, $\sigma_*:a\mapsto \sigma_a$, be any function such that for every $a\in R^*$, the set $\{\sigma_x:x\in H{\cdot}a\}\subseteq\Aut_H(R)$ is a singleton. For every $a\in R^*$ consider the $F$-submodule $\tilde S_a\defeq\{(x,\sigma_a(x){\cdot} a):x\in R\}$ of $R\times R$. 

\begin{proposition} If the $R$-module has finite $F$-dimension, then the family $$\tilde{\mathcal S}\defeq\{\tilde S_a:a\in R^*\}\cup\{\{0\}\times R,R\times\{0\}\}$$ is a direct spread of $F$-modules in the $F$-module $R\times R$.  The set $R\times R$ endowed with the family of lines $\mathcal L=\{S+x:S\in\tilde{\mathcal S},\;x\in R\times R\}$ is a Thalesian Playfair plane, called an \defterm{Andr\'e plane}.
\end{proposition}

\begin{proof} It is easy to see that for every $a\in R^*$ the $F$-submodule $\tilde S_a$ has $\dim_F(\tilde S_a)=\dim_F(R)$. To show that the family $\tilde {\mathcal S}$ is a direct spread of $F$-modules, it suffices to check that $\tilde S_a\cap\tilde S_b=\{(0,0)\}$ for any distinct elements $a,b\in R^*$.  In the opposite case, there exists a point $x\in R^*$ such that $(x,\sigma_a(x){\cdot}a)=(x,\sigma_b(x){\cdot}b)\in \tilde S_a\cap\tilde S_b$. It follows from $\sigma_a,\sigma_b\in\Aut_H(R)$ that $H{\cdot}\sigma_a(x)=H{\cdot}x=H{\cdot}\sigma_b(x)$ and hence $$x{\cdot}H{\cdot}a=H{\cdot}x{\cdot}a=H{\cdot}\sigma_a(x){\cdot}a=H{\cdot}\sigma_b(x){\cdot}b=H{\cdot}x{\cdot}a=x{\cdot}H{\cdot}a,$$
which implies $H{\cdot}a=H{\cdot}b$ and hence $\sigma_a=\sigma_b$. Then the equality $\sigma_a(x){\cdot}a=\sigma_b(x){\cdot} b$ implies $a=b$, which contradicts the choice of the (distinct) elements $a,b$. This contradiction shows that  $\tilde S_a\cap\tilde S_b=\{(0,0)\}$. Therefore, the family $\tilde{\mathcal S}$ is a direct spread of $F$-modules in $R\times R$. By Theorem~\ref{t:Andre2},  the set $R\times R$ endowed with the family of lines $\mathcal L=\{S+x:S\in\tilde{\mathcal S},\;x\in R\times R\}$ is a Thalesian Playfair plane.
\end{proof}


\chapter{Isomorphisms of Desarguesian proaffine spaces}

In this section we classify Desarguesian proaffine regular liners up to an isomorphism.
Observe that two liners $X,Y$ of rank $\|X\|=\|Y\|\le 2$ are isomorphic if and only if $|X|=|Y|$. So, we restrict our attention in this section to liners of rank $>2$. Such liners are called spaces.

More precisely, by a \index{space}\defterm{space} we shall understand any $3$-long regular liner $X$ of rank $\|X\|\ge 3$.

\section{Geometry of modules over corps}

By Theorem~\ref{t:paraD=>RX-module}, every Thalesian affine space $X$ is an $\IR_X$-module over the scalar corps $\IR_X$. This fact motivates studying the geometry of modules over corps in more details. First we present two important examples of $R$-modules.

\begin{example} For a corps $R$ and a set $B$, the set $R^B$ of functions $f:B\to R$ is an $R$-module with respect to the operation of addition of functions and left mutiplication of functions by elements of the corps $R$. The subset $R^{\oplus B}\defeq\{f\in R^B:|\supp(f)|<\w\}$ of finitely supported functions is an $R$-submodule of the $R$-module $R^B$.
\end{example}

For a function $f:B\to R$ with values in a corps $R$, its \defterm{support} is the set
$$\supp(f)\defeq\{x\in B:f(x)\ne 0\}.$$

Every $R$-module $X$ over a corps $R$ carries the \defterm{canonical line relation}
$$\Af=\{xyz\in X^3:\exists r\in R\;\;y-x=r\cdot(z-x)\},$$
turning $X$ into a Desarguesian affine regular liner. This important fact will be proved in Theorem~\ref{t:R-module=>Des-aff-reg}.

First we check that the canonical line relation on an $R$-module is indeed as line relation.

\begin{lemma} For every $R$-module $X$ over a corps $R$, the ternary relation 
$$\Af=\{xyz\in X^3:y\in x+R\cdot(z-x)\}$$
is a line relation on $X$.
\end{lemma}

\begin{proof} We have to verify the axioms {\sf(IL)}, {\sf(RL)}, {\sf(EL)} of Definition~\ref{d:liner}. Denote by $\mathbf 0$ the neutral element of the additive group of the $R$-module $X$.

To check the axiom {\sf(IL)}, take any point $x,y\in X$ with $xyx\in \Af$ and observe that $y\in x+R\cdot(x-x)=\{x\}$ and hence $y=x$.

To check the axiom {\sf(RL)}, take any points $x,y\in X$ and observe that $$\{x,y\}=\{x+0\cdot(y-x),y+1\cdot(y-x)\}$$ and hence $\{xxy,xyy\}\subseteq\Af$.

To check the Exchange Axiom {\sf(EL)}, take any points $a,b,x,y\in X$ such that $axb,ayb\in\Af$ and $x\ne y$. We have to prove that $xay,xby\in\Af$. By the definition of the ternary relation $\Af$, there exist elements $r,s\in R$ such that $x-a=r\cdot(b-a)$ and $y-a=s\cdot(b-a)$. It follows from $x\ne y$ that $r\ne s$.  Since $R$ is a corps, the nonzero element $r-s$ has an inverse $(r-s)^{-1}$ in the multiplicative group $R\setminus\{0\}$. Then $y-x=(s-r)\cdot(b-a)=(r-s)\cdot (a-b)$ and $b-a=(s-r)^{-1}\cdot(y-x)$. It follows from $x-a=r\cdot(b-a)$ that $a-x=r\cdot(a-b)=r\cdot(r-s)^{-1}\cdot(y-x)$ and hence $xay\in\Af$. Also $y-a=s\cdot(b-a)$ implies $y-b=a-b+s\cdot(b-a)=(1-s)\cdot(a-b)$ and 
$$
\begin{aligned}
b-x&=b-y+y-x=(1-s)\cdot(b-a)+y-x=(1-s)\cdot(s-r)^{-1}\cdot(y-x)+(y-x)\\
&=((1-s)\cdot(s-r)^{-1}+1)\cdot(y-x)\in R\cdot(y-x),
\end{aligned}
$$ witnessing that $xby\in\Af$.
\end{proof}

The definition of the canonical line relation on an $R$-module $X$ ensures that for any points $x,y\in X$,
$$\Aline xy=x+R\cdot (y-x),$$
and the family of all lines in the liner $X$ coincides with the family $\{x+R\cdot v:x\in X,\;v\in X\setminus\{\mathbf 0\}\}$, where $\mathbf 0$ is the neutral element of the additive group of the $R$-module $X$. 

\begin{lemma}\label{l:para-lines-in-R-module} Let $X$ be an $R$-module over a corps $R$ and  $aob$ be a triangle in the liner $X$. For every elements $r,s\in R\setminus\{0,1\}$ and points $x\defeq o+r\cdot(a-o)\in\Aline oa$ and $y\defeq o+s\cdot(b-o)\in\Aline ob$, we have $\Aline xy\cap\Aline ab=\varnothing$ if and only if $r=s$.
\end{lemma}

\begin{proof} If $r=s$, then $$\Aline xy=x+R\cdot(y-x)=o+r\cdot(a-o)+R\cdot r\cdot(b-a)\subseteq a+(1-r)\cdot(o-a)+R\cdot (b-a).$$Assuming that $\Aline xy\cap \Aline ab\ne\varnothing$, we conclude that $(a+(1-r)\cdot(o-a)+R\cdot(b-a))\cap(a+R\cdot(b-a))=\Aline xy\cap\Aline ab\ne\varnothing$ and hence $(1-r)\cdot(o-a)\in R\cdot(b-a)$ and $o=a+(1-r)^{-1}\cdot R\cdot (b-a)\subseteq a+R\cdot (b-a)=\Aline ab$, which contradicts the choice of the points $a,o,b$ that form a triangle in the liner $X$. This contradiction shows that $\Aline xy\cap\Aline ab=\varnothing$ and proves the ``if'' part of the lemma.
\smallskip

To prove the ``only if'' part, assume that $r\ne s$. Consider the element $t\defeq (1-r)\cdot(s-r)^{-1}$ of the corps $R$ and observe that the point
$$
\begin{aligned}
x+t\cdot(y-x)&=o+r\cdot(a-o)+t\cdot(s\cdot(b-o)-r\cdot(a-o))\\
&=a+(r-1)\cdot(a-o)+t\cdot(s\cdot(b-a)+s\cdot(a-o)-r\cdot (a-o))\\
&=a+t\cdot r\cdot(b-a)+(r-1)\cdot(a-o)+t\cdot(s-r)\cdot(a-o)\\
&=a+t\cdot r\cdot(b-a)+(r-1)\cdot(a-o)+(1-r)\cdot(s-r)^{-1}\cdot(s-r)\cdot(a-o)\\
&=a+t\cdot r\cdot(b-a)\in a+R\cdot(b-a)=\Aline ab
\end{aligned}$$belongs to the intersection $\Aline xy\cap\Aline ab$, witnessing that $\Aline xy\cap\Aline ab$.
\end{proof}

Next, we describe the algebraic structure of flats in modules over corps.

\begin{lemma}\label{l:flat<=>R-submodule} Let $X$ be an $R$-module over a corps $R\ne \{0,1\}$. For every nonempty set $A\subseteq X$, the following conditions are equivalent:
\begin{enumerate}
\item $A$ is flat in the liner $X$;
\item for every $o\in A$, the set $A-o$ is an $R$-submodule of the $R$-module;
\item for some $o\in A$, the set $A-o$ is an $R$-submodule of the $R$-module.
\end{enumerate}
\end{lemma}

\begin{proof} Let $A$ be a nonempty subset of $X$.
\smallskip

$(1)\Ra(2)$ Assume that the set $A$ is flat in the liner $X$. Given any point $o\in A$, we should check that $A-o$ is an $R$-submodule of the $R$-module $X$. 

First we show that $(x-o)+(y-o)\in A-o$ for every points $x,y\in A$. Since $R\ne\{0,1\}$, there exists an element $r\in R\setminus\{0,1\}$. Since $R$ is a corps, the nonzero elements $r,1-r$ have multiplicative inverses $r^{-1}$ and $(1-r)^{-1}$ in the multiplicative group $R\setminus\{0\}$. Since the set $A$ is flat, the points $a\defeq o+r^{-1}\cdot (x-o)\in\Aline ox$ and $b\defeq o+(1-r)^{-1}\cdot (y-o)\in\Aline oy$ belong to $A$ and so does the point 
$$b+r\cdot(a-b)=r\cdot a+(1-r)\cdot b=r\cdot(o+r^{-1}\cdot (x-o))+(1-r)\cdot(o+(1-r)^{-1}\cdot(y-o))=1\cdot o+ (x-o)+(y-o).$$
Then $(x-o)+(y-o)\in A-o$.

On the other hand, for every elements $x\in A$ and $r\in R$, the point $o+r\cdot(x-o)$ belongs to the line $o+R\cdot(x-o)=\Aline ox\subseteq A$ and hence $r\cdot(x-o)\in A-o$, witnessing that $A-o$ is an $R$-submodule of the $R$-module $X$.
\smallskip

The implication $(2)\Ra(3)$ is trivial.
\smallskip

$(3)\Ra(1)$ Assume that for some point $o\in A$, the set $A-o$ is an $R$-submodule of the $R$-module $X$. To see that $A$ is flat, take any points $x,y\in A$ and $z\in\Aline xy$. The definition of the line relation on $X$ ensures that $z=x+r\cdot(y-x)$ for some $r\in R$. Then 
$$
\begin{aligned}
z-o&=x+r\cdot(y-x)-o=(1-r)\cdot x+r\cdot y-(1-r+r)\cdot o\\
&=(1-r)\cdot(x-o)+r\cdot(y-o)\in R\cdot(A-o)+R\cdot(A-o)\subseteq A-o
\end{aligned}
$$ and hence $z\in A$.
\end{proof}

\begin{proposition}\label{p:flathull-in-R-module} Let $X$ be an $R$-module over a corps $R\ne\{0,1\}$. For every set $B\subseteq X$ its flat hull $\overline B$ in the liner $X$ coincides with the set
$$\textstyle\Sigma\defeq \{\sum_{b\in B}f(b)\cdot b:f\in R^{\oplus B},\;\sum_{b\in B}f(b)=1\}.$$
\end{proposition}

\begin{proof} If $B=\varnothing$, then $\overline B=\varnothing=\Sigma$ and we are done. So, assume that $B\ne\varnothing$ and fix any point $o\in B$. By Lemma~\ref{l:flat<=>R-submodule}, the set $\overline B-o$ is an $R$-submodule of the $R$-module $X$. The set $$
\begin{aligned}
\Sigma-o&\textstyle=\{\sum_{b\in B}f(b)\cdot (b-o):f\in R^{\oplus B},\;\sum_{b\in B}f(b)=1\}\\
&\textstyle=\{(1-\sum_{b\in B\setminus\{o\}}f(b))\cdot(o-o)+\sum_{b\in B\setminus\{o\}}f(b)\cdot (b-o):f\in R^{\oplus B\setminus\{o\}}\}\\
&\textstyle=\{\sum_{b\in B\setminus\{o\}}f(b)\cdot b:f\in R^{\oplus B\setminus\{o\}}\}
\end{aligned}
$$also is an $R$-submdule of the $R$-module $X$. Moreover, $\Sigma-o$ is the smallest $R$-submodule of $X$ that contains the set $B-o$. Then $\Sigma-o\subseteq\overline{B}-o$ and hence $\Sigma\subseteq\overline B$.

On the other hand, since the set $\Sigma-o$ is an $R$-submodule of the $R$-module, the set $\Sigma$ is flat in the liner $X$, by Lemma~\ref{l:flat<=>R-submodule}. Then the inclusion $B\subseteq \Sigma$ implies $\overline B\subseteq\Sigma$. 
\end{proof}

\begin{lemma}\label{l:A+b=A+R(b-o)} Let $X$ be an $R$-module over a corps $R\ne\{0,1\}$. For every flat $A\subseteq X$ and points $o\in A$ and $b\in X\setminus A$, the flat hull $\overline{A\cup\{b\}}$ coincides with the set $A+R\cdot (b-o)$.
\end{lemma}

\begin{proof} By Lemma~\ref{l:flat<=>R-submodule}, the sets $A-o$ and $\overline{A\cup\{b\}}-o$ are $R$-submodules of the $R$-module $X$. Then $(A-o)\cup\{b-o\}\subseteq \overline{A\cup\{b\}}-o$ implies $A-o+R\cdot(b-o)\subseteq\overline{A\cup\{b\}}-o$ and hence $A+R\cdot(b-o)\subseteq\overline{A\cup\{b\}}$.
On the other hand, since $A-o+R\cdot(b-o)$ is an $R$-submodule of $X$, Lemma~\ref{l:flat<=>R-submodule} ensures that the set $A+R\cdot(b-o)$ is flat in $X$. Taking into account that $A\cup\{b\}=(A+0\cdot(b-o))\cup(o+1\cdot(b-o))\subseteq A+R\cdot(b-o)$, we conclude that $\overline{A\cup\{b\}}\subseteq A+R\cdot(b-o)$ and hence $\overline{A\cup\{b\}}=A+R\cdot(b-o)$.
\end{proof}

\begin{exercise} Find examples showing that Lemmas~\ref{l:flat<=>R-submodule} and \ref{l:A+b=A+R(b-o)} do not hold for $R$-modules over the field $R=\{0,1\}$.
\end{exercise}

\begin{theorem}\label{t:R-module=>Des-aff-reg} Every $R$-module over a corps $R$ is a Desarguesian affine  regular liner.
\end{theorem}

\begin{proof} Let $X$ be an $R$-module over a corps $R$. We consider $X$ as a liner endowed with the canonical line relation.

\begin{lemma} The liner $X$ is regular.
\end{lemma}

\begin{proof} If $R=\{0,1\}$, then $|X|_2=|R|=2$ and hence the liner $X$ is (strongly) regular. So, assume that $R\ne\{0,1\}$ and take any flat $A\subseteq X$ and points $o\in A$ and $b\in X\setminus A$. Given any point $z\in\overline{A\cup\{b\}}$, we have to find points $x\in A$ and $y\in \Aline ob$ such that $z\in\Aline xy$. 
By Lemma~\ref{l:A+b=A+R(b-o)}, $\overline{A\cup\{b\}}=A+R\cdot(b-o)$ and hence $z=a+r\cdot (b-o)$ for some $a\in A$ and $r\in R$. Choose any element $s\in R\setminus\{0,1\}$ and observe that the point $y\defeq o+s^{-1}\cdot r\cdot (b-o)$ belongs to the line $o+R\cdot (b-o)=\Aline ob$ and the point $x\defeq o+(1-s)^{-1}\cdot (a-o)\subseteq o+R\cdot(a-o)=\Aline oa$ belongs to the flat $A$. Observe that $$
\begin{aligned}
x+s\cdot(y-x)&=(1-s)\cdot x+s\cdot y=(1-s)\cdot(o+(1-s)^{-1}(a-o))+s\cdot(o+s^{-1}\cdot r\cdot(b-o))\\
&=(1-s)\cdot o+(a-o)+s\cdot o+r\cdot(b-o)=a+r\cdot(b-o)=z
\end{aligned}
$$ and hence $z=x+s\cdot(y-x)\in \Aline xy$, witnessing that the liner $X$ is regular.
\end{proof}

\begin{lemma} The liner $X$ is affine.
\end{lemma}

\begin{proof} Given any points $o,x,y\in X$ and $p\in \Aline xy\setminus\Aline ox$, we have to prove that the set $I\defeq\{u\in\Aline yo:\Aline up\cap\Aline ox=\varnothing\}$ is a singleton. If $p\in \Aline oy$, then the set $I$ coincides with the singleton $\{p\}$ and we are done. So, assume that $p\notin\Aline oy$ and hence $o\ne x$ and $oyx$ is a triangle in the liner $X$. It follows from $p\in\Aline xy\setminus(\Aline ox\cup\Aline oy)\subseteq(y+R\cdot(x-y))\setminus\{x\}$ that $p=y+r\cdot(x-y)$ for some $r\in R\setminus\{0,1\}$. Lemma~\ref{l:para-lines-in-R-module} implies that $I=\{y+r\cdot(o-y)\}$ is a singleton.
\end{proof}

\begin{lemma} The liner $X$ is Desarguesian.
\end{lemma}

\begin{proof} Given two centrally perspective triangles $abc,a'b'c'\in X^3$ with $\Aline ab\cap\Aline{a'}{b'}=\varnothing=\Aline bc\cap\Aline{b'}{c'}$, we need to prove that $\Aline ac\cap\Aline{a'}{c'}=\varnothing$.

Let $o\in\Aline a{a'}\cap\Aline b{b'}\cap\Aline c{c'}$ be a unique perspector of the triangles $abc,a'b'c'$. It follows from $a'\in \Aline oa$, $b'\in\Aline ob$, $c'\in \Aline oc$ that $a'=o+r\cdot(a-o)$, $b'=o+s\cdot(b-o)$, $c'=o+t\cdot(c-o)$ for some $r,s,t\in R\setminus\{0,1\}$. By Lemma~\ref{l:para-lines-in-R-module}, the equalities $\Aline ab\cap\Aline{a'}{b'}=\varnothing=\Aline bc\cap\Aline{b'}{c'}$ imply $r=s=t$. By Lemma~\ref{l:para-lines-in-R-module}, the equality $r=t\in R\setminus\{0,1\}$ implies $\Aline ac\cap\Aline {a'}{c'}=\varnothing$.
\end{proof}
\end{proof}

\begin{proposition}\label{p:IR_X=R} Let $X$ be an $R$-module over a corps $R\ne\{0,1\}$. If $\|X\|\ge 3$, then the corps $\IR_X$ of the Desarguesian affine space $X$ is isomorphic to the corps $R$.
\end{proposition}

\begin{proof} Let $o$ be the neutral element of the additive group of the $R$-module $X$. For every $r\in R$ and $x\in X$, consider the point $r_x\defeq r\cdot x\in\Aline ox$. The axioms of an $R$-module imply that for every elements $r,s\in R$ and point $x\in X$ we have $$r_{s_x}=r\cdot(s\cdot x)=(r\cdot s)\cdot x=rs_x.$$

  Since $\|X\|\ge 3$, there exist a point $e\in X\setminus\{o\}$. By Theorem~\ref{t:R-module=>Des-aff-reg}, the liner $X$ is a Desarguesian affine space.  By Theorems~\ref{t:Desargues<=>3Desargues} and \ref{t:triple-unique}, $\IR_X=\{\overvector{oxe}:x\in \Aline oe\}$. Consider the function $F:R\to\IR_X$ assigning to every element $r\in R$ the scalar $\overvector{or_ee}$. Theorem~\ref{t:triple-unique} ensures that the function $F$ is bijective. 

Let us check that $F$ is an isomorphism of the corps $R$ and $\IR_X$. Observe that $F(0)=\overvector{ooe}=0$ and $F(1)=\overvector{oee}=1$. 

\begin{claim} For every $r,s\in R$ we have $F(r\cdot s)=F(r)\cdot F(s)$.
\end{claim}

\begin{proof} If $r=0$, then $F(r)=\overvector{ooe}=0$ and $F(r\cdot s)=F(0)=0=F(0)\cdot F(s)=F(r)\cdot F(s)$, by Theorem~\ref{t:multiplication-scalar}.
If $s=0$, then $F(s)=0$ and $F(r\cdot s)=F(0)=0=F(r)\cdot 0=F(r)\cdot F(s)$, by Theorem~\ref{t:multiplication-scalar}.

So, assume that $r\ne 0\ne s$. Since $\|X\|\ge 3$, there exists a point $u\in X\setminus R\cdot e$.  By Lemma~\ref{l:para-lines-in-R-module}, $\Aline eu\parallel \Aline{r_e}{r_u}$ and $\Aline u{s_e}\parallel \Aline{r_u}{r_{s_e}}=\Aline{r_u}{rs_e}$. Then $$F(r\cdot s)=\overvector{o(rs)_ee}=\overvector{or_{s_e}e}=\overvector{or_{s_e}{s_e}}\cdot\overvector{os_ee}=\overvector{or_uu}\cdot F(s)=\overvector{or_ee}\cdot F(s)=F(r)\cdot F(s).$$ 

\begin{picture}(120,85)(-120,-10)

\put(0,0){\line(1,0){120}}
\put(0,0){\line(2,1){120}}
\put(120,0){\line(0,1){60}}
\put(60,0){\line(0,1){30}}
\put(90,0){\line(1,2){30}}
\put(45,0){\line(1,2){15}}

\put(0,0){\circle*{3}}
\put(-2,-8){$o$}
\put(45,0){\circle*{3}}
\put(40,-8){$rs_e$}
\put(60,0){\circle*{3}}
\put(57,-8){$r_e$}
\put(90,0){\circle*{3}}
\put(87,-8){$s_e$}
\put(120,0){\circle*{3}}
\put(118,-8){$e$}
\put(60,30){\circle*{3}}
\put(55,35){$r_u$}
\put(120,60){\circle*{3}}
\put(117,64){$u$}
\end{picture}
\end{proof}

\begin{claim} For every $r,s\in R$, we have $F(r+s)=F(r)+F(s)$.
\end{claim}

\begin{proof} If $r=0$, then $F(r+s)=F(0+s)=F(s)=0+F(s)=F(0)+F(s)=F(r)+F(s)$.

If $s=0$, then $F(r+s)=F(r+0)=F(r)=F(r)+F(0)=F(r)+F(s)$.

So, assume that $r\ne 0\ne s$, which implies $r_e=r\cdot e\ne o\ne s\cdot e=s_e$. Since $\|X\|\ge 3$, there exists a point $u\in X\setminus (R\cdot e)$. Consider the point $v\defeq u+s_e\ne u$ and observe that the lines $\Aline os_e=R\cdot e$ and $\Aline uv=u+R\cdot s_e=u+R\cdot e$ are parallel. 

\begin{picture}(200,65)(-100,-15)

\put(0,0){\line(1,0){120}}
\put(30,30){\line(1,0){30}}
\put(0,0){\line(1,1){30}}
\put(30,0){\line(1,1){30}}
\put(30,30){\line(1,-1){30}}
\put(60,30){\line(1,-1){30}}

\put(0,0){\circle*{3}}
\put(-2,-8){$o$}
\put(30,0){\circle*{3}}
\put(27,-8){$s_e$}
\put(60,0){\circle*{3}}
\put(57,-8){$r_e$}
\put(90,0){\circle*{3}}
\put(78,-10){$(r{+}s)_e$}
\put(120,0){\circle*{3}}
\put(118,-8){$e$}
\put(30,30){\circle*{3}}
\put(28,34){$u$}
\put(60,30){\circle*{3}}
\put(58,34){$v$}

\end{picture}

Also the lines $\Aline ou=o+R\cdot u$ and $\Aline {s_e}v=s_e+R\cdot(v-s_e)=s_e+R\cdot(u+s_e-s_e)=s_e+R\cdot u$ are parallel. Then $os_e\# uv$ and hence $\overvector{os_e}=\overvector{uv}$. On the other hand, the lines $\Aline uv=u+R\cdot e$ and $\Aline {r_e}{(r{+}s)_e}=r_e+R\cdot((r{+}s)_e-r_e)=r_e+R\cdot e$ are parallel and the lines $\Aline u{r_e}=u+R\cdot(r_e-u)$ and $\Aline v{(r{+}s)_e}=v+R\cdot((r{+}s)_e-(u+s_e))=v+R\cdot(r_e-u)$ are parallel, which implies $\overvector{uv}=\overvector{r_e(r{+}s)_e}$. Then 
$$\overvector{o(r{+}s)_e}=\overvector{or_e}+\overvector{r_e(r{+}s)_e}=\overvector{or_e}+\overvector{uv}=\overvector{or_e}+\overvector{os_e}.$$ The definition of addition of scalars in the Desarguesian liner $X$ ensures that
$$F(r+s)=\overvector{o(r{+}s)_ee}=\overvector{or_ee}+\overvector{os_ee}=F(r)+F(s).$$
\end{proof}
\end{proof}


\section{Isomorphisms of Desarguesian affine spaces}

Since vectors and scalars in an affine space are defined using the line relation of the affine space, those notions are preserved by liner isomorphisms of affine spaces. More precisely, every isomorphism $F:X\to Y$ between affine spaces $X,Y$ induces  bijective maps $\vec F:\overvector X\to \overvector Y$ and $\dddot F:\IR_X\to\IR_Y$ assigning to every vector $\vec{\boldsymbol v}\in\overvector X$ the vector $\vec F(\vec{\boldsymbol v})\defeq\{Fxy:xy\in\vec{\boldsymbol v}\}$ and to every scalar $\sigma\in \IR_X$ the scalar $\dddot F(\sigma)\defeq\{Fxyz:xyz\in\sigma\}$. The bijection $\vec F:\overvector X\to\overvector Y$ preserves the zero vectors and the inversion of vectors, and the bijection $\dddot F:\IR_X\to\IR_Y$ is an isomorphism of the corps $\IR_X$ and $\IR_Y$. 

The definitions of the operations of addition of vectors, addition and muliplication of scalars, and multiplication of scalars and vectors imply the following proposition.

\begin{proposition}\label{p:iso-para-algebra} For every isomorphism $F:X\to Y$ between two affine spaces $X,Y$, the bijections $\vec F:\overvector X\to\overvector Y$ and $\dddot F:\IR_X\to\IR_Y$ have the following properties:
\begin{enumerate}
\item $\forall s,t\in \IR_X\;\;\dddot F(s\cdot t)=\dddot F(s)\cdot \dddot F(t)$;
\item $\forall s,t\in\IR_X\;\;\dddot F(s+t)=\dddot F(s)+\dddot F(t)$;
\item $\forall \vec{\boldsymbol x},\vec{\boldsymbol y}\in\overvector X\;\;\vec F(\vec{\boldsymbol x}+
\vec{\boldsymbol y})=\vec F(\vec {\boldsymbol x})+\vec F(\vec{\boldsymbol y})$;
\item $\forall s\in\IR_X\;\forall \vec{\boldsymbol v}\in\overvector X\;\;\vec F(s{\cdot}\vec {\boldsymbol v})=\dddot F(s){\cdot}\vec F(\vec{\boldsymbol v})$;
\item $\forall o\in X\;\forall \vec{\boldsymbol v}\in\overvector X\;\;F(o+\vec{\boldsymbol v})=F(o)+ \vec F(\vec {\boldsymbol v})$.
\end{enumerate} 
\end{proposition}

\begin{exercise} Write down a proof of Proposition~\ref{p:iso-para-algebra}.
\end{exercise}

For Desarguesian affine spaces, Proposition~\ref{p:iso-para-algebra} turns into a characterization of isomorphisms between Desarguesian affine spaces.

\begin{theorem}\label{t:Disomorphism<=>} A  function $F:X\to Y$ between two Desarguesian affine liners $X,Y$ is an isomorphism if and only if there exist bijective maps $\vec F:\overvector X\to\overvector Y$ and $\dddot F:\IR_X\to\IR_Y$  such that the following conditions are satisfied:
\begin{enumerate}
\item $\forall s,t\in \IR_X\;\;\dddot F(s\cdot t)=\dddot F(s)\cdot \dddot F(t)$;
\item $\forall s,t\in\IR_X\;\;\dddot F(s+t)=\dddot F(s)+\dddot F(t)$;
\item $\forall \vec{\boldsymbol x},\vec{\boldsymbol y}\in\overvector X\;\;\vec F(\vec{\boldsymbol x}+
\vec{\boldsymbol y})=\vec F(\vec {\boldsymbol x})+\vec F(\vec{\boldsymbol y})$;
\item $\forall s\in\IR_X\;\forall \vec{\boldsymbol v}\in\overvector X\;\;\vec F(s{\cdot}\vec {\boldsymbol v})=\dddot F(s){\cdot}\vec F(\vec{\boldsymbol v})$;
\item $\exists o\in X\;\forall \vec{\boldsymbol v}\in\overvector X\;\;F(o+\vec{\boldsymbol v})=F(o)+ \vec F(\vec {\boldsymbol v})$.
\end{enumerate}
\end{theorem}

\begin{proof} The ``only if'' part follows from Proposition~\ref{p:iso-para-algebra}. To prove the ``if'' part, assume that for a function $F:X\to Y$ there exist bijective maps $\vec F:\overvector X\to\overvector Y$ and $\dddot F:\IR_X\to\IR_Y$ satisfying the conditions (1)--(5). By Theorem~\ref{t:ADA=>AMA}, the Desarguesian affine space $X$ is Thalesian and by Corollary~\ref{c:Thalesian<=>vector=funvector}, the set $\overvector X$ of functional vectors in $X$ coincides with the set $X^2_{\#}$ of all vectors in $X$.

The conditions (1) and (2) ensure that the bijection $\dddot F:\IR_X\to\IR_Y$ is an isomorphism of the corps $\IR_X$ and $\IR_Y$. The condition (3) implies that the bijection $\vec F:\overvector X\to\overvector Y$ is an isomorphism of the Abelian groups $(\overvector X,+)$ and $(\overvector Y,+)$.  By the condition (5),  there exists a point $o\in X$ such that $F(o+\vec{\boldsymbol v})=F(o)+\vec F(\vec{\boldsymbol v})$ for every vector $\vec{\boldsymbol v}\in\overvector X$. Consider the point $o'\defeq F(o)\in Y$.

\begin{claim}\label{cl:ox->o'x'} For every point $x\in X$ and its image $x'\defeq F(x)\in Y$ we have
$\overvector{o'x'}=\vec F(\overvector{ox})$.
\end{claim}

\begin{proof} Corollary~\ref{c:vector-action}(1) and the condition (5) ensure that
$$o'+\overvector{o'x'}=x'=F(x)=F(o+\overvector{ox})=F(o)+\vec F(\overvector{ox})=o'+\vec F(\overvector{ox})$$and hence $\overvector{o'x'}=\vec F(\overvector{ox})$.
\end{proof}

\begin{claim}\label{cl:xy->x'y'} For any pair $xy\in X^2$ and its image $x'y'\defeq Fxy\in Y^2$ we have
$$\overvector{x'y'}=\vec F(\overvector{xy}).$$
\end{claim}

\begin{proof} By Claim~\ref{cl:ox->o'x'}, $\overvector{o'x'}=\vec F(\overvector{ox})$ and $\overvector{o'y'}=\vec F(\overvector{oy})$. By Theorem~\ref{t:vector-addition} and the condition (3),
$$\overvector{o'x'}+\overvector{x'y'}=\overvector{o'y'}=\vec F(\overvector{oy})=\vec F(\overvector{ox}+\overvector{xy})=\vec F(\overvector{ox})+\vec F(\overvector{xy})=\overvector{o'x'}+\vec F(\overvector{xy})$$ and hence 
$\overvector{x'y'}=\vec F(\overvector{xy})$, see Theorem~\ref{t:vector-addition}.
\end{proof}

 In the following three claims we shall prove that the map $F$ is an isomorphism of the liners $X,Y$.

\begin{claim} The map $F:X\to Y$ is injective.
\end{claim}

\begin{proof} Given any distinct points $x,y\in X$, we should prove that their images $x'\defeq F(x)$ and $y'\defeq F(y)$ are distinct. Since every vector in the Desarguesian affine space $X$ is functional, the inequality $x\ne y$ implies  $\overvector{ox}\ne\overvector{oy}$. Applying Claim~\ref{cl:ox->o'x'}, we conclude that
$$\overvector{o'x'}=\vec F(\overvector{ox})\ne\vec F(\overvector{oy})=\overvector{o'y'}$$and hence $x'\ne y'$.
\end{proof} 

\begin{claim} The map $F:X\to Y$ is surjective.
\end{claim}

\begin{proof} Given any point $y\in X$, we should find a point $x\in X$ such that $y=F(x)$. Since the bijection $\vec F:\overvector X\to\overvector Y$ is surective, there exists a vector $\vec{\boldsymbol v}$ such that $\vec F(\vec{\boldsymbol v})=\overvector{o'y}$. By Corollary~\ref{c:Thalesian<=>vector=funvector}, there exists a unique point $x\in X$ such that $\overvector{ox}=\vec{\boldsymbol v}$. Then $x=o+\overvector{ox}=o+\vec{\boldsymbol v}$ and hence
$$F(x)=F(o+\vec{\boldsymbol v})=F(o)+\vec F(\vec{\boldsymbol v})=o'+\overvector{o'y}=y.$$
\end{proof}

\begin{claim} The bijective map $F:X\to Y$ is a liner isomorphism.
\end{claim}

\begin{proof} Let $\Af_X,\Af_Y$ be the line relations of the liners $X,Y$, respectively.   Given any triple $xyz\in X^3$, we should prove that $xyz\in\Af_X$ if and only if $Fxyz\in\Af_Y$. Consider the triple $x'y'z'\defeq Fxyz$. Claim~\ref{cl:xy->x'y'} ensures that $\overvector{x'y'}=\vec F(\overvector{xy})$ and $\overvector{x'z'}=\vec F(\overvector{xz})$.  

First we show that $xyz\in\Af_X$ implies $x'y'z'\in\Af_Y$. If $x=z$, then $xyz\in\Af_X$ implies $x=y=z$ and hence $x'=y'=z'$ and $x'y'z'\in\Af_Y$, by the Identity Axiom {(\sf IL)} in Definition~\ref{d:liner}. So, assume that $x\ne z$. In this case, the bijectivity of $F$ ensures that $x'\ne z'$. Since $xyz\in\Af_X$, by  Theorem~\ref{t:paraD<=>Desargues}, there exists a scalar $s\in\IR_X$ such that $\overvector{xy}=s\cdot\overvector{xz}$. Applying the condition (4), we conclude that
$$\overvector{x'y'}=\vec F(\overvector{xy})=\vec F(s\cdot\overvector{xz})=\dddot F(s)\cdot\vec F(\overvector{xz})=\dddot F(s)\cdot\overvector{x'z'}.$$ By Proposition~\ref{p:triple-exists}, there exists a point $u\in \Aline{x'}{z'}\subseteq Y$ such that $x'uz'\in \dddot F(s)$. Theorem~\ref{t:scalar-by-vector}, $$\overvector{x'y'}=\dddot F(s)\cdot\overvector{x'z'}=\overvector{x'uz'}\cdot\overvector{x'z'}=\overvector{x'u}.$$ By Corollary~\ref{c:Thalesian<=>vector=funvector}, the equality $\overvector{x'y'}=\overvector{x'u}$ implies  $y'=u\in\Aline{x'}{z'}$ and hence $x'y'z'\in  \Af_Y$.

 By analogy we can prove that $x'y'z'\in\Af_Y$ implies $xyz\in\Af_X$.
\end{proof}
\end{proof}

\begin{proposition}\label{p:hull=sum} Let $M$ be a set in a Desarguesian affine space $X$, $o$ be a point of $M$ and $B\defeq M\setminus\{o\}$. The flat hull $\overline M$ of the set $M$ in $X$ is equal to the set
$$\textstyle\Sigma\defeq \{x\in X:\exists c\in\IR_X^{\oplus B}\;\;\overvector{ox}=\sum_{b\in B}c(b)\cdot\overvector{ob}\}.$$
\end{proposition}

\begin{proof} By Theorem~\ref{t:paraD=>RX-module}, the space of vector $\overvector X$ is an $\IR_X$-module. Consider the set $$\overvector{\Sigma}\defeq\{\overvector{ox}:x\in\Sigma\}=\textstyle\big\{\sum_{b\in B}c(b)\cdot\overvector{ob}:c\in\IR_X^{\oplus B}\big\}$$and observe that it is an $\IR_X$-submodule of the $\IR_X$-module $\overvector{X}$. Moreover, $\overvector{\Sigma}$ is the smallest $\IR_X$-submodule of $\overvector X$ that contains the set $\{\overvector{ob}:b\in B\}\subseteq \overvector{oM}$. By Proposition~\ref{p:flat=>R-submodule}, the set $$\overvector{oM}\defeq\{\overvector{ox}:x\in\overline M\}$$ is an $\IR_X$-submodule of the $\IR_X$-module $\overvector X$. Then $\{\overvector{ob}:b\in B\}\subseteq\overvector{\Sigma}$ implies $\overvector{\Sigma}\subseteq \overvector{oM}$. By the uniquness part of Corollary~\ref{c:Thalesian<=>vector=funvector}, $\overvector{\Sigma}\subseteq\overvector{oM}$ implies $\Sigma\subseteq\overline M$. Since $\overvector{\Sigma}$ is an $\IR_X$-submodule of $\overvector X$, the set $\Sigma$ is flat in $X$, by Theorem~\ref{t:flat<=>R-submodule}. Since $M\subseteq\Sigma$, the flat $\Sigma$ contains the flat hull $\overline M$ of the set $M$. Therefore, $\overline M=\Sigma$.
\end{proof}

\begin{theorem}\label{t:c-isomorphism} Let $M$ be a maximal independent set in a Desarguesian affine space $X$, $o$ be a point in $M$ and $B\defeq M\setminus\{o\}$. For every point $x\in X$, there exists a unique function $c_x\in \IR_X^{\oplus B}$ such that $\overvector{ox}=\sum_{b\in B}c_x(b)\cdot \overvector{ob}$. The map $c:X\to \IR_X^{\oplus B}$, $c:x\mapsto c_x$, is an isomorphism of liners.
\end{theorem}

\begin{proof} The maximality of $M$ and Proposition~\ref{p:add-point-to-independent} imply $\overline M=X$. By Proposition~\ref{p:hull=sum}, the set $X=\overline M$ is equal to the set 
$$\Sigma\defeq\big\{x\in X:\exists c\in \IR_X^{\oplus B}\;\;\big(\overvector{ox}=\sum_{b\in B}c(b)\cdot\overvector{ob}\big)\big\}.$$

Therefore, for every point $x\in X$ there exists a finitely supported function $c_x\in \IR_X^{\oplus B}$ such that $\overvector{ox}=\sum_{b\in B}c_x(b)\cdot \overvector{ob}$. We claim that the function $c_x$ is unique. Assume that $c_x'\in\IR_X^{\oplus B}$ is another finitely supported function such that $\overvector{ox}=\sum_{b\in B}c_x'(b)\cdot \overvector{ob}$. Then 
$$\vec{\mathbf 0}=\sum_{b\in B}(c_x(b)-c_x'(b))\cdot\overvector{ob}.$$
Assuming that $c_x\ne c_x'$, find a point $p\in B$ such that $c_x(p)\ne c_x'(p)$ and consider the set $B'\defeq B\setminus\{p\}$. Since $\IR_X$ is a corps and $c'_x(p)-c_x(p)\ne 0$, there exists a scalar $s\in \IR_X$ such that $s\cdot(c'_x(p)-c_x(p))=1$. 
Then $\vec{\mathbf 0}=s\cdot\vec{\mathbf 0}=\sum_{b\in B}s\cdot(c_x(b)-c_x'(b))\cdot\overvector{ob}$ implies
$$\overvector{op}=1\cdot \overvector{op}=s\cdot(c'_x(p)-c_x(p))\cdot\overvector{op}=\sum_{b\in B\setminus\{p\}}s\cdot(c_x(b)-c_x'(b))\cdot\overvector{ob}.$$Applying Proposition~\ref{p:hull=sum}, we conclude that $p\in \overline{M\setminus\{p\}}$, which contradicts the independence of $M$. This contradiction shows that the function $c_x$ is unique.   

Then the map $c_*:X\to\IR_X^{\oplus B}$, $c_*:x\mapsto c_x$, is well-defined. It is easy to see that this map is surjective and injective. To see that the map $c_*$ is an isomorphism of the liners $X$ and $\IR_X^{\oplus B}$, take any triple $xyz\in X^3$ with $y\in\Aline xz$. By Theorem~\ref{t:paraD<=>Desargues}, there exists a scalar $s\in\IR_X$ such that $\overvector{xy}=s\cdot\overvector{xz}$. Then 
$$
\begin{aligned}
\sum_{b\in B}c_y(b)\cdot b&=\overvector{oy}=\overvector{ox}+\overvector{xy}=\overvector{ox}+s\cdot\overvector{xz}=1\cdot\overvector{ox}+s\cdot(\overvector{oz}-\overvector{ox})=(1-s)\cdot\overvector{ox}+s\cdot\overvector{oz}\\
&=(1-s)\cdot\sum_{b\in B}c_x(b)\overvector{ob}+s\cdot\sum_{b\in B}c_z(b)\overvector{ob}=\sum_{b\in B}(1-s)\cdot c_x(b)+s\cdot c_z(b))\cdot\overvector{ob}
\end{aligned}
$$ and hence
$$c_y=(1-s)\cdot c_x+s\cdot c_y\in \Aline {c_x}{x_y}\subseteq \IR_X^{\oplus B},$$
witnessing that the points $c_x,c_y,c_z$ of the liner $\IR_X^{\oplus B}$ are collinear and the function $c_*:X\to\IR_X^{\oplus B}$ is a liner isomorphism, according to Theorem~\ref{t:liner-isomorphism<=>}.
\end{proof}

Theorem~\ref{t:c-isomorphism} and Corollary~\ref{c:Max=dim} imply the following important corollary.

\begin{corollary}\label{c:Desargues=>RXdim} Every Desarguesian affine space $X$ is isomorphic to the liner $\IR_X^{\oplus\dim(X)}$.
\end{corollary}

\begin{theorem}\label{t:aff-Des<=>} An affine space $X$ is Desarguesian if and only if $X$ is isomorphic to an $R$-module over a corps $R\ne\{0,1\}$.
\end{theorem}

\begin{proof} The ``if'' part follows from Theorem~\ref{t:R-module=>Des-aff-reg}.
To prove the ``only if'' part, assume that the affine space $X$ is Desarguesian. By Corollary~\ref{c:Desargues=>RXdim}, the Deasrguesian affine space $X$ is isomorphic to the $\IR_X$-module $\IR_X^{\oplus\dim(X)}$ over the corps $\IR_X$. It remains to prove that $\IR_X\ne\{0,1\}$. By Theorem~\ref{t:paraD=>RX-module}, the space of vectors $\overvector X$ is an $\IR_X$-module over the corps $\IR_X$. By Theorem~\ref{t:paraD<=>Desargues}, the liner relation $\Af$ on $X$ is equal to $$\{xyz\in X^3:\overvector{xy}\in R_X\cdot\overvector{xz}\}.$$ This implies that for any distinct points $x,z\in X$, the line $\Aline xz$ coincides with the set $x+\IR_X\cdot\overvector xz$. Since the affine space $X$ is $3$-long,
$3\le|\Aline xz|=|\IR_X|$ and hence $\IR_X\ne\{0,1\}$. 
\end{proof}

\begin{theorem}\label{t:extension-iso} Let $X,Y$ be two Desarguesian affine spaces and $I:\IR_X\to \IR_Y$ be an isomorphism between their corps. Let $M_X,M_Y$ be any maximal independent sets in the liners $X,Y$, respectively. Every bijection $f:M_X\to M_Y$ extends to a unique liner isomorphism $F:X\to Y$ such that $\dddot F=I$.
\end{theorem}

\begin{proof} Fix any point $o\in M_X$ and consider the set $B\defeq M_X\setminus \{o\}$. Let $o'\defeq f(o_X)$ and $B'\defeq f[B]$. For every point $b\in B$, let $b'$ denote the point $f(b)\in B'$. By Theorem~\ref{t:c-isomorphism}, for every point $x\in X$, there exists a unique function $c_x\in \IR_X^{\oplus B}$ such that $\overvector{ox}=\sum_{b\in B}s_x(b)\cdot\overvector{ob}$. Define a function $F:X\to Y$ assigning to every point $x\in X$ the unique point $y\in Y$ such that
$$\overvector{o'y}=\sum_{b\in B} I(c_x(b))\cdot \overvector{o'b'}.$$
It is easy to show that $F:X\to Y$ is a required liner isomorphism such that $F{\restriction}_B=f$ and $\dddot F=I$.
\end{proof}

Theorem~\ref{t:extension-iso} and Corollary~\ref{c:Max=dim} imply the following isomorphic classification of Desarguesian affine spaces.

\begin{corollary}
Two Desarguesian affine spaces $X,Y$ are isomorphic if and only if $\|X\|=\|Y\|$ and the corps $\IR_X$ and $\IR_Y$ are isomorphic.
\end{corollary}

\section{Projective spaces of modules over corps}

In this section we define and study an important example of a Desarguesian projective liner of algebraic origin.

Let $X$ be an $R$-module over a corps $R$ and $\mathbf 0$ be the neutral element of the additive group of the $R$-module $X$. An \index{$R$-line}\defterm{$R$-line} in the $R$-module $X$ is any $R$-submodule $L\subseteq X$ of form $R\cdot x$ for some nonzero element $x\in X$.

The \index{projective space}\index[note]{$\mathbb PX$}\defterm{projective space} $\mathbb PX$ of the $R$-module $X$ is the set of $R$-lines $\{R\cdot b:b\in X\setminus\{\mathbf 0\}\}$ endowed with the \defterm{canonical line relation}
$$\Af\defeq\{ABC\in (\mathbb P X)^3:B\subseteq A+C\},$$
where $A+C\defeq\{a+c:a\in A,\;c\in C\}$.
Let us show that the ternary relation $\Af$ is indeed a line relation on the projective space $\mathbb PX$.

\begin{lemma} For every $R$-module $X$ over a corps $R$, the ternary relation  
 $$\Af\defeq\{ABC\in (\mathbb P X)^3:B\subseteq A+C\}$$
in a line relation on $\mathbb PX$.
\end{lemma}

\begin{proof} We have to check that the ternary relation $\Af$ satisfies the axioms {\sf(IL)}, {\sf(RL)}, {\sf(EL)} in Definition~\ref{d:liner}.  
Consider the $R$-module $X$ as a liner endowed with its canonical liner relation. By Theorem~\ref{t:R-module=>Des-aff-reg}, the liner $X$ is affine and regular. 
\smallskip

To check the axiom {\sf(IL)}, take any $R$-lines $A,B\in \mathbb PX$ with $ABA\in \Af$. Choose any nonzero points $a\in A$ and $b\in B$ and observe that $B\subseteq A+A=R\cdot a+R\cdot a=(R+R)\cdot a=R\cdot a=A$ and hence $b=r\cdot a$ for some $r\in R\setminus\{0\}$. Since $R$ is a corps, $b=r\cdot a$ implies $a=1\cdot a=r^{-1}\cdot r\cdot a=r^{-1}\cdot b$ and hence $A=R\cdot a=R\cdot r^{-1}\cdot b\subseteq R\cdot b=B$ and finally $A=B$.
\smallskip

To check the axiom {\sf(RL)}, take any $R$-lines $A,B\in \mathbb PX$ and observe that $A\cup B\subseteq A+B$ and hence $\{AAB,ABB\}\subseteq\Af$.
\smallskip

To check the Exchange Axiom {\sf(EL)}, take any $R$-lines $A,B,U,V\in \mathbb PX$ such that $AUB,AVB\in\Af$ and $U\ne V$. We have to prove that $UAV,UBV\in\Af$. The definition of the ternary relation $\Af$ ensures that $U\cup V\subseteq A+B$. Observe that the sets $A+B$ and $U+V$ are $R$-submodules of the $R$-module $X$. Then $U\cup V\subseteq A+B$ implies $U+V\subseteq A+B$. By Propositions~\ref{t:flat<=>R-submodule} and \ref{p:flathull-in-R-module}, the $R$-submodules $U+V\subseteq A+B$ are flat in the liner $X$ such that $U+V=\overline{U\cup V}$ and $A+B=\overline{A\cup B}$. Since $U$ and $V$ are two concurrent lines in $X$, the flat $U+V=\overline{U\cup V}$ is a plane and so is the flat $A+B\subseteq{A\cup B}$ containing the plane $U+V$. By Proposition~\ref{p:k-regular<=>2ex}, the regular liner $X$ is ranked and hence $U+V=A+B$ and $A\cup B\subseteq U+V$, which means that $UAV,UBV\in\Af$.
\end{proof}

From now on, we consider projective spaces over $R$-modules as liners endowed with their canonical line relations.

\begin{lemma}\label{l:PX-3long} The projective space $\IP X$ of any $R$-module $X$ over a corps $R$ is a $3$-long liner.
\end{lemma}

\begin{proof} Take any line $\mathcal A\subseteq \IP X$ and choose any distinct $R$-lines $A,B\in\A$. Next, choose any nonzero elements $a\in A$ and $b\in B$. It is easy to see that $a+b\in (A+B)\setminus(A\cup B)$ and hence the $R$-line $C\defeq R\cdot(a+b)$ is an element of the line $\Aline AB=\mathcal A$ such that $A\ne C\ne B$. Then $|\mathcal A|\ge|\{A,B,C\}|=3$, which means that the liner $\IP X$ is $3$-long.
\end{proof}

\begin{proposition}\label{p:PX-Steiner} The projective space $\mathbb PX$ of any $R$-module $X$ over a two-element field $R=\{0,1\}$ is a Steiner projective liner.
\end{proposition}

\begin{proof} To see that the liner $\mathbb PX$ is Steinter, take any distinct $R$-lines $A,B\in\mathbb PX$ and fix unique nonzero points $a\in A$ and $b\in B$. It follows that $A+B=\{\mathbf 0,a\}+\{\mathbf 0,b\}=\{\mathbf 0,a,b,a+b\}$ and hence $\Aline AB=\{R\cdot a,R\cdot b,R\cdot (a+b)\}$ and $|\Aline AB|=3$, witnessing that the liner $\mathbb PX$ is Steiner. 

Next, take any triangle $ABC$ in the liner $\mathbb PX$. Let $a\in A\setminus\{\mathbf 0\}$, $b\in B\setminus\{\mathbf 0\}$, $c\in C\setminus\{\mathbf 0\}$ be unique nonzero points in the $R$-lines $A,B,C$. Since $A,B,C$ is a triangle, the $R$-subliner $A+B+C$ contains eight points $\{\mathbf 0,a,b,c,a+b,b+c,a+c,a+b+c\}$ which determine seven $R$-lines $R\cdot x$ where $x\in(A+B+C)\setminus\{\mathbf 0\}$ of the plane $\overline{\{A,B,C\}}$ in the liner $\mathbb PX$. Then $|\overline{\{A,B,C\}}|=7$ and by Theorem~\ref{t:Steiner<=>projective}, the Steiner liner $\mathbb PX$ is projective.
\end{proof}

\begin{proposition}\label{p:PX-projective} The projective space $\mathbb PX$ of any $R$-module $X$ over a corps $R$ is a projective liner.
\end{proposition}

\begin{proof}  If $R=\{0,1\}$, then the liner $\mathbb PX$ is projective, by Proposition~\ref{p:PX-Steiner}. So, assume that $R\ne\{0,1\}$.

Given any $R$-lines $A,B,C\in \mathbb PX$ and $R$-lines $P\in \Aline BC$ and $Q\in \Aline AB\setminus\{P\}$, we should prove that $\Aline PQ\cap\Aline AC\ne\varnothing$ in $\mathbb PX$. If $P\in \Aline AB$, then $P,Q\in\Aline AB$ and $P\ne Q$ implies $\Aline PQ=\Aline AB$ and $A\in \Aline AB\cap\Aline AC=\Aline PQ\cap \Aline AC$. 
So, assume that $P\notin\Aline AB$, which implies $A\ne C$. Then $P+Q$ and $A+C$ are $R$-submodules in the $R$-module $X$ and planes in the liner $X$, according to Lemma~\ref{l:flat<=>R-submodule}. Observe that $P+Q\subseteq (B+C)+(A+B)=A+B+C$ and hence  $\|(P+Q)\cup(A+C)\|\le\|A+B+C\|\le 4$. By Theorems~\ref{t:R-module=>Des-aff-reg} and \ref{t:w-modular<=>}, the liner $X$ is regular and weakly modular. Since $\mathbf 0\in(P+Q)\cap(A+C)$, the weak modularity of the liner $X$ ensures that $$\|(P+Q)\cap(A+C)\|=\|P+Q\|+\|A+C\|-\|A+B+C\|\ge 3+3-4=2.$$ Consequently, $(P+Q)\cap(A+C)\ne\{\mathbf 0\}$ and there exists an $R$-line $L\subseteq(P+Q)\cap(A+C)$. This $R$-line $L$ belongs to the intersection of the lines $\Aline PQ$ and $\Aline AC$ in the liner $\mathbb PX$, witnessing that the liner $\mathbb PX$ is projective.
\end{proof}

\begin{proposition} Let $Y$ be an $R$-module over a corps $R$. For every $R$-submodule $X\subseteq Y$, the projective space $\mathbb PX$ is a flat subliner of the liner $\mathbb PY$. If $X\ne Y$, then the projective liner $\mathbb PY$ is a projective completion of the Proclus liner $\mathbb PY\setminus\mathbb PX$.
\end{proposition}

\begin{proof} Since $X$ is an $R$-submodule of the $R$-module $Y$, every $R$-line in $X$ is an $R$-line in $Y$ and hence $\mathbb PX\subseteq\mathbb PY$. Moreover, for every $R$-lines $A,B\in\mathbb PX$ their sum $A+B$ in $X$ is an $R$-submodule of $Y$, which implies that $\mathbb PX$ is a subliner of the liner $\mathbb PY$. If $X\ne Y$, then the subliner $Z\defeq \mathbb PY\setminus \mathbb PX$ is not empty and the complement $\mathbb PY\setminus Z=\mathbb PX$ is a proper flat in $\mathbb Y$. By Proposition~\ref{p:PX-projective}, the liner $\mathbb PY$ is projective. By Theorem~\ref{p:projective-minus-proflat}, the liner $Z$ is Proclus and by Definition~\ref{d:procompletion}, the projective liner $\mathbb PY$ is a projective completion of the Proclus liner $Z=\mathbb PY\setminus\mathbb PX$.
\end{proof}

\begin{lemma}\label{l:flat-in-PX<=>R-submodule} Let $X$ be an $R$-module over a corps $R$ and $\mathbb PX$ be its projective space. A subset $\mathcal A\subseteq \mathbb PX$ is flat in the liner $\mathbb PX$ if and only if $\bigcup \mathcal A$ is an $R$-submodule of the $R$-module $X$.
\end{lemma}

\begin{proof} Assume that a set $\A\subseteq\mathbb PX$ is flat in the liner $\mathbb PX$. To show that the set $\bigcup\A$ is an $R$-submodule in $X$, we have to check that for every elements $a,b\in\bigcup\A$ and $r\in R$, we have $a+b\in\bigcup\A$ and $r\cdot a\in \bigcup A$. Since $a,b\in\bigcup\A$, there exist $R$-lines $A,B\in\A$ such that $a\in A$ and $b\in B$. Then $a+b\in A+B$. If $a+b$ coincides with the neutral element $\mathbf 0$, then $a+b=\mathbf 0\in A\subseteq\bigcup\A$. If $a+b\ne\mathbf 0$, then $C\defeq R\cdot(a+b)$ is an $R$-line in $X$ such that $C\subseteq A+B$ and hence $C\in\Aline AB\subseteq \A\subseteq\mathbb PX$ and $a+b\in C\subseteq\A$. Also $r\cdot a\in A\subseteq\bigcup\A$, witnessing that $\bigcup\A$ is an $R$-submodule of $X$.

Next, assuming that $\bigcup\A$ is an $R$-submodule of $X$, we shall prove that $\A$ is a flat subset of the liner $\mathbb PX$. Given any $R$-lines $A,B\in\A$, we should prove that any $R$-line $C\in\Aline AB\subseteq\mathbb PX$ is an element of $\A$. Fix any nonzero element $c\in C$. It follows from $C\in\Aline AB$ that $c\in C\subseteq A+B\subseteq \bigcup \A+\bigcup \A=\bigcup \A$ and hence $c\in C'$ for some $R$-line $C'\in\A$. Since $R$ is a corps, $c\in C\cap C'\setminus\{\mathbf 0\}$ implies $C=R\cdot c=C'\in\A$.
\end{proof}

Now we investigate the relation between the $R$-independence and independence in  $R$-modules and their projective spaces.

\begin{definition} Let $X$ be an $R$-module over a corps $R$. A subset $B\subseteq X$ is called \index{$R$-independent set}\defterm{$R$-independent} in the $R$-module $X$ if every finitely supported function $f:B\to R$ with $\mathbf 0=\sum_{b\in B}f(b)\cdot b$ has empty support $\supp(f)$ and hence is the zero function.
\end{definition} 

\begin{lemma}\label{l:independent<=>R-independent} A subset $B$ of an $R$-module $X$ over a corps $R$ is $R$-independent if and only if $\mathbf 0\notin B$, $R\cdot a\ne  R\cdot b$ for any distinct elements $a,b\in B$, and the set $\mathbb PB\defeq \{R\cdot b:b\in B\}$ is independent in the liner $\mathbb PX$.
\end{lemma}

\begin{proof} To prove the ``only if'' part, assume that a subset $B\subseteq X$ of the $R$-module $X$ is $R$-independent. Then $\mathbf 0\notin B$ and for every distinct elements $a,b\in B$ we have $R\cdot a\ne R\cdot b$. Assuming that the set $\mathbb PB\defeq \{R\cdot b:b\in B\}$ is not independent in the liner $\mathbb PX$, we could find a point $p\in B$ such that $R\cdot p\in\overline{\mathbb PB_p}$ where $B_p\defeq B\setminus\{p\}$ and $\mathbb PB_p=\{R\cdot b:b\in B_p\}=\mathbb PB\setminus\{R\cdot p\}$. Consider the set $\Sigma\defeq\big\{\sum_{b\in B_p}f(b)\cdot b:f\in R^{\oplus B_p}\big\}$ and observe that $\Sigma$  is an $R$-submodule of the $R$-module $X$. By Lemma~\ref{l:flat-in-PX<=>R-submodule}, the set $\mathbb P\Sigma=\{R\cdot x:x\in \Sigma\setminus\{\mathbf 0\}\}$ is flat in the liner $\mathbb PX$. Since the set $\mathbb P\Sigma$ contains the set $\mathbb PB_p$, it also contains its flat hull in the liner $\mathcal PX$ and hence $R\cdot p\in \overline{\mathbb PB_p}\subseteq\mathbb P\Sigma$ and $p\in\bigcup\mathbb P\Sigma=\Sigma$, which contradicts the $R$-independence of the set $B$.  This contradiction shows that the set $\mathbb P\cdot B$ is independent in the liner $\mathbb PX$.
\smallskip

To prove the ``if'' part, assume that $\mathbf 0\notin B$, $R\cdot a\ne R\cdot b$ for every distinct elements $a,b\in B$, and the set $\mathbb PB\defeq\{R\cdot b:b\in B\}$ is independent in the liner $\mathbb PX$. Assuming that the set $B$ is not $R$-independent, we can find a nonzero finitely supported function $f\in R^{\oplus B}$ such that $\sum_{b\in B}f(b)\cdot b=\mathbf 0$.
Since $f$ is not zero, there exists $p\in  B$ such that $f(p)\ne 0$. Replacing the function $f(p)$ by $-f(p)^{-1}\cdot f$, we can assume that $f(p)=-1$. Then $p=1\cdot p=\sum_{b\in B_p}f(b)\cdot b$, where $B_p\defeq B\setminus\{p\}$. Since $R\cdot a\ne R\cdot b$ for any distinct points $a,b\in B$, the set $\mathbb P B_b\defeq\{R\cdot x:x\in B_p\}$ is a proper subset of the set $\mathbb PB$. Let $\mathcal H$ be the flat hull of the set $\mathbb PB_b$ in the liner $\mathbb PX$. By Lemma~\ref{l:flat-in-PX<=>R-submodule}, the union $\bigcup\mathcal H$ is an $R$-submodule of the $R$-module $X$. Then $p=\sum_{p\in B_p}f(b)\cdot b\subseteq\bigcup\mathcal H$ and hence $R\cdot p\in \mathcal H=\overline{\mathbb PB\setminus\{R\cdot p\}}$, which contradicts the independence of the set $\mathbb PB$ in the liner $\mathbb PX$. This contradiction shows that the set $B$ is $R$-indepenent in the $R$-module $X$.
\end{proof}

\begin{definition} A set $B$ in an $R$-module $X$ over a corps $X$ is called an \index{$R$-base}\defterm{$R$-base} for $X$ if for every element $x\in X$ there exists a unique finitely supported function $f\in R^{\oplus B}$ such that $x=\sum_{b\in B}f(b)\cdot b$.
\end{definition}

\begin{proposition}\label{p:R-base<=>} For a set $B$ in an $R$-module $X$ over a corps $R$, the following conditions are equivalent:
\begin{enumerate}
\item $B$ is an $R$-base for the $R$-module $X$;
\item $B$ is a maximal $R$-independent set in the $R$-module $X$;
\item $\mathbf 0\notin B$, $R\cdot a\ne R\cdot b$ for any distinct points $a,b\in B$, and the set $\{R\cdot b:b\in B\}$ is maximal independent in the liner $\mathbb PX$.
\end{enumerate}  
\end{proposition}

\begin{proof} Let $B$ be a set in an $R$-module $X$ over a corps $R$.
\smallskip

$(1)\Ra(2)$ Assume that $B$ is an $R$-base for the $R$-module $X$. The $R$-independence of $B$ follows from the uniqueness of the representation of the zero element $\mathbf 0=\sum_{b\in B}0\cdot b$. Assuming that the $R$-independent set $B$ is not maximal, we can find a point $x\in X\setminus B$ such that the set $B\cup\{x\}$ is $R$-independent. Since $B$ is a base, there exists a finitely supported function $f\in R^{\oplus B}$ such that $x=\sum_{b\in B}f(b)\cdot b$. Consider the function $g:B\cup\{x\}\to R$ defined by $g{\restriction}_B=f$ and $g(x)=-1$. Then $\sum_{b\in B\cup\{x\}}g(b)\cdot b=\sum_{b\in B}f(b)-1\cdot x=\mathbf 0$, which contradicts the $R$-independence of the set $B\cup\{x\}$. This contradiction shows that the $R$-base $B$ is maximal $R$-independent in $X$.
\smallskip

$(2)\Ra(3)$ Assume that $B$ is a maximal $B$-independent set in the $R$-module $X$. The $R$-independence of $B$ implies that $\mathbf 0\notin B$ and $R\cdot a\ne R\cdot b$ for any distinct points $a,b\in B$. By Lemma~\ref{l:independent<=>R-independent}, the set $\mathbb PB\defeq \{B\cdot b:b\in B\}$ is independent in the liner $\mathbb PX$. Assuming that $\mathbb PB$ is not maximal independent, we can find an independent set $\mathcal I\subseteq \mathbb PX$ containing the independent set $\mathbb PB$ as a proper subset. Choose any $R$-line $A\in \mathcal I\setminus\mathbb PB$ and a point $a\in A\setminus\{\mathbf 0\}$. Since $A\notin\mathbb PB$, the point $a$ does not belong to the set $R\cdot B$. By Lemma~\ref{l:independent<=>R-independent}, the independence of the set  $\mathcal I=\{R\cdot x:x\in B\cup\{a\}\}$ implies the $R$-independence of the set $B\cup\{a\}$, which contradicts the maximal $R$-independence of the set $B$ in $X$. This contradiction shows that the set $\mathbb PB$ is maximal independent in the liner $\mathbb PX$.
\smallskip

$(3)\Ra(1)$ Assume that $\mathbf 0\notin B$, $R\cdot a\ne R\cdot b$ for any distinct points $a,b\in B$, and the set $\mathbb PB$ is maximal independent in the liner $\mathbb PX$. By Proposition~\ref{p:add-point-to-independent}, the flat hull $\overline{\mathbb PB}$ of the maximal independent set $\mathbb PB$ in the liner $\mathbb PX$ coincides with $\mathbb PX$. By Lemma~\ref{l:independent<=>R-independent}, the set $B$ is $R$-independent in the $R$-module $X$. To prove that $B$ is an $R$-base for $X$, consider the $R$-submodule $\Sigma\defeq\{\sum_{b\in B}f(b)\cdot b:f\in R^{\oplus B}\}$ of the $R$-module $X$. By Lemma~\ref{l:flat-in-PX<=>R-submodule}, the set $\{R\cdot x:x\in\Sigma\setminus\{\mathbf 0\}\}$ is flat in the liner $\mathbb PX$. Since this set contains the set $\mathbb PB$, it also contains its flat hull $\overline{\mathbb PB}$ in the liner $\mathbb PX$. Then $\mathbb PX =\overline{\mathbb PB}=\{R\cdot x:x\in\Sigma\setminus\{\mathbf 0\}\}$ and $X=\bigcup\mathbb PX=\bigcup\big\{R\cdot x:x\in \Sigma\setminus\{\mathbf 0\}\big\}=\Sigma$. Consequently, for every $x\in X=\Sigma$ there exists a finitely supported function $f\in\Sigma$ such that $x=\sum_{b\in B}f(b)\cdot b$. The uniqueness of the function $f$ follows from the $R$-independence of $B$. Therefore, $B$ is an $R$-base for the $R$-module $X$.
\end{proof}

\begin{corollary}\label{c:|base|=||PX||} Let $X$ be an $R$-module over a corps $R$. Every $R$-base $B$ in the $R$-module $X$ has cardinality $|B|=\|\mathbb PX\|$.
\end{corollary}

\begin{proof} Let $B$ be an $R$-base for the $R$-module $X$. By Proposition~\ref{p:R-base<=>}, $\mathbf 0\notin B$,  $R\cdot a\ne R\cdot b$ for any distinct points $a,b\in B$, and $\mathbb PB\defeq\{R\cdot b:b\in B\}$ is a maximal independent set in the liner $\mathbb PX$ such that $|\mathbb PB|=|B|$. By Proposition~\ref{p:PX-projective} the liner $\mathbb PX$ is projective and by Theorem~\ref{t:projective<=>}, the projective liner $\mathbb PX$ is strongly regular. By Corollary~\ref{c:proregular=>ranked}, the regular liner $\mathbb PX$ is ranked and has the Exchange Property. By Corollary~\ref{c:Max=dim}, 
$\|X\|=|\mathbb PB|=|B|$.  
\end{proof}

\begin{definition}\label{d:dimR(X)} For an $R$-module $X$ over a corps $R$,  its \index{$R$-dimension}\defterm{$R$-dimension} $\dim_R(X)$ is defined as the cardinality of any $R$-base for $X$.
\end{definition} 

Corollary~\ref{c:|base|=||PX||} ensures that any two $R$-bases in an $R$-module $X$ over a corps $R$ have the same cardinality, so the cardinal $\dim_R(X)$ is well-defined.

\begin{definition}
Let $X$ be an $R$-module over a corps $X$. A subset $H\subseteq X$ is called an \index{$R$-hyperplane}\defterm{$R$-hyperplane} in the $R$-module $X$ if $H$ is an $R$-submodule of $X$ such that $H\ne X=H+R\cdot x$ for every $x\in X\setminus H$.
\end{definition}

\begin{lemma} Let $X$ be an $R$-module over a corps $R$. A subset $\mathcal H\subseteq \mathbb PX$ is a hyperplane in the liner $\mathbb PX$ if and only if $\bigcup\mathcal H$ is an $R$-hyperplane in the $R$-module $X$.
\end{lemma}

\begin{proof} Let $\mathcal H$ be a subset of the projective space $\IP X$.

 If $\mathcal H\subseteq\mathbb PX$ is a hyperplane in the liner $\mathbb PX$, then $\mathcal H$ is proper flat in $\mathbb PX$ and its union $H\defeq\bigcup\mathcal H$ is a proper $R$-submodule of the $R$-module $X$, by Lemma~\ref{l:flat-in-PX<=>R-submodule}. To show that $H$ is an $R$-hyperplane in the $R$-module, take any point $b\in X\setminus H$ and consider the $R$-submodule $H+R\cdot b$ of the $R$-module $X$. By Lemma~\ref{l:flat-in-PX<=>R-submodule}, the set $\{R\cdot x:x\in H+R\cdot b\}$ is flat in $\IP X$. This flat contains the hyperplane $\mathcal H$ and also the $R$-line $R\cdot b\in \IP X\setminus \mathcal H$. Taking into account that $\mathcal H$ is a hyperplane in $\IP X$,  
we conclude that $\IP X=\overline{\mathcal H\cup\{R\cdot b\}}\subseteq \{R\cdot x:x\in H+R\cdot b\}$ and $X=\bigcup\IP X=\bigcup\{R\cdot x:x\in H+R\cdot b\}=H+R\cdot b$, witnessing that $H=\bigcup\mathcal H$ is an $R$-hyperplane in the $R$-module $X$.
\smallskip

To prove the ``if'' part, assume that the set $H\defeq\bigcup\mathcal H$ is an $R$-hyperplane in the $R$-module $X$. By Lemma~\ref{l:flat-in-PX<=>R-submodule}, $\mathcal H$ is a proper flat in the liner $\IP X$. To show that the flat $\mathcal H$ is a hyperplane in $\IP X$, take any $R$-line $L\in \IP X\setminus\mathcal H$ and consider the flat hull $\overline{\mathcal H\cup\{L\}}$ of the set $\mathcal H\cup\{L\}$ in the liner $\mathbb PX$. By Lemma~\ref{l:flat-in-PX<=>R-submodule}, the union $U\defeq\bigcup\overline{\mathcal H\cup \{L\}}$ is an $R$-submodule of the $R$-module $X$. This submodule contains the set $H\cup L$. Since $H$ is an $R$-hyperplane in $X$, $X=H+L\subseteq U$ and hence $\overline{\mathcal H\cup\{L\}}=\IP X$, witnessing that $\mathcal H$ is a hyperplane in the liner $\IP X$.
\end{proof}



\begin{proposition}\label{p:H=PX-PH} Let $X$ be an $R$-module over a corps $R\ne\{0,1\}$ and $H$ be an $R$-hyperplane in $X$. The liner $H$ is isomorphic to the subliner $\IP X\setminus \IP H$ of the projective liner $\IP X$, and hence $\IP X$ is isomorphic to the spread completion of the affine regular liner $H$.
\end{proposition}

\begin{proof} Choose any point $b\in X\setminus H$ and consider the function $F:H\to \IP X$, assigning to every point $x\in H$ the $R$-line $R\cdot(x+b)$. We claim that $F$ is an injective function with $F[H]=\IP X\setminus\IP H$. To see that $F$ is injective, take any distinct points $x,y\in H$. Assuming that $F(x)=F(y)$, we can find nonzero elements $r,s\in R$ such that $r\cdot(x+b)=s\cdot(y+b)$. Assuming that $r\ne s$, we conclude that $b=(s-r)^{-1}\cdot(r\cdot x-s\cdot y)\in H$, which contradicts the choice of $b\notin H$. This contradiction shows that $r=s$ and then $r\cdot(x+b)=s\cdot(y+b)$ implies $x=y$, which contradicts the choice of the points $x,y$. This contradiction shows that the function $F:H\to\IP X$ is injective. 

Since $H$ is an $R$-hyperplane in $X$, $X=H+R\cdot b$ and hence for every $a\in X=H+R\cdot b$, there exist elements $x\in H$ and $r\in R$ such that $a=x+r\cdot b$. If $a\notin H$, then $r\ne 0$ and hence $r^{-1}\cdot a=r^{-1}\cdot x+b$ and $R\cdot a=R\cdot r^{-1}\cdot a=F(r^{-1}\cdot x)\in F[H]$, witnessing that $\IP X\setminus\IP H\subseteq F[H]$. Assuming that $\IP X\setminus\IP H\ne F[H]$, we could find an $R$-line $L\in F[H]\setminus(\IP X\setminus\IP H)=\IP H\cap F[H]$. Then $L=R\cdot (x+b)$ for some $x\in H$ and hence $b=(b+x)-x\in L-x\subseteq H-H=H$, which contradicts the choice of $b$. This contradiction shows that $F[H]=\IP X\setminus\IP H$ and hence $F:H\to \IP X\setminus \IP H$ is a bijection.

To prove that $F$ is an isomorphism of the liners $H$ and $\IP X\setminus \IP H$, it suffices to show that for every line $L\subseteq H$, the image $F[L]$ is a line in the subliner $\IP X\setminus \IP H$ of $\IP X$. Fix any distinct points $u,v\in L$ and consider the $R$-lines $U\defeq R\cdot(u+b)$ and $V\defeq R\cdot(v+b)$ in the $R$-module $X$.  For every point $w\in\Aline uv=u+R\cdot(u-v)$ there exists an element $r\in R$ such that $w=u+r\cdot(u-v)=(1-r)\cdot u+r\cdot v$. Then $F[w]=R\cdot (w+b)=R\cdot((1-r)\cdot (u+b)+r\cdot (v+b))\subseteq R\cdot (u+b)+R\cdot(v+b)=U+V$ and hence $F(w)\in\Aline UV\subseteq \IP X$. This shows that $F[\Aline uv]\subseteq\Aline UV$. 

On the other hand, for every $Z\in \Aline UV\setminus\IP H$, we can fix any nonzero element $z\in Z\subseteq U+V$ and find elements $s,t\in R$ such that $z=s\cdot (u+b)+ t\cdot(v+b)$. If $s+t=0$, then $z=s\cdot u+t\cdot v+(s+t)\cdot b=s\cdot u+t\cdot v\in H$, which contradicts the choice of the $R$-line $Z\notin\IP  H$. Therefore, $s+t\ne 0$ and for the point $\tilde z\defeq(s+t)^{-1}\cdot s\cdot u+(s+t)^{-1}\cdot t\cdot v\in \Aline uv=L$ we have  $Z=R\cdot z=R\cdot (s+t)^{-1}\cdot z=R\cdot(((s+t)^{-1}\cdot s\cdot u+(s+t)^{-1}\cdot t\cdot v)+b)=F(\tilde z)\in F[L]$, witnessing that $F[L]=\Aline UV\setminus \IP H$. By Theorem~\ref{t:liner-isomorphism<=>}, the bijective function $F:H\to\IP X\setminus\IP H$ is an isomorphisms of the liners $H$ and $\IP X\setminus \IP H$.    

By Theorem~\ref{t:R-module=>Des-aff-reg}, the $R$-module $H$ endowed with the canonical line relation is an affine regular liner. By Corollary~\ref{c:affine-spread-completion}, the affine regular liner $H$ is completely regular, and by Corollary~\ref{c:spread-3long}, the spread completion $\overline H$ of $H$ is a projective $3$-long liner. Then $\overline H$ is a projective completion of $H$.  
By Theorem~\ref{t:extend-isomorphism-to-completions}, the isomorphism $F:H\to\IP X\setminus \IP H$ extends to an isomorphism $\bar F:\overline H\to \IP X$ of the projective completions of the liners $H$ and $\IP X\setminus\IP H$. Therefore, the projective space $\IP X$ is isomorphic to the spread completion of the affine regular liner $H$.
\end{proof}

\begin{theorem}\label{t:PX-Desarguesian} The projective space $\mathbb PX$ of any $R$-module $X$ over a corps $R$ is a $3$-long Desarguesian projective liner of rank $\|\mathbb PX\|=\dim_R(X)$.
\end{theorem}

\begin{proof} By Lemma~\ref{l:PX-3long} and Proposition~\ref{p:PX-projective}, $\mathbb PX$ is a projective liner. By Theorem~\ref{t:projective<=>}, the projective liner $\mathbb PX$ is strongly regular, and by Proposition~\ref{p:k-regular<=>2ex}, the regular liner $\IP X$ is ranked.  By the  Kuratowski--Zorn Lemma, the $R$-module $X$ contains a maximal $R$-independent set $B$. By Proposition~\ref{p:R-base<=>}, $\mathbf 0\notin B$, $R\cdot a\ne R\cdot b$ for any distinct elements $a,b\in B$, and $\mathbb PB\defeq\{R\cdot b:b\in B\}$ is a maximal independent set in the liner $\mathbb PX$. By Corollary~\ref{c:|base|=||PX||} and Definition~\ref{d:dimR(X)}, $$\|\mathbb PX\|=|\IP B|=|B|=\dim_R(X).$$ It remains to prove that the projective liner $\IP X$ is Desarguesian.

If $R=\{0,1\}$, then the projective liner $\IP X$ is Steiner, according to Proposition~\ref{p:PX-Steiner}. By Theorem~\ref{p:Steiner+projective=>Desargues}, the Steiner regular liner $\IP X$ is Desarguesian. So, assume that $R\ne\{0,1\}$.

If $\|\mathbb PX\|\ne 3$, then the $3$-long projective liner $\mathbb PX$ is Desarguesian, by Theorem~\ref{t:proaffine-Desarguesian}. So, assume that $3=\|\mathbb PX\|=\dim_R(X)=|B|$.  Choose any point $p\in B$ and observe that the set $B_p\defeq B\setminus\{p\}$ is an $R$-base for the $R$-submodule $H\defeq\big\{\sum_{x\in B_p}f(b)\cdot b:f\in R^{\oplus B_p}\big\}$ of $X$. Since $H+R\cdot p=X$, the set $H$ is an $R$-hyperplane in the $R$-module $X$. By Theorem~\ref{t:R-module=>Des-aff-reg}, the $R$-module $H$ endowed with the canonical line relation is a Desarguesian affine regular liner. By Theorem~\ref{t:Desargues-completion}, the spread completion $\overline H$ of the $3$-long Desarguesian affine regular liner $H$ is a Desarguesian $3$-long  projective liner. By Proposition~\ref{p:H=PX-PH}, the spread completion $\overline H$ of $H$ is isomorphic to the projective space $\IP X$. Then the projective liner $\IP X$ is Desarguesian, being isomorphic to the Desarguesian projective liner $\overline H$.
\end{proof}

\begin{theorem}\label{t:Des-projspace<=>} A projective space $Y$ is Desarguesian if and only if it is isomorphic to the projective space $\IP X$ of some $R$-module $X$ over a corps $R$. 
\end{theorem}

\begin{proof} The ``if'' part follows from Theorem~\ref{t:PX-Desarguesian}. To prove the ``only if'' part, assume that $Y$ is a Desarguesian projective space.  By the Kuratowski--Zorn Lemma, the liner $Y$ contains a maximal independent set $B$. Choose any point $o\in B$ and observe that the flat $H\defeq\overline{B\setminus\{o\}}$ is a hyperplane in $Y$ (this follows from the rankedness of the projective space $Y$).  By Theorem~\ref{t:affine<=>hyperplane}, the subliner $A\defeq Y\setminus H$ is affine and regular. Theorem~\ref{t:pD-minus-flat}, the affine subliner $A=Y\setminus H$ of the Desargusian projective liner $Y$ is Desarguesian.

If the liner $Y$ is $4$-long, then then the affine liner $A=Y\setminus H$ is $3$-long. By Corollary~\ref{c:procompletion-rank}, $\|A\|=\|Y\|\ge 3$ and hence $A$ is an affine space. By Theorem~\ref{t:aff-Des<=>}, the Desarguesian affine space $A$ is isomorphic to some $R$-module $M$ over a corps  $R\ne\{0,1\}$ (in fact, for $M$ we can take the vector space$\overvector A$ of $A$, which is an $\IR_A$-module over the corps of scalars $\IR_A$).  Consider the $R$-module $X\defeq M\times R$ and observe that $M_0\defeq M\times\{0\}$ is an $R$-hyperplane in the $R$-module $X$. By Proposition~\ref{p:H=PX-PH}, the hyperplane $M_0$ is isomorphic to the subliner $\IP X\setminus\IR M_0$ of the projective space $\IP X$ of the $R$-module $X$. Since the $R$-modules $M$ and $M_0$ are isomorphic, they also are isomorphic as liners. Then there exists a liner isomorphism $F:A\to \IP X\setminus\IP M_0$. Since $Y$ is a projective completion of $A$ and 
$\IP X$ is a projective completion of $\IP X\setminus\IP M_0$, the liner isomorphism $F$ extends to a liner isomorphism $\bar F:Y\to \IP X$, according to Theorem~\ref{t:extend-isomorphism-to-completions}. Therefore, the Desarguesian projective space $Y$ is isomorphic to the projective space $\IP X$ of the $R$-module $X$.
\smallskip

Next, assume that the projective space $Y$ is not $4$-long. By Corollary~\ref{c:Avogadro-projective}, the $3$-long projective liner $Y$ is $2$-balanced. Since $Y$ is not $4$-long, $|Y|_2=3$, which means that $Y$ is a Steiner projective liner. Then $Y$ carries the midpoint operation $\circ:Y\times Y\to Y$ assigning to any points $x,y\in Y$ the unique point $x\circ y\in Y$ such that $\Aline xy=\{x,y,x\circ y\}$.
It is clear that the midpoint operation $\circ$ on $Y$ is commutative. 

\begin{claim}\label{cl:center-Fano} For every points $x,y,z\in A$ with $x\ne y\ne z$ we have $x\circ(y\circ z)=(x\circ y)\circ z\in A$.
\end{claim}

\begin{proof} If $x=z$, then $x\circ(y\circ z)=z\circ(y\circ x)=(x\circ y)\circ z$, by the commutativity of the midpoint operation. So, assume that $x\ne z$. By Corollary~\ref{c:line-meets-hyperplane}, every line $L$ in $Y$ intersects the hyperplane $H$, which implies that $|L\cap A|\le 2$ and hence the $3$-element set $\{x,y,z\}\subseteq A$ is not a line in $Y$. Then $\overline{\{x,y,z\}}$ is a plane in $Y$. By Theorem~\ref{t:Steiner<=>projective}, $|\overline{\{x,y,z\}}|=7$. Since every line in $Y$ intersects the hyperplane $H$, $\{x\circ y,y\circ z\}\subseteq H$ and hence $(x\circ y)\circ(y\circ z)\in H$. 
Therefore, $|\overline{\{x,y,z\}}\cap A|=|\overline{\{x,y,z\}}\setminus\{x\circ y,y\circ z,(x\circ y)\circ(y\circ z)\}|=4$ and $\{x\circ(y\circ z),(x\circ y)\circ z\}=\overline{\{x,y,z\}}\setminus\{x\circ y,y\circ z,(x\circ y)\circ(y\circ z),x,y,z\}$ is a singleton, which implies that $x\circ(y\circ z)=(x\circ y)\circ z$.
\end{proof}

Let $+:A\times A\to A$ be the binary operation on $A$, assigning to every points $x,y\in A$ the point $$x+y\defeq \begin{cases}
x&\mbox{if $y=o$};\\
y&\mbox{if $x=o$};\\
x\circ(o\circ y)=(x\circ o)\circ y&\mbox{if $x\ne o\ne y$}.
\end{cases}
$$
We recall that $o$ is a unique point of the set $B\setminus H\subset X\setminus H=A$.
Claim~\ref{cl:center-Fano} ensures that the operation $+$ is well-defined.

\begin{claim}\label{cl:+Boolean} The binary operation $+:A\times A\to A$ has the following properties:
\begin{enumerate}
\item $\forall x,y\in A\;\;(x+y=y+x)$;
\item $\forall x,y\in A\;\;(x+y=o\;\Leftrightarrow\;x=y)$;
\item $\forall x,y,z\in A\;\;(x+y=x+z\;\Leftrightarrow y=z)$;
\item $\forall x,y,z\in A\;\;(x+y)+z=x+(y+z)$.
\end{enumerate}
\end{claim} 

\begin{proof} 1. The commutativity of the binary operation $+$ follows from the commutativity of the midpoint operation and Claim~\ref{cl:center-Fano}.
\smallskip

2. Given any points $x,y\in A$, we should check that $x+y=o$ if and only if $x=y$. If $x=y=o$, then $x+y=o$, by the definition of the binary operation $+$. If $x=y\ne o$, then $x+y=x+x=x\circ(o\circ x)=x\circ(x\circ o)=o$, by the commutativity and involutarity of the midpoint operation $\circ$, see Theorem~\ref{t:Steiner<=>midpoint}. Therefore, $x=y$ implies $x+y=o$.

Next, assume that  $x\ne y$. If $o=x$, then $x+y=o+y=y\ne x=o$. If $o=y$, then $x+y=x+o=x\ne y=o$. If $x\ne o\ne y$, then $x+y=x\circ(o\circ y)$. Observe that $\Aline oy=\{o,y,o\circ y\}$. It follows from $x\in A\setminus \{o,y\}$ and $o\circ y\in H=Y\setminus A$ that $x\notin\Aline oy$ and hence $x+y=x\circ (o\circ y)\notin \Aline oy$ and $x+y\ne o$.
\smallskip

3. Given any points $x,y,z\in A$ with $x+y=x+z$, we should prove that $y=z$.
If $x=o$, then $y=o+y=x+y=x+z=o+z=z$ and we are done.
 So, we assume that $x\ne o$. If $y=o\ne z$, then $x=x+o=x+y=x+z=x\circ(o\circ z)$ implies $x=o\circ z\in H$, which contradicts the choice of $x$. By analogy we can derive  a contradiction assuming that $y\ne o=z$. Therefore, $y\ne o\ne z$ and $o\circ y\ne o\ne o\circ z$. The definition of the binary operation $+$ ensures that
$x\circ(o\circ y)=x+y=x+z=x\circ(o\circ z)$. The involutarity of the midpoint operation $\circ$ implies that
$$o\circ y=x\circ(x\circ(o\circ y))=x\circ(x\circ(o\circ z))=o\circ z$$ and
$y=o\circ(o\circ y)=o\circ(o\circ z)=z.$
\smallskip

4. To prove that the binary operation $+$ is associative, choose any points $x,y,z\in A$. We have to prove that $(x+y)+z=x+(y+z)$.

If $x=o$, then $(x+y)+z=(o+y)+z=y+z=o+(y+z)=x+(y+z)$.

If $y=o$, then $(x+y)+z=(x+o)+z=x+z=x+(o+z)=x+(y+z)$.

If $z=o$, then $(x+y)+z=(x+y)+o=x+y=x+(y+o)=x+(y+z)$.

So, assume that $o\notin\{x,y,z\}$. 

If $\|\{o,x,y,z\}\|=2$, then $\overline{\{o,x,y,z\}}$ is a line in $Y$ such that $2\le|\{o,x,y,z\}\cap A\}|\le|\overline{\{o,x,y,z\}}\cap A|=|\overline{\{o,x,y,z\}}\setminus H|=3-1=2$, by Corollary~\ref{c:line-meets-hyperplane}. Then $x=y=z$ and $(x+y)+z=(x+x)+x=x+(x+x)=x+(y+z)$, by the commutativity of the operation $+$.

Next, assume that $\|\{o,x,y,z\}\|=3$, which implies that $P\defeq\overline{\{o,x,y,z\}}$ is  plane in the Steiner projective liner $Y$. By Theorem~\ref{t:Steiner<=>projective}, $|P|=|Y|_3=7$, which means that $P$ is a Steiner projective plane. Then $|P\setminus H|=7-3=4$. 

If $|\{x,y,z\}|=3$, then $\{o,x,y,z\}=P\setminus H=\{o,x,y,x+y\}=\{o,y,z,y+z\}$ and hence  $x+y=z$ and $y+z=x$. Then $(x+y)+z=z+z=o=x+x=x+(y+z)$.

So, assume that $|\{x,y,z\}|=2$. If $x=z$, then $(x+y)+z=(x+y)+x=x+(y+x)=x+(y+z)$, by the commutativity of the binary operation $+$.

If $x=y$, then $|\{x,y,z\}|=2$ implies $x=y\ne z$ and $\overline{\{o,x,y,z\}}\cap A=\overline{\{o,x,z\}}\cap A=\{o,x,z,x+z\}$.
Claim~\ref{cl:+Boolean}(2) ensures that $(x+y)+z=(x+x)+z=o+z=z$.
On the other hand, $x=y\ne z\ne o$ implies $y+z\ne o$, $x+(y+z)\ne x$, $y+z=x+z\ne x+o=x$ and $x+(y+z)\ne x+x=o$, by Claim~\ref{cl:+Boolean}(2,3). Also $y=x\ne o$ implies $x+(y+z)\ne y+z=x+z$. Then $$x+(y+z)=x+(x+z)\in (A\cap\overline{\{o,x,z\}})\setminus\{o,x,x+z\}=\{o,x,x+z,z\}\setminus\{o,x,x+z\}=\{z\}=\{(x+y)+z\}$$ and hence $(x+y)+z=x+(y+z)$.

By analogy we can prove that $y=z$ implies $(x+y)+z=x+(y+z)$.

Finally, assume that $\|\{o,x,y,z\}\|=4$. By Corollary~\ref{c:projective-order-n}, the Steiner projective space $\overline{\{o,x,y,z\}}$ has cardinality $2^3+2^2+2^1+1=15$.
The intersection $P=\overline{\{o,x,y,z\}}\cap H$ is a Steiner projective plane and hence contains exactly 7 points. Then the set $P\cup\{o,x,y,z,x+y,y+z,x+z\}$ contains $7+7=14$ points and the set $S\defeq\overline{\{o,x,y,z\}}\setminus(P\cup\{o,x,y,z,x+y,y+z,x+z\})$ is a singleton. Claim~\ref{cl:+Boolean}(2,3) implies that the points $(x+y)+z=(x+y)\circ(o\circ z)=(x\circ(o\circ y))\circ(o\circ z)$ and $x+(y+z)=x\circ(o\circ (y+z))=x\circ(o\circ(y\circ(o\circ z))$ belongs to the singleton
$S$, which implies that $(x+y)+z=x+(y+z)$.
\end{proof}

By Claim~\ref{cl:+Boolean}, $(A,+)$ is a commutative group such that $x+x=o$ for every $x\in X$. Consider the two-element field $\IF_2\defeq\{0,1\}$ and endow $A$ with the  operation $\cdot:\IF_2\times A\to A$ defined by $0\cdot x=o$ and $1\cdot x=x$ for every $x\in A$. It is easy to see that $A$ is an $\IF_2$-module. Consider the $\IF_2$-module $X\defeq A\times\IF_2$. We claim that the projective space $Y$ is isomorphic to the projective space $\IP X$ of the $\IF_2$-module $X=A\times\IF_2$. 

Consider the bijective function $F:Y\to \IP X$ assigning to every point $y\in Y$ the $\IF_2$-line
$$F(y)\defeq\begin{cases}
\{(o,0),(y,1)\}&\mbox{if $y\in A$};\\
\{(o,0),(o\circ y,0)\}&\mbox{if $y\in Y\setminus A$}.
\end{cases}
$$
It is easy to check that $F$ is an isomorphism of the projective liners $Y$ and $\IP X$. 
\end{proof} 

\begin{exercise} Check that the bijection $F$ defined at the end of the proof of Theorem~\ref{t:Des-projspace<=>} is indeed an isomorphism of the projective spaces $Y$ and $\IP X$.
\end{exercise}

\section{Automorphisms of Desarguesian projective spaces}

The following theorem establishes an important homogeneity properties of projective spaces of modules over corps.

\begin{theorem}\label{t:PX-autoextend} Let $X$ be an $R$-module over a corps $R$ and $\mathcal B,\mathcal B'$ be two maximal independent sets in the projective space $\IP X$. Every bijection $F:\mathcal B\to \mathcal B'$ extends to an automorphism $\bar F:\IP X\to \IP X$ of the liner $\IP X$.
\end{theorem}

\begin{proof} Choose sets $B,B'\subseteq X$ such that $\mathcal B=\{R\cdot x:x\in B\}$, $\mathcal B'=\{R\cdot x:x\in B'\}$ and $|B\cap L|=|B'\cap L'|=1$ for every $R$-lines $L\in \mathcal B$ and $L'\in \mathcal B'$. By Proposition~\ref{p:R-base<=>}, the sets $B$ and $B'$ are $R$-bases of the  $R$-module $X$. 

The bijective function $F:\mathcal B\to\mathcal B'$ induces a unique bijective function $f:B\to B'$ such that for every $b\in B$, the $R$-line $F(R\cdot b)$ coincides with the $R$-line $R\cdot f(b)$. Let $\bar  f:X\to X$ be the function assigning to every point $x\in X$ the point $\sum_{b\in B}c_x(b)\cdot f(b)\in X$ where $c_x\in R^{\oplus B}$ is a unique  finitely supported function such that $x=\sum_{b\in B}c_x(b)\cdot b$.  Taking into account that $B,B'$ are $R$-bases for the $R$-module $X$, we conclude that the function $\bar f$ is an automorphism of the $R$-module $X$ such that $\bar f(b)=f(b)$ for every $b\in B$. Then the function $\bar F:\IP X\to \IP X$, $\bar F:L\mapsto\bar f[L]$, is an automorphism of the projective liner $\IP X$ such that $\bar F(L)=F(L)$ for every $L\in\M$.
\end{proof}

Theorems~\ref{t:PX-autoextend} and \ref{t:Des-projspace<=>} imply the following important fact.

\begin{theorem}\label{t:Des-autoextend} Every bijection $\phi:M\to M'$ between maximal independent sets in a Desarguesian projective space $X$ extends to an automorphism of the liner $X$.
\end{theorem}

\begin{corollary}\label{c:hyperplanes-automorphism} For every hyperplanes $H,H'$ in a  Desarguesian projective space $X$ and every flat $\Pi\subseteq H\cap H'$, there exists an automorphism $\Phi:X\to X$ such that $\Phi[\Pi]=\Pi$, $\Phi[H]=H'$.
\end{corollary}

\begin{proof} Using the Kuratwski-Zorn Lemma, choose a maximal independent set $A$ in the flat $\Pi$. Next, enlarge $A$ to a maximal independent set $B$ in the hyperplane $H$ and also to a maximal independent set $B'$ in the hyperplane $H'$. Finally, chose any points $b\in X\setminus H$ and $b'\in X\setminus H'$. The maximal independence of the sets $B,B'$ implies the maximal $A$-independence of the sets $(B\cup\{b\})\setminus A$ and $(B'\cup\{b'\})\setminus A$. By Theorem~\ref{t:Max=codim}, $|(B\cup \{b\})\setminus A|=\|X\|_\Pi=|(B'\cup\{b'\})\setminus A|$. Then there exists a bijection $F:B\cup\{b\}\to B'\cup\{b'\}$ such that $F(b)=b'$ and $F(a)=a$ for every $a\in A$. By Theorem~\ref{t:Des-autoextend}, the bijection $F:B\cup\{b\}\to B'\cup\{b'\}$ enxtends to an automorphism $\Phi:X\to X$ of the Desarguesian projective space $X$. Then $\Phi[\Pi]=\Phi[\overline A]=\overline A=\Pi$ and $\Phi[H]=\Phi[\overline{B}]=\overline{B'}=H'$.
\end{proof}

\begin{corollary}\label{c:hyperplanes-are-isomorphic} For any hyperplanes $H,H'$ in a Desarguesian projective space $X$, the affine liners $X\setminus H$ and $X\setminus H'$ are isomorphic.
\end{corollary}

\begin{corollary}\label{c:affine-isomorphic<=>} Two Desarguesian affine spaces $X,Y$ are isomorphic if and only if their spread completions are isomorphic.
\end{corollary}

\begin{proof} The ``only if'' part follows from Theorem~\ref{t:extend-isomorphism-to-completions} and Corollary~\ref{c:affine-spread-completion}. To prove the ``if'' part, assume that the spread completions of two Desarguesian affine spaces $X,Y$ are isomorphic. Fix any isomorphism $\Psi:\overline X\to\overline Y$. By Theorem~\ref{t:Desargues-completion}, the spread completion $\overline Y$ of the Desarguesian affine space $Y$ is a Desarguesian projective space. By Theorem~\ref{t:affine<=>hyperplane}, the horizon $\partial X\defeq\overline X\setminus X$ is a hyperplane in the projective space $\overline X$ and hence its image $H\defeq \Phi[\partial X]$ is a hyperplane in the projective space $\overline Y$. By Corollary~\ref{c:hyperplanes-automorphism}, there exists an automorphism $\Psi:\overline Y\to\overline Y$ such that $\Psi[H]=\partial Y$. Then $F\defeq \Psi\circ\Phi:\overline X\to\overline Y$ is a liner isomorphism such that $F[X]=Y$, witnessing that the Desarguesian affine spaces $X,Y$ are isomorphic.  
\end{proof}

In fact, Corollary~\ref{c:affine-isomorphic<=>} can be generalized to Thalesian affine spaces.

\begin{theorem}\label{t:para-D-isomorphic<=>} Two Thalesian affine spaces are isomorphic if and only if their spread completions are isomorphic.
\end{theorem}

\begin{proof} Let $X,Y$ be Thalesian affine spaces. By Corollary~\ref{c:affine-spread-completion}, the spread completions $\overline X,\overline Y$ are $3$-long projective liners. If the liners $X,Y$ are isomorphic, then their spread completions are isomorphic by Theorem~\ref{t:extend-isomorphism-to-completions}. Now assume that the spread completions $\overline X$ and $\overline Y$ are isomorphic. Then $\|\overline X\|=\|\overline Y\|$ and hence $\|X\|=\|\overline X\|=\|\overline Y\|=\|Y\|$, by Corollary~\ref{c:procompletion-rank}. If $\|X\|=\|Y\|\ge 4$, then the affine spaces $X,Y$ are Desarguesian, by Theorem~\ref{c:affine-Desarguesian}. In this case we can apply Corollary~\ref{c:affine-isomorphic<=>} and conclude that the affine liners $X,Y$ are isomorphic. So, assume that $\|X\|=\|Y\|=3$.  Fix any isomorphism $F:\overline X\to\overline Y$ of the projective planes $\overline X,\overline Y$ and consider the lines $\partial Y\defeq \overline Y\setminus Y$, $\partial X\defeq\overline X\setminus X$. If $F[\partial X]=\partial Y$, then the restriction $F{\restriction}_X:X\to Y$ is an isomorphism of the affine liners $X,Y$ and we are done.

So, assume that $F[\partial X]\ne\partial Y$. Since $F[\partial X]$ and $\partial Y$ are distinct lines in the projective plane $\overline Y$, there exists a unique point $o\in F[\partial X]\cap\partial Y$.  Since the affine space $Y$ is $3$-long, its spread completion $\overline Y$ is $4$-long. By Proposition~\ref{p:cov-aff}, there exists a point $y\in\overline Y\setminus (\partial Y\cup F[\partial X])$. Then the line $\Lambda\defeq\Aline oy$ in $\overline Y$ is concurrent with the lines $\partial Y$ and $F[\partial X]$. 
Fix any point $z\in \partial Y\setminus\{o\}$. Since the affine space $X$ is Thalesian, so is its isomorphic copy $F[X]\subseteq \overline Y$. Since the points $y,z$ belong to the Thalesian affine space $F[X]$, we can consider the vector $\overvector{yz}\in\overvector{F[X]}$, which determines an automorphism $\Phi:F[X]\to F[X]$, $x\mapsto x+\overvector{yz}$, such that $\Phi[L]\parallel L$ for every line $L$ in $F[X]$. The automorphism $\Phi$ extends to an automorphism $\bar\Phi:\overline Y\to\overline Y$ such that $\bar\Phi(y)=\Phi(y)=y+\overvector{yz}=z$ and $\bar\Phi(o)=o$. Then $\bar\Phi[\Lambda]=\Phi[\Aline oy]=\Aline oz=\partial Y$. By analogy we can construct an automorphism $\bar\Psi:\overline Y\to\overline Y$ such that $\bar\Psi[\partial X]=\Lambda$. Then $I\defeq \bar\Phi\circ\bar\Psi\circ F:\overline X\to\overline Y$ is a liner isomorphism such that $I[\partial X]=\partial Y$ and $I[X]=I[\overline X\setminus\partial X]=\overline Y\setminus I[\partial X]=\overline Y\setminus \partial Y=Y$. Therefore, $I{\restriction}_X:X\to Y$ is a required isomorphism of the Thalesian affine spaces $X,Y$.
\end{proof}

Theorem~\ref{t:para-D-isomorphic<=>} allows us to define the \index{scalar corps}\index[note]{$\IR_X$}\defterm{scalar corps} $\IR_X$ of every projective space $X$ as follows. If for some hyperplane $H\subseteq X$, the affine liner $A\defeq X\setminus H$ is a  Thalesian affine space, then put $\IR_X\defeq\IR_A$. Theorem~\ref{t:para-D-isomorphic<=>} ensures that the isomorphic type of the corps $\IR_A$ does not depend on the choice of the hyperplane $H\subseteq X$ with Thalesian complement. If for every hyperplane $H\subseteq X$ the complement $X\setminus H$ is not a Thalesian affine space, then put $\IR_X\defeq\{0,1\}$ be the two-element field.
In particular, $\IR_X=\{0,1\}$ for any Steiner projective liner $X$.


\begin{proposition}\label{p:IR(PX)=R} Let $X$ be an $R$-module over a corps $R$. If $\dim_R(X)\ge 3$, then the scalar corps $\IR_{\IP X}$ of the projective space $\IP X$ is isomorphic to the corps $R$.
\end{proposition}

\begin{proof} If $R=\{0,1\}$, then the projective space $\IP X$ is Steiner and its scalar corps is equal to $\{0,1\}=R$. So, assume that $R\ne\{0,1\}$. Using the Kuratowski--Zorn Lemma, choose an $R$-hyperplane $H$ in the $R$-module $X$. By Theorem~\ref{t:R-module=>Des-aff-reg}, the $R$-module $H$ endowed with the canonical line relation is a Desargusian affine space. By Proposition~\ref{p:IR_X=R}, the scalar corps $\IR_H$ of the Desarguesian affine space $H$ is isomorphic to the corps $R$.  By Proposition~\ref{p:H=PX-PH}, the projective space $\IP X$ is isomorphic to the spread completion $\overline H$ of the Desarguesian affine space $H$ and hence the scalar corps $\IR_{\IP X}$ is isomorhic to the scalar corps $\IR_H$ of the Desarguesian affine space $H$, which is isomorphic to the corps $R$.
\end{proof}

\section{Isomorphic classification of Desarguesian proaffine spaces}


\begin{theorem}\label{t:Dp=PRX} Every Desarguesian projective space $X$ is isomorphic to the projective space $\IP\IR_X^{\oplus\|X\|}$ of the $\IR_X$-module $\IR_X^{\oplus\|X\|}$.
\end{theorem}

\begin{proof} By Theorem~\ref{t:Des-projspace<=>}, the Desarguesian projective space $X$ is isomorphic to the projective space $\IP Y$ of some $R$-modude $Y$ over a corps $R$. Using the Kuratowski--Zorn Lemma, fix an $R$-base $B\subseteq Y$ and consider the bijective function $c_*:Y\to R^{\oplus B}$ assigning to every point $y\in Y$ the unique finitely supported function $c_y\in R^{\oplus B}$ such that $y=\sum_{b\in B}c_y(b)\cdot b$. It is easy to see that $c_*:Y\to R^{\oplus B}$ is an isomorphism of the $R$-modules $Y$ and $R^{\oplus B}$. Then the projective spaces $X$ and $\IP Y$ are isomorphic to the projective space $\IP R^{\oplus B}$.

By Proposition~\ref{p:IR(PX)=R}, the scalar corps $\IR_X$ and $\IR_{\IP Y}$ of the (isomorphic) projective spaces $X$ and $\IP Y$ are isomorphic to the corps $R$. By Corollary~\ref{c:|base|=||PX||}, $|B|=\dim_R(Y)=\|\IP Y\|=\|X\|$. Then the projective spaces $X$ and  $\IP R^{\oplus B}$ are isomorphic to the projective space $\IP \IR_X^{\oplus\|X\|}$.
\end{proof}

\begin{corollary}\label{c:Desproj-isomorphic<=>} Two  Desarguesian projective spaces $X,Y$ are isomorphic if and only if $\|X\|=\|Y\|$ and the corps $\IR_X$ and $\IR_Y$ are isomorphic.
\end{corollary}


\begin{theorem}\label{t:Dproaff-iso<=>} Two Desarguesian proaffine spaces $X,Y$ are isomorphic if and only if $\|\partial X\|=\|\partial Y\|$, $\|\overline X\|_{\partial X}=\|\overline Y\|_{\partial Y}$ and the corps $\IR_X$ and $\IR_Y$ are isomorphic.
\end{theorem}

\begin{proof} By the definition of a space, the spaces $X,Y$ are $3$-long, regular, and have rank $\ge 3$.  By Theorem~\ref{t:Desargues-completion} and Corollary~\ref{c:procompletion-rank}, the spread completions $\overline X$ and $\overline Y$ of the Desarguesian proaffine spaces $X$ and $Y$  are Desarguesian projective spaces. 
By Proposition~\ref{p:nD=>pP}, the horizons $\partial X$ and $\partial Y$ are flat in the projective liners $\overline X$ and $\overline Y$, respectively.
\smallskip

To prove the ``only if'' part, assume that the liners $X,Y$ are isomorphic. Fix any isomorphism $F:X\to Y$. By Theorem~\ref{t:extend-isomorphism-to-completions}, the isomorphism $F$ extends to an isomorphism $\bar F:\overline X\to\overline Y$ of the spread completions. Then the projective spaces $\overline X$ and $\overline Y$ are isomorphic, and so are their scalar corps $\IR_{\overline X}=\IR_X$ and $\IR_{\overline Y}=\IR_Y$. 

The restriction $\bar F{\restriction}_{\partial X}:\partial X\to\partial Y$ is an isomorphism between the liners $\partial X$ and $\partial Y$, which implies that $\|\partial X\|=\|\partial Y\|$. Choose a maximal independent set $A$ in the flat $\partial X\subseteq\overline X$ and enlarge $A$ to a maximal independent set $B\subseteq \overline X$. Then the set $B\setminus A$ is maximal $\partial X$-independent in $\overline X$ and its image $\bar F[B]$ is a maximal $\partial Y$-independent set in the liner $\overline Y$. Applying Theorem~\ref{t:Max=codim}, we conclude that $\|\overline X\|_{\partial X}=|B\setminus A|=|\bar F[B]\setminus\bar F[A]|=\|\overline Y\|_{\partial Y}$.
This completes the proof of the ``only if'' part. 
\smallskip

To prove the ``if'' part, assume that  $\|\partial X\|=\|\partial Y\|$, $\|\overline X\|_{\partial X}=\|\overline Y\|_{\partial Y}$, and the scalar corps $\IR_X$ and $\IR_Y$ are 
 isomorphic. By Corollary~\ref{c:rank+},
$$\|\overline X\|=\|\partial X\|+\|\overline X\|_{\partial X}=\|\partial Y\|+\|\overline Y||_{\partial Y}=\|\overline Y\|,$$ and by Corollary~\ref{c:Desproj-isomorphic<=>}, the Desarguesian projective spaces $\overline X$ and $\overline Y$ are isomorphic. Fix an isomorphism $\Phi:\overline X\to\overline Y$. Using the Kuratowski--Zorn Lemma, find  maximal independent sets $A\subseteq \partial X$ and $A'\subseteq \partial X$. Enlarge the sets $A,A'$ to maximal independent sets $B\subseteq \overline X$ and $B'\subseteq \overline Y$. By Theorem~\ref{t:Max=codim},
$|A|=\|\partial X\|=\|\partial Y\|=|A'|$ and $|B\setminus A|=\|\overline X\|_{\partial X}=\|\overline Y\|_{\partial Y}=|B'\setminus A'|$. Then there exists a bijection $F:\Phi[B]\to B'$ such that $F[\Phi[A]]=A'$. Since $\Phi[B]$ and $B'$ are maximal independent sets in the Desarguesian projective liner $\overline Y$, the bijection $F:\Phi[B]\to B'$ extends to an automorphism $\bar F:\overline Y\to\overline Y$ of the Desarguesian projective liner $\overline Y$, by Theorem~\ref{t:Des-autoextend}. Then $\Psi\defeq\bar F\circ \Phi:\overline X\to\overline Y$ is a liner isomorphism such that $\Psi[\partial X]=\Psi[\overline A]=\overline{[\Psi[A]]}=\overline{F[\Phi[A]]}=\overline{A'}=\partial Y$. Then $\Psi[X]=\Psi[\overline X\setminus\partial X]=\overline Y\setminus\Psi[\partial X]=\overline Y\setminus\partial Y=Y$, which implies that the restriction $\Psi{\restriction}_X:X\to Y$ is an isomorphism of the liners $X,Y$.
\end{proof}

\section{Isomorphic classification of finite corps}

By Theorem~\ref{t:Dproaff-iso<=>}, two Desarguesian proaffine spaces $X,Y$ are isomorphic if and only if $\|\partial X\|=\|\partial Y\|$, $\|\overline X\|_{\partial X}=\|\overline Y\|_{\partial Y}$ and the corps $\IR_X$ and $\IR_Y$ are isomorphic. This characterization motivates the problem of classifying corps up to an isomorphism. In this section we prove that two finite corps are isomorphic if and only if they have the same cardinality. This characterization follows from two fundamental results:  Wedderburn's Theorem~\ref{t:Wedderburn-Witt} (saying that every finite corps is a field), and Moore's Theorem~\ref{t:Moore} (saying that every two finite fields of the same cardinality are isomorphic).

Wedderburn's Theorem is one of jewels of Mathematics, included to the famous collection of the most brilliant mathmatical proofs ``{\em The proof from THE BOOK}''\footnote{{\bf Proofs from THE BOOK} is a book of mathematical proofs by Martin Aigner and G\"unter M. Ziegler. The book is dedicated to the mathematician Paul Erd\H os, who often referred to ``The Book'' in which God keeps the most elegant proof of each mathematical theorem. During a lecture in 1985, Erd\H os said, ``{\em You don't have to believe in God, but you should believe in The Book}''. It has gone through six editions in English, and has been translated into Persian, French, German, Hungarian, Italian, Japanese, Chinese, Polish, Portuguese, Korean, Turkish, Russian and Spanish.} 

\index[person]{Wedderburn}Wedderburn\footnote{{\bf Joseph Henry Maclagan Wedderburn} (1882 -- 1948) was a Scottish mathematician, who taught at Princeton University for most of his career. A significant algebraist, he proved that a finite division algebra is a field, and part of the Artin–Wedderburn theorem on simple algebras. He also worked on group theory and matrix algebra. 

Joseph Wedderburn was the tenth of fourteen children of Alexander Wedderburn of Pearsie, a physician, and Anne Ogilvie. He was educated at Forfar Academy, then in 1895 his parents sent Joseph and his younger brother Ernest to live in Edinburgh with their paternal uncle, J R Maclagan Wedderburn, allowing them to attend George Watson's College. In 1898 Joseph entered the University of Edinburgh. In 1903, he published his first three papers, worked as an assistant in the Physical Laboratory of the University, obtained an MA degree with First Class Honours in mathematics, and was elected a Fellow of the Royal Society of Edinburgh. Aged only 21 he remains one of the youngest Fellows ever.

He then studied briefly at the University of Leipzig and the University of Berlin, where he met the algebraists Frobenius and Schur. A Carnegie Scholarship allowed him to spend the 1904--1905 academic year at the University of Chicago where he worked with Oswald Veblen, E. H. Moore, and most importantly, Leonard Dickson, who was to become the most important American algebraist of his day.

Returning to Scotland in 1905, Wedderburn worked for four years at the University of Edinburgh as an assistant to George Chrystal, who supervised his D.Sc, awarded in 1908 for a thesis titled ``On Hypercomplex Numbers''. He gained a PhD in algebra from the University of Edinburgh in 1908. From 1906 to 1908, Wedderburn edited the Proceedings of the Edinburgh Mathematical Society. In 1909, he returned to the United States to become a Preceptor in Mathematics at Princeton University; his colleagues included Luther P. Eisenhart, Oswald Veblen, Gilbert Ames Bliss, and George Birkhoff.}
himself gave three proofs of this theorem in 1905, and another proof was given by \index[person]{Dickson}Dickson\footnote{{\bf Leonard Eugene Dickson} (1874 -- 1954) was an American mathematician. He was one of the first American researchers in abstract algebra, in particular the theory of finite fields and classical groups, and is also remembered for a three-volume history of number theory, ``History of the Theory of Numbers.''} 
in the same year. 

More proofs were later given by \index[person]{Artin}Artin\footnote{{\bf Emil Artin} (1898 --1962) was an Austrian mathematician of Armenian descent. Artin was one of the leading mathematicians of the twentieth century. He is best known for his work on algebraic number theory, contributing largely to class field theory and a new construction of $L$-functions. He also contributed to the pure theories of rings, groups and fields.
Along with Emmy Noether, he is considered the founder of modern abstract algebra.},
\index[person]{Zassenhaus}Zassenhaus\footnote{{\bf Hans Julius Zassenhaus} (1912 -- 1991) was a German mathematician, known for work in many parts of abstract algebra, and as a pioneer of computer algebra.}, 
\index[person]{Bourbaki}Bourbaki\footnote{{\bf Nicolas Bourbaki} is the collective pseudonym of a group of mathematicians, predominantly French alumni of the \'Ecole normale sup\'erieure (ENS). Founded in 1934--1935, the Bourbaki group originally intended to prepare a new textbook in analysis. Over time the project became much more ambitious, growing into a large series of textbooks published under the Bourbaki name, meant to treat modern pure mathematics. The series is known collectively as the ``\'El\'ements de math\'ematique'', the group's central work. Topics treated in the series include set theory, abstract algebra, topology, analysis, Lie groups and Lie algebras.}, 
and many others (including the misterous ``unabomber'' 
\index[person]{Kaczynski}Teodor Kaczynski\footnote{{\bf Theodore John Kaczynski} (1942 -- 2023) was an American mathematician and domestic terrorist.  He was a mathematics prodigy, but abandoned his academic career in 1969 to pursue a primitive lifestyle. In 1971, Kaczynski moved to a remote cabin without electricity or running water near Lincoln, Montana, where he lived as a recluse while learning survival skills to become self-sufficient. After witnessing the destruction of the wilderness surrounding his cabin, he concluded that living in nature was becoming impossible and resolved to fight industrialization and its destruction of nature through terrorism. In 1979, Kaczynski became the subject of what was, by the time of his arrest in 1996, the longest and most expensive investigation in the history of the Federal Bureau of Investigation (FBI). The FBI used the case identifier UNABOM (University and Airline Bomber) before his identity was known, resulting in the media naming him the ``Unabomber". In context of Wedderburn's theorem, Kaczynski is known for his paper ``Another Proof of Wedderburn's Theorem'', published in American Mathematical Monthly in 1964.}. 

One proof of Wedderburn's Theorem stands out for its simplicity and elegance. It was found by \index[person]{Witt}Ernst Witt\footnote{{\bf Ernst Witt} (1911 -- 1991) was a German mathematician, one of the leading algebraists of his time. Witt was born on the island of Alsen, then a part of the German Empire. Shortly after his birth, his parents moved the family to China to work as missionaries, and he did not return to Europe until he was nine. 
After his schooling, Witt went to the University of Freiburg and the University of G\"ottingen. He joined the NSDAP (Nazi Party) and was an active party member. Witt was awarded a Ph.D. at the University of G\"ottingen in 1934 with a thesis titled: ``Riemann-Rochscher Satz und Zeta-Funktion im Hyperkomplexen'' that was supervised by Gustav Herglotz, with Emmy Noether suggesting the topic for the doctorate. He qualified to become a lecturer and gave guest lectures in G\"ottingen and Hamburg. He became associated with the team led by Helmut Hasse who led his habilitation. In June 1936, he gave his habilitation lecture. During World War II he joined a group of five mathematicians, recruited by Wilhelm Fenner, and which included Georg Aumann, Alexander Aigner, Oswald Teichm\"uller, Johann Friedrich Schultze and their leader professor Wolfgang Franz, to form the backbone of the new mathematical research department in the late 1930s that would eventually be called: Section IVc of Cipher Department of the High Command of the Wehrmacht. From 1937 until 1979, he taught at the University of Hamburg.
} in 1931.

\begin{theorem}[Wedderburn, 1905; Dickson, 1905; Witt, 1931]\label{t:Wedderburn-Witt} Every finite corps is a field.
\end{theorem}

\begin{proof} The first ingredient of Witt's proof of Wedderburn's theorem comes from a blend of linear algebra and
basic group theory. Let $R$ be a finite corps. For an arbitrary element $s\in R$, let $C_s$ be the set
$\{x\in R : xs = sx\}$ of elements which commute with $s$; $C_s$ is called the
\index{centralizer}\defterm{centralizer} of $s$. Clearly, $C_s$ contains $0$ and $1$ and is a subcorps of $R$. The \index{center}\defterm{center} $Z$ of $R$ is the set of elements which commute with all elements of $R$, thus $Z =
\bigcap_{s\in R}C_s$. In particular, all elements of $Z$ commute, $0$ and $1$
are in $Z$, and so $Z$ is a finite field. Let $q\defeq |Z|$ be the cardinality of the field $Z$. 
We can regard $R$ and $C_s$ as vector spaces over the field $Z$ and deduce that
$|R| = q^n$, where $n$ is the dimension of the vector space $R$ over $Z$, and
similarly $|C_s| = q^{n_s}$ for suitable integers $n_s\ge 1$.
Now let us assume that $R$ is not a field. This means that for some $s\in R$
the centralizer $C_s$ is not all of $R$, or, what is the same, $n_s < n$.
On the multiplicative group $R^*\defeq R\setminus\{0\}$ of the corps $R$ we consider the relation
$$r'\sim r\overset{\mathsf{def}}\Longleftrightarrow  r' \in\bigcup_{x\in R^*}x^{-1}rx.$$
It is easy to check that $\sim$ is an equivalence relation. Let
$A_s := \{x^{-1}sx : x \in R^*\}$
be the equivalence class containing $s$. We note that $|A_s| = 1$ precisely
when $s$ is in the center $Z$. So by our assumption, there are classes $A_s$ with
$|A_s| \ge 2$. Consider now for $s\in R^*$ the map $f_s : x \mapsto x^{-1}sx$ from $R^*$ onto $A_s$. For $x, y\in R^*$ we find
$$x^{-1}sx = y^{-1}sy\;\Leftrightarrow\;(yx^{-1})s = s(yx^{-1})\;\Leftrightarrow\;yx^{-1}\in C^*_s\;\Leftrightarrow\;y\in C^*_sx,$$
wehere $C^*_s \defeq C_s\setminus \{0\}$ and $C^*_sx = \{zx : z \in C^s \}$ has size $|C^*_sx|=|C^*_s|$. Hence any
element $x^{-1}sx$ is the image of precisely $|C_s^*| = q^{n_s}-1$ elements in $R^*$
under the map $f_s$, and we deduce $|R^*| = |A_s|\cdot |C_s |$. In particular, we note
that
$$
\frac{|R^*|}{|C^*_s |}
=\frac{q^n - 1}{q^{n_s}-1}= |A_s|$$
is an {\em integer} for all $s$.
We know that the equivalence classes partition $R^*$. We now group the
central elements $Z^*$ together and denote by $A_1,\dots,A_t$ the equivalence
classes containing more than one element. By our assumption we know
$t \ge 1$. Since $|R^*| = |Z^*| +\sum_{k=1}^t |A_k|$, we have proved the so-called
{\em class formula} 
\begin{equation}\label{eq:Witt1}q^n - 1 = q - 1 +\sum^t_{k=1}\frac{q^n-1}{q^{n_k}-1},
\end{equation}
where we have $1 < \frac{q^n-1}{q^{n_k}-1}\in \IN$ for all $k$. 

\begin{claim}\label{cl:nk|n} For every $k\in\{1,\dots,t\}$ the number $n_k$ divides the number $n$.
\end{claim}

\begin{proof} Given two positive integers $a,b$, we write $a|b$ if $a$ \defterm{divides} $b$, which means that $b=ac$ for some integer number $c\in\IZ$. Write $n = an_k+r$
with $a\in \IN$ and $0 \le r < n_k$. Observe that $$q^n-1=q^{an_k+r}-q^{n_k}+q^{n_k}-1=q^{n_k}(q^{(a-1)n_k+r}-1)+q^{n_k}-1.$$
 Then $q^{n_k}-1 | q^{an_k+r} -1$ implies $q^{n_k}-1 | q^{n_k} (q^{(a-1)n_k+r} -1),$ and thus $q^{n_k}-1 | q^{(a-1)n_k+r} - 1,$ 
since $q^{n_k}$ and $q^{n_k}-1$ are relatively
prime. Continuing in this way we find $q^{n_k}-1 | q^r - 1$ with $0 \le r < n_k$,
which is only possible for $r = 0$, that is, $n_k | n$.
\end{proof}

Now comes the second ingredient: the complex numbers $\IC$. Consider the
polynomial $x^n -1$. Its roots in $\IC$ form a multiplicative subgroup
$$\sqrt[n]{1}\defeq\{z\in\IC:z^n=1\}=\{e^{2\pi i k/n}:k\in\{0,1,\dots,n-1\}\}$$of the field of complex numbers $\IC$. Some of the roots $\lambda\in\sqrt[n]{1}$ satisfy $\lambda^d=1$ for $d < n$; for example, the
root $\lambda = -1$ satisfies $\lambda^2 = 1$. For a root $\lambda$, let $d$ be the smallest positive
exponent with $\lambda^d = 1$, that is, $d$ is the order of $\lambda$ in the group of the roots of unity. Then $d | n$, by Lagrange's Theorem~\ref{t:Lagrange} (saying that the order of any subgroup of a finite group divides the order of the group).

Now we group all roots of order $d$ together and set
$$\phi_d(x)\defeq \prod_{\mbox{\scriptsize$\lambda$ of order $d$}}(x-\lambda).$$
Note that the definition of $\phi_d(x)$ is independent of $n$. Since every root has
some order $d$, we conclude that
\begin{equation}\label{eq:Witt3}x^n -1 =
\prod_{d|n}\phi_d(x).
\end{equation}
Here is the crucial observation: The coefficients of the polynomials $\phi_d(x)$
are integers (that is, $\phi_d(x)\in\IZ[x]$ for all $d$), where in addition the constant
coefficient is either $1$ or $-1$.
Let us carefully verify this claim. For $n = 1$ we have $1$ as the only root,
and so $\phi_1(x)= x -1$. Now we proceed by induction, where we assume
$\phi_d(x)\in\IZ[x]$ for all $d < n$, and that the constant coefficient of $\phi_d(x)$ is $1$ or $-1$. By (\ref{eq:Witt3}),
\begin{equation}\label{eq:Witt4}
x^n-1 = p(x) \phi_n(x)
\end{equation}
where $p(x) =\sum_{j=0}^\ell p_jx^j$, $\phi_n(x) =
\sum^{n-\ell}_{k=0}a_kx^k$, with $p_0 = 1$ or $p_0 = -1$.
Since $-1 = p_0a_0$, we see $a_0\in\{1,-1\}$. Suppose we already know that
$a_0, a_1,\dots, a_{k-1}\in\IZ$. Computing the coefficient of $x^k$ on both sides of (\ref{eq:Witt4})
we find
$$\sum^k_{j=0}p_ja_{k-j} =\sum^k_{j=1}
p_ja_{k-j} + p_0a_k \in\IZ.$$
By assumption, all $a_0,\dots,a_{k-1}$ (and all $p_j$) are in $\IZ$. Thus $p_0a_k$ and hence $a_k$ must also be integers, since $p_0$ is $1$ or $-1$.

We are ready for the coup de gr\^ace. Let $n_k | n$ be one of the numbers
appearing in the class formula (\ref{eq:Witt1}). Then
$$x^n - 1 =\prod_{d|n}\phi_d(x) = (x^{n_k} - 1)\phi_n(x)
\prod_{d | n, d\nmid n_k, d\ne n}\phi_d(x).$$
We conclude that in $\IZ$ we have the divisibility relations
\begin{equation}
\label{eq:Witt5}
\phi_n(q) | q^n-1\quad\mbox{and}\quad\phi_n(q)|\frac{q^n -1}{q^{n_k}-1}.
\end{equation}
Since (\ref{eq:Witt5}) holds for all $k$, we deduce from the class formula (\ref{eq:Witt1}) that $\phi_n(q) | q - 1$, but this cannot be true. Why? We know $\phi_n(x) =\prod(x - \lambda)$ where $\lambda$ runs through all roots of $x^n-1$ of order $n$. Let $\tilde \lambda = a+ib$ be one of those roots. By $n > 1$ (because of $R \ne Z$) we have $\tilde \lambda\ne 1$, which implies that the real part $a$ of $\tilde L$ is smaller than 1.

\begin{picture}(200,100)(-150,-50)
{\linethickness{1pt}
\put(120,0){\color{red}\line(-105,26){105}}
\put(120,0){\color{blue}\line(-1,0){90}}
}
\put(0,0){\circle{60}}
\put(0,0){\line(1,0){30}}
\put(120,0){\circle*{3}}
\put(15,26){\circle*{3}}
\put(14,30){$\tilde \lambda$}
\put(15,0){\circle*{3}}
\put(13,-9){$a$}
\put(0,0){\circle*{3}}
\put(-2,-10){$0$}
\put(123,-2){$q$}
\put(30,0){\circle*{3}}
\put(31,-10){$1$}
\put(60,-10){$|q-1|$}
\put(55,22){\rotatebox{-14}{$|q-\tilde\lambda|$}}
\put(0,0){\line(15,26){15}}
\put(15,0){\line(0,1){26}}
\end{picture}

Taking into account that $|\tilde \lambda|^2 = a^2 + b^2 = 1$, we conclude that $$|q - \tilde \lambda|^2 = |q - a - ib|^2 = (q - a)^2 + b^2
= q^2 - 2aq + a^2 + b^2 = q^2 - 2aq + 1
> q^2 - 2q + 1
= (q - 1)^2,$$ and so $|q - \tilde \lambda| > q -1$ holds for all roots of order $n$. This implies
$$|\phi_n(q)| =\prod_{\lambda}|q -\lambda| > q - 1,$$
which means that $\phi_n(q)$ cannot be a divisor of $q -1$, contradiction and end
of proof.
\end{proof}

Theorem~\ref{t:Wedderburn-Witt} reduces the problem of isomorphic classification of finite corps to the problem of isomorphic classification of finite fields, which was resolved by \index[person]{Moore}E.H.~Moore\footnote{{\bf Eliakim Hastings Moore} (1862 -- 1932),  was an American mathematician.  
Moore first worked in abstract algebra, proving in 1893 the classification of the structure of finite fields (also called Galois fields). Around 1900, he began working on the foundations of geometry. He reformulated Hilbert's axioms for geometry so that points were the only primitive notion, thus turning David Hilbert's primitive lines and planes into defined notions. In 1902, he further showed that one of Hilbert's axioms for geometry was redundant. His work on axiom systems is considered one of the starting points for metamathematics and model theory. After 1906, he turned to the foundations of analysis. The concept of a closure operator first appeared in his 1910 ``{\em Introduction to a form of general analysis}''. He also wrote on algebraic geometry, number theory, and integral equations. 

At Chicago, Moore supervised 31 doctoral dissertations, including those of George Birkhoff, Leonard Dickson, Robert Lee Moore (no relation), and Oswald Veblen. Birkhoff and Veblen went on to lead departments at Harvard and Princeton, respectively. Dickson became the first great American algebraist and number theorist. Robert Moore founded American topology. According to the Mathematics Genealogy Project, as of March 2024, E. H. Moore had 30,199 known ``descendants."

Moore convinced the New York Mathematical Society to change its name to the American Mathematical Society, whose Chicago branch he led. He presided over the AMS, 1901--02, and edited the Transactions of the American Mathematical Society, 1899--1907. He was elected to the National Academy of Sciences, the American Academy of Arts and Sciences, and the American Philosophical Society. He was an Invited Speaker at the International Congress of Mathematicians in 1908 in Rome and in 1912 in Cambridge, England.
The American Mathematical Society established a prize in his honor in 2002.
} in 1893.

The proof of Moore's Theorem uses some standard facts from Group Theory. So, we first recall proofs of those standard facts in order to make our presentation self-contained.
 The first of them is the classical Largangue Theorem on the order of subgroups in finite groups. By the \index{order}\defterm{order} of a finite group $G$ we understand its cardinality $|G|$.

\index[person]{Largange}
\begin{theorem}[Lagrange, 1770/71]\label{t:Lagrange}\index[person]{Lagrange} The order of any subgroup $H$ of a finite group $G$ divides the order of $G$.
\end{theorem}

\begin{proof} Observe that $G/H\defeq\{Hx:x\in G\}$ is a partition of $G$ into pairwise disjoint sets of cardinality $|Hx|=|H|$. Then $|G|=|G/H|\cdot|H|$ and hence the number $|H|$ divides the number $|G|$.
\end{proof}

The \index{order}\defterm{order} of an element $x$ of a finite group $G$ is defined as the order of the cyclic subgroup of $G$, generated by $x$. The order of $x$ is equal to the smallest number $d$ such that $x^d$ is equal to the neutral element $1$ of the group $G$. 
Lagrange Theorem~\ref{t:Lagrange} implies the following corollary.

\begin{corollary}\label{c:Lagrange} The order of every element of a finite group $G$ divides the order of the group $G$.
\end{corollary}

The second standard fact is the following well-known theorem (whose proof is new).

\begin{theorem}\label{t:F*-cyclic} Let $F$ be a field and $F^*\defeq F\setminus\{0\}$ be  its multiplicative group. Every finite subgroup $G$ of $F^*$ is cyclic.
\end{theorem}

\begin{proof} Let $g$ be an element of the largest possible order $n$ in the group $G$. Let $H\defeq\{g^k:0\le k<n\}$ be the cyclic subgroup in $G$ generated by the element $g$. If $H=G$, then the group $G$ is cyclic and we are done. So, assume that $H\ne G$. Observe that every element $x\in H$ satisfies the equation $x^n-1=0$, which has at most $n$ roots in the field $F$. Therefore, $H=\{x\in F:x^n=1\}$. Let $p$ be the smallest positive number for which there exists an element $y\in G\setminus H$ such that $y^p\in H$. The minimality of $p$ implies that $p$ is a prime number. Let $k\in\{0,\dots,n-1\}$ be the smallest possible number such that $y^p=g^k$ for some element $y\in G\setminus H$.
We claim that $k<p$. Assuming that $k\ge p$, we can observe that the element $z=yg^{-1}\in G\setminus H$ has $(yg^{-1})^p=g^{k-p}$, which contradicts the minimality of $k$. 

If $k=0$, then $y^p=g^k=1$ and $y\notin H$ implies that $p$ does not divide $n$. In this case the element $yg\in G$ has order $pn>n$, which contradicts the maximality of $n$. This contradiction shows that $k>0$. Since $k<p$, the element $g^k$ has order $d>\frac{n}p$. Then the element $y$ has order $pd>n$, which contradicts the maximality of $n$. This is a final contradiction showing that the group $G=H$ is cyclic.
\end{proof}

Now we are ready to prove Moore's classification of finite fields. 

\begin{theorem}[Moore, 1893]\label{t:Moore} Two finite fields are isomorphic if and only if they have the same cardinality.
\end{theorem}

\begin{proof} The ``only if'' part of the theorem is trivial. To prove the ``if'' part, assume that  $X,Y$ are two finite fields of the same cardinality. Then $|X|=|Y|=p^n$ for some prime number $p$ and some $n\in\IN$. The minimal fields of the fields $X,Y$ both are isomorphic to the field $\IF_p=\IZ/p\IZ$. So, both fields $X$ and $Y$ are vectors spaces of dimension $n$ over the field $\IF_p$. If $n=1$, then the fields $X,Y$ are isomorphic to the field $\IF_p$ and we are done. So, assume that $n\ge 2$. A non-constant polynomial $P\in \IF_p[x]$ is \index{irreducible polynomial}\defterm{irreducible} if for any polynomias $A,B\in \IF_p[x]$ with $P=A\cdot B$ one of the polynomials has degree zero (i.e., is constant). Write the polynomial $P(x)=x^{p^n}-x$ as the product 
$$P=P_1\cdots P_m$$ of irreducible polynomials $P_1,\dots,P_m\in\IF_p[x]$. By the Lagrange Theorem, every point of $X$ is a root of the polynomial $P(x)=x^{p^n}-x$. So, this polynomial has exactly $p^n$ roots in the field $X$. Since $p^n=\deg(P)=\sum_{i=1}^m\deg(P_i)$, every irreducible polynomial $P_i$ has exactly $\deg(P_i)$ roots in the field $X$. The same argument show that every polynomial $P_i$ has exactly $\deg(P_i)$ roots in the field $Y$. 

By Theorem~\ref{t:F*-cyclic}, the multplicative group $X^*\defeq X\setminus\{0\}$ of the field $X$ is cyclic. So, we can fix a generator $g$ of $X^*$. By Corollary~\ref{c:Lagrange}, $0=g^{p^n}-g=P_1(g)\cdots P_m(g)$ and hence $P_i(g)=0$ for some $i$.

\begin{claim}\label{cl:Moore1} If $Q(g)=0$ for some nonzero polynomial $Q\in\IF_p[x]$, then $\deg(Q)\ge n$.
\end{claim}

\begin{proof} Assuming that $Q(g)=0$ for some polynomial $Q\in\IF_p[x]$ of degree $\deg(Q)<n$, we conclude that $g^k=a_0+a_1g+\dots+a_{k-1}g^{k-1}$ for some positive $k<n$. Then the set $\{x_0+x_1g+\dots+x_{k-1}g^{k-1}:x_0,\dots,x_{k-1}\in \IF_p\}$ is a proper subring of the ring $X$, which is not possible because this subring contains the generator $g$ of the multiplicative group $X^*$ and hence contains all elements of $X$.
\end{proof}

Claim~\ref{cl:Moore1} implies that the polynomial $P_i$ has degree $\deg(P_i)\ge n$. 
We claim that $\deg(P_i)=n$. Since the vector space $X$ over the field $\IF_p$ has dimension $n$, the set $\{1,g,g^2,\dots,g^n\}$ is not linearly independent and hence there exists $k\le n$ such that $g^{k}=a_0+a_1g+\dots+a_{k-1}g^{k-1}$ for some $a_0,\dots,a_{k-1}\in \IF_p$. Then the set $\{x_0+x_1g+\dots+x_{k-1}g^{k-1}:x_0,\dots,x_k\in \IF_p\}$ is a subring of the ring $X$. Since $g$ is a generator of the multiplicative group $X^*$, this subring equals $X$, which implies that $k=n$ and hence $A(g)=0$ for the polynomial $A(x)=a_0+a_1x+\dots+a_{k-1}x^{k-1}-x^{k}$ of degree $n$. Write the polynomial $P_i$ as $P_i=A\cdot Q+R$ for some polynomials $Q$ and $R$ with $\deg(R)<\deg(Q)=n$. Since $R(g)=P_i(g)-A(g)\cdot Q(g)=0$, Claim~\ref{cl:Moore1} ensures that the polynomial $R$ is zero. Then $P_i=A\cdot Q$ and the irreducibility of $P_i$ implies that the polynomial $A$ has degree zero and hence $\deg(P_i)=\deg(Q)=n$. 

 The Largange Theorem implies that the polynomial $y^{p^n}-y=P(y)=P_1(y)\cdots P_m(y)$ has exactly $p^n=|Y|$ roots in $Y$. Then every polynomial $P_j$ has exactly $\deg(P_j)$ roots in $Y$. In particular, the polynomial $P_i$ has $n$ roots in the field $Y$. Let $\tilde g$ be any root of the polynomial $P_i$ in the field $Y$. Consider the function $f:X\to Y$ that assigns to every element $x\in X$ the element $y\defeq x_0+x_1\tilde g+\cdots + x_{n-1}\tilde g^{n-1}$, where $x_0,\dots,x_{n-1}\in\IF_p$ are unique elements such that $x=x_0+x_1g+\cdots +x_{n-1}g^{n-1}$. Taking into account that $P_i(g)=0=P_i(\tilde g)$, one can check that the function $f:X\to Y$ is a required isomorphism of the fields $X$ and $Y$. 
\end{proof}


\section{Isomorphic classification of finite Desarguesian proaffine spaces}

In this section we shall present an isomorphic classification of finite Desarguesian proaffine spaces. We start with calculation of the cardinality of such spaces.

\begin{proposition}\label{p:card-compreg} Every $3$-long finite completely regular liner $X$ has cardinality $$|X|=\frac{(|\overline X|_2-1)^{\|X\|}-(|\overline X|_2-1)^{\|\partial X\|}}{|\overline X|_2-2}.$$
\end{proposition}

\begin{proof} Since the liner $X$ is $3$-long and completely regular, its spread completion $\overline X$ is projective and $3$-long, by Corollary~\ref{c:spread-3long}. Therefore, the projective liner $\overline X$ is a projective completion of the liner $X$. By Corollary~\ref{c:Avogadro-projective}, the $3$-long projective liner $\overline X$ is $2$-balanced, so the number $|\overline X|_2\ge 3$ is well-defined. 
By Corollary~\ref{c:procompletion-rank}, $\|\overline X\|=\|X\|$. By Corollary~\ref{c:projective-order-n}, the projective liner $\overline X$ has cardinality $$\overline X=\sum_{k=0}^{\|X\|-1}(|\overline X|_2-1)^k=\frac{(|\overline X|_2-1)^{\|X\|}-1}{(|\overline X|_2-1)-1}=\frac{(|\overline X|_2-1)^{\|X\|}-1}{|\overline X|_2-2}.$$

Taking into account that the horizon $\partial X=\overline X\setminus X$ of $X$ is flat in the $2$-balanced projective liner $\overline X$, we conclude that $\partial X$ is a $2$-balanced projective liner with $|\partial X|_2=|\overline X|_2$ and $$|\partial X|=\frac{(|\partial X|_2-1)^{\|\partial X\|}-1}{|\partial X|_2-2}=\frac{(|\overline X|_2-1)^{\|\partial X\|}-1}{|\overline X|_2-2}.$$ Then
$$|X|=|\overline X|-|\partial X|=\frac{(|\overline X|_2-1)^{\|X\|}-1}{|\overline X|_2-2}-\frac{(|\overline X|_2-1)^{\|\partial X\|}-1}{|\overline X|_2-2}=\frac{(|\overline X|_2-1)^{\|X\|}-(|\overline X|_2-1)^{\|\partial X\|}}{|\overline X|_2-2}.$$
\end{proof}

\begin{proposition}\label{p:card-Dproaff} Every $3$-long finite Desarguesian proaffine regular space $X$ has cardinality $$|X|=\frac{|\IR_X|^{\|X\|}-|\IR_X|^{\|\partial X\|}}{|\IR_X|-1}.$$
\end{proposition}

\begin{proof} By Theorem~\ref{t:Desargues-completion}, the spread completion $\overline X$ of the $3$-long Desarguesian proaffine regular space $X$ is a $3$-long Desarguesian projective space. By Corollary~\ref{c:procompletion-rank}, $\|\overline X\|=\|X\|$ and by Theorem~\ref{t:Dp=PRX}, the projective space $\overline X$ is isomorphic to the projective space $\IP\IR_X^{\oplus\|X\|}$. Then $|\overline X|_2=|\IP\IR_X^{\oplus\|X\|}|_2=|\IR_X|+1$. By Proposition\ref{p:card-compreg},
$$|X|=\frac{(|\overline X|_2-1)^{\|X\|}-(|\overline X|_2-1)^{\|\partial X\|}}{|\overline X|_2-2}=\frac{|\IR_X|^{\|X\|}-|\IR_X|^{\|\partial X\|}}{|\IR_X|-1}.$$
\end{proof}

\begin{theorem}\label{t:finite-proj-isomorphic<=>} Two $3$-long finite Desarguesian proaffine regular liners $X,Y$ are isomorphic if and only if $|X|=|Y|$ and $\|X\|=\|Y\|$. 
\end{theorem}

\begin{proof} The ``only if'' part is obvious. To prove the ``if'' part, assume that  $|X|=|Y|$ and $\|X\|=\|Y\|$. If $\|X\|=\|Y\|\le 2$, then any bijection $F:X\to Y$ is an isomorphism of the liners $X,Y$. So, we assume that $\|X\|=\|Y\|\ge 3$. In this case, the liners $X,Y$ are Desarguesian proaffine spaces.  By Proposition~\ref{p:card-Dproaff},
$$\frac{|\IR_X|^{\|X\|}-|\IR_X|^{\|\partial X\|}}{|\IR_X|-1}=|X|=|Y|=\frac{|\IR_Y|^{\|Y\|}-|\IR_Y|^{\|\partial Y\|}}{|\IR_Y|-1}.$$
First we show that $|\IR_X|=|\IR_Y|$. To derive a contradiction, assume that $|\IR_X|\ne|\IR_Y|$. We lose no generality assuming that $|\IR_X|<|\IR_Y|$. Let $r\defeq\|X\|=\|Y\|\ge 3$ and $\ell\defeq|\IR_X|\ge 2$. Then $\ell+1=|\IR_X|+1\le|\IR_Y|$ and  
$$
\begin{aligned}
(\ell+1)^{r-1}&\le |\IR_Y|^{r-1}=\frac{|\IR_Y|^{r}-|\IR_Y|^{r-1}}{|\IR_Y|-1}\le  \frac{|\IR_Y|^{r}-|\IR_Y|^{\|\partial Y\|}}{|\IR_Y|-1}=|Y|\\
&=|X|=\frac{|\IR_X|^{\|X\|}-|\IR_X|^{\|\partial X\|}}{|\IR_X|-1}\le\frac{|\IR_X|^{\|X\|}-1}{|\IR_X|-1}= \frac{\ell^{r}-1}{\ell-1}.
\end{aligned}
$$
On the other hand, $$(\ell+1)^{r-1}(\ell-1)=(\ell+1)^{r-3}(\ell^3+\ell(\ell-1)-1)> (\ell+1)^{r-3}\ell^3\ge \ell^r>\ell^r-1,$$
which contradicts $(\ell+1)^{r-1}\le \frac{\ell^r-1}{\ell-1}$.
This contradiction shows that $|\IR_X|=|\IR_Y|$. By Theorems~\ref{t:Wedderburn-Witt} and \ref{t:Moore}, the finite corps $\IR_X$ and $\IR_X$ are isomorphic.
The equalitities $|\IR_X|=|\IR_Y|$, $\|X\|=\|Y\|$ and
$$\frac{|\IR_X|^{\|X\|}-|\IR_X|^{\|\partial X\|}}{|\IR_X|-1}=|X|=|Y|=\frac{|\IR_Y|^{\|Y\|}-|\IR_Y|^{\|\partial Y\|}}{|\IR_Y|-1}$$
imply $\|\partial X\|=\|\partial Y\|$. By Corollary~\ref{c:rank+}, $$\|\partial X\|+\|\overline X\|_{\partial X}=\|X\|=\|Y\|=\|\partial Y|+\|\overline Y\|_{\partial Y}$$and hence $\|\overline X\|_{\partial X}=\|\overline Y\|_{\partial Y}$. By Theorem~\ref{t:Dproaff-iso<=>}, the Desarguesian proaffine spaces $X,Y$ are isomorphic.
\end{proof}

\begin{exercise} Find two non-isomorphic finite Desarguesian proaffine spaces $X,Y$ of the same cardinality. Show that such spaces necessarily are affine.
\end{exercise}

\section{Projective repers and frames in Desarguesian projective spaces}

In this section we establish some better homogeneity properties of Desarguesian projective spaces, using the notion of a projective reper.

\begin{definition} A \index{projective reper}\defterm{projective reper} in a projective space $X$ is a function $r\subseteq X\times X$ whose domain $\dom[r]$ is a maximal independent set in $X$ containing a unique point $o=r(o)$ such that $r(x)\in\Aline ox\setminus\{o,x\}$ for every $x\in \dom[r]\setminus\{o\}$. The point $o$ is called the \index{projective reper!origin of}\index{origin}\defterm{origin} of the projective reper $r$. 
\end{definition} 

\begin{picture}(100,90)(-200,-15)
\put(0,0){\line(1,1){60}}
\put(0,0){\line(-1,1){60}}

\put(0,0){\color{red}\circle*{3}}
\put(-2,-7){$o$}
\put(30,30){\color{blue}\circle*{3}}
\put(32,23){$a$}
\put(60,60){\color{violet}\circle*{3}}
\put(62,53){$r(a)$}
\put(-30,30){\color{blue}\circle*{3}}
\put(-38,23){$b$}
\put(-60,60){\color{violet}\circle*{3}}
\put(-78,53){$r(b)$}
\end{picture}

\begin{exercise} Show that for every isomorphism $A:X\to Y$ between projective spaces and every projective reper $r$ in $X$, the function $A\circ r\circ A^{-1}$ is a projective reper in the projective space $Y$.
\end{exercise}

\begin{proposition}\label{p:proj-reper} Every projective reper $r$ in a projective space $X$ has the following properties:
\begin{enumerate}
\item the function $r$ is injective;
\item the intersection $\dom[r]\cap\rng[r]$ is a singleton $\{o\}$ containing the origin $o$ of the projective reper $r$;
\item the set $\rng[r]$ is maximal independent in $X$;
\item the function $r^{-1}\defeq\{(y,x):(x,y)\in r\}$ is a projective reper in $X$;
\item the flat $H\defeq \overline{\dom[r]\setminus\rng[r]}$ is a hyperplane in $X$ and the subliner $X\setminus H$ is affine and regular.
\end{enumerate}
\end{proposition}

\begin{proof} Let $o=r(o)$ be the origin of the reper $r$.
\smallskip

1. Assuming that the function $r$ is not injective, we can find two distinct points $x,x'\in\dom[r]$ such that $r(x)=r(x')$. If $x'=o$, then $x\ne o$ and hence $o=r(o)=r(x')=r(x)\in\Aline ox\setminus\{o,x\}$, which is a contradiction showing that $x'\ne o$. By analogy we can prove that $x\ne o$. Then the point $b=r(x)=r(x')$ belongs to the intersection $(\Aline ox\setminus\{o,x\})\cap(\Aline o{x'}\setminus\{o,x'\})=\Aline ox\cap\Aline o{x'}\setminus\{o,x,x'\}$, which implies that $\{o,x,x'\}\subseteq \Aline ob$ and $x'\in\Aline o{x'}=\Aline ox$. But this contradicts the independence of the set $\dom[r]$.
\smallskip

2. It is clear that $o=r(o)\in \dom[r]\cap \rng[r]$. Assuming that $\dom[r]\cap\rng[r]\ne\{o\}$, we can find a point $y\in\dom[r]\cap\rng[r]\setminus\{o\}$. Since $y\in\rng[r]$, there exists a point $x\in\rng[r]$ such that $y=r(x)$. Assuming that $x=o$, we conclude that $y=r(x)=r(o)=o$, which contradicts the choice of the point $y$.
Then $x\ne o$ and hence $y\in\Aline ox\setminus\{o,x\}$ and $y\in \Aline ox\subseteq \overline{\dom[r]\setminus\{y\}}$, which contradicts the independence of the set $\dom[r]$.
\smallskip

3. Assuming that the set $\rng[r]$ is not independent, we can find a point $b\in \rng[r]$ such that $b\in\overline{\rng[r]\setminus\{b\}}$. Since $b\in\rng[r]$, there exists a point $a\in\dom[r]$ such that $b=r(a)$. If $a\ne b$, then $a\ne o\ne b$ and then $o=r(o)\in \rng[r]\setminus \{b\}$. Then $\rng[r]\setminus\{b\}\subseteq \bigcup_{x\in \dom[r]\setminus\{a\}}\Aline ox\subseteq \overline{\dom[r]\setminus\{a\}}$ and $a\in\Aline ob\subseteq \overline {\rng[r]\setminus\{b\}}\subseteq \overline{\dom[r]\setminus\{a\}}$, which contradicts the idependence of the set $\dom[r]$. This contradiction shows that $a=b=o$. Then $o\in\overline{\rng[r]\setminus\{o\}}$ and hence $o\in\overline{F}$ for some finite set $F\subseteq \rng[r]\setminus\{o\}$ of the smallest possible cardinality. In this case, the set $F$ is independent and $\|\overline F\|=|F|$, by Corollary~\ref{c:Max=dim}. For every point $y\in F\subseteq\rng[r]$, and the point $x\defeq r^{-1}(y)$, we have $y\in\Aline ox\setminus \{o,x\}$, which implies $x\in\Aline oy\subseteq \overline{\{o\}\cup F}=\overline F$. Then $r^{-1}[F]\subseteq \dom[r]\cap \overline F$ is an independent set of cardinality $|r^{-1}[F]|=|F|$, by the injectivity of the function $r$. The rankedness of the projective space $X$ ensures that $o\in \overline F=\overline{r^{-1}[F]}\subseteq\overline{\dom[r]\setminus \{o\}}$, which contradicts the independence of the set $\dom[r]$. This contradiction shows that the set $\rng[r]$ is independent. 

Since $r$ is a projective reper, for every point $x\in \dom[r]$ and its image $y\defeq r(x)$, we have $x\in\Aline oy \subseteq\overline{\rng[r]}$ and hence $\dom[r]\subseteq\overline{\rng[r]}$. The maximality of the independent set $\dom[r]$ implies $X=\overline{\dom[r]}\subseteq\overline{\rng[r]}\subseteq X$, which means that the independent set $\rng[r]$ is maximal independent in $X$.
\smallskip

4. The injectivity of the function $r$ (proved in the first statement) implies that the relation $r^{-1}\defeq\{(y,x):(x,y)\in r\}$ is a function. The equality $o=r(o)$ implies $r^{-1}(o)=o$. The preceding statement ensures that the set $\dom[r^{-1}]=\rng[r]$ is maximal independent in $X$. Finally, for every point $y\in \dom[r^{-1}]\setminus\{o\}=\rng[r]\setminus\{o\}$ and the point $x\defeq r^{-1}(y)$ we have $y\in\Aline ox\setminus\{o,x\}$, which implies $x\in\Aline oy\setminus\{o,y\}$, witnessing that the function $r^{-1}$ is a projective reper in $X$.
\smallskip

5. It follows from $\dom[r]\cap\rng[r]=\{o\}$ that $\dom[r]\setminus\rng[r]=\dom[r]\setminus\{o\}$. Since the set $\dom[r]$ is independent, the flat $$H\defeq \overline{\dom[r]\setminus\{o\}}=\overline{\dom[r]\setminus\rng[r]}$$ does not contain the point $o$. The maximal independence of $\dom[r]$ ensures $X=\overline{\dom[r]}=\overline{H\cup\{o\}}$. For every point $x\in X\setminus H=\overline{H\cup\{o\}}\setminus H$, the Exchange Property of the projective space $X$ ensures $o\in\overline{H\cup\{x\}}$ and hence $X=\overline{H\cup\{o\}}\subseteq \overline{H\cup\{x\}}$, which means that $H=\overline{\dom[r]\setminus\rng[r]}$ is a hyperplane in $X$. By Proposition~\ref{p:projective-minus-hyperplane}, the subliner $X\setminus H$ of $X$ is affine and regular. 
\end{proof}

\begin{theorem}\label{t:Des-iso-reper} Let $X,Y$ be two Desarguesian projective spaces and $r,r'$ be two projective repers in $X,Y$, respectively. Let $I:\IR_A\to\IR_B$ be any isomorphism between scalar corps of the affine liners $A\defeq X\setminus\overline{\dom[r]\setminus\rng[r]}$ and $B\defeq Y\setminus\overline{\dom[r']\setminus\rng[r']}$. For every bijective function $\psi:\dom[r]\cup\rng[r]\to\dom[r']\cup\rng[r']$ with $\psi\circ r=r'\circ \psi$, there exists a unique isomorphism $\Phi:X\to Y$ such that $\psi\subseteq\Phi$,  and the restriction $\phi\defeq\Phi{\restriction}_A$ is an isomorphism between the affine spaces $A$ and $B$ such that $\dddot\phi=I$.
\end{theorem} 

\begin{proof} The equality $\psi\circ r=r'\circ \psi$ implies that for every point  $e\in\dom[r]$, the point $\psi(e)$ belongs to $\dom[r']$.
Consequently, $\psi{\restriction}_{\dom[r]}$ is a bijection between the maximal independent sets $\dom[r]$ and $\dom[r']$ in the projective liners $X$ and $Y$, respectively. Corollary~\ref{c:Max=dim} implies that $$\|X\|=|\dom[r]|=|\dom[r']|=\|Y\|.$$ 
Let $o=r(o)$ and $o'=r(o')$ be the origins of the projective repers $r$ and $r'$, respectively. 

By Proposition~\ref{p:proj-reper}, the set $\rng[r]$ is maximal independent in $X$, the flat hull $H\defeq\overline{\rng[r]\setminus\dom[r]}$ is a hyperplane in $X$ whose complement $A\defeq X\setminus H$ is an affine regular liner. By Theorem~\ref{t:pD-minus-flat}, the affine liner $A$ is Desarguesian.

By analogy we can show that $H'\defeq\overline{\rng[r']\setminus\dom[r']}$ is a hyperplane in $X$ and $B\defeq Y\setminus H'$ is a Desarguesian affine regular liner. 

If the projective space $X$ is Steiner, then $\IR_A=\{0,1\}$ and hence $\IR_B=I[\IR_A]=\{0,1\}$, which implies that all lines in the Desarguesian affine liner $B$ have cardinality $2$ and the projective space $Y$ is Steiner. By (the proof of) Theorem~\ref{t:Des-projspace<=>}, the midpoint operation $\circ$ on the Steiner liner $X$ induces a group operation $+:A\times A\to A$,
$$x+y=\begin{cases} x&\mbox{if $y=o$};\\
y&\mbox{if $x=o$};\\
(x\circ o)\circ y&\mbox{if $x\ne o\ne y$};
\end{cases}
$$
turning $A$ into an $\IF_2$-module over the two-element field $\IF_2$. The point $o$ is the zero point of the $\IF_2$-module $A$. The maximal independence of the set $\dom[r]\subseteq A$ implies that $\dom[r]\setminus\{o\}$ is an $\IF_2$-basis of the $\IF_2$-module $A$. 

By analogy, the affine liner $B$ carries a canonical structure of an $\IF_2$-module such that $o'$ is the zero element of $B$ and $\dom[r']\setminus\{o'\}$ is the $\IF_2$-basis of the $\IF_2$-module $B$. Then the bijective function $\psi{\restriction}_{\dom[r]}:\dom[r]\to\dom[r']$ extends to a unique $\IF_2$-linear isomorphism $\phi:A\to B$ of the $\IF_2$-modules $A$ and $B$.  The $\IF_2$-linear isomorphism $\phi$ extends to a unique isomorphism $\Phi:X\to Y$ of the projective spaces $X,Y$. The isomorphism $\Phi$ assigns to every point $x\in H$ the point $o'\circ \phi(x\circ o)\in H'$. The isomorphism $\Phi$ has the required property:
$\psi\subseteq\Phi$ and $\dddot\phi=I$ is the unique identity isomorphism of the scalar corps $\IR_A=\{0,1\}=\IR_B$.
\smallskip

Next, assume that the projective spaces $X,Y$ are not Steiner and hence they are $4$-long. Then the affine regular liners $A,B$ are $3$-longs and hence they are affine spaces. By Lemma~\ref{l:trace-flat}, the maximal independence of the set $\dom[r]\subseteq A$ in $X$ implies the maximal independence of $\dom[r]$ in the affine space $A$. By analogy, the maximal independence of the set $\dom[r']\subseteq B$ in $Y$ implies the maximal independence of $\dom[r']$ in the affine space $B$. By Theorem~\ref{t:extension-iso}, there exists a unique isomorphism $\phi:A\to B$ such that $\phi(x)=\psi(x)$ for all $x\in\dom[r]$ and $\dddot\phi=I$. By Theorem~\ref{t:extend-isomorphism-to-completions}, the isomorphism $\phi:A\to B$ extends to a unique isomorphism $\Phi:X\to Y$. For every point $e\in \dom[r]\setminus\{o\}$ and its image $e'\defeq \psi(e)=\phi(e)=\Phi(e)$, we have $\{r(e)\}=\Aline oe\cap H$ and hence $$\Phi(r(e))=\Phi[\Aline oe]\cap\Phi[H]=\Aline{o'}{e'}\cap H'=\{r'(e')\}=\{r'(\psi(e))\}=\{\psi(r(e))\}.$$Therefore, $\Phi(x)=\psi(x)$ for all $x\in\dom[r]\cup\rng[r]$. The choice of the isomorphism $\phi=\Phi{\restriction}_A$ ensures that $\dddot\phi=I$.
\end{proof}

\begin{corollary}\label{c:Des-iso-reper} Let $r,r'$ be two projective repers in a Desarguesian projective space $X$. Every bijective function $\varphi:\dom[r]\cup\rng[r]\to\dom[r']\cup\rng[r']$ with $\varphi\circ r=r'\circ \varphi$ extends to an automorphism $\Phi:X\to X$.
\end{corollary} 

\begin{proof} Let $o$ and $o'$ be the origins of the repers $r$ and $r'$, respectively. By Proposition~\ref{p:proj-reper}, the flats $H\defeq\overline{\dom[r]\setminus\{o\}}$ and  $H'\defeq\overline{\dom[r']\setminus\{o'\}}$ are hyperplanes in $X$ and the subliners $A\defeq X\setminus H$ and $B\defeq X\setminus H'$ of $X$ are affine and regular. By Theorem~\ref{t:Des-autoextend}, there exists an automorphism $\Psi:X\to X$ such that $\Psi(x)=\varphi(x)$ for all $x\in\dom[r]$, which implies $\Psi[H]=H'$ and $\Psi[A]=B$. Then the affine liner $A,B$ and isomorphic and their scalar corps $\IR_A$ and $\IR_B$ are isomorphic. So, we can find an isomorphism $I:\IR_A\to\IR_B$. Now we can apply Theorem~\ref{t:Des-iso-reper} and find an automorphism $\Phi:X\to X$ such that $\Phi(x)=\varphi(x)$ for all $x\in\dom[r]\cup\rng[r]$ and moreover, the restriction $\phi\defeq\Phi{\restriction}_A$ is an isomorphism between the affine liners $A,B$ such that $\dddot\phi=I$.
\end{proof}

\begin{definition} Let $X$ be a projective space of finite rank $\|X\|$. A subset $M\subseteq X$ is called a \index{projective frame}\defterm{projective frame} in $X$ if $|M|=\|X\|+1$ and for every $x\in M$, the set $M\setminus\{x\}$ is independent in $X$.
\end{definition}



\begin{proposition}\label{p:reper-for-frame} Let $X$ be a projective space of finite rank and $M$ be a projective frame in $X$. For every distinct points $o,u\in M$ there exists a projective reper $r$ in $X$ having the following properties: 
\begin{enumerate}
\item $\dom[r]=M\setminus\{u\}$;
\item $r(o)=o$;
\item $\{r(e)\}=\Aline oe\cap\overline{M\setminus\{o,e\}}$ for every point $e\in M\setminus\{o,u\}$;
\item $\overline{M\setminus\{o,e\}}=\overline{\{r(e)\}\cup (M\setminus\{o,u,e\})}$ for every point $e\in M\setminus\{o,u\}$; 
\item $\displaystyle\{u\}=\bigcap_{e\in M\setminus\{o,u\}}\overline{M\setminus\{o,e\}}=\bigcap_{e\in M\setminus\{o,u\}}\overline{\{r(e)\}\cup (M\setminus\{o,u,e\})}$.
\end{enumerate}
The projective reper $r$ is uniquely determined by the properties \textup{(1), (2), (3)}.
\end{proposition}

\begin{proof} Since $M$ is a projective frame in $X$, for every $e\in M\setminus\{o,u\}$, the flat $H_e\defeq\overline{M\setminus\{o,e\}}$ has rank $\|H_e\|=|M|-2=\|X\|-1$ and hence  $H_e$ is a hyperplane in $X$. The independence of the sets $M\setminus\{o\}$ and $M\setminus\{e\}$ implies that $o,e\notin H_e$. By Corollary~\ref{c:line-meets-hyperplane}, the line $\Aline oe$ has a unique common point $r_e$ with the hyperplane $H_e$. Then $$r\defeq\{(o,o)\}\cup\big\{(e,r_e):e\in M\setminus\{o,u\}\big\}$$is a  projective reper in $X$ having the properties (1), (2), (3), which determine the reper $r$ uniquely.
\smallskip

It remains to prove that the projective reper $r$ has the properties (4) and (5).
\smallskip

Assuming that the property (4) fails, one can find a point  $e\in M\setminus\{o,u\}$ such that $$\overline{\{r(e)\}\cup(M\setminus\{o,u,e\})}\ne\overline{M\setminus\{o,e\}}.$$ The choice of the point $r(e)\in \Aline oe\cap\overline{M\setminus\{o,e\}}$ ensures that  
$$\overline{M\setminus\{o,u,e\}}\subseteq \overline{\{r(e)\}\cup (M\setminus\{o,u,e\})}\subsetneq\overline{M\setminus\{o,e\}}$$
and hence $\overline{M\setminus\{o,u,e\}}=\overline{\{r(e)\}\cup (M\setminus\{o,u,e\})}$, by the rankedness of the projective space $X$. Then $e\in \overline{\{o,r(e)\}}\subseteq\overline{\{o\}\cup (M\setminus\{o,u,e\})}=\overline{M\setminus\{u,e\}}$, which contradicts the independence of the set $M\setminus\{u\}$. 
This contradiction completes the proof of the condition (4).
\smallskip 

Since the set $M\setminus\{o\}$ is independent, Proposition~\ref{p:inter-indep-flats} and condition (4) ensure that
$$\{u\}=\overline{(M\setminus\{o\})\setminus(M\setminus\{u,o\})}=\bigcap_{e\in M\setminus\{o,u\}}\overline{M\setminus\{o,e\}}=\bigcap_{e\in M\setminus\{o,u\}}\overline{\{r(e)\}\cup (M\setminus\{o,u,e\})},$$witnessing that the condition (5) holds.
\end{proof}

\begin{theorem}\label{t:Papp-iso-frame} Let $X$ be a Desarguesian projective space of finite rank. Every bijection $\varphi:M\to M'$ between projective frames $M,M'$ in $X$ extends to an automorphism $\Phi:X\to X$ of $X$.
\end{theorem}

\begin{proof} Fix any distinct points $o,u\in M$ and consider the points $o'\defeq\varphi(o)$ and $u'\defeq\varphi(u)$ in the projective frame $M'$. By Proposition~\ref{p:reper-for-frame}, there exists a projective reper $r$ in the projective space $X$ such that $\dom[r]=M\setminus\{u\}$, $r(o)=o$, and $\{r(e)\}=\Aline oe\cap\overline{M\setminus\{o,e\}}$ for every $e\in\dom[r]\setminus\{o\}$. By analogy, there exists a projective reper $r'$ in the projective space $X$ such that $\dom[r']=M'\setminus\{u'\}$, $r(o')=o'$, and $\{r'(e)\}=\Aline {o'}{e'}\cap\overline{M'\setminus\{o',e'\}}$ for every $e'\in\dom[r']\setminus\{o'\}$. Consider the bijection $\psi:\dom[r]\cup\rng[r]\to\dom[r']\cup\rng[r']$ defined by
$$\psi(x)=\begin{cases}\varphi(x)&\mbox{if $x\in \dom[r]$};\\
r'\varphi r^{-1}(x)&\mbox{if $x\in\rng[r]\setminus\{o\}$}
\end{cases}
$$
The definition of the function $\psi$ ensures that $\psi\circ r=r'\circ \psi$. By Corollary~\ref{c:Des-iso-reper}, there exists an automorphism $\Phi:X\to X$ such that $\Phi(x)=\psi(x)$ for all $x\in\dom[r]\cup\rng[x]$. Then also $\Phi(x)=\psi(x)=\varphi(x)$ for all $x\in\dom[r]=M\setminus\{u\}$. It remains to show that $\Phi(u)=u'$. Applying Proposition~\ref{p:reper-for-frame}(5) to the projective repers $r,r'$, we conclude that
$$\{u\}=\bigcap_{e\in M\setminus\{o,u\}}\overline{\{r(e)\cup(M\setminus\{o,u,e\})}\quad\mbox{and}\quad \{u'\}=\bigcap_{e'\in M'\setminus\{o,u\}}\overline{\{r'(e')\cup(M'\setminus\{o',u',e'\}}),$$ which impy the equality
$$\{\Phi(u)\}=\bigcap_{e\in M\setminus\{o,u\}}\Phi[\overline{\{r(e)\}\cup(M\setminus\{o,u,e\})}]=\bigcap_{e'\in M'\setminus\{o',u'\}}\overline{\{r'(e')\cup(M'\setminus\{o',u',e'\})}=\{u'\}.$$
Therefore, $\Phi:X\to X$ is an automorphism of the Desarguesian projective space such that $\Phi{\restriction}_M=\varphi$.
\end{proof}

\section{Homogeneity of Desarguesian proaffine spaces}

We recall that a liner $X$ is called a {\em space} if $X$ is $3$-long, regular and has rank $\|X\|\ge 3$.

\begin{theorem}\label{t:D=>point-homogeneous} For any points $x,y$ of a Desarguesian proaffine space $X$, there exists an automorphism $F:X\to X$ such that $F(x)=y$.
\end{theorem}

\begin{proof}  By Corollary~\ref{c:Desarg=>compreg}, the Desarguesian proaffine space $X$ is completely regular, and by Theorem~\ref{t:Desargues-completion}, its spread completion $\overline X$ is a Desarguesian projective space. Choose any maximal independent set $M_\partial$ in the horizon $\partial X\defeq \overline X\setminus X$. Given two points $x,y\in X$, choose any maximal $\partial X$-independent sets $M_x,M_y$ in $\overline X$ such that $x\in M_x$ and $y\in M_y$. By Proposition~\ref{p:union-of-independent},  the sets $M_\partial \cup M_x$ and $M_\partial \cup M_y$ are maximal and independent in $\overline X$. By Corollary~\ref{c:Max=dim}, $|M_x\setminus M_\partial|=\|X\|_{\partial X}=|M_y\setminus M_\partial|$, so there exists a bijective function $\varphi:M_\partial\cup M_x\to M_\partial \cup M_y$ such that $\varphi(x)=y$ and $f(z)=z$ for all $z\in M_\partial$.  By Theorem~\ref{t:Des-autoextend}, there exists an automorphism $\Phi:\overline X\to\overline X$ such that $\Phi(z)=\varphi(z)$ for all $z\in \partial M\cup M_x$. Then $\Phi[\partial X]=\Phi[\overline{M_\partial}]=\overline{\Phi[M_\partial]}=\overline{M_\partial}=\partial X$ and hence $F\defeq\Phi{\restriction}_X$ is an automorphism of the liner $X$ such that $F(x)=y$.
\end{proof}

\begin{definition} A liner $X$ is called
\begin{itemize}
\item \index{point-homogeneous liner}\index{liner!point-homogeneous}\defterm{point-homogeneous} if for every points $x,y\in X$ there exists an automorphism $A:X\to X$ such that $A(x)=y$;
\item \index{line-homogeneous liner}\index{liner!line-homogeneous}\defterm{line-homogeneous} if for every lines $L,\Lambda$ in $X$ there exists an automorphism $A:X\to X$ such that $A[L]=\Lambda$.
\end{itemize}
\end{definition}

Theorem~\ref{t:D=>point-homogeneous} implies that Desarguesian proaffine spaces are point-homogeneous.

\begin{corollary}\label{c:Des-proaffine=>point-homogeneous} Every Desarguesian proaffine space is point-homogeneous.
\end{corollary}

\begin{example} A finite punctured projective plane is not line-homogeneous because it contains lines of different cardinalities.
\end{example} 

Nonetheless, we have the following weaker line-homogeneity property of Desarguesian proaffine spaces.

\begin{proposition}\label{p:line-homogeneous} Let $L,\Lambda$ be two lines in a Desargusian proaffine space $X$. If the lines $L,\Lambda$ are both spreading or both non-spreading, then for every distinct points $o,e\in L$ and any distinct points $o',e'\in \Lambda$, there exists an automorphism $A:X\to X$ such that $Aoe=o'e'$ and hence $A[L]=\Lambda$.
\end{proposition}

\begin{proof} By Corollary~\ref{c:Desarg=>compreg}, the Desarguesian proaffine space $X$ is completely regular and by Theorem~\ref{t:Desargues-completion}, its spread completion $\overline X$ is a Desarguesian projective space. Let $\partial X\defeq \overline X\setminus X$ be the horizon of $X$ in its spread completion $\overline X$. Let $\overline L,\overline{\Lambda}$ be the flat hulls of the lines $L,\Lambda$ in the projective space $\overline X$. If the lines $L,\Lambda\subseteq X$ are non-spreading, then $\overline L=L$, $\overline{\Lambda}=\Lambda$ and hence $\overline L\cap \partial X=\varnothing=\overline\Lambda\cap\partial X$. We claim that the set $\{o,e\}\subset L$ is $\partial X$-independent. In the opposite case, $e\in\overline{\{o\}\cup \partial X}$ and $\overline L\cap\partial X=\Aline oe\cap\partial X\ne\varnothing$, by Corollary~\ref{c:line-meets-hyperplane}. This is a contradiction showing that the doubleton $\{o,e\}$ is $\partial X$-independent. By analogy we can prove that the doubleton $\{o',e'\}\subseteq\Lambda$ is $\partial X$-indepednent. By the Kuratowski--Zorn Lemma, the $\partial X$-independent sets $\{o,e\}$ and $\{o',e'\}$ can be enlarged to maximal $\partial X$-independent sets $M,M'\subseteq X$, respectively. By Proposition~\ref{p:union-of-independent}, the maximal $\partial X$-independence of the sets $M,M'$ ensures that $M\cap\partial X=\varnothing=M'\cap\partial X$ and the unions $M_\partial\cup M$ and $M_\partial M'$ are maximal independent sets in $\overline X$. Choose any function $r$ such that $\dom[r]=M_\partial \cup M$, $r(o)=o$ and $r(x)\in\Aline ox\setminus\{o,x\}$ for every $x\in\dom[r]\setminus\{o\}$. 
Also choose any function $r'$ such that $\dom[r']=M_\partial \cup M'$, $r(o')=o'$, and $r(x)\in\Aline {o'}x\setminus\{o',x\}$ for every $x\in\dom[r']\setminus\{o'\}$. 

By Corollary~\ref{c:Max=dim}, $|M|=\|X\|_{\partial X}=|M'|$ and hence there exists a bijectiion $\varphi:\dom[r]\cup\rng[r]\to\dom[r']\cup\rng[r']$ such that $r'\circ\varphi=\varphi\circ r$, $\varphi(o)=o'$, $\varphi(e)=e'$ and $\varphi[M_\partial]=M_\partial$.  By Corollary~\ref{c:Des-iso-reper}, there exists an automorphism $\overline\varphi:\overline X\to\overline X$ such that $\overline\varphi(x)=\varphi(x)$ for all $x\in\dom[r]\cup\rng[r]$. Then $\overline\varphi[\partial X]=\overline\varphi[\overline{M_\partial}]=\overline{\overline\varphi[M_\partial]}=\overline{M_\partial}=\partial X$ and hence $\Phi\defeq\overline\varphi{\restriction}_X$ is an automorphism of the liner $X$ such that $\Phi oe=o'e'$ and hence $\Phi[L]=\Lambda$.
\smallskip

Next, assume that the lines $L,\Lambda$ both are spreading. In this case, $\overline L\cap\partial X=\{L_\parallel\}\ne\varnothing\ne\{\Lambda_\parallel\}=\overline \Lambda\cap\partial X$. Using the Kuratowski--Zorn Lemma, choose any maximal independent set $M_\partial$ in $\partial X$ such that $L_\parallel\in M_\partial$. By the Kuratowski--Zorn Lemma, the $M_\partial$-indepedent set $\{o\}$ can be enlarged to a maximal $M_\partial$-independent set $M$. By Proposition~\ref{p:union-of-independent}, the set $M\cup M_\partial$ is maximal independent in $X$. Choose any function $r$ such that $\dom[r]=M$, $r(o)=o$, $r(L_\parallel)=e$ and $r(x)\in\Aline ox\setminus\{o,x\}$ for all $x\in\dom[r]\setminus \{o\}$. The maximal independence of the set $M\cup M_\partial=\dom[r]$ implies that $r$ is a projective reper in $X$.

By analogy we can find a maximal independent set $M_\partial'$ in $\partial X$ such that $\Lambda_\parallel\in M'_\partial$, and choose a maximal $M_\partial'$-independent set $M'$ such that $o'\in M'$. Choose any projective reper $r'$ in $X$ such that $\dom[r']=M'\cup M'_\partial$, $r'(o')=o'$, and $r'(\Lambda_\parallel)=e'$.  Corollary~\ref{c:Max=dim} ensures that $|M_\partial|=\|\partial X\|=|M'_\partial|$ and $|M|=\|X\|_{\partial X}=|M'|$. Then there exists a bijective map $\varphi:\dom[r]\to\dom[r']\to\dom[r']\cup\rng[r']$ such that $\varphi(o)=o'$, $\varphi(L_\parallel)=\Lambda_\parallel$, $\varphi[M_\parallel]=M'_\parallel$, $\varphi[M]=M'$ and $r'\circ\varphi=\varphi\circ r$.

By Corollary~\ref{c:Des-iso-reper}, there exists a projective automorphism $\overline\varphi:\overline X\to\overline X$ such that $\overline\varphi(x)=\varphi(x)$ for all $x\in\dom[r]\cup\rng[r]$. Then $\overline\varphi[\partial X]=\overline\varphi[\overline{M_\partial}]=\overline{\overline\varphi[M_\partial]}=\overline{M_\partial'}=\partial X$ and hence $\Phi\defeq\overline\varphi{\restriction}_X$ is an automorphism of the liner $X$ such that $\Phi oe=o'e'$ and hence $\Phi[L]=\Lambda$.
\end{proof}

\begin{corollary} A Desarguesian proaffine space $X$ is line-homogeneous if and only if $X$ is affine or projective.
\end{corollary}

\begin{proof} By Corollary~\ref{c:Desarg=>compreg} and Theorem~\ref{t:spread=projective1}, the Desarguesian proaffine space $X$ is para-Playfair.

If $X$ is projective, then $X$ contains no spreading lines and then $X$ is line-homogeneous, by Proposition~\ref{p:line-homogeneous}. If $X$ is affine, then by Theorem~\ref{t:Playfair<=>}, every line in $X$ is spreading and $X$ is line-homogeneous, by Proposition~\ref{p:line-homogeneous}. 

Now assume that $X$ is line-homogeneous. Assuming that $X$ is not projective, we shall prove that $X$ is affine.  Since $X$ is not projective, $X$ contains two coplanar disjoint lines $L,L'$. By Proposition~\ref{p:lines-in-para-Playfair}, the line $L$ is spreading. Since $X$ is line-homogeneous, for every line $\Lambda$ there exists an automorphism $F:X\to X$ such that $F[L]=\Lambda$. Since the line $L$ is spreading, so is the line $\Lambda=F[L]$. Therefore, every line in $X$ is spreading the the para-Playfair liner $X$ is Playfair, by  Theorem~\ref{t:Playfair<=>spreading}. By Theorem~\ref{t:Playfair<=>}, the Playfair liner $X$ is affine.
\end{proof}






\chapter{Automorphisms of  affine spaces}

In this chapter we shall study automorphisms of (Desarguesian) affine spaces, and also shall calculate the cardinality of the automorphism group of a finite Desarguesian affine space. 

\section{Coordinate charts on lines in Desarguesian affine spaces}

Let $L$ be a line in an affine space $X$. Given any distinct points $z,u\in L$, consider the function $$\phi_{zu}:L\to \dddot X_{\!\Join},\quad \phi_{zu}:x\mapsto \overvector{zxu},$$ assigning to every point $x\in L$ the portion $\overvector{zxu}$. Proposition~\ref{p:triple-exists} implies that the function $\phi_{zu}:L\to \dddot X_{\!\Join}$ is surjective. If the affine space $X$ is Desarguesian, then $\dddot X_{\!\Join}=\IR_X$ and the function $\phi_{zu}:L\to\dddot X_{\!\Join}=\IR_X$ is bijective. In this case the function $\phi_{zu}:L\to\IR_X$ is called a \index{coordinate chart}\defterm{coordinate chart} on the line $L$. The coordinate chart $\phi_{zu}$  identifies the line $L$ with the corps $\IR_X$. Under this identification, we have $\phi_{zu}(z)=\overvector{zzu}=0$ and $\phi_{zu}(u)=\overvector{zuu}=1$. 

For four points $z,u,o,e\in L$ with $z\ne u$ and $o\ne e$, the bijective function
$$\phi_{zu}\circ\phi_{oe}^{-1}:\IR_X\to\IR_X$$is called the \index{transition function}\defterm{transition function} between the charts $\phi_{oe}$ and $\phi_{zu}$.

$$
\xymatrix{
&L\ar^{\phi_{zu}}[rd]\ar_{\phi_{oe}}[ld]\\
\IR_X\ar_{\phi_{zu}\circ\phi^{-1}_{oe}}[rr]&&\IR_X
}
$$

We are going to show that the transition function is affine.

\begin{lemma}\label{l:transition-translation} Let $L$ be a line in a Desarguesian affine space $X$, and $z,u,o,e\in L$ be points with $z\ne u$ and $o\ne e$. The following conditions are equivalent:
\begin{enumerate}
\item $\overvector{zu}=\overvector{oe}$;
\item $\forall s\in\IR_X\;\;\phi_{zu}\circ\phi_{oe}^{-1}(s)=s+\overvector{zou}$;
\item $\exists b\in\IR_X\;\forall s\in \IR_X\;\;\phi_{zu}\circ\phi_{oe}^{-1}(s)=s+b$.
\end{enumerate}
\end{lemma}

\begin{proof} $(1)\Ra(2)$ Assume that $\overvector{zu}=\overvector{oe}$. Then $zu=Toe$ for some line translation $T:L\to L$. Given any scalar $s\in\IR_X$, consider the points $x\defeq\phi_{oe}^{-1}(s)$ and $x'\defeq T(x)$. Then $zx'u=Toxe$ and hence $\overvector{zx'u}=\overvector{oxe}$. It follows from $zx'=Tox$ that $\overvector{zx'}=\overvector{ox}$ and hence
$$\overvector{zo}+\overvector{zx'}=\overvector{zo}+\overvector{ox}=\overvector{zx}.$$The definition of the scalar addition ensures that
$$\phi_{zu}\circ \phi_{oe}^{-1}(s)=\phi_{zu}(x)=\overvector{zxu}=\overvector{zou}+\overvector{zx'u}=\overvector{zou}+\overvector{oxe}=\overvector{oxe}+\overvector{zou}=s+\overvector{zou}.$$

The implication $(2)\Ra(3)$ is trivial.
\smallskip

$(3)\Ra(1)$ Assume that there exists $b\in\IR_X$ such that $\phi_{zu}\circ\phi_{oe}^{-1}(s)=s+b$ for every scalar $s\in\IR_X$. Then $\phi_{zu}(o)=\phi_{zu}\circ\phi_{oe}^{-1}(0)=0+b=b$ and $$\overvector{zeu}=\phi_{zu}(e)=\phi_{zu}\circ\phi_{oe}^{-1}(1)=1+b=1+\phi_{zu}(o)=\overvector{zuu}+\overvector{zou}.$$
By the definition of the scalar addition, the equality $\overvector{zeu}=\overvector{zuu}+\overvector{zou}$ implies $\overvector{ze}=\overvector{zu}+\overvector{zo}$ and hence $\overvector{zu}=\overvector{ze}-\overvector{zo}=\overvector{oz}+\overvector{ze}=\overvector{oe}$.
\end{proof}

\begin{lemma}\label{l:transition-homothety} Let $L$ be a line in a Desarguesian affine space and $o\in L$ and $u,e\in L\setminus\{o\}$ be any points. For every scalar $s\in\IR_X$ we have the equality
$$\phi_{ou}\circ\phi_{oe}^{-1}(s)=s\cdot \overvector{oeu}.$$
\end{lemma}

\begin{proof} Given any scalar $s\in \IR_X$, consider the point $x\defeq\phi_{oe}^{-1}(s)$ and observe that $$\phi_{ou}\circ\phi^{-1}_{oe}(s)=\phi_{ou}(x)=\overvector{oxu}=\overvector{oxe}\cdot\overvector{oeu}=s\cdot \overvector{oeu}.
$$
\end{proof}

\begin{theorem}\label{t:transition-affine} Let $L$ be a line in a Desarguesian affine space $X$, and $z,u,o,e\in L$ be points with $z\ne u$ and $o\ne e$. Then $$\phi_{oe}\circ\phi_{zu}^{-1}(s)=s\cdot (\overvector{zeu}-\overvector{zou})+\overvector{zou}$$for every scalar $s\in\IR_X$.
\end{theorem}

\begin{proof} Being Desarguesian, the affine space $X$ is Thalesian. Then every vector in $X$ is functional and hence there exists a unique point $v\in L$ such that $\overvector{ov}=\overvector{zu}$. Given any scalar $s\in\IR_X$, we can apply Lemmas~\ref{l:transition-translation}, \ref{l:transition-homothety}, and conclude that
$$\phi_{zu}\circ\phi_{oe}^{-1}(s)=\phi_{zu}\circ \phi^{-1}_{ov}\circ\phi_{ov}\circ\phi_{oe}^{-1}(s)=\phi_{zu}\circ\phi^{-1}_{ov}(s\cdot \phi_{ov}(e))=(s\cdot\phi_{ov}(e))+\phi_{zu}(o).$$In particular,
$\phi_{zu}(e)=\phi_{zu}\circ\phi_{oe}^{-1}(1)=(1\cdot\phi_{ov}(e))+\phi_{zu}(o)$, which implies $\phi_{ov}(e)=\phi_{zu}(e)-\phi_{zu}(o)$ and hence
$$\phi_{zu}\circ\phi_{oe}^{-1}(s)=s\cdot(\phi_{zu}(e)-\phi_{zu}(o))+\phi_{zu}(o)=s\cdot(\overvector{zeu}-\overvector{zou})+\overvector{zou}.$$
\end{proof}

Now we shall reveal the algebraic structure of line affinities in Desarguesian affine spaces.

\begin{lemma}\label{l:transition-identity} Let $F:L\to L'$ be a bijective function between two lines $L,L'$ in a Desarguesian affine space $X$. Given two distinct points $o,u\in L$, consider their images $o'\defeq F(o)$ and $u'\defeq F(u)$. The function $F$ is a line affinity if and only if $\phi_{o'u'}\circ F\circ \phi_{ou}^{-1}$ is the identity map of the corps $\IR_X$.
\end{lemma}

\begin{proof} To prove the ``only if'' part, assume that $F$ is a line affinity. Given any scalar $s\in\IR_X$, find a unique point $x\in L$ such that $\overvector{oxu}=s$. Then for the point $x'\defeq F(x)$ we have $o'x'u'=Foxu$ and hence $\overvector{o'x'u'}=\overvector{oxu}$ and 
$$\phi_{o'u'}\circ F\circ \phi_{ou}^{-1}(s)=\phi_{o'u'}\circ F(x)=\phi_{o'u'}(x')=\overvector{o'x'u'}=\overvector{oxu}=s,$$
witnessing that $\phi_{o'u'}\circ F\circ \phi_{ou}^{-1}$ is the identity function of the corps $\IR_X$. 
\smallskip

Now assume conversely that $\phi_{o'u'}\circ F\circ \phi_{ou}^{-1}$ is the identity function of the corps $\IR_X$. 
Then $\phi_{ou}=\phi_{o'u'}\circ F$. 
By Theorem~\ref{t:aff-trans}, there exists a line affinity $A:L\to L'$ such that $Aou=o'u'$. By the definition of the affine equivalence of line triples, for every $x\in L$, the line triple $Aoxu$ is affinely equivalent to the line triple $oxu$. Since $\phi_{ou}=\phi_{o'u'}\circ F$, for the point $x'\defeq F(x)$, we have $\overvector{oxu}=\phi_{ou}(x)=\phi_{o'u'}\circ F(x)=\phi_{o'u'}(x')=\overvector{o'x'u'}$.
Taking into account that $Aou=o'u'$ and $Aoxu\Join oxu\Join o'x'u'$,
we conclude that $A(x)=x'=F(x)$, by the defintion of a Desarguesian triple. Therefore, $F=A$ is a line affinity.  
\end{proof}

Now we are able to prove our main result on coordinate description of line affinities.

\begin{theorem}\label{t:lineaff=aff} Let $L,\Lambda$ be two lines in a Desarguesian affine space and $z,u\in L$, $o,e\in \Lambda$ be points with $z\ne u$ and $o\ne e$. For any function $F:L\to \Lambda$ and the function  $\Phi\defeq \phi_{oe}\circ F\circ\phi_{zu}^{-1}:\IR_X\to\IR_X$, the following conditions are equivalent:
\begin{enumerate}
\item The function $F:L\to \Lambda$ is a line affinity;
\item $\exists a\in\IR^*_X\;\exists b\in\IR_X\;\forall s\in\IR_X\;\;\Phi(s)=(s\cdot a)+b$.
\end{enumerate}
\end{theorem}

\begin{proof} $(1)\Ra(2)$ Assume that $F:L\to\Lambda$ is a line affinity.
Consider the pair $z'u'\defeq Fzu$ on the line $\Lambda$. By Lemma~\ref{l:transition-identity}, the transition function $\phi_{z'u'}\circ F\circ \phi_{zu}^{-1}$ is the identity function of the corps $\IR_X$. Consider the scalars $a\defeq \overvector{z'eu'}-\overvector{z'ou'}\ne 0$ and $b\defeq \overvector{z'ou'}$. By Theorem~\ref{t:transition-affine}, for every scalar $s\in\IR_X$, we have
$$\Phi(s)=\phi_{oe}\circ F\circ\phi_{zu}^{-1}=\phi_{oe}\circ\phi_{z'u'}^{-1}\circ\phi_{z'u'}\circ F\circ \phi^{-1}_{zu}(s)=\phi_{oe}\circ\phi_{z'u'}^{-1}(s)=s\cdot a+b.$$
\smallskip

$(2)\Ra(1)$ Assume that the condition (2) is satisfied and hence there exist scalars $a\in\IR_X^*$ and $b\in\IR_X$ such that $\Phi(s)=s\cdot a+b$ for every $s\in \IR_X$. It follows that  $b=\Phi(0)=\phi_{oe}\circ F\circ\phi_{zu}^{-1}(0)=\phi_{oe}\circ F(z)$ and $a+b=\Phi(1)=\phi_{oe}\circ F\circ\phi_{zu}^{-1}(1)=\phi_{oe}\circ F(u)$.  Taking into account that $a\ne 0$, we conclude that 
$\phi_{oe}\circ F(z)=b\ne a+b=\phi_{oe}\circ F(u)$ and hence $F(z)\ne F(u)$. By Theorem~\ref{t:aff-trans}, there exists a line affinity $A:L\to\Lambda$ such that $Fzu=Azu$. By the (already proved) implication $(1)\Ra(2)$, there exist scalars $a'\in\IR^*_X$ and  $b'\in\IR_X$ such that $\phi_{oe}\circ A\circ\phi_{zu}^{-1}(s)=s\cdot a'+b'$ for every $s\in\IR_X$. In particular, 
$$b'=0\cdot a'+b'=\phi_{oe}\circ A\circ\phi_{zu}^{-1}(0)=\phi_{oe}\circ A(z)=\phi_{oe}\circ F(z)=\phi_{oe}\circ F\circ\phi_{zu}^{-1}(0)=0\cdot a+b=b$$and
$$a'+b=a'+b'=1\cdot a'+b'=\phi_{oe}\circ A\circ \phi_{zu}^{-1}(1)=\phi_{oe}\circ A(u)=\phi_{oe}\circ F(u)=\phi_{oe}\circ F\circ\phi_{zu}^{-1}(1)=1\cdot a+b=a+b$$and hence $a'=a$.
Then $$\phi_{oe}\circ F\circ\phi_{zu}^{-1}(s)=s\cdot a+b=s\cdot a'+b'=\phi_{oe}\circ A\circ\phi_{zu}^{-1}(s)$$and hence the function
$$F=\phi_{oe}^{-1}\circ (\phi_{oe}\circ F\circ\phi_{zu}^{-1})\circ\phi_{zu}=\phi_{oe}^{-1}\circ (\phi_{oe}\circ A\circ\phi_{zu}^{-1})\circ\phi_{zu}=A$$is a line affinity.
\end{proof}

\section{Affinities, portionalities and scalarities in affine spaces}

\begin{definition}\label{d:affinity-portionality-scalarity} An automorphism $A:X\to X$ of an affine space $X$ is called
\begin{enumerate}
\item an \index{affinity}\index{automorphism!affinity}\defterm{affinity} if for every line $L\subset X$, the restriction $A{\restriction}_L$ is a line affinity of $X$;
\item a \index{portionality}\index{automorphism!portionality}\defterm{portionality} if $Axyz\Join xyz$ for every line triple $xyz\in\dddot X$ in $X$;
\item a \index{scalarity}\index{automorphism!scalarity}\defterm{scalarity} if $Axyz\Join xyz$ for every Desarguesian line triple $xyz\in\dddot X$ in $X$.
\end{enumerate}
\end{definition}

We recall that an automorphism $A: X\to X$ of an affine space $X$ is \defterm{paracentral} if for every $x,y\in X$ with $x'\defeq A(x)\ne x$ and $y'\defeq A(y)\ne y$, the lines $\Aline x{x'}$ and $\Aline y{y'}$ are parallel.  By Proposition~\ref{p:aff-paracentral<=>}, an automorphism of an affine space is paracentral if and only if it is a translation, a hypershear or a hyperscale.

\begin{proposition}\label{p:affinity=>portionality=>scalarity} For any automorphism $A:X\to X$ of an affine space $X$, the following conditions hold.
\begin{enumerate}
\item If $A$ is paracentral, then $A$ is an affinity.
\item If $A$ is an affinity, then $A$ is portionality.
\item If $A$ is a portionality, then $A$ is a scalarity.
\end{enumerate}
\end{proposition}

\begin{proof} 1. Assume that the automorphism $A$ is paracentral.  To prove that $A$ is an affinity, we have to check that for every line $L\subseteq X$ the restriction $A{\restriction}_L$ is a line affinity. This is trivially true, if $A$ is the identity automorphism of the affine space $X$. So, assume that $A$ is not identity and hence the family of parallel lines $\{\Aline xy:xy\in A\setminus 1_X\}$ is not empty and is contained in a unique direction $\delta\in\partial X$. 

Consider the line $L'\defeq A[L]$. If $L\notin\delta$, then $A{\restriction}_L=\{(x,y)\in L\times L':\Aline x\delta=\Aline y\delta\}$, which means that $A{\restriction}_L$ is a line projection and hence a line affinity between the lines $L$ and $L'$. So, assume that $L\in\delta$ and hence $L'=A[L]=L$. Take any line $\Lambda\notin \delta$ that has a common point $o$ with the line $L$. Take any points $a\in L\setminus \Lambda$, $b\in \Lambda\setminus L$ and concider the line projection $B\defeq\{xy\in L\times\Lambda:\Aline xy\subparallel \Aline ab\}$ from the line $L$ onto the line $\Lambda$. Consider the line $\Lambda'\defeq A[\Lambda]$ and observe that the restriction $A'\defeq A{\restriction}_\Lambda$ is a line projection from the line $\Lambda$ onto the line $\Lambda'$. Consider the map $C:\Lambda'\to L'$, $C:x\mapsto AB^{-1}A^{-1}(x)$ and observe that it is a line projection (by Exercise~\ref{ex:conj-lineprojection}). Then $A{\restriction}_L=CA'B$ is a composition of three line projections and hence it is a line affinity.
\smallskip

The statements 2 and 3 follow immediately from the definitions of portions and scalars.
\end{proof}

By Proposition~\ref{p:aff-paracentral<=>}, an automorphism of an affine space is paracentral if and only if it is a translation or a hypershear or a hyperscale. Therefore, for every automorphism of an affine space, we have the implications:
$$\mbox{translation or hypershear or hyperscale} \Leftrightarrow\mbox{paracentral}\Ra\mbox{affinity}\Ra\mbox{portionality}\Ra\mbox{scalarity}.$$
For Desarguesian affine spaces, the last three notions are equivalent.

\begin{theorem}\label{t:D=>a=p=s} For any automorphism $A:X\to X$ of a  Desarguesian affine space $X$, the following conditions are equivalent:
\begin{enumerate}
\item $A$ is an affinity;
\item $A$ is a portionality;
\item $A$ is a scalarity.
\end{enumerate}
\end{theorem}

\begin{proof} The implications $(1)\Ra(2)\Ra(3)$ follow from Proposition~\ref{p:affinity=>portionality=>scalarity}. To prove the implication $(3)\Ra(1)$, assume that $A$ is a scalarity. To prove that $A$ is an affinity, fix any line $L$. Choose any distinct points $o,u\in L$ and consider the pair $o'u'\defeq Aou$. By Theorem~\ref{t:paraproj-exists}, there exists a line affinity $B$ in $X$ such that $Bou=o'u'$. We claim that $A=B$. Given any point $x\in L$, we should prove that the points $a\defeq A(x)$ and $b\defeq B(x)$ coincide.  Since $A$ is a scalarity and $B$ is a line affinity, $o'au'=Aoxu\Join oxu\Join Boxu=o'bu'$ and hence $a=b$ because the line triples $o'au'$ and $o'bu'$ are Desarguesian and projectively equivalent. Therefore, $A(x)=a=b=B(x)$ and $A=B$.
\end{proof} 

\begin{exercise} Find ane axample of an affinity of a Desarguesian affine plane, which is not paracentral.
\end{exercise}

\begin{definition}\label{d:affine-aut} An automorphism $A:X\to X$ of a Desarguesian affine space $X$ is called \index{affine automorphism}\index{automorphism!affine}\defterm{affine} if $A$ satisfies the equivalent conditions of Theorem~\ref{t:D=>a=p=s}.
\end{definition}

For an affine space $X$, let $\Aff(X)$, $\Aff_P(X)$, $\Aff_S(X)$ be the subsets of the automorphism group $\Aut(X)$, consisting of all affinities, portionalities, and scalarities, respectively.

The following proposition is easy and is left to the reader as an exercise.

\begin{proposition} For every affine space $X$, the sets $\Aff(X),\Aff_P(X),\Aff_S(X)$ are normal subgroups of the automorphism group $\Aut(X)$ such that $\Aff(X)\subseteq \Aff_P(X)\subseteq\Aff_S(X)$. If the affine space $X$ is Desarguesian, then  $\Aff(X)=\Aff_P(X)=\Aff_S(X)$.
\end{proposition}

\begin{remark} By Corollary$^\dag$~\ref{c:portions=3}, $|\dddot X_{\Join}|=3$ for every non-Thalesian finite affine space $X$. This corollary$^\dag$ implies that $\Aff_P(X)=\Aff_S(X)=\Aut(X)$ for all non-Thalesian finite affine spaces $X$.
\end{remark}

\begin{remark} By Theorem~\ref{t:non-Thales-line-aff}, for every line $L$ in a non-Thalesian affine space of finite order $n\ne 24$, the group $\Sym^{\Join}_X(L)$ of line affinities of the line $L$ contains the alternating group $\Alt(L)$. Moreover, if $|X|_2\in(2+4\IZ)\cup(3+\IZ)$, then $\Sym^{\Join}_X(L)=\Sym(L)$. In the latter case, $\Aff(X)=\Aut(X)$. If $\Sym^{\Join}_X(L)=\Alt(L)$ for some (equivalently, any) line $L$, then 
$$\Aff(X)=\{A\in\Aut(X):\forall L\in\mathcal L_X\;A{\restriction}_L\in \Alt(L)\}.$$
\end{remark}

By Theorem~\ref{t:aff-trans}, for every distinct points $x,y$ of an affine space $X$, there exists a line affinity $A:\Aline xy\to\Aline xy$ such that $Axy=yx$ and hence $Axy\# yx$.

\begin{question} Let $L$ be a line in an affine space $X$. Is there a line affinity $A:L\to L$ such that $Axy\in\overvector{yx}$ for every point $x,y\in L$?
\end{question}

\begin{question}\label{q:centsymmetry} Let an affine space $X$. Is there a dilation $A:X\to X$ such that $Axy\in\overvector{yx}$ for every $x,y\in X$?
\end{question}

\begin{remark} The answer to Question~\ref{q:centsymmetry} is affirmative in Thalesian affine spaces $X$: for every point $o\in X$, the homothety $H_o:X\to X,\quad H_o:x\mapsto o+(-1)\cdot \overvector{ox}$, has the required property. Here $-1$ is the opposite scalar to the scalar $1\in\IR_X$.
\end{remark}




Let $X$ be an affine space and $\IR_X$ be its scalar corps. For every automorphism $A:X\to X$ and every line triple $xyz$ in $X$, the image $Axyz$ is a line triple in $X$. If the line triple $xyz$ is Desarguesian, then so is the line triple $Axyz$. This observation allows us to define the map $\dddot A:\IR_X\to\IR_X$ assigning to every scalar $s\in\IR_X$ the scalar $\dddot A(s)=\{Axyz:xyz\in s\}$. It is easy to see that the function $\dddot A:\IR_X\to\IR_X$ is an automorphism of the corps $\IR_X$. 
For two automorphisms $A,B:X\to X$ and their composition $C=A\circ B$, the automorphism $\dddot C:\IR_X\to\IR_X$ is equal to the composition $\dddot A\circ\dddot B$ of the automorphisms $\dddot A,\dddot B$ of the corps $\IR_X$. This means that the correspondence $T:\Aut(X)\to\Aut(\IR_X)$, $T:A\mapsto \dddot A$, is a homomorphism from the automorphism group $\Aut(X)$ of the liner $X$ into the automorphism group $\Aut(\IR_X)$ of the corps $\IR_X$. The definition of a scalarity implies that the automorphism $\dddot A:\IR_X\to\IR_X$ is the identity map of $\IR_X$ if and only if the automorphism $A$ is a scalarity of $X$. 
 

Let us write down this fact for future references.

\begin{proposition}\label{p:Aut(X)->Aut(RX)} For every affine space $X$, the function $$\Aut(X)\to \Aut(\IR_X),\quad A\mapsto \dddot A,$$ is a homomorphism from the automorphism group $\Aut(X)$ of the liner $X$ into the automorphism group $\Aut(\IR_X)$ of the corps $\IR_X$. The kernel of the homomorphism $\Aut(X)\to\Aut(\IR_X)$ coincides with the group $\Aff_S(X)$ of scalarities.
\end{proposition}

\begin{remark} By Proposition~\ref{p:aff-paracentral<=>}, translations, hypershears and hyperscales of affine spaces are affinities and hence scalarities. On the other hand, homotheties of Desarguesian affine spaces need not be affine. Indeed, by Proposition~\ref{p:central-homo}, for every non-zero element $s$ of a corps $R$, the map $H:\Pi\to \Pi$, $H:(x,y)\mapsto (s\cdot x,s\cdot y)$, is a homothety of the coordinate plane $\Pi=R\times R$ of $R$. The scalar corps $\IR_\Pi$ of the plane $\Pi=R\times R$ can be identified with the corps $R$. The homothety $H$ induces the inner automorphism $\dddot H:\IR_\Pi\to\IR_\Pi$, $\dddot H:x\mapsto sxs^{-1}$,  of the scalar corps $\IR_\Pi$. Therefore, the homothety $H$ is an affinity of the plane $\Pi$ if and only if $s$ belongs to the centre $\mathcal Z(R)=\{z\in R:\forall x\in R\;\;x\cdot z=z\cdot x\}$ of the corps $R$.
\end{remark}

\section{Affine automorphisms of Desarguesian affine spaces}

Propositon~\ref{p:Aut(X)->Aut(RX)} and Theorem~\ref{t:D=>a=p=s} imply the following characterization of affine automorphisms of Desarguesian affine spaces.

\begin{corollary} A self-map $F:X\to X$ on a Desarguesian affine space $X$ is an affine automorphism of $X$ if and only if there exists a bijective map $\vec F:\overvector X\to\overvector X$ satisfying the following conditions:
\begin{enumerate}
\item $\forall \vec{\boldsymbol x},\vec{\boldsymbol y}\in\overvector X\;\;\vec F(\vec{\boldsymbol x}+
\vec{\boldsymbol y})=\vec F(\vec {\boldsymbol x})+\vec F(\vec{\boldsymbol y})$;
\item $\forall s\in\IR_X\;\forall \vec{\boldsymbol v}\in\overvector X\;\;\vec F(s{\cdot}\vec {\boldsymbol v})=s{\cdot}\vec F(\vec{\boldsymbol v})$;
\item $\exists o\in X\;\forall \vec{\boldsymbol v}\in\overvector X\;\;F(o+\vec{\boldsymbol v})=F(o)+ \vec F(\vec {\boldsymbol v})$.
\end{enumerate}
\end{corollary}
 
Theorem~\ref{t:extension-iso} implies the following extension result.

\begin{corollary}\label{c:affine-extend-auto} For every bijection $B:M\to M'$ between two maximal independent sets $M,M'$ in a Desarguesian affine liner $X$, there exists a unique affine automorphism  $F:X\to X$ of $X$ such that $F{\restriction}_M=B$.
\end{corollary}

\begin{theorem}\label{t:Aut(X)=Aff(X)xAut(RX)} Let $X$ be a Desarguesian affine space, $M$ be a maximal independent subset in $X$, and $\Aut_M(X)\defeq\{A\in\Aut(X):A{\restriction}_M=1_M\}$. Then $\Aut_M(X)$ is a subgroup of the group $\Aut(X)$ such that $\Aut_M(X)\cap\Aff(X)=\{1_X\}$, $\Aut(X)=\Aff(X)\cdot\Aut_M(X)$ and the function $\Aut_M(X)\to\Aut(\IR_X)$, $A\mapsto\dddot A$, is a group isomorphism.
\end{theorem}

\begin{proof} It is clear that $\Aut_M(X)$ is a subgroup of the group $\Aut(X)$. The uniqueness part of Corollary~\ref{c:affine-extend-auto} implies that $\Aut_M(X)\cap\Aff(X)=\{1_X\}$.

To prove that $\Aut(X)=\Aff(X)\cdot\Aut_M(X)$, take any automorphism $A:X\to X$ and consider the maximal independent set $M'\defeq A[M]$ in $X$. By Corollary~\ref{c:affine-extend-auto}, there exists an affine automorphism $B:X\to X$ such that $B{\restriction}_M=A{\restriction}_M$. Then the automorphism $C=B^{-1}A$ belongs to the subgroup $\Aut_M(X)$ and hence $A=BC\in \Aff(X)\cdot\Aut_M(X)$.

It remains to prove that the function $I:\Aut_M(X)\to\Aut(\IR_X)$, $I:A\mapsto\dddot A$, is a group isomorphism. By Proposition~\ref{p:Aut(X)->Aut(RX)}, the function $T:\Aut(X)\to\Aut(\IR_X)$ is a group homomorphism whose kernel coincides with the subgroup $\Aff_S(X)$ of scalarities, which coincides with the group of affinities $\Aff(X)$, by Theorem~\ref{t:D=>a=p=s}. Since $\Aff(X)\cap\Aut_M(X)=\{1_X\}$, the restriction $I\defeq T{\restriction}_{\Aut_M(X)}$ is injective. 

It remains to prove that the homomorphism $I:\Aut_M(X)\to\Aut(\IR_X)$ is surjective. Fix any point $o\in M$ and let $B\defeq M\setminus\{o\}$. Consider the map $\Sigma:X\to\IR_X^{\oplus B}$ assigning to every $x\in X$ the unique finitely supported function $c_x\in \IR_X^{\oplus B}$ such that $\overvector{ox}=\sum_{b\in B}c_x(b)\cdot\overvector{ob}$.

Given any automorphism $A:\IR_X\to \IR_X$ of the corps $\IR_X$, consider the function $\tilde A:X\to X$ assigning to every $x\in X$ the unique point $y\in X$ such that
$$\overline{oy}=\sum_{b\in B}A(c_x(b))\cdot b.$$
By (the proof of) Theorem~\ref{t:extension-iso}, the function $\tilde A$ is an automorphism of the liner $X$. Since $A(0)=0$ and $A(1)=1$, the automorphism $\Lambda\defeq \tilde A$ belongs to the subgroup $\Aut_M(X)$. We claim that $\dddot \Lambda=A$. Indeed, take any scalar $s\in\IR_X$ and consider its image $t\defeq A(s)\in\IR_X$. Choose any point $b\in B$ and using Theorem~\ref{t:triple-unique} find a unique point $x\in \Aline ob$ such that $oxb\in s$. Then $\overvector{ox}=\overvector{oxb}\cdot\overvector{ob}=s\cdot\overvector{ob}$ and for the point $y\defeq\tilde A(x)=\Lambda(x)$ we have $\overvector{oy}=A(s)\cdot \overvector{ob}=t\overvector{ob}$. The definition of the automorphism $\dddot \Lambda$ ensures that $\dddot\Lambda(s)=\dddot \Lambda(\overvector {oxb})=\overvector{oyb}=t=A(s)$, witnessing that $\dddot \Lambda=A$. Therefore, $I:\Aut_M(X)\to\Aut(\IR_X)$ is a group isomorphism. 
\end{proof}

A group $G$ is called a \index{semidirect product of groups}\defterm{semidirect product} of groups $A,B$ if there exist a normal subgroup $\tilde A\subseteq G$ and a subgroup $\tilde B\subseteq G$ such that $\tilde A$ is isomorphic to $A$, $\tilde B$ is isomorphic to $B$, $\tilde A\cap\tilde B=\{1_G\}$ and $G=\tilde A\cdot\tilde B$. In this case we write that $G=A\rtimes B$.
This definition of a semidirect product and Theorem~\ref{t:Aut(X)=Aff(X)xAut(RX)} imply the following description of the automorphism group of a Desarguesian affine space. 

\begin{corollary}\label{c:Aut(X)=Aff(X)xAut(RX)} The automorphism group $\Aut(X)$ of a Desarguesian affine space $X$ is a semidirect product $\Aff(X)\rtimes\Aut(\IR_X)$ of the groups $\Aff(X)$ and $\Aut(\IR_X)$.
\end{corollary}

\begin{exercise}
Given a finite-dimensional Desarguesian affine space $X$, find a matrix representation of the group $\Aff(X)$ of affine automorphisms of $X$.
\end{exercise}

\section{Automorphisms of finite fields} 

Corollary~\ref{c:Aut(X)=Aff(X)xAut(RX)} motivates the problem of studying the automorphism groups of finite fields. 

Let $F$ be a finite field and $p$ be its characteristic, which is the smallest number $p$ such that $p\cdot 1=0$. Since the field $F$ has no divisors of zero, the characteristic $p$ is a prime number. 

\begin{lemma}\label{l:Frobenius} Let $F$ be a finite field of characteristic $p$. For every $x,y\in F$,$$(x+y)^p=x^p+y^p.$$
\end{lemma} 

\begin{proof} By the binomial formula,
$$(x+y)^p=\sum_{k=0}^p\tfrac{p!}{k!(p-k)!}x^ky^{p-k}.$$
Observe that for every $k\in\{1,\dots,p-1\}$, the prime number $p$ divides $p!$ but does not divide $k!(p-k)!$, which implies that $p$ divides the binomial coefficient $\frac{p!}{k!(p-k)!}$. Then $\frac{p!}{k!(p-k)!}x^ky^{p-k}$ is zero in the field $F$ and hence
$$(x+y)^p=x^p+\Big(\sum_{k=1}^{p-1}\tfrac{p!}{k!(p-k)!}x^ky^{p-k}\Big)+y^p=x^p+0+y^p=x^p+y^p.$$
\end{proof}

Lemma~\ref{l:Frobenius} implies that for every finite field $F$ of characteristic $p$, the function $\Phi:F\to F$, $\Phi:x\mapsto x^p$, is an endomorphism of the field $F$. Since $x^p=0$ iff $x=0$, the endomorphism $\Phi$ is injective. Since $F$ is finite, the injectivity of $\Phi$ implies its surjectivity and bijectivity. Therefore, the map $\Phi:F\to F$, $\Phi:x\mapsto x^p$ is an automorphism of the field $F$, called the \index{Frobenius automorphism}\index{automorphism!Frobenius}\defterm{Frobenius automorphism} of the field $F$, after \index[person]{Frobenius}Frobenius\footnote{{\bf Ferdinand Georg Frobenius} (1849 -- 1917) was a German mathematician, best known for his contributions to the theory of elliptic functions, differential equations, number theory, and to group theory. He is known for the famous determinantal identities, known as Frobenius--Stickelberger formulae, governing elliptic functions, and for developing the theory of biquadratic forms. He was also the first to introduce the notion of rational approximations of functions (nowadays known as Pad\'e approximants), and gave the first full proof for the Cayley--Hamilton theorem. He also lent his name to certain differential-geometric objects in modern mathematical physics, known as Frobenius manifolds.}, a known algebraist.

\begin{theorem}\label{t:Aut(F)} Let $F$ be a finite field of character $p$ and cardinality $p^n$ for some $n\in\IN$. The automorphism group $\Aut(F)$ of the field $F$ is cyclic of order $n$, generated by the Frobenius automorphism. 
\end{theorem}

\begin{proof} By Theorem~\ref{t:F*-cyclic}, the mutiplicative group $F^*\defeq F\setminus\{0\}$ of the finite field $F$ is cyclic, so has some generator $g$. Since $F$ is an $n$-dimension vector space over the minimal field $\IF_p=\{0,1,\dots,p-1\}\subseteq F$, the elements $1,g,g^2,\dots, g^n$ are linearly dependent over $\IF_p$ and hence $P(g)=0$ for some polynomial $P\in\IF_p[x]$ of degree at most $n$. Observe that every automorphism $A:F\to F$ of the field has $A(1)=1$ and hence $A(c)=c$ for every $c\in\IF_p=\{0,1,2,\dots,p-1\}$. Therefore, the automorphism $A$ does not change the coefficients of the polynomial $P$ and hence $0=A(0)=A(P(g))=P(A(g))$, which means that $A(g)$ is a root of the polynomial $P(x)$ in the field $F$. Since this polynomial has at most $n$ roots in $F$, the set $\{A(g):A\in\Aut(F)\}\subseteq\{x\in F:P(x)=0\}$ has cardinality $\le n$. Since $g$ is a generator of the multiplicative group $F^*$, every automorphism $A$ of the field $F$ is uniquely determined by its value $A(g)$. This implies that 
$$|\Aut(F)|=|\{A(g):A\in\Aut(F)\}|\le|\{x\in F:P(x)=0\}|\le n.$$On the other hand, the Frobenius automorphism $\Phi$ and its powers belong to the group $\Aut(F)$. It follows from $|F^*|=p^n-1$ that $g^{p^n-1}=1$ and hence all elements $g^p,g^{p^2},\cdots,g^{p^n}$ are distinct, which implies that $\Phi,\Phi^2,\dots,\Phi^n$ are distinct elements of the group $\Aut(F)$ and hence $\Aut(F)=\{\Phi,\Phi^2,\dots,\Phi^n\}$ is a cyclic group of order $n$, generated by the Frobenius automorphism $\Phi$.
\end{proof}

\section{Automorphisms of finite Desarguesian affine spaces}

Applying Corollary~\ref{c:Aut(X)=Aff(X)xAut(RX)} and Theorem~\ref{t:Aut(F)} we can calculate the cardinality of the automorphism group $\Aut(X)$ of a finite Desarguesian affine space $X$ and its normal subgroup $\Aff(X)$, consisting of affine automorphisms.

\begin{theorem}\label{t:|Aff(X)|=} For every finite Desarguesian affine space $X$, the group $\Aff(X)$ of affine automorphisms of $X$ has cardinality $$|\Aff(X)|=r^d\cdot\prod_{k=0}^{d-1}(r^d-r^k),$$
where  $d\defeq\dim(X)=\|X\|-1$ and $r\defeq|\IR_X|$.
\end{theorem}

\begin{proof} By Corollary~\ref{c:Max=dim}, every maximal independent set in $X$ has cardinality $\|X\|=d+1$. Fix any maximal independent set $M=\{x_0,x_1,\dots,x_{d}\}$. By Corollary~\ref{c:affine-extend-auto}, every affine automorphism $A\in\Aff(X)$ is uniquely determined by its restriction $A{\restriction}_M$. Moreover, for every maximal independent set $M'=\{y_0,\dots,y_{d}\}$ in $X$, there exists an automorphism $A\in\Aut(X)$ such that $A(x_i)=y_i$ for all $i\le d$. Therefore, the number of affine automorphisms of $X$ is equal to the number of $(d+1)$-tuples $(y_0,y_1,\dots,y_{d})\in X^{d+1}$ such that the set $\{y_0,y_1,\dots,y_{d}\}$ is maximal independent  in $X$. By Proposition~\ref{p:add-point-to-independent}, the number of such $(d+1)$-tuples in equal to the product
$$|X|\cdot(|X|-|X|_1)\cdot(|X|-|X|_2)\cdots (|X|-|X|_{d}).$$
Let $r\defeq|\IR_X|$ be the cardinality of the corps $\IR_X$ of the Desarguesian affine space $X$. By Theorem~\ref{t:paraD<=>Desargues}, every line in $X$ has cardinality $r\defeq|\IR_X|$. Therefore, $|X|_2=r$ and $|X|_k=|X|_2^{k-1}=r^{k-1}$ for every $k\in\{1,\dots,d\}$, by Corollary~\ref{c:affine-cardinality}. Then
$$|\Aff(X)|=|X|\cdot\prod_{k=1}^{d}(|X|-|X|_k)=r^d\cdot\prod_{k=1}^d(r^d-r^{k-1})=r^d\cdot\prod_{k=0}^{d-1}(s^d-s^k).$$
\end{proof}

\begin{theorem} For every finite Desarguesian affine space $X$, the group $\Aut(X)$ of  automorphisms of $X$ has cardinality $$|\Aut(X)|=|\Aut(\IR_X)|\cdot|\Aff(X)|=n\cdot s^r\cdot\prod_{k=0}^{r-1}(s^r-s^k),$$
where $r\defeq\|X\|$ and $s\defeq |\IR_X|=p^n$ for some prime number $p$ and some $n\in\IN$.
\end{theorem}

\begin{proof} By Corollary~\ref{c:Aut(X)=Aff(X)xAut(RX)}, the automorphism group $\Aut(X)$ of the Desarguesian affine space $X$ is a semidirect product $\Aff(X)\rtimes\Aut(\IR_X)$ and hence $|\Aut(X)|=|\Aff(X)|\cdot|\Aut(\IR_X)|$. Since the affine space $X$ is finite, so is its corps $\IR_X$. By Wedderburn Theorem~\ref{t:Wedderburn-Witt}, the finite corps $\IR_X$ is finite field. Then $|\IR_X|=p^n$ for some $n\in\IN$ and some prime number $p$, equal to the characteristic of the field $\IR_X$. By Theorem~\ref{t:Aut(F)}, $|\Aff(\IR_X)|=n$. Applying Theorem~\ref{t:|Aff(X)|=}, we conclude that
$$|\Aut(X)|=|\Aut(\IR_X)|\cdot|\Aff(X)|=n\cdot r^d\cdot\prod_{k=0}^{d-1}(r^d-r^k),$$where $d\defeq\dim(X)=\|X\|-1$ and $r\defeq|\IR_X|=p^n$.
\end{proof}

\begin{Exercise} Find the cardinality of the automorphism group $\Aut(X)$ of a Desarguesian projective space $X$.
\smallskip

{\em Hint:} This requires some extra-knowledge and will be done in Theorem~\ref{t:proj-|Aut|}.
\end{Exercise}

\chapter{Permutation groups of line affinities}\label{ch:permutation}

\rightline{\em Groups, as men, will be known by their actions.}

\rightline{Guillermo Moreno}
\bigskip

In this chapter we analyse the structure of the permutation groups 
$\Sym^{\Join}_X(L)$ acting on the lines $L$ of affine spaces 
$X$. We obtain a satisfactory description of these groups for Thalesian affine spaces and also for non-Thalesian finite affine spaces. The latter result is obtained by applying a non-trivial classification of finite $2$-transitive permutation groups. This classification is one of the principal achievements of permutation group theory and one of the most spectacular applications of the classification of finite simple groups, which was completed only at the end of the twentieth century.

\section{Groups of line affinities in Thalesian affine spaces}

In this section we describe the structure of the group $\Sym_X^{\Join}(L)$ of line affinities on a line $L$ in a Thalesian affine space $X$.

We recall that every line $L$ in a Thalesian affine space $X$, the set $\vec L\defeq\{\vec{\boldsymbol v}:\vec{\boldsymbol v}\subparallel L\}$ is an $\IR_X$-submodule of the $\IR_X$-module $\overvector X$.

\begin{proposition}\label{p:line-projection=alg-structure} For every line projection $A:L\to \Lambda$ in a Thalesian affine space $X$, there exists a unique isomorphism $\vec A:\vec L\to\vec \Lambda$ of the $\IR_X$-modules $\vec L$ and $\vec\Lambda$ such that $A(x+\vec {\boldsymbol v})=A(x)+\vec A(\vec {\boldsymbol v})$ for all $x\in L$ and $\vec {\boldsymbol v}\in \vec L$.
\end{proposition}

\begin{proof} By Proposition~\ref{p:flat=>R-submodule}, the subsets $\vec L\defeq\{\vec {\boldsymbol x}\in\overvector X:\vec{\boldsymbol x}\subparallel L\}$ and $\vec \Lambda\defeq\{\vec{\boldsymbol x}\in\overvector X:\vec{\boldsymbol x}\subparallel \Lambda\}$ are $\IR_X$-submodules of the $\IR_X$-module $\overvector X$.

 Define the map $\vec A:\vec L\to\vec\Lambda$ assignining to every vector $\overvector{xy}\in \vec L$ the vector $\overline{x'y'}$ where $x'y'\defeq Axy$. Let us show that the map $\vec A$ is well-defined. Given any pair $uv\in L^2$ with $\overvector{uv}=\overvector{xy}$, we need to show that $\overvector{u'v'}=\overvector{x'y'}$ where $u'v'\defeq Auv$. 

If $x=y$, then the equality $\overvector{uv}=\overvector{xy}$ implies $u=v$. Then $y'=A(y)=A(x)=x'$, $v'=A(v)=A(u)=u'$, and finally, $\overvector{u'v'}=\vec {\bf 0}=\overvector{x'y'}$.

So, assume that $x\ne y$, which implies $u\ne v$. Since the affine space $X$ is Thalesian, there exists a translation $T:X\to X$ such that $Txy=uv$. Then $T[L]=T[\Aline xy]=\Aline uv=L$. By Propositions~\ref{p:translation|=>line-translation} and \ref{p:RTR=T}, the function $B\defeq ATA^{-1}:\Lambda\to\Lambda$ is a line translation. Observe that $Bx'y'=ATA^{-1}x'y'=ATxy=Auv=u'v'$. Then $\overvector{x'y'}=\overvector{u'v'}$, by Definition~\ref{d:vector} of a vector. Therefore, the function $\vec A:\vec L\to\vec\Lambda$ is well-defined. The bijectivity of the line projection $A$ ensures that the function $\vec A:\vec L\to\vec \Lambda$ is bijective.

To see that $\vec A$ is an isomorphism of the $\IR_X$-modules $\vec L$ and $\vec \Lambda$, it suffices to check that $\vec A(\vec {\boldsymbol u}+\vec {\boldsymbol v})=\vec A(\vec{\boldsymbol u})+\vec A(\vec{\boldsymbol v})$ and $\vec A(s{\cdot}\vec{\boldsymbol u})=s{\cdot} A(\vec{\boldsymbol u})$ for all $\vec{\boldsymbol u},\vec{\boldsymbol v}\in\vec L$ and $s\in\IR_X$. By Corollary~\ref{c:Thalesian<=>vector=funvector}, there exist points $x,y,z\in L$ such that $\vec {\boldsymbol u}=\overvector{xy}$ and $\vec{\boldsymbol v}=\overvector{yz}$. By Theorem~\ref{t:vector-addition}, $\vec {\boldsymbol u}+\vec{\boldsymbol v}=\overvector{xz}$. Consider the points $x'\defeq A(x)$, $y'\defeq A(y)$ and $z'\defeq A(z)$ in the line $\Lambda$, and observe that $\vec A(\vec{\boldsymbol u})=\vec A(\overvector{xy})=\overvector{x'y'}$, $\vec A(\vec{\boldsymbol v})=\vec A(\overvector{yz})=\overvector{y'z'}$ and hence 
$$\vec A(\vec{\boldsymbol u})+\vec A(\vec{\boldsymbol v})=\overvector{x'y'}+\overvector{y'z'}=\overvector{x'z'}=\vec A(\overvector{xz})=\vec A(\vec{\boldsymbol u}+\vec{\boldsymbol v}).$$

To show that $\vec A(s\cdot\vec{\boldsymbol u})=s\cdot\vec A(\vec{\boldsymbol u})$, find a unique point $\lambda\in L$ such that $x\lambda y\in s$. Consider the point $\lambda'\defeq A(\lambda)\in\Lambda$ and observe that $\overvector{x'\lambda'y'}=\overvector{x\lambda y}=s$, by the definition of the portion $s$. Then $$\vec A(s\cdot\vec{\boldsymbol u})=\vec A(\overvector{x\lambda y}\cdot\overvector{xy})=\vec A(\overvector{x\lambda})=\overvector{x'\lambda'}=\overvector{x'\lambda'y'}\cdot\overvector{x'y'}=s\cdot\vec A(\overvector {xy})=s\cdot \vec A(\vec{\boldsymbol v}).$$
Therefore, $\vec A:\vec L\to\vec\Lambda$ is an isomorphism of the $\IR_X$-modules $\vec L$ and $\vec \Lambda$.

Now given any point $x\in L$ and vector $\vec{\boldsymbol v}\in\overvector L$, we shall prove that $A(x+\vec{\boldsymbol v})=A(x)+\vec A(\vec{\boldsymbol v})$. Find a point $y\in L$ such that $\vec {\boldsymbol v}=\overvector{xy}$ and consider the points $x'\defeq A(x)$ and $y'\defeq A(y)$. Then $$A(x+\vec{\boldsymbol v})=A(x+\overvector{xy})=A(y)=y'=x'+\overvector{x'y'}=A(x)+\vec A(\vec{\boldsymbol v}).$$ The $\IR_X$-linearity of the isomorphism $\vec A$ ensures that for every scalar $s\in\IR_X$, we
have
$$A(x+s\cdot \vec{\boldsymbol v})=A(x)+\vec A(s\cdot\vec{\boldsymbol v})=A(x)+s\cdot \vec A(\vec{\boldsymbol v}).$$

To see that the function $\vec A$ is unique, assume that $\vec B:\vec L\to\vec \Lambda$ is another function such that $A(x+\vec{\boldsymbol v})=A(x)+B(\vec{\boldsymbol v})$ for all points $x\in L$ and vectors $\vec{\boldsymbol v}\in\vec L$. Given any vector $\vec{\boldsymbol v}\in\vec L$, find two points $x,y\in L$ with $\vec {\boldsymbol v}=\overvector{xy}$. Consider the points $x'\defeq A(x)$ and $y'\defeq A(y)$ and observe that $$x'+\vec A(\vec {\boldsymbol v})=x'+\overvector{x'y'}=y'=A(y)=A(x+\overvector{xy})=A(x+\vec{\boldsymbol v})=A(x)+\vec B(\vec {\boldsymbol v})=x'+\vec B(\vec {\boldsymbol v})$$ and hence $\vec A(\vec {\boldsymbol v})=\overvector{x'y'}=\vec B(\vec{\boldsymbol v})$, by Corollary~\ref{c:vector-action}. 
\end{proof}

Since every line affinity in an affine space is the composition of finitely many line projections, Proposition~\ref{p:line-projection=alg-structure} implies the following corollary describing the algebraic structure of line affinities in a Thalesian affine space.

\begin{corollary}\label{c:line-affinity=alg-structure} Let $L$ be a line in a Thalesian affine space $X$. For any line affinity $A:L\to L$, there exists a unique automorphism $\vec A:\vec L\to\vec L$ of the $\IR_X$-module $\vec L$ such that $A(x+\vec{\boldsymbol  v})=A(x)+\vec A(\vec{\boldsymbol  v})$ for all $x\in X$, $s\in\IR_X$ and $\vec{\boldsymbol v}\in\overvector L$.
\end{corollary}

For a line $L$ in a Thalesian affine space $X$, let $\Aff_{\IR_X}(L)$ be the set of all permutations $A:L\to L$ for which there exists an automorphism $\vec A$ of the $\IR_X$-module $\vec L$ such that $A(x+\vec {\boldsymbol v})=A(x)+\vec A(\vec{\boldsymbol v})$ for all $x\in L$ and $\vec{\boldsymbol v}\in\vec L$. Such functions $A$ will be called \defterm{$\IR_X$-affine}. 
It is easy to see that $\Aff_{\IR_X}(L)$ is a subgroup of the symmetric group $\Sym(L)$ on $L$. 
The definitions of the groups $\Sym^\#_X(L),\Sym^{\Join}_X(L)$ and Corollary~\ref{c:line-affinity=alg-structure} ensure that 
$$\Sym^\#_X(L)\subseteq\Sym^{\Join}_X(L)\subseteq\Aff_{\IR_X}(L)\subseteq \Sym(L).$$
In the following proposition we describe the structure of the permutation groups $\Sym^{\Join}_X(L)$ and $\Aff_{\IR_X}(L)$ on a line $L$ is a Thalesian affine space $X$.

\begin{theorem}\label{t:grp-line-trans} Let $L$ be a line in a Thalesian affine space $X$ and $o\in L$.
\begin{enumerate}
\item The group $\Sym_X^{\#}(L)$ of line translations of $L$ is isomorphic to the additive group $(\vec L,+)$ of the $\IR_X$-module $\vec L$. 
\item The group $\Sym_X^\#(L)$ is a normal Abelian subgroup of the groups $\Sym^{\Join}_X(L)$ and $\Aff_{\IR_X}(L)$. 
\item The group $\mathsf{GL}_o(L)\defeq \{A\in\Aff_{\IR_X}(L):A(o)=o\}\subseteq\Aff_{\IR_X}(L)$ is isomorphic to the group $\Aut(\vec L)$ of all automorphisms of the $\IR_X$-module $\vec L$.
\item The group $\Aff_{\IR_X}(L)$ is the semidirect product $\Sym_X^{\#}(L)\rtimes \mathsf{GL}_o(L)$ of its subgroups $\Sym_X^\#(L)$ and $\mathsf{GL}_o(L)$.
\item The group $\Sym^{\Join}_X(L)$ of line affinities of the line $L$ is the semidirect product\\ $\Sym_X^{\#}(L)\rtimes \Sym^{\Join}_o(L)$ of its subgroups $\Sym_X^\#(L)$ and\\ $\Sym^{\Join}_o(L)\defeq\{A\in\Sym^{\Join}_X(L):A(o)=o\}=\Sym^{\Join}_X(L)\cap\mathsf{GL}_o(L)$.
\end{enumerate}
\end{theorem}

\begin{proof} By Theorem~\ref{t:paraD<=>translation}, the Thalesian affine space $X$ is translation. 
\smallskip

1. We recall that $\vec L\defeq\{\vec{\boldsymbol v}\in \overvector X:\vec{\boldsymbol v}\subparallel L\}$ and hence the restriction operator $R:\vec L\to \Sym(L)$, $R:\vec{\boldsymbol v}\mapsto \vec{\boldsymbol v}{\restriction}_L$, is a well-defined function from the additive group $(\vec L,+)$ to the symmetric group $\Sym(L)$. By Proposition~\ref{p:Ax=Bx=>A=B}, the function $R$ is injective. By Proposition~\ref{p:translation|=>line-translation} and Theorem~\ref{t:paraD<=>translation|}, $R[\vec L]=\Sym^\#_X(L)$.  By Corollary~\ref{c:Thalesian-vectors-commutative}, the group $(\vec L,+)$ is commutative and the bijective function $R:\vec L\to \Sym^\#_X(L)$ is an isomorphism between the groups $(\vec L,+)$ and $\Sym^\#_X(L)$. 
\smallskip

2. By the preceding item, the group $\Sym^\#_X(L)$ coincides with the additive group of the $\IR_X$-module $\vec L$, which is commutative by the definition of an $\IR_X$-module, see Exercise~\ref{ex:R-module+Abelian-group}. To see that the Abelian subgroup $\Sym^\#_X(L)$ is normal in the group $\Aff_{\IR_X}[L]$, take any line translation $T\in\Sym^\#_X(L)$ and any affine map $A\in\Aff_{\IR_X}(L)$. Find a unique automorphism $\vec A:\vec L\to\vec L$ of the $\IR_X$-module $\vec L$ such that $A(x+\vec{\boldsymbol v})=A(x)+\vec A(\vec{\boldsymbol v})$ for all $x\in L$ and $\vec{\boldsymbol v}\in \vec L$. By Theorem~\ref{t:paraD<=>translation|}, there exists a vector $\vec {\boldsymbol v}\in\vec L$ such that $T(x)=x+\vec {\boldsymbol v}$ for all $x\in L$. Consider the vector $\vec {\boldsymbol u}\defeq\vec A(\vec{\boldsymbol v})\in\vec L$. By Proposition~\ref{p:translation|=>line-translation}, the map $S:L\to L$, $S:x\mapsto x+\vec {\boldsymbol u}$, is a line translation. The equality
$$AT(x)=A(x+\vec {\boldsymbol v})=A(x)+\vec A(\vec{\boldsymbol v})=A(x)+\vec{\boldsymbol u}=SA(x)$$holding for all points $x\in L$ implies $AT=SA$ and $ATA^{-1}=S\in \Sym^\#_X(L)$ and witnesses that the subgroup $\Sym^\#_X(L)$ is normal in the group $\Aff_{\IR_X}(L)$ and also in its subgroup $\Sym^{\Join}_X(L)\subseteq\Aff_{\IR_X}(L)$.
\smallskip

3. Consider the function $\Phi:\mathsf{GL}_o(L)\to \mathsf{GL}(\vec L)$ assigning to every $\IR_X$-affine bijection $A\in\mathsf{GL}_o(L)$, the unique automorphism  $\vec A\in \mathsf{GL}[\vec L]$ such that $A(x+\vec{\boldsymbol v})=A(x)+\vec A(\vec {\boldsymbol v})$ for all $x\in L$. To see that the map $\Phi$ is injective, take any $\IR_X$-affine bijections $A,B\in \AGL[L]$ with $\vec A=\vec B$. Then for all $x\in L$ we have the equality $$A(x)=A(o+\overvector{ox})=A(o)+\vec A(\overvector{ox})=o+\vec B(\overvector{ox})=B(o)+\vec B(\overvector{ox})=B(x),$$ witnessing that $A=B$ and the function $\Phi$ is injective. To see that $\Phi$ is surjective, take any automorphism $\vec A\in \mathsf{GL}[\vec L]$ and consider the permutation $A:L\to L$, $A:x\mapsto o+\vec A(\overvector{ox})$. 
Given any point $x\in L$ and vector $\vec{\boldsymbol v}\in \vec L$, let $y\defeq x+\vec v$ and observe that $\vec v=\overvector{xy}$ according to Corollary~\ref{c:vector-action}. The associativity of the addition of vectors ensures that
$$
\begin{aligned}
A(x+\vec v)&=A(x+\overvector{xy})=A(y)=o+\vec A(\overvector{oy})=
o+\vec A(\overvector{ox}+\overvector{xy})\\
&=o+(\vec A(\overvector{ox})+\vec A(\overvector{xy}))=(o+\vec A(\overvector{ox}))+\vec A(\overvector{xy})=A(x)+\vec A(\vec v),
\end{aligned}
$$witnessing that $\Phi(A)=\vec A$ and the function $\Phi:\mathsf{GL}_o[L]\to\mathsf{GL}[\vec L]$ is bijective.  To see that $\Phi$ is a group isomorphism, take any $\IR_X$-affine bijections $A,B\in\mathsf{GL}_o[L]$ and  consider the automorphisms $\vec A\defeq\Phi(A)$, $\vec B\defeq\Phi(B)$, and $\vec C\defeq\Phi(AB)$ of the $\IR_X$-module $\vec L$. For every $\vec v\in\vec L$, we have the equalities
$$AB(o)+\vec C(\vec v)=AB(o+\vec v)=A(B(o)+\vec B(\vec v))=AB(o)+\vec A(\vec B(\vec v))=AB(o)+\vec A\vec B(\vec v),$$withessing that $\Phi(AB)=\vec C=\vec A\vec B=\Phi(A)\Phi(B)$ and $\Phi:\mathsf{GL}_o[L]\to\mathsf{GL}[\vec L]$ is a group isomorphism.
\smallskip

4. Theorem~\ref{t:unique-translation} ensures that $\Sym_X^\#(L)\cap\mathsf{GL}_o(L)=\{1_L\}$. By the item (2), the subgroup $\Sym^\#_X(L)$ is normal in the group $\AGL[L]$. It remains to prove that $\Aff_{\IR_X}(L)=\Sym^\#(L)\circ\mathsf{GL}_o(L)$. Given any permutation $A\in \Sym^{\Join}_X(L)$, apply Theorem~\ref{t:unique-translation}, and find a line translation $T\in\Sym^\#(L)$ such that $T(o)=A(o)$. Then the permutation $T^{-1}A\in\Aff_{\IR_X}[L]$ belongs to the subgroup $\mathsf{GL}_o(L)$, witnessing that $A=TT^{-1}A\in \Sym_X^\#(L)\circ\mathsf{GL}_o(L)$. Therefore, $\AGL(L)=\Sym^\#(L)\rtimes \mathsf{GL}_o(L)$ is the semidirect product of the normal subgroup $\Sym^\#_X(L)$ and $\mathsf{GL}_o(L)$.
\smallskip

5. By analogy we can prove that the subgroup $\Sym^{\Join}_X(L)$ of $\AGL(L)$ is the semidirect product $\Sym^{\Join}_X(L)=\Sym^\#(L)\rtimes \Sym^{\Join}_o(L)$ is the semidirect product of the normal subgroup $\Sym^\#_X(L)$ and $\Sym^{\Join}_o(L)=\{A\in \I^{\Join}_X[L]:A(o)=o\}=\Sym^{\Join}_X(L)\cap\mathsf{GL}_o(L)$.
\end{proof}

\begin{exercise} Show that for any lines $L$ and $\Lambda$ in an affine space $X$, the groups $\Sym^{\Join}_X(L)$ and $\Sym^{\Join}_X(\Lambda)$ are isomorphic.
\smallskip

{\em Hint:} Look at Proposition~\ref{p:Sym-Join-isomorpic}
\end{exercise}

\begin{remark} By computer calculations Ivan Hetman found a non-Thalesian affine plane of order $9$ containing two lines $L$ and $\Lambda$  such that the groups $\Sym^{\#}_X(L)$ and $\Sym^{\#}_X(\Lambda)$ are not isomorphic (and moreover, have different cardinality), see the table in Remark~\ref{r:Sym-Hetman}.
\end{remark}

\begin{theorem} For an affine space $X$, the following conditions are equivalent:
\begin{enumerate}
\item $X$ is Thalesian;
\item for every line $L$ in $X$, the group $\Sym^\#_X(L)$ is Abelian;
\item for some concurrent lines $L,\Lambda$ in $X$ the groups $\Sym^\#_X(L)$ and $\Sym^\#_X(\Lambda)$ are Abelian.
\end{enumerate}
\end{theorem}

\begin{proof} $(1)\Ra(2)$ Assume that the affine space $X$ is Thalesian. By Theorem~\ref{t:grp-line-trans}(2), for every line $L$ in $X$, the group $\Sym^\#_X(L)$ is Abelian.
\smallskip

The implication $(2)\Ra(3)$ is trivial.
\smallskip

$(3)\Ra(1)$ assume that $L$ and $\Lambda$ are two concurrent lines in $X$ such that the groups  $\Sym^\#_X(L)$ and $\Sym^\#_X(\Lambda)$ are Abelian. We have to prove that the afffine space $X$ is Thalesian. If $\|X\|\ne 3$, then the affine space $X$ is Desarguesian and Thalesian, by Corollary~\ref{c:affine-Desarguesian} and Theorem~\ref{t:ADA=>AMA}. So, assume that $\|X\|=3$, which means that $X$ is a Playfair plane. Let $o$ be the unique common point of the concurrent lines $L,\Lambda$. By Theorem~\ref{t:paraD<=>translation}, it suffices to show that for every point $p\in X$ there exists a translation $T:X\to X$ such that $T(o)=p$.

\begin{claim}\label{cl:T(o)=p} For every point $p\in L\cup\Lambda$ there exists a translation $T:X\to X$ such that $T(o)=p$.
\end{claim}

\begin{proof} First assume that $p\in L$. Choose any line $A$ in $X$ such that $A\cap L=\varnothing$. Choose any point $a\in A$ and consider the directions $\boldsymbol u\defeq(\Aline oa)_\parallel$ and $\boldsymbol v\defeq(\Aline ap)_\parallel$ in the Playfair plane $X$. By definition of the group $\Sym^\#_X(L)$, the permutation $P\defeq \boldsymbol v_{L,A}\boldsymbol u_{A,L}$ of the line $L$ belongs to the group $\Sym^\#_X(L)$. Observe that $P(o)=\boldsymbol v_{L,A}\boldsymbol u_{A,L}(o)=\boldsymbol v_{L,A}(a)=p$. Since the group $\Sym^\#_X(L)$ is Abelian, we can apply Gleason's Theorem~\ref{t:Gleason1} and conclude that $P=T{\restriction}_L$ for some translation $T:X\to X$. Then $T(o)=P(o)=p$. By analogy we can find a translation $T:X\to X$ with $T(o)=p$ if $p\in\Lambda$.
\end{proof}

Now we are ready to prove that for every point $p\in X$ there exists a translation $T:X\to X$ such that $T(o)=p$. If $p\in L\cup \Lambda$, then the translation $T$ exists, by Claim~\ref{cl:T(o)=p}. So, assume that $p\notin L\cup\Lambda$. Since $X$ is a Playfair plane, there exist points $x\in L\setminus\{o\}$ and $y\in\Lambda\setminus\{o\}$ such that $oxpy$ is a parallelogram. By Claim~\ref{cl:T(o)=p}, there exist translations $T_x,T_y$ of the plane $X$ such that $T_x(o)=x$ and $T_y(o)=y$. Then $T=T_xT_y$ is a translation of $X$ such that $T(o)=p$, witnessing that the Playfair plane $X$ is translation. By Theorem~\ref{t:paraD<=>translation}, the translation Playfair plane $X$ is Thalesian.
\end{proof}

\section{The holomorph of a group}

In this section we recall some basic notions of the theory of permutation groups and also will discuss the holomorphs of groups (consisting of affine transformations of groups).

By a \defterm{permutation group} we understand a subgroup $G$ of the group $\Sym(X)$ of all bijections of a set $X$. In this case we say that the permutation group $G$ \defterm{acts} on the set $X$. The cardinality of the set $X$ is called the \defterm{degree} of the permutation group $G$. In other to exclude trivial cases, we shall always assume that $|X|>1$. The group $\Sym(X)$ is endowed with the operation $\circ$ of composition of functions. We shall often omit the symbol operation writing $fg$ instead of $f\circ g$ for two permutations $f,g\in\Sym(X)$.

Two permutation groups $G\subseteq \Sym(X)$ and $H\subseteq \Sym(Y)$ are called \defterm{isomorphic} if there exists a bijection $f:X\to Y$ such that $H=\{fgf^{-1}:g\in G\}$. In this case the function $(\cdot)^f:G\to H$, $(\cdot)^f:g\mapsto g^f\defeq fgf^{-1}$, is an isomorphism of the groups $G$ and $H$. 

\begin{exercise} Find two permutation groups $G,H\subseteq \Sym(X)$ which are isomorphic as groups but not isomorphic as permutation groups.

{\em Hint:} Take any two involutions $g,h:X\to X$ with $$|\{x\in X:f(x)\ne x\}|\ne |\{x\in X:h(x)\ne x\}|$$ and consider the two-element subgroups $G\defeq\{1_X,g\}$ and $H\defeq\{1_X,h\}$ of $\Sym(X)$.\end{exercise}

\begin{proposition}\label{p:Sym-Join-isomorpic} For any lines $L,\Lambda$ in an affine space $X$, the permutation groups $\Sym^{\Join}_X(L)$ and $\Sym^{\Join}_X(\Lambda)$ are isomorphic.
\end{proposition}

\begin{proof} Choose any point $a\in \Lambda$ and find a unique line $I$ in $X$ such that $a\in I$ and $L\parallel I$. Observe that the parallel lines $L,I$ are coplanar and the intersecting ines $I,\Lambda$ are coplanar.
Then there exist planes $\Pi,\Pi'$ in $X$ such that $L\cup I\subseteq\Pi$ and $I\cup\Lambda\subseteq \Pi'$. Choose any directions $\boldsymbol u\in\partial \Pi\setminus L_\parallel$ and $\boldsymbol v\in\partial\Pi'\setminus(I_\parallel\cup\Lambda_\parallel)$, and consider the line affinity $A\defeq \boldsymbol v_{\Lambda,I}\boldsymbol u_{I,L}:L\to \Lambda$. The permutation $A$ determines an isomorphism $(\cdot)^A:\Sym^{\Join}_X(L)\to\Sym^{\Join}_X(\Lambda)$, $(\cdot)^A:F\mapsto AFA^{-1}$ of the permutation groups $\Sym^{\Join}_X(L)$ and $\Sym^{\Join}_X(\Lambda)$.
\end{proof}
 
Now we consider some special permutations of groups.

\begin{definition}
A permutation $A:G\to G$ of a group $G$ is called
\begin{enumerate}
\item a \defterm{left translation} if there exists an element $a\in G$ such that $A(x)=ax$ for all $x\in G$;
\item a \defterm{right translation} if there exists an element $b\in G$ such that $A(x)=xb$ for all $x\in G$;
\item a \defterm{two-sided translation} if there exist elements $a,b\in G$ such that $A(x)=axb$ for all $x\in G$;
\item an \defterm{automorphism} of $G$ if $A(x{\cdot}y)=A(x){\cdot}A(y)$ for all $x,y\in G$;
\item an \defterm{inner automorphism} of $G$ if there exists an element $a\in G$ such that $A(x)=axa^{-1}$ for all $x\in G$.
\end{enumerate}
\end{definition}

It is easy to see that the sets $L(G)$, $R(G)$, $T(G)$, and $\Aut(G)$ of all left translations, right translations, two-sided translations, and automorphisms are subgroups of the permutation group $\Sym(G)$. Moreover, $T(G)=L(G)\circ R(G)=R(G)\circ L(G)$. 

\begin{definition} Let $\Hol(G)$ be the subgroup of the symmetric group $\Sym(G)$, generated by the set $T(G)\cup\Aut(G)$. The permutation group $\Hol(G)$ is called the \defterm{holomorph} of the group $G$. Elements of the holomorph $\Hol(G)$ are called \defterm{holomorphisms} or else \defterm{affine transformations} of the group $G$.
\end{definition}

\begin{proposition} Let $G$ be a group.
\begin{enumerate}
\item $L(G)\cap \Aut(G)=\{1_G\}=\Aut(G)\cap R(G)$.
\item $T(G)\cap \Aut(G)$ is the subgroup of inner automorphisms of the group $G$.
\item The subgroups $L(G),R(G),T(G)$ are normal in the holomorph $\Hol(G)$ of $G$.
\item $\Hol(G)=L(G)\circ \Aut(G)=R(G)\circ \Aut(G)$.
\item $\Hol(G)=L(G)\rtimes \Aut(G)=R(G)\rtimes \Aut(G)$.
\end{enumerate}
\end{proposition}

\begin{proof} 1. The equality  $L(G)\cap \Aut(G)=\{1_G\}=\Aut(G)\cap R(G)$  follows immediately from the equality $A(1)=1$ holding for every automorphism $A$ of $G$.
\smallskip

2. The definitions imply that every inner automorphism of $G$ is an element of the set $T(G)\cap\Aut(G)$. On the other hand, for every automorphism $A\in T(G)\cap\Aut(G)$, there exist elements $a,b\in G$ such that $A(x)=axb$ for all $x\in G$. Since $A$ is an automorphism of the group $G$, $1=A(1)=ab$ and hence $b=a^{-1}$. Then $A(x)=axa^{-1}$ for all $x\in G$, witnessing that the automorphism $A$ is inner.
\smallskip

3. Since $T(G)=L(G)\circ R(G)$ and elements of the subgroups $L(G)$ and $R(G)$ commute, to prove that the subgroup $L(G)$ is normal in the holomorph $\Hol(G)$, it suffices to check that for every left shift $L_a:V\to V$, $L_a:x\mapsto a{\cdot}x$, and every automorphism $A\in \Aut(G)$ the permutation $AL_aA^{-1}$ is a left shift. Indeed, for every $x\in G$ we have
$AL_aA^{-1}(x)=A(a{\cdot}A^{-1}(x))=A(a){\cdot}AA^{-1}(x)=A(a){\cdot}x$, witnessing that $AL_aA^{-1}$ is a left translation of the group $G$, and hence the subgroup $L(G)$ in normal in the holomorph $\Hol(G)$ of $G$. By analogy we can prove the normality of the subgroup $R(G)$ of right translations in $\Hol(G)$. The normality of the subgroups $L(G)$ and $R(G)$ implies the normality of their product $T(G)=L(G)\circ R(G)$ in the group $\Hol(G)$.
\smallskip

4. The equality $\Hol(G)=L(G)\circ \Aut(G)$ will follow as soon as we check that $L(G)\circ\Aut(G)$ is a subgroup of $\Sym(G)$, containing the set $T(G)$.
To see that $L(G)\circ\Aut(G)$ is a subgroup, fix any elements $A,B\in L(G)\circ \Aut(G)$  and find elements $a,b\in G$ and automorphisms $\vec A,\vec B\in\Aut(G)$ such that $A(x)=a{\cdot}\vec A(x)$ and $B(x)=b{\cdot}\vec B(x)$ for all $x\in G$. Consider the element $c\defeq a\cdot\vec A(b)$ and the automorphism $\vec C\defeq \vec A\circ \vec B$ of $G$, and observe that
$$AB(x)=a\cdot\vec A(b\cdot\vec B(x))=a\cdot\vec A(b)\cdot \vec A(\vec B(x))=c\cdot\vec C(x)$$for all $x\in G$, witnessing that the permutation $AB$ of $G$ belongs to the set $L(G)\circ \Aut(G)$.

To see that the inverse permutation $A^{-1}$ belongs to $L(G)\circ\Aut(G)$, observe that for any $x\in G$ and $y=A(x)=a\cdot\vec A(x)$, we have $$A^{-1}(y)=x=\vec A^{-1}(a^{-1}\cdot y)=\vec A^{-1}(a^{-1})\cdot\vec A^{-1}(y),$$witnessing that $A^{-1}\in L(G)\circ \Aut(G)$. Therefore, $L(G)\circ \Aut(G)$ is a subgroup of $\Hol(G)$.
\smallskip

To see that $T(G)\subseteq L(G)\circ \Aut(G)$, take any two-sided shift $T\in R(G)$ and find elements $a,b\in G$ such that $T(x)=a{\cdot}x{\cdot}b$ for all $x\in G$. Consider the left shift $L:G\to G$, $L:x\mapsto (a{\cdot}b){\cdot}x$, and the inner automorphism $\vec A:G\to G$, $\vec A:x\mapsto b^{-1}{\cdot}x{\cdot}b$. Observe that $T=L\vec A\in L(G)\circ \Aut(G)$. 

Since $T(G)\subseteq L(G)\circ\Aut(G)$, the subgroup generated by the set $T(G)\cup \Aut(G)$ coincides with the subgroup $L(G)\circ \Aut(G)$ and hence $\Hol(G)=L(G)\circ \Aut(G)$. By analogy we can prove that $\Hol(G)=R(G)\circ \Aut(G)$.
\smallskip

5. Since $L(G)$ is a normal subgroup of the group $\Hol(G)=L(G)\circ\Aut(G)$ and $L(G)\cap\Aut(G)=\{1_G\}$, the product $\Hol(G)=L(G)\circ\Aut(G)$ is semidirect and hence $\Hol(G)=L(G)\rtimes \Aut(G)$. By analogy we can show that $\Hol(G)=R(G)\rtimes \Aut(G)$.
\end{proof}

\section{Permutation groups of affine type}

In this section we prove some basic facts about the structure of permutation groups of affine type and then apply these results in Theorem~\ref{t:Schleiermacher}, which characterizes Thalesian affine spaces via properties of the permutation groups of line affinities. We begin by introducing several notions of transitivity for sets in permutation groups.

\begin{definition} Let $n$ be any cardinal number. A subset $\mathcal S$ of the symmetric group $\Sym(X)$ on a set $X$ is called
\begin{itemize}
\item \defterm{$n$-transitive} if for any injective function $F$ with $|F|=n$ and $\dom[F]\cup\rng[F]\subseteq X$, there exists  a permutation $S\in \mathcal S$ such that $F\subseteq S$;
\item \defterm{$n$-sharp} if for any injective function $F$ with $|F|=n$ and $\dom[F]\cup\rng[F]\subseteq X$, there exists at most one permutation $S\in \mathcal S$ such that $F\subseteq S$;
\item  \defterm{sharply $n$-transitive} if for any injective function $F$ with $|F|=n$ and $\dom[F]\cup\rng[F]\subseteq X$, there exists a unique permutation $S\in \mathcal S$ such that $F\subseteq S$;
\item \defterm{transitive} if $\mathcal S$ is $1$-transitive;
\item \defterm{sharply transitive} if $\mathcal S$ is sharply $1$-transitive.
\end{itemize}
\end{definition}

\begin{exercise} Show that a subset of a permutation group is sharply $n$-transitive if and only if it is $n$-sharp and $n$-transitive.
\end{exercise}

\begin{exercise} Show that every $n$-sharp set in a permutation group is $k$-sharp for all $k\ge n$.
\end{exercise}

\begin{exercise} Show that for every group $G$, the permutation groups $L(G)$ and $R(G)$ of left and right translations of $G$ are sharply transitive.
\end{exercise}



\begin{theorem}\label{t:affine-type} If $N$ is a sharply transitive normal subgroup of a permutation group $G\subseteq\Sym(X)$ on a set $X$, and $o$ is any point of the set $X$, then 
\begin{enumerate}
\item the group $G$ is the semidirect product $N\rtimes G_o$ of the normal subgroup $N$ and the stabilizer $G_o\defeq\{g\in G:g(o)=o\}$ of the point $o$;
\item the bijective map $\alpha:N\to X$, $\alpha:v\mapsto v(o)$, determines the isomorphism $(\cdot)^\alpha:G\to G^\alpha$, $(\cdot)^\alpha:g\mapsto \alpha^{-1}g\alpha$, of the permutation group $G$ and the subgroup $G^\alpha\defeq \{\alpha^{-1}g\alpha:g\in G\}$ of the holomorph $\Hol(N)$  of the group $N$.
\end{enumerate}
\end{theorem}

\begin{proof} 


1. Consider the stabilizer $G_o\defeq\{g\in G:g(o)=o\}$ of the point $o$. The sharp transitivity of the group $N$ ensures that $N\cap G_o=\{1_X\}$. Given any element $g\in G$, find a permutation $h\in N$ such that $h(o)=g(o)$ and conclude that $h^{-1}g\in G_o$ and hence $g\in N\cdot G_o$ and $G=N\cdot G_o$. Since $N\cap G_o=\{1_X\}$, the product $G=N\cdot G_o$ is semidirect and hence $G=N\rtimes G_o$.
\smallskip

2. Consider the map $\gamma:G_o\to\Aut(N)$ assigning to each $g\in G_o$ the conjugating automorphism $\gamma_g:N\to N$, $\gamma_g:x\mapsto gxg^{-1}$, of the group $N$. Let us show that $\gamma:G_o\to\Aut(N)$ is a group homomorphism. Indeed, for any $g,h\in G_o$ and any $v\in N$, we obtain $\gamma_{gh}(v)=(gh)v(gh)^{-1}=g(hvh^{-1})g^{-1}=\gamma_g\gamma_h(v)$, withessing that $\gamma_{gh}=\gamma_g\circ\gamma_h$ and hence $\gamma:G_o\to \Aut(N)$ is a group homomorphism. 
\smallskip

Next, consider the map $\delta:G\to \Hol(N)$ assigining to every 
element $g=va\in G=N\rtimes G_o$ the holomorphism $\delta_g=\ell_v\gamma_a\in\Aff(N)$ of the group $N$. Here $v\in N$ and $a\in G_o$ are unique elements such that $g=va$, and $\ell_v:N\to N$, $\ell_v:x\mapsto vx$, is the left translation of the group $N$ by the element $v$. We claim that $\delta$ is a group homomorphism. Indeed, for every $va,ub\in N\rtimes G_o=G$ we have $$\delta_{vaub}(x)=\delta_{vaua^{-1}ab}(x)=\ell_{vaua^{-1}}\gamma_{ab}(x)=vaua^{-1}abx(ab)^{-1}=va(ubxb^{-1})a^{-1}=\delta_{va}\delta_{ub}(x),$$
witnessing that $\delta:G\to \Hol(N)$ is a group homomorphism.

Consider the bijective map $\alpha:N\to X$, $\alpha:v\mapsto v(o)$, and the injective homomorphism $(\cdot)^\alpha:G\to\Sym(N)$, $(\cdot)^\alpha:g\mapsto g^\alpha\defeq \alpha^{-1}g\alpha$. We claim that the subgroup $G^\alpha\defeq\{g^\alpha:g\in G\}=\{\alpha^{-1}g\alpha:g\in G\}$ coincides with the subgroup $\delta[G]\subseteq\Hol(N)\subseteq \Sym(N)$.
The equality $G^\alpha=\delta[G]$ will follow as soon as we check that $\alpha^{-1}g\alpha=\delta_g$ for all $g\in G$. The latter equality is equivalence to the equality $g\alpha=\alpha\delta_g$. Given any element $g\in G=N\rtimes G_o$, find unique elements $v\in N$ and $a\in G_o$ such that $g=va$. Taking into account that $a\in G_o$, we conclude that $a(o)=o$ and $a^{-1}(o)=o$. For every $u\in N$ we have the equality
$$\alpha\delta_g(u)=\alpha\delta_{va}(u)=\alpha(t_v\gamma_a(u))=\alpha(vaua^{-1})=vaua^{-1}(o)=gu(o)=g(\alpha(u)),$$ witnessing that $\alpha\delta_{g}=g\alpha$ and $g=\alpha^{-1}\delta_g\alpha=g^\alpha$. Therefore, $G^\alpha=\delta[G]\subseteq \Hol(N)$ and the permutation group $G$ is isomorphic to the permutation subgroup $G^\alpha=\delta[G]$ of the holomorph $\Hol(N)$ of the group $N$.
\end{proof}

Theorem~\ref{t:affine-type} motivates the following definition.

\begin{definition}\label{d:affine-type} A permutation group $G\subseteq \Sym(X)$ is defined to be \defterm{of affine type} if it contains a sharply transitive normal Abelian subgroup. 
\end{definition}

\begin{exercise}\label{ex:S4-of affine type} Show that for any set $X$ of cardinality $|X|\le 4$, the permutation group $\Sym(X)$ is of affine type.
\end{exercise}

\begin{remark} By Theorem~\ref{t:affine-type}, every permutation group $G$ of affine type is isomorphic to a subgroup of the holomorph $\Hol(N)$ of a sharply transitive normal Abelian subgroup $N$ of $G$. So, the structure of permutation groups of affine type is more-or-less understandable.
\end{remark}

\section{Groupable liners}

In this section we define groupable affine spaces and prove that an affine space is Thalesian if and only if it is groupable.

\begin{definition} A liner $X$ is called \defterm{groupable} if there exists a sharply transitive permutation group $G\subseteq \Sym(X)$ such that for every line $L$ in $X$ and every point $x\in L$, the set $\vec L_x\defeq \{g\in G:g(x)\in L\}$ is a normal subgroup of the group $G$.
\end{definition}

The following proposition shows that the group $G$ in the definition of a groupable liner often is (elementary) Abelian.

\begin{proposition}\label{p:groupable-commutative} Let $X$ be a liner and $G\subseteq \Sym(X)$ be a sharply transitive permutation group such that for some point $o\in X$ and any line $L\in\mathcal L_o$, the set $\vec L=\{g\in G:g(o)\in L\}$ is a normal subgroup of $G$. If $\|X\|\ne 2$ (and some nonidentity element of the group $G$ has finite order), then the group $G$ is Abelian (and elementary).
\end{proposition}

\begin{proof} If $\|X\|\le 1$, then $X=\{o\}$ and the group $G$ is trivial, by the sharp transitivity. So, assume that $\|X\|\ge 3$. By the sharp transitivity of the group $G$, for any distinct lines $L,\Lambda\in\mathcal L_o$, the normal subgroups $\vec L$ and $\vec \Lambda$ are transversal in the sense that $\vec L\cap\vec\Lambda=\{1_X\}$. The equality $X=\bigcup\mathcal L_o$ implies $G=\bigcup_{L\in\mathcal L_o}\vec L$. Given any elements $x,y\in G$, find lines $L,\Lambda\in \mathcal L_o$ such that $x\in\vec L$ and $y\in \vec\Lambda$. If $L\ne\Lambda$, then the normality of the subgroups $\vec L$ and $\vec\Lambda$ ensures that $xyx^{-1}y^{-1}\in \vec L\cap\vec\Lambda=\{1_X\}$ and hence $xy=yx$. So, asume that $L=\Lambda$. Since $\|X\|\ge 3$, there exists a line $A\in\mathcal L_o\setminus\{L\}$. Choose any element $a\in\vec A\setminus\{1_X\}$ and find a line $B\in\mathcal L_o$ such that $xa\in\vec B$. It follows from $a\in \vec A\setminus\{1_X\}=\vec A\setminus \vec L$ that $xa\notin\vec L$ and hence $ya=ay$ and $y(xa)=(xa)y$, which implies $(yx)a=y(xa)=(xa)y=x(ay)=x(ya)=(xy)a$ and $xy=yx$. This shows that the group $G$ is Abelian.

Now assume that some nonidentity element $g\in G$ has finite order $p$. We can assume that this order is the smallest possible and hence $p$ is a prime number. Find a line $L\in\mathcal L_o$ such that $g\in \vec L$. Given any element $x\in G\setminus\vec L$ find a line $\Lambda\in\mathcal L_o\setminus\{L\}$ such that $x\in\vec\Lambda$. The commutativity of the group $G$ ensures that $(gx)^p=g^px^p=x^p\in\vec\Lambda$. Assuming that $(gx)^p=x^p\ne 1_X$, we conclude that $gx\in \vec \Lambda\setminus\{1_X\}$ and hence $g=(gx)x^{-1}\in\vec \Lambda\cap\vec L=\{1_X\}$, which contradicts the choice of the element $g$. This contradiction shows that $x^p=1_X$ for all $x\in G\setminus\vec L$. Now take any elements $y\in \vec L$ and $x\in\vec G\setminus\vec L$, and observe that $xy\notin\vec L$ and hence $x^p=1=(xy)^p$, which implies $y^p=1$ and  witnesses that the Abelian group $G$ is elementary. 
\end{proof}

\begin{exercise} Find a groupable hyperbolic plane.
\smallskip

{\em Hint:} Consider the subliner $\mathbb Z^2$ of the Euclidean plane $\IR^2$.
\end{exercise}

\begin{problem} Calculate the number of groupable affine liners of order $64$, up to an isomorphism. 
\end{problem}

\begin{theorem}\label{t:Thalesian<=>groupable} An affine space is Thalesian if and only if it is groupable.
\end{theorem}

\begin{proof}  Let $\mathcal L$ be the family of lines in the affine space $X$. If $X$ is Thalesian, then the group $\Trans(X)=\vec X$ of translations is a sharply transitive subgroup of the permutation group $\Sym(X)$, by Propositions~\ref{p:Trans(X)isnormal}, \ref{p:Ax=Bx=>A=B} and Theorem~\ref{t:paraD<=>translation}. By Proposition~\ref{p:flat=>R-submodule}, for every line $L\in\mathcal L$ and point $x\in L$, the set $\vec L\defeq\{\boldsymbol v\in\vec X:x+\boldsymbol v\subseteq L\}$ is an $\IR_X$-submodule of the $\IR_X$-module $\vec X$. In particular, $\vec L$ is a normal subgroup of the additive group $\vec X=\Trans(X)\subseteq\Sym(X)$, witnessing that the Thalesian affine space $X$ is groupable.
\smallskip

To prove the ``if'' part, assume that the affine space $X$ is groupable.
We have to prove that the affine space $X$ is Thalesian. If $\|X\|\ne 3$, then $X$ is Desarguesian and Thalesian, by Corollary~\ref{c:affine-Desarguesian} and Theorem~\ref{t:ADA=>AMA}. So, assume that $\|X\|= 3$, which means that $X$ is a Playfair plane, by Theorem~\ref{t:Playfair<=>}. Since $X$ is groupable, there exists a sharply transitive subgroup $G\subseteq\Sym(X)$ such that for every line $L$ in $X$ and point $x\in L$, the set $\vec L_x\defeq \{g\in G:g(x)\in L\}$ is a normal subgroup of the group $G$.

\begin{claim}\label{cl:vecLx=vecLy} For every line $L$ in $X$ and all points $x,y\in L$, we have $\vec L_x=\vec L_y$.
\end{claim}

\begin{proof} By the transitivity of the permutation group $G$, there exists a  permutation $u\in G$ such that $u(x)=y\in L$. The definition of the set $\vec L_x$ ensures that $u\in \vec L_x$. Since $\vec L_x$ is a subgroup of $G$, for every $g\in\vec L_x$, the permutation $gu$ belongs to the subgroup $\vec L_x$ and hence $g(y)=g(u(x))=gu(x)\in L$, witnessing that $g\in\vec L_y$. Therefore, $\vec L_x\subseteq \vec L_y$. By analogy we can prove that $\vec L_y\subseteq \vec L_x$ and hence $\vec L_x=\vec L_y$.
\end{proof}

For every line $L$ in $X$, consider the subgroup $\vec L\defeq \vec L_x$ where $x$ is any point of $L$. Claim~\ref{cl:vecLx=vecLy} ensures that the subgroup $\vec L=\vec L_x$ does not depend on the choice of the point $x\in L$.

\begin{claim}\label{cl:vecLx=L} For every line $L\in\mathcal L$ and every point $x\in L$, the set $\vec L(x)\defeq \{g(x):g\in \vec L\}$ equals $L$.
\end{claim}

\begin{proof} Observe that $\vec L(x)=\{g(x):g\in \vec L\}\subseteq L$, by definition of the subgroup $\vec L=\vec L_x=\{g\in G:g(x)\in L\}$. On the other hand, for every $y\in L$, by the transitivity of the permutation group $G$, there exists a permutation $g\in G$ such that $g(x)=y$. The definition of the subgroup $\vec L=\vec L_x$ ensures that $g\in \vec L$ and hence $y=g(x)\in \vec L(x)$, witnessing that $L\subseteq\vec L(x)$. Therefore, $L=\vec L(x)$.
\end{proof} 

\begin{claim}\label{cl:L||L'=>vecL=vecL'} For any disjoint lines $L,L'\in\mathcal L$, we have $\vec L=\vec L'$.
\end{claim}

\begin{proof} By symmetry, it suffices to show that $\vec L\subseteq \vec L'$. Fix any points $x\in L$ and $y\in L'$. Let $\mathcal L_y\defeq\{\Lambda\in\mathcal L:y\in\Lambda\}$ and $\mathcal L_y'\defeq\mathcal L_y\setminus\{L'\}$. Since $X$ is a Playfair plane and $y\notin L$, the line $L'$ is a unique line that contains $y$ and is disjoint with the line $L$. Then for every line $\Lambda\in\mathcal L'_y$, the intersection $L\cap \Lambda$ is a singleton $\{z\}$. Claim~\ref{cl:vecLx=vecLy} ensures that $\vec L=\vec L_x=\vec L_z$ and $\vec \Lambda=\vec \Lambda_y=\vec \Lambda_z$. Then the intersection $$\vec L\cap\vec\Lambda=\vec L_z\cap\vec\Lambda_z=\{g\in G:g(z)\in L\cap \Lambda\}=\{g\in G:g(z)=z\}=\{1_X\},$$ by the sharp transitivity of the permutation group $G$. Since $\Lambda\cap L'=\{y\}$, the sharp transitivity of the group $G$ ensures that 
$$\vec \Lambda\cap\vec L'=\vec \Lambda_y\cap\vec L'_y=\{g\in G:g(y)\in \Lambda\cap L'\}=\{g\in G:g(y)=y\}=\{1_X\}.$$
Taking into account that $\vec L\cap \vec\Lambda=\{1_X\}$ for all $\Lambda\in\mathcal L'_y$, and $L'=\{y\}\cup(X\setminus\bigcup\mathcal L_y')$, we conclude that $$
\begin{aligned}
\vec L&\subseteq \{1_X\}\cup(G\setminus\bigcup_{\Lambda\in\mathcal L'_y}\vec \Lambda\}=\{1_X\}\cup\{g\in G:g(y)\notin \textstyle{\bigcup}\mathcal L'_y\}\\&=\{1_X\}\cup\{g\in G:g(y)\in L'\setminus\{y\}\}=\{g\in G:g(y)\in L'\}=\vec L'_y=\vec L'.
\end{aligned}$$
\end{proof}

\begin{claim}\label{cl:permutation=>dilation} Any permutation $g\in G$ is a dilation of the Playfair plane $X$.
\end{claim}

\begin{proof} Given any line $L$, we should prove that $g[L]$ is a line in $X$, parallel to the line $L$. Take any point $x\in L$ and consider the point $y\defeq g(x)\in g[L]$. Claim~\ref{cl:vecLx=L} ensures that $L=\vec L(x)$. If $y\in L$, then $g\in \vec L_x=\vec L$ and $g\vec L=\vec L$ because $\vec L$ is a subgroup of $G$. Then $g[L]=g[\vec L(x)]=g\vec L(x)=\vec L(x)=L$, so $g[L]=L$ is a line parallel to $L$.
\smallskip

Now assume that the point $y=g(x)$ does not belong to the line $L$. Since $X$ is a Playfair plane, there exists a line $\Lambda\in\mathcal L$ such that $y\in \Lambda$ and $\Lambda\parallel L$. Claim~\ref{cl:L||L'=>vecL=vecL'} ensures that $\vec L=\vec \Lambda$. Claim~\ref{cl:vecLx=L} ensures that $L=\vec L(x)$ and $\Lambda=\vec \Lambda(y)$.
Since $\vec L$ is a normal subgroup of the group $G$, $g\vec L=\vec Lg$. Then $$g[L]=g[\vec L(x)]=g\vec L(x)=\vec Lg(x)=\vec L(y)=\vec \Lambda(y)=\Lambda$$and hence $g[L]=\Lambda$ is a line in $X$, parallel to the line $L$.
\end{proof}

\begin{claim}\label{cl:permutation=>translation} The Playfair plane $X$ is translation.
\end{claim}

\begin{proof} Given any distinct points $x,y\in X$, we should find a translation $T:X\to X$ such that $T(x)=y$. By the sharp transitivity of the group $G$, there exists a unique permutation $T\in G$ such that $T(x)=y$. The sharp transitivity of $G$ ensures that the permutation $T$ has no fixed points. By Claim~\ref{cl:permutation=>dilation}, the permutation $T$ is a dilation of the Playfair plane $X$.
Since $T$ has no fixed points, the dilation $T$ is a translation of $X$, witnessing that the affine space $X$ is translation. 
\end{proof}

By Claim~\ref{cl:permutation=>translation}, the Playfair plane $X$ is translation, and by Theorem~\ref{t:paraD<=>translation}, $X$ is Thalesian.
\end{proof} 

A binary operation $\oplus:X\times X\to X$ on a set $X$ is called a \defterm{group operation} if $(X,\oplus)$ is a group.

A subset $A$ of a group $G$ is called a \defterm{coset of a normal subgroup} if $A=gN=Ng$ for some $g\in G$ and some normal subgroup $N$ of the group $G$. 

\begin{proposition}\label{p:groupable<=>} An affine space $X$ is groupable if and only if $X$ admits a (commutative) group operation $\oplus:X\times X\to X$ such that every line in $X$ is a coset of a normal subgroup in the group $(X,\oplus)$.
\end{proposition}

\begin{proof} If an affine space $X$ is groupable, then there exists a sharply transitive subgroup $G\subseteq \Sym(X)$ such that for every line $L$ in $X$ and every point $x\in L$, the set $\vec L_x\defeq\{g\in L:g(x)\in L\}$ is a normal subgroup of $G$. Fix any point $o\in X$. Since $G$ is sharply transitive, the map $\alpha:G\to X$, $\alpha:g\mapsto g(o)$, is bijective. Endow $X$ with a unique binary operation $\oplus:X\times X\to X$ such that the bijection $\alpha:G\to X$ is an isomorphism of the groups $(G,\circ)$ and $(X,\oplus)$. 

We claim that every line $L$ is a coset of some normal subgroup of the group $(X,\oplus)$. Choose any point $x\in L$ and consider the normal subgroup $\vec L_x\defeq\{g\in G:g(x)\in L\}$ of the group $G$. Since $\alpha:G\to X$ is an isomorphism of the groups $(G,\circ)$ and $(X,\oplus)$, the image $H\defeq\alpha[\vec L_x]=\vec L_x(o)$ is a normal subgroup of the group $(X,\oplus)$. We claim that $L=x\oplus H$. Consider the permutation $f\defeq \alpha^{-1}(x)$ and observe that $f(o)=\alpha(f)=x$.

To show that $L\subseteq x\oplus H$, take any point $y\in L$. By the sharp transitivty of the group $G$, there exists a unique permutation $g\in G$ such that $g(x)=y\in L$. Then $g\in\vec L_x$, by the definition of the group $\vec L_x$. The normality of the subgroup $\vec L_x$ in $G$ ensures that $f^{-1}gf\in f^{-1}\vec L_xf=\vec L_x$. Then $y=g(x)=gf(o)=ff^{-1}gf(o)=\alpha(f\circ f^{-1}gf)\in \alpha(f\circ \vec L_x)=\alpha(f)\oplus \alpha[\vec L_x]=x\oplus H$ and hence $L\subseteq x\oplus H$. 

On the other hand, for every $h\in H=\alpha[\vec L_x]$, the permutation $\hbar \defeq\alpha^{-1}(h)$ belongs to the subgroup $\vec L_x$. Since $\vec L_x$ is normal in $G$, $f\hbar f^{-1}\in\vec L_x$ and hence $x\oplus h=\alpha(f\circ \hbar)=f\circ \hbar (o)=f\circ \hbar f^{-1}f(o)=f\circ \hbar f^{-1}(x)\in L$. This proves that $x\oplus H\subseteq L$ and hence $L=x\oplus H$ is a coset of the normal subgroup $H$ in the group $(X,\oplus)$.

By Proposition~\ref{p:groupable-commutative}, the starply transitive permutation group $G$ is commutative and so is its isomorphic copy $(X,\oplus)$.
\smallskip

To prove the ``only if'' part, assume that $X$ admits a group operation $\oplus:X\times X\to X$ such that every line in $X$ is a coset of a normal subgroup in the group $(X,\oplus)$. For every element $a\in X$, consider the left translation $\ell_a:X\to X$, $\ell_a:x\mapsto a\oplus x$, of the group $(X,\oplus)$ and observe that $G\defeq\{\ell_a:a\in X\}$ is a sharply transitive subgroup of the permutation group $\Sym(X)$. We claim that for every line $L$ in $X$ and every point $x\in L$, the set $\vec L_x=\{g\in G:g(x)\in L\}$ is a normal subgroup of the group $G$. By the choice of the group operation $\oplus$, the group $(X,\oplus)$ contains a normal subgroup $H$ such that $L=x\oplus H$. Since the map $\ell_*:X\to G$, $\ell_*:a\mapsto \ell_a$, is a group isomorphism, the set $\ell_*[H]=\{\ell_a:a\in H\}$ is a normal subgroup of the group $G$. We claim that $\vec L_x=\ell_*[H]$. Indeed, for every element $\ell_a\in\vec L_x$, we have $a\oplus x=\ell_a(x)\in L=x\oplus H=H\oplus x$ and hence $a\in H$. This proves the inclusion $\vec L_x\subseteq \ell_*[H]$. On the other hand, for every $a\in H$, we have $a\oplus x\in H\oplus x=x\oplus H=L$ and hence $\ell_a(x)=a\oplus x\in L$, witnessing that $\ell_a\in \vec L_x$. Therefore, $\ell_*[H]\subseteq \vec L_x$ and finally $\vec L_x=\ell_*[H]$ is  normal subgroup of the group $G$. The sharply transitive group $G\subseteq\Sym(X)$ witnesses that the affine space $X$ is groupable.
\end{proof}

Theorem~\ref{t:Thalesian<=>groupable} and Proposition~\ref{p:groupable<=>} imply the following ``group'' characterization of Thalesian affine spaces.

\begin{corollary}\label{c:Thalesian<=>groupable} An affine space $X$ is Thalesian if and only if  $X$ admits a (commutative) group operation $\oplus:X\times X\to X$ such that every line in $X$ is a coset of a normal subgroup of the group $(X,\oplus)$.
\end{corollary}

\section{Groups of line affinities characterizing Thalesian liners}

In this section we characterize Thalesian affine spaces via properties of permutation groups of line affinities on lines in those affine spaces. At first we establish some transitivity properties of the permutations groups $\Sym^\#_X(\Delta)$ and $\Sym^{\Join}_X(\Delta)$ on a line $\Delta$ in an affine space $X$. 

\begin{exercise}\label{ex:line-aff-2-transitive} Let $X$ be a Playfair plane, $\Delta\in\mathcal L_X$ be a line in $X$, and $\boldsymbol h,\boldsymbol v\in \partial X\setminus \{\Delta_\parallel\}$ be two distinct directions on $X$. Show that
\begin{enumerate}
\item the subset $\Sym^\#_X[\Delta;\boldsymbol h,\boldsymbol v]\defeq\{\boldsymbol h_{\Delta,L}\boldsymbol v_{L,\Delta}:L\in \Delta_\parallel\}$ of the permutation group $\Sym^\#_X(\Delta)$ is sharply transitive;
\item the subset $\Sym^{\Join}_X[\Delta;\boldsymbol h,\boldsymbol v]\defeq\big\{\boldsymbol h_{\Delta,L}\boldsymbol v_{L,\Delta}:L\in \mathcal L_X\setminus(\boldsymbol h\cup\boldsymbol v)\big\}$ of the permutation group $\Sym^{\Join}_X(\Delta)$ is sharply $2$-transitive.
\end{enumerate}
\smallskip

{\em Hint:} Look at the pictures:

\begin{picture}(300,125)(-100,-15)

\put(0,0){\line(1,1){70}}
\put(-15,15){\color{red}\line(1,1){70}}
\put(20,20){\color{cyan}\vector(0,1){30}}
\put(20,50){\color{teal}\vector(1,0){30}}

\put(75,65){$\Delta$}
\put(57,87){\color{red}$L$}

\put(20,20){\circle*{3}}
\put(23,17){$x$}
\put(50,50){\circle*{3}}
\put(53,47){$x'$}

\put(110,90){\color{cyan}$\boldsymbol v$}
\put(110,45){\color{teal}$\boldsymbol h$}

\put(190,105){\color{red}$L$}
\put(240,95){$\Delta$}

\put(150,0){\line(1,1){90}}
\put(160,10){\color{cyan}\vector(0,1){45}}
\put(175,25){\color{cyan}\vector(0,1){55}}
\put(160,55){\color{teal}\vector(1,0){45}}
\put(175,80){\color{teal}\vector(1,0){55}}
\put(148,35){\color{red}\line(3,5){40}}

\put(160,10){\circle*{3}}
\put(162,5){$x$}
\put(175,25){\circle*{3}}
\put(178,22){$y$}
\put(205,55){\circle*{3}}
\put(208,52){$x'$}
\put(230,80){\circle*{3}}
\put(233,77){$y'$}
\end{picture}
\end{exercise}

\begin{definition} For a line $\Delta$ in an affine space $X$ and two distinct directions $\boldsymbol h,\boldsymbol v\in \partial X\setminus \{\Delta_\parallel\}$, let 
\begin{itemize}
\item $\Sym^{\#}_X(\Delta;\boldsymbol h,\boldsymbol v)$ be the subgroup of the symmetric group $\Sym(\Delta)$, generated by the sharply transitive set  $\Sym^{\#}_X[\Delta;\boldsymbol h,\boldsymbol v]\defeq\{\boldsymbol h_{\Delta,L}\boldsymbol v_{L,\Delta}:L\in \Delta_\parallel\}$;
\item $\Sym^{\Join}_X(\Delta;\boldsymbol h,\boldsymbol v)$ be the subgroup of the symmetric group $\Sym(\Delta)$, generated by the sharply $2$-transitive set  $\Sym^{\Join}_X[\Delta;\boldsymbol h,\boldsymbol v]\defeq\{\boldsymbol h_{\Delta,L}\boldsymbol v_{L,\Delta}:L\in \mathcal L_X\setminus(\boldsymbol h\cup\boldsymbol v)\}$.
\end{itemize}
\end{definition}

The permutation groups $\Sym^{\#}_X(\Delta;\boldsymbol h,\boldsymbol v)$ and $\Sym^{\Join}_X(\Delta;\boldsymbol h,\boldsymbol v)$ will be used in Theorem~\ref{t:Schleiermacher} characterizing Thalesian affine spaces. It is clear that $$\Sym^{\#}_X(\Delta;\boldsymbol h,\boldsymbol v)\subseteq\Sym^{\Join}_X(\Delta;\boldsymbol h,\boldsymbol v)\subseteq\Sym^{\Join}_X(\Delta)\subseteq\Sym(\Delta).$$

The following two characterizations of Thalesian and Desarguesian affine spaces are reformulations of Theorems~\ref{t:unique-translation} and Corollary~\ref{c:Desarg<=>unique-aff}, respectively.

\begin{theorem}\label{t:Thalesian<=>sharply-transitive} For an affine space $X$, the following conditions are equivalent:
\begin{enumerate}
\item The affine space $X$ is Thalesian.
\item For every line $\Delta$ in $X$, the permutation group $\Sym^\#_X(\Delta)$ is sharply transitive.
\item For every line $\Delta$ in $X$ and distinct directions $\boldsymbol h,\boldsymbol v\in\partial X\setminus \{D_\parallel\}$, the permutation group $\Sym^\#_X(\Delta)$ coincides with the sharply transitive set $\Sym^\#_X[\Delta;\boldsymbol h,\boldsymbol v]$.
\item For every line $\Delta$ in $X$, the permutation group $\Sym^\#_X(\Delta)$ coincides with the sharply transitive set $\Sym^\#_X[\Delta;\boldsymbol h,\boldsymbol v]$ for some distinct directions $\boldsymbol h,\boldsymbol v\in\partial X\setminus \{D_\parallel\}$.
\end{enumerate}
\end{theorem}

\begin{theorem}\label{t:Desarguesian<=>sharply-2transitive} For an affine space $X$, the following conditions are equivalent:
\begin{enumerate}
\item The affine space $X$ is Desarguesian.
\item For every line $\Delta$ in $X$ the permutation group $\Sym^{\Join}_X(\Delta)$ is sharply $2$-transitive.
\item For every line $\Delta$ in $X$ and distinct directions $\boldsymbol h,\boldsymbol v\in\partial X\setminus \{D_\parallel\}$, the permutation group $\Sym^{\Join}_X(\Delta)$ coincides with the sharply $2$-transitive set $\Sym^{\Join}_X[\Delta;\boldsymbol h,\boldsymbol v]$.
\item For some line $\Delta$ in $X$ and some distinct directions $\boldsymbol h,\boldsymbol v\in\partial X\setminus \{D_\parallel\}$, the permutation group $\Sym^{\Join}_X(\Delta)$ coincides with the sharply $2$-transitive set $\Sym^{\Join}_X[\Delta;\boldsymbol h,\boldsymbol v]$.
\end{enumerate}
\end{theorem}

A subset $B$ of a group $G$ is called 
\begin{itemize}
\item \defterm{normal in} $G$ if $x^{-1}Bx=B$ for all $x\in G$;
\item \defterm{Boolean} if $x=x^{-1}$ for all $x\in B$.
\end{itemize}

\begin{proposition}\label{p:1sharp=>Boolean} If an affine space $X$ has scalar corps $\IR_X=\{0,1\}$, then for any line $L$ in $X$, every $1$-sharp normal set in the permutation group $\Sym_X^{\Join}(L)$ is Boolean. 
\end{proposition}

\begin{proof} To derive a contradiction, assume that for some line $L$ in $X$,  the permutation group $\Sym_X^{\Join}(L)$ contains a non-Boolean $1$-sharp normal subset $\mathcal S$. Since $\mathcal S$ is not Boolean, there exists a permutation $S\in\mathcal S$ such that $S\circ S\ne 1_L$. Then there exists an element $o\in X$ whose image $a\defeq S(S(o))$ is distinct from $o$. We claim that the element $e\defeq S(o)$ is distinct from $o$ and $a$. If $e=o$, then $a=S(S(o))=S(e)=S(o)=e=o$, which contradicts the choice of $a$. If $e=a$, then $e=a=S(S(o))=S(e)$ and the bijectivity of $S$ ensures that $o=S^{-1}(e)=e=a$, which again contradicts the choice of $o$. Therefore, $e\ne o$ and $o\ne a\ne e$. We claim that the line triple $ose$ is Desarguesian. Given any permutation $P\in \Sym_X^{\Join}(L)$ with $Poe=oe$, we need to prove that $P(a)=a$.  
Observe that $P^{-1}SP(o)=P^{-1}S(o)=P^{-1}(e)=e=S(o)$ and hence $P^{-1}SP=S$, by the normality and $1$-sharpness of the set $\mathcal S=P^{-1}\mathcal SP$. Then $SP=PS$ and hence $P(a)=PS(e)=SP(e)=S(e)=a$, witnessing that the line triple $oae$ is Desarguesian and its portion $\overvector{oae}$ is a scalar, distinct from $0$ and $1$. Therefore, $\IR_X\ne\{0,1\}$, which contradicts the assumption.
\end{proof}

The following characterization of Thalesian affine spaces was essentially proved by \index[person]{Schleiermacher}Schleiermacher\footnote{{\bf Adolf Schleiermacher} was a German geometer active mainly in the late 1960s and early 1970s in the area of incidence geometry, with particular emphasis on affine and projective planes and their groups of collineations. He received his doctorate in 1969 from Goethe-Universit\"at Frankfurt am Main with the dissertation {\em Projektivit\"aten in projektiven und affinen Ebenen}, supervised by Ruth Moufang and Peter Dembowski. During this period he published several papers in Mathematische Zeitschrift, Archiv der Mathematik, and Abhandlungen aus dem Mathematischen Seminar der Universit\"at Hamburg, including joint work with Karl Strambach on projectivities and free extensions in affine geometry.
Schleiermacher’s contributions form part of the German school of axiomatic and group-theoretic geometry that developed around Frankfurt, Erlangen, and Hamburg in the postwar decades.} \cite{Schleiermacher1970}.

\begin{theorem}[Schleiermacher, 1970]\label{t:Schleiermacher} For a line $\Delta$ in an affine space $X$ and distinct directions $\boldsymbol h,\boldsymbol v\in\partial X\setminus\{\Delta_\parallel\}$, the following conditions are equivalent:
\begin{enumerate}
\item The affine space $X$ is Thalesian.
\item The permutation group $\Sym^{\#}_X(\Delta)$ is Abelian,  sharply transitive, and normal in $\Sym^{\Join}_X(\Delta)$.
\item The permutation group $\Sym^{\Join}_X(\Delta)$ contains a sharply transitive normal subgroup.
\item The sharply transitive set $\Sym^\#_X[\Delta;\boldsymbol h,\boldsymbol v]$ is a normal Abelian subgroup of the group $\Sym^{\Join}_X(\Delta;\boldsymbol h,\boldsymbol v)$;
\item The permutation group $\Sym^{\Join}_X(\Delta;\boldsymbol h,\boldsymbol v)$ is of affine type.
\item The sharply $2$-transitive set  $\Sym^{\Join}_X[\Delta;\boldsymbol h,\boldsymbol v]$ is contained in some permutation subgroup $G\subseteq\Sym(\Delta)$ of affine type.
\end{enumerate}
If the affine space $X$ is not Boolean, then the conditions \textup{(1)--(6)} are equivalent to
\begin{enumerate}
\item[\textup{(7)}] For some line $L$ in $X$, the permutation group $\Sym^{\Join}_X(L)$ contains a non-Boolean $1$-sharp normal subset. 
\item[\textup{(8)}] The permutation group $\Sym^{\Join}_X(\Delta)$ contains a non-Boolean $1$-sharp normal subset. 
\item[\textup{(9)}] $\IR_X\ne\{0,1\}$.
\end{enumerate}
\end{theorem}

\begin{proof} If $\|X\|>4$, then the affine space $X$ is Desarguesian and Thalesian, by Corollary~\ref{c:affine-Desarguesian} and Theorem~\ref{t:ADA=>AMA}. In this case,  Theorems~\ref{t:grp-line-trans} and \ref{t:Thalesian<=>sharply-transitive} ensure that the conditions (1)--(5) hold and hence are equivalent. So, assume that $\|X\|=3$ and hence $X$ is a Playfair plane.
\smallskip

The implication $(1)\Ra(2)$ follows from Theorems~\ref{t:grp-line-trans} and \ref{t:Thalesian<=>sharply-transitive}; $(2)\Ra(3)$ and $(4)\Ra(5)\Ra(6)$ are trivial.
\smallskip

$(2)\Ra(4)$ Assume that the permutation group $\Sym_X^\#(\Delta)$ is Abelian, sharply transitive and normal in the group $\Sym^{\Join}_X(\Delta)$. 
The sharp transitivity of the sets $\Sym_X^{\#}[\Delta;\boldsymbol h,\boldsymbol v]\subseteq \Sym^\#_X(\Delta)$ implies that the sharply transitive set $\Sym_X^{\#}[\Delta;\boldsymbol h,\boldsymbol v]=\Sym^\#_X(\Delta)$ is an Abelian normal subgroup in the group $\Sym_X^{\Join}(\Delta)$ and also in its subgroup $\Sym_X^{\Join}(\Delta;\boldsymbol h,\boldsymbol v)$.
\smallskip


$(3)\Ra(6)$ Assume that the permutation group $\Sym^{\Join}_X(\Delta)$ contains a sharply transitive normal subgroup $N$. If the group $N$ is Boolean, then it is an Abelian sharply transitive normal subgroup of $\Sym^{\Join}_X(\Delta)$, witnessing that the permutation group $\Sym^{\Join}_X(\Delta)$ is of affine type. If $N$ is not Boolean, then $N$ is a non-Boolean $1$-sharp normal set in $\Sym^{\Join}_X(\Delta)$ and we can apply Proposition~\ref{p:1sharp=>Boolean} to conclude that $\IR_X\ne\{0,1\}$. By Theorem~\ref{t:RXne01=>paraD}, the affine space $X$ is Thalesian. By Theorems~\ref{t:grp-line-trans} and \ref{t:Thalesian<=>sharply-transitive}, $\Sym^\#_X(\Delta)$ is a sharply transitive normal Abelian subgroup of $\Sym^{\Join}_X(\Delta)$, witnessing that the permutation group $\Sym^{\Join}_X(\Delta)$ is of affine type. 
\smallskip

$(6)\Ra(1)$ Assume that the sharply $2$-transitive set  $\Sym_X^{\Join}[\Delta;\boldsymbol h,\boldsymbol v]$ is contained in some permutation group $G\subseteq \Sym(\Delta)$ of affine type. By Definition~\ref{d:affine-type}, the group $G$ contains a sharply transitive normal Abelian subgroup $N\subseteq G$. 

Fix any point $o\in \Delta$. The sharp transitivity of the group $N$ ensures that the map $\alpha:N\to \Delta$, $\alpha:g\mapsto g(o)$, is bijective. Endow $\Delta$ with the binary operation $\oplus:\Delta\times \Delta\to \Delta$, defined by $x\oplus y=\alpha(\alpha^{-1}(x)\circ\alpha^{-1}(y))$. Then the bijection $\alpha:N\to \Delta$ is an isomorphism of the magmas $(N,\circ)$ and $(\Delta,\oplus)$, which implies that $(\Delta,\oplus)$ is an Abelian group. The bijective map $\alpha:N\to \Delta$ induces the group isomorphism $(\cdot)^\alpha:\Sym(\Delta)\to\Sym(N)$, $(\cdot)^\alpha:g\mapsto \alpha^{-1}g\alpha$. Theorem~\ref{t:affine-type} ensures that the group $G^\alpha\defeq\{\alpha^{-1}g\alpha:g\in G\}$ is a subgroup of the holomorph $\Hol(N)$ of the Abelian group $N$. Since $\alpha:N\to \Delta$ is an isomorpism of the Abelian groups $(N,\circ)$ and $(\Delta,\oplus)$, the group $G$ is a subgroup of the holomorph $\Hol(\Delta)$ of the Abelian group $(\Delta,\oplus)$.

Consider the bijective map $\gamma:\Delta\times \Delta\to X$, assigning to any pair $(x,y)\in \Delta\times \Delta$ the unique point $\gamma(x,y)\in \Aline x{\boldsymbol v}\cap\Aline y{\boldsymbol h}$. 
The bijection $\gamma:\Delta\times \Delta\to X$ allows us to identity the plane $X$ with the Cartesian square $\Delta\times \Delta$ of the line $\Delta$. For any point $p\in X$, the points $x,y\in \Delta$ forming the ordered pair $(x,y)=\gamma^{-1}(p)\in \Delta\times \Delta$ are called the {\em horizontal} and {\em vertical coordinates} of the point $p$. Since the line $\Delta$ carries the structure of an Abelian group, its square $\Delta\times \Delta$ also has the structure of an Abelian group. Then the bijecion $\gamma:\Delta\times \Delta\to X$ allows us to define a binary operation $\oplus:X\times X\to X$ on $X$ such that $\gamma$ is an isomorphism of the Abelian groups $\Delta\times \Delta$ and $(X,\oplus)$.

We claim that every line $\Lambda$ in the liner $X$ is a coset of some normal subgroup in the Abelian group $(X,\oplus)$. This is clear if $\Lambda\in\boldsymbol h\cup\boldsymbol v$. So, assume that $\Lambda\notin \boldsymbol h\cup\boldsymbol v$. Consider the permutation $g\defeq \boldsymbol h_{\Delta,\Lambda}\boldsymbol v_{\Lambda,\Delta}\in \Sym^{\Join}_X(\Delta;\boldsymbol u,\boldsymbol v)\subseteq G$ of the line $\Delta$. Observe that a point $p\in X$ with coordinates $(x,y)=\gamma^{-1}(p)\in \Delta\times \Delta$ belongs to the line $\Lambda$ if and only if $y=g(x)$. So, $\Lambda$ is the graph of the permutation $g$. Since $g\in G\subseteq\Hol(\Delta)$, the permutation is an affine transformation of the Abelian group $(\Delta,\oplus)$ and hence $g=\ell_b\vec A$ for some left translation $\ell_b:\Delta\to \Delta$, $\ell_b:x\mapsto b\oplus x$,  and some automorphism $A\in\Aut(\Delta)$ of the Abelian group $(\Delta,\oplus)$. The  translation $\ell_b$ of the group $\Delta$ induces the translation $T:X\to X$ of the group $X$, assigning to a point $p\in X$ with coordinates $(x,y)$, the point $T(p)$ with coordinates $(x,b\oplus y)$. On the other hand, the automorphism $A$ of the Abelian group $\Delta$ considered as a relation on $\Delta$ is a subgroup of the Abelian group $\Delta\times \Delta$. Then its preimage $H\defeq\gamma^{-1}[A]$ is a (normal) subgroup of the Abelian group $(X,\oplus)$. Therefore, $\Lambda=T[H]$ is a coset of the normal subgroup $H$ in the Abelian group $(X,\oplus)$. By Corollary~\ref{c:Thalesian<=>groupable}, the affine space $X$ is Thalesian.
\smallskip

$(1)\Ra(7)$ Assume that the affine space $X$ is Thalesian but not Boolean. Then $X$ contains a non-Boolean parallelogram $uowe$ in $X$. Consider the line $L\defeq\Aline oe$ and  find a line $\Lambda$ in the plane $\Pi\defeq\overline{\{u,o,w,e\}}$ such that $w\in\Lambda$ and $\Lambda\parallel L$. The directions $\boldsymbol u\defeq(\Aline ou)_{\parallel}$ and $\boldsymbol w\defeq(\Aline ow)_\parallel$ determine the line translation $T\defeq \boldsymbol u_{L,\Lambda}\boldsymbol w_{\Lambda,L}:L\to L$ such that $T(o)=e$. Since the parallelogram $uowe$ is not Boolean, $T(e)\ne o$, witnessing that $T\circ  T\ne 1_\Delta$ and the set  $\Sym^\#_X(L)$ is not Boolean. By Theorems~\ref{t:grp-line-trans} and \ref{t:Thalesian<=>sharply-transitive}, the non-Boolean set $\Sym^\#_X(L)$ is $1$-sharp and normal in the permutation group $\Sym^{\Join}_X(L)$. 
\smallskip

$(7)\Ra(8)$ Assume that for some line $L$ in $X$, the permutation group $\Sym^{\Join}_X(L)$ contains a non-Boolean $1$-sharp normal subset. By Proposition~\ref{p:Sym-Join-isomorpic}, the permutation groups $\Sym^{\Join}_X(L)$ and $\Sym^{\Join}_X(\Delta)$  are isomorphic. Then the permutation group $\Sym_X^{\Join}(\Delta)$ also contains a non-Boolean $1$-sharp normal subset.
\smallskip

The implication $(8)\Ra(9)$ follows from Proposition~\ref{p:1sharp=>Boolean}, and $(9)\Ra(1)$ follows from Theorem~\ref{t:RXne01=>paraD}.
\end{proof}

\section{Primitive permutation groups}

In this section we discuss an important notion of a primitive permutation group, which is intermediate between the notions of sharply transitive and transitive permutation groups. We prove some basic properties of primitive permutation groups, which will be applied in the proof of the fundamental Burnside's Dichotomy~\ref{t:Burnside2} for $2$-transitive finite permutation groups. Primitive permutation groups are defined with the help of blocks.

\begin{definition} A subset $B$ of a set $X$ is called a \defterm{$G$-block} for a permutation group $G\subseteq \Sym(X)$ if for every $g\in B$, either $g[B]=B$ or $g[B]\subseteq X\setminus B$.  A $G$-block $B$ is \defterm{trivial} if $B=X$ of $|B|\le 1$. A permutation group $G\subseteq \Sym(X)$ is called \defterm{primitive} if every $G$-block is trivial. 
\end{definition}

The following proposition establishes some elementary properties of primitive permutation groups.

\begin{proposition}\label{p:primitive-permutation-group} Let $G\subseteq S_X$ be a permutation group and $H$ be a nontrivial normal subgroup in $G$.
\begin{enumerate}
\item If $G$ is $2$-transitive, then it is primitive.
\item If $G$ is primitive, then the subgroup $H$ of $G$ is transitive.
\item If $H$ is transitive and Abelian, then $H$ is sharply transitive.
\item If $G$ is primitive and $H$ is Abelian, then $H$ is contained in every nontrivial normal subgroup of $G$.
\item If $G$ is $2$-transitive and $H$ is sharply transitive, then any two nonidentity elements of the group $H$ are conjugated in the group $G$ and hence they have the same order.
\end{enumerate}
\end{proposition}

\begin{proof} 1. Assuming that $G$ is $2$-transitive but not primitive, we can find a nontrivial $G$-block $B\subseteq X$. Since $B$ is not trivial, there exist distinct points $x,y\in B$ and $z\in X\setminus B$. Since $X$ is $2$-transitive, there exists an element $g\in G$ such that $g(x)=x$ and $g(y)=z$. Then $B\ne g[B]\not\subseteq X\setminus B$ witnessing that $B$ is not a block.
\smallskip

2. Assume that some nontrivial normal subgroup $H$ of $G$ is not transitive.
Since $H$ is not trivial, it contains a permutation $\hbar\in H$ such that $\hbar(o)\ne o$ for some point $o\in X$. Since $H$ is not transitive, the set $Ho\defeq\{h(o):h\in H\}$ is not equal to $X$. Also $|Ho|\ge|\{o,\hbar(o)\}|=2$. We claim that $Ho$ is a $G$-block. Given any $g\in G$ with $g[Ho]\cap Ho\ne \varnothing$, we should prove that $g[Ho]=Ho$. Since $g[Ho]\cap Ho\ne \varnothing$, there exist elements $\varphi,\psi\in H$ such that $g\varphi(o)=\psi(o)$. Then $o=\psi^{-1}g^{-1}\varphi(o)$. By the normality of the subgroup $H$ in $G$, $g\psi^{-1}  g^{-1}\in H$ and hence $g(o)=g\psi^{-1}g^{-1}\varphi(o)\in HHo=Ho$ and $Hg(o)=HHo=Ho$. The normality of the group $H$ in $G$ ensures that $g[Ho]=Hg(o)=Ho$,  witnessing that $Ho$ is a nontrivial $G$-block and the permutation group $G$ is not primitive. Therefore, all non-trivial normal subgroups of the primitive permutation group $G$ are transitive.
\smallskip

3. Given any non-trivial Abelian transitive permutation group $H\subseteq \Sym(X)$, we shall prove that $H$ is sharply transitive. Given any $g,h\in H$ with $g(x)=h(x)$ for some $x\in X$, we should show that $g(y)=h(y)$ for all $y\in X$. By the transitivity of $H$, there exists an element $f\in H$ such that $f(x)=y$. The commutativity of the group $H$ ensures that $g(y)=gf(x)=fg(x)=fh(x)=hf(x)=h(y)$ and hence $g=h$, witnessing that the Abelian group $H\subseteq G$ is sharply transitive.
\smallskip

4. Let $A$ be a nontrivial normal Abelian subgroup of a primitive group $G$ and let $H$ be a nontrivial normal subgroup of $G$. We have to prove that $A\subseteq H$. By the preceding items, the permutation group $A$ is sharply transitive. Fix any point $o\in X$. First we prove that $A\cap H\ne\{1_X\}$. To derive a contradiction, assume that $A\cap H=\{1_X\}$. By the normality of the subgroups $A,H$ in $G$, for any $a\in A$ and $h\in H$, the commutator $aha^{-1}h^{-1}$ belongs to the group $AA\cap HH=\{1_X\}$, which means that $ah=ha$.

Since the group $H$ is not trivial, it contains an element $h\ne 1_X$. By the transitivity of the permutation group $A$, there exists an element $a\in A$ such that $a(o)=h(o)$. We claim that $a=h$. Indeed, for any $x\in X$, by the transitivity of the subgroup $A$, there exists an element $b\in A$ such that $b(o)=x$. Then $h(x)=hb(o)=bh(o)=ba(o)=ab(o)=a(x)$, by the commutativity of subgroup $A$. Therefore, $h=a\in H\cap A=\{1_X\}$, which is a desired contradiction showing that $A\cap H\ne\{1_X\}$. By the item (2), the non-trivial normal subgroup $A\cap H$ of $G$ is transitive. Now the sharp transitivity of the group $A$ implies that $A\cap H=A$ and hence $A\subseteq H$.
\smallskip

5. Assume that $G$ is $2$-transitive and $H$ is sharply transitive. Let $h,\hbar$ be arbitrary nonidentity elements in $H$. Fix any element $o\in X$. Since $H$ is sharply transitive, $h(o)\ne o\ne \hbar(o)$. Since the permutation group $G$ is $2$-transitive, there exists an element $g\in G$ such that $g(o)=o$ and $g(h(o))=\hbar(o)$. Since $H$ is a normal subgroup of $G$, the element $g^{-1}\hbar g$ belongs to the group $H$. Observe that $g^{-1}\hbar g(o)=g^{-1}\hbar(o)=h(o)$ and hence $g^{-1}\hbar g=h$, by the sharp transitivity of the permutation group $H$. Since the elements $h,\hbar$ of the group $H$ are conjugated, they have the same order.
\end{proof}

Proposition~\ref{p:primitive-permutation-group} implies the following characterization.

\begin{proposition}\label{p:EAb<=>sharptrans} For a nontrivial finite normal subgroup $H$ of a $2$-transitive permutation group $G\subseteq \Sym(X)$, the following conditions are equivalent:
\begin{enumerate}
\item $H$ is Abelian;
\item $H$ is sharply transitive;
\item $H$ is Abelian and elementary.
\end{enumerate}
\end{proposition}

\begin{proof} By Proposition~\ref{p:primitive-permutation-group}(1,2), the $2$-transitive permutation group $G$ is primitive and the normal subgroup $H$ of $G$ is transitive.
\smallskip

$(1)\Ra(2)$ If $H$ is Abelian, then $H$ is sharply transitive, by Proposition~\ref{p:primitive-permutation-group}(3).
\smallskip

$(2)\Ra(3)$ If $H$ is sharply transitive, then by Proposition~\ref{p:primitive-permutation-group}(5), any nonidentity elements of the group $H$ have the same order $p$, which is a prime number (because the nontrivial finite group $H$ does contain an element of prime order). Since all nontrivial elements of the group $H$ have prime order $p$, the group $H$ is elementary. Moreover, the Sylow Theorem~\ref{t:Sylow} ensures that $|H|=p^n$ for some $n\in\IN$. By Theorem~\ref{t:Burnside-center}, the group $H$ has nontrivial center $\mathcal Z(H)$. We claim that $\mathcal Z(H)$ is a normal subgroup of the group $G$. Indeed, for every $z\in \mathcal Z(H)$, $g\in G$ and $h\in H$, we have 
$gzg^{-1}h=g(zg^{-1}hg)g^{-1}=g(g^{-1}hgz)g^{-1}=hgzg^{-1}$, witnessing that $gzg^{-1}\in \mathcal Z(H)$. Therefore, $\mathcal Z(H)$ is a normal subgroup of $G$. By Proposition~\ref{p:primitive-permutation-group}(2), the normal subgroup $\mathcal Z(H)$ of the primitive permutation group $G$ is transitive. Now the sharp transitivity of the group $H$ and the transitivity of the group $\mathcal Z(H)$ imply that $H=\mathcal Z(H)$, which means that the group $H$ is Abelian. Since every nontrivial element of $H$ has prime order $p$, the Abelian group $H$ is elementary Abelian.
\vskip3pt

The implication $(3)\Ra(1)$ is trivial.
\end{proof}

\begin{definition}
A subgroup $H$ of a group $G$ is called
\begin{itemize}
\item \defterm{nontrivial} if $|H|>1$;
\item \defterm{proper} if $H\ne G$;
\item \defterm{normal} if $xHx^{-1}=H$ for all $x\in G$;
\item \defterm{minimal} if $H$ is nontrivial and every nontrivial subgroup of $H$ coincides with $H$;
\item \defterm{maximal} if $H$ is a proper subgroup of $G$ and every proper subgroup $G$ that contains $H$ coincides with $H$;
\item \defterm{minimal normal} if $H$ is nontrivial, normal, and $H$ contains no non-trivial proper subgroup, which is normal in $G$;
\item \defterm{maximal normal} if $H$ is a proper normal subgroup of $G$ and every proper normal subgroup $G$ that contains $H$ coincides with $H$.
\end{itemize}
\end{definition}

\begin{proposition}\label{p:primitive<=>stab-max} A transitive permutation group $G\subseteq \Sym(X)$ is primitive if and only if for every point $x\in X$, its stabilizer $G_x\defeq\{g\in G:g(x)=x\}$ is a maximal subgroup of $G$.
\end{proposition}

\begin{proof} First assume that the permutation group $G$ is primitive. Take any element $x\in X$ and consider any subgroup $H$ of $G$ such that $G_x\subsetneq H$. It follows from $G_x\ne H$ that the set $Hx\defeq\{h(x):h\in H\}$ is not equal to the singleton $\{x\}$. Given any permutation $g\in G$ with $g[Hx]\cap Hx\ne\varnothing$, we shall prove that $g[Hx]=Hx$. Find elements $\varphi,\psi\in H$ such that $g\varphi(x)=\psi(x)$ and observe that $\psi^{-1}g\varphi(x)=x$ and hence $\psi^{-1}g\varphi\in G_x$ and $g\in \psi G_x\varphi^{-1}\subseteq H$. Then $g[Hx]\subseteq HHx=Hx$.
Applying the inverse permutation $g^{-1}$ to the inequality $g[Hx]\cap Hx\ne\varnothing$, we obtain the inequality $Hx\cap g^{-1}[Hx]\ne\varnothing$. Repeating the above argument, we can show that $g^{-1}[Hx]\subseteq Hx$, which is equivalent to $Hx\subseteq g[Hx]$. Since $Hx\subseteq g[Hx]\subseteq H_x$, the set $H_x$ and $g[H_x]$ are equal,  witnessing that $Hx$ is a $G$-block. Since $G$ is primitive, the $G$-block $Hx$ is trivial and hence $Hx=X$. Given any $g\in G$, we have $g(x)\in X=Hx$ and hence $g(x)=h(x)$ for some $h\in H$. Then $h^{-1}g\in G_x\subseteq H$ and hence $g\in hH$, witnessing that $H=G$ and hence the  stabilizer $G_x$ is a maximal nontrivial subgroup of $G$.

Now assume that the group $G$ is not primitive and find a non-trivial $G$-block $B\subseteq X$. Choose any distinct points $x,y\in B$ and observe that $G_x$ is a proper subgroup of the group $G_B\defeq\{g\in G:g[B]=B\}$ (because $G_B$ contains any permutation $g\in G$ with $g(x)=y$). Choose any point $z\in X\setminus B$ and find a permutation $g\in G$ such that $g(x)=z$. Since $g\notin G_B$, the proper subgroup $G_B$ of $G$ witnesses that the subgroup $G_x$ is not maximal.
\end{proof}

A group $G$ is defined to be \defterm{simple} if it contains no proper nontrivial normal subgroup.

\begin{proposition}\label{p:normal-simple=>minimal} Let $N$ be a nontrivial finite normal subgroup of a $2$-transitive permutation group $G\subseteq \Sym(X)$. If the group $N$ is simple, then $N$ is contained in any nontrivial normal subgroup $H$ of $G$.
\end{proposition}

\begin{proof} To derive a contradiction, assume that $N\not\subseteq H$ for some non-trivial normal subgroup $H$ of $G$. Since $N$ is simple, the proper normal subgroup $N\cap H$ of $N$ is trivial. The normality of the groups $N,H$ implies that $x{\cdot}y=y{\cdot}x$ for all $x\in N$ and $y\in H$. By Proposition~\ref{p:primitive-permutation-group}, the permutation groups $H$ and $N$ are transitive. We claim that the group $N$ is sharply transitive. Given any permutation $g\in N$ with $g(o)=o$ for some $o\in X$, we need to check that $g(x)=x$ for all $x\in X$. By the transitivity of the permutation group $H$, for every $x\in X$, there exists a permutation $h\in H$ such that $h(o)=x$. Then $g(x)=gh(o)=hg(o)=h(o)=x$, witnessing that $g=1_X$ and the group $N$ is sharply transitive. By Proposition~\ref{p:EAb<=>sharptrans}, the finite group $N$ is Abelian, and by Proposition~\ref{p:primitive-permutation-group}(4), $N\subseteq H$, which contradicts our assumption. 
\end{proof}

\begin{proposition}\label{p:trivial-center=>embedsinAut} Let $N$ be a nontrivial finite normal subgroup of a $2$-transitive permutation group $G\subseteq \Sym(N)$. If the group $N$ is not Abelian, then the subgroup $C_G(N)\defeq\{g\in G:\forall x\in N\;g{\cdot}x=x{\cdot}g\}$ is trivial  and hence $G$ is isomorphic to a subgroup of the automorphism group $\Aut(N)$.
\end{proposition}

\begin{proof} To derive a contradiction, assume that the subgroup $C_G(N)$ is not trivial. The normality of the subgroup $N$ implies the normality of the subgroup $C_G(N)$ in $G$. By Proposition~\ref{p:primitive-permutation-group}, the permutation groups $N$ and $C_G(N)$ are transitive. We claim that the group $N$ is sharply transitive. Given any permutation $g\in N$ with $g(o)=o$ for some $o\in X$, we need to check that $g(x)=x$ for all $x\in X$. By the transitivity of the permutation group $C_G(N)$, for every $x\in X$, there exists a permutation $h\in H$ such that $h(o)=x$. Then $g(x)=gh(o)=hg(o)=h(o)=x$, witnessing that $g=1_X$ and the group $N$ is sharply transitive. By Proposition~\ref{p:EAb<=>sharptrans}, the finite sharply transitive group $N$ is Abelian, which contradicts the assumption.

Consider the homomorphism $\gamma:G\to \Aut(N)$, assigining to every element $g\in G$ the conjugating automorphism $(\cdot)^g:N\to N$, $ 
(\cdot)^g:x\mapsto x^g\defeq gxg^{-1}$. It follows from $G_G(N)=\{1\}$ that the homomorphism $h$ is injective, which implies that the group $G$ is isomorphic to the subgroup $h[G]$ of the automorphism group $\Aut(N)$ of $N$. 
\end{proof}
 
\section{Classification of $2$-transitive finite permutation groups}

The following fundamental dychotomy is a classical result of Burnside \cite{Burnside1897} (see also \cite[Theorem 4.3]{Cameron-PG}). 

\begin{theorems}[Burnside, 1897]\label{t:Burnside2} Every $2$-transitive finite permutation group $G\subseteq\Sym(X)$ contains a unique minimal normal subgroup $N$. The minimal normal subgroup $N$ is either elementary Abelian and sharply transitive or simple and primitive.
\end{theorems}

\begin{proof} By Proposition~\ref{p:primitive-permutation-group}(1), the $2$-transitive permutation group $G$ is primitive. Fix any minimal normal subgroup $N$ of $G$. If $N$ is elementary Abelian, then $N$ is sharply transitive, by Proposition~\ref{p:primitive-permutation-group}(3). So, we assume that $N$ is not elementary Abelian.  In this case we prove that $N$ is primitive and simple. To derive a contradiction, assume that $N$ is not primitive. Then there exists a non-trivial $N$-block $B\subseteq X$ and we can assume $B$ has the smallest possible cardinality. 

\begin{claim} For every $g\in G$, the set $g[B]$ is an $N$-block. 
\end{claim}

\begin{proof} By the normality of $N$ in $G$, for every $f\in N$ the permutation $f^g\defeq g^{-1}fg$ is an element of $N$ and hence either $f^g[B]=B$ and then $f[g[B]]=gf^g[B]=g[B]$ or else $f^g[B]\cap B=\varnothing$ and then $f[g[B]]\cap g[B]=g[f^g[B]]\cap g[B]=g[f^g[B]\cap B]=g[\varnothing]=\varnothing$, witnessing that $g[B]$ is an $N$-block.
\end{proof}

Since $B$ is a minimal $N$-block and the intersection of two $N$-blocks is an $N$-block, for any $f,g\in B$, the $N$-blocks $f[B]$ and $g[B]$ either coincide or have at most one common point. Since the permutation group $G$ is $2$-transitive, for any distinct points $x,y\in X$ there exists a permutation $g\in G$ such that $\{g(x),g(y)\}\subseteq B$ and hence $\{x,y\}\subseteq g^{-1}[B]$. This shows that the family $\mathcal L\defeq\{g[B]:g\in X\}$ of ``lines'' turns the set $X$ into a liner. Moreover, every permutation $g\in G$ is an automorphism of the liner $(X,\mathcal L)$.  If $g\in N$, then for every $L\in \mathcal L$ either $g[L]=L$ or $g[L]\cap L=\varnothing$. This implies that for every line $L\in\mathcal L$ and point $x\in L$, every element $g\in N_x$ in the stabilizer  $N_x\defeq\{g\in N:g(x)=x\}$ of $x$ fixes $L$ in the sense that $g[L]=L$. 

\begin{claim}\label{cl:gNxg=Ny} For any points $x,y\in X$ and a permutation $g\in G$ with $y=g(x)$, we have $gN_xg^{-1}=N_y$.
\end{claim}

\begin{proof} Given any $f\in N_x$, observe that $gfg^{-1}(y)=gf(x)=g(x)=y$, witnessing that $gfg^{-1}\in N_y$ and hence $gN_xg^{-1}\subseteq N_y$. By analogy we can check that $g^{-1}N_yg\subseteq N_x$, which implies $N_y\subseteq gN_xg^{-1}$. Therefore, $gN_xg^{-1}=N_y$.
\end{proof}

\begin{claim}\label{clNxy-trivial} For any distinct points $x,y\in X$, the group $N_{xy}\defeq N_x\cap N_y$ is trivial.
\end{claim}

\begin{proof} Given any distinct points $x,y\in X$, find a unique line $\Aline xy\in\mathcal L$ containing both points $x$ and $y$.

Given any element $g\in N_{xy}$ we have to show that $g(z)=z$ for all $z\in X$. If $z\notin\Aline xy$, then $\{z\}=\Aline xz\cap\Aline zy$ implies $\{g(z)\}=g[\Aline xz]\cap g[\Aline yz]=\Aline xz\cap\Aline yz=\{z\}$. Therefore, $g\in N_p$ for all $p\in X\setminus\Aline xy$. 

If $z\in \Aline xy\setminus \{x,y\}$, then we can choose a point $p\in X\setminus \Aline xy$ (such point $p$ exists because the $N$-block $B$ is not trivial and hence $|\Aline xy|=|B|<|X|$) and observe that $g\in N_x\cap N_p$. Then $\{z\}=\Aline xz\cap \Aline zp$ and $\{g(z)\}=g[\Aline xz]\cap g[\Aline zp]=\Aline xz\cap\Aline zp=\{z\}$. 

Therefore, $g(z)=z$ for all $z\in X$, and the group $N_{x}\cap N_y$ is trivial.
\end{proof}

Taking into account that the group $N$ is not elementary Abelian, we can apply Proposition~\ref{p:EAb<=>sharptrans} and conclude that the permutation group $N$ is not sharply transitive. Then for some point $o\in X$, the stabilizer group $N_o$ is not trivial. 

\begin{claim} The subgroup $N_o$ is malnormal in $G$ in the sense that $gN_og^{-1}\cap N_o=\{1\}$ for all $g\in G\setminus N_o$.
\end{claim}

\begin{proof} It follows from $g\notin N_o$ that the point $y\defeq g(o)$ is distinct from $o$. Claim~\ref{cl:gNxg=Ny} ensures that $gN_og^{-1}=N_y$. Then $gN_og^{-1}\cap N_o=N_y\cap N_o=N_{yo}=\{1\}$, by Claim~\ref{clNxy-trivial}.
\end{proof}

By Frobenius Theorem$^\dag$~\ref{t:Frobenius-Kernel}, the set $$K\defeq\{1\}\cup (G\setminus \bigcup_{g\in G}gN_og^{-1})$$ is a normal subgroup of $G$ such that $G=K\rtimes N_o$. Then $N=(K\cap N)\rtimes N_o$ and hence $K\cap N$ is a nontrivial proper normal subgroup of $N$, which contradicts the minimality of the normal group $N$. This contradiction shows that the normal subgroup $N$ is primitive. 

\begin{claim} The group $N$ is simple.
\end{claim}

\begin{proof} 
Assuming that $N$ is not simple, we can choose a minimal normal subgroup $S\ne N$ in $N$. The minimality of the normal subgroup $N$ in $G$ ensures that the proper subgroup $S$ is not normal in $G$ but $N$ is generated by the union $\bigcup_{g\in G}gSg^{-1}$. Consider the normalizer $G_S\defeq\{g\in G:gSg^{-1}=S\}$ of the subgroup $S$ in $G$. Choose a set $A\subseteq G$ that intersects each coset $gG_S$, $g\in G$, at a single point. We can assume that $A\cap G_S=\{1\}$. Then $\bigcup_{g\in G}gSg^{-1}=\bigcup_{a\in A}aSa^{-1}$. Since the group $N\ne S$ is generated by the set $\bigcup_{a\in A}aSa^{-1}$, the set $A_1\defeq A\setminus\{1\}$ is not empty.

Observe that for every $a\in A$, the conjugation $(\cdot)^a:N\to N$, $(\cdot)^a:x\mapsto x^a\defeq axa^{-1}$, is an automorphism of the group $N$. Then the minimality of $S$ ensures that the subgroup $S^a\defeq aSa^{-1}$ is a minimal normal subgroup of $N$.

For any distinct elements $a,b\in A$, the normal subgroups $S^a$ and $S^b$ of $N$ are distinct and hence the normal subgroup $S^a\cap S^b$ of $N$ is a proper subgroup of the minimal normal group $S^a$, which implies $S^a\cap S^b=\{1\}$. Then for every $x\in S^a$ and $y\in S^b$ we have $xyx^{-1}y^{-1}\in S^bS^b\cap S^aS^a=S^b\cap S^a=\{1\}$ and hence the elements $x,y$ commute. Therefore, $N=\oplus_{a\in A}S^a$ is the direct product of the subgroups $S^a$, $a\in A$. Then every normal subgroup of $S$ is a normal subgroup of $N=\oplus_{a\in A}S^a$. The minimality of the normal subgroup $S$ ensures that $S$ contains no proper nontrivial normal subgroups, which means that the group $S$ is simple and so are its isomorphic copies $S^a$ for all $a\in A$. We claim that the simple group $S$ is Abelian. To derive a contradiction, assume that $S$ is not Abelian.
%

Consider the normal subgroup $\Pi\defeq\oplus_{a\in A_1}S^a$ of the group $N$ and observe that $N=S\times \Pi$. Let $\pi_S:N\to S$ and $\pi_\Pi:N\to\Pi$ be the projections of the group $N=S\times \Pi$ onto its factors $S$ and $\Pi$. 

By Proposition~\ref{p:primitive-permutation-group}(2), for every $a\in A$, the normal subgroup $S^a$ of the primitive group $N$ is transitive. We claim that $S^a$ is sharply transitive. Indeed, take any point $o\in X$ and consider its stabilizer $S^a_o=\{g\in S^a:g(o)=o\}$ in the group $S^a$. Take any permutation $g\in S^a_o$ and choose any $b\in A\setminus\{a\}$. By the transitivity of the group $S^b$, for every $x\in X$, there exists a permutation $f\in S^b$ such that $f(o)=x$. Taking into account that the elements $g\in S^a_o\subseteq S^a$ and $f\in S^b$ commute, we conclude that $g(x)=g(f(o))=fg(o)=f(o)=x$, witnessing that $g$ is the identity permutation of $X$ and hence $S^a_o=\{1_X\}$, which implies that that the transitive group $S^a$ is sharply transitive. 

Now fix any point $o\in X$ and consider its stabilizer $N_o\defeq\{g\in N:g(o)=o\}$ in the primitive permutation group $N$. By Proposition~\ref{p:primitive<=>stab-max}, $N_o$ is a maximal subgroup of $N$. On the other hand, the sharp transitivity of the subgroups $S^a$, $a\in A$, implies that $N_o\cap S^a=\{1_X\}$ for all $a\in A$. The normality of the subgroup $S$ in $N$ ensures that $SN_o$ is a subgroup of $N$. The maximality of the subgroup $N_o$ in $N$ ensures that $SN_o=N$ and hence $N=S\rtimes N_o$ is a semidirect product of the normal group $S$ and the group $N_o$. Then the restriction $\pi_{\Pi}{\restriction}_{N_o}:N_o\to\Pi$ of the projection $\pi_\Pi:N\to\Pi$ is an isomorphism, its inverse $\varphi:\Pi\to N_o$ also is an isomorphism, and the composition $\psi\defeq \pi_S\varphi:\Pi\to S$ is a homomorphism from the group $\Pi=\oplus_{a\in A_1}S^a$ to the group $S$. For every $a\in A_1$ the restriction $\psi_a\defeq\psi{\restriction}_{S^a}:S^a\to S$ is either injective or constant (because the group $S^a$ is simple). 

Assuming that for some $a\in A_1$, the homomorphism $\psi_a$ is constant, we conclude that $\{1\}=\psi_a[S^a]=\psi[S^a]=\pi_S\varphi[S]$ and hence  $\varphi[S^a]\subseteq N_o\cap\pi^{-1}_S(1)=N_o\cap\Pi$. Taking into account that projection $\pi_\Pi:N\to\Pi$ is identity on the set $\Pi$, we conclude that $S^a=\pi_\Pi\varphi[S^a]=\varphi[S^a]\subseteq N_o$, which contradicts $N_o\cap S^a=\{1\}$. This contradiction shows that for every $a\in A_1$, the homomorphism $\psi_a:S^a\to S$ is injective and hence bijective (because $|S^a|=|S|$).

We claim that $|A_1|=1$. To derive a contradiction, assume that the set $A_1$  contains two distinct elements $a,b$. Since the homomorphisms $\psi_a:S^a\to S$ and $\psi_b:S^b\to S$ are isomorphisms, for any elements $x,y\in S$ we can find elements $x_a\in S^a$ and $y_b\in S^b$ such that $\psi_a(x_a)=x$ and $\psi_b(y_b)=y$. Since the elements $x_a\in S^a$ and $y_b\in S^b$ commute in the group $N$,
$$x\cdot y=\psi(x_a)\cdot \psi(y_b)=\psi(x_a\cdot y_b)=\psi(y_b\cdot x_a)=\psi(y_b)\cdot\psi(x_a)=y\cdot x,$$
witnessing that the group $S$ is commutative, which contradicts our assumption. This contradiction shows that $|A_1|=1$ and hence $|A|=2$.

By Lemma${}^\dag$~\ref{l:3distinct-orders}, there are elements $s,t,u\in S\setminus\{1_X\}$ of pairwise distinct orders. Fix any point $o\in X$. Since the permutation group $S$ is sharply transitive, $o\notin\{s(o),t(o),u(o)\}$. Since the permutation group $G$ is 2-transitive, there exist elements $g,h\in G$ such that $g(o)=o=h(o)$, $g(s(o))=u(o)$ and $h(u(o))=t(o)$. Assuming that $g\in C_S$, we conclude that $gsg^{-1}\in S$ and $gsg^{-1}(o)=gs(o)=u(o)$ and hence $gsg^{-1}=u$, by the sharp transitivity of the permutation group $S$. Then the elements $s$ and $u$ have the same order, which contradicts the choice of $u$. This contradiction shows that $g\notin C_S$. By analogy we can prove that $h\notin C_S$. Since $|A|=2$, the subgroup $C_S$ has index $2$ in the group $G$ and hence the element $\hbar\defeq hg$ belongs to the subgroup $G_S$, which implies $\hbar s\hbar^{-1}\in S$. Observe that $\hbar s\hbar^{-1}(o)=hgsg^{-1}h^{-1}(o)=hgs(o)=hu(o)=t(o)$. The sharp transitivity of $S$ implies that $\hbar s\hbar^{-1}=t$ and hence the elements $s,t$ have the same order, which contradicts the choice of $s,t$. This contradiction shows that the simple group $S$ is Abelian. Then the group $N=\oplus_{a\in A}S^a$ is elementary Abelian, which contradicts our assumption. This contradiction show that the group $N$ is simple.
\end{proof}

By Propositions~\ref{p:primitive-permutation-group} and \ref{p:normal-simple=>minimal}, the Abelian or simple minimal normal subgroup $N$ of the primitive permutation group $G$ is contained in any nontrivial normal subgroup of $G$, which implies that $N$ is a unique minimal normal subgroup of the group $G$.
\end{proof}

\begin{lemmas}\label{l:3distinct-orders} Any non-commutative finite simple group $G$ contains elements $a,b,c\in G\setminus\{1\}$ of pairwise distinct orders.
\end{lemmas}

\begin{proof} We have to prove that the set $S$ of all possible orders of nonidentity elements of $G$ contains at least three elements. To derive a contradiction, assume that $|S|\le 2$. By Sylow Theorem~\ref{t:Sylow}, the set $S$ contains all prime divisors of the cardinal $|G|$. If $|G|$ has only one prime divisor, then $|G|=p^n$ for some prime number $p$ and some $n\in\IN$. By Theorem~\ref{t:Burnside-center}, the group $G$ has a non-trivial center $\mathcal Z(G)$. Taking into account that $\mathcal Z(G)$ is a normal subgroup of the simple group $G$, we conclude that $G=\mathcal Z(G)$, which contradicts the non-commutativity of the group $G$. This contradiction shows that $|G|$ has at least two prime divisors and hence $S=\{p,q\}$ for two distinct prime numbers $p,q$. Therefore, every nonidentity element of the group $G$ has a prime order (equal to $p$ or $q$), which means that the group $G$ is elementary. By Theorem${}^\dag$~\ref{t:elementary-group<=>}, the elementary group $G$ is not simple.
\end{proof}

Lemma${}^\dag$~\ref{l:3distinct-orders} can be also deduced from the following result, called the Burnside's $p$-$q$ Theorem, see \cite[8.5.3]{Robinson}.

\begin{theorems} A finite group is solvable if its order has at most two prime divisors.
\end{theorems}

\begin{remark} By Theorems~\ref{t:Burnside2}, \ref{t:affine-type}, and Propositions~\ref{p:trivial-center=>embedsinAut}, every $2$-transitive finite permutation group $G$ is isomorphic to a subgroup of the holomorph $\Hol(N)$ of its minimal normal subgroup $N$.
\end{remark} 

\begin{remark} The Classification Theorem for Finite Simple Groups  (CFSG) finally completed at the end of XXth century allowed to classify all finite $2$-transitive groups. The complete lists of finite simple groups  and of finite $2$-transitive permutation groups can be found in the book \cite{Cameron-PG} of \index[person]{Cameron}Cameron\footnote{{\bf Peter Jephson Cameron} (born 1947) is an Australian mathematician noted for his extensive contributions to group theory, combinatorics, coding theory, and model theory. After earning a B.Sc. from the University of Queensland and then a D.Phil. at the University of Oxford as a Rhodes Scholar under Peter M. Neumann, he held academic posts at Oxford and Bedford College before long-term positions at Queen Mary University of London and the University of St Andrews, where he is Emeritus Professor of Mathematics. Cameron has authored over 350 research papers and a series of influential texts on permutation groups, combinatorial designs, logic, and algebra, and his work includes the famous Cameron--Erd\H os conjecture with Paul Erd\H os. Cameron’s research bridges structural algebra and discrete mathematics, particularly through the study of permutation group actions on combinatorial structures. He has received several honours, including the London Mathematical Society’s Whitehead Prize (1979), Euler Medal (2003), Senior Whitehead Prize (2017), and Forder Lectureship, and in 2018 was elected Fellow of the Royal Society of Edinburgh.}.
\end{remark}

The following classification of finite permutation groups containing a sharply $2$-transitive subset is a combined result of \index[person]{Grundh\"ofer}Grundh\"ofer\footnote{{\bf  Theo Grundhöfer} (born in 1955) is a Professor at Institut f\"ur Mathematik  of  Julius-Maximilians-Universit\"at W\"urzburg.}, \index[person]{M\"uller}M\"uller\footnote{{\bf Peter M\"uller} is a Holder of the Chair of Mathematics I at Institut f\"ur Mathematik  of  Julius-Maximilians-Universit\"at W\"urzburg.}, and \index[person]{Nagy}Nagy\footnote{{\bf G\'abor P\'eter Nagy} (born 1972) is a Professor of Bolyai Institute of Mathematics at University of Szeged.}  \cite[Theorem 1.9]{GrundhoferMuller2009} and \cite{MullerNagy2011}.

\begin{theorems}[Grundh\"ofer--M\"uller--Nagy, 2010]\label{t:Grundhofer-Muller-Nagy} If a finite permutation group $G\subseteq \Sym(X)$ contains a sharply $2$-transitive subset, then it satisfies one of the following condition:
\begin{enumerate} 
\item $G$ is of affine type;
\item $G=\Alt(X)$;
\item $G=\Sym(X)$;
\item $|X|=24$ and $G$ is isomorphic to the Mathieu group $M_{24}$.
\end{enumerate}
Moreover, if $|X|\in (2+4\IZ)\cup(3+4\IZ)$, then $G\ne \Alt(X)$.
\end{theorems}

\begin{remark} The Mathieu group $M_{24}$ appearing in Theorem~\ref{t:Grundhofer-Muller-Nagy} is the largest group in the series of sporadic simple groups $M_{11}$, $M_{12}$, $M_{22}$, $M_{23}$, $M_{24}$, discovered by French mathematician \index[person]{Mathieu}\'Emile Mathieu
\footnote{{\bf \'Emile L\'eonard Mathieu} (1835 -- 1890) was a French mathematician best known for his 19th-century discovery of the first sporadic simple groups. In a series of papers published between 1861 and 1873, motivated by the study of multiply transitive permutation groups, he constructed five exceptional finite simple groups, now denoted 
$M_{11},M_{12},M_{22},M_{23},M_{24}$. 
These examples remained isolated curiosities for decades and later assumed central importance in the twentieth-century Classification of Finite Simple Groups. Mathieu spent most of his professional life teaching at the Lyc\'ee Louis-le-Grand in Paris and held no university professorship, publishing comparatively little; 
nonetheless, his work exerted lasting influence on group theory and permutation-group theory.}
in 1861--1873. The Mathieu group $M_{24}$ is the unique permutation group that acts $5$-transitively on a set of $24$ points and is the full automorphism group the \defterm{large Witt design} $S(5,8,24)$, which is a Steiner system that has $759$ blocks of length $8$ such that any $5$ points belong to a single block. The groups $M_{23}$, $M_{22}$, $M_{12}$, $M_{11}$ arise as stabilizers of points or blocks in the large Witt design $S(5,8,24)$. The Mathieu group $M_{23}$ is the stabilizer of a point in $M_{24}$; it is simple and acts 4-transitively on the remaining 23 points. $M_{22}$ is the stabilizer of two points in $M_{24}$; it is simple and acts 3-transitively on 22 points.
The Mathieu group $M_{12}$ is unique simple permutation subgroup of 
acting $5$-transitively on $12$ points, preserving the Steiner system 
$S(5,6,12)$, called the \defterm{small Witt design}. $M_{11}$ is the stabilizer of a point in $M_{12}$; it is simple and acts 4-transitively on 11 points. The Mathieu group $M_{24}$ has order $|M_{24}|=244823040=2^{10}\cdot 3^3\cdot 5\cdot 7\cdot 11\cdot 23$. 
\end{remark}

The classification of $2$-transitive finite permutation groups implies the  following classification of $3$-transitive finite permutation groups. This classification can be found in the books of Cameron \cite{Cameron-PG} and Huber \cite[\S2.2]{Huber2009}, see also the papers of Kantor \cite{Kantor1985}, \cite{Kantor1985h}.

\begin{theorems}\label{t:3-transitive-classification} For every $3$-transitive finite permutation group $G\subseteq \Sym(X)$, the minimal normal subgroup $N$ of $G$ is one of the following:
\begin{enumerate}
\item the Boolean group $\IZ_2^n$ for some $n\in\IN$;
\item the projective special linear group $\mathrm{PSL}(2,q)$ for some prime power $q=|X|-1\ge 4$;
\item the alternating group $\Alt(X)$;
\item one of the Mathieu groups $M_{11}$, $M_{12}$, $M_{22}$, $M_{23}$, $M_{24}$.  
\end{enumerate}
\end{theorems}

\begin{remark} The projective special liner group $\mathrm{PSL}(2,q)$ is the group of M\"obius transformations of the projective $\infty$-extension $F\cup\{\infty\}$ of a $q$-element field $F$, which can be identified with the projective line of order $q$. This group is studied in details in Section~\ref{s:Mobious}.
\end{remark}

\section{Groups of line affinities in non-Thalesian affine spaces}

In this section we apply the classification theorem for $2$-transitive finite permutation groups to classification of the groups of line translation on lines in non-Thalesian finite affine spaces.

\begin{theorems}\label{t:non-Thales-line-aff} Let $\Delta$ be a line in a non-Thalesian finite affine space $X$, and $\boldsymbol h,\boldsymbol v\in\partial X\setminus\{\Delta_\parallel\}$ be two distinct directions on $X$. For every subgroup $G\subseteq \Sym(\Delta)$ containing the sharply $2$-transitive set $\Sym_X^{\Join}[\Delta;\boldsymbol h,\boldsymbol v]$, one of the following conditions holds:
\begin{enumerate}
\item $G=\Alt(\Delta)$;
\item $G=\Sym(\Delta)$;
\item $|\Delta|=24$ and the permutation group $G$ is isomorphic to the Mathieu group $M_{24}$.
\end{enumerate}
Moreover, if $|\Delta|\in (2+4\IZ)\cup(3+4\IZ)$, then $G=\Sym(L)$.
\end{theorems}

\begin{proof} By Theorem~\ref{t:Schleiermacher}, the permutation group $G$ is not of affine type. Since the group $G$ contains the sharply $2$-transitive set 
$\Sym_X^{\Join}[\Delta;\boldsymbol h,\boldsymbol v]$, we can apply Theorem~\ref{t:Grundhofer-Muller-Nagy} to complete the proof.
\end{proof}

\begin{corollarys}\label{c:portions=3} If a finite affine space $X$ is not Thalesian, then  $|\dddot X_{\Join}|=3$ and $\IR_X=\{0,1\}$.
\end{corollarys}

\begin{proof} Assume that an affine space $X$ is not Thalesian. By Theorem~\ref{t:RXne01=>paraD}, $\IR_X=\{0,1\}$. By Corollary~\ref{c:4-Pappian}, $|X|_2\ge 7$. 
To see that $|\dddot X_{\Join}|=3$, it suffices to show that for any distinct points $o,e\in X$ and any points $x,y\in \Aline oe\setminus\{o,e\}$, there exists a line affinity $A$ such that $Aoxe=oye$. By Theorem~\ref{t:non-Thales-line-aff},  the permutation group  $\Sym^{\Join}_X(\Delta)$ on line $\Delta\defeq \Aline oe$ is equal to $\Sym(\Delta)$, $\Alt(\Delta)$ or $M_{24}$. All these permutation groups are $3$-transitive and hence $\Sym^{\Join}_X(\Delta)$ contain a permutation $A$ of the line $\Delta$ such that $Aoxe=oye$, witnessing that $\dddot X_{\Join}=\{0,1,\overvector{oxe}\}$ and hence $|\dddot X_{\Join}|=3$. 
\end{proof} 

\begin{remark}\label{r:Sym-Hetman} The following table (prepared by Ivan Hetman) presents the structure of the permutation groups $\Sym^\#_X(\Delta)$ and $\Sym^{\Join}_X(\Delta)$ on lines $\Delta$ in all seven affine planes $X$ of order $9$.

\begin{center}
\renewcommand{\arraystretch}{1.1}
\begin{tabular}{|c|c|c|}
\hline
$X$&$\Sym^\#_X(\Delta)$&$\Sym^{\Join}_X(\Delta)$\\
\hline
{\tt Desarg}&\phantom{\LARGE $|$}$C_3^2$\phantom{\LARGE $|$}&$C_3^2\rtimes C_8$\\
{\tt Thales}&$C_3^2$&$C_3^2\rtimes Q_8$\\
{\tt Hall}&$C_3^2$ or $A_9$&$A_9$\\
{\tt Hughes}&$(C_3^3\,{:}\,C_3)\,{:}\,C_2$ or $A_9$&$S_9$\\
{\tt dhall}&$A_9$ or $S_9$&$S_9$\\
{\tt hughes}&$A_9$ or $S_9$&$S_9$\\
{\tt hall}&$S_9$&$S_9$\\
\hline
\end{tabular}
\end{center}
\end{remark}

\chapter{Pappian liners}

In this section we study \index[person]{Pappus}Pappian\footnote{{\bf Pappus of Alexandria} (290 -- 350 AD) was a Greek mathematician of late antiquity known for his Synagoge ($\Sigma \upsilon \nu\alpha\gamma\omega\gamma\acute\eta$) or Collection and for Pappus's hexagon theorem in projective geometry. Almost nothing is known about his life except for what can be found in his own writings, many of which are lost. Pappus apparently lived in Alexandria, where he worked as a mathematics teacher to higher level students, such one named Hermodorus.
 The Collection, his best-known work, is a compendium of mathematics in eight volumes, the bulk of which survives. It covers a wide range of topics that were part of the ancient mathematics curriculum, including geometry, astronomy, and mechanics. Pappus was active in a period generally considered one of stagnation in mathematical studies, where he stands out as a remarkable exception. In many respects, his fate strikingly resembles that of Diophantus', originally of limited importance but becoming very influential in the late Renaissance and Early Modern periods.} liners and their relation to Desarguesian liners.

\section{The Pappus Axiom}

Theorem~\ref{t:RX-commutative} and Corollary~\ref{c:parallel-lines<=>} imply that the corps $\IR_X$ of a Desarguesian affine space $X$ is commutative if and only if $X$ satisfies the \index{Affine Pappus Axiom}\defterm{Affine Pappus Axiom}:
\begin{itemize}
\item[] {\em for any concurrent lines $L,L'$ in $X$ and any distinct points $a,b,c\in L\setminus L'$ and $a',b',c'\in L'\setminus L$ with $\Aline a{b'}\cap \Aline{a'}b=\varnothing=\Aline b{c'}\cap\Aline{b'}c$, we have $\Aline a{c'}\cap\Aline {a'}c=\varnothing$.}
\end{itemize}

\begin{picture}(140,90)(-150,-15)

{\linethickness{1pt}
\put(20,0){\color{red}\line(-1,1){20}}
\put(60,0){\color{red}\line(-1,1){60}}
}

\put(0,0){\line(1,0){70}}
\put(75,-4){$L$}
\put(0,0){\line(0,1){70}}
\put(4,65){$L'$}
\put(0,60){\color{blue}\line(1,-2){30}}
\put(20,0){\color{blue}\line(-1,2){20}}
\put(0,20){\color{cyan}\line(3,-2){30}}
\put(60,0){\color{cyan}\line(-3,2){60}}

\put(0,0){\circle*{2}}
\put(-5,-10){$o$}
\put(20,0){\circle*{2}}
\put(17,-10){$a$}
\put(30,0){\circle*{2}}
\put(27,-10){$b$}
\put(60,0){\circle*{2}}
\put(57,-10){$c$}
\put(0,20){\circle*{2}}
\put(-9,18){$c'$}
\put(0,40){\circle*{2}}
\put(-9,38){$b'$}
\put(0,60){\circle*{2}}
\put(-9,58){$a'$}
\end{picture}

Therefore, the Affine Pappus Axiom is responsible for the commutativity of the corps of scalars of a Desarguesian  affine space.  In Section~\ref{s:Pappian-completion} we shall prove that an affine space $X$ satisfies the Affine Pappus Axiom if and only if its spread completion $\overline X$ satisfies the \defterm{Projective Pappus Axiom}:
\begin{itemize}
\item[] {\em for every concurrent lines $L,L'$ in $\overline X$, every distinct points $a,b,c\in L\setminus L'$, $a',b',c'\in L'\setminus L$ the points $x\in \Aline a{b'}\cap\Aline {a'}b$, $y\in \Aline b{c'}\cap\Aline {b'}c$ and $z\in \Aline a{c'}\cap\Aline {a'}c$ are collinear.}
\end{itemize}

\begin{picture}(150,170)(-130,-15)

\put(-15,0){\line(1,0){155}}
\put(145,-5){$L$}
\put(0,0){\line(1,1){135}}
\put(140,135){$L'$}
\put(0,0){\line(-1,-1){10}}
{\linethickness{1pt}
\put(60,0){\line(1,0){60}}
\put(30,30){\line(1,1){90}}
}
\put(90,0){\color{blue}\line(-2,1){60}}
\put(90,0){\color{blue}\line(1,4){30}}
\put(120,0){\color{teal}\line(-1,1){60}}
\put(120,0){\color{cyan}\line(-3,1){90}}
\put(60,0){\color{blue}\line(0,1){60}}
\put(60,0){\color{cyan}\line(1,2){60}}

{\linethickness{1pt}
\put(60,15){\color{red}\line(4,1){36}}
}

\put(0,0){\circle*{3}}
\put(60,0){\circle*{3}}
\put(52,-8){$a$}
\put(90,0){\circle*{3}}
\put(87,-10){$b$}
\put(120,0){\circle*{3}}
\put(115,-8){$c$}

\put(30,30){\circle*{3}}
\put(25,33){$a'$}
\put(60,60){\circle*{3}}
\put(55,63){$b'$}
\put(120,120){\circle*{3}}
\put(115,123){$c'$}

\put(60,15){\color{blue}\circle*{3}}
\put(52,10){\color{blue}$z$}
\put(68.6,17.2){\color{cyan}\circle*{3}}
\put(64,22){\color{cyan}$y$}
\put(96,24){\color{teal}\circle*{3}}
\put(99,23){\color{blue}$x$}
\end{picture}

By analogy with the Desargues Axiom which unifies the Affine Desargues Axiom for affine spaces and the Projective Desargues Axiom for projective spaces, the Affine Pappus Axiom for affine spaces and the Projective Pappus Axiom for projective spaces can be unified in a single Pappus Axiom, introduced in the following definition.

\begin{definition} A liner $X$ is called \index{Pappian liner}\index{liner!Pappian}\defterm{Pappian} if it satisfies the \index{Pappus Axiom}\index{Axiom!Pappus}\defterm{Pappus Axiom}:
\begin{itemize}
\item[] {\em for any concurrent lines $L,L'$ in $X$, and distinct points $a,b,c\in L\setminus L'$, $a',b',c'\in L' \setminus L$, the set $T\defeq(\Aline a{b'}\cap\Aline {a'}b)\cup(\Aline a{c'}\cap\Aline {a'}c)\cup(\Aline b{c'}\cap\Aline {b'}c)$ has rank $\|T\|\in\{0,2\}$}.
\end{itemize}
\end{definition}

\begin{proposition}\label{p:Papp<=>4-Papp} A projective liner $X$ is Pappian if and only if every $4$-long plane in $X$ is Pappian.
\end{proposition}

\begin{proof} The ``only if'' part is trivial. To prove the ``if'' part, assume that every $4$-long plane in $X$ is Pappian. To prove that $X$ is Pappian, take any concurrent lines $L,L'\subseteq X$ and distinct points $a,b,c\in L\setminus L'$ and $a',b',c'\in L'\setminus L$. By the projectivity of the liner $X$, there exist points $x\in \Aline ab\cap\Aline {a'}{b'}$, $y\in \Aline bc\cap\Aline{b'}{c'}$ and $z\in \Aline ac\cap\Aline {a'}{c'}$. We have to prove that the points $x,y,z$ are collinear. 

Let $o$ be the unque common point of the lines $L,L'$, and let $M$ be a maximal $3$-long flat in $X$. By Corollary~\ref{c:Avogadro-projective}, the $3$-long projective liner $M$ is $2$-balanced and hence $|M|_2$ is a well-defined cardinal. By Lemma~\ref{l:ox=2}, for every points $m\in M$ and $n\in X\setminus M$, we have $\Aline mn=\{m,n\}$, which implies $\{o,a,b,c,a',b',c',x,y,z\}\subseteq M$. Then $|M|_2\ge |L|\ge |\{o,a,b,c,\}|=4$.  Since all $4$-long planes in $X$ are Pappian, the $4$-long plane $\overline{L\cup L'}\subseteq M\subseteq X$ is Pappian and hence the points $x,y,z$ are collinear.
\end{proof}

\begin{proposition}\label{p:Pappian-minus-flat} If $Y$ is a Pappian projective liner, then for every flat $H\subseteq X$, the subliner $X\defeq Y\setminus H$ is Pappian.
\end{proposition}

\begin{proof} To show that the liner $X$ is Pappian, take any concurrent lines $L,L'\subseteq X$ and distinct points $a,b,c\in L\setminus L'$ and $a',b',c'\in L'\setminus L$. We have to prove that for the set $$T\defeq(\Aline a{c'}\cap \Aline {a'}c)\cup(\Aline b{c'}\cap\Aline {b'}c)\cup(\Aline a{c'}\cap\Aline{a'}c),$$ the intersection $T\cap X$ has rank $\|T\cap X\|=\{0,2\}$ in the liner $X$. Since the lines $L,L'$ are concurrent, the set $P\defeq\overline{L\cup L'}$ is a plane in $Y$. Since the liner $Y$ is projective, the set $T$ has cardinality $|T|=3$. Since the liner $Y$ is Pappian, the nonempty set $T$ has rank $\|T\|=2$ in $Y$ and hence $\|T\cap X\|\le 2$ in $X$.

 If $|T\cap H|\ge 2$, then $T\subseteq H$ because $\|T\|=2$ and $H$ is a flat. In this case  the set $T\cap X=T\setminus H$ is empty and hence has rank $0$ in $X$. If $|T\cap H|\le 1$, then $|T\cap X|\ge 3-1=2$ and hence $\|T\cap X\|=|T\cap X|=2$. In both cases, we obtain $\|T\cap X\|\in\{0,2\}$, witnessing that the liner $X$ is Pappian.
\end{proof} 


\begin{proposition}\label{p:Pappus=>APA} Every Pappian liner $X$ satisfies the Affine Pappus Axiom.
\end{proposition}

\begin{proof} Given any concurrent lines $L,L'$ in $X$ and distinct points $a,b,c\in L\setminus L'$ and $a',b',c'\in L'\setminus L$ with $\Aline a{b'}\cap\Aline {a'}b=\varnothing=\Aline b{c'}\cap\Aline {b'}c$, we should prove that  $\Aline a{c'}\cap\Aline{a'}{c}=\varnothing$. Assuming that $\Aline a{c'}=\Aline{a'}c$, we conclude that $\{a\}=L\cap\Aline a{c'}=L\cap\Aline{a'}c=\{c\}$, which contradicts the choice of the distinct points $a,c$. This contradiction shows that $\Aline a{c'}\ne\Aline{a'}{c}$ and hence $|\Aline a{c'}\cap\Aline{a'}{c}|\le 1$. Since the intersections $\Aline a{b'}\cap\Aline {a'}b$ and $\Aline b{c'}\cap\Aline {b'}c$ are empty, the set 
$T\defeq (\Aline a{b'}\cap\Aline {a'}b)\cup(\Aline b{c'}\cap\Aline {b'}c)\cup(\Aline a{c'}\cap\Aline{a'}{c})$ coincides with the intersection $\Aline a{c'}\cap\Aline{a'}{c}$ and hence $\|T\|\le|T|=|\Aline a{c'}\cap\Aline {a'}c|\le 1$. The Pappus Axiom ensures that $\|T\|\in\{0,2\}$. Taking into account that $\|T\|\le 1$, we conclude that $\|T\|=0$ and hence the set $\Aline a{c'}\cap\Aline {a'}c=T$ is empty, witnessing that $X$ satisfies the Affine Pappus Axiom.
\end{proof}

Applying Theorem~\ref{t:4-long-affine} and Corollary~\ref{c:parallel-lines<=>}, it is easy to prove the following characterization of affine liners satisfying the Affine Pappus Axiom.

\begin{proposition}\label{p:Pappian<=>} An affine liner $X$ satisfies the Affine Pappus Axiom if and only if for every concurrent lines $L,L'\subseteq X$ and distinct points $a,b,c\in L\setminus L'$ and $a',b',c'\in L'\setminus L$, if $\Aline a{b'}\parallel \Aline {a'}b$ and $\Aline b{c'}\parallel \Aline {b'}c$, then  $\Aline a{c'}\parallel \Aline{a'}{c}$.
\end{proposition}







\begin{theorem}\label{t:Pappian=>para-Pappian} If a Playfair liner $X$ satisfies the Affine Pappus Axiom, 
then  for every disjoint lines $L,L'$ in $X$ and distinct points $a,b,c\in L$ and $a',b',c'\in L'$ with $\Aline a{b'}\cap \Aline {a'}b=\varnothing=\Aline b{c'}\cap\Aline {b'}c$, we have $\Aline a{c'}\cap \Aline {a'}c=\varnothing$.
\end{theorem}

\begin{picture}(200,60)(-140,-10)

\put(-15,0){\line(1,0){90}}
\put(80,-3){$L$}
\put(-15,40){\line(1,0){90}}
\put(80,37){$L'$}
\put(0,0){\color{blue}\line(1,2){20}}
\put(0,40){\color{cyan}\line(1,-1){40}}
\put(40,0){\color{blue}\line(1,2){20}}
\put(20,40){\color{cyan}\line(1,-1){40}}
{\linethickness{1pt}
\put(0,0){\color{red}\line(0,1){40}}
\put(60,0){\color{red}\line(0,1){40}}
}

\put(0,0){\circle*{3}}
\put(-3,-8){$a$}
\put(40,0){\circle*{3}}
\put(37,-10){$b$}
\put(60,0){\circle*{3}}
\put(57,-8){$c$}

\put(0,40){\circle*{3}}
\put(-2,44){$c'$}
\put(20,40){\circle*{3}}
\put(18,44){$b'$}
\put(60,40){\circle*{3}}
\put(58,44){$a'$}

\end{picture}

\begin{proof} Assume that a Playfair liner $X$ satisfies the Affine Pappus Axiom. Take any disjoint lines $L,L'\in X$ and distinct points $a,b,c\in L\setminus L'$ and $a',b',c'\in L'\setminus L$ such that $\Aline a{b'}\cap\Aline {a'}b=\varnothing=\Aline b{c'}\cap\Aline {b'}c$. We have to prove that $\Aline ac\cap\Aline{a'}{c'}=\varnothing$. To derive a contradiction, assume that $\Aline ac\cap\Aline{a'}{c'}\ne \varnothing$. Then $\|L\cup L'\|=\|\Aline ac\cup\Aline{a'}{c'}\|=3$, which means that the flat hull $\overline{L\cup L'}$ is a plane in $X$. Applying Theorem~\ref{t:parallel-char}, we can see that $L\parallel L'$,  $\Aline a{b'}\parallel\Aline {a'}b$, $\Aline b{c'}\parallel \Aline {b'}c$, and $\Aline a{c'}\nparallel \Aline{a'}c$. 

Since the liner $X$ is Playfair, there exists a unique line $\ell$ in $X$ such that $a'\in\ell$ and $\ell\parallel \Aline a{c'}$. Since $\Aline a{c'}\cap \Aline b{c'}=\{c'\}$ and $\Aline b{c'}\cap\Aline {b'}c=\varnothing$, Proposition~\ref{p:Proclus-Postulate} ensures that the lines $\Aline a{c'}$ and $\Aline {b'}c$ are concurrent. Since $|\Aline a{c'}\cap\Aline {b'}c|=1$ and $\Aline a{c'}\parallel \ell$, Proposition~\ref{p:Proclus-Postulate} ensures that the lines $\ell$ and $\Aline c{b'}$ have a unique common point $x$. Assuming that $x\in L$, we conclude that $x\in L\cap\Aline c{b'}=\{c\}$ and hence $\Aline {a'}c=\Aline{a'}v=\ell\parallel \Aline a{c'}$, which contradicts our assumption. This contradiction shows that $x\notin L$ and hence $\Aline bx\ne L$. Taking into account that $L$ is the unique line in the Playfair liner $X$ that contains the point $b$ and does not intersect the line $L'$, we conclude that the lines $\Aline bx$ and $L'$ are concurrent.

\begin{claim}\label{cl:para-Pappus} The points $a',b',c'$ do not belong to the line $\Aline bx$.
\end{claim}

\begin{proof} Assuming that $a'\in\Aline bx$, we conclude that $\Aline {a'}b=\Aline bx$ and hence $\Aline a{c'}$ and $\Aline a{b'}$ are two distinct lines that contain the point $a$ and are parallel to the line $\Aline {a'}b=\Aline bx=\ell$, which is not possible in a Playfair liner.

Assuming that $b'\in\Aline bx$, we conclude that $b'\in\Aline bx\cap \Aline {b'}c=\{x\}$ and hence the line $\ell=\Aline {a'}x=\Aline{a'}{b'}=L'$ is not parallel to the line $\Aline a{c'}$, which contradicts the choice of the line $\ell$.

Assuming that $c'\in \Aline bx$, we conclude that $x\in\Aline b{c'}\cap\Aline c{b'}$ and hence the parallel lines $\Aline b{c'}$ and $\Aline c{b'}$ coincide. Then $\{b\}=L\cap\Aline b{c'}=L\cap\Aline c{b'}=\{c\}$, which contradicts the choice of the (distinct) points $b,c$.
\end{proof}

\begin{picture}(400,110)(-20,-30)
\linethickness{0.8pt}
\put(0,0){\line(1,0){140}}
\put(143,-3){$L$}
\put(0,60){\line(1,0){360}}
\put(143,63){$L'$}
\put(0,0){\color{red}\line(0,1){60}}
\put(0,0){\color{cyan}\line(1,1){60}}
\put(0,60){\color{blue}\line(1,-1){60}}
\put(60,0){\color{cyan}\line(1,1){60}}
\put(60,60){\color{blue}\line(1,-1){60}}
\put(120,0){\line(0,1){60}}
\put(120,60){\color{red}\line(-1,-5){10}}
\put(60,0){\line(5,1){300}}
\put(60,0){\line(-5,-1){75}}
\put(0,0){\line(-1,-1){15}}
\put(0,60){\color{red}\line(-1,-5){15}}

\put(0,0){\circle*{3}}
\put(-2,-8){$a$}
\put(60,0){\circle*{3}}
\put(57,-9){$b$}
\put(120,0){\circle*{3}}
\put(118,-8){$c$}
\put(0,60){\circle*{3}}
\put(-2,63){$c'$}
\put(60,60){\circle*{3}}
\put(58,64){$b'$}
\put(120,60){\circle*{3}}
\put(118,63){$a'$}
\put(-15,-15){\circle*{3}}
\put(-18,-23){$y$}
\put(110,10){\circle*{3}}
\put(105,2){$x$}
\put(108,30){$\ell$}

\end{picture}

Since $\Aline a{b'}\cap \Aline b{a'}=\varnothing$ and $\Aline b{a'}\cap\Aline bx=\{b\}$, Proposition~\ref{p:Proclus-Postulate} implies that the lines $\Aline bx$ and $\Aline a{b'}$ have a unique common point $y$.

\begin{claim} The points $x,b,y$ are distinct and do not belong to the line $L'$.
\end{claim}

\begin{proof} Assuming that $x=b$, we conclude that $b=x\in L\cap \Aline c{b'}=\{c\}$, which contradicts the choice of the points $b\ne c$.  Assuming that $y=b$, we conclude that $b=y\in L\cap \Aline a{b'}=\{a\}$, which contradicts the choice of the points $a\ne b$.  Those contradictions show that $x\ne b\ne y$. Assuming that $x=y$, we conclude that $x=y\in\Aline a{b'}\cap\Aline{b'}c=\{b'\}$ and hence $b'=x\in\Aline bx$, which contradicts Claim~\ref{cl:para-Pappus}. Therefore, the points $x,b,y$ are distinct.

The points $b$ does not belong to the line $L'$ because $b\in L$ and the distinct parallel lines $L,L'$ are disjoint. Assuming that $y\in L'$, we conclude that $y\in L'\cap\Aline a{b'}=\{b'\}$ and hence $b'=y\in\Aline bx$, which contradicts Claim~\ref{cl:para-Pappus}.   Assuming that $x\in L'$, we conclude that $x\in L'\cap\Aline c{b'}=\{b'\}$ and hence $b'=x\in\Aline bx$, which contradicts Claim~\ref{cl:para-Pappus}.  
\end{proof}

Since the lines $L'=\overline{\{a',b',c'\}}$ and $\overline{\{x,b,y\}}$ are concurrent and $\Aline y{b'}\cap \Aline b{a'}=\varnothing=\Aline b{c'}\cap \Aline x{b'}$, we can apply the Affine Pappus Axiom and conclude that the lines $\Aline {a'}x$ and $\Aline{c'}y$ are disjoint and hence parallel, by Theorem~\ref{t:parallel-char}. On the other hand, the lines $\Aline {a'}x=\ell$ and $\Aline{c'}a$ are parallel and hence $\Aline{c'}a=\Aline{c'}y$. Then $\{y\}=\Aline{c'}y\cap\Aline a{b'}=\Aline {c'}a\cap\Aline a{b'}=\{a\}$ and hence the line $\Aline bx=\Aline yb=\Aline ab=L$ is parallel to the line $L'$, which is a contradiction completing the proof.
\end{proof}

\begin{corollary}\label{c:Pappian} A Playfair liner $X$ satisfies the Affine Pappus Axiom if and only if for any lines $L,L'\subset X$ and any points $a,b,c\in L\setminus L'$ and $a',b',c'\in L'\setminus L$, if $\Aline a{b'}\parallel \Aline{a'}b$ and $\Aline b{c'}\parallel \Aline {b'}c$, then $\Aline a{c'}\parallel \Aline {a'}c$.
\end{corollary}

\begin{proof} To prove the ``only if'' part, assume that the Playfair liner $X$ satisfies the Affine Pappus Axiom. Take any lines $L,L'\subset X$ and points $a,b,c\in L\setminus L'$ and $a',b',c'\in L'\setminus L$ such that $\Aline a{b'}\parallel \Aline{a'}b$ and $\Aline b{c'}\parallel \Aline {b'}c$. We have to prove that $\Aline a{c'}\parallel \Aline {a'}c$.

If $a=b$, then $\Aline a{b'}\parallel\Aline{a'}{b}$ implies $\Aline a{b'}=\Aline {a'}{b}$ and hence $\{a'\}=L'\cap\Aline {a'}b=L'\cap\Aline {a}{b'}=\{b'\}$ and $a'=b'$. Then $\Aline a{c'}=\Aline b{c'}\parallel \Aline {b'}c=\Aline{a'}c$ and we are done.

By analogy we can prove that $b=c$ implies $\Aline a{c'}\parallel \Aline {a'}c$.

If $a=c$, then $\Aline {a'}b\parallel\Aline a{b'}=\Aline c{b'}\parallel\Aline b{c'}$ and hence $\Aline{a'}b=\Aline b{c'}$ and $\{a'\}=L'\cap\Aline{a'}b=L'\cap\Aline b{c'}=\{c'\}$. Then $\Aline ac=\Aline {a'}{c'}$ and hence $\Aline ac\parallel \Aline {a'}{c'}$. 

Therefore, $|\{a,b,c\}|\le 2$ implies $\Aline a{c'}\parallel \Aline{a'}c$. By analogy we can prove that $|\{a',b',c'\}|\le 2$ implies $\Aline a{c'}\parallel \Aline{a'}c$. 

So, assume that the points $a,b,c,a',b',c'$ are distinct. It follows from $\Aline a{b'}\parallel \Aline {a'}b$ that $$\|\{a,b',a',b\}\|=\|\Aline a{b'}\cup\Aline {a'}b\|=3$$ and hence $P\defeq\overline{L\cup L'}=\overline{\{a,b,a',b'\}}$ is a plane in $X$. By Theorem~\ref{t:Playfair<=>}, the Playfair liner $X$ is $3$-regular, $3$-long and affine. By Proposition~\ref{p:k-regular<=>2ex}, the $3$-regular liner $X$ is $3$-ranked, and by Corollary~\ref{c:parallel-lines<=>}, any  lines in the plane $P$ are concurrent or parallel. In particular, the lines $L$ and $L'$ are concurrent or parallel. If the lines $L,L'$ are concurrent, then $\Aline a{c'}\parallel \Aline {a'}c$, by Proposition~\ref{p:Pappian<=>}. If $L,L'$ are parallel, then  $\Aline a{c'}\parallel \Aline {a'}c$, by Theorem~\ref{t:Pappian=>para-Pappian}.
\smallskip

To prove the ``if'' part, assume that if for any lines $L,L'\subset X$ and any points $a,b,c\in L\setminus L'$ and $a',b',c'\in L'\setminus L$ with $\Aline a{b'}\parallel \Aline{a'}b$ and $\Aline b{c'}\parallel \Aline {b'}c$, we have $\Aline a{c'}\parallel \Aline {a'}c$. To prove that $X$ satifies the Affine Pappus Axiom, take any concurrent lines $L,L'$ in $X$ and distinct points $a,b,c\in L\setminus L'$ and $a',b',c'\in L'\setminus L$ such that $\Aline a{b'}\cap\Aline {a'}b=\varnothing=\Aline b{c'}\cap\Aline {b'}c$. We have to prove that $\Aline a{c'}\cap\Aline{a'}c=\varnothing$. Since the lines $L,L'$ are concurrent, the flat $P\defeq\overline{L\cup L'}$ is a plane.  By Theorem~\ref{t:Playfair<=>} and Proposition~\ref{p:k-regular<=>2ex}, the Playfair liner $X$ is $3$-regular and $3$-ranked, and by Corollary~\ref{c:parallel-lines<=>}, the disjoint lines $\Aline a{b'},\Aline {a'}b
$ in the plane $P$ are parallel. By the same reason, the disjoint lines $\Aline b{c'},\Aline {b'}c$ are parallel. Our assumption implies that  $\Aline a{c'}\parallel \Aline {a'}c$. Assuming that  $\Aline a{c'}\cap \Aline {a'}c\ne\varnothing$, we conclude that $\Aline a{c'}=\Aline{a'}{c}$ and hence $\{a\}=L\cap\Aline a{c'}=L\cap\Aline{a'}c=\{c\}$, which contradicts the choice of the (distinct) points $a,c$. This contradiction shows that $\Aline a{c'}\cap\Aline{a'}c=\varnothing$.
\end{proof}

\section{Pappian affine liners are Desarguesian}

The following important theorem was proved by \index[person]{Hessenberg}Hessenberg
\footnote{{\bf Gerhard Hessenberg} (1874 -- 1925) was a German mathematician who worked in projective geometry, differential geometry, and set theory. His name is usually associated with projective geometry, where he is known for proving that Desargues' theorem is a consequence of Pappus's hexagon theorem, and differential geometry where he is known for introducing the concept of a connection.} \cite{Hessenberg1905} in 1905. Our proof follows the lines of the original Hessenberg's proof but with filling all necessary details, elaborated in the paper of Vlad Pshyk \cite{Pshyk}.

\begin{theorem}[Hessenberg, 1905]\label{t:Hessenberg-affine} If an affine liner $X$ satisfies the Affine Pappus Axiom, then $X$ is Desarguesian.
\end{theorem}

\begin{proof} If $|X|_2\le 3$ or $\|X\|\ne 3$, then the affine liner $X$ is Desarguesian, by Proposition~\ref{p:Steiner+affine=>Desargues} and  Corollary~\ref{c:affine-Desarguesian}. So, we assume that $|X|_2\ge 4$ and $\|X\|=3$. By Theorem~\ref{t:4-long-affine}, the $4$-long affine liner $X$ is regular, and by Theorem~\ref{t:Playfair}, the $4$-long regular affine liner $X$ is Playfair. Since $\|X\|=3$, any distinct lines in $X$ are either parallel or concurrent, according to Corollary~\ref{c:parallel-lines<=>}.

To show that $X$ is Desarguesian, we shall apply  Theorem~\ref{t:ADA<=>} and Lemma~\ref{l:ADA<=>}. Take any concurrent lines $A,B,C$ in $X$ and points $a,a'\in A\setminus(B\cup C)$, $b,b'\in B\setminus(A\cup C)$, $c,c'\in C\setminus(A\cup B)$ such that $\Aline ab\parallel\Aline{a'}{b'}$ and $\Aline bc\parallel \Aline {b'}{c'}$. We have to prove that $\Aline ac\parallel \Aline {a'}{c'}$. 
Let $o$ be the unique common point of the concurrent lines $A,B,C$. 

If $\Aline ab\parallel \Aline bc$, then $\Aline ac=\Aline ab=\Aline bc$. By Theorem~\ref{t:Proclus-lines}, the parallelity relations $\Aline {a'}{b'}\parallel \Aline ab$ and $\Aline {b'}{c'}\parallel \Aline{b}{c}$ imply $\Aline {a'}{b'}\parallel\Aline{b'}{c'}$ and hence the line $\Aline{a'}{c'}=\Aline {a'}{b'}$ is parallel to the line $\Aline ab=\Aline ac$. So, we assume that $\Aline ab\nparallel \Aline bc$, which implies $\Aline{a'}{b'}\nparallel \Aline{b'}{c'}$.

If $b=b'$, then $\Aline ab\parallel \Aline {a'}{b'}$ implies $\Aline ab=\Aline {a'}{b'}$ and hence $\{a\}=\Aline oa\cap\Aline ab=\Aline o{a'}\cap\Aline{a'}{b'}=\{a'\}$. By analogy we can prove that $c=c'$ and hence $\Aline ac=\Aline{a'}{c'}$ and $\Aline ac\parallel\Aline{a'}{c'}$. So, we assume that $b\ne b'$. In this case $a\ne a'$ and $c\ne c'$. 

If $\Aline ac\parallel B\parallel \Aline {a'}{c'}$, then $\Aline ac\parallel \Aline {a'}{c'}$, by Theorem~\ref{t:Proclus-lines}. So, we assume that either $\Aline ac\nparallel B$ or $B\nparallel \Aline {a'}{c'}$. We lose no generality assuming that 
\begin{equation}\label{as:P=>D1}
B\nparallel \Aline{a'}{c'}.
\end{equation} Since $X$ is Playfair, by Theorem~\ref{t:Playfair}, there exists a unique line $\ell_a\subseteq X$ such that $a\in\ell_a$ and $\ell_a\parallel B$. It follows from $\ell_a\parallel B\nparallel \Aline{a'}{c'}$ that the line $\ell_a$ is not parallel to the line $\Aline {a'}{c'}$. By Corollary~\ref{c:parallel-lines<=>}, the lines $\ell_a$ and $\Aline {a'}{c'}$ have a common point $x$.

\begin{claim}\label{cl:P=>D1} $a\ne x$ and $x\notin\Aline{a'}{b'}$.
\end{claim}

\begin{proof} Assuming that $x=a$, we conclude that $a=x\in\Aline {a'}{c'}$ and hence $a\in A\cap\Aline{a'}{c'}=\{a'\}$, which contradicts our assumption. This contradiction shows that $x\ne a$.

Assuming that $x\in \Aline {a'}{b'}$, we conclude that $x\in\Aline{a'}{b'}\cap\Aline{a'}{c'}=\{a'\}$ and hence $\ell_a=\Aline a{a'}=\Aline oa$. Then $\Aline oa=\ell_a\parallel \Aline ob$ and hence $A=\Aline oa=\Aline ob=B$, which contradicts the choice of the point $a\in A\setminus B$. This contradiction shows that $x\notin\Aline{a'}{b'}$.   
\end{proof}

Since, $\ell_a\parallel B$ and $B\cap C=\{o\}$, the lines $\ell_a$ and $C$ have a common point $y$, by Proposition~\ref{p:Proclus-Postulate}.
Assuming that $y=o$ and taking into account that $\ell_a\parallel B$, we conclude that $y=o\in \ell_a\cap B$ and hence $a\in\ell_a=B$, which contradicts the choice of the point $a$. This contradiction shows that 
\begin{equation}\label{eq:P=>Dy} y\ne o.
\end{equation}

Claim~\ref{cl:P=>D1} implies that $\Aline x{b'}\cap\Aline{a'}{b'}=\{b'\}$. Since $\Aline ab\parallel \Aline{a'}{b'}$, Proposition~\ref{p:Proclus-Postulate} ensures that $\Aline x{b'}\cap\Aline ab=\{z\}$ for some point $z$. 

\begin{picture}(200,125)(-130,-9)

{\linethickness{1pt}
\put(0,0){\line(1,0){158}}
\put(162,-3){$C$}
\put(0,0){\line(2,1){158}}
\put(160,77){$A$}
\put(0,0){\line(1,1){100}}
\put(103,100){$B$}

\put(60,0){\line(0,1){30}}
\put(60,0){\line(-25,35){25}}

\put(60,30){\line(-25,5){25}}

\put(120,60){\line(-25,5){50}}

\put(120,0){\line(0,1){60}}
\put(120,0){\line(-25,35){50}}

}

\put(30,0){\line(1,1){90}}
\put(62,45){$\ell_a$}

\put(120,0){\line(0,1){90}}

\put(60,30){\line(-25,5){60}}
\put(0,0){\line(0,1){42}}

\put(120,90){\line(-120,-48){120}}

\put(30,0){\line(-30,42){30}}

\put(0,0){\circle*{3}}
\put(-3,-8){$o$}
\put(30,0){\circle*{3}}
\put(27,-8){$y$}
\put(60,0){\circle*{3}}
\put(58,-8){$c$}
\put(120,0){\circle*{3}}
\put(122,3){$c'$}
\put(60,30){\circle*{3}}
\put(62,24){$a$}
\put(120,60){\circle*{3}}
\put(123,54){$a'$}
\put(120,90){\circle*{3}}
\put(123,88){$x$}

\put(0,42){\circle*{3}}
\put(-2,45){$z$}
\put(70,70){\circle*{3}}
\put(67,73){$b'$}
\put(35,35){\circle*{3}}
\put(32,39){$b$}
\end{picture}

\begin{claim}\label{cl:P=>D2} The points $x,b',z$ do not belong to the line $A$.
\end{claim}

\begin{proof} 
The point $b'$ does not belong to $A$ by the choice of $b'$. Assuming that $x\in A$, we obtain $x\in \ell_a\cap A=\{a\}$, which contradicts Claim~\ref{cl:P=>D1}. Assuming that $z\in A$, we conclude that $z\in\Aline ab\cap A=\{a\}$ and hence $a=z\in \Aline x{b'}$ and $b'\in \Aline a x=\ell_a$. Taking into account that $\ell_a\parallel B$, we conclude that $a\in\ell_a=B$, which contradicts the choice of the point $a$. This contradiction shows that $z\notin A$.
\end{proof}

\begin{claim}\label{cl:P=>D3} The points $o,a,a'$ do not belong to the line $\overline{\{z,b',x\}}$.
\end{claim}

\begin{proof} If $o\in\overline{\{z,b',x\}}$, then $x\in\ell_a\cap\Aline o{b'}=\ell_a\cap B$. Taking into account that $\ell_a\parallel B$, we conclude that $a\in\ell_a=B$, which contradicts the choice of the point $a$. 

If $a\in\overline{\{z,b',x\}}$, then $b'\in\Aline ax=\ell_a$. Taking into account that $\ell_a\parallel B$, we conclude that $a\in\ell_a=B$, which contradicts the choice of the point $a$. 

Assuming that $a'\in\overline{\{z,b',x\}}$, we conclude that 
$z\in \Aline{a'}{b'}\cap\Aline ab$ which is not possible because the lines $\Aline ab$ and $\Aline{a'}{b'}$ are parallel and disjoint.

Therefore, $o,a,a'\notin\overline{\{z,b',x\}}$.
\end{proof}

Claims~\ref{cl:P=>D2} and \ref{cl:P=>D3} ensure that $\{o,a,a'\}\cap\overline{\{z,b',x\}}=\varnothing=\{z,b',x\}\cap\overline{\{o,a,a'\}}$. Since $\Aline ax\parallel \Aline o{b'}$ and $\Aline az\parallel \Aline {a'}{b'}$, Corollary~\ref{c:Pappian} applied to the triples  $o,a,a'$ and $x,b',z$ ensures that 
\begin{equation}\label{eq:P=>D2}
\Aline oz\parallel \Aline {a'}x=\Aline {a'}{c'}.
\end{equation}

\begin{claim}\label{cl:P=>D4} $z\notin A\cup B\cup C$.
\end{claim}

\begin{proof} By Claim~\ref{cl:P=>D3}, $o\ne z$.  Assuming that $z\in A$, we conclude that $\Aline oz=A$. Then $A=\Aline oz\parallel \Aline {a'}{c'}$ and hence $c'\in \Aline{a'}{c'}=A$, which contradicts the choice of $c'$. 

Assuming that $z\in C$, we conclude that $\Aline oz=C$. Then  $C=\Aline oz\parallel \Aline {a'}{c'}$ and hence $a'\in\Aline{a'}{c'}=C$, which contradicts the choice of the point $a'$.

Assuming that $z\in B$, we conclude that $\Aline oz=B$ and hence $B=\Aline oz\parallel \Aline{a'}{c'}$, which contradicts the assumption~(\ref{as:P=>D1}).
\end{proof}

\begin{claim}\label{cl:P=>D5} The points $o,y,c$ do not belong to the line $\Aline ab=\overline{\{a,b,z\}}$.
\end{claim}

\begin{proof} 
Assuming that $o\in\Aline ab$, we conclude that $a\in\Aline ob=B$, which contradicts the choice of the point $a$. Assuming that $c\in\Aline ab$, we conclude that $\Aline ab=\Aline bc$, which contradicts our assumption $\Aline ab\nparallel \Aline bc$. Assuming that $y\in\Aline ab$, we conclude that $y\in\ell_a\cap\Aline ab=\{a\}$ and hence $a=y\in C$, which contradicst the choice of the point $a$.
Those contradictions show that $o,y,c\notin \Aline ab$.
\end{proof}

Two cases are possible.

1. First we assume that $x\notin C$. 

\begin{claim}\label{cl:P=>D6} The points $x,b',z$ do not belong to the line $\overline{\{o,y,c'\}}=C$.
\end{claim}

\begin{proof} 
By the assumption, $x\notin C$. The choice of the point $b'$ ensures that $b'\notin C$. By  Claim~\ref{cl:P=>D4}, $z\notin C$.
\end{proof}

\begin{claim}\label{cl:P=>D7} The points $o,y,c'$ do not belong to the line $\overline{\{x,b',z\}}$.
\end{claim}

\begin{proof} 
By Claim~\ref{cl:P=>D3}, $o\notin\overline{\{x,b',z\}}$. It follows from  $x\notin C=\Aline oy$ that $x\ne y$. Assuming that $y\in\overline{\{x,b',z\}}$ and taking into account that $\ell_a\parallel B$, we conclude that $b'\in\Aline yx=\ell_a$ and hence $a\in\ell_a=B$, which contradicts the choice of the point $a$.

Assuming that $c'\in\overline{\{x,b',z\}}$ and taking into account that $x\in\Aline{a'}{c'}\setminus  C$, we conclude that $b'\in \Aline x{c'}=\Aline {a'}{c'}$ and hence $\Aline {a'}{b'}\parallel \Aline{b'}{c'}$, which contradicts our assumption. 

Those contradictions show that the points $o,y,c'$ do not belong to the line $\overline{\{x,b',z\}}$.
\end{proof}

Claims~\ref{cl:P=>D6} and \ref{cl:P=>D7} ensure that $\{c',o,y\}\cap\overline{\{z,x,b'\}}=\varnothing=\{z,x,b'\}\cap\overline{\{c',o,y\}}$. Since $\Aline yx\parallel \Aline o{b'}$ and $\Aline zo\parallel \Aline{a'}{c'}=\Aline x{c'}$, Corollary~\ref{c:Pappian} implies $\Aline zy\parallel \Aline {b'}{c'}$. Since $\Aline bc\parallel\Aline{b'}{c'}\parallel \Aline zy$, the transitivity of the parallelity proved in Theorem~\ref{t:Proclus-lines}, implies 
that $$\Aline zy\parallel \Aline bc.$$

\begin{claim}\label{cl:P=>D8} The points $a,b,z$ do not belong to the line $C$.
\end{claim}

\begin{proof} The choice of the points $a\in A\setminus (B\cup C)$ and $b\in B\setminus (A\cup C)$ ensure that $a,b\notin C$.  By Claim~\ref{cl:P=>D4}, $z\ne y$. Assuming that $z\in C$ and taking into account that $y\in C$, we conclude that $C=\Aline zy$ and hence $C\parallel \Aline bc$ and $b\in \Aline bc=C$, which contradicts the choice of the point $b$. Therefore, $a,b,z\notin C$.
\end{proof}


 
Claims~\ref{cl:P=>D5} and \ref{cl:P=>D8} ensure that $\{o,y,c\}\cap\overline{\{a,b,z\}}=\varnothing=\{a,b,z\}\cap\overline{\{o,y,c\}}$. Since $\Aline ob\parallel\Aline ya$ and $\Aline yz\parallel \Aline bc$, Corollary~\ref{c:Pappian} implies $\Aline oz\parallel \Aline ac$. Taking into account that $\Aline oz\parallel \Aline {a'}{c'}$, we conclude that $\Aline ac\parallel \Aline {a'}{c'}$.
\smallskip

2. Next, consider the second case: $x\in C$. Then $x\in C\cap\Aline{a'}{c'}=\{c'\}$ and $y\in\ell_a\cap C=\{c'\}=\{x\}$. Then $$\Aline z{c'}=\Aline zx=\Aline {b'}x=\Aline{b'}{c'}\parallel \Aline bc$$ and $\Aline a{c'}=\ell_a\parallel \Aline bo$. The choice of the points $a,b\notin C$ and Claim~\ref{cl:P=>D4} ensure that the points $a,b,z$ do not belong to the line $C=\overline{\{o,c',c\}}$.

\begin{picture}(200,100)(-130,-10)

{\linethickness{1pt}
\put(0,0){\line(1,0){138}}
\put(142,-4){$C$}
\put(0,0){\line(2,1){140}}
\put(140,65){$A$}
\put(0,0){\line(1,1){70}}
\put(72,70){$B$}

\put(60,0){\line(0,1){30}}
\put(60,0){\line(-70,50){35}}

\put(120,0){\line(0,1){60}}
\put(120,0){\line(-70,50){70}}

\put(120,60){\line(-70,-10){70}}

\put(60,30){\line(-70,-10){35}}
}

\put(50,-10){\line(1,1){88}}
\put(140,80){$\ell_a$}

\put(120,60){\line(-70,-10){120}}

\put(60,0){\line(-70,50){60}}

\put(0,0){\line(0,1){43}}

\put(0,0){\circle*{3}}
\put(-3,-8){$o$}
\put(60,0){\circle*{3}}
\put(58,-10){$c'$}
\put(120,0){\circle*{3}}
\put(117,-8){$c$}
\put(60,30){\circle*{3}}
\put(62,24){$a'$}
\put(120,60){\circle*{3}}
\put(122,55){$a$}
\put(0,43){\circle*{3}}
\put(-2,47){$z$}
\put(50,50){\circle*{3}}
\put(46,53){$b$}
\put(25,25){\circle*{3}}
\put(22,29){$b'$}
\end{picture}

\begin{claim}\label{cl:P=>D10} The points $o,c',c$ do not belong to the line $\overline{\{a,b,z\}}=\Aline ab$.
\end{claim}

\begin{proof} 
Claim~\ref{cl:P=>D5} implies that $o,c\notin\overline{\{a,b,z\}}$. Assuming that $c'\in \Aline ab$ and taking into account that $z\in\Aline ab\setminus C$, we conclude that $\Aline ab=\Aline z{c'}\parallel \Aline bc$, which contradicts our assumption.
\end{proof}

Claims~\ref{cl:P=>D4} and \ref{cl:P=>D10} ensure that $\{a,b,z\}\cap\overline{\{o,c',c\}}=\varnothing=\{o,c',c\}\cap\overline{\{a,b,z\}}$. Also $\Aline a{c'}=\ell_a\parallel \Aline ob$ and $\Aline z{c'}\parallel \Aline bc$. Applying Corollary~\ref{c:Pappian} to the triples $a,b,z$ and $o,c',c$, we conclude that $\Aline zo\parallel \Aline ac$ and hence $\Aline ac\parallel \Aline zo\parallel \Aline{a'}{c'}$. 
\end{proof}

\section{Characterizations of Pappian affine liners}

\begin{theorem}\label{t:Papp<=>APA} An affine liner $X$ is Pappian if and only if it satisfies the Affine Pappus Axiom.
\end{theorem}

\begin{proof} The ``only if'' part has been proved in Proposition~\ref{p:Pappus=>APA}. To prove the ``if'' part, assume that an affine liner $X$ satisfies the Affine Pappus Axiom.  By Theorem~\ref{t:affine=>Avogadro}, the affine liner $X$ is $2$-balanced. If $|X|_2\le 3$, then $X$ is Pappian vacuously. So, assume that $|X|_2\ge 4$. In this case, the affine liner $X$ is regular and Playfair, by Theorems~\ref{t:4-long-affine} and \ref{t:Playfair<=>}. By Theorem~\ref{t:Hessenberg-affine}, the affine liner $X$ is Desarguesian.

To prove that $X$ is Pappian, take any concurrent lines $L,L'$ and distinct points $a,b,c\in L\setminus L'$ and $a',b',c'\in L'\setminus L$. We have to show that the set $T=(\Aline a{b'}\cap\Aline {a'}b)\cup(\Aline b{c'}\cap\Aline {b'}c)\cup(\Aline a{c'}\cap \Aline {a'}c)$ has rank $\|T\|\in\{0,2\}$. To derive a contradiction, assume that $\|T\|\notin\{0,2\}$. The Affine Pappus Axiom implies that $\|T\|\ne 1$ and hence $\|T\|=3$. Then the sets $T_{ab}\defeq \Aline a{b'}\cap\Aline {a'}b$, $T_{bc}\defeq \Aline b{c'}\cap\Aline {b'}c$, $T_{ac}\defeq\Aline a{c'}\cap \Aline {a'}c$ are not empty and hence $T_{ab}=\{t_{ab}\}$, $T_{bc}=\{t_{bc}\}$, $T_{ac}=\{t_{ac}\}$ for some points $t_{ab},t_{bc},t_{ac}$. Let $o$ be the unique point of the intersection $L\cap L'$. Since $\|\{t_{ab},t_{bc},t_{ac}\}\|=\|T\|=3$, the point $o$ does not belong to one of the lines $\Aline {t_{ab}}{t_{bc}}$, $\Aline {t_{bc}}{t_{ac}}$, $\Aline {t_{ab}}{t_{ac}}$. We lose no generality assuming that $o$ does not belong to the line $\Lambda\defeq \Aline{t_{ab}}{t_{bc}}$. Consider the plane $P\defeq\overline{L\cup L'}$. Since $X$ is Desarguesian, so is the plane $P$. By Theorem~\ref{t:Desargues-completion}, the spread completion $\overline P$ of the Desarguesian Playfair plane $P$ is a Desarguesian projective plane. Since $X$ is $4$-long, the projective plane $\overline P$ is $5$-long. By Corollary~\ref{c:hyperplanes-automorphism}, there exists an automorphism $A:\overline P\to\overline P$ such that $A[\overline\Lambda]=\overline P\setminus P$. For a point $x\in P$ denote the point $A(x)$ by $x_A$. Since $o\notin\Lambda$, the point $o_A=A(o)$ does not belong to the projective line $\overline P\setminus P$ and hence $o_A\in P$. It is easy to see that $\{a,b,c,a',b',c'\}\cap \Lambda=\varnothing$, which implies that $\{a_A,b_A,c_A,a_A',b_A',c_A'\}\subseteq P$. Since $A$ is an automorphism of the projective plane and $A(t_{ab})\in A[\Lambda]=\overline P\setminus P$, the lines $\Aline {a_A}{b_A'}$ and $\Aline {a_A'}{b_A}$ in $P$ are disjoint. By the same reason, the lines $\Aline {b_A}{c_A'}$ and $\Aline {b'_A}{c_A}$ in $P$ are disjoint. Since $P$ satisfies the Affine Pappus Axiom, the lines $\Aline {a_A}{c'_A}$ and $\Aline {a'_A}{c_A}$ in $P$ are disjoint. On the other hand, $\|T\|=3$ implies that $t_{ac}\notin\Lambda$ and hence $A(t_{ac})\notin A[\Lambda]=\overline P\setminus P$. Since $A(t_{ac})\in \Aline {a_A}{c'_A}\cap\Aline {a'_A}{c_A}$, we obtain a contradiction showing that the affine liner $X$ is Pappian.  
\end{proof}

\begin{theorem}\label{t:Papp<=>Des+RX} An affine space $X$ is Pappian if and only if $X$ is Desarguesian and the corps $\IR_X$ is a field.
\end{theorem}

\begin{proof} To prove the ``only if'' part, assume that an affine space $X$ is Pappian. By  Proposition~\ref{p:Pappus=>APA}, the affine liner $X$ satisfies the Affine Pappus Axiom and by Hessenberg's Theorem~\ref{t:Hessenberg-affine}, the Pappian affine space $X$ is Desarguesian. By Theorem~\ref{t:paraD=>hasRX}, the corps $\IR_X$ is well-defined.  Since affine space $X$ satisfies the Affine Pappus Axiom, Theorem~\ref{t:RX-commutative} ensures that the corps $\IR_X$ is commutative and hence is a field.
\smallskip

To prove the ``if'' part, assume that the affine space $X$ is Desarguesian and the corps $\IR_X$ of $X$ is a field. By Theorem~\ref{t:Desargues<=>3Desargues}, every line triple in the Desarguesian affine space $X$ is Desarguesian. Since the corps $\IR_X$ is commutative,  Theorem~\ref{t:RX-commutative} ensures that the affine space $X$ satisfies the Affine Pappus Axiom. Applying Theorem~\ref{t:Papp<=>APA}, we conclude that the affine space $X$ is Pappian. 
\end{proof}

\begin{theorem}\label{t:fin-affP<=>D} An affine liner $X$ with $|X|_2<\w$ is Pappian if and only if $X$ is Desarguesian.
\end{theorem}

\begin{proof} The ``only if'' part follows from Theorem~\ref{t:Papp<=>Des+RX}. To prove the ``if'' part, assume that $|X|_2<\w$ and the liner $X$ is Desarguesian. It follows from $|X|_2<\w$ that the corps $\IR_X$ is finite. By Wedderburn Theorem~\ref{t:Wedderburn-Witt}, the finite corps $\IR_X$ is a field. By Theorem~\ref{t:Papp<=>Des+RX}, the affine space $X$ is Pappian.
\end{proof}

\begin{proposition}\label{p:Pappian-xy} Let $L,L'$ be two coplanar lines in Pappian Playfair liner $X$ and $a,b,c\in L\setminus L'$, $a',b',c'\in L'\setminus L$ be distinct points such that 
$\Aline a{c'}\cap\Aline {a'}c=\varnothing$ and $\Aline a{b'}\cap\Aline {a'}b=\{x\}$ for some point $x\in X$. Then there exists a unique point $y\in X$ such that $\Aline b{c'}\cap\Aline {b'}c=\{y\}$ and $\Aline xy\cap\Aline a{c'}=\varnothing=\Aline xy\cap\Aline{a'}c$.
\end{proposition}

\begin{proof} By the Pappus Axiom, the non-empty set
$$T\defeq(\Aline a{b'}\cap\Aline{a'}b)\cup(\Aline b{c'}\cap\Aline {b'}c)\cup(\Aline a{c'}\cap\Aline {a'}c)=\{x\}\cup(\Aline b{c'}\cap\Aline{b'}c)\cup\varnothing$$has rank $\|T\|=2$, which implies that the set $\Aline b{c'}\cap\Aline {b'}c$ contains some point $y\ne x$. Assuming that $\Aline b{c'}\cap\Aline {b'}c\ne\{y\}$, we conclude that $\Aline b{c'}=\Aline {b'}c$ and hence $\{b\}=L\cap\Aline b{c'}=L\cap\Aline {b'}c=\{c\}$, which contradicts the choice of the (distinct) points $a,b,c$. This contradiction shows that $\Aline b{c'}\cap\Aline {b'}{c}=\{y\}$ and hence the point $y$ is unique. 

It remains to show that   $\Aline xy\cap\Aline a{c'}=\varnothing=\Aline xy\cap\Aline{a'}c$. By Theorems~\ref{t:Papp<=>APA} and \ref{t:Hessenberg-affine}, the Pappian affine liner $X$ is Desarguesian and so is the plane $P\defeq\overline{L\cup L'}$ in $X$. By Theorem~\ref{t:Desargues-completion}, the spread completion $\overline P$ of the Desarguesian Playfair plane $P$ is a Desarguesian projective plane. Since the Playfair plane $P$ is $3$-long, the projective plane $\overline P$ is $4$-long. By Corollary~\ref{c:hyperplanes-automorphism}, there exists an automorphism $A:\overline P\to\overline P$ such that $A[\overline{\{x,y\}}]=\overline P\setminus P$, where $\overline{\{x,y\}}$ is the line containing the points $x,y$ in the spread completion $\overline P$ of $P$. For a point $p\in P$, denote the point $A(p)$ by $p_A$.  It is easy to see that $\{a,b,c,a',b',c'\}\cap \Aline xy=\varnothing$, which implies that $\{a_A,b_A,c_A,a_A',b_A',c_A'\}\subseteq P$. Since $A$ is a collineation of the projective plane $\overline P$ and $\{A(x),A(y)\}\subseteq  A[\overline{\{x,y\}}]=\overline P\setminus P$, the lines $\Aline {a_A}{b_A'}$ and $\Aline {a_A'}{b_A}$ in $P$ are disjoint. By the same reason, the lines $\Aline {b_A}{c_A'}$ and $\Aline {b'_A}{c_A}$ in $P$ are disjoint. By Corollary~\ref{c:Pappian}, the lines $\Aline {a_A}{c'_A}$ and $\Aline {a'_A}{c_A}$ in $P$ are disjoint. Then $A[\overline{\{a,c'\}}\cap\overline{\{a',c\}}]= \overline{\{a_A,c'_A\}}\cap\overline{\{a'_A,c_A\}}\subseteq \overline P\setminus P=A[\overline{\{x,y\}}]$ and hence $\overline{\{a,c'\}}\cap\overline{\{a',c\}}\subseteq \overline{\{x,y\}}$ in $\overline P$. Taking into account that $\Aline a{c'}\cap \Aline {a'}c=\varnothing$, we conclude that $$\varnothing\ne \overline{\{a,c'\}}\cap\overline{\{a',c\}}\setminus P=\{(\Aline a{c'})_\parallel\}\cap\{(\Aline{a'}c)_\parallel\}=\Aline xy\setminus P=\{(\Aline xy)_\parallel\}$$ and hence
$(\Aline a{c'})_\parallel=(\Aline{a'}c)_\parallel=(\Aline xy)_\parallel$ and $\Aline a{c'}\parallel\Aline xy\parallel \Aline{a'}c$ in $P$. Assuming that the parallel lines $\Aline xy$ and $\Aline a{c'}$ are not disjoint, we conclude that $\Aline xy=\Aline a{c'}$ and then $b'\in \Aline ax\cap \Aline a{b'}\subseteq \Aline a{c'}\cap\Aline a{b'}=\{a\}$, which contradicts the choice of the point $a\in L\setminus L'\subseteq L\setminus\{b'\}$. This contradiction shows that $\Aline xy\cap\Aline a{c'}=\varnothing$. By analogy we can prove that $\Aline xy\cap\Aline {a'}c=\varnothing$.
\end{proof}

\begin{corollary} Let $L,L'$ be two concurrent lines in Pappian affine liner $X$ and $a,b,c\in L\setminus L'$, $a',b',c'\in L'\setminus L$ be distinct points such that 
$\Aline a{c'}\cap\Aline {a'}c=\varnothing$ and $\Aline a{b'}\cap\Aline {a'}b=\{x\}$ for some point $x\in X$. Then there exists a unique point $y\in X$ such that $\Aline b{c'}\cap\Aline {b'}c=\{y\}$ and $\Aline xy\cap\Aline a{c'}=\varnothing=\Aline xy\cap\Aline{a'}c$.
\end{corollary}

\begin{proof} Let $o$ be a unique common point of the concurrent lines $L,L'$. The choice of the points $a,b,c,a',b',c'$ ensures that $|\{o,a,b,c,a',b',c'\}|=7$. By Theorem~\ref{t:affine=>Avogadro}, the affine liner $X$ is $2$-balanced, so the cardinal $|X|_2$ is well-defined and $|X|_2\ge|L|\ge|\{o,a,b,c\}|=4$, which means that the affine liner $X$ is $4$-long and hence regular, by Theorem~\ref{t:4-long-affine}. By Theorem~\ref{t:Playfair<=>}, the $4$-long affine regular liner $X$ is Playfair.
So, we can apply Proposition~\ref{p:Pappian-xy} and find a unique point $y\in X$ such that $\Aline b{c'}\cap\Aline {b'}c=\{y\}$ and $\Aline xy\cap\Aline a{c'}=\varnothing=\Aline xy\cap\Aline{a'}c$.
\end{proof}

\begin{question} Is there any direct geometric proof of Proposition~\ref{p:Pappian-xy}, not involving the complicated machinery of homogeneity of Desarguesian projective planes?
\end{question}




\begin{exercise} Find an example of a Desarguesian affine regular liner, which is not Pappian.
\vskip5pt

{\em Hint:} Consider the quaternion plane $\mathbb H\times\mathbb H$.
\end{exercise}
\section{The complete regularity of Pappian liners}

In this section we shall prove an important and difficult Theorem~\ref{t:Pappian=>compreg} saying that every Pappian proaffine regular liner is completely regular. The proof of this theorem is preceded by six lemmas. In those lemmas we denote by $\mathcal L$ the family of lines in a liner $X$. For a point $x$ of a liner $X$ we denote by $\mathcal L_x\defeq\{L\in\mathcal L:x\in L\}$ the family of lines passing through the point $x$. For a line $L$ in a liner $X$ we denote by $L_\parallel\defeq \{\Lambda\in\mathcal L:\Lambda\parallel L\}$ the family of lines that are parallel to $L$. If $X$ is a $3$-ranked plane, then for every line $L\subseteq X$, the family $L_\parallel$ is equal to the family $\{L\}\cup\{\Lambda\in\mathcal L:\Lambda\cap L=\varnothing\}$, according to Corollary~\ref{c:parallel-lines<=>}.

\begin{lemma}\label{l:Vlad1} Let $X$ be an $\w$-long Pappian Proclus plane and $b,b'\in X$ be points such that the family $\F_b\defeq \{L\in\mathcal L_b:\mathcal L_{b'}\cap L_\parallel=\varnothing\}$ is finite. Then for every line $\Lambda\in \mathcal L\setminus(\mathcal L_{b}\cup\mathcal L_{b'})$, the family $\Lambda_\parallel$ is infinite.
\end{lemma}

\begin{proof} Since $X$ is $\w$-long and $b\notin\Lambda$, there exists a point $a'\in \Lambda\setminus (\Aline b{b'}\cup\bigcup\F_b)$. Consider the line $L'\defeq\Aline {a'}{b'}$ and the finite set $F\defeq \Aline{a'}{b'}\cap\bigcup\F_b$. Since the liner $X$ is Proclus, the set $I\defeq\{x\in L':\Aline bx\cap \Lambda=\varnothing\}$ contains at most one point, and the family $\mathcal L_b\cap L'_\parallel$ contains at most one line.

If $\mathcal L_b\cap L'_\parallel=\varnothing$, then put $F'\defeq \{a',b'\}\cup F\cup I$. If $\mathcal L_b\cap L'_\parallel\ne\varnothing$, then let $L$ be the unique line in the family $\mathcal L_b\cap L'_\parallel$. By Theorem~\ref{t:Proclus-lines}, there exists a unique point $\lambda\in L\cap\Lambda$. Since $X$ is Proclus, the set $J\defeq\{x\in L':\Aline xb\parallel \Aline\lambda {b'}\}$ contains at most one point. Put $F'\defeq \{a',b'\}\cup F\cup I\cup J$. 

\begin{picture}(100,100)(-140,-25)

\put(0,40){\line(1,0){100}}
\put(105,37){$L'$}
\put(0,40){\line(1,-1){40}}
\put(40,0){\line(1,1){40}}
\put(0,-10){\line(0,1){60}}
\put(0,10){\line(4,-1){80}}
\put(30,40){\line(1,-1){50}}
\put(0,10){\line(1,1){30}}
\put(80,40){\line(0,-1){50}}
\put(0,0){\line(1,0){40}}
\put(20,-10){$L$}

\put(-3,55){$\Lambda$}

\put(0,40){\circle*{3}}
\put(-10,38){$a'$}
\put(30,40){\circle*{3}}
\put(28,43){$b'$}
\put(80,40){\circle*{3}}
\put(78,43){$c'$}
\put(0,10){\circle*{3}}
\put(-8,8){$c$}
\put(0,0){\circle*{3}}
\put(-8,-3){$\lambda$}
\put(40,0){\circle*{3}}
\put(38,-10){$b$}
\put(80,-10){\circle*{3}}
\put(78,-18){$a$}

\end{picture}

We claim that for every $c'\in L'\setminus F'$, the intersection $\mathcal L_{c'}\cap \Lambda_\parallel$ is not empty. Since $c'\notin F$, the line $\Aline b{c'}$ does not belong to the family $\F_b$ and hence there exists a line $L_c\in\mathcal L_{b'}$ such that $L_{c}\parallel \Aline b{c'}$. Since $c'\notin I$, the parallel lines $L_{c}$ and $\Aline b{c'}$ are not parallel to the line $\Lambda$ and hence there exists a unique point $c\in L_c\cap\Lambda$. Assuming that $\Aline c{b}\parallel L'$, we conclude that $c=\lambda$ and $c'\in J$, which contradicts the choice of the point $c'$. Therefore, the lines $L'=\Aline {a'}{b'}$ and $\Aline cb$ are concurrent. The choice of the point $a'\notin F$ ensures that $\Aline {a'}b\notin\F_b$ and hence there exists a line $L_a\in\mathcal L_{b'}$ such that $L_a\parallel \Aline {a'}b$. Since the lines $\Aline {a'}{b}$ and $\Aline c{b}$ are concurrent, we can apply Theorem~\ref{t:Proclus-lines} and conclude that there exists a unique point $a\in \Aline {c}{b}\cap L_a$. It follows from $a'\notin\Aline b{b'}$ that the parallel lines $\Aline {a'}b$ and $\Aline a{b'}$ are disjoint and hence $a\ne b$. Assuming that $a=c$, we conclude that $\Aline b{a'}\parallel \Aline {b'}a=\Aline{b'}{c}\parallel \Aline b{c'}$ and hence $\Aline b{a'}=\Aline b{c'}$ and $\{a'\}=L'\cap\Aline b{a'}=L'\cap\ b{c'}=\{c'\}$, which contradicts the choice of the point $c'$. Therefore, the points $a,b,c$ are distinct. It is easy to see that $\{a,b,c\}\cap\Aline {a'}{b'}=\varnothing =\{a',b',c'\}\cap\Aline ab$. By the Proclus Axiom, the set
$$T\defeq(\Aline a{b'}\cap\Aline {a'}b)\cup(\Aline b{c'}\cap\Aline {b'}c)\cup(\Aline a{c'}\cap\Aline{a'}c)=\varnothing\cup\varnothing\cup(\Aline a{c'}\cap\Aline {a'}c)$$
has rank $\|T\|=\{0,2\}$. Since $|T|=|\Aline a{c'}\cap\Aline {a'}c|\le 1$, the set $T$ is empty, which means that the line $\Aline a{c'}$ is parallel to the line $\Aline {a'}c=\Lambda$. Therefore, $\Aline a{c'}\in\mathcal L_{c'}\cap \Lambda_\parallel\ne\varnothing$.

Now we see that the family $$\Lambda_\parallel\supseteq \bigcup_{c'\in L'\setminus F'}(\mathcal L_{c'}\cap\Lambda_\parallel)$$ is infinite.
\end{proof}

\begin{lemma}\label{l:Vlad2} Let $X$ be an $\w$-long Pappian Proclus plane and $b,b'\in X$ be points such that the family $\F_b\defeq \{L\in\mathcal L_b:\mathcal L_{b'}\cap L_\parallel=\varnothing\}$ is finite. Then for every line $\Lambda\in \mathcal L\setminus(\mathcal L_{b}\cup\mathcal L_{b'})$ and every point ${c'}\in X\setminus(\Aline b{b'}\cup\bigcup\F_b)$, the intersection  $\mathcal L_{c'}\cap \Lambda_\parallel$ is not empty.
\end{lemma}

\begin{proof} To derive a contradiction, assume that for some point  ${c'}\in X\setminus(\Aline b{b'}\cup\bigcup\F_b)$, the intersection  $\mathcal L_{c'}\cap \Lambda_\parallel$ is empty.  Then the lines $L'\defeq \Aline {b'}{c'}$ and $\Aline {b}{c'}$ are not parallel to the line $\Lambda$.  
 
Consider the finite set $F\defeq L'\cap\bigcup\F_b$. Since $X$ is a Proclus plane, the set $\mathcal L_{b}\cap L'_\parallel$ contains at most one line. If $\mathcal L_b\cap L'_\parallel=\varnothing$, then put $F'\defeq \{c',b,b'\}\cup F$. If $\mathcal L_b\cap L'_\parallel\ne\varnothing$, then let $L$ be the unique line in the family $\mathcal L_b\cap L'_\parallel$. Since $X$ is a Proclus plane, the set $I\defeq\{u\in L:\Aline u{b'}\parallel \Aline b{c'}\}$ contains at most one point. 
Put $F'=\{c',b,b'\}\cup F\cup I$. 

By Lemma~\ref{l:Vlad1}, the family $\Lambda_\parallel$ is infinite. So, we can find a line $\Lambda'\in\Lambda_\parallel$ such that $\Lambda'\cap F'=\varnothing$. Since $\Lambda'\parallel \Lambda\nparallel L'$, there exists a unique point $a'\in L'\cap\Lambda'$. 

\begin{picture}(100,100)(-140,-25)

\put(0,40){\line(1,0){100}}
\put(105,37){$L'$}
\put(0,40){\line(1,-1){40}}
\put(40,0){\line(1,1){40}}
\put(0,-10){\line(0,1){60}}
\put(0,10){\line(4,-1){80}}
\put(30,40){\line(1,-1){50}}
\put(0,10){\line(1,1){30}}
\put(80,40){\line(0,-1){50}}
\put(0,0){\line(1,0){40}}
\put(20,-10){$L$}

\put(-3,55){$\Lambda'$}

\put(0,40){\circle*{3}}
\put(-10,38){$a'$}
\put(30,40){\circle*{3}}
\put(28,43){$b'$}
\put(80,40){\circle*{3}}
\put(78,43){$c'$}
\put(0,10){\circle*{3}}
\put(-8,8){$c$}
\put(0,0){\circle*{3}}
\put(-8,-3){$\lambda$}
\put(40,0){\circle*{3}}
\put(38,-10){$b$}
\put(80,-10){\circle*{3}}
\put(78,-18){$a$}

\end{picture}

Since $a'\notin F$, the line $\Aline {a'}b$ does not belong to the family $\F_b$ and hence there exists a line $L_a\in\mathcal L_{b'}$ such that $L_a\parallel\Aline {a'}b$.
Since $c'\notin\bigcup\mathcal F_b$, there exists a line $L_c\in\mathcal L_{b'}$ such that $L_c\parallel \Aline b{c'}$. Since $\mathcal L_{c'}\cap\Lambda_\parallel'=\mathcal L_{c'}\cap\Lambda_\parallel=\varnothing$, the parallel lines $\Aline b{c'}$ and $L_c$ are concurrent with the line $\Lambda'$. Then there exists a unique point $c\in L_c\cap\Lambda'$. Assuming that $\Aline cb\parallel L'$, we conclude that $c\in \Lambda'\cap I\subseteq \Lambda'\cap F'=\varnothing$, which is a contradiction showing that the lines $\Aline cb$ and $L'=\Aline{a'}{b'}$ are concurrent. Since $\Aline cb\cap \Aline {a'}{b}=\{b\}$ and $L_a\parallel\Aline {a'}b$, we can apply Theorem~\ref{t:Proclus-lines} and conclude that there exists a unique point $a\in L_a\cap\Aline cb$. 

By the Proclus Axiom, the set
$$T\defeq(\Aline a{b'}\cap\Aline {a'}b)\cup(\Aline b{c'}\cap\Aline {b'}c)\cup(\Aline a{c'}\cap\Aline{a'}c)=\varnothing\cup\varnothing\cup(\Aline a{c'}\cap\Aline {a'}c)$$
has rank $\|T\|=\{0,2\}$. Since $|T|=|\Aline a{c'}\cap\Aline {a'}c|\le 1$, the set $T$ is empty, which means that the line $\Aline a{c'}$ is parallel to the lines $\Aline {a'}c=\Lambda'$ and $\Lambda$. Therefore, $\Aline a{c'}\in\mathcal L_{c'}\cap \Lambda_\parallel\ne\varnothing$, which is a desired contradiction, completing the proof of the lemma.
\end{proof}

\begin{lemma}\label{l:Vlad3} Let $X$ be an $\w$-long Pappian Proclus plane and $p,p'\in X$ be points such that the family $\F_p\defeq \{L\in\mathcal L_p:\mathcal L_{p'}\cap L_\parallel=\varnothing\}$ is finite. Then for every line $\Lambda\in \mathcal L\setminus \mathcal L_{p}$ and every point ${c}\in X\setminus\{p,p'\}$, the intersection  $\mathcal L_{c}\cap \Lambda_\parallel$ is not empty.
\end{lemma}

\begin{proof} To derive a contradiction, assume that $\mathcal L_{c}\cap\Lambda_\parallel=\varnothing$. In this case $c\notin \Lambda$. Since $p\notin \Lambda$, the set $\Lambda\cap\bigcup\F_p$ is finite, so we can choose a point $a\in \Lambda\setminus(\Aline {c}{p}\cup\Aline{c}{p'}\cup\bigcup\F_p)$. Choose any line $L_{c}\in\mathcal L_{c}\setminus(\mathcal L_p\cup\mathcal L_{p'}\cup\mathcal L_a)$. Since $\mathcal L_c\cap\Lambda_\parallel=\varnothing$,  the line $L_{c}$ is concurrent with the line $\Lambda$. 

Since $p\notin \Aline ac$, the set $\Aline {a}{c}\cap\F_p$ is finite, so we can find a point $b\in\Aline ac\setminus(\{a,c\}\cup\Aline p{p'}\cup\bigcup\F_p)$. By Lemma~\ref{l:Vlad2}, there exists a line $L_{c'}\in \mathcal L_b$ such that $L_{c'}\parallel  L_c$. Since $L_{c'}\parallel L_c\nparallel \Lambda$, there exists a point $c'\in L_{c'}\cap \Lambda$. It follows from $b\ne a$ and $L_{c'}\parallel L_c\nparallel \Aline ac$ that $c'\ne a$. Since $p\notin L_c$, the set $L_c\cap\bigcup\F_p$ is finite. Since $X$ is Proclus, the set $I\defeq\{u\in L_c:\Aline u{c'}\cap \Aline ac= \varnothing\}$ contains at most one point. Then there exists a point $$\textstyle b'\in L_c\setminus(\{c\}\cup\Lambda\cup I\cup\Aline ap\cup \Aline a{p'}\cup\Aline p{p'}\cup \bigcup\F_p).$$
By Lemma~\ref{l:Vlad2}, there exist a line $L_{a'}\in\mathcal L_{b}$ such that $L_{a'}\parallel \Aline a{b'}$. Since $L_{a'}\parallel \Aline a{b'}\nparallel \Aline {c'}{b'}$, there exists a unique point $a'\in L_{a'}\cap\Aline {c'}{b'}$.

\begin{picture}(100,110)(-140,-35)

\put(0,40){\line(1,0){80}}
\put(0,40){\line(1,-1){40}}
\put(80,40){\line(-1,-1){55}}
\put(17,-25){$L_c$}
\put(0,-10){\line(0,1){60}}
\put(0,10){\line(4,-1){80}}
\put(30,40){\line(1,-1){50}}
\put(0,10){\line(1,1){30}}
\put(80,40){\line(0,-1){50}}

\put(-3,55){$\Lambda$}

\put(0,40){\circle*{3}}
\put(-10,38){$a$}
\put(30,40){\circle*{3}}
\put(28,43){$b$}
\put(80,40){\circle*{3}}
\put(78,43){$c$}
\put(0,10){\circle*{3}}
\put(-8,8){$c'$}
\put(40,0){\circle*{3}}
\put(38,-11){$b'$}
\put(80,-10){\circle*{3}}
\put(78,-18){$a'$}
\put(55,31){\circle*{3}}
\put(58,28){$p'$}
\put(52,5){\circle*{3}}
\put(55,2){$p$}

\end{picture}

 It follows from $b'\notin\Aline ab$ and $a\ne b$ that the parallel lines $\Aline a{b'}$ and $L_{a'}$ are disjoint and hence $a'\ne b'$. Assuming that $a'=c'$, we conclude that $\Aline a{b'}\parallel \Aline {a'}{b}=\Aline{c'}b\parallel \Aline c{b'}$ and hence $\Aline a{b'}=\Aline c{b'}$ and $c=a\in\Lambda$, which is a contradiction showing that $a',b',c'$ are three distinct points of the line $\Aline {c'}{b'}$. The choice of the point $b'\notin I$ ensures that the lines $\Aline{c'}{b'}$ and $\Aline ac$ are concurrent. It is easy to see that $\{a',b',c'\}\cap\Aline ac=\varnothing=\{a,b,c\}\cap\Aline {a'}{c'}$. By the Pappus Axiom, the set
$$T\defeq(\Aline a{b'}\cap\Aline {a'}{b})\cup(\Aline b{c'}\cap\Aline{b'}c)\cup(\Aline a{c'}\cap\Aline {a'}c)=\varnothing\cup\varnothing\cup (\Aline a{c'}\cap\Aline {a'}{c})=\Lambda\cap\Aline {a'}c$$ has rank $\|T\|\in\{0,2\}$. Since $|\Lambda\cap\Aline {a'}c|\le 1$, we conclude that $\Aline c{a'}\cap\Lambda=\varnothing$ and hence $\Aline c{a'}\in\mathcal L_c\cap\Lambda_\parallel=\varnothing$, which is a final contradiction, completing the proof of the lemma.
\end{proof}

\begin{lemma}\label{l:Vlad4} Let $X$ be an $\w$-long Pappian Proclus plane and $b,b'\in X$ be points such that the family $\F_b\defeq \{L\in\mathcal L_b:\mathcal L_{b'}\cap L_\parallel=\varnothing\}$ is finite. Then the liner $X$ is Playfair.
\end{lemma}

\begin{proof} Given any line $\Lambda$ and point $x\in X\setminus\Lambda$, we should prove that the family $\mathcal L_x\cap\Lambda_\parallel $ contains a unique line. Since $X$ is Proclus, this family contains at most one line. Take any distinct points $p,p'\in X\setminus (\{x\}\cup\Lambda\cup\Aline b{b'})$. We claim that the family $\mathcal F_p\defeq\{L\in\mathcal L_p:\mathcal L_{p'}\cap L_\parallel=\varnothing\}$ is a subfamily of the family $\{\Aline pb,\Aline p{b'}\}$. Indeed, for every line $L\in\mathcal L_p\setminus\{\Aline pb,\Aline p{b'}\}$ we obtain $\{b,b'\}\cap L=\varnothing$. Applying Lemma~\ref{l:Vlad3}, we conclude that $\mathcal L_{p'}\cap L_\parallel\ne\varnothing$ and hence $L\notin \F_p$. Therefore, the set $\F_p$ is finite. Since $p,p'\notin \Lambda$ and $x\notin\{p,p'\}$, we can apply Lemma~\ref{l:Vlad3} and conclude that the family $\mathcal L_x\cap\Lambda_\parallel$ is not empty and hence it is a singleton, witnessing that the Proclus plane $X$ is Playfair.
\end{proof}

\begin{lemma}\label{l:Vlad5} Every $\w$-long Pappian Proclus plane $X$ is para-Playfair.
\end{lemma}

\begin{proof} To derive a contradiction, assume that $X$ is not para-Playfair. Then $X$ contains two disjoint lines $L,L'$ and a point $x\in X\setminus (L\cup L')$ such that $\mathcal L_x\cap L_\parallel=\varnothing$. By Proposition~\ref{p:cov-aff}, there exists a  $c\in X\setminus(\{x\}\cup L\cup L')$. We claim that the family $\F_x\defeq\{\Lambda\in\mathcal F_x:\mathcal L_c\cap\Lambda_\parallel=\varnothing\}$ is empty. 
In the opposite case, we can find a line $\Lambda\in\mathcal L_x$ such that $\mathcal L_c\cap\Lambda_\parallel=\varnothing$. This implies that $c\notin\Lambda$. 

Since $\mathcal L_x\cap L_\parallel=\varnothing=\mathcal L_x\cap L'_\parallel$,  there exist points $a\in\Lambda\cap L$ and $a'\in \Aline cx\cap L'$.  Then $\Lambda=\Aline ax$. By Theorem~\ref{t:Proclus-lines}, there exists a unique point $b\in \Aline ac\cap L'$. Assuming that $b=a'$, we conclude that $x\in \Aline c{a'}=\Aline cb=\Aline ca$ and $\Lambda=\Aline ax=\Aline ca$, which contradicts $c\in\Lambda$. This contradiction shows that $b\ne a'$. Assuming that $b\in \Aline ax$, we conclude that $c\in\Aline ab\cap\Aline {a'}x\subseteq \Aline ax\cap \Aline {a'}x=\{x\}$, which contradicts $c\notin\Lambda$. Therefore, $b\notin\{a'\}\cup \Aline ax$. Since the Proclus plane $X$ is not Playfair, we can apply Lemma~\ref{l:Vlad4} and conclude that the family $\F_{a'}=\{\Lambda'\in \mathcal L_{a'}:\mathcal L_b\cap\Lambda'_\parallel=\varnothing\}$ is infinite. Then we can choose a line $\Lambda'\in\F_{a'}\setminus\{\Aline {a'}x,\Aline {a'}a,L'\}$ such that $\Lambda'\nparallel \Lambda$. Since $\Lambda'\nparallel L'$, there exists a unique point $b'\in\Lambda'\cap L$. The inclusion $\Lambda'\in\F_{a'}$ implies that every line containing the point $b$ is not parallel to the line $\Aline{a'}{b'}$. In particular, the line $\Aline ab$ is concurrent with the line $\Lambda'=\Aline{a'}{b'}$. Since $\Aline ax=\Lambda \nparallel\Aline{a'}{b'}$, there exists a point $c'\in \Aline ax\cap\Aline{a'}{b'}$. It follows from $x\notin \Aline a{a'}$ that $c'\notin\{a',b'\}$. Also $c\in\Aline x{a'}$ and $b'\notin\Aline x{a'}$ imply $\{a',b',c'\}\cap \Aline ab=\varnothing=\{a,b,c\}\cap\Aline{a'}{b'}$.  

 Since $X$ is Pappian, the nonempty set $$T\defeq(\Aline a{b'}\cap\Aline {a'}b)\cup(\Aline b{c'}\cap\Aline {b'}c)\cup(\Aline a{c'}\cap\Aline{a'}c)=\varnothing\cup(\Aline b{c'}\cap\Aline {b'}c)\cup\{x\}$$ has rank $\|T\|=2$, which implies that the set $\Aline b{c'}\cap\Aline {b'}c$ is not empty and hence contains some point $y$. It is easy to see that $y\ne x$. The choice of the point $x$ ensures that the intersection $\Aline xy\cap L$ contains some point $z$. 

 The choice of the line $\Lambda'\in\F_{a'}$ ensures that the line $\Aline bz$ is not parallel to the line $\Aline{a'}{b'}=\Lambda'$ and hence there exists a point $\alpha'\in\Aline bz\cap\Aline{a'}{b'}$. 

We claim that $\alpha'\notin\{b',c'\}$. Assuming that $\alpha'=b'$, we conclude that $\alpha'=b'\in\Aline{a'}{b'}\cap L\cap\Aline {b}z=\{b'\}\cap\{z\}$ and hence $z=\alpha'=b'$ and $x\in\Aline yz\cap\Aline{a'}c\subseteq \Aline {b'}c\cap\Aline{a'}c=\{c\}$, which contradicts the choice of $c$. Assuming that $\alpha'=c'$, we conclude that $z\in \Aline b{\alpha'}=\Aline b{c'}$, $x\in\Aline yz\subseteq\Aline b{c'}$, $x\in\Aline a{c'}\cap\Aline b{c'}=\{c'\}$, and finally, $b'\in\Aline {a'}{c'}=\Aline {a'}x$, which contradicts the choice of the point $b'$. 

Therefore, $\alpha'\notin\{b',c'\}$ and  $\alpha',b',c'$ are three distinct points on the line $\Aline{a'}{c'}$. Assuming that $\alpha'\in\Aline ab$, we conclude that $z\in \Aline {\alpha'}b\cap L\subseteq\Aline ab\cap L=\{a\}$, $y\in \Aline xz\in\Aline xa=\Aline a{c'}$ and hence $y\in\Aline a{c'}\cap\Aline b{c'}\cap\Aline {b'}{c}=\{c'\}\cap\Aline{b'}c=\varnothing$, which is a contradition showing that $\alpha'\notin\Aline ab$. Also $\alpha'\ne a'$ because $z\in \Aline b{\alpha'}\cap L\ne \varnothing=L\cap L'=\Aline b{a'}\cap L'$. 

\begin{picture}(200,160)(-180,-20)

\put(-40,0){\line(1,3){40}}
\put(-40,0){\line(1,0){190}}
\put(90,-10){$L$}
\put(-30,30){\line(1,0){100}}
\put(75,27){$L'$}
\put(40,0){\line(-1,3){40}}
\put(-40,0){\line(1,1){60}}
\put(40,0){\line(-1,1){60}}
\put(-30,30){\line(5,3){50}}
\put(30,30){\line(-5,3){50}}
\put(-5,45){\line(1,0){85}}
\put(-30,30){\line(6,-1){180}}
\qbezier(80,45)(130,45)(150,0)

\put(5,45){\circle*{3}}
\put(3,37){$x$}
\put(-5,45){\circle*{3}}
\put(-15,44){$y$}
\put(-40,0){\circle*{3}}
\put(-45,-9){$a$}
\put(40,0){\circle*{3}}
\put(37,-10){$b'$}
\put(-30,30){\circle*{3}}
\put(-37,29){$b$}
\put(30,30){\circle*{3}}
\put(32,32){$a'$}
\put(-20,60){\circle*{3}}
\put(-28,59){$c$}
\put(20,60){\circle*{3}}
\put(22,58){$c'$}
\put(33.5,19.5){\circle*{3}}
\put(35,21){$\alpha'$}
\put(150,0){\circle*{3}}
\put(148,-8){$z$}

\end{picture}

Applying the Pappus Axiom to the line triples $abc$ and $\alpha'b'c'$, we conclude that  the set 
$$T'\defeq (\Aline a{b'}\cap\Aline {\alpha'}b)\cup(\Aline b{c'}\cap\Aline {b'}c)\cup(\Aline {\alpha'}c\cap\Aline a{c'})=\{z\}\cup \{y\}\cup(\Aline {\alpha'}c\cap\Aline a{c'})$$ has rank $\|T'\|=2$. Assuming that the intersection $\Aline {\alpha'}c\cap\Aline a{c'}$ contains some point $p$ and taking into account that $\|\{z,y,p\}\|=\|T'\|=2$, we conclude that $p\in \Aline zy=\Aline xy$. It follows from $a'\ne\alpha'$ that $\{p\}\cap\{x\}\subseteq\Aline {\alpha'}c\cap\Aline {a'}c\cap\{x\}=\{c\}\cap\{x\}=\varnothing$ and hence $p\ne x$. Then $p\in\Aline xy$ implies $y\in\Aline xp\subseteq\Aline a{c'}$ and hence $y\in\Aline a{c'}\cap\Aline b{c'}\cap\Aline {b'}c=\{c'\}\cap\Aline{b'}c=\varnothing$, which is a contradiction showing that $\Aline {\alpha'}c\cap\Aline a{c'}=\varnothing$ and hence $\Aline c{\alpha'}$ is a line that contains $c$ and is parallel to the line $\Aline a{c'}=\Aline ax=\Lambda$. Therefore, $\Aline c{\alpha}\in\mathcal L_c\cap\Lambda_\parallel\ne\varnothing$, which contradicts the choice of $\Lambda\in\F_c$. This contradiction shows that $\F_c=\varnothing$. Applying Lemma~\ref{l:Vlad4}, we conclude that the Proclus plane $X$ is Playfair and hence para-Playfair, which contradicts our assumption. This contradiction shows that $X$ is para-Playfair. 
\end{proof}

\begin{lemma}\label{l:Vlad6} For every $4$-long completion $Y$ of a Pappian liner $X$, the horizon $Y\setminus X$ is flat in $Y$.
\end{lemma}

\begin{proof}
To prove that the horizon $H := Y \setminus X$ of $X$ in $Y$ is flat, take any points $x,y\in H$ and a point $z\in\Aline xy$. If $z\in\{x,y\}$, then $z\in H$ and we are done. So, assume that $z\not\in\{x,y\}$, which implies $x\neq z\neq y$. Since $Y$ is a completion of $X$, there exists a point $a\in Y\setminus\overline H\subseteq X$. It follows from $x\in H$ and $a\not\in H$ that $\Aline ax\cap H = \{x\}$.  
Consider the plane $P\defeq\overline{\{a,x,y\}}$ in the projective liner $Y$. The $3$-rankedness of the regular projective liner $Y$ ensures that $P\cap\overline H=\Aline xy$. By Proposition~\ref{p:cov-aff}, there exist a point $o\in P\setminus(\Aline xy\cup\Aline ax)$. Choose any point $b'\in\Aline ax\setminus\{a,x\}$. Since the liner $Y$ is projective, the lines $\Aline oa$ and $\Aline y{b'}$ in the plane $P$ have a unique common point $c$. By the same reason, the lines $\Aline za$ and $\Aline o{b'}$ in the plane $P$ have a unique common point $c'$. Let $b$ be the unique common points of the lines $\Aline oa$ and $\Aline y{c'}$. Finally, the lines $\Aline bx$ and $\Aline o{b'}$ in the plane $P$ have a unique common point $a'$. 

\begin{picture}(150,160)(-100,-15)

\put(0,0){\line(1,0){150}}
\put(0,0){\line(1,1){135}}

\put(90,0){\color{blue}\line(-2,1){60}}
\put(90,0){\color{blue}\line(1,4){30}}
\put(120,0){\color{blue}\line(-1,1){60}}
\put(120,0){\color{blue}\line(-3,1){90}}
\put(60,0){\color{blue}\line(0,1){60}}
\put(60,0){\color{blue}\line(1,2){60}}

{\linethickness{1pt}
\put(60,15){\color{red}\line(4,1){36}}
}

\put(0,0){\circle*{3}}
\put(60,0){\circle*{3}}
\put(57,-10){$a$}
\put(90,0){\circle*{3}}
\put(87,-10){$b$}
\put(120,0){\circle*{3}}
\put(117,-10){$c$}

\put(30,30){\circle*{3}}
\put(36,30){$a'$}
\put(60,60){\circle*{3}}
\put(65,57){$b'$}
\put(120,120){\circle*{3}}
\put(123,115){$c'$}

\put(60,15){\color{red}\circle*{3}}
\put(52,10){\color{red}$x$}
\put(68.6,17.2){\color{red}\circle*{3}}
\put(64,22){\color{red}$z$}
\put(96,24){\color{red}\circle*{3}}
\put(99,23){\color{red}$y$}

\put(-3,-10){$o$}
\end{picture}  

The choice of the points $a,b,c,a',b',c'$ guarantees that we have triples $abc$ and $a'b'c'$ in $X$ such that $\Aline a{b'}\cap\Aline {a'}b=\{x\}$, $\Aline b{c'}\cap\Aline {b'}c=\{y\}$ and $z\in\Aline a{c'}$. Since the liner $X$ is Pappian, the set $$T:=(\Aline a{c'}\cap\Aline {a'}c\cap X)\cup(\Aline b{c'}\cap\Aline {b'}c\cap X)\cup(\Aline a{b'}\cap\Aline{a'}b\cap X)$$ has rank $\|T\|\in\{0,2\}$ in the liner $X$. Since $x,y\not\in X$, the set $T$ is equal to the set $\Aline a{c'}\cap\Aline {a'}c\cap X$ which contain at most one point $\{z\}$. Then $\|T\|=0$ and hence $\Aline a{c'}\cap\Aline {a'}c\cap X=\emptyset$. Since the liner $Y$ is projective, the set $\Aline a{c'}\cap\Aline {a'}c$ is not empty. Then $\emptyset\not=\Aline a{c'}\cap\Aline {a'}c\subseteq\Aline a{c'}\setminus X\subseteq\{z\}$ and hence $z\in Y\setminus X=H$, which means that the horizon $H$ of $X$ in $Y$ is flat.
\end{proof}

\begin{theorem}\label{t:Pappian=>compreg} Every Pappian proaffine regular liner $X$ is completely regular.
\end{theorem}

\begin{proof} By Theorem~\ref{t:Proclus<=>}, the proaffine regular liner $X$ is Proclus.
If $X$ is projective, then $X$ is strongly regular and hence completely regular.
So, assume that $X$ is not projective.  If $\|X\|\ne 3$, then $X$ is completely regular, by Theorem~\ref{t:proaffine3=>compregular}. So, assume that $\|X\|=3$. If $X$ is not $3$-long, then by Theorem~\ref{t:Proclus-not-3long}, $X$ has a projective completion with a flat horizon.  Theorem~\ref{t:spread=projective2}, the liner $X$ is completely regular. So, assume that the liner $X$ is $3$-long. If $X$ is finite, then $X$ has a projective completion $Y$, by Theorem~\ref{t:procompletion-finite}. By Corollary~\ref{c:Avogadro-projective}, the $3$-long projective liner $Y$ is $2$-balanced. Since the liner $X$ is $3$-long and $X\ne Y$, the $2$-balanced projective liner $Y$ is $4$-long. By  Lemma~\ref{l:Vlad6}, the horizon $Y\setminus X$ is flat, and by Theorem~\ref{t:spread=projective2}, the liner $X$ is completely regular.
\smallskip

So, assume that $X$ is $3$-long and infinite. Repeating the argument of Claim~\ref{cl:w-long}, we can prove that the infinite regular plane $X$ is $\w$-long. By Lemma~\ref{l:Vlad5}, the Pappian Proclus plane $X$ is para-Playfair. It remains to prove that for every concurrent Bolyai lines $A,B\subseteq X$, every line $L$ is the plane $\overline{A\cup B}=X$ is Bolyai. 

By Theorem~\ref{t:spreading<=>}, it suffices to show that every point $x\in X\setminus L$ is contained in some line $L_x\in L_{\parallel}$. If $L$ parallel to $A$ or $B$, then the existence of the line $L_x$ follows from the para-Playfair property of $X$. So, assume that $A\nparallel L\nparallel B$. Since the liner $X$ is Proclus, the set $I\defeq\{z\in L:\Aline xz\parallel B\}\cup\{z\in L:\Aline xz\parallel A\}$ contains at most two points. So, we can choose a point $z\in L\setminus I$. 
Since the lines $A$, $B$ are spreading, there exist lines $A_z\in A_{\parallel}$ and $B_{x}\in B_{\parallel}$ such that $z\in A_z$ and $x\in B_x$. Since the lines $A,B$ are concurrent, so are the lines  $A_z\in A_{\parallel}$ and $B_x\in B_{\parallel}$ in the plane $\overline{A\cup B}=X$. Let $y'$ be the unique common point of the concurrent lines $A_z$ and $B_x$. The choice of the point $z\in L\setminus I$ ensures that $y'\notin\{x,z\}$ and hence $y'\notin \Aline zx\cup L$. 

Since $X$ is Proclus, the set $J\defeq\{x'\in L:\Aline{x'}{y'}\cap\Aline xz=\varnothing\}$ contains at most one point. Then there exists a point $x'\in L\setminus(J\cup \{z\}\cup B_x)$. The choice of the point $x'\in L\setminus J$ ensures that the lines $\Aline xz$ and $\Aline{x'}{y'}$ have a unique common point $o$. Since $x',y'\notin \Aline xz$, the point $o$ is not equal to $x'$ or $y'$. The choice of the point $x'\notin \{z\}\cup B_x$ implies that $o\notin\{x,z\}$.

\begin{picture}(150,180)(-180,-35)
\put(0,60){\line(2,-1){150}}
\put(0,120){\line(1,0){150}}
\put(0,120){\line(-1,0){120}}
\put(0,60){\line(-2,1){120}}

{\linethickness{1pt}
\put(0,35){\color{red}\line(0,1){95}}
\put(2,40){$L$}
\put(150,-20){\color{red}\line(0,1){150}}
\put(152,30){$L_x$}

\put(-10,135){\color{blue}\line(2,-3){100}}
\put(82,0){$A_z$}
\put(50,135){\color{blue}\line(2,-3){105}}
\put(125,25){$A_y$}

\put(75,135){\color{cyan}\line(-1,-1){100}}
\put(-30,45){$B_{x'}$}
\put(165,135){\color{cyan}\line(-1,-1){130}}
\put(30,15){$B_x$}
}

\put(-120,120){\circle*{3}}
\put(-128,118){$o$}

\put(0,60){\circle*{3}}
\put(4,59){$x'$}
\put(60,30){\circle*{3}}
\put(58,36){$y'$}
\put(150,-15){\circle*{3}}
\put(152,-15){$z'$}

\put(0,120){\circle*{3}}
\put(2,122){$z$}
\put(60,120){\circle*{3}}
\put(59,126){$y$}
\put(150,120){\circle*{3}}
\put(142,122){$x$}
\end{picture}

Since the line $B$ is spreading, there exists a line $B_{x'}\in B_{\parallel}$ such that $x'\in B_{x'}$. The choice of the point $z\in L\setminus I$ ensures that the line $\Aline xz$ is not parallel to the parallel lines $B$ and $B_{x'}$. Then the lines $\Aline xz$ and $B_{x'}$ have a unique common point $\{y\}=\Aline xz\cap B_{x'}$. Since $z\ne x'\notin B_x$ and $L\nparallel B_{x'}$, the point $y$ is distinct from the points $z$ and $x$.  Assuming that $y=o$, we conclude that $B_{x'}=\Aline {x'}{y}=\Aline {x'}o=\Aline{x'}{y'}$ and hence $B_{x'}=B_x$, which contradicts the choice of the point $x'\notin B_x$. This contradiction shows that $y\ne o$ and hence $o\notin\{x,y,z\}$.

Since the line $A$ is spreading, there exists a line $A_y\in A_{\parallel}$ such that $y\in A_y$. By Proposition~\ref{p:Proclus-Postulate}, there exists a unique common point $z'=\Aline {x'}{y'}\cap A_y$. It follows from $y\notin A_z$ that $A_y\cap A_z=\varnothing$ and hence $z'\ne y'$. It follows from $\{y\}=B_{x'}\cap A_y$ that $z'\ne x'$. Since the lines $A_y$ and $\Aline xy$ are not parallel and $o\ne y$, the point $z'$ does not belong to the line $\Aline xz$ and hence $z'\ne o$. Therefore, $x',y',z'$ are three distinct collinear points such that $o\notin\{x',y',z'\}$.

Applying the Pappus Axiom to the triples $zyx$ and $x'y'z'$, we conclude that the set $$T:=(\Aline {x'}z\cap\Aline x{z'})\cup(\Aline {y'}z\cap\Aline y{z'})\cup(\Aline {x'}y\cap\Aline x{y'})=(L\cap\Aline x{z'})\cup(A_z\cap A_y)\cup(B_{x'}\cap B_x)=(L\cap\Aline x{z'})\cup\varnothing\cup\varnothing$$ has rank $\|T\|=0$, which implies that the intersection $L\cap\Aline x{z'}$ is empty. This means that $L\parallel\Aline x{z'}$ and hence $L_x:=\Aline x{z'}\in L_{\parallel}$.

By Theorem~\ref{t:spread=projective1}, $X$ is completely regular.
\end{proof}




Theorem~\ref{t:Pappian=>compreg} will help us to prove the following important extension of Proposition~\ref{p:Pappian-xy} to Pappian Proclus liners.

\begin{theorem}\label{t:Pappus=>xy||ac'} Let $L,L'$ be two coplanar lines in a Pappian Proclus liner $X$ and $a,b,c\in L\setminus L'$ and $a',b',c'\in L'\setminus L$ be distinct points such that $\Aline a{c'}\cap\Aline{a'}c=\varnothing$ and $\Aline a{b'}\cap\Aline{a'}b=\{x\}$ for some point $x\in X$. Then there exists a unique point $y\in X$ such that $\{y\}=\Aline b{c'}\cap\Aline {b'}c$ and $\Aline xy\cap \Aline a{c'}=\varnothing=\Aline xy\cap\Aline{a'}{c'}$.
\end{theorem}

\begin{proof} 
Consider the plane $\Pi=\overline{L\cup L'}$. This plane contains two disjoint lines $\Aline a{c'}$ and $\Aline{a'}c$ and hence $\Pi$ is not projective. Observe that $\Pi$ has cardinality $|\Pi|\ge|\{a,b,c,a',b',c',x\}|\ge 7$. Applying Theorem~\ref{t:Proclus-not-3long}, we conclude that the Proclus plane $\Pi$ is $3$-long. By Theorems~\ref{t:Pappian=>compreg} and \ref{t:Proclus<=>}, the Pappian Proclus plane $\Pi$ is completely regular and hence the spread completion $\overline\Pi$ of $\Pi$ is a $3$-long projective liner. By Corollary~\ref{c:procompletion-rank}, $\|\overline \Pi\|=\|\Pi\|=3$. By Corollary~\ref{c:Avogadro-projective}, the $3$-long projective liner $\overline \Pi$ is $2$-balanced. Since $\Pi$ is $3$-long and not projective, the projective liner $\overline\Pi$ is $4$-long. Since $\|\partial\Pi\|<\|\overline\Pi\|=3$, the horizon $\partial \Pi=\overline\Pi\setminus\Pi$ of $\Pi$ is either a singleton or a line.

If $\partial \Pi$ is a line, then $\Pi$ is a Playfair plane, and  Proposition~\ref{p:Pappian-xy} applies.
So, assume that $\partial\Pi$ is a singleton. In this case $\Pi$ is a punctured projective plane and $(\Aline a{c'})_\parallel$ is a unique spread of parallel lines in $\Pi$. This implies that any two disjoint lines in $\Pi$ are parallel to the lines $\Aline a{c'}$ and $\Aline {a'}c$. Since the lines $L,L'$ are not parallel to the line $\Aline a{c'}$, they are concurrent.

By the Pappus Axiom, the set 
$$T\defeq(\Aline a{b'}\cap\Aline{a'}b)\cup(\Aline b{c'}\cap\Aline{b'}c)\cup(\Aline a{c'}\cap\Aline{a'}c)=\{x\}\cup(\Aline b{c'}\cap\Aline{b'}c)\cup\varnothing$$ has rank $\|T\|=2$, which implies that the intersection $\Aline b{c'}\cap\Aline{b'}c$ contains some point $y$. Assuming that $\Aline b{c'}\cap\Aline{b'}c\ne\{y\}$, we conclude that $\Aline b{c'}=\Aline{b'}c$ and hence $\{b\}=L\cap\Aline b{c'}=L\cap\Aline {b'}c=\{c\}$, which contradicts the choice of the distinct points $b,c$. This contradiction shows that $\Aline b{c'}\cap\Aline{b'}c=\{y\}$. It remains to prove that $\Aline xy\cap\Aline ac=\varnothing=\Aline xy\cap\Aline{a'}{c'}$.

Assuming that $\Aline xy\cap\Aline a{c'}$ contains some point $z$, consider the line $\Aline cz$. Since the lines $\Aline cz$ and $L'$ are not parallel to the line $\Aline a{c'}$, there exists a point $\alpha'\in L'\cap \Aline cz$. It is easy to see that $\alpha'\notin\{b',c'\}\cup L$. 

\begin{picture}(200,160)(-180,-20)

\put(-40,0){\line(1,3){40}}
\put(-40,0){\line(1,0){190}}
\put(-20,90){$L$}
\put(-30,30){\line(1,0){60}}
\put(14,90){$L'$}
\put(40,0){\line(-1,3){40}}
\put(-40,0){\line(1,1){60}}
\put(40,0){\line(-1,1){60}}
\put(-30,30){\line(5,3){50}}
\put(30,30){\line(-5,3){50}}
\put(-5,45){\line(1,0){85}}
\put(-30,30){\line(6,-1){180}}
\qbezier(80,45)(130,45)(150,0)

\put(5,45){\circle*{3}}
\put(3,37){$x$}
\put(-5,45){\circle*{3}}
\put(-15,44){$y$}
\put(-40,0){\circle*{3}}
\put(-45,-9){$a$}
\put(40,0){\circle*{3}}
\put(37,-10){$c'$}
\put(-30,30){\circle*{3}}
\put(-37,29){$c$}
\put(30,30){\circle*{3}}
\put(32,32){$a'$}
\put(-20,60){\circle*{3}}
\put(-28,59){$b$}
\put(20,60){\circle*{3}}
\put(23,58){$b'$}
\put(33.5,19.5){\circle*{3}}
\put(35,21){$\alpha'$}
\put(150,0){\circle*{3}}
\put(148,-8){$z$}

\end{picture}

Applying the Proclus Axiom to the triples $abc$ and $\alpha'b'c'$, we conclude that the nonempty set 
$$T'\defeq(\Aline a{b'}\cap\Aline{\alpha'}b)\cup(\Aline b{c'}\cap\Aline{b'}c)\cup(\Aline a{c'}\cap\Aline{\alpha'}c)=(\Aline a{b'}\cap\Aline{\alpha'}b)\cup\{y\}\cup\{z\}$$ has rank $\|T'\|=2$. Since the lines $\Aline a{b'}$ and $\Aline{\alpha'}b$ are not parallel to the line $\Aline a{c'}$, there exists a point $x'\in\Aline a{b'}\cap\Aline{\alpha'}b$. It follows from $|\{x',y,z\}\|=\|T'\|=2$ that $x'\in\Aline yz=\Aline xy$ and $x'\in \Aline xy\cap \Aline a{b'}=\{x\}$ and $\alpha'\in\Aline b{x'}\cap L'=\Aline bx\cap L'=\Aline b{a'}\cap L=\{a'\}$. Then $z\in \Aline c{\alpha'}\cap \Aline a{c'}=\Aline c{a'}\cap\Aline a{c'}=\varnothing$, which is a contradiction that shows that $\Aline a{c'}\cap \Aline xy=\varnothing$. By the Proclus Axiom, $\Aline a{c'}\cap \Aline xy=\varnothing$ implies $\Aline xy\cap\Aline{a'}c=\varnothing$.
\end{proof}

\begin{corollary}\label{c:Pappus=>A-Pappus} Let $X$ be a  Pappian proaffine regular liner and let $H$ be a hyperplane in the spread completion $\overline X$ of $X$ such that $\overline X\setminus X\subseteq H$. Then the affine subliner $A\defeq \overline X\setminus H$ of $\overline X$ is Pappian.
\end{corollary}

\begin{proof} By Theorem~\ref{t:Pappian=>compreg}, the  Pappian proaffine regular liner $X$ is completely regular. By Theorem~\ref{t:spread=projective2}, the spread completion $\overline X$ of the completely regular liner is a projective liner. By Proposition~\ref{p:projective-minus-hyperplane}, the subliner $A=\overline X\setminus H$ of the projective liner $\overline X$ is affine and regular. 

Assuming that the affine liner $A$ is not Pappian, we can apply Theorem~\ref{t:Papp<=>APA} and conclude that the liner $A$ does not satisfy the Affine Pappus Axiom. Then there exist concurrent lines $L,L'$ in $A$ and distinct points $a,b,c\in L\setminus L'$ and $a',b',c'\in L'\setminus L$ such that $\Aline a{b'}\cap\Aline {a'}b\cap A=\varnothing=\Aline{b'}c\cap\Aline b{c'}\cap A$ but $\Aline a{c'}\cap\Aline{a'}c\cap A\ne \varnothing$. Consider the plane $P\defeq\overline{L\cup L'}$ in the projective liner $\overline X$. Since $\overline X$ is projective, any two lines in the plane $P$ have a common point. Therefore,  there exist unique points $x\in \Aline {a'}b\cap \Aline a{b'}\subseteq \overline X\setminus A=H$, $y\in \Aline {b'}c\cap \Aline b{c'}\subseteq \overline X\setminus A=H$,  and $z\in \Aline {a'}c\cap \Aline a{c'}\subseteq A$. It is easy to see that the points $x,y,z$ are distinct. Since the liner $X$ is Pappian, the nonempty set $\{x,y,z\}\cap X$ has rank $\|\{x,y,z\}\cap X\|=2$ in $X$. If $x,y\in X$, then the set $\{x,y,z\}=\{x,y,z\}\cap X$ has rank $2$ in $X$ and hence $z\in \Aline xy\subseteq H$, which contradicts the choice of the points $a,b,c,a',b',c'$. Therefore, $x\notin X$ or $y\notin X$. If $x\notin X$, then Theorem~\ref{t:Pappus=>xy||ac'} ensures that $y\in X$ and $(\Aline yz\cap X)\cap(\Aline a{b'}\cap X)=\varnothing$. In this case $(\Aline yz\cap X)_\parallel=(\Aline a{b'}\cap X)_\parallel=x$ and hence $x\in \Aline yz$ in $\overline X$. It follows from $\{x,y\}\subseteq (\overline X\setminus X)\cup(X\setminus A)\subseteq H$ that $z\in \Aline xy\subseteq H$, which contradicts $z\in A$. By analogy we can derive a contradiction assuming that $y\notin X$. Those contradictions show that the affine liner $A$ is Pappian.
\end{proof}

\section{Pappian proaffine regular liners are Desarguesian}

The following theorem generalizes the Hessenberg's Theorem~\ref{t:Hessenberg-affine} to arbitrary proaffine regular liners.

\begin{theorem}\label{t:Hessenberg-proaffine} Every Pappian proaffine regular liner is Desarguesian.
\end{theorem}

\begin{proof} Assume that $X$ is a Pappian proaffine regular liner.
By Theorem~\ref{t:Proclus<=>}, the proaffine regular liner $X$ is Proclus. By Theorem~\ref{t:Desargues<=>planeD}, to show that $X$ is Desarguesian, it suffices to check that every $3$-long plane in $X$ is Desarguesian. So, we lose no generality assuming that $X$ is a Pappian $3$-long Proclus plane. By Theorem~\ref{t:Pappian=>compreg}, the Pappian $3$-long Proclus plane $X$ is completely regular. By Theorem~\ref{t:spread=projective2}, the spread completion $\overline X$ of $X$ is a $3$-long projective plane. By Corollary~\ref{c:Avogadro-projective}, the $3$-long projective plane $\overline X$ is $2$-balanced. If $|\overline X|_2=3$, then the projective plane $\overline X$ is Desarguesian, by Proposition~\ref{p:Steiner+projective=>Desargues}. Since the subliner $X$ of $\overline X$ is $3$-long, $|\overline X|_2=3$ implies that $X=\overline X$ and hence the liner $X=\overline X$ is Desarguesian. So, assume that $|\overline X|_2\ge 4$. Let $H$ be any line in $\overline X$ that contains the flat $\partial X\defeq\overline X\setminus X$. By Theorem~\ref{t:affine<=>hyperplane}, $A\defeq \overline X\setminus H$ is a $3$-long affine regular liner in $\overline X$ and $\overline X$ is a projective completion of $A$. By Corollary~\ref{c:procompletion-rank}, $\|A\|=\|\overline X\|=3$, so $A$ is a Playfair plane, by Theorem~\ref{t:Playfair<=>}. By Corollary~\ref{c:Pappus=>A-Pappus}, the affine liner $A$ is Pappian. By Theorems~\ref{t:Papp<=>APA} and \ref{t:Hessenberg-affine}, the Pappian Playfair plane $A$ is Desarguesian. By Theorem~\ref{t:projDes<=>}, the projective completion $\overline X$ of the Desarguesian Playfair plane $A$ is Desarguesian. Since the horizon $\partial X$ of $X$ in $\overline X$ is flat in $\overline X$, the subliner $X=\overline X\setminus\partial  X$ of the Desarguesian projective liner $\overline X$ is Desarguesian, by Theorem~\ref{t:pD-minus-flat}.
\end{proof}

\section{Finite Desarguesian proaffine regular liners are Pappian}

\begin{theorem}\label{t:finite-Papp<=>Des} A finite proaffine regular liner is Pappian if and only if it is Desarguesian.
\end{theorem}

\begin{proof} The ``only if'' part follows from Theorem~\ref{t:Hessenberg-proaffine}. 
To prove the ``if'' part, assume that a finite proaffine regular liner $X$ is Desarguesian.
By Theorem~\ref{t:Proclus<=>}, the proaffine regular liner $X$ is Proclus. 
To prove that $X$ is Pappian, choose any concurrent lines $L,L'$ and distinct points $a,b,c\in L\setminus L'$ and $a',b',c'\in L'\setminus L$. Let $o$ be the unique point of the intersection $L\cap L'$. Consider the plane $P\defeq\overline{L\cup L'}$ in the liner $X$ and observe that $|P|\ge|\{o,a,b,c,a',b',c'\}|=7$. 

\begin{claim} The plane $P$ is $3$-long.
\end{claim}

\begin{proof} If the Proclus plane $P$ is not projective, then $P$ is $3$-long, by Theorem~\ref{t:Proclus-not-3long}.

 If $P$ is projective, then we can consider the maximal $3$-long flat $M$ in $P$ that contains the point $o$. By Lemma~\ref{l:ox=2}, $\Aline ox=\{o,x\}$ for every $x\in P\setminus M$, which implies that $\{o,a,b,c,a',b',c'\}\subseteq M$. Then $3\le\|M\|\le\|P\|=3$ and hence $M=P$ by the $3$-rankedness of the regular liner $X$. Therefore, the projective plane $P=M$ is $3$-long. 
\end{proof}

By Corollary~\ref{c:Desarg=>compreg}, the Desarguesian $3$-long proaffine regular liner $P$ is completely regular. 
By Theorem~\ref{t:Desargues-completion}, the spread completion $\overline P$ of $P$ is a Desarguesian $3$-long projective liner. By Corollary~\ref{c:procompletion-rank}, $\|\overline P\|=\|P\|=3$. Let $H$ be any line in the projective plane $\overline P$. By Theorem~\ref{t:affine<=>hyperplane}, the subliner $A\defeq\overline P\setminus H$ is affine  and regular. Since the liner $X$ is finite, so are the liners $P$ and $\overline P$. By Theorem~\ref{t:pD-minus-flat}, the subliner $A=\overline P\setminus H$ of the Desarguesian projective liner $\overline P$ is  Desarguesian. By Theorem~\ref{t:fin-affP<=>D}, the finite Desarguesian affine liner $A$ is Pappian. By Theorem~\ref{t:Pappian-completion}, the spread completion $\overline P$ of the Pappian affine regular liner $A$ is Pappian. Since the set $\overline P\cap(\overline X\setminus X)=\overline P\setminus X$ is flat in $\overline P$, the subliner $P=\overline P\setminus(\overline P\setminus X)$ of the Pappian projective liner $\overline P$ is Pappian, by Proposition~\ref{p:Pappian-minus-flat}. Then the set $$T\defeq(\Aline a{b'}\cap \Aline {a'}b)\cup(\Aline b{c'}\cap\Aline {b'}c)\cup(\Aline a{c'}\cap\Aline {a'}c)$$ has rank $\|T\|\in\{0,2\}$ in $P$ and also in $X$, witnessing that the finite Desarguesian liner $X$ is Pappian. 
\end{proof}

\begin{exercise} Find a Desarguesian projective liner, which is not Pappian.
\smallskip

{\em Hint:} Consider the projective plane $\mathbb{ PH}^3$ over the corps of quaternions $\mathbb H$.
\end{exercise}

\section{Projective completions of Pappian liners}\label{s:Pappian-completion}

\begin{theorem}\label{t:Pappian-completion} A completely regular liner $X$ is Pappian if and only if its spread completion $\overline X$ is a Pappian projective liner.\end{theorem}

\begin{proof} To prove the ``if'' part, assume that  the spread completion $\overline X$ of $X$ is a Pappian projective liner. Since the horizon $\partial X\defeq\overline X\setminus X$ of $X$ in $\overline X$ is flat in $\overline X$, the subliner $X$ of the Pappian projective liner $\overline X$ is Pappian, by Proposition~\ref{p:Pappian-minus-flat}. 

To prove the ``only if'' part, assume that the completely regular liner $X$ is Pappian. If $X$ is projective, then the spread completion $\overline X=X$ is Pappian and we are done. So, assume that the liner $X$ is not projective. 

\begin{claim} The spread completion $\overline X$ of $X$ is a $3$-long Desarguesian projective liner.
\end{claim}

\begin{proof} By Proposition~\ref{p:spread-3long}, the spread completion $\overline X$ of the non-projective completely regular liner $X$ is projective and $3$-long. 
By Corollary~\ref{c:Avogadro-projective}, the $3$-long projective liner $\overline X$ is $2$-balanced. If $|\overline X|_2=3$, then the Steiner projective liner $\overline X$ is Desarguesian, by Proposition~\ref{p:Steiner+projective=>Desargues}.  So, assume that $|\overline X|_2\ge 4$. Then the projective liner $\overline X$ is $4$-long and the liner $X=\overline X\setminus\partial X$ is $3$-long. By Theorem~\ref{t:spread=projective1}, the completely regular liner $X$ is para-Playfair and hence $X$ is Proclus and proaffine, see Theorem~\ref{t:Proclus<=>}. By Theorem~\ref{t:Hessenberg-proaffine}, the Pappian proaffine regular liner $X$ is Desarguesian, and by Theorem~\ref{t:spread=projective2}, the spread completion $\overline X$ of $X$ is a  Desarguesian projective liner.
\end{proof}

 Assuming that the projective liner $\overline X$ is not Pappian, we can find two concurrent lines $L,L'\subseteq \overline X$ and distinct points $a,b,c\in L\setminus L'$ and $a',b',c'\in L'\setminus L$ such that the set $$T\defeq(\Aline a{b'}\cap\Aline{a'}b)\cup(\Aline b{c'}\cap\Aline {b'}c)\cup(\Aline a{c'}\cap\Aline {a'}c)$$has rank $\|T\|\notin\{0,2\}$. Let $o$ be a unique point of the intersection $L\cap L'$. Since the $3$-long projective liner $\overline X$ is $2$-balanced, $|\overline X|_2\ge|L|\ge |\{0,a,b,c\}|=4$ and hence the projective liner $\overline X$ is $4$-long and the liner $X=\overline X\setminus\partial X$ is $3$-long.

Consider the plane $P\defeq\overline{L\cup L'}$ in the projective liner $\overline X$. 
Since $\overline X$ is projective, any lines in the plane $P$ have a common point. In particular, there exist points $x\in \Aline a{b'}\cap\Aline {a'}b$, $y\in \Aline b{c'}\cap\Aline {b'}c$ and $z\in \Aline a{c'}\cap\Aline {a'}c$. It is easy to see that the points $x,y,z$ are distinct, which implies that $\|T\|=\|\{x,y,z\}\|\ge 2$. Since $\|T\|\notin\{0,2\}$, the set $T=\{x,y,z\}$ has rank $\|T\|=3$. In this case $o\notin \varnothing=\Aline xy\cap\Aline yz\cap\Aline xz$ and we lose no generality assuming that $o\notin\Aline xy$. 

If $P\cap X\ne \varnothing$, then $P\cap(\overline X\setminus X)=P\setminus X$ is a proper flat in the projective plane $P$. Choose any line $H\subseteq P$ that contains the flat $P\setminus X$. By Theorem~\ref{t:affine<=>hyperplane}, the subliner $A\defeq P\setminus H$ is affine and regular. Since the projective liner $\overline X$ is $4$-long, the affine subliner $A$ of $X$ is $3$-long. It follows that $P$ is a projective completion of the affine regular liner $A$. By Corollary~\ref{c:procompletion-rank}, $\|A\|=\|P\|=3$. Since the liner $X$ is Pappian, so is the flat $P\cap X$ in $X$. Observe that $A=P\setminus H\subseteq P\setminus(P\setminus X)=P\cap X$. By Corollary~\ref{c:Pappus=>A-Pappus}, the affine subliner $A$ of the Pappian proaffine regular liner $P\cap X$ is Pappian. 

Since $P$ is a Desarguesian projective space, we can apply Corollary~\ref{c:hyperplanes-automorphism}, and find an automorphism $\Phi:P\to P$ such that $\Phi[\Aline xy]=H$. For a point $p\in P$ we denote by $p_\Phi$ the image $\Phi(p)$ of the point $p$ under the automorphism $\Phi$. Taking into account that $\{o,a,b,c,a',b',c'\}\cap\Aline xy=\varnothing$, we conclude that the points $o_\Phi,a_\Phi,b_\Phi,c_\Phi,a'_\Phi,b'_\Phi,c'_\Phi$ belong to $P\setminus\Phi[\Aline xy]=P\setminus H=A$. Then $\Phi[L]\cap A$ and $\Phi[L']\cap A$ are two concurrent lines in $A$, and $a_\Phi,b_\Phi,c_{\Phi}\in \Phi[L]\cap A$ and $a'_\Phi,b'_\Phi,c'_{\Phi}\in \Phi[L']\cap A$ are distinct points in the liner $A$. Since $A$ is Pappian, the set 
$$
\begin{aligned}
T_A&\defeq A\cap \big((\Aline {a_\Phi}{b'_\Phi}\cap\Aline {a'_\Phi}{b_\Phi})\cup(\Aline {b_\Phi}{c'_\Phi}\cap\Aline {b'_\Phi}{c_\Phi})\cup(\Aline {a_\Phi}{c'_\Phi}\cap\Aline {a'_\Phi}{c_\Phi})\big)\\
&=A\cap\Phi[\{x,y,z\}]\subseteq A\cap (H\cup \{\Phi(z)\})=A\cap\{\Phi(z)\}
\end{aligned}
$$ has rank $\|T_A\|\in\{0,2\}$. Since $|T_A|=|A\cap\{\Phi(z)\}|\le|\{\Phi(z)\}|=1$, the set $T_A=A\cap\{\Phi(z)\}$ has rank $\|T_A\|=0$, which implies that $\Phi(z)\in P\setminus A=H=\Phi[\Aline xy]$ and $z\in\Aline xy$. Then $\|T\|=\|\{x,y,z\}\|=2$, which contradicts the choice of the points $a,b,c,a',b',c'$.

This contradiction shows that $P\subseteq \overline X\setminus X$. Choose any point $v\in X$ and consider the flat $V\defeq\overline{P\cup\{v\}}$ in the projective liner $\overline X$. Corollary~\ref{c:Max=dim} ensures that $\|V\|=4$ and hence $P$ is a hyperplane in $V$. Since $v\in V\cap X$, the intersection $V\cap(\overline X\setminus X)=V\setminus X$ is a proper flat in $V$, containing the plane $P$. Since $3=\|P\|\le\|V\setminus X\|<\|V\|=4$, the rankedness of the regular liner $X$ ensures that $P=V\setminus X$. Then $V\setminus P=V\cap X$. By Theorem~\ref{t:affine<=>hyperplane}, the subliner $A\defeq V\setminus P=V\cap X$  of the projective liner $V$ is affine and regular. Since the liner $V$ is $4$-long, the affine subliner $A$ of $V$ is $3$-long. It follows that $V$ is a projective completion of the affine regular liner $A$. By Corollary~\ref{c:procompletion-rank}, $\|A\|=\|V\|=4$. Since the liner $X$ is Pappian, so is the flat $A=V\cap X$ in $X$.

By Proposition~\ref{p:add-point-to-independent}, the set $\{x,y,o,v\}$ is independent in the projective liner $V$. Since $\|V\|=4$, this set is maximal independent in $V$. By Theorem~\ref{t:Des-autoextend}, there exists an automorphism $\Psi:V\to V$ of the Desarguesian projective space $V$ such that $\Psi oxyv=vxyo$. Consider the plane $\Pi\defeq\overline{\{v,x,y\}}\subseteq V$ and observe that $\Psi[P]=\Psi[\overline{\{o,x,y\}}]=\overline{\{v,x,y\}}=\Pi$. Since $v\notin P$, the intersection $\Pi\cap P$ is a proper flat in $\Pi$. Since $\Aline xy\subseteq \Pi\cap P$, the flat $\Pi\cap P$ coincides with the line $\Aline xy$, by the rankedness of the regular liner $X$.

 For a point $p\in V$ we denote by $p_\Psi$ the image $\Psi(p)$ of the point $p$ under the automorphism $\Psi$. Taking into account that $\{o,a,b,c,a',b',c'\}\subseteq P\setminus\Aline xy$,  we conclude that the points $o_\Psi,a_\Psi,b_\Psi,c_\Psi,a'_\Psi,b'_\Psi,c'_\Psi$ belong to $\Psi[P\setminus \Aline xy]=\Psi[P]\setminus\Psi[\Aline xy]=\Pi\setminus \Aline xy\subseteq V\setminus P=A$. Then $\Psi[L]\cap A$ and $\Psi[L']\cap A$ are two concurrent lines in $A$, and $a_\Psi,b_\Psi,c_{\Psi}\in \Psi[L]$ and $a'_\Psi,b'_\Psi,c'_{\Psi}\in \Psi[L']$ are distinct points in the liner $A$. Since $A$ is Pappian, the set 
$$
\begin{aligned}
T_A&\defeq A\cap \big((\Aline {a_\Psi}{b'_\Psi}\cap\Aline {a'_\Psi}{b_\Psi})\cup(\Aline {b_\Psi}{c'_\Psi}\cap\Aline {b'_\Psi}{c_\Psi})\cup(\Aline {a_\Psi}{c'_\Psi}\cap\Aline {a'_\Psi}{c_\Psi})\big)\\
&=A\cap\Psi[\{x,y,z\}]\subseteq A\cap (H\cup \{\Psi(z)\})=A\cap\{\Psi(z)\}
\end{aligned}
$$ has rank $\|T_A\|\in\{0,2\}$. Since $|T_A|=|A\cap\{\Psi(z)\}|\le|\{\Phi(z)\}|=1$, the set $T_A=A\cap\{\Phi(z)\}$ has rank $\|T_A\|=0$, which implies that $\Phi(z)\in \Pi\setminus A=\Aline xy=\Phi[\Aline xy]$ and $z\in\Aline xy$. Then $\|T\|=\|\{x,y,z\}\|=2$, which contradicts the choice of the points $a,b,c,a',b',c'$. This is a final contradiction showing that the projective liner $\overline X$ is Pappian.
\end{proof}



\begin{corollary}\label{c:completion-Pappian} Any projective completion of a Pappian liner is a Pappian projective liner.
\end{corollary}

\begin{proof} Let $Y$ be a projective completion of a Pappian liner $X\subseteq Y$. Then $Y$ is a $3$-long projective liner and $\overline{Y\setminus X}\ne Y$. By Corollary~\ref{c:Avogadro-projective}, the $3$-long projective liner $Y$ is $2$-balanced and hence the cardinal number $|Y|_2\ge 3$ is well-defined. If $|Y|_2=3$, then the Steiner projective liner $Y$ is Pappian, because it contains no $4$-long lines. So, assume that the liner $Y$ is $4$-long. Since the liner $X$ is Pappian, its horizon $H=Y\setminus X$ is flat in the $4$-long liner $Y$, by Lemma~\ref{l:Vlad6}. By Theorem~\ref{t:proaffine<=>proflat} and Proposition~\ref{p:projective-minus-flat}, the liner $X=Y\setminus H$ is proaffine and regular. Since $H$ is flat in $Y$, $|Y|_2\ge 4$ implies that the liner $X$ is $3$-long. 
By Theorem~\ref{t:Pappian=>compreg}, the Pappian proaffine regular liner $X$ is completely regular, and by Theorem~\ref{t:Pappian-completion}, the spread completion $\overline X$ of $X$ is a Pappian projective liner. By Corollary~\ref{c:pcompletion=scompletion}, the projective completion $Y$ of $X$ is isomorphic to the Pappian projective liner $\overline X$ and hence $Y$ is Pappian.
\end{proof}

Let us recall that a liner $X$ is \defterm{completely $\mathcal P$} for some property $\mathcal P$ of projective liners if $X$ is completely regular and its spread completion $\overline X$ has property $\mathcal P$. Theorems~\ref{t:Pappian=>compreg} and \ref{t:Pappian-completion} imply the following corollary answering Problem~\ref{prob:inner-completely-P} for the Pappus property.

\begin{corollary}\label{c:completelyP<=>} A liner is completely Pappian if and only if it is Pappian, proaffine and regular.
\end{corollary}

Let us recall that a property $\mathcal P$ of liners is called \defterm{complete} if a completely regular liner has property $\mathcal P$ if and only if its spread completion has property $\mathcal P$. Theorem~\ref{t:Pappian-completion} implies the following important fact.

\begin{corollary}\label{c:Pappus=>complete} The property of a liner to be Pappian is complete.
\end{corollary}

Next, we prove that a (3-long) projective liner is Pappian if and only if it is everywhere Pappian (if and) only if it is somewhere Pappian.

\begin{theorem}\label{t:projPapp<=>} For a projective liner $X$, the following conditions are equivalent:
\begin{enumerate}
\item $X$ is Pappian;
\item for every flat $H$ in $X$, the liner $X\setminus H$ is Pappian;
\item for every hyperplane $H\subset X$, the affine liner $X\setminus H$ is Pappian.
\end{enumerate}
If the projective liner $X$ is $3$-long, then the conditions \textup{(1)--(3)} are equivalent to the conditions:
\begin{enumerate}
\item[(4)] for some hyperplane $H$ in $X$, the affine liner $X\setminus H$ is Pappian;
\item[(5)] for some flat $H\ne X$, the liner $X\setminus H$ is Pappian.
\end{enumerate}
\end{theorem}

\begin{proof} The implication $(1)\Ra(2)$ follows from Proposition~\ref{p:Pappian-minus-flat}, and the implication $(2)\Ra(3)$ is trivial.
\vskip5pt

$(3)\Ra(1)$ Assume that for every hyperplane $H\subseteq X$, the liner $X\setminus H$ is Pappian. By Proposition~\ref{p:Papp<=>4-Papp}, to prove that the liner $X$ is Pappian, it suffices to check that every $4$-long plane $P$ in $X$ is Pappian. Fix any independent set $a,b,c\in P$. By the Kuratowski--Zorn Lemma, there exists a maximal independent set $I\subseteq X$ such that $\{a,b,c\}\subseteq I$. The maximal independence of $I$ ensures that the flat $H\defeq\overline{I\setminus\{a\}}$ is a hyperplane in $X$. By the assumption, the liner $X\setminus H$ is Pappian and so is the liner $(X\setminus H)\cap P=P\setminus H$. Since $a\notin H$ and $\{b,c\}\subseteq P\cap H$, the flat $H\cap P$ is a hyperplane in $P$. Then $P$ is a projective completion of the pappian liner  $P\setminus H$. By Corollary~\ref{c:completion-Pappian}, the projective liner $P$ is Pappian.
\smallskip

The implications $(3)\Ra(4)\Ra(5)$ are trivial, and the implication $(5)\Ra(1)$ follows from Corollary~\ref{c:completion-Pappian}.
\end{proof}

Theorem~\ref{t:projPapp<=>} implies the following corollary (which can be also deduced from Corollary~\ref{c:Pappus=>complete} and Theorem~\ref{t:everywhereP<=>somewhereP}).

\begin{corollary} A ($4$-long) projective liner $X$ is Pappian if and only if $X$ is everywhere Pappian (if and) only if $X$ is somewhere Pappian.
\end{corollary}

\begin{exercise} Prove that the projective space $\mathbb PV$ of every $R$-module $V$ over a field $R$ is a Pappian projective liner.
\end{exercise} 


\section{The Dual Pappus Axiom} 


Lines $L_1,\dots,L_n$ in a liner $X$ are called \index{coplanar lines}\index{lines!coplanar}\defterm{coplanar} if $\|L_1\cup\dots\cup L_n\|\le 3$.

\begin{definition} A liner $X$ is defined to satisfy the
 \index{Dual Pappus Axiom}\index{Axiom!Dual Pappus}\defterm{\sf Dual Pappus Axiom} if\newline for any distinct coplanar lines $A,B,C,A',B',C'\subseteq X$ 
with\newline $(A\cap B\cap C)\setminus(A'\cup B'\cup C')\ne\varnothing\ne(A'\cap B'\cap C')\setminus(A\cup B\cup C)$, the intersection\newline $\overline{(A\cap B')\cup(A'\cap B)}\cap\overline{(A\cap C')\cup(A'\cap C)}\cap\overline{(B\cap C')\cup(B'\cap C)}$ is not empty.
\end{definition} 

\begin{picture}(150,80)(70,-20)
%






\put(200,0){\line(4,1){105}}
\put(310,25){$C$}
\put(200,0){\line(3,1){90}}
\put(290,33){$B$}
\put(200,0){\line(2,1){70}}
\put(273,35){$A$}
\put(320,0){\line(-4,1){105}}
\put(201,24){$A'$}
\put(320,0){\line(-3,1){90}}
\put(220,33){$B'$}
\put(320,0){\line(-2,1){70}}
\put(240,35){$C'$}

\put(272,24){\color{teal}\circle*{3}}
\put(248,24){\color{blue}\circle*{3}}
\put(268.6,17.15){\color{teal}\circle*{3}}
\put(251.4,17.15){\color{blue}\circle*{3}}

\put(260,30){\color{cyan}\circle*{3}}
\put(260,15){\color{cyan}\circle*{3}}

\put(272,24){\color{teal}\line(-1,-2){12}}
\put(248,24){\color{blue}\line(1,-2){12}}
\put(260,0){\color{cyan}\line(0,1){30}}

\put(200,0){\circle*{3}}
\put(320,0){\circle*{3}}
\put(260,0){\color{red}\circle*{4}}


\end{picture}

\begin{theorem} A projective liner $X$ is Pappian if and only if $X$ satisfies the Dual Pappus Axiom.
\end{theorem}

\begin{proof}  Assume that a projective liner $X$ is Pappian. To prove the Dual Pappus Axiom, take any distinct coplanar lines $A,B,C,A',B',C'$ in $X$ such that $b\in (A\cap B\cap C)\setminus(A'\cup B'\cup C')$ and $b'\in (A'\cap B'\cap C')\setminus(A\cup B\cup C)$ for some distinct points $b,b'\in X$. Since the lines $A,B,C,A',B',C'$ are coplanar, by the $0$-parallelity of projective liners, there exist unique points $a,c,a',c'\in X$ such that 
$$\{a\}=A\cap C',\quad\{c\}=A\cap B',\quad \{a'\}=A'\cap B,\quad \{c'\}=A'\cap C.$$
Observe that $\{a,b,c\}\subseteq A\setminus A'$ and $\{a',b',c'\}\subseteq A'\setminus A$.

\begin{picture}(240,120)(-80,-25)
\put(0,0){\line(4,1){210}}
\put(216,52){$C$}
\put(0,0){\line(3,1){180}}
\put(182,61){$B$}
\put(0,0){\line(2,1){140}}
\put(142,72){$A$}
\put(240,0){\line(-4,1){210}}
\put(16,50){$A'$}
\put(240,0){\line(-3,1){180}}
\put(48,61){$B'$}
\put(240,0){\line(-2,1){140}}
\put(90,72){$C'$}

\put(144,48){\color{teal}\circle*{3}}
\put(143,52){\color{teal}$z$}
\put(96,48){\color{blue}\circle*{3}}
\put(93,52){\color{blue}$c$}
\put(137.14,34.3){\color{teal}\circle*{3}}
\put(136,26.5){\color{teal}$x$}
\put(102.9,34.3){\color{blue}\circle*{3}}
\put(94,24.5){\color{blue}$a'$}

\put(120,60){\color{cyan}\circle*{3}}
\put(118,63){$a$}
\put(120,30){\color{cyan}\circle*{3}}
\put(112,20){$c'$}

\put(144,48){\color{teal}\line(-1,-2){24}}
\put(96,48){\color{blue}\line(1,-2){24}}
\put(120,0){\color{cyan}\line(0,1){60}}

\put(0,0){\circle*{3}}
\put(-8,-3){$b$}
\put(240,0){\circle*{3}}
\put(244,-3){$b'$}
\put(120,0){\color{red}\circle*{4}}
\put(117,-10){$y$}
\end{picture}

Consider the unique points $x,y,z\in X$ such that $$\{x\}=\Aline{b}{c'}\cap\Aline {b'}{c}=B'\cap C,\quad\{y\}=\Aline a{c'}\cap\Aline{a'}c\quad\mbox{and}\quad\{z\}=\Aline a{b'}\cap\Aline {a'}b=B\cap C'.$$
The Pappus Axiom ensures that the points $x,y,z$ are collinear and hence 
$$y\in\Aline a{c'}\cap\Aline {a'}c\cap \Aline xy=\overline{(A\cap C')\cup(A'\cap C)}\cap\overline{(A'\cap B)\cup(A\cap B')}\cap\overline{(B'\cap C)\cup(B\cap C')},$$
witnessing that the dual Pappus Axiom holds.
\smallskip

 Now assume that a projective liner $X$ satisfies the Dual Pappus Axiom. To prove the  Pappus Axiom for $X$, take two concurrent lines $A,A'$ and distinct points $a,b,c\in A\setminus A'$ and $a',b',c'\in A'\setminus A$.
Consider the lines $$B\defeq\Aline{a'}b,\quad C\defeq \Aline b{c'},\quad B'\defeq\Aline {b'}c,\quad C'\defeq\Aline a{b'}$$ and observe that $b\in A\cap B\cap C\setminus(A'\cup B'\cup C')$ and $b'\in A'\cap B'\cap C'\setminus(A\cup B\cup C)$.

\begin{picture}(150,190)(-130,-30)

\put(-15,0){\line(1,0){155}}
\put(145,-5){$A$}
\put(0,0){\line(1,1){135}}
\put(140,135){$A'$}
\put(0,0){\line(-1,-1){10}}
{\linethickness{1pt}
\put(60,0){\line(1,0){60}}
\put(30,30){\line(1,1){90}}
}
\put(90,0){\color{blue}\line(-2,1){70}}
\put(10,37){\color{blue}$B$}
\put(90,0){\color{teal}\line(1,4){35}}
\put(123,144){\color{teal}$C$}
\put(120,0){\color{teal}\line(-1,1){60}}
\put(120,0){\color{teal}\line(1,-1){15}}
\put(135,-23){\color{teal}$B'$}
\put(120,0){\color{cyan}\line(-3,1){90}}
\put(60,0){\color{blue}\line(0,1){60}}
\put(60,0){\color{blue}\line(0,-1){15}}
\put(56,-25){\color{blue}$C'$}

\put(60,0){\color{cyan}\line(1,2){60}}

{\linethickness{1pt}
\put(60,15){\color{red}\line(4,1){36}}
}

\put(60,0){\circle*{3}}
\put(52,-8){$a$}
\put(90,0){\circle*{3}}
\put(87,-10){$b$}
\put(120,0){\circle*{3}}
\put(115,-8){$c$}

\put(30,30){\circle*{3}}
\put(36,30){$a'$}
\put(60,60){\circle*{3}}
\put(65,57){$b'$}
\put(120,120){\circle*{3}}
\put(123,115){$c'$}

\put(60,15){\color{blue}\circle*{3}}
\put(52,10){\color{blue}$z$}
\put(68.6,17.2){\color{cyan}\circle*{3}}
\put(64,22){\color{cyan}$y$}
\put(96,24){\color{teal}\circle*{3}}
\put(99,23){\color{blue}$x$}
\end{picture}

By the $0$-parallelity of projective liners, there exist unique points $x,y,z$ such that
$$\{x\}=\Aline {b'}c\cap\Aline b{c'}=B'\cap C,\quad \{y\}=\Aline a{c'}\cap\Aline{a'}c\quad\mbox{and}\quad\{z\}=\Aline a{b'}\cap\Aline {a'}b=B\cap C'.$$
By the Dual Pappus Theorem,  
$$
\begin{aligned}
\varnothing&\ne \overline{(A\cap B')\cup(A'\cap B)}\cap\overline{(A\cap C')\cap(A'\cap C)}\cap\overline{(B\cap C')\cup(B'\cap C)}\\
&=\Aline {a'}c\cap \Aline a{c'}\cap \Aline zx\subseteq \Aline {a'}c\cap\Aline a{c'}=\{y\},
\end{aligned}
$$
and hence $y\in\Aline zx$ and 
$$\|(\Aline b{c'}\cap\Aline {b'c})\cup(\Aline a{c'}\cap\Aline{a'}c)\cup(\Aline a{b'}\cap\Aline{a'}b)\|=\|\{x,y,z\}\|=\|\{x,z\}\|=2\in\{0,2\},$$ witnessing that $X$ satisfies the Pappus Axiom and hence is Pappian.
\end{proof}

\section{$4$-Pappian liners}

Let us recall that a subset $A$ of a liner $X$ is \defterm{closed} if for any distinct points $a,b,c,d\in A$ with $\Aline ab\ne\Aline cd$, we have $\Aline ab\cap\Aline cd\subseteq A$. The \defterm{closure} $\langle A\rangle$ of a subset $A$ in a liner $X$ is the smallest closed subset of $X$ that contains $A$. It is equal to the intersection of all closed subsets of $X$ that contain the set $A$. 

\begin{proposition}\label{p:Pappian<=>5-Pappian} For a liner $X$, the following conditions are equivalent:
\begin{enumerate}
\item $X$ is Pappian;
\item for every subset $A\subseteq X$ its closure $\langle A\rangle$ is a Pappian subliner of $X$.
\item for every $5$-element subset $A\subseteq X$, its closure $\langle A\rangle$ is a Pappian subliner of $X$.
\end{enumerate}
\end{proposition} 

\begin{proof} The implications $(1)\Ra(2)\Ra(3)$ follow immediately from the corresponding definitions.
\smallskip

$(3)\Ra(1)$ Assume that the closure of any $5$-element subset of $X$ is Pappian. To prove that $X$ is Pappian, take any concurrent lines $L,L'\subseteq X$ and distinct points $a,b,c\in L\setminus L'$ and $a',b',c'\in L'\setminus L$. We have to prove that the set $T\defeq(\Aline a{b'}\cap\Aline {a'}b)\cup(\Aline a{c'}\cap\Aline{a'}c)\cup(\Aline b{c'}\cap\Aline {b'}c)$ has rank $\|T\|\in\{0,2\}$ in $X$.
If $T=\varnothing$, then $\|T\|=0$ and we are done. So, assume that $T$ contains some point $x$. We lose no generality assuming that $x\in \Aline a{b'}\cap\Aline {a'}b$. 

In this case consider the $5$-element set $A\defeq \{x,a,a',c,c'\}$ and its closure $\langle A\rangle$ in $X$. By our assumption, the liner $\langle A\rangle$ is Pappian. Observe that $L\cap L'=\Aline ac\cap\Aline{a'}{c'}\subseteq \langle A\rangle$ and $\{b,b'\}\subseteq(\Aline x{a'}\cap\Aline ac)\cup(\Aline xa\cap\Aline {a'}{c'})\subseteq \langle A\rangle$. Therefore, $\{a,b,c,a',b',c'\}\subseteq\langle A\rangle$ and hence $T\subseteq \langle A\rangle$. Since the liner $\langle A\rangle$ is Pappian, the set $T$ has rank $\{0,2\}$ in $\langle A\rangle$ and also in $X$.
\end{proof}

Proposition~\ref{p:Pappian<=>5-Pappian} motivates the following definition.

\begin{definition} A liner $X$ is called \defterm{$4$-Pappian} if the closure $\langle A\rangle$ of any $4$-element subset $A$ in $X$ is a Pappian subliner of $X$.
\end{definition} 

\begin{exercise} Find an example of a $4$-Pappian projective plane which is not Pappian.
\smallskip

{\em Hint:} Any Desarguesian projective plane is $4$-Pappian, see Theorem~\ref{t:invertible-add<=>4-Pappian}.
\end{exercise}

\chapter{Gallucci liners}

In this section we study skew lines and reguli in liners, introduce  
\index[person]{Gallucci}Gallucci\footnote{{\bf Generoso Gallucci} (1874--1942), an Italian geometer, a member of Academy of Sciences of Napoli (since 1913). His researches in projective geometry were summarized in the volume ``Complementi di geometria proiettiva'' (1928), which he himself considered his main work.} liners, which are tightly related to reguli, and reveal the interplay between Pappian and Gallucci liners.

\section{Reguli in liners}

Let us recall that two distinct lines $L,\Lambda$ in a liner $X$ are \index{skew lines}\index{lines!skew}\defterm{skew} if they are not coplanar. This happens if and only if $\|L\cup\Lambda\|=4$. Observe that two lines in a projective liner are skew if and only if they are disjoint.

A family of lines $\F$ is \index{skew family of lines}\defterm{skew} if any distinct lines $L,\Lambda\in\F$ are skew. 

Given a family of lines $\F$ in a liner $X$, let $\F^\pm$ be the family of lines which are coplanar with every line in the family $|F|$. So, $\F^\pm=\{L\in\mathcal L:\forall F\in\F\;\;(\|F\cup L\|\le 3)\}$.  

\begin{lemma}\label{l:skewAB} For every skew lines $A,B$ in a weakly modular liner $X$ and every point $x\in\overline{A\cup B}$, there exists a line $L_x\in\{A,B\}^\pm$ such that $x\in L_x$.
\end{lemma}

\begin{proof} If $x\in A\cup B$, then we lose no generality assuming that $x\in A$. In this case, take any point $y\in B$ and observe that the line $L_x\defeq\Aline xy$ belongs to the family $\{A,B\}^\pm$ and contains the point $x$. So, assume that $x\notin A\cup B$.

Since $x\notin \overline{A\cup B}\setminus(A\cup B)$, the flats $\overline{A\cup\{x\}}$ and $\overline{B\cup\{x\}}$ are planes in the flat $\overline{A\cup B}$. Taking into account that the lines $B,C$ are skew, we conclude that $4=\|A\cup B\|\le\|\overline{A\cup\{x\}}\cup\overline{B\cup\{x\}}\|\le \|\overline{A\cup\{x\}\cup B}\|=\|\overline{A\cup B}\|=4$. The weak modularity of $X$ ensures that the flat $L_x\defeq \overline{A\cup\{x\}}\cap\overline{B\cup\{x\}}$ has rank 
$$\|L_x\|=\|\overline{A\cup\{x\}}\|+\|\overline{B\cup\{x\}}\|-\|\overline{A\cup\{x\}}\cup\overline{B\cup\{x\}}\|=3+3-4=2$$ and hence is a line such that $L_x\in\{A,B\}^\pm$ and $x\in L_x$. 
\end{proof}

\begin{proposition}\label{p:regulus1} Let $A,B,C$ be skew lines in a liner $X$.
\begin{enumerate}
\item If the liner $X$ is $3$-ranked, then the family $\{A,B,C\}^\pm$ consists of pairwise disjoint lines.
\item If the liner $X$ is $4$-ranked, then $\bigcup\{A,B,C\}^\pm\subseteq\overline{A\cup B}\cap\overline{B\cup C}\cap\overline{A\cup C}$.
\item If the liner $X$ is weakly modular, then\newline $(A\cap\overline{B\cup C})\cup (B\cap\overline{A\cup C})\cup(C\cap\overline{A\cup C})\subseteq\bigcup\{A,B,C\}^\pm$.
\item If the liner $X$ is  weakly modular and $\|A\cup B\cup C\|=4$, then\newline  
$A\cup B\cup C\subseteq \bigcup\{A,B,C\}^\pm\subseteq\overline{A\cup B}=\overline{B\cup C}=\overline{A\cup C}=\overline{A\cup B\cup C}$.
\item If the liner $X$ is proaffine and regular, then the family $\{A,B,C\}^\pm$ consists of pairwise skew lines.
\item If the liner $X$ is proaffine and regular and $|\{A,B,C\}^\pm|\ge 2$, then $\|A\cup B\cup C\|=4$.
\end{enumerate}
\end{proposition}

\begin{proof} 1. Assume that the liner $X$ is $3$-ranked and take any distinct lines $L,\Lambda\in\{A,B,C\}^\pm$. To derive a contradiction, assume that the lines $L$ and $\Lambda$ have a common point $x$. Since the lines $A,B,C$ are skew, they are parwise disjoint and hence at least two lines among the lines $A,B,C$ do not contain the point $x$. We lose no generality assuming that $x\notin A\cup B$. Then $\overline{A\cup\{x\}}$ is a plane. Since $L\in\{A,B,C\}^\pm\subseteq\{A\}^\pm$, the lines $L$ and $A$ are coplanar and hence $\overline{A\cup \{x\}}\subseteq \overline{A\cup L}\subseteq P$ for some plane $P$. Since the liner $X$ is $3$-ranked, the planes $\overline{A\cup\{x\}}$ and $P$ coincide. Therefore, $L\subseteq P=\overline{A\cup\{x\}}$. By analogy we can prove that $\Lambda\subseteq \overline{A\cup\{x\}}$.

Since $L,\Lambda$ are two distinct lines in the plane $\overline{A\cup\{x\}}$, the flat hull $\overline{L\cup\Lambda}$ is a plane in the plane $\overline{A\cup\{x\}}$. The $3$-rankendess of $X$ ensures that $A\subseteq \overline{A\cup\{x\}}=\overline{L\cup\Lambda}$. By analogy we can prove that $B\subseteq\overline{L\cup\Lambda}$. Then the lines $A,B$ are coplanar and hence not skew. This contradiction show that distinct lines in the family $\{A,B,C\}^\pm$ are disjoint. 
\smallskip

2. Assume that the liner $X$ is $4$-ranked. Given any point $x\in\bigcup\{A,B,C\}^\pm$, we should prove that $x\in\overline{A\cup B}\cap\overline{A\cup C}\cap\overline{B\cup C}$. Since $x\in\bigcup\{A,B,C\}^\pm$, there exists a line $L\in\{A,B,C\}^\pm$ such that $x\in L$. Taking into account that the lines $A,B,C$ are skew and complanar to the line $L$, we conclude that $L\notin\{A,B,C\}$ and hence there exist points $a\in A\setminus L$, $b\in B\setminus L$ and $c\in C\setminus L$. Since the lines $A$ and $L$ are complanar, there exists a plane $P$ containing the lines $A$ and $L$. Since the $4$-ranked liner $X$ is $3$-ranked, the plane $\overline{L\cup\{a\}}\subseteq P$ coincides with the plane $P$ and hence $A\subseteq \overline{L\cup\{a\}}$. By analogy we can prove that $B\subseteq\overline{L\cup\{b\}}$ and $C\subseteq\overline{L\cup\{c\}}$. Then $A\cup B\subseteq\overline{L\cup\{a,b\}}$ and hence $4=\|A\cup B\|\le\|L\cup\{a,b\}\|\le 4$. The $4$-rankedness of the liner $X$ ensures that $x\in L\subseteq \overline{L\cup\{a,b\}}=\overline{A\cup B}$. By analogy we can prove that $x\in L\subseteq\overline{B\cup C}\cap\overline{A\cup C}$.
\smallskip

3. Assume that the liner $X$ is weakly modular, and take any point $$x\in (A\cap\overline{B\cup C})\cup (B\cap\overline{A\cup C})\cup(C\cap\overline{A\cup C}).$$ We lose no generality assuming that $x\in A\cap\overline{B\cup C}\subseteq\overline{B\cup C}\setminus(B\cup C)$. By Lemma~\ref{l:skewAB}, there exists a line such that $L_x\in\{B,C\}^\pm$ and $x\in L_x$. Then $x\in L_x\cap A$ and hence $x\in L_x\in \{A,B,C\}^\pm$ and $x\in \bigcup\{A,B,C\}^\pm$. 
\smallskip

4. Assume that the liner $X$ is weakly modular and $\|A\cup B\cup C\|=4$. By Theorem~\ref{t:w-modular<=>}, the weakly modular liner $X$ is ranked. Then the equality $\|A\cup B\cup C\|=4=\|A\cup B\|=\|B\cup C\|=\|A\cup C\|$ implies the equalities $\overline{A\cup B\cup C}=\overline{A\cup B}=\overline{B\cup C}=\overline{A\cup C}$. 
By Proposition~\ref{p:regulus1}(3),
$$
\begin{aligned}
A\cup B\cup C&=(A\cap\overline{B\cup C})\cup(B\cap\overline{A\cup C})\cap(C\cap\overline{A\cup B})\subseteq \textstyle\bigcup\{A,B,C\}^\pm\\
&\subseteq \overline{A\cup B}\cap\overline{B\cup C}\cap\overline{A\cup C}=\overline{A\cup B}=\overline{B\cup C}=\overline{A\cup C}=\overline{A\cup B\cup C}.
\end{aligned}
$$

5. Assume that the liner $X$ is proaffine and regular.  By Corollary~\ref{c:proregular=>ranked}, the regular  liner $X$ is ranked. Given any distinct lines $L,\Lambda\in\{A,B,C\}^\pm$, we should prove that the lines $L,\Lambda$ are skew. By Proposition~\ref{p:regulus1}(1) the lines $L,\Lambda$ are disjoint. Assuming that the lines $L,\Lambda$ are not skew, we conclude that $L,\Lambda$ are coplanar and hence parallel, by Corollary~\ref{c:parallel-lines<=>}. Since the lines $A,B,C$ are skew, at most one line among the lines $A,B,C$ is contained in the plane $\overline{L\cup\Lambda}$. We lose not generality assuming that the lines $A,B$ are not contained in the plane $\overline{L\cup\Lambda}$. Since the lines $L,\Lambda$ are coplanar with the line $A$, the flats $\overline{L\cup A}$ and $\overline{\Lambda\cup A}$ are planes in $X$, distinct from the plane $\overline{L\cup\Lambda}$. The $3$-rankedness of $X$ ensures that $\overline{L\cup A}\cap\overline{L\cup \Lambda}=L$, $\overline{L\cup B}\cap\overline{L\cup\Lambda}=\Lambda$, and $\overline{A\cup L}\cap\overline{A\cup\Lambda}=A$. Then 
$$A\cap L=(\overline{A\cup\Lambda}\cap\overline{A\cup L})\cap (\overline{A\cup L}\cap\overline{L\cup\Lambda})=(\overline{A\cup\Lambda}\cap\overline{L\cup\Lambda})\cap(\overline{A\cup L}\cap\overline{L\cup\Lambda})=\Lambda\cap L=\varnothing$$and hence the disjoint coplanar lines $L$ and $A$ are parallel, by Corollary~\ref{c:parallel-lines<=>}. By analogy we can prove that the lines $L$ and $B$ are parallel. By Theorem~\ref{t:Proclus-lines}, $A\parallel L\parallel B$ implies $A\parallel B$ and hence the lines $A,B$ are coplanar, which is a contradiction showing that any distinct lines in the family $\{A,B,C\}^\pm$ are skew. 
\smallskip

6. Assume that the liner $X$ is proaffine and regular and the family $\{A,B,C\}^\pm$ contains two distinct lines $L,\Lambda$. By Proposition~\ref{p:regulus1}(5), the lines $L\cup \Lambda$ are skew and hence $\|L\cup\Lambda\|=4$. By Proposition~\ref{p:regulus1}(2), $L\cup\Lambda\subseteq\bigcup\{A,B,C\}^\pm\subseteq \overline{A\cup B}$ and hence $4=\|L\cup\Lambda\|\le\|\overline{A\cup B}\|=4$, which implies $\overline{L\cup\Lambda}=\overline{A\cup B}$, by the $4$-rankedness of the regular liner $X$. By analogy we can prove that $\overline{B\cup C}=\overline{L\cup\Lambda}=\overline{A\cup C}$. Then $\overline{A\cup B\cup C}=\overline{L\cup\Lambda}$ and hence $\|A\cup B\cup C\|=\|L\cup\Lambda\|=4$.
\end{proof}

\begin{definition} A family of lines $\mathcal R$ in a liner $X$ is called a \index{regulus}\defterm{regulus} in $X$ if $|\mathcal R|\ge 3$ and $\mathcal R=\{A,B,C\}^\pm$ for some skew lines $A,B,C$ in $X$.
\end{definition}

By Proposition~\ref{p:regulus1}, any regulus $\mathcal R=\{A,B,C\}^\pm$ is a proaffine regular liners $X$ consists of skew lines. Given any lines $A',B',C',A'',B'',C''\in\{A,B,C\}^\pm$ with $|\{A',B',C'\}|=3=|\{A'',B'',C''\}|$, we can consider the reguli $\{A',B',C'\}^\pm$ and $\{A'',B'',C''\}^\pm$ and ask when those reguli coincide. This happens for reguli in Gallucci liners, introduced in the next section.

\section{Reguli in Gallucci liners}

In the following definition we identify the number $8$ with the set $\{0,1,2,3,4,5,6,7\}$ of numbers smaller than $8$.

\begin{definition} A family of eight distinct lines $(L_i)_{i\in 8}$ in a liner $X$ forms a \index{Gallucci configuration}\defterm{Gallucci configuration} if for every numbers $i,j\in 8$ with $|i-j|<7$, the flat $\overline{L_i\cup L_j}$ has rank $\|L_i\cup L_j\|\in|i-j|+2\IZ$. A liner $X$ is called \index{Gallucci liner}\index{liner!Gallucci}\defterm{Gallucci} if for every Gallucci configuration $(L_i)_{i\in 8}$ in $X$, the lines $L_0,L_7$ are coplanar.
\end{definition}

Observe that two lines $L,L'$ in a liner have 
$$\|L\cup L'\|=
\begin{cases}
2&\mbox{iff $L$ and $L'$ coincide};\\
3&\mbox{iff $L$ and $L'$ are coplanar and distinct};\\
4&\mbox{iff $L$ and $L'$ are skew.}
\end{cases}
$$
So, lines $L_i,L_j$ in a Gallucci configuration $(L_i)_{i\in 8}$ in a Gallucci line are skew if and only if $i-j$ is even.

\begin{exercise} Show that every liner $X$ of rank $\|X\|\le 3$ is Gallucci.
\smallskip

{\em Hint:} Observe that every liner $X$ containing a Gallucci configuration  $(L_i)_{i\in 8}$ has rank $\|X\|\ge\|L_0\cup L_2\|=4$.
\end{exercise}





\begin{theorem}\label{t:Gallucci} Let $X$ be a Gallucci proaffine regular liner. For every regulus $\mathcal R$ in $X$, the family $\mathcal R^\pm$ is a regulus in $X$, equal to the regulus $\{A,B,C\}^\pm$ for every distinct lines $A,B,C\in\mathcal R$. Moreover, $\mathcal R=(\mathcal R^\pm)^\pm$ and $\bigcup\mathcal R=\bigcup\mathcal R^\pm$.
\end{theorem}

\begin{proof} Since $\mathcal R$ is a regulus in $X$, there exist skew lines $A',B',C'$ in $X$ such that $\mathcal R=\{A',B',C'\}^\pm$. Since $|\mathcal R|\ge 3$, we can apply Proposition~\ref{p:regulus1}(6) and conclude that $\|A'\cup B'\cup C'\|=4$. Choose any distinct lines $A,B,C\in\mathcal R$. By Proposition~\ref{p:regulus1}(5), the lines $A,B,C$ are skew. The family $\{A,B,C\}^\pm$ contains three distinct lines $A',B',C'$ and hence is a regulus. It follows from $\{A,B,C\}\subseteq \mathcal R$ that $\mathcal R^\pm\subseteq\{A,B,C\}^\pm$. To show that $\{A,B,C\}^\pm\subseteq \mathcal R^\pm$, fix any line $D'\in\{A,B,C\}^\pm$. We have to prove that $D'\in\mathcal R^\pm$. This will follow as soon as we check that the line $D'$ is coplanar with every line $D\in\mathcal R$. 

If $D'\in\{A',B',C'\}$, then the line $D\in\mathcal R=\{A',B',C'\}^\pm$ is coplanar with the line $D'$, by the definition of the family $\{A',B',C'\}^\pm$. 
If $D\in\{A,B,C\}$, then the coplanarity of the lines $D$ and $D'$ follows from the choice of the line  $D'\in\{A,B,C\}^\pm$. So, assume that $D\notin\{A,B,C\}$ and $D'\notin \{A',B',C'\}$. Consider the lines
$$L_0\defeq D,\quad L_1\defeq A',\quad L_2\defeq A,\quad L_3\defeq B',\quad L_4\defeq B,\quad L_5\defeq C',\quad L_6\defeq C,\quad L_7\defeq D'$$ and observe that $(L_i)_{i\in 8}$ is a Gallucci configuration. Since the liner $X$ is Gallucci, the lines $L_0=D$ and $L_7=D'$ are coplanar. Therefore, $(\{A',B',C'\}^\pm)^\pm=\mathcal R^\pm=\{A,B,C\}^\pm$ is a regulus in $X$.
\smallskip

By analogy we can prove that $(\{A,B,C\}^\pm)^\pm=\{A',B',C'\}^\pm=\mathcal R$, which implies that $$(\mathcal R^\pm)^\pm=(\{A,B,C\}^\pm)^\pm=\{A',B',C'\}^\pm=\mathcal R.$$

Finally, we show that $\bigcup\mathcal R=\bigcup\mathcal R^\pm$. Given any point $x\in\bigcup\mathcal R$, find a line $A\in\mathcal R$ containing the point $x$. Since $|\mathcal R|\ge 3$, there exist distinct lines $B,C\in\mathcal R\setminus\{A\}$. By Proposition~\ref{p:regulus1}(5,6), $A,B,C$ are skew lines with $\|A\cup B\cup C\|=4$. We already know that $\mathcal R^\pm=\{A,B,C\}^\pm$. 
By Theorem~\ref{t:w-modular<=>}, the regular liner $X$ is weakly modular. By Proposition~\ref{p:regulus1}(4), $x\in A\cup B\cup C\subseteq\bigcup\{A,B,C\}^\pm=\bigcup\mathcal R^\pm$. Therefore, $\bigcup\mathcal R\subseteq\bigcup\mathcal R^\pm$. Applying the same argument to the regulus $\mathcal R^\pm$ and taking into account that $(\mathcal R^\pm)^\pm=\mathcal R$, we conclude that $\bigcup\mathcal R^\pm\subseteq\bigcup(\mathcal R^\pm)^\pm=\bigcup\mathcal R$ and hence $\bigcup\mathcal R=\bigcup\mathcal R^\pm$. 
\end{proof}

\begin{remark} By Theorem~\ref{t:Gallucci}, for every regulus $\mathcal R$ in a Gallucci liner $X$, the family $\mathcal R^\pm$ is a regulus, called to \index{dual regulus}\defterm{dual regulus} to the regulus $\mathcal R$. Dual reguli have the same unions. 
\end{remark}

\begin{exercise} Show that the hyperbolic paraboloid $H=\{(x,y,z)\in \IR^3:z=xy\}$ is equal to the union $\bigcup\mathcal R$ of some regulus $\mathcal R$ in $\IR^3$.  Generalize this fact to vector spaces over an arbitrary field.

{\em Hint:}  Observe that families lines 
$$
\begin{aligned}
L_a&\defeq\{(x,y,z)\in\IR^3:  z=ax\;\wedge\;y=a\},\\
\Lambda_b&\defeq\{(x,y,z)\in\IR^3:  x=b\;\wedge\;z=by\},
\end{aligned}$$determine two dual reguli $\{L_a:a\in \IR\setminus\{0\}\}$ and $\{\Lambda_b:b\in\IR^*\setminus\{0\}\}$ whose unions coincide with the hyperbolic paraboloid $H$.
\end{exercise}

\begin{exercise}[Bell\footnote{{\bf Robert J. T. Bell} (1876 -- 1963) was a Scottish mathematician. Bell was educated at Hamilton Academy from which he matriculated at the University of Glasgow, having won a high placement in the university’s Open Bursary Competition. Bell graduated in 1898 as M.A. (with First Class Honours) in Mathematics and Natural Philosophy. Appointed a William Ewing Fellow, Bell continued at the university as a tutorial assistant and in 1901 was promoted to junior assistant. In 1911 Bell was appointed Lecturer in Mathematics and awarded a D.Sc. by the university for his treatise on the geometry of three dimensions, published in book form in 1910 as ``An Elementary Treatise on Co-ordinate Geometry of Three Dimensions''. An instant success, this textbook was to be translated into other languages, including Japanese and three of the languages of the Indian sub-continent. The textbook ran to a third edition (1944) and, from 1938, chapters 1--9 were issued separately as Co-ordinate Solid Geometry. It has since been reprinted (BiblioBazaar, 2009.) In March 1899 Bell had become a member of Edinburgh Mathematical Society and from 1911–1920 he served as editor of the Society’s journal, the Proceedings.  On 6 March 1916, Bell was elected to the Royal Society of Edinburgh. 
In 1920 Bell was appointed Professor of Pure and Applied Mathematics at the University of Otago, Dunedin, New Zealand.}, 1910]\index[person]{Bell} Show that the one-sheeted hyperboloid 
$$H\defeq\{(x,y,z)\in \IR^3:x^2+y^2-z^2=1\}$$ is equal to the union $\bigcup\mathcal R$ of some regulus $\mathcal R$ in $\IR^3$.  Generalize this fact to vector spaces over  fields of characteristic $\ne 2$.
\smallskip

{\em Hint:}  Rewrite the equation $x^2+y^2-z^2=1$ as $(x-z)(x+z)=(1-y)(1+y)$ and observe that the families lines 
$$
\begin{aligned}
L_a&\defeq\{(x,y,z)\in\IR^3:  x-z=a(1-y)\;\wedge\;a(x+z)=1+y\},\\
\Lambda_b&\defeq\{(x,y,z)\in\IR^3:  x+z=b(1-y)\;\wedge\;b(x-z)=1+y\},
\end{aligned}$$determine two dual reguli $\{L_a:a\in \IR\setminus\{0\}\}$ and $\{\Lambda_b:b\in\IR\setminus\{0\}\}$ whose unions coincide with the hyperboloid $H$.
\end{exercise}

\section{Pappian proaffine regular liners are Gallucci}

\begin{theorem}[Gallucci, 1906]\label{t:Pappian=>Gallucci} Every Pappian projective liner is Gallucci.
\end{theorem}

\begin{proof}  Let $X$ be a Pappian projective liner. To show that $X$ is Gallucci, take any family of distinct lines $(L_i)_{i\in 8}$ in $X$ such that $\|L_i\cup L_j\|\in|i-j|+2\IZ$ for any numbers $i,j\in 8$ with $|i-j|<7$. We have to prove that the lines $L_0,L_7$ are coplanar. 

Since $L_0,L_2,L_4,L_6\in \{L_1,L_3,L_5\}^\pm$ and $L_1,L_3,L_5,L_7\in\{L_2,L_4,L_6\}^\pm$, there exist points 
$$
\begin{gathered}
a\in L_1\cap L_4,\quad b\in L_1\cap L_2,\quad c\in L_1\cap L_0,\quad u\in L_3\cap L_4\\
a'\in L_6\cap L_3\quad b'\in L_6\cap L_5,\quad c'\in L_6\cap L_7,\quad v\in L_2\cap L_5.
\end{gathered}
$$
Since $a,b,c,a',b',c'$ are points in the projective plane $\overline{L_1\cup L_6}$, there exist unique points $x,y,z\in X$ such that 
$$
\begin{aligned}
x&\in\Aline a{b'}\cap \Aline {a'}b=\overline{(L_1\cap L_4)\cup(L_5\cap L_6)}\cap\overline{(L_6\cap L_3)\cup(L_1\cap L_2)},\\
y&\in\Aline a{c'}\cap\Aline {a'}c=\overline{(L_1\cap L_4)\cup(L_6\cap L_7)}\cap\overline{(L_6\cap L_3)\cup(L_1\cap L_0)},\\
z&\in\Aline b{c'}\cap\Aline {b'}c=\overline{(L_1\cap L_2)\cup(L_6\cap L_7)}\cap\overline{(L_6\cap L_5)\cup(L_1\cap L_0)}.
\end{aligned}
$$
Since $X$ is Pappian, the points $x,y,z$ are collinear.




We claim that the points $x,u,v$ are collinear. For this consider the planes $\overline{L_2\cup L_3}$ and $\overline{L_4\cup L_5}$ and observe that $\{u,v\}\subseteq \overline{L_2\cup L_3}\cap\overline{L_4\cup L_5}$. Then 
$$x\in \overline{(L_1\cap L_4)\cup(L_5\cap L_6)}\cap\overline{(L_1\cap L_2)\cup(L_6\cap L_3)}\subseteq \overline{L_4\cup L_5}\cap\overline{L_2\cup L_3}=\Aline uv$$ and hence $v\in\Aline xu\subseteq \overline{\{y,z,u\}}$, which implies that $\overline{\{y,z,u,v\}}$ is a plane in liner $X$. Since the liner $X$ is projective, the lines $\Aline yu$ and $\Aline zv$ in the plane $\overline{\{y,z,u,v\}}$ have a common point $w$. We claim that $w\in L_0\cap L_7$. 

Observe that the points $y\in \overline{(L_1\cap L_0)\cup(L_6\cap L_3)}\cap\overline{(L_1\cap L_4)\cup(L_6\cap L_7)}$ and $u$ belong to the planes $\overline{L_3\cup L_0}$ and $\overline{L_4\cup L_7}$, which implies that $\overline{L_3\cup L_0}\cap\overline{L_4\cup L_7}=\Aline yu$. By analogy we can show that $\Aline zv=\overline{L_5\cup L_0}\cap\overline{L_2\cup L_7}$. Then 
\begin{multline*}
w\in \Aline yu\cap\Aline zv=(\overline{L_3\cup L_0}\cap\overline{L_4\cup L_7})\cap(\overline{L_5\cup L_0}\cap\overline{L_2\cup L_7})\\
=(\overline{L_3\cup L_0}\cap\overline{L_5\cup L_0})\cap(\overline{L_4\cup L_7}\cap\overline{L_2\cup L_7})=L_0\cap L_7,
\end{multline*}
and hence the lines $L_0,L_7$ are coplanar. 
\end{proof}

Theorem~\ref{t:Pappian=>Gallucci} admits the following self-generalization.

\begin{theorem}\label{t:Pappian=>Gallucci-my} Every Pappian proaffine regular liner $X$ is Gallucci.
\end{theorem}

\begin{proof} Let $(L_i)_{i\in 8}$ be a Gallucci configuration in the Pappian proaffine regular line $X$. Since $X$ contains the skew lines $L_0,L_2$, $\|X\|\ge \|L_0\cup L_2\|=4$. By Theorem~\ref{t:proaffine3=>compregular}, the proaffine regular liner $X$ of rank $\|X\|\ge 4$ is completely regular.  By Theorem~\ref{t:Pappian-completion}, the spread completion $\overline X$ of $X$ is a Pappian projective liner. For a line $L$ in $X$ we shall denote by $\overline L$ the flat hull of $L$ in $\overline X$. It follows that two lines $L,L'\subseteq X$ are coplanar in $X$ if and only if the lines $\overline L,\overline {L'}$ have a common point (and hence are coplanar) in the projective liner $\overline X$.  This implies that $(\overline L_i)_{i\in 8}$ is a Gallucci configuration in the Pappian projective liner $\overline X$. By Theorem~\ref{t:Pappian=>Gallucci}, the lines $\overline L_0$ and $\overline{L_7}$ have a common point in $\overline X$ and hence the lines $L_0$ and $L_7$ are coplanar in $X$.
\end{proof}

\begin{theorem}\label{t:Pappian<=>Gallucci} A $3$-long proaffine regular liner $X$ of rank $\|X\|\ne 3$ is Gallucci if and only if it is Pappian. 
\end{theorem}

\begin{proof} The ``only if'' follows from Theorem~\ref{t:Pappian=>Gallucci-my}. To prove the ``if'' part, assume that the liner $X$ is Gallucci. To prove that $X$ is Pappian, take any concurrent lines $L,L'$ in $X$ and distinct points $a,b,c\in L\setminus L'$ and $a',b',c'\in L'\setminus L$. Since $X$ contains distinct lines $L,L'$, it has rank $\|X\|\ge 3$. Taking into account that $\|X\|\ne 3$, we conclude that $\|X\|\ge 4$. By Theorem~\ref{t:proaffine3=>compregular}, the proaffine regular liner $X$ of rank $\|X\|\ge 4$ is completely regular. Then the spread completion $\overline X$ of $X$ is a $3$-long projective liner, according to Corollary~\ref{c:spread-3long}. For two distinct points $x,y$ in $\overline X$ we denote by $\Aline xy$ the unique line in $\overline X$ that contains the points $x,y$. 

Since $\overline X$ is projective, there exist points $x\in\Aline a{b'}\cap\Aline {a'}b$, $y\in\Aline a{c'}\cap\Aline {a'}c$ and $z\in\Aline b{c'}\cap\Aline{a'}c$. We shall prove that the points $x,y,z$ are collinear in the projective liner $\overline X$. Since $\|X\|\ge 4$, we can choose a point $u\in X\setminus\overline{L\cup L'}$ and consider the flat $V\defeq\overline{L\cup L'\cup\{u\}}$ of rank $\|V\|=4$. Since $X$ is $3$-long, there exists a point $v\in (X\cap\Aline xu)\setminus\{x,u\}$. Consider the lines 
$$L_1\defeq L,\quad L_2\defeq X\cap\Aline bv,\quad L_3\defeq X\cap\Aline u{a'},\quad L_4\defeq X\cap\Aline au,\quad L_5\defeq X\cap \Aline v{b'},\quad L_6\defeq L'.$$  
It is easy to see that $\|L_3\cup L_5\|=4=\|L_2\cup L_4\|$. 
By Lemma~\ref{l:skewAB}, there exists a line $L_0\in \{L_3,L_5\}^\pm$ such that $c\in L_0$. By analogy, we can find a line $L_7\in\{L_2,L_4\}^\pm$ such that $c'\in L_7$.
Then $(L_i)_{i\in 8}$ is Gallucci configuration in $X$. Since the liner $X$ is Gallucci, the  lines $L_0$ and $L_7$ are coplanar in $X$ and hence the lines $\overline{L_0}$ and $\overline{L_7}$ have a common point $w\in\overline X$. 

Observe that the points $y,u,w$ belong to the planes $\overline{L_4\cup L_7}$ and $\overline{L_0\cup L_3}$ and hence $\overline{\{y,u,w\}}=\overline{L_4\cup L_7}\cap\overline{L_0\cup L_3}$ and the points $y,u,w$ are collinear. By analogy we can prove that the points $z,v,w$ are collinear. Then the points $u,v$ belong to the plane $\overline{\{y,z,w\}}$ and hence $x\in\Aline uv\subseteq \overline{\{y,z,w\}}\cap\overline{L\cup L'}=\Aline yz$, witnessing that the points $x,y,z$ are collinear in the projective liner $\overline X$. Since the set $\overline X\setminus X$ is flat in $\overline X$, the set $$T\defeq X\cap\big((\Aline a{b'}\cap\Aline {a'}b)\cup(\Aline a{c'}\cap\Aline{a'}c)\cup(\Aline b{c'}\cap\Aline{b'}c)\big)=X\cap\{x,y,z\}$$has rank $\|T\|\in\{0,2\}$ in $X$, witnessing that the liner $X$ is Pappian.
\end{proof} 

\begin{remark} The ideas of the proofs of Theorems~\ref{t:Pappian=>Gallucci} and \ref{t:Pappian<=>Gallucci} are taken from the papers \cite{Horvath2019} and \cite{Horvath2023} of \index[person]{Horv\'ath}Horvath\footnote{{\bf \'Akos G. Horv\'ath} (born in 1960), a Hungarian mathematician, a head of Deptment of Geometry at Institute of Mathematics in Budapest University of Technology and Economics.}.
\end{remark}

\chapter{Central projections in projective spaces}

\section{Projective spaces} 

\begin{definition} A \index{projective space}\index{space!projective}\defterm{projective space} is any $3$-long projective liner $X$ of rank $\|X\|\ge 3$.
\end{definition}

Theorem~\ref{t:proj-1-axioms}  implies that projective spaces admit the following first-order axiomatization.

\begin{theorem} A set $X$ endowed with a ternary relation $\Af\subseteq X^3$ is a projective space if and only if the following four axioms are satisfied:
\begin{itemize}
\item[{\sf (L)}] {\sf Long Lines:} $\forall x,z\in X\;\;\big(x\ne z\;\Leftrightarrow\;\exists y\in X\;(\Af xyz\;\wedge\;x\ne y\ne z)\big)$;
\item[{\sf (R)}] {\sf Reflexivity:} $\forall x,y\in X\;\;(\Af xxy\;\wedge\;\Af xyy)$;
\item[{\sf (E)}] {\sf Exchange:} $\forall a,b,x,y\in X\;\big(( \Af axb\wedge \Af ayb\wedge x\ne y)\Rightarrow (\Af xay\wedge\Af xby)\big)$;
\item[{\sf (V)}] {\sf Veblen:} $\forall o,x,y,p,u\in X\;\big(\Af ypx\;\wedge\;p\ne y\;\wedge\; \Af ouy)\;\Ra\;\exists v\;(\Af ovx\;\wedge \Af upv)$;
\item[{\sf(D)}] {\sf Dimension:} $\exists x,y,z\in X \;(x\ne z\;\wedge\;\neg\Af xyz)$.
\end{itemize}
\end{theorem}

The following proposition is a projective counterpart of Proposition~\ref{p:flat-relation}. 

\begin{proposition}\label{p:central-projection} Let $X$ be a projective space. For two flats $A,B\subseteq X$ and a point $o\in X\setminus (A\cup B)$, the relation
$$F\defeq\{(x,y)\in A\times B:y\in \Aline ox\}$$ is a flat injective function having the following properties.
\begin{enumerate}
\item $F(x)=x$ for every $x\in A\cap B$;
\item $\dom[F]=A$ if and only if $A\subseteq \overline{\{o\}\cup B}$;
\item $\rng[F]=B$ if and only if $B\subseteq \overline{\{o\}\cup A}$;
\item If $A,B$ are two coplanar lines and $o\in\overline{A\cup B}$, then $\dom[F]=A$ and $\rng[F]=B$.
\end{enumerate}
\end{proposition}

\begin{proof} 0. First we show that $F$ is a function. Given any pairs $(x,y),(x,z)\in F$, we should prove that $y=z$. In the opposite case, $y,z\in\Aline ox$ implies $o\in\Aline ox=\Aline yz\subseteq B$, which contradicts the choice of the point $o$. By analogy we can prove that the function $F$ is injective.
\smallskip 

To prove that the function $F$ is flat, fix any flat $C\subseteq X$. We should prove that the sets $F[C]$ and $F^{-1}[C]$ are flat in $X$. To show that $F[C]$ is flat, take any distinct points $x,y\in F[C]$. We have to check that $\Aline xy\subseteq F[C]$. Fix any point $z\in\Aline xy$. Since $x,y\in F[C]$, there exist points $a,b\in C$ such that $(a,x),(b,y)\in F$. The definition of the function $F$ ensures that  $x\in\Aline oa$ and $y\in \Aline ob$. Then $z\in\Aline xy\subseteq \overline{\{o,a,b\}}$.  Since $z\in\Aline xy\subseteq B$, the point $z$ is distinct from the point $o\notin B$ and hence $\Aline oz$ is a line in the plane $\overline{\{o,a,b\}}$. By the $0$-parallelity of the projective space $X$, the lines $\Aline oz$ and $\Aline ab$ in the plane $\overline{\{o,a,b\}}\subseteq X$ have a common point $c$. It follows from $c\in\Aline ab\subseteq C\cap A\subseteq X\setminus\{o\}$ that $o\ne c$ and hence $z\in \Aline oz=\Aline oc$. The definition of the function $F$ ensures that $z=F(c)\in F[C]$, witnessing that the set $F[C]$ is flat.
By analogy, we can prove that $F^{-1}[C]$ is a flat in $X$.
\smallskip

Now we check that the function $F$ has the properties (1)--(4).
\smallskip

1. The definition of the relation $F$ ensures that for every $x\in A\cap B$, the pair $(x,x)$ belongs to $F$, which means that $F(x)=x$.
\smallskip

2. If $\dom[F]=A$, then for every point $a\in A$ and its image $b=F(a)\in B$ we have $b\in\Aline oa$ and hence $a\in\Aline oa=\Aline ob\subseteq\overline{\{o\}\cup B}$. On the other hand, assuming that $A\subseteq \overline{\{o\}\cup B}$, we can use the strong regularity of the projective space $X$ (see Theorem~\ref{t:projective<=>} and conclude that for every $a\in A$ there exists a point $b\in B$ such that $a\in\Aline ob$ and hence $b\in\Aline ob=\Aline oa$, witnessing that $b=F(a)$ and $a\in\dom[F]$.
Therefore, $A=\dom[F]$.
\smallskip

3. By analogy we can prove the equivalence $\rng[F]=B\;\Leftrightarrow\;B\subseteq\overline{\{o\}\cup A}$.
\smallskip

4. Finally, assume that $A,B$ are two coplanar lines and $o\in\overline{A\cup B}$. The latter condition and the choice of $o\notin A\cup B$ imply that $A\ne B$ and hence $\|A\cup B\|=3$. By Theorem~\ref{t:projective<=>}, the projective space $X$ is strongly regular and by Corollary~\ref{c:proregular=>ranked}, $X$ is ranked. Then $\overline{\{o\}\cup A}=\overline{A\cup B}=\overline{\{o\}\cup B}$ and hence $A\subseteq \overline{\{o\}\cup B}$ and $B\subseteq\overline{\{o\}\cup A}$. By the preceding two statements, $\dom[F]=A$ and $\rng[F]=B$.
\end{proof}

\section{Central projections in liners}

\begin{definition}\label{d:central-projection} A function $F\subseteq X\times X$ in a liner $X$ is called a \index{central projection}\defterm{central projection with center $o$} if $$F=\{(x,y)\in\dom[F]\times\rng[F]:y\in \Aline ox\}.$$
The point $o$ is called a \index{center}\defterm{center} of the central projection. A central projection $F$ whose domain and range are lines in $X$ is called a \defterm{perspectivity between lines} in $X$ or just a \index{line perspectivity}\defterm{line perspectivity} in $X$.  In this case, the point $o$ is called a \index{perspectivity center}\defterm{perspectivity center} of the line perspectivity $F$.
\end{definition}

Proposition~\ref{p:central-projection} implies the following example.

\begin{example} For every flats $A,B$ in $X$ and point $o\in X\setminus (A\cup B)$, the relation
$$F\defeq\{(x,y)\in A\times B: y\in\Aline xo\}$$ is a central projection with the center $o$. If $A,B$ are two coplanar lines and $o\in\overline{A\cup B}$, then $F$ is a line perspectivity with $\dom[F]=A$ and $\rng[F]=B$.
\end{example} 

\begin{example}\label{ex:identity} For every flat $A\ne X$ in a liner $X$, the identity bijection $1_A:A\to A$ is a central projection (with an arbitrary center $o\in X\setminus A$).
\end{example}

\begin{proposition}\label{p:central-projection2} Let $X$ be a liner and $P$ be a central projection with center $o$ such that the sets  $\dom[P]$ and $\rng[P]$ are flat in $X$. 
\begin{enumerate}
\item The point $o$ belongs to $\dom[P]\cup\rng[P]$ if and only if $\rng[P]=\{o\}$.
\item If $o\notin\dom[P]\cup\rng[P]$, then the function $P$ is injective and $P^{-1}$ is a central projection with center $o$.
\item $P$ is a flat function such that $P(x)=x$ for all $x\in\dom[P]\cap\rng[P]$.
\item If $|\dom[P]|>1$ and $\dom[P]\ne\rng[P]$, then the center $o$ of the central projection $P$ is uniquely determined by $P$.
\end{enumerate}
\end{proposition}

\begin{proof} Let $o$ be a center of the central projection $P$. Then $P=\{(x,y)\in\dom[P]\times\rng[P]:y\in \Aline ox\}$. 
\smallskip

1. Assume that $o\in\dom[P]\cup\rng[P]$. If $o\in\dom[P]$, then for the point $y=P(o)\in\rng[P]$, we have $y\in\Aline oo=\{o\}$ and hence $o=y\in\rng[P]$. Therefore, $o\in\dom[P]\cup\rng[P]$ implies $o\in\rng[P]$.  Assuming that $\rng[P]\ne\{o\}$, we can choose a point  $y\in\rng[P]\setminus\{o\}$ and find a point $x\in \dom[P]$ with $y=P(x)$. The definition of $P$ ensures that $y\in\Aline ox$ and hence $x\in\Aline ox=\Aline oy\subseteq \rng[P]$. For the point $x$ we have $\{o,y\}\in\Aline ox$ and hence $o=P(x)=y$, which contradicts the choice of $y$. This contradiction shows that $\rng[P]=\{o\}$.

On the other hand, if $\rng[P]=\{o\}$, then $o\in\rng[P]\subseteq \dom[P]\cup\rng[P]$.
\smallskip

2. Assume that $o\notin\dom[P]\cup\rng[P]$. To see that the function $P$ is injective, take any pairs $(x,y),(z,y)\in P$. Assuming that $x\ne z$ and taking into account that $y\in\Aline ox\cap\Aline oz$, we conclude that $o\in\Aline ox=\Aline oy=\Aline oz=\Aline xz\subseteq \dom[P]$, which contradicts $o\notin\dom[P]$. Since $o\notin\dom[P]\cup\rng[P]$, for every $x,y\in\dom[P]\cup\rng[P]$, the inclusion $y\in\Aline ox$ is equivalent to $x\in\Aline oy$. Then  $$P^{-1}=\{(y,x)\in \rng[P]\times \dom[P]:x\in\Aline oy\}=\{(y,x)\in\rng[P]\times\dom[P]:y\in\Aline ox\}$$ is a central projection, by Definition~\ref{d:central-projection} 
\smallskip

3. If $o\notin\dom[P]\cup\rng[P]$, then $P=\{(x,y)\in\dom[P]\times[\rng[P]:y\in\Aline ox\}$ is a flat function, by Proposition~\ref{p:central-projection}. If $o\in\dom[P]\cup\rng[P]$, then $\rng[P]=\{o\}$, by the preceding statement. In this case, for every flat $C\subseteq X$, the image $P[C]$ is empty if $C\cap\dom[P]=\varnothing$ and $P[C]=\{o\}$ if $C\cap\dom[P]\ne\varnothing$. In both cases, $P[C]$ is flat in $X$. On the other hand, $P^{-1}[C]=\varnothing$ if $o\notin C$ and $P^{-1}[C]=\dom[P]$ if $o\in C$. In both cases, $P^{-1}[C]$ is flat in $X$, witnessing that the function $P$ is flat.

For every $x\in\dom[P]\cap\rng[P]$, we have $x\in\Aline ox$ and hence $P(x)=x$.
\smallskip

4. Assume that $\dom[P]\ne\rng[P]$ and $|\dom[P]|>1$. To derived a contradiction, assume that $o,o'\in X$ are two distinct points such that $$\{(x,y)\in\dom[P]\times\rng[P]:y\in\Aline ox\}=P=\{(x,y)\in\dom[P]\times\rng[P]:y\in\Aline {o'}x\}.$$ If $\dom[P]\subseteq\rng[P]$, then for every $x\in \dom[P]$, we obtain $x=P(x)$ and hence $\rng[P]=\dom[P]$, we contradicts our assumption This contradiction shows that $\dom[P]\setminus \rng[P]$. Choose any point $x\in\dom[P]\setminus \rng[P]$. Since $|\dom[P]|>1$, there exists a point $x''\in\dom[P]\setminus\{x\}$. Since $\rng[P]$ is flat, the line $\Aline x{x''}$ contains at most one point of the flat $\rng[P]$. Since the projective space $X$ is $3$-long, the line $\Aline x{x''}$ contains some point $x'\in\Aline x{x'}\setminus(\{x\}\cup\rng[P])$. Consider the points $y\defeq P(x)$ and $y'\defeq P(x')$. It follows from $\{x,x'\}\cap\rng[P]=\varnothing$ that  $x\ne y$ and $x'\ne y'$. The choice of the points $o,o'$ ensures that $y\in \Aline ox\cap\Aline {o'}x$ and $y'\in\Aline o{x'}\cap\Aline {o'}{x'}$ and hence $o\ne x$ and $o'\ne x'$. Taking into account that $o\ne o'$, $x\ne x'$, $x\ne y$ and $x'\ne y'$, we conclude that $\Aline xy=\Aline o{o'}=\Aline {x'}{y'}=\Aline x{x'}\subseteq \dom[P]$ and hence $\{o,o'\}\subseteq \dom[P]$.
Then $P(o)\in\Aline oo=\{o\}$ and $P(o')\in\Aline {o'}{o'}$ and hence $\{o,o'\}\subseteq \rng[P]$ and $\Aline o{o'}\subseteq \rng[P]$. The equality  $\Aline o{o'}=\Aline ox$ implies $\{x\}\times\Aline o{o'}\subseteq P$ and hence $P$ is not a function. This is a final contradiction completing the proof.
\end{proof}

\begin{corollary}\label{c:line-perspectivity} Every line perspectivity $P$ in a liner $X$ is a line bijection such that the lines $\dom[P],\rng[P]$ are coplanar and $P(x)=x$ for all $x\in\dom[P]\cap\rng[P]$. If $\dom[P]\ne\rng[P]$, then the center $o$ of the central projection $P$ is uniquely determined by $P$ and belongs to $\overline{\dom[P]\cup\rng[P]}\setminus(\dom[P]\cup\rng[P])$.
\end{corollary}

\begin{proof} Since $P$ is a line perspectivity, $\dom[P]$ and $\rng[P]$ are lines. If $\dom[P]=\rng[P]$, then $P$ is the identity map of the line $\dom[P]=\rng[P]$, by Proposition~\ref{p:central-projection2}(3).
So, assume that $\dom[P]\ne\rng[P]$. By Proposition~\ref{p:central-projection2}(4), the the center $o$ of the perspectivity $P$ is uniquely determined by $P$. Since $\dom[P]\ne\{o\}$, Proposition~\ref{p:central-projection2}(1) ensures that $o\notin\dom[P]\cup\rng[P]$.  Applying Proposition~\ref{p:central-projection2}, we conclude that the function $F$ is bijective and $F(x)=x$ for all $x\in\dom[P]\cap\rng[P]$. It remains to show that the lines $\dom[P],\rng[P]$ are complanar and  $o\in\overline{\dom[P]\cup\rng[P]}$. Choose any point $x\in\dom[P]\setminus\rng[P]$ and consider its image $y\defeq P(x)\in\rng[P]$. The choice of $x\notin\rng[P]$ ensures that $x\ne y$. Then $y\in\Aline ox$ implies $o\in\Aline ox=\Aline xy\subseteq \overline{\dom[P]\cup\rng[P]}$. The definition of the function ensures that $\rng[P]\subseteq\overline{\dom[P]\cup\{o\}}$, so the lines $\dom[P],\rng[P]$ are subsets of the plane $\overline{\dom[P]\cup\{o\}}$.
\end{proof}

By Corollary~\ref{c:line-perspectivity}, every line perspectivity $P$ in a liner $X$ is a line bijection in $X$. Let $\mathcal{I}\projupind_X$ be the smallest submonoid of the inverse symmetric monoid $\I_X$, containing all line perspectivities in $X$. Line bijections that belong to the monoid $\I_X\projupind$ are called \index{line!perspectivity}\index{line perspectivity}\defterm{line projectivities}. Line projectivities can be equivalently defined without mentioning the inverse monoid $\I_X\projupind$.

\begin{definition} A \index{line!projectivity}\index{line projectivity}\defterm{line projectivity} in a liner $X$ is a bijection $P:L\to \Lambda$ between lines $L,\Lambda\subseteq X$ such that $P$ can be written as the composition of finitely many line perspectivities in $X$.
\end{definition}

\begin{exercise} Prove that for every line projectivity $P$ in an affine space $X$, the lines $\dom[P]$ and $\rng[P]$ are parallel.
\end{exercise}

\begin{Exercise} Prove that every line projectivity in a Desarguesian affine space is a line affinity.
\end{Exercise}

\begin{picture}(200,70)(-150,-40)

\put(0,0){\line(1,0){60}}
\put(0,0){\line(2,1){60}}
\put(0,0){\line(2,-1){60}}

\put(30,0){\line(0,1){15}}
\put(30,0){\line(0,-1){15}}
\put(60,0){\line(0,1){30}}
\put(60,0){\line(0,-1){30}}

\put(30,-15){\line(2,3){30}}
\put(30,0){\line(2,1){15}}
\put(60,0){\line(-2,1){15}}

\put(0,0){\color{red}\circle*{3}}
\put(-8,-2){$o$}
\put(30,0){\color{red}\circle*{3}}
\put(23,3){$x$}
\put(60,0){\color{red}\circle*{3}}
\put(62,-3){$z$}
\put(30,15){\circle*{3}}
\put(28,18){$b$}
\put(30,-15){\circle*{3}}
\put(27,-23){$a$}
\put(60,30){\circle*{3}}
\put(62,30){$b'$}
\put(60,-30){\circle*{3}}
\put(63,-33){$a'$}
\put(45,7.5){\circle*{3}}
\put(42,12){$y$}
\end{picture}

\section{Line projectivities in projective spaces}

A triple $xyz$ of points of a liner will be called a {\em colinear triple} if $|\{x,y,z\}|=3$ and $\|\{x,y,z\}\|=2$. So, a colinear triple consists of three colinear pairwise distinct points. 

\begin{theorem}\label{t:projective=>3-transitive} For any colinear triples $x_1x_2x_3$ and $z_1z_2z_3$ in a projective space $X$, there exists a line projectivity $P$ such that $P(x_i)=z_i$ for all $i\in\{1,2,3\}$ and  $P$ is the composition of two or three line perspectivities.  
\end{theorem}

\begin{proof} Consider the lines $L\defeq\overline{\{x_1,x_2,x_3\}}$ and $L'\defeq\overline{\{z_1,z_2,z_3\}}$. Three cases are possible.

If $L\ne L'$, then $|L\cap L'|\le 1$ and there exists a number $i\in\{1,2,3\}$ such that $x_i\notin L'$ and $z_i\notin  L$. Since the projective space $X$ is $3$-long, there exists a point $o\in \Aline {x_i}{z_i}\setminus\{x_i,z_i\}$. The choice of $i$ ensures that $o\notin L\cup L'$. By Proposition~\ref{p:cov-aff}, there exists a line $\Lambda$ in the plane $\overline{L\cup\{o\}}$ such that $z_i\in \Lambda$ and  $L'\ne \Lambda\ne \Aline {z_i}o$. By Proposition~\ref{p:central-projection}, the relation
$$P_1\defeq\{(x,y)\in L\times\Lambda:y\in \Aline ox\}$$ is a line perspectivity with $\dom[P_1]=L$ and $\rng[P_1]=\Lambda$. For every $k\in\{1,2,3\}$, let $y_k\defeq P_1(x_k)$. Since $z_i\in\Lambda\cap\Aline o{x_i}$, the point $y_i$ is equal to $z_i$. Let $j,k$ be any distinct numbers in the set $\{1,2,3\}\setminus\{i\}$. By the $0$-parallelity of the projective space $X$, there exists a unique point $o'\in \Aline {y_j}{z_j}\cap\Aline {y_k}{z_k}$. It follows from $\Lambda\ne L'$ that $o'\notin \Lambda\cup L'$. By Proposition~\ref{p:central-projection}, the relation
$$P_2\defeq\{y,z)\in \Lambda\times L':z\in \Aline {o'}y\}$$ is a line perspectivity such that $\dom[P_2]=\Lambda$ and $\rng[P_2]=L'$. The choice of the point $o'$ guarantees that $P_2(y_j)=z_j$ and $P_2(y_k)=z_k$. Since $y_i=z_i\in \Lambda\cap L'$, $P_2(y_i)=z_i$. Then the composition $P\defeq P_2P_1:L\to L'$ is a required line projectivity such that $P(x_n)=P_2(P_1(x_n))=P_2(y_n)=z_n$ for every $n\in\{1,2,3\}$.
\smallskip

Next, assume that $L=L'$ and $x_i=z_i$ for some $i\in\{1,2,3\}$. In this case, choose any line $\Lambda$ in $X$ such that $\Lambda\cap L=\{x_i\}=\{z_i\}$. Let $j,k$ be distinct numbers in the set $\{1,2,3\}\setminus\{i\}$. Since the projective space $X$ is $3$-long, the set $\Lambda\setminus\{x_i\}$ contains two distinct points $y_j$ and $y_k$. By the $0$-parallelity of the projective space $X$, there exist points $o\in\Aline{x_j}{y_j}\cap\Aline{x_k}{y_k}$ and $o'\in\Aline{y_j}{z_j}\cap\Aline {y_k}{z_k}$. By Proposition~\ref{p:central-projection}, the relations
$$P_1\defeq\{(x,y)\in L\times\Lambda:y\in\Aline ox\}\quad\mbox{and}\quad  P_2\defeq\{(y,z)\in \Lambda\times L':z\in\Aline{o'}y\}$$are line perspectivities such that $\dom[P_1]=L$, $\rng[P_1]=\Lambda=\dom[P_2]$, $\rng[P_2]=L'$, $P_1(x_n)=y_n$ and $P_2(y_n)=z_n$ for all $n\in\{1,2,3\}$. Then the composition $P\defeq P_1:P_2$ is a required line projectivity such that $P(x_n)=P_2(P_1(x_n))=P_2(y_n)=z_n$ for all $n\in\{1,2,3\}$.
\smallskip

Finally, assume that $L=L'$ and $x_i\ne z_i$ for all $i\in\{1,2,3\}$. Choose any point $p\in X\setminus L$.  Since the projective space $X$ is $3$-long, there exist points $o_1\in\Aline {x_1}{p}\setminus\{x_1,p\}$ and $o_3\in\Aline {p}{z_1}\setminus\{p,z_1\}$. Since $x_1\ne z_1$, the lines $\Lambda\defeq\Aline{o_3}{p}$ and $\Lambda'\defeq\Aline{o_1}{p}$ are distinct. Moreover, $o_1\notin L\cup \Lambda$ and $o_3\notin \Lambda'\cup L'$.  
 By Proposition~\ref{p:central-projection}, the relations
$$P_1\defeq\{(x,y)\in L\times\Lambda:y\in\Aline {o_1}x\}\quad\mbox{and}\quad  P_3\defeq\{(y,z)\in \Lambda'\times L':z\in\Aline{o_3}y\}$$are line perspectivities such that $\dom[P_1]=L$, $\rng[P_1]=\Lambda$, $\dom[P_2]=\Lambda'$, $\rng[P_2]=L'$, $P_1(x_1)=p$ and $P_3(p)=z_1$.
For every $i\in\{1,2,3\}$, consider the points $v_i\defeq P_1(x_i)$ and $w_i\defeq P_3^{-1}(z_i)$. Observe that $v_1=p=w_1$. By the $0$-parallelity of the projective space $X$, there exists a point $o_2\in\Aline {v_2}{w_2}\cap\Aline{v_3}{w_3}$.  By Proposition~\ref{p:central-projection}, the relation
$$P_2\defeq\{(v,w)\in \Lambda\times\Lambda':w\in\Aline {o_2}v\}$$is a  line perspectivity such that $\dom[P_2]=\Lambda$, $\rng[P_2]=\Lambda'$, and $P_2(v_i)=(w_i)$ for all $i\in\{1,2,3\}$.
Then the composition $P=P_3P_2P_1$ is a required line projectivity such that $P(x_i)=P_3P_2P_1(x_i)=P_3(P_2(v_i))=P_3(w_i)=z_i$ for all $i\in\{1,2,3\}$.
\end{proof}

\section{Line projectivities in Desarguesian projective spaces}

\begin{theorem}\label{t:Desargues-perspective} A projective space $X$ is Desarguesian if and only if for any concurrent lines $A,B,C$, with $A\ne C$ and line perspectivities $P:A\to B$ and $R:B\to C$, the composition $RP:A\to C$ is a line perspectivity.
\end{theorem} 

\begin{proof} To prove the ``if'' part, assume that the projective space $X$ is not Desarguesian. By Theorem~\ref{t:Desargues-projective}, the liner $X$ has rank $\|X\|=3$. Since $X$ is not Desarguesian, there exist disjoint centrally perspective triangles $abc$ and $a'b'c'$ in $X$ such that the set 
$$T=(\Aline ab\cap\Aline{a'}{b'})\cup(\Aline bc\cap\Aline{b'}{c'})\cup(\Aline ac\cap\Aline {a'}{c'})$$has rank $\|T\|=3$. Since the triangles $abc$ and $a'b'c'$ are disjoint and centrally perspective, the lines $A\defeq \Aline a{a'}$, $B\defeq \Aline b{b'}$  and $C\defeq\Aline c{c'}$ are pairwise distinct and have a unique common point $o$, the perspector of the triangles $abc$ and $a'b'c'$.  It follows from $\|T\|=3$ that the sets $T_a\defeq \Aline bc\cap\Aline {c'}{c'}$, $T_b\defeq \Aline ac\cap\Aline {a'}{c'}$ and $T_c\defeq \Aline ab\cap\Aline{a'}{b'}$ are singletons and hence $T_a=\{t_a\}$, $T_b=\{t_b\}$ and $T_c=\{t_c\}$ for unique points $t_a,t_b,t_c\in X$. Since $\|\{t_a,t_b,t_c\}\|=\|T\|=3$, the intersection $\Aline{t_a}{t_b}\cap\Aline{t_a}{t_c}\cap\Aline{t_b}{t_c}$ is empty. Therefore, the point $o$ does not belong to one of the lines  $\Aline{t_a}{t_b}$, $\Aline{t_a}{t_c}$ or $\Aline{t_b}{t_c}$. We lose no generality assuming that $o\notin\Aline{t_a}{t_c}$.

Taking into account that $a\ne a'$, $b\ne b'$ and $\Aline a{a'}\cap\Aline b{b'}=\{o\}\not\subseteq\{a,a',b,b'\}$, we can show that $T_c=\Aline ab\cap\Aline{a'}{b'}\not\subseteq A\cup B$. By analogy we can show that $T_a=\Aline bc\cap\Aline{b'}{c'}\not\subseteq B\cup C$. By Proposition~\ref{p:central-projection}, the relations
$$P\defeq\{(x,y)\in A\times B:y\in \Aline x{t_c}\}\quad\mbox{and}\quad R\defeq\{(y,z)\in B\times C:z\in\Aline y{t_a}\}$$are lines perspectivities with $\dom[P]=A$, $\rng[P]=B$, $\dom[R]=B$, $\rng[R]=C$. Assuming that the composition $RP$ is a line perspectivity we can find a point $o'\in X\setminus(A\cup C)$ such that $$RP=\{(x,z)\in A\times B:z\in\Aline x{o'}\}.$$Then $o'\in\Aline ac\cap\Aline{a'}{c'}=\{t_b\}$ and hence $o'=t_b$. Since $\|X\|=3=|\{t_a,t_b,t_c\}|$, the lines $\Aline {t_a}{t_c}$ and $A$ have a common point $x$, by Theorem~\ref{t:projective<=>}. Since $o\notin\Aline{t_a}{t_c}$, the point $x\in A\cap\Aline{t_a}{t_c}$ is not equal to the point $o$. Since $x\in\Aline{t_a}{t_c}$, the points $y\defeq P(x)\in\Aline x{t_c}$ and $z\defeq R(y)\in\Aline y{t_a}$ belong to the line $\Aline{t_a}{t_c}$. On the other hand, $z=RP(x)\in\Aline x{o'}$ and hence $t_b=o'\in\Aline xz\subseteq\Aline{t_a}{t_c}$, which contradicts $\|\{t_a,t_b,t_c\}\|=3$. This contradiction shows that the composition $RP$ is not a line perspectivity.
\smallskip

 To prove the ``only if'' part, assume that $X$ is a Desarguesian projective space and take any line perspectivities $P:A\to B$ and $R:B\to C$ between concurrent lines $A,B,C$ with $A\ne C$. We should prove that the line bijection $RP$ is a line perspectivity. If $A=B$, then by Corollary~\ref{c:line-perspectivity}, $P$ is the identity map of the line $A=B$ and then $RP=R$ is a line perspectivity.  If $B=C$, then $R$ is the identity map of the line $B=C$ and $RP=P$ is a line perspectivity. So, assume that $A\ne B\ne C$.

Let $o$ be the unique point of the intersection $A\cap B\cap C$. By Corollary~\ref{c:line-perspectivity}, $P(o)=o=R(o)$. 
Since $P,R$ are line perspectivities with $\dom[R]=A\ne B=\rng[R]$ and $\dom[R]=B\ne C=\rng[R]$, by Corollary~\ref{c:line-perspectivity}, there exist points $p\in \overline{A\cup B}\setminus(A\cup B)$ and $r\in\overline{B\cup C}\setminus(B\cup C)$ such that
$$P=\{(x,y)\in A\times B:y\in\Aline xp\}\quad\mbox{and}\quad R=\{(y,z)\in B\times C:z\in\Aline yr\}.$$

Two cases are possible. First we assume that $p=r$.

In this case, by Proposition~\ref{p:central-projection}(3), $\rng[R]=C$ and $\rng[P]=B$ imply $$C\subseteq\overline{B\cup\{r\}}\subseteq\overline{\overline{A\cup\{p\}}\cup\{r\}}=\overline{A\cup\{p,r\}}=\overline{A\cup p}.$$On the other hand, $\dom[P]=A$ and $\dom[R]=B$ imply 
$$A\subseteq\overline{B\cup\{p\}}\subseteq\overline{\overline{C\cup\{r\}}\cup\{p\}}=
\overline{C\cup\{r,p\}}=\overline{C\cup\{p\}}.$$ By Proposition~\ref{p:central-projection}, the central projection 
$$F\defeq\{(x,z)\in A\times C:z\in\Aline xp\}$$
is a line perspectivity with $\dom[F]=A$ and $\rng[F]=C$. 
We claim that $PR=F$. Given any point $x\in A$ and its image $z\defeq RP(x)$, we should check that $z\in\Aline xp$. If $x=o$, then $z=RP(x)=RP(o)=R(o)=o$ and hence $z=o\in\Aline op=\Aline xp$. So, we assume that $x\ne o$. In this case the points $y\defeq P(x)$ and $z=R(y)$ are not equal to $o=P(o)=RP(o)$, by the injectivity of the functions $P$ and $RP$.
It follows from $y=P(x)$, $z=R(y)$, and $p=r\notin A\cup B$ that $y\in\Aline xp$ and $z\in\Aline yr=\Aline yp=\Aline xp$, witnessing that $F(x)=z=RP(x)$ and $F=RP$. 
 
Next, assume that $p\ne r$. Since the projective space $X$ is $3$-long, there exist distinct points $a,a'\in A\setminus\{o\}$. Consider the points $b\defeq P(a)$, $b'\defeq P(a')$, and $c\defeq R(b)$, $c'\defeq R(b')$. Taking into account that  $P,R$ are line bijections with $P(o)=o=R(o)$, we conclude that $|\{o,b,b'\}|=3=|\{o,c,c'\}|$. 
Assuming that $o\in \Aline ac$, we conclude that $c\in \Aline ac=\Aline ao=A$ and hence $c\in A\cap C=\{o\}=\{RP(o)\}$, which contradicts the injectivity of the function $RP$. By analogy we can prove that $o\notin \Aline {a'}{c'}$.

By the $0$-parallelity of the projective space $X$, there exist a point $q\in \Aline ac\cap\Aline {a'}{c'}$. 
Assuming that $q=a$, we conclude that $a=q\in\Aline {a'}{c'}\setminus\{a'\}$ and hence $c'\in\Aline {a'}{c'}=\Aline a{a'}=A$ and hence $RP(a)=c'\in A\cap C=\{o\}=\{RP(o)\}$, which contradicts the bijectivity of the function $RP$. This contradiction shows that $q\ne a$. By analogy we can show that $q\notin\{a',c,c'\}$. Assuming that $q\in A$, we conclude that $c'\in\Aline {a'}{c'}=\Aline{a'}q\subseteq A$, which contradicts the bijectivity of the function $RP$. By analogy we can show that $q\notin C$. Then $\Aline ac\cap\Aline{a'}{c'}=\{q\}$ and $o\notin\{p,q,r\}$. Since the projective space $X$ is Desarguesian, $\|\{p,q,r\}\|\le 2$ and hence $q\in\Aline pr$. 

By Proposition~\ref{p:central-projection}, the relation
$$F\defeq\{(x,z)\in A\times C:z\in\Aline xq\}$$is a line projectivity with $\dom[F]=A$ and $\rng[F]=C$. We claim that $RP=F$. Given any point $x\in A$ and its image $z=RP(x)$, we should check that $z\in\Aline xq$. If $x=o\in A\cap B\cap C$, then $RP(o)=R(o)=o$, by Corollary~\ref{c:line-perspectivity}. Then $z=o\in\Aline oq$ and we are done. If $x=a$, then $z=RP(a)=R(b)=c$ and $z=c\in\Aline aq=\Aline xq$ by the choice of $q\in\Aline ac\setminus (A\cup C)$. If $x=a'$, then $z=RP(a')=R(b')=c'$ and hence $z=c'\in\Aline{a'}q=\Aline xq$ by the choice of $q\in\Aline {a'}{c'}\setminus (A\cup C)$. 

So, assume that $x\notin\{o,a,a'\}$ and consider the point $y\defeq P(x)\in B\setminus\{o,b,b'\}$. 
If $x\in\Aline pr$, then $y\in\Aline xp$, $z\in \Aline yr$, and $q\in \Aline pr\setminus\{x\}$ imply $y\in\Aline xp=\Aline pr$ and $z\in\Aline yr=\Aline pr=\Aline xr=\Aline xq$, witnessing that $F(x)=z=RP(x)$.

So, assume that $x\notin\Aline pr$. The choice of the points $p,r$ ensures that $p\in\Aline xy\cap\Aline ab\cap\Aline{a'}{b'}$ and $r\in\Aline yz\cap\Aline bc\cap\Aline{b'}{c'}$. By the $0$-parallelity of the projective space $X$, there exist points $s\in \Aline ac\cap\Aline xz$ and $s'\in\Aline{a'}{c'}\cap\Aline xz$. It follows from $o\notin\Aline ac\cup\Aline{a'}{c'}$ that $o\notin\{s,s'\}$ and hence $o\notin\{p,s,s',r\}$. Since the projective space $X$ is Desarguesian, $\|\{p,s,r\}\|\le  2$ and $\|p,s',r\|\le 2$. Since $p\ne r$, $\{s,s'\}\subseteq\Aline pr\cap\Aline xz$. Since $x\notin \Aline pr$, $|\{s,s'\}|\le|\Aline pr\cap\Aline xz|\le 1$ and hence $s=s'\in \Aline ac\cap\Aline{a'}{c'}=\{q\}$. Now $q=s\in\Aline xz\setminus A$ implies $z\in\Aline xq$, witnessing that $F(x)=z=RP(x)$. 
\end{proof} 

\begin{lemma}\label{l:change-center} Let $A,B,C,D$ be concurrent lines in a Desarguesian projective space $X$ such that $A\ne C\ne D$. For any line perspectivities $P:A\to B$ and $P':B\to C$ and any point $o\in X\setminus(A\cup D)$, there exist line perspectivities $R:A\to D$ and $R':D\to C$ such that $P'P=R'R$ and $o$ is the perspectivity center of the line perspectivity $R$.
\end{lemma}

\begin{proof} By Theorem~\ref{t:Desargues-perspective}, the composition $P'P$ is a line perspectivity. Given any point $o\in X\setminus(A\cup D)$, consider the line perspectivity $R:A\to D$  defined by $$R\defeq\{(x,y)\in A\times D:y\in\Aline ox\}.$$ 
 By Proposition~\ref{p:central-projection2}, the function $R^{-1}:D\to A$ is a line perspectivity and by Theorem~\ref{t:Desargues-perspective}, the composition $(P'P)R^{-1}:D\to C$ is a line perspectivity.  Then $R$ and $R'$ are required line perspectivities such that $P'P=P'PR^{-1}R=R'R$ and the perspectivity center of $R$ is $o$.
\end{proof}

\begin{lemma}\label{l:perspectivity-midchange} Let $X$ be a Desarguesian projective space, $L,\Lambda,L'$ be lines with $L\cap\Lambda\cap L'=\varnothing$, and $P:L\to \Lambda$, $P':\Lambda\to L'$ be two line perspectivities. For every line $\Lambda'$ such that $L\cap \Lambda'\ne\varnothing\ne \Lambda'\cap\Lambda\cap L'$ and $L\ne \Lambda'\ne L'$, there exist line perspectivities $R:L\to\Lambda'$ and $R':\Lambda'\to L'$ such that $P'P=R'R$. 
\end{lemma}

\begin{proof} By Corollaries~\ref{c:line-perspectivity} and the $0$-parallelity of the projective space $X$, $L\cap\Lambda\ne\varnothing \ne\Lambda\cap L'$. It follows from $L\cap \Lambda\cap L'=\varnothing$ that $L\ne \Lambda\ne L'$.   

If $\Lambda'=\Lambda$, then the line perspectivities $R\defeq P$ and $R'\defeq P'$ have the required properties. So, assume that $\Lambda'\ne \Lambda$.

By the assumptions of the lemma, there exists a point $\lambda\in \Lambda'\cap\Lambda\cap L'$. Taking into account that $L\cap \Lambda\cap L'=\emptyset$ and $\lambda\in \Lambda\setminus \Lambda'\setminus L$, we conclude that $L\cap\Lambda\cap\Lambda'=\varnothing$ and $\overline{L\cup\{\lambda\}}$ is a plane containing the center $o$ of the line projectivity $P$ and also the line $\Lambda$. Since $\Lambda'\cap L\ne\emptyset$, the plane $\overline{L\cup\{\lambda\}}$ contains also the line $\Lambda'$. It follows from $L\cap\Lambda\cap\Lambda'=\varnothing$ that there exists a point $s\in \overline{L\cup\{\lambda\}}\setminus(L\cup \Lambda\cup\Lambda')$. By Lemma~\ref{l:change-center}, there exist perspectivities $S:L\to\Lambda$ and $S':\Lambda\to L'$ such that $P'P=S'S$ and $s$ is the perspectivity center of the line perspectivity $S$.  

By Proposition~\ref{p:central-projection}, the relations
$$R\defeq\{(x,y)\in L\times\Lambda':y\in\Aline sx\}\quad\mbox{and}\quad T\defeq\{(y,z)\in \Lambda'\times\Lambda:z\in\Aline sy\}$$are line perspectivities such that $\dom[R]=L$,  $\rng[R]=\Lambda'=\dom[T]$, $\rng[T]=\Lambda$, and $S=TR$.  Since the lines $\Lambda',\Lambda,L'$ are concurrent and $\Lambda'\ne L$, Theorem~\ref{t:Desargues-perspective} ensures that the line projectivity $R'\defeq S'\circ T^{-1}:\Lambda'\to L'$ is a line perspectivity. Then $P'P=S'S=S'T^{-1}TS=(S'T^{-1})(TS)=R'R$ and hence $R,R'$ are required line perspectivities with $P'P=R'R$.
\end{proof}

By analogy we can prove the following symmetric version of Lemma~\ref{l:perspectivity-midchange}.

\begin{lemma}\label{l:perspectivity-midchange2} Let $X$ be a Desarguesian projective space, $L,\Lambda,L'$ be lines with $L\cap\Lambda\cap L'=\emptyset$, and $P:L\to \Lambda$, $P':\Lambda\to L'$ be two line perspectivities. For every line $\Lambda'$ such that $\Lambda'\cap L'\ne\varnothing\ne L\cap\Lambda\cap\Lambda'$ and $L\ne \Lambda'\ne L'$, there exist line perspectivities $R:L\to\Lambda'$ and $R':\Lambda'\to L'$ such that $P'P=R'R$. 
\end{lemma}

\begin{lemma}\label{l:perspective-concurrent} Let $P_1,P_2,P_3$ be three line perspectivities in a Desarguesian projective space $X$ such that $\rng[P_1]=\dom[P_2]$ and $\rng[P_2]=\dom[P_3]$. If $\dom[P_1]\cap\dom[P_2]\cap\dom[P_3]\ne\varnothing$ or $\rng[P_1]\cap\rng[P_2]\cap\rng[P_3]\ne\varnothing$, then $P_3P_2P_1=R'R$ for some line perspectivities $R$ and $R'$ in $X$.
\end{lemma}

\begin{proof} Consider the lines $L_1\defeq\dom[P_1]$, $L_2\defeq\rng[P_1]=\dom[P_2]$, $L_3\defeq\rng[P_2]=\dom[P_3]$ and $L_4\defeq\rng[P_3]$. 
If $L_i=L_{i+1}$ for some $i\in\{1,2,3\}$, then $P_i:L_i\to L_{i+1}$ is the identity map of the line $L_i$ and hence $P_3P_2P_1$ is the composition of two line perspectivities. So, we assume that $L_i\ne L_{i+1}$ for every $i\in\{1,2,3\}$. By Corollary~\ref{c:line-perspectivity} the $0$-parallelity of the projective space $X$, for every $i\in\{1,2,3\}$ the lines $L_i$ and $L_{i+1}$ are coplanar and hence concurrent.

By our assumption, $L_1\cap L_2\cap L_3$ or $L_2\cap L_3\cap L_4$ is not empty.
First assume that $L_1\cap L_2\cap L_3\ne\emptyset$. 

If $L_1\ne L_3$, then the composition $P_2P_1:L_1\to L_3$ is a line perspectivity, by Theorem~\ref{t:Desargues-perspective}. In this case the line projectivity $P_3P_2P_1$ is the composition of two line perspectivities $(P_2P_1)$ and $P_3$. 

If $L_1=L_3$, then using Proposition~\ref{p:cov-aff}, we can choose a line $L_2'$ in $X$ such that $L_3\ne L_2'\ne L_4$ and the lines $L_3,L_2',L_4$ are concurrent. Since $L_1=L_3$, the lines $L_1,L_2',L_3$ are concurent. Since $L_1\cap L_2\cap L_3\ne\varnothing$, the lines $L_1,L_2,L_2',L_3$ are concurrent. By Lemma~\ref{l:change-center}, there exist line perspectivities $P'_1:L_1\to L_2'$ and $P'_2:L_2'\to L_3$ such that $P_2P_1=P_2'P_1'$. Then $P_3P_2P_1=P_3P_2'P_1'$. Since the lines $L_2',L_3,L_4$ and concurrent and $L_2'\ne L_4$, the composition $P_3P_2'$ is a line perspectivity, by Theorem~\ref{t:Desargues-perspective}. Then the line projectivity $P_3P_2P_1$ is the composition $(P_3P_2')P_1'$ of two line perspectivities $P_1'$ and $P_3P_2'$.
\smallskip

By analogy we can prove that $P_3P_2P_1$ is the composition of two line perspectivities if $L_2\cap L_3\cap L_4\ne\varnothing$.
\end{proof}

\begin{lemma}\label{l:perspective-last}  Let $P_1,P_2,P_3$ be three line perspectivities in a Desarguesian projective space $X$ such that $\rng[P_1]=\dom[P_2]$ and $\rng[P_2]=\dom[P_3]$. If $\dom[P_1]\ne\rng[P_3]$, then  $P_3P_2P_1=R'R$ for some line perspectivities $R$ and $R'$ in $X$.
\end{lemma}

\begin{proof} Consider the lines $L_1\defeq\dom[P_1]$, $L_2\defeq\rng[P_1]=\dom[P_2]$, $L_3\defeq\rng[P_2]=\dom[P_3]$ and $L_4\defeq\rng[P_3]$. By our assumption, $L_1\ne L_4$.
If $L_i=L_{i+1}$ for some $i\in\{1,2,3\}$, then $P_i:L_i\to L_{i+1}$ is the identity map of the line $L_i$ and hence $P_3P_2P_1$ is the composition of two line perspectivities. So, we assume that $L_i\ne L_{i+1}$ for every $i\in\{1,2,3\}$. By Corollary~\ref{c:line-perspectivity}, for every $i\in\{1,2,3\}$ the lines $L_i$ and $L_{i+1}$ are coplanar and hence concurrent. 

If  $L_1\cap L_2\cap L_3\ne\varnothing$ or $L_2\cap L_3\cap L_4\ne\varnothing$, then by Lemma~\ref{l:perspective-concurrent}, $P_3P_2P_1$ is the composition of two line perspectivities. 
So, assume that $L_1\cap L_2\cap L_3=\varnothing=L_2\cap L_3\cap L_4$. In this case, the singetons $L_1\cap L_2$ and $L_3\cap L_4$ are distinct. Let $\Lambda$ be the unique line containing the singletons $L_1\cap L_2$ and $L_3\cap L_4$. 
Assuming that $\Lambda=L_2$, we conclude that $\emptyset =L_2\cap L_3\cap L_4=\Lambda\cap L_3\cap L_4=L_3\cap L_3\ne\varnothing$, which is a contradiction shows that $L_2\ne\Lambda$. By analogy we can prove that $\Lambda\ne L_3$. 

If $L_3\cap L_4\not\subseteq L_1$, then $L_1\ne \Lambda$.  By Lemma~\ref{l:change-center}, there exist line perspectivities $P'_1:L_1\to \Lambda$ and $P_2':\Lambda\to L_3$ such that $P'_2P'_1=P_2P_1$. Then $P_3P_2P_1=P_3P_2'P_1'$.
If $\Lambda=L_4$, then the lines $\Lambda,L_3,L_4$ are concurrent, and by Lemma~\ref{l:perspective-concurrent}, $P_3P_2'P_1'=R'R$ for some line perspectivities $R,R'$ in $X$. So, assume that $\Lambda\ne L_4$. In this case, the line projectivity $P_3P_2':\Lambda\to L_4$ is a line perspectivity, by Theorem~\ref{t:Desargues-perspective}. Then $P_3P_2P_1$ is the composition of two line perspectivites $P_1'$ and $P_3P_2'$.

By analogy we can prove that $L_1\cap L_2\not\subseteq L_4$ implies that $P_3P_2P_1$ is equal to the composition of two line perspectivities in $X$.

So, assume that $L_1\cap L_2\subseteq L_4$ and $L_2\cap L_3\subseteq L_1$. Then $(L_1\cap L_2)\cup(L_3\cup L_4)\subseteq L_1\cap L_4$ and hence $L_1=L_4$, which contradicts our assumption.
\end{proof} 

\index[person]{Hessenberg}
\begin{theorem}[Hessenberg\footnote{G. Hessenberg, Beweis des Desarguesschen Satzes aus dem Pascalschen. Math. Ann. 61
(1905), 161–172.}, 1905]\label{t:projectivity=3perspectivities} Every line projectivity $P$ in a Desarguesian projective space is the composition of three line perspectivities. Moreover, if $\dom[P]\ne\rng[P]$, then $P$ is the composition of two line perspectivities.
\end{theorem} 

\begin{proof} Given a line projectiving $P$ in a Desarguesian projective space $X$, write $P$ as the composition $P=P_n\cdots P_1$ of line perspectivities $P_1,P_2,\dots,P_n$. We can assume that the number $n$ is the smallest possible. We claim that $n\le 3$. To derive a contradiction, assume that $n\ge 4$.  Since the composition $P_1P_2\cdots P_n$ is a line bijection, $\rng[P_i]=\dom[P_{i+1}]$ for every $i<n$. The minimality of $n$ ensures that $\dom[P_i]\ne\rng[P_i]$ for all $i\le n$. Let $L_0\defeq\dom[P_1]$ and for every $i\in\{1,\dots,n\}$, consider the line $L_i\defeq\rng[P_i]$. The minimality of $n$ and Lemma~\ref{l:perspective-last} ensure that $L_0=L_3$ and $L_1=L_4$. By Proposition~\ref{p:cov-aff}, there exists a line $\Lambda$ such that $L_0\ne \Lambda\ne L_1$ and $L_0\cap \Lambda\cap L_1\ne\varnothing$. By Lemma~\ref{l:change-center}, there exist line perspectivities $R:L_0\to\Lambda$ and $R':\Lambda\to L_1$ such that $P_1=R'R$. Since $\Lambda\ne L_0=L_3$, the composition $P_3P_2R'$ is equal to the composition $T_2T_1$ of two line perspectivities, by Lemma~\ref{l:perspective-last}. Then $P_4P_3P_2R'=P_4T_2T_1$. Since $\dom[T_1]=\dom[R']=\Lambda\ne L_1=L_4=\rng[P_4]$, the composition $P_4T_2T_1$ is equal to the composition $S'S$ of two line perspectivities, by Lemma~\ref{l:perspective-last}. Then $P_4P_3P_2P_1=P_4P_3P_2R'R=P_4T_2T_1R=S'SR$ and $P=P_n\cdots P_5S'SR$ is represented as the composition of $n-1$ line perspectivities, which contradicts the minimality of $n$. This contradiction shows that $n\le 3$ and hence $P$ is the composition of three line projectivities (because the indentity map of a line is a line projectivity). If $\dom[P]\ne\rng[P]$, then Lemma~\ref{l:perspective-last} ensures that $P$ is the composition of two line projectivities. 
\end{proof}

\section{Line projectivities in Pappian projective spaces}

\begin{theorem}[Pappus, $\approx 320$ AD]\label{t:Pappus-perspectivity} A line projectivity $P$ in a Pappian projective space is a line perspectivity if and only if $\dom[P]\cap\rng[P]\ne\varnothing$ and $P(x)=x$ for every $x\in\dom[P]\cap\rng[P]$.
\end{theorem}

\begin{proof} The ``only if'' part follows from the definition of a line perspectivity and Corollaries~\ref{c:line-perspectivity}. To prove the ``if'' part, assume that $P$ is a line projectivity such that $\dom[P]\cap\rng[P]\ne\varnothing$ and $P(x)=x$ for all $x\in\dom[P]\cap\rng[P]$. If $\dom[P]=\rng[P]$, then $P$ is the identity map of the line $\dom[P]=\rng[P]$ and hence $P$ is a line perspectivity, by Example~\ref{ex:identity}. So, assume that $\dom[P]\ne\rng[P]$. By Theorem~\ref{t:projectivity=3perspectivities}, $P=R'R$ for some line perspectivities $R,R'$ in $X$. If the line $\Lambda\defeq\rng[R]=\dom[R']$ is concurrent with the (concurrent) lines $L\defeq \dom[P]=\dom[R]$ and $L'\defeq \rng[P]=\rng[R']$, then the line projectivity $P=R'R$ is a line perspectivity, by Theorem~\ref{t:Desargues-perspective}. 

So, we assume that $L\cap\Lambda\cap L'=\varnothing$. Let $o,b'$ be the  centers of the line perspectivities $R,R'$, respectively. Corollary~\ref{c:line-perspectivity} ensures that $o\notin L\cup\Lambda$ and $b'\notin\Lambda\cup L'$. Let $c,a',b\in X$ be unique points such that $\{c\}=L\cap L'$, $\{a'\}=L\cap\Lambda$, $\{b\}=\Lambda\cap L'$. We claim that $P(c)=c$ implies $o\in\Aline{b'}c$.

Indeed, consider the point $\lambda\defeq R(c)\in\Lambda$ and observe that $$R'(\lambda)=R'R(c)=P(c)=c\in L\cap L'=\dom[P]\cap\rng[P].$$ Then $\lambda\in\Aline oc\setminus\{c\}$ and $\lambda\in\Aline {b'}{c}$, which implies $o\in \Aline oc=\Aline c\lambda=\Aline {b'}{c}$ and $b'\in\Aline oc$ (as $o\notin L$ and hence $o\ne c$) .

By the $0$-parallelity of the projective space $X$, there exists a unique point $c'\in \Aline ob\cap\Aline{b'}{a'}$. Assuming that $c'\in\Lambda$, we conclude that $\{c'\}\in \Lambda\cap\Aline ob\cap\Aline {b'}{a'}=(\Lambda\cap\Aline ob)\cap(\Lambda\cap\Aline {b'}{a'})=\{b\}\cap\{a'\}$ and hence $c'=a'=b\in\Lambda\cap L\cap L'=\varnothing$, which is a contradiction showing that $c'\notin\Lambda$ and hence $a'\ne c'\ne b$. Then $c'\in\Aline ob$ implies $o\in\Aline b{c'}$ and hence
$$o\in\Aline b{c'}\cap\Aline {b'}{c}.$$
 Assuming that $c'\in L$ and taking into account that $c'\in\Aline{b'}{a'}$ and $a'\in L$, we conclude that $b'\in\Aline {a'}{c'}\subseteq L$ and then $o\in \Aline{b'}{c}\subseteq L$, which contradicts the choice of the point $o$ (as the center of the line perspectivity $R:L\to\Lambda$). Assuming that $c'\in L' $ and taking into account that $c'\in\Aline ob$ and $b\in L'$, we conclude that $o\in\Aline b{c'}\subseteq L'$ and then $b'\in \Aline oc\subseteq L'$, which contradicts the choice of the point $b'$ (as the center of the line perspectivity $R':\Lambda\to L'$).

Therefore, $c'\notin L\cup L'$ and we can consider the line perspectivity $P':L\to L'$ with center $c'$. So,
$$P'=\{(x,y)\in L\times L':y\in\Aline{c'}x\}.$$ We claim that the line projectivity $P=R'R$ is equal to the line perspectivity $P'$. Given any point $x\in L$, we should prove that $P(x)=P'(x)$. If $x=c$, then $P'(c)=c=P(c)$ by our assumption. So, we assume that $x\ne c$. Then for the point $a\defeq P'(x)\in L'$, we have $a\in\Aline x{c'}\setminus\{x\}$ and hence $x\in\Aline a{c'}\cap L=\Aline a{c'}\cap\Aline{a'}c$. By the $0$-parallelity of the projective space $X$, there exists a unique point $y\in \Lambda\cap\Aline {a}{b'}$.
Since $R':\Lambda\to L'$ is a line perspectivity with center $b'$, we have $a=R'(y)$.
 
$$
\begin{picture}(200,200)(0,-100)

\put(0,0){\line(2,1){180}}
\put(0,0){\line(2,-1){180}}
{\linethickness{1pt}
\put(0,0){\color{cyan}\line(5,-1){112.5}}
\put(60,30){\color{blue}\line(0,-1){105}}
\put(52,3){$\Lambda$}
\put(60,-30){\color{cyan}\line(1,1){120}}
\put(140,40){$L$}
\put(60,30){\color{teal}\line(1,-1){120}}
\put(140,-47){$L'$}
\put(90,-45){\color{teal}\line(0,1){90}}
\put(60,-75){\color{blue}\line(1,1){52.5}}
\put(90,45){\color{red}\line(-1,-4){30}}
}

\put(0,0){\circle*{3}}
\put(-8,-2){$c'$}
\put(75,-15){\circle*{3}}
\put(69,-11){$x$}
\put(90,45){\circle*{3}}
\put(88,48){$o$}
\put(60,-75){\circle*{3}}
\put(57,-83){$y$}
\put(90,-45){\circle*{3}}
\put(87,-56){$b'$}
\put(90,0){\circle*{3}}
\put(94,-2){$c$}
\put(60,30){\circle*{3}}
\put(58,33){$b$}
\put(60,-30){\circle*{3}}
\put(51,-40){$a'$}
\put(112.5,-22.5){\circle*{3}}
\put(113,-20){$a$}

\end{picture}
$$

By the Pappus Axiom, applied to the points $a,b,c,a',b',c'$, we obtain
$$\|\{y,x,o\}\|=\|(\Aline a{b'}\cap\Aline {a'}b)\cup(\Aline a{c'}\cap\Aline{a'c})\cup(\Aline b{c'}\cap\Aline{b'}c)\|\le 2$$and hence $y\in\Lambda\cap\Aline ox$. Since $R:L\to\Lambda$ is a line perspectivity with center $o$, we conclude that $y=R(x)$ and hence $P(x)=R'R(x)=R'(y)=a=P'(x)$, completing the proof of the equality $P=P'$.
\end{proof}

\begin{corollary}[Pappus, $\approx 320$ AD] Every line projectivity between distinct lines in a Pappian projective space is the composition of two line perspectivities.
\end{corollary}

\begin{proof} Let $P$ be a line projectivity between two distinct lines $\dom[P]$ and $\rng[P]$ in a Pappian projective space $X$. Choose any point $a\in \dom[P]\setminus\rng[P]$ and put  $b\defeq P(a)\in\rng[P]$. Since the projective space $X$ is $3$-long, there exists a point $o\in\Aline ab\setminus\{a,b\}$. By Proposition~\ref{p:cov-aff}, the plane $\overline{\dom[P]\cup\{b\}}$ contains a line $\Lambda$ such that $b\in\Lambda$ and $\Aline ab\ne \Lambda\ne\rng[P]$. By Proposition~\ref{p:central-projection}, the relation $$R\defeq \{(x,y)\in\dom[P]\times\Lambda:y\in\Aline ox\}$$ is a line perspectivity with center $o$ such that $R(a)=b$ and hence  $PR^{-1}:\Lambda\to \rng[P]$ is a line projectivity such that $PR^{-1}(b)=P(a)=b$. By Theorem~\ref{t:Pappus-perspectivity}, $PR^{-1}$ is a line perspectivity and hence $P=(PR^{-1})R$ is the composition of two line perspectivities: $R$ and $PR^{-1}$.
\end{proof}

\begin{corollary}\label{c:Pappus-Fix3} A line projectivity $P$ in a Pappian projective space $X$ is the identity map of the line $\dom[P]=\rng[P]$ if and only if the set $\Fix(P)\defeq\{x\in \dom[P]:P(x)=x\}$ contains at least three points.
\end{corollary}

\begin{proof} The ``only if'' is trivial (because the projective space $X$ is $3$-long). 
To prove the ``if'' part, assume that $|\Fix(P)|\ge 3$, choose three distinct points $a,b,c\in \Fix(P)$ and consider the line $L\defeq\overline{\{a,b,c\}}$. It follows from $\{a,b,c\}\subseteq \dom[P]\cap L\cap\rng[P]$ that $\dom[P]=L=\rng[P]$.  Since $\dim(X)>1$, there exists a line $L'\subseteq X$ such that $L'  \cap L=\{a\}$. Since $X$ is $3$-long, there exist distinct points $a',b',c'\in L'$ such that $a'=a$. By the $0$-parallelity of the projective space $X$, there exists a point $o\in\Aline b{b'}\cap\Aline c{c'}$. It is easy to see that $o\notin L\cup L'$. So, we can consider the line perspectivity
$$R\defeq\{(x,x')\in L\times L':x'\in\Aline ox\}.$$ 
Then $RP:L\to L'$ is a line projectivity such that $RP(a)=a$. By Theorem~\ref{t:Pappus-perspectivity}, the line projectivity $RP$ is a line perspectivity. Since $RP(b)=b'$ and $RP(c)=c'$, the center of the line perspectivity $RP$ is the point $o$ and hence $RP=R$ and $P=R^{-1}R$ is the identity map of the line $L=\dom[P]=\rng[P]$.
\end{proof}

Let us recall that a {\em colinear triple} in a liner $X$ is a triple $xyz\in X^3$ such that $|\{x,y,z\}|=3$ and $\|\{x,y,z\}\|=2$.

\begin{theorem}\label{t:proj-Papp<=>3-transitive} A projective space $X$ is Pappian if and only if for every colinear triples  $xyz$ and $x'y'z'$ in $X$ there exists a unique line projectivity $P$ such that $Pxyz=x'y'z'$.
\end{theorem}

\begin{proof} To prove the ``only if'' part, take any colinear triples $xyz$ and $x'y'z'$ in $X$. By Theorem~\ref{t:projective=>3-transitive}, there exists a line projectivity $P:\overline{\{x,y,z\}}\to\overline{\{x',y',z'\}}$ such that $Pxyz=x'y'z'$. To see that $P$ is unique, assume that $R$ is another line projectivity such that $Rxyz=x'y'z'$. Then $R^{-1}P$ is a line projectivity such that $\{x,y,z\}\subseteq\Fix(R^{-1}P)$. By Corollary~\ref{c:Pappus-Fix3}, $R^{-1}P$ is the identity map $I$ of the line $\overline{\{x,y,z\}}$. Then $P=(RR^{-1})P=R(R^{-1}P)=RI=R$, witnessing that the line projectivity $P$ with $Pxyz=x'y'z'$ is unique.
\smallskip

To prove the ``only if'' part, assume that for every colinear triples $abc,a'b'c'\in X^3$ there exists a unique line  projectovity $P$ in $X$ such that $Pabc=a'b'c'$. To show that the projective space $X$ is Pappian, take any lines $L,L'\subseteq X$ and distinct points $a,b,c\in L\setminus L'$ and $a',b',c'\in L'\setminus L$. By Theorem~\ref{t:projective<=>}, the projective space $X$ is $0$-parallel. Consequently, there exist unique points $x,y\in X$ such that $x\in\Aline b{c'}\cap\Aline {b'}c$ and $y\in\Aline a{c'}\cap\Aline {a'}c$.  It can be shown that the points $x,y$ are distinct and do not belong to the lines $L,L'$. Consider the line $\Lambda\defeq\Aline xy$ and let $o$ be the unique point of the intersection $\Lambda\cap L$. Let $o'$ be the unique point of the intersection $L\cap L'$. Choose using points $\beta\in\Lambda\cap\Aline b{b'}$ and $\gamma\in\Lambda\cap\Aline c{c'}$. Next, find (unique) points $z\in\Lambda\cap\Aline a{b'}$ and let $a''$ be the unique point of the intersection $L'\cap \Aline{b}z$. The Pappus Axiom will follow as soon as we check that $a''=a'$.

\begin{picture}(150,150)(-100,-25)

\put(-15,0){\line(1,0){155}}
\put(145,-4){$L$}
\put(0,0){\line(1,1){110}}
\put(115,110){$L'$}
\put(0,0){\line(-1,-1){10}}
\put(90,0){\color{black}\line(-2,1){60}}
\put(90,0){\color{black}\line(0,1){90}}
\put(120,0){\color{black}\line(-1,1){60}}
\put(120,0){\color{black}\line(-3,1){90}}
\put(60,0){\color{black}\line(0,1){60}}
\put(60,0){\color{black}\line(1,3){30}}
\put(120,0){\line(-1,3){30}}
\put(90,0){\line(-1,2){30}}
{
\put(10,-10){\color{red}\line(2,1){97}}
\put(6,-20){\color{red}$\Lambda$}
}

\put(78,24){\circle*{3}}
\put(78,28){$\beta$}
\put(107.14,38.57){\circle*{3}}
\put(109,39){$\gamma$}

\put(60,0){\circle*{3}}
\put(57,-10){$a$}
\put(90,0){\circle*{3}}
\put(87,-10){$b$}
\put(120,0){\circle*{3}}
\put(118,-10){$c$}

\put(30,30){\circle*{3}}
\put(21,32){$a''$}
\put(36,30){$a'$}
\put(60,60){\circle*{3}}
\put(65,57){$b'$}
\put(90,90){\circle*{3}}
\put(93,86){$c'$}

\put(60,15){\color{black}\circle*{3}}
\put(54,7){\color{black}$z$}
\put(66,18){\color{black}\circle*{3}}
\put(62,23){\color{black}$y$}
\put(90,30){\color{black}\circle*{3}}
\put(95,27){\color{black}$x$}
  
\put(0,0){\circle*{3}}
\put(-3,3){$o'$}
\put(30,0){\circle*{3}}
\put(27,-10){$o$}
\end{picture}

Consider the line perspectivities 
$$
\begin{aligned}
P_{b'}&\defeq\{(u,v)\in L\times \Lambda:v\in\Aline u{b'}\},&P_{c'}&\defeq\{(u,v)\in L\times\Lambda:v\in\Aline u{c'}\},\\
P_b&\defeq\{(v,w)\in \Lambda\times L':w\in\Aline vb\},
&P_c&\defeq\{(v,w)\in\Lambda\times L':w\in\Aline vc\}
\end{aligned}
$$
and observe that $P_{b}P_{b'}obc=P_bo\beta x=o'b'c'=P_cox\gamma=P_cP_{c'}obc$. 
Our assumption guarantees that $P_bP_{b'}=P_cP_{c'}$ and hence
$$a'=P_c(y)=P_cP_{c'}(a)=P_bP_{b'}(a)=P_b(z)=a'',$$
witnessing that $X$ satisfies the Pappus Axiom and hence is Pappian.
\end{proof}

\section{Line projectivities in non-Pappian finite projective planes}

Given any line $L$ in a projective liner $X$, let $\Sym\projupind_X(L)$ be the group of all line projectivities $P:L\to L$. 
 For Pappian liners, the structure of the (elements of the) group $\Sym\projupind_X(L)$ will be studied in Section~\ref{s:projchart}.

Now we describe the structure of the group $\Sym\projupind_X(L)$ for lines $L$ in   non-Pappian finite projective planes $X$. For lines $L$ of cardinality $|L|\ne 24$, the  following theorem was proved by Grundh\"ofer \cite{Grundhofer1988}. The case of $|L|=24$ was resolved by M\"uller M\"uller and Nagy in \cite{MN2007}.

\begin{theorems}[Grundh\"ofer, 1988; M\"uller--Nagy, 2007]\label{t:Grundhofer-Muller-Nagy-p} For every line $L$ in a non-Pappian finite projective plane $X$, the group of line projectivities $\Sym\projupind_X(L)$ equals $\Alt(L)$ or $\Sym(L)$.
\end{theorems}

\begin{proof} The proof is a bit complicated and refers to some results that will be proved (or merely cited) in the next chapters of the book. Since the projective plane $X$ is not Pappian, $|L|=|X|_2\ge 7$, by Corollary~\ref{c:p5-Pappian}. By Theorem~\ref{t:Moufang-finite<=>}, the finite non-Pappian projective plane $X$ is not Moufang. Fix any point $h\in L$. Since the projective plane $X$ is not Moufang, we can apply Theorem~\ref{t:Skornyakov-San-Soucie} (or Theorem~\ref{t:proj-Moufang<=>}) and find a line $H$ in $X$ such that $H\cap L=\{h\}$ and the affine plane $\Pi\defeq X\setminus H$ is not Thalesian.
By Theorem~\ref{t:non-Thales-line-aff}, for the line $L'\defeq L\setminus\{h\}$ in the non-Thalesian Playfair plane $\Pi$, the permutation group $\Sym^{\Join}_\Pi(L')$ either contains the alternating group $\Alt(L')$ or is isomorphic to the Mathieu group $M_{24}$. In the latter case, $|L'|=24$.

 By Lemma~\ref{l:line-affinity-extends-to-line-projectivity}, for every line affinity $A:L'\to L'$ of the line $L'$ in the Playfair plane $\Pi$, the function $\bar A\defeq A\cup\{(h,h)\}$ is a line projectivity on the line $L$ in the projective space $X$. The correspondence $A\mapsto\bar A$ is an isomorphic embedding of the group of line affinities $\Sym^{\Join}_X(L')$ into the (stabilizer of the point $h$ in the) group $\Sym\projupind_X(L)$ of line projectivities on the line $L$. 
 
 By Theorem~\ref{t:projective=>3-transitive}, the permutation group $\Sym\projupind_X(L)$ is $3$-transitive. Since, the $3$-transitive group $\Sym\projupind_X(L)$ contains an isomorphic copy of the group $\Sym^{\Join}_\Pi(L')$ which is equal to $\Alt(L')$, $\Sym(L')$ or $M_{24}$, we can apply Theorem$^\dag$~\ref{t:3-transitive-classification} and conclude that the minimal normal subgroup $N$ of the group $\Sym\projupind_X(L)$ equals $\Alt(L)$. By Theorem~\ref{p:trivial-center=>embedsinAut}, the  conjugating map $\Sym\projupind_X(L)\to\Aut(N)$ is injective.
Since $|L|\ne 6$, the automorphism group $\Aut(N)$ of the group $N\cong\Alt(L)$ can be identified with the symmetric group $\Sym(L)$, see H\"older's Theorem~\ref{t:HolderAn}. Therefore, the permutation group $\Sym\projupind_X(L)$ contains the group $N=\Alt(L)$ and canonically embeds into the group $\Aut(N)=\Sym(L)$, which implies that $\Sym\projupind_X(L)$ equals $\Alt(L)$ or $\Sym(L)$.
\end{proof}

\chapter{Proportions and proscalars in projective spaces}

\rightline{\em In more recent times one has … to make the geometry of position} 

\rightline{\em  an independent science that does not require measuring.}

\rightline{\index[person]{von Staudt}Karl von Staudt\footnote{{\bf Karl Georg Christian von Staudt} (1798--1867) was a German mathematician whose work laid the foundations of modern projective geometry. Born in Rothenburg ob der Tauber and educated at the University of G\"ottingen, where he attended Gauss’s lectures, he earned his doctorate in 1822 with a dissertation on number theory. After several years of teaching, von Staudt was appointed professor of mathematics at the University of Erlangen in 1835, a position he held until his death in 1867. His major work, ``Geometrie der Lage'' (1847), together with the later ``Beitr\"age zur Geometrie der Lage'' (1856--1860), established for the first time a purely synthetic foundation of projective geometry, entirely independent of metric and analytic notions. In these writings von Staudt introduced the concept of the Wurf (throw), the geometric counterpart of the modern cross ratio, and showed that arithmetic operations can be constructed within geometry itself by means of purely projective (incidence) constructions such as harmonic division. In this way he reversed the traditional analytic order of dependence, deriving algebra and coordinates from geometric principles rather than the reverse. Von Staudt’s work also contained the first synthetic formulation of the fundamental theorem of projective geometry and anticipated the inclusion of imaginary elements in geometric reasoning.  Though little appreciated during his lifetime---his style was extremely austere and abstract---his ideas were later recognized by Klein, Pasch, and Hilbert as decisive for the logical autonomy and axiomatic development of geometry. Von Staudt thus stands as one of the key figures in the nineteenth-century transformation of geometry from a metric to a purely structural, projective discipline.
}, ``Geometrie der Lage'', 1847}
\vskip30pt

By analogy with affine spaces 
$X$, which possess a canonical field of scalars 
$\IR_X$, projective spaces carry a canonical field of proscalars 
$\IF_X$, defined and studied in this section following and developing the pioneering ideas of Karl von Staudt.

\section{Projective equivalence of line quadruples}

Proscalars in projective spaces will be defined as equivalence classes of some special line quarduples.

\begin{definition} A \index{line quadruple}\defterm{line quadruple} in a projective space $X$ is a quadruple $zvui\in X^4$ such that  $\|\{z,v,u,i\}\|=2$ and $|\{z,u,i\}|=3$.
\end{definition}

For a projective space $X$, by \index[note]{$\ddddot X$}\defterm{$\ddddot X$} we denote the set of line quadruples in $X$.

In a line quadruple $zvui$, the points $z,v,u,i$ play non-symmetric roles. The notations $z,v,u,i$ are abbreviations of {\em zero}, {\em  variable}, {\em unit}, {\em infinity}, respectively.

\begin{picture}(100,35)(-100,-15)

\put(0,0){\line(1,0){200}}

\put(30,0){\circle*{3}}
\put(28,3){$0$}
\put(27,-8){$z$}
\put(60,0){\color{red}\circle*{3}}
\put(57,-8){$v$}
\put(80,0){\circle*{3}}
\put(77,-8){$u$}
\put(78,3){$1$}
\put(180,0){\circle*{3}}
\put(175,5){$\infty$}
\put(178,-10){$i$}
\end{picture}

Observe that every quadruple $zvui\in X^4$ in a set $X$ is a function $zvui:4\to X$. So we can consider its composition $Fzvui$ with any other function $F$. This composition is a quadruple if and only if $\{z,v,u,i\}\subseteq\dom[F]$.

\begin{exercise} Show that every line quadruple $zvui$ in a Steiner projective space $X$ is equal to $zzui$, $zuui$, or $ziui$.
\end{exercise}

\begin{definition} Given two line quadruples $zvui,oxej\in \ddddot X$ in a projective space $X$, we write $zvui\projeq oxej$ and say that the line quadruples $zvui,oxej$ are \index{projectively equivalent line quadruples}\defterm{projectively equivalent} if there exists a line projectivity $P$ in $X$ such that $Pzvui=oxej$.
\end{definition}

The following proposition shows that the projective equivalence is indeed an equivalence on the set of line quadruples in a projective space.

\begin{proposition}\label{p:projeq} Let $zvui,\hat z\hat v \hat u\hat \imath,\check z\check v\check u\check \imath$ be three line quadruples in a projective space $X$. Then
\begin{enumerate}
\item $zvui\projeq zvui$;
\item $\hat z\hat v\hat u\hat \imath\projeq \check z\check v\check u\check \imath\;\Ra\;
 \check z\check v\check u\check \imath\;\projeq\;\hat z\hat v\hat u\hat \imath$;
\item $(\,\hat z\hat v\hat u\hat \imath\projeq zvui\;\wedge\; zvui\projeq\check z\check v\check u\check \imath\,)\;\Ra\;
 \check z\check v\check u\check \imath\projeq\hat z\hat v\hat u\hat \imath$.
\end{enumerate}
\end{proposition}

\begin{proof} 1. Since $zvui$ is a line quadruple, $L\defeq\overline{\{z,v,u,i\}}$ is a line in $X$. By Exercise~\ref{ex:identity}, the identity map $I:L\to L$ is a line projectivity witnessing that $zvui\projeq zvui$.
\smallskip

2. If $\hat z\hat v\hat u\hat \imath\projeq \check z\check v\check u\check \imath$, then there exists a line projectivity $P$ in $X$ such that $P\hat z\hat v\hat u\hat \imath=\check z\check v\check u\check \imath$.
Write $P$ as the composition $P=P_n\cdots P_1$ of line perspectivities $P_1,\dots,P_n$ in $X$. By Proposition~\ref{p:central-projection2}, for every $i\in\{1,\dots,n\}$, the inverse relation $P_i^{-1}$ is a line projectivity. Then $P^{-1}=(P_n\cdots P_1)^{-1}=P_1^{-1}\cdots P_n^{-1}$ is the composition of line perspectivities and hence $P^{-1}$ is a line projectivity. It follows from $P\hat z\hat v\hat u\hat \imath=\check z\check v\check u\check \imath$ that $\hat z\hat v\hat u\hat \imath=P^{-1}\check z\check v\check u\check \imath$, witnessing that  $\check z\check v\check u\check \imath\;\projeq\;\hat z\hat v\hat u\hat \imath$.
\smallskip

3. If $\hat z\hat v\hat u\hat \imath\projeq zvui$ and $zvui\projeq\check z\check v\check u\check \imath$, then there exist line projectivities $P,P'$ in $X$ such that $P\hat z\hat v\hat u\hat \imath=zvui$ and $P'zvui\projeq\check z\check v\check u\check \imath$. It follows that $\rng[P]=\Aline zu=\dom[P']$ and hence the composition $P'P$ is a line projectivity with $P'P \check z\check v\check u\check \imath=\hat z\hat v\hat u\hat \imath$, witnessing that  $\check z\check v\check u\check \imath\projeq\hat z\hat v\hat u\hat \imath$.
\end{proof}

Theorems~\ref{t:projective=>3-transitive} and \ref{t:proj-Papp<=>3-transitive} imply the following corollary on the existence of projectively equivalent line quadruples.

\begin{corollary}\label{c:proj4-exist} For every line quadruple $zvui\in \ddddot X$ in a (Pappian) projective space $X$ and any colinear triple $oej\in X^3$ there exists a (unique) point $x\in  X$ such that $oxej\projeq zvui$.
\end{corollary}


\section{Pappian quadruples in projective spaces}

Definition~\ref{d:Desarg-triple} of a Desarguesian triple and Corollary~\ref{c:proj4-exist} motivate the following definition of a Pappian quadruple.

\begin{definition} A line quadruple $zvui\in \ddddot X$ in a projective space $X$ is called \index{Pappian line quadruple}\index{line quadruple!Pappian}\defterm{Pappian} if for every colinear triple $oej\in X^3$ there exists a unique point $x\in X$ such that $oxej\projeq zvui$.
\end{definition}

\begin{proposition}\label{p:Pappian-quadruple<=>} For a line quadruple $zvui\in\ddddot X$ in a projective space $X$, the following statements are equivalent:
\begin{enumerate}
\item the line quadruple $zvui$ is Pappian;
\item there exists a unique point $x\in X$ such that $zxui\projeq zvui$;
\item for every line projectivity $P$ with $Pzui=zui$ we have $P(v)=v$.
\end{enumerate}
\end{proposition}

\begin{proof} The implication $(1)\Ra(2)$ follows from the definition of a Pappian line quadruple, and $(2)\Ra(3)$ is trivial. 
\smallskip

$(3)\Ra(1)$ If the line quadruple $zvui$ is not Pappian, then there exist a colinear triple $oej$ and two distinct points $x,y\in X$ such that $oxej\projeq zvui\projeq oyej$. Since  $oxej\projeq zvui\projeq oyej$, there exist line projectivities $P,P'$ in $X$ such that $Poxej=zvui$ and $P'zvui=oyej$. Observe that $\dom[P]=\Aline oe=\rng[P']$ and hence $PP'$ is a line perspectivity such that $PP'zui=Poej=zui$. Since $x\ne y$, the point $w\defeq P(y)$ is distinct from the point $P(x)=v$. Taking into account that  $PP'zvui=zwui\ne zvui$, we conclude that the line quadruple $zvui$ is not Pappian.
\end{proof}

\begin{corollary} For every colinear triple $zui\in X^3$ in a projective space $X$ and every point $v\in\{z,u,i\}$, the line quadruple $zvui$ is Pappian.
\end{corollary}

\begin{exercise} Show that every line quadruple $zvui$ in a Steiner projective space is Pappian.
\end{exercise}

Proposition~\ref{p:Pappian-quadruple<=>} and Theorem~\ref{t:proj-Papp<=>3-transitive}  imply the following characterization of Pappian projective spaces via Pappian quadruples.

\begin{corollary}\label{c:Papp<=>allPapp4} A projective space $X$ is Pappian if and only if every line quadruple in $X$ is Pappian.
\end{corollary}

\begin{exercise} Let $zvui$ and $oxej$ be two projectively equivalent line quadruples in a projective space $X$. Show that the line quadruple $zvui$ is Pappian if and only if  $oxej$ is Pappian.
\end{exercise}

\begin{proposition}\label{p:Pappian-quadruples} If  $zvui$ is a Pappian quadruple in a projective space $X$, then
\begin{enumerate}
\item the line quadruple $uvzi$ is Pappian.
\item If $v\notin\{z,i\}$, then the line quadruple $zuvi$ is Pappian.
\end{enumerate}
\end{proposition}

\begin{proof} Assume that $zvui$ is a Pappian line quadruple in $X$.
\smallskip

1. To prove that the line quadruple $uvzi$ is Pappian, take any line projectivity $P$ with $Puzi=uzi$. We need to show that $P(v)=v$. The equality $Puzi=uzi$ implies $Pzui=zui$. Since the line quadruple $zvui$ is Pappian, the equality $Pzui=zui$ implies $P(v)=v$.
\smallskip

2. Assume that $v\notin\{z,i\}$. In this case $zuvi$ is a line quadruple. To prove that this line quadrpuple is Pappian, take any line projectivity $P$ with $Pzvi=zvi$. We have to prove that $P(u)=u$.  Consider the point $u'\defeq P(u)$. It follows that $3=|\{z,u,i\}|=|\{P(z),P(u),P(i)\}|=|\{z,u',i\}|$. Since $\|X\|\ge 3$, there exists a line $L\subseteq X$ such that $L\cap\Aline zi=\{z\}$. Since the projective space $X$  is $3$-long, there exists a point $c\in X\setminus(\Aline zi\cup L)$. Consider the line perspectivity $$C\defeq\{(x,y)\in \Aline zi\times L:y\in\Aline cx\}$$and the triple $jww'\defeq Ciuu'$ on the line $L$.

\begin{picture}(200,155)(-150,-15)

\put(0,0){\line(1,2){60}}
\put(0,0){\line(-1,2){60}}
\put(180,120){\line(-1,0){240}}
\put(180,120){\line(-6,-1){222}}
\put(180,120){\line(-5,-2){200}}

\put(-52,78){$w$}
\put(-41.5,83){\circle*{3}}
\put(-42,86){$w'$}
\put(48,90){$u'$}
\put(49,98){\circle*{3}}
\put(44,101){$u$}
\put(181,113){$c'$}
\put(-60,120){\circle*{3}}
\put(-62,125){$j$}
\put(60,120){\circle*{3}}
\put(59,123){$i$}
\put(180,120){\circle*{3}}
\put(175,123){$c$}
\put(0,0){\circle*{3}}
\put(-3,-8){$z$}
\put(30,60){\circle*{3}}
\put(32,55){$v$}
\put(-20,40){\circle*{3}}
\put(-17,35){$v'$}
\end{picture}

Assuming that $z\in \Aline ij\subseteq \Aline ic$, we conclude that $c\in\Aline zi$, which contradicts the choice of the point $c$. This contradiction shows that $z\notin\Aline ij=\Aline ic$ and hence $\Aline zi\cap\Aline ij=\{i\}$ and $u\notin \Aline ij$. Since the liner $X$ is projective, there exists a unique point $c'\in \Aline u{w'}\cap\Aline ij$. Assuming that $c'\in \Aline zi$, we conclude that $c'\in\Aline zi\cap \Aline ij=\{i\}$ and hence $w'\in\Aline u{c'}=\Aline ui$ and $c\in \Aline {u'}{w'}\subseteq \Aline ui=\Aline zi$, which contradicts the choice of the point $c$. This contradiction shows that $c'\notin \Aline zi$. 
Assuming that $c'\in L=\Aline zj$, we conclude that $c'\in L\cap\Aline u{w'}\cap \Aline ij=(L\cap\Aline u{w'})\cap(L\cap\Aline ij)=\{w'\}\cap\{j\}$ and hence $w'=c'=j$ and $u'=C^{-1}(j)=i$. Then $P(u)=u'=i=P(i)$ implies $u=i$, which contradicts the  choice of the line quadruple $zvui$. This contradiction shows that $c'\notin L$. Then we can consider the line perspectivity $C'\defeq\{(x,y)\in L\times\Aline zi:y\in \Aline x{c'}\}$, and observe that $C'zw'j=zui$. Consider the line projectivity $C'CP$ and observe that $C'CPzui=C'Czu'i=C'zw'j=zui$. Since the line quadruple $zvui$ is a Pappian, $C'CPzui=zui$ implies $C'CP(v)=v=P(v)$ and hene $C'C(v)=v$. Consider the point $v'\defeq C(v)\ne C(z)=z$ and observe that $C'(v')=C'C(v)=v$ and hence $c,c'\in\Aline v{v'}$. Assuming that $c\ne c'$, we conclude that $\Aline v{v'}=\Aline c{c'}=\Aline ij$ and then $v\in \Aline v{v'}\cap\Aline zi=\Aline ij\cap\Aline zi=\{i\}$, which contradicts our assumption. This contradiction shows that $c'=c$ and then $C'C$ is the identity map of the line $\Aline zi$. In particular, $C'(w')=u=C'C(u)=C'(w)$ and hence $C(u')=w'=w=C(u)$ and finally $P(u)=u'=u$, witnessing that the line quadruple $zuvi$ is Pappian.
\end{proof}

\begin{corollary}\label{c:Pappian-permutations} Let $zvui$ be a Pappian line quadruple in a projective space $X$ such that $|\{z,v,u,i\}|=4$. Then for every points $a,b,c,d\in X$ with $\{a,b,c,d\}=\{z,v,u,i\}$, the line quadruple $abcd$ is Pappian.
\end{corollary}

\begin{proof} First we prove the following partial case of this corollary.

\begin{claim}\label{cl:Pappian-permutations} For any Pappian line quadruple $zvui\in X^4$ with $|\{z,v,u,i\}|=4$ and any points $a,b,c,d\in X$ with $\{a,c,d\}=\{z,u,i\}$ and $b=v$, the line quadruple $abcd$ is Pappian.
\end{claim}

\begin{proof} Given any line projectivity $P$ with $Pacd=acd$, we have to prove that $P(b)=b$. Since $\{a,c,d\}=\{z,u,i\}$, the equality $Pacd=acd$ implies $Pzui=zui$. Since $zvui$ is Pappian, $Pzui=zui$ implies $P(b)=P(v)=v=b$. 
\end{proof}

By Claim~\ref{cl:Pappian-permutations}, the line quadruples $uvzi$ and $zviu$ are Pappian. By Proposition~\ref{p:Pappian-quadruples}(2) applied to the Pappian line quadruples $zvui,uvzi,zviu$, the line quadruples $zuvi$, $uzvi$ and $zivu$ are Pappian. 

Now we are able to prove that for every points $a,b,c,d\in X$ with $\{a,b,c,d\}=\{z,v,u,i\}$ the line quadruple $abcd$ is Pappian. If $b=v$, then $abcd$ is Pappian by Claim~\ref{cl:Pappian-permutations}. If $b=z$, then $abcd$ is Pappian by Claim~\ref{cl:Pappian-permutations} applied to the Pappian line quadruple $uzvi$. If $b=u$ or $b=i$, then $abcd$ is Pappian by Claim~\ref{cl:Pappian-permutations} applied to the Pappian line quadruple $zuvi$ and $zivu$.
\end{proof}

\begin{lemma}\label{l:Pappian-multiplication} Let $X$ be a projective space and $z,a,b,u,i,o,x,y,e,j$ be points in $X$ such that $zabi,oxyj$ are projectively equivalent line quadruples and $zbui,oyej$ are projectively equivalent line quadruples.  If one of the line quadruples $zabi$ or $zbui$ is Pappian, then $zaui$ and $oxej$ are projectively equivalent line quadruples.
\end{lemma}

\begin{proof} Taking into account that $zabi$ and $zbui$ are line quadruples, we conclude that $$|\{z,b,i\}|=3=|\{z,u,i\}|\quad\mbox{and}\quad\|\{z,a,b,i\}\|=2=\|\{z,b,u,i\}\|,$$ which implies that $zaui$ is a line quadruple. By analogy we can show that $oxej$ is a line quadruple. 
Since $zabi\projeq oxyj$ and $zbui\projeq oyej$, there exist line projectivities $P$ and $Q$ such that $Pzabi=oxyj$ and $Qzbui=oyej$. Then $Q^{-1}P$ is a line projectivity with $Q^{-1}Pzbi=zbi$. 

If the line quadruple $zabi$ is Pappian, then the equality $Q^{-1}Pzbi=zbi$ implies $Q^{-1}P(a)=a$ and hence $x=P(a)=Q(a)$. Then $Qzaui=oxej$, witnessing that $zaui\projeq oxej$.

If the line quadruple $zbui$ is Pappian, then so is the line quadruple $zubi$, see Proposition~\ref{p:Pappian-quadruples}(2). In this case, the equality $Q^{-1}Pzbi=zbi$ implies $Q^{-1}P(u)=u$ and hence $e=Q(u)=P(u)$. Then $Pzaui=oxej$, witnessing that  $zaui\projeq oxej$.
\end{proof}

%

\section{Proportions and proscalars in projective spaces}

By Proposition~\ref{p:projeq}, the relation of projective equivalence $\projeq$ is an equivalence relation on the set $\ddddot X$ of all line quadruples in $X$. So, we can consider the quotient set\index[note]{$\ddddot X\projeqind$} $\ddddot X\projeqind$ of $\ddddot X$ by the equivalence relation $\projeq$. Elements of the set $\ddddot X\projeqind$ are called \index{proportion}\defterm{proportions}\footnote{Abbreviated from ``projective portions''.} in $X$.

\begin{definition} For a line quadruple $zvui\in \ddddot X$ in a projective space $X$, its \index{proportion}\defterm{proportion} is the set $$\overvector{zvui}\defeq\{xyzj\in \ddddot X:xyzj\projeq zvui\}\in \ddddot X\projeqind.$$
\end{definition}

Corollary~\ref{c:proj4-exist} implies the following existence result.

\begin{corollary}\label{c:Pappian-exist-unique} Let $X$ be a (Pappian) projective space. For every colinear triple $zui$ in $X$ and every proportion $\alpha\in \ddddot X\projeqind$ there esists a (unique) point $v\in X$ such that $zvui\in\alpha$.
\end{corollary}

\begin{definition} Let $X$ be a projective space. A proportion $\alpha\in \ddddot X\projeqind$ is called a \index{proscalar}\defterm{proscalar} if $\alpha=\overvector{zvui}$ for some (equivalently, every) Pappian line quadruple $zvui\in \ddddot X$. The set of all proscalars in a projective space $X$ is denoted by \index[note]{$\bar\IF_X$}$\bar\IF_X$.
\end{definition} 

The set of proscalars $\bar\IF_X$ is a subset of the set $\ddddot X\projeqind$ of proportions. By Corollary~\ref{c:Papp<=>allPapp4}, a projective space if Pappian if and only if $\bar\IF_X=\ddddot X\projeqind$.

\begin{example} For every projective space $X$, the sets
$$0\defeq\{zvui\in \ddddot X:v=z\},\quad 1\defeq\{zvui\in\ddddot X:v=u\}, \quad\mbox{and}\quad \infty\defeq\{zvui\in \ddddot X:v=i\}$$
are proscalars and hence $0$, $1$, $\infty$ are elements of the set of $\bar \IF_X\subseteq\ddddot X\projeqind$.
\end{example}

\begin{proposition}\label{p:Steiner=>FX=01} Every Steiner projective space $X$ has $\bar\IF_X=\ddddot X\projeqind=\{0,1,\infty\}$.
\end{proposition}

\begin{proof} Fix any colinear triple $zui$ in $X$. Given any proportion $\alpha\in\ddddot X\projeqind$, find a point $v\in\overline{\{z,u,i\}}$ such that $\overvector{zvui}=\alpha$. Since the liner $X$ is Steiner, $v\in \overline{\{z,u,i\}}=\{z,u,i\}$. Then $\alpha=\overvector{zvui}\in\{\overvector{zzui},\overvector{zuui},\overvector{zuii}\}=\{0,1,\infty\}\subseteq\bar\IF_X$ and hence $\ddddot X\projeqind=\{0,1,\infty\}=\bar\IF_X$.
\end{proof}

The set of proscalars $\bar\IF_X$ of a projective space $X$ contains two important subsets:\index[note]{$\IF_X$}\index[note]{$\bar\IF_X^*$}
$$
\IF_X\defeq\bar\IF_X\setminus\{\infty\},\quad\mbox{and}\quad \IF_X^*\defeq\IF_X\setminus\{0\}=\bar\IF_X\setminus\{0,\infty\}.$$
The sets $\IF_X$ and $\IF_X^*$ consist of finite proscalars and finite nonzero proscalars,  respectively. 

The set $\IF_X^*$ of nonzero finite proscalars is a subset of the set\index[note]{$\ddddot X\projeqind^*$}
$$\ddddot X\projeqind^{*}\defeq\ddddot X\projeqind\setminus\{0,\infty\}$$of all nonzero finite proportions.

\section{The inversion and $01$-involution of proportions}

In this section we consider two natural unary operations on the set of proportions of a projective space. Those two operations are induced by suitable permutations of line quadruples. 

\begin{proposition}\label{p:Pappian-inverse} For every projective space $X$, there exists a unique unary operation $$(\cdot)^{-1}:\ddddot X\projeqind\to\ddddot X\projeqind,\quad (\cdot)^{-1}:p\mapsto p^{-1},$$such that 
$0^{-1}=\infty$, $\infty^{-1}=0$, and $(\overvector{zvui})^{-1}=\overvector{zuvi}$ for every line quadruple $zvui\in \ddddot X$ with $v\notin\{z,i\}$.
\end{proposition}

\begin{proof} Given any proportion $\alpha\in \ddddot X\projeqind$, define the proportion $\alpha^{-1}$ as follows. If $\alpha=0$, then put $\alpha^{-1}\defeq\infty$. If $\alpha=\infty$, then put $\alpha^{-1}\defeq 0$. If $\alpha\notin\{0,\infty\}$, then choose any line quadruple $zvui\in\alpha$. It follows from $\overvector{zvui}=\alpha\notin\{0,\infty\}$ that $v\notin\{z,i\}$ and hence $zuvi$ is a line quadruple. So, we can consider its proportion $\overvector{zuvi}$ and put $\alpha^{-1}\defeq\overvector{zuvi}$. Let us show that $\alpha^{-1}$ does not depend on the choice of the line quadruple $zvui\in\alpha$. Let $oaej$ be another line quadruple in $\alpha$. Then $\overvector{zvui}=\alpha=\overvector{oaej}$ and hence $Pzvui=oaej$ for some line projectivity $P$ in $X$. The equality $Pzvui=oaej$ implies $Pzuvi=oeaj$ and then $\overvector{zuvi}=\overvector{oeaj}$. Therefore, the unary operation $(\cdot)^{-1}:\ddddot X\projeqind\to\ddddot X\projeqind$, $(\cdot)^{-1}:\alpha\mapsto\alpha^{-1}$, is well-defined.
\end{proof} 

\begin{definition} The unary operation $(\cdot)^{-1}$ on $\ddddot X\projeqind$ (resp. $\bar\IF_X$) introduced in Proposition~\ref{p:Pappian-inverse} will be called the \index{inversion}\defterm{operation of inversion} of the set of proportions (resp. proscalars) on $X$. 
\end{definition} 

Propositions~\ref{p:Pappian-inverse} and \ref{p:Pappian-quadruples}(2) imply the following corollary.

\begin{corollary}\label{c:inversion=>involution} Let $X$ be a projective space. The operation of inversion on the set $\bar\IF_X$ has the following properties:
\begin{enumerate}
\item $(p^{-1})^{-1}=p$ for every $p\in \ddddot X\projeqind$;
\item for every proscalar $p\in\bar\IF_X$ its inverse $p^{-1}$ is a proscalar.
\end{enumerate}
\end{corollary}

By Corollary~\ref{c:inversion=>involution}(2), the operation of inversion is an involution on the sets $\ddddot X\projeqind$ and $\bar\IF_X$. Let us recall that a self-map $F:X\to X$ of a set $X$ is called an \index{involution}\defterm{involution} if $F(F(x))=x$ for every $x\in X$.

\begin{exercise} Show that $[\IF_X^*]^{-1}=\IF_X^*$ and $[\ddddot X\projeqind^*]^{-1}=\ddddot X\projeqind^*$ for every projective space $X$.
\end{exercise}

The definition of the proscalars $0,1,\infty$ and Corollary~\ref{c:Pappian-permutations} imply that the following proposition.  
 
\begin{proposition} For every projective space $X$, there exists a unique unary operation\index[note]{$1$-} $$1\mbox{-}:\ddddot X\projeqind\to\ddddot X\projeqind,\quad 1\mbox{-}:p\mapsto 1\mbox{-}\,p,$$such that
\begin{enumerate}
\item $1\mbox{-}\,\overvector{zvui}=\overvector{uvzi}$ for every line quadruple $zvui\in \ddddot X$;
\item $1\mbox{-}(1\mbox{-}\,p)=p$ for every proportion $p\in\ddddot X\projeqind$;
\item $1\mbox{-}\,0=1$, $1\mbox{-}\,1=0$, $1\mbox{-}\,\infty=\infty$;
\item for every proscalar $p\in\bar\IF_X$, the proportion $1\mbox{-}\,p$ is a proscalar.
\end{enumerate}
This unary operation is called the \index{$01$-involution}\defterm{$01$-involution} of the set of proportions.
\end{proposition} 


\section{Multiplication of proportions and proscalars}

In this section we introduce and study the binary operation of multiplication of proportions and proscalars. We start with multiplication of nonzero finite proportions and proscalars.

\begin{theorem} For every projective space $X$, there exists a binary operation $$\cdot:\dom[\cdot]\to \ddddot X\projeqind^*,\quad\cdot:(\alpha,\beta)\mapsto\alpha\cdot\beta,\quad\mbox{on the set}\quad\dom[\cdot]=(\IF_X^*\times\ddddot X\projeqind^*)\cup(\ddddot X\projeqind^*\times \IF_X^*),$$which is uniquely determined by the condition
\begin{enumerate}
\item $\forall (\alpha,\beta)\in\dom[\cdot]\;\;\forall z,a,b,u,i\in X\;\;(zabi\in\alpha\;\wedge\;zbui\in\beta)\;\Ra\;(zaui\in \alpha\cdot\beta)$. 
\end{enumerate}

The binary operation $\cdot$  has the following properties for every proportions $\alpha,\beta,\gamma\in\ddddot X\projeqind^*$:
\begin{enumerate}
\item[\textup{(2)}] $1\cdot\alpha=\alpha=\alpha\cdot 1$;
\item[\textup{(3)}] if $\alpha\in\IF^*_X$, then $\alpha\cdot\beta=\beta\cdot\alpha$;
\item[\textup{(4)}] if $\alpha\in\IF^*_X$, then $\alpha\cdot\alpha^{-1}=1=\alpha^{-1}\cdot\alpha$;
\item[\textup{(5)}] if $\alpha,\beta\in\IF^*_X$, then $\alpha\cdot\beta\in\IF_X^*$;
\item[\textup{(6)}] if $\alpha,\beta\in\IF_X^*$, then $(\alpha\cdot\beta)\cdot\gamma=\alpha\cdot(\beta\cdot\gamma)$.
\end{enumerate}
\end{theorem}

\begin{proof} Given any pair $(\alpha,\beta)\in \dom[\cdot]\defeq (\IF_X^*\times\ddddot X\projeqind^*)\cup(\ddddot X\projeqind^*\times \IF_X^*)$, choose a line quadruple $zbui\in\beta$. It follows from $\beta\in \ddddot X\projeqind^*=\ddddot X\projeqind\setminus\{0,\infty\}$ that $b\notin\{z,i\}$ and hence $zbi$ is a colinear triple. By Corollary~\ref{c:proj4-exist}, there exists a point $a\in \Aline zi$ such that $zabi\in\alpha$. 
Define the product $\alpha\cdot\beta$ letting $\alpha\cdot\beta\defeq\overvector{zaui}$.
It follows from $\alpha\notin\{0,\infty\}$ that $a\notin\{z,i\}$ and hence $\overvector{zaui}\notin\{0,\infty\}$.
 Lemma~\ref{l:Pappian-multiplication} ensures that the proportion $\alpha\cdot\beta=\overvector{zaui}$ does not depend on the choice of the quadruples $zbui\in\beta$ and $zabi\in\alpha$. Therefore, the binary operation $\cdot:\dom[\cdot]\to\ddddot X\projeqind^*$ is well-defined and uniquely defined by the condition (1). It remains to prove that this binary operation has the properties (2)--(6). Fix any nonzero finite proportions $\alpha,\beta,\gamma\in\ddddot X\projeqind^*$. 
\smallskip

2. To prove that $1\cdot \beta=\beta=\beta\cdot 1$, choose any line quadruple $zbui\in\beta$. It follows from $\beta\notin\{0,\infty\}$ that $b\notin\{z,i\}$. Then $zbbi\in 1$ and $$1\cdot\beta=\overvector{zbbi}\cdot\overvector{zbui}=\overvector{zbui}=\beta=\overline{zbui}=\overline{zbui}\cdot\overline{zuui}=\beta\cdot 1,$$by the condition (1). 
\smallskip

3. Assuming that $\alpha\in\IF_X^*$, we shall prove that $\alpha\cdot\beta=\beta\cdot\alpha$.

If $1\in\{\alpha,\beta\}$, then the equality $\alpha\cdot\beta=\beta\cdot\alpha$ follows from the preceding item. So, assume that $\alpha\ne 1\ne\beta$. Choose any triangle $zij$ in $X$ and any points $u\in\Aline zi\setminus\{z,i\}$ and $v\in\Aline zj\setminus\{z,j\}$. By Corollary~\ref{c:proj4-exist}, there exist points $a\in\Aline zi$ and $b\in\Aline zj$ such that $zaui\in\alpha$ and $zbvj\in\beta$.
Since the liner $X$ is projective, there exist unique points $c\in\Aline av\cap\Aline ij$,  $p\in\Aline cb\cap\Aline zi$, $d\in\Aline ub\cap\Aline ij$, and $q\in \Aline ad\cap\Aline zj$.

\begin{picture}(300,155)(-160,-15)

\put(0,0){\line(1,2){60}}
\put(0,0){\line(-1,2){60}}
\put(-129,120){\line(1,0){336}}
\put(-10,20){\line(13,6){217}}
\put(-40,80){\line(247,40){247}} 
\put(16,32){\line(-23,14){145}}
\put(-129,120){\line(89,-40){154}}
\put(-30,60){\line(77,34){136}}
\put(-10,20){\line(116,100){116}}

\put(207,120){\circle*{3}}
\put(210,118){$w$}
\put(60,120){\circle*{3}}
\put(58,125){$j$}
\put(-60,120){\circle*{3}}
\put(-62,123){$i$}
\put(-10,20){\circle*{3}}
\put(-18,17){$p$}
\put(106,120){\circle*{3}}
\put(104,123){$c$}
\put(47,94){\circle*{3}}
\put(41,96){$v$}
\put(25.4,50.8){\circle*{3}}
\put(22,54){$b$}
\put(-129,120){\circle*{3}}
\put(-131,123){$d$}
\put(16,32){\circle*{3}}
\put(17,26){$q$}
\put(-40,80){\circle*{3}}
\put(-39,83){$u$}
\put(-30,60){\circle*{3}}
\put(-37,55){$a$}
\put(0,0){\circle*{3}}
\put(-2,-8){$z$}
\end{picture}

 Consider the line perspectivities 
$$\begin{aligned}
C&\defeq\{(x,y)\in \Aline zj\times\Aline zi:y\in \Aline cx\},&D&\defeq\{(x,y)\in \Aline zi\times\Aline zj:y\in \Aline dx\},\\
U&\defeq\{(x,y)\in \Aline ij\times\Aline zj:y\in \Aline ux\},&
V&\defeq\{(x,y)\in \Aline zi\times\Aline ij:y\in \Aline vx\},\\
P&\defeq\{(x,y)\in \Aline ij\times\Aline zj:y\in \Aline px\},&
Q&\defeq\{(x,y)\in \Aline zi\times\Aline ij:y\in \Aline qx\},
\end{aligned}
$$
with centers $c,d,u,v,p,q$, respectively.

The equalities $Dzaui=zqbj$ and $Czbvj=zpai$ imply the equalities $\overvector{zqbj}=\overvector{zaui}$ and $\overvector{zbvj}=\overvector{zpai}$. Then $$\alpha\cdot\beta=\overvector{zaui}\cdot\overvector{zbvj}=\overvector{zqbi}\cdot\overvector{zbvj}=\overvector{zqvj}\quad\mbox{and}\quad\beta\cdot\alpha=\overvector{zbvj}\cdot\overvector{zaui}=\overvector{zpai}\cdot\overvector{zaui}=\overvector{zpui}.$$
The equality $\alpha\cdot\beta=\beta\cdot\alpha$ will follow as soon as we show that $\overvector{zqvj}=\overvector{zpui}$. 

Consider the line projectivities $PV$ and $UQ$ and observe that 
$$PVzai=Pjci=jbz=Ujdi=UQzai.$$
Since the quadruple $zaui\in\alpha\in \IF_X^*$ is Pappian, so is the quadruple $zuai$, see Proposition~\ref{p:Pappian-quadruples}(1). Then  the equality $PVzai=UQzai$ implies $PV(u)=UQ(u)=q$. Then for the point $w\defeq V(u)\in \Aline uv\cap\Aline ij$ we have $q=UQ(u)=PV(u)=P(w)$ and hence $Wzpui=zqvj$ for the line perspectivity $W\defeq\{(x,y)\in \Aline zi\times\Aline zj:y\in\Aline wx\}$ with center $w$.
Therefore, $\alpha\cdot\beta=\overvector{zqvj}=\overvector{zpui}=\beta\cdot\alpha$.
\smallskip

4. Assuming that $\alpha\in\IF_X^*$, we shall prove that $\alpha\cdot\alpha^{-1}=1=\alpha^{-1}\cdot \alpha$. Choose any line quadruples $zaui\in\alpha$. The definition of the inverse proportion $\alpha^{-1}$ ensures that $zuai\in\alpha^{-1}$. The condition (1) ensures that 
$$\alpha^{-1}\cdot \alpha=\overvector{zuai}\cdot \overvector{zaui}=\overvector{zuui}=1=\overvector{zaai}=\overvector{zaui}\cdot\overvector{zuai}=\alpha\cdot\alpha^{-1}.$$

5. Assuming that $\alpha,\beta\in\IF_X^*$, we need to show that $\alpha\cdot\beta\in\IF^*_X$. Choose any line quadruple $zbui\in\beta$. It follows from $\beta\ne 0$ that $b\ne z$. By Corollary~\ref{c:proj4-exist}, there exists a point $a\in \Aline zb$ such that $zabi\in \alpha$. It follows from $\alpha\notin\{0,\infty\}$ that $a\notin\{z,i\}$. Since $\alpha,\beta\in\IF_X^*$, the line quadruples $zabi\in\alpha$ and $zbui\in\beta$ are Pappian. We claim that the line quadruple $zaui$ is Pappian. Given any line projectivity $P$ with $Pzui=zui$, we need to show that $P(a)=a$. Since the line quadruple $zbui$ is Pappian, the equality $Pzui=zui$ implies $P(b)=b$ and hence $Pzbi=zbi$. Since the  line quadruple $zabi$ is Pappian, the equality $Pzbi=zbi$ implies $P(a)=a$, witnessing that the line quadruple $zaui$ is Pappian and the proportion $\alpha\cdot\beta=\overvector{zaui}$ belongs to the set $\bar\IF_X$ of proscalars. Since $a\notin\{z,i\}$, the proscalar $\alpha\cdot\beta=\overvector{zaui}$ belongs to the set $\IF_X^*=\bar\IF_X\setminus\{0,\infty\}$ of nonzero finite proscalars.
\smallskip

6. Assuming that $\alpha,\beta\in\IF_X^*$, we shall prove that $(\alpha\cdot\beta)\cdot\gamma=\alpha\cdot(\beta\cdot\gamma)$. By the preceding item, $\alpha\cdot\beta\in\IF_X^*$, so that the product $(\alpha\cdot\beta)\cdot\gamma$ is well-defined. Choose any line quadruple $zcui\in\gamma$. It follows from $\gamma\ne 0$ that $c\ne z$. By Corollary~\ref{c:proj4-exist}, there exists a line quadruple $zbci\in\beta$.  It follows from $\beta\ne 0$ that $b\ne z$.  By Corollary~\ref{c:proj4-exist}, there exists a line quadruple $zabi\in\alpha$. By the condition (1),
$$\alpha\cdot(\beta\cdot\gamma)=\overvector{zabi}\cdot(\overvector{zbci}\cdot\overvector{zcui})=\overvector{zabi}\cdot\overvector{zbui}=\overvector{zaui}=\overvector{zaci}\cdot\overvector{zcui}=(\overvector{zabi}\cdot\overvector{zbci})\cdot\overvector{zcui}=(\alpha\cdot\beta)\cdot\gamma.$$
\end{proof}

\begin{corollary}\label{c:FX*-group} For every projective space $X$, the set $\IF_X^*$ is a commutative group with respect to the binary operation $\cdot:\IF_X^*\times\IF_X^*\to\IF_X^*$, defined by the formula: 
$$\forall \alpha,\beta\in\IF_X^*\;\forall z,a,b,u,i\in X\;(zabi\in\alpha\;\wedge\;zbui\in\beta)\;\Ra\;zaui\in\alpha\cdot\beta).$$
\end{corollary}

For a projective space $X$, it is natural to extend the multiplication operation from the commutative group $\IF_X^*$ to the set of all finite proscalars $\IF_X=\IF_X\cup\{0\}$ letting $0\cdot \alpha=0=\alpha\cdot 0$ for all $\alpha\in\IF_X$. Then $\IF_X$ became a commutative monoid with zero, containing $\IF_X^*$ as a subgroup. 

 We can also extend the multiplicative operation further on to the set $\bar\IF_X$ of all proscalars letting 
$$\infty\cdot 0=1=0\cdot \infty\quad\mbox{and}\quad \forall \alpha\in\bar\IF_X\setminus\{0\}\;\;(\alpha\cdot\infty=\infty=\infty\cdot\alpha).$$ Then $\bar F_X$ is a commutative unital magma such that $\alpha\cdot\alpha^{-1}=1=\alpha^{-1}\cdot\alpha$ for all $\alpha\in\bar\IF_X$. The subset $\bar\IF_X\setminus\{0\}$ of $\bar\IF_X$ is a commutative monoid, isomorphic to the commutative monoid $\IF_X$. However, the magma $\bar\IF_X$ is not a semigroup because its binary operation is not associative:
$$(0\cdot 0)\cdot \infty=0\cdot\infty=1\ne 0=0\cdot 1=0\cdot (0\cdot\infty).$$On the other hand, for every $\alpha,\gamma\in\bar\IF_X$ and $\beta\in\IF_X^*$
the associative law $(\alpha\cdot\beta)\cdot\gamma=\alpha\cdot(\beta\cdot\gamma)$ does hold.

\section{Proscalars versus scalars}

In this section we study the interplay between scalars in an affine space  and proscalars in its  projective completion. First, we prove three lemmas on extension of line projections to line perspectivities and line affinities to line projectivities.
 

\begin{lemma}\label{l:line-projection-extends-to-line-persectivity} Let $X$ be a projective completion of an affine space $A$. For every line projection $P$ in $A$, the set $$\bar P\defeq P\cup\big((\overline{\dom[P]}\setminus A)\times(\overline{\rng[P]}\setminus A)\big)$$is a line perspectivity in $X$ extending the line projection $P$.
\end{lemma}

\begin{proof} By Theorem~\ref{t:affine<=>hyperplane}, the horizon $H\defeq X\setminus A$ of $A$ in $X$ is a hyperplane in $X$. This implies that for every line $L$ in $A$, its flat hull $\overline L$ in $X$ is a line in $X$ such that $|\overline L\cap H|\le 1$. By Corollary~\ref{c:line-meets-hyperplane}, the intersection $\overline L\cap H$ is not empty, which implies that $\overline L\setminus A$ is a singleton. Let $a\in\overline{\dom[P]}\cap H$ and $b\in \overline{\rng[P]}\cap H$ be unique common points of the lines $\overline{\dom[P]},\overline{\rng[P]}$ with the horizon $H$. Then the set $$\bar P\defeq P\cup\big((\overline{\dom[P]}\setminus A)\times(\overline{\rng[P]}\setminus A)\big)=P\cup\{(a,b)\}$$ is a function. 

Since $P$ is a line projection in $A$, there exists a line $\Lambda$ in $A$ such that $$P=\{(x,y)\in \dom[P]\times\rng[P]:\Aline xy\subparallel \Lambda\}$$ and $\Lambda$ is concurrent with the lines $\dom[P]$ and $\rng[P]$. Let $c$ be the unique point of the set $\overline \Lambda\cap H$. Since the line $\Lambda$ is concurrent with the lines $\dom[P]$ and $\rng[P]$, the point $c$ does not belong to the lines $\overline{\dom[P]}$ and $\overline{\rng[P]}$ in $X$. For every pair $(x,y)\in P$ with $x\ne y$, the line $\Aline xy$ is parallel to the line $\Lambda$ in $A$ and hence $\overline{\{x,y\}}\cap\overline\Lambda=\{c\}$, which implies $P=\{(x,y)\in \dom[P]\times\rng[P]:y\in\overline{\{c,x\}}\}$.

 Since $\Pi\defeq\overline{\dom[P]\cup\rng[P]}$ is a plane in the projective space $X$, the intersection $\Pi\cap H$ is a line in $\Pi$ containing the points $a,b$ and $c\notin\{a,b\}$. Then $\Aline ac=\Pi\cap H=\Aline ay$ and hence  $$\bar P=P\cup\{(a,b)\}=\{(x,y)\in\overline{\dom[P]}\times\overline{\rng[P]}:y\in\overline{\{c,x\}}\}$$ is a line perspectivity in $X$ extending the line projection $P$. 
\end{proof}

\begin{lemma}\label{l:line-affinity-extends-to-line-projectivity} Let $X$ be a projective completion of an affine space $A$. For every line affinity $P$ in $A$,  the set $$\bar P\defeq P\cup\big((\overline{\dom[P]}\setminus A)\times(\overline{\rng[P]}\setminus A)\big)$$is a line projectivity in $X$ extending the line affinity $P$.
\end{lemma}

\begin{proof} By definition, the line affinity $P$ is the composition $P_1\cdots P_n$ of some line projections $P_1,\dots,P_n$ in $A$. By Lemma~\ref{l:line-projection-extends-to-line-persectivity}, for every $i\in\{1,\dots,n\}$, the set $$\bar P_i\defeq P_i\cup\big((\overline{\dom[P_i]}\setminus A)\times(\overline{\rng[P_i]}\setminus A)\big)$$ is a line perspectivity in $X$, extending the line projection $P_i$. Then $P'\defeq \bar P_1\cdots \bar P_n$ is a line projectivity in $X$ extending the line affinity $P$. Observe that $\dom[P']=\dom[\bar P_n]=\overline{\dom[P_n]}=\overline{\dom[P]}$ and $\rng[P']=\rng[\bar P_1]=\overline{\rng[P_1]}=\overline{\rng[P]}$. Since the line projectivity $P':\overline{\dom[P]}\to\overline{\rng[P]}$ is a bijective map extending the bijective map $P:\dom[P]\to\rng[P]$, we conclude that $P'[\overline{\dom[P]}\setminus\dom[P]]=P'[\overline{\dom[P]}]\setminus P'[\dom[P]]=\overline{\rng[P]}\setminus \rng[P]$. Taking into account that the sets $\overline{\dom[P]}\setminus\dom[P]$ and $\overline{\rng[P]}\setminus \rng[P]$ are singletons, we conclude that $$
\bar P\defeq P\cup\big((\overline{\dom[P]}\setminus A)\times(\overline{\rng[P]}\setminus A)\big)=P'$$ and hence $\bar P=P'$ is a line projectivity in $X$, extending the line affinity $P$.
\end{proof}

\begin{lemma}\label{l:Pappian-quadruple=>Desarg-triple} Let $X$ be a projective completion of an affine space $A$, $H\defeq X\setminus A$ be the horizon of $A$ in $X$, and $\ddddot X_{\!\!A}\defeq \ddddot X\cap (A^3\times H)$.
 For every Pappian quadruple $zvui\in \ddddot X_{\!\!A}$, the triple $zvu$ is Desarguesian in the affine space $A$. Moreover, for any Pappian quadruples $\hat z\hat v\hat u\hat\imath, \check z\check v\check u\check \imath\in \ddddot X_{\!\!A}$ with $\hat z\hat v\hat u\hat\imath\projeq \check z\check v\check u\check \imath$ we have $\hat z\hat v\hat u\Join \check z\check v\check u$.
\end{lemma}

\begin{proof} Let $zvui\in\ddddot X_{\!\!A}$ be a Pappian line quadruple. To show that the line triple $zvu$ is Desarguesian, take any line affinity $F$ in $A$ such that $Fzu=zu$. By Lemma~\ref{l:line-affinity-extends-to-line-projectivity}, the set
$$\bar F\defeq F\cup\big((\overline{\dom[F]}\setminus A)\times(\overline{\rng[F]}\setminus A)\big)$$is a line projectivity in $X$ extending the line affinity $F$. Since $Fzu=zu$, $\dom[F]=\Aline zu=\rng[F]$ in $A$ and hence $\dom[\bar F]=\overline{\{z,u\}}=\rng[\bar F]$ in $X$. By Theorem~\ref{t:affine<=>hyperplane}, the horizon $H\defeq X\setminus A$ is a hyperplane in the projective space $X$. Consider the line $\dom[\bar F]=\rng[\bar F]$ in $X$ and let $i$ be the unique point of the intersection $\dom[\bar F]\cap H=\rng[\bar F]\cap H$ (which is not empty, by Corollary~\ref{c:line-meets-hyperplane}). Then $\bar F(i)\in \bar F[\dom[\bar F]\setminus A]=\rng[\bar F]\setminus A=\{i\}$ and hence $\bar Fzui=zui$. Since the line quadruple 
$zvui$ is Pappian, $\bar F(v)=v$ and hence $F(v)=\bar F(v)=v$, witnessing that the line triple $zvu$ is Desarguesian.
\smallskip

Now choose any Pappian line quadruples $\hat z\hat v\hat u\hat\imath, \check z\check v\check u\check \imath\in \ddddot X_{\!\!A}$ with $\hat z\hat v\hat u\hat\imath\projeq \check z\check v\check u\check \imath$. By Theorem~\ref{t:aff-trans}, there exists a line affinity $F$ in $X$ such that $F\hat z\hat u=\check z\check u$. By Lemma~\ref{l:line-affinity-extends-to-line-projectivity}, the set $$\bar F\defeq F\cup\big((\overline{\dom[F]}\setminus A)\times(\overline{\rng[F]}\setminus A)\big)$$is a line projectivity in $X$ extending the line affinity $F$. Then $\hat \imath\in\overline{\{\hat z,\hat u\}}=\overline{\dom[F]}=\dom[\bar F]$ and  $\bar F(\hat \imath)\in \bar F[\dom[\bar F]\setminus \dom[F]]=\rng[\bar F]\setminus \rng[F]=\overline{\{\check z,\check u\}}\setminus A=\{\check \imath\}$. Therefore, $\bar F$ is a line projectivity in $X$ such that $\bar F\hat z\hat u\hat\imath=\check z\check u\check\imath$. Since the quadruples $\hat z\hat v\hat u\hat\imath$ and $\check z\check v\check u\check \imath$ are projectively equivalent and Pappian,  $\bar F\hat z\hat u\hat\imath=\check z\check u\check\imath$ implies $F(\hat v)=\bar F(\hat v)=\check v$ and hence $F\hat z\hat v\hat u=\check z\check v\check u$, witnessing that the line triples  $\hat  z\hat v\hat u$ and $\check z\check v\check u$ are affinely equivalent (which is denoted by  $\hat  z\hat v\hat u\Join\check z\check v\check u$).
 \end{proof}

Let $X$ be a projective completion of an affine space $A$ and let $H\defeq X\setminus A$ be the horizon of $A$ in $X$. By Theorem~\ref{t:affine<=>hyperplane}, $H$ is a hyperplane in $X$. Consider the set of line quadruples 
$$\ddddot X_{\!\!A}\defeq\ddddot X\cap (A^3\times H)=\{zvui\in\ddddot X:z,v,u\in A,\;i\in H\}.$$


Consider the map $\Psi_A:\IF_X\to\IR_A$ assigning to every finite proscalar $\alpha\in\IF_X$ the scalar $\overvector{zau}$ where $zaui$ is any line quadruple in $\alpha\cap\ddddot X_{\!\!A}$ (such a line quadruple $zaui$ exists, by Corollary~\ref{c:proj4-exist}).  Lemma~\ref{l:Pappian-quadruple=>Desarg-triple} ensures that the value $\Psi_A(\alpha)=\overvector{zau}$ does not depend on the choice of the line quadruple $zaui\in\alpha\cap\ddddot X_{\!\!A}$. Therefore, $\Psi_A:\IF_{\!X}\to\IR_A$ is a well-defined function.

\begin{theorem}\label{t:PsiA} For any projective completion $X$ of an affine space $A$, the function $\Psi_A:\IF_X\to\IR_A$ has the following properties:
\begin{enumerate}
\item $\Psi_A$ is injective;
\item $\Psi_A(0)=0$ and $\Psi_A(1)=1$;
\item $\Psi_A(\alpha\cdot\beta)=\Psi_A(\alpha)\cdot\Psi_A(\beta)$ for any proscalars $\alpha,\beta\in\IF_X$;
\item $\Psi_A(1\mbox{-}\,\alpha)=1-\Psi_A(\alpha)$ for every proscalar $\alpha\in\IF_X$; 
\item $\Psi_A[\IF_X]\subseteq Z(\dddot A_{\Join})\defeq\{\alpha\in \IR_A:\forall \beta\in \dddot A_{\Join}\;\;(\alpha\cdot \beta=\beta\cdot \alpha)\}$.
\end{enumerate}
\end{theorem} 

\begin{proof} By Theorem~\ref{t:affine<=>hyperplane}, the horizon $H\defeq X\setminus A$ of $A$ in $X$ is a hyperplane in $X$. 
\smallskip

 1. To prove that the function $\Psi_A$ is injective, take any proscalars $\alpha,\beta\in\IF_X$ and assume that $\Psi_A(\alpha)=\Psi_A(\beta)$. Choose any points $i\in H$, $z\in A$, and $u\in\Aline zi\setminus\{z,i\}$. Using Corollary~\ref{c:proj4-exist}, find points $a,b\in\Aline zi$ such that $zaui\in \alpha$ and $zbui\in\beta$. Since $\overvector{zau}=\Psi_A(\alpha)=\Psi_A(\beta)=\overvector{zbu}$, the line triples $zau$ and $zbu$ are affinely equivalent in the affine space $A$. Then there exists a line affinity $F$ in $A$ such that $Fzau=Fzbu$. By Lemma~\ref{l:line-affinity-extends-to-line-projectivity}, the set
$$\bar F\defeq F\cup\big((\overline{\dim[F]}\setminus A)\times(\overline{\rng[F]}\setminus A)\big)=F\cup\{(i,i)\}$$is a line projectivity extending the line affinity $F$. Then $\bar Fzaui=zbui$, witnessing that the line quadruples $zaui$ and $zbui$ are projectively equivalent and hence $\alpha=\overvector{zaui}=\overvector{zbui}=\beta$. Therefore, the function $\Psi_A$ is injective.
\smallskip

2. Choose any points $i\in H$, $z\in A$, and $u\in\Aline zi\setminus\{z,i\}$. The definitions of the function $\Psi_X$ and the (pro)scalars $0$ and $1$ ensure that $\Psi_A(0)=\Psi_A(\overvector{zzui})=\overvector{zzu}=0$ and $\Psi_A(1)=\Psi_A(\overvector{zuui})=\overvector{zuu}=1$.
\smallskip

3. Given any finite proscalars $\alpha,\beta\in\IF_X$, we need to check that $\Psi(\alpha\cdot\beta)=\Psi_A(\alpha)\cdot\Psi_A(\beta)$. If $\beta=0$, then $\Psi_A(\alpha\cdot\beta)=\Psi_A(\alpha\cdot 0)=\Psi_A(0)=0=\Psi_A(\alpha)\cdot 0=\Psi_A(\alpha)\cdot\Psi_A(\beta)$, by the preceding item. So, assume that $\beta\ne 0$. Choose any line quadruple $zbui\in\beta\cap (A^3\times H)$. It follows from $0\ne\beta\ne \infty$ that $z\ne b\ne i$. By Corollary~\ref{c:proj4-exist}, there exists a point $a\in\Aline zi$ such that $zabi\in\alpha$. It follows from $\alpha\ne\infty$ that $a\in\Aline zi\setminus\{i\}\subseteq A$.  The condition (1) and Theorem~\ref{t:scalar-by-portion} ensure that\newline
$\Psi_A(\alpha\cdot\beta)=\Psi_A(\overvector{zabi}\cdot\overvector{zbui})=\Psi_A(\overvector{zaui})=\overvector{zau}=\overvector{zab}\cdot\overvector{zbu}=\Psi_A(\overvector{zabi})\cdot\Psi_A(\overvector{zbui})=\Psi_A(\alpha)\cdot\Psi_A(\beta)$.
\smallskip

4. Given any finite proscalar $\alpha\in\IF_X$, choose a line quadruple $zvui\in\alpha\cap (A^3\times H)$. The definition of the functions $1\mbox{-}:\bar\IF_X\to\bar\IR_X$ and $1\mbox{-}:\IR_A\to\IR_A$, Proposition~\ref{p:1-portion} and Theorem~\ref{t:1-=1-} imply that
$\Psi_A(1\mbox{-}\,\alpha)=\Psi_A(1\mbox{-}\,\overvector{zvui})=\Psi_A(\overvector{uvzi})=\overvector{uvz}=1\mbox{-}\,\overvector{zvu}=1-\Psi_A(\alpha).$ 

5. To prove that $\Psi_A[\IF_X]\subseteq Z(\dddot X_{\Join})$, take any finite proscalar $\alpha\in\IF_X$ and find a line quadruple $zaui\in\alpha\cap (A^3\times H)$. The definition of the function $\Psi_A$ ensures that $\Psi_A(\alpha)=\overvector{zau}$. To prove that $\overvector{zau}\in Z(\dddot A_{\Join})$, we should check that $\overvector{zau}\cdot\beta=\beta\cdot\overvector{zau}$ for every portion $\beta\in\dddot A_{\Join}$. If $a=z$, then
$$\overvector{zau}\cdot\beta=\overvector{zzu}\cdot\beta=0\cdot\beta=0=\beta\cdot 0=\beta\cdot\overvector{zzu}=\beta\cdot\overvector{zau}$$and we are done.
If $\beta=0$, then again $\overvector{zau}\cdot\beta=\overvector{zau}\cdot0=0=0\cdot\overvector{zau}=\beta\cdot\overvector{zau}$ and we are done.

So, assume that that $a\ne z$ and $\beta\ne0$. Choose any point $j\in H\setminus\{i\}$ and any point $v\in\Aline zj\setminus\{z,j\}$. By Corollary~\ref{c:proj4-exist}, there exists a point $b\in\Aline zv\cap A$ such that $zbv\in\beta$. It follows from $\beta\ne0$ that $b\ne z$. 
Since the liner $X$ is projective, there exist unique points $c\in\Aline av\cap\Aline ij$,  $p\in\Aline cb\cap\Aline zi$, $d\in\Aline ub\cap\Aline ij$, and $q\in \Aline ad\cap\Aline zj$. Since the horizon $H$ of $A$ in $X$ is flat, $c,d\in \Aline ij\subseteq H$.

\begin{picture}(300,155)(-160,-15)

\put(0,0){\line(1,2){60}}
\put(0,0){\line(-1,2){60}}
\put(-129,120){\line(1,0){336}}
\put(-10,20){\line(13,6){217}}
\put(-40,80){\line(247,40){247}} 
\put(16,32){\line(-23,14){145}}
\put(-129,120){\line(89,-40){154}}
\put(-30,60){\line(77,34){136}}
\put(-10,20){\line(116,100){116}}

\put(207,120){\circle*{3}}
\put(210,118){$w$}
\put(60,120){\circle*{3}}
\put(58,125){$j$}
\put(-60,120){\circle*{3}}
\put(-62,123){$i$}
\put(-10,20){\circle*{3}}
\put(-18,17){$p$}
\put(106,120){\circle*{3}}
\put(104,123){$c$}
\put(47,94){\circle*{3}}
\put(41,96){$v$}
\put(25.4,50.8){\circle*{3}}
\put(22,54){$b$}
\put(-129,120){\circle*{3}}
\put(-131,123){$d$}
\put(16,32){\circle*{3}}
\put(17,26){$q$}
\put(-40,80){\circle*{3}}
\put(-39,83){$u$}
\put(-30,60){\circle*{3}}
\put(-37,55){$a$}
\put(0,0){\circle*{3}}
\put(-2,-8){$z$}
\end{picture}

 Consider the line perspectivities 
$$\begin{aligned}
C&\defeq\{(x,y)\in \Aline zj\times\Aline zi:y\in \Aline cx\},&D&\defeq\{(x,y)\in \Aline zi\times\Aline zj:y\in \Aline dx\},\\
U&\defeq\{(x,y)\in \Aline ij\times\Aline zj:y\in \Aline ux\},&
V&\defeq\{(x,y)\in \Aline zi\times\Aline ij:y\in \Aline vx\},\\
P&\defeq\{(x,y)\in \Aline ij\times\Aline zj:y\in \Aline px\},&
Q&\defeq\{(x,y)\in \Aline zi\times\Aline ij:y\in \Aline qx\},
\end{aligned}
$$
in $X$ with centers $c,d,u,v,p,q$, respectively.

Since $c,d\in H$, the line perspectivities $C$ and $D$ resptricted to $A$ are line projections on $A$. The equalities $Dzau=zqb$ and $Czbv=zpa$ imply the equalities $\overvector{zqb}=\overvector{zau}$ and $\overvector{zbv}=\overvector{zpa}$. Then $$\overvector{zau}\cdot\beta=\overvector{zqb}\cdot\overvector{zbv}=\overvector{zqv}\quad\mbox{and}\quad \beta\cdot\overvector{zau}=\overvector{zbv}\cdot\overvector{zau}=\overvector{zpa}\cdot\overvector{zau}=\overvector{zpu}.$$
The equality $\overvector{zau}\cdot\beta=\beta\cdot\overvector{zau}$ will follow as soon as we show that $\overvector{zqv}=\overvector{zpu}$. 

Consider the line projectivities $PV$ and $UQ$ and observe that 
$$PVzai=Pjci=jbz=Ujdi=UQzai.$$
Since the quadruple $zaui\in\alpha\in \IF_X$ is Pappian, so is the quadruple $zuai$, see Proposition~\ref{p:Pappian-quadruples}(1). Then  the equality $PVzai=UQzai$ implies $PV(u)=UQ(u)=q$. For the point $w\defeq V(u)\in \Aline uv\cap\Aline ij\subseteq H$ we have $q=UQ(u)=PV(u)=P(w)$ and hence $Wzpu=zqv$ for the line perspectivity
$W\defeq\{(x,y)\in \Aline zi\times\Aline zj:y\in\Aline wx\}$ with center $w$. Since $w\in H$, the restriction $W{\restriction}_A$ is a line projection in the affine space, witnessing that $\overvector{zau}\cdot\beta=\overvector{zqv}=\overvector{zpu}=\beta\cdot\overvector{zau}$ and hence $\Psi_A(\alpha)=\overvector{zau}\in Z(\dddot A_{\Join})$.
\end{proof}

Theorem~\ref{t:PsiA}(5) motivates the following natural question.

\begin{question}\label{q:center} Let $X$ be a projective completion of an affine space $A$. Is $\Psi_A[\IF_X]=Z(\dddot A_{\Join})$?
\end{question}

In Theorem~\ref{t:D=>FX=Z(RX)} we shall prove that the answer to this question is affirmative for Desarguesian spaces. The proof exploits the following non-trivial lemma, proved by tools of analytic geometry.

\begin{lemma}\label{l:anal-geom} Let $\Pi$ be a plane in a Desarguesian affine space $A$ and $zuv$ and $opq$ be two triangles in $\Pi$ such that $z\in\Aline pq$ and $\Aline zu\cap\Aline op=\varnothing=\Aline zv\cap\Aline oq$. Let $\Lambda$ be a line in $\Pi$ such that $\{o\}=\Lambda\cap\Aline op=\Lambda\cap\Aline oq$ and  the lines $\Lambda,\Aline pu,\Aline qv$ are paraconcurrent in $A$. Let $a\in\Aline zu\setminus\{z\}$ and $b\in\Aline z{v}\setminus\{z\}$ be points such that the lines $\Lambda,\Aline pa,\Aline qb$ are paraconcurrent. 
If $$\overvector{zau}\in Z(\IR_A)\defeq\{\alpha\in\IR_A:\forall \beta\in\IR_A\;\;(\alpha\cdot\beta=\beta\cdot\alpha)\},$$ then the lines $\Aline uv$ and $\Aline ab$ are parallel.
\end{lemma}

\begin{proof} A picture visualizing the hypothese and conclusion of the lemma looks as follows.

\begin{picture}(120,130)(-160,-15)

{\linethickness{0.75pt}
\put(0,0){\color{blue}\vector(1,0){45}}
\put(0,0){\color{blue}\line(1,0){105}}
\put(0,0){\color{teal}\vector(0,1){45}}
\put(0,0){\color{teal}\line(0,1){100}}
\put(30,15){\color{blue}\line(1,0){75}}
\put(30,15){\color{teal}\line(0,1){85}}
\put(0,0){\line(1,1){105}}

\put(30,45){\color{red}\line(1,-2){15}}
\put(30,60){\color{red}\line(1,-2){22.5}}
}
\put(45,0){\line(0,1){45}}
\put(0,45){\line(1,0){45}}
\put(45,0){\line(1,2){45}}
\put(0,45){\line(2,1){90}}
\put(0,45){\line(1,-1){45}}

\put(100,90){$\Lambda$}
\put(0,0){\circle*{3}}
\put(-5,-8){$o$}
\put(45,15){\circle*{3}}
\put(38,8){$a$}
\put(30,15){\circle*{3}}
\put(25,8){$z$}
\put(52.5,15){\circle*{3}}
\put(53,8){$u$}
\put(0,45){\circle*{3}}
\put(-7,44){$q$}
\put(30,45){\circle*{3}}
\put(23,47){$b$}
\put(30,60){\circle*{3}}
\put(23,60){$v$}
\put(45,0){\circle*{3}}
\put(43,-8){$p$}
\put(45,45){\circle*{2}}
\put(90,90){\circle*{2}}
\end{picture}

Since $A$ is Desarguesian, $\overvector A$ is an $\IR_A$-module and for every point $x\in \Pi$, the vector $\overvector{ox}$ can be written uniquely as the linear combination $\overvector{ox}=x_1\cdot\overvector {op}+x_2\cdot\overvector{oq}$ for some scalars $x_1,x_2\in\IR_A$. Therefore, the points of the plane $\Pi$ can be identified with pairs of scalars $(x_1,x_2)\in\IR_A\times\IR_A$.

\begin{claim} $z\in \Aline pq$ implies $z_1+z_2=1$.
\end{claim}

\begin{proof}  It follows from $z\in\Aline pq$ that $\overvector{pz}=\overvector{pzq}\cdot \overvector{pq}$ and hence
$$
z_1\cdot\overvector{op}+z_2\cdot\overvector{oq}=\overvector{oz}=\overvector{op}+\overvector{pz}=\overvector{op}+\overvector{pzq}\cdot\overvector{pq}=\overvector{op}+\overvector{pzq}\cdot(\overvector{oq}-\overvector{op})=
(1-\overvector{pzq})\cdot\overvector{op}+\overvector{pzq}\cdot\overvector{oq}
$$
which implies $z_1=1-\overvector{pzq}$, $z_2=\overvector{pzq}$ and finally 
$z_1+z_2=(1-\overvector{pzq})+\overvector{pzq}=1$.
 \end{proof}

Since the line $\Lambda\subseteq \Pi$ is concurrent with the line $\Aline oq$, there exists a unique scalar $\lambda\in\IR_A$ such that the vector $\overvector{op}+\lambda\cdot\overvector{oq}$ is parallel to the line $\Lambda$. Then the line  $\Lambda$ coincides with the set $\{o+x\cdot \overvector{op}+x\cdot\lambda\cdot\overvector{oq}:x\in\IR_A\}$. 

\begin{claim}\label{cl:line-equation} For any distinct points $c,e\in\Pi$, the line $\Aline ce$ consists of points $x\in\Pi$ whose coordinates $(x_1,x_2)$ satisfy the equation
$$x_2=c_2+(x_1-c_1)\cdot(e_1-c_1)^{-1}\cdot (e_2-c_2)$$
in the corps $\IR_A$.
 \end{claim}

\begin{proof} A point $x\in\Pi$ belongs to the line $\Aline ce$ if and only if $\overvector{cx}=s\cdot \overvector{ce}$ for some scalar $s\in\IR_A$. Then 
\begin{multline*}
(x_1-c_1)\cdot\overvector{op}+(x_2-c_2)\cdot\overvector{oq}=
(x_1\cdot\overvector{op}+x_2\cdot\overvector{oq})-(c_1\cdot\overvector{op}+c_2\cdot\overvector{oq})\\
=\overvector{ox}-\overvector{oc}=\overvector{cx}=s\cdot\overvector{ce}=s\cdot((e_1-c_1)\cdot \overvector{op}+(e_2-c_2)\cdot\overvector{oq})
\end{multline*}
and hence
$$x_1-c_1=s\cdot(e_1-c_1)\quad\mbox{and}\quad x_2-c_2=s\cdot (e_2-c_2),$$
which is equivalent to $(x_1-c_1)\cdot (e_1-c_1)^{-1}=s=(x_2-c_2)\cdot(e_2-c_2)^{-1}$ and to
the required equation of the line $\Aline ce$: $x_2=c_2+(x_1-c_1)\cdot(e_1-c_1)^{-1}\cdot (e_2-c_2)$. 
\end{proof}

\begin{claim}\label{cl:anal-geom1} $v_2-z_2=z_1\cdot \lambda\cdot z_2^{-1}\cdot (u_1-z_1)$.
\end{claim}

\begin{proof} Since $\Aline zv\parallel \Aline oq$, the points $z$ and $v$ have the same first coordinates $z_1=v_1$. Since $\Aline zu\parallel \Aline op$, the points $z$ and $p$ have the same second coordinate $z_2=u_2$. If the lines $\Lambda,\Aline pu,\Aline {q}{v}$ are concurrent, then their common point has coordinates $(x,x\cdot \lambda)$ satisfying the equations of the lines $\Aline pu$ and $\Aline qv$:
$$x\cdot \lambda=(x-1)\cdot (u_1-1)^{-1}\cdot u_2\quad\mbox{and}\quad x\cdot\lambda=1+x\cdot v_1^{-1}\cdot (v_2-1),$$
by Claim~\ref{cl:line-equation}. 

Taking into account that $v_1=z_1$, $u_2=z_2$ and $z_1+z_2=1$, we conclude that 
$$
x-1=x\cdot\lambda\cdot u_2^{-1}\cdot (u_1-1)=x\cdot\lambda\cdot z_2^{-1}\cdot (u_1-z_1-z_2)=
x\cdot \lambda\cdot (z_2^{-1}(u_1-z_1)-1)$$and hence
\begin{equation}\label{eq:angeom1}
x\cdot (1+\lambda-\lambda\cdot z_2^{-1}\cdot (u_1-z_1))=1.
\end{equation}
On the other hand, the equation $$x\cdot\lambda=1+x\cdot v_1^{-1}\cdot (v_2-1)=1+x\cdot z_1^{-1}\cdot(v_2-z_2-z_1)=1+x\cdot(z_1^{-1}(v_2-z_2)-1)$$ implies
\begin{equation}\label{eq:angeom2}
1=x\cdot(1+\lambda-z_1^{-1}\cdot(v_2-z_2)).
\end{equation}
The equations (\ref{eq:angeom1}) and (\ref{eq:angeom2}) imply  
$\lambda\cdot z_2^{-1}\cdot(u_1-z_1)=z_1^{-1}\cdot(v_2-z_2)$ 
and finally,
\begin{equation}\label{eq:angeometry3}
v_2-z_2=z_1\cdot \lambda\cdot z_2^{-1}\cdot (u_1-z_1).
\end{equation}

If the lines $\Lambda,\Aline pu,\Aline {q}{v}$ are parallel, then the equations  (\ref{eq:angeom1}) and (\ref{eq:angeom2}) have no solutions, which happens if and only if 
$1+\lambda-\lambda\cdot z_2^{-1}\cdot (u_1-z_1)=0=1+\lambda-z_1^{-1}\cdot(u'_2-z_2)$ and only if the equality (\ref{eq:angeometry3}) holds.
\end{proof}

By analogy with Claim~\ref{cl:anal-geom1}, we can prove that the paraconcurrence of the lines $\Lambda,\Aline pa,\Aline qb$ implies the equality 
$$b_2-z_2=z_1\cdot \lambda\cdot z_2^{-1}\cdot (a_1-z_1).$$

Now assume that the  scalar $\overvector{zau}$ belongs to the center $Z(\IR_A)$ of the corps $\IR_A$. Then 
$$(a_1-z_1)\cdot\overvector{op}+(a_2-z_2)\cdot\overvector{oq}=\overvector{za}=\overvector{zau}\cdot\overvector{zu}=\overvector{zau}\cdot(u_1-z_1)\cdot\overvector{op}+\overvector{zua}\cdot(u_2-z_2)\cdot\overvector{oq},$$and hence
$$
\begin{aligned}
b_2-z_2&=z_1\cdot\lambda\cdot z_2^{-1}\cdot(a_1-z_1)=z_1\cdot\lambda\cdot z_2^{-1}\cdot\overvector{zau}\cdot (u_1-z_1)\\
&=\overvector{zau}\cdot z_1\cdot\lambda\cdot z_2^{-1}\cdot (u_1-z_1)=\overvector{zua}\cdot(v_2-z_2),
\end{aligned}
$$
see Claim~\ref{cl:anal-geom1}. 
Then $$
\begin{aligned}
\overvector{zbv}\cdot\overvector{zv}&=\overvector{zb}=(b_1-z_1)\cdot \overvector{op}+(b_2-z_2)\cdot\overvector{oq}=0\cdot\overvector{op}+\overvector{zua}\cdot(v_2-z_2)\cdot\overvector{oq}\\
&=\overvector{zua}\cdot(0\cdot\overvector{op}+(v_2-z_2)\cdot\overvector{oq})=\overvector{zua}\cdot((v_1-z_1)\cdot\overvector{op}+(v_2-z_2)\cdot\overvector{oq})=\overvector{zua}\cdot \overvector{zv}
\end{aligned}
$$
and hence $\overvector{zbv}=\overvector{zau}$. Since $$\overvector{ab}=\overvector{zb}-\overvector{za}=\overvector{zbv}\cdot\overvector{zv}-\overvector{zau}\cdot\overvector{zu}=\overvector{zau}\cdot\overvector{zv}-\overvector{zau}\cdot\overvector{zu}=\overvector{zau}\cdot(\overvector{zv}-\overvector{zu})=\overvector{zau}\cdot\overvector{uv},$$the lines $\Aline ab$ and $\Aline uv$ are parallel.
\end{proof}

\begin{problem} Find a geometric proof of Lemma~\textup{\ref{l:anal-geom}}.
\end{problem}

Now we are able to prove the promised equality $\Psi_A[\IF_X]=Z(\IR_A)$ for Desarguesian spaces.

\begin{theorem}\label{t:D=>FX=Z(RX)} If $X$ is a projective completion $X$ of a Desarguesian affine space $A$, then $$\Psi_A[\IF_X]=Z(\IR_A)=Z(\dddot A_{\Join}).$$
\end{theorem}

\begin{proof} Since the affine space $A$ is Desarguesian, $\IR_A=\dddot A_{\Join}$, by Theorem~\ref{t:Desargues<=>3Desargues}. Then $Z(\IR_A)=Z(\dddot A_{\Join})$. By Theorem~\ref{t:PsiA}(5), $\Psi_A[\IF_X]\subseteq Z(\dddot A_{\Join})=Z(\IR_A)$. To prove that $\Psi_A[\IF_X]=Z(\IR_A)$, take any scalar $\alpha\in Z(\IR_A)$ in the center of the corps $\IR_A$. Choose any points $i\in H$, $z\in A$ and $u\in\Aline zi\setminus\{z\}$. By Proposition~\ref{p:triple-exists}, there exists a point $a\in \Aline zu\cap A$ such that $\overvector{zau}=\alpha$. We claim that the quadruple $zaui$ is Pappian in the projective space $X$. Given any line projectivity $P$ in $X$ with $Pzui=zui$, we should prove that $P(a)=a$. Choose any point $j\in H\setminus\{i\}$ and any point $c\in \Aline ij\setminus\{i,j\}$. Consider the line perspectivity $C\defeq\{(x,y)\in \Aline zi\times\Aline zj:y\in \Aline xc\}$ and observe that $Czi=zj$. Let $v\defeq C(u)\in \Aline zj$ and observe that $Czui=zvj$ and $CP$ is a line projectivity such that $CPzui=zvj$. By Theorem~\ref{t:projectivity=3perspectivities}, $CP=R'R$ for some line perspectivities $R,R'$ in $X$. 

If the line $\Lambda\defeq\rng[R]=\dom[R']$ contains the point $z\in\dom[R]\cap\rng[R']$, then the line projectivity $CP=R'R$ is a line perspectivity, by Theorem~\ref{t:Desargues-perspective}. Since $CPui=vj$, the center of the perspectivity $CP$ is the unique common point $c$ of the lines $\Aline ij$ and $\Aline uv$. Therefore, the line perspectivities $CP$ and $C$ have the same perspectivity center  and the same domain and range. In this case $CP=C$ and $P=C^{-1}C$ is the identity map of the line $\Aline zi$. Therefore, $P(a)=a$.

\begin{picture}(200,220)(-160,-80)

\put(0,0){\color{blue}\line(1,0){105}}
\put(0,0){\color{teal}\line(0,1){105}}
\put(-15,120){\line(1,-1){165}}
\put(0,105){\color{teal}\line(-2,-5){70}}
\put(105,0){\color{blue}\line(-5,-2){175}}
\put(30,-30){\line(1,1){60}}
\put(90,30){\line(-8,-5){160}}
\put(90,30){\line(8,5){80}}
\put(-30,30){\line(4,1){200}}
\put(170,80){\line(-140,-110){140}}
\put(0,30){\color{red}\line(2,-1){150}}
\put(-30,30){\line(1,0){120}}
\put(30,-30){\line(-1,1){60}}
\multiput(150,-45)(-6,3.3){25}{\color{red}\line(-20,11){4}}

\put(0,0){\circle*{3}}
\put(-5,-8){$z$}
\put(60,0){\circle*{3}}
\put(48,-6){$u$}
\put(-30,30){\circle*{3}}
\put(-40,27){$q$}
\put(0,37.5){\circle*{3}}
\put(2,42){$b$}
\put(0,30){\circle*{3}}
\put(-8,22){$v$}
\put(90,30){\circle*{3}}
\put(85,33){$w$}
\put(30,-30){\circle*{3}}
\put(28,-38){$p$}
\put(105,0){\circle*{3}}
\put(109,-1){$i$}
\put(0,105){\circle*{3}}
\put(4,107){$j$}
\put(-70,-70){\circle*{3}}
\put(-77,-73){$o$}
\put(170,80){\circle*{3}}
\put(168,83){$d$}
\put(-40,-45){$\Lambda$}
\put(68.375,0){\circle*{3}}
\put(67,3){$a$}
\put(150,-45){\circle*{3}}
\put(154,-50){$c$}
\end{picture}

Next, assume that $z\notin\Lambda$. Let $p,q$ be the  centers of the line perspectivities $R,R'$, respectively. Then $p\notin \Aline zi\cup\Lambda$ and $q\notin \Lambda\cup\Aline zj$. Since $R'R(z)=CP(z)=z\notin\Lambda$, the point $z$ belongs to the line $\Aline pq$.
Consider the point $o\defeq R(i)\in \Lambda\cap\Aline pi$ and observe that $R'(o)=R'R(i)=CP(i)=C(i)=j$, which implies that $o\in \Aline j{q}$. Assuming that $o\in\Aline ij$, we conclude that $o\in\Aline ip\cap\Aline ij\cap \Aline j{q}=\{i\}\cap\{j\}$, which contradicts the choice of the point $j\ne i$. This contradiction shows that $o\notin\Aline ij$. Consider the points $w\defeq R(u)$, $d\defeq R(a)$, and $b\defeq R'(d)=R'R(a)$. Since $w\in \Aline pu\cap\Lambda\cap\Aline {q}{v}$, the lines $\Aline pu\cap A,\Lambda\cap A,\Aline{q}{v}\cap A$ in the affine space $A$ are paraconcurrent. Since $d\in\Aline pa\cap \Lambda\cap\Aline {q}{b}$, we obtain another triple of paraconcurrent lines $\Aline pa\cap A,\Lambda\cap A,\Aline{q}{b}$ in the affine space $A$. Applying Lemma~\ref{l:anal-geom}, we conclude that the lines $\Aline uv\cap A$ and $\Aline ab\cap A$ are parallel in the affine space $A$. Then  $\varnothing\ne\Aline uv\cap\Aline ab\cap \Aline ij\subseteq \Aline uv\cap \Aline ij=\{c\}$ and hence $c\in\Aline ab$ and $C(a)=b=R'R(a)=CP(a)$. Since the line perspectivity $C$ is injective, the equality $C(a)=CP(a)$ implies $P(a)=a$, witnessing that the line quadruple $zaui$ is Pappian. 

Then $\overvector{zaui}\in\IF_X$ and $\alpha=\overvector{zau}=\Psi_X(\overvector{zaui})\in \Psi_A[\IF_X]$, completing the proof of the equality $\Psi_A[\IF_X]=Z(\IR_A)=Z(\dddot A_{\Join})$.
\end{proof}

\begin{corollary}\label{c:FX!=01=>Moufang} If a projective space $X$ has $\IF_X\ne\{0,1\}$, then $X$ is everywhere Thalesian.
\end{corollary}

\begin{proof} Assume that $\IF_X\ne\{0,1\}$. We need to check that for every hyperplane $H\subset X$, the affine space $A\defeq X\setminus H$ is Thalesian. By Theorem~\ref{t:PsiA}(1), the function $\Psi_A:\IF_X\to \IR_A$ is injective, which implies that $\IR_A\ne\{0,1\}$. By Theorem~\ref{t:RXne01=>paraD}, the affine space $A$ is Thalesian, witnessing that the projective space $X$ is everywhere Thalesian.
\end{proof}

\begin{remark} Let $X$ be the non-Desarguesian Thalesian affine plane of order $9$ (considered in Section~\ref{s:J9xJ9}). For this plane, the corps of scalars $\IR_X$ is a $3$-element field $\{0,1,-1\}$ and hence $\IR_X=Z(\IR_X)=Z(\dddot X_{\Join})$, by Theorems~\ref{t:Papp<=>Des+RX} and \ref{t:Desargues<=>3Desargues}. Let $\overline X$ be the spread completion of the affine plane $X$. By Theorems~\ref{t:Desargues-completion} and \ref{t:Moufang-finite<=>}, the finite projective plane $\overline X$ is not Desarguesian and not Moufang, and hence $\IF_{\overline X}=\{0,1\}$, by Corollary~\ref{c:FX!=01=>Moufang}. Now consider the injective homomorphism $\Psi:\IF_{\overline X}\to \IR_X$ and observe that $\Psi[\IF_{\overline X}]=\{0,1\}\ne\{0,1,-1\}=\IR_X=Z(\dddot X_{\Join})$. This example shows that in general, Problem~\ref{q:center} has negative answer. 
\end{remark}

\section{Addition of proscalars} In this section we shall prove that for every projective space $X$ the operations of multiplication and $01$-involution on the set of finite proscalars $\IF_X$ induce a  unique operation of addition that turns $\IF_X$ into a field. Let us recall that the $01$-involution $1\mbox{-}:\bar\IF_X\to\bar\IF_X$ assigns to every proscalar $\overvector{zvui}$ the proscalar $1\mbox{-}\overvector{zvui}\defeq\overvector{uvzi}$.

\begin{theorem}\label{t:profield-addition} For every projective space $X$, there exists a unique binary operation $+:\IF_X\times \IF_X\to \IF_X$ such that $(\IF_X,+,\cdot,0,1)$ is a field, and  $1\mbox{-}\alpha=1-\alpha$ for every finite proscalar $\alpha\in\IF_X$.
\end{theorem}

\begin{proof} The uniqueness of the operation $+$ follows from Theorem~\ref{t:01*=>+}. It remains to prove its existence. 

If $\IF_X=\{0,1\}$, let $+:\IF_X\times\IF_X\to\IF_X$ be the operation of addition modulo 2. It is clear that $(\IF,+,\cdot,0,1)$ is a field and $1\mbox{-}\,\alpha=1-\alpha$ for every proscalar $\alpha\in\IF_X=\{0,1\}$.

So, assume that $\IF_X\ne\{0,1\}$. By Corollary~\ref{c:FX!=01=>Moufang} and Proposition~\ref{p:Steiner=>FX=01}, the projective space $X$ is Moufang and not Steiner. By Corollary~\ref{c:Avogadro-projective}, the $3$-long projective liner $X$ is $2$-balanced. Since $X$ is not Steiner, $X$ is $4$-long.

Fix any hyperplane $H$ in $X$, and consider the $3$-long subliner $A\defeq X\setminus H$, which is an affine space, by Theorem~\ref{t:affine<=>hyperplane} and Corollary~\ref{c:procompletion-rank}. By Theorem~\ref{t:PsiA}, the function $\Psi_A:\IF_X\to\IR_A$ is an injective homomorphism of multiplicative monoids $(\IF_X,\cdot)$ and $(\IR_A,\cdot)$ such that $\Psi_A(1\mbox{-}\alpha)=1-\Psi_A(\alpha)$ for every finite proscalar $\alpha\in\IF_A$. By Theorem~\ref{t:RX-corps}, $\IR_A$ is a corps and by Theorem~\ref{t:PsiA}, $\Psi_A[\IF_X]\subseteq Z(\dddot A_{\!\Join})\subseteq Z(\IR_A)$.

For every proscalar $\alpha\in\IF_X\setminus\{1\}$, define the proscalar $\mbox{-}\alpha$ by the formula
$$\mbox{-}\alpha\defeq (1\mbox{-}\,\alpha)\cdot(1\mbox{-}(1\mbox{-}\,\alpha)^{-1}).$$Consider the scalar $\beta\defeq\Phi_A(\alpha)\in\IR_A$ and observe that
\begin{multline*}
\Psi_A(\mbox{-}\alpha)=\Psi_A\big((1\mbox{-}\alpha)\cdot(1\mbox{-}(1\mbox{-}\alpha)^{-1})\big)=(1-\beta)\cdot(1-(1-\beta)^{-1})\\
=(1-\beta)-(1-\beta)\cdot(1-\beta)^{-1}=(1-\beta)-1=-\beta=-\Psi_A(\alpha).
\end{multline*}
The additive inverse $\mbox{-}1$ to $1$ can be defined as follows. Choose any proscalar $\alpha\in\IF_X\setminus\{0,1\}$ and put
$$\mbox{-}1\defeq \alpha^{-1}\cdot(\mbox{-}\alpha).$$ Then  
$$\Psi_A(\mbox{-1})=\Psi_A(\alpha^{-1}\cdot(\mbox{-}\alpha))=\Psi_A(\alpha^{-1})\cdot\Psi_A(\mbox{-}\alpha)=\Psi_A(\alpha)^{-1}\cdot (-\Psi_A(\alpha))=-1.$$
Therefore, for every proscalar $\alpha\in\IF_X$ we have the equality $$\Psi_A(\mbox{-}\alpha)=-\Psi_A(\alpha).$$

Now define the operation of addition $+:\IF_X\times\IF_X\to\IF_X$ by the formula 
$$\alpha+\beta=\begin{cases}
\beta&\mbox{if $\alpha=0$};\\
\alpha\cdot(1\mbox{-}(\alpha^{-1}\cdot (\mbox{-}\beta)))&\mbox{if $\alpha\ne 0$}.
\end{cases}
$$
If $\alpha\in\IF_X\setminus\{0\}$, then for every $\beta\in\IF_X$ we have the equalities
$$
\begin{aligned}
\Psi_A(\alpha+\beta)&=\Psi_A(\alpha\cdot(1\mbox{-}(\alpha^{-1}\cdot (\mbox{-}\beta)))=\Psi_A(\alpha)\cdot\Psi_A(1\mbox{-}(\alpha^{-1}\cdot(\mbox{-}\beta))=\Psi_A(\alpha)\cdot(1-\Psi_A(\alpha^{-1}\cdot(\mbox{-}\beta)))\\
&=
\Psi_A(\alpha)-\Psi_A(\alpha)\cdot\Psi_A(\alpha)^{-1}\cdot\Psi_A(\mbox{-}\beta)=\Psi_A(\alpha)-1\cdot(-\Psi_A(\beta))=\Psi_A(\alpha)+\Psi_B(\beta).
\end{aligned}
$$
Therefore, the injective function $\Psi_A$ preserves the operations of addition, multiplication and inversion (both additive and multiplicative). This implies that $\Psi_A[\IF_X]$ is a subcorps of the corps $\IR_X$ and hence $(\IF_X,+,\cdot,0,1)$ is a corps. Since the multiplication in $\IF_X$ is commutative, the corps $(\IF_X,+,\cdot,0,1)$ is a field.
\end{proof}
 



\begin{definition} For every projective space $X$, the set $\IF_X$ of finite proscalars endowed with the binary operations of addition and multiplication of scalars is called the \index{proscalar field}\defterm{proscalar field} of $X$. For a completely regular space $X$, its proscalar field $\IF_X$ is defined as the proscalar field $\IF_{\overline X}$ of the spread completion $\overline X$ of $X$.
\end{definition}

Theorems~\ref{t:PsiA}, \ref{t:D=>FX=Z(RX)}, \ref{t:profield-addition}, \ref{t:01*=>+} imply the following important fact.

\begin{corollary}\label{c:FX=Z(RA)} For every projective completion $X$ of a Desarguesian affine space $A$, the function $\Psi_A:\IF_X\to Z(\IR_A)$ is an isomorphism of the fields $\IF_X$ and $Z(\IR_A)$.
\end{corollary}

At the moment, for every projective space $X$, we have defined the multiplicative operation on the whole set of proscalars $\bar\IF_X$, but the additive operation is defined only on the set $\IF_X$ of finite proscalars. The additive operation can be extended to the whole set $\bar\IF_X$ letting
$$\mbox{$\infty+\infty\defeq 0\quad$ and $\infty+\alpha=\alpha+\infty\defeq \infty$ for all $\alpha\in \IF_X$.}$$ The set of proscalars $\bar \IF_X$ endowed with the extended addition and multiplication operations is called the \index{proscalar profield}\defterm{proscalar profield} of the projective space $X$. Proscalar profields of projective spaces motivate studying abstract profields and procorps, which will be done in Section~\ref{s:procorps}.

\section{Permutations of Pappian line quadruples}

In Corollary~\ref{c:Pappian-permutations} we have seen that for every Pappian line quadruple $abcd$ that consists of 4 distinct points of a projective space, every permutation of this quadruple remain a Pappian quadruple. In this section we shall show that among 24 possible permutations of the Pappian quadruple $abcd$, there are only six distinct proscalars.

 It suffices to understand the action of three transpositions on Pappian line quadruples, generating the 24-element permutation group.

\begin{theorem}\label{t:abcd=dcba}  For every Pappian line quadruple $zvui$ with $v\notin\{z,i\}$ in a projective space $X$, the quadruples $vuzi$, $zuvi$ and $iuvz$ are Pappian and
$$\overvector{uvzi}=1-\overvector{zvui},\quad \overvector{zuvi}=(\overvector{zvui})^{-1},\quad\overvector{iuvz}=\overvector{zvui}$$
in the proscalar field $\IF_X$.
\end{theorem}

\begin{proof} By Corollary~\ref{c:Pappian-permutations}, the line quadruples $vuzi$, $zuvi$ and $iuvz$ are Pappian. The equality $\overvector{uvzi}=1-\overvector{zvui}$ follows from the definition of the addition of proscalars in Theorem~\ref{t:profield-addition} and the definition of $01$-involution on $\IF_X$. By Corollary~\ref{c:FX*-group}, $\overvector{zuvi}\cdot\overvector{zvui}=\overvector{zuui}=1$ and hence $\overvector{zuvi}=\overvector{zvui}^{-1}$ in the multiplicative group $\IF_X^*$ of the field $\IF_X$. It remains to prove the equality $\overvector{iuvz}=\overvector{zvui}$. 

By Lemma~\ref{l:1/x} (proved below), there exists a line projectivity $P$ in $X$ such that $Pzui=iuz$ and for the point $w\defeq P(v)$ we have $\overvector{zwui}=(\overvector{zvui})^{-1}$. Then $$iuvz\projeq Piuvz=zuwi\in\overvector{zuwi}=(\overvector{zwui})^{-1}=(\overvector{zvui}^{-1})^{-1}=\overvector{zvui}$$ and hence $\overvector{iuvz}=\overvector{zvui}$.
\end{proof}

\begin{lemma}\label{l:1/x} For every line $L$ in a projective space $X$ and every distinct points $z,u,i\in L$,  there exists a line projectivity $P:L\to L$ such that $Pzui=iuz$ and for every point $v\in L$ with $\overvector{zvui}\in\bar \IF_X$ and its image $w\defeq P(v)$, we have $\overvector{zwui}=(\overvector{zvui})^{-1}$ in $\bar\IF_X$.
\end{lemma}

\begin{proof}
Choose any point $o\in X\setminus L$. Take any point $j\in\Aline uo\setminus\{u,o\}$ and consider the projective plane $\Pi\defeq\overline{\{z,u,o\}}$ in $X$ and the subliner $A\defeq\Pi\setminus\Aline ij$ in $\Pi$.  Consider the unique points $c\in \Aline jz\cap\Aline io\subseteq A$ and $k\in \Aline zo\cap\Aline ij$, and the line perspectivities
$$C\defeq \{(x,y)\in \Aline zi\times\Aline uo:y\in\Aline xc\}\quad\mbox{and}\quad K\defeq\{(x,y)\in \Aline uo\times \Aline iz:y\in\Aline xk\}$$with centers at the points $c$ and $k$. We claim that the line projectivity $P\defeq KC$ has the required properties. First observe that $Pzui=KCzui=Kjuo=iuz$.

\begin{picture}(200,145)(-160,-10)

\put(0,0){\line(1,0){120}}
\put(0,0){\line(0,1){120}}
\put(0,0){\line(-1,0){60}}
\put(-60,0){\line(1,2){60}}
\put(120,0){\line(-1,1){120}}
\put(120,0){\line(-4,1){160}}
\put(-60,0){\line(2,1){120}}

\put(0,0){\circle*{3}}
\put(-2,-7){$u$}
\put(-60,0){\circle*{3}}
\put(-65,-8){$z$}
\put(120,0){\circle*{3}}
\put(123,-2){$i$}
\put(0,30){\circle*{3}}
\put(2,22){$o$}
\put(0,120){\circle*{3}}
\put(-2,125){$j$}
\put(-40,40){\circle*{3}}
\put(-48,40){$c$}
\put(60,60){\circle*{3}}
\put(62,62){$k$}
\end{picture}

Now take any point $v\in L$ with $\overvector{zvui}\in\bar\IF_X$ and consider the points $x\defeq C(v)$ and $w\defeq K(x)=P(v)\in L$. We have to prove that $\overvector{zwui}=\overvector{zvui}^{-1}$. 

If $v=z$, then $w=P(v)=P(z)=i$ and hence
$\overvector{zwui}=\overvector{ziui}=\infty=0^{-1}=\overvector{zzui}^{-1}=\overvector{zvui}^{-1}$.

If $v=u$, then $w=P(v)=P(u)=u$ and hence
$\overvector{zwui}=\overvector{zuui}=1=1^{-1}=\overvector{zuui}^{-1}=\overvector{zvui}^{-1}$.

If $v=i$, then $w=P(v)=P(i)=z$ and hence
$\overvector{zwui}=\overvector{zzui}=0=\infty^{-1}=\overvector{ziui}^{-1}=\overvector{zvui}^{-1}$. 

So, assume that $v\notin\{z,u,i\}$.  In this case, $|X|_2=|L|\ge\{z,v,u,i\}|=4$  and the affine subliner $A=\Pi\setminus\Aline ij$ of $\Pi$ is $3$-long and hence $A$ is an affine plane. 

Taking into account that $Ciuvz=ouxj$ and $Kouxj=zuwi$, we conclude that
$$\overvector{iuvz}=\overvector{ouxj}=\overvector{zuwi}.$$ 
Since $zvui$ is a Pappian line quadruple in $X$, so are the quadruples $iuvz$ and $zuwi\projeq iuvz$, by Corollary~\ref{c:Pappian-permutations}. By Theorem~\ref{t:PsiA}(5), $\overvector{uvz},\overvector{zuw},\overvector{oux}\in \Psi_A[\IF_X]\subseteq Z(\IR_A)$ and hence the line triples $uvz,zuw,oux$ are Desarguesian in the affine plane $A$. 
Let $p\in A$ be the unique common point of the lines $\Aline zo$ and $\Aline wj$ in the projective plane $\Pi$.

\begin{picture}(200,190)(-180,-50)

\put(0,0){\line(1,0){120}}
\put(0,0){\line(0,1){120}}
\put(0,0){\line(-1,0){60}}
\put(0,0){\line(0,-1){40}}
\put(-60,0){\line(1,2){60}}
\put(120,0){\line(-1,1){120}}
\put(120,0){\line(-4,1){160}}
\put(-40,40){\line(1,-2){40}}
\put(-60,0){\line(2,1){120}}
\put(0,-40){\line(3,5){60}}
\put(0,120){\line(24,-120){24}}

\put(24,0){\circle*{3}}
\put(24,-8){$w$}
\put(0,0){\circle*{3}}
\put(2,-8){$u$}
\put(-60,0){\circle*{3}}
\put(-65,-8){$z$}
\put(-20,0){\circle*{3}}
\put(-25,-8){$v$}
\put(120,0){\circle*{3}}
\put(123,-2){$i$}
\put(0,30){\circle*{3}}
\put(2,22){$o$}
\put(0,120){\circle*{3}}
\put(-2,125){$j$}
\put(-40,40){\circle*{3}}
\put(-48,40){$c$}
\put(60,60){\circle*{3}}
\put(62,62){$k$}
\put(16.4,38.2){\circle*{3}}
\put(18,43){$p$}
\put(0,-40){\circle*{3}}
\put(-3,-49){$x$}
\end{picture}

\begin{claim} $\overvector{zvu}=\overvector{zop}=\overvector{zuw}\in \Psi_A[\IF_X]$.
\end{claim}

\begin{proof} Taking into account that the lines $\Aline ou\cap A$ and $\Aline wp\cap A$ are parallel in  the affine plane $A$, we conclude that $\overvector{zop}=\overvector{zuw}\in\Psi_A[\IF_X]\subseteq Z(\IR_A)$ and hence the line triple $zop$ in the affine plane $A$ is Desarguesian.
Taking into account that the lines $\Aline vu\cap A$ and $\Aline co\cap A$ are parallel in the affine plane $A$, we conclude that $\overvector{xvc}=\overvector{xuo}\in\Psi_{A}[\IF_X]\subseteq Z(\IR_A)$ and hence the line triple $xvz$ in the affine plane $A$ is Desarguesian. 

The equality $\overvector{zvu}=\overvector{zop}$ will follow as soon as we show that the lines $\Aline vo\cap A$ and $\Aline up\cap A$ are parallel in the affine space $A$. It follows from $\overvector{zvu}\in Z(\IR_A)\setminus\{0,1\}$ that $\IR_A\ne\{0,1\}$ and hence the affine plane $A$ is Thalesian, by Theorem~\ref{t:RXne01=>paraD}. By Theorem~\ref{t:paraD=>RX-module}, the set of vectors $\overvector A$ in $A$ is an $\IR_A$-module and also $\IF_A$-module over the subfield $$\IF_A\defeq \Psi_A[\IF_X]\subseteq Z(\IR_A)$$ of the corps $\IR_A$.
Endowed with the family of lines $\{a+\IF_A{\cdot} \vec{\boldsymbol v}:x\in A,\;\vec{\boldsymbol v}\in\overvector A\}$, the liner $A$ is Desarguesian, by Theorem~\ref{t:R-module=>Des-aff-reg}. Taking into account that the field $\IF_A$ contains the scalars $\overvector{zuw},\overvector{vuz}, \overvector{zop}, \overvector{oux},\overvector{xvc}$, we conclude that the set $\{z,v,u,w,o,p,x,c\}$ has rank $2$ in the $\IF_A$-module $A$ and hence it is contained in some plane $\Pi_A$ in the $\IF_A$-module $A$. 

Since $j\notin A$, the lines $z+\IF_A{\cdot}(z-c)$, $o+\IF_A{\cdot}(u-o)$ and $p+\IF_A{\cdot}(w-p)$ are parallel in the plane $\Pi_A$. Also $i\notin A$ imply $(o+\IF_A{\cdot}(c-o))\cap(z+\IF_A{\cdot}(z-u))\subseteq (\Aline co\cap A)\cap(\Aline zu\cap A)=\varnothing$. So, the lines $o+\IF_A{\cdot}(c-o)$ and $z+\IF_A{\cdot}(z-u)$ are parallel in $\Pi_A$. By analogy we can show that $k\notin A$ implies that the lines $o+\IF_A{\cdot}(z-o)$ and $x+\IF_A{\cdot}(w-x)$ are parallel in $P$.

\begin{picture}(150,200)(-150,-85)
{\linethickness{0.75pt}
\put(0,60){\color{orange}\line(1,0){60}}
\put(60,0){\color{orange}\line(1,0){30}}
\put(0,60){\color{blue}\line(2,-3){60}}
\put(60,0){\color{cyan}\line(0,1){60}}
\put(60,0){\color{cyan}\line(0,-1){30}}
\put(90,0){\color{cyan}\line(0,1){90}}
\put(90,0){\color{cyan}\line(0,-1){45}}
\put(40,0){\color{red}\line(1,3){20}}
\put(60,0){\color{red}\line(1,3){30}}
\put(60,0){\color{blue}\line(2,-3){30}}
\put(0,0){\color{cyan}\line(0,1){60}}
\put(0,0){\color{green}\line(1,1){60}}
\put(60,-30){\color{green}\line(1,1){30}}
}
\put(60,60){\line(1,1){30}}

\put(0,0){\line(1,0){60}}
\put(0,60){\line(1,-1){120}}

\put(0,0){\line(2,-1){120}}
\put(60,60){\line(1,-2){60}}

\put(0,0){\circle*{3}}
\put(-8,-3){$z$}
\put(40,0){\circle*{3}}
\put(35,-8){$v$}
\put(60,0){\circle*{3}}
\put(52,-8){$u$}
\put(90,0){\circle*{3}}
\put(93,-2){$w$}
\put(120,-60){\circle*{2}}
\put(0,60){\circle*{3}}
\put(-5,63){$c$}
\put(60,60){\circle*{3}}
\put(57,63){$o$}
\put(90,90){\circle*{3}}
\put(93,93){$p$}
\put(60,-30){\circle*{3}}
\put(57,-38){$x$}
\put(90,-45){\circle*{3}}
\put(86,-53){$q$}

\end{picture}

By the Inverse Desargues Theorem~\ref{t:ID1}, the lines $x+\IF_A{\cdot}(z-x)$, $u+\IF_A{\cdot}(c-u)$ and $o+\IF_A{\cdot}(w-o)$ are paraconcurrent in the plane $\Pi_A$. Since $(x+\IF_A{\cdot}(z-x))\cap(o+\IF_A{\cdot}(x-o))=\{x\}$ and $(o+\IF_A{\cdot}(x-o))\parallel (w+\IF_A{\cdot}(p-w))$, there exists a unique point $q\in (z+\IF_A{\cdot}(x-z))\cap(w+\IF_A{\cdot}(p-w))$. Since the triangles $xoc$ and $qwu$ are paraperspective,   $(q+\IF_A{\cdot}(w-q))\cap(x+\IF_A{\cdot}(o-x))=\varnothing=(w+\IF_A{\cdot}(u-w))\cap(o+\IF_A{\cdot}(c-o))$ imply $(c+\IF_A{\cdot}(x-c))\parallel(u+\IF_A{\cdot}(q-u))$. Since the triangles $vxo$ and $uqp$ are perspective from the point $z$,   $(v+\IF_A{\cdot}(x-v))\cap(u+\IF_A{\cdot}(q-u))=\varnothing=(x+\IF_A{\cdot}(o-x))\cap(q+\IF_A{\cdot}(p-q))$ imply $(v+\IF_A{\cdot}(o-v))\cap(u+\IF_A{\cdot}(p-u))=\varnothing$ and hence the vectors $\overvector{ov}$ and $\overvector{pu}$ are parallel and so are the lines $\Aline ov$ and $\Aline up$. Then  $\overvector{zvu}=\overvector{zop}=\overvector{zuw}$.
\end{proof}

It follows from $\overvector{zuw}=\overvector{zvu}$ that $\overvector{zuwi}=\overvector{zvui}$, and finally, $\overvector{iuvz}=\overvector{zuwi}=\overvector{zvui}$.
\end{proof}

\begin{problem} Is $\overvector{abcd}=\overvector{dcba}$ for every line quadruple $abcd$ consisting of four distinct points in a projective space?
\end{problem}

\begin{lemma}\label{l:lambda6} For every element $\lambda\notin\{0,1\}$ of a field $F$, we have the equality 
$$(1-\lambda^{-1})^{-1}=1-(1-\lambda)^{-1}.$$
\end{lemma}

\begin{proof}
Observe that $$1-(1-\lambda)^{-1}=(1-\lambda)^{-1}((1-\lambda)-1)=(\lambda-1)^{-1}\cdot\lambda=(1-\lambda^{-1})^{-1}\cdot\lambda^{-1}\cdot\lambda=(1-\lambda^{-1})^{-1}.
$$
\end{proof}

Theorem~\ref{t:abcd=dcba} and Lemma~\ref{l:lambda6} imply the following interesting corollary.

\begin{corollary} Let $a,b,c,d$ be four distinct points in a projective space $X$. If the quadruple $abcd$ is Pappian, then for the proscalar $\lambda\defeq\overvector{abcd}\in \bar \IF_X\setminus\{0,1,\infty\}$ we have the following identities:
$$
\begin{aligned}
&\overvector{abcd}=\overvector{badc}=\overvector{cdab}=\overvector{dcba}=\lambda\\
&\overvector{acbd}=\overvector{bdac}=\overvector{cadb}=\overvector{dbca}=\lambda^{-1}\\
&\overvector{cbad}=\overvector{dabc}=\overvector{adcb}=\overvector{bcda}=1-\lambda\\
&\overvector{bcad}=\overvector{adbc}=\overvector{dacb}=\overvector{bdca}=1-\lambda^{-1}\\
&\overvector{bacd}=\overvector{abdc}=\overvector{dcab}=\overvector{bcda}=(1-\lambda^{-1})^{-1}\\
&\overvector{cabd}=\overvector{dbac}=\overvector{acdb}=\overvector{bdca}=(1-\lambda)^{-1}
\end{aligned}
$$
in the field $\IF_X$.
\end{corollary}

%


\chapter{Automorphisms of projective spaces}

\section{Procorps and profields}\label{s:procorps}

\begin{definition}\label{d:procorps} A \index{procorps}\defterm{procorps} is a set $F$ endowed with two binary operations $+,\cdot: F\times F\to F$ and three distinct constants $0,1,\infty\in F$ that satisfy the following nine axioms:
\begin{enumerate}
\item $\forall x,z\in F\;\forall y\in F\setminus\{\infty\}\;\;\big(x+(y+z)=(x+y)+z\big)$;
\item $\forall x\in F\;\;(x+0=x=0+x)$;
\item $\forall x\in F\;\exists y\in F\;(x+y=0=y+x)$;
\item $\forall x,z\in F\;\forall y\in F\setminus\{0,\infty\}\;\;(x\cdot(y\cdot z)=(x\cdot y)\cdot z)$;
\item $\forall x\in F\; (x\cdot 1=x=1\cdot x)$;
\item $\forall x\in F\;\exists y\in F\;\;(x\cdot y=1=y\cdot x)$;
\item $\forall a\in F\setminus\{0,\infty\}\;\forall x,y\in F\;\;\big(a\cdot(x+y)=a\cdot x+a\cdot y)$;
\item $\forall x,y\in F\;\forall b\in F\setminus\{0,\infty\}\;\;\big((x+y)\cdot b=x\cdot b+y\cdot b\big)$;
\item $0\cdot 0=0$,  $\infty\cdot\infty=\infty$ and $1+\infty=\infty=\infty+1$.
\end{enumerate}
A procorps $F$ is called a \index{profield}\defterm{profield} if $x\cdot y=y\cdot x$ for all elements $x,y\in F$.
\end{definition} 

\begin{example} For every projective space $X$ its proscalar profield $\bar F_X$ is a profield.
\end{example}

\begin{exercise} Given any corps $F$, choose any element $\infty\notin F$,  consider the set $\bar F\defeq F\cup\{\infty\}$, and extend the operations of addition and multiplication from $F$ to $\bar F$ letting
$$
\begin{aligned}
&\infty+\infty=0,\quad\infty\cdot 0=1=0\cdot\infty,\\
&\forall x\in \bar F\setminus\{\infty\}\;\;(x+\infty=\infty=\infty+x),\\
&\forall x\in \bar F\setminus\{0\}\;\;(x\cdot\infty=\infty=\infty\cdot x).
\end{aligned}
$$ 
Check that $\bar F$ is a procorps. This procorps is called the \index{projective $\infty$-extension}\defterm{projective $\infty$-extension} of the corps $F$. If $F$ is a field, then its projective $\infty$-extension is a profield.
\end{exercise}

The following theorem shows that procorps are exactly projective $\infty$-extensions of corps.

\begin{theorem} For every procorps $\bar F$, the set $F\defeq \bar F\setminus\{\infty\}$ endowed with the induced operations of addition and multiplication is a corps and $\bar F$ is the projective $\infty$-extension of $F$.
\end{theorem} 

\begin{proof} The proof of this theorem is divided into a series of claims, which gradually reveal the structure of the procorps $\bar F$. In the proofs we shall refer to the axioms of a procorps as numbered in Definition~\ref{d:procorps}.

\begin{claim}\label{cl:+cancellative} For every $x\in \bar F\setminus\{\infty\}$ and $y,z\in \bar F$, the equality $x+y=x+z$ implies $y=z$.
\end{claim}

\begin{proof} Assume that $x+y=x+z$. By the axiom (3), there exists an element $x'\in \bar F$ such that $x'+x=0$. Since $x\ne\infty$, we can apply the axioms (2) and (1) and conclude that
$$y=0+y=(x'+x)+y=x'+(x+y)=x'+(x+z)=(x'+x)+z=0+z=z.$$
\end{proof}

By analogy we can prove 

\begin{claim}\label{cl:+cancellative2} For every $x,y\in\bar F$ and $z\in \bar F\setminus\{\infty\}$, the equality $x+z=y+z$ implies $x=y$.
\end{claim}

\begin{claim}\label{cl:x+x=x=>} Every element $x=x+x$ in $\bar F$ is equal to $0$ or $\infty$.
\end{claim}

\begin{proof}  Assuming that $x+x=x\ne\infty$, we can apply the axiom (2) and Claim~\ref{cl:+cancellative} and conclude that $x+x=x=x+0$ implies $x=0$.
\end{proof}

\begin{claim}\label{cl:0x0} For every $x\in \bar F\setminus\{\infty\}$, the elements $x\cdot 0$ and $0\cdot x$ belong to the set $\{0,\infty\}$.
\end{claim}

\begin{proof} If $x=0$, then $x\cdot 0=0\cdot x=0\cdot 0=0$, by the axiom (9). So, assume that $x\ne 0$. By the axioms (2) and (7), $$x\cdot 0=x\cdot(0+0)=x\cdot 0+x\cdot 0$$and hence $x\cdot 0\in\{0,\infty\}$, by Claim~\ref{cl:x+x=x=>}. By analogy we can prove that $0\cdot x\in\{0,\infty\}$. 
\end{proof}

\begin{claim}\label{cl:infty0=1} $0\cdot\infty=1=\infty\cdot 0$.
\end{claim}

\begin{proof}  By the axiom (6), there exists $x\in \bar F$ such that $0\cdot x=1$. Assuming that $x\ne\infty$, we can apply Claim~\ref{cl:0x0} and conclude that $1=0\cdot x\in\{0,\infty\}$, which is a contradiction showing that $x=\infty$ and hence $0\cdot \infty=0\cdot x=1$. By analogy we can prove that $\infty\cdot 0=1$.
\end{proof}









\begin{claim}\label{cl:a0=infty=>ainfty!=1} For every $a\in \bar F\setminus\{0,\infty\}$, if $a\cdot 0=\infty$, then $a\cdot\infty\ne 1$.
\end{claim}

\begin{proof} To derive a contradiction, assume that $a\cdot 0=\infty$ and $a\cdot\infty=1$. By the axioms (2), (7), and (9), 
$$1=a\cdot \infty=a\cdot(0+\infty)=a\cdot 0+a\cdot \infty=1+\infty=\infty,$$  which is a desired contradiction.
\end{proof} 

\begin{claim}\label{cl:a0=0} For every $\forall a\in \bar F\setminus\{\infty\}$ we have $a\cdot 0=0=0\cdot a$.
\end{claim}

\begin{proof}  To derive a contradition, assume that $a\cdot 0\ne 0$ for some $a\in \bar F\setminus\{\infty\}$.  By Claim~\ref{cl:0x0}, $a\cdot 0=\infty$. The axiom (9) ensures that $a\ne 0$. By the axiom (6), there exists $b\in \bar F$ such that $a\cdot b=1=b\cdot a$. It follows from $a\cdot b=1\ne\infty =a\cdot 0$ that $b\ne 0$. Claim~\ref{cl:a0=infty=>ainfty!=1}, implies $b\ne \infty$. 
By the axioms (4) and (5),
$$0=1\cdot 0=(b\cdot a)\cdot 0=b\cdot(a\cdot 0)=b\cdot\infty,$$
and by the axioms (3) and (8),
$$0+0=0=b\cdot\infty=b\cdot(0+\infty)=b\cdot 0+b\cdot\infty=b\cdot 0+0$$and hence $b\cdot 0=0$. Applying the axioms (4) and (5), we conclude that 
$$\infty=a\cdot 0=a\cdot(b\cdot 0)=(a\cdot b)\cdot 0=1\cdot 0=0,$$which is a desired contradiction showing that $a\cdot 0=0$. By analogy we can prove that $0\cdot a=0$.
\end{proof}

\begin{claim}\label{cl:a.i=1=>a=0} An element $a\in\bar F$ is equal to $0$ if and only if $a\cdot\infty=1$ or $\infty\cdot a=1$.
\end{claim}

\begin{proof} The ``only if'' part follows from Claim~\ref{cl:infty0=1}. To prove the ``if'' part, assume that $a\cdot\infty=1$.  Since $\infty\cdot\infty=\infty\ne 1=a\cdot\infty$, the element $a$ is not equal to $\infty$. Assuming that $a\ne 0$, we can apply the axiom (4), (5), and Claim~\ref{cl:a0=0} to obtain a contradiction $$0=0\cdot 1=0\cdot(a\cdot \infty)=(0\cdot a)\cdot \infty=0\cdot \infty=1,$$
showing that $a=0$. By analogy we can prove that $\infty\cdot a=1$ implies $a=0$.
\end{proof} 

\begin{claim}\label{cl:a.infty=infty} For every $a\in \bar F\setminus\{0,\infty\}$ we have $a\cdot \infty=\infty=\infty\cdot a$. 
\end{claim}

\begin{proof} To derive a contradiction, assume that for some element $a\in \bar F\setminus\{0\}$, the element $b:=a\cdot\infty$ is not equal to $\infty$. By the axiom (6), there exists $c\in \bar F$ such that $c\cdot b=1=b\cdot c$.  By Claim~\ref{cl:a0=0}, $c\ne 0$. Assuming that $c=\infty$, we conclude that $1=c\cdot b=\infty\cdot b$ and hence $b=0$, by Claim~\ref{cl:a.i=1=>a=0}. Then $a\cdot\infty=b=0$ and 
$$0=0\cdot 0=0\cdot(a\cdot\infty)=(0\cdot a)\cdot\infty=0\cdot\infty=1,$$
by the axiom (4) and Claims~\ref{cl:a0=0} and \ref{cl:infty0=1}.
This contradiction shows that $c\ne \infty$.  By Claim~\ref{cl:a.i=1=>a=0}, the equality
$$1=c\cdot b=c\cdot(a\cdot \infty)=(c\cdot a)\cdot\infty$$ implies $c\cdot a=0$. By the axioms (5), (4) and Claim~\ref{cl:a0=0},
$a=1\cdot a=(b\cdot c)\cdot a=b\cdot (c\cdot a)=b\cdot 0=0$, which contradicts the choice of $a$. By analogy we can prove that $\infty\cdot a=\infty$. 
\end{proof}

\begin{claim}\label{cl:x.yinF} For every $x,y\in \bar F\setminus\{\infty\}$, we have $x\cdot y\in \bar F\setminus\{\infty\}$.
\end{claim}

\begin{proof} To derive a contradiction, assume that $x\cdot y=\infty$ for some $x,y\in \bar F\setminus\{\infty\}$. Claim~\ref{cl:a0=0} implies that $x\ne 0\ne y$. By the axiom (5), there exists an element $z\in\bar F$ such that $z\cdot x=1$. By the axiom (4) and  Claim~\ref{cl:a.infty=infty}, we have
$$y=1\cdot y=(z\cdot x)\cdot y=z\cdot(x\cdot y)=z\cdot\infty=\infty,$$
which contradicts the choice of $y\in \bar F\setminus\{\infty\}$.
\end{proof}

\begin{claim}\label{cl:x+i=i} For every $x\in \bar F\setminus \{\infty\}$, we have $x+\infty=\infty=\infty+x$.
\end{claim}

\begin{proof} If $x=0$, then $x+\infty=0+\infty=\infty$ by the axiom (2). So, assume that $x\ne 0$. By Claim~\ref{cl:a.infty=infty} and the axiom (7),
$$x+\infty=x\cdot 1+x\cdot\infty=x\cdot(1+\infty)=x\cdot\infty=\infty.$$
On the other hand, Claim~\ref{cl:a.infty=infty} and the axiom (8) implies
$$\infty+x=\infty\cdot x+1\cdot x=(\infty+1)\cdot x=\infty\cdot x=\infty.$$
\end{proof}

\begin{claim}\label{cl:i+i=0} $\infty+\infty=0$.
\end{claim}

\begin{proof} By the axiom (3), there exists an element $x\in \bar F$ such that $x+\infty=0$. Claim~\ref{cl:x+i=i} ensures that $x=\infty$.
\end{proof}

\begin{claim}\label{cl:x+yinF} For every $x,y\in \bar F\setminus\{\infty\}$, $x+y\in\bar F\setminus\{\infty\}$.
\end{claim}

\begin{proof} To derive a contradiction, assume that $x+y=\infty$ for some elements $x,y\in\bar F\setminus\{\infty\}$. The axiom (2) ensures that $x\ne 0\ne y$. By the axiom (3), there exists an element $z\in\bar F$ such that $y+z=0$. Assuming that $z=\infty$, we can apply Claim~\ref{cl:x+i=i} and obtain a contradiction: $0=y+z=y+\infty=\infty$. Therefore, $z\ne\infty$ and $\infty+z=\infty$, by Claim~\ref{cl:x+i=i}. Applying the axiom (1) we obtain a contradiction $x=x+0=x+(y+z)=(x+y)+z=\infty+z=\infty$ with the choice of $x$.
\end{proof}

\begin{claim}\label{cl:x+y=y+x} For every $x,y\in\bar F$ we have $x+y=y+x$.
\end{claim}

\begin{proof} If $\infty\in\{x,y\}$, then $x+y=y+x$, by Claim~\ref{cl:x+i=i} and \ref{cl:i+i=0}. So, assume that $x,y\in \bar F\setminus\{\infty\}$. By Claim~\ref{cl:x+yinF}, $x+y\ne \infty\ne 1+1$. By Claim~\ref{cl:a0=0} and the axioms (7), (5), (1) and (8),
$$(1+1)\cdot (x+y)=1\cdot(x+y)+1\cdot(x+y)=(x+y)+(x+y)=x+(y+(x+y))=x+((y+x)+y)$$
and
$$(1+1)\cdot (x+y)=(1+1)\cdot x+(1+1)\cdot y=(1\cdot x+1\cdot x)+(1\cdot y+1\cdot y)=(x+x)+(y+y)=x+((x+y)+y).$$
Therefore, $x+((y+x)+y)=x+((x+y)+y)$. Applying Claims~\ref{cl:+cancellative} and \ref{cl:+cancellative2}, we conclude that $(y+x)+y=(x+y)+y$ and $y+x=x+y$.
\end{proof}

Finally, we are able to prove that the procorps $\bar F$ is a projective $\infty$-extension of the corps $F\defeq \bar F\setminus\{\infty\}$. Claims~\ref{cl:x+yinF} and \ref{cl:x.yinF} ensure that for any elements $x,y\in F$ we have $\{x+y,x\cdot y\}\subseteq F$. By the axioms (1) and (4), the induced operations of addition and multiplication on $F$ are associative. Claim~\ref{cl:x+y=y+x} ensures that the operation of addition on $F$ is commutative and the axioms (7), (8), (9), and Claim~\ref{cl:a0=0} imply that the distributive laws on $F$ also hold. The axiom (3) and Claim~\ref{cl:x+i=i} ensures that every element $x\in F$ has an additive inverse in $F$. The axiom (6) and Claim~\ref{cl:a0=0} imply that every non-zero element $a\in F$ has left and right multiplicative inverses in $F$.  Therefore, $F$ is a corps. 

Claims~\ref{cl:infty0=1}, \ref{cl:i+i=0}, \ref{cl:x+i=i}, \ref{cl:a.infty=infty} ensure that 
$$0\cdot\infty=1,\quad \infty+\infty=0,\quad\mbox{and}\quad
\forall x\in \bar F\setminus\{0\}\;\forall y\in \bar F\setminus\{\infty\}\;\;(\infty =x\cdot \infty=\infty\cdot x=y+\infty),
$$
which means that $\bar F$ is a projective $\infty$-extension of the corps $F$.
\end{proof}

\begin{example} The $3$-element set $\bar F_3\defeq\{0,1,\infty\}$, endowed with the addition and multiplication operations, defined by the following addition and multiplication tables
$$
\begin{array}{c|ccc}
+&0&1&\infty\\
\hline
0&0&1&\infty\\
1&1&\infty&0\\
\infty&\infty&0&1
\end{array}\qquad\mbox{and}\qquad
\begin{array}{c|ccc}
\cdot&0&1&\infty\\
\hline
0&0&0&1\\
1&0&1&\infty\\
\infty&1&\infty&\infty
\end{array}
$$
satisfies all the axioms of a profield except for the equality $1+\infty=\infty$ in the last axiom. The additive operation on $\bar F_3$ is associative, so $(\bar F_3,+)$ is a group.
\end{example}

\begin{example} The $4$-element set $\bar F_4\defeq\{0,1,\infty,1{+}\infty\}$, endowed with the addition and multiplication operations, defined by the following addition and multiplication tables
$$
\begin{array}{c|cccc}
+&0&1&\infty&1{+}\infty\\
\hline
0&0&1&\infty&1{+}\infty\\
1&1&\infty&1{+}\infty&0\\
\infty&\infty&1{+}\infty&0&1\\
1{+}\infty&1{+}\infty&0&1&\infty
\end{array}\qquad\mbox{and}\qquad
\begin{array}{c|cccc}
\cdot&0&1&\infty&1{+}\infty\\
\hline
0&0&0&1&0\\
1&0&1&\infty&1{+}\infty\\
\infty&1&\infty&\infty&\infty\\
1{+}\infty&0&1{+}\infty&\infty&1
\end{array}
$$
satisfies all the axioms of a profield except for the equality $1+\infty=\infty$ in the last axiom. The additive operation on $\bar F_4$ is associative, so $(\bar F_4,+)$ is a cyclic group.
\end{example}

\begin{example} The $4$-element set $\bar B_4\defeq\{0,1,\infty,1{+}\infty\}$, endowed with the addition and multiplication operations, defined by the following addition and multiplication tables
$$
\begin{array}{c|cccc}
+&0&1&\infty&1{+}\infty\\
\hline
0&0&1&\infty&1{+}\infty\\
1&1&0&1{+}\infty&\infty\\
\infty&\infty&1{+}\infty&0&1\\
1{+}\infty&1+\infty&\infty&1&0
\end{array}\qquad\mbox{and}\qquad
\begin{array}{c|cccc}
\cdot&0&1&\infty&1{+}\infty\\
\hline
0&0&0&1&0\\
1&0&1&\infty&1{+}\infty\\
\infty&1&\infty&\infty&\infty\\
1{+}\infty&0&1{+}\infty&\infty&1
\end{array}
$$
satisfies all the axioms of a profield except for the equality $1+\infty=\infty$ in the last axiom. The additive operation on $\bar B_4$ is associative, and $(\bar B_4,+)$ is a Boolean group.
\end{example}

By analogy with Theorem~\ref{t:01*=>+} we can prove the following helpful fact.

\begin{theorem} The operation of addition in a procorps $F$ is uniquely determined by  the unary operation of $01$-involution $1\mbox{-}:F\to F$, $1\mbox{-}:x\mapsto 1-x$ and the operation of multiplication in the group $F^*\defeq F\setminus\{0,\infty\}$.
\end{theorem}

\section{M\"obius transformations of profields}\label{s:Mobious}

In this section we study \index[person]{M\"obius}M\"obius\footnote{{\bf August Ferdinand M\"obius} (1790 -- 1868) was a German mathematician and theoretical astronomer. He is best known for his discovery of the M\"obius strip, a non-orientable two-dimensional surface with only one side when embedded in three-dimensional Euclidean space. It was independently discovered by Johann Benedict Listing a few months earlier. The M\"obius configuration, formed by two mutually inscribed tetrahedra, is also named after him. M\"obius was the first to introduce homogeneous coordinates into projective geometry. He is recognized for the introduction of the Barycentric coordinate system. Before 1853 and Schl\"afli's discovery of the 4-polytopes, M\"obius (with Cayley and Grassmann) was one of only three other people who had also conceived of the possibility of geometry in more than three dimensions.  Many mathematical concepts are named after him, including the M\"obius plane, the M\"obius transformations, important in projective geometry, and the M\"obius transform of number theory. His interest in number theory led to the important M\"obius function $\mu(n)$ and the M\"obius inversion formula. In Euclidean geometry, he systematically developed the use of signed angles and line segments as a way of simplifying and unifying results.} transformations of profields.

\begin{definition}\label{d:Mobius} Let $\bar F$ is a profield, $F\defeq\bar F\setminus\{\infty\}$ and $F^*\defeq F\setminus\{0\}$. A bijective function $\mu:\bar F\to\bar F$ is a \index{M\"obius transformation}\defterm{M\"obius transformation} of $\bar F$ if there exist elements $a\in F^*$, $b,c\in F$ and a sign $\pm1\in\{-1,+1\}$ such that  $\mu(x)=a\cdot (x+b)^{\pm1}+c$ for all $x\in\bar F$. 
\end{definition} 

In Definition~\ref{d:Mobius} we assume that for every element $x\in\bar F$, $x^{+1}=x$, and $x^{-1}$ is the unique element such that $x\cdot x^{-1}=1=x^{-1}\cdot x$. For example, $\infty^{+1}=\infty$ and $\infty^{-1}=0$.

\begin{exercise} Show that $(x\cdot y)^{-1}=y^{-1}\cdot x^{-1}$ for every elements $x,y$ of a procorps.
\end{exercise}

\begin{proposition} The set of M\"obius transformations of a profield $\bar F$ is a subgroup of the permutation group of the profield $\bar F$.
\end{proposition}

\begin{proof} First we show that for every M\"obius transformation $\mu:\bar F\to\bar F$ of $\bar F$, its inverse $\mu^{-1}:\bar F\to\bar F$ is M\"obius. Find constants $a\in F^*$, $b,c\in F$ and a sign $\pm1\in\{-1,+1\}$ such that $\mu(x)=a\cdot(x+b)^{\pm1}+c$. Then for every $y\in \bar F$ the element $x=\mu^{-1}(y)$ can be found by the formula
$$x=(a^{-1}\cdot(y-c))^{\pm1}-b=(a^{-1})^{\pm1}\cdot(y-c)^{\pm1}-b,$$witnessing that $\mu^{-1}$ is a M\"obius transformation of the profield $\bar F$.

Next, we prove that the composition $\mu\circ\eta$ of two M\"obius transformations $\mu,\eta:\bar F\to\bar F$ is a M\"obius transformation. Find elements $a,\alpha\in F^*$, $b,c,\beta,\gamma\in F$ and signs $s,\sigma\in\{-1,+1\}$ such that $$\mu(x)=a\cdot(x+b)^s+c\quad\mbox{and}\quad\eta(x)=\alpha\cdot(x+\beta)^\sigma+\gamma$$ for all $x\in\bar F$.

Then for every $x\in\bar F$, we have the equality
$$\mu\circ\eta(x)=\mu(\alpha\cdot(x+\beta)^\sigma+\gamma)=a\cdot(\alpha\cdot(x+\beta)^\sigma+\gamma+b)^s+c.$$
If $s=+1$, then 
$$\mu\circ\eta(x)=a\cdot(\alpha\cdot(x+\beta)^\sigma+\gamma+b)^s+c=(a\cdot\alpha)\cdot(x+\beta)^\sigma+(a\cdot(\gamma+b)+c)$$is a M\"obius transformation.

If $\sigma=+1$, then $$\mu\circ\eta(x)=a\cdot(\alpha\cdot(x+\beta)^\sigma+\gamma+b)^s+c=(a\cdot\alpha^s)\cdot(x+ (\beta+\alpha^{-1}\cdot(\gamma+b))^s+c$$is a M\"obius transformation.

Finally, assume that $s=\sigma=-1$. In this case $\mu\circ\eta(x)=a\cdot(\alpha\cdot(x+\beta)^{-1}+\gamma+b)^{-1}+c$. If $\gamma+b=0$, then 
$$\mu\circ\eta(x)=a\cdot(\alpha\cdot(x+\beta)^{-1}+\gamma+b)^{-1}+c=a\cdot(\alpha\cdot(x+\beta)^{-1})^{-1}+c=
(a\cdot\alpha^{-1})\cdot(x+\beta)+c$$
is a M\"obius transformation.
So, we assume that $\gamma+b\ne 0$. In this case
$$
\begin{aligned}
\mu\circ\eta(x)&=a\cdot(\alpha\cdot(x+\beta)^{-1}+\gamma+b)^{-1}+c=a\cdot\big((\alpha+(\gamma+b)(x+\beta))(x+\beta)^{-1}\big)^{-1}+c\\
&=a\cdot\big((\gamma+b)((x+\beta)+(\gamma+b)^{-1}\cdot\alpha\big)^{-1}\cdot(x+\beta)+c=\\
&=a\cdot(\gamma{+}b)^{-1}\cdot\big(x+(\beta+(\gamma{+}b)^{-1}\cdot\alpha)\big)^{-1}\cdot(x+\beta+(\gamma{+}b)^{-1}\cdot \alpha-(\gamma{+}b)^{-1}\cdot\alpha)+c\\
&=a\cdot(\gamma{+}b)^{-1}+a\cdot(\gamma{+}b)^{-1}\cdot\big(x+(\beta+(\gamma{+}b)^{-1}\cdot\alpha)\big)^{-1}\cdot(-(\gamma{+}b)^{-1}\cdot\alpha)+c\\
&=(-a\cdot(\gamma+b)^{-2}\cdot\alpha)\cdot\big(x+(\beta+(\gamma+b)^{-1}\cdot\alpha)\big)^{-1}+(a\cdot(\gamma+b)^{-1}+c)
\end{aligned}
$$is a M\"obius transformation of $\bar F$, too.
\end{proof}

\begin{theorem}\label{t:Mobius-coincide} Two M\"obius transformations $\mu,\eta:\bar F\to\bar F$ of a profield $\bar F$ coincide if and only if $\mu{\restriction}_{\{0,1,\infty\}}=\eta{\restriction}_{\{0,1,\infty\}}$.
\end{theorem}

\begin{proof} The ``only if'' part of this characterization is trivial. To prove the ``if'' part, assume that $\mu{\restriction}_{\{0,1,\infty\}}=\eta{\restriction}_{\{0,1,\infty\}}$. 
Find elements $a,\alpha\in F^*\defeq \bar F\setminus\{0,\infty\}$, $b,c,\beta,\gamma\in F\defeq\bar F\setminus\{\infty\}$ and signs $s,\sigma\in\{-1,+1\}$ such that $\mu(x)=a\cdot(x+b)^s+c$ and $\eta(x)=\alpha\cdot(x+\beta)^\sigma+\gamma$ for every $x\in \bar F$. If $s\ne \sigma$, then we lose no generality assuming that $s=+1$ and $\sigma=-1$. In this case
$$\mu(\infty)=a\cdot(\infty+b)^{+1}+c=a\cdot\infty+c=\infty+c=\infty$$and
$$\eta(\infty)=\alpha\cdot(\infty+\beta)^{-1}+\gamma=\alpha\cdot\infty^{-1}+\gamma=\alpha\cdot 0+\gamma=0+\gamma=\gamma\ne\infty=\mu(\infty),$$which is a contradiction showing that $s=\sigma$. 

If $s=\sigma={+1}$, then $$a\cdot b+c=\mu(0)=\eta(0)=\alpha\cdot\beta+\gamma$$and 
$$a=a\cdot(1+b)+c-(a\cdot b+c)=\mu(1)-\mu(0)=\eta(1)-\eta(0)=\alpha\cdot(1+\beta)+\gamma-(\alpha\cdot\beta+\gamma)=\alpha.$$ Then for every $x\in\bar F$ we have
$$\mu(x)=a\cdot (x+b)+c=a\cdot x+(a\cdot b+c)=\alpha\cdot x+(\alpha\cdot \beta+\gamma)=\alpha\cdot(x+\beta)+\gamma=\eta(x)$$and hence $\mu=\eta$.

Next, assume that $s=\sigma=-1$. In this case
$$\mu(\infty)=a\cdot(\infty+b)^{-1}+c=a\cdot\infty^{-1}+c=a\cdot 0+c=c$$and 
$$\gamma=\alpha\cdot(\infty+\beta)^{-1}=\eta(\infty)=\mu(\infty)=c.$$ 

If $b=0$, then $$\mu(0)=a\cdot 0^{-1}+c=a\cdot\infty+c=\infty.$$ Assuming that $\beta\ne 0$, we conclude that $$\eta(0)=\alpha\cdot(0+\beta)^{-1}+\gamma=\alpha\cdot\beta^{-1}+\gamma\ne\infty=\mu(0),$$ which is a contradiction showing that $b=0$ implies $\beta=0$. Then $\mu(1)=a\cdot(1+b)^{-1}+c=a+c$ and $\alpha+c=\alpha+\gamma=\eta(1)=\mu(1)=a+c$ implies the equality $\alpha=a$. Therefore, $(a,b,c)=(a,\beta,c)$ and hence $\mu=\eta$. By analogy we can prove that $\beta=0$ implies $\mu=\eta$.

If $b=-1$, then $\mu(1)=a\cdot (1+b)^{-1}+c=a\cdot 0^{-1}+c=a\cdot\infty+c=\infty+c=\infty$. Assuming that $\beta\ne b=-1$, we conclude that $1+\beta\ne 0$ and hence $$\eta(1)=\alpha\cdot(1+\beta)^{-1}+c\ne \infty=\mu(1),$$which is a contradiction showing that $b=-1$ implies $\beta=-1$. In this case
$$-a+c=a\cdot(0-1)^{-1}+c=\mu(0)=\eta(0)=\alpha\cdot(0-1)^{-1}+\gamma=-\alpha+\gamma$$and $c=\gamma$ imply $a=\alpha$ and hence $(a,b,c)=(\alpha,\beta,\gamma)$ and $\mu=\eta$.
By analogy we can prove that $\beta=-1$ implies $\mu=\eta$.

So, assume that $b,\beta\notin\{0,-1\}$. Then $$a\cdot b^{-1}=(a\cdot(0+b)^{-1}+c)-c=\mu(0)-\mu(\infty)=\eta(0)-\eta(\infty)=\alpha\cdot\beta^{-1}$$and
 $a\cdot (1+b)^{-1}=\mu(1)-\mu(\infty)=\eta(1)-\eta(\infty)=\alpha\cdot(1+\beta)^{-1}$ imply
 $a\cdot\beta=\alpha\cdot b$ and $a\cdot(1+\beta)=\alpha\cdot(1+b)$. Then 
 $$a=a\cdot(1+\beta)-a\cdot\beta=\alpha\cdot(1+b)-\alpha\cdot b=\alpha$$and
 $$b=a\cdot(a\cdot b^{-1})^{-1}=a\cdot\mu(0)^{-1}=\alpha\cdot\eta(0)^{-1}=\alpha\cdot(\alpha\cdot\beta^{-1})^{-1}=\beta.$$ Therefore, $(a,b,c)=(\alpha,\beta,\gamma)$ and hence $\mu=\eta$.
\end{proof}

\begin{theorem} For every distinct elements $z,u,i$ of a profield $\bar F$, there exists a unique M\"obius transformation $\mu:\bar F\to\bar F$ such that $\mu(0)=z$, $\mu(1)=u$, and $\mu(i)=\infty$.
\end{theorem}

\begin{proof} If $i=\infty$, then the M\"obius transformation $$\mu:\bar F\to\bar F,\quad \mu:x\mapsto (u-z)\cdot x+z,$$ has the required property $\mu01\infty=zui$.

If $z=\infty$, then the M\"obius transformation $$\mu:\bar F\to\bar F,\quad \mu:x\mapsto (u-i)\cdot x^{-1}+i,$$ has the required property $\mu01\infty=zui$.

If $u=\infty$, then the M\"obius transformation 
$$\mu:\bar F\to\bar F,\quad \mu:x\mapsto (i-z)\cdot (x-1)^{-1}+i,$$has the required property $\mu01\infty=zui$.

So, assume that $\infty\notin\{z,u,i\}$. In this case the M\"obius transformation
$$\mu:\bar F\to\bar F,\quad\mu:x\mapsto (z-i)\cdot(z-u)^{-1}\cdot (u-i)\cdot(x+(z-u)^{-1}\cdot (u-i))^{-1}+i,$$has the required property $\mu01\infty=zui$. Indeed,
$$
\begin{aligned}
\mu(0)&=(z-i)\cdot(z-u)^{-1}\cdot (u-i)\cdot((z-u)^{-1}\cdot (u-i))^{-1}+i\\
&=(z-i)\cdot(z-u)^{-1}\cdot (u-i)\cdot((u-i)^{-1}\cdot (z-u))+i=z,\\
\mu(1)&=(z-i)\cdot(z-u)^{-1}\cdot (u-i)\cdot(1+(z-u)^{-1}\cdot (u-i))^{-1}+i\\
&=(z-i)\cdot(z-u)^{-1}\cdot (u-i)\cdot((z-u)^{-1}((z-u)+(u-i))^{-1}+i\\
&=(z-i)\cdot(z-u)^{-1}\cdot (u-i)\cdot(z-u)\cdot (z-i)^{-1}+i=u,\\
\mu(\infty)&=(z-i)\cdot(z-u)^{-1}\cdot (u-i)\cdot(\infty+(z-u)^{-1}\cdot (u-i))^{-1}+i\\
&=(z-i)\cdot(z-u)^{-1}\cdot (u-i)\cdot\infty^{-1}+i=0+i=i.\\
\end{aligned}
$$
Therefore, for every distinct elements $z,u,i\in\bar F$, there exists a M\"obius transformation $\mu:\bar F\to\bar F$ with $\mu01\infty=zui$. The uniqueness of $\mu$ follows from Theorem~\ref{t:Mobius-coincide}. 
\end{proof} 

\section{Coordinate charts on lines in Pappian projective spaces}\label{s:projchart}

Let $L$ be a line in a projective space $X$. Given any distinct points $z,u,i\in L$, consider the function $$\phi_{zui}:L\to \ddddot X\projeqind,\quad \phi_{zui}:x\mapsto \overvector{zxui},$$ assigning to every point $x\in L$ the proportion $\overvector{zxui}$. Corollary~\ref{c:proj4-exist} implies that the function $\phi_{zui}:L\to \ddddot X\projeqind$ is surjective. If the projective space $X$ is Pappian, then $\ddddot X\projeqind=\bar\IF_X$ and the function $\phi_{zui}:L\to\ddddot X\projeqind=\bar\IF_X$ is bijective, by Corollary~\ref{c:Pappian-exist-unique}. In this case the function $\phi_{zui}:L\to\bar\IF_X$ is called a \index{coordinate chart}\defterm{coordinate chart} on the projective line $L$. The coordinate chart $\phi_{zui}$ identifies the projective line $L$ with the proscalar profield $\bar\IF_X$. Under this identification, we have 
$$\phi_{zui}(z)=\overvector{zzui}=0,\quad\phi_{zui}(u)=\overvector{zuui}=1,\quad\mbox{and}\quad \phi_{zui}(i)=\overvector{ziui}=\infty.$$ 

For six points $z,u,i,o,e,j\in L$ with $|\{z,u,i\}|=3=|\{o,e,j\}|$, the bijective function
$$\phi_{zui}\circ\phi_{oej}^{-1}:\bar\IF_X\to\bar\IF_X$$is called the \index{transition function}\defterm{transition function} between the charts $\phi_{oej}$ and $\phi_{zui}$.

$$
\xymatrix{
&L\ar^{\phi_{zui}}[rd]\ar_{\phi_{oej}}[ld]\\
\bar\IF_X\ar_{\phi_{zui}\circ\phi^{-1}_{oej}}[rr]&&\bar\IF_X
}
$$

We are going to show that the transition functions are M\"obius transformations of the proscalar profield $\bar\IF_X$.

\begin{lemma}\label{l:transition-affine-i} Let $L$ be a line in a Pappian projective space $X$, and $z,u,o,e,i\in L$ be points such that $|\{z,u,i\}|=3=|\{o,e,i\}|$. Then $$\phi_{zui}\circ\phi_{oei}^{-1}(p)=p\cdot (\overvector{zeui}-\overvector{zoui})+\overvector{zoui}$$for every proscalar $p\in\bar\IF_X$. Consequently, $\phi_{zui}\circ\phi_{oei}^{-1}$ is a M\"obius transformation of the profield $\bar F_X$.
\end{lemma}

\begin{proof} If $|L|=3$, then $L=\{z,u,i\}=\{o,e,i\}$ and hence $zu=oe$ or $zu=eo$. If $zu=oe$, then $$\phi_{zui}\circ\phi_{oei}^{-1}(p)=p=p\cdot 1+0=p\cdot (\overvector{zuui}-\overvector{zzui})+\overvector{zzui}=p\cdot (\overvector{zeui}-\overvector{zoui})+\overvector{zoui}$$for every proscalar $p\in\bar\IF_X$.

So, assume that $zu=eo$. In this case $$p\cdot (\overvector{zeui}-\overvector{zoui})+\overvector{zoui}=p\cdot (\overvector{zzui}-\overvector{zuui})+\overvector{zuui}=p\cdot (-1)+1$$
The equality $|L|=3$ implies $\bar\IF_X=\{0,1,\infty\}$. 

If $p=0$, then $\phi_{zui}\circ\phi_{oei}^{-1}(p)=\phi_{zui}\circ\phi_{oei}^{-1}(0)=\phi_{zui}(o)=\phi_{zui}(u)=1=0\cdot(-1)+1=p\cdot(-1)+1.$

If $p=1$, then $\phi_{zui}\circ\phi_{oei}^{-1}(p)=\phi_{zui}\circ\phi_{oei}^{-1}(1)=\phi_{zui}(e)=\phi_{zui}(z)=0=1\cdot(-1)+1=p\cdot(-1)+1.$

If $p=\infty$, then $\phi_{zui}\circ\phi_{oei}^{-1}(p)=\phi_{zui}\circ\phi_{oei}^{-1}(\infty)=\phi_{zui}(i)=\infty=\infty\cdot(-1)+1=p\cdot(-1)+1.$

Therefore, for every $p\in \bar\IF_X=\{0,1,\infty\}$ we have
$$\phi_{zui}\circ\phi_{oei}^{-1}(p)=p\cdot(-1)+1=p\cdot (\overvector{zeui}-\overvector{zoui})+\overvector{zoui}.$$

Next, assume that $|L|\ge 4$. In this case the projective space $X$ is $4$-long. Choose any hyperplane $H\subset X$ such that $L\cap H=\{i\}$, and consider the subliner $A\defeq X\setminus H$, which is an affine space, by Theorem~\ref{t:affine<=>hyperplane}.

By Proposition~\ref{p:Pappian-minus-flat}, the affine space $A$ is Pappian and by Theorem~\ref{t:Papp<=>Des+RX}, the corps $\IR_A$ is a field. Then $Z(\IR_A)=\IR_A$. By Corollary~\ref{c:FX=Z(RA)}, the function $$\Psi_A:\IF_X\to Z(\IR_A)=\IR_A,\quad \Psi_A:\overvector{zvui}\mapsto\overvector{zvu},$$ is an isomorphism of the fields $\IF_X$ and $\IR_A$.  By Theorem~\ref{t:lineaff=aff}, $$\phi_{zu}\circ\phi_{oe}^{-1}(s)=s\cdot (\overvector{zeu}-\overvector{zou})+\overvector{zou}$$ for every scalar $s\in\IR_A$.  

The definitions of the coordinate charts $\phi_{zu},\phi_{zui},\phi_{oe},\phi_{oei}$ imply $$\phi_{zu}=\Psi_A\circ \phi_{zui}{\restriction}_{L\setminus\{i\}}\quad\mbox{and}\quad\phi_{oe}=\Psi_A\circ \phi_{oei}{\restriction}_{L\setminus\{i\}}.$$ 
Given a finite proscalar $p\in\IF_X$, consider the scalar $s\defeq\Psi_A(p)\in\IR_A$ and observe that 
$$
\begin{aligned}
\Psi_A\circ \phi_{zui}\circ\phi_{oei}^{-1}(p)&=\phi_{zu}\circ\phi_{oe}^{-1}(s)=s\cdot(\overvector{zeu}-\overvector{zou})+\overvector{zou}\\
&=\Psi_A(p)\cdot(\Psi_A(\overvector{zeui})-\Psi_A(\overvector{zoui}))+\Psi_A(\overvector{zoui})=\Psi_A(p\cdot(\overvector{zeui}-\overvector{zoui})+\overvector{zoui}),
\end{aligned}
$$which implies
the desired equality $$\phi_{zui}\circ\phi_{oei}^{-1}(p)=p\cdot(\overvector{zeui}-\overvector{zoui})+\overvector{zoui},$$ by the injectivity of the function $\Psi_A$.

For $p=\infty$, we obtain the same equality as
$$\phi_{zui}\circ\phi_{oei}^{-1}(p)=\phi_{zui}\circ\phi_{oei}^{-1}(\infty)=\phi_{zui}(i)=\infty=\infty\cdot(\overvector{zeui}-\overvector{zoui})+\overvector{zoui}=p\cdot(\overvector{zeui}-\overvector{zoui})+\overvector{zoui}.$$
\end{proof}

For an element $s\in \bar F$ of a profield $\bar F$, let $s^{-1}$ denote its multiplicative inverse in $\bar F$.

\begin{lemma}\label{l:transition-inversion} Let $L$ be a line in a Pappian affine space and $z,u,i\in L$ be distinct points. Then $$\phi_{zui}\circ\phi_{iuz}^{-1}(p)=p^{-1}$$ for every proscalar $p\in\bar\IF_X$. Consequently, $\phi_{zui}\circ\phi_{iuz}^{-1}$ is a M\"obius transformation of the profield $\bar\IF_X$.
\end{lemma}

\begin{proof} By Lemma~\ref{l:1/x}, there exists a line projectivity $P:L\to L$ such that $Pzui=iuz$ and for every point $v\in L$ and its image $w\defeq P(v)$ we have $\overvector{zvui}=(\overvector{zwui})^{-1}$.  Given any proscalar $p\in \bar\IF_X$, consider the point $v\defeq\phi_{iuz}^{-1}(p)$ and its image $w\defeq P(v)$ under the line projectivity $P$. It follows from $Pivuz=zwui$ that $\overvector{ivuz}=\overvector{zwui}$. The choice of $P$ ensures that $$\phi_{zui}\circ\phi^{-1}_{iuz}(p)=\phi_{zui}(v)=\overvector{zvui}=(\overvector{zwui})^{-1}=(\overvector{ivuz})^{-1}=(\phi_{iuz}(v))^{-1}=p^{-1}.$$
\end{proof}

\begin{theorem}\label{t:transition-Mobius} Let $L$ be a line in a projective space $X$ and $z,u,i,o,e,i\in L$ be points such that $|\{z,u,i\}|=3=|\{o,e,j\}|$.
The transition function $\phi_{oej}\circ\phi_{zui}^{-1}$ is a M\"obius transformation of the proscalar profield $\bar\IF_X$.
\end{theorem}

\begin{proof} If $i=j$, then $\phi_{oej}\circ\phi_{zui}^{-1}$ is a M\"obius transformation of $\bar\IF_X$, by Lemma~\ref{l:transition-affine-i}. So, assume that $i\ne j$. Choose any point $w\in L\setminus\{i,j\}$. By Lemmas~\ref{l:transition-affine-i} and \ref{l:transition-inversion}, the transition functions $\phi_{jwi}\circ\phi_{zui}^{-1}$, $\phi_{iwj}\circ\phi_{jwi}^{-1}$, and $\phi_{oej}\circ\phi_{iwj}^{-1}$ are M\"obius transformations of the profield $\bar\IF_X$. Since the composition of M\"obius transformations is a M\"obius tranformation, the transition function
$$\phi_{oej}\circ\phi_{zui}^{-1}=(\phi_{oej}\circ\phi_{iwj}^{-1})\circ(\phi_{iwj}\circ\phi_{jwi}^{-1})\circ(\phi_{jwi}\circ\phi_{zui}^{-1})$$is a M\"obius transformation of $\bar\IF_X$.
\end{proof}

\begin{lemma}\label{l:ptransition-identity} Let $F:L\to \Lambda$ be a bijective function between two lines $L,\Lambda$ in a Pappian projective space $X$. Given three distinct points $z,u,i\in L$, consider the line triple $oej\defeq Fzui$. The function $F$ is a line projectivity if and only if $\phi_{oei}\circ F\circ \phi_{zui}^{-1}$ is the identity map of the profield $\bar\IF_X$.
\end{lemma}

\begin{proof} To prove the ``only if'' part, assume that $F$ is a line projectivity in $X$. Given any proscalar $p\in\bar\IF_X$, find a unique point $x\in L$ such that $\overvector{zxui}=p$. Then for the point $y\defeq F(x)$ we have $oyej=Fzxui$ and hence $\overvector{oyej}=\overvector{zxui}$ and 
$$\phi_{oej}\circ F\circ \phi_{zui}^{-1}(p)=\phi_{oej}\circ F(x)=\phi_{oej}(y)=\overvector{oyej}=\overvector{zxui}=p,$$
witnessing that $\phi_{oej}\circ F\circ \phi_{zui}^{-1}$ is the identity function of the profield $\bar \IF_X$. 
\smallskip

Now assume conversely that $\phi_{oej}\circ F\circ \phi_{zui}^{-1}$ is the identity function of the profield $\bar\IF_X$. 
Then $\phi_{zui}=\phi_{oej}\circ F$. 
By Theorem~\ref{t:projective=>3-transitive}, there exists a line projectivity $P:L\to \Lambda$ such that $Pzui=oej$. By the definition of the projective equivalence of line quadruples, for every point $x\in L$ and its image $y\defeq P(x)$, the line quadruple $Pzxui$ is affinely equivalent to the line quadruple $oyej$ and hence $\overvector{oyej}=\overvector{zxui}$. Since $\phi_{zui}=\phi_{oej}\circ F$, for the point $y'\defeq F(x)$, we have $\overvector{oyej}=\overvector{zxui}=\phi_{zui}(x)=\phi_{oej}\circ F(x)=\phi_{oej}(y')=\overvector{oy'ej}$ and hence $P(x)=y=y'=F(x)$, because the Pappian quadruples  $oyej$ and $oy'ej$ are projectively equivalent. Therefore, $F=P$ is a line projectivity.  
\end{proof}

Now we are able to prove our main result on coordinate description of line projectivities in Pappian projective spaces.

\begin{theorem} Let $L,\Lambda$ be two lines in a Pappian projective space, and $z,u,i\in L$, $o,e,j\in \Lambda$ be points such that $|\{z,u,i\}|=3=|\{o,e,j\}|$. A function $F:L\to \Lambda$ is a line projectivity if and only if the function  $\phi_{oej}\circ F\circ\phi_{zui}^{-1}$ is a M\"obius tranformation of the proscalar profield $\bar\IF_X$.
\end{theorem}

\begin{proof} The ``only if'' part is proved in the following lemma.

\begin{lemma}\label{l:lineproj=>Mobius} For every line projectivity $P:L\to \Lambda$, the function  $ \phi_{oej}\circ P\circ\phi_{zui}^{-1}$ is a M\"obius tranformation of the  profield $\bar\IF_X$.
\end{lemma}

\begin{proof} Given any line projectivity $P:L\to\Lambda$, consider the line triple $abc\defeq Pzui$ on the line $\Lambda$. By Lemma~\ref{l:ptransition-identity}, the transition function $\phi_{abc}\circ P\circ \phi_{zui}^{-1}$ is the identity function of the profield $\bar\IR_X$. By Theorem~\ref{t:transition-Mobius}, the transition function $\phi_{oej}\circ\phi_{abc}^{-1}$ is a M\"obius transformation of the profield $\bar\IF_X$. Then $
\phi_{oej}\circ P\circ\phi_{zui}^{-1}=(\phi_{oej}\circ \phi_{abc}^{-1})\circ(\phi_{abc}\circ P\circ\phi_{zui}^{-1})$ is a M\"obius transformation of $\bar\IR_X$.
\end{proof}

To prove the ``if'' part, assume that $F:L\to\Lambda$ is a function such that $\phi_{oej}\circ F\circ\phi_{zui}^{-1}$ is a M\"obius transformation of the profield $\bar \IF_X$. 
Since the M\"obius transformation $\phi_{oej}\circ F\circ\phi_{zui}^{-1}$ is a bijective function of $\bar\IF_X$ and the coordinate charts $\phi_{oej}:\Lambda\to\bar\IF_X$ and $\phi_{zui}:L\to\bar\IF_X$ are bijective functions, the function $F= \phi_{oej}^{-1}\circ(\phi_{oej}\circ F\circ\phi_{zui}^{-1})\circ \phi_{zui}$ is bijective. Then the triple $z'u'i'\defeq Fzui$ consists of three distinct points of the line $\Lambda$. By Theorem~\ref{t:projective=>3-transitive}, there exists a line projectivity $P:L\to\Lambda$ such that $Pzui=z'u'i'=Fzui$. By Lemma~\ref{l:lineproj=>Mobius}, the function $\phi_{oej}\circ P\circ\phi_{zui}^{-1}$ is a M\"obius transformation of the profield $\bar \IF_X$. Since
$$\phi_{oej}\circ P\circ\phi_{zui}^{-1}01\infty=\phi_{oej}\circ Pzui=\phi_{oej}\circ Fzui=\phi_{oej}\circ F\circ\phi_{zui}^{-1}01\infty,$$the M\"obius tranformations $\phi_{oej}\circ P\circ\phi_{zui}^{-1}$ and $\phi_{oej}\circ F\circ\phi_{zui}^{-1}$ coincide, according to Theorem~\ref{t:Mobius-coincide}. Then the function 
$F=\phi_{oej}^{-1}\circ(\phi_{oej}\circ F\circ\phi_{zui}^{-1})\circ\phi_{zui}=\phi_{oej}^{-1}\circ(\phi_{oej}\circ P\circ\phi_{zui}^{-1})\circ\phi_{zui}=P$ is a line projectivity.
\end{proof}

\section{Projectivities, proportionalities, and proscalarities}

\begin{definition} An automorphism $A:X\to X$ of a projective space $X$ is called
\begin{itemize}
\item a \index{projectivity}\defterm{projectivity} if for every line $L\subseteq X$, the restriction $A{\restriction}_L$ is a line projectivity;
\item a \index{proportionality}\defterm{proportionality} if for every line quadruple $zvui\in \ddddot X$ we have $Azvui\projeq zvui$;
\item a \index{proscalarity}\defterm{proscalarity} if for every Pappian line quadruple $zvui\in\ddddot X$ we have $Azvui\projeq zvui$.
\end{itemize}
\end{definition}

\begin{proposition}\label{p:proj=>pp=>ps} For any automorphism $A:X\to X$ of a projective space $X$, the following conditions hold.
\begin{enumerate}
\item If $A$ is a projectivity, then $A$ is proportionality.
\item If $A$ is a proportionality, then $A$ is a proscalarity.
\end{enumerate}
\end{proposition}

Therefore, for every automorphism of a projective spaces, we have the implications:
$$\mbox{projectivity}\Ra\mbox{proportionality}\Ra\mbox{proscalarity}.$$
For Pappian projective spaces, those three notions are equivalent.

\begin{theorem}\label{t:P=>p=pp=ps} For any automorphism $A:X\to X$ of a Pappian projective space $X$, the following conditions are equivalent:
\begin{enumerate}
\item $A$ is a projectivity;
\item $A$ is a proportionality;
\item $A$ is a proscalarity.
\end{enumerate}
\end{theorem}

\begin{proof} The implications $(1)\Ra(2)\Ra(3)$ follow from Proposition~\ref{p:proj=>pp=>ps}. To prove the implication $(3)\Ra(1)$, assume that $A$ is a proscalarity. To prove that $A$ is a projectivity, fix any line $L$. Choose any distinct points $z,u,i\in L$ and consider the pair $oej\defeq Azui$. By Theorem~\ref{t:projective=>3-transitive}, there exists a line projectivity $B$ in $X$ such that $Bzui=oej=Azui$. We claim that $A{\restriction}_L=B$. Given any point $x\in L$, we should prove that the points $a\defeq A(x)$ and $b\defeq B(x)$ coincide.  Since $A$ is a proscalarity and $B$ is a line projectivity, $oaej=Azxui\projeq zxui\projeq Bzxui=obej$ and hence $a=b$ because the line quadruples $oaej$ and $obej$ are Pappian and projectively equivalent. Therefore, $A(x)=a=b=B(x)$ and $A=B$.
\end{proof} 

\begin{definition} An automorphism $A:X\to X$ of a Pappian projective space $X$ is called \index{projective automorphism}\index{automorphism!projective}\defterm{projective} if $A$ satisfies the equivalent conditions of Theorem~\ref{t:P=>p=pp=ps}.
\end{definition}

For a projective space $X$, let $\Proj(X)$, $\Proj_{P}(X)$, $\Proj_{S}(X)$ be the subsets of the automorphism group $\Aut(X)$, consisting of all projectivities, proportionalities, and proscalarities, respectively.
Proposition~\ref{p:proj=>pp=>ps} and Theorem~\ref{t:P=>p=pp=ps} imply the following corollary.

\begin{corollary} For every projective space $X$, the sets $\Proj(X),\Proj_{P}(X),\Proj_{S}(X)$ are normal subgroups of the automorphism group $\Aut(X)$ such that $\Proj(X)\subseteq \Proj_{P}(X)\subseteq\Proj_{S}(X)$. If the projective space $X$ is Pappian, then  $\Proj(X)=\Proj_{P}(X)=\Proj_{S}(X)$.
\end{corollary}

\begin{remark} By Theorem$^\dag$~\ref{t:Grundhofer-Muller-Nagy-p}, for every line $L$ in a non-Pappian finite projective space $X$, the group $\Sym\projupind_X(L)$ of line projectivities of the line $L$ is either the symmetric group $\Sym(L)$ of $L$ or the alternative group $\Alt(L)$ of $L$. This implies that $\Proj_{P}(X)=\Proj_{S}(X)=\Aut(X)$. Moreover, $\Proj(X)=\Aut(X)$ if $\Sym\projupind_X(L)=\Sym(L)$, and $\Proj(X)=\{A\in\Aut(X):\forall L\in\mathcal L_X\;\;A{\restriction}_L\in\Alt(L)\}$ if $\Sym\projupind_X(L)=\Alt(L)$.
\end{remark}




Let $X$ be a projective space and $\bar\IF_X$ be its proscalar profield. For every automorphism $A:X\to X$ and every line quadruple $zvui$ in $X$, the image $Azvui$ is a line quadruple in $X$. If the line quadruple $zvui$ is Pappian, then so is the line quadruple $Azvui$. This observation allows us to define the map $\ddddot A:\bar\IF_X\to\bar\IF_X$ assigning to every proscalar $p\in\IR_X$ the proscalar $\ddddot A(p)=\{Azvui:zvui\in p\}$. It is easy to see that the function $\ddddot A:\bar\IF_X\to\bar \IF_X$ is an automorphism of the profield $\bar\IF_X$. 
For two automorphisms $A,B:X\to X$ and their composition $C=A\circ B$, the automorphism $\ddddot C:\bar\IF_X\to\bar\IF_X$ is equal to the composition $\ddddot A\circ\ddddot B$ of the automorphisms $\ddddot A,\ddddot B$ of the profield $\bar\IF_X$. This means that the correspondence $T:\Aut(X)\to\Aut(\bar\IF_X)$, $T:A\mapsto \ddddot A$, is a homomorphism from the automorphism group $\Aut(X)$ of the liner $X$ into the automorphism group $\Aut(\bar\IF_X)$ of the profield $\bar\IF_X$. The definition of a proscalarity implies that the automorphism $\ddddot A:\bar\IF_X\to\bar\IF_X$ is the identity map of $\bar\IF_X$ if and only if the automorphism $A$ is a proscalarity of $X$. 
 

Let us write down this fact for future references.

\begin{proposition}\label{p:Aut(X)->Aut(FX)} For every projective space $X$, the function $\Aut(X)\to \Aut(\bar\IF_X)$, \ $A\mapsto \ddddot A$, is a homomorphism from the automorphism group $\Aut(X)$ of the projective space $X$ into the automorphism group $\Aut(\bar\IF_X)$ of the profield $\bar\IF_X$. The kernel of the homomorphism $\Aut(X)\to\Aut(\bar \IF_X)$ coincides with the group $\Proj_S(X)$ of proscalarities of $X$.
\end{proposition}

\section{Projective repers and frames in Pappian projective spaces}

Let us recall that a \index{projective reper}\defterm{projective reper} in a projective space $X$ is a function $r\subseteq X\times X$ whose domain $\dom[r]$ is a maximal independent set in $X$ containing a unique point $o=r(o)$ such that $r(x)\in\Aline ox\setminus\{o,x\}$ for every $x\in \dom[r]\setminus\{o\}$. The point $o$ is called the \index{origin}\defterm{origin} of the projective reper $r$. 

\begin{picture}(100,90)(-200,-15)
\put(0,0){\line(1,1){60}}
\put(0,0){\line(-1,1){60}}

\put(0,0){\color{red}\circle*{3}}
\put(-2,-7){$o$}
\put(30,30){\color{blue}\circle*{3}}
\put(32,23){$a$}
\put(60,60){\color{violet}\circle*{3}}
\put(62,53){$r(a)$}
\put(-30,30){\color{blue}\circle*{3}}
\put(-38,23){$b$}
\put(-60,60){\color{violet}\circle*{3}}
\put(-78,53){$r(b)$}
\end{picture}

\begin{lemma}\label{l:autoproj=1} Let $r$ be a projective reper in a Pappian projective space $X$. A projective automorphism $A:X\to X$ is the identity if and only if $A(x)=x$ for every $x\in \dom[r]\cup\rng[r]$.
\end{lemma}

\begin{proof} The ``only if'' part is trivial. To prove the ``if'' part, assume that $A(x)=x$ for every $x\in \dom[r]\cup\rng[r]$. Let $o=r(o)$ be the origin of the reper $r$. By the projectivity of the automorphism $A$, for every $p\in\dom[r]$, the restriction $A{\restriction}_{L_p}$ of $A$ to the line $L_p\defeq\Aline op$ is a line projectivity such that $A(x)=x$ for every point $x\in\{o,p,r(p)\}$. Since $X$ is Pappian, we can apply Corollary~\ref{c:Pappus-Fix3} and conclude that $A{\restriction}_{L_p}$ is the identity map of the line $L_p$. We claim that $A$ is the identity map of $X$. Consider the set $E\defeq\{x\in X:A(x)=x\}$ and using the Kuratowski-Zorn Lemma, find a maximal flat $M$ in $X$ such that $o\in M\subseteq E$. We claim that $M=X$. The maximal independence of the set $\dom[r]$ in $X$ and Proposition~\ref{p:add-point-to-independent} imply the equality $X=\overline{\dom[r]}$. Assuming that $M\ne X=\overline{\dom[r]}$, we can find a point $u\in \dom[r]\setminus M$. The flatness of $M$ implies $M\cap\Aline ou=\{o\}$ and hence the point $v\defeq r(u)$ belongs to the set $\Aline ou\setminus\{o,u\}\subseteq E\setminus(M\cup\{u\})$. 

By the strong regularity of the projective space $X$, for every $x\in \overline{M\cup\{u\}}\setminus \Aline ou$ there exist points $a\in M\cap\Aline ux$ and $b\in M\cap\Aline vx$. Then $x\in \Aline ua\cap\Aline vb$. Assuming that $\{x\}\ne\Aline ua\cap\Aline vb$, we conclude that $\Aline ua=\Aline vb$ and hence $x\in\Aline ua=\Aline uv=\Aline ou$, which contradicts the choice of the point $x$. Taking into account that $A$ is an automorphism of $X$ with $\{u,v,a,b\}\subseteq \Aline {e}o\cup M\subseteq E$, we conclude that 
$A(x)\in A[\Aline ua\cap\Aline vb]=A[\Aline ua]\cap A[\Aline vb]=\Aline ua\cap\Aline vb=\{x\}$ and hence $x\in E$ and $\overline{M\cup\{e\}}\subseteq E$, which contradicts the maximality of the flat $M$. This contradiction shows that $E=X$ and hence $A$ is the identity automorphism of the projective space $X$.
\end{proof}

\begin{corollary}\label{c:isoPproj-unique} Let $r$ be a projective reper in a Pappian projective space $X$. Two isomorphisms $A,B:X\to Y$ to a projective space $Y$ coincide if and only if $\ddddot A=\ddddot B$ and  $A(x)=B(x)$ for every $x\in\dom[r]\cup\rng[r]$.
\end{corollary}

\begin{proof} The ``only if'' part is trivial. To prove the ``if'' part, assume that $\ddddot A=\ddddot B$ and $A(x)=B(x)$ for every point $x\in \dom[r]\cup\rng[r]$. The the automorphism $B^{-1}A$ of $X$ does not move points of the set $\dom[r]\cup\rng[r]$. The equality $\ddddot A=\ddddot B$ implies that the automorphim $B^{-1}A:X\to X$ is a proscalarity of $X$. By Theorem~\ref{t:P=>p=pp=ps}, $B^{-1}A$ is a projective automorphism of $X$. By Lemma~\ref{l:autoproj=1}, $B^{-1}A=1_X$ and hence $A=BB^{-1}A=B1_X=B$.
\end{proof}

\begin{theorem}\label{t:Papp-iso-reper} Let $X,Y$ be two Pappian projective spaces and $r,r'$ be two projective repers in $X,Y$, respectively. For every profield isomorphism $\bar I:\bar\IF_X\to\bar\IF_Y$ and every bijective function $\varphi:\dom[r]\cup\rng[r]\to\dom[r']\cup\rng[r']$ with $\varphi\circ r=r'\circ \varphi$, there exists a unique isomorphism $\Phi:X\to Y$ such that $\varphi\subseteq\Phi$ and $\ddddot \Phi=\bar I$.
\end{theorem} 

\begin{proof} By Hessenberg's Theorem~\ref{t:Hessenberg-proaffine}, the Pappian projective spaces $X,Y$ are Desarguesian. Let $o=r(o)$ and $o'=r'(o')$ be the origins of the repers $r,r'$, respectively. 
The maximal independence of the set $\dom[r]$ and the condition $r(x)\in\Aline ox\setminus\{o,x\}$ for all $x\in\dom[r]\setminus\{o\}$ imply that the set $\rng[r]$ is maximal independent in the projective space $X$ and hence the flat $H\defeq\overline{\rng[r]\setminus\{o\}}$ is a hyperplane and $A\defeq X\setminus H$ is an affine subliner of $X$. By analogy we can show that the flat $H'\defeq\overline{\rng[r']\setminus\{o'\}}$ is a hyperplane and $B\defeq Y\setminus H'$ is an affine subliner of $Y$. By Proposition~\ref{p:Pappian-minus-flat}, the affine liners $A,B$ are Pappian and by Theorem~\ref{t:Papp<=>Des+RX}, the scalar corps $\IR_A,\IR_B$ of the Pappian affine liners $A,B$ are fields.  Then $Z(\IR_A)=\IR_A$ and $Z(\IR_B)=\IR_B$. By Corollary~\ref{c:FX=Z(RA)}, the functions $\Psi_A:\IF_X\to Z(\IR_A)=\IR_A$ and $\Psi_B:\IF_Y\to Z(\IR_B)=\IR_B$ are field isomorphisms. Then $\Psi_B^{-1}\bar I\Psi_A$ is an isomorphism of the fields $\IR_A$ and $\IR_B$.

By Theorem~\ref{t:Des-iso-reper}, there exists an isomorphism $\Phi:X\to Y$ such that $\Phi{\restriction}_{\dom[r]\cup\rng[r]}=\varphi$ and the restriction $\phi\defeq\Phi{\restriction}_A$ is an isomorphism between the affine liners $A,B$ such that $\dddot\Phi= \Psi_B\bar I\Psi_A^{-1}$. The latter equality and the definition of the homomorphisms $\dddot \phi$ and $\ddddot \Phi$ imply the equality $\ddddot\Phi=\bar I$.

Therefore, $\Phi:X\to Y$ is an isomorphism of the projective spaces such that $\Phi{\restriction}_{\dom[r]\cup\rng[r]}=\varphi$ and $\ddddot {\,\Phi}=\bar I$. The uniqueness of the isomorphism $\Phi$ follows from Corollary~\ref{c:isoPproj-unique}.
\end{proof}

\begin{corollary}\label{c:Papp-iso-reper} Let $r,r'$ be two projective repers in a Pappian projective space $X$. Every bijective function $\varphi:\dom[r]\cup\rng[r]\to\dom[r']\cup\rng[r']$ with $\varphi\circ r=r'\circ \varphi$ extends to a unique projective isomorphism $\Phi:X\to Y$.
\end{corollary} 

\begin{proof} By Theorem~\ref{t:Papp-iso-reper}, the bijective function $\varphi$ extends to a unique automorphism $\Phi:X\to X$ such that $\varphi\subseteq\Phi$ and $\ddddot{\,\Phi}$ is the identity automorphism of the profield $\bar\IF_X$. By Theorem~\ref{t:P=>p=pp=ps}, the proscalarity $\Phi$ is a  projective automorphism of $X$.
\end{proof}
Let us recall that a subset $M\subseteq X$ of a projective space $X$ is called a {\em projective frame} in $X$ if $|M|=\|X\|+1<\w$ and for every $x\in M$, the set $M\setminus\{x\}$ is independent in $X$.

\begin{lemma}\label{l:unique-frame} Let $X$ be a Pappian projective space of finite rank $\|X\|$ and $M$ be a projective frame in $X$. A projective automorphism $A:X\to X$ is the identity map of $X$ if and only if $A(x)=x$ for every $x\in M$.
\end{lemma}

\begin{proof} The ``only if'' part of this characterization is trivial. To prove the ``only if'' part, assume the $A:X\to X$ is a projective automorphism such that $A(x)=x$ for every $x\in M$. Fix any distinct points $o,u\in M$. By Proposition~\ref{p:reper-for-frame}, there exists a projective reper $r$ in $X$ such that $\dom[r]=M\setminus\{u\}$, $r(o)=o$, and $\{r(e)\}=\Aline oe\cap\overline{M\setminus\{o,e\}}$ for every $e\in\dom[r]\setminus\{o\}$. It follows from $A{\restriction}_M=1_M$ that $A(x)=x$ for every $x\in\dom[r]$. On the other hand, for every $e\in \dom[r]\setminus\{o\}=M\setminus\{o,u\}$ we have
$$\{A(r(e))\}=A[\Aline oe]\cap A[\overline{M\setminus\{o,e\}}]=\Aline oe\cap\overline{M\setminus\{o,e\}}=\{r(e)\}$$ and hence $A(x)=x$ for all $x\in\dom[r]\cup\rng[r]$. Applying Lemma~\ref{l:autoproj=1}, we conclude that the projective automorphism $A$ is the identity function of $X$.
\end{proof} 

\begin{theorem}\label{t:Papp-iso-frame2} Let $X,Y$ be two Pappian projective spaces of finite rank $\|X\|=\|Y\|$ and $M,M'$ be two projective frames in the projective spaces $X,Y$, respectively. For every profield isomorphism $\bar I:\bar\IF_X\to\bar\IF_Y$ and every bijection $\varphi:M\to M'$ there exists a unique isomorphism $\Phi:X\to Y$ such that $\Phi{\restriction}_M=\varphi$ and $\ddddot{\, \Phi}=\bar I$.
\end{theorem}

\begin{proof} Fix any distinct points $o,u\in M$ and consider the points $o'\defeq\varphi(o)$ and $u'\defeq\varphi(u)$ in the projective frame $M'$. By Proposition~\ref{p:reper-for-frame}, there exists a projective reper $r$ in the projective space $X$ such that $\dom[r]=M\setminus\{u\}$, $r(o)=o$, and $\{r(e)\}=\Aline oe\cap\overline{M\setminus\{o,e\}}$ for every $e\in\dom[r]\setminus\{o\}$. By analogy, there exists a projective frame $r'$ in the projective space $Y$ such that $\dom[r']=M'\setminus\{u'\}$, $r(o')=o'$, and $\{r'(e)\}=\Aline {o'}{e'}\cap\overline{M'\setminus\{o',e'\}}$ for every $e'\in\dom[r']\setminus\{o'\}$. Consider the bijection $\psi:\dom[r]\cup\rng[r]\to\dom[r']\cup\rng[r']$ defined by
$$\psi(x)=\begin{cases}\varphi(x)&\mbox{if $x\in \dom[r]$};\\
r'\varphi r^{-1}(x)&\mbox{if $x\in\rng[r]\setminus\{o\}$}
\end{cases}
$$
The definition of the function $\psi$ ensures that $\psi\circ r=r'\circ \psi$. By Theorem~\ref{t:Papp-iso-reper}, there exists an isomorphism $\Phi:X\to Y$ such that $\ddddot{\,\Phi}=\bar I$ and $\Phi(x)=\psi(x)$ for all $x\in\dom[r]\cup\rng[x]$. Then also $\Phi(x)=\psi(x)=\varphi(x)$ for all $x\in M\setminus\{u\}$. It remains to show that $\Phi(u)=u'$. 
The independence of the sets $M\setminus\{o\}$ and $M'\setminus\{o'\}$ and Proposition~\ref{p:reper-for-frame}(5) imply the equalities
$$\{u\}=\bigcap_{e\in M\setminus\{o,u\}}\overline{\{r(e)\cup(M\setminus\{o,u,e\})}\quad\mbox{and}\quad \{u'\}=\bigcap_{e'\in M'\setminus\{o,u\}}\overline{\{r'(e')\cup(M'\setminus\{o',u',e'\}}),$$ which imply the equality
$$\{\Phi(u)\}=\bigcap_{e\in M\setminus\{o,u\}}\Phi[\overline{\{r(e)\}\cup(M\setminus\{o,u,e\})}]=\bigcap_{e'\in M'\setminus\{o',u'\}}\overline{\{r'(e')\cup(M'\setminus\{o',u',e'\})}=\{u'\}.$$
Therefore, $\Phi:X\to Y$ is an isomorphism of the Pappian projective spaces such that $\ddddot{\,\Phi}=\bar I$ and $\Phi{\restriction}_M=\varphi$. The uniqueness of the isomorphism $\Phi$ follows from Lemma~\ref{l:unique-frame}. 
\end{proof}

\begin{corollary}\label{c:Papp-iso-frame} Let $X$ be a Pappian projective space of finite rank. Every bijection $\varphi:M\to M'$ between two projective frames $M,M'\subseteq X$ extends to a unique projective automorphism of $X$.
\end{corollary} 

\begin{proof} By Theorem~\ref{t:Papp-iso-frame2}, the bijective function $\varphi:M\to M'$ extends to a unique automorphism $\Phi:X\to X$ such that $\ddddot{\,\Phi}$ is the identity automorphism of the profield $\bar\IF_X$. By Theorem~\ref{t:P=>p=pp=ps}, the proscalarity $\Phi$ is a  projective automorphism of $X$.
\end{proof}

\section{The automorphism groups of Pappian projective spaces}

\begin{theorem}\label{t:Aut(X)=Proj(X)xAut(FX)} Let $X$ be a Pappian projective space, $r$ be a projective reper in $X$, and $\Aut_r(X)\defeq\{A\in\Aut(X):\forall x\in \dom[r]\cup\rng[r]\;\;A(x)=x\}$. Then $\Aut_r(X)$ is a subgroup of the group $\Aut(X)$ such that $\Aut_r(X)\cap\Proj(X)=\{1_X\}$, $\Aut(X)=\Proj(X)\cdot\Aut_r(X)$ and the function $\Aut_r(X)\to\Aut(\bar\IF_X)$, $A\mapsto\ddddot A$, is a group isomorphism.
\end{theorem}

\begin{proof} It is clear that $\Aut_r(X)$ is a subgroup of the automorphism group $\Aut(X)$. Lemma~\ref{l:autoproj=1} implies $\Aut_r(X)\cap\Proj(X)=\{1_X\}$.

To prove that $\Aut(X)=\Proj(X)\cdot\Aut_M(X)$, take any automorphism $A:X\to X$ and consider the projective reper $r'\defeq ArA^{-1}$ in $X$. By Corollary~\ref{c:Papp-iso-reper}, there exists a projective automorphism $B:X\to X$ such that  $B(x)=A(x)$ for all $x\in\dom[r]\cup\rng[r]$. Then the automorphism $C=B^{-1}A$ belongs to the subgroup $\Aut_r(X)$ and hence $A=BC\in \Proj(X)\cdot\Aut_r(X)$.

It remains to prove that the function $I:\Aut_r(X)\to\Aut(\bar\IF_X)$, $I:A\mapsto\ddddot A$, is a group isomorphism. By Proposition~\ref{p:Aut(X)->Aut(FX)}, the function $T:\Aut(X)\to\Aut(\bar\IF_X)$ is a group homomorphism whose kernel coincides with the subgroup $\Proj_S(X)$ of proscalarities, which coincides with the group of projectivities $\Proj(X)$, by Theorem~\ref{t:P=>p=pp=ps}. Since $\Proj(X)\cap\Aut_r(X)=\{1_X\}$, the restriction $I\defeq T{\restriction}_{\Aut_r(X)}$ is injective. 

It remains to prove that the homomorphism $I:\Aut_r(X)\to\Aut(\bar \IF_X)$ is surjective. Take any automorphism $\bar J:\bar\IF_X\to\bar\IF_X$ of the profield $\bar\IF_X$. Then $J\defeq \bar J\setminus\{(\infty,\infty)\}$ is an automorphism of the proscalar field $\IF_X\defeq\bar\IF_X\setminus\{\infty\}$.  Let $o\in\dom[r]$ be the unique fixed point of the projective reper $r$. Consider the hyperplane $H\defeq \overline{\rng[r]\setminus\{o\}}$ and the affine subliner $A\defeq X\setminus H$. By Proposition~\ref{p:Pappian-minus-flat}, the affine subliner $A$ of $X$ is Pappian. By Proposition~\ref{p:Pappus=>APA} and  Theorem~\ref{t:Hessenberg-affine}, the Pappian affine liner $A$ is Desarguesian. By Theorem~\ref{t:Papp<=>Des+RX}, the scalar corps $\IR_A$ of the Pappian affine liner $A$ is a field and hence $Z(\IR_A)=\IR_A$. By Corollary~\ref{c:FX=Z(RA)}, the function $\Psi_A:\IF_X\to Z(\IR_A)=\IR_A$ is a field isomorphism. Then $\Psi_AJ\Psi_A^{-1}$ is an automorphism of the field $\IR_A$. Observe that $\dom[r]$ is a maximal independent set in the affine space $A$. By Theorem~\ref{t:extension-iso}, there exists an automorphism $\varphi:A\to A$ such that $\varphi(x)=x$ for all $x\in\dom[r]$ and $\dddot\varphi=\Psi_AJ\Psi_A^{-1}$. By Theorem~\ref{t:extend-isomorphism-to-completions}, the automorphism $\varphi$ extends to an automorphism $\Phi:X\to X$ of the projective completion $X$ of $A$. Observe that every point $e\in \dom[r]\setminus\{o\}$, the point $r(e)$ is a unique point of the intersection $\Aline oe\cap H$. Then $\{\Phi(e)\}=\Phi[\Aline oe\cap H]=\overline{\{\Phi(o),\Phi(e)\}}\cap \Phi[H]=\Aline oe\cap H=\{e\}$ and hence $\Phi(x)=x$ for all $x\in\dom[r]\cup\rng[r]$, which means that $\Phi\in\Aut_r(X)$. It follows from $\Phi{\restriction}_A=\varphi$ that $\Psi_A\!\!\ddddot{\,\Phi}\Psi_A^{-1}=\dddot\varphi=\Psi_AJ\Psi_A^{-1}$ and hence $\ddddot {\,\Phi}{\restriction}_{\IF_X}=J$ and $\ddddot {\,\Phi}=\bar J$, witnessing that the injective homomorphism $I:\Aut_r(X)\to\Aut(\bar\IF_X)$ is surjective and hence $I$ is a group isomorphism.
\end{proof}

Theorem~\ref{t:Aut(X)=Proj(X)xAut(FX)} implies the following description of the automorphism group of a Pappian projective  space. 

\begin{corollary}\label{c:Aut(X)=Proj(X)xAut(FX)} The automorphism group $\Aut(X)$ of a Pappian projective space $X$ is a semidirect product $\Proj(X)\rtimes\Aut(\bar\IF_X)$ of the groups $\Proj(X)$ and $\Aut(\bar\IF_X)$.
\end{corollary}

\begin{remark} By Theorem~\ref{t:Des-projspace<=>}, every Desarguesian projective space  $X$ of finite rank $\|X\|=n$ is isomorphic to the projective space $\mathbb P\IR_X^n$. For a field $F$, the group of projective automorphisms of the projective space $\mathbb P F^n$ is isomorphic to the \index{projective linear group}\defterm{projective linear group} $PGL(n,F)\defeq GL(n,F)/\mathcal Z(GL(n,F))$, where $GL(n,F)$ is the group of  invertible $n\times n$-matrices with coefficients in the field $F$, and $\mathcal Z(GL(n,F))=\{c \cdot I:c\in F\}$ is the center of the group $GL(n,F)$. If the field $F$ is finite, then the group $PGL(n,F)$ is denoted by $PGL(n,q)$ where $q=|F|$ (because fields of the same cardinality are isomorphic). In 1870 Camille \index[person]{Jordan}Jordan\footnote{{\bf Marie Ennemond Camille Jordan} (1838 -- 1922) was a French mathematician, known both for his foundational work in group theory and for his influential Cours d'Analyse. Jordan was born in Lyon and educated at the \'Ecole polytechnique. He was an engineer by profession; later in life he taught at the \'Ecole polytechnique and the Coll\`ege de France, where he had a reputation for eccentric choices of notation.
He is remembered now by name in a number of results:
(1)   The {\em Jordan curve theorem}, a topological result required in complex analysis; (2) The {\em Jordan normal form} and the {\em Jordan matrix} in linear algebra. (3) In mathematical analysis, {\em Jordan measure} is an area measure that predates measure theory. (4) In group theory, the {\em Jordan--H\"older theorem} on composition series is a basic result. (5) {\em Jordan's theorem} on finite linear groups. \newline
Jordan's work did much to bring Galois theory into the mainstream. He also investigated the Mathieu groups, the first examples of sporadic groups. His ``Trait\'e des substitution'', on permutation groups, was published in 1870; this treatise won for Jordan the 1870 prix Poncelet. He was an Invited Speaker of the ICM in 1920 in Strasbourg.} proved that for every finite field $F$ of cardinality $|F|\notin\{2,3\}$ and every number $n\ge 2$, the projective linear group $PGL(n,F)$ is simple (= has only two normal subgroups).
\end{remark}

\section{Automorphisms of finite Pappian projective spaces}

Applying Corollary~\ref{c:Aut(X)=Proj(X)xAut(FX)} and Theorem~\ref{t:Aut(F)} we can calculate the cardinality of the automorphism group $\Aut(X)$ of a finite Pappian projective space $X$ and its normal subgroup $\Proj(X)$, consisting of projective automorphisms.

\begin{theorem}\label{t:|Proj(X)|=} For every finite Pappian projective space $X$, the group $\Proj(X)$ of projective automorphisms of $X$ has cardinality $$|\Proj(X)|=\frac1{\ell-1}\prod_{k=0}^{\|X\|-1}(\ell^{\|X\|}-\ell^k),$$
where $\ell\defeq|\IF_X|=|X|_2-1$.
\end{theorem}

\begin{proof} Fix any projective reper $r$ in $X$ and let $d\defeq\dim(X)=\|X\|-1$. By Corollary~\ref{c:Max=dim}, every maximal independent set in $X$ has cardinality $\|X\|$, which implies $|\dom[r]|=\|X\|=d+1$.  By Corollary~\ref{c:Papp-iso-reper}, every projective automorphism $A\in\Proj(X)$ is uniquely determined by its restriction $A{\restriction}_{\dom[r]\cup\rng[M]}$. Moreover, for every projective reper $r'$ in $X$ and any bijective map $\varphi:\dom[r]\cup\rng[r]\to\dom[r']\cup\rng[r']$ with $\varphi\circ r=r'\circ\varphi$ there exists a projective automorphism $A\in\Proj(X)$ such that $\varphi\subseteq A$. Therefore, the number of projective automorphisms of $X$ is equal to the number of $(2d+1)$-tuples $(x_0,x_1,\dots,x_{d},y_1,\dots,y_d)\in X^{2d+1}$ such that the set $\{x_0,x_1,\dots,x_{d}\}$ is maximal independent  in $X$ and $y_i\in\Aline{x_0}{x_i}\setminus\{x_0,x_i\}$ for every $i\in\{1,\dots,d\}$. By Proposition~\ref{p:add-point-to-independent}, the number of such $(2d+1)$-tuples in equal to the product
$$|X|\cdot(|X|-|X|_1)\cdot(|X|-|X|_2)\cdots (|X|-|X|_{d})\cdot (|X|_2-2)^d.$$
Let $\ell\defeq|X|_2-1=|\IF_X|$ be the order of the Pappian projective space $X$. By Corollary~\ref{c:projective-order-n} $|X|_k=\frac{\ell^k-1}{\ell-1}=\sum_{i=0}^{k-1}\ell^i$ for every $k\in\{1,\dots,d\}$. Then
$$|\Proj(X)|=(\ell-1)^d\prod_{k=0}^{d}\frac{\ell^{d+1}-\ell^k}{\ell-1}=\frac1{\ell-1}\prod_{k=0}^{\|X\|}(\ell^{\|X\|}-\ell^k).$$
\end{proof}

\begin{theorem}\label{t:proj-|Aut|} For every finite Pappian projective space $X$, the group $\Aut(X)$ of  automorphisms of $X$ has cardinality $$|\Aut(X)|=|\Aut(\IF_X)|\cdot|\Proj(X)|=\frac{n}{\ell-1}\cdot \prod_{k=0}^{\|X\|-1}(\ell^{\|X\|}-\ell^k),$$
where $r\defeq\|X\|$ and $\ell\defeq |\IF_X|=p^n$ for some prime number $p$ and some $n\in\IN$.
\end{theorem}

\begin{proof} By Corollary~\ref{c:Aut(X)=Proj(X)xAut(FX)}, the automorphism group $\Aut(X)$ of the Pappian projective space $X$ is a semidirect product $\Proj(X)\rtimes\Aut(\bar\IF_X)$ and hence $|\Aut(X)|=|\Aff(X)|\cdot|\Aut(\bar\IF_X)|$. Since the affine space $X$ is finite, so is its profield $\bar\IF_X$ and the field $\IF_X\defeq\bar\IF_X\setminus\{\infty\}$. Then $|\IF_X|=p^n$ for some $n\in\IN$ and some prime number $p$, equal to the characteristic of the field $\IF_X$. By Theorem~\ref{t:Aut(F)}, $|\Aff(\IF_X)|=n$. Since every automorphism of the profield $\bar\IF_X$ is uniquely determined by its restriction to $\IF_X$, the automorphism groups $\Aut(\bar\IF_X)$ and $\Aut(\IF_X)$ are canonically isomorphic and hence $|\Aut(\bar\IF_X)|=|\Aut(\IF_X)|=n$. Applying Theorem~\ref{t:|Proj(X)|=}, we conclude that
$$|\Aut(X)|=|\Aut(\bar\IF_X)|\cdot|\Proj(X)|=\frac{n}{\ell-1}\cdot\prod_{k=0}^{d-1}(r^d-r^k),$$where $d\defeq\|X\|-1$ and $\ell\defeq|\IF_X|=p^n=|X|_2-1$.
\end{proof}

\section{Automorphisms of finite Desarguesian proaffine spaces}

By Corollary~\ref{c:Des-proaffine=>point-homogeneous}, every Desarguesian proaffine space is point-homogeneous. In this section we calculate the cardinality of the automorphism group $\Aut(X)$ of a finite Desarguesian proaffine space.

By Corollary~\ref{c:Desarg=>compreg}, every Desarguesian proaffine space is completely regular, and by Theorem~\ref{t:Desargues-completion}, its spread completion $\overline X$ is a Desarguesian projective space. By Theorem~\ref{t:extend-isomorphism-to-completions}, every automorphism $A:X\to X$ of a $X$ uniquely extends to an automorphism $\bar A:\overline X\to\overline X$ of the spread completion $\overline X$ of $X$. So, the automorphism group $\Aut(X)$ of $X$ can be identified with the subgroup $$\Aut(\overline X,\partial X)\defeq \{\bar A\in\Aut(\overline X): \bar A[\partial X]=\partial X\}.$$ This observation allows to calculate the cardinality of the automorphism group of any finite Desarguesian proaffine space.

\begin{theorem} The automorphism group $\Aut(X)$ of a finite Desarguesian proaffine space $X$ has cardinality 
$$|\Aut(X)|=\frac{n}{\ell-1}\cdot\prod_{k=0}^{\|\partial X\|-1}(\ell^{\|\partial X\|}-\ell^k)\cdot\prod_{k=\|\partial X\|}^{\|X\|-1}(\ell^{\|X\|}-\ell^k),$$
where $\ell\defeq |\IR_X|=p^n$ is the cardinality of the scalar corps $\IR_X\defeq\IR_{\overline X}$ of $X$, and $p$ is a prime number.
\end{theorem}

\begin{proof} By Corollary~\ref{c:Desarg=>compreg} and Theorem~\ref{t:Desargues-completion}, the Desarguesian proaffine space is completely regular and its spread completion $\overline X$ is a Desarguesian projective space. Since $X$ is finite, so is its spread completion $\overline X$. By Theorem~\ref{t:finite-Papp<=>Des}, the finite Desarguesian projective space $\overline X$ is Pappian, and by Theorem~\ref{t:Papp<=>Des+RX}, its scalar corps $\IR_{\overline X}$ is a field. So, the scalar corps $\IR_X\defeq\IR_{\overline X}$ of $X$ is a field, too. By Corollary~\ref{c:FX=Z(RA)}, the proscalar field $\IF_{\overline X}$ of the projective space $\overline X$ is isomorphic to its scalar field $\IR_{\overline X}$. The automorphism group $\Aut(X)$ of the proaffine space $X$ can be identified with the subgroup $$\Aut(\overline X,\partial X)\defeq \{\bar A\in\Aut(X):\bar A[\partial X]=\partial X\}$$ of the automorphism group $\Aut(\overline X)$ of the spread completion $\overline X$ of $X$. In the group $\Aut(\overline X,\partial X)$, consider the subgroup
$$\Proj(\overline X,\partial X)\defeq\Proj(\overline X)\cap\Aut(\overline X,\partial X).$$

Fix any projective reper $r$ in $\overline X$ such that the set $\dom[r]\cap\partial X$ is maximal independent in $\partial X$ and the origin $o=r(o)$ of $r$ does not belong to $\partial X$.  Consider the subgroup $$\Aut_r(\overline X)\defeq\{A\in\Aut(\overline X):\forall x\in\dom[r]\cup\rng[r]\;\;A(x)=x\}$$ and observe that $\Aut_r(\overline X)\subseteq \Aut(\overline X,\partial X)$. Theorem~\ref{t:Aut(X)=Proj(X)xAut(FX)} implies that $$\Aut(\overline X,\partial X)=\Aut_r(\overline X)\cdot \Proj(\overline X,\partial X)$$ and hence
$$|\Aut(X)|=|\Aut(\overline X,\partial X)|=|\Aut_r(\overline X)|\cdot|\Proj(\overline X,\partial X)|.$$ By Theorem~\ref{t:Aut(X)=Proj(X)xAut(FX)}, the group $\Aut_r(\overline X)$ is isomorphic to the groups $\Aut(\bar\IF_X)$, $\Aut(\IF_{\overline X})$ and $\Aut(\IR_X)$. By Theorem~\ref{t:Moore}, those groups have cardinality $n$ where $|\IR_X|=p^n$ for a prime number $p$. Let $\ell\defeq|\IR_X|=|\IF_{\overline X}|$.

By Corollary~\ref{c:Max=dim}, $|\dom[r]\cap\partial X|=\|\partial X\|$ and $|\dom[r]|=\|X\|$. Write the set $\dom[r]$ as $\{x_0,\dots,x_{\|\partial X\|-1},x_{\|\partial X\|},\dots,x_{\|X\|}\}$ where $\{x_0,\dots,x_{\|\partial X\|-1}\}=\dom[r]\cap\partial X$ and $x_{\|X\|}=o\notin\partial X$. For every $i\in\{0,\dots,\|X\|-1\}$, consider the point $y_i\defeq r(x_i)\in \Aline o{x_i}\setminus\{o,x_i\}$ and observe that $r(x_i)\notin\partial X$. Every projective automorphism $A\in \Proj(\overline X,\partial X)$ maps the set $\dom[r]\cup\rng[r]=\{x_0,\dots,x_{\|X\|},y_0,\dots,y_{\|X\|-1}\}$ onto the set $A[\dom[r]\cup\rng[r]]=\{x'_0,\dots,x'_{\|X\|},y'_0,\dots,y'_{\|X\|-1}\}$ such that $\{x'_0,\dots,x'_{\|\partial X\|-1}\}$ is a maximal independent set in the flat $\partial X$. Moreover, $A$ is uniquely determined by the points $x_0',\dots,x'_{\|X\|},y_0',\dots,y'_{\|X\|-1}$. This implies that
$$
\begin{aligned}
|\Proj(\overline X,\partial X)|&=|\partial X|\cdot(|\partial X|-1)\cdots (|\partial X|-|\overline X|_{\|\partial X\|-1})\cdot\\
&\hskip33pt\cdot(|X|-|X|_{\|\partial X\|})\cdots (|X|-|X|_{\|X\|-1})\cdot(|X|_2-2)^{\|X\|-1}\\
&=\prod_{k=0}^{\|\partial X\|-1}\frac{\ell^{\|\partial X\|}-\ell^k}{\ell-1}\cdot\prod_{k=\|\partial X\|}^{\|X\|-1}\frac{\ell^{\|X\|}-\ell^k}{\ell-1}\cdot(\ell-1)^{\|X\|-1}\\
&=\frac1{\ell-1}\cdot\prod_{k=0}^{\|\partial X\|-1}(\ell^{\|\partial X\|}-\ell^k)\cdot\prod_{k=\|\partial X\|}^{\|X\|-1}(\ell^{\|X\|}-\ell^k).
\end{aligned}  
$$
Then $$|\Aut(X)|=|\Aut_r(\overline X)|\cdot|\Proj(\overline X,\partial X)|=\frac{n}{\ell-1}\cdot \prod_{k=0}^{\|\partial X\|-1}(\ell^{\|\partial X\|}-\ell^k)\cdot\prod_{k=\|\partial X\|}^{\|X\|-1}(\ell^{\|X\|}-\ell^k).$$
\end{proof}

\chapter{The Bruck--Ryser Theorem}

In this section we present a (relatively simple) proof of the famous \index[person]{Ryser}Bruck--Ryser\footnote{{\bf Herbert John Ryser} (1923 -- 1985) was an American mathematician known for his influential work in combinatorics, finite geometry, and number theory. He earned his Ph.D. in 1948 from the University of Wisconsin--Madison under the supervision of Richard Bruck. Ryser made foundational contributions to the theory of Latin squares, incidence structures, and combinatorial designs. He is best known in finite geometry for the Bruck--Ryser theorem, which provides necessary conditions for the existence of certain finite projective planes. Ryser was also the author of the classic text 
``Combinatorial Mathematics'', which helped shape modern combinatorics as a field. He held academic positions at Ohio State University and the California Institute of Technology and was a key figure in the development of combinatorial mathematics in the mid-20th century.} Theorem \cite{BR1949}, which is a unique available theoretical result restricting possible orders of finite affine or projective planes. The proof of this theorem exploits the machinery of linear algebra and arithmetics, in particular, Euler's two-square and four-square theorems and the Lagrange's four-square Theorem.

\section{Euler's two-square theorem}

In this section we prove the Euler's two-square identity and its corollaries. This identity was used by \index[person]{Euler}Euler\footnote{{\bf Leonhard Euler} (1707 --1783) was a Swiss polymath who was active as a mathematician, physicist, astronomer, logician, geographer, and engineer. He founded the studies of graph theory and topology and made influential discoveries in many other branches of mathematics, such as analytic number theory, complex analysis, and infinitesimal calculus. He also introduced much of modern mathematical terminology and notation, including the notion of a mathematical function. He is known for his work in mechanics, fluid dynamics, optics, astronomy, and music theory.  He spent most of his adult life in Saint Petersburg, Russia, and in Berlin, then the capital of Prussia.
\\
Euler is credited for popularizing the Greek letter $\pi$ to denote the ratio of a circle's circumference to its diameter, as well as first using the notation 
$f(x)$ for the value of a function, the letter 
$i$ to express the imaginary unit $\sqrt{-1}$, the Greek letter 
$\Sigma$ to express summations. He gave the current definition of the constant $e$, the base of the natural logarithm, now known as Euler's number. Euler made contributions to applied mathematics and engineering, such as his study of ships which helped navigation, his three volumes on optics contributed to the design of microscopes and telescopes, and he studied the bending of beams and the critical load of columns.
\\
Euler is credited with being the first to develop graph theory (partly as a solution for the problem of the Seven Bridges of K\"onigsberg, which is also considered the first practical application of topology). He also became famous for, among many other accomplishments, solving several unsolved problems in number theory and analysis, including the famous Basel problem. Euler has also been credited for discovering that the sum of the numbers of vertices and faces minus the number of edges of a polyhedron equals 2, a number now commonly known as the Euler characteristic. In physics, Euler reformulated Isaac Newton's laws of motion into new laws in his two-volume work ``Mechanica'' to better explain the motion of rigid bodies. He contributed to the study of elastic deformations of solid objects. Euler formulated the partial differential equations for the motion of inviscid fluid, and laid the mathematical foundations of potential theory. 
} in his paper ``{\em Demonstratio theorematis Fermatiani omnem numerum primum formae $4n+1$ esse summam duorum quadratorum}'', written in 1749 and published in 1760.

The identity in the following theorem will be referred to as the \index{Euler's two-square identity}\defterm{Euler two-square identity}. It can be verified by routine calculations, which are left to the reader.

\begin{theorem}[Euler, 1749]\label{t:Euler-2square-id} For any real numbers $x_0,x_1,y_1,y_2$,
$$(x_0^2+x_1^2)\cdot(y_0^2+y_1^2)=(x_0y_0+x_1y_1)^2+(x_0y_1-x_1y_0)^2=(x_0y_0-x_1y_1)^2+(x_0y_1+x_1y_0)^2.$$
\end{theorem}

The Euler two-square identity has two important implications.

\begin{theorem}\label{t:Euler-2-squares} A number $n\in\IN$ belongs to the set $S=\{a^2+b^2:a,b\in\IZ\}$ if and only if $n\IZ\cap S\ne\varnothing$.
\end{theorem}

\begin{proof} The ``only if'' part is trivial. To prove the ``if'' part, assume that $n\IZ\cap S\ne\varnothing$ and take the smallest number $m\in\IN$ such that $nm\in S$ and hence $nm=x_0^2+x_1^2$ for some integer numbers $x_0,x_1$. We claim that $m=1$. To derive a contradiction, assume that $m>1$. For every $i\in\{0,1\}$, find a unique number $y_i\in(-\frac{m}2,\frac{m}2]\cap(x_i+m\IZ)$. Then $y_0^2+y_1^2\in x_0^2+x_1^2+m\IZ=nm+m\IZ=m\IZ$. It follows that $y_0^2+y_1^2\le 2\big(\frac{m}2\big)^2=\frac12m^2$ and hence $y_0^2+y_1^2=mr$ for some  number $r\le \frac12m<m$. By Euler's two-square identity, $nmmr=(x_0^2+x_1^2)\cdot (y_0^2+y_1^2)=z_0^2+z_1^2$, where the numbers
$$
\begin{aligned}
z_0&=x_0y_0+x_1y_1\in x_0^2+x_1^2+m\IZ=pm+m\IZ=m\IZ\quad\mbox{and}\\
z_1&=x_0y_1-x_1y_0=x_0x_1-x_1x_1+m\IZ=m\IZ
\end{aligned}
$$
are divisible by $m$ and hence $nr=\big(\frac {z_0}m\big)^2+\big(\frac{z_1}m\big)^2\in S$, which contradicts the minimality of $m$. This contradiction shows that $m=1$ and hence $n\in S$.
\end{proof}

\begin{corollary}\label{c:Euler-2-squares} A number $n\in \IN$ is the sum of two squares of integer numbers if and only if $n$ is the sum of two squares of rational numbers.
\end{corollary}

\begin{proof} The ``only if'' part is trivial. To prove the ``if'' part, assume that a number $n\in\IN$ is equal to the sum $a^2+b^2$ of the squares of two rational numbers $a,b$. Find a number $c\in\IN$ such that the rational numbers $ca$ and $cb$ are integer. Then $nc^2=(ca)^2+(cb)^2$. Applying Theorem~\ref{t:Euler-2-squares}, we conluclude that $n$ is the sum of two squares of integers.
\end{proof}

\section{Euler's four-square identity}

The identity in the following theorem will be refreed to as the \index{Euler's four-square identity}\defterm{Euler's four-square identity}. 

\begin{theorem}[Euler, 1749]\label{t:4square-id} For any real numbers $x_0,x_1,x_2,x_3,y_0,y_1,y_2,y_3$,
$$(x^2_0+x_1^2+x_2^2+x_3^2)\cdot(y_0^2+y_1^2+y_2^2+y_3^2)=z_0^2+z_1^2+z_2^2+z_3^2,$$
where 
$$
\begin{aligned}
z_0&=x_0y_0+x_1y_1+x_2y_2+x_3y_3,\\
z_1&=x_0y_1-x_1y_0+x_2y_3-x_3y_2,\\
z_2&=x_0y_2-x_1y_3-x_2y_0+x_3y_1,\\
z_3&=x_0y_3+x_1y_2-x_2y_1-x_3y_0.
\end{aligned}
$$
\end{theorem}

\begin{proof}Observe that
$$
\begin{aligned}
&\sum_{i\in 4}z_i^2-\sum_{i\in 4}x_i^2\cdot\sum_{j\in 4}y_j^2\\
&=(x_0y_0+x_1y_1+x_2y_2+x_3y_3)^2+
(x_0y_1-x_1y_0+x_2y_3-x_3y_2)^2\\
&+(x_0y_2-x_1y_3-x_2y_0+x_3y_1)^2+
(x_0y_3+x_1y_2-x_2y_1-x_3y_0)^2-\sum_{i\in 4}x_i^2\cdot\sum_{j\in 4}y_j^2\\
&=x_0^2y_0^2{+}x_1^2y_1^2{+}x_2^2y_2^2{+}x_3^2y_3^2+2(x_0y_0x_1y_1{+}x_0y_0x_2y_2{+}x_0y_0x_3y_3{+}x_1y_1x_2y_2{+}x_1y_1x_3y_3{+}x_2y_2x_3y_3)\\
&+x^2_0y^2_1{+}x_1^2y_0^2{+}x_2^2y_3^2{+}x_3^2y_2^2+
2({-}x_0y_1x_1y_0{+}x_0y_1x_2y_3{-}x_0y_1x_3y_2{-}x_1y_0x_2y_3{+}x_1y_0x_3y_2{-}x_2y_3x_3y_2)\\
&+x_0^2y_2^2{+}x_1^2y_3^2{+}x_2^2y_0^2{+}x_3^2y_1^2+2({-}x_0y_2x_1y_3{-}x_0y_2x_2y_0{+}x_0y_2x_3y_1{+}x_1y_3x_2y_0{-}x_1y_3x_3y_1{-}x_2y_0x_3y_1)\\
&+x_0^2y_3^2{+}x_1^2y_2^2{+}x_2^2y_1^2{+}x_3^2y_0^2+2(x_0y_3x_1y_2{-}x_0y_3x_2y_1{-}x_0y_3x_3y_0{-}x_1y_2x_2y_1{-}x_1y_2x_3y_0{+}x_2y_1x_3y_0)\\
&-\sum_{i,j\in 4}x_i^2y_j^2\\
&=2({\color{red}x_0y_0x_1y_1}+{\color{blue}x_0y_0x_2y_2}+{\color{cyan}x_0y_0x_3y_3}+{\color{magenta}x_1y_1x_2y_2}+{\color{brown}x_1y_1x_3y_3}+{\color{olive}x_2y_2x_3y_3})\\
&+2(-{\color{red}x_0y_1x_1y_0}+{\color{orange}x_0y_1x_2y_3}-{\color{purple}x_0y_1x_3y_2}-{\color{teal}x_1y_0x_2y_3}+{\color{gray}x_1y_0x_3y_2}-{\color{olive}x_2y_3x_3y_2})\\
&+2(-x_0y_2x_1y_3-{\color{blue}x_0y_2x_2y_0}+{\color{purple}x_0y_2x_3y_1}+{\color{teal}x_1y_3x_2y_0}-{\color{brown}x_1y_3x_3y_1}-{\color{violet}x_2y_0x_3y_1})\\
&+2(x_0y_3x_1y_2-{\color{orange}x_0y_3x_2y_1}-{\color{cyan}x_0y_3x_3y_0}-{\color{magenta}x_1y_2x_2y_1}-{\color{gray}x_1y_2x_3y_0}+{\color{violet}x_2y_1x_3y_0})\\
&=0.
\end{aligned}
$$
\end{proof}

\begin{remark} The Euler's four-square identity can be written in the matrix form as $$\sum_{i\in 4}x_i^2\cdot\sum_{j\in 4}y_j^2=\sum_{k\in 4}z_k^2,$$ where 
$$(z_0,z_1,z_2,z_3)=(x_0,x_1,x_2,x_3)\cdot
\left(\begin{array}{cccc}
y_0&y_1&y_2&y_3\\
y_1&-y_0&-y_3&y_2\\
y_2&y_3&-y_0&-y_1\\
y_3&-y_2&y_1&-y_0
\end{array}\right)
$$
\end{remark}

\section{The Lagrange four-square theorem}

In this section we present a proof of the classical \index[person]{Lagrange}Largange\footnote{{\bf Joseph-Louis Lagrange} (1736 -- 1813) was an Italian--French mathematician, physicist and astronomer, naturalized French. He made significant contributions to the fields of analysis, number theory, and both classical and celestial mechanics.\\
In 1766, on the recommendation of Leonhard Euler and d'Alembert, Lagrange succeeded Euler as the director of mathematics at the Prussian Academy of Sciences in Berlin, Prussia, where he stayed for over twenty years, producing many volumes of work and winning several prizes of the French Academy of Sciences. Lagrange's treatise on analytical mechanics ``M\'ecanique analytique'', written in Berlin and first published in 1788, offered the most comprehensive treatment of classical mechanics since Isaac Newton and formed a basis for the development of mathematical physics in the nineteenth century. \\In 1787, at age 51, he moved from Berlin to Paris and became a member of the French Academy of Sciences. He remained in France until the end of his life. He was instrumental in the decimalisation process in Revolutionary France, became the first professor of analysis at the \'Ecole Polytechnique upon its opening in 1794, was a founding member of the Bureau des Longitudes, and became Senator in 1799. \\
Lagrange was one of the creators of the calculus of variations, deriving the Euler–Lagrange equations for extrema of functionals. He extended the method to include possible constraints, arriving at the method of Lagrange multipliers. Lagrange invented the method of solving differential equations known as variation of parameters, applied differential calculus to the theory of probabilities and worked on solutions for algebraic equations. He proved that every natural number is a sum of four squares. His treatise Theorie des fonctions analytiques laid some of the foundations of group theory, anticipating Galois. In calculus, Lagrange developed a novel approach to interpolation and Taylor's theorem. He studied the three-body problem for the Earth, Sun and Moon (1764) and the movement of Jupiter's satellites (1766), and in 1772 found the special-case solutions to this problem that yield what are now known as Lagrangian points. Lagrange is best known for transforming Newtonian mechanics into a branch of analysis, Lagrangian mechanics. He presented the mechanical ``principles'' as simple results of the variational calculus.} \index{Lagrange's four-square Theorem}four-square Theorem.

\begin{theorem}[Lagrange, 1770]\label{t:4Lagrange} For every number $n\in\IN$, there exist integer numbers $a,b,c,d$ such that $n=a^2+b^2+c^2+d^2$.
\end{theorem}

\begin{proof} The proof is by induction on $n$. Let $S\defeq\{a^2+b^2+c^2+d^2:a,b,c,d\in\IZ\}$. The Euler four-square identity implies that $S\cdot S\defeq\{a\cdot b:a,b\in S\}\subseteq S$. It is clear that $\{1,2,3,4\}\subseteq S$. Assume that for some number $p\ge 5$ we know that every positive integer $k<p$ belongs to the set $S$. If $p=m\cdot k$ for some $m,k<p$, then $p=m\cdot k\in S\cdot S\subseteq S$, by the inductive assumption. So, assume that the number $p$ is prime.

Observe that for any distinct numbers $x,y\in\{0,1,\dots,\frac12(p-1)\}$, we have $0<|x-y|\le\frac12(p+1)<p$ and $0<x+y\le p-2<p$, which implies that  the difference $x^2-y^2=(x-y)(x+y)$ is not divisible  by the prime number $p$.
Then the set $A\defeq\{x^2+p\IZ:0\le x\le\frac12(p-1)\}$ in the field $\IZ_p=\IZ/p\IZ$ has cardinality $|A|=\frac12(p+1)$. Since $|A|+|-A-1|=p+1>|\IZ_p|$, there exists a point $x\in A\cap(-A-1)$ and hence $a^2+p\IZ=x+p\IZ=-b^2-1+p\IZ$ for some numbers $a,b\in\{0,\dots,\frac12(p-1)\}$. Therefore $a^2+b^2+1^2+0^2\in p\IZ\cap S\ne \varnothing$. It follows from $p\ge 5$ that $a^2+b^2+1\le 1+2(\frac{p+1}2)^2=1+\frac{1+2p+p^2}2<p^2$ and hence $pk=a^2+b^2+1\in S$ for some $k<p$. Let $m\le k<p$ be the smallest number such that $pm\in S$ and hence $pm=\sum_{i\in 4}x_i^2$ for some integer numbers $x_0,x_1,x_2,x_3$. We claim that $m=1$. To derive a contradiction, assume that $m>1$.

 Find unique integer numbers $y_0,y_1,y_2,y_3\in(-\frac12m,\frac12m]$ such that $y_i\in x_i+m\IZ$ for all $i\in 4$. It follows that $\sum_{i\in 4}y_i^2+m\IZ=\sum_{i\in 4}x_i^2+m\IZ=pm+m\IZ=m\IZ$ and $\sum_{i\in 4}y_i^2\le 4\cdot(\tfrac12m)^2=m^2$. Assuming that $\sum_{i\in 4}y_i^2=m^2$, we conclude that $y_i=\frac12m$ for all $i\in 4$. For every $i\in 4$ find $z_i\in\IZ$ such that $x_i=y_i+mz_i=\frac12m+mz_i=m(\frac12+z_i)$. Then $pm=\sum_{i\in 4}x_i^2=\sum_{i\in4}(\frac12m+mz_i)^2=m^2\big(1+\sum_{i\in 4}(z_i+z_i^2)\big)$, which implies that the prime number $p$ is divisible by $m>1$ and hence $m=p$, which contradicts the choice of $m\le k<p$. This contradiction shows that $\sum_{i\in 4}y_i^2<m^2$ and hence $\sum_{i\in 4}y_i^2=mr\in S$ for some $r<m$. Applying the Euler four-square identity, we conclude that $pmmr=\sum_{i\in 4}x_i^2\cdot\sum_{i\in 4}y_i^2=\sum_{i\in 4}z_i^2$, where the numbers
$$
\begin{aligned}
z_0&=x_0y_0+x_1y_1+x_2y_2+x_3y_3\in  x_0^2+x_1^2+x_2^2+x_3^2+m\IZ=mp+m\IZ=m\IZ;\\
z_1&=x_1y_2-x_2y_1+x_3y_4-x_4y_3\in x_1x_2-x_2x_1+x_3x_4-x_4x_3+m\IZ=m\IZ;\\
z_2&=x_1y_3-x_2y_4-x_3y_1+x_4y_2\in x_1x_3-x_2x_4-x_3x_1+x_4x_2+m\IZ=m\IZ;\\
z_3&=x_1y_4+x_2y_3-x_3y_2-x_4x_1\in x_1x_4+x_2x_3-x_3x_2-x_4x_1+m\IZ=m\IZ;\\
\end{aligned}
$$ 
are divisible by $m$, which implies $pr=\sum_{i\in 4}\big(\tfrac{z_i}m)^2\in S$ and contradicts the minimality of $m$. This contradiction shows that $m=1$ and hence $p\in S$.
\end{proof}

\section{Perturbations of invertible matrices}

In this section we prove one result on the preservation of invertibility of a matrix by a perturbation with another matrix.

\begin{lemma}\label{l:perturbation-matrix} For every number $v\in\w$ and an indexed family of real numbers $(c_S)_{S\subseteq v}$ with $c_\emptyset\ne\varnothing$, there exists a vector $(x_i)_{i\in v}\in\{-1,1\}^v$ such that $\sum_{S\subseteq v}c_S\cdot\prod_{i\in S}x_i\ne 0$.
\end{lemma}

\begin{proof} This lemma will be proved by induction on $v$. For $v=0$, the unique vector $(x_i)_{i\in v}=\emptyset\in \{\emptyset\}=\{-1,1\}^0$ has the required property $\sum_{S\subseteq v}c_S\cdot\prod_{i\in S}x_i=c_\emptyset\ne 0$. Assume that the lemma has been proved for some number $v\in \w$. Take any indexed family of real numbers $(c_S)_{S\subseteq v+1}$ with $c_\emptyset\ne 0$. By the inductive assumption, there exists a vector $(x_i)_{i\in v}\in\{-1,1\}^v$ such that $c\defeq\sum_{S\subseteq v}c_S\cdot\prod_{i\in S}x_i\ne 0$. For every set $S\subseteq v$, consider the real number $c'_S\defeq C_{S\cup\{v\}}$ and let $c'\defeq \sum_{S\subseteq v}c_S'\cdot\prod_{i\in S}x_i$. Since $c\ne 0$, there exists a number $x_v\in\{-1,1\}$ such that $c+x_v{\cdot}c'\ne 0$. Then for the vector $(x_i)_{i\in v+1}$, we have
$$\sum_{S\subseteq v+1}c_S\cdot\prod_{i\in S}x_i=\sum_{S\subseteq v}c_S\cdot\prod_{i\in S}x_i+x_v\cdot\sum_{S\subseteq v}c_S'\cdot\prod_{i\in S}x_i=c+x_i\cdot c'\ne 0.$$
\end{proof}

The following theorem is the main result of this section.

\begin{theorem}\label{t:BR-matrix} For any $v\in\IN$, any invertible real matrix $(a_{ij})_{i,j\in v}$ and any real matrix $(b_{ij})_{i,j\in v}$, there exists a vector $(x_i)_{i\in v}\in\{-1,1\}^v$ such that the matrix $(a_{ij}-x_jb_{ij})_{i,j\in v}$ is invertible.
\end{theorem}

\begin{proof} For every vector $\boldsymbol x=(x_0,\dots,x_{v-1})\in \IR^v$, consider the matrix $C_{\boldsymbol x}\defeq(a_{ij}-x_jb_{ij})_{i,j\in v}$ and observe that its  determinant $|C_{\boldsymbol x}|$ is a polynomial of the variables $x_0,\dots,x_{v-1}$ such that $|C_{\boldsymbol x}|=\sum_{S\subseteq v}c_S\cdot\prod_{i\in S}x_i$ for some real numbers $c_S$, $S\subseteq v$. Since the matrix $C_{\boldsymbol 0}=(a_{ij})_{i,j\in v}$ is invertible, $c_\emptyset=|C_{\boldsymbol 0}|\ne 0$. By Lemma~\ref{l:perturbation-matrix}, there exists a vector $(x_i)_{i\in v}\in\{-1,1\}^v$ such that $0\ne \sum_{S\subseteq v}c_S\prod_{i\in S}x_i=|C_{\boldsymbol x}|$, witnessing that the matrix $C_{\boldsymbol x}=(a_{ij}-x_jb_{ij})_{i,j\in v}$ is invertible.
\end{proof}

\section{The Bruck--Ryser Theorem}

Now we have all tools for presenting a proof of the \index{Bruck--Ryser Theorem} Bruck--Ryser Theorem (which is a detalized and simplified version of Cameron's proof \cite[\S 9.8]{Cameron}). 

\begin{theorem}[Bruck--Ryser, 1949]\label{t:Bruck-Ryser} If $n\in (1+4\IZ)\cup(2+4\IZ)$ is an order of a projective plane, then $n\in\{a^2+b^2:a,b\in\IZ\}$.
\end{theorem}

\begin{proof} Let $(\Pi,\mathcal L)$ be a projective plane of order $n$. Consider the number $v\defeq n^2+n+1$ and observe that $n\in(1+4\IZ)\cup(2+4\IZ)$ implies  $v\in 3+4\IZ$. Corollary~\ref{c:projective-order-n} ensures that $|\Pi|=|\mathcal L|=v$ and hence the sets $\Pi$ and $\mathcal L$ can be written as $\Pi=\{p_i\}_{i\in v}$ and $\mathcal L=\{\ell_i\}_{i\in v}$. Let $A=(a_{ij})_{i,j\in v}$ be the incidence matrix of the projective plane $(\Pi,\mathcal L)$. This matrix consists of the elements 
$$a_{ij}\defeq\begin{cases} 1&\mbox{if $p_j\in \ell_i$};\\
0&\mbox{otherwise}.
\end{cases}
$$Since every line in $\Pi$ contains exactly $n+1$ points, each row and column of the matrix $A$ contains exactly $n+1$ units.
Since any two lines in the projective plane $\Pi$ have exactly one common point, 
$AA^T=nI+J$, where $I$ is the identity matrix and $J$ is the matrix with all elements  equal to $1$.

For a vector  $\boldsymbol x=(x_0,\dots,x_{v-1})\in\IQ^v$, let $$\|\boldsymbol x\|^2\defeq \boldsymbol x\cdot\boldsymbol x^T=\sum_{i\in v}|x_i|^2$$be the squared norm of $\boldsymbol x$. Then 
$$\|\boldsymbol xA\|^2=\boldsymbol x AA^T\boldsymbol x^T=\boldsymbol x(nI+J)\boldsymbol x^T=n\cdot (\boldsymbol x\cdot\boldsymbol x^T)+\boldsymbol xJ\boldsymbol x^T=n\cdot\|\boldsymbol x\|^2+|\Sigma\boldsymbol x|^2,$$ where $\Sigma\boldsymbol x\defeq\sum_{i\in v}x_i$. 
If $\|\boldsymbol x\|^2\ne 0$, then $\|\boldsymbol xA\|^2=n\cdot\|\boldsymbol x\|^2+|\Sigma\boldsymbol x|^2>0$, witnessing that the matrix $A$ is invertible.


By the Lagrange four-square Theorem~\ref{t:4Lagrange}, $n=b_0^2+b_1^2+b_2^2+b_3^2$ for some integer numbers $b_0,b_1,b_2,b_3$. The non-zero vector $\boldsymbol b\defeq(b_0,b_1,b_2,b_3)$ determines the matrix 
$$B_4\defeq \left(\begin{array}{cccc}
b_0&b_1&b_2&b_3\\
b_1&-b_0&-b_3&b_2\\
b_2&b_3&-b_0&-b_1\\
b_3&-b_2&b_1&-b_0
\end{array}\right)
$$
such that $$\|\boldsymbol yB_4\|^2=\|\boldsymbol y\|^2\cdot\|\boldsymbol b\|^2=\|\boldsymbol b\|^2\cdot\|\boldsymbol y\|^2=n\cdot\|\boldsymbol y\|^2$$ for every vector $\boldsymbol y\in\IQ^4$, according to Euler's four-square identity \ref{t:4square-id}. Let $B$ be the matrix consisting of $\frac14(v+1)$ copies of the matrix $B_4$ along the diagonal. More precisely, $B=(b_{ij})_{i,j\in v+1}$ where $(b_{4k+i,4k+j})_{i,j\in 4}=B_4$ for every $k\in \frac14(v+1)$ (since $v\in 3+4\IZ$, the number $\frac14(v+1)$ is integer). The matrix $B$ has the property $\|\boldsymbol zB\|^2=n\cdot\|\boldsymbol z\|^2$ for every vector $\boldsymbol z\in \IQ^{v+1}$. 
Then 
$$\|\boldsymbol xA\|^2+n=n\cdot(\|\boldsymbol x\|^2+1)+|\Sigma\boldsymbol x|^2=\|(\boldsymbol x\hat{\;}1)B\|^2+|\Sigma\boldsymbol x|^2,$$
where $\boldsymbol x\hat{\;}1\defeq(x_0,\dots,x_{v-1},1)$.

By Theorem~\ref{t:BR-matrix}, there exists a vector of signs $(s_0,\dots,s_{v-1})\in\{-1,1\}^v$ such that the matrix $C\defeq (a_{ij}-s_jb_{ij})_{i,j\in v}$ is invertible. Consider the vector $\boldsymbol y=(y_0,\dots,y_{v-1})\in\IQ^v$ with $y_j\defeq s_jb_{mj}$ for all $j\in v$. Since the matrix $C$ is invertible, there exists a vector $\boldsymbol x=(x_0,\dots,x_{v-1})\in\IQ^v$ such that $\boldsymbol xC=\boldsymbol y$. Let $x_v=1$. Then for every $j\in v$ we have $\sum_{i\in v}x_i(a_{ij}-s_jb_{ij})=y_j=s_jb_{mj}$ and hence $\sum_{i\in v}x_ia_{ij}=s_j\cdot(b_{mj}+\sum_{i\in v}x_ib_{ij})=s_j\sum_{i\in v+1}x_ib_{ij}$ and
$$\|\boldsymbol xA\|^2=\sum_{j\in v}\Big(\sum_{i\in v}x_ia_{ij}\Big)^2=\sum_{j\in v}\Big(s_j\cdot\sum_{i\in v+1}x_ib_{ij}\Big)^2=\|(\boldsymbol x\hat{\;}1)B\|^2-q^2,$$ where $q\defeq\sum_{i\in m+1}x_ib_{im}\in\IQ$.
Then the equality
$$\|\boldsymbol xA\|^2+n=\|(\boldsymbol x\hat{\;}1)B\|^2+|\Sigma(\boldsymbol x)|^2=\|\boldsymbol x A\|^2+q^2+|\Sigma\boldsymbol x|^2$$implies that $n=q^2+|\Sigma\boldsymbol x|^2$ is the sum of two squares of rational numbers. By Corollary~\ref{c:Euler-2-squares}, $n$ is the sum of two squares of integer numbers.
\end{proof}

\begin{corollary} If $n\in(1+4\IZ)\cup(2+4\IZ)$ is an order of an affine liner $X$ of rank $\|X\|\ge 3$, then $n=a^2+b^2$ for some integer numbers $a,b$.
\end{corollary}

\begin{proof} If $n=2$, then $n=1^2+1^2$ and we are done. So, assume that $n\ge 3$. It follows from $n\in(1+4\IZ)\cup(2+4\IZ)$ that $n\ge 5$. By Theorem~\ref{t:4-long-affine}, the $5$-long affine liner $X$ is regular. Since $\|X\|\ge 3$, the affine liner $X$ contains a plane $\Pi$. By Corollary~\ref{c:affine-spread-completion}, the affine regular plane $\Pi$ is completely regular and hence the spread  completion $\overline\Pi$ is a projective plane of order $n$. By Theorem~\ref{t:Bruck-Ryser}, $n=a^2+b^2$ for some integer numbers $a,b$.
\end{proof} 

Bruck--Ryser Theorem implies the following corollary (that will applied in the proof of Corollary~\ref{c:p5-Pappian}).   

\begin{corollary}\label{c:no6order} There is no affine or projective planes of order $6$.
\end{corollary}

\chapter{Based affine planes and their ternars}

By Corollary~\ref{c:affine-Desarguesian}, proaffine spaces $X$ of rank $\|X\|\ge 4$ are Desarguesian, so their structure is well described by standard tools of Linear Algebra. The case of non-Desagusian planes is more complicated and requires special algebraic tools, developed in this chapter. Such  tools include affine bases, ternars, triloops, biloops,  etc.
\smallskip

\section{Affine bases in affine planes}

\begin{definition} An \index{affine plane}\defterm{affine plane} is any affine space $X$ of dimension $\dim(X)=2$, i.e., a $3$-long affine regular liner $X$ of rank $\|X\|=3$.
\end{definition}

Theorem~\ref{t:Playfair<=>} implies that a liner $X$ is an affine plane if and only if $X$ is a Playfair liner of dimension $\dim(X)=2$. So, affine planes will be also called \index{Playfair plane}\defterm{Playfair planes}.

\begin{definition} An \index{affine base}\defterm{affine base} in an affine plane $\Pi$ is any ordered triple of points $uow\in \Pi^3$ such that $\|\{u,o,w\}\|=\|\Pi\|=3$. The points $u,o,w$ are called the \index{affine base!unit}\index{unit}\defterm{unit}, \index{affine base!origin}\index{origin}\defterm{origin}, and the \index{biunit}\index{affine base!biunit}\defterm{biunit} of the affine base, respectively. An affine plane $\Pi$ endowed with an affine base $uow$ is called a \index{based affine plane}\defterm{based affine plane}.
\end{definition}

Let $\Pi$ be a based affine plane, endowed with an affine base $uow$. By Theorem~\ref{t:parallelogram3+1}, there exists a unique point $e\in\Pi$ such that $\Aline ue\parallel \Aline ow$ and $\Aline we\parallel \Aline ou$. Therefore, $uowe$ is a parallelogram. The point $e$ is called the \index{diunit}\index{affine base!diunit}\defterm{diunit} of the affine base $uow$. The lines $\Aline ou$, $\Aline ow$, and $\Delta\defeq\Aline oe$ are called the \index{horizontal axis}\index{affine base!horizontal axis}\defterm{horizontal axis}, the \index{vertical axis}\index{affine base!vertical axis}\defterm{vertical axis} and the \index{diagonal}\index{affine base!diagonal}\defterm{diagonal} of the affine based affine plane $(\Pi,uow)$. 
The directions $\boldsymbol h\defeq(\Aline ou)_\parallel$, $\boldsymbol v\defeq(\Aline ow)_\parallel$, and $\Delta_\parallel=(\Aline oe)_\parallel$ are called the \index{horizontal direction}\index{direction!horizontal}\defterm{horizontal}, \index{vertical direction}\index{direction!vertical}\defterm{vertical}, and \index{diagonal direction}\index{direction!diagonal}\defterm{diagonal directions} on the based affine  plane $(\Pi,uow)$.
\smallskip

\begin{picture}(100,140)(-50,-15)

{\linethickness{0.8pt}
\put(0,0){\color{teal}\vector(1,0){100}}
\put(0,0){\color{cyan}\vector(0,1){100}}
\put(0,0){\color{red}\line(1,1){90}}
}

\put(105,-3){\color{teal}$\boldsymbol h$}
\put(-3,105){\color{cyan}$\boldsymbol v$}
\put(93,93){\color{red}$\Delta$}

\put(0,0){\circle*{3}}
\put(-6,-8){$o$}
\put(40,0){\circle*{3}}
\put(37,-8){$u$}
\put(0,40){\circle*{3}}
\put(-10,38){$w$}
\put(40,40){\circle*{3}}
\put(42,34){$e$}

\put(150,105){$o$, the origin}
\put(150,93){$u$, the unit}
\put(150,81){$w$, the biunit}
\put(150,69){$e$, the diunit}
\put(150,57){$\Aline ou$, the horizontal axis}
\put(150,45){$\Aline ow$, the vertical axis}
\put(150,33){$\Delta\defeq\Aline oe$, the diagonal}
\put(150,21){$\boldsymbol h\defeq(\Aline ou)_\parallel$, the horizontal direction}
\put(150,9){$\boldsymbol v\defeq(\Aline ow)_\parallel$, the vertical direction}
\put(150,-3){$\Delta_\parallel=(\Aline oe)_\parallel$, the diagonal direction}

\end{picture}
\smallskip

\begin{remark} An affine base $uow$ in an affine plane $\Pi$ is uniquely determined by the quadruple $oe\boldsymbol{hv}$ consisting of the origin $o$, the diunit $e$, the horizontal direction $\boldsymbol h$ and the vertical direction $\boldsymbol v$. The unit $u$ and biunit $w$ can be recovered as unique points $u\in \Aline o{\boldsymbol h}\cap\Aline e{\boldsymbol v}$ and $w\in \Aline o{\boldsymbol v}\cap\Aline e{\boldsymbol h}$. Here for a point $p\in\Pi$ and a spread of parallel lines ${\boldsymbol \delta}\in\partial\Pi$, we denote by $\Aline p{\boldsymbol\delta}$, the unique line in the direction ${\boldsymbol\delta}$ that contain the point $p$.
\end{remark}

\pagebreak

Since $\Pi$ is a Playfair plane and $\Delta\notin {\boldsymbol h}\cup {\boldsymbol v}$, for every point $p\in \Pi$ there exist unique points $p'\in\Delta\cap\Aline p{\boldsymbol h}$ and $p''\in\Delta\cap\Aline p{\boldsymbol v}$, called the \index{horizontal coordinate}\defterm{horizontal} and \index{vertical coordinate}\defterm{vertical coordinates} of the point $p$, respectively. Conversely, for any points $x,y\in\Delta$ there exists a unique point $p\in\Aline x{\boldsymbol v}\cap\Aline y{\boldsymbol h}$, which implies that the map $C:\Pi\to\Delta^2$, $C:p\mapsto p'p''$, assigning to every point $p\in \Pi$ the pair of its coordinates $p'p''$ is a bijective map from the affine plane $\Pi$ onto the square $\Delta^2$ of the diagonal $\Delta$ of the affine base $uow$. This bijective map $C:\Pi\to \Delta^2$ is called the \index{coordinate chart}\defterm{coordinate chart} of the based affine plane. 
\smallskip

The set $\Delta^2$ endowed with the line relation $\{Cxyz:xyz\in \Pi^3\;\wedge\;y\in\Aline xz\}$ is a Playfair plane, isomorphic to the Playfair plane $\Pi$ via the coordinate chart $C:\Pi\to \Delta^2$. Moreover, the coordinate chart $C$ maps the affine base $(u,o,w)$ in $\Pi$ onto the affine base $(eo,oo,oe)$ in the plane $\Delta^2$. The based affine plane $\Delta^2$ endowed with the canonical base $(eo,oo,oe)$ is called the \index{coordinate plane}\defterm{coordinate plane} of the based affine plane $(\Pi,uew)$. The coordinate chart $C:\Pi\to \Delta^2$ is the canonical isomorphism of the based affine planes $\Pi$ and $\Delta^2$.
\smallskip

\begin{picture}(100,130)(-30,-15)


{\linethickness{0.8pt}
\put(0,0){\color{teal}\line(1,0){80}}
\put(0,0){\color{cyan}\line(0,1){80}}
\put(0,0){\color{red}\line(1,1){70}}
\put(20,20){\color{cyan}\line(0,1){40}}
\put(60,60){\color{teal}\line(-1,0){40}}
}

\put(85,-3){$\Aline o{\boldsymbol h}$}
\put(-5,85){$\Aline o{\boldsymbol v}$}
\put(73,73){$\Delta$}

\put(20,20){\color{teal}\circle*{3}}

\put(22,14){$p'$}
\put(60,60){\color{cyan}\circle*{3}}
\put(62,54){$p''$}
\put(20,60){\circle*{3}}
\put(18,65){$p$}

\put(0,0){\circle*{3}}
\put(-6,-8){$o$}

\put(40,0){\circle*{3}}
\put(37,-8){$u$}
\put(40,40){\circle*{3}}
\put(40,33){$e$}
\put(0,40){\circle*{3}}
\put(-10,38){$w$}

\put(100,40){\vector(1,0){70}}
\put(130,43){$C$}

\put(40,90){$\Pi$}
\put(240,90){$\Delta^2$}

\put(220,60){\color{teal}\line(-1,0){20}}
\put(220,20){\color{cyan}\line(0,-1){20}}

{\linethickness{0.8pt}
\put(200,0){\color{teal}\line(1,0){80}}
\put(200,0){\color{cyan}\line(0,1){80}}
\put(200,0){\color{red}\line(1,1){70}}
\put(220,20){\color{cyan}\line(0,1){40}}
\put(260,60){\color{teal}\line(-1,0){40}}
}

\put(285,-3){$\Delta o$}
\put(195,85){$o\Delta$}

{\linethickness{1pt}
\put(220,-2){\color{teal}\line(0,1){4}}
}
\put(216,-9){$p'o$}
\put(220,20){\color{teal}\circle*{3}}
\put(223,15){$p'p'$}

{\linethickness{1pt}
\put(198,60){\color{teal}\line(1,0){4}}
}
\put(184,58){$op''$}
\put(260,60){\color{cyan}\circle*{3}}
\put(261,53){$p''p''$}
\put(220,60){\circle*{3}}
\put(215,65){$p'p''$}

\put(200,0){\circle*{3}}
\put(190,-9){$oo$}

\put(240,0){\circle*{3}}
\put(236,-9){$eo$}
\put(240,40){\circle*{3}}
\put(240,33){$ee$}
\put(200,40){\circle*{3}}
\put(187,38){$oe$}

\end{picture}

\begin{exercise} Let $(\Pi,uow)$ be a based affine plane and $C:\Pi\to \Delta^2$ be its coordinate chart. Show that
$$\begin{aligned}
C[\Aline o{\boldsymbol h}]&=\Delta o\defeq\{xo:x\in \Delta\},\\
C[\Aline o{\boldsymbol v}]&=o\Delta\defeq \{oy:y\in \Delta\}\quad\mbox{and}\\
C[\Delta]&=\{xx:x\in\Delta\}.
\end{aligned}
$$
\end{exercise}

Using the coordinates we can write down the equations of lines in a based affine plane. 

A line $L$ in a based affine plane $(\Pi,uow)$ is called
\begin{itemize}
\item \index{horizontal line}\index{line!horizontal}\defterm{horizontal} if $L\in {\boldsymbol h}\defeq(\Aline ou)_\parallel$;
\item \index{vertical line}\index{line!vertical}\defterm{vertical} if $L\in {\boldsymbol v}\defeq(\Aline ow)_\parallel$.
\end{itemize} 

\begin{exercise} Let $\Pi$ be a based affine plane. 
\begin{enumerate}
\item Show that every horizontal line $L\subseteq\Pi$ is equal to the set $\{p\in\Pi:p''=y\}$, where $y$ is the unique point of the intersection $L\cap\Delta$.
\item Show that every vectical line $L\subseteq\Pi$ is equal to the set $\{p\in\Pi:p'=x\}$, where $x$ is the unique point of the intersection $L\cap \Delta$.
\end{enumerate}
\end{exercise}

The equations of non-vertical lines can be written down using a special ternary operation $T_{uow}:\Delta^3\to\Delta$, which assigns to every triple $xab\in\Delta^3$ the unique point $y\in\Delta$ denoted by $x_{\times}a_{+}b$ such that $\Aline {ob}{xy}\subparallel\Aline {oo}{ea}$. The point $y$ is the unique common point of the vertical line $L_x\defeq\{p:p'=x\}$ and the unique line $L_{a,b}$ which contains the point $ob\in \Aline ov$ (with coordinates $o,b$) and is parallel to the line $L_{o,a}\defeq \Aline {oo}{ea}$ connnecting the points $o$ (with coordinates $o,o$) and the point $ea$ (with coordinates $e,a$). 

\begin{picture}(240,250)(-70,-15)

{\linethickness{0.75pt}
\put(0,0){\color{teal}\vector(1,0){240}}
\put(0,0){\color{cyan}\vector(0,1){210}}
\put(0,0){\color{red}\line(1,1){210}}
\put(0,0){\color{blue}\line(2,1){230}}
\put(0,90){\color{blue}\line(2,1){230}}
\put(0,90){\color{blue}\line(-2,-1){20}}
\put(60,60){\color{cyan}\line(0,-1){30}}
\put(30,30){\color{teal}\line(1,0){30}}
\put(0,90){\color{teal}\line(1,0){90}}
\put(120,120){\color{cyan}\line(0,1){30}}
\put(120,150){\color{teal}\line(1,0){30}}
}
\put(0,60){\color{teal}\line(1,0){60}}
\put(60,0){\color{cyan}\line(0,1){30}}

\put(245,-3){\color{teal}${\boldsymbol h}$}
\put(-3,215){\color{cyan}${\boldsymbol v}$}
\put(213,213){\color{red}$\Delta$}
\put(233,205){\color{blue}$L_{a,b}$}
\put(233,115){\color{blue}$L_{a,o}$}
\put(0,0){\circle*{3}}
\put(-5,-7){$o$}
\put(60,0){\circle*{3}}
\put(58,-8){$u$}
\put(0,60){\circle*{3}}
\put(-10,58){$w$}
\put(30,30){\circle*{3}}
\put(22,29){$a$}
\put(60,30){\circle*{3}}
\put(62,24){$ea$}
\put(60,60){\circle*{3}}
\put(56,63){$e$}
\put(90,90){\circle*{3}}
\put(86,92){$b$}
\put(120,120){\circle*{3}}
\put(122,115){$x$}
\put(150,150){\circle*{3}}
\put(152,145){$y=x_{\times}a_{+}b$}
\put(120,150){\circle*{3}}
\put(122,142){$xy$}
\put(0,90){\circle*{3}}
\put(2,81){$ob$}
\end{picture}

The definition of the ternary operation $T_{uow}$ ensures that the line $L_{a,b}$ is determined  by the equation $y=x_{\times}a_{+}b$, more precisely, $$L_{a,b}=\{p\in\Pi:p''=p'_\times a_+b\}=C^{-1}[\{xy\in\Delta^2:y=x_\times a_+b\}].$$

\begin{theorem}\label{t:aff-reper=>ternary-operation} For every based affine plane $(\Pi,uow)$, the ternary operation
$$T_{uow}:\Delta^3\to\Delta,\quad T_{uow}:xab\mapsto x_\times a_+b,$$ has the following four properties:
\begin{enumerate}
\item $x_\times o_+b=o_\times x_+b=b$ and $x_\times e_+o=e_\times x_+o=x$ for every points $x,b\in\Delta$;
\item for every points $a,x,y\in\Delta$, there exists a unique point  $b\in\Delta$ such that $x_\times a_+b=y$;
\item for every points $a,b,c,d\in\Delta$ with $a\ne c$, there exists a unique point $x\in\Delta$ such that $x_\times a_+b=x_\times c_+d$;
\item for every points $\check x,\check y,\hat x,\hat y\in\Delta$ with $\check x\ne \hat x$, there exist unique points $a,b\in \Delta$ such that $\check x_\times a_+b=\check y$ and $\hat x_\times a_+b=\hat y$.
\end{enumerate}
\end{theorem}

\begin{proof} Let us recall that for every points $a,b\in\Delta$, the equation $y=x_\times a_+b$ determines the line $L_{a,b}$ that contains the point $ob$ and is parallel to the line $L_{a,o}=\Aline {oo}{ea}$.
\smallskip

1. Since the line $L_{o,b}$ is horizontal and contains the point $ob$ (with coordinates $o,b$), we conclude that $x_\times o_+b=b$ for all $x\in\Delta$. On the other hand, for every $x\in \Delta$, the definition of the point $y\defeq o_\times x_+b$ ensures that the point $oy$ belongs to the intesection $L_{x,b}\cap Y$ which coincides with the point $ob$. Therefore, $o_\times x_+b=y=b$.
\smallskip

 The equality $L_{e,o}=\Delta$ implies $x_\times e_+o=x$ for all $x\in\Delta$. On the other hand, the definition of $y\defeq e_\times x_+o$ ensures that $ey\in L_{x,o}=\Aline {oo}{ex}$ and hence $e_\times x_+o=y=x$.  
\smallskip

2. Given any points $a,x,y\in\Delta$, consider the unique line $L$ that contains the point $xy$ (with coordinates $x,y$) and is parallel to the line $L_{a,o}=\Aline {oo}{ea}$. Since the line $L_{a,o}$ is not vertical, the line $L\parallel L_{a,o}$ has a common point $\beta\in L\cap Y$ with the vertical line $Y$. For the point $\beta$ there exists a unique point $b\in\Delta$ such that $\Aline \beta b\subparallel X$. The definition of the ternary operation  $T_{uow}$ ensures that $L=L_{a,b}$ and hence $y=x_\times a_+b$.
\smallskip

3. For every points $a,b,c,d\in\Delta$ with $a\ne c$, the lines $L_{a,b}$ and $L_{c,d}$ are not parallel and hence have a unique common point whose coordinates $x,y$ satisfy the equations $y=x_\times a_+b$ and $y=x_\times c_+d$ of the lines $L_{a,b}$ and $L_{c,d}$.
\smallskip

4. Given four points $\check x,\check y,\hat x,\hat y\in\Delta$ with $\check x\ne \hat x$, find a unique line $L$ containing the points $\check x\check y$ and $\hat x\hat y$ (with coordinates $\check x,\check y$ and $\hat x,\hat y$). Since $\check x\ne \hat x$, the line $L$ is not vertical and hence has a unique common point $ob$ with the vectical line $Y$. Let $a\in \Delta$ be the unique point such that  
the line $L_{a,o}$ is parallel to the line $L$. Then $L=L_{a,b}$ and hence the points $a,b$ are unique elements of $\Delta$ that satisfy 
the equations  $\check x{\cdot}a{+}b=\check y$ and $\hat x{\cdot}a{+}b=\hat y$, which uniquely determine the line $L_{a,b}$.
\end{proof}

\section{Ternars and their coordinate planes}

The properties (1)--(4) of the ternary operation $T_{uow}:\Delta^3\to \Delta$ appearing in Theorem~\ref{t:aff-reper=>ternary-operation}  are taken as axioms of an algebraic structure, called a ternar. This algebraic structure was introduced by Hall \cite{Hall1943} in 1943, who called them \defterm{ternary rings}. The term \defterm{ternar} was suggested by Skornyakov \cite{Skornyakov} and we prefer to use the terminology of Skornyakov (because Hall's ternary rings of Hall are not rings at all). 

\begin{definition}[Hall, 1943; Skornyakov, 1951] A \index{ternar}\defterm{ternar} is a set $R$ endowed with a ternary operation 
$$T:R^3\to R,\quad T:xab\mapsto x_{\times}a_{+}b,$$satisfying the following four axioms:
\begin{itemize}
\item[\textup{\sf(T1)}] there exist distinct elements $0,1\in R$ such that $x_\times 0_+b=0_\times x_+b=b$ and\newline $x_\times 1_+0=1_\times x_+0=x$ for all $x,b\in R$;
\item[\textup{\sf(T2)}] for all elements $a,x,y\in R$, there exists a unique element $b\in R$ such that $x_\times a_+b=y$;
\item[\textup{\sf(T3)}] for all elements $a,b,c,d\in R$ with $a\ne c$, there exists a unique element $x\in R$ such that $x_\times a_+b=x_\times c_+d$;
\item[\textup{\sf(T4)}] for all elements $\check x,\check y,\hat x,\hat y\in R$ with $\check x\ne \hat x$, there exist unique elements $a,b\in R$ such that $\check x_\times a_+b=\check y$ and $\hat x_\times a_+b=\hat y$.
\end{itemize}
The elements $0,1$ appearing in axiom {\sf(T1)} are called the \index{ternar!zero}\defterm{zero} and the \index{ternar!unit}\defterm{unit} of the ternar $R$.
\end{definition}

The following proposition shows that the zero and the unit of a ternar are uniquely determined by the axioms {\sf(T1)} and {\sf(T3)}.
  
\begin{proposition}\label{p:TR-01-unique} For every ternar $R$, the elements $0$ and $1$ satisfying the axiom {\sf(T1)} are unique. 
\end{proposition}

\begin{proof} Take any points $0',1'\in R$ such that
$$x_\times 0'_+b=0'_\times x_+b=b\quad\mbox{and}\quad x_\times 1'_+0'=1'_\times x_+0'=x$$ for all $x,b\in R$.
Assuming that $0\ne 0'$, we conclude that $x_\times 0_+0=0=x_\times 0'_+0$ for all $x\in R$, which contradicts the axiom {\sf(T3)}. This contradiction shows that $0'=0$. Assuming that $1\ne 1'$, we conclude that $x_\times 1'_+0'=x=x_\times 1_+0$ for all $x\in R$, which contradicts the axiom {\sf(T3)}. This contradiction shows that $1'=1$. 
\end{proof}

\begin{remark}\label{rem:ternar-T2T3T4} Since for a ternar $R$ its zero and unit are uniquely determined by the ternary operation, they can be added to the structure of a ternar, which can be thought as an algebraic structure $(R,T,0,1)$ consisting the the underlying set $R$ of the ternar, the ternary operation $T:R^3\to R$, and two constants $0$ and $1$. In fact, any ternar $R$ has much more natural operations than just the ternary operation $T$ and two constants $0$ and $1$. In particular, the axiom {\sf(T2)}  determines a unique ternary operation $T_2:R^3\to R$ assigning to every triple $(a,x,y)\in R^3$ a unique element $b=T_2(a,x,y)\in R$ such that $T(x,a,b)=y$. The axiom {\sf(T3)} determines a unique quaternary operation $T_3:R^4_{\ne}\to R$ on the set $R^4_{\ne}\defeq\{(x_1,x_2,x_3,x_4)\in R^4:x_1\ne x_3\}$, assigning to every quadruple $(a,b,c,d)\in R^4_{\ne}$ a unique point $x=T_3(a,b,c,d)\in R$ such that $T(x,a,b)=T(x,c,d)$.  The axiom {\sf(T4)} determines two quaternary operations $T_4',T_4'':R^4_{\ne}\to R$ assigning to every quadruple $(x,y,x',y')\in R^4_{\ne}$ unique elements $a=T_4'(x,y,x',y')$ and $b=T_4''(x,y,x',y')$ such that $T(x,a,b)=y$ and $T(x',a,b)=y'$. The operations $T,T_2,T_3,T_4',T_4''$ satisfy the following identities
\begin{itemize}
\item[(T1)] $T(x,0,b)=b=T(0,x,b)$ and $T(x,1,0)=x=T(1,x,0)$,
\item[(T2)] $T(x,a,T_2(a,x,y))=y$,
\item[(T3)] $T(T_3(a,b,c,d),a,b)=T(T_3(a,b,c,d),c,d)$,
\item[(T4)] $T(x,T_4'(x,y,x',y'),T_4''(x,y,x',y'))=y$ and $T(x',T_4'(x,y,x',y'),T_4''(x,x',y,y'))=y'$,
\end{itemize}
for all elements $a,b,c,d,x,y,x',y'\in R$ such that $a\ne c$ and $x\ne x'$.
Such an extended representation of ternars is useful for analyzing the structure of their subternars.
\end{remark}

For finite ternars, the axiom {\sf(T4)} follows from the axioms {\sf(T1)}--{\sf(T3)}.

\begin{proposition}\label{p:T1-T3=>T4} A finite set $R$ endowed with a ternary operation $T:R^3\to R$ satisfying the axioms {\sf(T1), (T2), \sf(T3)} is a ternar.
\end{proposition}

\begin{proof} Given any distinct points $\check x,\check y,\hat x,\hat y\in R$ with $\check x\ne \hat x$, consider the function 
$$F:R^2\to R^2,\quad F:(a,b)\mapsto(\check x_\times a_+b,\hat x_\times a_+b).$$To show that the function $F$ is injective, take any   pairs $ab,\alpha\beta\in R^2$ such that $\check x_\times a_+b=\check{x}_\times\alpha_+\beta$ and $\hat x_\times a_+b=\hat{x}_\times \alpha_+\beta$. Since $\check x\ne\hat x$, those equalities and the axiom {\sf(T3)} imply $a=\alpha$. Then the axiom {\sf(T2)} and the equality $\check x_\times a_+b=\check x_\times \alpha_+\beta=\check x_\times a_+\beta$ implies $b=\beta$ and hence $ab=\alpha\beta$, witnessing that the function $F:R^2\to R^2$ is injective. Since $R$ is finite, the injective function $F$ is bijective. Then for the pair $\hat y\check y\in R^2$ there exists a unique pair $ab\in R^2$ such that $(\hat y,\check y)=F(\check x,\hat x)=(\check x_\times a_{+}b,\hat x_\times a_{+}b)$, confirming the axiom {\sf(T4)}.
\end{proof}

Given a ternar $R\ne\{0,1\}$, consider the affine plane whose set of point is $R^2$ and the family of lines is $$\mathcal L\defeq\{L_c:c\in R\}\cup\{L_{a,b}:a,b\in R\},$$
where $$L_c\defeq\{xy\in R^2:x=c\}\quad\mbox{and}\quad L_{a,b}\defeq\{xy\in R^2:y=x_\times a_{+}b\}\quad\mbox{for $a,b,c\in R$}.
$$The affine plane $R^2$ is endowed with the canonical affine base $(10,00,01)$ and hence is a based affine plane, called the \index{coordinate plane}\defterm{coordinate plane} of the ternar $R$.

The following theorem shows that the coordinate plane of a ternar is well-defined.

\begin{theorem}\label{t:coord-plane=>based-aff-plane} For every ternar $R\ne\{0,1\}$ its coordinate plane $R^2$ is a based affine plane.
\end{theorem}

\begin{proof} Since the ternar $R$ contains two distinct elements $0$ and $1$, every line $L\in\mathcal L\defeq\{L_c:c\in R\}\cup\{L_{a,b}:a,b\in R\}$ in its coordinate plane $R^2$ contains at least two points. Also the axiom {\sf(T4)} of a ternar implies that any distinct points $p,q\in R^2$ belong to a unique line $L\in\mathcal L$. Therefore, the family $\mathcal L$ satisfies the axioms {\sf(L1)}, {\sf(L2)} in Theorem~\ref{t:L1+L2} and uniquely determines the line relation of the liner $R^2$. Since $R^2\notin\mathcal L$, the liner $R^2$ has rank $\|R^2\|\ge 3$. 

To show that $\|R^2\|=3$, it suffices to show that any point $p=(x,y)\in R^2$ belongs to the flat hull of the $3$-element set $\{u,o,w\}$ consisting of the points $u\defeq(1,0)$, $o\defeq(0,0)$ and $w\defeq(0,1)$. If $x=0$, then $p\in \Aline ow\subseteq\overline{\{u,o,w\}}$. If $y=0$, then $p\in\Aline ou\subseteq\overline{\{u,o,w\}}$. So, assume that $x\ne 0\ne y$. Since $R\ne\{0,1\}$, there exists an element $x'\in R\setminus\{0,x\}$. By the axiom {\sf(T3)} of a ternar, there exist unique elements $a,b\in R$ such that $x_\times a_+b=y$ and $x'_\times a_+b=0$. Observe that the points $u'=(x',0)\in \Aline ou$, $w'=(o,b)\in\Aline ow$, and $p=(x,y)$ belong to the liner $L_{a,b}$ and hence $p\in \Aline {u'}{w'}\subseteq\overline{\{u,o,w\}}$, witnessing that the liner $R^2$ has rank $\|R^2\|=3$.

To show that the liner $R^2$ is an affine plane, it suffices to check that for every line $L\in\mathcal L$ and point $p=p'p''\in R^2\setminus L$, there exists a unique line $\Lambda\in\mathcal L$ such that $p\in\Lambda\subseteq R^2\setminus L$. If $L=L_c\defeq\{xy\in R^2:x=c\}$ for some element $c\in R$, then $\Lambda\defeq L_{p'}\defeq\{xy\in R^2:x=p'\}$  is a unique line $R^2$ with $p\in\Lambda\subseteq R^2\setminus L$.

If $L=L_{a,b}$ for some elements $a,b\in R$, then by the axiom {\sf(T2)}, there exists a unique element $\beta\in R$ with $p'_\times a_+\beta=p''$. Taking into account that $p\in L_{a,\beta}\setminus L_{a,b}$, we conclude that $\beta\ne b$ and hence $L_{a,\beta}\cap L_{a,b}=\varnothing$, by the axiom {\sf(T2)}. The axiom {\sf(T3)} implies that $\Lambda\defeq L_{a,\beta}\defeq\{xy\in R^2:y=x_\times a_+\beta\}$ is a unique line with $p\in\Lambda\subseteq R^2\setminus L$. Therefore, the liner $R^2$ is a Playfair plane. By Theorem~\ref{t:Playfair<=>}, the liner $R^2$ is an affine plane.

Since the points $10,00,01$ in $R^2$ are non-collinear, the triple $(10,00,01)$ is an affine base in $R^2$ and hence $R^2$ is a based affine plane.
\end{proof} 

\begin{proposition} For every based affine plane $(\Pi,uow)$ its coordinate plane $\Delta^2$ coincides with the coordinate plane of the ternar $(\Delta,T_{uow})$.
\end{proposition}

\begin{proof} Let $\Delta$ be the diagonal of the affine base $uow$ and let  $C:\Pi\to\Delta^2$ be the coordinate chart of the based affine plane $(\Pi,uow)$. By definition of the liner $\Delta^2$, the family of lines on the coordinate plane $\Delta^2$ is equal to the family $\mathcal L'\defeq\{C[L]:L\in\mathcal L\}$ where $\mathcal L$ is the family of lines in the affine plane $\Pi$. The definition of the ternary operation $T_{uow}:\Delta^3\to\Delta$ implies that the family $\mathcal L'$ coincides with the family of lines $\{\{c\}\times\Delta\}\cup\{\{xy\in\Delta^2:y=x_\times a_+b\}:a,b\in\Delta\}$ of the coordinate plane $\Delta^2$ of the ternar $(\Delta,T_{uow})$.  
\end{proof}

\section{Isomorphisms and isotopisms of ternars}

\begin{definition} Two ternars $(R,T)$ and $(R',T')$ are \index{ternars!isomorphic}\index{isomorphic ternars}\defterm{isomorphic} if there exists a bijective function $F:R\to R'$ such that  $$F(T(x,a,b))=T'(F(x),F(a),F(b))$$ for every elements $x,a,b\in R$. The function $F$ is called an \index{isomorphism of ternars}\defterm{isomorphism} of the ternars $(R,T)$ and $(R',T')$.
\end{definition}

\begin{proposition}\label{p:R=Delta} Let $R$ be a ternar and $R^2$ be its coordinate plane. The ternar $\Delta$ of the based affine plane $R^2$ is isomorphic to the ternar $R$ via the isomorphism $F:R\to\Delta$, $F:x\mapsto xx$.
\end{proposition}

\begin{proof} Consider the bijective map $F:R\to\Delta$, $F:x\mapsto xx$,  of $R$ onto the diagonal $\Delta\defeq\{xy\in R^2:x=y\}$ of the coordinate plane $R^2$ of the ternar $(R,T)$. The definition of the coordinate plane $R^2$ ensures that its family of lines $\mathcal L$ consists of the lines $$L_{a,b}\defeq\{xy\in R^2:y=T(x,a,b)\}\quad\mbox{and}\quad L_c\defeq\{xy\in R^2:x=c\}\quad\mbox{for \ $a,b,c\in R$}.$$ Let $(u,o,w)\defeq(10,00,01)$ be the affine base of the based affine plane $R^2$. Let $$C:R^2\to\Delta^2,\quad C:xy\mapsto (F(x),F(y)),$$ be the coordinate chart of the based affine plane $R^2$. The definition of the ternary operation $T_{uow}:\Delta^3\to \Delta$ ensures that $$C[L_{a,b}]=C[\{xy\in R^2:y=T(x,a,b)\}]=\{xy\in\Delta^2:y=T_{uow}(x,aa,bb)\}$$ for every $a,b\in R$. Then for every $x\in R$ we have $$(F(x),F(T(x,a,b))=C(x,T(x,a,b))=(xx,T_{uow}(xx,aa,bb))$$ and hence
 $$F(T(x,a,b))=T_{uow}(xx,aa,bb)=T_{uow}(F(x),F(a),F(b)),$$ witnessing that $F:R\to\Delta$ is an isomorphism of the ternars $(R,T)$ and $(\Delta,T_{uow})$.
\end{proof}

\begin{definition} Two based affine planes $(\hat\Pi,\hat u\hat o\hat w)$ and $(\check \Pi,\check u\check o\check w)$ are \index{isomorphic based affine planes}\index{affine planes!isomorphic}\defterm{isomorphic} if there exists a liner isomorphism $I:\hat \Pi\to \check \Pi$ such that $I\hat u\hat o\hat w=\check u\check o\check w$.
\end{definition}

\begin{theorem}\label{t:tring-iso<=>} Two ternars $\hat R,\check R$ are isomorphic if and only if their coordinate planes $\hat R^2,\check R^2$ are isomorphic as based affine planes.
\end{theorem}

\begin{proof} Assuming that the ternars $(\hat R,\hat T)$ and $(\check R,\check T)$ are isomorphic, fix an isomorphism $F:\hat R\to \check R$. This isomorphism induces a bijective map $\ddot F:\hat R^2\to \check R^2$, $\ddot F:(x,y)\mapsto(F(x),F(y))$, between the coordinate planes $\hat R^2$ and $\check R^2$ of the ternars $\hat R$ and $\check R$. We claim that $\ddot F$ is an isomorphism of the based affine planes $\hat R^2$ and $\check R^2$.

Given any line $L\subseteq\hat R^2$, we need to show that $\ddot F[L]$ is a line in the liner $\check R^2$. By the definition of coordinate plane $\hat R^2$, there exist elements $\hat a,\hat b,\hat c\in \hat R$ such that $L=L_{\hat c}\defeq\{xy\in\hat R:x=\hat c\}$ or $L=L_{\hat a,\hat b}\defeq\{xy\in\hat R^2:y=\hat T(x,\hat a,\hat b)\}$. Consider the elements $\check a\defeq F(\hat a)$, $\check b\defeq F(\check b)$, and $\check c\defeq F(\hat c)$ in the ternar $\check R$. If $L=L_{\hat c}$, then $\ddot F[L_{\hat c}]=L_{\check c}$ is a line in the liner $R^2$. Next, assume that $L=L_{\hat a,\hat b}$ and observe that for every element $\hat x\in\hat R$ and its image $\check x\defeq F(\hat x)\in \check R$, we have the equality
$$F\big(\hat T(\hat x,\hat a,\hat b)\big)=\check T(F(\hat x),F(\hat a),F(\hat b))=\check T(\check x,\check a,\check b),$$
which implies the equality 
$$
\begin{aligned}
\ddot F[L]&=\ddot F[L_{\hat a,\hat b}]=\{(F(\hat x),F(\hat y)):\hat x\hat y\in\hat R^2\;\wedge\; \hat y=\hat T(\hat x,\hat a,\hat b)\}\\
&=\{(F(\hat x),F(\hat y)):\hat x\hat y\in\hat R^2\;\wedge\;F(\hat y)=F(\hat T(\hat x,\hat a,\hat b))\}=\{(\check x,\check y)\in\check R^2:\check y=\check T(\check x,\check a,\check b)\}=L_{\check a,\check b},
\end{aligned}
$$ 
witnessing that $\ddot F[L]=\ddot F[L_{\hat a,\hat b}]=L_{\check a,\check b}$ is a line in the coordinate plane $\check R^2$. By analogy, we can prove that for every line $\Lambda$ in the liner $\check R^2$, the preimage $\ddot F^{-1}[\Lambda]$ is a line in the liner $\hat R^2$.
Therefore, $\ddot F:\hat R^2\to\check R^2$ is an isomorphism of the liners $\hat R^2$ and $\check R^2$.

Let $\hat 0,\hat 1$ be the zero and unit of the ternar $\hat R$, and $\check 0,\check 1$ be the   zero and unit of the ternar $\check R$. 
Proposition~\ref{p:TR-01-unique} implies $\check 0=F(\hat 0)$ and $\check 1=F(\hat 0)$. Then the isomorphism $\ddot F$ maps the affine base $(\hat 1\hat 0,\hat 0\hat 0,\hat 0\hat 1)$ of the coordinate plane $\hat R^2$ onto the affine base $(\check 1\check 0,\check 0\check 0,\check 0\check 1)$ of the coordinate plane $\check R^2$, and hence $\ddot F$ is an isomorphism of the based affine planes $\hat R^2$ and $\check R^2$.
\smallskip

Now assume conversely, that there exists an isomorphism $I:\hat R^2\to\check R^2$ of the based affine planes $\hat R^2$ and $\check R^2$. Denote by $$\hat\Delta\defeq\{\hat x\hat y\in\hat R^2:\hat x=\hat y\}\quad\mbox{and}\quad \check\Delta\defeq\{\check x\check y\in\check R^2:\check x=\check y\}$$the diagonals of the coordinate planes $\hat R^2$ and $\check R^2$, respectively.

We claim that $I[\hat \Delta]=\check \Delta$. Taking into account that  $I$ maps the affine base $(\hat 1\hat 0,\hat 0\hat 0,\hat 0\hat 1)$ of the coordinate plane $\hat R^2$ onto the affine base $(\check 1\check 0,\check 0\check 0,\check 0\check 1)$ of $\check R^2$, we conclude that
$$
I[L_{\hat 0,\hat 0}]=I[\Aline{\hat 0\hat 0}{\hat 1\hat 0}]=\Aline{I(\hat 0\hat 0)}{I(\hat 1\hat 0)}=\Aline{\check 0\check 0}{\check 1\check 0}=L_{\check 0,\check 0}\;\;\mbox{and}\;\;I[L_{\hat 0}]=I[\Aline{\hat 0\hat 0}{\hat 0\hat 1}]=\Aline{I(\hat 0\hat 0)}{I(\hat 0\hat 1)}=\Aline{\check 0\check 0}{\check 0\check 1}=L_{\check 0}.$$
It follows from $L_{\hat 0,\hat 1}\parallel L_{\hat 0,\hat 0}$ and $L_{\hat 0}\parallel L_{\hat 1}$ that $I[L_{\hat 0,\hat 1}]\parallel I[L_{\hat 0,\hat 0}]=L_{\check 0,\check 0}$ and $I[L_{\hat 1}]\parallel I[L_{\hat 0}]=L_{\check 0}$, which implies $I[L_{\hat 0,\hat 1}]=L_{\check 0,\check 1}$, $I[L_{\hat 1}]=L_{\check 1}$, and finally,  
$$\{I(\hat 1\hat 1)\}=I[L_{\hat 0,\hat 1}\cap L_{\hat 1}]=I[L_{\hat 0,\hat 1}]\cap I[L_{\hat 1}]=L_{\check 0,\check 1}\cap L_{\check 1}=\{\check 1\check 1\}.$$
Therefore, $I[\hat \Delta]=I[\Aline{\hat 0\hat 0}{\hat 1\hat 1}]=\Aline{I(\hat 0\hat 0)}{I(\hat 1\hat 1)}=\Aline{\check 0\check 0}{\check 1\check 1}=\check\Delta.$  Consider the bijections 
 $$\hat J:\hat R\to\hat\Delta,\quad\hat J:\hat x\mapsto\hat x\hat x,\quad\mbox{and}\quad\check J:\check R\to\check \Delta,\quad\check J:\check x\mapsto\check x\check x.$$  
 We claim that the bijection $F\defeq \check J^{-1}I\hat J:\hat R\to\check R$ is an isomorphism of the ternars $\hat R$ and $\check R$.
Given any elements $\hat x,\hat a,\hat b\in \hat R$ and the element $\hat y\defeq \hat T(\hat x,\hat a,\hat b)$, we have to show that $F(\hat y)=\check y\defeq \check T(F(\hat x),F(\hat a),F(\hat b))$. Consider the points $\check x\defeq F(\hat x)$, $\check a\defeq F(\hat a)$, $\check b\defeq F(\check b)$. Taking into account that $I=\check JF\hat J^{-1}$, we conclude that  $I(\hat x\hat x)=\check x\check x$, $I(\hat a\hat a)=\check a\check a$ and $I(\hat b\hat b)=\check b\check b$.

It follows from $\check a\check a=F(\hat a\hat a)\in F[L_{\hat 0,\hat a}]\parallel L_{\check 0,\check 0}$ that $F[L_{\hat 0,\hat a}]=L_{\check 0,\check a}$. By analogy we can prove that $F[L_{\hat 0,\hat b}]=L_{\check 0,\check b}$. Consider the point $\hat 1\hat a\in L_{\hat 1}\cap L_{\hat 0,\hat a}$ and its image $I(\hat 1\hat a)\in I[L_{\hat 1}]\cap I[L_{\hat 0,\hat a}]=L_{\check 1}\cap L_{\check 0,\check a}=\{\check 1\check a\}$. Then  $I[L_{\hat a,\hat 0}]=I[\Aline{\hat 0\hat 0}{\hat 1\hat a}]=\Aline{\check 0\check 0}{\check 1\check a}=L_{\check a,\check 0}$.

Observe that $\{F(\hat 0\hat b)\}=F[L_{\hat 0}\cap L_{\hat 0,\hat b}]=F[L_{\hat 0}]\cap F[L_{\hat 0,\hat b}]=L_{\check 0}\cap L_{\check 0,\check b}=\{\check 0\check b\}$.  It follows from $L_{\hat a,\hat b}\parallel L_{\hat a,\hat 0}$ that $I[L_{\hat a,\hat b}]\parallel I[L_{\hat a,\hat 0}]=L_{\check a,\check 0}$. Since $\check 0\check b=I(\hat 0\hat b)\in I[L_{\hat a,\hat b}]\cap L_{\check a,\check b}$, the parallelity relation $I[L_{\hat a,\hat b}]\parallel L_{\check a,\check b}$ implies $F[L_{\hat a,\hat b}]=L_{\check a,\check b}$. Then 
$$\{I(\hat x\hat y)\}=I[L_{\hat x}\cap L_{\hat a,\hat b}]=I[L_{\hat x}]\cap I[L_{\hat a,\hat b}]=L_{\check x}\cap L_{\check a,\check b}=\{\check x\check y\},$$ which implies $I[L_{\hat 0,\hat y}]=L_{\check 0,\check y}$,  
$$\{I(\hat y\hat y)\}=I[L_{\hat 0,\hat y}\cap\hat \Delta]=I[L_{\hat 0,\hat y}]\cap\check\Delta=L_{\check 0,\check y}\cap\check \Delta=\{\check y\check y\}$$and finally,
$F(\hat y)=\check J^{1}I\hat J(\hat y)=\check J^{-1}I(\hat y\hat y)=\check J^{-1}(\check y\check y)=\check y$, 
witnessing that $F:\hat R\to \check R$ is an isomorphism of the ternars $\hat R$ and $\check R$. 
\end{proof}

Since every based affine plane $(\Pi,uow)$ is isomorphic to its coordinate plane $(\Delta^2,(eo,oo,oe))$, Theorem~\ref{t:tring-iso<=>} implies the following characterization.

\begin{corollary}\label{c:basedafplane-iso<=>} Two based affine planes are isomorphic if and only if their ternars are isomorphic.
\end{corollary}

A bit weaker notion than an isomorphism is an isotopism of ternars.
 
\begin{definition}\label{d:ternars-isotopic} Two ternars $(R,T)$ and $(R',T')$ are defined to be  \index{ternars!isotopic}\index{isotopic ternars}\defterm{isotopic} if there exist bijective functions $F,G,H:R\to R'$ such that  $$H(T(x,a,b))=T'(F(x),G(a),H(b))$$ for all elements $x,a,b\in R$. The triple of bijections $(F,G,H)$ is called an \index{isotopism of ternars}\defterm{isotopism} of the ternars $(R,T)$ and $(R',T')$.
\end{definition}

Observe that a bijection $F:R\to R'$ is an isomorphism of ternars  if and only if the triple $(F,F,F)$ is an isotopism of the ternars. 

\begin{remark} If $(F,G,H)$ is an isotopism of ternars $R,R'$, then the functions $F,G$ are uniquely determined by the function $H$ and the values $F(1)$ or $G(1)$. Indeed,
$$F(1)_\times G(x)_+H(0)=H(1_\times x_+0)=H(x)=H(x_\times 1_+0)=F(x)_\times G(1)_+H(0).$$ So, the computation complexity of finding an isotopism between two ternars of cardinality $n$ is at most $O(n!{\cdot}n{\cdot}n^3)$.
\end{remark}

Like isomorphisms, isotopisms of ternars (partly) preserve the constant $0$.

\begin{proposition}\label{p:isotopic=>zero} If $(F,G,H)$ is an isotopism of ternars $(\check R,\check T)$ and $(\hat R,\hat T)$, then $F(0)=0=G(0)$.
\end{proposition}

\begin{proof} By the axioms {\sf(T1)} and {\sf(T3)}, the equation $T'(y,G(0),H(0))=T'(y,G(1),H(0))$ has a unique solution $y=0$. Since $$T'(F(0),G(0),H(0))=H(T(0,0,0))=H(0)=H(T(0,1,0))=T'(F(0),G(1),H(0)),$$
the element $F(0)$ equals zero in the ternar $(\hat R,\hat T)$. 

By the axioms {\sf(T1)} and {\sf(T3)}, for every $a\in \hat R\setminus\{0\}$, the equation $T'(x,a,0)=H(0)$ has a unique solution $x\in\hat R$. On the other hand
$$T'(F(0),G(0),H(0))=H(T(0,0,0))=H(0)=H(T(1,0,0))=T'(F(1),G(0),H(0)),$$ which implies $G(0)=0$.
\end{proof}

\begin{proposition}\label{p:isotopic=>isomorphic} If $(F,G,H)$ is an isotopism of ternars $(\check R,\check T)$ and $(\hat R,\hat T)$, then the bijection $\Psi:\check R^2\to \hat R^2,\quad \Psi:(x,y)\mapsto (F(x),H(y))$, 
is an isomorphism of the coordinate planes of the ternars $(\check R,\check T)$ and $(\hat R,\hat T)$.
\end{proposition}

\begin{proof} We need to check that for every line $L$ in the coordinate plane $\check R^2$ its image $\Psi[L]$ is a line in the coordinate plane $\check R^2$. If the line $L$ is vertical, then $L=\{c\}\times \check R$ for some $c\in \check R$ and its image $\Psi[L]=\{F(c)\}\times \hat R$ is a vertical line in the plane $\hat R^2$. If $L$ is not vertical, then $L=\{(x,T(x,a,b)):x\in\check R\}$ for some elements $a,b\in\check R$. Consider the elements $a'\defeq G(a)$ and $b'\defeq H(b)$ of the ternar $\hat R$. Since $(F,G,H)$ is an isotopism of the ternars $(\check R,\check T)$ and $(\hat R,\hat T)$, for every $x\in\check R$ we have the equality $H(T(x,a,b))=T'(F(x),a',b')$, which implies that $\Psi[L]=\{(x,T'(x,a',b')):x\in \check R\}$ is a line in the plane $\hat R^2$. Therefore, the bijection $\Psi$ is an isomorphism of the coordinate planes $\check R^2$ and $\hat R^2$, by Theorem~\ref{t:liner-isomorphism<=>}.
\end{proof}

Proposition~\ref{p:isotopic=>isomorphic} can be reversed as follows.

\begin{proposition}\label{p:isomorphic=>isotopic} Let $(\check R,\check T)$ and $(\hat R,\hat T)$ be two ternars and $F,H:\check R\to\hat R$ be two bijective functions such that $F(0)=0$ and the bijection  $\Psi:\check R^2\to\hat R^2$, $\Psi:(x,y)\mapsto(F(x),H(y))$, is an isomorphism of the coordinate planes of the ternars. Then there exists a unique bijection $G:\check R\to\hat R$ such that the triple $(F,G,H)$ is an isotopism of the ternars $(\check R,\check T)$ and $(\hat R,\hat T)$.
\end{proposition}

\begin{proof} Given any element $a\in\check R$, consider the line $\check L_{a,0}\defeq \{(x,T(x,a,0)):x\in\check R\}$ in the coordinate plane $\check R^2$ of the ternar $(\check R,\check T)$. Since the function $\Psi$ is an isomorphism of the coordinate planes $\check R^2$ and $\hat R^2$, there exist unique elements $a'$ and $b'$ such that $\Psi[\check L_{a,0}]=\hat L_{a',b'}\defeq\{(x,T'(x,a',b')):x\in \check R\}$.
 Consider the function $G:\check R\to\hat R$ assigning to every element $a\in \check R$ the unique element $a'$ such that $\Psi[\check L_{a,0}]=\hat L_{a',b'}$ for some $b'\in\hat R$. To see that the function $G$ is injective, observe that for two distinct elements $a,c\in\check R$, the lines $\check L_{a,0}$ and $\check L_{c,0}$ are concurrent and so are the lines $\Psi[\check L_{a,0}]$ and $\Psi[\check L_{c,0}]$, which implies $a'\ne c'$. To see that $G$ is surjective, take any element $\hat a\in \hat R$, consider the line  $\hat L_{\hat a,0}$ in the coordinate plane $\hat R^2$, and find unique elements $a,b\in\check R$ such that $\Psi^{-1}[\hat L_{\hat a,0}]=\check L_{a,b}$. Next, find unique elements $a',b'\in\hat R$ such that $\Psi[\check L_{a,0}]=\hat L_{a',b'}$. We claim that  $a'=\hat a$. Indeed, the parallelity of the lines $\check L_{a,0}$ and $\check L_{a,b}$ implies the parallelity of the lines $\Psi[\check L_{a,0}]=\hat L_{a',b'}$ and $\Psi[\check L_{a,b}]=\hat L_{\hat a,0}$, which implies $\hat a=a'$. Therefore, $G$ is a bijection.
 
It remains to show that the triple $(F,G,H)$ is an isotopism of the ternars $(\check R,\check T)$ and $(\hat R,\hat T)$. Given any elements $ a,b\in \check R$, we should check that $H(T(x,a,b))=T'(F(x),G(a),H(b))$ for every $x\in\check R$. Consider the line $\check L_{a,b}\defeq \{(x,T(x,a,b)):x\in\check R\}$ in the coordinate plane $\check R^2$ of the ternar $(\check R,\check T)$. Since $\Psi$ is an isomorphism of the coordinate planes $\check R^2$ and $\hat R^2$, there exist unique elements $a',b',a'',b''\in \hat R$ such that $\Psi[\check L_{a,b}]=\hat L_{a',b'}$ and $\Psi[\check L_{a,0}]=\hat L_{a'',b''}$. Since the lines $\check L_{a,0}$ and $\check L_{a,b}$ are parallel, their images $\Psi[\check L_{a,0}]=\hat L_{a',b'}$ and $\Psi[\check L_{a,b}]=\hat L_{a'',b''}$ are parallel, which implies $a'=a''\defeq G(a)$. 

For every $x\in \check R$, consider the vertical lines $\check L_x\defeq\{x\}\times \check R$ and $\check L_{F(x)}\defeq\{F(x)\}\times\hat R$ in the coordinate planes $\check R^2$ and $\hat R^2$, respectively. Observe that 
 $$
 \begin{aligned}
 (F(x),H(T(x,a,b))&=\Psi(x,T(x,a,b))\in \Psi[\check L_x\cap \check L_{a,b}]=\Psi[\check L_x]\cap\Psi[\check L_{a,b}]=\hat L_{F(x)}\cap\hat L_{a'',b''}\\
 &=\{(F(x),T'(F(x),a'',b''))\}=\{(F(x),T'(F(x),G(a),b'')\}
 \end{aligned}
 $$ and hence $H(T(x,a,b))=T'(F(x),G(a),b'')$ 
for every $x\in \check R$. In particular, for the element $x=0\in\check R$, we have
$H(b)=H(T(0,a,b))=T'(F(0),G(a),b''))=T'(0,G(a),b'')=b''$, and finally,
$$H(T(x,a,b))=T'(F(x),G(a),b'')=T'(F(x),G(a),H(b))$$for all $x\in\check R$,
 witnessing that $(F,G,H)$ is an isotopism of the ternars $(\check R,\check T)$ and $(\hat R,\hat T)$.
\end{proof}

\begin{proposition}\label{p:ternars-isotopic} The ternars $\Delta$ and $\Delta'$ of two affine bases $uow$ and $u'o'w'$ in an affine plane $\Pi$ are isotopic provided $\Aline ou=\Aline{o'}{u'}$ and $\Aline ow=\Aline{o'}{w'}$.
\end{proposition}

\begin{proof} Assume that $\Aline ou=\Aline{o'}{u'}$ and $\Aline ow=\Aline{o'}{w'}$, and consider the horizontal and vertical directions $\boldsymbol h\defeq(\Aline ou)_\parallel=(\Aline o{u'})_\parallel$ and $\boldsymbol v\defeq (\Aline ow)_\parallel=(\Aline o{w'})_\parallel$ of the affine bases $uow$ and $u'o'w'$. Observe that $o'\in\Aline {o'}{u'}\cap\Aline{o'}{w'}=\Aline ou\cap\Aline ow=\{o\}$. 

Consider the line projections $F\defeq\boldsymbol v_{\Delta',\Delta}:\Delta\to\Delta'$ and $H\defeq\boldsymbol h_{\Delta',\Delta}:\Delta\to\Delta'$. Also, consider the bijective function $G:\Delta\to\Delta'$ assigning to every $x\in\Delta$ the unique point $y\in\Delta'$ such that $$\|\{o\}\cup(\Aline x{\boldsymbol h}\cap\Aline u{\boldsymbol v})\cup(\Aline y{\boldsymbol h}\cap\Aline {u'}{\boldsymbol v})\|=2.$$We claim that the triple $(F,G,H)$ is an isotopism of the ternars $\Delta$ and $\Delta'$. Given any points $x,a,b\in\Delta$, we need to check that $H(x_\times a_+b)=F(x)_\times G(a)_+H(b)$.
Consider the points $x'\defeq F(x)$, $a'\defeq G(a)$, and $b'\defeq H(b)$. The definition of the functions $F,H$ ensures that $\Aline x{\boldsymbol v}=\Aline{x'}{\boldsymbol v}$ and  $\Aline {b}{\boldsymbol h}=\Aline{b'}{\boldsymbol h}$. 

The definition of the function $G$ ensures that the flat hull $L$ of the set $\{o\}\cup(\Aline a{\boldsymbol h}\cap\Aline u{\boldsymbol v})\cup(\Aline {a'}{\boldsymbol h}\cap\Aline {u'}{\boldsymbol v})$ is a line in the plane $\Pi$. Let $\Lambda$ be a unique line in $\Lambda$ such that $L'\parallel L$ and $\Aline ow\cap\Aline b{\boldsymbol h}=\Aline{o'}{w'}\cap\Aline{b'}{\boldsymbol h}\subseteq \Lambda$. Let $y$ be the unique point of the intersection $\Aline {x}{\boldsymbol v}\cap L'=\Aline {x}{\boldsymbol v}\cap L'$. The definition of the ternary operation on the ternars $\Delta$ and $\Delta'$ guarantee that $x_\times a_+b\in \Aline y{\boldsymbol h}\cap\Delta$ and $F(x)_{\times}G(a)_+H(b)\in  \Aline y{\boldsymbol h}\cap\Delta'=\{\boldsymbol h_{\Delta',\Delta}(x_\times a_+b)\}=\{H(x_\times a_+ b)\}$. Therefore, 
$H(x_\times a_+b)=F(x)_\times G(a)_+H(b)$ and the triple $(F,G,H)$ is an isotopy of the ternars $\Delta$ and $\Delta'$.
\end{proof}

\begin{remark}\label{r:trings-1693} By computer calculations, Ivan Hetman found that there exist exactly 1693 non-isomorphic and 33 non-isotopic ternars of order 9.
The numbers of isomorphy and isotopy classes of ternars whose coordinate planes are isomorphic to one of seven possible affine planes of order 9 are presented in the following table.
$$
\begin{array}{|l|c|c|c|c|c|c|c|c|}
\hline
\mbox{The affine plane:}&\mbox{\tt Desarg}&\mbox{\tt Thales}&\mbox{\tt Hall}&\mbox{\tt hall}&\mbox{\tt dhall}&\mbox{\tt Hughes}&\mbox{\tt hughes}&\mbox{All}\\
\hline
\mbox{Ternars:}&1&4&25&144&165&214&1140&1693\\
\mbox{up to isotopy:}&1&2&4&4&6&7&9&33\\
\hline
\end{array}
$$
\end{remark}

\begin{definition} We say that a ternar $R$ is a \defterm{ternar of an affine plane} $\Pi$ if $\Pi$ has an affine base $uow$ whose ternar $(\Delta,T_{uow})$ is isomorphic to the ternar $R$. More generally, a ternar $R$ is defined to be \defterm{a ternar of an affine space} $X$ if $R$ is a ternar of some affine plane in $X$.
\end{definition}

\begin{examples} The computer calculations presented in Remark~\ref{r:trings-1693} show that the (unique) Thalesian affine plane $\mathbb J_9\times\mathbb J_9$ of order 9 (described in Section~\ref{s:J9xJ9}) has 4 non-isomorphic ternars, whereas two non-isomorphic affine planes {\tt Hughes} and {\tt Hughes} in the Hughes projective plane have 214 and 1140 nonisomorphic ternars, respectively.
\end{examples}

\section{The plus, puls and dot operations in ternars}

Every ternar $(R,T)$ carries three binary operations 
$$
\begin{aligned}
\plus\!&:R\times R\to R,\quad \plus\!:(x,b)\mapsto x\plus b\defeq T(x,1,b)=x_\times 1_+b,\\
\puls\!&:R\times R\to R,\quad \puls\!:(x,b)\mapsto x\puls b\defeq T(1,x,b)=1_\times x_+b,\\
\cdot&:R\times R\to R,\quad\;\, \cdot:(x,a)\mapsto x\cdot a\defeq T(x,a,0)=x_{\times}a_+0,\\
\end{aligned}
$$
called the \index{ternar!plus operation}\defterm{plus operation} the \index{ternar!puls operation}\defterm{puls operation}, and the \index{ternar!dot operation}\defterm{dot operation} on $R$, respectively.

\begin{definition} A ternar $R$ is called \defterm{plus$=$puls} if the binary operations $\plus$ and $\puls$ coincide.
\end{definition}

Important examples of plus=puls ternars are linear ternars, in which the ternary operation is uniquely determined by the addition and multiplication operations.

\begin{definition}\label{d:linear-ternar} A ternar $R$ is \index{linear ternar}\index{ternar!linear}\defterm{linear} if $$x_\times a_+b=(x\,{\cdot}\, a)\plus b=(x\,{\cdot}\,a)\puls b$$ for every elements $x,a,b\in R$.
\end{definition}

The following characterization of linear ternars show that for defining linear ternars we can use any of two addition operations.

\begin{proposition}\label{p:ternar-linear<=>} For a ternar $R$, the following conditions are equivalent:
\begin{enumerate}
\item $R$ is linear;
\item $\forall x,a,b\in R\;\; x_\times a_+b=(x\,{\cdot}\, a)\plus b\defeq (x_\times a_+0)_\times 1_+b$;
\item $\forall x,a,b\in R\;\; x_\times a_+b=(x\,{\cdot}\, a)\puls b\defeq 1_\times(x_\times a_+0)_+b$;
\item $R$ is plus\textup{=}puls and  $\forall x,a,b\in R\;\; x_\times a_+b=(x\,{\cdot}\,a)+b$.
\end{enumerate}
\end{proposition}

\begin{proof} By Definition~\ref{d:linear-ternar}, $(1)\Leftrightarrow(4)\Ra(2\wedge 3)$.
\smallskip

$(2)\Ra(4)$ If (2) holds, then for all $x,b\in R$,
$$x\plus b=x_\times 1_+b=(x\,{\cdot}\,1)\plus b=(1\,{\cdot}\, x)\plus b=1_\times x_+b=x\puls b,$$
witnessing that the ternar $T$ is plus=puls. Then for all $x,a,b\in R$ we have
$$x_\times a_+b=(x\,{\cdot}\,a)\plus b=(x\,{\cdot}\,a)\puls b,$$
which means that the ternar $R$ is linear.
\smallskip

The implication $(3)\Ra(4)$ can be proved by analogy.
\end{proof}

\begin{example} Every corps $R$ endowed with the ternary operation $$T:R^3\to R,\quad T:(x,a,b)\mapsto (x\cdot a)+b,$$ is a linear ternar.
\end{example} 

\begin{remark}[Ivan Hetman, 2024]\label{rem:plus=puls11} Computer calculations show that among 1693 non-isomorphic ternars of seven affine planes of order 9 there are 11 ternars which are plus=puls but not linear.
\end{remark}

Let us recall that any set $M$ endowed with a binary operation $\cdot:M\times M\to M$ is called a \index{magma}\defterm{magma}. A magma $M$ is
\begin{itemize}
\item \index{unital magma}\index{magma!unital}\defterm{unital} if it contains a (necessarily unique) element $e\in M$ (called the \defterm{identity} or else the \index{magma!identity of}\index{magma!neutral element of}\defterm{neutral element} of $M$) such that $e\cdot x=x=x\cdot  e$ for all $x\in M$;
\item a \index{loop}\defterm{loop} if $M$ is a unital magma such that  for every $a,b\in M$, the equations $x\cdot a=b$ and $a\cdot y=b$ have unique solutions $x,y\in M$;
\item \index{$0$-magma}\defterm{$0$-magma} if $M$ contains a (necessarily unique) element $0\in M$ (called the \defterm{zero} of $M$) such that $x\cdot 0=0=0\cdot x$ for every $x\in M$;
\item a \index{$0$-loop}\defterm{$0$-loop} if $M$ is a unital $0$-magma such that for every $a\in M\setminus\{0\}$ and $b\in M$  the equations $x\cdot a=b$ and $a\cdot y=b$ have unique solutions $x,y\in M$.
\end{itemize}

\begin{exercise} Show that for every $0$-loop $X$ containing more than one element, the set $X^*\defeq X\setminus\{0\}$ is a loop, with respect to the binary operation inherited from $X$.
\end{exercise}

\begin{theorem}\label{t:tring-triloop} If $R$ is a ternar, then $(R,\plus)$ and $(R,\!\puls\!)$ are loops with neutral element $0$ and $(R,\cdot)$ is a $0$-loop such that the following conditions are satisfied:
\begin{enumerate}
\item $\forall x\in R\;\;(1\plus x=1\puls  x)$;
\item for every $a\in R\setminus\{1\}$ and $b\in R$, the equation $x\cdot a=x\plus b$ has a unique solution $x\in R$.
\end{enumerate}
\end{theorem}

\begin{proof} The axiom {\sf(T1)} ensures that for every $x\in R$ we have the equalities
$$0+x=0_\times 1_+x=x=x_\times 1_+0=x+0$$
and
$$0\puls x=1_\times 0_+x=x=1_\times x_+0=x\puls 0,$$
witnessing that $0$ is the neutral element of the magmas $(R,+)$ and $(R,\puls)$.  

The definitions of the plus and puls operations ensure that
$$1+x=1_\times 1_+x=1\puls x$$
for every $x\in X$.

By the axiom {\sf(T2)}, for every elements $\alpha,\beta\in R$, the equation $\beta=\alpha+y=\alpha_\times 1_+y$ has a unique solution $y\in R$. By the axioms {\sf(T3)} and {\sf(T1)}, the equation
$$x+\alpha=x_\times 1_+\alpha=x_\times 0_+\beta=\beta$$has a unique solution $x\in R$, witnessing that the unital magma $(R,+)$ is a quasigroup and hence a loop.

By the axiom {\sf(T2)}, for every elements $\alpha,\beta\in R$, the equation $\beta=\alpha\puls y=1_\times\alpha_+y$ has a unique solution $y\in R$. By the axioms {\sf(T4)} and {\sf(T1)}, the equation
$$x\puls \alpha=1_\times x_+\alpha=0_\times x_+\beta=\beta$$has a unique solution $x\in R$, witnessing that the unital magma $(R,\!\puls\!)$ is a quasigroup and hence a loop.

The axiom {\sf(T1)} ensures that 
$$\forall x\in R\;\;(x_\times 1_+0=x=1_\times x_+0)\;\wedge\;(x\cdot 0=x_\times 0_+0=0=0_\times x_+0=0\cdot x),$$
witnessing that $(R,\cdot)$ is a unital $0$-magma with the identity element $1$ and zero element $0$. To check that $(R,\cdot)$ is a $0$-loop, take any elements $\alpha\in R^*\defeq R\setminus\{0\}$ and $\beta\in R$. By the axioms {\sf(T3)}, the equation $x_\times\alpha_+0=x_\times 0_+\beta$ has a unique solution $x\in R$, which is also a unique solution of the equalition $x\cdot\alpha=\beta$, according to the axiom {\sf(T1)}. To solve the equation $\alpha\cdot y=\beta$, apply the axiom {\sf(T4)} and find unique elements $y,b\in R$ such that $\alpha_\times y_+b=\beta$ and $0_\times y_+b=0$.    The axiom {\sf(T1)} ensures that $b=0_\times y_+b=0$ and hence $\alpha\cdot y=\alpha_\times y_+0=\beta$. Therefore, $y$ is a unique solution of the equation $\alpha\cdot y=\beta$, witnessing that $(R,\cdot)$ is a $0$-loop.

Finally, for every $a\in R\setminus\{1\}$ and $b\in R\setminus\{0\}$, the equation
$$x\cdot a=x_\times a_+0=x_\times 1_+b=x+b$$has a unique solution, by the axiom {\sf(T3)}.
\end{proof}

\begin{lemma}\label{l:Boolean} Let $(\Pi,uow)$ be a based affine plane and $\Delta$ be its ternar. For every point $x\in\Delta$, the following conditions are equivalent:
\begin{enumerate}
\item $x+x=o$;
\item $\Aline{ox}{xo}\parallel \Aline {oo}{xx}$.
\end{enumerate}
\end{lemma}

\begin{proof} If $x=o$, then $x+x=o+o=o$ and $\Aline {xo}{ox}=\Aline {oo}{xx}$, so both conditions (1) and (2) hold and hence are equivalent. So, assume that $x\ne o$. For the point $y\defeq x+x$, the definition of plus operation ensures that $\Aline {ox}{xy}\subparallel\Delta=\Aline {oo}{xx}$ and hence $\Aline{ox}{xy}\parallel \Aline{oo}{xx}$, by Corollary~\ref{c:subparallel}.

$(1)\Ra(2)$ If $x+x=o$, then $y=x+x=o$ and hence $\Aline{ox}{xo}=\Aline{ox}{xy}\parallel\Aline{oo}{xx}$.

$(2)\Ra(1)$ If $\Aline{xo}{ox}\parallel \Aline{oo}{xx}$, then the parallelity relations
$\Aline{xo}{ox}\parallel \Aline{oo}{xx}$ and  $\Aline{ox}{xy}\parallel \Aline{oo}{xx}$ imply $\Aline{ox}{xo}=\Aline{ox}{xy}$ and hence $\{xo\}=(\{x\}\times\Delta)\cap\Aline{ox}{xo}=(\{x\}\times \Delta)\cap\Aline{ox}{xy}=\{xy\}$ and $x+x=y=o$.
\end{proof}

\begin{definition} Let $R$ be a ternar and $R^*\defeq R\setminus \{0\}$. The loops $(R,+)$, $(R,\puls)$, and $(R^*,\cdot)$ are called the \index{plus  loop}\index{ternar!plus loop of}\defterm{plus loop},  \index{puls loop}\index{ternar!puls loop of}\defterm{puls loop}, and the \index{dot loop}\index{ternar!dot loop of}\defterm{dot loop} of the ternar $R$, respectively.
\end{definition}

\begin{remark} By computer calculations, Ivan Hetman has established that among plus, puls and dot loops of all 1693 non-isomorphic ternars of cardinality 9 there exist exactly 1372 non-isomorphic plus loops, 1373 non-isomorphic puls loops, and 1424 non-isomorphic dot loops. More information on the number of non-isomorphic ternars and also plus, puls and dot loops of all seven affine planes of order 9 can be found in the following table.
$$
\begin{array}{|l|c|c|c|c|c|c|c|c|}
\hline
\mbox{The affine plane:}&\mbox{\tt Desarg}&\mbox{\tt Thales}&\mbox{\tt Hall}&\mbox{\tt hall}&\mbox{\tt dhall}&\mbox{\tt Hughes}&\mbox{\tt Hughes}&\mbox{All}\\
\hline
\mbox{Ternars:}&1&4&25&144&165&214&1140&1693\\
\mbox{Plus loops:}&1&1&10&108&145&142&1042&1372\\
\mbox{Puls loops:}&1&1&22&124&118&189&996&1373\\
\mbox{Dot loops:}&1&4&21&124&127&204&1054&1424\\
\hline
\end{array}
$$
\end{remark}

\section{Triloops}

In the preceding section we have seen that every ternar $R$ carries three binary operations $+$, $\puls$, ${\cdot}$  such that $(R,+)$ and $(R,\!\puls\!)$ are loops and $(R,\cdot)$ is a $0$-loop. This observation motivates the following definition of an algebraic structure called a triloop.

\begin{definition} A \index{triloop}\defterm{triloop} is a set $R$ endowed with three binary operations $$
\begin{aligned}
&+:R\times R\to R,\quad +:(x,y)\mapsto x+y,\\
&\!\puls\!:R\times R\to R,\quad \puls\!:(x,y)\mapsto x\puls y,\\
&{}\;\cdot\,:R\times R\to R,\quad \;\,\cdot\,:(x,y)\mapsto x\cdot y,
\end{aligned}
$$ and two distinct constants $0,1\in X$ such that the following axioms are satisfied:
\begin{itemize}
\item[{\sf(TL0)}] $\forall x\in R\;(x+0=x=0+x)$;
\item[{\sf(TL1)}] $\forall a,b\in R\;\exists ! x\in R\;\exists ! y\in R\;\;(x+a=b=a+y)$;
\smallskip
\item[{\sf(TL2)}] $\forall x\in R\;(x\puls 0=x=0\puls x)$;
\item[{\sf(TL3)}] $\forall a,b\in R\;\exists ! x\in R\;\exists ! y\in R\;\;(x\puls a=b=a\puls y)$;
\smallskip
\item[{\sf(TL4)}] $\forall x\in R\;\;(x\cdot 0=0=0\cdot x)\;\wedge\; (x\cdot 1=x=1\cdot x)$;
\item[{\sf(TL5)}] $\forall a\in R\setminus\{0\}\;\forall b\in R\;\exists! x\in X\;\exists !y\in R\;(x\cdot a=b=a\cdot y)$.
\end{itemize}
\end{definition} 

\begin{remark} A set $X$ endowed with three binary operations $+,\!\puls\!,\,{\cdot}\,$, and two distinct constants $0,1\in X$ is a triloop if and only if $(X,+)$ and $(X,\puls)$ are loops and $(X,\cdot)$ is a $0$-loop whose zero element $0$ is the neutral element of the loops $(X,+)$ and $(X,\!\puls\!)$, and whose identity is $1$.
\end{remark}

By Theorem~\ref{t:tring-triloop}, for a ternar $R$, the set $R$ endowed with the plus, puls and dot operations is a triloop, called  \index{ternar!triloop of}\defterm{the triloop of the ternar} $R$.

\begin{definition}\label{d:triloop-is:} A triloop $(R,+,\!\puls\!,\cdot)$ is defined to be
\begin{itemize}
\item \index{plus=puls triloop}\index{triloop!plus=puls}\defterm{plus\textup{=}puls} if $x+y=x\puls y$ for all $x,y\in R$;
\item \index{$1$-compatible triloop}\index{triloop!$1$-compatible}\defterm{plus-$1$-puls} if $1\,{+}\,y=1\puls y$ for every $x\in R$;
\item \index{planar triloop}\index{triloop!planar}\defterm{planar} if $R$ is the triloop of some ternar;
\item \index{triloop!concurrent}\index{concurrent triloop}\defterm{concurrent} if for all $a\in R\setminus\{1\}$ and $b\in R$ there exists a unique element $x\in R$ such that $x\cdot a=x+b$.
\end{itemize}
\end{definition}

It is clear that every plus=puls triloop is plus-$1$-puls. 
Theorem~\ref{t:tring-triloop} implies

\begin{proposition} Every planar triloop is puls-$1$-plus and concurrent.
\end{proposition}

\begin{remark} In Theorem~\ref{t:Hughes-Wilker}, we shall prove that an infinite triloop is planar if and only if it is plus-$1$-plus and concurrent.
\end{remark}

\begin{example}\label{ex:concurrent-non-planar3} The plus\textup{=}puls triloop $R\defeq\{0,1,a,b,c\}$ whose plus, puls and dot operations are defined by the tables
$$\begin{array}{c|ccccc}
{+}{=}\puls&0&1&a&b&c\\
\hline
0&0&1&a&b&c\\
1&1&0&c&a&b\\
a&a&b&0&c&1\\
b&b&c&1&0&a\\
c&c&a&b&1&0
\end{array}
\quad\mbox{and}\quad
\begin{array}{c|ccccc}
\cdot&0&1&a&b&c\\
\hline
0&0&0&0&0&0\\
1&0&1&a&b&c\\
a&0&a&1&c&b\\
b&0&b&c&1&a\\
c&0&c&b&a&1
\end{array}
$$
is plus-$1$-puls and concurrent but not planar (because planar triloops of order 5 are fields).
\end{example}

\begin{Exercise} Find an example of a planar triloop which is not plus=puls. Show that 9 is the smallest cardinality of a planar triloop, which is not plus=puls.
\smallskip

{\em Hint:} Look at Remark~\ref{rem:plus=puls11}.
\end{Exercise}           

\begin{remark}[Ivan Hetman]\label{r:1687biloops} Computer calculations show that two ternars of cardinality $\le 10$ are isomorphic if and only if they have isomorphic triloops. 
\end{remark}

\begin{Exercise} Find examples of two non-isomorphic ternars whose triloops are isomorphic.
\end{Exercise}

\section{Biloops}

In plus=puls triloops, the plus and puls operations coincide so the structure of a plus=puls triloop can be reduced to a simpler structure of a biloop.

\begin{definition} A \index{biloop}\defterm{biloop} is a set $R$ endowed with two binary operations $$
+:R\times R\to R\quad +:(x,y)\mapsto x+y,\quad\mbox{and}\quad
\cdot:R\times R\to R,\quad \cdot:(x,y)\mapsto x\cdot y,
$$ and two distinct constants $0,1\in X$ such that the following axioms are satisfied:
\begin{itemize}
\item[{\sf(BL0)}] $\forall x\in R\;(x+0=x=0+x)$;
\item[{\sf(BL1)}] $\forall a,b\in R\;\exists ! x\in R\;\exists ! y\in R\;\;(x+a=b=a+y)$;
\smallskip
\item[{\sf(BL2)}] $\forall x\in R\;\;(x\cdot 0=0=0\cdot x)\;\wedge\; (x\cdot 1=x=1\cdot x)$;
\item[{\sf(BL3)}] $\forall a\in R\setminus\{0\}\;\forall b\in R\;\exists! x\in X\;\exists !y\in R\;(x\cdot a=b=a\cdot y)$.
\end{itemize}
\end{definition} 

\begin{remark} A set $X$ endowed with two binary operations $+,\cdot:X\times X\to X$ and two distinct constants $0,1\in X$ is a biloop if and only if $(X,+)$ is a loop and $(X,\cdot)$ is a $0$-loop whose zero element $0$ coincides with the identity element of the loop $(X,+)$.
\end{remark}

\begin{remark} By the axiom {\sf(BL3)}, for every elements $a\in R\setminus\{0\}$ and $b\in R$ of a biloop $R$, there exist unique elements $a\backslash b$ and $b/a$ in $R$ such that $a\cdot(a\backslash b)=b=(b/a)\cdot a$. We shall often use the operations $a\backslash b$ and $b/a$ in the sequel without additional explanations.
\end{remark}

\begin{remark} By the axiom {\sf(BL1)}, for every elements $a,b\in R$ of a biloop $R$, there exist unique elements $b\mbox{-}a$ and $\mbox{-}a_+b$ in $R$ such that $(b\mbox{-}a)+a=b=a+(\mbox{-}a_+b)$. We shall often use the operations $b\mbox{-}a$ and $\mbox{-}a_+b$ in the sequel without additional  explanations.
\end{remark}

\begin{definition} For a triloop $(R,+,\!\puls\!,\cdot)$, the biloop $(R,+,\cdot)$ is called the  \defterm{biloop of the triloop} $(R,+,\!\puls\!,\cdot)$. Conversely, for every biloop $(R,+,\cdot)$, the plus=puls triloop $(R,+,+,\cdot)$ is called the \defterm{plus\textup{=}puls triloop} of the biloop $(R,+,\cdot)$.
\end{definition}

By Theorem~\ref{t:tring-triloop}, for a ternar $R$, the set $R$ endowed with the plus and dot operations is a biloop, called  \index{ternar!biloop of}\defterm{the biloop of the ternar} $R$.


\begin{definition}\label{d:biloop-is:} A biloop $R$ is defined to be
\begin{itemize}
\item \index{planar biloop}\index{biloop!planar}\defterm{planar} if $R$ is the biloop of some ternar;
\item \index{linear biloop}\index{biloop!linear}\defterm{linear} if the set $R$ endowed with the ternary operation
$$T:X^3\to X,\quad T:(x,a,b)\mapsto (x\cdot a)+b,$$is a (necessarily linear) ternar;
\item \index{biloop!concurrent}\index{concurrent biloop}\defterm{concurrent} if for every $a\in R\setminus\{1\}$ and $b\in R$ there exists a unique element $x\in R$ such that $x\cdot a=x+b$.
\end{itemize}
\end{definition}

Definitions~\ref{d:triloop-is:} and \ref{d:biloop-is:} imply the following propositions.

\begin{proposition} If a triloop $(R,+,\!\puls,\cdot)$ is planar, then so is the biloop $(R,+,\cdot)$.
\end{proposition}

\begin{proposition} A triloop $(R,+,\!\puls,\cdot)$ is concurrent if and only if the biloop $(R,+,\cdot)$ is concurrent.
\end{proposition} 

\begin{example}\label{ex:plus1puls7} The $7$-element triloop $(7,+,\!\puls\!,\cdot)$ whose binary operations are defined by the tables
$$\begin{array}{c|ccccccc}
+&0&1&2&3&4&5&6\\
\hline
0&0&1&2&3&4&5&6\\
1&1&2&3&4&5&6&0\\
2&2&3&4&5&6&0&1\\
3&3&4&5&6&0&1&2\\
4&4&5&6&0&1&2&3\\
5&5&6&0&1&2&3&4\\
6&6&0&1&2&3&4&5\\
\end{array}
\qquad
\begin{array}{c|ccccccc}
\puls&0&1&2&3&4&5&6\\
\hline
0&0&1&2&3&4&5&6\\
1&1&2&3&4&5&6&0\\
2&2&0&5&1&6&4&3\\
3&3&5&4&6&1&0&2\\
4&4&3&6&0&2&1&5\\
5&5&6&1&2&0&3&4\\
6&6&4&0&5&3&2&1\\
\end{array}
\qquad
\begin{array}{c|ccccccc}
\cdot&0&1&2&3&4&5&6\\
\hline
0&0&0&0&0&0&0&0\\
1&0&1&2&3&4&5&6\\
2&0&2&4&6&1&3&5\\
3&0&3&6&2&5&1&4\\
4&0&4&1&5&2&6&3\\
5&0&5&3&1&6&4&2\\
6&0&6&5&4&3&2&1\\
\end{array}
$$
is plus-$1$-puls and concurrent but it is not plus=puls and hence not planar (because the unique planar triloop of order 7 is plus=plus). Moreover, analyzing the list\footnote{https://users.cecs.anu.edu.au/$\sim$bdm/data/latin.html} of all 23746 loops of order 7, it can be shown that there exist exactly 129 concurrent puls-$1$-puls triloops of order 7 (one of which is planar). 
\end{example}

It is clear that every linear biloop is planar. By Theorem~\ref{t:tring-triloop}, every planar biloop is concurrent. Therefore, for every biloop $R$, the following implications hold:
$$\mbox{linear}\Ra\mbox{planar}\Ra\mbox{concurrent}.$$
 
\begin{example}\label{ex:concurrent-non-planar} The biloop $R\defeq\{0,1,a,b,c\}$ whose plus and dot operations are defined in Example~\ref{ex:concurrent-non-planar3} 
is concurrent but not planar (because planar biloops of order 5 are fields).
\end{example}

\begin{Exercise} Find an example of a non-linear planar biloop of cardinality 9. Show that 9 is the smallest cardinality of a non-linear planar biloop.
\smallskip

{\em Hint:} Look at Remark~\ref{r:1687biloops}.
\end{Exercise}           





\begin{proposition}\label{p:linear-biloop<=>} A biloop $R$ is linear if and only if it has the following two properties:
\begin{enumerate}
\item for every $a,b,c,d\in R$ with $a\ne c$, the equation $(x\cdot a)+b=(x\cdot c)+d$ has a unique solution $x\in R$;
\item for every elements $\check x,\hat x,\check y,\hat y\in R$ with $\check x\ne \hat x$, there exits unique elements $a,b\in R$ satisfying the equations $(\check x\cdot a)+b=\check y$ and $(\hat x\cdot a)+b=\hat y$.
\end{enumerate}
\end{proposition} 

\begin{proof} To prove the ``only if'' part, assume that the biloop $(R,+,\cdot)$ is linear. Then the ternary operation $T:R^3\to R$, $T:(x,a,b)\mapsto (x\cdot a)+b$, turns $R$ into a ternar. The axioms {\sf(T3)} and {\sf(T4)} of a ternar imply the conditions (1) and (2), respectively. 
\smallskip

To prove the ``if'' part, assume that a biloop $R$ satisfies the conditions (1) and (2). Consider the ternary operation $T:R^3\to R$, $T:(x,a,b)\mapsto x_\times a_+b\defeq (x\cdot a)+b$. The conditions (1) and (3) imply that the ternary operation $T$ satisfies the axioms {\sf(T3)} and {\sf(T4)} of a ternar. Since $(R,+,\cdot)$ is a biloop, there exist elements $0,1\in R$ such that $x+0=x=0+x$, $x\cdot 0=0=0\cdot x$, and $x\cdot 1=x=1\cdot x$ for every $x\in R$.  Then for every $x,b\in R$ we have the equalities
$$x_\times 0_+b=(x\cdot 0)+b=0+b=b=0+b=(0\cdot x)+b=0_\times x_+b$$ and
$$x_\times 1_+0=(x\cdot 1)+0=x+0=x=x+0=(1\cdot x)+0=1_\times x_+0,$$
witnessing that the axiom {\sf(T1)} holds.

To see that the axiom {\sf(T2)} hold, take any elements $x,a,y\in R$. Since $(R,+)$ is a loop, the equation $x_\times a_+b=(x\cdot a)+b=y$ has a unique solution $b\in R$.

Therefore, $(R,T)$ is a ternar and hence the biloop $(R,+,\cdot)$ is linear.
\end{proof}


\begin{definition}\label{d:biloops-isotopic}
Two biloops $R$ and $R'$ are defined to be \index{isotopic biloops}\index{biloops!isotopic}\defterm{isotopic} if there exist bijective functions $F,G,H:R\to R'$ such that $H((x{\cdot}y)+z)=(F(x){\cdot}G(y))+H(z)$ for all $x,y,z\in R$. The triple $(F,G,H)$ is called an \defterm{isotopism} of the biloops $R,R'$. If $F=G=H$, then the bijection $F=G=H$ is an \defterm{isomorphism} of the biloops $R,R'$.
\end{definition}

Definitions~\ref{d:ternars-isotopic} and \ref{d:biloops-isotopic} imply the following simple characterization.

\begin{proposition}\label{p:ternar-isotopic<=>biloops} Two linear ternars are isotopic if and only if their biloops are isotopic.
\end{proposition}

\begin{remark} It is clear that any isomorphic biloops are isotopic. The converse is not always true. On the oter hand, by a difficult Theorem~\ref{t:Schafer} of Schafer, isotopic alternative division rings are isomorphic.
\end{remark}

\begin{remark}\label{r:1687biloops} By computer calculations, Ivan Hetman has established that among the biloops of all 1693 non-isomorphic ternars of cardinality 9 there exist exactly 1687 non-isomorphic biloops. So, there exist exactly 1687 of non-isomorphic planar biloops of cardinality 9. The exact numbers of non-isomorphic biloops for all seven affine planes of order 9 are presented in the following table.
$$
\begin{array}{|l|c|c|c|c|c|c|c|c|}
\hline
\mbox{The affine plane:}&\mbox{\tt Desarg}&\mbox{\tt Thales}&\mbox{\tt Hall}&\mbox{\tt hall}&\mbox{\tt dhall}&\mbox{\tt Hughes}&\mbox{\tt hughes}&\mbox{All}\\
\hline
\mbox{Ternars:}&1&4&25&144&165&214&1140&1693\\
\mbox{Biloops:}&1&4&22&144&165&213&1140&1687\\
\hline
\end{array}
$$
The computer calculations show that non-isomorphic affine planes can have isomorphic biloops. More precisely, three non-isomorphic affine planes of order 9 (namely, {\tt Hall}, {\tt dhall}, and {\tt Hughes}) endowed with suitable affine bases have isomorphic biloops (which are anti-isomorphic\footnote{Two biloops $X,Y$ are \defterm{anti-isomorphic} if there exists a bijective function $F:X\to Y$ such that $F(x+y)=F(y)+F(x)$ and $F(x\cdot y)=F(y)\cdot F(x)$ for every $x,y\in X$.} to the near-field $\mathbb J_9$, considered in Section~\ref{s:J9xJ9}). This example also shows that there exist three non-isomorphic  ternars of cardinality 9 whose biloops are isomorphic (only one of these three non-isomorphic ternars is linear).
\end{remark}

\section{A characterization of infinite planar triloops}

In this section we consider the problem of recognizing planar triloops. For planar biloops this problem was first posed by \index[person]{Evans}Professor Trevor Evans\footnote{{\bf Trevor Evans} (1925 -- 1991) received bachelor's and master's degrees, and the Doctor of Science from Oxford University. He taught at Manchester and Wisconsin before joining the Emory faculty in 1951. Here, he served as departmental chair for the periods 1957-72 and 1976-82; he was named Fuller E. Callaway Professor of Mathematics in 1980. In the early 50's he was invited to spend a year at the Institute for Advanced Study in Princeton; he travelled extensively throughout his career as a research mathematician. An expert in algebra and its applications to combinatorial problems, Professor Evans is most well-known for his work on partial planes, nonassociative algebras, and word and decision problems. He had many PhD students during his career at Emory: Curt Lindner is a notable research mathematician, a Distinguished University Professor at Auburn, and the late Etta Falconer, who served as dean at Spelman University and who received the Hay Award from the Association of Women in Mathematics, to name two of his most successful students. Professor Evans was awarded the Emory Williams Distinguished Teaching Award for Graduate Education in 1972. 
In recognition of his lifelong service to mathematics, Professor Evans received the first distinguished service award from the Southeastern Section of the Mathematical Association of America. As well, prizes in his name are given by the MAA to recognize the best articles published in ``Math Horizons''.} to \index[person]{Hughes}Daniel Hughes\footnote{{\bf Daniel R.~Hughes} (1927--2012), was a British and American mathematician, working in Finite Geometries, Finite Groups, and Design Theory. Hughes defended his Ph.D. Thesis ``Planar Division Neo-Rings'' in 1955 in University of Wisconsin-Madison, under supervision of R.H. Bruck.  Since 1964 till 1992 he worked in (Queen Mary and) Westfield College in University of London (UK). In Geometry Hughes is best known for Hughes planes, non-Desarguesian planes of order $p^{2n}$ of every odd $p$ and every $n$, constructed by Hughes in 1957. Playfair planes with a given countable translation group were constructed by Hughes in  [D.R.~Hughes, {\em Additive and multiplicative loops of planar ternary rings}, Proc. Amer. Math. Soc. {\bf 6} (1955), 973--980.]} who solved it in 1955 for countable biloops. In 1964, the result of Hughes was extended to arbitrary infinite  biloops by \index[person]{Wilker}Wilker\footnote{{\bf Peter Wilker} (1922 -- 2002), born in Vienna, Austria, was a Swiss mathematician and amateur magician. He studied math and became a professor of mathematics at the University of Bern. He therefore was quite interested in mathematical principles and in particular magic squares. He was the editor of the Swiss magazine Hokus Pokus from 1976 until 1985, and he contributed regularly to other periodicals such as Abracadabra and Magische Welt. From 1986 to 1988 he published the periodical Magische Bl\"atter. In 1988 he moved to Cornwall England where he became a member of the Cornish Magic Society. Peter Wilker collected his entire life jokes. His daughter Jessica Wilker published his best 100 jokes (in German) in the book ``Peters Witzb\"uchlein''.}. 

\begin{theorem}[Hughes--Wilker, 1955--1964]\label{t:Hughes--Wilker} An infinite biloop is planar if and only if it is concurrent.
\end{theorem}

The Hughes--Wilker Theorem~\ref{t:Hughes--Wilker} is a corollary of the following more general theorem characterizing infinite planar triloops.

\begin{theorem}\label{t:Hughes-Wilker} An infinite triloop is planar if and only if it is plus-$1$-puls and concurrent.
\end{theorem}

\begin{proof} The ``only if'' part follows immediately from Theorem~\ref{t:tring-triloop}. To prove the ``if'' part, we first introduce some notations and prove some auxiliary results.

Given two families of sets $\A$ and $\mathcal B$, we write
\begin{itemize}
\item $\A\parallel \mathcal B$ if $A\cap B=\varnothing$ for every distinct sets $A\in\A$ and $B\in\mathcal B$;
\item $\A\perp\mathcal B$ if $|A\cap B|=1$ for every distinct sets $A\in\A$ and $B\in\mathcal B$;
\item $\A\pm\mathcal B$ if $|A\cap B|\le1$ for every distinct sets $A\in\A$ and $B\in\mathcal B$.
\end{itemize}
We recall that $S_X$ denotes the set of all bijections of a set $X$. Elements of the set $S_X$ are bijective functions $F:X\to X$, which are subsets of the square $X\times X$.

\begin{lemma}\label{l:Wilker1} Let $X$ be a infinite set and $\mathcal A_0\subseteq\mathcal A_1\subseteq S_X$ be two families of bijections of $X$ such that $\A_0\parallel\A_0$, $\A_1\pm\A_1$,  $|\mathcal A_0|<|X|$, and $\sup_{p\in X^2\setminus\bigcup\A_0}|\{A\in\A_1:p\in A\}|<|X|$. Let $\varphi\subset X\times X$ be an injective function such that $|\varphi|<|X|$, $\varphi\cap \bigcup\A_0=\varnothing$, and $\{\varphi\}\pm\A_1$. Then there exists a bijective function $\Phi:X\to X$ such that $\varphi\subset\Phi$, $\Phi\cap\bigcup\A_0=\varnothing$ and $\{\Phi\}\perp(\A_1\setminus\A_0)$.
\end{lemma}

\begin{proof} Consider the infinite cardinal $\kappa\defeq|X|$. Since $\A_1\pm\A_1$, $\A_0\subseteq\A_1$ and $|\A_0|<\kappa$, for every function $A\in\A_1\setminus\A_0$ the set $A\setminus\bigcup\A_0$ is not empty. Then the condition $\sup_{p\in X^2\setminus\bigcup\A_0}|\{A\in\A_1:p\in A\}|<|X|=\kappa$ implies $|\A_1|\le \kappa$.

Partition the set $\kappa$ as the union $\kappa=\{0\}\cup\Omega\cup \Omega_1\cup\Omega_2$ of four disjoint sets $\{0\},\Omega,\Omega_1,\Omega_2$ of cardinalities $|\Omega|=|\A_1\setminus\A_0|$ and $|\Omega_1|=|\Omega_2|=\kappa$. 
Write the sets $\A_1\setminus\A_0$, $\kappa\setminus\dom[\varphi]$ and $\kappa\setminus\rng[\varphi]$ as $\A_1\setminus\A_0=\{A_\alpha\}_{\alpha\in \Omega}$, $\kappa\setminus\dom[\varphi]=\{x_\alpha\}_{\alpha\in\Omega_1}$ and $\kappa\setminus\rng[\varphi]=\{y_\alpha\}_{\alpha\in\Omega_2}$. By transfinite induction we shall define an increasing transfinite sequence of injective functions $(\Phi_\alpha)_{\alpha\in\kappa}$ such that $\Phi_0=\varphi$ and for every ordinal $\alpha\in\kappa$ the following conditions are satisfied:
\begin{itemize}
\item[$(1_\alpha)$] $\Phi_{<\alpha}\defeq\varphi\cup\bigcup_{\beta<\alpha}\Phi_\beta\subseteq\Phi_\alpha\subset X\times X$;
\item[$(2_\alpha)$] $|\Phi_\alpha\setminus\Phi_{<\alpha}|\le 1$;
\item[$(3_\alpha)$] $\Phi_\alpha\cap\bigcup\A_0=\varnothing$ and  $\{\Phi_\alpha\}\pm \A_1$;
\item[$(4_\alpha)$] if $\alpha\in\Omega$, then $|\Phi_\alpha\cap A_\alpha|=1$;
\item[$(5_\alpha)$] if $\alpha\in\Omega_1$, then $x_\alpha\in \dom[\Phi_\alpha]$;
\item[$(6_\alpha)$] if $\alpha\in\Omega_2$, then $y_\alpha\in\rng[\Phi_\alpha]$.
\end{itemize}
We start the inductive construction letting $\Phi_0=\varphi$. The properties of the function $\varphi=\Phi_0$ ensure that the properies $(1_0)$--$(6_0)$ are satisfied. Assume that for some nonzero ordinal $\alpha\in\kappa$ an increasing sequence $(\Phi_\beta)_{\beta<\alpha}$ of injective functions satisfing the inductive conditions $(1_\beta)$--$(6_\beta)$ for $\beta<\alpha$ have been constructed. Consider the set $\Phi_{<\alpha}\defeq\bigcup_{\beta<\alpha}\Phi_\beta\subseteq X\times X$ and observe that it is an injective function.  The inductive conditions $(1_\beta)$--$(3_\beta)$ for $\beta<\alpha$ imply that $\Phi_{<\alpha}\cap\bigcup\A_0=\varnothing$ and $\{\Phi_{<\alpha}\}\pm\A_1$. The inductive condition $(2_\beta)$ for $\beta<\alpha$ guarantees that $$|\Phi_{<\alpha}|=|\varphi\cup\bigcup_{\beta<\alpha}(\Phi_{\beta}\setminus\Phi_{<\beta})|\le|\varphi|+|\alpha|<\kappa.$$ Since $\max\{|\Phi_{<\alpha}|,|\A_0|\}<\kappa$, $\Phi_{<\alpha}\cap\bigcup\A_0=\varnothing$, and $\sup_{p\in X^2\setminus\bigcup\A_0}|\{A\in \A_1:p\in A\}|<\kappa$, the set $\A\defeq\A_0\cup\{A\in\A_1:A\cap \Phi_{<\alpha}\ne\varnothing\}$ has cardinality $|\A|<\kappa$.

Depending on the location of the ordinal $\alpha$ in the set $\Omega\cup\Omega_1\cup\Omega_2$, three cases are possible.
\smallskip

1. If $\alpha\in \Omega$ and $\Phi_{<\alpha}\cap A_\alpha\ne\varnothing$, then put $\Phi_\alpha\defeq\Phi_{<\alpha}$. 

If $\alpha\in\Omega$ and $\Phi_{<\alpha}\cap A_\alpha=\varnothing$, then choose any pair $p\in A_\alpha\setminus \bigcup\A$ and put $\Phi_\alpha\defeq\Phi_{<\alpha}\cup\{p\}$. The set $A_\alpha
\setminus\bigcup\A$ is not empty because $\A\pm\A$ and $A_\alpha\notin\A$.
\smallskip

2. If $\alpha\in \Omega_1$ and $x_\alpha\in\dom[\Phi_{<\alpha}]$, then put $\Phi_\alpha\defeq\Phi_{<\alpha}$. If $\alpha\in \Omega_1$ and $x_\alpha\notin\dom[\Phi_{<\alpha}]$, then choose any point $y_\alpha\in X\setminus (\rng[\Phi_{<\alpha}]\cup\{A(x_\alpha):A\in\A\})$ and put $\Phi_\alpha\defeq\Phi_{<\alpha}\cup\{(x_\alpha,y_\alpha)\})$.
\smallskip

3. If $\alpha\in \Omega_2$ and $y_\alpha\in\rng[\Phi_{<\alpha}]$, then put $\Phi_\alpha\defeq\Phi_{<\alpha}$. If $\alpha\in \Omega_2$ and $y_\alpha\notin\rng[\Phi_{<\alpha}]$, then choose any point $x_\alpha\in X\setminus (\dom[\Phi_{<\alpha}]\cup\{A^{-1}(y_\alpha):A\in\A\})$ and put $\Phi_\alpha\defeq\Phi_{<\alpha}\cup\{(x_\alpha,y_\alpha)\}$.
\smallskip

This completes the inductive step of the inductive construction.
\smallskip

After completing the inductive construction, consider the injective function $\Phi\defeq\bigcup_{\alpha\in\kappa}\Phi_{\alpha}$. The inductive conditions $(3_\alpha)$ and $(4_\alpha)$, $\alpha\in\kappa$, ensure that $\Phi\cap\bigcup\A_0=\varnothing$ and $\{\Phi\}\perp (\A_1\setminus\A_0)$. The inductive conditions $(5_\alpha)$ and $(6_\alpha)$, $\alpha\in\kappa$, imply $\dom[\Phi]=\kappa=\rng[\Phi]$, so that $\Phi:\kappa\to\kappa$ is a bijective function.
\end{proof} 

\begin{lemma}\label{l:Wilker2b}   Let $(X,\puls\!)$ be an infinite loop with neutral element $0$,  let $a,1\in  X\setminus\{0\}$ be two  elements and $\mathcal A_0\subseteq\mathcal A_1\subseteq S_X$ be two families of bijections of $X$ satisfying the following conditions:
\begin{enumerate}
\item $\A_0\parallel\A_0$, $\A_1\pm\A_1$, $\A_0\perp(\A_1\setminus\A_0)$, and $|\mathcal A_0|<|X|$; 
\item $(0,0)\in\bigcup\A_0$ and $\sup_{p\in X^2\setminus\bigcup\A_0}|\{A\in\A_1:p\in A\}|<|X|$;
\item $\forall F\in\A_0\;\;F(1)=a\puls F(0)$;
\item $\forall F\in \A_1\setminus\A_0\;\;F(1)\in \{b\puls F(0):b\in X\setminus\{a\}\}$.
\end{enumerate}
Then there exists a family of bijections $\F\subseteq S_X$ such that $\A_0\subseteq \F$, $\F\parallel\F$, $\F\perp(\A_1\setminus\A_0)$, $\bigcup\F=X\times X$, and $F(1)=a\puls F(0)$ for every $F\in \F$.
\end{lemma} 

\begin{proof} Let $\kappa\defeq|X|$ and write the set $X\times X$ as $X\times X=\{p_\alpha\}_{\alpha\in\kappa}$ where $p_0=(0,0)\in\bigcup\A_0$. Let $\F_0\defeq\A_0$ and observe that $\F_0\parallel \F_0$ and $\F_0\perp(\A_1\setminus\F_0)$. For every nonzero ordinal $\alpha\in\kappa$, we shall construct a  family $\F_\alpha\subseteq S_X$ satisfying the following conditions:
\begin{itemize}
\item[$(1_\alpha)$] $\F_{<\alpha}\defeq\bigcup_{\beta<\alpha}\F_\beta\subseteq\F_\alpha\subseteq S_X$,  $|\F_\alpha\setminus \F_{<\alpha}|\le 1$, and $|\F_{\alpha}|<\kappa$;
\item[$(2_\alpha)$] $\F_\alpha\parallel\F_\alpha$ and $p_\alpha\in\bigcup\F_\alpha$;
\item[$(3_\alpha)$] $F(1)=a\puls F(0)$ for every $F\in \F_\alpha$;
\item[$(4_\alpha)$] $\F_\alpha\perp(\A_1\setminus\A_0)$.
\end{itemize}

Assume that for some ordinal $\alpha\in\kappa$, we have constructed families $\F_\beta$, $\beta<\alpha$, satisfying the inductive conditions $(1_\beta)$--$(4_\beta)$, $\beta<\alpha$. Then the family $\F_{<\alpha}\defeq\bigcup_{\beta<\alpha}\F_\beta$ has $\F_{<\alpha}\parallel\F_{<\alpha}$ and $\F_{<\alpha}\perp(\A_1\setminus\A_0)$. The inductive conditions $(1_\beta)$, $\beta<\alpha$, imply that $|\F_{<\alpha}|\le|\F_0|+|\alpha|<\kappa$. If $p_\alpha\in\bigcup\F_{<\alpha}$, then put $\F_\alpha\defeq\F_{<\alpha}$. Next, assume that $p_\alpha\notin\bigcup\F_{<\alpha}$ and hence $p_\alpha\ne (0,0)$ (because $(0,0)\in\bigcup\A_0\subseteq \bigcup\F_{<\alpha}$). 

\begin{claim}\label{cl:varphi-alpha} There exists an injective function $\varphi_\alpha\subseteq X\times X$ satisfying the following conditions:
\begin{enumerate}
\item[(i)] $p_\alpha\in \varphi_\alpha$;
\item[(ii)] $\{0,1\}\subseteq \dom[\varphi_\alpha]$ and $|\dom[\varphi_\alpha]|\le 3$;
\item[(iii)] $\varphi_\alpha(1)=a\puls \varphi_\alpha(0)$;
\item[(iv)] $\varphi_\alpha\cap \bigcup\F_{<\alpha}=\varnothing$.
\item[(v)] $\varphi_\alpha\pm\A_1$;
\end{enumerate}
\end{claim}

\begin{proof} If $p_\alpha=(0,y)$ for some $y\in X$, then consider the function $\varphi_\alpha\defeq\{(0,y),(1,a\puls y)\}$. Assuming that $\varphi_\alpha\cap\bigcup\F_{<\alpha}\ne\varnothing$, we can find an ordinal $\beta<\alpha$ and a function $F\in\F_\beta$ such that $\varphi_\alpha\cap F\ne\varnothing$. Taking into account that $(0,y)=p_\alpha\notin\bigcup\F_{<\alpha}$, we conclude that $(1,a\puls y)\in F$. The inductive condition $(3_\beta)$ ensures that $a \puls y=F(1)=a\puls F(0)$ and hence $y=F(0)$ by the cancellativity of the loop $(X,\puls\!)$. Then $p_{\alpha}=(0,y)=(0,F(0))\in F$, which contradicts the assumption $p_\alpha\notin\bigcup\F_{<\alpha}$. This contradiction shows that $\varphi_\alpha$ satisfies the condition (iv).

 If $p_\alpha=(1,y)$ for some $y\in X$, then find a unique element $x\in X$ such that $y=a\puls x$ and put $\varphi_\alpha\defeq\{(0,x),(1,a\puls x)\}$. Assuming that $\varphi_\alpha\cap\bigcup\F_{<\alpha}\ne\varnothing$, we can find an ordinal $\beta<\alpha$ and a function $F\in\F_\beta$ such that $\varphi_\alpha\cap F\ne\varnothing$. Taking into account that $(1,y)=p_\alpha\notin\bigcup\F_{<\alpha}$, we conclude that $(0,x)\in F$ and hence $x=F(0)$. The inductive condition $(3_\beta)$ ensures that $y=a \puls  x=a\puls F(0)=F(1)$ and hence $p_{\alpha}=(1,y)=(1,F(1))\in F$, which contradicts the assumption $p_\alpha\notin\bigcup\F_{<\alpha}$. This contradiction shows that $\varphi_\alpha$ satisfies the condition (iv).

Finally, consider the case $p_\alpha=(p_\alpha',p_\alpha'')\notin \{0,1\}\times X$. By the assumption (2) of the lemma, the family $\A'_1\defeq\{A\in\A_1:p_\alpha\in A\}$ has cardinality $|\A'_1|<\kappa$. Then there exist points $x,y\in R$ such that $y=a\puls x$ and $$x,y\notin\{p_\alpha''\}\cup\{F(k):k\in\{0,1\},\;F\in\A'_1\cup\F_{<\alpha}\}.$$ Put $\varphi_\alpha\defeq\{p_\alpha\}\cup\{(0,x),(1,y)\}$. The choice of $x,y\notin \{F(k):k\in\{0,1\},\;F\in\F_{<\alpha}\}$ ensures that $\varphi_\alpha\cap\bigcup\F_{<\alpha}=\varnothing$.

Therefore, the function $\varphi_\alpha$ satisfies the conditions (i)--(iv).
 Assuming that the condition (v) is not satisfied, we can find a function $F\in\A_1$ such that $|\varphi_\alpha\cap F|\ge 2$. The condition (4) of the lemma ensures that $F(1)=b\puls F(0)$ for some $b\in X\setminus\{a\}$. Assuming that $\{0,1\}\subseteq\dom[\varphi_\alpha\cap F]$, we conclude that $b\puls F(0)= F(1)=\varphi_\alpha(1)=a\puls\varphi_\alpha(0)=a\puls  F(0)$ and hence $b=a$ (by the cancellativity of the loop $(X,\puls\!)$), which contradicts $b\in X\setminus\{a\}$. Therefore, $\{0,1\}\not\subseteq\dom[\varphi_\alpha\cap F]$ and thus $\varnothing\ne (\varphi_\alpha\cap F)\setminus(\{0,1\}\times R)\subseteq \{p_\alpha\}$. Then $p_\alpha\notin\{0,1\}\times X$ and we have the third case of the construction of the function $\varphi_\alpha$. It follows from $p_\alpha\in \varphi_\alpha\cap F$ that $F\in\A'_1$. Then the choice of the elements $x,y\notin\{F(0),F(1)\}$ guarantees that $\varphi_\alpha\cap F=\{p_\alpha\}$, which contradicts $|\varphi_\alpha\cap F|\ge 2$. 
 \smallskip
 
 This contradiction completes the proof of the conditions (i)--(v) for the function $\varphi_\alpha$.
\end{proof}

Let $\varphi_\alpha$ be an injective function satisfying the conditions of Claim~\ref{cl:varphi-alpha}. Applying Lemma~\ref{l:Wilker1}, find a bijection $\Phi_\alpha:X\to X$ such that $p_\alpha\in\varphi_\alpha\subseteq \Phi_\alpha$, $\Phi_\alpha\cap \bigcup\F_{<\alpha}=\varnothing$, and $\{\Phi\}\perp(\A_1\setminus\A_0)$. Then the family $\F_\alpha\defeq\F_{<\alpha}\cup\{\Phi_\alpha\}$ satisfies the conditions $(1_\alpha)$--$(3_\alpha)$. This completes the inductive step.
\smallskip

After completing the inductive construction, consider the family of bijections $\F\defeq\bigcup_{\alpha\in\kappa}\F_\alpha$ and observe that it has the required properties.
\end{proof}

Having Lemmas~\ref{l:Wilker1} and \ref{l:Wilker2b} in our disposition, we are able to prove Theorem~\ref{t:Hughes-Wilker}. Fix a plus-1-puls concurrent triloop $(R,+,\!\puls\!,\cdot)$ of infinite cardinality $\kappa\defeq|R|$.

Fix a well-order $\le$ on the set $R$ such that 
\begin{itemize}
\item for every $a\in R$ the set ${\downarrow}a\defeq\{x\in R:x\le a\}$ has cardinality $|{\downarrow}a|<\kappa$;
\item $0=\min_{\le} R$ and $1=\min_{\le}(R\setminus\{0\})$.
\end{itemize}
Here for a non-empty set $A\subseteq R$ we denote by $\min_{\le}A$ the smallest element of the set $A$ in the well-order $\le$.
Let $a_0\defeq 0$ and $a_1\defeq 1$.

Let $\mathcal B$ be the family of all injective functions $\varphi\subset R\times R$ such that $|\varphi|=2$ and $(0,0)\notin \varphi$. Let $\preceq$ be any well-order on the set $\mathcal B$ such that for every $\beta\in \mathcal B$ the set ${\downarrow}\beta\defeq\{\alpha\in \mathcal B:\alpha\preceq\beta\}$ has cardinality $|{\downarrow}\beta|<\kappa$.

For every $a\in R$, consider the set
$$L_a\defeq\{(x,x{\cdot}a):x\in R\}.$$
Observe that $$L_{a_0}=L_0=\{(x,y)\in R\times R:y=x{\cdot}0=0\}\quad\mbox{and}\quad L_{a_1}=L_1=\{(x,y)\in R\times R:y=x{\cdot}1=x\}.$$ Next, consider the families 
$$
\begin{aligned}
&\F_0\defeq\big\{\{(x,y)\in R\times R:y=b\}:b\in R\big\},\\
&\F_1\defeq\big\{\{(x,y)\in R\times R:y=x+b\}:b\in R\big\},\\
&\F_{<0}\defeq\{L_a:a\in R\},\quad\mbox{and}\quad \F_{<1}\defeq \F_{<0}\cup\F_0.
\end{aligned}
$$ Taking into account that $R$ is a concurrent triloop, we can see that $\F_{<0}\perp\F_{<0}$, $\F_i\parallel\F_i$ and $\F_i\perp(\F_{<i}\setminus\{L_{a_i}\})$ for every $i\in\{0,1\}$.

Write the cardinal $\kappa$ as the disjoint union $\kappa=S\cup P$ of two sets $S,P$ of cardinality $|S|=\kappa=|P|$. In the future inductive construction, the set $S$ will be responsible for spreads of lines and $P$ for pairs. 
 
By induction we shall construct transfinite sequences $(\F_\alpha)_{\alpha\in\kappa}$, $(a_\alpha)_{\alpha\in \kappa}$, and $(\varphi_\alpha)_{\alpha\in\kappa}$ such that for every $\alpha\in\kappa\setminus\{0,1\}$ the following conditions are satisfied:
\begin{itemize}
\item[$(1_\alpha)$] If $\alpha\in S$, then $a_\alpha=\min_{\le}(R\setminus\{a_\beta:\beta<\alpha\})$;
\item[$(2_\alpha)$] $L_{a_\alpha}\in\F_\alpha\subseteq S_X$, $\F_\alpha\parallel\F_\alpha$, $\bigcup\F_\alpha=X\times X$, and 
$\F_\alpha\perp (\F_{<\alpha}\setminus\{L_{a_\alpha}\})$,\\ where $\F_{<\alpha}\defeq \F_{<0}\cup\bigcup_{\beta<\alpha}\F_\beta$;
\item[$(3_\alpha)$] for every $F\in\F_\alpha$, $F(1)=a_\alpha\puls F(0)$;
\item[$(4_\alpha)$] if $\alpha\in P$, then $\varphi_\alpha=\min_\preceq(\mathcal B\setminus\bigcup_{\beta<\alpha}\bigcup_{F\in\F_\beta}\{\varphi\in \mathcal B:\varphi\subseteq F\})$\\ and $\varphi_\alpha\subseteq \Phi_\alpha$ for some $\Phi_\alpha\in \F_\alpha$.
\end{itemize}

Assume that for some ordinal $\alpha\in \kappa\setminus\{0,1\}$, a transfinite sequence $(\F_\beta)_{\beta<\alpha}$ satisfying the inductive conditions $(1_\beta)$--$(4_\beta)$, $\beta<\alpha$, has been constructed. The cancellativity of the loop $(R^*,\cdot)$, the definition of the family $\F_{<0}=\{L_a:a\in R\}$, and the inductive conditions $(2_\beta)$, $\beta<\alpha$, imply that the family $$\F_{<\alpha}\defeq\F_{<0}\cup\bigcup_{\beta<\alpha}\F_\beta$$ has property $\F_{<\alpha}\pm\F_{<\alpha}$. 

Consider the set $$\mathcal B_{<\alpha}\defeq \bigcup_{\beta<\alpha}\bigcup_{F\in\F_\beta}\{\varphi\in \mathcal B:\varphi\subseteq F\}.$$ 

\begin{claim}\label{cl:BneB} $\mathcal B_{<\alpha}\ne \mathcal B$. 
\end{claim}

\begin{proof} Take any point $a\in R\setminus\{a_\beta:\beta<\alpha\}$ and observe that $a\ne a_0=0$ and hence $\varphi\defeq \{(0,a),(1,a\puls a)\}$ is an element of the set $\mathcal B$. To derive a contradiction, assume that $\varphi\in \mathcal B_{<\alpha}$. Taking into account that  $\varphi(0)=a\ne 0$, we conclude that $\varphi\not\subseteq F$ for any function $F\in\F_{<0}$, and hence there exist $\beta<\alpha$ and $F\in\F_\beta$ such that $\varphi\subseteq F$. The inductive condition $(3_\beta)$ ensures that
$$a\puls a=\varphi(1)=F(1)=a_\beta\puls F(0)=a_\beta\cdot \varphi(0)=a_\beta\puls a$$ and hence $a=a_\beta$, which contradicts the choice of $a$. This contradiction shows that $\mathcal B_{<\alpha}\ne \mathcal B$.
\end{proof}

Claim~\ref{cl:BneB} ensures that the set $\mathcal B\setminus\mathcal B_{<\alpha}$ is not empty, so we can consider the injective function $\varphi_\alpha\defeq\min_{\preceq}(\mathcal B\setminus \mathcal B_{<\alpha})$. The inductive conditions $(2_\beta)$, $\beta<\alpha$, ensure that the family $\F'_{<\alpha}\defeq\{F\in\F_{<\alpha}:\varphi_\alpha\cap F\ne\varnothing\}$ has cardinality $|\F'_{<\alpha}|<\kappa$. 

\begin{claim} There exist $\bar a_\alpha\in R\setminus\{a_\beta:\beta<\alpha\}$ and  a function $\bar\varphi_\alpha\subseteq R\times R$ satisfying the following conditions:
\begin{enumerate}
\item $\varphi_\alpha\subseteq\bar\varphi_\alpha$ and $\bar\varphi_\alpha\setminus(\{0,1\}\times R)\subseteq\varphi_\alpha$;
\item $\{0,1\}\subseteq\dom[\bar\varphi_\alpha]$;
\item $\bar \varphi_\alpha(1)=\bar a_\alpha\puls \,\bar\varphi_\alpha(0)$;
\item $\bar \varphi_\alpha\cap L_{\bar a_\alpha}=\varnothing$;
\item $\{\bar\varphi_\alpha\}\pm\F_{<\alpha}$.
\end{enumerate}
\end{claim}

\begin{proof} If $\dom[\varphi_\alpha]=\{0,1\}$, then let $\bar a_\alpha\in R$ be a unique element such that $\varphi_\alpha(1)=\bar a_\alpha\puls \varphi_\alpha(0)$. In this case put $\bar \varphi_\alpha\defeq \varphi_\alpha$.
We claim that $\bar a_\alpha\notin\{a_\beta:\beta<\alpha\}$. In the oppsite case, we could find $\beta<\alpha$ such that $\bar a_\alpha=a_\beta$. By the inductive condition $(2_\beta)$, $\bigcup\F_\beta=X\times X$ and hence there exists a function $F\in\F_\beta$ such that $(0,\varphi_\alpha(0))\in F$. The inductive condition $(3_\beta)$ ensures that 
$\varphi_\alpha(1)=\bar a_\beta\puls \varphi_\alpha(0)=a_\beta\puls F(0)=F(1)$ and hence $\varphi_\alpha\subseteq F$, which contradicts $\varphi_\alpha\notin \mathcal B_{<\alpha}$. It is clear that the function $\bar \varphi_\alpha$ satisfies the conditions (1)--(3) of the claim. It follows from $\bar\varphi_\alpha(0)=\varphi_\alpha(0)\ne 0=L_{\bar a_\alpha}(0)$ and $\bar \varphi_\alpha(1)=\bar a_\alpha\puls\bar \varphi_\alpha(0)\ne \bar a_{\alpha}\puls 0=\bar a_\alpha= L_{\bar a_\alpha}(1)$ that $\bar \varphi_\alpha\cap L_{\bar a_\alpha}=\varnothing$, so the condition (4) is satisfied. It follows from $\bar\varphi_\alpha=\varphi_\alpha\notin \mathcal B_{<\alpha}$ that the condition (5) is satisfied, too.
\smallskip

If $\dom[\varphi_\alpha]\cap\{0,1\}=\varnothing$, then let $\bar a_\alpha\defeq
\min_{\le}(R\setminus\{a_\beta:\beta<\alpha\})$ and take any element $x_\alpha\in R\setminus\{0\}$ such that $\{(0,x_\alpha),(1,\bar a_\alpha\puls x_\alpha)\}\cap\bigcup\F'_{<\alpha}=\varnothing$. It is clear that the function $$\bar \varphi_\alpha\defeq\varphi_\alpha\cup \{(0,x_\alpha),(1,\bar a_\alpha\puls x_\alpha)\}$$ satisfies the conditions (1)--(3) of the claim. 
It follows from $\bar\varphi_\alpha(0)=x_\alpha\ne 0=L_{\bar a_\alpha}(0)$ and $\bar \varphi_\alpha(1)=\bar a_\alpha\puls x_\alpha\ne \bar a_{\alpha}\puls 0=\bar a_\alpha= L_{\bar a_\alpha}(1)$ that $\bar \varphi_\alpha\cap L_{\bar a_\alpha}=\varnothing$, so the condition (4) is satisfied.
We claim that the condition (5) holds, too. In the opposite case, we could find a function $F\in\F_{<\alpha}$ such that $|\bar\varphi_\alpha\cap F|\ge 2$. If $F\in\F'_\alpha$, then the choice of the point $x_\alpha$ ensures that $\bar\varphi_\alpha\cap F=(\bar\varphi\cap F)\setminus(\{0,1\}\times R)\subseteq \varphi_\alpha\cap F$ and hence $|\varphi_\alpha\cap F|=2$ and $\varphi_\alpha\subseteq F$, which contradicts the choice of $\varphi_\alpha\notin \mathcal B_{<\alpha}$.
This contradiction shows that $F\notin\F_\alpha'$ and hence $\varphi_\alpha\cap F=\varnothing$. Then $|\bar \varphi_\alpha\cap F|\ge 2$ implies $\{(0,x_\alpha),(1,\bar a_\alpha\puls x_\alpha)\}\subseteq F$. Since $x_\alpha\ne 0$, $F\notin\F_{<0}$ and hence $F=\F_\beta$ for some $\beta<\alpha$. 
The inductive condition $(3_\beta)$ ensures that 
$$\bar a_\alpha\puls F(0)=\bar a_\alpha\puls \bar\varphi_\alpha(0)=\bar\varphi_\alpha(1)=F(1)=a_\beta\puls F(0)$$
and hence $\bar a_\alpha=a_\beta$, which contradicts the choice of $\bar a_\alpha$.
\smallskip

If $0\in\dom[\varphi_\alpha]$ and $1\notin\dom[\varphi_\alpha]$, then we can choose an element $\bar a_\alpha\in R\setminus\{a_\beta:\beta<\alpha\}$ such that $$\bar a_\alpha\puls \varphi_\alpha(0)\notin \rng[\varphi_\alpha]\cup \{F(1):F\in\F'_{<\alpha}\}.$$ Put $\bar\varphi_\alpha\defeq\varphi_\alpha\cup\{(1,\bar a_\alpha\puls \varphi_\alpha(0)\}$ and observe that the function $\bar\varphi_\alpha$ satisfies the conditions (1)--(3) of the claim. 
It follows from $\bar\varphi_\alpha(0)=\varphi_\alpha(0)\ne 0=L_{\bar a_\alpha}(0)$ and $\bar \varphi_\alpha(1)=\bar a_\alpha\puls \varphi_\alpha(0)\ne \bar a_{\alpha}\!\puls 0=\bar a_\alpha= L_{\bar a_\alpha}(1)$ that $\bar \varphi_\alpha\cap L_{\bar a_\alpha}=\varnothing$, so the condition (4) is satisfied. Assuming that 
 the condition (5) does not hold, we could find a function $F\in\F_{<\alpha}$ such that $|\bar\varphi_\alpha\cap F|\ge 2$. If $F\in\F'_{<\alpha}$, then the choice of the element $\bar a_\alpha$ ensures that $(1,\bar \varphi_\alpha(1))=(1,\bar a_\alpha\puls \varphi_\alpha(0))\notin F$. Then $|\varphi_\alpha\cap F|=|\bar\varphi_\alpha\cap F|\ge 2$ and hence $\varphi_\alpha\subseteq F$, which contradicts the choice of $\varphi_\alpha\in \mathcal B_\alpha$. Therefore, $F\notin\F'_{<\alpha}$. In this case $\varphi_\alpha\cap F=\varnothing$ and hence $\bar\varphi_\alpha\cap F=\{(0,\varphi_\alpha(0)),(1,\bar a_\alpha\puls \varphi_\alpha(0))\}$. Since $F(0)=\varphi_\alpha(0)\ne 0$, the function $F$ does not belong to the family $\F_{<0}$ and hence $F\in \F_\beta$ for some $\beta<\alpha$. The inductive condition $(3_\beta)$ ensures that
 $$a_\beta\puls F(0)=F(1)=\bar a_\alpha\puls \varphi_\alpha(0)=\bar a_\alpha\puls F(0)$$ and hence $\bar a_\alpha=a_\beta$, which contradicts the choice of $\bar a_\alpha$. This contradiction shows that the function $\bar \varphi_\alpha$ satisfies the conditions (1)--(5). 
 \smallskip

If $1\in\dom[\varphi_\alpha]$ and $0\notin\dom[\varphi_\alpha]$, then we can consider the sets $$D\defeq \{0\}\cup \rng[\varphi_\alpha]\cup\{F(0):F\in\F'_{<\alpha}\}\cup\bigcup_{\beta<\alpha}\{x\in R:a_\beta\puls x=\varphi_\alpha(1)\},$$ and observe that $|D|<\kappa$. Then we can choose any point $x_\alpha\in R\setminus D$ and find a unique element $\bar a_\alpha\in R$ such that $\varphi_\alpha(1)=\bar a_\alpha\puls x_\alpha$. It follows from $x_\alpha\notin D$ that $\bar a_\alpha\notin\{a_\beta:\beta<\alpha\}$. It is easy to see that the function $\bar\varphi_\alpha\defeq\varphi_\alpha\cup\{(0,x_\alpha),(1,\bar a_\alpha,x_\alpha)\}$ is injective and satisfies the conditions (1)--(3). It follows from $\bar\varphi_\alpha(0)=x_\alpha\ne 0=L_{\bar a_\alpha}(0)$ and $\bar \varphi_\alpha(1)=\bar a_\alpha\puls x_\alpha\ne \bar a_{\alpha}\!\puls 0=\bar a_\alpha= L_{\bar a_\alpha}(1)$ that $\bar \varphi_\alpha\cap L_{\bar a_\alpha}=\varnothing$, so the condition (4) is satisfied. Assuming that the condition (5) does not hold, we could find a function $F\in\F_{<\alpha}$ such that $|\bar\varphi_\alpha\cap F|\ge 2$. If $F\in\F'_{<\alpha}$, then the choice of the element $x_\alpha$ ensures that $(0,x_\alpha)\notin F$. Then $|\varphi_\alpha\cap F|=|\bar\varphi_\alpha\cap F|\ge 2$ and hence $\varphi_\alpha\subseteq F$, which contradicts the choice of $\varphi_\alpha\in \mathcal B_\alpha$. Therefore, $F\notin\F'_{\alpha}$ and $\bar \varphi_\alpha\cap F=\{(0,x_\alpha),(1,\bar a_\beta\puls x_\alpha)\}$. Since $F(0)=x_\alpha\ne 0$, the function $F$ does not belong to the family $\F_{<0}$ and hence $F\in \F_\beta$ for some $\beta<\alpha$. The inductive condition $(3_\beta)$ ensures that $a_\beta\puls x_\alpha=F(1)=\bar a_\alpha\puls x_\alpha$ and hence $\bar a_\alpha=a_\beta$, which contradicts the choice of $\bar a_\alpha$. This contradiction shows that the function $\bar \varphi_\alpha$ satisfies the conditions (1)--(5). 
\end{proof}

If $\alpha\in S$, then put $a_\alpha\defeq\min_{\le}(R\setminus\{a_\beta:\beta<\alpha\})$, $\A_0\defeq\{L_{a_\alpha}\}$ and $\A_1\defeq\F_{<\alpha}$.
The inductive conditions $(2_\beta)$ and $(3_\beta)$ for $\beta<\alpha$ ensure that the elements $1,a_\alpha$ and the families $\A_0$ and $\A_1$ satisfy the hypothesis of  Lemma~\ref{l:Wilker2b}, which produces a family of bijections $\F_\alpha\subseteq S_X$ satisfying the conditions $(1_\alpha)$--$(4_\alpha)$.

 If $\alpha\in P$, then put $a_\alpha\defeq\bar a_\alpha$, where $\bar a_\alpha\in R$ and $\bar\varphi_\alpha$ satisfy the conditions (1)--(4) of Claim~\ref{cl:varphi-alpha}. Applying Lemma~\ref{l:Wilker1}, to the families $\{L_{a_\alpha}\}$, $\{L_{a_\alpha}\}\cup\F_{<\alpha}$ and the injective function $\bar \varphi_\alpha$, find a bijective function $\Phi_\alpha\in S_X$ such that $\bar\varphi_\alpha\subseteq\Phi_\alpha$, $\Phi_\alpha\cap L_{a_\alpha}=\varnothing$ and $\{\Phi_\alpha\}\perp (\F_{<\alpha}\setminus\{L_{a_\alpha}\})$. 


Applying Lemma~\ref{l:Wilker2b}, to the elements $1,a_\alpha\in R$ and the families  $\A_0\defeq\{L_{a_\alpha},\Phi_\alpha\}$ and $\A_1\defeq\A_0\cup\F_{<\alpha}$, find a family $\F_\alpha\subseteq S_X$ such that $\A_0\subseteq\F_\alpha$, $\bigcup\F_\alpha=X\times X$, $\F_\alpha\parallel \F_\alpha$, $\F\perp(\A_1\setminus\A_0)$, and $F(1)=a_\alpha\puls F(0)$ for every $F\in\F_\alpha$. Since $\varphi_\alpha\subseteq\bar\varphi_\alpha\subseteq\Phi_\alpha\in\F_\alpha$, the  family $\F_\alpha$ satisfies the inductive conditions $(1_\alpha)$--$(4_\alpha)$. This completes the inductive step.
\smallskip

After completing the inductive construction, we obtain a transfinite family $(\F_\alpha)_{\alpha\in\kappa\setminus \{0\}}$ of pairwise transversal disjoint covers $\F_\alpha$ of the set $R\times R$ by (graphs of) bijective functions. We claim that the set $R^2$ endowed with the family of lines $$\mathcal L\defeq\{\{c\}\times R:c\in R\}\cup\F_{<0}\cup\bigcup_{\alpha\in\kappa}\F_\alpha$$ is a Playfair liner. Indeed, the inductive conditions $(2_\alpha)$ for $\alpha\in\kappa$ imply $\mathcal L\pm\mathcal L$, which means that any two distinct lines $L,\Lambda\in\mathcal L$ have at most one common point.

We claim that $\{a_\alpha:\alpha\in\kappa\}=R$. In the opposite case, we can find an element $a\in R\setminus\{a_\alpha:\alpha\in \kappa\}$. The inductive conditions $(1_\alpha)$ for $\alpha\in S$ imply that $\{a_\alpha:\alpha\in S\}\subseteq{\downarrow}a$ and hence $|{\downarrow}a|\ge|S|=\kappa$, which contradicts the choice of the well-order $\le$ on the set $R$.

Next, we show that any distinct points $p,q\subseteq R^2$ belong to some line $L\in\mathcal L$. To derive a contradiction, assume that $\{p,q\}\not\subseteq L$ for every $L\in\mathcal L$. Then $\{p,q\}\pm\mathcal L$. Since the family $\mathcal L$ contains all horizontal and all vertical lines in $R\times R$, $\{p,q\}$ is an injective function. If $(0,0)\in\{p,q\}$, then $\{p,q\}\subseteq L_a\in\F_{<0}$ for some $a\in R\setminus\{0\}$, which contradicts our assumption. 
Then $(0,0)\notin\{p,q\}$ and hence the set $\{p,q\}\in \mathcal B$.
It follows from $\{p,q\}\pm\mathcal L$ that $\{p,q\}\in \mathcal B\setminus\bigcup_{\beta<\alpha}\bigcup_{F\in\F_\beta}\{\varphi\in\mathcal B:\varphi\subseteq F\}$. The inductive conditions $(4_\alpha)$ for $\alpha\in P$ imply that $\{\varphi_\alpha:\alpha\in P\}\subseteq{\downarrow}\{p,q\}$. Since the functions $\varphi_\alpha$, $\alpha\in P$, are pairwise distinct,
$|{\downarrow}\{p,q\}|\ge|P|=\kappa$, which contradicts the choice of the well-order $\preceq$ on the set $\mathcal B$. This contradiction shows that any distinct points of the set $R^2$ belong to a unique line $L\in\mathcal L$. Therefore, the pair $(R^2,\mathcal L)$ satisfies the axioms {\sf(L1)}, {\sf(L2)} of Theorem~\ref{t:L1+L2} and hence $(R^2,\mathcal L)$ is a liner. 

\begin{claim} The liner $(R^2,\mathcal L)$ is a Playfair plane.
\end{claim}

\begin{proof} Given any line $L\in \mathcal L$ and point $p\in R^2\setminus L$, we should find a unique line $L'\in\mathcal L$ such that $p\in L'\subseteq R^2\setminus L$.

If $L=\{c\}\times R$ for some $c\in R$, then there exists a unique point $c'\in R$ such that $p\in L'\defeq \{c'\}\times R$. In this case, $L'$ is a unique line such that $p\in L'\subseteq R^2\setminus\mathcal L$. 

If $L=R\times \{b\}$ for some $b\in R$, then there exists a unique point $b'\in R$ such that $p\in L'\defeq R\times \{b'\}$. In this case, $L'$ is a unique line such that $p\in L'\subseteq R^2\setminus\mathcal L$. 

So, we assume that $L\notin\{\{c\}\times R:c\in R\}\cup\{R\times\{b\}:b\in R\}$.
Then $L\in\F_\alpha$ for some non-zero ordinal $\alpha\in \kappa$. Since the family $\F_\alpha$ is a disjoint cover of $R^2$, there exists a unique line $L'\in\F_\alpha$ containing the point $p$. The inductive conditions $(3_\alpha)$, $\alpha\in\kappa$ imply that $\F_\alpha\perp(\mathcal L\setminus\F_\alpha)$, which implies that all lines in the family $\mathcal L\setminus\F_\alpha$ intersect the line $L$ and hence $L'$ is a unique line in $\mathcal L$ such that $p\in L'\subseteq R^2\setminus L$. Therefore, the liner $(R^2,\mathcal L)$ is a Playfair plane. 
\end{proof}

Consider the points $u\defeq (1,0)$, $o\defeq(0,0)$ and $w\defeq(0,1)$ in $R^2$ and observe that $uow$ is an affine base in the Playfair plane $(R^2,\mathcal L)$  whose  diagonal $\Delta=\Aline oe$ coincides with the diagonal $\{(x,y)\in R^2:x=y\}$ of the square $R^2$. Identifying every element $x\in R$ with the pair $(x,x)\in \Delta$, we can identify the set $R$ with the diagonal $\Delta$.
The choice of the lines $L_a=\{(x,y)\in R^2:y=x\cdot a\}$ and the families $\F_0=\big\{\{(x,y)\in R^2:y=x+b\}:b\in R\big\}$ and $\F_1=\big\{\{(x,y)\in R^2:y=x+b\}:b\in R\big\}$ ensures that the plus and dot operations on the diagonal $\Delta$ (identified with $R$) induced by the ternary operation $T_{uow}$ coincide with the plus and dot operations of the triloop $R$. 

Finally, we show that the puls operation of the  ternar $(\Delta,T_{uow})$ coincides with the puls operation $\puls$ of the triloop. Given any points $a,b\in\Delta$, consider the unique line $L_{a,b}\in\mathcal L$ such that $(0,b)\in L_{a,b}$ and $L_{a,b}\parallel L_a=\{(x,x{\cdot}a):x\in R\}$.  Since $R=\{a_\alpha:\alpha\in \kappa\}$, there exists an ordinal $\alpha\in\kappa$ such that $a_\alpha=a$. Then $L_{a,b}\in\mathcal F_{a_\alpha}$ and hence $T_{uow}(1,a,b)=L_{a,b}(1)=a_\alpha\!\puls L_{a,b}(0)=a\puls b$. So, the puls operations on the ternar $(\Delta,T_{uow})$ and the triloop $(R,+,\!\puls\!,\cdot)$ coincide (after identification of the diagonal $\Delta$ with $R$).
\end{proof}

\begin{remark} Example~\ref{ex:concurrent-non-planar3} of a non-planar concurrent triloop of cardinality 5 shows that Theorem~\ref{t:Hughes-Wilker} does not generalize to all (not necessarily infinite) triloops.
\end{remark}

\begin{remark} Theorem~\ref{t:Hughes-Wilker} shows that there is no symmetry between properties of the plus and puls operations in a ternar: whereas the plus and dot operations are bounded by the (two-parametric) concurence condition, the only restriction on the interplay between the puls and dot operations is the (one-parameter) condition ($\forall a\in R\setminus \{0\}\;\exists !x\in R\;\; 1\puls x=x\cdot a$) which follows from the plus-$1$-plus and concurrence conditions.
\end{remark}

\begin{proposition}\label{p:loop+=>biloop} For every infinite loop $(R,+)$ with identity element $0$ and every element $1\in R\setminus \{0\}$, there exists a binary operation $\cdot:R\times R\to R$ such that $(R,\cdot)$ is a $0$-loop with identity $1$, and the algebraic structure $(R,+,\cdot)$ is a concurrent biloop.
\end{proposition}

\begin{proof} Let $(R,+)$ be an infinite loop with neutral element $0$ and $1\in R\setminus\{0\}$ be any element. Let $R^*\defeq R\setminus\{0\}$ and $R^\circ\defeq R^*\setminus\{1\}$. Fix any well-order $\preceq$ on the set $R^*$ such that $1=\min_\preceq R^*$ and  for every $a\in R^*$, the set $$\cev a\defeq\{x\in R^*:a\ne x\preceq a\}$$ has cardinality $|\cev a|<|R^*|=|R|$.

\begin{lemma}\label{l:enumeration-Wilker} There exists an enumeration $\{(x_a,y_a)\}_{a\in R^\circ}$ of the set $R^\circ\times R^*$ such that for every $a\in R^\circ$ we have  $a\ne y_a\ne x_a+a_{-1}$ where $a_{-1}\in R$ is a unique element such that $1+a_{-1}=a$.
\end{lemma}

\begin{proof} Fix an enumeration $\{(x_a',y'_a)\}_{a\in R^\circ}$ of the set $R^\circ\times R^*$ such that for every pair $(x,y)\in R^\circ\times R^*$ the set $\{a\in R^\circ:(x_a',y_a')=(x,y)\}$ has cardinality $|R|$. Since $(R^*,\cdot)$ is a loop, for every $(x,y)\in R^\circ\times R^*$, the set $\{a\in R^\circ:y=x+a_{-1}\}$ contains at most one element. Then we can choose an element $\alpha(x,y)\in \{a\in R^\circ:(x_a',y_a')=(x,y)\;\wedge\; a\ne y\ne x+a_{-1}\}$. Consider the set $A\defeq \{\alpha(x,y):(x,y)\in R^\circ\times R^*\}\subseteq R^*$. For every $a\in A$, put $(x_a,y_a)\defeq(x'_a,y_a')$. For every $a\in R^\circ\setminus A$, let $(x_a,y_a)$ be any point of the set $\{(x,y)\in R^\circ\times R^*:a\ne y\ne x+a_{-1}\}$. Then the enumeration $R^\circ\times R^*=\{(x_a,y_a)\}_{a\in R^\circ}$ has the required property. 
\end{proof}
 
 By Lemma~\ref{l:enumeration-Wilker}, there exists an enumeration $\{(x_a,y_a)\}_{a\in R^\circ}$ of the set $R^\circ\times R^*$ such that $y_a\notin \{a,x_a+a_{-1}\}$ for every $a\in R^\circ$. 
 
 For every $b\in R$, consider the ``diagonal'' function
 $$
 \Delta_b\defeq\{(x,y)\in R\times R:y=x+b\},$$and the family
$$\F_1\defeq\big\{\Delta_b:b\in R\}$$of such ``diagonal'' functions. 
We shall construct inductively a transfinite sequence $(\Phi_a)_{a\in R^*}$ of bijective maps of $R$ such that $\Phi_1=\Delta_0$ and for every $a\in R^*$, the following conditions are satisfied:
\begin{itemize}
\item[$(1_a)$] $\{(0,0),(1,a)\}\subseteq \Phi_a$;
\item[$(2_a)$] if $a\ne 1$ and $(x_a,y_a)\notin \bigcup_{c\in\cev a}\Phi_c$, then $(x_a,y_a)\in \Phi_a$;
\item[$(3_a)$] $\{\Phi_a\}\perp(\F_1\cup\{\Phi_\alpha\}_{c\in\cev a})$.
\end{itemize}
Assume that for some element $a\in R^\circ$ a transfinite sequence of bijections $(\Phi_c)_{c\in\cev a}$ satisfying the inductive conditions $(1_c)$--$(3_c)$ for $c\in  \cev a$, have been constructed. Consider the families $\A_0\defeq\varnothing$ and $\A_1\defeq\F_1\cup\{\Phi_c\}_{c\in\cev a}$, and observe that $\A_1\perp\A_1\setminus\F_1$ and $\sup_{p\in R^2}|\{A\in\A_1:p\in A\}|\le 1+|\cev a|<\kappa$. 

If $(x_a,y_a)\in \bigcup_{c\in\cev a}\Phi_c$, then put $\varphi_a\defeq\{(0,0),(1,a)\}$. In the other case, put $$\varphi_a\defeq\{(0,0),(1,a),(x_a,y_a)\}.$$ 

Since $(x_a,y_a)\in R^\circ\times R^*$ and $y_a\ne a$, the set $\varphi_a$ is an injective function. We claim that $\{\varphi_a\}\perp \A_1$. The choice of the function $\varphi_a\ni (1,a)$ ensures that $\{\varphi_a\}\perp\{\Phi_c\}_{c\in\vec a}$. Next we show that $|\varphi_a\cap\Delta_b|\le 1$ for every $b\in R$. To derive a contradiction, assume that $|\varphi_a\cap \Delta_b|\ge 2$ for some $b\in R$. If $b=0$, then $(1,a)\in \varphi_a\setminus\Delta_0$ implies $(x_a,y_a)\in \varphi_a\cap\Delta_1$, which contradicts the definition of the function $\varphi_a$. This contradiction shows that $b\ne 0$ and hence $(0,0)\notin \Delta_b$. Then $|\varphi_a\cap\Delta_b|\ge 2$ implies $\{(x_a,y_a),(1,a)\}\subseteq \Delta_b$. Then $a=1+b$ and hence $b=a_{-1}$ and $y_a=x_a+a_{-1}$, which contradicts Lemma~\ref{l:enumeration-Wilker}. This contradiction shows that $\{\varphi_a\}\pm\F_1$ and hence $\varphi_a\not\subseteq \Phi$ for every $\Phi\in\F_1$.

By Lemma~\ref{l:Wilker1}, there exists a bijective function $\Phi_a:X\to X$ such that $\varphi_a\subseteq\Phi_a$ and $\{\Phi_a\}\perp\A_1$. It follows that the function $\Phi_a$ satisfies the inductive conditions $(1_a)$--$(3_a)$. This completes the inductive step.
\smallskip

 After completing the inductive construction, we obtain a family of bijective functions $(\Phi_a)_{a\in R^*}$ satisfying the inductive conditions $(1_a)$--$(3_a)$ for every $a\in R^*$. 
 
Define a binary operation $\cdot:R\times R\to R$ letting $x\cdot a\defeq\Phi_a(x)$, where $\Phi_0=R\times\{0\}$ is the constant zero function. We claim that $(R,\cdot)$ is a $0$-loop with zero $0$ and identity $1$. Indeed, for every $x\in R$ we have $x\cdot 0=\Phi_0(x)=0=\Phi_x(0)=0\cdot x$ and $x\cdot 1=\Phi_1(x)=x=\Phi_x(1)=1\cdot x$, witnessing that $0$ and $1$ are the zero and the identity of the magma $(R,\cdot)$. 

Given any elements $a\in R^*$ and $b\in R$, observe that $x\defeq \Phi_a^{-1}(b)$ is the unqiue solution of the equation $b=x\cdot a=\Phi_a(x)$. Next, we prove that the equation $\Phi_y(a)=a\cdot y=b$ has a unique solution $y\in R$. To derive a contradiction, assume that $\Phi_y(a)\ne b$ for all $y\in R$, and hence $(a,b)\notin\bigcup_{y\in R}\Phi_y$. Then $b\ne\Phi_0(a)=0$ and $a\ne 1$ (because $1\cdot b=\Phi_{b}(1)=b$). Find an element $c\in R^*$ such that $(a,b)=(x_c,y_c)$. Since $(x_c,y_c)=(a,b)\notin\bigcup_{y\in R}\Phi_y$, the inductive condition $(2_c)$ ensures that $(x_c,y_c)\in \Phi_c$ and hence $b=y_c=\Phi_c(x_c)=\Phi_c(a)=a\cdot c$, which contrdaicts our assumption. This contradiction shows that the equation $a\cdot y=b$ has a solution $y\in R$. Assuming that $u\in R\setminus\{y\}$ in another element with $a\cdot u=b$, we conclude that $(a,b)\in \Phi_u\cap\Phi_y=\{(0,0)\}$, which contradicts the choice of $a\in R^*$. This contradiction shows that $(R,\cdot)$ is a $0$-loop and hence $(R,+,\cdot)$ is a biloop.

Finally, we prove that the biloop $(R,+,\cdot)$ is concurrent. Given any elements $a\in R^*$ and $b\in R$ we should prove that the equation $x\cdot a=x+b$ has a unique solution $x\in R$. This follows from the inductive condition $(3_a)$, which guarantees that $|\Phi_a\cap \Delta_b|=1$ and hence $x\cdot a=y=x+b$ for the unique pair $(x,y)\in \Phi_a\cap \Delta_b$.
\end{proof}

Theorem~\ref{t:Hughes-Wilker} and Proposition~\ref{p:loop+=>biloop} imply the following corollary, proved by Wilker in \cite{Wilker1964} and \cite{Wilker1965}.

\begin{corollary}[Wilker, 1964]\label{c:Wilker1} For any infinite loop $(R,+)$ with identity element $0$ and any element $1\in R\setminus \{0\}$, there exists a ternary operation $T:R^3\to R$ such that $(R,T)$ is a ternar and $x+y=T(x,1,y)$ for every $x,y\in R$. 
\end{corollary}

The following corollary of Theorem~\ref{t:Hughes-Wilker} and Proposition~\ref{p:loop+=>biloop} shows that there is no algebraic restrictions on the structure of the plus and puls operations in an infinite ternar, except for the puls-$1$-plus identity (for finite ternars the situation is different, see Example~\ref{ex:plus1puls7}).

\begin{corollary}\label{c:WilkerBanakh} Let $R$ be a set with distinguished elements $0\ne 1$ and $+,\!\puls\!:R\times R\to R$ be two loop operations on $R$ with neutral element $0$ such that $1+x=1\puls x$ for all $x\in R$. Then there exists a ternary operation $T:R^3\to R$ such that $(R,T)$ is a ternar and for every $x,y\in R$ we have $x+y=T(x,1,y)$ and $x\puls y=T(1,x,y)$.
\end{corollary}

\begin{proposition}\label{p:loopx=>biloop} For every infinite $0$-loop $(R,\cdot)$, there exists a binary operation $+:R\times R\to R$ such that the algebraic structure $(R,+,\cdot)$ is a concurrent biloop.
\end{proposition}

\begin{proof} Let $0$ and $1$ be the zero and the identity of the $0$-loop $(R,\cdot)$. Fix any well-order $\preceq$ on the set $R^*\defeq R\setminus\{0\}$ such that $1=\min_{\preceq}R^*$ and for every $a\in R^*$, the set $$\cev a\defeq\{x\in R^*:a\ne x\preceq a\}$$ has cardinality $|\cev a|<|R^*|$. By a suitable modification of Lemma~\ref{l:enumeration-Wilker}, we can find an enumeration $\{(x_a,y_a)\}_{a\in R^*}$ of the set $R^*\times R$ such that  $y_a\ne a$ for every $a\in R^*$.

Consider the functions $$\Phi_0\defeq\{(x,y)\in R^2:y=x\}\quad\mbox{and}\quad \Psi_a\defeq\{(x,y)\in R^2:y=x\cdot a\}\quad\mbox{for $a\in R$}.$$  
We shall construct inductively a transfinite sequence $(\Phi_a)_{a\in R^*}$ of bijective maps of $R$ such that for every $a\in R^*$, the following conditions are satisfied:
\begin{itemize}
\item[$(1_a)$] $(0,a)\in \Phi_a$;
\item[$(2_\alpha)$] If $(x_a,y_a)\notin \bigcup_{c\in\cev a}\Phi_c$, then $(x_a,y_a)\in \Phi_a$;
\item[$(3_a)$] $\Phi_a\cap\bigcup_{c\in\cev a}\Phi_c=\varnothing$;
\item[$(4_a)$] $\{\Phi_a\}\perp\{\Psi_c\}_{c\in R^*}$.
\end{itemize}
Assume that for some element $a\in R^*$ a transfinite sequence of bijections $(\Phi_c)_{c\in\cev a}$ satisfying the inductive conditions $(1_c)$--$(4_c)$ for $c\in\cev a$, have been constructed. Consider the families $\A_0\defeq\{\Phi_c\}_{c\in\cev a}$ and $\A_1\defeq\A_0\cup\{\Psi_c\}_{c\in R^*}$. The inductive conditions $(3_c)$ and $(4_c)$ for $c\in\cev a$ ensure that 
$$\A_0\parallel \A_0,\quad \A_1\pm\A_1\quad\mbox{and}\quad\sup_{p\in R^2\setminus\bigcup\A_0}|\{A\in\A_1:p\in A\}|\le 2<\kappa.$$ 

If $(x_a,y_a)\in \bigcup_{c\in\cev a}\Phi_c$, then put $\varphi_a\defeq\{(0,a)\}$. In the other case, put $\varphi_a\defeq\{(0,a),(x_a,y_a)\}$. 

Since $(x_a,y_a)\in R^*\times (R^*\setminus\{a\})$, the set $\varphi_a$ is an injective function such that $\varphi_a\cap\bigcup\A_0=\varnothing$ and $\varphi_a\not\subseteq \Phi$ for every $\Phi\in\A_1$. By Lemma~\ref{l:Wilker1}, there exists a bijective function $\Phi_a:X\to X$ such that $\varphi_a\subseteq\Phi_a$, $\Phi_a\cap\bigcup\A_0=\varnothing$ and $\{\Phi_a\}\perp(\A_1\setminus\A_0)$. It follows that the function $\Phi_a$ satisfies the inductive conditions $(1_a)$--$(4_a)$.
This completes the inductive step.
\smallskip

 After completing the inductive construction, we obtain a family of bijective functions $(\Phi_a)_{a\in R^*}$ satisfying the inductive conditions $(1_a)$--$(4_a)$ for every $a\in R^*$. 
 
Define a binary operation $+:R\times R\to R$ letting $x+ a\defeq\Phi_a(x)$. We claim that $(R,+)$ is a loop whose identity is the zero element $0$ of the $0$-loop $(R,\cdot)$. Indeed, for every $x\in R$ we have $$x+0=\Phi_0(x)=x=\Phi_x(0)=0+x$$witnessing that $0$ is the identity of the magma $(R,+)$. 

Given any elements $a\in R$ and $b\in R$, observe that $x\defeq \Phi_a^{-1}(b)$ is the unique solution of the equation $b=x+ a=\Phi_a(x)$. Next, we prove that the equation $\Phi_y(a)=a+y=b$ has a unique solution $y\in R$. To derive a contradiction, assume that $\Phi_y(a)\ne b$ for all $y\in R$, and hence $(a,b)\notin\bigcup_{y\in R}\Phi_y$. Then $b\ne\Phi_0(a)=a$ and $a\ne 0$ (because $0+b=b$). Find an element $c\in R^*$ such that $(a,b)=(x_c,y_c)$. Since $(x_c,y_c)=(a,b)\notin\bigcup_{y\in R}\Phi_y$, the inductive condition $(2_c)$ ensures that $(x_c,y_c)\in \Phi_c$ and hence $b=y_c=\Phi_c(x_c)=\Phi_c(a)=a+c$, which contrdaicts our assumption. This contradiction shows that the equation $a+y=b$ has a solution $y\in R$. Assuming that $u\in R\setminus\{y\}$ in another element with $a+ u=b$, we conclude that $(a,b)\in \Phi_u\cap\Phi_y=\varnothing$, which is a contradiction showing that $(R,+)$ is a loop and hence $(R,+,\cdot)$ is a biloop.

Finally, we prove that the biloop $(R,+,\cdot)$ is concurrent. Given any elements $a\in R^*$ and $b\in R$ we should prove that the equation $x\cdot a=x+b$ has a unique solution $x\in R$. This follows from the inductive condition $(3_b)$, which guarantees that $|\Phi_b\cap \Psi_a|=1$. Then $x\cdot a=y=x+b$ for the unique pair $(x,y)\in \Phi_b\cap \Psi_a$.
\end{proof}

Theorem~\ref{t:Hughes-Wilker} and Proposition~\ref{p:loopx=>biloop} imply the following corollary, proved by Wilker in \cite{Wilker1964} and \cite{Wilker1965}.

\begin{corollary}[Wilker, 1964]\label{c:Wilker2} For every infinite $0$-loop $(R,\cdot)$, there exists a ternary operation $T:R^3\to R$ such that $(R,T)$ is a ternar and $x\cdot y=T(x,y,0)$ for every $x,y\in R$. 
\end{corollary}

\begin{remark}Theorem~\ref{t:Hughes-Wilker} and Corollaries~\ref{c:WilkerBanakh} and \ref{c:Wilker2} show that there are no algebraic restrictions on the structure of the additive loops $(R,+)$, $(R,\!\puls\!)$, and the multiplicative $0$-loop $(R,\cdot)$ of an infinite planar triloop. However, this fact does not generalize to finite planar triloops because  every planar triloop $R$ of order $|R|<9$ is a field and any two finite fields of the same cardinality are isomorphic (see also Example~\ref{ex:plus1puls7}).
\end{remark}

\begin{problem} Find the number of non-isomorphic (concurent) biloops of cardinality $n$, and also triloops of cardinality $n$ (which are plus-$1$-puls and/or concurrent).
\end{problem}

\section{The permutation group of a ternar}

\begin{definition} For a ternar $(R,T)$, let 
$$\Sym[R,T]\defeq\{T_{a,b}:(a,b)\in R^*\times R\}$$ 
be the set of permutations $T_{a,b}:R\to R$, $T_{a,b}:x\mapsto T(x,a,b)$ for $a\in R^*\defeq R\setminus\{0\}$ and $b\in R$. The subgroup $\Sym(R,T)$ of $\Sym(R)$, generated by the set $\Sym[R,T]$ is called \index{permutation group of a ternar}\index{ternar!permutation group of}\defterm{permutation group} of the ternar $(R,T)$.
\end{definition}

The axiom {\sf (T4)} of a ternar implies that the set $\Sym[R,T]$ is sharply $2$-transitive. This set contains a sharply transitive subset $$\Sym[R,+]\defeq\{T_{1,b}:b\in R\}$$ consisting of the permutations $T_{1,b}:R\to R$, $T_{1,b}:x\mapsto x+b\defeq T(x,1,b)$ for all $b\in R$. The subgroup of $\Sym(R)$, generated by the set $\Sym[R,+]$ will be denoted by $\Sym(R,+)$.

We recall that a permutation group $G\subseteq \Sym(X)$ is \defterm{of affine type} if it contains a sharply transitive normal Abelian subgroup. A permutation group $G\subseteq\Sym(X)$ is \defterm{sharply transitive} if for all points $x,y\in X$ there exists a  unique permutation $g\in G$ such that $g(x)=y$.

\begin{theorem} For a based affine plane $(X,uow)$ and its ternar $(\Delta,T)$, the following conditions are equivalent:
\begin{enumerate}
\item The liner $X$ is Thalesian.
\item The permutation group $\Sym(\Delta,T)$ of the ternar $(\Delta,T)$ is of affine type.
\item The sharply transitive set $\Sym[\Delta,+]$ is a normal Abelian subgroup in the permutation group $\Sym(\Delta,T)$ of the ternar $(\Delta,T)$.
\end{enumerate}
If the affine plane $X$ has finite order $|X|_2\notin\{3,4,24\}$, then \textup{(1)--(3)} are equivalent to the condition:
\begin{enumerate}
\item[$(4^\dag)$] $\Alt(\Delta)\not\subseteq \Sym(\Delta,T)$.
\end{enumerate}
\end{theorem}

\begin{proof} Let $\boldsymbol h\defeq(\Aline ou)_\parallel$ and $\boldsymbol v\defeq(\Aline ow)_\parallel$ be the horizontal and vertical directions of the affine base $uow$, and let $e$ be its diunit. Observe that for any elements $a,b\in\Delta$ with $a\ne o$, and the line $L\subseteq X$ defined by the equation $y=T(x,a,b)$, the permutation $T_{a,b}:\Delta\to\Delta$, $T_{a,b}:x\mapsto x_\times a_+b$, coincides with the permutation $\boldsymbol h_{\Delta,L}\boldsymbol v_{L,\Delta}$. Then
$$\Sym^{\Join}_X[\Delta;\boldsymbol h,\boldsymbol v]\defeq\{\boldsymbol h_{\Delta,L}\boldsymbol v_{L,\Delta}:L\in \mathcal L_X\setminus(\boldsymbol h\cup\boldsymbol v\}=\{T_{a,b}:(a,b)\in\Delta^*\times\Delta\}=\Sym[R,T]$$ and $$\Sym^\#_X[\Delta;\boldsymbol h,\boldsymbol v]\defeq\{\boldsymbol h_{\Delta,L}\boldsymbol v_{L,\Delta}:L\in \Delta_\parallel\}=\{T_{e,b}:b\in  \Delta\}=\Sym[\Delta,+].$$ Therefore, the permutation group $\Sym(\Delta,T)$ of the ternar $(\Delta,T)$
coincides with the permutation group $\Sym(\Delta;\boldsymbol h,\boldsymbol v)$. By Theorem~\ref{t:Schleiermacher}, the liner $X$ is Thalesian if and only if the permutation group $\Sym^{\Join}_X(\Delta;\boldsymbol h,\boldsymbol v)=\Sym(\Delta,T)$ is of affine type if and only if the sharply transitive set $\Sym^\#_X[\Delta;\boldsymbol h,\boldsymbol v]=\Sym[\Delta,+]$ is an Abelian normal subgroup of the group $\Sym^{\Join}_X(\Delta;\boldsymbol h,\boldsymbol v)=\Sym(\Delta,T)$. This completes the proof of the equivalences $(1)\Leftrightarrow(2)\Leftrightarrow(3)$.
\smallskip

Now assume that the affine plane $X$ has finite order $|X|\notin\{3,4,24\}$. If $X$ is not Thalesian, then $\Alt(\Delta)\subseteq \Sym_X^{\Join}(\Delta;\boldsymbol h,\boldsymbol v)=\Sym(\Delta,T)$, by Theorem~\ref{t:non-Thales-line-aff}. If $X$ is Thalesian, then permutation group $\Sym(\Delta,T)$ of the ternar $(\Delta,T)$ is of affine type and the sharply transitive set $\Sym[\Delta,+]$ is a normal Abelian subgroup in $\Sym(\Delta,T)$, by the (alread proved) implication $(1)\Ra(2\wedge 3)$ of this theorem. By Proposition~\ref{p:EAb<=>sharptrans}, the sharply transitive normal Abelian subgroup $\Sym[\Delta,+]$ of the finite $2$-transitive group $\Sym(\Delta,T)$ is elementary Abelian and hence $\Sym[\Delta,+]$ is isomorphic to a power $C_p^n$ of some cyclic group $C_p$ of prime order $p$. By Theorem~\ref{t:affine-type}, the permutation group $\Sym(\Delta,T)$ is isomorphic to a subgroup of the holomorph $\Hol(\Sym[\Delta,+])$, which is isomorphic  to the holomorph $\Hol(C_p^n)=C_p^n\rtimes\Aut(C_p^n)$ of the elementary Abelian group $C_p^n$. The automorphism group of the elementary Abelian group $C_p^n$ has cardinality $|\Aut(C_p^n)|=\prod_{k=0}^{k-1}(p^n-p^k)$, which can be found by counting bases of the vector space $C_p^n$ over the $p$-element field. Then $$|\Sym(\Delta,T)|\le|\Hol(C_p^n)|=|C_p^n|\cdot|\Aut(C_p^n)|=p^n\cdot\prod_{k=0}^{k-1}(p^n-p^k)<\tfrac12(p^n)!=|\Alt(\Delta)|$$
and hence $\Alt(\Delta)\not\subseteq \Sym(\Delta,T)$.
\end{proof}

\chapter{Subternars and characteristics of ternars}

In this chapter we consider and characterize subternars of ternars, prove the existence of minimal subternars and introduce the charactacteristic of a ternars (as the cardinality of the minimal subternar). Also we define the characteristic range of a Playfair liner and the set of all characteristics of its ternars.

\section{Subternars of ternars}

\begin{definition}\label{d:subternar} A subset $S$ of a ternar $(R,T)$ is called a \index{subternar}\defterm{subternar} of $R$ if $(S,T{\restriction}_{S^3})$ is a ternar.
\end{definition}

In Remark~\ref{rem:ternar-T2T3T4} we have observed that the ternary operation $T:R^3\to R$ of a ternar $R$ induces four auxiliary operations  $T_2:R^3\to R$ and $T_3,T'_4,T_4'':R^4_{\ne}\to R$ satisfying the identities
\begin{itemize}
\item[(T2)] $T(x,a,T_2(a,x,y))=y$,
\item[(T3)] $T(T_3(a,b,c,d),a,b)=T(T_3(a,b,c,d),c,d)$,
\item[(T4)] $T(x,T_4'(x,y,x',y'),T_4''(x,y,x',y'))=y$ and $T(x',T_4'(x,y,x',y'),T_4''(x,x',y,y'))=y'$,
\end{itemize}
for all $a,b,c,d,x,y,x',y'\in R$ with $a\ne c$ and $x\ne x'$. The quaternary operations $T_3,T_4',T_4''$ are defined on the set $R^4_{\ne}\defeq\{(a,b,c,d)\in R^4:a\ne c\}$.

\begin{proposition}\label{p:subternar<=>} A subset $S$ of a ternar $(R,T)$ is a subternar of $R$ if and only if $\{0,1\}\subseteq S$ and $T[S^3]\cup T_2[S^3]\cup T_3[S^4_{\ne}]\cup T_4'[S^4_{\ne}]\cup T_4''[S^4_{\ne}]\subseteq S.$
\end{proposition}

\begin{proof} The ``if'' part follows from the axioms of a ternar and definitions of the operations $T_2,T_3,T_4',T_4''$. To prove the ``only if'' part, assume that $S$ is a subternar of the ternar $R$.  Then $S$ contains two distinct elements $0'$ and $1'$ such that $T(x,0',0')=0'$ and $T(x,1',0')=x$ for all $x\in S$. By the axiom {\sf(T1)}, the element $0$ of the ternar $(R,T)$  satisfies the identity $T(x,0,0')=0'$ for all $x\in R$. Then the set $\{x\in R:T(x,0',0')=0'=T(x,0,0')\}$ contains the subset $S$ and hence has cardinailty at least $|S|\ge|\{0',1'\}|=2$. Applying the axiom {\sf(T3)}, we conclude that $0'=0$. The equality $0'=0$ and the axiom {\sf(T1)} ensure that the set $\{x\in R:T(x,1,0)=x=T(x,1',0)\}$ contains the set $S$ and hence is not a singleton. Applying the axiom {\sf(T3)}, we conclude that $1=1'$. Therefore, $\{0,1\}=\{0',1'\}\subseteq S$. It remains to prove that $$T[S^3]\cup T_2[S^3]\cup T_3[S^4_{\ne}]\cup T_4'[S^4_{\ne}]\cup T_4''[S^4_{\ne}]\subseteq S.$$
The inclusion $T[S^3]\subseteq S$ follows from the definition of the subternar $S$. To prove that $T_2[S^3]\subseteq S$, take any elements $x,a,y\in S$. By the axiom {\sf(T2)} applied to the ternar $S$, there exists a unique point $b\in S$ such that $T(x,a,b)=y$. The definition of the operation $T_2:R^3\to R$ ensures that $T_2(x,a,y)=b\in S$ and hence $T_2[S^3]\subseteq S$. To prove that $T_3[S^4_{\ne}]\subseteq S$, take any quadruple $(a,b,c,d)\in S^4_{\ne}$. By the axiom {\sf(T3)} applied to the ternar $S$, there exists a unique point $x\in S$ such that $T(x,a,b)=T(x,c,d)$. The definition of the function $T_3:R^4_{\ne}\to R$ ensures that $T_3(a,b,c,d)=x\in S$. To prove that $T_4'[S^4_{\ne}]\cup T_4''[S^4_{\ne}]\subseteq S$, take any quadruple $(x,y,x',y')\in S^4_{\ne}$. By the axiom {\sf(T4)} applied to the ternar $S$, there exist unique points $a,b\in S$ such that  
$T(x,a,b)=y$ and $T(x',a,b)=y'$. Now the definition of the functions $T_4',T_4'':R^4_{\ne}\to R$ ensures that $T_4'(x,y,x',y')=a\in S$ and $T_4''(x,y,x',y')=b\in S$ and hence $T_4'[S^4_{\ne}]\cup T_4''[S^4_{\ne}]\subseteq S$.
\end{proof}

Finite subternars admit a much simpler characterization.

\begin{theorem}\label{t:fin-subternar<=>} A finite subset $S$ of a ternar $(R,T)$ is a subternar of $R$ if and only if $S\not\subseteq\{0\}$ and $T[S^3]\subseteq S$.
\end{theorem}

\begin{proof} The ``only if'' part follows from Proposition~\ref{p:subternar<=>}. To prove the ``if'' part, assume that $S\not\subseteq\{0\}$ and $T[S^3]\subseteq S$. By Proposition~\ref{p:T1-T3=>T4}, it suffices to check that $(S,T{\restriction}_{S^3})$ satisfies the axioms {\sf(T1)}, {\sf(T1)}, {\sf(T3)} of a ternar.

First we prove that the constants $0$ and $1$ of the ternar $R$ belong to the set $S$. Since $S\not\subseteq\{0\}$, there exists a point $s\in S\setminus\{0\}$. Consider the function $T_{s,s}:R\to R$, $T_{s,s}:x\mapsto T(x,s,s)$. Since $s\ne 0$, the axiom {\sf(T3)} ensures that the function $T_{s,s}$ is injective. Consider the sequence $(s_n)_{n\in\w}$ defined by the recursive formula $s_0\defeq s$ and $s_{n+1}\defeq T(s_n,s,s)$ for all $n\in\w$. By induction, it can be shown that $s_n\in S$ for all $n\in\w$. Since $0\ne s$, the injectivity of the function $T_{s,s}$ ensures that $s_1=T_{s,s}(s)\ne T_{s,s}(0)=s=s_0$. Since the set $S$ is finite, there exists a  number $n>2$ such that $s_n=s_i$ for some positive integer number $i<n$. We can assume that $n$ is the smallest number with this property. Then $s_i\ne s_j$ for all numbers $i<j<n$. We claim that $s_n=s_0$. In the opposite case,  $s_n=s_i$ for some positive $i<n$. The equality $T_{s,s}(s_{n-1})=s_n=s_i=T_{s,s}(s_{i-1})$ and the injectivity of the function $T_{s,s}$ imply $s_{n-1}=s_{i-1}$, which contradicts the minimality of $n$. This contradiction shows that $s_n=s_0$. The equality $T_{s,s}(s_{n-1})=s_n=s_0=T_{s,s}(0)$ and the injectivity of the function $T_{s,s}$ ensure that $0=s_{n-1}\in S$.

To prove that $1\in S$, consider the injective function $T_{s,0}:R\to R$, $T_{s,0}:x\mapsto T(x,s,0)$ and observe that $T_{s,0}[S]\subseteq T[S^3]\subseteq S$. Since the set $S$ is finite, the injective function $T_{s,0}{\restriction}_S:S\to S$ is bijective and hence there exists an element $z\in S$ such that $T_{s,0}(z)=s=T(1,s,0)=T_{s,0}(1)$. The injectivity of the function $T_{s,0}:R\to R$ ensures that $1=z\in S$. Therefore, $\{0,1\}\subseteq S$ and the algebraic structure $(S,T{\restriction}_S^3)$ satisfies the axiom {\sf(T1)} of a ternar. 
\smallskip

To check the axiom {\sf(T2)} for $S$, fix any elements $x,a\in S$ and consider the function $T_{x,a}:S\to S$, $T_{x,a}:b\mapsto T(x,a,b)$. The axiom {\sf(T2)} holding for the ternar $(R,T)$ ensures that the function  $T_{x,a}$ is injective. Since the set $S$ is finite, the injective function $T_{x,a}:S\to S$ is bijective and hence for every $y\in S$ there exists a unique element $b\in S$ such that $T(x,a,b)=y$, witnessing that $(S,T{\restriction}_{S^3})$ satisfies the axiom {\sf(T2)} of a ternar.
\smallskip

To check the axiom {\sf(T3)} for $S$, fix any points $a,b,c,d\in S$ with $a\ne c$. Given any elements $\alpha,\beta\in S$, consider the set $L_{\alpha,\beta}=\{(x,y)\in S\times S:y=T(x,\alpha,\beta)\}$ and observe that it has cardinality $|L_{\alpha,\beta}|=|S|$. The axiom {\sf(T2)} holding for the ternar $(R,T)$ ensures that $L_{a,\beta}\cap L_{a,\gamma}=\varnothing$ for any distinct elements $\beta,\gamma$ of $S$. Consequently, $(L_{a,\beta})_{\beta\in S}$ is a disjoint cover of the set $S\times S$ of size $|S|^2$ by $|S|$ many sets of size $|S|$. The axiom {\sf(T3)} holding for the ternar $(R,T)$ ensures that $|L_{a,\beta}\cap L_{c,d}|\le 1$ for all $\beta\in S$. Since $|L_{c,d}|=|S|$, for every $\beta\in S$ the intersection $L_{a,\beta}\cap L_{c,d}$ is a singleton. In particular, there exists an element $(x,y)\in L_{a,b}\cap L_{c,d}\subseteq S\times S$ witnessing that $T(x,a,b)=y=T(x,c,d)$. The uniqueness of $x$ follows from the axiom {\sf(T3)} holding for the ternar $(R,T)$. This completes the proof of the axiom {\sf(T3)} for $S$.
\smallskip

 By Proposition~\ref{p:T1-T3=>T4}, $(S,T{\restriction}_S^3)$ also satisfies 
the axiom {\sf(T4)} and hence is a ternar.
\end{proof}

\begin{corollary} The set $\{0,1\}\subseteq R$ is a subternar of $R$ if and only if $1+1=0$.
\end{corollary}

\begin{exercise} Let $(R,T)$ be a ternar and $2\defeq T(1,1,1)\in R$. Prove that the set $\{0,1,2\}\subseteq R$ is a subternar of $R$ if and only if $T(1,1,2)=T(1,2,1)=T(2,1,1)=R(2,2,2)=0$ and $T(2,2,0)=T(2,1,2)=T(2,1,2)=1$.
\end{exercise}

\begin{proposition}\label{p:order-subternar} Let $S$ be a subternar of a ternar $R$. If $|S|<|R|$, then either $|S|^2=|R|$ or $|S|^2+|S|\le |R|$.
\end{proposition}

\begin{proof} If $R$ is infinite, then the strict inequality $|S|<|R|$ implies $|S|^2\le|S|^2+|S|<|R|$ and we are done. So, assume that $|R|$ is finite.

Let $R^2$ be the coordinate plane of the ternar $R$ and $\overline{R^2}$ be its projective completion. The affine plane $R^2$ and its projective completion $\overline{R^2}$ have the same order $|R|$.  Observe that the coordinate plane $S^2$ of the ternar $S$ is a subplane in $R$ and the projective completion $\overline{S^2}$ of the affine plane $S^2$ can be identified with a projective subplane of the projective plane $\overline{R^2}$. Since the projective subplane $\overline{S^2}$ of $\overline{R^2}$ has order $|S|$, Bruck's  Theorem~\ref{t:Bruck55} ensures that $|S|^2=|R|$ or $|S|^2+|S|\le|R|$.
\end{proof}

\section{Minimal ternars and their isomorphisms}

Proposition~\ref{p:subternar<=>} implies the following corollary.

\begin{corollary} Let $R$ be a ternar and $\mathcal S$ be any nonempty family of subternars of $R$. Then the set $\bigcap S$ is a subternar of $R$.
\end{corollary}

\begin{definition} A ternar $R$ is called \index{minimal ternar}\index{ternar!minimal}\defterm{minimal} if every subternar of $R$ coincides with $R$.
\end{definition}

Proposition~\ref{p:subternar<=>} implies the following theorem description of minimal subternars (in finite ternars).

\begin{theorem}\label{t:min-ternar-structure} Every ternar $R$ contains a unique minimal subternar $\underline{R}$, which is the intersection of the family $\mathcal S$ of all subternars of $R$. The minimal subternar $\underline{R}$ is the union $\bigcup_{n\in\w}S_n$ of the sequence of finite sets $(S_n)_{n\in\w}$ defined by the recursive formula $S_0=\{0,1\}$ and $S_{n+1}\defeq T[S^3_n]\cup T_2[S^3_n]\cup T_3[(S_n)^4_{\ne}]\cup T_4'[(S_n)^4_{\ne}]\cup T_4''[(S_n)^4_{\ne}]$ for $n\in\w$. If $\underline{R}$ is finite, then $\underline{R}$ is the union $\bigcup_{n\in\w}F_n$ of the sequence of  sets $(F_n)_{n\in\w}$ defined by the recursive formula $F_0=\{0,1\}$ and $F_{n+1}=T[F_n^3]$ for all $n\in\w$.
\end{theorem} 

\begin{theorem}\label{t:autominimal=>trivial} Every minimal ternar has trivial automorphism group.
\end{theorem}

\begin{proof} Let $A:R\to R$ be an automorphism of a minimal ternar $(R,T)$. By Theorem~\ref{t:min-ternar-structure}, $R=\bigcup_{n\in\w}S_n$, where $S_0\defeq\{0,1\}$ and $$S_{n+1}=T[S_n^3]\cup T_2[S_n^3]\cup T_3[(S_n)^4_{\ne}]\cup T_4'[(S_n)^4_{\ne}]\cup T_4''[(S_n)^4_{\ne}]$$ for all $n\in\w$.  By induction on $n$, we shall prove that $A(x)=x$ for all $x\in S_n$ and all $n\in\w$. 

Proposition~\ref{p:TR-01-unique} ensures that $A(0)=0$ and $A(1)=1$, and hence $A(x)=x$ for all $x\in S_0=\{0,1\}$. Assume that for some $n\in\w$ we have proved that $A(x)=x$ for all $x\in S_n$. Given any $s\in S_{n+1}$, we shall prove that $A(s)=s$. By definition of the set $S_{n+1}$, four cases are possible.
\smallskip

1. If $s\in T[S_n^3]$, then $s=T(x,a,b)$ for some $x,a,b\in S_n$. The inductive hypothesis ensures that $A(x)=x$, $A(a)=a$ and $A(b)=b$. Taking into account that $A$ is an automorphism of the ternar $(R,T)$, we conclude that  $A(s)=A(T(x,a,b))=T(A(x),A(a),A(b))=T(x,a,b)=s$.
\smallskip

2. If $s\in T_2[S_n^3]$, then $s=T_2(x,a,y)$ for some $x,a,y\in S_n$.  The definition of the function $T_2$ ensures that $T(x,a,s)=y$. Taking into account that $A$ is an automorphism of the ternar $(R,T)$ and $A(x)=x$, $A(a)=a$, $A(y)=y$, we conclude that $$T(x,a,s)=y=A(y)=A(T(x,a,s))=T(A(x),A(a),A(s))=T(x,a,A(s))$$
and hence $s=A(s)$, by the axiom {\sf (T2)} of a ternar.
\smallskip

3. If $s\in T_3[(S_n)^4_{\ne}]$, then $s=T_3(a,b,c,d)$ for some $a,b,c,d\in S_n$ with $a\ne c$. The definition of the function $T_3$ ensures that $T(s,a,b)=T(s,c,d)$. Taking into account that $A$ is an automorphism of the ternar $(R,T)$ and $A(a)=a$, $A(b)=b$, $A(c)=c$, $A(d)=d$, we conclude that 
\begin{multline*}T(A(s),a,b)=T(A(s),A(a),A(b))=A(T(s,a,b))\\=
A(T(s,c,d))=T(A(s),A(c),A(d))=T(A(s),c,d)
\end{multline*}
and hence $A(s)=s$ by the uniqueness part of the axiom {\sf (T3)}.
\smallskip

4. If $s\in T_4'[(S_n)^4_{\ne}]\cup T_4''[(S_n)^4_{\ne}]$, then $s\in\{T_4'(x,y,x',y'),T_4''(x,y,x',y')\}$ for some $x,y,x',y'\in S_n$ with $x\ne x'$. Let $a\defeq T_4'(x,y,x',y')$ and $b\defeq T_4''(x,y,x',y')$, and observe that $s\in\{a,b\}$. The axiom {\sf(T4)} and the definition of the functions $T_4'$ and $T_4''$ ensure that $a,b$ are unique elements of the ternar $R$ such that $T(x,a,b)=y$ and $T(x',a,b)=y'$.
 Taking into account that $A$ is an automorphism of the ternar $(R,T)$ and $A(x)=x$, $A(y)=y$, $A(x')=x'$, $A(y')=y'$, we conclude that $$T(x,A(a),A(b))=T(A(x),A(a),A(b))=A(T(x,a,b))=A(y)=y=T(x,a,b)$$ and
$$T(x',A(a),A(b))=T(A(x'),A(a),A(b))=A(T(x',a,b))=A(y')=y'=T(x',a,b).$$
The uniqueness of the elements $a,b$ implies that $A(a)=a$ and $A(b)=b$ and hence $A(s)=s$. 
\smallskip

Therefore, $A(x)=x$ for all $x\in S_n$ and all $n\in\w$, and hence $A(x)=x$ for all $x\in\bigcup_{n\in\w}S_n=R$. So, $A$ is the trivial automorphism of the minimal ternar $R$.  
\end{proof}

\begin{remark}\label{rem:min-ternar-iso-test} Theorem~\ref{t:autominimal=>trivial} implies that for two minimal ternars $X,Y$, there exists at most one isomorphism $F:X\to Y$. This observation suggests a simple isomorphism test for minimal ternars. Given two minimal ternar $(X,T,0,1)$ and $(X',T',0',1')$ of the same finite cardinality $n=|X|=|X'|$, construct a sequence of relations $(F_n)_{n\in\w}$ defined by the recursive formulas $F_0\defeq \{(0,0'),(1,1')\}$ and $F_{n+1}\defeq\{(T(x,y,z),T'(x',y',z')):\{(x,x'),(y,y'),(z,z')\}\subseteq F_n\}$ for $n\in\w$. The calculation of the sets $F_n$ proceeds as long as $F_{n+1}\ne F_n$ and $F_n$ is an injective function. Therefore, the calculations terminates whenever either $F_n$ is not an injective function or $F_n$ is an injective function with $F_n=F_{n+1}$. In the first case, the ternars $X,X'$ are not isomorphic, in the second case the function $F_n=F_{n+1}$ is an isomorphism of the ternars $X,X'$. The computation complexity of this algorithm is $O(n^3)$ (under the assumption that the caclulations of the values of the operations $T$ and $T'$ has computation complexity $O(1)$).
\end{remark}

\section{The characteristic of a ternar}

\begin{definition} The \index{characteristic of a ternar}\index{ternar!characteristic of}\defterm{characteristic} $\har(R)$ of a ternar $R$ is defined as the cardinality $|\underline{R}|$ of the minimal subternar $\underline{R}$ of $R$. Therefore, $\har(R)\defeq|\underline{R}|$. 
\end{definition}

\begin{remark} A finite ternar $R$ is minimal if and only if $\har(R)=|R|$.
\end{remark}

\begin{proposition}\label{p:har-Bruck} The characteristic $p\defeq\har(R)$ of a ternar $R$ satisfies the following conditions:
\begin{enumerate}
\item If $p\le 10$, then $p\in\{2,3,5,7,9\}$;
\item $p=|R|$ or $p^2=|R|$ or $p^2+p\le|R|$.
\end{enumerate}
\end{proposition}

\begin{proof} 1. If $p\le 10$, then the minimal ternar $\underline R$ of $R$ has order $|\underline R|=p\le10$. The coordinate plane $\underline R^2$ is a Playfair plane of order $|\underline R|=p\le 10$. By a result of Lam, Thiel and Swiercz \cite{LTS1989}, no Playfair plane of order 10 exists, which implies $p\le 9$. If $p=9$, then $\underline R^2$ is a Playfair plane of order $9$. By Remark~\ref{rem:char9}, $p=\har(\underline R)\in\{2,3,9\}$. So, assume that $<9$. In this case the Playfair plane $\underline R^2$ is Desarguesian, by \cite{HSW1956} (see also Corollary~\ref{c:4-Pappian}). Then $(\underline R,+,\cdot)$ is a minimal field and its cardinality $p$ is a prime number, equal to its  characteristic. So, $p\in\{2,3,5,7\}$.
\smallskip

2. Since $p=|\underline R|$, Proposition~\ref{p:order-subternar} ensures that  $p=|R|$ or $p^2=|R|$ or $p^2+p\le|R|$.
\end{proof}

\begin{remark}
Proposition~\ref{p:order-subternar} suggests the following (relatively) fast algorithm of calculation of the minimal subternar of a finite ternar $(R,T)$: calculate the sequence of sets $(S_n)_{n\in\w}$ defined by the recursive formula $S_0\defeq\{0,1\}$ and $S_{n+1}\defeq T[S_n^3]$ for all $n\in\w$ satisfying $|S_{n+1}|>|S_{n}|$ and $|S_{n}|\le |R|^{\frac12}$. So, calculations will terminate if $|S_{n+1}|=|S_{n}|$ or $|S_{n+1}|>|R|^{\frac12}$. In the first case we put $\underline{R}\defeq S_{n+1}=S_n$. In the second case, put $\underline{R}\defeq R$. This algorithm has computational compexity $O(|R|^{3/2})$ (under the assumption that finding the value of the function $T(x,y,z)$ for a given triple $(x,y,z)\in R$ has computation complexity $O(1)$).
\end{remark}

\begin{remark} Computer calculations (by Ivan Hetman) show that every linear ternar of order $9$ has characteristic $2$ or $3$.
\end{remark}

\begin{problem} Is it true that the characteristic of every finite linear ternar is a prime number? 
\end{problem}

\section{The characteristic of a triangle in an affine liner}

Given a liner $X$, denote by $X^\vartriangle\defeq\{xyz\in X^3:\|\{x,y,z\}\|=3\}$ the family of all triangles in $X$. 

\begin{definition} A triangle $uow$ in a liner $X$ is called \index{Playfair triangle}\index{triangle!Playfair}\defterm{Playfair} if the flat hull $\overline{uow}$ of the set $\{u,o,w\}$ is a Playfair plane.
\end{definition}

\begin{remark} If an affine liner $X$ contains a Playfair triangle, then $|X|_2\ge 3$. 
\end{remark} 

Theorem~\ref{t:4-long-affine} implies the following proposition.

\begin{proposition} Every triangle in an affine liner $X$ of order $|X|_2\ge 4$ is Playfair.
\end{proposition}

\begin{exercise} Find an affine liner $X$ of order $|X|_2=3$ containing a non-Playfair triangle.
\smallskip

{\em Hint:} Look at Example~\ref{ex:Tao}.
\end{exercise}

\begin{definition} Let $uow$ be a triangle in an affine liner. If the triangle $uow$ is Playfaor, then its \index{ternar of a triangle}\index{triangle!ternar of}\defterm{ternar} $\Delta_{uow}$ is defined as the ternar of the based affine plane $(\overline{uow},uow)$. If $uow$ is not Playfair, then we put $\Delta_{uow}$ be the field $\IF_p$ of characteristic $p=|X|_2$, endowed with the ternary operation $x_\times y_+z\defeq (x{\cdot}y)+z$ for all $x,y,z\in \IF_p$. .
\end{definition}

\begin{definition} The \index{characteristic of a triangle}\index{triangle!characteristics of}\defterm{characteristic} $\har(uow)$ of a triangle $uow$ in an affine liner $X$ is defined as the characteristic $\har(\Delta_{uow})$ of its ternar $\Delta_{uow}$. Therefore, $$\har(uow)\defeq\har(\Delta_{uow}).$$ A triangle $uow\in X^\vartriangle$ is called \index{minimal triangle}\index{triangle!minimal}\defterm{minimal} if its ternar $\Delta_{uow}$ is minimal, i.e., contains no proper subternars. 
\end{definition}

\begin{remark} A triangle $uow$ in a finite affine liner $X$ is minimal if and only if $\har(uow)=|X|_2$.
\end{remark}

\begin{definition} A triangle $uow$ in a Playfair liner $X$ is \index{Boolean triangle}\index{triangle!Boolean}\defterm{Boolean} if it can be completed to a Boolean parallelogram $uowe$.
\end{definition}

\begin{proposition}\label{p:Boolean<=>har2} A triangle $uow$ in a Playfair liner $X$ is Boolean iff $\har(uow)=2$.
\end{proposition}

\begin{proof} Let $e$ be the diunit of the affine base $uow$ in the Playfair plane $\overline{\{u,o,w\}}$ and let $\Delta$ be the ternar of the based affine plane $(\overline{\{u,o,w\}},uow)$.
\smallskip

If the triangle $uow$ is Boolean, then the parallelogram $uowe$ is Boolean and hence $e+e=o$ in the ternar $\Delta_{uow}$ of the affine base $uow$, see Lemma~\ref{l:Boolean}. Then $\{o,e\}$ is a subternar of $\Delta_{uow}$, which implies that $\har(uow)=\har(\Delta_{uow})=|\{o,e\}|=2$.
\smallskip

If $\har(uow)=2$, then $\{o,e\}$ is the minimal subternar of the ternar $\Delta$ and hence $e+e=o$, witnessing that the parallelogram $uowe$ is Boolean, according to Lemma~\ref{l:Boolean}.
\end{proof} 

Observe that every triangle $uow$ in a Playfair plane $X$ is an affine base for the plane $\overline{\{u,o,w\}}$ in $X$. Let $\Delta_{uow}$ be the ternar of the based affine plane $(\overline{\{u,o,w\}},uow)$ and $\underline{\Delta}_{uow}$ be its minimal subternar. Let $\underline{\Pi}_{uow}$  be the set of all points $x\in\overline{uow}$ whose coordinates in the affine base $uow$ belong to the minimal subternar $\underline{\Delta}_{uow}$ of the ternar $\Delta_{uow}$. 

\begin{theorem}\label{t:minternar<=>minplane} Let $uow$ be a triangle in a Playfair liner $X$. If $\har(uow)>2$, then $\underline{\Pi}_{uow}$ is the Playfair subplane of $X$. Moreover, $\underline{\Pi}_{uow}$ is the smallest Playfair subliner of $X$ that contains the points $u,o,w$.
\end{theorem}

\begin{proof} Assume that $\har(uow)>2$. Let $e$ be the diunit of the affine base $uow$ and $\Delta_{uow}$ be the ternar of the affine base $uow$. 
Since $\har(uow)>2$, the minimal subternar $\underline{\Delta}_{uow}$ of the ternar $\Delta_{uow}$ contains more than two elements. Since the liner $\underline{\Pi}_{uow}$ is isomorphic to the coordinate plane of the ternar $\underline{\Delta}_{uow}$, it is a Playfair plane, by Theorem~\ref{t:coord-plane=>based-aff-plane}.
\smallskip

Now fix any Playfair subliner $\Pi\subseteq X$, containing the points $u,o,w$. Then $uow$ is an affine base for the liner $\Pi$. By Theorem~\ref{t:parallelogram3+1}, the diunit of the affine base $uow$ belongs to the Playfair subliner $\Pi$ of $X$. Then $\Delta_X\defeq\Delta_{uow}\cap \Pi$ is the underlying set of the ternar of the based affine plane $(\Pi,uow)$. Since $\Pi$ is a subliner of the liner $X$ and $e\in X$, the ternar of the based affine plane $(\Pi,uow)$ is a subternar of the ternar $\Delta_{uow}$. The minimality of the subternar $\underline{\Delta}_{uow}$ ensures that $\underline{\Delta}_{uow}\subseteq\Delta_X$ and hence $\underline{\Pi}_{uow}\subseteq X$ (by the Playfair property of $X$).
\end{proof}

\section{The characteristic range of an affine liner}

\begin{definition} For an affine liner $X$, the range $\har[X^\vartriangle]\defeq\{\har(uow):uow\in X^\vartriangle\}$ of the function $\har:X^\vartriangle\to \IN$ is called the \index{characteristic range}\index{Playfair liner!characteristic range of}\defterm{characteristic range} of  $X$. The function $$\har_X^-:\har[X^\vartriangle]\to \IN,\quad \har^-:n\mapsto|\{uow\in X^\vartriangle:\har(uow)=n\}|$$is called the \index{characteristic level function}
\index{Playfair liner!characteristic level function of}\defterm{characteristic level function} of $X$.
\end{definition} 

Proposition~\ref{p:har-Bruck} implies the following corollariey describing the structure of the characteristic range of an affine liner $X$.

\begin{corollary} Let $X$ be an affine liner of order $n$. Then $$\har[X^\vartriangle]\subseteq\{p\in\IN\cup\{\w\}:p=n\;\vee\;p^2=n\;\vee\;6\le p^2+p\le n\}.$$ 
\end{corollary}

\begin{corollary} Let $X$ be an affine liner and $p\in\har[X^\vartriangle]$. If $p\le 10$, then $p\in\{2,3,5,7,9\}$.
\end{corollary}

\begin{proposition} An affine liner $X$ is invertible-add if $\har[X^\vartriangle]\subseteq \{p\in\IN:p<9\}$
\end{proposition}

Proposition~\ref{p:Boolean<=>har2} implies the following two characterizations.

\begin{corollary}\label{c:Boolean<=>2har} A Playfair liner $X$ contains a Boolean parallelogram iff  $2\in\har[X^\vartriangle]$.
\end{corollary}

%

\begin{corollary}\label{c:Boolean<=>char2} A Playfair liner $X$ is Boolean if and only if $\har[X^\vartriangle]=\{2\}$.
\end{corollary}



\begin{remark}\label{rem:Thales=>unicharacteristic} In Theorem~\ref{t:BT=>f-elementary}  we shall prove that the characteristic range of any Thalesian Playfair liner $X$ is a singleton.
\end{remark}

Corollary~\ref{c:Boolean<=>char2} and Theorem~\ref{t:BT=>f-elementary}  motivate the following definition.

\begin{definition}\label{d:unicharacteristic-affine} An affine liner $X$ is called \index{unicharacteristic Playfair liner}\index{Playfaor liner!unicharacteristic}\defterm{unicharacteristic} if its characteristic range $\har[X^\vartriangle]$ is a singleton. The unique number in $\har[X^\vartriangle]$ is called \index{the characteristic}\defterm{the characteristic} of the unicharacteristic affine liner $X$.
\end{definition}


\begin{proposition}\label{p:char-divides-order} If $X$ is a unicharacteristic affine liner of order $n=|X|_2$, then the characteristic $p$ of $X$ divides $n$ and $p-1$ divides $n-1$.
\end{proposition}

\begin{proof} Since the liner $X$ is unicharacteristic, the characteristic range $\har[X^\vartriangle]$ is a singleton and hence is not empty. Then $X$  contains a triangle and hence has rank $\|X\|\ge 3$. If $X$ is not Playfair, then 
$n=|X|_2\le 3$, by Theorem~\ref{t:4-long-affine}. In this case, the ternar $\Delta_{uow}$ of every triangle $uow\in X^\vartriangle$ equals the minimal field $\IF_p$ of characteristic $p=n=|X|_2\le 3$. Then $p=n$ divides $n$ and $p-1=n-1$ divides $n-1$. 

So, assume that the affine liner $X$ is Playfair and hence $n=|X|_2\ge 3$, by Theorem~\ref{t:Playfair<=>}. For every triangle $uow$ in $X$, let $\underline{\Pi}_{uow}$ be the subliner of $X$, consisting of all points $x\in X$ whose both coordinates in the affine base $uow$ belong to the minimal subternar $\underline{\Delta}_{uow}$ of the ternar $\Delta_{uow}$. Since $\har(\Delta_{uow})=p$, the subplane $\underline{\Pi}_{uow}$ has cardinality $|\underline{\Pi}_{uow}|=p^2$.

\begin{claim}\label{cl:min-subplanes3=>coincide} Let $abc$ and $uvw$ be two triangles in $X$. If $\|\underline{\Pi}_{abc}\cap\underline{\Pi}_{uow}\|>2$, then $\underline{\Pi}_{abc}=\underline{\Pi}_{uow}$.
\end{claim}

\begin{proof} Assuming that $\|\underline{\Pi}_{abc}\cap\underline{\Pi}_{uow}\|>2$, find a triangle $xyz$ in the intersection $\underline{\Pi}_{abc}\cap\underline{\Pi}_{uow}$. Applying Theorem~\ref{t:minternar<=>minplane}, we conclude that $\underline{\Pi}_{xyz}\subseteq\underline{\Pi}_{abc}\cap\underline{\Pi}_{xyz}$. Taking into account that $|\underline{\Pi}_{abc}|=|\underline{\Pi}_{xyz}|=|\underline{\Pi}_{upw}|=p^2$, we conclude that $\underline{\Pi}_{abc}=\underline{\Pi}_{xyz}=\underline{\Pi}_{uow}$.
\end{proof}

Now we are able to prove the divisibility properties the numbers $p$ and $p-1$.
\smallskip
 
1. Fix any disjoint lines $L,\Lambda$ in $X$ and any distinct points $a,b\in \Lambda$. For every $x\in L$, consider the set $L_x\defeq L\cap \underline{\Pi}_{xab}$ and observe that $|L_x|=\har(\Delta_{xab})=p$. Claim~\ref{cl:min-subplanes3=>coincide} ensures that $\{L_x:x\in L\}$ is a partition of the set $L$ into  subsets  of cardinality $p$. Therefore, $p$ divides $|L|=n$.
\smallskip

2. Fix any lines $L,\Lambda$ in $X$ having a unique common point $o$. Fix any point $\lambda\in \Lambda\setminus \{o\}$. Let $L'\defeq L\setminus\{o\}$. For any point $x\in L'$, consider the set $L'_x\defeq L'\cap\underline{\Pi}_{xo\lambda}$ and observe that $|L'|=p-1$. Claim~\ref{cl:min-subplanes3=>coincide} ensures
$\{L_x':x\in L'\}$ is a partition of the set $L'$ into subsets of cardinality $p-1$, which implies that $p-1$ divides $|L'|=n-1$.
\end{proof}


\begin{remark}\label{rem:char9} The characteristic ranges and the characteristic level functions of all seven affine planes of order $9$ are presented in the following table,  calculated by Ivan Hetman.
$$
\begin{array}{c|c|c}
\hline
\mbox{affine plane}&\mbox{characteristic range}&\mbox{characteristic level function}\\
\hline
\mbox{\tt Desarg}&\{3\}&\{(3,466560)\}\\
\mbox{\tt Thales}&\{3\}&\{(3,466560)\}\\
\mbox{\tt Hall}&\{2,3,9\}&\{(2,124416),(3,186624),(9,155520)\}\\
\mbox{\tt hall}&\{2,3,9\}&\{(2,92160),(3,51840),(9,322560)\}\\
\mbox{\tt dhall}&\{2,3,9\}&\{(2,96768),(3,62208),(9,307584)\}\\
\mbox{\tt hughes}&\{2,3,9\}&\{(2,62208),(3,62208),(9,342144)\}\\
\mbox{\tt Hughes}&\{2,3,9\}&\{(2,62208),(3,93312),(9,311040)\}\\
\hline
\end{array}
$$
\end{remark}

\begin{remark} The above table shows that the characteristic level function distinguishes all non-Thalesian affine planes of order $9$. This is not true for affine planes of order 16: by computer calculations, Ivan Hetman has found four nonisomorphic affine planes (moreover, those affine planes have nonisomorphic projective completions: {\tt math, dmath, jowk, lmrh}) whose characteristic level functions coincide and are equal to $\{(2,7409664),(16,8257536)\}$. 
\end{remark}

The results of calulations show that every non-Thalesian affine plane of order $9$ has a minimal affine base. This suggests the following definition.

\begin{definition} An affine liner is called \index{miminal Playfair plane}\index{Playfair plane!minimal}\defterm{minimal} if it contains a minimal triangle.
\end{definition}

Let us recall that a triangle in a Playfair plane is \defterm{minimal} if its ternar is minimal (i.e., contains no proper subternars).

\begin{remark} A finite affine liner $X$ is minimal if and only if $|X|_2\in\har[X^\vartriangle]$.
\end{remark} 

\begin{problem} Is every non-Thalesian affine liner minimal?
\end{problem}

An affine liner  $X$ is defined to be \defterm{prime} if its order $|X|_2$ is a prime number. Proposition~\ref{p:char-divides-order} implies the following corollary.

\begin{corollary} Every unicharacteristic prime Playfair liner $X$ is minimal and every ternar of $X$ is minimal.
\end{corollary}

\begin{problem} Is a finite Playfair plane $X$ prime if all ternars of $X$ are minimal?
\end{problem}

\begin{problem}\label{prob:af-unichar=>prime} Is the characteristic of every unicharacteristic finite Playfair plane a prime number?
\end{problem}

\begin{problem}\label{prob:unichar=>Boolean-or-Thalesian} Is every unicharacteristic finite Playfair plane Boolean or Thalesian?
\end{problem}

\begin{problem}\label{prob:unichar-prime=>Pappian} Is every unicharacteristic prime Playfair plane Pappian?
\end{problem}

\begin{remark} Since the characteristic of any Thalesian liner is prime and every Thalesian Playfair plane of prime order is Pappian, the affirmative answer to Problem~\ref{prob:unichar=>Boolean-or-Thalesian} implies the affirmative answer to Problems~\ref{prob:unichar=>prime} and \ref{prob:unichar-prime=>Pappian}. 
\end{remark}

\begin{remark}\label{rem:min-Playfair-isomorphism-fast} The isomorphism test for minimal ternars described in Remark~\ref{rem:min-ternar-iso-test} suggests a simple isomorphism test for minimal finite Playfair planes. Given two minimal Playfair planes $X,Y$ of the same finite order $n$, fix any minimal affine base $uow$ in $X$ and calculate its ternar $\Delta_{uow}$. Then for every triangle $xyz$ in the Playfair plane $Y$, calculate its characteristic $\har(xyz)$ and if $\har(xyz)=n$ then using the isomorphism test for minimal ternars, find a (unique) isomorphism between the minimal ternars $\Delta_{uow}$ and $\Delta_{xyz}$, if it exists. If for all triangles $xyz$ in $Y$, the ternars $\Delta_{uow}$ and $\Delta_{xyz}$ are not isomorphic, then the Playfair planes $X,Y$ are not isomorphic, according to Theorem~\ref{t:tring-iso<=>}. This algorithm of recognizing (non)isomorphic minimal Playfair planes of order $n$ has computational complexity $O(n^9)$. 
\end{remark}



\chapter{Inversive, Bol and Moufang loops}

\rightline{\vbox{\hsize=260pt \em \noindent Algebra is the offer made by the devil to the mathematician.  The devil says: ``I will give you this powerful machine, it will answer any question you like.  All you need to do is give me your soul: give up geometry and you will have this marvellous machine.''}}

\rightline{Michael \index[person]{Atiyah}Atiyah\footnote{{\bf Sir Michael Atiyah} (1929--2019) was a British mathematician renowned for his groundbreaking work in geometry and topology. Born in London, he studied at Trinity College, Cambridge, and later held professorships at Oxford, Cambridge, and the University of Edinburgh. Atiyah made fundamental contributions to K-theory, the Atiyah--Singer Index Theorem, and topological quantum field theory, bridging pure mathematics with theoretical physics. He received numerous honors, including the Fields Medal (1966), the Abel Prize (2004), and a Knighthood (1983). Atiyah was known not only for his profound research but also for his eloquent reflections on the philosophy of mathematics, often emphasizing the interplay between geometry and algebra.}}
\vskip30pt

In this chapter we present basic information of some important classes of loops: alternative, inversive, \index[person]{Bol}Bol\footnote{{\bf Gerrit Door Bol} (1906 --1989) was a Dutch mathematician who specialized in geometry. He is known for introducing Bol loops in 1937, and Bol’s conjecture on sextactic points. Bol earned his PhD in 1928 at Leiden University under Willem van der Woude. In the 1930s, he worked at the University of Hamburg on the geometry of webs under Wilhelm Blaschke and later projective differential geometry. In 1931 he earned a habilitation. In 1933 Bol signed the Loyalty Oath of German Professors to Adolf Hitler and the National Socialist State. In 1942--1945 during World War II, Bol fought on the Dutch side, and was taken prisoner. On the authority of Blaschke, he was released. After the war, Bol became professor at the Albert-Ludwigs-University of Freiburg, until retirement there in 1971.} and Moufang. We shall meet such loops studying the algebraic structure of non-Desarguesian planes in Chapters~\ref{ch:PlusAff}, \ref{ch:PulsAff} and \ref{ch:dot-aff}.

\section{Alternative and flexible magmas}

\begin{definition} A magma $(X,\cdot)$ is called 
\begin{itemize}
\item \index{associative magma}\index{magma!associative}\defterm{associative} if $x{\cdot}(y{\cdot}z)=(x{\cdot}y){\cdot}z$ for every $x,y,z\in X$;
\item \defterm{left-alternative} if $\forall x,y\in X\;\;(x{\cdot} x){\cdot} y=x{\cdot}(x{\cdot} y)$;
\item \defterm{right-alternative} if $\forall x,y\in X\;\;x{\cdot}(y{\cdot} y)=(x{\cdot} y){\cdot} y$;
\item \defterm{alternative} if $X$ is left-alternative and right-alternative;
\item \defterm{flexible} if $\forall x,y\in X\;(x{\cdot} y){\cdot} x=x{\cdot}(y{\cdot} x)$.
\end{itemize}
\end{definition}

Associativity, alternativity and flexibility can be also expressed via properties of left and right shifts. Let us recall that for an element $a\in X$ of a magma $(X,\cdot)$, the maps
$$L_a:X\to X,\quad L_a:x\mapsto a{\cdot}x,\quad\mbox{and}\quad
R_a:X\to X,\quad R_a:x\mapsto x{\cdot}a,$$
are called the \defterm{left} and \defterm{right shifts} of $X$ by the element $a$.

\begin{proposition} A magma $(X,\cdot)$ is  
\begin{enumerate}
\item left-alternative if and only if $L_xL_x=L_{x{\cdot}x}$ for every $x\in X$;
\item right-alternative if and only if $R_xR_x=R_{x{\cdot}x}$ for every $x\in X$;
\item alternative if and only if $L_xL_x=L_{x{\cdot}x}$ and $R_xR_x=R_{x{\cdot}x}$ for every $x\in X$;
\item flexible if and only if $R_xL_x=L_xR_x$ for every $x\in X$;
\item associative if and only if $L_xL_y=L_{x{\cdot}y}$ and $R_xR_y=R_{y{\cdot}x}$ for every $x,y\in X$.
\end{enumerate}
\end{proposition}

We say that an element $a$ of a unital magma $X$ has \defterm{order $2$} if $a\ne e=a\cdot a$, where $e$ is the neutral element of the unital magma $X$.
 
\begin{proposition} Assume that a finite loop $X$ is left-alternative or right-alternative. If $X$ contains an element of order $2$, then the cardinal $|X|$ of $X$ is even.
\end{proposition}

\begin{proof} Let $a\in X$ be an element of order $2$ in the unital magma $X$. First assume that the loop $X$ is left-alternative. Consider the left shift $L_a:X\to X$, $L_a:x\mapsto a\cdot x$. The left-alternativity of $X$ ensures that $L_aL_a=L_{a{\cdot}a}=L_e$ is the identity map of $X$. Since $X$ is a quasigroup, for every $x\in X$, the inequality $a\ne e$ implies $L_a(x)=a\cdot x\ne e\cdot x=x$. Then $\{\{x,a{\cdot}x\}:x\in X\}$ is a partition of $X$ into doubletons, which implies that the cardinal $|X|$ is even. A similar argument shows that $X$ is even, if the loop $X$ is right-alternative. 
\end{proof}

\begin{definition}
A magma $(X,\cdot)$ is called
\begin{itemize}
\item \index{mono-associative magma}\index{magma!mono-associative}\defterm{mono-associative} if $(x{\cdot} x){\cdot} x=x{\cdot}(x{\cdot}x)$ for every $x\in X$;
\item \index{${\IN}$-associative magma}\index{magma!$\IN$-associative}\defterm{$\IN$-associative} if for every $x\in X$ there exists a sequence $(x^n)_{n\in\IN}$ such that $x^1=x$ and $x^n{\cdot} x^m=x^{n+m}$ for every $n,m\in\IN$;
\item  \index{$\IZ$-associative magma}\index{magma!$\IZ$-associative}\defterm{$\IZ$-associative} if for every $x\in X$ there exists a sequence $(x^n)_{n\in\IZ}$ such that $x^1=x$ and $x^n{\cdot} x^m=x^{n+m}$ for every $n,m\in\IZ$.
\end{itemize}
\end{definition}

For every magma the following implications hold:
$$
\xymatrix@C=-10pt{
&\mbox{associative}\ar@{=>}[d]&\\
&\mbox{alternative}\ar@{=>}[ld]\ar@{=>}[dr]\\
\mbox{left-alternative}\ar@{=>}[dr]&\mbox{flexible}\ar@{=>}[d]&\mbox{right-alternative}\ar@{=>}[dl]\\
&\mbox{mono-associative}\\
&\mbox{$\IN$-associative}\ar@{=>}[u]\\
&\mbox{$\IZ$-associative}\ar@{=>}[u]
}
$$

\begin{exercise} Prove the implications in the above diagram.
\end{exercise}

An $\IN$-associative loop $X$ is called
\begin{itemize} 
\item \index{$d$-divisible loop}\index{loop!$d$-divisible}\defterm{$d$-divisible} for $d\in\IN$ if for every element $x\in X$ there exists an element $y\in X$ such that $x=y^d$;
\item \index{locally cyclic group}\defterm{locally cyclic} if every finite subset of $X$ is contained in a cyclic subgroup of the loop $X$.
\end{itemize}

\begin{exercise} Prove that every locally cyclic loop is a commutative group.
\end{exercise}

\begin{proposition}\label{p:d-divisible-locally-cyclic} Let $d\in\IN$. Any element $x$ of a $d$-divisible $\IZ$-associative loop $X$ is contained in a $d$-divisible locally cyclic subgroup of $X$.
\end{proposition}

\begin{proof} Since $X$ is $d$-divisible, there exists a sequence $(x_n)_{n\in\w}$ such that $x_0=x$ and $x_n=x_{n+1}^d$ for all $n\in\w$. 
By the $\IZ$-associativity of $X$, for every $n\in\w$, the set $C_n\defeq\{x_n^m:m\in\IZ\}$ is a cyclic subgroup of $X$. Since $x_n=x_{n+1}^d$, the cyclic group $C_n$ is a subgroup of the cyclic group $C_{n+1}$. Then the union $C\defeq\bigcup_{n\in\w}C_n$ is a required $d$-divisible locally cyclic subgroup of the loop $X$.
\end{proof}

\begin{remark} In literature, $\IN$-associative and $\IZ$-associative magmas both are called \defterm{power-associative}. We have desided to modify the terminology in order to distinguish those two notions. However we shall also use the standard terminology calling $\IN$-associative magma power-associative magmas and $\IZ$-associative loops \defterm{power-associative loops}. An example of an infinite loop which is $\IN$-associative but not $\IZ$-associative was constructed by \index[person]{Smith}Warren Smith\footnote{{\bf Warren D. Smith} (born in 1964) is an Americal  mathematician and researcher whose work spans nonassociative algebra, computational mathematics, and combinatorics. He earned his Ph.D. in Applied Mathematics from Princeton University in 1989 under the supervision of John H. Conway and Robert E. Tarjan, with a dissertation ``Studies in Computational Geometry Motivated by Mesh Generation''. Smith is known for his contributions to loop theory, including studies of diassociativity, Moufang and Bol loops, and related structural properties, as well as the development of computational tools for exploring loop identities. Among his notable results is the demonstration that diassociativity in loops has no finite equational basis. In addition to algebra, his research interests include graph algorithms, computational complexity, and applications such as voting systems. Smith’s work combines rigorous algebraic theory with computational techniques, reflecting a distinctive interdisciplinary approach.} \cite[\S4.1]{WSmith}.
\end{remark}

\begin{remark} A more detail diagram of the interplay between various properties of loops can be found in the paper \cite{WSmith} of Warren Smith. The mentioned paper contains many examples of finite loops distinguishing various properties of loops. We present some examples of such distinguishing loops below.
\end{remark}

\begin{example}[Warren Smith] The $5$-element loop with multiplication table
$$
\begin{array}{c|ccccc}
*& 0& 1& 2& 3& 4\\
\hline
0&0&1&2&3&4\\
1&1&2&4&0&3\\
2&2&3&1&4&0\\
3&3&4&0&2&1\\
4&4&0&3&1&2
\end{array}
$$
is not mono-associative.
\end{example}

\begin{example}[Warren Smith] The $7$-element loop with multiplication table
$$
\begin{array}{c|ccccccc}
*&0&1&2&3&4&5&6\\
\hline
0&0&1&2&3&4&5&6\\
1&1&2&3&4&5&6&0\\
2&2&3&6&1&0&4&5\\
3&3&0&1&5&6&2&4\\
4&4&5&0&6&2&1&3\\
5&5&6&4&2&3&0&1\\
6&6&4&5&0&1&3&2
\end{array}
$$
is mono-associative but is neither left-alternative nor right-alternative.
\end{example}

\begin{example}[Warren Smith] The $6$-element loop with multiplication table
$$
\begin{array}{c|ccccccc}
*&0&1&2&3&4&5\\
\hline
0&0&1&2&3&4&5\\
1&1&0&3&4&5&2\\
2&2&3&0&5&1&4\\
3&3&5&4&0&2&1\\
4&4&2&5&1&0&3\\
5&5&4&1&2&3&0
\end{array}
$$
is mono-associative but is not left-alternative and not right-alternative. 
\end{example}

\begin{example}[Warren Smith] The $6$-element loop with multiplication table
$$
\begin{array}{c|cccccc}
*&0&1&2&3&4&5\\
\hline
0&0&1&2&3&4&5\\
1&1&2&3&4&5&0\\
2&2&3&4&5&0&1\\
3&3&0&5&2&1&4\\
4&4&5&0&1&2&3\\
5&5&4&1&0&3&2
\end{array}
$$
is left-alternative, but not right-alternative.
\end{example}

\section{Cancellative magmas}

\begin{definition} A magma $(X,\cdot)$ is called 
\begin{itemize}
\item \defterm{left-cancellative} if for every $a\in X$, the left shift $L_a:X\to X$, $L_a:x\mapsto a{\cdot}x$, is injective;
\item \defterm{right-cancellative} if for every $a\in X$, the right shift $R_a{:}X\to X$, $R_a{:}x\mapsto x{\cdot}a$, is injective;
\item \defterm{cancellative} if $X$ is left-cancellative and right-cancellative;
\item \defterm{semi-cancellative} if $X$ for every $a\in X$ at least one of the shifts $L_a$ or $R_a$ is injective.
\end{itemize}
\end{definition}

\begin{exercise} Show that a finite magma is a quasigroup if and only if it is cancellative.
\end{exercise}

\begin{proposition}\label{p:N-ass<=>Z-ass} A semi-cancellative finite magma is $\IN$-associative if and only if it is $\IZ$-associative.
\end{proposition}

\begin{proof} The ``if'' part is trivial. To prove the ``only if'' part, assume that a semi-cancellative finite magma $(X,\cdot)$ is power-$\IN$-associative.  Then for every $x\in X$, there exists a sequence $(x^n)_{n\in\IN}$ such that $x^1=x$ and $x^n\cdot x^m=x^{n+m}$ for all $n,m\in\IN$. Since $X$ is finite, there exist positive integer number $b$ and $p$ such that $x^b=x^{b+p}$.
We can assume that $b$ is the smallest numbers with those properties. More precisely, that $x^i\notin\{x^{i+p}:p\in\IN\}$ for all $i<b$.

We claim that $b=1$. Since the magma $X$ is semi-cancellative, one of the shifts $L_x$ or $R_x$ is injective. We lose no generality assuming that the left shift $L_x$ is injective. Assuming that $b>1$, we can observe that $L_x(x^{b-1})=x^{b}=x^{b+p}=L_x(x^{b-1+p})$ and hence $x^{b-1}=x^{b-1+p}$, by the injectivity of the left shift $L_x$. But the equality $x^{b-1}=x^{b-1+p}$ contradicts the minimality of $b$. This contradiction shows that $b=1$. 

We claim that $x^n=x^{n+p}$ for every $n\in\IN$. For $n=1$ this equality follows from the choice of $p$. Assume that $x^n=x^{n+p}$ for some $n\in\IN$. Then $x^{n+1}=x^n\cdot x=x^n\cdot x^{1+p}=x^{n+1+p}$. The Principle of Mathematical Induction ensures that $x^n=x^{n+p}$ for all $n\in\IN$, which implies $x^n=x^{n+kp}$ for every $n,k\in\IN$. For every integer $n<0$, let $x^n\defeq x^{n+kp}$, where $k\in\IN$ is the smallest number such that $n+kp\in\IN$. We claim that the sequence $(x^n)_{n\in\IZ}$ has the required property: $x^n\times x^m=x^{n+m}$ for every $n,m\in\IZ$. Given any numbers $n,m\in \IZ$, find smallest numbers $i,j,k\in\w$ such that the numbers $n+ip,m+jp,(n+m)+kp$ are positive. Then $n+ip+m+jp=(n+m)=(n+m)+(i+j)p\ge 1+1>0$ and hence $k\le i+j$ and $x^{(n+m)+kp}=x^{n+m+kp+(i+j-k)p}$. It follows that
$$x^n\cdot x^m=x^{n+ip}\cdot x^{m+jp}=x^{n+ip+m+jp}=x^{n+m+(i+j)p}=x^{(n+m)+kp+(i+j-k)p}=x^{(n+m)+kp}=x^{n+m},$$
witnessing that the magma $X$ is $\IZ$-associative.
\end{proof}

Since every quasigroup is cancellative, Proposition~\ref{p:N-ass<=>Z-ass} implies

\begin{corollary} A finite quasigroup is $\IZ$-associative if and only if it is $\IN$-associative.
\end{corollary}




\section{Inversive magmas}


\begin{definition} A magma $(X,\cdot)$ is called
\begin{itemize}
\item \index{left-inversive}\index{magma!left-inversive}\defterm{left-inversive} if $\forall a\in X\;\exists b\in X\;\forall x\in X\;\;b{\cdot}(a{\cdot} x)=x$;
\item  \index{right-inversive}\index{magma!right-inversive}\defterm{right-inversive} if $\forall a\in X\;\exists b\in X\;\forall x\in X\;\;(x{\cdot} a){\cdot} b=x$;
\item  \index{inversive}\index{magma!inversive}\defterm{inversive} if $\forall a\in X\;\exists b\in X\;\forall x\in X\;\; b{\cdot} (a{\cdot} x)=x=(x{\cdot} a){\cdot} b$.
\end{itemize}
\end{definition}

\begin{proposition}\label{p:inversive<=>lr-inversive} A unital magma is inversive if and only if it is left-inversive and right-inversive.
\end{proposition}

\begin{proof} The ``only if'' part is trivial. To prove the ``if'' part, assume that a unital magma $X$ is left-inversive and right-inversive. Then for every $a\in X$ there exist elements $b,c\in X$ such that $b{\cdot} (a{\cdot} x)=x=(x{\cdot} a){\cdot} c$ for every $x\in X$. The inversivity of $X$ will follow as soon as we check that $b=c$. Applying the identity $x=b{\cdot} (a{\cdot} x)$ to the neutral element $e$ of the unital magma $X$, we obtain $e=b{\cdot}(a{\cdot} e)=b{\cdot} a$ and hence
$c=e{\cdot} c=(b{\cdot} a){\cdot} c=b$.
\end{proof}

The inversivity can be also expressed via properties of left and right shifts. 
The following simple (but important) characterization follows immediately from the definitions. 

\begin{proposition} A quasigroup $X$ is 
\begin{enumerate}
\item left-inversive if and only if for every element $a\in X$ there exists an element $a^{-1}\in X$ such that $L_a^{-1}=L_{a^{-1}}$;
\item right-inversive if and only if for every element $a\in X$ there exists an element $a^{-1}\in X$ such that $R_a^{-1}=R_{a^{-1}}$;
\item inversive if and only if for every element $a\in X$ there exists an element $a^{-1}\in X$ such that $L_a^{-1}=L_{a^{-1}}$ and $R_a^{-1}=R_{a^{-1}}$.
\end{enumerate}
\end{proposition}

\begin{definition} A unital magma $(X,\cdot)$ with the identity $e$ is called
\begin{itemize}
\item \index{invertible}\index{magma!invertible}\defterm{invertible} if $\forall a\in X\;\exists b\in X \;\;a\cdot b=e=b\cdot a$;
\item  \index{Boolean}\defterm{Boolean} if $\forall x\in X\;\;x\cdot x=e$;
\end{itemize}
\end{definition}

It is clear that every Boolean unital magma is invertible.

\begin{proposition}\label{p:invertible<=lr-inv} A unital magma is invertible if it is left-inversive or right-inversive.
\end{proposition}

\begin{proof} Assume that a unital magma $X$ is left-inversive. Then for every $a\in X$ there exists an element $b\in X$ such that $b{\cdot}(a{\cdot} x)=x$ for all $x\in X$. Also, for the element $b$, there exists an element $c\in X$ such that $c{\cdot}(b{\cdot} x)=x$ for all $x\in X$. In particular, for the neutral element $e$ of the unital magma $X$, we have $e=b{\cdot} (a{\cdot} e)=b{\cdot} a$ and $e=c{\cdot}(b{\cdot} e)=c{\cdot} b$. 
Then $$a{\cdot} b=c{\cdot} (b{\cdot} (a{\cdot} b))=c{\cdot} b=e=b{\cdot} a,$$witnessing that the unital magma $X$ is invertible. 

By analogy we can prove that every right-inversive unital magma is invertible.
\end{proof}

\begin{proposition}[Yurii Yarosh, 2025] Every inversive unital magma is a loop.
\end{proposition}

\begin{proof} Let $X$ be an inversive unital magma. By the left-inversivity of $X$, for every $a\in X$ there exisa an element $b\in X$ such that $L_{b}L_a$ is  the identity map of $X$, which implies that the left shift $L_a$ is injective. By the left-inversivity of $X$, for the element $b$ there exists an element $c\in X$ such that $L_cL_b$ is the identity map of $X$, which implies that $L_c$ is surjective. We claim that $c=a$. Indeed, $e=b{\cdot}(a{\cdot}e)=b{\cdot}a$ and $a=c{\cdot}(b{\cdot}a)=c{\cdot}e=c$. Therefore, the left shift $L_a=L_c$ is biejctive. By analogy we can prove that for every $a\in X$, the right shift $R_a:X\to X$ is bijective, witnessing that the unital magma $X$ is a quasigroup and hence a loop.
\end{proof}

By definition of an invertible loop, for every element $x$ of an invertible loop $(X,\cdot)$, there exists an element $x^{-1}\in X$ such that $x{\cdot}x^{-1}=e=x^{-1}{\cdot}x$. Since $X$ is a loop, the element $x^{-1}$ is unique. This unique element $x^{-1}$ is called \defterm{the two-sided inverse} of $x$ in the loop $X$. The map $i:X\to X$, $i:x\mapsto x^{-1}$, is called the \defterm{inversion}. The uniqueness of the inverse element implies that $(x^{-1})^{-1}=x$ for every $x\in X$, which means that the inversion of an invertible loop is an involutive bijection of the loop. The following proposition shows that the inversion of an inversive loop is an anti-isomorphism of the loop.

\begin{proposition}\label{p:inversive=>anti} If $X$ is an inversive loop, then for every elements $x,y\in X$ we have $(x\cdot y)^{-1}=y^{-1}\cdot x^{-1}$.
\end{proposition}

\begin{proof} The left-inversivity of $X$ implies $x^{-1}{\cdot} (x{\cdot} y)=y$, and the right-inversivity of $X$ ensures that $$x^{-1}=(x^{-1}{\cdot}(x{\cdot} y)){\cdot}(x{\cdot} y)^{-1}=y{\cdot}(x{\cdot} y)^{-1}$$ and hence
$$(x{\cdot} y)^{-1}=y^{-1}{\cdot}(y{\cdot}(x{\cdot} y)^{-1})=y^{-1}{\cdot} x^{-1},$$by the left-inversivity of $X$.
\end{proof} 

\begin{proposition}\label{p:inversiveBoolean=>commutative} Every inversive Boolean loop is commutative.
\end{proposition}

\begin{proof} Let $X$ be an inversive Boolean loop with neutral element $e$.
For every $x\in X$, the equality $x\cdot x=e$ implies that $x^{-1}=x$. Applying Proposition~\ref{p:inversive=>anti}, we conclude that $$x{\cdot}y=(x{\cdot}y)^{-1}=y^{-1}{\cdot}x^{-1}=y{\cdot}x,$$
witnessing that $X$ is commutative.
\end{proof} 

The following simple (but useful) facts were observed by \index[person]{Asif Ali}Asif Ali\footnote{{\bf Asif Ali} currently works as a Full Professor at Department of Mathematics in Quaid-i-Azam University in Islamabad (Pakistan). His research interests include Algebra, Finite Group Theory, Finite Loop Theory, Group Theoretic Technique to solve non-linear differential equations, Functional Analysis, Teaching Mathematics, Coding Theory, Logic \&\ Automated Reasoning.} and \index[person]{Slaney}John Slaney\footnote{{\bf John Slaney}, Professor Emeritus in the School of Computing at the Australian National University. Research interests in the Computing Foundations and Intelligent Systems clusters. Former leader of the Logic and Computation Program in NICTA (National ICT Centre of Excellence, Australia). Author of a cool web site, providing LOGIC FOR FUN ({\tt https://logic4fun.cecs.anu.edu.au/about}). Instigator and convenor (1993--2019) of the annual LOGIC SUMMER SCHOOL. Other likes: travel, good food (enthusiastic but inexpert cook), classical music (ditto pianist). In his autobiography, John Slaney writes: ``I was born in England but escaped, taught logic in philosophy departments for several years, escaped again and moved to Canberra in 1988 where I have been automating reasoning ever since. I like doing this. The ANU is an idyllic place to be a researcher, Canberra is a better city to live in than you would believe from listening to Australians from anywhere else, and I actually get paid for thinking about logic and hacking code! That's as good as it gets.''} in \cite{AS2008}.

\begin{proposition}\label{p:Ali-Slaney} Any invertible loop $X$ of even order contains an element of order $2$.
\end{proposition}

\begin{proof} Since $X$ is an invertible loop, for every $x\in X$, there exists a unique element $x^{-1}\in X$ such that $x\cdot x^{-1}=e=x^{-1}\cdot x$. The uniqueness of the inverse element implies $\{\{x,x^{-1}\}:x\in X\}$ is a partition of $X$ into pairwise disjoints sets. Assuming that $x\ne x^{-1}$ for any $x\in X\setminus\{e\}$, we conclude that the set $X\setminus\{e\}=\bigcup_{x\in X\setminus\{e\}}\{x,x^{-1}\}$ has even cardinality and hence $|X|$ is odd, which contradicts our assumption. This contradiction shows that $a=a^{-1}$ for some $a\in X\setminus\{e\}$. Then $a\cdot a=a\cdot a^{-1}=e\ne a$, which means that $X$ is an element of order $2$.
\end{proof}

\begin{proposition}\label{p:Ali-Slaney2} Assume that a finite loop $X$ is left-inversive or right-inversive. The loop $X$ contains an element of order $2$ if and only if $|X|$ is even.
\end{proposition}

\begin{proof} The ``only if'' part follows from Propositions~\ref{p:invertible<=lr-inv} and \ref{p:Ali-Slaney}. 
To prove the ``if'' part, assume that $X$ contains an element $a\in X$ of order $2$. Then $a\ne e=a\cdot a$. If $X$ is left-inversive, then $a\cdot (a\cdot x)=a^{-1}\cdot (a\cdot x)=x$ for all $x\in X$, which implies that the left shift $L_a:X\to X$, $L_a:x\mapsto a\cdot x$, is an involution of $X$. Since $X$ is a loop, the inequality $a\ne e$ implies $L_a(x)=a\cdot x\ne e\cdot x=x$ for all $x\in X$. Then $\{\{a{\cdot}x,x\}:x\in X\}$ is a partition of $X$ into doubletons, which implies that the cardinal $|X|$ is even. By analogy we can prove that $|X|$ is even if $X$ is right-inversive. 
\end{proof}

\begin{corollary}[Ali, Slaney, 2008]\label{c:Ali-Slaney} A (left or right) inversive loop of odd order cannot contain a subloop of even order.
\end{corollary}

\begin{proof} To derive a contradition, assume that a (left or right) inversive loop $X$ of odd order contains a subloop $S$ of even order.  By Proposition~\ref{p:Ali-Slaney2}, the loop $S$ contains an element $a$ of order $2$, and by Proposition~\ref{p:Ali-Slaney2}, the loop $X$ has even order.
\end{proof}

The properties of loops, discussed in this section relate as follows:
$$
\xymatrix@C=-5pt{
&\mbox{associative}\ar@{=>}[d]\\
&\mbox{inversive}\ar@{=>}[ld]\ar@{=>}[rd]&\\
\mbox{left-inversive}\ar@{=>}[rd]&\mbox{Boolean}\ar@{=>}[d]
&\mbox{right-inversive}\ar@{=>}[ld]\\
&\mbox{invertible}\\
&\mbox{commutative}\ar@{=>}[u]
}
$$

\begin{example}[Warren Smith] The $6$-element loop with multiplication table
$$
\begin{array}{c|cccccc}
*&0&1&2&3&4&5\\
\hline
0&0&1&2&3&4&5\\
1&1&2&3&0&5&4\\
2&2&5&0&4&1&3\\
3&3&0&4&5&2&1\\
4&4&3&5&1&0&2\\
5&5&4&1&2&3&0\\
\end{array}
$$
is invertible but is not Boolean, not left-inversive, not right-inversive, not commutative and not power-associative.
\end{example}

\begin{example}[Warren Smith] The $6$-element loop with multiplication table
$$
\begin{array}{c|cccccc}
*&0&1&2&3&4&5\\
\hline
0&0&1&2&3&4&5\\
1&1&2&3&4&5&0\\
2&2&4&0&5&1&3\\
3&3&5&4&0&2&1\\
4&4&3&5&1&0&2\\
5&5&0&1&2&3&4\\
\end{array}
$$
is left-inversive, but not right-inversive.
\end{example}

\section{Bol magmas}

\begin{definition} A magma $(X,\cdot)$ is called
\begin{itemize} 
\item \index{left-Bol magma}\index{magma!left-Bol}\defterm{left-Bol} if $a{\cdot}(b{\cdot}(a{\cdot} x))=(a{\cdot} (b{\cdot} a)){\cdot} x$ for every elements $a,b,x\in X$;
\item \index{right-Bol magma}\index{magma!right-Bol}\defterm{right-Bol} if $((x{\cdot}a){\cdot} b){\cdot} a=x{\cdot} ((a{\cdot} b){\cdot} a)$ for every elements $x,a,b\in X$;
\item \index{Bol magma}\index{magma!Bol}\defterm{Bol} if $X$ is both left-Bol and right-Bol.
\end{itemize}
\end{definition}

Bol magmas can be characterized via left and right shifts as follows.

\begin{proposition} A magma $(X,\cdot)$ is 
\begin{enumerate} 
\item left-Bol if and only if $L_aL_bL_a=L_{a{\cdot}(b{\cdot}a)}$ for every elements $a,b\in X$;
\item right-Bol if and only if $R_aR_bR_a=R_{(a{\cdot}b){\cdot}a}$ for every elements $a,b\in X$.
\end{enumerate}
\end{proposition}

\begin{proposition}\label{p:left-Bol=>left-alt+flex} Every left-Bol unital magma is left-alternative. 
\end{proposition}

\begin{proof} Let $(X,\cdot)$ be a unital magma and $e$ be its neutral element. If the magma $X$ is left-Bol, then for every $x,y\in X$ we have the identity
$$(x{\cdot} x){\cdot} y=(x{\cdot}(e{\cdot}x)){\cdot}y=x{\cdot}(e{\cdot}(x{\cdot}y))=x{\cdot}(x{\cdot} y),$$
witnessing that $X$ is left-alternative.
\end{proof}

\begin{proposition}\label{p:right-Bol=>right-alt+flex} Every right-Bol unital magma is right-alternative.
\end{proposition}

\begin{proof} Let $(X,\cdot)$ be a right-Bol unital magma. Then $X$ endowed with the mirror binary operation $*:(x,y)\mapsto x*y\defeq y\cdot x$, is a left-Bol unital magma. By Proposition~\ref{p:left-Bol=>left-alt+flex}, the magma $(X,*)$ is left-alternative and hence its mirrow magma $(X,\cdot)$ is right-alternative.
\end{proof}

Propositions~\ref{p:left-Bol=>left-alt+flex} and \ref{p:right-Bol=>right-alt+flex} imply

\begin{corollary}\label{c:Bol=>alt} Every Bol unital magma is alternative.
\end{corollary}

The following lemma will help us to prove that every left-Bol unital magma is $\IN$-associative.

\begin{lemma}\label{l:power-as} Let $X$ be a left-Bol magma and $(x^n)_{n\in\IN}$ be a sequence such that $x^{n+1}=x^n\cdot x^1$ for every $n\in\IN$. If $x^3=x^1{\cdot}x^2$, $x^4=x^2{\cdot}x^2$, and $x^5=x^2{\cdot}x^3$, then $x^{n+m}=x^n{\cdot}x^m$ for every $n,m\in\IN$.
\end{lemma}

\begin{proof} The equality $x^{n+m}=x^n{\cdot}x^m$ will be proved by induction on $n+m$. Let $x\defeq x^1$. Our assumptions guarantee that $x^2=x{\cdot}x$ and  $x^3=x^2{\cdot}x=x{\cdot}x^2$. Therefore, the equality $x^{n+m}=x^n{\cdot}x^m$ holds if $n+m\le 3$. 
%

Assume that for some $k\ge 4$ and all numbers  $n,m\in\IN$ with $n+m<k$ the equality $x^n{\cdot} x^m=x^{n+m}$ holds. Choose any numbers $n,m\in\IN$ such that $n+m=k$. 

\begin{claim}\label{cl:xk=xxk-1} $x^k=x\cdot x^{k-1}$.
\end{claim}

\begin{proof} The left-Bol identity and the inductive assumption imply $$x^k\defeq x^{k-1}\cdot x=(x\cdot (x^{k-3}\cdot x))\cdot x=x\cdot (x^{k-3}\cdot (x\cdot x))=x\cdot (x^{k-3}\cdot x^2)=x\cdot x^{k-1}.$$ 
\end{proof}

If $n=1$, then $x^n\cdot x^m=x\cdot x^{k-1}=x^k=x^{n+m}$, by Claim~\ref{cl:xk=xxk-1}. The assumptions $x^2{\cdot}x^2=x^4$ and $x^2{\cdot}x^3=x^5$ imply that the equaliy $x^n\cdot x^m=x^{n+m}$ holds if $n=2$ and $m\in\{2,3\}$.

If $n=2$ and $m\ge 4$, then the left-Bol identity and the inductive assumption ensure that
$$x^n{\cdot}x^m=x^2{\cdot}(x^{m-3}{\cdot}(x^{2}{\cdot}x))=(x^2{\cdot}(x^{m-3}{\cdot}x^2)){\cdot}x=(x^2{\cdot}x^{m-1}){\cdot}x=x^{m+1}\cdot x=x^{m+2}=x^{n+m}.$$

If $n\ge 3$, then $$x^n\cdot x^m=(x\cdot (x^{n-2}\cdot x))\cdot x^m=x\cdot (x^{n-2}\cdot(x\cdot x^m))=x\cdot x^{n-1+m}=x^{n+m},$$
by the left-Bol identity, the inductive assumption, and Claim~\ref{cl:xk=xxk-1}.
\end{proof}

\begin{theorem}\label{t:left-Bol=>power-associative} Every left-Bol unital magma $X$ is $\IN$-associative.
\end{theorem}

\begin{proof} Let $e$ be the neutral element of the unital magma $X$. Given any element $x\in X$, consider the sequence $(x^n)_{n\in\w}$ defined by the recursive formula $x^0\defeq e$ and $x^{n+1}\defeq x^n{\cdot} x$ for every $n\in\w$. 

We claim that $x^3=x^1\cdot x^2$, $x^4=x^2\cdot x^2$ and $x^5=x^2\cdot x^3$.

The left-Bol identity ensures that 
$$
\begin{aligned}
x^1{\cdot}x^2&=x{\cdot}(x{\cdot}x)=x{\cdot}(e{\cdot}(x{\cdot}x))= (x{\cdot}(e{\cdot}x)){\cdot}x=(x{\cdot}x){\cdot}x=x^3,\\
x^2{\cdot}x^2&=(x{\cdot}(e{\cdot}x)){\cdot}x^2=x{\cdot}(e{\cdot}(x{\cdot}x^2))=x{\cdot}(x{\cdot}(x{\cdot}x))=(x{\cdot}(x{\cdot}x)){\cdot} x=x^3{\cdot}x=x^4,\\
x^2{\cdot} x^3&=x^2{\cdot}(e{\cdot}(x^2{\cdot}x))=(x^2{\cdot}(e{\cdot}x^2)){\cdot}x=(x^2{\cdot}x^2){\cdot}x=x^4{\cdot}x=x^5.
\end{aligned}
$$
Applying Lemma~\ref{l:power-as}, we conclude that $x^n\cdot x^m=x^{n+m}$ for all $n,m\in\IN$, which means that the magma $X$ is $\IN$-associative.
\end{proof}

\begin{corollary}\label{c:right-Bol=>power-associative} Every right-Bol unital magma $X$ is $\IN$-associative.
\end{corollary}

\begin{proof} Let $(X,\cdot)$ be a right-Bol unital magma. Then the magma $(X,*)$ endowed with the mirror operation $x*y=y{\cdot}x$ is a left-Bol unital magma. By Theorem~\ref{t:left-Bol=>power-associative}, the left-Bol magma $(X,*)$ is $\IN$-associative and so is its mirror magma $(X,\cdot)$.
\end{proof} 

\begin{proposition}\label{p:left-Bol=>power-inversive} Let $X$ be a left-Bol unital magma and $(x^n)_{n\in\IN}$ be a sequence such that $x^{n+m}=x^n\cdot x^m$ for all $n,m\in\IN$. Then $x^n\cdot (x^m\cdot y)=x^{n+m}\cdot y$ for all $y\in X$ and $n,m\in\IN$.
\end{proposition}

\begin{proof} The equality $x^n\cdot (x^m\cdot y)=x^{n+m}{\cdot}y$ will be proved by induction on $n+m$. Let $x^0$ be the neutral element of the unital magma $X$. Then the equality $x^nx^m=x^{n+m}$ holds for all $n,m\in\w$. 

Assume that for some $k\in\w$ and all $n,m\in\w$ with $n+m<k$ the equality $x^n{\cdot}(x^m{\cdot}y)=x^{n+m}{\cdot}y$ holds.
Take any numbers $n,m\in\w$ with $n+m=k$. 
If $n=0$, then 
$$x^n{\cdot}(x^m{\cdot}y)=x^0{\cdot}(x^m{\cdot}y)=x^m{\cdot}y=x^{0+m}{\cdot}y=x^{n+m}{\cdot}y$$ and we are done.
If $m=0$, then 
$$x^n{\cdot}(x^m{\cdot}y)=x^n{\cdot}(x^0{\cdot}y)=x^n{\cdot}y=x^{n+m}{\cdot}y$$ and we are done.

So, assume that $n,m>0$, which implies $\max\{n,m\}<m+n=k$.
 If $n\le m$, the  by the inductive assumption and the left-Bol identity,
$$x^n{\cdot}(x^m{\cdot}y)=x^n{\cdot}((x^{m-n}{\cdot}x^n){\cdot}y)=x^n{\cdot}(x^{m-n}{\cdot}(x^n{\cdot}y))=(x^n{\cdot}(x^{m-n}{\cdot}x^n)){\cdot}y=x^{n+m}{\cdot}y.$$
If $m< n$, then $n\ge 2$ and by the left-Bol identity and the inductive assumption,
$$
\begin{aligned}
x^n{\cdot}(x^m{\cdot}y)&=(x{\cdot}(x^{n-2}{\cdot}x))(x^m{\cdot}y)=x{\cdot}(x^{n-2}{\cdot}(x{\cdot}(x^m{\cdot} y)))\\
&=x{\cdot}(x^{n+m-1}{\cdot}y)=x{\cdot}(x^{n+m-2}{\cdot}(x{\cdot}y))=(x{\cdot}(x^{n+m-2}{\cdot}x)){\cdot}y=x^{n+m}{\cdot}y.
\end{aligned}
$$
\end{proof}

By the mirror argument, Proposition~\ref{p:left-Bol=>power-inversive} implies its ``right'' version.

\begin{proposition} Let $X$ be a right-Bol unital magma and $(x^n)_{n\in\IN}$ be a sequence such that $x^{n+m}=x^n{\cdot} x^m$ for all $n,m\in\IN$. Then $(y{\cdot} x^n){\cdot} x^m=y{\cdot} x^{n+m}$ for every $y\in X$ and $n,m\in\IN$.
\end{proposition}

\section{Bol loops}

In this section we establish some properties of (left or right) Bol loops.

\begin{proposition}\label{p:left-Bol=>left-inversive} Every left-Bol loop $X$ is left-inversive.
\end{proposition}

\begin{proof} Given any element $a\in X$, we need to find an element $b\in X$ such that $b{\cdot} (a{\cdot} x)=x$ for all $x\in X$. Since $X$ is a loop, there exists a unique element $b\in X$ such that $b{\cdot} a=e$, where $e$ is the neutral element of the loop $X$. Since the magma $X$ is left-Bol, for every $x\in X$, we have the equality 
$$a{\cdot}(b{\cdot}(a{\cdot}x))=(a{\cdot}(b{\cdot}a)){\cdot}x=(a{\cdot} e){\cdot} x=a{\cdot} x,$$which implies
$b{\cdot} (a{\cdot} x)=x$, by the right-cancellativity of the quasigroup $X$.
\end{proof}

By analogy we can prove 

\begin{proposition}\label{p:right-Bol=>right-inversive} Every right-Bol loop $X$ is right-inversive.
\end{proposition}

Propositions~\ref{p:left-Bol=>left-inversive} and \ref{p:left-Bol=>left-inversive} imply

\begin{corollary}\label{c:Bol=>inversive} Every Bol loop is inversive.
\end{corollary}

The main result of this section is the following characterization of Bol loops.

\begin{theorem}\label{t:Bol<=>} For a loop $X$,  the following conditions are equivalent:
\begin{enumerate}
\item $X$ is Bol;
\item $X$ is inversive and left-Bol;
\item $X$ is inversive and right-Bol;
\item $X$ is right-inversive and left-Bol;
\item $X$ is left-inversive and right-Bol;
\item $X$ is flexible and left-Bol;
\item $X$ is flexible and right-Bol.
\end{enumerate}
\end{theorem}

\begin{proof} Let $e$ be the neutral element of the loop $X$.
First we prove the implications\\ $(1)\Ra(2)\Ra(4)\Ra(6)\Ra(2)\Ra(1)$.
\smallskip

The implication $(1)\Ra(2)$ follows from Corollary~\ref{c:Bol=>inversive}, and $(2)\Ra(4)$ is trivial.
\smallskip

$(4)\Ra(6)$. Assume that the loop $X$ is right-inversive and left-Bol. Since $X$ is a right-inversive loop, for every $x\in X$  there exists a unique element $x^{-1}\in X$ such that $(y{\cdot}x){\cdot}x^{-1}=y$ for all $y\in X$. Then for every $y\in X$, the following identities hold (the last equality follows from the left-Bol property):
$$((x{\cdot}y){\cdot}x){\cdot}x^{-1}=x{\cdot}y=x{\cdot}(y{\cdot}e)=x{\cdot}(y{\cdot}(x{\cdot}x^{-1}))=(x{\cdot}(y{\cdot}x)){\cdot}x^{-1}$$
and hence $(x{\cdot}y){\cdot}x=x{\cdot}(y{\cdot}x)$, by the cancellativity of the binary operation in the loop $X$. The identity $(x{\cdot}y){\cdot}x=x{\cdot}(y{\cdot}x)$ witnesses that the loop $X$ is flexible.
\smallskip

$(6)\Ra(2)$ Assume that $X$ is flexible and left-Bol. By Proposition~\ref{p:left-Bol=>left-inversive}, $X$ is left-inversive. By Proposition~\ref{p:invertible<=lr-inv}, the left-inversive loop $X$ is invertible. Therefore, for every $x\in X$, there exists a unique element $x^{-1}\in X$ such that $x{\cdot} x^{-1}=x^{-1}{\cdot} x=e$. To prove that the loop $X$ is right-inversive, it suffices to prove that $(z{\cdot}x){\cdot}x^{-1}=z$ for every $x,z\in X$. Take any elements $x,z\in X$. Since $X$ is a quasigroup, there exists a unique element $y\in X$ such that $z=x{\cdot}y$. Since $X$ is a flexible left-Bol loop, 
$$z=x{\cdot}y=x{\cdot}(y{\cdot}e)=x{\cdot}(y{\cdot}(x{\cdot}x^{-1}))=(x{\cdot}(y{\cdot}x)){\cdot}x^{-1}=((x{\cdot}y){\cdot}x){\cdot}x^{-1}=(z{\cdot}x){\cdot}x^{-1},$$ witnessing that the loop $X$ is right-inversive. By Proposition~\ref{p:inversive<=>lr-inversive}, the loop $X$ is inversive.
\smallskip

$(2)\Ra(1)$ Assume that the loop $X$ is inversive and left-Bol. By Proposition~\ref{p:inversive=>anti}, the inversion $X\to X$, $x\mapsto x^{-1}$, is an involutive anti-isomorphism of the loop $X$. Then for every $x,y,z\in X$, the left-Bol property of $X$ ensures that
$$
\begin{aligned}
((z{\cdot}x){\cdot}y){\cdot}x&=\big(\big(((z{\cdot}x){\cdot}y){\cdot}x\big)^{-1}\big)^{-1}=
\big(x^{-1}{\cdot}(y^{-1}{\cdot}(x^{-1}{\cdot}z^{-1}))\big)^{-1}\\
&=
\big((x^{-1}{\cdot}(y^{-1}{\cdot}x^{-1})){\cdot}z^{-1}\big)^{-1}=z{\cdot}((x{\cdot}y){\cdot}x),
\end{aligned}
$$
witnessing that $X$ is right-Bol and hence Bol.
\smallskip

Applying the (already proved) implications $(1)\Ra(2)\Ra(4)\Ra(6)\Ra(2)\Ra(1)$, to the loop $(X,*)$ endowed with the binary operation $x*y\defeq y\cdot x$, we can see that  $(1)\Ra(3)\Ra(5)\Ra(7)\Ra(3)\Ra(1)$.
\end{proof}

Now we are going to prove that left-Bol loops are power-associative.

\begin{theorem}\label{t:Z-ass-y} For every element $x$ of a left-Bol loop $X$, there exists a sequence $(x^n)_{n\in\IZ}$ such that $x^1=x$ and $x^n\cdot(x^m\cdot y)=x^{n+m}\cdot y$ for all $y\in X$ and $n,m\in\IZ$.
\end{theorem}

\begin{proof} By Proposition~\ref{p:left-Bol=>left-inversive}, the left-Bol loop $X$ is left-inversive. By Proposition~\ref{p:invertible<=lr-inv}, $X$ is invertible and hence for every $x\in X$ there exists a unique element $x^{-1}\in X$ such that $x\cdot x^{-1}=e=x^{-1}\cdot x$, where $e$ is the neutral element of the loop $X$.
Moreover, since $X$ is left-inversive, $x^{-1}\cdot (x\cdot y)=x\cdot(x^{-1}\cdot y)$ for all $x,y\in X$. Fix any element $x\in X$.

\begin{claim}\label{cl:xnx-n} There exist sequences $(x^n)_{n\in\w}$ and $(x^{-n})_{n\in\w}$ in $X$ such that for all $n,m\in\IN$ and $y\in X$, the following conditions are satisfied:
\begin{enumerate}
\item $x^1=x$ and $x\cdot x^{-1}=x^0=e$;
\item $x^n\cdot x^m=x^{n+m}$ and $x^{-n}\cdot x^{-m}=x^{-(n+m)}$;
\item $x^n\cdot (x^m\cdot y)=x^{n+m}\cdot y$ and $x^{-n}\cdot (x^{-m}\cdot y)=x^{-(n+m)}\cdot y$;
\item $x^{n}\cdot x^{-n}=x^0=x^{-n}\cdot x^{n}$.
\end{enumerate}
\end{claim}

\begin{proof} Let $x^0\defeq e$ be the neutral element of the loop $X$. By Theorem~\ref{t:left-Bol=>power-associative}, the left-Bol loop is $\IN$-associative and hence for the elements $x$ and $x^{-1}$ there exist  sequences $(x^n)_{n\in\IN}$ and $(x^{-n})_{n\in\IN}$ satisfying the conditions (1) and (2) of the Claim. Proposition~\ref{p:left-Bol=>power-inversive} ensures that these two sequences satisfy the condition (3).
\smallskip

By induction, we shall prove that $x^n\cdot x^{-n}=x^0$ for all $n\in\w$. Since $x^0$ is the neutral element of the loop $X$, the equality $x^0\cdot x^0=x^0$ holds. Assume that for some number $n\in\w$ the equality $x^{n}\cdot x^{-n}=x^0$ holds. The conditions (2,3),  the left-inversivity of $X$, and the inductive assumption ensure that
$$x^{n+1}\cdot x^{-(n+1)}=(x^{n}\cdot x)\cdot x^{-(n+1)}=
x^n\cdot(x\cdot x^{-(n+1)})=x^n\cdot (x\cdot (x^{-1}\cdot x^{-n}))=x^n\cdot x^{-n}=x^0.$$
The equality $x^n\cdot x^{-n}=x^0$ implies $x^{-n}\cdot x^n=(x^n)^{-1}\cdot x^n=x^0$.
\end{proof}

Let $(x^n)_{n\in\w}$ and $(x^{-n})_{n\in\w}$ be two sequences satisfying the conditions (1)--(4) of Claim~\ref{cl:xnx-n}.

\begin{claim}\label{cl:xn-m} For all positive integer numbers $n\le m$ and all $y\in X$, the following equalities hold:
\begin{enumerate}
\item $x^{-n}\cdot (x^m\cdot y)=x^{-n+m}\cdot y$;
\item $x^n\cdot (x^{-m}\cdot y)=x^{n-m}\cdot y$;
\item $x^{-m}\cdot (x^n\cdot y)=x^{-m+n}\cdot y$;
\item $x^m\cdot (x^{-n}\cdot y)=x^{m-n}\cdot y$.
\end{enumerate}
\end{claim}

\begin{proof} Claim~\ref{cl:xnx-n}(3) and the left-inversivity of $X$ imply
$$
\begin{aligned}
&x^{-n}\cdot (x^m\cdot y)=x^{-n}\cdot(x^n\cdot (x^{m-n}\cdot y))=x^{m-n}\cdot y,\\
&x^n\cdot (x^{-m}\cdot y)=x^n\cdot (x^{-n}\cdot (x^{n-m}\cdot y))=x^{n-m}\cdot y,\\
&x^m\cdot (x^{-n}\cdot y)=(x^{m-n}\cdot x^n)\cdot (x^{-n}\cdot y)=x^{m-n}\cdot (x^n\cdot (x^{-n}\cdot y)))=x^{m-n}\cdot y,\\
&x^{-m}\cdot (x^n\cdot y)=(x^{-m+n}\cdot x^{-n})\cdot (x^n\cdot y)=x^{m-n}\cdot(x^{-n}\cdot (x^n\cdot y))=x^{m-n}\cdot y.
\end{aligned}
$$
\end{proof}
Claims~\ref{cl:xnx-n} and \ref{cl:xn-m} imply that $x^n\cdot (x^m\cdot y)=x^{n+m}\cdot y$ for all $n,m\in\IZ$ and all $y\in Y$.
\end{proof}

Theorem~\ref{t:Z-ass-y} implies the following corollary improving Theorem~\ref{t:left-Bol=>power-associative} in case of left-Bol loops.

\begin{corollary}\label{c:left-Bol=>Z-ass} Every left-Bol loop is $\IZ$-associative.
\end{corollary}

Next, we prove that left-Bol loops satisfy a weak version of the Lagrange Theorem. For a subset $C$ of a magma $(X,\cdot)$ and a point $x\in X$ we denote by $C{\cdot}x\defeq\{c\cdot x:c\in C\}$ the right shift of $C$ by the element $x$.

\begin{proposition}\label{p:left-Bol=>Lagrange} Let $C$ be a cyclic subgroup of a left-Bol loop $X$.
\begin{enumerate}
\item For every $x\in X$ and $z\in C{\cdot}x$ we have $C{\cdot}z=C{\cdot}x$. 
\item For every $x,y\in X$, the sets $C{\cdot}x$ and $C{\cdot}y$ either coincide or are disjoint.
\item $|X|=|C|\cdot\kappa$ for some cardinal $\kappa$.
\end{enumerate}
\end{proposition}

\begin{proof} Let $c$ be a generator of the cyclic group $C$. By Theorem~\ref{t:Z-ass-y}, there exists a sequence $(c^n)_{n\in\IZ}$ such that $c^1=c$ and $c^n\cdot(c^m\cdot x)=c^{n+m}\cdot x$ for all $x\in X$. Then $\{c^n\}_{n\in\IZ}=C$.
\smallskip

1. Let $x,z\in X$ be points such that $z\in C{\cdot}x$. Then $z=c^n{\cdot}x$ for some $n\in\IZ$. By the left-inversivity of the left-Bol loop $X$, $z=c^n\cdot x$ implies $x=c^{-n}\cdot z$. Then for every $m\in\IZ$ we have $c^m\cdot z=c^m\cdot (c^n\cdot x)=c^{m+n}\cdot x\in C\cdot x$ and hence $C\cdot z\subseteq C\cdot x$. By analogy we can show that $x=c^{-n}\cdot z$ implies $C\cdot x\subseteq C\cdot z$. Therefore, $C\cdot x=C\cdot z$.
\smallskip

2. Take any points $x,y\in X$. Assuming that the sets $C{\cdot}x$ and $C{\cdot}y$ have a common element $z$, apply the preceding item and conclude that $C{\cdot}x=C{\cdot}z=C{\cdot}y$.
\smallskip

3. By the preceding item, the family $\mathcal F\defeq\{C{\cdot}x\}_{x\in X}$ consists of pairwise disjoint sets. By the Axiom of Choice, there exists a set $S\subseteq X$ such that for every $x\in X$, the intersection $S\cap (C{\cdot}x)$ is a singleton $\{x\}$. The first item ensures that $C{\cdot}x=C{\cdot}s$.  Therefore, $X=\bigcup_{x\in X}C{\cdot}x=\bigcup_{x\in S}C{\cdot}x$ and $\mathcal F$ is a partition of $X$ into pairwise disjoint sets $C{\cdot}x$, $x\in S$. Since the loop $X$ is cancellative, for every $s\in S$, the right shift $R_s:X\to X$, $R_s:x\mapsto x{\cdot}s$, is injective and hence $|C{\cdot}s|=|R_s[C]|=|C|$. Therefore, $\mathcal F$ is a partition of $X$ into pairwise disjoint sets of cardinality $|C|$, which implies $|X|=|C|\cdot|S|$.
\end{proof}

Proposition~\ref{p:left-Bol=>Lagrange} implies

\begin{corollary}\label{c:Bol=>Lagrange} Let $C$ be a cyclic subgroup of a loop $X$. If the loop $X$ is left-Bol or right-Bol, then the cardinal $|C|$ divides the cardinal $|X|$.
\end{corollary}

\begin{corollary}\label{c:left-Bol-prime=>cyclic} If a left-Bol loop $X$ has prime order $|X|$, then it is a cyclic group.
\end{corollary}

\begin{proof} Let $e$ be the neutral element of the loop $X$. Since $|X|$ is a prime number, $X\ne \{e\}$ and hence $X$ contains some element $x\ne e$.  By Corollary~\ref{c:left-Bol=>Z-ass}, the left-Bol loop $X$ is $\IZ$-associative. Then the element $x$ generates a cyclic subgroup $C$ in the loop $X$. By Corollary~\ref{c:Bol=>Lagrange}, the cardinal $|C|>1$ divides the prime number $|X|$ and hence $|C|=|X|$, which implies that $X=C$ is a cyclic group.
\end{proof}

Corollaries~\ref{c:Bol=>Lagrange} and \ref{c:left-Bol-prime=>cyclic} do not generalize to inversive loops.

\begin{example}\label{ex:Ali-Slaney} The $7$-element loop with multiplication table
$$
\begin{array}{c|cccccccc}
+&1&2&3&4&5&6&7\\ 
\hline
1&1&2&3&4&5&6&7\\
2&2&3&1&6&7&5&4\\
3&3&1&2&7&6&4&5\\
4&4&7&6&5&1&2&3\\
5&5&6&7&1&4&3&2\\
6&6&4&5&3&2&7&1\\
7&7&5&4&2&3&1&6
\end{array}
$$
is inversive but is neither commutative nor associative. 
This loop contains $3$-element subgroups $\{1,2,3\}$, $\{1,4,5\}$ and $\{1,6,7\}$ (whose order $3$ does not divide the order $7$ of the loop). This example shows that the Lagrange Theorem~\ref{t:Lagrange} fails for finite inversive loops. Nonetheless, by Corollary~\ref{c:Ali-Slaney}, a (left or right) inversive loop of odd order cannot contain an subloop of even order.
\end{example}

\section{Moufang loops}

\begin{definition}\label{d:Moufang-magma} A magma $(X,\cdot)$ is \defterm{Moufang} if it satisfies two Moufang identities:
\begin{enumerate}
\item[\textup{(M1)}] $\forall x,y,z\in X\;\;(x{\cdot} y){\cdot}(z{\cdot}x)=(x{\cdot}(y{\cdot}z)){\cdot}x;$
\item[\textup{(M2)}] $\forall x,y,z\in X\;\;(x{\cdot} y){\cdot}(z{\cdot}x)=x{\cdot}((y{\cdot}z){\cdot}x).$
\end{enumerate}
\end{definition}

Let us show that for unital magmas, the Moufang identities (M1) and (M2) in Definition~\ref{d:Moufang-magma} are equivalent.

\begin{proposition}\label{p:M12=>flexible} If a unital magma $X$ satisfies one of the Moufang identities \textup{(M1)} or \textup{(M2)}, then $X$ is flexible.
\end{proposition}

\begin{proof}  Let $e$ be the identity of the unital magma $X$. If $X$ satisfies the Moufang identity (M1), then for every $x,z\in X$, we have
$x{\cdot}(y{\cdot}x)=(x{\cdot}e){\cdot}(y{\cdot}x)=(x{\cdot} (e{\cdot} y)){\cdot} x=(x{\cdot}y){\cdot} x$, witnessing that the magma $X$ is flexible. 

If $X$ satisfies the Moufang identity (M2), then for every $x,y\in X$, we have
$(x{\cdot}y){\cdot}x=(x{\cdot}y){\cdot}(e{\cdot}x)=x{\cdot}( (y{\cdot} e){\cdot} x)=x{\cdot}(y{\cdot} x)$, witnessing that the magma $X$ is flexible. 
\end{proof}

\begin{corollary}\label{c:uM=>M1=M2} For every unital magma $X$, the Moufang identities \textup{(M1)} and \textup{(M2)} are equivalent.
\end{corollary}

\begin{proof} If a unital magma $X$ satisfies (M1) or (M2), then it is flexible, by Proposition~\ref{p:M12=>flexible}. Then for every $x,y,z\in X$ we have the equality $$(x{\cdot}(y{\cdot}z)){\cdot} x=x{\cdot}((y{\cdot}z){\cdot}x),$$which implies that the the Moufang identities (M1) and (M2) are equivalent.
\end{proof}

Since every loop is a unital magma, Corollary~\ref{c:uM=>M1=M2} implies another corollary.

\begin{corollary}\label{c:loop=>M1=M2} For every loop $X$, the Moufang identities \textup{(M1)} and \textup{(M2)} are equivalent.
\end{corollary}

\begin{proposition}\label{p:Moufang=>inversive} Every Moufang loop is inversive.
\end{proposition}

\begin{proof} Let $e$ be the identity of the loop $X$. Since loop $X$ is a quasigroup, for every element $x\in X$, there exists a unique element $x^{-1}\in X$ such that $x^{-1}{\cdot} x=e$. By the Moufang identity (M1), for every $x,y\in X$, we have
$$y{\cdot}x^{-1}=e{\cdot}(y{\cdot}x^{-1})=(x^{-1}{\cdot} x){\cdot}(y{\cdot} x^{-1})=(x^{-1}{\cdot}(x{\cdot}y)){\cdot}x^{-1},$$which implies $y=x^{-1}{\cdot}(x{\cdot}y)$, by the cancellativity of the binary operation in loops. The latter identity witnesses that $X$ is left-inversive. 

By Proposition~\ref{p:invertible<=lr-inv}, the left-inversive loop $X$ is invertible and hence $x{\cdot} x^{-1}=e=x^{-1}{\cdot} x$ for every $x\in X$. 

By the Moufang identity (M2), for every $x,y\in X$ we have
$$x^{-1}{\cdot}y=(x^{-1}{\cdot}y){\cdot} e=(x^{-1}{\cdot}y){\cdot}(x{\cdot}x^{-1})=x^{-1}((y{\cdot}x){\cdot} x^{-1}),$$
which implies $y=(y{\cdot}x){\cdot}x^{-1}$, by the cancellativity of the binary opertaion in loops. The latter identity witnesses that the loop $X$ is right-inversive.
\end{proof}


\begin{theorem}[Bol, 1937; Bruck, 1946]\label{t:Bol<=>Moufang} For a loop $X$, the following conditions are equivalent:
\begin{enumerate}
\item $X$ is Moufang;
\item $\forall x,y,z\in X\;\;x{\cdot}(y{\cdot}(x{\cdot}z))=((x{\cdot}y){\cdot}x){\cdot}z$;
\item $\forall x,y,z\in X\;\;((z{\cdot}x){\cdot}y){\cdot}x=z{\cdot}(x{\cdot}(y{\cdot}x))$;
\item $\forall x,y\in X\;\;L_xL_yL_x=L_{(x{\cdot}y){\cdot}x}$;
\item $\forall x,y\in X\;\;R_xR_yR_x=R_{x{\cdot}(y{\cdot}x)}$;
\item $X$ is flexible and left-Bol;
\item $X$ is flexible and right-Bol;
\item $X$ is Bol.
\end{enumerate}
\end{theorem} 

\begin{proof} Let $e$ be the neutral element of the loop $X$.
First, we prove the implications\\ $(1)\Ra(2)\Ra(4)\Leftrightarrow(6)\Ra(8)\Ra(1)$.
\smallskip

 $(1)\Ra(2)$ Assume that the loop $X$ is Moufang. By Proposition~\ref{p:Moufang=>inversive}, $X$ is inversive. By Proposition~\ref{p:inversive=>anti}, the inversion $X\to X$, $x\mapsto x^{-1}$, is an involutive anti-isomorphism of the loop $X$. 

Let $a\defeq (x{\cdot}z)^{-1}=z^{-1}{\cdot}x^{-1}$ and $b\defeq y{\cdot}(x{\cdot}z)$. Observe that $a{\cdot}x=(z^{-1}{\cdot}x^{-1}){\cdot}x=z^{-1}$, by the right-inversivity of $X$. The inversivity of $X$ implies $y=(y{\cdot}(x{\cdot}z)){\cdot}(x{\cdot}z)^{-1}=b{\cdot}a$. By the inversivity and the Moufang identity $(M2)$, $$x{\cdot}(y{\cdot}(x{\cdot}z))=x{\cdot}b=((x{\cdot}b){\cdot}(a{\cdot}x))(a{\cdot}x)^{-1}=((x{\cdot}(b{\cdot}a)){\cdot}x)(z^{-1})^{-1}=((x{\cdot}y){\cdot}x){\cdot}z.$$
\smallskip

$(2)\Ra(4)$ Applying the condition (2)  to elements $x,y,e$, we obtain the identity
$$x{\cdot}(y{\cdot}x)=x{\cdot}(y{\cdot}(x{\cdot}e))=((x{\cdot}y){\cdot}x){\cdot}e=(x{\cdot}y){\cdot}x,$$
witnessing that the loop $X$ is flexible. Then for every $x,y,z\in X$, the condition (2) implies
$$x{\cdot}(y{\cdot}(x{\cdot}z))=((x{\cdot}y){\cdot}x){\cdot}z=(x{\cdot}(y{\cdot}x)){\cdot}z,$$
witnessing that the loop $X$ is left-Bol.
\smallskip

The equivalence $(4)\Leftrightarrow(6)$ is obvious, and the implication $(6)\Ra(8)$ has been proved in Theorem~\ref{t:Bol<=>}.
\smallskip

$(8)\Ra(1)$ Assume that the loop $X$ is Bol. By Theorem~\ref{t:Bol<=>}, $X$ is flexible, inversive, and hence invertible. By Proposition~\ref{p:inversive=>anti}, the inversion $X\to X$, $x\mapsto x^{-1}$, is an involutive anti-isomorphism of the loop $X$. We claim that the Bol loop $X$ satisfies the Moufang identity (M1). Given any elements $x,a,b\in X$, we need to check that $(x{\cdot}b){\cdot}(a{\cdot}x)=(x{\cdot}(b{\cdot}a)){\cdot}x$. Consider the element $y\defeq b{\cdot}a$. Since $X$ is a loop, there exist a unique element $z\in X$ such that $a=z^{-1}{\cdot}x^{-1}=(x{\cdot}z)^{-1}$. The inversivity of $X$ ensures that $a{\cdot}x=(z^{-1}{\cdot}x^{-1}){\cdot}x=z^{-1}$ and $b=(b{\cdot}a){\cdot}a^{-1}=y{\cdot}a^{-1}=y{\cdot}(x{\cdot}z)$. Since the loop $X$ is flexible, inversive, and Bol,
$$
\begin{aligned}
((x{\cdot}b){\cdot}(a{\cdot}x)){\cdot}z&=((x{\cdot}b){\cdot}(a{\cdot}x)){\cdot}(a{\cdot}x)^{-1}=x{\cdot}b\\
&=x{\cdot}(y{\cdot}(x{\cdot}z))=(x{\cdot}(y{\cdot}x)){\cdot}z=((x{\cdot}y){\cdot}x){\cdot}z=((x{\cdot}(b{\cdot}a)){\cdot}x){\cdot}z,
\end{aligned}
$$
and hence $(x{\cdot}b){\cdot}(a{\cdot}x)=(x{\cdot}(b{\cdot}a)){\cdot}x$, by the cancellativity of the binary operation in the loop $X$. Therefore, the Bol loop $X$ satisfies the Moufang identity (M1) and is Moufang, by Corollary~\ref{c:loop=>M1=M2}.
\smallskip

Applying the (already proved) implications $(1)\Ra(2)\Ra(4)\Ra(6)\Ra(8)\Ra(1)$, to the loop $(X,*)$ endowed with the binary operation $x*y\defeq y\cdot x$, we can see that  $(1)\Ra(3)\Ra(5)\Ra(7)\Ra(8)\Ra(1)$.
\end{proof}

\begin{proposition}\label{p:Boolean+Moufang=>commutative-group} Every Boolean Moufang loop is a commutative group.
\end{proposition}

\begin{proof} Let $X$ be a Boolean Moufang loop. By Theorem~\ref{t:Bol<=>Moufang}, the Moufang loop $X$ is Bol, by Corollary~\ref{c:Bol=>inversive}, the Bol loop $X$ is inversive, and by Proposition~\ref{p:inversiveBoolean=>commutative}, the inversive Boolean loop $X$ is commutative. To see that $X$ is associative, take any elements $x,y,z\in X$ and applying the Bol identity, observe that
$$((x{\cdot}y){\cdot}z){\cdot}y=x{\cdot}((y{\cdot}z){\cdot}y)=x{\cdot}((z{\cdot}y){\cdot}y^{-1}))=x{\cdot}z.$$
On the other hand, the commutativity and the Bol identity imply
$$(x{\cdot}(y{\cdot}z)){\cdot}y=((z{\cdot}y){\cdot}x){\cdot}y=z{\cdot}((y{\cdot}x){\cdot}y)=z{\cdot}((x{\cdot}y){\cdot}y^{-1})=z{\cdot}x.$$
Therefore, $$((x{\cdot}y){\cdot}z){\cdot}y=x{\cdot}z=z{\cdot}x=(x{\cdot}(y{\cdot}z)){\cdot}y$$ and hence $(x{\cdot}y){\cdot}z=x{\cdot}(y{\cdot}z)$, by the cancellativity of the binary operation in the loop $X$.
\end{proof}

\section{Autotopisms and pseudo-automorphisms of loops}

In this section we establish some basic facts on autotopisms and pseudo-automorphisms of magmas and loops. Those facts will be essentially used in the Dr\'apal's proof of the Moufang Theorem, presented in the next section. 

\begin{definition} Let $(M,\cdot)$ be a magma. A triple $(\alpha,\beta,\gamma)$ of permutations of $M$ is called an \defterm{autotopism} of the magma $M$ if 
$\alpha(x)\cdot \beta(y) =
\gamma(x\cdot y)$ for all $x, y\in M$.
\end{definition}

\begin{remark} A permutation $\alpha:M\to M$ of a magma $M$ is an \defterm{automorphism} of $M$ if and only if the triple $(\alpha,\alpha,\alpha)$ is an autotopism of $M$.
\end{remark}

Like automorphisms, the autotopisms of a magma form a group.

\begin{lemma}\label{l:autotopies-group}  For every autotopisms $(\alpha,\beta,\gamma)$ and $(\alpha',\beta',\gamma')$ of a magma, the triples $$(\alpha\alpha',\beta\beta',\gamma\gamma')\quad\mbox{and}\quad(\alpha^{-1},\beta^{-1},\gamma^{-1})$$ are autotopisms the magma.
\end{lemma}

\begin{proof} Observe that for every elements $x,y$ of the magma $M$, we have
$$(\alpha\alpha'(x)\cdot\beta\beta'(y))=\gamma(\alpha'(x)\cdot\beta'(y))=\gamma\gamma'(x\cdot y),$$
which means that $(\alpha\alpha',\beta\beta',\gamma\gamma')$ is an autotopism of $M$.

Since $(\alpha,\beta,\gamma)$ is an autotopism of $M$,
$$\gamma(\alpha^{-1}(x)\cdot \beta^{-1}(y))=\alpha(\alpha^{-1}(x))\cdot\beta(\beta^{-1}(y))=x\cdot y=\gamma\gamma^{-1}(x\cdot y),$$ which implies $\alpha^{-1}(x)\cdot\beta^{-1}(y)=\gamma^{-1}(x\cdot y)$ and witnesses that $(\alpha^{-1},\beta^{-1},\gamma^{-1})$ is an autotopism of the magma $M$.
\end{proof}

Lemma~\ref{l:autotopies-group} implies 

\begin{corollary}\label{c:autotopisms-group} For every magma $M$, the set of autotopisms of $M$ is a group with respect to the operation of composition $\circ$, defined by $(\alpha,\beta,\gamma)\circ(\alpha',\beta',\gamma')\defeq(\alpha\alpha',\beta\beta',\gamma\gamma')$.
\end{corollary}

\begin{theorem}\label{t:autotopism3} Let $M$ be an invetible unital magma, and $i:M\to M$, $i:x\mapsto x^{-1}$, be the inversion of $M$. If $(\alpha,\beta,\gamma)$ is an autotopism of $M$ and 
\begin{enumerate}
\item $M$ is right-inversive, then $(\gamma,i\beta i,\alpha)$ is an autotopism of $M$;
\item $M$ is left-inversive, then $(i\alpha i,\gamma,\beta)$ is an autotopism of $M$;
\item $M$ is inversive, then $(\beta,i\gamma i,i\alpha i)$ is an autotopism of $M$.
\end{enumerate}
\end{theorem}

\begin{proof} We divide the proof of this theorem into three claims.

\begin{claim}\label{cl:autotop1} If $(\alpha,\beta,\gamma)$ is an autotopism of $M$ and $M$ is right-invarsive, then $(\gamma,i\beta i,\alpha)$ is an autotopism of $M$.
\end{claim}

\begin{proof} Given any elements $x,y\in M$, consider the elements $a=x\cdot y$ and $b=y^{-1}$. Since the magma $M$ is right-inversive, $a\cdot b=(x\cdot y)\cdot y^{-1}=x$. Since $(\alpha,\beta,\gamma)$ is an autotopism of $X$, we have $\alpha(a)\cdot \beta(b)=\gamma(a\cdot b)$ and by the right-inversivity of $X$, $$\alpha(x\cdot y)=\alpha(a)=\gamma(a\cdot b)\cdot(\beta(b))^{-1}=\gamma(x)\cdot i\beta i(y),$$
witnessing that $(\gamma,i\beta i,\alpha)$ is an autotopism of $M$.
\end{proof}

\begin{claim}\label{cl:autotop2} If $(\alpha,\beta,\gamma)$ is an autotopism of $M$ and $M$ is left-inversive, then $(i\alpha i,\gamma,\beta)$ is an autotopism of $M$.
\end{claim}

\begin{proof} Given any elements $x,y\in M$, consider the elements $c\defeq x^{-1}$ and $d\defeq x\cdot y$, and observe that $c\cdot a=x^{-1}\cdot(x\cdot y)=y$, by the left-inversivity of $M$. By the left-inversivity of $M$, the identity $\alpha(c)\cdot\beta(d)=\gamma(c\cdot d)$ implies
$$\beta(x\cdot y)=\beta(d)=\alpha(c)^{-1}\cdot \gamma(c\cdot d)=i\alpha i(x)\cdot \gamma(y),$$
witnessing that $(i\alpha i,\gamma,\beta)$ is an autotopism of $M$. 
\end{proof}

\begin{claim} If $(\alpha,\beta,\gamma)$ is an autotopism of $M$ and $M$ is inversive, then $(\beta,i\gamma i,i\alpha i)$ is an autotopism of $M$.
\end{claim}

\begin{proof} By Claim~\ref{cl:autotop2}, $(i\alpha i,\gamma,\beta)$ is an autotopism of $M$. Applying Claim~\ref{cl:autotop1} to the autotopism
$(i\alpha i,\gamma ,\beta)$, we conclude that 
  and $(\beta,i\gamma i,i\alpha i)$ is an autotopism of $M$.
\end{proof}
\end{proof}

Recall that for an element $a$ of a magma $(M,\cdot)$, the maps
$$L_a:M\to M,\quad L_a:x\mapsto a{\cdot}x,\quad\mbox{and}\quad
R_a:M\to M,\quad R_a:x\mapsto x{\cdot}a,
$$are the \defterm{left} and \defterm{right shifts} of $M$ by the 
element $a$.

\begin{definition} A self-map $t:M\to M$ of a magma $M$ is called a \defterm{translation} of $M$ if $t$ is a finite composition of left or right shifts. A translation $t:M\to M$ of a unital magma $M$ is called an \defterm{inner mapping} if $t(e)=e$.
\end{definition}

Translations of an inversive loop form a group, called the \defterm{multiplication group} of the loop.  Inner mappings form a normal subgroup of the multplication group.
 
\begin{proposition}\label{p:Moufang<=>} A quasigroup $M$ is Moufang if and only if for every $a\in M$, the triples $$(L_a,R_a,R_aL_a)\quad\mbox{and}\quad(L_a,R_a,L_aR_a)$$ are autotopisms of the magma $M$.
\end{proposition}

\begin{proof} Since $M$ is a quasigroup, the left and right shifts of $M$ are bijections. If $M$ is Moufang, then for every $a,x,y\in M$, we have the equalities
$$
\begin{aligned}
L_a(x)\cdot R_a(y)&=(a{\cdot}x)\cdot(y{\cdot}a)=(a\cdot(x\cdot y))\cdot a=R_aL_a(x\cdot y)\quad\mbox{and}\\
L_a(x)\cdot R_a(y)&=(a{\cdot}x)\cdot(y{\cdot}a)=a\cdot ((x\cdot y)\cdot a)=L_aR_a(x\cdot y),
\end{aligned}
$$
witnessing that the triples $(L_a,R_a,R_aL_a)$ and $(L_a,R_a,L_aR_a)$ are autotopisms of the magma $M$.
\smallskip

Now assume that for every $a\in M$, the triples   $(L_a,R_a,R_aL_a)$ and $(L_a,R_a,L_aR_a)$ are autotopisms of $M$.  Then for every $x,y\in M$ we have the equalities
$$
\begin{aligned}
(a{\cdot}x){\cdot}(y{\cdot}a)&=L_a(x)\cdot R_a(y)=R_aL_a(x{\cdot}y)=(a{\cdot}(x{\cdot}y)){\cdot}a\quad\mbox{and}\\
(a{\cdot}x){\cdot}(y{\cdot}a)&=L_a(x)\cdot R_a(y)=L_aR_a(x{\cdot}y)=a{\cdot}((x{\cdot}y){\cdot}a),
\end{aligned}
$$
witnessing that the quasigroup $M$ is Moufang.
\end{proof}

Let us recall that a magma $(M,\cdot)$ is called \defterm{unital} if it contains a (necessarily unique) element $e$ such that $x\cdot e=x=e\cdot x$ for all $x\in M$. The element $e$ is called the \defterm{neutral element} of the unital magma $M$.

\begin{lemma}\label{l:autotopism-e} Let $(\alpha,\beta,\gamma)$ be an autotopism of an unital magma $M$ with neutral element $e\in M$. If $\alpha(e)=e$, then $\beta=\gamma=R_c\alpha$, where $c=\beta(e)=\gamma(e)$.
\end{lemma}

\begin{proof} For every $x\in M$ we have $\gamma(x)=\gamma(e\cdot x)=\alpha(e)\cdot\beta(x)=e\cdot\beta(x)=\beta(x)$ and hence $\gamma=\beta$. Consider the element $c\defeq\beta(e)=\gamma(e)$ and observe that $$\gamma(x)=\gamma(x\cdot e)=\alpha(x)\cdot\beta(e)=\alpha(x)\cdot c=R_c\alpha(x)$$and hence $\beta=\gamma=R_c\alpha$.
\end{proof}

Lemma~\ref{l:autotopism-e} implies the following characterization.

\begin{proposition}\label{p:pseudo-automorphism} Let $\alpha:M\to M$ be a permutation of a loop $M$ such that $\alpha(e)=e$, where $e$ is the neutral element of $M$. The following conditions are equivalent:
\begin{enumerate}
\item for some permutations $\beta,\gamma$ of $M$, the triple $(\alpha,\beta,\gamma)$ is an autotopism of $M$;
\item for some permutation $\beta$ of $M$, the triple $(\alpha,\beta,\beta)$ is an autotopism of $M$;
\item for some element $c\in M$, the triple $(\alpha,R_c\alpha,R_c\alpha)$ is an autotopism of $M$.
\end{enumerate}
\end{proposition}

Proposition~\ref{p:pseudo-automorphism} motivates the following definition.

\begin{definition} Let $M$ be a unital magma and $e$ be the neutral element of $M$. A permutation $\alpha:M\to M$ is called 
\begin{itemize}
\item a \defterm{pseudo-automorphism} of $M$ if $\alpha(e)=e$ and for some permutations $\beta,\gamma$ of $M$, the triple $(\alpha,\beta,\gamma)$ is an autotopism of $M$;
\item a \defterm{pseudo-automorphism of $M$ with a companion $c\in M$} if $\alpha(e)=e$ and the triple $(\alpha,R_c\alpha,R_c\alpha)$ is an autotopism of $M$.
\end{itemize}
\end{definition}

Proposition~\ref{p:pseudo-automorphism} implies the following characterization of pseudo-automorphisms.

\begin{corollary} A permutation $\alpha$ of a loop $M$ is a pseudo-automorphism of $M$ if and only if $\alpha$ is a pseudo-autotopism of $M$ with a companion $c\in M$.
\end{corollary}

Corollary~\ref{c:autotopisms-group} implies the following corollary.

\begin{corollary} The set of pseudo-automorphisms of a loop is a subgroup of the permutation group of the loop.
\end{corollary}

\begin{remark} Pseudo-automorphisms were introduced by Bruck  \cite{Bruck1958}. In the textbook of Pflugfelder \cite{Pflugfelder1990}, pseudo-automorphisms are called \defterm{left pseudo-automorphisms}.
\end{remark}

\begin{definition} A function $\alpha:X\to X$ on a unital magma $M$ with neutral element $e$ is called a \defterm{semi-endomorphism} of $M$ if $\alpha(e)=e$ and $\alpha(x{\cdot}(y{\cdot}x))=\alpha(x){\cdot}(\alpha(y){\cdot}\alpha(x))$ for every $x,y\in M$. If $\alpha$ is bijective, then $\alpha$ is called a \defterm{semi-automorphism} of $M$.
\end{definition}

\begin{example} For every flexible inversive loop $M$, the function $i:M\to M,\quad i:x\mapsto x^{-1},$ 
is a semi-automorphism of $M$.
\end{example}

It is clear that semi-automorphism of a unital magma $M$ form a subgroup of the symmetric group $\mathcal S_M$ of the set $M$.

\begin{theorem}\label{t:pseudo-auto=>semi-auto} Every pseudo-automorphism $\alpha:M\to M$ of a Moufang loop $M$ is a semi-automorphism of $M$.
\end{theorem}

\begin{proof} By Proposition~\ref{p:pseudo-automorphism}, there exists $c\in M$ such that $(\alpha,R_c\alpha,R_c\alpha)$ is an autotopism of $M$. Then for every $x,y\in M$, we have
$$
\begin{aligned}
\alpha(x{\cdot}(y{\cdot}x)){\cdot}c&=R_c\alpha(x{\cdot}(y{\cdot}x))=\alpha(x)\cdot R_c\alpha(y{\cdot}x)=\alpha(x){\cdot}(\alpha(y){\cdot} R_c\alpha(x))\\
&=\alpha(x){\cdot}(\alpha(y){\cdot}(\alpha(x){\cdot}c))=(\alpha(x){\cdot}(\alpha(y){\cdot}\alpha(x))){\cdot}c,
\end{aligned}
$$by the Bol identity. The right-cancellativity of $M$ ensures that 
$\alpha(x{\cdot}(y{\cdot}x))=\alpha(x){\cdot}(\alpha(y){\cdot}\alpha(x))$, witnessing that $\alpha$ is a semi-automorphism of the loop $M$.
\end{proof}

\begin{theorem}\label{t:inner=>pseudo-semi-auto} Every inner mapping of a Moufang loop $M$ is a pseudo-automorphism and a semi-automorphism of the loop $M$.
\end{theorem}

\begin{proof} By Proposition~\ref{p:Moufang<=>}, for every $x\in M$, the triple $(L_x,R_x,R_xL_x)$ is an autotopism of the loop $M$. Appying Theorem~\ref{t:autotopism3}(3) to the autotopism $(L_x,R_x,R_xL_x)$, we conclude that the triple $(R_x,iR_xL_xi,iR_xi)$ is an autotopism of $M$. Since autotopisms form group, every translation $\alpha$ of $M$ can be completed to an autotopism $(\alpha,\beta,\gamma)$. If $\alpha$ is an inner mapping, then $\alpha(e)=e$ and hence $\alpha$ is a pseudo-automorphism and a semi-automorphism of $M$, according to Theorem~\ref{t:pseudo-auto=>semi-auto}.
\end{proof}

For an element $a$ of a diassociative loop $X$, the map
$$T_a:X\to X,\quad T_a:x\mapsto a^{-1}{\cdot}x{\cdot}a,$$
is called the \defterm{conjugation} of $X$ by the element $a$. 

\begin{theorem} For every element $a$ of a Moufang loop $X$, the conjugation $T_a:X\to X$, $T_a:x\mapsto a^{-1}{\cdot}x{\cdot}a$, is a pseudo-automorphism of $X$ with companion $c=a^{-3}$. Consequently,
$$a^{-1}{\cdot}(x{\cdot}y){\cdot}a^{-2}=R_cT_a(x{\cdot}y)=T_a(x)\cdot(R_cT_a(y))=(a^{-1}{\cdot}x{\cdot}a)\cdot (a^{-1}{\cdot}x{\cdot}a^{-2})$$for all $x,y\in X$.
\end{theorem}

\begin{proof} By Proposition~\ref{p:Moufang<=>}, the triples $(L_a^{-1},R_{a^{-1}},R_{a^{-1}}L_{a^{-1}})$ and $(L_a,R_a,R_aL_a)$ are autotopisms of the loop $X$. Applying Theorem~\ref{t:autotopism3}(3) to the autotopism $(L_a,R_a,R_aL_a)$, we conclude that the triple $(R_a,iR_aL_ai,iR_ai)$ is an autotopism of $M$. Then the composition 
$$(L_{a^{-1}}R_a,R_{a^{-1}}iR_aL_ai,R_{a^{-1}}L_{a^{-1}}iR_ai)$$is also an autotopism of the loop $X$. Since $T_a(e)=L_{a^{-1}}R_a(e)=e$, the translation $T_a=L_{a^{-1}}R_a$ is a pseudo-automorphism of $X$ with companion $c=R_{a^{-1}}iR_aL_ai(e)=a^{-3}$, by Lemma~\ref{l:autotopism-e}. This lemma also implies that 
$a^{-1}{\cdot}(x{\cdot}y){\cdot}a^{-2}=R_cT_a(x{\cdot}y)=T_a(x)\cdot(R_cT_a(y))=(a^{-1}{\cdot}x{\cdot}a){\cdot}(a^{-1}{\cdot}y{\cdot}a^{-2})$ for all $x,y\in X$.
\end{proof}

\begin{corollary}\label{c:conjugation3=>automorphism}  For every element $a$ of order $3$ in a Moufang loop $X$, the conjugation $T_a:X\to X$, $T_a:x\mapsto a^{-1}{\cdot}x{\cdot}a$, is an automorphism of the loop $X$.
\end{corollary}

Semi-automorphisms of Moufang loops share many common properties with automorphisms.

\begin{lemma}\label{l:pseudo=>inverse} Let $\alpha:M\to M$ be a semi-endomorphism of an invertible loop $M$.  For every element $x$ of $M$, the equality $\alpha(x)=x$ implies $\alpha(x^{-1})=x^{-1}$.
\end{lemma}

\begin{proof} Observe that
$x{\cdot}e=x=\alpha(x)=\alpha(x{\cdot}(x^{-1}{\cdot}x))=\alpha(x){\cdot}(\alpha(x^{-1}){\cdot}\alpha(x))=x{\cdot}(\alpha(x^{-1})\cdot x)$ 
and hence $x^{-1}{\cdot} x=e=\alpha(x^{-1}){\cdot}x$ and $x^{-1}=\alpha(x^{-1})$, by the cancellativity of the loop $M$.
\end{proof}

Now we present some important examples of pseudo-automorphisms of Moufang loops.
We shall use the following two identities holding in Moufang loops.

\begin{proposition}\label{p:Moufang-pauto} For any elements $x,y$ of a Moufang loop $M$, 
\begin{enumerate}
\item $L_{x{\cdot}y}^{-1}L_xL_y=R_xL_y^{-1}R_x^{-1}L_y$, and
\item $R_{y{\cdot}x}^{-1}R_xR_y=L_xR_y^{-1}L_x^{-1}R_y=L_xL_{x^{-1}y}L_y^{-1}$.
\end{enumerate}
\end{proposition}

\begin{proof} 1. For every $z\in M$, the Moufang identity $(xy)(zx)=x((yz)x)$ implies $L_{xy}R_x=L_xR_xL_y,$ which is equivalent to $L_{xy}^{-1}L_xL_y=R_xL_y^{-1}R_x^{-1}L_y$.
\smallskip

2. Consider the loop $(M,*)$ endowed with the binary operation $*$ defined by $x*y=y\cdot x$. Since the loop $(M,\cdot)$ is Moufang, for every $a,b,c\in M$,
$$(a*b)*(c*a)=(a\cdot c)\cdot(b\cdot a)=(a\cdot (c\cdot b))\cdot a=a*((b*c)*a),$$ witnessing the the loop $(M,*)$ is Moufang. 
Observe that for every $a\in M$ the left and right shifts $L_a,R_a$ of the loop $(M,\cdot)$ are respectively the right and left shifts of the loop $(M,*)$. Taking this correspondence cin mind, we can see that the identity $L_{x{*}y}^{-1}L_xL_y=R_xL_y^{-1}R_x^{-1}L_y$ in the loop $(M,*)$ is equivalent to the identity $R_{y\cdot x}^{-1}R_xR_y=L_xR_y^{-1}L_x^{-1}R_y$ in the loop $(M,\cdot)$. Finally, observe that 
$$R_{y{\cdot}x}^{-1}R_xR_y=L_xR_y^{-1}L_x^{-1}R_y=L_x(R_{y^{-1}}L_x^{-1}R_{y^{-1}}^{-1}L_x)L_x^{-1}=L_x(L_{y^{-1}x}^{-1}L_{y^{-1}}L_x)L_x^{-1}=L_xL_{x^{-1}y}L_y^{-1}.$$
\end{proof}

Now we apply semi-automorphisms to prove the following important associativity property of Moufang loops.

\begin{theorem}\label{t:semi-Moufang} Let $x,y,z$ be elements of a Moufang loop $M$. If $(x{\cdot}y){\cdot}z=x{\cdot}(y{\cdot}z)$, then for every points $a,b,c\in\{x,y,z,x^{-1},y^{-1},z^{-1}\}$ we have $(a{\cdot}b){\cdot}c=a{\cdot}(b{\cdot}c)$.
\end{theorem}

\begin{proof} We divide the proof of Theorem~\ref{t:semi-Moufang} into four lemmas. Let $M$ be a Moufang loop.

\begin{lemma}\label{l:Moufang-as1} For all $x,y,z\in M$, the equality $(x{\cdot}y){\cdot}z=x{\cdot}(y{\cdot}z)$ implies $(x{\cdot}y){\cdot}z^{-1}\,{=}\,x{\cdot}(y{\cdot}z^{-1})$.
\end{lemma}

\begin{proof} By Theorem~\ref{t:inner=>pseudo-semi-auto}, the permutation $L_{x{\cdot}y}^{-1}L_xL_y$ is a semi-automorphism of the Moufang loop $M$. Assuming that $x{\cdot}(y{\cdot}z)=(x{\cdot}y){\cdot}z$, we conclude that  
$$
L_{x{\cdot}y}^{-1}L_xL_y(z)=(x{\cdot} y)^{-1}{\cdot}(x{\cdot}(y{\cdot} z))=(x{\cdot} y)^{-1}((x{\cdot} y){\cdot} z)=z.$$ Applying Lemma~\ref{l:pseudo=>inverse}, we obtain the equality $L_{xy}^{-1}L_xL_y(z^{-1})=z^{-1},$ equivalent to the desired equality $x{\cdot}(y{\cdot}z^{-1})=(x{\cdot}y){\cdot}z^{-1}$.
\end{proof}

\begin{lemma}\label{l:Moufang-as2} For all $x,y,z\in M$, the equality $(x{\cdot}y){\cdot}z=x{\cdot}(y{\cdot}z)$ implies $(y{\cdot}z){\cdot}x^{-1}\,{=}\,y{\cdot}(z{\cdot}x^{-1})$.
\end{lemma}

\begin{proof} It follows from $L_xL_y(z)=x{\cdot}(y{\cdot}z)=(x{\cdot}y){\cdot}z=L_{x{\cdot}y}(z)$ that $z=L_{x{\cdot}y}^{-1}L_xL_y(z)$. 
Applying the identity $L_{x{\cdot}y}^{-1}L_xL_y=R_xL^{-1}_yR_{x}^{-1}L_y$ from Proposition~\ref{p:Moufang-pauto}(1)  to the element $z$ and taking into account that $L_y^{-1}=L_{y^{-1}}$ and $R_x^{-1}=R_{x^{-1}}$, we obtain the equality
$$z=L_{x{\cdot}y}^{-1}L_xL_y(z)=R_xL_y^{-1}R_{x}^{-1}L_{y}(z)=
(y^{-1}{\cdot}((y{\cdot}z){\cdot}x^{-1}){\cdot}x,$$
which implies the desired equality $y{\cdot}(z{\cdot}x^{-1})=(y{\cdot}z){\cdot}x^{-1}$, by the inversivity of the loop $M$.
\end{proof}

\begin{lemma}\label{l:Moufang-as3}  Let $x,y,z$ be any elements of $M$. If $(x{\cdot}y){\cdot}z=x{\cdot}(y{\cdot}z)$, then 
\begin{enumerate}
\item $(y{\cdot}z){\cdot}x=y{\cdot}(z{\cdot}x)$;
\item $(z{\cdot}x){\cdot}y=z{\cdot}(x{\cdot}y)$;
\item $(y{\cdot}x){\cdot}z=y{\cdot}(x{\cdot}z)$;
\item $(z{\cdot}y){\cdot}x=z{\cdot}(y{\cdot}x)$;
\item $(x{\cdot}z){\cdot}y=x{\cdot}(z{\cdot}y)$.
\end{enumerate}
\end{lemma}

\begin{proof} 1. By Lemma~\ref{l:Moufang-as2}, the equality $(x{\cdot}y){\cdot}z=x{\cdot}(y{\cdot}z)$ implies $(y{\cdot}z){\cdot}x^{-1}=y{\cdot}(z{\cdot}x^{-1})$, and by Lemma~\ref{l:Moufang-as1}, the latter equality implies $(y{\cdot}z){\cdot}x=y{\cdot}(z{\cdot}x)$.
\smallskip

2. By the preceding item, the equality $(x{\cdot}y){\cdot}z=x{\cdot}(y{\cdot}z)$ implies $(y{\cdot}z){\cdot}x=y{\cdot}(z{\cdot}x)$ and the latter equality implies $(x{\cdot}z){\cdot}y=z{\cdot}(x{\cdot}y)$.
\smallskip

3. By Lemma~\ref{l:Moufang-as2}, $(x{\cdot}y){\cdot}z=x{\cdot}(y{\cdot}z)$ implies $(y{\cdot}z){\cdot}x^{-1}=y{\cdot}(z{\cdot}x^{-1})$, and  the latter equality implies $(z{\cdot}x^{-1}){\cdot}y^{-1}=z{\cdot}(x^{-1}{\cdot}y^{-1})$. Since the inversion of $M$ is an anti-isomorphism of $M$, the equality $(z{\cdot}x^{-1}){\cdot}y^{-1}=z{\cdot}(x^{-1}{\cdot}y^{-1})$ implies the equalities
$$(y{\cdot}x){\cdot}z^{-1}=(z{\cdot}(x^{-1}{\cdot}y^{-1}))^{-1}=((z{\cdot}x^{-1}){\cdot}y^{-1})^{-1}=y{\cdot}(x{\cdot}z^{-1}).$$
By Lemma~\ref{l:Moufang-as1}, the equality $(y{\cdot}x){\cdot}z^{-1}=y{\cdot}(x{\cdot}z^{-1})$ implies the equality $(y{\cdot}x){\cdot}z=y{\cdot}(x{\cdot}z)$.
\smallskip

4,5. Applying Lemma~\ref{l:Moufang-as3}(3) to the equalities $(y{\cdot}z){\cdot}x=y{\cdot}(z{\cdot}x)$ and $(z{\cdot}x){\cdot}y=z{\cdot}(x{\cdot}y)$ proved in Lemma~\ref{l:Moufang-as3}(1,2), we obtain the equalities $(z{\cdot}y){\cdot}x=z{\cdot}(y{\cdot}x)$ and $(x{\cdot}z){\cdot}y=x{\cdot}(z{\cdot}y)$. 
\end{proof}

\begin{lemma}\label{l:Moufang-as4} Let $x,y,z$ be any elements of $M$. If $(x{\cdot}y){\cdot}z=x{\cdot}(y{\cdot}z)$, then $(a{\cdot}b){\cdot}c=a{\cdot}(b{\cdot}c)$ for every $a,b,c\in\{x,y,z\}$.
\end{lemma}

\begin{proof} If $b\in \{a,c\}$, then the equality $(a{\cdot}b){\cdot}c=a{\cdot}(b{\cdot}c)$ follows from the alternativity of the Moufang loop $M$. So, we assume that $b\notin\{a,c\}$. If $a=c$, then $(a{\cdot}a){\cdot}b=a{\cdot}(a{\cdot}b)$ by the alternativity of $M$. By Lemma~\ref{l:Moufang-as3}, the equality $(a{\cdot}a){\cdot}b=a{\cdot}(a{\cdot}b)$ implies the equality $(a{\cdot}b){\cdot}a=a{\cdot}(b{\cdot}a)$ which coincides with the equality $(a{\cdot}b){\cdot}c=a{\cdot}(b{\cdot}c)$ because $a=c$. So, assume that $a\ne c$. Then $a,b,c$ are tree distinct points of the set $\{x,y,z\}$. In this case the equality $(a{\cdot}b){\cdot}c=a{\cdot}(b{\cdot}c)$ follows from Lemma~\ref{l:Moufang-as3}.
\end{proof}

Now we can complete the proof of Theorem~\ref{t:semi-Moufang}. Assume that $(x{\cdot}y){\cdot}z=x{\cdot}(y{\cdot}z)$. Lemma~\ref{l:Moufang-as4} implies that for every $a,b,c\in \{x,y,z\}$ we have $(a{\cdot}b){\cdot}c=a{\cdot}(b{\cdot}c)$. By Lemma~\ref{l:Moufang-as2} the latter equality implies $(b{\cdot}c){\cdot}a^{-1}=b{\cdot}(c{\cdot}a^{-1})$, $(c{\cdot}a^{-1}){\cdot}b^{-1}=c{\cdot}(a^{-1}{\cdot}b^{-1})$ and $(a^{-1}{\cdot}b^{-1}){\cdot}c^{-1}=a^{-1}{\cdot}(b^{-1}{\cdot}c^{-1})$.
Those three equalities and Lemma~\ref{l:Moufang-as3} imply that $(a{\cdot}b){\cdot}c=a{\cdot}(b{\cdot}c)$ for every elements $a,b,c\in\{x,y,z,x^{-1},y^{-1},z^{-1}\}$. 
\end{proof}

\begin{corollary}\label{c:diass} Let $x,y$ be any elements of a Moufang loop $M$. Then $(a{\cdot}b){\cdot}c=a{\cdot}(b{\cdot}c)$ for every elements $a,b,c\in\{x,y,x^{-1},y^{-1}\}$.
\end{corollary}

\begin{proof} By Proposition~\ref{p:M12=>flexible}, the Moufang loop $M$ is flexible and hence $(x{\cdot}y){\cdot}x=x{\cdot}(y{\cdot}x)$. By Theorem~\ref{t:semi-Moufang}, the equality $(a{\cdot}b){\cdot}c=a{\cdot}(b{\cdot}c)$ holds for every elements $a,b,c\in\{x,y,x^{-1},y^{-1}\}$.
\end{proof}

\section{The Moufang Theorem}

In this section we present a proof of the following fundamental and difficult theorem due to Ruth Moufang. Her proof was rather complicated and long (about eight pages). A shorter (but still difficult) proof was found by Richard Bruck in 1951. Bruck's proof is  reproduced in thee main textbooks in quasigroup theory \cite{Bruck1958}, \cite{Pflugfelder1990}, \cite{Belousov}. In this section we reproduce a relatively simple proof of the Moufang Theorem, found by \index[person]{Dr\'apal}Ale\v s Dr\'apal\footnote{{\bf Ale\v s Dr\'apal} is a prominent Czech mathematician, specializing in algebra, matematic logic and computational mathematics. He is a professor at the Department of Algebra at the Faculty of Mathematics and Physics, Charles University, Prague. Professor Drapal graduated from the Charles University of Prague, where he defended his RNDr Thesis ``Latinsk\'e \v ctverce a parci\'aln\'\i\ grupoidy'' in 1989 under supervision of Karle Drbohlav and Tom\'a\v s Kepka. One of Ph.D. students of Dr\'apal is Petr Vojt\v echovsk\'y.} \cite{Drapal2011} in 2011.

\begin{theorem}[Moufang, 1935; Bruck, 1951; Dr\'apal, 2011]\label{t:Moufang-xyz} Let $x,y,z$ be elements of a  Moufang loop $M$. If $(x{\cdot}y){\cdot}z=x{\cdot}(y{\cdot}z)$, then the set $\{x,y,z\}$ is contained in an associative subloop of the loop $M$.
\end{theorem}
 
\begin{proof} Let $M$ be a Moufang loop. We say that a subset $X\subseteq M$ \defterm{generates} $M$ if there is no proper subloop of $M$ that contains $X$. Since the Moufang loop $M$ has the inverse property, $X$ generates $M$ if and only if each $u\in M$ can be written as the product of elements of the set $X^\pm\defeq X\cup X^{-1}$. We shall give notation to some of such products. 
Namely, for a sequence $u_1,\dots,u_k$ of elements of the loop $M$, we define their left and right products $\ell(u_1,\dots,u_k)$ and $r(u_1,\dots,u_k)$ by the recursive formulas $$\ell(u_1,\dots,u_k)\defeq u_1{\cdot}\ell( u_2,\dots,u_k)\quad\mbox{and}\quad r(u_1,\dots,u_k)\defeq r(u_1,\dots,u_{k-1}){\cdot} u_k$$where $\ell(u_1)=r(u_1)=u_1$.

\begin{lemma}\label{l:Drapal1} Let $M$ be a Moufang loop generated by a set $X$ such that $\ell(x_1,\dots,x_k)=r(x_1,\dots,x_k)$ for all finite sequences $x_1,\dots,x_k$ over $X^\pm$. Then $M$ is a group.
\end{lemma}

\begin{proof} First we shall prove, by induction on $n$, that
$$\ell(u_1,\dots,u_n)\cdot\ell(v_1,\dots,v_m) = \ell(u_1,\dots,u_n, v_1,\dots,v_m)$$
for all $u_1,\dots,u_n \in X^\pm$ and $v_1,\dots,v_m \in X^\pm$. The case $n=1$ is clearly true. Assume $n\ge 2$, and put $x \defeq u_1$, $s \defeq\ell(u_2,\dots,u_n)$ and 
$t \defeq \ell(v_1,\dots,v_m)=r(v_1,\dots,v_m)$.

Applying the Moufang identity, express $(x{\cdot}s){\cdot} t$ as $(x{\cdot} s){\cdot} (t{\cdot}x^{-1}{\cdot} x) = x{\cdot}(s {\cdot}( t{\cdot}x^{-1})){\cdot}x$. Now, $tx^{-1} = r(v_1,\dots,v_m)x^{-1}=r(v_1,\dots,v_m, x^{-1}) = \ell(v_1,\dots,v_m, x^{-1})$, and so, by the induction assumption, 
$$
\begin{aligned}
(x{\cdot}s){\cdot} t &= x{\cdot}(s {\cdot}( t{\cdot}x^{-1})){\cdot}x = x{\cdot}(s{\cdot}\ell(v_1,\dots,v_m, x^{-1})){\cdot}x= x{\cdot}\ell(u_2,\dots,u_n, v_1,\dots,v_m, x^{-1}){\cdot}x\\
& =  x{\cdot}r(u_2,\dots,u_n, v_1,\dots,v_m, x^{-1}){\cdot}x=x{\cdot}(r(u_2,\dots,u_n, v_1,\dots,v_m){\cdot}x^{-1}){\cdot}x\\
& = x{\cdot}\ell(u_2,\dots,u_n, v_1,\dots,v_m) = \ell(u_1, u_2,\dots,u_n, v_1,\dots,v_m).
\end{aligned}$$ 

Set $a \defeq \ell(u_1,\dots,u_n)$, $b \defeq \ell(v_1,\dots,v_m)$ and $c \defeq\ell(w_1,\dots,w_p)$. The proved equality makes clear that both $a{\cdot} (b{\cdot}c)$ and $(a{\cdot}b){\cdot} c$ are equal to $\ell(u_1,\dots,u_n, v_1,\dots,v_m, w_1,\dots,w_p)$.
Hence $L = \{\ell(u_1,\dots,u_n):u_1,\dots,u_n \in X^\pm\}$ is a subsemigroup of $M$ that is generated by $X^\pm$. The semigroup $L$ is a group since all generating elements possess inverses. Therefore $L = M$.
\end{proof}

\begin{lemma}\label{l:Drapal2} Let $M$ be a Moufang loop generated by $X = \{x, y, z\}$. If $x{\cdot}(y{\cdot}z) =
(x{\cdot}y) {\cdot} z$, then $u_0 {\cdot}  (\ell(u_1,\dots,u_k){\cdot}u_{k+1}) = (u_0{\cdot}\ell(u_1,\dots,u_k) ){\cdot}  u_{k+1}$ and $\ell(u_0,\dots,u_{k+1}) = r(u_0,\dots,u_{k+1})$ for any sequence $u_0,\dots,u_{k+1}$ of elements of $X^\pm$, $k\ge 1$.
\end{lemma}

\begin{proof}  There are two equalities to be proved. We shall proceed by a common
induction on $k$. The case $k = 1$ follows from the assumption and from Lemma~\ref{l:Drapal1}.
In the induction step we shall first prove the equality
$u_0 {\cdot}  (\ell(u_1,\dots,u_k){\cdot}u_{k+1}) = (u_0{\cdot}\ell(u_1,\dots,u_k)){\cdot}  u_{k+1}$.
The case $u_{k+1}\in\{u_0,u_0^{-1}\}$ follows from Corollary~\ref{c:diass}. We can thus assume that, say, $u_0 = x$ and $u_{k+1} = y$.
Put $s\defeq \ell(u_2,\dots,u_k)$. We wish to prove the equality 
$
x{\cdot} ((u_1{\cdot}s){\cdot} y)=(x{\cdot}(u_1{\cdot}s)){\cdot}y
$. 
If $u_1 = x^{-1}$, then this is equivalent to $(x^{-1}{\cdot}s) {\cdot}  y = x^{-1} {\cdot}  (s{\cdot}y)$, and that holds by the induction assumption. 

If $u_1=x$, then apply the inductive assumption and conclude that $(x{\cdot}s){\cdot} y=x{\cdot}(s{\cdot}y)$ and hence $x^{-1}{\cdot}((x{\cdot}s){\cdot}y)=x^{-1}\cdot(x{\cdot}(s{\cdot}y))=s{\cdot}y=(x^{-1}{\cdot}(x{\cdot}s))\cdot y$. 
By Theorem~\ref{t:semi-Moufang}, the equality $x^{-1}{\cdot}((x{\cdot}s){\cdot}  y = (x^{-1} {\cdot} (x{\cdot}s)){\cdot}y$ implies  $x{\cdot}((x{\cdot}s){\cdot}  y = (x{\cdot} (x{\cdot}s)){\cdot}y$, which coincides with the equality $x{\cdot}((u_1{\cdot}s){\cdot}  y = (x{\cdot} (u_1{\cdot}s)){\cdot}y$. Thus we have solved the case $u_1\in\{x,x^{-1}\}$. 

The inductive assumption implies $(x{\cdot}s){\cdot}y=x{\cdot}(s{\cdot}y)$ and Theorem~\ref{t:semi-Moufang} implies $(x{\cdot}y){\cdot}s=x{\cdot}(y{\cdot}s)$ and 
$(x{\cdot}y^{-1}){\cdot}s=x{\cdot}(y^{-1}{\cdot}s)$. 
 Since the Moufang loop $M$ is right-Bol, $$x{\cdot}((y{\cdot}s){\cdot}y)=((x{\cdot}y){\cdot}s){\cdot}y=(x{\cdot}(y{\cdot}s)){\cdot}y.$$This proves the case $u_1=y$.
 
Since the Moufang loop $M$ is right-Bol, $x{\cdot}((y^{-1}{\cdot}s){\cdot}y^{-1})=((x{\cdot}y^{-1}){\cdot}s){\cdot}y^{-1}=(x{\cdot}(y^{-1}{\cdot}s)){\cdot}y^{-1}$. By Lemma~\ref{l:Moufang-as1}, the equality $x{\cdot}((y^{-1}{\cdot}s){\cdot}y^{-1})=(x{\cdot}(y^{-1}{\cdot}s)){\cdot}y^{-1}$ implies $$x{\cdot}((y^{-1}{\cdot}s){\cdot}y)=(x{\cdot}(y^{-1}{\cdot}s)){\cdot}y,$$ which solves the case $u_1=y^{-1}$.

 Therefore, the equality $$x{\cdot}(\ell(u_1,\dots,u_k){\cdot}y)=(x{\cdot}\ell(u_1,\dots,u_k){\cdot}y)$$ is proved if $u_1\in\{x,x^{-1},y,y^{-1}\}$. Since $\ell(u_1,\dots,u_k)=r(u_1,\dots,u_k)$, the mirror argument proves this equality in case $u_k\in\{x,x^{-1},y,y^{-1}\}$.
 
 It remains to consider the case $u_1,u_k\in\{z,z^{-1}\}$. It suffices
to consider the cases $u_1 = u_k = z$ and $u_1=u_k^{-1}=z$. 

Put $w = \ell(u_2,\dots,u_{k-1})=r(u_2,\dots,u_{k-1})$. In case $u_1=u_k=z$, we need to check that $x{\cdot} ((z{\cdot}w{\cdot}z) {\cdot}  y) = (x {\cdot}(z{\cdot}w{\cdot}z)){\cdot}y$, i.e. that
$R_y^{-1}L_x^{-1}R_yL_x(zwz) = zwz$. The induction assumption yields $R_y^{-1}L_x^{-1}R_yL_x(w) = w$. Since $(x{\cdot}y){\cdot}z=x{\cdot}(y{\cdot}z)$, Theorem~\ref{t:semi-Moufang} implies $R_y^{-1}L_x^{-1}R_yL_x(z)=z$. By Theorem~\ref{t:inner=>pseudo-semi-auto}, the inner mapping $R_y^{-1}L_x^{-1}R_yL_x=R_{y^{-1}}L_x^{-1}R_{y^{-1}}^{-1}L_x$ is a pseudo-automor\-phism and a semi-automorphism of the Moufang loop $M$. Therefore, $R_y^{-1}L_x^{-1}R_yL_x(z{\cdot}w{\cdot}z)=z{\cdot}w{\cdot}z$.

In case $u_1=u_k^{-1}=z$, we need to show that $x{\cdot}(z{\cdot}(w{\cdot}z^{-1})) {\cdot}  y) = (x {\cdot}  (z{\cdot}(w{\cdot}z^{-1})){\cdot} y$, i.e. that
$\alpha(z^{-1}) = z^{-1}$ where $\alpha = L^{-1}_w L^{-1}_z (R_y^{-1} L_x^{-1}R_yL_x)L_zL_w$. By the induction assumption,
$(x{\cdot}(z{\cdot}w)) {\cdot}  y=x {\cdot}((z{\cdot}w){\cdot}y)$, which implies $\alpha(e) =e$. Theorem~\ref{t:inner=>pseudo-semi-auto} implies that the inner mapping $\alpha$ is a pseudo-automorphism of the loop $M$. The already-proved equality $(x{\cdot}(z{\cdot}w{\cdot}z)) {\cdot}  y = x{\cdot} ((z{\cdot}w{\cdot}z){\cdot}y)$ can be expressed as $\alpha(z) = z$. By Lemma~\ref{l:pseudo=>inverse}, $\alpha(z^{-1})=z^{-1}$.
\smallskip

To finish the induction step observe that 
$$
\begin{aligned}
\ell&(u_0,\dots,u_{k+1}) = u_0{\cdot}\ell(u_1,\dots,u_{k+1}) =
u_0{\cdot}r(u_1,\dots,u_{k+1}) = u_0{\cdot}( r(u_1,\dots,u_k){\cdot}u_{k+1})\\
&=(u_0{\cdot}\ell(u_1,\dots,u_k)) {\cdot}  u_{k+1} =
\ell(u_0, u_1,\dots,u_k){\cdot}u_{k+1} = r(u_0,\dots,u_k){\cdot}u_{k+1} = r(u_0,\dots,u_{k+1}).
\end{aligned}
$$
\end{proof}

Now we are aple to complete the proof of Theorem~\ref{t:Moufang-xyz}. Let $x,y,z$ be elements of a Moufang loop $M$ such that $(x{\cdot}y){\cdot}z=x{\cdot}(y{\cdot}z)$. 
Lemmas~\ref{l:Drapal1} and \ref{l:Drapal2} imply that the set $X=\{x,y,z\}$ generates a subgroup in the loop $M$.
\end{proof}

\section{Diassociative loops}

In this section we discuss some properties of diassociative loops.

\begin{definition} A loop $X$ is called \defterm{diassociative} if any elements $x,y\in X$ are contained in a subgroup of the loop $X$.
\end{definition}

Theorems~\ref{t:Moufang-xyz} and Proposition~\ref{p:M12=>flexible} imply 

\begin{corollary}\label{c:Moufang=>diassociative} Every Moufang loop is diassociative.
\end{corollary}

The following proposition shows that diassociative loops satisfy a weak form of the Lagrange property.

\begin{proposition}\label{p:di-Lagrange} Let $C$ be a cyclic subgroup of a diassociative loop $X$. For every elements $x,y\in X$, the set $x{\cdot}C$ and $y{\cdot}C$ are either disjoint or coincide. Consequently, the cardinality of the cyclic group $C$ divides the cardinality of $X$.
\end{proposition}

\begin{proof} Let $c$ be a generator of the cyclic group $C$. Assume that the cosets $x{\cdot}C$ and $y{\cdot}C$ contain a common point $z$. Since the loop $X$ is diassociative, the set $\{z,c\}$ generates a subgroup $H$ of the loop $X$. Since the diassociative loop $X$ is inversive, $z\in (x{\cdot}C)\cap(y{\cdot}C)$ implies  
 $x,y\in z{\cdot}C\subseteq H$. The associativity of the group $H$ ensures that $x{\cdot}C\subseteq (z{\cdot}C){\cdot}C=z{\cdot}(C{\cdot}C)=z{\cdot}C\subseteq H$. On the other hand, $z{\cdot}C\subseteq (x{\cdot}C){\cdot}C=x{\cdot}(C{\cdot}C)=x{\cdot}C$ and hence $x{\cdot}C=z{\cdot}C$. By analogy we can prove that $y{\cdot}C=z{\cdot}C$ and hence $x{\cdot}C=y{\cdot}C$. 

Since $\{g{\cdot}C:g\in X\}$ is a partition of the set $X$ into paiwise dosjoint subsets of cardinality $|C|$, the cardinality of the cyclic group $C$ divides the cardinality of $X$.
\end{proof}

Since Moufang loops are diassociative, Proposition~\ref{p:di-Lagrange}  implies that the cardinality of any cyclic subgroup of a Moufang  loop divides the cardinality of the loop. In fact, Moufang loops possess a much stronger property.

\begin{theorem}[Gagola--Hall and Grishkov--Zavarnitsine, 2005] The cardinality of any subloop of a Moufang loop divides the cardinality of the loop.
\end{theorem}

The proof of this theorem is difficult and can be found in the papers \cite{GagolaHall2005} and \cite{Grish-Zava2005}.

\section{Commutative Moufang loops}

In this section we present some additional results on the structure of commutative Moufang loops. Commutative Moufang loops have been already considered in Section~\ref{s:comMoufang}. 

We recall that the \defterm{centre} of a Moufang loop $X$ is the set
$$\{z\in X:\forall x,y\in X\;\;([x,y,z]=e=[x,z])\}.$$
For commutative Moufang loops, the centre coincides with the set $$\{z\in X:\forall x,y\in X\;[x,y,z]=e\}=\{z\in X:\forall x,y\in X\;\;(x\cdot y)\cdot z=x\cdot(y\cdot z)\}.$$ 

\begin{theorem}\label{t:cubic-Moufang} For every element $a$ of a commutative Moufang loop $M$ its cube $a^3$ belongs to the centre of $M$.
\end{theorem}

\begin{proof} Let $\bar a\defeq a^{-1}$. By Proposition~\ref{p:Moufang<=>}, the triples $(L_a,R_a,R_aL_a)$ and $(L_{\bar a},R_{\bar a},R_{\bar a}L_{\bar a})$ are autotopisms of $M$. Applying Theorem~\ref{t:autotopism3}(3) to the autotopism $(L_{\bar a},R_{\bar a},R_{\bar a}L_{\bar a})$, we conclude that the triple $(R_{\bar a},iR_{\bar a}L_{\bar a}i,iL_{\bar a}i)$ is an autotopism of $M$. Since autotopisms form a group, the triple $(L_{a}R_{\bar a},R_{a}iR_{\bar a}L_{\bar a}i,R_{a}L_{a}iL_{\bar a}i)$ is an autotopism of $M$. By Lemma~\ref{l:autotopism-e}, the inner mapping $T_a\defeq L_{a}R_{\bar a}$ is a pseudo-automorphism with the companion 
$$c=R_{a}iR_{\bar a}L_{\bar a}i(e)=a(\bar ae\bar a)^{-1}=a^{3}.$$ Then $(T_a,R_cT_a,R_cT_a)$ is an autotopism of $M$ and by the commutativity and diassociativity of the Moufang magma $M$, for every $x,y\in M$, we have the identities
$$x{\cdot} (y{\cdot} a^3)=(a{\cdot}x{\cdot}\bar a)\cdot((a{\cdot}y{\cdot}\bar a){\cdot}a^3)=T_a(x)\cdot R_cT_a(y)=R_cT_a(x{\cdot} y)=(a{\cdot}(x{\cdot}y){\cdot}\bar a){\cdot}a^3=(x{\cdot}y){\cdot}a^3,$$
witnessing that the element $a^3$ belongs to the centre of the commutative Moufang loop $M$. 
\end{proof} 

\begin{theorem}[Bruck--Slaby, 1958]\label{t:Bruck-Slaby} Every  non-trivial finitely generated commutative Moufang loop has nontrivial centre.
\end{theorem}

\begin{proof} The proof of this theorem is difficult and can be found in \cite[10.1]{Bruck1958}. Here we present a short proof for finite commutative Moufang loops. So, let $M$ be a finite commutative Moufang loop. If $M$ does not have exponent 3, then it contains an element $a$ whose cube $a^3$ is not trivial. By Theorem~\ref{t:cubic-Moufang}, $a^3$ belongs to the centre of $M$, witnessing that it is not trivial. So, assume that the commutative Moufang loop has exponent $3$. By  Proposition~\ref{p:Ali-Slaney2}, $|M|$ is an odd number, and by Ss\"org\H o Theorem~\ref{t:Scorgo}, $M$ has non-trivial nucleus, which coincides with the centre of $M$, by the commutativity of $M$. 
\end{proof}

\section{The interplay between various properties of loops}

The properties of loops, discussed in the preceding sections relate as follows:
$$
\xymatrix@C=-5pt@R=20pt{
&\mbox{associative}\ar@{=>}[d]\\
&\mbox{Moufang}\ar@{<=>}[d]\\
\mbox{left-Bol}\ar@{=>}[ddd]&\mbox{Bol}\ar@{=>}[d]\ar@{=>}[l]\ar@{=>}[r]&\mbox{right-Bol}\ar@{=>}[ddd]\\
&\mbox{diassociative}\ar@{=>}[d]&\\
&\mbox{inversive}\ar@{=>}[ld]\ar@{=>}[rd]&\\
\mbox{left-inversive}\ar@{=>}[rd]&\mbox{Boolean}\ar@{=>}[d]
&\mbox{right-inversive}\ar@{=>}[ld]\\
&\mbox{invertible}\\
&\mbox{commutative}\ar@{=>}[u]
}
$$


\begin{example}[Warren Smith] The 8-element loop with multiplication table
$$
\begin{array}{c|cccccccc}
*&0&1&2&3&4&5&6&7\\
\hline
0&0&1&2&3&4&5&6&7\\
1&1&2&3&4&5&6&7&0\\
2&2&5&4&7&6&3&0&1\\
3&3&6&5&0&7&2&1&4\\
4&4&3&6&1&0&7&2&5\\
5&5&4&7&6&1&0&3&2\\
6&6&7&0&5&2&1&4&3\\
7&7&0&1&2&3&4&5&6\\
\end{array}
$$
is inversive but not power-associative and hence not diassociative.
\end{example}

\begin{example}[Warren Smith] The 10-element loop with multiplication table
$$
\begin{array}{c|cccccccccc}
*&0&1&2&3&4&5&6&7&8&9\\ 
\hline
0&0&1&2&3&4&5&6&7&8&9\\
1&1&0&3&2&7&9&8&4&6&5\\
2&2&3&0&1&9&8&7&6&5&4\\
3&3&2&1&0&8&7&9&5&4&6\\
4&4&7&9&8&0&6&5&1&3&2\\
5&5&9&8&7&6&0&4&3&2&1\\
6&6&8&7&9&5&4&0&2&1&3\\
7&7&4&6&5&1&3&2&0&9&8\\
8&8&6&5&4&3&2&1&9&0&7\\
9&9&5&4&6&2&1&3&8&7&0\\
\end{array}
$$
is diassociative, commutative, and Boolean, but not Moufang.
\end{example}

\begin{example}[Warren Smith]\label{ex:Smith14} The 14-element loop with multiplication table
$$
\begin{array}{c|cccccccccccccc}
*&0&1&2&3&4&5&6&7&8&9&A&B&C&D\\
\hline
0&0&1&2&3&4&5&6&7&8&9&A&B&C&D\\
1&1&2&3&4&5&6&0&A&7&8&D&C&9&B\\
2&2&3&4&5&6&0&1&D&A&7&B&9&8&C\\
3&3&4&5&6&0&1&2&8&9&C&7&D&B&A\\
4&4&5&6&0&1&2&3&C&B&D&9&7&A&8\\
5&5&6&0&1&2&3&4&B&D&A&C&8&7&9\\
6&6&0&1&2&3&4&5&9&C&B&8&A&D&7\\
7&7&8&9&A&B&C&D&0&1&2&3&4&5&6\\
8&8&9&C&7&D&B&A&3&0&1&6&5&2&4\\
9&9&C&B&8&A&D&7&6&3&0&4&2&1&5\\
A&A&7&8&D&C&9&B&1&2&5&0&6&4&3\\
B&B&D&A&C&8&7&9&5&4&6&2&0&3&1\\
C&C&B&D&9&7&A&8&4&6&3&5&1&0&2\\
D&D&A&7&B&9&8&C&2&5&4&1&3&6&0\\
\end{array}
$$
is alternative and power-associative but not left-inversive and not right-inversive.
\end{example}

\begin{example}[Warren Smith] The 12-element loop with multiplication table
$$
\begin{array}{c|cccccccccccc}
*&0&1&2&3&4&5&6&7&8&9&A&B\\
\hline
0&0&1&2&3&4&5&6&7&8&9&A&B\\
1&1&2&0&5&3&4&8&6&7&B&9&A\\
2&2&0&1&4&5&3&7&8&6&A&B&9\\
3&3&4&5&0&1&2&9&B&A&6&8&7\\
4&4&5&3&2&0&1&B&A&9&8&7&6\\
5&5&3&4&1&2&0&A&9&B&7&6&8\\
6&6&7&8&9&B&A&0&1&2&3&5&4\\
7&7&8&6&B&A&9&2&0&1&5&4&3\\
8&8&6&7&A&9&B&1&2&0&4&3&5\\
9&9&A&B&6&8&7&3&5&4&0&1&2\\
A&A&B&9&8&7&6&5&4&3&2&0&1\\
B&B&9&A&7&6&8&4&3&5&1&2&0\\
\end{array}
$$
is Moufang but not associative (this loop is the Chein extension $M(S_3,2)$ of the symmetric group $\mathcal S_3$).
\end{example}


\chapter{Alternative rings}\label{ch:alt-rings}

\rightline{\em Civilization advances by extending the number of important}

\rightline{\em operation which we can perform without thinking about them.}

\rightline{Alfred North Whitehead}
\bigskip

In this chapter we analyze the algebraic stucture of alternative (division) rings and present several characterizations of such ring, which will be applied in subsequent chapters for studying the geometry of Moufang planes. 

\section{Rings}

\begin{definition} A \index{ring}\defterm{ring} is a set $R$ endowed with two binary operations $+,\cdot:R\times R\to R$ and two distinct elements $0,1\in R$ satisfying the following axioms:
\begin{itemize}
\item $\forall x,y,z\in R\;\;x+(y+z)=(x+y)+z$;
\item $\forall x\in R\;\;x+0=x=0+x$;
\item $\forall x\in R\;\exists y\in R\;\;x+y=0=y+x$;
\item $\forall x\in R\;\;x\cdot 1=x=1\cdot x$;
\item $\forall a,x,y\in R\;\;a\cdot(x+y)=a\cdot x+a\cdot y$;
\item $\forall x,y,b\in R\;\;(x+y)\cdot b=x\cdot b+y\cdot b$.
\end{itemize}
\end{definition}

Let us note the the operation of addition in a ring is always commutative.
\index[person]{Hankel}

\begin{theorem}[Hankel\footnote{
{\bf Hermann Hankel} (1839 -- 1873) was a German mathematician. Having worked on mathematical analysis during his career, he is best known for introducing the Hankel transform and the Hankel matrix.
Hankel was born on 14 February 1839 in Halle, Germany. His father, Wilhelm Gottlieb Hankel, was a physicist. Hankel studied at Nicolai Gymnasium in Leipzig before entering Leipzig University in 1857, where he studied with Moritz Drobisch, August Ferdinand Möbius and his father. In 1860, he started studying at University of G\"ottingen, where he acquired an interest in function theory under the tutelage of Bernhard Riemann. Following the publication of an award winning article, he proceeded to study under Karl Weierstrass and Leopold Kronecker in Berlin. He received his doctorate in 1862 at Leipzig University. Receiving his teaching qualifications a year after, he was promoted to an associate professor at Leipzig University in 1867. At the same year, he received his full professorship in University of Erlangen--Nuremberg and spent his last four years in University of T\"ubingen. In 1867, he published ``Theorie der Complexen Zahlensysteme'', a treatise on complex analysis. His works on the theory of functions include 1870's ``Untersuchungen \"uber die unendlich oft oscillirenden und unstetigen functionen'' and his 1871 article ``Grenze'' for the Ersch-Gruber Encyklop\"adie. His work for Mathematische Annalen has highlighted the importance of Bessel functions of the third kind, which were later known as Hankel functions. His 1867 exposition on complex numbers and quaternions is particularly memorable. For example, Fischbein notes that he solved the problem of products of negative numbers by proving the following theorem: ``The only multiplication in $\mathbb R$ which may be considered as an extension of the usual multiplication in $\mathbb R_+$ by respecting the law of distributivity to the left and the right is that which conforms to the rule of signs.'' Furthermore, Hankel draws attention to the linear algebra that Hermann Grassmann had developed in his ``Extension Theory'' in two publications. This was the first of many references later made to Grassmann's early insights on the nature of space.}, 1867]\label{t:Hankel} For every ring $R$, the group $(R,+)$ is commutative.
\end{theorem}

\begin{proof} Given any elements $x,y\in R$, apply the left and right distributivites to obtain
$$(1+1)\cdot (x+y)=1\cdot (x+y)+1\cdot(x+y)=1\cdot x+1\cdot y+1\cdot x+1\cdot y=x+y+x+y$$and
$$(1+1)\cdot (x+y)=(1+1)\cdot x+(1+1)\cdot y=1\cdot x+1\cdot x+1\cdot y+1\cdot y=x+x+y+y,$$
which implies $x+y=y+x$ because $(R,+)$ is a group.
\end{proof}

\begin{definition} A ring $R$ is called
\begin{itemize}
\item \index{commutative ring}\index{ring!commutative}\defterm{commutative} if $\forall x,y\in R\;\;x\cdot y=y\cdot x$;
\item \index{associative ring}\index{ring!associative}\defterm{associative} if $\forall x,y,z\in R\;\;(x\cdot y)\cdot z=x\cdot(y\cdot z)$;
\item \index{alternative ring}\index{ring!alternative}\defterm{alternative} if $\forall x,z\in R\;\forall y\in\{x,z\}\;\;(x\cdot y)\cdot z=x\cdot(y\cdot z)$;
\item \index{division ring}\index{ring!division}\defterm{division} if $\forall a\in R\setminus\{0\}\;\forall b\in R\;\;\exists! x,y\in R\;\;a\cdot x=b=y\cdot a$;
\item \index{invertible ring}\index{ring!invertible}\defterm{invertible} if $\forall a\in R\setminus\{0\}\;\exists! a^{-1}\in R\;\;a\cdot a^{-1}=1=a^{-1}\cdot a$.
\end{itemize}
\end{definition}

Observe that every associative ring is alternative and corps are just invertible  associative rings.
 
\begin{remark} In classical (associative) algebra all rings are assumed to be associative. However in geometry nonassociative ring appear quite often, so it is more convenient to define ring as non-associative rings and then add this adjective if necessary. Such an approach also is adopted in the known textbook ``The role of Nonassociative Algebra in Projective Geometry'' \cite{Faulkner} of Faulkner. 
\end{remark}

\begin{exercise} Find an example of an alternative ring which is not associative.
\smallskip

{\em Hint:} The octonion ring $\mathbb O$.
\end{exercise}

\section{Some identities in alternative rings}\label{s:alt-identities}

In this section we shall prove some useful identities holding in (alternative) rings.

\begin{definition} For two elements $x,y$ of a ring $R$, the element
$$[x,y]\defeq x{\cdot} y-y{\cdot} x$$is called the \index{commutator}\defterm{commutator} of the elements $x,y$.
\end{definition}

The following elementary properties of the commutators follow immediately from the definition.

\begin{proposition} For two elements $x,y,z$ of a ring $R$,
\begin{enumerate}
\item $[x,y]=0$ if and only if $x{\cdot} y=y{\cdot} x$;
\item $[x,y]+[y,x]=0$;
\item $[x+y,z]=[x,z]+[y,z]$ and $[x,y+z]=[x,y]+[x,z]$.
\end{enumerate}
\end{proposition} 

Observe that a ring is commutative if and only if $[x,y]=0$ for every $x,y\in R$.

\begin{definition} For three elements $x,y,z$ of a ring $R$, the element
$$[x,y,z]\defeq (x{\cdot}y){\cdot} z-x{\cdot}(y{\cdot} z)$$is called the \index{associator}\defterm{associator} of the elements $x,y,z$.
\end{definition}

A ring $R$ is associative if and only if $[x,y,z]=0$ for any elements $x,y,z\in R$.

The distributivity of the multiplication over addition in a ring implies the following additivity property of the associator.

\begin{lemma}\label{l:associator-additive} Any elements $a,b,x,y$ of a ring satisfy the following identities:
$$
\begin{aligned}
&[x+y,a,b]=[x,a,b]+[y,a,b],\\
&[a,x+y,b]=[a,x,b]+[a,y,b],\\
&[a,b,x+y]=[a,b,x]+[a,b,y].
\end{aligned}
$$
\end{lemma}

\begin{lemma}[Teichm\"uller\footnote{{\bf Paul Julius Oswald Teichm\"uller} (1913 -- 1943) was a German mathematician who made contributions to complex analysis. He introduced quasiconformal mappings and differential geometric methods into the study of Riemann surfaces. Teichm\"uller spaces are named after him.
\newline 
Born in Nordhausen, Teichmüller attended the University of G\"ottingen, where he graduated in 1935 under the supervision of Helmut Hasse. His doctoral dissertation was on operator theory, though this was his only work on functional analysis. His next few papers were algebraic, but he switched his focus to complex analysis after attending lectures given by Rolf Nevanlinna. In 1937, he moved to the University of Berlin to work with Ludwig Bieberbach. Bieberbach was the editor of Deutsche Mathematik and much of Teichm\"uller's work was published in the journal, which made his papers hard to find in modern libraries before the release of his collected works.
\newline
A member of the Nazi Party (NSDAP) and Sturmabteilung (SA), the military wing of the NSDAP, from 1931, Teichmüller agitated against his Jewish professors Richard Courant and Edmund Landau in 1933. He was drafted into the Wehrmacht in July 1939 and took part in the invasion of Norway in 1940 before being recalled to Berlin to undertake cryptographic work with the Cipher Department of the High Command of the Wehrmacht. In 1942, he was released from his military duties and returned to teach at the University of Berlin. After the German defeat at Stalingrad in February 1943, he gave up his position in Berlin to volunteer for combat on the Eastern Front. He was killed in action in September 1943.}] Any elements $a,b,c,d$ of a ring satisfy the \index[person]{Teichm\"uller}\index{Teichm\"uller identity}\index{identity!Teichm\"uller}\defterm{Teichm\"uller identity}\textup{:}

$$a{\cdot}[b, c, d] - [a{\cdot}b, c, d] + [a, b{\cdot}c, d] - [a, b, c{\cdot}d] + [a, b, c]{\cdot}d = 0.$$
\end{lemma}

\begin{proof} The distributivity of the multiplication over addition in the ring $R$ implies
\begin{multline*}
a{\cdot}[b, c, d] - [a{\cdot}b, c, d] + [a, b{\cdot}c, d] - [a, b, c{\cdot}d] + [a, b, c]{\cdot}d\\
=a{\cdot}((b{\cdot}c){\cdot}d)-a{\cdot}(b{\cdot}(c{\cdot}d))-((a{\cdot}b){\cdot}c)\cdot d+(a{\cdot}b){\cdot}(c{\cdot} d) +(a{\cdot}(b{\cdot}c)){\cdot}d\\
-a{\cdot}((b{\cdot}c){\cdot}d)
 - (a{\cdot}b){\cdot}(c{\cdot}d)+a{\cdot}(b{\cdot}(c{\cdot}d))
+ ((a{\cdot}b){\cdot}c){\cdot}d-(a{\cdot}(b{\cdot} c)){\cdot}d=0
\end{multline*}
\end{proof}

Observe that a ring $R$ is alternative if and only if $$[x,x,y]=0=[x,y,y]$$ for every elements $x,y\in R$.

In the following proposition we establish some identities in alternative rings, which will be used in the proof of the Artin Theorem~\ref{t:Artin}.

\begin{proposition}\label{p:alternative9} For every elements $x,y,z$ of an alternative ring $R$ the following identities hold:
\begin{enumerate}
\item $[x,y,z]=-[y,x,z]=-[x,z,y]$;
\item $[x,y,x]=0$;
\item $[x,x{\cdot}y,y]=0$;
\item $[x,x{\cdot}y,z]+[x,x{\cdot}z,y]=0$;
\item $[x,y{\cdot}z,z]+[y,x{\cdot}z,z]=0$;
\item $x{\cdot}(y{\cdot}(x{\cdot}z))=(x{\cdot}(y{\cdot}x)){\cdot}z$;
\item $((x{\cdot}y){\cdot}z){\cdot}y=x{\cdot}((y{\cdot}z){\cdot}y)$;
\item $(x{\cdot}y){\cdot}(z{\cdot}x)=(x{\cdot}(y{\cdot}z)){\cdot}x$.
\end{enumerate}
\end{proposition} 

\begin{proof} 1. The alternativity of the ring $R$ and Lemma~\ref{l:associator-additive} imply 
$$[x,y,z]+[y,x,z]=[x{+}y,x{+}y,z]-[x,x,z]-[y,y,z]=0$$
and 
$$[x,y,z]+[x,z,y]=[x,y{+}z,y{+}z]-[x,y,y]-[x,z,z]=0.$$

2. The preceding statement and the alternativity of the ring $R$  imply
$$[x,y,x]=-[y,x,x]=0.$$

3. The alternativity of the ring $R$ implies
$$[x,x{\cdot}y,y]=(x{\cdot}(x{\cdot}y)){\cdot}y-x{\cdot}((x{\cdot}y){\cdot}y)=
((x{\cdot}x){\cdot}y){\cdot}y-x{\cdot}(x{\cdot}(y{\cdot}y))=(x{\cdot}x){\cdot}(y{\cdot}y)-(x{\cdot}x){\cdot}(y{\cdot}y)=0.
$$

4, 5. The preceding statement, the alternativity and the distributivity of the multiplication over addition in the ring $R$ imply 
$$[x,x{\cdot}y,z]+[x,x{\cdot}z,y]=[x,x{\cdot}(y{+}z),y{+}z]-[x,x{\cdot}y,y]-[x,x{\cdot}z,z]=0$$
and
$$[x,y{\cdot}z,z]+[y,x{\cdot}z,z]=[x{+}y,(x{+}y){\cdot}z,z]-[x,x{\cdot}z,z]-[y,y{\cdot}z,z]=0.$$

6. The statements (4), (1), (2) imply
$$
\begin{aligned}
0&=[x,x{\cdot}y,z]+[x,x{\cdot}z,y]=-[x{\cdot}y,x,z]-[x,y,x{\cdot}z]\\
&=-((x{\cdot}y){\cdot}x){\cdot}z+(x{\cdot}y){\cdot}(x{\cdot}z)-(x{\cdot}y){\cdot}(x{\cdot}z)+x{\cdot}(y{\cdot}(x{\cdot}z))\\
&=-((x{\cdot}y){\cdot}x){\cdot}z+x{\cdot}(y{\cdot}(x{\cdot}z))=-(x{\cdot}(y{\cdot}x)){\cdot}z+x{\cdot}(y{\cdot}(x{\cdot}z)).
\end{aligned}
$$

7. By analogy, the identity (7) can be deduced from (5), (1), and (2).

8. The identities (1) and (7) imply
$$
\begin{aligned}
&(x{\cdot}y){\cdot}(z{\cdot}x)-(x{\cdot}(y{\cdot}z)){\cdot}x=
-[x{\cdot}y,z,x]+[x,y,z]{\cdot}x\\
&=[z,x{\cdot}y,x]+[z,x,y]{\cdot}x=
-z{\cdot}((x{\cdot}y){\cdot}x)+((z{\cdot}x){\cdot}y){\cdot}x=0.
\end{aligned}
$$
\end{proof}

The identities (6) and (7) in Proposition~\ref{p:alternative9} have special names and will play an important role in studying the algebraic structure of Moufang planes.

\begin{definition} A ring $R$ is defined to be 
\begin{itemize}
\item \index{left-Bol ring}\index{ring!left-Bol}\defterm{left-Bol} if $a{\cdot}(b{\cdot}(a{\cdot}x))=(a{\cdot}(b{\cdot}a)){\cdot}x$ for every elements $a,b,x\in R$;
\item \index{right-Bol ring}\index{ring!right-Bol}\defterm{right-Bol} if $((x{\cdot}a){\cdot}b){\cdot}a=x{\cdot}((a{\cdot}b){\cdot}a)$ for every elements $x,a,b\in R$;
\item \index{Bol ring}\index{ring!Bol}\defterm{Bol} if $R$ is left-Bol and right-Bol.
\end{itemize}
\end{definition}

\begin{remark} A ring is (left or ring) Bol iff so is its multiplicative magma. By Example~\ref{ex:Smith14}, an alternative loop need not to be Bol. However, for rings we have the following characterization.
\end{remark}

\begin{proposition}\label{p:Bol<=>alternative} A ring is alternative if and only if it is Bol.
\end{proposition}

\begin{proof} The ``only if'' part follows from Proposition~\ref{p:alternative9}.
The prove the ``if'' part, assume that a ring $R$ is Bol. Given any elements $a,b\in R$, apply the left-Bol identity to the elements $a,1,a,b$ and obtain the left-alternative identity
$$a{\cdot}(a{\cdot}b)=a{\cdot}(1{\cdot}(a{\cdot}b))=(a{\cdot}(1{\cdot}a)){\cdot}b=(a{\cdot}a){\cdot}b.$$ Applying the right-Bol identity to the elements $a,b,1,b$, we obtain the right-alternative idenity
$$(a{\cdot}a){\cdot}b=(a{\cdot}(1{\cdot}a)){\cdot}b=a{\cdot}(1{\cdot}(a{\cdot}b))=a{\cdot}(a{\cdot}b).$$
\end{proof}

 In the following proposition we establish some useful properties of left-Bol rings. A deeper study of left-Bol rings will be continued in Section~\ref{s:Bol}.

\begin{proposition}\label{p:left-Bol} For every elements $a,b,c,d$ of a left-Bol ring the following identities hold:
\begin{enumerate}
\item $[a,a,b]=0$;
\item $[a,b,c]+[b,a,c]=0$;
\item $a{\cdot}[b,a,d]+[a,b{\cdot}a,d]=0$;
\item $a{\cdot}[b,c,d]+c{\cdot}[b,a,d]+[a,b{\cdot}c,d]+[c,b{\cdot}a,d]=0$;
\item $[a,b,c{\cdot}d]-[a,b,c]{\cdot}d-c{\cdot}[a,b,d]+[[a,b],c,d]=0$.
\end{enumerate}
\end{proposition}

\begin{proof} 1. The identity $[a,a,b]=0$ was actually proved in Proposition~\ref{p:Bol<=>alternative}.
\smallskip

2. Observe that
$$[a,b,c]+[b,a,c]=[a{+}b,a{+}b,c]-[a,a,c]-[b,b,c]=0.$$
by the preceding statement.
\smallskip

3. Applying the left-Bol identity, we obtain the equality
$$
\begin{aligned}
&a{\cdot}[b,a,d]+[a,b{\cdot}a,d]=a{\cdot}((b{\cdot}a){\cdot}d)-a{\cdot}(b{\cdot}(a{\cdot}d))+(a{\cdot}(b{\cdot}a)){\cdot}d-a{\cdot}((b{\cdot}a){\cdot}d)\\
&=-(a{\cdot}(b{\cdot}a)){\cdot}d+(a{\cdot}(b{\cdot}a)){\cdot}d=0.
\end{aligned}
$$
\smallskip

4. The preceding statement and Lemma~\ref{l:associator-additive} imply
$$
\begin{aligned}
0&=(a{+}c){\cdot}[b,a{+}c,d]+[a{+}c,b{\cdot}(a{+}c),d]\\
&=a{\cdot}[b,a,d]+a{\cdot}[b,c,d]+c{\cdot}[b,a,d]+c{\cdot}[b,c,d]+[a,b{\cdot}a,d]+[a,b{\cdot}c,d]+[c,b{\cdot}a,d]+[c,b{\cdot}c,d]\\
&=a{\cdot}[b,a,d]+[a,b{\cdot}a,d]+a{\cdot}[b,c,d]+c{\cdot}[b,a,d]+[a,b{\cdot}c,d]+[c,b{\cdot}a,d]+c{\cdot}[b,c,d]+[c,b{\cdot}c,d]\\
&=a{\cdot}[b,c,d]+c{\cdot}[b,a,d]+[a,b{\cdot}c,d]+[c,b{\cdot}a,d].
\end{aligned}
$$

5. Subtracting the Teichm\"uller identity 
$$a{\cdot}[b, c, d] - [a{\cdot}b, c, d] + [a, b{\cdot}c, d] - [a, b, c{\cdot}d] + [a, b, c]{\cdot}d = 0$$
from the identity proved in the preceding statement, and then apply the identity (2) we obtain
$$
\begin{aligned}
0&=a{\cdot}[b,c,d]+c{\cdot}[b,a,d]+[a,b{\cdot}c,d]+[c,b{\cdot}a,d]\\
&-a{\cdot}[b, c, d] + [a{\cdot}b, c, d] - [a, b{\cdot}c, d] + [a, b, c{\cdot}d] - [a, b, c]{\cdot}d\\
&=c{\cdot}[b,a,d]+[c,b{\cdot}a,d]+ [a{\cdot}b, c, d] + [a, b, c{\cdot}d] - [a, b, c]{\cdot}d\\
&= [a, b, c{\cdot}d]- [a, b, c]{\cdot}d-c{\cdot}[a,b,d]+ [a{\cdot}b, c, d]-[b{\cdot}a,c,d] \\
&= [a, b, c{\cdot}d]- [a, b, c]{\cdot}d-c{\cdot}[a,b,d]+ [[a,b]c, d].
\end{aligned}
$$
\end{proof}

\section{The diassociativity of alternative rings}

In this section we shall prove an important theorem of Artin on the diassociativity of alternative rings. A ring $R$ is \index{diassociative ring}\index{ring!diassociative}\defterm{diassociative} if every $2$-generated subring of $R$ is associative. Artin's Theorem implies the famous Artin--Zorn Theorem~\ref{t:Artin-Zorn} saying that every finite alternative ring is a field.

For a ring $R$, the set 
 $$\mathcal N(R)\defeq\{a\in R:\forall x,y\in R\;\;[a,x,y]=[x,a,y]=[x,y,a]=0\}$$is called the \index{nucleus}\index{ring!nucleus}\defterm{nucleus} of the ring.

\begin{theorem}[Artin, $<$1930]\label{t:Artin} For any elements $a,b$ of an alternative ring $R$, the subring of $R$, generated by the set $\mathcal N(R)\cup\{a,b\}$ is associative.
\end{theorem}

\begin{proof} Consider the set of monomials $M\defeq\bigcup_{i=1}^\infty M_i$ where $M_1\defeq\mathcal N(R)\cup\{a,b\}$ and $$M_{n}=\bigcup_{i=1}^{n-1} \{x\cdot y:x\in M_i\;\wedge\; y\in M_{n-i}\}.$$ Since $-1\in \mathcal N(R)$, the subring $S$ of $R$, generated by the set $\mathcal N(R)\cup\{a,b\}=M_1$ coincides with the set of finite sums of elements of the set $M$. Assuming that $S$ is not associative, we can find three elements $m_1,m_2,m_3\in M$ such that $[m_1,m_2,m_3]\ne 0$. Find numbers $d_1,d_2,d_3\in\IN$ such that $m_i\in M_{d_i}$ for $i\in\{1,2,3\}$. We can assume that the number $d\defeq d_1+d_2+d_3$ is the smallest possible and hence for every numbers $d_1',d_2',d_3'\in\IN$ with $d_1'+d_2'+d_3'<d$ and monomials $m'_1\in M_{d_1'}$, $m'_2\in M_{d_2}'$, $m_3'\in M_{d_3'}$ we have $[m_1',m_2',m_3']=0$. Also we can assume that the number $d_1+d_3$ is the smallest possible. 

We claim that $d_1=d_3=1$. To derive a contradiction, assume that $d_1\ne 1$. Then $m_1=m_1'{\cdot}m_1''$ for some monomials $m_1',m_1''$ such that $m_1'\in M_{d_1'}$ and $m_1''\in M_{d_1''}$ for some numbers $d_1',d_1''\in\IN$ with $d_1=d_1'+d_1''$. By Proposition~\ref{p:alternative9}(6), the alternative ring $R$ is left-Bol. By Propositions~\ref{p:alternative9}(1), \ref{p:left-Bol}(4) and the minimality of $d$,   
$$
\begin{aligned}
0&\ne [m_1,m_2,m_3]=[m_1'{\cdot}m_1'',m_2,m_3]=-[m_2,m_1'{\cdot}m_1'',m_3]\\
&=m_2{\cdot}[m_1',m_1'',m_2]+m_1''{\cdot}[m_1',m_2,m_3]+[m_1'',m_1'{\cdot}m_2,m_3]=[m_1'',m_1'{\cdot}m_2,m_3]
\end{aligned}
$$
which contradicts the minimality of $d_1+d_3$ (because $d_1''+d_3<d_1+d_3$ and $d_1''+d_1'+d_2+d_3=d_1+d_2+d_3=d$). This contradiction shows that $d_1=1$. 

Assuming that $d_3\ne 1$, we can write the monomial $m_3$ as $m_3=m_3'{\cdot}m_3''$ for some monomials $m_3',m_3''$ such that $m_3'\in M_{d_3'}$ and $m_3''\in M_{d_3''}$ for some numbers $d_3',d_3''\in\IN$ with $d_3=d_3'+d_3''$. By Propositions~\ref{p:alternative9}(1), \ref{p:left-Bol}(4) and the minimality of $d$,  
$$
\begin{aligned}
0&\ne [m_1,m_2,m_3]=-[m_2,m_1,m_3]=[m_2,m_3,m_1]=[m_2,m_3'{\cdot}m_3'',m_1]\\
&=-m_2{\cdot}[m_3',m_3'',m_1]-m_3''{\cdot}[m_3',m_2,m_1]-[m_3'',m_3'{\cdot}m_2,m_1]=-[m_3'',m_3'{\cdot}m_2,m_1]
\end{aligned}
$$
which contradicts the minimality of $d_1+d_3$ (because $d_3''+d_1<d_3+d_1$ and $d_3''+d_3'+d_2+d_1=d_1+d_2+d_3=d$). This contradiction shows that $d_3=1$. 

Therefore, $m_1,m_3\in M_1=\mathcal N(R)\cup\{a,b\}$. Taking into account that $[m_1,m_2,m_3]\ne 0$, we conclude that $m_1,m_3\subseteq M_1\setminus\mathcal N(R)\subseteq\{a,b\}$ and $m_1\ne m_3$. Since $0\ne [m_1,m_2,m_3]=-[m_3,m_2,m_1]$, we lose no generality assuming that $m_1=a$ and $m_2=b$. 

\begin{claim}\label{cl:Mi=*M1} For every number $i=\{2,\dots,d-1\}$, the following inclusion holds: $$M_i\subseteq M_{i-1}\cdot M_1\defeq\{x{\cdot}y:x\in M_{i-1},\;y\in M_1\}.$$
\end{claim}

\begin{proof} For $i=2$, the definition of the set $M_2$ ensures that $M_2=M_1\cdot M_1$. Assume that for some $i\in \{3,\dots,d-1\}$ and all $j\in\{2,\dots,i-1\}$ the inclusion $M_j\subseteq M_{j-1}{\cdot} M_1$ holds. By the definition of the set $M_i=\bigcup_{j=1}^{i-1}M_j\cdot M_{i-j}$, every element $m\in M_i$ can be written as $s{\cdot}t$ for some number $j\in\{1,\dots,i-1\}$ and some elements $s\in M_j$ and $t\in M_{i-j}$. If $j=i-1$, then $m=s{\cdot}t\in M_{i-1}\cdot M_1$ and we are done. If $j<i-1$, then the inductive assumption ensures that $t=x{\cdot}y$ for some elements $x\in M_{i-j-1}$ and $y\in M_1$. Then $m=s{\cdot}t=s{\cdot}(x{\cdot}y)$. By the minimality of the number $d$, the strict inequality $j+(i-j-1)+1=i<d$, implies $m=s{\cdot}(x{\cdot}y)=(s{\cdot}x){\cdot}y\in (M_j\cdot M_{i-j-1})\cdot M_1\subseteq M_{i-1}\cdot M_1$, witnessing that $M_i\subseteq M_{i-1}{\cdot} M_1$.
\end{proof}

Since $d_2<d$, Claim~\ref{cl:Mi=*M1} ensures that $m_2=s{\cdot}t$ for some $s\in M_{d_2-1}$ and $t\in M_1=\mathcal N(R)\cup\{a,b\}$. If $t\in\mathcal N(R)$, then by Proposition~\ref{p:left-Bol}(4) and the minimality of $d$,
$$
0\ne -[m_1,m_2,m_3]=-[a,s{\cdot}t,b]=a{\cdot}[s,t,b]+t{\cdot}[s,a,b]+[t,s{\cdot}a,b]=0+0+0,
$$
which is a contradiction showing that $t\in M_1\setminus\mathcal N(R)\subseteq\{a,b\}$. If $t=a$, then Proposition~\ref{p:left-Bol}(3) and the minimality of $d$  imply
$$0\ne[m_1,m_2,m_2]=[a,s{\cdot}t,b]=[a,s{\cdot}a,b]=-a{\cdot}[s,a,b]=0.
$$
If $t=b$, then Propositions~\ref{p:alternative9}(1), \ref{p:left-Bol}(3), and the minimality of $d$ imply
$$0\ne [m_1,m_2,m_2]=[a,s{\cdot}t,b]=[a,s{\cdot}b,b]=-[b,s{\cdot}b,a]=b{\cdot}[s,b,a]=0.
$$
In all cases, we obtain a contradiction showing that $[m_1,m_2,m_3]=0$ for all monomial $m_1,m_2,m_3\in M$, which implies that the subring generated by the set $\mathcal N(R)\cup\{a,b\}$ is associative.
\end{proof}

Artin's Theorem~\ref{t:Artin} implies the following generalization of the Wedderburn Theorem~\ref{t:Wedderburn-Witt}, due to Artin and \index[person]{Zorn}Zorn\footnote{{\bf Max August Zorn} (1906 -- 1993) was a German mathematician. He was an algebraist, group theorist, and numerical analyst. He is best known for Zorn's lemma, which was first postulated by Kazimierz Kuratowski in 1922, and then independently by Zorn in 1935.
\newline
Zorn was born in Krefeld, Germany. He attended the University of Hamburg. He received his PhD in April 1930 for a thesis on alternative algebras, written under supervision of Emil Artin. Max Zorn was appointed to an assistant position at the University of Halle. However, he did not have the opportunity to work there for long as he was forced to leave Germany in 1933 because of policies enacted by the Nazis. Zorn immigrated to the United States and was appointed a Sterling Fellow at Yale University. While at Yale, Zorn wrote his paper ``A Remark on Method in Transfinite Algebra'' that stated his Maximum Principle, later called Zorn's lemma. In 1936 he moved to UCLA and remained until 1946. While at UCLA Zorn revisited his study of alternative rings and proved the existence of the nilradical of certain alternative rings. In 1946 Zorn became a professor at Indiana University, where he taught until retiring in 1971.}.

\begin{theorem}[Artin--Zorn, 1930]\label{t:Artin-Zorn} Every finite alternative divisible ring $R$ is a field.
\end{theorem}

\begin{proof} We claim that the ring $R$ is associative. Given any elements $a,b,c\in R$, we should prove that $[a,b,c]=0$. By Artin's Theorem~\ref{t:Artin}, the subring $A$ of $R$, generated by the set $\{1,a,b\}$ is associative. Then $A^*\defeq A\setminus\{0\}$ is a finite semigroup with respect to the multiplication. The divisibility of the ring $R$ ensures that for every element $\alpha,\beta\in A^*$ the two-sided shift $s:A^*\to A^*$, $s:x\mapsto \alpha{\cdot}x{\cdot}\beta$, is injective and hence bijective, witnessing that the finite semigroup $A^*$ is a group and the finite ring $A$ is a corps. By the Wedderburn Theorem~\ref{t:Wedderburn-Witt}, the finite corps $A$ is a field. By Theorem~\ref{t:F*-cyclic}, the multiplicative group $A^*=A\setminus\{0\}$ of the field $A$ is cyclic and hence it is generated by a single element $g\in A^*$. Then the ring $A$ is also generated by the element $g$. By Artin's Theorem~\ref{t:Artin}, the subring $C$ of $R$, generated by the set $\{g,c\}$ is associative. Since $a,b,c\in S\cup\{c\}\subseteq C$, the associator $[a,b,c]$ is equal to zero, witnessing that the division ring $R$ is associative and hence $R$ is a corps. By Wedderburn Theorem~\ref{t:Wedderburn-Witt}, the finite corps $R$ is a field.
\end{proof}

\section{The alternativity of left-Bol rings}\label{s:Bol}

In this section we continue the study of left-Bol rings, started in Section~\ref{s:alt-identities}. The proof of the following theorem of Mikheev is copied from the monograph \cite[4.2]{Faulkner}. We recall that $[a,b]$ and $[a,b,c]$ denote the commutator and associator of the elements $a,b,c$ in a ring.

\begin{theorem}[Mikheev, 1969]\label{t:Mikheev} Every elements $a,b,c$ of a left-Bol ring $R$ satisfy the following identities:
\begin{enumerate}
\item $[a, b, [a, b, c]] =[a, b]\cdot[a, b, c]$;
\item $[a, b, [a, b]]^2 = 0$;
\item $[a, b, a]^4= 0$.
\end{enumerate}
\end{theorem}

\begin{proof} We shall exploit the natural representation of the ring $R$ in the endomorphism ring $\End(R)$ of the additive group $(R,+)$. Elements of the endomorphism ring $\End(R)$ are endomorphisms of the additive group $(R,+)$, i.e., functions $A:R\to R$ such that $A(x+y)=A(x)+A(y)$ for every $x,y\in R$. Since the group $(R,+)$ is commutative, for every endomorphisms $A,B\in \End(R)$, the function $A+B:x\mapsto A(x)+B(x)$, is an element of the endomorphism ring.  So, $\End(R)$ is a commutative group with respect to the operation of addition of endmorphisms. Endowed wih the operations of addition and composition of endomorphism, $\End(R)$ is an associative ring.

The distributivity of the multiplication over addition in the ring $R$ implies that for every $a\in R$, the left shift 
$$l_a:R\to R,\;\;l_a:x\mapsto a\cdot x,$$is well-defined endomorphisms of the group $(R,+)$, and  the map $$l_*:R\to\End(R),\;\;l_*:a\mapsto l_a,$$ is a homomorphism of the additive groups $(R,+)$ and $\End(R)$. Observe that every endomorphism $A\in l_*[R]$ is equal to the left shift $l_{A(1)}$.

Since the multiplicative magma $(R,\cdot)$ is not necessarily associative, the additive subgroup  $l_*[R]$ of $\End(R)$ is not necessarily a subring of the endomorphism ring $\End(R)$. Nonetheless,  the left-Bol identity implies the set $l_*[R]$ is closed under the binary operation $(a,b)\mapsto aba\defeq a\circ b\circ a$, called the \index{Jordan operation}\defterm{Jordan operation}. Indeed, for every elements $a,b,x\in R$, we have 
$$l_al_bl_a(x)=a{\cdot}(b{\cdot}(a{\cdot}x))=(a{\cdot}(b{\cdot}a)){\cdot x}=l_{a{\cdot}(b{\cdot}a)}(x).$$

Therefore, the subset $l_*[R]$ of the endomorphism ring $\End(R)$ is closed under the operations:
$$
\begin{aligned}
&aba\defeq a\circ b\circ a\\
&\{abc\}\defeq abc+cba=(a+c)b(a+c)-aba-cbc\\
&a^2\defeq a1a\\
&a{\dotplus} b\defeq ab+ba=(a+b)^2-a^2-b^2.
\end{aligned}
$$
Also for the polynomial 
$$f:\End(R)^3\to\End(R),\;\;f:(x,y,z)\mapsto z-xy,$$
the set $l_*[R]\subseteq\End[R]$ is closed under the following polynomials:
$$
\begin{aligned}
p(x, y, z)&=f(x, y, z)f(y, x, z) = z^2 - \{xyz\} + xy^2x,\\
q(x, y, z)&=f(x, y, z)xf(y, x, z) = zxz -((xyx)\dotplus z)+ x(yxy)x,\\
r(x, y, z,w)&=f(x, y, z)w + wf(y, x, z) = z{\dotplus}w - \{xyw\},\\
s(x, y, z,w)&=xf(y,x,z)w+wf(x,y,z)x = \{xzw\}- xyx{\dotplus}w.
\end{aligned}
$$
Given any elements $a,b,c\in R$, consider the elements $A\defeq f(l_a, l_b, l_{a{\cdot}b})$ and $B\defeq f(l_b, l_a, l_{a{\cdot}b})$ of the endomorphism ring $\End(R)$. Observe that $B(1)=l_{a{\cdot}b}(1)-l_bl_a(1)=a{\cdot}b-b{\cdot}a=[a,b]$. By  Proposition~\ref{p:left-Bol}(2) for every $x\in R$ we have $A(x) =[a, b, x]=-[b,a,x]$, which implies $A=f(l_a,l_b,l_{a{\cdot}b})=-f(l_b,l_a,l_{b{\cdot}a})$ and $A(1) = 0$, Then 
$$A+B = -f(l_b, l_a, l_{b{\cdot}a}) + B = l_{[a,b]}.$$
We see that
$$l_{[a,b]}A - A^2 =BA= p(l_b, l_a, l_{a{\cdot}b}) = l_{BA(1)} = 0,$$
giving the identity (1). 

Also, 
$$AB=p(l_a,l_b,l_{ab})=l_{AB(1)}=l_{[a,b,[a,b]]}$$
so $BA = 0$ implies
$[a, b, [a, b]]^2 = ABAB(1) = 0$,
giving the identity (2). 

Finally, consider the element $d\defeq [a, b, a]= A(a)$ and observe that $B(d) =BA( a)=
0$. Thus,
\begin{itemize}
\item[(i)]  $Bl_d + l_dA= r( l_b, l_a, l_{a{\cdot}b}, l_d) = l_{(Bl_d+l_dA)(1)} = l_0=0$,
\item[(ii)]  $l_dABl_d= l_dp(l_a, l_b, l_{a{\cdot}b}) l_d = 0$,
\item[(iii)] $l_dAl_aBl_d= l_d q(l_a, l_b, l_{a{\cdot}b}) l_d = 0$,
\item[(iv)] $l_aBl_d + l_dAl_a = s(l_a, l_b, l_{a{\cdot}b},l_d)= l_{l_dAl_a(1)} = l_{d^2} = l_d^2$.
\end{itemize}
We now see that
$$\begin{aligned}
l^3_dA&=-l_d^2Bl_d&&\mbox{by (i)}\\
&=-(l_aBl_d + l_dAl_a)Bl_d&&\mbox{by (iv)}\\
&=l_al_dABl_d - l_dAl_aBl_d&&\mbox{by (i)}\\
&=0&&\mbox{by (ii) and (iii)}.
\end{aligned}
$$
Thus, $d^4 = l^3_dA(a) = 0$, giving Mikheev's identity (3).
\end{proof}

A ring $R$ is defined to be \index{reduced ring}\index{ring!reduced}\defterm{reduced} if $$\forall a\in R\;(a=0\;\Leftrightarrow\;a{\cdot}a=0).$$

\begin{theorem}\label{t:alternative<=>left-Bol} A reduced ring  is alternative if and only if it is left-Bol.
\end{theorem}

\begin{proof} The ``only if'' part of this characterization has been proved in Proposition~\ref{p:alternative9}(6).

To prove the ``if'' part, assume that a reduced ring $R$ is left-Bol. Fix any elements $a,b\in R$.  Proposition~\ref{p:left-Bol}(1) ensures that $[a,a,b]=0$. By Theorem~\ref{t:left-Bol=>power-associative}, the left-Bol ring $R$ is power-associative. By Theorem~\ref{t:Mikheev}(3), $[a,b,a]^4=0$. Since the ring $R$ is reduced and power-associative, the equality $[a,b,a]^4=0$ implies $[a,b,a]^2=0$ and $[a,b,a]=0$. Applying Proposition~\ref{p:left-Bol}(2), we conclude that
$$[b,a,a]=-[a,b,a]=0,$$witnessing that the left-Bol ring $R$ is alternative. 
\end{proof}

\begin{corollary} A reduced ring is alternative if and only if it is right-Bol.
\end{corollary}

\begin{proof} The ``only if'' part of this characterization has been proved in Proposition~\ref{p:alternative9}(7). To prove the ``if'' part, assume that a reduced ring $R$ is right-Bol. Denote by $R^{op}$ the additive group $(R,+)$ of the ring $R$, endowed with the ``opposite'' multiplication $$*:R\times R\to R,\quad *:(x,y)\mapsto x*y\defeq y\cdot x.$$ The right-Bol identity of the ring $R$ implies the left-Bol identity of the opposite ring $R^{op}$. Since $R$ is reduced, so is its opposite ring $R^{op}$. Applying Theorem~\ref{t:alternative<=>left-Bol} to the reduced left-Bol ring $R^{op}$, we conclude that $R^{op}$ is alternative. Then for every elements $a,b\in R$, we have the equalities
$$(a\cdot a) \cdot b=b*(a*a)=(b*a)*a=a\cdot(a\cdot b)$$
and
$$b\cdot(a\cdot a)=(a*a)*b=a*(a*b)=(b\cdot a)\cdot a,$$
witnessing that the ring $R$ is alternative.
\end{proof}

\section{The alternativity of division rings}

In this section we discuss a \index[person]{Hua}Hua\footnote{{\bf Hua Luogeng or Hua Loo-Keng} (1910 -- 1985) was a Chinese mathematician and politician famous for his important contributions to number theory and for his role as the leader of mathematics research and education in the People's Republic of China. He was largely responsible for identifying and nurturing the renowned mathematician Chen Jingrun who proved Chen's theorem, the best known result on the Goldbach conjecture. In addition, Hua's later work on mathematical optimization and operations research made an enormous impact on China's economy. He was elected a foreign associate of the US National Academy of Sciences in 1982. He was elected a member of the standing Committee of the first to sixth National people's Congress, Vice-chairman of the sixth National Committee of the Chinese People's Political Consultative Conference (April 1985) and vice-chairman of the China Democratic League (1979). He joined the Chinese Communist Party in 1979.
\newline 
Hua did not receive a formal university education. Although awarded several honorary PhDs, he never got a formal degree from any university. In fact, his formal education only consisted of six years of primary school and three years of secondary school. For that reason, Xiong Qinglai, after reading one of Hua's early papers, was amazed by Hua's mathematical talent, and in 1931 Xiong invited him to study mathematics at Tsinghua University.} identity in associative rings and its application to characterizing alternative division rings.

 Let us recall that an element $a$ of a ring $R$ is \index{invertible element}\index{element!invertible}\defterm{invertible} if there exists an element $a^{-1}\in R$ such that $a\cdot a^{-1}=1=a^{-1}\cdot a$. 
 
\begin{exercise} Show that any invertible element in an associative ring has a unique inverse.
\end{exercise}

\begin{lemma}[Hua, 1949]\label{l:Hua} Let $a,b$ be two  invertible elements in an associative ring $R$. If the element $a-b^{-1}$ is invertible, then $a^{-1}-(a-b^{-1})^{-1}$ is invertible and
$$a-(a^{-1}-(a-b^{-1})^{-1})^{-1}=aba.$$
\end{lemma}

\begin{proof} Since the ring $R$ is associative, its set of invertible elements is a group with respect to the multiplication operation. Then $(x\cdot y)^{-1}=y^{-1}\cdot x^{-1}$ for any invertible elements $x,y$ in $R$. By the associativity and distributivity of the multiplication over addition in the ring $R$, we obtain  
$$
\begin{aligned}
a^{-1}-(a-b^{-1})^{-1}&=(a^{-1}{\cdot}(a-b^{-1})-1)\cdot (a-b^{-1})^{-1}\\
&=-a^{-1}{\cdot}b^{-1}{\cdot}(a-b^{-1})^{-1}=-((a-b^{-1}){\cdot}b{\cdot}a)^{-1}=(a-aba)^{-1},
\end{aligned}
$$
which implies  \index{identity!Hua}\index{Hua's identity}\defterm{Hua's identity}  
$a-(a^{-1}-(a-b^{-1})^{-1})^{-1}=aba$.
\end{proof}

\begin{definition} A ring $R$ is defined to be 
\begin{itemize}
\item \index{ring!left-inversive}\defterm{left-inversive} if for every $a\in R^*\defeq R\setminus\{0\}$ there exists $b\in R$ such that $b{\cdot}(a{\cdot}x)=x$ for every $x\in R$;
\item \index{ring!right-inversive}\defterm{right-inversive} if for every $a\in R^*\defeq R\setminus\{0\}$ there exists $b\in R$ such that $x=(x{\cdot}a){\cdot}b$ for every $x\in R$.
\end{itemize}
\end{definition}
Observe that a division ring $R$ is left-inversive (resp. right-inversive) if and only if so is the multiplicative loop $(R^*,\cdot)$.

Let us recall that a ring $R$ is called \index{division ring}\index{ring!division}\defterm{division} if for any elements $a\in R^*\defeq R\setminus\{0\}$ and $b\in R$, there exist unique elements $b/a$ and $a\backslash b$ in $R$ such that
$$(b/a){\cdot} a=b=a{\cdot} (a\backslash b).$$ Observe that $0/a=0\ne a$ and hence $a{\cdot}a\ne (0/a){\cdot}a=0$, witnessing that every division ring is reduced. 

\begin{theorem}\label{t:alternative<=>left-inverse} For a division ring $R$, the following conditions are equivalent:
\begin{enumerate}
\item the ring $R$ is alternative;
\item the ring $R$ is left-Bol;
\item the ring $R$ is left-inversive;
\item for any elements $a,b\in R^*$ with $a\backslash 1+b\backslash 1\ne 0$, there exists an element $c\in R^*$ such that\newline $a\backslash x+b\backslash x=c\backslash x$ for all $x\in R$.
\end{enumerate}
\end{theorem}

\begin{proof} The equivalence $(1)\Leftrightarrow(2)$ follows from Theorem~\ref{t:alternative<=>left-Bol} and the reducibility of divisible rings.
\smallskip

$(2)\Ra(3)$ Assume that the ring $R$ is left-Bol. Given any element $a\in R^*\defeq R\setminus\{0\}$, we should find an element $b\in R^*$ such that $b{\cdot}(a{\cdot}x)=x$ for all $x\in R$. Since $(R^*,\cdot)$ is a loop, there exists an element $a^{-1}\in R^*$ such that $a^{-1}{\cdot}a=1$. For every $x\in R$, the left-Bol identity implies
$$a{\cdot}(a^{-1}{\cdot}(a{\cdot}x))=(a{\cdot}(a^{-1}{\cdot}a)){\cdot}x=(a{\cdot}1){\cdot}x=a{\cdot}x$$and hence $a^{-1}{\cdot}(a{\cdot}x)=x$, by the injectivity of the left shift $l_{a}:R\to R$.
\smallskip

$(3)\Ra(4)$ Assume that the ring $R$ has the left-inverse property. Then for every element $a\in R^*\defeq R\setminus\{0\}$ there exists an element $a^{-1}\in R$ such that $a^{-1}{\cdot}(a{\cdot}x)=x$. In particular, $a^{-1}{\cdot}a=a^{-1}{\cdot}(a{\cdot}1)=1$ and $a^{-1}{\cdot}(a{\cdot}a^{-1})=a^{-1}=a^{-1}{\cdot}1$, which implies $a\cdot a^{-1}=1$. So, $a^{-1}$ is a two-sided inverse to $a$ in the loop $(R^*,\cdot)$, and hence $a=(a^{-1})^{-1}$. For every $x\in R$, the equality  $a{\cdot}(a^{-1}{\cdot}x)=x$ implies $a\backslash x=a^{-1}{\cdot}x$.

Now take any elements $a,b\in R^*$ with $a\backslash 1+b\backslash 1\ne 0$ and find an element $c\in R^*$ such that $c\backslash 1=a\backslash 1+b\backslash 1$. Then for every $x\in R$, we have 
$$a\backslash x+b\backslash x=a^{-1}{\cdot}x+b^{-1}{\cdot}x=(a^{-1}+b^{-1}){\cdot}x=c^{-1}\cdot x=c\backslash x.$$
\smallskip

$(4)\Ra(2)$ Assume that the condition (4) holds. For every $a\in R$, consider the left shift $l_a:R\to R$, $l_a:x\mapsto a{\cdot}x$. Then $l_a$ is an element of the endomorphism ring $\End(R)$ of the additive group $(R,+)$. The distribution of the group $(R,+)$ implies that the map $l_*:R\to\End(R)$, $l_*:a\mapsto l_a$, is a homomorphism of the additive groups of the rings $R$ and $\End(R)$. The left distributivity of $R$ implies that for every $a\in R^*$ and $x\in R$ we have $0=(a+(-a))\cdot x=l_a(x)+l_{-a}(x)$ and hence $-l_a=l_{-a}$. This observation  and the condition (4) imply that the set $$G\defeq \{0\}\cup\{l_a^{-1}:a\in R^*\}=\{0\}\cup\{a^{-1}:a\in l_*[R^*]\}$$ is a subgroup of the additive group of the endomorphism ring $\End(R)$. Observe that $G$ contains the identity element $\mathbf 1=l_1$ of the endomorphism ring $\End(R)$. We claim that $l_a^{-1}\in l_*[R^*]$ for every nonzero element $a\in R^*$. If $a=1$, then $l_a^{-1}=l_1^{-1}=l_1\in l_*[R^*]$. 
If $a\ne 1$, then $0\ne \mathbf 1-l_a^{-1}=l_1^{-1}-l_a^{-1}\in G$, by the condition (4). Then $(\mathbf 1-l_a^{-1})^{-1}\in l_*[R]$ and $b\defeq \mathbf 1-(\mathbf 1-l_a^{-1})^{-1}\in l_*[R]$ because the set $l_*[R]$ is closed under addition.  Applying Hua's identity from Lemma~\ref{l:Hua}, we conclude that
$$l_a=\mathbf 1l_a\mathbf 1=\mathbf 1-(\mathbf 1^{-1}-(\mathbf 1-l^{-1}_a)^{-1})^{-1}=\mathbf 1-b^{-1}\in G$$because $G$ is a group. Therefore, $l_a=l_{a'}^{-1}$ for some element $a'\in R^*$. Since $l_{a'}l_a=\mathbf 1$,  $a'{\cdot}a=a'{\cdot}(a{\cdot}1)=l_{a'}l_a(1)=\mathbf 1(1)=1$ and hence $a'=a^{-1}$. For every $x\in R$, the equality $l_{a'}l_a=\mathbf 1$ implies 
$$(a^{-1}\cdot a)\cdot x=1\cdot x=x=\mathbf 1(x)=l_{a^{-1}}l_a(x)=a^{-1}{\cdot}(a{\cdot}x),$$witnessing that the ring $R$ is left-inversive and moreover $l_a^{-1}=l_{a^{-1}}$ for every $a\in R^*$. 

Next, we show that the ring $R$ is left-Bol. Given any elements $a,b,x\in R$, we need to show that
$$a{\cdot}(b{\cdot}(a{\cdot}x))=(a{\cdot}(b{\cdot}a)){\cdot}x.$$
If $0\in\{a,b\}$, then 
$$a{\cdot}(b{\cdot}(a{\cdot}x))=0=(a{\cdot}(b{\cdot}a)){\cdot}x$$
and we are done. So, assume that $a,b\in R^*$. If $b=a^{-1}$, then 
$$a{\cdot}(b{\cdot}(a{\cdot}x))=a{\cdot}(a^{-1}{\cdot}(a{\cdot}x))=a{\cdot}x=(a{\cdot}(a^{-1}{\cdot}a)){\cdot}x=(a{\cdot}(b{\cdot}a)){\cdot}x$$by the left-inversivity. So, assume that $b\ne a^{-1}$. In this case, the Hua identity \ref{l:Hua} implies 
$$l_al_bl_a=l_a-(l_a^{-1}-(l_a-l_b^{-1})^{-1})^{-1}\in l_*[R^*]$$ and hence $l_al_bl_a=l_d$ for some element $d\in R^*$ such that $$d=l_d(1)=l_al_bl_a(1)=a{\cdot}(b{\cdot}(a{\cdot 1}))=a{\cdot}(b{\cdot}a).$$ For every $a,b\in R^*$ and $x\in R$, the equality $l_al_bl_a=l_d$ implies the identity
$$a{\cdot}(b{\cdot}(a{\cdot}x))=l_al_bl_a(x)=l_d(x)=d{\cdot}x=(a{\cdot}(b{\cdot}a)){\cdot}x,$$
witnessing that the ring $R$ is left-Bol. 
\end{proof}

Aplying Theorem~\ref{t:alternative<=>left-inverse} to the opposite ring of a division ring, we obtain the following ``right'' characterization of alternative division rings.

\begin{theorem}\label{t:alternative<=>right-Bol} For a division ring $R$, the following conditions are equivalent:
\begin{enumerate}
\item the ring $R$ is alternative;
\item the ring $R$ is right-Bol;
\item the ring $R$ is right-inversive;
\item for any distinct elements $a,b\in R^*$, there exists an element $c\in R^*$ such that\newline $x/a-x/b=x/c$ for all $x\in R$.
\end{enumerate}
\end{theorem}

For a ring $R$, the set 
$$\mathcal N(R)\defeq\{a\in R:\forall x,y\in R\;[a,x,y]=[x,a,y]=[x,y,a]=0\}$$ is called the \defterm{nucleus} of $R$. Here $[x,y,z]\defeq(x{\cdot}y){\cdot}z-x{\cdot}(y{\cdot}z)$ is the \defterm{associator} of elements $x,y,z\in R$. On the other hand, $[x,y]\defeq x{\cdot}y-y{\cdot}x$ is the \defterm{commutator} of $x,y$.
The set $$\mathcal Z(R)\defeq\{a\in \mathcal N(R):\forall x\in R\;\;[a,x]=0\}$$is the \defterm{centre} of the ring $R$. 

\begin{exercise} Check that the nucleus (and the centre) of a ring $R$ is an associative (and commutative) subring of $R$.
\end{exercise}

\begin{proposition}\label{p:NZ-corps-field} The nucleus $\mathcal N(R)$ of a division ring $R$ is a corps and the center $\mathcal Z(R)$ of $R$ is a field.
\end{proposition}

\begin{proof} Using the distributivity of the multiplication over addition in the ring $R$, it is easy to show that the nucleus $\mathcal N(R)$ is an associative subring of $R$. To show that $\mathcal N(R)$ is a corps, it suffices to check that every non-zero element $a\in \mathcal N(R)$ has a multiplicative inverse $a^{-1}$ in $\mathcal N(R)$. Since $R$ is a division ring, there exists a unique element $a^{-1}\in R$ such that $a{\cdot}a^{-1}=1$. Since $a\in\mathcal N(R)$, for every $x,y\in R\setminus\{0\}$ we have $[a,x,y]=0$. By Theorems~\ref{t:alternative<=>left-inverse} and \ref{t:alternative<=>right-Bol}, the alternative division ring $R$ is Bol and hence its multiplicative loop $R^*\defeq R\setminus\{0\}$ is Bol and Moufang, by Theorem~\ref{t:Bol<=>Moufang}. By Moufang Theorem~\ref{t:Moufang-xyz}, $[a,x,y]=0$ implies $[a^{-1},x,y]=[x,a^{-1},y]=[x,y,a^{-1}]=0$, which ensures that $a^{-1}\in \mathcal N(R)$, witnessing that $\mathcal N(R)$ is an associative division ring and hence a corps.

Next, we show that the centre $\mathcal Z(R)$ of $R$ is a field. Since $\mathcal Z(R)$ is a commutative subring of the corps $\mathcal N(R)$, it suffices to check that for every $z\in\mathcal Z(R)\setminus\{0\}$ its multiplicative inverse $z^{-1}\in\mathcal N(R)$ belongs to $\mathcal Z(R)$. Since $R$ is a division ring and $z\ne 0$, for every $x\in R$, there exists a unique element $y\in R$ such that $z{\cdot}y=x$. It follows from $z\in \mathcal Z(R)$ that $y{\cdot}z=z{\cdot}y=x$. Observe that $z{\cdot}(z^{-1}{\cdot}x)=(z{\cdot}z^{-1}){\cdot}x=1\cdot x=x=z{\cdot}y$ and hence $z^{-1}{\cdot}x=y$. On the other hand, $(x{\cdot}z^{-1}){\cdot} z=x{\cdot}(z^{-1}{\cdot}z)=x{\cdot}1=x=y{\cdot}z$ implies $x{\cdot}z^{-1}=y$. Therefore, $z^{-1}{\cdot}x=y=x{\cdot}z^{-1}$ and $z^{-1}\in\mathcal Z(R)$.
\end{proof}

\section{The associativity of commutative alternative rings}

In this section we shall prove that commutative division alternative rings are associative.

We shall exploit the following identity holding in every ring.

\begin{lemma}\label{l:com-alt} Every elements $a,b,c,d$ of a ring $R$ satisfy the following identity:
$$[a,b]{\cdot}c-[a,b{\cdot}c]+b{\cdot}[a,c]=[a,b,c]-[b,a,c]+[b,c,a].$$
\end{lemma}

\begin{proof} Indeed,
$$\begin{aligned}
[a,b]{\cdot}c-[a,b{\cdot}c]+b{\cdot}[a,c]&=\big((a{\cdot}b){\cdot}c-(b{\cdot}a){\cdot}c\big)-\big(a{\cdot}(b{\cdot}c)-(b{\cdot} c){\cdot} a\big)+\big(b{\cdot}(a{\cdot}c)-b{\cdot}(c{\cdot}a)\big)\\
&=\big((a{\cdot}b){\cdot}c-a{\cdot}(b{\cdot}c)\big)-\big((b{\cdot}a){\cdot}c-b{\cdot}(a{\cdot}c)\big)+\big((b{\cdot}c){\cdot}a-b{\cdot}(c{\cdot}a)\big)\\
&=[a,b,c]-[b,a,c]+[b,c,a].
\end{aligned}
$$
\end{proof}

\begin{theorem}\label{t:3[a,b,c]=0} Every elements $a,b,c$ of a comutative alternative ring $R$ satisfy the identities
$$3[a,b,c]=[a,b,c]^2=0.$$
\end{theorem}

\begin{proof} Since $R$ is commutative, Lemma~\ref{l:com-alt} and Proposition~\ref{p:alternative9}(1) imply
$$\begin{aligned}
0&=[a,b]{\cdot}c-[a,b{\cdot}c]+b{\cdot}[a,c]=[a,b,c]-[b,a,c]+[b,c,a]=3[a,b,c].
\end{aligned}
$$
By the commutativity of $R$ and Proposition~\ref{p:left-Bol}(5), for every $a,b\in R$, the function $A:R\to R$, $A:x\mapsto [a,b,x]$, is a derivation in the sense that 
$$A(x{\cdot}y)=A(x){\cdot}y+x{\cdot}A(y)$$
for every $x,y\in R$. Moreover, Theorem~\ref{t:Mikheev}(1) ensures that $A\cdot A=0$.

Then $$
\begin{aligned}
0&=A(A(c{\cdot}c))=A(A(c){\cdot}c+c{\cdot}A(c))=A(2c{\cdot}A(c))=2A(c{\cdot}A(c))\\
&=2(A(c){\cdot}A(c)+c\cdot A(A(c)))=2A(c)^2=2[a,b,c]^2
\end{aligned}
$$
and finally,
$$[a,b,c]^2=3[a,b,c]{\cdot}[a,b,c]-2[a,b,c]^2=0.$$
\end{proof}

Theorem~\ref{t:3[a,b,c]=0} implies the following corollary.

\begin{corollary}\label{c:reduced+com+alt=>ass} Every reduced commutative alternative ring is associative.
\end{corollary}

Since every division ring is reduced, Corollary~\ref{c:reduced+com+alt=>ass} implies another corollary, which will be applied in the proof of Theorem~\ref{t:commutative-dot<=>}.

\begin{corollary}\label{c:div+com+alt=>ass} Every divisible commutative alternative ring is associative.
\end{corollary}

\chapter{Involutive rings and the Cayley--Dickson process}\label{ch:Cayley-Dickson}

In this chater we study involutive rings, define and analyzse their Cayley--Dickson extensions, and establish the structure  and isotopies of alternative division rings. The principal results are difficult Skornyakov--Bruck--Kleinfeld Theorem~\ref{t:SBK} on the octonion structure of non-associative alternative division rings, and Schafer Theorem~\ref{t:Schafer} on isotopies of alternative division rings.

\section{Involutive rings}

In this section we introduce and study involutive rings, i.e., rings endowed with a central involution.

\begin{definition} An \index{involution}\defterm{involution} on a ring $R$ is a unary operation $j:R\to R$, $j:x\mapsto \bar x$, satisfying the following properties:
\begin{itemize}
\item $\overline{x+y}=\overline{x}+\overline y$;
\item $\overline{x{\cdot}y}=\overline y{\cdot}\overline x$;
\item $\overline{\overline x}=x$
\end{itemize}
for all $x,y\in R$. 

If $\bar x=x$ for all $x\in R$, then the involution $j$ is called \index{trivial involution}\index{involution!trivial}\defterm{trivial}.
\end{definition}

Therefore, an involution is an anti-isomorphism of $R$ of order $\le 2$.

\begin{definition} An involution $j:R\to R$, $j:x\mapsto \bar x$, is called \index{central involution}\index{involution!central}\defterm{central} if every element $x\in R$ its \index{norm}\defterm{norm} $\mathsf n(x)\defeq x{\cdot}\bar x$ and \index{trace}\defterm{trace} $\mathsf t(x)\defeq x+\bar x$ both belong to the  \index{centre}\defterm{centre}$$\mathcal Z(R)\defeq\{z\in R:\forall x,y\in R\;\;x{\cdot}(y{\cdot}z)=(x{\cdot}y){\cdot}z=(y{\cdot}x){\cdot}z=y{\cdot}(x{\cdot}z)\}$$ of the ring $R$.
\end{definition}

\begin{definition} An \index{involutive ring}\index{ring!involutive}\defterm{involutive ring} is ring $R$ endowed with a central involution $j:R\to R$, $j:x\mapsto\bar x$. 
\end{definition}

\begin{definition} For an involutive ring $R$, the set $$\Re(R)\defeq\{x\in \mathcal Z(R):\bar x=x\}$$ is called the \index{real axis}\index{involutive ring!real axis}\defterm{real axis} of $R$. Elements of the real axis are called \defterm{real elements} of the involutive ring $R$.
\end{definition}

\begin{proposition}\label{p:CD-real-axis} For any involutive ring $R$, its real axis $\Re(R)$ is a commutative associative subring of $R$. If $R$ is an alternative division ring, then its  real axis $\Re(R)$ is a field.
\end{proposition}

\begin{proof}  Since every anti-isomorphism of the ring $R$ preserves $0$ and $1$, the set $\Re(R)$ contains $0$ and $1$. For every $x,y\in \Re(R)\subseteq \mathcal Z(R)$ we have $\overline{x-y}=\bar x-\bar y=x-y$ and $\overline{x{\cdot}y}=\bar y{\cdot}\bar x=y{\cdot} x=x{\cdot} y$, witnessing that $\Re(R)$ is a subring of the ring $R$. Since $\Re(R)\subseteq\mathcal Z(R)$, the ring $\Re(R)$ is commutative and associative.  

Now assume that the ring $R$ is a division ring. By Theorem~\ref{p:NZ-corps-field}, the centre $\mathcal Z(R)$ of $R$ is a field. To see that $\Re(R)$ is a field, take any non-zero element $x\in \Re(R)\subseteq\mathcal Z(R)$. Since $\mathcal Z(R)$ is a field, there exists a unique element $x^{-1}\in \mathcal Z(R)$ such that $x{\cdot}x^{-1}=1$. Applying to the latter equality the central involution of $R$, we obtain $$x^{-1}{\cdot}x=1=\bar 1=\overline{x{\cdot}x^{-1}}=\overline{x^{-1}}\cdot\bar x=\overline{x^{-1}}\cdot x$$and hence $x^{-1}=\overline{x^{-1}}$, by the cancellativity of the multiplicative group of the field $\mathcal Z(R)$. Therefore $x^{-1}\in \Re(R)$, witnessing that the commutative associative ring $\Re(R)$ is invertible and hence $\Re(R)$ is a field.
\end{proof}

\begin{proposition}\label{p:tn-bar} For any element $x$ of an involutive ring $R$, 
$$\overline{\mathsf t(x)}=\mathsf t(\bar x)=\mathsf t(x)\in \Re(R)\quad\mbox{and}\quad\overline{\mathsf n(x)}=\mathsf n(\bar x)=\mathsf n(x)\in\Re(R).
$$
\end{proposition}

\begin{proof} Observe that
$$\overline{\mathsf t(x)}=\overline{x+\bar x}=\bar x+\bar{\bar x}=\mathsf t(\bar x)=\bar x+x=x+\bar x=\mathsf t(x)\in\mathcal Z(R)$$ and hence $\mathsf t(x)\in \Re(R)$. The identity $x{\cdot}(x+\bar x)=x{\cdot}\mathsf t(x) = \mathsf t(x){\cdot}x=(x+\bar x){\cdot} x$ implies $\mathsf n(x) = x{\cdot}\bar x =\bar x{\cdot }x= \mathsf n(\bar x)$. Also $\overline{\mathsf n(x)}=\overline {x{\cdot}\bar x}=\bar{\bar x}\cdot\bar x=x{\cdot}\bar x=\mathsf n(x)\in\mathcal Z(R)$ implies $\mathsf n(x)\in\Re(R)$.
\end{proof}



\begin{proposition}\label{p:CD-1gen} Let $R$ be an involutive ring and $Z$ be a subring of $\mathcal Z(R)$ such that $\Re(R)\subseteq Z$. Then for every $x\in R$, the set $Z+x{\cdot}Z\defeq\{a+x{\cdot}b:a,b\in Z\}$ is a commutative associative subring of the ring $R$, which implies that the involutive ring $R$ is $\IN$-associative. If for every $x\in R\setminus\{0\}$ the element  $\mathsf n(x)\defeq x{\cdot}\bar x$ is invertible in the ring $Z$, then the ring $Z+x{\cdot}Z$ is a field and the ring $R$ is $\IZ$-associative.
\end{proposition}

\begin{proof} Observe that $x^2=x{\cdot}\mathsf t(x)-\mathsf n(x)$. Taking into account that $\mathsf t(x)$ and $\mathsf n(x)$ belong to the centre of the ring $R$, we conclude that
$$
\begin{aligned}
x^2{\cdot}x&=(x{\cdot}\mathsf t(x)-\mathsf n(x)){\cdot}x=(x{\cdot}\mathsf t(x)){\cdot}x-\mathsf n(x){\cdot}x\\
&=x{\cdot}(\mathsf t(x){\cdot}x)-x{\cdot}\mathsf n(x)=
x{\cdot}(x{\cdot}\mathsf t(x)-\mathsf n(x))=x{\cdot}x^2,
\end{aligned}
$$
which implies that the set $Z+x{\cdot}Z$ is a commutative associative subring of $R$. 

Let $x^0\defeq 1$ and $x^{k+1}=x^k{\cdot}x$ for every $k\in\w$. Since $Z+x{\cdot}Z$ is a subring of $R$, the sequence $(x^k)_{k\in\w}$ belongs to $Z+x{\cdot}Z$. The associativity of $Z+x{\cdot}Z$ implies that $x^{m+k}=x^m{\cdot}x^k$ for all $m,k\in\w$, witnessing that the ring $R$ is $\IN$-associative.

Now assume that the element $\mathsf n(x)\defeq x{\cdot}\bar x$ has inverse $\mathsf n(x)^{-1}$ in the ring $Z$. Then consider the element $x^{-1}\defeq \mathsf n(x)^{-1}{\cdot}\mathsf t(x)-\mathsf n(x)^{-1}{\cdot}x\in Z+x{\cdot}Z$ and observe that 
$$
\begin{aligned}
x^{-1}{\cdot}x&=(\mathsf n(x)^{-1}{\cdot}\mathsf t(x)-\mathsf n(x)^{-1}{\cdot}x){\cdot}x=\mathsf n(x)^{-1}{\cdot}\mathsf t(x){\cdot}x-\mathsf n(x)^{-1}{\cdot}x^2\\
&=\mathsf n(x)^{-1}{\cdot}\mathsf t(x){\cdot}x-\mathsf n(x)^{-1}{\cdot}(x{\cdot}\mathsf t(x)-\mathsf n(x))\\
&=\mathsf n(x)^{-1}{\cdot}\mathsf t(x){\cdot}x-\mathsf n(x)^{-1}{\cdot}x{\cdot}\mathsf t(x)+\mathsf n(x)^{-1}{\cdot}\mathsf n(x)=1.\\
\end{aligned}
$$
Therefore, $x^{-1}\in Z+x{\cdot}Z$ is the multiplicative inverse to $x$ in the commutative associative ring $Z+x{\cdot}Z$. Then there exists a sequence $(x^n)_{n\in\w}$ in $Z+x{\cdot}Z$ such that $x^1=x$ and $x^m{\cdot}x^k=x^{m+k}$ for all $m,k\in\IZ$, witnessing that the ring $R$ is $\IZ$-associative.
\end{proof}

\begin{definition} Two involutive rings $X,Y$ are \index{involutive rings!isomorphic}\defterm{isomorphic} if there exists a bijective function $F:X\to Y$ such that 
\begin{itemize}
\item $F(x+y)=F(x)+F(y)$;
\item $F(x\cdot y)=F(x)\cdot F(y)$;
\item $\overline{F(x)}=F(\overline x)$;
\end{itemize}
for all $x,y\in X$
\end{definition}

\section{The Cayley--Dickson process}

In this section we describe the Cayley--Dickson process  extending involutive rings to larger involutive rings with worse algebraic properties. 

\begin{theorem}[Zorn, 1930]\label{t:Cayley-Dickson} Suppose that $R$ is an involutive ring and $\mathfrak z\in \Re(R)$ is a fixed real element of $R$. Then set 
$R'\defeq CD(R,\mathfrak z) = R\times R$ endowed with the operations of sum, product and involution
$$
\begin{aligned}
&(a,\alpha)+(b,\beta)\defeq (a+b,\alpha+\beta)\\
&(a,\alpha)\cdot(b,\beta)\defeq(a{\cdot}b+{\mathfrak z}{\cdot}\beta{\cdot}\bar\alpha,\bar a{\cdot}\beta+b{\cdot}\alpha)\\
&\overline{(a,\alpha)}\defeq(\bar a,-\alpha),
\end{aligned}
$$
is an involutive ring. Moreover, the ring $R'$ is
\begin{enumerate}
\item alternative if and only if $R$ is associative,
\item associative if and only if $R$ is associative and commutative,
\item commutative if and only if the involution of $R$ is trivial (i.e., $\bar x=x$ for all $x\in R$),
\item invertible if $R$ is an alternative division ring and $\mathfrak z\in\Re(R)\setminus\{\mathsf n(x):x\in R\}$.
\end{enumerate}
\end{theorem}

\begin{proof} The definition of the addition and multiplication in $R'$ shows that $R$ is a ring (i.e., associative-plus distributive biloop). It remains to show that the function $$j':R'\to R',\;\;j':(a,\alpha)\mapsto \overline{(a,\alpha)}\defeq (\bar a,-\alpha),$$ is a central involution of the ring $R'$.  For every $(a,\alpha)\in R'$, we have
$$\overline{\overline{(a,\alpha)}}=\overline{(\bar a,-\alpha)}=(\bar{\bar a},-(-\alpha))=(a,\alpha),$$
so $j'$ is an involutive map on $R'$. Next, we show that $j'$ is an anti-isomorphism of $R'$.

Indeed, for every $(a,\alpha),(b,\beta)\in R'$, we have
$$\overline{(a,\alpha)+(b,\beta)}=\overline{(a+b,\alpha+\beta)}=(\overline{a+b},-(\alpha+\beta))=(\bar a,-\alpha)+(\bar b,-\beta)=\overline{(a,\alpha)}+\overline{(b,\beta)}$$
and
$$
\begin{aligned}
\overline{(b,\beta)}\cdot \overline{(a,\alpha)}&=(\bar b,-\beta)\cdot(\bar a,-\alpha)=(\bar b{\cdot}\bar a+\mathfrak z(-\alpha)(\overline{-\beta}),\bar{\bar b}{\cdot}(-\alpha)+\bar a{\cdot}(-\beta))\\
&=(\bar b{\cdot}\bar a+\mathfrak z{\cdot}\alpha{\cdot}\overline{\beta},-b{\cdot}\alpha-\bar a{\cdot}\beta)=j'(a{\cdot} b+\mathfrak z{\cdot}\beta{\cdot}\bar\alpha,\bar a{\cdot}\beta+b{\cdot}\alpha)=\overline{(a,\alpha)\cdot(b,\beta)}.
\end{aligned}
$$
Therefore, $j'$ is an anti-isomorphism of $R'$.

To see that the involution $j'$ is central, observe that $$
\begin{aligned}
\mathsf n'(a,\alpha)&=(a,\alpha)\cdot \overline{(a,\alpha)}=(a,\alpha)\cdot(\bar a,-\alpha)=(a{\cdot}\bar a+\mathfrak z{\cdot}(-\alpha){\cdot}\bar \alpha,\bar a{\cdot}(-\alpha)+\bar a{\cdot}\alpha)\\
&=(\mathsf n(a)-\mathfrak z{\cdot}\mathsf n(\alpha),0)\in \mathcal Z(R)\times\{0\}\subseteq \mathcal Z(R')
\end{aligned}
$$
and
$$\mathsf t(a,\alpha)=(a,\alpha)+\overline{(a,\alpha)}=(a,\alpha)+(\bar a,-\alpha)=(\mathsf t(a),0)\in \mathcal Z(R)\times\{0\}\subseteq \mathcal Z(R').$$
\smallskip

1. If the ring $R$ is associative, then for every $(a,\alpha),(b,\beta)\in R'$ we have
$$
\begin{aligned}
((a,\alpha)\cdot(a,\alpha))\cdot(b,\beta)&=(a{\cdot}a+\mathfrak z{\cdot}\alpha{\cdot}\bar\alpha,\bar a{\cdot}\alpha+a{\cdot}\alpha){\cdot}(b,\beta)=(a{\cdot}a+\mathfrak z{\cdot}\mathsf n(\alpha),\mathsf t(a){\cdot}\alpha)\cdot(b,\beta)\\
&=((a{\cdot}a+\mathfrak z{\cdot}\mathsf n(\alpha)){\cdot}b+\mathfrak z{\cdot}\beta{\cdot}\overline{\mathsf t(a){\cdot}\alpha},\overline{a{\cdot}a+\mathfrak z{\cdot}\mathsf n(\alpha)}{\cdot}\beta+b{\cdot}(\mathsf t(a){\cdot}\alpha))\\
&=(a{\cdot}a{\cdot}b+\mathfrak z{\cdot}\mathsf n(\alpha){\cdot}b+\mathfrak z{\cdot}\beta{\cdot}\bar\alpha{\cdot}\overline{\mathsf t(a)},\bar a{\cdot}\bar a{\cdot}\beta+\mathsf n(\alpha){\cdot}\bar{\mathfrak z}{\cdot}\beta+\mathsf t(a){\cdot}b{\cdot}\alpha))
\end{aligned}
$$
On the other hand,
$$
\begin{aligned}
(a,\alpha)\cdot((a,\alpha)\cdot(b,\beta))&=(a,\alpha)\cdot(a{\cdot}b+\mathfrak z{\cdot}\beta{\cdot}\bar\alpha,\bar a{\cdot}\beta+b{\cdot}\alpha)\\
&=(a{\cdot}(a{\cdot}b+\mathfrak z{\cdot}\beta{\cdot}\bar\alpha)+\mathfrak z{\cdot}(\bar a{\cdot}\beta+b{\cdot}\alpha){\cdot}\bar\alpha,\bar a{\cdot}(\bar a{\cdot}\beta+b{\cdot}\alpha)+(a{\cdot}b+\mathfrak z{\cdot}\beta{\cdot}\bar\alpha){\cdot}\alpha)
\\
&=(a{\cdot}a{\cdot}b+\mathfrak z{\cdot}a{\cdot}\beta{\cdot}\bar\alpha+\mathfrak z{\cdot}\bar a{\cdot}\beta{\cdot}\bar\alpha+\mathfrak z{\cdot}b{\cdot}\alpha{\cdot}\bar\alpha,\bar a{\cdot}\bar a{\cdot}\beta+\bar a{\cdot}b{\cdot}\alpha+a{\cdot}b{\cdot}\alpha+\mathfrak z{\cdot}\beta{\cdot}\bar\alpha{\cdot}\alpha)
\\
&=(a{\cdot}a{\cdot}b+\mathfrak z{\cdot}\mathsf t(a){\cdot}\beta{\cdot}\bar\alpha+\mathfrak z{\cdot}b{\cdot}\mathsf n(\alpha),\bar a{\cdot}\bar a{\cdot}\beta+\mathsf t(a){\cdot}b{\cdot}\alpha+\mathfrak z{\cdot}\beta{\cdot}\mathsf n(\alpha)).
\\
\end{aligned}
$$
Taking into account that $\mathfrak z,\mathsf t(a),\mathsf n(a)\in \Re(R)$, we see that $$((a,\alpha)\cdot(a,\alpha))\cdot(b,\beta)=(a,\alpha)\cdot((a,\alpha)\cdot(b,\beta)),$$so the ring $R'$ is left-alternative. Since $j'$ is an anti-isomorphism of the ring $R'$, the left-alternativity of $R'$ implies the right-alternativity of $R'$. Therefore, the ring $R'$ is alternative.

Now assume that the ring $R'$ is alternative. To prove that $R$ is associative, take any elements $a,b,\alpha\in R$ and consider the products
$$
\begin{aligned}
((a,\alpha)\cdot(a,\alpha))\cdot(b,0)&=(a{\cdot}a+\mathfrak z{\cdot}\alpha{\cdot}\bar\alpha,\bar a{\cdot}\alpha+a{\cdot}\alpha){\cdot}(b,0)=(a{\cdot}a+\mathfrak z{\cdot}\mathsf n(\alpha),\mathsf t(a){\cdot}\alpha)\cdot(b,0)\\
&=((a{\cdot}a+\mathfrak z{\cdot}\mathsf n(\alpha)){\cdot}b,b{\cdot}(\mathsf t(a){\cdot}\alpha))=((a{\cdot}a){\cdot}b+\mathfrak z{\cdot}\mathsf n(\alpha){\cdot}b,\mathsf t(a){\cdot}b{\cdot}\alpha)
\end{aligned}
$$
and 
$$
(a,\alpha)\cdot((a,\alpha)\cdot(b,0))=(a,\alpha)\cdot(a{\cdot}b,b{\cdot}\alpha)
=(a{\cdot}(a{\cdot}b)+\mathfrak z{\cdot}(b{\cdot}\alpha){\cdot}\bar\alpha,\bar a{\cdot}(b{\cdot}\alpha)+(a{\cdot}b){\cdot}\alpha).
$$
The left-alternativity of $R'$ implies 
$$
{\bar a}{\cdot}(b{\cdot}\alpha)+a{\cdot}(b{\cdot}\alpha)=\mathsf t(a){\cdot}b{\cdot}\alpha=\bar a{\cdot}(b{\cdot}\alpha)+(a{\cdot}b){\cdot}\alpha
$$and hence $a{\cdot}(b{\cdot}\alpha)=(a{\cdot}b){\cdot}\alpha$, witnessing that the ring $R$ is associative.
\smallskip

2. Assume that the ring $R$ is associative and commutative. Then for every elements $(a,\alpha),(b,\beta),(c,\gamma)\in R'$, we have
$$
\begin{aligned}
((a,\alpha)\cdot(b,\beta))\cdot(c,\gamma)&=(a{\cdot}b+\mathfrak z{\cdot}\beta{\cdot}\bar\alpha,\bar a{\cdot}\beta+b{\cdot}\alpha)\cdot(c,\gamma)\\
&=(a{\cdot}b{\cdot}c+\mathfrak z{\cdot}\beta{\cdot}\bar\alpha{\cdot}c+\mathfrak z{\cdot}\gamma\cdot\overline{\bar a{\cdot}\beta+b{\cdot}\alpha},\overline{a{\cdot}b+\mathfrak z{\cdot}\beta{\cdot}\bar\alpha}{\cdot}\gamma+c{\cdot}\bar a{\cdot}\beta+c{\cdot}b{\cdot}\alpha)\\
&=(a{\cdot}b{\cdot}c+\mathfrak z{\cdot}\beta{\cdot}\bar\alpha{\cdot}c+\mathfrak z{\cdot}\gamma{\cdot}\bar\beta{\cdot}a+\mathfrak z{\cdot}\gamma{\cdot}\bar{\alpha}{\cdot}\bar b,\bar b{\cdot}\bar a{\cdot}\gamma+\alpha{\cdot}\bar\beta{\cdot}\bar{\mathfrak z}{\cdot}\gamma+c{\cdot}\bar a{\cdot}\beta+c{\cdot}b{\cdot}\alpha)\\
\end{aligned}
$$and
$$
\begin{aligned}
(a,\alpha)\cdot((b,\beta)\cdot(c,\gamma))&=(a,\alpha)\cdot(b{\cdot}c+\mathfrak z{\cdot}\gamma{\cdot}\bar\beta,\bar b{\cdot}\gamma+c{\cdot}\beta)\\
&=(a{\cdot}b{\cdot}c+a{\cdot}\mathfrak z{\cdot}\gamma{\cdot}\bar\beta+\mathfrak z{\cdot}(\bar b{\cdot}\gamma+c{\cdot}\beta){\cdot}\bar\alpha,\bar a{\cdot}(\bar b{\cdot}\gamma+c{\cdot}\beta)+(b{\cdot}c+\mathfrak z{\cdot}\gamma{\cdot}\bar\beta){\cdot}\alpha).
\end{aligned}
$$
The commutativity and associativity of $R$ ensure that $((a,\alpha)\cdot(b,\beta))\cdot(c,\gamma)=(a,\alpha)\cdot((b,\beta)\cdot(c,\gamma))$ witnessing that the ring $R'$ is associative.

Now assume that the ring $R'$ is associative. Since the associative ring $R'$ is alternative, Theorem~\ref{t:Cayley-Dickson}(1) ensures that $R$ is associative. For every elements $b,c\in R$, the equality 
$$
\begin{aligned}
(0,c{\cdot}b)&=(0,b)\cdot(c,0)=((0,1)\cdot(b,0))\cdot(c,0)=(0,1)\cdot((b,0)\cdot(c,0))=(0,1)\cdot(b{\cdot}c,0)=(0,b{\cdot}c),
\end{aligned}
$$
witnesses that the ring $R$ is commutative.
\smallskip

3. If $R$ is commutative and the involution $j:R\to R$ is trivial, then for every $(a,\alpha),(b,\beta)\in R'$ we have
$$
\begin{aligned}
(a,\alpha)\cdot(b,\beta)&=(a{\cdot}b+\mathfrak z{\cdot}\beta{\cdot}\bar\alpha,\bar a{\cdot}\beta+b{\cdot}\alpha)=(a{\cdot}b+\mathfrak z{\cdot}\beta{\cdot}\alpha,a{\cdot}\beta+b{\cdot}\alpha)\\
&=(b{\cdot}a+\mathfrak z{\cdot}\alpha{\cdot}\beta,b{\cdot}\alpha+a{\cdot}\beta)=(b{\cdot}a+\mathfrak z{\cdot}\alpha{\cdot}\bar\beta,\bar b{\cdot}\alpha+a{\cdot}\beta)=(b,\beta)\cdot(a,\alpha),
\end{aligned}
$$
witnessing that the ring $R'$ is commutative. 

Now assume that the ring $R'$ is commutative.  Then for every $a,b\in R$ we have $(a{\cdot}b,0)=(a,0)\cdot(b,0)=(b,0)\cdot(a,0)=(b{\cdot}a,0)$, witnessing that the ring $R$ is commutative. To see that the involution of $R$ is trivial, observe that 
$$
(0,\bar a)=(0,\bar a{\cdot}1)=(a,0)\cdot(0,1)=(0,1)\cdot(a,0)=(0,0{\cdot}\bar 1)=(0,a)$$ and hence $\bar a=a$ for all $a\in R$.
\smallskip

4. Assume that $R$ is an alternative division ring and $\mathfrak z\in \Re(R)\setminus \{\mathsf n(x):x\in R\}$. By Proposition~\ref{p:CD-real-axis}, the real axis $\Re(R)$ of the alternative division ring $R$ is a field.

\begin{claim}\label{cl:c-nonzero} For any nonzero element $(a,\alpha)\in R'$, the element $c\defeq a{\cdot}\bar a-\mathfrak z{\cdot}\alpha{\cdot}\bar\alpha\in \Re(R)$ is nonzero.
\end{claim}
 
\begin{proof} By Artin's Theorem~\ref{t:Artin}, the subring $A$ of $R$, generated by the set $\mathcal Z(R)\cup\{a,\alpha\}$ is associative. By Proposition~\ref{p:CD-1gen}, for every nonzero element $b\in A$, the set $\mathcal Z(R)+b{\cdot}\mathcal Z(R)\subseteq A$ is a field, which implies that the associative ring $A$ is a corps. Observe also that $\bar b=b^{-1}\cdot\mathsf n(b)\in A\cdot\mathcal Z(R)\subseteq A$. Therefore, the corps $A$ is closed under the involution of the involutive ring $R$.

If $\alpha=0$, then $a\ne 0\ne\bar a$ and hence $c=a{\cdot}\bar a\ne 0$ because $A$ is a corps. So, assume that $\alpha\ne 0$, which implies $\bar\alpha\ne 0$ and $\mathsf n(\alpha)=\alpha{\cdot}\bar\alpha\ne 0$. Since $\alpha{\cdot}\bar\alpha=\mathsf n(\alpha)\in \Re(R)$, $(\alpha{\cdot}\bar\alpha)^{-1}\in \Re(R)$, by Proposition~\ref{p:CD-real-axis}. Since the involution of the field $C$ is an anti-isomorphism of $C$, $\overline{\alpha^{-1}}=\bar\alpha^{-1}$. Consider the element $b\defeq a{\cdot}\bar\alpha^{-1}$ and observe that $\bar b=\overline{\bar\alpha^{-1}}{\cdot}\bar a=\alpha^{-1}{\cdot}\bar a$. Assuming that $c=0$, we conclude that $$\mathfrak z=(a{\cdot}\bar a)\cdot(\alpha{\cdot}\bar\alpha)^{-1}=a{\cdot}( \alpha{\cdot} \bar\alpha)^{-1}{\cdot}\bar a=
a{\cdot}\bar\alpha^{-1}{\cdot}\alpha^{-1}{\cdot} \bar a=b\cdot\bar b=\mathsf n(b),
$$by the associativity of the corps $A$. But the equality $\mathfrak z=\mathsf n(b)$ contradicts the assumption $\mathfrak z\notin\{\mathsf n(x):x\in R\}$.
\end{proof}

By Claim~\ref{cl:c-nonzero}, for any nonzero element $(a,\alpha)\in R'$, the element  $c\defeq a{\cdot}\bar a-\mathfrak z{\cdot}\alpha{\cdot}\bar\alpha\in \Re(R)$ is nonzero and hence has multiplicative inverse $c^{-1}$ in the field $\Re(R)$. We claim that the element 
$$(b,\beta)\defeq (\bar a{\cdot}c^{-1},-\alpha{\cdot}c^{-1})$$ is inverse to the element $(a,\alpha)$ in the ring $R'$. Indeed,
$$
\begin{aligned}
(a,\alpha){\cdot}(b,\beta)&=(a{\cdot} b+\mathfrak z{\cdot}\beta{\cdot}\bar \alpha,\bar a{\cdot}\beta+b{\cdot}\alpha)\\
&=(a{\cdot}\bar a{\cdot}c^{-1}-\mathfrak z{\cdot}\alpha{\cdot}c^{-1}{\cdot}\bar\alpha,\bar a{\cdot}(-\alpha{\cdot}c^{-1})+(\bar a{\cdot}c^{-1}){\cdot}\alpha)=((a{\cdot}\bar a-\mathfrak z{\cdot}\alpha{\cdot}\bar\alpha){\cdot}c^{-1},0)=(1,0).
\end{aligned}
$$
On the other hand,
$$
\begin{aligned}
(b,\beta){\cdot}(a,\alpha)&=(b{\cdot} a+\mathfrak z{\cdot}\alpha{\cdot}\bar \beta,\bar b{\cdot}\alpha+a{\cdot}\beta)\\
&=(\bar a{\cdot}c^{-1}{\cdot}a-\mathfrak z{\cdot}\alpha{\cdot}c^{-1}{\cdot}{\bar\alpha},c^{-1}{\cdot}\bar{\bar a}{\cdot}\alpha+a{\cdot}(-\alpha{\cdot}c^{-1}))=((\bar a{\cdot}a-\mathfrak z{\cdot}\alpha{\cdot}\bar\alpha){\cdot}c^{-1},0)=(1,0).
\end{aligned}
$$
In the last equality we have used the identity $\mathsf n(\bar a)=\mathsf n(a)$, proved in Proposition~\ref{p:tn-bar}. Therefore, the ring $R'$ is invertible.
\end{proof}

For every involutive ring $R$ and an element $\mathfrak z\in \Re(R)$, we can view the involutive ring $R'=CD(R,\mathfrak z)$ as an extension of $R$ as follows. Note that $a\mapsto (a,0)$ is
an isomorphism of $R$ with a subring of $R'$. Identifying each element $y\in R$ with $(y, 0)$ and
setting ${s}= (0, 1)$, we see that ${s}{\cdot}y = (0,1){\cdot}(y, 0) = (0,y)$. Thus, we have that
$R' = CD(R,{\mathfrak z}) = \{x+{s}{\cdot}y:x,y\in R\}$ and
$(a+{s}{\cdot}\alpha){\cdot}(b+{s}{\cdot}\beta)= ab + {\mathfrak z}{\cdot}\beta{\cdot}\bar\alpha+{s}{\cdot}(\bar a{\cdot}\beta+b{\cdot}\alpha)$. 
We say that the involutive ring $R'=CD(R,\mathfrak z)$ is a \index{Cayle-Dickson extension}\index{involutive ring!Cayley--Dickson extension}\defterm{Cayley--Dickson extension} of the involutive ring $R$, obtained from $R$ by the \index{Cayley--Dickson process}\defterm{Cayley--Dickson process}.  

We now will show that the multiplication in the Cayley-Dickson extension is natural within an alternative ring.

\begin{proposition}\label{p:CD-natural} Let $R$ be an involutive ring. If $R$ is a subring of an alternative ring $S$ and if
$s\in S$ is invertible with
\begin{itemize}
\item $\mathfrak z\defeq s{\cdot}s\in \mathcal N(R)\defeq\{x\in R:\forall a,b\in R\;\;[x,a,b]=0\}$,
\item $x{\cdot}s=s{\cdot}\bar x$ for all $x\in R$, and
\item $R\cap sR=\{0\}$, 
\end{itemize}
then the set $R+sR\defeq\{x+s{\cdot}y:x,y\in R\}$ is a subring of $S$, isomorphic to $CD(R,\mathfrak z)$. Moreover, for all $x,y\in R$,
$$x{\cdot}(s{\cdot}y)=s{\cdot}(\bar x{\cdot}y),\quad (s{\cdot}x){\cdot} y=s{\cdot}(y{\cdot}x),\quad\mbox{and}\quad (s{\cdot}x){\cdot}(s{\cdot}y)=s^2{\cdot}(y{\cdot}\bar x).$$
\end{proposition}

\begin{proof} By Artin's Theorem~\ref{t:Artin}, the alternative ring $S$ is diassociative, and hence $x{\cdot}s^2 = s{\cdot}\bar x{\cdot}s = s^2{\cdot}\bar{\bar x} =s^2{\cdot}x$ for all $x\in R$. Now we see that $\mathfrak z\defeq s^2\in \mathcal Z(R)$. Also 
$s{\cdot}s^2=s^2{\cdot}s=s{\cdot}\overline{s^2}$ implies $s^2=\overline{s^2}$ and hence $s^2\in\Re(R)$ (by the invertibility of $s$). Since $s$ is invertible, the map $R\times R\to R+sR$, $(x,y)\mapsto x+s{\cdot}y$, is bijective. 

Fix any elements $x,y\in R$. Using the left Bol identity, we have
$$s{\cdot}(x{\cdot}(s{\cdot}y)) = (s{\cdot}x{\cdot}s){\cdot}y = (s^2{\cdot}\bar x){\cdot}y =s^2{\cdot}(\bar x{\cdot}y)=s{\cdot}(s{\cdot}(\bar x{\cdot}y))$$
and $x{\cdot}(s{\cdot}y) = s{\cdot}(\bar x{\cdot}y)$, since $s$ is invertible. Similarly, the right Bol identity
gives
$$((s{\cdot}x){\cdot}y){\cdot}s = ((\bar x{\cdot}s){\cdot}y){\cdot}s = \bar x{\cdot}(s{\cdot}y{\cdot}s) = \bar x{\cdot}(\bar y{\cdot}s^2)
=(\bar x{\cdot}\bar y)(s{\cdot}s)=((\bar x{\cdot}\bar y){\cdot}s){\cdot}s = (s{\cdot}(y{\cdot}x)){\cdot}s$$
so $(s{\cdot}x){\cdot}y = s{\cdot}(y{\cdot}x)$. Finally, the Moufang identity shows
$$(s{\cdot}x){\cdot}(s{\cdot}y)=(s{\cdot}x){\cdot}(\bar y{\cdot}s)=(s{\cdot}(x{\cdot}\bar y)){\cdot}s=((y{\cdot}\bar x){\cdot}s){\cdot}s= (y{\cdot}\bar x){\cdot}s^2=s^2{\cdot}(y{\cdot}\bar x).$$

Therefore, for every $a,\alpha,b,\beta\in R$ we have the equalities $$(a+s{\cdot}\alpha)\cdot(b+s{\cdot}\beta)=a{\cdot}b+a{\cdot}(s{\cdot}\beta)+(s{\cdot}\alpha){\cdot}b+(s{\cdot}\alpha){\cdot}(s{\cdot}\beta)=
(a{\cdot}b+\mathfrak z{\cdot}\beta{\cdot}\bar\alpha))+s{\cdot}(\bar a{\cdot}\beta+b{\cdot}\alpha),
$$
witnessing that the ring $R+sR$ is isomorphic to the Cayley--Dickson extension $CD(R,\mathfrak z)$ of $R$.
\end{proof}

\section{Complex fields, quaternion corps, and octonion rings}

\begin{definition} An involutive division ring $R$ is called 
\begin{itemize}
\item a \index{real field}\index{involutive ring!real field}\defterm{real field} if $R$ is associative and the involution of $R$ is trivial;
\item a \index{complex field}\index{involutive ring!complex field}\defterm{complex field} if $R$ is commutative, associative and the involution of $R$ is not trivial;
\item a \index{quaternion corps}\index{involutive ring!quaternion corps}\defterm{quaternion corps} if $R$ is non-commutative and associative;
\item an \index{octonion ring}\index{involutive ring!octonion ring}\defterm{octonion ring} if $R$ is non-associative and alternative.
\end{itemize}
\end{definition}

\begin{exercise} Show that every real field is commutative.
\end{exercise}

A ring $R$ is called \index{Boolean ring}\index{ring!Boolean}\defterm{Boolean} if so is its additive group. The latter happens if and only if $1+1=0$ in $R$.

\begin{example} Let $X$ is an involutive division ring, $\Re(X)\defeq\{x\in\mathcal Z(X):\bar x=x\}$ be its real axis, $\mathfrak z\in \Re(X)\setminus\{\mathsf n(x):x\in X\}$, and $Y\defeq CD(X,\mathfrak z)$ be the Cayley--Dickson extension of the involutive ring $X$.
\begin{enumerate} 
\item If $X$ is a Boolean real field, then $Y$ is a real field.
\item If $X$ is a real field and $X$ is not Boolean, then  $Y$ is a complex field.
\item If $X$ is a complex field, then $Y$ is a quaternion corps.
\item If $X$ is a quaternion corps, then $Y$ is an octonion ring.
\end{enumerate}
\end{example}

\begin{proof} 1. If $X$ is a Boolean real field, then for every $(a,\alpha)\in Y=CD(X,\mathfrak z)$ we have $\overline{(a,\alpha)}=(\bar a,-\alpha)=(a,\alpha)$, witnessing that the involution on the involutive ring $Y$ is trivial and hence $Y$ is a real field.
\smallskip

2. If $X$ is a real field and $X$ is not Boolean, then $1+1\ne 0$ and $\overline{(0,1)}=(\bar 0,-1)\ne (0,1)$. So, the involution of the involutive ring $Y$ is non-trivial. By Theorem~\ref{t:Cayley-Dickson}(3), the ring $Y$ is commutative and associative, and hence the involutive ring $Y$ is a complex field.
\smallskip

3. If $X$ is a complex field, then the ring $Y$ is non-commutative and associative, by Theorem~\ref{t:Cayley-Dickson}(2,3),  and hence $Y$ is a quaterion corps.
\smallskip

4. If $X$ is a quaternion corps, then the ring $Y$ is non-associative and alternative, by Theorem~\ref{t:Cayley-Dickson}(1,2),  and hence $Y$ is an octonion ring.
\end{proof}

\begin{example}\begin{enumerate}
\item The field of real numbers $\IR$ endowed with the trivial involution is a real field.
\item The Cayley--Dickson extension $\IC\defeq CD(\IR,-1)$ of the real field $\IR$ is a complex field, isomorphic to the field of complex numbers.
\item The Cayley--Dickson extension $\IH\defeq CD(\IC,-1)$ of the complex field $\IC$ is a quaternion corps, isomorphic to the corps of quaternions.
\item The Cayley--Dickson extension $\mathbb O\defeq CD(\IH,-1)$ of the quaternion corps $\IH$ is an octonion ring, isomorphic to the alternative ring of octonions.
\end{enumerate}
\end{example}

If $D$ is an involutive division ring, then its real axis $R\defeq\Re(D)$ is a field (by Proposition~\ref{p:CD-real-axis} and $D$ can be considered as an $R$-vector space over the field $R$. The dimension $\dim_R(D)$ of this vector space is called the \defterm{$R$-dimension} of $D$.

\begin{lemma}\label{l:CD-orthogonal} Let $D$ be an involutive division ring, $R\defeq\Re(D)$ be its real axis, and $A$ be a division subring of $D$ such that $R\subseteq A\ne D$ and $A$ has finite $R$-dimension. If $D$ is not Boolean or $A$ contains an element $a\ne \bar a$, then there exists an element $d\in D\setminus A$ with $d^2\in R$ such that $x{\cdot}\bar d+d{\cdot}\bar x=0$ and $x{\cdot}d=d{\cdot}\bar x$ for all $x\in A$.
\end{lemma}

\begin{proof} Since $A\ne D$, there exists an element $b\in D\setminus A$. Observe that the left shift $l_b:D\to D$, $l_b:x\mapsto b{\cdot}x$, is an automorphism of the $R$-vector space $D$ and hence $\dim_R(bA)=\dim_R(A)$.

Assuming that $A\cap bA\ne\{0\}$, we can find a non-zero elements $x,y\in A$ such that $y=b{\cdot}x$. Since $A$ is a division ring, there exists an element $\beta\in A$ such that $y=\beta{\cdot}x$. Since $D$ is a division ring, the equality $b{\cdot}x=y=\beta{\cdot}x$ implies $b=\beta\in A$, which contradicts the choice of $b$. This contradiction shows that $A\cap bA=\{0\}$ and hence $\dim_R(D)\ge 2k$ where $k\defeq\dim_R(A)$. 

Let $e_1,\dots,e_k$ be a basis of the $R$-vector space $A$ such that $e_1=1\in R\subseteq A$. Observe that for every $i\in\{1,\dots,k\}$,  the function $$\mathsf n_i:D\to R,\;\mathsf n_i:x\mapsto \mathsf n(e_i+x)-\mathsf n(e_i)-\mathsf n(x)=e_i{\cdot}\bar x+x{\cdot}e_i,$$is an $R$-linear functional on the $R$-vector space $D$. Its kernel $\Ker(\mathsf n_i)\defeq\mathsf n_i^{-1}(0)$ has codimension $\dim_R(D/\Ker(n_i))$ at most $1$ in the $R$-vector space $D$.

Then the $R$-vector subspace $A^\perp\defeq\bigcap_{i=1}^k\Ker(\mathsf n_i)$ of $D$ has $R$-dimension $\dim_R(A^\perp)\ge \dim_R(D)-k\ge 2k-k=k$. Assuming that $A^\perp\subseteq A$, we conclude that $A^\perp=A$ and hence $\mathsf n_i(x)=0$ for all $x\in A$ and $i\in\{1,\dots,k\}$. In particular $0=\mathsf n_1(1)=1{\cdot}\bar 1+1{\cdot}\bar 1=1+1$ and $x+x=0=\mathsf n_1(x)=1{\cdot}\bar x+x{\cdot}\bar 1=\bar x+x$, which implies $\bar x=x$ and contradicts the assumption (that $D$ is not Boolean or $A$ contains an element $a\ne\bar a$). This contradiction shows that $A^\perp\ne A$ and hence there exists an element $d\in A^\perp\setminus A$. Then $0=\mathsf n_i(d)=e_i{\cdot}\bar d+d{\cdot}\bar e_i$ for all $i\in\{1,\dots,n\}$. In particular, $0=\mathsf n_1(d)=1{\cdot}\bar d+d{\cdot}\bar 1=d+\bar d$ and hence $\bar d=-d$ and $d^2=-d{\cdot}\bar d\in R$.

Given any element $x\in A$, find elements $x_1,\dots,x_k\in R$ such that $x=\sum_{i=1}^k x_i{\cdot}e_i$. Then 
$$x{\cdot}\bar d+d{\cdot}\bar x=\Big(\sum_{i=1}^kx_i{\cdot}e_i\Big){\cdot}\bar d+d\cdot\sum_{i=1}^k\bar e_i{\cdot}\bar x_i=\sum_{i=1}^kx_i{\cdot}(e_i{\cdot}\bar d+d{\cdot}\bar e_i)=\sum_{i=1}^kx_i{\cdot}\mathsf n_i(d)=0,$$
which implies $x{\cdot}d=x{\cdot}(-\bar d)=-x{\cdot}\bar d=d{\cdot}\bar x$.
\end{proof}

\begin{theorem}\label{t:248CD}Let $D$ be an involutive alternative division ring and $R\defeq\Re(D)$ be its real axis.
\begin{enumerate}
\item If $D$ is a complex field, then $\dim(D)=2$. Moreover, if $D$ is not Boolean, then $D$ is isomorphic to a Cayley--Dickson extension of some real field.
\item If $D$ is a quaternion corps, then $\dim_R(D)=4$ and $D$ is isomorphic to a Cayley--Dickson extension of some complex field.  
\item If $D$ is an octonion ring, then $\dim_R(D)=8$ and $D$ is isomorphic to a Cayley--Dickson extension of some quaternion corps.
\end{enumerate}
\end{theorem}

\begin{proof} By Proposition~\ref{p:CD-real-axis}, the real axis $R\defeq\Re(D)=\{x\in \mathcal Z(D):\bar x=x\}$ of the division ring $D$ is a subfield of $D$.
\smallskip

1. Assume that $D$ is a complex field. Then the involution on $D$ is not trivial and hence $\bar c\ne c$ for some $c\in C$. By Proposition~\ref{p:CD-1gen}, the set $C\defeq R+c{\cdot}R$ is a subfield of the ring $D$. Assuming that $C\ne D$, we can apply Lemma~\ref{l:CD-orthogonal} and find an $q\in D\setminus C$ with $q^2\in R$ such that $x{\cdot}\bar q+q{\cdot}\bar x=0$ for all $x\in C$. By Proposition~\ref{p:CD-natural}, the subset $C+q{\cdot}C$ of the ring $D$ is a ring, isomorphic to the Cayley--Dickson extension $CD(C,q^2)$ of the complex field $C$. Since the involution on $C$ is not trivial, the ring $C+q{\cdot}C\cong CD(C,q^2)$ is not commutative, by Theorem~\ref{t:Cayley-Dickson}(3). But this contradicts the commutativity of the complex field $D$. This contradiction shows that $D=C$ and hence $\dim_R(D)=\dim_R(C)=2$. 

If the ring $D$ is not Boolean, then by Lemma~\ref{l:CD-orthogonal}, there exists an element $d\in D\setminus R$ such that $d^2\in R$ and $x{\cdot}d=d{\cdot}x$ for all $x\in R$. It follows from $\dim_R(D)=2$ that $D=R+d{\cdot}R$. By Proposition~\ref{p:CD-natural}, the field $R+d{\cdot}R=D$ is isomorphic to the Cayley--Dickson extension $CD(R,\gamma^2)$ of the real field $R$.
\smallskip

2. Now assume that $D$ is a quaternion corps. Since $D$ is not commutative, the involution of $D$ is not trivial. In the opposite case, we would obtain that $x{\cdot}y=\bar x{\cdot}\bar y=\overline{y{\cdot} x}=y{\cdot}x$ for all $x,y\in D$, which contradicts the non-commutativity of $D$. Choose any element $c\in D$ with $\bar c\ne c$. Repeating the argument from the preceding item, show that $C\defeq R+c{\cdot}C$ is a complex field and there exists an element $q\in D\setminus C$ such that the set $Q\defeq C+q{\cdot}C$ is a subring of $D$, isomorphic to the Cayley--Dickson extension of the complex field $C$.  Assuming that $D\ne Q$, we can apply Lemma~\ref{l:CD-orthogonal} and find an element $p\in D\setminus Q$ such that $p^2\in R$ and $x{\cdot} p=p{\cdot}\bar x$ for all $x\in Q$. By Proposition, the set $Q+p{\cdot}Q$ is isomorphic to the Cayley-Dickson extension $CD(Q,p^2)$ of the quaternion corps $Q=C+d{\cdot}C$. Since the corps $Q$ is not commutative, the ring $Q+p{\cdot}Q\cong CD(Q,p^2)$ is not associative, which contradicts the associativity of the corps $D$. This contradiction shows that $D=Q=C+q{\cdot}C$ and hence $D$ is isomorphic to the Cayley--Dickson extension $CD(C,q^2)$ of the complex field $C$.
\smallskip

3. Finally, assume that $D$ is an octonion ring. Then $D$ is not associative and not commutative, by Corollary~\ref{c:div+com+alt=>ass}. Repeating the argument from the preceding item, find an element $c\ne \bar c$ and show that $C\defeq R+c{\cdot}C$ is a complex field. Since $D$ is not commutative, $C\ne D$, so we can find an element $q\in D\setminus C$ such that the set $Q\defeq C+q{\cdot}C$ is an quaternion corps, isomorphic to the Cayley--Dickson extension $CD(C,d^2)$ of the complex field $C$. Since $Q$ is associative and $D$ is not, we can apply Lemma~\ref{l:CD-orthogonal} and Proposition~\ref{p:CD-natural}, and find an element $p\in D\setminus Q$ such that $p^2\in R$ and the set $A\defeq Q+p{\cdot}Q$ is a subring of $D$, isomorphic to the Cayley--Dickson extension $CD(Q,p^2)$ of the quaternion corps $Q$. Since $Q$ is not commutative, the ring $A$ is not associative. Assuming that $D\ne A$, we can apply Lemma~\ref{l:CD-orthogonal} and find an element $b\in D\setminus A$ such that $b^2\in R$ and $x{\cdot}b=b{\cdot}\bar x$ for every $x\in A$. By Proposition, the set $A+b{\cdot}A$ is isomorphic to the Cayley-Dickson extension $CD(A,b^2)$ of the non-associative alternative ring $A$. Since the ring $A$ is not associative, the ring $A+b{\cdot}A\cong CD(A,p^2)$ is not alternative (by Theorem~\ref{t:Cayley-Dickson}(1)), which contradicts the alternativity of the ring $D$. This contradiction shows that $D=A$ and hence $D$ is isomorphic to the Cayley--Dickson extension $CD(Q,p^2)$ of the quaternion corps $Q$.
\end{proof}

\begin{corollary}\label{c:dim=1248} Let $D$ be an involutive alternative division ring and $R\defeq \Re(D)$ be its real axis. The $R$-dimension $\dim_R(D)$ of $D$ is equal to $1,2,4$ or $8$. More precisely, 
\begin{enumerate}
\item $\dim_R(D)=1$ if and only if $D$ is a real field;
\item $\dim_R(D)=2$ if and only if $D$ is a complex field;
\item $\dim_R(D)=4$ if and only if $D$ is a quaternion corps;
\item $\dim_R(D)=8$ if and only if $D$ is an octonion ring.
\end{enumerate}
\end{corollary}  
 
\section{Non-associative alternative division rings are octonion rings}

In this section we prove an important and difficult \index[person]{Kleinfeld}Skornyakov--Bruck--Kleinfeld\footnote{{\bf Erwin Kleinfeld} (1929--2022) was an american mathematician, born in Austria. In 1951 he defended the Ph.D. Thesis ``Alternative rings'' in University of Wisconsin-Madison, under supervision of Richard Bruck. Subsequently, Kleinfeld held lecturer positions at the University of Chicago, Yale University, Ohio State University, Syracuse University, and finally at the University of Iowa (since 1968). His career rapidly took off due to fundamental contributions in algebra, and more specifically in non-associative algebras, in the 50’s and early 1960’s. His work on alternative rings, and later Jordan algebras was featured at several International Congress of Mathematicians, and was essential in many subsequent landmark contributions by top mathematicians such as M. Artin, N. Jacobson, M. Zorn, and of Fields Medalist Efim Zelmanov, and lead to later work related to that of Fields Medalist Sergei Novikov. 
In 2017, together with his wife, and shortly after her death, the Kleinfelds made the single largest donation to the Department of Mathematics in its entire history. The very generous donation, followed by an additional endowment from the Kleinfeld estate, created the Kleinfeld fellowship as well as several other awards for graduate students in mathematics.} Theorem on the octonion structure of non-associative alternative division rings.

Let us recall that for a ring $R$, the sets 
$$
\begin{aligned}
\mathcal N(R)&\defeq\{a\in R:\forall x,y\in R\;[a,x,y]=[x,a,y]=[x,y,a]=0\}\quad\mbox{and}\\
\mathcal Z(R)&\defeq\{a\in \mathcal N(R):\forall x\in R\;\;[a,x]=0\}
\end{aligned}
$$are called the {\em nucleus} and the {\em centre} of $R$. Here $$[x,y,z]\defeq(x{\cdot}y){\cdot}z-x{\cdot}(y{\cdot}z)\quad\mbox{and}\quad
[x,y]\defeq x{\cdot}y-y{\cdot}x$$are the {\em associator} and the {\em commutator} of elements $x,y,z$ and $x,y$, respectively.
By Proposition~\ref{p:NZ-corps-field}, the nucleus $\mathcal N(R)$ of a division ring $R$ is a corps, and the center $\mathcal Z(R)$ of $R$ is a field.

\begin{theorem}[Skornyakov, 1950; Bruck, Kleinfeld, 1951]\label{t:SBK} Every non-associative alternative division ring $R$ admits a central involution turning $R$ into an octonion ring  with real axis $\Re(R)=\mathcal Z(R)=\mathcal N(R)$.
\end{theorem}

\begin{proof}  Let $R$ be a non-associative alternative division ring. Given any elements $a,b\in R$, consider the operators 
$$
\begin{aligned}
&l_a:R\to R,\;\;l_a:x\mapsto a{\cdot}x,\\
&A_{a,b}:R\to R,\;\;A_{a,b}:x\mapsto [a,b,x],\quad\mbox{and}\\
&B_{a,b}\defeq l_{[a,b]}-A_{a,b}.
\end{aligned}
$$

\begin{claim}\label{cl:4.7-4.10} For every elements $a,b,x,y\in R$ the following identities hold:
\begin{enumerate}
\item $l_aA_{a,b}(x)= A_{a,b{\cdot}a}(x)$;
\item $A_{a,b}(x{\cdot}y) = A_{a,b}(x){\cdot}y + x{\cdot}A_{a,b}(y) - [[a, b], x, y]$;
\item $A_{a,b}^2(x) = l_{[a,b]}A_{a,b}(x)$;
\item $B_{a,b}(x{\cdot}y) = B_{a,b}(x){\cdot}y - x{\cdot}A_{a,b}(y)$;
\item $[a,b{\cdot}a]=[a,b]{\cdot}a$.
\end{enumerate}
\end{claim}

\begin{proof} By Theorems~\ref{t:alternative<=>left-inverse} and \ref{t:alternative<=>right-Bol}, the alternative division ring is Bol. 
\smallskip

1. By Proposition~\ref{p:left-Bol}(2,3), for every $x\in R$ we have
$$l_aA_{a,b}(x)=a{\cdot}[a,b,x]=-a{\cdot}[b,a,x]=[a,b{\cdot}a,x]=A_{a,b{\cdot}a}(x).$$ 

2. By Proposition~\ref{p:left-Bol}(5), 
$$A_{a,b}(x{\cdot}y) =[a,b,x{\cdot}y]=[a,b,x]{\cdot}y+x{\cdot}[a,b,y]-[[a,b],x,y]=
A_{a,b}(x){\cdot}y+x{\cdot}A_{a,b}(y)- [[a, b], x, y].$$

3. By Theorem~\ref{t:Mikheev}(1), 
$A_{a,b}^2(x)=[a,b,[a,b,x]]=[a,b]{\cdot}[a,b,x]=l_{[a,b]}A_{a,b}(x)$.
\smallskip

4. Applying the second item and the definition of the operator $B_{a,b}$, we obtain
$$
\begin{aligned}
B_{a,b}(x{\cdot}y)&=l_{[a,b]}(x{\cdot}y)-A_{a,b}(x{\cdot}y)=
[a,b]{\cdot}(x{\cdot}y)-A_{a,b}(x){\cdot}y-x{\cdot}A_{a,b}(y)+[[a,b],x,y]\\
&=[a,b]{\cdot}(x{\cdot}y)-A_{a,b}(x){\cdot}y-x{\cdot}A_{a,b}(y)+([a,b]{\cdot}x){\cdot}y-[a,b]{\cdot}(x{\cdot}y)\\
&=(l_{[a,b]}(x)-A_{a,b}(x)){\cdot}y-x{\cdot}A_{a,b}(y)=B_{a,b}(x){\cdot}y-x{\cdot}A_{a,b}(y).
\end{aligned}
$$

5.  By Theorem~\ref{t:Bol<=>}, the Bol ring $R$ is flexible and hence  
$$[a,b{\cdot}a]=a{\cdot}(b{\cdot}a)-(b{\cdot}a){\cdot}a=(a{\cdot}b){\cdot}a-(b{\cdot}a){\cdot}a=(a{\cdot}b-b{\cdot}a){\cdot}a=[a,b]{\cdot}a.$$
\end{proof}
 
\begin{claim}\label{cl:4.11} If $A_{a,b} = 0$ for some elements $a,b\in R$, then $[a, b] = 0$ or $a\in\mathcal N(R)$.
\end{claim}

\begin{proof}  By Claim~\ref{cl:4.7-4.10}(2,5), $[a, b]\in\mathcal N\defeq\mathcal  N(R)$ and $[a,b{\cdot}a]=[a,b]{\cdot}a$. Since $A_{a,b{\cdot}a} = 0$ by Claim~\ref{cl:4.7-4.10}(1), we also have $[a,b{\cdot}a]\in\mathcal N$. If $[a,b]\ne 0$, then $a= [a, b]^{-1}{\cdot}([a,b]{\cdot}a)=[a,b]^{-1}\cdot [a, b{\cdot}a]\in\mathcal N$ (because $\mathcal N$ is a corps, according to Proposition~\ref{p:NZ-corps-field}).
\end{proof}

\begin{claim}\label{cl:4.12} $\mathcal N(R) = \mathcal Z(R)$.
\end{claim}

\begin{proof}
Let $b\in\mathcal N\defeq\mathcal N(R)$, so $A_{a,b} = 0$ for all $a\in R$. If $a\notin \mathcal N$, then $[a, b] = 0$ according to Claim~\ref{cl:4.11}. On the other hand, if $a\in\mathcal N$, there is some $c\notin\mathcal N$, since $R$ is not associative.
Since $c, a+ c\notin\mathcal N$, we have $[a, b] = [a+ c, b] - [c, b] = 0$, by Claim~\ref{cl:4.11}. Therefore, $b\in\mathcal N$ implies $[a, b] = 0$ for all $a\in R$ and hence $b\in\mathcal Z(R)$.
\end{proof}

Claims~\ref{cl:4.11} and \ref{cl:4.12} imply the following improvement of Claim~\ref{cl:4.11}.

\begin{claim}\label{cl:4.13} If $A_{a,b} = 0$ for some elements $a,b\in R$, then $[a, b] = 0$.
\end{claim}

We next show that for every $a,b\in R$, the kernel $\mathcal A_{a,b}\defeq A_{a,b}^{-1}(0)$ of the operator $A_{a,b}$ is a subalgebra of $R$. If $[a, b] = 0$, this is immediate from Claim~\ref{cl:4.7-4.10}(2). If $[a, b]\ne 0$, we shall show even more. Let $\mathcal B_{a,b}\defeq B_{a,b}^{-1}(0)$ be the kernel of the operator $B_{a,b}\defeq l_{[a,b]}-A_{a,b}$.

\begin{claim}\label{cl:4.14-4.21} If $a,b\in R$ are two elements with $[a,b]\ne 0$, then
\begin{enumerate}
\item $R = \mathcal A_{a,b} \oplus \mathcal B_{a,b}$,
\item $A_{a,b} = \rng[B_{a,b}] \ne \{0\}$,
\item $\mathcal B_{a,b} = \rng[A_{a,b}]\ne \{0\}$,
\item $\mathcal A_{a,b}$ is a division subalgebra of $R$,
\item  $\mathcal B_{a,b}{\cdot}\mathcal A_{a,b} + \mathcal A_{a,b}{\cdot}\mathcal B_{a,b} \subseteq \mathcal B_{a,b}$,
\item  $\mathcal B_{a,b}{\cdot}\mathcal B_{a,b} \subseteq \mathcal A_{a,b}$,
\item  $v{\cdot}\mathcal A_{a,b} = \mathcal B_{a,b}$ for any $v\in\mathcal  B_{a,b}\setminus\{0\}$,
\item $u\in\mathcal N(\mathcal A_{a,b})$ and $v\in \mathcal B_{a,b}\setminus\{0\}$ with $[u, v, \mathcal A_{a,b}] = 0$ imply $u\in\mathcal N(R)$,
\item $s^2\in \mathcal N(R)$ for every $s\in \mathcal B_{a,b}$ with $s^2\in\mathcal N(\mathcal A_{a,b})$.
\end{enumerate}
\end{claim}

\begin{proof} Observe that $[a,b]\ne 0$ implies $a,b\in R^*\defeq R\setminus\{0\}$. By Theorems~\ref{t:alternative<=>left-inverse} and \ref{t:alternative<=>right-Bol}, the alternative division ring $R$ is Bol, and hence its  multiplicative loop $(R^*,\cdot)$ is Bol and Moufang.
\smallskip

 1--3. Since $A_{a,b}([a, b]) = [a, b, [a, b]] = 0$ by Artin's Theorem~\ref{t:Artin}, we see that for $x = [a, b]$ in Claim~\ref{cl:4.7-4.10}(2) we obtain
$$
\begin{aligned}
A_{a,b}l_{[a,b]}(y)&=A_{a,b}([a,b]{\cdot}y)=A_{a,b}([a,b]){\cdot}y+[a,b]{\cdot}A_{a,b}(y)-[[a,b],[a,b],y]\\
&=0{\cdot}y+[a,b]{\cdot}[a,b,y]-0=l_{[a,b]}A_{a,b}(y),
\end{aligned}
$$
so the operators $A_{a,b}$ and $l_{[a,b]}$ commute. Thus, the operator $P\defeq l^{-1}_{[a,b]}A_{a,b}$ has $P^2 = P$ by Claim~\ref{cl:4.7-4.10}(3). Thus, $R = \rng[P] \oplus \rng[Id - P]$. Note that $Id - P = l^{-1}_{[a,b]}B_{a,b} = B_{a,b}l^{-1}_{[a,b]}$, so
$$
\begin{aligned}
\mathcal A_{a,b} &= \Ker(P) = \rng[Id - P] = \rng[B_{a,b}],\\
\mathcal B_{a,b} & = \Ker[Id - P] = \rng[P] = \rng[A_{a,b}].
\end{aligned}
$$
Since $0 \ne a\in \mathcal A_{a,b}$ and $\mathcal A_{a,b}\ne R$ by Claim~\ref{cl:4.13}, the conditions (1)--(3) of Claim~\ref{cl:4.14-4.21} hold.
\smallskip 

4. If $z, y \in \mathcal A_{a,b}=\Ker(A_{a,b})=\rng[B_{a,b}]$, we find an element $x\in R$ with $z = B_{a,b}(x)$ and use Claim~\ref{cl:4.7-4.10}(4) to get $z{\cdot}y =B_{a,b}(x){\cdot}y=
B_{a,b}(x{\cdot}y)+x{\cdot}A_{a,b}(y)=B_{a,b}(x{\cdot}y)+x{\cdot}0\in\rng[B_{a,b}]=\mathcal A_{a,b}$. Thus, $\mathcal A_{a,b}$ is a subalgebra. 
If $0\ne x\in \mathcal A_{a,b}$, then $[a,b,x]=0$ and $[a,b,x^{-1}]=0$, by Theorem~\ref{t:semi-Moufang}. Therefore, $x^{-1}\in\mathcal A_{a,b}$, witnessing that the condition (4) of Claim~\ref{cl:4.14-4.21} holds.
\smallskip

5. Taking $x\in\mathcal B_{a,b}$ and $y\in\mathcal A_{a,b}$ in Claim~\ref{cl:4.7-4.10}(4) shows $x{\cdot}y\in \Ker(B_{a,b}) = \mathcal B_{a,b}$; i.e.,
$\mathcal B_{a,b}{\cdot}\mathcal A_{a,b}\subseteq \mathcal B_{a,b}$. Consider the opposite ring $R^{op}$ endowed with the mirror multiplication $x*y=y\cdot x$. Since $R^{op}$ is an alternative division ring, $\mathcal B^{op}_{a,b}{\cdot}\mathcal A^{op}_{a,b}\subseteq \mathcal B^{op}_{a,b}$, by the preceding paragraph. In the opposite ring $R^{op}$, the associator is
$$[a, b, c]^{op} = c{\cdot}(b{\cdot}a) - (c{\cdot}b){\cdot}a = -[c, b, a] = [a, b, c],$$according to Proposition~\ref{p:alternative9}(1).
Thus, $A^{op}_{a,b} = A_{a,b}$, $\mathcal A^{op}_{a,b}=\mathcal A_{a,b}$, and $\mathcal B^{op}_{a,b}=\rng[A_{a,b}] = \mathcal B_{a,b}$, so $$\mathcal A_{a,b}{\cdot}\mathcal B_{a,b}=\mathcal B_{a,b}*\mathcal A_{a,b}=\mathcal B^{op}_{a,b}*\mathcal A^{op}_{a,b}\subseteq\mathcal B^{op}_{a,b}=\mathcal B_{a,b}.$$ Thus, the condition (5) of Claim~\ref{cl:4.14-4.21} holds. 
\smallskip

6. Take any $x, z \in\mathcal B_{a,b}$ and find $y\in R$ with $z = A_{a,b}(y)$. By Claim~\ref{cl:4.7-4.10}(4), we get $x{\cdot}z =x{\cdot}A_{a,b}(y)=B_{a,b}(x){\cdot}y-B_{a,b}(x{\cdot}y)=0{\cdot}y-B_{a,b}(x{\cdot}y)=\rng[B_{a,b}]=\mathcal A_{a,b}$. Therefore, $\mathcal B_{a,b}{\cdot}\mathcal B_{a,b}\subseteq \mathcal A_{a,b}$, witnessing that the condition (6) holds.
\smallskip

7.  Given any elements $v\in \mathcal B_{a,b}\setminus\{0\}$ and $x\in\mathcal B_{a,b}$, write the element $l^{-1}_v(x)\in R=\mathcal A_{a,b}\oplus\mathcal B_{a,b}$ as $l^{-1}_v(x) = y + z$ for some $y\in\mathcal A_{a,b}$ and $z\in\mathcal B_{a,b}$. By Claim~\ref{cl:4.14-4.21}(5), 
$$v{\cdot}z=x-v{\cdot}y\in \mathcal B_{a,b}-\mathcal B_{a,b}\cdot\mathcal A_{a,b}=\mathcal B_{a,b}-\mathcal B_{a,b}=\mathcal B_{a,b}$$and by Claim~\ref{cl:4.14-4.21}(6), $v{\cdot}z\in \mathcal B_{a,b}\cap(\mathcal B_{a,b}{\cdot}\mathcal B_{a,b})\subseteq\mathcal B_{a,b}\cap\mathcal A_{a,b}=\{0\}$, which implies $x=v{\cdot}y\in v{\cdot}\mathcal A_{a,b}$. Therefore, $\mathcal B_{a,b}\subseteq v{\cdot}\mathcal A_{a,b}$. On the other hand, the condition (5) ensures that $v{\cdot}\mathcal A_{a,b}\subseteq \mathcal B_{a,b}{\cdot}\mathcal A_{a,b}\subseteq \mathcal B_{a,b}$. Therefore, $v{\cdot}\mathcal A_{a,b}=\mathcal B_{a,b}$, showing the condition (7).
\smallskip

8. Fix any elements $u\in \mathcal N(\mathcal A_{a,b})$ and $v\in\mathcal B_{a,b}\setminus\{0\}$ with $[u,v,\mathcal A_{a,b}]\defeq \{[u,v,x]:x\in\mathcal A_{a,b}\}=\{0\}$. Then for every $w\in\mathcal A_{a,b}$, we have $[u, w, v] = -[u, v, w] = 0$. Since
$u\in\mathcal N(\mathcal A_{a,b})$, the subalgebra $\mathcal A_{u,w}\defeq\{x\in R:[u,w,x]=0\}$ contains $\mathcal A_{a,b}$ and $v$. Thus,
$$R = \mathcal A_{a,b}\oplus \mathcal B_{a,b}= \mathcal A_{a,b}\oplus v{\cdot}\mathcal A_{a,b}\subseteq \mathcal A_{u,w},$$
for all $w \in\mathcal A_{a,b}$. Thus, $[u, R,\mathcal A_{a,b}] = [u, \mathcal A_{a,b}, R] = 0$. Now if $z\in\mathcal B_{a,b}\setminus\{0\}$,
the subalgebra $\mathcal A_{u,z}$ contains $\mathcal A_{a,b}$ and $z$. As before, this gives $R\subseteq \mathcal A_{u,z}$, so $[u, \mathcal B_{a,b}, R] = 0$. Hence, $[u, R, R] = 0$; i.e., $u\in\mathcal N(R) =\mathcal Z(R)$.
\smallskip

9. Given any $s\in\mathcal B_{a,b}$ with $s^2\in\mathcal N(\mathcal A_{a,b})$, we have $[s^2,s,R]=0$, by Artin's Theorem~\ref{t:Artin}. Applying the preceding item, we get $s^2\in \mathcal N(R)=\mathcal Z(R)$.
\end{proof}


\begin{claim}\label{cl:pqZ} There exists elements $a,b\in R$ such that $[a,b]\notin \mathcal Z(R)$.
\end{claim}

\begin{proof} Since the ring $R$ is not commutative, there exist elements $p,q\in R$ such that $[p,q]\ne 0$, which implies $p,q\notin \mathcal Z(R)$. By Claim~\ref{cl:4.7-4.10}(5), $[p,q{\cdot}p]=[p,q]{\cdot}p$.
If $[p,q]\notin \mathcal Z(R)$, then put $a\defeq p$ and $b\defeq q$. If $[p,q]\in \mathcal Z(R)$, then $p\notin \mathcal Z(R)$ implies $[p,q{\cdot}p]=[p,q]\cdot p\notin \mathcal Z(R)$. In this case, put $q\defeq p$ and $b\defeq q{\cdot}p$.
\end{proof}

By Claim~\ref{cl:pqZ}, there exist elements $a,b\in \mathcal Z(R)$ such that $[a,b]\notin \mathcal Z(R)$. Consider the division ring $\mathcal A\defeq \mathcal A_{a,b}$ and the set $\mathcal B \defeq \mathcal B_{a,b}$. 

\begin{claim} The division ring $\mathcal A$ is associative.
\end{claim}

\begin{proof} By Claim~\ref{cl:4.7-4.10}(2), 
$a,b\in\mathcal A$ and $[a,b]\in \mathcal N(\mathcal A)$.  Since $A_{a,b} = l_{[a,b]}P$ and $B_{a,b} = l_{[a,b]}(Id - P)$ 
where $P$ is the projection onto $\mathcal B$, Claims~\ref{cl:4.14-4.21}(5) and \ref{cl:4.7-4.10}(4) with $x\in \mathcal A$ and $y \in \mathcal B$ give  $x{\cdot}y\in \mathcal A\cdot\mathcal B\subseteq\mathcal B$ and
$$
\begin{aligned}
0 &=B_{a,b}(x{\cdot}y)=B_{a,b}(x){\cdot}y-x{\cdot}A_{a,b}(y)= l_{[a,b]}(x){\cdot}y -A_{a,b}(x){\cdot}y-x{\cdot}l_{[a,b]}P(y)\\
&=([a,b]{\cdot}x){\cdot}y-0{\cdot}y-x{\cdot}l_{[a,b]}(y)=
\big[[a,b],x\big]{\cdot}y+(x{\cdot}[a,b]){\cdot}y-x{\cdot}([a,b]{\cdot}y)\\
&= \big[[a,b], x\big]{\cdot}y + [x,[a,b], y].
\end{aligned}
$$
Assuming that $[a,b]\in \mathcal Z(\mathcal A)$, we conclude that $[\mathcal A,[a,b],y]=0=[[a,b],y,\mathcal A]$ and $[a,b]\in \mathcal Z(R)$ by Claim~\ref{cl:4.14-4.21}(8), which is a contradiction. This shows that $[a,b]\in \mathcal N(\mathcal A)\ne\mathcal Z(\mathcal A)$. Claim~\ref{cl:4.12} (applied to the alternative divisible ring $\mathcal A$) implies that the ring $\mathcal A$ is associative.
\end{proof}

By Claim~\ref{cl:4.14-4.21}(3), there exists a non-zero element $s\in \mathcal B\setminus\{0\}$. Then $s^2\in \mathcal B{\cdot}\mathcal B\subseteq \mathcal A=\mathcal N(\mathcal A)$. Claim~\ref{cl:4.14-4.21}(9) guarantees that $\mathfrak z\defeq  s^2\in \mathcal Z(R)$. 

We will show
that $j : x\mapsto \bar x\defeq s{\cdot}x{\cdot}s^{-1}$ is a central involution of $\mathcal A$. By Artin's Theorem~\ref{t:Artin}, for every $x\in\mathcal A$, the subring  of $R$, generated by the set $\{x,s\}\cup \mathcal Z(R)$ is associative. This subring contains the element $s^{-1}=\mathfrak z^{-1}{\cdot}s$. 
Then
$$s{\cdot}x=s{\cdot}x{\cdot}s^{-1}{\cdot}s=\bar x{\cdot}s.$$
Also 
$$\bar{\bar x}= s{\cdot}(s{\cdot}x{\cdot}s^{-1}){\cdot}s^{-1}=\mathfrak z{\cdot}x{\cdot}\mathfrak z^{-1} = x.$$ By Claim~\ref{cl:4.14-4.21}(5,6,9), $x{\cdot}s\in \mathcal A{\cdot}\mathcal B\subseteq\mathcal B$, $(x{\cdot}s)^2\in\mathcal B{\cdot}\mathcal B\subseteq \mathcal A=\mathcal N(\mathcal A)$, and $(x{\cdot}s)^2\in \mathcal Z(R)$.  Then  
$$\mathsf n(x)\defeq x{\cdot}\bar x = x{\cdot}s{\cdot}x{\cdot}s^{-1}= x{\cdot}s{\cdot}x{\cdot}\mathfrak z^{-1}{\cdot}s = \mathfrak z^{-1}{\cdot}(x{\cdot}s)^2\in \mathcal Z(R)$$
and $\mathsf t(x)=\mathsf n(x+1)-\mathsf n(x)-1\in \mathcal Z(R)$. For every $z\in \mathcal Z(R)$, we have $\bar z=s{\cdot}z{\cdot}s^{-1}=z{\cdot}s{\cdot}s^{-1}=z$ and hence $\mathcal R(R)=\mathcal Z(R)=\mathcal N(R)$. Then, $[\bar x,y,z] + [x,y,z] =[\mathsf t(x),y,z]=0$ for all $y,z\in R$. 


Now we can prove that $\bar x\in\mathcal A$ for all $x\in\mathcal A$. If $x=0$, then $\bar x=s{\cdot}x{\cdot}s^{-1}=0\in \mathcal A$. If $x\ne 0$, then the division ring $\mathcal A$ contains a multiplicative inverse $x^{-1}$ of $x$. Artin's Theorem~\ref{t:Artin} ensures that $ \bar x=x^{-1}{\cdot}(x{\cdot}\bar x)\in \mathcal A\cdot\mathcal Z(R)\subseteq \mathcal A$.
Therefore, $j:\mathcal A\to\mathcal A$, $j:x\mapsto \bar x\defeq s{\cdot}x{\cdot}s^{-1}$, is a well-defined involutive map on $\mathcal A$. The distributivity of the mutiplication over addition in the ring $R$ ensures that $j(x+y)=j(x)+j(y)$ for all $x\in\mathcal A$.

If $x,y\in\mathcal A$, then
$$x{\cdot}(s{\cdot}y)=(x{\cdot}s){\cdot}y - [x, s, y]=(s{\cdot}\bar x){\cdot}y - [s,\bar x, y]= s{\cdot}(\bar x{\cdot}y).$$
Applying the left-Bol identity, we obtain
$$(x{\cdot}y{\cdot}x){\cdot}s = x{\cdot}(y{\cdot}(x{\cdot}s)) = x{\cdot}(y{\cdot}(s{\cdot}\bar x))
= x{\cdot}(s{\cdot}(\bar y{\cdot}\bar x)) = s{\cdot}(\bar x{\cdot}\bar y{\cdot}\bar x).$$
Then $j(x{\cdot}y{\cdot}x)=s{\cdot}(x{\cdot}y{\cdot}x){\cdot}s^{-1}=(\bar x{\cdot}\bar y{\cdot}\bar x){\cdot}s{\cdot}s^{-1}=\bar x{\cdot}\bar y{\cdot}\bar x$, witnessing that that $j$ is a Jordan homomorphism of $\mathcal A$; i.e., $\overline{x{\cdot}y{\cdot}x} = \bar x{\cdot}\bar y{\cdot}\bar x$. As we noted in
the proof of Theorem~\ref{t:Mikheev},
$$p(x, y, z) \defeq (z - x{\cdot}y)\cdot (z - y{\cdot}x)$$
is a Jordan polynomial. Thus, for $x,y,z\in\mathcal A$, we have $\overline{p(x,y,z)} = p(\bar x, \bar y, \bar z)$.
Taking $z\defeq x{\cdot}y$, we have
$$0 = \overline{p(x,y,x{\cdot}y)}= p(\bar x, \bar y, \overline{x{\cdot}y})= (\overline{x{\cdot}y} - \bar x{\cdot}\bar y)\cdot(\overline{x{\cdot}y} -\bar y{\cdot}\bar x).$$
We see that for each pair $x,y\in\mathcal A$, either $\overline{x{\cdot}y} = \bar x{\cdot}\bar y$ or $\overline{x{\cdot}y} =\bar y{\cdot}\bar x$. If $\overline{x{\cdot}y} = \bar x{\cdot}\bar y$, then 
$$[\bar x,\bar y, s] =[x, y, s] =(x{\cdot}y){\cdot}s-x{\cdot}(y{\cdot}s)=s\cdot\overline{x{\cdot}y}-x{\cdot}(s{\cdot}\bar y)=s\cdot\overline{x{\cdot}y}-s(\bar x{\cdot}\bar y)=s{\cdot}(\overline{x{\cdot}y}-\bar x{\cdot}\bar y)=0.$$
Since $\mathcal A$ is associative, the subalgebra $\mathcal A_{\bar x,\bar y}$  contains $\mathcal A$ and $s$ and hence $R$.
We see that $A_{\bar x,\bar y}=0$ and $[\bar a,\bar b] = 0$ by Claim~\ref{cl:4.13}. Thus, $\overline{a{\cdot}b}= \bar b{\cdot}\bar a$ in any case.
Therefore, $j$ is a central involution on the non-commutative corps $\mathcal A$, and hence $\A$ is a quaternion corps.

Claim~\ref{cl:4.14-4.21}(1,7) implies $R=\A+\mathcal B=\A+s{\cdot}\A$. 
By Proposition~\ref{p:CD-natural}, the ring $R=\mathcal A+s{\cdot}\mathcal A$ is  isomorphic to the Cayley--Dickson extension $CD(\mathcal A,s^2)$ of the quaternion corps $\A$. So, $R$ admits a central involution, turning $R$ into an octonian ring with real axis $\Re(R)=\{x\in \mathcal Z(R):\bar x=x\}=\{x\in\mathcal Z(R):s{\cdot}x{\cdot}s^{-1}=x\}=\mathcal Z(R)=\mathcal N(R)$.

\end{proof}

Theorems~\ref{t:SBK} and Corollary~\ref{c:dim=1248} imply

\begin{corollary}\label{c:dimR=8} Every non-associative alternative division ring $R$ is a $8$-dimensional vector space over the centre $\mathcal Z(R)$.
\end{corollary}

\section{Orderable alternative division rings are associative}

\begin{definition} A ring $R$ is defined to be \index{orderable ring}\index{ring!orderable}\defterm{orderable} if there exists a subset $P\subseteq R$ such that $$(P+P)\cup(P\cdot P)\subseteq P,\quad R=(-P)\cup P\quad\mbox{and}\quad P\cap (-P)=\{0\}.$$
\end{definition}

\begin{theorem}\label{t:ab=-ba} Every non-commutative involutive alternative division ring $D$ contain two elements $a,b\in D$ with $a{\cdot}b=-b{\cdot}a\ne 0$ and hence $D$ cannot be orderable.
\end{theorem}

\begin{proof} By Proposition~\ref{p:CD-real-axis}, the real axis $R\defeq\{z\in\mathcal Z(D):\bar z=z\}$ is a subfield of $D$. If the field $R$ is Boolean, then the elements $a\defeq 1$ and $b\defeq 1$ have the required property: $a{\cdot}b=1=-1=-b{\cdot}a\ne 0$. So, assume that $R$ is not Boolean. Since $R\ne D$, we can apply Lemma~\ref{l:CD-orthogonal} and find an element $a\in D\setminus R$ such that $a+\bar a=a{\cdot}\bar 1+1{\cdot}\bar a=0$. By Proposition~\ref{p:CD-1gen}, the set $C\defeq R+a{\cdot}R$ is a subfield of the ring $D$. Since $D$ is not commutative, $D\ne C$. Applying Lemma~\ref{l:CD-orthogonal}, find an element $b\in D\setminus C$ such that $a{\cdot} b=b{\cdot}\bar a$. Then $a{\cdot} b=b{\cdot}\bar a=b{\cdot}(-a)=-a{\cdot}b$.

It remains to show that the ring $D$ is not orderable. To derive a contradiction, assume that $D$ is orderable. Then $D=P\cup(-P)$ for some subset $P\subseteq D$ such that $(P+P)\cup(P\cdot P)\subseteq P$ and $P\cap(-P)=\{0\}$. Replacing $a,b$ by $-a,-b$, if necessary, we can assume that $a,b\in P$ and hence $0\ne a{\cdot}b=-b{\cdot}a\in P\cap(-P)=\{0\}$, which is a desired contradiction showing that the ring $D$ cannot be orderable.
\end{proof}

\begin{corollary}[Bruck, Kleinfeld, 1951]\label{c:orderable-ring} Every orderable involutive alternative division ring $R$ is a field.
\end{corollary}

\begin{proof} By Theorem~\ref{t:ab=-ba}, the orderable involutive alternative division ring $R$ is commutative, and by Corollary~\ref{c:div+com+alt=>ass}, the commutative alternative division ring $R$ is associative, witnessing that $R$ is a field.
\end{proof}

\begin{corollary}[Bruck, Kleinfeld, 1951] Every orderable alternative division ring $R$ is associative.
\end{corollary}

\begin{proof} Assuming that the orderable alternative division ring $R$ is not associative, we can apply Theorem~\ref{t:SBK} and Corollary~\ref{c:orderable-ring} to conclude that $R$ admits a central involution and hence is associative.
\end{proof}

\begin{exercise} Find an example of an orderable complex field.
\smallskip

{\em Hint:} Look at the field $\IQ+\sqrt{2}\IQ$ endowed with the involution $j:x+\sqrt{2}y\mapsto x-\sqrt{2}y$.
\end{exercise}

\begin{Exercise} Find an example of an orderable non-commutative corps.
\smallskip

{\em Hint:} Such an example was constructed by David Hilbert in \cite[\S33]{Hilbert}.
\end{Exercise}

\section{Isotopic non-Boolean alternative division rings are isomorphic}

We recall that two biloops $X,Y$ are \index{isotopic biloops}\index{biloops!isotopic}\defterm{isotopic} if there exist three bijections $F,G,H:X\to Y$ such that $H(x{\cdot} y+z)=F(x){\cdot}G(y)+H(z)$ for all $x,y,z\in X$. The triple $(F,G,H)$ is called an \index{isotopism}\defterm{isotopism} of the biloops $X,Y$.

\begin{proposition}\label{p:nice-isotopy} Let $X,Y$ be isotopic biloops. If the biloop $X$ is associative-plus, then there exists an isotopism $(F,G,H):X\to Y$ such that
\begin{enumerate}
\item $H(x{\cdot}y)=F(x){\cdot}G(y)$ and $H(x+y)=H(x)+H(y)$ for all $x,y\in X$, and
\item $0=F(0)=G(0)=H(0)$.
\end{enumerate}
\end{proposition}

\begin{proof} Since the biloops $X,Y$ are isotopic, there exists an isotopism $(F,G,H'):X\to Y$ between $X$ and $Y$. Let $b\defeq H'(0)\in Y$ and $R_b:Y\to Y$, $R_b:y\mapsto y+b$, be the right shift of the plus loop of the biloop $Y$. Consider the bijection $H\defeq R_b^{-1}H'$ and observe that $H(0)=R_b^{-1}H'(0)=R_b^{-1}(b)=0$ and $H'=R_bH$. We claim that $(F,G,H)$ is an isotopism of $X$ onto $Y$. Indeed, for any $x,y,z\in X$ we have
$$
\begin{aligned}
H(x{\cdot}y+z)+b&=R_bH(x{\cdot}y+z)=H'(x{\cdot}y+z)=F(x){\cdot} G(y)+H'(z)=F(x){\cdot}G(y)+R_bH(z)\\
&=F(x){\cdot}G(y)+(H(z)+b)=(F(x){\cdot}G(y)+H(z))+b,
\end{aligned}
$$ by the associativity-plus and hence $H(x{\cdot}y+z)=F(x){\cdot}G(y)+H(z)$, by the cancellativity of the plus operation in the biloop $X$. The latter equality witnesses that the triple $(F,G,H)$ is an isotopism of the biloops $X,Y$.

Now observe that for any $x,y\in X$ we have
$$H(x{\cdot}y)=H(x{\cdot}y+0)=F(x){\cdot}G(y)+H(0)=F(x){\cdot}G(y)+0=F(x){\cdot} G(y)$$and 
$$H(x+y)=H((x{\cdot}1)+y)=F(x)\cdot G(1)+H(y)=H(x{\cdot}1)+H(y)=H(x)+H(y).$$

The definition of $H$ implies $H(0)=0$. Assuming that $F(0)\ne 0$, we conclude that $H(0)=H(0{\cdot}0)=F(0){\cdot} G(0)\ne F(0){\cdot }G(1)=H(0{\cdot} 1)=H(0)$, which is a contradiction showing that $F(0)=0$. By analogy we can prove that $G(0)=0$.
\end{proof}

We recall that a  ring $R$ is called \index{Boolean ring}\index{ring!Boolean}\defterm{Boolean} if $x+x=0$ for all $x\in R$. This happens if and only if $2\defeq 1+1=0$ in $R$, i.e., $R$ has characteristic $2$.

The main result of this section is the following difficult theorem of \index[person]{Schafer}Schafer\footnote{{\bf Richard Donald Schafer} (1918 -- 2014) was an American mathematician. Richard Schafer studied at the University at Buffalo, where he received his bachelor's degree in 1938 and his master's degree in 1940. He received in 1942 from the University of Chicago his PhD under Abraham Adrian Albert with dissertation ``Alternative Algebras over an Arbitrary Field''. After service in the U.S. Naval Reserve from 1942 to 1945, he was an instructor at the University of Michigan for the academic year 1945–1946. From 1946 to 1948 he was at the Institute for Advanced Study. From 1948 to 1953 he was a professor at the University of Pennsylvania. From 1953 to 1958 he was at University of Connecticut as professor and head of the mathematics department. He spent the academic year 1958–1959 at the Institute for Advanced Study. From 1959 until his retirement in 1988, he was a professor at the Massachusetts Institute of Technology. In 2012 he was elected was a Fellow of the American Mathematical Society. Schafer did research on algebra, specifically on Jordan algebras and Lie algebras. He is best known for his textbook ``An Introduction to Nonassociative Algebras'', first published in 1966. He also studied the Cayley--Dickson construction.} \cite{Schafer1943}.

\begin{theorem}[Schafer, 1943]\label{t:Schafer} Two non-Boolean alternative division rings are isotopic if and only if they are isomorphic.
\end{theorem} 

\begin{proof} Let us observe that alternative division rings are  alternative-dot associative-plus commutative-plus distributive biloops. The ``if'' part of Schafer's theorem is trivial. So, assume that $X,Y$ are two isotopic non-Boolean alternative division rings. Since every ring is associative-plus, we can apply Proposition~\ref{p:nice-isotopy} and find an isotopism $(F,G,H):X\to Y$ such that  $H(x\cdot y)=F(x)\cdot G(y)$ and $H(x+y)=H(x)+H(y)$ for all $x,y\in H$, and $F(0)=G(0)=H(0)=0$.  
Since $F,G$ are bijections, there exist elements $f,g\in X$ such that $F(f)=1=G(g)$, where $1$ is the multiplicative unit of the ring $Y$. Since $F(0)=0=G(0)$, the elements $f,g\in X$ are distinct from zero and hence $f,g$ have multiplicative inverses $f^{-1},g^{-1}$ in the alternative division ring $X$.

\begin{claim}\label{cl:Hxy=} For any $x,y\in X$, 
\begin{enumerate}
\item $H(x{\cdot}y)=H(x{\cdot}g)\cdot H(f{\cdot}y)$;
\item $H(x)\cdot H(y)=H\big((x{\cdot}g^{-1}){\cdot}(f^{-1}{\cdot}y)\big).$
\end{enumerate}
\end{claim}

\begin{proof} 1. Observe that  $H(f{\cdot}y)  =  F(f){\cdot}G(y)= 1{\cdot} G(y)=G(y)$  and 
$H(x{\cdot}g)  = F(x){\cdot}G(g)=F(x){\cdot}1=F(x)$ and hence $H(x{\cdot}y)  = F(x){\cdot}G(y)=H(x{\cdot}g)\cdot H(f{\cdot}y)$.
\smallskip

2. Since the alternative ring $X$ is inversive-dot, for $x'\defeq x{\cdot}g^{-1}$ and $y'\defeq f^{-1}{\cdot}y$, the preceding item implies $$H\big((x{\cdot}g^{-1}){\cdot}(f^{-1}{\cdot}y)\big)=H(x'{\cdot}y')=H(x'{\cdot}g){\cdot}H(f{\cdot}y')=H((x{\cdot}g^{-1}){\cdot}g){\cdot}H(f{\cdot}(f^{-1}{\cdot}y))=H(x){\cdot}H(y).$$
\end{proof}

Since the elements $f,g$ of the division ring $X$ are distinct from zero, their product $e\defeq f{\cdot}g\in X$ also is not equal to zero. Then $H(e)\ne H(0)=0$.
 


Consider the bijection $\Phi:X\to Y$, $\Phi:x\mapsto H(x{\cdot}e)$.  The following claim shows that the bijection $\Phi$ is an isomorphism of the plus loops of the biloops $X,Y$.

\begin{claim}\label{cl:Phi+} $\Phi(x+y)=\Phi(x)+\Phi(y)$ for all $x,y\in X$.
\end{claim}

\begin{proof} Observe that $\Phi(x+y)=H((x+y){\cdot}e)=H(x{\cdot}e+y{\cdot}e)=H(x{\cdot}e)+H(y{\cdot}e)=\Phi(x)+\Phi(y)$.
\end{proof}

\begin{claim}\label{cl:Phi-multhom} The equality $\Phi(x{\cdot}y)=\Phi(x){\cdot} \Phi(y)$ holds for all $x,y$ that belong to some associative subring $A\subseteq X$ that contains the elements $f^{-1},g^{-1}$ and $e$.
\end{claim}

\begin{proof} The inclusions $f^{-1},g^{-1}\in A$ imply $e^{-1}=(f{\cdot}g)^{-1}=g^{-1}{\cdot}f^{-1}\in A$. The associativity of the multiplication in the ring $A$ ensures that $(x{\cdot}g^{-1}){\cdot}(f^{-1}{\cdot}y)  =  x{\cdot}(g^{-1}{\cdot}f^{-1}){\cdot} y=x{\cdot}e^{-1}{\cdot}y$. By Claim~\ref{cl:Hxy=}, 
$$H(x)\cdot H(y)=H((x{\cdot}g^{-1}){\cdot}(f^{-1}{\cdot}y))=H(x{\cdot}e^{-1}{\cdot}y)$$for all $x,y\in A$. 
Then
$$\Phi(x{\cdot}y)=H(x{\cdot}y{\cdot}e)=H(x{\cdot}e{\cdot}e^{-1}{\cdot}y{\cdot}e)=H(x{\cdot}e){\cdot}H(y{\cdot}e)=\Phi(x)\cdot \Phi(y)$$ for all $x,y\in A$.
\end{proof}

If the ring $X$ is associative, then Claims~\ref{cl:Phi+} and \ref{cl:Phi-multhom} ensure that $\Phi$ is an isomorphism of the rings $X,Y$ and we are done. So, assume that the ring $X$ is not associative. 
By Theorem~\ref{t:SBK}, the non-associative alternative division ring $X$ admits a central involution $j:X\to X$, $j:x\mapsto\bar x,$ whose real axis $R\defeq\{z\in\mathcal Z(X):\bar z=z\}$ coinsides with the centre $\mathcal Z(X)$ of the ring $X$. So, we can think of $X$ as an involutive ring with real axis $R=\mathcal Z(X)$.
 
\begin{claim}\label{cl:fginQ} The elements $f,g$ are contained in some quaternion subring of the involutive ring $X$. 
\end{claim}

\begin{proof} By Artin's Theorem~\ref{t:Artin}, the set $\{f,g\}\cup R$ generates an associative subring $A$ of the alternative ring $X$. Proposition~\ref{p:CD-1gen} implies that the ring $A$ is divisible. Observe that for every $x\in A$ we have $\bar x=x^{-1}{\cdot}(x\cdot\bar x)\in A^{-1}\cdot R\subseteq A$, which implies that $A$ is an involutive subring of the involutive ring $X$. If $A$ is not commutative, then $A$ is a required quaternion subring of the ring $X$. So, assume that $A$ is commutative, which implies $A\ne X$. If $A$ contains an element $a\ne \bar a$, then $A$ is a complex subfield of $X$. Applying Lemma~\ref{l:CD-orthogonal}, we can find an element $q\in X\setminus A$ such that $q^2\in R$ and $x{\cdot}q=q{\cdot}\bar x$ for all $x\in X$. By Proposition~\ref{p:CD-natural} and Theorem~\ref{t:Cayley-Dickson}, the set $Q=A+q{\cdot}A$ is a quaternion subring of $X$, isomorphic to the Cayley--Dickson extension $CD(A,q^2)$ of the complex field $A$. It is clear that $f,g\in A\subseteq Q$.

Finally, assume that $\bar x=x$ for all $x\in A$. Then $2{\cdot}x=x+x=x+\bar x\in R$ and $x\in R$ for all $x\in A$ (because the ring $X$ is not Boolean and hence $2\defeq 1+1\ne 0$ in the field $R$). Therefore, $A=R$. By Theorem~\ref{t:248CD}, the involutive ring $X$ is isomorphic to the Cayley--Dickson extension of some quaternion corps and hence $X$ contains a quaternion subcorps $Q\subseteq X$ such that $f,g\in A=R\subseteq Q$.
\end{proof}

By Claim~\ref{cl:fginQ}, the elements $f,g$ are contained in some quaternion subcorps $Q$ of the division ring $X$. Repeating the argument of the proof of Theorem~\ref{t:248CD}, we can show that $X=Q+i{\cdot}Q$ for some element $i\in X\setminus Q$ such that $\bar i=-i$ and $i{\cdot}x=\bar x{\cdot}i$ for all $x\in Q$. Then $i^2=-i{\cdot}\bar i\in R$.
 Applying Proposition~\ref{p:CD-natural}, we obtain the following identities which will help us to make calculations in the alternative ring $Q+i{\cdot}Q=X$.

\begin{claim}\label{cl:Q-operations} If $x,y\in Q$, then
\begin{enumerate}
\item $(i{\cdot}x){\cdot} y=i{\cdot}(y{\cdot}x)$,
\item $x{\cdot}(i{\cdot}y)=i{\cdot}(\bar x{\cdot}y)$, and
\item $(i{\cdot}x){\cdot}(i{\cdot}y)=i^2{\cdot}(y\cdot\bar x)$.
\end{enumerate}
\end{claim}

Consider the function $\Psi:X\to Y$, defined by the formula $\Psi(x+i{\cdot}y)=\Phi(x)+\Phi(i){\cdot}\Phi(y)$ for $x,y\in Q$. We claim that $\Psi$ is an isomorphism of the rings $X$ and $Y$. 

 Claim~\ref{cl:Phi-multhom} implies that the following multiplicativity property of the map $\Psi$.

\begin{claim} $\Psi(x{\cdot}y)=\Psi(x)\cdot\Psi(y)$ for every $x,y\in Q$.
\end{claim}

\begin{claim}\label{cl:Psi+} $\Psi(x+y)=\Psi(x)+\Psi(y)$ for all $x,y\in X$.
\end{claim}

\begin{proof}  Given any elements $x,y\in X$, find elements $x',x'',y',y''\in Q$ such that $x=x'+i{\cdot}x''$ and $y=y'+i{\cdot}y''$. Claim~\ref{cl:Phi+} implies 
$$\begin{aligned}
\Psi(x+y)&=\Psi((x'+y')+i{\cdot}(x''+y''))=\Phi(x'+y')+\Phi(i){\cdot}\Phi(x''+y'')\\
&=(\Phi(x')+\Phi(y'))+\Phi(i)\cdot(\Phi(x'')+\Phi(y''))=(\Phi(x')+i{\cdot}\Phi(x''))+(\Phi(y')+i{\cdot}\Phi(y''))\\
&=\Phi(x)+i{\cdot}\Phi(y).
\end{aligned}
$$ 
\end{proof}

\begin{claim} The function $\Psi$ is injective.
\end{claim}

\begin{proof} Since $\Psi$ is a homomorphism of the additive groups of the rings $X,Y$, it suffices to show that every element $x+i{\cdot}y\in \Psi^{-1}(0)$ is zero. The choice of $i\in X\setminus Q$ implies $\Phi(i)\notin \Phi[Q]$. Claims~\ref{cl:Phi+} and \ref{cl:Phi-multhom} ensure that $\Phi[Q]$ is an associative subring of the alternative ring $Y$. It follows from $\Phi(i)\notin\Phi[Q]$ that $\Phi[Q]\cap (\Phi(i){\cdot}\Phi[Q])=\{0\}$. Then  $0=\Psi(x+i{\cdot}y)=\Phi(x)+\Phi(i){\cdot}\Psi(y)$ implies $\Phi(x)=0=\Phi(y)$ and hence $x=0=y$ because $\Phi$ is an isomorphism of the groups $(X,+)$, $(Y,+)$, ccording to Claim~\ref{cl:Phi+}. Therefore, $x+i{\cdot}y=0$ and the homomorphism $\Psi$ is injective.
\end{proof}

\begin{claim} The function $\Psi:X\to Y$ is bijective.
\end{claim}

\begin{proof} Since the function $\Phi:X\to Y$ is bijective, for every element $y\in Y$, there exist elements $x,z\in Q$ such that $\Phi(x+i{\cdot}z)=y$. Consider the element $a=e{\cdot}z{\cdot}e^{-1}{\cdot}g\in Q$. Claim~\ref{cl:Q-operations} implies 
$$((i{\cdot}e){\cdot}g^{-1}){\cdot}a=(i{\cdot}(g^{-1}{\cdot}e)){\cdot}a=i{\cdot}(a{\cdot}g^{-1}{\cdot}e)=i{\cdot}(e{\cdot}z{\cdot}e^{-1}{\cdot}g{\cdot}g^{-1}{\cdot}e)=i{\cdot}(e{\cdot}z)=(i{\cdot}z){\cdot}e.$$ Applying Claim~\ref{cl:Hxy=}, we conclude that $$\Phi(i{\cdot}z)=H((i{\cdot}z){\cdot}e)=H(((i{\cdot}e){\cdot}g^{-1}){\cdot} a)=H(i{\cdot}e)\cdot H(f{\cdot}a)=\Phi(i){\cdot}\Phi(f{\cdot}a{\cdot}e^{-1}).$$ Then for the element $x+i{\cdot}(f{\cdot}a{\cdot}e^{-1})\in Q+i{\cdot}Q=X$, we have 
$$\Psi(x+i{\cdot}(f{\cdot}a{\cdot}e^{-1}))\defeq\Phi(x)+\Phi(i){\cdot}\Phi(f{\cdot}a{\cdot}e^{-1})=\Phi(x)+\Phi(i{\cdot}z)=\Phi(x+i{\cdot}z)=y,$$
witnessing that the injective function $\Psi$ is surjective and hence bijective.
\end{proof}

\begin{claim}\label{cl:Phi((ix)y)} For any $x,y\in Q$, we have 
$$(\Phi(i){\cdot}\Phi(x))\cdot\Phi(y)=\Phi(i)\cdot\Phi(y{\cdot}x).$$
\end{claim}

\begin{proof} By Claim~\ref{cl:Hxy=},$$
\begin{aligned} 
(\Phi(i){\cdot}\Phi(x)){\cdot}\Phi(y)&=(H(i{\cdot}e){\cdot}H(x{\cdot}e)){\cdot}H(y{\cdot}e)=H((i{\cdot}e){\cdot}g^{-1}){\cdot}(f^{-1}{\cdot}x{\cdot}e)){\cdot}H(y{\cdot}e)\\
&=H((i{\cdot}(g^{-1}{\cdot}e)){\cdot}(f^{-1}{\cdot}x{\cdot}e)){\cdot}H(y{\cdot}e)=H(i{\cdot}(f^{-1}{\cdot}x{\cdot}e{\cdot}g^{-1}{\cdot}e)){\cdot}H(y{\cdot}e)\\
&=H((i{\cdot}(f^{-1}{\cdot}x{\cdot}f{\cdot}e)){\cdot}g^{-1}){\cdot}(f^{-1}{\cdot}y{\cdot}e))=H((i{\cdot}(g^{-1}{\cdot}f^{-1}{\cdot}x{\cdot}f{\cdot}e)){\cdot}(f^{-1}{\cdot}y{\cdot}e))\\
&=H(i{\cdot}(f^{-1}{\cdot}y{\cdot}e{\cdot}g^{-1}{\cdot}f^{-1}{\cdot}x{\cdot}f{\cdot}e))=H(i{\cdot}(f^{-1}{\cdot}y{\cdot}x{\cdot}f{\cdot}e))
\end{aligned}
$$
and 
$$
\begin{aligned} 
\Phi(i){\cdot}\Phi(y{\cdot} x)&=H(i{\cdot}e){\cdot}H(y{\cdot}x{\cdot}e)=
H(((i{\cdot}e){\cdot}g^{-1}){\cdot}(f^{-1}{\cdot}y{\cdot}x{\cdot}e))=H((i{\cdot}(g^{-1}{\cdot}e)){\cdot}(f^{-1}{\cdot}y{\cdot}x{\cdot}e))\\
&=H(i{\cdot}(f^{-1}{\cdot}y{\cdot}x{\cdot}e{\cdot}g^{-1}{\cdot}e))=
H(i{\cdot}(f^{-1}{\cdot}y{\cdot}x{\cdot}f{\cdot}e)).
\end{aligned}
$$
\end{proof}

\begin{claim}\label{cl:Phi(x(iy))} For any $x,y\in Q$, we have 
$$\Phi(x){\cdot}(\Phi(i){\cdot}\Phi(y))=\Phi(i){\cdot}\Phi(\bar x{\cdot}y).$$
\end{claim}

\begin{proof} Since the involution is an anti-isomorphism of the corps $Q$, we have $\overline{f^{-1}}=\bar f^{-1}$. By Claim~\ref{cl:Hxy=} and the inclusion $\bar f^{-1}f^{-1}=\overline{f^{-1}}f^{-1}\in\mathcal Z(Q)$, we have $$
\begin{aligned} 
\Phi(x){\cdot}(\Phi(i){\cdot}\Phi(y))&=H(x{\cdot}e){\cdot}(H(i{\cdot}e){\cdot}H(y{\cdot}e))=H(x{\cdot}e){\cdot}H\big(((i{\cdot}e){\cdot}g^{-1}){\cdot}(f^{-1}{\cdot}y{\cdot}e)\big)\\
&=H(x{\cdot}e){\cdot}H\big((i{\cdot}(g^{-1}{\cdot}e)){\cdot}(f^{-1}{\cdot}y{\cdot}e)\big)=H(x{\cdot}e){\cdot}H\big(i{\cdot}(f^{-1}{\cdot}y{\cdot}e{\cdot}g^{-1}{\cdot}e)\big)\\
&=H(x{\cdot}e){\cdot}H(i{\cdot}(f^{-1}{\cdot}y{\cdot}f{\cdot}e))=H\big((x{\cdot}e{\cdot}g^{-1}){\cdot}(f^{-1}{\cdot}(i{\cdot}(f^{-1}{\cdot}y{\cdot}f{\cdot}e)))\big)=\\
&=H\big((x{\cdot}f){\cdot}(i{\cdot}(\overline{f^{-1}}{\cdot}f^{-1}{\cdot}y{\cdot}f{\cdot}e)\big)=H\big(i{\cdot}(\overline{x{\cdot}f}{\cdot}\bar f^{-1}{\cdot}f^{-1}{\cdot}y{\cdot}f{\cdot}e)\big)\\
&=H\big(i{\cdot}(\bar f{\cdot}\bar x{\cdot}\bar f^{-1}{\cdot}f^{-1}{\cdot}y{\cdot}f{\cdot}e)\big)=H\big(i{\cdot}(\bar f{\cdot}\bar f^{-1}{\cdot}f^{-1}{\cdot}\bar x{\cdot}y{\cdot}f{\cdot}e)\big)\\
&=H\big(i{\cdot}(f^{-1}{\cdot}\bar x{\cdot}y{\cdot}f{\cdot}e)\big)
\end{aligned}
$$
and 
$$
\begin{aligned} 
\Phi(i){\cdot}\Phi(\bar x{\cdot}y)&=H(i{\cdot}e){\cdot}H(\bar x{\cdot}y{\cdot}e)=H\big(((i{\cdot}e){\cdot}g^{-1}){\cdot}(f^{-1}{\cdot}\bar x{\cdot}y{\cdot}e)\big)\\
&=H\big((i{\cdot}(g^{-1}{\cdot}e))(f^{-1}{\cdot}\bar x{\cdot}y{\cdot}e)\big)=H\big(i{\cdot}(f^{-1}{\cdot}\bar x{\cdot}y{\cdot}e{\cdot}g^{-1}{\cdot}e)\big)\\
&=H\big(i{\cdot}(f^{-1}{\cdot}\bar x{\cdot}y{\cdot}f{\cdot}e)\big).
\end{aligned}
$$
\end{proof}

\begin{claim}\label{cl:Phi:ixiy} For any $x,y\in Q$, we have 
$$(\Phi(i){\cdot}\Phi(x))\cdot(\Phi(i){\cdot}\Phi(y))=\Phi(i^2{\cdot}y{\cdot}\bar x).$$
\end{claim}

\begin{proof} By Claims~\ref{cl:Q-operations}, \ref{cl:Hxy=} and the inclusion $e{\cdot}\bar e,f{\cdot}\bar f\in\mathcal Z(Q)$, we have $$
\begin{aligned} 
(\Phi(i){\cdot}\Phi(x))\cdot(\Phi(i){\cdot}\Phi(y))&=(H(i{\cdot}e){\cdot}H(x{\cdot}e))\cdot(H(i{\cdot}e){\cdot}H(y{\cdot}e))\\
&=H\big(((i{\cdot}e){\cdot}g^{-1}){\cdot}(f^{-1}{\cdot}x{\cdot}e)\big)\cdot H\big(((i{\cdot}e){\cdot}g^{-1}){\cdot}(f^{-1}{\cdot}y{\cdot}e)\big)\\
&=H\big((i{\cdot}(g^{-1}{\cdot}e)){\cdot}(f^{-1}{\cdot}x{\cdot}e)\big)\cdot H\big(((i{\cdot}(g^{-1}{\cdot}e)){\cdot}(f^{-1}{\cdot}y{\cdot}e)\big)\\
&=H\big(i{\cdot}(f^{-1}{\cdot}x{\cdot}e{\cdot}g^{-1}{\cdot}e)\big)\cdot H\big(i{\cdot}(f^{-1}{\cdot}y{\cdot}e{\cdot}g^{-1}{\cdot}e)\big)\\
&=H\big((i{\cdot}(f^{-1}{\cdot}x{\cdot}f{\cdot}e)){\cdot}g^{-1}){\cdot}(f^{-1}{\cdot}(i{\cdot}(f^{-1}{\cdot}y{\cdot}f{\cdot}e))\big)\\
&=H\big((i{\cdot}(g^{-1}{\cdot}f^{-1}{\cdot}x{\cdot}f{\cdot}e)){\cdot}(i{\cdot}(\bar f^{-1}{\cdot}f^{-1}{\cdot}y{\cdot}f{\cdot}e)\big)\\
&=H\big(i^2{\cdot}\bar f^{-1}{\cdot}f^{-1}{\cdot}y{\cdot}f{\cdot}e{\cdot}\bar e{\cdot}\bar f{\cdot}\bar x{\cdot}\bar e^{-1}\big)\\
&=H\big(i^2{\cdot}\bar f^{-1}{\cdot}f^{-1}{\cdot}f\bar f{\cdot}y{\cdot}\bar x{\cdot}e{\cdot}\bar e{\cdot}\bar e^{-1}\big)=H\big(i^2{\cdot}y{\cdot}\bar x{\cdot}e\big)=\Phi(i^2{\cdot}y{\cdot}\bar x).
\end{aligned}
$$
\end{proof}

\begin{claim}\label{cl:Psi-dot} For every $x,y\in X$ we have $\Psi(x{\cdot}y)=\Psi(x){\cdot}\Psi(y)$.
\end{claim}

\begin{proof} Since $x,y\in X=Q+i{\cdot}Q$, there exist elements $x_1,x_2,y_1,y_2\in Q$ such that $x=x_1+i{\cdot}x_2$ and $y=y_1+i{\cdot}y_2$. Applying Claims~\ref{cl:Phi((ix)y)}, \ref{cl:Phi(x(iy))}, \ref{cl:Phi:ixiy}, \ref{cl:Phi+}, \ref{cl:Phi-multhom}, we conclude that
$$
\begin{aligned}
x{\cdot}y&=(x_1+i{\cdot}x_2){\cdot}(y_1+i{\cdot}y_2)=x_1{\cdot}y_1+(i{\cdot}x_2){\cdot}(i{\cdot}y_2)+(i{\cdot}x_2){\cdot}y_1+x_1{\cdot}(i{\cdot}y_2)\\
&=x_1{\cdot}y_1+i^2{\cdot}y_2{\cdot}\overline{x_2}+i{\cdot}(y_1{\cdot}x_2+\overline{x_1}{\cdot}y_2)
\end{aligned}
$$
and
$$
\begin{aligned}
\Psi(x)\cdot\Psi(y)&=\Psi(x_1+i{\cdot}x_2)\cdot\Psi(y_1+i{\cdot}y_2)=(\Phi(x_1)+\Phi(i){\cdot}\Phi(x_2))\cdot(\Phi(y_1)+\Phi(i){\cdot}\Phi(y_2))\\
&=\Phi(x_1){\cdot}\Phi(y_1)+(\Phi(i){\cdot}\Phi(x_2)){\cdot}(\Phi(i){\cdot}\Phi(y_2))+(\Phi(i){\cdot}\Phi(x_2)){\cdot}\Phi(y_1)+\Phi(x_1){\cdot}(\Phi(i){\cdot}\Phi(y_2))\\
&=\Phi(x_1{\cdot}y_1)+\Phi(i^2{\cdot}y_2{\cdot}\overline{x_2})+\Phi(i){\cdot}\Phi(y_1{\cdot}x_2)+\Phi(i){\cdot}\Phi(\overline{x_1}{\cdot}y_2)\\
&=\Phi(x_1{\cdot}y_1+i^2{\cdot}y_2{\cdot}\overline{x_2})+\Phi(i){\cdot}\Phi(y_1{\cdot}x_2+\overline{x_1}{\cdot}y_2)\\
&=\Psi((x_1{\cdot}y_1+i^2{\cdot}y_2{\cdot}\overline{x_2})+i{\cdot}(y_1{\cdot}x_2+\overline{x_1}{\cdot}y_2))=\Psi(x{\cdot}y).
\end{aligned}
$$
\end{proof}

Claims~\ref{cl:Psi+} and \ref{cl:Psi-dot} show that $\Psi:X\to Y$ is an isomorphism of the rings $X,Y$.
\end{proof}

\begin{problem}\label{prob:Schafer-Boolean} Are any isotopic Boolean alternative division rings isomorphic?
\end{problem}

\begin{remark} One could guess that Problem~\ref{prob:Schafer-Boolean} was resolved already in the original paper of Schafer \cite{Schafer1943}, but it seems that the standard proof of this theorem (presented also in the textbook  \cite[14.3.5]{Stevenson}) is not complete in the Boolean case.
\end{remark}

\section{Cayley-Dickson octonion rings of characteristic $2$}

In this section we study the structure of Boolean octonion rings .

\begin{proposition}\label{p:Boolean-CD1} If $X$ is a Boolean octonion ring with real axis $R$, then there exist  elements $a,b,c\in X$ such that 
\begin{enumerate}
\item $\bar a=a$, $\bar b=b$ and $c+\bar c=1$;
\item $C\defeq R+cR$ is a complex field;
\item $Q\defeq C+bQ$ is a quatenion corps;
\item $A\defeq R+aR+bR+abR$ is a real subfield of $X$;
\item $bx=\bar xb$ for all $x\in C$;
\item $ax=\bar xa$ for all $x\in Q$;
\item $Q+aQ=X=A+cA$.
\end{enumerate}
\end{proposition}

\begin{proof} Since the involutive ring $X$ is octonion, it is not real and hence there exists an element $c\in X$ such that $\bar c\ne c$. Since the field $R$ has characteristic $2$, $c+c=0$ and $\bar c\ne c$ imply $\mathsf{t}(c)=\bar c+c\ne 0$. Multiplying $c$ by $\mathsf{t}(c)^{-1}$, we can assume that $\bar c+c=1$.  Propositions~\ref{p:CD-1gen} ensure that $C\defeq R+cR$ is a complex subfield of the involutive ring $X$. Applying Lemma~\ref{l:CD-orthogonal}, find an element $b\in X\setminus C$ such that $\bar b=b$ and $x\bar b=bx$ for all $x\in C$. Since $cb=b\bar c\ne bc$, the subring $Q=C+bC$ is a quaternion corps in $X$.  Applying Lemma~\ref{l:CD-orthogonal} again, find an element $a\in X\setminus Q$ such that $\bar a=a$ and $xa=a\bar x$ for all $x\in Q$. In particular, $ba=a\bar b=ab$, which measn that the elements $a,b$ commute. Taking into account that $\dim_R(X)=8=\dim_R(Q+aQ)$, we conclude that $X=Q+aQ$.
Observe that $A\defeq R+bR+aR+abR$ is a subfield of the ring $X$ such that $\bar x=x$ for all $x\in A$, which implies that $c\notin A$ and $A\cap cA=\{0\}$. Taking into account that $\dim_R(X)=8=\dim_R(A+cA)$, we conclude that $X=A+cA$.
\end{proof}

Now we elaborate some tools simplifying calculations in a Boolean octonion ring.
Given a Boolean octonion ring $X$, find elements $a,b,c\in X$ satisfying the conditions of Proposition~\ref{p:Boolean-CD1}. Consider the real subfield $A=R+aR+bR+abR$ of $X$. Any element $x\in A$ can be uniquely written in the form $x=\Re(x)+\vec x$ for some elements $\Re(x)\in R$ and $\vec x\in aR+bR+abR$, called the \defterm{real} and \defterm{vector} parts of the element $x$, respectively. For two elements $x,y\in A$, the real part $\Re(x{\cdot}y)\in R$ of the product $x{\cdot}y\in A$ is called the \defterm{scalar product} of the elements $x,y$ and the element
$$x\times y\defeq x{\cdot}y+\Re(x{\cdot}y)+x{\cdot}\Re(y)+\Re(x){\cdot}y$$is called the \defterm{vector product} of the elements $x,y$.

\begin{exercise} Show that for any distinct basic vectors $x,y\in\{a,b,ab\}$, we have $x\times y=x{\cdot}y$.
\end{exercise}

For two elements $x,y\in X$, let $[x,y]\defeq x{\cdot}y-y{\cdot}x=x{\cdot}y+y{\cdot}x$ be the \defterm{commutator} of the elements $x,y$ in the Boolean ring $X$. 

\begin{proposition}\label{p:commutator-Boolean}The equality $[c,x]=x+\Re(x)$ holds for all $x\in A$.
\end{proposition}

\begin{proof} It suffices to check that this equality holds for all basic elements $x\in\{1,a,b,ab\}$ in the $R$-submodule $A$ of $X$. For $x=1$ this equality follows from the equality $\Re(1)=1$. For $x\in a,b$, we have $\Re(x)=0$ and $cx=x\bar c=x(c+1)=xc+x=xc+x+\Re(x)$, by Proposition~\ref{p:Boolean-CD1}. For $x=ab$, we have
$c(ab)=a(\bar cb)=a((c+1)b)=a(cb)+ab=(ab)c+ab+\Re(ab)$, 
by Propositions~\ref{p:Boolean-CD1} and \ref{p:CD-natural}
\end{proof} 

For three elements $x,y,z\in X$ we denote by $[x,y,z]\defeq(x{\cdot}y){\cdot}z-x{\cdot}(y{\cdot}z)=(x{\cdot}y){\cdot}z+x{\cdot}(y{\cdot}z)$ their associator in the Boolean ring $X$. By Proposition~\ref{p:alternative9}, $$[x,y,z]=[x,z,y]=[y,x,z]=[y,z,x]=[z,x,y]=[z,y,x].$$

\begin{proposition}\label{p:associator-Boolean} For any elements $x,y\in A$, we have the equality
$[c,x,y]=x\times y$.
\end{proposition}

\begin{proof} Since the functions $[c,x,y]$ and $x\times y$ are bilinear (on the arguments $x,y$), it suffices to check the equality $(c{\cdot}x){\cdot}y+c{\cdot}(x{\cdot}y)=x\times y$ for all basic vectors $x,y\in\{1,a,b,ab\}$. This equality trivialy holds if $1\in\{x,y\}$. If $x=y$, then $[c,x,y]=0=x\times y$. It remains to check that $(c{\cdot}x){\cdot}y+c{\cdot}(x{\cdot}y)=x\times y$ for distinct basic vectors $x,y\in\{a,b,ab\}$. Since $[c,x,y]=[c,y,x]$ and $x\times y=y\times x$, it suffices to consider the case of $(x,y)\in\{(a,b),(a,ab),(b,ab)\}$. In the following formulas we omit the notations of multiplication operation in the ring $X$. Applying Propositions~\ref{p:Boolean-CD1} and \ref{p:CD-natural}, we conclude that
$$
\begin{aligned}
&(ca)b+c(ab)=(a\bar c)b+a(\bar cb)=a(b\bar c)+a(bc)=a(b(\bar c+c))=ab=a\times b,\\
&(ca)(ab)+c(a(ab))=(a\bar c)(ab)+cba^2=a^2(bc)+a^2(b\bar c)=a^2(b(c+\bar c))=a^2b=a\times ab\\
&(cb)(ab)+c(b(ab))=a(\overline{cb}b)+cab^2=a(cb^2)+cab^2=\bar cab^2+cab^2=ab^2=b\times ab.
\end{aligned}
$$
\end{proof}

The following proposition allows us to fulfil the multiplication in the Boolean octonion ring $X=A+cA$.

\begin{proposition} For any elements $x,y\in A$, we have the equalities
\begin{enumerate}
\item $(c{\cdot}x){\cdot}y=c{\cdot}(x{\cdot}y)+(x\times y)$;
\item $x{\cdot}(c{\cdot}y)=c{\cdot}(x{\cdot}y)+x{\cdot}y+\Re(x){\cdot}y$;
\item $(c{\cdot}x){\cdot}(c{\cdot}y)=(x{\cdot}y{\cdot}\mathsf n(c)+x\times y)+c{\cdot}(x\times y+\Re(x){\cdot}y)$.
\end{enumerate}
\end{proposition}

\begin{proof} 1. The first item has been proved in Proposition~\ref{p:associator-Boolean}.
\smallskip

2. By Propositions~\ref{p:associator-Boolean} and \ref{p:commutator-Boolean}, we have
$$
\begin{aligned}
x(cy)&=(xc)y+x\times y=(cx+x+\Re(x))y+x\times y=(cx)y+xy+\Re(x)y+x\times y\\
&=c(xy)+x\times y+xy+\Re(x)y+x\times y=c(xy)+xy+\Re(x)y.
\end{aligned}
$$

3. By Propositions~\ref{p:associator-Boolean}, \ref{p:commutator-Boolean}, and the Moufang identity,
$$
\begin{aligned}
(cx)(cy)&=(cx)(yc+y+\Re(y))=(cx)(yc)+(cx)y+cx\Re(y)
\\
&=c((xy)c)+c(xy)+x\times y+cx\Re(y)\\
&=c((xy)c+xy)+x\times y+cx\Re(y)\\
&=c(c(xy)+\Re(xy))+x\times y+cx\Re(y)\\
&=c^2(xy)+c\Re(xy)+x\times y+cx\Re(y)\\
&=(c(\bar c+1)){\cdot} (xy)+x\times y+c\Re(xy)+cx\Re(y)\\
&=\mathsf{n}(c){\cdot}(xy)+x\times y+c(xy+\Re(xy)+x\Re(y))\\
&=xy{\cdot}\mathsf{n}(c)+x\times y+c(x\times y+\Re(x)y).
\end{aligned}
$$
\end{proof}

\chapter{Affine transformations of based affine planes}

In this section we study the interplay between homogeneity properties of a based affine plane and the algebraic structure of its ternar.
The homogeneity properties are expressed in terms of (para)central automorphisms of affine planes. 

Let us recall that an automorphism $A:X\to X$ of an affine space $X$ is   
\begin{enumerate}
\item \index{automorphism!central}\index{central automorphism}\defterm{central} if there exists a point $c\in X$ such that $A[\Aline xc]=\Aline xc$ for all $x\in X$;
\item \index{automorphism!paracentral}\index{paracentral automorphism}\defterm{paracentral} if for every pairs $xx',yy'\in A\setminus 1_A$, the lines $\Aline x{x'}$ and $\Aline y{y'}$ are parallel;
\item \index{automorphism!hyperfixed}\index{hyperfixed automorphism}\defterm{hyperfixed} if the set $\Fix(A)\defeq\{x\in X:A(x)=x\}$ contains some hyperplane $H$ in $X$.
\end{enumerate}
By Theorem~\ref{t:central<=>homothety}, an automorphism $A$ of an affine plane $X$ is central if and only if $A$ is a homothety. By Proposition~\ref{p:hyperfixed=>paracentral}, for every hyperfixed automorphism $A:X\to X$ of an affine  space, the set $\Fix(A)$ is flat (more precisely, $\Fix(A)=X$ or $\Fix(A)$ is a hyperplane in $X$).  By Proposition~\ref{p:aff-paracentral<=>}, an automorphism $A:X\to X$ of an affine space $X$ is paracentral if and only if $A$ is a translation, a hypershear, or a hyperscale. In Sections~\ref{s:tring-trans}, \ref{s:tring-shears},   \ref{s:tring-scales}, \ref{s:tring-homos}, we shall study translations, hypershears, hyperscales, and homotheties of based affine planes. In Section~\ref{s:tring-homos} we also prove that the scalar corps of a Thalesian affine plane can be canonically identified with the kernel of the ternar of any affine base in the plane. All ternars appearing in this chapter are linear, so their structure is uniquely determined by their dot and plus operation. Therefore, there is no need to consider the puls operation (which is equal to the plus operation in linear ternars). 

\section{Some algebraic properties of biloops}

In this section we define some algebraic properties of biloops that will appear in the subsequent algebraic characterizations of various geometric properties of (based) affine planes. For a biloop $R$ we denote by $R^*$ the set $R\setminus\{0\}$ (which is a loop with respect to the dot operation inherited from the biloop $R$).

\begin{definition}\label{d:biloops-properties} A biloop $R$ is called
\begin{itemize}
\item \index{commutative-plus biloop}\index{biloop!commutative-plus}\defterm{commutative-plus} if $\forall x,y\in R\;\;x+y=y+x$;
\item \index{commutative-dot biloop}\index{biloop!commutative-dot}\defterm{commutative-dot} if $\forall x,y\in R\;\;x\cdot y=y\cdot x$;
\item \index{commutative biloop}\index{biloop!commutative}\defterm{commutative} if $R$ is commutative-plus and commutative-dot;
\smallskip
\item \index{associative-plus biloop}\index{biloop!associative-plus}\defterm{associative-plus} if $\forall x,y,z\in R\;\;x+(y+z)=(x+y)+z$;
\item \index{associative-dot biloop}\index{biloop!associative-dot}\defterm{associative-dot} if $\forall x,y,z\in R\;\;x\cdot(y\cdot z)=(x\cdot y)\cdot z$;
\item \index{associative biloop}\index{biloop!associative}\defterm{associative} if $R$ is associative-plus and associative-dot;
\smallskip
\item \index{Moufang-plus biloop}\index{biloop!Moufang-plus}\defterm{Moufang-plus} if $\forall x,y,z\in R\;\;(x+y)+(z+x)=(x+(y+z))+x=x+((y+z)+x)$;
\item \index{Moufang-dot biloop}\index{biloop!Moufang-dot}\defterm{Moufang-dot} if $\forall x,y,z\in R\;\;(x\cdot y)\cdot(z\cdot x)=(x\cdot (y\cdot z))\cdot x=x\cdot((y\cdot z)\cdot x)$;
\item \index{Moufang biloop}\index{biloop!Moufang}\defterm{Moufang} if $R$ is Moufang-plus and Moufang-dot;
\smallskip
\item \index{alternative-plus biloop}\index{biloop!alternative-plus}\defterm{alternative-plus} if $\forall x,z\in R\;\;\forall y\in\{x,z\}\;x+(y+ z)=(x+y)+z$;
\item \index{alternative-dot biloop}\index{biloop!alternative-dot}\defterm{alternative-dot} if $\forall x,z\in R\;\;\forall y\in\{x,z\}\;x\cdot(y\cdot z)=(x\cdot y)\cdot z$;
\item \index{alternative biloop}\index{biloop!alternative}\defterm{alternative} if $R$ is alternative-plus and alternative-dot;
\smallskip
\item \index{inversive-plus biloop}\index{biloop!inversive-plus}\defterm{inversive-plus} if $\forall x\in R\;\exists y\in R\;\forall z\in R\;\;y+(x+z)=z=(z+x)+y$;
\item \index{inversive-dot biloop}\index{biloop!inversive-dot}\defterm{inversive-dot} if $\forall x\in R^*\;\exists y\in R^*\;\forall z\in R\;\;y\cdot(x\cdot z)=z=(z\cdot x)\cdot y$;
\item \index{inversive biloop}\index{biloop!inversive}\defterm{inversive} if $R$ is inversive-plus and inversive-dot;
\smallskip
\item \index{invertible-plus biloop}\index{biloop!invertible-plus}\defterm{invertible-plus} if $\forall x\in R\;\exists y\in R\;x+y=0=y+x$;
\item \index{invertible-dot biloop}\index{biloop!invertible-dot}\defterm{invertible-dot} if $\forall x\in R^*\;\exists y\in R^*\;x\cdot y=1=y\cdot x$;
\item \index{invertible biloop}\index{biloop!invertible}\defterm{invertible} if $R$ is invertible-plus and invertible-dot;
\smallskip
\item \index{left-distributive biloop}\index{biloop!left-distributive}\defterm{left-distributive} if $\forall a,x,y\in R\;\;a\cdot(x+y)=(a\cdot x)+(a\cdot y)$;
\item \index{right-distributive biloop}\index{biloop!right-distributive}\defterm{right-distributive} if $\forall x,y,b\in R\;\;(x+y)\cdot b=(x\cdot b)+(y\cdot b)$;
\item \index{distributive biloop}\index{biloop!distributive}\defterm{distributive} if $R$ is left-distributive and right-distributive;
\smallskip
\item a \index{corps}\index{biloop!corps}\defterm{corps} if $R$ is distributive and associative;
\item a \index{field}\index{biloop!field}\defterm{field} if $R$ is a commutative corps.
\end{itemize}
\end{definition}

The first  18 properties of biloops introduced in Definition~\ref{d:biloops-properties} are special instances of the following general definition.

\begin{definition} Let $\mathcal P$ be a property of loops. A biloop $(R,+,\cdot)$ is defined to have the property
\begin{itemize}
\item \defterm{$\mathcal P$-plus} if the loop $(R,+)$ has property $\mathcal P$;
\item \defterm{$\mathcal P$-dot} if the loop $(R^*,\cdot)$ has property $\mathcal P$;
\item \defterm{$\mathcal P$} if the loops $(R,+)$ and $(R^*,\cdot)$ have property $\mathcal P$.
\end{itemize}
\end{definition}

The following proposition generalizes Hankel's Theorem~\ref{t:Hankel} to distributive Moufang-plus biloops.

\begin{proposition}\label{p:Moufang-plus=>commutative-plus} Every distributive Moufang-plus biloop is commutative-plus.
\end{proposition}

\begin{proof} Let $R$ be a distributive Moufang-plus biloop. By Corollary~\ref{c:Moufang=>diassociative}, any two elements $x,y\in R$ generate an associative subloop in the Moufang loop $(R,+)$. By the left-distributivity of the biloop $R$, 
$$(1+1)\cdot(x+y)=(1+1)\cdot x+(1+1)\cdot y=(1\cdot x+1\cdot x)+(1\cdot y+1\cdot y)=(x+x)+(y+y).$$
On the other hand, the right-distributivity of the biloop $R$ ensures that
$$(1+1)\cdot(x+y)=1\cdot(x+y)+1\cdot(x+y)=(1\cdot x+1\cdot y)+(1\cdot x+1\cdot y)=(x+y)+(x+y).$$
The associativity of the subloop generated by the elements $x,y$ implies
$$(x+(x+y))+y=(x+x)+(y+y)=(1+1)\cdot(x+y)=(x+y)+(x+y)=(x+(y+x))+y.$$By the cancellativity of the loop $(R,+)$, the latter equality implies $x+y=y+x$. Therefore, the biloop $R$ is commutative-plus.
\end{proof}

Proposition~\ref{p:Moufang-plus=>commutative-plus} motivates the following (open) problem.

\begin{problem} Is every distributive Moufang(-plus) biloop associative-plus?
\end{problem}

\begin{proposition}\label{p:Moufang-plus-not-associative-plus} If a distributive Moufang-plus biloop $R$ is not associative-plus, then $(R,+)$ is an infinite commutative Moufang loop of exponent $3$, with trivial centre.
\end{proposition}

\begin{proof} By Proposition~\ref{p:Moufang-plus=>commutative-plus}, the Moufang loop $(R,+)$ is commutative. Assuming that $(R,+)$ has non-trivial centre, we can find a nonzero element $c\in R$ such that $(x+y)+c=x+(y+c)$ for all $x,y\in R$. We claim that $(x+y)+z=x+(y+z)$ for all elements $x,y,z\in R$. If $z=0$, then $$(x+y)+z=(x+y)+0=x+y=x+(y+0)=x+(y+z).$$So, assume that $z\ne 0$. Since $(X,\cdot)$ is a $0$-loop, there exists a nonzero element $a\in R$ such that $z=a\cdot c$, and also there exist elements $x',y'\in X$ such that $x=a\cdot x'$ and $y=a\cdot y'$. By the right-distributivity of the biloop $R$,
$$
\begin{aligned}
(x+y)+z&=(a\cdot x'+a\cdot y')+a\cdot c=a\cdot (x'+y')+a\cdot c=a\cdot((x'+y')+c)\\
&=a\cdot (x'+(y'+c))=a\cdot x'+(a\cdot x'+a\cdot y')=x+(y+z).
\end{aligned}
$$
Therefore, the loop $(R,+)$ is associative, which contradicts our assumption. This contradiction shows that the centre of the Moufang loop $(R,+)$ is trivial.
Then the Bruck--Slaby Theorem~\ref{t:Bruck-Slaby} ensures that the commutative Moufang loop $(R,+)$ is infinite.

By Theorem~\ref{t:cubic-Moufang}, for every element $a\in R$, the element $a+a+a$ belongs to the (trivial) centre of $R$ and hence $a+a+a=0$, which means that the Moufang loop $(R,+)$ has exponent 3.
\end{proof}

The following corollary of Proposition~\ref{p:Moufang-plus-not-associative-plus} generalizes the Artin--Zorn Theorem~\ref{t:Artin-Zorn} to finite distributive Moufang biloops.

\begin{corollary} Every finite distributive Moufang biloop is a field.
\end{corollary}

\begin{proof} Let $F$ be a finite distributive Moufang biloop. By Propositions~\ref{p:Moufang-plus=>commutative-plus} and \ref{p:Moufang-plus-not-associative-plus}, $F$ is associative-plus and commutative-plus, so, $F$ is a divisible ring. Being Moufang-dot, the Moufang biloop $F$ is alternative-dot. By the Artin--Zorn Theorem~\ref{t:Artin-Zorn}, the finite alternative divisible ring $F$ is a field.
\end{proof}

\begin{remark} By Proposition~\ref{p:Moufang-plus=>commutative-plus},  distributive Moufang-plus biloops are commutative-plus. On the other hand, distributive Moufang (even associative) biloops (i.e., corps) need not be commutative-dot.
\end{remark}  

We recall that a biloop $R$ is \index{linear biloop}\index{biloop!linear}\defterm{linear} if $R$ is a ternar with respect to the ternary operation $$T:R^3\to R,\quad T:(x,a,b)\mapsto x_\times a_+b\defeq (x\cdot a)+b.$$ In this case, $R$ is the biloop of the ternar $(R,T)$. Indeed, $$x\cdot a\defeq x_\times a_+0=(x\cdot a)+0=x\cdot a\quad\mbox{and}\quad x+b\defeq x_\times 1_+b=(x\cdot 1)+b=x+b.$$
Linear biloops are characterized in Proposition~\ref{p:linear-biloop<=>}. For right-distributive associative-plus biloops this characterization can be simplified as follows.

\begin{proposition}\label{p:quasi-field} A right-distributive associative-plus biloop $R$ is linear if and only if for any elements $a,b,c\in R$ with $a\ne c$, the equation $x\cdot a=(x\cdot c)+b$ has a unique solution $x\in R$.
\end{proposition}

\begin{proof} The ``only if'' part follows immediately from Proposition~\ref{p:linear-biloop<=>}. To prove the ``if'' part, assume that for any elements $a,b,c\in R$ with $a\ne c$, the equation $(x\cdot a)+b=x\cdot c$ has a unique solution $x\in R$.
 To show that $R$ is a linear biloop, it suffices to check the conditions (1) and (2) of Proposition~\ref{p:linear-biloop<=>}. 

To check the condition (1), take any elements $a,b,c,d\in R$ with $a\ne c$ and consider the equality $(x\cdot a)+b=(x\cdot c)+d$. Since the biloop $R$ is associative-plus, the loop $(R,+)$ is a group. Then there exists a unique element $b'\in R$ such that $b+b'=0$. Then the equation $(x\cdot a)+b=(x\cdot c)+d$ is equivalent to the equation
$$x\cdot a=(x\cdot a)+0=(x\cdot a)+(b+b')=((x\cdot a)+b)+b'=((x\cdot c)+d)+b'=(x\cdot c)+(d+b'),$$
which has a unique solution by our assumption.

To check the condition (2), take any elements $\hat x,\hat y,\check x,\check y\in R$ with $\hat x\ne\check x$. We have to show that the equation $(\hat x\cdot a)+b=\hat y$ and $(\check x\cdot a)+b=\check y$ has a unique solution $a,b\in R$. Assume that elements $a,b\in R$ satisfy the equations $(\hat x\cdot a)+b=\hat y$ and $(\check x\cdot a)+b=\check y$.

Since  the additive loop $(R,+)$ is a group, there exist unique elements $\mbox{-}\hat x$ and $\mbox{-}\hat y$ in $R$ such that $(\mbox{-}\hat x)+\hat x=0=\hat x+(\mbox{-}\hat x)$ and $(\mbox{-}\hat y)+\hat y=0=\hat y+(\mbox{-}\hat y)$. Using the associativity of the addition operation and the right-distributivity of the quasi-field $R$, we conclude that 
$$(\mbox{-}\hat x)\cdot a+\hat y=(\mbox{-}\hat x)\cdot a+((\hat x\cdot a)+b)=((\mbox{-}\hat x)\cdot a+(\hat x\cdot a))+b=(((\mbox{-}\hat x)+\hat x)\cdot a)+b=(0\cdot a)+b=0+b=b$$ and hence
$$
\begin{aligned}
\check y+(\mbox{-}\hat y)&=((\check x\cdot a)+b)+(\mbox{-}\hat y)=((\check x\cdot a)+((\mbox{-}\hat x)\cdot a)+\hat y))+(\mbox{-}\hat y)\\
&=(((\check x\cdot a)+(\mbox{-}\hat x)\cdot a))+\hat y)+(\mbox{-}\hat y)=((\check x+(\mbox{-}\hat x))\cdot a)+(\hat y+(\mbox{-}\hat y))\\
&=(\check x+(\mbox{-}\hat x))\cdot a
\end{aligned}
$$
Assuming that $\check x+(\mbox{-}\hat x)=0$, we conclude that 
$$\hat x=0+\hat x=(\check x+(\mbox{-}\hat x))+\hat x=\check x+((\mbox{-}\hat x)+\hat x))=\check x+0=\check x,$$
which contradicts the choice of the points $\hat x\ne\check x$. This contradiction shows that $\check x+(\mbox{-}\hat x)\ne 0$. Since $(R,\cdot)$ is a $0$-loop, the equation $\check y+(\mbox{-}\hat y)=(\check x+(\mbox{-}\hat x))\cdot a$ has a unique solution $a\in R$. Since $(R,+)$ is a loop and $a$ is unique, the equation $(\hat x\cdot a)+b=\hat y$ has a unique solution $b\in R$, which is also a solution of the equation $(\check x\cdot a)+b=\check y$.
\end{proof}

\begin{proposition}\label{p:semifield=>quasifield} Every distributive associative-plus biloop $R$ is linear.
\end{proposition}

\begin{proof} By Proposition~\ref{p:quasi-field}, it suffices to show that for every elements $a,b,c\in R$ with $a\ne c$, the equation $x\cdot a=(x\cdot c)+b$ has a unique solution $x\in R$. Since $(R,+)$ is a group, this equation is equivalent to $-(x\cdot c)+(x\cdot a)=b$. The left-distributivity of $R$ implies $0\cdot(c-c)=x\cdot(c-c)=(x\cdot c)+x\cdot(-c)$ and hence $-(x\cdot c)=x\cdot(-c)$. Then $$b=-(x\cdot c)+(x\cdot a)=x\cdot(-c)+x\cdot a=x\cdot (-c+a).$$ Since $(R,\cdot)$ is a $0$-loop, the equation $x\cdot (-a+c)=b$ has a unique solution $x\in R$, which is also a unique solution of the equation $(x\cdot a)=(x\cdot c)+b$.
\end{proof}

\begin{proposition}\label{p:fin-near-field=>quasi-field} Every finite right-distributive associative-plus biloop $R$ is linear.
\end{proposition}

\begin{proof} By Proposition~\ref{p:quasi-field}, it suffices to check that for every elements $a,b,c\in R$ with $a\ne c$, there exists a unique solution $x\cdot a=(x\cdot c)+b$.  Since $(R,+)$ is a group, for every element $y\in R$ there exists a unique element $y'\in R$ such that $y+y'=0=y'+y$.
Since $R$ is right-distributive, $0=0\cdot c=(y+y')\cdot c=y\cdot c+y'\cdot c$ and hence $(y\cdot c)'=y'\cdot c$. We claim that the function $F:R\to R$, $F:x\mapsto x'\cdot c+x\cdot a$, is injective. Indeed, for any $x,y\in R$ with $F(x)=F(y)$, we have $x'\cdot c+x\cdot a=F(x)=F(y)=y'\cdot c+y\cdot a$. Since the loop $(R,+)$ is a group and $R$ is right-distributive, the equality $x'\cdot c+x\cdot a=y'\cdot c+y\cdot a$ implies 
$(y+x')\cdot c=(y+x')\cdot c$ and hence $y+x'=0$ and $y=x$, witnessing the injectivity of the function $F:R\to R$. Since $R$ is finite, the injective function $F:R\to R$ is bijective and hence there exists a unique element $x$ such that $x'\cdot c+x\cdot a=F(x)=b$. Adding to both sides of this equality $x\cdot c$ from the left, we obtain that $x$ is a unique element of $R$ satisfying the equation  $x\cdot a=(x\cdot c)+b$.
\end{proof}

\begin{definition}\label{d:quasi-semi-near-field} A biloop $R$ is called 
\begin{itemize}
\item a \index{cartesian-group}\index{biloop!cartesian-group}\defterm{cartesian-group} if $R$ is linear and associative-plus;
\item a \index{quasi-field}\index{biloop!quasi-field}\defterm{quasi-field} if $R$ is linear, right-distributive and associative-plus;
\item a \index{semi-field}\index{biloop!semi-field}\defterm{semi-field} if $R$ is a distributive and associative-plus;
\item a \index{near-field}\index{biloop!near-field}\defterm{near-field} if $R$ is right-distributive and associative.
\end{itemize}
\end{definition}

\begin{remark} In Theorem~\ref{t:VW-Thalesian<=>quasifield} we shall prove that a based affine plane  is translation if and only if its ternar is linear and its biloop is a quasi-field.
\end{remark}

By Propositions~\ref{p:semifield=>quasifield}, \ref{p:fin-near-field=>quasi-field}, for every biloop the following implications hold:
$$
\xymatrix{
\mbox{field}\ar@{=>}[r]\ar@{=>}[d]&\mbox{semi-field}\ar@{=>}[r]&\mbox{quasi-field}\ar@{=>}[r]&\mbox{cartesian}\atop\mbox{group}\ar@{=>}[d]\\
\mbox{corps}\ar@{=>}[r]&\mbox{linear}\atop\mbox{near-field}\ar@{=>}[ru]&\mbox{finite}\atop\mbox{near-field}\ar@{=>}[l]&\mbox{linear}\atop\mbox{biloop}
}
$$

\begin{Exercise} Construct an example of a non-linear near-field.
\end{Exercise}

\begin{definition} Let $\mathcal P$ be a property of biloops. A ternar $R$ is defined to have the property $\mathcal P$ if its biloop has the property $\mathcal P$.
\end{definition}





\section{Linear ternars}

In this section we find a geometric characterization of linear ternars.


Let us recall that a ternar $R$ is \index{linear ternar}\index{ternar!linear}\defterm{linear} if $$x_\times a_+b=(x\cdot a)+b\defeq (x_\times a_+0)_\times 1_+b$$ for every elements $x,a,b\in R$. In this case the ternary operation of the ternar is uniquely determined by the multiplication and addition.

\begin{definition} A based affine plane $(\Pi,uow)$ is defined to be \index{linear based affine plane}\index{based affine plane!linear}\defterm{linear} if for any distinct parallel lines $A,B,C\subseteq \Pi$ with $C=\Aline ow$ and any distinct points $a,a'\in A$, $b,b'\in B$, $c,c'\in C$ with $\Aline ab\parallel \Aline ou\parallel \Aline{a'}{b'}$ and $\Aline bc\parallel \Aline oe\parallel {b'}{c'}$, we have $\Aline ac\parallel \Aline{a'}{c'}$.
\end{definition}

Therefore, linear based affine planes satisfy a weak form of the Thales Axiom.

\begin{theorem}\label{t:tring-linear<=>}  A based affine plane is linear if and only if its ternar is linear.
\end{theorem}

\begin{proof} Let $(\Pi,uow)$ be a based affine plane and $\Delta$ be its ternar. 

Assume that the based affine plane $(\Pi,uow)$ is linear. To prove that the ternar $\Delta$ is linear, we need to check that $x_\times p_+q=(x{\cdot} p){+}q$ for every points $x,p,q\in \Delta$. The equality $x_\times p_+q=(x{\cdot} p){+}q$ trivially holds if $p=e$ or $q=o$ or $x=o$. So, assume that $p\ne e$ and $q\ne o\ne x$.  Consider the unique point $p'\in\Aline ev\cap\Aline ph$. Next, find unique points $a\in \Aline o{p'}\cap \Aline xv$ and $b\in \Delta\cap\Aline ah$. Consider the points $c\defeq o$ and $c'\in\Aline ov\cap \Aline qh$. Find a unique point $b'\in \Aline bv$ such that $\Aline{c'}{b'}\parallel \Aline cb=\Delta$, and unique points $a'\in\Aline av\cap \Aline{b'}h$, and $s\in\Delta\cap\Aline {a'}{b'}$. The definitions of the multiplicaton and addition in the ternar $\Delta$ ensure that $x\times p=b$ and $b+q=s$. The linearity of the based affine plane $(\Pi,uow)$ implies $\Aline ac\parallel \Aline{a'}{c'}$ and hence $x_\times p_+b=s=b{+}q=(x{\times} p){+}q$. 

\begin{picture}(180,185)(-130,-15)
\linethickness{0.8pt}
\put(0,0){\color{teal}\vector(1,0){150}}
\put(0,0){\color{red}\line(1,1){150}}
\put(0,0){\color{blue}\line(1,2){45}}
\put(0,0){\color{cyan}\vector(0,1){150}}
\put(0,60){\color{blue}\line(1,2){45}}
\put(0,60){\color{red}\line(1,1){90}}
\put(15,15){\color{cyan}\line(0,1){15}}
\put(30,30){\color{teal}\line(-1,0){15}}
\put(45,45){\color{cyan}\line(0,1){105}}
\put(60,60){\color{teal}\line(-1,0){60}}
\put(90,90){\color{cyan}\line(0,1){60}}
\put(90,90){\color{teal}\line(-1,0){45}}
\put(150,150){\color{teal}\line(-1,0){105}}

\put(0,0){\circle*{3}}
\put(1,-8){$o$}
\put(-7,-8){$c$}
\put(153,-3){$\color{teal}\boldsymbol h$}
\put(-2,153){$\color{cyan}\boldsymbol v$}
\put(15,15){\circle*{3}}
\put(15,7){$e$}
\put(30,30){\circle*{3}}
\put(32,24){$p$}
\put(45,45){\circle*{3}}
\put(47,38){$x$}
\put(60,60){\circle*{3}}
\put(62,53){$q$}
\put(90,90){\circle*{3}}
\put(92,84){$b$}
\put(150,150){\circle*{3}}
\put(153,148){$s$}

\put(0,60){\circle*{3}}
\put(-8,58){$c'$}
\put(15,30){\circle*{3}}
\put(8,30){$p'$}
\put(45,90){\circle*{3}}
\put(47,93){$a$}
\put(45,150){\circle*{3}}
\put(42,153){$a'$}
\put(90,150){\circle*{3}}
\put(88,153){$b'$}
\end{picture}

Now assume that the ternar $\Delta$ is linear. To prove that the based affine plane $(\Pi,uow)$ is linear, take any distinct parallel lines $A,B,C$ with $C=\Aline ov$, and distinct points $a,a'\in A$, $b,b'\in B$, $c,c'\in C$ such that $\Aline ab\parallel \Aline ou\parallel\Aline{a'}{b'}$ and $\Aline bc\parallel\Delta\parallel \Aline{b'}{c'}$. Consider the unique points $x\in A\cap\Delta$ and $y\in B\cap\Delta$. It follows from $A\ne C=\Aline ov\ne B$ that the points $o,x,y$ are distinct. Find a unique point $p'\in\Aline xv\cap\Aline yh$. Find a unique point $p\in\Delta$ such that $y=x\cdot p$. Find unique points $\gamma\in \Delta\cap\Aline ch$, $\gamma'\in\Delta\Aline{c'}h$, $z\in\Delta\cap\Aline ab$ and $z'\in\Delta\cap \Aline {a'}{b'}$. The definition of the operations of addition and multiplication, and the linearity of the ternar $\Delta$ imply $$z=y+\gamma=(x{\cdot}p)+\gamma=x_\times p_+\gamma\quad\mbox{and}\quad z'=y+\gamma'=(x{\cdot}p)+\gamma'=x_\times p_+{\gamma'}.$$ Then the lines $\Aline ca=L_{p,\gamma}$ and $\Aline{c'}{a'}=L_{p,\gamma'}$ are parallel, witnessing that the based affine plane $(\Pi,uow)$ is linear.

\begin{picture}(150,170)(-130,-15)
\linethickness{0.8pt}
\put(0,0){\color{red}\line(1,1){135}}
\put(0,0){\color{cyan}\line(0,1){60}}
\put(15,15){\color{teal}\line(-1,0){15}}
\put(45,45){\color{cyan}\line(0,1){90}}
\put(75,75){\color{cyan}\line(0,1){60}}
\put(0,15){\color{blue}\line(3,5){45}}
\put(0,15){\color{red}\line(1,1){75}}
\put(0,60){\color{red}\line(1,1){75}}
\put(0,60){\color{blue}\line(3,5){45}}
\put(60,60){\color{teal}\line(-1,0){60}}
\put(90,90){\color{teal}\line(-1,0){45}}
\put(135,135){\color{teal}\line(-1,0){90}}
\put(75,75){\color{teal}\line(-1,0){30}}
\put(0,0){\color{blue}\line(3,5){45}}

\put(45,75){\circle*{3}}
\put(47,78){$p'$}
\put(0,0){\circle*{3}}
\put(0,-7){$o$}
\put(15,15){\circle*{3}}
\put(16,6){$\gamma$}
\put(45,45){\circle*{3}}
\put(45,37){$x$}
\put(60,60){\circle*{3}}
\put(60,53){$\gamma'$}
\put(75,75){\circle*{3}}
\put(75,68){$y$}
\put(90,90){\circle*{3}}
\put(90,83){$z$}
\put(0,15){\circle*{3}}
\put(-8,13){$c$}
\put(135,135){\circle*{3}}
\put(135,138){$z'$}
\put(75,90){\circle*{3}}
\put(76,92){$b$}
\put(45,90){\circle*{3}}
\put(47,92){$a$}
\put(45,135){\circle*{3}}
\put(42,138){$a'$}
\put(75,135){\circle*{3}}
\put(72,138){$b'$}
\put(0,60){\circle*{3}}
\put(-10,58){$c'$}
\end{picture}
\end{proof}

Theorem~\ref{t:tring-linear<=>} implies the following characterization of Thalesian affine planes.

\begin{theorem}\label{t:Thalesian<=>linear} An affine plane $\Pi$ is Thalesian if and only if for every affine base $ouw$ in $\Pi$ the based affine plane $(\Pi,uow)$ is linear.
\end{theorem} 

Another characterization of Thalesian planes involves the plus=puls property of its ternars. We recall that a ternar $R$ is called a \defterm{ternar of an affine space} $X$ if it is isomorphic to the ternar $\Delta_{uow}$ of some triangle $uow$ (which is an affine base for the plane $\overline{\{u,o,w\}}$ in $X$).  A ternar  $(X,T)$ is \defterm{plus\textup{=}puls} if $x+y=T(x,1,y)=T(1,x,y)=x\puls y$ for all $x,y\in R$.

\begin{theorem}\label{t:plus=puls<=>} An affine space $X$ is Thalesian if and only if every ternar of $X$ is plus\textup{=}puls.
\end{theorem}

\begin{proof}  The ``only if'' part follows from Theorem~\ref{t:Thalesian<=>linear} and  Proposition~\ref{p:ternar-linear<=>}.

To prove the ``if'' part, assume that every ternar of an affine space $X$ is plus=puls. If $\|X\|\ne 3$, then the affine space $X$ is Desarguesian and hence Thalesian, by Corollaries~\ref{c:affine-Desarguesian} and Theorem \ref{t:ADA=>AMA}. So, assume that $\|X\|=3$, which means that $X$ is a Playfair plane. To check that $X$ is Thalesian, take any distinct parallel lines $A,B,C\subset X$ and any points $a,a'\in A$, $b,b'\in B$, $c,c'\in C$ such that $\Aline ab\cap\Aline{a'}{b'}=\varnothing=\Aline bc\cap\Aline{b'}{c'}$. We need to prove that $\Aline ac\cap\Aline {a'}{c'}=\varnothing$. Corollary~\ref{c:parallel-lines<=>} ensures that $\Aline ab\parallel \Aline {a'}{b'}$ and $\Aline bc\parallel \Aline{b'}{c'}$.

If $b\in \Aline ac$, then $\Aline {a'}{b'}\parallel \Aline ab=\Aline bc\parallel \Aline{b'}{c'}$ implies $\Aline {a'}{b'}=\Aline{b'}{c'}=\Aline{a'}{c'}$ and hence $\Aline ac\cap\Aline {a'}{c'}=\Aline ab\cap\Aline{a'}{b'}=\varnothing$. By analogy we can prove that $b'\in\Aline{a'}{c'}$ implies $\Aline ac\cap \Aline{a'}{c'}=\varnothing$. 

So, assume that $abc$ and $a'b'c'$ are triangles. Since $A,B$ are two parallel lines in the Playfair (and hence Proclus) plane $X$ and $\Aline ac\cap A=\{a\}$, there exists a unique point $e\in B\cap \Aline ac$. Since $X$ is Playfair, there exist unique points $w\in A$, $u\in B$, and $y\in\Aline ac$ such that the lines $\Aline ou$, $\Aline we$ and $\Aline {a'}y$ are parallel to the line $\Aline bc$.

\begin{picture}(180,200)(-130,-30)
\linethickness{0.8pt}
\put(0,0){\color{teal}\line(1,0){150}}
\put(0,0){\color{red}\line(1,1){150}}
\put(0,0){\color{blue}\line(1,2){45}}
\put(0,-15){\color{cyan}\line(0,1){165}}
\put(0,60){\color{blue}\line(1,2){45}}
\put(0,60){\color{red}\line(1,1){90}}
\put(0,45){\color{teal}\line(1,0){45}}
\put(45,-15){\color{cyan}\line(0,1){165}}
\put(60,60){\color{teal}\line(-1,0){60}}
\put(90,-15){\color{cyan}\line(0,1){165}}
\put(90,90){\color{teal}\line(-1,0){45}}
\put(150,150){\color{teal}\line(-1,0){105}}

\put(0,0){\circle*{3}}
\put(1,-8){$o$}
\put(-6,-8){$a$}
\put(0,45){\circle*{3}}
\put(-11,42){$w$}
\put(-4,-25){\color{cyan}$A$}
\put(41,-25){\color{cyan}$B$}
\put(87,-25){\color{cyan}$C$}
\put(45,0){\circle*{3}}
\put(46,-8){$u$}
\put(45,45){\circle*{3}}
\put(47,38){$e$}
\put(60,60){\circle*{3}}
\put(62,53){$y$}
\put(90,90){\circle*{3}}
\put(92,84){$c$}
\put(150,150){\circle*{3}}
\put(153,148){$s$}

\put(0,60){\circle*{3}}
\put(-9,58){$a'$}
\put(45,90){\circle*{3}}
\put(47,93){$b$}
\put(45,150){\circle*{3}}
\put(42,153){$b'$}
\put(90,150){\circle*{3}}
\put(88,153){$c'$}
\end{picture}

Then $uow$ is an affine base in the affine plane $X$ whose diagonal $\Delta$ coincides with the line $\Aline ac$. By our assumption, the ternar $(\Delta,T_{uow})$ is plus=puls and hence $c+y=c\puls y=s$ for some point $s\in \Delta=\Aline ac$. Since $\Aline {a'}{b'}\parallel \Aline ab$ and $\Aline {a'}y\parallel \Aline ou$, the equality $s=c\puls y=T_{uow}(1,c,y)$ ensures that $s\in \Aline {b'}{c'}$. On the other hand, the equality $s=c+y=T(x,1,y)$ ensures that $\Aline{a'}{c'}\parallel \Delta=\Aline ac$. Therefore, the lines $\Aline ac$ and $\Aline{a'}{c'}$ are parallel. Assuming that these parallel lines have common point, we conclude that they coincide and hence $\{a\}=\Aline ac\cap A=\Aline{a'}{c'}\cap A=\{a'\}$, which contradicts $\Aline ab\cap\Aline{a'}{b'}=\varnothing$. This contradiction completes the proof of the equality $\Aline ac\cap\Aline{a'}{c'}=\varnothing$.
\end{proof}

\section{Translations of based affine planes}\label{s:tring-trans} 

In this section we study the interplay between horizontal and vertical translations of based affine planes and algebraic properties of their ternars.

\begin{definition}\label{d:trans-hv} Let $(\Pi,uow)$ be a based affine plane. A translation $T:\Pi\to\Pi$ is called 
\begin{itemize}
\item \index{horizontal translation}\index{translation!horizontal}\defterm{horizontal} if $T(o)\in \Aline ou$;
\item \index{vertical translation}\index{translation!vertical}\defterm{vertical} if $T(o)\in \Aline ow$;
\item \index{diagonal translation}\index{translation!diagonal}\defterm{diagonal} if $T(o)\in \Aline oe$.
\end{itemize}
\end{definition}

Based planes admitting all possible vertical (resp. horizontal, diagonal) translations are called vertical-translation (resp. horizontal-translation, diagonal-translation). A precise definition follows.

\begin{definition} A based affine plane $(\Pi,uow)$ is defined to be 
\begin{enumerate}
\item \index{based affine plane!horizontal-translation}\defterm{horizontal-translation} if $\forall x\in \Aline ou$ there exists a translation $T:\Pi\to \Pi$ such that $T(o)=x$;
\item \index{based affine plane!vertical-translation}\defterm{vertical-translation} if $\forall y\in \Aline ow$ there exists a translation $T:\Pi\to \Pi$ such that $T(o)=y$;
\item \index{based affine plane!diagonal-translation}\defterm{diagonal-translation} if $\forall z\in \Aline oe$ there exists a translation $T:\Pi\to \Pi$ such that $T(o)=z$.
\end{enumerate}
\end{definition}

\begin{exercise} Show that a based affine plane $(\Pi,uow)$ is horizontal-translation (resp. vertical-translation, diagonal-translation) if and only if for every point $x$ in the line $\Aline ou$ (resp. $\Aline ow$, $\Aline oe$), the vector $\overvector{ox}$ is functional.
\end{exercise}

\begin{proposition}\label{p:trans<=>hv-trans} For a based affine plane $(\Pi,uow)$, the following conditions are equivalent:
\begin{enumerate}
\item the affine plane $\Pi$ is a translation plane;
\item the based affine plane $(\Pi,uow)$ is horizontal-translation and vertical-translation;
\item the based affine plane $(\Pi,uow)$ is horizontal-translation and diagonal-translation;
\item the based affine plane $(\Pi,uow)$ is vertical-translation and diagonal-translation.
\end{enumerate}
\end{proposition}

\begin{proof} It is clear that $(1)$ implies the conditions $(2),(3),(4)$.
\smallskip

$(2)\Ra(1)$. Assume that a based affine plane $(\Pi,uow)$ is horizontal-translation and vertical-translation. To see that $\Pi$ is a translation plane, it suffices to show that for every $z\in \Pi$ there exists a translation $T:\Pi\to\Pi$ such that $T(o)=z$. If $z\in\Aline ou\cup\Aline ow$, then the translation $T$ with $T(o)=z$ exists because $(\Pi,uow)$ is horizontal-translation and vertical-translation. So, assume that $z\in \Pi\setminus(\Aline ou\cup\Aline ow)$. Since $\Pi$ is a Playfair plane, 
there exist unique points $x\in \Aline ou\cap\Aline zv$ and $y\in\Aline ow\cap\Aline zh$. Since the based affine plane $(\Pi,uow)$ is horizontal-translation and vertical-translation, there exist translations $T_x,T_y$ of the plane $\Pi$ such that $T_x(o)=x$ and $T_y(o)=y$. 

By Proposition~\ref{p:Trans(X)isnormal}, the composition $T\defeq T_y\circ T_x$ is a translation of the affine plane $\Pi$. We claim that the point $z'\defeq T(o)=T_y(T_x(o))=T_y(x)$ coincides with the point $z$. Since $T_y(o)=y\in\Aline ow\setminus\{o\}$, Proposition~\ref{p:Trans-spread} implies $\Aline x{z'}\parallel \Aline oy=\Aline ow\parallel \Aline xz$ and hence $\Aline x{z'}=\Aline xz$. Since $T_y$ is a dilation, $\Aline yz\parallel \Aline ou\parallel T_y[\Aline ou]=T_y[\Aline ox]=\Aline y{z'}$ and hence $\Aline yz=\Aline y{z'}$. Then $z'\in \Aline x{z'}\cap\Aline y{z'}=\Aline xz\cap\Aline yz=\{z\}$ and hence $T(o)=z'=z$.
\smallskip

The implications $(3)\Ra(1)$ and $(4)\Ra(1)$ can be proved by analogy with the implication $(2)\Ra(1)$.
\end{proof}

Applying Definition~\ref{d:trans-hv} to the coordinate plane $R^2$ of a ternar $R$, we conclude that a translation $T:R^2\to R^2$ is  
\begin{itemize}
\item \index{translation!vertical}\index{vertical translation}\defterm{vertical} iff $T(0,0)\in \{0\}\times R$;
\item \index{translation!horizontal}\index{horizontal translation}\defterm{horizontal} iff $T(0,0)\in R\times\{0\}$;
\item \index{translation!diagonal}\index{diagonal translation}\defterm{diagonal} iff $T(0,0)\in \Delta\defeq\{(x,y)\in R^2:x=y\}$.
\end{itemize}

\begin{proposition}\label{p:tring-vshifts} Let $R$ be a ternar and $R^2$ be its coordinate plane. For every element $s\in R$, the following conditions are equivalent:
\begin{enumerate}
\item there exists a vertical translation $S:R^2\to R^2$ such that $S(0,0)=(0,s)$;
\item the function $S:R^2\to R^2$, $S:(x,y)\mapsto (x,y+s)$, is a translation of the affine plane $R^2$;
\item $\forall x,a,b\in R\;\;x_\times a_+(b+s)=(x_\times a_+b)+s$.
\end{enumerate}
The equivalent conditions \textup{(1)--(3)} imply the condition
\begin{enumerate}
\item[(4)] $\forall x\in R\;\;x\puls s=x+ s$.
\end{enumerate}
\end{proposition}

\begin{proof} $(1)\Ra(2)$ Assume that $S:R^2\to R^2$ is a translation of the plane $R^2$ such that $S(0,0)=(0,s)$. Then $S[L_0]=L_0$ and $S[L_c]=L_c$ for every  $c\in R$, by Proposition~\ref{p:Trans-spread}. Since $S$ is a translation, $S[L_{1,0}]\parallel L_{1,0}$ and hence $S[L_{1,0}]=L_{1,s}$. Then for every elements $x,y\in R$ we have that $S(y,y)\in S[L_{1,0}\cap L_y]=S[L_{1,0}]\cap S[L_y]=L_{1,s}\cap L_y$, which implies $S(y,y)=(y,y_\times 1_+s)=(y,y+s)$. Then $(y,y+s)=S(y,y)\in S[ L_{0,y}]=L_{0,y+s}$ and finally $S(x,y)\in S[L_x\cap L_{0,y}]=S[L_x]\cap S[L_{0,y}]=L_x\cap L_{0,y+s}=\{(x,y+s)\}$. 
\smallskip

$(2)\Ra(3)$ Assume that the function $S:R^2\to R^2$, $S:(x,y)\mapsto (x,y+s)$, is a translation of the affine plane $R^2$. Given any elements $x,a,b\in R$, we need to check that $x_\times a_+(b+s)=(x_\times a_+b)+s$. Since $S^{-1}$ is a dilation, $S[L_{a,b}]\parallel L_{a,0}$ and hence $S[L_{a,b}]=L_{a,c}$ for some $c\in R$. It follows from $(0,b+s)=S[L_{a,b}]=L_{a,c}$ that $c=b+s$. Then $(x,(x_\times a_+b)+s)=S(x,x_\times a_+b)\in S[L_{a,b}]=L_{a,c}$ and hence $$(x_\times a_+b)+s=x_\times a_+(b+s).$$
\smallskip 

$(3)\Ra(1)$ Assume that the condition (3) holds. We claim that the function $S:R^2\to R^2$, $S:(x,y)\mapsto (x,y+s)$ is a translation. Given any elements $a,b\in R$, it suffices to show that the set $S[L_{a,b}]$ is a line, parallel to the line $L_{a,b}$. We claim that $S[L_{a,b}]=L_{a,c}$ where $c\defeq b+s$. Indeed, for every point $(x,y)\in L_{a,b}$, we have 
$y=x_\times a_+b$ and hence $S(x,y)=(x,y+s)=(x,(x_\times a_+b)+s)=(x,x_\times a_+(b+s))\in L_{a,c}$, witnessing that $S[L_{a,b}]\subseteq L_{a,c}$. Since both sets $S[L_{a,b}]$ and $L_{a,c}$ project bijectively onto the set $R\times\{0\}$, the inclusion $S[L_{a,b}]\subseteq L_{a,c}$ implies the equality  $S[L_{a,b}]=L_{a,c}$, witnessing that  the function $S:R^2\to R^2$ is a translation of the affine plane $R^2$. It is clear that $S(0,0)=(0,0+s)=(0,s)$.
\smallskip

$(3)\Ra(4)$ For every $x\in R$, the condition (3) implies
$$x\puls s=1_\times x_+s=1_\times x_+(o+s)=(1_\times x_+o)+s=x+s.$$
\end{proof}
\smallskip

\begin{proposition}\label{p:tring-hshifts} Let $R$ be a ternar and $R^2$ be its coordinate plane. For every element $s\in R$, the following conditions are equivalent:
\begin{enumerate}
\item there exists a horizontal translation $S:R^2\to R^2$ such that $S(0,0)=(s,0)$;
\item the function $S:R^2\to R^2$, $S:(x,y)\mapsto (x+s,y)$, is a translation of the affine plane $R^2$;
\item $\forall x,a,b\in R\;\;(x+s)_\times a_+b=x_\times a_+(s_\times a_+b)$.
\end{enumerate}
\end{proposition}

\begin{proof} $(1)\Ra(2)$ Assume that there exists a translation $S:R^2\to R^2$ of the plane $R^2$ such that $S(0,0)=(s,0)$. We will prove that $S(x,y)=(x+s,y)$ for every $(x,y)\in R^2$. It follows from $S(0,0)=(s,0)$ that $S[L_{0,0}]=L_{0,0}$. Proposition~\ref{p:Trans-spread} implies that $S[L_{0,b}]=L_{0,b}$ for every $b\in R$. Since $S$ is a translation, $S[L_{1,0}]\parallel L_{1,0}$ and hence $(s,0)=S(0,0)\in S[L_{1,0}]=L_{1,b}$ for some $b\in R$ such that $(s,0)\in L_{1,b}$, which implies $s+b=s_\times 1_+b=0$ and $b=(\mbox{-}s)$, where $(\mbox{-}s)$ is the unique element of the loop $(R,+)$ such that $s+(\mbox{-}s)=0$. Then for every element $x\in R$ we have 
 $S(x,x)\in S[L_{1,0}\cap L_{0,x}]=S[L_{1,0}]\cap S[L_{0,x}]=L_{1,\mbox{-}s}\cap L_{0,x}$, which implies $S(x,x)=(z,x)$ for a unique element $z\in R$ such that $z+(\mbox{-}s)=x$ and hence $z=x{-}(\mbox{-}s)$. Here $x{-}(\mbox{-}s)$ is the unique element of the loop $(R,+)$ such that $(x{-}(\mbox{-}s))+(\mbox{-}s)=x$. Therefore, $S(x,x)=(z,x)=(x{-}(\mbox{-}s),x)$ and $S(x,y)\in S[L_x\cap L_{0,y}]=S[L_x]\cap S[L_{0,y}]=L_{x{-}(\mbox{-}s)}\cap L_{0,y}=\{(x{-}(\mbox{-}s),y)\}$. 

Since $S^{-1}$ is a translation, $S^{-1}[L_{1,0}]\parallel L_{1,0}$ and hence $S^{-1}[L_{1,0}]=L_{1,c}$ for some $c\in R$. Observe that $s+(\mbox{-}s)=0=(0{-}(\mbox{-}s))+(\mbox{-}s)$ and hence $s=0{-}(\mbox{-}s)$. Then $$(0,s)=S^{-1}(0{-}(\mbox{-}s),s_\times 1_+0)=S^{-1}(s,s_\times 1_+0)\in S^{-1}[L_{1,0}]=L_{1,c}$$ and hence $s=0_\times 1_+c=c$. Then $(x,x{-}(\mbox{-}s))=S^{-1}(x{-}(\mbox{-}s),(x{-}(\mbox{-}s))_\times 1_+0)\in S^{-1}[L_{1,0}]=L_{1,c}$ and hence $x{-}(\mbox{-}s)=x_\times 1_+c=x+s$ and 
$$S(x,y)=(x{-}(\mbox{-}s),y)=(x+s,y)$$for every $x,y\in R$.
\smallskip

2. Assume that the function $S:R^2\to R^2$, $S:(x,y)\mapsto (x+s,y)$, is a translation of the affine plane $R^2$. Given any elements $x,a,b\in R$, we need to check that $(x+s)_\times a_+b=x_\times a_+(s_\times a_+b)$. Since $S^{-1}$ is a dilation, $S^{-1}[L_{a,b}]\parallel L_{a,0}$ and hence $S^{-1}[L_{a,b}]=L_{a,c}$ for some $c\in R$. Then $$(0,s_\times a_+b)=S^{-1}(s,s_\times a_+b)\in S^{-1}[L_{a,b}]=L_{a,c}$$ and hence $s_\times a_+b=0_\times a_+c=c$. Then $(x,(x+s)_\times a_+b)=S^{-1}(x+s,(x+s)_\times a_+b)\in S^{-1}[L_{a,b}]=L_{a,c}$ and hence $$(x+s)_\times a_+b=x_\times a_+c=x_\times a_+(s_\times a_+b).$$

3. Assume that the condition (3) is satisfied. We claim that the  function $S:R^2\to R^2$, $S:(x,y)\mapsto (x+s,y)$ is a translation. Given any elements $a,b\in R$, it suffices to show that the set $S^{-1}[L_{a,b}]$ is a line, parallel to the line $L_{a,b}$. We claim that $S^{-1}[L_{a,b}]=L_{a,c}$ where $c\defeq s_\times a_+b$. Indeed, for every point $(x,y)\in L_{a,c}$, we have 
$$y=x_\times a_+c=x_\times a_+(s_\times a_+b)=(x+s)_\times a_+b$$ and hence $S(x,y)=(x+s,y)\in L_{a,b}$, witnessing that $S[L_{a,c}]\subseteq L_{a,b}$. Since both sets $S[L_{a,b}]$ and $L_{a,c}$ projects bijectively onto the set $\{0\}\times R$, the inclusion $S[L_{a,c}]\subseteq L_{a,b}$ implies the equalities $S[L_{a,c}]=L_{a,b}$ and $S^{-1}[L_{a,b}]=L_{a,c}$. Therefore, the function $S^{-1}:R^2\to R^2$ is a translation of the affine plane $R^2$ and so is the function $S:R^2\to R^2$, $S:(x,y)\mapsto (x+s,y)$. 
\end{proof}

The characterization of vertical and horizontal translation proved in Propositions~\ref{p:tring-vshifts} and \ref{p:tring-hshifts}  imply the following characterizations of vertical-translation and horizontal-translation based affine planes suggest the following definition.

\begin{definition} A ternar $R$ is called
\begin{itemize}
\item  \index{vertical-translation ternar}\index{ternar!vertical-translation}\defterm{vertical-translation} if
$\forall x,a,b,s\in R\;\;x_\times a_+(b+s)=(x_\times a_+b)+s;$
\item  \index{horizontal-translation ternar}\index{ternar!horizontal-translation}\defterm{horizontal-translation} if $\forall x,a,b,s\in R\;\;(x+s)_\times  a_+b=x_\times a_+(s_\times a_+ b)$;
\item  \index{translation ternar}\index{ternar!translation}\defterm{translation} if it is both vertical-translation and horizontal-translation.
\end{itemize}
\end{definition}

\begin{proposition}\label{p:vtrans=>linear+plus=puls} If a ternar $R$ is  vertical-translation, then it is linear and plus\textup{=}puls.
\end{proposition}

\begin{proof}  If a ternar $R$ is vertical-translation, then for every points $x,a,s\in R$ we have
$$x_\times a_+s=x_\times a_+(0+s)=(x_\times a_+0)+s=(x\cdot a)+s,$$
witnessing that the ternar $R$ is linear. By Proposition~\ref{p:ternar-linear<=>}, the linear ternar $R$ is plus=puls.
\end{proof}

\begin{theorem}\label{t:cart-group<=>} For a based affine plane $(\Pi,uow)$ and  its ternar $\Delta$, the following conditions are equivalent:
\begin{enumerate}
\item the based affine plane $(\Pi,uow)$ is vertical-translation;
\item the ternar $\Delta$ is vertical-translation;
\item the ternar $\Delta$ is linear and associative-plus;
\item the ternar $\Delta$ is linear and associative-puls.
\end{enumerate}
 \end{theorem}
 
\begin{proof} Observe that the coordinate chart $C:\Pi\to\Delta^2$ is an isomorphism of based affine planes $\Pi$ and $\Delta^2$. 
\smallskip

The equivalence $(1)\Leftrightarrow(2)$ follows from Proposition~\ref{p:tring-vshifts}.
\smallskip 

$(2)\Ra(3)$ If the ternar $\Delta$ is vertical-translation, then it is linear, by Proposition~\ref{p:vtrans=>linear+plus=puls}. To see that the loop $(\Delta,+)$ is associative, take any points $x,b,s\in \Delta$ and observe that the vertical-translation property of $\Delta$ implies
$$x+(b+s)=x_\times 1_+(b+s)=(x_\times 1_+b)+s=(x+b)+s.$$ 
\smallskip

The equivalence $(3)\Leftrightarrow(4)$ follows from Proposition~\ref{p:ternar-linear<=>} (implying that linear ternars are plus=puls).
\smallskip

$(3)\Ra(2)$ If $\Delta$ is linear and associative-puls, then for every $x,a,b,s\in \Delta$ we have
$$x_\times a_+(b+s)=(x\cdot a)+(b+s)=((x\cdot a)+b)+s=(x_\times a_+b)+s,$$witnessing that the ternar $\Delta$ is vertical-translation.
\end{proof}

Proposition~\ref{p:tring-hshifts} implies the following characterization of horizontal-translation ternars and based affine planes.

\begin{theorem}\label{t:right-distributive<=>} A based affine plane is horizontal-translation if and only if its ternar is horizontal-translation.
\end{theorem} 

\begin{proposition}\label{p:tring-rd=>+ass} If a ternar $R$ is horizontal-translation, then it is associative-plus.
\end{proposition}

\begin{proof} For every $x,y,z\in R$, the horizontal-translativity of the ternar $R$ implies the equality
$$(x+y)+z=(x+y)_\times 1_+z=x_\times 1_+(y_\times 1_+z)=x+(y+z),$$
witnessing that the addition operation on $R$ is associative.
\end{proof}

\begin{proposition}\label{p:ltring-rdist<=>} A linear ternar is  horizontal-translation $R$ if and only if it is right-distributive and associative-plus.
\end{proposition}

\begin{proof} To prove the ``only if'' part, assume that a linear ternar $R$ is horizontal-translation. By Proposition~\ref{p:tring-rd=>+ass}, $R$ is associative-plus.  To see that the ternar $R$ is right-distributive, observe that for every $x,y,a\in R$ we have 
$$(x+y)\cdot a=(x+y)_\times a_+0=x_\times a_+(y_\times a_+0)=(x\cdot a)+(y\cdot a),$$by the linearity of $R$.

To prove the ``if'' part, assume that the linear ternar $R$ is right-distributive and associative-plus. Then for every $x,a,b,s\in R$ we have
$$(x+s)_\times a_+b=(x+s)\cdot a+b=((x\cdot a)+(s\cdot a))+b=(x\cdot a)+((s\cdot a)+b)=x_\times a_+(s_\times a_+b),$$ witnessing that the linear ternar $R$ is horizontal-translation.
\end{proof}

Theorems~\ref{t:cart-group<=>}, \ref{t:right-distributive<=>}, Proposition~\ref{p:trans<=>hv-trans} and Definition~\ref{d:quasi-semi-near-field} imply the following important characterization of translation affine planes.

\begin{theorem}[Veblen, Wedderburn, 1907]\label{t:VW-Thalesian<=>quasifield} For a based affine plane $\Pi$ and its ternar $\Delta$, the following conditions are equivalent:
\begin{enumerate}
\item the affine plane $\Pi$ is translation;
\item the based affine plane $(\Pi,uow)$ is horizontal-translation and vertical-translation;
\item the ternar $\Delta$ is translation;
\item the ternar $\Delta$ is linear and horizontal-translation;
\item the ternar $\Delta$ is linear, right-distributive and associative-plus;
\item the ternar $\Delta$ is linear and its biloop is a quasi-field.
\end{enumerate}
\end{theorem}

\begin{proposition}\label{p:trans=s+t} Let $R$ be a translation ternar and $R^2$ be its coordinate plane. A function $F:R^2\to R^2$ is a translation of the affine plane $R^2$ if and only if there exists a pair $(s,t)\in R^2$ such that $F(x,y)=(x+s,y+t)$ for every $(x,y)\in R^2$.
\end{proposition}

\begin{proof} To prove the ``if'' part, assume that there exists a pair $(s,t)\in R^2$ such that $F(x,y)=(x+s,y+t)$ for all $(x,y)\in R^2$. By Propositions~\ref{p:tring-hshifts} and \ref{p:tring-vshifts}, the functions $$F_s:R^2\to R^2,\;\;F_s:(x,y)\mapsto(x+s,t)\quad\mbox{and}\quad F_t:R^2\to R^2,\;F_t(x,y)\mapsto (x,y+t),$$ are translations of the affine plane $R^2$, and so is their composition $F_s\circ F_t$. Observe that for every pair $(x,y)\in R^2$, we have the equality
$$F_t\circ F_s(x,y)=F_t(x+s,y)=(x+s,y+t)=F(x,y),$$
which implies that $F=F_t\circ F_s$ is a translation of the affine plane $R^2$.
\smallskip

To prove the ``only if'' part, take any translation $F:R^2\to R^2$ of the affine plane $R^2$ and consider the pair $(s,t)\defeq F(0,0)$. By Propositions~\ref{p:tring-hshifts} and \ref{p:tring-vshifts}, the functions $$F_s:R^2\to R^2,\;\;F_s:(x,y)\mapsto(x+s,t)\quad\mbox{and}\quad F_t:R^2\to R^2,\;F_t(x,y)\mapsto (x,y+t),$$ are translations of the affine plane $R^2$, and so is their composition $F_s\circ F_t$. Observe that $$F_t\circ F_s(0,0)=F_t(0+s,0)=F_t(s,0)=(s,0+t)(s,t)=F(0,0).$$By Proposition~\ref{p:Ax=Bx=>A=B}, the translations $F$ and $F_t\circ F_s$ coincide, and hence
$$F(x,y)=F_t\circ F_s(x,y)=F_t(x+s,y)=(x+s,y+t)$$ for every pair $(x,y)\in R^2$.
\end{proof}

In the following corollary, for a ternar $R$, we endow its coordinate plane $R^2$ with the binary operation ${+}\hskip-6pt{+}:R^2\times R^2$ defined by $(x,y){+}\hskip-6pt{+}(s,t)=(x+s,y+t)$, where $+:R\times R\to R$, $+:(x,y)\mapsto x_\times 1_+ y$, is the addition operation on the ternar $R$. To simplify the notations, we shall denote the operation ${+}\hskip-6pt{+}$ on $R^2$ by $+$.

\begin{corollary}\label{c:Trans(R)=R2} If $R$ is a translation ternar, then the map $$I:\Trans(R^2)\to R^2,\quad I:F\mapsto F(0,0),$$ is an isomorphism of the groups $\Trans(R^2)$ and $(R^2,+)$.
\end{corollary}

\begin{proof} By Theorem~\ref{t:VW-Thalesian<=>quasifield}, the coordinate plane $R^2$ of the ternar $R$ is translation and hence the map  $I:\Trans(R^2)\to R^2$, $I:F\mapsto F(0,0)$, is bijective.  To show that $I$ is a group isomorphism, fix any translations $F,G\in\Trans(R^2)$. Let $(x_F,y_F)\defeq F(0,0)=I(F)$ and $(x_G,y_G)=G(0,0)=I(G)$. Applying Proposition~\ref{p:trans=s+t}, we conclude that
$$I(G\circ F)=G\circ F(0,0)=G(x_F,y_F)=(x_F+x_G,y_F+y_G)=(x_F,y_F)+(x_G,y_G)=I(F)+I(G),$$which means that the bijective map $I:\Trans(R^2)\to R^2$ is an isomorphism of the groups $\Trans(R^2)$ and $(R^2,+)$.
\end{proof}

\begin{corollary}\label{c:Trans=vectors=+} Let $\Pi$ be a based affine plane and $\Delta$ be its ternar. If the affine plane $\Pi$ is Thalesian, then the groups $\Trans(\Pi)$, $\overvector \Pi$ and $(\Delta^2,+)$ are isomorphic.
\end{corollary}

Corollary~\ref{c:Trans(R)=R2} implies that for a translation based affine plane $(\Pi,uow)$, the plus loop $(\Delta,+)$ of its ternar is isomorphic to the subgroup $\{T\in\Trans(\Pi):T(o)\in\Delta\}$ of the translation group $\Trans(\Pi)$ of the affine plane $\Pi$. We are going to show that the same conclusion holds for any diagonal-translation based affine plane. Let us recall that a based affine plane $(\Pi,uow)$ is \defterm{diagonal-translation} if for every point $z\in\Delta$, there exists a translation $T:\Pi\to\Pi$ such that $T(o)=z$.  In this case the translation $T$ is equal to the vector $\overvector{oz}$. Therefore, a based affine plane $(\Pi,uow)$ is diagonal-translation if and only if for every point $z\in\Delta$, the vector $\overvector{oz}$ is functional.

\begin{lemma}\label{l:vector-plus} Let $(\Pi,uow)$ be a based affine plane and, $\Delta$ be its ternar, and $x,y,z\in\Delta$ be any points.
\begin{enumerate}
\item If $z=x+y$, then $\overvector{ox}=\overvector{yz}$.
\item If $\overvector{ox}=\overvector{yz}\in\overvector X$, then $z=x+y$.
\item If $\overvector{ox},\overvector{oy}\in\overvector X$, then $(z=x+y\;\Leftrightarrow\; \overvector{oz}=\overvector{oy}+\overvector{ox})$.
\item If $\overvector{ox},\overvector{oy},\overvector{oz}\in\overvector X$, then $(x+y)+z=x+(y+z)$.
\end{enumerate}
\end{lemma} 

\begin{proof} 1. Assume that $z=x+y$.  Consider the points $a,b\in\Pi$ with coordinates $oy$ and $xz$, respectively. By the definition of the plus operation, the equlity $z=x+y$ implies $\Aline ab=\Aline{oy}{xz}\subparallel \Delta$. Then $ox\# ab\# yz$ and hence $\overvector{ox}=\overvector{yz}$. 
\smallskip

2. Assume that $\overvector{ox}=\overvector{yz}\in\overvector X$. Consider the point $z'\defeq x+y\in\Delta$. The preceding item implies $\overvector{yz}=\overvector{ox}=\overvector{yz'}$. Since the vector $\overvector{yz}=\overvector{yz'}$ is functional, $z=z'=x+y$.
\smallskip

3. Assume that $\overvector{ox},\overvector{oy}\in\overvector X$. If 
$z=x+y$, then the first item implies $\overvector{yz}=\overvector{ox}\in\overvector X$. By Theorem~\ref{t:vector-addition}, $\overvector{oz}=\overvector{oy}+\overvector{yz}=\overvector{oy}+\overvector{ox}$. Now assume that $\overvector{oz}=\overvector{oy}+\overvector{ox}$. By Theorem~\ref{t:vector-addition},
$$\overvector{oy}+\overvector{ox}=\overvector{oz}=\overvector{oy}+\overvector{yz}$$ and $\overvector{ox}=\overvector{yz}$. By the second item, the equality  $\overvector{ox}=\overvector{yz}$ implies $z=x+y$.
\smallskip

4. Assume that the vectors $\overvector{ox},\overvector{oy},\overvector{oz}$ are functional. By Theorem~\ref{t:vector-addition}, $(\overvector{oz}+\overvector{oy})+\overvector{ox}=\overvector{oz}+(\overvector{oy}+\overvector{ox})$, which is equivalent to the equality $(x+y)+z=x+(y+z)$, according to the preceding item.
\end{proof}

\begin{theorem}\label{t:diagonal-trans=>ass-plus} If a based affine plane $(\Pi,uow)$ is diagonal-translation, then its ternar $\Delta$ is associative-plus. More precisely, the plus loop $(\Delta,+)$ is isomorphic to the subgroup $\{T\in\Trans(\Pi):T(o)\in\Delta\}$ of the translation group $\Trans(\Pi)$. If the affine plane $\Pi$ is translation, then the ternar $\Delta$ is commutative-plus.
\end{theorem}

\begin{proof} Since the based affine plane $(\Pi,uow)$ is diagonal-translation, for every $x\in \Delta$ the vector $\overvector{ox}$ is functional and hence the function $\Phi:\Delta\to\Trans(\Pi)$, $\Phi:x\mapsto\overvector{ox}$, is well-defined. Since $x=\overvector{ox}(o)$, the function $\Phi$ is injective. 

Observe that $\Phi[\Delta]=H\defeq\{T\in\Trans(\Pi):T(o)\in\Delta\}$. Lemma~\ref{l:vector-plus} and the definition of addition of functional vectors ensure that for every points $x,y,z\in\Delta$ with $z=x+y$, we have
$$\Phi(x+y)=\Phi(z)=\overvector{oz}=\overvector{oy}+\overvector{ox}=\overvector{ox}\circ\overvector{oy}=\Phi(x)\circ\Phi(y),$$
which means that $\Phi:\Delta\to \Trans(\Pi)$ is an injective homomorphism of the loops $(\Delta,+)$ and $\Trans(\Pi)$. This implies that $H$ is a subgroup of the group $\Trans(\Pi)$ and $\Phi:\Delta\to H$ is an isomorphism of the loops $(\Delta,+)$ and $(H,\circ)$. Since the loop $(H,\circ)$ is associative, so is the loop $(\Delta,+)$. Therefore, the ternar $\Delta$ is associative-plus.

If the affine plane $\Pi$ is translation, then the group $\Trans(\Pi)$ is commutative, by Corollary~\ref{c:Trans-commutative}, and so are its subgroup $H$ and the isomorphic copy $(\Delta,+)$ of the group $H$. Since the loop $(\Delta,+)$ is commutative, the ternar $\Delta$ is commutative-plus. 
\end{proof}

Finally, we extend Theorem~\ref{t:01*=>+}  to a class of biloops containing all quasi-fields.
 
\begin{definition} A biloop $R$ is called \index{$1$-commutative biloop}\index{biloop!$1$-comutative}\defterm{$1$-commutative} if $1+x=x+1$ for every $x\in R$.
\end{definition}

\begin{theorem} Let $R$ be a right-distributive $1$-commutative biloop whose additive loop $(R,+)$ is right-inversive. The addittion operation on $R$ is uniquely determined by the multiplication operation and the unary operation $1\mbox{-}:R\to R$ assigning to every $x\in R$ a unique element $1\mbox{-}x\in R$ such that $(1\mbox{-}x)+x=1$. 
\end{theorem}

\begin{proof} Since $R$ is a $0$-loop, for every $x\in R$ and $y\in R^*\defeq R\setminus\{0\}$, there exits a unique element $x/y\in R$ such that $(x/y)\cdot y=x$. Since the loop $(R,+)$ has the right-inverse property, for every element $x\in X$ there exists an element $x'\in R$ such that $y=(y+x)+x'$ for every $y\in R$. Applying the latter equality for $y\in \{0,x'\}$, we obtain the equalities $0=(0+x)+x'=x+x'$ and $0+x'=x'=(x'+x)+x'$ implying $x'+x=0$. Therefore, $x'$ is the two-sided inverse to $x$ in the loop $(R,+)$.

For every $y\in R$, the right-distributivity of the biloop $R$ implies $$(x\cdot y)+(x\cdot y)'=0=0\cdot y=(x+x')\cdot y=(x\cdot y)+(x'\cdot y)$$
and hence $x'\cdot y=(x\cdot y)'$.
 
The definition of the unary operation $1\mbox{-}:R\to R$ and the right-inversivity of the loop $(R,+)$ ensure that $(1\mbox{-}x)+x=1=(1+x')+x$ and hence $$1\mbox{-}x=1+x'\quad\mbox{for every \;\; $x\in X$}.$$ 
 We claim that $1'$ is the unique element of the set $M\defeq R\setminus \{(1\mbox{-}(1/x))\cdot x:x\in R^*\}$. Indeed, assuming that $1'=(1\mbox{-}(1/x))\cdot x$ for some $x\in R^*$, we conclude that 
 $$0+1'=1'=(1\mbox{-}(1/x))\cdot x=(1+(1/x)')\cdot x=1\cdot x+((1/x)'\cdot x)=x+((1/x)\cdot x)'=x+1'$$and hence $x=0$, which contradicts the choice of $x$. This contradiction shows that $1'\notin M$. On the other hand, for every $z\in R\setminus\{1'\}$, the equation $z=x+1'$ has a unique solution $x\in R^*$ and then 
$$z=x+1'=1\cdot x+((1/x)\cdot x)'=1\cdot x+((1/x)'\cdot x)=(1+(1/x)')\cdot x\in M,$$witnessing that $\{1'\}=M=R\setminus\{(1\mbox{-}(1/x))\cdot x:x\in R^*\}$ and hence the element $1'$ is uniquely determined by the multiplication and the unary operation $1\mbox{-}:R\to R$.

On the other hand, for every $x\in R\setminus\{1\}$, we have $0+x=x\ne 1=(1\mbox{-}x)+x$ and hence $1\mbox{-}x\ne 0$. So, we can consider the element $1/(1\mbox{-}x)\in R^*$ and the element 
$$
\begin{aligned}
\big(1\mbox{-}(1/(1\mbox{-}x))\big)\cdot(1\mbox{-}x)&=\big(1+(1/(1\mbox{-}x))')\cdot (1\mbox{-}x)=(1\cdot (1\mbox{-}x))+((1/(1\mbox{-}x))'\cdot(1\mbox{-}x))\\
& =(1\mbox{-}x)+((1/(1\mbox{-}x))\cdot (1\mbox{-}x))'=(1+x')+1'=(x'+1)+1'\\
&=x'+(1+1')=x'+0=x'.
\end{aligned}
$$
Now we see that for every element $x\in R$, its two-sided inverse $x'$ in the loop $(R,+)$ is uniquely determined by the multiplication operation and the unary operation $1\mbox{-}:R\to R$. 

Finally, for every elements $x\in R^*$ and $y\in R$, observe that
$$x+y=(1\cdot x)+(y/x)\cdot x=(1+(y/x))\cdot x=(1\mbox{-}(y/x)')\cdot x,$$which means that the addition operation of the biloop $R$ is uniquely determined by the multiplication operation and the unary operation $1\mbox{-}:R\to R$.
\end{proof}

\begin{corollary} The addition operation on any right-distributive associative-plus biloop $R$ is uniquely determined by the multiplication operation and the unary operation $1\mbox{-}:R\to R$ assigning to every $x\in R$ a unique element $1\mbox{-}x\in R$ such that $(1\mbox{-}x)+x=1$. 
\end{corollary}

\section{Shears of based affine planes}\label{s:tring-shears}

\begin{definition}\label{d:hv-shear} Let $(\Pi,uow)$ be a based affine plane and $e$ be its diunit. An automorphism $A:\Pi\to\Pi$ is called
\begin{itemize}
\item a \index{horizontal shear}\index{shear!horizontal}\defterm{horizontal shear} if $\Aline ou\subseteq \Fix(A)$ and $A(w)\in\Aline we$; 
\item a \index{vertical shear}\index{shear!vertical}\defterm{vertical shear} if $\Aline ow\subseteq \Fix(A)$ and $A(u)\in\Aline ue$.
\end{itemize} 
\end{definition}

\begin{definition}\label{d:shear-hv} A based affine plane $(\Pi,uow)$ is called
\begin{itemize}
\item \index{horizontal-shear based affine plane}\index{based affine plane!horizontal-shear}\defterm{horizontal-shear} if for every $x\in \Aline wh$ there exists a horizontal shear $A:\Pi\to\Pi$ such that $A(w)=x$;
\item \index{vertical-shear based affine plane}\index{based affine plane!vertical-shear}\defterm{vertical-shear} if for every $y\in \Aline uv$ there exists a vertical shear $A:\Pi\to\Pi$ such that $A(u)=y$.
\end{itemize}
\end{definition}

Applying Definition~\ref{d:shear-hv} to the coordinate plane $R^2$ of a ternar $R$, we conclude that an automorphism $A:R^2\to R^2$ is  
\begin{itemize}
\item a \index{horizontal shear}\defterm{horizontal shear} iff $A(0,1)\in R\times\{1\}$ and $R\times\{0\}\subseteq\Fix(A)$;
\item a \index{vertical shear}\defterm{vertical shear} iff $A(1,0)\in \{1\}\times R$ and $\{0\}\times R\subseteq \Fix(A)$.
\end{itemize}

\begin{proposition}\label{p:(x,y)->(x,xc+y)} For a ternar $R$, its coordinate plane $R^2$, and an element $c\in R$, the following conditions are equivalent:
\begin{enumerate}
\item there exists a vertical shear $A:R^2\to R^2$ such that $A(1,0)=(1,c)$;
\item the function $A:R^2\to R^2$, $A:(x,y)\mapsto (x,x_\times c_+y)$, is an automorphism of the affine plane $R^2$;
\item $\forall x,a,b\in R\;\;x_\times c_+(x_\times a_+b)=x_\times(1_\times c_+a)_+b$.
\end{enumerate}
\end{proposition} 

\begin{proof} $(1)\Ra(2)$ It follows from $A(1,0)=(1,c)$ and $A(0,0)=(0,0)$ that $A[L_{0,0}]=L_{c,0}$. For every $y\in R$ we have $L_{0,y}\parallel L_{0,0}$ and hence the line $A[L_{0,y}]$ is parallel to the line $A[L_{0,0}]=L_{c,0}$. Since $(0,y)\in \{0\}\times R\subseteq \Fix(A)$, $(0,y)\in L_{c,y}\cap A[L_{0,y}]$ and hence $A[L_{0,y}]=L_{c,y}$.

Since $\{0\}\times R\subseteq\Fix(A)$, the automorphism $A$ is hyperfixed. If $c=0$, then $A(1,0)=(1,c)=(1,0)\in\Fix(A)$ and the hyperfixed automorphism $A$ is the identity, by Proposition~\ref{p:hyperfixed=>paracentral}. In this case $A(x,y)=(x,y)=(x,x_\times 0_+y)=(x,x_\times c_+y)$ for every $x\in R$. So, assume that $c\ne 0$. In this case, $A$ is non-identity hyperfixed automorphism of the affine plane $R^2$. By Proposition~\ref{p:hyperfixed=>paracentral}, the automorphism is paracentral with $\Fix(A)=\{0\}\times R$, which implies that for every point $x\in R\setminus\{0\}$ the line $\overline{\{(x,0),A(x,0)\}}$ is parallel to the line $\overline{\{(1,0),A(1,0)\}}=\overline{\{(1,0),(1,c)\}}=L_1$. Then $L_x=\overline{\{(x,0),A(x,0)\}}$ and hence $A[L_x]=L_x$ for every $x\in R$. Finally, $$A(x,y)\in A[L_x\cap L_{0,y}]=A[L_x]\cap A[L_{0,y}]=L_x\cap L_{c,y}=\{(x,x_\times c_+y)\}.$$

$(2)\Ra(3)$ Assume that the map $A:R^2\to R^2$, $A:(x,y)\mapsto (x,x_\times c_+y)$, is an automorphism of the affine plane $R^2$. Then for every $a,b\in R$, the line $A[L_{a,0}]$ coincides with the line $L_{1_\times c_+a,0}$ (because $A(0,0)=(0,0_\times c_+0)=(0,0)$ and $A(1,a)=(1,1_\times c_+a)$), and the line $A[L_{a,b}]$ coincides with the line $L_{1_\times c_+a,b}$ (because $A[L_{a,b}]=L_{1_\times c_+a,b}$ is a unique line which contains the points $A(0,b)=(0,0_\times c _+b)=(0,b)$ and is parallel to the line $A[L_{a,0}]=L_{1_\times c_+a,0}$). The equality $A[L_{a,b}]=L_{1_\times c_+a,b}$ implies that for every $x\in R$ we have the equality 
$$(x,x_\times c_+(x_\times a_+b))=A(x,x_\times a_+b)\in A[L_{a,b}]=L_{1_\times c_+a,b}$$ and hence $x_\times c_+(x_\times a_+b)=x_\times (1_\times c_+a)_+b.$
\smallskip

$(3)\Ra(1)$ Consider the map $A:R^2\to R^2$, $A:(x,y)\mapsto (x,x_\times c_+y)$. The axioms of a ternar ensure that the map $A$ is bijective, $A(1,0)=(1,1_\times c_+0)=(1,c)$ and $A(0,y)=(0,0_\times c_+y)=(0,y)$ for every $y\in R$. The condition (3) ensures that   $A[L_{a,b}]=L_{1_\times c_+a,b}$. Since $A[L]=L$ for every vertical line $L$, the function $A$ is an automorphism of the affine plane $R^2$. Since $\{0\}\times R\subseteq \Fix(A)$, $A$ is a vertical shear with $A(1,0)=(1,c)$.
\end{proof}

Proposition~\ref{p:(x,y)->(x,xc+y)} motivates the following definition.

\begin{definition}\label{d:tring-ldistributive} A ternar $R$ is defined to be \index{vertical-shear ternar}\index{ternar!vertical-shear}\defterm{vertical-shear} if 
$$\forall x,a,b,c\in R\;\;x_\times c_+(x_\times a_+b)=x_\times(1_\times c_+a)_+ b.$$
\end{definition}

Proposition~\ref{p:(x,y)->(x,xc+y)} and Definition~\ref{d:tring-ldistributive} imply the following characterization of vertical-shear based affine planes.

\begin{theorem}\label{t:vshear<=>vshear} A based affine plane  is vertical-shear if and only if its ternar is vertical-shear.
\end{theorem}

\begin{proposition}\label{p:ltring-ldist<=>} A linear ternar $R$ is vertical-shear if and only if it is left-distributive and associative-plus.
\end{proposition}

\begin{proof} Assume that a linear ternar $R$ is vertical-shear. We have to check that the biloop $(R,+,\cdot)$ is left-distributive and the loop $(R,+)$ is a group. To see that the biloop $(R,+,\cdot)$ is left-distributive, observe that for every $x,a,c\in R$, the vertical-shear property of the ternar $R$ implies
$$x\cdot (c+a)=x\cdot((1\cdot c)+a)=x_\times (1_\times c_+a)_+0=x_\times c_+(x_\times a_+0)=(x\cdot c)+(x\cdot a).$$

To see that the loop $(R,+)$ is associative, observe that for every $a,b,c\in R$, the  linearity and the vertical-shear property of the ternar $R$ imply
$$(c+a)+b=1_\times(1_\times c_+a)_+b=1_\times c_+(1_\times a_+b)=(1\cdot c)+((1\cdot a)+b)=c+(a+b).$$

To prove the ``if'' part, assume that the ternar $R$ is left-distributive and associative-plus. Then for every $x,a,b,c\in R$ the linearity of the ternar $R$ implies
\begin{multline*}
x_\times(1_\times c_+a)_+b=(x\cdot ((1\cdot c)+a))+b=(x\cdot(c+a))+b\\
=((x\cdot c)+(x\cdot a))+b=(x\cdot c)+((x\cdot a)+b)=x_\times c_+(x_\times a_+b),
\end{multline*}
witnessing that the ternar $R$ is vertical-shear.
\end{proof}

For a based affine plane $(\Pi,uow)$, a line $L\subseteq\Pi$ is called \index{vertical line}\index{line!vertical}\defterm{vertical} if $L\in v\defeq(\Aline ow)_\parallel$.

\begin{theorem}\label{t:vtrans+vshear=vshears} For a based affine plane $(\Pi,uow)$ and its ternar $\Delta$,  the following conditions are equivalent:
\begin{enumerate}
\item the based affine plane $\Pi$ is vertical-translation and vertical-shear;
\item the ternar $\Delta$ is vertical-translation and vertical-shear;
\item the ternar $\Delta$ is linear, left-distributive and associative-plus;
\item for every vertical line $L\subseteq\Pi$ and points $x,y\in \Pi\setminus L$ with $\Aline xy\cap L=\varnothing$, there exists a hypershear $A:\Pi\to\Pi$ such that $A(x)=y$ and $L\subseteq\Fix(A)$;
\item the based affine plane $\Pi$ is vertical-shear and for every point $y\in\Aline ow\setminus\{o\}$ there exists an automorphism $A:\Pi\to\Pi$ such that $A(o)=y$ and $\Aline ue\subseteq\Fix(A)$.
\end{enumerate}
\end{theorem}

\begin{proof} The equivalences $(1)\Leftrightarrow(2)\Leftrightarrow(3)$ follow from Theorems~\ref{t:cart-group<=>}, \ref{t:vshear<=>vshear} and Proposition~\ref{p:ltring-ldist<=>}.
\smallskip

$(1)\Ra(4)$ Assume that based affine plane $(\Pi,uow)$ is vertical-translation and vertical-shear. Given any vertical line $L\subseteq\Pi$ and points $x,y\in \Pi\setminus L$ with $\Aline xy\cap L=\varnothing$, we have to find an automorphism $A:\Pi\to\Pi$ such that $A(x)=y$ and $L\subseteq\Fix(A)$. If $x=y$, the the identity automorphism of $\Pi$ has the required property. So, we assume that $x\ne y$. If $L=\Aline ow$, then the automorphism $A$ exists because the based affine plane $(\Pi,uow)$ is vertical-shear. So, assume that $L\ne\Aline ow$. By Proposition~\ref{p:lines3+}, there exists a line $L''\subseteq\Pi\setminus(L\cup\Aline ow)$. Choose any point $z\in L$ and consider the unique points $x'\in \Aline xz\cap \Aline ow$, $x''\in \Aline xz\cap L''$, $y'\in\Aline yz\cap\Aline ow$, and $y''\in \Aline yz\cap L''$. Since the based affine plane $(\Pi,uow)$ is vertical-translation, there exists a translation $T:\Pi\to\Pi$ such that $T(x')=y'$. Since $(\Pi,uow)$ is vertical-shear, there exists a hypershear $S:\Pi\to\Pi$ such that $S(T(x''))=y''$ and $\{x',y'\}\subseteq \Aline ow\subseteq\Fix(S)$. Then the composition $A\defeq S\circ T:\Pi\to\Pi$ is an automorphism such that $A(x'')=S\circ T(x'')=y''$ and $A(x')=S\circ T(x')=S(y')=y'$. 

\begin{picture}(200,115)(-140,-35)

\put(0,-20){\line(0,1){80}}
\put(-3,-30){$L$}
\put(30,-20){\line(0,1){80}}
\put(60,-20){\line(0,1){80}}
\put(90,-20){\line(0,1){80}}
\put(87,-30){$L''$}
\put(0,0){\line(1,0){90}}
\put(0,0){\line(2,1){90}}
\put(30,15){\line(1,0){60}}
{\linethickness{0.8pt}
\put(90,0){\color{blue}\vector(0,1){15}}
\put(30,0){\color{blue}\vector(0,1){15}}
\put(90,15){\color{cyan}\vector(0,1){30}}
}

\put(0,0){\circle*{3}}
\put(-8,-2){$z$}
\put(30,0){\circle*{3}}
\put(32,2){$x'$}
\put(60,0){\circle*{3}}
\put(62,2){$x$}
\put(90,0){\circle*{3}}
\put(92,-2){$x''$}
\put(80,4){\color{blue}$T$}
\put(30,15){\color{cyan}\circle*{3}}
\put(32,21){$y'$}
\put(60,30){\circle*{3}}
\put(62,36){$y$}
\put(90,45){\circle*{3}}
\put(93,45){$y''$}
\put(90,15){\circle*{3}}
\put(93,13){$T(x'')$'}
\put(81,25){\color{cyan}$S$}
\put(30,-20){\circle*{3}}
\put(27,-29){$o$}
\put(30,60){\circle*{3}}
\put(27,63){$w$}

\end{picture}

Since $T$ is a vertical translation and $S$ is a vertical shear, for every line $\Lambda\parallel \Aline ow$, we have
$$A[\Lambda]=S[T[\Lambda]]=S[\Lambda]=\Lambda,$$
which implies that the automorphism $A$ is paracentral. In particular, 
$A[L]=S[T[L]]=S[L]=L$ and hence $A(z)\in A[\Aline {x'}{x''}\cap L]=A[\Aline{x'}{x''}]\cap A[L]=\Aline{y'}{y''}\cap L=\{z\}$, which implies that $A$ is not a translation. By Proposition~\ref{p:aff-paracentral<=>}, the paracentral automorphism $A$ is hyperfixed and $\Fix(A)$ is a line such that $z\in\Fix[A]\parallel \Aline {x'}{x''}=\Aline ow\parallel L$ and hence $\Fix[A]=L$. Since $A$ is paracentral and $\Aline xy\cap \Fix[A]=\Aline xy\cap L=\varnothing$, $A[\Aline xy]=\Aline xy$. Then 
$$A(x)\in A[\Aline z{x'}\cap\Aline {x}{y}]=A[\Aline z{x'}]\cap \Aline xy=\Aline z{y'}\cap \Aline xy=\{y\}.$$ 

The implication $(4)\Ra(5)$ is obvious.
\smallskip

$(5)\Ra(1)$ Assume that the condition (5) is satisfied. Then the based affine plane $(\Pi,uow)$ is vertical-shear. To prove that $(\Pi,uow)$ is vertical-translation, take any point $y\in \Aline ow\setminus\{o\}$. By the condition (5), there exists an automorphism $A:\Pi\to\Pi$ such that $A(o)=y$ and $\Aline ue\subseteq\Fix(A)$. Since the affine plane $\Pi$ is Playfair, there exists  a unique point $z\in \Aline ue$ such that $\Aline yz\parallel \Aline ou$. Since the based affine plane $(\Pi,uow)$ is vertical-shear, there exists a hypershear $S:\Pi\to\Pi$ such that $S(u)=z$ and $\Aline ow\subseteq\Fix(S)$. By Lemma~\ref{l:hypershear2=translation}, the automorphism $T\defeq SA=AS$ is a translation of the affine plane $\Pi$ such that $T(o)=SA(o)=S(y)=y$, witnessing that the based affine plane $(\Pi,uow)$ is vertical-translation.
\end{proof}

\begin{theorem}\label{t:semifield<=>} For a based affine plane $(\Pi,uow)$ and its ternar $\Delta$,  the following conditions are equivalent:
\begin{enumerate}
\item the ternar $\Delta$ is translation and vertical-shear;
\item the ternar $\Delta$ is linear, horizontal-translation and vertical-shear;
\item the ternar $\Delta$ is linear, distributive and associative-plus\footnote{Let us recall that distributive associative-plus biloops are called \defterm{semi-fields}.};
\item the based affine plane $(\Pi,uow)$ is translation and vertical-shear;
\item the based affine plane $(\Pi,uow)$ is horizontal-translation and vertical-shear.
\end{enumerate}
\end{theorem}

\begin{proof} $(1)\Ra(2)$ If the ternar $\Delta$ is translation and vertical-shear, then $\Delta$ is horizontal-translation and linear, by Theorem~\ref{t:cart-group<=>}.
\smallskip

$(2)\Ra(3)$ If the ternar $\Delta$ is linear, horizontal-translation and vertical-shear, then by Propositions~\ref{p:ltring-rdist<=>} and \ref{p:ltring-ldist<=>},the biloop $(\Delta,+,\cdot)$ is distributive and associative-plus.
\smallskip

The implication $(3)\Ra(1)$ follows from Propositions~\ref{p:ltring-rdist<=>} and \ref{p:ltring-ldist<=>}.
\smallskip

The equivalence $(1)\Leftrightarrow(4)$ follows from Theorems~\ref{t:VW-Thalesian<=>quasifield}, \ref{t:vshear<=>vshear}, and $(4)\Ra(5)$ is trivial.
\smallskip

$(5)\Ra(4)$ Assume that the based affine plane $(\Pi,uow)$ is horizontal-translation and vertical-shear. We claim that the based affine plane $(\Pi,uow)$ satisfies the condition (5) of Theorem~\ref{t:vtrans+vshear=vshears}. Given any point $y\in\Aline ow\setminus\{o\}$, we have to find an automorphism $A:\Pi\to\Pi$ such that $A(o)=y$ and $\Aline ue\subseteq\Fix(A)$. Since the based affine plane $(\Pi,uow)$ is horizontal-translation, there exists a horizontal translation $H:\Pi\to\Pi$ such that $H(o)=u$. Then for the points $y'\defeq H^{-1}(y)$, and $u'\defeq H^{-1}(o)\in \Aline ou$ the line $\Aline{u'}{y'}=H^{-1}[\Aline oy]$  is parallel to the line $\Aline oy=\Aline ow$. Since the based affine plane  $(\Pi,uow)$ is vertical-shear, there exist vertical shear $S:\Pi\to\Pi$ such that $S(u')=y'$. 
Consider the automorphism $A\defeq HSH^{-1}:\Pi\to\Pi$ and observe that $S$ is a hypershear of the affine plane $\Pi$ such that $\Aline u{e}=H[\Aline ow]=H[\Fix(S)]=\Fix(A)$ and $A(o)=HSH^{-1}(o)=HS(u')=H(y')=y$. By Theorem~\ref{t:vtrans+vshear=vshears} and the based affine plane $(\Pi,uow)$ is vertical-translation, and by  Proposition~\ref{p:trans<=>hv-trans}, the vertical-translation horizontal-translation based affine plane $(\Pi,uow)$ is translation.
\end{proof}

In the following proposition, for every elements $x\in R^*$ and $y\in R$ of a biloop $R$, we denote by $x\backslash y$ and $y/x$ the unique elements of $R$ such that $x{\cdot}(x\backslash y)=y=(y/x){\cdot} x$.

\begin{proposition}\label{p:horizontal-shear} For a translation ternar $R$, its coordinate plane $R^2$ and an element $c\in R^*\defeq R\setminus\{0\}$, the following conditions are equivalent:
\begin{enumerate}
\item there exists a horizontal shear $A:R^2\to R^2$ such that $A(0,1)=(1/c,1)$;
\item the function $A:R^2\to R^2$, $A:(x,y)\mapsto (x+(y/c),y)$, is an automorphism of the affine plane $R^2$;
\item for every $a\in R^*$ either $y/a+y/c=0$ for all $y\in R$ or there exists an element $b\in R^*$ such that $y/a+y/c=y/b$ for all $y\in R$.
\end{enumerate}
\end{proposition}

\begin{proof} By Theorems~\ref{t:VW-Thalesian<=>quasifield} and \ref{t:diagonal-trans=>ass-plus} , the translation ternar $R$ is linear, right-distributive, associative-plus, and commutative-plus.
\smallskip

$(1)\Ra(2)$ Assume that $A:R^2\to R^2$ is a horizontal shear with $A(0,1)=(1/c,1)$. Given any pair $(x,y)\in R^2$, we have to prove that $A(x,y)=(x+(y/c),y)$. If $y=0$, then $(x,y)=(x,0)\in \Fix(A)$ and hence $A(x,y)=(x,y)=(x+0/c,y)=(x+y/c,y)$. So, assume that $y\ne 0$. Since $A$ is a horizontal shear with $A(0,1)=(c,1)\ne(0,1)$, Proposition~\ref{p:hyperfixed=>paracentral} ensures that the hypershear $A$ is paracentral and has $\Fix(A)=R\times\{0\}$. Then $A(x,y)\ne (x,y)$ and $\overline{\{(x,y),A(x,y)\}}\parallel \overline{\{(0,1),A(0,1)\}}=\overline{\{(0,1),(1/c,1)\}}=L_{0,1}$, which implies $\overline{\{(x,y),A(x,y)\}}=L_{0,y}$.

Observe that $A[L_x]\parallel A[L_0]=L_{c,0}$, because $A(0,0)=(0,0)$ and $A(0,1)=(1/c,1)\in L_{c,0}$. Then $A[L_x]=L_{c,b}$ for a unique element $b\in R$ such that $(x,0)=A(x,0)\in L_{c, b}$. Since $R$ is linear and $(R,+)$ is a group, $0=x_\times c_+b=x\cdot c+b$ implies $b=-(x{\cdot}c)$. Then $A(x,y)\in L_{0,y}\cap L_{c,b}=\{(z,y)\}$ where $z\in R$ is a unique element such that $$y=z_\times c_+b=(z\cdot c)-(x\cdot c)=(z-x)\cdot c$$and hence $z=x+(y/c)$ and  $A(x,y)=(z,y)=(x+(y/c),y)$.
\smallskip

$(2)\Ra(3)$ If the function $A:R^2\to R^2$, $A:(x,y)\mapsto (x+(y/c),y)$, is an automorphism of the affine plane $R^2$, then for every $a\in R^*$, the set $A[L_{a,0}]$ is a line in the affine plane $R^2$. Observe that $(1/a,1)\in L_{a,0}$ and hence $((1/a)+(1/c),1)=A(1/a,1)\in A[L_{a,0}]$. If $(1/a)+(1/c)=0$, then the line $A[L_{a,0}]$ is vertical and hence for every $y\in R$ the pair $(y/a,y)$ belongs to the line $L_{a,0}$ and its image $((y/a)+(y/c),y)=A(y/a,y)\in A[L_{a,0}]=L_0$ has $(y/a)+(y/c)=0$. If $(1/a)+(1/c)\ne 0$, then the line $A[L_{a,0}]$ is not vertical and hence it is equal to a line $L_{d,0}$ for some $d\in R^*$. For every $y\in R$, the point $(y/a,y)\in R^2$ belongs to the line $L_{a,0}$ and hence $((y/a)+(y/c),y)=A(y/a,y)\in L_{d,0}$, which implies $y/a+y/c=y/d$.
\smallskip

$(3)\Ra(1)$ Assume that the condition (3) holds. Given any element $c\in R^*$, consider the function $A:R^2\to R^2$, $A:(x,y)\mapsto (x+(y/c),y)$. It is clear that $A(x,0)=(x,0)$ and $A(0,1)=(1/c,1)$. We claim that $A$ is an automorphism of the affine plane $R^2$. Given any line $L\subseteq R^2$, we should prove that $A[L]$ is a line. If the line $L$ is horizontal, then $A[L]=L$. If $L$ is vertical, then $L=L_{d/c}$ for some $d\in R$. Then for every $(x,y)\in L=L_{d/c}$, we have $A(x,y)=(d/c+(y/c),y)\in L_{c,-d}$ because $(d/c+y/c)\cdot c+(-d)=y$. Therefore, $A[L]=L_{c,-d}$. So, assume that $L$ is not horizontal and not vertical. Then $L=L_{a,b}$ for some $a\in R^*$ and $b\in R$. If $(1/a)+(1/c)=0$, then the condition (3) ensures that $(y/a)+(y/c)=0$ for all $y\in R$. Then for every $(x,y)\in L=L_{a,b}$ we have $y=(x\cdot a)+b=(x+b/a)\cdot a$ and hence $x=(y/a)-(b/a)$ and $A(x,y)=((y/a)-(b/a)+(y/c),y)=((y/a)+(y/c)-(b/c),y)=(-(b/c),y)$. Then $A[L]=L_{-(b/c)}$ is a vertical line. If $(1/a)+(1/c)\ne0$, then by the condition (3), there exists an element $d\in R^*$ such that $(y/a)+(y/c)=(y/d)$ for all $y\in R$. In this case,  for every $(x,y)\in L=L_{a,b}$ we have $y=(x\cdot a)+b=(x+b/a)\cdot a$ and hence $x=(y/a)-(b/a)$ and $A(x,y)=((y/a)-(b/a)+(y/c),y)=((y/a)+(y/c)-(b/c),y)=(y/d-(b/c),y)\in L_{d,(b/c){\cdot}d}$, so $A[L]=L_{d,(b/c){\cdot}d}$ is a line, witnessing that the map $A$ is an automorphism of the afine plane $R^2$ such that $R\times\{0\}=\Fix(A)$ and $A(0,1)=(1/c,1)$. 
\end{proof}

Let us recall that an affine plane $\Pi$ is called \index{shear affine plane}\index{affine plane!shear}\defterm{shear} if for every line $L\subseteq \Pi$ and points $x,y\in \Pi$ with $\Aline xy\cap L=\varnothing$, there exists an automorphism $A:\Pi\to\Pi$ such that $A(x)=y$ and $L\subseteq \Fix(A)$. By Theorem~\ref{t:shear=>translation}, every shear affine plane is translation.

\begin{theorem}\label{t:shear<=>alternative} For any based affine plane $(\Pi,uow)$ and its ternar $\Delta$, the following conditions are equivalent:
\begin{enumerate}
\item the ternar $\Delta$ is linear, distributive, associative-plus and alternative-dot;
\item the based affine plane $(\Pi,uow)$ is vertical-shear, horizontal-shear, and translation;
\item the based affine plane $(\Pi,uow)$ is vertical-shear, horizontal-shear, and $A(o)\ne o$ for some automorphism $A:\Pi\to\Pi$;
\item the affine plane $\Pi$ is shear.
\end{enumerate}
\end{theorem}

\begin{proof} $(1)\Ra(2)$ Assume that the ternar $\Delta$ is linear, distributive, associative-plus and alternative-dot. By Theorem~\ref{t:semifield<=>}, the based affine plane $\Pi$ is translation and vertical-shear. To show that the based affine plane $(\Pi,uow)$ is horizontal-shear, we shall apply Proposition~\ref{p:horizontal-shear} and Theorem~\ref{t:alternative<=>right-Bol}. Observe that the coordinate chart $C:\Pi\to\Delta^2$ is an isomorphism of based affine planes $\Pi$ and $\Delta^2$. Thus it suffices to check that the based affine plane $\Delta^2$ is horizontal-shear.

By Proposition~\ref{p:horizontal-shear}, it suffices to check that for every elements $a,c\in \Delta^*\defeq\Delta\setminus\{o\}$ either $y/a+y/c=0$ for all $y\in\Delta$ or there exists an element $d\in\Delta$ such that $y/a+y/c=y/d$ for all $d\in\Delta$. Observe that the diunit $e$ of the affine base $uow$ is the identity of the $0$-loop $(\Delta,\cdot)$. 

If $e/a+e/c=o$, then $$o=o\cdot a=(e/a+e/c)\cdot a=(e/a)\cdot a+(e/c)\cdot a=e+(e/c)\cdot a,$$ by the right-distributivity of the biloop $(\Delta,+,\cdot)$. Since $(\Delta,+)$ is a loop, there exists a unique element $c'\in\Delta$ such that $c+c'=o$. The left-distributivity of the biloop $(\Delta,+,\cdot)$ implies
$$e+(e/c)\cdot a=o=(e/c)\cdot o=(e/c)\cdot(c+c')=(e/c)\cdot c+(e/c)\cdot c'=e+(e/c)\cdot c'$$and hence $a=c'$ (because $(\Delta,+)$ and $(\Delta,\cdot)$ are quasigroups). For every $y\in\Delta$, the distributivity of the biloop $(\Delta,+,\cdot)$ ensures that
$$o=c+c'=(y/c)\cdot (c+c')=(y/c)\cdot c+(y/c)\cdot c'=y+(y/c)\cdot a=(y/a)\cdot a+(y/c)\cdot a=((y/a)+(y/c))\cdot a$$and hence $(y/a)+(y/c)=o$ because $a\ne o$ and $(\Delta,+)$ is a $0$-loop.

Next, assume that $e/a+e/c\ne o$. Since the biloop $(\Delta,+,\cdot)$ is an alternative ring, we can apply Theorem~\ref{t:alternative<=>right-Bol} and find an element $d\in \Delta$ such that $y/a+y/c=y/d$ for all $y\in\Delta$. Applying Proposition~\ref{p:horizontal-shear}, we conclude that the based affine plane $\Delta^2$ is horizontal-shear and so is its isomorphic copy $(\Pi,uow)$.
\smallskip

$(2)\Ra(1)$ Assume that the based affine plane $(\Pi,uow)$ is vertical-shear, horizontal-shear and translation. By Theorem~\ref{t:semifield<=>}, the ternar $\Delta$ of the based affine plane $(\Pi,uow)$ is linear, distributive and associative-plus, so the biloop $(R,+,\cdot)$ is a division ring. Since the based affine plane $\Pi$ is horizontal-shear, so is its isomorphic copy $\Delta^2$. Applying Proposition~\ref{p:horizontal-shear}, we conclude that for every elements $a,b\in\Delta^*$ with $(e/a)+(e/c)\ne o$, there exists an element $d\in \Delta^*$ such that $(y/a)+(y/c)=y/d$ for all $y\in\Delta$. Applying Theorem~\ref{t:alternative<=>right-Bol}, we conclude that the biloop (which is a division ring) $\Delta$ is alternative-dot and so is the ternar $\Delta$. 
\smallskip

The implications $(2)\Ra(3)$ is obvious.
\smallskip

The implication $(3)\Ra(4)$ will be deduced from the following lemma.

\begin{lemma}\label{l:Lshear} If a based affine plane $(\Pi,uow)$ is vertical-shear and horizontal-shear, then for every line $L\subseteq \Pi$ with $o\in L$ and every points $x,y\in\Pi$ with $\Aline xy\cap L=\varnothing$, there exists a hypershear $S:\Pi\to\Pi$ such that $S(x)=y$ and $L\subseteq\Fix(S)$.
\end{lemma}

\begin{proof}  If the line $L$ is vertical, then the existence of the hypershear $S$ follows from the vertical-shear property of the based affine plane. So, assume that $L$ is not vertical and  find a unique point $u'\in L$ such that $\Aline u{u'}\parallel \Aline ow$. Since the based affine plane $(\Pi,uow)$ is vertical-shear, there exists a vertical shear $A:\Pi\to\Pi$ such that $A(u)=u'$ and $\Aline ow\subseteq\Fix(A)$. It follows that $A[\Aline ou]=\Aline o{u'}=L$. Consider the points $x'\defeq A^{-1}(x)$ and $y'\defeq A^{-1}(y)$, and observe that
$$\Aline{x'}{y'}\cap \Aline ou=A^{-1}[\Aline xy\cap L]=A[\varnothing]=\varnothing.$$ Since the based affine plane $(\Pi,uow)$ is horizontal-shear, there exists a horizontal shear $H:\Pi\to\Pi$ such that $H(x')=y'$ and $\Aline ou\subseteq\Fix(H)$. Then the automorphism $S\defeq AHA^{-1}$ of $\Pi$ is a hypershear such that $S(x)=AHA^{-1}(x)=AH(x')=A(y')=y$ and $L=A[\Aline ou]\subseteq A[\Fix[H]]=\Fix[S]$.
\end{proof}

Now we can present the proof of the implication $(3)\Ra(4)$. Assume that a based affine plane $(\Pi,uow)$ is horizontal-shear, vertical-shear and there exists an automorphism $A:\Pi\to\Pi$ such that $A(o)\ne o$. Then $A(o)\notin\Aline ou$ or $A(o)\notin\Aline ow$. We lose no generality assuming that $A(o)\notin\Aline ow$. In this case we shall prove that the based affine plane $(\Pi,uow)$ is vertical-translation. Given any point $y\in\Aline ow\setminus\{o\}$, we have to find a translation $T$ of $\Pi$ such that $T(o)=y$. Since $\Pi$ is Playfair, there exists a unique line $L\subseteq\Pi$ such that $A(o)\in L\subseteq \Pi\setminus \Aline ow$.  Consider the unique point $x\in \Aline ou\cap L$ and observe that $x\ne o$. Since the plane $\Pi$ is Playfair, there exists a unique point $z\in L\cap\Aline yh$. Since the based affine plane $(\Pi,uow)$ is vertical-shear, there exists a vertical shear $S:\Pi\to\Pi$ such that $S(x)=z$ and $\Aline ow\subseteq\Fix(S)$. Consider the line $A^{-1}[L]$ and the points $o'\defeq A^{-1}(o)$ and $y'\defeq A^{-1}(y)$. Observe that $o\in A^{-1}[L]$ and $\Aline {o'}{y'}\cap A^{-1}[L]=A^{-1}[\Aline oy\cap L]=\varnothing$. By Lemma~\ref{l:Lshear}, there exist a hypershear $H:\Pi\to\Pi$ such that $H(o')=y'$ and $A^{-1}[L]\subseteq\Fix(H)$. Then $S'=AHA^{-1}$ is a hypershear of $\Pi$ such that $S'(o)=AHA^{-1}(o)=AH(o')=A(y')=y$ and $L\subseteq \Fix(S')$. Proposition~\ref{p:hyperfixed=>paracentral} implies that $\Fix(S)=\Aline ow$ and $\Fix(S')=L$. Then $\Fix(S)\cap \Fix(S')=\Aline ow\cap L=\varnothing$. Since $\Aline ox\parallel \Aline yz$, Lemma~\ref{l:hypershear2=translation} ensures that the map $T\defeq S'S$ is a translation with $T(o)=SS'(o)=S(y)=y$, witnessing that the based affine plane $(\Pi,uow)$ is vertical-translation. Then the based affine plane $(\Pi,wou)$ is horizontal-translation. Since the based affine plane $(\Pi,uow)$ is horizontal-shear, the based affine plane $(\Pi,wou)$ is vertical-shear. By Theorem~\ref{t:semifield<=>}, the affine plane $\Pi$ is translation. 

To show that $\Pi$ is shear, take any line $L$ and points $x,y\in \Pi$ with $\Aline xy\cap L=\varnothing$. Since $\Pi$ is a translation plane, there exists a translation $T:\Pi\to\Pi$ such that $T(o)\in L$. Consider the points $x'\defeq T^{-1}(x)$ and $y'\defeq T^{-1}(y)$ and observe that $\Aline {x'}{y'}\cap T^{-1}[L]=T^{-1}[\Aline xy\cap L]=\varnothing$. By Lemma~\ref{l:Lshear}, there exists a hypershear $H:\Pi\to\Pi$ such that $H(x')=y'$ and $T^{-1}[L]\subseteq\Fix(H)$. Then $S\defeq THT^{-1}$ if a hypershear of the plane $\Pi$ such that $S(x)=THT^{-1}(x)=TH(x')=T(y')=y$ and $L\subseteq H[\Fix(H)]=\Fix(S)$, whitnessing that the affine plane $\Pi$ is shear.
\end{proof}

\begin{Exercise}[Yaqub, 1960]\label{ex:vh-shear} Find a based affine plane $(\Pi,uow)$ which is horizontal-shear, vertical-shear and not shear.
\smallskip

{\em Hint:} Look at the complement to a line in the projective completion of the Multon plane.
\end{Exercise}

\begin{remark} By a deep group-theoretic result of Hering and Kantor \cite{HK1971}, every {\em finite} vertical-shear horizontal-shear based affine plane is Desarguesian. So, the non-shear vertical-shear horizontal-shear based affine plane in Exercise~\ref{ex:vh-shear} is necessarily infinite.
\end{remark}

\section{Hyperscales of based affine planes}\label{s:tring-scales}

\begin{definition}\label{d:hvd-shear} Let $(\Pi,uow)$ be a based affine plane. An automorphism $A:\Pi\to\Pi$ is called
\begin{itemize}
\item a \index{horizontal scale}\index{scale!horizontal}\defterm{horizontal scale} if $A(u)\in\Aline ou$ and $\Aline ow\subseteq \Fix(A)$; 
\item a \index{vertical scale}\index{scale!vertical}\defterm{vertical scale} if $A(w)\in\Aline ow$ and $\Aline ou\subseteq \Fix(A)$;
\item a \index{diagonal scale}\index{scale!diagonal}\defterm{diagonal scale} if $A(o)\in\Aline oe$ and $\Aline ue\subseteq \Fix(A)$.
\end{itemize} 
\end{definition}

\begin{definition} A based affine plane $(\Pi,uow)$ is called
\begin{itemize}
\item \index{horizontal-scale based affine plane}\index{based affine plane!horizontal-scale}\defterm{horizontal-scale} if for every $x\in \Aline ou\setminus\{o\}$ there exists a horizontal scale $A:\Pi\to\Pi$ such that $A(u)=x$;
\item \index{vertical-scale based affine plane}\index{based affine plane!vertical-scale}\defterm{vertical-scale} if for every $y\in \Aline ow\setminus\{o\}$ there exists a vertical scale $A:\Pi\to\Pi$ such that $A(w)=y$.
\item \index{diagonal-scale based affine plane}\index{based affine plane!diagonal-scale}\defterm{diagonal-scale} if for every $z\in \Aline oe\setminus\{e\}$ there exists a diagonal scale $A:\Pi\to\Pi$ such that $A(o)=z$.
\end{itemize}
\end{definition}

Applying Definition~\ref{d:trans-hv} to the coordinate plane $R^2$ of a ternar $R$, we conclude that an automorphism $A:R^2\to R^2$ is  
\begin{itemize}
\item a \index{horizontal scale}\defterm{horizontal scale} iff $A(1,0)\in R\times\{0\}$ and $\{0\}\times R\subseteq\Fix(A)$;
\item a \index{vertical scale}\defterm{vertical scale} iff $A(0,1)\in \{0\}\times\ R$ and $R\times \{0\}\subseteq \Fix(A)$;
\item a \index{diagonal scale}\defterm{diagonal scale} iff $A(0,0)\in \Delta$ and $\{1\}\times R\subseteq \Fix(A)$.
\end{itemize}

The following proposition provides an algebraic description of vertical scales in the coordinate planes of ternars.

\begin{proposition}\label{p:vscale} For a ternar $R$, and an element $c\in R^*\defeq R\setminus\{0\}$, the following conditions are equivalent:
\begin{enumerate}
\item there exists a vertical scale $A:R^2\to R^2$ such that $A(0,1)=(0,c)$;
\item the function $A:R^2\to R^2$, $A:(x,y)\mapsto (x,y\cdot c)$, is an automorphism of the affine plane $R^2$;
\item $\forall x,a,b\in R\;\; (x_\times a_+b)\cdot c=x_\times (a\cdot c)_+(b\cdot c)$.
\end{enumerate}
\end{proposition}

\begin{proof} $(1)\Ra(2)$ Let $A:R^2\to R^2$ be a vertical scale of $R^2$ such that $A(0,1)=(0,c)$. Then $$A(1,1)=(1,c)\in A[L_{1,0}]=A[\overline{\{(0,0),(1,1)\}}]=\overline{\{A(0,0),A(1,1)\}}=\overline{\{(0,0),(1,c)\}}=L_{c,0}$$ and for every $x,y\in R^2$ we have $$A(y,y)\in A[L_y\cap L_{1,0}]=A[L_y]\cap A[L_{1,0}]=L_y\cap L_{c,0}=\{(y,y\cdot c)\}$$ and $$(y,y\cdot c)=A(y,y)\in A[L_{0,y}]\parallel A[L_{0,0}]=L_{0,0},$$ which implies $A[L_{0,y}]=L_{0,y{\cdot}c}$, and finally, $$A(x,y)\in A[L_x\cap L_{0,y}]=A[L_x]\cap A[L_{0,y}]=L_x\cap L_{0,y{\cdot}c}=\{(x,y\cdot c)\}.$$

$(2)\Ra(3)$ Assume that the function $A:R^2\to R^2$, $A:(x,y)\mapsto (x,y\cdot c)$, is an automorphism of the affine plane $R^2$. Then for every $a,b\in R$, the set $A[L_{a,b}]$ is a line, parallel to the line $A[L_{a,0}]$, which contains the points $(0,0)$ and $A(1,a)=(1,a\cdot c)$. Then $A[L_{a,0}]=L_{a{\cdot}c}$. Observe that $A(0,b)=(0,b\cdot c)$ and hence $A[L_{a,b}]=L_{a{\cdot}c,b{\cdot}c}$, which implies
$$(x,(x_\times a_+b)\cdot c)=A(x,x_\times a_+b)\in A[L_x\cap L_{a,b}]=A[L_x]\cap A[L_{a,b}]=L_x\cap L_{{a\cdot}c,b{\cdot}c}=\{(x,x_\times(a{\cdot}c)_+(b\cdot c))\}$$for every $x\in R$.
\smallskip

$(3)\Ra(1)$ Assume that the condition (3) is satisfied. Consider the function $A:R^2\to R^2$, $A(x,y)\mapsto (x,y\cdot c)$ and observe that $A(0,1)=(0,1\cdot c)=(1,c)$. We claim that $A$ is an automorphism of the affine plane $R^2$. Given any line $L\subseteq R^2$, we should prove that $A[L]$ is a line in $R^2$. If the line $L$ is vertical, then $A[L]=L$, by the definition of $A$. So, assume that $L$ is not vertical and find elements $a,b\in R$ such that $L=L_{a,b}$. The condition (3) ensures that 
\begin{multline*}
A[L]=A[L_{a,b}]=\{A(x,x_\times a_+b):x\in R\}=\{(x,(x_\times a_+b)\cdot c):x\in R\}\\
=\{(x,x_\times(a\cdot c)_+{b\cdot c}):x\in R\}=L_{a{\cdot}c,b{\cdot}c}
\end{multline*}
is a line in $R^2$.  Therefore, $A$ is an automorphism of the affine plane $A^2$. It is clear that $A(0,1)=(0,1\cdot c)=(0,c)$ and $R\times\{0\}\subseteq \Fix(A)$, witnessing that $A$ is a vetrical scale.
\end{proof}

Now we consider horizontal scales of based affine planes.

\begin{proposition}\label{p:hscale} For a ternar $R$, and an element $c\in R^*$, the following conditions are equivalent:
\begin{enumerate}
\item there exists a horizontal scale $A:R^2\to R^2$ such that $A(1/c,0)=(1,0)$; 
\item the function $A:R^2\to R^2$, $A:(x,y)\mapsto (x\cdot c,y)$, is an automorphism of the affine plane $R^2$;
\item $\forall x,a,b\in R\;\; x_\times(c\cdot a)_+b=(x\cdot c)_\times a_+b$.
\end{enumerate}
\end{proposition}

\begin{proof} $(1)\Ra(2)$ Assume that $A:R^2\to R^2$  is a horizontal scale with $A(1/c,0)=(1,0)$. By Proposition~\ref{p:hyperfixed=>paracentral}, the hyperfixed automorphism $A$ is paracentral and hence $$A(1/c,1)=(1,1)\in A[L_{c,0}]=A[\overline{\{(0,0),(1/c,1)\}}]=\overline{\{A(0,0),A(1/c,1)\}}=\overline{\{(0,0),(1,1)\}}=L_{1,0}.$$ Then for every $x,y\in R^2$ we have $$A(x,x{\cdot}c)\in A[L_{c,0}\cap L_{0,x{\cdot}c}]=A[L_{c,0}]\cap  A[L_{0,x{\cdot}c}]=L_{1,0}\cap L_{0,x{\cdot}c}=\{(x{\cdot}c,x{\cdot}c)\}$$ and $$(x{\cdot}c,x{\cdot}c)=A(x,x{\cdot}c)\in A[L_{x}]\parallel A[L_{0}]=L_{0},$$ which implies $A[L_{x}]=L_{x{\cdot}c}$, and finally, $$A(x,y)\in A[L_x\cap L_{0,y}]=A[L_x]\cap A[L_{0,y}]=L_{x{\cdot}c}\cap L_{0,y}=\{(x\cdot c,y)\}.$$

$(2)\Ra(3)$ Assume that the function $A:R^2\to R^2$, $A:(x,y)\mapsto(x\cdot c,y)$, is an automorphism of the affine plane $R^2$. Given any elements $a,b\in R$, observe that $A(1,c\cdot a)=(1\cdot c,c\cdot a)=(c,c\cdot a)\in L_{a,0}$ and hence $A[L_{c\cdot a,0}]=L_{a,0}$. Since $A(0,b)=(0\cdot c,b)=(0,b)$, the line $A[L_{c\cdot a,b}]$ coincides with the line $L_{a,b}$. Then for every $x\in R$ we have $(x\cdot c,x_\times (c\cdot a)_+b)=A(x,x_\times (c\cdot a)_+b)\in A[L_{c\cdot a,b}]=L_{a,b}$
and hence $x_\times (c\cdot a)_+b=(x\cdot c)_\times a_+b$.
\smallskip

$(3)\Ra(1)$ Assume that the condition (3) holds. Observe that the function $A:(x,y)\mapsto (x\cdot c,y)$ is bijective and has $A(1/c,0)=(1,0)$ and $\{0\}\times R\subseteq \Fix(A)$. To prove that $A$ is an automorphism of the affine plane $R^2$, take any line $L\subseteq R^2$. If the line $L$ is vertical, then $L=L_d$ for some $d\in R$ and hence $A[L]=A[L_d]$ is the vertical line $L_{d\cdot c}$. If $L$ is not vertical, then $L=L_{a,b}$ for some elements $a,b\in R$. Since $(R,\cdot)$ is a $0$-loop, there exists a unique element $\alpha\in R$ such that $c\cdot \alpha=a$. Then $A[L]=A[L_{a,b}]=L_{\alpha,b}$. Indeed, for every $(x,y)\in L_{a,b}$, the condition (3) ensures that $$A(x,y)=A(x,x_\times a_+b)=(x\cdot c,x_\times a_+b)=(x\cdot c,x_\times(c\cdot \alpha)_+b)=(x\cdot c,(x\cdot c)_\times \alpha_+b)\in L_{\alpha,b}.$$ Therefore, $A[L]=L_{\alpha,b}$ is a line in $R^2$ and $A$ is a horizontal scale with $A(1/c,0)=(1,0)$.
\end{proof}

Proposition~\ref{p:vscale} and \ref{p:hscale} motivate the following definitions.

\begin{definition} A ternar $R$ is defined to be 
\begin{itemize}
\item \index{vertical-scale ternar}\index{ternar!vertical-scale}\defterm{vertical-scale} if $\forall x,a,b,c\in R\;\; (x_\times a_+b)\cdot c=x_\times (a\cdot c)_+(b\cdot c)$;
\item  \index{horizontal-scale ternar}\index{ternar!horizontal-scale}\defterm{horizontal-scale} if 
$\forall x,a,b,c\in R\;\;x_\times(c\cdot a)_+b=(x\cdot c)_\times a_+b$.
\end{itemize}
 \end{definition}
 
Proposition~\ref{p:vscale} implies the following characterization of vertical-scale based affine planes and their ternars.
 
 \begin{theorem}\label{t:vscale<=>vscale} A based affine plane is vertical-scale if and only if its ternar is vertical-scale.
 \end{theorem}
 
 \begin{proposition}\label{p:vscale=>multass} Every vertical-scale ternar $R$ is associative-dot.
 \end{proposition}
 
 \begin{proof} For every $x,a,c\in R$ the vertical-scale property of the ternar $R$ implies
 $$(x\cdot a)\cdot c=(x_\times a_+0)\cdot c=x_\times(a\cdot c)_+(0\cdot c)=x\cdot(a\cdot c),$$
 witnessing that the ternar $R$ is associative-dot.
 \end{proof}
 
\begin{proposition}\label{p:lv-scale=>rdist+ass-dot} A linear ternar $R$ is vertical-scale if and only if it is  right-distributive and associative-dot. 
\end{proposition}

\begin{proof} Assume that the liner ternar $R$ is vertical-scale. By Proposition~\ref{p:vscale=>multass}, the biloop $R$ is associative-dot. To see that  $R$ is right-distributive, observe that for every $x,b,c\in R$, we have
$$(x+b)\cdot c=(x_\times 1_+b)\cdot c= x_\times(1\cdot c)_+(b\cdot c)=(x\cdot c)+(b\cdot c),$$
by the linearity of the vertical-scale ternar $R$.

Now assume that $R$ is right-distributive and associative-dot. Then for every $x,a,b,c\in R$ we have
$$(x_\times a_+b)\cdot c=((x{\cdot} a){+}b)\cdot c=((x{\cdot} a){\cdot}c)+(b{\cdot }c)=(x{\cdot}(a{\cdot}c))+(b{\cdot}c)=x_\times(a{\cdot}c)_+(b{\cdot}c),$$
witnessing that the ternar $R$ is vertical-scale.
\end{proof}

Next, we characterize horizontal-scale based affine planes and their ternars.

\begin{theorem}\label{t:horizontal-scale} For a based affine plane $\Pi$ and its ternar $\Delta$, the following conditions are equivalent:
\begin{enumerate}
\item the based affine plane $\Pi$ is horizontal-scale;
\item the ternar $\Delta$ is horizontal-scale;
\item the ternar $\Delta$ is linear and associative-dot.
\end{enumerate}
\end{theorem}

\begin{proof} The equivalence $(1)\Leftrightarrow(2)$ follows from Proposition~\ref{p:hscale}.
\smallskip 

$(2)\Ra(3)$ Assume that the ternar $\Delta$ is horizontal-scale. To see that $\Delta$ is linear, observe that for every $x,c,b\in \Delta$, the horizontal-scale property of $\Delta$ implies
$$x_\times c_+b=x_\times (c\cdot 1)_+b=(x\cdot c)_\times 1_+b=(x\cdot c)+b.$$
To see that the magma $(\Delta,\cdot)$ is associative, observe that for every $x,a,c\in \Delta$ we have  
$$x\cdot(c\cdot a)=x_\times(c\cdot a)_+0=(x\cdot c)_\times a_+0=(x\cdot c)\cdot a.$$

$(3)\Ra(2)$ Assume that the ternar $\Delta$ is linear and associative-dot. Observe that for every $x,a,b,c\in \Delta$, the linearity of $\Delta$ and associativity of $(\Delta,\cdot)$ imply
$$x_\times(c{\cdot} a)_+b=(x{\cdot} (c{\cdot} a))+b=((x{\cdot} c){\cdot} a)+b=(x{\cdot} c)_\times a_+b.$$
\end{proof}

Theorems~\ref{t:horizontal-scale}, \ref{t:vscale<=>vscale} and Proposition~\ref{p:lv-scale=>rdist+ass-dot} imply the following characterization.

\begin{corollary}\label{c:hv-scale<=>} For a based affine plane $\Pi$ and its ternar $\Delta$, the following conditions are equivalent:
\begin{enumerate}
\item the based affine plane $\Pi$ is horizontal-scale and vertical-scale;
\item the ternar $\Delta$ is horizontal-scale and vertical-scale;
\item the ternar $\Delta$ is linear, right-distributive, and associative-dot.
\end{enumerate}
\end{corollary}



Theorem~\ref{t:cart-group<=>}, Proposition~\ref{p:lv-scale=>rdist+ass-dot} and Theorems~\ref{t:horizontal-scale}, \ref{t:vscale<=>vscale},  Theorems~\ref{t:VW-Thalesian<=>quasifield} imply the following characterization of based affine planes with linear right-distributive associative ternars.

\begin{corollary}\label{c:vtrans+vscale} For a based affine plane $\Pi$ and its ternar $\Delta$, the following conditions are equivalent:
\begin{enumerate}
\item the based affine plane $\Pi$ is vertical-translation and vertical-scale;
\item the ternar $\Delta$ is vertical-translation and vertical-scale;
\item the ternar $\Delta$ is translation, horizontal-scale and vertical-scale;
\item the ternar $\Delta$ is linear, right-distributive and associative.
\end{enumerate}
\end{corollary}

Theorems~\ref{t:horizontal-scale}, \ref{t:vshear<=>vshear} and Proposition~\ref{p:ltring-ldist<=>} imply the following characterization of based affine planes with linear left-distributive associative ternars.

\begin{corollary}\label{c:vertical-shear+horizontal-scale} For a based affine plane $\Pi$ and its ternar $\Delta$, the following conditions are equivalent:
\begin{enumerate}
\item the ternar $\Delta$ is linear, left-distributive and associative;
\item the ternar $\Delta$ is vertical-shear and horizontal-scale;
\item the based affine plane $\Pi$ is vertical-shear and horizontal-scale;
\end{enumerate}
\end{corollary}

Corollary~\ref{c:vertical-shear+horizontal-scale} and Proposition~\ref{p:lv-scale=>rdist+ass-dot}, \ref{p:ltring-rdist<=>} imply

\begin{corollary}\label{c:vshear+vh-scale+htrans} For a based affine plane $\Pi$ and its ternar $\Delta$, the following conditions are equivalent:
\begin{enumerate}
\item the ternar $\Delta$ is linear, distributive and associative\footnote{Linear distributive associative biloops are \defterm{corps}.};
\item the based affine plane $\Pi$ is horizontal-scale, vertical-shear, and vertical-scale;
\item the based affine plane $\Delta$ is horizontal-scale, vertical-shear and horizontal-translation.
\end{enumerate}
\end{corollary}

Next, we find an algebraic characterization of diagonal scales in coordinate planes of linear ternars. For two elements $x,y$ of a loop $(X,+)$, we shall denote by $x-y$ and $x\mp y$ the unique elements of the loop $X$ such that $(x-y)+y=x$ and $x+(x\mp y)=y$. By the cancellativity, the identities $((x+y)-y)+y=x+y$ and $x+(x\mp (x+y))=x+y$ imply $(x+y)-y=x$ and $x\mp(x+y)=y$. So, the binary operations $\mp$ and $-$ satisfy the identities
$$(x-y)+y=x=(x+y)-y\quad\mbox{and}\quad x\mp(x+y)=y=x+(x\mp y).$$ The dot operation in a ternar has priority over the operations $+,-,\mp$,  and the operations $\mp,-$ have priority over the addition, 
i.e., the expressions $$x{\cdot}y+z,\;x{\cdot}y-z,\;x{\cdot}y\mp z,\;x\mp y+z,\;x+y\mp z,\;x-y+z,\;x+y-z$$ should be read as $$(x{\cdot}y)+z,\;(x{\cdot}y)-z,\;(x{\cdot}y)\mp z,\; (x\mp y)+z,\;x+(y\mp z),\;(x-y)+z,\;x+(y-z),$$ respectively.

\begin{proposition}\label{p:diagonal-scale<=>} For a linear ternar $R$ and an element $c\in R\setminus\{1\}$, the following conditions are equivalent:
\begin{enumerate}
\item there exists a diagonal scale $A:R^2\to R^2$ such that $A(0,0)=(c,c)$;
\item the function $A:R^2\to R^2$, $A:(x,y)\mapsto (c+x{\cdot}c\mp x,c+x{\cdot}c\mp y)$, is an automorphism of the affine plane $R^2$;
\item $\forall a,b,x\in R\;\;c+x{\cdot}c\mp (x{\cdot}(a+c{\cdot}a\mp c)+b)=(c+x{\cdot}c\mp x){\cdot}a+(c{\cdot}a\mp(c+b))$.
\end{enumerate}
\end{proposition}

\begin{proof} For two elements $a,b\in R$ let $L_{a,0}$ be the line connecting the points $(0,0)$ and $(1,a)$ in the coordinate plane $R^2$, and $L_{a,b}$ be the unique line in $R^2$ that is parallel to the line $L_{a,0}$ and contains the point $(0,b)$. The linearity of the ternar $R$ ensures that the line $L_{a,b}$ has equation $y=x_\times a_+b=(x{\cdot}a)+b$ in the coordinate plane $R^2$. Also, let $V_a\defeq\{a\}\times R$ be the vertical line containing the point $(a,0)$. 
\smallskip 

$(1)\Ra(2)$ Assume that $A:R^2\to R^2$ is a diagonal scale such that $A(0,0)=(c,c)$. Since $A$ is a diagonal scale, $R\times\{1\}\subseteq\Fix(A)$ and hence $A(1,1)=(1,1)$ and $A[\Delta]=\Delta$, where $\Delta\defeq\{(x,y)\in R^2:x=y\}$ is the diagonal of the coordinate plane $R^2$. Choose any point $(x,y)\in R^2$ and and consider the image $(x',x')\in A(x,x)$ of the point $(x,x)\in \Delta$. Observe that the lines $L_{0,c}$ and $L_{0,x'}$ are parallel and hence the lines $A^{-1}[L_{0,c}]=L_{c,0}$ and $A^{-1}[L_{0,x'}]$ are also parallel. The line $A^{-1}[L_{0,x'}]$ contains the point $(x,x)$ and hence it coincides with the line $L_{c,b}$ for a unique element $b\in R$ such that $x{\cdot}c+b=x$.
The definition of the binary operation $\mp:R\times R\to R$ ensures that $b=x{\cdot}c\mp x$.  Therefore, $A^{-1}[L_{0,x'}]=L_{c,x{\cdot}c\mp x}$. Taking into account that $A(1,x)=(1,x')\in L_{0,x'}\cap L_{c,b}$, we conclude that $x'=1{\cdot}c+b=c+x{\cdot}c\mp x$. 

Observe that $(x,y)\in V_x\cap L_{1,x\mp y}$. Since $A$ is a diagonal scale, $$A(x,y)\in A[V_{x}\cap L_{1,x\mp y}]=A[V_x]\cap A[L_{1,x\mp y}]=V_{x'}\cap L_{1,x\mp y}$$ and hence $A(x,y)=(x',x'+x\mp y)=(c+x{\cdot}c\mp x,(c+x{\cdot}c\mp x)+x\mp y)$.

Observe that $A[L_{c,0}]$ is the line containing the points $(1,c)$ and $(c,c)=A(0,0)$, which implies that the line $A[L_{c,0}]=L_{0,c}$ is horizontal. Then for $\beta\defeq x{\cdot}c\mp y$,  the line $A[L_{c,\beta}]\parallel A[L_{c,0}]=L_{0,c}$ is horizontal, too, and the point $(x,y)$ belongs to the line $L_{c,\beta}$. Since the horizontal line $A[L_{c,\beta}]$ contains the point $A(0,\beta)=(c,c+\beta)$, it coincides with the horizontal line $L_{0,c+\beta}$. Then 
$$A(x,y)=(c+x{\cdot}c\mp x,(c+x{\cdot}c\mp x)+x\mp y)=
(c+x{\cdot}c\mp x,c+\beta)=(c+x{\cdot}c\mp x,c+x\cdot c\mp y).$$
\smallskip

$(2)\Ra(3)$ Assume that the function  $$A:R^2\to R^2,\quad A:(x,y)\mapsto (c+x{\cdot} c\mp x,c+x{\cdot} c\mp y),$$ is an automorphism of the affine plane $R^2$. Then for every $\alpha\in R$, the image $A[L_{\alpha,0}]$ is a line in $R^2$ and hence $A[L_{\alpha,0}]=L_{a,\beta}$ for some $a,\beta\in R$. For every $x\in R$, the point 
$$(c+x{\cdot}c\mp x,c+x{\cdot}c\mp x{\cdot}\alpha)=A(x,x{\cdot}\alpha)$$ belongs to the line $A[L_{\alpha,o}]=L_{a,\beta}$ and hence satisfies its equation
$$
c+x{\cdot}c\mp x{\cdot}\alpha=(c+x{\cdot}c\mp x){\cdot}a+\beta.$$
In particular, for $x=0$ we have the identity $$c=c+0{\cdot}c\mp 0=(c+0{\cdot}c\mp 0){\cdot}a+\beta=c{\cdot}a+\beta$$ implying $\beta=c{\cdot}a\mp c$. On the other hand, for $x=1$ we have the identity
$$\alpha=c+1{\cdot c}\mp 1{\cdot}\alpha=(c+1{\cdot}c\mp 1){\cdot}a+\beta=a+\beta=a+c{\cdot}a\mp c.$$

For any $b\in R$, the line $L_{\alpha,b}$ is parallel to the line $L_{\alpha,0}$ and hence the line $A[L_{\alpha,b}]$ is parallel to the line $A[L_{\alpha,b}]=L_{a,\beta}$. Then $A[L_{\alpha,b}]=L_{a,\gamma}$ for some $\gamma\in R$. For every $x\in R$, the point $(c+x{\cdot} c\mp x,c+x{\cdot} c\mp (x{\cdot}\alpha+b))=A(x,x{\cdot}\alpha+b)\in A[L_{\alpha,b}]=L_{a,\gamma}$ satisfies the equation of the line $L_{a,\gamma}$ and hence
$c+x{\cdot} c\mp (x{\cdot}\alpha+b)=(c+x{\cdot} c\mp x){\cdot}a+\gamma$. 
In particular, for $x=0$ we have the equality
$$c+b=(c+0{\cdot}c\mp (0{\cdot}\alpha+b)=(c+0{\cdot} c\mp 0){\cdot}a+\gamma=c{\cdot}a+\gamma$$which implies $\gamma=c{\cdot}a\mp(c+b)$. Since $\alpha=a+c{\cdot}a\mp c$, we obtain the identity
$$(c+x{\cdot} c\mp (x{\cdot}(a+c{\cdot}a\mp c)+b)=(c+x{\cdot} c\mp x){\cdot}a+c{\cdot}a\mp(c+b)$$
holding for all $x,a,b\in R$ and appearing in the condition (3). 
\smallskip

$(3)\Ra(1)$ Assume that the identity $(3)$ holds. We have to check that the function $A:R^2\to R^2$, $A:(x,y)\mapsto (c+x{\cdot} c\mp x,c+x{\cdot} c\pm y)$, is an automorphism of the affine plane $R^2$. To see that the function $A$ is bijective, we have to check that for every $a,b\in R^2$, the equation $(a,b)=(c+x{\cdot} c\mp x,c+x{\cdot} c\mp y)$ has a unique solution $(x,y)$. Indeed, the equation $a=c+x{\cdot}c\mp x$ is equivalent to the equations $c\mp a=x{\cdot}c\mp x$ and $x{\cdot}c+c\mp a=x$, which has a unique solution by the axiom {\sf(T3)} of a ternar. Since $x$ is uniquely determined, the equation $b=c+x{\cdot}c\mp y$ has a unique solution $y=x\cdot c+c\mp b$. This completes the proof of the bijectivity of the function $A$. To see that $A$ is an automorphism of the plane $R^2$, it suffices to show that for every  line $L\subseteq R^2$ the  preimage $A^{-1}[L]$ is a line. This is clear if the line $L$ is vertical. If $L$ is not vertical, then $L=L_{a,\gamma}$ for some elements $a,\gamma\in R$. Consider the elements $\alpha\defeq a+c{\cdot}a\mp c$ and $b\defeq c\mp(c{\cdot}a+\gamma)$ and observe that $\gamma=c{\cdot}a\mp(c+b)$. We claim that $A^{-1}[L]=L_{\alpha,b}$, which is equivalent to $A[L_{\alpha,b}]=L_{a,\gamma}$, by the bijectivity of $A$. Indeed, for every point $(x,y)\in L_{\alpha,b}$, we have $y=x{\cdot}\alpha+b$ and hence
$$A(x,y)=(c+x{\cdot} c\mp x,c+x{\cdot} c\mp y)= (c+x{\cdot} c\mp x,c+x{\cdot} c\mp (x{\cdot}\alpha+b)).$$
The identity (3) ensures that
$$
\begin{aligned}
&c+x{\cdot} c\mp (x{\cdot}\alpha+b)=c+x{\cdot} c\mp (x{\cdot}(a+c{\cdot}a\mp c)+b)\\
&=
(c+x{\cdot}c\mp x){\cdot}a+c{\cdot}a\mp(c+b)=(c+x{\cdot}c\mp x){\cdot}a+\gamma
\end{aligned}
$$ witnessing that $A(x,y)\in L_{a,\gamma}=L$ and hence $A[L_{\alpha,b}]=L$.
\end{proof}




Proposition~\ref{p:diagonal-scale<=>} motivates the following definition.

\begin{definition}\label{d:diagonal-scale} A ternar $R$ is defined to be \index{ternar!diagonal-scale}\index{diagonal-scale ternar}\defterm{diagonal-scale} if 
$$c+x{\cdot}c\mp (x{\cdot}(a+c{\cdot}a\mp c)+b)=(c+x{\cdot}c\mp x){\cdot}a+c{\cdot}a\mp(c+b)$$for all $a,b,c,x\in R$.
\end{definition} 

Theorem~\ref{t:horizontal-scale} and Proposition~\ref{p:diagonal-scale<=>} imply the following characterization.

\begin{theorem}\label{t:diagonal-scale<=>} A based affine plane is horizontal-scale and diagonal-scale if and only if its ternar is horizontal-scale and diagonal-scale if and only if its ternar is linear, associative-dot and diagonal-scale.
\end{theorem}

For two elements $x,y$ of a ternar $R$, their \index{diagonal product}\defterm{diagonal product}  $x\circ y$ is defined by the formula $$x\circ y\defeq x+y{\cdot}x\mp y.$$ Observe that $x\circ 0=x=0\circ x$ and $x\circ 1=1=1\circ x$, which means that $0$ and $1$ are the unit and zero of the operation $\circ$, which will be called the \index{dit operation}\defterm{dit operation} (abbreviation from {\bf di}agonal produc{\bf t}).

The dit operation $\circ$ has priority over the operations $+,-,\mp$, i.e.
$$x\circ y+z,\;x\circ y-z,\;x\circ y\mp z,\;x+y\circ z,\;x-y\circ z,\;x\mp y\circ z,$$should be read as 
$$(x\circ y)+z,\;(x\circ y)-z,\;(x\circ y)\mp z,\;\;x+(y\circ z),\;x-(y\circ z),\;x\mp (y\circ z),$$
respectively.

The dit operation $\circ$ was introduced by Pickert \cite{Pickert1959} (who denoted it by $*$) and Spencer\footnote{{\bf Jill Courtaney Donaldson Yaqub} (born {\bf Spencer}) (1931 -- 2002) was an american mathematician. She defended her Ph.D. ``Non-Desarguesian Geometries'' in 1961 under supervision of Philip D.~Watson in University of Oxford, USA. Later she worked in the Ohio State University. Jill Yaqub studied the geometry of projective and M\"obius planes and has published several influential papers on the Lenz--Barlotti classification.} \cite{Spencer1960}.

\begin{definition} A ternar $R$ is defined to be
\index{ternar!associative-dit}\index{associative-dit ternar}\defterm{associative-dit} if its dit operation is associative, i.e. $(x\circ y)\circ z=x\circ(y\circ z)$ for all $x,y,z\in R$.
\end{definition}

For associative-plus linear ternars, the following characterization was proved by Pickert \cite{Pickert1959} and Spencer~\cite{Spencer1960}.

\begin{theorem} A linear (associative-plus) ternar is associative-dit if (and only if) it is diagonal-scale.
\end{theorem}

\begin{proof} Let $o\defeq(0,0)$ and $e\defeq(1,1)$ be the origin and the diunit of the coordinate plane $R^2$ of the ternar $R$. Let $\Delta\defeq\Aline oe$ be the diagonal of the coordinate plane $R^2$ and $d\defeq\Delta_\parallel$ be its direction. Also, let $E\defeq\{e\}\times R$ be the vertical line containing the diunit of the coordinate plane $R^2$.
Let $\Aut_{d,E}(R^2)$ be the group of diagonal scales of the affine plane $R^2$. By Proposition~\ref{p:hyperfixed=>paracentral}, the map $A:\Aut_{d,E}(R^2)\to \Delta\setminus\{e\}$, $A:\Phi\mapsto \Phi(o)$, is injective. This map is surjective if and only if the coordinate plane $R^2$ is diagonal-scale.

If the linear ternar $R$ is diagonal-scale, then its coordinate plane $R^2$ is diagonal-scale, by Proposition~\ref{p:diagonal-scale<=>} and hence the map $A:\Aut_{d,E}(R^2)\to\Delta\setminus\{e\}$ is bijective. Then for every $c\in R\setminus\{1\}$, there exists a unique diagonal scale $\Phi_c$ of $R^2$ such that $\Phi_c(0,0)=(c,c)$. By Proposition~\ref{p:diagonal-scale<=>}, for every $a,b\in R\setminus\{1\}$ we have $$\Phi_a\Phi_b(0,0)=\Phi_a(b,b)=(a+b{\cdot}a\mp b,a+b{\cdot}a\mp b)=(a\circ b,a\circ b)=\Phi_{a\circ b}(0,0),$$which implies $\Phi_{a\circ b}=\Phi_a\circ\Phi_b$. Then for any $a,b\in R\setminus\{1\}$ and $c\in R$, $$\Phi_{a\circ(b\circ c)}=\Phi_a\Phi_{b\circ c}=\Phi_a(\Phi_b\Phi_c)=(\Phi_a\Phi_b)\Phi_c=\Phi_{a\circ b}\Phi_c=\Phi_{(a\circ b)\circ c}$$ and hence $(a\circ b)\circ c=a\circ (b\circ c)$. If $a=1$, then $(a\circ b)\circ c=(1\circ b)\circ c=1\circ c=1=1\circ (b\circ c)$. If $b=1$, then  $(a\circ b)\circ c=(a\circ 1)\circ c=1\circ c=1=a\circ 1=a\circ (1\circ c)=a\circ (b\circ c)$. Therefore, the dit operation on $R$ is associative.
\smallskip

Now assume that $R$ is associative-plus and associative-dit. To prove that $R$ is diagonal-scale, take any elements $a,b,c,x\in R$ and observe that the associativity of the plus and dit operations imply the equalities
$$
\begin{aligned}
&c+x{\cdot}c\mp (x{\cdot}(a+c{\cdot}a\mp c)+b)=
c-x{\cdot}c+ x{\cdot}(a\circ c)+b\\
&=c-x{\cdot}c+x-(a\circ c-x\cdot (a\circ c)+x)+a\circ c+b=c\circ x-(a\circ c)\circ x+a\circ c+b=\\
&=c\circ x-a\circ(c\circ x)+a\circ c +b=c\circ x-(a-(c\circ x)\cdot a+c\circ x)+a\circ c+b=\\
&=c\circ x-c\circ x+(c\circ x)\cdot a-a+a-c{\cdot}a+c+b=(c+x{\cdot}c\mp x){\cdot}a+c{\cdot}a\mp(c+b),
\end{aligned}
$$
witnessing that the ternar $R$ is diagonal-scale.
\end{proof}

The following problem is a major unsolved problem on diagonal-scale based affine planes. 

\begin{problem}\label{prob:dh-scale=>Desarg} Is every diagonal-scale horizontal-scale based affine plane Desarguesian? 
\end{problem}

This problem has an algebraic reformulation.

\begin{problem}\label{prob:sh-scale=>corps} Is every diagonal-scale linear associative-dot ternar a corps?
\end{problem}

\begin{remark} We shall return to Problems~\ref{prob:dh-scale=>Desarg}  and \ref{prob:sh-scale=>corps} in Section~\ref{s:discale}.
\end{remark}

\section{Central homotheties in based affine planes}\label{s:tring-homos}

\begin{definition} Let $(\Pi,uow)$ be a based affine plane. An automorphism $A:\Pi\to\Pi$ is called a \index{central homothety}\defterm{central homothety} of the based affine plane $(\Pi,uow)$ if $A$ is a dilation of $\Pi$ such that $A(o)=o$.
\end{definition}

\begin{proposition}\label{p:central-homo} For a ternar $R$, its coordinate plane $R^2$, and an element $c\in R^*\defeq R\setminus\{0\}$, the following conditions are equivalent:
\begin{enumerate}
\item there exists a homothety $A:R^2\to R^2$ such that $A(1,1)=(c,c)$ and $A(0,0)=(0,0)$;
\item the function $A:R^2\to R^2$, $A:(x,y)\mapsto (c\cdot x,c\cdot y)$, is a central homothety of the based affine plane $R^2$;
\item $\forall x,a,b\in R\;\;c{\cdot} (x_\times a_+b)=(c{\cdot} x)_\times a_+(c{\cdot }b)$.
\end{enumerate}
\end{proposition}

\begin{proof} $(1)\Ra(2)$ Assume that $A:R^2\to R^2$ is a homothety with $A(1,1)=(c,c)$ and $A(0,0)=(0,0)$. Then $A[L_1]=L_c$,  $A[L_{0,1}]=L_{0,c}$, and  $A[L]=L$ for every line $L\subseteq R^2$ containing the origin $(0,0)$ of the based affine plane $R^2$. Observe that for every $(x,y)\in R^2$ we have $$A(1,y)=A[L_1\cap L_{y,0}]=A[L_1]\cap A[L_{y,0}]=L_c\cap L_{y,0}=\{(c,c\cdot y)\}$$ and hence $A[L_{0,y}]=L_{0,c{\cdot}y}$. Then $$A(x,y)\in A[L_{x\backslash y,0}\cap L_{0,y}]=L_{x\backslash y,0}\cap L_{0,c{\cdot}y}=\big\{\big((c{\cdot}y)/(x\backslash y),c\cdot y\big)\big\}.$$ 

On the other hand, $$A(x,1)\in A[L_{x\backslash 1,0}\cap L_{0,1}]=A[L_{x\backslash 1,0}]\cap A[L_{0,1}]=L_{x\backslash 1,0}\cap L_{0,c}=\{(c/(x\backslash 1),c)\}$$ and hence $A[L_x]=L_{c/(x\backslash 1)}$ and $$A(x,y)\in A[L_x\cap L_{x\backslash y,0}]=A[L_x]\cap A[L_{x\backslash y,0}]= L_{c/(x\backslash 1)}\cap L_{x\backslash y,0}=\big\{\big(c/(x\backslash 1),(c/(x\backslash 1))\cdot(x\backslash y)\big)\}.$$

Therefore, $$A(x,y)=\big((c{\cdot}y)/(x\backslash y),c\cdot y\big)=\big(c/(x\backslash 1),(c/(x\backslash 1))\cdot(x\backslash y)\big)=(c/(x\backslash 1),c\cdot y).$$
Since $A[L_{1,0}]=L_{1,0}$, we have $A(x,x)=(c/(x\backslash 1),c\cdot x)\in L_{1,0}$ and hence $c/(x\backslash 1)=c\cdot x$ and finally,
$$A(x,y)=(c/(x\backslash 1),c\cdot y)=(c\cdot x,c\cdot y).$$

$(2)\Ra(3)$ Assume that the function $A:R^2\to R^2$, $A:(x,y)\mapsto (c\cdot x,c\cdot y)$, is a homothety of the affine plane $R^2$. Then for every $a,b\in R$, the set $A[L_{a,b}]$ is a line, parallel to the line $A[L_{a,0}]=L_{a,0}$ and hence $A[L_{a,b}]=L_{a,c\cdot b}$ (because $A(0,b)=(0,c\cdot b)\in A[L_{a,b}]$).
Then for every $x\in R$ we have the equality
$$A(x,x_\times a_+b)=(c\cdot x,c\cdot (x_\times a_+b))\in A[L_{a,b}]=L_{a,c{\cdot}b},$$implying the desired identity
$$c\cdot (x_\times c_+b)=(c\cdot x)_\times a_+(c\cdot b).$$
\smallskip

$(3)\Ra(1)$ Assume that the condition (3) is satisfied. Consider the function $A:R^2\to R^2$, $A(x,y)\mapsto (c\cdot x,c\cdot y)$, and observe that $A(1,1)=(c,c)$. We claim that $A$ is a dilation of the affine plane $R^2$. Given any line $L\subseteq R^2$, we should prove that $A[L]$ is a line in $R^2$, parallel to the line $L$. If the line $L$ is vertical, then $L=L_a$ for some $a\in R$ and $A[L]=A[L_a]=L_{c\cdot a}\parallel L_a$. So, assume that $L$ is not vertical and find elements $a,b\in R$ such that $L=L_{a,b}$. The condition (3) ensures that 
\begin{multline*}
A[L]=A[L_{a,b}]=\{A(x,x_\times a_+b):x\in R\}=\{(c\cdot x,c\cdot (x_\times a_+b)):x\in R\}\\
=\{(c\cdot x,(c\cdot x)_\times a_+(c\cdot b)):x\in R\}=L_{a,c{\cdot}b}
\end{multline*}
is a line in $R^2$, parallel to the line $L_{a,b}=L$. Since $A(0,0)=(0,0)$, $A$ is a central homothety of $R^2$ with $A(1,1)=(c,c)$.
\end{proof}

Proposition~\ref{p:central-homo} motivates the following definitions.

\begin{definition} A based affine plane $(\Pi,uow)$ is called \index{homocentral based affine plane}\index{based affine plane!homocentral}\defterm{homocentral} if for every point $x\in \Aline oe\setminus\{o\}$, there exists a central homothety $A:\Pi\to\Pi$ such that $A(e)=x$, where $e$ is the diunit of the affine base $uow$.
\end{definition}

\begin{definition} A ternar $R$ is called \index{homocentral ternar}\index{ternar!homocentral}\defterm{homocentral} if 
 $$\forall x,a,b,c\in R\;\; c{\cdot}(x_\times a_+b)=(c{\cdot}x)_\times a_+(c{\cdot} b).$$
 \end{definition}

Proposition~\ref{p:central-homo} implies the following characterization of homocentral based affine planes.
 
 \begin{theorem}\label{t:homocentral<=>} A based affine plane is homocentral if and only if its ternar is homocentral.
 \end{theorem}
 
 \begin{proposition}\label{p:tring-homocentral<=>} A (linear) ternar is homocentral (if and) only if it is left-distributive and associative-dot.
 \end{proposition}

\begin{proof} Assume that a ternar $R$ is homocentral. To prove that the magma $(R,\cdot)$ is associative, observe that for every elements $x,a,c\in R$, the homocentrality of $R$ implies
$$c\cdot (x\cdot a)=c\cdot (x_\times a_+0)=(c\cdot x)_\times a_+(c\cdot 0)=(c\cdot a)_\times a_+0=(c\cdot x)\cdot a.$$
To see that $R$ is left-distributive, observe that for every $x,c,b\in R$  the homocentrality of $R$ implies
$$c\cdot (x+b)=c\cdot (x_\times 1_+b)=(c\cdot x)_\times 1_+(c\cdot b)=(c\cdot x)+(c\cdot b).$$

Now assume that the ternar $R$ is linear, left-distributive and associative-dot. To see that the ternar $R$ is homocentral, observe that for every $x,a,b\in R$ we have
$$c\cdot (x_\times a_+b)=c\cdot ((x\cdot a)+b)=(c\cdot (x\cdot a))+(c\cdot b)=((c\cdot x)\cdot a)+(c\cdot b)=(c\cdot x)_\times a_+(c\cdot b).$$
\end{proof}

Theorems~\ref{t:horizontal-scale}, \ref{t:homocentral<=>}, and Proposition~\ref{p:tring-homocentral<=>} imply the following characterization.

\begin{theorem}\label{t:hscale+homocentral} For a based affine plane $\Pi$ and its ternar $\Delta$, the following conditions are equivalent:
\begin{enumerate}
\item the based affine plane $\Pi$ is horizontal-scale and homocentral;
\item the ternar $\Delta$ is horizontal-scale and homocentral;
\item the ternar $\Delta$ is linear, left-distributive and associative-dot.
\end{enumerate}
\end{theorem}

\begin{definition} A point $c$ of an affine plane $\Pi$ is called \index{central point}\index{affine plane!central point of}\defterm{central} if for every points $x\in \Pi\setminus\{c\}$ and $y\in\Aline cx\setminus\{c\}$, there exists a homothety $H\in\Dil_c(\Pi)$ such that $H(x)=y$.
\end{definition}

Observe that an affine plane $\Pi$ is dilation if and only if every point $c\in\Pi$ is central. A based affine plane $(\Pi,uow)$ is homocentral if and only if the origin $o$ is a central point of $\Pi$. 

\begin{theorem}\label{t:leftcorps} For a based affine plane $(\Pi,uow)$ and its ternar $\Delta$, the following conditions are equivalent:
\begin{enumerate}
\item the ternar $\Delta$ is linear, left-distributive and associative;
\item the based affine plane $(\Pi,uow)$ is vertical-translation and homocentral;
\item the based affine plane $(\Pi,uow)$ is vertical-shear and horizontal-scale;
\item every point $c\in \Aline ow$ is central;
\item the points $o$ and $w$ are central.
\end{enumerate}
\end{theorem}

\begin{proof} The equivalence $(1)\Leftrightarrow(2)$ follows from Theorems~\ref{t:cart-group<=>}, \ref{t:homocentral<=>}, and Proposition~\ref{p:tring-homocentral<=>}.
\vskip3pt

The equivalence $(1)\Leftrightarrow(3)$ follows from Theorems~\ref{t:horizontal-scale}, \ref{t:vshear<=>vshear} and Proposition~\ref{p:ltring-ldist<=>}.
\smallskip

$(2)\Ra(4)$ Assume that the based affine plane $(\Pi,uow)$ is vertical-translation and homocentral. Given any points $c\in\Aline ow$, $x\in\Pi\setminus\{c\}$ and $y\in \Aline cx\setminus\{c\}$, we should find a homothety $H\in\Dil_w(\Pi)$ such that $H(x)=y$. Since $(\Pi,uow)$ is vertical-translation, there exists a translation $T:\Pi\to\Pi$ such that $T(o)=c$. Consider the points $x'\defeq T^{-1}(x)$ and $y'\defeq T^{-1}(y)$. It follows  $o=T^{-1}(c)$, $x\in\Pi\setminus\{c\}$ and $y\in\Aline cx\setminus\{c\}$ that $x'\in \Pi\setminus\{o\}$ and $y'\in\Aline o{x'}\setminus\{o\}$. Since the based affine plane $(\Pi,uow)$ is homocentral, there exists a homothety $H_o\in\Dil_o(\Pi)$ such that $H_o(x')=y'$. Then $H\defeq TH_oT^{-1}:\Pi\to\Pi$ is a dilation of $\Pi$ such that $H(c)=TH_oT^{-1}(w)=TH_o(o)=T(o)=c$ and $H(x)=TH_oT^{-1}(x)=TH_o(x')=T(y')=y$. 
\smallskip

The implication $(4)\Ra(5)$ is obvious.
\smallskip

$(5)\Ra(2)$ Assume that the points $o$ and $w$ are central. Since the point $o$ is central, the based affine plane $(\Pi,uow)$ is homocentral. To prove that it is vertical-translation, take any point $y\in \Aline ow\setminus\{o\}$. We have to construct a translation $T:\Pi\to\Pi$ such that $T(o)=y$. Let $e$ be the diunit of the affine base $uow$ and $z\in\Aline ue$ be a unique point such that $\Aline yz\parallel \Aline ou$. Two cases are possible.
\smallskip

(i) First assume that the lines $\Aline wu$ and $\Aline oz$ have a common point $c$. Since the points $o$ and $w$ are central, there exist homotheties $H_w\in\Dil_w(\Pi)$ and $H_o\in\Dil_o(\Pi)$ such that $H_w(u)=c$ and $H_o(c)=z$. Then $T\defeq H_oH_w$ is a dilation of the plane $\Pi$ such that $T(u)=H_oH_w(u)=H_o(c)=z$ and $T[\Aline ow]=H_oH_w[\Aline ow]=H_o[\Aline ow]=\Aline ow$. Assuming that the dilation $T$ is not a translation, we can find a fixed point $p=T(p)$, which is the centre of the homothety $T$. Since $T(u)=z\ne u$, the centre $p$ belongs to the line $\Aline uz=\Aline ue$. Since $T(o)\in \Aline op\cap \Aline ow=\{o\}$, the dilation $T$ has two fixed points $o$ and $p$ and hence is the identity map of $\Pi$, by Theorem~\ref{t:dilation-atmost1}. On the other hand $T(u)=z\ne u$ witneses that $T$ is not identity. This contradiction shows that $T$ is a translation. Since $T$ is a dilation, $z=T(u)\in T[\Aline ou]\parallel \Aline ou$ and hence $T[\Aline ou]=\Aline zy$ and $T(o)\in T[\Aline ou\cap\Aline ow]=T[\Aline ou]\cap T[\Aline ow]=\Aline zy\cap\Aline ow=\{y\}$.
\smallskip

(ii) Next, assume that the lines $\Aline wu$ and $\Aline oz$ are disjoint and hence parallel. Choose any point $p\in \Aline uz\setminus\{u,z\}$. Since the affine plane $\Pi$ is Proclus and $\Aline wu\cap\Aline oz=\varnothing$, there exist points $c\in \Aline wu\cap\Aline op$ and $d\in \Aline wp\cap\Aline oz$. Since the points $o$ and $w$ are central, there exist homotheties $H_o,H_o'\in\Dil_o(\Pi)$ and $H_w,H_w'\in\Dil_w(\Pi)$ such that $H_w(u)=c$, $H_o(c)=p$, $H_w'(p)=d$, $H_o'(d)=z$. Then $T\defeq H_o'H_w'H_oH_w$ is a dilation of the plane $\Pi$ such that $T(u)=z$ and $T[\Aline ow]=\Aline ow$. Since the lines $\Aline ow=T[\Aline ow]$ and $\Aline uz=T[\Aline uz]$ are disjoint, the dilation $T$ cannot be a homothety and hence $T$ is a translation. Since $z=T(u)\in T[\Aline ou]\parallel\Aline ou$, the lines $\Aline yz$ and $T[\Aline ou]$ coincide and hence $T(o)\in T[\Aline ou\cap\Aline ow]=T[\Aline ou]\cap T[\Aline ow]=\Aline yz\cap\Aline ow=\{y\}$. Therefore $T$ is a required translation with $T(o)=y$.  
\end{proof}

\begin{theorem}\label{t:corps<=>} For a based affine plane $(\Pi,uow)$ and its ternar $\Delta$, the following conditions are equivalent:
\begin{enumerate}
\item the affine plane $\Pi$ is Desarguesian;
\item the affine plane $\Pi$ is dilation;
\item every point of $\Pi$ is central;
\item the points $o,u,w$ are central;
\item the based affine plane $(\Pi,uow)$ is translation and homocentral;
\item the ternar $\Delta$ is linear, distributive and associative\footnote{Linear distributive associative biloops are \defterm{corps}.}.
\end{enumerate}
\end{theorem}

\begin{proof} The equivalence $(1)\Leftrightarrow(2)$ is proved in Theorem~\ref{t:Des<=>dilation}, the equivalence $(2)\Leftrightarrow(3)$ follows from the definition of a dilation plane, the implications $(5)\Ra(3)\Ra (4)$ is obvious, and the equivalence $(5)\Leftrightarrow(6)$ follows from Theorems~\ref{t:VW-Thalesian<=>quasifield}, \ref{t:homocentral<=>}, and Proposition~\ref{p:tring-homocentral<=>}. 
\smallskip

$(4)\Ra(5)$ Assume that the points $u,o,w$ are central. Then the based affine planes $(\Pi,uow)$ and $(\Pi,wou)$ satisfy the condition (3) of Theorem~\ref{t:leftcorps}, which implies that those two based affine planes are homocentral and vertical-translation. Since the based affine plane $(\Pi,wou)$ is vertical-translation, the based affine plane $(\Pi,uow)$ is horizontal-translation. Therefore, $(\Pi,uow)$ is both vertical-translation and horizontal-translation and hence translation.
\end{proof}

\section{Ternars of Thalesian affine planes}

By Theorem~\ref{t:RX-corps}, every affine plane $\Pi$ has a canonical corps of scalars $\IR_X$. In this section we shall prove that for a Thalesian plane $\Pi$ this corps can be identified with the kernel $\Ker(\Delta)$ of the ternar $\Delta$ of any affine base in $\Pi$. 

For a ternar $R$, the set 
$$\Ker(R)\defeq\bigcap_{x,a,b\in R}\{s\in R:s\cdot(x_\times a_+b)=(s\cdot x)_\times a_+(s\cdot b)\}$$is called the \index{ternar!kernel}\index{kernel}\defterm{kernel} of the ternar $R$.
Let also $\Ker^*(R)\defeq\Ker(R)\setminus\{0\}$. 

\begin{proposition}\label{p:Ker-tring} For every ternar $R$, the following statements hold.
\begin{enumerate}
\item $\forall s\in\Ker(R)\;\forall x,a\in R\;\;\big(s\cdot(x\cdot a)=(s\cdot x)\cdot a\big)$.
\item $\forall s\in \Ker(R)\;\forall x,b\in R\;\;\big(s\cdot(x+b)=(s\cdot x)+(s\cdot b)\big)$.
\item $\forall s,t\in\Ker(R)\;\;\big(s\cdot t\in\Ker(R)\big)$.
\item $0,1\in\Ker(R)$.
\end{enumerate}
\end{proposition} 
 
\begin{proof} 1,2. For every $s\in\Ker(R)$ and $x,a,b\in R$ we have the equalities
$$s\cdot(x\cdot a)=s\cdot(x_\times a_+0)=(s\cdot x)_\times a_+(s\cdot 0)=(s\cdot x)_\times a_+0=(s\cdot x)\cdot a$$
and
$$s\cdot(x+b)=s\cdot (x_\times 1_+b)=(s\cdot x)_\times 1_+(s\cdot b)=(s\cdot x)+(s\cdot b).$$
\smallskip

3. For every $s,t\in\Ker(R)$ and every $x,a,b\in R$ we can apply the first statement and conclude that
$$
\begin{aligned}
(s\cdot t)\cdot(x_\times a_+b)&=s\cdot(t\cdot(x_\times a_+b))=s\cdot((t\cdot x)_\times a_+(t\cdot b))=(s\cdot(t\cdot x))_\times a_+(s\cdot (t\cdot b))\\
&=((s\cdot t)\cdot x)_\times a_+((s\cdot t)\cdot b),
\end{aligned}
$$witnessing that $s\cdot t\in\Ker(R)$.
\smallskip

4. Observe that for every $x,a,b\in R$ we have the equalities
$$0\cdot(x_\times a_+b)=0=(0\cdot x)_\times a_+(0\cdot b)\quad\mbox{and}\quad
1\cdot(x_\times a_+b)=x_\times a_+b=(1\cdot x)_\times a_+(1\cdot b),
$$
witnessing that $0,1\in\Ker(R)$.
\end{proof}

\begin{corollary} If $R$ is a linear ternar, then
$$\Ker(R)=\bigcap_{x,y\in R}\{s\in R:s\cdot(x\cdot y)=(s\cdot x)\cdot y\;\wedge\;s\cdot(x+y)=(s\cdot x)+(s\cdot y)\}.$$
\end{corollary}

\begin{proof} We have to prove that $\Ker(R)$ is equal to the set 
$$\ker(R)\defeq \bigcap_{x,y\in R}\{s\in R:s\cdot(x\cdot y)=(s\cdot x)\cdot y\;\wedge\;s\cdot(x+y)=(s\cdot x)+(s\cdot y)\}.$$The inclusion $\Ker(R)\subseteq\ker(R)$ follows from Proposition~\ref{p:Ker-tring}(1,2). To prove that $\ker(R)\subseteq\Ker(R)$, take any elements $s\in \ker(R)$ and $x,a,b\in R$.
The linearity of the ternar $R$ ensures that 
$$s\cdot(x_\times a_+ b)=s\cdot((x\cdot a)+b)=(s\cdot (x\cdot a))+(s\cdot b)=((s\cdot x)\cdot a)+(s\cdot b)=(s\cdot x)_\times a_+(s\cdot b),$$witnessing that $s\in\Ker(R)$. 
\end{proof}

Proposition~\ref{p:Ker-tring} implies $\Ker(R)$ is a submonoid of the unital magma $(R,\cdot)$, and $\Ker^*(R)\defeq\Ker(R)\setminus\{0\}$ is a submonoid of the monoid $(R^*,\cdot)$. Proposition~\ref{p:central-homo} implies the following characterization of the kernel.

\begin{proposition}\label{p:Ker<=>homothety} Let $R$ be a ternar and $R^2$ be its coordinate plane. A  nonzero element $s\in R^*$ belongs to the kernel $\Ker(R)$ of $R$ if and only if the function $H_s:R^2\to R^2$, $H_s:(x,y)\mapsto(s\cdot x,s\cdot y)$, is a homothety of the affine plane $R^2$.
\end{proposition}

\begin{theorem}\label{t:Ker-group} Let $R$ be a ternar $R$ and $o\defeq(0,0)$ be the origin of its coordinate plane $R^2$. The monoid $\Ker^*(R)$ is a group and the function $$H_*:\Ker^*(R)\to \Dil_o(R^2),\quad H_*:s\mapsto H_s,$$ is a group isomorphism.
\end{theorem}

\begin{proof} Propositions~\ref{p:Ker<=>homothety} and \ref{p:central-homo} imply that $H_*:\Ker^*(R)\to \Dil_{o}(R^2)$, $H_*:s\mapsto H_s$, is a well-defined bijective monoid homomorphism.
%
 Since the monoid $\Dil_o(R^2)$ is a group, so is the monoid $\Ker^*(R)$.
\end{proof}

\begin{lemma}\label{l:Ker-submonoid+} For every translation ternar $R$, its kernel $\Ker(R)$ is a submonoid of the commutative group $(R,+)$.
\end{lemma}

\begin{proof} By Theorem~\ref{t:VW-Thalesian<=>quasifield}, the coordinate plane $R^2$ of the translation ternar $R$ is translation. By Corollary~\ref{c:Trans-commutative}, the translation group $\Trans(R^2)$ of the affine plane $R^2$ is commutative. By Corollary~\ref{c:Trans(R)=R2}, the function $\Trans(R^2)\to R^2$, $F\mapsto F(0,0)$, is an isomorphism of the groups $\Trans(X)$ and $(R^2,+)$. The commutativity of the group $\Trans(R^2)$ implies the commutativity of the groups $(R^2,+)$ and $(R,+)$. We claim that the kernel $\Ker(R)$ of $R$ is a submonoid of the commutative group $(R,+)$. By Proposition~\ref{p:Ker-tring}, the set $\Ker(R)$ contains the identity $0$ of the group $(R,+)$. 

To prove that $\Ker(R)$ is a submonoid of the group $(R,+)$, take any elements  $s,t\in\Ker(R)$. For every elements $x,a,b\in R$, by the right-distributivity and linearity of $R$, and the commutativity of the group $(R,+)$, we obtain
$$
\begin{aligned}
(s+t)\cdot (x_\times a_+b)&=s\cdot (x_\times a_+b)+t\cdot(x_\times a_+b)=(s\cdot x)_\times a_+(s\cdot b)+(t\cdot x)_\times a_+(t\cdot b)\\
&=(s\cdot x)\cdot a+s\cdot b+(t\cdot x)\cdot a+t\cdot b=(s\cdot x)\cdot a+(t\cdot x)\cdot a+s\cdot b+t\cdot b\\
&=\big((s\cdot x)+(t\cdot x)\big)\cdot a+(s+t)\cdot b=\big((s+t)\cdot x\big)\cdot a+(s+t)\cdot b\\
&=((s+t)\cdot x)_\times a_+((s+t)\cdot b),
\end{aligned}
$$
witnessing that $s+t\in\Ker(R)$. 
\end{proof}

Theorem~\ref{t:Ker-group} and Lemma~\ref{l:Ker-submonoid+} ensure that the kernel $\Ker(R)$ of a translation ternar $R$ is closed with respect to the operations of addition and multiplication. 

\begin{theorem}\label{t:Ker(R)=RPi} Let $(\Pi,uow)$ be a based affine plane, $e$ be its diunit, and $\Delta$ be its ternar. If the affine plane $\Pi$ is translation, then the set $\Ker(\Delta)$ endowed with the operations of addition and multiplication is a corps. Moreover, the function $F:\Ker(\Delta)\to \IR_\Pi$, $F:s\mapsto\overvector{ose}$, is a well-defined  isomorphism between the corps $\Ker(\Delta)$ and $\IR_\Pi$.
\end{theorem}

\begin{proof} Assume that the affine plane $\Pi$ is Thalesian. By Proposition~\ref{p:Ker<=>homothety}, for every element $s\in\Ker^*(\Delta)$, the map $\ddot H_s:\Delta^2\to \Delta^2$, $\ddot H_s:(x,y)\mapsto(s\cdot x,s\cdot y)$, is a homothety of the affine plane $\Delta^2$. Moreover, by Theorem~\ref{t:Ker-group}, the function $\Ker^*(\Delta)\to\Dil_{00}(\Delta^2)$, $s\mapsto \ddot H_s$, is an isomorphism of the groups $\Ker^*(\Delta)$ and $\Dil_{00}(\Delta^2)$.

Let $C:\Pi\to\Delta^2$ be the coordinate chart of the based affine plane $\Pi$. Since $C$ is an isomorphism of the based affine planes $\Pi$ and $\Delta^2$, the function $I:\Ker^*(\Delta)\to\Dil_o(\Pi)$, $I:s\mapsto H_s\defeq C^{-1}\ddot H_sC$, is an isomorphism of the groups $\Ker^*(\Delta)$ and $\Dil_{o}(\Pi)$.

For every non-zero scalar $\sigma\in\IR^*_\Pi$, consider the map $ H_\sigma:\Pi\to\Pi$ assigning to every point $x\in\Pi$ the unique point $y\in\Aline ox$ such that $oyx\in\{ooo\}\cup \sigma$. By Theorem~\ref{t:RX=Dilo}, the map $J:\IR^*_\Pi\to\Dil_o(\Pi)$, $J:\sigma\mapsto H_\sigma$, is an isomorphism of the groups $\IR^*_\Pi$ and $\Dil_o(\Pi)$. Then $J^{-1}I:\Ker^*(\Delta)\to\IR^*_\Pi$ is a group isomorphism.  We claim that $I(s)=JF(s)$ for every $s\in\Ker^*(\Delta)$.

Given any nonzero element $s\in\Ker^*(\Delta)$, consider the non-zero scalar $\sigma=\overvector{ose}=F(s)$ and the homotheties $H_s$ and $H_\sigma$ of the affine plane $\Pi$.  Observe that 
$$H_s(e)=C^{-1}\ddot H_sC(e)=C^{-1}\ddot H_s(e,e)=C^{-1}(s\cdot e,s\cdot e)=C^{-1}(s,s)=s=H_\sigma(e)$$ and 
$$H_s(o)=C^{-1}\ddot H_sC(o)=C^{-1}\ddot H_s(o,o)=C^{-1}(s\cdot o,s\cdot o)=C^{-1}(o,o)=o=H_\sigma(o).$$ Since $H_s$ and $H_\sigma $ are two dilations with $H_soe=os=H_\sigma oe$, Theorem~\ref{t:dilation-atmost1} ensures that $H_s=H_\sigma$ and hence
$I(s)=H_s=H_\sigma =J(\sigma)=JF(s)$ and $F(s)=J^{-1}I(s)$.
Now we see that the restriction $F{\restriction}_{\Ker^*(\Delta)}=J^{-1}I$ is an isomorphism of the groups $\Ker^{*}(\Delta)$ and $\IR^*_\Pi$. Since $F(o)=\overvector{ooe}=0$, the function $F:\Ker(\Delta)\to\IR_\Pi$ is an isomorphism of the monoids $\Ker(\Delta)$ and $\IR_\Pi$.

It remains to prove that the map $F$ preserves the operation of addition. Given any points $a,b\in\Ker(\Delta)\subseteq \Delta$, consider their sum $c=a+b\in R$. Lemma~\ref{l:Ker-submonoid+} ensures that $c=a+b\in\Ker(R)$. The equality $c=a+b$ implies the equality $cc=(c,c)=(a+b,a+b)=(a,a)+(b,b)=aa+bb$ in the group $(\Delta^2,+)$. By Proposition~\ref{p:trans=s+t}, for every pair $(s,t)\in \Delta^2$, the function $T_{st}:\Delta^2\to \Delta^2$, $T_{st}:(x,y)\mapsto (x+s,y+t)$, is a translation of the affine plane $\Delta^2$. Moreover, by Theorem~\ref{c:Trans(R)=R2}, the function $\Delta^2\to\Trans(\Delta^2)$, $st\mapsto T_{st}$, is an isomorphism of the groups  $(\Delta^2,+)$ and $\Trans(\Delta^2)$. Then the equality $cc=aa+bb$ implies $T_{cc}=T_{aa}\circ T_{bb}$. It follows that for every $s\in\Delta$, $T_s\defeq C^{-1}T_{ss}C$ is a unique translation of the affine plane $\Pi$ such that $T_s(o)=s$. Observe that 
$$T_c=C^{-1}T_{cc}C=C^{-1}T_{aa}T_{bb}C=C^{-1}T_{aa}CC^{-1}T_{bb}C=T_aT_b.$$By Theorems~\ref{t:Trans(X)=vecX} and \ref{t:vector-addition}, the group $\overvector \Pi=\Trans(\Pi)$ is commutative and hence 
$$\overvector{oc}=T_c=T_a\circ T_b=T_b\circ T_a=\overvector{oa}+\overvector{ob}.$$By Theorem~\ref{t:scalar-addition}, the equality $\overvector{oa}+\overvector{ob}=\overvector{oc}$ implies $\overvector{oce}=\overvector{oae}+\overvector{obe}$ and hence
$$F(a+b)=F(c)=\overvector{oce}=\overvector{oae}+\overvector{obe}=F(a)+F(b),$$witnessing that the bijection $F:\Ker(R)\to \IR_\Pi$ preserves the operation of addition. 

Therefore, the bijective function $F:\Ker(\Delta)\to\IR_\Pi$ preserves the operations of addition and multiplications. Taking into account that $(\IR_\Pi,+,\cdot)$ is a corps, we conclude that its isomorphic copy $(\Ker(\Delta),+,\cdot)$ is a corps and the map $F$ is an isomorphism of the corps $\Ker(\Delta)$ and $\IR_\Pi$.
\end{proof}

\begin{corollary}\label{c:R=RPi} Let $(\Pi,uow)$ be a based affine plane, $e$ be its diunit, and $\Delta$ be its ternar. If the affine plane $\Pi$ is Desarguesian, then the ternar $\Delta$ is a corps. Moreover, the function $F:\Delta\to \IR_\Pi$, $F:s\mapsto\overvector{ose}$, is a well-defined  isomorphism between the corps $\Delta$ and $\IR_\Pi$.
\end{corollary}

\begin{proof} By Theorem~\ref{t:corps<=>}, the ternar $\Delta$ is a corps and hence $\Ker(\Delta)=\Delta$. By Theorem~\ref{t:Ker(R)=RPi}, the function $F:\Delta\to \IR_\Pi$, $F:s\mapsto\overvector{ose}$, is a well-defined  isomorphism between the corps $\Delta$ and $\IR_\Pi$.
\end{proof}

\begin{corollary}\label{c:coordinates=vectors} Let $(\Pi,uow)$ be a Thalesian based affine plane, $e$ be its diunit, and $\Delta=\Aline oe$ be its ternar. For any point $p\in\Pi$ with coordinates $xy\in \Ker(\Delta)^2\subseteq \Delta^2$,  the equality $\overvector{op}=\overvector{oxe}\cdot\overvector{ou}+\overvector{oye}\cdot\overvector{ow}$ holds.
\end{corollary}

\begin{proof} Let $a\in\Aline ou$ and $b\in\Aline ow$ be the points with coordinates $xo$ and $oy$, respectively. Since $x,y\in\Ker(\Delta)$, the portions $\overvector{oxe}$ and $\overvector{oye}$ are scalars, by Theorem~\ref{t:Ker(R)=RPi}. It follows from $\Aline {x}{a}\subparallel \Aline ue$ and $\Aline{b}{y}\subparallel \Aline ew$ that $\overvector{oxe}=\overvector{oau}$ and $\overvector{oye}=\overvector{obw}$. Applying Corollary~\ref{c:Thalesian<=>vector=funvector} and Theorems~\ref{t:vector-addition}, \ref{t:scalar-by-vector}, we conclude that $\overvector{op}=\overvector{oa}+\overvector{ob}=\overvector{oau}\cdot\overvector{ou}+\overvector{obw}\cdot\overvector{ow}=\overvector{oxe}\cdot\overvector{ou}+\overvector{oye}\cdot\overvector{ow}$. 

\begin{picture}(100,120)(-150,-15)
\linethickness{0.6pt}
\put(0,0){\color{teal}\line(1,0){100}}
\put(0,0){\color{cyan}\line(0,1){100}}
\put(0,0){\line(1,1){90}}
\put(50,0){\color{cyan}\line(0,1){80}}
\put(0,80){\color{teal}\line(1,0){80}}
\put(30,0){\color{cyan}\line(0,1){30}}
\put(0,30){\color{teal}\line(1,0){30}}
{\linethickness{1pt}
\put(0,0){\color{red}\vector(5,8){50}}
\put(0,0){\color{teal}\vector(1,0){50}}
\put(0,0){\color{cyan}\vector(0,1){80}}
}

\put(0,0){\circle*{3}}
\put(-7,-7){$o$}
\put(30,0){\circle*{3}}
\put(27,-8){$u$}
\put(0,30){\circle*{3}}
\put(-10,28){$w$}
\put(30,30){\circle*{3}}
\put(32,26){$e$}
\put(50,0){\circle*{3}}
\put(52,3){$a$}
\put(0,80){\circle*{3}}
\put(-8,78){$b$}
\put(50,80){\circle*{3}}
\put(48,83){$p$}
\put(50,50){\circle*{3}}
\put(52,46){$x$}
\put(80,80){\circle*{3}}
\put(82,76){$y$}
\end{picture}
\end{proof}

\begin{theorem}\label{t:Pappian<=>tring-field} For a based affine plane $(\Pi,uow)$ and its ternar $\Delta$, the following conditions are equivalent:
\begin{enumerate}
\item the affine plane $\Pi$ is Pappian;
\item the ternar $\Delta$ is linear and its biloop is a field.
\end{enumerate}
\end{theorem}

\begin{proof} $(1)\Ra(2)$ Assume that the affine plane $\Pi$ is Pappian. 
By Hessenberg's Theorem~\ref{t:Hessenberg-affine}, the Pappian affine plane $\Pi$ is Desarguesian, by Theorem~\ref{t:ADA=>AMA}, the Desarguesian affine plane $\Pi$ is Thalesian, and by Theorem~\ref{t:VW-Thalesian<=>quasifield}, the ternar $\Delta$ is translation. By Theorem~\ref{t:Ker(R)=RPi}, the map $\Ker(\Delta)\to \IR_\Pi$, $s\mapsto \overvector{ose}$, is an isomorphism of the corps $\Ker(\Delta)$ and $\IR_\Pi$. Since $\Pi$ is Desarguesian, for every point $s\in\Delta$, the portion $\overvector{ose}$ is a scalar and hence $s\in\Ker(\Delta)$. Therefore, $\Delta=\Ker(\Delta)$ is a corps. Since $\Pi$ is Pappian, by Theorem~\ref{t:Papp<=>Des+RX}, the corps $\IR_\Pi$ is a field and so is its isomorphic copy $\Delta=\Ker(\Delta)$.
\smallskip


$(2)\Ra(1)$ Assume that the ternar $\Delta$ is linear and its biloop is a field. Then $\Ker(\Delta)=\Delta$. By Theorem~\ref{t:VW-Thalesian<=>quasifield}, the affine plane $\Pi$ is Thalesian and the map $\Ker(\Delta)\to\IR_\Pi$, $s\mapsto \overvector{ose}$, is a corps isomorphism. Since $\Ker(\Delta)=\Delta$, for every $s\in\Delta$ the portion $\overvector{ose}$ is a scalar and hence $\IR_\Pi={\dddot{\, \Pi}}_{\!\Join}$ and the affine plane $\Pi$ is Desarguesian, by Theorem~\ref{t:Desargues<=>3Desargues}. Since  the corps $\Delta=\Ker(\Delta)$ is a field, so is its isomorphic copy $\IR_\Pi$. By Theorem~\ref{t:Papp<=>Des+RX}, the affine plane $\Pi$ is Pappian. 
\end{proof}


We shall say that a ternar $R$ is a ternar of an affine space $X$ if $R$ is isomorphic to the ternar $\Delta_{uow}$ of some triangle $uow$ in $X$, which is an affine base for the plane $\overline{uow}$. In this case we shall say that the plane $\overline{uow}$ can be \defterm{coordinatized} by the ternar $R$.

\begin{theorem}\label{t:Desarg=>tring=corps} Every ternar $R$ of a Desarguesian affine space $X$ is isomorphic to the scalar corps $\IR_X$ endowed with the ternary operation $x_\times a_+b\defeq (x\cdot a)+b$.
\end{theorem}

\begin{proof} Given any ternar $R$ of $X$, find a plane $\Pi\subseteq X$ and an affine base $uow$ in $\Pi$ whose ternar $\Delta$ is isomorphic to the ternar $R$. Since $X$ is Desarguesian, the plane $\Pi$ is Desarguesian, too.  By Theorem~\ref{t:corps<=>}, the ternar $\Delta$ is a linear corps and hence $\Ker(\Delta)=\Delta$. By Theorem~\ref{t:Ker(R)=RPi}, the map $F:\Delta\to\IR_\Pi$, $F:x\mapsto\overvector{oxe}$, is an isomorphism of the corps (here $e$ is the diunit of the affine base $uow$).  Since every line triple in $X$ is Desarguesian, the functions $I_\Pi:\Delta\to\IR_\Pi$, $I_\Pi:x\mapsto \overvector{oxe}$, and $I_X:\Delta\to \IR_X$, $I_X:x\mapsto\overvector{oxe}$, are bijective. This implies that the ``identity'' function $I:\IR_\Pi\to\IR_X$, $I:\overvector{oxe}\mapsto\overvector{oxe}$, is an isomorphism of the scalar corps $\IR_\Pi$ and $\IR_X$. Then the function $IF:\Delta\to\IR_X$ also is an isomorphism of the corps.
 Since the ternar $\Delta$ is linear, the function $IF$ is an isomorphism between the ternar $\Delta$ and the corps $\IR_X$, endowed with the ternary operation $x_\times a_+b\defeq (x\cdot a)+b$.
\end{proof}



\section{Algebra versus Geometry in based affine planes}

The results of the preceding sections imply the following equivalences between  geometric properties of based affine plane and algebraic properties of its ternar.

\begin{corollary}\label{c:Algebra-vs-Geometry-in-trings:affine} For a based affine plane $\Pi$, its ternar is
\begin{enumerate}
\item linear and associative-plus iff $\Pi$ is vertical-translation;
\item linear, right-distributuve, and associative-plus iff $\Pi$ vertical-translation and horizontal-translation;
\item linear, left-distributive, and associative-plus iff $\Pi$ is vertical-translation and vertical-shear;
\item linear, distributive, and associative-plus iff $\Pi$ is vertical-translation, horizontal-translation, and vertical-shear;
\smallskip
\item linear and associative-dot iff $\Pi$ is horizontal-scale;
\item linear, right-distributive, and associative-dot iff $\Pi$ is horizontal-scale and vertical-scale;
\item linear, associative-dot, and diagonal-scale iff $\Pi$ is horizontal-scale and diagonal-scale;
\item linear, left-distributive, and associative-dot iff $\Pi$ is horizontal-scale and homocentral;
\item linear, distributive, and associative-dot iff $\Pi$ is horizontal-scale, vertical-scale and homocentral;
\smallskip
\item linear and associative iff $\Pi$ is vertical-translation and horizontal-scale;
\item linear, left-distributive, and associative iff $\Pi$ is vertical-shear and horizontal-scale iff $\Pi$ is vertical-translation and homocentral;
\item linear, right-distributive, and associative  iff $\Pi$ is vertical-translation and vertical-scale iff $\Pi$ is translation, vertical-scale, and horizontal-scale;
\item linear, distributive, and associative iff $\Pi$ Desarguesian iff $\Pi$ is dilation iff $\Pi$ is translation and homocentral iff $\Pi$ is vertical-translation, vertical-shear and vertical-scale iff $\Pi$ is horizontal-scale, vertical-scale, and vertical-shear iff $\Pi$ is vertical-shear, horizontal-scale and horizontal-translation iff $\Pi$ is shear and vertical-scale iff $\Pi$ is shear and horizontal-scale;
\smallskip
\item linear, distributive, associative-plus, and alternative-dot iff $\Pi$ is shear iff $\Pi$ is translation, vertical-shear and horizontal-shear;
\smallskip
\item linear, distributive, associative, and commutative iff $\Pi$ is Pappian.
\smallskip
\end{enumerate}
\end{corollary}

Therefore, for every linear based affine plane $\Pi$ we have the following correspondence between geometric properties of the based affine plane $\Pi$ and algebraic properties of the biloop of its linear ternar. 
\newpage

{\small
\begin{center}
\begin{tabular}{|c|c|}
\hline
{\bf \!Geometric property of a based affine plane\!}&{\bf Algebraic property of its ternar}\\
\hline
linear&linear\\
\hline
\hline
vertical-translation&linear associative-plus\\
\hline translation $\Leftrightarrow$&linear right-distributive associative-plus\\
vertical-translation and horizontal-translation&\\
\hline
 vertical-translation and vertical-shear&linear left-distributive associative-plus\\
\hline
translation vertical-shear&linear distributive associative-plus\\
\hline
\hline
 horizontal-scale&linear associative-dot\\
\hline
\end{tabular}
\end{center}
}

\bigskip
The following table collects the algebraic identities defining various properties of ternars.
\smallskip

\begin{center}
\begin{tabular}{|c|c|}
\hline
{\bf Algebraic property of a ternar}&{\bf The algebraic identity}\\
\hline
linear&$x_\times a_+b=(x\cdot a)+b$\\
\hline
vertical-translation&$x_\times a_+(b+s)=(x_\times a_+b)+s$\\
\hline
horizontal-translation & $(x+s)_\times a_+b=x_\times a_+(s_\times a_+b)$\\
\hline
vertical-shear&$x_\times c_+(x_\times a_+b)=x_\times (1_\times c_+a)_+b$\\
\hline
vertical-scale& $(x_\times a_+b)\cdot c=x_\times(a\cdot c)_+(b\cdot c)$\\
\hline
horizontal-scale& $x_\times (c\cdot a)_+b=(x\cdot c)_\times a_+b$ \\
\hline
diagonal-scale& $c+x{\cdot}c\mp (x{\cdot}(a\circ c)+b)=(c\circ x){\cdot}a+c{\cdot}a\mp(c+b)$ \\
\hline
homocentral& $c{\cdot}(x_\times a_+b)=(c{\cdot}x)_\times a_+(c{\cdot}b)$\\
\hline
commutative-plus&$a+b=b+a$\\
\hline
commutative-dot&$a\cdot b=b\cdot a$\\
\hline
commutative&$a+b=b+a\;\wedge\;a\cdot b=b\cdot a$\\
\hline
associative-plus&$(a+b)+c=a+(b+c)$\\
\hline
associative-dot&$(a\cdot b)\cdot c=a\cdot(b\cdot c)$\\
\hline
associative&$(a{+}b){+}c=a{+}(b{+}c)\;\wedge\;(a{\cdot} b){\cdot} c=a{\cdot}(b{\cdot} c)$\\
\hline
alternative-dot&$(a{\cdot} a){\cdot} b=a{\cdot}(a{\cdot} b)\;\wedge\;(a{\cdot} b){\cdot} b=a{\cdot}(b{\cdot} a)$\\
\hline
left-distributive&$a\cdot (x+y)=a{\cdot}x+a{\cdot}y$\\
\hline
right-distributive&$(x+y)\cdot b=x{\cdot}b+y\cdot b$\\
\hline
distributive&$a{\cdot}(x+y)=a{\cdot}x+a{\cdot}y\;\wedge\;(x+y){\cdot} b=x{\cdot}b+y{\cdot}b$\\
\hline
\end{tabular}
\end{center}

\chapter{Various classification schemes for affine planes}

In this section we present several classifications of affine planes, describing their homogeneity properties and also by algebraic properties of their coordinate ternars. 

We shall say that an affine space $X$ is \index{coordinatized affine space}\defterm{coordinatized} by a ternar (or biloop) $R$ if $X$ contains a plane $\Pi\subseteq X$ that has an affine base $uow\in \Pi^3$ whose ternar (or its biloop) is isomorphic to $R$. In this case we shall also say that the ternar $R$ is \index{affine space!ternar of}\defterm{a ternar} of the affine space $X$. If $X$ is an affine space of rank $\|X\|\ge 4$, then it is Desarguesian and every ternar of $X$ is isomorphic to the scalar corps $\IR_X$ of $X$. 

By Theorem~\ref{t:Pappian<=>tring-field}, an affine plane is Pappian if and only if it is coordinatized by a field. For an affine plane $\Pi$, we denote by $\mathcal L_\Pi$ the family of all lines in $\Pi$, and by $\partial \Pi$ the boundary of $\Pi$ (i.e., the set of all spreads of parallel lines in $\Pi$). Elements of the boundary $\partial\Pi$ will be called \index{direction}\defterm{directions} on the affine plane $\Pi$.

For labeling various types of affine planes we shall extensively use ranks of sets in affine planes, their families of lines and their boundaries. 
\smallskip

For a subset $A$ of an affine plane $\Pi$, the number
\smallskip

$\|A\|\defeq\begin{cases}
0&\mbox{if $A=\varnothing$};\\
1&\mbox{if $A=\{p\}$ for some point $p\in\Pi$};\\
2&\mbox{if $A\subseteq L$ for a unique line $L\in\mathcal L_\Pi$};\\
3&\mbox{otherwise}
\end{cases}
$ 
\smallskip

\noindent is called the \index{rank of a set of points}\defterm{rank of the set of points} $A$. Observe that $\|A\|$ coincides with the usual rank of the set $A$ in the plane $\Pi$.
\smallskip

For a family of lines $\mathcal B$ in $\Pi$, the number
\smallskip

$\|\mathcal B\|\defeq\begin{cases}
0&\mbox{if $\mathcal B=\varnothing$};\\
1&\mbox{if $\mathcal B=\{L\}$ for some line $L\in\mathcal L_\Pi$};\\
2&\mbox{if $\mathcal B$ is paraconcurrent and $|\mathcal B|>1$;}\\
3&\mbox{otherwise}
\end{cases}
$ 

\noindent is called the \index{rank of a family of lines}\defterm{rank of the family of lines} $\mathcal B$. Let us recall that a family of line $\mathcal B$ is called \defterm{paraconcurrent} if either $\bigcap\mathcal B\ne\varnothing$ or any two lines in $\mathcal B$ are parallel.
\smallskip

For a family of directions $\mathfrak D\subseteq\partial\Pi$, the number
\smallskip

$\|\mathfrak D\|\defeq\begin{cases}
0&\mbox{if $\mathfrak D=\varnothing$};\\
1&\mbox{if $\mathfrak D=\{\delta\}$ for some direction $\delta\in\partial\Pi$};\\
2&\mbox{otherwise}
\end{cases}
$

\noindent is called the \index{rank of a set of directions}\defterm{rank of the family of directions} $\mathfrak D$.

\section{The translation classification of affine planes} 

\begin{definition} For an affine plane $\Pi$, a direction $\delta\in\partial\Pi$ is called  a \index{translation direction}\defterm{translation direction} of $\Pi$ if for every points $x,y\in \Pi$ with $\Aline xy\in \delta$, there exists a translation $T:\Pi\to\Pi$ such that $T(x)=y$. The set $\TT_\Pi\subseteq\partial \Pi$ of all translation directions of $\Pi$ is called is the \index{trace!translation}\index{translation trace}\defterm{translation trace} of the affine plane $\Pi$, and its rank $\|\TT_\Pi\|\in\{0,1,2\}$ is called the \index{translation rank}\defterm{translation rank} of the affine plane $\Pi$. 
\end{definition}

\begin{example} A based affine plane $(\Pi,uow)$ is
\begin{itemize}
\item vertical-translation if and only if $(\Aline ow)_\parallel\in \TT_\Pi$;
\item horizontal-translation if and only if $(\Aline ou)_\parallel\in \TT_\Pi$;
\item translation if and only if $\TT_\Pi=\partial\Pi$.
\end{itemize}
\end{example}

\begin{theorem}\label{t:trans-clas} For every affine plane $\Pi$ one of the following $3$ cases holds:
\begin{enumerate}
\item[$(0)$] $\|\TT_\Pi\|=0$ and $\TT_\Pi=\varnothing$;
\item[$(1)$] $\|\TT_\Pi\|=1$ and $\TT_\Pi=\{\delta\}$ for some direction $\delta\in\partial\Pi$;
\item[$(2)$] $\|\TT_\Pi\|=2$ and $\TT_\Pi=\partial \Pi$.
\end{enumerate}
\end{theorem}

\begin{proof} It suffices to prove that $|\TT_\Pi|\ge 2$ implies $\TT_\Pi=\partial\Pi$.  Assuming that $|\TT_\Pi|\ge 2$, fix two distinct translation directions $h,v\in\partial\Pi$. Choose an affine base $uow$ in $\Pi$ such that $\Aline ou\in h$ and $\Aline ow\in v$. Since $h$ and $v$ are translation directions, the based affine plane $(\Pi,uow)$ is horizontal-translation and vertical-translation. By Proposition~\ref{p:trans<=>hv-trans}, the plane $\Pi$ is translation and hence $\TT_\Pi=\partial\Pi$.
\end{proof}

\begin{theorem}\label{t:trans-coordinate} An affine plane $\Pi$ is coordinatized by 
\begin{enumerate}
\item a linear associative-plus ternar\footnote{Linear  associative-plus biloops are called {\em cartesian groups}.} iff $\|\TT_\Pi\|\ge 1$;
\item a linear right-distributive associative-plus ternar\footnote{Linear right-distributive associative-plus biloops are called {\em quasi-fields}.} iff $\|\TT_\Pi\|=2$.
\end{enumerate}
\end{theorem}
 
\begin{proof} 1. If the affine plane $\Pi$ is coordinatized by a linear associative-plus ternar, then $\Pi$ has an affine base $uow$ whose ternar $\Delta$ is linear and associative-plus. By Theorem~\ref{t:cart-group<=>}, the based affine plane $(\Pi,uow)$ is vertical-translation, which means that $(\Aline ow)_\parallel\in \TT_\Pi$ and hence $\|\TT_\Pi\|\ge 1$.

If $\|\TT_\Pi\|\ge 1$, then the affine plane $\Pi$ has a translation direction $v\in\partial \Pi$. Choose any affine base $uow$ in $\Pi$ such that $\Aline ow\in v$. Then the based affine plane $(\Pi,uow)$ is vertical-translation and by Theorem~\ref{t:cart-group<=>} its ternar is linear and associative-plus, witnessing that the affine plane $\Pi$ is coordinatized by a linear associative-plus ternar.
\smallskip

2. If  the affine plane $\Pi$ is coordinatized by a linear right-distributive associative-plus ternar, then $\Pi$ has an affine base $uow$ whose ternar $\Delta$ is linear, right-distributive, and associative-plus. By Theorem~\ref{t:VW-Thalesian<=>quasifield}, the based affine plane $(\Pi,uow)$ is translation, which means that $\TT_\Pi=\partial\Pi$ and hence $\|\TT_\Pi\|=2$.

If $\|\TT_\Pi\|\ge 2$, then $\TT_\Pi=\partial\Pi$, by Theorem~\ref{t:trans-clas}, and hence the affine plane $\Pi$ is translation. Fix any affine base $uow$ in $\Pi$. By Theorem~\ref{t:VW-Thalesian<=>quasifield}, the ternar $\Delta$ of the based affine plane $(\Pi,uow)$ is linear, right-distributive and associative-plus, witnessing that the affine plane $\Pi$ is coordinatized by a linear right-distributive associative-plus ternar.
\end{proof}

Theorem~\ref{t:aa'-Thalesian<=>} implies the following characterization of translation directions, which can be helpful for computational purposes.

\begin{proposition} For an affine plane $\Pi$, a direction $\delta\in\partial\Pi$ is  translation if and only if for every distinct lines $A,B,C\in \delta$ and distinct points $a,a'\in A$, $b,b'\in B$, $c,c'\in C$ with $\Aline ab\cap\Aline {a'}{b'}=\varnothing=\Aline bc\cap\Aline{b'}{c'}$, we have $\Aline ac\cap\Aline{a'}{c'}=\varnothing$.
\end{proposition}

\begin{exercise} Show that the algorithmic complexity of finding the translation trace of an affine plane of order $n$ is at most $O(n^4)$.
\smallskip

{\em Hint:} Apply Theorem~\ref{t:cart-group<=>}.
\end{exercise}

\section{The shear classification of affine planes} 

\begin{definition} For an affine plane $\Pi$, a line $L\subseteq \Pi$ is called a \index{shear line}\defterm{shear line} of $\Pi$ if for every points $x,y\in \Pi$ with $\Aline xy\cap L=\varnothing$, there exists an automorphism $A:\Pi\to\Pi$ such that $A(x)=y$ and $L\subseteq \Fix(A)$. The set $\NT_\Pi$ of all shear lines of $\Pi$ is called is the \index{shear trace}\index{trace!shear}\defterm{shear trace} of the affine plane $\Pi$, and its rank $\|\NT_\Pi\|\in\{0,1,2,3\}$ is called the \index{shear rank}\defterm{shear rank} of $\Pi$. 
\end{definition}

\begin{example} A based affine plane $(\Pi,uow)$ is
\begin{itemize}
\item vertical-shear if and only if $\Aline ow\in\NT_\Pi$;
\item horizontal-shear if and only if $\Aline ou\in\NT_\Pi$;
\item shear if and only if $\NT_\Pi$ coincides with the family $\mathcal L_\Pi$ of all lines of $\Pi$.\end{itemize}
\end{example}

For a point $p$ of an affine plane $\Pi$, let $\mathcal L_p\defeq\{L\in\mathcal L_\Pi:p\in L\}$ be the family of all lines in $\Pi$ that contain the point $p$ (the family $\mathcal L_p$ is called the \index{pencil of lines}\defterm{pencil of lines} passing through the point $p$).

\begin{theorem}\label{t:shear-class} For every affine plane $\Pi$ one of the following $5$ cases holds:
\begin{enumerate}
\item[$(0)$] $\|\NT_\Pi\|=0$ and $\NT_\Pi=\varnothing$;
\item[$(1)$] $\|\NT_\Pi\|=1$ and $\NT_\Pi=\{L\}$ for some line $L\in\mathcal L_\Pi$;
\item[$(\dot 2)$] $\|\NT_\Pi\|=2$ and $\NT_\Pi=\mathcal L_p$ for some point $p\in\Pi$;
\item[$(\bar 2)$] $\|\NT_\Pi\|=2$ and $\NT_\Pi=\delta$ for some direction $\delta\in\partial\Pi$;
\item[$(3)$] $\|\NT_\Pi\|=3$ and $\NT_\Pi=\mathcal L_\Pi$.
\end{enumerate} 
\end{theorem}

\begin{proof} If the set $\NT_\Pi$ is empty, then the case (0) holds. If the set $\NT_\Pi$ is a singleton, then the case (1) holds. So, assume that the set $\NT_\Pi$ contains two distinct shear lines $L'$ and $L''$. Two cases are possible.

(i) If the shear lines $L',L''$ have a common point $o\in\Pi$, then every line in the family $\mathcal L_o$ is a shear line, by Lemma~\ref{l:Lshear}. In this case $\mathcal L_o\subseteq \NT_\Pi$. If $\mathcal L_o=\NT_\Pi$, then the case $(\dot 2)$ holds. 

If $\mathcal L_o\ne\NT_\Pi$, then there exists a shear line $\Lambda\in\mathcal L_\Pi\setminus\mathcal L_o$. Choose a line $L_o\in\mathcal L_o$ such that $L_o\parallel \Lambda$. Take any points $w\in L_o\setminus\{o\}$ and $u\in \Lambda$. Then $(\Pi,uow)$ is a vertical-shear based affine plane, satisfying the condition (5) of Theorem~\ref{t:vtrans+vshear=vshears}, which implies that the based affine plane $(\Pi,uow)$ is vertical-translation. Then there exists a translation $T:\Pi\to\Pi$ such that $T(o)=w$. Since $T$ is an automorphism of the affine plane $\Pi$, the image $T[\Aline ou]=\Aline we$ of the shear line $\Aline ou\in \mathcal L_o\subseteq\NT_\Pi$ is a shear line in $\Pi$. Then the based affine plane $(\Pi,wou)$ satisfies the condition (5) of Theorem~\ref{t:vtrans+vshear=vshears}, which implies that the based affine plane $(\Pi,wou)$ is vertical-translation and hence the based affine plane $(\Pi,uow)$ is horizontal-translation. By Proposition~\ref{p:trans<=>hv-trans}, the horizontal-translation vertical-translation based affine plane $(\Pi,uow)$ is translation. In this case $\mathcal L_o\subseteq \NT_\Pi$ implies $\bigcup_{p\in\Pi}\mathcal L_p=\mathcal L_\Pi=\NT_\Pi$. Therefore, the case (3) holds.

(ii) Next, assume that the shear lines $L'$ and $L''$ are disjoint. Choose any distinct points $o,w\in L'$ and any point $u\in L''$. Then the triple $uow$ is an affine base for the affine plane $\Pi$. Since $\Aline ow=L'\in\NT_\Pi$ and $\Aline ue=L''\in\NT_\Pi$, the based affine space $(\Pi,uow)$ is vertical-shear and satisfies the condition (5) of Theorem~\ref{t:vtrans+vshear=vshears}, which implies that the based affine plane $(\Pi,uow)$ is vertical-translation and $(\Aline ow)_\parallel \subseteq \NT_\Pi$. If $(\Aline ow)_\parallel=\NT_\Pi$, then the case $(\bar 2)$ holds. If $(\Aline ow)_\parallel \ne \NT_\Pi$, then the family $\NT_\Pi$ contains two concurrent lines and we have the case (i) implying that the case $(\dot 2)$ or $(3)$ holds.
\end{proof}

The following characterization of shear lines can be proved by analogy with Theorem~\ref{t:aa'-Thalesian<=>}. It can be helpful for computation of the shear trace of an affine plane.

\begin{proposition} A line $L$ in an affine plane $\Pi$ is a shear line if and only if for every distinct lines $A,B,C\in L_\parallel\setminus\{L\}$ and distinct points $a,a'\in A$, $b,b'\in B$, $c,c'\in C$ with $\varnothing\ne\Aline ab\cap\Aline {a'}{b'}\subseteq L$ and $\varnothing\ne\Aline bc\cap\Aline{b'}{c'}\subseteq L$, we have $\Aline ac\cap\Aline{a'}{c'}\subseteq L$.
\end{proposition}

\begin{exercise} Show that the algorithmic complexity of finding the shear trace of an affine plane of order $n$ is at most $O(n^6)$.
\smallskip

{\em Hint:} Apply Theorem~\ref{t:vshear<=>vshear}.
\end{exercise}

\section{The central classification of affine planes} 

\begin{definition} A point $o$ of an plane $\Pi$ is called \index{central point}\defterm{central} if for every points $x,y\in \Pi\setminus\{o\}$ with $o\in\Aline xy$, there exists a central automorphism $C:\Pi\to\Pi$ such that $C(x)=y$ and $C(o)=o$. The set $\CT_\Pi$ of all central points of $\Pi$ is called is the \index{central trace}\index{trace!central}\defterm{central trace} of the affine plane $\Pi$ and its rank $\|\CT_\Pi\|\in\{0,1,2,3\}$ is called the \index{central rank}\defterm{central rank} of the affine plane $\Pi$.
\end{definition}

\begin{remark} By Theorem~\ref{t:central<=>homothety} an automorphism of an affine space is central if and only if it is a homothety. So, the central trace also could be called the homothety trace of an affine plane. The notation $\CT_\Pi$ for denoting the central trace is chosen becase the letter $\mathsf X$ has a central point and resembles an action of a central automorphism on an affine plane. 
\end{remark} 

\begin{example} A based affine plane $(\Pi,uow)$ is homocentral if and only if $o\in \CT_\Pi$.
\end{example}

\begin{example} An affine plane $\Pi$ is dilation if and only if $\CT_\Pi=\Pi$.
\end{example}

\begin{theorem}\label{t:central-class} For every affine plane $\Pi$ one of the following $4$ cases holds:
\begin{itemize}
\item[$(0)$] $\|\CT_\Pi\|=0$ and $\CT_\Pi=\varnothing$;
\item[$(1)$] $\|\CT_\Pi\|=1$ and $\CT_\Pi=\{p\}$ for some point $p\in\Pi$;
\item[$(2)$] $\|\CT_\Pi\|=2$ and $\CT_\Pi=L$ for some line $L\in\mathcal L_\Pi$;
\item[$(3)$] $\|\CT_\Pi\|=3$ and $\CT_\Pi=\Pi$.
\end{itemize} 
\end{theorem}

\begin{proof} If $\CT_\Pi=\varnothing$, then the case $(0)$ holds. If $\CT_\Pi$ is a singleton, then the case $(1)$ holds. So, assume that the set $\CT_\Pi$ contains two distinct points $o,w$. Choose any point $u\in\Pi\setminus\Aline ow$. We can additionally assume that $u\in \CT_\Pi$ if $\CT_\Pi\not\subseteq\Aline ow$.
Since the points $o,w$ are central, the based affine plane $(\Pi,uow)$ satisfies the condition (4) of Theorem~\ref{t:leftcorps}, which implies that $\Aline ow\subseteq \CT_\Pi$. If $\Aline ow=\CT_\Pi$, then the case (2) holds. 

If $\Aline ow\ne\CT_\Pi$, then the choice of the point $u$ ensures that $u\in\CT_\Pi$. In this case the  the based affine plane $(\Pi,uow)$ satisfies the condition (3) of Theorem~\ref{t:leftcorps}, which implies that $\CT+\Pi=\Pi$ and hence the case (3) holds.
\end{proof}

Theorems~\ref{t:leftcorps} and \ref{t:corps<=>} imply the following theorem describing the inverplay between the central rank $\|\CT_\Pi\|$ of an affine plane and its coordinatizations by nice ternars.

\begin{theorem}\label{t:coord-central} An affine plane $\Pi$ is coordinatized by 
\begin{enumerate}
\item a homocentral ternar iff $\|\CT_\Pi\|\ge 1$;
\item a linear left-distributive associative ternar iff $\|\CT_\Pi\|\ge 2$;
\item a linear distributive associative ternar iff $\|\CT_\Pi\|=3$.
\end{enumerate}
\end{theorem}

For computational purposes, the following characterization of central points can be helpful. This characterization can be proved by analogy with 
Theorems~\ref{t:aa'-translation<=>} and \ref{t:local-homothety}.

\begin{proposition} A point $o$ of an affine plane $\Pi$ is central if and only if for every distinct lines $A,B,C\subset \Pi$ with $o\in A\cap B\cap C$ and any  distinct points $a,a'\in A\setminus\{o\}$, $b,b'\in B\setminus\{o\}$, $c,c'\in C\setminus\{o\}$ with $\Aline ab\cap\Aline {a'}{b'}=\varnothing=\Aline bc\cap\Aline{b'}{c'}$, we have $\Aline ac\cap\Aline{a'}{c'}=\varnothing$.
\end{proposition}

\begin{exercise} Show that the algorithmic complexity of finding the central trace of an affine plane of order $n$ is at most $O(n^6)$.
\smallskip

{\em Hint:} Apply Theorem~\ref{t:homocentral<=>}.
\end{exercise}

\section{The translation-shear classification of affine planes}

\begin{definition} For an affine plane $\Pi$, the set  
$$\TS_\Pi\defeq\TT_\Pi\cup\NT_\Pi$$
is called the \index{translation-shear trace}\index{trace!translation-shear}\defterm{translation-shear trace} of $\Pi$. 

The \index{translation-shear rank}\index{rank!translation-shear}\defterm{translation-shear rank} $\|\TS_\Pi\|$ of $\Pi$ is the pair $in$ of the numbers $i\defeq\|\TT_\Pi\|$ and $n\defeq \|\NT_\Pi\|$.
\end{definition}

\begin{theorem}\label{t:transhear-class} For every affine plane $\Pi$, exactly one of the following $8$ cases hold:
\begin{enumerate}
\item $\|\mathsf{TS}_\Pi\|=00$ and $\TS_\Pi=\varnothing$;   
\item $\|\mathsf{TS}_\Pi\|=01$ and $\TS_\Pi=\{L\}$ for some line $L$ in $\Pi$;
\item $\|\mathsf{TS}_\Pi\|=02$ and $\TS_\Pi=\mathcal L_p$ for some point $p\in\Pi$;
\item $\|\mathsf{TS}_\Pi\|=10$ and $\TS_\Pi=\{\delta\}$ for some direction $\delta\in\partial \Pi$;
\item $\|\mathsf{TS}_\Pi\|=12$ and $\TS_\Pi=\{\delta\}\cup \delta$, for some direction $\delta\in\partial \Pi$;
\item $\|\mathsf{TS}_\Pi\|=20$ and $\TS_\Pi=\partial \Pi$;
\item $\|\mathsf{TS}_\Pi\|=22$ and $\TS_\Pi=\partial \Pi\cup \delta$ for some direction $\delta\in\partial\Pi$;
\item $\|\mathsf{TS}_\Pi\|=23$ and $\TS_\Pi=\partial \Pi\cup\mathcal L_\Pi$.
\end{enumerate}
\end{theorem}

\begin{proof} By Theorems~\ref{t:trans-clas} and \ref{t:shear-class}, $\TT_\Pi\in\{\varnothing\}\cup\{\{\delta\}:\delta\in\partial\Pi\}\cup\{\partial\Pi\}$ and
$$\NT_\Pi\in\{\varnothing\}\cup\{\{L\}:L\in\mathcal L_\Pi\}\cup\{\mathcal L_p:p\in\Pi\}\cup\{\delta:\delta\in\partial\Pi\}\cup\{\mathcal L_\Pi\}.$$

Depending on the structure of the translation trace $\TT_\Pi$, we consider three cases.
\smallskip

0. First assume that $\TT_\Pi=\varnothing$. If $\NT_\Pi=\varnothing$, then $\TS_\Pi=\TT_\Pi\cup\NT_\Pi=\varnothing$ and $\|\TS_\Pi\|=00$. If $\NT_\Pi=\{L\}$ for some line $L\in\mathcal L_\Pi$, then $\TS_\Pi=\TT_\Pi\cup\NT_\Pi=\{L\}$ and $\|\TS_\Pi\|=01$. If $\NT_\Pi=\mathcal L_p$ for some $p\in\Pi$, then $\TS_\Pi=\TT_\Pi\cup\NT_\Pi=\mathcal L_p$ and $\|\TS_\Pi\|=02$. If $$\NT_\Pi\notin\{\varnothing\}\cup\{\{L\}:L\in\mathcal L_\Pi\}\cup\{\mathcal L_p:p\in\Pi\},$$
then $\NT_\Pi\in\{\{\delta:\delta\in \partial\Pi\}\cup\{\mathcal L_\Pi\}$ and hence $\delta\subseteq\NT_\Pi$ for some direction $\delta\in\partial\Pi$. Choose an affine base $uow$ for the affine plane $\Pi$ such that $\Aline ow\in \delta$. It follows from $\delta\subseteq\NT_\Pi$ that every vertical line in the based affine plane $(\Pi,uow)$ is shear. By Theorem~\ref{t:vtrans+vshear=vshears}, the based affine plane $(\Pi,uow)$ is vertical-traslation and hence $\Aline ow\in \TT_\Pi=\varnothing$, which is a contradiction showing that the case $\NT_\Pi\in\{\{\delta:\delta\in \partial\Pi\}\cup\{\mathcal L_\Pi\}$ is impossible.
\smallskip

1. Next assume that $\TT_\Pi=\{\delta\}$ for some direction $\delta\in\partial \Pi$. We claim that $\NT_\Pi\subseteq\delta$. 
In the opposite case, we can find a shear line $L\in\NT_\Pi\setminus\delta$ and  choose an affine base $uow$ in $\Pi$ such that $\Aline ow=L$ and $\Aline ou\in \delta$. Then the based affine plane $(\Pi,uow)$ is vertical-shear (because $\Aline ow=L\in\NT_\Pi$) and horizontal-translation (because $\Aline ou \in\delta\in\TT_\Pi$). Applying Theorem~\ref{t:semifield<=>}, we conclude that the affine plane $\Pi$ is translation and hence $\TT_\Pi=\partial\Pi$, which contradicts our assumption. This contradiction shows that $\NT_\Pi\subseteq \delta$.

If $\NT_\Pi=\varnothing$, then $\ST_\Pi=\TT_\Pi\cup\NT_\Pi=\{\delta\}$ and $\|\ST_\Pi\|=10$. If $\NT_\Pi\ne\varnothing$, then we can fix a line $L\in \NT_\Pi\subseteq \delta$ and choose an affine base $uow$ for the affine plane $\Pi$ such that $\Aline ow=L\in\delta$. In this case the based affine plane $(\Pi,uow)$ is vertical-translation and vertical-shear. Applying Theorem~\ref{t:vtrans+vshear=vshears}, we conclude that $\NT_\Pi=\delta$ and hence $\ST_\Pi=\TT_\Pi\cup\NT_\Pi=\{\delta\}\cup\delta$ and $\|\ST_\Pi\|=12$.
\smallskip

2. Finally, assume that $\TT_\Pi=\partial\Pi$ and hence $\Pi$ is a translation plane. 
If $\NT_\Pi=\varnothing$, then $\ST_\Pi=\TT_\Pi\cup\NT_\Pi=\partial\Pi$ and $\|\ST_\Pi\|=20$. Next, assume that the set $\NT_\Pi$ is not empty and take any line $L\in\NT_\Pi$. Since the plane $\Pi$ is translation, $L\in\NT_\Pi$ implies $L_\parallel\subseteq \NT_\Pi$. If $L_\parallel =\NT_\Pi$, then $\ST_\Pi=\TT_\Pi\cup\NT_\Pi=\partial\Pi\cup\{\delta\}$ for the direction $\delta\defeq L_\parallel$ and $\|\ST_\parallel\|=22$. If $L_\parallel\ne\NT_\Pi$, then $\NT_\Pi=\mathcal L_\Pi$, by Theorem~\ref{t:shear-class}. In this case $\ST_\Pi=\TT_\Pi\cup\NT_\Pi=\partial\Pi\cup\mathcal L_\Pi$ and $\|\ST_\Pi\|=23$. 
\end{proof}

\begin{remark}\label{r:HeKantor} By a difficult result of Hering and Kantor \cite{HK1971}, a finite affine plane $\Pi$  cannot have translation-shear rank $\|\ST_\Pi\|=02$. The proof of this fact is not elementary (it involves some deep group-theoretic considerations), so will be not reproduced here.
\end{remark}

Endow the 8-element set of translation-shear ranks $$\{00,01,02,10,12,20,22,23\}$$ with the partial order $\preceq$ defined by
$$
in\preceq jk\;\Leftrightarrow\;i\le j\wedge  n\le r \wedge \neg (0=i<j\wedge n=k=2).
$$
 
The Hasse diagram of this partial order is drawn in the following diagram:
$$
\xymatrix{
&23\\
&22\ar@{-}[u]&02\ar@{-}[ul]\\
20\ar@{-}[ur]&12\ar@{-}[u]\\
10\ar@{-}[ur]\ar@{-}[u]&&01\ar@{-}[uu]\ar@{-}[ul]\\
&00\ar@{-}[ul]\ar@{-}[ur]
}
$$ 



The partial order $\preceq$ on the set of all translation-shear ranks allows us to formulate the following coordinatization result.

\begin{theorem}\label{t:coord-transhear} An affine plane $\Pi$ is coordinatized by 
\begin{itemize}
\item a vertical-shear ternar iff $01\preceq \|\TS_\Pi\|$;
\item a linear associative-plus ternar ($=$ cartesian  group) iff $10\preceq\|\TS_\Pi\|$;
\item a linear left-distributive associative-plus ternar  iff $12\preceq\|\TS_\Pi\|$;
\item  a linear right-distributive associative-plus ternar ($=$ quasi-field) iff $20\preceq\|\TS_\Pi\|$;
\item a linear distributive associative-plus ternar ($=$ semi-field) iff $22\preceq\|\TS_\Pi\|$;
\item a linear distributive alternative-dot associative-plus ternar iff $23=\|\TS_\Pi\|$.
\end{itemize}
\end{theorem}
   
\begin{proof} 01. If the affine plane $\Pi$ is coordinatized by a vertical-shear ternar, then there exists an affine base $uow$ in $\Pi$ whose ternar is vertical-shear. By Theorem~\ref{t:vshear<=>vshear}, the based affine plane $(\Pi,uow)$ is vertical-shear and hence $\Aline ow\in\NT_\Pi\subseteq \TS_\Pi$, which implies $01\preceq\|\TS_\Pi\|$. 

On the other hand, if $01\preceq\|\TS_\Pi\|$, then $\NT_\Pi\ne\varnothing$ and we can choose an affine base $uow$ for the affine plane $\Pi$ so that $\Aline ow\in \NT_\Pi$ and hence the based affine plane $(\Pi,uow)$ is vertical-shear and its ternar $\Delta$ is vertical-shear, by Theorem~\ref{t:vshear<=>vshear}. Therefore, the affine plane is coordinatized by the vertical-shear ternar $\Delta$.
\smallskip

10.  If the affine plane $\Pi$ is coordinatized by a linear associative-plus ternar,  then there exists an affine base $uow$ in $\Pi$ whose ternar is linear and associative-plus. By Theorem~\ref{t:cart-group<=>}, the based affine plane $(\Pi,uow)$ is vertical-translation, and hence $\|TT_\Pi\|\ge 1$ and $10\preceq\|\ST_\Pi\|$, by the definition of the partial order $\preceq$. 

On the other hand, if $10\preceq\|\ST_\Pi\|$, then the affine plane $\Pi$ has a translation direction $\delta$. Take any affine base $uow$ for the affine plane $\Pi$  with $\Aline ow\in\delta$. Then the based affine plane $(\Pi,uow)$ is vertical-translation and its ternar $\Delta$ is linear and associative-plus, by Theorem~\ref{t:cart-group<=>}. Therefore, the affine plane is coordinatized by the linear associative-plus ternar $\Delta$.
\smallskip 

12.  If the affine plane $\Pi$ is coordinatized by a linear left-distributive associative-plus ternar,  then there exists an affine base $uow$ in $\Pi$ whose ternar is linear, left-distributive, and associative-plus. By Theorem~\ref{t:vtrans+vshear=vshears}, the based affine plane $(\Pi,uow)$ is vertical-translation and vertical-shear, moreover, $(\Aline ow)_\parallel\in\TT_\Pi$ and $(\Aline ow)_\parallel \subseteq \NT_\Pi$, which implies $12\preceq\|\TS_\Pi\|$.

On the other hand, if $12\preceq\|\TS_\Pi\|$, then Theorem~\ref{t:transhear-class} ensures that $\{\delta\}\cup\delta\subseteq\TS_\Pi$ for some direction $\delta\in\partial\Pi$. Take any affine base $uow$ for $\Pi$  with $\Aline ow\in\delta$. Then the based affine plane $(\Pi,uow)$ is vertical-translation and vertical-shear and its ternar $\Delta$ is linear, left-distrubitive and associative-plus, by Theorem~\ref{t:vtrans+vshear=vshears}. Therefore, the affine plane is coornatized by the linear left-distributive associative-plus ternar  $\Delta$.
\smallskip

20.  If the affine plane $\Pi$ is coordinatized by a linear right-distributive associative-plus ternar,  then there exists an affine base $uow$ in $\Pi$ whose ternar is linear, right-distributive, and associative-plus. By Theorem~\ref{t:trans-coordinate}, the based affine plane $(\Pi,uow)$ is translation and hence $\|\TT_\Pi\|=2$ and $20\preceq\|\TS_\Pi\|$.

On the other hand, if $20\preceq\|\TS_\Pi\|$, then $\|\TT_\Pi\|=2$ and Theorem~\ref{t:trans-coordinate} ensures that the affine plane $\Pi$ is translation and is coordinatized by a linear right-distributive associative-plus ternar $\Delta$.
\smallskip

22.  If the affine plane $\Pi$ is coordinatized by a linear distributive associative-plus ternar,  then there exists an affine base $uow$ in $\Pi$ whose ternar is linear, distributive, and associative-plus. By Theorem~\ref{t:semifield<=>}, the based affine plane $(\Pi,uow)$ is translation and vertical-shear. Moreover, $\TT_\Pi=\partial\Pi$ and $(\Aline ow)_\parallel \subseteq\NT_\Pi,$ which implies the inequality $22\preceq \|\TS_\Pi\|$.

On the other hand, if $22\preceq\|\TS_\Pi\|$, then Theorem~\ref{t:transhear-class} ensures that $\partial\Pi\cup\delta\subseteq \TS_\Pi$ for some direction $\delta\in\partial\Pi$. Choose any affine base $uow$ for the affine plane $\Pi$ with $\Aline ow\in\delta$. By Theorem~\ref{t:semifield<=>}, the based affine plane $(\Pi,uow)$ is translation and vertical-shear, and its ternar $\Delta$ is linear, distrubutive and associative-plus.
\smallskip

23.  If the affine plane $\Pi$ is coordinatized by a linear distributive alternative-dot associative-plus ternar,  then there exists an affine base $uow$ in $\Pi$ whose ternar is linear, distributive, associative-plus and alternative-dot. By Theorem~\ref{t:shear<=>alternative}, the based affine plane $(\Pi,uow)$ is translation and shear, which means that $\TS_\Pi=\TT_\Pi\cup\NT_\Pi=\partial\Pi\cup\mathcal L_\Pi$ and hence $\|\TS_\Pi\|=23$.

On the other hand, if $23=\|\TS_\Pi\|$, then Theorem~\ref{t:transhear-class} ensures that $\TS_\Pi=\partial\Pi\cup\mathcal L_\Pi$. Then for every affine base $uow$ for the affine plane $\Pi$, the based affine plane $(\Pi,uow)$ is shear. By Theorem~\ref{t:shear<=>alternative}, the ternar of the shear based affine plane $(\Pi,uow)$ is linear, distributive, aletrnative-dot and associative-plus.
\end{proof}  
  
\section{The translation-shear-central classification of affine planes}  

\begin{definition} For an affine plane $\Pi$, the set $$\TSH_\Pi\defeq \TT_\Pi\cup\NT_\Pi\cup\CT_\Pi$$ is called the $\TSH$-trace of the affine plane $\Pi$. 

The \index{translation-shear-central rank}\index{rank!translation-shear-central}\defterm{translation-shear-central rank} $\|\TSH_\Pi\|$ of the affine plane $\Pi$ is the triple $inx$ of the numbers $i\defeq \|\TT_\Pi\|$, $n\defeq\|\NT_\Pi\|$, $x\defeq\|\CT_\Pi\|$.
\end{definition}

\begin{theorem}\label{t:INX-class} For every affine plane $\Pi$, exactly one of the following $13$ cases hold:
\begin{enumerate}
\item $\|\TSH_\Pi\|=000$ and $\TSH_\Pi=\varnothing$;   
\item $\|\TSH_\Pi\|=001$ and $\TSH_\Pi=\{p\}$ for some point $p\in\Pi$;   
\item $\|\TSH_\Pi\|=010$ and $\TSH_\Pi=\{L\}$ for some line $L$ in $\Pi$;
\item $\|\TSH_\Pi\|=011$ and $\TSH_\Pi=\{L,p\}$ for some line $L$ in $\Pi$ and point $p\in L$;
\item $\|\TSH_\Pi\|=020$ and $\TSH_\Pi=\mathcal L_p$ for some point $p\in\Pi$;
\item $\|\TSH_\Pi\|=021$ and $\TSH_\Pi=\mathcal L_p\cup\{p\}$ for some point $p\in\Pi$;
\item $\|\TSH_\Pi\|=100$ and $\TSH_\Pi=\{\delta\}$ for some direction $\delta\in\partial \Pi$;
\item $\|\TSH_\Pi\|=120$ and $\TSH_\Pi=\{\delta\}\cup \delta$, for some direction $\delta\in\partial \Pi$;
\item $\|\TSH_\Pi\|=122$ and $\TSH_\Pi=\{\delta\}\cup \delta\cup L$, for some direction $\delta\in\partial \Pi$ and some line $L\in \delta$;
\item $\|\TSH_\Pi\|=200$ and $\TSH_\Pi=\partial \Pi$;
\item $\|\TSH_\Pi\|=220$ and $\TSH_\Pi=\partial \Pi\cup \delta$ for some direction $\delta\in\partial\Pi$;
\item $\|\TSH_\Pi\|=230$ and $\TSH_\Pi=\partial \Pi\cup\mathcal L_\Pi$;
\item $\|\TSH_\Pi\|=233$ and $\TSH_\Pi=\partial \Pi\cup\mathcal L_\Pi\cup \Pi$.
\end{enumerate}
\end{theorem}
  
\begin{proof} We follow the translation-shear classification, taking into account the structure of the central trace. By Theorem~\ref{t:transhear-class}, the following 8 types of the translation-shear trace $\TS_\Pi$ exist.
\vskip3pt

00. Assume that $\TS_\Pi=\varnothing$. We claim that $|\CT_\Pi|\le 1$. Assuming that $|\CT_\Pi|\ge 2$, we can find an affine base $uow$ for $\Pi$ such that $o,w\in \CT_\Pi$ and hence the points $o$ and $w$ are central. Applying Theorem~\ref{t:leftcorps}, we conclude that the based affine plane $(\Pi,uow)$ is vertical-translation and hence $|\TT_\Pi|\ge 1$, which contradicts $\TS_\Pi=\varnothing$. This contradiction shows that $|\CT_\Pi|\le 1$.
If $|\CT_\Pi|=0$, then $\TSH_\Pi=\TS_\Pi\cup\CT_\Pi=\varnothing$ and $\|\TSH_\Pi\|=000$. If $|\CT_\Pi|=1$, then $\TSH_\Pi=\CT_\Pi=\{p\}$ for some point $p\in\Pi$ and $\|\TSH_\Pi\|=001$.
\smallskip

01. Assume that $\TS_\Pi=\{L\}$ for some line $L$ in $\Pi$ and hence $\TT_\Pi=\varnothing$ and $\NT_\Pi=\{L\}$. Repeating the argument from the preceding case, we can show that $|\CT_\Pi|\le 1$. 

If $|\CT_\Pi|=0$, then $\TSH_\Pi=\TS_\Pi\cup\CT_\Pi=\{L\}$ and $\|\TSH_\Pi\|=010$. 

If $|\CT_\Pi|=1$, then $\TSH_\Pi=\TS_\Pi\cup\CT_\Pi=\{L,p\}$ for some point $p\in\Pi$. Assuming that $p\notin L$ and taking into account that $L$ is a shear line, we can choose any point $q\in \Pi\setminus\{p\}$ with $\Aline pq\parallel L$ and find an automorphism $A:\Pi\to\Pi$ such that $A(p)=q$ and conclude that $q$ is a central point and hence $\{p,q\}\subseteq\CT_\Pi$, which contradicts $|\CT_\Pi|\le 1$. This contradiction shows that $p\in L$ and hence $\TSH_\Pi=\{L,p\}$ for a line $L\subseteq \Pi$ and point $p\in L$. In this case $\|\TSH_\Pi\|=011$.
\smallskip

02. Assume that $\TS_\Pi=\mathcal L_p$ for some point $p\in\Pi$. In this case $\TT_\Pi=\varnothing$ and $\NT_\Pi=\mathcal L_p$. Repeating the argument from the case 00, we can show that $|\CT_\Pi|\le 1$. 

If $|\CT_\Pi|=0$, then $\TSH_\Pi=\TS_\Pi\cup\CT_\Pi=\mathcal L_p$ and $\|\TSH_\Pi\|=020$. 

If $|\CT_\Pi|=1$, then $\TSH_\Pi=\TS_\Pi\cup\CT_\Pi=\mathcal L_p\cup\{c\}$ for some point $c\in\Pi$. Repeating the argument from the preceding case, we can show that $c=p$ and hence $\TSH_\Pi=\TS_\Pi\cup\CT_\Pi=\mathcal L_p\cup\{p\}$ and $\|\TSH_\Pi\|=021$.
\smallskip

10. Assume that $\TS_\Pi=\{\delta\}$ for some direction $\delta\in\Pi$. We claim that $\CT_\Pi=\varnothing$. Assuming that the set $\CT_\Pi$ is not empty, we can choose an affine base $uow$ for the plane $\Pi$ such that $o\in\CT_\Pi$ and $\Aline ow\in\delta$. Then the based affine plane $(\Pi,uow)$ is homocentral and vertical-translation. By Theorem~\ref{t:leftcorps}, the ternar $\Delta$ of the affine base $uow$ is linear, left-distributive and associative. Since the ternar $\Delta$ is linear, left-distributive and associative-plus, we can apply Theorem~\ref{t:vtrans+vshear=vshears} and conclude that the based affine space $(\Pi,uow)$ is vertical-shear and hence $\Aline ow\in\NT_\Pi\subseteq \TS_\Pi$, which contradicts $\TS_\Pi=\{\delta\}$. This contradiction shows that $\CT_\Pi=\varnothing$ and hence $\TSH_\Pi=\TS_\Pi\cup\CT_\Pi=\{\delta\}$ and $\|\TSH_\Pi\|=100$.
\smallskip

12. Assume that $\TS_\Pi=\{\delta\}\cup\delta$ for some direction $\delta\in \partial\Pi$. Then $\TT_\Pi=\{\delta\}$ and $\NT_\Pi=\delta$. We claim that $\|\CT_\Pi\|\ne 3$. Assuming that $\|\CT_\Pi\|=3$, we can apply Theorem~\ref{t:central-class} and conclude that $\CT_\Pi=\Pi$. By Theorems~\ref{t:corps<=>}, the affine plane $\Pi$ is translation, which contradicts $\TT_\Pi=\{\delta\}$. This contradiction shows that $\|\CT_\Pi\|\ne 3$. 
If $\CT_\Pi=\varnothing$, then $\TSH_\Pi=\TS_\Pi\cup\CT_\Pi=\{\delta\}\cup\delta$ and $\|\TSH_\Pi\|=120$. So, assume that $\CT_\Pi\ne\varnothing$. In this case we can choose an affine base $uow$ for the plane $\Pi$ such that $o\in\CT_\Pi$ and $\Aline ow\in\delta$. Then  the based affine plane $(\Pi,uow)$ is vertical-translation and homocentral. Applying Theorem~\ref{t:leftcorps}, we conclude that $\Aline ow\subseteq\CT_\Pi$. Since $\|\CT_\Pi\|\ne 3$, Theorem~\ref{t:central-class} ensures that $\Aline ow=\CT_\Pi$. Then $\TSH_\Pi=\TS_\Pi\cup\CT_\Pi=\{\delta\}\cup\delta\cup L$ for the line $L\defeq\Aline ow\in\delta$ and $\|\TSH_\Pi\|=122$.  
\vskip3pt

20. Assume that $\TS_\Pi=\partial \Pi$. Then $\TT_\Pi=\partial\Pi$ and $\NT_\Pi=\varnothing$. Assuming that $\CT_\Pi\ne\varnothing$ and taking into account that $\Pi$ is a translation plane, we conclude that $\CT_\Pi=\Pi$. Fix any affine base $uow$ for the plane $\Pi$. By Theorem~\ref{t:corps<=>}, the ternar $\Delta$ of the base $uow$ is linear, distributive and associative. In particular, $\Delta$ is linear, distributive, associative-plus and alternative-dot. By Theorem~\ref{t:shear<=>alternative}, $\Pi$ is a shear plane, which contradicts $\NT_\Pi=\varnothing$. This contradiction shows that $\CT_\Pi=\varnothing$ and hence $\TSH_\Pi=\TS_\Pi\cup\CT_\Pi=\partial\Pi$ and $\|\TSH_\Pi\|=200$.
\smallskip

22. Assume that $\TS_\Pi=\partial\Pi\cup\delta$ for some direction $\delta\in\partial\Pi$. Then $\TT_\Pi=\partial\Pi$ and $\NT_\Pi=\delta\ne\mathcal L_\Pi$. Repeating the argument from the preceding case, we can show that $\CT_\Pi=\varnothing$ and hence $\TSH_\Pi=\TS_\Pi\cup\CT_\Pi=\partial\Pi\cup\delta$ and $\|\TSH_\Pi\|=220$.

23. Assume that $\TS_\Pi=\partial\Pi\cup\mathcal L_\Pi$. If $\CT_\Pi=\varnothing$, then $\TSH_\Pi=\TS_\Pi\cup\CT_\Pi=\partial\Pi\cup\mathcal L_\Pi$ and $\|\TSH_\Pi\|=230$. If $\CT_\Pi\ne\varnothing$, then $\CT_\Pi=\Pi$ (because the affine plane $\Pi$ is translation). In this case $\TSH_\Pi=\TS_\Pi\cup\CT_\Pi=\partial\Pi\cup\mathcal L_\Pi\cup\Pi$ and $\|\TSH_\Pi\|=233$.
\end{proof}

Endow the 13-element set $$\{000,001,010,011,020,021,100,120,122,200,220,230,233\}$$ with the partial order $\preceq$ defined by
$$
inx\le jky\;\Leftrightarrow\;i\le j\wedge  n\le k\wedge x\le y \wedge \neg (0=i<j\wedge n=k=2).
$$
 
The Hasse diagram of this partial order  is drawn in the following diagram:
$$
\xymatrix{
&&233\\
&230\ar@{-}[ur]&\\
&220\ar@{-}[u]&122\ar@{-}[uu]&021\ar@{-}[uul]\\
200\ar@{-}[ur]&120\ar@{-}[u]\ar@{-}[ur]&020\ar@{-}[ur]\ar@{-}[uul]&011\ar@{-}[u]\ar@{-}[lu]\\
&100\ar@{-}[ul]\ar@{-}[u]&010\ar@{-}[lu]\ar@{-}[ur]\ar@{-}[u]&001\ar@{-}[u]\\
&&000\ar@{-}[ul]\ar@{-}[u]\ar@{-}[ur]
}
$$



The partial order $\preceq$ on the set of all translation-shear-central ranks allows us to formulate the following coordinatization result, extending Theorem~\ref{t:coord-transhear}.

\begin{theorem}\label{t:coord-INX} An affine plane $\Pi$ is coordinatized by 
\begin{itemize}
\item a homocentral ternar iff $001\preceq\|\TSH_\Pi\|$;
\item a vertical-shear ternar iff $010\preceq\|\TSH_\Pi\|$;
\item a vertical-shear homocentral ternar iff $011\preceq\|\TSH_\Pi\|$;
\item  a linear associative-plus ternar $(=$ cartesian group$)$ iff $100\preceq\|\TSH_\Pi\|$;
\item a linear left-distributive associative-plus ternar iff $120\preceq\|\TSH_\Pi\|$;
\item a linear left-distributive associative ternar iff $122\preceq\|\TSH_\Pi\|$;
\item a linear right-distributive associative-plus ternar iff $200\preceq\|\TSH_\Pi\|$;
\item a linear distributive associative-plus ternar iff $220\preceq\|\TSH_\Pi\|$;
\item a linear distributive associative-plus alternative-dot ternar iff $230\preceq\|\TSH_\Pi\|$;
\item a linear distributive associative ternar $(=$ linear corps$)$ iff $233=\|\TSH_\Pi\|$.
\end{itemize}
\end{theorem}

\begin{proof} 001.  By Theorem~\ref{t:coord-central}, the affine plane $\Pi$ is coordinatized by a homocentral ternar iff $1\le\|\TT_\Pi\|$, which is equivalent to $001\preceq\|\TSH_\Pi\|$.
\smallskip

010. By Theorem~\ref{t:coord-transhear}, the affine plane $\Pi$ is coordinatized by a vertical-shear ternar iff $01\preceq\|\TS_\Pi\|$, which is equivalent to $010\preceq\|\TSH_\Pi\|$.
\smallskip

011. If the affine plane $\Pi$ is coordinatized by a vertical-shear homocentral ternar, then there exists an affine base $uow$ whose ternar $\Delta$ is vertical-shear and homocentral. By Theorems~\ref{t:vshear<=>vshear} and \ref{t:homocentral<=>}, the based affine plane $(\Pi,uow)$ is vertical-shear and homocentral and hence $\Aline ow\in\NT_\Pi$ and $o\in\CT_\Pi$, which implies $011\preceq \|\TSH_\Pi\|$.

On the other hand, if $011\preceq\|\TSH_\Pi\|$, then the definition of the partial order $\preceq$ and Theorem~\ref{t:INX-class}(4) imply that $\{L,p\}\subseteq \TSH_\Pi$ for some line $L$ and point $p\in L$. Choose an affine base $uow$ for the plane $\Pi$ such that $o=p$ and $\Aline ow=L$. Then the based affine plane $(\Pi,uow)$ is vertical-shear and homocentral. By Theorems~\ref{t:vshear<=>vshear} and \ref{t:homocentral<=>}, the ternar $\Delta$ of the base $uow$ is vertical-shear and homocentral, implying that the affine plane $\Pi$ is coordinatized by a vertical-shear homocentral ternar.
\smallskip

100. By Theorem~\ref{t:coord-transhear}, the affine plane $\Pi$ is coordinatized by a linear associative-plus ternar iff $10\preceq\|\TS_\Pi\|$, which is equivalent to  $100\preceq\|\TSH_\Pi\|$.
\smallskip

120. By Theorem~\ref{t:coord-transhear}, the affine plane $\Pi$ is coordinatized by a linear left-distributive associative-plus ternar iff $12\preceq\|\TS_\Pi\|$, which is equivalent to $120\preceq\|\TSH_\Pi\|$.
\smallskip

122. If the affine plane $\Pi$ is coordinatized by a linear left-distrubutive associative ternar, then there exists an affine base $uow$ whose ternar $\Delta$ is linear, left-distributive, and associative. By Theorems~\ref{t:leftcorps} and \ref{t:vtrans+vshear=vshears}, $(\Aline ow)_\parallel\in\TT_\Pi$, $(\Aline ow)_\parallel\subseteq\NT_\Pi$ and $\Aline ow\subseteq\CT_\Pi$, which implies  $122\preceq \|\TSH_\Pi\|$.

On the other hand, if $122\preceq\|\TSH_\Pi\|$, then the definition of the partial order $\preceq$ and Theorem~\ref{t:INX-class}(9) imply that $\{\delta\}\cup\delta\cup L\subseteq \TSH_\Pi$ for some direction $\delta\in\partial\Pi$ and line $L\in\delta$. Choose an affine base $uow$ for the plane $\Pi$ such that $\Aline ow=L$. Then the based affine plane $(\Pi,uow)$ is vertical-translation, vertical-shear and homocentral. By Theorem~\ref{t:leftcorps}, the ternar $\Delta$ of the based affine plane  $(\Pi,uow)$ is linear, left-distrubitive and associative, implying that the affine plane $\Pi$ is coordinatized by a linear left-distributive associative ternar.
\smallskip

200. By Theorem~\ref{t:coord-transhear}, the affine plane $\Pi$ is coordinatized by a linear right-distributive associative-plus ternar iff $20\preceq\|\TS_\Pi\|$, which is equivalent to $200\preceq\|\TSH_\Pi\|$.
\smallskip

220. By Theorem~\ref{t:coord-transhear}, the affine plane $\Pi$ is coordinatized by a linear distributive associative-plus ternar iff $22\preceq\|\TS_\Pi\|$, which is equivalent to $220\preceq\|\TSH_\Pi\|$.
\smallskip

230. By Theorem~\ref{t:coord-transhear}, the affine plane $\Pi$ is coordinatized by a linear distributive associative-plus alternative-dot ternar iff $23\preceq\|\TS_\Pi\|$, which is equivalent to $230\preceq\|\TSH_\Pi\|$.
\smallskip

233. If the affine plane $\Pi$ is coordinatized by a linear distributive associative ternar, then there exists an affine base $uow$ whose ternar $\Delta$ is linear, distributive, and associative. By Theorems~\ref{t:corps<=>} and \ref{t:VW-Thalesian<=>quasifield}, $\TSH_\Pi=\TT_\Pi\cup\NT_\Pi\cup\CT_\Pi=\partial\Pi\cup\mathcal L_\Pi\cup\Pi$, which implies  $233=\|\TSH_\Pi\|$.

On the other hand, if $233=\|\TSH_\Pi\|$, then $\TSH_\Pi=\partial\Pi\cup\mathcal L_\Pi\cup\Pi$. Choose any affine base $uow$ for the plane $\Pi$ and obseve that the based affine plane $(\Pi,uow)$ is translation and homocentral.  By Theorem~\ref{t:corps<=>}, the ternar $\Delta$ of the based affine plane  $(\Pi,uow)$ is linear, distributive and associative, witnessing that the affine plane $\Pi$ is coordinatized by a linear distributive associative ternar.
\end{proof}

\begin{proposition} If an affine plane $\Pi$ is finite, then $\|\TSH_\Pi\|\notin\{020,021,230\}$.
\end{proposition}

\begin{proof} The result of Hering and Kantor mentioned in Remark~\ref{r:HeKantor} implies that $\|\TSH_\Pi\|\notin\{020,021\}$. Assuming that $\|\TSH_\Pi\|=230$, we can apply Theorem~\ref{t:coord-INX} and conclude that the affine plane $\Pi$ has an affine base $uow$ whose ternar is linear distributive associative-plus, and associative-dot.  Since the plane $\Pi$ is finite, so is the ternar $\Delta$. Then the biloop of the ternar $\Delta$ is a finite alternative ring.  By the Artin--Zorn Theorem~\ref{t:Artin-Zorn}, the finite alternative ring is a corps and hence the ternar $\Delta$ is distributive and associative. Therefore, the affine plane $\Pi$ is coordinatized by a liner distributive associative ternar. Applying Theorem~\ref{t:coord-INX}, we conclude that $\|\TSH_\Pi\|=233$, which contradicts our assumption $\|\TSH_\Pi\|=230$.
\end{proof}


\begin{remark} In the following table we present the translation-shear-central ranks of all seven affine planes of order $9$ (calculated by Ivan Hetman).
$$
\begin{array}{l|c|c|c|c|c|c|c}
\mbox{The affine plane:}&\mbox{\tt Desarg}&\mbox{\tt Thales}&\mbox{\tt Hall}&\mbox{\tt dhall}&\mbox{\tt hall}&\mbox{\tt Hughes}&\mbox{\tt hughes}\\
\hline
\mbox{The rank:}&233&200&122&011&000&000&000\\
\end{array}
$$
\end{remark}

\section{Hyperscale classification of affine planes}\label{s:hyperscale-class}

\begin{definition} Let $\Pi$ be an affine plane and $\mathcal L_\Pi$ be the family of all lines in $\Pi$. A pair $(\delta,L)\in\partial \Pi\times\mathcal L_\Pi$ is called a \index{hyperscale pair}\defterm{hyperscale pair} if $L\notin\delta$ and for every line $D\in\delta$ and points $x,y\in D\setminus L$, there exists an automorphism $A:\Pi\to\Pi$ such that $A(x)=y$ and $L\subseteq\Fix(A)$.
\end{definition}


 
\begin{definition} For an affine plane $\Pi$, the set $\ST_\Pi$ of all hyperscale pairs is called the \index{hyperscale trace}\index{trace!hyperscale}\defterm{hyperscale trace} of the affine plane $\Pi$. 
\end{definition}

The hyperscale trace $\ST_\Pi$ is a subset of the Cartesian product $\partial\Pi\times\mathcal L_\Pi$ and so $\ST_\Pi$ is a relation between directions and lines in the affine plane $\Pi$. Let 
$$
\begin{aligned}
\dom[\ST_\Pi]&\defeq\{\delta\in\partial\Pi:\exists L\in\mathcal L_\Pi\;\;(\delta,L)\in\ST_\Pi\}\quad\mbox{and}\\
\rng[\ST_\Pi]&\defeq\{L\in\mathcal L_\Pi:\exists \delta\in\partial \Pi\;\;(\delta,L)\in\ST_\Pi\}
\end{aligned}
$$be the domain and range of the relation $\ST_\Pi$.

The complexity of the hyperscale trace can be measures by its birank, defined as follows.

\begin{definition} For an affine plane $\Pi$, the \index{birank}\defterm{birank} $\|\ST_\Pi\|$ of the relation $\ST_\Pi$ is the pair $dr$ of the numbers $d\defeq\|\dom[\ST_\Pi]\|$ and $r\defeq\|\rng[\ST_\Pi]\|$.
\end{definition}

 In Theorem~\ref{t:scale-trace} we shall classify the hyperscale traces of affine planes into 10 possible types. To approach smoothly this classification, we first consider a somewhat simpler problem of classifying the relation
$$\PSC_\Pi=\{(\delta,\delta'):\exists L\in\delta'\;(\delta,L)\in\SC_\Pi\}\subseteq\partial\Pi\times\partial\Pi,$$called the \defterm{hyperscale figure} of an affine plane $\Pi$. The hyperscale figure is the image of the hyperscale trace $\SC_\Pi$ under the projection $$\partial\Pi\times\mathcal L_\Pi\to\partial\Pi\times\partial\Pi,\quad(\delta,L)\mapsto(\delta,L_\parallel).$$ 
The complexity of the hyperscale figure $\PSC_\Pi$ can be measures by its birank, defined as follows.

 \begin{definition} For an affine plane $\Pi$, the \defterm{birank} $\|\PSC_\Pi\|$ of the relation $\PSC_\Pi$ is the pair $dr$ of the numbers $d\defeq\|\dom[\PSC_\Pi]\|$ and $r\defeq\|\rng[\PSC_\Pi]\|$.
\end{definition}

Observe that for every affine plane $\Pi$, the relation $\PSC_\Pi$ is irreflexive. i.e., $(\delta,\delta)\notin\PSC_\Pi$ for all $\delta\in \partial\Pi$. A relation $F$ is a \defterm{function} if for any $x\in\dom[F]$ there exists a unique element $y\in \rng[F]$ (denoted by $F(x)$) such that $(x,y)\in F$. A relation $F$ is called an \defterm{involution} if $F$ is a function such that $F=F^{-1}$. 
 
\begin{theorem}\label{t:scale-figure} For every affine plane $\Pi$, exactly one of the following $7$ cases holds:
\begin{enumerate}
\item[{$(\,\circ\,)$}]\index[note]{$(\,\circ\,)$}  $\|\PSC_\Pi\|=00$ and $\PSC_\Pi=\varnothing$;
\item[{$(\,\bullet\,)$}] \index[note]{$(\,\bullet\,)$} $\|\PSC_\Pi\|=11$ and $\PSC_\Pi=\{(h,v)\}$ for two distinct directions $h,v\in\partial\Pi$;
\item[{$(\hskip2.5pt\vert\hskip2.5pt)$}]\index[note]{$(\hskip2.5pt\vert\hskip2.5pt)$}  $\|\PSC_\Pi\|=12$ and $\PSC_\Pi=\{\delta\}\times (\partial\Pi\setminus\{\delta\})$ for some direction $\delta\in\partial\Pi$;
\item[{$(-)$}]\index[note]{$(-)$} $\|\PSC_\Pi\|=21$ and  $\PSC_\Pi=(\partial \Pi\setminus\{\delta\})\times\{\delta\}$ for some direction $\delta\in\partial\Pi$;
\item[$(\didots)$]\index[note]{$(\didots)$} $\|\PSC_\Pi\|=22$ and $\PSC_\Pi=\{(\delta,\delta'),(\delta',\delta)\}$ for two distinct directions $\delta,\delta'\in\partial\Pi$;
\item[{$(\invol)$}]\index[note]{$(\invol)$} $\|\PSC_\Pi\|=22$ and $\PSC_\Pi$ is an involution with $|\PSC_\Pi|=10=|\partial\Pi|$;
\item[{$(\Dscale)$}]\index[note]{$(\Dscale)$} $\|\PSC_\Pi\|=22$ and $\PSC_\Pi=\{(\delta,\delta')\in\partial\Pi\times\partial\Pi:\delta\ne\delta'\}$.
\end{enumerate}
\end{theorem}

\begin{proof} 
 The proof of Theorem~\ref{t:scale-figure} is divided into 19 lemmas, which will be used also in the proof of Theorem~\ref{t:scale-trace}.

First we show that the hyperscale trace and the hyperscale figure both are invariant under automorphisms of a projective plane. 

Given an automorphism $\Phi$ of an affine plane $\Pi$, let $\partial\Phi:\partial\Pi\to\partial\Pi$ be the bijection assigning to each direction $\delta\in\partial\Pi$, the direction $\partial\Phi(\delta)\defeq\{\Phi[L]:L\in\delta\}$, and let $\overline \Phi\defeq\Phi\cup\partial\Phi$ be the spread completion of the automorphism $\Phi$.

\begin{lemma}\label{l:SC-auto} For any automorphism $\Phi$ of an affine plane $\Pi$, $$\{(\bar\Phi(\delta),\Phi[L]):(\delta,L)\in\SC_\Pi\}=\SC_\Pi\quad\mbox{and}\quad\{(\bar\Phi(\delta),\bar\Phi(\delta')):(\delta,\delta')\in\PSC_\Pi\}=\PSC_\Pi.$$
\end{lemma}

\begin{proof} First we prove that  $(\bar\Phi(\delta),\Phi[L])\in\SC_\Pi$ for all $(\delta,L)\in\SC_\Pi$. Given any distinct points $y,y'\in \Pi\setminus\Phi[L]$ with $\Aline y{y'}\in\bar\Phi(\delta)$, we should find an automorphism $A$ of the affine plane $\Pi$ such that $A(y)=y'$ and $\Phi[L]\subseteq \Fix(A)$.

Consider the points $x\defeq\Phi^{-1}(y)\in\Pi\setminus L$ and $x'\defeq\Phi^{-1}(y')\in \Pi\setminus L$ and observe that $\Aline y{y'}\in\bar\Phi(\delta)$ implies $\Aline x{x'}\in\delta$. Since $(\delta,L)\in\SC_\Pi$, there exists an automorphism $B$ of the affine plane $\Pi$ such that $B(x)=x'$ and $L\subseteq\Fix(B)$. Then the automorphism $A\defeq \Phi B\Phi^{-1}$ of $\Pi$ has the required properties: $A(y)=\Phi B\Phi^{-1}(y)=\Phi(B(x))=\Phi(x')=y'$ and $\Phi[L]\subseteq\Fix(A)$. Indeed, for any $z\in \Phi[L]$, the point $\Phi^{-1}(z)$ belongs to the line $L\subseteq\Fix(B)$ and hence $A(z)=\Phi B\Phi^{-1}(z)=\Phi(\Phi^{-1}(z))=z$, witnessing that $\Phi[L]\subseteq\Fix(A)$. Therefore, $\{(\bar\Phi(\delta),\Phi[L]):(\delta,L)\in\SC_\Pi\}\subseteq\SC_\Pi$. By analogy we can prove that for the inverse automorphism $\Psi\defeq\Phi^{-1}$ of the plane $\Pi$, we have $\{(\bar\Psi(\delta),\Psi[L]):(\delta,L)\in\SC_\Pi\}\subseteq\SC_\Pi$. Since $\Phi\circ\Psi$ is the identity automorphism of the affine plane $\Pi$, $\bar\Phi\circ\bar\Psi$ is the identity automorphism of its spread extension.
Then $$
\SC_\Pi=\{(\bar\Phi(\bar\Psi(\delta)),\Phi[\Psi[L]]):(\delta,L)\in\SC_\Pi\}\subseteq \{(\bar\Phi(\delta),\Phi[L]):(\delta,L)\in\SC_\Pi\}\subseteq\SC_\Pi$$and finally, $\{(\bar\Phi(\delta),\Phi[L]):(\delta,L)\in\SC_\Pi\}=\SC_\Pi$.
\smallskip

Next, we show that $\{(\bar \Phi(\delta),\bar\Phi(\delta')):(\delta,\delta')\in\PSC_\Pi\}=\PSC_\Pi$. Given any pair $(\delta,\delta')\in\PSC_\Pi$, find a line $L\in \delta'$ such that $(\delta,L)\in \SC_\Pi$. Taking into account that $(\bar\Phi(\delta),\Phi[L])\in\SC_\Pi$ and $\Phi[L]\in\bar\Phi(\delta')$, we conclude that $(\bar\Phi(\delta),\bar\Phi(\delta'))\in \PSC_\Pi$. Therefore, $\{(\bar \Phi(\delta),\bar\Phi(\delta')):(\delta,\delta')\in\PSC_\Pi\}\subseteq\PSC_\Pi$. Applying the same argument to the inverse automorphism $\Phi^{-1}$ we can prove that $\PSC_\Pi=\{(\bar \Phi(\delta),\bar\Phi(\delta')):(\delta,\delta')\in\PSC_\Pi\}$.
\end{proof}


\begin{lemma}\label{l:PSC-not-function}Let $\Pi$ be an affine plane. If $\{(\delta,h),(\delta,v)\}\subseteq \PSC_\Pi$ for some  distinct directions $\delta,h,v\in\partial\Pi$, then $\{\delta\}\times(\partial\Pi\setminus\delta)\subseteq\PSC_\Pi$ and $\{\delta\}\times(\mathcal L_p\setminus\delta)\subseteq\SC_\Pi$ for some point $p\in\Pi$.
\end{lemma}

\begin{proof} Since $\{(\delta,h),(\delta,v)\}\subseteq \PSC_\Pi$, there are lines $H\in h$ and $V\in v$ such that $\{(\delta,H),(\delta,V)\}\subseteq\SC_\Pi$. Since the directions $h,v$ are distinct, the lines $H\in h$ and $v\in V$ are concurrent and hence have a common point $p\in H\cap V$.
We claim that $\{\delta\}\times (\mathcal L_p\setminus\delta)\subseteq\SC_\Pi$. Given any line $L\in \mathcal L_p\setminus\{\delta\}$, we have to prove that $(\delta,L)\in\SC_\Pi$. If $L\in\{H,V\}$, then $(\delta,L)\in\{\delta,H),(\delta,V)\}\subseteq \SC_\Pi$ and we are done. So, assume that $H\ne L\ne V$. Choose any point $x\in V\setminus \{p\}$. Since $L\notin\delta$, there exists a unique point $y\in L$ such that $\Aline xy\in \delta$. Since $(\delta,H)\in\SC_\Pi$, there exists an automorphism $\Phi$ of the plane $\Pi$ such that $\Phi(x)=y$ and $p\in H\subseteq\Fix(\Phi)$. 
Let $\bar\Phi$ be the extension of the automorphism $\Phi$ to the spread completion $\overline\Pi=\Pi\cup\partial\Pi$ of $\Pi$. 

Observe that $\Phi[V]=\Phi[\Aline px]=\Aline py=L$. Since $\Aline xy\in\delta\ne H_\parallel$, there exists a unique point $o\in \Aline xy\cap H$. Then $\Phi[\Aline xy]=\Phi[\Aline ox]=\Aline oy=\Aline xy$ and hence $\bar\Phi(\delta)=\bar\Phi(\Aline xy_\parallel)=\Phi[\Aline xy]_\parallel=\Aline xy_\parallel=\delta$. Lemma~\ref{l:SC-auto} ensures that $(\delta,L)=(\bar\Phi(\delta),\Phi[V])\in\SC_\Pi$.  Therefore, $\{\delta\}\times(\mathcal L_p\setminus\delta)\subseteq\SC_\Pi$, which implies $\{\delta\}\times(\partial\Pi\setminus\delta)\subseteq\PSC_\Pi$. 
\end{proof}

\begin{lemma}\label{l:h-rozmaznia} Let $\Pi$ be an affine plane. If $\{(\delta,V),(\delta,\Lambda)\}\subseteq\SC_\Pi$ for some direction $\delta\in\partial\Pi$ and disjoint lines $V,\Lambda$, then $\{\delta\}\times V_\parallel\subseteq\SC_\Pi$.
\end{lemma}

\begin{proof} Choose an affine base $uow$ in the plane $\Pi$ such that $\Aline ow=V$ and $\Aline ue=\Lambda$, where $e$ is the diunit of the affine base $uow$. Given any line $L\in V_\parallel\setminus\{V,\Lambda\}$, we should prove that the pair $(\delta,L)$ belongs to the set $\SC_\Pi$. Given any line $D\in\delta$ and points $x,y\in D\setminus L$, we should find an automorphism $A$ of $\Pi$ such that $A(x)=y$ and $L\subseteq\Fix(A)$. Since $(\delta,V)\in \SC_\Pi$, there exists an automorphism $B$ of $\Pi$ such that $B(e)\in \Aline we\cap L$ and $V\subseteq \Fix(B)$. It follows from $\Aline ue\parallel V$ that $B[\Aline ue]\parallel B[V]=V\parallel L$ and hence $B[\Aline ue]=L$. Find a unique point $d\in D\cap V\subseteq\Fix(B)$ and observe that $D\parallel\Aline we$ imples $B[D]\parallel B[\Aline we]=\Aline we\parallel D$. Then $d\in D\cap B[D]$ implies $B[D]=D$ and hence $\{B^{-1}(x),B^{-1}(y)\}\subseteq B^{-1}[D\setminus L]=D\setminus \Lambda$. Since the pair $(\delta,\Lambda)$ belongs to $\SC_\Pi$, there exists an automorphism $C$ of the plane $\Pi$ such that $C(B^{-1}(x))=B^{-1}(y)$ and $\Lambda\subseteq \Fix(C)$. Then $A\defeq BCB^{-1}$ is an automorphism of the affine plane $\Pi$ such that $A(x)=BCB^{-1}(x)=BB^{-1}(y)=y$. Also for every $z\in L$ we have $B^{-1}(z)\in B^{-1}[L]=\Lambda\subseteq \Fix(C)$ and hence $A(z)=BCB^{-1}(z)=BB^{-1}(z)=z$, so $L\subseteq \Fix(A)$, witnessing that $(\delta,L)\in\SC_\Pi$.
\end{proof}

\begin{lemma}\label{l:Phi-v-trans} Let $\Pi$ be an affine plane and $v,h\in\partial\Pi$ be two distinct directions such that $\{v\}\times h\subseteq \SC_\Pi$. Then for any  line $V\in v$ and distinct points $x,x'\in V$, there exists a translation $T:\Pi\to\Pi$ such that $T(x)=x'$.
\end{lemma}

\begin{proof} If $|\Pi|_2\le 4$, then the affine plane $\Pi$ is Desarguesian (by Corollary~\ref{c:p5-Pappian}), and hence Thalesian and translation, by Theorem~\ref{t:ADA=>AMA}. So, assume that $|\Pi|_2\ge 5$.

Given any line $V\in v$ and distinct points $x,x'\in V$, choose distinct lines $A,B,C\in h$ such that $\{x,x'\}\cap(A\cup B\cup C)=\varnothing$. Such lines exist because $|\Pi|_2\ge 5$. Next, choose any disjoint lines $\Lambda,\Lambda'\notin h\cup v$ such that $x\in \Lambda$ and $x'\in \Lambda'$. Since $\Lambda,\Lambda'\notin h$, there exist unique points $a\in \Lambda\cap A$, $b\in \Lambda\cap B$, $c\in \Lambda\cap C$, $a'\in \Lambda'\cap A$, $b'\in \Lambda'\cap B$ and $c'\in \Lambda'\cap C$. Since the lines $\Aline a{b'}$ and $\Aline a{c'}$ are concurrent and $a\notin V$, one of the sets  $\Aline a{b'}\cap V$ or $\Aline a{c'}\cap V$ is not empty. We lose no generality assuming that $\Aline a{b'}\cap V$ is not empty and hence contains some point $o$. It follows from $\{x,x'\}\cap (A\cup B)=\varnothing$ that $o\notin A\cup B\cup\{x,x'\}$. 

\begin{picture}(100,150)(-100,-15)

\put(0,0){\color{teal}\line(1,0){170}}
\put(175,-3){$A$}
\put(10,0){\color{blue}\line(1,1){115}}
\put(128,120){$\Lambda'$}
\put(50,0){\color{blue}\line(1,1){115}}
\put(168,120){$\Lambda$}
\put(0,80){\color{teal}\line(1,0){170}}
\put(175,77){$B$}
\put(0,100){\color{teal}\line(1,0){170}}
\put(175,97){$C$}
\put(70,0){\color{cyan}\line(0,1){115}}
\put(68,120){$V$}
\put(50,0){\line(1,2){40}}

\put(10,0){\circle*{3}}
\put(8,-9){$a'$}
\put(50,0){\circle*{3}}
\put(46,-8){$a$}
\put(70,0){\circle*{3}}
\put(70,20){\circle*{3}}
\put(72,15){$x$}
\put(70,40){\circle*{3}}
\put(72,36){$o$}
\put(70,60){\circle*{3}}
\put(61,60){$x'$}
\put(70,80){\circle*{3}}
\put(90,80){\circle*{3}}
\put(86,83){$b'$}
\put(130,80){\circle*{3}}
\put(126,82){$b$}
\put(110,100){\circle*{3}}
\put(107,103){$c'$}
\put(150,100){\circle*{3}}
\put(145,103){$c$}

\end{picture}

Since $\{(v,A),(v,B)\}\subseteq\SC_\Pi$, there exist automorphisms $\Phi_A$ and $\Phi_B$ of the plane $\Pi$ such that $\Phi_A(x)=o$, $\Phi_B(o)=x'$,  $A\subseteq\Fix(\Phi_A)$, and $B\subseteq\Fix(\Phi_B)$. 
Then the automorphism $\Phi\defeq\Phi_B\circ\Phi_A$ of the plane $\Pi$ maps the point $x$ to the point $x'$. We claim that $\Phi$ is a translation of the plane $\Pi$. 

\begin{claim}\label{cl:Phi-h} $\{\Phi_A[L]:L\in h\}\cup\{\Phi_B[L]:L\in h\}\cup\{\Phi[L]:L\in h\}\subseteq h$.
\end{claim}

\begin{proof} For every $L\in h$ we have $A\parallel L\parallel B$ and hence $L\parallel A=\Phi_A[A]\parallel \Phi_A[L]$ and $L\parallel B=\Phi_B[B]\parallel \Phi_B[L]$, which implies $\{\Phi_A[L]:L\in h\}\cup\{\Phi_B[L]:L\in h\}\subseteq h$ and finally, $\{\Phi[L]:L\in h\}=\{\Phi_B\Phi_A[L]:L\in h\}\subseteq h$.
\end{proof}

\begin{claim}\label{cl:Phi-v} For every line $L\in v$, we have $L=\Phi_A[L]=\Phi_B[L]=\Phi[L]$.
\end{claim}

\begin{proof} Find unique points $\alpha\in V\cap A$ and $\beta\in V\cap B$ and observe that $\Phi_A[V]=\Phi_A[\Aline \alpha x]=\Aline\alpha o=V$ and $\Phi_B[V]=\Phi_B[\Aline o\beta]=\Aline {x'}\beta=V$. Then for every $L\in v$ we have $L\parallel V=\Phi_A[V]$ and hence $\Phi_A[L]\parallel \Phi_A[V]=V\parallel L$. By analogy we can prove that $\Phi_B[L]\parallel L$. Next, find unique ponts $\alpha'\in L\cap A$ and $\beta'\in L\cap B$ and observe that $\alpha'=\Phi_A(\alpha')\in L\cap\Phi_A[L]$ and  $\beta'=\Phi_B(\beta')\in L\cap\Phi_A[L]$. Since $\Phi_A[L]\parallel L\parallel \Phi_B[L]$, this implies $\Phi_A[L]=L=\Phi_B[L]$ and finally, $\Phi[L]=\Phi_B[\Phi_A[L]]=\Phi_B[L]=L$.
\end{proof}

\begin{claim}\label{cl:Phi[L]neL} For every line $L\in \Lambda_\parallel$, the line $\Phi[L]$ is disjoint with the line $L$.
\end{claim}

\begin{proof} Given any line $L\in\Lambda_\parallel$, find a unique point $\alpha\in L\cap A$. It follows from $L\parallel \Lambda$ that $\Phi_A[L]\parallel\Phi_A[\Lambda]=\Phi_A[\Aline ax]=\Aline ao\nparallel L$. By the Proclus Postulate~\ref{p:Proclus-Postulate}, there exists a point $\beta\in B\cap\Phi_A[L]$. It follows from $\alpha\in A\cap L$ and $A\subseteq\Fix(\Phi_A)$ that $\Phi_A(\alpha)=\alpha$ and hence $\Phi_A[L]=\Aline \alpha\beta$. Since $\Phi_A[L]\nparallel L$, $\beta\notin L$. It follows from $\Phi_A[L]=\Aline \alpha\beta\parallel \Aline ao=\Aline o{b'}$ that $\Phi[L]=\Phi_B[\Phi_A[L]]\parallel \Phi_B[\Aline o{b'}]=\Aline{x'}{b'}=L'\parallel L$. Since $\beta=\Phi_B(\beta)\in \Phi_B[\Phi_A[L]]=\Phi[L]$ and $\beta\notin L$, the parallel lines $\Phi[L]$ and $L$ are disjoint.
\end{proof} 

\begin{claim}\label{cl:Phi-no-fixed-points} The automorphism $\Phi$ of the plane $\Pi$ has no fixed points.
\end{claim}

\begin{proof} Given any point $p\in\Pi$, find a unique line $L\in\Lambda_\parallel$ such that $p\in L$. Claim~\ref{cl:Phi[L]neL} ensures that $\Phi(p)\in \Phi[L]\subseteq \Pi\setminus L\subseteq \Pi\setminus\{p\}$ and hence $\Phi(p)\ne p$.  
\end{proof}

Now we are able to prove that the automorphism $\Phi$ is a translation of the plane $\Pi$. By Claim~\ref{cl:Phi-no-fixed-points}, the automorphism $\Phi$ has no fixed points. Assuming that $\Phi$ is not a translation, we conclude that there exists a line $L$, which is concurrent to its image $\Phi[L]$. Claims~\ref{cl:Phi-h} and \ref{cl:Phi-v} ensure that $L\notin h\cup v$. Let $p$ be the unique point of the intersection $L\cap \Phi[L]$.  Let $V'\in v$ be a unique ``vertical'' line containing the point $p$. It follows from $L\notin v$ that $L\cap V'=\{p\}$. Claim~\ref{cl:Phi-v} ensures that $\Phi[V']=V'$ and hence   $p\in \Phi[L]\cap V'=\Phi[L\cap V']=\Phi[\{p\}]$ and hence $p=\Phi(p)$, which contradicts Claim~\ref{cl:Phi-no-fixed-points}. This contradiction shows that the automorphism $\Phi$ is a translation.
\end{proof}

\begin{lemma}\label{l:trans-swap} Let $\Pi$ be an affine plane and $v,h\in\partial\Pi$ be two directions such that $\{v\}\times h\subseteq\SC_\Pi$. Then $\Pi$ is a translation plane and $\{h\}\times v\subseteq \SC_\Pi$.
\end{lemma}

\begin{proof} Given any line $V\in v$, choose an affine base $uow$ in the affine plane $\Pi$ such that $\Aline ou\in h$ and $\Aline ow=V$. Lemma~\ref{l:Phi-v-trans} ensures that the based affine plane $(\Pi,uow)$ is vertical-translation. Since $(v,\Aline ou)\in \SC_\Pi$, the based affine plane $(\Pi,uow)$ is vertical-scale. By Corollary~\ref{c:vtrans+vscale}, the based affine plane $(\Pi,uow)$ is translation and horizontal-scale, which implies $\{h\}\times v\subseteq \SC_\Pi$.
\end{proof}

\begin{lemma}\label{l:SC=>translation} An affine plane $\Pi$ is a translation plane if $\{(\delta,L),(\delta,\Lambda)\}\subseteq\SC_\Pi$ for some direction $\delta$ and disjoint lines $L,\Lambda$.
\end{lemma}

\begin{proof} Assume that  $\{(\delta,L),(\delta,\Lambda)\}\subseteq\SC_\Pi$ for some direction $\delta\in\partial \Pi$ and disjoint lines $L,\Lambda$ in $\Pi$. Lemma~\ref{l:h-rozmaznia} ensures that $\{\delta\}\times L_\parallel\subseteq\SC_\Pi$. By Lemma~\ref{l:trans-swap}, $\Pi$ is a translation plane.
\end{proof}

\begin{lemma}\label{l:scale=>shear} A line $L$ in an affine plane $\Pi$ is a shear line if $\{(v,L),(v',L)\}\subseteq \SC_\Pi$ for two distinct directions $v,v'\in\partial\Pi$.
\end{lemma}

\begin{proof} Given two distinct points $x,x'\in \Pi$ with $\Aline x{x'}\cap L=\varnothing$, we should find an automorphism $A$ of the affine plane $\Pi$ such that $A(x)=x'$ and $L\subseteq\Fix(A)$. Find unique lines $\Lambda\in\mathcal L_x\cap v$ and $\Lambda'\in\mathcal L_{x'}\cap v'$.
Since $v\ne v'$, there exists a point $y\in \Lambda\cap\Lambda'$.

If $y\notin L$, then we can find automorphisms $\Phi,\Psi$ of the affine plane $\Pi$ such that $\Phi(x)=y$, $\Psi(y)=x'$ and $L\subseteq \Fix(\Phi)\cap\Fix(\Psi)$. Then the automorphism $A\defeq\Psi\circ\Phi$ has the desired properties: $A(x)=x'$ and $L\subseteq\Fix(A)$.

If $y\in L$, then choose any point $o\in \Aline x{x'}\setminus\{x,x'\}$ and find points $z\in \Lambda$ and $z'\in\Lambda'$ such that $\Aline oz\in v'$ and $\Aline o{z'}\in v$. It follows from $y\in L$ that $z,z'\notin L$. Repeating the above argument, find two automorphisms $A,A'$ of the plane $\Pi$ such that $A(x)=o$, $A'(o)=x'$, and $L\subseteq\Fix(A)\cap\Fix(A')$. Then the automorphism $A'A$ has the required properties: $A'A(x)=x'$ and $L\subseteq\Fix(A'A)$. 
\end{proof}  

\begin{lemma}\label{l:SC<=Desarg} An affine plane $\Pi$ is Desarguesian if there exist distinct directions $v,h\in\partial \Pi$ such that $(h,L)\in\SC_\Pi$ and  $\{v\}\times (\mathcal L_o\setminus v)\subseteq \SC_\Pi$ for some line $L\in\mathcal L_\Pi$ and point $o\in\Pi$.
\end{lemma}

\begin{proof} Separately we consider four possible cases.
\smallskip

1. First we assume that $o\in L\in v$. The definition of the set $\SC_\Pi\ni (h,L)$ ensures that $v=L_\parallel\ne h$. Consider the unique line $H\in\mathcal L_o\cap h$ and take any points $u\in H\setminus\{o\}$ and $w\in L\setminus\{o\}$. Then $uow$ is an affine base such that the based affine plane $(\Pi,uow)$ is horizontal-scale and vertical-scale (because $\{(h,L),(v,H)\}\subseteq\SC_\Pi$). Let $e$ be the diunit of the affine base $uow$. Since $(v,\Aline oe)\in \SC_\Pi$, there exists an automorphism $\Phi$ of the plane $\Pi$ such that $\Phi(u)\in \Aline ue\setminus\{u\}$ and $\Aline oe\subseteq\Fix(\Phi)$. It follows from $L\parallel \Aline ue$ that $\Phi[L]\parallel \Phi[\Aline ue]=\Aline ue\parallel L$ and hence $\Phi[L]=L$. Lemma~\ref{l:SC-auto} ensures that $(\Phi[\Aline ou]_\parallel,L)=(\Phi[\Aline ou]_\parallel,\Phi[V])\in\SC_\Pi$. Since $\{(h,L),(\Phi[\Aline ou],L)\}\subseteq\SC_\Pi$, we can apply Lemma~\ref{l:scale=>shear} and conclude that the line $L$ is shear and hence the based affine plane $(\Pi,uow)$ is vertical-shear, horizontal-scale and vertical-scale. By Corollary~\ref{c:vshear+vh-scale+htrans}, the ternar $\Delta$ of the based affine plane $(\Pi,uow)$ is linear, distributive and associative. By Theorem~\ref{t:corps<=>}, the affine plane $\Pi$ is Desarguesian.
\smallskip

2. Next, assume that $o\in L\notin v$. Consider the unique lines $H\in\mathcal L_o\cap  h$ and $V\in\mathcal L_o\cap v$. The definition of the relation $\SC_\Pi\ni (h,L)$ ensures that $L\notin h$.  Choose any point $x\in L\setminus\{o\}$ and observe that $L_\parallel \ne h=H_\parallel$ and $L\notin v$ imply $x\notin H\cup V$. Choose any point $y\in \Pi\setminus H$ such that $\Aline xy\in v$. Since $(v,H)\in\SC_\Pi$, there exists an automorphism $\Phi$ of  the plane $\Pi$ such that $\Phi(x)=y$ and $o\in H\subseteq\Fix(\Phi)$. It is easy to show $\Phi[V]\parallel \Phi[\Aline xy]=\Aline xy\parallel V$ and hence $\Phi[V]=V$. Find a unique point $z\in V$ such that $\Aline yz\in h$. Since $(h,L)\in\SC_\Pi$, there exists an automorphism $\Psi$ of the plane $\Pi$ such that $\Psi(y)=z$ and $L\subseteq\Fix(\Psi)$. Consider the automorphism $A\defeq\Psi\circ\Phi$ and its extension $\bar A$ to the spread completion $\overline\Pi=\Pi\cup\partial\Pi$. Observe that $A[L]=\Psi[\Phi[L]]=\Psi[\Phi[\Aline ox]]=\Psi[\Aline oy]=\Aline oz=V$. Lemma~\ref{l:SC-auto} ensures that $(\bar A(h),V)=(\bar A(h),A[L])\in \SC_\Pi$ and hence $\bar A(h)\ne v$. Since $o\in V\in v$, the first case ensures that the plane $\Pi$ is Desarguesian.
\smallskip

3. Next assume that $o\notin L\in v$. Choose any point $x\in \Pi\setminus L$ such that $\Aline ox\in h$. Since $(h,L)\in\SC_\Pi$, there exists an automorphism $\Phi$ of the plane $\Pi$ such that $\Phi(o)=x$ and $L\subseteq \Fix(\Phi)$. Let $\bar\Phi$ be the extension of $\Pi$ to an automorphism of the spread completion of the plane $\Pi$. It follows from $\Phi[L]=L\in v$ that $\bar\Phi(v)=v$. Lemma~\ref{l:SC-auto} ensures that $\{v\}\times(\mathcal L_x\setminus \{v\})=\{(\bar\Phi(v),\Phi[\Lambda]):\Lambda\subseteq \mathcal L_o\setminus v\}\subseteq\SC_\Pi$.
Choose any disjoint lines $\Lambda\in\mathcal L_o\setminus v$ and $\Lambda'\in\mathcal L_x\setminus v$. Since $\{(v,\Lambda),(v,\Lambda')\}\subseteq\SC_\Pi$, we can apply Lemma~\ref{l:SC=>translation} and conclude that $\Pi$ is a translation palne. Then there exists a translation $T$ of the plane $\Pi$ such that $o\in T[L]$. Let $\bar T$ be the extension of the translation $T$ to a unique automorphism of the spread completion of the plane $\Pi$. Since translations do no change directions, $T[L]\in \bar T(v)=v$ and $\bar T(h)=h$. Lemma~\ref{l:SC-auto} ensures that $(h,T[L])=(\bar T(h),T[L])\in \SC_\Pi$. Since $o\in T[L]\in v$, we can apply the first case and conclude that the plane $\Pi$ is Desarguesian.
\smallskip

4. Finally, assume that $o\notin L\notin v$. Since $|\mathcal L_o|=|\Pi|_2+1\ge 4$, there exists a line $\Lambda\in\mathcal L_o\setminus(v\cup h\cup L_\parallel)$. Consider a unique point $x\in \Lambda\cap L$ and find a unique point $y\in \Pi$ such that $\Aline oy\in h$ and $\Aline xy\in v$. Since $(h,L)\in\SC_\Pi$, there exists an automorphism $\Phi$ of the plane $\Pi$ such that $\Phi(o)=y$ and $L\subseteq\Fix(\Phi)$. Then $\Phi[\Lambda]=\Phi[\Aline xo]=\Aline xy\in v$. Let $\bar\Phi$ be the extension of $\Phi$ to an automorphism of the spread completion of the plane $\Pi$. It follows from $\Lambda\notin v$ that $\Aline xy=\Phi[\Lambda]\notin \bar\Phi(v)$ and hence $v=\Aline xy_\parallel\ne\bar\Phi(v)$. By Lemma~\ref{l:SC-auto}, $(v,\Lambda)\in \SC_\Pi$ implies $(\bar\Phi(v),\Phi[\Lambda])\in\SC_\Pi$. Since $\Phi[\Lambda]\in v$, we can apply the preceding case and conclude that the plane $\Pi$ is Desarguesian.
\end{proof}

\begin{lemma}\label{l:PSC-not-function2} Let $\Pi$ be a non-Desarguesian affine plane. If $\PSC_\Pi$ is not a function, then $\PSC_\Pi=\{\delta\}\times(\partial\Pi\setminus\{\delta\})$ for some direction $\delta\in\partial\Pi$.
\end{lemma} 

\begin{proof} If $\PSC_\Pi$ is not a function, then there exist three distinct directions $\delta,\delta',\delta''$ such that $\{(\delta,\delta'),(\delta,\delta'')\}\subseteq \PSC_\Pi$. Lemma~\ref{l:PSC-not-function} ensures that $\{\delta\}\times(\mathcal L_p\setminus \delta)\subseteq \SC_\Pi$ and hence $\{\delta\}\times(\partial\Pi\setminus\{\delta\})\subseteq\PSC_\Pi$. Assuming that $\{\delta\}\times(\partial\Pi\setminus\{\delta\})\ne\PSC_\Pi$, we can find a pair $(h,v)\in\PSC_\Pi$ such that $(h,v)\notin\{\delta\}\times(\partial\Pi\setminus\{\delta\})$. Assuming that $h=\delta$ and taking into account that $v\ne h=\delta$, we conclude that $(h,v)\in\{\delta\}\times(\partial\Pi\setminus \delta)$, which contradicts the choice of $(h,v)$. This contradiction shows that $h\ne \delta$. Since $(h,v)\in\PSC_\Pi$, there exists a line $L\in v$ such that $(h,L)\in\SC_\Pi$. By Lemma~\ref{l:SC<=Desarg}, the plane $\Pi$ is Desarguesian, which contradicts our assumption. This contradiction shows that $\{\delta\}\times(\partial\Pi\setminus\{\delta\})=\PSC_\Pi$.
\end{proof}

\begin{lemma}\label{l:PCS-not-function} Let $\Pi$ be an affine plane such that $\{(v,\delta),(h,\delta)\}\subseteq\PSC_\Pi$ for some distinct directions $v,h,\delta\in\partial\Pi$. Then $(\partial\Pi\setminus\{\delta\})\times\{\delta\}\subseteq\PSC_\Pi$.
\end{lemma}

\begin{proof} Since $\{(v,\delta),(h,\delta)\}\subseteq\PSC_\Pi$, there exist lines $L,\Lambda\in\delta$ such that $\{(v,L),(h,\Lambda)\}\subseteq\SC_\Pi$. To check that $(\partial\Pi\setminus\{\delta\})\times\{\delta\}\subseteq\PSC_\Pi$, take any direction $d\in \partial\Pi\setminus\{\delta\}$. We have to show that $(d,\delta)\in \PSC_\Pi$. This is clear if $d\in\{v,h\}$. So, assume that $d\notin\{h,v\}$. Fix any point $o\in L$. Since $L\in\delta\ne h$, there exists a point $x\in\Pi\setminus L$ such that $\Aline ox\in h$. Find unique lines $D\in\mathcal L_o\cap d$ and $V\in\mathcal L_x\cap v$. Since $d\ne v$, there exists a unique point $y\in D\cap V$. It is easy to see that $y\notin L\cup\{x\}$. Since $(v,L)\in\SC_\Pi$, there exists an automorphism $\Phi$ of the plane $\Pi$ such that $\Phi(x)=y$ and $L\subseteq\Fix(\Phi)$. Then $\Phi[\Aline ox]=\Aline oy=D\in d$.
Let $\bar\Phi$ be the extension of the automorphism $\Phi$ to an automorphism of the spread completion $\overline \Pi=\Pi\cup\partial\Pi$ of $\Pi$. Observe that $\bar\Phi(h)=\bar\Phi(\Aline ox_\parallel)=\Phi[\Aline ox]_\parallel=D_\parallel=d$. Lemma~\ref{l:SC-auto} ensures that $(d,\Phi[\Lambda])=(\bar\Phi(h),\Phi[\Lambda])\in \SC_\Pi$. It follows from $\Lambda\parallel L=\Phi[L]$ that $\Phi[\Lambda]\parallel \Phi[L]=L\in \delta$ and hence $\Phi[\Lambda]\in\delta$. Then $(d,\Phi[L])\in\SC_\Pi$ implies $(d,\delta)\in\PSC_\Pi$, witnessing that $(\partial\Pi\setminus\delta)\times\{\delta\}\subseteq\PSC_\Pi$.
\end{proof}

\begin{lemma}\label{l:PCS-not-function2} Let $\Pi$ be a non-Desarguesian affine plane. If $\PSC_\Pi^{-1}$ is not a function, then $\PSC_\Pi=(\partial\Pi\setminus\{\delta\})\times\{\delta\}$ for some direction $\delta\in\partial\Pi$.
\end{lemma} 

\begin{proof} Assume that $\PSC^{-1}_\Pi$ is not a function. Then there exist three distinct directions $\delta,\delta',\delta''$ such that $\{(\delta',\delta),(\delta'',\delta)\}\subseteq\PSC_\Pi$.  Lemma~\ref{l:PCS-not-function} ensures that $(\partial\Pi\setminus\{\delta\})\times\{\delta\}\subseteq\PSC_\Pi$.
Assuming that $(\partial\Pi\setminus\{\delta\})\times\{\delta\}\ne\PSC_\Pi$, we can find a pair $(h,v)\in \PSC_\Pi$ such that $(h,v)\notin (\partial\Pi\setminus\{\delta\})\times\{\delta\}$ and hence $v\ne\delta$. Since $\{(h,v),(h,\delta)\}\subseteq\SC_\Pi$, the relation $\PSC_\Pi$ is not a function. By Lemma~\ref{l:PSC-not-function2}, $\PSC_\Pi=\{h\}\times(\partial\Pi\setminus\{h\})$, which contradicts $(\partial\Pi\setminus\{\delta\})\times\{\delta\}\subseteq\PSC_\Pi$. This contradiction shows that $(\partial\Pi\setminus\{\delta\})\times\{\delta\}=\PSC_\Pi$.
\end{proof}

\begin{lemma}\label{l:SC-symmetry} Let $\Pi$ be an affine plane such that $\PSC_\Pi$ is an injective function. If $(\delta,L),(\delta',L')\in\SC_\Pi$ are two pairs with $L'\in\delta$, then $L\in\delta'$.
\end{lemma}

\begin{proof} To derive a contradiction, assume that $L'\in\delta$ but $L\notin \delta'$. The definition of the relation $\SC_\Pi$ ensures that $L\notin \delta=L'_\parallel \ne\delta'$. Then the lines $L,L'$ are concurrent and hence have a common point $o$.  Choose any point $x\in \Pi\setminus L'$ such that $\Aline ox\in \delta'$. It follows from $o\in L\notin\delta'$ that $x\notin L$. Choose any point $y\in\Pi\setminus L$ such that $\Aline xy\in \delta$. It follows from $\Aline ox_\parallel=\delta'\ne\delta$ that $\Aline oy\notin \delta'$. 

Since $(\delta,L)\in\SC_\Pi$, there exists an automorphism $\Phi$ of the plane $\Pi$ such that $\Phi(x)=y$ and $o\in L\subseteq\Fix(\Phi)$. Then $\Phi[\Aline ox]=\Aline oy$. Let $\bar\Phi$ be the extension of the automorphism to the spread completion of the affine plane $\Pi$. Then $\bar\Phi(\delta')=\bar\Phi(\Aline ox_\parallel)=\Phi[\Aline ox]_\parallel=\Aline oy_\parallel\ne \delta'$. 

We claim that $\Phi[L']=L'$. Since $L\notin\delta=\Aline xy_\parallel$, there exists a unique point $p\in L\cap \Aline xy$. Then $\Phi[\Aline xy]=\Phi[\Aline px]=\Aline py=\Aline xy$ and $L'\parallel\Aline xy$ imply $\Phi[L']\parallel\Phi[\Aline xy]=\Aline xy\parallel L'$ and hence $\Phi[L']=L'$ (because $o\in L'\cap \Phi[L']$). 
 Lemma~\ref{l:SC-auto} ensures that $(\bar\Phi(\delta'),L')=(\bar\Phi(\delta'),\Phi[L'])\in\SC_\Pi$ and hence $\{(\delta',L'_\parallel),(\bar\Phi(\delta'),L'_\parallel)\}\subseteq\PSC_\Pi$, which contradicts the injectivity of the function $\PSC_\Pi$. This contradiction shows that $L\in\delta'$.
\end{proof}

\begin{lemma}\label{l:SC-symmetry2} Let $\Pi$ be an affine plane and $(v,h),(v',h')\in\PSC_\Pi$ be two pairs. If $\PSC_\Pi$ is an injective function and $\{v,h\}\cap\{v',h'\}\ne\varnothing$, then $\{v,h\}=\{v',h'\}$.\end{lemma}

\begin{proof} Assume that $\{v,h\}\cap\{v',h'\}\ne\varnothing$. If $(v,h)=(v',h')$, then $\{v,h\}=\{v',h'\}$ and we are done. So, assume that $(v,h)\ne(v',h')$. Since $\PSC_\Pi$ is an injective function, the inequality $(v,h)\ne (v',h')$ implies $v\ne v'$ and $h\ne h'$.

Then $\{v,h\}\cap\{v',h'\}\ne\varnothing$ implies $v=h'$ or $h=v'$. We lose no generality assuming that $v=h'$.  Since $\{(v,h),(v',h')\}\subseteq\PSC_\Pi$, there exist lines $H\in h$ and $H'\in h'$ such that $\{(v,H),(v',H')\}\subseteq\SC_\Pi$. It follows from $v=h'$ that $H'\in h'=v$ and hence $h=H_\parallel =v'$, by Lemma~\ref{l:SC-symmetry}. Then $\{h',v'\}=\{v,h\}$. 
\end{proof}

\begin{lemma}\label{l:PSC-involution} Let $\Pi$ be an affine plane. If $\PSC_\Pi$ is an injective function with $\dom[\PSC_\Pi]=\rng[\PSC_\Pi]$, then $\PSC_\Pi$ is an involution.
\end{lemma}

\begin{proof} First we show that $\PSC_\Pi^{-1}\subseteq\PSC_\Pi$. Given any pair $(\delta,\delta')\in \PSC_\Pi$, we need to check that $(\delta',\delta)\in\PSC_\Pi$. Since $(\delta,\delta')\in \PSC_\Pi$, there exists a line $L'\in\delta'$ such that $(\delta,L')\in\SC_\Pi$. Since $\delta'\in\rng[\PSC_\Pi]=\dom[\PSC_\Pi]=\dom[\SC_\Pi]$, there exists a line $L$ such that $(\delta',L)\in\SC_\Pi$.  Since $L'\in\delta'$, we can apply Lemma~\ref{l:SC-symmetry} and  conclude that $L\in\delta$. Then $(\delta',\delta)=(\delta',L_\parallel)\in\PSC_\Pi$, witnessing that $\PSC_\Pi^{-1}\subseteq\PSC_\Pi$. 

To see that $\PSC_\Pi\subseteq\PSC_\Pi^{-1}$, take any pair $(\delta,\delta')\in\PSC_\Pi$.
Since $\delta\in\dom[\PSC_\Pi]=\rng[\PSC_\Pi]$, there exists a direction $\delta''$ such that $(\delta'',\delta)\in\PSC_\Pi$. Since $\PSC_\Pi$ is a function, $\{(\delta,\delta'),(\delta,\delta'')\}\subseteq\PSC_\Pi\cup\PSC_\Pi^{-1}=\PSC_\Pi$ implies $\delta'=\delta''$. Then $(\delta,\delta')=(\delta,\delta'')\in \PSC_\Pi^{-1}$, witnessing that $\PSC_\Pi\subseteq\PSC_\Pi^{-1}$ and hence $\PSC_\Pi^{-1}=\PSC_\Pi$. 
\end{proof}

\begin{lemma}\label{l:SC-largedom} Let $\Pi$ be an affine plane whose hyperscale figure $\PSC_\Pi$ is an injective function. If $|\dom[\PSC_\Pi]\cup\rng[\PSC_\Pi]|>2$, then $\dom[\PSC_\Pi]=\partial\Pi=\rng[\PSC_\Pi]$. Moreover, there exists a point $p\in\Pi$ such that $\dom[\SC_\Pi\cap(\partial\Pi\times\mathcal L_p)]=\partial\Pi$ and $\rng[\SC_\Pi\cap(\partial\Pi\times\mathcal L_p)]=\mathcal L_p$.
\end{lemma}

\begin{proof} Fix any pair $(v,h)\in \PSC_\Pi$. Since $|\dom[\PSC_\Pi]\cup\rng[\PSC_\Pi]>2|$, there exists a pair $(v',h')\in\PSC_\Pi$ such that $\{v',h'\}\not\subseteq\{v,h\}$. Applying Lemma~\ref{l:SC-symmetry2}, we obtain that $\{v',h'\}\cap\{v,h\}=\varnothing$. 

Since $\{(v,h),(v',h')\}\subseteq\PSC_\Pi$, there exist lines $H\in h$ and $H'\in h'$ such that $\{(v,H),(v',H')\}\subseteq\SC_\Pi$.
Since $H_\parallel=h\ne h'=H'_\parallel$, the lines $H,H'$ are concurrent and hence contain a common point $p$. We claim that
$\dom[\SC_\Pi\cap(\partial\Pi\times\mathcal L_p)]=\partial\Pi$ and $\rng[\SC_\Pi\cap(\partial\Pi\times\mathcal L_p)]=\mathcal L_p$.

To prove the equality $\dom[\SC_\Pi\cap(\partial\Pi\times\mathcal L_p)]=\partial\Pi$, it suffices for every direction $\delta\in\partial\Pi$ to find a line $L\in\mathcal L_p$ such that $(\delta,L)\in\SC_\Pi$. If $\delta=v$ (resp. $\delta=v'$), then the line $L\defeq H$ (resp. $L\defeq H'$) has the required property. So, assume that $\delta\notin\{v,v'\}$. 
Since $h\ne h'$, either $\delta\ne h$ or $\delta\ne h'$. We lose no generality assuming that $\delta\ne h$. Since $H\in h\ne v'$, there exists a point $x\in \Pi\setminus H$ such that $\Aline px\in v'$. Since $v\ne\delta\ne v'$, there exists  a unique point $y\in\Pi$ such that $\Aline xy\in v$ and $\Aline py\in \delta$. Since $(v,H)\in\SC_\Pi$, there exists an automorphism $\Phi$ of the plane $\Pi$ such that $\Phi(x)=y$ and $H\subseteq\Fix(\Phi)$.  Let $\bar\Phi:\overline\Pi\to\overline\Pi$ be the extension of the automorphism $\Phi$ to an automorphism of the spread complestion of the affine plane $\Pi$. Then $\bar\Phi(v')=\bar\Phi(\Aline px_\parallel)=\Phi[\Aline px]_\parallel=\Aline py_\parallel=\delta$. By Lemma~\ref{l:SC-auto}, $(\delta,\Phi[H'])=(\bar\Phi(\delta'),\Phi[H'])\in\SC_\Pi$ and hence $(\delta,\Phi[H'])\in \SC_\Pi\cap(\partial\Pi\times\mathcal L_p)$. Therefore, $\partial\Pi=\dom[\SC_\Pi\cap(\partial\times\mathcal L_p)]$ and $\partial\Pi=\dom[\SC_\Pi]=\dom[\PSC_\Pi]$.
\smallskip

Next, we prove the equality $\rng[\SC_\Pi\cap(\partial\Pi\times\mathcal L_p)]=\mathcal L_p$. Given any line $L\in\mathcal L_p$, we need to find a direction $\delta$ such that $(\delta,L)\in\SC_\Pi$. If $L=H$ (resp. $L=H'$), then the direction $\delta\defeq v$ (resp. $\delta\defeq v'$) has the required property. So, assume that $H\ne L\ne H'$ and hence $h\ne L_\parallel\ne h'$. Since $h\ne h'$, either $L_\parallel\ne h$ or $L_\parallel \ne h'$. We lose no generality assuming that $L_\parallel\ne h$. Choose any point $x\in H'\setminus\{p\}$ and find a unique point $y\in L$ such that $\Aline xy\in v$. Since $(v,H)\in\SC_\Pi$, there exists an automorphism $\Psi$ of the plane $\Pi$ such that $\Psi(x)=y$ and $p\in H\subseteq\Fix(\Psi)$. Then $\Psi[H']=\Psi[\Aline px]=\Aline py\in L$.  Let $\bar\Phi$ be the extension of the automorphism to the spread completion $\overline\Pi=\Pi\cup\partial\Pi$ of the affine plane $\Pi$. Lemma~\ref{l:SC-auto} ensures that $(\bar\Psi(v'),L)=(\bar\Psi(v'),\Psi[H'])\in\SC_\Pi$, witnessing that $L\in \rng[\SC_\Pi\cap(\partial\Pi\times \mathcal L_p)]$. Therefore, $\mathcal L_p=\rng[\SC_\Pi\cap(\partial\Pi\times\mathcal L_p)]\subseteq\rng[\SC_\Pi]$ and $\partial \Pi=\{L_\parallel:L\in\mathcal L_p\}=\rng[\PSC_\Pi]$. 
\end{proof}

\begin{lemma}\label{l:Desarg=>SC=} If $\Pi$ is a Desarguesian affine plane, then $$\SC_\Pi=\{(\delta,L)\in\partial\Pi\times\mathcal L_\Pi:L\notin\delta\}\quad\mbox{and}\quad\PSC_\Pi=\{(v,h)\in\partial\Pi\times\partial\Pi:v\ne h\}.$$
\end{lemma}

\begin{proof} By Theorem~\ref{t:corps<=>}, every ternar of the Desarguesian plane $\Pi$ is linear and associative. The definition of the set $\SC_\Pi$ ensures that  $\PSC_\Pi\subseteq\{(\delta,L)\in\partial\Pi\times\mathcal L_\Pi:L\notin\delta\}$. To prove that this inclusion is an equality, take any direction $h\in\partial\Pi$ and line $V\in\mathcal L_\Pi\setminus\delta$. Choose an affine base $uow$ in the plane $\Pi$ such that $\Aline ow=V$ and $\Aline ou\in h$. By Theorem~\ref{t:horizontal-scale}, the based affine plane $(\Pi,uow)$ is horizontal-scale, which implies $(h,V)\in\SC_\Pi$. Therefore, $\ SC_\Pi=\{(\delta,L)\in\partial\Pi\times\mathcal L_\Pi:L\notin\delta\}$ and hence $\PSC_\Pi=\{(v,h)\in\partial\Pi\times\partial\Pi:v\ne h\}$.
\end{proof}

\begin{lemma}\label{l:PSC=2} Let $\Pi$ be an affine plane. If $|\PSC_\Pi|=2$, then $\PSC_\Pi=\{(\delta,\delta'),(\delta',\delta)\}$ for some distinct directions $\delta,\delta'\in\partial\Pi$.
\end{lemma}

\begin{proof} Since $|\PSC_\Pi|=2$, the affine plane $\Pi$ is not Desarguesian, by Lemma~\ref{l:Desarg=>SC=}. If $\PSC_\Pi$ is not a function, then by Lemma~\ref{l:PSC-not-function}, for some $\delta\in\partial\Pi$ we have $3\le|\{\delta\}\times(\partial\Pi\setminus\{\delta\})|\le|\PSC_\Pi|$, which contradicts $|\PSC_\Pi|=2$.  If $\PSC_\Pi^{-1}$ is not a function, then by Lemma~\ref{l:PCS-not-function}, for some $\delta\in\partial\Pi$ we have $3\le|(\partial\Pi\setminus\{\delta\})\times\{\delta\}|\le|\PSC_\Pi|$, which contradicts $|\PSC_\Pi|=2$. Therefore, $\PSC_\Pi$ is an injective function and hence $|\dom[\PSC_\Pi]|=|\rng[\PSC_\Pi]|=|\PSC_\Pi|=2$. If $|\dom[\PSC_\Pi]\cup\rng[\PSC_\Pi]|>2$, then Lemma~\ref{l:SC-largedom} ensures that $4\le|\partial\Pi|=|\dom[\PSC_\Pi]|=|\PSC_\Pi|$, which contradicts $|\PSC_\Pi|=2$. This contradiction shows that $|\dom[\PSC_\Pi]\cup\rng[\PSC_\Pi]|=2$ and hence $\dom[\PSC_\Pi]=\rng[\PSC_\Pi]=\{\delta,\delta'\}$ for some distinct directions $\delta,\delta'$. Since the relation $\PSC_\Pi$ is irreflexive, $\PSC_\Pi=\{(\delta,\delta'),(\delta',\delta)\}$. 
\end{proof}

\begin{lemma}\label{l:PSC-all} Let $\Pi$ be an affine plane, $o\in\Pi$ be a point and $H,V\in\mathcal L_o$ be two distinct lines such that $\{(H_\parallel,V),(V_\parallel,H)\}\subseteq\SC_\Pi$. Let $\Phi:\Pi\to\Pi$ be an automorphism of the affine plane $\Pi$ such that $H\subseteq\Fix(\Phi)$, $\Phi[V]=V$ and $(A_\parallel,\Phi[A])\in\SC_\Pi$ for some line $A\in\mathcal L_o\setminus\{H,V\}$. Then $(B_\parallel,\Phi[B])\in\SC_\Pi$ for all lines $B\in\mathcal L_o\setminus\{H,V\}$.
\end{lemma} 

\begin{proof} Given any line $B\in\mathcal L_o\setminus\{H,V,A\}$, we should prove that $(B_\parallel,\Phi[B])\in\SC_\Pi$. Fix any point $a\in A\setminus\{o\}$ and consider the point $a'\defeq\Phi(a)$ on the line $A'\defeq\Phi[A]\ne A$. It follows from $o\in H\subseteq\Fix(\Phi)$ that $A'=\Aline o{a'}$. Since $B\cap A=\{o\}$, there exists a unique point $b\in B$ such that $\Aline ab\parallel H$. Since the plane $\Pi$ is Playfair, there exist unique points $\alpha,\beta\in H$ such that $\Aline a\alpha\parallel V\parallel b\beta$.  
It follows from $H\subseteq\Fix(\Phi)$ and $V=\Phi[V]$ that $\Phi[\Aline a\alpha]=\Aline a\alpha$ and $\Phi[\Aline b\beta]=\Aline b\beta$. In particular, $a'\in\Phi(a)\in \Aline a\alpha$ and $b'\defeq\Phi(b)\in\Aline b\beta$. Moreover, $\Aline{a'}{b'}=\Phi[\Aline ab]\parallel \Phi[H]=H$. 
Since $\Aline ab\in H_\parallel$ and $(H_\parallel,V)\in\SC_\Pi$, there exists an automorphism $\Psi$ of the plane $\Pi$ such that $\Psi(a)=b$ and $V\subseteq\Fix(\Psi)$. Then $\Psi[\Aline a{a'}]\parallel \Psi[V]=V$ and hence $\Psi[\Aline a{a'}]=\Aline {b}{b'}$. It follows from $V\subseteq\Fix(\Psi)$ and $H=\Psi[H]$ that $\Psi[\Aline {a'}{b'}]=\Aline{a'}{b'}$. Then $\Psi(a')\in \Aline{b}{b'}\cap\Aline{a'}{b'}=\{b'\}$, $\Psi[A]=\Psi[\Aline oa]=\Aline ob=B$ and $\Psi[A']=\Psi[\Aline o{a'}]=\Aline o{b'}=\Phi[B]$. Since $(A_\parallel, A')=(A_\parallel ,\Phi[A])\in \SC_\Pi$, we can apply Lemma~\ref{l:SC-auto} and conlcude  that $(B_\parallel,\Phi[B])=(\Psi[A]_\parallel,\Psi[A'])\in \SC_\Pi$.
\end{proof}

\begin{lemma}[Ostrom, 1957]\label{l:PSC-involution=>9} If for an affine plane $\Pi$, the relation $\PSC_\Pi$ is an involution with $\dom[\PSC_\Pi]=\partial\Pi=\rng[\PSC_\Pi]$, then $|\Pi|_2=9$ and $|\PSC_\Pi|=|\partial\Pi|=10$.
\end{lemma}

\begin{proof} Since $\PSC_\Pi$ is an involution, $\PSC_\Pi\ne \{(h,v)\in\partial\Pi\times\partial\Pi:h\ne v\}$ and hence the affine plane $\Pi$ is not Desarguesian, by Lemma~\ref{l:Desarg=>SC=}. By Corollary~\ref{c:4-Pappian}, $|\Pi|_2\ge 7$. Choose any pair of directions $(h,v)\in\PSC_\Pi$. Since $\PSC_\Pi$ is an involution, $(v,h)\in\PSC_\Pi$. Since $\{(h,v),(v,h)\}\subseteq\PSC_\Pi$, there exist lines $V\in v$ and $H\in h$ such that $\{(v,H),(h,V)\}\subseteq \SC_\Pi$. Since $v\ne h$, the lines $V$ and $H$ are concurrent and hence have a unique common point $o$. 

\begin{claim}\label{cl:SC-Delta} $\mathcal L_o\subseteq\rng[\SC_\Pi]$.
\end{claim}

\begin{proof} Given any line $D\in \mathcal L_o$, we should find a direction $\delta\in\partial\Pi$ such that $(\delta,D)\in \SC_\Pi$. If $D\in\{V,H\}$, then the direction $\delta\in\{v,h\}\setminus \{D_\parallel\}$ has the required property. So, assume that $V\ne D\ne H$. Since $D_\parallel\in\partial\Pi=\dom[\PSC_\Pi]$, there exists a direction $\delta\in\partial\Pi$ such that $(\delta,D_\parallel)\in \PSC_\Pi$. Since $(\delta,D_\parallel)\in\PSC_\Pi$, there exists a line $\Delta\in D_\parallel$ such that $(\delta,\Delta)\in\SC_\Pi$. If $\Delta=D$, then $(\delta,D)=(\delta,\Delta)\in\SC_\Pi$ and we are done. So, assume that $\Delta\ne D$ and hence $o\notin \Delta$. Since $D_\parallel \ne h=H_\parallel$, there exist a unique point $x\in H\cap \Delta$.  Since  $D_\parallel\ne h$, the direction $\delta=\PSC_\Pi(D_\parallel)$ is not equal to the direction $v=\PSC_\Pi(h)$. Then there exists a unique point $y\in\Pi$ such that $\Aline xy\in v$ and $\Aline oy\in\delta$. It follows from $o\notin\Delta$ that $y\notin\Delta$.  Since $(\delta,\Delta)\in\SC_\Pi$, there exists an automorphism $\Phi$ of the plane $\Pi$ such that $\Phi(o)=y$ and $\Delta\subseteq\Fix(\Phi)$. Let $\bar\Phi$ be the unique extension of $\Pi$ to an automorphism of the spread completion $\overline\Pi$ of  the plane $\Pi$. Lemma~\ref{l:SC-auto} ensures that $(\bar\Phi(v),\Phi[H])\in \SC_\Pi$. Observe that $\Phi[H]=\Phi[\Aline xo]=\Aline xy\in v$ and hence $(\bar\Phi(v),v)\in\PSC_\Pi$. Since $\PSC_\Pi$ is an injective  function, $\bar\Phi(v)=\PSC_\Pi^{-1}(v)=h$. Then $(h,\Aline xy)=(\bar\Phi(v),\Phi[H])\in\SC_\Pi$. Since $o\notin\Delta$, the point $x\in H\cap\Delta$ does not belong to the line $V$, which implies $V\cap\Aline xy=\varnothing$. Since  $\{(h,V),(h,\Aline xy)\}\subseteq \SC_\Pi$, we can apply Lemma~\ref{l:SC=>translation} and conclude that $\Pi$ is a translation plane. Then there exists a translation $T$ of $\Pi$ such that $o\in T[\Delta]$ and hence $T[\Delta]=D$. Let $\bar T$ be the unique extension of $T$ to an automorphism of the spread completion $\overline\Pi$ of  the plane $\Pi$. Since $T$ is a translation, $\bar T(\delta)=\delta$. Lemma~\ref{l:SC-auto} ensures that $(\delta,D)=(\bar T(\delta),T[\Delta])\in\SC_\Pi$.
\end{proof}

Choose any points $u\in H\setminus\{o\}$, $w\in V\setminus\{o\}$, and observe that $uow$ is an affine base in the affine plane $\Pi$. Let $e$ be the diunit of the affine base $uow$. Since $\{(v,\Aline ou),(h,\Aline ow)\}=\{(v,H),(h,V)\}\subseteq\SC_\Pi$, the based affine plane $(\Pi,uow)$ is vertical-scale and horizontal-scale. By Theorem~\ref{c:Algebra-vs-Geometry-in-trings:affine}(6), the ternar $\Delta$ of the based affine plane $(\Pi,\Delta)$ is linear, right-distributive and associative-dot.

By Claim~\ref{cl:SC-Delta}, for the diagonal $\Delta\defeq\Aline oe\in\mathcal L_o$ of the affine base $uow$, there exists a direction $\bar\delta$ such that $(\bar\delta,\Delta)\in\SC_\Pi$. Let $\delta=\Delta_\parallel$ be the direction of the diagonal $\Delta$. It follow from $\delta\ne h$ that $\bar\delta=\PSC_\Pi(\delta)\ne\PSC_\Pi(h)=v$. Then there exists a unique point $w'\in V$ such that $\Aline{u}{w'}\in\bar\delta$. Since $(\bar\delta,\Delta)\in\SC_\Pi$, there exists an automorphism $S$ of the plane $\Pi$ such that $S(u)=w'$ and $\Delta\subseteq\Fix(S)$. Let $\bar S$ be the extension of $S$ to an automorphism of the spread completion $\overline\Pi=\Pi\cup\partial\Pi$ of the plane $\Pi$. Observe that $S[\Aline ou]=\Aline o{w'}=V$ and hence $\bar S(h)=v$. Lemma \ref{l:SC-auto} ensures that $(v,S[V])=(\bar S(h),S[V])\in \SC_\Pi$ and $\bar S(v)=S[V]_\parallel=\PSC_\Pi(v)=h$. 

\begin{claim} The automorphism $S$ maps each point $p\in\Pi$ with coordinates $(x,y)$ in the base $uow$ onto the point $S(p)$ with coordinates $(y,x)$. In particular, $w'=S(u)=w$ and hence $\Aline uw\in\bar\delta$.
\end{claim}

\begin{proof} If $p\in\Delta$, then $x=p=y$ and $S(p)=p$ has coordinates $(y,x)=(x,y)=(p,p)$. So, assume that $p\notin \Delta$. The definition of the coordinates of $p$ in the base $uow$ ensures that $x,y\in \Delta$ are two points such that $\Aline px\in v$ and $\Aline py\in h$. Then for the image $q\defeq S(p)$, we have $\Aline xq=S[\Aline xp]\in S[\Aline xp]_\parallel=\bar S(\Aline xp_\parallel)=\bar S(h)=v$ and $\Aline yq=S[\Aline yp]\in S[\Aline yp]_\parallel=\bar S(\Aline yp_\parallel)=\bar S(v)=h$, which implies that $(y,x)$ are coordinates of the point $q=S(p)$ in the base $uow$.

Since the point $u$ has coordinates $(e,o)$, its image $w'=S(w)$ has coordinates $(o,e)$ and hence coincides with the point $w$. Then $\Aline uw=\Aline u{w'}\in\bar\delta$.
\end{proof}

Since the ternar $\Delta$ of the based affine plane $(\Pi,uow)$ is associative-dot, its dot loop $\Delta^*\defeq\Delta\setminus\{o\}$ is a group. So, for every $a\in \Delta^*$ there exists a unique element $a^{-1}\in\Delta^*$ such that $a\cdot a^{-1}=e=a^{-1}\cdot a$.

For every $a\in\Delta$, let $L_{a,o}\in\mathcal L_o$ be the unique line that contains the point with coordinates $(e,a)$. The line $L_{a,o}$ is defined by the equation $y=x\cdot a=x_\times a_+o$ in the coordinate plane of the based affine plane $(\Pi,uow)$. 

\begin{claim} $S[L_{a,o}]=L_{a^{-1},o}$ for every $a\in\Delta^*$.
\end{claim}

\begin{proof} Since the line $L_{a,o}$ contains the point $p$ with coordnates $(e,a)$, the line $S[L_{a,o}]$ contains the point $q=S(p)$ with coordinates $(a,e)$. Since $a\cdot a^{-1}=a_\times a^{-1}_+o=e$, the point $q$ belongs to the line $L_{a^{-1},o}$ and hence $S[L_{a,o}]=L_{a^{-1},o}$.
\end{proof} 

Let $\bar e\in\Delta$ be the unique point such that $L_{\bar e,o}\in\bar\delta$. Since the automorphism $S$ preserves the direction $\bar\delta$, $L_{\bar e^{-1},o}=S[L_{\bar e,o}]=L_{\bar e,o}$ and hence $\bar e^{-1}=\bar e$ and $\bar e{\cdot}\bar e=e$.

\begin{claim}\label{cl:CS-unique2} Any point $a\in\Delta$ with $a{\cdot}a=e$ is equal to $e$ or $\bar e$.
\end{claim}

\begin{proof} Take any point $a\in\Delta$ with $a{\cdot}a=e$. Then $a^{-1}=a$ and hence $L_{a,o}=L_{a^{-1},o}=S[L_{a,o}]$. By Proposition~\ref{p:Fix(A)=L+p}, $\Fix(\bar S)=\Delta\cup\{\bar\delta\}$ and hence the direction $\lambda$ of the line $L_{a,o}=S[L_{a,o}]$ is equal to $\delta=\Delta_\parallel$ or $\bar \delta$. Then $L_{a,o}\in\{L_{e,o},L_{\bar e,o}\}$ and hence $a\in\{e,\bar e\}$.
\end{proof}

\begin{claim} $x\cdot\bar e=\bar e\cdot x$ for all $x\in \Delta$.
\end{claim}

\begin{proof} If $x=o$, then $x\cdot \bar e=o=\bar e\cdot x$. If $x\ne o$, then the element $x^{-1}{\cdot}\bar e{\cdot}x\ne e$ of the group $(\Delta^*,\cdot)$ has order $2$ and hence is equal to $\bar e$, by Claim~\ref{cl:CS-unique2}. The equality $x^{-1}{\cdot}\bar e{\cdot}x=\bar  e$ is equivalent to $\bar e{\cdot}x=x{\cdot}\bar e$.
\end{proof} 

\begin{claim}\label{cl:aa=bare} For any point $a\in \Delta\setminus\{o,e,\bar e\}$ we have $a{\cdot}a=\bar e$ and hence $\bar e{\cdot}a=a^{-1}$.
\end{claim}

\begin{proof} Fix any point $a\in\Delta\setminus\{o,e,\bar e\}$. Since the based affine plane $(\Delta,uow)$ is vertical-scale, we can apply  Proposition~\ref{p:vscale} and conclude that the function $A:\Pi\to \Pi$ mapping a point with coordinates $(x,y)$ to the point with coordinates $(x,y{\cdot}\bar e)$, is an automorphism of the plane $\Pi$ such that $H\subseteq\Fix(A)$ and $A[V]=V$. Since $(\Delta_\parallel,A[\Delta])=(\delta,L_{\bar e,o})\in \SC_\Pi$, we can apply Lemma~\ref{l:PSC-all} and conclude that $((L_{a,o})_\parallel,A[L_{a,o}])\in\SC_\Pi$. On the other hand $(H_\parallel,S[H])=(h,V)\in \SC_\Pi$ implies $((L_{a,o})_\parallel,S[L_{a,0}])\in \SC_\Pi$, by the same Lemma~\ref{l:PSC-all}. Since $\PSC_\Pi$ is a function, $A[L_{a,o}]\parallel S[L_{a,o}]=L_{a^{-1},o}$ and hence $L_{a{\cdot}\bar e,o}=A[L_{a,o}]=L_{a^{-1},o}$. Then $a{\cdot}\bar e=a^{-1}$ and hence $a{\cdot}a=\bar e$.
\end{proof}

\begin{claim}\label{cl:yx=exy} For any elements $x,y\in \Delta^*$ with $\{x,y,x{\cdot}y\}\cap\{e,\bar e\}=\varnothing$, we have $y{\cdot}x=\bar e{\cdot} x{\cdot}y$.
\end{claim}

\begin{proof} Since $x,y,x{\cdot}y\notin \{e,\bar e\}$, Claim~\ref{cl:aa=bare} ensures that $$x{\cdot}(y{\cdot}x){\cdot}y=(x{\cdot}y){\cdot}(x{\cdot}y)=\bar e=\bar e{\cdot}\bar e{\cdot}\bar e=\bar e{\cdot}(x{\cdot}x){\cdot}(y{\cdot}y)=x{\cdot}(\bar e{\cdot}x{\cdot}y){\cdot}y$$ and hence $y{\cdot}x=\bar e{\cdot}x{\cdot}y$, by the cancellativity of the group $\Delta^*$.
\end{proof}

Consider the two-element central (and hence normal) subgroup $E\defeq\{e,\bar e\}$ in the group $(\Delta^*,\cdot)$. Since $x^2\in E$ for all $x\in \Delta^*$, the quotient group $G=\Delta^*/E$ is Boolean. Let $q:\Delta^*\to G$ be the quotient homomorphism.

Choose any element $a\in\Delta^*\setminus\{e,\bar e\}$. Since $|\Delta^*|=|\Pi|_2-1\ge 6$, there exists an element $b\in\Delta^*\setminus\{e,\bar e,a^{-1},a^{-1}\bar e\}$. Such a choice of $b$ ensures that $a{\cdot}b\notin \{e,\bar e\}$ and hence $b{\cdot}a=\bar e{\cdot} b{\cdot}a\ne b{\cdot}a$, by Claim~\ref{cl:yx=exy}. 
Since the elements $a,b$ do not commute, their images in the quotient group $\Delta^*/E$ are distinct and hence generate a subgroup of order $4$ in the Boolean group $\Delta^*/E$. Then $a,b$ generate a subgroup $G$ of order $8$ in the group $\Delta^*$. Assuming that $|\Delta^*|>8$, we can find an element $c\in \Delta^*\setminus G$. For every $x\in G\setminus E$ we have $c{\cdot}x{\cdot}c^{-1}=\bar e{\cdot}x=x^{-1}$, by Claims~\ref{cl:yx=exy} and \ref{cl:aa=bare}.
In particular, $$(a{\cdot}b)^{-1}=c{\cdot}(a{\cdot}b){\cdot}c^{-1}=(c{\cdot}a{\cdot}c^{-1}){\cdot}(c{\cdot}b{\cdot}c^{-1})=a^{-1}{\cdot}b^{-1}=(b{\cdot}a)^{-1}$$ and hence $a{\cdot}b=b{\cdot}a$, which contradicts Claim~\ref{cl:yx=exy}. This contradiction shows that $\Delta^*=G$ and hence $|\Pi|_2=|\Delta|=|\Delta^*|+1=|G|+1=8+1=9$ and $|\PSC_\Pi|=|\partial\Pi|=10$.
\end{proof}

Now we are able to present {\em a proof of Theorem~\ref{t:scale-figure}}. Let $\Pi$ is an affine plane.
\smallskip

If $\Pi$ is Desarguesian, then Lemma~\ref{l:Desarg=>SC=} ensures that $\PSC_\Pi=\{(\delta,L)\in\partial\Pi\times\mathcal L_\Pi:L\notin \delta\}$ and $\|\PSC_\Pi\|=22$, which is the case $(\Dscale)$. So, assume that the plane $\Pi$ is not Desarguesian.
\smallskip

If $|\PSC_\Pi|\le 1$, then one of the cases $(\circ)$ or $(\bullet)$ holds.
\smallskip

If $|\PSC_\Pi|=2$, then we have the case $(\didots)$, by Lemma~\ref{l:PSC=2}.

So, assume that $|\PSC_\Pi|\ge 3$. If $\PSC_\Pi$ is not a function, then the case $(\,\vert\,)$ holds by Lemma~\ref{l:PSC-not-function2}. If $\PSC_\Pi^{-1}$ is not a function, then the case $({-})$ holds by Lemma~\ref{l:PCS-not-function2}. So, assume that $\PSC_\Pi$ and $\PSC_\Pi^{-1}$ both are functions. Applying Lemma~\ref{l:SC-largedom}, we conclude that $\dom[\PSC_\Pi]=\partial\Pi=\rng[\PSC_\Pi]$. By Lemma~\ref{l:PSC-involution}, $\PSC_\Pi$ is an involution and by Lemma~\ref{l:PSC-involution=>9}, $|\PSC_\Pi|=|\partial\Pi|=10$. Therefore, the case $(\invol)$ holds.
\end{proof}

Now we are able to classify hyperscale traces of affine planes into 10 types (including two exceptional).

\begin{theorem}\label{t:scale-trace} For every affine plane $\Pi$, exactly one of the following $10$ cases holds:
\begin{enumerate}
\item[{$(\,\circ\,)$}]\index[note]{$(\,\circ\,)$}
$\|\SC_\Pi\|=00$ and $\SC_\Pi=\varnothing$;
\item[{$(\,\bullet\,)$}]\index[note]{$(\,\bullet\,)$} $\|\SC_\Pi\|=11$ and $\SC_\Pi=\{(\delta,L)\}$ for some direction $\delta\in\partial\Pi$ and line $L\in\mathcal L_\Pi\setminus\delta$;
\item[{$(\hskip2.5pt\vert\hskip2.5pt)$}]\index[note]{$(\hskip2.5pt\vert\hskip2.5pt)$}   $\|\SC_\Pi\|=12$ and $\SC_\Pi=\{\delta\}\times (\mathcal L_p\setminus\delta)$ for some direction $\delta\in\partial\Pi$ and point $p\in\Pi$;
\item[{$(-)$}]\index[note]{$(-)$} $\|\SC_\Pi\|=21$ and  $\SC_\Pi=(\partial \Pi\setminus\{\delta\})\times\{L\}$ for some direction $\delta\in\partial\Pi$ and line $L\in\delta$;
\item[$(\discale)$]\index[note]{$(\discale)$} $\|\SC_\Pi\|=22$ and $\SC_\Pi$ is an injective function with $\dom[\SC_\Pi]=\partial \Pi\setminus\{\delta\}$ and $\rng[\SC_\Pi]=\delta$;
\item[{$(\didots)$}]
\index[note]{$(\didots)$} 
$\|\SC_\Pi\|=22$ and $\SC_\Pi=\{(L_\parallel,\Lambda),(\Lambda_\parallel,L)\}$ for two concurent lines $L,\Lambda\in\partial\Pi$;
\item[{$(\reshitka)$}]
\index[note]{$(\reshitka)$}
 $\|\SC_\Pi\|=23$ and $\SC_\Pi=(\{h\}\times v)\cup(\{v\}\times h)$ for distinct directions $v,h\in\partial\Pi$;
\item[{$(\invol)$}]\index[note]{$(\invol)$} $\|\SC_\Pi\|=22$, $\PSC_\Pi$ is an involution with $|\PSC_\Pi|=|\partial\Pi|=10$ and\\ \phantom{m}\hskip62pt $\SC_\Pi=\{(\delta,L)\in\partial\Pi\times\mathcal L_p:(\delta,L_\parallel)\in\PSC_\Pi\}$ for some point $p\in\Pi$;
\item[$(\Invol)$]\index[note]{$(\Invol)$} $\|\SC_\Pi\|=23$, $\PSC_\Pi$ is an involution with $|\PSC_\Pi|=|\partial\Pi|=10$ and\\
\phantom{m}\hskip62pt $\SC_\Pi=\{(\delta,L)\in\partial\Pi\times\mathcal L_\Pi:(\delta,L_\parallel)\in\PSC_\Pi\}$;
\item[{$(\Dscale)$}]\index[note]{$(\Dscale)$} $\|\SC_\Pi\|=23$ and $\SC_\Pi=\{(\delta,L)\in\partial\Pi\times\mathcal L_\Pi:L\notin \delta\}$.
\end{enumerate}
\end{theorem}

\begin{proof} Let $\Pi$ be an affine plane. If  $\Pi$ is Desarguesian, then $\SC_\Pi=\{(\delta,L)\in\partial\Pi\times\mathcal L_\Pi:L\notin\delta\}$, by Lemma~\ref{l:Desarg=>SC=}, witnessing that the case $(\Dscale)$ holds. So, assume that the affine plane $\Pi$ is not Desarguesian.
\smallskip

If $|\SC_\Pi|=0$, then $\SC_\Pi=\varnothing$ and we have the case $(\circ)$.

If $|\SC_\Pi|=1$, then $\SC_\Pi=\{(\delta,L)\}$ for some direction $\delta\in\partial\Pi$ and line $L\in\mathcal L_\Pi\setminus\delta$, so the case $(\bullet)$ holds.

If $|\dom[\SC_\Pi]|=1<|\SC_\Pi|$, then there exist a direction $\delta\in\partial\Pi$ and two distinct lines $L,\Lambda\in\mathcal L_\Pi\setminus\delta$ such that $\{(\delta,L),(\delta,\Lambda)\}\subseteq\SC_\Pi$. If the lines $L,\Lambda$ are disjoint, then Lemmas~\ref{l:h-rozmaznia} and \ref{l:trans-swap} ensure that $(\{\delta\}\times L_\parallel)\cup(\{L_\parallel\}\times \delta)\subseteq \SC_\Pi$ and hence $\{\delta,L_\parallel\}\subseteq\dom[\SC_\Pi]$, which contradicts $|\dom[\SC_\Pi]|=1$. This contradiction shows that the lines $L,\Lambda$ are concurrent and hence $L_\parallel\ne\Lambda_\parallel$. Since $\{(\delta,L_\parallel),(\delta,\Lambda_\parallel)\}\subseteq\PSC_\Pi$, the relation $\PSC_\Pi$ is not a function. By Lemmas~\ref{l:PSC-not-function} and \ref{l:PSC-not-function2}, $\{\delta\}\times(\mathcal L_p\setminus\delta)\subseteq\SC_\Pi$ and $\{\delta\}\times(\partial\Pi\setminus\{\delta\})=\PSC_\Pi$. Assuming that $\{\delta\}\times(\mathcal L_p\setminus\delta)\ne\SC_\Pi$, we can find a pair $(\delta,A)\in \SC_\Pi$ such that $A\notin\mathcal L_p$. There exists a line $B\in\mathcal L_p$ such that $A\cap B=\varnothing$. Since $\{(\delta,A),(\delta,B)\}\subseteq\SC_\Pi$, we can apply Lemma~\ref{l:h-rozmaznia} and \ref{l:trans-swap} conclude that $\{A_\parallel\}\times\delta\subseteq\SC_\Pi$ and hence $\{\delta,A_\parallel\}\subseteq\dom[\SC_\Pi]$, which contradicts $|\dom[\SC_\Pi]|=1$. This contradiction shows that $\{\delta\}\times(\mathcal L_p\setminus\delta)=\SC_\Pi$ and hence $\|\SC_\Pi\|=12$ and the case $(\,\vert\,)$ holds.
\smallskip

If $|\dom[\SC_\Pi]|>1$ and $\rng[\PSC_\Pi]=1$, then by Theorem~\ref{t:scale-figure}, there exists a direction  $\delta\in\partial\Pi$ such that $\PSC_\Pi=(\partial\Pi\setminus\{\delta\})\times\{\delta\}$ and hence $\dom[\SC_\Pi]=\dom[\PSC_\Pi]=\partial\Pi\setminus\{\delta\}$.
We claim that $\SC_\Pi$ is a function. In the opposite case, there exist 
a direction $d\in\partial\Pi\setminus\{\delta\}$ and two distinct lines $L,\Lambda\in \delta$ such that $\{(d,L),(d,\Lambda)\}\subseteq\SC_\Pi$.  Applying Lemma~\ref{l:h-rozmaznia} and \ref{l:trans-swap}, we conclude that $d\in \rng[\PSC_\Pi]$, which contradicts $|\PSC_\Pi|=1$. This contradiction shows that $\SC_\Pi$ is a function. If the function $\SC_\Pi$ is constant, then $\SC_\Pi=(\partial\Pi\setminus\{\delta\}) \times\{L\}$ for some line $L$. Then $\|\SC_\Pi\|= 21$ and the case $({-})$ holds. 

So, assume that the function $\SC_\Pi$ is not constant. We claim that $\SC_\Pi$ is an injective
function with $\dom[\SC_\Pi] = \partial\Pi\setminus\{\delta\}$ and $\rng[\SC_\Pi] = \delta$.
Assuming that the function $\SC_\Pi$ is not injective, we can find two distinct directions $v, h\in\partial\Pi\setminus\{\delta\}$ and a line $L\in\delta$ such that $\{(v, L),(h, L)\}\subseteq \SC_\Pi$. We claim that $(\partial\Pi\setminus\{\delta\})\times\{L\}\subseteq\SC_\Pi$. 
Indeed, given any direction $w\in \partial \Pi\setminus\{\delta,v,h\}$, fix any point $o\in L$ and choose any point
$x\in\Pi\setminus L$ such that $\Aline ox\in h$. Next, find a unique point $y\in \Pi\setminus\{x\}$ such that $\Aline oy\in w$ and
$\Aline xy \in v$. Since $L\in \delta\ne w$, the point $y$ does not belong to the line $L$. Since $(v, L)\in\SC_\Pi$, there
exists an automorphism $\Phi$ of the plane $\Pi$ such that $\Phi(x) = y$ and $L\subseteq\Fix(\Phi)$. Let $\bar\Phi:\overline\Pi\to\overline\Pi$ be the spread extension of the automorphism $\Phi$. It follows from $\Aline o x\in h$ and $\Phi[\Aline o x] = \Aline o y\in w$ that $\bar\Phi (h) = \Aline o y_\parallel = w$. By Lemma~\ref{l:SC-auto}, $(w, L) = (\bar\Phi(h),\Phi[L]) \in\SC_\Pi$, witnessing that $(\partial\Pi\setminus\{\delta\})\times\{L\}\subseteq\SC_\Pi$. Taking into account that $\SC_\Pi$ is a function with $\dom[\SC_\Pi] = \dom[\PSC_\Pi] = \partial\Pi\setminus\{\delta\}$, we conclude that $\SC_\Pi=(\partial\Pi\times\{\delta\})\times\{L\}$ is a
constant function, which contradicts our assumption. This constradiction shows that the function $\SC_\Pi$ is injective.

It remains to show that $\rng[\SC_\Pi] = \delta$. Since $\SC_\Pi$ is an injective function defined on the set $\partial \Pi\setminus\{\delta\}$ of cardinality $\ge 2$, there exist two distinct pairs $(v, L),(v',L')\in\SC_\Pi$. Then $\{L, L'\}\subseteq \rng[\SC_\Pi]$. It remains to show that every line $\Lambda\in\delta\setminus\{L,L'\}$ belongs to the set $\rng[\SC_\Pi]$. Choose
any point $x \in L'$. Since $\Lambda\in\delta\ne v$, there exists a unique point $y\in\Lambda$ such that $\Aline x y\in v$. Since
$(v, L)\in\SC_\Pi$ and $\{x, y\} \subseteq L'\cup\Lambda \subseteq \Pi\setminus\{L\}$, there exists an automorphism $\Psi$ of the plane $\Pi$  such that $\Psi(x) = y$ and $L \subseteq\Fix(\psi)$. Then $y = \Psi(x)\in \Lambda\cap\Psi[L']$ and $\Psi[L']\parallel\Psi[L]= L\parallel\Lambda$ imply $\Psi[L'] = \Lambda$. Let $\bar\Psi$ be the spread completion of the automorphism $\Psi$. Applying Lemma~\ref{l:SC-auto}, we conclude that $(\bar\Psi(v),\Lambda) = (\bar\Psi(v),\Psi[L])\in\SC_\Pi$ and hence $\Lambda\in\rng[\SC_\Pi]$, witnessing that $\rng[\SC_\Pi] = \delta$. Then $\|\SC_\Pi\|= 22$ and the case $(\discale)$ holds. 
\smallskip

Next, assume that $|\dom[\SC_\Pi]|>1$ and $|\rng[\PSC_\Pi]|>1$. If $|\PSC_\Pi|=2$, then by Theorem~\ref{t:scale-figure}, $\PSC_\Pi=\{(v,h),(h,v)\}$ for some distinct directions $v,h$. Then there exist lines $H\in h$ and $V\in V$ such that $\{(v,H),(h,V)\}\subseteq \SC_\Pi$. If $\{(v,H),(h,V)\}=\SC_\Pi$, then $\|\SC_\Pi\|=22$ and the case $({}^\centerdot{\;}_\centerdot)$ holds. 
\vskip3pt

So, assume that $\{(v,H),(h,V)\}\ne\SC_\Pi$. The equality $\PSC_\Pi=\{(v,h),(h,v)\}$ implies $\SC_\Pi\subseteq(\{v\}\times h)\cup(\{h\}\times v)$. Since $\{(v,H),(h,V)\}\subseteq\SC_\Pi\subseteq(\{v\}\times h)\cup(\{h\}\times v)$, there exists a line $L\in v\cup h$ such that $(v,L)$ or $(h,L)$ belongs to the set $\SC_\Pi\setminus\{(v,H),(h,V)\}$. We lose no generality assuming that $(v,L)\in\SC_\Pi\setminus\{(v,H)\}$. Then $L,H\in h$ are two disjoint lines with $\{(v,L),(v,H)\}\subseteq\SC_\Pi$. Applying Lemma~\ref{l:SC=>translation}, we conclude that $\Pi$ is a translation plane. Then for every line $H'\in h$ we can find a translation $T$ of the plane $\Pi$ such that $T[H]=H'$. Since translations do no change directions, $(v,H)\in\SC_\Pi$ implies $(v,H')=(\bar T(h),T[H])\in \SC_\Pi$, witnessing that $\{v\}\times h\subseteq\SC_\Pi$. By analogy, we can prove that $\{h\}\times v\subseteq\SC_\Pi$ and hence $\SC_\Pi=(\{v\}\times h)\cup(\{h\}\times v)$. Then $\|\SC_\Pi\|=23$ and the case $(\reshitka)$ holds.
\smallskip

If $|\PSC_\Pi|>2$, then Theorem~\ref{t:scale-figure} ensures that $\PSC_\Pi$ is an involution with $|\PSC_\Pi|=|\partial\Pi|=10$. 
Choose any pair $(v,h)\in\PSC_\Pi$. Since $\PSC_\Pi$ is an involution, $(h,v)\in \PSC_\Pi$. Then there exist lines $H\in h$ and $V\in v$ such that $\{(h,V),(v,H)\}\subseteq\SC_\Pi$. Since the directions $h,v$ are distinct, the lines $H,V$ are concurrent and hence have a common point $p\in H\cap V$. Claim~\ref{cl:SC-Delta} ensures that $\mathcal L_p\subseteq\rng[\SC_\Pi]$. Observe that $\SC_\Pi^{-1}$ is a function assigning to every line $L\in\rng[\SC_\Pi]$ the direction $\PSC^{-1}_\Pi(L_\parallel)$. Since the map $\mathcal L_p\to\partial \Pi$, $L\mapsto L_\parallel$, is bijective, the restriction $\SC_\Pi^{-1}{\restriction}_{\mathcal L_p}$ is a bijective function between the sets $\mathcal L_p$ and $\partial\Pi$. If $\rng[\SC_\Pi]=\mathcal L_p$, then $\SC_\Pi=\{(\delta,L)\in\partial\Pi\times\mathcal L_p:(\delta,L_\parallel)\in \PSC_\Pi\}$ is a bijective function from $\partial\Pi$ into $\mathcal L_p$, $\|\SC_\Pi\|=22$ and the case $(\invol)$ holds.
\smallskip

Finally, assume that $\rng[\SC_\Pi]\ne\mathcal L_p$. Then there exists a line $L\in\rng[\SC_\Pi]\setminus\mathcal L_p$. Find a unique line $L'\in\mathcal L_p$ with $L'\parallel L$ and observe that $\delta\defeq\SC_\Pi^{-1}(L)=\SC_\Pi^{-1}(L')$. Then $\{(\delta,L),(\delta,L')\}\subseteq\SC_\Pi$ and we can apply Lemma~\ref{l:h-rozmaznia} and \ref{l:trans-swap} to conclude that the plane $\Pi$ is translation. In this case, for every pair $(v,H)\in\partial \Pi\times\mathcal L_\Pi$ with $(\delta,L_\parallel)\in\PSC_\Pi$, we can find a line $H'\in\mathcal L_p\cap L_\parallel$ and conclude that $(v,H')\in \SC_\Pi$. Choose any translation $T$ of $\Pi$ with $T[H']=H$ and conclude that $(v,H)=(\bar T(v),T[H'])\in \SC_\Pi$, witnessing that $\{(\delta,L)\in\partial\Pi\times\mathcal L_\Pi:(\delta,L_\parallel)\in \PSC_\Pi\}=\SC_\Pi$. Then $\|\SC_\Pi\|=23$ and the final case $(\Invol)$ holds.
\end{proof}


\begin{remark} Among the 10 classes in the classification Theorem~\ref{t:scale-trace}, the classes $(\invol)$, $(\Invol)$ and $(\discale)$ are somewhat special: up to an isomorphism there is just one plane in each of the classes $(\invol)$ and $(\Invol)$. There are no known examples of affine planes of scale type $(\discale)$.
\end{remark}


\begin{definition}\label{d:scale-trace} An affine plane $\Pi$ has scale type {\em at least}
\begin{enumerate}
\item[{$(\,\circ\,)$}]\index[note]{$(\,\circ\,)$} if $\varnothing\subseteq \SC_\Pi$;
\item[{$(\,\bullet\,)$}]\index[note]{$(\,\bullet\,)$} if $\{(\delta,L)\}\subseteq \SC_\Pi$ for some direction $\delta\in\partial \Pi$ and line $L\in\mathcal L_\Pi$;
\item[{$(\hskip2.5pt\vert\hskip2.5pt)$}]\index[note]{$(\hskip2.5pt\vert\hskip2.5pt)$} if $\{\delta\}\times(\mathcal L_p\setminus\delta)\subseteq\SC_\Pi$ for some direction $\delta\in\partial\Pi$ and point $p\in\Pi$;
\item[{$(-)$}]\index[note]{$(-)$} if $(\partial\Pi\setminus\{L_\parallel\})\times \{L\}\subseteq \SC_\Pi$ for some line $L\in\mathcal L_\Pi$;
\item[$(\discale)$]\index[note]{$(\discale)$} if $F\subseteq\SC_\Pi$ for some direction $\delta\in\partial\Pi$ and some bijective function $F:\partial\Pi\setminus\{\delta\}\to\delta$;
\item[{$(\didots)$}]
\index[note]{$(\didots)$}  $\{(L_\parallel,\Lambda),(\Lambda_\parallel,L)\}\subseteq\SC_\Pi$ for two concurent lines $L,\Lambda\in\partial\Pi$;
\item[{$(\reshitka)$}]
 if $(\{h\}\times v)\cup(\{v\}\times h)\subseteq \SC_\Pi$ for distinct directions $v,h\in\partial\Pi$;
 \item[{$(\invol)$}] if  $\{(\delta,L)\in\partial\Pi\times\mathcal L_p:(\delta,L_\parallel)\in F\}\subseteq\SC_\Pi$ for some involution $F:\partial \Pi\to \partial\Pi$ and some point $p\in\Pi$;
\item[$(\Invol)$] $\{(\delta,L)\in\partial\Pi\times\mathcal L_\Pi:(\delta,L_\parallel)\in F\}\subseteq\SC_\Pi$ for some involution $F:\partial\Pi\to\partial\Pi$;
\item[{$(\Dscale)$}] if $\{(\delta,L)\in\partial\Pi\times\mathcal L_\Pi:L\notin \delta\}\subseteq\SC_\Pi$.
\end{enumerate}
\end{definition}

Definition~\ref{d:scale-trace} actually involves a natural partial order  
on the $10$-element set of scale types $$\{(\circ), (\bullet), (-), (\,|\,), (\discale), (\didots), (\reshitka), (\invol), (\Invol), (\Dscale)\}.$$
The Hasse diagram of this partial order looks as follows (the class $(\discale)$ is colored in gray because it can disappear if someone will succeed to prove that it is empty):
$$
\xymatrix@C=0pt@R=12pt{
&&&(\Dscale)&&&&&&&\\
&&&(\Invol)\ar@{-}[u]&&&&&\\
(-)\ar@{-}[uurrr]&&(\invol)\ar@{-}[ur]&{\color{gray}(\discale)}\ar@/^12pt/@{..}[uu]&(\reshitka)\ar@{-}[ul]&&(\,\vert\,)\ar@{-}[uulll]\\
&&&(\didots)\ar@{-}[ur]\ar@{-}[ul]&&&&\\
&&&(\bullet)\ar@{-}[u]\ar@{-}[uulll]\ar@/_12pt/@{..}[uu]\ar@{-}[uurrr]\\
&&&(\circ)\ar@{-}[u]
}
$$

\begin{theorem}\label{p:at-least-scale} An affine plane $\Pi$ has scale type {\em at least}
\begin{enumerate}
\item[{$(\,\bullet\,)$}]\index[note]{$(\,\bullet\,)$} iff some ternar of $\Pi$ is linear and associative-dot;
\item[{$(-)$}]\index[note]{$(-)$} iff some ternar of $\Pi$ is linear, left-distributive, and associative;
\item[$(\discale)$]\index[note]{$(\discale)$} iff some ternar of $\Pi$ is linear, associative-dot and diagonal-scale;
\item[{$(\didots)$}]\index[note]{$(\didots)$} iff some ternar of $\Pi$ is linear right-distributive and associative-dot;
\item[{$(\reshitka)$}]\index[note]{$(\reshitka)$}
 iff some ternar of $\Pi$ is linear, right-distributive and associative;
\item[{$(\Dscale)$}]\index[note]{$(\Dscale)$} iff some ternar of $\Pi$ is linear, distributive and associative.
\end{enumerate}
\end{theorem}

\begin{proof} $(\bullet)$ If the affine plane $\Pi$ has the scale type at least $(\bullet)$, then $(h,V)\in\SC_\Pi$ for some direction $h\in\partial\Pi$ and line $V\in\mathcal L_\Pi\setminus h$. Choose an affine base $ouw$ in $\Pi$ such that $\Aline ow=V$ and $\Aline ou\in h$. Then the based affine plane $(\Pi,uow)$ is horizontal-scale and by Theorem~\ref{t:horizontal-scale}, its ternar is linear and associative-dot.

Now assume that some ternar of the affine plane $\Pi$ is linear and associative-dot. Then $\Pi$ has an affine base $uow$ whose ternar is linear and associative-dot. By Theorem~\ref{t:horizontal-scale}, the based affine plane $(\Pi,uow)$ is horizontal-scale and hence $((\Aline ow)_\parallel,\Aline ou)\in\SC_\Pi$, witnessing that $\Pi$ is of scale type at least $(\bullet)$.
\smallskip

$(-)$ If the affine plane $\Pi$ has  the scale type at least $(-)$, then $(\partial\Pi\setminus\{V_\parallel\})\times\{V\}$ for some line $V$ in $\Pi$.
Choose an affine base $uow$ in $\Pi$ such that $V=\Aline ow$. 
Then the based affine plane $(\Pi,uow)$ is horizonthal-scale and vertical-shear, by Lemma~\ref{l:scale=>shear}. By Corollary~\ref{c:vertical-shear+horizontal-scale}, its ternar is linear, associative-dot and left-distributive.

Now assume that some ternar of the affine plane $\Pi$ is linear, left-distributive and associative-dot. Then $\Pi$ has an affine base $uow$ whose ternar is linear, left-distributive and associative-dot. By Corollary~\ref{c:vertical-shear+horizontal-scale}, the based affine plane $(\Pi,uow)$ is horizontal-scale and vertical-shear. Then for the horizontal direction $h\defeq(\Aline ou)_\parallel$ and the vertical line $V\defeq\Aline ow$, the pair $(h,V)$ belongs to the scale trace $\SC_\Pi$. Since the line $V$ is shear, for every line $L\in\mathcal L_o\setminus\{V\}$ there exists an automorphism $\Phi$ of the plane $\Pi$ such that $\Phi[\Aline ou]=L$ and $V\subseteq\Fix(\Phi)$. Then $(L_\parallel,V)=(\bar\Phi(h),V)\in \SC_\Pi$, by Lemma~\ref{l:SC-auto} and hence $\{(\delta,V):\delta\in\partial\Pi\setminus\{V_\parallel\}\}\subseteq\SC_\Pi$, witnessing that the plane $\Pi$ has scale type at least $(-)$.
\smallskip

$(\discale)$ If the affine plane $\Pi$ has the scale type at least $(\discale)$, then $F\subseteq\SC_\Pi$ for some direction $v\in\partial\Pi$ and some bijective function $F:\partial\Pi\setminus\{v\}\to V$. Choose any distinct pairs $(h,V),(d,E)$ in $\SC_\Pi$ and observe that the distinct lines $V,E\in v$ are parallel and hence disjoint. Choose any points $o\in V$ and $e\in E$ such that $\Aline oe\in d$ and then find unique points $u\in E$ and $w\in V$ such that $\Aline ou\in h$ and $\Aline we\in h$. It follows from $\{((\Aline ou)_\parallel,\Aline ow),((\Aline oe)_\parallel,\Aline ue)\}=\{(h,V),(d,E)\}\subseteq\SC_\Pi$ that the based affine pane $(\Pi,uow)$ is horizontal-scale and diagonal-scale. By Theorem~\ref{t:diagonal-scale<=>}, its ternar is linear, associative-dot and diagonal-scale.  

Now assume that some ternar of the affine plane $\Pi$ is linear, associative-dot, and  diagonal-scale. Then $\Pi$ has an affine base $uow$ whose ternar $\Delta$ is linear, associative-dot, and diagonal-scale.  By Theorem~\ref{t:diagonal-scale<=>}, the based affine plane $(\Pi,uow)$ is horizontal-scale and diagonal-scale. Consider the horizontal and vertical directions $h\defeq(\Aline ou)_\parallel$ and $v\defeq(\Aline ow)_\parallel$. For every $x\in\Delta$, consider the horizontal and vertical lines $H_x\in  h\cap \mathcal L_x$ and $V_x\in v\cap\mathcal L_x$. Let $x^{-1}$ be the unique element of the ternar $\Delta$ such that $x{\cdot}x^{-1}{\cdot}x=x^{-1}$ and $x^{-1}\cdot x\cdot x^{-1}$. So, $x^{-1}$ is the inverse element of $x$ in the group $(\Delta\setminus\{o\},\cdot)$ if $x\ne o$ and $x^{-1}=o^{-1}=o$ if $x=o$. Since the based affine plane $(\Pi,uow)$ is horizontal-scale and diagonal-scale, $\{(x^{-1}_\parallel,V_x):x\in\{0,1\}\}\subseteq\SC_\Pi$.  Since $\Pi$ is horizontal-scale, for every $c\in \Delta\setminus\{o\}$, there exists an automorphism $\Phi_c:\Pi\to \Pi$ such that $\Phi_c[V_e]=V_c$, $V_o\subseteq\Fix(\Phi_c)$, and $\bar\Phi_c(h)=h$, where $\bar\Phi:\overline\Pi\to\overline\Pi$ is the spread extension of the automorphism $\Phi$. Then $\Phi_c(e)\in \Phi[H_e\cap V_e]=H_e\cap V_c$ and $\Phi[\Delta]=L_{c^{-1},o}$ is the line with slope $c^{-1}$.  Lemma~\ref{l:SC-auto} ensures that $(c^{-1}_\parallel,V_c)=(\bar\Phi_c(e_\parallel),\Phi_c[V_e])\in\SC_{\Pi}$ and hence $F\defeq \{(c^{-1}_\parallel,V_c):c\in R\}\subseteq\SC_{\Pi}$. Taking into account that $F$ is an injective function with domain $\partial\Pi\setminus\{v\}$ and range $\rng[F]=v$, we conclude that the plane $\Pi$ is of scale type at least $(\discale)$.
\smallskip

$(\didots)$ If the ffine plane $\Pi$ is of scale type at least $(\didots)$, then $\{(V_\parallel,H),(H_\parallel,V)\}\subseteq\SC_\Pi$ for two concurrent lines $V,H$ in $\Pi$. Choose an affine base $uow$ in the plane $\Pi$ such that $\Aline ou=H$ and $\Aline ow=V$. Then the based affine plane $(\Pi,ouw)$ is horizontal-scale and vertical-scale and its ternar is linear,  right-distributive, and associative-dot, by Corollary~\ref{c:hv-scale<=>}.

Now assume that some ternar of the affine plane $\Pi$ is linear, right-distributive and associative-dot. Then $\Pi$ has an affine base $uow$ whose ternar is linear, right-distributuve and associative-dot. By Corollary~\ref{c:hv-scale<=>}, the based affine plane $(\Pi,uow)$ is horizontal-scale and vertical-scale and hence $\{(\Aline ou_\parallel,\Aline ow),(\Aline ow_\parallel,\Aline ou)\}\subseteq\SC_\Pi$, witnessing that the affine plane $\Pi$ is of scale tyle at least $(\didots)$.
\smallskip

$(\reshitka)$ If the ffine plane $\Pi$ is of scale type at least $(\reshitka)$, then $(\{v\}\times h)\cup(\{h\}\times v)\subseteq\SC_\Pi$ for distinct directions $h,v\in\partial\Pi$. Choose an affine base $uow$ in the plane $\Pi$ such that $\Aline ou\in h$ and $\Aline ow\in v$. Then the based affine plane $(\Pi,uow)$ is horizontal-scale and vertical-scale. By Corollary~\ref{c:hv-scale<=>}, its ternar $\Delta$ is linear, associative-dot and right-distributive. Since $\{v\}\times h\subseteq\SC_\Pi$, the  based affine plane $(\Pi,uow)$ is translation, by Lemma~\ref{l:trans-swap} . By Theorem~\ref{t:VW-Thalesian<=>quasifield}, the ternar $\Delta$ is associative-plus. Therefore, the ternar $\Delta$ is linear right-distributive and associative.

Now assume that some ternar of the affine plane $\Pi$ is linear, right-distributive and associative. Then $\Pi$ has an affine base $uow$ whose ternar is linear, right-distributuve and associative. Consider the horizontal and vertical directions $h\defeq(\Aline ou)_\parallel$ and $v\defeq(\Aline ow)_\parallel)$ on the based affine plane $(\Pi,uow)$. By Corollary~\ref{c:vtrans+vscale}, the based affine plane $(\Pi,uow)$ is translation, horizontal-scale and vertical-scale and hence $(\{v\}\times h)\cup(\{h\}\times v)\subseteq\SC_\Pi$, witnessing that the plane $\Pi$ is of scale type at least $(\reshitka)$.
\smallskip

$(\Dscale)$ Assume that the plane $\Pi$ is of scale type at least $(\Dscale)$. Then $\{(\delta,L)\times\partial\Pi\times \mathcal L_\Pi:\delta\notin L\}\subseteq\SC_\Pi$. Choose any affine base $uow$ for the affine plane $\Pi$ and consider its horizontal and vertical directions  $h\defeq(\Aline ou)_\parallel$ and $v\defeq(\Aline ow)_\parallel$. Since $(\{h\}\times v)\cup(\{v\}\times h)\subseteq\SC_\Pi$, the ternar $\Delta$ of the based affine plane $(\Pi,uow)$ is linear, right-distributive and associative. Since $\{(\delta,\Aline ow):\delta\in\partial\Pi\setminus\{\Aline ow_\parallel\}\}\subseteq\SC_\Pi\}$, the based affine plane $(\Pi,uow)$ is vertical-shear, by  Lemma~\ref{l:scale=>shear}. By Proposition~\ref{p:ltring-ldist<=>}, the tenar $\Delta$ of the based affine plane $(\Pi,uow)$ is left-distributive and hence $\Delta$ is linear, distributive and associative.

Now assume conversely that some ternar of the plane $\Pi$ is linear, distributive and associative. Then $\Pi$ of scale type at least $(-)$ and $(\reshitka)$. The classification Theorem~\ref{t:scale-trace} ensures that $\Pi$ is of scale type (at least) $(\Dscale)$.
\end{proof}

\begin{problem} Find an algebraic characterization of affine planes of scale class at least $(\,\vert\,)$.
\end{problem}

\section{Affine planes of scale type $({\diagdown})$}\label{s:discale}

\rightline{\vbox{\hsize210pt 
\em \noindent\small Jak wiadomo, smok\'ow nie ma. Prymitywna ta konstatacja wystarczy mo\.ze umysłowi prostackiemu, ale nie nauce, poniewa\.z Wy\.zsza Szko\l a Neantyczna tym, co istnieje, wcale si\c e nie zajmuje; banalno\'s\'c istnienia zosta\l a ju\.z udowodniona zbyt dawno, by warto jej po\'swi\c eca\'c cho\'cby jedno jeszcze s\l owo.}}

\rightline{\index[person]{Lem}Stanis\l aw Lem\footnote{{\bf Stanis\l aw Lem} (1921 – 2006) was a Polish science fiction writer, philosopher, and futurist, born in Lw\'ow (now Lviv, Ukraine). Trained in medicine, he turned to writing after World War II and became one of the most influential speculative thinkers of the 20th century. His works, including Solaris, The Cyberiad, and His Master's Voice, blend science fiction with deep philosophical inquiry into consciousness, technology, and the limits of human understanding. Lem's writing, often satirical and intellectually rigorous, has been translated into over 40 languages and sold millions of copies. Though working under communist censorship, he earned global acclaim for his unique fusion of science, literature, and philosophy. He died in Kraków in 2006.}, Cyberiada}

\vskip20pt

In spite of the absence of known examples of affine planes of scale type $(\discale)$, we shall spend a lot of time and efforts studying properties of such (rather phantom) planes. Among many other results we shall prove Theorem~\ref{t:hdiscale<=>corps} characterizing corps as diagonal-scale horizontal-scale ternars which are commutative-plus or diassociative-plus. This theorem will be essentially used in the Lenz--Barlotti classification~\ref{t:Lenz-Barlotti} of projective planes. The results of this section are inspired by the papers of Pickert~\cite{Pickert1959}, Spencer-Yaqub \cite{Spencer1960}, \cite{Yaqub1967}  and J\'onsson\footnote{{\bf Wilbur Jacob J\'onsson} (1936 -- 2022) was a distinguished German‑Canadian mathematician born in Winnipeg by Icelandic immigrant parents. He earned his B.Sc. (1958) and M.Sc. (1959) in mathematics from the University of Manitoba, even competing on their Putnam team which placed fourth nationally. J\'onsson completed his doctorate (Dr. rer. nat.) at the Eberhard Karls Universit\"at T\"ubingen in 1963 under the supervision of G\"unter Pickert and Helmut Wielandt, with a thesis ``{\em Transitivit\"at und Homogenit\"at projektiver Ebenen}'' on transitivity and homogeneity in projective planes. Upon graduation, J\"onsson worked for a time at the University of Manitoba and at the University of Birmingham before joining McGill University as an Assistant Professor in 1966. He was promoted to the rank of Associate Professor in 1969 and retired after 48 years of service in 2014.} \cite{Jonsson1963}.

First we define ternars, coordinatizing affine planes of scale type $(\discale)$.

\begin{definition}\label{d:discale-ternar} A ternar is called a \defterm{$(\discale)$-ternar} if it is horizontal-scale and diagonal-scale but is not a corps.
\end{definition}

\begin{proposition}\label{p:discale<=>discale} An affine plane is of scale type $(\discale)$ if and only if some its ternar is a $(\discale)$-ternar.
\end{proposition}

\begin{proof} If an affine plane $\Pi$ is of scale type $(\discale)$, then its scale trace $\SC_\Pi$ is an injective function with $\dom[\SC_\Pi]=\partial\Pi\setminus\{v\}$ and $\rng[\SC_\Pi]=v$ for some direction $v\in\partial \Pi$. Choose any distinct pairs $(h,V)$ and $(d,E)$ in the scale trace $\SC_\Pi$ of $\Pi$. Then $V,E$ are two distinct lines in the direction $v$. Chose any points $o\in V$ and $e\in E$ such that $\Aline oe\in d$. Next, find unique points $u\in E$ and $w\in V$ such that $\{\Aline oh,\Aline we\}\subseteq h$. Then the based affine plane $(\Pi,uow)$ is horizontal-scale and diagonal-scale and so is its ternar $\Delta=\Aline oe$, by Theorem~\ref{t:diagonal-scale<=>}. Since $\SC_\Pi\ne \{(\delta,L)\in \partial\Pi\times\mathcal L_\Pi:\delta\notin L\}$, the plane $\Pi$ is not Desarguesian and hence its ternar $\Delta$ is not a corps.
\smallskip

Now assume that some ternar $R$ of an affine plane $\Pi$ is a $(\discale)$-ternar. Find an affine base $uow\in\Pi^3$ whose ternar  is isomorphic to the ternar $R$. Let $e$ be the diunit of the affine base $uow$. By Theorem~\ref{t:diagonal-scale<=>}, the based affine plane $(\Pi,uow)$ is horizontal-scale and diagonal-scale. Then the scale trace $\SC_\Pi$ of $\Pi$ contains the pairs $((\Aline ou)_\parallel,\Aline ow)$ and $((\Aline oe)_\parallel,\Aline ue)$. Since $R$ is not a corps, the plane $\Pi$ is not Desarguesian and $\Pi$ is not of type $(\Dscale)$. Since $\SC_\Pi$ contains two pairs $((\Aline ou)_\parallel,\Aline ow)$ and $((\Aline oe)_\parallel,\Aline ue)$ with $(\Aline ou)_\parallel\ne (\Aline oe)_\parallel$ and $\Aline ow\cap\Aline ue=\varnothing$, it is not of scale types $(\circ)$, $(\bullet)$, $(\didots)$, $(\,\vert\,)$, $(-)$, $(\reshitka)$, $(\invol)$, $(\Invol)$, so the unique possibility for $\Pi$ is to have the scale type $(\discale)$, by the classification Theorem~\ref{t:scale-trace}.
\end{proof}

For every element $x$ of an associative-dot ternar $R$, denote by $x^{-1}\in R$ the unique element such that $x{\cdot}x^{-1}{\cdot}x=x$ and $x^{-1}{\cdot}x{\cdot}x^{-1}=x^{-1}$. If $x\in R\setminus\{0\}$, then $x^{-1}$ is the inverse of $x$ in the group $(R\setminus\{0\},\cdot)$. If $x=0$, then $x^{-1}=0^{-1}=0$. 

For any element $a$ of a ternar $R$, we denote by $V_a$ the vertical line $\{a\}\times R$ in the coordinate plane $R^2$ of $R$ and by $a_\parallel$ the direction of the line $L_{a,0}=\{(x,y)\in R^2:y=x{\cdot}a\}$ with slope $a$. 

\begin{proposition}\label{p:scale-trace-of-diagonal} If $R$ is a $(\discale)$-ternar, then the scale trace $\SC_{R^2}$ of its coordinate plane $R^2$ is equal to the set $\{(x^{-1}_\parallel, V_{x}):x\in R\}$.
\end{proposition}

\begin{proof} Since $R$ is horizontal-scale and diagonal-scale, $\{(0_\parallel, V_0),(1_\parallel,V_1)\}\subseteq\SC_{R^2}$. By Theorem~\ref{t:horizontal-scale}, the coordinate plane $R^2$ of the ternar $R$ is horizontal-scale. Then for every $c\in R\setminus\{0\}$, there exists an automorphism $\Phi_c:R^2\to R^2$ such that $\Phi_c(1,0)=(c,0)$ and $V_0\subseteq\Fix(\Phi_c)$. By Proposition~\ref{p:hscale}, $\Phi_c(x,y)=(x{\cdot}c,y)$ for all $(x,y)\in R^2$. Observe that the line $\Phi_c[L_{1,0}]$ contains the points $\Phi_c(0,0)=(0,0)$ and $\Phi_c(1,1)=(c,1)$ and hence $\Phi_c[L_{1,0}]=L_{c^{-1},0}$ and $\bar\Phi_c(1_\parallel)=c^{-1}_\parallel$.  Lemma~\ref{l:SC-auto} ensures that $(c^{-1}_\parallel,V_c)=(\bar\Phi_c(1_\parallel),\Phi_c[V_1])\in\SC_{R^2}$ and hence $\{(c^{-1}_\parallel,V_c):c\in R\}\subseteq\SC_{R^2}$. Taking into account that $\SC_{R^2}$ is an injective function with domain $\partial\Pi\setminus\{v\}$ for some direction $v\in\partial R^2$, we conclude that $v\notin\{c^{-1}_\parallel:c\in R\}$ is the vertical direction and $\SC_{R^2}=\{(c^{-1}_\parallel,V_c):c\in R\}$.
\end{proof}

\begin{proposition}\label{p:diagonal11} If $R$ is a $(\discale)$-ternar, then
\begin{enumerate}
\item $c+x{\cdot}c\mp(x{\cdot}(a\circ c)+b)=(c\circ x)\cdot a+c{\cdot}a\mp(c+b)$ for all $a,b,c,x\in R$;
\item $c+x\cdot c\mp (x+b)=(c+x{\cdot}c\mp x)+b=c\circ x+b$ for all $c,b,x\in R$;
\item $c+c^{-1}=1$ and $c\circ c^{-1}=0$ for all $c\in R\setminus\{0,1\}$;
\item $c\circ c^{-1}\circ c=c$ and $c^{-1}\circ c\circ c^{-1}=c^{-1}$ for all $c\in R$;
\item $1+a=a+1$ and $(a+b)+1=a+(b+1)$ for all $a,b\in R$;
\item $c+(c^{-1}+b)=1+b$ for all $c\in R\setminus\{0,1\}$ and $b\in R$.
\item $c+x{\cdot}c\mp(x{\cdot}\alpha+b)=(c\circ x){\cdot} (\alpha\circ c^{-1})+c^{-1}\mp(c^{-1}{\cdot}\alpha+b)$ for all $\alpha,b,c,x\in R$ with $c\notin\{0,1\}$;
\item $c\cdot((c\cdot p)\circ c^{-1})=c+p\mp c^{-1}$ for all $c\in R\setminus\{0,1\}$ and $p\in R$.
\end{enumerate}
\end{proposition}

\begin{proof} Let $R$ be a $(\discale)$-ternar and $R^2$ be its coordinate plane. Since $R$ is horizontal-scale and diagonal-scale, so is its coordinate plane. Since $R$ is not a corps, the plane $R^2$ is not Desarguesian and hence $|R|>4$, by Corollary~\ref{c:4-Pappian}. Let $\Delta\defeq\{(x,y)\in R^2:x=y\}$ be the diagonal of the coordinate plane $R^2$. By Proposition~\ref{p:scale-trace-of-diagonal}, $\SC_{R^2}=\{(x^{-1}_\parallel,V_{x}):x\in R\}$.
\smallskip

\smallskip

1. The identity (1) follows from the Definition~\ref{d:diagonal-scale} of a diagonal-scale ternar.
\smallskip

2. For $a=1$, the identity (1) implies $c+x{\cdot}c\mp(x+b)=(c\circ x)\cdot 1+c\mp(c+b)=c\circ x+b=(c+x{\cdot}c\mp x)+b.$
\smallskip

3. Take any $a\in R\setminus\{0,1\}$ and observe that the line $L_{a^{-1},0}$ contains the points $(0,0)$, $(1,a^{-1})$ and $(a,1)$. Consider the line $L_{1,a\mp 1}$ containing the point $(a,1)$. Find a unique element $x\in R$ such that $(x,a^{-1})\in L_{1,a\mp 1}$. Since $L_{a^{-1},0}\nparallel L_{1,a\mp 1}$, the element $x$ is not equal to $1$.  Since $(1_\parallel,V_1)\in\SC_{R^2}$, there exists an automorphism $\Phi:R^2\to R^2$ such that $\Phi(a,1)=(x,a^{-1})$ and $V_1\subseteq\Fix(\Phi)$. Assumimg that $x\ne 0$ and taking into account that $(0_\parallel,V_0)\in\SC_{R^2}$, we can find an automorphism $\Psi:R^2\to R^2$ such that $\Psi(x,a^{-1})=(a,a^{-1})$ and $V_0\subseteq\Fix(\Psi)$. Then $\Psi\Phi[V_a]=\Psi[V_x]=V_a$. By Lemma~\ref{l:SC-auto}, $(\bar\Psi\bar\Phi(a^{-1}_\parallel),V_a)=(\bar\Psi\bar\Phi(a^{-1}_\parallel),\Psi[\Phi[V_a]])\in\SC_{R^2}$ and hence $\bar\Psi\bar\Phi(a^{-1}_\parallel)=a^{-1}_\parallel$. On the other hand, the line $\Phi[L_{a^{-1},0}]$ contains the points $\Phi(1,a^{-1})=(1,a^{-1})$ and $\Phi(a,1)=(x,a^{-1})$ and hence is horizontal. Since the automorphism $\Psi$ preserves the horizontal direction, the line $\Psi\Phi[L_{a^{-1},0}]$ belongs to the horizontal direction $0_\parallel$ and hence cannot belong to the direction $a^{-1}_\parallel$. This contradiction shows that $x=0$ and hence the point $(0,a^{-1})=(x,a^{-1})$ belongs to the line $L_{1,a\mp 1}$, which implies $a^{-1}=0{\cdot}1+(a\mp 1)=a\mp 1$ and finally, $a+a^{-1}=1$.
\smallskip

Next, we prove that $c\circ c^{-1}=0$. Writing the identity (2)  for $b=c$ and $x=c^{-1}$, we obtain 
$(c\circ c^{-1})+c=c+1\mp (c^{-1}+c)=c$ and hence $c\circ c^{-1}=0$, by the cancellativity of the plus operation.
\smallskip

4. If $c\in R\setminus\{0,1\}$, then the identities $c\circ c^{-1}\circ c=c$ and $c^{-1}\circ c\circ c^{-1}=c^{-1}$ follows from Proposition~\ref{p:diagonal11}(3). For $c\in\{0,1\}$ those identities follow from $0\circ 0\circ 0=0$ and $1\circ 1\circ 1=1$.
\smallskip

5. First we check that $c+(1+d)=1+(c+d)$ for all $c,d\in R$. This equality is clear if $c\in\{0,1\}$. So, assume that $c\notin\{0,1\}$. The item (3) ensures that $c^{-1}\circ c=0$. For $x\defeq a\defeq c^{-1}$ and  $b\defeq 1+d$, the identity in item (1) yields $c+d=1\mp (c+b)$, which is equivalent to $1+(c+d)=c+(1+d)$.  For $d=0$ the latter identity implies $1+c=c+1$. Then 
$(c+d)+1=1+(c+d)=c+(1+d)=c+(d+1)$.
\smallskip

6. Writing down the identity in item (2) for $x=c^{-1}$ and taking into account the item (5), we obtain the identity $b=c\circ c^{-1}+b=c+1\mp(c^{-1}+b)$, which implies $1+b=b+1=(c+1\mp(c^{-1}+b))+1=c+(1\mp(c^{-1}+b)+1)=c+(1+1\mp(c^{-1}+b))=c+(c^{-1}+b)$.
\smallskip



7. Given any elements $\alpha,b,c,x\in R$ with $c\notin \{0,1\}$, consider the element $a\defeq \alpha\circ c^{-1}$ and observe that $\alpha=a\circ c$. Writing down the identity (1) in Proposition~\ref{p:diagonal11} for $x=c^{-1}$, we obtain the equality
$c+1\mp(c^{-1}{\cdot}\alpha+b)=c{\cdot}a\mp(c+b)$. Adding $c^{-1}$ to both sides of this equality, we obtain
$$c^{-1}+c{\cdot}a\mp (c+b)=c^{-1}+(c+1\mp(c^{-1}{\cdot}\alpha+b)=1+1\mp(c^{-1}{\cdot}\alpha+b)=c^{-1}{\cdot}\alpha+b.$$ Then the identity (1) yields the identity $c+x{\cdot}c\mp(x{\cdot}\alpha+b)=(c\circ x){\cdot} (\alpha\circ c^{-1})+c^{-1}\mp(c^{-1}{\cdot}\alpha+b)$ appearing in the condition (7).
\smallskip

8. Choose any elements $c\in R\setminus\{0,1\}$ and $p\in R$. Writing down the identity (7) of  Proposition~\ref{p:diagonal11} for $x=0$, $\alpha=c\cdot p$ and $b=p\mp c^{-1}$, we obtain the identity
$$
\begin{aligned}
c+p\mp c^{-1}&=c+b=c\cdot((c\cdot p)\circ c^{-1})+c^{-1}\mp(c^{-1}{\cdot}\alpha+b)\\
&=c\cdot((c\cdot p)\circ c^{-1})+c^{-1}\mp(p+p\mp c^{-1})\\
&=c\cdot((c\cdot p)\circ c^{-1})+c^{-1}\mp c^{-1}=c\cdot((c\cdot p)\circ c^{-1})
\end{aligned}
$$appearing in the condition (8).
\end{proof}

\begin{remark} By Proposition~\ref{p:diagonal11}(4), for each element $x$ of a $(\discale)$-ternar $R$, its inverse $x^{-1}$ in the dot semigroup $(R,\cdot)$ coincides with its inverse in the dit semigroup $(R,\circ)$. So, we can use the same notation $x^{-1}$ for those two (coinciding) inverse elements.
\end{remark}


\begin{proposition}\label{p:diagonal-alpha} Let $R$ be a $(\discale)$-ternar. For every $c\in R\setminus\{0,1\}$, there exists an isomorphism $\alpha_c:(R,\cdot)\to (R,\circ)$ of the semigroups $(R,\cdot)$ and $(R,\circ)$ such that  
$$
\begin{aligned}
\alpha_c(x)&\defeq (c\circ (c^{-1}\cdot x^{-1}))\cdot c^{-1}=((c\circ x^{-1})\cdot c^{-1})\circ c\\
&=c\cdot ((x^{-1}\cdot c)\circ c^{-1})=c^{-1}\circ (c\cdot ((x^{-1}\circ c^{-1}))\\
&=((c\circ ((c^{-1}\circ (z\cdot c))\cdot x^{-1}))\cdot c^{-1})\circ z^{-1}\\
&=z\circ (c\cdot ((x^{-1}\cdot ((c^{-1}\cdot z^{-1})\circ c)\circ c^{-1}))
\end{aligned}
$$
for all $x\in R$ and $z\in R\setminus\{1\}$.
\end{proposition}

\begin{proof} For every $c\in R$, consider the group $$\Aut_c(R^2)=\{\Phi\in\Aut(R_2):\bar\Phi(c^{-1}_\parallel)=c^{-1}_\parallel\;\wedge\;V_{c}\subseteq\Fix(\Phi)\}$$ of hyperscales with vertical axis $V_{c}\defeq\{c\}\times R$ and central direction $c^{-1}_\parallel\defeq(L_{c^{-1},0})_\parallel$. Since $(c^{-1}_\parallel,V_c)\in \SC_\Pi$, the group $\Aut_c(R^2)$ acts transitively on each punctured line $L\setminus V_c$ of direction $c^{-1}_\parallel$. 

Since the ternar $R$ is horizontal-scale, for every $a\in R\setminus\{0\}$, the map $\varphi_a:R^2\to R^2$, $\varphi_a:(x,y)\mapsto (x{\cdot}a,y)$, is an automorphism of the plane $R^2$, according to Proposition~\ref{p:hscale}. Moreover, the map $\varphi:R\setminus\{0\}\to\Aut_0(R^2)$ assigning to each element $a\in R\setminus \{0\}$ the automorphism $\varphi_{a^{-1}}:R^2\to R^2$, $\varphi_{a^{-1}}:(x,y)\mapsto (x{\cdot}a^{-1},y)$, is an isomorphism of the dot group $(R\setminus\{0\},\cdot)$ onto the group $\Aut_0(R^2)$. Observe that $\varphi_{a^{-1}}[V_a]=V_1$. 
Since $\varphi$ is a group isomorphism, $\varphi_{a}^{-1}=\varphi_{a^{-1}}$. 

Since the ternar $R$ is diagonal-scale, the map $\psi:R\setminus\{1\}\to\Aut_1(R^2)$ assigning to each element $a\in R\setminus \{1\}$ the automorphism $\psi_a:R^2\to R^2$, $\psi_a:(x,y)\mapsto (a\circ x,a+x{\cdot}a\mp y)$, is an isomorphism of the dit group $(R\setminus\{1\},\circ)$ onto the group $\Aut_1(R^2)$, by Proposition~\ref{p:diagonal-scale<=>} and the definition of the dit operation. 

For every $c\in R\setminus\{0,1\}$, the isomorphisms $\varphi_{c}:R^2\to R^2$ and $\psi_c:R^2\to R^2$ induce the group isomorphisms 
$$
\begin{aligned}
&\Phi_{c}:\Aut_1(R^2)\to \Aut_c(R^2),\;\Phi_c:F\mapsto \varphi_{c}F\varphi_c^{-1},\quad\mbox{and}\\
&\Psi_c:\Aut_0(R^2)\to \Aut_c(R^2),\;\Psi_c:F\mapsto \psi_cF\psi_c^{-1}.
\end{aligned}
$$
 Then the composition $\psi^{-1}\Phi_c^{-1}\Psi_c \varphi:R\setminus\{0\}\to R\setminus\{1\}$ is an isomorphism of the groups $(R\setminus\{0\},\,\cdot\,)$ and $(R\setminus\{1\},\circ)$. Let $\alpha_c:(R,\cdot)\to (R,\circ)$ be the unique extension of the isomorphism $\psi^{-1}\Phi_c^{-1}\Psi_c \varphi$ such that $\alpha_c(0)=1$. 
 
 For every $x\in R\setminus\{0\}$, let $y\defeq\alpha_c(x)=\psi^{-1}\Phi_c^{-1}\Psi_c(\varphi_{x^{-1}})$ and observe that $\Phi_c(\psi_y)=\Psi_c(\varphi_{x^{-1}})$ and hence
$$\varphi_c\psi_y\varphi_c^{-1}=\Phi_c(\psi_y)=\Psi_c(\varphi_{x^{-1}})=\psi_c\varphi_{x^{-1}}\psi_{c^{-1}}.$$ Then for every $z\in R$ we have $\varphi_c\psi_y\varphi_{c^{-1}}[V_{z{\cdot}c}]=\psi_{c}\varphi_{x^{-1}}\psi_{c^{-1}}[V_{z{\cdot}c}]$, which implies 
$$V_{(y\circ z){\cdot}c}=\varphi_c\psi_y\varphi_{c^{-1}}[V_{z{\cdot}c}]=\psi_{c}\varphi_{x^{-1}}\psi_{c^{-1}}[V_{z{\cdot}c}]=V_{c\circ ((c^{-1}\circ (z\cdot c))\cdot x^{-1})}$$ and $(y\circ z)\cdot c=c\circ ((c^{-1}\circ (z\cdot c))\cdot x^{-1})$. Then 
$$
\alpha_c(x)=y=((c\circ ((c^{-1}\circ (z\cdot c))\cdot x^{-1}))\cdot c^{-1})\circ z^{-1},
$$ 
for all $z\in R\setminus\{1\}$. For $z=0$ and $z=c^{-1}$, the latter identity implies  
$$\alpha_c(x)=(c\circ (c^{-1}\cdot x^{-1}))\cdot c^{-1}=((c\circ x^{-1})\cdot c^{-1})\circ c.$$
For $x=0$, the identity
$$\alpha_c(0)=1=((c\circ ((c^{-1}\circ (z\cdot c))\cdot 0^{-1}))\cdot c^{-1})\circ z^{-1}$$ also holds.

Since $\alpha_c:(R,\cdot)\to (R,\circ)$ is an isomorphism of inverse semigroups, 
$$\alpha_c(x)=(\alpha_c(x^{-1}))^{-1}=(((c\circ ((c^{-1}\circ (z\cdot c))\cdot x))\cdot c^{-1})\circ z^{-1})^{-1}=
z\circ (c\cdot ((x^{-1}\cdot ((c^{-1}\cdot z^{-1})\circ c)\circ c^{-1}))$$
for all $x\in R$ and $z\in R\setminus\{1\}$. For $z=0$ and $z=c^{-1}$, the latter identity implies 
$$\alpha_c(x)=c\cdot ((x^{-1}\cdot c)\circ c^{-1})=
c^{-1}\circ (c\cdot ((x^{-1}\circ c^{-1})).$$
\end{proof}

For any elements $c,x$ of a $(\discale)$-ternar $R$, let
$$x^c\defeq c+x^{-1}\mp c^{-1}.$$

\begin{proposition}\label{p:di-conjugation} Let $R$ be a $(\discale)$-ternar and $c\in R\setminus\{0,1\}$. Then 
\begin{enumerate}
\item $x^c=c\cdot((c\cdot x^{-1})\circ c^{-1})$ for all $x\in R$;
\item $(x\cdot y)^c=x^c\circ y^c$  and $(x\circ y)^c=x^c\cdot y^c$ for all $x,y\in R$;
\item $(x^c)^{c^{-1}}=x$ and $(x^c)^{-1}=(x^{-1})^c$ for all $x\in R$;
\item $x+c=c+x^{c^{-1}}$ for all $x\in R\setminus\{0,1\}$;
\item $x^c+c=c+x$, and $x^c=(c+x)-c$ for all $x\in R\setminus\{0,1\}$;
\item $c^c=c^{c^{-1}}=c$.
\end{enumerate} 
\end{proposition}

\begin{proof} 1. By Propositions~\ref{p:diagonal11}(8) and \ref{p:diagonal-alpha}, for every $x\in R$, we have $$x^c=c+x^{-1}\mp c^{-1}=c\cdot ((c\cdot x^{-1})\circ c^{-1})=\alpha_{c}(c\cdot x\cdot c^{-1}).$$ 
\smallskip

2. It is clear that the map $\beta_c:R\to R$, $\beta_c:x\mapsto c\cdot x\cdot c^{-1}$, is an automorphism of the semigroup $(R,\cdot)$. Since $\alpha_{c}:(R,\cdot)\to (R,\circ)$ is an isomorphism of the semigroups $(R,\cdot)$ and $(R,\circ)$, so is the map $\gamma_c\defeq \alpha_{c}\beta_c:(R,\cdot)\to(R,\circ)$. Then $$(x\cdot y)^c=\gamma_c(x\cdot y)=\gamma_c(x)\circ \gamma_c(y)=x^c\circ y^c$$ for all $x,y\in R$. 

Since the map $\gamma_{c}:(R,\cdot)\to (R,\circ)$ is an isomorphism of the semigroups $(R,\cdot)$ and $(R,\circ)$, the inverse map $\gamma_{c}^{-1}$ is an isomorphism of the semigroups $(R,\circ)$ and $(R,\cdot)$. For every $x\in R$ and $y\defeq\gamma_{c}(x)=x^{c}=c\cdot ((c\cdot x^{-1})\circ c^{-1})$, we have $x=\gamma_{c}^{-1}(y)=(c^{-1}\cdot((c^{-1}\cdot y)\circ c))^{-1}=(c^{-1}\circ (y^{-1}\cdot c))\cdot c=\alpha_{c^{-1}}(c^{-1}{\cdot} y{\cdot} c)=\alpha_{c^{-1}}\beta_{c^{-1}}(y)=\gamma_{c^{-1}}(y)=y^{c^{-1}}$, by Proposition~\ref{p:diagonal-alpha}. Taking into account that  $\gamma_{c}^{-1}:(R,\circ)\to (R,\cdot)$ is an isomorphism of the semigroups $(R,\circ)$ and $(R,\cdot)$, we conclude that
$$(x\circ y)^{c^{-1}}=\gamma_{c^{-1}}(x\circ y)=\gamma_{c^{-1}}(x)\cdot\gamma_{c^{-1}}(y)=x^{c^{-1}}\cdot y^{c^{-1}}$$ for all $x,y\in R$.
\smallskip

3. Observe also that $(x^c)^{c^{-1}}=\gamma_c^{-1}(\gamma_c(x))=x$ for all $x\in R$.

Since $\gamma_c:(R,\cdot)\to(R,\circ)$ is an isomorphism, for every 
$x\in R$ we have
$$(x^c)^{-1}=(c+x^{-1}\mp c^{-1})^{-1}=\gamma_c(x)^{-1}=\gamma_c(x^{-1})=c+x\mp c^{-1}=(x^{-1})^c.$$
\smallskip

4. Now assume that $x\in R\setminus\{0,1\}$. By definition, $x^{c^{-1}}=c^{-1}+x^{-1}\mp c$. Adding $c$ to both sides of this equality, we obtain the equality
$$c+x^{c^{-1}}=c+(c^{-1}+x^{-1}\mp c)=1+x^{-1}\mp c=(x+x^{-1})+x^{-1}\mp c=x+(x^{-1}+x^{-1}\mp c)=x+c,$$
by Proposition~\ref{p:diagonal11}(6,3).

5. For every $x\in R\setminus\{0,1\}$, two preceding items imply $c+x=c+(x^c)^{c^{-1}}=x^c+c$ and hence $x^c=(c+x)-c$.
\smallskip

6. It follows from $c^c+c=c+c=c+c^{c^{-1}}$ that $c^c=c=c^{c^{-1}}$.
\end{proof}

For any element $x$ of a ternar $R$, define its powers by the recursive formula: $x^0=1$ and $x^{n+1}=x^n\cdot x$ for $n\in\w$. Also put $x^{-n}=(x^{-1})^n$ for all $n\in\IN$.

\begin{proposition}\label{p:diagonal-n-power} For any element $c\in R\setminus\{0,1\}$ of a $(\discale)$-ternar and every $n\in\w$ we have the identity
$c^{n+1}=c\circ ((((c^{-1}\circ c^n)\cdot c)\circ c^{-1})\cdot c^{-1})$.
\end{proposition}

\begin{proof} Proposition~\ref{p:diagonal-alpha} implies that 
$(c^{-1}\circ (c\cdot x))\cdot c=((c^{-1}\circ x)\cdot c)\circ c^{-1}$ for all $x\in R$. In particular, for $x=c^n$, we obtain the identity
$(c^{-1}\circ c^{n+1})\cdot c=((c^{-1}\circ c^n)\cdot c)\circ c^{-1}$ implying
$$c^{n+1}=c\circ ((((c^{-1}\circ c^n)\cdot c)\circ c^{-1})\cdot c^{-1}).$$
\end{proof}

Let $R$ be a $(\discale)$-ternar. For every $x\in R$ define its powers in the semigroup $(R,\circ)$ by the recursive formula: $x^{\circ 0}=0$ and $x^{\circ(n+1)}=x^{\circ n}\circ x$ for $n\in\w$. Also put $x^{\frac1n}\defeq x^{\circ n}$ and $x^{\frac{m}n}\defeq (x^m)^{\circ n}$ for all $m\in\IZ$ and $n\in\IN$.

\begin{proposition}\label{p:discale=>divisible} If $R$ is a $(\discale)$-ternar, then $$(c^{\frac1n})^n=(c^{\circ n})^n=c=(c^n)^{\circ n}=(c^n)^{\frac1n}$$for all $n\in\IN$ and $c\in R$. Moreover $x^n\ne 1$ for all $x\in R\setminus\{1\}$.
\end{proposition}

\begin{proof} This theorem will be  proved by induction on $n\in\IN$. For $n=1$, the identity $(c^{\circ n})^{n}=c=(c^n)^{\circ n}$ is trivial. Assume that for some $n\in \IN$ and all positive $k\le n$ we know that $(c^{\circ k})^k=c=(c^k)^{\circ k}$ for all $c\in R$, and $x^k\ne 1$ for all $x\in R\setminus\{1\}$.  We have to prove that $(c^{\circ(n+1)})^{n+1}=c=(c^{n+1})^{\circ(n+1)}$. If $c\in\{0,1\}$, then this equality is clear.  So, assume that $c\notin\{0,1\}$. By Proposition~\ref{p:diagonal-n-power} and the inductive assumption,
$$
\begin{aligned}
c^{n+1}&=c\circ ((((c^{-1}\circ c^n)\cdot c)\circ c^{-1})\cdot c^{-1})=
c\circ (((((c^{-n})^{\circ n}\circ c^n)\cdot c)\circ c^{-1})\cdot c^{-1})\\
&=c\circ ((((c^{-n})^{\circ (n-1)}\cdot c)\circ c^{-1})\cdot c^{-1})=
c\circ (((c^{\frac{-n}{n-1}}\cdot c)\circ c^{-1})\cdot c^{-1})\\
&=c\circ ((c^{\frac{-1}{n-1}}\circ c^{-1})\cdot c^{-1})=
c\circ (((c^{-1})^{\circ(n-1)}\circ c^{-1})\cdot c^{-1})\\
&=c\circ ((c^{-1})^{\circ n}\cdot c^{-1})=
c\circ (c^{\frac{-1}n}\cdot c^{-1})=c\circ c^{\frac{-(n+1)}n}=c\circ (c^{-(n+1)})^{\circ n}
\end{aligned}
$$
which implies
$(c^{n+1})^{\circ(n+1)}=c$. 

Applying to the equality $(c^{n+1})^{\circ(n+1)}=c$ the homomorphism $x^c$ and taking into account that $c^c=c$, we conclude that $c=c^c=((c^{n+1})^{\circ(n+1)})^c=((c^c)^{\circ(n+1)})^{n+1}=(c^{\circ(n+1)})^{n+1}$.

For any $x\in R\setminus\{0\}$ with $x^n=1$ we have
$x=(x^n)^{\circ n}=1^{\circ n}=1$, which means that the group $(R\setminus\{0\},\cdot)$ has no elements of order $n$.
\end{proof}

Let us recall that a group $G$ is \defterm{divisible} if for every $x\in G$ and $n\in\IN$ there exists an element $y\in G$ such that $y^n=x$. Proposition~\ref{p:discale=>divisible} imply the following corollary.

\begin{corollary}\label{c:finite-discale=>divisible} The dot and dit groups of every $(\discale)$-ternar are divisible and have no elements of finite order.
\end{corollary}

Since nontrivial divisible groups without elements of finite order are infinite, Corollary~\ref{c:finite-discale=>field} and Wedderburn---Witt Theorem~\ref{t:Wedderburn-Witt} imply the following known characterization of finite fields.  

\begin{corollary}\label{c:finite-discale=>field} A finite ternar is a field if and only if it is horizontal-scale and diagonal-scale. 
\end{corollary}

\begin{remark} Corollary~\ref{c:finite-discale=>field} is known. It was proved by complicated Group Theoretic tools 60 years ago by combined efforts of J\'onsson \cite{Jonsson1963}, L\"uneburg \cite{Luneburg1964}, Cofman \cite{Cofman1966}, Yaqub \cite{Yaqub1967} and Dembowski \cite[\S4.37]{Dembowski}. Our proof of Corollary~\ref{c:finite-discale=>field} is much simple and natural.
\end{remark}

\begin{proposition}\label{p:discale-fraction} Let $R$ be a $(\discale)$-ternar and $c\in R\setminus\{0,1\}$ be any element.
Then
\begin{enumerate}
\item $(c^{\frac nm})^{c^{-1}}=(c^{\frac nm})^c=c^{\frac mn}\ne c^{\frac nm}$ and $c^{\frac nm}+c=c+c^{\frac mn}$ for any distinct numbers $n,m\in\IN$;
\item $c\cdot c\ne c\circ c$.
\end{enumerate}
\end{proposition}

\begin{proof} 1. Fix any distinct numbers $n,m\in\IN$. By Proposition~\ref{p:di-conjugation}(1), $$
\begin{aligned}
(c^\frac nm)^{c^{-1}}&=c^{-1}\cdot ((c^{-1}\cdot c^{\frac{-n}{m}})\circ c)=
c^{-1}\cdot(c^{\frac{-(m+n)}{m}}\circ c)\\
&=c^{-1}\cdot (c^{-(m+n)})^{\circ m}\circ (c^{(m+n)})^{\circ (m+n)})\\
&=c^{-1}\cdot (c^{(m+n)})^{\circ n}=c^{-1}\cdot c^{\frac{(m+n)}n}=c^{\frac mn}.
\end{aligned}
$$ 
Applying Proposition~\ref{p:di-conjugation}(6) and the (already proved) equality $(c^{\frac mn})^{c^{-1}}=c^{\frac nm}$, we obtain the identity
$$c^{\frac mn}=((c^{\frac mn})^{c^{-1}})^c=(c^{\frac nm})^c.$$

Assuming that $c^{\frac nm}=c^{\frac mn}$, we conclude that $c^{n^2m}=(c^{\frac nm})^{nm}=(c^{\frac mn})^{nm}=c^{m^2n}$ and hence $c^{(n-m)nm}=1$, which contradicts Proposition~\ref{p:discale=>divisible}.

By Proposition~\ref{p:di-conjugation}(4), $c^{\frac nm}+c=c+(c^{\frac nm})^{c^{-1}}=c+c^{\frac mn}$.
\smallskip

2. To derive a contradiction, assume that $c\cdot c=c\circ c$. Then $c^2=c\circ c=c^{\frac12}$ and hence $c^3=1$, which contradicts Proposition~\ref{p:discale=>divisible}.
\end{proof}

Now we deduce some properties of the inverse element $\bar 1$ to $1$ in the plus loop of a $(\discale)$-ternar.

\begin{proposition}\label{p:discale-1} For any $(\discale)$-ternar $R$, the element $\bar 1\defeq 1\mp 0$ has the following properties:
\begin{enumerate}
\item $\bar 1+1=0=1+\bar 1$;
\item $(x+y)+\bar 1=x+(y+\bar 1)$ for all $x,y\in R$;
\item $c\mp x=(c^{-1}+x)+\bar 1=c^{-1}+(x+\bar 1)$ for all $c\in R\setminus\{0,1\}$ and $x\in R$;
\item $c+\bar 1=c^{-1}\mp 0$ for all $c\in R\setminus\{0,1\}$;
\item $a^c=c+a^{-1}\mp c^{-1}=c+(a+c\mp 0)$ for all $a,c\in R\setminus\{0,1\}$.
\item $1\mp x=x+\bar 1$ for all $x\in R$;
\end{enumerate}
\end{proposition}

\begin{proof} 1. The definition of $\bar 1\defeq 1\mp 0$ and Proposition~\ref{p:diagonal11}(5) imply $\bar 1+1=1+\bar 1=0$,
\smallskip

2. By Proposition~\ref{p:diagonal11}(5), $$((x+y)+\bar 1)+1=(x+y)+(\bar 1+1)=x+y=x+(y+(\bar 1+1))=x+((y+\bar 1)+1)=(x+(y+\bar 1))+1,$$ which implies $(x+y)+\bar 1=x+(y+\bar 1)$, by the cancellativity of the plus loop $(R,+)$.
\smallskip
 
3. Given any $c\in\{0,1\}$ and $x\in R$, consider the element $y\defeq c\mp x$ and observe that $c+y=x$. Adding to both sides $c^{-1}$ and applying Proposition~\ref{p:diagonal11}(5,6), we conclude that $$y+1=1+y=c^{-1}+(c+y)=c^{-1}+x$$ and hence $c\mp x=y=(y+1)+\bar 1=(c^{-1}+x)+\bar 1=c^{-1}+(x+\bar 1)$.
\smallskip

4. By the preceding item, $c^{-1}\mp 0=c+(0+\bar 1)=c+\bar 1$ for all $c\in R\setminus\{0,1\}$.
\smallskip

5. By the preceding items, for every $a,c\in R\setminus\{0,1\}$, we have
$$a^c=c+a^{-1}\mp c^{-1}=c+(a+(c^{-1}+\bar 1))=c+(a+c\mp 0).$$
\smallskip

6. For every  $x\in R$, Propositions~\ref{p:discale-1}(1,2) and \ref{p:diagonal11}(5) imply
$$1\mp x=1\mp x+(1+\bar 1)=(1\mp x+1)+\bar 1=(1+1\mp x)+\bar 1=x+\bar 1.$$
\end{proof}





\begin{proposition}\label{p:xab=b+a+x-(b+a)} If $R$ is a $(\discale)$-ternar, then
\begin{enumerate}
\item $c+(a+b)=(c+(a+c\mp 0))+(c+b)$ for all $a,b,c\in R$;
\item $c+(a+b)=a^c+(c+b)$ for all $a,c\in R\setminus\{0,1\}$ and $b\in R$.
\item $(x^a)^b=(b+(a+x))-(b+a)$ for all $a,b\in R\setminus\{0,1\}$ and $x\in R$.
\end{enumerate}
\end{proposition}

\begin{proof} 1. Let $a,b,c\in R$ be any elements. We have to prove the equality $$c+(a+b)=(c+(a+c\mp 0))+(c+b).$$ This equality is trivially true if $a=0$ or $c=0$. So, assume that $a\ne 0\ne c$. If $a=1$, then 
$$c+(1+c\mp 0)=c+(c\mp 0+1)=(c+c\mp 0)+1=0+1=1,$$ by Proposition~\ref{p:diagonal11}(5), and hence
$$(c+(1+c\mp 0))+(c+b)=1+(c+b)=(c+b)+1=c+(b+1)=c+(1+b).$$
If $c=1$, then 
$$c+(a+c\mp 0)=1+(a+1\mp 0)=(a+1\mp 0)+1=a+(1\mp 0+1)=a+(1+1\mp 0))=a+0=a,$$ by Proposition~\ref{p:diagonal11}(5), and hence 
$$(c+(a+c\mp 0))+(c+b)=a+(1+b)=a+(b+1)=(a+b)+1=1+(a+b)=c+(a+b).$$
So, assume that $a,c\notin\{0,1\}$.
Writing down the identity (7) in Proposition~\ref{p:diagonal11} for $\alpha=0$, and applying Proposition~\ref{p:diagonal11}(3), we obtain the identity
$$c+x{\cdot}c\mp b=(c\circ x)\cdot c^{-1}+c^{-1}\mp b.$$ 
By Propositions~\ref{p:discale-1}(1,3), \ref{p:diagonal11}(5), for $x\defeq a^{-1}{\cdot}c^{-1}$ we have 
$$(c+x{\cdot} c\mp b)+1=(c+a^{-1}\mp b)+1=(c+(a+b+\bar 1))+1=c+(a+b).$$ On the other hand, Propositions~\ref{p:di-conjugation}(1,3) and \ref{p:discale-1}(5) imply
$$
\begin{aligned}
&((c\circ x)\cdot c^{-1}+c^{-1}\mp b)+1=
(c\circ (a^{-1}\cdot c^{-1}))\cdot c^{-1}+(c+b+\bar 1))+1\\
&=(c\cdot ((c\cdot a)\circ c^{-1})^{-1}+((c+b+\bar 1)+1)
=((a^{-1})^c)^{-1}+((c+b+\bar 1)+1)\\
&=a^c+((c+b+\bar 1)+1)=(c+a^{-1}\mp c^{-1})+(c+b)=(c+(a+c\mp 0))+(c+b).
\end{aligned}
$$
Then $c+(a+b)=(c+x\cdot c\mp b)+1=((c\circ x)\cdot c^{-1}+c^{-1}\mp b)+1=(c+(a+c\mp 0))+(c+b)$.
\smallskip

2. For any $a,c\in R\setminus\{0,1\}$ and $b\in R$, Propositions~\ref{p:Yaqub1967}(1) and \ref{p:discale-1}(5) imply 
$$c+(a+b)=(c+(a+c\mp 0))+(c+b)=a^c+(c+b).$$

3. For every $a,b,x\in R\setminus\{0,1\}$, the preceding item implies 
$(x^a)^b+(b+a)=b+(x^a+a)=b+(a+x)$ and hence 
$(x^a)^b=(b+(a+x))-(b+a)$. Let us show that this equality also holds for $x\in\{0,1\}$. If $x=0$, then $(x^a)^b=(0^a)^b=1^b=0=(b+(a+0))-(b+a)=(b+(a+x))-(b+a)$.
If $x=1$, then 
$$
\begin{aligned}
(x^a)^b&=(1^a)^b=0^b=1=(1+(b+a))-(b+a)=((b+a)+1)-(b+a)\\
&=(b+(a+1))-(b+a)=(b+(a+x))-(b+a),
\end{aligned}
$$
by Proposition~\ref{p:diagonal11}(5).
\end{proof}

Now our aim is to prove that $1+1=0$ in $(\discale)$-ternars. This (nontrivial) fact was proved by Yaqub (Spencer) in \cite{Yaqub1967}.

\begin{proposition}[Yaqub, 1967]\label{p:Yaqub1967} If $R$ is a $(\discale)$-ternar, then
\begin{enumerate}
\item for any $a,b\in R\setminus\{0,1\}$, the function $$\psi_{b,a}:R^2\to R^2,\quad \psi_{b,a}:(x,y)\mapsto ((b+(a+x))-(b+a),b+(a+y)),$$ is an automorphism of the plane $R^2$;
\item $1+1=0$ in $R$;
\item the commutative center $\{x\in R:\forall a\in R\;\;x+a=a+x\}$ of the loop $(R,+)$ equals $\{0,1\}$;
\item the right nucleus  $\{x\in R:\forall a,b\in R\;\;a+(b+x)=(a+b)+x\}$ of $(R,+)$ equals $\{0,1\}$.
\end{enumerate}
\end{proposition}

\begin{proof} 1. For every $c\in R\setminus\{0,1\}$, consider the horizontal scale $$h_c:R^2\to R^2,\quad h_c:(x,y)\mapsto (x\cdot c^{-1},y),$$and the diagonal scale
$$d_c:R^2\to R^2,\quad d_c:(x,y)\mapsto (c\circ x,c+x{\cdot}c\mp y).$$
By Propositions~\ref{p:hscale} and \ref{p:diagonal-scale<=>}, the maps $h_c$ and $d_c$ are automorphisms of the affine plane $R^2$. Then the composition $\varphi_c\defeq h_cd_ch_c$ also is an automorphism of the plane $R^2$. By Proposition~\ref{p:di-conjugation}(1),
$$
\begin{aligned}
\varphi_c(x,y)&=\big((c\circ (x\cdot c^{-1}))\cdot c^{-1},c+x\mp y\big)=
\big((c\cdot(c\cdot x^{-1})\circ c^{-1}))^{-1},c+x\mp y\big)\\
&=((x^{c})^{-1},c+x\mp y)=((x^{-1})^{c},c+x\mp y)=(c+x\mp c^{-1},c+x\mp y).
\end{aligned}
$$

Given any elements $a,b\in R\setminus\{0,1\}$, consider the automorphism $\varphi_{b,a}\defeq \varphi_b\varphi_a$ of the plane $R^2$ and observe that 
$$
\varphi_{b,a}(x,y)=(b+(a+x\mp a^{-1})\mp b^{-1},b+(a+x\mp a^{-1})\mp (a+x\mp y)).$$

Proposition~\ref{p:xab=b+a+x-(b+a)}(3) ensures that $b+(a+x\mp a^{-1})\mp b^{-1}=b+(x^{-1})^a\mp b^{-1}=b+(x^a)^{-1}\mp b^{-1}=(x^a)^b=(b+(a+x))-(b+a)$. Next, we prove that $b+(a+x\mp a^{-1})\mp (a+x\mp y)=(b+(a+y))+\bar 1$ for all $x,y\in R$. If $x\notin\{0,1\}$, then  $a+x^{-1}\mp a^{-1}=x^a\notin\{0,1\}$ and hence 
$$
\begin{aligned}
&b+(a+x\mp a^{-1})\mp (a+x\mp y)=b+(x^{-1})^a\mp (a+x\mp y)=b+\big((x^a+(a+x\mp y))+\bar 1\big)\\
&=\big(b+(a+(x+x\mp y)\big)+\bar 1=
(b+(a+y))+\bar 1,
\end{aligned}
$$by Propositions~\ref{p:discale-1}(3) and \ref{p:xab=b+a+x-(b+a)}(2). 
If $x=0$, then 
$$b+(a+x\mp a^{-1})\mp(a+x\mp y)=b+(a+a^{-1})\mp (a+y)=b+1\mp(a+y)=b+((a+y)+\bar 1)=(b+(a+y))+\bar 1,$$by Propositions~\ref{p:diagonal11}(3) and \ref{p:discale-1}(2,6).
If $x=1$, then
\begin{multline*}b+(a+1\mp a^{-1})\mp(a+1\mp y)=b+(a\circ a^{-1})\mp (a+(y+\bar 1))=b+((a+y)+\bar 1)=(b+(a+y))+\bar 1,
\end{multline*}
by Propositions~\ref{p:diagonal11}(3) and \ref{p:discale-1}(2,6).
Therefore, 
$$\varphi_{b,a}(x,y)=((b+(a+x))-(b+a),(b+(a+y))+\bar 1)$$for all $x,y\in R$.

By Proposition~\ref{p:diagonal11}(5), the function $T:R^2\to R^2$, $T:(x,y)\mapsto (x,y+1)$ is a translation of the coordinate plane $R^2$. Then the function $$\psi_{b,a}=T\varphi_{b,a}=T\varphi_b\varphi_a:(x,y)\mapsto ((b+(a+x))-(b+a),b+(a+y)),$$ is an automorphism of the plane $R^2$.
\smallskip


2. To derive a contradiction, assume that $2\defeq 1+1\ne 0$. Let $R'$ be the smallest subternar of $R$ containing the elements $0$ and $1$. Then $2=1+1\in R'$. Since the multiplicative group $(R\setminus\{0\},\cdot)$ contains no elements of finite order, the subternar $R'$ is infinite. Given any element $b\in R\setminus\{0,1,\bar 1\}$, consider the element $a\defeq b\mp 0$ and observe that $a\in R\setminus\{0,1\}$. Since $b+a=b+b\mp 0=0$, the automorphism $\psi_{b,a}$ has very simple form: $$\psi_{b,a}(x,y)=(b+(a+x),b+(a+y))\quad\mbox{for any $(x,y)\in R^2$}.$$Since $b+(a+x)=x$ for $x\in\{0,1\}$, the set $\Fix(\psi_{b,a})$ contains the affine subplane $R'\times R'$ of $R$. 
Then, $b+(b\mp 0+x)=x$ for all $x\in R'$ and $b\in R\setminus\{0,1,\bar 1\}$. Since $1+1\ne 0$, for every $b\in R\setminus\{\bar 1\}$ we have $b\mp 0\ne 1$ and hence $b+(1+x)\ne b+(b\mp 0+x)=x$. Then $\bar 1+(1+x)=x=(x+1)+\bar 1$ for all $x\in R'$. Therefore, the element $c\defeq \bar 1\in R\setminus\{0,1\}$ belongs to the commutative center of the loop $(R',+)$. By Proposition~\ref{p:di-conjugation}(5), $x+c=c+x=x^c+c$ and hence $x^c=x$ for all $x\in R'\setminus\{0,1\}$. In particular, 
$$c\cdot c=(c\cdot c)^c=c^c\circ c^c=c\circ c,$$which contradicts Proposition~\ref{p:discale-fraction}(2). 
\smallskip

3. Proposition~\ref{p:diagonal11}(5) implies that $\{0,1\}\subseteq \mathcal Z\defeq\{x\in R:\forall a\in R\;\;x+a=a+x\}$. On the other hand, by Proposition~\ref{p:discale-fraction}(1), for every $c\in R\setminus\{0,1\}$, we have $c^2+c=c+c^{\frac12}\ne c+c^2$ and hence $c\notin \mathcal Z$. Therefore, $\mathcal Z=\{0,1\}$.
\smallskip

4. We have to prove that the right nucleus $\mathcal N_r=\{x\in R:\forall a,b\in R\;\;a+(b+x)=(a+b)+x\}$ of the loop $(R,+)$ coincides with the doubleton $\{0,1\}$. Proposition~\ref{p:diagonal11}(5) implies that $\{0,1\}\subseteq\mathcal N_r$. Assuming that $\{0,1\}\ne\mathcal N_r$, we can find an element $c\in \mathcal N_r\setminus\{0,1\}$. We claim that the function $T_c:\Pi\to\Pi$, $T_c:(x,y)\mapsto (x,y+c)$, is a translation of the plane $\Pi$. Given any line $L\subseteq\Pi$ we should prove that $T_c[L]$ is a line, parallel to the line $L$. If $L=\{a\}\times R$ for some $a\in R$, then $T_c[L]=L$ and we are done. So, assume  $L=L_{a,b}\defeq\{(x,x{\cdot}a+b):x\in R\}$ for some $a,b\in R$. In this case $T_c[L]=\{T(x,x{\cdot}a+b):x\in R\}=\{(x,(x{\cdot}a+b)+c):x\in R\}=\{(x,x{\cdot}a+(b+c):x\in R\}=L_{a,b+c}$. 

Let $R'$ be the smallest subternar of the ternar $R$, containing the set $\{0,1,c\}$. By Proposition~\ref{p:Yaqub1967}(1), for any elements $a,b\in R\setminus\{0,1\}$ with $a+b=c$, the map $$\varphi\defeq T_c^{-1}\psi_{a,b}:R^2\to R^2,\quad \varphi:(x,y)\mapsto (a+(b+x))-c,(a+(b+y))-c)$$ is an automorphism of the plane $\Pi=R\times R$. Observe that for every $x\in\{0,1,c\}$ we have $(a+(b+x))-c=x$. Then the set $\Fix(\varphi)$ contains the set $R'\times R'$ and hence $a+(b+x)-c=x$ for all $x\in R'$. 

Since $c^4\ne c^9$, there exists $n\in\{2,3\}$ such that $1\mp c\ne c^{n^2}$. Put $b\defeq c^{n^2}\in R'$ and find a unique $a\in R'$ such that $a+b=c$. It follows from $c^{n^2}\notin\{c,1\mp c\}$ that $a,b\notin\{0,1\}$. Applying Propositions~\ref{p:discale=>divisible} and \ref{p:di-conjugation}(2,6), we conclude that
$$c^{c^n}=((c^n)^{\frac1n})^{c^n}=((c^n)^{\circ n})^{c^n}=((c^n)^{c^n})^n=(c^n)^n=c^{n^2}.$$

Then for $x\defeq c^n\in R'$, we can apply Proposition~\ref{p:di-conjugation}(5) and obtain the equality 
$$a+(c^{n^2}+c^n)=a+(b+x)=x+c=c^n+c=c^{c^n}+c^n=c^{n^2}+c^n,$$which implies $a=0$ and contradicts the choice of $b=c^{n^2}\ne c$.
\end{proof}

A ternar $R$ is defined to be {\em diassociative-plus} if its plus loop $(R,+)$ is diassociative.

\begin{theorem}\label{t:hdiscale<=>corps} For a diagonal-scale horizontal-scale ternar $R$ the following conditions are equivalent:
\begin{enumerate}
\item $R$ is a corps;
\item $R$ is commutative-plus;
\item $R$ is associative-plus;
\item $R$ is diassociative-plus. 
\end{enumerate}
\end{theorem}

\begin{proof} The implications $(1)\Ra(2)$ and $(1)\Ra(3)\Ra(4)$ are trivial.
The implication $(2)\Ra(1)$ follows from Proposition~\ref{p:Yaqub1967}(3).\smallskip

$(4)\Ra(1)$ To derive a contradiction, assume that $R$ is diassociative-plus, but not a corps. Then $R$ is a $(\discale)$-ternar. Since the plus loop $(R,+)$ is diassociative, it is inversive. Then for every element $c\in R$ there exists a unique element $\bar c\in R$ such that $\bar c+(c+x)=x=(x+c)+\bar c$ for all $x\in R$. Proposition~\ref{p:inversive=>anti} ensures that $\overline{x+y}=\overline{y}+\overline x$ for all $x,y\in R$. Assuming that $x=\bar x$ for all $x\in R$, we conclude that $x+y=\overline{x+y}=\overline y+\overline x=y+x$ for all $x,y\in R$ and hence $R$ is commutative-plus, which contradicts Proposition~\ref{p:Yaqub1967}(3). This contradiction shows that $c\ne \bar c$ for some $c\in R$, which implies $c+c\ne 0$. Since $R$ is diassociative-plus, the elements $c,c^2$ and contained in some subgroup $G$ of the plus loop $(R,+)$. Proposition~\ref{p:di-conjugation}(5) and the associativity of the group $(G,+)$ imply 
$$x^{c+c}=((c+c)+x)-(c+c)=(c+((c+x)-c))-c=(x^c)^c$$ for all $x\in G$. 
By Proposition~\ref{p:di-conjugation}(2), $c=c^c$ and hence $c^{c+c}=(c^c)^c=c^c=c$. By Proposition~\ref{p:di-conjugation}(5,6), 
$(c\cdot c)^{c+c}=c^{c+c}\circ c^{c+c}=c\circ c$.
On the other hand,
$$(c\cdot c)^{c+c}=((c\cdot c)^c)^c=(c^c\circ c^c)^c=(c\circ c)^c=c^c\cdot c^c=c\cdot c.$$
Therefore, $c\circ c=c\cdot c$, which contradicts Proposition~\ref{p:discale-fraction}(2).
\end{proof}

\begin{remark} The equivalence $(1)\Leftrightarrow(2\wedge 3)$ in Theorem~\ref{t:hdiscale<=>corps} was proved by Pickert \cite{Pickert1959} and $(1)\Leftrightarrow(3)$ by  Spencer \cite{Spencer1960} and J\'onsson \cite{Jonsson1963}.
\end{remark}
   
\begin{problem} Is every inversive-plus diagonal-scale horizontal-scale ternar a corps?
\end{problem}

\section{The translation-shear-central-scale classification of affine planes} 
  
\begin{definition} For an affine plane $\Pi$, the set $$\INXT_\Pi\defeq \TT_\Pi\cup\NT_\Pi\cup\CT_\Pi\cup\SC_\Pi$$ is called the \defterm{translation-shear-central-scale trace} of the affine plane $\Pi$. 

The \index{translation-shear-central-scale rank}\index{rank!translation-shear-central-scale}\defterm{translation-shear-central-scale rank} $\|\INXT_\Pi\|$ of the affine plane $\Pi$ is the quintuple $inxdr$ of the numbers $i\defeq \|\TT_\Pi\|$, $n\defeq\|\NT_\Pi\|$, $x\defeq\|\CT_\Pi\|$, $d\defeq\|\dom[\SC_\Pi]\|$, $r\defeq\|\rng[\SC_\Pi]\|$.
\end{definition}

\begin{theorem}\label{t:INXT24} For every affine plane $\Pi$, exactly one of the following $25$ cases holds: 
\begin{enumerate}
\item $\|\INXT_\Pi\|=00000$ and $\INXT_\Pi=\varnothing$;
\item  $\|\INXT_\Pi\|=00011$ and $\INXT_\Pi=\{(\delta,L)\}$ for some direction $\delta\in\partial\Pi$ and line $L\in\mathcal L_\Pi\setminus \delta$;
\item  $\|\INXT_\Pi\|=00022$ and $\INXT_\Pi=\{(L_\parallel,\Lambda),(\Lambda_\parallel,L)\}$ for some concurrent lines $L,\Lambda\in\mathcal L_\Pi$;
\item  $\|\INXT_\Pi\|=00022$, $|\Pi|_2=9$, and $\INXT_\Pi=\{(L_\parallel,\phi(L)):L\in\mathcal L_p\}$ for some point $p\in\Pi$ and involution $\phi:\mathcal L_p\to\mathcal L_p$;
\item  $\|\INXT_\Pi\|=00022$, and $\INXT_\Pi$ is an injective function with $\dom[\INXT_\Pi]=\partial\Pi\setminus\{\delta\}$ and $\rng[\INXT_\Pi]=\delta$ for some direction $\delta\in\partial\Pi$;
\item $\|\INXT_\Pi\|=00100$ and $\INXT_\Pi=\{p\}$ for some point $p\in\Pi$;
\item $\|\INXT_\Pi\|=00111$  and $\INXT_\Pi=\{p\}\cup\{(\delta,L)\}$ for some direction $\delta\in\partial\Pi$, line $L\in\mathcal L_\Pi\setminus\{\delta\}$ and point $p\in L$;
\item  $\|\INXT_\Pi\|=00122$ and $\INXT_\Pi=(L\cap\Lambda)\cup\{(L_\parallel,\Lambda),(\Lambda_\parallel,L)\}$ for some concurrent lines $L,\Lambda\in\mathcal L_\Pi$;
\item  $\|\INXT_\Pi\|=01000$ and $\INXT_\Pi=\{L\}$ for some line $L\in\mathcal L_\Pi$;
\item  $\|\INXT_\Pi\|=01100$ and $\INXT_\Pi=\{L\}\cup\{p\}$ for some line $L\in\mathcal L_\Pi$ and point $p\in L$;
\item $\|\INXT_\Pi\|=01112$ and $\INXT_\Pi=\{V\}\cup\{p\}\cup (\{v\}\times (\mathcal L_p\setminus v))$ for some direction $v\in\partial \Pi$, line $V\in v$ and point $p\in V$;
\item  $\|\INXT_\Pi\|=02000$ and $\INXT_\Pi=\delta$ for some direction $\delta\in\partial\Pi$;
\item $\|\INXT_\Pi\|=02000$  and $\INXT_\Pi=\mathcal L_p$  for some point $p\in\Pi$;
\item $\|\INXT_\Pi\|=02100$ and $\INXT_\Pi=\mathcal L_p\cup\{p\}$ for some point $p\in\Pi$;
\item  $\|\INXT_\Pi\|=10000$ and $\INXT_\Pi=\{\delta\}$ for some direction $\delta\in \partial\Pi$;
\item  $\|\INXT_\Pi\|=10011$ and $\INXT_\Pi=\{v\}\cup\{(h,V)\}$ for some distinct directions $h,v\in\partial\Pi$ and line $V\in v$;
\item  $\|\INXT_\Pi\|=12000$ and $\INXT_\Pi=\{\delta\}\cup\delta$ for some direction $\delta\in\partial\Pi$;
\item  $\|\INXT_\Pi\|=12200$ and $\INXT_\Pi=\{\delta\}\cup\delta\cup L$ for some direction $\delta\in\partial\Pi$ and line $L\in\delta$;
\item $\|\INXT_\Pi\|=12221$ and $\INXT_\Pi=\{v\}\cup v\cup V\cup\{(\delta,V):\delta\in\partial\Pi\setminus\{V_\parallel\}\}$ for some direction $v\in\partial\Pi$ and line $V\in v$;
\item  $\|\INXT_\Pi\|=20000$ and $\INXT_\Pi=\partial\Pi$;
\item  $\|\INXT_\Pi\|=20023$ and $\INXT_\Pi=\partial\Pi\cup(\{v\}\times h)\cup(\{h\}\times v)$ for two distinct directions $h,v\in\partial\Pi$.
\item  $\|\INXT_\Pi\|=20023$, $|\Pi|_2=9$ and $\INXT_\Pi=\partial\Pi\cup\{(\delta,L):\delta\in\partial\Pi\;\wedge\;L\in \phi(\delta)\}$ for some  involution $\phi:\partial\Pi\to\partial\Pi$;
\item  $\|\INXT_\Pi\|=22000$ and $\INXT_\Pi=\partial\Pi\cup\delta$ for some direction $\delta\in\partial\Pi$;
\item  $\|\INXT_\Pi\|=23000$  and $\INXT_\Pi=\partial\Pi\cup\mathcal L_\Pi$;
\item  $\|\INXT_\Pi\|=23323$ and $\INXT_\Pi=\partial\Pi\cup\mathcal L_\Pi\cup\Pi\cup\{(\delta,L)\in\partial\Pi\times\mathcal L_\Pi:L\notin\delta\}$.
\end{enumerate}
\end{theorem}

\begin{proof} This theorem can be easily deduced from the Lenz--Barlotti classification~\ref{t:Lenz-Barlotti} by analyzing possible positions of the horizon line on a projective plane. 
\end{proof}

\begin{exercise} Deduce Theorem~\ref{t:INXT24} from Theorem~\ref{t:Lenz-Barlotti}.
\smallskip

{\em Hint:} Prove the items of Theorem~\ref{t:INXT24} in order of the items of Theorem~\ref{t:Lenz-Barlotti}.
\end{exercise}

\begin{remark} In the following table we present the translation-shear-central-scale ranks of all (seven) affine planes of order $9$ (calculated by Ivan Hetman).
$$
\begin{array}{l|c|c|c|c|c|c|c}
\mbox{The affine plane:}&\mbox{\tt Desarg}&\mbox{\tt Thales}&\mbox{\tt Hall}&\mbox{\tt dhall}&\mbox{\tt hall}&\mbox{\tt Hughes}&\mbox{\tt hughes}\\
\hline
\mbox{The rank:}&23323&20023&12221&01112&00022&00000&00000\\
\end{array}
$$
\end{remark}

\chapter{Based projective planes and their ternars}

In this chapter we consider projective bases in projective spaces, the corresponding ternars of projective planes. 

\section{Duality in projective planes}

\begin{definition} A \index{projective plane}\defterm{projective plane} is a $3$-long projective liner $X$ of rank $\|X\|=3$.
\end{definition}

Observe that the points and lines of any projective plane satisfy the following symmetric axioms:
\begin{itemize}
\item any two distinct points belong to a unique line;
\item any two distinct lines have a unique common point;
\item no line contains all points;
\item no point belongs to all lines.
\end{itemize}

Observe that those axioms remain true if we replace the words ``points'' by ``lines'' and vice versa. This fact implies that all statements about projective planes remain true if we exchange the words ``points'' in their formulations by ``lines'' and vice versa. This fact is known in Projective Geometry as the \index{Duality Principle}\defterm{Duality Principle}. 

For every projective plane $\Pi$, the \index{dual projective plane}\defterm{dual projective plane} has the family of lines $\mathcal L_\Pi$ of $\Pi$ as the set of points. Every point $p\in\Pi$ determines the ``line'' $\mathcal L_p\defeq\{L\in\mathcal L:p\in L\}$ in the dual projective plane. In the plane $\Pi$ this ``line'' $\mathcal L_p$ is called the \index{pencil of lines}\defterm{pencil of lines} determined by the point $p$.

A projective plane $\Pi$ is called \index{self-dual projective plane}\index{projective plane!self-dual}\defterm{self-dual} if it is isomorphic to its dual plane.

\begin{theorem}\label{t:Pappian-self-dual} Every Pappian projective plane is self-dual.
\end{theorem}    

\begin{proof} Let $X$ be a Pappian projective plane and let $\IR_X$ be its field of scalars. By Theorems~\ref{t:Dp=PRX} and \ref{t:Papp<=>Des+RX}, $X$ is isomorphic to the projective space $\mathbb P\IR_X^3$ of the $3$-dimensional vector space $\IR_X^3$ over the field $\IR_X$.
So, $X$ be can be identified with the projective space $\mathbb P\IR_X^3$. 
Let $\mathcal L$ be the family of lines in $X$. Consider the map $F:X\to\mathcal L$, $F:V\mapsto V^\perp$ assigning to every point $V\in X$ (which is a $1$-dimensional vector subspace of $\IR_X^3$) the line $$V^\perp\defeq\{L\in X:\forall (x,y,z)\in L\;\forall (a,b,c)\in V\;a{\cdot}x+b{\cdot}y+c{\cdot}z=0\}.$$It can be shown that $F$ is an isomorphism of the projective plane $X$ and its dual projective plane.
\end{proof}

\begin{exercise} Find a Desarguesian projective plane which is self-dual but not Pappian.
\smallskip

{\em Hint:} Look at the projective plane over the corps of quaternions.
\end{exercise}

\begin{exercise} Find an example of a projective plane of order 9, which is not self-dual.
\smallskip

{\em Hint:} Look at the right-Hall and left-Hall projective planes.
\end{exercise}

\begin{Exercise} Find an example of a Desarguesian projective plane which is not self-dual.
\smallskip

{\em Hint:} See {\tt https://math.stackexchange.com/a/45086}. 
\end{Exercise}

\section{Projective bases in projection planes}

\begin{definition} A quadruple $uowe$ of four distinct points of a projective plane $\Pi$ is called a \index{projective base}\defterm{projective base} in $\Pi$ if $uow$ is a triangle in $\Pi$ and $e\notin\Aline uo\cup\Aline ow\cup\Aline uw$. A projective plane $\Pi$ endowed with a projective base $uowe$ is called a \index{based projective plane}\index{projective plane!based}\defterm{based projective plane}. 
\end{definition}

The points $u,o,w,e$ of a projective base $uowe$ are called respectively  the \index{projective base!unit of}\defterm{unit}, the \index{projective base!origin of}\defterm{origin}, the \index{projective base!biunit of}\defterm{biunit}, and the \index{projective base!diunit of}\defterm{diunit} of the projective base $uowe$. The unique points $h\in\Aline ou\cap\Aline we$ and $v\in\Aline ow\cap\Aline ue$ are called the \index{horizontal infinity point}\index{projective base!horizontal infinity point}\defterm{horizontal infinity point} and the \index{vertical infinity point}\index{projective base!vertical infinity point}\defterm{vertical infinity point}, the line $\Aline hv$ is called the \index{infinite line}\index{projective base!infinite line of}\defterm{infinite line} of the projective base, the affine plane $\Pi\setminus\Aline hv$ is called the \index{projective base!affine plane of}\defterm{affine plane} of the projective base, and the punctured line $\Delta\defeq\Aline oe\setminus\Aline hv$ is called the \index{diagonal}\index{projective base!diagonal of}\defterm{diagonal} of the projective base $uowe$.

\begin{picture}(100,120)(-90,-15)

{\linethickness{0.8pt}
\put(0,0){\color{teal}\line(1,0){90}}
\put(0,0){\color{cyan}\line(0,1){90}}
\put(0,0){\color{red}\line(1,1){60}}
\put(0,90){\color{violet}\line(1,-1){90}}
}

\put(0,90){\line(1,-2){45}}
\put(90,0){\line(-2,1){90}}

\put(90,0){\circle*{3}}
\put(95,-3){$h$}
\put(0,90){\circle*{3}}
\put(-3,95){$v$}
\put(63,63){\color{red}$\Delta$}

\put(45,45){\color{red}\circle*{3}}
\put(45,45){\color{white}\circle*{2}}
\put(0,0){\circle*{3}}
\put(-6,-8){$o$}
\put(45,0){\circle*{3}}
\put(42,-8){$u$}
\put(0,45){\circle*{3}}
\put(-10,43){$w$}
\put(30,30){\circle*{3}}
\put(29,34){$e$}

\put(150,87){$o$, the origin}
\put(150,77){$u$, the unit}
\put(150,67){$w$, the biunit}
\put(150,57){$e$, the diunit}
\put(150,47){$h$, the horizontal infinity point}
\put(150,37){$v$, the vertical infinity point}
\put(150,27){$\Delta\defeq\Aline oe\setminus\Aline hv$, the diagonal}
\put(150,17){$\Aline ou=\Aline oh$, the horizontal axis}
\put(150,7){$\Aline ow=\Aline ov$, the vertical axis}
\put(150,-3){$\Aline hv$, the infinity line}
\end{picture}

If the projective plane $\Pi$ has $|\Pi|_2=3$, then $|\Delta|=\{o,e\}$ and $\Delta$ carries a unique structure of a linear ternar (whose biloop is the two-element field consisting of zero $o$ and identity $e$). If $|\Pi|_2\ge 4$, then the affine subliner $\Pi\setminus\Aline hv$ of $\Pi$ is $3$-long and has rank $\|\Pi\setminus\Aline hv\|=\|\Pi\|$, by Corollary~\ref{c:procompletion-rank}. In this case, $\Pi\setminus\Aline hv$ is an affine plane.

 The choice of the points $h$ and $v$ ensures that $uow$ is an affine base in the affine plane $\Pi\setminus\Aline hv$ and the diunit $e$ of the projective base $uowe$ coincides with the diunit of the affine base $uow$ in the affine plane $\Pi\setminus\Aline hv$. The points $h$ and $v$ can be canonically indentified with the horizontal and vertical directions $\boldsymbol h$ and $\boldsymbol v$ on the affine plane $\Pi\setminus\Aline hv$. The ternar $(\Delta,T_{uow})$ of the based affine plane $(\Pi\setminus\Aline hv,uow)$ is called \index{based projective plane!ternar of}\defterm{the ternar} of the based projective plane $(\Pi,uowe)$.

A ternar $R$ is called \index{ternar of a based projective plane}\defterm{a ternar of a projective plane} $\Pi$ if there exists a projective base $uowe$ such that $R$ is isomorphic to the ternar of the based projective plane $(\Pi,uowe)$.

For every ternar $R$, its \index{coordinate projective plane}\defterm{coordinate projective plane} $\overline{R^2}$ is defined as the spread completion of the coordinate plane $R^2$ of the ternar $R$. The projective plane $\overline{R^2}$ carries the canonical projective base $(10,00,01,11)$.

\begin{definition} Two based projective planes $(\hat\Pi,\hat u\hat o\hat w\hat e)$ and $(\check \Pi,\check u\check o\check w\check e)$ are defined to be \index{isomorphic based projective planes}\index{based projective planes!isomorphic}\defterm{isomorphic} if there exists an isomorphism $I:\hat \Pi\to\check \Pi$ of the projective planes $\hat \Pi$, $\check \Pi$ such that $I\hat u\hat o\hat w\hat e=\check u\check o\check w\check e$.
\end{definition}

For every based projective plane $(\Pi,uowe)$, consider the coordinate chart $C:\Pi\setminus\Aline hv\to \Delta^2$, which is an isomorphism of the based affine planes $(\Pi\setminus\Aline hv,uow)$ and $(\Delta^2,(eo,oo,oe))$. By Theorem~\ref{t:extend-isomorphism-to-completions}, the isomorphism $C$ extends to an isomorphism $\bar C:\Pi\to\overline{\Delta^2}$ between the projective planes $\Pi$ and $\overline{\Delta^2}$. Since $\bar C(e)=C(e)=ee$, the bijection $\bar C$ is an isomorphism of the based projective planes $(\Pi,uowe)$ and $(\overline{\Delta^2},(eo,oo,oe,ee))$. 

\begin{theorem}\label{t:p-isomorphic<=>ternar-isomorphic} Two based projective planes are isomorphic if and only if their ternars are isomorphic.
\end{theorem}

\begin{proof} Let $(\hat \Pi,\hat u\hat o\hat w\hat e)$ and $(\check \Pi,\check u\check o\check w\check e)$ be two based projective planes. Let $\hat h\in\Aline{\hat o}{\hat u}\cap \Aline{\hat w}{\hat e}$ and $\hat v\in\Aline{\hat o}{\hat w}\cap \Aline{\hat u}{\hat e}$ be the horizontal and vertical infinity points in the based projective plane $(\hat \Pi,\hat u\hat o\hat w\hat e)$, and $\check h\in\Aline{\check o}{\check u}\cap \Aline{\check w}{\check e}$ and $\check v\in\Aline{\check o}{\check w}\cap \Aline{\check u}{\check e}$ be the horizontal and vertical infinity points in the based projective plane $(\check \Pi,\check u\check o\check w\check e)$.
Let $\hat \Delta\defeq\Aline{\hat o}{\hat e}\setminus\Aline {\hat h}{\hat v}$ and $\check \Delta\defeq\Aline{\check o}{\check e}\setminus\Aline {\check h}{\check v}$ be the diagonals of the based affine planes 
$(\hat\Pi\setminus\Aline{\hat h}{\hat v},\hat u\hat o\hat w)$ and $(\check \Pi\setminus\Aline{\check h}{\check v},\check u\check o\check w)$, respectively.

If the based projective planes $(\hat \Pi,\hat u\hat o\hat w\hat e)$ and $(\check \Pi,\check u\check o\check w\check e)$ are isomorphic, then there exists an isomorphism $I:\hat\Pi\to\check \Pi$ such that $I\hat u\hat o\hat w\hat e=\check u\check o\check w\check e$. If one of the projective planes $\hat\Pi$ or $\check\Pi$ is Steiner, then the isomorphiness of $\hat\Pi$ and $\check\Pi$ implies that both planes are Steiner. In this case, the their ternars are two-elements fields, which are isomorphic. So, assume that both projective planes $\hat\Pi$ and $\check\Pi$ are not Steiner. The choice of the points $\hat h,\hat v,\check h,\check v$ ensures that $I(\hat h)\in I[\Aline {\hat o}{\hat u}\cap\Aline {\hat w}{\hat e}]=\Aline{\check o}{\check u}\cap\Aline{\check w}{\check e}=\{\check h\}$ and $I(\hat v)=\check v$. Then $I[\Aline{\hat h}{\hat v}]=\Aline{\check h}{\check v}$ and hence $I{\restriction}_{\hat \Pi\setminus\Aline{\hat h}{\hat v}}$ is an isomorphism of the based affine planes $(\hat\Pi\setminus\Aline{\hat h}{\hat v},\hat u\hat o\hat w)$ and $(\check \Pi\setminus\Aline{\check h}{\check v},\check u\check o\check w)$.  By Corollary~\ref{c:basedafplane-iso<=>}, the isomorphness of the based affine planes 
 $(\hat\Pi\setminus\Aline{\hat h}{\hat v},\hat u\hat o\hat w)$ and $(\check \Pi\setminus\Aline{\check h}{\check v},\check u\check o\check w)$ implies the isomorphness of their ternars $(\hat\Delta,T_{\hat u\hat o\hat w})$ and $(\check \Delta,T_{\check u\check o\check w})$, which are also the ternars of the based projective planes  $(\hat\Pi,\hat u\hat o\hat w\hat)$ and $(\check \Pi,\check u\check o\check w\check e)$. This completes the proof of the ``only if'' part of the theorem.
\smallskip

To prove the ``if'' part, assume that the ternars  $(\hat\Delta,T_{\hat u\hat o\hat w})$ and $(\check \Delta,T_{\check u\check o\check w})$ of the based projective planes  $(\hat\Pi,\hat u\hat o\hat w\hat e)$ and $(\check \Pi,\check u\check o\check w\check e)$ are isomorphic. By Corollary~\ref{c:basedafplane-iso<=>}, there exists an isomorphism $I:\hat \Pi\setminus\Aline{\hat h}{\hat v}\to\check \Pi\setminus\Aline{\check h}{\check v}$ of the based affine planes  $(\hat\Pi\setminus\Aline{\hat h}{\hat v},\hat u\hat o\hat w)$ and $(\check \Pi\setminus\Aline{\check h}{\check v},\check u\check o\check w)$. Since $\hat e$ is a unique point in $\hat \Pi$ such that $\Aline {\hat e}{\hat u}\parallel \Aline {\hat w}{\hat o}$ and $\Aline {\hat e}{\hat w}\parallel \Aline{\hat u}{\hat o}$ and the point $\check e$ has the analogous uniqueness property, $I(\hat e)=\check e$. Since the projective planes $\hat \Pi$ and $\check\Pi$ are projective completions of the affine planes $\hat\Pi\setminus\Aline{\hat h}{\hat v}$ and $\check \Pi\setminus\Aline{\check h}{\check v}$, respectively, we can apply Theorem~\ref{t:extend-isomorphism-to-completions}, and conclude that the ismorphism $I$ extends to an isomorphism $\bar I:\hat\Pi\to\check\Pi$ of the projective planes $\hat \Pi$ and $\check \Pi$. Since $\bar I\hat u\hat o\hat w\hat e=I\hat u\hat o\hat w\hat e=\check u\check o\check w\check e$, $\bar I$ is an isomorphism of the based projective planes  $(\hat\Pi,\hat u\hat o\hat w\hat e)$ and $(\check \Pi,\check u\check o\check w\check e)$.
\end{proof}

We recall that a two ternars $(R,T)$ and $ (R',T')$ are \defterm{isotopic} if there exist bijections $F,G,H:R\to R'$ such that $T'(F(x),G(y),H(z))=H(T(x,y,z))$ for all $x,y,z\in R$. 

\begin{theorem} The ternars of two based projective planes $(\check \Pi,\check u\check o\check w\check e)$ and $(\hat \Pi,\hat u\hat o\hat w\hat e)$ are isotopic if and only if there exists a liner isomorphism $\bar\Phi:\check \Pi\to\hat \Pi$ such $\bar\Phi(\check h)=\hat h$, $\bar\Phi(\check v)=\hat v$ and $\bar\Phi(\check o)\in \Aline{\hat o}{\hat v}$, where $\check h\in\Aline{\check o}{\check u}\cap\Aline{\check w}{\check e}$, $\hat h\in\Aline{\hat o}{\hat u}\cap\Aline{\hat w}{\hat e}$,   $\check v\in\Aline{\check o}{\check w}\cap\Aline{\check u}{\check e}$, $\hat v\in\Aline{\hat o}{\hat w}\cap\Aline{\hat u}{\hat e}$ are the horizontal and vertical infinity points of the projective bases $\check u\check o\check w\check e$ and $\hat u\hat o\hat w\hat e$.  
\end{theorem}

\begin{proof} By definition, the ternars $\check \Delta$ and $\hat \Delta$ of the based projective planes $(\check \Pi,\check u\check o\check w\check e)$ and $(\hat \Pi,\hat u\hat o\hat w\hat e)$ coincide with the ternars of the based affine planes $(\check \Pi\setminus\Aline{\check h}{\check v},\check u\check o\check w)$ and $(\hat \Pi\setminus\Aline{\hat h}{\hat v},\hat u\hat o\hat w)$. In its turn, those based affine planes are isomorphic to the coordinate planes of the corresponding ternars. More precisly, there exist besed affine plane isomorphisms $\check C:\check\Pi\setminus\Aline {\check h}{\check v}\to\check\Delta^2$ and  $\hat C:\hat\Pi\setminus\Aline {\hat h}{\hat v}\to\hat\Delta^2.$
\smallskip

If the ternars of the based projective planes $(\check \Pi,\check u\check o\check w\check e)$ and $(\hat \Pi,\hat u\hat o\hat w\hat e)$ are isotopic, then so are the ternars of the based affine planes $(\check \Pi\setminus\Aline{\check h}{\check v},\check u\check o\check w)$ and $(\hat \Pi\setminus\Aline{\hat h}{\hat v},\hat u\hat o\hat w)$. By Proposition~\ref{p:isotopic=>isomorphic}, there exist two bijections $F,H:\check \Delta\to\hat\Delta$ such that the function $\Psi:\check\Delta^2\to\hat\Delta^2$, $\Psi:(x,y)\mapsto (F(x),H(y))$, is an isomorphism of the coordinate planes of the ternars $\check \Delta$ and $\hat \Delta$. Then the function $\Phi\defeq \hat C^{-1}\circ\Psi\circ\check C:\check\Pi\setminus\Aline{\check h}{\check v}\to\hat\Pi\setminus\Aline{\hat h}{\hat v}$ is an isomorphism of the affine planes that uniquely extends to an isomorphism $\bar\Phi:\check\Pi\to\hat\Pi$ of the projective planes, by Theorem~\ref{t:extend-isomorphism-to-completions}. Taking into account that the isomorphisms $\Psi$, $\hat C$, $\check C$ map vertical lines to vertical lines and horizontal lines, we conclude that $\bar\Phi(\check h)=\hat h$ and $\bar\Phi(\check v)=\hat v$. Proposition~\ref{p:isotopic=>zero} ensures that $F(\check o)=\hat o$, which implies $\bar\Phi(\check o)=\Phi(\check o)\in \Aline{\hat o}{\hat v}$. 
\smallskip

Now assume conversely that there exists a liner isomorphism $\bar \Phi:\check \Pi\to\hat\Pi$ such that  $\bar\Phi(\check h)=\hat h$, $\bar\Phi(\check v)=\hat v$ and $\bar\Phi(\check o)\in \Aline{\hat o}{\hat v}$. Then its restriction $\Phi\defeq\bar\Phi{\restriction}_{\check\Pi\setminus\Aline{\check h}{\check v}}$ is an isomorphism of the affine planes $\check\Pi\setminus\Aline{\check h}{\check v}$ and $\hat\Pi\setminus\Aline{\hat h}{\hat v}$ that maps horizontal lines to horizontal lines and vertical lines to vertical lines. This implies that the isomorphism $\Psi\defeq \hat C\circ\Phi\circ\check C^{-1}:\check\Delta^2\to\hat\Delta^2$ of the coordinate planes of the ternars $\check\Delta$ and $\hat \Delta$ also maps horizontal lines to horizontal lines and vertical lines to vertical lines. Then $\Psi$ is of the form $\Psi:(x,y)\mapsto (F(x),H(y))$ for some bijections $F,G:\check\Delta\to\hat H$. The inclusion $\bar\Phi(\check o)\in\Aline{\check o}{\check v}$ implies $F(\check o)=\hat o$. By Proposition~\ref{p:isomorphic=>isotopic}, the ternars $\check\Delta$ and $\hat \Delta$ are isotopic. 
\end{proof} 

\begin{remark} By Remark~\ref{r:trings-1693}, there exist exactly 1693 non-isomorphic and 33 non-isotopic ternars of order 9.
The numbers of isomorpy and isotopy classes of ternars whose coordinate projective planes are isomorphic to one of four possible projective planes of order 9 are presented in the following table (composed from the table in Remark~\ref{r:trings-1693}).
$$
\begin{array}{l|c|c|c|c}

\mbox{Projective plane:}&\mbox{Desarg}&\mbox{right-Hall}&\mbox{left-Hall}&\mbox{Hughes}\\
\hline
\mbox{Ternars:}&1&169&169&1354\\
\mbox{up to isotopy:}&1&8&8&16\\
\end{array}
$$
\end{remark}

\section{A duality between plus and puls operations}

In this section we discuss the duality between the plus and puls operations in projective planes and their duals. 

\begin{theorem}\label{t:puls<=>plus} For every projective base $uowe$ in a projective plane $(\Pi,\mathcal L)$ there exists a projective base $UOWE$ in the dual projective plane $(\mathcal L,\Pi)$ whose plus loop $(\Delta^*,+)$ is isomorphic to the puls loop $(\Delta,\!\puls\!)$ of the based projective plane $(\Pi,uowe)$.
\end{theorem}

\begin{proof} Let $uowe$ be a projective base in the projective plane $(\Pi,\mathcal L)$. Consider the points $h\defeq\Aline ou\cap\Aline we$, $v\defeq \Aline ow\cap\Aline ue$, and $d\defeq\Aline oe\cap \Aline hv$. 

Recall that the puls operation on the diagonal $\Delta\defeq\Aline oe\setminus\{d\}$ is defined as follows. Given two points $x,y\in \Delta$, we find unique points $s\defeq \Aline yh\cap\Aline ov$, $\bar x\in \Aline xh\cap\Aline ev$, $t\in \Aline o{\bar x}\cap\Aline hv$, $\bar z\in \Aline st\cap\Aline ev$ and $z\in \Aline {\bar z}h\cap\Delta$. The point $z$ is equal to the sum $x\puls y$.

Recall also the definition of the plus operation on $\Delta$. Given two points $x,y\in \Delta$, we find unique points $s\defeq \Aline yh\cap\Aline ov$, $t\in \Aline sd\cap\Aline xv$ and $z\in \Aline th\cap \Aline od$. The point $z$ is equal to the sum $x+y$.

\begin{picture}(450,210)(0,-10)

\put(220,0){\line(1,0){180}}
\put(220,0){\line(0,1){180}}
\put(220,0){\line(1,1){90}}
\put(220,180){\line(1,-1){180}}
\put(220,30){\line(6,-1){180}}
\put(220,30){\line(3,2){90}}
\put(220,180){\line(1,-6){25.7}}
\put(220,180){\line(3,-8){49}}
\put(265,60){\line(9,-4){135}}

\put(220,0){\circle*{3}}
\put(214,-6){$o$}
\put(265,60){\circle*{3}}
\put(265,64){$t$}
\put(400,0){\circle*{3}}
\put(403,-3){$h$}
\put(220,180){\circle*{3}}
\put(218,184){$v$}
\put(310,90){\circle*{3}}
\put(312,92){$d$}
\put(220,30){\circle*{3}}
\put(212,28){$s$}
\put(245.7,25.7){\circle*{3}}
\put(245,19){$y$}
\put(269,49){\circle*{3}}
\put(267,40){$x$}
\put(275.3,55.3){\circle*{3}}
\put(280,54){$z=x{+}y$}

\put(0,0){\line(1,0){180}}
\put(0,0){\line(0,1){180}}
\put(0,0){\line(1,1){90}}
\put(0,180){\line(1,-1){180}}
\put(0,20){\line(9,-1){180}}
\put(0,0){\line(2,3){72}}
\put(0,180){\line(1,-6){30}}
\put(0,20){\line(9,11){72}}
\put(24,36){\line(39,-9){156}}
\put(180,0){\line(-57,17){158}}

\put(0,0){\circle*{3}}
\put(-6,-6){$o$}
\put(72,108){\circle*{3}}
\put(75,108){$t$}
\put(180,0){\circle*{3}}
\put(183,-3){$h$}
\put(0,180){\circle*{3}}
\put(-2,184){$v$}
\put(90,90){\circle*{3}}
\put(92,92){$d$}
\put(0,20){\circle*{3}}
\put(-8,18){$s$}
\put(25.7,25.7){\circle*{3}}
\put(19,24){$e$}
\put(24,36){\circle*{3}}
\put(16,34){$\bar x$}
\put(22.15,47.1){\circle*{3}}
\put(15,47){$\bar z$}
\put(18,18){\circle*{3}}
\put(16,10){$y$}
\put(33.6,33.6){\circle*{3}}
\put(31,26){$x$}
\put(41.3,41.3){\circle*{3}}
\put(46,41){$z=x\puls y$}
\end{picture}

Now we shall dualize the procedure of finding the sum $z=x+y$. First we define a suitable projective base in the dual projective plane $(\mathcal L,\Pi)$. Consider the lines $O\defeq \Aline oh$, $H\defeq \Aline ue$, $D\defeq \Aline ov$, $V\defeq \Aline vh$. The lines $O,D$ are points of the dual projective plane $(\mathcal L,\Pi)$. The line in $(\mathcal L,\Pi)$ passing through these points can be identified with the pencil of lines $\mathcal L_o\defeq\{L\in\mathcal L:o\in L\}$. Let $\Delta^*\defeq\mathcal L_o\setminus\{D\}$. Consider the lines $E\defeq\Aline oe$, $U\defeq \overline{(V\cap E)\cup (O\cap H)}$ and $W\defeq \overline{(H\cap E)\cup(O\cap V)}$. Then $UOWE$ is a projective base in the dual projective plane $(\mathcal L,\Pi)$, and $\Delta^*=\mathcal L_o\setminus\{D\}$ is the diagonal of the based projective plane $(\mathcal L,UOW\!E)$.  Given any elements $X,Y\in\Delta^*$, let us find their sum $Z\defeq X+Y$ in the plus loop $(\Delta^*,+)$. For this find unique lines $S\defeq \overline{(Y\cap H)\cup (O\cap V)}$, $T\defeq \overline{(S\cap D)\cup (X\cap V)}$ and $Z\defeq \overline{(T\cap H)\cup (O\cap D)}$. Then $Z$ is equal to the sum $X+Y$ in the plus loop $(\Delta^*,+)$ of the based projective plane $(\mathcal L,UOW\!E)$.

\begin{picture}(250,285)(-50,-30)

\put(240,0){\line(-1,0){255}}
\put(-25,-3){$O$}
\put(0,240){\line(0,-1){255}}
\put(-4,-25){$D$}
\put(0,0){\line(1,1){130}}
\put(0,240){\line(1,-1){250}}
\put(250,-20){$V$}
\put(0,240){\line(1,-6){42}}
\put(38,-24){$H$}
\put(120,120){\line(-4,-3){135}}
\put(-25,14){$T$}
\put(240,0){\line(-8,1){255}}
\put(0,0){\line(7,5){150}}
\put(0,0){\line(7,12){95}}
\put(0,0){\line(3,1){190}}
\put(195,62){$E$}
\put(240,0){\line(-6,1){205}}
\put(240,0){\line(-47,12){209}}

\put(37.8,12.6){\circle*{3}}
\put(33,5){$e$}
\put(65.4,21.8){\circle*{3}}
\put(62,14){$y$}
\put(79.8,26.6){\circle*{3}}
\put(78,30){$x$}
\put(104.1,34.7){\circle*{3}}
\put(114,32.8){$z=x\puls y$}
\put(0,0){\circle*{3}}
\put(-7,-7){$o$}
\put(31.15,53.1){\circle*{3}}
\put(24,52){$\bar z$}
\put(240,0){\circle*{3}}
\put(243,-1){$h$}
\put(0,240){\circle*{3}}
\put(-2,244){$v$}
\put(120,120){\circle*{3}}
\put(0,30){\circle*{3}}
\put(-25,30){$S$}
\put(-6,32){$s$}
\put(35.7,25.5){\circle*{3}}
\put(30.5,15.5){$\bar y$}
\put(155,108){$Y$}
\put(133,133){$X$}
\put(98,168){$Z=X+Y$}
\put(34.3,34.3){\circle*{3}}
\put(26,33){$\bar x$}

\end{picture}

We claim that the plus loop $(\Delta^*,+)$ is isomorphic to the puls loop $(\Delta,\!\puls\!)$ of the projective plane $(\Pi,uowe)$. Given any line $X\in \Delta^*$, find unique points $\bar x\in X\cap H$ and $x\in \Delta\cap \Aline {\bar x}h$. We claim that the map $\varphi:\Delta^*\to\Delta$, $\varphi:X\mapsto x$, is an isomorphism of the loops $(\Delta^*,+)$ and $(\Delta,\!\puls\!)$ such that $\varphi(E)=e$. It is clear that the map $\varphi:\Delta^*\to\Delta$ is bijective, and $\varphi(E)=\bar e=e$. To prove that it preserves the operations, take any elements $X,Y\in \Delta^*$, find their sum $Z\defeq X+Y$ in the plus loop $(\Delta^*,+)$ and consider the points $\bar x\in X\cap H$, $\bar y\in Y\cap H$, $\bar z\in Z\cap H$, $x\in \Delta\cap\Aline {\bar x}h$, $y\in\Delta\cap\Aline {\bar y}h$, and $z=\Delta\cap\Aline {\bar z}h$. 
Consider the unique point $s\in S\cap D=\Aline {\bar y}h\cap\Aline ov$. Since $\Aline s{\bar z}\cap \Aline o{\bar x}\subseteq V=\Aline vh$, the lines $\Aline s{\bar z}\setminus\Aline hv$ and $\Aline o{\bar x}\setminus\Aline hv$ are parallel in the affine plane $\Pi\setminus\Aline hv$. Also the lines $\Aline s{\bar y}\setminus\{h\}$ and $\Aline oh\setminus\{h\}$ are parallel in the affine plane $\Pi\setminus\Aline hv$. The definition of the puls operation in the based affine plane $(\Pi\setminus\Aline hv,uow)$ ensures that $z= x\puls  y$ and hence
$$\varphi(X+Y)=\varphi(Z)=z= x\puls y=\varphi(X)\puls\varphi(Y).$$
Therefore, the map $\varphi:\Delta^*\to\Delta$ is an  isomorphism of the loops $(\Delta^*,+)$ and $(\Delta,\!\puls\!)$ such that $\varphi(E)=e$. Its inverse $\psi\defeq\varphi^{-1}:\Delta\to\Delta^*$  is an  isomorphism of the loops $(\Delta,\!\puls\!)$ and $(\Delta^*,+)$ such that $\varphi(E)=e$.
\end{proof}
\begin{definition} Let $\mathcal P$ be a property of loops. A projective plane $\Pi$ is defined to be 
\begin{itemize}
\item \defterm{$\mathcal P$-plus} if for every ternar $R$ of $\Pi$, the plus loop $(R,+)$ has property $\mathcal P$;
\item \defterm{$\mathcal P$-puls} if for every ternar $R$ of $\Pi$, the puls loop $(R,\!\puls\!)$ has property $\mathcal P$.
\end{itemize}
\end{definition}

\begin{theorem}\label{t:plus-puls-duality} Let $\mathcal P$ be a property of loops. A projective plane is $\mathcal P$-puls if and only if its dual projective plane is $\mathcal P$-plus.  
\end{theorem}

\begin{proof} Let $(\Pi,\mathcal L)$ be a projective plane and $(\mathcal L,\Pi)$ be its dual projective plane.

First assume that the dual projective plane $(\mathcal L,\Pi)$ is $\mathcal P$-plus. To prove that the projective plane $(\Pi,\mathcal L)$ is $\mathcal P$-puls, take any ternar $R$ of $(\Pi,\mathcal L)$, and find a projective base $uowe$ in $\Pi$ whose ternar $\Delta$ is isomorphic to the ternar $R$. Consider the points $h\in \Aline ou\cap\Aline we$ and $v\in \Aline ow\cap\Aline ue$, and the lines $O\defeq \Aline oh$, $E\defeq \Aline oe$, $H\defeq \Aline eh$, $V\defeq\Aline vh$,  $U\defeq \overline{(O\cap H)\cup(V\cap E)}$ and $W\defeq\overline{(O\cap V)\cup (E\cap H)}$. Then $UOWE$ is a projective base in the dual projective plane $(\mathcal L,\Pi)$ whose plus loop $(\Delta^*,+)$ is isomorphic to the puls loop $(\Delta,\!\puls\!)$ of the ternar of the based projective plane $(\Pi,\mathcal L)$, according to (the proof of) Theorem~\ref{t:puls<=>plus}. Since the projective plane $(\mathcal L,\Pi)$ is $\mathcal P$-plus, the plus loop $(\Delta^*,+)$ has the property $\mathcal P$ and so does its isomorphic copy $(\Delta,\!\puls\!)$, which is isomorphic to the puls loop $(R,\!\puls\!)$ of the ternar $R$, witnessing that the pojective plane $(\Pi,\mathcal L)$ is $\mathcal P$-puls.
\smallskip
 
Now assume that a projective plane $(\Pi,\mathcal L)$ is $\mathcal P$-puls. To prove that the dual projective plane $(\mathcal L,\Pi)$ is $\mathcal P$-plus, take any ternar $R$ of $(\mathcal L,\Pi)$ and find a projective base $UOWE\in\mathcal L^4$ in $(\mathcal L,\Pi)$ whose ternar $\Delta^*$ is isomorphic to the ternar $R$. Consider the lines $H\defeq\overline{(O\cap U)\cup(W\cap E)}$, $V\defeq\overline{(O\cap W)\cup(U\cap E)}$ in $\Pi$ and the unique points $o\in O\cap E$, $h\in O\cap V$, $v\in H\cap V$, $e\in H\cap E$, $u\in \Aline oh\cap\Aline ev=O\cap H$ and $w\in \Aline ov\cap\Aline eh$. Then $uowe$ is a projective base for the projective plane $(\Pi,\mathcal L)$. The proof of Theorem~\ref{t:puls<=>plus} ensures that the puls loop $(\Delta,\!\puls\!)$ of the based projective plane $(\Pi,uowe)$ is isomorphic to the plus loop $(\Delta^*,+)$ if the based projective plane $(\mathcal L,UOW\!E)$. Since the projective plane $\Pi$ is $\mathcal P$-puls, the puls loop $(\Delta,\!\puls\!)$ has property $\mathcal P$ and so does its isomorphic copy $(\Delta^*,+)$, which is isomorphic to the plus loop $(R,+)$ of the ternar $R$. Therefore, the dual projective plane $(\mathcal L,\Pi)$ is $\mathcal P$-plus.  
\end{proof}

\begin{remark} In Theorem~\ref{t:puls<=>plus} we have proved that for every projective base $uowe$ in a projective plane $(\Pi,\mathcal L)$, there exist a projective base $UOW\!E$ in the dual projective plane $(\mathcal L,\Pi)$ and a bijection $\psi:\Delta\to\Delta^*$ between the ternars of the based projective planes  $(\Pi,uowe)$ and $(\mathcal L,UOW\!E)$ such that  $\psi(e)=E$ and $$\psi(e_\times x_+y)=\psi(x\puls y)=\psi(x)+\psi(y)=\psi(x)_\times E_+\psi(y)=\psi(x)_\times \psi(e)_+\psi(y)\mbox{\quad for all $x,y\in \Delta$}.$$ Looking at this identity, one may suggest that the bijection $\psi$ has a stronger property:
 $$\psi(a_\times x_+y)=\psi(x)_\times \psi(a)_+\psi(y)\mbox{\quad for all $a,x,y\in \Delta$}.$$
 However, this is not true: by computer calculations Ivan Hetman has found that every non-Desarguesian projective plane $\Pi$ of order $9$ has a ternar $(R,T)$ such that the ternary function $T^*:R^3\to R$, $T^*:(x,y,z)\mapsto T(y,x,z)$, does not satisfy the axiom $\mathsf{(T3)}$ of a ternar and hence cannot be identified with a ternar of the dual projective plane $\Pi^*$ of $\Pi$.
\end{remark}

\begin{exercise} Show for every ternar $(R,T)$ of a  Thalesian affine plane, the function $T^*:R^3\to R$, $T^*:(x,y,z)\mapsto T(y,x,z)$, satisfies the axioms $\mathsf{(T1)}$--$\mathsf{(T4)}$ of a ternar.
\end{exercise}

\section{Central and axial automorphisms of projective planes}


\begin{definition}\label{d:central-axial} An automorphism $A:\Pi\to\Pi$ of a projective plane $\Pi$ is called 
\begin{itemize}
\item \index{central automorphism}\index{automorphism!central}\defterm{central} if there exists a point $c\in \Pi$ (called a \defterm{centre} of $A$) such that $A[\Aline xc]=\Aline xc$ for all $x\in \Pi$;
\item \index{axial automorphism}\index{automorphism!axial}\defterm{axial} if there exists a line $L\subset \Pi$ (called an \defterm{axis} of $A$) such that $A(x)=x$ for all $x\in L$.
\end{itemize}
\end{definition}

Observe that an automorphism $A$ of a projective plane $\Pi$ is axial if and only if it is hyperfixed (because hyperplanes in planes are just lines). Theorem~\ref{t:central<=>hyperfixed} implies that both notions defined in Definition~\ref{d:central-axial} are equivalent.

\begin{corollary}\label{c:central<=>axial} An automorphism $A$ of a projective plane $\Pi$ is central iff it is axial.
\end{corollary}

By Proposition~\ref{p:cenralauto-has1-center} and Corollary~\ref{c:central<=>axial}, a non-identity automorphism $A$ of a projective plane can have at most one center (and one axis), which is called {\em the} center (and {\em the} axis) of the automorphism $A$. 

Let $\Pi$ be a projective plane and $\mathcal L_\Pi$ be its family of lines. 
A pair $(p,L)\in\Pi\times\mathcal L_\Pi$ will be called a \index{point-line pair}\defterm{point-line pair} in $\Pi$.

For any point-line pair $(p,L)$ in a projective plane $\Pi$, consider the set $\Aut_{p,L}(\Pi)$  of all automorphisms $A:\Pi\to\Pi$ with centre $p$ and axis $L$. It is easy to see that $\Aut_{p,L}(\Pi)$ is a subgroup of the automorphism group $\Aut(\Pi)$. Let us recall that $$\Fix(F)\defeq\{x\in\dom[F]:F(x)=x\}$$ denotes the set of fixed point of a function $F$.

\begin{proposition}\label{p:Fix(A)=L+p} For any point-line pair $(p,L)$ in a projective plane $\Pi$ and any non-identity automorphism $A\in\Aut_{p,L}(\Pi)$ we have $\Fix(A)=\{p\}\cup L$.
\end{proposition}

\begin{proof} The inclusion $\{p\}\cup L\subseteq \Fix(A)$ is trivial. The inclusion $\Fix(A)\subseteq \{p\}\cup L$ will be proved with the help of the following lemma.

\begin{lemma}\label{l:Fix(A)=L+p} If $a,b\in\Pi$ are two distinct points with $b=A(a)$, then $A(x)\ne x$ for every point $x\in \Pi\setminus(\{p\}\cup L\cup\Aline ab)$.
\end{lemma}

\begin{proof} Consider the unique point $\lambda\in L\cap\Aline ax$ and observe that $y\defeq A(x)\in A[\Aline \lambda a\cap\Aline xp]=\Aline\lambda b\cap\Aline xp$. Assuming that $y=x$, we conclude that $b\in\Aline \lambda y\cap\Aline ap=\Aline \lambda x\cap\Aline ap=\{a\}$, which contradicts the choice of the points $a,b$.
\end{proof}

Assuming that $\Fix(A)\not\subseteq\{p\}\cup L$, we can find a point $x\in \Fix(A)\setminus(\{x\}\cup L)$. Since $A$ is a non-identity automorphism, there exists a point $a\in\Pi$ such that $b\defeq A(a)\ne a$. Since $A(x)=x$, Lemma~\ref{l:Fix(A)=L+p} ensures that $x\in \Aline ab$. Since $L\subseteq \Fix(A)\subseteq \Pi\setminus\{a,b\}$, the line $\Aline ab$ is distinct from the line $L$ and hence $L\cap\Aline ab=\{\lambda\}$ for a unique point $\lambda$. Since $p$ is a centre of $A$, $\{\lambda,a,b,p\}\subseteq\Aline ab$ and hence $|\Pi|_2\ge|\Aline ab|\ge 4$. By Proposition~\ref{p:cov-aff}, there exists a point $c\in \Pi\setminus (\{p\}\cup L\cup\Aline ab)$. Lemma~\ref{l:Fix(A)=L+p} ensures that $d\defeq A(c)\ne c$. Since the automorphism $A$ has centre $p$, the choice of the point $c\notin\Aline ab=\Aline ap$ ensures that $\Aline ab\cap\Aline cd=\Aline ap\cap\Aline cd=\{p\}$ and hence $x\in\Aline ap\setminus \Aline cd$. Applying Lemma~\ref{l:Fix(A)=L+p}, we conclude that $A(x)\ne x$, which contradicts the choice of the point $x$.
\end{proof}

\section{Desarguesian point-pairs and transitivity in projective planes}

In this section we discuss the notion of a Desarguesian point-line pair in a projective plane and the related notion of point-line transitivity. For a liner $X$ we denote by $\mathcal L_X$, the family of all lines in $X$.

\begin{definition} A projective plane $\Pi$ is called
\index{$(p,L)$-transitive projective plane}\index{projective plane!$(p,L)$-transitive}\defterm{$(p,L)$-transitive} for a point-line pair $(p,L)\in\Pi\times\mathcal L_\Pi$ if for every points $x\in \Pi\setminus(\{p\}\cup L)$ and $y\in \Aline px\setminus(\{p\}\cup L)$, there exists an automorphism $A\in\Aut_{p,L}(\Pi)$ such that $A(x)=y$.
\end{definition}

The point-line transitivity is preserved by automorphisms in the following sense.

\begin{proposition}\label{p:pL-Auto} If a projective plane $\Pi$ is $(p,L)$-transitive for a point-line pair $(p,L)$, then for every automorphism $A\in\Aut(\Pi)$, the plane $\Pi$ is $(q,\Lambda)$-transitive for the point $q\defeq A(p)$ and line $\Lambda\defeq A[L]$.
\end{proposition}

\begin{proof} Given any points $x\in\Pi\setminus(\{q\}\cup\Lambda)$ and $y\in \Aline xq\setminus(\{q\}\cup\Lambda)$, consider the points $x'\defeq A^{-1}(x)$ and $y'\defeq A^{-1}(y)$ and observe that $x'\in A^{-1}[\Pi\setminus(\{q\}\cup\Lambda)]=\Pi\setminus(\{p\}\cup L)$ and $y'=A^{-1}(y)\in A^{-1}[\Aline xq\setminus(\{q\}\cup\Lambda)]=\Aline {x'}p\setminus (\{p\}\cup L)$. By the $(p,L)$-transitivity of $\Pi$, there exists an automorphism $\Phi\in\Aut_{p,L}(\Pi)$ such that $\Phi(x')=y'$. Then the automorphism $F\defeq A\Phi A^{-1}\in \Aut_{q,\Lambda}(\Pi)$ has $F(x)=y$, witnessing that the plane $\Pi$ is $(q,\Lambda)$-transitive.
\end{proof}

The $(p,L)$-transitivity can be characterized by a geometric condition, which is a ``point-line'' specification of the Desargues Axiom.

\begin{definition} A point-line pair $(p,L)\in\Pi\times\mathcal L_\Pi$ is called \index{Desarguesian point-line pair}\index{point-line pair!Desarguesian}\defterm{Desarguesian} if for every distinct lines $A,B,C\in\mathcal L_p\setminus\{L\}$  and distinct points $a,a'\in A\setminus(\{p\}\cup L)$, $b,b'\in B\setminus(\{p\}\cup L)$, $c,c'\in C\setminus(\{p\}\cup L)$, $(\Aline ab\cap\Aline {a'}{b'})\cup(\Aline bc\cap\Aline {b'}{c'})\subseteq L$ implies  $\Aline ac\cap\Aline {a'}{c'}\subseteq L$. 
\end{definition}

\begin{remark} A projective plane is Desarguesian if and only if every point-line pair $(p,L)\in\Pi\times\mathcal L_\Pi$ is Desarguesian. 
\end{remark}

The following theorem was first proved by Baer in \cite{Baer1942}.

\begin{theorem}[Baer, 1942]\label{t:Baer-pL-Des<=>}  A point-line pair $(p,L)$ in a projective plane $\Pi$ is Desarguesian if and only if $\Pi$ is $(p,L)$-transitive.
\end{theorem}

\begin{proof} Let $(p,L)$ be a point-line pair in a projective plane $\Pi$. 
\smallskip

 To prove the ``if'' part, assume that the projective plane $\Pi$ is $(p,L)$-transitive. 
 To prove that the point-pair $(p,L)$ is Desarguesian, take any distinct lines $A,B,C\in\mathcal L_p\setminus\{L\}$ and distinct points $a,a'\in A\setminus(\{p\}\cup L)$, $b,b'\in B\setminus(\{p\}\cup L)$, $c,c'\in C\setminus(\{p\}\cup L)$ such that $(\Aline ab\cap\Aline {a'}{b'})\cup(\Aline bc\cap\Aline {b'}{c'})\subseteq L$. We have to prove that $\Aline ac\cap \Aline {a'}{c'}\subseteq L$. Since $\Pi$ is $(p,L)$-transitive,  there exists an automorphism $F\in\Aut_{p,L}(\Pi)$ such that $F(b)=b'$. Since the plane $\Pi$ is projective, there exist unique points $x,y,z\in \Pi$ such that $\{x\}=\Aline ab\cap\Aline{a'}{b'}$, $\{y\}=\Aline bc\cap\Aline {b'}{c'}$ and $\{z\}=\Aline ac\cap L$. The choice of the points $a,b,c,a',b',c'$ ensures that $x,y\in L$. Taking into account that $L$ is the axis and $p\in A\cap B\cap C$ is the centre of the automorphism $F$, we conclude that $F(z)=z\in L$ and $F(a)\in F[A\cap \Aline ab]=F[A]\cap F[\Aline xb]=A\cap \Aline x{b'}=\{a'\}$. By analogy we can prove that $F(c)=c'$. Then $z=F(z)\in F[\Aline ac]=\Aline{a'}{c'}$ and hence $\Aline ac\cap\Aline {a'}{c'}=\{z\}\subseteq L$.
\smallskip

The proof of the ``only if'' part is more difficult. Assume that the point-line pair $(p,L)$ is Desarguesian. Given any points $a\in \Pi\setminus(\{p\}\cup L)$ and $a'\in\Aline ap\setminus(\{p\}\cup L)$, we have to construct an automorphism $F\in\Aut_{p,L}(\Pi)$ such that $F(a)=a'$. If $a'=a$, then the identity automorphism of $\Pi$ has the required property. So, assume that $a'\ne a$.

If $|\Pi|_2=3$, then $\Pi$ is a Desarguesian (Steiner) plane, by Proposition~\ref{p:Steiner+projective=>Desargues}. Then $\Aline ap=\{a,a',p\}$ and $\Aline ap\cap L=\{p\}$. Choose any point $\lambda\in L\setminus\{p\}$ and observe that the sets $\{p,\lambda,a\}$ and $\{p,\lambda,a'\}$ are maximal independent in $\Pi$. By Theorem~\ref{t:Des-autoextend}, there exists an automorphism $F:\Pi\to\Pi$ of the Desarguesian projective plane $\Pi$ such that $Fp\lambda a=p\lambda a'$. Then $F[L]=F[\Aline p\lambda]=\Aline p\lambda=L$. Since the set  $L\setminus\{p,\lambda\}$ consists of a unique line, $$F[L\setminus\{p,\lambda\}]=F[L]\setminus\{F(p),F(\lambda)\}=L\setminus\{p,\lambda\}$$ and hence $F(x)=x$ for every $x\in L$. Therefore, $L$ is the axis of the automorphism $F$. Since the set $\mathcal L_p\setminus\{L,\Aline pa\}$ contains a unique line and $F[L]=L$, $F[\Aline ap]=\Aline {a'}p=\Aline ap$, the point $p$ is a centre of the automorphism $F$. Therefore, $F$ is a required automorphism with center $p$, axis $L$ and $F(a)=a'$. 

Next, assume that $|\Pi|_2\ge 4$, which means that the projective plane $\Pi$ is $4$-long. Given any points $x\in \Aline ap\setminus(\{p\}\cup L)$ and $b\in \Pi\setminus(L\cup \Aline ap)$, find a unique point $\lambda_b \in L\cap\Aline ab$, a unique point $b'\in \Aline pb\cap\Aline {\lambda_b} {a'}$, a unique point $\mu_b\in L\cap\Aline xb$, and a unique point $x_b\in \Aline ap\cap \Aline {\mu_b}{b'}$. 

\begin{picture}(100,155)(-190,-15)
{\linethickness{=0.7pt}
\put(0,0){\color{teal}\line(-1,2){60}}
}
\put(-55,85){$L$}
\put(0,0){\line(0,1){120}}
\put(0,0){\line(1,2){60}}
\put(-60,120){\line(1,0){120}}
\put(-60,120){\line(5,-2){100}}

\put(60,120){\line(-3,-2){90}}
\put(40,80){\line(-7,-2){70}}
\put(0,120){\vector(0,-1){24}}
\put(0,80){\vector(0,-1){11}}
\put(60,120){\vector(-1,-2){20}}

\put(0,0){\color{red}\circle*{3}}
\put(-3,-9){$p$}
\put(60,120){\circle*{3}}
\put(63,117){$b$}
\put(-60,120){\circle*{3}}
\put(-68,113){$\lambda_b$}
\put(0,120){\circle*{3}}
\put(-2,123){$a$}
\put(40,80){\circle*{3}}
\put(43,76){$b'$}
\put(-30,60){\circle*{3}}
\put(-38,53){$\mu_b$}
\put(0,96){\circle*{3}}
\put(2,96){$a'$}
\put(0,80){\circle*{3}}
\put(-8,79){$x$}
\put(0,68.6){\circle*{3}}
\put(2,62){$x_b$}
\end{picture}

We are going to show that the point $x_b$ depends only on $x$ and not on $b$.

\begin{claim}\label{cl:Baer1} For any points $b\in \Pi\setminus(L\cup \Aline ap)$ and $c\in \Pi\setminus(L\cup \Aline ap\cup\Aline bp\cup \Aline b{\lambda_b})$, the points $x_b$ and $x_c$ coincide.
\end{claim}

\begin{proof} Since $\Aline ab\cap\Aline {a'}{b'}=\{\lambda_b\}\subset L$, $\Aline ac\cap\Aline{a'}{c'}=\{\lambda_c\}\subset L$, and the point-line pair $(p,L)$ is Desarguesian, $\Aline bc\cap\Aline{b'}{c'}\subseteq L$. Since
$\Aline bc\cap\Aline{b'}{c'}\subseteq L$ and $\{\mu_b\}=\Aline bx\cap\Aline {b'}{x_b}\subseteq L$, and the point-line pair $(p,L)$ is Desarguesian, $\Aline cx\cap\Aline {c'}{x_b}\subseteq L$. Taking into account that $\mu_c$ is the unique point of the intersection $L\cap \Aline cx$, we conclude that 
$\mu_c\in \Aline {c'}{x_b}\cap \Aline {c'}{x_c}$ and hence $\{x_b,x_c\}\in \Aline xp\cap\Aline{c'}{\mu_c}$ and finally, $x_b=x_c$.
\end{proof}

\begin{claim}\label{cl:Baer2} For any points $b\in\Pi\setminus(L\cup \Aline ap)$ and  $d\in \Aline bp\setminus(\{p\}\cup L)$, the points $x_b$ and $x_d$ coincide.
\end{claim}

\begin{proof} Since the projective plane $\Pi$ is $4$-long, there exist a line $C\in\mathcal L_p\setminus\{L,\Aline ap,\Aline bp\}$ and a point $c\in C\setminus(L\cup\Aline b{\lambda_b}\cup\Aline d{\lambda_d})$. Claim~\ref{cl:Baer1} ensures that $x_b=x_c=x_d$.
\end{proof}

\begin{claim}\label{cl:Baer3}  For any points  $b\in\Pi\setminus(L\cup \Aline ap)$ and $c\in \Pi\setminus(L\cup \Aline ap)$, the points $x_b$ and $x_c$ coincide.
\end{claim}

\begin{proof} If $c\in \Aline bp$, then $x_b=x_c$, by Claim~\ref{cl:Baer2}. So, assume that $c\notin\Aline bp$. If $c\notin\Aline b{\lambda_b}$, then $x_b=x_c$, by Claim~\ref{cl:Baer1}. So, assume that $c\in \Aline b{\lambda_b}$. Choose any point $d\in\Aline bp\setminus(\{p,b\}\cup L)$ and observe that $\Aline ad\ne \Aline ab$, which implies $\Aline b{\lambda_b}\cap\Aline d{\lambda_d}=\{a\}$ and hence $c\in \Aline b{\lambda_b}\setminus A=\Aline b{\lambda_b}\setminus \Aline d{\lambda_d}$. By Claims~\ref{cl:Baer2} and \ref{cl:Baer1}, $x_b=x_d=x_c$.
\end{proof}

Now we are able to define the automorphism $F\in\Aut_{p,L}(\Pi)$ such that $F(a)=a'$. Given any point $x\in \Pi\setminus \Aline ap$, find a unique point $\lambda\in \Aline xa\cap L$ and let $F(x)$ be the unique point of the intersection $\Aline xp\cap\Aline \lambda {a'}$. If $x\in L$, then $F(x)=\lambda=x$. For every point $x\in \Aline ap\cap(\{p\}\cup L)$ put $F(x)\defeq x$, and for every point $x\in \Aline ap\setminus(\{p\}\cup L)$, put $F(x)\defeq x_b$, where $b\in \Pi\setminus(L\cup \Aline ap)$ is any point. Claim~\ref{cl:Baer3} ensures that the value $F(x)\defeq x_b$ does not depend on the choice of the point $b$, so the function $F$ is well-defined. It is easy to see that $F$ is a bijection of $\Pi$ such that $\{p\}\cup L\subseteq\Fix(F)$ and $F[\Aline xp]=\Aline xp$ for every $x\in\Pi$. 

It remains to prove that $F$ is an automorphism of the projective plane $\Pi$. Given any line $\Lambda\subset \Pi$, we should check that $F[\Lambda]$ is a line in $\Pi$. If $\Lambda=L$ or $L\in\mathcal L_p$, then $F[\Lambda]=\Lambda$ by the definition of the function $F$. So, assume that $L\ne \Lambda\notin\mathcal L_p$. Let $\lambda$ be the unique common point of the lines $\Lambda$ and $L$. 

If $\lambda\notin \Aline ap$, then consider the unique point $x\in \Lambda\cap \Aline ap$ and its image $y\defeq F(x)$. Assuming that $x\in L$, we conclude that $x\in L\cap \Lambda=\{\lambda\}$, which contradicts the assumption $\lambda\notin \Aline ap$. So, $x\notin L$ and hence $y=F(x)=x_b$ for every point $b\in \Lambda\setminus(L\cup \Aline ap)=\Lambda\setminus\{\lambda,x\}$. The definition of the function $F$ at the point $b$ ensures that $F(b)=b'\in \Aline {a'}{\lambda_b}$, where $\lambda_b$ is a unique common point of the lines $\Aline ab$ and $L$. On the other hand, the definition of the point $y=F(x)=x_b$ ensures that $y\in \Aline ap\cap\Aline \lambda{b'}$. Then $F(b)=b'\in \Aline \lambda y$ and hence $F[\Lambda]=\Aline \lambda y$ is a line in the plane $\Pi$.

If $\lambda\in \Aline ap$, then choose any point $b\in\Lambda\setminus\{\lambda\}$ and consider its image $F(b)=b'$. The definition of the function $F$ ensures that $b'\in \Aline {a'}{\lambda_b}$, where $\lambda_b\in \Aline ab\cap L$. We claim that $F[\Lambda]=\Aline \lambda{b'}$. Indeed, for any point $c\in \Lambda\setminus\{b,\lambda\}$, the definiton of $F(c)\defeq c'$ ensures that $c'\in \Aline {a'}{\lambda_c}$, where $\lambda_c$ is a unique common point of the lines $L$ and $\Aline ac$. Taking into account that  
$$(\Aline ab\cap\Aline {a'}{b'})\cup(\Aline ac\cap\Aline{a'}{c'})=\{\lambda_b\}\cup\{\lambda_c\}\subseteq L$$ and the point-line pair $(p,L)$ is Desarguesian, we conclude that $\varnothing\ne\Aline bc\cap\Aline{b'}{c'}\subseteq L$. Then $\varnothing\ne \Aline bc\cap\Aline {b'}{c'}\cap L\subseteq \Aline bc\cap L=\{\lambda\}$ and $F(c)=c'\in\Aline \lambda{b'}$, witnessing that $F[\Lambda]=\Aline \lambda{b'}$ is a line in $\Pi$.

Therefore, in both cases we have proven that the image $F[\Lambda]$ of a line $\Lambda\subset \Pi$ is a line, which means that $F$ is a required automorphism of the projective plane $\Pi$ with centre $p$, axis $L$, and $F(a)=a'$.
\end{proof}

Now we prove some inheritance properties of the $(p,L)$-transitivity that will be essentially used in the Lenz and Lenz--Barlotti classifications of projective planes.

\begin{theorem}\label{t:pL-transitive-dual} A projective plane $\Pi$ is $(p,L)$-transitive for some point-line pair $(p,L)$ if and only if the dual projective plane $\mathcal L_\Pi$ is $(L,\mathcal L_p)$-transitive.
\end{theorem}

\begin{proof} Assume that $\Pi$ is $(p,L)$-transitive for some point-line pair $(p,L)\in\Pi\times\mathcal L_\Pi$. To prove that the dual projective plane $\mathcal L_\Pi$ is $(L,\mathcal L_p)$-transitive, take any distinct lines $A,B\in \mathcal L_\Pi\setminus(\{L\}\cup\mathcal L_p)$ with $A,B,L\in\mathcal L_o$ for some point $o\in \Pi$. Since $A,B\notin(\{L\}\cup\mathcal L_p)$, there exists points $x\in A\setminus(L\cup\{p\})$ and $y\in \Aline px\cap B$. By the $(p,L)$-transitivity of $\Pi$, there exists an automorphism $F\in\Aut_{p,L}(\Pi)$ such that $F(x)=y$. It follows from $o\in L$ that $F(o)=o$ and hence $F[A]=F[\Aline ox]=\Aline oy=B$. The automorphism $F$ induces the automorphism $[F]:\mathcal L_\Pi\to\mathcal L_\Pi$, $[F]:\Lambda\mapsto F[\Lambda]$, of the dual projective plane $\mathcal L_\Pi$ such that $[F](L)=F[L]=L$ and $[F](\Lambda)=F[\Lambda]=\Lambda$ for every $\Lambda\in\mathcal L_p$. Therefore, $[F]\in \Aut_{L,\mathcal L_p}(\mathcal L_\Pi)$ and $[F](A)=F[A]=B$, witnessing that the dual projective plane $\mathcal L_\Pi$ is $(L,\mathcal L_p)$-transitive. 

Since the projective plane $\Pi$ can be identified with the dual projective plane to $\mathcal L_\Pi$, the $(L,\mathcal L_p)$-transitivity of $\mathcal L_\Pi$ implies the $(p,L)$-transitivity of $\Pi$, by the preceding argument.
\end{proof}

\begin{theorem}\label{t:(pq,L)-transitive} If a projective plane $\Pi$ is $(p,L)$-transitive and $(q,L)$-transitive for some  points $p,q\in\Pi$ and line $L\subset \Pi$, then $\Pi$ is $(c,L)$-transitive for every point $c\in\Aline pq$.
\end{theorem}

\begin{proof} If $|\Pi|_2\le 4$, then the projective plane $\Pi$ is Desarguesian (by Proposition~\ref{p:Steiner+projective=>Desargues} and Corollary~\ref{c:4-Pappian}) and hence $(c,L)$-transitive, by Baer's Theorem~\ref{t:Baer-pL-Des<=>}. So, assume that $|\Pi|_2\ge 5$. If $c\in\{p,q\}$, then the $(c,L)$-transitivity of $\Pi$ follows from the $(p,L)$-transitivity and $(q,L)$-transitivity of $\Pi$. So, assume that $c\in\Aline pq\setminus\{p,q\}$.

\begin{claim}\label{cl:pLqL=>cL} For every points $x\in \Pi\setminus(L\cup \Aline pq)$ and $y\in \Aline xc\setminus(\{x,c\}\cup L)$ with $\Aline px\cap\Aline qy\not\subseteq L$, there exists an automorphism $A\in\Aut_{c,L}(\Pi)$ such that $A(x)=y$.
\end{claim}

\begin{proof} Consider the unique points $\lambda\in L\cap\Aline xc$ and $z\in\Aline xp\cap\Aline yq$. Since $p,q\notin \Aline xy$, the points $x,z,y$ are distinct. Our assumption ensures that $z\notin L$. By the $(p,L)$- and $(q,L)$-transitivities of $\Pi$, there exist automorphisms $P\in\Aut_{p,L}(\Pi)$ and $Q\in\Aut_{q,L}(\Pi)$ such that $P(x)=z$ and $Q(z)=y$. Then the automorphism $A\defeq QP:\Pi\to\Pi$ has axis $L$, and  $A(x)=QP(x)=Q(z)=y$. Also $A[\Aline pq]=QP[\Aline pq]=Q[\Aline pq]=\Aline pq$ and hence $A(c)\in A[\Aline pq\cap\Aline x\lambda]=\Aline pq\cap\Aline y\lambda=\{c\}$. So, $c\in \Fix(A)$. By Proposition~\ref{p:Fix(A)=L+p}, $\Fix(P)=\{p\}\cup L$ and $\Fix(Q)=\{q\}\cup L$, which implies that $\Fix(A)\subseteq \Aline pq\cup L$ and hence $\{c\}=\Aline xy\cap \Fix(A)$. By Corollary~\ref{c:central<=>axial}, the axial automorphism $A$ is central with centre $c'\in \Fix(A)\cap\Aline xy=\{c\}$. Therefore, $A$ is a required automorphism of $\Pi$ with axis $L$ and centre $c$ such that $A(x)=y$.
\end{proof}

\begin{claim}\label{cl:pLqL=>cL2} For any distinct points $x\in \Pi\setminus(L\cup \Aline pq)$ and $y\in \Aline xc\setminus(\{x,c\}\cup L)$, there exists an automorphism $A\in\Aut_{c,L}(\Pi)$ such that $A(x)=y$.
\end{claim}

\begin{proof} If $\Aline px\cap\Aline qy\not\subseteq L$,  then the automorphism $A$ exists by Claim~\ref{cl:pLqL=>cL}. If $\Aline py\cap\Aline qx\not\subseteq L$, then by Claim~\ref{cl:pLqL=>cL}, there exists an automorphism $B\in\Aut_{c,L}(\Pi)$ such that $B(y)=x$. Then the automorphism $A\defeq B^{-1}\in\Aut_{c,L}(\Pi)$ has $A(x)=y$.

So, assume that $(\Aline px\cap\Aline qy)\cup(\Aline py\cap\Aline qx)\subseteq L$. Since $|\Pi|_2\ge 5$, there exists a point $z\in \Aline xy\setminus(\{c,x,y\}\cup L)$. Since $x\ne z\ne y$, the inclusion $(\Aline px\cap\Aline qy)\cup(\Aline py\cap\Aline qx)\subseteq L$ implies $\Aline px\cap\Aline qz\not\subseteq L$ and $\Aline pz\cap\Aline qy\not\subseteq L$. By Claim~\ref{cl:pLqL=>cL}, there exist automorphisms $A,B\in\Aut_{c,L}(\Pi)$ such that $A(x)=z$ and $B(z)=y$. Then the automorphism $C\defeq BA$ has the required properties: it has axis $L$, center $c$, and $C(x)=y$.
\end{proof}
  
The following claim completes the proof of the theorem.

\begin{claim} For every points $x\in \Pi\setminus(\{c\}\cup L)$ and $y\in \Aline xc\setminus(\{c\}\cup L)$, there exists an automorphism $A\in\Aut_{c,L}(\Pi)$ such that $A(x)=y$. 
\end{claim}

\begin{proof} If $y=x$, the the identity automorphism of $\Pi$ has the required property. So, assume that $x\ne y$.  If $x\notin\Aline pq$, then the automorphism $A$ exists by Claim~\ref{cl:pLqL=>cL2}. So, assume that $x\in \Aline pq$ and hence $y\in \Aline xc=\Aline pq$. Choose any points $\lambda\in L\setminus\Aline pq$, $x'\in \Aline x\lambda\setminus\{x,\lambda\}$ and consider the unique point $y'\in \Aline \lambda y\cap\Aline c{x'}$. By Claim~\ref{cl:pLqL=>cL2}, there exists an automorphism $A\in\Aut(\Pi;c,L)$ such that $A(x')=y'$. Taking into account that $c\in\Aline pq=\Aline xy$ is a center of $A$ and $\lambda\in L\subseteq\Fix(A)$, we conclude that $A(x)\in A[\Aline pq\cap\Aline {x'}\lambda]=\Aline pq\cap \Aline {y'}\lambda=\{y\}$.
\end{proof}
\end{proof}

\begin{corollary}\label{c:(p,q]-transitive} If a projective plane $\Pi$ is $(p,L)$-transitive and $(p,L')$-transitive for a point $p$ and two distinct lines $L,L'$, then $\Pi$ is $(p,\Lambda)$-transitive for every line $\Lambda\in\mathcal L_o$, where $o$ is a unique common point of the lines $L,L'$.
\end{corollary}

\begin{proof} Let $o$ be the unique common point of the lines $L,L'$. By Theorem~\ref{t:pL-transitive-dual}, the dual projective plane $\mathcal L_\Pi$ is $(L,\mathcal L_p)$-transitive and $(L',\mathcal L_p)$-transitive, and by Theorem~\ref{t:(pq,L)-transitive}, $\mathcal L_\Pi$ is $(\Lambda,\mathcal L_p)$-transitive for every line $\Lambda\in \mathcal L_o$. By Theorem~\ref{t:pL-transitive-dual}, $\Pi$ is $(p,\Lambda)$-transitive for every line $\Lambda\in\mathcal L_o$.
\end{proof}

Theorem~\ref{t:(pq,L)-transitive} and Corollary~\ref{c:(p,q]-transitive} motivate the following extensions of the point-line transitivity.

\begin{definition} A projective plane $\Pi$ is 
\begin{itemize}
\item \index{$(p,q]$-transitive projective plane}\index{projective plane!$(p,q]$-transitive}\defterm{$(p,q]$-transitive} for two points $p,q\in \Pi$ if $\Pi$  is $(p,L)$-transitive for every line $L\in\mathcal L_q$;
\item \index{$[\Lambda,L)$-transitive projective plane}\index{projective plane!$[\Lambda,L)$-transitive}\defterm{$[\Lambda,L)$-transitive} for two lines $\Lambda,L\subset \Pi$  if $\Pi$ is $(p,L)$-transitive for every point $p\in \Lambda$.
\end{itemize}
\end{definition}

\begin{theorem}\label{t:pq<=>qp} Let $p,q$ be two distinct points in a projective plane $\Pi$. The projective plane $\Pi$ is $(p,q]$-transitive if and only if $\Pi$ is $(q,p]$-transitive.
\end{theorem}

\begin{proof} It suffices to show that the $(p,q]$-transitivity of $\Pi$ implies the $(q,p]$-transitivity of $\Pi$. So, assume that the projective plane $\Pi$ is $(p,q]$-transitive. Choose a projective base $uowe$ in $\Pi$ so that its horizontal infinity point $h\in\Aline ou\cap\Aline we$ coinsides with $q$ and the vertical infinity point $v\in\Aline ow\cap\Aline ue$ coincides with the point $p$. Since $(p,q)=(v,h)$, the $(p,q]$-transitivity of $\Pi$ implies the $(v,\Aline hv)$-transitivity of $\Pi$ and the $(v,\Aline oh)$-transitivity of $\Pi$. Then the based affine plane $(\Pi\setminus\Aline hv,uow)$ is vertical-translation and vertical-scale. By Theorem~\ref{t:cart-group<=>} and Proposition~\ref{p:lv-scale=>rdist+ass-dot}, the ternar $\Delta=\Aline oe\setminus\Aline hv$ of the based affine plane $(\Pi,uow)$ is linear, right-distributive and associative. By  Theorems~\ref{t:VW-Thalesian<=>quasifield} and \ref{t:horizontal-scale}, the based affine plane $(\Pi\setminus\Aline hv,uow)$ is translation and horizontal-scale. Then the projective plane $\Pi$ is $(h,\Aline hv)$-transitive and $(h,\Aline ov)$-transitive. By Corollary~\ref{c:(p,q]-transitive}, the projective plane $\Pi$ is $(h,v]$-transitive. Since $(h,v)=(q,p)$, the projective plane $\Pi$ is $(q,p]$-transitive.
\end{proof}

\begin{corollary}\label{c:LL'<=>L'L} Let $L,\Lambda$ be two distinct lines in a projective plane $\Pi$. The projective plane $\Pi$ is $[L,\Lambda)$-transitive if and only if $\Pi$ is $[\Lambda,L)$-transitive.
\end{corollary}

\begin{proof} If $\Pi$ is $[L,\Lambda)$-transitive, then by Theorem~\ref{t:pL-transitive-dual}, the dual projective plane $\mathcal L_\Pi$ is $(\Lambda,L]$-transitive. By Theorem~\ref{t:pq<=>qp}, $\mathcal L_\Pi$ is $(L,\Lambda]$-transitive and by Theorem~\ref{t:pL-transitive-dual}, $\Pi$ is $[\Lambda,L)$-transitive. By analogy we can prove that the $[\Lambda,L)$-transitivity of $\Pi$ implies its $[L,\Lambda)$-transitivity.
\end{proof}

\begin{theorem}\label{t:circledast} If a projective plane $\Pi$ is $(h,\Aline oh)$-transitive and $(v,\Aline ov)$-transitive for three distinct points $h,o,v\in\Pi$, then $\Pi$ is $(p,\Aline op)$-transitive for every point $p\in\Aline hv\setminus\{o\}$.
\end{theorem} 

\begin{proof} If $o\in \Aline hv$, then the lines $\Aline vo$ and $\Aline oh$ coincide and hence $\Pi$ is $(h,L)$-transitive and $(v,L)$-transitive for the line $L\defeq \Aline hv$. By Theorem~\ref{t:(pq,L)-transitive}, the projective plane $\Pi$ is $(p,L)$-transitive for every point $p\in \Aline hv$. If $p\in\Aline hv\setminus\{o\}$, then $\Aline op=L$ and the $(p,L)$-transitivity implies the $(p,\Aline op)$-transitivity.

So, assume that $o\notin \Aline hv$. In this case we can choose any point $e\in\Pi\setminus(\Aline hv\cup\Aline oh\cup\Aline ov)$, find unique points $u\in \Aline oh\cap \Aline ve$ and $w\in \Aline ov\cap\Aline he$, and observe that $uowe$ is a projective base for the projective plane $\Pi$. Observe that the line $\Aline hv$ can be identified with the boundary $\partial(\Pi\setminus\Aline hv)$ of the affine plane $\Pi\setminus\Aline hv$. The $(h,\Aline oh)$-transitivity and  $(v,\Aline ov)$-transitivity of the projective plane $\Pi$ implies that the affine plane $\Pi\setminus\Aline hv$ is horizontal-shear and vertical-shear. By Lemma~\ref{l:Lshear}, for every point $p\in \Aline hv$, the set $\Aline op\setminus\{p\}$ is a shear line in the affine plane $\Pi\setminus\Aline hv$, which implies that the projective plane $\Pi$ is $(p,\Aline op)$-transitive.
\end{proof}

\section{Algebra versus Geometry in based projective planes}

The main result of this section is the following characterization theorem describing the interplay between point-line transitivelation properties of a based projective plane and algebraic properties of its ternar.

\begin{theorem}\label{t:Algebra-vs-Geometry-proj} Let $(\Pi,uowe)$ be a based projective plane, $h,v$ be its horizontal and vertical infinity points, and $\Delta=\Aline oe\setminus\Aline hv$ be the ternar of the based affine plane $\IA=(\Pi\setminus\Aline hv,uow)$ and the based projective plane $(\Pi,uowe)$. The ternar $\Delta$ is
\begin{enumerate}
\item linear and associative-plus iff $\Pi$ is $(v,\Aline hv)$-transitive\\ iff $\IA$ is vertical-translation;
\item linear, right-distributive and associative-plus iff $\Pi$ is $[\Aline hv,\Aline hv)$-transitive\\ iff $\IA$ is translation;
\item linear, left-distributive and associative-plus iff $\Pi$ is $(v,v]$-transitive\\ iff $\IA$ is vertical-translation and vertical-shear;
\item linear, distributive and associative-plus iff $\Pi$ is $[\Aline hv,\Aline hv)$-transitive and $(v,v]$-transitive iff $\IA$ is translation and vertical-shear;
\smallskip
\item linear and associative-dot iff $\Pi$ is $(h,\Aline ov)$-transitive iff $\IA$ is horizontal-scale;
\item linear, right-distributive and associative-dot iff $\Pi$ is $(h,\Aline ov)$-transitive and $(v,\Aline oh)$-transitive iff $\IA$ is horizontal-scale and vertical-scale;
\item linear, left-distributive and associative-dot iff $\Pi$ is $(h,\Aline ov)$-transitive and $(o,\Aline hv)$-transitive iff $\IA$ is horizontal-scale and homocentral;
\item linear, distributive and associative-dot iff $\Pi$ is $(h,\Aline ov)$-transitive, $(v,\Aline oh)$-transitive and $(o,\Aline hv)$-transitive iff $\IA$ is horizontal-scale, vertical-scale and homocentral;
\smallskip
\item linear and associative iff $\Pi$ is $(v,\Aline hv)$-transitive and $(h,\Aline ov)$-transitive iff $\IA$ is vertical-translation and horizontal-scale;
\item linear, left-distributive and associative iff $\Pi$ is $[\Aline ov,\Aline hv)$-transitive iff $\Pi$ is $[\Aline hv,\Aline ov)$-transitive iff $\IA$ is vertical-translation and homocentral iff $\IA$ is vertical-shear and horizontal-scale;
\item linear, right-distributive and associative iff $\Pi$ is $(v,h]$-transitive iff $\Pi$ is $(h,v]$-transitive iff the based affine plane $\IA$ is vertical-translation and vertical-scale iff $\IA$ is horizontal-translation and horizontal-scale iff $\IA$ is translation, vertical-scale and horizontal-scale;
\item linear, distributive and associative iff $\Pi$ is $(p,L)$-transitive for every point-line pair $(p,L)$ iff $\Pi$ is $(p,\Aline hv)$-transitive for every point $p\in\{o,h,v\}$ iff $\Pi$ is Desarguesian iff $\IA$ is Desarguesian iff $\IA$ is a dilation plane;
\smallskip
\item linear, distributive, associative, and commutative iff $\Pi$ Pappian iff $\IA$ is Pappian;
\smallskip
\item linear, distributive, associative-plus and alternative-dot iff $\Pi$ is  $[\Aline hv,\Aline hv)$-transitive,\break $(h,\Aline oh)$-transitive, and $(v,\Aline ov)$-transitive iff the based affine plane $\IA$ is translation, vertical-shear and horizontal-shear iff $\IA$ is shear;
\item linear, distributive, associative-plus and alternative-dot iff\/ $\Pi$ is $[\Aline hv,\Aline hv)$-transitive,\break  $(v,\Aline ov)$-transitive, and $(o,\Aline ov)$-transitive. 
\end{enumerate}
\end{theorem}

\begin{proof} Observe that all point-line translation properties of the projective plane $(\Pi,uowe)$ with respect to points $o,h,v$ and lines $\Aline ov,\Aline oh,\Aline hv$ (except for the $(o,\Aline oh)$-translation and $(o,\Aline ov)$-translation) can be rewritten as translation-shear-scale-homothety properties of the affine plane $\Pi\setminus\Aline hv$. This allows us to deduce the statements $(1)$--$(14)$ from the corresponding statements in Corollary~\ref{c:Algebra-vs-Geometry-in-trings:affine} (for the statements (10) and (11),  one should also apply Theorems~\ref{t:pq<=>qp} and Corollary~ \ref{c:LL'<=>L'L}).  It remains to prove the statement (15).

To prove the ``if'' part of (15), assume that the based projective plane $(\Pi,uowe)$ is $[\Aline hv,\Aline hv)$-transitive, $(v,\Aline ov)$-transitive, and $(o,\Aline ov)$-transitive. Then the based affine plane $\IA=(\Pi\setminus\Aline hv,uow)$ is translation and vertical-shear (we have not use the $(o,\Aline ov)$-transitivity yet). By  Theorem~\ref{t:semifield<=>}, the ternar $\Delta=\Aline oe\setminus\Aline hv$ of the affine plane $\IA$ is linear, distributive and associative-plus. By Hankel's Theorem~\ref{t:Hankel}, the ring $\Delta$ is commutative-plus. Let $I\defeq\Aline hv$ be the line at infinity for the affine plane $\IA$. Let $f\in \Aline ou\setminus I$ and $e'\in \Delta$ be unique points such that the line $\Aline fw\setminus I$ is parallel to the diagonal $\Delta$ and the line $\Aline f{e'}\setminus I$ is parallel to the vertical axis $\Aline ow\setminus I$ in the affine plane $\Pi\setminus I$. The definition of addition in the ternar $\Delta$ ensures that $e'+e=o$ and hence $e'$ is the additive inverse to $e$ in $\Delta$.

\begin{picture}(60,70)(-160,-15)
\linethickness{0.7pt}
\put(0,0){\color{red}\line(1,1){40}}
\put(0,0){\color{cyan}\line(0,1){20}}
\put(0,20){\color{red}\line(1,1){20}}
\put(20,0){\color{cyan}\line(0,1){50}}
\put(20,40){\color{teal}\line(1,0){20}}
\put(-10,20){\color{teal}\line(1,0){60}}

\put(0,0){\circle*{2.5}}
\put(-2,-8){$e'$}
\put(20,20){\circle*{2.5}}
\put(22,13){$o$}
\put(0,20){\circle*{2.5}}
\put(-8,12){$f$}
\put(20,40){\circle*{2.5}}
\put(10,38){$w$}
\put(40,20){\circle*{2.5}}
\put(37,12){$u$}
\put(40,40){\circle*{2.5}}
\put(43,38){$e$}
\end{picture}

Since the projective plane $\Pi$ is $(o,\Aline ov)$-transitive and $o\in \Aline fh$, there exists an automorphism $\Phi\in \Aut_{o,\overline{o\,v}}(\Pi)$ such that $\Phi(h)=f$.
Let us write down the action of $\Phi$ in coordinates of the based affine plane $\IA=(\Pi\setminus I,uow)$ (we use standard notations $L_b$ and $L_{a,b}$ for lines in the based affine plane $\IA$). For a line $L$ in $\Pi\setminus I$ we denote by $\overline L$ the flat hull of the line $L$ in the projective plane $\Pi$.

Since $o$ is the centre of the automorphism $\Phi$, $\Phi[\Aline oe]=\Aline oe$. For every $b\in \Delta\setminus\{e\}$, let $b'\defeq\Aline {b}h\cap \Aline ov$ and observe that $L_{o,b}=\Aline {b}h\setminus\{h\}=\Aline {b'}h\setminus\{h\}$.  Taking into account that ${b'}\in \Aline ov\subseteq\Fix(\Phi)$ and $\Phi(h)=f$, we conclude that $\Phi[L_{o,b}]\subseteq \Phi[\Aline {b'}h]=\Aline {b'}f$. Then for the point $b\in L_{e,o}\cap L_{o,b}$ with coordinates $(b,b)$, its image $\Phi(b)$ belongs to the intersection of the lines $\Phi[L_{e,o}]=L_{e,o}$ and $\Phi[L_{o,b}]=\Aline f{b'}$. 

\begin{picture}(100,160)(-160,-50)

\linethickness{0.7pt}
\put(0,0){\line(1,0){120}}
\put(0,0){\color{cyan}\line(0,1){80}}
\put(0,0){\line(1,1){85}}
\put(0,0){\line(-1,0){60}}
\put(0,0){\color{cyan}\line(0,-1){30}}
\put(0,0){\line(-1,-1){30}}
\put(-35,60){\color{teal}\vector(1,0){155}}
\put(-20,0){\color{red}\line(1,3){32}}
\put(-20,0){\color{red}\line(-1,-3){10}}
\put(-22,-6){\color{red}\vector(1,3){2}}
\put(-18,6){\color{red}\vector(-1,-3){2}}

\qbezier(-30,65)(-25,90)(5,90)
\put(0,90){\vector(1,0){6}}
\put(-25,86){$\Phi$}

\put(87,87){$L_{e,o}$}
\put(125,58){\color{teal}$h$}
\put(12,100){\color{red}$L_{b,b}$}
\put(-50,53){\color{teal}$L_{o,b}$}

\put(0,0){\circle*{2.5}}
\put(2,-7){$o$}
\put(20,20){\circle*{2.5}}
\put(20,14){$e$}
\put(-20,-20){\circle*{2.5}}
\put(-20,-27){$e'$}
\put(-20,0){\color{red}\circle*{3}}
\put(-26,3){\color{red}$f$}
\put(0,20){\circle*{2.5}}
\put(2,18){$w$}
\put(0,60){\color{red}\circle*{3}}
\put(4,62){$b'$}
\put(60,60){\circle*{3}}
\put(57,63){\color{teal}$b$}
\put(-30,-30){\color{red}\circle*{3}}
\put(-40,-42){\color{red}$\Phi(b)$}

\qbezier(60,55)(60,-30)(-20,-30)
\put(-15,-30){\vector(-1,0){8}}
\put(40,-18){$\Phi$}
\end{picture}

Since the line $\Aline f{b'}\setminus I$ is not vertical and $b'\in \Aline ov\cap\Aline bh$, there exist an element $\beta\in \Delta$ such that $\Aline f{b'}\setminus I=L_{\beta,b}$. Taking into account that the point $f$ has coordinates $(e',o)$ in the based affine plane $(\Pi\setminus I,uow)$, we conclude that $e'{\cdot} \beta+b=o$. On the other hand, the distributivity of the ternar $R$ ensures that $$e'{\cdot} \beta+b=o=(e'+e){\cdot} b=e'{\cdot} b+e{\cdot} b=e'{\cdot} b+b$$ and hence $\beta=b$ by the divisibility of the loops $(\Delta,+)$ and $(\Delta^*,\cdot)$. Therefore, the line $\Aline {b'}f\setminus I$ coincides with the line $L_{b,b}$. Since $b\ne e$,  the lines $L_{e,o}$ and $L_{b,b}$ have a unique common point whose coordinates $(x,x)$ satisfy the equation 
$x=x{\cdot} b+b$, which has a unique  solution $x=b/(e-b)$.
Therefore, for every $b\in\Delta\setminus\{e\}$, the automorphism $\Phi$ maps the point with coordinates $(b,b)$ to the point with coordinates $(b/(e-b),b/(e-b))$. Since the line $\Aline ov$ is the axis of $\Phi$, the automorphism $\Phi$ maps the vertical line $L_b$ to the vertical line $L_{b/(e-b)}$ in the affine plane $\Pi\setminus I$.

Since $o$ is the centre of the automorphism, $\Phi[\overline L_{c,o}]=\overline L_{c,o}$ for every $c\in \Delta$. This imples that for every point $p\in \Pi\setminus (I\cup L_e)$ with coordinates $(x,y)$ its image $\Phi(p)$ has coordinates $(x/(e-x),(x/(e-x))\cdot(x\backslash y))$.

Now observe that for every $b\in \Delta$ the line $L_{e,b}$ contains the points $b'$ and $a'$ with coordinates $(o,b)$ and $(e,e+b)$. It follows from $b'\in \Aline ov\subseteq \Fix(A)$ that $\Phi(b')=b'$. On the other hand, the point $a'$ belongs to the line $L_{e+b,o}$ and its image $\Phi(a')$ belongs to the intersection $\overline L_{e+b,o}\cap I$, which implies that $\Phi[L_{e,b}]\subseteq \overline L_{e+b,b}$.
Then for every point $p\in L_{e,b}\setminus L_e$ with coordinates $$(x,x +b)=(x,x{\cdot}(e+x\backslash b))$$ its image $\Phi(p)$ has coordinates $$\big(x/(e-x),(x/(e-x)){\cdot} (e+x\backslash b)\big),$$which satisfy the equation of the line $L_{e+b,b}$ and hence
$$(x/(e-x)){\cdot} (e+x\backslash b)=(x/(e-x)){\cdot}(e+b)+b.$$Since the ternar $R$ is distributive, associative-plus and commutative-plus, this equation is equivalent to
$$(x/(e-x)){\cdot} (x\backslash b)=(x/(e-x)){\cdot} b+b.$$
Now take any elements $z,s\in \Delta\setminus\{o\}$ and apply the above identity to the elements $x\defeq e-z$ and $b\defeq x{\cdot} s$. Then we obtain the identity
$$((e-z)/z){\cdot} s=((e-z)/z){\cdot} ((e-z){\cdot} s)+(e-z){\cdot} s,$$
which is transformed to
$$(e/z){\cdot} s-s=(e/z-e){\cdot}(s-z{\cdot}s)+s-z{\cdot} s=
(e/z){\cdot} s-(e/z){\cdot} (z{\cdot} s)-s+z{\cdot} s+s-z{\cdot} s
$$
and after cancellations, to the identity
$$(e/z){\cdot} (z{\cdot} s)=s,$$
witnessing that the ternar $\Delta$ has the left inverse property. By Theorem~\ref{t:alternative<=>left-inverse}, $\Delta$ is alternative-dot. Therefore, the ternar $\Delta$ is linear, distributive, associative-plus and alternative-dot.
\smallskip

To prove the ``only if'' part, assume that the ternar $\Delta$ of the based projective plane $(\Pi,uowe)$ is  linear, distributive, associative-plus and alternative-dot. Since $\Delta$ is the ternar of the based affine plane $\IA=(\Pi\setminus\Aline hv,uow)$, we can apply Theorems~\ref{t:VW-Thalesian<=>quasifield}, \ref{t:shear<=>alternative} and conclude that the based affine plane $\IA$ is translation and shear. In particular, $\IA$ is translation and vertical-shear, which implies that the projective plane $\Pi$ is $[\Aline hv,\Aline hv)$-transitive and $(v,\Aline ov)$-transitive. It remains to prove that $\Pi$ is $(o,\Aline ov)$-transitive. 

Let $e'\in \Delta$ be a unique point such that $e+e'=o=e'+e$. Consider the unique point $f\in \Aline {e'}v\cap\Aline oh$. The point $f$ has coordinates $(e',o)$ in the affine base $uow$ of the affine plane $\IA=\Pi\setminus\Aline hv$.

Consider the bijective function $\Phi:\Pi\to\Pi$ defined as follows. Every point $p\in \Pi\setminus\Aline hv$ with coordinates $(x,y)\in (\Delta\setminus\{o,e\})\times\Delta$ in the affine base $uow$ is mapped to the point $\Phi(p)$ with coordinates $(x/(e-x),(x/(e-x))\cdot(x\backslash y))$. Each point $p\in\Aline ov$ is mapped to itself, each point $p\in \Aline ev\setminus\{v\}$ is mapped to the unique common point of the lines $\Aline op$ and $\Aline hv$, and each point $p\in \Aline hv$ is mapped to the unique common point of the lines $\Aline op$ and $\Aline fv$. Then $\Phi(h)=f$. We claim that $\Phi$ is an automorphism of the projective plane. Given any line $L\subseteq \Pi$, we should check that $\Phi[L]$ is a line in $\Pi$. If $L=\Aline hv$, then $\Phi[L]=\Aline fv$ is a line in $\Pi$. If $L=\Aline ov$, then $\Phi[L]=\Phi[\Aline ov]=\Aline ov$. If $L=\Aline ev$, then $\Phi[L]=\Aline hv$ is a line. If $L=\Aline xv$ for some $x\in \Delta\setminus\{o,e\}$, then $\Phi[L]=\Aline yv$ where $y=x/(e-x)\in \Delta$. In all other cases $L=\overline {L_{a,b}}$ for some $a,b\in\Delta$. The line $\overline L_{a,b}$ is the flat hull of the line $L_{a,b}$ in the affine plane $\IA$. In the coordinates of the based affine plane $(\Pi\setminus\Aline hv,uow)$, the line $L_{a,b}$ has the equation $y=x\cdot a+b$. We claim that $\Phi[L]=\Phi[\overline L_{a,b}]=\overline L_{a+b,b}$. Indeed, take any point $p\in L=\overline L_{a,b}$. If $p\in \Aline hv$, then $p\in \bar L_{a,b}\cap \bar L_{a,o}$ and $\Phi(p)\in \Aline v{e'}\cap\Aline op$, so the point $\Phi(p)$ has coordinate $(-e,-a)$. Since $-e\cdot(a+b)+b=-e\cdot a=-a$, the point $\Phi(p)$ belongs to the line $L_{a+b,b}$. If $p\in \Aline ev\cap L_{a,b}$, then $p$ has coordinates $(e,a+b)$ and $\Phi(p)\in \Aline op\cap\Aline hv=\overline L_{a+b,o}\cap\Aline hv$. If $p\in \Aline ov\cap L_{a,b}$, then $p$ has coordinates $(o,b)$. Since $o\cdot(a+b)+b=b$, the point $\Phi(p)=p$ belongs to the line $L_{a+b,b}$. Finally, assume that $p\in L_{a,b}\setminus(\Aline ov\cup\Aline ev)$. In this case, the point has coordinates $(x,x{\cdot} a+b)$ for some $x\in\Delta\setminus\{o,e\}$ and its image $\Phi(p)$ has coordinates 
$$(x/(e{-}x),(x/(e{-}x))\cdot(x\backslash (x{\cdot}a+b))).$$

By Artin's Theorem~\ref{t:Artin}, the alternative division ring $\Delta$ is diassociative and hence every element $z\in \Delta\setminus\{o\}$ has a two-sided inverse $z^{-1}\in\Delta\setminus\{o\}$ such that $z\backslash s=z^{-1}{\cdot} s$ and $s/z=s{\cdot} z^{-1}$ for every element $s\in \Delta$. By the distributivity and diassociativity of the ternar $\Delta$, we obtain 
$$
\begin{aligned}
x/(e{-}x)&=x{\cdot}(e{-}x)^{-1}=((e{-}x){\cdot}x^{-1})^{-1}=((e{\cdot}x^{-1})-(x{\cdot}x^{-1}))^{-1}\\
&=((x^{-1}{\cdot}e)-(x^{-1}{\cdot}x))^{-1}=(x^{-1}\cdot(e{-}x))^{-1}=(e{-}x)^{-1}\cdot x.
\end{aligned}
$$ 
By Artin's Theorem~\ref{t:Artin}, the subring of $\Delta$, generated by the elements $e,x,b$ is associative. Then
$$
\begin{aligned}
(x/(e{-}x))\cdot(x\backslash b)&=\big((e{-}x)^{-1}{\cdot} x\big)\cdot (x^{-1}{\cdot} b)=(e{-}x)^{-1}{\cdot} (x\cdot x^{-1}){\cdot} b=(e{-}x)^{-1}\cdot e\cdot b\\
&=(e{-}x)^{-1}\cdot(e-x+x)\cdot b=(e{-}x)^{-1}{\cdot}(e{-}x){\cdot}b+(e{-}x)^{-1}{\cdot}x{\cdot}b\\
&=b+(x/(e{-}x)){\cdot}b=(x/(e{-}x)){\cdot} b+b
\end{aligned}
$$and 
$$
\begin{aligned}
(x/(e{-}x))\cdot(x\backslash(x{\cdot} a+b))&=(x/(e{-}x))\cdot(a+(x\backslash b))=
(x/(e{-}x)){\cdot} a+(x/(e{-}x)){\cdot} (x\backslash b)\\
&=(x/(e{-}x)){\cdot} a+(x/(e{-}x)){\cdot} b+b=(x/(e{-}x)){\cdot}(a+b)+b,
\end{aligned}
$$
witnessing that the point $\Phi(p)$ belongs to the line $\overline {L_{a+b,b}}$. Therefore, $\Phi[L_{a,b}]=L_{a+b,b}$ is a line in the projective plane $\Pi$, and $\Phi$ is an automorphism of $\Pi$, see Theorem~\ref{t:liner-isomorphism<=>}. Since $\Phi(h)=f$ and $\Phi(v)=v$, the $(h,\Aline hv)$-transitivity of $\Pi$ implies the $(f,\Aline fv)$-transitivity of $\Pi$. Since $\Pi$ is $(h,\Aline hv)$-transitive, there exists an automorphism $F\in\Aut_{h,\overline{h\,v}}(\Pi)$ such that $F(f)=o$. By Proposition~\ref{p:pL-Auto}, the $(f,\Aline fv)$-transitivity of $\Pi$ implies the $(o,\Aline ov)$-transitivity of $\Pi$. Therefore, the projective plane $\Pi$ is $[\Aline hv,\Aline hv)$-transitive, $(v,\Aline ov)$-transitive, and $(o,\Aline ov)$-transitive.
\end{proof}

\section{Moufang projective planes}\label{s:Moufang-proj1}

In this section we characterize Moufang projective planes. By Theorem~\ref{t:Moufang<=>everywhere-Thalesian}, a projective plane $Y$ is Moufang if and only if for every line $H$ in $Y$, the affine subliner $Y\setminus H$ of $Y$ is Thalesian.  By Theorem~\ref{t:paraD<=>translation}, an affine space is Thalesian if and only if it is translation.  In Theorem~\ref{t:Skornyakov-San-Soucie} we shall prove that a projective plane $\Pi$ is Moufang if and only if every line $L\subseteq \Pi$ is \defterm{translation} in the sense that for every distinct points $x,y\in X\setminus H$ there exists an automorphism $A:X\to X$ such that $A(x)=y$ and $\Fix(A)=H$.

The main result of this section is the following deep characterization of Moufang planes, proved by combined efforts of Moufang (who proved the equivalences of the conditions (1), (2), (3), (5), (6)), \index[person]{Skornyakov}Skornyakov\footnote{{\bf Lev Anatolievich Skornyakov} (1924 -- 1989), a Soviet algebraist, defended his Ph.D. ``Alternative skew-fields and alternative planes'' in 1950 under supervision of A.G. Kurosh. Since 1960 he worked at the Chair of Higher Algebra in Lomonosov Moscow State University.} (1951)  and \index[person]{San-Soucie}San-Soucie\footnote{{\bf Robert Louis San Soucie} (1927 -- 2017), an American mathematician, defended his Ph.D. Thesis  ``Right Alternative Division Rings of Characteristic Two'' under supervision of R.H. Bruck 1953. After receiving his PhD, Robert taught Mathematics at the University of Oregon, in Eugene.
\newline 
San Soucie, a graduate of Adams high school with the class of 1944, received the Bachelor of Arts degree in Mathematics, Summa cum Laude. He is a member of Phi Kappa Phi, national honorary society; Sigma Xi, national honorary society of research scientists, and the Phi Beta Kappa association of Massachusetts. He received the Phi Kappa Phi scholarship award, was Phi Beta scholar of the senior class; was a member of Kappa Sigma fraternity, Mathematics and Newman clubs, summer session social committee, German club, chairman of the Connecticut Valley student science conference, solicitation chairman for the campus community chest and chairman of the inter-fraternity council constitution committee. He has accepted an appointment to the research faculty of the department of mathematics at the University of Wisconsin where he will do research and graduate work for his doctor's degree. The Wisconsin fellowship award is for \$1,100. He also won a \$1,000 fellowship to McGill university at Montreal, which he relinquished because of the Wisconsin award. He served in the U. S. Navy during World War II.} (1955) who proved (independently) the equivalence of the condition (4) to other conditions of Theorem~\ref{t:Skornyakov-San-Soucie}.

\begin{theorem}[Moufang, 1933; Skornyakov, San-Soucie; 1955]\label{t:Skornyakov-San-Soucie} For a projective plane $\Pi$, the following conditions are equivalent:
\begin{enumerate}
\item $\Pi$ is Moufang;
\item $\Pi$ is everywhere Thalesian;
\item every line in $\Pi$ is translation;
\item at least two lines in $\Pi$ are translation;
\item some ternar of $\Pi$ is linear, distributive, associative-plus, and alternative-dot;
\item every ternar of $\Pi$ is linear, distributive, associative-plus, and alternative-dot.
\end{enumerate}
\end{theorem}

\begin{proof} The equivalence $(1)\Leftrightarrow(2)$ has been proved in Theorem~\ref{t:Moufang<=>everywhere-Thalesian}.
\smallskip

$(2)\Ra(3)$ Assume that the projective plane $\Pi$ is everywhere Thalesian and take any line $L$ in $\Pi$. Given any distinct points $x,y\in\Pi\setminus L$, we have to find an automorphism $A:\Pi\to\Pi$ such that $A(x)=y$ and $\Fix(A)=L$. Since $\Pi$ is a projective plane, there exists a unique point $p\in L\cap\Aline xy$. If the point-line pair $(L,p)$ is Desarguesian, then an automorphism $A$ of $\Pi$ with $A(x)=y$ and $\Fix(A)=L$ exists, by the definition of a Desarguesian point-line pair. If the pair $(L,p)$ is not Desarguesian, then the projective plane $\Pi$ is not Desarguesian and hence $\Pi$ is $4$-long, by Proposition~\ref{p:Steiner+projective=>Desargues}. Then the affine subliner $X\defeq\Pi\setminus L$ is a $3$-long Playfair plane, by Propositions~\ref{p:projective-minus-flat}, \ref{p:projective-minus-hyperplane}, and Theorem~\ref{t:Playfair<=>}. Since $\Pi$ is everywhere Thalesian, the Playfair plane $X$ is Thalesian and by Theorem~\ref{t:paraD<=>translation}, $X$ is translation. Then there exists a translation $T:X\to X$ such that $T(x)=y$. By Theorem~\ref{t:extend-isomorphism-to-completions}, the translation $T$ can be extended to an automorphism $A$ of the projective plane $\Pi$. Since the translation $T$ has no fixed points, $\Fix(A)\subseteq L$. To prove that $L=\Fix(A)$, take any point $z\in L$ and choose any  line $\Lambda$ in $\Pi$ such that $\Lambda\cap L=\{z\}$. Since $T$ is a translation of the affine plane $X$, the line $T[\Lambda\cap X]$ is parallel to the line $\Lambda\cap X$ in the affine plane $X$. If $T[\Lambda\cap X]=\Lambda\cap X$, then $A[\Lambda]=\Lambda$ and $A(z)\in A[\Lambda\cap L]=A[\Lambda]\cap A[L]=\Lambda\cap L=\{z\}$ and hence $z\in \Fix(A)$. If $T[\Lambda\cap X]\ne\Lambda\cap X$, then 
$T[\Lambda\cap X]\cap (\Lambda\cap X)=\varnothing$. Since $\Pi$ is a projective plane, the lines $A[\Lambda]$ and $\Lambda$ have unique common point $o$. Assuming that $o\in X$, we conclude that $$o\in X\cap A[X]\cap A[\Lambda]\cap \Lambda=A[\Lambda\cap X]\cap(\Lambda\cap X)=T[\Lambda\cap X]\cap(\Lambda\cap X)=\varnothing,$$ which is a contradiction showing that 
$o\in \Pi\setminus X=L=A[L]$ and hence $o\in (L\cap \Lambda)\cap (A[L]\cap A[\Lambda]=\{z\}\cap A[L\cap\Lambda]=\{z\}\cap \{A(z)\}$ and $A(z)=o=z\in \Fix(A)$.
\smallskip

$(3)\Ra(2)$ Assume that every line in $\Pi$ is translation. To prove that the projective plane $\Pi$ is everywhere Thalesian, we should prove that for every line $L\subset\Pi$, the affine liner $X\defeq \Pi\setminus L$ is Thalesian. If $|\Pi|_2=3$, then $|X|_2=2$ and hence the liner $|X|_2$ contains no disjoint parallel lines and hence is Thalesian vacuously. So, assume that $|\Pi|_2\ge 4$. In this case the affine liner $X=\Pi\setminus L$ is regular and Playfair, by Proposition~\ref{p:projective-minus-hyperplane} and Theorem~\ref{t:Playfair<=>}. Since the line $L$ is translation, the Playfair liner $X$ is translation and Thalesian, by Theorem~\ref{t:paraD<=>translation}. 
\smallskip

The implication $(3)\Ra(4)$ is trivial. 
\smallskip

$(4)\Ra(5)$ Assume that the projective plane $\Pi$ contains at least two distinct translation lines. In this case we can choose a projective base $uowe$ for $\Pi$ so that the lines $\Aline hv$ and $\Aline ov$ are translation. Then the projective plane $\Pi$ is $[\Aline hv,\Aline hv)$-transitive, $(v,\Aline ov)$-transitive and $(o,\Aline ov)$-transitive. By Theorem~\ref{t:Algebra-vs-Geometry-proj}(15), the ternar $\Delta$ of the projective base $uowe$ is linear, distributive, associative-plus and alternative-dot.
\smallskip

$(5)\Ra(3)$ Assume that some  ternar of $\Pi$ is linear, distributive, associative-plus, and alternative-dot. Then there exists a projective base $uowe$ for $\Pi$ whose ternar is  linear, distributive, associative-plus, and alternative-dot. By Theorem~\ref{t:Algebra-vs-Geometry-proj}(15), the projective plane $\Pi$ is $[\Aline hv,\Aline hv)$-transitive and $(o,\Aline ov)$-transitive, and by Theorem~\ref{t:Algebra-vs-Geometry-proj}(14), the affine plane $\Pi\setminus\Aline hv$ is shear. Now we can prove that every line $L$ in the projective plane $\Pi$ is translation. If $L=\Aline hv$, then $L$ is translation because the line $\Aline hv$ is translation. So, assume that $L\ne \Aline hv$. In this case $L\setminus\Aline hv$ is a line in the affine plane $\Pi\setminus\Aline hv$. By Proposition~\ref{p:shear=>2-homogen}, there exists an automorphism $A$ of the affine plane $\Pi\setminus\Aline hv$ such that $A[\Aline ov\setminus\Aline hv]=L\setminus\Aline hv$. By Theorem~\ref{t:extend-isomorphism-to-completions}, the automorphism $A$ can be extended to an automorphism $\bar A$ of the projective plane $\Pi$. It follows from $A[\Aline ov\setminus\Aline hv]=L\setminus\Aline hv$ that $\bar A[\Aline ov]=L$. Since the line $\Aline ov$ is translation, so is its image $L=\bar A[\Aline ov]$ under the automorphism $\bar A$. 
\smallskip

The implication $(6)\Ra(5)$ is trivial.
\smallskip

$(3)\Ra(6)$ Assume that every line in the projective plane $\Pi$ is translation. Given any ternar $R$ of the projective plane $\Pi$, find a projective base $uowe$ of $\Pi$ whose ternar $\Delta$ is isomorphic to the ternar $R$.  Let $h\in \Aline ou\cap\Aline we$ and $v\in \Aline ow\cap \Aline ue$ be the horizontal and vertical infinity points of the projective base $uowe$. Since every line in $\Pi$ is translation, the lines $\Aline hv$ and $\Aline ov$ are translation. Then the projective plane $\Pi$ is $[\Aline hv,\Aline hv)$-transitive and $[\Aline ov,\Aline ov)$-transitive. Applying Theorem~\ref{t:Algebra-vs-Geometry-proj}(15), we conclude that the ternar $\Delta$ of the based projective plane $(\Pi,uowe)$ is linear, distributive, associative-plus and alternative-dot, and so is the isomorphic copy $R$ of $\Delta$. 
\end{proof}

\begin{remark} Theorem~\ref{t:Algebra-vs-Geometry-proj} does not contain an algebraic characterization of based projective planes which are $(v,\Aline ov)$-transitive and $(h,\Aline oh)$-transitive.
By Lemma~\ref{l:Lshear}, such based projective planes are $(p,\Aline op)$-transitive for every point $p\in\Aline hv$.  The ternars of such based projective planes are not linear (in general) but are at least vertical-shear, see Definition~\ref{d:tring-ldistributive}. Such based projective planes are dual to $[\Aline hv,\Aline hv)$-transitive projective planes (which are coordinatized by linear right-distributive associative-plus ternars).
\end{remark}

\begin{problem} Find an algebraic characterization of based projective planes which are $(v,\Aline ov)$-transitive and $(h,\Aline oh)$-transitive. This problem can be reformulated as the problem of finding an algebraic characterization of vertical-shear and horizontal-shear based affine planes.
\end{problem}

\section{The Lenz classification of projective planes}\label{s:Lenz}

The \index{Lenz--Barlotti figure}\index{projective plane!Lenz--Barlotti figure}\defterm{Lenz--Barlotti figure} $\LB_\Pi$ of a projective plane $\Pi$ is the set of all Desarguesian point-line pairs $(p,L)\in \Pi\times\mathcal L_\Pi$.  The \index{Lenz figure}\index{projective plane!Lenz figure}\defterm{Lenz figure} $\Lenz_\Pi$ of a projective plane $\Pi$ is the subset of the Lenz--Barlotti figure of $\Pi$, consisting of all Desarguesian point-line pairs $(p,L)$ with $p\in L$. By Baer's Theorem~\ref{t:Baer-pL-Des<=>}, a point-line pair $(p,L)$ belongs to the Lenz--Barlotti figure if and only if the projective plane $\Pi$ is $(p,L)$-transitive.

Lenz figures of projective planes are classified by the following theorem due to \index[person]{Lenz}Lenz\footnote{{\bf Hanfried Lenz}  (1916 -- 2013) was a German mathematician, who is mainly known for his work in geometry and combinatorics.
\newline 
Hanfried Lenz was the eldest son of Fritz Lenz an influential German geneticist, who is associated with Eugenics and hence also with the Nazi racial policies during the Third Reich. He was also the older brother of Widukind Lenz, a geneticist. He started to study mathematics and physics at the University of T\"ubingen, but interrupted his studies from 1935 to 1937 to do a (at this time, in Weimar Republic voluntary) military service. After that he continued to study in Munich, Berlin and Leipzig. In 1939 when World War II broke out in Europe, he became a soldier in the western front and during a vacation he passed the exams for his teacher certification. He married Helene Ranke in 1943 and 1943--45 he worked on radar technology in a laboratory near Berlin.
\newline 
After World War II Hanfried Lenz was classified as a ``follower'' by the denazification process. He started to work as a math and physics teacher in Munich and in 1949 he became an assistant at the Technical University of Munich. He received his PhD in 1951 and his Habilitation in 1953. He worked as a lecturer until he became an associate professor in 1959. In 1969 he finally became a full professor at the Free University of Berlin and worked there until his retirement in 1984.
\newline
He was also politically active and in connection with his opposition to the rebuilding of the German army in the early 50s, he became a member of the Social Democratic Party (SPD) in 1954. Later, partially due to being alienated by the student movement of the '60s, his leanings became more conservative again and in 1972 he left the SPD to join the Christian Democratic Union.
\newline
Hanfried Lenz is known for his work on the classification of projective planes and in 1954 he showed how one can introduce affine spaces axiomatically without constructing them from projective spaces or vector spaces. This result is now known as the theorem of Lenz. During his later years he also worked in the area of combinatorics and published a book on design theory (together with Dieter Jungnickel and Thomas Beth).
\newline
In 1995 the Institute of Combinatorics and its Applications awarded the Euler Medal to Hanfried Lenz.} \cite{Lenz1954}.

\begin{theorem}[Lenz, 1954]\label{t:Lenz} The Lenz figure $\Lenz_\Pi$ of any projective plane $\Pi$ is equal to one of the following sets:
\begin{itemize}
\item[{$(\,\circ\,)$}]\index[note]{$(\,\circ\,)$} the empty set $\varnothing$;
\item[{$(\,\bullet\,)$}]\index[note]{$(\,\bullet\,)$} $\{(p,L)\}$ for some point $p\in\Pi$ and line $L\in\mathcal L_p$;
\item[{$(\;\vert\;)$}]\index[note]{$(\,\vert\,)$} $\{p\}\times\mathcal L_p$ for some point $p\in\Pi$;
\item[{$({-})$}]\index[note]{$({-})$} $L\times\{L\}$ for some line $L\in\mathcal L_\Pi$;
\item[{$({+})$}]\index[note]{$({+})$} $(L\times\{L\})\cup(\{p\}\times\mathcal L_p)$ for some line $L\in\mathcal L_\Pi$ and point $p\in L$;
\item[{$(\circledast)$}]\index[note]{$(\circledast)$} $\{(x,\Aline xp):x\in L\}$ for some line $L\in\mathcal L_\Pi$ and point $p\in \Pi\setminus L$;
\item[{$(\square)$}]\index[note]{$(\square)$} $\{(p,L)\in\Pi\times\mathcal L_\Pi:p\in L\}$.
\end{itemize}
\end{theorem}

\begin{proof} Consider the projections $$\dom[\Lenz_\Pi]\defeq\{p\in \Pi:\exists L\in\mathcal L_\Pi \;(p,L)\in\Lenz_\Pi\}\quad\mbox{and}\quad\rng[\Lenz_\Pi]\defeq\{L\in\mathcal L_\Pi:\exists p\in\Pi \;(p,L)\in\Lenz_\Pi\}$$ of the Lenz figure $\Lenz_\Pi\subset\Pi\times\mathcal L_\Pi$ onto the factors. Depending on the rank of the set $\dom[\Lenz_\Pi]$ (and also the set $\rng[\Lenz_\Pi]$) in the projective plane $\Pi$ (and the dual projective plane $\mathcal L_\Pi$), we consider four cases (with subcases).
\smallskip

0. If $\|\dom[\Lenz_\Pi]\|=0$, then $\Lenz_\Pi=\varnothing$ and we have the case $(\,\circ\,)$.
\smallskip

1. If $\|\dom[\Lenz_\Pi]=1\|$, then $\dom[\Lenz_\Pi]=\{p\}$ for some point $p\in\Pi$ and hence $\Lenz_\Pi\subseteq \{p\}\times\mathcal L_p$. If $\Lenz_\Pi=\{(p,L)\}$ for some point-line pair $(p,L)$, then we have the case $(\,\bullet\,)$. If $\Lenz_\Pi$ contains two distinct point-line pairs $(p,L)$, $(p,L')$, then by Corollary~\ref{c:(p,q]-transitive}, $\Lenz_\Pi=\{p\}\times\mathcal L_p$ and we have the case $(\,|\,)$.
\smallskip

21. Next, assume that $\|\dom[\Lenz_\Pi]\|=2$ and $\|\rng[\Lenz_\Pi]\|=1$. In this case  $\dom[\Lenz_\Pi]\subseteq \Aline pq$ for some distinct points $p,q\in \dom[\Lenz_\Pi]$ and $\rng[\Lenz_\Pi]=\{L\}$ for some line $L$. Applying Theorem~\ref{t:(pq,L)-transitive}, we conclude that $\Lenz_\Pi=\Aline pq\times \{L\}=L\times\{L\}$, which means that the case $(-)$ holds.
\smallskip

22. Assume that $\|\dom[\Lenz_\Pi]\|=2=\|\rng[\Lenz_\Pi]\|$. This case has three subcases. 
\smallskip

22a. First assume that there exists a point $p\in\Pi$ and two distinct lines $L_1,L_2$ such that $(p,L_1),(p,L_2)\in\Lenz_\Pi$. In this case we can apply Corollary~\ref{c:(p,q]-transitive} and conclude that $\{p\}\times\mathcal L_p\subseteq \Lenz_\Pi$ and hence $\rng[\Lenz_\Pi]=\mathcal L_p$ (because $\|\rng[\Lenz_\Pi]\|=2$). Since $\|\dom[\Lenz_\Pi]\|=2$, there exists a point $q\in \dom[\Lenz_\Pi]$ such that $q\ne p$. Find a line $L\in\mathcal L_q$ such that $(q,L)\in \Lenz_\Pi$. Since $L\in\rng[\Lenz_\Pi]=\mathcal L_p$, the line $L\in \mathcal L_p\cap\mathcal L_q$ is equal to the line $\Aline pq$. Then $\{(p,\Aline pq),(q,\Aline pq)\}\subseteq\Lenz_\Pi$ and hence $\Aline pq\times\{\Aline pq\}\subseteq \Lenz_\Pi$, by Theorem~\ref{t:(pq,L)-transitive}. It follows from $\|\dom[\Lenz_\Pi]\|=2$ that $\dom[\Lenz_\Pi]=\Aline pq$. Then for every $(c,L)\in \Lenz_\Pi\setminus(\{p\}\times\mathcal L_p)$, we have $c\in \dom[\Lenz_\Pi]=\Aline pq$ and $L\in \mathcal L_c\cap\mathcal L_p$, which implies $L=\Aline cp=\Aline pq$. Therefore, $$\Lenz_\Pi=(\Aline pq\times\{\Aline pq\})\cup(\{p\}\times\mathcal L_p)$$ and we have the case $(+)$.
\smallskip

22b. Next, assume that there exists a line $L\in \mathcal L_\Pi$ and two distinct points $p',p''\in\Pi$ such that $(p',L),(p'',L)\in\Lenz_\Pi$. By Theorem~\ref{t:(pq,L)-transitive}, $\Aline {p'}{p''}\times\{L\}\subseteq\Lenz_\Pi$ and hence $L=\Aline {p'}{p''}$. It follows from $\|\dom[\Lenz_\Pi]\|=2$ that $\dom[\Lenz_\Pi]=\Aline {p'}{p''}=L$. Since $\|\rng[\Lenz_\Pi]\|=2$, there exists a line $\Lambda\in\rng[\Lenz_\Pi]\setminus\{L\}$. Find a point $q\in\Pi$ such that $(q,\Lambda)\in\Lenz_\Pi$ and observe that $q\in\dom[\Lenz_\Pi]=L$. Since $L,\Lambda\in \mathcal L_q$, we have the case $22a$, which implies that $\Lenz_\Pi=(L\times\{L\})\cup(\{q\}\times\mathcal L_q)$ and hence the case $(+)$ holds.
\smallskip

22c. The cases 22a and 22b do not hold. Choose any distinct point-line pairs $(a,A),(b,B)\in\Lenz_\Pi$. Since the cases $22a$ and $22b$ do not hold, $a\ne b$ and $A\ne B$. Let $p$ be the unique common point of the lines $A$ and $B$. Since $\|\dom[\Lenz_\Pi]\|=2=\|\rng[\Lenz_\Pi]\|$, $\dom[\Lenz_\Pi]\subseteq \Aline ab$ and $\rng[\Lenz_\Pi]\subseteq \mathcal L_p$. 

If $a=p$, then by the $(b,B)$-transitivity, for every line $C\in \mathcal L_p\setminus\{A,B\}$, there exists an automorphism $F\in\Aut_{b,B}(\Pi)$ such that $F[A]=C$. Since $F(a)=F(p)=p$, the $(a,A)$-transitivity of $\Pi$ implies that $\{(p,C):C\in\mathcal L_p\setminus\{B\}\}\subseteq \Lenz_\Pi$, which means that the case 22a holds. But this contradicts our assumption. This contradiction shows that $a\ne p$. By analogy we can prove that $b\ne p$.
Assuming that $p\in\Aline ab$ and taking into account that $a\in A\in\mathcal L_p$ and $b\in B\in\mathcal L_p$, we conclude that $A=\Aline ap=\Aline pb=B$ and hence the case 22b holds. But this contradicts out assumption.

This contradiction shows that $p\notin \Aline ab$. Then $(a,A)=(a,\Aline ap)$ and $(b,B)=(b,\Aline bp)$. Applying Theorem~\ref{t:circledast}, we conclude that $\{(x,\Aline xp):x\in \Aline ab\}\subseteq \Lenz_\Pi$ and hence $\Aline ab=\dom[\Lenz_\Pi]$ and $\mathcal L_p=\rng[\Lenz_\Pi]$. Taking into account that the cases 22a and 22b do not hold, we conclude that $\Lenz_\Pi=\{(x,\Aline xp):x\in\Aline ab\}$, which means that the case $(\circledast)$ holds.
\smallskip

3. Finally, assume that $\|\dom[\Lenz_\Pi]\|=3$. In this case we shall apply the following claims to prove that the Lenz figure $\Lenz_\Pi$ is equal to the set $\square\defeq\{(p,L)\in\Pi\times\mathcal L_\Pi:p\in L\}$ and hence the case $(\square)$ holds. 

\begin{claim}\label{cl:Lenz1} If the projective plane $\Pi$ contains a translation line, then every line in $\Pi$ is translation and $\Lenz_\Pi=\square$.
\end{claim}

\begin{proof} Assume that $\Lambda$ is a translation line in $\Pi$. If $|\Pi|_2=3$, then $\Pi$ is a Steiner projective plane and every line in $\Pi$ is translation, by Theorem~\ref{t:Skornyakov-San-Soucie}. 
So, assume that $|\Pi|_2\ge 4$. Since $\|\dom[\Lenz_\Pi]\|=3>2=\|\Lambda\|$, there exists a point-line pair $(p,L)\in\Lenz_\Pi$ such that $p\notin \Lambda$. Since $p\in L\setminus\Lambda$, there exists a point $\lambda\in \Lambda\setminus L$. Since $|\Pi|_2\ge 4$, there exists a point $\lambda'\in \Aline p\lambda\setminus(\{p\}\cup L\cup\Lambda)$. Since the plane $\Pi$ is $(p,L)$-transitive, there exists an automorphism $A\in\Aut_{p,L}(\Pi)$ such that $\lambda=\lambda'$. Then $\Lambda'\defeq A[\Lambda]$ is a translation line, distinct from the translation line $\Lambda$. By Theorem~\ref{t:Skornyakov-San-Soucie}, every line in the projective plane $\Pi$ is translation and hence $\Lenz_\Pi=\square$.
\end{proof}

\begin{claim}\label{cl:Lenz2} If there exist distinct points $a,b\in\dom[\Lenz_\Pi]$ such that $(a,\Aline ab)\in \Lenz_\Pi$, then $\Pi$ contains a translation line and hence $\Lenz_\Pi=\square$.
\end{claim}

\begin{proof} Since $b\in\dom[\Lenz_\Pi]$, there exists a line $B\in\mathcal L_\Pi$ such that $(b,B)\in\Lenz_\Pi$. If $B=\Aline ab$, then the plane $\Pi$ is $(a,\Aline ab)$-transitive and $(b,\Aline ab)$-transitive. By Theorem~\ref{t:(pq,L)-transitive}, the plane $\Pi$ is $[\Aline ba,\Aline ab)$-transitive, which means that the line $\Aline ab$ is translation.

So, assume that $B\ne \Aline ab$ and hence $a\notin B$. Choose any point $c\in \Aline ab\setminus\{a,b\}$. Since $\Pi$ is $(b,B)$-transitive, there exists an automorphism $A\in\Aut_{b,B}(\Pi)$ such that $A(a)=c$. Since $A[\Aline ab]=\Aline cb=\Aline ab$, the $(a,\Aline ab)$-transitivity of $\Pi$ implies the $(c,\Aline ab)$-transitivity of $\Pi$. By Theorem~\ref{t:(pq,L)-transitive}, the plane $\Pi$ is $[\Aline ac,\Aline ab)$-transitive, which implies that the line $\Aline ab=\Aline ac$ is translation. Now we can apply Claim~\ref{cl:Lenz1} and conclude that $\Lenz_\Pi=\square$.
\end{proof}

Now we are able to finish the analysis of the case 3. Take any point-line pair $(a,A)\in\Lenz_\Pi$. Since $\|\dom[\Lenz_\Pi]\|=3$, there exists a point-line pair $(b,B)\in \Lenz_\Pi$ such that $b\notin A$. If $a\in B$, then $B=\Aline ab$ and by Claim~\ref{cl:Lenz2}, $\Lenz_\Pi=\square$.

So, assume that $a\notin B$. Let $o$ be the unique common point of the lines $A$ and $B$. The assumption $a\notin B$ implies $o\notin\Aline ab$ (otherwise $a\in A=\Aline ao=\Aline pb=B$). Since $(a,A)=(a,\Aline ao)$ and $(b,B)=(b,\Aline bo)$, we can apply Theorem~\ref{t:circledast} and conclude that $\{(x,\Aline xo):x\in \Aline ab\}\subseteq\Lenz_\Pi$. Since $\|\dom[\Lenz_\Pi]\|=3>2=\|\Aline ab\|$, there exists a point-line pair $(c,C)\in\Lenz_\Pi$ such that $c\notin\Aline ab$. Since $\Pi=\bigcup_{x\in \overline{a\,b}}\Aline xo$, there exists a point $x\in \Aline ab$ such that $c\in \Aline xo$. Since $(x,\Aline xc)=(x,\Aline xo)\in\Lenz_\Pi$, we can apply Claim~\ref{cl:Lenz2} and conclude that $\Lenz_\Pi=\square$. 
\end{proof}

\begin{remark}\label{r:Lenz} The cases $(\circ),(\bullet),(\circledast),(-),(\,|\,),(+),(\square)$ in the Lenz classification are usually numbered by the roman numerals I, II, III, IVa, IVb, V, VII, respectively. The number VI is missed because in the initial classification of Lenz, it denoted two dual classes VI.a and VI.b of projective planes, which were later excluded by  Theorem~\ref{t:Skornyakov-San-Soucie} of Skornyakov and San-Soucie. Projective planes of Lenz class $(-)$ are called \defterm{translation}. By Theorem~\ref{t:paraD<=>translation}, a projective plane is translation if and only if it is somewere Thalesian. By a result of Hering and Kantor \cite{HK1971}, projective planes of Lenz type $(\circledast)$ are infinite.
\end{remark}

\begin{definition}\label{d:po-Lenz} A projective plane $\Pi$ has the Lenz type at least
\begin{itemize}
\item[{$(\,\circ\,)$}]\index[note]{$(\,\circ\,)$} if $\varnothing\subseteq\Lenz_\Pi$;
\item[{$(\,\bullet\,)$}]\index[note]{$(\,\bullet\,)$} if $\Lenz_\Pi\ne\varnothing$;
\item[{$(\;\vert\;)$}]\index[note]{$(\,\vert\,)$} if $\{p\}\times\mathcal L_p\subseteq\Lenz_\Pi$ for some point $p\in\Pi$;
\item[{$({-})$}]\index[note]{$({-})$} if $L\times\{L\}\subseteq\Lenz_\Pi$ for some line $L\in\mathcal L_\Pi$;
\item[{$({+})$}]\index[note]{$({+})$} if $(L\times\{L\})\cup(\{p\}\times\mathcal L_p)\subseteq\Lenz_\Pi$ for some line $L\in\mathcal L_\Pi$ and point $p\in L$;
\item[{$(\circledast)$}]\index[note]{$(\circledast)$} if $\{(x,\Aline xp):x\in L\}\subseteq\Lenz_\Pi$ for some line $L\in\mathcal L_\Pi$ and point $p\in \Pi\setminus L$;
\item[{$(\square)$}]\index[note]{$(\square)$} if $\{(p,L)\in\Pi\times\mathcal L_\Pi:p\in L\}=\Lenz_\Pi$.
\end{itemize}
\end{definition}

Definition~\ref{d:po-Lenz} actually involves a natural partial order  
on the set of Lenz types $$\{(\circ), (\bullet), (-), (\,|\,), (+), (\circledast), (\square)\}.$$
The Hasse diagram of this partial order looks as follows.
$$
\xymatrix@C=8pt@R=12pt{
&(\square)\\
&(+)\ar@{-}[u]\\
(\,|\,)\ar@{-}[ur]&(-)\ar@{-}[u]&(\circledast)\ar@{-}[uul]\\
&(\bullet)\ar@{-}[u]\ar@{-}[ul]\ar@{-}[ur]\\
&(\circ)\ar@{-}[u]
}
$$

\begin{theorem}\label{t:Lenz-algebra} A projective plane $\Pi$ has  Lenz type at least
\begin{enumerate}
\item[$(\bullet)$] iff some ternar of $\Pi$ is linear and associative-plus;
\item[$(\,|\,)$] iff some ternar of $\Pi$ is linear, left-distributive, and associative-plus;
\item[$(-)$] iff some ternar of $\Pi$ is linear, right-distributive, and associative-plus;
\item[$(+)$] iff some ternar of $\Pi$ is iff linear, distributive, and associative-plus;
\item[$(\square)$] iff some ternar of $\Pi$ is linear, distributive, associative-plus, and alternative-dot.
\end{enumerate}
\end{theorem}

\begin{proof} 1. If $\Pi$ has Lenz type at least $(\bullet)$, then $(p,L)\in \Lenz_\Pi$ for some point $p\in \Pi$ and line $L\in\mathcal L_p$. In this case we can choose a projective base $uowe$ in $\Pi$ such that $p=v$ and $L=\Aline hv$, where $h\in\Aline ou\cap\Aline we$ and $v\in\Aline ow\cap\Aline ue$ are the horizontal and vertical infinity points of the projective base $uowe$. The $(p,L)$-transitivity of $\Pi$ ensures that the based projective plane $(\Pi,uowe)$ is $(v,\Aline hv)$-transitivive. By Theorem~\ref{t:Algebra-vs-Geometry-proj}(1), the ternar $\Delta=\Aline oe\setminus\Aline hv$ of the projective base $uowe$ is linear and associative-plus. 

Now assume that some ternar of $\Pi$ is linear and associative-plus. This means that $\Pi$ has a projective base $uowe$ whose ternar is linear and associative-plus. By Theorem~\ref{t:Algebra-vs-Geometry-proj}(1), the based projective plane $(\Pi,uowe)$ is $(v,\Aline hv)$-transitive. Then $(v,\Aline hv)\in\Lenz_\Pi$ and the projective plane $\Pi$ has Lenz type at least $(\bullet)$.
\smallskip

2. If $\Pi$ has Lenz type at least $(\,|\,)$, then $\{p\}\times\mathcal L_p\subseteq \Lenz_\Pi$ for some point $p\in \Pi$. Choose a projective base $uowe$ in $\Pi$  whose vertical infinity point $v\in\Aline ow\cap\Aline ue$ coincides with $p$. The inclusion $\{p\}\times\mathcal L_p\subseteq \Lenz_\Pi$ implies that the projective plane $\Pi$ is $(v,v]$-transitive. By Theorem~\ref{t:Algebra-vs-Geometry-proj}(3), the ternar of the projective base $uowe$ is linear, left-distributive, and associative-plus. 

Now assume that some ternar of $\Pi$ is linear, left-distributive, and associative-plus. This means that $\Pi$ has a projective base $uowe$ whose ternar is linear, left-distributive and associative-plus. By Theorem~\ref{t:Algebra-vs-Geometry-proj}(3), the based projective plane $(\Pi,uowe)$ is $(v,v]$-transitive. Then $\{v\}\times\mathcal L_v\subseteq\Lenz_\Pi$ and the projective plane $\Pi$ has Lenz type at least $(\,|\,)$.
\smallskip 

3. If $\Pi$ has Lenz type at least $(-)$, then $L\times\{L\}\subseteq \Lenz_\Pi$ for some line $L\in\mathcal L_\Pi$. Choose a projective base $uowe$ in $\Pi$  such that $L=\Aline hv$, where $h\in\Aline ou\cap\Aline we$ and $v\in\Aline ow\cap\Aline ue$ are the horizontal and vertical infinity points of the projective base $uowe$. The inclusion $L\times\{L\}\subseteq\Lenz_\Pi$ implies that the projective plane $\Pi$ is $[\Aline hv,\Aline hv)$-transitive. By Theorem~\ref{t:Algebra-vs-Geometry-proj}(2), the ternar of the projective base $uowe$ is linear, right-distributive, and associative-plus. 

Now assume that some ternar of $\Pi$ is linear, right-distributive, and associative-plus. This means that $\Pi$ has a projective base $uowe$ whose ternar is linear, right-distributive and associative-plus. By Theorem~\ref{t:Algebra-vs-Geometry-proj}(2), the based projective plane $(\Pi,uowe)$ is $[\Aline hv,\Aline hv)$-transitive. Then $\Aline hv\times\{\Aline hv\}\subseteq\Lenz_\Pi$ and the projective plane $\Pi$ has Lenz type at least $(-)$.
\smallskip

4. If $\Pi$ has Lenz type at least $(+)$, then $(L\times\{L\})\cup(\{p\}\times\mathcal L_p)\subseteq \Lenz_\Pi$ for some point $p\in\Pi$ and line $L\in\mathcal L_p$. Choose a projective base $uowe$ in $\Pi$  such that $L=\Aline hv$ and $p=v$, where $h\in\Aline ou\cap\Aline we$ and $v\in\Aline ow\cap\Aline ue$ are the horizontal and vertical infinity points of the projective base $uowe$. The inclusion $(L\times\{L\})\cup(\{p\}\times\mathcal L_p)\subseteq\Lenz_\Pi$ implies that the projective plane $\Pi$ is $[\Aline hv,\Aline hv)$-transitive and $(v,v]$-transitive. By Theorem~\ref{t:Algebra-vs-Geometry-proj}(4), the ternar $\Delta=\Aline oe\setminus\Aline hv$ of the projective base $uowe$ is linear, distributive, and associative-plus. 

Now assume that some ternar of $\Pi$ is linear, distributive, and associative-plus. This means that $\Pi$ has a projective base $uowe$ whose ternar is linear, distributive and associative-plus. By Theorem~\ref{t:Algebra-vs-Geometry-proj}(4), the based projective plane $(\Pi,uowe)$ is $[\Aline hv,\Aline hv)$-transitive and $(v,v]$-transitive. Then $(\Aline hv\times\{\Aline hv\})\cup(\{v\}\times\mathcal L_v)\subseteq\Lenz_\Pi$ and the projective plane $\Pi$ has Lenz type at least $(+)$.
\smallskip

5. If $\Pi$ has Lenz type at least $(\square)$, then the projective plane $\Pi$ is $(p,L)$-transitive for every point $p\in\Pi$ and line $L\in\mathcal L_p$. By Theorem~\ref{t:Algebra-vs-Geometry-proj}(14) or Theorem~\ref{t:Algebra-vs-Geometry-proj}(15), the ternar of any projective base in $\Pi$ is linear, distributive, associative-plus, and alternative-dot. 

Now assume that some ternar of $\Pi$ is linear, distributive, associative-plus, and alternative-dot. This means that $\Pi$ has a projective base $uowe$ whose ternar is linear, distributive, associative-plus, and alternative-dot. By Theorem~\ref{t:Algebra-vs-Geometry-proj}(15), the based projective plane $(\Pi,uowe)$ is $[\Aline hv,\Aline hv)$-transitive, $(v,\Aline ov)$-transitive, and $(o,\Aline ov)$-transitive. By Theorem~\ref{t:(pq,L)-transitive}, the $(o,\Aline ov)$-transitivity of the plane $\Pi$ implies that $\Pi$ is $[\Aline ov,\Aline ov)$-transitive. Therefore, the projective plane $\Pi$ contains two translation lines $\Aline hv$ and $\Aline ov$. By Theorem~\ref{t:Skornyakov-San-Soucie}, every line in $\Pi$ is translation, which implies that $\Lenz_\Pi=\{(p,L)\in\Pi\times\mathcal L_\Pi:p\in L\}$ and $\Pi$ has the Lenz type (at least) $(\square)$.
\end{proof}

\begin{problem} Find an algebraic characterization of projective planes of Lenz type at least $(\circledast)$.
\end{problem}

\begin{remark} The Lenz figures of projective planes can be also classified by the ranks of their domain and range. More precisely, for a projective plane $\Pi$, the \index{birank}\index{Lenz figure!birank of}\defterm{birank} of its Lenz figure $\Lenz_\Pi$ is the pair $dr$ of the ranks $d\defeq\|\dom[\Lenz_\Pi]\|$ and $r\defeq\|\rng[\Lenz_\Pi]\|$ of the projections 
$$\dom[\Lenz_\Pi]\defeq\{p\in\Pi:\exists L\in\mathcal L_\Pi\;(p,L)\in\Lenz_\Pi\}\quad\mbox{and}\quad \rng[\Lenz_\Pi]\defeq\{L\in\mathcal L_\Pi:\exists p\in\Pi\;(p,L)\in\Lenz_\Pi\}$$of the Lenz figure $\Lenz_\Pi$ onto the projective plane $\Pi$ and its dual projective plane $\mathcal L_\Pi$. 

The Lenz figures of projective planes of Lenz type $$(\circ),(\bullet),(\,|\,),(-),(+),(\circledast),(\square)$$ have biranks
$$00,\;\,11,\;\,12,\;\,21,\;\,22,\;\,22,\;\,33,$$
respectively. 
\end{remark}
 
\section{The Lenz--Barlotti classification of projective planes}

Let us recall that for a projective plane $\Pi$, its Lenz--Barlotti figure is the set $\LB_\Pi$ of all Desarguesian point-line pairs $(p,L)\in \Pi\times\mathcal L_\Pi$. Equivalently, $\LB_\Pi$ is the set of all point-line pairs $(p,L)\in \Pi\times\mathcal L_\Pi$ for which the plane $\Pi$ is $(p,L)$-transitive. The Lenz figure $\Lenz_\Pi$ of a projective plane coincides with the subset $\{(p,L)\in\LB_\Pi:p\in L\}$ of the Lenz--Barlotti figure $\LB_\Pi$.

The Lenz--Barlotti figures were classified by \index[person]{Barlotti}Adriano Barlotti\footnote{{\bf Adriano Barlotti} (1923 -- 2008) was an Italian mathematician who worked on geometry and combinatorics.
\newline
Barlotti studied at the University of Florence and then turned to finite geometry under the influence of Guido Zappa. In 1963/64 he was with Raj Chandra Bose at the University of North Carolina at Chapel Hill. He was a professor at the University of Perugia from 1968, at the University of Bologna from 1972 and at the University of Florence from 1982.
The current classification of projective planes is named after Barlotti and Hanfried Lenz. The Napoleon--Barlotti theorem (the centers of regular $n$-gons constructed over the sides of an $n$-gon $P$ form a regular $n$-gon if and only if $P$ is an affine image of a regular $n$-gon) is also named after him.
} in 1957. This classification extends the classification of Lenz figures, given in Theorem~\ref{t:Lenz}.
 
\begin{theorem}\label{t:Lenz-Barlotti} For every projective plane $\Pi$, its Lenz--Barlotti figure $\LB_\Pi$ is equal to one of the following $18$ sets: 
\begin{itemize}
\item[{$(\circ\circ)$}]\index[note]{$(\circ\circ)$} 
the empty set $\varnothing$;
\item[{$(\circ\bullet)$}]\index[note]{$(\circ\bullet)$} 
$\{(p,L)\}$ for some line $L\subseteq \Pi$ and point $p\in \Pi\setminus L$;

\item[{$(\cdidots)$}]\index[note]{$(\didots)$} 
$\{(p,L),(q,\Lambda)\}$ for some lines $L,\Lambda$ and points $p\in \Lambda\setminus L$ and $q\in L\setminus\Lambda$;
\item[{$(\circ\!\therefore$)}]\index[note]{$(\circ\hskip-2pt\therefore)$} 
$\{(a,\Aline bc),(b,\Aline ac),(c,\Aline ab)\}$ for some non-collinear points $a,b,c\in \Pi$;
\item[{$(\circ\backslash)$}]\index[note]{$(\circ\backslash)$} 
an injective function $F$ with $\dom[F]=L\setminus\{p\}$ and $\rng[F]=\mathcal L_p\setminus\{L\}$ for line $L\in\mathcal L_\Pi$ and point $p\in L$;
\item[{$(\bullet\circ)$}]\index[note]{$(\bullet\circ)$} 
$\{(p,L)\}$ for some line $L$ and point $p\in L$;
\item[{$(\bullet\bullet)$}]\index[note]{$(\bullet\bullet)$} 
$\{(p,L),(q,\Lambda)\}$ for some lines $L,\Lambda$ and points $p\in L\cap\Lambda$ and $q\in L\setminus\Lambda$;
\item[{$({-}\circ)$}]\index[note]{$({-}\circ)$} 
$\{(p,L):p\in L\}=L\times\{L\}$ for some line $L$;
\item[{$({-}\vert\vert)$}]\index[note]{$({-}\vert\vert)$} 
$(L\times\{L\})\cup(\{p\}\times\mathcal L_q)\cup(\{q\}\times\mathcal L_p)$ for some line $L$ and distinct points $p,q\in L$;
\item[{$({-}/\!/)$}]\index[note]{$(-/{\hskip-2pt}/)$}
$(L\times\{L\})\cup\bigcup_{p\in L}\{p\}\times\mathcal L_{\phi(p)}$ for some line $L$ and some irreflexive involution $\phi:L\to L$; moreover, $\Pi$ is a projective plane of order $9$;
\item[{$(\,|\,\circ)$}]\index[note]{$(\,\vert\circ)$} 
$\{p\}\times\mathcal L_p$ for some point $p$;
\item[{$(\,|{=})$}]\index[note]{$(\vert{=})$} 
$(\{p\}\times \mathcal L_p)\cup(L\times \{\Lambda\})\cup(\Lambda\times\{L\})$ for some point $p$ and distinct lines $L,\Lambda\in \mathcal L_p$;
\item[{$(\,|/\!/)$}]\index[note]{$(\,\vert/\hskip-2pt/)$} 
$(\{p\}\times\mathcal L_p)\cup\bigcup_{L\in \mathcal L_p}L\times \{\phi(L)\}$ for some point $p$ and some irreflexive involution $\phi:\mathcal L_p\to \mathcal L_p$; moreover, $\Pi$ is a projective plane of order $9$;
\item[{$({+}\circ)$}]\index[note]{$(+\circ)$}  
$(L\times\{L\})\cup(\{p\}\times\mathcal L_p)$ for some line $L$ and point $p\in L$;
\item[{$(\circledast\circ)$}]\index[note]{$(\circledast\circ)$} 
$\{(x,\Aline xp):x\in L\}$ for some line $L$ and point $p\notin L$;
\item[{$(\circledast\bullet)$}]\index[note]{$(\circledast\bullet)$} 
$\{(x,\Aline xp):x\in L\}\cup\{(p,L)\}$ for some line $L$ and point $p\notin L$;
\item[{$(\square\,\circ)$}]\index[note]{$(\square\,\circ)$} 
$\{(p,L):p\in L\in\mathcal L_\Pi\}$;
\item[{$(\square\blacksquare)$}]\index[note]{$(\square\blacksquare)$} 
$\Pi\times\mathcal L_\Pi$.
\end{itemize}
\end{theorem}

\begin{proof} By Theorem~\ref{t:Lenz}, for the Lenz figure $\Lenz_\Pi$ seven cases are possible.
\smallskip

$(\circ)$ The Lenz figure $\Lenz_\Pi$ is empty. In this case, the plane $\Pi$ is non-Desarguesian and hence $|\Pi|_2\ge 6$, by Corollary~\ref{c:p5-Pappian}.

\begin{claim}\label{cl:LB1} For any distinct pairs $(p,L),(q,\Lambda)\in\LB_\Pi$ we have $p\ne q$ and $L\ne \Lambda$. Consequently, $\LB_\Pi$ is an injective function.
\end{claim}

\begin{proof} It follows from $\Lenz_\Pi=\varnothing$ that $p\notin L$. Assuming that $p=q$, we conclude that $L\ne \Lambda$. Since $\Pi$ is a projective plane, the distinct lines $L,\Lambda$ have a unique common point $o$. Since $p\notin L$, the points $p$ and $o$ are distinct. Applying Corollary~\ref{c:(p,q]-transitive}, we conclude that the projective plane $\Pi$ is $(p,\mathcal L_o)$-transitive. In particular, $\Pi$ is $(p,\Aline op)$-transitive and hence $(p,\Aline op)\in \Lenz_\Pi$, which contradicts $\Lenz_\Pi=\varnothing$. This contradiction shows that $p\ne q$. 

Assuming that $L=\Lambda$, we can apply Theorem~\ref{t:(pq,L)-transitive} and conclude that the plane $\Pi$ is $(x,L)$-transitive for every point $x\in \Aline pq$. In particular, $\Pi$ is $(x,L)$-transitive for the unique point $x\in\Aline pq\cap L$, which implies $(x,L)\in\Lenz_\Pi$ and contradicts the assumption $\Lenz_\Pi=\varnothing$. This contradiction shows that $L\ne\Lambda$. 

Therefore, for any distinct pairs $(p,L),(q,\Lambda)\in\LB_\Pi$ we have  $p\ne q$ and $L\ne\Lambda$, which means that $\LB_\Pi$ is an injective function.
\end{proof}

\begin{claim}\label{cl:LB2} Let $(a,A),(b,B)\in\LB_\Pi$ be two distint point-line pairs.  If $a\in B$, then $b\in A$.
\end{claim}

\begin{proof} By Claim~\ref{cl:LB1}, $a\ne b$ and $A\ne B$. Let $o$ be the unique common point of the lines $A$ and $B$. It follows from $\Lenz_\Pi=\varnothing$ that $a\ne o\ne b$. To derive a contradiction, assume that $a\in B$ but $b\notin A$. Since $\Pi$ is $(b,B)$-transitive and $b\notin A\ne B$, there exists an automorphism $F\in\Aut_{b,B}(\Pi)$ such that $A'\defeq  F[A]\ne A$. Since $a\in B\subseteq \Fix(F)$, the $(a,A)$-transitivity of $\Pi$ implies $(a,A')$-transitivity of $\Pi$, see Proposition~\ref{p:pL-Auto}. By Corollary~\ref{c:(p,q]-transitive}, the projective plane $\Pi$ is $(a,L)$-transitive for every $L\in\mathcal L_o$. In particular, $\Pi$ is $(a,\Aline oa)$-transitive, which contradicts $\Lenz_\Pi=\varnothing$.
\end{proof}



\begin{claim}\label{cl:LB2n} Let $(a,A),(b,B)\in \LB_\Pi$ be two distinct point-line pairs,  and let $o$ be a common point of the lines $A,B$. If $b\notin A$, then $\Aline ab\setminus A\subseteq \dom[\LB_\Pi]$ and
 $\LB_\Pi[\Aline ab\setminus A]=\mathcal L_o\setminus\{\Aline ao\}$.
\end{claim}

\begin{proof} It follows from Claim~\ref{cl:LB1} that $a\ne b$ and $A\ne B$.  Assume that $b\notin A$. Applying Claim~\ref{cl:LB2}, we conclude that $a\notin B$. To check that $\Aline ab\setminus A\subseteq \dom[\LB_\Pi]$ and  $\LB_\Pi[\Aline ab\setminus A]\subseteq\mathcal L_o\setminus\{\Aline ao\}$, take any point $x\in \Aline ab\setminus A$. We have to show that $x\in\dom[\LB_\Pi]$ and $\LB_\Pi(x)\in\mathcal L_o\setminus\{\Aline ao\}$. This is clear if $x\in \{a,b\}$. So, assume that $x\notin \{a,b\}$. Since $(a,A)\in\LB_\Pi$, there exists an automorphism $\Phi\in\Aut_{a,A}(\Pi)$ such that $\Phi(b)=x$. Since $\{a,o\}\subseteq \{a\}\cup A\subseteq\Fix(\Phi)$, $\Phi[\Aline oa]=\Aline oa$ and hence $o\in \Phi[B]\ne\Phi[\Aline oa]=\Aline oa$. By Proposition~\ref{p:pL-Auto}, $(x,\Phi[B])=(\Phi(b),\Phi[B])\in\LB_\Pi$ and hence $x=\Phi(b)\in \dom[\LB_\Pi]$ and $\LB_\Pi[x]=\Phi[B]\in \mathcal L_o\setminus\{\Aline oa\}$, witnessing that $\Aline ab\setminus \{o\}\subseteq\dom[\LB_\Pi]$ and $\LB_\Pi[\Aline ab\setminus A]\subseteq\mathcal L_o\setminus\{\Aline oa\}$. 

To show that $\LB_\Pi[\Aline ab\setminus A]=\mathcal L_o\setminus\{\Aline ao\}$,  take any line $C\in\mathcal L_o\setminus\{\Aline ao\}$. We have to show that $C\in \LB_\Pi[\Aline ab\setminus A]$. If $C\in\{A,B\}$, then $C\in \{A,B\}=\LB_\Pi[\{a,b\}]\subseteq \LB_\Pi[\Aline ab\setminus A]$ and we are done. So, assume that $C\notin\{A,B\}$. Choose any point $\alpha\in A\setminus \{o\}$ and consider the unique points $\beta\in B\cap\Aline a\alpha$ and $\gamma\in C\cap\Aline a\alpha$. By the $(a,A)$-transitivity of the projective plane $\Pi$, there exists an automorphism $\Psi\in\Aut_{a,A}(\Pi)$ such that $\Psi(\beta)=\gamma$. Then $\Psi[B]=\Psi[\Aline o\beta]=\Aline o\gamma=C$ and $(\Psi(b),C)=(\Psi(b),\Psi[B])\in\LB_\Pi$, by Proposition~\ref{p:pL-Auto}. Since $a$ is the centre of the automorphism $\Psi$, we have $\Psi(b)\in \Aline ab\setminus A$ and hence $C=\LB_\Pi(\Psi(b))\in \LB_\Pi[\Aline ab\setminus A]$.
\end{proof}

\begin{claim}\label{cl:LB3} Let $(a,A)$ and $(b,B)$ be two distint point-line pairs in the Lenz--Barlotti figure $\LB_\Pi$ of $\Pi$. If $b\notin A$, then $A\cap B\subseteq \Aline ab$. 
\end{claim}

\begin{proof} By Claim~\ref{cl:LB1}, $\LB_\Pi$ is an injective function and hence $a\ne b$ and $A\ne B$. Let $o$ be the unique common point of the lines $A,B$.  If $b\notin A$, then $a\notin B$, by Claim~\ref{cl:LB2}. To derive a contradiction, assume that $\{o\}=A\cap B\not\subseteq\Aline ab$. By Claim~\ref{cl:LB2n}, 
$$
\begin{aligned}
\Aline ab&=\Aline ab\setminus(A\cap B)=(\Aline ab\setminus A)\cup(\Aline ab\setminus B)\subseteq\dom[\LB_\Pi]\quad\mbox{and}\\
\LB_\Pi[\Aline ab]&=\LB_\Pi[\Aline ab\setminus A]\cup\LB_\Pi[\Aline ba\setminus B]=(\mathcal L_o\setminus\{\Aline oa\})\cup(\mathcal L_o\setminus\{\Aline ob\})=\mathcal L_o.
\end{aligned}
$$
 Consider the affine plane $\mathbb A\defeq\Pi\setminus \Aline ab$ and identify the line $\Aline ab$ with the horizon line of the affine plane $\IA$. More precisely, each point $x\in\Aline ab$ is identified with the spread of parallel lines $\{L\setminus\Aline ab:L\in \mathcal L_x\setminus\{\Aline ab\}\}$ in $\IA$.  Let $\SC_\IA$ be the hyperscale trace of $\IA$. It is easy to show that a pair $(x,L)\in\Aline ab\times (\mathcal L_\Pi\setminus\{\Aline ab\})$ belongs to the Lenz-Barlotti figure $\LB_\Pi$ of the plane $\Pi$ if and only if the pair $(x,L\setminus\Aline ab)$ belongs  to the scale trace $\SC_\IA$ of the affine plane $\IA=\Pi\setminus\Aline ab$. This fact can be used to show that the scale trace $\SC_\IA$ of $\IA$ is an injective function (because $\LB_\Pi$ is an injective function) with $\dom[\SC_\IA]=\Aline ab=\partial\IA$. Lemma~\ref{l:PSC-involution=>9} ensures that $\IA$ is an affine plane of order $9$ and hence $\Pi$ is projective plane of order $9$. The unique projective plane of order $9$ with empty Lenz figure is the Hughes plane which also has empty Lenz-Barlotti figure. Since $\LB_\Pi\ne\varnothing$, we obtain a desired contradiction showing that $A\cap B=\{o\}\subseteq\Aline ab$.
\end{proof}

\begin{claim}\label{cl:LB5} Assume that the plane $\Pi$ contains a line $L\in\mathcal L_\Pi$ and a point $o\in L$ such that $L\setminus\{o\}\subseteq \dom[\LB_\Pi]$ and $\LB_\Pi[L\setminus\{o\}]=\mathcal L_o\setminus\{L\}$.  If $\dom[\LB_\Pi]\not\subseteq L$, then  $o\in\dom[\LB_\Pi]$.
\end{claim}

\begin{proof} 
Assuming that $\dom[\LB_\Pi]\not\subseteq L$, find a point-line pair $(c,C)\in\LB_\Pi$ such that $c\notin L$. Assuming that $C=L$ we can take any distinct points $x,y\in L\setminus \{o\}=C\setminus\{o\}$ and applying Claim~\ref{cl:LB2}, conclude that $c\in \LB_\Pi(x)\cap\LB_\Pi(y)=\{o\}\subseteq L$, which contradicts the choice of the pair $(c,C)$. This contradiction shows that $C\ne L$ and hence the lines $C,L$ have a unique common point $\gamma\in C\cap L$.

Assuming that $\gamma=o$, we conclude that $C\in\mathcal L_o\setminus\{L\}=\LB_\Pi[L\setminus\{o\}]$ and hence $C=\Phi(x)$ for some $x\in L\setminus\{o\}$. Then $(x,C),(c,C)\in \LB_\Pi$, which contradicts Claim~\ref{cl:LB1}. This contradiction shows that $o\ne\gamma$ and hence $o\notin C$. Since $|L|=|\Pi|_2\ge 4$, there exists a point $a\in L\setminus\{o,\gamma,\LB_\Pi^{-1}(\Aline oc)\}$. Consider the line $A\defeq\LB_\Pi(a)\in\mathcal L_o\setminus\{L\}$ and observe that $A=\LB_\Pi(a)\ne \Aline oc$. Since $|\Pi|_2\ge 4$, there exists a point $x\in \Aline ca\setminus(\{c,a\}\cup C)$. It follows from $a\ne o$ and $c\notin L$ that $\Aline ca\cap\Aline co=\{c\}$ and hence $x\notin\Aline oc$.  By the $(c,C)$-transitivity of $\Pi$, there exists an automorphism $\Phi\in\Aut_{c,C}(\Pi)$ such that $\Phi(a)=x$. Consider the line $X\defeq \Phi[A]$ and observe that $(x,X)=(\Phi(a),\Phi[A])\in\LB_\Pi$, according to Proposition~\ref{p:pL-Auto}.  Assuming that $o\in X$, we conclude that $\Phi^{-1}(o)\in A\cap \Aline oc=\{o\}$ and hence $o\in\Fix(\Phi)=C\cup\{c\}$, which is a desired contradiction showing that $o\notin X$.

 Since $|\Pi|_2\ge 6$, there exists a point $y\in \Aline ox\setminus(\{o,x\}\cup X\cup C\cup\Aline c\gamma)$. Let $b$ be the unique common point of the lines $L$ and $\Aline cy$. Assuming that $b=o$, we conclude that $y\in \Aline co\cap \Aline ox=\{o\}$, which contradicts the choice of $y$. This contradiction shows that $b\in L\setminus\{o\}$ and hence the line $B\defeq \LB_\Pi(b)\in\mathcal L_o\setminus\{L\}$ is well-defined. The choice of the point $y$ ensures that $b,y\notin \{c\}\cup C$ and $y\notin \{x\}\cup X$. By the $(c,C)$-transitivity, there exists an automorphism $\Psi\in\Aut_{c,C}(\Pi)$ such that $\Psi(b)=y$. Consider the line $Y\defeq \Psi[B]$ and observe that $(y,Y)=(\Psi(b),\Psi[B])\in\LB_\Pi$, according to Proposition~\ref{p:pL-Auto}. 
 By the $(x,X)$-transitivity of $\Pi$, there exists an automorphism $\varphi\in\Aut_{x,X}(\Pi)$ such that $\varphi(y)=o$. Then $(o,\varphi[Y])=(\varphi(y),\varphi[Y])\in \LB_\Pi$, witnessing that $o\in\dom[\LB_\Pi]$.
\end{proof}

\begin{claim}\label{cl:LB6} Assume that the plane $\Pi$ contains a line $L\in\mathcal L_\Pi$ and a point $o\in L$ such that $L\setminus\{o\}\subseteq \dom[\LB_\Pi]$ and $\LB_\Pi[L\setminus\{o\}]=\mathcal L_o\setminus\{L\}$. Then $\dom[\LB_\Pi]=L\setminus\{o\}$.
\end{claim}

\begin{proof} If $\dom[\LB_\Pi]\ne L\setminus\{o\}$, then $o\in\dom[\LB_\Pi]$, by Claim~\ref{cl:LB5}. Consider the line $\Lambda\defeq\LB_\Pi(o)$ and observe that $o\notin \Lambda$ (because $\Lenz_\Pi=\varnothing$) and hence $\Lambda\ne L$. Choose any point $x\in L\setminus(\{o\}\cup \Lambda)$. By the $(o,\Lambda)$-transitivity of $\Pi$, there exists an automorphism $F\in\Aut_{o,\Lambda}(\Pi)$ such that $x'\defeq F(x)\ne x$. Consider the line $X\defeq\LB_\Pi(x)\in\mathcal L_o\setminus\{L\}$. Since $o\in X$ is a centre of $F$, $F[X]=X$. By Proposition~\ref{p:pL-Auto}, $(x',X)=(F(x),F[X])\in\LB_\Pi$, and by Claim~\ref{cl:LB1}, $x'=x$, which contradicts the choice of the automorphism $F$. This contradiction shows that $\dom[\LB_\Pi]=L\setminus\{o\}$.
\end{proof}

Now we can finish the analysis of the case $(\circ)$. By Claim~\ref{cl:LB1}, $\LB_\Pi$ is an injective function. If $\LB_\Pi=\varnothing$, then we have the case $(\circ\circ)$ of the Lenz--Barlotti classification. So, assume that $\LB_\Pi\ne\varnothing$.
If $\LB_\Pi$ contains a unique point-line pair, then we have the case $(\circ\bullet)$. So, assume that $\LB_\Pi$ contains at least two distinct point-line pairs. 

If there exists two distinct point-line pairs $(h,V),(d,E)\in\LB_\Pi$ such that $h\notin E$ or $d\notin V$, then Claim~\ref{cl:LB1} ensures that $h\ne d$ and $V\ne E$. Claim~\ref{cl:LB3} ensures that the unique common point $v$ of the lines $V,E$ belongs to the line $L\defeq\Aline hd$, and Claims~\ref{cl:LB2n} and \ref{cl:LB6} ensure that $L\setminus\{v\}= \dom[\LB_\Pi]$ and $\LB_\Pi[\Aline hd\setminus\{v\}]=\mathcal L_v\setminus\{L\}$, so we have the case $(\circ\backslash)$.

So, assume that that for any distinct point-line pairs $(a,A),(b,B)\in \LB_\Pi$ we have $a\in B$ and $b\in A$. By our assumption, the Lenz--Barlotti figure $\LB_\Pi$ contains two distinct pairs $(a,A)$ and $(b,B)$. Our assumption guarantees that $a\in B$ and $b\in A$. If $\LB_\Pi=\{(a,A),(b,B)\}$, then we have the case $(\cdidots)$.
If $\LB_\Pi\ne\{(a,A),(b,B)\}$, then there exist a point-line pair $(c,C)\in \LB_\Pi\setminus\{(a,A),(b,B)\}$. Claim~\ref{cl:LB1} ensures that $a\ne c\ne b$ and $A\ne C\ne B$, and our assumption guarantees that $c\in A\cap B$, $a\in B\cap C$ and $c\in A\cap B$. Since the lines $A,B,C$ are distinct, the points $a,b,c$ are not collinear. Assuming that $\{(a,A),(b,B),(c,C)\}\ne\LB_\Pi$, we can find a point-line pair $(d,D)\in\LB_\Pi\setminus\{(a,A),(b,B),(c,C)\}$. Claim~\ref{cl:LB1} ensures that $d\notin\{a,b,c\}$ and $D\notin\{A,B,C\}$. Our assumption guarantees that $d\in A\cap B=\{c\}$ and hence $d=c$, which is a contradiction showing that $\LB_\Pi=\{(a,A),(b,B),(c,C)\}=\{(a,\Aline bc),(b,\Aline ac),(c,\Aline ab)\}$ and the case $(\circ\!\therefore)$ holds. This completes the analysis of the case $(\circ)$ of the Lenz classification.
\smallskip

$(\bullet)$ Next, assume that  $\Lenz_\Pi=\{(p,L)\}$ for some point $p\in\Pi$ and line $L\in\mathcal L_p$. If $\LB_\Pi=\{(p,L)\}$, then we have the case $(\bullet\circ)$. So, assume that $\LB_\Pi\ne\{(p,L)\}$.

\begin{claim}\label{cl:bullet12} For any point-line pair $(q,\Lambda)\in\LB_\Pi\setminus\{(p,L)\}$, we have $q\in L\setminus\Lambda$ and $p\in L\cap \Lambda$.
\end{claim}

\begin{proof} Since $(q,\Lambda)\notin\{(p,L)\}=\Lenz_\Pi$, $q\notin\Lambda$. If $p=q$, then $(p,L),(p,\Lambda)\in\LB_\Pi$ and Corollary~\ref{c:(p,q]-transitive} implies $\{p\}\times\mathcal L_o\subseteq \LB_\Pi$, where $o$ is the unique common point of the lines $L,\Lambda$. The inclusion $\{p\}\times\mathcal L_o$ implies that the projective plane $\Pi$ is $(p,o]$-transitive. By Theorem~\ref{t:pq<=>qp}, $\Pi$ is $(o,p]$-transitive. Then $(o,\Aline op)\in \Lenz_\Pi$, which contradicts the equality $\Lenz_\Pi=\{(p,L)\}$. This contradiction shows that $q\ne p$. If $p\notin \Lambda$, then by the $(q,\Lambda)$-transitivity of $\Pi$, there exists an automorphism $F\in\Aut_{q,\Lambda}(\Pi)$ such that $F(p)\ne p$. By Proposition~\ref{p:pL-Auto}, $(F(p),F[L])\in \Lenz_\Pi$, which contradicts $\Lenz_\Pi=\{(p,L)\}$. This contradiction shows that $p\in\Lambda$.

Assuming that $L=\Lambda$ and applyng Theorem~\ref{t:(pq,L)-transitive}, we conclude that the plane $\Pi$ is $[\Aline pq,L)$-transitive. By Corollary~\ref{c:LL'<=>L'L}, $\Pi$ is $[L,\Aline pq)$-transitive.
Then $(p,\Aline pq)\in \Lenz_\Pi$, which contradicts $\Lenz_\Pi=\{(p,L)\}$. This contradiction shows that $L\ne \Lambda$. Assuming that $q\notin L$, we can find a point $x\in L\setminus\Lambda$. Since $\Lenz_{\Pi}=\{(p,L)\}$, the projective plane $\Pi$ is not Desarguesian and hence $|\Pi|_2\ge 4$. Then there exists a point $y\in \Aline qx\setminus(\{q,x\}\cup \Lambda)$. Since $\Aline qx\cap L=\{x\}$, the point $y$ does not belong to the line $L$.  By the $(q,\Lambda)$-transitivity of $\Pi$, there exists an automorphism $F\in\Aut_{q,\Lambda}(\Pi)$ such that $F(x)=y$. Since $y\in F[L]\setminus L$, the line $F[L]$ is distinct from the line $L$ and hence the pair $(F(p),F[L])$ is distinct from the pair $(p,L)$. Proposition~\ref{p:pL-Auto} ensures that $(F(p),F[L])\in\Lenz_\Pi$, which contradicts $\Lenz_\Pi=\{(p,L)\}$. This contradiction shows that $q\in L\setminus\Lambda$.
\end{proof}

If $|\LB_\Pi|>2$, then we can find two distinct pairs $(h,V),(d,E)\in \LB_\Pi\setminus\{(p,L)\}$.  Claim~\ref{cl:bullet12} ensures that $p\in V\cap E$ and $\{h,d\}\subseteq L\setminus\{p\}$. Assuming that $h=d$, we can apply Corollary~\ref{c:(p,q]-transitive} and conclude that the projective plane $\Pi$ is $(h,p]$-transitive and hence $(h,\Aline pq)\in \Lenz_\Pi=\{(p,L)\}$, which is a contradiction showing that $h\ne d$. Assuming that $V=E$, we can apply Theorem~\ref{t:(pq,L)-transitive} and conclude that the projective plane $\Pi$ is $[L,V)$-transitive, which implies $(p,V)\in\Lenz_\Pi$ and contradicts the equality $\Lenz_\Pi=\{(p,L)\}$. Therefore, $h\ne d$ and $V\ne E$. Consider the affine plane $\IA\defeq\Pi\setminus L$ and choose any point $o\in V\setminus L\subseteq\IA$. Next, find unique points $e\in E\cap\Aline od$, $u\in E\cap\Aline pe$ and $w\in V\cap\Aline he$. It follows from $\{(p,L),(h,V),(d,E)\}\subseteq\LB_\Pi$ that the based affine plane $(\IA,uow)$ is vertical-translation, horizontal-scale and diagonal-scale. By Theorems~\ref{t:cart-group<=>} and \ref{t:diagonal-scale<=>}, its ternar $\Delta=\Aline oe$ is horizontal-scale, diagonal-scale and associative-plus. By Theorem~\ref{t:hdiscale<=>corps}, the ternar $\Delta$ is a corps, so is distributive. By Theorem~\ref{t:corps<=>}, the affine plane $\IA$ is Desarguesian and so is its projective completion $\Pi$, (by Theorem~\ref{t:projDes<=>}), which contradicts $\Lenz_\Pi=\{(p,L)\}$. This contradiction shows that $|\LB_\Pi|=2$ and hence $\LB_\Pi=\{(p,L),(q,\Lambda)\}$ for some pair $(q,\Lambda)\in (L\setminus\{p\})\times(\mathcal L_p\setminus\{L\})$, witnessing that the case  $(\bullet\bullet)$ holds. 
\smallskip

$(-)$ Next, assume that $\Lenz_\Pi=L\times\{L\}$ for some line $L$ in $\Pi$. 

\begin{claim}\label{cl:LB7} For every $(q,\Lambda)\in \LB_\Pi\setminus\Lenz_\Pi$ we have $q\in L\setminus \Lambda$ and $(\{q\}\times \mathcal L_o)\cup(\{o\}\times\mathcal L_q)\subseteq \LB_\Pi$, where $o$ is a unique common point of the lines $L$ and $\Lambda$.
\end{claim}

\begin{proof} Since $(q,\Lambda)\in\LB_\Pi\setminus\Lenz_\Pi$, $q\notin\Lambda$. Assuming that $L=\Lambda$, we can apply Theorem~\ref{t:(pq,L)-transitive} and conclude that the projective plane $\Pi$ is $[\Aline pq,L)$-transitive for every point $p\in L$. By Theorem~\ref{t:pq<=>qp}, $\Pi$ is $[L,\Aline pq)$-transitive and hence $\{(p,\Aline pq):p\in L\}\subseteq \Lenz_\Pi$, which contradicts the assumption $\Lenz_\Pi=L\times\{L\}$. This contradiction shows that $L\ne\Lambda$. Let $o$ be the unique common point of the lines $L$ and $\Lambda$.
Assuming that $q\notin L$, we can apply the $(q,\Lambda)$-transitivity of $\Pi$ and find an automorphism $\Phi\in\Aut_{q,\Lambda}(\Pi)$ such that $F[L]\ne L$. By Proposition~\ref{p:pL-Auto}, for every $p\in L$, we have $(F(p),F[L])\in \Lenz_\Pi\setminus(L\times\{L\})=\varnothing$, which is a contradiction showing that $q\in L$. By Corollary~\ref{c:(p,q]-transitive}, the $(q,L)$-transitivity and $(q,\Lambda)$-transitivity of $\Pi$ imply the $(q,o]$-transitivity of $\Pi$. By Theorem~\ref{t:pq<=>qp}, the $(q,o]$-transitivity of $\Pi$ implies $(o,q]$-transitivity of $\Pi$. Since $\Pi$ is $(q,o]$-transitive and $(o,q]$-transitive, we have $(\{q\}\times\mathcal L_o)\cup(\{o\}\times\mathcal L_q)\subseteq \LB_\Pi$.
\end{proof}

With Claim~\ref{cl:LB7} in our disposition we can now complete the analysis of the case $(-)$.
If $\LB_\Pi=L\times\{L\}$, then the case $(-\circ)$ holds. So, assume that $\LB_\Pi\ne L\times\{L\}$ and choose any pair $(a,A)\in\LB_\Pi\setminus(L\times\{L\})$. Claim~\ref{cl:LB7} ensures that $a\in L\setminus A$ and $${+}\!\!{+}\defeq (L\times\{L\})\cup (\{a\}\times\mathcal L_{a'})\cup(\{a'\}\times\mathcal L_a)\subseteq \LB_\Pi$$ where $a'$ is the unique common point of the lines $L$ and $A$. If $\LB_\Pi={+}\!\!{+}$, then we have the case $({-}||)$. So, assume that $\LB_\Pi\ne{+}\!\!{+}$ and choose any pair $(b,B)\in \LB_\Pi\setminus {+}\!\!{+}$. Claim~\ref{cl:LB7} ensures that $b\in L\setminus B$ and $(\{b\}\times\mathcal L_{b'})\cup(\{b'\}\times\mathcal L_b)\subseteq\LB_\Pi$ where $b'$ is the unique common point of the lines $L$ and $B$. 

Since $(b,B)\in (\{b\}\times\mathcal L_{b'})\setminus{+}\!\!{+}$, the pair $(b,b')$ is not equal to the pair $(a,a')$ or $(a',a)$. We claim that $\{a,a'\}\cap\{b,b'\}=\varnothing$. 
Indeed, assuming that $a=b$, we conclude that $a'\ne b'$. For every point $o\in \Pi\setminus L$, consider the unique lines $A'\in \mathcal L_{a'}\cap\mathcal L_o$ and $B'\in\mathcal L_{b'}\cap\mathcal L_o$. Since the point-line pairs $(a,A')$ and $(a,B')=(b,B')$ belong to the Lenz--Barlotti figure, Corollary~\ref{c:(p,q]-transitive} ensures that $\{a\}\times\mathcal L_o\subseteq\LB_\Pi$, and Theorem~\ref{t:pq<=>qp} ensures that $\{o\}\times\mathcal L_a\subseteq\LB_\Pi$. Then $(o,\Aline oa)\in\Lenz_\Pi=L\times\{L\}$, which is a contradiction showing that $a\ne b$. By analogy we can prove that $a\ne b'$ and $a'\notin\{b,b'\}$. Therefore, $\{a,a'\}\cap\{b,b'\}=\varnothing$.

For every $x\in L\setminus\{a,a'\}$, by the $(a,A)$-transitivity of $\Pi$, there exists an automorphism $F_x\in\Aut_{a,A}(\Pi)$ such that $F_x(b)=x$. Let $x'\defeq F_x(b')$. It follows from $(\{b\}\cup\mathcal L_{b'})\cup(\{b'\}\times\mathcal L_b)\subseteq \LB_\Pi$ that $(\{x\}\cup\mathcal L_{x'})\cup (\{x'\}\times\mathcal L_x)\subseteq \LB_\Pi$, see Proposition~\ref{p:pL-Auto}. Define a map $\phi:L\to L$ letting $\phi(a)=a'$, $\phi(a')=a$, and $\phi(x)=x'$ for $x\in L\setminus\{a'\}$.  Then $$\bigcup_{x\in L}(\{x\}\times \mathcal L_{\phi(x)})\subseteq \LB_\Pi.$$
Let us show that $\phi$ is an irreflexive involutive bijection of the line $L$. The irreflexivity of $\phi$ follows from the inequalities $a\ne a'$ and $b\ne b'$. It suffices to show that $\phi(\phi(x))=x$ for every $x\in L$. For $x\in\{a,a'\}$ this follows from the definition of $\phi(a)=a'$ and $\phi(a')=a$. Now take any point $x\in L\setminus\{a,a'\}$ and consider the points $x'=\phi(x)$ and $x''\defeq\phi(x')$. The definition of $\phi$ ensures that
$$(\{x\}\times\mathcal L_{x'})\cup(\{x'\}\times\mathcal L_x)\cup(\{x'\}\times\mathcal L_{x''})\cup(\{x''\}\times\mathcal L_{x'})\subseteq\LB_\Pi.$$
Assuming that $x''\ne x$, we conclude that for every line $\Lambda\in\mathcal L_{x'}\setminus\{L\}$, the plane $\Pi$ is $(x,\Lambda)$-transitive and $(x'',\Lambda)$-transitive. By Theorem~\ref{t:(pq,L)-transitive}, $\Pi$ is $[\Aline x{x''},\Lambda)$-transitive. Since $x'\in \Lambda\cap L=\Lambda\cap\Aline x{x''}$, the plane $\Pi$ is $(x',\Lambda)$-transitive and hence $(x',\Lambda)\in \Lenz_\Pi=L\times\{L\}$, which is a contradiction showing that $x''=x$ and hence $\phi$ is an involutive bijection.

Assuming that $\LB_\Pi\ne  \bigcup_{x\in L}(\{x\}\times \mathcal L_{\phi(x)})$, we can find a point-line pair $(c,C)\in\LB_\Pi\setminus \bigcup_{x\in L}(\{x\}\times \mathcal L_{\phi(x)}).$ Claim~\ref{cl:LB7} ensures that $c\in L\setminus C$ and $(c,C)\in \{c\}\times\mathcal L_{c'}\subseteq \mathcal LB_\Pi$, where $c'$ is the unique common point of the lines $L$ and $C$. The choice of $(c,C)$ ensures that $c'\ne\phi(c)$. For every $o\in \Pi\setminus L$, consider the unique lines $\Lambda\in \mathcal L_o\cap\mathcal L_{\phi(c)}$ and $\Lambda'\in \mathcal L_o\cap\mathcal L_{c'}$. Since $\{(c,\Lambda),(c,\Lambda')\}\subseteq(\{c\}\times\mathcal L_{\phi(c)})\cup(\{c\}\times\mathcal L_{c'})\subseteq\LB_\Pi$, Corollary~\ref{c:(p,q]-transitive} implies that $\{c\}\times\mathcal L_o\subseteq \LB_\Pi$, and Theorem~\ref{t:pq<=>qp} implies that $\{o\}\times\mathcal L_c\subseteq\LB_\Pi$. In particular, $(o,\Aline co)\in\Lenz_\Pi=L\times \{L\}$, which is a contradiction showing that $\LB_\Pi=\bigcup_{x\in L}\{x\}\times\mathcal L_{\phi(x)}$ and hence the case $({-}/\!/)$ holds. It remains to show that the plane $\Pi$ has order 9. Consider the affine plane $\IA\defeq\Pi\setminus L$ and observe that $\SC_\IA=\bigcup_{x\in L}\{\delta\}\times \phi(x)$, where we identity the line $L$ with the horizon line of the affine plane $\IA$. Applying Theorem~\ref{t:scale-trace}$(\invol)$, we conclude that $\IA$ is an affine plane of order $9$ and hence $\Pi$ is a projective plane of order $9$.  
\smallskip

$(\,|\,)$ Assume that $\Lenz_\Pi=\{p\}\times\mathcal L_p$ for some point $p\in\Pi$. This case can be deduced from the case $(-)$ by duality, but we present a direct proof. 

\begin{claim}\label{cl:LB8} For every $(q,\Lambda)\in \LB_\Pi\setminus (\{p\}\times\mathcal L_p)$ we have $q\notin\Lambda\in \mathcal L_p$ and $$(\Aline pq\times \{\Lambda\})\cup(\Lambda\times\{\Aline pq\})\subseteq \LB_\Pi.$$
\end{claim}

\begin{proof} Since $(q,\Lambda)\in\LB_\Pi\setminus(\{p\}\times\mathcal L_p)=\LB_\Pi\setminus\Lenz_\Pi$, $q\notin\Lambda$. Assuming that $p=q$, we conclude that $\Lambda\notin \mathcal L_p$. For every point $x\in \Lambda$, we have $\{(p,\Lambda),(p,\Aline xp)\}=\{(q,\Lambda),(p,\Aline xp)\}\subseteq \LB_\Pi$. By Corollary~\ref{c:(p,q]-transitive}, $\{p\}\times\mathcal L_x\subseteq \LB_\Pi$, and by Theorem~\ref{t:pq<=>qp}, $\{x\}\times\mathcal L_p\subseteq\LB_\Pi$. In particular, $(x,\Aline xp)\in\Lenz_\Pi=\{p\}\times\mathcal L_p$ and hence $q=p=x\in \Lambda$, which contradicts the choice of $(q,\Lambda)\notin\Lenz_\Pi$. This contradiction shows that $p\ne q$.

Assuming that $p\notin \Lambda$, we can apply the $(q,\Lambda)$-transitivity of $\Pi$ and find an automorphism $\Phi\in\Aut_{q,\Lambda}(\Pi)$ such that $F(p)\ne p$. By Proposition~\ref{p:pL-Auto}, $\{p\}\times\mathcal L_p\subseteq\LB_\Pi$ implies $\{F(p)\}\times \mathcal L_{F(p)}\subseteq\Lenz_\Pi=\{p\}\times\mathcal L_p$, which is a contradiction showing that $p\in\Lambda$ and hence $\Lambda\in \mathcal L_p$. By Theorem~\ref{t:(pq,L)-transitive}, the $(p,\Lambda)$-transitivity and $(q,\Lambda)$-transitivity of $\Pi$ imply the $[\Aline pq,\Lambda)$-transitivity of $\Pi$. By Theorem~\ref{t:pq<=>qp}, the $[\Aline pq,\Lambda)$-transitivity of $\Pi$ implies $[\Lambda,\Aline pq)$-transitivity of $\Pi$. Then $(\Aline pq\times\{\Lambda\})\cup(\{\Lambda\}\times\{\Aline pq\})\subseteq \LB_\Pi$.
\end{proof}

With Claim~\ref{cl:LB8} in our disposition we can now complete the analysis of the case $(\,|\,)$.
If $\LB_\Pi=\{p\}\times\mathcal L_p$, then the case $(\,|\circ)$ holds. So, assume that $\LB_\Pi\ne \{p\}\times \mathcal L_p$ and choose any pair $(a,A)\in\LB_\Pi\setminus(\{p\}\times\mathcal L_p)$. Claim~\ref{cl:LB8} ensures that $a\notin A\in\mathcal L_p$ and $${+}\!\!{+}\defeq (\{p\}\times\mathcal L_p)\cup (\Aline ap\times\{A\})\cup (A\times\{\Aline ap\})\subseteq \LB_\Pi.$$ If $\LB_\Pi={+}\!\!{+}$, then we have the case $(\,|{=})$. So, assume that $\LB_\Pi\ne{+}\!\!{+}$ and choose any pair $(b,B)\in \LB_\Pi\setminus {+}\!\!{+}$. Claim~\ref{cl:LB8} ensures that $b\notin B\in \mathcal L_p$ and $(\Aline bp\times\{B\})\cup(B\times\{\Aline bp\})\subseteq\LB_\Pi$. The choice of $(b,B)\notin\{p\}\times\mathcal L_p$ ensures that $b\ne p$. 

Assuming that $B=A$, we conclude that $b\notin \Aline ap$. Since $\{(a,A),(b,A)\}=\{(a,A),(b,B)\}\subseteq\LB_\Pi$, the plane $\Pi$ is $[\Aline ab,A)$-transitive and $[A,\Aline ab)$-transitive, by Theorems~\ref{t:(pq,L)-transitive} and \ref{t:pq<=>qp}. In particular, for the unique point $x\in A\cap\Aline ab$, the plane $\Pi$ is $(x,\Aline ab)$-transitive and hence $(x,\Aline ab)\in\Lenz_\Pi=\{p\}\times\mathcal L_p$ and $p=x\in \Aline ab$, which contradicts $b\notin \Aline ap$. This contradiction shows that $B\ne A$. 

Next, we show that $\Aline bp\ne\Aline ap$. In the opposite case, $(\Aline ap\times\{A\})\cup(\Aline bp\times\{B\})\subseteq\LB_\Pi$ imply $\Aline ap\times\{\Aline ap\}\subseteq \Aline ap\times\mathcal L_p\subseteq\LB_\Pi$, according to Corollary~\ref{c:(p,q]-transitive}. Then $\Aline ap\times\{\Aline ap\}\subseteq\Lenz_\Pi$, which contradicts the equality $\Lenz_\Pi=\{p\}\times\mathcal L_p$. This contradiction shows that $\Aline bp\ne\Aline ap$ and hence $p\notin\Aline ab$.

For every point $x\in \Aline ab\setminus(\{a\}\cup A)$, by the $(a,A)$-transitivity of $\Pi$, there exists an automorphism $F_x\in\Aut_{a,A}(\Pi)$ such that $F_x[b]=x$. It follows from $(\Aline bp\times \{B\})\cup(B\times\{\Aline bp\})\subseteq \LB_\Pi$ that $(\Aline xp\times\{F_x[B]\})\cup (F_x[B]\times\{\Aline xp\})\subseteq \LB_\Pi$, see Proposition~\ref{p:pL-Auto}. Define a map $\phi:\mathcal L_p\to \mathcal L_p$ letting $\phi(\Aline ap)\defeq A$, $\phi(A)\defeq \Aline ap$, and $\phi(\Aline xp)=F_x[B]$ for $x\in \Aline ab\setminus(\{a\}\cup A)$. Then $$\bigcup_{L\in\mathcal L_p}(L\times \{\phi(L)\}\subseteq \LB_\Pi.$$
Let us show that $\phi$ is an irreflexive involutive bijection of the pencil of lines $\mathcal L_p$. The irreflexivity of $\phi$ follows from $\Aline bp\ne B$. To show that the function $\phi$ is involutive, it suffices to check that $\phi(\phi(L))=L$ for every line $L\in \mathcal L_p$. If $L\in\{\Aline ap,A\}$, then $\phi(\phi(L))=L$ follows from $\phi(\Aline ap)\defeq A$ and $\phi(A)\defeq\Aline ap$. Now take any line $L\in\mathcal L_p\setminus\{\Aline ap,A\}$ and consider the lines $L'=\phi(L)$ and $L''\defeq\phi(L')$. The definition of $\phi$ ensures that
$$(L\times\{L'\})\cup(L'\times\{L\})\cup(L''\times\{L'\})\cup(L'\times\{L''\})\subseteq\LB_\Pi.$$
Assuming that $L''\ne L$, we conclude that for every point $x\in L'\setminus\{p\}$, the plane $\Pi$ is $(x,L)$-transitive and $(x,L'')$-transitive. By Corollary~\ref{c:(p,q]-transitive} and Theorem~\ref{t:pq<=>qp}, $\Pi$ is $(x,p]$-transitive and hence $(x,\Aline xp)\in\Lenz_\Pi$, which contradicts $\Lenz_\Pi=\{p\}\times\mathcal L_p$. This contradiction shows that $L''=L$ and hence $\phi$ is an involutive bijection of the set $\mathcal L_p$.

Assuming that $\LB_\Pi\ne  \bigcup_{L\in\mathcal L_p}(L\times\{\phi(L)\})$, we can find a point-line pair $(c,C)\in\LB_\Pi\setminus \bigcup_{L\in \mathcal L_p}(L\times\{\phi(L)\})$. Claim~\ref{cl:LB8} ensures that $c\notin C\in\mathcal L_p$ and $\Aline cp\times\{C\}\subseteq \LB_\Pi$. The choice of $(c,C)$ ensures that $c\ne p$ and $C\ne \phi(\Aline cp)$. By Corollary~\ref{c:(p,q]-transitive}, $\{(c,C),(c,\phi(\Aline cp)\}\subseteq\LB_\Pi$ imply $(c,\Aline cp)\in \{c\}\times\mathcal L_p\subseteq\LB_\Pi$, which contradicts $\Lenz_\Pi=\{p\}\times\mathcal L_p$. This contradiction shows that $\LB_\Pi=\bigcup_{L\in\mathcal L_p}(L\times\{\phi(L)\})$. Applying Theorem~\ref{t:scale-trace}$(\begin{picture}(10,10)(0,2)\linethickness{=0.7pt}
\put(5,0){\line(1,1){5}}
\put(0,5){\line(1,1){5}}
\end{picture})$ we can show that $\Pi$ is a projective plane of order 9. Therefore, the case $(\,|/\!/)$ holds.
\smallskip

$(+)$ Assume that $\Lenz_\Pi=(L\times\{L\})\cup(\{p\}\times\mathcal L_p)$ for some point $p\in\Pi$ and line $L\in\mathcal L_p$.  Assuming that $\LB_\Pi\ne\Lenz_\Pi$, find a  point-line pair $(q,\Lambda)\in\LB_\Pi\setminus\Lenz_\Pi$ and observe that $q\notin\Lambda$.

If $\Lambda=L$, then $q\notin \Lambda=L$. Fix a projective base $uowe$ in the plane $\Pi$ such that $o=q$,  $L=\Lambda=\Aline hv$ where $h\in \Aline ou\cap\Aline we$, $v\in\Aline ow\cap\Aline ue$ and $v=p$. It follows from $ (\Aline hv\times\{\Aline hv\})\cup\{(o,\Aline hv)\}=(L\times\{L\})\cup\{(q,\Lambda)\}\subseteq \LB_\Pi$ that the based affine plane $(\Pi\setminus\Aline hv,uow)$ is translation and homocentral. By Theorem~\ref{t:corps<=>}, the affine plane $\Pi\setminus\Aline hv$ is Desarguesian and so is its projective completion $\Pi$, which implies that $\Pi$ is of Lenz type $(\square)$, which contradicts our assumption. This contradiction implies that $\Lambda\ne L$. 

If $q\notin L$, then the $(q,\Lambda)$-transitivity of $\Pi$ implies the existence of an automorphism $\Phi\in\Aut_{q,\Lambda}(\Pi)$ such that $\Phi[L]\ne L$. Since $L$ is a translation line in $\Pi$, so is the line $\Phi[L]$. By Theorem~\ref{t:Skornyakov-San-Soucie}, the plane $\Pi$ is Moufang and hence it has the Lenz type $(\square)$, which contradicts our assumption. This contradiction shows that $q\in L$. 

If $q=p$, then choose a projective base $uowe$ in the projective plane $\Pi$ such that $\Aline ou=\Lambda$, $\Aline ou\cap\Aline we\in L$ and $p=q\in\Aline ow\cap\Aline ue$. It follows from $\{(q,\Lambda),(p,\lambda)\}\cup(L\times\{L\})\subseteq \LB_\Pi$ that the affine plane $(\Pi\setminus L,uow)$ is translation, vertical-scale and vertical-shear. By Corollary~\ref{c:Algebra-vs-Geometry-in-trings:affine}(13), the affine plane $\Pi\setminus L$ is Desarguesian and so is the projective plane $\Pi$. Then it is of Lenz type $(\square)$, which contradicts our assumption. This contradiction shows that $q\ne p$. If $p\notin\Lambda$, then we can find an automorphism $\Phi\in\Aut_{q,\Lambda}(\Pi)$ such that $p'\defeq \Phi(p)\ne p$. Then $\{p\}\times\mathcal L_p\subseteq\Lenz_\Pi$ implies  $\{p'\}\times\mathcal L_{p'}\subseteq \Lenz_\Pi=(L\times\{L\})\cup\{p\}\times\mathcal L_p)$, which is a contradiction implying $p\in\Lambda$. 

Choose a projective base $uowe$ in the projective plane $\Pi$ such that $p\in \Aline ow\cap\Aline ue$, $q\in \Aline ou\cap\Aline we$, and $\Aline ow=\Lambda$. It follows from $\{(q,\Lambda)\}\cup(L\times\{L\})\cup(\{p\}\times\mathcal L_p)\subseteq \LB_\Pi$ that the affine plane $(\Pi\setminus L,uow)$ is translation, vertical-shear, and horizontal-scale. By Corollary~\ref{c:vshear+vh-scale+htrans}, the ternar $\Delta=\Aline oe\setminus\Aline hv$ of the affine plane $\Pi\setminus\Aline hv$ is a corps and by Theorem~\ref{t:corps<=>}, the affine plane $\Pi\setminus\Aline hv$ is Desarguesian and so is its projective completion $\Pi$. Then $\Pi$ is of Lenz type $(\square)$, which contradicts our assumption. This contradiction implies that $\LB_\Pi=\Lenz_\Pi$ and hence $\Pi$ is of Lenz--Barlotti type $(+\circ)$. 

\smallskip

$(\circledast)$ Assume that $\Lenz_\Pi=\{(x,\Aline xp):x\in L\}$ for some line $L\in\mathcal L_\Pi$ and point $p\in \Pi\setminus L$. 

\begin{claim}\label{cl:LB9} For every $(q,\Lambda)\in\LB_\Pi\setminus\Lenz_\Pi$, 
$$q=p\;\Leftrightarrow\;\Lambda=L.$$
\end{claim}

\begin{proof} If $p=q$ and $\Lambda\ne L$, then we can choose a point $x\in L\setminus(\Lambda\cup\{q\})$ and using the $(q,\Lambda)$-transitivity of $\Pi$, find an an automorphism $F\in\Aut_{q,\Lambda}(\Pi)$ such that $y\defeq F(x)\ne x$. By Proposition~\ref{p:pL-Auto}, $(x,\Aline xp)\in\Lenz_\Pi$ implies $(y,\Aline yp)\in\Lenz_\Pi$ and hence $y\in L$ and $p=q\in \Aline xy=L$, which is a contradiction showing that $p=q$ implies $\Lambda=L$.

If $\Lambda=L$ and $p\ne q$, then by the $(q,\Lambda)$-transitivity, there exists an automorphism $\Phi\in\Aut_{q,\Lambda}(\Pi)$ such that $p'\defeq\Phi(p)\ne p$. Taking into account that $L=\Lambda\subseteq\Fix(\Phi)$, we can apply Proposition~\ref{p:pL-Auto} and conclude that $\{(x,\Aline x{p'}):x\in L\}\subseteq \Lenz_\Pi$, which contradicts $\Lenz_\Pi=\{(x,\Aline xp):x\in L\}$. This contradiction shows that $\Lambda=L$ implies $p=q$.
\end{proof}

\begin{claim}\label{cl:LB10} Every point-line pair $(q,\Lambda)\in \LB_\Pi\setminus\Lenz_\Pi$ it equal to the point-line pair $(p,L)$.
\end{claim}

\begin{proof} Assuming that $(q,\Lambda)\ne (p,L)$, we can apply Claim~\ref{cl:LB9} and conclude that $q\ne p$ and $\Lambda\ne L$. If $p\notin \Lambda$, then by the $(q,\Lambda)$-transitivity of $\Pi$, there exists an automorphism $F\in\Aut_{q,\Lambda}(\Pi)$ such that $F(p)\ne p$. Consider the point $p'\defeq F(p)\ne p$ and line $L'\defeq\Phi[L]$. It follows from $p\notin L$ that $p'\notin L'$. Then there exists a point $x'\in L'\setminus\Aline p{p'}$. Applying Proposition~\ref{p:pL-Auto}, we conclude that $(x',\Aline {x'}{p'})\in \{(x,\Aline x{p'}):x\in L'\}\subseteq\Lenz_\Pi$ and hence $(x',\Aline{x'}{p'})=(y,\Aline yp)$ for some $y\in L$. Then $\Aline {x'}{p'}=\Aline yp$ and hence $x'\in \Aline p{p'}$, which contradicts the choice of $x'$. This contradiction shows that $p\in\Lambda$. 

Consider the unique point $x\in L\cap \Lambda$ and observe that $(x,\Lambda)=(x,\Aline xp)\in \Lenz_\Pi$. Since $q\notin\Lambda$, the $(x,\Lambda)$-transitivity and $(q,\Lambda)$-transitivty of $\Pi$ implies the $[\Aline xq,\Lambda)$-transitivity and $[\Lambda,\Aline xq)$-transitivity of $\Pi$. Then $(x,\Aline xq)\in\Lenz_\Pi$ and hence $\Aline xq=\Aline xp=\Lambda$, which contradicts $q\notin\Lambda$. This contradiction shows that $(q,\Lambda)=(p,L)$.
\end{proof}

With Claim~\ref{cl:LB10} in our disposition, we can now complete the analysis of the case $(\circledast)$. If $\LB_\Pi=\Lenz_\Pi$, then we have the case $(\circledast\circ)$. If $\LB_\Pi\ne\Lenz_\Pi$, then Claim~\ref{cl:LB10} ensures that $\LB_\Pi=\Lenz_\Pi\cup\{(p,L)\}$ and hence we have the case $(\circledast\bullet)$.
\smallskip

$(\square)$ Assume that $\Lenz_\Pi=\{(p,L)\in\Pi\times\mathcal L_\Pi:p\in L\}$. If $\LB_\Pi=\Lenz_\Pi$, then the case $(\square\circ)$ holds. So, assume that $\LB_\Pi\ne\Lenz_\Pi$ and choose a point-line pair $(p,L)\in\LB_\Pi\setminus\Lenz_\Pi$. Then $p\notin L$. Choose a projective base $uowe$ in $\Pi$ such that $L=\Aline ow=\Aline ov$ and $p=h\in\Aline ou\cap\Aline we$ is the horizontal infinity point of the base $uowe$. It follows from $(h,\Aline ov)=(p,L)\in\LB_\Pi$ that the projective plane $\Pi$ is $(h,\Aline ov)$-transitive. By Theorem~\ref{t:Algebra-vs-Geometry-proj}(5,14), the ternar $\Delta=\Aline oe\setminus\Aline hv$ is linear, distributive, and associative, and by Theorem~\ref{t:Algebra-vs-Geometry-proj}(12), $\LB_\Pi=\Pi\times\mathcal L_\Pi$ and hence the case $(\square\blacksquare)$ holds.
\end{proof}

Now let us consider the interplay between the Lenz-Barlotti class of a projective plane and the algebraic properties of its ternars.

\begin{definition}\label{d:po-LB} A projective plane $\Pi$ is defined to have the Lenz--Barlotti type at least
\begin{itemize}
\item[{$(\circ\circ)$}] if $\varnothing\subseteq\LB_\Pi$;
\item[{$(\circ\bullet)$}] if $(p,L)\in\LB_\Pi$ for some line $L\subseteq \Pi$ and point $p\in \Pi\setminus L$;

\item[{$(\cdidots)$}] if $\{(p,L),(q,\Lambda)\}\subseteq\LB_\Pi$ for some lines $L,\Lambda$ and points $p\in \Lambda\setminus L$ and $q\in L\setminus\Lambda$;
\item[{$(\circ\!\therefore$)}] if $\{(a,\Aline bc),(b,\Aline ac),(c,\Aline ab)\}\subseteq\LB_\Pi$ for some non-collinear points $a,b,c\in \Pi$;
\item[{$(\circ\backslash)$}] if $F\subseteq\LB_\Pi$ for some line $L$, point $p\in L$ and a bijective function $F:L\setminus\{p\}\mapsto \mathcal L_p\setminus\{L\}$;
\item[{$(\bullet\circ)$}] if $(p,L)\in\LB_\Pi$ for some line $L$ and point $p\in L$;
\item[{$(\bullet\bullet)$}] if $\{(p,L),(q,\Lambda)\}\subseteq\LB_\Pi$ for some lines $L,\Lambda$ and points $p\in L\cap\Lambda$ and $q\in L\setminus\Lambda$;
\item[{$({-}\circ)$}] if $L\times\{L\}\subseteq\LB_\Pi$ for some line $L$;
\item[{$({-}\vert\vert)$}] if $(L\times\{L\})\cup(\{p\}\times\mathcal L_q)\cup(\{q\}\times\mathcal L_p)\subseteq\LB_\Pi$ for some line $L$ and distinct points $p,q\in L$;
\item[{$({-}/\!/)$}] if $(L\times\{L\})\cup\bigcup_{p\in L}(\{p\}\times\mathcal L_{\phi(p)})\subseteq\LB_\Pi$ for some line $L$ and some irreflexive involution $\phi:L\to L$;
\item[{$(\,|\,\circ)$}] if $\{p\}\times\mathcal L_p\subseteq\LB_\Pi$ for some point $p$;
\item[{$(\,|{=})$}] if $(\{p\}\times \mathcal L_p)\cup(L\times \{\Lambda\})\cup(\Lambda\times\{L\})\subseteq\LB_\Pi$ for some point $p$ and distinct lines $L,\Lambda\in \mathcal L_p$;
\item[{$(\,|\,/\!/)$}] if $(\{p\}\times\mathcal L_p)\cup\bigcup_{L\in \mathcal L_p}L\times \{\phi(L)\}\subseteq\LB_\Pi$ for some point $p$ and some  irreflexive involution $\phi:\mathcal L_p\to \mathcal L_p$;
\item[{$({+}\circ)$}]  if $(L\times\{L\})\cup(\{p\}\times\mathcal L_p)\subseteq\LB_\Pi$ for some line $L$ and point $p\in L$;
\item[{$(\circledast\circ)$}] if $\{(x,\Aline xp):x\in L\}\subseteq\LB_\Pi$ for some line $L$ and point $p\notin L$;
\item[{$(\circledast\bullet)$}] if $\{(x,\Aline xp):x\in L\}\cup\{(p,L)\}\subseteq\LB_\Pi$ for some line $L$ and point $p\notin L$;
\item[{$(\square\,\circ)$}] if $\{(p,L):p\in L\in\mathcal L_\Pi\}\subseteq\LB_\Pi$;
\item[{$(\square\blacksquare)$}] if $\LB_\Pi=\Pi\times\mathcal L_\Pi$.
\end{itemize}
\end{definition}

Definition~\ref{d:po-LB} actually involves a natural partial order  
on the $18$-element set of Lentz-Barlotti types 
$$\{(\circ\circ),\!(\circ\bullet),\!(\cdidots),\!({\circ}{\therefore}),\!(\circ\backslash),\!(\bullet\circ),\!(\bullet\bullet),\!({-}\circ),\!({-}\vert\vert),\!({-}/\!/),\!(\vert\circ),\!(\vert{=}),\!(\vert/\!/),\!(+\circ),\!(\circledast\circ),\!(\circledast\bullet),\!(\square\circ),\!(\square\blacksquare)\}.$$
This partial order is inherited from the partial orders of the product of the set of Lenz types
$$\{(\circ),(\bullet),(\circledast),(\,|\,),({-}),(+),(\square)\}$$ 
and the set of Barlotti types
$$\{(\circ),(\bullet),(\didots),(\therefore),(\,\backslash\,),(\,||\,),({=}),(/\!/),
(\blacksquare)\}.$$ The partial order on the set of Barlotti types is determined by the Hasse diagram

$$
\xymatrix@C=15pt@R=18pt{
&(\blacksquare)\\
&(/\!/)\ar@{-}[u]\\
(\,||\,)\ar@{-}[ru]&(\therefore)\ar@{-}[u]&({=})\ar@{-}[lu]&(\,\backslash\,)\ar@{-}[uull]\\
&(\didots)\ar@{-}[ul]\ar@{-}[u]\ar@{-}[ur]\\
&(\bullet)\ar@{-}[u]\ar@{-}[rruu]\\
&(\circ)\ar@{-}[u]
}
$$ 

The Hasse diagram of the partial order on the set of Lenz--Barlotti types looks as follows.
$$
\xymatrix@C=10pt{
&&(\square\blacksquare)\\
(\,\vert/\!/)\ar@{-}[urr]&{(\circ\backslash)}\ar@{-}[ur]&(\square\circ)\ar@{-}[u]&(\circledast\bullet)\ar@{-}[ul]&({-}/\!/)\ar@{-}[ull]\\
(\,\vert{=})\ar@{-}[u]&({\circ}{\therefore})\ar@{-}[ul]\ar@{-}[urrr]&(+\circ)\ar@{-}[u]&(\circledast\circ)\ar@{-}[u]&({-}||)\ar@{-}[u]\\
(\,\vert\circ)\ar@{-}[u]\ar@{-}\ar@{-}[urr]&(\cdidots)\ar@{-}[ul]\ar@{-}[u]\ar@{-}[urrr]&(\bullet\bullet)\ar@{-}[ull]\ar@{-}[urr]\ar@{-}[uur]&&({-}\circ)\ar@{-}[u]\ar@{-}[ull]\\
&(\circ\bullet)\ar@{-}[u]\ar@{-}[ur]\ar@/^20pt/@{-}[uuu]&(\bullet\circ)\ar@{-}[ull]\ar@{-}[u]\ar@{-}[uru]\ar@{-}[urr]\\
&(\circ\circ)\ar@{-}[u]\ar@{-}[ur]
}
$$

\begin{theorem}\label{t:LB>=} A projective plane $\Pi$ has Lenz--Barlotti type at least
\begin{itemize}
\item[$(\circ\bullet)$] iff some ternar of $\Pi$ is linear and associative-dot;
\item[$(\cdidots)$] iff some ternar of $\Pi$ is linear, left-distributive and associative-dot\newline iff some ternar of $\Pi$ is linear, right-distributive and associative-dot;
\item[$(\circ\!\therefore)$] iff some ternar of $\Pi$ is linear, distributive and associative-dot;
\item[$(\circ\backslash)$] iff some ternar of $\Pi$ is linear, associative-dot and diagonal-scale;
\item[$(\bullet\circ)$] iff some ternar of $\Pi$ is linear and associative-plus;
\item[$(\bullet\bullet)$] iff some ternar of $\Pi$ is linear and associative;
\item[$({-}\circ)$] iff some ternar of $\Pi$ is linear, right-distributive and associative-plus;
\item[$({-}||)$] iff some ternar of $\Pi$ is linear, right-distributive and associative;
\item[$(\,|\,\circ)$] iff some ternar of $\Pi$ is linear, left-distributive and associative-plus;
\item[$(\,|{=})$] iff some ternar of $\Pi$ is linear, left-distributive and associative;
\item[$(+\circ)$] iff some ternar of $\Pi$ is linear, distributive and associative-plus;
\item[$(\square\,\circ)$] iff some ternar of $\Pi$ is linear, distributive, associative-plus, and alternative-dot;
\item[$(\square\blacksquare)$] iff some ternar of $\Pi$ is linear, distributive and associative.
\end{itemize}
\end{theorem}

\begin{proof} In the subsequent proofs, for a projective base $uowe$ in the projective plane $\Pi$ we denote by $h\in\Aline ou\cap\Aline we$ and $v\in\Aline ow\cap\Aline ue$ the horizontal and vertical infinity point of the projective base $uowe$, respectively.
\smallskip

 $(\circ\bullet)$ If the projective plane $\Pi$ has Lenz--Barlotti type at least $(\circ\bullet)$, then it is $(p,L)$-transitive for some line $L\in\mathcal L_\Pi$ and point $p\in\Pi\setminus L$. Choose a projective base $uowe$ in $\Pi$ such that $p=h$ and $L=\Aline ov$. The $(p,L)$-transitivity of $\Pi$ implies the $(h,\Aline ov)$-transitivity of $\Pi$. By Theorem~\ref{t:Algebra-vs-Geometry-proj}(5), the ternar of the based projective plane $(\Pi,uowe)$ is linear and associative-dot.

Now assume that some ternar of $\Pi$ is linear and associative-dot, which means that so is the ternar of some projective base $uowe$ for $\Pi$. By Theorem~\ref{t:Algebra-vs-Geometry-proj}(5), the projective plane $\Pi$ is $(h,\Aline ov)$-transitive. Since $(h,\Aline ov)\in\LB_\Pi$, the projective plane $\Pi$ has Lenz--Barlotti class at least $(\circ\bullet)$.
\smallskip

$(\cdidots)$ If the projective plane $\Pi$ has Lenz--Barlotti type at least $(\cdidots)$, then $\{(p,L),(q,\Lambda)\}\subseteq\LB_\Pi$ for some lines $L,\Lambda$ and points $p\in\Lambda\setminus L$ and $q\in L\setminus\Lambda$.  Choose a projective base $uowe$ in $\Pi$ such that $\Aline ow=L$, $\Aline ou=\Lambda$, $p=h$ and $q=v$. The inclusion $\{(h,\Aline ov),(v,\Aline oh)\}=\{(p,L),(q,\Lambda)\}\subseteq \LB_\Pi$ implies that the projective plane $\Pi$ is $(h,\Aline ov)$-transitive and $(v,\Aline oh)$-transitive.  By Theorem~\ref{t:Algebra-vs-Geometry-proj}(6), the ternar of the based projective plane $(\Pi,uowe)$ is linear, right-distributive, and associative-dot. 

We can also choose a projective base $uowe$ so that $p=h$, $q=o$, $L=\Aline ov$ and $\Lambda=\Aline hv$. In this case, the inclusion $\{(h,\Aline ov),(o,\Aline oh)\}=\{(p,L),(q,\Lambda)\}\subseteq \LB_\Pi$ implies that the projective plane $\Pi$ is $(h,\Aline ov)$-transitive and $(o,\Aline hv)$-transitive.  By Theorem~\ref{t:Algebra-vs-Geometry-proj}(7), the ternar of the based projective plane $(\Pi,uowe)$ is linear, left-distributive, and associative-dot. 

Now assume that some ternar of $\Pi$ is linear, right-distributive and associative-dot, which means that so is the ternar of some projective base $uowe$ for $\Pi$. By Theorem~\ref{t:Algebra-vs-Geometry-proj}(6), the projective plane $\Pi$ is $(h,\Aline ov)$-transitive and $(v,\Aline oh)$-transitive. Since $\{(h,\Aline ov),(v,\Aline oh)\}\subseteq\LB_\Pi$, the projective plane $\Pi$ has Lenz--Barlotti class at least $(\cdidots)$.

If some ternar of $\Pi$ is linear, left-distributive and associative-dot, then by Theorem~\ref{t:Algebra-vs-Geometry-proj}(7), $\Pi$ has a projective base $uowe$ such that  $\Pi$ is $(h,\Aline ov)$-transitive and $(o,\Aline hv)$-transitive. Since $\{(h,\Aline ov),(o,\Aline hv)\}\subseteq\LB_\Pi$, the projective plane $\Pi$ has Lenz--Barlotti class at least $(\cdidots)$.
\smallskip

$(\circ\!\therefore)$ If the projective plane $\Pi$ has Lenz--Barlotti type at least $(\circ\!\therefore)$, then $$\{(a,\Aline bc),(b,\Aline ac),(c,\Aline ab)\}\subseteq\LB_\Pi$$ for some non-collinear points $a,b,c$ in $\Pi$. Choose a projective base $uowe$ in $\Pi$ such that $a=h$, $b=v$, $c=o$. The inclusion $\{(h,\Aline ov),(v,\Aline oh),(o,\Aline hv)\}=\{(a,\Aline bc),(b,\Aline ac),(c,\Aline ab)\}\subseteq\LB_\Pi$ implies that the projective plane $\Pi$ is $(h,\Aline ov)$-transitive, $(v,\Aline oh)$-transitive, and $(o,\Aline hv)$-transitive.  By Theorem~\ref{t:Algebra-vs-Geometry-proj}(8), the ternar of the based projective plane $(\Pi,uowe)$ is linear, distributive, and associative-dot.

Now assume that some ternar of $\Pi$ is linear, distributive and associative-dot, which means that so is the ternar of some projective base $uowe$ for $\Pi$. By Theorem~\ref{t:Algebra-vs-Geometry-proj}(8), the projective plane $\Pi$ is $(h,\Aline ov)$-transitive, $(v,\Aline oh)$-transitive, and $(o,\Aline hv)$-transitive. Since $$\{(h,\Aline ov),(v,\Aline oh),(o,\Aline hv)\}\subseteq\LB_\Pi,$$ the projective plane $\Pi$ has Lenz--Barlotti class at least $(\circ\!\therefore)$.
\smallskip

$(\circ\backslash)$ If the projective plane $\Pi$ has Lenz--Barlotti type at least $(\circ\backslash)$, then $F\subseteq \LB_\Pi$ for some line $L\in\mathcal L_\Pi$, point $p\in L$ and a bijective function $F:L\setminus\{p\}\to\mathcal L_p\setminus\{L\}$.  Choose two distinct pairs $(h,V),(d,E)\in F$. Since $F$ is a bijective function, $h\ne v$ and $V\ne E$. Choose any point $o\in V\setminus \{p\}$ and find unique points $e\in E\cap\Aline od$, $u\in E\cap\Aline oh$ and $w\in V\cap\Aline eh$. Then $uowe$ is a projective base whose horizon line $\Aline hv$ coincides with the line $L$ and the vertical infinity point coincides with the point $p$. Since $\{(h,\Aline ow),(d,\Aline ue)\}=\{(h,V),(d,E)\}\subseteq \LB_\Pi$, the based affine plane $(\Pi\setminus L,uow)$ is horizontal-scale and diagonal-scale. By Theorem~\ref{t:diagonal-scale<=>}, its ternar $\Delta=\Aline oe\setminus\{d\}$ is linear associative-dot and diagonal-scale. Then the projective plane $\Pi$ has a linear associative-dot diagonal-scale ternar $\Delta$.

Now assume that some ternar $R$ of the projective plane $\Pi$ is linear, associative-dot and diagonal-scale. Then there exists a projective base $uowe$ whose ternar $\Delta$ is linear, associative-dot and diagonal-scale. Let $h\in\Aline ou\cap\Aline we$, $v\in\Aline ow\cap\Aline ue$, and $d\in\Aline oe\cap\Aline hv$ be the horizontal, vertical and diagonal infinite points of the projective base $uowe$. By Theorem~\ref{t:diagonal-scale<=>}, the based affine plane $(\Pi\setminus\Aline hv,uow)$ is horizontal-scale and diagonal-scale, which implies that the projective plane $\Pi$ is $(h,\Aline vo)$-transitive and $(d,\Aline ve)$-transitive. Since $\Pi$ is $(h,\Aline vo)$-transitive, for every point $x\in \Aline ou\setminus \{o,h\}$, there exists an automorphism $\Phi_x$ of the projective plane $\Pi$ such that $\Phi_x(u)=x$ and $\Aline ov\subseteq\Fix(\Phi_x)$. By Proposition~\ref{p:pL-Auto}, the pair $(\Phi_x(d),\Phi[\Aline ev])$ belongs to the Lenz--Barlotti figure $\LB_\Pi$. It can be shown that $F\defeq\{(h,\Aline vo)\}\cup\{(\Phi_x(d),\Phi_x[\Aline ve]):x\in \Aline oh\setminus\{o,h\}\}$ is an injective function with $\dom[F]=\partial \Pi\setminus\{v\}$ and $\rng[F]=\mathcal L_v\setminus\{\Aline hv\}$. Since $F\subseteq\LB_\Pi$, the projective plane is of Lenz--Barlotti type at least $(\circ\backslash)$.
\smallskip

$(\bullet\circ)$ If the projective plane $\Pi$ has Lenz--Barlotti type at least $(\bullet\circ)$, then $\Pi$ is $(p,L)$-transitive for some line $L$ and point $p\in L$. Then $\Pi$ is $(v,\Aline ov)$-transitive, and  by Theorem~\ref{t:Algebra-vs-Geometry-proj}(1), the ternar of the based projective plane $(\Pi,uowe)$ is linear and associative-plus.

Now assume that some ternar of $\Pi$ is linear and associative-plus, which means that so is the ternar of some projective base $uowe$ for $\Pi$. By Theorem~\ref{t:Algebra-vs-Geometry-proj}(1), the projective plane $\Pi$ is $(v,\Aline ov)$-transitive. Since $(v,\Aline ov)\in\LB_\Pi$, the projective plane $\Pi$ has Lenz--Barlotti class at least $(\bullet\circ)$.
\smallskip

$(\bullet\bullet)$ If the projective plane $\Pi$ has Lenz--Barlotti type at least $(\bullet\bullet)$, then $\{(p,L),(q,\Lambda)\}\subseteq\LB_\Pi$ for some lines $L,\Lambda$ and points $p\in L\cap\Lambda$ and $q\in L\setminus\Lambda$.  Choose a projective base $uowe$ in $\Pi$ such that $L=\Aline hv$, $\Lambda=\Aline ov$, $v=p$, and $h=q$. The inclusion $\{(v,\Aline hv),(h,\Aline ov)\}=\{(p,L),(q,\Lambda)\}\subseteq\LB_\Pi$ implies that the projective plane $\Pi$ is $(v,\Aline hv)$-transitive and $(h,\Aline ov)$-transitive. By Theorem~\ref{t:Algebra-vs-Geometry-proj}(9), the ternar of the based projective plane $(\Pi,uowe)$ is linear and associative.

Now assume that some ternar of $\Pi$ is linear and associative, which means that so is the ternar of some projective base $uowe$ for $\Pi$. By Theorem~\ref{t:Algebra-vs-Geometry-proj}(9), the projective plane $\Pi$ is $(v,\Aline ov)$-transitive and $(h,\Aline ov)$-transitive. Since $\{(v,\Aline hv),(h,\Aline ov)\}\subseteq\LB_\Pi$, the projective plane $\Pi$ has Lenz--Barlotti class at least $(\bullet\bullet)$.
\smallskip

$({-}\circ)$ If the projective plane $\Pi$ has Lenz--Barlotti type at least $({-}\circ)$, then $L\times\{L\}\subseteq\LB_\Pi$ for some line $L$. Choose a projective base $uowe$ in $\Pi$ such that $L=\Aline hv$. The inclusion $\Aline hv\times\{\Aline hv\}=L\times\{L\}\subseteq\LB_\Pi$ implies that the projective plane $\Pi$ is $[\Aline hv,\Aline hv)$-transitive. By Theorem~\ref{t:Algebra-vs-Geometry-proj}(2), the ternar of the based projective plane $(\Pi,uowe)$ is linear, right-distributive, and associative-plus.

Now assume that some ternar of $\Pi$ is linear, right-distributive and associative-plus, which means that so is the ternar of some projective base $uowe$ for $\Pi$. By Theorem~\ref{t:Algebra-vs-Geometry-proj}(2), the projective plane $\Pi$ is $[\Aline hv,\Aline hv)$-transitive and hence $\Aline hv\times\{\Aline hv\}\subseteq\LB_\Pi$, which means that the projective plane $\Pi$ has Lenz--Barlotti class at least $({-}\circ)$.
\smallskip

$({-}||)$ If the projective plane $\Pi$ has Lenz--Barlotti type at least $({-}||)$, then $$(L\times\{L\})\cup(\{p\}\times\mathcal L_q)\cup(\{q\}\times\mathcal L_p)\subseteq\LB_\Pi$$ for some line $L$ and distinct points $p,q\in L$. Choose a projective base $uowe$ in $\Pi$ such that $h=p$and $v=q$. The inclusion 
$\{v\}\times\mathcal L_h=\{q\}\times\mathcal L_p\subseteq\LB_\Pi$ implies that the projective plane $\Pi$ is $(v,h]$-transitive. By Theorem~\ref{t:Algebra-vs-Geometry-proj}(11), the ternar of the based projective plane $(\Pi,uowe)$ is linear, right-distributive and associative.

Now assume that some ternar of $\Pi$ is linear, right-distributive and associative, which means that so is the ternar of some projective base $uowe$ for $\Pi$. By Theorem~\ref{t:Algebra-vs-Geometry-proj}(2,11), the projective plane $\Pi$ is $[\Aline hv,\Aline hv)$-transitive, $(v,h]$-transitive, $(h,v]$-transitive, and hence $$(\Aline hv\times\{\Aline hv\})\cup(\{h\}\times\mathcal L_v)\cup(\{v\}\times\mathcal L_h)\subseteq\LB_\Pi,$$ which means that the projective plane $\Pi$ has Lenz--Barlotti class at least $({-}||)$.
\smallskip

$(\,|\circ)$ If the projective plane $\Pi$ has Lenz--Barlotti type at least $(\,|\circ)$, then $\{p\}\times\mathcal L_p\subseteq\LB_\Pi$ for some point $p\in\Pi$. Choose a projective base $uowe$ in $\Pi$ such that $v=p$. Then $\{v\}\times\mathcal L_v=\{p\}\times\mathcal L_p\subseteq\LB_\Pi$, and by Theorem~\ref{t:Algebra-vs-Geometry-proj}(3), the ternar of the based projective plane $(\Pi,uowe)$ is linear, left-distributive, and associative-plus.

Now assume that some ternar of $\Pi$ is linear, left-distributive and associative-plus, which means that so is the ternar of some projective base $uowe$ for $\Pi$. By Theorem~\ref{t:Algebra-vs-Geometry-proj}(3), the projective plane $\Pi$ is $(v,v]$-transitive and hence $\{v\}\times\mathcal L_v\subseteq\LB_\Pi$, which means that the projective plane $\Pi$ has Lenz--Barlotti class at least $(\,|\circ)$.
\smallskip

$(\,|{=})$ If the projective plane $\Pi$ has Lenz--Barlotti type at least $(\,|{=})$, then $$(\{p\}\times\mathcal L_p)\cup(L\times\{\Lambda\})\cup(\Lambda\times\{L\})\subseteq\LB_\Pi$$ for some point $p\in\Pi$ and distinct lines $L,\Lambda\in \mathcal L_p$. Choose a projective base $uowe$ in $\Pi$ such that $v=p$, $\Aline ov=L$ and $\Aline hv=\Lambda$. The inclusion 
$\Aline ov\times\{\Aline hv\}=L\times\{\Lambda\}\subseteq\LB_\Pi$ implies that the projective plane $\Pi$ is $[\Aline ov,\Aline hv)$-transitive. By Theorem~\ref{t:Algebra-vs-Geometry-proj}(10), the ternar of the based projective plane $(\Pi,uowe)$ is linear, left-distributive and associative.

Now assume that some ternar of $\Pi$ is linear, left-distributive and associative, which means that so is the ternar of some projective base $uowe$ for $\Pi$. By Theorem~\ref{t:Algebra-vs-Geometry-proj}(3,10), the projective plane $\Pi$ is $(v,v]$-transitive, $[\Aline ov,\Aline hv)$-transitive, $[\Aline hv,\Aline ov)$-transitive, and $(h,v]$-transitive. Then $$(\{v\}\times\mathcal L_v)\cup (\Aline ov\times\{\Aline hv\})\cup(\Aline hv\times\{\Aline ov\})\subseteq\LB_\Pi,$$ which means that the projective plane $\Pi$ has Lenz--Barlotti class at least $(\,|\,{=})$.
\smallskip

$(+\circ)$ If the projective plane $\Pi$ has Lenz--Barlotti type at least $(+\circ)$, then $(L\times\{L\})\cup(\{p\}\times\mathcal L_p)\subseteq\LB_\Pi$ for some line $L$ and point $p\in L$. Choose a projective base $uowe$ in $\Pi$ such that $\Aline hv=L$ and $v=p$. The inclusion 
$(\Aline hv\times\{\Aline hv\})\cup(\{v\}\times\mathcal L_v)=(L\times\{L\})\cup(\{p\}\times\mathcal L_p)\subseteq\LB_\Pi$ implies that the projective plane $\Pi$ is $[\Aline hv,\Aline hv)$-transitive and $(v,v]$-transitive. By Theorem~\ref{t:Algebra-vs-Geometry-proj}(4), the ternar of the based projective plane $(\Pi,uowe)$ is linear, distributive and associative-plus.

Now assume that some ternar of $\Pi$ is linear, distributive and associative-plus, which means that so is the ternar of some projective base $uowe$ for $\Pi$. By Theorem~\ref{t:Algebra-vs-Geometry-proj}(4), the projective plane $\Pi$ is $[\Aline hv,\Aline hv)$-transitive and $(v,v]$-transitive. Then $(\Aline hv\times\{\Aline hv\})\cup (\{v\}\times\mathcal L_v)\subseteq\LB_\Pi$, which means that the projective plane $\Pi$ has Lenz--Barlotti class at least $(+\circ)$.
\smallskip

$(\square\,\circ)$ If the projective plane $\Pi$ has Lenz--Barlotti type at least $(\square\circ)$, then $\{(p,L)\in\Pi\times\mathcal L_\Pi:p\in L\}\subseteq\LB_\Pi$.  By Theorem~\ref{t:Algebra-vs-Geometry-proj}(14 or 15), the ternar of any projective base in $\Pi$ is linear, distributive, associative-plus, and alternative-dot.

Now assume that some ternar of $\Pi$ is linear, distributive, associative-plus, and alternative-dot. By Theorem~\ref{t:Skornyakov-San-Soucie}, every line in $\Pi$ is translation and hence $\{(p,L)\in\Pi\times\mathcal L_\Pi:p\in L\}\subseteq\LB_\Pi$, which means that the projective plane $\Pi$ has Lenz--Barlotti class at least $(\square\circ)$.
\smallskip

$(\square\blacksquare)$ If the projective plane $\Pi$ has Lenz--Barlotti type at least $(\square\blacksquare)$, then $\LB_\Pi=\Pi\times\mathcal L_\Pi$.  By Theorem~\ref{t:Algebra-vs-Geometry-proj}(12), the ternar of any projective base in $\Pi$ is linear, distributive and associative.

Assume that some ternar of $\Pi$ is linear, distributive, and associative. By Theorem~\ref{t:Algebra-vs-Geometry-proj}(12), $\LB_\Pi=\Pi\times\mathcal L_\Pi$, which means that the projective plane $\Pi$ has Lenz--Barlotti class at least $(\square\blacksquare)$.
\end{proof}

\begin{corollary}\label{c:Desarg-self-dual} A projective plane is Desarguesian if and only if its dual projective plane is Desarguesian.
\end{corollary}

\begin{proof} If a projective plane $X$ is Desarguesian, then any ternar $R$ of $X$ is linear, distributive, and associative, by Theorem~\ref{t:Algebra-vs-Geometry-proj}(12). By Theorem~\ref{t:LB>=}, $X$ has Lenz-Barlotti type $(\square\blacksquare)$. The definition of Lenz--Barlotti figure implies that the dual projective plane $X^*$ to $X$ also has Lenz-Barlotti type $(\square\blacksquare)$. By Theorem~\ref{t:LB>=}, some ternar of the projective plane $X^*$ is linear, distributive, and associative. By Theorem~\ref{t:Algebra-vs-Geometry-proj}(12), the projective plane $X^*$ is Desarguesian. By analogy we can prove that the Desarguesianity of the dual projective plane $X^*$ implies the Desarguesianity of $X$.
\end{proof}


\begin{remark} In standard textbooks in Projective Geometry, the 18 cases 
$$(\circ\circ),(\circ\bullet),(\cdidots),({\circ}{\therefore}),({\circ}{\backslash}),(\bullet\circ),
(\bullet\bullet),(\circledast\circ),(\circledast\bullet),({-}\circ),({-}\vert\vert),({-}/\!/),(\vert\circ),(\vert{=}),(\vert/\!/),(+\circ),(\square\circ),(\square\blacksquare)$$ in the Lenz--Barlotti classification are usually numbered by the roman.arabic numerals 
$$\mbox{I.1,\;I.2,\;I.3,\;I.4,\;I.6,\;II.1,\;II.2,\;III.1;\,III.2;\,IVa.1;\,IVa.2,\;IVa.3,\;IVb.1,\;IVb.2,\;IVb.3,\;V.1,\;VII.1,\;VII.2},$$ respectively. Theorem~\ref{t:Lenz-Barlotti} shows that all other types, initially mentioned by Barlotti \cite{Barlotti1957} do not occur among projective planes. 
\end{remark}

\begin{remark} According to the survey \cite{LB}, there exist projective planes of all Lenz--Barlotti types, with a possible exception of class I.6, i.e., $(\circ\backslash)$. By Theorem~\ref{t:Lenz-Barlotti}, there is only one plane of the Lenz--Barlotti type $(-/\!/)$ IV.a.3 (resp. $(\,|\,/\!/)$ IV.b.3) and it is the right-Hall (resp. left-Hall) projective plane of order 9. So, among projective planes of order $\ne 9$ the Lenz-Barlotti types $(-/\!/)$ and $(\,|\,/\!/)$ (i.e. IV.a.3 and IV.b.3) do not appear. 
The Lenz-Barlotti types $(\circ\backslash)$, $(\circledast\circ)$, $(\circledast\bullet)$, $(\square\circ)$ (numbered as I.6, III.1, III.2, VII.1) do not appear among finite projective planes. The existence of a projective plane of Lenz--Barlotti type $(\circ\backslash)$ and {\em finite} projective planes of types $(\circ\bullet)$, $(\cdidots)$, $(\circ\!\therefore)$ and $(\bullet\bullet)$ (numbered as II.2, I.3, I.4 and II.2) are long-standing open problems.
\end{remark}

\begin{remark} The Lenz--Barlotti figures of projective planes can be also classified by the ranks of the domain and range of their Lenz and Barlotti parts. More precisely, for a projective plane $\Pi$, the \index{tetrarank}\index{Lenz--Barlotti figure!tetrarank of}\defterm{tetrarank} of its Lenz--Barlotti figure $\LB_\Pi$ is the quadruple $dr\delta\rho$ of the ranks 
\vskip3pt

\centerline{$d\defeq\|\dom[\Lenz_\Pi]\|,\quad r\defeq\|\rng[\Lenz_\Pi]\|,\quad\delta\defeq\|\dom[\LB_\Pi\setminus\Lenz_\Pi]\|,\quad\rho\defeq\|\rng[\LB_\Pi\setminus\Lenz_\Pi]\|$}
\vskip3pt

\noindent of the projections of the Lenz figure $\Lenz_\Pi$ and the Barlotti figure $\LB_\Pi\setminus\Lenz_\Pi$ onto the projective plane $\Pi$ and its dual projective plane $\mathcal L_\Pi$. 

The Lenz--Barlotti figures of projective planes of Lenz--Barlotti type 
\vskip3pt

\centerline{$(\circ\circ),(\circ\bullet),(\cdidots),({\circ}{\therefore}),(\circ\backslash),(\bullet\circ),
(\bullet\bullet),({-}\circ),({-}\vert\vert),({-}/\!/),(\vert\circ),(\vert{=}),(\vert/\!/),(+\circ),(\circledast\circ),(\circledast\bullet),(\square\circ),(\square\blacksquare)$}
\vskip3pt

\noindent
have tetraranks, respectively:
\vskip2pt

\centerline{$\scalebox{0.97}{0000,\,0011,\,0022,\,0033,\,0022,\,1100,\,1111,\,2100,\,2123,\,2123,\,1200,\,1232,\,1232,\,2200,\,2200,\,2211,\,3300,\,3333.}$}
\end{remark}

\begin{remark} The projective planes of order 9 have the following Lenz--Barlotti types and tetraranks:
$$
\begin{array}{|l|c|c|c|c|}
\hline
\mbox{Projective plane:}&\mbox{Desarg}&\mbox{right-Hall}&\mbox{left-Hall}&\mbox{Hughes}\\
\hline
\mbox{Lenz--Barlotti type:}&(\square\blacksquare)&(-/\!/)&(\,|/\!/)&(\circ\circ)\\
\mbox{in roman.arabic:}&$VII.2$&$IV.a.3$&$IV.b.3$&$I.1$\\
\mbox{tetraranks:}&3333&2123&1232&0000\\
\hline
\end{array}
$$
\end{remark}
 
\begin{remark} All 22 known projective planes of order 16 \cite{Moorhouse} have the following Lenz--Barlotti types and tetraranks (calculated by Ivan Hetman).
$$
\begin{array}{|c|c|c|c|}
\hline
\mbox{\tt Projective plane}&\mbox{Lenz--Barlotti type}&\mbox{in roman.arabic}&\mbox{tetrarank}\\
\hline
\mbox{\tt bbh1,john,bbs4,dbbs4,bhh2,dbbh2}&(\circ\circ)&$I.1$&0000\\
\mbox{\tt math,dmath}&(\bullet\circ)&$II.1$&1100\\
\mbox{\tt hall, jowk,demp,lmrh,dsfp}&(\,|\,\circ)&$IV.b.1$&1200\\
\mbox{\tt dhall,djowk,ddemp,dlmrh,ddsfp}&({-}\circ)&$IV.a.1$&2100\\
\mbox{\tt semi2,semi4}&(+\circ)&$V.1$&2200\\
\mbox{\tt desarg}&(\square\blacksquare)&$VII.2$&3333\\
\hline
\end{array}
$$
\end{remark}

\begin{remark} The Lenz--Barlotti types of all 193 known projective planes of order 25 \cite{Moorhouse} were calculated by Ivan Hetman. The results of calculations are summarized in the following table displaying the (nonzero) number of projective planes of order 25 of a given Lenz--Barlotti type.
$$
\begin{array}{|c|c|c|c|}
\hline
\mbox{Number of projective}\atop\mbox{planes of order 25}&\mbox{Lenz--Barlotti type}&\mbox{in roman.arabic}&\mbox{tetrarank}\\
\hline
144&(\circ\circ)&$I.1$&0000\\
8&(\bullet\circ)&$II.1$&1100\\
18&(\,|\,\circ)&$IV.a.1$&1200\\
18&({-}\circ)&$IV.b.1$&2100\\
2&(|{=})&$IV.a.2$&1223\\
2&({-}||)&$IV.b.2$&2123\\
1&(\square\blacksquare)&$VII.2$&3333\\
\hline
\end{array}
$$
\end{remark}

\begin{remark} All 13 known projective planes of order 27 \cite{Moorhouse} have the following Lenz--Barlotti types and tetraranks (calculated by Ivan Hetman).
$$
\begin{array}{|c|c|c|c|}
\hline
\mbox{\tt Projective plane}&\mbox{Lenz--Barlotti type}&\mbox{in roman.arabic}&\mbox{tetrarank}\\
\hline
\mbox{\tt fig}&(\circ\circ)&$I.1$&0000\\
\mbox{\tt andre,hering,sherk,flag4,flag6,}&(\,|\,\circ)&$IV.a.1$&1200\\
\mbox{\tt dandre,dhering,dsherk,dflag4,dflag6}&({-}\circ)&$IV.b.1$&2100\\
\mbox{\tt twisted}&(+\circ)&$V.1$&2200\\
\mbox{\tt desarg}&(\square\blacksquare)&$VII.2$&3333\\
\hline
\end{array}
$$
\end{remark}

\begin{remark} All 12 known projective planes of order 32 \cite{Moorhouse} have the following Lenz--Barlotti types and tetraranks (calculated by Ivan Hetman).
$$
\begin{array}{|c|c|c|c|}
\hline
\mbox{\tt Projective plane}&\mbox{Lenz--Barlotti type}&\mbox{in roman.arabic}&\mbox{tetrarank}\\
\hline
\mbox{\tt e,f,g}&(\,|\,\circ)&$IV.a.1$&1200\\
\mbox{\tt eD,fD,gD}&({-}\circ)&$IV.b.1$&2100\\
\mbox{\tt a,aD,b,bD,d}&(+\circ)&$V.1$&2200\\
\mbox{\tt desarg}&(\square\blacksquare)&$VII.2$&3333\\
\hline
\end{array}
$$
\end{remark}




\section{The central matrix of a projective plane}

\begin{definition} For a projective plane $\Pi$, the indexed family of groups 
$$\big(\Aut_{p,L}(\Pi)\big)_{(p,L)\in\Pi\times\mathcal L_\Pi}$$is called the \index{autocentral matrix}\index{projective plane!autocentral matrix}\defterm{autocentral matrix} of $\Pi$, and the indexed family of cardinals
$$\big(|\Aut_{p,L}(\Pi)|\big)_{(p,L)\in\Pi\times\mathcal L_\Pi}$$
is called the \index{central matrix}\index{projective plane!central matrix}\defterm{central matrix} of $\Pi$. 

The central matrix is uniquely determined by the indexed family $(\overline *_n)_{n\in|\Aut_{**}(\Pi)|}$ of its \index{level figure}\defterm{level figures}
$$\overline *_n\defeq\{(p,L)\in\Pi\times\mathcal L_\Pi:|\Aut_{p,L}(\Pi)|=n\},$$
indexed by the numbers in the set
$$|\Aut_{**}(\Pi)|\defeq\{|\Aut_{p,L}(\Pi)|:(p,L)\in\Pi\times\mathcal L_\Pi\},$$
called the \index{central spectrum}\index{projective plane!central spectrum}\defterm{central spectrum} of the projective plane $\Pi$.
\end{definition}

\begin{proposition}\label{p:Aut(p,L)-action} For every point-line pair $(p,L)$ in a projective plane $\Pi$ and every point $x\in \Pi\setminus(\{p\}\cup L)$, the map $\Aut_{p,L}(\Pi)\to \Pi$, $A\mapsto A(x)$, is injective.
\end{proposition}

\begin{proof} Assume that $F,G\in \Aut_{p,L}(\Pi)$ are two automorphism such that $F(x)=G(x)$. Then $G^{-1}F\in\Aut_{p,L}(\Pi)$ is an automorphism of $\Pi$ such that $L\cup\{p,x\}\subseteq\Fix(G^{-1}F)$ and hence $G^{-1}F$ is an identity automorphism of $\Pi$, by Proposition~\ref{p:Fix(A)=L+p}. Then $F=G$, witnessing that the map $\Aut_{p,L}(\Pi)\to X$, $A\mapsto A(x)$, is injective.
\end{proof}

\begin{corollary}\label{c:|Aut|-divides} Let $\Pi$ be a finite projective plane and $(p,L)$ be a point-line pair in $\Pi$.
\begin{enumerate}
\item If $p\in L$, then $|\Aut_{p,L}(\Pi)|$ divides the number $|\Pi|_2-1$;
\item If $p\notin L$, then $|\Aut_{p,L}(\Pi)|$ divides the number $|\Pi|_2-2$.
\end{enumerate}
\end{corollary}

\begin{proof} By Proposition~\ref{p:Aut(p,L)-action}, for every $x\in \Pi\setminus(\{p\}\cup L)$, the set $\{F(x):F\in\Aut_{p,L}(\Pi)\}$ has cardinality $|\Aut_{p,L}(\Pi)|$. Fix any point $a\in \Pi\setminus(\{p\}\cup L)$ and observe that  
$$\big\{\{F(x):F\in\Aut_{p,L}(\Pi)\}:x\in \Aline pa\setminus(\{p\}\cup L)\big\}$$ is a partition of the set $\Aline pa\setminus(\{p\}\cup L)$ into pairwise disjoint sets of cardinality $|\Aut_{p,L}(\Pi)|$, which implies that $|\Aut_{p,L})\Pi)|$ divides $|\Aline pa\setminus(\{p\}\cup L)|$. 

If $p\in L$, then $|\Aut_{p,L}(\Pi)|$ divides $|\Aline pa\setminus(\{p\}\cup L)|=|\Pi|_2-1$. 

If $p\notin L$, then $|\Aut_{p,L}(\Pi)|$ divides $|\Aline pa\setminus(\{p\}\cup L)|=|\Pi|_2-2$. 
\end{proof}

\begin{corollary}
The central spectrum of a finite projective plane $\Pi$ is a subset of divisors of the numbers $|\Pi|_2-1$ or $|\Pi|_2-2$.
\end{corollary}

\begin{proposition} A finite projective plane $\Pi$ is Pappian if and only if its central spectrum coincides with the set $\{|\Pi|_2-1,|\Pi|_2-2\}$. 
\end{proposition}

\begin{remark} The Lenz figure $\Lenz_\Pi$ and the Lenz--Barlotti figure $\LB_\Pi$ of a finite projective plane $\Pi$ can be recovered from the central matrix of $\Pi$ because they are the highest level figures:
$$\Lenz_\Pi=\overline *_{|\Pi|_2-1}\quad\mbox{and}\quad\LB_\Pi=\overline *_{|\Pi|_2-1}\cup\overline *_{|\Pi|_2-2}.$$
\end{remark}

Since the central matrix is too large for handling manually, it is convenient to have some ``smaller'' characteristics of this matrix, which are introduced in the following definition.

\begin{definition} For a projective plane $\Pi$ its
\begin{itemize}
\item \index{central levelsize}\index{projective plane!central levelsize}\defterm{central levelsize} is the indexed family of cardinals $(|\overline *_n|)_{n\in|\Aut_{**}(\Pi)|\setminus\{1\}};$
\item \index{central levelrank}\index{projective plane!central multirank}\defterm{central levelrank} is the indexed family of pairs of ranks $(\|\dom[\overline *_n]\|,\|\rng[\overline *_n]\|)_{n\in|\Aut_{**}(\Pi)|\setminus\{1\}}.$
\end{itemize}
\end{definition}

\begin{example} A Pappian projective plane $\Pi$ of order $n$ has central spectrum $\{n,n-1\}$, central levelsize $$\big(n^3+2n^2+2n+1)_n(n^4+n^3+n^2)_{n-1}$$ 
and central levelrank $$(3,3)_n(3,3)_{n-1},$$ which will be written as $33_n33_{n-1}$.
\end{example} 

\begin{examples} The central spectra, central levelsizes and central levelranks of all projective planes of order $9$ are written in the following table (calculated by Ivan Hetman).
$$
\begin{array}{|l|c|c|c|c|}
\hline
\mbox{Projective plane:}&\mbox{Desargues}&\mbox{right-Hall}&\mbox{left-Hall}&\mbox{Hughes}\\
\hline
\mbox{Lenz--Barlotti type}&(\square\blacksquare)&(-/\!/)&(\,|/\!/)&(\circ\circ)\\
\mbox{Central spectrum}&9,8&9,8,2,1&9,8,2,1&3,2,1\\
\mbox{Central levelsize:}&910_9  7371_8&10_990_881_2&10_9 90_8  81_2&52_3 117_2\\
\mbox{Central levelrank:}&33_933_8&12_932_813_2&21_923_831_2&33_333_2\\
\hline
\end{array}
$$
\end{examples}

\begin{exercise} Show that the computation complexity of calculating the central matrix  of a projective plane of order $n$ equals $O(n^8\ln n)$.
\end{exercise}


\begin{definition} For a projective plane $\Pi$, its \index{Lenz graph}\index{projective plane!Lenz graph}\defterm{Lenz graph} $\LG(\Pi)$ has the set $$V_\Pi\defeq\{(p,L)\in \Pi\times\mathcal L_\Pi:p\in L\;\wedge\;|\Aut_{p,L}(\Pi)|>1\}$$ as the set of vertices and the set 
$$E_\Pi\defeq\{\{(p,L),(q,\Lambda)\}\in [V]^2: p=q\;\mbox{or}\;L=\Lambda\}$$as the set of edges. We say that the Lenz graph $\LG(\Pi)$ of a projective plane $\Pi$ is 
\begin{itemize}
\item \index{Lenz graph!line-wide}\defterm{line-wide} if for any line $L\in\mathcal L_\Pi$, the set $\{x\in L:(x,L)\in V_\Pi\}$ contains at least two distinct points;
\item \index{Lenz graph!point-wide}\defterm{point-wide} if for any point $p\in\mathcal L_\Pi$, the set $\{L\in\mathcal L_p:(p,L)\in V_\Pi\}$ contains at least two distinct lines;
\item \index{Lenz graph!wide}\defterm{wide} if $\LG(\Pi)$ is line-wide and point-wide;
\item \index{Lenz graph!connected}\defterm{connected} if for any $(p,L),(q,\Lambda)\in V_\Pi$, there exists a sequence\newline $(p_0,L_0),(p_1,L_1),\dots,(p_n,L_n)\in V_\Pi$ such that $(p_0,L_0)=(p,L)$, $(p_n,L_n)=(q,\Lambda)$,\newline and $\{(p_{i-1},L_{i-1}),(p_i,L_i)\}\in E_\Pi$ for all $i\in\{1,\dots,n\}$.
\end{itemize}
\end{definition}

The following theorem of Gleason implies that the autocentral matrix of a non-Desargueasian finite projective plane necessarily contains trivial groups.

\index[person]{Gleason}
\begin{theorem}[Gleason, 1956]\label{t:Gleason56} For any finite projective plane $\Pi$ of order $n\ge 4$, the following conditions are equivalent:
\begin{enumerate}
\item $\Pi$ is Pappian;
\item $\Pi$ is Desarguesian;
\item $|\Aut_{**}(\Pi)|=\{n,n-1\}$;
\item $\overline*_{n-1}=\{(p,L)\in\Pi\times\mathcal L_\Pi:p\notin L\}$;
\item $|\overline *_{n-1}|=(n^2+n+1)\cdot n^2$;
\item $\overline *_n=\{(p,L)\in \Pi\times\mathcal L_\Pi:p\in L\}$;
\item $|\overline *_n|=(n^2+n+1)\cdot (n+1)$;
\item for every line $L\subset \Pi$, the affine plane $\Pi\setminus L$ is  translation;
\item for every line $L\subset \Pi$, the affine plane $\Pi\setminus L$ is $\partial$-translation;
\item $|\Aut_{x,L}(\Pi)|>1$ for every line $L\subset \Pi$ and point $x\in L$;
\item the Lenz graph $\LG(\Pi)$ of $\Pi$ is wide and connected;
\item the Lenz graph $\LG(\Pi)$ of $\Pi$ is wide and $\{|\Aut_{x,L}(\Pi)|:(x,L)\in V_\Pi\}\subseteq p\IN$ for some prime number $p$;
\item the Lenz graph $\LG(\Pi)$ of $\Pi$ is wide and $n$ is a prime power.
\end{enumerate}
\end{theorem}

\begin{proof} The equivalence $(1)\Leftrightarrow(2)$ follows from Theorem~\ref{t:finite-Papp<=>Des}.
\smallskip

The implication $(2)\Ra(3)$ follows from Baer's Theorem~\ref{t:Baer-pL-Des<=>}.
\smallskip

$(3)\Ra(4)$ Assume that $|\Aut_{**}(\Pi)|=\{n,n-1\}$. By Corollary~\ref{c:|Aut|-divides}, for every $p\in \Pi$ and $L\in\mathcal L_\Pi$, the cardinal $|\Aut_{p,L}(\Pi)|\in \{n,n-1\}$ divides $n$ if $p\in L$ and $n-1$ if $p\notin L$. This implies that $\overline *_{n-1}=\{(p,L)\in\Pi\times\mathcal L_\Pi:p\notin L\}$.
\smallskip

$(4)\Ra(5)$ If $\overline *_{n-1}=\{(p,L)\in\Pi\times\mathcal L_\Pi:p\notin L\}$, then $$|\overline *_{n-1}|=|\{(p,L)\in \Pi\times\mathcal L_\pi:p\notin L\}|=|\mathcal L_\Pi|\cdot(|\Pi|-|\Pi|_2)=(n^2+n+1)\cdot n^2.$$
\smallskip

$(5)\Ra(4)$ Assume that $|\overline *_{n-1}|=(n^2+n+1)\cdot n^2$. Consider the set $S\defeq\{(p,L)\in\Pi\times\mathcal L_\Pi:p\notin L\}$ and observe that $|S|=(n^2+n+1)\cdot n^2$. By Corollary~\ref{c:|Aut|-divides}, for every pair $(p,L)\in (\Pi\times\mathcal L_\Pi)\setminus S$, the cardinal $|\Aut_{p,L}(\Pi)|$ divides $|\Pi|_2-1=n$ and hence $(p,L)\notin \overline*_{n-1}$. Then $*_{n-1}\subseteq S$ and hence
$$(n^2+n+1)\cdot n^2=|\overline *_{n-1}|\le|S|=(n^2+n+1)\cdot n^2,$$
which implies $\overline *_{n-1}=S$.
\smallskip

$(4)\Ra(6)$ Assume that $\overline *_{n-1}=\{(p,L)\in\Pi\times\mathcal L_\Pi:p\notin L\}$. We have to show that $\{(p,L)\in\Pi\times\mathcal L_\Pi:p\in L\}=\overline *_{n}$. It follows that $\overline *_n\subseteq (\Pi\times\mathcal L_\Pi)\setminus\overline *_{n-1}=\{(p,L)\in\Pi\times\mathcal L_\Pi:p\in L\}$. On the other hand, for every point-line pair $(p,L)\in\Pi\times\mathcal L_\Pi$ with $p\in L$, we can choose two distinct points  $x,y\in \Pi\setminus L$ with $p\in\Aline xy$. The equality  $\overline *_{n-1}=\{(p,L)\in\Pi\times\mathcal L_\Pi:p\notin L\}$ implies that the plane $\Pi$ is $(x,L)$-transitivite and $(y,L)$-transitive. By Theorem~\ref{t:(pq,L)-transitive}, $\Pi$ is $[\Aline xy,L)$-transitive and hence $(p,L)$-transitive, which implies $|\Aut_{p,L}(\Pi)|=n$ and $(p,L)\in\overline *_n$.
\smallskip

The equivalence $(6)\Leftrightarrow(7)$ follows from Corollary~\ref{c:|Aut|-divides} and can be proved by analogy with the equivalence $(4)\Leftrightarrow(5)$.
\smallskip

The implications $(6)\Ra(8)\Ra(9)\Ra(10)\Ra(11)$ are trivial.
\smallskip

$(11)\Ra(12)$ Assume that the Lenz graph $\LG(\Pi)$ of $\Pi$ is wide and connected. Observe that for every line $L\subseteq \Pi$, the set $$\Aut_{L}(\Pi)\defeq\bigcup_{x\in L}\Aut_{x,L}(\Pi)$$ can be identified with the translation group $\Trans(\Pi\setminus L)$ of the affine plane $\Pi\setminus L$. Since the Lenz grah $\LG(\Pi)$ of $\Pi$ is line-wide, the group $\Aut_L(\Pi)$ is elementary Abelian, by Theorem~\ref{t:Trans-commutative}. Consequently, there exists a prime number $p_L$ such that every non-identity element of the group $\Aut_L(\Pi)$ has order $p_L$. In particular, for every $x\in L$, every non-identity element of the group $\Aut_{x,L}(\Pi)$ has order $p_L$. 

Applying the same argument to the dual projective plane, we can prove that for every point $x\in \Pi$ the group $\Aut_x(\Pi)\defeq\bigcup_{L\in \mathcal L_p}\Aut_{x,L}(\Pi)$ is elementary Abelian and hence there exists a prime number $p_x$ such that every non-identity element of the group $\Aut_x(\Pi)$ has order $p_x$. 

Since for every point-line pair $(x,L)\in V_\Pi$, the group $\Aut_{x,L}(\Pi)\subseteq \Aut_x(\Pi)\cap\Aut_L(\Pi)$ is not trivial, we have the equality $p_x=p_L$. The connectedness of the Lenz graph $\LG(\Pi)$ implies  $\{p_x:x\in\Pi\}=\{p_L:L\in\mathcal L_\Pi\}=\{p\}$ for some prime number $p$. Therefore, $p$ divides $|\Aut_{x,L}(\Pi)$ for all $(x,L)\in V_\Pi$.
\smallskip

The implication $(11)\Ra(12)$ is trivial.
\smallskip

$(12)\Ra(1)$ Assume that the Lenz graph $\LG(\Pi)$ is wide and 
$\{|\Aut_{x,L}(\Pi)|:(x,L)\in V_\Pi\}\subseteq p\IN$ for some prime number $p$. For every line $L$, consider the group 
$$\Aut(\Pi,L)\defeq\{A\in\Aut(\Pi):A[L]=L\}.$$

\begin{claim}\label{cl:Gleason-trans} For every line $L\subset \Pi$ and points $x,y\in L$, there exists an automorphism $A\in\Aut(\Pi,L)$ such that $A(x)=y$.
\end{claim}

\begin{proof} To derive a contradiction, assume that for some $x\in L$, the set $O(x)\defeq \{A(x):A\in\Aut(\Pi,L)\}$ is not equal to $L$, and hence $y\notin O(x)$ for some point $y\in L$. Since the Lenz graph $\LG(\Pi)$ of $\Pi$ is point-wide, there exist lines $L_x\in\mathcal L_x\setminus\{L\}$ and $L_y\in\mathcal L_y\setminus\{L\}$ such that the groups $\Aut_{x,L_x}(\Pi)$ and $\Aut_{y,L_y}(\Pi)$ are non-trivial and hence their orders are divided by the prime number $p$.  Then there exist automorphims $\alpha\in\Aut_{x,L_x}(\Pi)$ and $\beta\in\Aut_{y,L_y}(\Pi)$ of order  $p$. Since the point $x\in L$ is a centre of the automorphism $\alpha$, $\alpha[L]=L$ and hence $\alpha\in \Aut(\Pi,L)$. By analogy we can see that $\beta\in\Aut(\Pi,L)$. Proposition~\ref{p:Fix(A)=L+p} ensures that $\Fix(\alpha)=\{x\}\cup L_x$ and $\Fix(\beta)=\{y\}\cup L_y$ and hence $O(x)\cap \Fix(\alpha)=\{x\}$ and $O(x)\cap\Fix(\beta)=\varnothing$. Let $A$ and $B$ be the cyclic subgroups of the group $\Aut(\Pi,L)$, generated by the automorphisms $\alpha$ and $\beta$, respectively. Since $\alpha$ and $\beta$ have order $p$, the cyclic groups $A$ and $B$ have cardinality $p$. Then for every points $a\in O(x)\setminus\{x\}$ and $b\in O(x)$, the orbits $A(a)\defeq\{F(a):F\in A\}$ and $B(b)\defeq\{F(b):F\in B\}$ have cardinality $p$. Observe that $\{A(a):a\in O(x)\setminus\{x\}\}$ is a partition of the set $O(x)\setminus\{x\}$ into pairwise disjoint sets of cardinality $p$, which implies that $p$ divides the number $|O(x)\setminus\{x\}|$. On the other hand $\{B(b):b\in O(x)\}$ is a partition of the set $O(x)$ into pairwise disjoint sets of cardinality $p$, which implies that $p$ divides $|O(x)|$. But $p$ cannot divide both numbers $|O(x)|$ and $|O(x)|-1$. This is a contradiction showing that $O(x)=L$.
\end{proof}

\begin{claim}\label{cl:AutxL=AutyL} For every line $L\subset \Pi$ and points $x,y\in L$, the subgroups $\Aut_{x,L}(\Pi)$ and $\Aut_{y,L}(\Pi)$ are conjugated in the group $\Aut(\Pi,L)$ and hence $|\Aut_{x,L}(\Pi)|=|\Aut_{y,L}(\Pi)|$.
\end{claim}

\begin{proof} By Claim~\ref{cl:Gleason-trans}, there exists an automorhism $\alpha\in\Aut(\Pi,L)$ such that $\alpha(x)=y$. It is easy to see that for every automorphism $\varphi\in\Aut_{y,L}(\Pi)$, the automorphism $\alpha^{-1}\varphi\alpha$ has centre $x$ and axis $L$, which implies that $\alpha^{-1}\Aut_{y,L}(\Pi)\alpha=\Aut_{x,L}(\Pi)$ and $|\Aut_{y,L}(\Pi)|=|\Aut_{x,L}(\Pi)|$.
\end{proof}

Applying Claim~\ref{cl:AutxL=AutyL} to the dual plane, we obtain that for every point $x\in\Pi$ and lines $L,\Lambda\in\mathcal L_x$, the groups $\Aut_{x,L}(\Pi)$ and $\Aut_{x,\Lambda}(\Pi)$ have the same cardinality. This implies that $V_\Pi=\{(x,L)\in \Pi\times\mathcal L_\Pi:x\in L\}$ and all groups $\Aut_{x,L}(\Pi)$, $(x,L)\in V_\Pi$, have the same cardinality $k>1$. 

Observe that for any line $L\subset \Pi$ and distinct points $x,y\in L$, the groups $\Aut_{x,L}(\Pi)$ and $\Aut_{y,L}(\Pi)$ intersect by the identity element, which implies that the group $\Aut_L(\Pi)=\bigcup_{x\in L}\Aut_{x,L}(\Pi)$ has cardinality $$|\Aut_L(\Pi)|=1+|L|\cdot (k-1).$$ By Proposition~\ref{p:Fix(A)=L+p}, for every point $a\in \Pi\setminus L$, the map $\Aut_L(\Pi)\to \Pi\setminus L$, $A\mapsto A(a)$, is injective, which implies that the number $|\Aut_L(\Pi)|=1+n\cdot(k-1)$  divides the cardinality $|\Pi\setminus L|=(|L|-1)^2$ of the affine plane $\Pi\setminus L$. 
Then $$|L|^2-2|L|+1=(1+|L|\cdot(k-1))\cdot m$$ for some $m\in\IN$, and hence $m-1$ is divisible by $|L|$. Assuming that $m>1$, we conclude that $m-1\ge |L|$ and hence
$$(|L|-1)^2=(1+|L|\cdot(k-1))\cdot m\ge (1+|L|)\cdot(|L|+1),$$
which is a contradiction showing that $m=1$. Then $|L|^2-2|L|+1=1+|L|\cdot (k-1)$ implies $k=|L|-1$ and $$|\Aut_L(\Pi)|=1+|L|\cdot(|L|-2)=(|L|-1)^2=|\Pi\setminus L|.$$ Then for every $a\in\Pi\setminus L$, the injective map $\Aut_L(\Pi)\to\Pi\setminus L$, $A\mapsto A(a)$, is surjective, which means that $L$ is a translation line in the projective plane $\Pi$. Therefore, every line in the projective plane $\Pi$ is translation and $\Pi$ is Moufang and Pappian, by Theorem~\ref{t:Skornyakov-San-Soucie}. 
\smallskip

The implication $(1)\Ra(13)$ follows from the implication $(1)\Ra(12)$ and Theorem~\ref{t:cardinality-pD}. The implication $(13)\Ra(12)$ follows from Corollary~\ref{c:|Aut|-divides}(1).
\end{proof} 

\begin{definition} For a projective plane $\Pi$, the cardinal numbers
$$\Lenz_\Sigma(\Pi)\defeq\sum_{p\in L\in\mathcal L_\Pi}(|\Aut_{x,L}(\Pi)|-1)$$
and 
$$\LB_\Sigma(\Pi)\defeq\sum_{(p,L)\in\Pi\times\mathcal L_\Pi}(|\Aut_{x,L}(\Pi)|-1)$$
are called the \index{Lenz weight}\index{projective plane!Lenz weight}\defterm{Lenz weight} and the \index{Lenz--Barlotti weight}\index{projective plane!Lenz--Barlotti weight}\defterm{Lenz--Barlotti weight} of $\Pi$, respectively.
\end{definition}

\begin{remark} The Lenz--Barlotti weight of a projective plane is equals to the number of non-trivial central automorphisms of the plane. It is strictly smaller that the cardinality of the automorphism group of the plane.
\end{remark}
 
\begin{proposition} Every Pappian projective plane $\Pi$ of finite order $n$ has
$$\Lenz_\Sigma(\Pi)=(n+1)(n^3-1)\quad\mbox{and}\quad\LB_\Sigma(\Pi)=(n^2+n+1)(n^3-n^2-1).$$
\end{proposition}

\begin{proof} Observe that $$\Lenz_\Sigma(\Pi)=\sum_{p\in L\in\mathcal L_\Pi}(|\Aut_{x,L}(\Pi)|-1)=(n-1)(n+1)(n^2+n+1)=(n+1)(n^3-1)$$and
$$
\begin{aligned}
\LB_\Sigma(\Pi)&=\Lenz_\Sigma(\Pi)+\sum_{L\in\mathcal L_\Pi}\sum_{p\in\Pi\setminus L}(|\Aut_{p,L}(\Pi)|-1)=(n^2-1)(n^2+n+1)+(n^2+n+1)n^2(n-2)\\
&=(n^2+n+1)(n^2-1+n^3-2n^2)=(n^2+n+1)(n^3-n^2-1).
\end{aligned}
$$
\end{proof}

\begin{exercise} Show that a finite projective plane $\Pi$ has Lenz--Barlotti weight (resp. Lenz weight) zero if and only if for every point-line pair $(p,L)\in\Pi\times\mathcal L_\Pi$ (with $p\in L$), the group $\Aut_{p,L}(\Pi)$ is trivial.
\end{exercise}

\begin{question} Is there a finite projective plane with Lenz weight zero?
\end{question}


\begin{definition} For a projective plane $\Pi$ of finite order $n$, the rational numbers
$$\Lenz_\sigma(\Pi)\defeq\frac{\Lenz_\Sigma(\Pi)}{(n+1)(n^3-1)}\quad\mbox{and}\quad\LB_\sigma(\Pi)\defeq\frac{\LB_\Sigma(\Pi)}{(n^2+n+1)(n^3-n^2-1)}$$
are called the \index{Lenz norm}\index{projective plane!Lenz norm}\defterm{Lenz norm} and the \index{Lenz--Barlotti norm}\index{projective plane!Lenz--Barlotti norm}\defterm{Lenz--Barlotti norm} of $\Pi$, respectively.
\end{definition}

\begin{examples} The Lenz and Lenz-Barlotti weights and norms of all projective planes of order $9$ are presented in the following table.
$$
\begin{array}{|l|c|c|c|c|}
\hline
\mbox{Projective plane:}&\mbox{Desargues}&\mbox{right-Hall}&\mbox{left-Hall}&\mbox{Hughes}\\
\hline
\mbox{Lenz weight}&7280&80&80&104\\
\mbox{Lenz--Barlotti weight}&58877&791&791&221\\
\mbox{Lenz norm}&1&1/91&1/91&1/70\\
\mbox{Lenz--Barlotti norm}&1&113/8433&113/8411&17/4529\\
\hline
\end{array}
$$
\end{examples}

\begin{exercise} Show that a finite projective plane $\Pi$ is Pappian if and only if $\Lenz_\sigma(\Pi)=1$ if and only if $\LB_\sigma(\Pi)=1$.
\end{exercise}

\begin{proposition}\label{p:Lenz-degree} Let $\Pi$ be a projective plane of order $n$, $p$ be the smallest prime divisor of $n$, and $d\defeq n/p$ be the largest proper divisor of $n$. If the projective plane $\Pi$ is of Lenz type
\begin{itemize}
\item[] $(\circ)$, then $\Lenz_\sigma(\Pi)< \frac 1p$;
\item[] $(\bullet)$, then $\frac {n-1}{(n+1)(n^3-1)}\le\Lenz_\sigma(\Pi)\le \frac{2nd-n-1}{(n+1)(n^3-1)}<\frac 2{n(n+1)p}$;
\item[] $(-)$ or $(\,|\,)$, then $\frac 1{n^2+n+1}\le\Lenz_\sigma(\Pi)\le \frac{nd-1}{n^3-1}<\frac1{np}$;
\item[] $(\circledast)$, then $\Lenz_\sigma(\Pi)=\frac{1}{n^2+n+1}<\frac1{n^2}$;
\item[] $(+)$, then $\Lenz_\sigma(\Pi)=\frac{2n+1}{(n+1)(n^2+n+1)}<\frac{2}{n^2}$;
\item[] $(\square)$, then $\Lenz_\sigma(\Pi)=1$.
\end{itemize}
\end{proposition}

\begin{proof} For every line $L\in\mathcal L_\Pi$, the union $\Aut_L(\Pi)=\bigcup_{x\in L}\Aut_{x,L}(\Pi)$ can be identified with the translation group of the affine plane $\Pi\setminus L$. 

\begin{lemma}\label{l:Lenz-graph} If $|\Aut_L(\Pi)|>n(d-1)>0$ for some line $L\subseteq\Pi$, then either $|\Aut_{x,L}(\Pi)|=n$ for some $x\in L$ or $|\Aut_{x,L}|>1$ for all $x\in L$.
\end{lemma}

\begin{proof} To derive a contradiction, assume that $|\Aut_{x,L}(\Pi)|<n$ for all $x\in L$ and $|\Aut_{x,L}(\Pi)|=1$ for some $x\in L$. By Corollary~\ref{c:|Aut|-divides}, $|\Aut_{x,L}(\Pi)|$ divides $n$ and hence $|\Aut_{x,L}|\le d$ for all $x\in L$. Fix any point $o\in\Pi\setminus L$ and consider the its orbit  $O\defeq\{A(o):A\in\Aut_L(\Pi)\}$ under the (free) action of the group $\Aut_L(\Pi)$. Proposition~\ref{p:Ax=Bx=>A=B} implies that $|O|=|\Aut_L(\Pi)|$. For every $x\in L$, the intersection $O\cap\Aline ox$ has cardinaility $|O\cap\Aline ox|=|\Aut_{x,L}(\Pi)|\le d$. By our assumption, $|\Aut_{x,L}(\Pi)|=1$ for some $x\in L$. Then $$|\Aut_L(\Pi)|=|O|=1+\bigcup_{x\in L}(|\Aut_{x,L}(\Pi)-1)\le 1+(|L|-1)\cdot(d-1)=1+n(d-1).$$ Since $\Aut_L(\Pi)$ acts freely on the affine plane $\Pi\setminus L$, $|\Aut_L(\Pi)|\le 1+n(d-1)$ divides $|\Pi\setminus L|=n^2$, which implies $|\Aut_L(\Pi)|\le n(d-1)$ and contradicts our assumption.
\end{proof}
 
 Now we are ready to complete the proof of the proposition. Let $\LG(\Pi)=(V_\Pi,E_\Pi)$ be the Lenz graph of the plane $\Pi$.
\smallskip

 $(\circ)$ Assume that the projective plane $\Pi$ is of Lenz type $(\circ)$.  Then $|\Aut_{x,L}(\Pi)|<n$ for all point-line pairs $(x,L)\in V_\Pi$. By Corollary~\ref{c:|Aut|-divides}, $|\Aut_{x,L}(\Pi)|$ divides $n=|\Pi|_2-1$ and hence $|\Aut_{x,L}(\Pi)|\le d$ for all $(x,L)\in V_\Pi$. If $n$ is prime, then $d=1$ and hence $\L(\Pi)=0<\frac1p$. So, assume that $n$ is not prime and hence $d>1$.
To derive a contradiction, assume that $\L(\Pi)\ge\frac1p$.

We claim that the set $\mathcal L'\defeq\{L\in\mathcal L_\Pi:L\times\{L\}\subseteq V_\Pi\}$ contains at least two distinct lines. By Lemma~\ref{l:Lenz-graph}, for every $L\in\mathcal L_\Pi\setminus\mathcal L'$, the group $\Aut_L(\Pi)\defeq\bigcup_{x\in L}\Aut_{x,L}(\Pi)$ has cardinality $|\Aut_L(\Pi)|\le n(d-1)$. On the other hand, $\sum_{x\in L}(|\Aut_{x,L}(\Pi)-1)\le |L|\cdot (d-1)=(n+1)(d-1)$ for all $L\in\mathcal L'$. Assuming that $|\mathcal L'|\le 1$, we conclude that 
$$\sum_{x\in L\in\mathcal L_\Pi}(|\Aut_{x,L}(\Pi)|-1)\le (n+1)(d-1)+
(n^2+n)(n(d-1)-1)=(n+1)(n^2d-n^2-n+d-1)$$and hence
$n^3-1\le (n^3-1)\cdot p\cdot \Lenz_\sigma(\Pi)\le p(n^2d-n^2-n+d-1)=n^3-pn^2-np+n-p< n^3-1$, 
which is a contradiction showing that $|\mathcal L'|\ge 2$. Applying the same argument to the dual plane, we conclude that the set $\{x\in \Pi:\{x\}\times\mathcal L_x\subseteq V_\Pi\}$ contains at least two distinct ponts. This implies that the Lenz graph of $\Pi$ is wide and connected. By Theorem~\ref{t:Gleason56}, the projective plane $\Pi$ is Pappian and hence cannot be of the Lenz class $(\circ)$. This contradiction shows that $\Lenz_\sigma(\Pi)<\frac1p$.  
\smallskip

$(\bullet)$ Assume that the projective plane $\Pi$ is of Lenz class $(\bullet)$. Then its Lenz figure $\Lenz_\Pi$ is equal to $\{(q,\Lambda)\}$ for some point-line pair $(q,\Lambda)\in V_\Pi$. Assuming that some pair $(x,L)\in V_\Pi$ has $x\ne q\notin L$, we can take any non-trivial automorphism $\Phi\in \Aut_{x,L}(\Pi)$ and conclude that the pair $(\Phi(q),\Phi(\Lambda))\ne (q,\Lambda)$ belongs to the Lenz figure $\Lenz_\Pi=\{(q,\Lambda)\}$, which is a contradiction showing that each pair $(x,L)\in V_\Pi$ has $x=q$ or $L=\Lambda$. Then
$$
\begin{aligned}
\sum_{x\in L\in\mathcal L_\Pi}|\Aut_{x,L}(\Pi)-1|&=(|\Aut_{q,\Lambda}(\Pi)|-1)+\sum_{x\in \Lambda\setminus\{q\}}(|\Aut_{x,\Lambda}(\Pi)|-1)+\sum_{L\in\mathcal L_q\setminus\{\Lambda\}}(|\Aut_{q,L}(\Pi)|-1)\\
&=
(n-1)+n(d-1)+n(d-1)\le 2nd-n-1
\end{aligned}
$$and hence 
$$\frac{n-1}{(n+1)(n^3-1)}\le\Lenz_\sigma(\Pi)\le\frac{2nd-n-1}{(n+1)(n^3-1)}<\frac2{n(n+1)p}.$$
\smallskip

$(-)$ Assume that the projective plane $\Pi$ is of Lenz class $(-)$, which means that $\Lenz_\Pi=\Lambda\times\{\Lambda\}$ for some line $\Lambda$. The uniqueness of $\Lambda$ implies that for any point $(x,L)\in V_\Pi$ we have $x\in \Lambda$ or $L=\Lambda$. Then 
$$
\begin{aligned}
\sum_{x\in L\in\mathcal L_\Pi}(|\Aut_{x,L}(\Pi)|-1)&=\sum_{x\in \Lambda}(|\Aut_{x,\Lambda}(\Pi)|-1)+\sum_{x\in\Lambda}\sum_{L\in\mathcal L_x\setminus\{\Lambda\}}|\Aut_{x,\Lambda}(\Pi)-1|\\
&\le (n+1)(n-1)+(n+1)n(d-1)=(n+1)(nd-1).
\end{aligned}
$$
Then
$$\frac{1}{n^2+n+1}=\frac{(n+1)(n-1)}{(n+1)(n^3-1)}=\sum_{x\in\Lambda}\frac{|\Aut_{x,\Lambda}(\Pi)|-1}{(n+1)(n^3-1)}\le \Lenz_\sigma(\Pi)\le \frac{nd-1}{n^3-1}<\frac1{pn}.$$
\smallskip

$(\,|\,)$ If $\Pi$ has Lenz type $(\,|\,)$, then the dual projective plane $\Pi^*$ has Lenz type $(-)$ and the preceding item ensures that $
\frac{1}{n^2+n+1}\le \Lenz_\sigma(\Pi)=\Lenz_\sigma(\Pi^*)\le \frac{nd-1}{n^3-1}<\frac1{pn}.$
\smallskip

$(\circledast)$ Assume that the projective plane $\Pi$ is of Lenz class $(\circledast)$. Then its Lenz figure $\Lenz_\Pi$ is equal to the set $\{(x,\Aline ox):x\in \Lambda\}$ for some line $\Lambda\subset \Pi$ and point $o\in \Pi\setminus\Lambda$. Assuming that $V_\Pi\ne\Lenz_\Pi$, we could find a pair $(q,L)\in V_\Pi\setminus\Lenz_\Pi$ and a non-trivial automorphism $\Phi\in\Aut_{q,L}(\Pi)$. Consider the point $o'\defeq\Phi(o)$ and observe that $\{(x,\Aline {o'}x):x\in\Phi[\Lambda]\}\subseteq \Lenz_\Pi=\{(x,\Aline ox):x\in\Lambda\}$, which implies $o'=o$ and $\Phi[\Lambda]=\Lambda$. Then $o'=o\in\Fix(\Phi)=L$ and $q\in \Lambda$. Then $(q,L)=(q,\Aline oq)\in \Lenz_\Pi$, which contradicts the choice of $(q,L)\notin \Lenz_\Pi$. This contradiction shows that $\Lenz_\Pi=V_\Pi$ and $$\Lenz_\sigma(\Pi)=\sum_{(x,L)\in\Lenz_\Pi}\frac{|\Aut_{x,L}(\Pi)|-1}{(n+1)(n^3-1)}=\frac{(n+1)(n-1)}{(n+1)(n^3-1)}=\frac{1}{n^2+n+1}<\frac{1}{n^2}.$$
\smallskip

$(+)$ Assume that the projective plane $\Pi$ is of Lenz class $(+)$. Then its Lenz figure $\Lenz_\Pi$ is equal to the set $(\Lambda\times\{\Lambda\})\cup(\{q\}\times\mathcal L_q)$ for some line $\Lambda\in\mathcal L_\Pi$ and point $q\in \Lambda$. Every automorphism of $\Pi$ fixes the line $L$ and the point $q$, which implies that $V_\Pi=\Lenz_\Pi$. Then $$\Lenz_\sigma(\Pi)=\sum_{(x,\in L)\in\Lenz_\Pi}\frac{|\Aut_{x,L}(\Pi)|-1}{(n+1)(n^3-1)}=\frac{(2n+1)(n-1)}{(n+1)(n^3-1)}=\frac{2n+1}{(n+1)(n^2+n+1)}<\frac{2}{n^2}.$$
\smallskip

$(\square)$ If the projective plane $\Pi$ has Lenz type $(\square)$, then $V_\Pi=\Lenz_\Pi=\{(x,L)\in\Pi\times\mathcal L_\Pi:x\in L\}$ and 
$$\Lenz_\sigma(\Pi)=\sum_{x\in L\in\mathcal L_\Pi}\frac{|\Aut_{x,L}(\Pi)|-1}{(n+1)(n^3-1)}=\frac{(n+1)(n^2+n+1)(n-1)}{(n+1)(n^3-1)}=1.$$
\end{proof}

Lenz Theorem~\ref{t:Lenz} and  Proposition~\ref{p:Lenz-degree} imply the following upper bound for the Lenz norm of finite non-Desarguesian projective planes.

\begin{corollary}\label{c:Lenz-degree} Every finite non-Desarguesian projective plane $\Pi$ has Lenz norm\break $\Lenz_\sigma(\Pi)<\frac1p$ where $p$ is the largest prime divisior of the order of $\Pi$.
\end{corollary}

\begin{remark} By a deep result of Hering and Kantor \cite{HK1971}, no finite projective plane has Lenz type $(\circledast)$, so the case $(\circledast)$ in Proposition~\ref{p:Lenz-degree} is in fact obsolete. However, we have included this (non-existing case) to Proposition~\ref{p:Lenz-degree} in order to present a self-contained proof of Corollary~\ref{c:Lenz-degree}.
\end{remark}

\begin{exercise} Show that any finite non-Desarguesian projective plane $\Pi$ has Lenz--Barlotti norm $\LB(\Pi)<\frac12$. Can $\frac12$ be replaced by $\frac1p$ where $p$ is the smallest prime divisor of the order of $\Pi$?
\end{exercise}


\begin{remark} The levelsize, levelrank and weight of all known projective planes \cite{Moorhouse} of orders 16, 25, 27 and 32 were calculated by Ivan Hetman. The results of the calculations are presented in the following tables. They show that the central levelsize and levelrank provide a finer classification of finite projective planes, comparing to the Lenz--Barlotti classification.
\end{remark} 

\newcolumntype{C}[1]{>{\centering\arraybackslash}m{#1}}

\begin{table}[h]
\caption{Levelsize, levelrank and weight of known projective planes of order 16}\label{t:16}
\begin{tabular}{|C{60pt}|c|C{40pt}|C{72pt}|C{60pt}|C{35pt}|C{60pt}|}
\hline
{\tt Projective plane}&type&central \!spectrum\!&central levelsize&central levelrank&Lenz weight&\!Lenz-Barlotti\! weight\\
\hline
{\tt bbs4}&$(\circ\circ)$&4,3,1&$5_48_3$&$12_432_3$&15&31\\
{\tt dbbs4}&$(\circ\circ)$&4,3,1&$5_48_3$&$21_423_3$&15&31\\
{\tt john}&$(\circ\circ)$&4,2,1&$5_412_2$&$12_432_2$&27&27\\
{\tt djohn}&$(\circ\circ)$&4,2,1&$5_412_2$&$21_423_2$&27&27\\
{\tt dbbh2}&$(\circ\circ)$&4,2,1&$5_420_2$&$21_423_2$&35&35\\
{\tt bhh2}&$(\circ\circ)$&4,2,1&$5_420_2$&$12_433_2$&35&35\\
{\tt bhh1,dbbh1}&$(\circ\circ)$&8,4,2,1&$1_88_416_324_2$&$11_822_422_322_2$&55&87\\
\hline
{\tt math,dmath}&$(\bullet\circ)$&16,2,1&$1_{16}32_2$&$11_{16}22_2$&47&47\\
\hline
{\tt demp}&$(\,|\,\circ)$&16,3,1&$17_{16}16_3$&$12_{16}21_3$&255&287\\
{\tt ddemp}&$({-}{\circ})$&16,3,1&$17_{16}16_3$&$21_{16}12_3$&255&287\\
{\tt dsfp}&$(\,|\,\circ)$&16,3,1&$17_{16}32_3$&$12_{16}23_3$&255&319\\
{\tt ddsfp}&$({-}{\circ})$&16,3,1&$17_{16}32_3$&$21_{16}32_3$&255&319\\
{\tt lmrh,jowk}&$(\,|\,\circ)$&16,2,1&$17_{16}48_2$&$12_{16}32_2$&303&303\\
{\tt dlmrh,djowk}&$({-}{\circ})$&16,2,1&$17_{16}48_2$&$21_{16}23_2$&303&303\\
{\tt hall}&$(\,|\,\circ)$&16,5,3,2,1&$17_{16}192_5256_380_2$&$12_{16}32_513_332_2$&335&1615\\
{\tt dhall}&$({-}{\circ})$&16,5,3,2,1&$17_{16}192_5256_380_2$&$21_{16}23_531_323_2$&335&1615\\
\hline
{\tt semi2}&$(+\circ)$&16,1&$33_{16}$&$22_{16}$&495&495\\
{\tt semi4}&$(+\circ)$&$16,3,1$&$33_{16}768_3$&$22_{16}33_3$&495&2031\\
\hline
{\tt desarg}&$(\square\blacksquare)$&16,15&$69615_{16}69888_{15}$&$33_{16}33_{15}$&69615&1048047\\
\hline
\end{tabular}
\end{table}

\begin{table}\label{t:25a}
\caption{Levelsize, levelrank, and weight of known projective planes of order 25}
\begin{tabular}{|C{100pt}|c|C{40pt}|C{42pt}|C{45pt}|C{30pt}|C{60pt}|}
\hline
{\tt Projective plane}&type&central \!spectrum\!&central levelsize&central levelrank&Lenz weight&\!Lenz-Barlotti\! weight\\
\hline
{\tt a1b,a4d,a5b,a5c,
a5d,a5e,a5f,a5g,
a6b,a7c,a8d,a8e,
a8f,b1a,b1b,b1c,
b2a,b2b,b2c,b3a,
b3b,b3c,b3d,b3e,
b6d,b6e,b7d,b8b,
b8e,b8f, b8g}
&$(\circ\circ)$&5,1&$6_5$&$12_5$&24&24\\
\hline
{\tt a1bD,a4dD,a5bD,a5cD, 
a5dD,a5eD,a5fD,a5gD, 
a6bD,a7cD,a8eD,a8fD, 
b1aD,b1bD,b1cD,b2aD, 
b2bD,b2cD,b3aD,b3bD,
b3cD,b3dD,b3eD,b6dD,
b6eD,b7dD,b8bD,b8dD, 
b8eD, b8fD, b8gD}&$(\circ\circ)$&5,1&$6_5$&$21_5$&24&24\\
\hline
{\tt a1c,a3b,a4b,a4c,
a5a,a6a,b6b,a6c,
b6c,b7c,b8d,s3b}&$(\circ\circ)$&5,2,1&$6_55_2$&$12_521_2$&24&29\\
\hline
{\tt a1cD,a3bD,a4bD,a4cD, 
a5aD,a6aD,b6bD,a6cD, 
b6cD,b7cD,b8dD,s3bD}&$(\circ\circ)$&5,2,1&$6_55_2$&$21_512_2$&24&29\\
\hline
{\tt b4b,s5b,s5d}&$(\circ\circ)$&5,2,1&$6_515_2$&$12_532_2$&24&39\\
\hline
{\tt b4bD,s5bD,s5dD}&$(\circ\circ)$&5,2,1&$6_515_2$&$21_523_2$&24&39\\
\hline
{\tt b4a,b5b}&$(\circ\circ)$&5,2,1&$6_525_2$&$12_522_2$&24&49\\
\hline
\end{tabular}
\end{table}

\begin{table}\label{t:25b}
\begin{tabular}{|C{55pt}|c|C{45pt}|C{75pt}|C{65pt}|C{30pt}|C{60pt}|}
\hline
{\tt Projective plane}&type&central spectrum&central\break levelsize&central levelrank&Lenz weight&\!Lenz-Barlotti\! weight\\
\hline
\mbox{\tt h1a}&$(\bullet\circ)$&25,4,2,1&$1_{25}25_426_2$&$11_{25}12_432_2$&24&125\\
\mbox{\tt h1aD}&$(\bullet\circ)$&25,4,2,1&$1_{25}25_426_2$&$11_{25}21_423_2$&24&125\\
\mbox{\tt h2a}&$(\bullet\circ)$&25,4,2,1&$1_{25}26_425_2$&$11_{25}32_412_2$&24&127\\
\mbox{\tt h2aD}&$(\bullet\circ)$&25,4,2,1&$1_{25}26_425_2$&$11_{25}23_421_2$&24&127\\
\mbox{\tt w1}&$(\bullet\circ)$&25,8,4,2,1&$1_{25}1_825_425_2$&$11_{25}11_821_412_2$&24&132\\
\mbox{\tt w1D}&$(\bullet\circ)$&25,8,4,2,1&$1_{25}1_825_425_2$&$11_{25}11_812_421_2$&24&132\\
\mbox{\tt a3a}&$(\bullet\circ)$&25,5,2,1&$1_{25}25_525_2$&$11_{25}12_522_2$&24&149\\
\hline
\end{tabular}
\end{table}

\begin{table}
\caption{Levelsize, levelrank and weight of known projective planes of order 27}\label{t:27}
\begin{tabular}{|C{60pt}|c|C{40pt}|C{70pt}|C{45pt}|C{30pt}|C{60pt}|}
\hline
{\tt Projective plane}&type&central \!spectrum\!&central levelsize&central levelrank&Lenz weight&\!Lenz-Barlotti\! weight\\
\hline
{\tt fig}
&$(\circ\circ)$&3,2,1&$52_3117_2$&$33_333_2$&104&121\\
\hline
{\tt  flag4, flag6, hering}
&$(\,|\,\circ)$&27,2,1&$28_{27}729_2$&$12_{27}13_2$&728&1457\\
\hline
{\tt  dflag4, dflag6, dhering}
&$({-}\circ)$&27,2,1&$28_{27}729_2$&$21_{27}31_2$&728&1457\\
\hline
{\tt andre}
&$(\,|\,\circ)$&27,13,2,1&$28_{27}54_{13}729_2$&$12_{27}32_313_2$&728&2105\\
{\tt dandre}
&$({-}\circ)$&27,13,2,1&$28_{27}54_{13}729_2$&$21_{27}23_331_2$&728&2105\\
\hline
{\tt sherk}
&$(\,|\,\circ)$&27,3,2,1&$28_{27}27_3729_2$&$12_{27}21_313_2$&782&1511\\
{\tt dsherk}
&$({-}\circ)$&27,3,2,1&$28_{27}27_3729_2$&$21_{27}12_331_2$&782&1511\\
\hline
{\tt twisted}
&$(+\circ)$&27,2,1&$55_{27}2187_2$&$22_{27}33_2$&1430&3617\\
\hline
{\tt desarg}
&$(\square\blacksquare)$&27,26&$21196_{27}551853_{26}$&$33_{27}33_{26}$&551096&14347421\\
\hline
\end{tabular}
\vskip20pt
\end{table}

\begin{table}\label{t:32}
\caption{Levelsize, levelrank and weight of known projective planes of order 32}
\begin{tabular}{|C{60pt}|c|C{40pt}|C{75pt}|C{45pt}|C{35pt}|C{60pt}|}
\hline
{\tt Projective plane}&type&central \!spectrum\!&central levelsize&central levelrank&Lenz weight&\!Lenz-Barlotti\! weight\\
\hline
{\tt e,f,g}
&$(\,|\,\circ)$&32,1&$33_{32}$&$12_{32}$&1023&1023\\
{\tt eD,fD,gD}
&$({-}\circ)$&32,1&$33_{32}$&$21_{32}$&1023&1023\\
\hline
{\tt a,b,d,aD,bD}
&$(+\circ)$&32,1&$65_{32}$&$22_{32}$&2015&2015\\
\hline
{\tt desarg}
&$(\square\blacksquare)$&32,31&$34881_{32}1082368_{31}$&$33_{32}33_{31}$&1081311&33552351\\
\hline
\end{tabular}
\end{table}

\chapter{Geometry of operations in ternars}

In this chapter we study the interplay between geometric properties of based affine planes and algebraic properties of their ternars. 

\section{Geometry of the plus operation}

In this section we find geometric counterparts of some algebraic properties of the plus operation of a ternar. Let us recall that for a based affine plane $(\Pi,uow)$ and any elements $a,b,c\in\Delta$ on its diagonal $\Delta$, the equality $a+b=c$ holds if and only if $\Aline{ac}{ob}\subparallel\Delta$. 

%

Here for two points $x,y\in \Delta$ we denote by $xy$ the unique point with cordinates $x,y$. So, $\Aline x{xy}\subparallel \Aline ow$ and $\Aline {xy}y\subparallel\Aline ou$. Therefore, we identify the based plane $\Pi$ with the coordinate plane $\Delta^2$ of the ternar $\Delta$. We denote by $\Delta_\parallel$ the spread of lines in $\Pi$, parallel to the diagonal $\Delta$.

\begin{proposition}\label{p:two-sided} For a based affine plane $(\Pi,uow)$ and its ternar $\Delta$,\newline the following conditions are equivalent:
\begin{enumerate}
\item the plus loop $(\Delta,+)$ is invertible;
\item $\forall x,y\in \Delta\;\;(\Aline {ox}{yo}\subparallel\Delta\;\Ra\;\Aline {xo}{oy}\subparallel \Delta)$;
\item $\forall x,y\in \Delta\;\;(\Aline {ox}{yo}\subparallel\Delta\;\Leftrightarrow \;\Aline {xo}{oy}\subparallel \Delta)$. 


\end{enumerate}
\end{proposition}

\begin{proof}  $(1)\Ra(2)$ Assume that the loop $(\Delta,+)$ is invertible.  Given any points $x,y\in \Delta$ with $\Aline {ox}{yo}\subparallel\Delta$, we should prove that $\Aline {xo}{oy}\subparallel \Delta$. The subparallelity relation  $\Aline {ox}{yo}\subparallel\Delta$ implies $y+x=o$, by the definition of the plus operation in the ternar $\Delta$. Since the loop $(\Delta,+)$ is invertible, $x+y=o$, which implies $\Aline {xo}{oy}\subparallel \Delta$. 

%

The implication $(2)\Ra(3)$ follows from the equalities $\Aline {ox}{yo}=\Aline{yo}{ox}$ and $\Aline {xo}{oy}=\Aline {yo}{xo}$. 
\smallskip

 $(3)\Ra(1)$ Assume that the condition (3) holds. To prove that the loop $(\Delta,+)$ is invertible, take any point $x\in \Delta$. Since $(\Delta,+)$ is a loop, there exists an element $y\in\Delta$ such that $y+x=o$. If $o\in \{x,y\}$, then $y+x=o=o+o$ implies $x=o=y$ and hence $y=o$ is the two-sided inverse to $x=o$ in the loop $(\Delta,+)$, witnessing that this loop is invertible. So, assume that $o\notin\{x,y\}$, which implies that $\Aline{ox}{yo}$ is a line in $R^2$. By the definition of the plus operation in the ternar $\Delta$, the equality $y+x=o$ implies $\Aline {ox}{yo}\subparallel \Delta$. Applying the condition (3), we conclude that $\Aline{xo}{oy}\subparallel\Delta$. Then $x+y=o$, by the  definition of the plus operation in the ternar $\Delta$. Therefore, $y$ is the two-sided inverse to $x$ in the loop $(\Delta,+)$.
\end{proof}

\begin{proposition}\label{p:rip<=>}  For a based affine plane $(\Pi,uow)$ and its ternar $\Delta$,\newline the following conditions are equivalent:
\begin{enumerate}
\item the plus loop $(\Delta,+)$ is right-inversive;
\item $\forall x,s,t\in \Delta\;\;(\Aline {ox}{st}\subparallel\Delta\;\Ra\;\Aline {xo}{ts}\subparallel \Delta)$; 
\item $\forall x,y,s,t\in\Delta\;\;(\Aline {xy}{st}\subparallel\Delta\;\Ra\;\Aline{yx}{ts}\subparallel\Delta)$.
\begin{picture}(100,0)(-180,-20)
\linethickness{0.5pt}
\put(-20,-20){\color{teal}\line(1,0){20}}
\put(-20,-20){\color{cyan}\line(0,1){20}}
\put(0,0){\color{teal}\line(-1,0){20}}
\put(0,0){\color{cyan}\line(0,-1){20}}

\put(5,5){\color{teal}\line(1,0){20}}
\put(5,5){\color{cyan}\line(0,1){20}}
\put(25,25){\color{teal}\line(-1,0){20}}
\put(25,25){\color{cyan}\line(0,-1){20}}
\put(-20,0){\color{red}\line(1,1){25}}
\put(0,-20){\color{red}\line(1,1){25}}
{\linethickness{0.75pt}
\put(0,0){\color{red}\line(1,1){25}}
\put(0,0){\color{red}\line(-1,-1){20}}
}

\put(0,0){\circle*{2}}
\put(5,5){\circle*{2}}
\put(25,5){\circle*{2}}
\put(-20,0){\circle*{2}}
\put(5,25){\circle*{2}}
\put(0,-20){\circle*{2}}
\put(25,25){\circle*{2}}
\put(-20,-20){\circle*{2}}
\end{picture}
\end{enumerate}
\end{proposition}

\smallskip

\begin{proof} $(1)\Ra(2)$ Assume that the loop $(\Delta,+)$ is right-inversive. Given any points $x,s,t\in \Delta$ with $\Aline {ox}{st}\subparallel\Delta$, we should prove that $\Aline {xo}{ts}\subparallel \Delta$. By the geometric definition of the plus operation, the subparallelity relation  $\Aline {ox}{st}\subparallel\Delta$ implies $s+x=t$. Since the loop $(\Delta,+)$ is right-inversive, there exists a point $y\in \Delta$ such that $z=(z+x)+y$ for every $z\in\Delta$. In particular, $o=(o+x)+y=x+y$ and $t+y=(s+x)+y=s$ and hence $\Aline {oy}{ts}\subparallel\Delta$. The equality $x+y=o$ implies $\Aline {oy}{xo}\subparallel\Delta$. By the Proclus Axiom, the subparallelity relations $\Aline {oy}{ts}\subparallel\Delta$ and $\Aline {oy}{xo}\subparallel\Delta$ imply $\Aline {xo}{ts}\subparallel\Delta$. 

\begin{picture}(150,145)(-160,-40)

\linethickness{0.75pt}
\put(0,0){\color{teal}\line(1,0){30}}
\put(0,-30){\color{cyan}\line(0,1){60}}
\put(90,90){\color{teal}\line(-1,0){30}}
\put(90,90){\color{cyan}\line(0,-1){30}}
{\linethickness{1pt}
\put(-30,-30){\color{red}\line(1,1){120}}
}
\put(0,30){\color{teal}\line(1,0){30}}
\put(30,0){\color{cyan}\line(0,1){30}}
\put(60,60){\color{teal}\line(1,0){30}}
\put(60,60){\color{cyan}\line(0,1){30}}
\put(0,30){\color{red}\line(1,1){60}}
\put(0,-30){\color{red}\line(1,1){90}}
\put(-30,-30){\color{teal}\line(1,0){30}}

\put(0,0){\circle*{3}}
\put(-8,-1){$o$}
\put(30,0){\circle*{3}}
\put(32,-5){$xo$}
\put(60,60){\circle*{3}}
\put(60,53){$s$}
\put(0,30){\circle*{3}}
\put(-9,33){$ox$}
\put(90,60){\circle*{3}}
\put(93,56){$ts$}
\put(30,30){\circle*{3}}
\put(32,25){$x$}
\put(90,90){\circle*{3}}
\put(93,90){$t$}
\put(60,90){\circle*{3}}
\put(56,93){$st$}
\put(-30,-30){\circle*{3}}
\put(-38,-32){$y$}
\put(0,-30){\circle*{3}}
\put(2,-35){$oy$}
\end{picture}
\smallskip

$(2)\Ra(1)$ Assume that the condition (2) holds. To prove that the loop $(X,+)$ is right-inversive, take any point $x\in\Delta$. Since $(X,+)$ is a quasigroup, there exists a point $y\in\Delta$ such that $x+y=o$. Given any point $s\in \Delta$, we need to check that $(s+x)+y=s$.  By the definition of the plus operation in the ternar $\Delta$, the equality $x+y=o$ implies $\Aline {oy}{xo}\subparallel\Delta$. Also, for the point  $t\defeq s+x$, we have the subparallelity relation $\Aline {ox}{st}\subparallel \Delta$, which implies the subparallelity relation $\Aline{xo}{ts}\subparallel\Delta$, by the condition (2). By the Proclus Axiom, the subparallelity relations  $\Aline {oy}{xo}\subparallel\Delta$ and $\Aline{xo}{ts}\subparallel\Delta$ imply  $\Aline {oy}{ts}\subparallel\Delta$ and hence $(s+x)+y=t+y=s$, by the definition of plus operation in the ternar $\Delta$.
\smallskip

 The implication $(3)\Ra(2)$ is trivial. To prove that $(2)\Ra(3)$, take any points $y,z,s,t\in\Delta$ with $\Aline {yz}{st}\subparallel \Delta$. We have to prove that $\Aline{zy}{ts}\subparallel\Delta$. Find a unique point $x\in\Delta$ such that $\Aline {ox}{yz}\subparallel \Delta$. The subparallelity relations $\Aline {ox}{yz}\subparallel \Delta$ and $\Aline {yz}{st}\subparallel \Delta$ imply $\Aline {ox}{st}\subparallel \Delta$. Applying the condition (2), we obtain $\Aline {xo}{ts}\subparallel \Delta$ and $\Aline{xo}{zy}\subparallel \Delta$. By the Proclus Axiom, the latter subparallelity relations imply $\Aline {zy}{ts}\subparallel \Delta$.


 \end{proof}

\begin{proposition}\label{p:lip<=>}   For a based affine plane $(\Pi,uow)$ and its ternar $\Delta$,\newline the following conditions are equivalent:
\begin{enumerate}
\item the plus loop $(\Delta,+)$ is left-inversive;
\item $\forall x,y,z,s\in \Delta\;\;(\Aline {xo}{oy}\subparallel\Delta\;\wedge\; \Aline{oz}{ys}\subparallel\Delta)\;\Ra\;\Aline {os}{xz}\subparallel \Delta)$;
\item $\forall a,b,x,y,s,t\in\Delta\;\;(\Aline {ax}{oy},\Aline {os}{bt},\Aline {os}{bt}\in\Delta_\parallel\;\Ra\;\Aline {as}{ot}\in\Delta_\parallel)$.
%
\end{enumerate}
\end{proposition}

\smallskip

\begin{proof} $(1)\Ra(2)$ Assume that the loop $(\Delta,+)$ is left-inversive.  To prove the condition (2), take any points $x,y,z,s\in \Delta$ with $\Aline {xo}{oy}\subparallel\Delta$ and $\Aline{oz}{ys}\subparallel\Delta$. By the definition of the plus operation in the ternar $\Delta$, the subparallelity relations $\Aline {xo}{oy}\subparallel\Delta$ and $\Aline{oz}{ys}\subparallel\Delta$ imply $x+y=o$ and $y+z=s$. Since the loop $(\Delta,+)$ is left-inversive, $x+y=o$ and $y+z=s$ imply $z=x+(y+z)=x+s$ and hence $\Aline {os}{xz}\subparallel\Delta$.

\begin{picture}(100,120)(-130,-15)

{\linethickness{1pt}
\put(0,0){\color{red}\line(1,1){100}}
}
\linethickness{0.75pt}
\put(0,0){\color{cyan}\line(0,1){80}}
\put(0,80){\color{red}\line(1,1){20}}
\put(0,80){\color{teal}\line(1,0){80}}
\put(0,20){\color{teal}\line(1,0){20}}
\put(20,20){\color{cyan}\line(0,1){80}}
\put(20,40){\color{teal}\line(1,0){20}}
\put(20,100){\color{teal}\line(1,0){80}}
\put(0,20){\color{red}\line(1,1){20}}
\put(20,80){\color{red}\line(1,1){20}}
\put(40,40){\color{cyan}\line(0,1){60}}

\put(0,0){\circle*{3}}
\put(0,-7){$x$}
\put(20,20){\circle*{3}}
\put(20,13){$o$}
\put(40,40){\circle*{3}}
\put(40,33){$y$}
\put(80,80){\circle*{3}}
\put(80,73){$z$}
\put(100,100){\circle*{3}}
\put(100,93){$s$}

\put(0,20){\circle*{3}}
\put(-14,18){$xo$}
\put(20,40){\circle*{3}}
\put(7,40){$oy$}
\put(0,80){\circle*{3}}
\put(-14,78){$xz$}
\put(20,100){\circle*{3}}
\put(16,103){$os$}
\put(40,100){\circle*{3}}
\put(35,104){$ys$}
\put(60,52){\color{red}$\Delta$}
\put(20,80){\circle*{3}}
\put(22,73){$oz$}

\end{picture}

$(2)\Ra(1)$ Assume that the condition (2) is satisfied. To prove that the loop $(\Delta,+)$ is left-inversive, take any element $y\in\Delta$. Since $(\Delta,+)$ is a quasigroup, there exists an element $x\in\Delta$ such that $x+y=o$ and hence $\Aline{oy}{xo}\subparallel\Delta$. We claim that $x+(y+z)=z$ for every $z\in \Delta$. Consider  the point $s\defeq y+z$ and observe that $\Aline {oz}{ys}\subparallel\Delta$. By the condition (2), the subparallelity relations $\Aline{oy}{ox}\subparallel\Delta$ and $\Aline {oz}{ys}\subparallel\Delta$ imply $\Aline{xz}{os}\subparallel\Delta$ and hence $x+(y+z)=x+s=z$.
\smallskip

$(2)\Ra(3)$ Assume that the condition (2) holds and take any points  $a,b,x,y,s,t\in\Delta$ with $\Aline {ax}{oy},\Aline {ox}{by},\Aline {os}{bt}\in\Delta_\parallel$. We need to show that $\Aline {as}{ot}\parallel \Delta$.  Our assumption ensures that $\Aline {ax}{oy}$, $\Aline {ox}{by}$, $\Aline {as}{bt}$ are lines, parallel to the line $\Delta$. This implies that $a\ne o\ne b$, $x\ne y$ and $s\ne t$.

\begin{picture}(100,120)(-130,-15)

\linethickness{0.5pt}
\put(0,0){\color{cyan}\line(0,1){75}}
\put(0,75){\color{teal}\line(1,0){75}}
\put(0,15){\color{teal}\line(1,0){15}}
\put(15,15){\color{cyan}\line(0,1){75}}
\put(15,30){\color{teal}\line(1,0){15}}
\put(15,90){\color{teal}\line(1,0){75}}
\put(0,15){\color{red}\line(1,1){15}}
\put(30,30){\color{cyan}\line(0,1){60}}
\put(15,60){\color{teal}\line(1,0){45}}
\put(0,45){\color{teal}\line(1,0){45}}

{\linethickness{1pt}
\put(0,0){\color{red}\line(1,1){90}}
\put(0,75){\color{cyan}\line(0,-1){30}}
\put(15,90){\color{cyan}\line(0,-1){45}}
\put(30,90){\color{cyan}\line(0,-1){30}}
\put(0,45){\color{red}\line(1,1){15}}
\put(15,45){\color{red}\line(1,1){15}}
\put(15,75){\color{red}\line(1,1){15}}
\put(0,75){\color{red}\line(1,1){15}}
\put(15,90){\color{teal}\line(1,0){15}}
\put(0,75){\color{teal}\line(1,0){15}}
\put(15,60){\color{teal}\line(1,0){15}}
\put(0,45){\color{teal}\line(1,0){15}}
}

\put(0,0){\circle*{3}}
\put(1,-7){$a$}
\put(-10,0){$\alpha$}
\put(15,15){\circle*{3}}
\put(16,8){$o$}
\put(30,30){\circle*{3}}
\put(31,22){$b$}
\put(45,45){\circle*{3}}
\put(46,37){$x$}
\put(60,60){\circle*{3}}
\put(61,53){$y$}
\put(75,75){\circle*{3}}
\put(76,68){$s$}
\put(90,90){\circle*{3}}
\put(91,83){$t$}

\put(0,15){\circle*{3}}
\put(-14,13){$ao$}
\put(0,45){\circle*{3}}
\put(-14,43){$ax$}
\put(0,75){\circle*{3}}
\put(-14,73){$as$}

\put(15,30){\circle*{3}}
\put(3,29){$ob$}
\put(15,45){\circle*{3}}
\put(16,38){$ox$}
\put(15,60){\circle*{3}}
\put(3,60){$oy$}
\put(15,75){\circle*{3}}
\put(16,68){$os$}
\put(15,90){\circle*{3}}
\put(10,93){$ot$}

\put(30,60){\circle*{3}}
\put(32,62){$by$}
\put(30,90){\circle*{3}}
\put(27,93){$bt$}

\end{picture}

We claim that $\Aline {ao}{ob}\parallel \Delta$. Let $\alpha\in\Delta$ be a unique point such that $\Aline {\alpha o}{ob}\parallel\Delta$. Since $\Aline ob\nparallel \Delta$, the point $\alpha$ is not equal to the point $o$. Applying the condition (2), we conclude that $\Aline{\alpha x}{oy}\parallel\Delta$. Since $\Aline {ax}{oy}\parallel \Delta$, the Proclus Axiom implies $\Aline{ax}{oy}=\Aline{\alpha x}{oy}$ and hence $\{ax\}=\Aline{ax}{oy}\cap\Aline {ox}x=\Aline{\alpha x}{oy}\cap\Aline {ox} x=\{\alpha x\}$. Then $\alpha=a$ and $\Aline {ao}{ob}=\Aline{\alpha o}{ob}\in\Delta_\parallel$. Applying the condition (2) to the points $a,b,s,t$, we conclude that $\Aline {as}{ot}\subparallel \Delta$. Since $a\ne o$, the points $as$ and $ot$ are distinct. In this case $\Aline{as}{ot}\subparallel\Delta$ implies $\Aline{as}{ot}\parallel\Delta$, by Corollary~\ref{c:subparallel}.  
\smallskip

$(3)\Ra(2)$ Assume that the condition (3) holds. To prove the condition (2), take any points $x,y,z,s\in\Delta$ such that the flats $\Aline{xo}{oy}$ and $\Aline {oz}{ys}$ are subparallel to the line $\Delta$. We have to show that $\Aline {os}{xz}\subparallel\Delta$. If $x=o$, then $\Aline {oo}{oy}=\Aline {xo}{oy}\subparallel \Delta$ implies $y=o$. In this case $\Aline {oz}{os}=\Aline {oz}{ys}\subparallel\Delta$ implies $z=s$ and hence $\Aline{oz}{ys}=\{oz\}\subparallel \Delta$. If $x\ne o$, then $\Aline {xo}{oy}\subparallel \Delta$ implies $o\ne y$ and hence $\Aline {oz}{ys},\Aline {xo}{oy}\in\Delta_\parallel$, by Corollary~\ref{c:subparallel}. Applying the condition (2) to the points $x,y,o,y,z,s$, we conclude that $\Aline {os}{xz}\parallel \Delta$.
\end{proof}

\begin{proposition}\label{p:add-ass<=>}  For a based affine plane $(\Pi,uow)$ and its ternar $\Delta$,\\ the following conditions are equivalent:
\begin{enumerate}
\item the plus loop $(\Delta,+)$ is associative;
\item $\forall x,y,p,q,s\in \Delta\;\;(\Aline {oy}{xp}\subparallel\Delta\;\wedge\; \Aline {yq}{ps}\subparallel\Delta) \;\Ra\;(\Aline {oq}{xs}\subparallel\Delta)$;
\item $\forall a,b,c,d,\alpha,\beta,\gamma,\delta\in \Delta\;\;(\Aline {ac}{\alpha\gamma},\Aline{bc}{\beta\gamma},\Aline {bd}{\beta\delta}\in\Delta_\parallel\;\Ra\;\Aline{ad}{\alpha\delta}\in\Delta_\parallel)$.
\begin{picture}(100,0)(-30,0)
\linethickness{0.5pt}
\put(0,0){\color{teal}\line(1,0){30}}
\put(0,0){\color{cyan}\line(0,1){15}}
\put(30,15){\color{cyan}\line(0,-1){15}}

\put(20,20){\color{teal}\line(1,0){30}}
\put(20,20){\color{cyan}\line(0,1){15}}
\put(50,35){\color{teal}\line(-1,0){30}}
\put(50,35){\color{cyan}\line(0,-1){15}}

\put(0,0){\color{red}\line(1,1){20}}
\put(30,15){\color{red}\line(1,1){20}}
\put(0,15){\color{red}\line(1,1){20}}
\put(30,0){\color{red}\line(1,1){20}}
\put(30,15){\color{teal}\line(-1,0){30}}

\put(0,0){\circle*{2}}
\put(30,0){\circle*{2}}
\put(0,15){\circle*{2}}
\put(30,15){\circle*{2}}

\put(20,20){\circle*{2}}
\put(50,20){\circle*{2}}
\put(20,35){\circle*{2}}
\put(50,35){\circle*{2}}
\end{picture}

\end{enumerate}
\end{proposition}

\begin{proof} $(1)\Ra(2)$ Assume that the loop $(\Delta,+)$ is associative and take any points $x,y,p,q,s\in \Delta$ with $\Aline {oy}{xp}\subparallel \Delta$ and $\Aline {yq}{ps}\subparallel\Delta$. By the definition of the plus operation in the ternar $\Delta$, the subparallelity relations $\Aline {oy}{xp}\subparallel \Delta$ implies $x+y=p$. 

Since $(\Delta,+)$ is a quasigroup, there exist unique points $z\in \Delta$ such that $s=p+z$ and hence $s=p+z=(x+y)+z=x+(y+z)$. 

The definition of the plus operation ensures that $\Aline {oz}{ps}\subparallel \Delta$. By the Proclus Axiom, the subparallelity relations $\Aline {yq}{ps}\subparallel\Delta$ and $\Aline {oz}{ps}\subparallel \Delta$ imply $\Aline {oz}{yq}\subparallel\Delta$ and hence $y+z=q$ and $s=x+(y+z)=x+q$, which implies $\Aline {oq}{xs}\subparallel\Delta$. 

\begin{picture}(150,175)(-130,-10)
\linethickness{0.75pt}
\put(113,105){\color{red}$\Delta$}
{\linethickness{1pt}
\put(0,0){\color{red}\line(1,1){150}}
}
\put(0,0){\color{cyan}\line(0,1){60}}
\put(0,15){\color{red}\line(1,1){135}}
\put(0,45){\color{red}\line(1,1){90}}
\put(0,60){\color{red}\line(1,1){90}}
\put(0,15){\color{teal}\line(1,0){15}}
\put(0,45){\color{teal}\line(1,0){45}}
\put(0,60){\color{teal}\line(1,0){60}}
\put(90,135){\color{teal}\line(1,0){45}}
\put(90,150){\color{teal}\line(1,0){60}}
\put(90,90){\color{cyan}\line(0,1){60}}
\put(45,45){\color{cyan}\line(0,1){15}}
\put(135,135){\color{cyan}\line(0,1){15}}

\put(45,60){\circle*{3}}
\put(37,65){$yq$}
\put(0,0){\circle*{3}}
\put(-8,-3){$o$}
\put(15,15){\circle*{3}}
\put(15,8){$z$}
\put(45,45){\circle*{3}}
\put(45,37){$y$}

\put(60,60){\circle*{3}}
\put(60,53){$q$}
\put(90,90){\circle*{3}}
\put(90,82){$x$}
\put(90,135){\circle*{3}}
\put(92,138){$xp$}
\put(90,150){\circle*{3}}
\put(85,153){$xs$}

\put(135,135){\circle*{3}}
\put(135,128){$p$}
\put(135,150){\circle*{3}}
\put(129,154){$ps$}

\put(150,150){\circle*{3}}
\put(150,142){$s$}
\put(0,15){\circle*{3}}
\put(-13,13){$oz$}
\put(0,45){\circle*{3}}
\put(-13,43){$oy$}
\put(0,60){\circle*{3}}
\put(-13,58){$oq$}

\end{picture}

$(2)\Ra(1)$ Now assume that the condition (2) is satisfied. To prove that the loop $(\Delta,+)$ is associative, take any points $x,y,z\in \Delta$ and consider the points $p\defeq x+y$, $q\defeq y+z$ and $s\defeq p+z=(x+y)+z$. 

The geometric definition of the plus operation ensures that the lines $\Aline {oy}{xp}$, $\Aline {oz}{yq}$ and $\Aline {oz}{ps}$ are subparallel to the diagonal $\Delta$. By the Proclus Axiom, the subparallelity relations 
$\Aline {oz}{yq}\subparallel\Delta$ and $\Aline {oz}{ps}\parallel\Delta$ imply $\Aline {yq}{ps}\subparallel\Delta$. Applying the condition (2), we conclude $\Aline {oq}{xs}\subparallel \Delta$ and hence $(x+y)+z=s=x+q=x+(y+z)$, witnessing that the loop $(\Delta,+)$ is associative.
\smallskip

$(2)\Ra(3)$ Assume that the condition (2) holds.

\begin{claim}\label{cl:add-ass} For any points $b,c,d,\alpha,\beta,\gamma,\delta\in\Delta$ with $\Aline {oc}{\alpha\gamma},\Aline{bc}{\beta\gamma},\Aline {bd}{\beta\delta}\in\Delta_\parallel$, we have $\Aline{od}{\alpha\delta}\in\Delta_\parallel$.
\end{claim}

\begin{proof}  First we show that $\Aline {cd}{\gamma\delta}\parallel \Delta$.
Find unique points $x,y\in\Delta$ such that $\Aline{xo}{bc}\subparallel\Delta$ and $\Aline {xy}{bd}\subparallel\Delta$. The condition (2) applied to the points $c,x,b,y,d$ implies $\Aline{oy}{cd}\subparallel\Delta$. Applying the condition (2) to the points $\gamma,x,\beta,y,\delta$, we conclude that $\Aline {oy}{\gamma\delta}\subparallel\Delta$. The subparallelity relations $\Aline{oy}{cd}\subparallel\Delta$ and $\Aline {oy}{\gamma\delta}\subparallel\Delta$ imply $\Aline {cd}{\gamma\delta}\subparallel\Delta$. Applying the condition (2) to the points $c,d,\alpha,\gamma,\delta$, we conclude that $\Aline{od}{\alpha\delta}\subparallel \Delta$. It follows from $\Aline{oc}{\alpha\gamma}\in\Delta_\parallel$ that $o\ne\alpha$. Then $od\ne\alpha\delta$ and $\Aline{od}{\alpha\delta}\subparallel \Delta$ implies $\Aline{od}{\alpha\delta}\parallel \Delta$, by Corollary~\ref{c:subparallel}.


\end{proof}

Now we can prove the condition (3). Take any points $a,b,c,d,\alpha,\beta,\gamma,\delta\in\Delta$ with $\Aline {ac}{\alpha\gamma}$, $\Aline{bc}{\beta\gamma}$, $\Aline {bd}{\beta\delta}\in\Delta_\parallel$. We have to prove that $\Aline{ad}{\alpha\delta}\in\Delta_\parallel$. If $a=o$ or $\alpha=o$, then $\Aline{ad}{\alpha\delta}\in\Delta_\parallel$, by Claim~\ref{cl:add-ass}. So, we assume that $a\ne o\ne \alpha$. In this case we can find unique points $y,x,z\in\Delta$ such that $\Aline {oy}{ac}\parallel \Delta$, $\Aline {xy}{bd}\parallel\Delta$ and $\Aline {xz}{cd}\parallel\Delta$. Claim~\ref{cl:add-ass} applied to the points $x,y,z,a,b,c,d$ implies $\Aline {oz}{ad}\parallel \Delta$. Claim~\ref{cl:add-ass} applied to the points $x,y,z,\alpha,\beta,\gamma,\delta$ implies $\Aline {oz}{\alpha\delta}\parallel \Delta$.  The parallelity of $\Aline{ac}{\alpha\gamma}$ to the line $\Delta$ implies $a\ne\alpha$. By the Proclus Axiom, the parallelity relations $\Aline {oz}{ad}\parallel \Delta$, $\Aline {oz}{\alpha\delta}\parallel \Delta$ and the inequality $ad\ne \alpha\delta$ imply $\Aline{ad}{\alpha\delta}\parallel\Delta$.
\smallskip

$(3)\Ra(2)$ Assume that the condition (3) holds. To check the condition (2), take any points $x,y,p,q,s\in\Delta$ with $\Aline {oy}{xp}\subparallel \Delta$ and $\Aline {yq}{ps}\subparallel \Delta$. We have to prove that $\Aline {oq}{xs}\subparallel \Delta$. If $o=x$, then the subparallelity relation $\Aline {oy}{xp}\subparallel\Delta$ implies $y=p$ and then the subparallelity $\Aline {yq}{ps}\subparallel \Delta$ implies $q=s$. Then $\Aline {oq}{xs}=\Aline{oq}{oq}=\{oq\}\subparallel \Delta$. By analogy, the equality $y=p$ implies $o=x$, $q=s$, and   $\Aline {oq}{xs}=\{oq\}\subparallel \Delta$.  So, we assume that $o\ne x$ and $y\ne p$. In this case the subparallelity relations $\Aline {oy}{xp}\subparallel \Delta$ and $\Aline {yq}{ps}\subparallel \Delta$ imply $\Aline {oy}{xp}\parallel \Delta$ and $\Aline {yq}{ps}\parallel \Delta$, by Corollary~\ref{c:subparallel}. Applying the condition (4) to the points $o,y,y,q,x,p,p,s$, we conclude that $\Aline {oq}{xs}\parallel\Delta$.
\end{proof}

\begin{proposition}\label{p:add-com<=>} For a based affine plane $\Pi$ and its ternar $\Delta$,\newline the following conditions are equivalent:
\begin{enumerate}
\item the plus loop $(\Delta,+)$ is commutative;
\item $\forall x,y,z\in \Delta\;\;(\Aline {ox}{yz}\subparallel\Delta\;\Leftrightarrow \;\Aline {oy}{xz}\subparallel \Delta)$. 
\begin{picture}(40,0)(-165,0)
\linethickness{0.5pt}
\put(0,0){\color{red}\line(1,1){40}}
\put(10,0){\color{red}\line(1,1){30}}
\put(0,0){\color{cyan}\line(0,1){30}}
\put(10,0){\color{cyan}\line(0,1){40}}
\put(0,30){\color{red}\line(1,1){10}}
\put(0,30){\color{teal}\line(1,0){40}}
\put(10,40){\color{teal}\line(1,0){30}}
\put(0,0){\color{teal}\line(1,0){10}}
\put(40,30){\color{cyan}\line(0,1){10}}

\put(10,10){\circle*{2}}

\put(0,0){\circle*{2}}
\put(40,40){\circle*{2}}
\put(40,30){\circle*{2}}
\put(0,30){\circle*{2}}
\put(10,40){\circle*{2}}
\put(10,0){\circle*{2}}

\end{picture}
\end{enumerate}
\end{proposition}

\begin{proof} $(1)\Ra(2)$ Assume that the loop $(\Delta,+)$ is commutative. Given any points $x,y,z\in \Delta$, we need to check that $\Aline {ox}{yz}\subparallel\Delta\;\Leftrightarrow \;\Aline {oy}{xz}\subparallel \Delta$. If $\Aline {ox}{yz}\subparallel\Delta$, then $y+x=z$, by the definition of plus operation in the ternar $\Delta$. Since the loop $(\Delta,+)$ is commutative, $z=y+x=x+y$ and hence $\Aline {oy}{xz}\subparallel\Delta$. By analogy we can prove that $\Aline {oy}{xz}\subparallel\Delta$ implies $\Aline {ox}{yz}\subparallel\Delta$.

\begin{picture}(150,100)(-150,-10)
\linethickness{0.75pt}
{\linethickness{1pt}
\put(0,0){\color{red}\line(1,1){80}}
}
\put(0,0){\color{cyan}\line(0,1){60}}
\put(0,60){\color{teal}\line(1,0){80}}
\put(0,0){\color{teal}\line(1,0){20}}
\put(20,0){\color{cyan}\line(0,1){80}}
\put(80,80){\color{teal}\line(-1,0){60}}
\put(0,60){\color{red}\line(1,1){20}}
\put(80,60){\color{cyan}\line(0,1){20}}
\put(20,0){\color{red}\line(1,1){60}}

\put(0,0){\circle*{3}}
\put(-8,-3){$x$}
\put(20,20){\circle*{3}}
\put(22,16){$o$}
\put(60,60){\circle*{3}}
\put(60,53){$z$}
\put(80,80){\circle*{3}}
\put(83,78){$y$}
\put(20,0){\circle*{3}}
\put(23,-3){$ox$}
\put(0,60){\circle*{3}}
\put(-13,58){$xz$}
\put(20,80){\circle*{3}}
\put(7,79){$oy$}
\put(80,60){\circle*{3}}
\put(83,58){$yz$}
\end{picture}

$(2)\Ra(1)$ Assume that the condition (2) holds. To prove that the loop $(\Delta,+)$ is commutative, take any points $x,y\in\Delta$ and consider the point $z\defeq x+y$. The definition of plus operation in the ternar $\Delta$ ensures that $\Aline {oy}{xz}\subparallel \Delta$ and hence $\Aline{ox}{yz}\subparallel \Delta$, by the condition (2). The latter subparallelity relation implies $z=y+x$ and hence $x+y=z=y+x$. 
\end{proof}

Lemma~\ref{l:Boolean} implies the following characterization.

\begin{proposition}\label{p:Boolean+<=>} For a based affine plane $\Pi$ and its ternar $\Delta$, the following conditions are equivalent:
\begin{enumerate}
\item the plus loop $(\Delta,+)$ is Boolean;
\item $\forall x\in\Delta \;\;\Aline {ox}{xo}\parallel\Aline{oo}{xx}$.
\end{enumerate}
\end{proposition}

Proposition~\ref{p:Ali-Slaney} and Lemma~\ref{l:Boolean} imply the following sufficient condition of existence of Boolean parallelograms in affine planes.

\begin{corollary}\label{c:inv-plus=>Boolean-paralelogram} Let $(\Pi,uow)$ be a based affine plane of even order, and $\Delta$ be its ternar. If the plus loop $(\Delta,+)$ is invertible, then there is an element $a\in \Delta\setminus\{o\}$ such that $a+a=o$ and hence the plane $\Pi$ contains a Boolean parallelogram.
\end{corollary}

\section{Group properties of the plus operation}

In this section we study the interplay between algebraic properties of the plus loop $(\Delta,+)$ of the ternar $\Delta$ of a based affine plane $(\Pi,uow)$ and properties of the subset $\Sym_\Pi^\#[\Delta;{\boldsymbol h},{\boldsymbol v}]$ in the group $\Sym_\Pi^\#(\Delta)$ of line translations of the diagonal $\Delta$ of the based affine plane $(\Pi,uow)$.

We recall that for a based affine plane $(\Pi,uow)$ we denote by $${\boldsymbol v}\defeq(\Aline ow)_\parallel\in\partial\Pi\quad\mbox{and}\quad \boldsymbol h\defeq(\Aline ou)_\parallel\in\partial\Pi$$the vertical and horizontal directions on the plane $\Pi$. For a non-diagonal direction $\boldsymbol \delta\in\partial\Pi\setminus\{\Delta_\parallel\}$ and two lines $D,D'\in\Delta_\parallel$, consider the line shift $$\boldsymbol \delta_{D',D}\defeq\{(x,y)\in D\times D':y\in\Aline x{\boldsymbol \delta}\}$$where $\Aline x{\boldsymbol \delta}$ denotes the unique line in the spread ${\boldsymbol \delta}$ that contains the point $x$. Given two non-diagonal directions ${\boldsymbol \delta},{\boldsymbol \delta'}\in \partial\Pi\setminus\{\Delta_\parallel\}$, consider the subset $$\Sym_\Pi^\#[\Delta;\boldsymbol \delta,\boldsymbol \delta']\defeq\{{\boldsymbol \delta}_{\Delta,D}{\boldsymbol \delta'}_{D,\Delta}:D\in\Delta_\parallel\}$$ of the group $\Sym_\Pi^\#(\Delta)$, introduced in Section~\ref{s:IX[L;u,v]}. 

\begin{theorem}\label{t:plus<=>group} Let $(\Pi,uow)$ be a based affine plane and $(\Delta,+)$ be the plus loop of its ternar. Let $\boldsymbol h\defeq(\Aline ou)_\parallel$ and $\boldsymbol v\defeq(\Aline ow)_\parallel$ be the horizontal and vertical directions of the affine base $uow$. 
\begin{enumerate}
\item For every point $a\in \Delta$, the right shift $R_a:\Delta\to\Delta$, $R_a:x\mapsto x+a$, of the loop $(\Delta,+)$ is equal to the permutation $\boldsymbol h_{\Delta,D}\boldsymbol v_{D,\Delta}$ where $D\in\Delta_\parallel$ is the unique line such that $\Aline o{\boldsymbol v}\cap\Aline{a}{\boldsymbol h}\subset D$.
\item For every line $L\in\Delta_\parallel$, the permutation $\boldsymbol h_{\Delta,L}\boldsymbol v_{L,\Delta}$ is equal to the right shift $R_a$ by the element $a\defeq\boldsymbol h_{\Delta,L}\boldsymbol v_{L,\Delta}(o)$. 
\item $\Sym^\#_\Pi[\Delta;\boldsymbol h,\boldsymbol v]=\{R_a:a\in\Delta\}$.
\item The plus loop $(\Delta,+)$ is right-inversive if and only if\\ $\Sym^\#_\Pi[\Delta;\boldsymbol h,\boldsymbol v]=\Sym^\#_\Pi[\Delta;\boldsymbol h,\boldsymbol v]^{-1}=\Sym^\#_\Pi[\Delta;\boldsymbol v,\boldsymbol h]$.
\item The loop $(\Delta,+)$ is associative if and only if $\Sym^\#_\Pi[\Delta;\boldsymbol h,\boldsymbol v]$ is a subgroup of the group $\Sym^\#_\Pi[\Delta]$ if and only if the function $\Phi:\Delta\to \Sym^\#_\Pi[\Delta;\boldsymbol h,\boldsymbol v]$, $\Phi:a\mapsto R_a^{-1}$, is an isomomorphism of the loops $(\Delta,+)$ and $\Sym^\#_\Pi[\Delta;\boldsymbol h,\boldsymbol v]$.
\end{enumerate}
\end{theorem}

\begin{proof} 1. Given any point $a\in \Delta$, consider the point $\alpha\in\Aline o{\boldsymbol v}\cap\Aline a{\boldsymbol h}$ with coordinates $oa$ and let $D\in\Delta_\parallel$ be the unique line containing the point $\alpha$. The definition of the plus operation ensures that $R_a=\boldsymbol h_{\Delta,D}\boldsymbol v_{D,\Delta}$.

\begin{picture}(60,100)(-160,-15)
{\linethickness{1pt}
\put(0,0){\color{red}
\line(1,1){70}}
}\linethickness{0.75pt}
\put(0,20){\color{red}\line(1,1){50}}
\put(0,0){\color{cyan}\vector(0,1){20}}
\put(0,20){\color{teal}\vector(1,0){20}}
\put(50,50){\color{cyan}\vector(0,1){20}}
\put(50,70){\color{teal}\vector(1,0){20}}
\put(35,25){\color{red}$\Delta$}
\put(20,50){\color{red}$D$}
\put(0,0){\color{cyan}\vector(0,1){70}}
\put(0,0){\color{teal}\vector(1,0){75}}
\put(80,-2){\color{teal}$\boldsymbol h$}
\put(-2,73){\color{cyan}$\boldsymbol v$}

\put(0,0){\circle*{3}}
\put(-2,-8){$o$}
\put(20,20){\circle*{3}}
\put(22,14){$a$}
\put(50,50){\circle*{3}}
\put(50,43){$x$}
\put(70,70){\circle*{3}}
\put(74,68){$x+a$}
\put(0,20){\circle*{3}}
\put(-9,18){$\alpha$}
\put(50,70){\circle*{3}}
\end{picture}

2. Given any line $L\in\Delta_\parallel$, consider the point $a\defeq\boldsymbol h_{\Delta,L}\boldsymbol v_{L,\Delta}(o)$. The preceding item implies that $\boldsymbol h_{\Delta,L}\boldsymbol v_{L,\Delta}=R_a$.
\smallskip

3. Two preceding items imply the equality $\Sym^\#_\Pi[\Delta;\boldsymbol h,\boldsymbol v]=\{R_a:a\in\Delta\}$.
\smallskip

4. If the loop $(\Delta,+)$ is right-inversive, then for every $a\in\Delta$ there exists a unique element $a^{-}\in\Delta$ such that $x=(x+a)+a^{-}=R_{a^{-}}R_a(x)$ for every $x\in \Delta$. Then $R_a^{-1}=R_{a^{-}}$ in the permutation group $\Sym^\#_\Pi(\Delta)\subseteq S_\Delta$. Then for the set $\Sym^\#_\Pi[\Delta;\boldsymbol h,\boldsymbol v]=\{R_a:a\in\Delta\}$, we have 
$$
\begin{aligned}
\Sym^\#_\Pi[\Delta;\boldsymbol h,\boldsymbol v]^{-1}&=\{R_a:a\in\Delta\}^{-1}=\{R^{-1}_a:a\in\Delta\}=\{R_{a^{-}}:a\in\Delta\}\\
&=\{R_a:a\in\Delta\}= \Sym^\#_\Pi[\Delta;\boldsymbol h,\boldsymbol v]
\end{aligned}
$$in the group $\Sym^\#_\Pi(\Delta)$.
Observe also that $$\Sym^\#_\Pi[\Delta;\boldsymbol h,\boldsymbol v]^{-1}=\{(\boldsymbol h_{\Delta,D}\boldsymbol v_{D,\Delta})^{-1}:D\in\Delta_\parallel\}=\{\boldsymbol v_{\Delta,D}\boldsymbol h_{D,\Delta}:D\in\Delta_\parallel\}=\Sym^\#_\Pi[\Delta;\boldsymbol v,\boldsymbol h].
$$

If $\Sym^\#_\Pi[\Delta;\boldsymbol h,\boldsymbol v]^{-1}=\Sym^\#_\Pi[\Delta;\boldsymbol h,\boldsymbol v]$, then 
$$\{R_a^{-1}:a\in\Delta\}=\Sym^\#_\Pi[\Delta;\boldsymbol h,\boldsymbol v]^{-1}=\Sym^\#_\Pi[\Delta;\boldsymbol h,\boldsymbol v]=\{R_a:a\in\Delta\}$$ and hence for every $a\in\Delta$, there exists an element $b\in\Delta$ such that $R_a^{-1}=R_b$ and hence $x=R_bR_a(x)=(x+a)+b$ for every $x\in\Delta$, witnessing that the plus loop $(\Delta,+)$ is right-inversive.
\smallskip

5. If the plus loop $(\Delta,+)$ is associative, then it is right-inversive and hence $\Sym^\#_\Pi[\Delta;\boldsymbol h,\boldsymbol v]^{-1}=\Sym^\#_\Pi[\Delta;\boldsymbol h,\boldsymbol v]$, by the preceding item. By the associativity of $(\Delta,+)$, for every elements $a,b,x\in\Delta$ we have $$R_bR_a(x)=(x+a)+b=x+(a+b)=R_{a+b}(x).$$
Then the bijective function $\Phi:\Delta\to \Sym^\#_\Pi[\Delta;\boldsymbol h,\boldsymbol v]$, $\Phi:a\mapsto R_a^{-1}$, is an isomorphism.
Indeed,
$$\Phi(a+b)=R_{a+b}^{-1}=(R_bR_a)^{-1}=R_a^{-1}R_b^{-1}=\Phi(a)\Phi(b).$$
Since $\Phi$ is a homomorphism, its image $\Phi[\Delta]=\Sym^\#_\Pi[\Delta;\boldsymbol h,\boldsymbol v]$ is a subgroup of the group $\Sym^\#_\Pi(\Delta)$.

If  $\Sym^\#_\Pi[\Delta;\boldsymbol h,\boldsymbol v]=\{R_a:a\in\Delta\}$ is a   subgroup of the group $\Sym^\#_\Pi(\Delta)$, then for every $a,b\in \Delta$ there exist an element $c\in \Delta$ such that $R_bR_a=R_c$ and hence $$x+c=R_c(x)=R_bR_a(x)=(x+a)+b.$$ In particular, $c=o+c=(o+a)+b=a+b$. Then $(x+a)+b=x+c=x+(a+b)$, witnessing that the loop $(\Delta,+)$ is associative.
\end{proof}

\section{Geometry of the puls operation}

In this section we find geometric counterparts of some algebraic properties of the puls operation $\puls$ of a ternar. Let us recall that for a based affine plane $(\Pi,uow)$ and any elements $a,b,c\in\Delta$ on its diagonal $\Delta$, the equality $a\puls b=c$ holds if and only if $\Aline{oo}{ea}\parallel\Aline {ob}{ec}$. 

\begin{picture}(120,155)(-150,-15)
\linethickness{0.8pt}
\put(0,0){\line(1,1){100}}
\put(0,0){\color{teal}\vector(1,0){120}}
\put(0,0){\color{cyan}\vector(0,1){120}}
\put(20,20){\color{cyan}\line(0,1){80}}
\put(40,40){\color{teal}\line(-1,0){20}}
\put(60,60){\color{teal}\line(-1,0){60}}
\put(100,100){\color{teal}\line(-1,0){80}}
\put(0,0){\color{red}\line(1,2){20}}
\put(0,60){\color{red}\line(1,2){20}}


\put(-3,123){\color{cyan}$\boldsymbol v$}
\put(124,-3){\color{teal}$\boldsymbol h$}

\put(0,0){\circle*{3}}
\put(-3,-8){$o$}
\put(20,20){\circle*{3}}
\put(21,13){$e$}
\put(40,40){\circle*{3}}
\put(41,33){$a$}
\put(60,60){\circle*{3}}
\put(64,55){$b$}
\put(20,40){\circle*{3}}
\put(8,41){$ea$}
\put(0,60){\circle*{3}}
\put(-13,58){$ob$}
\put(20,100){\circle*{3}}
\put(15,103){$ec$}
\put(100,100){\circle*{3}}
\put(104,99){$c=a\puls b$}

\end{picture}

Here for two points $x,y\in \Delta$ we denote by $xy$ the unique point with cordinates $x,y$. Therefore, we identify the based affine plane $\Pi$ with the coordinate plane $\Delta^2$ of the ternar $\Delta$. 

\begin{proposition}\label{p:two-sided+} For a based affine plane $(\Pi,uow)$ and its ternar $\Delta$,\newline the following conditions are equivalent:
\begin{enumerate}
\item the puls loop $(\Delta,\!\puls\!)$ is invertible;
\item $\forall x,y\in \Delta\;\;(\Aline {oo}{ex}\parallel\Aline{oy}{eo}\;\Ra\;\Aline {oo}{ey}\parallel \Aline{ox}{eo})$;
\item $\forall x,y\in \Delta\;\;(\Aline {oo}{ex}\parallel\Aline{oy}{eo}\;\Leftrightarrow\; \Aline{oo}{ey}\parallel\Aline {ox}{eo})$.

\end{enumerate}
\end{proposition} 
\smallskip

\begin{proof}  $(1)\Ra(2)$ Assume that the loop $(\Delta,\!\puls\!)$ is invertible.  Given any points $x,y\in \Delta$ with $\Aline {oo}{ex}\parallel\Aline {oy}{eo}$, we should prove that $\Aline {oo}{ey}\parallel \Aline {ox}{eo}$. The parallelity relation  $\Aline {oo}{ex}\parallel\Aline{oy}{eo}$ implies $x\puls y=o$, by the definition of the puls operation in the ternar $\Delta$. Since the loop $(\Delta,\!\puls\!)$ is invertible, $y\puls x=o$, which implies $\Aline {oo}{oy}\parallel \Aline{ox}{eo}$. 

%

The implication $(2)\Ra(3)$ is trivial.
\smallskip

 $(3)\Ra(1)$ Assume that the condition (3) holds. To prove that the loop $(\Delta,!\puls!)$ is invertible, take any point $x\in \Delta$. Since $(\Delta,\!\puls\!)$ is a loop, there exists an element $y\in\Delta$ such that $x\puls y=o$. By the definition of the puls operation in the ternar $\Delta$, the equality $x\puls y=o$ implies $\Aline {oo}{ex}\parallel \Aline{oy}{eo}$. Applying the condition (3), we conclude that $\Aline{oo}{oy}\parallel\Aline{ox}{eo}$. Then $y\puls x=o$, by the  definition of the plus operation in the ternar $\Delta$. Therefore, $y$ is the two-sided inverse to $x$ in the loop $(\Delta,\!\puls\!)$, and the loop $(\Delta,\!\puls\!)$ is invertible.
\end{proof}
\smallskip

\begin{proposition}\label{p:lip+<=>}   For a based affine plane $(\Pi,uow)$ and its ternar $\Delta$,\newline the following conditions are equivalent:
\begin{enumerate}
\item the puls loop $(\Delta,\!\puls\!)$ is left-inversive;
\item $\forall x,z,s\in \Delta\;\;(\Aline {oo}{ex}\parallel\Aline {oz}{es}\;\Ra\;\Aline {ox}{eo}\parallel \Aline{os}{ez})$.  
\begin{picture}(100,0)(-160,0)
\linethickness{0.5pt}
\put(0,0){\color{cyan}\line(0,1){40}}
\put(20,0){\color{cyan}\line(0,1){40}}
\put(0,0){\color{teal}\line(1,0){20}}
\put(0,10){\color{teal}\line(1,0){20}}
\put(0,30){\color{teal}\line(1,0){20}}
\put(0,40){\color{teal}\line(1,0){20}}
\put(0,0){\color{blue}\line(2,1){20}}
\put(0,30){\color{blue}\line(2,1){20}}
\put(0,10){\color{red}\line(2,-1){20}}
\put(0,40){\color{red}\line(2,-1){20}}

\put(0,0){\circle*{2}}
\put(20,0){\circle*{2}}
\put(0,10){\circle*{2}}
\put(20,10){\circle*{2}}

\put(0,30){\circle*{2}}
\put(20,30){\circle*{2}}
\put(0,40){\circle*{2}}
\put(20,40){\circle*{2}}

\end{picture}

\end{enumerate}
\end{proposition}

\begin{proof} $(1)\Ra(2)$ Assume that the loop $(\Delta,\!\puls\!)$  is left-inversive.  To prove the condition (2), take any points $x,z,s\in \Delta$ with $\Aline {oo}{ex}\parallel\Aline{oz}{es}$. We have to prove that $\Aline {ox}{eo}\parallel \Aline{os}{ez}$. Since the loop $(\Delta,\!\puls\!)$ is left-inversive, there exists an element $y\in\Delta$ such that $y\puls (x\puls a)=a$ for every $a\in \Delta$. In particular, $y\puls x=y\puls(x\puls o)=o$ and $y\puls (x\puls z)=z$. By definition of the puls operation, the equality $y\puls x=o$ implies $\Aline {oo}{ey}\parallel \Aline{ox}{eo}$. By definition of the puls operation, the parallelity relation $\Aline{oo}{ex}\parallel \Aline{oz}{es}$ implies $x\puls z=s$. Then $y\puls s=y\puls (x\puls z)=z$. By  definition of the plus operation, the equality $y\puls s=z$ implies $\Aline {oo}{ey}\parallel \Aline {os}{ez}$. Then $\Aline {ox}{eo}\parallel \Aline {oo}{ey}\parallel \Aline{os}{ez}$ and hence $\Aline {ox}{eo}\parallel \Aline{os}{ez}$.

\begin{picture}(140,135)(-140,-35)
\linethickness{0.8pt}
\put(-20,-20){\line(1,1){100}}
\put(-20,-20){\color{teal}\line(1,0){60}}
\put(0,0){\color{cyan}\line(0,1){80}}
\put(40,-20){\color{cyan}\line(0,1){100}}
\put(0,0){\color{teal}\line(1,0){40}}
\put(0,20){\color{teal}\line(1,0){40}}
\put(0,60){\color{teal}\line(1,0){60}}
\put(0,80){\color{teal}\line(1,0){80}}
\put(0,0){\color{blue}\line(2,1){40}}
\put(0,60){\color{blue}\line(2,1){40}}
\put(0,0){\color{red}\line(2,-1){40}}
\put(0,80){\color{red}\line(2,-1){40}}
\put(0,20){\color{red}\line(2,-1){40}}

\put(-20,-20){\circle*{3}}
\put(-23,-28){$y$}
\put(40,-20){\circle*{3}}
\put(35,-28){$ey$}
\put(0,0){\circle*{3}}
\put(-8,-2){$o$}
\put(40,0){\circle*{3}}
\put(43,-2){$eo$}
\put(0,20){\circle*{3}}
\put(-12,18){$ox$}
\put(20,20){\circle*{3}}
\put(16,23){$x$}
\put(40,20){\circle*{3}}
\put(42,18){$ex$}
\put(40,40){\circle*{3}}
\put(43,36){$e$}
\put(60,60){\circle*{3}}
\put(63,56){$z$}
\put(0,60){\circle*{3}}
\put(-12,58){$oz$}
\put(40,60){\circle*{3}}
\put(42,62){$ez$}
\put(80,80){\circle*{3}}
\put(79,83){$s$}
\put(0,80){\circle*{3}}
\put(-5,83){$os$}
\put(40,80){\circle*{3}}
\put(36,83){$es$}
\end{picture}

$(2)\Ra(1)$ Assume that the condition (2) is satisfied. To prove that the loop $(\Delta,\!\puls\!)$ is left-inversive, take any element $x\in\Delta$. Since $(\Delta,\!\puls\!)$ is a quasigroup, there exists an element $y\in\Delta$ such that $y\puls x=o$ and hence $\Aline{oo}{ey}\parallel\Aline{ox}{eo}$. We claim that $y\puls (x\puls z)=z$ for every $z\in \Delta$. Consider  the point $s\defeq x\puls z$ and observe that $\Aline {oo}{ex}\parallel\Aline{oz}{es}$. Applying the condition (2), we conlcude that $\Aline{ox}{eo}\parallel\Aline{os}{ez}$. 
Then $\Aline {os}{ez}\parallel \Aline{ox}{eo}\parallel\Aline{oo}{ey}$ and hence $y\puls (x\puls z)=y\puls s=z$, by definition of the puls operation.
\end{proof}

\begin{proposition}\label{p:rip+<=>}  For a based affine plane $(\Pi,uow)$ and its ternar $\Delta$,\newline the following conditions are equivalent:
\begin{enumerate}
\item the puls loop $(\Delta,\!\puls\!)$ is right-inversive;
\item $\forall x,y,z,s\in \Delta\;\;(\Aline {oo}{ex}\parallel\Aline {oy}{eo}\;\wedge\;\Aline {oo}{ez}\parallel \Aline {ox}{es} \;\Ra\;\Aline {oo}{es}\parallel \Aline{oy}{ez})$.
\begin{picture}(100,0)(-75,-20)
\linethickness{0.5pt}
\put(0,20){\color{teal}\line(1,0){20}}
\put(0,-20){\color{cyan}\line(0,1){40}}
\put(20,-10){\color{cyan}\line(0,1){30}}
\put(0,0){\color{teal}\line(1,0){20}}

\put(0,0){\color{violet}\line(2,-1){20}}
\put(0,20){\color{violet}\line(2,-1){20}}
\put(0,0){\color{red}\line(2,1){20}}
\put(0,-20){\color{red}\line(2,1){20}}

\put(0,-20){\color{blue}\line(1,1){20}}
\put(0,0){\color{blue}\line(1,1){20}}

\put(0,0){\circle*{2}}
\put(20,0){\circle*{2}}
\put(0,20){\circle*{2}}
\put(20,20){\circle*{2}}
\put(0,-20){\circle*{2}}
\put(20,-10){\circle*{2}}
\put(20,10){\circle*{2}}

\end{picture}
\end{enumerate}
\end{proposition}

\begin{proof} $(1)\Ra(2)$ Assume that the loop $(\Delta,\!\puls\!)$ is right-inversive. Given any points $x,y,z,s\in\Delta$ with 
$\Aline {oo}{ex}\parallel\Aline {oy}{eo}$ and $\Aline {oo}{ez}\parallel \Aline {ox}{es}$, we should prove that $\Aline {oo}{es}\parallel \Aline{oy}{ez}$.
By definition of the puls operation, $\Aline {oo}{ex}\parallel\Aline {oy}{eo}$ implies $x\puls y=o$, and  $\Aline {oo}{ez}\parallel \Aline {ox}{es}$ implies $z\puls x=s$. Since the loop $(\Delta,\!\puls\!)$ is right-inversive, 
$s\puls y=(z\puls x)\puls y=z$ and hence $\Aline {oo}{es}\parallel \Aline {oy}{ez}$.

\begin{picture}(150,145)(-170,-72)
\linethickness{0.8pt}
\put(-60,-60){\line(1,1){120}}
\put(0,-60){\color{cyan}\line(0,1){120}}
\put(40,-40){\color{cyan}\line(0,1){100}}
\put(-60,-60){\color{teal}\line(1,0){60}}
\put(-40,-40){\color{teal}\line(1,0){80}}
\put(0,0){\color{teal}\line(1,0){40}}
\put(0,60){\color{teal}\line(1,0){60}}
\put(0,-60){\color{red}\line(2,1){40}}
\put(0,-60){\color{blue}\line(2,3){40}}
\put(0,0){\color{red}\line(2,1){40}}
\put(0,0){\color{blue}\line(2,3){40}}
\put(0,0){\color{violet}\line(1,-1){40}}
\put(0,60){\color{violet}\line(1,-1){40}}
\put(20,20){\color{teal}\line(1,0){20}}

\put(-60,-60){\circle*{3}}
\put(-65,-68){$y$}
\put(0,-60){\circle*{3}}
\put(-5,-68){$oy$}
\put(-40,-40){\circle*{3}}
\put(-48,-41){$z$}
\put(40,-40){\circle*{3}}
\put(43,-42){$ez$}
\put(0,0){\circle*{3}}
\put(-8,-2){$o$}
\put(40,0){\circle*{3}}
\put(43,-2){$eo=u$}
\put(20,20){\circle*{3}}
\put(18,13){$s$}
\put(40,20){\circle*{3}}
\put(43,18){$es$}
\put(40,40){\circle*{3}}
\put(43,37){$e$}
\put(0,60){\circle*{3}}
\put(-5,63){$ox$}
\put(40,60){\circle*{3}}
\put(35,63){$ex$}
\put(60,60){\circle*{3}}
\put(58,63){$x$}
\end{picture}

$(2)\Ra(1)$ Assume that the condition (2) holds. To prove that the loop $(\Delta,\puls)$ is right-inversive, take any $x\in\Delta$. Since $(\Delta,\puls)$ is a loop, there exists a unique point $y\in\Delta$ such that $x\puls y=o$. By the definition of the puls operation, $x\puls y=o$ implies $\Aline {oo}{ox}\parallel\Aline{oy}{eo}$. We claim that $(z\puls x)\puls y=z$ for every $z\in \Delta$. Consider the point $s\defeq z\puls x$. By the definition of the puls operation, the equality $s=z\puls x$ implies $\Aline {oo}{ez}\parallel \Aline {ox}{es}$. By the conditions (2), the parrallelity relations $\Aline {oo}{ox}\parallel \Aline{oy}{eo}$ and $\Aline {oo}{ez}\parallel \Aline{ox}{es}$ imply $\Aline {oo}{es}\parallel\Aline{oy}{ez}$ and hence $(z\puls x)\puls y=s\puls y=z$, by definition of the puls operation.
\end{proof}

\begin{proposition}\label{p:add-ass<=>+}  For a based affine plane $(\Pi,uow)$ and its ternar $\Delta$,\\ the following conditions are equivalent:
\begin{enumerate}
\item the puls loop $(\Delta,\!\puls\!)$ is associative;
\item $\forall x,y,a,b,c\in \Delta\;\;(\Aline {ox}{ey}\parallel\Aline {oa}{eb}\;\wedge\; \Aline {oo}{ex}\parallel\Aline {oc}{ea}) \;\Ra\;(\Aline {oo}{ey}\parallel\Aline{oc}{eb})$;
\item $\forall a,b,c,\alpha,\beta,\gamma\in \Delta\;\;(\Aline {oa}{eb}\parallel\Aline {o\alpha}{e\beta}\;\wedge\; \Aline {oc}{ea}\parallel\Aline {o\gamma}{e\alpha}) \;\Ra\;(\Aline {oc}{eb}\parallel\Aline{o\gamma}{e\beta})$.
\begin{picture}(20,0)(-45,-3)
\linethickness{0.5pt}
\put(0,0){\color{cyan}\line(0,1){50}}
\put(20,0){\color{cyan}\line(0,1){40}}

\put(0,0){\color{teal}\line(1,0){20}}
\put(0,0){\color{blue}\line(2,1){20}}
\put(0,20){\color{violet}\line(1,-1){20}}
\put(0,20){\color{red}\line(2,-1){20}}

\put(0,30){\color{teal}\line(1,0){20}}
\put(0,30){\color{blue}\line(2,1){20}}
\put(0,50){\color{violet}\line(1,-1){20}}
\put(0,50){\color{red}\line(2,-1){20}}

\put(0,0){\circle*{2}}
\put(20,0){\circle*{2}}
\put(0,20){\circle*{2}}
\put(20,10){\circle*{2}}

\put(0,30){\circle*{2}}
\put(20,30){\circle*{2}}
\put(0,50){\circle*{2}}
\put(20,40){\circle*{2}}

\end{picture}
\end{enumerate}
\end{proposition}

\begin{proof} $(1)\Ra(2)$ Assume that the loop $(\Delta,\!\puls\!)$ is associative and take any points $x,y,a,b,c\in \Delta$ with $\Aline {ox}{ey}\parallel\Aline {oa}{eb}$ and $\Aline {oo}{ex}\parallel\Aline {oc}{ea}$. We have to prove that $\Aline {oo}{ey}\parallel\Aline{oc}{eb}$. Since $(\Delta,\!\puls\!)$ is a quasigroup, there exists a unique element $z\in\Delta$ such that $z\puls x=y$. By definition of the puls operation, the equality $z\puls x=y$ implies $\Aline {oo}{ez}\parallel \Aline {ox}{ey}$. Then $\Aline{oa}{eb}\parallel \Aline{ox}{ey}\parallel \Aline{oo}{ez}$, which implies $z\puls a=b$. On the other hand, the parallelity relation  $\Aline {oo}{ex}\parallel\Aline {oc}{ea}$ implies $x\puls c=a$. The associativity of the loop $(\Delta,\!\puls\!)$ ensures that $b=z\puls a=z\puls(x\puls c)=(z\puls x)\puls c=y\puls c$ and hence $\Aline{oo}{ey}\parallel\Aline{oc}{eb}$. 

\begin{picture}(150,170)(-160,-55)
\linethickness{0.75pt}
\put(-40,-40){\line(1,1){140}}
\put(-40,-40){\color{teal}\line(1,0){80}}
\put(-20,-20){\color{teal}\line(1,0){60}}
\put(20,20){\color{teal}\line(1,0){20}}
\put(60,60){\color{teal}\line(-1,0){60}}
\put(80,80){\color{teal}\line(-1,0){40}}
\put(100,100){\color{teal}\line(-1,0){100}}
\put(0,-40){\color{cyan}\line(0,1){140}}
\put(40,-40){\color{cyan}\line(0,1){120}}
\put(0,-40){\color{blue}\line(2,1){40}}
\put(0,0){\color{blue}\line(2,1){40}}
\put(0,60){\color{blue}\line(2,1){40}}
\put(0,0){\color{violet}\line(1,-1){40}}
\put(0,100){\color{violet}\line(1,-1){40}}
\put(0,0){\color{red}\line(2,-1){40}}
\put(0,100){\color{red}\line(2,-1){40}}

\put(-40,-40){\circle*{3}}
\put(-43,-48){$x$}
\put(0,-40){\circle*{3}}
\put(-5,-48){$ox$}
\put(40,-40){\circle*{3}}
\put(35,-48){$ex$}
\put(-20,-20){\circle*{3}}
\put(-28,-20){$y$}
\put(40,-20){\circle*{3}}
\put(43,-22){$ey$}
\put(0,0){\circle*{3}}
\put(-7,0){$o$}
\put(20,20){\circle*{3}}
\put(13,20){$z$}
\put(40,20){\circle*{3}}
\put(43,18){$ez$}
\put(40,40){\circle*{3}}
\put(33,40){$e$}
\put(60,60){\circle*{3}}
\put(63,55){$a$}
\put(0,60){\circle*{3}}
\put(-13,58){$oa$}
\put(40,60){\circle*{3}}
\put(42,62){$ea$}
\put(80,80){\circle*{3}}
\put(83,75){$b$}
\put(40,80){\circle*{3}}
\put(40,83){$eb$}
\put(100,100){\circle*{3}}
\put(102,95){$c$}
\put(0,100){\circle*{3}}
\put(-5,103){$oc$}
\end{picture}

$(2)\Ra(1)$ Assume that the condition (2) holds. Given any elements $x,z,c\in\Delta$, we should check that $(z\puls x)\puls c=z\puls(x\puls c)$. Consider the points $a\defeq x\puls c$, $y\defeq z\puls x$,  and $b\defeq z\puls a$. The definition of the puls operation ensures that $\Aline {oo}{ex}\parallel \Aline {oc}{ea}$, $\Aline {oo}{ez}\parallel \Aline {ox}{ey}$ and $\Aline{oo}{ez}\parallel \Aline{oa}{eb}$. Then $\Aline{ox}{ey}\parallel \Aline {oo}{ez}\parallel \Aline {oa}{eb}$.  Applying the condition (2), we conclude that $\Aline {oo}{ey}\parallel \Aline {oc}{eb}$, which implies $y\puls c=b$, by definition of the puls operation. Then $(z\puls x)\puls c=y\puls c=b=z\puls a=z\puls(x\puls z)$.
\smallskip

$(2)\Ra(3)$ Assume that the condition (2) holds. Given any points $a,b,c,\alpha,\beta,\gamma\in \Delta$ with $\Aline {oa}{eb}\parallel\Aline {o\alpha}{e\beta}$ and $\Aline {oc}{ea}\parallel\Aline {o\gamma}{e\alpha}$, we should prove $\Aline {oc}{eb}\parallel\Aline{o\gamma}{e\beta}$. Since the affine plane $\Pi$ is Playfair, there exist unique points $x,y\in\Delta$ such that $\Aline{oo}{ex}\parallel \Aline {oc}{ea}\parallel \Aline{o\gamma}{e\alpha}$ and $\Aline{ox}{ey}\parallel \Aline {oa}{eb}\parallel \Aline{o\alpha}{e\beta}$. Applying the condition (2), we conclude that $\Aline {oc}{eb}\parallel \Aline {oo}{ey}\parallel \Aline{o\gamma}{e\beta}$ and hence $\Aline{oc}{eb}\parallel \Aline{o\gamma}{e\beta}$, by the transitivity of the parallelity relation.

\begin{picture}(150,175)(-160,-55)
\linethickness{0.75pt}
\put(-40,-40){\line(1,1){140}}
\put(-40,-40){\color{teal}\line(1,0){80}}
\put(-20,-20){\color{teal}\line(1,0){60}}
\put(60,60){\color{teal}\line(-1,0){60}}
\put(80,80){\color{teal}\line(-1,0){40}}
\put(100,100){\color{teal}\line(-1,0){100}}
\put(0,10){\color{teal}\line(1,0){40}}
\put(30,30){\color{teal}\line(1,0){10}}
\put(0,50){\color{teal}\line(1,0){50}}
\put(0,-40){\color{cyan}\line(0,1){140}}
\put(40,-40){\color{cyan}\line(0,1){120}}
\put(0,-40){\color{blue}\line(2,1){40}}
\put(0,60){\color{blue}\line(2,1){40}}
\put(0,10){\color{blue}\line(2,1){40}}
\put(0,0){\color{violet}\line(1,-1){40}}
\put(0,100){\color{violet}\line(1,-1){40}}
\put(0,50){\color{violet}\line(1,-1){40}}
\put(0,0){\color{red}\line(2,-1){40}}
\put(0,100){\color{red}\line(2,-1){40}}
\put(0,50){\color{red}\line(2,-1){40}}

\put(-40,-40){\circle*{3}}
\put(-43,-48){$x$}
\put(0,-40){\circle*{3}}
\put(-5,-48){$ox$}
\put(40,-40){\circle*{3}}
\put(35,-48){$ex$}
\put(-20,-20){\circle*{3}}
\put(-28,-20){$y$}
\put(40,-20){\circle*{3}}
\put(43,-22){$ey$}
\put(0,0){\circle*{3}}
\put(-7,0){$o$}
\put(40,40){\circle*{3}}
\put(33,40){$e$}
\put(60,60){\circle*{3}}
\put(63,55){$a$}
\put(0,60){\circle*{3}}
\put(-13,58){$oa$}
\put(40,60){\circle*{3}}
\put(42,62){$ea$}
\put(80,80){\circle*{3}}
\put(83,75){$b$}
\put(40,80){\circle*{3}}
\put(42,83){$eb$}
\put(100,100){\circle*{3}}
\put(102,95){$c$}
\put(0,100){\circle*{3}}
\put(-13,98){$oc$}

\put(10,10){\circle*{3}}
\put(12,3){$\alpha$}
\put(0,10){\circle*{3}}
\put(-14,9){$o\alpha$}
\put(40,10){\circle*{3}}
\put(42,8){$e\alpha$}
\put(30,30){\circle*{3}}
\put(22,28){$\beta$}
\put(40,30){\circle*{3}}
\put(43,28){$e\beta$}
\put(50,50){\circle*{3}}
\put(52,45){$\gamma$}
\put(0,50){\circle*{3}}
\put(-13,47){$o\gamma$}
\end{picture}

\smallskip

$(3)\Ra(2)$ Assume that the condition (3) holds. Given any points $x,y,a,b,c\in \Delta$ with $\Aline {ox}{ey}\parallel\Aline {oa}{eb}$ and $\Aline {oo}{ex}\parallel\Aline {oc}{ea}$, we need to prove that $\Aline {oc}{eb}\parallel\Aline{oo}{ey}$. Consider the points $\alpha\defeq x$, $\beta\defeq y$, $\gamma\defeq o$, and observe that  $\Aline{o\alpha}{e\beta}=\Aline {ox}{ey}\parallel\Aline {oa}{eb}$ and $\Aline{o\gamma}{e\alpha}=\Aline {oo}{ex}\parallel\Aline {oc}{ea}$. Applying the condition (3) to the points $a,b,c,\alpha,\beta,\gamma$, we conclude that $\Aline{oc}{eb}\parallel \Aline{o\gamma}{e\beta}=\Aline{oo}{ey}$.
\end{proof}

\begin{proposition}\label{p:commutative-puls<=>} For a based affine plane $\Pi$ and its ternar $\Delta$,\newline the following conditions are equivalent:
\begin{enumerate}
\item the puls loop $(\Delta,\!\puls\!)$ is commutative;
\item $\forall x,y,z\in \Delta\;\;(\Aline{oo}{ex}\parallel\Aline {oy}{ez}\;\Ra\;\Aline {oo}{ey}\parallel \Aline{ox}{ez})$.
\begin{picture}(40,0)(-158,0)
\linethickness{0.5pt}
\put(0,0){\color{blue}\line(1,1){20}}
\put(0,10){\color{blue}\line(1,1){20}}
\put(0,0){\color{cyan}\line(0,1){30}}
\put(20,0){\color{cyan}\line(0,1){30}}
\put(0,10){\color{red}\line(2,-1){20}}
\put(0,30){\color{red}\line(2,-1){20}}
\put(0,0){\color{teal}\line(1,0){20}}
\put(0,30){\color{teal}\line(1,0){20}}

\put(0,10){\circle*{2}}
\put(20,0){\circle*{2}}
\put(0,0){\circle*{2}}
\put(0,30){\circle*{2}}
\put(20,20){\circle*{2}}
\put(20,30){\circle*{2}}
\end{picture}
\end{enumerate}
\end{proposition}

\begin{proof} $(1)\Ra(2)$ Assume that the puls loop $(\Delta,\!\puls\!)$ is commutative. Given any points $x,y,z\in \Delta$ with $\Aline{oo}{ox}\parallel\Aline {oy}{ez}$, we need to check that $\Aline {oo}{ey}\parallel \Aline{ox}{ez}$. By the definition of the puls operation, $\Aline{oo}{ox}\parallel\Aline {oy}{ez}$ implies $x\puls y=z$. Since the loop $(\Delta,\!\puls\!)$ is commutative, $z=x\puls y=y\puls x$ and hence  $\Aline {oo}{ey}\parallel \Aline{ox}{ez}$.

\begin{picture}(150,125)(-150,-50)
\linethickness{0.8pt}
\put(-40,-40){\line(1,1){100}}
\put(0,-40){\color{cyan}\line(0,1){100}}
\put(40,-40){\color{cyan}\line(0,1){100}}
\put(0,-40){\color{blue}\line(2,3){40}}
\put(0,0){\color{blue}\line(2,3){40}}
\put(-40,-40){\color{teal}\line(1,0){80}}
\put(0,60){\color{teal}\line(1,0){60}}
\put(0,0){\color{red}\line(1,-1){40}}
\put(0,60){\color{red}\line(1,-1){40}}

\put(-40,-40){\circle*{3}}
\put(-44,-48){$y$}
\put(0,-40){\circle*{3}}
\put(-5,-48){$oy$}
\put(40,-40){\circle*{3}}
\put(35,-48){$ey$}
\put(0,0){\circle*{3}}
\put(-8,-2){$o$}
\put(20,20){\circle*{3}}
\put(22,15){$z$}
\put(40,20){\circle*{3}}
\put(43,18){$ez$}
\put(40,40){\circle*{3}}
\put(42,35){$e$}
\put(0,60){\circle*{3}}
\put(-5,63){$ox$}
\put(40,60){\circle*{3}}
\put(35,63){$ex$}
\put(60,60){\circle*{3}}
\put(60,63){$x$}
\end{picture}

$(2)\Ra(1)$ Assume that the condition (2) holds. To prove that the loop $(\Delta,\!\puls\!)$ is commutative, take any points $x,y\in\Delta$ and consider the point $z\defeq x\puls y$. The definition of puls operation ensures that $\Aline {oo}{ox}\parallel \Aline{oy}{ez}$ and hence $\Aline{oo}{ey}\parallel \Aline{ox}{ez}$, by the condition (2). The latter parallelity relation implies $z=y\puls x$ and hence $x\puls y=z=y\puls x$. 
\end{proof}

Now we present a geometric characterization of Boolean puls loops of  a ternar. The definition of the puls operation implies the following lemma.

\begin{lemma}\label{l:Boolean-puls} Let $(\Pi,uow)$ be a based affine plane and $\Delta$ be its ternar.\\ For every point $x\in\Delta$, the following conditions are equivalent:
\begin{enumerate}
\item $x\puls x=o$;
\item $\Aline{oo}{ex}\parallel \Aline {ox}{eo}$.
\end{enumerate}
\end{lemma}

Lemma~\ref{l:Boolean-puls} implies the following characterization.

\begin{proposition}\label{p:Boolean-puls<=>} For a based affine plane $\Pi$ and its ternar $\Delta$,\\ the following conditions are equivalent:
\begin{enumerate}
\item the puls loop $(\Delta,\!\puls\!)$ is Boolean;
\item $\forall x\in\Delta \;\;\Aline {oo}{ex}\parallel\Aline{ox}{eo}$.
\begin{picture}(20,0)(-265,0)
\put(0,0){\color{teal}\line(1,0){20}}
\put(0,20){\color{teal}\line(1,0){20}}
\put(0,0){\color{cyan}\line(0,1){20}}
\put(20,0){\color{cyan}\line(0,1){20}}
\put(0,0){\color{red}\line(1,1){20}}
\put(0,20){\color{red}\line(1,-1){20}}

\put(0,0){\circle*{2}}
\put(20,0){\circle*{2}}
\put(0,20){\circle*{2}}
\put(20,20){\circle*{2}}
\put(10,10){\color{red}\circle*{2.5}}
\put(10,10){\color{white}\circle*{1.5}}
\end{picture}
\end{enumerate}
\end{proposition}

Proposition~\ref{p:Ali-Slaney} and Lemma~\ref{l:Boolean-puls} imply the following sufficient condition of existence of Boolean parallelograms in affine planes.

\begin{corollary}\label{c:inv-puls=>Boolean-paralelogram} Let $(\Pi,uow)$ be a based affine plane of even order, and $\Delta$ be its ternar. If the puls loop $(\Delta,\!\puls\!)$ is invertible, then there is an element $a\in \Delta\setminus\{o\}$ such that $a\puls a=o$ and hence the plane $\Pi$ contains a Boolean parallelogram.
\end{corollary}

\section{Group properties of the puls operation}

In this section we study the interplay between algebraic properties of the puls loop $(\Delta,\!\puls\!)$ of the ternar $\Delta$ of a based affine plane $(\Pi,uow)$ and properties of the subset $\Sym_\Pi^\#[\Aline u{\boldsymbol v};{\boldsymbol h},\Aline o{\boldsymbol v}]$ in the group $\Sym_\Pi^\#(\Aline u{\boldsymbol v})$ of line translations of the vertical line $\Aline u{\boldsymbol v}$ in the based affine plane $(\Pi,uow)$.

We recall that for a based affine plane $(\Pi,uow)$ we denote by $${\boldsymbol v}\defeq(\Aline ow)_\parallel\in\partial\Pi\quad\mbox{and}\quad \boldsymbol h\defeq(\Aline ou)_\parallel\in\partial\Pi$$the vertical and horizontal directions on the plane $\Pi$. It will be convenient to denote the vertical lines $\Aline o{\boldsymbol v}$ and $\Aline u{\boldsymbol v}$ by $Y$ and $U$, respectively. Fot two lines $L,\Delta\notin\boldsymbol h$, let $\boldsymbol h_{\Lambda,L}:L\to\Lambda$ be the bijective map assigning to every point $x\in L$ the unique point $y\in\Lambda\cap\Aline x{\boldsymbol h}$.

In this section we study algebraic properties of the subset
$$\Sym^\#_\Pi[U;\boldsymbol h,Y]\defeq\{\boldsymbol x_{U,Y}\boldsymbol h_{Y,U}:\boldsymbol x\in\partial X\setminus \{\boldsymbol v\}\}$$
of the group $\Sym_\Pi^\#(U)$ of line translations of the line $U$. 
This subset was also considered in Section~\ref{s:IX[L;u,v]}.


\begin{theorem}\label{t:puls-group} Let $(\Pi,uow)$ be a based affine plane and $(\Delta,\!\puls\!)$ be the plus loop of its ternar.  
\begin{enumerate}
\item For every point $a\in \Delta$, the left shift $L_a:\Delta\to\Delta$, $L_a:x\mapsto a\puls x$, of the loop $(\Delta,\!\puls\!)$ is equal to the permutation $\boldsymbol h_{\Delta,U}\boldsymbol a_{U,Y}\boldsymbol h_{Y,\Delta}$ where $\boldsymbol a\defeq(\Aline o\alpha)_\parallel$ and $\alpha\defeq\boldsymbol h_{U,\Delta}(a)$.
\item For every direction $\boldsymbol a\in\partial \Pi\setminus \{\boldsymbol v\}$, the permutation $\boldsymbol h_{\Delta,U}\boldsymbol a_{U,Y}\boldsymbol h_{Y,\Delta}$ is equal to the left shift $L_a$ by the element $a\defeq\boldsymbol h_{\Delta,U}\boldsymbol a_{U,Y}(o)$, and hence the permutation 
$\boldsymbol a_{U,Y}\boldsymbol h_{Y,U}$ is equal to $\boldsymbol h_{U,\Delta}L_a\boldsymbol h_{\Delta,U}$. 
\item $\Sym^\#_\Pi[U;\boldsymbol h,Y]=\{\boldsymbol h_{U,\Delta}L_a\boldsymbol h_{\Delta,U}:a\in\Delta\}$.
\item The puls loop $(\Delta,\!\puls\!)$ is left-inversive if and only if\\ $\Sym^\#_\Pi[U;\boldsymbol h,Y]=\Sym^\#_\Pi[U;\boldsymbol h,Y]^{-1}$.
\item The loop $(\Delta,\!\puls\!)$ is associative if and only if $\Sym^\#_\Pi[U;\boldsymbol h,Y]$ is a subgroup of the group $\Sym^\#_\Pi(U)$ if and only if the function $\Phi:\Delta\to \Sym^\#_\Pi[U;\boldsymbol h,Y]$, $\Phi:a\mapsto\boldsymbol h_{\Delta,U}L_a\boldsymbol h_{\Delta,U}$, is an isomorphsim of the loops $(\Delta,\puls)$ and $\Sym^\#_\Pi[U;\boldsymbol h,Y]$. 
\end{enumerate}
\end{theorem}

\begin{proof} 1. Given any point $a\in \Delta$, consider the point $\alpha\in\Aline e{\boldsymbol v}\cap\Aline a{\boldsymbol h}$ with coordinates $ea$ and let $\boldsymbol a\defeq(\Aline o\alpha)_\parallel$ be the direction of the line $\Aline o\alpha$. The definition of the puls operation ensures that $L_a=\boldsymbol h_{\Delta,U}\boldsymbol a_{U,Y}\boldsymbol h_{Y,\Delta}$.

\begin{picture}(120,155)(-150,-15)
\linethickness{0.8pt}
\put(0,0){\line(1,1){120}}
\put(123,123){$\Delta$}
\put(0,0){\color{teal}\vector(1,0){120}}
\put(0,0){\color{cyan}\line(0,1){120}}
\put(20,0){\color{cyan}\line(0,1){120}}
\put(40,40){\color{teal}\line(-1,0){20}}
\put(60,60){\color{teal}\vector(-1,0){60}}
\put(20,100){\color{teal}\vector(1,0){80}}
\put(0,0){\color{red}\line(1,2){20}}
\put(0,60){\color{red}\vector(1,2){20}}


\put(-3,123){\color{cyan}$Y$}
\put(17,123){\color{cyan}$U$}

\put(124,-3){\color{teal}$\boldsymbol h$}

\put(0,0){\circle*{3}}
\put(-3,-8){$o$}
\put(20,20){\circle*{3}}
\put(21,13){$e$}
\put(40,40){\circle*{3}}
\put(41,33){$a$}
\put(60,60){\circle*{3}}
\put(64,55){$x$}
\put(20,40){\circle*{3}}
\put(8,41){$\alpha$}
\put(0,60){\circle*{3}}
\put(20,100){\circle*{3}}
\put(100,100){\circle*{3}}
\put(103,97){$a\puls x$}

\end{picture}

2. Given any direction $\boldsymbol a\in\partial \Pi\setminus\{\boldsymbol v\}$,  consider the points $\alpha\defeq\boldsymbol a_{U,Y}(o)$ and $a\defeq\boldsymbol h_{\Delta,U}(\alpha)$. The preceding item implies that $\boldsymbol h_{\Delta,U}\boldsymbol a_{U,Y}\boldsymbol h_{Y,\Delta}=L_a$ and hence
$$\boldsymbol a_{U,Y}\boldsymbol h_{Y,U}=\boldsymbol h_{U,\Delta}\boldsymbol h_{\Delta,U}\boldsymbol a_{U,Y}\boldsymbol h_{Y,\Delta}\boldsymbol h_{\Delta,U}
=\boldsymbol h_{U,\Delta}L_a\boldsymbol h_{\Delta,U}.$$

\smallskip

3. Two preceding items imply the equality $\Sym^\#_\Pi[U;\boldsymbol h,Y]=\{\boldsymbol h_{U,\Delta}L_a\boldsymbol h_{\Delta,U}:a\in\Delta\}$.
\smallskip

4. If the loop $(\Delta,+)$ is left-inversive, then for every $a\in\Delta$ there exists a unique element $a^{-}\in\Delta$ such that $x=a^{-}\puls (a\puls x)=L_{a^{-}}L_a(x)$ for every $x\in \Delta$. Then $L_a^{-1}=L_{a^{-}}$ in the permutation group $S_\Delta$. Then for the set $\Sym^\#_\Pi[U;\boldsymbol h,Y]=\{\boldsymbol h_{U,\Delta}L_a\boldsymbol h_{\Delta,U}:a\in\Delta\}$, we have 
$$
\begin{aligned}
\Sym^\#_\Pi[U;\boldsymbol h,Y]^{-1}&=\{(\boldsymbol h_{U,\Delta}L_a\boldsymbol h_{\Delta,U})^{-1}:a\in\Delta\}=\{\boldsymbol h_{U,\Delta}L^{-1}_a\boldsymbol h_{\Delta,U}:a\in\Delta\}\\
&=\{\boldsymbol h_{U,\Delta}L_{a^{-}}\boldsymbol h_{\Delta,U}:a\in\Delta\}=\{\boldsymbol h_{U,\Delta}L_a\boldsymbol h_{\Delta,U}:a\in\Delta\}=\Sym^\#_\Pi[U;\boldsymbol h,Y]
\end{aligned}
$$
in the group $\Sym^\#_\Pi(U)$.

If $\Sym^\#_\Pi[U;\boldsymbol h,Y]^{-1}=\Sym^\#_\Pi[U;\boldsymbol h,Y]$, then 
$$\{\boldsymbol h_{U,\Delta}L^{-1}_a\boldsymbol h_{\Delta,U}:a\in\Delta\}=\Sym^\#_\Pi[U;\boldsymbol h,Y]^{-1}=\Sym^\#_\Pi[U;\boldsymbol h,Y]=\{\boldsymbol h_{U,\Delta}L_a\boldsymbol h_{\Delta,U}:a\in\Delta\}$$ and hence for every $a\in\Delta$, there exists an element $b\in\Delta$ such that $\boldsymbol h_{U,\Delta}L_a^{-1}\boldsymbol h_{\Delta,U}=\boldsymbol h_{U,\Delta}L_b\boldsymbol h_{\Delta,U}$ and hence $L_a^{-1}=L_b$ and $x=L_bL_a(x)=b\puls (a\puls x)$ for every $x\in\Delta$, witnessing that the puls loop $(\Delta,\!\puls\!)$ is left-inversive.
\smallskip

5. If the puls loop $(\Delta,\!\puls\!)$ is associative, then 
for every elements $a,b,x\in\Delta$ we have $$L_aL_b(x)=a\puls(b\puls x)=(a\puls b)\puls x=L_{a\puls b}(x)$$ and hence
$$\boldsymbol h_{U,\Delta}L_{a\puls b}\boldsymbol h_{\Delta,U}=\boldsymbol h_{U,\Delta}L_aL_b\boldsymbol h_{\Delta,U}
=\boldsymbol h_{U,\Delta}L_a(\boldsymbol h_{\Delta,U}\boldsymbol h_{U,\Delta})L_b\boldsymbol h_{\Delta,U}=(\boldsymbol h_{U,\Delta}L_a\boldsymbol h_{\Delta,U})(\boldsymbol h_{U,\Delta}L_b\boldsymbol h_{\Delta,U}),
$$which means that the function 
$$\Phi:\Delta\to\Sym^\#_\Pi(U),\quad \Phi:a\mapsto \boldsymbol h_{U,\Delta}L_a\boldsymbol h_{\Delta,U},$$
is a homomorphism from the loop $(\Delta,\!\puls\!)$ to the group $\Sym^\#_\Pi(U)$. Then its image $\Phi[\Delta]=\Sym^\#_\Pi[U;\boldsymbol h,Y]$ is a subgroup of $\Sym^\#_\Pi(U)$ and $\Phi$ is an isomorphism of the puls loop $(\Delta,\!\puls\!)$ and the group $\Sym^\#_\Pi[U;\boldsymbol h,Y]$.

If $\Sym^\#_\Pi[U;\boldsymbol h,Y]=\{\boldsymbol h_{U,\Delta}L_a\boldsymbol h_{\Delta,U}:a\in\Delta\}$ is a subgroup of the group $\Sym^\#_\Pi(U)$, then for every $a,b\in \Delta$ there exist an element $c\in \Delta$ such that $$(\boldsymbol h_{U,\Delta}L_a\boldsymbol h_{\Delta,U})(\boldsymbol h_{U,\Delta}L_b\boldsymbol h_{\Delta,U})=\boldsymbol h_{U,\Delta}L_c\boldsymbol h_{\Delta,U}$$ and hence $L_aL_b=L_c$ and $c\puls x=L_c(x)=L_aL_b(x)=a\puls(b\puls x).$ In particular, $c=c\puls o=a\puls (b\puls o)=a\puls b$. Then $a\puls (b\puls x)=c\puls x=(a\puls b)\puls x$, witnessing that the loop $(\Delta,\!\puls\!)$ is associative.
\end{proof}

\section{Geometry of the dot operation}

In this section we find geometric counterparts of some algebraic properties of the multiplicative loop of a ternar. Let us recall that for any based affine plane $(\Pi,uow)$ and points $x,y,z\in\Delta$ on its diagonal, the equality $x\cdot y=z$ holds if and only if $xz\in\Aline o{ey}$.

\begin{picture}(60,115)(-150,-10)

\linethickness{0.75pt}
\put(0,0){\color{red}\line(1,1){90}}
\put(0,0){\color{orange}\line(2,3){60}}
\put(30,30){\color{cyan}\line(0,1){15}}
\put(45,45){\color{teal}\line(-1,0){15}}
\put(60,60){\color{cyan}\line(0,1){30}}
\put(90,90){\color{teal}\line(-1,0){30}}

\put(0,0){\color{red}\circle*{3}}
\put(0,0){\color{white}\circle*{2}}

\put(-2,-8){$o$}
\put(30,30){\color{red}\circle*{3}}
\put(30,23){$e$}
\put(45,45){\circle*{3}}
\put(46,38){$y$}
\put(60,60){\circle*{3}}
\put(61,53){$x$}
\put(90,90){\circle*{3}}
\put(93,88){$z=x\cdot y$}
\put(30,45){\circle*{3}}
\put(18,45){$ey$}
\put(60,90){\circle*{3}}
\put(55,93){$xz$}

\end{picture}

\begin{proposition}\label{p:dot-dot-inverses<=>} For a based affine plane $(\Pi,uow)$ and its ternar $\Delta$,\newline the following conditions are equivalent:
\begin{enumerate}
\item the loop $(\Delta^*,\cdot)$ is invertible;
\item $\forall x,y\in\Delta^*\;\;(o\in\Aline{ex}{ye}\;\Ra\;o\in \Aline{xe}{ey})$;
\item $\forall x,y\in\Delta^*\;\;(o\in\Aline{xe}{ey}\;\Ra\;o\in \Aline{ex}{ye})$.
\begin{picture}(50,0)(-160,0)
\linethickness{0.5pt}
\put(0,0){\color{red}\line(1,1){40}}
\put(0,0){\color{orange}\line(2,1){40}}
\put(0,0){\color{magenta}\line(1,2){20}}
\put(10,10){\color{teal}\line(1,0){10}}
\put(10,10){\color{cyan}\line(0,1){10}}
\put(20,20){\color{teal}\line(-1,0){10}}
\put(20,20){\color{cyan}\line(0,-1){10}}
\put(20,20){\color{teal}\line(1,0){20}}
\put(20,20){\color{cyan}\line(0,1){20}}
\put(40,40){\color{teal}\line(-1,0){20}}
\put(40,40){\color{cyan}\line(0,-1){20}}

\put(0,0){\color{red}\circle*{2}}
\put(10,10){\circle*{2}}
\put(20,20){\circle*{2}}
\put(40,40){\circle*{2}}
\put(10,20){\circle*{2}}
\put(20,10){\circle*{2}}
\put(40,20){\circle*{2}}
\put(20,40){\circle*{2}}
\end{picture}
\end{enumerate}
\end{proposition}

\begin{proof}  $(1)\Ra(2)$ Assume that the loop $(\Delta^*,\cdot)$ is invertible. Given any points $x,y\in\Delta^*$ with $o\in\Aline{ex}{ye}$, we should prove that $o\in \Aline{xe}{ey}$. The inclusion $o\in \Aline{ex}{ye}$ implies $ex\ne ye$ and $ye\in \Aline o{ex}$. The latter inclusion implies the equality $y\cdot x=e$, by the geometric definition of multiplication. Since the loop $(\Delta^*,\cdot)$ is invertible, the equality $y\cdot x=e$ implies $x\cdot y=e$ and $xe\in \Aline{o}{ey}$, by the definition of multiplication in the 0-loop $(\Delta,\cdot)$. It follows from $ex\ne ye$ that $xe\ne ey$. In this case $xe\in \Aline o{ey}$ implies $o\in\Aline {xe}{ey}$.

\begin{picture}(200,110)(-140,-15)
\linethickness{0.75pt}
\put(0,0){\color{red}\line(1,1){80}}
\put(0,0){\color{orange}\line(2,1){80}}
\put(0,0){\color{magenta}\line(1,2){40}}
\put(20,20){\color{teal}\line(1,0){20}}
\put(20,20){\color{cyan}\line(0,1){20}}
\put(40,40){\color{teal}\line(-1,0){20}}
\put(40,40){\color{cyan}\line(0,-1){20}}
\put(40,40){\color{teal}\line(1,0){40}}
\put(40,40){\color{cyan}\line(0,1){40}}
\put(80,80){\color{teal}\line(-1,0){40}}
\put(80,80){\color{cyan}\line(0,-1){40}}

\put(0,0){\color{red}\circle*{3}}
\put(0,0){\color{white}\circle*{2}}
\put(-2,-7){$o$}
\put(20,20){\circle*{3}}
\put(12,18){$x$}
\put(40,40){\color{red}\circle*{3}}
\put(41.5,33){$e$}
\put(80,80){\circle*{3}}
\put(83,80){$y$}
\put(20,40){\circle*{3}}
\put(8,42){$xe$}
\put(40,20){\circle*{3}}
\put(40,12){$ex$}
\put(80,40){\circle*{3}}
\put(83,38){$ye$}
\put(40,80){\circle*{3}}
\put(37,83){$ey$}

\end{picture}

$(2)\Ra(3)$ Assume that the condition (2) holds. To prove the condition (3), take any points $x,y\in\Delta^*$ with $o\in \Aline {xe}{ey}$. Then $o\in\Aline {xe}{ey}=\Aline{ey}{xe}$. Applying the condition (2), we obtain $o\in \Aline {ye}{ex}=\Aline{ex}{ye}$, witnessing that (3) holds.
\smallskip

$(3)\Ra(1)$ Assume that the condition (3) holds. To prove that the loop $(\Delta^*,\cdot)$ is invertible, take any points $x\in \Delta^*$. Since $(\Delta^*,\cdot)$ is a loop, there exists a unique element $y\in \Delta^*$ such that $x\cdot y=e$. If $x=e$, then $y=e$ is a two-sided inverse to $e$. So, assume that $x\ne e$. By the geometric definition of multiplication, the equlity $x\cdot y=e$ implies $xe\in\Aline o{ey}$. Since $xe\ne ey$, the inclusion $xe\in\Aline o{ey}$ is equivalent to $o\in\Aline{ey}{xe}$.  The condition (3) implies $o\in\Aline{ye}{ex}$, which implies $ye\in\Aline o{ex}$ and $y\cdot x=e$, by the geometric definition of multiplication. Therefore, $x\cdot y=e=y\cdot x$, witnessing that the loop $(\Delta^*,\cdot)$ is invertible.
\end{proof}

\begin{proposition}\label{p:right-inversive-dot<=>} For a based affine plane $(\Pi,uow)$ and its ternar $\Delta$,\newline the following conditions are equivalent:
\begin{enumerate}
\item the loop $(\Delta^*,\cdot)$ has the right-inverse property;
\item $\forall x,y,z\in\Delta^*\;\;(o\in\Aline{ex}{yz}\;\Ra\;o\in \Aline{xe}{zy})$;
\item $\forall a,b,\alpha,\beta\in\Delta^*\;\;(o\in\Aline{ab}{\alpha\beta}\;\Ra\;o\in \Aline{ba}{\beta\alpha})$.
\begin{picture}(200,0)(-140,5)
\linethickness{0.6pt}
\put(0,0){\color{red}\line(1,1){50}}
\put(0,0){\color{orange}\line(2,1){50}}
\put(0,0){\color{magenta}\line(1,2){25}}
\put(10,10){\color{teal}\line(1,0){10}}
\put(10,10){\color{cyan}\line(0,1){10}}
\put(20,20){\color{teal}\line(-1,0){10}}
\put(20,20){\color{cyan}\line(0,-1){10}}
\put(25,25){\color{teal}\line(1,0){25}}
\put(25,25){\color{cyan}\line(0,1){25}}
\put(50,50){\color{teal}\line(-1,0){25}}
\put(50,50){\color{cyan}\line(0,-1){25}}

\put(0,0){\color{red}\circle*{2.5}}
\put(0,0){\color{white}\circle*{1.5}}
\put(10,10){\circle*{2}}
\put(20,20){\circle*{2}}
\put(25,25){\circle*{2}}
\put(50,50){\circle*{2}}
\put(10,20){\circle*{2}}
\put(20,10){\circle*{2}}
\put(50,25){\circle*{2}}
\put(25,50){\circle*{2}}
\end{picture}
\end{enumerate}
\end{proposition}
\smallskip

\begin{proof}  $(1)\Ra(2)$ Assume that the loop $(\Delta^*,\cdot)$ has the right-inverse property. Given any points $x,y,z\in\Delta^*$ with $o\in\Aline{ex}{yz}$, we should prove that $o\in \Aline{xe}{zy}$. The inclusion $o\in \Aline{ex}{yz}$ implies $e\ne y$ and $yz\in\Aline o{ex}$. The latter inclusion implies $y\cdot x=z$, by the geometric definition of multiplication. Since $(\Delta^*,\cdot)$ has the right-inverse property, there exists an element $x'\in\Delta^*$ such that $a=(a\cdot x)\cdot x'$ for every $a\in \Delta^*$. In particular, $e=(e\cdot x)\cdot x'=x\cdot x'$ and $y=(y\cdot x)\cdot x'=z\cdot x'$. 
The equalities $e=x\cdot x'$ and $y=z\cdot x'$ imply $xe,zy\in \Aline o{ex'}$ and hence $o\in \Aline {xe}{zy}$ (because $xe\ne zy$).
\smallskip

$(2)\Ra(1)$ Assume that the condition (2) holds. To prove that the loop $(\Delta^*,\cdot)$ has the right-inverse property, take any point $x\in \Delta^*$. We have to find a point $x'\in \Delta^*$ such that $y\cdot (x\cdot x')=y=(y\cdot x)\cdot x'$. 

Since $(\Delta^*,\cdot)$ is a quasigroup, there exists a unique element $x'\in\Delta^*$ such that $x\cdot x'=e$ and hence $xe\in \Aline o{ex'}$, by the definition of the multiplication.  It remains to prove that $(y\cdot x)\cdot x'=y$ for every $y\in\Delta^*$. If $x=e$, then $x'=e$ and $(y\cdot x)\cdot x'=(y\cdot e)\cdot e=y$.  So, assume that $x\ne e$. In this case $xe\in\Aline o{ex'}$ implies $o\in \Aline {ex'}{xe}$. Applying the condition (2), we conclude that $o\in\Aline {x'e}{ex}$ and hence $\Aline o{x'e}=\Aline o{ex}$ and $x'\cdot x=e$. If $y=x'$, then $(y\cdot x)\cdot x'=(x'\cdot x)\cdot y=e\cdot y=y$ and we are done. So, assume that $y\ne x'$. In this case, consider the element $z\defeq y\cdot x$ and observe that 
$yz\in\Aline o{ex}=\Aline o{x'e}$, by the definition of multiplication. Since $y\ne x'$, the inclusion $yz\in\Aline o{x'e}$ implies $o\in\Aline {x'e}{yz}$. Applying the condition (2), we obtain $o\in \Aline {ex'}{zy}$ and hence $zy\in\Aline o{ex'}$ and $z\cdot x'=y$, by the definition of multiplication. Then $(y\cdot x)\cdot x'=z\cdot x'=y$, witnessing that the loop $(\Delta,+)$ has the right-inverse property.

\begin{picture}(200,150)(-150,-15)
\linethickness{0.75pt}
\put(0,0){\color{red}\line(1,1){120}}
\put(0,0){\color{orange}\line(2,1){120}}
\put(0,0){\color{magenta}\line(1,2){60}}
\put(20,20){\color{teal}\line(1,0){20}}
\put(20,20){\color{cyan}\line(0,1){20}}
\put(40,40){\color{teal}\line(-1,0){20}}
\put(40,40){\color{cyan}\line(0,-1){20}}
\put(60,60){\color{teal}\line(1,0){60}}
\put(60,60){\color{cyan}\line(0,1){60}}
\put(120,120){\color{teal}\line(-1,0){60}}
\put(120,120){\color{cyan}\line(0,-1){60}}
\put(20,20){\color{cyan}\line(0,-1){10}}
\put(20,10){\color{teal}\line(-1,0){10}}

\put(20,10){\circle*{3}}
\put(16,2){$ex'$}
\put(10,10){\circle*{3}}

\put(0,0){\color{red}\circle*{3}}
\put(0,0){\color{white}\circle*{2}}
\put(-2,-7){$o$}
\put(20,20){\circle*{3}}
\put(13,20){$e$}
\put(40,40){\circle*{3}}
\put(42,36){$x$}
\put(60,60){\circle*{3}}
\put(52,60){$y$}
\put(120,120){\circle*{3}}
\put(123,118){$z$}
\put(20,40){\circle*{3}}
\put(8,40){$ex$}
\put(40,20){\circle*{3}}
\put(40,12){$xe$}
\put(120,60){\circle*{3}}
\put(118,52){$zy$}
\put(60,120){\circle*{3}}
\put(54,123){$yz$}

\end{picture}

$(2)\Ra(3)$ Assume that the condition (2) holds. To prove the condition (3), take any points $a,b,\alpha,\beta$ with $o\in\Aline {ab}{\alpha\beta}$. We have to prove that $o\in\Aline{ba}{\beta\alpha}$. The inclusion $o\in\Aline {ab}{\alpha\beta}$ implies $a\ne\alpha$, $b\ne\beta$, and $\Aline o{ab}=\Aline o{\alpha\beta}$. If $e\in\{a,\alpha\}$, then $o\in\Aline {ab}{\alpha\beta}$ implies $o\in\Aline{ba}{\beta\alpha}$, by the condition (2). So, assume that $a\ne e\ne\alpha$. Find a unique point $x\in\Delta$ such that $ex\in\Aline o{ab}=\Aline o{\alpha\beta}$ and hence $o\in\Aline{ex}{ab}\cap\Aline{ex}{\alpha\beta}$ (because $a\ne e\ne\alpha$). The condition (2) ensures that $o\in\Aline{xe}{ba}\cap\Aline{xe}{\beta\alpha}$ and hence $o\in\Aline{ba}{\beta\alpha}$ (because $b\ne\beta$).
\smallskip

The implication $(3)\Ra(2)$ is trivial.
\end{proof}

\begin{proposition}\label{p:left-inversive-dot<=>} For a based affine plane $(\Pi,uow)$ and its ternar $\Delta$, the following conditions are equivalent:
\begin{enumerate}
\item the loop $(\Delta^*,\cdot)$ has the left-inverse property;
\item $\forall x,x',y,z\in\Delta^*\;\;(o\in\Aline{ex}{x'e}\cap \Aline {xz}{ey}\;\Ra\;o\in \Aline{x'y}{ez})$.
\end{enumerate}
\end{proposition}

\begin{proof} $(1)\Ra(2)$ Assume that the loop $(\Delta^*,\cdot)$ has the left-inverse property. Given any points $x,x',y,z\in\Delta^*$ with $o\in\Aline{ex}{x'e}\cap \Aline {xz}{ey}$, we have to prove that $o\in \Aline{x'y}{ez}$. 

\begin{picture}(100,115)(-150,-12)

\linethickness{0.7pt}

\put(15,15){\color{cyan}\line(0,1){30}}
\put(30,30){\color{cyan}\line(0,1){60}}
\put(60,60){\color{cyan}\line(0,1){30}}
\put(30,30){\color{teal}\line(-1,0){15}}
\put(45,45){\color{teal}\line(-1,0){30}}
\put(60,60){\color{teal}\line(-1,0){30}}
\put(60,60){\color{cyan}\line(0,1){30}}
\put(90,90){\color{teal}\line(-1,0){60}}
{\linethickness{1pt}
\put(0,0){\color{red}\line(1,1){90}}
}
\put(0,0){\color{orange}\line(2,3){60}}
\put(0,0){\color{magenta}\line(1,2){30}}
\put(0,0){\color{purple}\line(1,3){30}}

\put(0,0){\color{red}\circle*{3}}
\put(0,0){\color{white}\circle*{2}}
\put(1,-7){$o$}
\put(15,15){\circle*{3}}
\put(16,8){$x'$}
\put(30,30){\circle*{3}}
\put(31,23){$e$}
\put(45,45){\circle*{3}}
\put(46,38){$y$}
\put(60,60){\circle*{3}}
\put(61,53){$x$}
\put(90,90){\circle*{3}}
\put(91,83){$z$}
\put(15,30){\circle*{3}}
\put(15,45){\circle*{3}}
\put(30,45){\circle*{3}}
\put(30,60){\circle*{3}}
\put(30,90){\circle*{3}}
\put(60,90){\circle*{3}}

\end{picture}

By the definition of the multiplication, the inclusion $o\in\Aline {ex}{x'e}$ implies $x'e\in\Aline o{ex}$ and $x'\cdot x=e$. The left-inverse property ensures that $x'\cdot (x\cdot y)=y$. The inclusion $o\in \Aline{xz}{ey}$ implies $x\cdot y=z$. Then $y=x'\cdot (x\cdot y)=x'\cdot z$ and hence $x'y\in \Aline o{ez}$. It follows from $o\in \Aline {ex}{x'e}$ that $e\ne x'$ and hence $x'y\in\Aline o{ez}$ implies $o\in\Aline{x'y}{ez}$.

$(2)\Ra(1)$ Assume that the condition (2) holds. To show that the loop $(\Delta^*,\cdot)$ has the left-inverse property, take any element $x\in \Delta^*$. Since $(\Delta^*,\cdot)$ is a loop, there exists a unique element $x'$ such that $x'\cdot x=e$. To show that $x'\cdot (x\cdot y)=y$, consider the point $z\defeq x\cdot y$ and observe tat $xz\in\Aline o{ey}$. If $x=e$, then $x'=e$ and $x'\cdot (x\cdot y)=e\cdot(e\cdot y)=y$. So, assume that $x\ne e$. In this case the inclusions $x'e\in\Aline e{ex}$ and  $xz\in\Aline o{ey}$ imply $o\in\Aline{x'e}{ex}\cap\Aline{xz}{ey}$. The condition (2) ensures that $o\in\Aline{x'y}{ez}$ and hence $x'y\in\Aline o{ez}$ and $x'\cdot z=y$, by the definition of the multiplication. Then $y=x'\cdot z=x'\cdot(x\cdot y)$, witnessing that the loop $(\Delta^*,\cdot)$ has the left-inverse property.
\end{proof}

\begin{proposition}\label{p:ass-dot<=>} For a based affine plane $(\Pi,uow)$ and its ternar $\Delta$,\newline the following conditions are equivalent:
\begin{enumerate}
\item the loop $(\Delta^*,\cdot)$ is associative;
\item $\forall x,y,z,p,q,s\in\Delta^*\;\;(o\in\Aline{xp}{ey}\cap\Aline {yq}{ez}\cap \Aline {ps}{ez}\;\Ra\;o\in \Aline{xs}{eq})$;
\item $\forall a,b,c,d,\alpha,\beta,\gamma,\delta\in\Delta^*\;\;(o\in\Aline{ac}{\alpha\gamma}\cap\Aline{ad}{\alpha\delta}\cap\Aline{bc}{\beta\gamma}\;\Ra\;o\in \Aline{bd}{\beta\delta})$.
\begin{picture}(50,0)(-35,0)
\linethickness{0.64pt}

\put(16,8){\color{teal}\line(1,0){8}}
\put(24,24){\color{teal}\line(-1,0){8}}
\put(32,16){\color{teal}\line(1,0){16}}
\put(48,48){\color{teal}\line(-1,0){16}}

\put(0,0){\color{orange}\line(3,1){48}}
\put(0,0){\color{red}\line(2,1){32}}
\put(0,0){\color{magenta}\line(1,1){48}}
\put(0,0){\color{purple}\line(2,3){32}}
\put(24,8){\color{cyan}\line(0,1){16}}
\put(16,8){\color{cyan}\line(0,1){16}}
\put(32,16){\color{cyan}\line(0,1){32}}
\put(48,48){\color{cyan}\line(0,-1){32}}

\put(0,0){\color{red}\circle*{2}}
\put(16,8){\circle*{2}}
\put(24,8){\circle*{2}}
\put(16,24){\circle*{2}}
\put(24,24){\circle*{2}}
\put(32,16){\circle*{2}}
\put(48,16){\circle*{2}}
\put(32,48){\circle*{2}}
\put(48,48){\circle*{2}}

\end{picture}
\end{enumerate}
\end{proposition}

\begin{proof} $(1)\Ra(2)$ Assume that the multiplicative loop $(\Delta^*,\cdot)$ is associative. To prove the condition (2), take any points $x,y,z,p,q,s\in\Delta^*$ with $o\in\Aline{xp}{ey}\cap\Aline {yq}{ez}\cap \Aline {ps}{oz}$. We have to prove that $o\in \Aline{xs}{eq}$. The inclusions  $o\in\Aline{xp}{ey}\cap\Aline {yq}{ez}\cap \Aline {ps}{ez}$ imply $e\notin\{x,y,p\}$ and $x\cdot y=p$, $y\cdot z=q$ and $p\cdot z=s$. The associativity of the loop $(\Delta^*,\cdot)$ ensures that $s=p\cdot z=(x\cdot y)\cdot z=x\cdot(y\cdot z)=x\cdot q$ and hence $xs\in \Aline o{eq}$. Since $x\ne e$, the inclusion $xs\in\Aline o{eq}$ implies $o\in\Aline{xs}{eq}$, proving the condition (2).

\begin{picture}(150,150)(-130,-15)
\linethickness{0.64pt}

\put(0,0){\color{orange}\line(3,2){96}}
{\linethickness{1pt}
\put(0,0){\color{red}\line(1,1){128}}
}
\put(0,0){\color{magenta}\line(3,4){96}}
\put(0,0){\color{purple}\line(1,2){64}}
\put(16,16){\color{cyan}\line(0,1){16}}
\put(16,16){\color{teal}\line(1,0){8}}
\put(24,24){\color{cyan}\line(0,-1){8}}
\put(24,24){\color{cyan}\line(0,1){24}}
\put(32,32){\color{teal}\line(-1,0){16}}
\put(48,48){\color{teal}\line(-1,0){24}}
\put(64,64){\color{teal}\line(1,0){32}}
\put(64,64){\color{cyan}\line(0,1){64}}
\put(96,64){\color{cyan}\line(0,1){64}}
\put(128,128){\color{teal}\line(-1,0){64}}

\put(0,0){\color{red}\circle*{2.5}}
\put(0,0){\color{white}\circle*{1.5}}
\put(0,-7){$o$}
\put(16,16){\circle*{2.5}}
\put(13,9){$y$}
\put(24,24){\circle*{2.5}}
\put(25,20){$e$}
\put(32,32){\circle*{2.5}}
\put(34,28){$q$}
\put(48,48){\circle*{2.5}}
\put(49,42){$z$}
\put(64,64){\circle*{2.5}}
\put(62,56){$p$}
\put(96,96){\circle*{2.5}}
\put(98,90){$x$}
\put(128,128){\circle*{2.5}}
\put(129,122){$s$}
\put(24,16){\circle*{2.5}}
\put(22,8){$ey$}
\put(16,32){\circle*{2.5}}
\put(4,30){$yq$}
\put(24,32){\circle*{2.5}}
\put(12,48){$ez$}
\put(24,48){\circle*{2.5}}
\put(48,48){\circle*{2.5}}
\put(64,128){\circle*{2.5}}
\put(60,131){$ps$}
\put(96,64){\circle*{2.5}}
\put(98,60){$xp$}
\put(92,131){$xs$}
\put(96,128){\circle*{2.5}}
\end{picture}

$(2)\Ra(1)$ Assume that the condition (2) holds. 

\begin{claim}\label{cl:Delta*-ass} The equality $(x\cdot y)\cdot z=x\cdot (y\cdot z)$ holds for every $x,y,z\in\Delta^*$ with $x\cdot y\ne e$.
\end{claim}

\begin{proof} Consider the points $p\defeq x\cdot y\ne e$, $q\defeq y\cdot z$ and $s\defeq (x\cdot y)\cdot z=p\cdot z$. The definition of multiplication in the loop $(\Delta^*,\cdot)$ ensures that $xp\in\Aline o{ey}$, $yq\in\Aline o{ez}$ and $ps\in\Aline o{ez}$.

If $x=e$, then $(x\cdot y)\cdot z=(e\cdot y)\cdot z=y\cdot z=e\cdot(y\cdot z)=x\cdot (y\cdot z)$. 

If $y=e$, then $(x\cdot y)\cdot z=(x\cdot e)\cdot z=x\cdot z=x\cdot(e\cdot z)=x\cdot (y\cdot z)$.  

So, we assume that $x\ne e\ne y$. Since $p\ne e$, the inclusions  $xp\in\Aline o{ey}$, $yq\in\Aline o{ez}$ and $ps\in\Aline o{ez}$ imply $o\in\Aline{xp}{ey}\cap\Aline {yq}{ez}\cap\Aline {ps}{ez}$. Applying the condition (2), we conclude that $o\in\Aline {xs}{eq}$ and hence $xs\in\Aline o{eq}$ and $s=x\cdot q$, by the definition of multiplication. Then $(x\cdot y)\cdot z=p\cdot z=s=x\cdot q=x\cdot (y\cdot z).$
\end{proof}

Assuming that the loop $(\Delta^*,\cdot)$ is not associative, we can find points $x,y,z\in \Delta^*$ with $(x\cdot y)\cdot z\ne x\cdot(y\cdot z)$. In this case  Claim~\ref{cl:Delta*-ass} ensures that $x\cdot y=e$. Since $(\Delta^*,\cdot)$ is a loop, there exists an element $x'\in \Delta^*$ such that $(x\cdot y)\cdot z=x'\cdot(y\cdot z)$. It follows from $x'\cdot (y\cdot z)=(x\cdot y)\cdot z\ne x\cdot(y\cdot z)$ that $x'\ne x$ and hence $x'\cdot y\ne x\cdot y=e$. So, we can apply Claim~\ref{cl:Delta*-ass} and conclude that $(x'\cdot y)\cdot z=x'\cdot (y\cdot z)=(x\cdot y)\cdot z$. Since $(\Delta^*,\cdot)$ is a quasigroup, the equality $(x'\cdot y)\cdot z=(x\cdot y)\cdot z$ implies $x'\cdot y=x\cdot y$ and $x'=x$, which is a contradiction showing that the loop $(\Delta^*,\cdot)$ is associative.   
\smallskip

$(2)\Ra(3)$ Assume that the condition (2) holds. 

\begin{claim}\label{cl:Delta*-ass2} $\forall a,b,c,\alpha,\beta,\gamma\in\Delta^*\;\;(o\in\Aline{ab}{\alpha\beta}\cap\Aline{bc}{\beta\gamma}\;\Ra\;o\in \Aline{ac}{\alpha\gamma}$.
\end{claim} 

\begin{proof} Given any points  $\forall a,b,c,\alpha,\beta,\gamma\in\Delta^*\;\;o\in\Aline{ab}{\alpha\beta}\cap\Aline{bc}{\beta\gamma}$, we should prove that $o\in \Aline{ac}{\alpha\gamma}$. The inclusion $o\in \Aline{ab}{\alpha\beta}\cap\Aline{bc}{\beta\gamma}$ implies $a\ne\alpha$, $\Aline o{ab}=\Aline o{\alpha\beta}$, and $\Aline o{bc}=\Aline o{\beta\gamma}$. Find unique points $y,q,z\in\Delta^*$ such that $ye\in \Aline o{ab}=\Aline o{\alpha\beta}$ and $yq,ez\in\Aline o{bc}=\Aline o{\beta\gamma}$. The condition (2) implies $o\in\Aline {bc}{eq}\cap\Aline{\alpha\gamma}{eq}$ and hence $o\in \Aline {ac}{\alpha\gamma}$ (because $b\ne\beta$). 
\end{proof}

\begin{claim}\label{cl:Delta*-ass3} $\forall a,b,c,\alpha,\beta,\gamma\in\Delta^*\;\;(o\in\Aline{ab}{\alpha\beta}\cap\Aline{ac}{\alpha\gamma}\;\Ra\;o\in \Aline{bc}{\beta\gamma}$.
\end{claim}

\begin{proof} Given any points $a,b,c,\alpha,\beta,\gamma\in\Delta^*\;\;(o\in\Aline{ab}{\alpha\beta}\cap\Aline{ac}{\alpha\gamma}$, we should prove that $o\in \Aline{bc}{\beta\gamma}$. The inclusion $o\in\Aline {ab}{\alpha\beta}$ implies $a\ne\alpha$ and $b\ne \beta$. Then there exists a unique point $\delta\in\Delta^*$ such that $o\in \Aline{bc}{\beta\delta}$. Applying Claim~\ref{cl:Delta*-ass2} to the points $a,b,c,\alpha,\beta,\delta$, we conclude that $o\in\Aline{ac}{\alpha\delta}\cap\Aline{ac}{\alpha\gamma}$ and hence $\{\alpha\gamma,\alpha\delta\}\subseteq \Aline o{ac}$ and $\gamma=\delta$. Then $o\in \Aline {bc}{\beta\delta}=\Aline{bc}{\beta\gamma}$.
\end{proof}

Now we can prove the condition (3). Given any points $a,b,c,d,\alpha,\beta,\gamma,\delta\in\Delta^*$ with $o\in\Aline {ac}{\alpha\gamma}\cap\Aline{ad}{\alpha\delta}\cap\Aline {bc}{\beta\gamma}$, we should prove that $o\in\Aline {bd}{\beta\delta}$. By Claim~\ref{cl:Delta*-ass3}, $o\in\Aline {ac}{\alpha\gamma}\cap\Aline{ad}{\alpha\delta}$ implies $o\in \Aline{cd}{\gamma\delta}$. By Claim~\ref{cl:Delta*-ass2}, $o\in \Aline {bc}{\beta\gamma}\cap  \Aline{cd}{\gamma\delta}$ implies  $o\in \Aline{bd}{\beta\delta}$.
\smallskip

$(3)\Ra(2)$ Assume that the condition (3) holds. To prove the condition (2), take any points  $x,y,z,p,q,s\in\Delta^*$ with $o\in\Aline{xp}{ey}\cap\Aline {yq}{ez}\cap \Aline {ps}{ez}$. We have to prove that $o\in \Aline{xs}{eq}$. The inclusion $o\in \Aline{xp}{ey}$ implies $p\ne y$ and hence $yq\ne ps$. Then $o\in\Aline {yq}{ez}\cap \Aline {ps}{ez}$ implies $o\in \Aline {yq}{ps}$.

For the sequence of points $$(a,b,c,d,\alpha,\beta,\gamma,\delta)\defeq (p,x,p,s,y,e,y,q),$$ we have $o\in \Aline{pp}{yy}\cap\Aline {ps}{yq}\cap \Aline{xp}{ey}=\Aline {ac}{\alpha\gamma}\cap\Aline{ad}{\alpha\delta}\cap\Aline {bc}{\beta\gamma}$. Applying the condition (3), we obtain $o\in\Aline {bd}{\beta\delta}=\Aline{xs}{eq}$.
\end{proof}

\begin{proposition}\label{p:com-dot<=>} For a based affine plane $(\Pi,uow)$ and its ternar $\Delta$,\newline the following conditions are equivalent:
\begin{enumerate}
\item the loop $(\Delta^*,\cdot)$ is commutative;
\item $\forall x,y,z\in\Delta^*\setminus\{e\}\;\;(o\in\Aline {xz}{ey}\;\Ra\;o\in \Aline {yz}{ex})$.

\end{enumerate}
\end{proposition}

\begin{proof} $(1)\Ra(2)$ Assume that the loop $(\Delta,+)$ is commutative. Given any points $x,y,z\in\Delta^*\setminus\{e\}$ with $o\in\Aline {xz}{ey}$, we have to check that $o\in \Aline {yz}{ex}$. Since $x\ne e$, the inclusion $o\in\Aline{xz}{ey}$ implies $xz\in\Aline o{ey}$ and hence $z=x\cdot y$, by the definition of multiplication. The commutativity of the loop $(\Delta^*,\cdot)$ implies $y\cdot x=x\cdot y=z$ and hence $yz\in\Aline o{ex}$. Since $y\ne e$, the latter inclusion implies $o\in \Aline{yz}{ex}$.

%

$(2)\Ra(1)$ Assume that the condition (2) holds. Assuming that the loop $\Delta^*$ is not commutative, we can find two points $x,y\in \Delta^*$ such that $x\cdot y\ne y\cdot x$. In this case $x\ne e\ne y$. Since $x\cdot y\ne y\cdot x$, one of the points $x\cdot y$ or $y\cdot x$ is not equal to $e$. We lose no generality assuming that $x\cdot y\ne  e$. The definition of the multiplication ensures that for the product $z\defeq x\cdot y$ we have  $xz\in\Aline o{ey}$ and hence $o\in\Aline {xz}{ey}$ (because $x\ne e$). Applying the condition (2), we obtain $o\in\Aline {yz}{ex}$ and hence $yz\in \Aline o{ex}$ and $y\cdot x=z=x\cdot y$, which contradicts the choice of the points $x,y$. This contradiction shows that the loop $(\Delta^*,\cdot)$ is commutative.
\end{proof}

\chapter{Plus properties of affine spaces}\label{ch:PlusAff}

In this chapter we provide geometric characterizations of some properties of the plus operation in ternars of affine spaces.

Let us recall that an affine space $X$ is \defterm{coordinatized} by a ternar $R$ if $X$ contains a plane $\Pi$ that has an affine base $uow$ such that the ternar of the based affine plane $(\Pi,uow)$ is isomorphic to the ternar $R$. In this case we shall also say that $R$ is \index{affine space!ternar of}\defterm{a ternar} of the affine space $X$.

\section{Invertible-plus affine spaces}

\begin{definition}\label{d:a++inversive} An affine liner $X$ is defined to be \index{invertible-plus affine space}\index{affine space!invertible-plus}\defterm{invertible-plus} if for every points $a,b,c,\alpha,\beta,\gamma\in X$ and $o\in\Aline a\alpha\cap\Aline b\beta\cap\Aline c\gamma$ 
$$\big((\Aline ac\parallel \Aline ob\parallel \Aline o\beta\parallel\Aline\alpha\gamma\nparallel \Aline b\gamma\parallel\Aline o\alpha\parallel \Aline oa\parallel\Aline c\beta)\;\wedge\;(\Aline ab\parallel\Aline co\parallel \Aline o\gamma)\big)\;\Ra\;\Aline co\parallel\Aline \alpha\beta.$$
\end{definition}

\begin{picture}(150,80)(-180,-45)
\linethickness{0.75pt}
\put(0,0){\color{teal}\line(1,0){30}}
\put(0,0){\color{cyan}\line(0,1){30}}
\put(0,0){\color{teal}\line(-1,0){30}}
\put(0,0){\color{cyan}\line(0,-1){30}}

\put(0,30){\color{teal}\line(1,0){30}}
\put(30,0){\color{cyan}\line(0,1){30}}
\put(-30,-30){\color{teal}\line(1,0){30}}
\put(-30,-30){\color{cyan}\line(0,1){30}}
\put(-30,0){\color{red}\line(1,1){30}}
\put(0,-30){\color{red}\line(1,1){30}}
{\linethickness{1pt}
\put(0,0){\color{red}\line(1,1){30}}
\put(0,0){\color{red}\line(-1,-1){30}}
}
\put(0,0){\circle*{3}}
\put(2,-7){$o$}
\put(30,0){\circle*{3}}
\put(33,-2){$\alpha$}
\put(-30,0){\circle*{3}}
\put(-38,-2){$a$}
\put(0,30){\circle*{3}}
\put(-7,31){$b$}
\put(0,-30){\circle*{3}}
\put(2,-36){$\beta$}
\put(30,30){\circle*{3}}
\put(32,31){$\gamma$}
\put(-30,-30){\circle*{3}}
\put(-37,-36){$c$}
\end{picture}

Let us recall that a loop $(X,\cdot)$ is called \defterm{invertible} if for every $x\in X$ there exists an element $y\in X$ such that $x\cdot y=e=y\cdot x$, where $e$ is the identity of the loop $X$.

\begin{theorem}\label{t:invertible-plus<=>} An affine space $X$ is invertible-plus if and only if for every ternar $R$ of $X$, the plus loop $(R,+)$ is invertible.
\end{theorem}

\begin{proof} Assume that an affine plane $X$ is invertible-plus. Given any plane $\Pi\subseteq X$ and an affine base $uow$ in $\Pi$, we need to check that the plus loop $(\Delta,+)$ of the ternar $\Delta$ of the based affine plane $(\Pi,uow)$ is invertible. By Proposition~\ref{p:two-sided}, it suffices to prove that for every points $x,y\in\Delta$ with $\Aline {ox}{yo}\subparallel \Delta$, we have $\Aline{xo}{oy}\subparallel\Delta$. If $o\in\{x,y\}$ then $\Aline {ox}{yo}\subparallel \Delta$ implies $x=o=y$ and hence $\Aline {xo}{yo}=\Aline{oo}{oo}\subparallel \Delta$. So, we assume that $o\notin\{x,y\}$. If $x=y$, then $\Aline {xo}{oy}=\Aline {yo}{ox}=\Aline {ox}{yo}\subparallel \Delta$. 

So, assume that $x\ne y$ and consider the points $a\defeq yo$, $b\defeq ox$, $c\defeq yy$, $\alpha\defeq xo$, $\beta\defeq oy$, $\gamma\defeq xx$ of the coordinate plane $\Delta^2$, identified with the plane $\Pi$. Observe that the points $a,b,c,\alpha,\beta,\gamma$
satisfy the parallelity relations $$\big((\Aline ac\parallel \Aline o\beta\parallel\Aline ob\parallel\Aline\alpha\gamma\nparallel \Aline b\gamma\parallel\Aline o\alpha\parallel \Aline oa\parallel\Aline \beta c)\;\wedge\;(\Aline ab=\Aline {yo}{ox}\subparallel \Delta=\Aline{yy}{xx}=\Aline co=\Aline o\gamma)\big).$$Since the affine plane $\Pi$ is invertible-plus, $\Aline {xo}{oy}=\Aline \alpha\beta\parallel\Aline co=\Delta$, witnessing that the plus loop $(\Delta,+)$ is invertible.

\begin{picture}(150,90)(-180,-45)
\linethickness{0.75pt}
\put(0,0){\color{teal}\line(1,0){30}}
\put(0,0){\color{cyan}\line(0,1){30}}
\put(0,0){\color{teal}\line(-1,0){30}}
\put(0,0){\color{cyan}\line(0,-1){30}}

\put(0,30){\color{teal}\line(1,0){30}}
\put(30,0){\color{cyan}\line(0,1){30}}
\put(-30,-30){\color{teal}\line(1,0){30}}
\put(-30,-30){\color{cyan}\line(0,1){30}}
\put(-30,0){\color{red}\line(1,1){30}}
\put(0,-30){\color{red}\line(1,1){30}}
{\linethickness{1pt}
\put(0,0){\color{red}\line(1,1){30}}
\put(0,0){\color{red}\line(-1,-1){30}}
}
\put(0,0){\circle*{3}}
\put(2,-7){$o$}
\put(30,0){\circle*{3}}
\put(33,-2){$\alpha$}
\put(-30,0){\circle*{3}}
\put(-38,-2){$a$}
\put(0,30){\circle*{3}}
\put(-7,31){$b$}
\put(0,-30){\circle*{3}}
\put(2,-36){$\beta$}
\put(30,30){\circle*{3}}
\put(32,31){$\gamma$}
\put(-30,-30){\circle*{3}}
\put(-37,-36){$c$}
\end{picture}

Next, assume that for every ternar $R$ of the affine space $X$, the plus loop $(R,+)$ is invertible. Given any points $a,b,c,\alpha,\beta,\gamma\in X$ and $o\in\Aline a\alpha\cap\Aline b\beta\cap\Aline c\gamma$ with 
$$(\Aline ac\parallel \Aline o\beta\parallel \Aline ob\parallel\Aline\alpha\gamma\nparallel \Aline b\gamma\parallel\Aline o\alpha\parallel\Aline oa\parallel\Aline c\beta)\;\wedge\;(\Aline ab\parallel\Aline co\parallel\Aline o\gamma),$$ we need to show that $\Aline \alpha\beta\parallel\Aline oc$. 
 
If $a=c$, then the parallelity relations $\Aline ac\parallel \Aline o\beta\parallel\Aline ob\parallel\Aline\alpha\gamma$ imply $b=o=\beta$ and $\alpha=\gamma$. It follows from $\{\alpha\}=\Aline\alpha\gamma\nparallel \Aline b\gamma\parallel\Aline o\alpha\parallel\Aline oa\parallel\Aline c\beta$ that $b\ne\gamma$, $a\ne o\ne\alpha$, $c\ne \beta$.  Then $\Aline \alpha\beta=\Aline \alpha o=\Aline ao=\Aline co$.
By analogy we can prove that $c=\beta$ implies $\Aline \alpha\beta=\Aline co$. 

So, assume that $a\ne c\ne \beta$. In this case, $\Aline ac\parallel \Aline o\beta\parallel \Aline ob\parallel\Aline\alpha\gamma\nparallel \Aline b\gamma\parallel \Aline o\alpha\parallel \Aline oa\parallel\Aline c\beta$ implies $b\ne\gamma\ne\alpha$, $a\ne o\ne \alpha$, and $b\ne o\ne\alpha$. 
It follows that $uow\defeq \alpha o b$ is an affine base for the affine plane $\Pi\defeq \overline{\{a,b,c,o,\alpha,\beta,\gamma\}}$, $\gamma$ is the diunit of the affine base $uow$ and $\Delta\defeq\Aline o\gamma$ is its diagonal. Identify the affine plane $\Pi$ with the coordinate plane $\Delta^2$ of the ternar $\Delta$. The parallelity relations $$(\Aline ac\parallel \Aline o \beta\parallel\Aline ob\parallel\Aline\alpha\gamma)\;\wedge\; (\Aline b\gamma\parallel\Aline o\alpha\parallel\Aline oa\parallel\Aline c\beta)\;\wedge\;(\Aline ab\parallel\Aline co\parallel\Aline o\gamma)$$ imply that $\Aline {o\gamma}{co}=\Aline ba\parallel \Delta$.  Since the plus loop $(\Delta,+)$ is invertible, we can apply Proposition~\ref{p:two-sided} and conclude that $\Aline {\alpha}{\beta}=\Aline {\gamma o}{oc}\subparallel\Delta=\Aline co$. Since $\alpha\ne\beta$, the subparallelity $\Aline{\alpha}{\beta}\subparallel \Aline c{\gamma}$ implies $\Aline\alpha\beta\parallel\Aline co$, witnessing that the affine space $X$ is  invertible-plus.
\end{proof}

Let us recall that a loop $(X,\cdot)$ is called \defterm{$\IZ$-associative} (or else \defterm{power-associative}) if for every element $x\in X$ there exists a sequence $(x^n)_{n\in\IZ}$ in $X$ such that $x^1=x$ and $x^{n+m}=x^n\cdot x^m$ for every integer numbers $n,m\in\IZ$. 

\begin{theorem}\label{t:invertible-plus<=>power-associative} An affine space $X$ is invertible-plus if and only if for every ternar $R$ of $X$, the plus loop $(R,+)$ is power-associative.
\end{theorem}

\begin{proof} The ``if'' part follows from the ``if'' part of Theorem~\ref{t:invertible-plus<=>} (because every power-associative loop is invertible). To prove the ``only if'' part, assume that an affine space $X$ is invertible-plus. Given any ternar $R$ of $X$, we need to check that the additive loop $(R,+)$ is power-associative. Find an affine base $uow$ whose ternar $\Delta$ is isomorphic to the ternar $R$. Given any element $x\in\Delta$, we need construct a sequence $(x^n)_{n\in\IZ}$ such that $x^1=x$ and $x^{n+m}=x^n+x^m$ for every integer numbers $n,m$. If $x=o$, then the constant sequence $(x^n)_{n\in\IZ}\defeq \{(n,o):n\in\IZ\}$ has the required property. So we assume that $x\ne o$. Put $x^0\defeq o$ and $x^1\defeq x$. By Theorem~\ref{t:invertible-plus<=>}, the loop $(\Delta,+)$ is invertible, so there exist a unique element $x^{-1}\in X$ such that $x^{-1}+x=o=x+x^{-1}$.  For every $n\in \IN$, let $x^{n+1}\defeq x^n+x$ and $x^{-(n+1)}\defeq x^{-n}+x^{-1}$. We claim that the sequence $(x^n)_{n\in\IZ}$ has the required property. 

Let $\boldsymbol h\defeq(\Aline ou)_\parallel\in\partial X$ and $\boldsymbol v\defeq(\Aline ov)_\parallel\in\partial X$ be the horizontal and vertical directions determined by the affine base $uow$. For every integer numbers $i,j\in\IZ$, find a unique point $x_i^j\in X$ such that $\Aline {x_i^j}{\boldsymbol h}=\Aline {x^i}{\boldsymbol h}$ and $\Aline {x_i^j}{\boldsymbol v}=\Aline {x^j}{\boldsymbol v}$.



 For every $n\in\omega$, the equality $x^{n+1}=x^n+x$ and the definition of addition operation in the ternar $\Delta$ ensures that $\Aline {x_0^1}{x_n^{n+1}}\subparallel \Delta$, which implies $\Aline{x_{n}^{n+1}}{x_{n+1}^{n+2}}\parallel \Delta$ for every $n\in\omega$. By the invertibility-plus of $X$, the parallelity relation  $\Aline{x_{n}^{n+1}}{x_{n+1}^{n+2}}\parallel \Delta$ implies $\Aline {x_n^{n-1}}{x_{n+1}^n}\parallel \Delta$.  By analogy, $x^{-(n+1)}=x^{-n}+x^{-1}$ implies $\Aline {x_0^{-1}}{x_{-n}^{-(n+1)}}\subparallel \Delta$, $\Aline {x_{-n}^{-(n+1)}}{x_{-(n+1)}^{-(n+2)}}\parallel \Delta$ and $\Aline {x_{-(n+1)}^{-n}}{x_{-(n+2)}^{-(n+1)}}\parallel \Delta$.
Therefore, $\Aline {x_n^{n+1}}{x_{n+1}^{n+2}}\parallel \Delta\parallel \Aline{x_n^{n-1}}{x_{n+1}^n}$ for all $n\in\IZ$. 

Applying the invertibility-plus, we can prove that for every $k\in \w$ and $n\in\IZ$, the lines $\Aline {x_n^{n+k}}{x_{n+1}^{n+1+k}}$ and $\Aline {x_n^{n-k}}{x_{n+1}^{n+1-1}}$ are parallel to the diagonal $\Delta$ of the affine base $uow$.

Then for every $n,m\in\IZ$ we have $\Aline {x_0^n}{x_m^{m+n}}\subparallel \Delta$, which implies $x^{n+m}=x^n+x^m$, witnessing that the plus loop $(\Delta,+)$ is power-associative and so is its isomorphic copy $(R,+)$.
 \end{proof}
 
\begin{exercise}Show that the Moulton plane is not invertible-plus.






{\em Hint:} Consider the points 
$$a=(-8,4),\;\; b=(0,4),\;\; c=(-12,0),\;\; o=(-4,0),\;\; \alpha=(-1,-3),\;\; \beta=(-8,-4),\;\;\gamma=(2,0)$$ and observe that
$\Aline ac\parallel \Aline ob\parallel \Aline o\beta\parallel\Aline\alpha\gamma\nparallel \Aline \beta\gamma\parallel\Aline o\alpha\parallel \Aline oa\parallel\Aline c\beta$ and $\Aline ab\parallel\Aline co\parallel\Aline o\gamma\nparallel\Aline \alpha\beta$.

\begin{picture}(200,100)(-210,-45)
\linethickness{0.8pt}
\put(-120,0){\color{blue}\line(1,1){40}}
\put(-120,0){\color{cyan}\line(1,-1){40}}
\put(-80,-40){\color{blue}\line(1,1){80}}
\put(-80,40){\color{cyan}\line(1,-1){70}}
\put(0,40){\color{cyan}\line(1,-2){20}}
\put(-10,-30){\color{blue}\line(1,1){30}}

\put(-80,40){\color{red}\line(1,0){80}}
\put(-120,0){\color{red}\line(1,0){140}}
\put(-80,-40){\color{orange}\line(7,1){70}}

\put(-120,0){\circle*{3}}
\put(-128,-3){$c$}
\put(-80,40){\circle*{3}}
\put(-83,43){$a$}
\put(-80,-40){\circle*{3}}
\put(-83,-50){$\beta$}
\put(-40,0){\circle*{3}}
\put(-42,-8){$o$}
\put(-10,-30){\circle*{3}}
\put(-12,-38){$\alpha$}
\put(0,40){\circle*{3}}
\put(2,43){$b$}
\put(20,0){\circle*{3}}
\put(23,-2){$\gamma$}
\end{picture}

\end{exercise}

\section{Inversive-plus affine spaces}

\begin{definition}\label{d:a+inversive} An affine liner $X$ is defined to be \index{inversive-plus affine space}\index{affine space!inversive-plus}\defterm{inversive-plus} if for every line $D\subset X$ and points $a,d,\alpha,\delta\in D$ and $b,c,\beta,\gamma\in X\setminus D$ $$\big[(\Aline ab\parallel\Aline cd\parallel\Aline \alpha\beta\parallel\Aline\gamma\delta\nparallel 
\Aline ac\parallel\Aline bd\parallel\Aline \alpha\gamma\parallel\Aline\beta\delta)\;\wedge\;(\Aline b\beta\parallel D)\big]\;\Ra\;\Aline c\gamma\parallel D.$$
\end{definition}

\begin{picture}(150,100)(-170,-40)

\linethickness{0.75pt}
\put(0,0){\color{teal}\line(-1,0){30}}
\put(0,0){\color{cyan}\line(0,-1){30}}
\put(-30,-30){\color{teal}\line(1,0){30}}
\put(-30,-30){\color{cyan}\line(0,1){30}}

{\linethickness{1pt}
\put(-30,-30){\color{red}\line(1,1){75}}
}
\put(-30,0){\color{red}\line(1,1){45}}
\put(0,-30){\color{red}\line(1,1){45}}

\put(15,15){\color{teal}\line(1,0){30}}
\put(15,15){\color{cyan}\line(0,1){30}}
\put(45,45){\color{teal}\line(-1,0){30}}
\put(45,45){\color{cyan}\line(0,-1){30}}

\put(-30,-30){\circle*{3}}
\put(-36,-37){$a$}
\put(0,0){\circle*{3}}
\put(2,-5){$d$}
\put(-30,0){\circle*{3}}
\put(-37,-1){$b$}
\put(0,-30){\circle*{3}}
\put(0,-38){$c$}
\put(15,45){\circle*{3}}
\put(11,48){$\beta$}
\put(45,45){\circle*{3}}
\put(47,46){$\delta$}
\put(45,15){\circle*{3}}
\put(47,10){$\gamma$}
\put(15,15){\circle*{3}}
\put(6,14){$\alpha$}

\end{picture}

\begin{theorem}\label{t:+inversive} For an affine space $X$ the following statements are equivalent:
\begin{enumerate}
\item the affine space $X$ is inversive-plus;
\item for every ternar $R$ of $X$, the plus loop $(R,+)$ is inversive;
\item for every ternar $R$ of $X$, the plus loop $(R,+)$ is left-inversive;
\item for every ternar $R$ of $X$, the plus loop $(R,+)$ is right-inversive;
\smallskip
\item for every ternar $R$ of $X$, the plus loop $(R,+)$ is left-Bol;
\item for every ternar $R$ of $X$, the plus loop $(R,+)$ is right-Bol;
\item for every ternar $R$ of $X$, the plus loop $(R,+)$ is Bol;
\item for every ternar $R$ of $X$, the plus loop $(R,+)$ is Moufang.
\end{enumerate}
\end{theorem}

\begin{proof} The equivalence $(1)\Leftrightarrow(4)$ follows from Proposition~\ref{p:rip<=>} .
\smallskip

$(1)\Ra(3)$ Assume that the affine plane is inversive-plus. Given any plane $\Pi\subseteq X$, an affine base $uow$ in $\Pi$ and its diagonal $\Delta$, we need to show that the additive loop $(\Delta,+)$ has the left-inverse property. By Proposition~\ref{p:lip<=>}, it suffices to check that for every points $x,y,z,s\in\Delta$ with $\Aline {xo}{oy}\subparallel \Delta$ and $\Aline {oz}{ys}\subparallel \Delta$, we have $\Aline {os}{xz}\subparallel \Delta$. Here we identify the based affine plane $(\Pi,uow)$ with its coordinate plane $\Delta^2$. If $x=o$, then $\Aline {xo}{oy}\subparallel \Delta$ and $\Aline {oz}{ys}\subparallel \Delta$ imply $y=o$ and $z=s$. In this case $\Aline{oz}{ys}=\{oz\}\subparallel \Delta$ and we are done. So, assume that $x\ne o$. In this case the subparallelity relations $\Aline {xo}{oy}\subparallel \Delta$ and $\Aline {oz}{ys}\subparallel \Delta$ imply $o\ne y$ and $z\ne s$.

\begin{picture}(100,135)(-130,-15)

{\linethickness{1pt}
\put(0,0){\color{red}\line(1,1){100}}
}
\linethickness{0.75pt}
\put(0,0){\color{cyan}\line(0,1){80}}
\put(0,80){\color{red}\line(1,1){20}}
\put(0,80){\color{teal}\line(1,0){80}}
\put(0,20){\color{teal}\line(1,0){20}}
\put(20,20){\color{cyan}\line(0,1){80}}
\put(20,40){\color{teal}\line(1,0){20}}
\put(20,100){\color{teal}\line(1,0){80}}
\put(0,20){\color{red}\line(1,1){20}}
\put(20,80){\color{red}\line(1,1){20}}
\put(40,40){\color{cyan}\line(0,1){60}}

\put(0,0){\circle*{3}}
\put(0,-7){$x$}
\put(20,20){\circle*{3}}
\put(20,13){$o$}
\put(13,22){$a$}
\put(40,40){\circle*{3}}
\put(40,33){$y$}
\put(34,42){$b$}
\put(80,80){\circle*{3}}
\put(80,73){$z$}
\put(100,100){\circle*{3}}
\put(100,93){$s$}

\put(0,20){\circle*{3}}
\put(-14,18){$xo$}
\put(2,13){$c$}
\put(20,40){\circle*{3}}
\put(7,40){$oy$}
\put(22,42){$d$}
\put(0,80){\circle*{3}}
\put(-14,78){$xz$}
\put(2,73){$\gamma$}
\put(20,100){\circle*{3}}
\put(16,103){$os$}
\put(22,92){$\delta$}
\put(40,100){\circle*{3}}
\put(35,104){$ys$}
\put(34,89){$\beta$}
\put(60,52){\color{red}$\Delta$}
\put(20,80){\circle*{3}}
\put(22,73){$oz$}
\put(12,82){$\alpha$}

\end{picture}

Consider the points $a\defeq oo$, $d\defeq oy$, $\alpha\defeq oz$, $\delta\defeq os$ on the line $V\defeq \Aline ow$, and the points $b\defeq yy$, $\beta\defeq ys$ and $c\defeq xo$. Let $\gamma\in \Aline {oz}{z}$ be a unique point such that $\Aline {\gamma\delta}\subparallel \Delta$. Observe that the lines $\Aline ab=\Aline oy=\Delta$, 
$\Aline cd=\Aline {xo}{oy}\subparallel\Delta$, $\Aline \alpha\beta=\Aline {oz}{ys}\subparallel\Delta$, and $\Aline \gamma\delta$ are (sub)parallel to the line $\Delta=\Aline oe$.
On the other hand, the lines $\Aline ac=\Aline {oo}{xo}$, $\Aline bd=\Aline {yy}{oy}$, $\Aline \alpha\gamma=\Aline{oz}{zz}$, $\Aline \beta\delta=\Aline {ys}{os}$ are parallel to the horizontal line $H\defeq \Aline ou$. Since the line $\Aline b\beta=\Aline{yy}{ys}$ is parallel to the line $V$, the additive inverse property of the affine plane $\Pi$ ensures that $\Aline c\gamma\parallel V$and hence $\Aline c\gamma=\Aline{xx}{xo}$ and hence $\gamma\in \Aline {xx}{xo}\cap\Aline {oz}{zz}=\{xz\}$. Then $\Aline {os}{xz}=\Aline{\delta}{\gamma}\parallel \Delta$, witnessing that the loop $(\Delta,+)$ has the left-inverse property.
\smallskip

$(3)\Ra(1)$ Assume that the condition (3) is satisfied. Given any line $D\subset\Pi$ and  points $a,d,\alpha,\delta\in D$ and $b,c,\beta,\gamma\in X\setminus D$ with $\Aline ab\parallel\Aline cd\parallel\Aline \alpha\beta\parallel\Aline\gamma\delta\nparallel\Aline ac\parallel\Aline bd\parallel\Aline \alpha\gamma\parallel\Aline\beta\delta$, and $\Aline b\beta\parallel D$, we have to prove that $\Aline c\gamma\parallel D$. Observe that the flat hull $\Pi$ of the set $\{a,b,c,d,\alpha,\beta,\gamma,\delta\}$ is a plane in the affine space $X$. Consider the points $o\defeq a$, $w\defeq d$, $e\defeq b$ and find a unique point $u\in\Pi$ such that $\Aline ou\parallel \Aline we$ and $\Aline ue\parallel \Aline ow$. It follows from $\Aline ad\nparallel \Aline db$ that $a\ne d$. Then $uow$ is an affine base for $\Pi$ and the point $e=b$ is the diunit of the affine base $uow$. Let $\Delta$ be the diagonal of the affine base $uow$. Identify the affine plane $\Pi$ with the coordinate plane $\Delta^2$ of the affine base $uow$. Let $y\defeq e=b$ and find a unique point $x\in\Delta$ such that $\Aline xc\parallel V=\Aline ow=D$. Also find unique points $z,s\in\Delta$ such that $\Aline \alpha z\parallel \Aline ou\parallel \Aline \beta s$.
In this case the points $\alpha,\beta$ have coordinates $oz$ and $ys$, respectively. The parallelity relations $\Aline ab\parallel\Aline cd\parallel\Aline \alpha\beta\parallel\Aline\gamma\delta$ imply  $\Aline {xo}{oy}=\Aline cd\parallel \Aline ab=\Delta$ and $\Aline {oz}{y\beta}=\Aline \alpha\beta\parallel \Delta$. The condition (3) ensures that the plus loop $(\Delta,+)$ has the left-inverse property. By Proposition~\ref{p:lip<=>}, $\Aline {xz}{os}\subparallel \Delta$. Taking into account that $\Aline \gamma{os}=\Aline \gamma\delta\parallel \Delta$, we conclude that $\Aline \gamma{os}=\Aline{xz}{os}$ and $\gamma=xz$. Then $\Aline c\gamma=\Aline {xo}{xz}\parallel \Aline ow=D$.
\smallskip

By the definition of the inverse property, $(2)\Leftrightarrow ((3)\wedge (4))$. Taking into account that $(3)\Leftrightarrow(1)\Leftrightarrow(4)$, we conclude that 
$(3)\Leftrightarrow(4)\Leftrightarrow(1)\Leftrightarrow((3)\wedge(4))\Leftrightarrow (2).$ Therefore, the conditions $(1)$--$(4)$ are equivalent.
\smallskip

$(1)\Ra(5)$ Assume that the affine space $X$ is inversive-plus. We have to prove that for every ternar $R$ of $X$, the plus loop $(R,+)$ is left-Bol.
Find a plane $\Pi\subseteq X$ and an affine base $uow$ in $\Pi$ whose ternar $\Delta$ is isomorphic to the ternar $R$. Now it suffices to check that the plus loop $(\Delta,+)$ is left-Bol. Given any elements $x,y,z\in\Delta$, we need to check that $x+(y+(x+z))=(x+(y+x))+z$.

 Consider the points $$
\begin{gathered}
a\defeq x+z,\quad b\defeq y+a=y+(x+z),\quad s\defeq x+b=x+(y+(x+z)),\\
c\defeq y+x\quad\mbox{and}\quad
d\defeq x+c=x+(y+x).
\end{gathered}$$ 
We need to show that $s=d+z$. 
By definition of the plus operation, the equalities $a=x+z$, $b=y+a$, $c=y+x$, $d=x+c$, $s=x+b$ imply that the flats $\Aline{oz}{xa},\Aline{oa}{yb},\Aline{ox}{yc},\Aline{oc}{xd},\Aline{ob}{xs}$ are subparallel to the diagonal $\Delta\defeq\Aline oe$ of the affine base $uow$. By the (already proved) implication (1)$\Ra$(2), the loop $(\Delta,+)$ is inversive and hence  there exists a point $\bar y\in \Delta$ such that $\bar y+(y+p)=p=(p+y)+\bar y$ for every $p\in\Delta$. In particular, $x=\bar y+(y+x)=\bar y+c$ and $a=\bar y+(y+a)=\bar y+b$. 
By definition of the plus operation, the equalities $x=\bar y+c$ and $a=\bar y+b$ imply that the flats $\Aline {oc}{\bar yx}$ and $\Aline{ob}{\bar ya}$ are subparallel to the diagonal $\Delta$ of the affine base $uow$. 

If $\bar y=x$, then $\bar y=x=\bar y+(y+x)$ implies $y+x=o$ and hence
$$x+(y+(x+z))=\bar y+(y+(x+z))=x+z=(x+o)+z=(x+(y+x))+z$$and we are done. So, assume that $\bar y\ne x$. 

The subparallelity of the flats $\Aline{oc}{xd},\Aline{oc}{\bar yx}$, $\Aline{ob}{xs},\Aline{ob}{\bar yx}$ to $\Delta$ implies that the flines $\Aline{xd}{\bar yx}$ and $\Aline{xs}{\bar ya}$ are parallel to $\Delta$. By Proposition~\ref{p:lip<=>}(3), the left-inversivity of the loop $(\Delta,+)$ ensures that $\Aline {xa}{ds}\subparallel \Aline {\bar yx}{xd}\parallel \Delta$ and hence $\Aline{xa}{ds}\subparallel \Delta$. Taking into account that $\Aline{xa}{oz}\subparallel \Delta$, we conclude that $\Aline {oz}{ds}\subparallel \Delta$ and hence $s=d+z$, by definition of the plus operation.



\smallskip

$(5)\Ra(7)$ Assume that for every ternar $R$ of $X$, the plus loop $(R,+)$ is left-Bol. By Proposition~\ref{p:left-Bol=>left-inversive}, every left-Bol loop is left-inversive, which implies that the condition (3) holds. Since $(3)\Leftrightarrow(4)$, for every ternar $R$ of $X$, the plus loop $(R,+)$ is right-inversive. By Theorem~\ref{t:Bol<=>}, every right-inversive left-Bol loop is Bol. Therefore, for every ternar $R$ of $X$, the plus loop $(R,+)$ is Bol.
\smallskip

The equivalence $(7)\Leftrightarrow(8)$ follows from Theorem~\ref{t:Bol<=>Moufang}; (7)$\Ra$(6) is trivial and (6)$\Ra$(4) follows from Proposition~\ref{p:right-Bol=>right-inversive}. 
\end{proof}

Theorem~\ref{t:+inversive} and Corollary~\ref{c:left-Bol-prime=>cyclic} imply the following corollary.

\begin{corollary}\label{c:inversive-plus-prime=>cyclic} If an inversive-plus affine space $X$ has prime order $p$, then for every ternar $R$ of $X$, the additive loop $(R,+)$ is a cyclic group of order $p$.
\end{corollary}

\begin{problem} Can Corollary~\textup{\ref{c:inversive-plus-prime=>cyclic}} be generalized to invertible-plus (instead of inversive-plus) affine spaces?
\end{problem}

\begin{example}[Hetman, 2024] There exists a based affine plane $X$ of order 9 such that
\begin{enumerate}
\item $X$ is not invertible-plus, 
\item the ternar $\Delta$ of $X$ is linear,
\item the multiplicative $0$-loop $(\Delta,\cdot)$ is associative (but not commutative),
\item the plus loop $(\Delta,+)$ is invertible-plus,
\item the plus loop $(\Delta,+)$ contains two elements of order 3 and six elements of order 2, and hence is not a group.
\end{enumerate} The linear ternar $\Delta$ is uniquely determined by its addition and multiplication tables:
$$
\begin{array}{c|ccccccccc}
+&0&1&2&3&4&5&6&7&8\\
\hline
0&0&1&2&3&4&5&6&7&8\\
1&1&8&6&7&5&4&2&3&0\\
2&2&4&0&6&1&8&3&5&7\\
3&3&5&7&0&8&1&4&2&6\\
4&4&7&3&8&0&6&5&1&2\\
5&5&6&8&2&7&0&1&4&3\\
6&6&3&4&1&2&7&0&8&5\\
7&7&2&1&5&6&3&8&0&4\\
8&8&0&5&4&3&2&7&6&1
\end{array}
\quad\mbox{and}
\quad
\begin{array}{c|ccccccccc}
\cdot&0&1&2&3&4&5&6&7&8\\
\hline
0&0&0&0&0&0&0&0&0&0\\
1&0&1&2&3&4&5&6&7&8\\
2&0&2&8&4&7&3&1&5&6\\
3&0&3&5&8&2&6&4&1&7\\
4&0&4&3&6&8&1&7&2&5\\
5&0&5&7&2&1&8&3&6&4\\
6&0&6&1&5&3&7&8&4&2\\
7&0&7&4&1&6&2&5&8&3\\
8&0&8&6&7&5&4&2&3&1\\
\end{array}.
$$
\end{example}

\section{Associative-plus affine spaces}

\begin{definition}\label{d:a+associative} An affine liner $X$ is defined to be \index{associative-plus affine space}\index{affine space!associative-plus}\defterm{associative-plus} if for every distinct points $a,b,c,d,\alpha,\beta,\gamma,\delta\in X$,$$\big[(\Aline ab\parallel\Aline cd\parallel\Aline \alpha\beta\parallel\Aline\gamma\delta\nparallel
\Aline bc\parallel\Aline ad\parallel\Aline \beta\gamma\parallel\Aline\alpha\delta)\;\wedge\;(\Aline a\alpha\parallel b\beta\parallel \Aline c\gamma)\big]\;\Ra\;\Aline d\delta\parallel\Aline a{\alpha}.$$
\end{definition}

\begin{picture}(100,90)(-150,-10)

\linethickness{0.7pt}
\put(0,0){\color{teal}\line(1,0){60}}
\put(0,0){\color{cyan}\line(0,1){30}}
\put(60,30){\color{cyan}\line(0,-1){30}}

\put(40,40){\color{teal}\line(1,0){60}}
\put(40,40){\color{cyan}\line(0,1){30}}
\put(100,70){\color{teal}\line(-1,0){60}}
\put(100,70){\color{cyan}\line(0,-1){30}}

\put(0,0){\color{red}\line(1,1){40}}
\put(60,30){\color{red}\line(1,1){40}}
\put(0,30){\color{red}\line(1,1){40}}
\put(60,0){\color{red}\line(1,1){40}}
\put(60,30){\color{teal}\line(-1,0){60}}

\put(0,0){\circle*{3}}
\put(-5,-8){$a$}
\put(60,0){\circle*{3}}
\put(60,-8){$d$}
\put(0,30){\circle*{3}}
\put(-7,28){$b$}
\put(60,30){\circle*{3}}
\put(62,25){$c$}

\put(40,40){\circle*{3}}
\put(31,39){$\alpha$}
\put(100,40){\circle*{3}}
\put(103,38){$\delta$}
\put(40,70){\circle*{3}}
\put(32,70){$\beta$}
\put(100,70){\circle*{3}}
\put(103,70){$\gamma$}
\end{picture}


\begin{theorem}\label{t:add-ass<=>} An affine space $X$ is associative-plus if and only if every ternar of $X$ is associative-plus.
\end{theorem}

\begin{proof} Assume that an affine space $X$ is associative-plus. Given any plane $\Pi\subseteq X$ and an affine base $uow$ in $\Pi$, we should prove that the additive loop $(\Delta,+)$ is associative. By Proposition~\ref{p:add-ass<=>}, it suffices to show that for every points $a,b,c,d,\alpha,\beta,\gamma,\delta\in \Pi$ with $\Aline{ac}{\alpha\gamma}\parallel \Aline{bc}{\beta\gamma}\parallel \Aline {bd}{\beta\delta}\parallel \Delta$, we have $\Aline {ad}{\alpha\delta}\parallel \Delta$. Those assumptions imply $ac\ne\alpha\gamma$, $bc\ne\beta\gamma$, $bd\ne\beta\delta$. 

If $a=b$, then $\Aline {ac}{\alpha\gamma}\parallel \Aline{bc}{\beta\gamma}\parallel \Delta$ implies $\alpha=\beta$ and then $\Aline {ad}{\alpha\delta}=\Aline{bd}{\beta\delta}\parallel \Delta$. 

If $c=d$, then $\Aline{bc}{\beta\gamma}\parallel \Aline {bd}{\beta\delta}\parallel \Delta$ implies $\gamma=\delta$ and hence $\Aline {ad}{\alpha\delta}=\Aline{ac}{\alpha\gamma}\parallel\Delta$.

So, we assume that $a\ne b$ and $c\ne d$. In this case $\Aline{ac}{\alpha\gamma}\parallel \Aline{bc}{\beta\gamma}\parallel \Aline {bd}{\beta\delta}\parallel \Delta$ implies $\alpha\ne \beta$ and $\gamma\ne \delta$. Since $\Aline {ac}{ad}\parallel \Aline {bc}{bd}\parallel \Aline {\alpha\gamma}{\alpha\delta}\parallel \Aline {\beta\gamma}{\beta\delta}\nparallel \Aline {ac}{bc}\parallel \Aline {ad}{bd}\parallel \Aline {\alpha\gamma}{\beta\gamma}\parallel \Aline{\alpha\delta}{\beta\delta}$ and $\Aline{ac}{\alpha\gamma}\parallel \Aline{bc}{\beta\gamma}\parallel \Aline {bd}{\beta\delta}\parallel \Delta$, we can apply the additive associativity of $\Pi$ and conclude that $\Aline {ad}{\alpha\delta}\parallel \Delta$. This completes the proof of the ``only if'' part.
\smallskip

In the proof of the ``if'' part we shall use the following lemma (which will be also applied in the proof of Theorem~\ref{t:add-com=>add-ass}.

\begin{lemma}\label{l:add-ass<=>} Let $\Pi$ be an affine plane and $a,b,c,d,\alpha,\beta,\gamma,\delta$ be points in $\Pi$ such that $$\Aline ab\parallel\Aline cd\parallel\Aline \alpha\beta\parallel\Aline\gamma\delta\nparallel
\Aline bc\parallel\Aline ad\parallel\Aline \beta\gamma\parallel\Aline\alpha\delta\quad\mbox{and}\quad\Aline a\alpha\parallel b\beta\parallel \Aline c\gamma.$$
If $(\Aline ab\subparallel\Aline a\alpha)\;\vee\;(\Aline a\alpha\subparallel \Aline ab)\;\vee\;(\Aline ad\subparallel\Aline a\alpha)\;\vee\;(\Aline a\alpha\subparallel \Aline ad)$, then $\Aline d\delta\parallel \Aline a\alpha$.
\end{lemma}

\begin{proof} We divide the proof into four claims.

\begin{claim}\label{cl:add-ass4} If $\Aline ab\subparallel \Aline a\alpha$, then $\Aline d\delta\parallel \Aline a\alpha$.
\end{claim}

\begin{proof} If $a=b$, then the parallelity relations $\Aline ab\parallel\Aline cd\parallel\Aline \alpha\beta\parallel\Aline\gamma\delta$ imply $c=d$, $\alpha=\beta$, and $\gamma=\delta$. Then $\Aline d\delta=\Aline c\gamma\parallel \Aline a\alpha$ and we are done. 

So, assume that $a\ne b$. In this case the subparallelity relation $\Aline ab\subparallel \Aline a\alpha$ implies $a\ne\alpha$ and $\Aline ab=\Aline a\alpha$, by Corollary~\ref{c:subparallel}. The parallelity relations $$\Aline c\gamma\parallel \Aline b\beta\parallel \Aline a\alpha=\Aline ab\parallel\Aline cd\parallel\Aline \alpha\beta\parallel\Aline\gamma\delta$$
imply $c\ne d$, $\alpha\ne\beta$, $\gamma\ne\delta$, $b\ne\beta$, $c\ne\gamma$ and $\Aline cd=\Aline  c\gamma=\Aline\gamma\delta\parallel \Aline ab=\Aline a\alpha$.  Assuming that $d=\delta$ and taking into account that $\Aline ad\parallel \Aline \alpha\delta\nparallel \Aline ab$, we conclude that $\Aline ad=\Aline \alpha\delta$ and $\{a\}=\Aline ad\cap\Aline ab=\Aline\alpha\delta\cap\Aline a{\alpha}=\{\alpha\}$, which contradicts $a\ne b$ and $\Aline ab\subparallel \Aline a\alpha$. This contradiction shows that $d\ne\delta$. In this case $d,\delta\in \Aline c\gamma$ implies $\Aline d\delta=\Aline c\gamma\parallel \Aline a\alpha$.
\end{proof}

\begin{claim}\label{cl:add-ass5} If $\Aline a\alpha\subparallel \Aline ab$, then $\Aline d\delta\parallel \Aline a\alpha$.
\end{claim}

\begin{proof} If $a=\alpha$, then the parallelity relations $\Aline a\alpha=\Aline b\beta=\Aline c\gamma$ imply $b=\beta$ and $c=\gamma$. Since $a=\alpha$ and $b=\beta$, the parallelity relations $$\Aline ad\parallel \Aline \alpha\delta\parallel \Aline bc\parallel \Aline\beta\gamma\nparallel \Aline cd\parallel \Aline \gamma\delta$$ imply $\Aline ad=\Aline \alpha\delta$, $\Aline bc=\Aline\beta\gamma$, $\Aline cd=\Aline \gamma\delta$,
$$\{d\}=\Aline ad\cap\Aline cd=\Aline\alpha\delta\cap\Aline \gamma\delta=\{\delta\},$$and finally,
$\Aline d{\delta}=\{d\}\parallel \{a\}=\Aline a\alpha$. 

So, assume that $a\ne\alpha$. In this case the subparallelity relation $\Aline a\alpha\subparallel\Aline ab$ implies $a\ne b$ and $\Aline ab\subparallel \Aline a\alpha$, by Corollary~\ref{c:subparallel}. Applying Claim~\ref{cl:add-ass4}, we conclude that $\Aline d\delta\parallel\Aline a\alpha$.
\end{proof}

By analogy with Claims~\ref{cl:add-ass4} and \ref{cl:add-ass5}, we can prove the following two claims completing the proof of the lemma.

\begin{claim} If $\Aline ad\subparallel \Aline a\alpha$, then $\Aline d\delta\parallel \Aline a\alpha$.
\end{claim}

\begin{claim} If $\Aline a\alpha\subparallel \Aline ad$, then $\Aline d\delta\parallel \Aline a\alpha$.
\end{claim}
\end{proof}

With Lemma~\ref{l:add-ass<=>} in our disposition, we now can present the proof of the ``if'' part of Theorem~\ref{t:add-ass<=>}. Assume that every ternar of an affine space $X$ is associative-plus. To prove that the affine space $X$ is associative-plus, take any points 
$a,b,c,d,\alpha,\beta,\gamma,\delta\in X$ with $\Aline ab\parallel\Aline cd\parallel\Aline \alpha\beta\parallel\Aline\gamma\delta\nparallel
\Aline bc\parallel\Aline ad\parallel\Aline \beta\gamma\parallel\Aline\alpha\delta$ and $\Aline a\alpha\parallel b\beta\parallel \Aline c\gamma$. We have to show that $\Aline a\alpha\parallel\Aline d{\delta}$. 

If the flat hull $\Pi$ of the set $\{a,b,c,d,\alpha,\beta,\gamma,\delta\}$ has rank $\|\Pi\|\ge 4$, then the parallel lines $\Aline ab$, $\Aline cd$, $\Aline \alpha\beta$, $\Aline\gamma\delta$ are distinct. Since $\|X\|\ge\|\Pi\|\ge 4$, the affine space $X$ is Desarguesian and Thalesian, by Corollary~\ref{c:affine-Desarguesian} and Theorem~\ref{t:ADA=>AMA}. In this case, the parallelity relations $\Aline ad\parallel\Aline bc$ and $\Aline a\alpha\parallel \Aline b\beta$ imply $\Aline \alpha d\parallel \Aline \beta c$. Applying the Thalesian Axiom to the triangles $d\alpha\delta$ and $c\beta\gamma$, we conclude that $\Aline d\alpha\parallel\Aline c\beta$ and $\Aline \alpha\delta\parallel \Aline\beta\gamma$ imply $\Aline d\delta\parallel \Aline c\gamma$.

So, assume that $\|\Pi\|= 3$, which means that $\Pi$ is a plane. 
If $(\Aline ab\subparallel\Aline a\alpha)\;\vee\;(\Aline a\alpha\subparallel \Aline ab)\;\vee\;(\Aline ad\subparallel\Aline a\alpha)\;\vee\;(\Aline a\alpha\subparallel \Aline ad)$, then $\Aline d\delta\parallel \Aline a\alpha$, by Lemma~\ref{l:add-ass<=>}. So, assume that $\Aline ab,\Aline ad,\Aline a\alpha$ are non-parallel lines. In this case we can choose an affine base $uow$ in the plane $\Pi$ such that $\Aline ab\parallel \Aline ou$, $\Aline ad\parallel \Aline ow$ and $\Aline a\alpha\parallel \Delta\defeq\Aline oe$, where $e$ is the diunit of the affine base $uow$. By our assumption, the additive loop $(\Delta,+)$ is associative. Applying Proposition~\ref{p:add-ass<=>}, we conclude that $\Aline d\delta\parallel \Delta\parallel \Aline a\alpha$. 
\end{proof}



Let us recall that an affine space $X$ is called \index{$\partial$-translation}\index{affine space!$\partial$-translation}\defterm{$\partial$-translation} if for every direction $\boldsymbol \delta\in \partial X$ there exists a non-identity translation $T:X\to X$ such that $\boldsymbol\delta=\{\Aline xy:xy\in T\}$.

\begin{theorem}\label{t:add-ass=>partialT} Let $X$ be an affine space of  order $|X|_2=p^n$ for some prime number $p$ and some $n\in\IN$. If $X$ is associative-plus, then $X$ is $\partial$-translation.
\end{theorem}

\begin{proof} Assume that $X$ is associative-plus and $|X|_2$ is a prime power. If $\|X\|\ge 4$, then $X$ is Desarguesian, Thalesian, translation, and $\partial$-translation, by Corollary~\ref{c:affine-Desarguesian} and  Theorems~\ref{t:ADA=>AMA} and \ref{t:paraD<=>translation}.  So, assume that $\|X\|=3$, which means that $X$ is a Playfair plane. To show that $X$ is $\partial$-translation, we shall apply Gleason's Theorem~\ref{t:partial-translation<=}. Given any line $\Delta\subseteq X$ and two distinct directions ${\boldsymbol h},{\boldsymbol v}\in\partial X\setminus\{\Delta_\parallel\}$, we should show that the set $\I_X^\#[\Delta;{\boldsymbol h},{\boldsymbol v}]$ is a subgroup of the group $\I_X^\#[\Delta]$. Choose any distinct points $o,e\in\Delta$ and consider the unique points $u\defeq \Aline o{\boldsymbol h}\cap\Aline e{\boldsymbol v}$ and $w\defeq\Aline o{\boldsymbol v}\cap\Aline e{\boldsymbol h}$. Then $uow$ is an affine base such that ${\boldsymbol h}=(\Aline ou)_\parallel$, ${\boldsymbol v}=(\Aline ow)_\parallel $, and $\Delta=\Aline oe$. The affine base $uow$ determines a canonical structure of a ternar on the diagonal $\Delta=\Aline oe$ of the based affine plane $(X,uow)$. Since $X$ is associative-plus, so is the ternar $\Delta$. By Theorem~\ref{t:plus<=>group}, the set $\I_X^\#[\Delta;{\boldsymbol h},{\boldsymbol v}]$ is a subgroup of the group $\I_X^\#[\Delta]$ of all line translations of the line $\Delta$. By Theorem~\ref{t:partial-translation<=}, the affine plane $X$ is $\partial$-translation.
\end{proof}

\section{Commutative-plus affine spaces}

\begin{definition}\label{d:a+commutative} An affine liner $X$ is defined to be \index{commutative-plus affine space}\index{affine space!commutative-plus}\defterm{commutative-plus} if for every points $a,b,c,\alpha,\beta,\gamma\in X$,$$\big((\Aline ab\parallel\Aline \alpha c\parallel\Aline \beta\gamma\nparallel \Aline \alpha\beta\parallel\Aline a\gamma\parallel\Aline bc)\;\wedge\;(\Aline a\alpha\parallel b\beta)\big)\;\Ra\;\Aline b\beta\parallel\Aline c\gamma.$$
\end{definition}

\begin{picture}(150,100)(-150,-10)
\linethickness{0.75pt}
{\linethickness{1pt}
\put(0,0){\color{red}\line(1,1){80}}
}
\put(0,0){\color{cyan}\line(0,1){60}}
\put(0,60){\color{teal}\line(1,0){80}}
\put(0,0){\color{teal}\line(1,0){20}}
\put(20,0){\color{cyan}\line(0,1){80}}
\put(80,80){\color{teal}\line(-1,0){60}}
\put(0,60){\color{red}\line(1,1){20}}
\put(80,60){\color{cyan}\line(0,1){20}}
\put(20,0){\color{red}\line(1,1){60}}

\put(0,0){\circle*{3}}
\put(-7,-4){$b$}
\put(80,80){\circle*{3}}
\put(83,78){$\beta$}
\put(20,0){\circle*{3}}
\put(23,-4){$a$}
\put(0,60){\circle*{3}}
\put(-8,57){$c$}
\put(20,80){\circle*{3}}
\put(12,81){$\gamma$}
\put(80,60){\circle*{3}}
\put(83,57){$\alpha$}
\end{picture}


\begin{theorem}\label{t:add-com<=>} An affine space $X$ is commutative-plus if and only if every ternar of $X$ is commutative-plus.
\end{theorem}

\begin{proof} Assume that an affine space $X$ is commutative-plus. Given any plane $\Pi\subseteq X$ and an affine base $uow$ in $\Pi$, we have to prove that the ternar $\Delta$ of the based affine plane $(\Pi,uow)$  is commutative-plus. By Proposition~\ref{p:add-com<=>}, it suffices to prove that for every points  $x,y,z\in \Delta$ with $\Aline {ox}{yz}\subparallel \Delta$, we have $\Aline {oy}{xz}\subparallel \Delta$.

If $o=y$ or $x=z$, then $\Aline {ox}{yz}\subparallel \Delta$ implies $o=y$ and $x=z$ and hence $\Aline {oy}{xz}=\Aline {oo}{xx}\subparallel\Delta$. If $o=x$ ot $y=z$, then the subparallelity relation $\Aline {ox}{yz}\subparallel \Delta$ implies $o=x$ and $y=z$. In this case $\Aline {oy}{xz}=\{oy\}\subparallel\Delta$. 
If $x=y$, then $\Aline{oy}{xz}=\Aline{ox}{yz}\subparallel\Delta$. 

So, assume that $x\ne o\ne y$, $x\ne z\ne y$, and $x\ne y$. 
Then for the points $a\defeq ox$, $b\defeq xx$, $\alpha\defeq yz$, $\beta\defeq yy$, $\gamma\defeq oy$, $c\defeq xz$, we have $\Aline ab\parallel \Aline \alpha c\parallel \Aline \beta\gamma$, $\Aline\alpha\beta\parallel\Aline a\gamma\parallel\Aline bc$, and $\Aline a\alpha=\Aline{ox}{yz}\parallel \Delta=\Aline b\beta$.
Since the affine plane $\Pi$ is commutative-plus, $\Aline {oy}{xz}=\Aline c\gamma\parallel \Aline b\beta=\Delta$, witnessing that the additive loop $(\Delta,+)$ is commutative.
\smallskip

Next, assume that every ternar of $X$ is commutative-plus. To prove that the affine space $X$ is commutative-plus, take any points 
$a,b,c,\alpha,\beta,\gamma\in X$ with $\Aline ab\parallel\Aline \alpha c\parallel\Aline \beta\gamma\nparallel \Aline \alpha\beta\parallel\Aline a\gamma\parallel\Aline bc$ and $\Aline a\alpha\parallel \Aline b\beta$. We have to prove that $\Aline c\gamma\parallel \Aline b\beta$.  Observe that the flat hull $\Pi$ of the set $\{a,b,c,\alpha,\beta,\gamma\}$ is a plane in $X$.

If $a=b$, then the parallelity relations  $\Aline ab\parallel\Aline \alpha c\parallel\Aline \beta\gamma\nparallel \Aline \alpha\beta\parallel\Aline a\gamma\parallel\Aline bc$ imply $\alpha=c$, $\beta=\gamma$, $\alpha\ne\beta$, $a\ne\gamma$, $b\ne c$ and $\Aline c\gamma=\Aline \alpha\beta=\Aline a\gamma=\Aline b\beta$.

By analogy we can prove that $\alpha=\beta$ implies $\Aline c\gamma=\Aline ba=\Aline c\alpha=\Aline b\beta$.

So, assume that $a\ne b$ and $\alpha\ne\beta$. If $\Aline a\alpha\subparallel \Aline ab$, then $\Aline a\alpha\parallel \Aline b\beta$ implies $\alpha,\beta\in \Aline ab$, which contradicts $\Aline ab\nparallel \Aline \alpha\beta$. This contradiction shows that $a\ne\alpha$ and $\Aline a\alpha\nparallel \Aline ab$. By analogy we can show that $\Aline a\alpha\nparallel \Aline \alpha\beta$. Then we can choose an affine base $uow$ in the affine plane $\Pi$ such that $o\in \Aline b\beta\cap\Aline a\gamma$, $\Aline ou\parallel \Aline ab$, $\Aline ow=\Aline a\gamma$ and $\Delta=\Aline b\beta$, where $\Delta$ is the diagonal of the affine base $uow$. Then for the points $x\defeq b$, $y\defeq \beta$, and $z\in \Aline c\alpha\cap \Delta$, we have $a=ox$, $\alpha=yz$, $c=xz$, $\gamma=oy$ and $\Aline {ox}{yz}=\Aline a\alpha\parallel \Delta$. Applying Proposition~\ref{p:add-com<=>}, we conclude that $\Aline c{\gamma}=\Aline {xz}{oy}=\Aline{oy}{xz}\subparallel \Delta=\Aline b\beta$.
\end{proof}

\begin{lemma}\label{l:add-com} If an affine space $X$ is commutative-plus, then for every points $a,b,c,\alpha,\beta,\gamma\in X$,$$\big((\Aline ab\parallel\Aline \alpha c\parallel\Aline \beta\gamma\nparallel \Aline \alpha\beta\parallel\Aline a\gamma\parallel\Aline bc)\;\wedge\;\big((\Aline a\alpha\subparallel b\beta)\;\vee\;(\Aline b\beta\subparallel \Aline a\alpha)\big)\big)\;\Ra\;\big((\Aline c\gamma\parallel \Aline a\alpha)\;\vee\;(\Aline c\gamma\parallel\Aline b\beta)\big).$$
\end{lemma}

\begin{proof} Assumming that $$(\Aline ab\parallel\Aline \alpha c\parallel\Aline \beta\gamma\nparallel \Aline \alpha\beta\parallel\Aline a\gamma\parallel\Aline bc)\;\wedge\;\big((\Aline a\alpha\subparallel b\beta)\;\vee\;(\Aline b\beta\subparallel \Aline a\alpha)\big),$$
we should prove that the flat $\Aline c\gamma$ is subparallel to $\Aline a\alpha$ or $\Aline b\beta$. 

If $a=b$, then the parallelity relations $\Aline ab\parallel\Aline \alpha c\parallel\Aline \beta\gamma\nparallel \Aline \alpha\beta\parallel\Aline a\gamma\parallel\Aline bc$ imply $\alpha=c$, $\beta=\gamma$, $\alpha\ne\beta$, $a\ne\gamma$, $b\ne c$, and $\Aline \alpha\beta=\Aline a\gamma=\Aline b\beta$. Then $\Aline c\gamma=\Aline \alpha\beta=\Aline b\beta$ and hence $\Aline c\gamma\parallel \Aline b\beta$.

By an analogous argument, we can prove that $\alpha=\beta$ implies $\Aline c\gamma\parallel \Aline b\beta$. 

So, assume that $a\ne b$ and $\alpha\ne\beta$. 

If $a=\alpha$, then $\Aline ab\parallel \Aline \alpha c$ implies $\Aline ab=\Aline \alpha c$ and hence $\{b\}=\Aline cb\cap\Aline ab=\Aline cb\cap\Aline\alpha c=\{c\}$. On the other hand, $a=\alpha$ and $\Aline a\gamma\parallel \Aline\alpha\beta$ imply $\Aline a\gamma=\Aline \alpha\beta$ and hence $\{\gamma\}=\Aline \gamma\beta\cap\Aline a\gamma=\Aline\gamma\beta\cap\Aline \alpha\beta=\{\beta\}$. Finally, $\Aline c\gamma=\Aline b\beta$.
By analogy we can prove that $b=\beta$ implies $\Aline c\gamma=\Aline a\alpha$. 

So, assume that $a\ne\alpha$ and $b\ne\beta$. In this case $(\Aline a\alpha\subparallel b\beta)\;\vee\;(\Aline b\beta\subparallel \Aline a\alpha)$ implies $\Aline a\alpha\parallel \Aline b\beta$. Applying the commutativity-plus of $X$, we conclude that $\Aline c\gamma\parallel \Aline b\beta$.
\end{proof}

The following theorem is a ``plus'' counterpart of Hessenberg's  Theorem~\ref{t:Hessenberg-affine}.

\begin{theorem}\label{t:add-com=>add-ass} Every  commutative-plus affine space $X$ is associative-plus.
\end{theorem}

\begin{proof} Given any points $a,b,c,d,\alpha,\beta,\gamma,\delta\in X$ with $$\Aline ab\parallel\Aline cd\parallel\Aline \alpha\beta\parallel\Aline\gamma\delta\nparallel\Aline bc\parallel\Aline ad\parallel\Aline \beta\gamma\parallel\Aline\alpha\delta\quad\mbox{and}\quad\Aline a\alpha\parallel b\beta\parallel \Aline c\gamma,$$ we have to prove that $\Aline a\alpha\parallel\Aline d{\delta}.$

\begin{picture}(100,95)(-150,-15)

\linethickness{0.7pt}
\put(0,0){\color{teal}\line(1,0){60}}
\put(0,0){\color{cyan}\line(0,1){30}}
\put(60,30){\color{cyan}\line(0,-1){30}}

\put(40,40){\color{teal}\line(1,0){60}}
\put(40,40){\color{cyan}\line(0,1){30}}
\put(100,70){\color{teal}\line(-1,0){60}}
\put(100,70){\color{cyan}\line(0,-1){30}}

\put(0,0){\color{red}\line(1,1){40}}
\put(60,30){\color{red}\line(1,1){40}}
\put(0,30){\color{red}\line(1,1){40}}
\put(60,0){\color{red}\line(1,1){40}}
\put(60,30){\color{teal}\line(-1,0){60}}

\put(0,0){\circle*{3}}
\put(-5,-8){$a$}
\put(60,0){\circle*{3}}
\put(60,-8){$d$}
\put(0,30){\circle*{3}}
\put(-7,28){$b$}
\put(60,30){\circle*{3}}
\put(62,25){$c$}

\put(40,40){\circle*{3}}
\put(31,39){$\alpha$}
\put(100,40){\circle*{3}}
\put(103,38){$\delta$}
\put(40,70){\circle*{3}}
\put(32,70){$\beta$}
\put(100,70){\circle*{3}}
\put(103,70){$\gamma$}
\end{picture}

If the flat hull $\Pi$ of the set $\{a,b,c,d,\alpha,\beta,\gamma,\delta\}$ is not a plane, then $\|X\|\ge\|\Pi\|\ge 4$ and  the affine space $X$ is associative-plus, by Theorem~\ref{t:corps<=>}. So, assume that $\Pi$ is a plane, which is commutative-plus (being a flat in the commutative-plus affine space $X$).
\smallskip

If $(\Aline ab\subparallel\Aline a\alpha)\vee(\Aline a\alpha\subparallel \Aline ab)\vee(\Aline ad\subparallel\Aline a\alpha)\vee(\Aline a\alpha\subparallel \Aline ad)$, then $\Aline d\delta\parallel \Aline a\alpha$, by Lemma~\ref{l:add-ass<=>}. In the opposite case, $b\ne a\ne d$, $a\ne\alpha$, and the lines $\Aline ab$, $\Aline ad$, $\Aline a\alpha$ are not parallel.

\begin{claim}\label{cl:add-com=>add-ass1} If $\gamma=a$ and $c=\alpha$, then $\Aline d\delta=\Aline\beta b\parallel \Aline a\alpha$.
\end{claim}

\begin{proof}  It follows from $\gamma=a$ and $\Aline ab\parallel \Aline \gamma\delta\nparallel \Aline ad\parallel \Aline \gamma\beta$ that $\Aline ab= \Aline \gamma\delta\nparallel \Aline ad=\Aline \gamma\beta$. On the other hand, $c=\alpha$ and $\Aline \alpha\beta\parallel\Aline cd\nparallel \Aline \alpha\delta\parallel \Aline cb$ imply $\Aline cd=\Aline \alpha\beta\nparallel \Aline cb=\Aline \alpha\delta$. Then $\{\delta\}=\Aline \alpha\delta\cap \Aline\gamma\delta=\Aline cb\cap \Aline ab=\{b\}$ and $\{d\}=\Aline ad\cap\Aline cd=\Aline\gamma\beta\cap\Aline \alpha\beta=\{\beta\}$. Then $\Aline d\delta=\Aline \beta b\parallel \Aline a\alpha$.
\end{proof}

 \begin{claim}\label{cl:add-com=>add-ass2} If $c=\alpha$ and $a\ne\gamma$, then $\Aline d\delta\parallel \Aline a\alpha$.
\end{claim}

\begin{proof} It follows from $c=\alpha$, $a\ne\gamma$ and $\Aline \gamma\beta\parallel \Aline da\nparallel \Aline \gamma c\parallel \Aline \alpha a\nparallel \Aline ab\parallel \Aline\gamma\delta$ that $\Aline a\gamma=\Aline \gamma c=\Aline \alpha a$, $\beta\ne d$ and $b\ne\delta$. Observe that $$\Aline \gamma\delta\parallel \Aline \beta\alpha=\Aline \beta d=\Aline cd\parallel\Aline ba\nparallel ad\parallel \Aline bc=\Aline b\delta=\Aline \alpha\delta\parallel\Aline \beta\gamma.$$ 

\begin{picture}(100,85)(-150,-15)

\linethickness{0.7pt}
\put(0,0){\color{teal}\line(1,0){30}}
\put(0,0){\color{cyan}\line(0,1){30}}
\put(30,30){\color{cyan}\line(0,-1){30}}

\put(30,30){\color{teal}\line(1,0){30}}
\put(30,30){\color{cyan}\line(0,1){30}}
\put(60,60){\color{teal}\line(-1,0){30}}
\put(60,60){\color{cyan}\line(0,-1){30}}

\put(0,0){\color{red}\line(1,1){30}}
\put(30,30){\color{red}\line(1,1){30}}
\put(0,30){\color{red}\line(1,1){30}}
\put(30,0){\color{red}\line(1,1){30}}
\put(30,30){\color{teal}\line(-1,0){30}}

\put(0,0){\circle*{3}}
\put(-5,-8){$a$}
\put(30,0){\circle*{3}}
\put(30,-8){$d$}
\put(0,30){\circle*{3}}
\put(-7,28){$b$}
\put(30,30){\circle*{3}}
\put(32,23){$c$}

\put(30,30){\circle*{3}}
\put(22,32){$\alpha$}
\put(60,30){\circle*{3}}
\put(63,28){$\delta$}
\put(30,60){\circle*{3}}
\put(26,63){$\beta$}
\put(60,60){\circle*{3}}
\put(60,64){$\gamma$}
\end{picture}

Since $\Aline b\beta\parallel \Aline \gamma a$, we can apply the commutative-plus property of the affine plane $\Pi$ and conclude that $\Aline d\delta\parallel \Aline \gamma a=\Aline \gamma c\parallel \Aline a\alpha$.
\end{proof}

By analogy we can prove   

 \begin{claim}\label{cl:add-com=>add-ass3} If $a=\gamma$ and $c\ne\alpha$, then $\Aline d\delta\parallel \Aline a\alpha$.
\end{claim}

\begin{claim} If $\Aline c\gamma=\Aline a\alpha$, then $\Aline d\delta\parallel \Aline a\alpha$.
\end{claim}

\begin{proof} If $a=\gamma$ or $c=\alpha$, then the parallelity relation $\Aline d\delta\parallel \Aline a\alpha$ follows from Claims~\ref{cl:add-com=>add-ass1}, \ref{cl:add-com=>add-ass2}, \ref{cl:add-com=>add-ass3}. So, assume that $a\ne\gamma$ and $c\ne \alpha$. Consider the unique points $x\in \Aline bc\cap\Aline \gamma\delta$ and $y\in \Aline ad\cap \Aline\alpha\beta$.

\begin{picture}(150,95)(-170,-45)

\linethickness{0.75pt}
\put(0,0){\color{teal}\line(-1,0){30}}
\put(0,0){\color{teal}\line(1,0){45}}
\put(0,0){\color{cyan}\line(0,-1){30}}
\put(-30,-30){\color{teal}\line(1,0){45}}
\put(-30,-30){\color{cyan}\line(0,1){30}}

{\linethickness{1pt}
\put(-30,-30){\color{red}\line(1,1){75}}
}
\put(-30,0){\color{red}\line(1,1){45}}
\put(0,-30){\color{red}\line(1,1){45}}

\put(15,15){\color{teal}\line(1,0){30}}
\put(15,-30){\color{cyan}\line(0,1){75}}
\put(45,45){\color{teal}\line(-1,0){30}}
\put(45,45){\color{cyan}\line(0,-1){45}}
\put(15,-30){\color{red}\line(1,1){30}}

\put(-30,-30){\circle*{3}}
\put(-36,-37){$a$}
\put(0,0){\circle*{3}}
\put(2,-6){$c$}
\put(-30,0){\circle*{3}}
\put(-37,-1){$b$}
\put(0,-30){\circle*{3}}
\put(-2,-39){$d$}
\put(15,45){\circle*{3}}
\put(11,48){$\beta$}
\put(45,45){\circle*{3}}
\put(47,46){$\gamma$}
\put(45,15){\circle*{3}}
\put(47,10){$\delta$}
\put(15,15){\circle*{3}}
\put(6,14){$\alpha$}
\put(45,0){\circle*{3}}
\put(48,-2){$x$}
\put(15,-30){\circle*{3}}
\put(12,-38){$y$}
\end{picture}

The commutativity-plus of the affine pace $X$ implies $\Aline xy\parallel \Aline a\gamma=\Aline \alpha c$. Applying the commutative-plus property to the points $x,c,\delta,y,\alpha,d$, we conclude that $\Aline d\delta\parallel \Aline c\alpha=\Aline a\alpha$.
\end{proof}

So, assume that $\Aline c\gamma\ne\Aline a\alpha$. 
Since the lines $\Aline ad$ and $\Aline c\gamma$ are not parallel, there exists a unique point $x\in \Aline ad\cap\Aline c\gamma$. It follows from $\Aline a\alpha\ne \Aline c\gamma\ne \Aline \gamma\delta$ that $x\notin\Aline ab\cup\Aline cd\cup\Aline \alpha\beta\cup\Aline \gamma\delta$.

Let $y\in \Aline b\beta$ be a unique point such that $\Aline xy\parallel \Aline\gamma\delta$. Let $s$ be the unique common point of the lines $\Aline ad$ and $\Aline\alpha\beta$. Since $x\notin\Aline \gamma\delta$, the line $\Aline xy$ is disjoint with the line $\Aline \gamma\delta$ and hence there exists a unique point $t\in \Aline \gamma\delta$ such that $\Aline yt\parallel \Aline ad$. Applyng the  commutative-plus property of $\Pi$ to the points $y,x,s,\beta,\gamma,t$, we conclude that $\Aline st\parallel \Aline x\gamma$.

\begin{picture}(120,135)(-120,-15)
\linethickness{0.75pt}
\put(0,0){\color{cyan}\line(0,1){15}}
\put(0,0){\color{teal}\line(1,0){90}}
\put(0,15){\color{teal}\line(1,0){105}}
\put(0,0){\color{red}\line(1,1){105}}
\put(0,15){\color{red}\line(1,1){90}}
\put(45,0){\color{cyan}\line(0,1){60}}
\put(60,0){\color{cyan}\line(0,1){60}}
\put(90,0){\color{cyan}\line(0,1){105}}
\put(105,15){\color{cyan}\line(0,1){90}}
\put(45,0){\color{red}\line(1,1){45}}
\put(60,0){\color{red}\line(1,1){45}}
\put(90,45){\color{teal}\line(1,0){15}}
\put(45,60){\color{teal}\line(1,0){15}}
\put(90,0){\color{red}\line(1,1){15}}
\put(90,105){\color{teal}\line(1,0){15}}

\put(0,0){\circle*{3}}
\put(-3,-8){$x$}
\put(45,0){\circle*{3}}
\put(43,-8){$a$}
\put(60,0){\circle*{3}}
\put(59,-9){$d$}
\put(90,0){\circle*{3}}
\put(88,-8){$s$}
\put(0,15){\circle*{3}}
\put(-4,19){$y$}
\put(60,15){\circle*{3}}
\put(54,17){$z$}
\put(105,15){\circle*{3}}
\put(108,13){$t$}
\put(90,45){\circle*{3}}
\put(81,45){$\alpha$}
\put(105,45){\circle*{3}}
\put(108,43){$\delta$}
\put(45,60){\circle*{3}}
\put(41,63){$b$}
\put(60,60){\circle*{3}}
\put(56,63){$c$}
\put(90,105){\circle*{3}}
\put(85,108){$\beta$}
\put(105,105){\circle*{3}}
\put(105,108){$\gamma$}

\end{picture}

Let $z$ be the unique common point of the lines $\Aline cd$ and $\Aline yt$. Applying the commutative-plus property of $\Pi$ to the points $y,x,a,b,c,z$ we conclude that $\Aline az\parallel \Aline xc=\Aline c\gamma\parallel \Aline a\alpha$, which implies $\Aline az=\Aline a\alpha=\Aline z\alpha$. Applying the commutative-plus property to the points $z,d,s,\alpha,\delta,t$, we conclude that $\Aline st\parallel \Aline d\delta$ and hence $\Aline d\delta\parallel \Aline st\parallel\Aline x\gamma=\Aline c\gamma$.
\end{proof}

Theorems~\ref{t:diagonal-trans=>ass-plus} and \ref{t:paraD<=>translation}  imply the following important corollary. 

\begin{corollary}\label{c:Thales=>commutative-plus} Every Thalesian affine space is commutative-plus.
\end{corollary}

\section{Plus loops of Boolean affine spaces}

We recall that a liner $X$ is called \defterm{Boolean} if for any points 
$a,b,c,d\in X$ with $\Aline ab\cap\Aline cd=\varnothing=\Aline bc\cap\Aline ad$, we have $\Aline ac\cap\Aline bd=\varnothing)$. By Theorem~\ref{t:Boolean<=>}, a $3$-ranked liner $X$ is Boolean if and only if every parallelogram in $X$ is Boolean (i.e., has parallel diagonals).

Proposition~\ref{p:Boolean+<=>} implies the following characterization of Boolean affine spaces.

\begin{theorem}\label{t:Boolean<=>Boolean-plus} An affine space $X$ is Boolean if and only if for every ternar $R$ of $X$ its plus loop $(R,+)$ is Boolean.
\end{theorem}

\begin{theorem}\label{t:Boolean=>commutative-plus} Every Boolean affine space is commutative-plus.
\end{theorem}

\begin{proof} Given any distinct points  $a,b,c,\alpha,\beta,\gamma\in X$ with $\Aline ab\parallel\Aline \alpha c\parallel\Aline \beta\gamma\nparallel \Aline \alpha\beta\parallel\Aline a\gamma\parallel\Aline bc$ and $\Aline a\alpha\parallel b\beta$, we need to prove that $\Aline b\beta\parallel\Aline c\gamma$. Since $\Aline ab\parallel\Aline \alpha c\parallel\Aline \beta\gamma\nparallel \Aline \alpha\beta\parallel\Aline a\gamma\parallel\Aline bc$, there exist unique points $x\in\Aline bc\cap\Aline \beta\gamma$, $y\in\Aline \alpha c\cap\Aline a\gamma$, and $z\in\Aline ab\cap\Aline\alpha\beta$. Since $X$ is Boolean, the parallelograms $x\beta zb$ and $y\alpha za$ have parallel diagonals. Therefore, $\Aline xz\parallel \Aline b\beta\parallel \Aline a\alpha\parallel \Aline yz$ and hence $\Aline xz=\Aline yz=\Aline xy$. Since $X$ is Boolean, the parallelogram $x\gamma yc$ has parallel diagonals and hence $\Aline c\gamma\parallel \Aline xy\parallel  \Aline b\beta$.

\begin{picture}(150,100)(-150,-10)
\linethickness{0.8pt}
\put(0,0){\color{red}\line(1,1){80}}
\put(0,80){\color{red}\line(1,-1){80}}
\put(0,0){\color{cyan}\line(0,1){80}}
\put(0,60){\color{teal}\line(1,0){80}}
\put(0,0){\color{teal}\line(1,0){80}}
\put(20,0){\color{cyan}\line(0,1){80}}
\put(80,80){\color{teal}\line(-1,0){80}}
\put(0,60){\color{red}\line(1,1){20}}
\put(80,0){\color{cyan}\line(0,1){80}}
\put(20,0){\color{red}\line(1,1){60}}

\put(0,0){\circle*{3}}
\put(-6,-7){$b$}
\put(80,80){\circle*{3}}
\put(83,78){$\beta$}
\put(20,0){\circle*{3}}
\put(18,-7){$a$}
\put(0,60){\circle*{3}}
\put(-8,57){$c$}
\put(20,80){\circle*{3}}
\put(22,74){$\gamma$}
\put(80,60){\circle*{3}}
\put(83,57){$\alpha$}

\put(0,80){\circle*{3}}
\put(-8,78){$x$}
\put(20,60){\circle*{3}}
\put(14,53){$y$}
\put(80,0){\circle*{3}}
\put(82,-7){$z$}
\end{picture}

\end{proof}

\begin{corollary}\label{c:Boolean-order} If a Boolean affine space $X$ has finite order, then $|X|_2=2^n$ for some $n\ge 2$.
\end{corollary}

\begin{proof} Fix any plane $\Pi$ and an affine base $uow$ in $\Pi$. By Theorems~\ref{t:Boolean<=>Boolean-plus}, \ref{t:Boolean=>commutative-plus} and \ref{t:add-com=>add-ass}, the plus loop $(\Delta,+)$ of the ternar of the based affine plane $(\Pi,uow)$ is a Boolean group. Since $X$ has finite order $|X|_2\ge 2$, the Boolean group $(\Delta,+)$ is finite and hence has cardinality $2^n$ for some $n\in\IN$ (being a finite vector space over the two-element field). Then $3\le |X|_2=|\Delta|=2^n$ and hence $n\ge 2$.
\end{proof}

Corollary~\ref{c:inv-plus=>Boolean-paralelogram} implies the following sufficient condition of the existence of Boolean parallelograms in invertible-plus affine spaces.

\begin{corollary}\label{c:inv-plus=>Boolean-parallelogram2} Every invertible-plus affine space of even order contains a Boolean parallelogram.
\end{corollary}

For inversive-puls affine spaces, Corollary~\ref{c:inv-plus=>Boolean-parallelogram2} can be reversed.

\begin{proposition}\label{p:inversive-plus<=>even<=>Boole} Let $X$ be an inversive-plus affine space of finite order $|X|_2$. The liner $X$ contains a Boolean parallelogram if and only if $|X|_2$ is even.
\end{proposition}

\begin{proof} The ``if'' part follows from Corollary~\ref{c:inv-plus=>Boolean-parallelogram2}. To prove the ``only if'' part, assume that the liner $X$ contains a Boolean parallelogram $uowe$. Consider the affine base $uow$ in the plane $\Pi\defeq\overline{\{u,o,w,e\}}$ and its ternar $\Delta=\Aline oe$. Since $uowe$ is a Boolean parallelogram, $e+e=o$ in the ternar $\Delta$, according to Lemma~\ref{l:Boolean}. Since the affine space $X$ is inversive-plus, the plus loop $(\Delta,+)$ is inversive. Applying Proposition~\ref{p:Ali-Slaney}, we conclude that $|X|_2=|\Delta|$ is even.
\end{proof}

\section{Plus properties of prime affine spaces}

An affine space is defined to be \defterm{prime} if its order is a prime number.

\begin{theorem}\label{t:prime-affine<=>} For every prime affine space $X$, the following conditions are equivalent:
\begin{enumerate}
\item $X$ is Pappian;
\item $X$ is Desarguesian;
\item $X$ is Thalesian;
\item $X$ is translation;
\item $X$ is $\partial$-translation;
\item $X$ is commutative-plus;
\item $X$ is associative-plus;
\item $X$ is inversive-plus.
\end{enumerate}
\end{theorem}

\begin{proof} Since the affine space $X$ has prime (and hence finite) order, every plane in $X$ is finite. Using this fact, we can deduce the equivalence $(1)\Leftrightarrow(2)$ from Theorem~\ref{t:finite-Papp<=>Des}. The equivalence of the conditions (2), (3), (4), (5) was proved in Theorem~\ref{t:partial-translation<=>}. The implication $(3)\Ra(6)$ follows from Corollary~\ref{c:Thales=>commutative-plus}, $(6)\Ra(7)$ follows from Theorem~\ref{t:add-com=>add-ass}, and $(7)\Ra(8)$ is trivial.

It remains to prove that $(8)\Ra(5)$. Assume that a prime affine  pace $X$ is inversive-plus. By Corollary~\ref{c:inversive-plus-prime=>cyclic}, for every ternar $R$ of $X$, the additive loop $(R,+)$ is a cyclic group. By Theorem~\ref{t:add-ass<=>}, $X$ is associative-plus, and by 
Theorem~\ref{t:add-ass=>partialT}, $X$ is $\partial$-translation. 
\end{proof}

\begin{remark}\label{diag:plus-additive} By Theorems~\ref{t:paraD<=>translation}, \ref{t:Boolean=>commutative-plus},  \ref{t:prime-affine<=>}, and Corollary~\ref{c:Thales=>commutative-plus}, for any affine plane we have the implications:
$$
\xymatrix@C=40pt@R=19pt{
\mbox{Pappian}\ar@{=>}[d]\\
\mbox{Desarguesian}\ar@{=>}[d]\ar@/_10pt/_{+\mbox{\em \scriptsize prime}}[u]\\
\mbox{Thalesian}\ar@{=>}[d]\ar@{<=>}[r]\ar@/_10pt/_{+\mbox{\em \scriptsize prime}}[u]&\mbox{translation}\ar@{=>}[ddd]\\
\mbox{Boolean or Thalesian}\ar@{=>}[d]&\\
\mbox{commutative-plus}\ar@{=>}[d]&\\
\mbox{associative-plus}\ar@{=>}[d]\ar_{\phantom{.}+\mbox{\em \scriptsize prime-power}}[r]&\mbox{$\partial$-translation}\ar@/_10pt/_{+\mbox{\em \scriptsize prime}}[uuu]\\
\mbox{inversive-plus}\ar@/_10pt/_{+\mbox{\em \scriptsize prime}}[u]\ar@{=>}[d]\\
\mbox{invertible-plus}&
}
$$
\end{remark}

This diagram motivates the following problems.

\begin{problem} Is every associative-plus affine plane commutative-plus?
\end{problem}

\begin{problem} Is every (finite) invertible-plus affine plane inversive-plus?
\end{problem}

\begin{problem} Is every (finite) inversive-plus affine plane associative-plus?
\end{problem}

\begin{problem} Is every (finite) Boolean affine space Thalesian?
\end{problem}


\begin{problem} Is every  prime (invertible) affine space Desarguesian?
\end{problem}

\chapter{Puls properties of affine spaces}\label{ch:PulsAff}

In this chapter we provide geometric characterizations of some properties of the puls operation in ternars of affine spaces, analogous to the properties of the plus operation, proved in the preceding chapter.

\section{Invertible-puls affine spaces}

\begin{definition}\label{d:a-invertible-puls} An affine liner $X$ is defined to be \index{invertible-puls affine space}\index{affine space!invertible-puls}\defterm{invertible-puls} if for any distinct parallel lines $L,L'$ in $X$ and points $a,b,c\in L$ and $a',b',c'\in L'$,
$$\big(\Aline {a}{a'}\parallel \Aline b{b'}\parallel \Aline c{c'}\;\wedge\;\Aline{a}{b'}\parallel\Aline b{c'}\big)\Ra\Aline {a'}b\parallel \Aline {b'}c.$$
\end{definition}

\begin{picture}(150,50)(-150,-15)
\linethickness{0.75pt}
\put(0,0){\color{teal}\line(1,0){60}}
\put(0,30){\color{teal}\line(1,0){60}}
\put(0,0){\color{blue}\line(1,1){30}}
\put(30,0){\color{blue}\line(1,1){30}}
\put(0,30){\color{red}\line(1,-1){30}}
\put(30,30){\color{red}\line(1,-1){30}}
\put(0,0){\color{cyan}\line(0,1){30}}
\put(30,0){\color{cyan}\line(0,1){30}}
\put(60,0){\color{cyan}\line(0,1){30}}

\put(0,0){\circle*{3}}
\put(-2,-9){$a$}
\put(30,0){\circle*{3}}
\put(26,-9){$b$}
\put(60,0){\circle*{3}}
\put(58,-9){$c$}
\put(0,30){\circle*{3}}
\put(-2,33){$a'$}
\put(30,30){\circle*{3}}
\put(28,33){$b'$}
\put(60,30){\circle*{3}}
\put(58,33){$c'$}
\end{picture}

Let us recall that a loop $(X,\cdot)$ is called \defterm{invertible} if for every $x\in X$ there exists an element $y\in X$ such that $x\cdot y=e=y\cdot x$, where $e$ is the neutral element of the loop $X$.

\begin{theorem}\label{t:invertible-puls<=>} An affine space $X$ is invertible-puls if and only if for every ternar $R$ of $X$ its puls loop $(R,\!\puls\!)$ is invertible.
\end{theorem}

\begin{proof} Assume that the affine plane is invertible-puls. Given any ternar $R$ of $X$, we need to check that its puls loop $(R,\!\puls\!)$ is invertible. Find a plane $\Pi\subseteq X$ and an affine base $uow$ in $\Pi$ whose ternar $\Delta$ is isomorphic to the ternar $R$. It suffices to check that the puls loop $(\Delta,+)$ of the ternar $\Delta$ is invertible. By Proposition~\ref{p:two-sided+}, this will follow as soon as we check that that for every points $x,y\in\Delta$, the parellelity relation $\Aline {oo}{ex}\parallel\Aline{oy}{eo}$ implies $\Aline {oo}{ey}\parallel \Aline{ox}{eo}$.

\begin{picture}(150,120)(-180,-60)
\linethickness{0.8pt}
\put(0,0){\color{teal}\line(1,0){20}}
\put(0,40){\color{teal}\line(1,0){40}}
\put(-40,-40){\color{teal}\line(1,0){60}}
\put(0,-40){\color{cyan}\line(0,1){80}}
\put(20,-40){\color{cyan}\line(0,1){80}}
\put(-40,-40){\line(1,1){80}}
\put(0,0){\color{blue}\line(1,2){20}}
\put(0,-40){\color{blue}\line(1,2){20}}
\put(0,0){\color{red}\line(1,-2){20}}
\put(0,40){\color{red}\line(1,-2){20}}

\put(0,0){\circle*{3}}
\put(-7,-1){$o$}
\put(20,0){\circle*{3}}
\put(23,-2){$eo$}
\put(-40,-40){\circle*{3}}
\put(-43,-48){$y$}
\put(0,-40){\circle*{3}}
\put(-5,-48){$oy$}
\put(20,-40){\circle*{3}}
\put(15,-48){$ey$}
\put(40,40){\circle*{3}}
\put(40,42){$x$}
\put(0,40){\circle*{3}}
\put(-5,43){$ox$}
\put(20,40){\circle*{3}}
\put(15,43){$ex$}
\put(20,20){\circle*{3}}
\put(22,14){$e$}
\end{picture} 

Consider the points $a\defeq oy$, $b\defeq oo$, $c\defeq ox$ on the line $L\defeq \Aline ow$ and the $a'\defeq ey$, $b'\defeq eo$ and $c'\defeq ex$ on the line $\Aline ue$. It is clear that the line $\Aline a{a'}=\Aline {oy}{ey}$, $\Aline b{b'}=\Aline {oo}{eo}$ and $\Aline {c}{c'}=\Aline {ex}{ox}$ are parallel. Since $X$ is inversive-puls, $\Aline a{b'}=\Aline {oy}{eo}\parallel \Aline {oo}{ex}=\Aline {b}{c'}$ implies $\Aline {oy}{oo}=\Aline {a'}b\parallel \Aline {b'}c=\Aline{eo}{ox}$.  
\smallskip
\pagebreak

Now assume that for every ternar $R$ of $X$, its puls loop $(R,\!\puls\!)$ is invertible. Given any distinct parallel lines $L,L'$ in $X$ and points $a,b,c\in L$ and $a',b',c'\in L'$ with $\Aline a{a'}\parallel \Aline b{b'}\parallel \Aline {c}{c'}$ and $\Aline a{b'}\parallel \Aline b{c'}$, we need to check that $\Aline {a'}b\parallel \Aline {b'}c$. If $a=b$, then the parallelity relations $\Aline a{b'}\parallel \Aline b{c'}$ and $\Aline a{a'}\parallel \Aline b{b'}\parallel \Aline c{c'}$ imply $b'=c'$, $a'=b'$ and $c=b$. Then $\Aline {a'}b=\Aline {b'}b=\Aline {b'}c$ and we are done. So, assume that $a\ne b$. 

\begin{picture}(60,80)(-150,-15)
\linethickness{0.75pt}
\put(0,0){\color{teal}\line(1,0){60}}
\put(0,0){\color{red}\line(1,1){60}}
\put(30,0){\color{cyan}\line(0,1){60}}
\put(60,0){\color{cyan}\line(0,1){60}}
\put(30,0){\color{red}\line(1,1){30}}
\put(30,30){\color{teal}\line(1,0){30}}
\put(30,60){\color{teal}\line(1,0){30}}
\put(30,60){\color{violet}\line(1,-1){30}}
\put(30,30){\color{violet}\line(1,-1){30}}

\put(0,0){\circle*{3}}
\put(-8,0){$c''$}
\put(30,0){\circle*{3}}
\put(23,2){$c$}
\put(60,0){\circle*{3}}
\put(63,-1){$c'$}
\put(30,30){\circle*{3}}
\put(23,30){$b$}
\put(60,30){\circle*{3}}
\put(63,28){$b'$}
\put(30,60){\circle*{3}}
\put(22,58){$a$}
\put(60,60){\circle*{3}}
\put(63,58){$a'$}
\end{picture}

Consider the plane $\Pi\defeq\overline{L\cup L'}$ and observe that $uow=b'ba$ is an affine base in $\Pi$. The parallelity relation $\Aline a{a'}\parallel \Aline b{b'}$ implies that the point $a'$ coincides with the diunit $e$ of the base $uow$. Then the line $\Aline b{a'}$ coincides with the diagonal $\Delta$ of the base $uow$. Since $\Aline c{c'}\parallel \Aline b{b'}$ and $\Aline b{b'}\cap \Delta=\{b\}$, there exists a unique point $c''\in \Delta\cap\Aline c{c'}$. Since $\Aline a{b'}\parallel \Aline {b}{c'}$, the definition of the puls operation ensures that $c''\puls a'=b=o$. Since the loop $(\Delta,\!\puls\!)$ is  invertible, $a'\puls c''=o=b'$ and hence $\Aline c{b'}\parallel \Aline {b}{a'}$, by the definition of the puls operation.
\end{proof}

We recall that  a loop $(X,\cdot)$ is \index{power-associative loop}\index{loop!power-associative}\defterm{power-associative} if for every element $x\in X$ there exists a sequence $(x^n)_{n\in\IZ}$ in $X$ such that $x^1=x$ and $x^{n+m}=x^n\cdot x^m$ for every integer numbers $n,m\in\IZ$.

\begin{theorem}\label{t:invertible-puls<=>power-associative} An affine space $X$ is invertible-puls if and only if for every ternar $R$ of $X$, its puls loop $(R,\puls)$ is power-associative.
\end{theorem}

\begin{proof} The ``if'' part follows from the ``if'' part of Theorem~\ref{t:invertible-puls<=>} (because every power-associative loop is invertible). To prove the ``only if'' part, assume that an affine space $X$ is invertible-puls. Given any ternar $R$ of $X$, we need to check that its puls loop $(R,\!\puls\!)$ is power-associative. Find an affine base $uow$ whose ternar $\Delta$ is isomorphic to the ternar $R$. Given any element $x\in\Delta$, we need construct a sequence $(x^n)_{n\in\IZ}$ such that $x^1=x$ and $x^{n+m}=x^n\puls x^m$ for every integer numbers $n,m$. If $x=o$, then the constant sequence $(x^n)_{n\in\IZ}\defeq \{(n,o):n\in\IZ\}$ has the required property. So we assume that $x\ne o$. Put $x^0\defeq o$ and $x^1\defeq x$. By Theorem~\ref{t:invertible-puls<=>}, the loop $(\Delta,\!\puls\!)$ is invertible, so there exists a unique element $x^{-1}\in X$ such that $x^{-1}\puls x=o=x\puls x^{-1}$.  For every $n\in \IN$, let $x^{n+1}\defeq x\puls x^n$ and $x^{-(n+1)}\defeq x^{-1}\puls x^{-n}$. 

\begin{picture}(100,170)(-120,-15)
\linethickness{0.75pt}
\put(60,0){\color{cyan}\line(0,1){140}}
\put(90,0){\color{cyan}\line(0,1){140}}
\put(0,0){\color{red}\line(1,1){140}}
\put(60,60){\color{teal}\line(1,0){30}}
\put(90,0){\color{teal}\line(-1,0){90}}
\put(60,20){\color{violet}\line(3,-2){30}}
\put(90,20){\color{teal}\line(-1,0){70}}
\put(60,40){\color{violet}\line(3,-2){30}}
\put(60,40){\color{blue}\line(3,2){30}}
\put(90,40){\color{teal}\line(-1,0){50}}
\put(60,60){\color{blue}\line(3,2){30}}
\put(60,60){\color{violet}\line(3,-2){30}}
\put(60,60){\color{teal}\line(1,0){30}}
\put(60,80){\color{blue}\line(3,2){30}}
\put(60,80){\color{violet}\line(3,-2){30}}
\put(60,80){\color{teal}\line(1,0){30}}
\put(60,100){\color{blue}\line(3,2){30}}
\put(60,100){\color{teal}\line(1,0){40}}
\put(60,120){\color{blue}\line(3,2){30}}
\put(60,120){\color{teal}\line(1,0){60}}
\put(60,140){\color{teal}\line(1,0){80}}

\put(0,0){\circle*{3}}
\put(-15,0){$x^{-3}$}
\put(20,20){\circle*{3}}
\put(5,20){$x^{-2}$}
\put(40,40){\circle*{3}}
\put(25,40){$x^{-1}$}
\put(60,60){\circle*{3}}
\put(52,60){$o$}
\put(80,80){\circle*{3}}
\put(75,83){$x$}
\put(90,90){\circle*{3}}
\put(92,86){$e$}
\put(100,100){\circle*{3}}
\put(103,96){$x^2$}
\put(120,120){\circle*{3}}
\put(123,116){$x^3$}
\put(140,140){\circle*{3}}
\put(143,136){$x^4$}
\put(60,0){\circle*{3}}
\put(60,20){\circle*{3}}
\put(60,40){\circle*{3}}
\put(60,60){\circle*{3}}
\put(60,80){\circle*{3}}
\put(46,79){$ox$}
\put(60,100){\circle*{3}}
\put(44,99){$ox^2$}
\put(60,120){\circle*{3}}
\put(44,119){$ox^3$}
\put(60,140){\circle*{3}}
\put(44,139){$ox^4$}
\put(90,0){\circle*{3}}
\put(93,-2){$ex^{-3}$}
\put(90,20){\circle*{3}}
\put(93,18){$ex^{-2}$}
\put(90,40){\circle*{3}}
\put(93,38){$ex^{-1}$}
\put(90,60){\circle*{3}}
\put(93,58){$eo$}
\put(90,80){\circle*{3}}
\put(93,77){$ex$}
\put(90,100){\circle*{3}}
\put(90,120){\circle*{3}}
\put(90,140){\circle*{3}}
\end{picture}

We claim that the sequence $(x^n)_{n\in\IZ}$ has the required property: $x^{n+m}=x^n\puls x^m$ for all $n,m\in\IZ$.

\begin{claim}\label{cl:par-1} For every $m\in\IZ$, the lines $\Aline {ox^m}{ex^{m-1}}$ and $\Aline {oo}{ex^{-1}}$ are parallel.
\end{claim}

\begin{proof} For every $m\in\w$, the definition of $x^{-(m+1)}\defeq x^{-1}\puls x^{-m}$ ensures that $\Aline {ox^{-m}}{ex^{-m-1}}\parallel \Aline {oo}{ex^{-1}}$.  Therefore, $\Aline {ox^{m}}{ex^{m-1}}\parallel \Aline {oo}{ex^{-1}}$ for all $m\le 0$.  

For $m=1$, the parallelity of the lines $\Aline {ox^m}{ex^{m-1}}=\Aline{ox}{eo}$ and $\Aline{oo}{ex^{-1}}$ follows from the equality $x^{-1}\puls x=o$ and the definition of the puls operation. 

Assume that for some $m\in\IN$ the lines $\Aline {ox^m}{ex^{m-1}}$ and $\Aline {oo}{ex^{-1}}$ are parallel.
By the definition of the puls operation, the equalities $x^{1+m}=x\puls x^m$ and $x^m=x\puls x^{m-1}$ imply $$\Aline {ox^m}{ex^{1+m}}\parallel \Aline {oo}{ex}\parallel \Aline {ox^{m-1}}{ex^m}.$$ Since $\Aline {ox^m}{ex^{m+1}}\parallel \Aline{ox^{m-1}}{ex^m}$ and $\Aline {ox^{m-1}}{ex^{m-1}}\parallel \Aline {ox^m}{ex^m}\parallel \Aline {ox^{m+1}}{ex^{m+1}}$, the invertibility-puls of $X$ ensures that $\Aline {ox^{m+1}}{ex^m}\parallel \Aline {ox^m}{ex^{m-1}}\parallel \Aline{oo}{ex^{-1}}$.

Therefore, $\Aline {ox^m}{ex^{m-1}}\parallel\Aline {oo}{ex^{-1}}$ for all $m\in\IZ$.   
\end{proof}

By analogy we can prove 

\begin{claim}\label{cl:par+1} For every $m\in\IZ$, the lines $\Aline {ox^m}{ex^{m+1}}$ and $\Aline {oo}{ex}$ are parallel.
\end{claim}

\begin{claim}\label{cl:par:n+m} For every $n\in\w$ and $m\in\IZ$, the lines $\Aline {ox^{m-1}}{ex^{n-1+m}}$ and $\Aline o{x^{m}}{ex^{n+m}}$ are parallel.
\end{claim}

\begin{proof} For $n=0$ the parallelity of the lines $\Aline {ox^{m-1}}{ex^{n-1+m}}$ and $\Aline o{x^{m}}{ex^{n+m}}$ is trivial and for $n=1$ it follows from Claim~\ref{cl:par+1}. 

Assume that for some $n\ge 1$ and every $k\in\w$ with $k\le n$ we know that $\Aline {ox^{m-1}}{ex^{k-1+m}}\parallel \Aline {ox^{m}}{ex^{k+m}}$ for every $m\in\IZ$.  In particular, for every $m\in\IZ$ we have $\Aline {ox^{m-1}}{ex^{n+m-1}}\parallel \Aline {ox^{m}}{ex^{n+m}}\parallel \Aline {ox^{m+1}}{e^{n+m+1}}$ and $\Aline{ox^m}{ex^{n-1+m}}\parallel \Aline {ox^{m+1}}{ex^{n+m}}$.

\begin{picture}(60,90)(-150,-15)
\linethickness{0.75pt}
\put(30,0){\color{orange}\line(1,0){30}}
\put(30,30){\color{red}\line(1,1){30}}
\put(30,0){\color{cyan}\line(0,1){60}}
\put(60,0){\color{cyan}\line(0,1){60}}
\put(30,0){\color{red}\line(1,1){30}}
\put(30,30){\color{orange}\line(1,0){30}}
\put(30,60){\color{orange}\line(1,0){30}}
\put(30,60){\color{violet}\line(1,-1){30}}
\put(30,30){\color{violet}\line(1,-1){30}}

\put(30,0){\circle*{3}}
\put(0,-2){$ox^{m-1}$}
\put(60,0){\circle*{3}}
\put(63,-1){$ex^{n+m-1}$}
\put(30,30){\circle*{3}}
\put(5,28){$ox^m$}
\put(60,30){\circle*{3}}
\put(63,28){$ex^{n+m}$}
\put(30,60){\circle*{3}}
\put(0,58){$ox^{m+1}$}
\put(60,60){\circle*{3}}
\put(63,58){$ex^{n+m+1}$}
\end{picture}

By the invertibility-puls of $X$, those parallelity relations imply $\Aline {ox^{m-1}}{ex^{n+m}}\parallel\Aline {ox^m}{ex^{n+1+m}}$.
\end{proof}

By analogy we can prove 

\begin{claim}\label{cl:par-n+m} For every $n\in\w$ and $m\in\IZ$, the lines $\Aline {ox^m}{ex^{-n+m}}$ and $\Aline o{x^{m+1}}{ex^{-n+m+1}}$ are parallel.
\end{claim}

Claims~\ref{cl:par:n+m} and \ref{cl:par-n+m} imply that for every $n,m\in\IZ$ the lines $\Aline{ox^m}{e^{n+m}}$ and $\Aline {oo}{ex^n}$ are parallel. By the definition of the puls operation, $\Aline{ox^m}{e^{n+m}}\parallel\Aline {oo}{ex^n}$ imply the desired equality $x^n\puls x^m=x^{n+m}$, witnessing that the puls loop $(\Delta,\!\puls\!)$ is power-associative and so is its isomorphic copy $(R,\!\puls\!)$.
 \end{proof}
 
\begin{exercise}Show that the Moulton plane is not invertible-puls.

{\em Hint:} Consider the points 
$$a=(-2,0),\;\; b=(0,0),\;\; c=(1,0),\;\; a'=(2,2),\;\; b'=(0,2),\;\; c'=(-2,2)$$ and observe that
$\Aline a{b'}\parallel \Aline{a'}b$, $\Aline b{c'}\parallel \Aline {b'}c$ and $\Aline a{c'}\parallel \Aline b{b'}$, but $\Aline {a'}c\nparallel\Aline b{b'}$.

\begin{picture}(200,75)(-200,-15)
\linethickness{0.75pt}
\put(-40,0){\color{teal}\line(1,0){60}}
\put(-40,40){\color{teal}\line(1,0){80}}
\put(-40,0){\color{cyan}\line(1,1){40}}
\put(0,0){\color{cyan}\line(1,1){40}}
\put(-40,40){\color{violet}\line(1,-1){40}}
\put(0,40){\color{violet}\line(1,-2){20}}
\put(-40,0){\color{red}\line(0,1){40}}
\put(0,0){\color{red}\line(0,1){40}}
\put(20,0){\color{orange}\line(1,2){20}}

\put(-40,0){\circle*{3}}
\put(-43,-8){$a$}
\put(0,0){\circle*{3}}
\put(-3,-9){$b$}
\put(20,0){\circle*{3}}
\put(18,-8){$c$}
\put(-40,40){\circle*{3}}
\put(-42,43){$c'$}
\put(0,40){\circle*{3}}
\put(-2,43){$b'$}
\put(40,40){\circle*{3}}
\put(38,43){$a'$}
\end{picture}

\end{exercise}

\section{Inversive-puls affine spaces}

\begin{definition}\label{d:left-inversive-puls} An affine space $X$ is defined to be \index{left-inversive-puls affine space}\index{affine space!left-inversive-puls}\defterm{left-inversive-puls} if for all distinct parallel lines $L,L'\subseteq X$ and points $a,b,\alpha,\beta\in L$ and $a',b',\alpha',\beta'\in L'$,
$$\big(\Aline a{a'}\parallel \Aline b{b'}\parallel \Aline \alpha{\alpha'}\parallel \Aline \beta{\beta'}\;\wedge\;\Aline a{b'}\parallel \Aline \alpha{\beta'}\big)\Ra\Aline {a'}b\parallel \Aline {\alpha'}\beta.$$
\end{definition}

\begin{picture}(100,60)(-150,-15)

\linethickness{0.75pt}
\put(0,0){\color{teal}\line(1,0){80}}
\put(0,30){\color{teal}\line(1,0){80}}
\put(0,0){\color{cyan}\line(0,1){30}}
\put(30,0){\color{cyan}\line(0,1){30}}
\put(50,0){\color{cyan}\line(0,1){30}}
\put(80,0){\color{cyan}\line(0,1){30}}
\put(0,0){\color{blue}\line(1,1){30}}
\put(50,0){\color{blue}\line(1,1){30}}
\put(0,30){\color{red}\line(1,-1){30}}
\put(50,30){\color{red}\line(1,-1){30}}

\put(0,0){\circle*{3}}
\put(-3,-9){$a$}
\put(30,0){\circle*{3}}
\put(26,-10){$b$}
\put(50,0){\circle*{3}}
\put(47,-9){$\alpha$}
\put(80,0){\circle*{3}}
\put(79,-9){$\beta$}
\put(0,30){\circle*{3}}
\put(-3,33){$a'$}
\put(30,30){\circle*{3}}
\put(28,33){$b'$}
\put(50,30){\circle*{3}}
\put(47,33){$\alpha'$}
\put(80,30){\circle*{3}}
\put(77,33){$\beta'$}

\end{picture}

\begin{theorem}\label{t:left-inversive-puls<=>} For an affine space $X$ the following statements are equivalent:
\begin{enumerate}
\item the affine space $X$ is left-inversive-puls;
\item for every ternar $R$ of $X$, its puls loop $(R,\!\puls\!)$ is left-inversive;
\item for every ternar $R$ of $X$, its puls loop $(R,\!\puls\!)$ is left-Bol.
\end{enumerate}
\end{theorem}

\begin{proof} The equivalence $(1)\Leftrightarrow(2)$ follows from Proposition~\ref{p:lip+<=>}, and the implication $(3)\Ra(2)$ follows from Proposition~\ref{p:left-Bol=>left-inversive}.
\smallskip

$(2)\Ra(3)$ Assume that the affine plane is  left-inversive-puls. Given any plane $\Pi\subseteq X$, an affine base $uow$ in $\Pi$ and its diagonal $\Delta$, we need to show that the puls loop $(\Delta,\!\puls\!)$ is left-Bol. Take any points $x,y,z\in\Delta$. We have to show that $(x\puls (y\puls x))\puls z=x\puls (y\puls (x\puls z))$. Consider the points $a\defeq y\puls x$, $b\defeq x\puls a=x\puls(y\puls x)$, $c\defeq x\puls z$, $d\defeq y\puls c=y\puls(x\puls z)$, and $s\defeq x\puls d=x\puls(y\puls(x\puls z))$. We have to prove that $s=b\puls z$.

\begin{picture}(200,183)(-150,-45)
\linethickness{0.7pt}
\put(0,-30){\color{cyan}\line(0,1){150}}
\put(-30,-30){\line(1,1){150}}
\put(40,-30){\color{cyan}\line(0,1){150}}
\put(0,0){\color{blue}\line(4,5){40}}
\put(0,0){\color{violet}\line(4,-3){40}}
\put(0,-30){\color{blue}\line(4,5){40}}
\put(-30,-30){\color{teal}\line(1,0){70}}
\put(0,20){\color{teal}\line(1,0){40}}
\put(0,20){\color{violet}\line(4,-3){40}}
\put(-10,-10){\color{teal}\line(1,0){50}}
\put(0,0){\color{red}\line(4,-1){40}}
\put(0,100){\color{teal}\line(1,0){100}}
\put(0,70){\color{teal}\line(1,0){70}}
\put(0,100){\color{violet}\line(4,-3){40}}
\put(0,70){\color{blue}\line(4,5){40}}
\put(0,120){\color{teal}\line(1,0){120}}
\put(0,120){\color{violet}\line(4,-3){40}}
\put(0,100){\color{red}\line(4,-1){40}}
\put(40,90){\color{teal}\line(1,0){50}}
\put(40,50){\color{teal}\line(1,0){10}}
\put(0,20){\color{orange}\line(4,-5){40}}
\put(0,120){\color{orange}\line(4,-5){40}}

\put(-30,-30){\circle*{3}}
\put(-37,-33){$x$}
\put(-10,-10){\circle*{3}}
\put(-17,-12){$b$}
\put(0,0){\circle*{3}}
\put(-7,-1){$o$}
\put(20,20){\circle*{3}}
\put(18,13){$a$}
\put(40,40){\circle*{3}}
\put(42,36){$e$}
\put(50,50){\circle*{3}}
\put(52,45){$y$}
\put(70,70){\circle*{3}}
\put(72,65){$c$}
\put(90,90){\circle*{3}}
\put(92,85){$s$}
\put(100,100){\circle*{3}}
\put(102,95){$z$}
\put(120,120){\circle*{3}}
\put(123,118){$d$}

\put(0,-30){\circle*{3}}
\put(-5,-38){$ox$}
\put(40,-30){\circle*{3}}
\put(43,-33){$ex$}
\put(40,-10){\circle*{3}}
\put(43,-12){$eb$}
\put(0,20){\circle*{3}}
\put(-13,17){$oa$}
\put(40,20){\circle*{3}}
\put(43,18){$ea$}
\put(40,50){\circle*{3}}
\put(28,48){$ey$}
\put(0,70){\circle*{3}}
\put(-13,67){$oc$}
\put(40,70){\circle*{3}}
\put(42,72){$ec$}
\put(40,90){\circle*{3}}
\put(42,92){$es$}
\put(0,100){\circle*{3}}
\put(-13,98){$oz$}
\put(0,120){\circle*{3}}
\put(-5,123){$od$}
\put(40,120){\circle*{3}}
\put(36,123){$ed$}

\end{picture}

 It follows from $a=y\puls x$ and $d=y\puls c$ that $\Aline {ox}{ea}\parallel\Aline{oo}{ey}\parallel \Aline {oc}{ed}$. By the left-invertibility-puls, $\Aline {ox}{ea}\parallel\Aline {oc}{ed}$ implies $\Aline {oa}{ex}\parallel\Aline {od}{ec}$. It follows from $b=x\puls a$ , $c=x\puls z$, $s=x\puls d$ that the line $\Aline{oo}{ex}$ is parallel to the lines $\Aline {oa}{eb}$, $\Aline {oz}{ec}$, $\Aline{od}{es}$. By the left-inversivity-puls, the parallelity relations $\Aline{oo}{ex}\parallel\Aline {oa}{eb}\parallel \Aline {oz}{ec}\parallel\Aline{od}{es}$ and $\Aline{oa}{ex}\parallel \Aline{od}{ec}$ imply $\Aline {oo}{eb}\parallel \Aline{oz}{es}$. By the definition of the puls operation, the latter paralleity relation implies $s=b\puls z$ and hence $$
x\puls(y\puls (x\puls z))=s=b\puls z=(x\puls(y\puls x))\puls z,$$
witnessing that the loop $(\Delta,\!\puls\!)$ is left-Bol.
\end{proof}

Theorem~\ref{t:left-inversive-puls<=>} and Corollary~\ref{c:left-Bol-prime=>cyclic} imply 

\begin{corollary}\label{c:left-inversive-puls-prime=>cyclic} If a left-inversive-puls affine space $X$ has prime order $p$, then for every ternar $R$ of $X$, the additive loop $(R,\!\puls\!)$ is a cyclic group of order $p$.
\end{corollary}

\begin{problem} Can Corollary~\textup{\ref{c:inversive-plus-prime=>cyclic}} be generalized to invertible-puls (instead of left-inversive-puls) affine spaces?
\end{problem}

\begin{definition}\label{d:right-inversive-puls} An affine space $X$ is defined to be \index{right-inversive-puls affine space}\index{affine space!right-inversive-puls}\defterm{right-inversive-puls} if for all distinct parallel lines $L,L'\subseteq X$ and points $x,y,z\in L$ and $a,b,a',b'\in L'$
$$\big(\Aline xa\parallel \Aline yb\;\wedge\;\Aline ya\parallel \Aline zb\;\wedge\; \Aline x{a'}\parallel \Aline y{b'}\big)\Ra\Aline y{a'}\parallel \Aline z{b'}.$$
\end{definition}

\begin{picture}(100,60)(-150,-15)

\linethickness{0.75pt}
\put(0,0){\color{teal}\line(1,0){60}}
\put(0,30){\color{teal}\line(1,0){75}}
\put(0,0){\color{cyan}\line(0,1){30}}
\put(30,0){\color{cyan}\line(0,1){30}}
\put(30,0){\color{orange}\line(-1,1){30}}
\put(60,0){\color{orange}\line(-1,1){30}}
\put(0,0){\color{blue}\line(3,2){45}}
\put(30,0){\color{blue}\line(3,2){45}}
\put(30,0){\color{red}\line(1,2){15}}
\put(60,0){\color{red}\line(1,2){15}}

\put(0,0){\circle*{3}}
\put(-3,-8){$x$}
\put(30,0){\circle*{3}}
\put(27,-8){$y$}
\put(60,0){\circle*{3}}
\put(57,-8){$z$}
\put(0,30){\circle*{3}}
\put(-3,33){$a$}
\put(30,30){\circle*{3}}
\put(28,33){$b$}
\put(45,30){\circle*{3}}
\put(42,33){$a'$}
\put(75,30){\circle*{3}}
\put(72,33){$b'$}

\end{picture}

\begin{theorem}\label{t:right-inversive-puls<=>} An affine space  $X$ is right-inversive-puls if and only if for every ternar $R$ of $X$, its puls loop $(R,\!\puls\!)$ is right-inversive.
\end{theorem}

\begin{proof} Assume that an affine space $X$ is right-inversive-puls. Given any ternar $R$ of $X$, we need to show that its puls loop $(R,\!\puls\!)$ is right-inversive. Find a plane $\Pi\subseteq X$ and an affine base $uow$ in $\Pi$ whose ternar $\Delta$ is isomorphic to the ternar $R$. It suffices to prove that the loop $(\Delta,\!\puls\!)$ is right-inversive. The right-inversivity of the loop $(\Delta,\!\puls\!)$ will follow from Proposition~\ref{p:rip+<=>} as soon as we show that for every points $x,y,z,s\in\Delta$, the parallelity relations $\Aline {oo}{ex}\parallel \Aline{oy}{eo}$ and $\Aline{oo}{ez}\parallel \Aline{ox}{es}$ imply $\Aline {oo}{es}\parallel \Aline {oy}{ez}$.

\begin{picture}(150,140)(-170,-72)
\linethickness{0.8pt}
\put(-60,-60){\line(1,1){120}}
\put(0,-60){\color{cyan}\line(0,1){120}}
\put(40,-40){\color{cyan}\line(0,1){100}}
\put(-60,-60){\color{teal}\line(1,0){60}}
\put(-40,-40){\color{teal}\line(1,0){80}}
\put(0,0){\color{teal}\line(1,0){40}}
\put(0,60){\color{teal}\line(1,0){60}}
\put(0,-60){\color{red}\line(2,1){40}}
\put(0,-60){\color{blue}\line(2,3){40}}
\put(0,0){\color{red}\line(2,1){40}}
\put(0,0){\color{blue}\line(2,3){40}}
\put(0,0){\color{violet}\line(1,-1){40}}
\put(0,60){\color{violet}\line(1,-1){40}}
\put(20,20){\color{teal}\line(1,0){20}}

\put(-60,-60){\circle*{3}}
\put(-65,-68){$y$}
\put(0,-60){\circle*{3}}
\put(-5,-68){$oy$}
\put(-40,-40){\circle*{3}}
\put(-48,-41){$z$}
\put(40,-40){\circle*{3}}
\put(43,-42){$ez$}
\put(0,0){\circle*{3}}
\put(-8,-2){$o$}
\put(40,0){\circle*{3}}
\put(43,-2){$eo=u$}
\put(20,20){\circle*{3}}
\put(18,13){$s$}
\put(40,20){\circle*{3}}
\put(43,18){$es$}
\put(40,40){\circle*{3}}
\put(43,37){$e$}
\put(0,60){\circle*{3}}
\put(-5,63){$ox$}
\put(40,60){\circle*{3}}
\put(35,63){$ex$}
\put(60,60){\circle*{3}}
\put(58,63){$x$}
\end{picture}

Consider the parallel lines $L\defeq \Aline ow$ and $L'\defeq \Aline ue$, where $e$ is the diunit of the affine base $uow$. Also consider the points $x'\defeq ox$, $y'\defeq oo$, $z'\defeq oy$ on the line $L$, and the points $a\defeq ex$, $b\defeq eo$, $a'\defeq es$ and $b'\defeq ez$ on the line $L'$. Since $\Aline {x'}a=\Aline {ox}{ex}\parallel \Aline{oo}{eo}=\Aline {y'}b$, $\Aline {y'}a=\Aline{oo}{ex}\parallel \Aline{z'}b=\Aline{oy}{eo}$ and $\Aline {x'}{a'}=\Aline {ox}{es}\parallel \Aline{oo}{ez}$, we can apply the right-inversivity-puls of $X$ and conclude that $\Aline{oo}{es}=\Aline{y'}{a'}\parallel \Aline{z'}{b'}=\Aline{oy}{ez}$.
\smallskip

Now assuming that for every ternar $R$ of $X$ its puls loop is right-inversive, we shall prove that the affine plane $X$ is right-inversive-puls. Given any disjoint parallel lines $L,L'$ in $X$ and points $\alpha,\beta,\gamma\in L$, $a,b,a',b'\in L'$ with $\Aline \alpha a\parallel \Aline \beta b$, $\Aline \beta a\parallel \Aline \gamma b$ and $\Aline \alpha{a'}\parallel \Aline 
\beta{b'}$, we need to show that $\Aline \beta{a'}\parallel \Aline \gamma{b'}$. If $\alpha=\beta$, then the parallelity relation $\Aline \alpha a\parallel\Aline \beta b$ implies $\Aline \alpha a=\Aline \beta b$ and hence $a\in L'\cap\Aline \alpha a=L'\cap\Aline \beta b=\{b\}$. The equality $a=b$ and the parallelity relation $\Aline \beta a\parallel \Aline \gamma b$ imply $\Aline \beta a=\Aline \gamma b$ and hence $\beta\in L\cap\Aline \beta a=L\cap\Aline \gamma b=\{\gamma\}$. Therefore, $\alpha=\beta=\gamma$. Then $\Aline \beta{a'}=\Aline \alpha{a'}\parallel \Aline \beta{b'}=\Aline \gamma{b'}$. So, we assume that $\alpha\ne \beta$.

\begin{picture}(150,100)(-170,-50)
\linethickness{0.8pt}
\put(-40,-40){\color{blue}\line(1,1){80}}
\put(0,-40){\color{cyan}\line(0,1){80}}
\put(40,-20){\color{cyan}\line(0,1){60}}
\put(-40,-40){\color{teal}\line(1,0){40}}
\put(-20,-20){\color{teal}\line(1,0){60}}

\put(0,0){\color{teal}\line(1,0){40}}
\put(0,40){\color{teal}\line(1,0){40}}
\put(0,-40){\color{red}\line(2,1){40}}
\put(0,-40){\color{blue}\line(1,1){40}}
\put(0,0){\color{red}\line(2,1){40}}
\put(0,0){\color{violet}\line(2,-1){40}}
\put(0,40){\color{violet}\line(2,-1){40}}
\put(20,20){\color{teal}\line(1,0){20}}

\put(-40,-40){\circle*{3}}
\put(-48,-41){$y$}
\put(0,-40){\circle*{3}}
\put(-3,-48){$\gamma$}
\put(-20,-20){\circle*{3}}
\put(-28,-21){$z$}
\put(40,-20){\circle*{3}}
\put(43,-22){$b'$}
\put(0,0){\circle*{3}}
\put(-27,-1){$o=\beta$}
\put(40,0){\circle*{3}}
\put(43,-2){$b=u$}
\put(20,20){\circle*{3}}
\put(18,13){$s$}
\put(40,20){\circle*{3}}
\put(43,18){$a'$}
\put(40,40){\circle*{3}}
\put(43,37){$a=e=x$}
\put(0,40){\circle*{3}}
\put(-30,38){$w=\alpha$}
\end{picture}

Consider the points $u\defeq b$, $o\defeq \beta$, $w\defeq \alpha$ and observe that $uow$ is an affine base in the plane $\Pi\defeq\overline{L\cup L'}$ whose diunit $e$ coincides with the point $a$. Let $\Delta\defeq\Aline oe$ be the diagonal and $\boldsymbol h\defeq (\Aline ou)_\parallel$ be the horizontal direction of the affine base $uow$. Consider the unique points $x\defeq a=e$, $y\in\Delta\cap\Aline\gamma{\boldsymbol h}$, $z\in \Delta\cap\Aline {b'}{\boldsymbol h}$, $s\in \Delta\cap\Aline {a'}{\boldsymbol h}$. Observe that $\Aline{oo}{ex}=\Aline \beta a\parallel \Aline \gamma b=\Aline{oy}{eo}$ and $\Aline {ox}{es}=\Aline{\alpha}{a'}\parallel \Aline{\beta}{b'}=\Aline {oo}{ez}$. Since the loop $(\Delta,\!\puls\!)$ is right-inversive, Proposition~\ref{p:rip+<=>} ensures that $\Aline\beta {a'}=\Aline {oo}{es}\parallel \Aline {oy}{ez}=\Aline \gamma{b'}$. 
\end{proof}


\begin{definition}\label{d:inversive-puls} An affine space $X$ is defined to be \index{inversive-puls}\defterm{inversive-puls} if $X$ is left-inversive-puls and right-inversive-puls.
\end{definition}

\begin{theorem}\label{t:inversive-puls} For an affine space $X$ the following statements are equivalent:
\begin{enumerate}
\item the affine space $X$ is inversive-puls;
\item for every ternar $R$ of $X$, its puls loop $(R,\!\puls\!)$ is inversive;
\item for every ternar $R$ of $X$, its plus loop $(R,\!\puls\!)$ is Bol;
\item for every ternar $R$ of $X$, its plus loop $(R,\!\puls\!)$ is Moufang.
\end{enumerate}
\end{theorem}

\begin{proof} The equivalence $(1)\Leftrightarrow(2)$ follows from Theorems~\ref{t:left-inversive-puls<=>} and \ref{t:right-inversive-puls<=>}. The equivalence $(3)\Leftrightarrow(4)$ follows from Theorem~\ref{t:Bol<=>Moufang}. The implication $(3)\Ra(2)$ follows from Corollary~\ref{c:Bol=>inversive}. It remain to prove that $(1)\Ra(3)$ Assume that the affine space $X$ is inversive-plus and take any ternar $R$ of $X$. By Theorems~\ref{t:left-inversive-puls<=>} and \ref{t:right-inversive-puls<=>}, the loop $(R,\!\puls\!)$ is left-Bol and right-inversive. By Theorem~\ref{t:Bol<=>}, the loop $(R,\!\puls\!)$ is Bol.
\smallskip

\end{proof}

\section{Associative-puls affine spaces}

\begin{definition}\label{d:a-associative-puls} An affine liner $X$ is defined to be \index{associative-puls affine space}\index{affine space!associative-puls}\defterm{associative-puls} if for every disjoint parallel lines $L,L'$ in $X$ and points $a,b,\alpha,\beta\in L$ and $a',b',\alpha',\beta'\in L'$, $$\big(\Aline a{a'}\parallel\Aline \alpha{\alpha'}\;\wedge\;\Aline a{b'}\parallel\Aline\alpha{\beta'}\;\wedge\;\Aline {a'}{b}\parallel\Aline{\alpha'}\beta\big)\;\Ra\;\Aline b{b'}\parallel\Aline \beta{\beta'}.$$
\end{definition}

\begin{picture}(100,50)(-140,-10)

\linethickness{0.7pt}
\put(10,0){\color{teal}\line(1,0){80}}
\put(0,30){\color{teal}\line(1,0){100}}
\put(10,0){\color{cyan}\line(-1,3){10}}
\put(10,0){\color{violet}\line(1,1){30}}
\put(30,0){\color{blue}\line(-1,1){30}}
\put(30,0){\color{red}\line(1,3){10}}

\put(70,0){\color{cyan}\line(-1,3){10}}
\put(70,0){\color{violet}\line(1,1){30}}
\put(90,0){\color{blue}\line(-1,1){30}}
\put(90,0){\color{red}\line(1,3){10}}

\put(10,0){\circle*{3}}
\put(8,-9){$a$}
\put(30,0){\circle*{3}}
\put(27,-9){$b$}
\put(0,30){\circle*{3}}
\put(-3,33){$a'$}
\put(40,30){\circle*{3}}
\put(38,33){$b'$}

\put(70,0){\circle*{3}}
\put(68,-9){$\alpha$}
\put(90,0){\circle*{3}}
\put(87,-9){$\beta$}
\put(60,30){\circle*{3}}
\put(57,33){$\alpha'$}
\put(100,30){\circle*{3}}
\put(98,33){$\beta'$}
\end{picture}


\begin{theorem}\label{t:associative-puls<=>} An affine space $X$ is associative-plus if and only if for every ternar $R$ of $X$ its puls loop $(R,\!\puls\!)$ is associative.
\end{theorem}

\begin{proof} Assume that an affine space $X$ is associative-puls. Given any plane $\Pi\subseteq X$ and an affine base $uow$ in $\Pi$, we should prove that the puls loop $(\Delta,\!\puls\!)$ is associative. By Proposition~\ref{p:add-ass<=>+}, it suffices to show that for every points $x,y,p,q,r\in \Delta$, the parallelity relations $\Aline {ox}{ey}\parallel \Aline {op}{eq}$ and $\Aline{oo}{ex}\parallel \Aline{or}{ep}$ imply $\Aline {oo}{ey}\parallel \Aline{or}{eq}$.

Consider the points $a\defeq ex$, $b\defeq ey$, $\alpha\defeq ep$, $\beta\defeq eq$ on the line $L\defeq \Aline ue$ and the points $a'\defeq ox$, $b'\defeq oo$, $\alpha'\defeq op$ and $\beta'\defeq or$ on the line $L'\defeq \Aline ow$. Since $\Aline a{a'}=\Aline {ex}{ox}\parallel \Aline {ep}{op}=\Aline \alpha{\alpha'}$, $\Aline a{b'}=\Aline{ex}{oo}\parallel \Aline {ep}{or}=\Aline\alpha{\beta'}$ and $\Aline {a'}b=\Aline {ox}{ey}\parallel \Aline {op}{eq}=\Aline {\alpha'}{\beta}$, we can apply the associativity-puls of $X$ and conclude that $\Aline {ey}{oo}=\Aline b{b'}\parallel \Aline\beta{\beta'}=\Aline{eq}{or}$. This completes the proof of the ``only if'' part.

\begin{picture}(150,170)(-160,-55)
\linethickness{0.75pt}
\put(-40,-40){\line(1,1){140}}
\put(-40,-40){\color{teal}\line(1,0){80}}
\put(-20,-20){\color{teal}\line(1,0){60}}
\put(60,60){\color{teal}\line(-1,0){60}}
\put(80,80){\color{teal}\line(-1,0){40}}
\put(100,100){\color{teal}\line(-1,0){100}}
\put(0,-40){\color{cyan}\line(0,1){140}}
\put(40,-40){\color{cyan}\line(0,1){120}}
\put(0,-40){\color{blue}\line(2,1){40}}
\put(0,60){\color{blue}\line(2,1){40}}
\put(0,0){\color{violet}\line(1,-1){40}}
\put(0,100){\color{violet}\line(1,-1){40}}
\put(0,0){\color{red}\line(2,-1){40}}
\put(0,100){\color{red}\line(2,-1){40}}

\put(-40,-40){\circle*{3}}
\put(-43,-48){$x$}
\put(0,-40){\circle*{3}}
\put(-5,-48){$ox$}
\put(1,-36){$a'$}
\put(40,-40){\circle*{3}}
\put(30,-48){$ex$}
\put(43,-40){$a$}
\put(-20,-20){\circle*{3}}
\put(-28,-20){$y$}
\put(40,-20){\circle*{3}}
\put(43,-25){$ey$}
\put(43,-18){$b$}
\put(0,0){\circle*{3}}
\put(-7,0){$o$}
\put(1,6){$b'$}
\put(40,40){\circle*{3}}
\put(33,40){$e$}
\put(60,60){\circle*{3}}
\put(63,55){$p$}
\put(0,60){\circle*{3}}
\put(-13,58){$op$}
\put(1,64){$\alpha'$}
\put(40,60){\circle*{3}}
\put(42,62){$ep$}
\put(33,52){$\alpha$}
\put(80,80){\circle*{3}}
\put(83,75){$q$}
\put(40,80){\circle*{3}}
\put(40,83){$eq$}
\put(33,70){$\beta$}
\put(100,100){\circle*{3}}
\put(103,95){$r$}
\put(0,100){\circle*{3}}
\put(-12,95){$or$}
\put(-3,103){$\beta'$}
\end{picture}

Now assume that for every ternar $R$ of the affine space $X$, the puls loop $(R,\!\puls\!)$ is associative. To prove that $X$ is associative-puls, take any distinct parallel lines $L,L'\in X$ and points $a,b,\alpha,\beta\in L$ and $a',b',\alpha',\beta'\in L'$ with  $\Aline a{b'}\parallel \Aline\alpha{\beta'}$, $\Aline{a'}b\parallel \Aline{\alpha'}\beta$ and  $\Aline a{a'}\parallel \Aline \alpha{\alpha'}$. We have to show that $\Aline b{b'}=\Aline\beta{\beta'}$. In the affine plane $\Pi\defeq\overline{L\cup L'}$, chose an affine base $uow$ such that $o=b'$, $\Aline ow=L'$ and $\Aline ou\parallel \Aline a{a'}\parallel \Aline\alpha{\alpha'}$. Let $e$ be the diunit of the affine base $uow$, $\Delta\defeq\Aline oe$ be the diagonal, and $\boldsymbol h\defeq(\Aline ou)_\parallel$ be the horizontal direction of the affine base $uow$.

Consider the unique points $x\in\Delta\cap\Aline a{a'}$, $y\in\Delta\cap\Aline b{\boldsymbol h}$, $p\in\Delta\cap\Aline\alpha{\boldsymbol h}$, $q\in\Delta\cap\Aline\beta{\boldsymbol h}$, $r\in\Delta\cap\Aline{\beta'}{\boldsymbol h}$. 
Observe that $\Aline{ox}{ey}=\Aline{a'}b\parallel \Aline{\alpha'}\beta=\Aline{op}{eq}$ and $\Aline {oo}{ex}=\Aline{b'}a\parallel \Aline{\beta'}\alpha=\Aline{or}{ep}$. Since the loop $(\Delta,\!\puls\!)$ is associative, we can apply Proposition~\ref{p:add-ass<=>+} and conclude that $\Aline b{b'}=\Aline {ey}{oo}\parallel \Aline {eq}{or}=\Aline{\beta}{\beta'}$, witnessing that the affine space $X$ is associative-puls.
 \end{proof}

Let us recall that an affine space $X$ is called \index{$\partial$-translation}\index{affine space!$\partial$-translation}\defterm{$\partial$-translation} if for every direction $\boldsymbol \delta\in \partial X$ there exists a non-identity translation $T:X\to X$ such that $\boldsymbol\delta=\{\Aline xy:xy\in T\}$.

\begin{theorem}\label{t:ass-puls=>partialT} Let $X$ be an affine space of  order $|X|_2=p^n$ for some prime number $p$ and some $n\in\IN$. If $X$ is associative-puls, then $X$ is $\partial$-translation.
\end{theorem}

\begin{proof} Assume that $X$ is associative-plus and $|X|_2$ is a prime power. If $\|X\|\ge 4$, then $X$ is Desarguesian, Thalesian, translation, and $\partial$-translation, by Corollary~\ref{c:affine-Desarguesian} and  Theorems~\ref{t:ADA=>AMA} and \ref{t:paraD<=>translation}.  So, assume that $\|X\|=3$, which means that $X$ is a Playfair plane. To show that $X$ is $\partial$-translation, we shall apply Theorem~\ref{t:partial-translation<=+}. Given any distinct directions $\boldsymbol h,\boldsymbol v\in\partial X$ and any distinct lines $L,\Lambda\in \boldsymbol v$, we should show that the set $\Sym_X^\#[L;{\boldsymbol h},\Lambda]$ is a subgroup of the group $\Sym_X^\#(L)$. Choose any distinct points $o,w\in \Lambda$ and consider the unique points $u\defeq L\cap \Aline o{\boldsymbol h}$ and $e\defeq L\cap\Aline w{\boldsymbol h}$. Then $uow$ is an affine base for the affine plane $X$ such that $\Aline ou\in {\boldsymbol h}$ and  $\Aline ow=\Lambda\in {\boldsymbol v}$. Consider the ternar $\Delta=\Aline oe$ of the  based affine plane $(X,uow)$. Since $X$ is associative-puls, the puls loop $(\Delta,\!\puls\!)$ is a group. 

Consider the line projection $P:\Delta\to L$ assigning to every point $x\in\Delta$ the unique point $P(x)\in \Aline x{\boldsymbol h}\cap L$. Endow the set $L$ with a unique binary operation $*:L\times L\to L$ such that the function $P$ is an isomorphism of the magmas $(L,*)$ and $(\Delta,\puls)$. Since the magma $(\Delta,\puls)$ is a group, so is its isomorphic copy $(L,*)$. 

Consider the function $T:L\to \Sym^\#_X[L;\boldsymbol h,\Lambda]\subseteq\Sym^\#(L)$ assigning to every point $a\in L$ the line translation $T_a\defeq {\boldsymbol a}_{L,\Lambda}{\boldsymbol h}_{\Lambda,L}:L\to L$ where ${\boldsymbol a}\defeq (\Aline oa)_\parallel$.  It is easy to see that the function $T$ is injective and $T[L]=\Sym^\#_X[L;\boldsymbol h,\Lambda]$.

\begin{claim}\label{cl:Ta-shift} $T_a(x)=a*x$ for all points $a,x\in L$.
\end{claim}

\begin{proof} Consider the points $y\defeq a*x$, $x'\defeq P^{-1}(x)$, $a'\defeq P^{-1}(a)$, and $y'\defeq P^{-1}(y)$. Taking into account that $P$ is an isomorphism of the groups $(L,*)$ and $(\Delta,\puls)$, we conclude that $y'=a'\puls x'$. Consider the point $x''\in {\boldsymbol h}_{\Lambda,L}(x)\in\Lambda$. By the definition of the puls operation, the equality $y'=a'\puls x'$ implies $\Aline{x''}y\parallel \Aline oa$.
Then 
$$a*x=y={\boldsymbol a}_{\Lambda,L}(x'')={\boldsymbol a}_{L,\Lambda}{\boldsymbol h}_{\Lambda,L}(x)=T_a(x).$$
\end{proof} 

The associativity of the binary operation $*$ and Claim~\ref{cl:Ta-shift} ensure that 
$$T_bT_a(x)=T_b(a*x)=b*(a*x))=(b*a)*x=T_{b*a}(x)$$
for every $a,b,x\in L$. Therefore, $T_bT_a=T_{b*a}$ and $T:L\to \Sym^\#_X(L)$ is a homomorphism of the groups $(L,*)$ into the group $\Sym^\#_X(L)$. Then the image $\Sym^\#_X[L;\boldsymbol h,\Lambda]$ of the group $(L,*)$ under the homomorphism $T$ is a subgroup of the group $\Sym^\#_X(L)$.

Now we can apply Theorem~\ref{t:partial-translation<=+} and conclude that the affine plane $X$ is $\partial$-translation.
\end{proof}




\section{Commutative-puls affine spaces}

\begin{definition}\label{d:a-commutative-puls} An affine liner $X$ is defined to be \index{commutative-puls affine space}\index{affine space!commutative-puls}\defterm{commutative-puls} if for any distinct parallel lines $L,L'\subseteq X$ and points $a,b,c\in L$ and $a',b',c'\in L'$,$$\big(\Aline a{b'}\parallel\Aline {a'}b\;\wedge \;\Aline b{c'}\parallel \Aline{b'}c\big)\;\Ra\;\Aline a{c'}\parallel\Aline {a'}c.$$
\end{definition}

\begin{picture}(60,60)(-150,-10)
\linethickness{0.75pt}
\put(0,0){\color{teal}\line(1,0){60}}
\put(0,40){\color{teal}\line(1,0){60}}

\put(0,0){\color{blue}\line(1,1){40}}
\put(20,0){\color{blue}\line(1,1){40}}
\put(0,40){\color{violet}\line(1,-2){20}}
\put(40,40){\color{violet}\line(1,-2){20}}
\put(0,0){\color{red}\line(0,1){40}}
\put(60,0){\color{red}\line(0,1){40}}

\put(0,0){\circle*{3}}
\put(-2,-9){$a$}
\put(20,0){\circle*{3}}
\put(17,-9){$b$}
\put(60,0){\circle*{3}}
\put(57,-9){$c$}
\put(0,40){\circle*{3}}
\put(-2,43){$c'$}
\put(40,40){\circle*{3}}
\put(38,43){$b'$}
\put(60,40){\circle*{3}}
\put(58,43){$a'$}
\end{picture}


\begin{theorem}\label{t:commutative-puls<=>} An affine space $X$ is commutative-puls if and only if for every ternar $R$ of $X$, the puls loop $(R,\!\puls\!)$ is commutative.
\end{theorem}

\begin{proof} Assume that an affine space $X$ is commutative-puls. Given any plane $\Pi\subseteq X$ and an affine base $uow$ in $\Pi$, we should prove that the puls loop $(\Delta,\!\puls\!)$ is commutative. By Proposition~\ref{p:commutative-puls<=>}, it suffices to show that for every points $x,y,z\in \Delta$, the parallelity relation $\Aline {oo}{ox}\parallel \Aline {oy}{ez}$ implies  $\Aline{oo}{ey}\parallel \Aline{ox}{ez}$.

Consider the points $a\defeq ox$, $b\defeq oy$, $c\defeq oo$ on the line $L\defeq \Aline ow$ and the points $a'\defeq ey$, $b'\defeq ex$ and $c'\defeq ez$ on the line $L'\defeq \Aline ow$. Since $\Aline a{b'}=\Aline {ox}{ex}\parallel \Aline {oy}{ey}=\Aline b{a'}$ and $\Aline b{c'}=\Aline{oy}{ez}\parallel \Aline {oo}{ex}=\Aline{b'}c$, we can apply the commutativity-puls of $X$ and conclude that $\Aline {ox}{ez}=\Aline a{c'}\parallel \Aline{a'}c=\Aline{ey}{oo}$. This completes the proof of the ``only if'' part.

\begin{picture}(150,135)(-170,-55)
\linethickness{0.8pt}
\put(-40,-40){\line(1,1){100}}
\put(0,-40){\color{cyan}\line(0,1){100}}
\put(40,-40){\color{cyan}\line(0,1){100}}
\put(0,-40){\color{blue}\line(2,3){40}}
\put(0,0){\color{blue}\line(2,3){40}}
\put(-40,-40){\color{teal}\line(1,0){80}}
\put(0,60){\color{teal}\line(1,0){60}}
\put(0,0){\color{red}\line(1,-1){40}}
\put(20,20){\color{teal}\line(1,0){20}}
\put(0,60){\color{red}\line(1,-1){40}}

\put(-40,-40){\circle*{3}}
\put(-44,-48){$y$}
\put(0,-40){\circle*{3}}
\put(-5,-48){$oy$}
\put(-7,-37){$b$}
\put(40,-40){\circle*{3}}
\put(35,-48){$ey$}
\put(43,-38){$a'$}
\put(0,0){\circle*{3}}
\put(-8,-2){$o$}
\put(4,-2){$c$}
\put(20,20){\circle*{3}}
\put(22,15){$z$}
\put(40,20){\circle*{3}}
\put(43,22){$c'$}
\put(43,14){$ez$}
\put(40,40){\circle*{3}}
\put(42,35){$e$}
\put(0,60){\circle*{3}}
\put(-5,63){$ox$}
\put(-8,56){$a$}
\put(40,60){\circle*{3}}
\put(28,63){$ex$}
\put(42,63){$b'$}
\put(60,60){\circle*{3}}
\put(60,63){$x$}
\end{picture}

Now assume that for every ternar $R$ of the affine space $X$, the puls loop $(R,\!\puls\!)$ is commutative. To prove that $X$ is commutative-puls, take any distinct parallel lines $L,L'\in X$ and points $a,b,c\in L$ and $a',b',c'\in L'$ such that $\Aline a{b'}\parallel \Aline{a'}b$ and $\Aline b{c'}\parallel \Aline{b'}c$. We have to show that $\Aline a{c'}=\Aline {a'}c$. In the affine plane $\Pi\defeq\overline{L\cup L'}$, choose an affine base $uow$ such that $o=c$, $u\in L'$, $\Aline ow=L$ and $\Aline ou\parallel \Aline a{b'}\parallel \Aline{a'}b$. Let $e$ be the diunit of the affine base $uow$, $\Delta\defeq\Aline oe$ be the diagonal, and $\boldsymbol h\defeq(\Aline ou)_\parallel$ be the horizontal direction of the affine base $uow$.

Consider the unique points $x\in\Delta\cap\Aline a{b'}$, $y\in\Delta\cap\Aline {a'}b$, $z\in\Delta\cap\Aline {c'}{\boldsymbol h}$. Observe that $\Aline{oo}{ex}=\Aline c{b'}\parallel\Aline b{c'}=\Aline{oy}{ez}$. Since the loop $(\Delta,\!\puls\!)$ is commutative, we can apply Proposition~\ref{p:commutative-puls<=>} and conclude that $\Aline c{a'}=\Aline {oo}{ey}\parallel \Aline {ox}{ez}=\Aline a{c'}$, witnessing that the affine space $X$ is commutative-puls.
\end{proof}

The following theorem is an ``puls'' counterpart of Hessenberg's  Theorem~\ref{t:Hessenberg-affine}.

\begin{theorem}\label{t:com-puls=>ass-puls} Every  commutative-puls affine space is associative-puls.
\end{theorem}

\begin{proof} Given any disjoint parallel lines $L,L'$ in a commutative-puls affine space $X$ and points $a,b,\alpha,\beta\in L$ and $a',b',\alpha',\beta'\in L'$ with $\Aline a{a'}\parallel \Aline \alpha{\alpha'}$, $\Aline a{b'}\parallel \Aline\alpha{\beta'}$ and $\Aline {a'}b\parallel \Aline{\alpha'}{\beta}$ we should prove that $\Aline b{b'}\parallel \Aline \beta{\beta'}$. Since the affine space $X$ is Playfair, there exist points $c\in L$ and $c'\in L'$ such that $\Aline c{b'}\parallel \Aline a{a'}\parallel \Aline \alpha{\alpha'}$ and $\Aline c{c'}\parallel \Aline b{a'}\parallel \Aline \beta{\alpha'}$.
Since $\Aline a{a'}\parallel \Aline {b'}c$ and $\Aline c{c'}\parallel \Aline {a'}b$, we can apply the commutativity-puls of $X$ to the triples $acb$ and $b'a'c'$ and conclude that $\Aline a{c'}\parallel \Aline {b'}b$. Since $\Aline\alpha {\alpha'}\parallel \Aline{b'}c$ and $\Aline c{c'}\parallel \Aline {\alpha'}\beta$, we can apply the commutativity-puls to the triples $\alpha c\beta$ and $b'\alpha'c'$, and conclude that $\Aline{\alpha}{c'}\parallel \Aline {b'}b$. 
Since $\Aline a{b'}\parallel \Aline {\beta'}\alpha$ and $\Aline \alpha{c'}\parallel\Aline {b'}\beta$, we can apply the commutativity-plus to the triples $a\alpha\beta$ and $\beta'b'c'$, and conclude that $\Aline a{c'}\parallel \Aline \beta{\beta'}$ and hence $\Aline b{b'}\parallel \Aline a{c'}\parallel \Aline\beta{\beta'}$.

\begin{picture}(100,55)(-140,-10)

\linethickness{0.7pt}
\put(10,0){\color{teal}\line(1,0){80}}
\put(0,30){\color{teal}\line(1,0){100}}
\put(10,0){\color{cyan}\line(-1,3){10}}
\put(40,30){\color{cyan}\line(1,-3){10}}
\put(10,0){\color{violet}\line(1,1){30}}
\put(30,0){\color{blue}\line(-1,1){30}}
\put(10,0){\color{red}\line(1,3){10}}
\put(30,0){\color{red}\line(1,3){10}}

\put(70,0){\color{cyan}\line(-1,3){10}}
\put(70,0){\color{violet}\line(1,1){30}}
\put(90,0){\color{blue}\line(-1,1){30}}
\put(90,0){\color{red}\line(1,3){10}}

\put(40,30){\color{olive}\line(5,-3){50}}
\put(20,30){\color{olive}\line(5,-3){50}}
\put(20,30){\color{blue}\line(1,-1){30}}

\put(20,30){\circle*{3}}
\put(17,33){$c'$}
\put(50,0){\circle*{3}}
\put(47,-9){$c$}
\put(10,0){\circle*{3}}
\put(8,-9){$a$}
\put(30,0){\circle*{3}}
\put(27,-9){$b$}
\put(0,30){\circle*{3}}
\put(-3,33){$a'$}
\put(40,30){\circle*{3}}
\put(38,33){$b'$}

\put(70,0){\circle*{3}}
\put(68,-9){$\alpha$}
\put(90,0){\circle*{3}}
\put(87,-9){$\beta$}
\put(60,30){\circle*{3}}
\put(57,33){$\alpha'$}
\put(100,30){\circle*{3}}
\put(98,33){$\beta'$}
\end{picture}
\end{proof}

Theorems~\ref{t:diagonal-trans=>ass-plus} and \ref{t:paraD<=>translation}  imply the following important corollary. 

\begin{corollary}\label{c:Thales=>commutative-puls} Every Thalesian affine space is commutative-puls.
\end{corollary}

\section{Puls loops of Boolean affine spaces}

We recall that a liner $X$ is \defterm{Boolean} if 
$$\forall a,b,c,d\in X\;(\Aline ab\cap\Aline cd=\varnothing=\Aline bc\cap\Aline ad\;\Ra\;\Aline ac\cap\Aline bd=\varnothing).$$

Proposition~\ref{p:Boolean-puls<=>} implies the following characterization of Boolean affine spaces.

\begin{theorem}\label{t:Boolean<=>Boolean-puls} An affine space is Boolean $X$ if and only if for every ternar $R$ of $X$ its puls loop $(R,\!\puls\!)$ is Boolean.
\end{theorem}

\begin{theorem}\label{t:Boolean=>commutative-puls} Every Boolean affine space is commutative-puls.
\end{theorem}

\begin{proof} Let $X$ be a Boolean affine space. To show that $X$ is commutative-puls, fix any distinct parallel lines $L,L'\subseteq X$ and points  $a,b,c\in L$ and $a',b',c'\in L'$ such that $\Aline a{b'}\parallel \Aline {a'}b$ and $\Aline b{c'}\parallel \Aline {b'}c=\varnothing$. Since $X$ is Boolean, the parallelogram $ab'b'b$ has parallel diagonals and hence $\Aline a{a'}\parallel\Aline b{b'}$. Also the parallelogram $bc'b'c$ has parallel diagonals $\Aline c{c'}\parallel \Aline b{b'}$. Then $\Aline a{a'}\parallel \Aline b{b'}\parallel \Aline c{c'}$ and hence $aa'c'c$ is a parallelogram, which has parallel diagonals $\Aline a{c'}$ and $\Aline {a'}c$.

\begin{picture}(60,65)(-150,-10)
\linethickness{0.75pt}
\put(0,0){\color{teal}\line(1,0){60}}
\put(0,40){\color{teal}\line(1,0){60}}

\put(0,0){\color{blue}\line(1,1){40}}
\put(20,0){\color{blue}\line(1,1){40}}
\put(0,40){\color{violet}\line(1,-2){20}}
\put(40,40){\color{violet}\line(1,-2){20}}
\put(0,0){\color{cyan}\line(0,1){40}}
\put(60,0){\color{cyan}\line(0,1){40}}
\put(0,0){\color{red}\line(3,2){60}}
\put(0,40){\color{red}\line(3,-2){60}}
\put(20,0){\color{red}\line(1,2){20}}
\put(30,20){\color{white}\circle*{3}}

\put(0,0){\circle*{3}}
\put(-2,-9){$a$}
\put(20,0){\circle*{3}}
\put(17,-9){$b$}
\put(60,0){\circle*{3}}
\put(57,-9){$c$}
\put(0,40){\circle*{3}}
\put(-2,43){$c'$}
\put(40,40){\circle*{3}}
\put(38,43){$b'$}
\put(60,40){\circle*{3}}
\put(58,43){$a'$}
\end{picture}
\end{proof}

Corollary~\ref{c:inv-puls=>Boolean-paralelogram} implies the following sufficient condition of the existence of Boolean parallelograms in invertible-puls affine spaces.

\begin{corollary}\label{c:inv-puls=>Boolean-parallelogram2} Every invertible-puls affine space of even order contains a Boolean parallelogram.
\end{corollary}

For inversive-puls affine spaces, Corollary~\ref{c:inv-puls=>Boolean-parallelogram2} can be reversed.

\begin{proposition}\label{p:inversive-puls<=>even<=>Boole} Let $X$ be an inversive-puls affine space of finite order. The liner $X$ contains a Boolean parallelogram if and only if $|X|_2$ is even.
\end{proposition}

\begin{proof} The ``if'' part follows from Corollary~\ref{c:inv-puls=>Boolean-parallelogram2}. To prove the ``only if'' part, assume that the liner $X$ contains a Boolean parallelogram $uowe$. Consider the affine base $uow$ in the plane $\Pi\defeq\overline{\{u,o,w,e\}}$ and its ternar $\Delta=\Aline oe$. Since $uowe$ is a Boolean parallelogram, $e\puls e=o$ in the ternar $\Delta$, by Lemma~\ref{l:Boolean-puls}. Since $X$ is inversive-puls, the puls loop $(\Delta,\puls)$ is inversive. Applying Proposition~\ref{p:Ali-Slaney}, we conclude that the number $|X|_2=|\Delta|$ is even.
\end{proof}

\section{Puls properties of prime affine spaces}

An affine space is defined to be \defterm{prime} if its order is a prime number.

\begin{theorem}\label{t:prime-affine<=>+} For every prime affine space $X$, the following conditions are equivalent:
\begin{enumerate}
\item $X$ is Pappian;
\item $X$ is Desarguesian;
\item $X$ is Thalesian;
\item $X$ is translation;
\item $X$ is $\partial$-translation;
\item $X$ is commutative-plus;
\item $X$ is associative-plus;
\item $X$ is inversive-plus;
\item $X$ is commutative-puls;
\item $X$ is associative-puls;
\item $X$ is inversive-puls;
\item $X$ is left-inversive-puls.
\end{enumerate}
\end{theorem}

\begin{proof} The equivalence of the conditions (1)--(8) was proved in Theorem~\ref{t:prime-affine<=>}. The implication $(3)\Ra(9)$ follows from Corollary~\ref{c:Thales=>commutative-puls}, $(9)\Ra(10)$ follows from Theorem~\ref{t:com-puls=>ass-puls}, and $(10)\Ra(11)\Ra(12)$ are trivial.

It remains to prove that $(12)\Ra(5)$. Assume that a prime affine  space $X$ is left-inversive-puls. 
 By Corollary~\ref{c:left-inversive-puls-prime=>cyclic}, for every ternar $R$ of $X$, the additive loop $(R,\!\puls\!)$ is a cyclic group. By Theorem~\ref{t:associative-puls<=>}, $X$ is associative-puls, and by 
Theorem~\ref{t:ass-puls=>partialT}, $X$ is $\partial$-translation. 
\end{proof}

\begin{remark}\label{diag:puls-additive} By Theorems~\ref{t:paraD<=>translation}, \ref{t:Boolean=>commutative-puls},  \ref{t:prime-affine<=>+}, and Corollary~\ref{c:Thales=>commutative-puls}, for any affine plane we have the implications:
$$
\xymatrix@C=40pt@R=19pt{
\mbox{Pappian}\ar@{=>}[d]\\
\mbox{Desarguesian}\ar@{=>}[d]\ar@/_10pt/_{+\mbox{\em \scriptsize prime}}[u]\\
\mbox{Thalesian}\ar@{=>}[d]\ar@{<=>}[r]\ar@/_10pt/_{+\mbox{\em \scriptsize prime}}[u]&\mbox{translation}\ar@{=>}[ddd]\\
\mbox{Boolean or Thalesian}\ar@{=>}[d]&\\
\mbox{commutative-puls}\ar@{=>}[d]&\\
\mbox{associative-puls}\ar@{=>}[d]\ar_{\phantom{.}+\mbox{\em \scriptsize prime-power}}[r]&\mbox{$\partial$-translation}\ar@/_10pt/_{+\mbox{\em \scriptsize prime}}[uuu]\\
\mbox{inversive-puls}\ar@/_10pt/_{+\mbox{\em \scriptsize prime}}[u]\ar@{=>}[d]\\
\mbox{left-inversive-puls}
\ar@/_10pt/_{+\mbox{\em \scriptsize prime}}[u]\ar@{=>}[d]
\\
\mbox{invertible-puls}
}
$$
\end{remark}

This diagram motivates the following problems.

\begin{problem} Is every associative-puls affine plane commutative-puls?
\end{problem}

\begin{problem} Is every (finite) left-invertible-puls affine plane inversive-puls?
\end{problem}

\begin{problem} Is every (finite) inversive-puls affine plane associative-puls?
\end{problem}

\begin{problem} Is every right-inversive-puls affine plane right-Bol?
\end{problem}
\chapter{Dot properties of affine spaces}\label{ch:dot-aff}

In this chapter we shall consider some affine configuration axioms describing properties of the dot operation in ternars of affine spaces.

\section{Invertible-dot affine spaces}

\begin{definition}\label{d:a-invertible-dot} An affine liner $X$ is defined to be \index{invertible-dot affine space}\index{affine space!invertible-dot}\defterm{invertible-dot} if for every points\\ $o,a,b,d,\beta,\gamma,\delta\in X$ $$\big[(\Aline ab\parallel \Aline d\beta\parallel \Aline \delta\gamma\nparallel \Aline ad\parallel \Aline b\delta\parallel \Aline \beta\gamma)\;\wedge\;(\Aline b\delta\cap\Aline \beta d\cap\Aline a\gamma\ne\varnothing)\;\wedge\;(o\in \Aline \alpha\gamma\cap\Aline b\beta)\big]\;\Ra\;(\,o\in \Aline d\delta\,).$$
\end{definition}

\begin{picture}(200,95)(-140,-5)
\linethickness{0.75pt}
\put(0,0){\color{red}\line(1,1){80}}
\put(0,0){\color{orange}\line(2,1){80}}
\put(0,0){\color{magenta}\line(1,2){40}}
\put(20,20){\color{blue}\line(1,0){20}}
\put(20,20){\color{cyan}\line(0,1){20}}
\put(40,40){\color{blue}\line(-1,0){20}}
\put(40,40){\color{cyan}\line(0,-1){20}}
\put(40,40){\color{blue}\line(1,0){40}}
\put(40,40){\color{cyan}\line(0,1){40}}
\put(80,80){\color{blue}\line(-1,0){40}}
\put(80,80){\color{cyan}\line(0,-1){40}}

\put(0,0){\color{red}\circle*{3}}
\put(20,20){\circle*{3}}
\put(12,18){$a$}
\put(40,40){\circle*{3}}
\put(80,80){\circle*{3}}
\put(83,80){$\gamma$}
\put(20,40){\circle*{3}}
\put(14,42){$b$}
\put(40,20){\circle*{3}}
\put(40,12){$d$}
\put(80,40){\circle*{3}}
\put(83,38){$\delta$}
\put(40,80){\circle*{3}}
\put(37,83){$\beta$}

\end{picture}

\begin{theorem}\label{t:invertible-dot<=>} An affine space $X$ is  invertible-dot if and only if for every ternar $R$ of $X$, its dot loop $(R^*,\cdot)$ is invertible.
\end{theorem}

\begin{proof} Assume that the affine space $X$ is invertible-dot. Given any ternar $R$ of $X$, find a plane $\Pi\subseteq X$ and an affine base $uow$ in $\Pi$ whose ternar $\Delta$ is isomorphic to the ternar $R$. Since the ternars $R$ and $\Delta$ are isomorphic, it suffices to prove that the loop $(\Delta^*,\cdot)$ is invertible. By Proposition~\ref{p:dot-dot-inverses<=>}, it suffices to show that for every points $x,y\in\Delta^*$ with $o\in \Aline {xe}{ey}$ we have $o\in\Aline {ex}{ye}$. Here we identify the plane $\Pi$ with the coordinate plane $\Delta^2$ of the ternar $\Delta$, and denote by $e$ the diunit of the affine base $uow$.

\begin{picture}(200,110)(-140,-15)
\linethickness{0.75pt}
\put(0,0){\color{red}\line(1,1){80}}
\put(0,0){\color{orange}\line(2,1){80}}
\put(0,0){\color{magenta}\line(1,2){40}}
\put(20,20){\color{teal}\line(1,0){20}}
\put(20,20){\color{cyan}\line(0,1){20}}
\put(40,40){\color{teal}\line(-1,0){20}}
\put(40,40){\color{cyan}\line(0,-1){20}}
\put(40,40){\color{teal}\line(1,0){40}}
\put(40,40){\color{cyan}\line(0,1){40}}
\put(80,80){\color{teal}\line(-1,0){40}}
\put(80,80){\color{cyan}\line(0,-1){40}}

\put(0,0){\color{red}\circle*{3}}
\put(0,0){\color{white}\circle*{2}}
\put(-2,-7){$o$}
\put(20,20){\circle*{3}}
\put(12,18){$x$}
\put(20,14){$a$}
\put(40,40){\color{red}\circle*{3}}
\put(41.5,33){$e$}
\put(80,80){\circle*{3}}
\put(83,78){$y$}
\put(74,83){$\gamma$}
\put(20,40){\circle*{3}}
\put(8,42){$xe$}
\put(21,32){$b$}
\put(40,20){\circle*{3}}
\put(40,12){$ex$}
\put(34,21){$d$}
\put(80,40){\circle*{3}}
\put(83,38){$ye$}
\put(74,42){$\delta$}
\put(40,80){\circle*{3}}
\put(37,84){$ey$}
\put(41,72){$\beta$}

\end{picture}

 The inclusion $o\in\Aline {xe}{ey}$ implies $x\ne e\ne y$.  Consider the points $a\defeq xx$, $b\defeq xe$, $d\defeq ex$, $\beta\defeq ey$, $\gamma=yy$ and $\delta\defeq ye$. Observe that $e\in \Aline {xe}{ey}\cap\Aline {ex}{ey}\cap\Delta=\Aline b\delta\cap\Aline d\beta\cap\Aline a\gamma$ and $o\in\Aline {xe}{ey}\cap\Delta=\Aline b\beta\cap\Aline a\gamma$. The invertibility-dot of $X$ ensures that $o\in \Aline d\delta=\Aline {xe}{ey}$. By Proposition~\ref{p:dot-dot-inverses<=>}, the loop $(\Delta^*,\cdot)$ is invertible and so is the loop $(R^*,\cdot)$.
\smallskip

Now assume that for every ternar $R$ of $X$, its multiplicative loop $(R^*,\cdot)$ is invertble. To prove that $X$ is invertible-dot, take any points $o,a,b,d,\beta,\gamma,\delta\in X$ with $$(\Aline ab\parallel \Aline d\beta\parallel \Aline \delta\gamma\nparallel \Aline ad\parallel \Aline b\delta\parallel \Aline \beta\gamma)\;\wedge\;(\Aline b\delta\cap\Aline \beta d\cap\Aline a\gamma\ne\varnothing)\;\wedge\;(o\in \Aline \alpha\gamma\cap\Aline b\beta).$$
We have to prove that $o\in \Aline d\delta$.  It follows from $\Aline b\delta\nparallel\Aline \beta  d$ that the nonempty set $\Aline b\delta\cap\Aline \beta d\cap\Aline a\gamma$ is a singleton containing a unique point $e$.

If $a=b$, then $\Aline ab\parallel \Aline d\beta\parallel \Aline \delta\gamma\nparallel \Aline ad\parallel \Aline b\delta\parallel \Aline \beta\gamma$ implies $d=\beta$, $\delta=e=\gamma$, $a\ne d$, $b\ne\delta$, $\beta\ne \gamma$. Then $o\in\Aline b\beta=\Aline e\gamma=\Aline d\delta$. 

If $a=d$, then $\Aline ab\parallel \Aline d\beta\parallel \Aline \delta\gamma\nparallel \Aline ad\parallel \Aline b\delta\parallel \Aline \beta\gamma$ implies $b=e=\delta$, $\beta=\gamma$, $a\ne b$, $d\ne\beta$, $\delta\ne \gamma$. Then $o\in\Aline b\beta=\Aline ae=\Aline d\delta$. 

So, assume that $b\ne a\ne d$. In this case $d\ne\beta\ne\gamma$, $b\ne\delta\ne\gamma$, $e\notin \{a,b,d,\beta,\gamma,\delta\}$, and $o\in \Aline a\gamma\setminus \{a,e,\gamma\}$. Let $u\in \Aline d\beta$ and $w\in\Aline b\delta$ be unique points such that $\Aline ou\parallel \Aline b\delta$ and $\Aline ow\parallel\Aline\beta d$. Then  $uow$ is an affine base for the plane $\Pi\defeq\overline{\{o,a,b,d,e,\beta,\gamma,\delta\}}$. By our asumption the multiplicative loop $(\Delta^*,\cdot)$ of the ternary ring $\Delta=\Aline oe=\Aline \alpha\gamma$ is invertible. Consider the points $x\defeq a$ and $y\defeq\gamma$ in $\Delta$ and observe that $o\in\Aline b\beta=\Aline {xe}{ey}$. Here we identify the plane $\Pi$ with the coordinate plane $\Delta^2$ of the ternar $\Delta$. Since the multiplicative loop $(\Delta^*,\cdot)$ is invertible, we can apply Proposition~\ref{p:dot-dot-inverses<=>} and conclude that $o\in \Aline {ex}{ye}=\Aline d\delta$.
\end{proof}

\begin{proposition} Every invertible-dot affine space is invertible-plus.
\end{proposition}

\begin{proof} Let $X$ be an invertible-dot affine space. To prove that $X$ is invertible-plus, take any points $a,b,c,\alpha,\beta,\gamma\in X$ and $o\in\Aline a\alpha\cap\Aline b\beta\cap\Aline c\gamma$ such that
$$(\Aline ac\parallel \Aline ob\parallel \Aline o\beta\parallel\Aline\alpha\gamma\nparallel \Aline b\gamma\parallel\Aline o\alpha\parallel \Aline oa\parallel\Aline c\beta)\;\wedge\;(\Aline ab\parallel\Aline co).$$ We have to prove that $\Aline co\parallel\Aline \alpha\beta$. 

\begin{picture}(150,80)(-180,-45)
\linethickness{0.75pt}
\put(0,0){\color{blue}\line(1,0){30}}
\put(0,0){\color{cyan}\line(0,1){30}}
\put(0,0){\color{blue}\line(-1,0){30}}
\put(0,0){\color{cyan}\line(0,-1){30}}

\put(0,30){\color{blue}\line(1,0){30}}
\put(30,0){\color{cyan}\line(0,1){30}}
\put(-30,-30){\color{blue}\line(1,0){30}}
\put(-30,-30){\color{cyan}\line(0,1){30}}
\put(-30,0){\color{red}\line(1,1){30}}
\put(0,-30){\color{red}\line(1,1){30}}
{\linethickness{1pt}
\put(0,0){\color{red}\line(1,1){30}}
\put(0,0){\color{red}\line(-1,-1){30}}
}
\put(0,0){\circle*{3}}
\put(2,-7){$o$}
\put(30,0){\circle*{3}}
\put(33,-2){$\alpha$}
\put(-30,0){\circle*{3}}
\put(-38,-2){$a$}
\put(0,30){\circle*{3}}
\put(-7,31){$b$}
\put(0,-30){\circle*{3}}
\put(2,-36){$\beta$}
\put(30,30){\circle*{3}}
\put(32,31){$\gamma$}
\put(-30,-30){\circle*{3}}
\put(-37,-36){$c$}
\end{picture}

If $c=a$, then $\Aline ac\parallel \Aline ob\parallel \Aline o\beta\parallel\Aline\alpha\gamma\nparallel \Aline b\gamma\parallel\Aline o\alpha\parallel \Aline oa\parallel\Aline c\beta$ implies $b=o=\beta$, $\alpha\ne \gamma\ne b$, $a\ne o\ne \alpha$, $c\ne\beta=o$. Then $\Aline \alpha\beta=\Aline co$ and we are done.

If $c=\beta$, then $\Aline ac\parallel \Aline ob\parallel \Aline o\beta\parallel\Aline\alpha\gamma\nparallel \Aline b\gamma\parallel\Aline o\alpha\parallel \Aline oa\parallel\Aline c\beta$ implies $a\ne c$, $b\ne o\ne\beta$, $\alpha\ne\gamma$, $\beta\ne\gamma$ and $\alpha=o=a$. Then $\Aline \alpha\beta=\Aline oc$ and we are done.

So, assume that $a\ne c\ne \beta$. Assuming that $\Aline\alpha\beta\nparallel \Aline oc$, we can find a unique point $o'\in\Aline oc\cap\Aline \alpha\beta$. The invertibility-dot of $X$ implies $o'\in \Aline ab$, which contradicts $\Aline ab\parallel \Aline oc$ and $a\ne o$.
\end{proof}

\section{Inversive-dot affine spaces}

\begin{definition} An affine liner $X $ is called \index{left-inversive-dot affine space}\index{affine space!left-inversive-dot}\defterm{left-inversive-dot} if\newline for every points $o,a,b,c,d,\alpha,\beta,\gamma,\delta\in X$ $$\big[(\Aline ab\parallel \Aline cd\parallel \Aline \alpha\beta\parallel \Aline \gamma\delta\nparallel \Aline ad\parallel \Aline bc\parallel \Aline \beta\gamma\parallel \Aline \alpha\delta)\;\wedge\;(\Aline bc=\Aline\alpha\delta)\;\wedge\;(o\in\Aline a\alpha\cap\Aline b\beta\cap\Aline c\gamma)\big]\;\Ra\;(\,o\in\Aline d\delta\,).$$
\end{definition}

\begin{picture}(100,150)(-150,-15)

\linethickness{0.7pt}

\put(60,30){\color{cyan}\line(0,1){30}}
\put(60,60){\color{blue}\line(1,0){120}}
\put(180,120){\color{cyan}\line(0,-1){60}}
\put(60,30){\color{blue}\line(1,0){30}}
\put(90,60){\color{cyan}\line(0,-1){30}}
\put(120,120){\color{cyan}\line(0,-1){60}}
\put(120,120){\color{blue}\line(1,0){60}}

\put(0,0){\color{red}\line(1,1){120}}
\put(0,0){\color{red}\line(3,2){180}}
\put(0,0){\color{red}\line(2,1){120}}
\put(0,0){\color{red}\line(3,1){180}}

\put(0,0){\color{red}\circle*{3}}
\put(-9,-5){$o$}
\put(60,60){\circle*{3}}
\put(57,63){$b$}
\put(120,120){\circle*{3}}
\put(115,123){$\beta$}
\put(60,30){\circle*{3}}
\put(57,22){$a$}
\put(90,30){\circle*{3}}
\put(90,22){$d$}
\put(90,60){\circle*{3}}
\put(87,63){$c$}
\put(120,60){\circle*{3}}
\put(119,53){$\alpha$}
\put(180,60){\circle*{3}}
\put(183,55){$\delta$}
\put(180,120){\circle*{3}}
\put(183,119){$\gamma$}
\end{picture}

Proposition~\ref{p:left-inversive-dot<=>} implies the following algebraic characterization of left-inversive-dot affine spaces.

\begin{theorem}\label{t:left-inversive-dot<=>} An affine space $X$ is left-inversive-dot if and only if every ternar $R$ of $X$, the multiplicative loop $(R^*,\cdot)$ is left-inversive.
\end{theorem}

\begin{exercise} Write down the proof of Theorem~\ref{t:left-inversive-dot<=>}.
\end{exercise}

Since every left-inversive loop is invertible, Theorems~\ref{t:invertible-dot<=>} and \ref{t:left-inversive-dot<=>} imply

\begin{corollary} Every left-inversive-dot affine space is invertible-dot.
\end{corollary} 


\begin{definition}\label{d:r-inversive-dot} An affine liner $X$ is called \index{right-inversive-dot affine space}\index{affine space!right-inversive-dot}\defterm{right-inversive-dot} if for every line $D\subset X$ and points $o,a,c,\alpha,\gamma\in D$ and $b,d,\beta,\delta\in X\setminus D$ $$\big[(\Aline ab\parallel \Aline cd\parallel \Aline \alpha\beta\parallel \Aline \gamma\delta\nparallel \Aline ad\parallel \Aline bc\parallel \Aline \beta\gamma\parallel \Aline \alpha\delta)\;\wedge\;(o\in \Aline b\beta)\big]\;\Ra\;(o\in \Aline d\delta).$$
\end{definition}

\begin{picture}(200,150)(-150,-15)
\linethickness{0.75pt}
\put(0,0){\color{red}\line(1,1){120}}
\put(0,0){\color{orange}\line(2,1){120}}
\put(0,0){\color{magenta}\line(1,2){60}}
\put(20,20){\color{blue}\line(1,0){20}}
\put(20,20){\color{cyan}\line(0,1){20}}
\put(40,40){\color{blue}\line(-1,0){20}}
\put(40,40){\color{cyan}\line(0,-1){20}}
\put(60,60){\color{blue}\line(1,0){60}}
\put(60,60){\color{cyan}\line(0,1){60}}
\put(120,120){\color{blue}\line(-1,0){60}}
\put(120,120){\color{cyan}\line(0,-1){60}}


\put(0,0){\color{red}\circle*{3}}
\put(-3,-7){$o$}
\put(20,20){\circle*{3}}
\put(13,20){$a$}
\put(40,40){\circle*{3}}
\put(42,36){$c$}
\put(60,60){\circle*{3}}
\put(52,60){$\alpha$}
\put(120,120){\circle*{3}}
\put(123,120){$\gamma$}
\put(20,40){\circle*{3}}
\put(13,40){$b$}
\put(40,20){\circle*{3}}
\put(40,12){$d$}
\put(120,60){\circle*{3}}
\put(123,58){$\delta$}
\put(60,120){\circle*{3}}
\put(51,118){$\beta$}

\end{picture}

Proposition~\ref{p:right-inversive-dot<=>} implies the following algebraic characterization of right-inversive-dot affine spaces.

\begin{theorem}\label{t:right-inversive-dot<=>} An affine space $X$ is right-inversive-dot if and only if for every ternar $R$ of $X$, the multiplicative loop $(R^*,\cdot)$ is right-inversive.
\end{theorem}

\begin{exercise} Write down the proof of Theorem~\ref{t:right-inversive-dot<=>}.
\end{exercise}

\begin{proposition}\label{p:right-inv-dot=>inv-plus} Every-right inversive-dot affine space is inversive-plus.
\end{proposition}

\begin{proof} Let $X$ be a right inversive-dot affine space. To prove that $X$ is inversive-plus, take any line $D\subset X$ and points $a,d,\alpha,\delta\in D$ and $b,c,\beta,\gamma\in X\setminus D$ such that $$(\Aline ab\parallel\Aline cd\parallel\Aline \alpha\beta\parallel\Aline\gamma\delta\nparallel 
\Aline ac\parallel\Aline bd\parallel\Aline \alpha\gamma\parallel\Aline\beta\delta)\;\wedge\;(\Aline b\beta\parallel D).$$
We have to prove that $\Aline c\gamma\parallel D$. The above parallelity relations ensure that $\Pi\defeq\overline{\{a,b,c,d,\alpha,\beta,\gamma,\delta\}}$ is a plane.

\begin{picture}(150,90)(-170,-40)

\linethickness{0.75pt}
\put(0,0){\color{blue}\line(-1,0){30}}
\put(0,0){\color{cyan}\line(0,-1){30}}
\put(-30,-30){\color{blue}\line(1,0){30}}
\put(-30,-30){\color{cyan}\line(0,1){30}}

{\linethickness{1pt}
\put(-30,-30){\color{red}\line(1,1){75}}
}
\put(-30,0){\color{red}\line(1,1){45}}
\put(0,-30){\color{red}\line(1,1){45}}
\put(0,0){\color{red}\line(-1,-1){30}}
\put(15,15){\color{blue}\line(1,0){30}}
\put(15,15){\color{cyan}\line(0,1){30}}
\put(45,45){\color{blue}\line(-1,0){30}}
\put(45,45){\color{cyan}\line(0,-1){30}}
\put(-30,-30){\circle*{3}}
\put(-36,-37){$a$}
\put(0,0){\circle*{3}}
\put(2,-5){$d$}
\put(-30,0){\circle*{3}}
\put(-37,-1){$b$}
\put(0,-30){\circle*{3}}
\put(0,-38){$c$}
\put(15,45){\circle*{3}}
\put(11,48){$\beta$}
\put(45,45){\circle*{3}}
\put(47,46){$\delta$}
\put(45,15){\circle*{3}}
\put(47,10){$\gamma$}
\put(15,15){\circle*{3}}
\put(6,14){$\alpha$}

\end{picture}

Assuming that $\Aline c\gamma\nparallel D$ and taking into account that $\gamma\in \Pi\setminus D$, we conclude that $\Aline c\gamma\cap D=\{o'\}$ for some point $o'$. Since $X$ is right-inversive-dot, $o'\in \Aline b\beta\cap D$, which contradicts $\Aline b\beta\parallel D$ and $b,\beta\notin D$.
\end{proof}

Let us recall that a liner $X$ is \index{hypersymmetric liner}\index{liner!hypersymmetric}\defterm{hypersymmetric} if for every triangle $abc$ in $X$, there exists an involutive hyperfixed automorphism $A:X\to X$ such that $Aabc=cba$.

\begin{proposition}\label{p:rinvdot=>hypersymmetric} Every right-inversive-dot affine space $X$ is hypersymmetric.
\end{proposition}

\begin{proof} Given any triangle $abc$ in $X$, we need to construct an involutive hyperfixed automorphism $A$ of $X$ such that $Aabc=cba$. By Theorem~\ref{t:parallelogram3+1}, there exists a unique point $d\in X$ such that $abcd$ is a parallelogram. 

If $\|X\|=3$, then consider  the line $\Delta \defeq \Aline ad$ in $X$ and define an involutive bijection $A:X\to X$ as follows. Given any point $x\in X$, find unique points $b',d'\in \Delta$ such that $\Aline x{b'}\subparallel \Aline ab$ and $\Aline x{d'}\subparallel \Aline ad$. Next, find a unique point $y\in X$ such that $\Aline y{b'}\subparallel \Aline cb$ and $\Aline y{d'}\subparallel \Aline cd$. Put $A(x)\defeq y$. Thus we have defined an involutive bijection $A:X\to X$ such that $Aabc=cba$ and $\Aline bd=\Fix(A)$. The definition of a right-inversive-dot affine space and Proposition~\ref{p:right-inv-dot=>inv-plus} imply that the bijection $A$ is an automorphism of the Playfair plane $X$. Since $\Fix(A)=\Aline bd$, the automorphism $A$ is hyperfixed.

\begin{picture}(200,160)(-150,-15)
\linethickness{0.75pt}
\put(0,0){\color{red}\line(1,1){140}}
\put(20,20){\color{blue}\line(1,0){20}}
\put(20,20){\color{cyan}\line(0,1){20}}
\put(40,40){\color{blue}\line(-1,0){20}}
\put(40,40){\color{cyan}\line(0,-1){20}}
\put(60,60){\color{blue}\line(1,0){60}}
\put(60,60){\color{cyan}\line(0,1){60}}
\put(120,120){\color{blue}\line(-1,0){60}}
\put(120,120){\color{cyan}\line(0,-1){60}}


\put(20,20){\circle*{3}}
\put(13,19){$b$}
\put(40,40){\circle*{3}}
\put(43,35){$d$}
\put(60,60){\circle*{3}}
\put(52,58){$b'$}
\put(120,120){\circle*{3}}
\put(123,116){$d'$}
\put(20,40){\circle*{3}}
\put(13,42){$a$}
\put(40,20){\circle*{3}}
\put(43,18){$c$}
\put(120,60){\circle*{3}}
\put(123,56){$y$}
\put(60,120){\circle*{3}}
\put(51,122){$x$}

\end{picture}

Now assume that $\|X\|\ge 4$. By Corollary~\ref{c:affine-Desarguesian}, the affine space $X$ is Desarguesian. Using the Kuratowski--Zorn Lemma, enlarge the independent set $\{a,b,c\}$ to a maximal independent set $M$ in $X$. By Corollary~\ref{c:affine-extend-auto}, there exists a unique affine automorphism $A:X\to X$ such that $A(a)=c$, $A(c)=a$ and $A(x)=x$ for all $x\in M\setminus\{a,c\}$. Since $Aabc=bca$, the uniqueness of the vertex $d$ of the parallelogram $d$ implies that $A(d)=d$ and hence $\{b,d\}\subseteq \Fix(A)$.
Consider the set $M'\defeq (M\setminus\{a,c\})\cup\{d\}$. The Exchange Property of the regular liner $X$ implies that the flat $H\defeq \overline{M'}$ is a hyperplane in $X$. Since the affine automorphism $A$ is a scalarity (see Theorem~\ref{t:D=>a=p=s} and Definition~\ref{d:affine-aut}), the inclusion $M'\subseteq \Fix(A)$ implies $H=\overline{M'}\subseteq \Fix(A)$, which means that the automorphism $A$ is hyperfixed. Since the restriction $AA{\restriction}_M$ of the affine  automorphism $AA$ is the identity map of the maximal independent set $M$, the uniqueness part of Corollary~\ref{c:affine-extend-auto} ensures that the automorphism $AA$ is the identity map of $X$ and hence $A$ is an involutive hyperfixed automorphism of $X$, witnessing that the Desarguesian affine space $X$ is hypersymmetric.
\end{proof}


\begin{definition}\label{d:a.inversive} An affine liner $X$ is called \index{inversive-dot affine space}\index{affine space!inversive-dot}\defterm{inversive-dot} if $X$ is right-inversive-dot and left-inversive-dot.
\end{definition}

\begin{theorem}\label{t:inversive-dot} For an affine space $X$, the following conditions are equivalent:
\begin{enumerate}
\item $X$ is inversive-dot;
\item $X$ is right-inversive-dot;
\item $X$ is hypersymmetric;
\item $X$ is shear;
\item  every ternar of $X$ is linear, distributive, commutative-plus and alternative-dot;
\item some ternar of $X$ is linear, distributive, associative-plus and alternative-dot;
\item every ternar of $X$ is Moufang;
\item every ternar of $X$ is Moufang-dot;
\item every ternar of $X$ is inversive-dot.
\end{enumerate}
If the order of $X$ is finite, then \textup{(1)--(9)} are equivalent to
\begin{enumerate}
\item[\textup{(10)}] some ternar of $X$ is linear, distributive, and Moufang.
\end{enumerate}
\end{theorem}

\begin{proof} If $\|X\|\ge 4$, then the affine space $X$ is Desarguesian. In this case all conditions (1)--(9) hold and thus they are equivalent (to show that Desarguesian affine spaces are shear and hypersymmetric, one should apply Corollary~\ref{c:affine-extend-auto}). So, assume that $\|X\|=3$, which means that $X$ is a Playfair plane.
\smallskip

The implication $(1)\Ra(2)$ is trivial, the implication $(2)\Ra(3)$ follows from Proposition~\ref{p:rinvdot=>hypersymmetric}, and $(3)\Ra(4)$ was proved in Theorem~\ref{t:hypersymmetric=>shear}. The equivalences $(4)\Leftrightarrow(5)\Leftrightarrow(6)$ were proved in  Theorem~\ref{t:shear<=>alternative}.
The implication $(5)\Ra(7)$ follows from Proposition~\ref{p:Bol<=>alternative} and Theorem~\ref{t:Bol<=>Moufang}; $(7)\Ra(8)$ is trivial, and $(8)\Ra(9)$ follows from Proposition~\ref{p:Moufang=>inversive}. The implication $(9)\Ra(1)$ follows from Theorems~\ref{t:left-inversive-dot<=>} and \ref{t:right-inversive-dot<=>}. 
\smallskip

$(10)\Ra(6)$ Assume that some ternar $R$ of $Y$ is linear, distributive, and Moufang. If $Y$ has finite order, then $|R|=|Y|_2<\w$ and hence the ternar $R$ is finite. By Theorem~\ref{t:Bol<=>Moufang} and Corollary~\ref{c:Bol=>alt}, the Moufang-dot ternar $R$ is alternative-dot. By Proposition~\ref{p:Moufang-plus-not-associative-plus}, the finite distributive Moufang-plus biloop $(R,+,\cdot)$ is associative-plus, and hence the ternar $R$ is linear, associative-plus and alternative-dot.
\smallskip

The implication $(6\wedge 7)\Ra(10)$ is trivial.
\end{proof}

By Theorem~\ref{t:inversive-dot}, every right-inversive-dot affine plane is inversive-dot.

\begin{problem} Is every left-inversive-dot affine plane inversive-dot?
\end{problem}

\section{Associative-dot affine spaces}

\begin{definition}\label{d:ass-dot} An affine liner $\Pi$ is defined to be \index{associative-dot affine space}\index{affine space!associative-dot}\defterm{associative-dot} if for every points\\ $o,a,b,c,d,\alpha,\beta,\gamma,\delta\in\Pi$ $$\big[(\Aline ab\parallel \Aline cd\parallel \Aline \alpha\beta\parallel \Aline \gamma\delta\nparallel \Aline ad\parallel \Aline bc\parallel \Aline \beta\gamma\parallel \Aline \alpha\delta)\;\wedge\;(o\in\Aline a\alpha\cap\Aline b\beta\cap\Aline c\gamma)\big]\;\Ra\;(\,o\in\Aline d\delta\,).$$
\end{definition}

\begin{picture}(100,150)(-120,-15)
\linethickness{0.64pt}

\put(40,20){\color{blue}\line(1,0){20}}
\put(60,60){\color{blue}\line(-1,0){20}}
\put(80,40){\color{blue}\line(1,0){40}}
\put(120,120){\color{blue}\line(-1,0){40}}

\put(0,0){\color{red}\line(3,1){120}}
\put(0,0){\color{red}\line(2,1){80}}
\put(0,0){\color{red}\line(1,1){120}}
\put(0,0){\color{red}\line(2,3){80}}
\put(60,20){\color{cyan}\line(0,1){40}}
\put(40,20){\color{cyan}\line(0,1){40}}
\put(80,40){\color{cyan}\line(0,1){80}}
\put(120,120){\color{cyan}\line(0,-1){80}}

\put(0,0){\color{red}\circle*{3}}
\put(-7,-6){$o$}
\put(40,20){\circle*{3}}
\put(32,20){$a$}
\put(60,20){\circle*{3}}
\put(60,13){$d$}
\put(40,60){\circle*{3}}
\put(33,59){$b$}
\put(60,60){\circle*{3}}
\put(62,55){$c$}
\put(80,40){\circle*{3}}
\put(71,39){$\alpha$}
\put(120,40){\circle*{3}}
\put(123,35){$\delta$}
\put(80,120){\circle*{3}}
\put(71,119){$\beta$}
\put(120,120){\circle*{3}}
\put(123,120){$\gamma$}
\end{picture}




\index[person]{Klingenberg}
\begin{theorem}[Klingenberg\footnote{
{\bf Wilhelm Paul Albert Klingenberg} (1924 -- 2010) was a German mathematician who worked on differential geometry and in particular on closed geodesics. Klingenberg was born in 1924 as the son of a Protestant minister. In 1934 the family moved to Berlin; he joined the Wehrmacht in 1941. After the war, he studied mathematics at the University of Kiel, where he finished his Ph.D. in 1950 with Karl-Heinrich Weise, with a thesis in affine differential geometry. After some time as an assistant of Friedrich Bachmann, he worked in the group of Wilhelm Blaschke at the University of Hamburg, where he defended his Habilitation in 1954. He then visited Sapienza University of Rome, working in the group of Francesco Severi and Beniamino Segre, after which he obtained a faculty position at the University of Göttingen (with Kurt Reidemeister), where he stayed until 1963. In 1954–55 Klingenberg spent a year at Indiana University Bloomington; during this time he also visited Marston Morse at Princeton University. In 1956–58 he accepted invitations to the Institute for Advanced Study in Princeton, New Jersey. In 1962 he visited the University of California, Berkeley as a guest of Shiing-Shen Chern, who he knew from his time in Hamburg. Later he became a full professor at the University of Mainz, and in 1966 a full professor at the University of Bonn, a position he kept till his retirement in 1989.}, 1955]\label{t:ass-dot<=>} For an affine space $X$, the following conditions are equivalent:\footnote{In his paper \cite{Wagner1981}, Wagner writes that Prof. Klingenberg told him that ``Ruth Moufang's initial reaction [to Theorem~\ref{t:ass-dot<=>}] was one of incredulity''.}\label{t:Klingenberg-ass}
\begin{enumerate}
\item $X$ is associative-dot;
\item $X$ is Desarguesian;
\item $X$ is dilation;
\item every ternar of $X$ is associative-dot; 
\item every ternar of $X$ is linear, distributive, and associative; 
\item some ternar of $X$ is linear, distributive, and associative.
\end{enumerate}
\end{theorem}

\begin{proof} If $\|X\|\ge 4$, then $X$ is Desarguesian, by Corollary~\ref{c:affine-Desarguesian}. Theorems~\ref{t:corps<=>} and  \ref{t:Des<=>dilation} imply that the conditions (1)--(9) hold and hence are equivalent. So, assume that $\|X\|=3$.
\smallskip

The implication $(1)\Ra(4)$ follows from Definition~\ref{d:ass-dot} and Proposition~\ref{p:ass-dot<=>}.
\smallskip

$(4)\Ra(5)$ Assume that every ternar of $X$ is associative-dot. Since every associative loop is inversive, every ternar of $X$ is inversive-dot. By Theorem~\ref{t:inversive-dot}, every ternar $R$ of $X$ is linear, distributive, and associative-plus. Since $R$ is associative-dot, $R$ is linear, distributive, and associative.
\smallskip

The implication $(5)\Ra(6)$ is trivial and $(6)\Ra(2)$  was proved in Theorem~\ref{t:corps<=>}. The equivalence $(2)\Leftrightarrow(3)$ was proved in Theorem~\ref{t:Des<=>dilation}, and $(2)\Ra(1)$ follows from Theorem~\ref{t:corps<=>}.


%
%
\end{proof}

\begin{exercise} Find an example of a based affine plane $X$ which is right-distributive and associative-dot, but the affine space $X$ is not associative-dot.
\smallskip

{\em Hint:} Consider the non-Desarguesian Thalesian affine plane of order $9$.
\end{exercise}

\section{Commutative-dot affine spaces}

\begin{definition}\label{d:a.commutative} An affine liner $X$ is called \index{commutative-dot affine space}\index{affine space!commutative-dot}\defterm{commutative-dot} if for every points $a,b,c,\alpha,\beta,\gamma\in X$ with
$$\Aline ab\parallel \Aline \alpha c\parallel \Aline \beta\gamma\nparallel \Aline bc\parallel \Aline a\gamma\parallel\Aline \alpha\beta,$$the lines $\Aline a\alpha,\Aline b\beta,\Aline c\gamma$ are paraconcurrent.
\end{definition}

\begin{picture}(50,100)(-150,0)

\put(0,0){\color{red}\line(1,1){90}}
\put(0,0){\color{orange}\line(3,2){90}}
\put(0,0){\color{magenta}\line(1,2){45}}
{\linethickness{0.7pt}
\put(45,30){\color{blue}\line(-1,0){15}}
\put(45,30){\color{cyan}\line(0,1){60}}
\put(30,30){\color{cyan}\line(0,1){30}}
\put(30,60){\color{blue}\line(1,0){60}}
\put(45,90){\color{blue}\line(1,0){45}}
\put(90,90){\color{cyan}\line(0,-1){30}}
\put(30,30){\color{red}\line(1,1){60}}
\put(46,30){\color{orange}\line(3,2){45}}
\put(30,60){\color{magenta}\line(1,2){15}}
}

\put(0,0){\color{red}\circle*{3}}
\put(0,0){\color{white}\circle*{2}}
\put(45,30){\circle*{3}}
\put(46,23){$a$}
\put(30,60){\circle*{3}}
\put(22,59){$c$}
\put(45,90){\circle*{3}}
\put(38,92){$\gamma$}
\put(90,60){\circle*{3}}
\put(93,57){$\alpha$}
\put(90,90){\circle*{3}}
\put(93,90){$\beta$}
\put(30,30){\circle*{3}}
\put(23,28){$b$}
\end{picture}

\begin{remark} The commutativity-dot axiom can be considered as a dual form of the Affine Pappus Axiom. The configuration appearing in this axiom is known in Geometry as the Thomsen's figure.
\end{remark}

\begin{lemma}\label{l:com-dot=>right-inv-dot} Every commutative-dot affine space $X$ is right-inversive-dot.
\end{lemma}

\begin{proof} To prove that $X$ is right-inversive-dot, take any line $D\subset X$ and points $o,a,c,\alpha,\gamma\in D$ and $b,d,\beta,\delta\in X\setminus D$ such that $$(\Aline ab\parallel \Aline cd\parallel \Aline \alpha\beta\parallel \Aline \gamma\delta\nparallel \Aline ad\parallel \Aline bc\parallel \Aline \beta\gamma\parallel \Aline \alpha\delta)\;\wedge\;(o\in \Aline b\beta).$$
We have to show that $o\in \Aline d\delta$. Observe that $o\in D\cap\Aline b\beta$ and $b\notin D$ implies $b\ne \beta$. 

If $a=\alpha$, then $\Aline ab=\Aline \alpha\beta$ and $o\in \Aline b\beta\cap D=\{a\}$. The parallelity $\Aline ad\parallel \Aline \alpha\delta$ implies $\Aline ad=\Aline \alpha\delta$. It follows from $b\ne\beta$ that $c\ne\gamma$ and $d\ne\delta$ and hence $o=a\in \Aline ad=\Aline d\delta$.
By analogy we can check that $o\in \Aline d\delta$ if $c=\gamma$. So, assume that $a\ne\alpha$ and $c\ne\gamma$. The choice of the points $b,\delta\notin D$ and $\Aline ab\nparallel \Aline bc\parallel \Aline \alpha\delta\nparallel \Aline \delta\gamma$ imply $a\ne c$ and $\alpha\ne\gamma$. If $c=\alpha$, then the commutativity-dot axiom implies that the lines $\Aline b\beta,\Aline a\gamma,\Aline d\delta$ are paraconcurrent and hence $\{o\}=\Aline b\beta\cap\Aline a\gamma\subset\Aline d\delta$.

So, assume that $c\ne \alpha$. Consider the unique points $x\in \Aline ad\cap\Aline c\beta$ and $y\in \Aline bc\cap\Aline \gamma\delta$. The commutativity-dot axiom ensures that the lines $\Aline b\beta,\Aline a\gamma,\Aline xy$ are paraconcurrent and hence $\{o\}=\Aline b\beta\cap\Aline a\gamma\subset\Aline xy$. Applying the commutativity-dot axiom to the points $x,y,d,\delta,c,\alpha$, we conclude that the lines $\Aline xy,\Aline d\delta,\Aline c\alpha$ are paraconcurrent and hence $\{o\}=\Aline xy\cap\Aline c\alpha\subset\Aline d\delta$.

\begin{picture}(200,140)(-150,-5)
\linethickness{0.75pt}
\put(0,0){\color{red}\line(1,1){120}}
\put(0,0){\color{magenta}\line(2,1){120}}
\put(0,0){\color{purple}\line(1,2){60}}
\put(20,20){\color{blue}\line(1,0){40}}
\put(20,20){\color{cyan}\line(0,1){20}}
\put(20,40){\color{blue}\line(1,0){100}}
\put(40,40){\color{cyan}\line(0,-1){20}}
\put(60,60){\color{blue}\line(1,0){60}}
\put(60,60){\color{cyan}\line(0,1){60}}
\put(120,120){\color{blue}\line(-1,0){60}}
\put(120,120){\color{cyan}\line(0,-1){80}}
\put(60,60){\color{cyan}\line(0,-1){40}}
\put(0,0){\color{orange}\line(3,1){120}}


\put(0,0){\color{red}\circle*{3}}
\put(-3,-9){$o$}
\put(20,20){\circle*{3}}
\put(13,20){$a$}
\put(40,40){\circle*{3}}
\put(36,42){$c$}
\put(60,60){\circle*{3}}
\put(52,60){$\alpha$}
\put(120,120){\circle*{3}}
\put(123,120){$\gamma$}
\put(20,40){\circle*{3}}
\put(13,40){$b$}
\put(40,20){\circle*{3}}
\put(33,22){$d$}
\put(120,60){\circle*{3}}
\put(123,58){$\delta$}
\put(60,120){\circle*{3}}
\put(51,118){$\beta$}
\put(60,20){\circle*{3}}
\put(58,10){$x$}
\put(120,40){\circle*{3}}
\put(123,38){$y$}
\end{picture}
\end{proof}

\begin{theorem}\label{t:commutative-dot<=>} For any affine space $X$, the following conditions are equivalent:
\begin{enumerate}
\item $X$ is commutative-dot;
\item $X$ is Pappian;
\item every ternar of $X$ is commutative-dot;
\item every ternar of $X$ is linear, distributive, associative, and commutative;
\item some ternar of $X$ is linear, distributive, associative, and commutative.
\end{enumerate}
\end{theorem} 

\begin{proof} The equivalence $(1)\Leftrightarrow(3)$ follows from Proposition~\ref{p:com-dot<=>}, and $(4)\Ra(3)$ is obvious.
\smallskip

$(1)\Ra(4)$ If $X$ is commutative-dot, then $X$ is right-inversive-dot, by Lemma~\ref{l:com-dot=>right-inv-dot}. By Theorem~\ref{t:inversive-dot}, every ternar $R$ of $X$ is linear, distributive, associative-plus, and alternative-dot. Since $X$ is commutative-dot, the ternar $R$ is also commutative-dot. Then $R$ is a commutative alternative division ring. By Corollary~\ref{c:div+com+alt=>ass}, the ring $R$ is associative and hence the ternar  $R$ is associative and commutative.
\smallskip

The implication $(4)\Ra(5)$ is trivial, and the implications $(5)\Ra(2)\Ra(4)$ follow from Theorem~\ref{t:corps<=>}. 
\end{proof}

\section{Dot loops of Boolean affine spaces}

In Theorems~\ref{t:Boolean<=>Boolean-plus} and \ref{t:Boolean<=>Boolean-puls} we have characterized Boolean affine liners via properties of the plus plus and puls operations. In this section we characterize Boolean affine liners via properties of the dot operation.

\begin{definition}An element $x$ of a unital magma $(X,\cdot)$ if called an \defterm{involutive} if $x\cdot x=e\ne x$, where $e$ is the identity element of the unital magma.
\end{definition}

\begin{theorem} An affine space $X$ is Boolean if and only if for every ternar $R$ of $X$, the dot loop $(R^*,\cdot)$ has no involutive elements.
\end{theorem}

\begin{proof} Assume that $X$ is Boolean. Given any ternar $R$ of $X$, find a plane $\Pi\subseteq X$ and an affine base $uow$ in $\Pi$ whose ternar $\Delta$ is isomorphic to the ternar $R$. Let $e\in\Delta$ be the diunit of the affine base $uew$. We claim that the dot loop $(\Delta^*,\cdot)$ has no involutive elements. Indeed, given any element $a\in \Delta^*\setminus\{e\}$, consider the point $y\defeq a\cdot a\in\Delta$ and the points $b,d$ with coordinates $ay$ and $ea$ in the based plane $(\Pi,uow)$. By the definition of the dot operation, the equality $y=a\cdot a$ implies $b\in \Aline od$. Assuming that $y=e$, we conclude that $abed$ is a parallelogram. Since the affine space $X$ is Boolean, the parallelogram $abed$ has disjoint diagonals $\Aline ae=\Delta$ and $\Aline bd=\Aline od$, which contradicts $o\in \Delta\in \Aline od$. This contradiction shows that $y\ne e$ and hence the element $a$ is not involutive. 
\smallskip

Assuming that the affine space $X$ is not Boolean, we can find a parallelogram $abcd$ whose diagonals $\Aline ac$ and $\Aline bd$ have a common point $o$. Consider the directions $\boldsymbol h=(\Aline bc)_\parallel =(\Aline ad)_{\parallel}$ and $\boldsymbol v=(\Aline ab)_{\parallel}=(\Aline cd)_{\parallel}$. Let $H\in\boldsymbol h$ and $V\in\boldsymbol v$ be unique lines containing the point $o$. By Proposition~\ref{p:Proclus-Postulate} (the Proclus Postulate), there exist unique points $u\in H\cap\Aline cd$ and $w\in V\cap\Aline bc$. Then $uow$ is an affine base in the plane $\Pi$ and the point $c$ coincides with the diunit $e$ of the affine base $uow$. The definition of the dot operation in the ternar $\Delta$ of the affine base $uow$ ensures that $a\cdot a=e$ and hence $a\ne e$ is an involutive element of the loop $(\Delta^*,\cdot)$.  
\end{proof}

\section{Affine spaces of finite order}

In this section we prove that for affine spaces of finite order many algebraic and geometric properties are equivalent.

\begin{theorem}\label{t:finite-shear<=>Pappian} For an affine space $X$ of finite order, the following conditions are equivalent:
\begin{enumerate}
\item $X$ is commutative-dot;
\item $X$ is associative-dot;
\item $X$ is inversive-dot;
\item $X$ is Pappian;
\item $X$ is Desarguesian;
\item $X$ is dilation;
\item $X$ is shear.
\end{enumerate}
If the order of $X$ is prime, then the conditions \textup{(1)--(7)} are equivalent to the conditions:
\begin{enumerate}
\item[\textup{(8)}]  $X$ is Thalesian;
\item[\textup{(9)}]  $X$ is translation;
\item[\textup{(10)}]  $X$ is $\partial$-translation.
\item[\textup{(11)}] $X$ is commutative-plus;
\item[\textup{(12)}]  $X$ is associative-plus;
\item[\textup{(13)}]  $X$ is inversive-plus;
\item[\textup{(14)}] $X$ is commutative-puls;
\item[\textup{(15)}]  $X$ is associative-puls;
\item[\textup{(16)}]  $X$ is inversive-puls;
\item[\textup{(17)}]  $X$ is left-inversive-puls.
\end{enumerate}
\end{theorem}

\begin{proof} Theorems~\ref{t:commutative-dot<=>}, \ref{t:ass-dot<=>}, \ref{t:inversive-dot} ensure $(1)\Leftrightarrow(4)\Ra(2)\Leftrightarrow(5)\Leftrightarrow(6)\Ra(3)\Leftrightarrow(7)$. It remains to prove $(7)\Ra(1)$. Assume that $X$ is shear. By Theorem~\ref{t:inversive-dot}, some ternar $R$ of $X$ is linear, distributive, associative-plus and alternative-dot, and hence the biloop of $R$ is an alternative division ring. Since $|R|=|X|_2<\w$, the division ring $R$ is finite. By the Artin--Zorn Theorem~\ref{t:Artin-Zorn}, the finite altenative division ring $R$ is a field and hence the ternar $R$ is commutative-dot. Therefore, the condition (1)--(7) are equivalent.
\smallskip

If the order of the affine space $X$ is prime, then the condition (5) is equivalent to the condition  (8)--(17), by Theorem~\ref{t:prime-affine<=>+}.
\end{proof} 

Although the following corollary can be proved by more elementary means, Theorem~\ref{t:finite-shear<=>Pappian} provides a short proof.

\begin{corollary}\label{c:4-Pappian} Every affine space $X$ of order $|X|_2\le 6$ is Pappian and Desarguesian.
\end{corollary}

\begin{proof} By Corollary~\ref{c:no6order}, $|X|_2\ne 6$ and hence $|X|_2\le 6$ implies $|X|_2\le 5$. So, it suffices to check that every affine space $X$ of order $|X|_2\le 5$ is Pappian. Given any plane $\Pi\subseteq X$ and an affine base $uow$ in $\Pi$, we shall prove that the ternar $\Delta$ of the based affine plane $(\Pi,uow)$ is commutative-dot. Let $e$ be the diunit of the affine base $uow$. 

If $|X|_2=3$, then $\Delta=\{o,e,a\}$ for some $a\in \Delta$. It is clear that the loop $(\Delta,\cdot)$ is commutative. 

If $|X|_2=4$, then there exist elements $a,b\in\Delta$ such that $\Delta=\{o,e,a,b\}$.  Since $(\Delta,\cdot)$ is a $0$-loop, $a,b\notin\{o,e\}$ implies $a\cdot b\notin\{a\cdot o,a\cdot e,o\cdot b,e\cdot b\}=\{o,a,b\}$ and hence $a\cdot b=e$. By analogy we can show that $b\cdot a=e$. This shows that the multiplication operation in the ternar $\Delta$ is commutative and hence the affine space $X$ is commutative-dot. 

Finally, assume that $|X|_2=5$ and find elements $a,b,c\in\Delta$ such that $\Delta=\{o,e,a,b,c\}$. If $e\notin\{a{\cdot}a,b{\cdot}b,c{\cdot}c\}$, then we lose no generality assuming that $a{\cdot} a=b$, $b{\cdot} b=c$ and $c{\cdot} c=a$. Then $a{\cdot} b\notin\{a,b,b{\cdot} b\}=\{a,b,c\}$ and hence $a{\cdot} b=e$. Also $a{\cdot} c\notin\{a,c,a\cdot a\}=\{a,c,b\}$ and again $a{\cdot} c=e=a{\cdot}b$, which is a contradiction showing that $e\in \{a{\cdot}a,b{\cdot}b,c{\cdot}c\}$. We lose no generality assuming that $a{\cdot}a=e$. In this case there exist exactly two ways to fill the multiplication table of the loop $(\Delta^*,\cdot)$:
$$
\begin{array}{c|cccc}
\cdot&e&a&b&c\\
\hline
e&e&a&b&c\\
a&a&e&c&b\\
b&b&c&e&a\\
c&c&b&a&e\
\end{array}\quad\mbox{and}\quad \begin{array}{c|cccc}
\cdot&e&a&b&c\\
\hline
e&e&a&b&c\\
a&a&e&c&b\\
b&b&c&a&e\\
c&c&b&e&a\
\end{array}.
$$
In both cases the multiplication is commutative and hence the ternar $\Delta$ is commutative-dot.
By Theorem~\ref{t:finite-shear<=>Pappian}, $X$ is Pappian, and by Theorem~\ref{t:Hessenberg-affine}, the Pappian affine space is Desarguesian.
\end{proof} 

\begin{corollary}\label{c:p5-Pappian} Every projective space $Y$ of order $|Y|_2-1\le 6$ is Pappian and Desarguesian.
\end{corollary}

\begin{proof} By Theorem~\ref{t:projPapp<=>}, the Pappianity of the space $Y$ will follow as soon as we check that for every hyperplane $H\subseteq Y$, the affine liner $X\defeq Y\setminus H$ is Pappian. If $|Y|_2\le 4$, then all lines in the affine liner have length $3$ and hence $X$ is Pappian vacuously. So, assume that $5\le |Y|_2$. By Proposition~\ref{p:projective-minus-hyperplane} and Corollary~\ref{c:procompletion-rank}, $X$ is an affine regular liner of order $4\le |X|_2=|Y|_2-1\le 6$ and rank $\|X\|=\|Y\|$. By Corollary~\ref{c:4-Pappian}, the affine space $X$ is Pappian. Therefore, the projective space $Y$ is Pappian, by Theorem~\ref{t:projPapp<=>}. By Hessenberg's Theorem~\ref{t:Hessenberg-proaffine}, the Pappian projective liner $Y$ is Desarguesian.
\end{proof}

\begin{Exercise} Prove that every affine space $X$ of order $|X|_2\le 8$ is Pappian.
\smallskip

{\em Hint:} This is difficult and cannot be done without a computer.
\end{Exercise}

\begin{remark}\label{rem:Des-algorithm} Theorems~\ref{t:Klingenberg-ass} and \ref{t:finite-shear<=>Pappian} suggest a fast algorithm for recognizing Desarguesian finite projective planes. Given a projective plane $X$ of order $n$, choose any quadrangle $uowe$ in $X$ and calculate the ternary operation $T:\Delta^3\to \Delta$ on the ternar $\Delta$ of the based projective plane $(X,uowe)$. Finding the ternar $(\Delta,T)$ has computational complexity $O(n^3)$ and memory usage $O(n^3)$.  The projective plane $(\Delta,T)$ is Desarguesian and Pappian if and only if its ternar $(\Delta,T)$ is linear, distributive, and associative, according to Theorems~Theorems~\ref{t:Klingenberg-ass} and \ref{t:finite-shear<=>Pappian}. Checking the linearity, distributivity and associativity of the ternar $(\Delta,T)$ has computation complexity $O(n^3)$. Therefore, this algorithm of recognizing Pappian finite projective planes has complexity $O(n^3)$. On the other hand, the direct verification of Desargues Axiom has computation complexity $O(n^{11})$ and the direct verification of Pappus Axiom has computation complexity $O(n^{10})$. \end{remark} 
\newpage

By Remarks~\ref{diag:plus-additive} and \ref{diag:puls-additive},  Theorems~\ref{t:inversive-dot}, \ref{t:ass-dot<=>}, \ref{t:commutative-dot<=>}, \ref{t:finite-shear<=>Pappian}, \ref{t:add-ass=>partialT}, for any affine space the following implications hold:
$$
\xymatrix{
\mbox{commutative-dot}\ar@{=>}[d]&\mbox{Pappian}\ar@{<=>}[l]\ar@{=>}[d]\\
\mbox{associative-dot}\ar@{=>}[d]&\mbox{Desarguesian}\ar@{<=>}[l]\ar@{=>}[d]\ar@/_10pt/_{\mbox{\scriptsize finite}}[u]\\ 
\mbox{inversive-dot}\ar@{=>}[d]&\mbox{shear}\ar@{<=>}[l]\ar@{=>}[d]\ar@/_10pt/_{\mbox{\scriptsize finite}}[u]\\
\mbox{invertible-dot}&\mbox{Thalesian}\ar@/_10pt/_{\mbox{\scriptsize prime}}[u]\\
\mbox{commutative-puls}\ar@{=>}[d]&\mbox{translation}\ar@{<=>}[u]\ar@{=>}[r]\ar@{=>}[l]\ar@{=>}[d]&\mbox{commutative-plus}\ar@{=>}[d]\\
\mbox{associative-puls}\ar@{=>}[d]\ar^{\mbox{\scriptsize prime-}}_{\mbox{\scriptsize power}}[r]&\mbox{$\partial$-translation}\ar@/_10pt/_{\mbox{\scriptsize prime}}[u]&\mbox{associative-plus}\ar@{=>}[d]\ar_{\mbox{\scriptsize prime-}}^{\mbox{\scriptsize power}}[l]\\
\mbox{inversive-puls}\ar@{=>}[d]\ar@/_10pt/_{\mbox{\scriptsize prime}}[u]&&\mbox{inversive-plus}\ar@{=>}[d]\ar@/_10pt/_{\mbox{\scriptsize prime}}[u]\\
\mbox{invertible-puls}&&\mbox{invertible-plus}\\
}
$$

\chapter{Elementary finite Moufang loops}\label{ch:elementary-Moufang}

In this section we explore the structure of elementary finite Moufang loops. The obtained information will be applied in Chapter~\ref{ch:invertible-add} devoted to invertible-add affine liners.

\section{Elementary loops}

We recall that a nonempty subset $A$ of a loop $(X,\cdot)$ is a \defterm{subloop} if for every elements $a,b\in A$, there exist elements $x,y,z\in A$ such that $a\cdot x=b=y\cdot a$ and $a\cdot b=z$. In this case $A$ endowed with the binary operation inherited from $X$ is a loop sharing with $X$ the same neutral element.

For every subset $A$ of a loop $X$, the intersection $\langle A\rangle$ of all subloops of $X$ that contain $A$ is called \defterm{the subloop} of $X$ generated by the set $A$. 

A subloop $B$ of a loop $X$ is called \defterm{$n$-generated} for some $n\in\IN$ if the smallest cardinality of a subset generating the loop $B$ equals $n$. In particular, the trivial subloop $\{e\}$ of $X$ is $0$-generated. 

\begin{definition}\label{d:elementary-loop} A loop $X$ is defined to be \defterm{elementary} if every $1$-generated subloop of $X$ is a cyclic group of prime order.
\end{definition}

Definition~\ref{d:elementary-loop} implies that every elementary loop is power-associative. 

\begin{definition} The \defterm{exponent} of a finite power-associative loop $X$ is the smallest number $n\in\IN$ such that for every $x\in X$ we have $x^n=e$, where $x^1\defeq x$ and $x^{k+1}\defeq x^k\cdot x$ for all $k\in\IN$.
\end{definition}

A number $n\in\IN$ is called \defterm{square-free} if $n$ is not divisible by the square $p^2$ of any prime number $p$. 

\begin{exercise} Prove that every finite elementry loop $X$ has square-free exponent.
\end{exercise}

\section{Normal subloops in loops}

For a subset $A$ of a loop $X$ and an element $x\in X$, let
$$xA\defeq\{x{\cdot}a:a\in A\}\quad\mbox{and}\quad Ax\defeq\{a{\cdot}x:a\in A\}.$$

A subset $A$ of a loop $X$ is called \defterm{normal} in $X$ if $xA=Ax$, $(xA)y=x(Ay)$ and $x(yA)=(x{\cdot}y)A$ for all $x,y\in X$. 

\begin{proposition}\label{p:normal-subloop} If $H$ is a normal subloop of a loop $X$, then 
\begin{enumerate}
\item for all $x\in X$ and $y\in xH$, we have $xH=yH$;
\item for all $x,y\in X$, either $xH=yH$ or $xH\cap yH=\varnothing$;
\item for all $x,y\in H$, $(xH)\cdot(yH)=(x{\cdot}y)H$;
\item for all $a,x,y\in X$, $a(xH)=a(yH)$ or $(xH)a=(yH)a$ implies $xH=yH$.
\end{enumerate}
\end{proposition}

\begin{proof} 1. For any $x\in X$ and $y\in xH$, we can find an element $h\in H$ such that $y=x{\cdot}h$. The normality of $H$ ensures that $yH=(x{\cdot}h)H=x(hH)=xH$ and hence $yH\subseteq xH$. 

Since $X$ is a loop, there exists an element $z\in X$ such that $y{\cdot} z=x$. Then $x=y{\cdot}z=(x{\cdot}h){\cdot}z\in (xH)z=x(Hz)$ and hence $e\in Hz$ and $e=uz$ for some $u\in H$. Since $H$ is a loop, there exists $v\in H$ such that $e=uv$. Since $X$ is a quasigroup, $uz=e=uv$ implies $z=v\in H$ and hence $x=y{\cdot}z\in yH$. Repeating the above argument, we conclude that $xH\subseteq yH$ and hence $xH=yH$.
\smallskip

2. Assuming that for $x,y\in X$, the cosets $xH$ and $yH$ have a common point $z$, we can apply the preceding item and conclude that $xH=zH=yH$.
\smallskip

3. Given any elements $x,y\in X$, we shall check that $(xH){\cdot}(yH)=(x{\cdot}y)H$. By the normality of $H$, we have
$$(xH){\cdot}(yH)=((xH)y)H=(x(Hy))H=(x(yH))H=((x{\cdot}y)H)H=(x{\cdot}y)(H{\cdot}H)=(x{\cdot}y)H.$$ 
\smallskip

4. Assume that $a(xH)=a(yH)$ for some $a,x,y\in X$. Then $a{\cdot}x=a{\cdot}(y{\cdot}h)$ for some $h$. Since loops are left cancellative, $x=y{\cdot}h\in yH$ and the first item implies $xH=yH$. If $(xH)a=(yH)a$, then $x{\cdot}a=(y{\cdot}\hbar){\cdot}a$ for some $\hbar\in H$ and the right cancellativity of $X$ ensures that $x\in y{\cdot}\hbar\in yH$. Applying the first statement, we conclude that $xH=yH$.
\end{proof}

For a normal subloop $H$ of a loop $X$, consider the set $X/H\defeq \{xH:x\in X\}$ and the quotient map $q:X\to X/H$, $q:x\mapsto xH$. The set $X/H$ is endowed with the binary operation $(xH){\cdot}(yH)=(x{\cdot}y)H$ making the function $q:X\to X/H$ a homomorphism of magmas. Proposition~\ref{p:normal-subloop} implies that the magma $X/H$ is a loop whose neutral element is $H=eH$. This implies the following external characterization of normal subloops.

\begin{proposition} A subloop $H$ of a loop $X$ is normal iff there exists a homomorphism $h:X\to Y$ to some loop $Y$ such that $H=q^{-1}(e)$, where $e$ is the neutral element of the loop $Y$.
\end{proposition}

\begin{proposition}\label{p:index2=>normal} Let $H$ be a subloop of an inversive loop $X$. If $|X|=2|H|$, then the subloop $H$ is normal in $X$.
\end{proposition}

\begin{proof} Since the left and right shifts in the loop $X$ are bijections of $X$, $x\in X\setminus H$, the sets $xH$ and $Hx$ have cardinality $|H|$. The inversivity of the loop $X$ ensures that the sets $xH$ and $Hx$ are disjoint with the set $H$, which implies $xH=X\setminus H=Hx$.  To complete the proof of the normality of the subloop $H$, we need to check that $x(Hy)=(xH)y$ and $x(yH)=(x{\cdot}y)H$ for all $x,y\in X$. Four cases are possible.
\smallskip

1. If $x,y\in H$, then $x(Hy)=xH=H=Hy=(xH)y$ and $x(yH)=xH=H=(x{\cdot}y)H$.
\smallskip

2. If $x\in H$ and $y\in X\setminus H$, then $x(Hy)=x(X\setminus H)=X\setminus xH=X\setminus H=Hy=(xH)y$ and $x(yH)=x(X\setminus H)=X\setminus H=(x{\cdot}y)H$ because $x{\cdot}y\in Hy\in X\setminus H$.

3. If $x\in X\setminus H$ and $y\in H$, then $x(Hy)=xH=X\setminus H=X\setminus Hy=(X\setminus H)y=(xH)y$ and $x(yH)=xH=X\setminus H=(x{\cdot}y)H$ because $x{\cdot}y\in xH=X\setminus H$.

4. If $x,y\in X\setminus H$, then $x^{-1}\in X\setminus H$ and hence  $x^{-1}H=X\setminus H=yH$ which implies $y\in x^{-1}H$ and $x{\cdot}y\in x(x^{-1}H)=H$, by the inversivity of the loop $X$. Also $xH=X\setminus H$ implies $x(X\setminus H)=X\setminus xH=X\setminus (X\setminus H)=H$. By analogy we can show that $(X\setminus H)y=H$. Then $x(Hy)=x(X\setminus H)=H=(X\setminus H)y=(xH)y$ and $x(yH)=x(X\setminus H)=H=(x{\cdot}y)H$.
\end{proof}

\begin{definition} A subset $S$ of a magma $X$  is \defterm{characteristic} if $\alpha[S]=S$ for every automorphism $\alpha$ of the magma $X$.
\end{definition}

\begin{exercise} Show that every characteristic subgroup of a group $X$ is normal in $Y$.
\end{exercise}

\begin{Exercise} Find an example of a finite Moufang loop $X$ containing a characteristic subloop, which is not normal in $X$.
\smallskip

{\em Hint:} Look at Remark~\ref{r:com-center}.
\end{Exercise}

\begin{exercise} Show that for every unital magma $X$ with neutral element $e$, the set of involutions $I\defeq\{x\in X\setminus\{e\}:x{\cdot}x=e\}$  is characteristic in $X$.
\end{exercise}

\begin{proposition}\label{p:I-normal} For every Moufang loop $X$, the set of  
involutions $I\defeq\{x\in X\setminus\{e\}:x{\cdot}x=e\}$ is normal in $X$.
\end{proposition}

\begin{proof} Observe that for every semi-automorphism $\alpha$ of $M$ and every $x\in I$, we have $\alpha(x)\ne\alpha(e)=e$ and 
$\alpha(x){\cdot}e=\alpha(x)=\alpha(x{\cdot}(x{\cdot}x))=\alpha(x){\cdot}(\alpha(x){\cdot}\alpha(x))$, which implies $\alpha(x){\cdot}\alpha(x)=e$, by the cancellativity of $X$. Therefore, $\alpha(x)\in I$ and $\alpha[I]=I$.

For every $a\in X$ let $L_a:X\to X$, $L_a:x\mapsto a{\cdot}x$, and $R_a:X\to X$, $R_a:x\mapsto x{\cdot}a$, be the left and right shifts on the loop $X$ by the element $a$. By Theorem~\ref{t:inner=>pseudo-semi-auto}, for every $a,b\in X$, the function $L_{a{\cdot}b}^{-1}L_aL_b=R_aL_b^{-1}R_a^{-1}L_b:X\to X$ is a semi-automorphism of $X$. 

Then $L_{a{\cdot}b}^{-1}L_aL_b[I]=I$ and $R_{a}L_b^{-1}R_a^{-1}L_b[I]=I$, which implies $$a(bI)=L_aL_b[I]=L_{a{\cdot}b}[I]=(a{\cdot}b)I\quad\mbox{and}\quad(bI)a^{-1}=R_{a}^{-1}L_b[I]=L_bR_a^{-1}[I]=b(Ia^{-1}),$$ witnessing that the set $I$ is normal in $X$.  
\end{proof}

\section{Finite simple Moufang loops}

In this section we survey some known results on the structure of finite simple Moufang loops.

\begin{definition} A loop $X$ is  \defterm{simple} if $X$ has exactly two normal subloops: $\{e\}$ and $X$. 
\end{definition}

In particular, a group $G$ is simple if it has exactly two normal subgroups: $\{e\}$ and $X$. The classification of finite simple groups is one of the main achievements of Mathematics in XX century. According to this classification, there exist 18 infinite series of finite simple groups: cyclic groups $C_p$  of prime order, alternating groups $A_n$ for $n\ge 5$, 16 infinite series of simple groups of Lie type and 27 sporadic groups. The proof of this theorem consists of about 5000 pages. The classification of finite simple groups is illustrated by the following diagram, borrowed from Wikimedia.




For us it is important to know that cyclic groups of prime order are simple and so are alternating groups $A_n$ for $n\ge 5$.

The classification of finite simple Moufang loops was completed soon after completing the classification of finite simple groups. Namely, Liebeck proved in 1987 that every non-associative finite simple Moufang loop is a Paige loop $M(\IF)$ over some finite field $\IF$. 
The Paige loops are non-associative counterparts of the projective linear groups $PGL(n,q)$ and are named in honour of \index[person]{Paige}Lowell J. Paige\footnote{{\bf Lowell J. Paige}  (December 10, 1919 – December 10, 2010) was an American mathematician renowned for his contributions to abstract algebra, particularly in the study of non-associative algebraic structures. He is best known for constructing the finite simple Moufang loops, now known as Paige loops. Paige earned his undergraduate degree from the University of Wyoming in 1939 and his PhD in 1947 from the University of Wisconsin-Madison under the supervision of Richard Hubert Bruck. He joined the faculty of the University of California, Los Angeles (UCLA) in 1947, later serving as chair of the mathematics department (1964–1968) and playing a key role in establishing the department’s reading room during the construction of the Mathematical Sciences Building. Paige subsequently held senior academic and administrative positions, including chairman of the UCLA Academic Senate, Dean of the Division of Physical Sciences (1968–1973), and Assistant Director of the National Science Foundation (1973–1975). He retired from UCLA in 1983, having also spent research periods at the Institute for Advanced Study and Yale University and held a Fulbright Award in Australia.} who proved in 1956 that such loops are simple. 
The Paige loops are defined with the help of Zorn matrices.

A \defterm{Zorn matrix} over a field $\IF$ is a $2\times 2$ matrix of the form
$$\left(\begin{array}{cc}
a&\vec{\boldsymbol \alpha}\\
\vec{\boldsymbol \beta}&d
\end{array}
\right)
$$
where $a,b\in\IF$ and $\vec{\boldsymbol \alpha},\vec{\boldsymbol\beta}\in\IF^3$. The addition and multiplication of two Zorn matrices is defined by the formulas:
$$\left(\begin{array}{cc}
a&\vec{\boldsymbol \alpha}\\
\vec{\boldsymbol\beta}&b
\end{array}
\right)+ \left(\begin{array}{cc}
c&\vec{\boldsymbol \gamma}\\
\vec{\boldsymbol\delta}&d\end{array}\right)\defeq \left(\begin{array}{cc}
a+c&\vec{\boldsymbol \alpha}+\vec{\boldsymbol \gamma}\\
\vec{\boldsymbol\beta}+\vec{\boldsymbol\delta}&b+d
\end{array}
\right)$$
and
$$\left(\begin{array}{cc}
a&\vec{\boldsymbol\alpha}\\
\vec{\boldsymbol\beta}&b
\end{array}
\right)\cdot \left(\begin{array}{cc}
c&\vec{\boldsymbol \gamma}\\
\vec{\boldsymbol\delta}&d\end{array}\right)\defeq \left(\begin{array}{cc}
a{\cdot}c+\vec{\boldsymbol\alpha}{\circ}\vec{\boldsymbol\delta}&a{\cdot}\vec{\boldsymbol\gamma}+d{\cdot}\vec{\boldsymbol \alpha}\\
c{\cdot}\vec{\boldsymbol\beta}+b{\cdot}\vec{\boldsymbol\delta}&\vec{\boldsymbol\beta}{\circ}\vec{\boldsymbol\gamma}+b{\cdot}d
\end{array}
\right)+\left(\begin{array}{cc}
0&-\vec{\boldsymbol\beta}{\times}\vec{\boldsymbol\delta}\\
\vec{\boldsymbol\alpha}{\times}\vec{\boldsymbol\gamma}&0
\end{array}
\right),$$
suggested by Max Zorn in 1933\footnote{Zorn, Max. Alternativkörper und quadratische systeme. (German) Abh. Math. Sem. Univ. Hamburg 9 (1933), no. 1, 395--402.}. In Zorn's formula, for two vectors $\vec{\boldsymbol\alpha}=(\alpha_0,\alpha_1,\alpha_2)$ and $\vec{\boldsymbol\beta}=(\beta_0,\beta_1,\beta_2)$ in $\IF^3$, 
$$
\begin{aligned}
\vec{\boldsymbol\alpha}{\circ}\vec{\boldsymbol\beta}&\defeq\alpha_0{\cdot}\beta_0+\alpha_1{\cdot}\beta_1+\alpha_2{\cdot}\beta_2\quad\mbox{and}\\
\vec{\boldsymbol\alpha}{\times}\vec{\boldsymbol\beta}&\defeq\left(\left|\begin{array}{cc}
\alpha_1&\alpha_2\\
\beta_1&\beta_2
\end{array}\right|,-\left|\begin{array}{cc}
\alpha_0&\alpha_2\\
\beta_0&\beta_2
\end{array}\right|,\left|\begin{array}{cc}
\alpha_0&\alpha_1\\
\beta_0&\beta_1
\end{array}\right|\right)\\
&=\big(\alpha_1{\cdot}\beta_2-\alpha_2{\cdot}\beta_1,\alpha_2{\cdot}\beta_0-\alpha_0{\cdot}\beta_2,\alpha_0{\cdot}\beta_1-\alpha_1{\cdot}\beta_0\big)
\end{aligned}
$$
are their standard scalar and vector products.

For a Zorn matrix 
$A=\left(\begin{array}{cc}
a&\vec{\boldsymbol \alpha}\\
\vec{\boldsymbol \beta}&d
\end{array}
\right)
$ its determinant is defined by the formula
$$\det(A)\defeq a{\cdot}d-\vec{\boldsymbol\alpha}\circ\vec{\boldsymbol\beta}.$$

\begin{propositions}[Paige, 1956]\label{p:Paige-det} $\det(A\cdot B)=\det(A)\cdot\det(B)$ for any Zorn matrices $A,B$ over a field $\IF$.
\end{propositions}


\begin{definition} A Zorn matrix $A$ over a field $\IF$ is called \defterm{special} if $\det(A)=1$.
\end{definition}

Proposition~\ref{p:Paige-det} implies that the set $SLZ(\IF)$ of all special Zorn matrices over a field $\IF$ is a magma with respect to the operation of multiplication.

\begin{propositions}[Paige, 1956]\label{p:Paige-loop} For every field $\IF$, the magma $SLZ(\IF)$ is a Moufang loop.
\end{propositions}


\begin{proposition}[Paige, 1956]\label{p:|SLZ(F)|} For every finite field $\IF$ of cardinality $q$, the Moufang loop $SLZ(\IF)$ has cardinality $q^3(q^4-1)$.
\end{proposition}


\begin{propositions}[Paige, 1956]\label{p:Paige-center} For every field $\IF$, the center $\mathcal Z(SLZ(F))$ of the Moufang loop $SLZ(\IF)$ coincides with the set
$$\left\{\left(\begin{array}{cc}
1&\vec{\boldsymbol 0}\\
\vec{\boldsymbol 0}&1\end{array}\right),\left(\begin{array}{cc}
-1&\vec{\boldsymbol 0}\\
\vec{\boldsymbol 0}&-1\end{array}\right)\right\},
$$
which is one-element if and only if the field $\IF$ has characteristic $2$.
\end{propositions}


The quotient loop $M(\IF)\defeq SLZ(\IF)/\mathcal Z(SLZ(\IF))$ is called the \defterm{Paige loop} over the field $\IF$. If the field $\IF$ is finite of cardinality $q$, then the Paige loop $M(\IF)$ is denoted by $M(q)$  (because all fields of cardinality $q$ are isomorphic).

\begin{theorems}[Paige, 1956] For every field $\IF$, the Paige loop $M(\IF)$ is a non-associative simple Moufang loop.
\end{theorems}

Propositions~\ref{p:|SLZ(F)|} and \ref{p:Paige-center} imply the following corollary on the cardinality of the Paige loops.

\begin{corollary} For every finite field $\IF$ of cardinality $q$, the Paige loop $M(q)$ has cardinality $|M(q)|=\frac1{d}q^3(q^4-1)$, where $d$ the the greatest common divisor of $2$ and $q-1$.
\end{corollary}

In 1987 \index[person]{Liebeck}Liebeck\footnote{{\bf Martin Liebeck} (born 23 September 1954) is a professor of Pure Mathematics at Imperial College London whose research interests include group theory and algebraic combinatorics. Martin Liebeck studied mathematics at the University of Oxford earning a First Class BA in 1976, an MSc in 1977, and a D.Phil. in 1979, with the Dissertation ``Finite Permutation Groups'' under supervision of Peter M. Neumann. In January 1991 he was appointed Professor at Imperial College London and became Head of the Pure Mathematics section there in 1997. Liebeck has published over 150 research articles and 10 books. His research interests include algebraic combinatorics, algebraic groups, permutation groups, and finite simple groups. 
He was elected Fellow of the American Mathematical Society (AMS) in 2019, and was awarded the London Mathematical Society’s Pólya Prize in 2020.} \cite{Liebeck1987} has proved that finite Paige loops are the only non-associative finite simple Moufang loops.

\begin{theorems}[Liebeck, 1987] Every finite non-associative simple Moufang loop is isomorphic to the Paige loop $M(\IF)$ over some finite field $\IF$.
\end{theorems} 

All possible orders of elements in the Paige loops were calculated by Grishkov and Zavarnitsine \cite[Lemma 22]{Grish-Zava2006} in 2006.

\begin{theorems}[Grishkov and Zavarnitsine, 2006]\label{t:ordersM(q)} Let $\IF$ be a finite field of cardinality $q$ and characteristic $p$, and $d$ be the greatest common divisor of the numbers $2$ and $p-1$. The set of orders of elements of the Paige loop $M(q)$ is the set of all divisors of the numbers $p$, $\frac1d(q-1)$ and $\frac1d(q+1)$.
\end{theorems}

Theorem~\ref{t:ordersM(q)} implies the following classification of elementary non-associative finite simple Moufang loops.

\begin{theorems}\label{t:elementary-simple<=>} A finite simple Moufang loop $X$ is elementary and not commutative if and only if $X$ is isomorphic to the alternating group $A_5$ of the Paige loop $M(q)$ for $q\in\{2,3,4,5\}$.
\end{theorems}

\begin{proof} The ``if'' part follows from Theorem~\ref{t:ordersM(q)} and known properties of the alternating group $A_5$. So, assume that $X$ is an elementary finite non-commutative simple Mougang loop. If $X$ is a group, then one can look at the known list of all finite simple groups and check that the unique elementary non-commutative group in this list is the alternating group $A_5$. If $X$ is not associative, then $X$ is isomorphic to the simple Paige loop $M(q)$ for some prime-power number $q$. Let $d$ be the greatest common divisor of the numbers $2$ and $q-1$. The elementarity of $X$ and Theorem~\ref{t:ordersM(q)} implies that the numbers $\frac1d(q-1)$, $\frac1d(d+1)$ are prime, which implies $q\in\{2,3,4,5\}$.  
\end{proof}

The following universal property of the Paige loop $M(2)$ was proved by Grishkov and Zavarnitsine in \cite{Grish-Zava2006}.

\begin{propositions}[Grishkov and Zavarnitsine, 2006]\label{p:M(2)-universal} For every $q\in \{2,3,4,5\}$, the simple Paige loop $M(q)$ contains an isomorphic copy of the Paige loop $M(2)$. 
\end{propositions}

\section{The nucleus and the center of a magma}

\begin{definition}\label{d:nucleus-center} For a magma $(X,\cdot)$ its
\begin{itemize}
\item \defterm{nucleus} is the set $\mathcal N(X)$ of all elements $x\in X$ such that\newline  $(x\cdot y)\cdot z=x\cdot (y\cdot z)$,  $(y\cdot x)\cdot z=y\cdot (x\cdot z)$, and $(y\cdot z)\cdot x=y\cdot (z\cdot x)$ for all $y,z\in X$;
\item \defterm{commutant} is the set $\mathcal C(X)$ of all elements $c\in X$ such that $x\cdot c=c\cdot x$ for all $x\in X$;
\item \defterm{center} is the set $\mathcal Z(X)\defeq \mathcal N(X)\cap \mathcal C(X)$. 
\end{itemize}
\end{definition}

Definition~\ref{d:nucleus-center} implies that for every magma $X$,
\begin{enumerate}
\item  $\mathcal Z(X)\subseteq \mathcal N(X)$;
\item $\mathcal N(X)$ is an associative submagma of $X$;
\item $\mathcal Z(X)$ is a commutative associative submagma of $X$.
\end{enumerate}

\begin{Exercise} Show that the commutant of a Moufang loop is a commutative subloop of the loop.
\end{Exercise}

\begin{Exercise} Find an example of a loop whose commutant is not a subloop.
\end{Exercise}

\begin{exercise} Show that for every loop $X$, 
\begin{enumerate}
\item the nucleus $\mathcal N(X)$ is a normal characteristic subgroup of $X$;
\item the commutant $\mathcal C(X)$ is a characteristic subset of $X$;
\item the center $\mathcal Z(X)$ is a normal characteristic subgroup of $X$.
\end{enumerate}
\end{exercise}

\begin{remark}\label{r:com-center} In \cite{Grish-Zava2021} Grishkov and Zavarnitsine have found a Moufang loop $X$ of order $|X|=3^8$ whose commutant $\mathcal C(X)$ is not normal in $X$.  
\end{remark}

The following important theorem was recently proved by \index[person]{Cs\"og\H o}Piroska Cs\"org\H o\footnote{{\bf Piroska Cs\"org\H o}, is a Hungarian mathematician specializing in algebra, particularly in the study of loops and non-associative algebraic structures. She has been affiliated with the Institute of Mathematics at Eszterh\'azy K\'aroly University in Eger, Hungary. Cs\"org\H o earned her Candidate of Science degree in 1989 and became a Doctor of the Hungarian Academy of Sciences in 2011.} in \cite{Csorgo2022}.

\begin{theorems}[Cs\"org\H o, 2022]\label{t:Scorgo} Every non-trivial Moufang loop of odd order has non-trivial nucleus.
\end{theorems}

\begin{definition} A magma $X$ is defined to be \defterm{prime-power} if $|X|=p^n$ for some prime number $p$ and some $n\in\w$.
\end{definition}

The following important theorem is a culmination of research of many mathematicians starting from Bruck, Glauberman, Wright and many others.

\begin{theorems}[Dr\'apal, 2023]\label{t:prime-power=>center} Every non-trivial prime-power Moufang loop has non-trivial center.
\end{theorems}

\begin{definition}\label{d:centrally-nilpotent} A loop $X$ is \defterm{centrally-nilpotent} if there exists a finite sequence $(X_k)_{k\in n}$ of normal subloops of $X$ such that $X_0=\{e\}$, $X_n=X$ and $X_k/X_{k-1}\subseteq \mathcal Z(X/X_{k-1})$ for all positive $k\le n$.
\end{definition}

Dr\'apal's Theorem~\ref{t:prime-power=>center} (combined with the Lagrange Theorem~\ref{t:Lagrange-Moufang} for Moufang loops)  implies the following corollary.

\begin{corollarys} Every prime-power Moufang loop is centrally-nipotent.
\end{corollarys}

\section{Solvable, nucleo-solvable and group-solvable loops}

In this section we introduce and study solvable loops and their two modifications: nucleo-solvable and group-solvable loops.

\begin{definition}\label{d:solvable} A loop $X$ is \defterm{solvable} if there exist $n\in \IN$ and an increasing sequence of normal subloops $(X_n)_{0\le k\le n}$ of $X$ such that $X_0=\{e\}$, $X_n=X$ and for every positive $k\in n$, the quotient loop $X_k/X_{k-1}$ is a commutative group.
\end{definition}

\begin{definition}\label{d:group-solvable} A loop $X$ is \defterm{group-solvable} if there exist $n\in \IN$ and an increasing sequence of normal subloops $(X_n)_{0\le k\le n}$ of $X$ such that $X_0=\{e\}$, $X_n=X$ and for every positive $k\in n$,  the quotient loop $X_k/X_{k-1}$ is a group.
\end{definition}

\begin{definition}\label{d:nucleo-solvable} A loop $X$ is \defterm{nucleo-solvable} if there exist $n\in \IN$ and an increasing sequence of normal subloops $(X_n)_{0\le k\le n}$ such that $X_0=\{e\}$, $X_n=X$ and for every positive $k\in n$, the quotient loop $X_k/X_{k-1}$ is a commutative subgroup of the nucleus $\mathcal N(X/X_{k-1})$ of the quotient loop $X/X_{k-1}$.
\end{definition}

Definitions~\ref{d:centrally-nilpotent} and \ref{d:solvable}--\ref{d:nucleo-solvable} imply that for every loop $X$, we have the implications:
$$\mbox{centrally-nilpotent}\Ra\mbox{nucleo-solvable}\Ra\mbox{solvable}\Ra\mbox{group-solvable}.$$

\begin{exercise} Find examples of finite Moufang loops showing the none of the above implications can be reversed.
\end{exercise}

One of crucial theorem in Group Theory is the Odd-Order Theorem, proved by Feit (1962) and Thomson (1963). 

\begin{theorems}[Odd Order Theorem]\label{t:Odd-Order} Every finite group of odd order is solvable.
\end{theorems}

The Odd Order Theorem combined with Cs\"org\H o Theorem~\ref{t:Scorgo} implies the following valuable fact.

\begin{theorems}\label{t:odd=>nucleo-solvable} Every Moufang loop of odd order is nucleo-solvable.
\end{theorems}

\begin{proof} The proof is by induction on the order of a Moufang loop. Assume that for some $n\in \IN$, all Moufang loops of odd order $<n$ are nucleo-solvable. Let $X$ be a Moufang loop of odd order $|X|=n$. By Cs\"org\H o Theorem~\ref{t:Scorgo}, the nucleus $G\defeq \mathcal N(X)$ of $X$ is not trivial, and by Lagrange Theorem~\ref{t:Lagrange-Moufang}, $|\mathcal N(X)|$ divides $|X|$ and hence is odd. By Odd Order Theorem~\ref{t:Odd-Order}, $G$ is a solvable group. 

Let $G_0=G$ and for every $n\in\w$, let $G_{n+1}$ be the subgroup of $G_n$ generated by commutators $xyx^{-1}y^{-1}$ of elements $x,y\in G_n$.  By induction we can prove that for every $n\in\w$ the subgroup $G_n$  is characteristic in $G$.
Since $G$ is solvable, there exists $m\in\IN$ such that $\{e\}=G_m\ne G_{m-1}$. Then $H\defeq G_{m-1}$ is a nontrivial characteristic Abelian subgroup of the group $G$. Since $H\subseteq G=\mathcal N(X)$, for every $a\in X$, the map $\alpha:G\to G$, $x\mapsto x^a\defeq a^{-1}xa$, is an automorphism of $G$. Since the subgroup $H$ is characteristic in $G$, $a^{-1}Ha=\alpha[H]=H$ and hence $aH=Ha$. On the other hand, $H\subseteq G=\mathcal N(X)$ implies $a(bH)=(a{\cdot}b)H$ for every $a,b\in X$.
Therefore, the Abelian subgroup $H$ of $G$ is normal in $X$. So, we can consider the quotient loop $X/H$. Since $|X/H|=|X|/|H|<n$ is odd, the inductive assumption ensures that the Moufang loop $X/H$ is nucleuo-solvable and hence contains an increasing sequence of normal subloops $(Y_k)_{1\le k\le l}$ such that $Y_1=\{H\}$, $Y_l=X/H$, and for every $k\in\{2,\dots,l\}$ the quotient subloop $Y_k/Y_{k-1}$ is an Abelian subgroup of the nucleus $\mathcal N(Y_l/Y_{k-1})$ of the loop $Y_l/Y_{k-1}$. Let $\pi:X\to X/H$ be the quotient homomorphism. Let $X_0\defeq \{e\}$ and for every positive $k\le l$, let $X_k\defeq\pi^{-1}[Y_k]$. Then $(X_k)_{0\le k\le l}$ is a sequence of normal subloops in $X$, witnessing that the Moufang loop $X$ is nucleo-solvable.
\end{proof}

We finish this section with the following important theorem of Grishkov and Zavarnitsine, describing the structure of finite Moufang loops.

\begin{theorems}[Grishkov, Zavarnitsine, 2009]\label{t:structure-Moufang} Every finite Moufang loop $X$ contains two normal subloops $G\subseteq P$ such that the loop $G$ is group-solvable, the quotient loop $P/G$ is isomorphic to the direct product of simple Paige loops, and $X/P$ is a Boolean group. 
\end{theorems}

\section{Cauchy, Largange, Sylow and Hall properties of loops}

\begin{definition} A finite loop $X$ is defined to have 
\begin{itemize}
\item the \defterm{Lagrange property} if the order of every subloop divides the order of $X$;
\item the \defterm{Cauchy property} if for every prime divisor $p$ of the order $|X|$ there exists a subloop $H\subseteq X$ of order $|H|=p$;
\item the \defterm{Sylow property} if for every maximal prime-power divisor $p^n$ of the order $|X|$, there exists a subloop $H\subseteq X$ of order $|H|=p^n$;
\item the \defterm{Hall property} if for every coprime numbers $n,m$ with $|X|=n\cdot m$, there exists a subloop $H\subseteq X$ of order $|H|=n$.
\end{itemize}
\end{definition}

%

The following important theorem was proved by Grishkov and Zavarnitsine \cite{Grish-Zava2005}, and independently by Gagola III and Hall \cite{GagolaHall2005} in 2005.

\begin{theorems}[Grishkov--Zavarnitsine and Hall--Gagola-III, 2005]\label{t:Lagrange-Moufang} Every finite Moufang loop has the Largange property.
\end{theorems}

The following characterization of Sylow Moufang loops follows from the results of Grishkov and Zavarnitsine \cite{Grish-Zava2009}.
 
\begin{theorems}[Grishkov--Zavarnitsine, 2009]\label{t:Sylow-Moufang} A finite Moufang loop has the Sylow property if and only if it has the Cauchy property if and only if it is group-solvable.
\end{theorems}

Moufang loops with the Hall property were characterized by Gagola III in \cite[6.10]{Gagola2011}. 

\begin{theorems}[Gagola III, 2011]\label{t:Hall-Moufang} A finite Moufang loop has the Hall property if and only if it is solvable.
\end{theorems}

Therefore, for every finite Moufang loop we have the following implications.
$$\xymatrix{
\mbox{of odd order}\ar@{=>}[d]&\mbox{prime-power}\ar@{=>}[d]\\
\mbox{nucleo-solvable}\ar@{=>}[d]&\mbox{centrally-nilpotent}\ar@{=>}[l]\ar@{=>}[d]\\
\mbox{solvable}\ar@{=>}[d]&\mbox{Hall property}\ar@{<=>}[l]\ar@{=>}[d]\\
\mbox{group-solvable}\ar@{=>}[d]&\mbox{Sylow property}\ar@{<=>}[l]\ar@{=>}[d]&\mbox{Cauchy property}\ar@{<=>}[l]\\
\mbox{arbitrary}\ar@{=>}[r]&\mbox{Lagrange property}
}
$$

The fact that every finite group has the Sylow property was first proved in 1872 by \index[person]{Sylow}Ludwig Sylow\footnote{{\bf Peter Ludvig Meidell Sylow} (1832--1918) was a Norwegian mathematician who proved foundational results in group theory. Sylow processed and further developed the innovative works of mathematicians Niels Henrik Abel and \'Evariste Galois in algebra. Sylow theorems and $p$-groups, known as Sylow subgroups, are fundamental in finite groups. By profession, Sylow was a teacher at the Frederiksborg Latin School for 40 years from 1858 to 1898, and then a professor at the University of Oslo for 20 years from 1898 to 1918. Despite the isolation in Frederiksborg, Sylow was an active member of the mathematical world. He wrote a total of approximately 25 mathematical and biographical works, corresponded with many of the leading mathematicians of the time, and was an able co-editor of Acta Mathematica from the journal's start in 1882. He was also elected into the Norwegian Academy of Science and Letters in 1868, a corresponding member of the Academy of Sciences in G\"ottingen and the University of Copenhagen awarded him an honorary doctorate in 1894.}. 

Let us recall that a group $X$ is called a \defterm{$p$-group} for a prime number $p$ if $|X|=p^n$ for some $n\in\IN$. A $p$-subgroup $S$ of a finite group $G$ is called a \defterm{Sylow $p$-subgroup} of $G$ if $|S|$ is the maximal power of $p$ dividing the order of the group $G$. 

\begin{theorem}[Sylow, 1872]\label{t:Sylow} Let $p$ be a prime number dividing the order $|G|$ of a finite group $G$.
\begin{enumerate}
\item Every $p$-subgroup of $G$ is contained in some Sylow $p$-subgroup of $G$.
\item The number of Sylow $p$-subgroups in $G$ equals $1$ modulo $p$.
\item Any two Sylow $p$-subgroups of $G$ are conjugated.
\end{enumerate}
\end{theorem}

\begin{proof} 1. The proof is by induction on the order of the group $G$. Assume that for some $n\in\IN$ any group of order $<n$ contains a $p$-Sylow subgroup. Let $G$ be a group of order $|G|=n$ and $p$ be any prime number dividing $n$. 

If $p$ divides the order of the center $\mathcal Z(G)$ of $G$, then the Abelian group $\mathcal Z(G)$ contains an element $z$ of order $p$ (by the Cauchy property of finite Abelian groups). Let $C$ be the cyclic subgroup generated by the element $z$, and $\pi:G\to G/C$ be the quotient homomorphism to the quotient group $G/C$. By the inductive assumption, the group $G/C$ contains a Sylow $p$-subgroup $S$. Then the preimage $\pi^{-1}[S]$ is a Sylow $p$-subgroup of $G$. 

Next, assume that $p$ does not divide $|\mathcal Z(G)|$. For every $x\in G$, let $x^G\defeq\{g^{-1}xg:g\in G\}$ be the conjugacy class of $G$. It is easy to see that two conjugacy classes either coincide or are disjoint. Moreover, $x\in\mathcal Z(G)$ if and only if $x^G=\{x\}$. Choose a set $A\subseteq G$ that has one-point intersection with each conjugacy class $x^G$. 
Then $$|G|=\bigcup_{a\in A}|a^G|=|\mathcal Z(G)|+\sum_{a\in A\setminus \mathcal Z(G)}|a^G|.$$Since $p$ divides $|G|$ but not does not divide $|\mathcal Z(G)|$, there exists $a\in A\setminus \mathcal Z(G)$ such that $|a^G|$ is not divided by $p$. Consider the stabilizer $H=\{g\in G:a^g=a\}$ of the element $a$ and observe that $|H|=|G|/|a^G|<|G|$. By the inductive assumption, the subgroup $H$ of $G$ contains a Sylow $p$-subgroup $S$. Since $|G|/|H|$ is not dividible by $p$, the subgroup $S$ is a Sylow $p$-subgroup in $G$.
\smallskip

2. Given any prime number $p$ that divides the order $|G|$ of the group $G$, consider the family $\Syl_p(G)$ of all Sylow $p$-subgroups of $G$. By the preceding item, the family $\Syl_p(G)$ is not empty. So, we can fix a Sylow $p$-subgroup $S\in\Syl_p(G)$. Consider the action $\alpha:S\times\Syl_p(G)\to\Syl_p(G)$, $(s,H)\mapsto H^s\defeq s^{-1}Hs$, of the group $S$ on the set $\Syl_p(G)$.
For a group $H\in\Syl_p(G)$, let $H^S\defeq\{H^s:s\in S\}$ be the orbit of the group $H$ under the action of the group $S$. We claim that $H^S=\{H\}$ if any only if $H=S$. Indeed, if $H\ne S$, then we can choose any element $s\in S\setminus H$. By the Lagrange Theorem, the element $s$ has order $p^k$ in the $p$-group $S$ for some $k\in\IN$. Assuming that $H^s=H$, we conclude that the subgroup of $G$ generated by $H\cup\{s\}$ has cardinality $|H|\cdot p^k$, which is not possible as $|H|$ is the largest power of $p$ that divides $|G|$. Therefore, $H^s\ne H$ and $H^S\ne\{H\}$. Taking into account that $|H^S|$ divides $|S|$, we conclude that $|H^S|$ is divided by $p$. Now we see that the set $\Syl_p(G)=\{S\}\cup\bigcup\{H^S:H\in\Syl_p(G)\setminus\{S\}\}$ is the  union of $S$-orbits all of which except for $S^S=\{S\}$ have orders divisible by $p$. This implies that $|\Syl_p(G)|=1 \mod p$. Moreover, for every $S$-invariant subset $\mathcal I\subseteq \Syl_p(X)$, we have 
$$|\mathcal I|=
\begin{cases}
1\mod p,&\mbox{ if $S\in\mathcal I$;}\\
0\mod p,&\mbox{ if $S\notin \mathcal I$}.
\end{cases}
$$
A subset $\I\subseteq\Syl_p(G)$ is called {\em $S$-invariant} if $H^s\in \mathcal I$ for every $H\in\I$ and $s\in S$. 
\smallskip

3. To see that all Sylow $p$-subgroups of $G$ are conjugated, fix any Sylow $p$-group $S\in\Syl_p(X)$ and consider the $S$-invariant sets $S^G=\{S^g:g\in G\}$ and $\mathcal I\defeq\Syl_p(G)\setminus S^G$ in $\Syl_p(G)$. By the preceding item $|S^G|=1\mod p$ and $|\mathcal I|=0\mod p$. Assuming that $S^G\ne\Syl_p(G)$, we can choose a Sylow $p$-subgroup $S'\in\Syl_p(G)\setminus S^G$ and observe that $\mathcal I$ is an $S'$-invariant subset of $\Syl_p(G)$. Since $S'\in\mathcal I$, $|\mathcal I|=1\mod p$, by the preceding idem, which contradicts $|\mathcal I|=0\mod p$. This contradiction shows that $\Syl_p(G)=S^G$ and hence all Sylow $p$-subgroups of $G$ are conjugated.
\end{proof}


The following proposition on the solvability of finite groups with cyclic Sylow $2$-subgroups was suggested to the author by Derek Holt and Peter Mueller\footnote{\tt https://mathoverflow.net/a/487000/61536}.

\begin{propositions}\label{p:Sylow-cyclic=>solvable} If Sylow $2$-subgroups of a finite group $G$ are cyclic, then $G$ is solvable.
\end{propositions}

\begin{proof} The proof is by induction on $|G|$. Assume that for some $n\in\IN$ we know that every group of order $<n$ with cyclic Sylow $2$-subgroups is solvable. Take any group $G$ of order $|G|=n$ with cyclic Sylow $2$-subgroups. If $n\le 3$, then the group $G$ is cyclic and hence solvable. So, assume that $n>3$. 

\begin{claim}\label{cl:2-Sylow=>normal} The group $G$ contains a non-trivial proper normal subgroup.
\end{claim}

\begin{proof} If $n$ is odd, this follows from the Odd Order Theorem. So, assume that $n$ is even and write $n$ as $n=2^sm$ for some $s\in\IN$ and some odd number $m$. Consider the injective homomorphism $h:G\to \Sym(G)$ assigning to every $g\in G$ the left shift $L_g:G\to G$, $L_g:x\mapsto g{\cdot}x$. Let $\Alt(G)$ be the group of even permutations of the set $G$. It is a normal subgroup of index 2 in the symmetric group $\Sym(G)$. By Sylow's Theorem, the group $G$ contains a Sylow $2$-subgoup $S$ of order $|S|=2^s$. By our assumption, the group $S$ is cyclic and hence contains an element $a$ of order $2^s$. Then the permutation $L_a:G\to G$ consists of $m$ cycles of length $2^s$ each. Every cycle of length $2^s$ is the product of $2^s-1$ transpositions. So, $L_a$ is a product of $m(2^s-1)$ transposition, which means that the permutation $L_a$ is odd and $L_a\notin \Alt(G)$. Then $N\defeq h^{-1}[\Alt(G)]$ is a non-trivial normal subgroup of index 2 in $G$.
\end{proof}

By Claim~\ref{cl:2-Sylow=>normal}, the group $G$ contains a non-trivial proper normal subgroup $N$. Let $S$ be a Sylow $2$-subgroup of $G$ and observe that $S\cap N$ is a Sylow $2$-subgroup of the group $N$, and $S/(S\cap N)$ is a Sylow $2$-subgroup of the quotient group $G/N$. Since $S$ is cyclic, so are the Sylow $2$-subgroups $S\cap N$ and $S/(S\cap N)$ in the groups $N$ and $S/N$, respectively. By the inductive hypotheses, the groups $N$ and $G/N$ are solvable, and then $G$ is solvable.   
\end{proof}
 
 
\begin{corollarys}\label{c:2Odd-Order} Every group $X$ of order $|X|=2m$ for some odd number $m$ is solvable. 
\end{corollarys}

For Moufang loops we have a weaker result.

\begin{theorems}\label{t:2m=>normal-m} Let $X$ be a diasociative loop with the Lagrange property. If $|X|=2m$ for some odd number $m$, then the set $H$ of all elements of odd order in $X$ is a normal subloop of order $|H|=m$ in $X$.
\end{theorems}

\begin{proof} Assuming that the set $H$ is not a subloop of $X$, we can find two elements $x,y\in H$ such that the element $x{\cdot}y$ has even order. 
Since $X$ is diassociative, the subloop $G$, generated by the elements $x,y$ is a group. Since the element $x{\cdot}y\in G$ has even order, the group $G$ has even order. The Lagrange property of $X$ ensures that $|G|$ divides $|X|=2m$ and hence $|G|=2d$ for some divisor $d$ of $m$. By Corollary~\ref{c:2Odd-Order}, the group $G$ is solvable and hence has the Hall property. Then $G$ contains a subgroup $D$ of order $|D|=d$. Proposition~\ref{p:index2=>normal} guarantees that the subgroup $D$ is normal in $G$. It follows that any element $z\in G\setminus D$ has even order, which implies that $x,y\in D$ and hence $x{\cdot} y\in D$. Since the order $|D|=d$ is odd, the element $x{\cdot}y$ has odd order that divides the odd number $d$. But this contradicts the choice of the elements $x,y$. This contradiction shows that $H$ is a submagma of $X$. Since every element of $X$ has finite order, the submagma $H$ of the inversive loop $X$ is a subloop of $X$.

By Proposition~\ref{p:Ali-Slaney}, the invertible loop $X$ contains an element $b$ of order $2$. We claim that $H\cup bH=X$. Indeed, take any element $x\in X\setminus H$ and consider the subgroup $B$ of $X$, generated by the elements $x$ and $b$. The group $B$ contains the element $b$ of order $2$ and hence $B$ has even order. By the Largange property of the loop $X$, $|B|$ divides $|X|=2m$ and hence $|B|=2k$ for some odd number $k$. By Corollary~\ref{c:2Odd-Order}, the group $B$ is solvable and by the Hall property, $B$ contains a subgroup $K$ of cardinality $|K|=k$. Since $k$ is odd, $K\subseteq H$ and $b,x\in B\setminus K$. Proposition~\ref{p:index2=>normal} guarantees that the subgroup $K$ is normal in $B$ and hence $x\in B\setminus K=bK\subseteq bH$, witnessing that $X=H\cup bH$ and hence $|X|=2|H|=2m$ and $|X|=m$. Proposition~\ref{p:index2=>normal} ensures that the subloop $H$ is normal in $X$.
\end{proof}

In spite of the fact that not every finite Moufang loop is Cauchy, we have the following important fact proved by Grishkov and Zavarnitsine in 
\cite{Grish-Zava2009}.

\begin{theorems}[Grishkov and Zavarnitsine, 2009]\label{t:23Sylow} Any finite Moufang loop contains Sylow $p$-subloops for every $p\in\{2,3\}$.
\end{theorems} 


 
\section{Semidirect products of loops}

\begin{definition} A loop $X$ is called a \defterm{product} of two subloops $K,H\subseteq X$ if $X=K{\cdot}H\defeq\{x{\cdot}y:x\in K,\;y\in H\}$ and $K\cap H=\{e\}$. The product $X=K{\cdot} H$ is called
\begin{itemize}
\item \defterm{direct} if both subloops $K$ and $H$ are normal in $X$;
\item \defterm{semidirect} if one of the subgroups $K$ or $H$ is normal in $X$.
\end{itemize}
\end{definition}

\begin{exercise} Find a finite inversive loop $X$ and two subloops $K,H\subseteq X$ such that $K\cap H=\{e\}$ and $X=K{\cdot}H$, but $|X|<|K|\cdot|H|$.
\smallskip

{\em Hint:} Look at Example~\ref{ex:Ali-Slaney}.
\end{exercise}

\begin{problem} Assume that a finite Moufang loop $X$ is a product of two subloops $K,H\subseteq X$ with $K\cap H=\{e\}$. Is $|X|=|K|\cdot|H|$?
\end{problem}

\begin{proposition} If a loop $X=A{\cdot}B$ is a direct product of two normal subloops $A,B$, then $a{\cdot}b=b{\cdot}a$ for all $a\in A$ and $b\in B$.
\end{proposition}

\begin{proof} Since $X$ is a loop, for every $c\in X$ the left and right shifts $$L_c:X\to X,\;L_c:x\mapsto c{\cdot}x,\quad\mbox{and}\quad R_c:X\to X,\;R_c:x\mapsto x{\cdot}c,$$are bijective maps of $X$.
Then $$
\begin{aligned}
\{a{\cdot}b\}&=L_aR_b[\{e\}]=L_aR_b[A\cap B]=L_aR_b[A]\cap L_aR_b[B]=(a(Ab))\cap (a(Bb))\\
&=((aA)b)\cap(aB)=(Ab)\cap (aB)=(bA)\cap (Ba)
\end{aligned}
$$ and 
$$
\begin{aligned}
\{b{\cdot}a\}&=L_bR_a[\{e\}]=L_bR_a[A\cap B]=L_bR_a[A]\cap L_bR_a[B]=(b(Aa))\cap(b(Ba))\\
&=(bA)\cap((bB)a)=(bA)\cap (Ba)=\{a{\cdot}b\}.
\end{aligned}
$$ 
\end{proof}

If $X=K\cdot H$ is a semidirect product of two subloops $K,H$ and the loop $K$ (resp. $H$) is normal in $X$, then we write $X=K\rtimes H$ (resp. $X=K\ltimes H$).

\begin{proposition}\label{p:=>cyclic-divides-kernel-1} If an elementary  diassociative loop $X$ is a semidirect product $K\rtimes H$ of two subloops $K,H$, then the order of any finite cyclic subgroup of $H$ divides $|K|$ or $|K\setminus\{e\}|$. 
\end{proposition}

\begin{proof} Let $C$ be any nontrivial cyclic subgroup of the loop $H$ and let $c\in C$ be its generator. Since $X$ is elementary, the cyclic group $C$ has prime order $|C|$. Assume the prime number $p\defeq |C|$ does not divide the cardinality $|K|$ of the group $K$. 

For every element $x\in K\setminus\{e\}$, consider the function $\gamma_x:C\to K\setminus\{e\}$, $\gamma_x:c\mapsto cxc^{-1}$.
We claim that the function $\gamma_c$ is injective. In the opposite case, we could find two distinct elements $a,b\in C$ such that $axa^{-1}=bxb^{-1}$. Since $C$ is a cyclic group of prime order $p$, any nonidentity element of $C$ has order $p$. In particular, the element $d\defeq b^{-1}a$ has order $p$. Since $X$ is diassociative, the elements $c,x$ are contained in a subgroup $G$ of the loop $X$. It is clear that $a,b\in C\subseteq G$. The equality $axa^{-1}=bxb^{-1}$ implies $dx=(b^{-1}a)x=x(b^{-1}a)=xd$. Then the elements $d$ and $x$ generate a commutative subgroup of the elementary group $G$ and hence the elements $d$ and $x$ have the same prime order $p$. By Proposition~\ref{p:di-Lagrange}, the order $p=|C|$ of the element $x$ divides $|K|$, which contradicts our assumption. This contradiction shows that the function $\gamma_x$ is injective. 

Next we show that for any $x,y\in K\setminus\{e\}$, the sets $\gamma_x[C]$ and $\gamma_y[C]$ either coincide or are disjoint. Assume that the intersection $\gamma_x[C]\cap\gamma_y[C]$ contains some element $z$. Then $c^nxc^{-n}=z=c^myc^{-m}$ for some $n,m\in\IZ$ and then $x=c^{m-n}yc^{n-m}$ and $$\gamma_x[C]=\{c^kxc^{-k}:k\in\IZ\}=\{c^{k+m-n}yc^{n-m-k}:k\in\IZ\}=\gamma_y[C].$$ Now we see that $\{\gamma_x[C]:x\in K\setminus\{e\}\}$ is a partition of the set $K\setminus\{e\}$ into subsets of cardinality $|C|=p$, which implies that $p$ divides $|K\setminus\{e\}|=|K|-1$. 
\end{proof}

\begin{definition} A subloop $H$ of a finite loop $X$ is called \defterm{Hall} if the order $|H|$ of $H$ is coprime with the index $|X/H|$ of $X$.
\end{definition}

The following fundamental theorem was proved by Issai Shur in 1904 (for solvable groups) and Zassenhaus in 1937 for the general case.

\begin{theorems}[Shur, 1904; Zassenhaus, 1937] For every normal Hall subgroup $H$ of a finite group $X$, there exists a subgroup $K\subseteq X$ such that $X=H\rtimes K$.
\end{theorems}

\section{Frobenius loops}

For a set $A\subseteq X$ of a loop $X$ and an element $x\in X$, let $$A^x\defeq x{\backslash}(Ax)\defeq\{x{\backslash}(a{\cdot}x):a\in A\}$$be the conjugation of the set $A$ in $X$. We recall that for two elements $a,b$ of a loop $X$, $a{\backslash}b$ is the unique element of $X$ such that $a{\cdot}(a{\backslash}b)=b$. If the loop $X$ is left-inversive, then  $a{\backslash}b=a^{-1}{\cdot}b$, where $a^{-1}$ is a unique element of $X$ such that $a^{-1}{\cdot}a=e=a{\cdot}a^{-1}$.

\begin{definition} A subloop $H$ of a loop $X$ is called \defterm{malnormal} if $H^x\cap H=\{e\}$ for every $x\in X\setminus H$.
\end{definition}

\begin{definition} A loop $X$ is defined to be \defterm{Frobenius} if $X$ is a semidirect product $K\rtimes M$ of a normal subloop $K$ and a malnormal subloop $M$ of $X$. The normal subloop $K$ is called the  \defterm{kernel} of $X$ and the malnormal subloop $M$ is called a \defterm{complement} of $X$. 
\end{definition}

\begin{example}\label{ex:Frobenius-affine} For every field $\IF$ and any subgroup $M$ of its multiplicative group $\IF^*$, the semidirect product $\IF\rtimes M$ endowed with the operation
$$(x,a)\cdot(y,b)=(x+ay,ab)$$is a Frobenius group with kernel $\IF\times\{0\}$ and complement $\{1\}\times M$.
\end{example}

\begin{exercise} Check that the group $\IF\rtimes M$ is indeed Frobenius.
\end{exercise}




\begin{Exercise} Show that the kernel $K$ of a Frobenius group $X$ is uniquely determined and the complement $M$ of $K$ in $X$ is isomorphic to the quotient group $X/K$.
\end{Exercise}

\begin{proposition}\label{p:=>Frobenius} If an elementary finite diassociative loop $X$ contains a normal subloop $K\subseteq X$ of prime index $p$ that does not divide $|K|$, then $X$ is a Frobenius loop with kernel $K$.
\end{proposition}

\begin{proof} Fix any element $c\in X\setminus K$. Since $X$ is elementary, $c$ is contained in a cyclic subgroup $C\subseteq X$ of prime order. Consider the quotient homomorphism $\pi: X\to X/K$ and observe that  $\pi[C]$ is a cyclic subgroup of prime order in the loop $X/K$ whose order is equal to the index $p$ of the normal subloop $K$ in $X$. The diassociativity of the loop $X$ implies the diassociativity of the quotient loop $X/K$. By Proposition~\ref{p:di-Lagrange}, the order $|\pi[C]|=|C|$ divides the prime number $p$ and hence $|C|=p$ and $\pi[C]=X/K$. Observe that $C\cap K$ is a subgroup of the cyclic group $C$. Since $c\in C\setminus K$ and $|C|$ is prime, $K\cap C=\{e\}$. The equality $\pi[C]=X/K$ implies $K\cdot C=X$. It remains to prove that the subgroup $C$ is malnormal in $X$. In the opposite case, we can find an element $x\in X\setminus C$ such that $C\cap C^x$ contains some element $y\ne\{e\}$. Since $y\in  C^x$, there exists an element $z\in C$ such that $y=x\backslash(z{\cdot}x)$ 
Then $x{\cdot}y=z{\cdot}x$ and $$\pi(x){\cdot}\pi(y)=\pi(x{\cdot}y)=\pi(z{\cdot}x)=\pi(z){\cdot}\pi(x)=\pi(x){\cdot}\pi(z)$$ and hence $\pi(y)=\pi(z)$ and $y=z$, by the commutativity of the cyclic group $\pi[C]=X/K$, and the bijectivity of the restriction $\pi{\restriction}_C:C\to X/K$. Then $x{\cdot}y=z{\cdot}x=y{\cdot}x$. Since $X$ is diassociative, the subloop $G$ generated by the commuting elements $x,y$ is a commutative group. Since $x\in X=K\cdot C$, there exist elements $u\in K$ and $v\in C$ such that $x=u{\cdot}v$. Since $C$ is a cyclic group of prime order, the element $y\in C\setminus\{e\}$ generates $C$ and hence $v\in C\subseteq G$ and $u=x{\cdot}v^{-1}\in G$. Since $X$ is elementary, the element $u\in K$ has a prime order $q$. By Proposition~\ref{p:di-Lagrange}, $q$ divides the order $|K|$ of $K$. Since $p$ does not divide $|K|$, $q\ne p$. Then the element $u{\cdot}y$ of the commutative group $G$ has non-prime order $pq$, which contradicts the elementarity of the loop $X$. This contradiction shows that the group $C$ is malnormal in $X$ and hence $X=K\rtimes C$ is a Frobenius loop.   
\end{proof}

\begin{proposition}\label{p:complement-divides} Let $X=K\rtimes M$ be a finite Frobenius loop with kernel $K$ and complement $M$. If $X$ is Moufang and $K\subseteq\mathcal N(X)$, then $|M|$ divides $|K|-1$. 
\end{proposition}

\begin{proof} We claim that for every $x\in K\setminus\{e\}$, the function $x^*:M\mapsto K$, $x^*:y\mapsto x^y\defeq y^{-1}xy$, is injective. Assume that $x^y=x^z$ for two elements $y,z\in M$.  It follows from $x\in K\subseteq \mathcal N(X)$ that $(x{\cdot}y){\cdot}z=x{\cdot}(y{\cdot}z)$. Applying the Moufang Theorem~\ref{t:Moufang-xyz}, we conclude that the elements $x,y,z$ generate a subgroup $G$ in $X$. By the associativity of the operation in the group $G$, the equality $y^{-1}xy=z^{-1}xz$ implies $x^{-1}zy^{-1}x=zy^{-1}\in M^x\cap M=\{e\}$ and hence $z=y$. Therefore, the set $x^M\defeq\{x^y:y\in M\}$ has cardinality $|M|$.

Next, we show that for any $x,y\in K\setminus\{e\}$ the orbits $x^M$ and $y^M$ either coincide or are disjoint. Assume that $x^M\cap y^M$ contains some point $z$ and find an element $s\in M$ such that $s^{-1}xs=z$. For every $t\in M$, the inclusion $x\in K\subseteq\mathcal N(X)$ and Moufang Theorem~\ref{t:Moufang-xyz} imply that the elements $x,s,t$ generate a subgroup $H$ in $X$. Then the associativity of the operation in the group $H$ ensures that $t^{-1}zt=t^{-1}(s^{-1}xs)t=(st)^{-1}x(st)\in x^M$ and hence $z^M\subseteq x^M$. On the other hand, $t^{-1}xt=t^{-1}(szs^{-1})t=(s^{-1}t)^{-1}z(s^{-1}t)\in z^M$. Therefore, $z^M=x^M$. By analogy we can prove that $z^M=y^M$ and hence $x^M=y^M$. 

 Therefore, $\{x^M:x\in K\setminus\{1\}\}$ is a partition of the set $K\setminus\{e\}$ into subsets of cardinality $|M|$ and hence $|M|$ divides $|K\setminus\{e\}|=|K|-1$.
\end{proof}

\begin{question} Let $X=K\rtimes M$ be a finite Frobenius Moufang loop with kernel $K$ and complement $M$. Does $|M|$ divide $|K|-1$?
\end{question}

The answer is affirmative if the loop $M$ is a cyclic group.

\begin{proposition}\label{p:Frobenius-prime-divides} Let $X=K\rtimes H$ be a diassociative Frobenius loop with kernel $K$ and complement $H$. The order of any cyclic subgroup of $H$ divides $|K|-1$.
\end{proposition}

\begin{proof} Let $C$ be a cyclic subgroup of $H$ and $c$ be a generator of $C$ and $p\defeq|C|$. Consider the function $\gamma:K\to K$, $\gamma:x\mapsto x^c\defeq c{\backslash}(c{\cdot}x)$. Let $\gamma^0:K\to K$ be the identity function of $X$ and for every $n\in\w$, let $\gamma^{n+1}\defeq\gamma^n\gamma$. Since the loop $X$ is diassociative, $\gamma^n(x)=c^{-n}xc^n=\gamma^{n+p}(x)$ for every $n\in\IN$ and $x\in K$. In particular, $\gamma^p(x)=x$, which implies that $\gamma$ is a bijective function of $K$. 

For every $x\in K\setminus \{e\}$, the malnormality of the group $C$ ensures that $x^{-1}Cx\cap C=\{e\}$. Then for every positive $k<p$, we have $x^{-1}c^{-k}x\ne c^{-k}$ and $\gamma^k(x)=c^{-k}xc^k\ne x$. This implies $\gamma^n(x)=\gamma^{n-k}(\gamma^k(x))\ne \gamma^k(x)$ for any numbers $0\le k<n<p$. Therefore, the set $x^C\defeq\{\gamma^k:0\le k<p\}=\{\gamma^k(x):k\in \IN\}$ has cardinality $p$. If for some elements $x,y\in K$ the sets $x^C$ and $y^C$ have a common point $z$, then $x^K=z^C=y^C$. Therefore, $\{x^C:x\in K\setminus\{e\}\}$ is a partition of the set $K\setminus \{e\}$ into pairwise disjoint sets of cardinality $p$ and hence $p$ divides $|K\setminus\{e\}|=|K|-1$.
\end{proof}

The following deep theorem was proved by Frobenius in 1901. The proof is not elementary and essentially uses Representation Theory of groups, see the post\footnote{\tt https://terrytao.wordpress.com/tag/frobenius-theorem/} in the blog of Terence Tao.

\begin{theorems}[Frobenius, 1901]\label{t:Frobenius-Kernel} If $H$ is a malnormal subgroup of a finite group $X$, then $X$ is a Frobenius group with kernel $$K\defeq\{e\}\cup(X\setminus\bigcup_{x\in X}H^x)$$ and complement $H$.
\end{theorems}

\section{Elementary nucleo-solvable Moufang loops}

In this section we describe the structure of elementary finite nucleo-solvable Moufang loops. We start with the following key lemma, adapted from the paper \cite{Wagner1981} of Wagner. 


\begin{lemma}\label{l:key-elementary-loop}
Let $X$ be a loop with neutral element $e$, $A$ be a commutative subgroup of $X$ and $B$ be a finite subloop of $X$ such that for all $x,y\in A$ and $a,b\in B$, the following conditions are satisfied:
\begin{enumerate}
\item $x^{a}\defeq a\backslash(x{\cdot}a)\in A$;
\item $(x\cdot y)^a=x^a\cdot y^a$;
\item $(x^a)^b=x^{a\cdot b}$.
\end{enumerate}
Let $\mathcal C$ be a family of subloops of $B$ such that $\bigcup\mathcal C=B$, $C\cap C'=\{e\}$ for every distinct loops $C,C'\in\mathcal C$, and for every $a\in A\setminus\{e\}$ and $C\in\mathcal C$ there exists an element $c\in C$ such that $c{\cdot} a\ne a{\cdot }c$.  Then $x^{|\mathcal C|-1}=e$ for every $x\in A$. 
\end{lemma}

\begin{proof} For any subloop $C\subseteq B$ and element $x\in A$, the condition (1) ensures that the finite product $x^C\defeq\prod_{c\in C}x^c$ is a well-defined element of the Abelian group $A$. For every element $a\in C$, the conditions (2) and (3) ensure that for 
$$(x^C)^a=\big(\prod_{c\in C}x^c)^a=\prod_{c\in C}(x^c)^a=\prod_{c\in C}x^{c\cdot a}=\prod_{c\in C}x^c=x^C,$$
which implies $a{\cdot}x^C=x^C\cdot a$. Then our assumptions guarantee that $x^B=e$ and $x^C=e$ for every $C\in\mathcal C$.
Taking into account that $\{C\setminus\{e\}:C\in\mathcal C\}$ is a disjoint cover of the set $B\setminus\{e\}$, we conclude that 
$$e=\prod_{C\in\mathcal C}x^C=\prod_{C\in\mathcal C}(x^e\cdot\prod_{c\in C\setminus\{e\}}x^c)=x^{|\mathcal C|}\cdot\prod_{c\in B\setminus\{e\}}x^c=x^{|\mathcal C|-1}\cdot\prod_{c\in B}x^c =x^{|\mathcal C|-1}x^B=x^{|\mathcal C|-1}\cdot e$$
and $x^{|\mathcal C|-1}=e$. 
\end{proof}

\begin{corollary}\label{c:cover-power}
Let $X$ be a Moufang loop with neutral element $e$, $A\subseteq \mathcal N(X)$ be a commutative normal subgroup of $X$, $B$ be a finite subloop of $X$, and $\mathcal C$ be a family of subloops of $B$ such that $\bigcup\mathcal C=B$, $C\cap C'=\{e\}$ for every distinct loops $C,C'\in\mathcal C$, and for every $a\in A\setminus\{e\}$ and $C\in\mathcal C$ there exists an element $c\in C$ such that $c{\cdot} a\ne a{\cdot }c$.  Then $x^{|\mathcal C|-1}=e$ for every $x\in A$. 
\end{corollary}

\begin{proof} The corollary will follow from Lemma~\ref{l:key-elementary-loop} as soon as we check that the subloops $A,B$ satisfy the conditions (1)--(3) of Lemma~\ref{l:key-elementary-loop}. 
\smallskip

1. Given any elements $x\in A$ and $a\in B$ we should prove that the element $x^a\defeq a\backslash(x{\cdot} a)$ belongs to the subgroup $A$. Since the Moufang loop $X$ is inversive, $a\backslash(x{\cdot}a)=a^{-1}{\cdot}(x{\cdot}a)$. The normality of $A$ in the inversive loop $X$ implies $Aa=aA$ and $a^{-1}{\cdot}(x{\cdot}a)\in a^{-1}{\cdot}Aa=a^{-1}\cdot aA=A$.
\smallskip

2. Given any elements $x,y\in A$ and $a \in B$, we should prove that $(x{\cdot}y)^a=x^a\cdot y^a$. Since $A\subseteq \mathcal N(X)$, $(x{\cdot}y){\cdot}a=x{\cdot}(y{\cdot}a)$. By Moufang Theorem~\ref{t:Moufang-xyz}, the elements $x,y,a$ generate a subgroup $G\subseteq X$. The associativity of the binary operation in the group $G$ ensures that $$(x{\cdot}y)^a=a^{-1}\cdot (x\cdot y)\cdot a=a^{-1}\cdot x\cdot a\cdot a^{-1}\cdot y\cdot a=x^a\cdot y^a.$$

3. Given any elements $x\in A$ and $a,b \in B$, we should prove that $x^{a{\cdot}b}=(x^a)^b$. Since $A\subseteq \mathcal N(X)$, $(x{\cdot}a){\cdot}b=x{\cdot}(a{\cdot}b)$. By Moufang Theorem~\ref{t:Moufang-xyz}, the elements $x,a,b$ generate a subgroup $G\subseteq X$. The associativity of the binary operation in the group $G$ ensures that 
$$x^{a\cdot b}=(a\cdot b)^{-1}\cdot x\cdot (a\cdot b)=b^{-1}\cdot a^{-1}\cdot x\cdot a\cdot b=b^{-1}\cdot x^a\cdot b=(x^a)^b.$$
\end{proof}

\begin{lemmas}\label{l:prime-not-prime} Let $A\subseteq \mathcal N(X)$ be a commutative normal subloop of an elementary finite Moufang loop $X$. If $|X/A|$ is a prime-power but not prime, then $X$ is prime-power.
\end{lemmas}

\begin{proof} Assume that $|X/A|$ is a prime-power but not prime. 
Then $|X/A|=p^n$ for some prime number $p$ and some $n\ge 2$. The elementarity of the loop $X$ implies the elementarity of its Abelian subgroup $A$. Then $|A|=q^m$ for some prime number $q$ and some $m\in\w$. If $m=0$, then the group $A$ is trivial and $|X|=|X/A|$ is a prime power. So, we assume that $m\ne 0$.

We claim that $q=p$. To derive a contradiction, assume that $q\ne p$. By Dr\'apal's Theorem~\ref{t:prime-power=>center}, the prime-power Moufang loop $X/A$ is solvable and so is the Moufang loop $X$. By Theorem~\ref{t:Sylow-Moufang}, the solvable Moufang loop $X$ has the Sylow property. Then there exists a subgroup $H\subseteq X$ of order $|H|=|X/A|=p^n$. By Dr\'apal's Theorem~\ref{t:prime-power=>center}, the center $\mathcal Z(H)$ is not trivial. So, we can fix an element $z\in \mathcal Z(H)\setminus\{e\}$. Let $C(z)$ be the cyclic subgroup of $H$, generated by $z$. Taking into account that $X$ is elementary and $|H|=p^n$, we conclude that $|C(z)|=p<p^n=|H|$ and hence there exists an element $y\in H\setminus C$. Let $B$ be a subloop generated by the elements $y$ and $z$. By Corollary~\ref{c:Moufang=>diassociative},  the Moufang loop $X$ is diassociative and hence $B$ is a subgroup of $H$. Since $z\in \mathcal Z(H)$, the subgroup $B$ is commutative. Being elementary, the commutative group $B$ has order $|B|=p^2$. Let $\mathcal C$ be the family of all non-trivial cyclic subgroups of $B$. Since the group $G$ is elementary, $\mathcal C$ is a cover of $B$ by cyclic groups of the same prime order $p$ such that $C\cap C'=\{e\}$ for any distinct groups $C,C'\in\mathcal C$. This implies that $|\mathcal C|=(|G\setminus\{e\}|)/(p-1)=p+1$. Observe that for every $a\in A\setminus \{e\}$ and $c\in B\setminus \{e\}$ we have $a{\cdot}c\ne c{\cdot}a$  (because  otherwise the element $a{\cdot}c$ would have non-prime order $pq$). Applying Corollary~\ref{c:cover-power}, we conclude that $x^p=x^{|\C|-1}=e$ for every $x\in A$, which contradicts $|A|=q^m$. This contradiction shows that $p=q$ and then $|X|=|A|\cdot |X/A|=q^m\cdot p^n=p^{m+n}$ is a prime power.
\end{proof}

\begin{lemmas}\label{l:prime-power-prime=>prime-power} Let $A\subseteq B$ be normal subloops of an elementary finite Moufang loop $X$. If $A$ is Abelian and $A\subseteq \mathcal N(X)$, $|X/B|$ is prime, and $|B/A|$ is prime-power, then $|B|$ is prime-power.
\end{lemmas}

\begin{proof} Assume that $|X/B|=p$ is prime. If $|B/A|$ is not prime, then $|B|$ is prime-power, by Lemma~\ref{l:prime-not-prime}. So, assume that $|B/A|=q$ is prime. If $p=q$, then $|X/A|=|X/B|\cdot|B/A|=p\cdot q=p^2$ is a prime-power but not prime. In this case, $X$ is prime-power, by Lemma~\ref{l:prime-not-prime}. By Lagrange Theorem~\ref{t:Lagrange-Moufang}, $|B|$ divides $|X|$ and hence is prime-power. So, assume that $p\ne q$. Since the Abelian group $A$ is elementary, $|A|=r^n$ for some prime number $r$ and some $n\in\w$. If $n=0$, then the group $A$ is trivial and $|B|=|B/A|$ is a prime-power. So, assume that $n>0$. 

We claim that $r=q$. To derive a contradiction, assume that $r\ne q$. Two cases should be considered. First we assume that $r\ne p$.  In this case the numbers $r^n$ and $pq$ are coprime. Applying the Hall Theorem~\ref{t:Hall-Moufang} to the solvable Moufang loop $X$, we can find a subloop $H\subseteq X$ of cardinality $|H|=pq$. Since $r$ does not divide $|H|=pq$, the intersection $A\cap H$ is trivial. 
Then the restriction $\pi{\restriction}_H:H\to X/A$ of the quotient homomorphism $\pi:X\to X/A$ is bijective and hence $|B\cap H|=q$, every element $x\in H\cap B\setminus \{e\}$ has order $q$ and every element $y\in H\setminus B$ has order $p$. Let $\mathcal C$ be the family of all cyclic subgroups of order $p$ in the group $H$. Since $p$ is prime, for every distinct groups $C,C'\in \mathcal C$ we have $C\cap C'=\{e\}=C\cap (B\cap H)$. Then $\C\cup\{B\cap H\}$ is a cover of the loop $H$ such that for every $a\in A\setminus\{e\}$ and $c\in H\setminus\{e\}$ we have $a\cdot c\ne c\cdot a$. Observe that $|\mathcal C\cup\{B\cap H\}|=1+(|H\setminus B|/(p-1))=1+(pq-q)/(p-1)=1+q$. Applying Corollary~\ref{c:cover-power}, we conclude that $x^q=e$ for all $x\in A$ which is impossible because every non-zero element of $A$ has prime order $r\ne q$. This contradiction shows that $r=p$.

Therefore, $|X|=r^nqp=p^{n+1}q$. Choose any elements $b\in B\setminus A$ and $c\in X\setminus B$. Since $|B/A|=q$ and $|X/B|=p$, the elementarity of $X$ ensures that $b^q=e=c^p$. Let $H$ be the subloop of $X$ generated by the elements $b,c$. Since the Moufang loop is diassociative, $H$ is a subgroup of $X$. The Lagrange Theorem~\ref{t:Lagrange-Moufang} ensures that $|H|$ divides $|X|=p^{n+1}q$ and hence $|H|=p^{m+1}q$ for some $m\in\w$. It follows that the cyclic group $Q$ generated by the element $b$ is Sylow. Observe that $H\cap A$ is a normal Abelian subgroup of the group $H$ and $H\cap B$ is a normal subgroup of $H$ containing the group $Q$. Then an element $x\in H$ has order $q$  if and only if $x\in H\cap (B\setminus A)$, which implies that the number of cyclic subgroups of order $q$ in $H$ is equal to $|H\cap B\setminus A|/(q-1)=(p^mq-p^m)/(q-1)=p^m$. By Sylow's Theorem~\ref{t:Sylow}, all cyclic subgroups of order $q$ in $H$ are conjugated. This implies that the subgroup $S=\{h\in H:h^{-1}Qh=Q\}$ has cardinality $|H|/p^m=pq$. It is clear that $Q\subseteq S$. We claim that $S\cap A=\{e\}$. Assuming that $S\cap A$ contains some element $a\ne e$, we conclude that $a^{-1}Qa=Q$ and hence $a^{-1}ba\in Q\cap Ab=\{b\}$. Then $ba=ab$ is an element of order $pq$, which contradicts the elementarity of the loop $X$. This contradiction shows that $S\cap A=\{e\}$ and hence $Q=S\cap B$. Let $\mathcal C$ be the family of all non-trivial cyclic subgroups of the group $S$. The elementarity of $X$ ensures that every cyclic group $C\in\mathcal C$ has order $p$ or $q$. More precisely, $Q$ is a unique cyclic group of order $q$ in $H$ and all other $(pq-q)/(p-1)=q$ cyclic groups have order $p$. Therefore, $|\mathcal C|=q+1$.

For every $x\in A$, consider the element $x^S\defeq\prod_{s\in S}x^s\in A$ where $x^s\defeq s^{-1}xs$. Also for every cyclic group $C\in\mathcal C$, consider the element $x^C\defeq\prod_{c\in C}x^c\in A$. Taking into account that $A\subseteq \mathcal N(X)$, we can repeat the argument of the proof of Lemma~\ref{l:key-elementary-loop} and show that $x^S=e$ and $x^Q=e$. We claim that $x^C=e$ for every $C\in\mathcal C\setminus\{Q\}$. Let $c$ be a generator of the cyclic group $C$. Then 
$$x^C=\prod_{k=0}^{p-1}c^kxc^{-k}=(xc)^{p}c^{-p}=(xc)^p=e$$because $xc,c\in cB$ and $(cB)^p=B$ in the group $X/B$.
Now we see that $$e=\prod_{C\in\mathcal C}x^C=\prod_{C\in\mathcal C}\prod_{c\in C}x^c=\prod_{C\in\mathcal C}x^e\prod_{c\in C\setminus\{e\}}x^c=x^{|\mathcal C|}\prod_{s\in S\setminus\{e\}}x^s=x^{|\mathcal C|-1}\cdot\prod_{s\in S}x^s=x^{|\mathcal C|-1}x^S=x^q,$$ which is impossible because $x\in A\setminus \{e\}$ has prime order $p\ne q$. This contradiction shows that $p=q$ and hence $|B|=|A|\cdot|B/A|=p^n\cdot q=p^{n+1}$ is prime-power.
\end{proof}

Now we can prove the main result of this section, describing the structure of elementary finite nucleo-solvable Moufang loops.

\begin{theorems}\label{t:nucleo-solvable} If $X$ is an elementary finite nucleo-solvable Moufang loop, then either $|X|$ is prime-power or else $X$ is a Frobenius loop whose kernel $K$ has prime-power order and complement has prime order $p=|X/K|$ that divides $|K|-1$.
\end{theorems}

\begin{proof} The proof is by induction on the cardinality of $X$. Assume that for some number $m\in\IN$ we know that every elementary finite nucleo-solvable Moufang loop $X$ of order $|X|<m$ is either  prime-power or else it contains a normal subloop $K$ of prime-power order and prime index. Take any elementary finite nucleo-solvable Moufang loop $X$ of order $|X|=m$. If $m=1$, then $|X|=1=2^0$ is prime-power and we are done. So, assume that $m>1$. Since $X$ is nucleo-solvable, there exists an increasing sequence of normal subloops $(X_k)_{k\le n}$ in $X$ such that $X_0=\{e\}$,  $X_n=X$, and for every positive $k\le n$, the quotient loop $X_k/X_{k-1}$ is a non-trivial commutative subgroup of the nucleus $\mathcal N(X/X_{k-1})$. In particular, $X_1$ is a non-trivial Abelian subgroup of the nucleus $\mathcal N(X)$ of $X$. Since $X$ is elementary, the Abelian group $X_1$ is elementary and hence prime-power. The increasing sequence of loops $(X_k/X_1)_{1\le k\le n}$ witnesses that the quotient loop $X/X_1$ is nucleo-solvable. By the inductive assumption, either $X/X_1$ is prime-power or else $X/X_1$ contains a normal subloop $H$ of prime-power order and prime index.

If $|X/X_1|$ is prime-power but not prime, then $|X|$ is prime-power,  by Lemma~\ref{l:prime-not-prime}. 

If $|X/X_1|$ is prime, then $X_1$ is a normal commutative subgroup of $X$ that has prime-power order and prime index $p$ in $X$. If $p$ divides the order of $X_1$, then $|X|$ is prime-power.

If $X/X_1$ is not prime-power, then by the inductive hypothesis, $X/X_1$ contains a normal subloop $H$ of prime-power order and prime index $p$. Let $\pi:X\to X_1$ be the quotient homomorphism. Then $K\defeq \pi^{-1}[H]$ is a normal subloop of prime index $p$ in $X$ such that $|K/X_1|=|H|$ is prime-power. By Lemma~\ref{l:prime-power-prime=>prime-power}, the subloop $K$ has prime-power order. 

In both cases, we have found a normal subloop $K\subseteq X$ of prime-power order and index $p$. If $p$ divides $|K|$, then $X$ is prime-power. If $p$ does not divide $|K|$, then by Proposition~\ref{p:=>Frobenius}, $X$ is a Frobenius loop with kernel $K$ and compement $C$ of prime order $p$. By Proposition~\ref{p:Frobenius-prime-divides}, $p$ divides $|K|-1$.
\end{proof}



Since every Moufang loop of odd order is nucleo-solvable, Theorem~\ref{t:nucleo-solvable} implies the following corollary.

\begin{corollarys}\label{c:elementary-odd<=>} If $X$ is an elementary Moufang loop of odd order, then either $|X|$ is prime-power or else $X$ is a Frobenius loop whose kernel $K$ has prime-power order and complement has  prime order $|X/K|$ that divides $|K|-1$.
\end{corollarys}

\section{The Chein extensions of loops}

For an invertible loop $X$, its \defterm{Chein extension} $M(X,2)$ is the magma $X\times\{-1,1\}$ endowed with the binary operation
$$(x,i)\cdot(y,j)=\big((x^j\cdot y^{ij})^j,ij).$$
The Chein extensions were introduced in 1974 by \index[person]{Chein}Orin Chein\footnote{{\bf Orin N.~Chein} is an American mathematician known for his work in nonassociative algebra, particularly in loop theory. He earned his Ph.D. in 1968 from New York University under Wilhelm Magnus for the Ph.D.  Thesis ``Some I-A Authomorphisms of a Free Group'' and spent most of his academic career at Temple University, where he taught for over 40 years. Chein is best known for the ``Chein extension'', a method for constructing nonassociative Moufang loops from groups. In addition to his research, Chein co-authored the book ``Problem Solving Through Recreational Mathematics'', promoting mathematics education through engaging problems. He received the Temple University Great Teacher Award in 1995, recognizing his excellence in teaching and mentorship.} \cite{Chein1974}.

\begin{proposition} The Chein extension of any (elementary) invertible loop is an (elementary) invertible loop.
\end{proposition}

\begin{proof} Let $X$ be an invertible loop with neutral element $e$. First we check that $(e,1)$ is the neutral element of the magma $M(X,2)$ and hence $M(X,2)$ is a unital magma.

Indeed,  for every $(x,k)\in M(X,2)$ we have
$(e,1)\cdot(x,k)=((e^k\cdot x^{1k})^k,1k)=((e\cdot x^k)^k,k)=((x^k)^k,k)=(x,k)$ and $(x,k)\cdot(e,1)=((x^1\cdot e^{k1})^1,k1)=(x,k)$.

Next, we show that $M(X,2)$ is a quasigroup. Given any elements $(a,n),(b,m)\in M(X,2)$ we should show that there exist unique elements $(x,i),(y,j)\in M(X,2)$ such that $(x,i)\cdot (a,n)=(b,m)$ and $(a,n)\cdot (y,j)=(b,m)$.

Since $\{-1,1\}$ is a multiplicative group, there exists a unique element $k\in\{-1,1\}$ such that $nk=m$. For $x\in X$, the equation $(x,k)\cdot(a,n)=(b,m)$ is equivalent to the equation $(x^n\cdot a^{kn})^n=b$, which is equivalent to $x^n\cdot a^m=b^n$ and has a unique solution $x=(b^n/a^m)^n$ in the invertible loop $X$. 
On the other hand, for every $y\in X$, the equation $(a,n)\cdot (y,k)=(b,m)$ is equivalent to $(a^k\cdot y^{nk})^k=b$, which is equivalent to $a^k\cdot y^m=b^k$ and has a unique solution $y=(a^k\backslash y^m)^k$ in the invertible loop $X$.

Therefore, $M(X,2)$ is a loop. To show that the loop $M(X,2)$ is invertible, fix any element $(x,k)\in M(X,2)$. If $k=1$, then the element $(x^{-1},1)$ is inverse to $(x,1)=(x,k)$ in the loop $M(X,2)$. If $k=-1$, then the element $(x,-1)$ in inverse to $(x,-1)=(x,k)$. Indeed,
$$(x,-1)\cdot(x,-1)=((x^{-1}\cdot x^{1})^{-1},1)=(e^{-1},1)=(e,1).$$
 
If the loop $X$ is elementary, then every element $x\in X\setminus\{e\}$ is contained in a cyclic group $C\subseteq X$ of prime order and hence the element $(x,1)$ is contained in the cyclic subgroup $C\times\{1\}$ of prime order in $M(X,2)$. On the other hand, for every $x\in X$, the element $(x,-1)$ is contained in the cyclic subgroup $\{(x,-1),(e,1)\}$ of order two. This witnesses that the loop $M(X,2)$ is elementary.
\end{proof}

\begin{proposition} The Chein extension of an inversive loop is an inversive loop.
\end{proposition}

\begin{proof} Let $X$ be an inversive loop. Proposition~\ref{p:inversive=>anti} ensures that $(x\cdot y)^{-1}=y^{-1}\cdot x^{-1}$ for all $x,y\in X$. To see that the loop $M(X,2)$ is inversive, take any elements $(x,i),(y,j)\in M(X,2)$. We should prove that $(x,i)^{-1}\cdot((x,i)\cdot(y,j))=(y,j)=((y,j)\cdot(x,i))\cdot(x,i)^{-1}$.
\smallskip

If $i=1$ and $j=1$, then the required equalities follow from the inversivity of the loop $X$.
\smallskip

If $i=1$ and $j=-1$, then 
$(x,i)^{-1}\cdot((x,i)\cdot (y,j))=(x,1)^{-1}\cdot((x,1)\cdot (y,-1))
=(x^{-1},1)\cdot((x^{-1}\cdot y^{-1})^{-1},-1)=
(x^{-1},1)\cdot(y{\cdot} x,-1)
=((x\cdot(y{\cdot} x)^{-1})^{-1},-1)=((x\cdot(x^{-1}\cdot y^{-1}))^{-1},-1)=(y,-1)=(y,j)
$ 
and 
$
((y,j)\cdot (x,i))\cdot(x,i)^{-1}=((y,-1)\cdot(x,1))\cdot(x,1)^{-1}=(y\cdot x^{-1},-1)\cdot (x^{-1},1)=((y\cdot x^{-1})\cdot x,-1)=(y,j),
$
by the inversivity of the loop $X$.
\smallskip

If $i=-1$ and $j=1$, then  $(x,i)^{-1}\cdot((x,i){\cdot} (y,j))=(x,-1)^{-1}\cdot((x,-1){\cdot} (y,1))
=(x,-1)\cdot(x{\cdot} y^{-1},-1)=
((x^{-1}{\cdot}(x{\cdot} y^{-1}))^{-1},1)=(y,j)
$ and $((y,j){\cdot }(x,i)){\cdot}(x,i)^{-1}=((y,1){\cdot}(x,-1)){\cdot}(x,-1)^{-1}=((y^{-1}{\cdot} x^{-1})^{-1},-1){\cdot}(x,-1)
=(x{\cdot} y,-1){\cdot} (x,-1)=(((x{\cdot} y)^{-1}{\cdot} x)^{-1},1)=(((y^{-1}{\cdot} x^{-1}){\cdot} x)^{-1},1)=(y,j)$.
\smallskip

If $i=-1$ and $j=-1$, then 
$((y,j)\cdot(x,i))\cdot(x,i)^{-1}=((y,-1)\cdot (x,-1))\cdot (x,-1)^{-1}=((y^{-1}{\cdot} x)^{-1},1)\cdot (x,-1)=(x^{-1}{\cdot} y,1)\cdot(x,-1)=(((x^{-1}{\cdot} y)^{-1}{\cdot} x^{-1})^{-1},-1)=((y^{-1}{\cdot }x){\cdot} x^{-1})^{-1},-1)=(y,j)$
and $((y,j)\cdot (x,i))\cdot(x,i)^{-1}=((y,-1)\cdot(x,-1))\cdot(x,-1)^{-1}=((y^{-1}\cdot x)^{-1},1)\cdot (x,-1)
=(x^{-1}\cdot y,1)\cdot (x,-1)=(((x^{-1}\cdot y)^{-1}\cdot x^{-1})^{-1},-1)=(((y^{-1}\cdot x)\cdot x^{-1})^{-1},-1)=(y,j)$, 
witnessing that the loop $M(X,2)$ is inversive.
\end{proof}

\begin{theorem}[Chein, 1974]\label{t:M(G,2)-Moufang<=>G-group} For an invertible loop $X$, the loop $M(X,2)$ is Moufang if and only if $X$ is a group.
\end{theorem}

\begin{proof} First assume that $X$ is a group. To prove that the loop $M(X,2)$ is Moufang, it suffices to check that $\big((x,i){\cdot}(y,j)\big){\cdot}\big((z,k){\cdot}(x,i)\big)=\big((x,i){\cdot}\big((y,j){\cdot}(z,k)\big)\big){\cdot}(x,i)$ for every elements $(x,i),(y,j),(z,k)$ of the loop $M(X,2)$.

It follows that
$$
\begin{aligned}
&\big((x,i){\cdot}(y,j)\big){\cdot}\big((z,k){\cdot}(x,i)\big)=\big((x^j\cdot y^{ij})^j,ij\big)\cdot\big((z^i\cdot x^{ki})^i,ki\big)\\
&=
\big(((x^j\cdot y^{ij})^{jki}\cdot (z^i\cdot x^{ki})^{iijki})^{ki},ijki\big)=
\big(((x^j\cdot y^{ij})^{jki}\cdot (z^i\cdot x^{ki})^{jki})^{ki},jk\big)
\end{aligned}$$
and
$$
\begin{aligned}
&\big((x,i){\cdot}\big((y,j){\cdot}(z,k)\big)\big){\cdot}(x,i)=
\big((x,i)\cdot((y^k\cdot z^{jk})^k,jk)\big)\cdot(x,i)\\
&=
\big(((x^{jk}\cdot (y^k\cdot z^{jk})^{kijk})^{jk},ijk\big)\cdot(x,i)=\big(((x^{jk}\cdot (y^k\cdot z^{jk})^{ij})^{jk},ijk\big)\cdot(x,i)\\
&=\big(((x^{jk}\cdot (y^k\cdot z^{jk})^{ij})^{jki}\cdot x^{ijki})^i,ijki\big)=\big(((x^{jk}\cdot (y^k\cdot z^{jk})^{ij})^{jki}\cdot x^{jk})^i,jk\big)
\end{aligned}
$$
If $k=1$, then
$$((x^j\cdot y^{ij})^{jki}\cdot (z^i\cdot x^{ki})^{jki})^{ki}=((x^j\cdot y^{ij})^{ji}\cdot (z^i\cdot x^{i})^{ji})^{i}$$ and
$$((x^{jk}\cdot (y^k\cdot z^{jk})^{ij})^{jki}\cdot x^{jk})^i=((x^{j}\cdot (y\cdot z^{j})^{ij})^{ji}\cdot x^{j})^i,$$
so, it suffices to check that 
\begin{equation}\label{eq:M2}
(x^j\cdot y^{ij})^{ji}\cdot (z^i\cdot x^{i})^{ji}=(x^{j}\cdot (y\cdot z^{j})^{ij})^{ji}\cdot x^{j}.
\end{equation}

If $ij=1$, then $i=j$ and the latter equality is equivalent to 
 $$(x^i\cdot y)\cdot (z^i\cdot x^{i})=(x^{i}\cdot (y\cdot z^{i}))\cdot x^{i},$$
 which holds because the group $X$ is a Moufang loop.
 
If $ij=-1$, then $i=-j$ and the equation (\ref{eq:M2}) is equivalent to
$$(x^{-i}\cdot y^{-1})^{-1}\cdot (z^i\cdot x^{i})^{-1}=(x^{-i}\cdot (y\cdot z^{-i})^{-1})^{-1}\cdot x^{-i}$$ and to 
$$(y\cdot x^i)\cdot (x^{-i}\cdot z^{-i})=((y\cdot z^{-i})\cdot x^i)\cdot x^{-i},$$which holds by the associativity of $X$.

If $k=-1$, then 
$$((x^j\cdot y^{ij})^{jki}\cdot (z^i\cdot x^{ki})^{jki})^{ki}=((x^j\cdot y^{ij})^{-ji}\cdot (z^i\cdot x^{-i})^{-ji})^{-i}$$ and
$$((x^{jk}\cdot (y^k\cdot z^{jk})^{ij})^{jki}\cdot x^{jk})^i=((x^{-j}\cdot (y^{-1}\cdot z^{-j})^{ij})^{-ji}\cdot x^{-j})^i,$$
so, it suffices to check that 
\begin{equation}\label{eq:M2a}
(x^j\cdot y^{ij})^{-ji}\cdot (z^i\cdot x^{-i})^{-ji}=
((x^{-j}\cdot (y^{-1}\cdot z^{-j})^{ij})^{-ji}\cdot x^{-j})^{-1}=
x^j\cdot (x^{-j}\cdot(y^{-1}\cdot z^{-j})^{ij})^{ji}.
\end{equation}
If $ij=1$, then $i=j$ and the equation (\ref{eq:M2a}) is equivalent to the equations
$$(x^i\cdot y)^{-1}\cdot (z^i\cdot x^{-i})^{-1}=
x^i\cdot (x^{-i}\cdot(y^{-1}\cdot z^{-i}))$$
and 
$$(y^{-1}\cdot x^{-i})\cdot(x^i\cdot z^{-i})=
x^i\cdot (x^{-i}\cdot(y^{-1}\cdot z^{-i})),$$
which hold by the associativity of $X$.

If $ij=-1$, then $j=-i$ and the equation (\ref{eq:M2a}) is equivalent to the equations
$$(x^{-i}\cdot y^{-1})\cdot (z^i\cdot x^{-i})=
x^{-i}\cdot (x^{i}\cdot(y^{-1}\cdot z^{i})^{-1})^{-1}=
x^{-i}\cdot ((y^{-1}\cdot z^{i})\cdot x^{-i}),
$$ 
which holds by the Moufang property of the group $X$. This completes the proof of the Moufang property of the loop $M(X,2)$.
\smallskip

Now assume that the loop $M(X,2)$ is Moufang. Then the loop $X$ also is Moufang, being isomorphic to the subloop $X\times\{1\}$ of the Moufang loop $M(X,2)$. To prove that the Moufang loop $X$ is associative, take any elements $x,y,z\in X$, put $a\defeq (x\cdot y)^{-1}$ and observe that 
$$
\begin{aligned}
&\big((y,-1)\cdot (z^{-1},1)\big)\cdot\big((a,1)\cdot(y,-1)\big)=(y\cdot z,-1)\cdot((a^{-1}\cdot y^{-1})^{-1},-1)=(y\cdot z,-1)\cdot(y\cdot a,-1)\\
&=(((y\cdot z)^{-1}\cdot(y\cdot a))^{-1},1)=((y\cdot a)^{-1}\cdot (y\cdot z),1)=\big((a^{-1}\cdot y^{-1})\cdot (y\cdot z),1\big)\\
&=\big((x{\cdot} y)\cdot y^{-1})\cdot (y\cdot z),1\big)=(x\cdot (y\cdot z),1).
\end{aligned}
$$ 
By the Moufang identity,
$$\begin{aligned}
\big(x\cdot (y\cdot z),1\big)&=\big((y,-1)\cdot (z^{-1},1)\big)\cdot\big((a,1)\cdot(y,-1)\big)=\big((y,-1)\cdot ((z^{-1},1)\cdot(a,1))\big)\cdot(y,-1)\\
&=\big((y,-1)\cdot ((z^{-1}\cdot a,1)\big)\cdot(y,-1)=
(y\cdot (z^{-1}\cdot a)^{-1},-1))\cdot(y,-1)\\
&=(y\cdot (a^{-1}\cdot z),-1))\cdot(y,-1)=\big(((y\cdot (a^{-1}\cdot z))^{-1}\cdot y)^{-1},1\big)\\
&=\big(y^{-1}\cdot (y\cdot (a^{-1}\cdot z)),1\big)=(a^{-1}\cdot z,1)=\big((x\cdot y)\cdot z,1\big),
\end{aligned}
$$
witnessing that the loop $X$ is associative and hence is a group.
\end{proof} 

\begin{theorem}\label{t:M(X,2)-associative<=>X-commutative} For an invertible loop $X$, the loop $M(X,2)$ is associative if and only if $X$ is a commutative group.
\end{theorem}

\begin{proof} If  $X$ is a commutative group, then the binary operation in $M(X,2)$ can be defined by a simpler formula:
$$(x,i)\cdot (y,j)\defeq ((x^j\cdot y^{ij})^j,ij)=(x\cdot y^i,ij).$$
To show that the loop $M(X,2)$ is associative, take any elements $(x,i),(y,j),(z,k)$ of $M(X,2)$ and observe that
$$
\big((x,i)\cdot(y,j)\big)\cdot(z,k)=
(x\cdot y^i,ij)\cdot (z,k)=(x\cdot y^i\cdot z^{ij},ijk)
$$
and
$$(x,i)\cdot\big((y,j)\cdot (z,k)\big)=
(x,i)\cdot(y\cdot z^j,jk)=(x\cdot(y\cdot z^j)^i,ijk)=(x\cdot y^i\cdot z^{ji},ijk),$$
witnessing that $\big((x,i)\cdot (y,j)\big)\cdot(z,k)=(x,i)\cdot\big((y,j)\cdot (z,k)\big)$.
\smallskip

Now assume that the loop $M(X,2)$ is associative. By Theorem~\ref{t:M(G,2)-Moufang<=>G-group}, the loop $X$ is a group. Moreover, for all $x,y\in X$, we have
$$((x,-1)\cdot(y^{-1},-1))\cdot(e,-1)=((x^{-1}\cdot y^{-1})^{-1},1)\cdot (e,-1)=(y\cdot x,1)\cdot (e,-1)=(y\cdot x,-1)$$
and $$(x,-1)\cdot ((y^{-1},-1)\cdot (e,-1))=(x,-1)\cdot (y^{-1},1)=(x\cdot y,-1).$$
Since $M(X,2)$ is associative,
$$(x\cdot y,-1)=(x,-1)\cdot ((y^{-1},-1)\cdot (e,-1))=((x,-1)\cdot (y^{-1},-1))\cdot(e,-1)=(y\cdot x,-1),$$
witnessing that the group $X$ is commutative.
\end{proof}

An unitary magma $X$ is called \defterm{Boolean} if $x\cdot x=e$ for all $x\in X$. Observe that an invertible loop is Boolean if and only if $x^{-1}=x$ for all $x\in X$. This characterization and the definition of the Chain extension implies the following characterization. 

\begin{proposition}\label{p:Boolean<=>Chein-Boolean} An invertible loop $X$ is Boolean if and only if its Chein extension $M(X,2)$ is Boolean.
\end{proposition}

\begin{theorem}\label{t:M(X,2)-commutative<=>X-Boolean}For an invertible loop $X$, the following conditions are equivalent:
\begin{enumerate}
\item the loop $M(X,2)$ is commutative;
\item the loop $X$ is commutative and Boolean;
\item the loop $M(X,2)$ is commutative and Boolean.
\end{enumerate}
\end{theorem}

\begin{proof} $(1)\Ra(2)$ Assume that the loop $M(X,2)$ is commutative. Observe that for every $x\in X$, we have
$(x,1)\cdot(e,-1)=((x^{-1}\cdot e^{-1})^{-1},-1)=(x,-1)$ and $(e,-1)\cdot(x,1)=(x^{-1},-1)$. The commutativity of the loop $M(X,2)$ ensures that
$(x,-1)=(x,1)\cdot(e,-1)=(e,-1)\cdot(x,1)=(x^{-1},-1)$ and hence $x=x^{-1}$, which means that the invertible loop $X$ is Boolean.  
\smallskip

$(2)\Ra(3)$ Assume that the loop $X$ is commutative and Boolean. Then the loop $M(X,2)$ is Boolean, by Proposition~\ref{p:Boolean<=>Chein-Boolean}. Given any elements $(x,i),(y,j)$ of $M(X,2)$, observe that
$$
(x,i)\cdot(y,j)=((x^j\cdot y^{ij})^j,ij)=(x\cdot y,ij)=(y\cdot x,ji)=(y,j)\cdot (x,i),$$witnessing that the loop $M(X,2)$ is commutative.
\smallskip

The implication $(3)\Ra(1)$ is trivial.
\end{proof}

\begin{theorem}\label{t:Chein-index} Let $X$ is a Moufang loop and let $G$ be a proper subloop of $X$ containing all elements of order $>2$ in $X$. If $X$ is not Boolean, then $G$ is a normal subgroup of index $|X/G|\in\{2,4\}$ in $X$. Moreover,
\begin{enumerate}
\item if $G$ has index $2$ in $X$, then $X$ is isomorphic to the Chein extension $M(G,2)$ of $G$;
\item if $G$ has index $4$ in $X$, then $G$ is commutative, and $X$ is isomorphic to the Chein extension $M(M(G,2),2)$ of the Chein extension $M(G,2)$ of $G$.
\end{enumerate}
\end{theorem}

\begin{proof} By Corollary~\ref{c:Moufang=>diassociative}, the Moufang loop $X$ is diassociative and hence power-associative. So, every element of $X$ generates a cyclic subgroup whose cardinality is equal to the order of the element in $X$. Since the loop $X$ is not Boolean, some element of $X$ has order $>2$ and belongs to the subloop $G$. Since $G$ contains all elements of order $>2$, every element $b\in X\setminus G$ has order $2$.  

\begin{claim}\label{cl:GGx=0} For every $x\in X\setminus G$, $Gx\cap G=\varnothing=xG\cap G$.
\end{claim}

\begin{proof} Assuming that $Gx\cap G$ contains some element $a$, we conclude that $a=b{\cdot}x$ for some $b\in G$ and hence $x=b^{-1}{\cdot}(b{\cdot} x)=b^{-1}{\cdot} a\in G$, by the inversivity of the Moufang loop $X$. But $x\in G$ contradicts the choice of $x$. This contradiction shows that $Gx\cap G=\varnothing$.
By analogy we can show that $xG\cap G=\varnothing$.
\end{proof}

\begin{claim}\label{cl:xg=gx} For every $g\in G$ and $x\in X\setminus G$ we have $x{\cdot}g=g^{-1}{\cdot}x$.
\end{claim}

\begin{proof} Claim~\ref{cl:GGx=0} ensures that $xG\cap G=\varnothing$ and hence the element $x{\cdot}g\in xG\subseteq X\setminus G\subseteq X\setminus T$ has order 2. Then $(x{\cdot}g)\cdot(x{\cdot}g)=e$, which implies $x\cdot g\cdot x=g^{-1}$ and $x{\cdot} g=g^{-1}{\cdot} x$, by the diassociativity of $X$. \end{proof}

Claims~\ref{cl:xg=gx} implies

\begin{claim} For every $a\in X$ we have $Ga=aG$.
\end{claim}

\begin{claim}\label{cl:M(G,2)->X}  For every $a\in X\setminus G$, the function 
$$h:M(G,2)\to X,\quad
h(x,k)\defeq\begin{cases}
x&\mbox{if $k=1$};\\
x{\cdot}a&\mbox{if $k=-1$};
\end{cases}
$$ is an injective homomorphism of the Chein loop $M(G,2)$ onto the subloop $G\cup Ga$ of $X$. 
\end{claim}

\begin{proof} Given any elements $(x,i),(y,j)\in M(G,2)$, we should check that $h((x,i)\cdot (y,j))=h(x,i)\cdot h(y,j)$.
If $i=j=1$, then $$h((x,1)\cdot(y,1))=h(x\cdot y,1)=x\cdot y=h(x,1)\cdot h(y,1).$$
If $i=1$ and $j=-1$, then  $h((x,i)\cdot (y,j))=h((x^{-1}\cdot y^{-1})^{-1},-1))=h(y\cdot x,-1)=(y\cdot x)\cdot a$. 
On the other hand, Claim~\ref{cl:xg=gx} the Bol identity imply 
$$
\begin{aligned}
&h(x,i)\cdot h(y,j)=x\cdot(y\cdot a)=x\cdot(a\cdot y^{-1})=a\cdot(a\cdot (x\cdot (a\cdot y^{-1})))=a\cdot((a\cdot (x\cdot a))\cdot y^{-1})\\
&=a\cdot((a\cdot (a\cdot x^{-1}))\cdot y^{-1})=a\cdot (x^{-1}\cdot y^{-1})=(x^{-1}\cdot y^{-1})^{-1}a=(y\cdot x)\cdot a=h((x,i)\cdot (y,j)).
\end{aligned}
$$
If $i=-1$ and $j=1$, then $h((x,i)\cdot (y,j))=h((x\cdot y^{-1}),-1))=(x\cdot y^{-1})\cdot a$. 
On the other hand, Claim~\ref{cl:xg=gx} and the Bol identity imply
$$
\begin{aligned}
h(x,i)\cdot h(y,j)&=(x\cdot a)\cdot y=(((x\cdot a)\cdot y)\cdot a)\cdot a=(x\cdot((a\cdot y)\cdot a))\cdot a\\
&=(x\cdot ((y^{-1}\cdot a)\cdot a))\cdot a=(x\cdot y^{-1})\cdot a=h((x,i)\cdot (y,j)).
\end{aligned}
$$
If $i=-1$ and $j=-1$, then $h((x,i)\cdot (y,j))=h((x^{-1}\cdot y)^{-1},1))=y^{-1}\cdot x$.
On the other hand,  Claim~\ref{cl:xg=gx} and the Moufang identity imply 
$$
\begin{aligned}
h(x,i)\cdot h(y,j)&=(x\cdot a)\cdot (y\cdot a)=(a\cdot x^{-1})\cdot(y\cdot a)=(a\cdot(x^{-1}\cdot y))\cdot a\\
&=((x^{-1}\cdot y)^{-1}\cdot a)\cdot a)=y^{-1}\cdot x=h((x,i)\cdot(y,j))
\end{aligned}
$$
Thefore $h:M(G,2)\to X$ is an injective homomorphism of the loop $M(G,2)$ into the loop $X$.
\end{proof}

\begin{claim}\label{cl:G-subgroup} If $G\ne X$, then $G$ is a subgroup of $X$.
\end{claim}

\begin{proof} Take any element $a\in X\setminus G$ and consider the injective homomorphism $h:M(G,2)\to X$, defined in Claim~\ref{cl:M(G,2)->X}. Since $X$ is Moufang, so is its subloop $h[M(G,2)]$ and the isomorphic copy $M(G,2)$ of $h[M(G,2)]$. By Theorem~\ref{t:M(G,2)-Moufang<=>G-group}, $G$ is a group. 
\end{proof} 

Now we are able to complete the proof of Theorem~\ref{t:Chein-index}. 
Choose any element $a\in X\setminus G$ and consider the injective homomorphism $h_a:M(G,2)\to X$ defined in Claim~\ref{cl:M(G,2)->X}. Then $H\defeq h_a[M(G,2)]=G\cup Ga$ is a subloop of the Moufang loop $X$. By Claim~\ref{cl:G-subgroup}, $G$ is a group. If $H=X$, then $G$ is a normal subgroup of index $2$ in $X$ and $X$ is isomorphic to the Chein extension $M(G,2)$ of the group $G$.

So, assume that $H\ne X$. In this case $H\cong M(G,2)$ is a group, by Claim~\ref{cl:G-subgroup}, and $G$ is a commutative group, by Theorem~\ref{t:M(X,2)-associative<=>X-commutative}. Choose any element $b\in X\setminus H$ and consider the injective homomorphism $h_b:M(H,2)\to X$ defined in Claim~\ref{cl:M(G,2)->X}. Then $C\defeq h_b[M(H,2)]=H\cup Hb$ is a subloop of the Moufang loop $X$. Assuming that $C\ne X$, we can apply Claim~\ref{cl:G-subgroup} and conclude that $C\cong M(H,2)$ is a group. By Theorem~\ref{t:M(X,2)-associative<=>X-commutative}, $H\cong M(T,2)$ is a commutative group, and by Theorem~\ref{t:M(X,2)-commutative<=>X-Boolean}, $C$ is a Boolean subgroup of $X$, which is not possible because $C$ contains an element of order $>2$. This contradiction shows that $C=X$ and hence the commutative group $G$ is normal and has index 4 in the Moufang loop $X=C$, which is isomorhic to $M(H,2)\cong M(M(G,2),2)$.
\end{proof}

\begin{corollary}\label{c:ind8=>Boolean}  If an elementary Moufang loop $X$ contains a normal subloop $G$ of index $8$, then the loop $G$ is a Boolean group.
\end{corollary}

\begin{proof} By the normality of $G$ in $X$, every element $a\in X\setminus G$ has even order and by the elementarity of $X$, this order is equal to $2$. Therefore, the subloop $G$ contains all elements of order $>2$ in $X$. If $X$ is not Boolean, then $G$ contains an element of order $>2$ and hence the loop $G$ is not trivial. By Theorem~\ref{t:Chein-index}, the loop $G$ is normal and has index $|X/G|\in\{2,4\}$, which contradicts the asumption. This contradiction shows that the elementary Moufang loop $X$ is Boolean. By Proposition~\ref{p:Boolean+Moufang=>commutative-group}, the Boolean Moufang loop $X$ is a Boolean group.
\end{proof}

We shall also need the following fact proved by Chein in \cite{Chein1978}.

\begin{theorems}[Chein, 1978]\label{t:Chein63} Let $X$ be a non-associative Moufang loop such that every element $x\in X$ has order $\le 3$. If $|X|\le 63$, then $X$ is isomorphic to the Chein extension $M(G,2)$ of some group $G$. 
\end{theorems}

In the proof of Theorem~\ref{t:elementary-Moufang<=>} we shall use the following fact, proved by \index[person]{Vojt\v echovsk\'y}Petr Vojt\v echovsk\'y\footnote{{\bf Petr Vojt\v echovsk\'y} is a Czech‑American mathematician specializing in nonassociative algebra, particularly the theory of loops, quasigroups, and related combinatorial structures. He received his M.S. in Mathematics from Charles University in Prague in 1998 and his Ph.D. in Mathematics from the same university in 2003 under the supervision of Ale\v s Dr\'apal, with a dissertation on connections between codes, groups, and loops. He also holds a Ph.D. from Iowa State University, where he worked on finite simple Moufang loops. Vojt\v echovsk\'y is currently Professor of Mathematics at the University of Denver. His research encompasses the structure theory of Moufang and Bol loops, isotopy and nuclei, automorphism groups, and computational approaches in nonassociative algebra. He is the principal developer of the LOOPS package for GAP, a standard computational tool in loop and quasigroup theory. His work includes classification of small Moufang loops, investigations of Bruck loops and Bol loops of odd order, and connections between nonassociative algebra, combinatorics, and geometry. 
} in his Ph.D. Thesis \cite[\S5.2.5]{VojtechovskyPhD}.

\begin{propositions}[Vojt\v echovsk\'y, 2001]\label{p:M(S3,2)inM(2)} The simple Paige loop $M(2)$ contains exactly $112$ subloops, isomorphic to the Chein extension $M(S_3,2)$ of the symmetric group $S_3$.
\end{propositions}



\section{Elementary finite groups}

In this section we prove a classification of elementary finite groups, due to L\"uneburg \cite{Luneburg1961} and later reproved by Deaconescu \cite{Deaconescu1989} in 1989 and with corrections by Cheng Kai Nah, M. Deaconescu, Lang Mong Lung and Shi Wujie \cite{CDeaconescu1993} in 1993.

\begin{theorems}[L\"uneburg; 1961]\label{t:elementary-group<=>} If $X$ is an elementary finite group, then one of the following conditions holds:
\begin{enumerate}
\item $X$ has prime-power order;
\item $X$ is a Frobenius group whose kernel $K$ has prime-power order and the complement has prime order that divides $|K|-1$;
\item $X$ is isomorphic to the alternating group $A_5$.
\end{enumerate}
\end{theorems}

\begin{proof} If $X$ is solvable, then by Theorem~\ref{t:nucleo-solvable}, $X$ has prime-power order or else $X$ is a Frobenius group whose kernel $K$ is has prime-power order and complement has prime order $q$ that divides $|K|-1$.

Next, assume that $X$ is not solvable. By induction on the order $X$, we shall prove that $X$ is isomorphic to the alternating group $A_5$. Assume that for some $n\in\IN$ we know that all elementary non-solvable groups of order $<n$ are isomorphic to the group $A_5$. Let $G$ be a group of order $n$. If $G$ is simple, then by Theorem~\ref{t:elementary-simple<=>}, $G$ is isomorphic to $A_5$ and we are done. So, assume that $G$ is not simple. Then $G$ contains a non-trivial proper normal subgroup $N$. Since $G$ is not solvable, either $N$ is not solvable or $G/N$ is not solvable. Let $\pi:G\to G/N$ be the quotient homomorphism.

First assume that $G/N$ is not solvable. By the inductive assumption, $G/N$ is isomorphic to the alternating group $A_5$. Then $G/N$ contains subgroups isomorphic to the Chein extensions $M(C_3,2)$ and $M(C_5,2)$ of cyclic groups of order $3$ and $5$. Also $G/N\cong A_5$ contains a Boolean subgroup $B$ of order 4.  Then $\pi^{-1}[B]$ is an elementary group containing the normal subgroup $N$ of index 4. By Theorem~\ref{t:Chein-index}, the group $N$ is commutative and being elementary, has  prime-power order $|N|=p^n$. Let $k\in\{3,5\}$ be any number such that $p\ne k$. Since the alternating group $A_5$ contains a subgroup isomorphic to the dihedral group $M(C_k,2)$, we can find a subgroup $E\subseteq G/N$ isomorphic to the Chein extension $M(C,2)$ of a cyclic group $C\subseteq E$ of order $|C|=k$. Then $\pi^{-1}[C]$ is a normal subgroup of order $kp^n$ and index $2$ in the group $\pi^{-1}[E]$. By Theorem~\ref{t:Chein-index}, the group $\pi^{-1}[E]$ is isomorphic to the Chein extension $M(\pi^{-1}[C],2)$ of the group $\pi^{-1}[C]$. Since $M(\pi^{-1}[C],2)$ is a group, $\pi^{-1}[C]$ is a commutative group, by Theorem~\ref{t:M(X,2)-associative<=>X-commutative}. Since the commutative group $\pi^{-1}[C]$ contains elements of order $p$ and $k$, it contains an element of non-prime order $pk$, which contradicts the elementarity of the group $X$. This contradiction shows that the group $G/N$ is solvable. Then the group $N$ is not solvable and by the inductive assumption, $N$ is isomorphic to the group $A_5$. Consider the homomorphism $\alpha:G\to \Aut(N)$ assigning to every $y\in G$ the automorphism $\alpha_y:N\to N$, $\alpha_y:x\mapsto y^{-1}xy$, of the group $N\cong A_5$. The elementarity of $X$ ensures that every element $y\in G\setminus\{e\}$ has prime order. Choose an element $x\in N$ whose order is distinct from the order of $y$. By the elementarity of $X$, the elements $x$ and $y$ do not commute, which implies that $\alpha(x)=y^{-1}xy\ne x$. So, the automorphism $\alpha_y$ is not trivial, and the homomorphism $\alpha:G\to\Aut(N)$ is injective.
It is known that the automorphism group of the group $A_5$ is isomorphic to the symmetric group $S_5$. Then the image $\alpha[N]$ has index 2 in the group $\Aut(N)$, which implies that the group $N$ has index 2 in $G$. By Theorem~\ref{t:Chein-index}, $G$ is isomorphic to the Paige extension $M(N,2)$ of the group $N\cong A_5$. Since $G\cong M(N,2)$ is associative, we can apply Theorem~\ref{t:M(X,2)-associative<=>X-commutative} and conclude that the group $N\cong A_5$ is commutative, which is not true. This contradiction completes the proof.
\end{proof}

We shall deduce from Theorem~\ref{t:elementary-group<=>} the following corollary that will be used in the proof of Theorem~\ref{t:elementary-group-solvable<=>}. 

\begin{corollarys}\label{c:A5} If an elementary Moufang loop $X$ contains a normal proper subgroup $G$, isomorphic to the group $A_5$, then $G$ has index $2$ in $X$ and $X$ is isomorphic to $M(A_5,2)$.
\end{corollarys}

\begin{proof} If $G$ has index $2$ in $X$, then $X$ is isomorphic to $M(A_5,2)$, by Theorem~\ref{t:Chein-index}.

So, assume that the quotient loop $X/G$ has cardinality $>2$. If $|X/G|$ is divisible by $4$, then by Theorem~\ref{t:23Sylow}, the Sylow $2$-group $S$ of $X/G$ has cardinality $\ge 4$. The elementarity of the Moufang loop $X$ implies the elementarity of the quotient loop $X/G$. By the elementarity of $X/G$, the $2$-group $S$ is Boolean. Since $|S|\ge 4$, $S$ contains a Boolean subgroup $B\subseteq S$ of order $4$. Let $\pi:X\to X/G$ be the quotient homomorphism. Then $\pi^{-1}[B]$ is an elementary Moufang loop containing the normal subloop $G$ of index 4. By Theorem~\ref{t:Chein-index}, the group $G\cong A_5$ is commutative, which is not true. This contradiction shows that $|X/G|$ is not divisible by 4.  

Since $X/G$ is an elementary loop of order $|X/G|>2$ and and $|X/G|$ is not divisible by 4, $X/G$ has an element $y$ of prime order $p>2$. Choose any point $x\in\pi^{-1}(y)$ and observe that the order of $x$ equals the order $p$ of $y$, by the elementarity of $X$.  Let $P$ be the cyclic subgroup of order $p$ generated by $x$ in the Moufang loop $X$. If $p\ne 5$, then take any cyclic subgroup $C$ of order $5$ in $G\cong A_5$ and let $c$ be a generator of the group $C$. Consider the subloop $H\subseteq X$ generated by the elements $c$ and $x$. Since the Moufang loop $X$ is diassociative, the subloop $H$ is a subgroup of $X$. Assuming that $H\cap G=C$, we conclude that $C$ is a normal subgroup of $H$. Consider the homomorphism $\alpha:P\to\Aut(C)$ assigning to every element $z\in P$ the automorphism $\alpha_z:C\to C$, $\alpha_z:\gamma\mapsto z^{-1}\gamma z$. The automorphism group of the cyclic group $C$ is isomorphic to the multiplicative group $\IF _5^*$ of the field $\IF_5$. For every $z\in P\setminus \{e\}$ and $\gamma\in C\setminus\{e\}$ the orders of $z$ and $\gamma$ are distinct. By the elementarity, the elements $z$ and $\gamma$ do not commute, which implies that $\alpha_z(\gamma)\ne\gamma$ and hence the homomorphism $\alpha:P\to\Aut(C)$ is injective, which is not possible because the odd order $p$ of the cyclic group $P$ does not divide the order $|\Aut(C)|=4$ of the automorphism group $\Aut(C)\cong\IF_5^*$. This contradiction shows that $H\cap G\ne C$ and hence $xCx^{-1}$ are two distinct subgroups of order $5$ in the subgroup $H\cap G$ of the group $G\cong A_5$.  Since no proper subgroup of the group $A_5$ contains two distinct $5$-elements subgroups, $G\cong A_5$ is a proper subgroup of the elementary group $H$, which contradicts Theorem~\ref{t:elementary-group<=>}.

This contradiction shows that $p=5$. In this case, take any cyclic subgroup $C$ of order $3$ in $G\cong A_5$ and let $c$ be a generator of the group $C$. Consider the subloop $H\subseteq X$ generated by the elements $c$ and $x$. Since the Moufang loop $X$ is diassociative, the subloop $H$ is a subgroup of $X$. Let $z\in P\setminus\{e\}$ be any element. Since the elements $z$ and $c$ have distinct orders and the loop $X$ is elementary, $zc\ne cz$. Assuming that $zcz^{-1}=c^{-1}$, we conclude that
$c=z^pcz^{-p}=c^{(-1)^5}=c^{-1}$, which is a contradiction showing that $zcz^{-1}\ne c^{-1}$ and hence $C$ and $zCz^{-1}$ are two distinct subgroup of order $3$ in $H$. This implies that the subgroup $H\cap G$ contains five pairwise distinct subgroups $zCz^{-1}$, $z\in P$, of order 3. Since no proper subgroup of $A_5$ contains five pairwise distinct subgroups of order 3, $G\cong A_5$ is  proper subgroup of the elementary group $H$, which contradicts Theorem~\ref{t:elementary-group<=>}.
\end{proof}

\section{Elementary solvable finite Moufang loops}

In this section we prove a classification theorem for elementary solvable finite Moufang loops.

First we prove two lemmas.

\begin{theorems}\label{t:solvable-elementary<=>} If a finite solvable Moufang loop $X$ is elementary, then one of the following conditions holds.
\begin{enumerate}
\item $X$ has prime-power order.
\item $X$ is a Frobenius loop of odd order whose kernel $K$ has  prime-power order and complement has prime order that divides $|K|-1$.
\item $X$ is isomorphic to the Chein extension $M(G,2)$ of some elementary solvable group $G$.
\item $X$ is a semidirect product $B\rtimes E$ of a normal Boolean subgroup $B\subset X$ and an elementary Moufang loop $E\subset X$ of odd order.
\item $X$ is a semidirect product $A\rtimes B$ of an elementary Abelian group $A\subset X$ of odd order and a subgroup $B\subset X$, isomorphic to the alternating group $A_4$.
\end{enumerate}
\end{theorems}

\begin{proof} The proof is by induction on the cardinality of $X$. Assume that for some $n\in\IN$, any elementary finite solvable Moufang loop $X$ of cardinality $|X|<n$ satisfies one of the conditions (1)--(5). 

Let $X$ be an elementary finite solvable Moufang loop of order $n$.
If $n$ is odd, then the Moufang loop $X$ is nucleo-solvable, by Theorem~\ref{t:odd=>nucleo-solvable}. Applying Theorem~\ref{t:nucleo-solvable}, we conclude that $X$ satisfies condition (1) or (2). So, assume that $n$ is even. If $n$ is a power of $2$, then $X$ is a Boolean group, by Proposition~\ref{p:Boolean+Moufang=>commutative-group}. So, assume that $n$ has an odd prime divisor.

Let $N$ be a normal subloop of maximal order in $X$ such that $|N|$ is either odd or a power of 2. Since the loop $X$ is solvable and elementary, the subloop $N$ is not trivial. Then the quotient loop $X/N$ is elementary, solvable, and has order $|X/N|<|X|=n$. By the inductive assumption, $X/N$ satisfies one of the conditions (1)--(5). Let $\pi:X\to X/N$ be the quotient homomorphism.

If $X/N$ satisfies the conditions (1), (2) or (4), then $X/N$ is a semidirect product $B\rtimes C$ of a Boolean group $B$ and an elementary Moufang loop $C$ of odd order (if the condition (1) holds, then $B$ or $C$ is trivial, if (2) holds, then $B$ is trivial). 

 If $|N|$ is a power of $2$, then $K\defeq \pi^{-1}[B]$ is a normal subloop of $X$ with $|K|=|B|\cdot|N|$ a power of $2$. The elementarity of the loop $X$ ensures that the loop $K$ is a Boolean group. Since $|K|$ is coprime with $|X/K|=|C|$, we can apply the Hall's Theorem~\ref{t:Hall-Moufang} and find a subloop $E\subseteq X$ of odd order $|E|=|X/K|$.
Since $|K|$ and $|E|$ are coprime, $K\cap E=\{e\}$ and hence $X=K\rtimes E$ is a semidirect product of the Boolean subgroup $K$ and the elementary Moufang loop $E$ of odd order. 

So, assume that $|N|$ is not a power of $2$ and hence $|N|$ is odd. Since $|X|=|N|\cdot|B|\cdot|C|$ is even, $|B|\ge 2$.  Applying Proposition~\ref{p:Frobenius-prime-divides}, we can show that the order of any element $c\in C\setminus B$ divides $|B|-1$, which implies $|B|\ge 4$. 
Assuming that $|B|>4$,  we can find a Boolean subgroup $D\subseteq B$ of order $8$. Then the normal subloop $N$ has index $8$ in the elementary Moufang loop $\pi^{-1}[D]$, and hence $N$ is Boolean, by Corollary~\ref{c:ind8=>Boolean}. But this contradicts our assumption. This contradiction shows that $|B|=4$.
Since the Boolean group $B$ has odd index $|C|$ in the Moufang loop $X/N$, we can apply Lemma~\ref{l:4n=>3} and conclude that $|C|=3$ and the loop $X/N$ is isomorphic to the alternating group $A_4$. By Lemma~\ref{l:4=>A45}, $X$ is a semidirect product $N\rtimes A$ of the normal loop of odd order $N$ and a subgroup $A$, isomorphic to the alternating group $A_4$. Therefore, the loop $X$ satisfies the condition (5) of the theorem.
\smallskip

If the loop $X/N$ satisfies the condition (3), then $X/N$ contains a normal subgroup $H$ of index $2$.  Then $G\defeq\pi^{-1}[H]$ is a normal subloop of index $2$ in the loop $X$. By the elementarity of $X$ every element in the set $X\setminus G$ has order $2$. By Theorem~\ref{t:Chein-index}, the subloop $G$ is a group and $X$ is isomorphic to the Chein extension $M(G,2)$. The group $G$ is elementary and solvable, being a subloop of the elementary solvable loop $X$. Therefore, the loop $X$ satisfies the condition (3).

If the loop $X/N$ satisfies the condition (5), then $X$ is the semidirect product $X/N=A\rtimes B$ of a normal elementary Abelian group $A$ of odd order and a subgroup $B$, isomorphic to $A_4$. Then the group $B$ contains a normal subgroup $H\subseteq B$ of  order $4$ and index 3 in $B$. The loop $A$ has index $4$ in the normal subloop $A{\cdot}H$ of the loop $X/N$, and the normal subloop $\pi^{-1}[A]$ has index $4$ in the elementary Moufang loop $\pi^{-1}[A{\cdot}H]$. Applying Theorem~\ref{t:Chein-index}, we conlcude that the loop $\pi^{-1}[A]$ is an elementary Abelian group. If this group is Boolean, then the loop $K\defeq \pi^{-1}[A{\cdot}H]$ also is Boolean. Since the group $H$ has index 3 in $B$, the loop $A{\cdot}H$ has index 3 in the loop $X/N$, and the Boolean group $K$ has index 3 in the elementary Moufang loop $X$. By the elementarity of $X$, any element $x\in X\setminus  K$ has order 3 and hence generates a cyclic group $C$ of order $3$ in $X$. Since $K$ is Boolean, $C\cap K=\{e\}$ and $X=K\rtimes C$ is a semidirect product of the Boolean group $K$ and the elementary group $C$ of odd order. This shows that the condition (3) is satisfied. So, assume that the elementary Abelian group $E\defeq \pi^{-1}[A]$ is not Boolean. Then $|E|$ is odd. Since $X/E\cong B\cong A_4$, we can apply Lemma~\ref{l:4=>A45} and conclude that $X$ is the semidirect product $E\rtimes \Gamma$ of the elementary Abelian group $E$ of odd order and a subgroup $\Gamma\subseteq X$, isomorphic to the alternating group $A_4$. Therefore, $X$ satisfies the condition (5).
\end{proof}

It remains to prove two lemmas, used in the proof of Theorem~\ref{t:solvable-elementary<=>}.

\begin{lemmas}\label{l:4n=>3} Let $X$ be an elementary Moufang loop containing a normal subloop $B$ of order $|B|=4$ and odd index $|X/B|\ge 3$. Then $|X/B|=3$ and $X$ is isomorphic to the group $A_4$.
\end{lemmas}

\begin{proof} By the elementarity of $X$, the subloop $B\subseteq X$ of order $|B|=4$ is a Boolean group. 

We claim that every element $c\in X/B$ has order 3. Since $X$ is elementary and $|X/B|$ is odd, the order of $x$ is an odd prime $p$. Fix any element $b\in B$ and consider the subloop $Y$, generated by the elements $c$ and $b$. The diassociativity of the Moufang loop $X$ ensures that $Y$ is a group. Since $b$ and $c$ have distinct prime order, $c{\cdot}b\ne b{\cdot}c$ and hence $b, c^{-1}bc$ and $cbc^{-1}$ are three distinct elements of order 2 in the Boolean group $B$. Since $|B|=4$, $B=\{e,b,cbc^{-1},c^{-1}bc\}\subseteq Y$ and hence $Y$ is an elementary solvable group of cardinality $|Y|=4p$. By Theorem~\ref{t:elementary-group<=>}, $p$ divides $|B|-1=3$ and hence $p=3$.

Therefore, every element of the Moufang loop $X/B$ has order $3$ and by the Cauchy property of the solvable Moufang loop $X/B$, $|X/B|=3^n$ for some $n\in\IN$. By Theorem~\ref{t:prime-power=>center}, the loop $X/B$ has non-trivial center.
If $|X/B|>3$, then $X/B$ contains a commutative subgroup $T$ of order $|T|=9$. Let $\pi:X\to X/B$ be the quotient homomorphism, and observe that $E\defeq \pi^{-1}[T]$ is an elementary Moufang loop of order 36. Theorem~\ref{t:elementary-group<=>} ensures that the elementary Moufang loop $E$ is not associative. Since $|E|=36=2^2\cdot 3^2$, the elementarity of $E$ ensures that every element of $E$ has order $\le 3$. By Theorem~\ref{t:Chein63}, the loop $E$ is isomorphic to the Chein extension $M(G,2)$ of some subgroup $G$ of $E$. Then $G$ is a normal subgroup of index $2$ in $E$ and $B\cap G$ is a normal subgroup of cardinality $|B\cap G|=2$ in $E$. Then for the unique non-identity element $b'\in B\cap G$ we have $c{\cdot}b'{\cdot}c^{-1}=b'$ and hence $c{\cdot}b'=b'{\cdot} c$ is an element of order 6 in $E$, which contradicts the elementarity of the loop $X$. This contradiction shows that $|X/B|=3$ and $X$ is an elementary group of order 12, containing a normal subgroup of order 4. The unique group with those properties is the alternating group $A_4\cong(C_2\times C_2)\rtimes C_3$.
\end{proof}

\begin{lemmas}\label{l:4=>A45} Let $X$ be an elementary solvable finite Moufang loop and $H$ be a normal subloop of odd oder in $X$. If the quotient loop $X/H$ is isomorphic to the alternating group $A_4$, then $X$ is the semidirect product $H\rtimes A$ of the subloop $H$ and a subgroup $A\subseteq X$, which is isomorphic to $A_4$.
\end{lemmas}

\begin{proof} Let $\pi:X\to X/H$ be the quotient homomorphism. Since the group $A_4$ is 2-generated, there exist elements $b,c\in X$ whose images $\pi(b),\pi(c)$ generate the group $X/N\cong A_4$. Let $A$ be the subloop of $X$, generated by the elements $b$ and $c$ in $X$. The diassociativity of the Moufang loop $X$ ensures that $A$ is a subgroup of $X$. Taking into account that the elements $\pi(b)$ and $\pi(c)$ generate the group $X/N\cong A_4$, we conlcude that $\pi[A]=X/N$. We claim that $A\cap H=\{e\}$. To derive a contradiction, assume that $A\cap H\ne \{e\}$. Then $A\cap H$ is a non-trivial normal subgroup in the group $A$, witnessing the the group $A$ is not simple. Since $|\pi[A]|=|X/N|=|A_4|=12$, the order $|A|$ is not prime-power. By Theorem~\ref{t:elementary-group<=>}, the group $A$ is a Frobenius group whose kernel  $K$ has prime-power order $p^n$ and complement $C$ has prime order $q$ that divides $p^n-1$. Taking into account that $12=|\pi[A]|$ divides $|A|=p^nq$, we conclude that $p=2$ and $q=3$. Assuming that $A\cap H\not\subseteq K$, we conclude that $A\subseteq K\cdot H$ and then $|\pi[A]|$ is a power of $2$, which contradicts $|\pi[A]|=12$. This contradiction shows that $A\cap H\subseteq K$. Since $|H|$ is odd and $|K|$ is a power of $2$, the group $A\cap K$ is trivial. Then the restriction $\pi{\restriction}_A:A\to X/H$ is an isomorphisms and hence the subgroup $A\subseteq X$ is isomorphic to the group $A_4$. Now we see that the group $X=H\rtimes A$ is a semidirect product of the groups $H$ and $A$.
\end{proof}

\begin{lemma}\label{l:X/N=A5=>N=1} Let $N$ be a normal subloop of an elementary Moufang loop $X$. If the quotient loop $X/N$ is isomorphic to the alternative group $A_5$, then the subloop $N$ is trivial.
\end{lemma}


\section{Elementary group-solvable finite Moufang loops}

\begin{theorems}\label{t:elementary-group-solvable<=>} If a group-solvable finite Moufang loop $X$ is elementary, then one of the following conditions holds.
\begin{enumerate}
\item $X$ has prime-power order.
\item $X$ is a Frobenius loop whose kernel $K$ has odd prime-power order and complement has prime order $q\ge 3$ that divides $|K|-1$.
\item $X$ is isomorphic to the Chein extension $M(G,2)$ of some elementary group $G$.
\item $X$ is the semidirect product $K\rtimes H$ if a Boolean group $K$ and an elementary Moufang loop $H$ of odd order.
\item $X$ is the semidirect product $A\rtimes B$ of an elementary Abelian group $A\subseteq X$ of odd order and a subgroup $B$, isomorphic to the alternating groups $A_4$ or $A_5$.
\end{enumerate}
\end{theorems}

\begin{proof} The proof is by induction on the cardinality of $X$. Assume that for some $n\in\IN$, any elementary finite group-solvable Moufang loop $X$ of cardinality $|X|<n$ satisfies one of the conditions (1)--(5). 

Let $X$ be an elementary finite group-solvable Moufang loop of order $n$.
If $n$ is odd, then the Moufang loop $X$ is nucleo-solvable, by Theorem~\ref{t:odd=>nucleo-solvable}. Applying Theorem~\ref{t:nucleo-solvable}, we conclude that $X$ satisfies condition (1) or (2). So, assume that $n$ is even. If $n$ is a power of $2$, then $X$ is a Boolean group. So, assume that $n$ has an odd divisor.

Since $X$ is group-solvable, $X$ contains a non-trivial normal subgroup $G\subseteq X$. By Theorem~\ref{t:elementary-group<=>}, $G$ is prime-power, or Frobenius loop whose kernel $K$ is prime power and complement has prime order dividing $|K|-1$ or else $G$ is isomorphic to the alternating group $A_5$. 

First we assume that $G$ is isomorphic to $A_5$. If $G=X$, then the condition (5) holds. If $G\ne X$, then $X$ is isomorphic to the Chein extension $M(G,2)$ of the elementary group $G\cong A_5$, by Corollary~\ref{c:A5}, and hence the condition (3) holds.

So, assume that $G$ is not isomorphic to $A_5$. 

\begin{claim}\label{cl:odd-or2} The group $G$ contains a non-trivial subgroup $K$ such that $K$ is normal in $X$ and $|K|$ is odd or a power of $2$.
\end{claim}

\begin{proof} If $G$ is prime-power or $|G|$ is odd, then the subgroup $K\defeq G$ has the required property. So, assume that $G$ is not prime-power and $|G|$ is even. Then $|G|\ge 4$. By Theorem~\ref{t:elementary-group<=>}, $G$ is a Frobenius group whose kernel $K$ is prime-power and the complement has prime order $p$ that divides $|K|-1$. By Proposition~\ref{p:I-normal}, the set of involutions $I\defeq\{x\in X\setminus\{e\}:x{\cdot}x=e\}$ is normal in the loop $X$.

 If $p=2$, then $|K|$ is odd because $2=p$ divides $|K|-1$. By the elementarity of the loop $X$, every element $x\in G\setminus K$ has order $2$ and hence $G\setminus K=G\cap I$ is a normal set in $X$ and then $K=G\setminus I$ is a normal subgroup of $X$.

If $p\ne 2$, then $|K|$ is a power of $2$ because the order $|G|=p|K|$ is even. The elementarity of $X$ ensures that the group $K$ is Boolean and hence the group $K=G\cap(I\cup\{e\})$ is normal in $X$.
\end{proof}

Let $N$ be a maximal normal subgroup of $X$ such that the order $|N|$ is either odd or a power of $2$. Claim~\ref{cl:odd-or2}  ensures that the group $N$ is not trivial. Since the elementary Moufang loop $X$ is group-solvable, the quotient Moufang loop $X/N$ is elementary and group-solvable. By the inductive assumption, $X/N$ satisfies one of the conditions (1)--(5).
Let $\pi:X\to X/N$ be the quotient homomorphism.

If the loop $X/N$ contains a normal subgroup $H$ of index $2$, then $E\defeq\pi^{-1}[H]$ is a normal subgroup of index $2$ in the elementary Moufang loop $X$. By Theorem~\ref{t:Chein-index}, $E$ is a group and $X$ is isomorphic to the Chein extension $M(E,2)$ of the elementary Moufang group $E$, witnessing that the loop $X$ satisfies the condition (3). So, assume that the quotient loop $X/N$ does not contain normal subloops of index $2$. In this case, the group $X/N$ does not satisfy the condition (3) of Theorem~\ref{t:elementary-group-solvable<=>}.

If $X/N$ satisfies the conditions (1), (2) or (4), then $X/N$ is the  semidirect product $B\rtimes C$ of a Boolean group $B$ and an elementary Moufang loop $C$ of odd order (in cases (1) and (2) the Boolean group $B$ is trivial). Since $X/N$ contains no normal subloops of index $2$, the loop $C$ of odd order is not trivial.  If $|N|$ is a power of $2$, then $K\defeq \pi^{-1}[B]$ is a normal subloop of $X$ with $|K|=|B|\cdot|N|$ a power of $2$. The elementarity of the loop $X$ ensures that the loop $K$ is a Boolean group and so is its normal subloop $N$. Then the loop $X$ is solvable. Since $|K|$ is coprime with $|X/K|=|C|$, we can apply the Hall Theorem~\ref{t:Hall-Moufang} and find a subgroup $E\subseteq X$ of odd order $|E|=|X/K|$.
Since $|K|$ and $|E|$ are coprime, $K\cap E=\{e\}$ and hence $X=K\rtimes E$ is a semidirect product of the Boolean subgroup $K$ and the elementary Moufang loop $E$ of odd order, witnessing that the condition (4) holds.

So, assume that $|N|$ is not a power of $2$ and hence $|N|$ is odd. Since $|X|=|N|\cdot|B|\cdot|C|$ is not odd, $|B|\ge 2$. Given any element $c\in C\setminus B$, consider the subgroup $\langle \{c\}\cup B\rangle$ of $X/N$, generated by the set $\{c\}\cup B$. Since $|C|$ odd, the elementarity of $X$ ensures that the order of $c$ is an odd prime number. Applying Theorem~\ref{t:elementary-group<=>}, we conclude that the order of $c$ divides $|B|-1$, which implies $|B|\ge 4$. 
Assuming that $|B|>4$,  we can find a Boolean subgroup $D\subseteq B$ of order $8$. Then the normal subloop $N$ has index $8$ in the elementary Moufang loop $\pi^{-1}[D]$, and hence $N$ is Boolean, by Corollary~\ref{c:ind8=>Boolean}. But this contradicts our assumption. This contradiction shows that $|B|=4$.

Since the normal subloop $N$ has index $4$ in the elementary Moufang loop $\pi^{-1}[B]$, the loop $N$ is an elementary Abelian group, by Theorem~\ref{t:Chein-index}. 
Since the Boolean group $B$ has odd index $|C|$ in the Moufang loop $X/N$, we can apply Lemma~\ref{l:4n=>3} and conclude that $|C|=3$ and the loop $X/N$ is isomorphic to the alternating group $A_4$. By Lemma~\ref{l:4=>A45}, $X$ is a semidirect product $N\rtimes A$ of the elementary Abelian group $N$ and the subgroup $A$, isomorphic to the alternating group $A_4$. Therefore, the loop $X$ satisfies the condition (5) of the theorem.
\smallskip


If the loop $X/N$ satisfies the condition (5), then $X$ is the semidirect product $X/N=A\rtimes B$ of a normal elementary Abelian group $A$ and a subgroup $B$, isomorphic to $A_4$ of $A_5$. Since the group $B$ is isomorphic to $A_4$ or $A_5$, it contains a Boolean subgroup  $D\subset B$ of  order $4$.  Then the normal group $A$ has index $4$ in the  subloop $A{\cdot}D$ of the loop $X/N$, and the normal subloop $K\defeq\pi^{-1}[A]$ has index $4$ in the elementary Moufang loop $E\defeq \pi^{-1}[A{\cdot}D]$. Applying Theorem~\ref{t:Chein-index}, we conlcude that the loop $K$ is an elementary Abelian group such that the quotient loop $X/K$ is isomorphic to the group $B$, which is isomorphic to $A_4$ or $A_5$.

First assume that the group $B$ is isomorphic to the alternating group $A_4$. 
If the elementary Abelian group $K$ is Boolean, then the normal loop $E$ has order a power of $2$ and hence $E$ is a Boolean group, by the elementarity of $X$. Since the group $D$ has index $3$ in $B\cong A_4$, the Boolean group $E$ has  index 3 in $X$. By the elementarity, the group $X$ contains a cyclic subgroup $C$ of order $3$. Then $X=E\rtimes C$ is a semidirect product of the Boolean group $E$ and the cyclic group $C$ of odd order, witnessing that the condition (3) is satisfied. If the elementary Abelian group $K$ is not Boolean, then $|K|$ is an odd prime power and hence the subgroup $N$ of $K$ has odd order. Applying Lemma~\ref{l:4=>A45}, we conclude that that $X$ is the semidirect product $K\rtimes \Gamma$ of the elementary Abelian group $K$ and a subgroup $\Gamma\subseteq X$, which is isomorphic to the alternating group $A_4$. Therefore, $X$ satisfies the condition (5).

Finally, that the group $B$ is isomorphic to the group $A_5$. Then the quotient loop $X/K\cong B$ is isomorphic to the group $A_5$. Consider te quotient homomorphism $\tilde \pi:X\to X/K$. Since the group $A_5$ is $2$-generated, we can chose two elements $a,b\in X$ whose images $\tilde \pi(a)$ and $\tilde\pi(b)$ generate the simple group $A_5$. The diassociativity of the Moufang loop $X$ ensures that the subloop $\Gamma\subseteq X$ generated by the elements $a,b$ is a group. Since $\tilde\pi[\Gamma]=X/K\cong A_5$, the group $\Gamma$ is not solvable. By Theorem~\ref{t:elementary-group<=>}, $\Gamma$ is isomorphic to $A_5$. Then the surjective homomorphism $\tilde\pi{\restriction}_\Gamma:\Gamma\to X/K$ is an isomorphism, and hence $K\cap\Gamma=\{e\}$ and $X$ is the semidirect product $K\rtimes\Gamma$ of the elementary Abelian group $K$ and the group $\Gamma\cong A_5$, witnessing that the condition (5) is satisfied.   
\end{proof}

\section{Elementary finite Moufang loops}

Now we are able to prove the main classification theorem.

\begin{theorems}\label{t:elementary-Moufang<=>} If a finite Moufang loop $X$ is elementary, then one of the following conditions holds.
\begin{enumerate}
\item $X$ has prime-power order.
\item $X$ is a Frobenius loop whose kernel $K$ has odd prime-power order and complement has prime order $q\ge 3$ that divides $|K|-1$.
\item $X$ is isomorphic to the Chein extension $M(G,2)$ of some elementary group $G$.
\item $X$ is the semidirect product $K\rtimes H$ of a Boolean group $K$ and an elementary Moufang loop $H$ of odd order.
\item $X$ is the semidirect product $A\rtimes B$ of an elementary Abelian group  $A$ of odd order and a subgroup $B\subseteq X$, isomorphic to the alternating groups $A_4$ or $A_5$.
\item $X$ is isomorphic to the simple Paige loop $M(q)$ for some $q\in \{2,3,4,5\}$.
\end{enumerate}
\end{theorems}

\begin{proof} By Theorem~\ref{t:structure-Moufang}, the finite Moufang loop $X$ contains two normal subloops $H\subseteq \Pi$ such that the quotient loop $X/\Pi$ is Boolean, the quotient loop $\Pi/H$ is the direct product of some family $\mathcal P$ of simple Paige loops, and $H$ is a group-solvable Moufang loop.

If the Boolean group $X/\Pi$ is not trivial, then $X$ contains a normal subloop $H$ of index $2$. In this case $H$ is an elementary group and $X$ is isomorphic to $M(H,2)$, by Theorem~\ref{t:Chein-index}. Therefore, the case (3) holds. 

So, assume that the Boolean group $X/\Pi$ is trivial and hence $\Pi=X$. If $H\ne \Pi$, then the family $\mathcal P$ of simple Paige loops is not empty. Since the direct product of Paige loops $\Pi/H$ is elementary and every Paige loop contains elements of two distinct prime orders, the family $\mathcal P$ contains a single Paige loop $M(q)$ for some number $q$ that belongs to the set $q\in\{2,3,4,5\}$, by Theorem~\ref{t:elementary-simple<=>}. Then the loop $X/H$ is isomorphic to the Paige loop $M(q)$. We claim that the loop $H$ is trivial. To derive a contradiction, assume that $H$ is not trivial. 

By Proposition~\ref{p:M(2)-universal}, the loop $X/H\cong M(q)$ contains a subloop $P$, isomorphic to the simple Paige loop $M(2)$. By  Theorem~\ref{t:23Sylow}, the Paige loop $P$ contains a Sylow $2$-subloop $S$ of order $8$. Since $X$ is elementary, the loop $S$ is a Boolean group. Let $\pi:X\to X/H$ be the quotient homomorphism, and observe that $H$ is a normal subloop of index 8 in the elementary Moufang loop $\pi^{-1}[S]$. By Corollary~\ref{c:ind8=>Boolean}, the loop $H$ is a Boolean group. By Proposition~\ref{p:M(S3,2)inM(2)}, the simple Paige loop $M(2)$ contains a subloop isomorphic to the Chein extension $M(S_3,2)$ of the symmetric group $S_3$. Then we can choose a subloop $D\subseteq P$ containing a normal cyclic subgroup $C$ of order 3 and index 4 in $D$. Then $\pi^{-1}[C]$ is a normal subloop of order $3|H|$ and index 4 in the elementary Moufang loop $\pi^{-1}[D]$. By Theorem~\ref{t:Chein-index}, the loop $\pi^{-1}[C]$ is a commutative group. Since $H$ is a Boolean subgroup of index 3 in $\pi^{-1}[C]$, the commutative group $\pi^{-1}[C]$ contains an element of order 6, which contradicts the elementarity of $X$. This contradiction shows that the loop $H$ is trivial and hence the loop $X$ is isomorphic to the Paige loop $M(q)$. Therefore the condition (7) holds.

Finally, assume that $H=\Pi=X$. In this case the elementary Moufang loop $X=H$ is group-solvable, and by Theorem~\ref{t:elementary-group-solvable<=>}, one of the conditions (1)--(5) holds.
\end{proof}

Applying~Theorem~\ref{t:elementary-Moufang<=>} and   Proposition~\ref{p:=>cyclic-divides-kernel-1}, we obtain the following classification of possible cardinalities of elementary finite Moufang loops. 

\begin{corollarys}\label{c:elementary-Moufang=>cardinality} The cardinality of any elementary finite Moufang loop is equal to one of the numbers:
\begin{enumerate}
\item $p^n$ for some prime number $p$ and some positive integer $n$;
\item $p^nq$ for some $n\in\IN$ and some odd primes $p,q$ such that $q$ divides $p^n-1$;
\item $2^np^m$ for $n\ge 3$, $m\in\IN$ and some odd prime number $p$ that divides $2^n-1$;
\item $2^np^mq$ for some $n,m\in\IN$ and some odd prime numbers $p,q$ such that $q$ divides $p^m-1$ and $pq$ divides $2^n-1$;
\item $2p^n$ for some odd prime number $p$ and some $n\in\IN$;
\item $4p^n$ for some odd prime number $p$ and some $n\in\IN$;
\item $2p^nq$ for some $n\in\IN$ and some prime number $p$ and odd prime number $q$ that divides $p^n-1$;
\item $12p^n$ for some $n\in\IN$ and some odd prime number $p\ge 5$ such that $3$ divides $p^n-1$;
\item $60 p^n$ for some $n\in\IN$ and some prime number $p\ge 7$ such that $30$ divides $p^n-1$;
\item $60,120,1080,16320$ or $39000$.
\end{enumerate}The numbers in cases \textup{(1)--(10)} are pairwise distinct.
\end{corollarys}

\begin{proof} By Theorem~\ref{t:elementary-Moufang<=>}, one of the following six cases holds:
\begin{enumerate}
\item[(i)] $X$ has prime-power order;
\item[(ii)] $X$ is a Frobenius loop whose kernel $K$ has odd prime-power order and complement has prime order $q\ge 3$ that divides $|K|-1$;
\item[(iii)] $X$ is isomorphic to the Chein extension $M(G,2)$ of some elementary group $G$;
\item[(iv)] $X$ is the semidirect product $K\rtimes H$ of a Boolean group $K$ and an elementary Moufang loop $H$ of odd order;
\item[(v)] $X$ is the semidirect product $K\rtimes H$ of an elementary Abelian group $A$ of odd order and a subgroup $B\subseteq X$, isomorphic to the alternating groups $A_4$ or $A_5$;
\item[(vi)] $X$ is isomorphic to $M(q)$ for some $q\in \{2,3,4,5\}$.
\end{enumerate}

(i-ii) The cases (i) and (ii) imply the cases (1) and (2) of Corollary~\ref{c:elementary-Moufang=>cardinality}, respectively.
\smallskip

(iii) If the case (iii) holds, then the loop $X$ is isomorphic to the Chein extension $M(G,2)$ of some elementary group $G$. By Theorem~\ref{t:elementary-group<=>}, the elementary group $G$ satisfies the conditions (1)--(3) of Theorem~\ref{t:elementary-group<=>}. If $G$ satisfies the condition (1), then $|G|=p^n$ for some prime number $p$ and some number $n\in\IN$. In this case $|X|=|M(G,2)|=2|G|=2p^n$. If $p=2$, then $|X|=2^{n+1}$ and the case (1) of Corollary~\ref{c:elementary-Moufang=>cardinality} holds. If $p\ne 2$, then the case (5) of Corollary~\ref{c:elementary-Moufang=>cardinality} holds. 

If the condition (2) of Theorem~\ref{t:elementary-group<=>} holds, then $|G|=p^nq$ for some $n\in\IN$, some prime numbers $p,q$ such that $q$ divides $p^n-1$. If $q=2$, then $p$ is odd and $|X|=|M(G,2)|=2|G|=4p^n$.
So, the case (6) of Corollary~\ref{c:elementary-Moufang=>cardinality} holds. If $q\ge 3$, then $|X|=|M(G,2)|=2p^nq$ and the case (7)  of Corollary~\ref{c:elementary-Moufang=>cardinality} holds. If the condition (3)  of Theorem~\ref{t:elementary-group<=>} holds, then $|X|=|M(G,2)|=2|G|=2|A_5|=120$ and the condition (10) of Corollary~\ref{c:elementary-Moufang=>cardinality} holds.
\smallskip

(iv) If (iv) holds, then  $X$ is the semidirect product $K\rtimes H$ of a Boolean group $K$ and an elementary Moufang loop $H$ of odd order. Since $K$ is Boolean, $|K|=2^n$ for some $n\in\IN$, by Corollary~\ref{c:Boolean-order}. By Theorem~\ref{t:elementary-Moufang<=>}, the elementary Moufang loop $H$ satisfies the condition (1) or (2) of Theorem~\ref{t:elementary-Moufang<=>}. 

If $H$ satisfies the condition (1), then $|H|=p^m$ for some $m\in\IN$ and some odd prime number $p$. By Proposition~\ref{p:di-Lagrange}, the order of any nonidentity element of the elementary Moufang loop divides $|H|=p^m$ and hence equals $p$. By Proposition~\ref{p:=>cyclic-divides-kernel-1}, $p$ divides $|K|-1=2^n-1$, which implies $n\ge 2$. If $n\ge 3$, then the case (3) of Corollary~\ref{c:elementary-Moufang=>cardinality} holds. If $n=2$, then the case (6) holds.

If $H$ satisfies the condition (2) of Theorem~\ref{t:elementary-Moufang<=>}, then $H$ is a Frobenius loop whose kernel has odd prime-power order $p^m$ and has prime index $q\ge 3$ that divides $p^m-1$ and hence is not equal to $p$. Since the Frobenius loop $H$ contains elements of orders $p$ and $q$, the prime numbers $p$ and $q$ divide the cardinality $|K|-1=p^n-1$, by Proposition~\ref{p:=>cyclic-divides-kernel-1}. Since $p\ne q$, the product $pq$ divides $p^n-1$ and we have the case (4) of Corollary~\ref{c:elementary-Moufang=>cardinality}.
\smallskip

(v) If the case (v) holds, then $X$ is the semidirect product $K\rtimes H$ of an elementary Abelian group $K$ of odd order and a subgroup $H\subseteq X$, isomorphic to the alternating groups $A_4$ or $A_5$. Find we assume that $H\cong A_4$. Since the abelian group $K$ is elementary and has odd order, it has cardinality $p^n$ for some $n\in\IN$ and some prime number $p\ge 3$. If $p=3$, then $|X|=|H|\cdot|K|=12\cdot 3^n=4\cdot 3^{n+1}$ and hence the condition (6) of  Corollary~\ref{c:elementary-Moufang=>cardinality} holds. So, we assume that $p\ge 5$. Since $H$ is isomorphic to  the alternating group $A_4$, $H$ is the semidirect product $B\rtimes T$ of a Boolean group $B\subseteq H$ of cardinality $|B|=4$ and a $3$-element group $T\subseteq H$.  Proposition~\ref{p:=>cyclic-divides-kernel-1} ensures that $|T|=3$ divides  $p^n-1$. Therefore, $|X|=|H|\cdot|K|=12p^n$ and the case (8) of Corollary~\ref{c:elementary-Moufang=>cardinality} holds. 
\smallskip

Next, assume that the group $H$ is isomorphic to $A_5$. Since the group $A_5$ contains elements of order $2$, $3$ and $5$, the numbers $2,3,5$ divide the cardinal $|H|-1=p^n-1$ (by Proposition~\ref{p:=>cyclic-divides-kernel-1}), which implies that $p\notin\{2,3,5\}$ and hence $p\ge 7$. Then the number $30=2\cdot 3\cdot 5$ divides the number $p^n-1$ and the case (9) of Corollary~\ref{c:elementary-Moufang=>cardinality} holds.
\smallskip

(vi) If $X$ is isomorphic to the simple Page loop $M(q)$ for some $q\in\{2,3,4,5\}$. In this case $|X|\in\{|M(q)|:q\in\{2,3,4,5\}\}=\{120,1080,16320,39000\}$ and the case (10) of Corollary~\ref{c:elementary-Moufang=>cardinality} holds.
\smallskip

Now we show that the numbers in the numbers in cases (1)--(10) Corollary~\ref{c:elementary-Moufang=>cardinality} all are distinct. 
\smallskip

(1) The numbers in the case (1) have a unique prime divisior, whereas the numbers in cases (2)--(10) have at least two distinct prime divisiors. 
\smallskip

(2) The numbers in the case (2) are odd whereas the numbers in cases (3)--(10) are even.
\smallskip

(3) Take any number $x$ that satisfying the condition (3). Then $x=2^np^m$ for some $n\ge 3$, $m\in\IN$ and an odd prime $p$ that divides $2^n-1$. The number $2^np^m$ has exactly two prime divisors whereas the numbers in cases (4), (9), (10) have at least three distinct prime divisors. Since $n\ge 3$, the number $x=2^np^m$ does not satisfy the conditions (5), (6), (8). Assuming that $x$ satisfies the condition (7), we conclude that $m=1$ and $p$ divides $2^{n-1}-1$. Thenn $p$ also divides $(2^n-1)-(2^{n-1}-1)=2^{n-1}$, which is impossible as $p$ is odd. 
\smallskip

(4) Take any number $x$ satisfying the condition (4). Then $x=2^np^mq$ for some numbers $n,m\in\IN$ and odd primes $p,q$ such that $q$ divides $p^m-1$ and $pq$ divides $2^n-1$, which implies $n\ge 4$. Therefore, the number $2^np^mq$ is divisible by $16$. Since the number  $x$ has three distinct prime divisors, it cannot be equal to numbers that satisfy the conditions (5) or (6) (and have at most two prime divisiors) or the numbers that satisfy the condition (7) (which are divisible by $4$ only when they have exactly two prime divisors) or the numbers in cases (8), (9) (which are not divisible by $16$). It remains to check that the number $x$ is not equal to numbers  $60,120,1080,16320,39000$ in case (9). The numbers $60=2^2\cdot 3\cdot 5$, $120=2^3\cdot 3\cdot 5$, $1080=2^3\cdot 3^3\cdot 5$, and $39000=2^3\cdot 3\cdot 5\cdot 13$ are not divisible by 16 and hence are not equal to $x$.  The number $x$ has three prime divisors and cannot be equal to the number $16320=2^6\cdot 3\cdot 5\cdot 17$ that has four prime divisors.
\smallskip

(5)  Any number $x=2p^n$ satisfying the condition (5) has two prime divisiors and cannot be equal to the numbers that satisfy one of the cases (8)--(10) and have at least three distinct prime divisors. The number $x=2p^n$ is not divisible by $4$ and hence cannot be equal to a number satisfying the condition (6). If a number satisfying the condition (7) has only two prime divisors (like $x$), then it is divisible by $4$ and hence cannot be equal to the number $x$. 
\smallskip

(6)  Any number $x=4p^n$ satisfying the condition (6) has two prime divisiors and cannot be equal to the numbers that satisfy one of the cases (8)--(10) and have at least three distinct prime divisors. If a number satisfying the condition (7) has only two prime divisors (like $x$), then it is divisible by $8$ and hence cannot be equal to the number $x$. 
\smallskip

(7) Take any number $x$ that satisfies the condition (7), which means that  $x=2p^nq$ for some $n\in\IN$, some prime number $p$, and some odd prime number $q$ that divides $p^n-1$. If $p=2$, then the number $x$ has only two prime divisors and hence $x$ cannot be equal to the numbers that satisfy the conditions (8)--(10) and have at least three distinct divisors. So, assume that $p\ge 3$. In this case the number $x$ is not divisible by $4$ and cannot be equal to the numbers that satisfy the conditions (8)--(10) and are divisible by $4$.
\smallskip

(8) Take any number $x$ that satisfies the condition (8). Then $x=12p^n=2^2\cdot 3\cdot p^n$ for some $n\in\IN$ and some prime number $p\ge 5$ such that $3$ divides $p^n-1$. Since the number $x$ has exactly three prime divisors, it cannot be equal to the numbers in the condition (9) that have four prime divisiors. 
Assuming that $x=60$, we conclude that $p=5$ and $n=1$. Since $3$ does not divide $p^1-1=4$, the equality $x=60$ is impossible. Since $x=12p^n$ is not divisible by $8$ and the numbers $120=2^3\cdot 3\cdot 5$, $1080=2^3\cdot 3^3\cdot 5$, $16320=2^6\cdot 3\cdot 5\cdot 17$ and $39000=2^3\cdot 3\cdot 5\cdot 13$ are divisible by $8$, the number $x$ cannot be equal to numbers in the condition (10).
\smallskip

(9) Take any number $x$ that satisfies the condition (9). Then $x=60 p^n=2^2\cdot 3\cdot 5\cdot p^n$ for some $n\in\IN$ and some prime number $p\ge 7$ such that $30$ divides $p^n-1$. The number $x$ is not equal to the numbers in the condition (10) since $x>60$ and $x=2^2\cdot 3\cdot 5\cdot p^n$ is not divisible by $8$ whereas the numbers  $120,1080,16320,39000$ are are divisible by $8$.
\end{proof}


\begin{exercise} Let $\Pi$ be the set of all odd prime numbers. For every number in the set
$$
\begin{aligned}
&\{p^n,2p^n,4p^n:n\in\w,\;p\in \{2\}\cup\Pi\}\cup\{p^nq,2p^nq:n\in\IN,\;p\in\{2\}\cup \Pi,\;q\in \Pi,\;q|(p^n-1)\}\\
&\cup\{60,120,1080, 16320,39000\}
\end{aligned}
$$
construct an elementary Moufang loop $X$ of cardinality $|X|=n$.
\smallskip

{\em Hint:} Look at Exercise~\ref{ex:Frobenius-affine}.
\end{exercise}

\begin{problem} Are there elementary Moufang loops whose cardinality satisfy  the conditions \textup{(3), (4), (8), (9)} of Corollary~{\em\ref{c:elementary-Moufang=>cardinality}}?
\end{problem}

Theorem~\ref{t:elementary-Moufang<=>} implies the following characterization of the alternative group $A_5$. 

\begin{corollarys}\label{c:A5<=>} Every elementary Moufang loop $X$ of order $|X|=60$ is isomorphic to the simple group $A_5$.
\end{corollarys}

\begin{proof}  By Theorem~\ref{t:elementary-Moufang<=>}, the elementary Moufang loop $X$ satisfies one of the following conditions:
\begin{enumerate}
\item[(3)] $X$ is isomorphic to the Chein extension $M(G,2)$ of some elementary group $G$.
\item[(4)] $X$ is the semidirect product $K\rtimes H$ of a Boolean group $K$ and an elementary Moufang loop of odd order.
\item[(5)] $X$ is the semidirect product $K\rtimes H$ of an elementary Abelian group $K$ of odd order and a subgroup $H\subseteq X$, isomorphic to the alternating groups $A_4$ or $A_5$.
\end{enumerate}

If $X$ is isomorphic to the Chein extension $M(G,2)$ of some elementary group $G$, then the elementary group $G$ has cardinality $|G|=|X|/2=30$, which contradicts Theorem~\ref{t:elementary-group<=>}. Therefore, the case (3)  is impossible.

If $X$ satisfies the condition (4), then $X$ is the semidirect product $K\rtimes H$ of a Boolean group $K$ and an elementary Moufang loop $H$ of odd order. In this case $|X|=60$ implies $|K|=4$ and $|H|=15$. By Corollary~\ref{c:elementary-odd<=>}, $|H|$ is either prime-power (which is not true) or $|H|=qp$ for a prime-power number $q$ and a prime number $p$ that divides $q$. But neither $3$ divides $5-1$ nor $5$ divides $3-1$. Therefore, the case (4) is also impossible.

Assume that $X$ is the semidirect product $K\rtimes H$ of an elementary Abelian group $K$ of odd order and a subgroup $H\subseteq X$, isomorphic to the alternating group $A_4$ or $A_5$. If the group $H$ is isomorphic to the alternating group $A_4$, then the elementary Abelian group $K$ has cardinality $|K|=|X|/|A_4|=60/12=5$ and hence is cyclic. Since the group $H$ is isomorphic to $A_4$, it contains a cyclic subgroup $C$ of order $|C|=3$. Consider the subloop $K\rtimes C$ of $X$, generated by the set $K\cup C$. Since the groups $K$ and $C$ are cyclic, the loop $K\rtimes C$ is $2$-generated and hence is a group, by the diassociativity of the Moufang loop $X$. The group $K\rtimes C$ is elementary, being a subgroup of the elementary group $X$. By Proposition~\ref{p:=>Frobenius} and \ref{p:Frobenius-prime-divides}, $|C|=3$ divides $|K|-1=5-1=4$, which is not true. Therefore, the group $H$ is isomorphic to the group $A_5$ and then the group $X$ of cardinality $|X|=60=|A_5|=|H|$ coincides with the group $H$ and is isomorphic to the simple group $A_5$.
\end{proof}

\begin{corollarys}\label{c:M120<=>} Every elementary Moufang loop $X$ of order $|X|=120$ is isomorphic to the simple Paige loop $M(2)$ or to the Chein extension $M(A_5,2)$ of the simple group $A_5$.
\end{corollarys}

\begin{proof}  By Theorem~\ref{t:elementary-Moufang<=>}, the elementary Moufang loop $X$ satisfies one of the following conditions:
\begin{enumerate}
\item[(3)] $X$ is isomorphic to the Chein extension $M(G,2)$ of some elementary group $G$.
\item[(4)] $X$ is the semidirect product $K\rtimes H$ of a Boolean group $K$ and an elementary Moufang loop of odd order.
\item[(5)] $X$ is the semidirect product $K\rtimes H$ of an elementary Abelian group $K$ of odd order and a subgroup $H\subseteq X$, isomorphic to the alternating groups $A_4$ or $A_5$.
\item[(6)] $X$ is isomorphic to the simple Paige loop $M(2)$.
\end{enumerate}

If $X$ is isomorphic to the Chein extension $M(G,2)$ of some elementary group $G$, then the elementary group $G$ has cardinality $|G|=|X|/2=60$ and is isomorphic to the alternating group $A_5$, by Corollary~\ref{c:A5<=>}. It remains to show that the cases (4) and (5) are impossible.

If $X$ satisfies the condition (4), then $X$ is the semidirect product $K\rtimes H$ of a Boolean group $K$ and an elementary Moufang loop $H$ of odd order. In this case $|X|=120$ implies $|K|=8$ and $|H|=15$. By Corollary~\ref{c:elementary-odd<=>}, $|H|$ is either prime-power (which is not true) or $|H|=qp$ for a prime-power number $q$ and a prime number $p$ that divides $q$. But neither $3$ divides $5-1$ nor $5$ divides $3-1$. Therefore, the case (4) is impossible.

Assume that $X$ is the semidirect product $K\rtimes H$ of an elementary Abelian group $K$ of odd order and a subgroup $H\subseteq X$, isomorphic to the alternating group $A_4$ or $A_5$. If the group $H$ is isomorphic to the alternating group $A_4$, then the elementary Abelian group $K$ has cardinality $|K|=|X|/|A_4|=120/12=10$, which is not possible. Therefore, the group $H$ is isomorphic to the group $A_5$ and then the elementary Abelian group $K$ has cardinality $|K|=|X|/|A_5|=120/60=2$. The normal $2$-element group $K$ is central in $X=K\rtimes H$ and hence the group $X$ contain cyclic subgroups of orders $6$ and $10$, which contradicts the elementarity of the loop $X$.
\end{proof}

By analogy we can deduce from  Theorem~\ref{t:elementary-Moufang<=>}  the following characterization of the simple Paige loops $M(q)$ for $q\in\{3,4,5\}$.

\begin{corollarys}\label{c:M(q)<=>} An elementary Moufang loop $X$ is isomorphic to the simple Paige loop $M(q)$ for $q\in\{3,4,5\}$ if and only if $|X|=|M(q)|$.
\end{corollarys}

Theorem~\ref{t:elementary-Moufang<=>}, \ref{t:Sylow-Moufang} and Corollaries~\ref{c:A5<=>}, \ref{c:M120<=>}, \ref{c:M(q)<=>} imply the following corollary.

\begin{corollarys}\label{c:Sylow<=>notM(q)} If a finite elementary Moufang loop $X$ has cardinality $$|X|\notin \{120,1080,16320,39000\}=\{|M(q)|:2\le q\le 5\},$$ then $X$ is group-solvable and has the Sylow property.
\end{corollarys}

Theorem~\ref{t:elementary-Moufang<=>} motivates the following open problems.

\begin{problem} Is a normal subloop $H$ of an elementary Moufang loop $X$ Boolean if the quotient loop $X/H$ is isomorphic to the alternating group $A_4$?
\end{problem}

\begin{problem} Is a normal subloop $H$ of an elementary Moufang loop $X$ trivial if the quotient loop $X/H$ is isomorphic to the alternating group $A_5$?
\end{problem}

\begin{problem} Is a normal Boolean subgroup $B$ of an elementary Moufang loop $X$ trivial if the index of $B$ is $X$ is odd and not prime?
\end{problem}

\begin{remark} In \cite{Rajah2001} Andrew Rajah constructed non-associative Moufang loops of order $p^3q$ for all odd prime numbers $p,q$ such that $q$ divides $p-1$. Computer calculations show that the Rajah loops are elementary for all $p<100$.
\end{remark}

\chapter{Invertible-add affine liners}\label{ch:invertible-add}

An affine liner is defined to be
\begin{itemize}
\item \index{invertible-add affine liner}\index{affine liner!invertible-add}\defterm{invertible-add} if it is invertible-plus and invertible-puls;
\item \index{inversive-add affine liner}\index{affine liner!inversive-add}\defterm{inversive-add} if it is inversive-plus and inversive-puls;
\item \index{associative-add affine liner}\index{affine liner!associative-add}\defterm{associative-add} if it is associative-plus and associative-puls;
\item \index{commutative-add affine liner}\index{affine liner!commutative-add}\defterm{commutative-add} if it is commutative-plus and commutative-puls.
\end{itemize}
 In this chapter we shall prove that a finite affine liner is invertible-add if and only if every triangle in $X$ is contained in an affine Pappian subplane of prime order. This implies that for every ternar $R$ of a finite invertible-add Playfair liner the loops $(R,+)$ and $(R,\!\puls\!)$ are elementary in the sense that they are unions of subgroups of prime order.

\section{Elementary ternars}


\begin{definition} Let $(R,+,\cdot,0,1)$ be a ring. A ternar $X$ is defined to be \index{$R$-elementary ternar}\index{ternar!$R$-elementary}\defterm{$R$-elementary} if there exists a function $h:R\to X$ such that $h(1)=1$ and $h((x{\cdot}y)+z)=h(x)_\times h(y)_+h(z)$ for every $x,y,z\in R$. If, moreover, the function $h$ is bijective, then the ternar $X$ is called \defterm{$R$-isomorphic}.\index{$R$-isomorphic ternar}\index{ternar!$R$-isomorphic}
\end{definition}

\begin{exercise} Let $R$ be a corps endowed with the ternary operation $x_\times y_+z\defeq (x\cdot y)+z$. Show that the ternar $R$ is $\IZ$-elementary.
\end{exercise}

\begin{definition} A ternar $X$ if defined to be \index{field-elementary ternar}\index{ternar!field-elementary}\defterm{field-elementary} if $X$ is $F$-elementary for some field $F$.
\end{definition}

In the following proposition, for a prime number $p$ we denote by $\IF_p$ the $p$-element field $\IZ/p\IZ$. By $\IF_\w$ we denote the field $\IQ$ of rational numbers. Therefore, for every $p\in\IP\cup\{\w\}$, $\IF_p$ is a minimal field of characteristic $p$. We denote by $\IP$  the set of all prime numbers.

\begin{theorem}\label{t:Zelementary1-6} For a ternar $X$, the following conditions are equivalent:
\begin{enumerate}
\item $X$ is $\IF_p$-elementary for some $p\in\IP\cup\{\w\}$;
\item $X$ is field-elementary;
\item the minimal subternar of $X$ is $\IF_p$-isomorphic for some $p\in\IP\cup\{\w\}$.
\end{enumerate}
If the minimal subternar of $X$ is finite, then the conditions \textup{(1)--(3)} are equivalent to
\begin{enumerate}
\item[\textup{(4)}] $X$ is $R$-elementary for some ring $R$.
\item[\textup{(5)}] $X$ is $\IZ$-elementary;
\item[\textup{(6)}] $X$ is $\IF_p$-elementary for some $p\in\IP$.
\end{enumerate}
\end{theorem}

\begin{proof} The implications $(3)\Ra(1)\Ra(2)$ are trivial. To see that $(2)\Ra(1)$, observe that the minimal subfield of every field is isomorphic to the  field $\IF_p$ for some $p\in\IP\cup\{\w\}$. 
\smallskip

$(1)\Ra(3)$ Assuming that the ternar $X$ is $\IF_p$-elementary for some $p\in\IP\cup\{\w\}$, find a function $h:\IF_p\to X$ such that $h(1)=1$ and $h((x{\cdot}y)+z)=h(x)_\times h(y)_+h(z)$ for all $x,y,z\in \IF_p$. 

We claim that map $h:\IF_p\to X$ is injective. To derive a contradiction, assume that $h(x)=h(y)$ for certain distinct elements $x,y\in \IF_p$. Since $\IF_p$ is a field, there exists an element $u\in \IF_p\setminus\{0\}$ such that $x=u+y$. Then $0+h(y)=h(x)=h(u{\cdot}1+y)=h(u)_\times 1_+h(y)=h(u)+h(y)$, which implies $0=h(u)$, by the cancellability of the plus loop $(X,+)$. Since $u$ is a non-zero element of the field $\IF_p$, there exist an element $v\in\IF_p$ such that $u\cdot v=1$. Then $1=h(1)=h(u{\cdot}v+0)=h(u)_\times h(v)_+h(0)=0_\times h(v)_+0$, by the axiom {\sf(T1)}. But the equality $1=0$ in $X$ contradicts the definition of a ternar.

Next, we prove that $H\defeq h[\IF_p]$ is a subternar of the ternar $X$. By Proposition~\ref{p:subternar<=>}, it suffices to check that $\{0,1\}\subseteq H$ and 
$$T[H^3]\cup T_2[H^3]\cup T_3[H^4_{\ne}]\cup T_4'[H^4_{\ne}]\cup T''_4[H^4_{\ne}]\subseteq H,$$
where $T:X^3\to X$, $T:xyz\mapsto x_\times y_+z$, is the ternary operation of the ternar $X$, and the operation $T_2:X^3\to X$ and $T_3,T_4',T_4'':X^4_{\ne}\to X$ are defined in Remark~\ref{rem:ternar-T2T3T4}.

The axiom {\sf (T2)} of a ternar ensures that $1_\times 1_+1\ne 1_\times 1_+0=1$ and hence $h(2)=h(1{\cdot}1+1)=h(1)_\times h(1)_+h(1)=1_\times 1_+1\ne 1=h(1)$. Then the set $H\defeq h[\IF_2]$ contains at least two elements. For every $x\in \IF_p$, we have the equality $h(x)_\times 0_+1=1=h(1)=h((x{\cdot}0)+1)=h(x)_\times h(0)_+h(1)=h(x)_\times h(0)_+1$, which implies that the set $E\defeq\{x\in X:x_\times 0_+1=x_\times h(0)_+1\}$ contains the subset $H=h[\IF_p]$ and hence has cardinality $|E|\ge|H|\ge 2$. Applying the axiom {\sf(T3)}, we conclude that $0=h(0)\in h[\IF_p]$.

To see that $T[H]^3\subseteq H$, take any elements $x,y,z\in H$ and find elements $x',y',z'\in \IF_p$ such that $x=h(x)$, $y=h(y')$, $z=h(z')$. Then $T(x,y,z)=x_\times y_+z=h(x')_\times h(y')_+h(z')=h((x'{\cdot}y')+z')\in h[\IF_p]=H$.

To see that $T_2[H^3]\subseteq H$, take any elements $x,a,y\in H$ and find elements $x',a',y'\in \IF_p$. Since $\IF_p$ is a field, there exists a unique element $b'\in\IF_p$ such that $(x'{\cdot}a')+b'=y'$. Then the element $b\defeq h(b')\in h[\IF_p]=H$ satisfies the equalities $x_\times a_+b=h(x')_\times h(a')_+h(b')=h((x'{\cdot}a')+b')=h(y')=y$ and hence $T_2(x,a,y)=b\in H$ by the uniqueness part of the axiom {\sf(T2)} and the definition of the function $T_2$.

To see that $T_3[H^4_{\ne}]\subseteq H$, take any quadruple $abcd\in H^4_{\ne}$ and find elements $a',b',c',d'\in \IF_p$ such that $h(a')=a$, $h(b')=b$, $h(c')=c$, $h(d')=d$. The inequality $a\ne c$ implies $a'\ne c'$. Since $\IF_p$ is a field, there exists an element $x'\in\IF_p$ such that $(x'{\cdot}a')+b'=(x'{\cdot}c')+d'$.  Then the element $x\defeq h(x')\in h[\IF_p]=H$ satisfies the equalities 
$$x_\times a_+b=h(x')_\times h(a')_+h(b')=h((x'{\cdot}a')+b')=
h((x'{\cdot}c')+d')=h(x')_\times h(c')_+h(d')=x_\times c_+d$$
and hence $T_3(a,b,c,d)=x\in H$, by the uniqueness part of the axiom {\sf(T3)} and the definition of the function $T_3$.

To see that $T_4'[H^4_{\ne}]\cup T_4''[H^4_{\ne}]\subseteq H$, take any quadruple $xyuv\in H^4$ and find elements $x',y',u',v'\in \IF_p$ such that $h(x')=x$, $h(y')=y$, $h(u')=u$, $h(v')=v$. The inequality $x\ne u$ implies $x'\ne u'$. Since $\IF_p$ is a field, there exist unique elements $a',b'\in\IF_p$ satisfying the equations $x'{\cdot}a'+b'=y'$ and $u'{\cdot}a'+b'=v'$. Then the elements $a\defeq h(a')\in H$ and $b\defeq h(b')\in H$ satisfy the equations 
$x_\times a_+b=h(x')_\times h(a')_+h(b')=h(x'{\cdot}a'+b')=h(y')=y$ and
$u_\times a_+b=h(u')_\times h(a')_+h(b')=h(u'{\cdot}a'+b')=h(v')=v$. 
The uniqueness part of the axiom {\sf(T4)} and the definition of the functions $T_4'$ and $T''_4$ ensure that $T_4'(x,y,u,v)=a\in H$ and $T_4''(x,y,u,v)=b\in H$.

Therefore, $T[H^3]\cup T_2[H^3]\cup T_3[H^4_{\ne}]\cup T_4'[H^4_{\ne}]\cup T_4''[H^4_{\ne}]\subseteq H$ and $H$ is a subternar of the ternar $X$, by Proposition~\ref{p:subternar<=>}. We claim that $H$ is contained in the minimal subternar $\underline X$ of $X$. Let $x_1\defeq 1\in\IF_p$ and $y_1\defeq 1\in \underline{X}$. For every $n\in\IN$, let $x_{n+1}\defeq 1+x_n=1{\cdot}1+x_n\in\IF_p$ and $y_{n+1}\defeq 1_\times 1_+y_n\in \underline{X}$. By induction it can be shown that $h(x_n)=y_n$ for all $n\in\IN$.

 If $p<\w$, then $\IF_p=\{x_n:1\le n\le p\}$ and hence $H=h[\IF_p]=\{y_n:1\le n\le p\}\subseteq\underline X$. If $p=\w$, then $\IF_p=\IF_\w$ is the field of rational numbers and hence for every element $x\in\IF_\w$ there exist numbers $n,m,k\in\IN$ such that $x_n$ and $x{\cdot}x_n+x_m=x_k$.
Since $\IF_\w$ is the field of rational numbers, $x_n\ne 0$ and hence there exists an element $x_n^{-1}\in\IF_\w$ such that $x_n\cdot x_{n}^{-1}=1$.
Assuming that $y_n=h(x_n)=0$ in $X$, we conclude that 
$$1=h(1)=h(x_n{\cdot} x_n^{-1}+0)=h(x_n)_\times h(x_n^{-1})_+h(0)=0_\times h(x_n^{-1})_+0=0,$$
which is a contradiction showing that $y_n\ne 0$ n $X$.
Applying the axiom {\sf(T3)} to the ternar $\underline X$, find an element $y\in\underline{X}$ such that $T(y,y_n,y_m)=y_k=T(y,0,y_k)$.  Then $$T(y,y_n,y_m)=y_k=h(x_k)=h(x{\cdot}x_n+x_m)=T(h(x),h(x_n),h(x_m))=T(h(x),y_n,y_m)$$and hence $h(x)=y\in H$ by the uniquness part of the axiom {\sf(T3)}. This completes the proof of the inclusion $h[\IF_p]\subseteq\underline X$. The minimality of the ternar $\underline X$ ensures that $H=h[\IF_p]=\underline X$ and hence $h:\IF_p\to\underline X$ is a well-defined bijective map witnessing that the mininal subternar $\underline X$ of $X$ is $\IF_p$-isomorphic.
\smallskip

Therefore, the conditions (1), (2), (3) are equivalent. The implications $(6)\Ra(1)\Ra(4)$ are trivial.

$(4)\Ra(5)$ Assume that the ternar $X$ is $R$-elementary for some ring $R$. Then there exists a function $h:R\to X$ such that $h(1)=1$ and $h(x{\cdot}y+z)=h(x)_\times h(y)_+h(z)$ for all $x,y,z\in R$. Since $(R,+)$ is a group, there exists a unique map $g:\IZ\to R$ such that $g(1)=1$ and $g(n+m)=g(n)+g(m)$ for all $n,m\in\IZ$. Using the distributivity of the multiplication over addition in the ring $R$, it is easy to show by induction that $g(n\cdot m)=g(n)\cdot g(m)$ for all $n,m\in\IZ$.
Then the map $h\circ g:\IZ\to X$ witnesses that the trenar $X$ is $\IZ$-elementary.
\smallskip

Finally, assuming that the minimal ternar $\underline X$ of $X$ is finite, we shall prove that $(5)\Ra(6)$. Assume that the ternar $X$ is $\IZ$-elementary. Then there exists a function $h:\IZ\to X$ such that  $h(1)=1$ and $h((x\cdot y)+z)=h(x)_\times h(y)_+h(z)$ for all $x,y,z\in \IZ$. It follows from $$h(0)+0=h(0)=h((0\cdot 1)+0)=h(0)_\times h(1)_+h(0)=h(0)_\times 1_+h(0)=h(0)+h(0)$$ that $h(0)=0$, by the cancellativity of the plus operation in the ternar $R$. Repeating the argument from the proof of the implication $(1)\Ra(3)$, one can show that $h[\w]\subseteq\underline X$.

\begin{lemma}\label{l:Z-elementary=>p} There exists a prime number $p$ such that $h(p)=0$ and $h(i)\ne h(j)$ for all distinct numbers $i,j\in\{0,\dots,p-1\}$.
\end{lemma}

\begin{proof} Since $\underline X\supset h[\IN]$ is finite, there exist two non-negative numbers $x<y$ such that $h(x)=h(y)$. We can choose the numbers $x,y$ so that the difference $p\defeq y-x$ is the smallest possible. Then $$h(x)+0=h(x)=h(y)=h((x\cdot 1)+p)=h(x)_\times h(1)_+h(p)=h(x)_\times  1_+h(p)=h(x)+h(p)$$and hence $h(p)=0$ because $(X,+)$ is a loop with neutral element $0$. Since $h(1)=1\ne 0=h(p)$, the number $p$ is not equal to $1$.

The minimality of $p$ ensures that $h(x)\ne 0$ for all $x\in\{0,\dots,p-1\}$. Assuming that the number $p$ is composite, we can find two positive numbers $n,m<p$ such that $p=n\cdot m$. Then $h(n)\ne 0\ne h(m)$ and $0=h(p)=h(n\cdot m+0)=h(n)_\times h(m)_+h(0)=h(n)\cdot h(m)$, which is not possible because $(X^*,\cdot)$ is a loop. This contradiction shows that the number $p$ is prime. The minimality of $p$ ensures that  $h(i)\ne h(j)$ for all  distinct numbers $i,j\in\{0,\dots,p-1\}$.
\end{proof}

By Lemma~\ref{l:Z-elementary=>p}, there exists a prime number $p$ such that $h(p)=0$. Since $p$ is prime, the quotient ring $\IF_p\defeq \{p \IZ+k:k\in\IZ\}$ is a field. Let $\pi:\IZ\to\IF_p$, $\pi:k\mapsto p\cdot\IZ+k$, be the quotient homomorphism. Observe that for every $n,k\in\IZ$, we have $h((p\cdot n)+k)=h(p)_\times h(n)_+h(k)=0_\times h(n)_+h(k)=h(k)$ and hence $h[p\IZ+k]=\{h(k)\}$. Therefore, there exists a function $\tilde h:\IF_p\to X$ such that $h=\tilde h\circ \pi$. Since $\tilde h(p\IZ+1)=h(1)=1$ and 
$$
\begin{aligned}
&\tilde h\big(\big((p\IZ+x)\cdot (p\IZ+y)\big)+(p\IZ+z)\big)=\tilde h(p\IZ+((x\cdot y)+z))=h((x\cdot y)+z)\\
&=h(x)_\times h(y)_+h(z)=\big(\tilde h(p\IZ+x)\cdot\tilde h(p\IZ+y)\big)+\tilde h(p\IZ+z),
\end{aligned}
$$
for all $x,y,z\in \IZ$, the function $\tilde h:\IF_p\to X$ witnesses that the ternar $X$ is $\IF_p$-elementary.
\end{proof}

\begin{corollary}\label{c:Z-elementary=>prime-har} Every $\IZ$-elementary ternar $X$ has characteristic $\har(X)\in\IP\cup\{\w\}$.
\end{corollary}

\begin{proof} Let $(X,T)$ be a $\IZ$-elementary ternar and let $\underline{X}$ be the minimal subternar of $X$. Proposition~\ref{p:subternar<=>} implies that $\har(X)=|\underline{X}|\le\w$. If $\underline{X}$ is infinite, then $\har(X)=|\underline{X}|=\w$ and we are done. So, assume that $\underline{X}$ is finite. By Theorem~\ref{t:Zelementary1-6}, the ternar $\underline X$ if $\IF_p$-isomorphic for some prime number $p$. Then $\har(R)=|\underline R|=|\IF_p|=p\in \IP$.
\end{proof}

\begin{problem} Is every $\IZ$-elementary minimal ternar a field?
\end{problem}

\section{Elementary Playfair liners}

\begin{definition} A Playfair liner $X$ is defined to be
\begin{itemize}
\item \index{Playfair liner!$R$-elementary} \defterm{$R$-elementary} for a ring $R$ if every ternar of $X$ is $R$-elementary;
\item \index{Playfair liner!field-elementary} \defterm{field-elementary} if every ternar of $X$ is field-elementary.
\end{itemize}
\end{definition}

Theorem~\ref{t:Zelementary1-6}  implies the following corollary.

\begin{corollary}\label{c:Z-elementary<=>field-elementary} A Playfair liner (of finite order) is $\IZ$-elementary if (and only if) it is field-elementary.
\end{corollary}

For an affine space $X$ we denote by $\underline{\IR}_X$ the minimal subfield of the scalar corps $\IR_X$ of $X$. 

\begin{theorem}\label{t:BT=>f-elementary} If an affine space $X$ is Boolean or Thalesian, then it is $\IR_X$-elementary and $\underline{\IR}_X$-elementary. Moreover, $X$ is unicharacteristic with characteristic range $\har[X^\vartriangle]=\{|\underline{\IR}_X|\}\subseteq\IP\cup\{\w\}$.
\end{theorem}

\begin{proof} It suffices to show that for every triangle $uow$ in $X$, its ternar $\Delta\defeq\Delta_{uow}$ is $\IR_X$-elementary, $\underline{\IR}_X$-elementary, and has characteristic $\har(\Delta)=|\underline{\IR}_X|$.

Let $e$ be the diunit of the affine base $ouw$ in the plane $\overline{\{u,o,w\}}$. Let 
$$\Ker(\Delta)\defeq\bigcap_{x,a,b\in \Delta}\{s\in R:s\cdot(x_\times a_+b)=(s\cdot x)_\times a_+(s\cdot b)\}$$be the kernel of the ternar $\Delta$. 

If $X$ is Thalesian, then by Theorems~\ref{t:paraD<=>translation}, \ref{t:Thalesian<=>linear} and \ref{t:Ker(R)=RPi}, the ternar $\Delta$ is linear and its kernel $\Ker(\Delta)$ endowed with the plus and dot operations is a corps, isomorphic to the scalar corps $\IR_X$ by the isomorphism $h:\Ker(\Delta)\to \IR_X$, $h:x\mapsto \overvector{oxe}$. Then the homomorphism $h^{-1}:\IR_X\to \Ker(\Delta)\subseteq\Delta$ witnesses that the ternar $\Delta$ is $\IR_X$-elementary. The restriction of $h^{-1}$ to to the minimal subfield $\underline{\IR}_X$ of the scalar corps $\IR_X$ witnesses that the ternar $\Delta$ is $\underline{\IR}_X$-elementary. 
The minimality of the subfield $\underline{\IR}_X$ ensures that its image $h^{-1}[\underline{\IR}_X]$ coincides with the minimal subternar $\underline \Delta$ of the ternar $\Delta$. Then $\har(\Delta)=|\underline{\Delta}|=|h^{-1}[\underline{\IR}_X]|=|\underline{\IR}_X\}|$.

Now assume that the Playfair liner $X$ is non-Thalesian but Boolean. By Theorem~\ref{t:RXne01=>paraD}, $\IR_X=\{0,1\}$. Since $X$ is Boolean, $e+e=o$ in the ternar $\Delta$ and hence the function $h:\IR_X\to \Delta$ defined by $h(0)=0$ and $h(1)=e$ witnesses that the ternar $\Delta$ is $\IR_X$-elementary. Since $\underline{\IR}_X=\IR_X=\{0,1\}$, the ternar $\Delta$ is $\underline{\IR}_X$-elementary. Taking into account that $h[\IR_X]=\{o,e\}$ is a subternar of the ternar $\Delta$, we conclude that $\underline{\Delta}=\{o,e\}$ and  $\har(\Delta)=|\underline{\Delta}|=2=|\IR_X|=|\underline{\IR}_X|$.

Since $\underline{\IR}_X$ is a minimal field, it has cardinality $|\underline{\IR}_X|\in\IP\cup\{\w\}$. Therefore, $\har[X^\vartriangle]=\{|\underline{\IR}_X|\}\subseteq\IP\cup\{\w\}$. 
\end{proof}

\begin{corollary}\label{c:BT=>f-elementary} If a Playfair liner $X$ is Boolean or Thalesian, then it is field-elementary.
\end{corollary}

\begin{proof} Given any ternar $R$ of $X$, find a triangle $uow$ whose ternar $\Delta_{uow}$ is isomorphic to $R$. The triangle $uow$ is an affine base for the plane $\Pi\defeq\overline{\{u,o,w\}}$ in $X$. Since $X$ is Boolean or Thalesian, so is the plane $\Pi$. By Theorem~\ref{t:BT=>f-elementary}, the Playfair plane $\Pi$ is $\underline{\IR}_{\Pi}$-elementary for the minimal subfield $\underline{\IR}_{\Pi}$ of the scalar corps $\IR_\Pi$ of the Playfair plane $\Pi$. Then the ternar $R\cong \Delta_{uow}$ is $\underline{\IR}_\Pi$-elementary and hence field-elementary.
\end{proof}

The following theorem is the principal result of this chapter.

\begin{theorem}\label{t:invertible-add<=>Z-elementary} A Playfair liner is invertible-add if and only if its is $\IZ$-elementary.
\end{theorem}

\begin{proof} To prove the ``if'' part, assume that every ternar $R$ of a Playfair liner $X$ is $\IZ$-elementary. 

To prove that $X$ is invertible-plus, take any points $a,b,c,\alpha,\beta,\gamma\in X$ and $o\in\Aline a\alpha\cap\Aline b\beta\cap\Aline c\gamma$ such that $\Aline ac\parallel \Aline ob\parallel \Aline o\beta\parallel\Aline\alpha\gamma\nparallel \Aline b\gamma\parallel\Aline o\alpha\parallel \Aline oa\parallel\Aline c\beta$ and $\Aline ab\parallel\Aline co\parallel \Aline o\gamma$. We have to prove that $\Aline co\parallel\Aline \alpha\beta$. 

\begin{picture}(150,80)(-180,-45)
\linethickness{0.75pt}
\put(0,0){\color{teal}\line(1,0){30}}
\put(0,0){\color{cyan}\line(0,1){30}}
\put(0,0){\color{teal}\line(-1,0){30}}
\put(0,0){\color{cyan}\line(0,-1){30}}

\put(0,30){\color{teal}\line(1,0){30}}
\put(30,0){\color{cyan}\line(0,1){30}}
\put(-30,-30){\color{teal}\line(1,0){30}}
\put(-30,-30){\color{cyan}\line(0,1){30}}
\put(-30,0){\color{red}\line(1,1){30}}
\put(0,-30){\color{red}\line(1,1){30}}
{\linethickness{1pt}
\put(0,0){\color{red}\line(1,1){30}}
\put(0,0){\color{red}\line(-1,-1){30}}
}
\put(0,0){\circle*{3}}
\put(2,-7){$o$}
\put(30,0){\circle*{3}}
\put(33,-2){$\alpha$}
\put(-30,0){\circle*{3}}
\put(-38,-2){$a$}
\put(0,30){\circle*{3}}
\put(-7,31){$b$}
\put(0,-30){\circle*{3}}
\put(2,-36){$\beta$}
\put(30,30){\circle*{3}}
\put(32,31){$\gamma$}
\put(-30,-30){\circle*{3}}
\put(-37,-36){$c$}
\end{picture}

If $\alpha=o$, then $\Aline o\alpha\parallel \Aline c\beta$ implies $c=\beta$  and $\Aline \alpha\beta=\Aline oc$. If $b=o$, then $\Aline ob\parallel \Aline o\beta\parallel \alpha\gamma$ implies $o=\beta$ and $\alpha=\gamma$. Then again $\Aline\alpha\beta=\Aline \gamma o\parallel \Aline co$. 

So, assume that $b\ne o\ne\alpha$. Then $\Aline \alpha\gamma\parallel\Aline bo\nparallel \Aline o\alpha\parallel \Aline b\gamma$ implies that $bo\alpha$ is an affine base in the plane $\Pi\defeq\overline{\{b,o,\alpha\}}$ and $\gamma$ is the diunit of this affine base $bo\alpha$. Let $\Delta=\Aline o\gamma$ be the ternar of the affine base $bo\alpha$. By our assumption, the ternar $\Delta$ is $\IZ$-elementary and hence there exists a function $h:\IZ\to \Delta$ such that $h(1)=\gamma$ and $h((n\cdot m)+k)=h(n)_\times h(m)_+ h(k)$ for all $n,m,k\in\IZ$.
Since $\Aline ac\parallel \Aline ob$, $a\in\Aline o\alpha$ and $\Aline ab\parallel \Delta$, the definition of the plus operation ensures that $c+\gamma=o=h(0)=h((-1)+h(1)=h(-1)+\gamma$ and hence $c=h(-1)$. Then $\gamma+c=h(1)+h(-1)=h(1+(-1))=h(0)=o$ and hence $\alpha\beta\parallel \Delta=\Aline oc$, by the definition of the plus operation. Therefore, the Playfair liner $X$ is invertible-plus.
\smallskip

To prove that $X$ is invertible-puls, take any distinct parallel lines $L,L'$ in $X$ and points $a,b,c\in L$ and $a',b',c'\in L'$ such that $\Aline {a}{a'}\parallel \Aline b{b'}\parallel \Aline c{c'}$ and $\Aline{a}{b'}\parallel\Aline b{c'}$. We have to prove that $\Aline {a'}b\parallel \Aline {b'}c$. 

If $b=c$, then $\Aline b{b'}\parallel \Aline c{c'}$ implies $\Aline b{b'}=\Aline c{c'}$ and $b'=c'$. Then $\Aline b{b'}=\Aline b{c'}\parallel \Aline a{b'}$ implies $\Aline b{b'}=\Aline a{b'}$ and $a=b$. Then $\Aline{a'}b=\Aline {a'}a\parallel \Aline{b'}b=\Aline {b'}c$ and we are done.  

So, assume that $c\ne b$. In this case $cbb'$ is an affine base whose diunit coincides with the point $c'$ and diagonal $\Delta$ with the line $\Aline b{c'}$. Since $\Delta\nparallel \Aline a{a'}$, there exists a unique point $\gamma\in \Delta\cap\Aline a{a'}$. 

\begin{picture}(150,95)(-180,-15)
\linethickness{0.75pt}
\put(0,0){\color{cyan}\line(0,1){60}}
\put(30,0){\color{cyan}\line(0,1){60}}
\put(0,0){\color{blue}\line(1,1){30}}
\put(-30,0){\color{blue}\line(1,1){60}}
\put(30,0){\color{red}\line(-1,1){30}}
\put(30,30){\color{red}\line(-1,1){30}}
\put(-30,0){\color{teal}\line(1,0){60}}
\put(0,30){\color{teal}\line(1,0){30}}
\put(0,60){\color{teal}\line(1,0){30}}

\put(-30,0){\circle*{3}}
\put(-34,-8){$\gamma$}
\put(0,0){\circle*{3}}
\put(-3,-8){$a$}
\put(0,30){\circle*{3}}
\put(-7,29){$b$}
\put(0,60){\circle*{3}}
\put(-3,63){$c$}
\put(30,0){\circle*{3}}
\put(29,-8){$a'$}
\put(30,30){\circle*{3}}
\put(33,28){$b'$}
\put(30,60){\circle*{3}}
\put(30,63){$c'$}
\end{picture}

By the definition of the plus operation, the parallelity relation $\Aline a{b'}\parallel \Delta$ implies $c'+\gamma=b$. Since the ternar $\Delta$ of the affine base $cbb'$ is $\IZ$-elementary, there exists a function $h:\IZ\to \Delta$ such that $h(1)=c'$ and $h((x\cdot y)+z)=h(x)_\times h(y)_+h(z)$ for all $x,y,z\in\IZ$. Since $b$ is the neutral element of the ternar $\Delta$,
$$ c'+\gamma=b=h(0)=h((1\cdot 1)+(-1))=h(1)_\times h(1)_+h(-1)=c'_\times c'_+h(-1)=c'+h(-1)$$ and hence $h(-1)=\gamma$. Then 
$$b=h(0)=h((1\cdot(-1))+1)=h(1)_\times h(-1)_+h(1)=c'_\times \gamma_+c'=\gamma\puls c'.$$
By the definition of the puls operation, the equality $\gamma\puls c'=b$ implies $\Aline b{a'}\parallel \Aline c{b'}$. This completes the proof of the ``if'' part of the theorem.
\smallskip
 
To prove the ``only if'' part, assume that a Playfair liner $X$ is invertible-add. Given any ternar $R$ of $X$, we should prove that $R$ is $\IZ$-elementary. Find a plane $\Pi\subseteq X$ and an affine base $uow$ in $\Pi$ whose ternar $\Delta$ is isomorphic to the ternar $R$. It suffices to prove that the ternar $\Delta$ is $\IZ$-elementary. Let $e$ be the diunit of the affine base $uow$, and $\boldsymbol h\defeq(\Aline ou)_{\parallel}$ and $\boldsymbol v\defeq(\Aline ow)_{\parallel}$ be the horizontal and vertical directions in the based affine plane $(\Pi,uow)$, respectively.

By Theorem~\ref{t:invertible-plus<=>power-associative}, the loop $(\Delta,+)$ is power-associative and hence there exists a function $f:\IZ\to\Delta$ such that $f(1)=e$ and $f(x+y)=f(x)+f(y)$, for all $x\in\IZ$. For all integer numbers $n,m\in\IZ$ consider the unique point $x_n^m\in\Aline {f(n)}{\boldsymbol v}\cap\Aline {f(m)}{\boldsymbol h}$ (with coordinates $f(n)$ and $f(m)$).
For every $n,m\in\IZ$, the equality $f(n+m)=f(n)+f(m)$ implies $\Aline {x_0^{m}}{x_n^{n+m}}\subparallel \Delta$. 



We divide the further proof into a series of claims.

\begin{claim}\label{cl:inv-puls=>} For all $n,m,k\in\IZ$, 
$\Aline {x_n^n}{x_{n+1}^{n+k}}\parallel \Aline {x_n^m}{x_{n+1}^{m+k}}.$
\end{claim}

\begin{proof} This claim can be proved by analogy with the proof of Theorem~\ref{t:invertible-puls<=>}, using the invertibility-puls of $X$.
\end{proof} 

\begin{claim}\label{cl:inv-add=1} For every $n,m,k\in\IZ$, 
$$\Aline {x_n^m}{x_{n+1}^{m+k}}=\Aline {x_{n+1}^{m+k}}{x_{n+2}^{m+2k}}.$$
\end{claim}

\begin{proof} Indeed, for $k=0$ we have $\Aline{x_n^m}{x_{n+1}^m}=\Aline{f(m)}{\boldsymbol v}=\Aline{x_{n+1}^m}{x_{n+2}^m}$.

 For $k=1$, consider the unique line $L\subset \Pi$ such that $x_0^{m-n}\in L$ and $L\parallel \Delta$. For every $i\in\{0,1,2\}$, the equality $$f(m+i)=f((n+i)+(m-n))=f(n+i)+f(m-n),$$ implies $x_{n+i}^{m+i}\in L$ and hence $\Aline{x_n^m}{x_{n+1}^{m+1}}=L=\Aline{x_{n+1}^{m+1}}{x_{n+2}^{m+2}}$.
For $k=-1$, we use the invertibility-puls of $X$ and the parallelity relations $\Aline {x_n^{m-1}}{x_{n+1}^{m}}\parallel \Aline {x_{n+1}^{m-1}}{x_{n+2}^{m}}$ to deduce that $\Aline {x_n^{m}}{x_{n+1}^{m-1}}\parallel \Aline {x_{n+1}^m}{x_{n+2}^{m-1}}$. 

\begin{picture}(100,120)(-150,-30)
\linethickness{0.8pt}

\put(0,0){\color{teal}\line(1,0){60}}
\put(0,30){\color{teal}\line(1,0){60}}
\put(0,0){\color{blue}\line(1,1){30}}
\put(30,0){\color{blue}\line(1,1){30}}
\put(0,60){\color{red}\line(1,-1){60}}
\put(0,30){\color{red}\line(1,-1){30}}
\put(0,0){\color{cyan}\line(0,1){60}}
\put(30,0){\color{cyan}\line(0,1){30}}
\put(60,0){\color{cyan}\line(0,1){30}}

\put(0,0){\circle*{3}}
\put(-8,-10){$x_n^{m-1}$}
\put(30,0){\circle*{3}}
\put(20,-10){$x_{n+1}^{m-1}$}
\put(60,0){\circle*{3}}
\put(63,-10){$x_{n+2}^{m-1}$}
\put(0,30){\circle*{3}}
\put(-15,27){$x_n^{m}$}
\put(30,30){\circle*{3}}
\put(27,35){$x_{n+1}^{m}$}
\put(60,30){\circle*{3}}
\put(63,30){$x_{n+2}^{m}$}
\put(0,60){\circle*{3}}
\put(-5,65){$x_n^{m+1}$}
\end{picture}

Claim~\ref{cl:inv-puls=>} ensures that 
$\Aline {x_n^m}{x_{n+1}^{m-1}}\parallel \Aline{x_n^n}{x_{n+1}^{n-1}}\parallel \Aline{x_n^{m+1}}{x_{n+1}^m}$. Then $\Aline {x_n^{m+1}}{x_{n+1}^m}\parallel \Aline {x_n^m}{x_{n+1}^{m-1}}\parallel \Aline{x_{n+1}^m}{x_{n+2}^{m-1}}$ and $\Aline {x_n^{m+1}}{x_{n+1}^m}=\Aline{x_{n+1}^m}{x_{n+2}^{m-1}}$. Therefore, the equality $\Aline {x_n^m}{x_{n+1}^{m+k}}=\Aline{x_{n+1}^{m+k}}{x_{n+2}^{m+2k}}$ holds for all $n,m\in\IZ$ and $k\in\{-1,0,1\}$.
 
Now assume that for some $k\ge 2$ we have proved that  for all $n,m,i\in \IZ$ with $|i|<k$ we have $\Aline {x_n^m}{x_{n+1}^{m+i}}=\Aline{x_{n+1}^{m+i}}{x_{n+2}^{m+2i}}$.

\begin{picture}(100,110)(-150,-30)
\linethickness{0.8pt}

\put(0,0){\color{blue}\line(1,0){60}}
\put(0,30){\color{blue}\line(1,0){60}}
\put(0,0){\color{red}\line(1,1){60}}
\put(30,0){\color{red}\line(1,1){30}}
\put(0,60){\color{violet}\line(1,-1){60}}
\put(0,30){\color{violet}\line(1,-1){30}}
\put(0,0){\color{cyan}\line(0,1){60}}
\put(30,0){\color{cyan}\line(0,1){30}}
\put(60,0){\color{cyan}\line(0,1){60}}

\put(0,0){\circle*{3}}
\put(-8,-10){$x_n^m$}
\put(30,0){\circle*{3}}
\put(16,-11){$x_{n+1}^{m+k-1}$}
\put(60,0){\circle*{3}}
\put(60,-10){$x_{n+2}^{m+2k-2}$}
\put(0,30){\circle*{3}}
\put(-25,28){$x_n^{m+1}$}
\put(30,30){\circle*{3}}
\put(18,43){$x_{n+1}^{m+k}$}
\put(60,30){\circle*{3}}
\put(63,28){$x_{n+2}^{m+2k-1}$}
\put(60,60){\circle*{3}}
\put(63,58){$x_{n+2}^{m+2k}$}
\put(0,60){\circle*{3}}
\put(-10,65){$x_{n}^{m+2}$}
\end{picture}

Then $\Aline{x_n^m}{x_{n+1}^{m+k-1}}=\Aline{x_{n+1}^{m+k-1}}{x_{n+2}^{m+2k-2}}$, $\Aline{x_n^{m+1}}{x_{n+1}^{m+k}}=\Aline{x_{n+1}^{m+k}}{x_{n+2}^{m+2k-1}}$, and $\Aline {x_n^{m+2}}{x_{n+1}^{m+k}}=\Aline {x_{n+1}^{m+k}}{x_{n+2}^{m+2k-2}}$. Claim~\ref{cl:inv-puls=>} ensures that $\Aline{x_n^{m+2}}{x_{n+1}^{m+k}}\parallel\Aline{x_n^{m+1}}{x_{n+1}^{m+k-1}}$ and hence $\Aline {x_n^{m+1}}{x_{n+1}^{m+k-1}}\parallel \Aline {x_n^{m+2}}{x_{n+1}^{m+k}}=\Aline{x_{n+1}^{m+k}}{x_{n+2}^{m+2k-2}}$. Applying the invertibility-puls of $X$ and Claim~\ref{cl:inv-puls=>}, we conclude that $\Aline {x_n^m}{x_{n+1}^{m+k}}\parallel \Aline {x_{n+1}^{m+k-1}}{x_{n+2}^{m+2k-1}}\parallel \Aline{x_{n+1}^{m+k}}{x_{n+2}^{m+2k}}$ and finally, $\Aline {x_n^m}{x_{n+1}^{m+k}}=\Aline {x_{n+1}^{m+k}}{x_{n+2}^{m+2k}}$.

By analogy we can prove the equality $\Aline {x_n^m}{x_{n+1}^{m-k}}=\Aline {x_{n+1}^{m-k}}{x_{n+2}^{m-2k}}$.
\end{proof}

By induction, Claim~\ref{cl:inv-add=1} implies the following claim.

\begin{claim}\label{cl:inv-add=2} For every $n,m,k\in\IZ$, $\Aline{x_0^k}{x_n^{mn+k}}\subparallel\Aline {x_0^k}{x_1^m}.$
\end{claim}

%

For every $n,m,k\in\IZ$, consider the lines $L_{f(m),o}\defeq\Aline{x_0^0}{x_1^m}$ and $L_{f(m),f(k)}\defeq\Aline {x_0^k}{x_1^{m+k}}$. Claim~\ref{cl:inv-add=1} implies $L_{f(m),f(k)}=\Aline {x_0^k}{x_1^{m+k}}\parallel \Aline{x_0^0}{x_1^m}=L_{f(m),o}$, and Claim~\ref{cl:inv-add=2} ensures that $x_n^{mn+k}\in L_{f(m),f(k)}$ and hence $f(n)_\times f(m)_+f(k)=f(nm+k)$, witnessing that the ternar $\Delta$ is $\IZ$-elementary.
\end{proof}

Theorem~\ref{t:invertible-add<=>Z-elementary} and Corollary~\ref{c:Z-elementary<=>field-elementary} imply the following corollary.

\begin{corollary}\label{c:invertible-add<=>field-elementary} A finite Playfair liner is invertible-add if and only if it is field-elementary.
\end{corollary}

\section{Characteristics of elementary Playfair liners}

In this section we analyze the structure of the characteristic range $\har[X^\vartriangle]$ of a $\IZ$-elementary Playfair liner. Many structural properties of $\har[X^\vartriangle]$ will be deduced from the following lemma,  where $\IP$ stands for the set of all prime numbres.

\begin{lemma}\label{l:har=order-plus-puls} Let $(X,uow)$ be a based affine plane, $e$ be its diunit, and $\Delta=\Aline oe$ be its ternar. If the affine plane $X$ is invertible-add, then the characteristic $p\defeq\har(\Delta)$ of the ternar $\Delta$ belongs to the set $\IP\cup\{\w\}$ and is equal to the order of the element $e$ in the plus and puls loops of the ternar $\Delta$.
\end{lemma}

\begin{proof} By Theorem~\ref{t:invertible-add<=>Z-elementary}, the ternar $\Delta$ is $\IZ$-elementary. By Corollary~\ref{c:Z-elementary=>prime-har}, $\har(\Delta)\in\IP\cup\{\w\}$. 

Let $e^+$ and $e^\heartsuit$ be the subloops of the plus and puls loops of the ternar $\Delta$, generated by the element $e$. We have to show that $\har(\Delta)=|e^+|=|e^\heartsuit|$.

By the $\IZ$-elementarity of the ternar $\Delta$, there exists a function $h:\IZ\to \Delta$ such that $h(1)=e$ and $h(x{\cdot}y+z)=h(x)_\times h(y)_+h(z)$ for all $x,y,z\in\IZ$. Observe that $h(1+n)=h(1{\cdot}1+n)=h(1)_\times h(1)_+h(n)=e_\times e_+h(n)=e+h(n)=e\puls h(n)$ for every $n\in\IZ$. By induction it can be shown that $h(n)\in E^+\cap E^\heartsuit\subseteq E^+\cup E^\heartsuit\subseteq\underline\Delta$ for all $n\in\w$. If the function $h{\restriction}_\w$ is injective, then $\har(\Delta)=|\underline{\Delta}|=|e^+|=|e^\heartsuit|=|h[\w]|=\w$ and we are done. If $h{\restriction}_\w$ is not injective, then (the proof of) Lemma~\ref{l:Z-elementary=>p} ensures that there exists a prime number $p$ and a bijective function $\tilde h: \IF_p\to \underline X$ such that $h=\tilde h\circ \pi$, where $\pi:\IZ\to\IF_p$, $\pi:x\mapsto x+p\IZ$, is the quotient homomorphism. Then $\underline\Delta=e^+=e^\heartsuit=\tilde h[\IF_p]$ and hence $\har(\Delta)=|\underline\Delta|=|e^+|=|e^\heartsuit|=p$. 
\end{proof}

The following theorem shows that the characteristic range of an invertible-add Playfair liner is equal to the set of all possible orders of elements in the plus and puls loops of ternars of the Playfair liner.

\begin{theorem}\label{t:har[X]<=>orders} Let $X$ be an invertible-add Playfair liner $X$. For a number   $p\in\IN\cup\{\w\}$, the following conditions are equivalent:
\begin{enumerate}
\item $p\in\har[X^\vartriangle]$;
\item \mbox{$p$ is equal to the order of some element in the plus and puls loops of some ternar $R$ of $X$;}
\item $p$ is equal to the order of some element in the plus or puls loop of some ternar $R$ of $X$.
\end{enumerate}
\end{theorem} 

\begin{proof} The implication $(1)\Ra(2)$ follows from Lemma~\ref{l:har=order-plus-puls}, and $(2)\Ra(3)$ is trivial. 
\smallskip

$(3)\Ra(1)$ Assume that the number $p$ is equal to the order of some element in the plus or puls loop of some ternar $R$ of $X$. Find a triangle $uow$ in $X$ whose ternar $\Delta_{uow}$ is isomorphic to the ternar $R$. Let $e$ be the diunit of the affine base $uow$ in the plane $\overline{\{u,o,w\}}$. It follows that for some operation $*\in\{+,\puls\!\}$, the loop $(\Delta_{uow},*)$ contains an element $x$ of order $p$.
 
If $*=+$, then find unique points $u'\in\Aline ou$ and $w'\in\Aline ow$ such that $u'ow'x$ is a parallelogram. Observe that the plus operations in the ternars $\Delta_{uow}$ and $\Delta_{u'ow'}$ coincide and hence $x$ has order $p$ in the loop $\Delta_{u'ow'}$. Since $x$ is the diunit of the affine base $u'ow'$, the order $p$ of $x$ equals the charactertistic $\har(\Delta_{u'ow'})$ of the ternar $\Delta_{u'ow'}$, by Lemma~\ref{l:har=order-plus-puls}. Therefore, $p=\har(\Delta_{u'ow'})\in\har[X^\vartriangle]$.

If $*=\puls$, then find unique points $w'\in\Aline ow$ and $e'\in \Aline ue$ such that $x\in \Aline {w'}{e'}\parallel\Aline ou$. Observe that the puls loops $(\Delta_{uow},\!\puls\!)$ and $(\Delta_{uow'},\!\puls\!)$ are isomorphic and the order of the point $x$ in the ternar $\Delta_{uow}$ is equal to the order of the point $e'$ in the ternar $\Delta_{uow'}$. Since $e'$ is the diunit of the affine base $uow'$, the order $p$ of $e'$ equals the charactertistic $\har(\Delta_{uow'})$ of the ternar $\Delta_{uow'}$, by Lemma~\ref{l:har=order-plus-puls}. 
Then $p=\har(\Delta_{uow'})\in\har[X^\vartriangle]$.
\end{proof}

Lemma~\ref{l:har=order-plus-puls}  implies the following two theorems.

\begin{theorem}\label{t:characteristic=>prime} If a Playfair liner $X$ is invertible-add, then $\har[X^\vartriangle]\subseteq\IP\cup\{\w\}$.
\end{theorem}

Let us recall that a loop $(X,\cdot)$ is \defterm{elementary} if every element of $X$ is contained in a cyclic subgroup of prime order in $X$.

\begin{theorem}\label{t:invertible-add=>pls-elementary} If a Playfair liner $X$ of finite order is invertible-add, then for every ternar $R$ of $X$, the plus and puls loops $(R,+)$ and $(R,\!\puls\!)$ of the ternar $R$ are elementary.
\end{theorem}

Corollaries~\ref{c:inv-plus=>Boolean-paralelogram}, \ref{c:inv-puls=>Boolean-paralelogram} and \ref{c:Boolean<=>2har} imply the following corollary.

\begin{corollary}\label{c:invertible=>2har} Let $X$ be a Playfair plane of even order. If $X$ is invertible-plus or invertible-puls, then $X$ contains a Boolean parallelogram and $2\in \har[X^\vartriangle]$.
\end{corollary}

\begin{theorems}\label{t:char-divides-order} If a Playfair liner $X$ is inversive-plus or inversive-puls, then any number $p\in\har[X^\vartriangle]$ divides the order $|X|_2$ of $X$.
\end{theorems}

\begin{proof} Given any cardinal number $p\in \har[X^\vartriangle]\subseteq\IN\cup\{\w\}$, find a triangle $uow$ in $X$ such that $p=\har(uow)$. If $|X|_2$ is infinite, then it is divisible by the cardinal number $p\in \IN\cup\{\w\}$. So, assume that $|X|_2$ is finite. Consider the ternar $\Delta\defeq\Delta_{uow}$ and its minimal subternar $\underline{\Delta}$. Let $T:\Delta^3\to\Delta$, $T:xyz\mapsto x_\times y_+ z$, be the ternary operation of the ternar $\Delta$. Proposition~\ref{p:subternar<=>} ensures that $\{0,1\}\subseteq\underline{\Delta}$ and $T[\underline{\Delta}]\subseteq\underline{\Delta}$. Observe that for every $x,y\in \underline{\Delta}$, we have $x+y=T(x,1,y)\in \underline{\Delta}$ and $x\puls y=T(1,x,y)\in \underline{\Delta}$. Therefore, $\underline{\Delta}$ is a submagma of the loops $(\Delta,+)$ and $(\Delta,\!\puls\!)$.

If $X$ is inversive-plus, then the plus loop $(\Delta,+)$ of the ternar $\Delta$ is Moufang, by Theorem~\ref{t:+inversive}. Since $\Delta$ is finite, the submagma $\underline{\Delta}$ is a subloop of the Moufang loop $(\Delta,+)$. By Theorem~\ref{t:Lagrange-Moufang}, the Moufang loop $(\Delta,+)$ has the Lagrange property and hence the order $|\underline{\Delta}|=\har(uow)=p$ of the subloop $\underline{\Delta}$ divides the order $|\Delta|=|X|_2$ of the Moufang loop $(\Delta,+)$.

If $X$ is inversive-puls, then the puls loop $(\Delta,\!\puls\!)$ of the ternar $\Delta$ is Moufang, by Theorem~\ref{t:inversive-puls}. Since $\Delta$ is finite, the submagma $\underline{\Delta}$ is a subloop of the Moufang loop $(\Delta,\!\puls\!)$. By Theorem~\ref{t:Lagrange-Moufang}, the Moufang loop $(\Delta,\!\puls\!)$ has the Lagrange property and hence the order $|\underline{\Delta}|=\har(uow)=p$ of the subloop $\underline{\Delta}$ divides the order $|\Delta|=|X|_2$ of the Moufang loop $(\Delta,\!\puls\!)$.

Therefore, in both cases, the number $p$ divides the order $|X|_2$ of the Playfair liner $X$.
\end{proof}

\begin{theorems}\label{t:har[X]=divisors} Let $X$ be a finite invertible-add Playfair plane and assume that $X$ is inversive-plus or inversive-puls.
\begin{enumerate}
\item If $|X|_2=120$, then $\har[X^\vartriangle]=\{2,3\}$ or $\har[X^\vartriangle]=\{2,3,5\}$.
\item If $|X|_2=1080$, then $\har[X^\vartriangle]=\{2,3\}$;
\item If $|X|_2\in\{16320,39000\}$, then $\har[X^\vartriangle]=\{2,3,5\}$;
\item If $|X|_2\notin\{120,1080,16320,39000\}$, then $\har[X^\vartriangle]$ is equal to the set of all prime divisors of the order $|X|_2$ of $X$.
\end{enumerate}
\end{theorems} 

\begin{proof} Assume that $X$ is inversive-plus or inversive-puls. By Theorems~\ref{t:har[X]<=>orders} and \ref{t:char-divides-order}, a number $p\in\IN$ belongs to the characteristic range $\har[X^\vartriangle]$ of $X$ if and only if $p$ is equal to the order of an element of the plus and/or puls loops of a suitable ternar $R$ of the Playfair plane $X$ only if $p$ is a divisor of the number $|X|_2$. Let $R$ be any ternar of the Playfair plane $X$. Since $X$ is inversive-plus or inversive-puls, for some operation $*\in\{+,\puls\}$, the loop $(R,*)$ is Moufang, according to Theorems~\ref{t:+inversive} and \ref{t:inversive-puls}. By Theorem~\ref{t:invertible-add=>pls-elementary}, the finite Moufang loop $(R,*)$ is elementary.
\smallskip

1. If $|R|=|X|_2=120=2^3{\cdot}3{\cdot}5$, then the elementary Moufang loop $(R,*)$ is isomorphic either to the simple Paige loop $M(2)$ or to the Chein extension $M(A_5,2)$ of the alterating group $A_5$, see Corollary~\ref{c:M120<=>}. Taking into account that  $M(2)$ and $M(A_5,2)$ contain elements of orders $2$ and $3$, we conclude that $\{2,3\}\subseteq \har[X^\vartriangle]\subseteq \{2,3,5\}$, by Theorems~\ref{t:har[X]<=>orders} and \ref{t:char-divides-order}.
\smallskip

2. If $|R|=|X|_2=1080=|M(3)|$, then the elementary Moufang loop $(R,*)$ is isomorphic to the simple Paige loop $M(3)$, according to Corollary~\ref{c:M(q)<=>}. By Theorem~\ref{t:ordersM(q)}, $\{2,3\}$ is the set of orders of elements in the loop $M(3)$. Then $\har[X^\vartriangle]=\{2,3\}$, by Theorem~\ref{t:har[X]<=>orders}.
\smallskip

3. If $|R|=|X|_2\in\{16320,39000\}=\{|M(4)|,|M(5)|\}$, then the elementary Moufang loop $(R,*)$ is isomorphic to the simple Moufang loop $M(q)$ for some $q\in\{4,5\}$, according to Corollary~\ref{c:M(q)<=>}. By Theorem~\ref{t:ordersM(q)}, $\{2,3,5\}$ is the set of orders of elements in the loop $M(q)$. Then $\har[X^\vartriangle]=\{2,3,5\}$, by Theorem~\ref{t:har[X]<=>orders}.
\smallskip

4. If $|R|=|X|_2\notin\{120,1080,16320, 39000\}$, then the elementary Moufang loop $(R,*)$ has the Sylow property, by Corollary~\ref{c:Sylow<=>notM(q)}. Then the set of orders of elements of the elementary loop $(R,*)$ coincides with the set of $D$ of prime divisors of the number $|R|=|X|_2$, and $\har[X^\vartriangle]=D$, by Theorem~\ref{t:har[X]<=>orders}.
\end{proof}






Theorems~\ref{t:+inversive}, \ref{t:inversive-puls}, \ref{t:invertible-add=>pls-elementary} and Corollary~\ref{c:elementary-Moufang=>cardinality} imply the following corollary describing possible orders of invertible Playfar planes which are inversive-plus or inversive-puls.

\begin{corollarys}\label{c:inversive=>order} If an invertible-add finite Playfair plane $X$ is inversive-plus or inversive-puls, then its order $|X|_2$ is equal to one of the numbers:
\begin{enumerate}
\item $p^n$ for some prime number $p$ and some positive integer $n$;
\item $p^nq$ for some $n\in\IN$ and some odd primes $p,q$ such that $q$ divides $p^n-1$;
\item $2^np^m$ for $n\ge 3$, $m\in\IN$ and some odd prime number $p$ that divides $2^n-1$;
\item $2^np^mq$ for some $n,m\in\IN$ and some odd prime numbers $p,q$ such that $q$ divides $p^m-1$ and $pq$ divides $2^n-1$;
\item $2p^n$ for some odd prime number $p$ and some $n\in\IN$;
\item $4p^n$ for some odd prime number $p$ and some $n\in\IN$;
\item $2p^nq$ for some $n\in\IN$ and some prime number $p$ and odd prime number $q$ that divides $p^n-1$;
\item $12p^n$ for some $n\in\IN$ and some odd prime number $p\ge 5$ such that $3$ divides $p^n-1$;
\item $60 p^n$ for some $n\in\IN$ and some prime number $p\ge 7$ such that $30$ divides $p^n-1$;
\item $60,120,1080,16320$ or $39000$.
\end{enumerate}
\end{corollarys} 

\begin{remark} In Theorem~\ref{t:not60} we shall prove that  an inversive-plus invertible-puls Playfair liner $X$ has order $|X|_2\ne 60$. 
\end{remark}

Theorem~\ref{t:har[X]=divisors} and Corollary~\ref{c:inversive=>order} imply the following classification of characteristic ranges of inversive-add Playfair liners.

\begin{corollarys} Let $X$ be an invertible-add Playfair liner $X$ of finite order. If $X$ is inversive-plus or inversive-puls, then $\har[X^\vartriangle]\in\big\{\varnothing,\{p,q\},\{2,p,q\},\{2,3,5,p\}:p,q\in\IP\big\}$.
\end{corollarys}

\begin{remark} A well-known (and still unproven) Prime-Power Conjecture states that the order of any finite Playfair plane is a prime power. So, if this conjecture is true (at least for inversive-plus or inversive-puls invertible-add Playfair planes), then the cases  \textup{(2)--(10)} of Corollary~\ref{c:inversive=>order} do not occur and $\har[X^\vartriangle]\in\{\varnothing,\{p\}:p\in\IP\}$ for any invertible-add Playfair liner $X$, which is inversive-plus or inversive-puls.
\end{remark}

\begin{propositions}\label{p:unichar+inversive=>prime-power} Let $X$ be an  invertible-add finite Playfair plane $X$, and assume that $X$ is inversive-plus or inversive-puls. The Playfair plane $X$ is unicharacteristic if and only if $|X|_2$ is a prime power.
\end{propositions}

\begin{proof} By Theorem~\ref{t:characteristic=>prime}, $\har[X^\vartriangle]\subseteq\IP$. 
If $X$ is unicharacteristic, then $\har[X^\vartriangle]=\{p\}$ for some prime number $p$. If $p=2$, then $X$ is Boolean and has order $|X|_2=2^n$ for some $n\in\IN$, by Corollaries~\ref{c:Boolean<=>char2} and \ref{c:Boolean-order}. So, assume that $p>2$. Corollary~\ref{c:invertible=>2har} ensures that $X$ has odd order $|X|_2$. By Theorem~\ref{t:har[X]=divisors}(4), $|X|_2=p^n$ for some $n\in\IN$.

To prove the ``if'' part, assume that $|X|_2=p^n$ for some $p\in\IP$ and $n\in\IN$. By Theorem~\ref{t:char-divides-order}, every number $c\in\har[X^\vartriangle]$ is a prime divisor of $|X|_2=p^n$ and hence $c=p$ and $\har[X^\vartriangle]=\{p\}$, witnessing that the Playfair plane $X$ is unicharacteristic.
\end{proof}

Proposition~\ref{p:unichar+inversive=>prime-power} and 
Theorems~\ref{t:add-ass=>partialT}, \ref{t:ass-puls=>partialT} imply the following corollary.

\begin{corollarys}\label{c:unichar+ass=>delta-trans} If a unicharacteristic invertible-add Playfair liner $X$ has finite order and is associative-plus or associative-puls, then $X$ is $\partial$-translation.
\end{corollarys}

\begin{proposition}\label{p:commutative-pls=>prime-power} If a finite invertible-add Playfair plane $X$ is commutative-plus or commutative-puls, then $X$ is unicharacteristic, $\partial$-translation, and has a prime-power order.
\end{proposition}

\begin{proof} By Theorems~\ref{t:add-com=>add-ass} and \ref{t:com-puls=>ass-puls}, the commutative-plus or commutative-puls Playfair plane $X$ is associative-plus or associative-puls and inversive-plus or inversive-puls. Let $R$ be any ternar of the Playfair plane $X$. Since $X$ is commutative-plus or commutative-puls, for some operation $*\in\{+,\puls\}$, the loop $(R,*)$ is a commutative and associative, according to Theorems~\ref{t:add-ass<=>},  \ref{t:add-com<=>}, \ref{t:add-com=>add-ass}, \ref{t:associative-puls<=>}, \ref{t:commutative-puls<=>}, \ref{t:com-puls=>ass-puls}. By Theorem~\ref{t:invertible-add=>pls-elementary}, the commutative group $(R,*)$ is elementary and hence has order $|R|=p^n$ for some $p\in\IP$ and $n\in\IN$. By Theorems~\ref{t:characteristic=>prime} and \ref{t:char-divides-order}, $\har[X^\vartriangle]=\{p\}$, which means that the Playfair plane $X$ is unicharacteristic.  By Corollary~\ref{c:unichar+ass=>delta-trans}, the Playfair plane $X$ is $\partial$-translation.
\end{proof} 

\begin{propositions}\label{p:a-Pappian-prime-order} For a finite Playfair plane $X$, the following conditions are equivalent:
\begin{enumerate}
\item $X$ is Pappian of prime order;
\item $X$ has prime order and is inversive-plus or inversive-puls;
\item $X$ is inversive-add of prime order;
\item $X$ is invertible-add and minimal.
\end{enumerate}
\end{propositions} 

\begin{proof} The equivalences $(1)\Leftrightarrow(2)\Leftrightarrow(3)$ follow from Theorem~\ref{t:finite-shear<=>Pappian}. 
\smallskip

$(3)\Ra(4)$ Assume that $X$ is inversive-add and has prime order $p$. Theorem~\ref{t:char-divides-order} ensures that every number $c\in \har[X^\vartriangle]$ divides the prime number $|X|_2=p$ and hence $c=p$ and $|X_2|\in\{p\}=\har[X^\vartriangle]$, witnessing that the Playfair plane $X$ is minimal. Since every inversive loop is invertible, the inversive-add Playfair plane $X$ is invertible-add.
\smallskip

$(4)\Ra(1)$ Assume that $X$ is invertible-add and minimal. Then $|X|_2\in\har[X^\vartriangle]$. Theorem~\ref{t:characteristic=>prime} ensures that $|X|_2\in\har[X^\vartriangle]\subseteq\IP$ and hence $|X|_2$ is a prime number. 
By Theorems~\ref{t:invertible-add<=>Z-elementary} and Corollary~ \ref{c:Z-elementary<=>field-elementary}, the Playfair plane $X$ is field-elementary. Since $X$ is minimal, for some triangle $uow$ in $X$, the ternar $\Delta_{uow}$ is minimal. By Theorem~\ref{t:Zelementary1-6}, the minimal ternar $\Delta_{uow}$ if $\IF_p$-isomorphic for some $p\in \IP$, and hence $\Delta_{uow}$ is linear, distributive, associative, and commutative.
By Theorem~\ref{t:commutative-dot<=>}, the Playfair plane $X$ is Pappian.
\end{proof}

\section{Some properties of the alternating group $A_5$}

In this section we discuss some properties of the alternating group $A_5$ which will be used in the proof of Kegel--L\"uneburg Theorem~\ref{t:Kegel-Luneburg-main} in Section~\ref{s:Kegel-Luneburg}. We recall that for every $n\in\IN$, the \defterm{alternating group} $A_n$ is a subgroup of the symmetric group $S_n$ consisting of all permutations that are compositions of even numbers of transpositions. A permutation $\pi$ of a set $X$ is called a \defterm{transposition} if $|\{x\in X:\pi(x)\ne x\}|=2$.

The automorphism groups of the alternating groups were described by \index[person]{H\"older}Otto H\"older\footnote{{\bf Otto H\"older} (1859 -- 1937) was a German mathematician best known for his foundational contributions to group theory, algebra, and analysis. He is most famous for H\"older's inequality, a key result in functional analysis, and for the Jordan--H\"older theorem, a fundamental result in group theory concerning composition series. H\"older studied under Leopold Kronecker and Karl Weierstrass and earned his doctorate from the University of Berlin in 1882. He spent much of his academic career at the University of Leipzig, where he influenced the development of modern algebra. His work helped lay the groundwork for the rigorous structure theory of groups and contributed to the arithmetization of analysis.} in \cite{Holder1895}.

\begin{theorems}[H\"older, 1895]\label{t:HolderAn} For every $n\in\IN$, let $\alpha:S_n\to\Aut(A_n)$ be the homomorphism assigning to every $a\in S_n$ the automorphism $\alpha^a:A_n\to A_n$, $\alpha^a:x\mapsto a^{-1}xa$. If $n\ne 6$, then the homomorphism $\alpha:S_n\to\Aut(A_n)$ is an isomorphism. If $n=6$, then $\alpha[S_n]$ is a subgroup of index $2$ in the automorphism group $\Aut(A_n)$.
\end{theorems} 

\begin{exercise}\label{ex:orders} Given an element $c\in S_5$, consider the automorphism $\alpha^c:A_5\to A_5$, $\alpha^c:x\mapsto c^{-1}xc$.
Show that the set $\Fix(\alpha^c)\defeq\{x\in A_5:c^{-1}{\cdot}x{\cdot}c=x\}$ has cardinality
\begin{enumerate}
\item $60$ if $c$ is the identity permutation;
\item $6$ if $c\in S_5\setminus A_5$;
\item $5$ is $c$ has order $5$;
\item $3$ if $c$ has order $3$;
\item $4$ if $c\in A_5$ has order $2$.
\end{enumerate}
\end{exercise}

\begin{lemmas}\label{l:A5xA5}  Let $A,B$ be two proper subgroups of a group $G$ such that $G=A{\cdot}B$ and the groups $A,B$ are isomorphic to $A_5$. If the subgroup $A$ is normal in $G$, then $G$ is isomorphic to $A_5\times A_5$.
\end{lemmas}

\begin{proof} Since $A\cap B$ is a normal subgroup of the simple group $B\cong A_5$, either $A\cap B=\{e\}$ or $A\cap B=B$. In the latter case, $G=A=B$, which contradicts the choice of the {\em proper} subgroups $A,B$. Therefore, $A\cap B=\{e\}$ and hence $G=A\rtimes B$ is the semidirect product of the subgroups $A$ and $B$.

Consider the homomorphism $\alpha:B\to\Aut(A)$ assigning to every $b\in B$ the automorphism $\alpha^b:A\to A$, $\alpha^b:a\mapsto a^b\defeq b^{-1}{\cdot}a{\cdot}b$, of $A$. Since the group $B$ is simple, the homomorphism $\alpha$ is either trivial or injective.
If the homomorphism $\alpha$ is trivial, then $G=A\times B$ is the direct sum of the subgroups $A,B$ and we are done. So, assume that the homomorphism $\alpha$ is injective. By H\"older's Theorem~\ref{t:HolderAn}, the automorphism group $\Aut(A)\cong\Aut(A_5)$ is isomorphic to the symmetric group $S_5$. Since the group $B$ is simple of cardinality $>2$, the image $\alpha[B]$ is contained in the normal subgroup of inner automorphisms $\mathsf{Inn}(A)=\{\alpha^a:a\in A\}$ in the automorphism group $\Aut(A)$. Since $|B|=60=|\mathsf{Inn}(A)|$, the homomorphism $\alpha:B\to\mathsf{Inn}(A)$ is an isomorphism. Then for every $b\in B$ there exists a unique element $a\in A$ such that $b^{-1}xb=a^{-1}xa$ for all $x\in A$. The latter equality implies $ab^{-1}x=xab^{-1}$ and hence $ab^{-1}\in\mathcal Z_G(A)\defeq\{z\in G:\forall x\in A\;z{\cdot}x=x{\cdot}z\}$. This implies $A{\cdot} Z_G(A)=A{\cdot}B=G$. Since the group $A\cong A_5$ has no central elements, $A\cap Z_G(A)=\{e\}$. Consider the quotient homomorphism $\pi:G\to G/A$ and observe that both $\pi{\restriction}_{Z_G(A)}:Z_G(A)\to G/A$ and $\pi{\restriction}_B:B\to G/A$ are  isomorphisms. Then $G=A\times Z_G(A)$ is isomorphic to the group $A_5\times A_5$.
\end{proof}

\begin{definition} A group homomorphism $h:D\to G$ is called  
\begin{itemize}
\item a \defterm{cover} of $G$ if its kernel is a subgroup of $\mathcal Z(D)\cap[D,D]$;
\item a \index{Schur cover}\defterm{Schur cover} of $G$ if $h$ is a cover and for every cover $\hbar:E\to G$ of $G$ there exists a unique homomorphism $\pi:D\to E$ such that $h=\hbar\circ\pi$.
\end{itemize}
\end{definition}

\begin{theorems}[Schur, 1904] For every (finite) group $G$, there exists a (finite) Schur cover of $G$. Moreover, any two Schur covers of a group $G$ are isomorphic.
\end{theorems}

\begin{definition} The \index{Schur multiplier}\defterm{Schur multiplier} $H^2(G,\IZ)$ of a group $G$ is the 
kernel of the Schur cover of $G$.
\end{definition}

\begin{theorems}[Schur, 1911]\label{t:Shur-multipliers} Let $n\in\IN$. The Schur multiplier $H^2(A_n,\IZ)$ 
\begin{enumerate}
\item is trivial if $n\le 3$;
\item is a cyclic group of order $6$ if $n\in\{6,7\}$;
\item has order $2$, otherwise.
\end{enumerate}
\end{theorems}
 
\begin{theorems}[Schur, 1907]\label{t:ShurPSL-SL} Let $p$ be an odd prime and $n\in\IN$. If a cover $h:G\to PSL(2,p^n)$ of the group $PSL(2,p^n)$ has kernel of size $2$, then the group $G$ is isomorphic to the group $SL(2,p^n)$.
\end{theorems}

With the help of \index[person]{Schur}Schur\footnote{{\bf Issai Schur} (1875--1941) was a Belarusian-born mathematician who became one of the leading figures in early 20th-century algebra and representation theory. Born in Mohyliv, then part of the Russian Empire, Schur studied mathematics in Berlin under Georg Frobenius and Lazarus Fuchs, earning his doctorate in 1901. He held professorships in Berlin and later at the University of Bonn, where he developed many foundational results in group theory and representation theory, including the Schur lemma, Schur decomposition, Schur multipliers, and the celebrated Schur--Weyl duality connecting group representations with linear algebra. His work influenced the development of character theory, the theory of symmetric and general linear groups, and the broader study of linear representations. After the rise of the Nazis, Schur, being Jewish, faced persecution and emigrated to Palestine in 1939. He continued his mathematical work until his death in Tel Aviv in 1941. Schur’s contributions remain central to modern algebra, combinatorics, and the representation theory of finite and Lie groups.} multipliers we now can prove the following theorem that will be used in the proof of Kegel-L\"uneburg Theorem~\ref{t:Kegel-Luneburg-main}.
 
\begin{theorems}[Kegel--L\"undeburg, 1963]\label{t:A5xA5=} Let $A,B$ be two proper subgroups of a group $G$ such that $G=A{\cdot}B$. If the groups $A,B$ are isomorphic to $A_5$, then the group $G$ is isomorphic to the group $A_6$ or $A_5\times A_5$.
\end{theorems}

\begin{proof} To derive a contradiction, assume that there exists a group $G$ containing two proper subgroups $A,B$ such that $G=A{\cdot}B$ and $A,B$ are isomorphic to $A_5$ but $G$ is not isomorphic to the groups $A_6$ or $A_5\times A_5$. We can assume that the group $G$ has the smallest possible cardinality. 

\begin{claim} The group $G$ coincides with its commutator subgroup $[G,G]$.
\end{claim}

\begin{proof} Consider the quotient homomorphism $h:G\to G/[G,G]$ and observe that the restrictions $h{\restriction}_A$ and $h{\restriction}_B$ of $h$ to the simple groups $A,B$ are constant homomorphisms. Then $A,B\subseteq [G,G]$ and hence $G=A{\cdot}B\subseteq[G,G]$, which implies $G=[G,G]$.
\end{proof}

Let $H$ be a normal subgroup of $G$ of smallest possible cardinality $|H|>1$. If the group $G$ is not simple, then $|H|<|G|$ and the quotient group $G/H$ has cardinality $|G/H|<|G|$. Consider the quotient homomorphism $\pi:G\to G/H$.

\begin{claim} $|H|\le 60$.
\end{claim}

\begin{proof} Since $H\subseteq G=A{\cdot}B$, for every $x\in H$, there exist elements $a_x\in A$ and $b_x\in B$ such that $x=a_x{\cdot}b_x$. Assuming that $|H|>60$, we can apply the Pigeonhole Principle and find two distinct elements $x,y\in H$ such that $a_x=a_y$. Then $x b_x^{-1}=a_x=a_y=yb_y^{-1}$ and hence $e\ne y^{-1}x=b_y^{-1}b_x\in H\cap B$. Since the group $B$ is simple, $H\cap B=B$ and hence $B\subseteq H$. By analogy we can show that $A\subseteq H$ and hence $G=A{\cdot}B\subseteq H$, which contradicts the choice of $H$ as a proper subgroup of $G$. This contradiction shows that $|H|\le 60$.
\end{proof} 

\begin{claim} The group $H$ is solvable.
\end{claim}

\begin{proof} If $H$ is not solvable, then $|H|\le 60$ implies that $H$ is isomorphic to the alternating group $A_5$, which is a unique non-solvable group of order 60. Since $H\cap A$ is a normal subgroup of the simple group $A$, either $A\cap H=\{e\}$ or $A\cap H=A$. In the latter case, $A\subseteq H$ and $|H|\le 60=|A|$ implies $A=H$ is a normal subgroup of $G$. Applying Lemma~\ref{l:A5xA5}, we conclude that $G$ is isomorphic to $A_5\times A_5$, which contradicts our assumptions. This contradiction shows that $A\cap H=\{e\}$. Then $$60\cdot 60=|A|\cdot|H|=|A\cdot H|\le |G|=|A\cdot B|\le|A|\cdot|B|=60\cdot 60,$$ and hence $G=A\cdot H$. Applying Lemma~\ref{l:A5xA5}, we conclude that $G$ is isomorphic to $A_5\times A_5$, which contradicts our assumption. 
\end{proof}

Since the solvable group $H$ is a minimal normal subgroup of $G$, it contains no proper non-trivial characteristic subgroups and hence $H$ is elementary Abelian. Then $H\cap A$ is a normal elementary Abelian subgroup of the simple group $A\cong A_5$ and hence $H\cap A=\{e\}$. Then $|H{\cdot}A|=|H|\cdot|A|$ divides $|G|$, which divides $|A|\cdot|B|$ and hence the prime-power number $|H|$ divides $|B|=60$. This implies that $H\in\{C_2,C_2^2,C_3,C_5\}$. Consider the automorphism group $\Aut(H)$ of the elementary Abelian group $H$ and the homomorphism $\alpha:G\to \Aut(H)$ assigning to every $g\in G$ the automorphism $\alpha^g:H\to H$, $\alpha^g:x\mapsto x^g\defeq g^{-1}xg$, of the normal group $H$. Observe that the kernel of this automormophism coincides with the centralizer $\mathcal Z_G(H)$ of $H$ in $G$, which implies that $\mathcal Z_G(H)$ is a normal subgroup of $G$. Assuming that $A\not\subseteq \mathcal Z_G(H)$, we conclude that $A\cap\mathcal Z_G(H)=\{e\}$ and hence the restriction $\alpha{\restriction}_A$ is injective. Then $|\Aut(H)|$ is divisible by $|A|=60$. On the other hand, $H\in\{C_2,C_2^2,C_3,C_5\}$ implies $|\Aut(H)|\in\{|\IF_2^*|,|GL(2,2)|,|\IF_3^*|,|\IF_5^*|\}=\{1,6,2,4\}$, which contradicts the divisibility of $|\Aut(H)|$ by 60. This contradiction shows that $A\subseteq\mathcal Z_G(H)$. By analogy we can prove that $B\subseteq \mathcal Z_G(H)$. Then $G\subseteq \mathcal Z_G(H)$ and $H\subseteq\mathcal Z(G)$. Since $H$ is a minimal normal subgroup in $G$, $H\subseteq\mathcal Z(G)$ implies that $|H|$ is prime and belongs to the set $\{2,3,5\}$ of prime divisors of $|A|\cdot|B|=3600$.

Since $A\cap H=\{e\}=B\cap H$, the restrictions $\pi{\restriction}_A:A\to G/H$ and $\pi{\restriction}_B:B\to G/H$ are injective. Then the quotient group $G/H$ is the product $\pi[A]\cdot\pi[B]$ of two subgroups, isomorphic to the alternating group $A_5$. The minimality of $|G|$ ensures  that $G/H$ is isomorphic to $A_5\cong PSL(2,5)$ or $A_6\cong PSL(2,9)$ (the quotient group $G/H$ cannot be isomorphic to $A_5\times A_5$ because $|G/H|<|G|\le|A|\cdot|B|=3600=|A_5\times A_5|$).  
If $G/H\cong A_5$, then $A\cdot H=G$, and $H\subseteq\mathcal Z(G)$ implies that the subgroup $A$ is normal in $G=A{\cdot}H$. By Lemma~\ref{l:A5xA5}, $G$ is isomorphic to $A_5\times A_5$, which contradicts our assumption. 

This contradiction shows that $G/H\cong A_6\cong PSL(2,9)$. Then $$360\cdot p=|G/H|\cdot |H|=|G|=|A\cdot B|=|A|\cdot|B/(A\cap B)|=(|A|\cdot|B|)/|A\cap B|=3600/|A\cap B|$$implies
$p\cdot|A\cap B|=10$ and hence $|H|=p\in\{2,5\}$. Schur Theorem~\ref{t:Shur-multipliers} implies that $p$ divides $6$ and hence $p=2$. By Theorem~\ref{t:ShurPSL-SL}, the group $G$ is isomorphic to the group $SL(2,9)$ and hence has a unique element of order 2. On the other hand, $G$ contains the subgroup $A$ which is isomorphic to the group $A_5$, which contains 15 elements of order 2. This contradiction shows that the group $G$ should be simple. 

Observe that the $|G|$ divides $|A|\cdot|B|=3600$ and is divided by 60. The unique simple group of order $>60$ dividing 3600 and divided by 60 is the alternating group $A_6$. Therefore, $G\cong A_6$, which contradicts the choice of $G$. This contradiction completes the proof.
\end{proof}

\section{Inversive-plus invertible-puls affine spaces cannot have order $60$}

In this section we prove that inversive-plus invertible-puls affine spaces cannot have order $60$. This result will be essentially used in the proof of important (and difficult) Kegel--L\"uneburg Theorem~\ref{t:Kegel-Luneburg-main} on Desarguesianity of finite projective spaces satisfying the little Reidemeister condition.

\begin{theorems}\label{t:not60} If an affine space $X$ is inversive-plus and invertible-puls, then $|X|_2\ne 60$.
\end{theorems}

\begin{proof} To derive a contradiction, assume that $X$ is inversive-plus and invertible-puls, but $|X|_2=60$.  If $\|X\|>3$, then the affine space $X$ is Desarguesian and Thalesian, by Corollary~\ref{c:affine-Desarguesian} and Theorem~\ref{t:ADA=>AMA}. By Corollary~\ref{c:Thales=>commutative-plus}, $X$ is commutative-plus, and by Proposition~\ref{p:commutative-pls=>prime-power}, $|X|_2=60$ is a prime-power, which is not true. Therefore, $\|X\|=3$ and hence $X$ is a Playfair plane. By Theorems~\ref{t:+inversive} and \ref{t:invertible-add=>pls-elementary}, for every ternar $R$ of $X$, the plus loop $(R,+)$ is an elementary Moufang loop. Since $|R|=|X|_2=60$, the elementary Moufang loop $(R,+)$ is isomorphic to the alternating group $A_5$, according to Corollary~\ref{c:A5<=>}. Therefore, the loop $(R,+)$ is associative, which means that the Playfair plane $X$ is associative-plus.

Let us recall that for every parallel lines $L,\Lambda$ in $X$ and any direction $\boldsymbol \delta\in\partial X\setminus\{L_\parallel\}$, we denote by $\boldsymbol\delta_{\Lambda,L}:L\to\Lambda$ the line translation assigning to every point $x\in L$ the unique point $y\in\Lambda\cap\Aline x{\boldsymbol \delta}$. 
Since $X$ is associative-plus, we can apply Theorem~\ref{t:plus<=>group} and conclude that for every line $\Delta\subset X$ and distinct directions $\boldsymbol h,\boldsymbol v\in\partial X\setminus\{\Delta_\parallel\}$, the set
$$\Sym^\#_X[\Delta;\boldsymbol h,\boldsymbol v]\defeq\{\boldsymbol h_{\Delta,L}\boldsymbol v_{L,\Delta}:L\in\Delta_\parallel\}$$ is a subgroup of the group $\Sym^\#_X(\Delta)$ of line translations of the line $\Delta$. Moreover, the group 
$\Sym^\#_X[\Delta;\boldsymbol h,\boldsymbol v]$ is isomorphic to the plus loop $(\Delta,+)$, which is isomorphic to the alternating group $A_5$.

\begin{claim}\label{cl:KL-distinct-groups} There exist a line $\Delta$ in $X$ and three distinct directions $\boldsymbol a,\boldsymbol b,\boldsymbol c\in \partial X\setminus\{\Delta_\parallel\}$ 
such that the groups $\Sym^\#_X[\Delta;\boldsymbol a,\boldsymbol b]$ and $\Sym^\#_X[\Delta;\boldsymbol b,\boldsymbol c]$ are distinct. 
\end{claim}

\begin{proof} To derive a contradiction, assume that for any line $\Delta$ in $X$ and any distinct directions $\boldsymbol a,\boldsymbol b,\boldsymbol c\in \partial X\setminus\{\Delta_\parallel\}$ the groups $\Sym^\#_X[\Delta;\boldsymbol a,\boldsymbol b]$ and $\Sym^\#_X[\Delta;\boldsymbol b,\boldsymbol c]$ coincide.

We claim that the Playfair plane $X$ is Thalesian. Given three distinct parallel lines $A,B,C$ in $X$ and points $a,a'\in A$, $b,b'\in B$, $c,c'\in C$ with $\Aline ab\parallel \Aline {a'}{b'}$ and $\Aline bc\parallel \Aline {b'}{c'}$, we need to check that $\Aline ac\parallel \Aline{a'}{c'}$. If $\Aline ab=\Aline bc$, then $\Aline ac=\Aline ab\parallel \Aline{a'}{b'}=\Aline {a'}{c'}$ and we are done. So, assume that $\Aline ab\ne\Aline bc$. In this case, the directions $\boldsymbol a\defeq(\Aline bc)_\parallel$, $\boldsymbol b\defeq(\Aline ac)_\parallel$ and $\boldsymbol c\defeq(\Aline ab)_\parallel$  are pairwise distinct. 

By our assumption, $\Sym^\#_X[A;\boldsymbol b,\boldsymbol a]=\Sym^\#_X[A;\boldsymbol a,\boldsymbol c]$ and hence
$$\boldsymbol b_{A,C}\boldsymbol a_{C,B}\boldsymbol c_{B,A}=(\boldsymbol b_{A,C}\boldsymbol a_{C,A})(\boldsymbol a_{A,B}\boldsymbol c_{B,A})\in 
\Sym^\#_X[A;\boldsymbol b,\boldsymbol a]\cdot \Sym^\#_X[A;\boldsymbol a,\boldsymbol c]=\Sym^\#_X[A;\boldsymbol b,\boldsymbol a].$$ Since the group $\Sym^\#_X[A;\boldsymbol b,\boldsymbol a]$ acts on the line $A$ freely, the equality $a=\boldsymbol b_{A,C}\boldsymbol a_{C,B}\boldsymbol c_{B,A}(a)$ implies $$a'=\boldsymbol b_{A,C}\boldsymbol a_{C,B}\boldsymbol c_{B,A}(a')=\boldsymbol b_{A,C}\boldsymbol a_{C,B}(b')=\boldsymbol b_{A,C}(c'),$$ which means that $\Aline {a'}{c'}\in\boldsymbol b=(\Aline ac)_\parallel$.

Therefore, the Playfair plane $X$ is Thalesian. By Corollary~\ref{c:Thales=>commutative-plus}, $X$ is commutative-plus and by Proposition~\ref{p:commutative-pls=>prime-power}, $|X|_2$ is a prime-power and hence $|X_2|\ne 60$, which is a desired contradiction completing the proof.
\end{proof}

By Claim~\ref{cl:KL-distinct-groups},  there exist a line $\Delta$ in $X$ and three distinct directions $\boldsymbol a,\boldsymbol b,\boldsymbol c\in \partial X\setminus\{\Delta_\parallel\}$ 
such that the groups $A\defeq\Sym^\#_X[\Delta;\boldsymbol b,\boldsymbol c]$ and  $B\defeq\Sym^\#_X[\Delta;\boldsymbol a,\boldsymbol c]$ are distinct.
By Lemmas~\ref{l:IX-commute} and \ref{IXuv=IXud+IXdv}, $G\defeq A\cdot B=B\cdot A$ is a subgroup of $\Sym^\#_X(\Delta)$ containing the group $C\defeq\Sym^\#_X[\Delta;\boldsymbol a,\boldsymbol b]$. By the same reason, $A\cdot C=C\cdot A$ is a group containing the group $B$ and hence containing the group $G$, which implies $G=A\cdot C=C\cdot A$ and $C\ne A$. By the same reason, $G=B\cdot C=C\cdot B$ and $C\ne B$. 
Since the groups $A,B$ are isomorphic to $A_5$, we can apply Theorem~\ref{t:A5xA5=} and conclude that the group $G=A{\cdot}B$ is isomorphic either to the alternating group $A_6$ or to the group $A_5\times A_5$. 

First assume that $G\cong A_6$. Let $\mathcal A_5$ be the family of all subgroups of the group $A_6$, which are isomorphic to the group $A_5$. The family $\mathcal A_5$ has the following known properties:
\begin{enumerate}
\item $|\mathcal A_5|=12$;
\item $|\{g^{-1}Ag:g\in A_6\}|=6$ for every $A\in\A_5$;
\item if two distinct groups $A,B\in\mathcal A_5$ are conjugated in $A_6$, then $|A\cap B|=12$, if $A,B$ are not conjugated, then $|A\cap B|=10$.
\end{enumerate}

Identify the group $G$ with the alternating group $A_6$. The subgroups $A,B,C$ are isomorphic to the alternating group $A_5$ and hence belong to the family $\A_{5}$. It follows from $A{\cdot}B=G$ that 
$$360=|A_6|=|G|=|A|\cdot(|B|/|A\cap B|)=3600/|A\cap B|$$and hence $|A\cap B|=10$. The property (3) of the family $\A_{15}$ ensures that the groups $A,B$ are not conjugated in $A_6$. By analogy, we can prove that $A,C$ are not conjugated. Since the group $A_6$ contains exactly two conjugacy classes of subgroups isomorphic to $A_5$, $B$ is congugated to $C$, which implies $|B\cap C|=12$ and hence $G\ne B\cdot C$, which is a desired contradiction showing that the group $G$ cannot be isomorphic to $A_6$.
\smallskip

Therefore, $G$ is isomorphic to $A_5\times A_5$. In this case $|G|=3600$ and $A\cap B=A\cap C=B\cap C=\{e\}$. Let $\mathcal A_5$ be the family of all subgroups of $A_5\times A_5$, isomorphic to the group $A_5$. The family $\A_5$ contains two normal subgroups of $A_5\times A_5$, namely the ``horizontal'' subgroup $A_5\times\{e\}$ and the ``vertical'' subgroup $\{e\}\times A_5$. All other subgroups in $\A_5$ are ``diagonal''. More precisely, they are graphs of automorphisms of $A_5$. It is well-known that the automorphism group of the group $A_5$ is isomorphic to the group $S_5$ and any automorphism of $A_5$ has at least 3 fixed points. This implies that among the groups $A,B,C$ ony one can be diagonal. More precisely, two groups are normal and one is diagonal. We lose no generality assuming that the groups $A,B$ are normal and $C$ is diagonal. Then $G=A\times B$.  
Fix any point $o\in\Delta$ and consider its stabilizer $S_o\defeq\{g\in G:g(o)=o\}$. Since the group $G$ acts transitively on $\Delta$, $|S_o|=|G|/|\Delta|=60$. Since the groups $A,B$ act freely on $\Delta$, $A\cap S_o=B\cap S_o=\{e\}$. This implies that the subgroup $S_o$ projects bijectively on $A$ and $B$ and hence $S_o$ is a diagonal subgroup of $G=A_5\times B_5$, isomorphic to $A_5$. Since any two diagonal subgroups of $A_5\times A_5$ have at least $3$ common points, $|S_o\cap C|\ge 3$, which is impossible because the action of the group $C$ on $\Delta$ is free.   
\end{proof}

Theorem~\ref{t:not60} combined with Theorems~\ref{t:add-ass<=>},  \ref{t:invertible-add=>pls-elementary}, and \ref{t:elementary-group<=>} implies the following corollary, which should be attributed to \index[person]{Kegel}Kegel\footnote{{\bf Otto Helmut Kegel} (Born 1934, in Bethlehem, Pennsylvania) is a German mathematician specializing in group theory. Kegel completed secondary school  in Frankfurt am Main in 1954. He earned his PhD in 1961 at Johann Wolfgang Goethe University Frankfurt am Main under Reinhold Baer, with a dissertation titled ``{\em Commutativity of Subgroups and Composition Structure in Finite Groups}''. In 1966, he completed his habilitation in Frankfurt. In 1967, he was a lecturer at the University of Cologne. In 1968, he became a Reader at the University of London, where he was appointed Professor in 1970. Starting in 1975, he served as a full professor at Albert Ludwig University of Freiburg, where he remained until his retirement in 1999. From 1975 to 1989, Kegel was a co-editor of the journal Archiv f\"ur Mathematik. He also served as an editor for Journal of the London Mathematical Society and Proceedings of the London Mathematical Society (1970–1972), Communications in Algebra (1974–1985), Mathematische Zeitschrift (1978–1985). In 1962, he gave a talk at the International Congress of Mathematicians (ICM) in Stockholm, titled ``{\em Locally Finite Groups with a Partition}.''} and 
 and \index[person]{L\"uneburg}L\"uneburg\footnote{{\bf Heinrich L\"uneburg} (1935 -- 2009) was a German mathematician who worked in combinatorics, geometry, and algebra.
 L\"uneburg earned his PhD in 1962 at Johann Wolfgang Goethe University Frankfurt under Reinhold Baer (with Ruth Moufang also involved in his supervision). His dissertation was titled ``{\em Affine Hjelmslev Planes with a Transitive Translation Group}''. In 1963, he moved to Johannes Gutenberg University Mainz, where he completed his habilitation in 1964. From 1970 onward, he was a professor at the newly founded University of Kaiserslautern, where he helped establish the mathematics department. He declined an offer for a professorship in Bayreuth in 1975 and retired in 2003. L\"uneburg focused on finite geometries, which he studied using group-theoretic methods. Later in his career, his interests shifted toward algorithmic algebra and the history of mathematics. His work in the history of mathematics arose from intensive study of primary sources (such as Leonardo da Pisa, Cardano, and Tartaglia) and was driven by a deep curiosity about linguistic, chronological, and cultural-historical details. This curiosity was reflected in his writing style, which included numerous digressions, side notes, and personal observations. His book on Leonardo da Pisa's Liber Abbaci, though explicitly described by him as merely a ``mathematical reading pleasure'' rather than a scholarly historical study, is notable for its meticulous and uninhibited engagement with the text. It remains a unique work in Leonardo research and the most significant German-language contribution to the subject.}  
 \cite{KegelLuneburg1963}.

\begin{corollarys}\label{c:Kegel-Luneburg} If a finite affine space $X$ is associative-plus and invertible-puls, then for every ternar $R$ of $X$, the plus loop $(R,+)$ satisfies one of the following conditions:
\begin{enumerate}
\item $(R,+)$ is a nilpotent group of prime-power order;
\item $(R,+)$ is a Frobenius group whose kernel $K$ has prime-power order and complement has prime order that divides $|K|-1$.
\end{enumerate}
\end{corollarys}

For inversive-plus invertible-puls affine spaces of odd order, Theorems~\ref{t:+inversive}, \ref{t:inversive-puls}, \ref{t:invertible-add=>pls-elementary}, \ref{t:odd=>nucleo-solvable} and \ref{t:nucleo-solvable} imply the following corollary.

\begin{corollarys} If a finite affine space $X$ of odd order is inversive-plus and invertible-puls, then for every ternar $R$ of $X$, the plus loop $(R,+)$ satisfies one of the following conditions:
\begin{enumerate}
\item $(R,+)$ is a centrally-nilpotent Moufang loop of prime-power order;
\item $(R,+)$ is a Frobenius Moufang loop whose kernel $K$ has prime-power order and complement has prime order  that divides $|K|-1$.
\end{enumerate}
\end{corollarys}

\begin{corollarys} If a finite affine space $X$ of odd order is inversive-puls and invertible-plus, then for every ternar $R$ of $X$, the puls loop $(R,\!\puls\!)$ satisfies one of the following conditions:
\begin{enumerate}
\item $(R,\!\puls\!)$ is a centrally-nilpotent Moufang loop of prime-power order;
\item $(R,\!\puls\!)$ is a Frobenius Moufang loop whose kernel $K$ has prime-power order and complement has prime order  that divides $|K|-1$.
\end{enumerate}
\end{corollarys}

We do not know if the puls-counterpart of Theorem~\ref{t:not60} holds.

\begin{problem} Can an associative-puls invertible-plus affine space have order $60$?
\end{problem}

\chapter{Commutative-add affine liners}\label{ch:commutative-add}

An affine liner is defined to \defterm{commutative-add} if it is commutative-plus and commutative-puls. In this chapter we study commutative-add Playfair liners and show that a Playfair liner is commutative-add if and only if it is Boolean or Thalesian if and only if it is uno-Thalesian. Those two (difficult) characterizations will be proved in Theorems~\ref{t:uno-Thales<=>} and \ref{t:commutative-add<=>}. The proof motivates a study of Playfair liners satisfying some weaker versions of the Thales Axiom with one or two additional incidences between  vertices and sides of the triangles appearing in the Thales Axiom.  The liners satisfying those weaker versions of Thales Axioms are studied in Sections~\ref{s:uno-Thalesian}, \ref{s:bi-Boolean}, \ref{s:bi-Thalesian}, \ref{s:quadratic}.

\section{Uno-Thalesian liners}\label{s:uno-Thalesian}

\begin{definition}\label{d:uno-Thalesian} A liner $X$ is defined to be \index{liner!uno-Thalesian}\index{uno-Thalesian liner}\defterm{uno-Thalesian} if for any distinct parallel lines $A,B,C\subseteq X$ and distinct points $a,a'\in A$, $b,b'\in B$ and $c,c'\in C$,
$$\big(\,\Aline ab\cap\Aline{a'}{b'}=\varnothing=\Aline bc\cap\Aline{b'}{c'}\;\wedge\;b'\in \Aline ac\,\big)\;\Ra\;\big(\,\Aline ac\cap\Aline{a'}{c'}=\varnothing\,\big).$$

\begin{picture}(200,80)(-150,-10)
\put(-10,0){\line(1,0){100}}
\put(95,-3){$C$}
\put(-10,30){\line(1,0){100}}
\put(95,27){$B$}
\put(-10,60){\line(1,0){100}}
\put(95,57){$A$}

{\linethickness{1pt}
\put(10,0){\color{red}\line(0,1){60}}
\put(10,0){\color{blue}\line(1,1){30}}
\put(10,60){\color{cyan}\line(1,-1){30}}

\put(40,0){\color{red}\line(0,1){60}}
\put(40,0){\color{blue}\line(1,1){30}}
\put(40,60){\color{cyan}\line(1,-1){30}}
}

\put(10,0){\circle*{3}}
\put(8,-9){$c'$}
\put(10,60){\circle*{3}}
\put(8,63){$a'$}
\put(70,30){\circle*{3}}
\put(71,32){$b$}
\put(40,0){\circle*{3}}
\put(38,-9){$c$}
\put(40,60){\circle*{3}}
\put(38,63){$a$}
\put(40,30){\color{red}\circle*{3}}
\put(42,32){$b'$}
\end{picture}
\end{definition}

Definitions~\ref{d:para-Desargues} and \ref{d:uno-Thalesian} imply 

\begin{proposition} Every Thalesian liner is uno-Thalesian.
\end{proposition}

Let us recall that a \defterm{parallelogram} in a liner $X$ is any quadruple $abcd\in X^4$ such that $\Aline ab\parallel\Aline cd\ne\Aline ab$ and $\Aline bc\parallel\Aline ad\ne\Aline bc$. A parallelogram $abcd$ is \defterm{Boolean} if its diagonals are parallel. 

A liner is called \defterm{Boolean}  for every points $a,b,c,d\in X$, $\Aline ab\cap\Aline cd=\varnothing=\Aline bc\cap\Aline ad$ implies $\Aline ac\cap\Aline bd=\varnothing$. By Theorem~\ref{t:Boolean<=>}, a $3$-ranked liner $X$ is Boolean if and only if every parallelogram in $X$ is Boolean.

\begin{proposition}\label{p:Boole=>uno-Thales} Every Boolean Proclus liner is uno-Thalesian.
\end{proposition}

\begin{proof}  Assume that $X$ is a Boolean Proclus liner. To prove that $X$ is uno-Thalesian, take any three distinct parallel lines $A,B,C\subset X$ and points $a,a'\in A$, $b,b'\in B$ and $c,c'\in C$ such that $\Aline ab\cap\Aline{a'}{b'}=\varnothing=\Aline bc\cap\Aline {b'}{c'}$ and $b'\in \Aline ac$. We have to prove that $\Aline {a'}{c'}\cap\Aline {a}{c}=\varnothing$. 

\begin{picture}(200,82)(-180,-10)
\linethickness{0.8pt}
\put(0,0){\color{teal}\line(1,0){30}}
\put(0,30){\color{teal}\line(1,0){30}}
\put(0,60){\color{teal}\line(1,0){30}}
\put(30,0){\color{red}\line(0,1){60}}
\put(0,0){\color{blue}\line(1,1){30}}
\put(0,60){\color{violet}\line(1,-1){30}}
\multiput(0,0)(0,3){20}{\color{red}\line(0,1){2}}
\put(30,0){\color{blue}\line(-1,1){30}}
\put(30,60){\color{violet}\line(-1,-1){30}}

\put(30,0){\circle*{3}}
\put(33,-2){$c$}
\put(30,60){\circle*{3}}
\put(33,58){$a$}
\put(0,30){\circle*{3}}
\put(-8,28){$b$}
\put(0,0){\circle*{3}}
\put(-9,-2){$c'$}
\put(0,60){\circle*{3}}
\put(-9,58){$a'$}
\put(30,30){\color{red}\circle*{3}}
\put(33,28){$b'$}
\put(15,15){\color{violet}\circle*{3.5}}
\put(15,15){\color{white}\circle*{2.5}}
\put(15,45){\color{blue}\circle*{3.5}}
\put(15,45){\color{white}\circle*{2.5}}
\end{picture}

It follows that $baa'b'$ and $bcc'b'$ are parallelograms, which have disjoint diagonals $\Aline b{a'}\cap\Aline a{b'}=\varnothing=\Aline b{c'}\cap\Aline c{b'}$ because the liner $X$ is Boolean.  
Taking into account that $A,B,C$ are parallel lines and $b'\in\Aline ac\setminus(A\cup B)$, we conclude that $\Aline a{b'}=\Aline ac=\Aline {b'}c$ and $\Pi\defeq\overline{A\cup B\cup C}$ is a plane. Since the liner $X$ is Proclus, the equalities $
\Aline ac\cap\Aline {a'}b=\Aline a{b'}\cap\Aline {a'}b=\varnothing=\Aline c{b'}\cap\Aline {c'}b=\Aline ca\cap\Aline {c'}b$ imply $\Aline {a'}b=\Aline{c'}b=\Aline {a'}{c'}$ and $\Aline {a'}{c'}\cap\Aline ac=\Aline {a'}b\cap\Aline{a}{b'}=\varnothing$.
\end{proof}

\begin{proposition}\label{p:uno-Thalesian1} A Playfair liner $X$ is uno-Thalesian if and only if for any distinct parallel lines $A,B,C\subseteq X$ and  points $a,a'\in A$, $b,b'\in B$ and $c,c'\in C$ the conditions $\Aline ab\cap\Aline{a'}{b'}=\varnothing=\Aline ac\cap\Aline{a'}{c'}$ and $b'\in \Aline ac$ imply $\Aline bc\cap\Aline{b'}{c'}=\varnothing$.

\begin{picture}(200,85)(-150,-10)
\put(-10,0){\line(1,0){100}}
\put(95,-3){$C$}
\put(-10,30){\line(1,0){100}}
\put(95,27){$B$}
\put(-10,60){\line(1,0){100}}
\put(95,57){$A$}

{\linethickness{1pt}
\put(10,0){\color{blue}\line(0,1){60}}
\put(10,0){\color{red}\line(1,1){30}}
\put(10,60){\color{cyan}\line(1,-1){30}}

\put(40,0){\color{blue}\line(0,1){60}}
\put(40,0){\color{red}\line(1,1){30}}
\put(40,60){\color{cyan}\line(1,-1){30}}
}

\put(10,0){\circle*{3}}
\put(8,-9){$c'$}
\put(10,60){\circle*{3}}
\put(8,63){$a'$}
\put(70,30){\circle*{3}}
\put(71,32){$b$}
\put(40,0){\circle*{3}}
\put(38,-9){$c$}
\put(40,60){\circle*{3}}
\put(38,63){$a$}
\put(40,30){\color{red}\circle*{3}}
\put(42,32){$b'$}
\end{picture}
\end{proposition}

\begin{proof} Assume that $X$ is uno-Thalesian. Given any distinct parallel lines $A,B,C\subset X$ and points $a,a'\in A$, $b,b'\in B$, $c,c'\in C$ with $\Aline ab\cap\Aline{a'}{b'}=\varnothing=\Aline ac\cap\Aline {a'}{c'}$ and $b'\in\Aline ac$, we have to prove that $\Aline bc\cap\Aline {b'}{c'}=\varnothing$. Since the lines $A,B,C$ are parallel and $b'\in B\cap\Aline ac$, the flat hull $\Pi\defeq\overline{A\cup B\cup C}$ is a plane in $X$. Assuming that $b'\in \Aline bc$, we conclude that $\Aline bc=\Aline ac=\Aline ab$ and $b'\in \Aline {a'}{b'}\cap \Aline ab=\varnothing$, which is a contradiction showing that $b'\notin\Aline bc$. Since $X$ is Playfair, the plane $\Pi$ contains a unique line $L$ such that $b'\in L\subset\Pi\setminus \Aline bc$. Since $L\parallel\Aline bc$ and $\Aline bc\cap C=\{c\}$, we can apply Proposition~\ref{p:Proclus-Postulate} and find a unique point $c''\in L\cap C$. Since $X$ is uno-Thalesian, $\Aline ab\cap\Aline{a'}{b'}=\varnothing=\Aline bc\cap\Aline {b'}{c''}$ and $b'\in\Aline ac$ imply $\Aline ac\cap\Aline {a'}{c''}=\varnothing$. By the Proclus Axiom, $\Aline {a'}{c'}\cap\Aline ac=\varnothing=\Aline{a'}{c''}\cap\Aline ac$ imply $\Aline {a'}{c'}=\Aline{a'}{c''}$ and hence $c''\in C\cap\Aline{a'}{c''}=C\cap\Aline{a'}{c'}=\{c'\}$ and finally, $\Aline bc\cap\Aline {b'}{c'}=\Aline bc\cap\Aline{b'}{c''}=\Aline bc\cap L=\varnothing$.
\smallskip

Now assume that a Playfair liner $X$ satisfies the ``if'' condition of the proposition. To prove that $X$ is uno-Thalesian, fix any distinct parallel lines $A,B,C\subset X$ and any points $a,a'\in A$, $b,b'\in B$, $c,c'\in C$ such that $\Aline ab\cap\Aline {a'}{b'}=\varnothing=\Aline bc\cap\Aline{b'}{c'}$ and $b'\in \Aline ac$. We have to prove that $\Aline ac\cap\Aline{a'}{c'}=\varnothing$. 
 Since the lines $A,B,C$ are parallel and $b'\in B\cap\Aline ac$, the flat hull $\Pi\defeq\overline{A\cup B\cup C}$ is a plane in $X$. Assuming that $a'\in \Aline ac$, we conclude that $a'\in \Aline ac\cap A=\{a\}$, which contradicts $\Aline ab\cap\Aline{a'}{b'}=\varnothing$. This contradiction shows that $a'\notin \Aline ac$. Since $X$ is Playfair, the plane $\Pi$ contains a unique line $L$ such that $a'\in L\subset\Pi\setminus \Aline ac$. Since $L\parallel\Aline ac$ and $\Aline ac\cap C=\{c\}$, we can apply Proposition~\ref{p:Proclus-Postulate} and find a unique point $c''\in L\cap C$. Since $\Aline ab\cap\Aline{a'}{b'}=\varnothing=\Aline ac\cap L=\Aline ac\cap\Aline{a'}{c''}$ and $b'\in\Aline ac$, the ``if'' assumption ensures that $\Aline bc\cap\Aline{b'}{c''}=\varnothing$. By the Proclus Axiom, $\Aline bc\cap\Aline {b'}{c''}=\varnothing=\Aline bc\cap\Aline {b'}{c'}$ imply $\Aline {b'}{c'}=\Aline{b'}{c''}$ and hence $c''\in C\cap\Aline{b'}{c''}=C\cap\Aline{b'}{c'}=\{c'\}$ and finally, $\Aline {a'}{c'}\cap\Aline ac=\Aline{a'}{c''}\cap\Aline ac=L\cap \Aline ac=\varnothing$.
\end{proof}

\begin{proposition}\label{p:uno-Thalesian2} Let $abc$ and $a'b'c'$ be two triangles in a uno-Thalesian Playfair liner $X$ such that $\Aline ab\parallel \Aline{a'}{b'}$, $\Aline bc\parallel\Aline{b'}{c'}$ and $\Aline ac\parallel\Aline{a'}{c'}$. If $\{a',b',c'\}\cap(\Aline ab\cup\Aline bc\cup\Aline ac)\ne\varnothing$ and $\Aline a{a'}\cap\Aline b{b'}=\varnothing$, then $\Aline b{b'}\parallel\Aline c{c'}$.
\end{proposition}

\begin{picture}(400,80)(-60,-15)
\linethickness{1pt}
\put(0,0){\color{red}\line(1,0){30}}
\put(0,0){\color{blue}\line(1,1){30}}
\put(0,0){\color{cyan}\line(0,1){60}}
\put(30,0){\color{blue}\line(1,1){30}}
\put(30,0){\color{cyan}\line(0,1){60}}
\put(0,60){\color{red}\line(1,0){30}}
\put(0,60){\color{violet}\line(1,-1){30}}
\put(30,30){\color{red}\line(1,0){30}}
\put(30,60){\color{violet}\line(1,-1){30}}

\put(0,0){\circle*{3}}
\put(-2,-9){$c'$}
\put(30,0){\circle*{3}}
\put(28,-9){$c$}
\put(30,30){\circle*{3}}
\put(31,32){$a'$}
\put(0,60){\circle*{3}}
\put(-2,63){$b'$}
\put(30,60){\circle*{3}}
\put(28,63){$b$}
\put(60,30){\circle*{3}}
\put(63,28){$a$}

\put(100,0){\color{red}\line(1,0){30}}
\put(100,0){\color{blue}\line(1,1){30}}
\put(100,0){\color{cyan}\line(0,1){60}}
\put(130,0){\color{blue}\line(1,1){30}}
\put(130,0){\color{cyan}\line(0,1){60}}
\put(100,60){\color{red}\line(1,0){30}}
\put(100,60){\color{violet}\line(1,-1){30}}
\put(130,30){\color{red}\line(1,0){30}}
\put(130,60){\color{violet}\line(1,-1){30}}

\put(100,0){\circle*{3}}
\put(98,-9){$c'$}
\put(130,0){\circle*{3}}
\put(128,-9){$c$}
\put(130,30){\circle*{3}}
\put(131,32){$b'$}
\put(100,60){\circle*{3}}
\put(98,63){$a'$}
\put(130,60){\circle*{3}}
\put(128,63){$a$}
\put(160,30){\circle*{3}}
\put(163,28){$b$}

\put(200,0){\color{red}\line(1,0){30}}
\put(200,0){\color{blue}\line(1,1){30}}
\put(200,0){\color{cyan}\line(0,1){60}}
\put(230,0){\color{blue}\line(1,1){30}}
\put(230,0){\color{cyan}\line(0,1){60}}
\put(200,60){\color{red}\line(1,0){30}}
\put(200,60){\color{violet}\line(1,-1){30}}
\put(230,30){\color{red}\line(1,0){30}}
\put(230,60){\color{violet}\line(1,-1){30}}

\put(200,0){\circle*{3}}
\put(197,-9){$b'$}
\put(230,0){\circle*{3}}
\put(227,-9){$b$}
\put(230,30){\circle*{3}}
\put(231,32){$c'$}
\put(200,60){\circle*{3}}
\put(198,63){$a'$}
\put(230,60){\circle*{3}}
\put(228,63){$a$}
\put(260,30){\circle*{3}}
\put(263,28){$c$}
\end{picture}

\begin{proof} Consider the planes $\overline{\{a,b,c\}}$ and $\overline{\{a',b',c'\}}$. By Theorem~\ref{t:subparallel-via-base}, $\Aline ab\parallel \Aline{a'}{b'}$ and $\Aline bc\parallel \Aline{b'}{c'}$ imply $\overline{\{a,b,c\}}\parallel \overline{\{a',b',c'\}}$. Since $\{a',b',c'\}\cap(\Aline ab\cup \Aline bc\cup\Aline ac)\ne\varnothing$, the parallel planes $\overline{\{a,b,c\}}$ and $\overline{\{a',b',c'\}}$ coincide.


If $a=a'$, then $\Aline ab\parallel \Aline {a'}{b'}$ and $\Aline ac\parallel \Aline{a'}{c'}$ imply $\Aline ab=\Aline{a'}{b'}$ and $\Aline ac=\Aline{a'}{c'}$. Since $\Aline ab=\Aline{a'}{b'}$, the condition $\Aline a{a'}\cap\Aline b{b'}=\varnothing$ implies $b=b'$. Then $\Aline bc\parallel \Aline {b'}{c'}$ implies $\Aline bc=\Aline{b'}{c'}$ and hence $\{c\}=\Aline ac\cap\Aline bc=\Aline{a'}{c'}\cap\Aline{b'}{c'}=\{c'\}$, which implies $\Aline c{c'}=\{c\}\parallel \{b\}=\Aline b{b'}$. By analogy we can prove that $b=b'$ implies $\Aline c{c'}=\{c\}\parallel\{b\}=\Aline b{b'}$.

So, assume that $a\ne a'$ and $b\ne b'$. In this case, $\Aline a{a'}$ and $\Aline{b}{b'}$ are two disjoint lines in the plane $\Pi\defeq\overline{\{a,b,c\}}=\overline{\{a',b',c'\}}$. It follows from $a\ne a'$ that the parallel lines $\Aline ab$, $\Aline {a'}{b'}$ are disjoint. If $c=c'$, then $\Aline ac\parallel \Aline {a'}{c'}$ implies $\Aline ac=\Aline{a'}{c'}=\Aline a{a'}$ and then $\Aline bc\parallel \Aline {b'}{c'}$ implies $\Aline bc=\Aline {b'}{c'}$ and $\{b\}=\Aline bc\cap\Aline b{b'}=\Aline{b'}{c'}\cap\Aline b{b'}=\{b'\}$, which contradicts $b\ne b'$. This contradiction shows that $c\ne c'$. 

 If $\{c,c'\}\cap\Aline a{a'}\ne\varnothing$ or $\{a,a'\}\cap\Aline c{c'}\ne\varnothing$, then $\Aline ac\parallel \Aline {a'}{c'}$ implies $\Aline ac=\Aline{a'}{c'}$ and hence $\Aline c{c'}=\Aline a{a'}\parallel \Aline b{b'}$. If $\{c,c'\}\cap\Aline b{b'}\ne\varnothing$ or $\{b,b'\}\cap\Aline c{c'}\ne\varnothing$, then $\Aline bc\parallel \Aline {b'}{c'}$ imply $\Aline c{c'}=\Aline bc=\Aline{b'}{c'}=\Aline b{b'}$.

So, assume that $\{c,c'\}\cap(\Aline a{a'}\cup\Aline b{b'})=\varnothing=\{a,a',b,b'\}\cap\Aline c{c'}$. Then  $\{a',b',c'\}\cap(\Aline ab\cup\Aline bc\cup\Aline ac)\ne\varnothing$ implies $a'\in \Aline bc$ or $b'\in \Aline ac$ or $c'\in \Aline ab$. 

If $a'\in \Aline bc$ or $b'\in \Aline ac$, then consider the unique line $C\subseteq\Pi$ such that $c\in C$ and $C\parallel\Aline a{a'}\parallel \Aline b{b'}$. The assumption $c\notin \Aline a{a'}\cup\Aline b{b'}$ ensures that $\Aline a{a'},\Aline b{b'},C$ are three distinct parallel lines in the plane $\Pi$. Taking into account that $a\ne a'$, $b\ne b'$, $\Aline ab\parallel \Aline{a'}{b'}$, $\Aline bc\parallel \Aline{b'}{c'}$, and $\Aline ac\parallel \Aline{a'}{c'}$, we conclude that $\varnothing=\Aline ab\cap\Aline {a'}{b'}=\Aline bc\cap\Aline{b'}{c'}=\Aline ac\cap\Aline{a'}{c'}$.  By Proposition~\ref{p:Proclus-Postulate}, there exists a unique point $c''\in C\cap\Aline{b'}{c'}$. If $b'\in\Aline ac$, then $\Aline ab\cap\Aline{a'}{b'}=\varnothing=\Aline bc\cap\Aline{b'}{c'}=\Aline bc\cap\Aline {b'}{c''}$ implies $\Aline ac\cap\Aline{a'}{c''}$ by the uno-Thalesian property of $X$. If $a'\in\Aline bc$, then $\Aline ab\cap\Aline{a'}{b'}=\varnothing=\Aline bc\cap\Aline {b'}{c''}$ implies $\Aline ac\cap\Aline{a'}{c''}=\varnothing$, by Proposition~\ref{p:uno-Thalesian1}. In both cases, we conclude that the lines $\Aline ac$ and $\Aline{a'}{c''}$ in the plane $\Pi$ are disjoint  and hence parallel. Then $\Aline {a'}{c'}\parallel \Aline ac\parallel\Aline{a'}{c''}$ implies $\Aline{a'}{c'}=\Aline{a'}{c''}$ and $c''\in \Aline {a'}{c''}\cap\Aline{b'}{c''}=\Aline {a'}{c'}\cap\Aline{b'}{c'}=\{c'\}$ and finally, $\Aline c{c'}=\Aline c{c''}=C\parallel \Aline b{b'}$.

If $c'\in\Aline ab$, then consider the unique line $C$ in the plane $\Pi$ such that $c'\in C$ and $C\parallel \Aline a{a'}\parallel \Aline b{b'}$. The assumption $c'\notin \Aline a{a'}\cup\Aline b{b'}$ ensures that $\Aline a{a'},\Aline b{b'},C$ are three distinct parallel lines in the plane $\Pi$. Since $\Aline {a'}{c'}\cap C=\{c'\}$ and $\Aline ac\parallel \Aline{a'}{c'}$, there exists a unique point $\gamma\in \Aline ac\cap C$. By Proposition~\ref{p:uno-Thalesian1}, $c'\in \Aline ab$ and $\Aline a\gamma\cap\Aline {a'}{c'}=\Aline ac\cap\Aline{a'}{c'}=\varnothing=\Aline ab\cap\Aline{a'}{b'}$ imply $\Aline \gamma{b'}\cap \Aline cb=\varnothing$. It follows that $\Aline \gamma{b}\parallel \Aline {c'}{b'}\parallel \Aline{c}{b}$ and hence $\Aline \gamma{b}=\Aline{c}{b}$ and $\gamma\in \Aline \gamma{b}\cap\Aline a\gamma=\Aline{c}{b}\cap\Aline ac=\{c\}$. 
Then $\Aline c{c'}=\Aline \gamma {c'}=C\parallel \Aline a{a'}\parallel \Aline b{b'}$.
\end{proof}

\section{Bi-Boolean, di-Boolean, and by-Boolean liners}\label{s:bi-Boolean}

In this section we consider three weaker versions of the Thales axiom in which every triangle has a vertex of a side of the other triangle. Since liners satisfying such weakenings of the Thales axiom are closely related to Boolean liners, we call them bi-Boolean, di-Boolean and by-Boolean.

\begin{definition} A liner $X$ is defined to be 
\begin{itemize}
\item \defterm{bi-Boolean} if for any disjoint lines $L,L'\subseteq X$\\

and distinct points $a,b,c\in L$, $a',b',c'\in L$, the equalities\\ $\Aline a{b'}\cap\Aline {a'}b=\Aline b{c'}\cap\Aline{b'}c=\Aline a{a'}\cap\Aline c{c'}=\varnothing$ imply $\Aline a{a'}\cap\Aline b{b'}=\varnothing$;
\smallskip
\item \defterm{di-Boolean} if for any disjoint lines $L,L'\subseteq X$\\
%
and distinct points $a,b,c\in L$, $a',b',c'\in L$, the equalities\\ $\Aline a{b'}\cap\Aline {a'}b=\Aline b{c'}\cap\Aline{b'}c=\Aline a{a'}\cap\Aline b{b'}=\varnothing$ imply $\Aline b{b'}\cap\Aline c{c'}=\varnothing$;
\item \defterm{by-Boolean} if for any lines $L,L'\subseteq X$\\
%
and distinct points $a,b,c\in L$, $a',b',c'\in L'$, the equalities\\ $\Aline a{b'}\cap\Aline {a'}b=\Aline b{c'}\cap\Aline{b'}c=\Aline a{a'}\cap\Aline b{b'}=\Aline b{b'}\cap\Aline c{c'}=\varnothing$ imply $L\cap L'=\varnothing$.
\end{itemize}
\end{definition}

\begin{proposition}\label{p:Boole=>bi-di-by-Boole} Every Boolean liner is bi-Boolean, di-Boolean and by-Boolean.
\end{proposition}

\begin{proof} Let $X$ be a Boolean liner.
\smallskip

(bi) To prove that $X$ is bi-Boolean, take any disjoint lines $L,L'\subset X$ and distinct points $a,b,c\in L$, $a',b',c'\in L'$ such that $\Aline a{b'}\cap\Aline {a'}b=\Aline b{c'}\cap\Aline {b'}c=\Aline a{a'}\cap\Aline c{c'}=\varnothing$. We should prove that $\Aline a{a'}\cap\Aline b{b'}=\varnothing$. To derive a contradiction, assume that $\Aline a{a'}\cap\Aline b{b'}\ne \varnothing$. Then the lines $\Aline a{a'}$ and $\Aline b{b'}$ are coplanar and $\Pi\defeq\overline{\{a,a',b,b'\}}$ is a plane containing the disjoint lines $L=\overline{\{a,b,c\}}$ and $L'=\overline{\{a',b',c'\}}$. Then $ab'ba'$ is a parallelogram in the plane $\Pi$. Since $X$ is Boolean, this parallelograms has disjoint diagonals: $\Aline a{a'}\cap\Aline b{b'}=\varnothing$, which contradicts our assumption.   
\smallskip

(di) To prove that $X$ is di-Boolean, take any disjoint lines $L,L'\subset X$ and distinct points $a,b,c\in L$, $a',b',c'\in L'$ such that $\Aline a{b'}\cap\Aline {a'}b=\Aline b{c'}\cap\Aline {b'}c=\Aline a{a'}\cap\Aline b{b'}=\varnothing$. We should prove that $\Aline b{b'}\cap\Aline c{c'}=\varnothing$. To derive a contradiction, assume that $\Aline b{b'}\cap\Aline c{c'}\ne \varnothing$. Then the lines $\Aline b{b'}$ and $\Aline c{c'}$ are coplanar and $\Pi\defeq\overline{\{b,b',c,c'\}}$ is a plane containing the disjoint lines $L=\overline{\{a,b,c\}}$ and $L'=\overline{\{a',b',c'\}}$. Then $bc'cb'$ is a parallelogram in the plane $\Pi$. Since $X$ is Boolean, this parallelograms has disjoint diagonals: $\Aline b{c'}\cap\Aline c{c'}=\varnothing$, which contradicts our assumption.
\smallskip

(by) To prove that $X$ is by-Boolean, take any lines $L,L'\subset X$
and distinct points $a,b,c\in L$ and $a',b',c'\in L'$ such that $\Aline a{b'}\cap\Aline {a'}b=\Aline b{c'}\cap\Aline{b'}c=\Aline a{a'}\cap\Aline b{b'}=\Aline b{b'}\cap\Aline c{c'}=\varnothing$. We have to prove that $L\cap L'=\varnothing$. To derive a contradiction, assume that $L\cap L'\ne \varnothing$. Assuming that $L=L'$, we conlude that $a,b,c,a',b',c'\in L=L'$. In this case $\Aline a{b'}\cap\Aline b{a'}=\Aline a{a'}\cap\Aline b{b'}=\varnothing$ implies $a=b'\ne b=a'=a$, which is a contradiction showing that $L\ne L'$. Then $\Pi\defeq\overline{L\cup L'}$ is a plane and $aa'bb'$ is a  paralelogram in the plane $\Pi$. Since the liner $X$ is Boolean, this paralellogram has disjoint diagonals and hence $L\cap L'=\Aline ab\cap \Aline b{b'}=\varnothing$, which contradicts our assumption.
\end{proof}

\begin{proposition}\label{p:uno-Thalesian=>bi-di-by-Boolean} Every uno-Thalesian Playfair liner is bi-Boolean, di-Boolean, and by-Boolean.
\end{proposition}

\begin{proof} Let $X$ be a uno-Thalesian Playfair liner.
\smallskip

(bi) To prove that $X$ is bi-Boolean, take any disjoint lines $L,L'\subset X$ and distinct points $a,b,c\in L$, $a',b',c'\in L'$ such that $\Aline a{b'}\cap\Aline {a'}b=\Aline b{c'}\cap\Aline {b'}c=\Aline a{a'}\cap\Aline c{c'}=\varnothing$. We should prove that $\Aline a{a'}\cap\Aline b{b'}=\varnothing$. To derive a contradiction, assume that $\Aline a{a'}\cap\Aline b{b'}\ne \varnothing$. Then the lines $\Aline a{a'}$ and $\Aline b{b'}$ are coplanar and $\Pi\defeq\overline{\{a,a',b,b'\}}$ is a plane containing the disjoint lines $L=\overline{\{a,b,c\}}$ and $L'=\overline{\{a',b',c'\}}$. It follows from Corollary~\ref{c:parallel-lines<=>} that $ab'c$ and $a'bc'$ are two triangles in the plane $\Pi$ such that $\Aline a{b'}\parallel \Aline {a'}b$, $\Aline {b'}c\parallel \Aline b{c'}$, $\Aline ac\parallel \Aline{a'}{c'}$, and $\Aline a{a'}\cap\Aline c{c'}=\varnothing$. Applying Proposition~\ref{p:uno-Thalesian2}, we conclude that $\Aline b{b'}\parallel \Aline a{a'}$, which contradicts the assumptions $\Aline a{a'}\cap\Aline b{b'}\ne\varnothing$ and $a\ne b$. This contradiction shows that $\Aline a{a'}\cap\Aline b{b'}=\varnothing$.
\smallskip

(di) To prove that $X$ is di-Boolean, take any disjoint lines $L,L'\subset X$ and distinct points $a,b,c\in L$, $a',b',c'\in L'$ such that $\Aline a{b'}\cap\Aline {a'}b=\Aline b{c'}\cap\Aline {b'}c=\Aline a{a'}\cap\Aline b{b'}=\varnothing$. We should prove that $\Aline b{b'}\cap\Aline c{c'}=\varnothing$. To derive a contradiction, assume that $\Aline b{b'}\cap\Aline c{c'}\ne \varnothing$. Then the lines $\Aline b{b'}$ and $\Aline c{c'}$ are coplanar and $\Pi\defeq\overline{\{b,b',c,c'\}}$ is a plane containing the disjoint lines $L=\overline{\{a,b,c\}}$ and $L'=\overline{\{a',b',c'\}}$. It follows from Corollary~\ref{c:parallel-lines<=>} that $ab'c$ and $a'bc'$ are two triangles in the plane $\Pi$ such that $\Aline a{b'}\parallel \Aline {a'}b$, $\Aline {b'}c\parallel \Aline b{c'}$, $\Aline ac\parallel \Aline{a'}{c'}$, and $\Aline a{a'}\cap\Aline b{b'}=\varnothing$. Applying Proposition~\ref{p:uno-Thalesian2}, we conclude that $\Aline b{b'}\parallel \Aline c{c'}$, which contradicts the assumptions $\Aline b{b'}\cap\Aline c{c'}\ne\varnothing$ and $b\ne c$. This contradiction shows that $\Aline b{b'}\cap\Aline c{c'}=\varnothing$.
\smallskip

(by) To prove that $X$ is by-Boolean, take any lines $L,L'\subset X$
and distinct points $a,b,c\in L$ and $a',b',c'\in L'$ such that $\Aline a{b'}\cap\Aline {a'}b=\Aline b{c'}\cap\Aline{b'}c=\Aline a{a'}\cap\Aline b{b'}=\Aline b{b'}\cap\Aline c{c'}=\varnothing$. We have to prove that $L\cap L'=\varnothing$. To derive a contradiction, assume that $L\cap L'\ne \varnothing$. Assuming that $L=L'$, we conlude that $a,b,c,a',b',c'\in L=L'$. In this case $\Aline b{c'}\cap\Aline{b'}c=\Aline b{b'}\cap\Aline c{c'}=\varnothing$ implies $b=c'\ne b'=c=c'$, which is a contradiction showing that $L\ne L'$. Then $\Pi\defeq\overline{L\cup L'}$ is a plane containing three  lines $A\defeq\Aline a{a'}$, $B\defeq\Aline b{b'}$ and $C\defeq\Aline c{c'}$ such that $A\cap B=\varnothing=B\cap C$. By Corollary~\ref{c:parallel-lines<=>}, $A\parallel B\parallel C$.
 Assuming that $A\cap C\ne\varnothing$, we conclude that $\Aline a{a'}=A=C=\Aline c{c'}$ and hence $L=\Aline ac=\Aline{a'}{c'}=L'$. Then $\Aline a{b'}\cap\Aline {a'}b=\Aline b{c'}\cap\Aline{b'}c=\Aline a{a'}\cap\Aline b{b'}=\Aline b{b'}\cap\Aline c{c'}=\varnothing$ implies $a=b'\ne a'=b=b'$, which is a contradiction, showing that $A\cap C=\varnothing$. Therefore, $A,B,C$ are three distinct parallel lines in $X$ and $a,a'\in A$, $b,b'\in B$, $c,c'\in C$ are distinct points such that $\Aline a{b'}\cap\Aline {a'}b=\varnothing=\Aline b{c'}\cap\Aline {b'}c$ and $b'\in \Aline {a'}{c'}$. Since $X$ is uno-Thalesian, $L\cap L'=\Aline ac\cap\Aline {a'}{c'}=\varnothing$, which contradicts the assumption $L\cap L'\ne\varnothing$. This contradiction shows that $L\cap L'=\varnothing$.
\end{proof}

Di-Boolean Playfair liners admit the following characterization. 

\begin{proposition}\label{p:di-Boolean<=>} A Playfair liner $X$ is di-Boolean if and only if for any disjoint lines $L,L'\subseteq X$ and distinct points $a,b,c\in L$, $a',b',c'\in L'$, the equalities $\Aline a{b'}\cap\Aline {a'}b=\Aline a{a'}\cap\Aline b{b'}=\Aline b{b'}\cap\Aline c{c'}=\varnothing$ imply $\Aline b{c'}\cap\Aline {b'}c=\varnothing$.

\begin{picture}(100,70)(-160,-15)
\linethickness{1pt}
\put(0,0){\color{cyan}\line(0,1){40}}
\put(0,0){\color{teal}\line(1,0){60}}
\put(0,40){\color{teal}\line(1,0){60}}
\put(20,0){\color{cyan}\line(0,1){40}}
\put(0,0){\color{blue}\line(1,2){20}}
\put(0,40){\color{blue}\line(1,-2){20}}
\put(20,0){\color{red}\line(1,1){40}}
\put(20,40){\color{red}\line(1,-1){40}}
\put(60,0){\color{cyan}\line(0,1){40}}
\put(0,0){\circle*{3}}
\put(-2,-9){$a$}
\put(20,0){\circle*{3}}
\put(17,-9){$b$}
\put(60,0){\circle*{3}}
\put(58,-9){$c$}
\put(0,40){\circle*{3}}
\put(-2,43){$a'$}
\put(20,40){\circle*{3}}
\put(18,43){$b'$}
\put(60,40){\circle*{3}}
\put(58,43){$c'$}
\put(10,20){\color{blue}\circle*{4}}
\put(10,20){\color{white}\circle*{3}}
\put(40,20){\color{red}\circle*{4}}
\put(40,20){\color{white}\circle*{3}}
\end{picture}
\end{proposition}

\begin{proof} Assume that a Playfair liner $X$ is di-Boolean. Let $L,L'\subset X$ be two disjoint lines and $a,b,c\in L$ and $a',b',c'\in L'$ be distinct points such that $\Aline a{b'}\cap\Aline{a'}b=\Aline a{a'}\cap\Aline b{b'}=\Aline b{b'}\cap\Aline c{c'}=\varnothing$. We have to prove that $\Aline b{c'}\cap\Aline{b'}c=\varnothing$. 
To derive a contradiction, assume that $\Aline b{c'}\cap\Aline {b'}c\ne\varnothing$. Then the lines $\Aline b{c'},\Aline {b'}c$ are coplanar and $\Pi\defeq\overline{\{b,b',c,c'\}}$ is a plane containing the lines $L=\Aline bc$ and $L'=\Aline {b'}{c'}$. Since the liner $X$ is Playfair, there exists a unique point $c''\in L'$ such that $\Aline b{c''}\parallel \Aline {b'}c$ and hence $\Aline b{c''}\cap\Aline{b'}c=\varnothing$. Since $X$ is di-Boolean, $\Aline {a'}b\cap\Aline a{b'}=\varnothing=\Aline {b'}c\cap\Aline b{c''}$ implies 
$\Aline b{b'}\cap \Aline c{c''}=\varnothing$ and hence $\Aline c{c'}\parallel \Aline b{b'}\parallel \Aline c{c''}$ and $\Aline c{c'}=\Aline c{c''}$. Then $c''\in \Aline c{c''}\cap L'=\Aline c{c'}\cap L'=\{c'\}$ and $\Aline b{c'}\cap\Aline {b'}c=\Aline b{c''}\cap\Aline {b'}c=\varnothing$, which contradicts the assumption. This contradiction shows that $\Aline b{c'}\cap\Aline{b'}c=\varnothing$.
\smallskip

Assume that a Playfair liner $X$ satisfies the ``if'' condition of the proposition. To prove that $X$ is di-Thalesian, take any disjoint lines $L,L'\subseteq X$ and distinct points $a,b,c\in L$ and $a',b',c'\in L$ such that $\Aline a{b'}\cap\Aline {a'}b=\varnothing=\Aline a{a'}\cap\Aline b{b'}$. 
We have to prove that $\Aline b{b'}\cap\Aline{c}{c'}=\varnothing$.
To derive a contradiction, assume that $\Aline b{b'}\cap\Aline{c}{c'}\ne\varnothing$.  Then the lines $\Aline b{b'},\Aline c{c'}$ are coplanar and $\Pi\defeq\overline{\{b,b',c,c'\}}$ is a plane containing the lines $L=\Aline bc$ and $L'=\Aline {b'}{c'}$. Since the liner $X$ is Playfair, there exists a unique point $c''\in L'$ such that $\Aline b{b'}\cap\Aline c{c''}=\varnothing$. The ``if'' condition ensures that $\Aline b{c''}\cap\Aline {b'}c=\varnothing$ and hence $\Aline b{c''}\parallel \Aline {b'}c\parallel \Aline b{c'}$ and $\Aline b{c''}=\Aline b{c'}$. Then $c''\in \Aline b{c''}\cap L'=\Aline b{c'}\cap L'=\{c'\}$ and $\Aline b{b'}\cap\Aline c{c'}=\Aline b{b'}\cap\Aline c{c''}=\varnothing$, which contradicts the assumption. This contradiction shows that $\Aline b{b'}\cap\Aline c{c'}=\varnothing$.
\end{proof}

\begin{proposition} A commutative-puls Playfair liner $X$ is bi-Boolean if and only if it is di-Boolean.
\end{proposition}

\begin{proof} Asume that $X$ is bi-Boolean. To prove that $X$ is di-Boolean, take any disjoint lines $L,L'\subseteq X$ and distinct points $a,b,c\in L$ and $a',b',c'\in L$ such that $\Aline a{b'}\cap\Aline {a'}b=\Aline b{c'}\cap\Aline {b'}c=\Aline a{a'}\cap\Aline b{b'}=\varnothing$. We have to prove that $\Aline b{b'}\cap\Aline c{c'}=\varnothing$.  Applying the commutativity-puls axiom to the triples $abc$ and $a'b'c'$, we conclude that $\Aline a{c'}\cap\Aline {a'}c=\varnothing$. Since $X$ is bi-Boolean, $\Aline {c'}a\cap\Aline c{a'}=\Aline a{b'}\cap\Aline {b'}a=\Aline c{c'}\cap\Aline a{a'}=\varnothing$ implies $\Aline c{c'}\cap\Aline b{b'}=\varnothing$.
\smallskip

 Asume that $X$ is di-Boolean. To prove that $X$ is bi-Boolean, take any disjoint lines $L,L'\subseteq X$ and distinct points $a,b,c\in L$ and $a',b',c'\in L$ such that $\Aline a{b'}\cap\Aline {a'}b=\Aline b{c'}\cap\Aline {b'}c=\Aline a{a'}\cap\Aline c{c'}=\varnothing$. We have to prove that $\Aline a{a'}\cap\Aline b{b'}=\varnothing$.  Applying the commutativity-puls axiom to the triples $abc$ and $a'b'c'$, we conclude that $\Aline a{c'}\cap\Aline {a'}c=\varnothing$. Since $X$ is di-Boolean, $\Aline {c'}a\cap\Aline c{a'}=\Aline a{b'}\cap\Aline {b'}a=\Aline c{c'}\cap\Aline a{a'}=\varnothing$ implies $\Aline a{a'}\cap\Aline b{b'}=\varnothing$.
\end{proof}

\begin{theorem}\label{t:Boolean<=>di-Boolean+} An affine space $X$ is Boolean if and only if $X$ is di-Boolean and contains a Boolean parallelogram.
\end{theorem}

\begin{proof} The ``only if'' part follows from Proposition~\ref{p:Boole=>bi-di-by-Boole}. To prove the ``if'' part, assume that $X$ is a di-Boolean affine space, containing a Boolean parallelogram. If $\|X\|\ge 4$, then $X$ is Desarguesian and Thalesian, by Theorems~\ref{t:proaffine-Desarguesian} and \ref{t:ADA=>AMA}. 
Let $oexy$ be a Boolean parallelogram in $X$. Then $\Aline oe\parallel \Aline xy$, $\Aline ox\parallel \Aline ey$ and $\Aline oy\parallel \Aline ex$, which implies $\overvector{oe}=\overvector{xy}=\overvector{eo}$ and hence $\overvector{oe}+\overvector{oe}=\overvector{oe}+\overvector{eo}=\overvector{oo}$ and  $1+1=\overvector{oee}+\overvector{oee}=\overvector{ooe}=0$ in the scalar corps $\IR_X$ of the Thalesian space $X$. 
For every points $a,b\in X$ we have $$\overvector{ab}+\overvector{ab}=1\cdot\overvector{ab}+1\cdot\overvector{ab}=(1+1)\cdot\overvector{ab}=0\cdot\overvector{ab}=\overvector{aa}=\overvector{ab}+\overvector{ba}$$ and hence $\overvector{ab}=\overvector{ba}$. Then for every parallelogram $abcd$, we have $$\overvector{ac}=\overvector{ab}+\overvector{bc}=\overvector{ab}+\overvector{cb}=\overvector{dc}+\overvector{cb}=\overvector{db}$$
and hence $\Aline ac\parallel \Aline bd$, witnessing that the parallelogram $abcd$ is Boolean, and so is the liner $X$. 
\smallskip

So, assume that $\|X\|=3$, which means that $X$ is a di-Boolean Playfair plane. In this case, the proof is divided into a series of claims.

\begin{claim}\label{cl:Boolean1} If $abcd$ is a Boolean parallelogram in $X$, then for every points $x\in \Aline ab$ and $y\in\Aline cd$ with $\Aline bc\cap\Aline xy=\varnothing$, the parallelogram $xbcy$ is Boolean.
\end{claim}

\begin{picture}(100,65)(-160,-15)
\linethickness{1pt}
\put(0,0){\color{cyan}\line(0,1){40}}
\put(0,0){\color{teal}\line(1,0){60}}
\put(0,40){\color{teal}\line(1,0){60}}
\put(20,0){\color{cyan}\line(0,1){40}}
\put(0,0){\color{blue}\line(1,2){20}}
\put(0,40){\color{blue}\line(1,-2){20}}
\put(20,0){\color{red}\line(1,1){40}}
\put(20,40){\color{red}\line(1,-1){40}}
\put(60,0){\color{cyan}\line(0,1){40}}
\put(0,0){\circle*{3}}
\put(-2,-9){$a$}
\put(20,0){\circle*{3}}
\put(17,-9){$b$}
\put(60,0){\circle*{3}}
\put(58,-8){$x$}
\put(0,40){\circle*{3}}
\put(-2,43){$d$}
\put(20,40){\circle*{3}}
\put(18,43){$c$}
\put(60,40){\circle*{3}}
\put(58,44){$y$}
\put(10,20){\color{blue}\circle*{4}}
\put(10,20){\color{white}\circle*{3}}
\put(40,20){\color{red}\circle*{4}}
\put(40,20){\color{white}\circle*{3}}
\end{picture}

\begin{proof} By Corollary~\ref{c:parallel-lines<=>}, the lines $\Aline ad$ and $\Aline xy$ are parallel to the line $\Aline bc$. If $\Aline ad=\Aline xy$, then the parallelogram $xbcy=abcd$ is Boolean. So, assume that $\Aline ad\ne\Aline xy$.
In this case $\Aline ad,\Aline bc,\Aline xy$ are three distinct parallel lines in $X$. 
Since the parallelogram $abcd$ is Boolean, $\Aline ac\cap\Aline bd=\varnothing$. Applying Proposition~\ref{p:di-Boolean<=>}, we conclude that $\Aline by\cap\Aline cx=\varnothing$, witnessing that the parallelogram $xbcy$ is Boolean.
\end{proof}

\begin{claim}\label{cl:Boolean2} If $abcd$ is a Boolean parallelogram in $X$, then for every distinct points $a',b'\in \Aline ab$ and $c',d'\in \Aline cd$ with $\Aline {a'}{d'}\parallel \Aline ad$ and $\Aline {b'}{c'}\parallel \Aline bc$, the parallelogram $a'b'c'd'$ is Boolean.
\end{claim}

\begin{picture}(100,65)(-160,-15)
\linethickness{1pt}
\put(0,0){\color{cyan}\line(0,1){40}}
\put(0,0){\color{teal}\line(1,0){80}}
\put(0,40){\color{teal}\line(1,0){80}}
\put(20,0){\color{cyan}\line(0,1){40}}
\put(40,0){\color{cyan}\line(0,1){40}}

\put(0,0){\color{blue}\line(1,2){20}}
\put(0,40){\color{blue}\line(1,-2){20}}
\put(40,0){\color{red}\line(1,1){40}}
\put(40,40){\color{red}\line(1,-1){40}}
\put(80,0){\color{cyan}\line(0,1){40}}
\put(0,0){\circle*{3}}
\put(-2,-9){$a$}
\put(20,0){\circle*{3}}
\put(17,-9){$b$}
\put(40,0){\circle*{3}}
\put(37,-9){$a'$}
\put(80,0){\circle*{3}}
\put(77,-9){$b'$}
\put(0,40){\circle*{3}}
\put(-2,43){$d$}
\put(20,40){\circle*{3}}
\put(18,43){$c$}
\put(40,40){\circle*{3}}
\put(38,43){$d'$}
\put(80,40){\circle*{3}}
\put(78,43){$c'$}
\put(10,20){\color{red}\circle*{4}}
\put(10,20){\color{white}\circle*{3}}
\put(60,20){\color{red}\circle*{4}}
\put(60,20){\color{white}\circle*{3}}
\end{picture}

\begin{proof} If $\Aline {a'}{d'}=\Aline bc$, then the parallelogram $a'b'c'd'=bb'c'c$ is Boolean by Claim~\ref{cl:Boolean1}. So, assume that $\Aline{a'}{d'}\ne \Aline bc$. In this case Claim~\ref{cl:Boolean1} ensures that the parallelogram $ba'd'c$ is Boolean and so is the parallelogram $a'b'c'd'$.
\end{proof}

\begin{claim}\label{cl:Boolean3} If $abcd$ is a Boolean parallelogram in $X$, then for every parallelogram $a'b'c'd'$ with $\Aline ab\parallel \Aline {a'}{b'}$ and $\Aline bc\parallel \Aline {b'}{c'}$ is Boolean.
\end{claim}

\begin{proof} Consider the unique points $a''\in\Aline ab\cap\Aline {a'}{d'}$, $b''\in \Aline ab\cap\Aline {b'}{c'}$, $c''\in \Aline cd\cap\Aline {b'}{c'}$ and $d''\in \Aline cd\cap\Aline {a'}{d'}$. Since $\Aline ab=\Aline {a''}{b''}$ and $\Aline cd=\Aline{c''}{d''}$, we can apply Claim~\ref{cl:Boolean2} and conclude that the parallelogram $a''b''c''d''$ is Boolean. Since $\Aline {a''}{d''}=\Aline {a'}{d'}$ and $\Aline {b''}{c''}=\Aline {b'}{c'}$, we can apply Claim~\ref{cl:Boolean2} and conclude that the parallelogram $a'b'c'd'$ is Boolean.

\begin{picture}(100,110)(-180,-50)

\put(-15,-40){\line(0,1){80}}
\put(15,-40){\line(0,1){80}}
\put(-40,-15){\line(1,0){80}}
\put(-40,15){\line(1,0){80}}
\put(-15,-40){\line(1,0){30}}
\put(-15,40){\line(1,0){30}}
\put(-40,-15){\line(0,1){30}}
\put(40,-15){\line(0,1){30}}

\put(-40,-15){\circle*{2.5}}
\put(-49,-17){$a'$}
\put(-40,15){\circle*{2.5}}
\put(-48,13){$b'$}
\put(-15,-40){\circle*{2.5}}
\put(-17,-49){$a$}
\put(15,-40){\circle*{2.5}}
\put(14,-49){$d$}
\put(40,-15){\circle*{2.5}}
\put(43,-17){$d'$}
\put(40,15){\circle*{2.5}}
\put(43,13){$c'$}
\put(-15,40){\circle*{2.5}}
\put(-17,43){$b$}
\put(15,40){\circle*{2.5}}
\put(13,43){$c$}

\put(-15,-15){\circle*{2.5}}
\put(-26,-25){$a''$}
\put(-15,15){\circle*{2.5}}
\put(-25,17){$b''$}
\put(15,-15){\circle*{2.5}}
\put(17,-25){$d''$}
\put(15,15){\circle*{2.5}}
\put(17,17){$c''$}
\end{picture}
\end{proof} 

\begin{claim}\label{cl:Boolean4} If $abcd$ is a Boolean parallelogram in $X$, then for every points $x,y\in \Aline bc$ with $\Aline ax\cap\Aline dy=\varnothing$, the parallelogram $axyd$ is Boolean.
\end{claim}

\begin{proof}
 If $x=b$, then $\Aline dy\parallel \Aline ax=\Aline ab\parallel \Aline cd$ implies $\Aline dy=\Aline dc$ and $y=c$. In this case the parallelogram $axyd=abcd$ is Boolean by the assumption. So, assume that $x\ne b$. Find a unique point $z\in\Aline ad$ such that $\Aline xz\parallel\Aline ab$ and observe that $x\ne b$ implies $z\ne a$. By Claim~\ref{cl:Boolean3}, the parallelogram $abxz$ is Boolean and hence $\Aline ax\cap \Aline bz=\varnothing$. 
 
\begin{picture}(100,65)(-150,-10)

\put(0,0){\line(1,0){60}}
\put(0,0){\line(0,1){40}}
\put(60,0){\line(0,1){40}}

\put(20,0){\line(0,1){40}}
\put(0,0){\color{blue}\line(1,2){20}}
\put(0,40){\line(1,0){80}}
\put(60,0){\line(1,2){20}}
\put(0,40){\color{blue}\line(1,-2){20}}
\put(20,40){\color{red}\line(1,-1){40}}
\put(0,0){\color{red}\line(2,1){80}}

\put(0,0){\circle*{2.5}}
\put(-4,-8){$a$}
\put(20,0){\circle*{2.5}}
\put(17,-8){$z$}
\put(60,0){\circle*{2.5}}
\put(60,-8){$d$}
\put(0,40){\circle*{2.5}}
\put(-4,43){$b$}
\put(20,40){\circle*{2.5}}
\put(18,43){$x$}
\put(60,40){\circle*{2.5}}
\put(58,43){$c$}
\put(80,40){\circle*{2.5}}
\put(82,43){$y$}

\put(10,20){\circle*{3}}
\put(10,20){\color{white}\circle*{2}}
\put(40,20){\circle*{3}}
\put(40,20){\color{white}\circle*{2}}

\end{picture}

By Proposition~\ref{p:di-Boolean<=>}, $\Aline xz\cap\Aline ba=\varnothing=\Aline bz\cap\Aline xa=\Aline xa\cap\Aline yd$ imply $\Aline xd\cap\Aline ya=\varnothing$, witnessing that the parallelogram $axyd$ is Boolean.
\end{proof} 

The following claim completes the proof of the ``if'' part of Theorem~\ref{t:Boolean<=>di-Boolean+}.

\begin{claim}\label{cl:Boolean5} If the liner $X$ contains a Boolean parallelogram, then every parallelogram in $X$ is Boolean.
\end{claim}

\begin{proof} Let $abcd$ be a Boolean parallelogram in $X$ and let $a'b'c'd'$ be any parallelogram in $X$. Choose a parallelogram $a_1b_1c_1d_1$ such that $\Aline {a_1}{b_1}\parallel \Aline ab$, $\Aline {a_1}{d_1}\parallel \Aline ad$ and $a_1=a'$. By Claim~\ref{cl:Boolean3}, the parallelogram $a_1b_1c_1d_1$ is Boolean. If $\Aline{a_1}{d_1}=\Aline{a'}{b'}$, then find a parallelogram $a_2b_2c_2d_2$ such that $a_2d_2=a'b'$, $\Aline {a_2}{b_2}=\Aline{a_1}{b_1}$ and $\Aline{b_2}{c_2}=\Aline {c'}{d'}$. By Claim~\ref{cl:Boolean3}, the parallelogram $a_2b_2c_2d_2$ is Boolean, and by Claim~\ref{cl:Boolean4}, the parallelogram  $a'b'c'd'$ is Boolean.

So, assume that $\Aline{a_1}{d_1}\ne\Aline {a'}{b'}$. In this case we can consider a parallelogram $a_2b_2c_2d_2$ such that $\Aline {a_2}{b_2}=\Aline{a_1}{b_1}$, $\Aline{a_2}{d_2}=\Aline {a_1}{d_1}$, $d_2\in \Aline {c'}{d'}$ and $b'\in\Aline {b_2}{c_2}$.  By Claim~\ref{cl:Boolean3}, the parallelogram $a_2b_2c_2d_2$ is Boolean.

\begin{picture}(100,110)(-150,-40)

\put(0,0){\line(1,0){75}}
\put(0,0){\line(0,1){60}}
\put(0,30){\line(1,0){40}}
\put(40,0){\line(0,1){30}}
\put(0,0){\line(1,2){30}}
\put(0,0){\line(2,-1){60}}
\put(60,-30){\line(1,2){45}}
\put(0,60){\line(1,0){105}}
\put(75,0){\line(0,1){60}}
\put(30,60){\line(2,-1){60}}

\put(0,0){\circle*{2.5}}
\put(-10,-10){$a'$}
\put(0,-10){$a_1$}
\put(-11,4){$a_2$}
\put(4,3){$a_3$}
\put(40,0){\circle*{2.5}}
\put(38,-10){$d_1$}
\put(40,30){\circle*{2.5}}
\put(42,32){$c_1$}
\put(0,30){\circle*{2.5}}
\put(-10,28){$b_1$}
\put(0,60){\circle*{2.5}}
\put(-10,60){$b_2$}
\put(30,60){\circle*{2.5}}
\put(23,63){$b'$}
\put(32,63){$b_3$}
\put(75,0){\circle*{2.5}}
\put(77,-5){$d_3$}
\put(65,4){$d_2$}
\put(75,60){\circle*{2.5}}
\put(73,64){$c_2$}
\put(90,30){\circle*{2.5}}
\put(92,25){$c'$}
\put(105,60){\circle*{2.5}}
\put(106,63){$c_3$}
\put(60,-30){\circle*{2.5}}
\put(62,-35){$d'$}

\end{picture}

Consider the parallelogram $a_3b_3c_3d_3$ such that ${a_3}{d_3}={a_2}{d_2}$ and ${a_3}{b_3}={a_2}{b'}={a'}{b'}$. Since $b_3=b'\in\Aline{b_2}{c_2}$, the parallel lines $\Aline {b_2}{c_2}$ and $\Aline {b_3}{c_3}$ coincide. Since the parallelogram $a_2b_2c_2d_2$ is Boolean, we can apply Claim~\ref{cl:Boolean4} and conclude that the parallelogram $a_3b_3c_3d_3$ is Boolean.

Since $d_3=d_2\in \Aline{c'}{d'}$, the parallel lines $\Aline{c_3}{d_3}\parallel \Aline{a_3}{b_3}=\Aline{a'}{b'}$ and $\Aline{c'}{d'}$ coincide. 
Since  the parallelogram $a_3b_3c_3d_3$ is Boolean, we can apply Claim~\ref{cl:Boolean4} and conclude that the parallelogram $a'b'c'd'$ is Boolean.
\end{proof}
\end{proof} 
%

\begin{exercise} Show that the Moulton plane is bi-Boolean, di-Boolean, and by-Boolean, but not Boolean.

{\em Hint:} Use the topological properties of the Moulton plane.
\end{exercise}

\section{Uno-Thalesian Playfair liners are Boolean or Thalesian}

The main result of this section is the following characterization of uno-Thalesian affine spaces.

\begin{theorem}\label{t:uno-Thales<=>} A Playfair liner is uno-Thalesian if and only if it is Thalesian or Boolean.
\end{theorem}

\begin{proof} The ``if'' part of this theorem follows from Definitions~\ref{d:para-Desargues}, \ref{d:uno-Thalesian} and Proposition~\ref{p:Boole=>uno-Thales}. The ``only if'' of Theorem~\ref{t:uno-Thales<=>} is proved in Lemma~\ref{l:uno-Thalesian=>BorT} (whose proof uses  five other lemmas).

\begin{lemma}\label{l:nT=>two-sided} Every uno-Thalesian Playfair liner is invertible-plus.
\end{lemma}

\begin{proof} Assume that a Playfair liner $X$ is uno-Thalesian. To show that $X$ is invertible-plus, take any points $a,b,c,\alpha,\beta,\gamma\in X$ and $o\in\Aline a\alpha\cap\Aline b\beta\cap\Aline c\gamma$ with 
$$(\Aline ac\parallel \Aline o\beta\parallel\Aline ob\parallel\Aline\alpha\gamma\nparallel \Aline b\gamma\parallel\Aline o\alpha\parallel \Aline oa\parallel\Aline c\beta)\;\wedge\;(\Aline ab\parallel\Aline co).$$
We have to show that $\Aline \alpha\beta\parallel\Aline co$.

\begin{picture}(150,80)(-180,-45)
\linethickness{0.75pt}
\put(0,0){\color{blue}\line(1,0){30}}
\put(0,0){\color{cyan}\line(0,1){30}}
\put(0,0){\color{blue}\line(-1,0){30}}
\put(0,0){\color{cyan}\line(0,-1){30}}

\put(0,30){\color{blue}\line(1,0){30}}
\put(30,0){\color{cyan}\line(0,1){30}}
\put(-30,-30){\color{blue}\line(1,0){30}}
\put(-30,-30){\color{cyan}\line(0,1){30}}
\put(-30,0){\color{red}\line(1,1){30}}
\put(0,-30){\color{red}\line(1,1){30}}
{\linethickness{1pt}
\put(0,0){\color{red}\line(1,1){30}}
\put(0,0){\color{red}\line(-1,-1){30}}
}
\put(0,0){\circle*{3}}
\put(2,-7){$o$}
\put(30,0){\circle*{3}}
\put(33,-2){$\alpha$}
\put(-30,0){\circle*{3}}
\put(-38,-2){$a$}
\put(0,30){\circle*{3}}
\put(-7,31){$b$}
\put(0,-30){\circle*{3}}
\put(2,-36){$\beta$}
\put(30,30){\circle*{3}}
\put(32,31){$\gamma$}
\put(-30,-30){\circle*{3}}
\put(-37,-36){$c$}
\end{picture}

If $a=c$, then $\Aline ac\parallel\Aline o\beta\parallel \Aline ob\parallel\Aline\alpha\gamma\nparallel \Aline b\gamma\parallel\Aline o\alpha\parallel\Aline oa\parallel\Aline c\beta$ implies $b=o=\beta$, $\alpha=\gamma$, and $a\ne o\ne \alpha$. Then $\Aline \alpha\beta=\Aline \alpha o=\Aline oa=\Aline co$ and hence $\Aline \alpha\beta\parallel \Aline co$. By analogy we can show that $c=\beta$ implies $\Aline \alpha\beta\parallel \Aline co$. 

So, assume that $a\ne c\ne\beta$. In this case $\Aline b\gamma,\Aline ao,\Aline c\beta$ are distinct parallel lines and $\Aline ab\cap\Aline o\gamma=\varnothing=\Aline ac\cap\Aline o\beta$. Since $b\in\Aline o\beta$, we can apply Proposition~\ref{p:uno-Thalesian1} to the triangles $abc$ and $o\gamma\beta$ and conclude that $\Aline cb\cap \Aline \beta\gamma=\varnothing$.



Since $\beta\in\Aline ob$, we can apply Proposition~\ref{p:uno-Thalesian1} to the triangles $ocb$ and $\alpha\beta\gamma$ and conclude that  $\Aline\alpha\beta\parallel\Aline oc$.

%

\end{proof}

Given any triangle $xyz$ in a Playfair liner $X$, let $\diamond xyz$ be a unique point $p\in X$ completing the triangle $xyz$ to a parallelogram such that $\Aline xp\parallel \Aline yz$ and $\Aline pz\parallel \Aline xy$. It is clear that $\diamond xyz=\diamond zyx$.

%

\begin{lemma}\label{l:uno-Thalesian3} For any triangle $xoe$ in a Playfair liner $X$, the point $\diamond ex(\diamond xoe)$ belongs to the set $\Aline oe\setminus\{e\}$.
\end{lemma}

\begin{proof} Consider the points $y\defeq\diamond xoe$ and $z\defeq \diamond ex(\diamond xoe)=\diamond exy$, and observe that $\Aline oe\parallel \Aline xy\parallel \Aline ez$ and hence $x\ne y$ and $z\ne e$ (because $o\ne e$). Therefore, $\Aline oe$ and $\Aline ez$ are two lines that contain the point $e$ and are parallel to the line $\Aline xy$. Since $X$ is Playfair, $\Aline oe=\Aline ez$ and hence $z\in \Aline oe\setminus\{e\}$. 

\begin{picture}(60,50)(-150,-5)

\linethickness{0.75pt}
\put(0,0){\color{red}\line(1,0){80}}
\put(10,30){\color{red}\line(1,0){40}}
\put(40,0){\color{blue}\line(1,3){10}}
\put(0,0){\color{blue}\line(1,3){10}}
\put(10,30){\color{violet}\line(1,-1){30}}
\put(50,30){\color{violet}\line(1,-1){30}}

\put(0,0){\circle*{3}}
\put(-5,-8){$o$}
\put(10,30){\circle*{3}}
\put(2,32){$x$}
\put(40,0){\circle*{3}}
\put(38,-8){$e$}
\put(50,30){\circle*{3}}
\put(53,32){$y$}
\put(80,0){\circle*{3}}
\put(80,-8){$z$}
\end{picture}

\end{proof}

\begin{lemma}\label{l:uno-Thalesian4} If a Playfair liner $X$ is uno-Thalesian, then for any triangle $xoe\in X^3$ and any point $y\in \Aline ox\setminus\{o\}$ the equality $\diamond ex(\diamond xoe)=\diamond ey(\diamond yoe)$ holds.
\end{lemma}

\begin{proof} Given any point $y\in \Aline ox\setminus\{o\}$, consider the points $x'\defeq\diamond xoe$, $y'\defeq yoe$, $p\defeq \diamond ex(\diamond xoe)=\diamond exx'$ and $q\defeq\diamond ey(\diamond yoe)=\diamond eyy'$.
We have to prove that the points $p$ and $q$ are equal. Lemma~\ref{l:uno-Thalesian3} ensures that $p,q\in\Aline oe$.

\begin{picture}(100,90)(-180,-15)

\linethickness{0.75pt}
\put(0,0){\color{cyan}\line(0,1){60}}
\put(20,0){\color{cyan}\line(0,1){60}}
\put(0,0){\color{red}\line(1,2){20}}
\put(20,0){\color{red}\line(1,2){20}}
\put(0,40){\color{blue}\line(1,0){40}}
\put(0,60){\color{teal}\line(1,-1){20}}
\put(20,60){\color{teal}\line(1,-1){20}}
\put(0,0){\color{blue}\line(1,0){20}}
\put(0,60){\color{blue}\line(1,0){20}}

\put(0,0){\circle*{2.5}}
\put(-4,-8){$y$}
\put(0,60){\circle*{2.5}}
\put(-4,63){$x$}
\put(20,0){\circle*{2.5}}
\put(20,-8){$y'$}
\put(20,40){\circle*{2.5}}
\put(22,42){$e$}
\put(20,60){\circle*{2.5}}
\put(20,63){$x'$}
\put(40,40){\circle*{2.5}}
\put(42,34){$q$}
\put(42,43){$p$}
\put(0,40){\circle*{2.5}}
\put(-7,38){$o$}
\end{picture}

If $x=y$, then $y'=x'$ and $p=q$, by the definition of the parallelogram operation $\diamond$. So, assume that $y\ne x$. In this case, $\Aline oe$, $\Aline x{x'}$ and $\Aline y{y'}$ are three distinct parallel lines in the plane $\Pi\defeq\overline{\{o,e,x\}}$. Taking into account that $\Aline e{x'}\parallel \Aline ox=\Aline oy\parallel \Aline e{y'}$, we conclude that $e\in \Aline e{x'}=\Aline e{y'}=\Aline {x'}{y'}$. Since $\Aline xe\parallel \Aline {x'}{p}$ and $e\in \Aline {x'}{y'}\parallel \Aline xy$, we can apply Proposition~\ref{p:uno-Thalesian1} and  conclude that $\Aline ye\parallel \Aline {y'}p$. On the other hand, the definition of the point $q=\diamond eyy'$ ensures that $\Aline {y'}{q}\parallel \Aline ye\parallel \Aline {y'}{p}$ and hence $\Aline {y'}p=\Aline {y'}q$. Then $\{p\}=\Aline oe\cap\Aline {y'}p=\Aline oe\cap\Aline {y'}q=\{q\}$ and $p=q$.
\end{proof}

\begin{lemma}\label{l:uno-Thalesian5} If a Playfair plane $X$ is uno-Thalesian and not Boolean, then for any triangle $xoe\in X^3$ and any point $y\in X$ with $\Aline xy\parallel \Aline oe$, the equality
$\diamond ex(\diamond xoe)=\diamond ey(\diamond yoe)$ holds.
\end{lemma}

\begin{proof}  Consider the points $x'\defeq\diamond xoe$ and $y'\defeq\diamond yoe$. We have to prove that the points $t\defeq\diamond exx'=\diamond ex(\diamond xoe)$ and $\tau\defeq \diamond eyy'=\diamond ey(\diamond y oe)$ are equal. The equality $t=\tau$ will be proved as soon as we show that the lines $\Aline ey$ and $\Aline t{y'}$ are parallel.

\begin{picture}(100,70)(-150,15)
\put(0,30){\line(1,0){60}}
\put(0,60){\line(1,0){150}}
\put(0,60){\line(1,-1){30}}
\put(30,60){\line(1,-1){30}}
\put(30,30){\line(0,1){30}}
\put(0,30){\line(0,1){30}}

\put(0,30){\line(4,1){120}}
\put(30,30){\line(4,1){120}}
\put(30,30){\color{red}\line(3,1){90}}
\put(60,30){\color{red}\line(3,1){90}}


\put(120,60){\circle*{3}}
\put(118,65){$y$}
\put(150,60){\circle*{3}}
\put(148,65){$y'$}
\put(0,30){\circle*{3}}
\put(-3,22){$o$}
\put(30,30){\circle*{3}}
\put(27,22){$e$}
\put(60,30){\color{red}\circle*{3}}
\put(54,22){$t$}
\put(62,22){$\tau$}
\put(0,60){\circle*{3}}
\put(-2,64){$x$}
\put(30,60){\circle*{3}}
\put(28,64){$x'$}
\end{picture}

Consider the points $z\defeq \diamond oex$, $z'\defeq \diamond oxz$, $e'\defeq exo$, and $o'\defeq \diamond oee'$.  By Lemma~\ref{l:nT=>two-sided}, the uno-Thalesian plane $\Pi$ is  invertible-plus, which implies $\Aline {z'}{o'}\parallel \Aline zo=\Aline o{e'}$ and $\Aline {z'}{o'}\cap\Aline o{e'}=\varnothing$. The choice of the point $z=\diamond oex$ ensures that $\Aline zo\parallel \Aline xe$. Assuming that $\Aline zx=\Aline{o'}{e'}$, we conclude that $e'=x'$ and hence $\Aline o{x'}=\Aline o{e'}=\Aline zo\parallel \Aline xe$, which implies that the parallelogram $oxx'e$ is Boolean. By Proposition~\ref{p:uno-Thalesian=>bi-di-by-Boolean}, the uno-Thalesian Playfair liner $X$ is di-Boolean, and Theorem~\ref{t:Boolean<=>di-Boolean+}, the di-Boolean Playfair liner $X$ is Boolean (because it contains the Boolean parallelogram $oxx'e$), which contradicts our assumption. This contradiction shows that $\Aline zx\ne\Aline{o'}{e'}$ and hence $\Aline zx$, $\Aline oe$, $\Aline {o'}{e'}$ are three distinct parallel lines in the Playfair plane $X$.

\begin{picture}(100,95)(-150,-15)
\put(0,0){\line(1,0){30}}
\put(0,0){\line(0,1){30}}
\put(0,0){\line(-1,1){30}}
\put(0,30){\line(1,-1){30}}
\put(30,0){\line(0,1){30}}
\put(-30,30){\line(1,1){30}}
\put(0,30){\line(1,1){30}}
\put(0,0){\line(1,1){30}}
\put(30,0){\line(1,1){30}}

\put(-30,30){\line(1,0){90}}
\put(-30,60){\line(1,0){180}}
\put(0,60){\line(1,-1){30}}
\put(30,60){\line(1,-1){30}}
\put(30,30){\line(0,1){30}}
\put(0,30){\line(0,1){30}}
\put(-30,30){\line(0,1){30}}
\put(0,30){\line(-1,1){30}}
\put(0,30){\line(4,1){120}}
\put(30,30){\line(4,1){120}}
\put(30,30){\line(3,1){90}}


\put(-30,30){\circle*{3}}
\put(-39,28){$z'$}
\put(-30,60){\circle*{3}}
\put(-32,64){$z$}


\put(120,60){\circle*{3}}
\put(118,65){$y$}
\put(150,60){\circle*{3}}
\put(148,65){$y'$}
\put(0,0){\circle*{3}}
\put(-3,-10){$o'$}
\put(0,30){\circle*{3}}
\put(-6,22){$o$}
\put(30,30){\circle*{3}}
\put(25,21){$e$}
\put(60,30){\circle*{3}}
\put(60,22){$t$}
\put(0,60){\circle*{3}}
\put(-2,64){$x$}
\put(30,60){\circle*{3}}
\put(28,64){$x'$}
\put(30,0){\circle*{3}}
\put(27,-10){$e'$}
\end{picture}

 Since the affine plane $\Pi$ is uno-Thalesian, $\Aline {z'}{o'}\cap\Aline o{e'}=\varnothing=\Aline x{o'}\cap\Aline {x'}{e'}$ and $o\in \Aline x{o'}$ imply $\Aline x{z'}\cap\Aline {x'}o=\varnothing$, see Proposition~\ref{p:uno-Thalesian1}. By the same Proposition~\ref{p:uno-Thalesian1}, $\Aline zo\cap\Aline {z'}{o'}=\varnothing=\Aline {z'}e\cap\Aline z{x'}$ and $o\in \Aline {z'}e$ imply $\Aline o{x'}\parallel \Aline {o'}e$, and $\Aline xe\cap \Aline {x'}t=\varnothing=\Aline {x'}{e'}\cap\Aline x{o'}$ and $e\in\Aline {x'}{e'}$ imply $\Aline {o'}e\parallel \Aline {e'}t$.

Since $\Pi$ is Playfair, there exist a unique point $u\in\Aline oe$ such that $\Aline y{u}\parallel \Aline ox$. By Proposition~\ref{p:uno-Thalesian1}, $\Aline zo\cap\Aline {z'}{o'}=\varnothing=\Aline {z'}{u}\cap\Aline zy$ and $u\in\Aline {z'}o$ imply $\Aline {o'}{u}\cap\Aline oy=\varnothing$ and hence $\Aline{o'}u\parallel \Aline {o}{y}\parallel \Aline e{y'}$. By Proposition~\ref{p:uno-Thalesian1}, $\Aline xe\cap\Aline e{e'}=\varnothing=\Aline o{u}\parallel \Aline xy$ and $u\in\Aline oe$ imply $\Aline ey\cap \Aline {e'}{u}=\varnothing$.

\begin{picture}(100,95)(-150,-15)
\put(0,0){\line(1,0){30}}
\put(0,0){\line(0,1){30}}
\put(0,0){\line(-1,1){30}}
\put(0,30){\line(1,-1){30}}
\put(30,0){\line(0,1){30}}
\put(-30,30){\line(1,1){30}}
\put(0,30){\line(1,1){30}}
\put(0,0){\line(1,1){30}}
\put(30,0){\line(1,1){30}}

\put(-30,30){\line(1,0){150}}
\put(-30,60){\line(1,0){180}}
\put(0,60){\line(1,-1){30}}
\put(30,60){\line(1,-1){30}}
\put(30,30){\line(0,1){30}}
\put(0,30){\line(0,1){30}}
\put(-30,30){\line(0,1){30}}
\put(0,30){\line(-1,1){30}}
\put(0,30){\color{blue}\line(4,1){120}}
\put(0,0){\color{blue}\line(4,1){120}}
\put(30,30){\color{blue}\line(4,1){120}}
\put(30,30){\color{teal}\line(3,1){90}}
\put(30,0){\color{teal}\line(3,1){90}}
\put(120,30){\line(0,1){30}}

\put(-30,30){\circle*{3}}
\put(-39,28){$z'$}
\put(-30,60){\circle*{3}}
\put(-32,64){$z$}


\put(120,60){\circle*{3}}
\put(118,64){$y$}
\put(120,30){\circle*{3}}
\put(122,28){$u$}
\put(150,60){\circle*{3}}
\put(148,64){$y'$}
\put(0,0){\circle*{3}}
\put(-3,-10){$o'$}
\put(0,30){\circle*{3}}
\put(-6,22){$o$}
\put(30,30){\circle*{3}}
\put(25,21){$e$}
\put(60,30){\circle*{3}}
\put(60,22){$t$}
\put(0,60){\circle*{3}}
\put(-2,64){$x$}
\put(30,60){\circle*{3}}
\put(28,64){${x'}$}
\put(30,0){\circle*{3}}
\put(27,-10){$e'$}
\end{picture}

By Theorem~\ref{t:Boolean<=>di-Boolean+}, the non-Boolean Playfair plane $X$ contains no Boolean parallelograms, which guarantees that $z'\ne e$ and hence the parallel lines $\Aline {z'}x$ and $\Aline {o'}e$ are disjoint. If $x=y'$, then $\Aline{o'}u\parallel \Aline e{y'}=\Aline ex\parallel\Aline {o'}{z'}$ implies $u={z'}$ and hence $\Aline u{y'}=\Aline{z'}{x}\parallel \Aline {o'}e$. By analogy we can show that $u=z'$ implies $\Aline u{y'}\parallel \Aline {o'}e$. So, assume that $x\ne y'$ and $z'\ne u$. In this case $o\in \Aline{z'}u\parallel \Aline x{y'}$. Since $\Aline xe\parallel\Aline {z'}{o'}$, $\Aline e{y'}\parallel \Aline {o'}{u}$ and $\Aline {z'}x\cap\Aline {o'}e=\varnothing$, we can apply Proposition~\ref{p:uno-Thalesian2} and conclude that the line $\Aline u{y'}$ is parallel to the line $\Aline{o'}e$, which is parallel to the lines $\Aline o{x'}$ and $\Aline{e'}t$.



Since $X$ is not Boolean, the parallelogram $oxx'e$ is not Boolean and hence $o\ne t$ and $\Aline o{x'}\cap\Aline {e'}t=\varnothing$. It follows from $\Aline xy\parallel \Aline oe$ that $x\ne y$ and hence $x'\ne y'$. In this case the parallelity $\Aline u{y'}\parallel \Aline o{x'}$ implies $o\ne u$.
Applying Proposition~\ref{p:uno-Thalesian1} to the triangles $x'ty'$ and $oe'u$, we conclude that $\Aline t{y'}\parallel \Aline {e'}{u}\parallel \Aline ey$ and hence $\Aline t{y'}\parallel \Aline ey$ and $t=\tau$.

%
\end{proof}

\begin{lemma}\label{l:uno-Thalesian6} If a Playfair plane $X$ is uno-Thalesian, then for any distinct points $o,e\in X$, the set $\{\diamond ex(\diamond xoe):x\in X\setminus\Aline oe\}$ is a singleton.
\end{lemma}

\begin{proof} Given any points $x,y\in X\setminus\Aline oe$, we need to show that $
\diamond ex(\diamond xoe)=\diamond ey(\diamond yoe)$. If $X$ is Boolean, then for the points $x'\defeq \diamond xoe$ and $y'\defeq\diamond yoe$, the parallelograms $xoex'$ and $yoey'$ are Boolean and hence 
$\diamond ex(\diamond xoe)=\diamond exx'=o=\diamond eyy'=\diamond ey(\diamond yoe)$.

So, assume that $X$ is not Boolean. Since $X$ is a Playfair plane, there exists a point $z\in \Aline ox$ such that $\Aline yz\parallel \Aline oe$. Lemmas~\ref{l:uno-Thalesian4} and \ref{l:uno-Thalesian5} ensure that
$$\diamond ex(\diamond xoe)=\diamond ez(\diamond zoe)=\diamond ey(\diamond yoe).$$
\end{proof}

\begin{lemma}\label{l:uno-Thalesian=>BorT} Every uno-Thalesian Playfair liner is Boolean or Thalesian.
\end{lemma}

\begin{proof} Let $X$ be a uno-Thalesian Playfair liner. If $\|X\|\ne 3$, then $X$ is Desarguesian and Thalesian, by Corollary~\ref{c:affine-Desarguesian} and Theorem~\ref{t:ADA=>AMA}. So, assume that $\|X\|=3$, which means that $X$ is a Playfaor plane. Assuming that $X$ is not Boolean, we shall prove that $X$ is Thalesian. By Lemma~\ref{l:uno-Thalesian6}, for any distinct points $o,e\in X$, the set $$\{\diamond ex(\diamond xoe):x\in X\setminus\Aline oe\}$$ contains a unique point, which will be denoted by $\two_{oe}$ (because $2_{oe}$ is equal to $e+e=e\puls e$ in the ternary-ring $\Delta=\Aline oe$ of any affine base $uow$ with origin $o$ and diunit $e$). Lemma~\ref{l:uno-Thalesian3} ensures that $\two_{oe}\in \Aline oe\setminus\{e\}$.

\begin{claim}\label{cl:two-shift} Let $S$ be a line shift between two parallel lines $\dom[S]$ and $\rng[S]$ in  $X$. For any distinct points $o,e\in\dom[S]$ and their images $o'\defeq S(o)$ and $e'\defeq S(e')$ we have $S(\two_{oe})=\two_{o'e'}$.
\end{claim}

\begin{proof} Since $S$ is a line shift between parallel lines, $e'o'oe$ is a parallelogram and hence $e'=\diamond o'oe$. Lemma~\ref{l:uno-Thalesian6} ensures that $$t\defeq \two_{oe}=\diamond eo'(\diamond o'oe)=\diamond ee'(\diamond e'oe)=\diamond ee't',$$ where $t'\defeq \diamond e'oe$. The equality $t=\diamond ee't'$ implies $\Aline t{t'}\parallel \Aline e{e'}$ and hence $t'=S(t)=S(\two_{oe})$. On the other hand, 
$S(\two_{oe})=t'=\diamond e'oe=\diamond e'o(\diamond oo'e')=\two_{o'e'}.$

\begin{picture}(100,65)(-180,-5)

\put(0,0){\line(1,0){30}}
\put(0,20){\line(1,0){30}}
\put(0,40){\line(1,0){30}}
\put(0,0){\line(0,1){40}}
\put(30,0){\line(0,1){40}}
\put(0,0){\line(3,2){30}}
\put(0,20){\line(3,-2){30}}
\put(0,20){\line(3,2){30}}
\put(0,40){\line(3,-2){30}}


\put(0,0){\circle*{2.5}}
\put(-8,-2){$o$}
\put(0,20){\circle*{2.5}}
\put(-8,18){$e$}
\put(0,40){\circle*{2.5}}
\put(-5,43){$t$}
\put(30,0){\circle*{2.5}}
\put(33,-2){$o'$}
\put(30,20){\circle*{2.5}}
\put(33,18){$e'$}
\put(30,40){\circle*{2.5}}
\put(30,43){$t'$}
\end{picture}

\end{proof}

\begin{claim}\label{cl:two-shear} Let $P$ be a line projection between concurrent lines $\dom[P]$ and $\rng[P]$ in $X$. Let $o\in\dom[P]\cap\rng[P]$ the center of the line projection $P$. For any point $e\in\dom[P]\setminus\{o\}$  and its image $e'\defeq P(e)$ we have $P(\two_{oe})=\two_{oe'}$.
\end{claim}

\begin{proof} Consider the point $p\defeq \diamond eoe'=\diamond e'oe$ and observe that $$t\defeq \two_{oe}=\diamond ee'(\diamond e'oe)=\diamond ee'p\quad\mbox{and}\quad t'\defeq \two_{o'e'}=\diamond e'e(\diamond eo'e')=\diamond e'ep.$$ The equalities $t=\diamond ee'p$ and $t'=\diamond e'ep$ imply $\Aline tp\parallel \Aline e{e'}\parallel \Aline p{t'}$ and hence $\Aline t{t'}\subseteq \Aline tp=\Aline p{t'}\parallel \Aline e{e'}$ and $\two_{o'e'}=t'=P(t)=P(\two_{ot})$, by the definition of the line projection $P$.

\begin{picture}(100,85)(-180,-10)

\put(0,0){\line(1,2){30}}
\put(0,0){\line(-1,2){30}}
\put(15,30){\line(-1,2){15}}
\put(-15,30){\line(1,2){15}}
\put(-30,60){\line(1,0){60}}
\put(-15,30){\line(1,0){30}}

\put(0,0){\circle*{2.5}}
\put(-3,-8){$o$}
\put(15,30){\circle*{2.5}}
\put(18,27){$e'$}
\put(-15,30){\circle*{2.5}}
\put(-23,27){$e$}
\put(30,60){\circle*{2.5}}
\put(31,63){$t'$}
\put(-30,60){\circle*{2.5}}
\put(-33,63){$t$}
\put(0,60){\circle*{2.5}}
\put(-2,65){$p$}
\end{picture}

\end{proof}

\begin{claim}\label{cl:two-projection} Let $P$ be a line projection between lines $\dom[P]$ and $\rng[P]$ in $X$. For any distinct points $o,e\in\dom[P]$  and their images $o'\defeq P(o)$ and $e'\defeq P(e)$ we have $P(\two_{oe})=\two_{o'e'}$.
\end{claim}

\begin{proof} If $\dom[P]=\rng[P]$, then $P$ is the identity map of the line $\dom[P]=\rng[P]$ and hence $o'e'=oe$ and $P(\two_{oe})=\two_{oe}=\two_{o'e'}$. So, assume that $\dom[P]\ne\rng[P]$. If $\dom[P]\cap\rng[P]=\varnothing$, then $P$ is a line shift and $P(\two_{oe})=\two_{o'e'}$, by Claim~\ref{cl:two-shift}. So, assume that $\dom[P]\cap \rng[P]$ are concurrent lines. Since the liner $X$ is Playfair, there exists a unique line $L\subseteq X$ such that $o\in L$ and $L\parallel \rng[P]$. Since $P$ is a line projection, there exists a direction $\boldsymbol \delta\in \partial X$ such that $$P=\{(x,y)\in\dom[P]\times\rng[P]:\Aline y\in \Aline x{\boldsymbol \delta}\}.$$ By Proposition~\ref{p:flat-relation}, the lines $\dom[P]$ and $\rng[P]$ do not belong to the direction $\boldsymbol\delta$.

Consider the line projection $$S\defeq\{(x,y)\in \dom[P]\times L:y\in \Aline x{\boldsymbol\delta}\}$$ and the line shift $$T\defeq\{(x,y)\in L\times\rng[P]:y\in\Aline x{\boldsymbol\delta}\}.$$ Observe that $P=T\circ S$ and $S(o)=o$. Consider the point $s\defeq S(e)\in L$ and observe that $e'=P(e)=TS(e)=T(s)$ and $o'=P(o)=TS(o)=T(o)$. By Claims~\ref{cl:two-shift} and \ref{cl:two-shear}, $$P(\two_{oe})=TS(\two_{oe})=T(\two_{os})=\two_{o'e'}.$$
\end{proof}

\begin{claim}\label{cl:two-affinity} Let $n\in\IN$ and $P_1,\dots,P_n$ be line projections such that their composition $A\defeq P_n\cdots P_1$ is a line affinity. For every distinct points $o,e\in \dom[A]$ and their images $o'\defeq A(o)$ and $e'\defeq A(e)$ we have $A(\two_{oe})=\two_{o'e'}$.
\end{claim} 

\begin{proof} This claim will be proved by induction on $n$. For $n=1$, the claim follows from Claim~\ref{cl:two-projection}. Assume that the claim is true for some $n\in\IN$. Choose any line projections $P_1,\dots,P_n,P_{n+1}$ such that $A\defeq P_{n+1}\cdots P_1$ is a line affinity. Then for every $k\in\{1,\dots,n+1\}$ the composition $A_k\defeq P_k\cdots P_1$ is a line affinity. Fix any distinct points $o,e\in\dom[A]=\dom[P_1]$ and consider their images $o'\defeq A(o)$ and $e'\defeq A_n(e)$. Also consider the points $o''\defeq A_n(o)$ and $e''\defeq A_n(e)$. The equality $A=P_{n+1}A_n$ implies $o'e'=P_{n+1}o''e''$. The inductive hypothesis ensures that $A_n(\two_{oe})=\two_{o''e''}$. On the other hand, Claim~\ref{cl:two-projection} guarantees that $P_{n+1}(\two_{o''e''})=\two_{o'e'}$. Then $$A(\two_{oe})=P_{n+1}A_n(\two_{oe})=P_{n+1}(\two_{o''e''})=\two_{o'e'}.$$
\end{proof}

\begin{claim} The plane $X$ is Thalesian.
\end{claim}

\begin{proof} Claim~\ref{cl:two-affinity} implies that for every line affinity $A$ in $X$ and every distinct points $o,e\in \dom[A]$ and their images $o'\defeq A(o)$ and $e'\defeq A(e)$ we have $A(\two_{oe})=\two_{o'e'}$. In particular, for any line affinity $A$ with $Aoe=oe$ we have $A(\two_{oe})=\two_{oe}$, which implies that the triple $o\two_{oe}e$ is Desarguesian. Then $s\defeq\overvector{o\two_{oe}e}$ is a scalar. Since $X$ is not Boolean, $\two_{oe}\ne o$ and hence $s\ne 0$. On the other hand, Lemma~\ref{l:uno-Thalesian3} ensures that $\two_{oe}\ne e$ and hence $s\ne 1$. Therefore, $s\in\IR_X\setminus\{0,1\}$ and hence $\IR_X\ne\{0,1\}$. By Theorem~\ref{t:RXne01=>paraD}, the Playfair plane $X$ is Thalesian.
\end{proof}
\end{proof}

Now we are able to complete the proof of Theorem~\ref{t:uno-Thales<=>}. Let $X$ be an uno-Thalesian Plyfair liner. By Theorem~\ref{t:affine=>Avogadro}, the affine liner $X$ is $2$-balanced. If $|X|_2\le 3$, then $X$ is Desarguesian by Proposition~\ref{p:Steiner+affine=>Desargues}, and Thalesian by Theorem~\ref{t:ADA=>AMA}. So, assume that $|X|_2\ge 4$. By Theorem~\ref{t:4-long-affine}, the $4$-long affine liner $X$ is regular. If $\|X\|\ne 3$, then the affine regular line $X$ is Desarguesian by Theorem~\ref{t:proaffine-Desarguesian}, and Thalesian, by Theorem~\ref{t:ADA=>AMA}. If $\|X\|=3$, then $X$ is a Playfair plane. In this case $X$ is Boolean or Thalesian, by Lemma~\ref{l:uno-Thalesian=>BorT}.
\end{proof}

\begin{corollary}\label{c:Thales=nT+B} An affine regular liner $X$ is uno-Thalesian if and only if it is Thalesian or Boolean.
\end{corollary}

\begin{proof} To prove the ``if'' part, assume that $X$ is Thalesian or Boolean. To check that $X$ is uno-Thalesian, take any distinct parallel lines $A,B,C\subseteq X$ and any distinct points $a,a'\in A$, $b,b'\in B$, $c,c'\in C$ such that  $\Aline ab\cap\Aline{a'}{b'}=\varnothing=\Aline bc\cap\Aline{b'}{c'}$ and $b'\in\Aline ac$. By Theorem~\ref{t:affine=>Avogadro}, the affine liner $X$ is $2$-balanced. Since $X$ contains distinct parallel lines, $|X|_2\ge 3$ and hence $X$ is $3$-long. By Theorem~\ref{t:Playfair<=>}, the $3$-long affine regular liner $X$ is Playfair. By Theorem~\ref{t:uno-Thales<=>}, the (Thalesian or Boolean) Playfair liner $X$ is uno-Thalesian. 
\smallskip

To prove the ``only if'' part, assume that an affine regular liner $X$ is uno-Thalesian. We need to check that $X$ is Thalesian or Boolean. If $X$ is not Thalesian, then there exist three distinct parallel lines $A,B,C\subseteq X$ and any distinct points $a,a'\in A$, $b,b'\in B$, $c,c'\in C$ such that $\Aline ab\cap\Aline{a'}{b'}=\Aline bc\cap\Aline{b'}{c'}=\varnothing\ne\Aline ac\cap\Aline {a'}{c'}$. Since the affine liner $X$ contains distinct parallel lines, $|X|_2\ge 3$ and hence $X$ is $3$-long. By Theorem~\ref{t:Playfair<=>}, the $3$-long affine regular liner $X$ is Playfair. By Theorem~\ref{t:uno-Thales<=>}, the non-Thalesian uno-Thalesian Playfair liner $X$ is Boolean.
\end{proof}

\begin{problem} Is every Boolean Playfair plane Thalesian?
\end{problem}


\section{Bi-Thalesian and by-Thalesian liners}\label{s:bi-Thalesian}

In this section we study two weaker versions of Thales Axiom in which two vertices of one triangle lie on sides of the other triangle.

\begin{definition}\label{d:bi-Thalesian} A liner $X$ is defined to be \defterm{bi-Thalesian} if for every distinct parallel lines $A,B,C$ and points $a,a'\in A$, $b,b'\in B$, $c,c'\in C$,
$$\big(\,\Aline ab\cap\Aline {a'}{b'}=\varnothing=\Aline bc\cap\Aline{b'}{c'}\;\wedge\; a'\in \Aline bc\;\wedge\; c'\in \Aline ba\,\big)\;\Ra\;(\Aline ac\cap\Aline{a'}{c'}=\varnothing).$$

\begin{picture}(100,80)(-150,-10)
\put(30,30){\color{blue}\line(1,1){30}}
\put(30,30){\color{violet}\line(1,-1){30}}
\put(0,0){\color{teal}\line(1,0){60}}
\put(0,60){\color{teal}\line(1,0){60}}
\put(30,30){\color{teal}\line(1,0){60}}

{\linethickness{1pt}
\put(0,0){\color{blue}\line(1,1){30}}
\put(0,60){\color{violet}\line(1,-1){30}}
\put(0,0){\color{red}\line(0,1){60}}
\put(60,0){\color{red}\line(0,1){60}}
\put(60,0){\color{blue}\line(1,1){30}}
\put(60,60){\color{violet}\line(1,-1){30}}
}

\put(0,0){\circle*{3}}
\put(-3,-8){$c$}
\put(60,0){\circle*{3}}
\put(59,-8){$c'$}
\put(0,60){\circle*{3}}
\put(-3,63){$a$}
\put(60,60){\circle*{3}}
\put(59,63){$a'$}
\put(30,30){\circle*{3}}
\put(21,27){$b$}
\put(90,30){\circle*{3}}
\put(93,27){$b'$}
\end{picture}
\end{definition}

\begin{proposition}\label{p:bi-Thalesian<=>} A Playfair liner $X$ is bi-Thalesian if and only if for any triangles $abc$ and $a'b'c'$ in $X$,
$$\big(\,\Aline ab\cap\Aline{a'}{b'}=\Aline bc\cap\Aline{b'}{c'}=\Aline ac\cap\Aline {a'}{c'}=\Aline a{a'}\cap\Aline {b}{b'}=\varnothing\;\wedge\;a'\in\Aline bc\;\wedge\;c'\in \Aline ab\,\big)\;\Ra\;(\Aline b{b'}\cap \Aline c{c'}=\varnothing).$$

\begin{picture}(100,80)(-150,-10)
\put(30,30){\color{blue}\line(1,1){30}}
\put(30,30){\color{violet}\line(1,-1){30}}

{\linethickness{1pt}
\put(60,0){\color{cyan}
\line(0,1){60}}\put(0,0){\color{red}\line(1,0){60}}
\put(0,60){\color{red}\line(1,0){60}}
\put(30,30){\color{red}\line(1,0){60}}
\put(0,0){\color{blue}\line(1,1){30}}
\put(0,60){\color{violet}\line(1,-1){30}}
\put(0,0){\color{cyan}\line(0,1){60}}

\put(60,0){\color{blue}\line(1,1){30}}
\put(60,60){\color{violet}\line(1,-1){30}}
}

\put(0,0){\circle*{3}}
\put(-3,-8){$c$}
\put(60,0){\circle*{3}}
\put(59,-8){$c'$}
\put(0,60){\circle*{3}}
\put(-3,63){$a$}
\put(60,60){\circle*{3}}
\put(59,63){$a'$}
\put(30,30){\circle*{3}}
\put(21,27){$b$}
\put(90,30){\circle*{3}}
\put(93,27){$b'$}
\end{picture}
\end{proposition}

\begin{proof} Assume that $X$ is bi-Thalesian. Given any triangles $abc$ and $a'b'c'$ in $X$ with $\Aline ab\cap\Aline{a'}{b'}=\Aline bc\cap\Aline{b'}{c'}=\Aline ac\cap\Aline {a'}{c'}=\Aline a{a'}\cap\Aline {b}{b'}=\varnothing$ and $a'\in\Aline bc$ and $c'\in \Aline ab$, we have to prove that $\Aline b{b'}\cap \Aline c{c'}=\varnothing$. Consider the plane $\Pi\defeq\overline{\{a,b,c\}}$ and observe that $a',c'\in \Aline bc\cup\Aline ab\subseteq\Pi$. Assuming that $b\in \Aline c{c'}$, we conclude that $a'\in \Aline bc\subseteq\Aline c{c'}$. It follows from $\Aline bc\cap\Aline{b'}{c'}=\varnothing$ that $b\ne c'$. Then $c'\in\Aline ab$ implies $a\in \Aline b{c'}\subseteq \Aline c{c'}$. Therefore, $a,b,c\in\Aline c{c'}$ and hence $abc$ is not a triangle, which is a contradiction showing that $b\notin\Aline c{c'}$. If $b'\notin\Pi$, then $\Aline b{b'}\cap \Aline c{c'}=\Aline b{b'}\cap(\Pi\cap\Aline c{c'})=(\Aline b{b'}\cap\Pi)\cap\Aline c{c'}=\{b\}\cap\Aline c{c'}=\varnothing$ and we are done. So, assume that $b'\in \Pi$. It follows from $\Aline a{a'}\cap\Aline b{b'}=\varnothing=\Aline ab\cap\Aline {a'}{b'}$ that $abb'a'$ is a parallelogram in $\Pi$, and $\Aline a{a'},\Aline b{b'}$ are two distinct parallel lines in the plane $\Pi$.
Since $X$ is Playfair, there exists a unique line $C\subseteq\Pi$ such that $c'\in C$ and $C\parallel \Aline a{a'}\parallel \Aline b{b'}$. Assuming that $C=\Aline a{a'}$, we conclude that $c'\in \Aline a{c'}\cap \Aline a{a'}\subseteq \Aline a{b}\cap \Aline a{a'}=\{a\}$, which contradicts $\Aline ac\cap\Aline {a'}{c'}=\varnothing$.
 Assuming that $B=\Aline b{b'}$, we conclude that $c'\in \Aline ab\cap\Aline {b}{b'}=\{b\}$, which contradicts $\Aline bc\cap\Aline{b'}{c'}=\varnothing$. These contradictions show that $\Aline a{a'},\Aline b{b'}$ and $C$ are three distinct parallel lines in the plane $\Pi$. Since $\Aline {a'}b\cap\Aline b{b'}=\{b\}$ and $C\parallel \Aline b{b'}$, there exists a unique point $\gamma\in C\cap\Aline {a'}b$. Applying the bi-Thalesian property of $X$ to the triangles $ab\gamma$ and $a'b'c'$, we conclude that $\Aline a\gamma\parallel\Aline {a'}{c'}\parallel \Aline ac$ and hence $\gamma\in \Aline {a'}\gamma\cap \Aline a\gamma=\Aline {a'}c\cap\Aline ac=\{c\}$ and finally, $\Aline c{c'}\cap\Aline b{b'}=\Aline \gamma{c'}\cap\Aline b{b'}=C\cap\Aline b{b'}=\varnothing$.
\smallskip

Now assume that a Playfair liner satisfies the ``if'' part of Proposition~\ref{p:bi-Thalesian<=>}. To prove that $X$ is bi-Thalesian, take any distinct parallel lines $A,B,C$ in $X$ and points $a,a'\in A$, $b,b'\in B$, $c,c'\in C$ such that  $a'\in \Aline bc$, $c'\in\Aline ab$, and $\Aline ab\cap\Aline {a'}{b'}=\varnothing=\Aline bc\cap\Aline{b'}{c'}$. We have to prove that $\Aline ac\cap\Aline{a'}{c'}=\varnothing$. It follows from $A\parallel B\parallel C$ and $b'\in\Aline ac$ that $\Pi\defeq\overline{A\cup B\cup C}$ is a plane. Since $X$ is Playfair, there exists a unique line $L$ in $\Pi$ such that $a\in L$ and $L\parallel\Aline{a'}{c'}$. 
Assuming that $b\in \Aline {a'}{c'}$ and taking into account that $a'\in\Aline bc$, we conclude that $c\in\Aline {a'}b\subseteq\Aline {a'}{c'}$ and hence $c\in \Aline {a'}{c'}\cap C=\{c'\}$, which contradicts $\Aline bc\cap\Aline{b'}{c'}=\varnothing$. This contradiction shows that $b\notin\Aline {a'}{c'}$ and hence the lines $\Aline {a'}{b}$ and $\Aline {a'}{c'}$ are concurrent. Then the lines $\Aline {a'}b$ and $L$ are concurrent, too. Hence $\Aline {a'}b\cap L=\{\gamma\}$ for some point $\gamma\notin\{a,b\}$.
It follows that $ab\gamma$ and $a'b'c'$ are two triangles in the plane $\Pi$ such that $a'\in\Aline bc=\Aline b\gamma$, $c'\in\Aline ab$ and $\Aline ab\cap\Aline {a'}{b'}=\Aline b\gamma\cap\Aline{b'}{c'}=\Aline a\gamma\cap\Aline{a'}{c'}=\varnothing$. The ``if'' condition ensures that $\Aline \gamma{c'}\cap \Aline b{b'}=\varnothing$ and hence $\Aline \gamma{c'}\parallel \Aline b{b'}\parallel \Aline{c}{c'}$ and hence $\gamma\in\Aline \gamma{c'}\cap\Aline {a'}b=\Aline c{c'}\cap\Aline cb=\{c\}$ and finally, $\Aline ac\cap\Aline{a'}{c'}=\Aline a\gamma\cap\Aline{a'}{c'}=L\cap\Aline{a'}{c'}=\varnothing$.
\end{proof}

Propositions~\ref{p:uno-Thalesian1} or \ref{p:bi-Thalesian<=>} imply the following corollary.

\begin{corollary} Every uno-Thalesian Playfair liner is bi-Thalesian.
\end{corollary}

\begin{definition}\label{d:by-Thalesian} A liner $X$ is defined to be \defterm{by-Thalesian} if for every distinct parallel lines $A,B,C$ and points $a,a'\in A$, $b,b'\in B$, $c,c'\in C$,
$$\big(\,\Aline ab\cap\Aline {a'}{b'}=\varnothing=\Aline ac\cap\Aline{a'}{c'}\;\wedge\; a'\in \Aline bc\;\wedge\; c'\in \Aline ba\,\big)\;\Ra\;(\Aline bc\cap\Aline{b'}{c'}=\varnothing).$$

\begin{picture}(100,80)(-150,-10)
\put(30,30){\color{red}\line(1,1){30}}
\put(30,30){\color{violet}\line(1,-1){30}}
\put(0,0){\color{teal}\line(1,0){60}}
\put(0,60){\color{teal}\line(1,0){60}}
\put(30,30){\color{teal}\line(1,0){60}}

{\linethickness{1pt}
\put(0,0){\color{red}\line(1,1){30}}
\put(0,60){\color{violet}\line(1,-1){30}}
\put(0,0){\color{cyan}\line(0,1){60}}
\put(60,0){\color{cyan}\line(0,1){60}}
\put(60,0){\color{red}\line(1,1){30}}
\put(60,60){\color{violet}\line(1,-1){30}}
}

\put(0,0){\circle*{3}}
\put(-3,-8){$c$}
\put(60,0){\circle*{3}}
\put(59,-8){$c'$}
\put(0,60){\circle*{3}}
\put(-3,63){$a$}
\put(60,60){\circle*{3}}
\put(59,63){$a'$}
\put(30,30){\circle*{3}}
\put(21,27){$b$}
\put(90,30){\circle*{3}}
\put(93,27){$b'$}
\end{picture}
\end{definition}

\begin{proposition}\label{p:by-Thalesian<=>} A Playfair liner $X$ is by-Thalesian if and only if for any triangles $abc$ and $a'b'c'$ in $X$,
$$\big(\,\Aline ab\cap\Aline{a'}{b'}=\Aline bc\cap\Aline{b'}{c'}=\Aline ac\cap\Aline {a'}{c'}=\Aline a{a'}\cap\Aline {c}{c'}=\varnothing\;\wedge\;a'\in\Aline bc\;\wedge\;c'\in \Aline ab\,\big)\;\Ra\;(\Aline b{b'}\cap \Aline a{a'}=\varnothing).$$

\begin{picture}(100,80)(-150,-10)
\put(30,30){\color{blue}\line(1,1){30}}
\put(30,30){\color{violet}\line(1,-1){30}}

{\linethickness{1pt}
\put(60,0){\color{cyan}
\line(0,1){60}}\put(0,0){\color{red}\line(1,0){60}}
\put(0,60){\color{red}\line(1,0){60}}
\put(30,30){\color{red}\line(1,0){60}}
\put(0,0){\color{blue}\line(1,1){30}}
\put(0,60){\color{violet}\line(1,-1){30}}
\put(0,0){\color{cyan}\line(0,1){60}}

\put(60,0){\color{blue}\line(1,1){30}}
\put(60,60){\color{violet}\line(1,-1){30}}
}

\put(0,0){\circle*{3}}
\put(-3,-8){$c$}
\put(60,0){\circle*{3}}
\put(59,-8){$c'$}
\put(0,60){\circle*{3}}
\put(-3,63){$a$}
\put(60,60){\circle*{3}}
\put(59,63){$a'$}
\put(30,30){\circle*{3}}
\put(21,27){$b$}
\put(90,30){\circle*{3}}
\put(93,27){$b'$}
\end{picture}
\end{proposition}

\begin{proof} Assume that $X$ is by-Thalesian. Given any triangles $abc$ and $a'b'c'$ in $X$ with $\Aline ab\cap\Aline{a'}{b'}=\Aline bc\cap\Aline{b'}{c'}=\Aline ac\cap\Aline {a'}{c'}=\Aline a{a'}\cap\Aline {c}{c'}=\varnothing$ and $a'\in\Aline bc$ and $c'\in \Aline ab$, we have to prove that $\Aline a{a'}\cap \Aline b{b'}=\varnothing$. Consider the plane $\Pi\defeq\overline{\{a,b,c\}}$ and observe that $a',c'\in \Aline bc\cup\Aline ab\subseteq\Pi$. Assuming that $b\in \Aline a{a'}$, we conclude that $c'\in \Aline ab\subseteq\Aline a{a'}$. It follows from $\Aline ab\cap\Aline{a'}{b'}=\varnothing$ that $a'\ne b$. Then $a'\in\Aline bc$ implies $c\in \Aline b{a'}\subseteq \Aline a{a'}$. Therefore, $a,b,c\in\Aline a{a'}$ and hence $abc$ is not a triangle, which is a contradiction showing that $b\notin\Aline a{a'}$. If $b'\notin\Pi$, then $\Aline b{b'}\cap \Aline a{a'}=\Aline b{b'}\cap(\Pi\cap\Aline a{a'})=(\Aline b{b'}\cap\Pi)\cap\Aline a{a'}=\{b\}\cap\Aline a{a'}=\varnothing$ and we are done. So, assume that $b'\in \Pi$. It follows from $\Aline a{a'}\cap\Aline c{c'}=\varnothing=\Aline ac\cap\Aline {a'}{c'}$ that $acc'a'$ is a parallelogram in $\Pi$, and $\Aline a{a'},\Aline c{c'}$ are two distinct parallel lines in the plane $\Pi$.
Since $X$ is Playfair, there exists a unique line $B\subseteq\Pi$ such that $b\in B$ and $B\parallel \Aline a{a'}\parallel \Aline c{c'}$. Assuming that $B=\Aline a{a'}$, we conclude that $c'\in \Aline ab\subseteq \Aline a{a'}$, which contradicts $\Aline a{a'}\cap\Aline {c}{c'}=\varnothing$.
 Assuming that $B=\Aline c{c'}$, we conclude that $a'\in \Aline bc\subseteq\Aline c{c'}$, which contradicts $\Aline a{a'}\cap\Aline c{c'}=\varnothing$. These contradictions show that $\Aline a{a'},B$ and $\Aline c{c'}$ are three distinct parallel lines in the plane $\Pi$. Since $\Aline ab\cap B=\{b\}$ and $\Aline ab\parallel \Aline {a'}{b'}$, there exists a unique point $\beta'\in B\cap\Aline {a'}{b'}$. Applying the by-Thalesian property of $X$ to the triangles $abc$ and $a'\beta'c'$, we conclude that $\Aline bc\cap\Aline {\beta'}{c'}=\varnothing$ and hence $\Aline{b'}{c'}\parallel\Aline bc\parallel\Aline{\beta'}{c'}$ and $\Aline {b'}{c'}=\Aline{\beta'}{c'}$. Then $\beta'\in \Aline {a'}{\beta'}\cap\Aline {\beta'}{c'}=\Aline{a'}{b'}\cap\Aline{b'}{c'}=\{b'\}$ and finally, $\Aline b{b'}\cap \Aline a{a'}=\Aline b{\beta'}\cap\Aline a{a'}\subseteq B\cap \Aline a{a'}=\varnothing$.
\smallskip

Now assume that a Playfair liner satisfies the ``if'' part of Proposition~\ref{p:by-Thalesian<=>}. To prove that $X$ is by-Thalesian, take any distinct parallel lines $A,B,C$ in $X$ and points $a,a'\in A$, $b,b'\in B$, $c,c'\in C$ such that  $a'\in \Aline bc$, $c'\in\Aline ab$, and $\Aline ab\cap\Aline {a'}{b'}=\varnothing=\Aline ac\cap\Aline{a'}{c'}$. We have to prove that $\Aline bc\cap\Aline{b'}{c'}=\varnothing$. It follows from $A\parallel B\parallel C$ and $b'\in\Aline ac$ that $\Pi\defeq\overline{A\cup B\cup C}$ is a plane. Since $X$ is Playfair, there exists a unique line $L$ in $\Pi$ such that $c'\in L$ and $L\parallel\Aline bc$. It follows from $c'\in \Aline ab$ and $\Aline ac\cap\Aline {a'}{c'}=\varnothing$ that $\Aline ab\cap\Aline bc=\{b\}$. Then  $L\cap\Aline {a'}{b'}=\{\beta'\}$ for some point $\beta'$. 
It follows that $abc$ and $a'\beta'c'$ are two triangles in the plane $\Pi$ such that $a'\in\Aline bc$, $c'\in\Aline ab$ and $\Aline ab\cap\Aline {a'}{\beta'}=\Aline bc\cap\Aline{\beta'}{c'}=\Aline a\gamma\cap\Aline{a'}{c'}=\varnothing$. The ``if'' condition ansures that $\Aline a{a'}\cap \Aline b{\beta'}=\varnothing$ and hence $\Aline b{\beta'}\parallel A\parallel B=\Aline b{b'}$ and hence $\beta'\in\Aline{a'}{\beta'}\cap \Aline b{\beta'}=\Aline{a'}{b'}\cap\Aline b{b'}=\{b'\}$ and finally, $\Aline bc\cap\Aline{b'}{c'}=\Aline bc\cap\Aline{\beta'}{c'}=\Aline bc\cap L=\varnothing$.
\end{proof}

Definitions~\ref{d:uno-Thalesian} and \ref{d:by-Thalesian} imply the following proposition.

\begin{proposition} Every uno-Thalesian liner is by-Thalesian.
\end{proposition}

\begin{proposition}\label{p:inv-pls<=>inv-puls+biT} An invertible-puls Playfair liner is bi-Thalesian if and only if it is invertible-plus.
\end{proposition}

\begin{proof} First assume that an invertible puls Playfair liner $X$ is invertible-plus. To prove that $X$ is bi-Thalesian, we shall apply Proposition~\ref{p:bi-Thalesian<=>}. Given any triangles $abc$ and $a'b'c'$ in $X$ with $a'\in\Aline bc$, $c'\in \Aline ab$ and $\Aline ab\cap\Aline {a'}{b'}=\Aline bc\cap\Aline{b'}{c'}=\Aline ac\cap\Aline{a'}{c'}=\Aline a{a'}\cap\Aline b{b'}=\varnothing$, we have to prove that $\Aline b{b'}\cap\Aline c{c'}=\varnothing$. Consider the plane $\Pi\defeq\overline{\{a,b,c\}}$ and observe that $a',c'\in \Aline bc\cup\Aline ab\subseteq\Pi$. Assuming that $b\in \Aline c{c'}$, we conclude that $a'\in \Aline bc\subseteq\Aline c{c'}$. It follows from $\Aline bc\cap\Aline{b'}{c'}=\varnothing$ that $b\ne c'$. Then $c'\in\Aline ab$ implies $a\in \Aline b{c'}\subseteq \Aline c{c'}$. Therefore, $a,b,c\in\Aline c{c'}$ and hence $abc$ is not a triangle, which is a contradiction showing that $b\notin\Aline c{c'}$. If $b'\notin\Pi$, then $\Aline b{b'}\cap \Aline c{c'}=\Aline b{b'}\cap(\Pi\cap\Aline c{c'})=(\Aline b{b'}\cap\Pi)\cap\Aline c{c'}=\{b\}\cap\Aline c{c'}=\varnothing$ and we are done. So, assume that $b'\in \Pi$.  
By Theorem~\ref{t:parallelogram3+1}, there exist points $\alpha,\beta,\gamma\in X$ such that $\alpha aba'$, $\beta abc$ and $\gamma cbc'$ are parallelograms.

\begin{picture}(100,110)(-200,-55)
{\linethickness{1pt}
\put(0,0){\line(-1,-1){20}}
\put(0,0){\line(-1,1){20}}
\put(0,0){\line(1,0){40}}
\put(-20,20){\line(0,-1){40}}
\put(20,20){\line(0,-1){40}}
\multiput(-20,-20)(3,0){13}{\line(1,0){2}}
\put(-20,20){\line(1,0){40}}
\put(20,20){\line(1,-1){20}}
\put(20,-20){\line(1,1){20}}
}
\put(0,0){\line(1,1){20}}
\put(0,0){\line(1,-1){20}}

\put(0,0){\line(-1,0){40}}
\put(0,0){\line(0,1){40}}
\put(0,0){\line(0,-1){40}}
\put(-40,0){\line(1,1){40}}
\put(-40,0){\line(1,-1){40}}
\put(40,0){\line(-1,1){40}}
\put(40,0){\line(-1,-1){40}}

\put(0,0){\circle*{3}}
\put(1,5){$b$}
\put(40,0){\circle*{3}}
\put(43,-3){$b'$}
\put(20,20){\circle*{3}}
\put(22,22){$a'$}
\put(0,40){\circle*{3}}
\put(-3,43){$\alpha$}
\put(-20,20){\circle*{3}}
\put(-27,22){$a$}
\put(-40,0){\circle*{3}}
\put(-49,-3){$\beta$}
\put(-20,-20){\circle*{3}}
\put(-28,-26){$c$}
\put(0,-40){\circle*{3}}
\put(-3,-48){$\gamma$}
\put(20,-20){\circle*{3}}
\put(22,-26){$c'$}
\end{picture}

Since $X$ is invertible-puls, $\Aline b{b'}\parallel \Aline a{a'}$ implies $\Aline \alpha b\parallel \Aline {a'}{c'}\parallel \Aline ac$, which implies $\Aline b\beta\parallel \Aline a{a'}$. Since $\Aline b{b'}\parallel \Aline a{a'}\parallel \Aline \beta b$, the invertibility-plus of $X$ implies $\Aline c{c'}\parallel \Aline b{b'}$. Assuming that $\Aline c{c'}\cap\Aline b{b'}\ne\varnothing$, we conclude that $\Aline c{c'}=\Aline b{b'}$ and hence $a'\in \Aline cb\subseteq \Aline b{b'}$, which contradicts $\Aline a{a'}\cap\Aline b{b'}=\varnothing$. This contradiction shows that $\Aline c{c'}\cap\Aline b{b'}=\varnothing$.
\smallskip

Now assume that an invertible puls Playfair liner $X$ is bi-Thalesian. To show that $X$ is invertible-plus, take any points $a,b,c,\beta,a',b',c'\in X$ such that 
$\Aline {\beta}c\parallel \Aline ab=\Aline b{c'}\parallel \Aline {a'}{b'}\nparallel \Aline \beta a\parallel\Aline cb=\Aline b{a'}\parallel \Aline{c'}{b'}$ and $\Aline {a}{a'}\parallel \Aline b{b'}\parallel \Aline \beta b$. We have to prove that $\Aline c{c'}\parallel \Aline b{b'}$. If $\beta=c$, then
$\Aline {\beta}c\parallel \Aline ab=\Aline b{c'}$ implies $a=b=c'$. Then $\Aline c{c'}=\Aline\beta b\parallel \Aline b{b'}$. If $a=\beta$, then  $ \Aline \beta a\parallel\Aline cb=\Aline b{a'}$ implies $c=b=a'$. In this case $\Aline c{c'}=\Aline b{b'}$. So, assume that $a\ne\beta\ne c$. By Theorem~\ref{t:parallelogram3+1}, there exists a unique point $\alpha\in X$ such that $\alpha aba'$ is a parallelogram.

\begin{picture}(100,90)(-160,-15)
\put(0,30){\line(0,1){30}}
\put(0,60){\line(1,0){30}}

{\linethickness{1pt}
\put(0,0){\line(1,0){30}}
\put(0,0){\line(0,1){30}}
\put(30,0){\line(0,1){60}}
\put(60,30){\line(0,1){30}}
\put(0,30){\line(1,0){60}}
\put(30,60){\line(1,0){30}}
\put(0,0){\line(1,1){60}}
\multiput(30,0)(3,3){10}{\line(1,1){2}}
\put(0,30){\line(1,1){30}}
}
\put(0,60){\line(1,-1){30}}

\put(30,0){\line(-1,1){30}}
\put(0,30){\line(1,-1){30}}
\put(30,60){\line(1,-1){30}}

\put(0,0){\circle*{3}}
\put(-8,-8){$\beta$}
\put(30,0){\circle*{3}}
\put(28,-8){$c$}
\put(0,30){\circle*{3}}
\put(-9,27){$a$}
\put(30,30){\circle*{3}}
\put(31,35){$b$}
\put(60,30){\circle*{3}}
\put(63,28){$c'$}
\put(0,60){\circle*{3}}
\put(-9,60){$\alpha$}
\put(30,60){\circle*{3}}
\put(28,63){$a'$}
\put(60,60){\circle*{3}}
\put(60,63){$b'$}
\end{picture}

 Since $X$ is invertible-puls, the parallelity relations $\Aline b{b'}\parallel\Aline a{a'}\parallel \Aline\beta b$ imply $\Aline {a'}{c'}\parallel \Aline \alpha b\parallel \Aline ac$. By Proposition~\ref{p:bi-Thalesian<=>}, $a'\in \Aline bc$, $c'\in \Aline ba$, $\Aline ab\cap\Aline{a'}{b'}=\Aline bc\cap\Aline{b'}{c'}=\Aline ac\cap \Aline{a'}{c'}=\Aline a{a'}\cap\Aline b{b'}=\varnothing$ imply $\Aline c{c'}\cap\Aline b{b'}=\varnothing$ and hence $\Aline c{c'}\parallel \Aline b{b'}$, witnessing that $X$ is invertible-plus.
\end{proof}

A liner is called \defterm{invertible-add} if it is invertible-plus and invertible-puls. Proposition~\ref{p:inv-pls<=>inv-puls+biT} implies

\begin{corollary}\label{c:invertible-add=>bi-Thalesian} Every invertible-add Playfair liner is bi-Thalesian.
\end{corollary}

\begin{question} Is every invertible-add Playfair liner by-Thalesian?
\end{question}

\begin{exercise} Show that a bi-Thalesian by-Thalesian Playfair liner $X$ is  invertible-plus if and only if it is invertible-puls.
\end{exercise}

\begin{proposition}\label{p:by-Thalesian<=>by-Boolean} Let $X$ be a Playfair liner.
\begin{enumerate}
\item If $X$ is by-Thalesian, then $X$ is by-Boolean;
\item If $X$ is bi-Thalesian and by-Boolean, then $X$ is by-Thalesian.
\end{enumerate}
\end{proposition}

\begin{proof} 1. Assuming that a Playfair liner $X$ is by-Thalesian. Assuming that $X$ is not by-Boolean, we can two lines $L,L'\subseteq X$ and distinct points $a,b,c\in L$ and $a',b',c'\in L'$ such that $\Aline a{b'}\cap \Aline {a'}b=\Aline b{c'}\cap\Aline {b'}c=\Aline a{a'}\cap\Aline b{b'}=\Aline b{b'}\cap\Aline c{c'}=\varnothing\ne L\cap L'$. Assuming that $L=L'$ and taking into account that $\Aline a{b'}\cap \Aline {a'}b=\Aline a{a'}\cap\Aline b{b'}$, we conclude that $a=b'\ne a'=b=b'$, which is a contradiction showing that $L\ne L'$. Then $L\cap L=\{\beta\}$ for some point $\beta$. We claim that $\beta\notin\{a,a',b,b',c,c'\}$.
Indeed, if $\beta=a$, then $a\ne a'$ or $a\ne b'$ and hence $\Aline a{a'}=L'$ or $\Aline a{b'}=L'$, which contradicts $\Aline a{a'}\cap\Aline b{b'}=\varnothing=\Aline a{b'}\cap\Aline {a'}b$. By analogy we can show that $\beta\notin\{a',b,b',c,c'\}$. By Theorem~\ref{t:parallelogram3+1}, the triangle $b\beta b'$ can be completed to a parallelogram $b\beta b'\beta'$. Applying Proposition~\ref{p:by-Thalesian<=>} to the triangles $a\beta a'$ and $b'\beta'b$, we conclude that $\Aline\beta{\beta'}\cap\Aline a{b'}=\varnothing$. Applying Proposition~\ref{p:by-Thalesian<=>} to the triangles $c\beta c'$ and $b'\beta'b$, we conclude that $\Aline\beta{\beta'}\cap\Aline c{b'}=\varnothing$. Then $\Aline a{b'}\parallel \Aline\beta{\beta'}\parallel\Aline c{b'}$ and hence $\Aline a{b'}=\Aline b{b'}$ and $a\in \Aline a{b'}\cap L=\Aline c{b'}\cap L=\{c\}$, which contradicts the choice of the distinct points $a,c$. This contradiction shows that the by-Thalesian Playfair liner $X$ is by-Boolean.

\begin{picture}(100,100)(-180,-50)
\linethickness{0.8pt}
\put(0,0){\color{violet}\line(1,-1){20}}
\put(0,0){\color{blue}\line(1,1){20}}
\put(0,0){\color{violet}\line(-1,1){40}}
\put(0,0){\color{blue}\line(-1,-1){40}}
\put(20,20){\color{teal}\line(-1,0){40}}
\put(20,20){\color{cyan}\line(0,-1){40}}
\put(20,-20){\color{teal}\line(-1,0){40}}
\put(-20,-20){\color{cyan}\line(0,1){40}}
\put(-40,-40){\color{cyan}\line(0,1){80}}
\put(20,-20){\color{red}\line(-3,-1){60}}
\put(20,20){\color{red}\line(-3,1){60}}
\put(20,-20){\color{blue}\line(1,1){20}}
\put(20,20){\color{violet}\line(1,-1){20}}
\put(0,0){\color{teal}\line(1,0){40}}

\put(-40,-40){\circle*{3}}
\put(-49,-42){$a'$}
\put(-20,-20){\circle*{3}}
\put(-28,-21){$c'$}
\put(20,-20){\circle*{3}}
\put(22,-27){$b$}
\put(0,0){\circle*{3}}
\put(-11,-4){$\beta$}
\put(-20,20){\circle*{3}}
\put(-28,17){$c$}
\put(-40,40){\circle*{3}}
\put(-48,38){$a$}
\put(20,20){\circle*{3}}
\put(21,21){$b'$}
\put(40,0){\circle*{3}}
\put(43,-3){$\beta'$}
\end{picture}

Now assume that a Playfair liner $X$ is bi-Thalesian and by-Boolean. To show that $X$ by-Thalesian, we shall apply Proposition~\ref{p:by-Thalesian<=>}. Given two triangles $abc$ and $a'b'c'$ in $X$ with $a'\in \Aline bc$, $c'\in \Aline ba$ and $\Aline ab\cap\Aline{a'}{b'}=\Aline bc\cap\Aline{b'}{c'}=\Aline ac\cap\Aline{a'}{c'}=\Aline a{a'}\cap\Aline c{c'}=\varnothing$, we need to show that $\Aline a{a'}\cap\Aline b{b'}=\varnothing$.
To derive a contradiction, assume that $\Aline a{a'}\cap\Aline b{b'}\ne\varnothing$.

 Consider the plane $\Pi\defeq\overline{\{a,b,c\}}$ and observe that $a',c'\in\Aline bc\cup\Aline ba\subseteq \Pi$.  Assuming that $b\in \Aline a{a'}$, we conclude that $c'\in \Aline ab\subseteq\Aline a{a'}$. It follows from $\Aline ab\cap\Aline{a'}{b'}=\varnothing$ that $a'\ne b$. Then $a'\in\Aline bc$ implies $c\in \Aline b{a'}\subseteq \Aline a{a'}$. Therefore, $a,b,c\in\Aline a{a'}$ and hence $abc$ is not a triangle, which is a contradiction showing that $b\notin\Aline a{a'}$. By analogy we can show that $b\notin\Aline c{c'}$. It follows from $b\notin\Aline a{a'}$ and $\Aline a{a'}\cap\Aline b{b'}\ne\varnothing$ that $\Aline b{b'}\subseteq \Pi$. Since $X$ is Playfair, there exist points $\alpha\in \Aline ab$ and $\gamma\in\Aline bc$ such that $\Aline \alpha{a'}\parallel \Aline b{b'}\parallel \Aline \gamma{c'}$. It follows from $\Aline b{c'}\cap\Aline {a'}{b'}=\varnothing=\Aline b{a'}\cap\Aline {c'}{b'}$ that $a',c'\notin\Aline b{b'}$ and hence the parallel lines $\Aline\alpha {a'},\Aline b{b'},\Aline\gamma{c'}$ are distinct. Taking into account that $X$ is bi-Thalesian, we conclude that $\Aline\alpha \gamma\cap\Aline{a'}{c'}=\varnothing$. Then $L\defeq\Aline ab$ and $L'\defeq\Aline bc$ are two lines and $a,c',\alpha\in L$, $c,a',\gamma\in L'$ are distinct points such that $\Aline a{a'}\cap\Aline c{c'}=\Aline {c'}\gamma\cap\Aline {a'}\alpha=\Aline ac\cap\Aline {c'}{a'}=\Aline {c'}{a'}\cap\Aline \alpha\gamma=\varnothing$. Since $X$ is by-Boolean, $L\cap L'=\varnothing$, which contradicts $b\in L\cap L'$. This contradiction shows that $\Aline b{b'}\cap\Aline a{a'}=\varnothing$ and hence the liner $X$ is by-Thalesian.



\end{proof}

\begin{exercise} Show that the Moulton plane is not bi-Thalesian and not by-Thalesian.
\end{exercise}


\section{Quadratic Playfair liners}\label{s:quadratic}

\begin{definition} A liner $X$ is called \defterm{quadratic} if for every points $a,b,c,d\in X$, $o\in \Aline ac\cap\Aline bd$, $h\in \Aline ab\setminus\Aline bc$, $v\in \Aline bc\setminus \Aline ba$,
$$\big(\,\Aline ab\parallel\Aline ov\parallel\Aline cd\;\wedge\;\Aline ad\parallel \Aline oh\parallel\Aline bc\,\big)\;\Ra\;(\Aline hv\parallel \Aline ac).$$

%
\end{definition}

\begin{proposition}\label{p:by-Thalesian=>quadratic} Every by-Thalesian Playfair liner is quadratic.
\end{proposition}

\begin{proof} Assume that a Playfair liner $X$ is by-Thalesian. To prove that $X$ is quadratic, take any points $a,b,c,d\in X$, $o\in\Aline ac\cap\Aline bd$, $h\in \Aline ab\setminus \Aline bc$, $v\in \Aline bc\setminus\Aline ba$ such that $\Aline ab\parallel\Aline ov\parallel\Aline cd$ and  $\Aline ad\parallel \Aline oh\parallel\Aline bc$. We have to prove that $\Aline hv\parallel \Aline ac$. It follows from $h\in \Aline ab\setminus \Aline bc$ and $v\in \Aline bc\setminus\Aline ba$ that $b\notin \{a,h,v,c\}$ and hence $abcd$ is a parallelogram whose centre $o\in\Aline ac\cap\Aline bd$ does not belong to the set $\Aline ab\cup\Aline bc\cup\Aline cd\cup\Aline da$. 
By Theorem~\ref{t:parallelogram3+1}, there exists a unique point $\alpha\in X$ such that $oa\alpha b$ is a parallelogram. 

\begin{picture}(60,60)(-150,-15)
\linethickness{0.7pt}
\put(0,0){\line(1,0){40}}
\put(0,0){\line(0,1){40}}
\put(0,0){\line(1,1){40}}
\put(0,40){\line(1,-1){40}}
\put(0,0){\line(-1,1){20}}
\put(-20,20){\line(1,1){20}}
\put(-20,20){\line(1,0){40}}
\put(20,20){\line(0,1){20}}
\put(0,20){\line(1,1){20}}
\put(0,40){\line(1,0){40}}
\put(40,0){\line(0,1){40}}

\put(0,0){\circle*{3}}
\put(-2,-8){$a$}
\put(40,0){\circle*{3}}
\put(39,-8){$d$}
\put(-20,20){\circle*{3}}
\put(-29,18){$\alpha$}
\put(0,20){\circle*{3}}
\put(-7,22){$h$}
\put(20,20){\circle*{3}}
\put(17,12){$o$}
\put(0,40){\circle*{3}}
\put(-2,43){$b$}
\put(20,40){\circle*{3}}
\put(17,43){$v$}
\put(40,40){\circle*{3}}
\put(38,43){$c$}

\end{picture}

Applying Proposition~\ref{p:by-Thalesian<=>} to the triangles $cod$ and $b\alpha a$, we conclude that $\Aline  o\alpha\parallel \Aline ad\parallel \Aline oh$. Applying Proposition~\ref{p:by-Thalesian<=>}  to the triangles $\alpha ha$ and $bvo$, we conclude that $\Aline  hv\parallel \Aline ao=\Aline ac$.
\end{proof}

\begin{proposition}\label{p:quadratic=>by-Thalesian} Any quadratic invertible-puls Playfair liner $X$ is by-Thalesian.
\end{proposition}

\begin{proof} Assume a quadratic invertible-puls Playfair liner $X$ is invertible-plus or inversive-puls. 
To prove that $X$ is by-Thalesian, we shall apply Proposition~\ref{p:by-Thalesian<=>}. Take any triangles $abc$ and $a'b'c'$ in $X$ such that $a'\in \Aline bc$, $c'\in \Aline ba$  and $\varnothing=\Aline ab\cap\Aline{a'}{b'}=\Aline bc\cap\Aline{b'}{c'}=\Aline ac\cap\Aline {a'}{c'}=\Aline a{a'}\cap\Aline c{c'}$. We have to prove that 
$\Aline a{a'}\cap\Aline{b}{b'}=\varnothing$. Consider the plane $\Pi\defeq\overline{\{a,b,c\}}$ and observe that $a',c'\in \Aline bc\cup \Aline ba\subseteq\Pi$. It follows from $\Aline ac\cap\Aline{a'}{c'}=\varnothing=\Aline a{a'}\cap\Aline c{c'}$ that  $aa'c'c$ is a parallelogram in the plane $\Pi$ and $b\in \Aline a{c'}\cap\Aline c{a'}$. If $b'\notin\Pi$, then $\Aline a{a'}\cap\Aline b{b'}=\Aline a{a'}\cap(\Pi\cap\Aline b{b'})=\Aline a{a'}\cap\{b\}=\varnothing$ and we are done. So, assume that $b'\in\Pi$.

Since $X$ is Playfair, there exist lines $H,V\subseteq \Pi$ such that $b\in H\cap V$, $H\parallel \Aline a{a'}\parallel\Aline c{c'}$ and $V\parallel \Aline ac\parallel \Aline {a'}{c'}$. Consider the unique points $\alpha\in V\cap \Aline a{a'}$, $\gamma\in V\cap\Aline c{c'}$,  $\beta\in H\cap\Aline ac$ and $\beta'\in H\cap\Aline {a'}{c'}$, which exist by the Proclus Postulate~\ref{p:Proclus-Postulate}. 

\begin{picture}(100,110)(-150,-15)
{\linethickness{1pt}
\put(0,0){\color{teal}\line(1,0){80}}
\put(0,80){\color{teal}\line(1,0){80}}
\put(0,0){\color{blue}\line(1,1){40}}
\put(0,80){\color{violet}\line(1,-1){40}}
\put(0,0){\color{cyan}\line(0,1){80}}
\put(80,0){\color{cyan}\line(0,1){80}}
\put(80,0){\color{blue}\line(1,1){40}}
\put(80,80){\color{violet}\line(1,-1){40}}

}
\put(0,0){\color{blue}\line(1,1){80}}
\put(0,40){\color{blue}\line(1,1){40}}
\put(0,40){\color{violet}\line(1,-1){40}}
\put(0,80){\color{violet}\line(1,-1){80}}
\put(0,0){\color{teal}\line(1,0){80}}
\put(0,80){\color{teal}\line(1,0){80}}
\put(0,40){\color{teal}\line(1,0){80}}
\put(40,0){\color{cyan}\line(0,1){80}}
\put(40,0){\color{blue}\line(1,1){40}}
\put(40,80){\color{violet}\line(1,-1){40}}

\put(0,0){\circle*{3}}
\put(-2,-8){$c$}
\put(40,0){\circle*{3}}
\put(37,-8){$\gamma$}
\put(80,0){\circle*{3}}
\put(78,-9){$c'$}
\put(0,40){\circle*{3}}
\put(-9,37){$\beta$}
\put(0,80){\circle*{3}}
\put(-3,83){$a$}
\put(40,80){\circle*{3}}
\put(38,84){$\alpha$}
\put(80,80){\circle*{3}}
\put(78,83){$a'$}
\put(40,40){\circle*{3}}
\put(41,45){$b$}
\put(80,40){\circle*{3}}
\put(82,37){$\beta'$}
\put(120,40){\circle*{3}}
\put(123,41){$b'$}
\end{picture}

 Since $X$ is quadratic, $\Aline \alpha \beta\parallel \Aline c{a'}\parallel \gamma{\beta'}$ and $\Aline \beta\gamma \parallel \Aline a{c'}\parallel \Aline\alpha{\beta'}$. The Proclus Postulate~\ref{p:Proclus-Postulate} ensures the existence of unique points $u\in \Aline \beta\gamma\cap\Aline cb$, $v\in \gamma{\beta'}\cap\Aline b{c'}$, $u'\in \Aline \alpha\beta\cap \Aline ab$, $v'\in \Aline \alpha{\beta'}\cap\Aline b{a'}$. Since $X$ is quadratic, $\Aline uv\parallel H\parallel \Aline{u'}{v'}$ and $\Aline u{u'}\parallel V\parallel \Aline v{v'}$. By the Proclus Postulate~\ref{p:Proclus-Postulate}, there exist unique points $w\in\Aline uv\cap\Aline \alpha{\beta'}$ and $w'\in \Aline {u'}{v'}\cap\Aline \gamma {\beta'}$.

\begin{picture}(100,110)(-150,-15)
{\linethickness{1pt}
\put(0,0){\color{teal}\line(1,0){80}}
\put(0,80){\color{teal}\line(1,0){80}}
\put(0,0){\color{blue}\line(1,1){40}}
\put(0,80){\color{violet}\line(1,-1){40}}
\put(0,0){\color{cyan}\line(0,1){80}}
\put(80,0){\color{cyan}\line(0,1){80}}
\put(80,0){\color{blue}\line(1,1){40}}
\put(80,80){\color{violet}\line(1,-1){40}}

}
\put(0,0){\color{blue}\line(1,1){80}}
\put(0,40){\color{blue}\line(1,1){40}}
\put(0,40){\color{violet}\line(1,-1){40}}
\put(0,80){\color{violet}\line(1,-1){80}}
\put(0,0){\color{teal}\line(1,0){80}}
\put(0,80){\color{teal}\line(1,0){80}}
\put(0,40){\color{teal}\line(1,0){80}}
\put(40,0){\color{cyan}\line(0,1){80}}
\put(40,0){\color{blue}\line(1,1){60}}
\put(40,80){\color{violet}\line(1,-1){60}}
\put(20,20){\color{teal}\line(1,0){80}}
\put(20,60){\color{teal}\line(1,0){80}}
\put(20,20){\color{cyan}\line(0,1){40}}
\put(60,20){\color{cyan}\line(0,1){40}}
\put(100,20){\color{cyan}\line(0,1){40}}

\put(0,0){\circle*{3}}
\put(-2,-8){$c$}
\put(40,0){\circle*{3}}
\put(37,-8){$\gamma$}
\put(80,0){\circle*{3}}
\put(78,-9){$c'$}
\put(0,40){\circle*{3}}
\put(-9,37){$\beta$}
\put(0,80){\circle*{3}}
\put(-3,83){$a$}
\put(40,80){\circle*{3}}
\put(38,84){$\alpha$}
\put(80,80){\circle*{3}}
\put(78,83){$a'$}
\put(40,40){\circle*{3}}
\put(41,45){$b$}
\put(80,40){\circle*{3}}
\put(83,37){$\beta'$}
\put(120,40){\circle*{3}}
\put(123,41){$b'$}

\put(20,20){\circle*{2.5}}
\put(17,12){$u$}
\put(60,20){\circle*{2.5}}
\put(57,12){$v$}
\put(100,20){\circle*{2.5}}
\put(102,15){$w$}
\put(20,60){\circle*{2.5}}
\put(17,63){$u'$}
\put(60,60){\circle*{2.5}}
\put(57,63){$v'$}
\put(100,60){\circle*{2.5}}
\put(102,60){$w'$}
\end{picture}

Consider the parallel lines $\overline{\{u,v,w\}}$ and $\overline{\{u',v',w'\}}$. 
Since $X$ is invertible-puls, the parallelity relations $\Aline u{v'}\parallel \Aline v{w'}$, $\Aline {u'}v\parallel\Aline {v'}w$ and $\Aline u{u'}\parallel \Aline v{v'}$ imply $\Aline w{w'}\parallel \Aline v{v'}\parallel \Aline u{u'}$. 

Next, consider the parallel lines $\overline{\{b,v,c'\}}$ and $\overline{\{v',\beta',w\}}$. Since $X$ is invertible-puls, the parallelity relations $\Aline b{\beta'}\parallel \Aline vw$, $\Aline v{v'}\parallel \Aline {c'}{\beta'}$ and $\Aline b{v'}\parallel\Aline v{\beta'}$ imply $\Aline {c'}w\parallel \Aline b{v'}=\Aline cb\parallel\Aline {c'}{b'}$ and hence $w\in \Aline {c'}{b'}$. By analogy we can prove that $w'\in \Aline{a'}{b'}$.

Finally, consider the parallel lines $\overline{\{c',w,b'\}}$ and $\overline{\{v,\beta',w'\}}$. Since $X$ is invertible-puls, the parallelity relations $\Aline{c'}{\beta'}\parallel \Aline w{w'}$ and $\Aline v{c'}\parallel \Aline w{\beta'}\parallel \Aline{b'}{w'}$ imply $\Aline {\beta'}{b'}\parallel \Aline vw\parallel H=\Aline b{\beta'}$ and hence $\Aline b{b'}=H\parallel \Aline a{a'}\parallel b{b'}$.
\end{proof}

Quadratic Playfair liners admit the following algebraic characterization.
  
\begin{proposition}\label{p:quadratic<=>} For a Playfair liner $X$, the following conditions are equivalent:
\begin{enumerate}
\item $X$ is quadratic;
\item for every ternar $R$ of $X$ and $\forall x\in R\setminus\{1\}\;\;(x\cdot x=1\;\Ra\;1+x=0)$;
\item for every ternar $R$ of $X$ and $\forall x\in R\setminus\{1\}\;\;(x\cdot x=1\;\Ra\;x+1=0)$;
\item for every ternar $R$ of $X$ and $\forall x\in R\setminus\{1\}\;\;(x\cdot x=1\;\Ra\;x\puls 1=0)$.
\end{enumerate}
\end{proposition}

\begin{proof} \vbox{\hsize335pt $(1)\Ra(2\wedge 3\wedge 4)$ Assume that a Playfair liner $X$ is quadratic. Given any ternar $R$ of $X$, find a plane $\Pi$ and affine base $uow$ in $\Pi$ whose ternar $\Delta$ is isomorphic to the ternar $R$. We need to prove that for every $x\in\Delta\setminus\{e\}$, $x\cdot x=e$ implies $e+x=x+e=x\puls e=o$, where $e$ is the diunit of the affine base $uow$. Consider the points $b,d\in\Pi$ with coordinates $xe$ and $ex$, respectively. Then $xbed$ is a parallelogram in the plane $\Pi$ such that $w\in \Aline be$ and $u\in \Aline ed$. Consider also the unique points $h\in \Aline ou\cap\Aline xb$ and $s\in \Aline ow\cap \Aline xd$.  By definition of the dot operation, $x\cdot x=e$ implies $b=\Aline od$ and hence $o\in \Aline xe\cap\Aline bd$. Since $X$ is quadratic, $\Aline hw\parallel \Aline oe\parallel \Aline su$ and $\Aline wu\parallel \Aline od$. By definition of the plus and puls operations, the latter parallelity relations imply $e+x=o=e+x$ and $x\puls e=o$.}


$(2\vee 3\vee 4)\Ra(1)$ Assume that $X$ satisfies one of the conditions (2), (3), (4). To prove that $X$ is quadratic, take any points $a,b,c,d\in X$, $o\in\Aline ac\cap\Aline bd$, $h\in \Aline ab\setminus\Aline bc$, $v\in \Aline bc\setminus\Aline ba$ such that $\Aline ab\parallel \Aline ov\parallel \Aline cd$ and $\Aline ad\parallel \Aline oh\parallel \Aline bc$. We have to prove that $\Aline hv\parallel \Aline ac$. 
\vskip-5pt

%

\vbox{\hsize350pt If $X$ satisfies condition (2), then consider the affine base $uow$ with $u\in \Aline ho\cap\Aline cd$ and $w=v$.  The diunit $e$ of the affine base $uow$ coincides with the point $c$. By definition of the dot operation in the ternar $\Delta=\Aline oe=\Aline oc$, the inclusion $b\in\Aline od$ implies $a\cdot a=c=e$. The condition (2) implies $e+a=o$. By definition of the plus operation, the equality $e+a=o$ implies $\Aline hv=\Aline hw\parallel \Delta=\Aline ac$.} 
\vskip-5pt

%

\vbox{\hsize350pt If $X$ satisfies condition (3), then consider the affine base $uow$ with $u=h$ and $w\in\Aline ov\cap\Aline ad$. The diunit $e$ of the affine base $uow$ coincides with the point $a$. By definition of the dot operation in the ternar $\Delta=\Aline oe=\Aline oa$, the inclusion $d\in\Aline ob$ implies $c\cdot c=a=e$. The condition (3) implies $c+e=o$. By definition of the plus operation, the equality $c+e=o$ implies $\Aline hv=\Aline uv\parallel \Delta=\Aline ac$.
}
\vskip-5pt

%

\vbox{\hsize350pt If $X$ satisfies condition (4), then consider the affine base $uow$ with $u=h$ and $w=v$.  The diunit $e$ of the base $uow$ coincides with the point $b$. By definition of the dot operation in the ternar $\Delta=\Aline oe=\Aline ob$, the inclusion $c\in\Aline oa$ implies $d\cdot d=b=e$. The condition (4) implies $d\puls e=o$. By definition of the puls operation, the equality $d\puls e=o$ implies $\Aline hv=\Aline uw\parallel \Aline ao=\Aline ac$.
}\end{proof}

\begin{definition} A Playfair liner $X$ is called \defterm{invertible-half} if for any disjoint parallel lines $L,L'$ in $X$ and any distinct points $a,b\in L$, $a',b',c'\in L'$, $x\in\Aline a{b'}\cap\Aline {a'}b$, $y\in\Aline b{b'}\cap\Aline a{c'}$, if $\Aline a{a'}\parallel \Aline b{b'}$ and $\Aline a{b'}\parallel \Aline b{c'}$, then $\Aline xy\parallel L$.

\begin{picture}(200,90)(-120,-15)
\linethickness{=0.8pt}
\put(0,0){\color{red}\line(1,0){60}}
\put(0,60){\color{red}\line(1,0){120}}
\put(0,0){\color{cyan}\line(0,1){60}}
\put(60,0){\color{cyan}\line(0,1){60}}
\put(0,0){\color{blue}\line(1,1){60}}
\put(60,0){\color{blue}\line(1,1){60}}
\put(0,0){\line(2,1){120}}
\put(0,60){\line(1,-1){60}}
\put(30,30){\color{red}\line(1,0){30}}

\put(0,0){\circle*{3}}
\put(-3,-9){$a$}
\put(60,0){\circle*{3}}
\put(60,-9){$b$}
\put(30,30){\circle*{3}}
\put(20,28){$x$}
\put(60,30){\circle*{3}}
\put(62,24){$y$}
\put(0,60){\circle*{3}}
\put(-3,63){$a'$}
\put(60,60){\circle*{3}}
\put(58,63){$b'$}
\put(120,60){\circle*{3}}
\put(118,63){$c'$}
\end{picture}
\end{definition}

\begin{proposition}\label{p:by-Thalesian<=invertible-half} Every invertible-half invertible-puls Playfair liner is by-Thalesian.
\end{proposition}

\begin{proof}  By Proposition~\ref{p:by-Thalesian<=>}, the by-Thalesian property of $X$ will follow as soon as we show that for any disjoint parallel lines $L,L'\subseteq X$ and distinct points $a,b\in L$, $a',b'\in L'$, $x\in \Aline a{b'}\cap\Aline {a'}b$ and $z\in X$ with $\Aline ax\parallel \Aline bz$, $\Aline {a'}x\parallel \Aline{b'}z$, and $\Aline a{a'}\parallel \Aline b{b'}$, we have $\Aline xz\parallel L$. By Proclus Postulate~\ref{p:Proclus-Postulate}, there exist unique points $c\in L\cap\Aline {b'}z$ and $c'\in L'\cap\Aline bz$. The invertible-puls property of $X$ ensures that $\Aline c{c'}\parallel \Aline b{b'}$. Since the lines $\Aline a{b'}$ and $\Aline {a'}b$ are concurrent, the point $c$ is not equal to the point $a$. Then the lines $\Aline a{c'}$ and $\Aline c{c'}$ are concurrent and so are the lines $\Aline a{c'}$ and $\Aline b{b'}$ (by the Proclus Postulate). Then there exists a unique point $y\in \Aline a{c'}\cap\Aline b{b'}$.
Since $X$ is invertible-half, $\Aline xy\parallel L$ and $\Aline yz\parallel L'$. Then the lines $\Aline xy$ and $\Aline yz$ are parallel and hence coincide with the line $\Aline xz$. Therefore, $\Aline xz=\Aline xy\parallel L$, witnessing that   the liner $X$ is by-Thalesian.

\begin{picture}(200,90)(-120,-15)
\linethickness{=0.8pt}
\put(0,0){\color{teal}\line(1,0){120}}
\put(0,60){\color{teal}\line(1,0){120}}
\put(0,0){\color{cyan}\line(0,1){60}}
\put(60,0){\color{cyan}\line(0,1){60}}
\put(120,0){\color{cyan}\line(0,1){60}}
\put(0,0){\color{blue}\line(1,1){60}}
\put(60,0){\color{blue}\line(1,1){60}}
\put(0,0){\color{red}\line(2,1){120}}
\put(0,60){\color{violet}\line(1,-1){60}}
\put(60,60){\color{violet}\line(1,-1){60}}
\put(30,30){\color{teal}\line(1,0){60}}

\put(0,0){\circle*{3}}
\put(-3,-9){$a$}
\put(60,0){\circle*{3}}
\put(60,-9){$b$}
\put(120,0){\circle*{3}}
\put(118,-9){$c$}
\put(30,30){\circle*{3}}
\put(27,21){$x$}
\put(60,30){\circle*{3}}
\put(62,23){$y$}
\put(90,30){\circle*{3}}
\put(87,21){$z$}
\put(0,60){\circle*{3}}
\put(-3,63){$a'$}
\put(60,60){\circle*{3}}
\put(58,63){$b'$}
\put(120,60){\circle*{3}}
\put(118,63){$c'$}
\end{picture}

\end{proof}

Let us recall that a liner is called {\em invertible-add} if it is invertible-plus and invertible-puls. 

\begin{theorem}\label{t:by-Thalesian<=>by-Boolean<=>quadratic} If a Playfair liner $X$ is invertible-add, then the following conditions are equivalent:
\begin{enumerate}
\item $X$ is by-Thalesian;
\item $X$ is by-Boolean;
\item $X$ is quadratic;
\item $X$ is invertible-half;
\item for every ternar $R$ of $X$, if $1+1\ne 0$ in $R$, then there exists an element $\frac12\in R$ such that $\frac12+\frac12=1=\frac12\puls\frac12$;
\item for every ternar $R$ of $X$, if $1+1\ne 0$ in $R$, then there exists an element $\frac12\in R$ such that $\frac12+\frac12=1$.
\end{enumerate}
\end{theorem}

\begin{proof} The equivalence $(1)\Leftrightarrow(2)$ follows from 
Proposition~\ref{p:by-Thalesian<=>by-Boolean} and Corollary~\ref{c:invertible-add=>bi-Thalesian}. The equivalence $(1)\Leftrightarrow(3)$ follows from Propositions~\ref{p:by-Thalesian=>quadratic} and \ref{p:quadratic=>by-Thalesian}.
Therefore, the conditions (1), (2), (3) are equivalent. The implication $(4)\Ra(1)$ follows from Proposition~\ref{p:by-Thalesian<=invertible-half}.
\smallskip

$(1)\Ra(4)$ Assume that the invertible-add Playfair liner $X$ is by-Thalesian. To prove that $X$ is invertible-half, fix disjoint parallel lines $L,L'$ in $X$ and distinct points $a,b\in L$, $a',b',c'\in L'$, $x\in\Aline a{b'}\cap\Aline {a'}b$, $y\in\Aline b{b'}\cap\Aline a{c'}$ such that $\Aline a{a'}\parallel \Aline b{b'}$ and $\Aline a{b'}\parallel \Aline b{c'}$.  We have to prove that $\Aline xy\parallel L$. Since $X$ is Playfair, there exists a unique point $c\in L$ such that $\Aline c{c'}\parallel \Aline b{b'}$. Since $X$ is invertible-puls, $\Aline a{b'}\parallel \Aline b{c'}$ implies $\Aline {a'}b\parallel \Aline{b'}c$. Since the lines $\Aline a{b'}$ and $\Aline{a'}b$ are not parallel, the lines $\Aline b{c'}$ and $\Aline {b'}c$ are not parallel and hence have a common point $z$. Applying the by-Thalesian property of $X$ to the triangles $axa'$ and $bzb'$, we conclude that $\Aline xz\parallel L$. By Proclus Postulate~\ref{p:Proclus-Postulate}, there exist unique points $\alpha\in \Aline a{a'}\cap\Aline xz$ and $\beta\in\Aline b{b'}\cap\Aline xz$. Since $X$ is Playfair, there exists a point $u,v\in L'$ such that $\Aline xu\parallel \Aline b{b'}\parallel \Aline zv$.



 By Proposition~\ref{p:by-Thalesian=>quadratic}, the by-Thalesian Playfair liner $X$ is quadratic, which implies $\Aline \alpha u\parallel \Aline ax$. By Corollary~\ref{c:invertible-add=>bi-Thalesian}, the invertible-add Playfair liner $X$ is bi-Thalesian. The bi-Thalesian property of $X$, applied to the triangles $ax\beta$ and $\alpha ub'$ implies $\Aline a{\beta}\parallel \Aline\alpha {b'}$.
 
%
 
 
Since $X$ is invertible-puls, the parallelity relations $\Aline \alpha u\parallel \Aline x{b'}$ and $\Aline \alpha{a'}\parallel \Aline xu\parallel \Aline\beta {b'}$ imply  $\Aline {a'}x\parallel \Aline u\beta$. Since $X$ is invertible-puls, the parallelity relations $\Aline u\beta\parallel \Aline {b'}z$ and $\Aline xu\parallel\Aline\beta b'\parallel zv$ imply $\Aline x{b'}\parallel \Aline \beta v$. Since $X$ is invertible-puls, the parallelity relations $\Aline \alpha u\parallel \Aline x{b'}\parallel \Aline \beta v\parallel \Aline z{c'}$ and $\Aline xu\parallel \Aline\beta{b'}\parallel \Aline zv$ imply $\Aline \alpha{b'}\parallel \Aline xv\parallel \Aline \beta{c'}$. Therefore, $\Aline \beta{c'}\parallel \Aline \alpha{b'}\parallel \Aline a\beta$ and hence $\beta\in \Aline a{c'}\cap b{b'}=\{y\}$ and finally $\Aline xy=\Aline x\beta=\Aline xz\parallel L$.
\smallskip

$(1)\Ra(5)$ Let $R$ be a ternar of $X$ such that $1+1\ne 0$ in $R$. Find a plane $\Pi$ in $X$ and an affine base $uow$ in $\Pi$ whose ternar $\Delta$ is isomorphic to the ternar $R$. Observe that $\Delta=\Aline oe$, where $e$ is the diunit of the affine base $uow$. Since $1+1\ne 0$ in $R$, $e+e\ne o$ in the ternar $\Delta$. It follows from $e+e\ne o$ that the lines $\Aline uw$ and $\Aline oe$ are not parallel and have a unique common point $c\in\Delta$. Since $X$ is Playfair, there exist unique points $a\in \Aline ow$, $b\in \Aline we$, $\alpha\in \Aline eu$, $\beta\in\Aline ou$ such that $\Aline ac\parallel \Aline ou\parallel \Aline c\alpha$ and $\Aline bc\parallel \Aline ow\parallel\Aline c\beta$. 

\begin{picture}(60,78)(-150,-10)
{\linethickness{0.8pt}
\put(0,0){\color{teal}\line(1,0){60}}
\put(0,60){\color{teal}\line(1,0){60}}
\put(0,0){\color{blue}\line(1,1){60}}
\put(0,60){\color{violet}\line(1,-1){60}}
\put(0,30){\color{violet}\line(1,-1){30}}
\put(0,0){\color{cyan}\line(0,1){60}}
\put(60,0){\color{cyan}\line(0,1){60}}
\put(30,0){\color{cyan}\line(0,1){60}}
\put(0,30){\color{teal}\line(1,0){60}}

}
\put(0,30){\color{blue}\line(1,1){30}}
\put(30,0){\color{blue}\line(1,1){30}}
\put(30,60){\color{violet}\line(1,-1){30}}
\put(0,0){\color{red}\line(2,1){60}}
\put(0,30){\color{red}\line(2,1){60}}

\put(0,0){\circle*{3}}
\put(-7,-7){$o$}
\put(30,30){\circle*{3}}
\put(24,20){$c$}
\put(30,60){\circle*{3}}
\put(28,63){$b$}
\put(60,30){\circle*{3}}
\put(63,28){$\alpha$}
\put(0,30){\circle*{3}}
\put(-8,28){$a$}
\put(30,0){\circle*{3}}
\put(27,-10){$\beta$}
\put(60,0){\circle*{3}}
\put(62,-6){$u$}
\put(0,60){\circle*{3}}
\put(-8,62){$w$}
\put(60,60){\circle*{3}}
\put(63,60){$e$}
\end{picture}

The implication $(1)\Ra(3)$ ensures that the by-Thalesian liner $X$ is quadratic. Then $\Aline ab\parallel \Aline oe$ and hence $c+c=e$, by the definition of the plus operation in the ternar $\Delta$. Since $X$ is invertible-plus, $\Aline ab\parallel \Delta$ implies $\Aline\alpha\beta\parallel\Delta$. Since $X$ is invertible-puls, $\Aline ab\parallel \Delta\parallel\Aline\alpha\beta$ implies $\Aline wc\parallel \Aline b\alpha\parallel \Aline cu\parallel \Aline a\beta$. Consider the triangles $ace$ and $o\beta\alpha$, and observe that $\Aline ac\parallel \Aline ou$,  $\Aline ce\parallel \Aline \beta\alpha$, and $\Aline ao\parallel \Aline c\beta\parallel \Aline e\alpha$. Since $X$ is by-Thalesian, $\alpha\in \Aline ac$ and $ o\in\Aline ce$ imply $\Aline ae\parallel \Aline o\alpha$ and hence $c\puls c=e$, by the definition of the puls operation.  Let $\frac12$ be the image of the point $c$ under the isomorphism $\Delta\to R$. The equalities $c+c=e=c\puls c$ imply $\frac12+\frac12=1=\frac12\puls\frac12$.
\smallskip

The implication $(5)\Ra(6)$ is trivial. 
\smallskip

$(6)\Ra(3)$ Assume that the Playfair liner $X$ satisfies the condition (6). To show that $X$ is quadratic, take any points $u,o,w,e\in X$ and $c\in\Aline oe\cap \Aline uw$, $a\in \Aline ow\setminus \Aline we$, $b\in\Aline we\setminus\Aline ow$ such that $\Aline ow\parallel \Aline cb\parallel \Aline ue$ and $\Aline ou\parallel \Aline ac\parallel \Aline we$. We have to prove that $\Aline ab\parallel \Aline oe$. It follows from $a\in \Aline ow\setminus \Aline we$ and $b\in\Aline we\setminus\Aline ow$ that $o\ne w\ne e$, $\Aline ow\cap\Aline we=\{w\}$. Then the points $o,w,e$ are not collinear and  $\Pi\defeq\overline{\{o,w,e\}}$ is a plane in $X$. It follows from $c\in\Aline oe\cap \Aline uw$ that $o,u,e,w,c,a,b\in\Pi$ and hence $uow$ is an affine base in the plane $\Pi$ whose diunit coincides with the point $e$. By the condition (6), there exists a point $\gamma\in\Aline oe$ such that $\gamma+\gamma=e$. Since $X$ is Playfair, there exist unique points $\alpha\in\Aline ow$, $\beta\in\Aline we$, $\delta\in\Aline eu$, $\e\in\Aline uo$ such that $\Aline \alpha\gamma\parallel \Aline ou\parallel\Aline \gamma\delta$ and $\Aline \beta\gamma\parallel \Aline ow\parallel\Aline \gamma\e$. By definition of the plus operation, $\gamma+\gamma=e$ implies $\Aline \alpha\beta\parallel \Aline oe$. Since $X$ is invertible-plus, $\Aline \alpha\beta\parallel \Aline oe$ implies $\Aline \delta\e\parallel\Aline oe$. Since $X$ is invertible-puls, $\Aline \alpha\beta\parallel \Aline oe\parallel\Aline \delta\e$ implies $\Aline w\gamma\parallel \Aline\beta\delta\parallel\Aline u\gamma$ and hence $\gamma\in \Aline oe\cap\Aline uw=\{c\}$. Then $\alpha=a$ and $\beta=b$ and $\Aline ab=\Aline\alpha\beta\parallel \Aline oe$.

\begin{picture}(60,80)(-150,-5)
{\linethickness{0.8pt}
\put(0,0){\color{teal}\line(1,0){60}}
\put(0,60){\color{teal}\line(1,0){60}}
\put(0,0){\color{blue}\line(1,1){60}}
\put(0,60){\color{violet}\line(1,-1){60}}
\put(0,30){\color{violet}\line(1,-1){30}}
\put(0,0){\color{cyan}\line(0,1){60}}
\put(60,0){\color{cyan}\line(0,1){60}}
\put(30,0){\color{cyan}\line(0,1){60}}
\put(0,30){\color{teal}\line(1,0){60}}

}
\put(0,30){\color{blue}\line(1,1){30}}
\put(30,0){\color{blue}\line(1,1){30}}
\put(30,60){\color{violet}\line(1,-1){30}}

\put(0,0){\circle*{3}}
\put(-7,-7){$o$}
\put(30,30){\circle*{3}}
\put(31,38){$\gamma$}
\put(24,20){$c$}
\put(30,60){\circle*{3}}
\put(26,63){$b$}
\put(32,63){$\beta$}
\put(60,30){\circle*{3}}
\put(63,28){$\delta$}
\put(0,30){\circle*{3}}
\put(-8,31){$a$}
\put(-8,22){$\alpha$}
\put(30,0){\circle*{3}}
\put(27,-10){$\e$}
\put(60,0){\circle*{3}}
\put(62,-6){$u$}
\put(0,60){\circle*{3}}
\put(-8,62){$w$}
\put(60,60){\circle*{3}}
\put(63,60){$e$}
\end{picture}
\end{proof}

\begin{exercise} Show that the Moulton plane is by-Boolean but not quadratic. 
\end{exercise}


\section{Commutative-add Playfair liners are Boolean or Thalesian}

A liner is {\em commutative-add} if it is commutative-plus and commutative-puls.

\begin{theorem}\label{t:commutative-add=>quadratic} If a Playfair liner is commutative-add, then it is quadratic, invertible-half, by-Thalesian, by-Boolean and di-Boolean.
\end{theorem}

\begin{proof} Let $X$ be a commutative-add Playfair liner. Then $X$ is invertible-add. First, we prove that the liner $X$ is quadratic. 
Given any points $a,b,c,d\in X$, $o\in\Aline ac\cap\Aline bd$, $h\in\Aline ab\setminus\Aline bc$ and $v\in \Aline bc\setminus\Aline ab$ with $\Aline ab\parallel \Aline ov\parallel \Aline cd$ and $\Aline ad\parallel \Aline ho\parallel \Aline bc$, we should prove that $\Aline hv\parallel \Aline ac$. It follows that $abcd$ is a parallelogram, $o$ is the intersection of its diagonals, $h\notin\{a,b\}$ and $v\notin\{b,c\}$. Since $X$ is Playfair, there exist unique points $\delta\in \Aline ad$, $d'\in \Aline bd$, and $\delta'\in\Aline cd$ such that $\Aline h\delta\parallel \Aline bd$, $\Aline\delta{d'}\parallel \Aline ov$, and $\Aline {d'}{\delta'}\parallel \Aline oh$. Since the commutative-add liner $X$ is inversive-plus, $\Aline v{\delta'}\parallel \Aline bd$. Since $X$ is Playfair, there exist unique points $\alpha\in \Aline ad\cap\Aline vo$, $\alpha'\in \Aline ab\cap\Aline {d'}{\delta'}$, $a'\in \Aline vo\cap\Aline {\delta'}{d'}$, $\gamma\in \Aline bc\cap\Aline \delta{d'}$, $\gamma'\in \Aline ho\cap\Aline cd$ and $c'\in \Aline ho\cap\Aline\delta{d'}$. Since $X$ is commutative-plus, $\Aline \alpha{\alpha'}\parallel \Aline h\delta\parallel \Aline v{\delta'}\parallel \Aline \gamma{\gamma'}$. Since $X$ is inversive-puls, $\Aline hv\parallel \Aline a{a'}\parallel \Aline \delta{\delta'}\parallel \Aline {c'}c$.

\begin{picture}(100,110)(-150,-15)
\linethickness{0.8pt}
\put(0,0){\color{teal}\line(1,0){80}}
\put(0,20){\color{teal}\line(1,0){80}}
\put(0,60){\color{teal}\line(1,0){80}}
\put(0,80){\color{teal}\line(1,0){80}}

\put(0,0){\color{cyan}\line(0,1){80}}
\put(20,0){\color{cyan}\line(0,1){80}}
\put(60,0){\color{cyan}\line(0,1){80}}
\put(80,0){\color{cyan}\line(0,1){80}}

\put(0,80){\color{violet}\line(1,-1){80}}
\put(20,80){\color{violet}\line(1,-1){60}}
\put(0,60){\color{violet}\line(1,-1){60}}

\put(0,20){\color{violet}\line(1,-1){20}}
\put(60,80){\color{violet}\line(1,-1){20}}

\put(0,0){\color{blue}\line(1,1){20}}
\put(0,60){\color{blue}\line(1,1){20}}
\put(60,60){\color{blue}\line(1,1){20}}
\put(60,0){\color{blue}\line(1,1){20}}

\put(0,0){\circle*{3}}
\put(-8,-8){$a$}
\put(20,0){\circle*{3}}
\put(18,-8){$\alpha$}
\put(60,0){\circle*{3}}
\put(58,-9){$\delta$}
\put(80,0){\circle*{3}}
\put(82,-8){$d$}
\put(0,20){\circle*{3}}
\put(-10,18){$\alpha'$}
\put(20,20){\circle*{3}}
\put(22,22){$a'$}
\put(60,20){\circle*{3}}
\put(62,22){$d'$}
\put(80,20){\circle*{3}}
\put(83,18){$\delta'$}

\put(0,60){\circle*{3}}
\put(-8,58){$h$}
\put(20,60){\circle*{3}}
\put(12,54){$o$}
\put(60,60){\circle*{3}}
\put(52,51){$c'$}
\put(80,60){\circle*{3}}
\put(82,58){$\gamma'$}
\put(0,80){\circle*{3}}
\put(-6,82){$b$}
\put(20,80){\circle*{3}}
\put(18,83){$v$}
\put(60,80){\circle*{3}}
\put(59,84){$\gamma$}
\put(80,80){\circle*{3}}
\put(78,83){$c$}

\end{picture}

Since $X$ is commutative-puls, $\Aline hv\parallel \Aline a{a'}$ implies $\Aline ao\parallel \Aline{\alpha'}v$. Also $\Aline bo\parallel \Aline \gamma{\gamma'}$ implies $\Aline oc\parallel \Aline h\gamma$; $\Aline b{d'}\parallel \Aline v{\delta'}$ implies $\Aline{d'}c\parallel \Aline{\alpha'}v\parallel \Aline ao$; and finally $\Aline b{d'}\parallel \Aline h\delta$ implies $\Aline a{d'}\parallel \Aline h\gamma\parallel \Aline oc$. Now recall that $\Aline ao=\Aline oc$, which implies $\Aline a{d'}\parallel \Aline oc=\Aline oa\parallel \Aline {d'}c$ and hence $\Aline a{d'}=\Aline {d'}c$ and $d'\in \Aline ac\cap \Aline bd=\{o\}$. The equality $o=d'$ implies $a'=c'=o=d$ and hence $\Aline hv\parallel \Aline a{a'}=\Aline {c'}c=\Aline ac$, witnessing that the liner $X$ is quadratic.

\begin{picture}(100,110)(-150,-15)
\linethickness{0.8pt}
\put(0,0){\color{red}\line(1,3){20}}
\put(0,20){\color{red}\line(1,3){20}}
\put(20,60){\color{orange}\line(3,1){60}}
\put(0,60){\color{orange}\line(3,1){60}}

\put(60,20){\color{red}\line(1,3){20}}
\put(0,0){\color{orange}\line(3,1){60}}

\put(0,0){\color{teal}\line(1,0){80}}
\put(0,20){\color{teal}\line(1,0){80}}
\put(0,60){\color{teal}\line(1,0){80}}
\put(0,80){\color{teal}\line(1,0){80}}

\put(0,0){\color{cyan}\line(0,1){80}}
\put(20,0){\color{cyan}\line(0,1){80}}
\put(60,0){\color{cyan}\line(0,1){80}}
\put(80,0){\color{cyan}\line(0,1){80}}

\put(0,80){\color{violet}\line(1,-1){80}}
\put(20,80){\color{violet}\line(1,-1){60}}
\put(0,60){\color{violet}\line(1,-1){60}}

\put(0,20){\color{violet}\line(1,-1){20}}
\put(60,80){\color{violet}\line(1,-1){20}}

\put(0,0){\color{blue}\line(1,1){20}}
\put(0,60){\color{blue}\line(1,1){20}}
\put(60,60){\color{blue}\line(1,1){20}}
\put(60,0){\color{blue}\line(1,1){20}}

\put(0,0){\circle*{3}}
\put(-8,-8){$a$}
\put(20,0){\circle*{3}}
\put(18,-8){$\alpha$}
\put(60,0){\circle*{3}}
\put(58,-9){$\delta$}
\put(80,0){\circle*{3}}
\put(82,-8){$d$}
\put(0,20){\circle*{3}}
\put(-10,18){$\alpha'$}
\put(20,20){\circle*{3}}
\put(22,22){$a'$}
\put(60,20){\circle*{3}}
\put(62,22){$d'$}
\put(80,20){\circle*{3}}
\put(83,18){$\delta'$}

\put(0,60){\circle*{3}}
\put(-8,58){$h$}
\put(20,60){\circle*{3}}
\put(13,54){$o$}
\put(60,60){\circle*{3}}
\put(52,51){$c'$}
\put(80,60){\circle*{3}}
\put(82,58){$\gamma'$}
\put(0,80){\circle*{3}}
\put(-6,82){$b$}
\put(20,80){\circle*{3}}
\put(18,83){$v$}
\put(60,80){\circle*{3}}
\put(59,84){$\gamma$}
\put(80,80){\circle*{3}}
\put(78,83){$c$}
\end{picture}

 By Theorem~\ref{t:by-Thalesian<=>by-Boolean<=>quadratic}, the quadratic invertible-add liner $X$ is invertible-half, by-Thalesian and by-Boolean.
 
Assuming that $X$ is not di-Boolean, we can find two disjoint lines $V,E$ and distinct points $o,w,a\in V$ and $u,e,b\in E$ such that $\varnothing=\Aline oe\cap\Aline uw=\Aline ou\cap\Aline we=\Aline we\cap\Aline ab$ but $\Aline wb\cap\Aline ae\ne\varnothing$. The latter inequality implies that $\Pi=\overline{V\cup E}=\overline{\{w,b,a,e\}}$ is a plane in $X$, $uow$ is an affine base in the plane $\Pi$, and $e$ is the diunit of the affine base $uow$. It follows from $\Aline oe\cap\Aline uw=\varnothing$ that $e+e=o=e\puls e$ in the ternar $\Delta=\Aline oe$ of the affine base $uow$. Consider the unique point $c\in\Aline ab\cap \Delta$. Since $X$ is commutative-puls, the puls loop $(\Delta,\!\puls\!)$ of the ternar $\Delta$ is a commutative group, by Theorem~\ref{t:com-puls=>ass-puls}.  Then there exists a unique element $\gamma\in \Delta$ such that $c=e\puls \gamma=\gamma\puls e$. Find unique points $\alpha\in V$ and $\beta\in E$ such that $\Aline \alpha\beta\parallel \Aline ou$ and $\gamma\in \Aline \alpha\beta$. It follows from $c=\gamma\puls e$ that $\Aline o\beta\parallel \Aline wb$. Since $X$ is inversive-puls, the latter parallelity relation implies $\Aline u\alpha\parallel \Aline ea$. Since the lines $\Aline wb$ and $\Aline ea$ are concurrent, so are the lines $\Aline o\beta$ and $\Aline u\alpha$. Then $\gamma\puls \gamma\ne o$ and $c\puls c=(e\puls \gamma)\puls (e\puls \gamma)=(e\puls e)\puls(\gamma\puls \gamma)=o\puls(\gamma\puls\gamma)=\gamma\puls\gamma\ne o$, which implies $\Aline ob\cap \Aline ua\ne\varnothing$. Since the invertible-add liner $X$ is quadratic, we can apply Theorem~\ref{t:by-Thalesian<=>by-Boolean<=>quadratic} and find elements $c',\gamma'\in \Delta$ such that $c'\puls c'=c$ and $\gamma'\puls\gamma'=\gamma$. Since $(\Delta,\!\puls\!)$ is a group, there exists a unique element $e'\in\Delta$ such that $c'=\gamma'\puls e'$. Observe  that
$$\gamma\puls e=c=c'\puls c'=(\gamma'\puls e')\puls(\gamma'\puls e')=(\gamma'\puls\gamma')\puls(e'\puls e')=\gamma\puls(e'\puls e')$$ and hence $e=e'\puls e'$, by the cancellativity of the group operation $\puls$. Then the element $e'$ has order $4$ in the group $(\Delta,\!\puls\!)$.  

Find a unique element $w'\in L$ such that $\Aline {w'}{e'}\parallel \Aline ou$ and observe that $uow'$ is an affine base in the plane $\Pi$ such that $e'$ is the diunit of the affine base $uow'$. The equality $e'\puls e'\puls e'\puls e'=o$ in the ternar $\Delta=\Aline oe$ implies that $e'\puls e'\puls e'\puls e'=o$ in the ternar $\Delta'=\Aline o{e'}$ of the affine base $uow'$. 
By Theorem~\ref{t:invertible-add<=>Z-elementary}, the ternar $\Delta'$ is $\IZ$-elementary. Then there exists a function $h:\IZ\to \Delta'$ such that $h(1)=e'$ and $h(x{\cdot} y+z)=h(x)_\times h(y)_+h(z)$ for all numbers $x,y,z\in \IZ$. In particular, $h(2)=h(1{\cdot}1+1)=h(1)_\times h(1)_+h(1)=e'+e'=e'\puls e'$ and $h(4)=h(1\cdot 2+2)=h(1)_\times h(2)_+h(2)=h(2)\puls h(2)=(e'\puls e')\puls(e'\puls e')=o$. Then $o=h(4)=h(2\cdot 2)=h(2)_\times h(2)_+h(0)=(e'\puls e')\cdot(e'\puls e')=e\cdot e=e$, which is a contradiction showing that the liner $X$ is di-Boolean.
\end{proof}

Now we can prove the main result of this chapter.

\begin{theorem}\label{t:commutative-add<=>} For a Playfair liner $X$ the following conditions are equivalent:
\begin{enumerate}
\item $X$ is commutative-add;
\item $X$ is uno-Thalesian;
\item $X$ is Boolean or Thalesian.
\end{enumerate}
\end{theorem}

\begin{proof} $(1)\Ra(2)$. Assume that a Playfair liner $X$ is commutative-add. By Theorem~\ref{t:commutative-add=>quadratic}, $X$ is by-Thalesian and di-Boolean. To prove that $X$ is uno-Thalesian, we shall apply Proposition~\ref{p:uno-Thalesian1}. Given any distinct parallel lines $A,B,C$ in $X$ and points $a,a'\in A$, $b,b'\in B$, $c,c'\in C$ such that $b\in\Aline {a'}{c'}$ and $\Aline ab\cap\Aline {a'}{b'}=\varnothing=\Aline ac\cap\Aline{a'}{c'}$, we should prove that $\Aline bc\cap\Aline{b'}{c'}=\varnothing$. Consider the plane $\Pi\defeq\overline{A\cup B\cup C}$. Since $X$ is Playfair, there exists a unique point $a''\in A$ such that $\Aline {a''}{b'}\parallel \Aline ac\parallel \Aline{a'}{c'}$. The Proclus property of the Playfair liner $X$ ensures that there exist unique points $c''\in \Aline {a''}{b'}\cap\Aline c{c'}$ and $b''\in \Aline b{b'}\cap\Aline {a}{c}$. 

If $b'\in \Aline ac$, then $b''=b'$, $c''=c$ and the parallelogram $aa'bb''$ is Boolean. Since $X$ is di-Boolean, the parallelogram $cb''bc'$ is Boolean, too, and hence $\Aline bc\cap\Aline {b'}{c'}=\Aline bc\cap\Aline {b''}{c'}=\varnothing$. 

If $c''\in \Aline ab$, then the invertibility-plus of $X$ implies $\Aline {b''}{c'}\parallel \Aline ab\parallel \Aline b{c''}$. Since $X$ is invertible-puls, the parallelity relations $\Aline {b''}c\parallel \Aline b{c'}\parallel \Aline{b'}{c''}$ and $\Aline {b''}{c'}\parallel \Aline b{c''}$ imply $\Aline bc\parallel \Aline{b'}{c'}$ and we are done.

\begin{picture}(100,70)(-170,-15)
\linethickness{0.7pt}

\put(0,0){\color{teal}\line(1,0){40}}
\put(0,20){\color{teal}\line(1,0){40}}
\put(0,40){\color{teal}\line(1,0){40}}

\put(0,0){\color{cyan}\line(0,1){40}}
\put(20,0){\color{cyan}\line(0,1){40}}
\put(40,0){\color{cyan}\line(0,1){40}}
\put(0,0){\color{red}\line(1,1){20}}
\put(20,0){\color{red}\line(1,1){20}}
\put(0,20){\color{violet}\line(1,-1){20}}
\put(0,40){\color{violet}\line(1,-1){40}}
\put(20,40){\color{violet}\line(1,-1){20}}

\put(40,0){\circle*{3}}
\put(38,-9){$c''$}
\put(40,40){\circle*{3}}
\put(38,43){$a''$}
\put(0,0){\circle*{3}}
\put(-2,-9){$c$}
\put(20,0){\circle*{3}}
\put(18,-9){$c'$}
\put(0,40){\circle*{3}}
\put(-3,43){$a$}
\put(20,40){\circle*{3}}
\put(18,43){$a'$}
\put(20,20){\circle*{3}}
\put(21,21){$b$}
\put(40,20){\circle*{3}}
\put(41,21){$b'$}
\put(0,20){\circle*{3}}
\put(-10,19){$b''$}
\end{picture}

So, assume that $b',c''\notin \Aline ab$. 
Since the lines $\Aline ab$ and $\Aline {a'}{b'}$ are disjoint, the unique points $\gamma\in \Aline {a}{b}\cap\Aline c{c'}$ and  $\gamma'\in \Aline{a'}{b'}\cap\Aline c{c'}$ are distinct. Since $X$ is Playfair, there exist a unique point $\gamma''\in \Aline c{c'}$ such that $\Aline {b''}{\gamma''}\parallel \Aline ab$. Since $X$ is Playfair, there exist unique points $\beta,\beta',\beta''\in \Aline b{b'}$ such that the lines $\Aline{\gamma}{\beta}$, $\Aline {\gamma'}{\beta'}$, $\Aline{\gamma''}{\beta''}$  are parallel to the line $\Aline ac$. By the Proclus property of the Playfair liner $X$, there exist a unique points $\alpha\in \Aline {\beta}{\gamma}\cap\Aline a{a'}$, $\alpha''\in\Aline{\beta''}{\gamma''}\cap\Aline a{a'}$, $\delta\in \Aline{b''}{\gamma''}\cap \Aline{\gamma}{\beta}$, $\delta'\in\Aline a\gamma\cap\Aline{\gamma'}{\beta'}$, and $\varepsilon\in \Aline ab\cap\Aline{\gamma''}{\beta''}$ and $\varepsilon'\in \Aline {a'}{b'}\cap\Aline{\gamma}{\beta}$.


\begin{picture}(100,120)(-150,-30)
\linethickness{0.7pt}

\put(0,0){\color{teal}\line(1,0){90}}
\put(0,60){\color{teal}\line(1,0){90}}
\put(0,75){\color{teal}\line(1,0){75}}
\put(60,15){\color{teal}\line(1,0){15}}
\put(75,-15){\color{teal}\line(1,0){15}}



\put(0,0){\color{cyan}\line(0,1){75}}
\put(15,0){\color{cyan}\line(0,1){75}}
\put(30,0){\color{cyan}\line(0,1){75}}
\put(60,0){\color{cyan}\line(0,1){75}}
\put(75,-15){\color{cyan}\line(0,1){90}}
\put(90,-15){\color{cyan}\line(0,1){75}}
\put(0,60){\color{violet}\line(1,-1){75}}
\put(0,75){\color{violet}\line(1,-1){90}}
\put(15,75){\color{violet}\line(1,-1){75}}
\put(60,75){\color{violet}\line(1,-1){15}}


\put(30,0){\circle*{3}}
\put(28,-9){$c''$}
\put(30,75){\circle*{3}}
\put(28,78){$a''$}
\put(0,0){\circle*{3}}
\put(-2,-9){$c$}
\put(15,0){\circle*{3}}
\put(13,-9){$c'$}
\put(0,75){\circle*{3}}
\put(-3,78){$a$}
\put(15,75){\circle*{3}}
\put(13,78){$a'$}
\put(15,60){\circle*{3}}
\put(9,52){$b$}
\put(30,60){\circle*{3}}
\put(31,62){$b'$}
\put(0,60){\circle*{3}}
\put(-10,59){$b''$}
\put(90,0){\circle*{3}}
\put(93,-2){$\gamma'$}
\put(60,0){\circle*{3}}
\put(51,-10){$\gamma''$}
\put(90,60){\circle*{3}}
\put(92,57){$\beta'$}
\put(60,75){\circle*{3}}
\put(58,78){$\alpha''$}

\put(75,0){\circle*{3}}
\put(68,-7){$\gamma$}
\put(60,60){\circle*{3}}
\put(48,50){$\beta''$}
\put(75,60){\circle*{3}}
\put(67,51){$\beta$}
\put(75,75){\circle*{3}}
\put(73,78){$\alpha$}
\put(75,-15){\circle*{3}}
\put(72,-25){$\delta$}
\put(90,-15){\circle*{3}}
\put(88,-25){$\delta'$}
\put(60,15){\circle*{3}}
\put(52,12){$\varepsilon$}
\put(75,15){\circle*{3}}
\put(78,14){$\varepsilon'$}


\end{picture}

We claim that the line $\Aline \varepsilon{\varepsilon'}$ is distinct from the lines $\Aline a{a'}$ and $\Aline c{c'}$. Assuming that 
$\Aline \varepsilon{\varepsilon'}=\Aline a{a'}$, we conclude that 
$\varepsilon'\in \Aline {a'}{b'}\cap \Aline \varepsilon{\varepsilon'}=\Aline {a'}{b'}\cap\Aline a{a'}=\{a'\}$ and then $\Aline \alpha\gamma=\Aline {a'}b$ and hence $\gamma\in \Aline \alpha\gamma\cap\Aline ab=\Aline {a'}b\cap\Aline ab=\{b\}$,
which contradicts $\Aline b{b'}\cap\Aline c{c'}=\varnothing$. 
Assuming that
$\Aline \varepsilon{\varepsilon'}=\Aline b{b'}$, we conclude that $\e=b$ and $\e'=b'$. Then $c''\in \Aline {c''}{b'}\cap\Aline c{c'}=\Aline {\e'}\gamma\cap\Aline c{c'}=\{\gamma\}\subseteq\Aline ab$, which contradicts the assumption $c''\in \Aline ab$. This contradiction shows that $\Aline \varepsilon{\varepsilon'}\ne\Aline b{b'}$. Since $X$ is commutative-plus, the parallelity relations $\Aline {a}{b''}\parallel \Aline {\alpha''}{\gamma''}\parallel \Aline \gamma\beta\nparallel \Aline {a}{\alpha''}\parallel \Aline {b''}{\beta}\parallel \Aline{\gamma''}{\gamma}$ and $\Aline a{\gamma}\parallel \Aline {b''}{\gamma''}$ imply $\Aline {\alpha''}{\beta}\parallel \Aline {a}{\gamma}=\Aline ab$.
By the same reason, the parallelity relations $\Aline {a'}{b}\parallel \Aline {\alpha''}{\varepsilon}\parallel \Aline \beta{\varepsilon'}\nparallel \Aline {b}{\varepsilon}\parallel \Aline {a'}{\varepsilon'}\parallel \Aline{\alpha''}{\beta}$ and $\Aline {a'}{\alpha''}\parallel \Aline b{\beta}$ imply $\Aline {\varepsilon}{\varepsilon'}\parallel \Aline {b}{\beta}$. Taking into account that $X$ is invertible-plus, we can see that $\Aline\delta{\delta'}\parallel \Aline \varepsilon{\varepsilon}\parallel \Aline c{c'}$.



By the Proclus property of the Playfair liner $X$, there exist unique points $d\in \Aline {\alpha''}{\beta}\cap\Aline{a'}b$, $\phi\in \Aline {\alpha''}\beta\cap\Aline{\beta'}{\gamma'}$, and $f\in \Aline {a'}b\cap\Aline {b''}{\gamma''}$. Since $X$ is commutative-plus, the parallelity relations $\Aline bf\parallel \Aline \beta\delta\parallel \Aline \phi{\delta'}\nparallel \Aline f\delta\parallel\Aline b{\delta'}\parallel \Aline \beta\phi$ and $\Aline \delta{\delta'}\parallel \Aline b{\beta'}$ imply $\Aline f\phi\parallel \Aline b{\beta'}$.
Since $X$ is inversive-puls, $\Aline {\alpha''}{\beta}\parallel \Aline\varepsilon{\gamma}$ implies $\Aline {\beta''}{\alpha}\parallel \Aline{\gamma''}{\varepsilon'}$. 

Applying the commutativity-puls to the triples $\beta'bb'$ and $c'\gamma'\gamma$, we conclude that $\Aline {c'}{b'}\parallel \Aline \gamma{\beta'}$. Applying the commutativity-puls to the triples $\delta\beta\gamma$ and $\beta'\delta'\phi$, we conclude that $\Aline \delta\phi\parallel \Aline\gamma{\beta'}\parallel \Aline {c'}{b'}$.  Applying the commutativity-puls to the triples $\delta f b''$ and $d\beta\phi$, we obtain $\Aline d{b''}\parallel\Aline\delta{\phi}\parallel \Aline{c'}{b'}$. 



Applying the commutativity-puls to the triples $\beta b''b$ and $c\gamma\gamma''$, we conclude that $\Aline {c}{b}\parallel \Aline {\gamma''}{\beta}$.  Since $X$ is inversive-puls, $\Aline {c}{b}\parallel \Aline {\gamma''}{\beta}$ implies $\Aline {b''}{c'}\parallel \Aline {\beta''}\gamma$.
Applying the commutativity-puls of $X$ to the triples $\gamma \beta\varepsilon'$ and $\alpha''\gamma''\beta''$, we conclude that $\Aline {\alpha''}{\varepsilon'}\parallel \Aline \gamma{\beta''}\parallel \Aline {c'}{b''}$.  Applying the commutativity-puls to the triples $\varepsilon'a'b'$ and $d\beta\alpha''$,  we obtain $\Aline d{b'}\parallel\Aline{\varepsilon'}{\alpha''}\parallel \Aline {c'}{b''}$. Since the liner $X$ is by-Thalesian, the parallelity relations $\Aline {b''}{c'}\parallel \Aline d{b'}$ and $\Aline d{b''}\parallel \Aline {b'}{c'}$ imply $\Aline bc\parallel \Aline {b'}{c'}$, witnessing that the Playfair liner $X$ is uno-Thalesian.

$(2)\Ra(3)$ If $X$ is uno-Thalesian, then $X$ is Boolean or Thalesian, by Theorem~\ref{t:uno-Thales<=>}.

$(3)\Ra(1)$ If $X$ is Boolean or Thalesian, then $X$ is commutative-add, by Theorems~\ref{t:Boolean=>commutative-plus}, \ref{t:Boolean=>commutative-puls} and Corollaries~\ref{c:Thales=>commutative-plus}, \ref{c:Thales=>commutative-puls}.
\end{proof}



\chapter{Associativity and commutativity of Moufang loops}\label{ch:AM-numbers}

In this chapter we discuss some results on the commutativity and associativity of finite Moufang loops depending on arithmetic properties of their orders, and then will apply those results to Playfair liners whose order has some specific arithmetic properties.

\section{Associativity of Moufang loops}

In this section we present some sufficient conditions of associativity of Moufang loops, which are based on a nice identity discovered by \index[person]{Gagola III}Stephen Gagola III\footnote{{\bf Stephen M. Gagola III} is an American algebraist and specialist in loop theory, particularly Moufang loops. He completed his Ph.D. at Michigan State University in 2005 with a dissertation titled ``Subloops of the Unit Octonions'' (Department of Mathematics, Michigan State University), in which he studied structural properties of loops arising from octonion algebras and proved consequences such as Lagrange’s theorem for finite Moufang loops.
Gagola has held faculty positions in the United States (including at Kent State University) and has published numerous articles on Moufang loops, Sylow $p$-subloops, Lagrange’s theorem in loops, and related non-associative algebraic structures.}  in \cite{GagolaIII2014}. In that identity, for elements $c,x$ of a Moufang loop $X$ we denote by $x^c$ the conjugated element $c^{-1}{\cdot}x{\cdot}c$ to $x$. Also for any $n\in \IZ$ we denote by $x^n$ the $n$-power of $x$ in $X$.

\begin{theorem}[Gagola III, 2013]\label{t:GagolaIII} Let $X$ be a Moufang loop. The identity
$$(x{\cdot}c^{3m})\cdot(y{\cdot}c^{3n})=(x^{c^{2m+n}}\cdot y^{c^{-m+n}})^{c^{-2m-n}}\cdot c^{3m+3n}$$
holds for all $x,y,c\in X$ and $m,n\in\IZ$.
\end{theorem}

\begin{proof} The required identity follows from the Bol and Moufang identities and the diassociativity of the Moufang loop $X$:
$$
\begin{aligned}
(x{\cdot}c^{3m})\cdot (y{\cdot}c^{3n})&=\big(\big(((x{\cdot}c^{2m+n}){\cdot}c^{m-n}){\cdot}(y{\cdot}c^{3n})\big){\cdot}c^{m-n}\big){\cdot}c^{-m+n}\\
&=\big((x{\cdot}c^{2m+n}){\cdot}(c^{m-n}{\cdot}y{\cdot}c^{m+2n})\big){\cdot}c^{-m+n}\\
&=\big((c^{2m+n}{\cdot}(c^{-2m-n}{\cdot}x{\cdot}c^{2m+n}))\cdot ((c^{m-n}{\cdot}y{\cdot}c^{-m+n}){\cdot}c^{2m+n})\big){\cdot}c^{-m+n}\\
&=c^{2m+n}{\cdot}\big((c^{-2m-n}{\cdot}x{\cdot}c^{2m+n})\cdot(c^{m-n}{\cdot}y{\cdot}c^{-m+n})\big)\cdot c^{-2m-n}{\cdot}c^{3m+3n}\\
&=(x^{c^{2m+n}}y^{-m+n})^{c^{-2m-n}}\cdot c^{3m+3n}.
\end{aligned}
$$
\end{proof}

\begin{corollary}\label{c:GagolaI} Let $N$ be a normal subgroup of a Moufang loop $X$, $g$ be an element of $X$, and $g^{3\IZ}\defeq\{g^{3n}:n\in\IZ\}$ be the cyclic subgroup generated by the element $g^3$. If the conjugation map $(\cdot)^{g}:N\to N$, $(\cdot)^{g}:x\mapsto x^g$, is an automorphism of the group $N$, then the subloop $Ng^{\IZ}\defeq\{x{\cdot}g^{3n}:x\in N,\;n\in\IZ\}$ of $X$ is associative.
\end{corollary}

\begin{proof} If the conjugation map $(\cdot)^g$ is an automorphism of the group $N$, then so are its iterations $(\cdot)^{g^n}$ for all  $n\in\IZ$. Applying Theorem~\ref{t:GagolaIII} and the diassociativity of the Moufang loop $X$, we obtain the identity
\begin{equation}\label{eq:GagolaIII}
\begin{aligned}
(x{\cdot} g^{3m})\cdot(y{\cdot}g^{3n})&=
(x^{g^{2m+n}}\cdot y^{g^{-m+n}})^{g^{-2m-n}}\cdot g^{3m+3n}\\
&=
((x^{g^{2m+n}})^{g^{-2m-n}}\cdot (y^{g^{-m+n}})^{g^{-2m-n}})\cdot g^{3m+2n}\\
&=
(x\cdot y^{g^{-3m}})\cdot g^{3m+3n}
\end{aligned}
\end{equation}
holding for all $x,y\in N$ and  $m,n\in\IZ$.
\smallskip

 To prove that the loop $Ng^{3\IZ}$ is associative, fix any elements $a,b,c\in Ng^{3\IZ}$ and find elements $x,y,z\in N$ and integer numbers $m,n,k$ such that $a=x{\cdot}g^{3m}$, $b=y{\cdot}g^{3n}$, and $c=z{\cdot}g^{3k}$. 
By the identity (\ref{eq:GagolaIII}) and the associativity of $N$, we obtain 
$$
\begin{aligned}
(a\cdot b)\cdot c&=
\big((x{\cdot}g^{3m})\cdot(y{\cdot}g^{3n})\big){\cdot}(z{\cdot}g^{3k})\\
&= 
\big((x\cdot y^{g^{-3m}}){\cdot}g^{3m+3n}\big)\cdot (z{\cdot}g^{3k})=
\big((x\cdot y^{g^{-3m}})\cdot z^{g^{-3m-3n}}\big)\cdot g^{3m+3n+3k}
\end{aligned}
$$
and
$$
\begin{aligned}
a\cdot (b\cdot c)&=
(x{\cdot}g^{3m})\cdot\big((y{\cdot}g^{3n}){\cdot}(z{\cdot}g^{3k})\big)= 
(x{\cdot}g^{3m})\cdot \big((y{\cdot} z^{g^{-3n}}){\cdot}g^{3n+3k}\big)\\
&=
\big(x\cdot (y{\cdot}z^{g^{-3n}})^{g^{-3m}}\big)\cdot g^{3m+3n+3k}=\big(x\cdot (y^{g^{-3m}}\cdot z^{g^{-3n-3m}})\big)\cdot g^{3m+3n+3k}\\
&=\big((x\cdot y^{g^{-3m}})\cdot z^{g^{-3m-3n}}\big)\cdot g^{3m+3n+3k}=(a\cdot b)\cdot c.
\end{aligned}
$$
Therefore, the loop $Ng^{3\IZ}$ is associative.
\end{proof}



A Moufang loop $X$ is called \defterm{$n$-divisible} for $n\in\IN$ if for every $x\in X$ there exist an element $y\in X$ such that $x=y^n$.

 \begin{corollary}\label{c:GagolaIII} Let $N$ be a normal Abelian subgroup of a Moufang loop $X$, $g$ be any element of $X$, and $g^{3\IZ}\defeq\{g^{3n}:n\in\IZ\}$ be the cyclic subgroup, generated by the element $g^3$. If $N$ is $2$-divisible, then the subloop $N{\cdot}g^{3\IZ}$ is associative.
 \end{corollary}
 
 \begin{proof} By Corollary~\ref{c:GagolaI}, it suffices to show that the conjugation map $(\cdot)^g:N\to N$, $(\cdot)^g:x\mapsto x^g$, is an automorphism of the Abelian group $N$. Fix any elements $x,y\in N$. Since $N$ is $2$-divisible, there exists an element $z\in N$ such that $x=z^2$. By Theorem~\ref{t:inner=>pseudo-semi-auto}, the conjugation map $(\cdot)^g:N\to N$ is a semi-automorphism of $N$ and hence
 $$(x\cdot y)^g=(z^2\cdot y)^g=(z\cdot y\cdot z)^g=z^g\cdot y^g\cdot z^g=(z^g)^2\cdot y^g=(z^2)^g\cdot y^g=x^g\cdot y^g,$$
 witnessing that the conjugation $(\cdot)^g$ is an automorphism of the group $N$. By Corollary~\ref{c:GagolaI}, the loop $Ng^{3\IZ}$ is associative.
 \end{proof}

A loop $X$ is called \defterm{central-by-abelian} if the quotient loop $X/\mathcal Z(X)$ of $X$ by its center $\mathcal Z(X)$ is an Abelian group. It is clear that every central-by-abelian loop is centrally nilpotent. In Theorem~\ref{t:central-by-abelian} we shall prove that every $6$-divisible central-by-abelian Moufang loop is associative. We precede the proof of this theorem with three lemmas.

\begin{lemma}\label{l:commuting-conjugations} Let $X$ be a central-by-abelian inversive loop. Then 
$$x^{y^z}=x^y\quad\mbox{and}\quad (x^y)^z=(x^z)^y$$ for all $x,y,z\in X$.
\end{lemma}

\begin{proof} Let $Z\defeq\mathcal Z(X)$ be the centre of the inversive loop $X$ and $\pi:Z\to X/Z$ be the quotient homomorphism. Fix any elements $x,y,z\in X$. The commutativity of the quotient loop $X/Z$ ensures that $\pi(x^z)=\pi(x)$ and hence $y^z=c{\cdot}y$ for some $c\in Z$. Taking into account that the inversion is an anti-automorphism of the inversive loop $X$, we conclude that $c^{-1}\in \mathcal Z(X)\subseteq \mathcal N(X)$. Then 
$$
\begin{aligned}
x^{y^z}&=x^{c{\cdot}y}=(c{\cdot}y)^{-1}(x{\cdot}(c{\cdot}y))=(y^{-1}{\cdot}c^{-1}){\cdot}((x{\cdot}c){\cdot}y)\\
&=y^{-1}{\cdot}(c^{-1}{\cdot}((c{\cdot}x){\cdot}y))=y^{-1}{\cdot}\big(c^{-1}{\cdot}(c{\cdot}(x{\cdot}y))\big)=y^{-1}{\cdot}(x{\cdot}y)=x^y.
\end{aligned}
$$

Since the quotient loop $X/Z$ is a commutative group, $x^y=a{\cdot}x$ and $x^z=b{\cdot}x$ for some $a,b\in Z$. Then 
$$
\begin{aligned}
(x^y)^z&=(a{\cdot}x)^z=z^{-1}{\cdot}((a{\cdot}x){\cdot}z)=z^{-1}{\cdot}((x{\cdot}a){\cdot}z)=z^{-1}{\cdot}(x{\cdot}(a{\cdot}z))=z^{-1}{\cdot}(x{\cdot}(z{\cdot}a))\\
&=z^{-1}{\cdot}((x{\cdot}z){\cdot}a)=(z^{-1}{\cdot}(x{\cdot}z)){\cdot}a=x^z{\cdot}a=(b{\cdot}x){\cdot}a=a{\cdot}(b{\cdot}x)=(a{\cdot}b){\cdot}x.
\end{aligned}
$$
By analogy we can show that $(x^z)^y=(b{\cdot}a){\cdot}x$, which implies $(x^y)^z=(a{\cdot}b){\cdot}x=(b{\cdot}a){\cdot}x=(x^y)^z$.
\end{proof}

\begin{lemma}\label{l:central-by-abelian=>automorphism} For every normal $2$-divisible  subloop $N$ of a  central-by-abelian Moufang loop $X$ and every element $c\in X$, the conjugation map $(\cdot)^c:N\to N$, $(\cdot)^c:x\mapsto x^c$, is an automorphism of the loop $N$.
\end{lemma}

\begin{proof} Given any elements $x,y\in N$, we need to check that $(x{\cdot}y)^c=x^c{\cdot}y^c$. Since the loop $N$ is $2$-divisible, there exists an element $z\in N$ such that $y=z^2$. Then $x{\cdot}y=x{\cdot}z^2=z{\cdot}x^{z}{\cdot}z$. By Theorem~\ref{t:inner=>pseudo-semi-auto}, the conjugation map $(\cdot)^c$ is a semi-automorphism of $N$, which implies $$(x\cdot y)^c=(z\cdot x^{z}\cdot z)^c=z^c\cdot (x^{z})^c\cdot z^c=((x^z)^c)^{(z^{-1})^c}\cdot (z^c)^2=((x^c)^z)^{z^{-1}}\cdot y^c=x^c\cdot y^c,$$
according to Lemma~\ref{l:commuting-conjugations}.
\end{proof}

Let us recall that a magma $X$ is \defterm{locally cyclic} if every finite subset of $X$ is contained in a cyclic subgroup of $X$.

\begin{lemma}\label{l:induction-associative} Let $Z$ be the centre of a Moufang loop $X$ and $\pi:X\to X/Z$ be the quotient homomorphism. Assume the group $Z$ is $2$-divisible, the quotient loop $X/Z$ is Abelian, $D\subseteq X/Z$ is a subgroup of $X/Z$, and $A$ is a $3$-divisible locally cyclic subgroup in $X/Z$.
\begin{enumerate}
\item If the group $D$ is $2$-divisible, then the subloop $\pi^{-1}[D]$ is $2$-divisible, too.
\item If the loop $\pi^{-1}[D]$ is $2$-divisible and associative, then the loop $\pi^{-1}[DA]$ is associative.
\end{enumerate}
\end{lemma}

\begin{proof} 1. Assume that the group $D$ is $2$-divisible. To show the the loop $\pi^{-1}[D]$ is $2$-divisible, take any element $x\in \pi^{-1}[D]$ and consider its image $\pi(x)\in D$. Since $D$ is $2$-divisible, there exists an element $d\in D$ such that $d^2=\pi(x)$. Choose any element $y\in \pi^{-1}(d)$ and observe that $\pi(y^2)=d^2=\pi(x)$ and hence $x=z{\cdot}y^2$ for some element $z\in Z$. Since the group $Z$ is $2$-divisible, there exists an element $c\in Z$ such that $c^2=z$. Consider the element $c{\cdot}y\in \pi^{-1}(\pi(y))=\pi^{-1}[D]$ and observe that $(c{\cdot}y)^2=c^2{\cdot}y^2=z{\cdot}y^2=x$, witnessing that the loop $\pi^{-1}[D]$ is $2$-divisible.
\smallskip

2. Now assume that the loop $\pi^{-1}[D]$ is $2$-divisible and associative. To show that the loop $\pi^{-1}[DA]$ is associative, take any elements $x_1,x_2,x_3\in \pi^{-1}[DA]$. For every $i\in\{1,2,3\}$, find  elements $d_i\in D$ and $a_i\in A$ such that $\pi(x_i)=d_i{\cdot}a_i$. Since $A$ is locally cyclic and $3$-divisible, there exists an element $\alpha\in A$ such that $\{a_1,a_2,a_3\}\subseteq \{\alpha^{3n}:n\in\IZ\}$. Since the loop $X/Z$ is an Abelian group, the subgroup $D$ is normal in $X/Z$ and hence its preimage $N=\pi^{-1}[D]$ is a normal subloop of the Moufang loop $X$.  Choose any element $a\in\pi^{-1}(\alpha)$. By Lemma~\ref{l:central-by-abelian=>automorphism}, the conjugation map $(\cdot)^a:N\to N$ is an automorphism of the loop $N$. By Corollary~\ref{c:GagolaIII}, the loop $N a^{3\IZ}$ is associative. It follows from $\{\pi(x_1),\pi(x_2),\pi(x_3)\}\subseteq D\alpha^{3\IZ}$ and $Z\subseteq N$ that $\{x_1,x_2,x_3\}\subseteq Na^{3\IZ}$ and hence $(x_1\cdot x_2)\cdot x_3=x_1\cdot(x_2\cdot x_3)$, witnessing that the loop $X$ is associative.
\end{proof}

Now we are able to prove the promised associativity of $6$-divisible abelian-by-central Moufang loops.

\begin{theorem}\label{t:central-by-abelian} A central-by-abelian Moufang loop $X$ is associative if its centre $Z=\mathcal Z(X)$ is $2$-divisible and the quotient loop $X/Z$ is $6$-divisible.
 \end{theorem}

\begin{proof} Given any elements $x_1,x_2,x_3\in X$, we need to check that $(x_1\cdot x_2)\cdot x_3=x_1\cdot(x_2\cdot x_3)$. Let $Z\defeq\mathcal Z(X)$ be the centre of the Moufang loop and $\pi:X\to X/Z$ be the homomorphism to the quotient loop $X/Z$. Since $X$ is central-by-abelian, the quotient loop $X/Z$ is an Abelian group. Since the group $X/Z$ is $6$-divisible, for every $i\in\{1,2,3\}$, there exists a $6$-divisible locally cyclic subgroup $A_i\subseteq X/Z$ containing the image $\pi(x_i)$ of $x_i$, see Proposition~\ref{p:d-divisible-locally-cyclic}. Since $A_1,A_2,A_3$ are $6$-divisible subgroups of the Abelian group $X/Z$, so are the subgroups $A_1A_2$ and $A_1A_2A_3$.

Applying Lemma~\ref{l:induction-associative} (three times), we conclude that subloops $N_1\defeq\pi^{-1}[A_1]$, $N_2\defeq\pi^{-1}[A_1A_2]$ and $N_3\defeq\pi[A_1A_2A_3]$ are $2$-divisible and associative. 
Since $x_1,x_2,x_3$ are elements of the group $N_3$, $(x_1\cdot x_2)\cdot x_3=x_1\cdot(x_2\cdot x_3)$, witnessing that the Moufang loop $X$ is associative.
\end{proof}

\section{Abelian numbers}

\begin{definition}\label{d:Abelian-num} A positive integer number $n$ is called 
\begin{itemize}
\item  \defterm{square-free} if $n$ is not divisible by the square $p^2$ of a prime number $p$;
\item \defterm{cube-free} if $n$ is not divisible by the cube $p^3$ of a prime number $p$;
\item \defterm{Abelian} if it is cube-free and for any prime-power divisor $d$ of $n$, the numbers $d-1$ and $n$ are coprime.
\end{itemize}
\end{definition}

The main result of this section is the following characterization of Abelian numbers.

\begin{theorems}\label{t:Abelian-num} For a number $n\in\IN$, the following conditions are equivalent:
\begin{enumerate}
\item the number $n$ is Abelian;
\item every group of order $n$ is Abelian;
\item every Moufang loop of order $n$ is an Abelian group.
\end{enumerate}
\end{theorems}

\begin{proof} The implication $(3)\Ra(2)$ is trivial.
\smallskip

$(2)\Ra(1)$ Assume that every group of order $n$ is Abelian.
The following two claims imply that the number $n$ is Abelian.

\begin{claim} For any prime divisor $p$ of $n$, the numbers $p-1$ and $n$ are coprime.
\end{claim}

\begin{proof} Consider the $p$-element field $\IF_p$ and observe that its multiplicative group $\IF_p^*$ has cardinality $|\IF_p^*|=p-1$. To derive a contradiction, assume that the numbers $p-1$ and $n$ have a common prime divisor $q$. It follows that $p\ne q$ and hence $pq$ is a divisor of $n$. Let $C_m$ be the cyclic group of order $m\defeq n/(pq)$.  Since $q$ divides $p-1=|\IF_p^*|$, the multiplicative group $\IF_p^*$ contains a cyclic subgroup $C_q$ of order $q$. Then the semidirect product $\IF_p\rtimes C_q$ endowed with the operation $(x,y)\cdot(x',y')=(x+y{\cdot}x',y{\cdot} y')$ is a noncommutative group of order $pq$, and $(\IF_p\rtimes C_q)\times C_m$ is a noncommutative group of order $n=(pq)m$, witnessing that the number $n$ is not Abelian. 
\end{proof}

\begin{claim} If for some prime number $p$, the number $p^2$ divides $n$, then $p^3$ does not divide $n$ and $p^2-1$ is coprime with $n$. 
\end{claim}

\begin{proof} In the opposite case, we can find prime numbers $p,q$ such that $q$ divides $p(p^2-1)$ and $p^2q$ divides $n$. Let $m=n/(p^2q)$. Observe that the group $\Aut(C_p^2)\cong GL(2,p)$ of automorphisms of the elementary Abelian group $C_p^2=C_p\times C_p$ has cardinality $(p^2-1)(p^2-p)=(p-1)^2p(p+1)$. Since $q$ divides $p(p^2-1)$, it divides the order of the group $\Aut(C_p^2)$ and hence the group $\Aut(C_p^2)$ contains a cyclic subgroup $C_q$ of order $q$. Then the semidirect product $C_p^2\rtimes C_q$ endowed with the operation $(x,y)\cdot(x',y')=(x+y(x'),y\circ y')$ is a noncommutative groups of order $p^2q$ and the group $(C_p^2\rtimes C_q)\times C_m$ is a noncommutative group of order $(p^2q)m=n$, witnessing that the number $n$ is not Abelian.
\end{proof}
\smallskip

$(1)\Ra(3)$ By induction we shall prove that every Moufang loop of Abelian order is a commutative group. Assume that for some number $n\in\IN$ we have proved that every Moufang loop of Abelian order $<n$ is a  commutative group. Take any Moufang loop $X$ of Abelian order $n$. If $n$ is even, then $n\in\{2,4\}$, by the definition of an Abelian number. If $X$ is Boolean, then it is an Abelian group, by Proposition~\ref{p:Boolean+Moufang=>commutative-group}. If $X$ is not Boolean, then it contains a cyclic subgroup $C$ of order $|C|>2$. By Corollary~\ref{c:Bol=>Lagrange}, $|C|$ divides $|X|\in\{2,4\}$ and hence $|C|=4=|X|$ and $X=C$ is a cyclic group. So, assume that the Abelian number $n$ is odd.

By Theorem~\ref{t:odd=>nucleo-solvable}, the Moufang loop $X$ is nucleo-solvable and hence contains a nontrivial normal Abelian subloop $N$ such that $N\subseteq\mathcal N(X)$. Let $p$ be any prime divisor of $|N|$ and $P$ be the unique $p$-Sylow subgroup of the Abelian group $N$. 
Since $P\subseteq N\subseteq\mathcal N(X)$, for every $x,y\in X$ we have $(xP)y=x(Py)$ and $x(yP)=(x{\cdot}y)P$.
Also $N\subseteq\mathcal N(X)$ implies that for every $a\in X$, the conjugation map $(\cdot)^a:N\to N$, $(\cdot)^a:x\mapsto x^a$, is an automorphism of the group $N$. Since $P$ is a characteristic subgroup of the Abelian group, $P^a=P$ and hence $Pa=aP$ for all $a\in X$, witnessing that the $p$-Sylow subgroup $P$ is normal in the Moufang loop $X$.

Consider the quotient Moufang loop $X/P$. Since the order $|X/P|<n$ of the loop $X/P$ divides the Abelian number $n=|X|$, $|X/P|$ is an Abelian number, too. By the inductive hypothesis, the loop $X/P$ is an Abelian group. Therefore, the Moufang loop $X$ is central-by-abelian. Let $\pi:X\to X/P$ be the quotient homomorphism. Let $S_P$ be the Sylow $p$-subgroup the Abelian group $X/P$. Since $S_P$ is normal in the Abelian group $X/P$, the preimage $S\defeq\pi^{-1}[S_P]$ is a normal $p$-Sylow subloop of $X$. We claim that the Moufang loop $S$ is an Abelian group. This is clear if $P=S$. So assume that $P\ne S$. Since $n$ is cube-free, $|S|=p^2$ and both groups $P$ and $S/P$ are cyclic groups of order $p$. Then the Moufang loop $S$ is $2$-generated and hence is a group, by the diassociativity of the Moufang loop $X$. Then for every $s\in S$, the conjugation $(\cdot)^s:P\to P$, $(\cdot)^s:x\mapsto x^s$, is an automorphism of the cyclic group $P$ whose order divides the order of $s$, which divides $|S|=p^2$. Since the automorphism group of the cyclic group $P$ has cardinality $p-1$, the automorphism $(\cdot)^s$ is trivial and hence the group $S$ is commutative. 

By Theorem~\ref{t:Hall-Moufang}, there exists a subgroup $H\subseteq X$ of cardinality $|H|=|X|/|S|$. Since the numbers $|H|$ and $|S|$ are coprime, the intersection  $S\cap H$ is trivial and hence $X=S\rtimes H$. The group $S$ has odd order and hence is $2$-divisible. By Lemma~\ref{l:central-by-abelian=>automorphism}, for every $c\in X$, the conjugating map $(\cdot)^c:S\to S$ is an automorphism of the Abelian group $S$. 

Since the Abelian number $n=|X|$ is cube-free, the commutative $p$-Sylow subgroup $S$ of $X$ has cardinality $p$ or $p^2$ and hence $S$ is isomorphic to one of the groups $C_p$, $C_{p^2}$ or $C_p\times C_p$. Then the automorphism group $\Aut(S)$ of the group $S$ has cardinality $p-1$, $p(p-1)$ or $(p^2-1)(p^2-p)$, respectively. Since the number $n$ is Abelian, the number $|H|=n/|S|$ is coprime with the order $|\Aut(S)|\in \{p-1,p(p-1),(p^2-1)(p^2-p)\}$ of the automorphism group $\Aut(S)$ of $S$. 
This implies that for every  $h\in H$, the conjugating automorphism $(\cdot)^h:S\to S$, $(\cdot)^h:x\mapsto x^h$, is trivial, which means that $h{\cdot}s=s{\cdot}h$ for all $s\in S$. 

Let $T$ be the unique $3$-Sylow subgroup of the group $H$. 

\begin{claim} The subloop $S\rtimes T$ of $X$ is a commutative group.
\end{claim}

\begin{proof} This is clear if the group $T$ is trivial. So, assume that $T$ is not trivial. Since $S$ is a $p$-Sylow subgroup of $X$, $p\ne 3$. Since the Abelian number $n$ is cube-free, $|T|\in\{3,9\}$ and $|S|\in\{p,p^2\}$. Then there exist elements $s_1,s_2\in S$ that generate the Abelian group $S$ and there exist elements $t_1,t_2\in T$ that generate the Abelian group $T$. Since the Moufang loop $X$ is diassociative, the commuting elements $s_1$ and $t_1$ generate a commutative group $C_1$. Since the orders of $s_1$ and $t_1$ are coprime, the commutative group $C_1$ is cyclic and is generated by the element $s_1t_1$. By the same reason, the group $C_2$ generated by the elements $s_2,t_2$ is cyclic and is generated by the element $s_2t_2$. By the diassociativity of the Moufang loop $X$, the loop $G$, generated by the elements $s_1t_1$ and $s_2t_2$ is a group. Since $\{s_1,t_1,s_2,t_2\}\subseteq C_1\cup C_2\subseteq G$, the semidirect product $S\rtimes T$ coincides with the group $G$ and hence is associative.
Since all elements of the set $S\cup T$ commute, the group $S\rtimes T$ is commutative.
\end{proof}

If $T=H$, then $X=S\rtimes H=S\rtimes T$ is a commutative group. So, assume that $T\ne H$. Since the Abelian group $H$ is a direct sum of its Sylow subgroups and each Sylow subgroup of $H$ is either cyclic or elementary Abelian, the group $H$ contains two subgroups $L$ and $C$ such that $H=L\times C$, $T\subseteq L$ and $C$ is a cyclic group  of order $q$ or $q^2$ for some prime number $q\ne 3$. The inductive hypothesis ensures that the subloop $A=S\rtimes L$ of $X$ is an Abelian group. 
Since $|C|$ is not divisible by $3$, the cyclic group $C$ is $3$-divisible. Applying Corollary~\ref{c:GagolaIII}, we conclude that the semidirect product $X=S\rtimes H=S\rtimes(L\times C)=(S\rtimes L)\rtimes C=A\rtimes C$ is a group.

It remains to show that the group $X=A\rtimes C$ is Abelian. Let $c$ be a generator of the cyclic group $C$ and $\Pi$ be the set of all prime divisors of $|A|$.
For every $p\in\Pi$, let $S_p$ be the (unique) Sylow $p$-subgroup of the Abelian group $A$. The order $|A|$ of the group $A$ divides the number $n=|X|$ and hence is cube-free. Then for every $p\in\Pi$ the Sylow $p$-group $S_p$ is isomorphic to one of the groups: $C_p$, $C_{p^2}$ or $C_p\times C_p$, and then the automorphism group of $S_p$ has order $p-1$, $p(p-1)$ or $(p^2-1)(p^2-p)$. If $q=|C|\in\Pi$, then the Sylow $q$-subgroup $S_q$ of $A$ is cyclic of order $q$ (because $n$ is cube-free) and $|\Aut(S_q)|=q-1$ is not divisible by $|C|$. Observe that for every $p\in\Pi$ the Sylow $p$-subgroup $S_p$ is characteristic in the group $A$ and hence $S_p^c=S_p$. Taking into account that $n$ is an Abelian number, we conclude that the order of $c$ does not divide the order $|\Aut(S_q)|\in\{q-1,q^2-q,(q^2-1)(q^2-q)\}$ of the automorphism group $\Aut(S_p)$ and hence $x^c=x$ for all $x\in \bigcup_{p\in\Pi}S_p$. Since the set $\bigcup_{p\in\Pi}S_p$ generates the group $A$, $x^c=x$ for all $x\in A$. Therefore, the group $X$ is generated by the set $A\cup\{c\}$ consisting of commuting elements and hence $X$ is a commutative group.
\end{proof}

\begin{remark} The equivalence $(1)\Leftrightarrow(2)$ in Theorem~\ref{t:Abelian-num} was proved by Dickson \cite{Dickson1905} in 1905. Other proofs of this characterization can be found in the paper \cite{CheinRajah2000} of Chein and Rajah, and also in the Bachelor Thesis  \cite{LCrew} of Logan Crew (which contains a lot of informations on some other numbers: nilpotent, ordered Sylow, supersolvable). The implications $(1)\Ra(3)$ in Theorem~\ref{t:Abelian-num} can be also deduced from  far more general results of Ademola and Rajah \cite{AdemolaRajah2016}, \cite{Rajah2018}.
\end{remark}
\pagebreak

\section{Moufang numbers}

\begin{definition}\label{d:Moufang-num} A number $n$ is \defterm{Moufang} if every Moufang loop of order $n$ is associative.
\end{definition}

Even Moufang numbers were characterized by Chein and \index[person]{Rajah}Rajah\footnote{{\bf Andrew Rajah Balasingam Gnanaraj} is an Associate Professor in the School of Mathematical Sciences at Universiti Sains Malaysia (USM). He completed his B.Sc. in Mathematics at USM in 1993, then entered the M.Sc. program the same year before transferring into the Ph.D. program in 1994. He received his Ph.D. in 1997 with a dissertation entitled Which Moufang Loops are Associative, and subsequently joined the USM faculty as a lecturer in the same year. Promoted to Associate Professor in 2009, he has focused his research on determining the orders for which all Moufang loops must be associative and on constructing minimally nonassociative Moufang loops in cases where nonassociativity can be demonstrated.} in  \cite{CheinRajah2000}.

\begin{theorem}[Chein--Rajah, 2000] A number $n$ is Abelian iff the number $2n$ is Moufang.
\end{theorem}

\begin{proof} If a number $n$ is not Abelian, then there exists a noncommutative group $G$ of order $n$. By Theorem~\ref{t:M(G,2)-Moufang<=>G-group}, the Chein extension $M(G,2)$ of $G$ is a nonassociative Moufang loop of order $2n$, witnessing that the number $2n$ is not Moufang.
\smallskip

Now assume that the number $n$ is Abelian. If $n$ is even, then $n\in\{1,2,4\}$, by Theorem~\ref{t:Abelian-num}. Since every Moufang loop of order $<12$ is associative, the number $2n\in\{2,4,8\}$ is Moufang. So, assume that the Abelian number $n$ is odd. Given any Moufang loop $X$ of order $|X|=2n$, we should prove that $X$ is a group. By Proposition~\ref{p:Ali-Slaney}, the invertible loop $X$ contains a  subgroup $D$  of order $|D|=2$. By Theorem~\ref{t:2m=>normal-m}, the set $H$ of all elements of odd order is a normal subloop in $X$ of odd order $|H|=n$.
Then $X=H\rtimes D$. By Theorem~\ref{t:Abelian-num}, the Moufang loop $H$ is an Abelian group. 
 Since $n=|A|$ is odd and $|D|=2\notin 3\IZ$, the group $H$ is $2$-divisible and the cyclic group $D$ is $3$-divisible. By Corollary~\ref{c:GagolaIII}, the loop $X=H\rtimes D$ is a group.
\end{proof}

In contrast of even Moufang numbers, odd Moufang numbers do not have a reasonable arithmetic characterization yet.

\begin{problem}\label{prob:Moufang-num} Find an arithmetic characterization of Moufang numbers.
\end{problem}

We shall characterize Moufang numbers among prime-powers. 

\begin{remark}\label{r:23Moufang} The $16$-element Moufang loop of Cayley octonion units shows that the number $2^4$ is not Moufang. The nonassociative Moufang loop algebraizing the Hall plane in Example~\ref{ex:HTS} has $81$ elements and shows that the number $3^4$ is not Moufang. On the other hand, for all prime numbers $p\notin\{2,3\}$, the number $p^4$ is Moufang, which was proved by Fook Leong\footnote{{\bf Fook Leong} is a Malaysian mathematician associated with the School of Mathematical Sciences at Universiti Sains Malaysia, known for his contributions to finite Moufang loop theory. His research activity was concentrated in the 1990s, during which he produced several technical monographs and reports at USM, including advanced work on the structure and classification of Moufang loops of small and mixed orders. He collaborated with leading algebraists at USM, among them Teh Pang Eng and Andrew Rajah, and helped develop a systematic approach to understanding associativity conditions and extensions in Moufang loops. His publications and internal research reports form a significant part of the Malaysian school of non-associative algebra during that period.} in \cite{Leong1974}.
\end{remark}

\begin{theorems}[Leong, 1974]\label{t:Leong} For every prime number $p\ge 5$, the number $p^4$ is Moufang.
\end{theorems}

\begin{proof} We should prove that every Moufang loop $X$ of order $|X|=p^4$ is associative. By Theorem~\ref{t:prime-power=>center}, the loop $X$ has non-trivial centre $Z=\mathcal Z(X)$. Since $p\ge 5$, the group $Z$ is $2$-divisible and the quotient loop $X/Z$ is $6$-divisible. 

If the quotient loop $X/Z$ is an Abelian group, then the loop $X$ is central-by-Abelian. By Theorem~\ref{t:central-by-abelian}, the $6$-divisible central-by-abelian Moufang loop $X$ is associative. 

So, assume that the quotient loop $Y=X/Z$ is not an Abelian group. Since the number $p^2$ is Abelian, the loop $Y$ has cardinality $|X/Z|=p^3$ and the centre $Z$ is a cyclic group of order $p$. Let $z$ be a generator of the cyclic group $Z$. By Theorem~\ref{t:prime-power=>center}, the Moufang loop $Y$ has nontrivial centre $C$. Since the order $|Y/C|\in\{1,p,p^2\}$ is Abelian, the quotient loop $Y/C$ is an Abelian group. By Theorem~\ref{t:central-by-abelian}, the central-by-Abelian Moufang loop $Y$ is a group. Since the group $Y$ is not commutative, there exists an element $a\in Y\setminus C$. Then the subgroup $N$ generated by the set $C\cup\{a\}$ is a commutative subgroup of $Y$. Since the quotient group $Y/C$ is Abelian, the subgroup $N$ is normal in $Y$. Since the group $Y$ is not commutative, $N\ne Y$, which implies $|C|=p$, $|N|=p^2$ and $|Y|=p^3$. Let $A\subseteq Y$ be the cyclic subgroup generated by the element $a$. Let $c$ be a generator of the cyclic group $C$, and $b$ be any element in $Y\setminus A$. Since the group $Y$ is generated by the elements $a,b,c$ and $c\in C=\mathcal Z(Y)$, the elements $a,b$ do not commute. Since the order $p$ of the element $b$ does not divide the order $p-1$ of the automorphism group of the cyclic group $A$, the element $a^b\in N$ does not belong to the cyclic group $A$ and hence the group $N$ is generated by the elements $a,a^b$. Moreover, the group $Y$ is generated by the elements $a,b$. Then the group $X$ is generated by the elments $z,a,b$. Since $z\in \mathcal Z(X)\subseteq\mathcal N(X)$, the loop $X$ is associative, by the Moufang Theorem~\ref{t:Moufang-xyz}.
\end{proof}

\begin{remarks}\label{r:p5Moufang} In \cite{Wright1965} Wright constructed a non-associative Moufang loop of of order $p^5$ for every prime number $p\ge 5$. In \cite{NagyVal2007} Nagy and Valsecchi proved that up to an isomorphism, there are exactly four non-associative Moufang loops of order $p^5$. These four Moufang loops  have presentation 
$$
\begin{aligned}
\langle w,g_1,g_2,g_3,g_4\mid& g_2^p=g_3^p=g_4^p=1,\; g_1^p=g_4^\alpha,\;w^p=g_4^\beta,\\
& [g_1,g_4]=[g_2,g_3]=[g_2,g_4]=[g_1,w]=[g_2,w]=[g_n,w]=1,\\
&[g_2,g_1]=g_3,\;[g_1,g_3]=g_4,\;[g_3,w]=g_4^6,\;[g_2,g_1,w]=g_4\rangle
\end{aligned}
$$where $(\alpha,\beta)\in\{(0,0),(0,1),(1,0),(1,6)\}$.
\end{remarks}

\begin{remarks} In \cite{Rajah2001} for every odd prime number $p$ and any odd prime number $q$ dividing $p-1$, Andrew Rajah constructed a non-associative Moufang loop of order $p^3q$. Stephen Gagola III noticed \cite{GagolaIII2014} that for $q>3$, the Rajah loop is the semidirect product $\IF_p^3\rtimes C_q$ where the generator of the cyclic group $C_q$ acts on the Abelian group $\IF^3_p$ by the semi-automorphism 
$$
\left(\begin{array}{ccc}1&x&z\\
0&1&y\\
0&0&1
\end{array}
\right)\mapsto \left(\begin{array}{ccc}1&xk^{-1}&zk+xy\tfrac{k^{-2}-k}2\\
0&1&yk^{-1}\\
0&0&1
\end{array}
\right)$$
with $k$ being any element of order $q$ in the multiplicative group $\IF^*_p$ of the field $\IF_p$.
\end{remarks} 

Theorem~\ref{t:Leong} and Remarks~\ref{r:23Moufang} and \ref{r:p5Moufang} imply the following characterization of Moufang prime-powers.

\begin{theorems}\label{t:Moufang-prime-power} For a prime number $p$ and a positive integer $n$, the number $p^n$ is Moufang if and only if $p\in\{2,3\}$ and $n\le 3$ or $p\ge 5$ and $n\le 4$.
\end{theorems}

The latest information on partial solutions of Problem~\ref{prob:Moufang-num} can be found in the paper \cite{Rajah2018} of Rajah and Ademola. Their paper \cite{AdemolaRajah2016} contains the following interesting theorem.

\begin{theorems}[Ademola--Rajah, 2016] The cube $n^3$ of every square-free Abelian number $n$ is a Moufang number.
\end{theorems}

\section{Abelian and Moufang numbers in Geometry} 

In this section we present some applications of Abelian and Moufang numbers in Linear Geometry.


Definitions~\ref{d:Abelian-num}, \ref{d:Moufang-num}, Theorems~\ref{t:+inversive}, \ref{t:add-ass<=>}, \ref{t:add-com<=>}, \ref{t:inversive-puls}, \ref{t:associative-puls<=>}, \ref{t:commutative-puls<=>} imply the following two characterizations.

\begin{proposition}\label{p:Moufang=>inversive<=>ass} If the order of a Playfair liner $X$ is a Moufang number, then
\begin{enumerate}
\item $X$ is inversive-plus if and only if $X$ is associative-plus;
\item $X$ is inversive-puls if and only if $X$ is associative-puls.
\item $X$ is inversive-add if and only if $X$ is associative-add.
\end{enumerate}
\end{proposition}

\begin{proposition}\label{p:inv-add<=>comm-add} If the order of a Playfair liner $X$ is an Abelian number, then
\begin{enumerate}
\item $X$ is inversive-plus if and only if $X$ is commutative-plus;
\item $X$ is inversive-puls if and only if $X$ is commutative-puls;
\item $X$ is inversive-add if and only if $X$ is commutative-add.
\end{enumerate}
\end{proposition}

\begin{proposition} Assume that the order of a Playfair liner $X$ is Moufang and prime-power. If $X$ is inversive-plus or inversive-puls, then $X$ is $\partial$-translation and $|X|_2\in\{p,p^2,p^3,p^4\}$ for some prime number $p$.
\end{proposition}

\begin{proof} Theorem~\ref{t:Moufang-prime-power} ensures that $|X|_2\in\{p,p^2,p^3,p^4\}$ for some prime number $p$. Since $|X|_2$ is a Moufang number, the inversivity-plus or inversivity-puls of $X$ implies the associativity-plus or associatovity-puls of $X$. Since the order $|X|_2$ is a prime-power, we can apply Theorems~\ref{t:add-ass=>partialT} and \ref{t:ass-puls=>partialT} and conclude that the Playfair liner $X$ is $\partial$-translation.
\end{proof}

\begin{proposition} Assume that the order of a Playfair liner $X$ is Abelian. If $X$ is inversive-add, then $X$ is Thalesian of order $|X|_2\in\{p,p^2\}$ for some prime number $p$.
\end{proposition}

\begin{proof}  If $X$ is inversive-add, then $X$ is commutative-add (by Proposition~\ref{p:inv-add<=>comm-add}) and Boolean or Thalesian, by Theorem~\ref{t:commutative-add<=>}. If $X$ is Boolean then $|X|_2=2^n$ for some number $n\ge 2$, by Corollary~\ref{c:Boolean-order}. Since the number $|X|_2$ is Abelian, it is cube-free, which implies $|X|_2=4$ and hence $X$ is Desarguesian and Thalesian, by Corollary~\ref{c:4-Pappian}. If $X$ is Thalesian, then $|X|_2=p^n$ for some $n\in\IN$, by Theorem~\ref{t:cardinality-pD}. Since the Abelian number $|X|_2=p^n$ is cube-free, $n\le 2$ and hence $|X|_2\in\{p,p^2\}$.
\end{proof} 

\section{Geometric and arithmetic properties of affine spaces}
Results of Chapters~\ref{ch:PlusAff}, \ref{ch:PulsAff}, \ref{ch:dot-aff}, \ref{ch:invertible-add}, \ref{ch:commutative-add}, \ref{ch:AM-numbers} imply that for every affine space, the following implications hold:
$$
\xymatrix@C=45pt@R=20pt{
\mbox{commutative-dot}\ar@{<=>}[r]\ar@{=>}[d]&\mbox{Pappian}\ar@{=>}[d]\\
\mbox{associative-dot}\ar@{<=>}[r]\ar@{=>}[d]&\mbox{Desarguesian}\ar@/_12pt/_{\mbox{\tiny finite}}[u]\ar@{=>}[d]\\
\mbox{inversive-dot}\ar@{<=>}[r]\ar@{=>}[d]\ar@/_12pt/_{\mbox{\tiny Boolean}}[u]&\mbox{Moufang}\ar@/_12pt/_{\mbox{\tiny finite}}[u]\ar@{<=>}[r]\ar@{=>}[d]&\mbox{shear}\ar@{=>}[d]\\
\mbox{invertible-dot}\ar@{=>}[dd]&\mbox{Thalesian}\ar@/_12pt/_{\mbox{\tiny prime}}[u]\ar@{<=>}[r]\ar@{=>}[d]&\mbox{translation}\ar@{=>}[d]\\
&\mbox{Thalesian}\atop\mbox{or Boolean}\ar@/_12pt/_{\mbox{\tiny prime}}[u]\ar@{<=>}[d]&\mbox{$\partial$-translation}\ar@/_12pt/_{\mbox{\tiny prime}}[u]\\
\mbox{invertible-plus}&\mbox{commutative-add}\ar@{=>}[ld]\ar@{=>}[rd]\ar@{=>}[r]&\mbox{field-elementary}\\
\mbox{commutative-plus}\ar@{=>}[u]\ar@{=>}[d]\ar^{\mbox{\tiny finite}}_{\mbox{\tiny invertible-puls}}[r]&\mbox{prime-power}\atop\mbox{$\partial$-translation}\ar@{=>}[d]&\mbox{commutative-puls}\ar@{=>}[d]\ar_{\mbox{\tiny finite}}^{\mbox{\tiny invertible-plus}}[l]\\
\mbox{associative-plus}\ar@/_12pt/_{\mbox{\tiny Abelian}\atop\mbox{\tiny order}}[u]\ar@{=>}[d]\ar^{\mbox{\tiny prime-power}}[r]&\mbox{$\partial$-translation}&\mbox{associative-puls}\ar@/_12pt/_{\mbox{\tiny Abelian}\atop\mbox{\tiny order}}[u]\ar_{\mbox{\tiny prime-power}}[l]\ar@{=>}[d]\\
\mbox{inversive-plus}\ar@/_12pt/_{\mbox{\tiny Moufang}\atop\mbox{\tiny order}}[u]\ar@{=>}[d]\ar^{\mbox{\tiny finite}}_{\mbox{\tiny invertible-puls}}[r]&\mbox{order not 60}&\mbox{inversive-puls}\ar@/_12pt/_{\mbox{\tiny Moufang}\atop\mbox{\tiny order}}[u]\ar@{=>}[d]\\
\mbox{invertible-plus}&\mbox{invertible-add}\ar@{=>}[l]\ar@{=>}[r]\ar@{<=>}[d]&\mbox{invertible-puls}\\
\mbox{Pappian}&\mbox{$\IZ$-elementary}\ar_{\mbox{\tiny finite and}}^{\mbox{\tiny minimal}}[l]\ar@/_12pt/^{\mbox{\tiny finite}}[r]&\mbox{field-elementary}\ar@{=>}[l]
}
$$

\chapter{Fano liners}\label{ch:Fano}

This chapter is devoted to \index[person]{Fano}Fano\footnote{{\bf Gino Fano} (1871 -- 1952) was an Italian mathematician, best known as the founder of finite geometry. Fano made various contributions on projective and algebraic geometry. His work in the foundations of geometry predates the similar, but more popular, work of David Hilbert by about a decade.} and \index[person]{Boole}Boolean\footnote{{\bf George Boole} (1815--1864) was a British mathematician and logician, best known as the founder of Boolean algebra, the algebraic system that underlies modern digital logic and computer science. Largely self-taught, Boole became the first professor of mathematics at Queen’s College, Cork, in Ireland. His work on the algebra of logic, published in ``The Mathematical Analysis of Logic'' (1847) and ``An Investigation of the Laws of Thought'' (1854), laid the foundations for symbolic logic and influenced later developments in set theory, probability, and theoretical computer science. Boole’s ideas are central to the design of digital circuits and programming languages.} liners. In the last Section~\ref{s:Gleason-Fano}, we prove the famous theorem of Gleason on the Pappianity of finite Fano projective planes and also its generalization to finite proaaffine Fano liners. The results of the Sections~\ref{s:Fano1}--\ref{s:Fano3}  are taken from the Bachelor Thesis \cite{Skygar} of Oksana Skygar. 

\section{Fano and Boolean liners}\label{s:Fano1}

\begin{definition} A liner $X$ is \index{Fano liner}\index{liner!Fano}\defterm{Fano} if for any distinct points $a,b,c,d\in \Pi$, the set 
$$D\defeq(\Aline ab\cap\Aline cd)\cup(\Aline bc\cap\Aline ad)\cup(\Aline ac\cap\Aline bd)$$ has rank $\|D\|\in\{0,2\}$.
\end{definition}

Fano liners can be equivalently defined by requiring all quadrangles to be Fano. 

\begin{definition} A \index{quadrangle}\defterm{quadrangle} in a liner $X$ is a quadruple $abcd$ of distinct points $a,b,c,d\in X$ such that $\|\{a,b,c,d\}\|=3$ and $\|\{x,y,z\}\|=3$ for any distinct points $x,y,z\in\{a,b,c,d\}$. A quadrangle $abcd$ in a liner $X$ is called \index{complete quadrangle}\index{quadrangle!complete}\defterm{complete} if the set $$D\defeq(\Aline ab\cap\Aline cd)\cup(\Aline bc\cap\Aline ad)\cup(\Aline ac\cap\Aline bd)$$  has cardinality $|D|=3$.
\end{definition}

\begin{exercise} Show that every quadrangle in a projective liner is complete.
\end{exercise}

\begin{definition} A quadrangle $abcd$ in a liner is defined to be \index{Fano quadrangle}\index{quadrangle!Fano}\defterm{Fano} if the set $$D\defeq(\Aline ab\cap\Aline cd)\cup(\Aline bc\cap\Aline ad)\cup(\Aline ac\cap\Aline bd)$$ has rank $\|D\|\in\{0,2\}$.
\end{definition}

\begin{proposition} A liner $X$ is Fano if and only if every quadrangle in $X$ is Fano. 
\end{proposition}

\begin{proof}   The ``only if'' part of this characterization follows immediately from the definitions. To prove the ``if'' part, assume that every quadrangle in a liner $X$ is Fano. To prove that $X$ is Fano, take any distinct points $a,b,c,d\in X$ and consider the set $D\defeq(\Aline ab\cap\Aline cd)\cup(\Aline ac\cap\Aline bd)\cup(\Aline ad\cap\Aline bc)$. We have to prove that $\|D\|\in\{0,2\}$.

If $\|\{a,b,c,d\}\|=4$, then the set $$D=(\Aline ab\cap\Aline cd)\cup(\Aline ac\cap\Aline bd)\cup(\Aline ad\cap\Aline bc)$$ is empty and hence $\|D\|=0$. If $\|\{a,b,c,d\}\|=2$, then the set $D=\overline{\{a,b,c,d\}}$ has rank $\|D\|=\|\{a,b,c,d\}\|=2$. So, assume that $\|\{a,b,c,d\}\|=3$. If $\|\{x,y,z\}\|=3$ for all distinct points $x,y,z\in \{a,b,c,d\}$, then $abcd$ is a quadrangle in $X$ and hence $\|D\|\in\{0,2\}$, by the assumption. So, assume that $\|\{x,y,z\}\|<3$ for some distinct points $x,y,z\in\{a,b,c,d\}$. We lose no generality assuming that $\{x,y,z\}=\{a,b,c\}$. Since $\|\{a,b,c,d\}\|=3>2=\|\{a,b,c\}\|$, the point $d$ does not belong to the line $\overline{\{a,b,c\}}$. Then the set $$D=(\Aline ab\cap\Aline cd)\cup(\Aline ac\cap\Aline bd)\cup(\Aline ad\cap\Aline bc)=\{c\}\cup\{b\}\cup\{a\}$$ has rank $\|D\|=\|\{c,b,a\}\|=2$.
\end{proof}

Recall that a liner is \defterm{Steiner} if each of its lines contains exactly three points.

\begin{example}\label{ex:Fano3} Every Steiner projective liner $X$ is Fano.
\end{example}

\begin{proof} Given any quadrangle $abcd$ in $X$, we should prove that the set $$D\defeq(\Aline ab\cap\Aline cd)\cup(\Aline ac\cap\Aline bd)\cup(\Aline ad\cap\Aline bc)$$ has rank $0$ or $2$. Since $abcd$ is a quadrangle, its flat hull $\Pi\defeq\overline{\{a,b,c,d\}}$ in $X$ is a plane. By Corollary~\ref{c:Steiner-projective<=>}, the Steiner projective plane $\Pi$ has cardinality $|\Pi|=7$. 

Since $\Pi$ is projective, there exist points 
$x\in \Aline ab\cap\Aline cd$, $y\in \Aline ac\cap\Aline bd$ and $z\in \Aline ad\cap\Aline bc$.  It is easy to see that all points $a,b,c,d,x,y,z$ are pairwise distinct and moreover, $\{a,b,c,d\}\cap\Aline xy=\varnothing$. Taking into account that $|\Pi|=7$, we conclude that $\Pi=\{a,b,c,d,x,y,z\}$ and $\Aline xy=\Pi\setminus\{a,b,c,d\}=\{x,y,z\}$ witnesing that $\|D\|=\|\{x,y,z\}\|=\|\Aline xy\|=2$.

\begin{picture}(100,75)(-200,-25)
{\linethickness{1pt}
\put(0,0){\color{orange}\line(0,1){33}}
\put(0,0){\color{orange}\line(0,-1){16}}
\put(0,0){\color{brown}\line(7,4){14}}
\put(0,0){\color{brown}\line(-7,-4){28}}
\put(0,0){\color{olive}\line(7,-4){28}}
\put(0,0){\color{olive}\line(-7,4){14}}

\put(28,-16){\color{magenta}\line(-4,7){28}}
\put(-28,-16){\color{cyan}\line(4,7){28}}
\put(28,-16){\color{blue}\line(-1,0){56}}

\put(0,0.3){\color{red}\circle{32}}
}

\put(0,0){\circle*{3}}
\put(-7,-11){$d$}
\put(14,8){\circle*{3}}
\put(16,9){$z$}
\put(28,-16){\circle*{3}}
\put(31,-19){$y$}
\put(-14,8){\circle*{3}}
\put(-20,9){$b$}
\put(-28,-16){\circle*{3}}
\put(-36,-20){$a$}
\put(0,33){\circle*{3}}
\put(-3,36){$x$}
\put(0,-16){\circle*{3}}
\put(-3,-24){$c$}
\end{picture}
\end{proof}

Let us recall that a liner $X$ is \defterm{Boolean} if for any distinct points $a,b,c,d\in X$ with $\Aline ab\cap\Aline cd=\varnothing=\Aline bc\cap\Aline ad$ we have $\Aline ac\cap\Aline cd=\varnothing$. 

\begin{proposition}\label{p:Fano=>Boolean} Every Fano liner is Boolean.
\end{proposition}

\begin{proof} Let $X$ be a Fano liner. Given any distinct points $a,b,c,d\in X$ with $\Aline ab\cap\Aline cd=\varnothing=\Aline bc\cap\Aline ad$ we need to show that $\Aline ac\cap\Aline cd=\varnothing$. Since $X$ is Fano, the set 
$$D\defeq(\Aline ab\cap\Aline cd)\cup(\Aline bc\cap\Aline ad)\cup(\Aline ac\cap\Aline bd)=\Aline ac\cap\Aline bd$$has rank $\|D\|\in\{0,2\}$. If $\|D\|=0$, then $D=\Aline ac\cap\Aline bd$ is empty and we are done. So, assume that $\|D\|=2$. In this case $\Aline ac=\Aline bd$ and hence $\Aline ab=\Aline  ac=\Aline bd=\Aline cd$, which contradicts our assumption. 
\end{proof}

\begin{remark} Proposition~\ref{p:Fano=>Boolean} cannot be reversed: just take any non-Pappian Thalesian Playfair liner of order $2^n$ for some $n\ge 4$. It is Boolean, by Proposition~\ref{p:Boolean<=>1+1=0} and not Fano, by Theorem~\ref{t:proaffFano=>Pappian}. 
\end{remark}  

The following theorem shows that a projective liner is Fano if and only if it is everywhere Boolean.

\begin{theorem}\label{t:Fano<=>} For a projective liner $X$, the following conditions are equivalent:
\begin{enumerate}
\item $X$ is Fano;
\item for every flat $H$ in $X$, the subliner $X\setminus H$ of $X$ is Fano;
\item for every flat $H$ in $X$, the subliner $X\setminus H$ of $X$ is Boolean;
\item for every hyperplane $H$ in $X$, the subliner $X\setminus H$ of $X$  is Boolean.
\end{enumerate}
\end{theorem}

\begin{proof} $(1)\Ra(2)$ Assume that a projective liner $X$ is Fano. Let $H$ be any flat in $X$. To prove that the subliner $X\setminus H$ of $X$ is Fano, choose any quadrangle $abcd$ in $X\setminus H$. Since $\|\{a,b,c,d\}\|=3$ and $X$ is projective, there exist points $x\in \Aline ab\cap\Aline cd$, $y\in\Aline ac\cap\Aline bd$ and $z\in \Aline ad\cap\Aline bc$.  Since the liner $X$ is Fano, the set $$D\defeq(\Aline ab\cap\Aline cd)\cup(\Aline ac\cap\Aline bd)\cup(\Aline ad\cap\Aline bc)=\{x, y, z\}$$ has rank $\|D\|=2$ in $X$. We need to check that the set $D\setminus H$ has rank $\|D\setminus H\|\in\{0,2\}$ in the liner $X\setminus H$.

    If $|\{x,y,z\}\cap H|\le1$, then $|D\setminus H|\ge 2$ and hence $\|D\setminus H\|=2$. So, assume that $|\{x,y,z\}\cap H|\ge 2$. We lose no generality assuming that $x,y\in H$. Since $\|\{x,y,z\}\|=2$ and $H$ is flat, $z\in\Aline xy\subseteq H$ and hence $\|D\setminus H\|=\|\varnothing\|=0$.        
\smallskip

The implication $(2)\Ra(3)$ follows from Proposition~\ref{p:Fano=>Boolean}, and $(3)\Ra(4)$ is trivial.
\smallskip

$(4)\Ra(1)$ Assume that for every hyperplane $H$ in $X$ the subliner $X\setminus H$ is Boolean. To prove that $X$ is Fano, take any quadrangle $abcd$ in $X$. Consider the plane $\Pi\defeq\overline{\{a,b,c,d\}}$ in $X$. Since $X$ is projective, there exist points  $x\in \Aline ab\cap\Aline cd$, $y\in\Aline ac\cap\Aline bd$ and $z\in \Aline ad\cap\Aline bc$. It is easy to see that $\{a,b,c,d\}\cap\Aline xz=\varnothing$. By the Kuratowski--Zorn Lemma, there exists a maximal flat $H\subseteq X$ such that $x,z\in H$ and $\{a,b,c,d\}\cap H=\varnothing$. We claim that $H$ is a hyperplane in $X$. Since the projective space $X$ is strongly regular and hence has the Exchange Property, it suffices to check that the flat $F\defeq \overline{H\cup\{a\}}$ coincides with $X$. Observe that $\{b,d\}\subseteq \Aline xa\cup\Aline za\subseteq F$ and $c\in \Aline bz\subseteq F$.  Assuming that $F \ne X$, we can choose a point $p\in X\setminus F$ and consider the flat $H'\defeq\overline{H\cup\{p\}}$. The strong regularity of the projective liner $X$ ensures that $H'=\bigcup_{h\in H}\Aline hp$. Then for every $h'\in H'\setminus H$ we can find a point $h\in H$ such that $h'\in\Aline ph$. It follows from $h\in H\subseteq F$ and $p\notin F$ that $\Aline ph\cap F=\{h\}$ and hence $h'\notin F$. Then  $\{a,b,c,d\}\cap H'\subseteq (H'\cap F)\setminus H=H\setminus H=\varnothing$, which contradicts the maximality of the flat $H$. This contradiction shows that the flat $H$ is a hyperplane in $X$.  By (4), the liner $B\defeq X\setminus H$ is Boolean, and hence
$$(\Aline ab\cap B\cap\Aline cd)\cup(\Aline ad\cap B\cap\Aline bc)=B\cap\{x,z\}=\varnothing$$ implies $\Aline ac\cap B\cap \Aline bd=\varnothing$. Then   
$\{z\}\cap B=(\Aline ac\cap\Aline bd)\cap B=\varnothing$ and $\{x,y,z\}\subseteq \Pi\setminus B=\Pi\cap H$. Since the projective plane $\Pi$ is ranked, the intersection $\Pi\cap H$ is a line in $\Pi$ and hence $\|\{x,y,z\}\|=\|\Pi\cap H\|=2$, witnessing that the projective liner $X$ is Fano. 
\end{proof}

For $4$-long Desarguesian projective spaces we can extend the characterization Theorem~\ref{t:Fano<=>} as follows.

\begin{proposition}\label{p:DesFano<=>} For a $4$-long Moufang projective space $X$, the following conditions  are equivalent:
\begin{enumerate}
\item $X$ is Fano;
\item for some hyperplane $H$ in $X$, the subliner $X\setminus H$ of $X$ is Boolean;
\item for some hyperplane $H$ in $X$, the affine subliner $X\setminus H$ of $X$ contains a Boolean parallelogram;
\item $1+1=0$ in the scalar corps $\IR_X$ of $X$. 
\end{enumerate}
\end{proposition}

\begin{proof} The implication $(1)\Ra(2)$ follows from the implication $(1)\Ra(4)$ of Theorem~\ref{t:Fano<=>}.
\smallskip

$(2)\Ra(3)$ Assume that for some hyperplane $H$ in $X$, the subliner $X\setminus H$ of $X$ is Boolean. By Proposition~\ref{p:projective-minus-hyperplane}, the subliner $A\defeq X\setminus H$ is affine and regular. Since the projective liner $X$ is $4$-long, the affine subliner $A$ of $X$ is $3$-long and hence $A$ is a Playfair liner, by Theorem~\ref{t:Playfair<=>}. By Theorem~\ref{t:parallelogram3+1}, the Playfair liner $A$ contains a parallelogram. Since the liner $A$ is Boolean, this parallelogram is Boolean.
\smallskip

$(3)\Ra(4)$ Assume that for some hyperplane $H$ in $X$, the affine subliner $A\defeq X\setminus H$ of $X$ contains a Boolean parallelogram.  Since the projective space $X$ is $4$-long and Moufang, the affine subliner $A$ is $3$-long, affine, regular, and Thaleasian, according to Proposition~\ref{p:projective-minus-hyperplane} and Theorem~\ref{t:Moufang<=>everywhere-Thalesian}. By
Proposition~\ref{p:Boolean<=>1+1=0} (proved below), $1+1=0$ in the scalar corps $\IR_A$ of the Thalesian affine space $A$. Since the scalar corps $\IR_X$ of the projective space $X$ is isomorphic to $\IR_A$, $1+1=0$ in $\IR_X$.
\smallskip

$(4)\Ra(1)$ Assume that $1+1=0$ in the scalar corps $\IR_X$ of $X$. To prove that $X$ is Fano, it suffices to show that for every hyperplane $H$, the affine subspace $A\defeq X\setminus H$ is Boolean.  By Proposition~\ref{p:projective-minus-hyperplane}  and Theorem~\ref{t:Moufang<=>everywhere-Thalesian}, the affine space $A$ is  Thalesian. Since the scalar corps $\IR_X$ of the Moufang projective space $X$ is isomorphic to the scalar corps $\IR_A$ of the Thalesian affine space $A$, $1+1=0$ in $\IR_A$. By Proposition~\ref{p:Boolean<=>1+1=0}, the Thalesian affine space $A=X\setminus H$ is Boolean.
\end{proof}

\begin{proposition}\label{p:Boolean<=>1+1=0} For a Thalesian affine space $X$, the following conditions are equivalent:
\begin{enumerate}
\item $X$ is Boolean;
\item $X$ contains a Boolean parallelogram;
\item $1+1=0$ in the scalar corps $\IR_X$ of $X$.
\end{enumerate}
\end{proposition}

\begin{proof} By Theorem~\ref{t:paraD=>RX-module}, the Thalesian affine space $X$ carries the structure of a vector space over the corps of scalars $\IR_X$, so we can apply Vector Calculus in $X$.
\smallskip

The implication $(1)\Ra(2)$ is trivial (more precisely, it follows from Theorem~\ref{t:parallelogram3+1} and $\|X\|\ge 3$). 
\smallskip

$(2)\Ra(3)$ Assume that $X$ contains a Boolean parallelogram $abcd$. By the definition of a vector, the equalities $$\varnothing=\Aline ab\cap\Aline cd=\Aline ac\cap\Aline bd=\Aline ad\cap\Aline bc$$ imply $\overvector {ac}=\overvector {bd}=\overvector{ca}$. Then $(1+1)\cdot\overvector{ac}=\overvector{ac}+\overvector{ac}=\vec 0$ and $1+1=0$ in the scalar corps $\IR_X$ of the affine space $X$.
\smallskip

$(3)\Ra(1)$ Assume that $1+1=0$ in the scalar corps $\IR_X$ of $X$. To prove that $X$ is Boolean choose any parallelogram $abcd$ in $X$. Let $\Pi\defeq\overline{\{a,b,c,d\}}$ be the plane of the parallelogram $abcd$. It follows from $\Aline ab\parallel \Aline cd$ and $\Aline ad\parallel\Aline cd$ that $\overvector{ba}=\overvector{cd}$ and $\overvector{ad}=\overvector{bc}$.
Since $1+1=0$ in the scalar corps $\IR_X$, $$\overvector{ab}+\overvector{ab}=(1+1)\cdot\overvector{ab}=0\cdot\overvector{ab}=\vec{\boldsymbol 0}=\overvector{ab}+\overvector{ba},$$
and $\overvector{ab}=\overvector{ba}$. Then $$\overvector{ac}=\overvector{ab}+\overvector{bc}=\overvector{ba}+\overvector{bc}=\overvector{bc}+\overvector{cd}=\overvector{bd}$$ and hence $\Aline ac\parallel \Aline bd$ and $\Aline ac\cap\Aline bd=\varnothing$, witnessing that the parallelogram $abcd$ is Boolean.
\end{proof}



\begin{theorem}\label{t:Boolean<=>forder} For an affine space $X$ of finite order, the following conditions are equivalent:
\begin{enumerate}
\item $X$ is Boolean;
\item $X$ is di-Boolean, invertible-plus, and $|X|_2$ is even;
\item $X$ is di-Boolean, invertible-puls, and $|X|_2$ is even;
\item $X$ is commutative-plus, invertible-puls, and $|X|_2$ is even;
\item $X$ is commutative-puls, invertible-plus, and $|X|_2$ is even;
\item $X$ is inversive-plus, invertible-puls, and $|X|_2=2^n$ for some $n\ge 2$;  
\item $X$ is inversive-puls, invertible-plus, and $|X|_2=2^n$ for some $n\ge 2$.
\end{enumerate}
\end{theorem}

\begin{proof} $(1)\Ra(2,3)$ If $X$ is Boolean, then it is di-Boolean, invertible-add, and has even order, by Propositions~\ref{p:Boole=>bi-di-by-Boole}, Theorem~\ref{t:add-com=>add-ass}, \ref{t:Boolean=>commutative-plus},  \ref{t:com-puls=>ass-puls}, \ref{t:Boolean=>commutative-puls}, and Corollary~\ref{c:Boolean-order}.
\smallskip

$(2,3)\Ra(1)$ Assume that $X$ is di-Boolean and has even order. If $X$ is invertible-plus or invertible-puls, then $X$ contains a Boolean parallelogram, by Corollaries~\ref{c:inv-plus=>Boolean-parallelogram2}  and \ref{c:inv-puls=>Boolean-parallelogram2}. By Theorem~\ref{t:Boolean<=>di-Boolean+}, the di-Boolean Playfair plane $X$ is Boolean.
\smallskip

$(1)\Ra(4,5)$ If $X$ is Boolean, then $X$ is commutative-add and hence invertible-add, by Theorems~\ref{t:Boolean=>commutative-plus} and \ref{t:Boolean=>commutative-puls}. By Corollary~\ref{c:Boolean-order}, $|X|_2=2^n$ for some $n\ge 2$. 
\smallskip

$(4)\Ra(6)$ Assume that $X$ is commutative-plus, invertible-puls and $|X|_2$ is even. 
By Proposition~\ref{p:commutative-pls=>prime-power}, $|X|_2=p^n$ for some prime number $p$ and some $n\in\IN$. Since $X$ is $3$-long and $|X|_2$ is even, $p=2$ and $n\ge 2$.
\smallskip 

$(5)\Ra(7)$ Assume that $X$ is commutative-puls, invertible-plus and $|X|_2$ is even. 
By Proposition~\ref{p:commutative-pls=>prime-power},  $|X|_2=p^n$ for some prime number $p$ and some $n\in\IN$. Since $X$ is $3$-long and $|X|_2$ is even, $p=2$ and $n\ge 2$.
\smallskip

$(6)\Ra(1)$ Assume that $X$ is inversive-plus, invertible-puls, and $|X|_2=2^n$ for some $n\ge 2$. Let $R$ be any ternar of $X$. By Theorems~\ref{t:+inversive} and \ref{t:invertible-add=>pls-elementary}, the plus loop $(R,+)$ of the ternar $R$ is an elementary Moufang loop. Then every element $x\in R$ generates a cyclic subgroup $C_x$  of prime order in the elementary Moufang loop $(R,+)$. By  Corollary~\ref{c:Bol=>Lagrange}, $|C_x|$ divides $|X|_2=2^n$ and hence $|C_x|\le 2$, which implies $x+x=0$, witnessing that the plus loop $(X,+)$ is Boolean. By Theorem~\ref{t:Boolean<=>Boolean-plus}, the affine space $X$ is Boolean.
\smallskip

$(7)\Ra(1)$ Assume that $X$ is inversive-puls, invertible-plus, and $|X|_2=2^n$ for some $n\ge 2$. Let $R$ be any ternar of $X$.  By Theorems~\ref{t:inversive-puls} and \ref{t:invertible-add=>pls-elementary}, the puls loop $(R,\!\puls\!)$ of the ternar $R$ is an elementary Moufang loop. Then every element $x\in R$ generates a cyclic subgroup $C_x$  of prime order in the elementary Moufang loop $(R,\!\puls\!)$. By  Corollary~\ref{c:Bol=>Lagrange}, $|C_x|$ divides $|X|_2=2^n$ and hence $|C_x|\le 2$, which implies $x+x=0$, witnessing that the puls loop $(X,\!\puls\!)$ is Boolean. By Theorem~\ref{t:Boolean<=>Boolean-plus}, the affine space $X$ is Boolean.
\end{proof}

\begin{corollary}\label{c:Fano-projective-order} Every $3$-long Fano projective liner $X$ of finite order and  rank $\|X\|\ge 3$ has order $|X|_2-1=2^n$ for some $n\in\IN$.
\end{corollary}

\begin{proof} By Corollary~\ref{c:Avogadro-projective}, the $3$-long projective liner $X$ is $2$-balanced.  Choose any plane $\Pi$ in $X$ and any line $H$ in $\Pi$. Since $X$ is Fano, so is the plane $\Pi$. By Proposition~\ref{p:projective-minus-hyperplane}  and \ref{t:Fano<=>}, the subliner $A\defeq \Pi\setminus H$ is affine and Boolean. By Theorem~\ref{t:Boolean<=>forder}, $|A|_2=2^n$ for some $n\in\IN$ and hence $|X|_2=|\Pi|_2=1+|A|_2=1+2^n$.
\end{proof}

\begin{proposition}\label{p:Fano-flat-horizon} If $Y$ be a projective completion of a Fano liner $X$, then the horizon $H\defeq Y\setminus X$ is flat in $X$.
\end{proposition}

\begin{proof} Assuming that the horizon $H$ is not flat in $X$, we can find distinct points $x,y\in H$ and $z\in\Aline xy\setminus H$. Since $\overline H\ne Y$, there exists a point $a\in Y\setminus\overline H\subseteq X$. Since the liner $Y$ if $3$-long, there exists a point $b\in \Aline ax\setminus\{a,x\}$. It follows from $x\in \overline H$ and $a\notin H$ that $\Aline ax\cap\overline H=\{x\}$ and hence $b\in Y\setminus\overline H\subseteq X$. It follows from $z\in \Aline xy\subseteq \overline H$ and $a\notin\overline H$ that $\Aline az\cap\overline H=\{z\}$ and hence $\Aline az=(\Aline az\setminus\overline H)\cup\{z\}\subseteq X$. Since the liner $Y$ is projective, there exist points $d\in \Aline by\cap \Aline az\subseteq Y$and $c\in \Aline xd\cap\Aline ay\subseteq Y$.

\begin{picture}(100,95)(-180,-15)
\put(0,0){\line(1,1){60}}
\put(0,0){\line(-1,1){60}}
\put(0,0){\line(0,1){60}}
\put(60,60){\line(-3,-1){90}}
\put(-60,60){\line(3,-1){90}}
\put(-60,60){\line(1,0){120}}

\put(0,0){\circle*{3}}
\put(-2,-8){$a$}
\put(-60,60){\circle*{3}}
\put(-69,58){$x$}
\put(-30,30){\circle*{3}}
\put(-38,25){$b$}
\put(30,30){\circle*{3}}
\put(33,25){$c$}
\put(60,60){\circle*{3}}
\put(63,58){$y$}
\put(0,60){\circle*{3}}
\put(-2,63){$z$}
\put(0,40){\circle*{3}}
\put(1,43){$d$}
\end{picture}

 It is easy to see that $c,d\in Y\setminus\overline H\subseteq X$ and $abcd$ is a  quadrangle in $X$. Since the liner $X$ is Fano, the set 
\begin{multline*}
(\Aline ab\cap X\cap\Aline cd)\cup(\Aline ac\cap X\cap\Aline bd)\cup(\Aline ad\cap X\cap\Aline bc)\\=(\{x\}\cap X)\cup(\{y\}\cap X)\cup(\Aline ad\cap X\cap\Aline bc)=\Aline az\cap X\cap \Aline bc=\Aline az\cap\Aline bc
\end{multline*} is empty, which contradicts the projectivity of the liner $Y$. This contradiction shows that the horizon $H$ of $X$ in $Y$ is flat.
\end{proof}

\begin{theorem}\label{t:Fano-dual} A projective plane is Fano is and only if its dual projective plane is Fano.
\end{theorem}

\begin{proof} It suffices to prove that for every Fano projective plane $(\Pi,\mathcal L)$, the dual projective plane $(\mathcal L,\Pi)$ is Fano. Take any quadrangle $L_0L_1L_2L_3$ in the dual projective plane $(\mathcal L,\Pi)$.
We have to prove that the lines $\overline{(L_0\cap L_1)\cup(L_2\cap L_3)}$, 
$\overline{(L_0\cap L_2)\cup(L_1\cap L_3)}$, $\overline{(L_0\cap L_3)\cup(L_1\cap L_2)}$ are concurrent (and hence belong to some pencil of lines).

Endow the set $4=\{0,1,2,3\}$ with the operation $\oplus$ of addition modulo $4$. For every $i\in 4$, consider the unique point $p_i\in L_i\cap L_{i\oplus 1}\subset \Pi$. Since the lines $L_0L_1,L_2,L_3$ are distinct, the quadruple $p_0p_1p_2p_3$ is a quadrangle in the projective plane $\Pi$. Since $\Pi$ is Fano, the points $x\in \Aline{p_0}{p_1}\cap\Aline{p_2}{p_3}$, $y\in \Aline{p_0}{p_2}\cap\Aline{p_1}{p_3}$ and $z\in \Aline {p_0}{p_3}\cap\Aline{p_1}{p_2}$ are colinear. Observe that 
$$(L_0\cap L_2)\cup(L_1\cap L_3)=(\Aline {p_0}{p_3}\cap\Aline {p_1}{p_2})\cup(\Aline{p_0}{p_1}\cap\Aline{p_2}{p_3})=\{z,x\}.$$  
Then the lines
$$\overline{(L_0\cap L_1)\cup(L_2\cap L_3)}=\Aline {p_0}{p_2},\quad\overline{(L_0\cap L_2)\cup(L_1\cap L_3)}=\Aline xz,\quad\overline{(L_0\cap L_3)\cup(L_1\cap L_2)}=\Aline {p_3}{p_1}$$
contain the point $y$ and hence are concurrent, witnessing that the dual projective plane $(\mathcal L,\Pi)$ is Fano.
\end{proof}

\section{Regularity of proaffine Fano liners}\label{s:Fano2}

The main result of this section is Theorem~\ref{t:proaffine-Fano=>3-regular} on the $3$-regularity of proaffine Fano liners. The proof of this theorem is long and is preceded by four lemmas.

\begin{lemma}\label{l:Fano1} Let $a,b,c,d,e,f$ be distinct points of a proaffine Fano liner $X$ such that $e\in\Aline cd$, $f\in \Aline ac\cap \Aline be$ and $\Aline ab\cap\Aline cd=\varnothing=\Aline bc\cap\Aline ad$. Then $\Aline be\cap\Aline ad\ne\varnothing$.   

\begin{picture}(80,55)(-180,-10)
\linethickness{=0.6pt}
\put(0,0){\line(1,0){60}}
\put(0,0){\line(1,1){30}}
\put(30,0){\line(0,1){30}}
\put(60,0){\line(0,1){30}}
\put(30,30){\line(1,0){30}}
\put(30,30){\line(1,-1){30}}
\put(0,0){\line(2,1){60}}

\put(0,0){\circle*{3}}
\put(-2,-9){$e$}
\put(30,0){\circle*{3}}
\put(27,-10){$d$}
\put(60,0){\circle*{3}}
\put(58,-9){$c$}
\put(30,30){\circle*{3}}
\put(27,33){$a$}
\put(60,30){\circle*{3}}
\put(58,33){$b$}
\put(40,20){\circle*{3}}
\put(36,10){$f$}
\put(30,15){\color{red}\circle*{3}}
\end{picture}
\end{lemma}

\begin{proof} To derive a contradiction, assume that $\Aline be\cap\Aline ad=\varnothing$. Consider the quadrangle $abce$. Since $X$ is Fano, the nonempty set $(\Aline ab\cap\Aline ce)\cup(\Aline ac\cap\Aline be)\cup(\Aline ae\cap\Aline bc)=\varnothing\cup \{f\}\cup(\Aline ae\cap\Aline bc)$ has rank 2 and hence there exists a point $g\in \Aline ae\cap\Aline bc$. Consider the quadrangle $afbg$.  Since $X$ is Fano, the nonempty set $(\Aline af\cap\Aline bg)\cup(\Aline ab\cap\Aline fg)\cup(\Aline ag\cap\Aline fb)=\{c\}\cup(\Aline ab\cap\Aline fg)\cup \{e\}$ has rank 2 and hence $\Aline ab\cap\Aline fg\subseteq \Aline ce\cap\Aline ab=\varnothing$. 
Consider the quadrangle $aedf$.  Since $X$ is Fano, the nonempty set $(\Aline ae\cap\Aline df)\cup(\Aline ad\cap\Aline ef)\cup(\Aline af\cap\Aline ed)=(\Aline ae\cap\Aline df)\cup\varnothing\cup\{c\}$ has rank 2 and hence there exists a point $h\in \Aline ae\cap\Aline df$. 

\begin{picture}(120,150)(-150,-15)
\linethickness{=0.8pt}
\put(0,0){\color{teal}\line(1,0){120}}
\put(0,0){\line(1,1){120}}
\put(60,0){\color{cyan}\line(0,1){60}}
\put(120,0){\color{cyan}\line(0,1){120}}
\put(60,60){\color{teal}\line(1,0){60}}
\put(60,60){\line(1,-1){60}}
\put(0,0){\line(2,1){120}}
\put(60,0){\color{blue}\line(1,2){60}}
\put(60,0){\line(1,1){60}}

\put(0,0){\circle*{3}}
\put(-2,-9){$e$}
\put(60,0){\circle*{3}}
\put(57,-10){$d$}
\put(120,0){\circle*{3}}
\put(118,-9){$c$}
\put(60,60){\circle*{3}}
\put(57,63){$a$}
\put(120,60){\circle*{3}}
\put(122,58){$b$}
\put(80,40){\circle*{3}}
\put(76,29.5){$f$}
\put(60,30){\circle*{3.5}}
\put(60,30){\color{white}\circle*{3}}
\put(120,120){\circle*{2.5}}
\put(123,113){\color{cyan}$g$}
\put(111,120){\color{blue}$h$}
\put(90,60){\circle*{3.5}}
\put(90,60){\color{white}\circle*{3}}
\put(90,30){\circle*{3.5}}
\put(90,30){\color{white}\circle*{3}}

\end{picture}

Consider the quadrangle $abcd$. Since $X$ is Fano, the set $(\Aline ab\cap\Aline cd)\cup(\Aline ac\cap\Aline bd)\cup(\Aline ad\cap bc)=\varnothing\cup(\Aline ac\cap\Aline bd)\cup\varnothing$ is empty. Consider the quadrangle $abdf$. Since $X$ is Fano, the set $$(\Aline ab\cap\Aline df)\cup(\Aline ad\cap\Aline bf)\cup(\Aline af\cap bd)=(\Aline ab\cap\Aline df)\cup\varnothing\cup\varnothing$$ is empty. Then the lines $\Aline df=\Aline dh=\Aline fh$ and $\Aline fg$ do not intersect the line $\Aline ab$. Since $b\in\Aline ef$ and $h,g\in \Aline ea$, the proaffinity of $X$ ensures that $g=h$.
Finally, consider the quadrangle $ecfg$. Since $X$ is Fano, the set $(\Aline ec\cap\Aline fg)\cup(\Aline ef\cap\Aline cg)\cup(\Aline eg\cap\Aline cf)=\{d\}\cup\{b\}\cup\{a\}$ has rank $2$. Then $a\in \Aline bd$, which contradicts $\Aline ac\cap\Aline bd=\varnothing$. This contradiction shows that $\Aline ad\cap\Aline be\ne\varnothing$ and completes the proof of the lemma.
\end{proof}

\begin{lemma}\label{l:Fano2} Let $a,b,d,c,e,f$ be distinct points of a proaffine Fano liner $X$ such that $e\in\Aline cd$, $f\in \Aline ad\cap \Aline be$ and $\Aline ab\cap\Aline cd=\varnothing=\Aline bc\cap\Aline ad$. Then $\Aline ae\cap\Aline bc\ne\varnothing$.   

\begin{picture}(60,70)(-150,-15)
\linethickness{0.7pt}
\put(0,0){\color{teal}\line(1,0){90}}
\put(0,30){\color{teal}\line(1,0){30}}
\put(0,45){\color{red}\line(2,-1){90}}
\put(0,0){\color{cyan}\line(0,1){45}}
\put(30,0){\color{cyan}\line(0,1){30}}
\put(0,45){\line(2,-3){30}}

\put(0,0){\circle*{3}}
\put(-5,-8){$c$}
\put(30,0){\circle*{3}}
\put(26,-9){$b$}
\put(0,30){\circle*{3}}
\put(-8,28){$d$}
\put(30,30){\circle*{3}}
\put(32,30){$a$}
\put(10,30){\circle*{3}}
\put(5,20){$f$}
\put(0,45){\circle*{3}}
\put(-8,45){$e$}
\put(90,0){\color{red}\circle*{3}}
\end{picture}
\end{lemma}

\begin{proof} To derive a contradiction, assume that $\Aline ae\cap\Aline bc=\varnothing$. Consider the quadrangle $abce$. Since $X$ is Fano, the set $(\Aline ab\cap\Aline ce)\cup(\Aline ac\cap\Aline be)\cup(\Aline ae\cap\Aline bc)=\varnothing\cup (\Aline ac\cap\Aline be)\cup \varnothing=\Aline ac\cap\Aline be$ is empty. Consider the quadrangle $abde$ and observe that 
$(\Aline ab\cap de)\cup(\Aline  ad\cap\Aline be)\cup(\Aline ae\cap\Aline bd)=\varnothing\cup\{f\}\cup(\Aline ae\cap\Aline bd)$. Since $X$ is Fano, the nonempty set 
$\{f\}\cup(\Aline ae\cap\Aline bd)$ has rank $2$ and hence $\Aline ae\cap\Aline bd$ contains some point $g$. Consider the quadrangle $abcf$. Since $X$ is Fano, the set $(\Aline ab\cap\Aline cf)\cup(\Aline ac\cap\Aline bf)\cup(\Aline bc\cap\Aline af)=(\Aline ab\cap\Aline cf)\cup\varnothing\cup\varnothing=\Aline ab\cap\Aline cf$ is empty.

\begin{picture}(120,150)(-200,-15)
\linethickness{0.6pt}
\put(-40,0){\color{teal}\line(1,0){220}}
\put(60,60){\color{teal}\line(-1,0){240}}
\put(60,100){\line(-6,-1){240}}

\put(0,90){\color{red}\line(2,-1){180}}
\put(0,90){\color{red}\line(-2,1){60}}
\put(0,0){\color{cyan}\line(0,1){90}}
\put(60,0){\color{cyan}\line(0,1){100}}
\put(0,90){\line(2,-3){60}}
\put(60,0){\line(-1,1){120}}
\put(0,0){\line(1,1){60}}
\put(0,0){\line(1,3){26}}
\put(0,60){\line(3,2){60}}
\put(60,100){\line(-1,-1){100}}

\put(0,0){\circle*{3}}
\put(-5,-8){$c$}
\put(60,0){\circle*{3}}
\put(56,-9){$b$}
\put(0,60){\circle*{3}}
\put(-8,51){$d$}
\put(60,60){\circle*{3}}
\put(62,61){$a$}
\put(20,60){\circle*{3}}
\put(16.5,47.5){$f$}
\put(0,90){\circle*{3}}
\put(-2,93){$e$}
\put(180,0){\circle*{3.5}}
\put(180,0){\color{white}\circle*{3}}
\put(-60,120){\circle*{3}}
\put(-68,120){$g$}
\put(36,36){\circle*{3.5}}
\put(36,36){\color{white}\circle*{3}}
\put(25.7,77.14){\circle*{3}}
\put(23,81){$h$}
\put(60,100){\circle*{3}}
\put(58,103){$i$}
\put(33.3,73.3){\circle*{3}}
\put(31,64){$j$}
\put(-180,60){\circle*{3}}
\put(-189,58){$k$}
\put(13.8,69.2){\circle*{3.5}}
\put(13.8,69.2){\color{white}\circle*{3}}
\put(0,40){\circle*{3}}
\put(-6,39){$l$}
\put(-40,0){\circle*{3}}
\put(-52,-3){$m$}

\end{picture}

Consider the quandrangle $acef$. Since $X$ is Fano, the set $(\Aline ac\cap\Aline ef)\cap(\Aline ae\cap\Aline cf)\cap(\Aline af\cap\Aline ce)=\varnothing\cup(\Aline ae\cap\Aline cf)\cup\{d\}$ has rank $2$ and hence there exists a point $h\in \Aline ae\cap\Aline cf$. By the proaffinity, $\Aline ab\cap\Aline de=\varnothing\ne\Aline ae\cap\Aline bd$ implies $\Aline dh\cap\Aline ab\ne\varnothing$. Therefore, there exists a point $i\in \Aline dh\cap\Aline ab$. Consider the quadrangle $afhi$. Since $X$ is Fano, the nonempty set $$(\Aline af\cap\Aline hi)\cup(\Aline ah\cap\Aline fi)\cap(\Aline ai\cap\Aline fh)=\{d\}\cup(\Aline ah\cap\Aline fi)\cap(\Aline ab\cap\Aline cf)=\{d\}\cup(\Aline ah\cap\Aline fi)\cap\varnothing$$ has rank 2 and hence there exists a point $j\in \Aline ah\cap\Aline fi$. Consider the quadrangle $adei$. Since $X$ is Fano, the nonempty set $(\Aline ad\cap\Aline ei)\cup(\Aline ae\cap\Aline di)\cup(\Aline ai\cap\Aline de)= (\Aline ad\cap\Aline ei)\cup\{h\}\cup\varnothing$ has rank $2$ and hence there exists a point $k\in \Aline ad\cap\Aline ei$. Consider the quadrangle $dehf$. Since $X$ is Fano, the set $(\Aline de\cap\Aline hf)\cup(\Aline dh\cap\Aline ef)\cup(\Aline df\cap\Aline eh)=\{c\}\cup (\Aline dh\cap\Aline ef)\cup\{a\}$ has rank $2$ and hence $\Aline dh\cap\Aline ef\subseteq \Aline ac$ and $\Aline di\cap\Aline ef=\Aline dh\cap\Aline ef=\Aline dh\cap\Aline ef\cap \Aline ac\subseteq \Aline ef\cap\Aline ac=\Aline eb\cap\Aline ac=\varnothing$. Consider the quadrangle $deif$. Since $X$ is Fano, the nonempty set $(\Aline de\cap\Aline if)\cup(\Aline di\cap\Aline ef)\cup(\Aline df\cap\Aline ei)=(\Aline de\cap\Aline if)\cup\varnothing\cup \{k\}$ has rank $2$ and hence there exists a point $l\in \Aline de\cap\Aline if$. It is easy to see that $l\ne c$.

 Since $X$ is proaffine, $\Aline fd\cap\Aline bc=\varnothing$, $\Aline dc\cap\Aline fb=\{e\}$ and $l\in \Aline ec\setminus\{d\}$ imply $\Aline  fl\cap\Aline bc\ne\varnothing$ and hence there exists a point $m\in \Aline  fl\cap\Aline bc$, which is not equal to the points $b,c$. Now observe that $a\in \Aline ib$, $f,j\in\Aline im$, and $\Aline af\cap\Aline bm=\varnothing =\Aline aj\cap\Aline bm$, which contradicts the proaffinity of $X$. This is the final contradiction showing that $\Aline ae\cap\Aline bc\ne\varnothing$.
\end{proof}

\begin{lemma}\label{l:Fano3} Let $a,b,c,\alpha,\beta,\lambda$ be distinct points of a proaffine Fano liner $X$ such that $\Aline ab\cap\Aline \alpha\beta=\varnothing$, and $\alpha\in\Aline \lambda a$, $\beta\in\Aline\lambda b$, and $c\in\Aline ab$. Then $\Aline \lambda c\cap\Aline \alpha\beta\ne\varnothing$.


\end{lemma}

\begin{proof} If $\Aline a\alpha \cap\Aline c\beta=\varnothing$, then $\Aline \lambda c\cap\Aline \alpha\beta\ne\varnothing$, by Lemma~\ref{l:Fano2}. So, assume that $\Aline a\alpha \cap\Aline c\beta$ contains some point $o$. Then $\Aline \lambda c\cap\Aline\alpha\beta\ne\varnothing$, by the proaffinity of $X$.

%
\end{proof}

\begin{lemma}\label{l:Fano4} Let $a,b,c,d,e,f$ be distinct points of a Fano liner $X$ such that $\Aline ab\cap\Aline cd=\varnothing=\Aline ac\cap\Aline bd$, $e\in \Aline cd$, $f\in \Aline ac$, $b\in \Aline ef$, and $\Aline ad\cap\Aline fe=\varnothing$. Then $|\Aline ac|\ge 4$.

%
\end{lemma}

\begin{proof} Consider the quadrangle $abcd$. Since $X$ is Fano, the set
$$(\Aline ab\cap \Aline cd)\cup(\Aline ac\cap \Aline bd)\cup(\Aline ad\cap \Aline bc)=
\varnothing\cup\varnothing\cup (\Aline ad\cap\Aline bc)=\Aline ad\cap\Aline bc$$ has rank $0$ and hence the intersection $\Aline ad\cap\Aline bc$ is empty. Next, consider the quadrangle $afbd$. Since $X$ is Fano, the set
$(\Aline af\cap \Aline bd)\cup(\Aline ab\cap \Aline fd)\cup(\Aline ad\cap \Aline fb)=
\varnothing\cup (\Aline ab\cap\Aline fd)\cup \varnothing$ has rank $0$ and hence the intersection $\Aline ab\cap\Aline fd$ is empty. 

\begin{picture}(80,150)(-150,-15)

\put(0,0){\line(1,0){120}}
\put(0,0){\line(0,1){120}}
\put(0,60){\line(1,0){60}}
\put(60,0){\line(0,1){60}}
\put(0,120){\line(1,-1){120}}
\put(0,60){\line(1,-1){60}}
\put(0,120){\line(1,-2){60}}
\put(0,0){\line(1,1){60}}
\qbezier(0,60)(70,45)(120,0)
\put(60,60){\line(-5,-2){60}}

\put(0,36){\circle*{3}}
\put(-9,34){$h$}
\put(0,0){\circle*{3}}
\put(-2,-9){$c$}
\put(60,0){\circle*{3}}
\put(58,-9){$d$}
\put(0,0){\circle*{3}}
\put(-2,-9){$c$}
\put(120,0){\circle*{3}}
\put(118,-9){$e$}
\put(0,60){\circle*{3}}
\put(-8,58){$a$}
\put(60,60){\circle*{3}}
\put(62,62){$b$}
\put(0,120){\circle*{3}}
\put(-7,119){$f$}
\put(30,30){\circle*{3.5}}
\put(30,30){\color{white}\circle*{3}}

\put(30,60){\circle*{3.5}}
\put(30,60){\color{white}\circle*{3}}
\put(35,50){\circle*{3}}
\put(30,42){$g$}
\end{picture}

Next, consider the quadrangle $afed$. Since $X$ is Fano, the nonempty set $$(\Aline af\cap \Aline ed)\cup(\Aline ae\cap \Aline fd)\cup(\Aline ad\cap \Aline fe)=
\{c\}\cup (\Aline ae\cap\Aline fd)\cup \varnothing$$ has rank $2$ and hence 
$\Aline ae\cap\Aline fd$ contains some point $g$. Assuming that $g\in \Aline bc$, we can consider the quadrangle $cgfe$ and conclude that the set $(\Aline cg\cap \Aline fe)\cup(\Aline cf\cap \Aline ge)\cup(\Aline ce\cap \Aline gf)=
\{b\}\cup\{a\}\cup \{d\}$ has rank $2$ and hence $d\in \Aline ab\cap\Aline cd=\varnothing$, which is a contradiction showing that $g\notin \Aline bc$. Since $\Aline ab\cap\Aline fd=\varnothing$, the point $g\in \Aline fd$ does not belong to the line $\Aline ab$. Then $afbg$ is a quadrangle and the Fano property of $X$ ensures that the nonempty set $$(\Aline af\cap \Aline bg)\cup(\Aline ab\cap \Aline fg)\cup(\Aline ag\cap \Aline fb)=
(\Aline af\cap\Aline bg)\cup \varnothing\cup\{e\}$$ has rank $2$ and hence the set $\Aline af\cap\Aline bg$ contains some point $h$. It is easy to see that $h\notin \{a,c,f\}$ and hence $|\Aline ac|\ge|\{a,c,f,h\}|=4$.
\end{proof}

\begin{theorem}\label{t:proaffine-Fano=>3-regular} Every proaffine Fano liner  $X$ is $3$-regular.
\end{theorem}

\begin{proof} To prove the $3$-regularity of a proaffine Fano liner $X$, fix any concurrent lines $A,\Lambda$ in $X$. We need to check that the set $\Aline A\Lambda\defeq\bigcup_{x\in A}\bigcup_{y\in \Lambda}\Aline xy$ is flat. To derive a contradiction, assume that $\Aline A\Lambda$ is not flat and find points $a,b\in \Aline A\Lambda$ such that $\Aline ab\not\subseteq \Aline A\Lambda$. Then $\{a,b\}\notin A\cup\Lambda$ and we lose no generality assuming that $b\notin A\cup\Lambda$.
Let $o$ be the unique common point of the lines $A,\Lambda$. 
Since $\Aline ab\not\subseteq\Aline A\Lambda$, there exists a point $c\in \Aline ab\setminus\Aline A\Lambda$. 

If $\Aline ab\cap (A\cup \Lambda) \ne\varnothing$, then we lose no generality assuming that $a\in A\cup\Lambda$. Since $b\in\Aline A\Lambda\setminus(A\cup\Lambda)$, there exit points $\beta\in A\setminus \{o\}$ and $\lambda\in\Lambda\setminus\{o\}$ such that $b\in\Aline\beta\lambda$. It follows from $c\notin \Aline A\Lambda$ that $\Aline c\lambda\cap A=\varnothing=\Aline c\beta\cap\Lambda$ and hence $\Aline oc\cap\Aline \lambda\beta=\varnothing$, by the Fano property of $X$. Taking into account that $b\in \Aline ac\cap\Aline\lambda \beta$, we conclude that $o\ne a$. Then $A=\Aline oa$ and $\Lambda=\Aline o\beta$.

If $a\in A$, then we can apply Lemma~\ref{l:Fano1} and find a point $d\in \Aline ac\cap\Aline o\lambda\subseteq\Lambda$. Then $c\in\Aline ab=\Aline ad\subseteq\Aline A\Lambda$, which contradicts the choice of $c$. 



If $a\in \Lambda$, then we can apply Lemma~\ref{l:Fano2} and conclude that $\varnothing\ne\Aline \lambda c\cap \Aline o\beta\subseteq \Aline\lambda c\cap A$ and hence $c\in \Aline \Lambda A$, which contradicts the choice of the point $c$. 

%

Those two contradictions show that $\Aline ab\cap(A\cap\Lambda)=\varnothing$.
Since $a,b\in \Aline A\Lambda$, there exist points $\alpha,\beta\in A$ and $u,v\in\Lambda$ such that $a\in \Aline \alpha u$ and $b\in \Aline \beta v$. It follows from $a,b\notin A$ that $u\ne o\ne v$.

%

 If $u=v$, then we can apply  Lemma~\ref{l:Fano3}, and conclude that $\Aline uc\cap A=\Aline uc\cap \Aline{\alpha}{\beta}\ne\varnothing$ and hence $c\in\Aline \Lambda A$, which contradicts our assumptions. Therefore, $u\ne v$. 

 Since $X$ is proaffine, the sets $I\defeq\{x\in \Aline ou:\Aline ax\cap \Aline \alpha o=\varnothing\}$ and $J\defeq\{x\in \Aline ov:\Aline bx\cap \Aline \beta o=\varnothing\}$ have cardinality $\max\{|I|,|J|\}\le 1$. Assuming that the set $\Lambda\setminus(\{o\}\cup I\cup J)$ contains some point $\lambda$, we can find points $\alpha'\in \Aline \lambda a\cap\Aline \alpha o\subseteq A$ and $\beta'\in\Aline\lambda b\cap\Aline \beta o\subseteq A$.

\begin{picture}(150,75)(-150,-15)
\linethickness{=0.6pt}
\put(0,0){\color{teal}\line(1,0){120}}
\put(0,0){\line(1,1){40}}
\put(20,20){\color{teal}\line(1,0){60}}
\put(40,40){\line(1,-1){40}}
\put(40,40){\line(2,-1){80}}
\put(40,0){\line(0,1){40}}
\put(40,30){\line(2,-1){60}}
\put(20,20){\line(2,-1){40}}

\put(40,0){\circle*{3}}
\put(38,-8){$o$}

\put(120,0){\color{red}\circle*{3}}
\put(0,0){\circle*{3}}
\put(-4,-10){$\alpha'$}
\put(80,0){\circle*{3}}
\put(76,-10){$\beta'$}
\put(20,20){\circle*{3}}
\put(12,20){$a$}
\put(60,20){\circle*{3}}
\put(56,10){$b$}
\put(80,20){\circle*{3}}
\put(82,20){$c$}
\put(40,40){\circle*{3}}
\put(37,43){$\lambda$}
\put(40,10){\circle*{3}}
\put(32,5){$u$}
\put(40,30){\circle*{3}}
\put(32,26){$v$}
\put(60,0){\circle*{3}}
\put(58,-9){$\alpha$}
\put(100,0){\circle*{3}}
\put(98,-10){$\beta$}
\end{picture}

 Applying Lemma~\ref{l:Fano3}, we conclude that $\Aline\lambda c\cap A=\Aline\lambda c\cap \Aline{\alpha'}{\beta'}\ne\varnothing$ and hence $c\in\Aline \Lambda A$, which contradicts our assumptions. Therefore, $\{o,u,v\}\subseteq \Lambda=\{o\}\cup I\cup J$ and hence $I=\{v\}$, $J=\{u\}$, and $\Lambda=\{o,u,v\}$. 
 
 We claim that $\Aline av\cap\Aline bu=\varnothing$. In the opposite case, the set $(\Aline av\cap\Aline bu)\cup(\Aline ab\cap\Aline vu)\cup(\Aline au\cap\Aline bv)=(\Aline av\cap\Aline bu)\cup\varnothing\cup(\Aline au\cap\Aline bv)$ is not empty and has rank $2$, by the Fano property of $X$. Then the intersection $\Aline au\cap\Aline bv=\Aline a\alpha\cap\Aline b\beta$ is not empty and contains some point $\gamma$. By the proaffinity of $X$, $\Aline ub\cap\Aline \alpha\beta=\varnothing$ and $a\in\Aline \gamma u\setminus\{u\}$ imply $\Aline ab\cap\Aline \alpha\beta\ne\varnothing$, which contradicts our assumption. This contradiction shows that $\Aline av\cap\Aline bu=\varnothing$ and hence  $\Aline au\cap\Aline bv=\varnothing$, by the Fano property of $X$.
 
 \begin{picture}(80,90)(-150,-15)

\put(0,0){\line(1,0){60}}
\put(0,0){\line(0,1){60}}
\put(0,30){\line(1,0){30}}
\put(30,0){\line(0,1){30}}
\put(0,60){\line(1,-1){60}}
\put(-30,60){\line(1,-1){60}}
\put(-30,60){\line(1,0){30}}
\put(-30,60){\line(2,-1){60}}

\put(-30,60){\circle*{3}}
\put(-38,58){$a$}
\put(0,0){\circle*{3}}
\put(-5,-7){$o$}
\put(30,0){\circle*{3}}
\put(28,-9){$\alpha$}
\put(60,0){\circle*{3}}
\put(58,-10){$\beta$}
\put(0,30){\circle*{3}}
\put(-8,26){$u$}
\put(30,30){\circle*{3}}
\put(31,31){$b$}
\put(0,60){\circle*{3}}
\put(-2,63){$v$}

\put(0,45){\circle*{3.5}}
\put(0,45){\color{white}\circle*{3}}

\end{picture}

Since $\Lambda=\{o,u,v\}$ and $\alpha\notin \Aline ob\cup\Aline ub\cup\Aline vb$, the line $\Aline \alpha b$ does not intersect the line $\Lambda$. Applying Lemma~\ref{l:Fano4}, we conclude that $|\Lambda|\ge 4$, which contradicts the equality $\Lambda=\{o,u,v\}$. This is the final contradiction completing the proof of the $3$-regularity of the proaffine Fano liner $X$.
\end{proof}


\begin{corollary}\label{c:proaffineFano=>para-Playfair} Every proaffine Fano liner $X$ is Proclus and para-Playfair.
\end{corollary}

\begin{proof} By Theorem~\ref{t:proaffine-Fano=>3-regular}, the proaffine Fano liner $X$ is $3$-regular and hence $3$-proregular. By Theorem~\ref{t:Proclus<=>}, the $3$-proregular proaffine liner $X$ is Proclus. 

To prove that $X$ is para-Playfair, fix any plane $\Pi\subseteq X$ and two disjoint lines $L,L'$ in $\Pi$. Given any point $x\in \Pi\setminus L$, we need to find a line $\Lambda$ such that $x\in \Lambda\subseteq\Pi\setminus L$. If $x\in L'$, then the line $\Lambda\defeq L'$ has the required property. So, assume that $x\notin L'$ and choose any distinct points $a,b\in L$. Since $X$ is Proclus, there exist points $\alpha\in \Aline xa\cap L'$ and $\beta\in\Aline xb\cap L'$. Consider the quadrangle $a\alpha\beta b$. Since $X$ is Fano, the nonempty set
$$(\Aline a\alpha\cap\Aline\beta b)\cup(\Aline a\beta\cap\Aline \alpha b)\cup(\Aline ab\cap\Aline \alpha\beta)=\{x\}\cup(\Aline a\beta\cap\Aline \alpha b)\cup\varnothing$$ has rank $2$ and hence the intersection $\Aline a\beta\cap\Aline \alpha b$ contains some point $y$. 

\begin{picture}(90,90)(-150,-15)
\linethickness{=0.6pt}
\put(0,0){\line(1,2){30}}
\put(60,0){\line(-1,2){30}}
\put(15,30){\color{red}\line(1,0){65}}
\put(85,27){$L'$}
\put(0,0){\color{red}\line(1,0){80}}
\put(85,-3){$L$}
\put(0,0){\line(3,2){45}}
\put(60,0){\line(-3,2){45}}
\put(30,60){\color{red}\line(0,-1){40}}

\put(0,0){\circle*{3}}
\put(-2,-9){$a$}
\put(60,0){\circle*{3}}
\put(57,-9){$b$}
\put(15,30){\circle*{3}}
\put(8,33){$\alpha$}
\put(45,30){\circle*{3}}
\put(46,33){$\beta$}
\put(30,60){\circle*{3}}
\put(27,63){$x$}
\put(30,20){\circle*{3}}
\put(27,12){$y$}
\put(30,30){\color{red}\circle*{3.5}}
\put(30,30){\color{white}\circle*{3}}
\end{picture}

Consider the quadrangle $xayb$. Since $X$ is Fano, the nonempty set
$(\Aline xa\cap\Aline yb)\cup(\Aline xy\cap\Aline ab)\cup(\Aline xb\cap\Aline ay)=\{\alpha\}\cup (\Aline xy\cap\Aline ab)\cup\{\beta\}$ has rank $2$ and hence $\Aline xy\cap\Aline ab\subseteq\Aline \alpha\beta$. Then
$$\Aline xy\cap L=\Aline xy\cap\Aline ab=\Aline xy\cap\Aline ab\cap\Aline \alpha\beta\subseteq \Aline ab\cap\Aline \alpha\beta\subseteq L'\cap L=\varnothing$$ and hence the line $\Aline xy$ has the required property: $x\in\Aline xy\subseteq\Pi\setminus L$.
\end{proof}

\begin{example} For any Steiner projective liner $X$ of rank $\|X\|=4$ and any line $H$ in $X$ the Fano subliner $X\setminus H$ of $X$ is not regular.
\end{example}

Now we are going to show that every regular proaffine Fano liner is completely regular. In the proof we shall use the following proposition on the complete regularity of Boolean para-Playfair regular liners.

\begin{proposition}\label{p:Boolean-para-Playfair=>compreg} Every Boolean para-Playfair regular liner $X$ is completely regular.
\end{proposition}

\begin{proof} By Theorem~\ref{t:spread=projective2}, it suffices to show that for every concurrent spreading lines $A,B$ in $X$, every line $L$ in the plane $\Pi\defeq \overline{A\cup B}$ is spreading. This is clear if $L\in A_\parallel\cup B_\parallel$. So, assume that $L\notin A_\parallel \cup B_\parallel$. Given any point $x\in  \Pi\setminus L$, we should find a line $L_x\subseteq \Pi\setminus L$ that contains the point $x$. Since the lines $A,B$ are spreading, there exist lines $A_x\in A_\parallel$ and $B_x\in B_\parallel$ such that $x\in A_x\cap B_x$. Since $L\notin A_\parallel\cup B_\parallel$, there exist unique points $a\in A_x\cap L$ and $b\in B_x\cap L$, see the Proclus Postulate~\ref{p:Proclus-Postulate}. Find unique lines $A_b\in A_\parallel$ and $B_a\in B_\parallel$ such that $b\in A_b$ and $a\in B_a$. The liner $X$ is para-Playfair and hence Proclus. By the Proclus Postulate~\ref{p:Proclus-Postulate}, the lines $A_b$ and $B_a$ in the plane $\Pi$ have a unique common point $c$. Then $xacb$ is a parallelogram in $X$. Since $X$ is Boolean, the diagonals $\Aline xc$ and $\Aline ab=L$ of the parallelogram $xacb$ are disjoint. Then the line $L_x\defeq \Aline xc$ has the required property: $x\in L_x\subseteq \Pi\setminus L$. 

\begin{picture}(100,100)(-150,-45)
\linethickness{=0.6pt}

\put(0,0){\color{red}\line(1,0){80}}
\put(0,0){\color{blue}\line(1,1){30}}
\put(0,0){\color{violet}\line(1,-1){30}}
\put(60,0){\color{violet}\line(-1,1){30}}
\put(60,0){\color{blue}\line(-1,-1){30}}
\put(30,-30){\color{red}\line(0,1){70}}

\put(85,-3){\color{red}$L$}
\put(26,45){\color{red}$L_x$}
\put(3,18){\color{blue}$A_x$}
\put(44,-22){\color{blue}$A_b$}
\put(44,19){\color{violet}$B_x$}
\put(3,-22){\color{violet}$B_a$}

\put(0,0){\circle*{3}}
\put(-8,-3){$a$}
\put(60,0){\circle*{3}}
\put(61,2){$b$}
\put(30,30){\circle*{3}}
\put(22,29){$x$}
\put(30,-30){\circle*{3}}
\put(27,-38){$c$}
\put(30,0){\color{red}\circle*{3.5}}
\put(30,0){\color{white}\circle*{3}}

\end{picture}
\end{proof}

\begin{theorem}\label{t:wreg-proFano=>compreg}    
Every weakly regular proaffine Fano liner $X$ is completely regular.
\end{theorem}
    
\begin{proof} By Theorem~\ref{t:proaffine-Fano=>3-regular} and Corollary~\ref{c:proaffineFano=>para-Playfair}, the proaffine Fano liner $X$ is $3$-regular and para-Playfair, and by Proposition~\ref{p:Fano=>Boolean}, the Fano liner $X$ is Boolean. By Propositions~\ref{p:reg<=>wreg+3-reg}, the $3$-regular weakly regular liner $X$ is regular. By Proposition~\ref{p:Boolean-para-Playfair=>compreg}, the Boolean para-Playfair regular liner $X$ is completely regular.
\end{proof}

In the following corollary, by a \defterm{proaffine Fano plane} we understand any proaffine Fano liner of rank $3$.

\begin{corollary}\label{c:proafFano=>compreg} Every proaffine Fano plane is completely regular.
\end{corollary}

\begin{proof} Let $X$ be a proaffine Fano plane. By Theorem~\ref{t:proaffine-Fano=>3-regular}, $X$ is $3$-regular and hence $\|X\|$-regular. By Corollary~\ref{c:reg<=>||X||-reg}, the $\|X\|$-regular liner $X$ is regular. By Theorem~\ref{t:wreg-proFano=>compreg}, the (weakly) regular proaffine Fano liner $X$ is completely regular. 
\end{proof}

\begin{example} Let $Y$ be a Steiner projective liner of rank $\|X\|=4$. For every line $L\subset Y$, the subliner $X\defeq Y\setminus L$ is proaffine and Fano, by Theorems~\ref{t:proaffine<=>proflat} and \ref{t:Fano<=>}. By  Theorem~\ref{t:regular-horizon3}, the liner $X=Y\setminus L$ is not regular. This example shows that Corollary~\ref{c:proafFano=>compreg} cannot be extended to arbitrary proaffine Fano liners. On the other hand, in Theorem~\ref{t:completelyF<=>} we shall prove that every finite $3$-long proaffine Fano liner is completely regular.
\end{example}

\begin{problem} Is every $\w$-long proaffine Fano liner regular?
\end{problem}

\section{The spread completions of Fano liners}\label{s:Fano3}

In section we explore the structure of the spread completions of completely regular  Fano liners. First, we consider the spread completions of proaffine Fano planes. 

\begin{theorem}\label{t:Fano-completion-plane} The spread completion of any proaffine Fano plane is a projective Fano plane.
\end{theorem} 

\begin{proof} Let $X$ be a proaffine Fano plane. 
By Corollary~\ref{c:proafFano=>compreg}, the liner $X$ is completely regular and hence its spead completion $\overline X$ is a projective liner. If $X$ is projective, then $\overline X=X$ is a projective Fano plane. So, assume that $X$ is not projective. By Corollary~\ref{c:proaffineFano=>para-Playfair}, the proaffine Fano plane $X$ is Proclus. If $X$ is not $3$-long, then we can apply Theorem~\ref{t:Proclus-not-3long} and conclude that $|X|=6$, $\overline X\setminus X$ is a singleton and $\overline X$ is a Steiner (and hence Fano) plane.

So, assume that the Fano liner $X$ is $3$-long. Then the spread completion $\overline X$ of $X$ is a $3$-long projective plane, by Corollary~\ref{c:spread-3long} and \ref{c:procompletion-rank}. It remains to prove that the projective plane $\overline X$ is Fano. Since $X$ is not projective, the horizon $H\defeq\overline X\setminus X$ is either a line of a singleton in $\overline X$. 

Let $abcd$ be any quadrangle in $\overline X$. Since the liner $\overline X$ is projective, there exist unique points $x\in\Aline ab\cap\Aline cd$, $y\in \Aline ac\cap\Aline bd$ and $z\in \Aline ad\cap bd$. We have to show that the set $\{x,y,z\}$ has rank $2$ in $\overline X$. Depending on the cardinality of the intersection $\{a,b,c,d,x,y,z\}\cap H$ we shall consider three cases.
\smallskip

0. First assume that $\{a,b,c,d,x,y,z\}\cap H=\varnothing$. In this case $\{a,b,c,d,x,y,z\}\subseteq X$. Since $X$ is Fano, the set $\{x,y,z\}$ has rank $2$ in $X$ and also in $\overline X$.
\smallskip

1. Next, assume that $\{a,b,c,d,x,y,z\}\cap H$ is a singleton. If this singleton is contained in the set $\{x,y,z\}$, then we lose no generality assuming that this singleton is $\{z\}$. In this case, $\Aline ad\cap X=\Aline az\cap X$ and $\Aline bc\cap X=\Aline bz\cap X$ are two distinct parallel lines in the liner $X$. Consider the quadrangle $bxcy$ in the liner $X$. Since $X$ is Fano, the nonempty set
$$(\Aline bx\cap X\cap\Aline cy)\cup(\Aline bc\cap X\cap\Aline xy)\cup(\Aline by\cap X\cap \Aline xc)= \{a\}\cup (\Aline bc\cap X\cap \Aline xy)\cup \{d\}$$ has rank 2 in $X$ and hence $\Aline bc\cap X\cap \Aline xy\subseteq \Aline ad$ and
$\Aline bc\cap X\cap\Aline xy=\Aline ad\cap\Aline bc\cap X\cap\Aline xy=\varnothing$.

\begin{picture}(200,90)(-140,-15)
\put(0,0){\color{red}\vector(1,0){150}}
\put(0,0){\line(0,1){60}}
\put(0,30){\color{red}\vector(1,0){150}}
\put(0,30){\line(2,-1){60}}
\put(20,20){\color{red}\vector(1,0){130}}
\put(0,0){\line(1,1){30}}
\put(0,60){\line(1,-1){60}}
\put(0,60){\color{red}\line(1,-2){20}}
\put(160,15){\color{red}$z$}

\put(0,0){\circle*{3}}
\put(-3,-8){$a$}
\put(60,0){\circle*{3}}
\put(57,-9){$d$}
\put(20,20){\circle*{3}}
\put(17,10){$y$}
\put(30,30){\circle*{3}}
\put(31,32){$c$}
\put(0,30){\circle*{3}}
\put(-8,28){$b$}
\put(0,60){\circle*{3}}
\put(-2,63){$x$}
\put(15,30){\color{red}\circle*{3.5}}
\put(15,30){\color{white}\circle*{3}}
\end{picture}

It follows that the lines $\Aline xy\cap X$ and $\Aline yz\cap X$ both are disjoint with the line $\Aline bc\cap X$ in the liner $X$. Since $X$ is Proclus, $\Aline xy\cap X=\Aline xz\cap X$ and hence the set $\{x,y,z\}$ has rank $2$ in the projective liner $\overline X$. 
\smallskip

 Next, assume that the singleton $\{a,b,c,d,x,y,z\}\cap H$  is contained in the set $\{a,b,c,d\}$. In this case we lose no generality assuming that $\{a,b,c,d,x,y,z\}\cap H=\{d\}$. Then $\Aline az\cap X=\Aline ad\cap X,\Aline yb\cap X= \Aline bd\cap X$, and $\Aline xc\cap X=\Aline cd\cap X$ are three distinct parallel lines in the liner $X$. Consider the quadrangle $azxc$ in the liner $X$. Since $X$ is Fano, the nonempty set
$$(\Aline az\cap X\cap\Aline xc)\cup(\Aline ax\cap X\cap\Aline zc)\cup(\Aline ac\cap X\cap \Aline zx)= \varnothing\cup\{b\}\cup (\Aline ac\cap X\cap \Aline zx)$$ has rank 2 in $X$ and hence there exists a point $e\in \Aline ac\cap X\cap \Aline zx$.

\begin{picture}(200,90)(-140,-15)
\put(0,0){\color{red}\vector(1,0){150}}
\put(0,0){\line(0,1){60}}
\put(0,30){\color{red}\vector(1,0){150}}
\put(0,30){\line(2,-1){60}}
\put(0,20){\color{red}\vector(1,0){150}}
\put(0,0){\line(1,1){30}}
\put(0,60){\line(1,-1){60}}
\put(0,60){\color{red}\line(1,-2){20}}
\put(160,15){\color{red}$d$}

\put(0,0){\circle*{3}}
\put(-3,-8){$c$}
\put(60,0){\circle*{3}}
\put(57,-8){$x$}
\put(20,20){\circle*{3}}
\put(17,10){$b$}
\put(30,30){\circle*{3}}
\put(31,32){$z$}
\put(0,20){\circle*{3}}
\put(-8,18){$y$}
\put(0,30){\circle*{3}}
\put(-8,28){$a$}
\put(0,60){\circle*{3}}
\put(-2,63){$e$}
\put(15,30){\color{red}\circle*{3.5}}
\put(15,30){\color{white}\circle*{3}}
\end{picture}

Consider the quadrangle $aezb$ in $X$. Since $X$ is Fano, the nonempty set 
$$(\Aline ae\cap X\cap\Aline zb)\cup(\Aline az\cap X\cap\Aline eb)\cup(\Aline ab\cap X\cap \Aline ez)= \{c\}\cup (\Aline az\cap X\cap \Aline eb)\cup\{x\}$$ has rank $2$ in $X$ and hence  $\Aline az\cap X\cap \Aline eb\subseteq \Aline cx$. Then $$\Aline az\cap X\cap \Aline eb=\Aline cx\cap \Aline az\cap X\cap \Aline eb=\varnothing,$$which means that the lines $\Aline eb\cap X$ and $\Aline az\cap X$ are disjoint in the plane $X$. Since the lines  $\Aline yb\cap X$ and $\Aline az\cap X$ also are disjoint,  
the Proclus property of the plane $X$ ensures that $\Aline be\cap X=\Aline by\cap X$ and hence $y\in \Aline ac\cap\Aline by=\Aline ac\cap \Aline be=\{e\}$. Then $y=e\in \Aline xz$ and hence the set $\{x,y,z\}$ has rank $2$ in $X$ and $\overline X$.
\smallskip

2. Finally, assume that the intersection $D\defeq \{a,b,c,d,x,y,z\}\cap H$ contains at least two points.  If $D$ has two points in the set $\{x,y,z\}$, then we lose no generality assuming that $x,y\in D$. Assuming that $a\in H$, we conclude that $b\in \Aline ax\in H$ and $c\in \Aline ay\subseteq H$, and hence the points $a,b,c$ are collinear, which contradicts the choice of the quadrangle $abcd$. Terefore, $a\notin H$. By analogy we can show that $b,c\notin H$. Then $abcd$ is a quadrangle in the liner $X$. Since $X$ is Fano, the set 
$$(\Aline ab\cap X\cap\Aline cd)\cup(\Aline ac\cap X\cap\Aline bd)\cup(\Aline ad\cap X\cap \Aline bc)= (\{x\}\cap X)\cup (\{y\}\cap X)\cup (\Aline ad\cap X\cap\Aline bc)=\Aline ad\cap X\cap \Aline bc$$
has rank $0$ in $X$ and hence $\{z\}=\Aline ad\cap \Aline bc\subseteq H$. Then the set $\{x,y,z\}\subseteq H$ has rank $\|\{x,y,z\}\|=\|H\|=2$ and we are done.

Next, assume that $D$ contains at most one point in the set $\{x,y,z\}$. In this case $D\cap\{x,y,z\}$ is a singleton. Indeed, assuming that $D\cap \{x,y,z\}$ is empty, we conclude that $D\cap \{a,b,c,d\}$ contains two distinct points. We lose no generality assuming that $a,b\in D$. Then $x\in\Aline ab\subseteq H$, which contradicts $D\cap\{x,y,z\}=\varnothing$. This contradiction shows that $D\cap \{x,y,z\}$ is a singleton. We lose no generality assuming that this singleton equals $\{y\}$. Since $D$ contains at least two points, the set $D\cap\{a,b,c,d\}$ is not empty and we lose no generality assuming that $d\in D$. In this case $b\in \Aline yd\subseteq \overline D\subseteq H$. Since $abcd$ is a quadrangle, $\{a,b,c,d,x,y,z\}\cap H=\{b,d,y\}$ and hence $azcx$ is a quadrangle in $X$ such that $\Aline az\cap\Aline xc=\{d\}$ and $\Aline ax\cap\Aline zc=\{b\}$.

\begin{picture}(100,115)(-140,-15)
\linethickness{=0.6pt}
\put(0,0){\color{teal}\vector(1,0){80}}
\put(90,13){\color{teal}$b$}
\put(0,0){\color{cyan}\vector(0,1){75}}
\put(14,85){\color{cyan}$d$}
\put(0,0){\color{red}\vector(1,1){65}}
\put(70,68){\color{red}$y$}
\put(0,30){\color{teal}\vector(1,0){80}}
\put(30,0){\color{cyan}\vector(0,1){75}}
\put(0,30){\color{red}\line(1,-1){30}}

\put(0,0){\circle*{3}}
\put(-3,-8){$a$}
\put(30,0){\circle*{3}}
\put(27,-8){$x$}
\put(0,30){\circle*{3}}
\put(-8,28){$z$}
\put(30,30){\circle*{3}}
\put(31,23){$c$}
\put(15,15){\color{red}\circle*{3.5}}
\put(15,15){\color{white}\circle*{3}}
\end{picture}

Since $X$ is Fano, the set $$(\Aline az\cap X\cap\Aline cx)\cup(\Aline ac\cap X\cap\Aline zx)\cup(\Aline ax\cap X\cap \Aline zc)= (\{d\}\cap X)\cup (\Aline ac\cap X\cap \Aline zx)\cup(\{b\}\cap X)=\Aline ac\cap X\cap \Aline zx$$ has rank $0$ in $X$ and hence  $\varnothing\ne\Aline ac\cap \Aline zx\subseteq \Aline ac\setminus X=\{y\}$, which implies $y\in \Aline zx$ and hence $\|\{x,y,z\}\|=2$.
\end{proof}

\begin{corollary}\label{c:SteinerFano=>projective} Every Steiner Fano liner is projective.
\end{corollary}

\begin{proof} Let $X$ be a Steiner Fano liner. To prove that $X$ is projective, it suffices to check that every distinct coplanar lines $L,\Lambda$ in $X$ are concurent. Consider the plane $\Pi\defeq\overline{L\cup \Lambda}$ in $X$. 
Since every Steiner liner is proaffine, the Steiner Fano liner $X$ is proaffine. 
By Corollary~\ref{c:proafFano=>compreg} and Theorem~\ref{t:Fano-completion-plane}, the proaffine Fano plane $\Pi$ is completely regular and its spread completion $\overline\Pi$ is a Fano projective plane. By Corollary~\ref{c:Fano-projective-order}, $|\overline\Pi|_2=1+2^k$ for some $k\in\IN$. Since $X$ is Steiner, $3=|\Pi|_2\in \{|\overline\Pi|_2,|\overline\Pi|_2-1\}=\{1+2^k,2^k\}$ implies $k=1$. Assuming that $\Pi\ne\overline\Pi$, we can fix points $x\in\Pi$ and $y\in\overline\Pi\setminus\Pi$ and conclude that the line $\Aline xy\cap\Pi=\Aline xy\setminus\{y\}$ has length $|\Aline xy\setminus\{y\}|=|\overline \Pi|_2-1=3-1=2$, which contradicts the Steiner property of $X$. This contradiction shows that $\Pi=\overline\Pi$ and hence the lines $L, \Lambda$ is the projective plane $\Pi=\overline\Pi$ are concurrent.  By Theorem~\ref{t:projective<=>}, the liner $X$ is projective.
\end{proof}

\begin{theorem}\label{t:compreg-Fano=>comp-Fano} The spread completion $\overline X$ of a completely regular Fano liner $X$ is a projective Fano liner.
\end{theorem}

\begin{proof} 
Since the liner $X$ is completely regular, its spread completion $\overline X$ is a projective liner. It remains to prove that $\overline X$ is Fano.  If $\|X\|\le 2$, then $\|\overline X\|\le 2$ and hence $\overline X$ is Fano (because $\overline X$ contains no quadrangles). If $\|X\|=3$, then $\overline X$ is a projective Fano liner, by Theorem~\ref{t:Fano-completion-plane}. So, assume that $\|X\|\ge 4$. If $X$ is projective, then $\overline X=X$ and hence $\overline X$ is Fano. So, assume that $X$ is not projective.
By Theorems~\ref{t:spread=projective1} and \ref{t:Proclus<=>}, the completely regular liner $X$ is para-Playfair and hence Proclus and proaffine.
 By Proposition~\ref{p:spread-3long}, the projective liner $\overline X$ is $3$-long. By Corollary~\ref{c:Avogadro-projective}, the $3$-long projective liner $\overline X$ is $2$-balanced. If $|\overline X|_2=3$, then $\overline X$ is Fano, by Example~\ref{ex:Fano3}. So, assume that $|\overline X|_2\ge 4$. By 
 
Choose any hyperplane $H$ in $\overline X$ containing the horizon $\overline X\setminus X$ of $X$ in $\overline X$.  By Propositions~\ref{p:projective-minus-flat}(2) and \ref{p:projective-minus-hyperplane}, the subliner $A=\overline X\setminus H$ of $\overline X$ is $3$-long, affine, regular, and has rank $\|A\|=\|\overline X\|\ge 4$. Choose any plane $\Pi$ in $\overline X$ such that $\Pi\cap A\ne\varnothing$. Then $\Pi\cap X$ is a plane in the Fano liner $X$ and hence $\Pi\cap X$ is Fano. Since the set $\Pi\setminus X$ is flat in $\Pi$, the liner $\Pi\cap X$ is proaffine, by Theorem~\ref{t:proaffine<=>proflat}.
By Corollary~\ref{c:proafFano=>compreg} and Theorem~\ref{t:Fano-completion-plane}, the proaffine Fano plane $\Pi\cap X$ is completely regular and its spread completion $\overline{\Pi\cap X}$ is a Fano projective plane. By Corollary~\ref{c:pcompletion=scompletion}, the projective completion $\Pi$ of the $3$-long liner $\Pi\cap X$ is isomorphic to the spread completion $\overline{\Pi\cap X}$ of $\Pi\cap X$. Then the projective plane $\Pi$ is Fano. Since the set $\Pi\cap H$ is flat in $\Pi$, the affine liner $A\cap\Pi=\Pi\setminus H$ is Boolean, by Theorem~\ref{t:Fano<=>}. Then  $A\cap\Pi$ contains a Boolean parallelogram, which remains Boolean in the affine space $A$. By Theorem~\ref{p:DesFano<=>}, the Desarguesian projective liner $\overline X$ is Fano.  
\end{proof}

The following theorem characterizes completely regular $3$-long Fano liners.

\begin{theorem}\label{t:completelyF<=>} For a $3$-long  liner $X$, the following conditions are equivalent:
\begin{enumerate}
\item $X$ is Fano, proaffine, and weakly regular;
\item $X$ is Fano and completely regular;
\item $X$ is completely regular and its spread completion is Fano;
\item $X$ is $3$-ranked and its spread completion is a Fano projective liner;
\item $X$ has a Fano projective completion with flat horizon;
\item $X$ is Fano and has a projective completion;
\item $X$ is Fano, proaffine, and regular.
\end{enumerate}
If the liner $X$ is finite, then \textup{(1)--(7)} are equivalent to the condition
\begin{enumerate}
\item[\textup{(8)}] $X$ is proaffine and Fano.
\end{enumerate}
\end{theorem}

\begin{proof} The implications $(1)\Ra(2)\Ra(3)$ follow from Theorems~\ref{t:wreg-proFano=>compreg} and \ref{t:compreg-Fano=>comp-Fano}, respectively. The implication $(3)\Ra(4)\Ra(5)$ follows from Definition~\ref{d:comp-regular} of a completely regular liner and Definition~\ref{d:procompletion} of a projective completion. The implication $(5)\Ra(6)$ follows from Theorem~\ref{t:Fano<=>}.
\smallskip

$(6)\Ra(7)$ Assume that $X$ is Fano and has a projective completion $Y$. By Proposition~\ref{p:Fano-flat-horizon}, the horizon $H\defeq Y\setminus X$ is flat in $Y$. By Theorem~\ref{t:proaffine<=>proflat} and Proposition~\ref{p:projective-minus-flat}, the $3$-long liner $X=Y\setminus H$ is proaffine and regular.
\smallskip

The implications $(7)\Ra(1\wedge 8)$ are trivial.
\smallskip

$(8)\Ra(1)$ Assume that $X$ is a finite proregular Fano liner. By Corollary~\ref{c:proaffineFano=>para-Playfair}, the proaffine Fano liner $X$ is para-Playfair. By Theorem~\ref{t:spread=projective2}, the finite $3$-long para-Playfair liner $X$ is completely regular. 
 \end{proof}

By Theorem~\ref{t:Fano<=>}, a $3$-long projective liner is Fano if and only if it is everywhere Fano if and only if it is everywhere Boolean. The following corollary implies that a $3$-long projective liner is Fano if and only if it is somewhere Fano.

\begin{corollary}\label{c:projDes<=>} For a $3$-long projective liner $Y$, the following conditions are equivalent:
\begin{enumerate}
\item $Y$ is Fano;
\item for every flat $H$ in $Y$, the liner $Y\setminus H$ is Fano;
\item for every hyperplane $H$ in $Y$, the liner $Y\setminus H$ is Boolean;
\item for some hyperplane $H$ in $Y$, the liner $Y\setminus H$ is Fano;
\item for some proper flat $H$ in $Y$, the liner $Y\setminus H$ is Fano;
\item $Y$ is a projective completion of some affine regular Fano liner;
\item $Y$ is a projective completion of some proaffine Fano liner.
\end{enumerate}
\end{corollary}

\begin{proof} Theorem~\ref{t:Fano<=>} implies that the conditions $(1)$--$(3)$ are equivalent. The implications $(2)\Ra(4)\Ra(5)$ and $(6)\Ra(7)$ are trivial. 
\smallskip

$(4)\Ra(6)$ Assume that  for some proper flat $H$ in $Y$, the liner $X\defeq Y\setminus H$ is Fano. By Theorem~\ref{t:proaffine<=>proflat}, the liner $Y$ is proaffine. Therefore, $Y$ is a projective completion of the proaffine Fano liner $X$.
\smallskip

$(5)\Ra(7)$ Assume that  for some hyperplane $H$ in $Y$, the liner $X\defeq Y\setminus H$ is Fano. By Proposition~\ref{p:projective-minus-hyperplane}, the liner $X$ is affine and regular. 
Thereore, $Y$ is a projective completion of the affine regular Fano liner $X$.
\smallskip

$(7)\Ra(1)$ Assume that $Y$ is a projective completion of some proaffine Fano liner $X\subseteq Y$. If $|Y|_2=3$, then $Y$ is Fano, by Example~\ref{ex:Fano3}. So, assume that $|Y|_2\ge 4$. Then $|X|_2\ge |Y|_2-1\ge 3$. By Proposition~\ref{p:Fano-flat-horizon}, the horizon $H\defeq Y\setminus X$ of the Fano liner $X$ in its projective completion $Y$ is flat. By Theorem~\ref{t:proaffine<=>proflat} and Proposition~\ref{p:projective-minus-flat}, the subliner $X=Y\setminus H$ is proaffine and regular. By Theorems~\ref{t:wreg-proFano=>compreg} and \ref{t:compreg-Fano=>comp-Fano}, the $3$-long regular proaffine Fano liner $X$ is completely regular and its spread completion $\overline X$ is a projective Fano liner. By Corollary~\ref{c:pcompletion=scompletion}, the projective completion $Y$ of $X$ is isomorphic to the spread completion $\overline X$ of $X$, which implies that the projective liner $Y$ is Fano. 
\end{proof}

Let us recall that a liner $X$ is \index{completely $\mathcal P$ liner}\index{liner!completely $\mathcal P$}\defterm{completely $\mathcal P$} for some property $\mathcal P$ of projective liners if $X$ is completely regular and its spread completion $\overline X$ has property $\mathcal P$. Theorem~\ref{t:completelyF<=>} implies the following corollary answering Problem~\ref{prob:inner-completely-P} for the Moufang property.

\begin{corollary}\label{c:completelyF<=>} A liner is completely Fano if and only if it is Fano, proaffine, and (weakly) regular.
\end{corollary}

\section{Finite proaffine Fano liners}\label{s:Gleason-Fano}

In this section we explore the structure of finite Fano liners. 
\begin{theorem}[Gleason, 1956]\label{t:Gleason-Fano} Every finite Fano projective liner is Pappian.
\end{theorem}

\begin{proof} Let $X$ be a finite Fano projective liner. By Proposition~\ref{p:Papp<=>4-Papp}, the projective liner $X$ is Pappian if and only if every $4$-long plane in $X$ is Pappian. So, fix any $4$-long projective plane $\Pi$ in $X$. The Pappian property of $\Pi$ will follow from the Gleason Theorem~\ref{t:Gleason56} as soon as we check that for every line $L\subseteq \Pi$ and point $p\in L$, the group $\Aut_{p,L}(\Pi)$ of automorphisms of $\Pi$ with centre $p$ and axis $L$ is not trivial. Consider the affine plane $A\defeq \Pi\setminus L$ and the spread of parallel lines $\boldsymbol \delta\defeq\{\Aline xp\setminus\{p\}:x\in A\}$ in $A$. Since the projective plane $\Pi$ is Fano, the affine plane $A=\Pi\setminus L$ is Boolean, by Theorem~\ref{t:Fano<=>}. By Corollary~\ref{c:Boolean-order}, $|A|_2=2^n$ for some $n\in \IN$. Since $\Pi$ is $4$-long, the affine plane $A$ is $3$-long and hence $3\le |A|_2=2^n$ and $n\ge 2$. By Theorem~\ref{t:4-long-affine} the $4$-long affine plane $A$ is regular, and by Theorem~\ref{t:Playfair<=>}, the $4$-long regular affine plane $A$ is Playfair. By Theorems~\ref{t:Boolean=>commutative-plus} and \ref{t:add-com=>add-ass}, the Boolean Playfair plane $A$ is commutative-plus and associative-plus.  By Theorem~\ref{t:add-ass=>partialT}, the Playfair plane $A$ is $\partial$-translation, which implies the existence of a non-trivial translation $T:A\to A$ such that $\{\Aline xy:xy\in T\}=\boldsymbol \delta$. By Theorem~\ref{t:extend-isomorphism-to-completions}, the translation $T$ entends to a (non-trivial) automorphism $\bar T\in\Aut_{p,L}(\Pi)$, witnessing that the group $\Aut_{p,L}(\Pi)$ is not trivial. The Gleason Theorem~\ref{t:Gleason56} ensures that the projective plane $\Pi$ is Pappian. By Proposition~\ref{p:Papp<=>4-Papp}, the projective liner $X$ is Pappian.
\end{proof}

The following theorem extends Gleason's Theorem~\ref{t:Gleason-Fano} to finite proaffine Fano liners.

\begin{theorem}\label{t:proaffFano=>Pappian} Every finite proaffine Fano liner is Pappian.
\end{theorem} 

\begin{proof} Let $X$ be a finite proaffine Fano liner. To prove that $X$ is Pappian, fix any concurrent lines $L,L'\subseteq X$ and any distinct points $a,b,c\in L\setminus L'$ and $a',b',c'\in L'\setminus L$. We have to prove that the set $T\defeq(\Aline a{b'}\cap\Aline {a'}b)\cup (\Aline a{c'}\cap\Aline{a'}c)\cup(\Aline b{c'}\cap\Aline{b'}c)$ has rank $\|T\|\in\{0,2\}$. Consider the plane $\Pi\defeq\overline{L\cup L'}$ in $X$. Since the lines $L,L'$ are concurrent, there exists a point $o\in L\cap L'$. By Theorem~\ref{t:proaffine-Fano=>3-regular}, the proaffine Fano liner $X$ is $3$-regular, which implies $\Pi=\bigcup_{x\in L}\bigcup_{y\in L'}\Aline xy$. If $\Pi$ is projective, then by Theorem~\ref{t:Gleason-Fano}, the finite Fano projective plane $\Pi$ is Pappian and hence $\|T\|=2$.

Now assume that the liner $\Pi$ is not projective. By Corollary~\ref{c:proaffineFano=>para-Playfair}, the proaffine Fano liner $X$ is para-Playfair and Proclus, and so is the plane $\Pi$ in $X$.  Since the Proclus plane $\Pi$ contains the line $L$ of length $|L|\ge|\{o,a,b,c\}|=4$, the liner $\Pi$ is $3$-long, by Theorem~\ref{t:Proclus-not-3long}. By Theorem~\ref{t:completelyF<=>}, the $3$-long finite proaffine Fano plane $\Pi$ is completely regular and its spread completion $\overline \Pi$ is a finite Fano projective liner.  By Gleason's Theorem~\ref{t:Gleason-Fano}, the finite Fano projective plane $\overline \Pi$ is Pappian. Since the horizon $\overline\Pi\setminus\Pi$ of $\Pi$ in $\overline\Pi$ is flat, the subliner $\Pi$ of the Pappian projective liner $\overline\Pi$ is Pappian, by Proposition~\ref{p:Pappian-minus-flat}. Then $\|T\|\in\{0,2\}$ in $\Pi$ and also in $X$. 
\end{proof}

\chapter{Little-Desarguesian and Moufang liners}\label{ch:Moufang}

In this chapter we introduce little-Desarguesian liners (as liners satisfying the little Desargues Axiom) and study the interplay between little-Desarguesian and Moufang liners. In particular, we show that a completely regular liner is Moufang if and only if it is little-Desarguesian if and only if its spread completion is a Moufang projective liner. Also we introduce para-Desarguesian liners and prove that a Playfair liner is Moufang if and only of it is para-Desarguesian.  

\section{Little-Desargues liners}

\begin{definition} A liner $X$ is defined to be \defterm{little-Desarguesian} if it satisfies the \defterm{Little Desargues Axiom}:
\begin{itemize}
\item[{\sf(lD)}] for any distinct lines $A,B,C,D\subset X$ and distinct points $o\in A\cap B\cap C\cap D$, $a,a'\in A$, $b,b'\in B$, $c,c'\in C$, if $\Aline ab\cap \Aline{a'}{b'}\cap D\ne\varnothing\ne \Aline bc\cap\Aline{b'}{c'}\cap D$, then $D\cap \Aline ac\subseteq\Aline{a'}{c'}$.
\end{itemize}

\begin{picture}(200,165)(-120,-30)
\put(0,0){\color{teal}\line(1,0){180}}
\put(120,120){\line(-1,-1){130}}
\put(120,120){\line(-1,-2){65}}
\put(120,120){\line(0,-1){130}}
\put(120,120){\line(1,-2){65}}
\put(80,80){\color{red}\line(1,0){60}}
\put(0,0){\color{red}\line(7,4){140}}
\put(80,80){\color{teal}\line(5,-4){100}}
\put(120,80){\color{cyan}\line(25.7,11.4){16}}
\put(60,0){\color{cyan}\line(7,8){76.4}}
{\linethickness{1pt}
\put(120,80){\color{cyan}\line(-25.7,-11.4){25.7}}
\put(60,0){\color{cyan}\line(7,8){60}}
\put(0,0){\color{teal}\line(1,0){60}}
\put(80,80){\color{teal}\line(5,-4){14}}
\put(80,80){\color{red}\line(1,0){40}}
\put(0,0){\color{red}\line(7,4){120}}
}

\put(-18,-20){$A$}
\put(50,-20){$B$}
\put(117,-20){$C$}
\put(184,-20){$D$}
\put(120,120){\circle*{3}}
\put(118,123){$o$}
\put(0,0){\circle*{3}}
\put(-6,3){$a'$}
\put(60,0){\circle*{3}}
\put(55,3){$b'$}
\put(180,0){\color{teal}\circle*{3}}
\put(80,80){\circle*{3}}
\put(76,83){$a$}
\put(120,80){\circle*{3}}
\put(113,83){$c$}
\put(94.3,68.6){\circle*{3}}
\put(86,63){$b$}
\put(120,68.6){\circle*{3}}
\put(122,63){$c'$}
\put(136.4,87.3){\color{cyan}\circle*{3}}
\put(140,80){\color{red}\circle*{3}}
\end{picture}
\end{definition}

The following two characterizations of little-Desarguesian liners and are counterparts of Propositions~\ref{p:Moufang<=>planeM} and \ref{p:Moufang<=>3plane} characterizing Moufang liners.

\begin{proposition}\label{p:lDes<=>plane} A liner $X$ is little-Desarguesian if and only if every plane in $X$ is little-Desarguesian.
\end{proposition}

\begin{proof} The ``only if'' part is trivial. To prove the ``if'' part, assume that every  plane in $X$ is little-Desarguesian. To prove that $X$ is little-Desarguesian, take any distinct lines $A,B,C,D$ in $X$  and distinct points $o\in A\cap B\cap C\cap D$, $a,a'\in A$, $b,b'\in B$, $c,c'\in C$ such that $\Aline ab\cap\Aline{a'}{b'}\cap D\ne \varnothing\ne \Aline bc\cap\Aline{b'}{c'}\cap D$. We have to prove that $\Aline ac\cap D\subseteq \Aline {a'}{c'}$. 

By the assumption, there exist points $x\in \Aline ab\cap\Aline{a'}{b'}\cap D$ and $y\in \Aline bc\cap\Aline{b'}{c'}\cap D$. It is easy to see that $x\notin \{o,a,b,a',b'\}$ and $y\notin \{o,b,c,b',c'\}$. Assuming that $x=y$, we conclude that $\Aline ab=\Aline xb=\Aline yb=\Aline cb$ and $\Aline {a'}{b'}=\Aline x{b'}=\Aline y{b'}=\Aline {c'}{b'}$. Then $\Aline ac\cap D=\Aline ab\cap D=\{x\}\subseteq \Aline{a'}{b'}=\Aline {a'}{c'}$. So, assume that $x\ne y$.

Observe that the plane $\Pi\defeq\overline{B\cup D}$ contains the points $a\in\Aline xb$, $a'\in \Aline x{b'}$, $c\in\Aline yb$ and $c'\in\Aline y{b'}$. By our assumption, the plane $\Pi$ is little-Desarguesian and hence $\Aline ac\cap D\subseteq\Aline{a'}{c'}$.
\end{proof}

For Proclus liners, Proposition~\ref{p:lDes<=>plane} can be improved as follows.

\begin{proposition}\label{p:lDes<=>3plane} A Proclus liner $X$ is little-Desargue if and only if every $3$-long plane in $X$ is little-Desarguesian.
\end{proposition}

\begin{proof} The ``only if'' part is trivial. To prove the ``only if'' part, assume that every $3$-long plane in a Proclus liner $X$ is little-Desarguesian. By Theorem~\ref{t:Proclus<=>}, the Proclus liner $X$ is $3$-proregular, and by Proposition~\ref{p:k-regular<=>2ex}, the $3$-proregular liner $X$ is $3$-ranked. By Proposition~\ref{p:lDes<=>plane}, the little-Desarguesian property of $X$ will follow as soon as we show that every plane $P$ in $X$ is little-Desarguesian. So, take any distinct lines $A,B,C,D$ in $P$ and distinct points $o\in A\cap B\cap C\cap D$, $a,a'\in A$, $b,b'\in B$, $c,c\in C$ such that $\Aline ab\cap\Aline{a'}{b'}\cap D\ne \varnothing\ne \Aline bc\cap\Aline{b'}{c'}\cap D$. We have to prove that $\Aline ac\cap D\subseteq \Aline {a'}{c'}$. Observe that $|P|\ge|\{o,a,b,c,a',b',c'\}|\ge 7$.

If the Proclus plane $P$ is not projective, then we can apply Theorem~\ref{t:Proclus-not-3long} and conclude that the plane $P$ is $3$-long. 
If  the plane $P$ is projective, then we can consider the maximal $3$-long flat $M\subseteq P$ that contains the point $o$. By Lemma~\ref{l:ox=2}, $\Aline op=\{o,p\}$ for all points $p\in P\setminus M$, which implies that $\{a,a',b,b',c,c'\}\subseteq M$ and hence $3=\|A\cup B\cup C\cup D\|\le\|M\|\le\|P\|=3$. The $3$-rankedness of the Proclus liner $X$ ensures that $M=P$ and hence the plane $P=M$ is $3$-long.

Therefore, in both cases the Proclus plane $P$ is $3$-long. By our assumption, the $3$-long plane $P$ is little-Desarguesian and hence $\Aline ac\cap D\subseteq\Aline {a'}{c'}$.
\end{proof}

\begin{proposition}\label{p:Moufang=>little-Desarguesian} Every Moufang liner is little-Desarguesian.
\end{proposition}

\begin{proof} Let $X$ be  a Moufang liner. To prove that $X$ is little-Desarguesian, take any distinct lines $A,B,C,D\subset X$ and distinct points $o\in A\cap B\cap C\cap D$, $a,a'\in A$, $b,b'\in B$, $c,c'\in C$ such that $\Aline ab\cap \Aline {a'}{b'}\cap D\ne\varnothing\ne \Aline bc\cap\Aline{b'}{c'}\cap D$. We have to prove that $\Aline ac\cap D\subseteq \Aline{a'}{c'}$. By our assumption, there exist points $x\in \Aline ab\cap\Aline{a'}{b'}\cap D$ and $y\in \Aline bc\cap\Aline{b'}{c'}\cap D$. 

If $b\in \Aline ac$, then $\Aline ab=\Aline ac=\Aline bc$ and $\{x\}=\Aline ab\cap D=\Aline bc\cap D=\{y\}$. Then $a',c'\in \Aline x{b'}=\Aline y{b'}$ and hence $\Aline ac\cap D=\Aline ab\cap D=\{x\}\subseteq \Aline {a'}{b'}=\Aline{a'}{c'}$. 
By analogy we can show that $b'\in\Aline{a'}{c'}$ implies $D\cap\Aline ac\subseteq\Aline{a'}{c'}$. 

So, assume that $abc$ and $a'b'c'$ are two triangles in $X$. In this case, the points $x,y\in D$ are distinct and hence $D=\Aline xy$. 

\begin{picture}(300,150)(-100,-15)
\linethickness{=0.6pt}
\put(0,0){\line(1,0){240}}
\put(60,120){\color{violet}\line(1,-2){60}}
\put(60,120){\color{blue}\line(-1,-2){60}}
\put(60,120){\line(0,-1){120}}
\put(30,60){\color{red}\line(7,-2){210}}
\put(30,60){\line(1,-2){30}}
\put(100,40){\line(-1,-1){40}}
\put(60,40){\color{violet}\line(3,-2){60}}
\put(60,40){\color{blue}\line(-3,-2){60}}
\put(45,30){\color{red}\line(13,-2){195}}

\put(84,24){\circle*{3}}
\put(81,28){$c'$}
\put(45,30){\circle*{3}}
\put(44,34){$a'$}
\put(0,0){\color{blue}\circle*{3}}
\put(-4,-8){\color{blue}$x$}
\put(60,0){\circle*{3}}
\put(57,-10){$o$}
\put(120,0){\color{violet}\circle*{3}}
\put(118,-8){\color{violet}$y$}
\put(240,0){\color{red}\circle*{3}}
\put(30,60){\circle*{3}}
\put(22,60){$a$}
\put(100,40){\circle*{3}}
\put(101,42){$c$}
\put(60,40){\circle*{3}}
\put(62,40){$b'$}
\put(60,120){\circle*{3}}
\put(58,123){$b$}

\end{picture}

Consider the triangles $xb'y$ and $aoc$, and observe that they are Moufang from the point $b$. moreover, the vertex $o$ of the triangle $aoc$ belongs to the side $\Aline xy$ of the triangle $xb'y$.  Since the liner $X$ is Moufang, the nonempty set $$T\defeq(\Aline x{b'}\cap\Aline {a}o)\cup(\Aline {b'}y\cap\Aline o{c})\cup(\Aline xy\cap\Aline {a}{c})=\{a',c'\}\cup(\Aline xy\cap\Aline{a}{c})$$ has rank $\|T\|=2$ and hence $D\cap\Aline ac=\Aline xy\cap\Aline{a}{c}\subseteq\Aline {a'}{c'}$, witnessing that the liner $X$ is little-Desarguesian.
\end{proof}

\begin{theorem}\label{t:projMoufang<=>lDes} A projective liner is Moufang if and only if it is little-Desarguesian.
\end{theorem}

\begin{proof} The ``only if'' part follows from Proposition~\ref{p:Moufang=>little-Desarguesian}. To prove the ``if' part, assume that $X$ is a little-Desarguesian projective liner. To prove that $X$ is Moufang, take any disjoint centrally perspective triangles $abc$ and $a'b'c'$ with $b'\in \Aline ac$. We have to prove that the set 
$$T\defeq(\Aline a{b}\cap\Aline {a'}{b'})\cup(\Aline {b}c\cap\Aline {b'}{c'})\cup(\Aline ac\cap\Aline {a'}{c'})$$ has rank $\|T\|\in\{0,2\}$. 
Since the triangles $abc$ and $a'b'c'$ are centrally perspective, there exists a unique point $o\in \Aline a{a'}\cap\Aline b{b'}\cap\Aline c{c'}$. Consider the plane $\Pi\defeq\overline{\{a,o,c\}}$ and observe that it contains the point $b'\in\Aline ac$ and also the set $\{a',b,c'\}\subseteq\Aline oa\cup\Aline o{b'}\cup \Aline oc$. By the projectivity of the plane $\Pi$, there exist unique points $x\in\Aline ab\cap\Aline{a'}{b'}$, $y\in\Aline bc\cap\Aline{b'}{c'}$ and $z\in\Aline ac\cap\Aline{a'}{c'}$.

\begin{picture}(300,150)(-100,-15)
\linethickness{=0.6pt}
\put(0,0){\line(1,0){240}}
\put(60,120){\color{violet}\line(1,-2){60}}
\put(60,120){\color{blue}\line(-1,-2){60}}
\put(60,120){\line(0,-1){120}}
\put(30,60){\color{red}\line(7,-2){210}}
\put(30,60){\line(1,-2){30}}
\put(100,40){\line(-1,-1){40}}
\put(60,40){\color{violet}\line(3,-2){60}}
\put(60,40){\color{blue}\line(-3,-2){60}}
\put(45,30){\color{red}\line(13,-2){195}}

\put(84,24){\circle*{3}}
\put(81,28){$y$}
\put(45,30){\circle*{3}}
\put(44,34){$x$}
\put(0,0){\color{blue}\circle*{3}}
\put(-4,-8){\color{blue}$a$}
\put(60,0){\circle*{3}}
\put(57,-10){$b'$}
\put(120,0){\color{violet}\circle*{3}}
\put(118,-8){\color{violet}$c$}
\put(240,0){\color{red}\circle*{3}}
\put(240,-8){\color{red}$z$}
\put(30,60){\circle*{3}}
\put(22,60){$a'$}
\put(100,40){\circle*{3}}
\put(101,42){$c'$}
\put(60,40){\circle*{3}}
\put(62,40){$b$}
\put(60,120){\circle*{3}}
\put(58,123){$o$}

\end{picture}

Now consider the distinct lines $A\defeq \Aline {b'}{a'}$, $B\defeq\Aline {b'}b$, $C\defeq\Aline {b'}{c'}$, $D\defeq\Aline ac$ containing the point $b'$. These lines contain the distinct points $a',x\in A$, $o,b\in B$, $c',y\in D$. Observe that
$$\Aline {a'}o\cap \Aline xb\cap D=\{a\}\ne \varnothing\ne \{c\}=\Aline o{c'}\cap\Aline by\cap D.$$ Since $X$ is little-Desarguesian, $\Aline {a'}{c'}\cap D\subseteq\Aline xy$ and hence $z\in \Aline ac\cap\Aline{a'}{c'}=D\cap\Aline{a'}{c'}\subseteq \Aline xy$. Since $z\in \Aline xy$, the set $T=\{x,y,z\}$ has rank $\|T\|=2$, witnessing that the projective liner $Y$ is Moufang.
\end{proof}

\begin{proposition}\label{p:lDes-subliner} Any subliner $X$ of a little-Desarguesian liner $Y$ is little-Desarguesian. 
\end{proposition}

\begin{proof}  To show that the liner $X$ is little-Desarguesian, take any distinct lines $A,B,C,D\subseteq X$ and distinct points $o\in A\cap B\cap C\cap D$, $a,a'\in A$, $b,b'\in B$, $c,c'\in C$ such that $\Aline ab\cap\Aline {a'}{b'}\cap D\ne\varnothing\ne \Aline bc\cap\Aline{b'}{c'}\cap D$. We have to show that $\Aline ac\cap D\subseteq\Aline {a'}{c'}$. To distinguish between lines in the lines $X$ and $Y$, for distinct points $x,y\in Y$, we shall denote by $\overline{\{x,y\}}$ the line in the liner $Y$ containing the points $x,y$.  The assumption $\Aline ab\cap\Aline {a'}{b'}\cap D\ne\varnothing\ne \Aline bc\cap\Aline{b'}{c'}\cap D$ implies $\overline{\{a,b\}}\cap\overline{\{a',b'\}}\cap \overline D\ne\varnothing\ne \overline{\{b,c\}}\cap\overline{\{b',c'\}}\cap \overline D$ and hence $\overline D\cap\overline{\{a,c\}}\subseteq\overline{\{a',c'\}}$, by the Little Desarges Axiom, holding for the little-Desarguesian liner $Y$. Then $$D\cap\Aline ac=\overline D\cap\overline{\{a,c\}}\cap X\subseteq\overline{\{a',c'\}}\cap X=\Aline{a'}{c'},$$
witnessing that the liner $X$ is little-Desarguesian.
\end{proof}

\begin{corollary} A liner $X$ is little-Desarguesian if it has a little-Desarguesian projective completion.
\end{corollary}

By Propositions~\ref{p:Desarg=>Moufang} and \ref{p:Moufang=>little-Desarguesian}, for every liner we have the implications:
$$\mbox{Desarguesian}\Ra\mbox{uno-Desarguesian}=\mbox{Moufang}\Ra\mbox{little-Desarguesian}.
$$

None of those implications can be reversed.

\begin{Exercise}[Moufang, 1931] Find an example of a Playfair plane, which is Moufang but not Desarguesian.
\smallskip

{\em Hint:} Consider the octonion plane $\mathbb O\times\mathbb O$.
\end{Exercise}

\begin{exercise} Find an example of a Proclus plane, which is little-Desarguesian but not Moufang.
\smallskip

{\em Hint:} Take any $4$-long Moufang projective plane $Y$ and choose any line $L\subset Y$ and point $p\in L$. Applying Theorem~\ref{t:proaffine<=>proflat} and Proposition~\ref{p:nD=>pP}, show that the subliner $X\defeq Y\setminus (L\setminus\{p\})$ is Proclus, little-Desarguesian, but not Moufang.
\end{exercise}



\begin{corollary}\label{c:Steiner=>little-Desarg} Every Steiner $3$-regular liner is Moufang and little-Desarguesian.
\end{corollary}

\begin{proof} Let $X$ be a Steiner $3$-regular liner. By Theorem~\ref{t:Steiner-3reg},  $X$ is Proclus (more precisely, $X$ is $p$-parallel for some $p\in\{0,1\}$). By  Proposition~\ref{p:Steiner+projective=>Desargues}, the Steiner $3$-regular liner $X$ is Desarguesian. By Propositions~\ref{p:Desarg=>Moufang} and  \ref{p:Moufang=>little-Desarguesian}, the Desarguesian Proclus liner $X$ is Moufang and little-Desarguesian.
\end{proof}

\begin{example}[Ivan Hetman, 2025] The cyclic group $\IZ_{25}=\{0,1,\dots,24\}$ endowed with the family of lines $\mathcal L\defeq\big\{x+L:x\in \IZ_{25},\;L\in\big\{\{0, 1,12\},\{0,2,18\},\{0,3,8\},\{0,4,19\}\big\},$
is a Steiner ranked plane, which is not little-Desarguesian. For the concurrent lines
$$A\defeq \{0,3,8\},\quad B\defeq\{0,5,22\},\quad C\defeq\{0,11,24\},\quad D\defeq\{0,15,21\}$$and points $$a\defeq 3,\;a'\defeq 8,\quad b\defeq 5,\;b'\defeq 22,\quad c\defeq 11,\;c'\defeq 24$$we have $$21\in \Aline ab\cap\Aline {a'}{b'}\cap D\quad\mbox{and}\quad 15\in \Aline bc\cap\Aline {b'}{c'}\cap D,\quad\mbox{but}\quad \Aline ac\cap\Aline{a'}{c'}=\{6\}\not\subseteq D,$$ witnessing that the Steiner plane $(\IZ_{25},\mathcal L)$ is not little-Desarguesian.
\end{example}

\begin{exercise} Find an example of a little-Desarguesian Playfair liner which is not regular.
\vskip5pt

{\em Hint:} Take the non-regular Hall liner from Example~\ref{ex:HTS} and apply Corollary~\ref{c:Steiner=>little-Desarg}.
\end{exercise}

\section{Projective completions of little-Desarguesian liners}

The main result of this section is the following (difficult) theorem on the preservation of little-Desarguesianity by projective completions. 

\begin{theorem}\label{t:lDes-projcompletion} A projective completion of any little-Desarguesian liner is Moufang.
\end{theorem}

\begin{proof} Let $Y$ be a projective completion of a little-Desargusian liner $X$.
If $\|Y\|\ne 3$, then the $3$-long projective liner $Y$ is Desarguesian and Moufang, by Theorem~\ref{t:proaffine-Desarguesian} and Proposition~\ref{p:Desarg=>Moufang}. So, assume that $\|Y\|=3$ and hence $Y$ is a projective plane. By Corollary~\ref{c:Avogadro-projective}, the $3$-long projective liner $Y$ is $2$-balanced. If $|Y|_2\le 5$, then $Y$ is Desarguesian, by Corollary~\ref{c:p5-Pappian}. So, assume that $|Y|_2\ge 6$. Since no projective plane of order $6$ exists, $|Y|_2\ge 7$.  Since $Y$ is a projective completion of $X$, $\overline{Y\setminus X}\ne Y$ and hence $Y\setminus X\subseteq H$ for some line $H$ in $Y$. By Theorems~\ref{t:affine<=>hyperplane} and \ref{t:Playfair<=>}, the subliner $Y\setminus H$ is Playfair. By Proposition~\ref{p:lDes-subliner}, the subliner $X\setminus H$ of the little-Desarguesian liner $X$ is little-Desarguesian. Replacing the liner $X$ by its Playfair subliner $X\setminus H$, we can assume that the liner $X$ is Playfair and its horizon $H\defeq Y\setminus Y$ is a line in $Y$. Since $Y$ is $7$-long, the Playfair plane $X$ is $6$-longs.

\begin{lemma}\label{l:lDes=>quadratic} Let $A,B,C,D$ be distinct lines in the little-Desarguesian Playfair plane $X$ and $o\in A\cap B\cap C\cap D$, $a\in A$, $b,b'\in B$, $c\in C$, $x\in D\cap \Aline ab$, $y\in D\cap\Aline bc$ be distinct points such that $\Aline xb\parallel \Aline oc\parallel \Aline {b'}y$ and $\Aline {b'}x\parallel \Aline oa$. Then $\Aline ac\parallel \Aline xy$.

\begin{picture}(200,85)(-190,-10)
\linethickness{=0.6pt}
\put(0,0){\line(0,1){60}}
\put(0,0){\color{blue}\line(1,1){30}}
\put(0,0){\color{violet}\line(-1,1){30}}
\put(-30,30){\color{blue}\line(1,1){30}}
\put(-30,30){\color{red}\line(1,0){60}}
\put(30,30){\line(-1,1){30}}
\put(0,30){\color{blue}\line(1,1){15}}
\put(0,30){\color{violet}\line(-1,1){15}}
\put(-15,45){\color{red}\line(1,0){30}}

\put(0,0){\circle*{3}}
\put(-3,-10){$b'$}
\put(-30,30){\circle*{3}}
\put(-38,28){$x$}
\put(0,30){\circle*{3}}
\put(2,23){$o$}
\put(30,30){\circle*{3}}
\put(32,28){$y$}
\put(-15,45){\circle*{3}}
\put(-22,47){$a$}
\put(15,45){\circle*{3}}
\put(16,47){$c$}
\put(0,60){\circle*{3}}
\put(-2,63){$b$}
\end{picture}

\end{lemma}

\begin{proof} To derive a contradiction, assume that $\Aline ac\nparallel D$. Since $X$ is Playfair, and $c\notin D$, there exists a unique point $\alpha\in \Aline xb$ such that $\Aline\alpha c\cap D=\varnothing$. Since $\Aline ac\nparallel D\parallel \Aline \alpha c$, the point $\alpha$ is not equal to the point $a$ and hence the line $\Aline \alpha o$ is not parallel to the lines $\Aline x{b'}\parallel \Aline oa$. Then there exists a point $a'\in \Aline x{b'}\cap \Aline \alpha o$. Since $\alpha\ne x$, the points $x$ and $a'$ are distinct and hence $a'\notin \Aline xy$. Since $X$ is Playfair, there exists a unique line $L'$ such that $a'\in L'\subseteq X\setminus D$. Since $\Aline oc\cap\Aline xy=\{o\}$, there exists a unique point $c'\in L'\cap\Aline oc$.
Assuming that $b'\in \Aline {a'}{c'}$, we conclude that $x\in \Aline {a'}{b'}\subseteq \Aline {a'}{c'}$ and hence the parallel lines $D$ and $\Aline {a'}{c'}$ coincide.
Then $b'\in D$, which contradicts the choice of the point $b'\in B\setminus\{o\}=B\setminus D$. This contradiction shows that $b'\notin\Aline{a'}{c'}$. By the Proclus Postulate~\ref{p:Proclus-Postulate}, there exists a point $y'\in \Aline {b'}{c'}\cap D$. It follows from $a'\ne c'$ that $y'\ne x$ and hence the lines $\Aline b{y'}$ and $\Aline bx$ are concurrent. Since $\Aline bx\parallel \Aline oc$, the lines $\Aline b{y'}$ and $\Aline oc$ are concurent and hence there exists a point $\gamma\in \Aline oc\cap\Aline b{y'}$. It follows from $\Aline {b'}y\cap \Aline oc=\varnothing\ne\Aline{b'}{y'}\cap\Aline oc=\{c'\}$ that $y'\ne y$ and hence $\gamma\ne c$ and $\Aline \alpha \gamma\ne\Aline \alpha c\parallel D$. Then $\Aline \alpha\gamma\cap D\ne\varnothing$.  

\begin{picture}(200,180)(-220,-15)
\linethickness{=0.7pt}
\put(0,0){\line(0,1){120}}
\put(0,0){\color{blue}\line(1,1){60}}
\put(0,0){\color{violet}\line(-1,1){150}}
\put(-60,60){\color{blue}\line(1,1){60}}
\put(-60,60){\color{red}\line(1,0){120}}
\put(60,60){\line(-1,1){60}}
\put(0,60){\color{blue}\line(1,1){90}}
\put(0,60){\line(-5,3){150}}
\put(-150,150){\color{red}\line(1,0){240}}
\put(0,0){\line(3,5){90}}
\put(0,120){\line(3,-5){36}}
\put(0,60){\color{violet}\line(-1,1){30}}
\put(30,90){\color{red}\line(-9,-1){67.5}}

\put(0,0){\circle*{3}}
\put(-3,-10){$b'$}
\put(-60,60){\circle*{3}}
\put(-67,53){$x$}
\put(0,60){\circle*{3}}
\put(-8,53){$o$}
\put(60,60){\circle*{3}}
\put(61,53){$y$}
\put(36,60){\circle*{3}}
\put(35,51){$y'$}
\put(-37.5,82.5){\circle*{3}}
\put(-48,80){$\alpha$}
\put(22.5,82.5){\circle*{3}}
\put(13.5,80.5){$\gamma$}
\put(-30,90){\circle*{3}}
\put(-37,92){$a$}
\put(30,90){\circle*{3}}
\put(28,94){$c$}
\put(0,120){\circle*{3}}
\put(-2,123){$b$}
\put(-150,150){\circle*{3}}
\put(-153,153){$a'$}
\put(90,150){\circle*{3}}
\put(89,153){$c'$}
\end{picture}

Consider the lines $\Lambda\defeq \Aline o\alpha$, $B\defeq \Aline ob$, $C\defeq\Aline oc$, $D\defeq \Aline xy$ and observe that $o\in \Lambda\cap B\cap C\cap D$, $\alpha,a'\in \Lambda$, $b,b'\in B$, $\gamma,c'\in C$ are distinct points such that 
$$\Aline \alpha b\cap\Aline {a'}{b'}\cap D=\{x\}\ne\varnothing\ne\{y'\}=\Aline b\gamma\cap\Aline {b'}{c'}\cap D.$$ Since $X$ is little-Desarguesian, this implies $\varnothing\ne \Aline \alpha\gamma\cap D\subseteq \Aline {a'}{c'}=L'$, which contradicts $L'\cap\Aline xy=\varnothing$.  This is the final contradiction showing that $\Aline ac\parallel\Aline xy$.
\end{proof}

For a line $L$ in the projective plane $Y$ and distinct points $a,a'\in Y\setminus L$, let  $o$ be the unique point of the intersection $L\cap\Aline a{a'}$ and $\Phi_L^{aa'}:Y\setminus\Aline a{a'}\to Y\setminus\Aline a{a'}$ be the bijective map assigning to each point $x\in Y\setminus\Aline a{a'}$ the unique point $y\in \Aline ox$ such that $\Aline y{b'}\cap \Aline xb\cap L\ne \varnothing$. This definition implies $\Phi_L^{aa'}(x)=x$ for every $x\in L\setminus\Aline a{a'}$, and $\Phi_L^{aa'}[\Lambda]=\Lambda$ for every line $\Lambda\subset Y$ containing the point $o$.

In the following lemmas we shall assume that $o,a,a'\in X$, and in this case we shall prove that the bijective map $\Phi_L^{aa'}$ can be extended to a central automorphism of the projective plane $Y$. Since $o\in L\cap X$, the set $L\cap H=L\setminus X$ is a singleton whose unique point will be denoted by $h$. First we prove that the bijection $\Phi_L^{aa'}$ preserves the collinearity of points in the set $Y\setminus (\Aline ha\cup\Aline oa)$.

\begin{lemma}\label{l:lD-lines1} For any line $\Lambda$ in $Y$ with $h\notin \Lambda$, there exists a line $\Lambda'\subset Y$ such that $h\notin \Lambda'$ and $\Phi_L^{aa'}[\Lambda\setminus(\Aline ha\cup \Aline oa)]\subseteq \Lambda'$.
\end{lemma}

\begin{proof} If $o\in \Lambda$, then the line $\Lambda'\defeq \Lambda$ has the required property. So, assume that $o\notin\Lambda$ and hence $L\ne\Lambda\ne\Aline oa$. By the projectivity of $Y$, there exists a  point $\lambda\in\Lambda\cap L\subseteq L\setminus\{h\}\subseteq X$. If $a\in\Lambda$, then the definition of the function $\Phi_L^{aa'}$ ensures that $\Phi_L^{aa'}[\Lambda\setminus\Aline a{a'}]=\Phi_L^{aa'}[\Aline\lambda a\setminus\Aline a{a'}]=\Aline \lambda{a'}\setminus\Aline a{a'}$. Therefore, the line $\Lambda'\defeq\Aline \lambda{a'}$ has the required property. So, assume that $a\notin \Lambda$. 

\begin{claim}\label{cl:lD-exist-bxb'} There exist points $b\in (X\cap \Lambda)\setminus(L\cup\Aline oa)$, $x\in \Aline ab\cap L\cap X$, and $b'\in \Aline x{a'}\cap\Aline ob\cap X$.
\end{claim}

\begin{proof} Since $X$ is a Playfair plane and $\Lambda\cap L=\{\lambda\}\subseteq X$, there exist unique points $\mu\in L\cap X$ such that the lines $\Aline a\mu\cap X$ and $\Lambda\cap X$ in $X$ are parallel.  Since $X$ is Proclus, the set $$I\defeq\{b\in \Lambda:\Aline ab\cap L\cap X=\varnothing\}\cup\{b\in\Lambda:\Aline ab\cap \Aline {a'}\mu\cap X=\varnothing\}\cup \{b\in\Lambda:\Aline ab\cap \Aline {a'}{\lambda}\cap X=\varnothing\}$$ contains at most three points. Since $X$ is $6$-long, there exists a point $b\in (X\cap \Lambda)\setminus(I\cup\{\lambda\}\cup\Aline oa)$. Since $b\notin I$, the unique point $x\in \Aline ab\cap L$ belongs to the set $X$. By the projectivity of the plane $Y$, there exists a unique point $b'\in \Aline ob\cap \Aline x{a'}\subseteq Y$.  If $b'\in X$, then we are done. So, assume that $b'\in Y\setminus X$.

\begin{picture}(200,150)(-200,-15)

\linethickness{=0.7pt}
\put(0,0){\line(0,1){120}}
\put(0,0){\color{blue}\line(1,1){75}}
\put(0,0){\color{violet}\vector(-1,1){90}}
\put(0,60){\color{violet}\vector(-1,1){60}}
\put(-60,60){\color{blue}\line(1,1){60}}
\put(-60,60){\color{teal}\line(1,0){150}}
\put(-35,63){\color{teal}$L$}
\put(0,120){\line(1,-1){60}}
\put(0,60){\color{blue}\line(1,1){30}}
\put(30,90){\color{teal}\line(-1,0){60}}
\put(0,60){\color{red}\line(5,1){75}}
\put(90,60){\color{cyan}\line(-4,1){150}}

\put(-90,115){\color{violet}$b'$}

\put(0,0){\circle*{3}}
\put(-3,-8){$a'$}
\put(50,70){\circle*{3}}
\put(50,73){$c$}
\put(75,75){\circle*{3}}
\put(78,75){$c'$}

\put(-60,60){\circle*{3}}
\put(-67,53){$x$}
\put(0,60){\circle*{3}}
\put(-7,53){$o$}
\put(60,60){\circle*{3}}
\put(57,51){$y$}
\put(90,60){\circle*{3}}
\put(93,58){$\lambda$}
\put(-70,95){\color{cyan}$\Lambda$}
\put(-30,90){\circle*{3}}
\put(-32,93){$b$}
\put(30,90){\circle*{3}}
\put(31,92){$\gamma$}
\put(0,120){\circle*{3}}
\put(-2,123){$a$}
\end{picture}

Since the liner $X$ is Playfair, there exists a unique point $y\in L\cap X$ such that $\Aline y{a'}\cap X\parallel \Aline xa\cap X$. The choice of $b\notin I$ ensures that $\Aline ab\cap X\nparallel \Aline{a'}{\lambda}\cap X$ and hence $y\ne \lambda$. The choice of $b\notin\Aline oa$ ensures that $x\notin \Aline a{a'}$ and hence the parallel lines $\Aline y{a'}$ and $\Aline xa$ are disjoint. Since $X$ is Playfair, there exists a unique point $\gamma\in \Aline ay\cap X$ such that $\Aline o\gamma\cap X\parallel \Aline xa\cap X\parallel \Aline {a'}{y}\cap X$. By Lemma~\ref{l:lDes=>quadratic}, $\Aline b\gamma\cap X\parallel L\cap X$. Assuming that $\Aline ay\cap X \parallel \Lambda\cap X$, we conclude that $y=\mu$ and $\Aline ab\cap X\parallel  \Aline {a'}y\cap X=\Aline {a'}\mu\cap X$, which contradicts the choice of the point $b\notin I$. This contradiction shows that $\Aline ay\cap X\nparallel \Lambda\cap X$ and hence there exists a point $c\in \Aline ay\cap\Lambda\cap X$. It follows from $y\ne\lambda$ that $c\notin L$. Since $a\notin \Lambda$, the point $c$ does not belong to the line $\Aline oa$. Assuming that $c\in\Aline ob$, we conclude that $y\in \Aline ac\cap L=\Aline ab\cap L=\{x\}$, which contradicts $\Aline {a'}y\cap \Aline xa\cap X=\varnothing$. Therefore, the lines $\Aline oc$, $\Aline oa$, $\Aline ob$ and $L$ are distinct. Since $\Aline b\gamma\cap X\parallel L\cap X\nparallel \Lambda \cap X=\Aline b c\cap X$, the points $\gamma$ and $c$ are distinct. Then $\Aline oc\cap X\nparallel \Aline o\gamma\cap X\parallel \Aline {a'}y\cap X$ and there exists a point $c'\in \Aline oc\cap \Aline {a'}y\cap X$. In this case, the points $c,y,c'$ have the required property.
\end{proof}

By Claim~\ref{cl:lD-exist-bxb'}, there exist points $b\in \Lambda\cap X\setminus(L\cup\Aline oa)$, $x\in \Aline ab\cap L\cap X$ and $b'\in\Aline x{a'}\cap \Aline ob\cap X$.  We claim that the line $\Lambda'\defeq\Aline \lambda{b'}$ has the required property: $\Phi_L^{aa'}[\Lambda\setminus(\Aline ha\cup\Aline oa)]\subseteq\Lambda'$. The definition of the function $\Phi_L^{aa'}$ ensures that $\Phi_L^{aa'}(b)=b'$. 
Given any pair $cc'\in \Phi_L^{aa'}$ with $c\in \Lambda\setminus (\Aline h{a}\cup \Aline oa)$, we should check that $c'\in \Lambda'$. This is clear if $c=b$. So, assume that $c\ne b$. By the projectivity of $Y$, there exists a unique point $y\in \Aline ac\cap L$. Since $c\notin \Aline ha$, the point $y$ is not equal to $h$ and hence $y\in L\setminus\{h\}\subseteq X$. The definition of the function $\Phi_L^{aa'}$ensures that $c'\in \Aline oc\cap\Aline y{a'}$. Then 
$$\Aline ab\cap \Aline {a'}{b'}\cap L\cap X=\{x\}\ne \varnothing\ne \{y\}=\Aline ac\cap\Aline{a'}{c'}\cap L\cap X.$$ 

Since $o\in \Aline c{c'}\cap X$, the set $\Aline c{c'}\setminus X$ is a singleton and hence $c\in X$ or $c'\in X$. 

If $c,c'\in X$, then the Little Desargues Axiom
ensures that $\{\lambda\}=\Aline bc\cap L\cap X\subseteq \Aline {b'}{c'}\cap X$ and hence $c'\in \Aline \lambda{b'}=\Lambda'$. 

If $c\in X$ and $c'\notin X$, then $\Aline y{a'}\cap X$ and $\Aline oc\cap X$ are disjoint parallel lines in the Playfair liner $X$. Assuming that the lines $\Aline \lambda{b'}\cap X$ and $\Aline oc\cap X$ are not parallel, we can find a point $c''\in \Aline \lambda{b'}\cap\Aline oc\cap X$. The inequlities $b\ne b'\ne x$ imply $c\ne c''\ne o$. Observe that $a,a',b,b',c,c''\in X$ and $\Aline ab\cap\Aline {a'}{b'}\cap L\cap X=\{x\}\ne \varnothing\ne\{\lambda\}=\Aline bc\cap\Aline{b'}{c''}\cap L\cap X$. Since $X$ is little-Desarguesian, this implies $\{y\}=L\cap\Aline ac\subseteq \Aline{a'}{c''}$ and hence $c''\in \Aline y{a'}\cap X$, which contradicts $\Aline y{a'}\cap X\cap \Aline oc=\varnothing$. This contradiction shows that $\Aline \lambda{b'}\cap X\parallel \Aline oc\cap X\parallel\Aline y{a'}\cap X$ and hence $c'\in \Aline \lambda{b'}=\Lambda'$.

By analogy we can show that $c'\in \Lambda'$ whenever $c'\in X$ and $c\notin X$.
\end{proof}

Lemma~\ref{l:lD-lines1} admits the following self-improvements.

\begin{lemma}\label{l:lD-lines2} For every line $\Lambda$ in $Y$ with $h\notin \Lambda$, there exists a line $\Lambda'\subset Y$ such that $$\Phi_L^{aa'}[\Lambda\setminus(\Aline ha\cup \Aline oa)]=\Lambda'\setminus(\Aline h{a'}\cup\Aline o{a'}).$$
\end{lemma}

\begin{proof} The definition of the functions $\Phi_L^{aa'}$ and $\Phi_L^{a'\!a}$ ensures that $(\Phi_L^{aa'})^{-1}=\Phi_L^{a'\!a}$ and $\Phi_L^{aa'}[\Aline ha\setminus\Aline a{a'}]=\Aline h{a'}\setminus\Aline a{a'}$. 
If $o\in \Lambda$, then $\varphi[\Lambda\setminus\Aline a{a'}]=\Lambda\setminus\Aline a{a'}$ and hence the line $\Lambda'\defeq\Lambda$ has the required property. So, assume that $o\notin\Lambda$. Since $h\notin\Lambda$, we have also that $\Lambda\ne\Aline ha$. By Lemma~\ref{l:lD-lines1} (applied to functions $\Phi_L^{aa'}$ and $\Phi_L^{a'\!a}$), there exist lines $\Lambda'$ and $\Lambda''$ such that  $\Phi_L^{aa'}[\Lambda\setminus(\Aline ha\cup \Aline oa)]\subseteq\Lambda'$ and $\Phi_L^{a'\!a}[\Lambda'\setminus(\Aline h{a'}\cup\Aline o{a'})]\subseteq\Lambda''$. Since the liner $Y$ is $5$-long and $\Aline ha\ne \Lambda\ne \Aline oa$, there exist two distinct points $x,y\in \Lambda\setminus(L\cup\Aline \lambda a\cup\Aline oa)$. Observe that  $\{\Phi_L^{aa'}(x),\Phi_L^{aa'}(y)\}\subseteq \Lambda'\setminus(\Aline \lambda{a'}\cup\Aline o{a'})$ and $\{x,y\}=\{\Phi_L^{a'\!a}(\Phi_L^{aa'}(x)),\Phi_L^{a'\!a}(\Phi_L^{aa'}(y))\}\subseteq \Lambda\cap\Lambda''$. Then $\Lambda''=\Lambda$ and hence $$\Phi_L^{a'\!a}[\Lambda'\setminus(\Aline h{a'}\cup\Aline o{a'})]\subseteq \Lambda''\setminus(\Aline h{a}\cup\Aline oa)]=\Lambda\setminus(\Aline ha\cup\Aline oa).$$ Applying to this inclusion the function $\Phi_L^{aa'}$ and taking into account that $\Phi_L^{aa'}\circ\Phi_L^{a'\!a}$ is the identity function, we conclude that  $\Lambda'\setminus(\Aline h{a'}\cup\Aline o{a'})\subseteq \Phi_L^{aa'}[\Lambda\setminus(\Aline h{a}\cup\Aline oa)]\subseteq \Lambda'\setminus(\Aline h{a'}\cup\Aline o{a'})$, 
which yields the required equality $\Phi_L^{aa'}[\Lambda\setminus(\Aline h{a}\cup\Aline oa)]=\Lambda'\setminus(\Aline h{a'}\cup\Aline o{a'})$.
\end{proof}

\begin{lemma}\label{l:lD-lines3} For every line $\Lambda$ in $Y$ with $h\notin \Lambda$, there exists a line $\Lambda'$ such that $$\Phi_L^{aa'}[\Lambda\setminus \Aline a{a'}]=\Lambda'\setminus\Aline a{a'}.$$
\end{lemma}

\begin{proof} If $o\in \Lambda$, then the line $\Lambda'=\Lambda$ has the required property. So, assume that $o\notin\Lambda$. By the projectivity of $Y$, there exist  unique points $\alpha\in \Lambda\cap\Aline a{a'}$ and $\lambda\in \Lambda\cap L$. It follows from $o,h\notin \Lambda$ that $\Lambda\cap L=\{\lambda\}\ne\{h\}$. By Lemma~\ref{l:lD-lines2}, for the line $\Lambda$, there exists a line $\Lambda'$  such that $\Phi_L^{aa'}[\Lambda\setminus(\Aline ha\cup\Aline a{a'})]=\Lambda'\setminus(\Aline h{a'}\cup\Aline a{a'})$. Observe that $\lambda\in \Lambda\setminus(\Aline ha\cup\Aline a{a'})$ and hence $\lambda=\Phi_L^{aa'}(\lambda)\in \Phi_L^{aa'}[\Lambda\setminus(\Aline ha\cup\Aline a{a'})]\subseteq\Lambda'$. 
Assuming that $\Phi_L^{aa'}[\Lambda\setminus\Aline a{a'}]\not\subseteq\Lambda'$, we can find a point $x\in (\Lambda\cap\Aline ha)\setminus\Aline a{a'}$ such that $y\defeq \Phi_L^{aa'}(x)\notin\Lambda'$. Then $y=\Phi_L^{aa'}(x)\in \Aline ox\cap\Aline h{a'}$, by definition of the function $\Phi^{aa'}_L$. By the projectivity of $Y$, there exists a unique point $z\in \Aline ox\cap\Lambda'$. Assuming that $z\in \Aline a{a'}$, we conclude that $z\in \Aline a{a'}\cap \Aline ox=\{o\}$ and hence $\Lambda'=\Aline o\lambda=L$ and $\Lambda\setminus\{x,\alpha\}=\Lambda\setminus(\Aline ha\cup\Aline a{a'})\subseteq \Phi_L^{a'\!a}[\Lambda'\setminus\{o\}]=\Phi^{a'\!a}_L[L\setminus\{o\}]\subseteq L$, which is impossible because the liner $Y$ is $4$-long and $x\in \Lambda\setminus L$. Therefore, $z\in \Lambda'\setminus(\Aline h{a'}\cup\Aline a{a'})=\Phi_L^{aa'}[\Lambda\setminus(\Aline ah\cup\Aline aa')]=\Phi_L^{aa'}[\Lambda\setminus\{x,\alpha\}]$. Then $\Phi_L^{a'\!a}(z)\in \Aline oz\cap (\Lambda\setminus\{x,\alpha\})=\Aline ox\cap(\Lambda\setminus\{x,\alpha\})=\varnothing$. This contradiction shows that $\Phi_L^{aa'}[\Lambda\setminus\Aline a{a'}]\subseteq\Lambda'\setminus\Aline a{a'}$. 

To see that $\Phi_L^{aa'}[\Lambda\setminus\Aline a{a'}]=\Lambda'\setminus\Aline a{a'}$, take any point $c'\in \Lambda'\setminus\Aline a{a'}$ and using the projectivity of the plane $Y$, find a unique point $c\in\Lambda\cap\Aline o{c'}$. The definition of the function $\Phi_L^{aa'}$ ensures that $\Phi_L^{aa'}(c)\in \Aline oc\cap\Phi_L^{aa'}[\Lambda\setminus\Aline a{a'}]\subseteq \Aline o{c'}\cap\Lambda'=\{c'\}$, and hence $\Phi_L^{aa'}[\Lambda\setminus\Aline a{a'}]=\Lambda'\setminus\Aline a{a'}$.
\end{proof}

\begin{lemma}\label{l:lD-lines4} For every line $\Lambda$ in $Y$, there exists a line $\Lambda'$ in $Y$ such that $\Psi_L^{aa'}[\Lambda\setminus \Aline a{a'}]=\Lambda'\setminus \Aline a{a'}$.
\end{lemma}

\begin{proof} If $o\in \Lambda$, then $\Lambda'\defeq\Lambda$ is a unique line with the property $\Phi_L^{aa'}[\Lambda\setminus\Aline a{a'}]=\Lambda'\setminus\Aline a{a'}$. 

So, assume that $o\notin\Lambda$ and hence $\Lambda\ne L$. 
If $h\notin \Lambda$, then by Lemma~\ref{l:lD-lines3}, there exists a line $\Lambda'$ such that $\Phi_L^{aa'}[\Lambda\setminus\Aline a{a'}]=\Lambda'\setminus\Aline a{a'}$. 

So, assume that $h\in \Lambda$. If $a\in \Lambda$, then the line $\Lambda'\defeq\Lambda$ has the desired property $\Phi_L^{aa'}[\Lambda\setminus\Aline a{a'}]=\Lambda'\setminus\Aline a{a'}$, by the definition of the function $\Phi_L^{aa'}$. So, assume that $a\notin \Lambda$ and hence $\Lambda\cap\Aline ha=\{h\}$.

Choose any point $b\in \Lambda\setminus(\{h\}\cup\Aline a{a'})$ and observe that  $b\notin\Aline ha$ implies $b'\defeq\Phi_L^{aa'}(b)\notin \Phi_L^{aa'}[\Aline ha]=\Aline h{a'}$. We claim that the line $\Lambda'\defeq\Aline h{b'}$ has the desired property.  To derive a contradiction, assume that for some point $c\in \Lambda\setminus \Aline a{a'}$ its image $c'\defeq\Phi_L^{aa'}(c)$ does not belong to the line $\Lambda'$. Then the unique point $z\in \Aline {b'}{c'}\cap L$ is not equal to $h$ and hence $z\in L\setminus\{h\}\subseteq X$. By Lemma~\ref{l:lD-lines3}, for the line $\Omega\defeq \Aline zc$, there exists a line $\Omega'$ in $Y$ such that $\Phi_L^{aa'}[\Omega\setminus\Aline a{a'}]=\Omega'\setminus\Aline a{a'}$. Then $\{z,c'\}\subseteq \Phi_L^{aa'}[\Omega\setminus\Aline a{a'}]\subseteq\Omega'$ and hence $\Omega'=\Aline z{c'}$. Then $b=\Phi_L^{a'\!a}(b')\in \Phi_L^{a'\!a}[\Omega'\setminus\Aline a{a'}]=\Omega\setminus \Aline oa$ and hence $z\in \Omega\cap L=\Aline bc\cap L=\Lambda\cap L=\{h\}$, which is a contradiction showing that $\Phi_L^{aa'}[\Lambda\setminus\Aline a{a'}]\subseteq \Lambda'$. To see that $\Phi_L^{aa'}[\Lambda\setminus\Aline a{a'}]=\Lambda'\setminus\Aline a{a'}$, take any point $c'\in \Lambda'\setminus \Aline a{a'}$ and find a unique point $c\in \Lambda\cap \Aline o{c'}$. Then $\Phi_L^{aa'}(c)\in \Aline o{c'}\cap\Phi_L^{aa'}[\Lambda\setminus\Aline oa]\subseteq\Lambda'\cap\Aline o{c'}=\{c'\}$, witnessing that $\Phi_L^{aa'}[\Lambda\setminus\Aline a{a'}]=\Lambda'\setminus\Aline a{a'}$.
%
 \end{proof}

By Lemma~\ref{l:lD-lines4}, for every point $x\in Y\setminus\{h\}$, the flat hull $$\bar\Phi_L^{aa'}[\Aline ha]\defeq\overline{\Phi_L^{aa'}[\Aline hx\setminus\Aline a{a'}]}$$ is a line containing the point $h$. This allows us to extend the bijective function $\Phi_L^{aa'}:Y\setminus\Aline a{a'}\to Y\setminus \Aline a{a'}$ to a bijective function $\bar\Phi_L^{aa'}:Y\to Y$ assigning to every point $x\in\Aline a{a'}$ the unique point $\bar\Phi_L^{aa'}(x)$ of the intersection $\bar\Phi_L^{aa'}[\Aline ha]\cap\Aline ox$. The inverse function to $\bar\Phi_L^{aa'}$ is the function $\bar\Phi_L^{a'\!a}$.

\begin{lemma}\label{l:lD-lines5} The function $\bar\Phi_L^{aa'}$ is an automorphism of the projective plane $Y$.  
\end{lemma}

\begin{proof} By Theorem~\ref{t:liner-isomorphism<=>}, it suffices to show that for every line $\Lambda\subseteq Y$, the set $\bar\Phi_L^{aa'}[\Lambda]$ is a line. By Lemma~\ref{l:lD-lines4}, there exists a line $\Lambda'$ such that $\Phi_{L}^{aa'}[\Lambda\setminus\Aline a{a'}]=\Lambda'\setminus\Aline a{a'}$. We claim that $\bar\Phi_{L}^{aa'}[\Lambda]=\Lambda'$. If $\Lambda=\Aline a{a'}$, then $\Lambda'=\Lambda$ and $\bar\Phi_L^{aa'}[\Lambda]=\Aline a{a'}=\Lambda'$. So, we assume that $\Lambda\ne\Aline a{a'}$.  In this case $\Lambda'=\overline{\Phi_L^{aa'}[\Lambda\setminus\Aline a{a'}]}$. 
Let $x$ be the unique common point of the lines $\Lambda$ and $\Aline a{a'}$. Assuming that $\bar\Phi_L^{aa'}[\Lambda]\ne\Lambda'$, we conclude that $y\defeq \bar\Phi_L^{aa'}(x)\notin\Lambda'$. In this case, the definition of the extension $\bar\Phi_L^{aa'}$ guarantees that $h\notin\Lambda$ and hence $h\notin\Lambda'$. By the projectivity of $Y$, there exists a unique point $c'\in \Lambda'\cap \Aline hy$. Assuming that $c'\in \Aline a{a'}$, we conclude that $c'\in \Aline hy\cap\Aline a{a'}=\{y\}$ and hence $y=c'\in \Lambda'$, which contradicts the choice of the points $x$ and $y$. This contradiction shows that $c'\in \Lambda'\setminus \Aline a{a'}=\Phi_L^{aa'}[\Lambda\setminus\Aline a{a'}]$.
On the other hand, the definition of $y=\bar\Phi_L^{aa'}(x)$ ensures that $\Phi_L^{aa'}[\Aline hx\setminus\Aline a{a'}]=\Aline hy\setminus \Aline a{a'}$. 
Then $\varnothing=\Phi_L^{aa'}[(\Lambda\cap \Aline hx)\setminus\Aline a{a'}]=(\Lambda'\cap\Aline hy)\setminus\Aline a{a'}=\{c'\}$, which is a desired contradiction showing that $\bar\Phi_L^{cc'}[\Lambda]=\Lambda'$, and hence $\bar\Phi_L^{aa'}$ is an automorphism of the projective plane $Y$.
\end{proof}

\begin{lemma}\label{l:lD-point-line} The point-line pair $(o,L)$ in $Y$ is Desarguesian.
\end{lemma}

\begin{proof} Given any points $b\in Y\setminus L$ and $b'\in \Aline ob\setminus\{o,b\}$, we need to find a central automorphism $\Phi:Y\to Y$ with center $o$ and axis $L$ such that $\Phi(b)=b'$. If $b,b'\in X$, then for $\Phi$ we can take the automorphism $\bar\Phi_L^{b,b'}$. If $b\in X$ and $b'\notin X$, then take any line $A$ in $Y$ such that $o\in A$ and $L\ne A\ne \Aline ob$. By the projectivity of $Y$, there exist a unique point $\alpha\in A\cap H=A\setminus X$. Since the liner $Y$ is $4$-long, there exists a point $x\in L\setminus(H\cup \{o\}\cup \Aline\alpha{b})$. By the projectivity of $Y$, there exist points $a\in A\cap \Aline xb$ and $a'\in A\cap\Aline x{b'}$. The choice of the point $x\notin H\cup\Aline \alpha{b}$ ensures that $a,a'\in X$. By Lemma~\ref{l:lD-lines5}, the bijective map $\Phi\defeq\bar\Phi_L^{aa'}$ is a central automorphism of $Y$ with center $o$, axis $L$ such that $\Phi(a)=a'$. Then $\Phi(b)\in \Phi[\Aline xa\cap \Aline ob]=\Phi[\Aline xa]\cap\Phi[\Aline ob]=\Aline x{a'}\cap\Aline o{b'}=\{b'\}$.   By analogy we can construct a central automorphism $\Phi\in \Aut_{o,L}(Y)$ with $\Phi(b)=b'$ in case $b\notin X$ and $b'\in X$. The case $b,b'\notin X$ is impossible because $o\in \Aline b{b'}\cap X$.
\end{proof}

Let $\mathcal L$ be the family of lines in the liner $Y$. By Lemma~\ref{l:lD-point-line}, for every line $L\subseteq Y$ and point $o\in L\cap X$, the point-line pair $(o,L)$ is Desarguesian. This implies that the Lenz figure of the projective plane $Y$ contains the set $\{(p,L)\in X\times \mathcal L:p\in L\}$ and hence $Y$ is of Lenz type $(\square)$. By Theorems~\ref{t:Lenz} and \ref{t:Skornyakov-San-Soucie}, the projective plane $Y$ is Moufang, and by Proposition~\ref{p:Moufang=>little-Desarguesian}, $Y$ is little-Desarguesian. 
\end{proof}


%


\section{Moufang projective spaces}

In this section we characterize Moufang projective spaces, extending the characterization of Moufang projective planes given in Theorem~\ref{t:Skornyakov-San-Soucie}.

\begin{theorem}\label{t:proj-Moufang<=>} For any $3$-long projective liner $Y$ of rank $\|Y\|\ge 3$, the following conditions are equivalent:
\begin{enumerate}
\item $Y$ is Moufang;
\item $Y$ is little-Desarguesian;
\item every subliner of $X$ is little-Desarguesian;
\item for some set $H\subseteq Y$ with $\overline H\ne Y$, the subliner $Y\setminus X$ is little-Desarguesian;
\item for every flat set $H$ in $Y$, the subliner $Y\setminus H$ is Moufang;
\item for every hyperplane $H$ in $Y$, the subliner $Y\setminus H$ is Moufang;
\item for every hyperplane $H$ in $Y$, the subliner $Y\setminus H$ is Thalesian;
\item for two distinct hyperplanes $H,H'$ in $Y$, the liners $Y\setminus H$ and $Y\setminus H'$ are Thalesian;
\item some ternar of $Y$ is linear, distributive, associative-plus and alternative-dot;
\item every ternar of $Y$ is linear, distributive, associative-plus and alternative-dot;
\item every ternar of $Y$ is Moufang;
\item every ternar of $Y$ is Moufang-dot;
\item every ternar of $Y$ is inversive-dot;
\end{enumerate}
If the projective space $Y$ has finite order, then \textup{(1)--(13)} are equivalent to the condition
\begin{enumerate}
\item[\textup{(14)}] some ternar of $Y$ is linear, distributive and Moufang;
\item[\textup{(15)}] some ternar of $Y$ is inear, distributive and associative;
\item[\textup{(16)}] $Y$ is Desarguesian.
\end{enumerate}
If $\|Y\|=3$, then the conditions \textup{(1)--(13)} are equivalent to the condition
\begin{enumerate}
\item[\textup{(17)}] $Y$ has Lenz type $(\square)$;
\item[\textup{(18)}] the dual projective plane to $Y$ is Moufang.
\end{enumerate}
\end{theorem}

\begin{proof} The implications $(1)\Ra(2)\Ra(3)$ follow from Propositions~\ref{p:Moufang=>little-Desarguesian} and \ref{p:lDes-subliner} and Theorem~\ref{t:lDes-projcompletion}, respectively. The implication $(3)\Ra(4)$ is trivial. 
\smallskip

$(4)\Ra(1)$ Assume that for some set $H\subseteq Y$ with $\overline H\ne Y$, the subliner $X\defeq Y\setminus H$ is little-Desarguesian. Then $Y$ is a projective completion and by Theorem~\ref{t:lDes-projcompletion}, the projective completion $Y$ of the little-Desarguesian liner $X$ is little-Desarguesian. By Theorem~\ref{t:projMoufang<=>lDes}, the projective little-Desarguesian liner $Y$ is Moufang. 
\smallskip

The implication $(1)\Ra(5)$ is proved in Proposition~\ref{p:Moufang-minus-flat}, the implication $(5)\Ra(6)$ is trivial, and $(6)\Ra(4)$ follows from Proposition~\ref{p:Moufang=>little-Desarguesian}. The equivalence $(7)\Leftrightarrow(1)$ has been proved in Theorem~\ref{t:Moufang<=>everywhere-Thalesian}. The implication $(7)\Ra(8)$ is trivial.
\smallskip

$(8)\Ra(1)$  Assume that the liner $Y$ contains two distinct hyperplanes $H_1,H_2$ such that the subliners $Y\setminus H_1$ and $Y\setminus H_2$ are Talesian. We have to prove that the projective space $Y$ is Moufang. If $\|Y\|\ge 4$, then the projective liner $Y$ is Desarguesian and Moufang, by Theorems~\ref{t:proaffine-Desarguesian} and Proposition~\ref{p:Desarg=>Moufang}.
If $\|Y\|=3$, then $Y$ is a projective plane and the hyperplanes $H_1,H_2,H$ are lines in $Y$. In this case we can apply Theorem~\ref{t:Skornyakov-San-Soucie} and conclude that the projective plane $Y$ is Moufang.
\smallskip

\smallskip

The implication $(1)\Ra(9)$ follows from Theorem~\ref{t:Skornyakov-San-Soucie}.
\smallskip

$(9)\Ra(10)$ Assume that some ternar of $Y$ is linear, distributive, associative-plus and alternative-dot. We should prove that every ternar $R$ of $Y$ is linear, distributive, associative-plus and alternative-dot. By definition, the ternar $R$ is isomorphic to the ternar $\Delta$ of some based projective plane $(P,uowe)$ in $Y$. If $|Y|_2=3$, then $\Delta=\{o,e\}$ coincides with the ternar of the two-element field and hence the ternar $R$ is linear, distributive and associative. So, assume that $|Y|_2\ge 4$. If $\|Y\|\ge 4$, then the projective space $Y$ is Desarguesian, by Theorem~\ref{t:proaffine-Desarguesian}. Then the plane $P$ in $Y$ is Desarguesian and so is the Playfair plane $P\setminus\Aline hv$, where $\{h\}=\Aline ou\cap\Aline we$ and $\{v\}=\Aline ow\cap\Aline ue$. By Theorem~\ref{t:Desarg=>tring=corps}, the ternar $\Delta=\Aline oe\setminus\Aline hv$ of the affine plane $P\setminus\Aline hv$ is linear, distributive and associative and so is its isomorphic copy $R$. So, assume that $\|Y\|=3$. In this case we can apply Theorem~\ref{t:Skornyakov-San-Soucie} and conclude that every ternar of the projective plane $Y$ is linear, distributive, associative-plus and associative-dot. In particular, so is the ternar $R$.
\smallskip

The implication $(10)\Ra(11)$ follow from Proposition~\ref{p:Bol<=>alternative} and Theorem~\ref{t:Bol<=>Moufang}; $(11)\Ra(12)$ is trivial, and $(12)\Ra(13)$ follows from Proposition~\ref{p:Moufang=>inversive}.
\smallskip

$(13)\Ra(7)$ Assume that every ternar of $Y$ is inversive-dot. If $\|Y\|\ge 4$ or $|Y|_2\le 3$, then the projective space $Y$ is Desarguesian by Theorems~\ref{t:proaffine-Desarguesian} and Proposition~\ref{p:Steiner+projective=>Desargues}. By Proposition~\ref{p:Desarg=>Moufang}, the Desarguesian liner $Y$ is Moufang. So, assume that $\|Y\|=3$ and $|Y|_2\ge 4$. We need to prove that for every line $H\subseteq Y$ the affine liner $X\defeq Y\setminus H$ is Thalesian. By Theorem~\ref{t:affine<=>hyperplane} and \ref{t:Playfair<=>}, the liner $X$ is a Playfair plane. By definition of ternars of the projective space $Y$, every ternar of the Playfair plane $X$ is a ternar of $Y$. Consequently, every ternar of the Playfair plane $X$ is inversive-dot. By Theorem~\ref{t:inversive-dot}, every ternar of the affine plane $X$ is linear, distributive, associative-plus, and alternative-dot. In particular, every ternar of $X$ is linear, right-ristributive and associative-plus. By Veblen--Weddenburn Theorems~\ref{t:VW-Thalesian<=>quasifield} and \ref{t:paraD<=>translation}, the Playfair plane $X$ is translation and Thalesian. By Theorem~\ref{t:Skornyakov-San-Soucie}, the projective plane $Y$ is Moufang.
\smallskip

Therefore, the conditions (1)--(13) are equivalent. It is clear that the conditions (10) and (11) imply (14). 
\smallskip

$(14)\Ra(15)$ Assume that some ternar $R$ of $Y$ is linear, distributive, and Moufang. Then there exists a plane $\Pi\subseteq Y$ and a projective base $uowe$ in $\Pi$ whose ternar is isomorphic to $R$ and hence  is linear, distributive, and Moufang. Let $h\in \Aline ou\cap\Aline we$ and $v\in \Aline ow\cap\Aline ue$ be the horizontal and vertical infinity points and $\Aline hv$ be the horizon line of the projective base $uowe$. By definition, the ternar of the based projective plane $(\Pi,uowe)$ coincides with the ternar $\Delta=\Aline oe\setminus\Aline hv$ of the based affine plane $(\Pi\setminus\Aline hv,uow)$. Then the ternar $\Delta$ is linear, distributive, and Moufang. 

If the projective liner $Y$ has finite order, then the affine plane $\Pi\setminus\Aline hv$ is finite. Applying Theorems~\ref{t:inversive-dot}, we conclude that the ternar $\Delta$ is linear, distributive, associative-plus and alternative-dot, and hence $\Delta$ is an alternative division ring. By Artin--Zorn Theorem~\ref{t:Artin-Zorn}, the finite alternative ring $\Delta$ is associative.
\smallskip

The implication $(15)\Ra(16)\Ra(1)$ follow from Theorem~\ref{t:Algebra-vs-Geometry-proj}(12) and Proposition~\ref{p:Desarg=>Moufang}, respectively. 
\smallskip

If $\|Y\|=3$, then the equivalence $(1)\Leftrightarrow(17)$ follows from Theorems~\ref{t:Skornyakov-San-Soucie}, \ref{t:Lenz} and the definition of the Lenz type $(\square)$. 

The equivalence $(17)\Leftrightarrow(18)$ follows from the equivalence $(1)\Leftrightarrow(17)$ and the definition of the Lenz figure of a projective plane.
\end{proof}

We do not know if the condition (14) is equivalent to other conditions of Theorem~\ref{t:proj-Moufang<=>} for arbitrary projective spaces. 

\begin{problem} Is a projective space Moufang if some its ternar is linear, distributive, and Moufang?
\end{problem}

Let us recall that a projective liner $Y$ is {\em somewhere $\mathcal P$} (resp. {\em everywhere $\mathcal P$}) if for some (resp. every) hyperplane $H\subseteq Y$, the affine subliner $Y\setminus H$ of $Y$ has property $\mathcal P$. 

Theorem~\ref{t:proj-Moufang<=>} implies the following characterization of Moufang projective spaces in terms of somewhere/everywhere $\mathcal P$ properties.

\begin{corollary}\label{c:Moufang<=>everywhereThales} For a projective space $X$, the following conditions are equivalent:
\begin{enumerate}
\item $X$ is Moufang;
\item $X$ is somewhere Moufang;
\item $X$ is everywhere Moufang;
\item $X$ is little-Desarguesian;
\item $X$ is somewhere little-Desarguesian;
\item $X$ is everywhere little-Desarguesian;
\item $X$ is everywhere Thalesian.
\end{enumerate}
\end{corollary}

\section{Completely regular Moufang liners}

In this section we prove that the uno-Desargues and Little Desargues Axioms are equivalent for completely regular liners.

\begin{theorem}\label{t:compreg-Moufang<=>} For a completely regular liner $X$, the following conditions are equivalent:
\begin{enumerate}
\item $X$ is Moufang;
\item $X$ is little-Desarguesian;
\item the spread completion of $X$ is Moufang;
\item some projective completion of $X$ is Moufang.
\end{enumerate}
\end{theorem}

\begin{proof} Let $X$ be a completely regular liner and $\overline X$ be its spread completion. By definition of the spread completion, the horizon $\partial X\defeq\overline X\setminus X$ of $X$ in $\overline X$ is flat.
\smallskip

The implications $(1)\Ra(2)\Ra(3)$ follow from Proposition~\ref{p:Moufang=>little-Desarguesian} and Theorem~\ref{t:lDes-projcompletion}, respectively. The implication $(3)\Ra(4)$ is trivial.
\smallskip

$(4)\Ra(1)$ Assume that some projective completion $Y$ of the completely regular liner $X$ is Moufang. By Theorem~\ref{t:spread=projective1}, the completely regular liner $X$ is para-Playfair, and by Theorem~\ref{t:flat-horizon}, the horizon $H\defeq Y\setminus X$ of the para-Playfair liner $X$ in $Y$ is flat. By Proposition~\ref{p:Moufang-minus-flat}, the subliner $X$ of the Moufang liner $Y$ is Moufang.
\end{proof}

Let us recall that a property $\mathcal P$ of liners is called \index{complete property}\index{property!complete}\defterm{complete} if a completely regular liner has property $\mathcal P$ if and only if its spread completion has property $\mathcal P$. Theorem~\ref{t:compreg-Moufang<=>}  implies the following important facts.

\begin{corollary}\label{c:Moufang=>complete} The properties of a liner to be Moufang (resp. little-Desarguesian) are complete.
\end{corollary}

Let us recall that a liner $X$ is \index{completely $\mathcal P$ liner}\index{liner!completely $\mathcal P$}\defterm{completely $\mathcal P$} for some property $\mathcal P$ of projective liners if $X$ is completely regular and its spread completion $\overline X$ has property $\mathcal P$. The following corollary answers Problem~\ref{prob:inner-completely-P} for the Moufang and little-Desarguesian properties.

\begin{corollary}\label{c:completelyM<=>} For a $3$-long liner $X$, the following conditions are equivalent:
\begin{enumerate}
\item $X$ is regular and Moufang;
\item $X$ is completely regular and little-Desarguesian;
\item $X$ is $3$-ranked and the spread completion of $X$ is Moufang and projective;
\item $X$ has a Moufang projective completion with flat horizon.
\item $X$ is completely little-Desarguesian;
\item $X$ is completely Moufang.
\end{enumerate}
\end{corollary}

\begin{proof} The implication $(1)\Ra(2)$ follows from Theorem~\ref{t:Moufang=>compreg} and Proposition~\ref{p:Moufang=>little-Desarguesian}; $(2)\Ra(3)$ follows from Theorem~\ref{t:lDes-projcompletion}; the implication $(3)\Ra(4)$ is trivial, and $(4)\Ra(1)$ follows from Theorems~\ref{t:flat-horizon} and \ref{p:Moufang-minus-flat}. The equivalence $(5)\Leftrightarrow(6)$ follows from Theorem~\ref{t:projMoufang<=>lDes}; $(3)\Leftrightarrow(6)$ from definition of completely Moufang spaces.
\end{proof}

\section{Para-Desarguesian liners}

In this section we introduce para-Desarguesian liners satisfying  a ``parallel'' version of the Desargues Axiom and establish some basic property of such liners.

\begin{definition}\label{d:para-Desarguesian} A liner $X$ is called \index{para-Desarguesian liner}\index{liner!para-Desarguesian}\defterm{para-Desarguesian} if satisfies the \defterm{para-Desargues Axiom}:
\begin{itemize}
\item[{\sf(pDA)}] for any plane $\Pi\subseteq X$, disjoint lines $A,B,C,D\subseteq \Pi$, and distinct points $a,a'\in A$, $b,b'\in B$, $c,c'\in C$, if
$\Line ab\cap\Line {a'}{b'}\cap D\ne\varnothing\ne\Line bc\cap\Line {b'}{c'}\cap D$, then $\Line ac\cap\Line{a'}{c'}\subseteq D$.
\end{itemize}

\begin{picture}(300,110)(-100,0)
\linethickness{=0.6pt}
\put(0,10){\color{teal}\line(1,0){160}}
\put(165,7){\color{teal}$D$}
\put(0,50){\color{teal}\line(1,0){160}}
\put(165,47){\color{teal}$C$}
\put(0,70){\color{teal}\line(1,0){160}}
\put(165,67){\color{teal}$B$}
\put(0,90){\color{teal}\line(1,0){160}}
\put(165,87){\color{teal}$A$}
\put(60,10){\color{red}\line(0,1){80}}
\put(140,10){\color{violet}\line(-1,2){40}}
\put(60,90){\color{violet}\line(1,-1){80}}
\put(60,10){\color{red}\line(1,2){40}}
\put(20,10){\color{blue}\line(1,1){60}}
\put(20,10){\color{blue}\line(3,2){90}}

\put(20,10){\color{blue}\circle*{3}}
\put(60,10){\color{red}\circle*{3}}
\put(140,10){\color{violet}\circle*{3}}
\put(60,50){\circle*{3}}
\put(53,52){$c$}
\put(80,50){\circle*{3}}
\put(75,52){$c'$}
\put(80,70){\circle*{3}}
\put(81,72){$b$}
\put(110,70){\circle*{3}}
\put(111,72){$b'$}
\put(60,90){\circle*{3}}
\put(58,93){$a$}
\put(100,90){\circle*{3}}
\put(98,93){$a'$}
\end{picture}

\end{definition}

Definition~\ref{d:para-Desarguesian} implies the following simple (but useful) characterization of para-Desarguesian liners.

\begin{proposition}\label{p:para-Des<=>plane-pD} A liner $X$ is para-Desarguesian if and only if every plane in $X$ is para-Desarguesian.
\end{proposition}

Many examples of para-Desarguesian liners can be constructed using the following proposition.

\begin{proposition}\label{p:subliner-Moufang} Any proaffine subliner of a projective Moufang liner is para-Desarguesian.
\end{proposition}

\begin{proof} Let $X$ be a proaffine subliner of a projective Moufang liner $Y$. By Proposition~\ref{p:Moufang=>little-Desarguesian}, the Moufang projective liner $Y$ is little-Desarguesian. By Theorem~\ref{t:proaffine<=>proflat}, the horizon $H\defeq Y\setminus X$ of the proaffine liner $X$ in $Y$ is proflat in $Y$. For every line $L\subseteq X$, let $\overline L$ be its flat hull in the projective liner $Y$ and $\partial L\defeq\overline L\setminus L$ be its horizon. Since $H$ is proflat, $|\partial L|\le 1$. 

 To prove that $X$ is para-Desarguesian, take any plane $\Pi\subseteq X$, disjoint lines $A,B,C,D\subset \Pi$ and distinct points  $a,a'\in A$, $b,b'\in B$, $c,c'\in C$ such that $\Line ab\cap\Line {a'}{b'}\cap D\ne\varnothing\ne\Line bc\cap\Line {b'}{c'}\cap D$. We have to prove that $\Line ac\cap\Line{a'}{c'}\subseteq D$. By our assumption, there exist points $x\in \Line ab\cap\Line {a'}{b'}\cap D$ and $y\in \Line bc\cap\Line {b'}{c'}\cap D$. 
 
If $b\in \Aline ac$, then $\{x\}=\Aline ab\cap D=\Aline bc\cap D=\{y\}$ and $\{a',c'\}\subseteq \Aline x{b'}=\Aline y{b'}$. Then $\Aline {a}{c}\cap\Aline {a'}{c'}=\Aline ab\cap\Aline {a'}{b'}=\{x\}\subseteq D$. By analogy we can show that $b'\in\Aline {a'}{c'}$ implies $\Aline ac\cap\Aline{a'}{c'}=\{x\}=\{y\}\subseteq D$. 

So, assume that $abc$ and $a'b'c$ are triangles. Consider the plane $\overline{\{b,b',x\}}$ in the projective liner $Y$ and observe that $\{a,a'\}\subseteq\Aline xb\cup\Aline x{b'}\subseteq \overline{\{b,b',x\}}$. By the projectivity of $Y$, the lines $\overline A=\overline{\{a,a'\}}$ and $\overline B=\overline{\{b,b'\}}$ in the projective plane $\overline{\{b,b',x\}}$ have a common point $o$, which is the unique horizon point of the lines $A$ and $B$. Considering the plane $\overline{\{b,b',y\}}$ we can prove that the point $o$ also coincides with the horizon point of the line $C$.

Let $M$ be the maximal $3$-long flat in $Y$, containing the point $o$. By Lemma~\ref{l:ox=2}, for every point $p\in Y\setminus M$, the line $\overline{\{o,p\}}$ coincides with the doubleton $\{o,p\}$. Taking into account that the lines $\overline A,\overline B,\overline C$ are $3$-long, we conclude that $\overline A\cup\overline B\cup\overline C\subseteq M$ and also $\{x,y\}\subseteq \Aline ab\cup\Aline bc\subseteq M$ and $\overline D\subseteq M$.

If $\|M\|\ge 4$, then the $3$-long projective space $M$ is Desarguesian (by Theorem~\ref{t:proaffine-Desarguesian}) and hence 
$$\Aline ac\cap \Aline {a'}{c'}=\overline{\{a,c\}}\cap\overline{\{a',c'\}}\cap X\subseteq \overline{\{x,y\}}\cap X=D.$$ So, assume that $\|M\|=3$. Then $M$ is the projective plane containing the lines $\overline A$ and $\overline D$ which have a unique common point, which coincides with the horizon point $o$ of the disjoint lines $A$ and $D$. Now we see that $\overline A,\overline B,\overline C,\overline D$ are concurrent lines in the plane $M$. 
By Proposition~\ref{p:Moufang=>little-Desarguesian}, the Moufang liner $Y$ is little-Desarguesian, which implies $\varnothing\ne \overline{\{a,c\}}\cap\overline D\subseteq\overline{\{a',c'\}}$ and 
$$\Aline ac\cap\Aline {a'}{c'}=\overline{\{a,c\}}\cap\overline{\{a',c'\}}\cap X\subseteq \overline D\cap X=D,$$witnessing that the liner $X$ is para-Desarguesian.
\end{proof}

\begin{exercise} Find an example of a para-Desarguesian liner, which is not regular, not para-Playfair, and not Moufang. 
\smallskip

{\em Hint:} Consider the complement to a line with removed point in a $4$-long Moufang projective liner $Y$ of rank $\|Y\|\ge4$ and apply Theorems~\ref{t:regular<=>flat4}, \ref{t:flat-horizon} and Proposition~\ref{p:nD=>pP}.
\end{exercise}

\begin{proposition}\label{p:Moufang=>paraD} Every Moufang Proclus liner is para-Desarguesian.
\end{proposition}

\begin{proof} Let $X$ be a Moufang Proclus liner. To show that $X$ is para-Desarguesian, it suffices to check that every plane $\Pi$ in $X$ is para-Desarguesian. Since $X$ is Moufang, so is the plane $\Pi$ in $X$. If $\Pi$ is projective, then $\Pi$ is para-Desarguesian, by Proposition~\ref{p:subliner-Moufang}. So, assume that $\Pi$ is not projective. If $\Pi$ is not $3$-long, then by Theorem~\ref{t:Proclus-not-3long}, $|\Pi|\le 6$ and $\Pi$ is para-Desarguesian because it contains no four disjoint lines. So, assume that $\Pi$ is $3$-long. By Theorem~\ref{t:Proclus<=>}, the $3$-long Proclus plane $\Pi$ is $3$-regular and hence regular. By Theorem~\ref{t:Moufang=>compreg}, the regular Moufang liner $\Pi$ is completely regular. By Theorem~\ref{t:compreg-Moufang<=>}, the spread completion $\overline\Pi$ of the Moufang plane $\Pi$ is Moufang. Since the horizon $\overline \Pi\setminus \Pi$ of $\Pi$ in $\overline \Pi$ is flat, the subliner $\Pi$ of the Moufang liner $\overline \Pi$ is para-Desarguesian, by Proposition~\ref{p:subliner-Moufang}. Therefore, every plane in $X$ is para-Desarguesian and hence the liner $X$ is para-Desarguesian, according to Proposition~\ref{p:para-Des<=>plane-pD}.
\end{proof}

\begin{corollary}\label{c:Steiner=>paraD} Every Steiner $3$-regular liner is para-Desarguesian.
\end{corollary}

\begin{proof} Let $X$ be a Steiner $3$-regular liner. Since Steiner liners are proaffine, the Steiner proaffine $3$-regular liner is Proclus, by Theorem~\ref{t:Proclus<=>}. By  Propositions~\ref{p:Steiner+projective=>Desargues} and \ref{p:Desarg=>Moufang}, the Steiner $3$-regular liner $X$ is Desarguesian and Moufang. By Proposition~\ref{p:Moufang=>paraD}, the Moufang Proclus liner $X$ is para-Desarguesian.
\end{proof}

The $3$-regularity in Corollary~\ref{c:Steiner=>paraD} is essential as shows the following example, found by Ivan Hetman.

\begin{example}[Ivan Hetman, 2025]\label{ex:Z15pM} The  cyclic group $X:=C_{15}=\{0,1,\dots,15\}$ endowed with the family of lines $$\mathcal L=\big\{x+L:x\in C_{15},\;L\in\{\{0,1,4\},\{0,2,9\},\{0,5,10\}\big\}$$ is an affine Steiner ranked plane, which is Desarguesian and Moufang, but not Thalesian and not para-Desarguesian.
\end{example}

\begin{proof} It is well-known that the group $C_{15}$ is isomorphic to the product $C_5\times C_3$, which can be visualized as a $5\times 3$ rectangle. Elements of $C_5\times C_3$ can be labeled by ordered pairs $xy$ of numbers $x\in 5=\{0,1,2,3,4\}$ and $y\in 3=\{0,1,2\}$. The family of lines $\mathcal L$ coincides with all possible shifts of the disjoint lines 
$$L_1\defeq\{02,11,22\},\quad L_2\defeq\{01,20,41\},\quad L_3\defeq\{30,31,32\}.$$

\begin{picture}(120,90)(-150,-15)
\multiput(0,0)(0,30){3}{\multiput(0,0)(30,0){5}{\color{gray}\circle*{3}}}
\put(0,60){\color{cyan}\line(1,-1){30}}
\put(30,30){\color{cyan}\line(1,1){30}}
\put(0,30){\color{blue}\line(2,-1){60}}
\put(60,0){\color{blue}\line(2,1){60}}
\put(90,0){\color{magenta}\line(0,1){60}}
\put(87,64){\color{magenta}$L_3$}
\put(125,26){\color{blue}$L_2$}
\put(-15,56){\color{cyan}$L_1$}
\put(0,60){\color{cyan}\circle*{3}}
\put(30,30){\color{cyan}\circle*{3}}
\put(60,60){\color{cyan}\circle*{3}}
\put(0,30){\color{blue}\circle*{3}}
\put(60,0){\color{blue}\circle*{3}}
\put(120,30){\color{blue}\circle*{3}}
\put(90,0){\color{magenta}\circle*{3}}
\put(90,30){\color{magenta}\circle*{3}}
\put(90,60){\color{magenta}\circle*{3}}
\end{picture}

 By Example~\ref{ex:Z15D}, the Steiner liner $X$ is affine, ranked, Desarguesian and not Thalesian. To see that the liner $X$ is not para-Desarguesian, consider the lines
$$A\defeq\{01,10,21\},\quad B\defeq\{22,31,42\},\quad C\defeq\{00,11,41\},\quad D\defeq\{20,30,40\}$$and the points 
$a\defeq 10$, $a'\defeq 21$, \ $b\defeq 42$, $b'\defeq 22$, \ $c\defeq 41$, $c'\defeq 00.$
Observe that\break $20\in\Aline ab\cap\Aline {a'}{b'}\cap D$, $40\in\Aline bc\cap\Aline{b'}{c'}\cap D$, but $\Aline ac\cap\Aline {a'}{c'}=\{31\}\not\subseteq D$, 
witnessing that the liner $X$ is not para-Desarguesian.

\begin{picture}(120,95)(-140,-15)
\linethickness{0.7pt}
\multiput(0,0)(0,30){3}{\multiput(0,0)(30,0){5}{\color{gray}\circle*{2}}}
\put(30,0){\color{teal}\line(1,1){30}}
\put(0,30){\color{teal}\line(1,-1){30}}
\put(60,60){\color{teal}\line(1,-1){30}}
\put(90,30){\color{teal}\line(1,1){30}}
\put(0,0){\color{teal}\line(1,1){30}}
{\color{teal}
\qbezier(30,30)(75,52)(120,30)
}
{\linethickness{1pt}
\put(60,0){\color{teal}\line(1,0){60}}
\put(30,0){\color{blue}\line(1,0){30}}
{\color{blue}
\qbezier(60,0)(105,15)(120,60)
}\put(60,0){\color{cyan}\line(0,1){60}}
\put(120,0){\color{brown}\line(0,1){60}}
{\color{olive}
\qbezier(0,0)(10,50)(60,60)
\qbezier(60,60)(110,50)(120,0)
}
\put(0,0){\color{magenta}\line(2,1){60}}
\put(60,30){\color{magenta}\line(1,0){30}}
\put(30,0){\color{red}\line(2,1){60}}
\put(90,30){\color{red}\line(1,0){30}}
}
\put(60,0){\color{teal}\circle*{3}}
\put(120,0){\color{teal}\circle*{3}}

\put(30,0){\circle*{3}}
\put(27,-8){$a$}
\put(60,30){\circle*{3}}
\put(62,32){$a'$}

\put(120,60){\circle*{3}}
\put(118,63){$b$}
\put(60,60){\circle*{3}}
\put(58,63){$b'$}

\put(120,30){\circle*{3}}
\put(123,28){$c$}
\put(0,0){\circle*{3}}
\put(-2,-8){$c'$}
\put(90,30){\color{red}\circle*{4}}

\put(-10,32){\color{teal}$A$}
\put(100,50){\color{teal}$B$}
\put(42,40){\color{teal}$C$}
\put(88,-10){\color{teal}$D$}
\end{picture}
\end{proof}

\section{Para-Desarguesian para-Playfair liners}

In this section we study the interplay between the para-Desarguesian and para-Playfair liners. In particular, in Theorem~\ref{t:para-Desarg<=>shear} we prove that a para-Playfair plane is para-Desarguesian if and only if it is shear; in Theorem~\ref{t:para-Desarg-compreg<=>} we prove that a para-Desarguesian liner is completely regular if and only if it is regular and para-Playfair.

First, we show that for para-Playfair liners, the para-Desargues Axiom holds in a stronger form.

\begin{proposition}\label{p:paraMouf=>z} Let $X$ be a  para-Desarguesian para-Playfair liner. For any plane $\Pi\subseteq X$, disjoint  lines $A,B,C,D\subseteq \Pi$ and distinct points $a,a'\in A$, $b,b'\in B$, $c,c'\in C$, if $\Line ab\cap\Line {a'}{b'}\cap D\ne\varnothing\ne\Line bc\cap\Line {b'}{c'}\cap D$ then $\varnothing\ne\Line ac\cap\Line{a'}{c'}\subseteq D$.
\end{proposition}

\begin{proof} By Proposition~\ref{p:Playfair=>para-Playfair=>Proclus}, the para-Playfair liner is Proclus; by Theorem~\ref{t:Proclus<=>}, the Proclus liner $X$ is $3$-proregular, and by Proposition~\ref{p:k-regular<=>2ex}, the $3$-proregular liner $X$ is $3$-ranked. Then any disjoint coplanar lines in $X$ are parallel, by Corollary~\ref{c:parallel-lines<=>}. Since $X$ is para-Desarguesian, $\Line ac\cap\Line {a'}{c'}\subseteq D$. It remains to prove that $\Line ac\cap\Line {a'}{c'}\ne\varnothing$. To derive a contradiction, assume that  $\Line ac\cap\Line {a'}{c'}=\varnothing$. 
By our assumption, $\Line ab\cap\Line {a'}{b'}\cap D=\{x\}$ and $\Line bc\cap\Line {b'}{c'}\cap D=\{y\}$ for some points $x,y$. 
By Proclus Postulate~\ref{p:Proclus-Postulate}, $\Line ac\cap D=\{z\}$ and $\Line {a'}{c'}\cap D=\{z'\}$ for some points $z,z'$.

\begin{picture}(300,125)(-100,-15)
\linethickness{=0.6pt}
\put(-40,0){\color{teal}\line(1,0){220}}
\put(185,-3){\color{teal}$D$}
\put(-40,50){\color{teal}\line(1,0){220}}
\put(185,47){\color{teal}$C$}
\put(-40,70){\color{teal}\line(1,0){220}}
\put(185,67){\color{teal}$B$}
\put(-40,90){\color{teal}\line(1,0){220}}
\put(185,87){\color{teal}$A$}
\put(10,0){\color{blue}\line(0,1){70}}
\put(60,0){\color{red}\line(0,1){90}}
\put(150,0){\color{violet}\line(0,1){90}}
\put(152,25){\color{red}$L$}
\put(60,90){\color{violet}\line(1,-1){90}}
\put(60,0){\color{red}\line(1,1){90}}
\put(10,0){\color{blue}\line(1,1){70}}
\put(10,0){\color{blue}\line(2,1){140}}
\put(150,0){\color{violet}\line(-2,1){180}}
\put(10,50){\color{red}\line(-1,1){40}}

\put(10,0){\color{blue}\circle*{3}}
\put(7,-9){$y$}
\put(60,0){\color{red}\circle*{3}}
\put(57,-9){$z$}
\put(150,0){\color{violet}\circle*{3}}
\put(147,-9){$x$}
\put(10,50){\circle*{3}}
\put(11,52){$c'$}
\put(60,50){\circle*{3}}
\put(53,52){$c$}
\put(110,50){\circle*{3}}
\put(106,52){$c''$}
\put(10,70){\circle*{3}}
\put(8,72){$b'$}
\put(80,70){\circle*{3}}
\put(81,72){$b$}
\put(150,70){\circle*{3}}
\put(152,72){$b''$}
\put(-30,90){\circle*{3}}
\put(-33,93){$a'$}
\put(60,90){\circle*{3}}
\put(58,93){$a$}
\put(150,90){\circle*{3}}
\put(147,93){$a''$}
\end{picture}

Since $\Aline ac$ and $\Aline{a'}{c'}$ are disjoint lines in the para-Playfair liner $X$,  there exists a line $L\subseteq \Pi$ such that $x\in L$ and $L\parallel \Aline ac$. By the Proclus Axiom, there exist points $a''\in A\cap L$, $b''\in B\cap L$, and $c''\in C\cap\Line y{b''}$. Then $x\in \Line ab\cap\Line {a'}{b'}\cap\Line {a''}{b''}\cap D$ and $y\in \Line bc\cap\Line{b'}{c'}\cap\Line {b''}{c''}\cap D$. Since $X$ is para-Desarguesian, $(\Line {a}{c}\cap\Line {a''}{c''})\cup(\Line {a''}{c''}\cap\Line {a'}{c'})\subseteq D$. 
By Proclus Postulate~\ref{p:Proclus-Postulate}, there exists a unique point $z''\in \Line {a''}{c''}\cap D$. Assuming that the line $\Line{a''}{c''}$ is not parallel to the parallel lines $\Line ac$ and $\Line{a'}{c'}$, we conclude that 
$$\varnothing\ne \Line ac\cap\Line{a''}{c''}=\Line ac\cap\Line{a''}{c''}\cap D=(\Line ac\cap D)\cap(\Line{a''}{c''}\cap D)=\{z\}\cap\{z''\}$$ and hence $z=z''$. By analogy we can prove that $z''=z'$ and hence $z''=z=z'\in\Line ac\cap\Line {a'}{c'}=\varnothing$, which is a contradiction showing that $\Line{a''}{c''}\parallel \Line ac\parallel L=\Line {a''}{b''}$ and hence $\Line {a''}{c''}=\Line {a''}{b''}=\Line {b''}{c''}$ and $y\in \Line{b''}{c''}\cap D=\Line{a''}{c''}\cap D=\{x\}$.

Then $a,c\in \Line xb=\Line yb$, $a',c'\in \Line x{b'}=\Line y{b'}$ and hence $\Line ac\cap\Line {a'}{c'}=\Line xb\cap\Line x{b'}=\{x\}$, which contradicts the assumption $\Line ac\cap\Line{a'}{c'}=\varnothing$.
\end{proof}

A plane $X$ is \defterm{shear} if every line $L$ in $X$ is \index{shear line}\index{line!shear}\defterm{shear}, which means that for every distinct points $a,a'\in X$ with $\Aline xy\cap L=\varnothing$, there exists a liner automorphism $\Phi:X\to X$ such that $\Phi(a)=a'$ and the fixed point set $\Fix(\Phi)\defeq\{z\in X:\Phi(z)=z\}$ of $\Phi$ equals $L$. 

\begin{proposition}\label{p:shear=>para-Desarguesian} Every shear Proclus plane is para-Desarguesian.
\end{proposition}

\begin{proof} By Theorem~\ref{t:Proclus<=>}, the Proclus liner $X$ is $3$-proregular and by Proposition~\ref{p:k-regular<=>2ex}, the $3$-proregular liner $X$ is $3$-ranked. To prove that the shear liner $X$ is para-Desarguesian, take any plane $\Pi$ , disjoint lines $A,B,C,D\subset \Pi$ and distinct points $a,a'\in A$,  $b,b'\in B$, $c,c'\in C$ such that $\Aline ab\cap\Aline{a'}{b'}\cap D\ne\varnothing\ne\Aline bc\cap\Aline{b'}{c'}\cap D$. We have to prove that $\Aline ac\cap\Aline {a'}{c'}\subseteq D$.
Since the line $D$ is shear, there exists a liner automorphism $\Phi:X\to X$ such that $\Phi(b)=b'$ and $\Fix(\Phi)=D$. 

\begin{claim} For any line $L\subseteq X$ with $L\cap D=\varnothing$, its image $\Phi[L]$ coincides with $L$.
\end{claim}

\begin{proof} It follows from $L\cap D=\varnothing$ and $\Phi[D]=D$ that $$\Phi[L]\cap D=\Phi[L]\cap\Phi[D]=\Phi[L\cap D]=\Phi[\varnothing]=\varnothing$$and hence $L$ and $\Phi[L]$ are two lines, disjoint with the line $D$. Assuming that $L\cap\Phi[L]\ne\varnothing$, we conclude that $L=\Phi[L]$, by the Proclus Axiom (holding for the Proclus plane $X$). So, assume that $L\cap \Phi[L]=\varnothing$. Choose a point $x\in L$ and consider the point $y\defeq\Phi(x)$. Since $L\cap\Phi(L)=\varnothing$ and $\Phi(L)\cap D=\varnothing=L\cap D$, the Proclus Postulate~\ref{p:Proclus-Postulate} ensures that the line $\Aline xy$ has a common point $z$ with the line $D$. Then $\Phi[\Aline xz]=\Aline{\Phi(x)}{\Phi(z)}=\Aline yz=\Aline xz$. Since $A\cap D=\varnothing$, Proclus Postulate~\ref{p:Proclus-Postulate} ensures that $\Aline xz\cap A$ contains some point $\alpha$. It follows from $a'\in A\cap\Phi[A]$ and $A\cap D=\varnothing=\Phi[A]\cap D$ that $A=\Phi[A]$ and hence $\Phi(\alpha)\in\Phi[A\cap \Aline xz]=\Phi[A]\cap\Phi[\Aline xz]=A\cap\Aline xz=\{\alpha\}$, which contradicts $\Fix(\Phi)=D$. This contradiction shows that $\Phi[L]=L$.
\end{proof}

By the assumption, there exist points $x\in\Aline ab\cap\Aline{a'}{b'}\cap D$ and $y\in\Aline bc\cap\Aline{b'}{c'}\cap D$. By Proclus Postulate~\ref{p:Proclus-Postulate}, there exists a point $z\in \Aline ac\cap D$. It follows from $x,y\in\Fix(\Phi)=D$ and $\Phi(b)=b'$ that $\Phi(a)\in\Phi[\Aline bx\cap A]=\Phi[\Aline bx]\cap \Phi[A]=\Aline{b'}x\cap A=\{a'\}$ and $\Phi(c)\in\Phi[\Aline by\cap C]=\Phi[\Aline by]\cap \Phi[C]=\Aline{b'}y\cap C=\{c'\}$, which implies $z=\Phi(z)\in \Phi[\Aline ac]=\Aline{a'}{c'}$ and hence $\Aline ac\cap \Aline{a'}{c'}=\{z\}\subseteq D$.
\end{proof}

\begin{proposition}\label{p:Pappian=>shear} Every Pappian para-Playfair plane is shear.
\end{proposition}

\begin{proof} Let $X$ be a Pappian para-Playfair plane.  By Proposition~\ref{p:Playfair=>para-Playfair=>Proclus}, Theorem~\ref{t:Proclus<=>}, and Proposition~\ref{p:k-regular<=>2ex}, the para-Playfair plane is Proclus, proaffine, $3$-proregular and $3$-ranked. Since $X$ is a plane, it suffices to prove that every line $D$ in $X$ is shear. Given any distinct points $a,a'\in X$ with $\Aline a{a'}\cap D=\varnothing$, we have to find an automorphism $\Phi:X\to X$ such that $\Phi(a)=a'$ and $\Fix(\Phi)=D$. The plane $X$ is not projective because it contains two disjoint lines $D$ and $\Aline a{a'}$. 

\begin{claim}\label{cl:Pap-has-Pcompletion} The liner $X$ has a Pappian projective completion $Y$ with flat horizon $H\defeq Y\setminus X$.
\end{claim}

\begin{proof} If $X$ is not $3$-long, then by Theorems~\ref{t:Proclus-not-3long} and \ref{t:flat-horizon}, $X=Y\setminus H$ for some Steiner projective plane $Y$ and some singleton $H$ in $Y$. Since $Y$ contains no $4$-long lines, the Steiner projective plane $Y$ is Pappian. 

If $X$ is $3$-long, then the $3$-proregular plane $X$ is $3$-regular and hence regular. By Theorem~\ref{t:Pappian=>compreg}, the $3$-long Pappian proaffine regular liner $X$ is completely regular, and by Theorem~\ref{t:Pappian-completion}, the spread completion of $X$ is a Pappian projective completion of $X$ with flat horizon. 
\end{proof}

By Claim~\ref{cl:Pap-has-Pcompletion}, the plane $X$ has a Pappian projective completion $Y$ with flat horizon $H\defeq Y\setminus X$. Since $X$ is a plane, its projective completion $Y$ is a projective plane. Then the lines $\overline D$ and $\overline{\{a,a'\}}$ in $Y$ have a comon point $h\in H$. Since $\|H\|<\|Y\|=3$, the horizon $H$ is either a singleton or a line in $Y$. Fix any distinct points $o,d\in D$. If $H$ is a line in $Y$, then by the projectivity of the plane $Y$, there exist points $v\in \overline{\{o,a\}}\cap H$ and $v'\in \overline{\{v,h\}}\subseteq H$. If $H=\{h\}$, then choose any point $v\in \overline{\{o,a\}}\setminus\{o,a\}$ and by the projectivity of $Y$, find a point $v'\in\overline{\{o,a'\}}\cap \overline{\{v,h\}}$. By the projectivity of the plane $Y$, there exists a point $d\in \overline{\{v,a'\}}\cap\overline D$.

\begin{picture}(100,125)(-180,-15)
\put(0,0){\line(1,0){90}}
\put(0,0){\line(0,1){90}}
\put(0,90){\line(1,-1){90}}
\put(0,0){\line(1,1){45}}
\put(0,45){\line(2,-1){90}}
\put(0,0){\line(1,1){45}}
\put(0,90){\line(1,-2){45}}

\put(0,0){\circle*{3}}
\put(-6,-8){$o$}
\put(45,0){\circle*{3}}
\put(43,-10){$d$}

\put(90,0){\circle*{3}}
\put(93,-2){$h$}
\put(0,90){\circle*{3}}
\put(-3,93){$v$}
\put(0,45){\circle*{3}}
\put(-8,43){$a$}
\put(45,45){\circle*{3}}
\put(45,47){$v'$}
\put(30,30){\circle*{3}}
\put(28.5,34){$a'$}
\end{picture}

 Then $r\defeq\{oo,dh,av\}$ and $r'\defeq\{oo,dh,a'v'\}$ are two projective repers in the projective plane $Y$. Consider the injective function $\varphi\defeq\{oo,dd,hh,aa',vv'\}$ mapping the reper $r$ onto the reper $r'$. By Corollary~\ref{c:Papp-iso-reper}, the function $\varphi$ extends to a unique projective automorphism $\overline \Phi:Y\to Y$. For every $x\in\{o,d,h\}$ we have $\overline \Phi(x)=\varphi(x)=x$, which implies $\overline \Phi[\overline D]=\overline D$. By Corollary~\ref{c:Pappus-Fix3}, the projective automorphism $\overline\Phi$ is identity on the line $\overline D$. If $H=\{h\}$, then $\overline\Phi(h)=\varphi(h)=h$ implies $\overline\Phi[H]=H$. If $H$ is a line, then $H=\Aline vh$ and $\overline\Phi[\Aline vh]=\Aline{v'}h=H$. In both cases, we obtain $\overline \Phi[H]=H$ and hence $\overline \Phi[X]=\overline\Phi[Y\setminus H]=Y\setminus\overline\Phi[H]=Y\setminus H=X$. Then $\Phi\defeq\overline\Phi{\restriction}_X$ is an automorphism of the liner $X$ such that $\Phi(a)=a'$ and $D\subseteq \Fix(\Phi)$.

 Assuming that $D\ne\Fix(\Phi)$, we can find a fixed point $c\in \Fix(\Phi)\setminus D\subseteq X$. By the projectivity of the plane $Y$, there exists a point $\hbar\in \overline D\cap\overline{\{a,c\}}$. Then $a\in \overline{\{\hbar,c\}}$ and hence $a'=\Phi(a)\in \overline \Phi[\overline{\{\hbar,c\}}]=\overline{\{\hbar,c\}}$, $\overline{\{\hbar,c\}}=\overline{\{a,a'\}}$, and $\hbar\in \overline D\cap\overline{\{a,c\}}=\overline D\cap\overline{\{a,a'\}}=\{h\}$. By the projectivity of the plane $Y$, there exists a point $\gamma\in \overline D\cap\overline{\{v,c\}}$. Then $v'=\overline\Phi(v)=\overline \Phi[\overline{\{\gamma,c\}}]=\overline{\{\gamma,c\}}=\overline{\{v,c\}}$. 
 It follows from $v\notin \overline{\{a,a'\}}$ and $c\ne h$ that $\gamma\ne h$ and hence $v'\in \overline{\{\gamma,v\}}\cap\overline{\{h,v\}}=\{v\}$, which contradicts $v'\in\overline{\{o,a'\}}$ and $a\ne a'$. This contradiction shows that $\Fix(\Phi)=D$, witnessing that the line $D$ is shear.
\end{proof}

The following important theorem was first proved by Oleksi\u\i\ Ilchuk in his Bachelor Thesis \cite{Ilchuk}.
 
\begin{theorem}\label{t:para-Desarg<=>shear} A para-Playfair plane is para-Desarguesian if and only if it is shear.
\end{theorem}

\begin{proof} The ``if'' part follows from Propositions~\ref{p:shear=>para-Desarguesian} and \ref{p:Playfair=>para-Playfair=>Proclus}. To prove the ``only if'' part, assume that $X$ is a para-Desarguesian para-Playfair plane.
 By Proposition~\ref{p:Playfair=>para-Playfair=>Proclus}, the para-Playfair liner $X$ is Proclus; by Theorem~\ref{t:Proclus<=>}, the Proclus liner $X$ is proaffine and $3$-proregular, and by Proposition~\ref{p:k-regular<=>2ex}, the $3$-proregular liner $X$ is $3$-ranked. By Corollary~\ref{c:parallel-lines<=>}, any disjoint lines in the $3$-ranked plane $X$ are parallel.

We have to prove that the plane $X$ is shear. This is vacuously true if the liner $X$ is projective (because projective planes contains no disjoint lines). So, assume that the liner $X$ is not projective. If $X$ is not $3$-long, then by Theorem~\ref{t:Proclus-not-3long}, every line in $X$ contains at most three points, which implies that $X$ is Pappian.
By Proposition~\ref{p:Pappian=>shear}, the Pappian para-Playfair liner is shear. So, assume that $X$ is $3$-long. In this case the $3$-proregular plane $X$ is $3$-regular and hence regular. 

If $X$ is not $4$-long, then every line in $X$ contains at most $4$-points, by Theorem~\ref{t:card-lines-in-Proclus}. By the regularity, the plane $X$ is finite, and by Theorem~\ref{t:procompletion-finite}, the finite proaffine regular plane $X$ has a projective completion $Y$. Since lines in $X$ contain at most $4$ points, $|Y|_2\le 5$. By Corollary~\ref{c:p5-Pappian}, the projective plane $Y$ with $|Y|_2\le 5$ is Pappian. By Theorem~\ref{t:flat-horizon}, the horizon $Y\setminus X$ of the para-Playfair liner $X$ in $Y$ is flat, and by Proposition~\ref{p:Pappian-minus-flat}, the subliner $X$ of the Pappian liner $Y$ is Pappian. 
By Proposition~\ref{p:Pappian=>shear}, the Pappian para-Playfar liner $X$ is shear.

So, assume that the liner $X$ is $4$-long and fix any line $D$ in $X$. We shall prove that the line $D$ is shear in $X$.

\begin{lemma}\label{l:part-automorphism} For every line $A\subseteq X\setminus D$ and distinct points $a,a'\in A$, there exists an injective function $\varphi$ such that $\dom[\varphi]=(X\setminus A)\cup\{a\}$, $\rng[\varphi]=(X\setminus A)\cup\{a'\}$, $aa'\in\varphi$, and for every pairs $xx',yy'\in\varphi$ the following conditions are satisfied:
\begin{enumerate}
\item if $x\in D$, then $x'=x$;
\item if $x\notin D$, then $\Aline x{x'}\cap D=\varnothing$;
\item if $y\in \Aline ox$ for some point $o\in D$, then $y'\in \Aline o{x'}$.
\end{enumerate}
\end{lemma}

\begin{proof} Since $X$ is a para-Playfair plane and the lines $D,A$ are disjoint, for every point $x\in X$ there exists a unique line $L_x$ such that $x\in L_x$ and $L_x\parallel D$. Since $X$ is Proclus, for every point $x\in X\setminus A$, there exists a unique point $d_x\in D\cap \Aline ax$ and a unique point $x'\in \Aline {d_x}{a'}\cap L_x$. We claim that the function 
$$\varphi\defeq\{aa'\}\cup\{xx':x\in X\setminus A\}$$
has the required properties.

It is easy to see that $\varphi$ is an injective function with $\dom[\varphi]=(X\setminus A)\cup\{a\}$ and $\rng[\varphi]=(X\setminus A)\cup\{a'\}$. If $x\in D$, then $d_x=x$ and $x'\in \Aline {d_x}{a'}\cap L_x=\Aline {d_x}{a'}\cap D=\{d_x\}=\{x\}$. If $x\notin D$, then $\Aline x{x'}\cap D\subseteq L_x\cap D=\varnothing$.

Finally, assume that $xx',yy'\in\varphi$ are two pairs such that $y\in\Aline ox$ for some point $o\in D$. We have to prove that $y'\in \Aline o{x'}$. If $x=y$, then $y'=y\in \Aline ox=\Aline o{x'}$ and we are done. So, assume that $x\ne y$. In this case $x\ne o$. If $x\in D$, then $y\in\Aline ox=D$ and $y'=y\in\Aline ox=\Aline o{x'}$. So, assume that $x\notin D$. If $x=a$, then $d_y\in \Aline ya\cap D=\Aline yx\cap D=\{o\}$ and $y'\in \Aline o{a'}=\Aline o{x'}$. If $y=a$, then $d_x\in \Aline xa\cap D=\Aline xy\cap D=\{o\}$ and $x'\in \Aline o{a'}=\Aline o{y'}$, which implies $y'\in\Aline o{x'}$. So, assume that $a\notin\{x,y\}$. Then the lines $D,L_x,L_a,L_y$ are distinct and parallel. 

\begin{picture}(200,135)(-120,-20)
\linethickness{=0.6pt}
\put(0,0){\color{teal}\line(1,0){160}}
\put(165,-3){\color{teal}$D$}
\put(0,30){\color{teal}\line(1,0){160}}
\put(165,27){\color{teal}$L_x$}
\put(0,60){\color{teal}\line(1,0){160}}
\put(165,57){\color{teal}$L_a$}
\put(0,90){\color{teal}\line(1,0){160}}
\put(165,87){\color{teal}$L_y$}

\put(105,90){\color{violet}\line(1,-2){45}}
\put(15,90){\color{violet}\line(3,-2){135}}
\put(30,0){\color{blue}\line(1,2){30}}
\put(30,0){\color{blue}\line(3,2){90}}
\put(60,0){\color{red}\line(1,2){45}}
\put(60,0){\color{red}\line(-1,2){45}}

\put(150,0){\color{violet}\circle*{3}}
\put(148,-10){$d_y$}
\put(30,0){\color{cyan}\circle*{3}}
\put(27,-10){$d_c$}
\put(60,0){\color{red}\circle*{3}}
\put(57,-8){$o$}
\put(45,30){\circle*{3}}
\put(49,32){$x'$}
\put(75,30){\circle*{3}}
\put(69,32){$x$}
\put(60,60){\circle*{3}}
\put(61,62){$a'$}
\put(120,60){\circle*{3}}
\put(121,62){$a$}
\put(15,90){\circle*{3}}
\put(12,95){$y'$}
\put(105,90){\circle*{3}}
\put(103,95){$y$}

\end{picture}

Consider the points $x,x'\in L_x$, $a,a'\in L_a$, $y,y'\in L_y$, and observe that $d_y\in \Aline ya\cap\Aline {y'}{a'}\cap D$ and $d_x\in \Aline xa\cap\Aline{x'}{a'}\cap D$. Since $X$ is para-Playfair and para-Desarguesian, we can apply Proposition~\ref{p:paraMouf=>z} and conclude that $\varnothing\ne \Aline xy\cap\Aline {x'}{y'}\subseteq\Aline xy\cap D=\{o\}$ and hence $o\in \Aline {x'}{y'}$ and finally $y'\in\Aline o{x'}$.  
\end{proof}

Now we are ready to prove that the line $D$ is shear. Given any distinct points $a,a'\in X$ with $\Aline a{a'}\cap D=\varnothing$, we should construct an automorphism $\Phi:X\to X$ such that $\Phi(a)=a'$ and $\Fix(\Phi)=D$. Since $X$ is a para-Playfair plane containing the disjoint lines $D$ and $\Aline a{a'}$, for every point $x\in X$ there exists a unique line $L_x\subseteq X$ such that $x\in L_x$ and $L_x\parallel D$. 

Fix any point $o\in D$. Since $X$ is $3$-long, there exists a point $b\in \Aline oa\setminus\{o,a\}$. By Lemma~\ref{l:part-automorphism}, there exists an injective function $\varphi$ such that $aa'\in\varphi$, $\dom[\varphi]=(X\setminus\Aline a{a'})\cup\{a\}$, $\rng[\varphi]=(X\setminus L_a)\cup\{a'\}$, 
and for every pairs $xx',yy'\in\varphi$ the conditions (1)--(3) of the lemma  are satisfied. The conditions (2) and (3) ensure that the point $b'\defeq \varphi(b)$ is a unique point of the set $\Aline o{a'}\cap L_b$ and hence $b'\ne b$. By Lemma~\ref{l:part-automorphism}, there exists an injective function $\psi$ such that $bb'\in\psi$, $\dom[\psi]=(X\setminus L_b)\cup\{b\}$, $\rng[\psi]=(X\setminus\Aline b{b'})\cup\{b'\}$, 
and for every pairs $xx',yy'\in\psi$ the conditions (1)--(3) of the lemma are satisfied.

\begin{claim}\label{cl:phi=psi} For every point $x\in X\setminus(L_a\cup L_b)$, we have $\varphi(x)=\psi(x)$.
\end{claim}

\begin{proof} If $x\in D$, then $\varphi(x)=x=\psi(x)$, by the condition (1) holding for the maps $\varphi$ and $\psi$. So, we assume that $x\notin D$. By the Proclus Postulate~\ref{p:Proclus-Postulate}, there exist points $d_x\in \Aline xa\cap D$ and $\delta_x\in \Aline xb\cap B$. The condition (3) holding for the maps $\varphi$ and $\psi$ imply that $x'\defeq \varphi(x)\in \Aline {d_x}{a'}$ and $\psi(x)\in\Aline{\delta_x}{b'}$.

\begin{picture}(200,130)(-120,-20)
\linethickness{=0.6pt}
\put(0,0){\color{teal}\line(1,0){160}}
\put(165,-3){\color{teal}$D$}
\put(0,30){\color{teal}\line(1,0){160}}
\put(165,27){\color{teal}$L_x$}
\put(0,60){\color{teal}\line(1,0){160}}
\put(165,57){\color{teal}$L_a$}
\put(0,90){\color{teal}\line(1,0){160}}
\put(165,87){\color{teal}$L_b$}

\put(105,90){\color{violet}\line(1,-2){45}}
\put(15,90){\color{violet}\line(3,-2){135}}
\put(30,0){\color{blue}\line(1,2){30}}
\put(30,0){\color{blue}\line(3,2){90}}
\put(60,0){\color{red}\line(1,2){45}}
\put(60,0){\color{red}\line(-1,2){45}}

\put(150,0){\color{violet}\circle*{3}}
\put(148,-8){$o$}
\put(30,0){\color{cyan}\circle*{3}}
\put(27,-10){$d_x$}
\put(60,0){\color{red}\circle*{3}}
\put(57,-10){$\delta_x$}
\put(45,30){\circle*{3}}
\put(49,32){$x'$}
\put(75,30){\circle*{3}}
\put(69,32){$x$}
\put(60,60){\circle*{3}}
\put(61,62){$a'$}
\put(120,60){\circle*{3}}
\put(121,62){$a$}
\put(15,90){\circle*{3}}
\put(12,95){$b'$}
\put(105,90){\circle*{3}}
\put(103,95){$b$}

\end{picture} 

Consider the disjoint parallel lines $L_a,L_b,L_x,D$ and distinct points $a,a'\in L_a$, $b,b'\in L_b$, $x,x'\in L_x$. Observe that $o\in \Aline ab\cap\Aline {a'}{b'}\cap D$ and $d_x\in \Aline ax\cap\Aline {a'}{x'}\cap D$. By Proposition~\ref{p:paraMouf=>z}, $\varnothing \ne\Aline bx\cap\Aline {b'}{x'}\cap D\subseteq \Aline bx\cap D=\{\delta_x\}$ and hence $\delta_x\in \Aline {b'}{x'}$. Then $\varphi(x)=x'\in \Aline {b'}{\delta_x}\cap L_x=\{\psi(x)\}$ and hence $\varphi(x)=\psi(x)$.
\end{proof}

Claim~\ref{cl:phi=psi} implies that $\Phi\defeq\varphi\cup\psi$ is a well-defined bijection of the set $X$ such that $aa'\in \Phi$. We claim that $\Phi$ is an automorphism of the liner $X$. By Theorem~\ref{t:liner-isomorphism<=>}, it suffices to check that for any line $\Lambda$ in $X$, its image in a line in $X$. If $\Lambda\cap D=\varnothing$, then $\Phi[\Lambda]=\Lambda$, by the properties (2) of the functions $\varphi$ and $\psi$. So assume that $\Lambda\cap D$ contains some point $o$. Since $X$ is $4$-long, there exists a point $c\in \Lambda\setminus(D\cup L_a\cup L_b)$. Let  $c'\defeq\Phi(c)=\varphi(c)=\psi(c)$ be its image. We claim that $\Phi[\Lambda]\subseteq \Aline o{c'}$. Take any $x\in \Lambda$.  If $x\notin L_a$, then $\Phi(x)=\varphi(x)\in \Aline o{c'}$, by the condition (3) for the function $\varphi$. If $x\in L_a$, then $\Phi(x)=\psi(x)\in \Aline o{c'}$ by the condition (3) for the function $\psi$. Therefore, $\Phi[\Lambda]\subseteq\Aline o{c'}$. Since $\Phi$ is bijective, $\Phi[\Lambda]=\Aline o{c'}$ is a line in $X$, witnessing that $\Phi$ is an automorphism of the liner $X$ and the line $D$ in $X$ is shear.
\end{proof}

\begin{theorem}\label{t:para-Desarg-compreg<=>}  A para-Desarguesian liner is completely regular if and only if it is regular and para-Playfair.
\end{theorem}

\begin{proof} The ``only if'' part follows from Theorem~\ref{t:spread=projective1}. To prove the ``if'' part, assume that $X$ is para-Desarguesian, para-Playfair, and regular.
By Proposition~\ref{p:Playfair=>para-Playfair=>Proclus} and Theorem~\ref{t:Proclus<=>}, the para-Playfair plane $X$ is Proclus, proaffine, $3$-proregular and $3$-ranked. By Theorem~\ref{t:spread=projective1}, the complete regularity of $X$ will follow as soon as we check that for every concurrent Bolyai lines $A,B$ in $X$, every line $L$ in the plane $P\defeq\overline{A\cup B}$ is Bolyai. This is clear if $L\in A_\parallel \cup B_\parallel$. So, assume that $L\notin A_\parallel\cup B_\parallel$.

By Theorem~\ref{t:spreading<=>}, the Bolyai lines $A,B$ are spreading. Then for any point $x\in X$ there exists unique lines $A_x\in\A_\parallel$ and $B_x\in B_\parallel$ such that $x\in A_x\cap B_x$. Fix any distinct points $o,e\in L$. It follows from $L\notin A_\parallel\cup B_\parallel$ that $A_o\cap A_e=\varnothing=B_o\cap B_e$. By Proclus Postulate~\ref{p:Proclus-Postulate}, there exist unique points $u\in A_o\cap B_e$ and $v\in A_e\cap B_o$.

By Theorem~\ref{t:para-Desarg<=>shear}, every line in the para-Desarguesian para-Playfair plane $P$ is shear. Then there exist two automorphisms $\Phi_o$ and $\Phi_v$ of the plane $P$ such that $\Phi_o(o)=u$, $\Fix(\Phi_o)=A_e$, and $\Phi_v(v)=e$, $\Fix(\Phi_v)=A_o$.

\begin{claim}\label{cl:Phi[Ax]=Ax} For every $x\in X$ we have $\Phi_o[A_x]=A_x=\Phi_v[A_x]$.
\end{claim}

\begin{proof} If $x\in A_e$, then the equality $\Phi_o[A_x]=A_x$ follows from $A_x=A_e=\Fix(\Phi_o)$. If $x\in A_o$, then $\Phi_o[A_x]\cap A_e=\Phi_o[A_o]\cap\Phi[A_e]=\Phi_o[A_o\cap A_e]=\Phi_0[\varnothing]=\varnothing$ and hence $A_x=A_o$ and $\Phi_o[A_x]=\Phi_o[A_o]$ are two lines that contain a common point $e=\Phi_o(o)$ and are disjoint with the line $A_e$. The Proclus Axiom ensures that $A_x=\Phi_o[A_x]$. So, assume that $x\notin A_o\cup A_e$. It follows from $x\notin A_e=\Fix(\Phi_o)$ that the point $x'\defeq \Phi_o(x)\notin\Phi[A_e]=A_e$ is distinct from $x$ and hence $\Aline x{x'}$ is a line in $X$. Assuming that $\Aline x{x'}\ne A_x$, we can apply the Proclus Postulate~\ref{p:Proclus-Postulate} and find  points $a_o\in \Aline x{x'}\cap A_o$ and $a_e\in \Aline x{x'}\cap A_e\subseteq\Fix(\Phi_o)$. Then $$\Phi_o(a_o)\in \Phi_o[\Aline x{x'}]\cap A_o=\Phi_o[\Aline x{a_e}]\cap A_o=\Aline {x'}{a_e}\cap A_o=\Aline x{x'}\cap A_o=\{a_o\},$$ which contradicts $\Fix(\Phi_o)=A_e$. This contradiction shows that $\Aline x{x'}=A_x$. Then $A_x$ and $\Phi_o[A_x]$ are two lines that contain the point $x'$ and are disjoint with the line $A_e$. The Proclus Axiom ensures that $\Phi_o[A_x]=A_x$. By analogy we can prove that $\Phi_e[A_x]=A_x$.
\end{proof}

\begin{picture}(100,145)(-140,-15)
\linethickness{0.6pt}
\put(-30,0){\color{teal}\line(1,0){140}}
\put(115,-3){\color{teal}$A_o$}
\put(-30,30){\color{teal}\line(1,0){140}}
\put(115,27){\color{teal}$A_e$}
\put(0,0){\color{cyan}\line(0,1){110}}
\put(-4,115){\color{cyan}$B_o$}
\put(30,0){\color{cyan}\line(0,1){110}}
\put(26,115){\color{cyan}$B_e$}
\put(0,0){\color{red}\line(1,1){100}}
\put(105,105){\color{red}$L$}
\put(30,90){\line(1,0){60}}
\put(-15,0){\line(1,2){45}}

\put(-15,0){\circle*{3}}
\put(-17,-9){$a$}
\put(0,0){\circle*{3}}
\put(-2,-9){$o$}
\put(30,0){\circle*{3}}
\put(28,-9){$u$}
\put(0,30){\circle*{3}}
\put(-8,32){$v$}
\put(30,30){\circle*{3}}
\put(23,32){$e$}
\put(30,90){\circle*{3}}
\put(32,92){$c'$}
\put(90,90){\circle*{3}}
\put(86,93){$c$}

\end{picture}

Consider the automorphism $\Phi\defeq \Phi_v\circ\Phi_o$ of the plane $P$. We claim that $\Phi[L]\cap L=\varnothing$. To derive a contradiction, assume that $L\cap\Phi[L]$ contains some point $c$. It follows that $\{c\}\in L\cap \Phi[L]\cap L_c$ and hence $\Phi[L]\cap L_c=\{c\}$. 
Consider the point $c'\defeq \Phi_o(c)$. Claim~\ref{cl:Phi[Ax]=Ax} ensures that $c'=\Phi_o(c)\in \Phi_o[L_c]=L_c$ and $\Phi(c)=\Phi_v(c')\in \Phi_v[L_c]\cap \Phi[L]=L_c\cap\Phi[L]=\{c\}$.  

Assuming that $c\in A_e$, we conclude that $c=\Phi(c)=\Phi_v\circ \Phi_o(c)=\Phi_v(c)$ and hence $c\in \Fix(\Phi_v)=A_o$, which contradicts $c\in A_e$. This contradiction shows that $c\notin A_e$ and hence $c'= \Phi_o(c)\notin \Phi_o[A_e]=A_e$.  Observe that $c'=\Phi_o(c)\in \Phi_o[L]=\Phi_o[\Aline oe]=\Aline ue=B_e$. By Proclus Postulate~\ref{p:Proclus-Postulate}, there exists a point $a\in \Aline {c'}v\cap A_o$. Assuming that $a=o$, we conclude that $c'\in \Aline av=\Aline ov=B_o$, which contradicts $c'\in B_e$. This contradiction shows that $a\ne o$, which implies $\Aline ae\cap\Aline oe=\{e\}$. Then $c=\Phi(c)=\Phi_v(c')\in \Phi_v[\Aline av]=\Aline ae$ and hence $c\in \Aline ae\cap \Aline oe=\{e\}$, which contradicts $c\notin A_e$. This contradiction shows that the lines $L$ and $\Phi[L]$ is the plane $P$ are disjoint. By Theorem~\ref{p:lines-in-para-Playfair}, the line $L$ is Bolyai. Therefore, the para-Playfair regular liner $X$ is bi-Bolyai. By Theorem~\ref{t:spread=projective1}, the bi-Bolyai para-Playfair regular liner $X$ is completely regular.
\end{proof}

\begin{corollary} A $4$-long para-Desarguesian liner is completely regular if and only if it is para-Playfair.
\end{corollary}

\begin{proof} The ``only if'' part follows immediately from Theorem~\ref{t:spread=projective1}. To prove the ``if'' part, assume that $X$ is a $4$-long para-Desarguesian para-Playfair liner. We claim that $X$ is bi-Bolyai. Given any concurrent Bolyai lines $A,B$ in $X$, we need to check that every line $L$ in the plane $P\defeq\overline{A\cup B}$ is Bolyai in $X$. Since the plane $P$ contains disjoint lines, it is not projective. By Proposition~\ref{p:Playfair=>para-Playfair=>Proclus} and Theorem~\ref{t:Proclus<=>}, the para-Playfair plane $P$ is Proclus and $3$-proregular. Being a $4$-long plane, the $3$-proregular liner $P$ is regular. By Theorem~\ref{t:para-Desarg-compreg<=>} and \ref{t:spread=projective1}, the regular para-Desarguesian para-Playfair plane $P$ is completely regular and bi-Bolyai. Then every line $L$ in the liner $P$ is Bolyai and hence there exists a line $\Lambda\subseteq P\setminus L$. By Proposition~\ref{p:Bolyai-in-para-Playfair}, the line $L$ is Bolyai in $X$, witnessing that the liner $X$ is bi-Bolyai. By 
Theorem~\ref{t:4long+pP+bB=>regular}, the $4$-long para-Playfair bi-Bolyai liner $X$ is regular, and by Theorem~\ref{t:spread=projective1}, the regular para-Playfair bi-Bolyai liner $X$ is completely regular.
\end{proof}

\begin{corollary}\label{c:para-Desarguesian=>Thalesian} Every para-Desarguesian para-Playfair liner $X$ is Thalesian.
\end{corollary} 

\begin{proof} By Proposition~\ref{p:Playfair=>para-Playfair=>Proclus}, Theorem~\ref{t:Proclus<=>}, and Proposition~\ref{p:k-regular<=>2ex}, the para-Playfair liner $X$ is Proclus, proaffine, $3$-proregular and $3$-ranked. By Proposition~\ref{p:Thales<=>planeT}, the Thalesian property of the proffine regular liner $X$ will follow as soon as we check that every $3$-long plane $\Pi$ in $X$ is Thalesian. To prove that $\Pi$ is Thalesian, take any disjoint lines $A,B,C\subset \Pi$ and points $a,a'\in A$, $b,b'\in B$, $c,c'\in C$ such that $\Aline ab\cap\Aline {a'}{b'}=\varnothing$ and $\Aline bc\cap\Aline {b'}{c'}=\varnothing$. We have to prove that $\Aline ac\cap\Aline {a'}{c'}=\varnothing$. By Corollary~\ref{c:parallel-lines<=>}, $\Aline ab\parallel \Aline {a'}{b'}$ and $\Aline bc\parallel \Aline {b'}{c'}$. If $\Aline ab=\Aline bc$, then $\Aline {a'}{b'}\parallel \Aline ab=\Aline bc\parallel\Aline {b'}{c'}$ implies $\Aline {a'}{b'}=\Aline{b'}{c'}=\Aline{a'}{c'}$ and hence $\Aline ac\cap\Aline{a'}{c'}=\Aline ab\cap\Aline{a'}{b'}=\varnothing$. So, assume that $\Aline ab\ne\Aline bc$.

Since $\Pi$ is para-Playfair and $\Aline ab\cap\Aline {a'}{b'}=\varnothing=\Aline bc\cap\Aline {b'}{c'}$, the lines $\Aline ab$ and $\Aline bc$ are spreading and hence Bolyai. By Theorem~\ref{t:para-Desarg-compreg<=>} and \ref{t:spread=projective1}, the para-Desarguesian para-Playfair regular liner $\Pi$ is completely regular and bi-Bolyai. Then every line in the plane $\Pi=\overline{A\cup B}$ is Bolyai and hence the Proclus liner $\Pi$ is Playfair. By Theorem~\ref{t:para-Desarg<=>shear}, the para-Desarguesian para-Playfair plane $\Pi$ is shear. By Theorems~\ref{t:shear<=>alternative} and \ref{t:paraD<=>translation}, the shear Playfair plane $\Pi$ is translation and Thalesian. Then 
$\Aline ab\cap\Aline {a'}{b'}=\varnothing=\Aline bc\cap\Aline {b'}{c'}=\varnothing$ implies $\Aline ac\cap\Aline {a'}{c'}=\varnothing$.
\end{proof}

\begin{example}[Ivan Hetman, 2025]\label{ex:Thales13} The $13$-element field $\IF_{13}$ endowed with the family of lines
$\mathcal L\defeq\big\{x+\{0,3,4\},x+\{0,5,7\}:x\in\IF_{13}\big\}$ is an example of a non-regular para-Desarguesian Steiner (and hence proaffine) ranked plane, which is not Thalesian.
\end{example}

\begin{proof} Computer calculations show that the liner $(\IF_{13},\mathcal L)$ is a  para-Desarguesian plane. To see that it is not Thalesian, consider the disjoint lines
$A\defeq\{0,1,4\}$, $B\defeq\{2,3,6\}$, $C\defeq\{5,7,12\}$, and points $a\defeq 0$, $a'\defeq 1$, $b\defeq 2$, $b'\defeq 3$, $c\defeq 7$, $c'\defeq 5$. Observe that $\Aline ab\cap\Aline {a'}{b'}=\{0,2,7\}\cap\{1,3,8\}=\varnothing$, $\Aline bc\cap\Aline{b'}{c'}=\{0,2,7\}\cap\{3,5,10\}=\varnothing$, but $\Aline ac\cap\Aline{a'}{c'}=\{0,2,7\}\cap\{1,2,5\}=\{2\}\ne\varnothing$, witnessing that the liner is not Thalesian.
\end{proof}

\begin{exercise} Show that a non-projective completely regular plane is para-Desarguesian if and only if its spread completion has Lenz type at least $(\,|\,)$.
\end{exercise}

\section{Moufang affine spaces}

In this section we present a characterization of Moufang affine spaces. We recall that an {\em affine space} is any $3$-long affine regular liner of rank at least $3$. An affine space $X$ is {\em hypersymmetric} if for any triangle $abc$ there exists a hyperfixed involutive automorphism $A$ of $X$ such that $Aabc=cba$.

\begin{theorem}\label{t:affine-Moufang<=>} For any affine space $X$, the following conditions are equivalent:
\begin{enumerate}
\item $X$ is Moufang;
\item $X$ is little-Desarguesian;
\item $X$ is para-Desarguesian;
\item $X$ is shear;
\item $X$ is hypersymmetric;
\item $X$ is inversive-dot;
\item every ternar of $X$ is Moufang;
\item every ternar of $X$ is Moufang-dot;
\item every ternar of $X$ is inversive-dot;
\item every ternar of $X$ is linear, distributive, associative-plus, and alternative-dot; 
\item some ternar of $X$ is linear, distributive, associative-plus, and alternative-dot;
\item the spread completion of $X$ is Moufang.
\end{enumerate}
\end{theorem}

\begin{proof} By Corollary~\ref{c:affine-spread-completion}, the affine space is completely regular.
\smallskip

The equivalences $(12)\Leftrightarrow(1)\Leftrightarrow(2)$ follow from Theorem~\ref{t:compreg-Moufang<=>}. The implications $(1)\Ra(3)\Ra(4)$ follow from Proposition~\ref{p:Moufang=>paraD} and Theorem~\ref{t:para-Desarg<=>shear}, respectively. The equivalence of the conditions $(4)$--$(11)$  was proved in Theorem~\ref{t:inversive-dot}.
\smallskip

$(10)\Ra(12)$ Since the affine space $X$ is completely regular, its spread completion $\overline X$ is a $4$-long projective liner of rank $\|\overline X\|=\|X\|\ge 3$, according to Corollary~\ref{c:procompletion-rank}. If $\|\overline X\|\ge 4$, then the projective liner $X$ is Desarguesian and Moufang, by Theorem~\ref{t:proaffine-Desarguesian} and Proposition~\ref{p:Desarg=>Moufang}. So, assume that $\|X\|=\|\overline X\|=3$.  Chose any affine base $uow$ in the affine plane $X$. The condition (10) ensures that the ternar of the based affine plane $(X,uow)$ is linear, distributive, associative-plus, and alternative-dot. Let $e$ be the diunit of the affine base $uow$. Then $(\overline X,uowe)$ is a based projective plane whose ternar coincides with the ternar of the based affine plane $(X,uow)$ and hence is linear, distributive, associative-plus and alternative-dot. By Theorem~\ref{t:proj-Moufang<=>}, the projective plane $\overline X$ is Moufang.
\end{proof}

\section{Finite  Moufang  liners}

\begin{theorem}\label{t:Moufang-finite<=>} A finite Proclus liner is Moufang if and only if it is Pappian.
\end{theorem}

\begin{proof} If $X$ is Pappian, then $X$ is Desarguesian and Moufang, by Theorems~\ref{t:Hessenberg-proaffine} and Proposition~\ref{p:Desarg=>Moufang}. Now assume that a Proclus liner $X$ is Moufang.  By Proposition~\ref{p:Papp<=>4-Papp}, it suffices to check that every $4$-long plane $P$ in $X$ is Pappian.  Since $X$ is Moufang, so is the plane $P$ in $X$. By Theorem~\ref{t:Proclus<=>}, the $4$-long finite Proclus plane $P$ is regular. By Theorem~\ref{t:compreg-Moufang<=>}, the $4$-long Moufang regular liner $P$ is completely regular. By Theorem~\ref{t:compreg-Moufang<=>}, the spread  completion $\overline P$ of the completely regular Moufang liner $P$ is projective and Moufang. Since $X$ is finite, so is the plane $P$ and its spread completion $\overline P$. Choose any line $H$ in $\overline P$ and consider the affine plane $A\defeq \overline P\setminus H$. Since $P$ is $4$-long, the projective plane $\overline P$ is $4$-long and the affine plane $A$ is $3$-long. By Theorem~\ref{t:proj-Moufang<=>}, the Playfair plane $A=\overline P\setminus H$ is Moufang. By Theorem~\ref{t:affine-Moufang<=>}, some ternar $R$ of the the Moufang affine plane $A$ is linear, distributive, associative-plus and alternative-dot. Since $A$ is finite, so is the ternar $R$. By Artin--Zorn Theorem~\ref{t:Artin-Zorn}, the ternar $R$ is associative and commutative. By Theorem~\ref{t:commutative-dot<=>}, the affine plane $A$ is Pappian. By Theorem~\ref{t:projPapp<=>}, the projective completion $\overline P$ of the Pappian affine plane $A$ is Pappian. Since the completely regular plane $P$ has Pappian spread completion, it is Pappian, by Theorem~\ref{t:Pappian-completion}, the plane $P$ is Pappian, and by Proposition~\ref{p:Papp<=>4-Papp}, the liner $X$ is Pappian.   
\end{proof}

\section{Homogeneity of Moufang liners}

We recall that a \defterm{projective reper} in a projective space $X$ is a function $r\subseteq X\times X$ whose domain $\dom[r]$ is a maximal independent set in $X$ containing a unique point $o=r(o)$ such that $r(x)\in\Aline ox\setminus\{o,x\}$ for all  $x\in \dom[r]\setminus\{o\}$. 

By Corollary~\ref{c:Des-iso-reper}, for any projective repers $r,r'$ in a Desarguesian projective space $X$, every bijective function $\varphi:\dom[r]\cup\rng[r]\to\dom[r']\cup\rng[r']$ with $\varphi\circ r=r'\circ \varphi$ extends to an automorphism $\Phi:X\to X$ of the projective space $X$. We shall prove that this homogeneity property holds, more generally, for non-Fano Moufang projective spaces.

We shall derive this fact from $4$-homogeneity of non-Fano Moufang planes, proved by Bruck and Kleinfeld \cite{BruckKleinfeld1951}.

\begin{definition} A liner $X$ is defined to be
\begin{itemize} 
\item \index{1-homogeneous liner}\index{liner!1-homogeneous}\defterm{$1$-homogeneous} if for any points $a,b\in X$ there exists an automorphism $\Phi:X\to X$ such that $\Phi(a)=b$;
\item \index{2-homogeneous liner}\index{liner!2-homogeneous}\defterm{$2$-homogeneous} if for any doubletons\footnote{A {\em doubleton} is an ordered pair $ab$ of distinct points.} $ab,a'b'$ in $X$, there exists an automorphism $\Phi:X\to X$ such that $\Phi ab=a'b'$;
\item \index{3-homogeneous liner}\index{liner!3-homogeneous}\defterm{$3$-homogeneous} if for any triangles $abc,a'b'c'$ in $X$, there exists an automorphism $\Phi:X\to X$ such that $\Phi abc=a'b'c'$;
\item \index{4-homogeneous liner}\index{liner!4-homogeneous}\defterm{$4$-homogeneous} if for any quadrangles\footnote{A {\em quadrangle} is a quadruple $abcd$ of distinct points of a liner such that $\|\{x,y,z\}\|=\|\{a,b,c,d\}\|=3$ for any distinct points $x,y,z\in\{a,b,c,d\}$.} $abcd,a'b'c'd'$ in $X$, there exists an automorphism $\Phi:X\to X$ such that $\Phi abcd=a'b'c'd'$.
\end{itemize}
\end{definition}

\begin{theorem}[Hall, 1943]\label{t:Hall-3hom} Every $3$-long Moufang projective plane is $3$-homo\-geneous.
\end{theorem}

\begin{proof} Given any triangles $hov$ and $h'o'v'$ in a $3$-long Moufang plane $\Pi$, we need to construct an automorphism $A:X\to X$ of $X$ such that $Ahow=h'o'v'$. Since the liner $\Pi$ is $3$-long and has rank $\|\Pi\|=3$, there exists a line $L_1\subseteq \Pi$ such that $\{h,h'\}\cap L_1=\varnothing$. Since the projective plane $\Pi$ is Moufang, the line $L_1$ is translation and hence there exists an automorphism $\Phi_1$ of $\Pi$ such that $\Phi_1(h)=h'$. Consider the triangle $h_1o_1v_1\defeq\Phi_1 hov$. Since $\Pi$ is $3$-long, there exists a line $L_2\subset \Pi$ such that $h'=h_1\in L_2$ and $v_1,v'\notin L_2$. Since $\Pi$ is Moufang, the line $L_2$ is translation and hence there exists an automorphism $\Phi_2$ of $\Pi$ such that $L_2\subseteq \Fix(\Phi_2)$ and $\Phi_2(v_1)=v'$. Consider the triangle $h_2o_2v_2\defeq\Phi_2 h_1o_1v_1$. Since $h_2o_2v_2$ and $h'o'v'$ are triangles, the points $o_2$ and $o'$ do not belong to the line $L_3\defeq\Aline {h_2} {v_2}=\Aline {h'}{v'}$.  Since $\Pi$ is Moufang, the line $L_3$ is translation and hence there exists an automorphism $\Phi_3$ of $\Pi$ such that $L_3\subseteq \Fix(\Phi_3)$ and $\Phi_2(o_2)=o'$. Then the automorphism $\Phi=\Phi_3\Phi_2\Phi_1$ has the required property: $\Phi hov=\Phi_3\Phi_2\Phi_1hov=h'o'v'$.
\end{proof}

Non-Fano Moufang projective liners have even better homogeneity properties.

\begin{theorem}[Bruck, Kleinfeld, 1951]\label{t:Moufang=>4-homogen} Every $3$-long non-Fano Moufang projective plane is $4$-homo\-geneous.
\end{theorem}

\begin{proof} Given two quadrangles $uowe$ and $u'o'w'e'$ in a non-Fano Moufang plane $\Pi$, we need to find an automorphism $\Phi$ of $\Pi$ such that $\Phi uowe=u'o'w'e'$. By projectivity of $\Pi$, there exist unique points $h\in \Aline ou\cap\Aline we$, $v\in \Aline ow\cap\Aline ue$, $h'\in \Aline {o'}{u'}\cap\Aline {w'}{e'}$, $v'\in \Aline {o'}{w'}\cap\Aline {u'}{e'}$. By Theorem~\ref{t:Hall-3hom}, there exists an automorphism $\Phi:\Pi\to\Pi$ such that $\Phi hov=h'o'v'$.

Since the projective plane $\Pi$ is Moufang, the ternars $\Delta,\Delta'$ are linear, distributive, associative-plus, and alternative-plus, according to Theorem~\ref{t:proj-Moufang<=>}. In other terminology, the biloops of these ternars are alternative division rings. Since the projective plane $\Pi$ is not Fano, these alternative division rings are not Boolean, by Proposition~\ref{p:DesFano<=>}.
By Theorem~\ref{t:Hall-3hom}, we lose no generality assuming that $hov=h'o'v'$. Then $uow$ and $u'o'w'$ are two affine bases in the affine plane $\Pi\setminus\Aline hv$ such that $\Aline ou=\Aline{o'}{u'}$ and $\Aline ow=\Aline{o'}{w'}$. By Proposition~\ref{p:ternars-isotopic}, the ternars $\Delta=\Aline oe\setminus\Aline hv$ and $\Delta'=\Aline{o'}{e'}\setminus\Aline hv$ of these affine bases are isotopic. Then the biloops of these (linear) ternars also are isotopic, by Proposition~\ref{p:ternar-isotopic<=>biloops}. 

 By Schafer's Theorem~\ref{t:Schafer}, the isotopic biloops of the non-Boolean linear ternars $\Delta,\Delta'$ are isomorphic. Then the ternars $\Delta,\Delta'$, being linear, are isomorphic, too. By Theorem~\ref{t:p-isomorphic<=>ternar-isomorphic}, there exists an automorphism $\Psi:\Pi\to\Pi$ such that $\Psi uowe=u'o'w'e'$.
\end{proof}

\begin{theorem}\label{t:Moufang-reper-hom} For any projective repers $r,r'$ in a non-Fano Moufang projective space $X$, every bijective function $\varphi:\dom[r]\cup\rng[r]\to\dom[r']\cup\rng[r']$ with $\varphi\circ r=r'\circ \varphi$ extends to an automorphism $\Phi:X\to X$ of the projective space $X$.
\end{theorem}

\begin{proof} If $\|X\|\ne 3$, then the projective space $X$ is Desarguesian (by Theorem~\ref{t:proaffine-Desarguesian}) and the existence of the automorphism $\Phi$ follows from Corollary~\ref{c:Des-iso-reper}. So, assume that $\|X\|=3$ and hence $X$ is a projective plane. Since $\|X\|=3$, $r=\{oo,uh,wv\}$ and $r'=\{o'o',u'h',w'v'\}$ for some points $o,u,h,w,v,o',u',h',w',v'\in\Pi$ with $\varphi: ouhwv\mapsto o'u'h'w'v'$. By the projectivity of $X$, there exist unique points $e\in \Aline uv\cap\Aline wh$ and $e'\in \Aline {u'}{v'}\cap\Aline {w'}{h'}$. By Theorem~\ref{t:Moufang=>4-homogen}, there exists an automorphism $\Phi:\Pi\to\Pi$ such that $\Phi uowe=u'o'w'e'$. Then $\Phi(h)\in \Phi[\Aline ow\cap\Aline ue]=\Aline{o'}{w'}\cap\Aline {u'}{e'}=\{h'\}$. By analogy we can prove that $\Phi(v)=v'$. Therefore, $\Phi$ is a required automorphism of $\Phi$ that extends the bijection $\varphi$.
\end{proof}

\begin{corollary}\label{c:Moufang=>4hom} Any non-Fano Moufang projective space is $4$-homogeneous.
\end{corollary}

\begin{proof} Given two quadrangles $uowe$ and $u'o'w'e'$ in a non-Fano Moufang projective space $X$, we need to find an automorphism $\Phi$ on $X$ such that $\Phi uowe=u'o'w'e'$. Since $uowe$ and $u'o'w'e'$ are quadrangles in $X$, their flat hulls $\Pi\defeq\overline{\{o,u,w,e\}}$ and $\Pi'\defeq\overline{\{u',o',w',e'\}}$ in the projective space $X$ are projective planes. By the projectivity of the planes $\Pi,\Pi'$, there exist unique points $h\in\Aline ou\cap\Aline we$, $v\in\Aline ow\cap\Aline ue$, $h'\in\Aline {o'}{u'}\cap\Aline {w'}{e'}$, $h'\in\Aline {o'}{u'}\cap\Aline {w'}{e'}$. Then $\{oo,uh,wv\}$ and $\{o'o',u'h',w'v'\}$ are two projective repers in the projective planes $\Pi,\Pi'$, respectively. By Kuratowski--Zorn Lemma, these projective repers can be enlarged to projective repers $r,r'$ of the Moufang projective space $X$. Since $\dom[r]$ and $\dom[r']$ are maximal independent sets in $X$,  Corollary~\ref{c:Max=dim} ensures that $|\dom[r]|=\|X\|=|\dom[r']|$, and hence there exists a bijective map $\varphi:\dom[r]\cup\rng[r]\to\dom[r']\cup\rng[r']$ such that $\varphi\circ r=r'\circ\varphi$ and $\varphi:ouhwv\mapsto o'u'h'w'v'$. By Theorem~\ref{t:Moufang-reper-hom}, the bijection $\varphi$ extends to an automorphism $\Phi:X\to X$ of the projective space. Then $\Phi uohwv=u'o'h'w'v'$ and $\Phi(e)\in \Phi[\Aline uv\cap\Aline wh]=\Aline {u'}{v'}\cap\Aline {w'}{h'}=\{e'\}$ and hence $\Phi$ is a required automorphism of the projective space $X$ with $\Phi uowe=u'o'w'e'$.
\end{proof}

\begin{corollary}\label{c:ternars-in-non-Fano-Moufang-plane-are-isomorphic} Any two ternars of a non-Fano Moufang projective space are isomorphic.
\end{corollary} 

\begin{proof} Let $R,R'$ be two ternars of a non-Fano Moufang projective space $X$. Then there exist two planes $\Pi,\Pi'$ in $X$ and projective bases $uowe$ and $u'o'w'e'$ in the projective planes $\Pi,\Pi'$, such that $R,R'$ are isomorphic to the ternars of the based projective planes $(\Pi,uowe)$, $(\Pi',u'o'w'e')$, respectively. By Corollary~\ref{c:Moufang=>4hom}, there exists an automorphism $\Phi$ of the projective space $X$ such that $\Phi uowe=u'o'w'e'$. Then $\Phi{\restriction}_\Pi$ is an isomorphism of the based projective planes $(\Pi,uowe)$ and $(\Pi',u'o'w'e')$. By Theorem~\ref{t:p-isomorphic<=>ternar-isomorphic}, the ternars of these based projective planes are isomorphic and hence  the ternars $R$ and $R'$ are isomorphic, too.
\end{proof}

\begin{corollary}\label{c:Moufang=>3hom} Any non-Fano Moufang affine space is $3$-homogeneous.
\end{corollary}

\begin{proof} Let $X$ be a Moufang affine space. Given two triangles $uow$ and $u'o'w'$ in $X$, we have to find an automorphism $\Phi$ of $X$ such that $\Phi uow=u'o'w'$. Using the Kuratowski--Zorn Lemma, find two maximal independent sets $M,M'$ in $X$ such that $\{u,o,w\}\subseteq M$ and $\{o',o',w'\}\subseteq M'$. By Corollary~\ref{c:Max=dim}, $|M|=\|X\|=|M'|$, so there exists a bijection $\varphi:M\to M'$ such that $\varphi:uow\mapsto o'o'w'$. If $\|X\|\ne 3$, then the affine space $X$ is Desarguesian, by Theorem~\ref{t:proaffine-Desarguesian}. In this case, we can apply Theorem~\ref{t:extension-iso} and extend the bijection $\varphi$ to an automorphism $\Phi$ of the affine space $X$. Then $\Phi$ is a required automorphism of $X$ with $\Phi uow=u'o'w'$. 

Now assume that $\|X\|=3$ and hence $X$ is a Playfair plane. By Theorem~\ref{c:affine-spread-completion}, the Playfair plane $X$ is completely regular, and by Theorem~\ref{t:compreg-Moufang<=>}, its spread completion $\overline X$ is a Moufang projective plane. Since $X$ is not fano, we can apply Theorem~\ref{t:Fano<=>} and conclude that the projective liner $\overline X$ is not Fano. Let $\partial X\defeq\overline X\setminus X$ be the horizon of the affine space $X$ in its spread completion. Consider the directions $h\defeq(\Aline ou)_\parallel$,  $v\defeq(\Aline ow)_\parallel$, $h'\defeq(\Aline ou')_\parallel$,  $v'\defeq(\Aline ow')_\parallel$ in the horison $\partial X$ of the affine space $X$, and observe that $r\defeq\{oo,uh,wv\}$ and $r'\defeq\{o'o',u'h',w'v'\}$ are two projective repers in the projective plane $\overline X$. By Theorem~\ref{t:Moufang-reper-hom}, there exists an automorphism $\Psi$ of the projective plane  $X$ such that $\Psi:ouhwv\mapsto o'u'h'w'v'$. Then $\Phi[\partial X]=\Phi[\Aline hv]=\Aline{h'}{v'}=\partial X$ and the restriction $\Phi\defeq\Psi{\restriction}_X$ is a required automorphism of the affine space $X$ with $\Phi uow=u'o'w'$.
\end{proof}

\begin{corollary}\label{c:ternars-in-non-Fano-Moufang-afplane-are-isomorphic} Any two ternars of a non-Fano Moufang affine space are isomorphic.
\end{corollary} 

\begin{proof} Let $R,R'$ be two ternars of a non-Fano Moufang affine space $X$. Then there exist two planes $\Pi,\Pi'$ in $X$ and affine bases $uow$ and $u'o'w'$ in the affine planes $\Pi,\Pi'$, such that $R,R'$ are isomorphic to the ternars of the based affine planes $(\Pi,uow)$, $(\Pi',u'o'w')$, respectively. By Corollary~\ref{c:Moufang=>3hom}, there exists an automorphism $\Phi$ of the affine space $X$ such that $\Phi uow=u'o'w'$. Then $\Phi{\restriction}_\Pi$ is an automorphism of the based affine planes $(\Pi,uow)$ and $(\Pi',u'o'w')$. By Theorem~\ref{t:tring-iso<=>}, the ternars of these based affine planes are isomorphic and hence  the ternars $R$ and $R'$ are isomorphic, too.
\end{proof}

\begin{remark} Since we have found a proof of Schafer's Theorem~\ref{t:Schafer} only for non-Boolean alternative division rings, we formulated (and proved)  Theorems~\ref{t:Moufang=>4-homogen}, \ref{t:Moufang-reper-hom} and Corollaries~\ref{c:Moufang=>4hom}, 
\ref{c:ternars-in-non-Fano-Moufang-plane-are-isomorphic}, \ref{c:Moufang=>3hom}, 
\ref{c:ternars-in-non-Fano-Moufang-afplane-are-isomorphic} only for non-Fano Moufang liners. If Schafer's Theorem~\ref{t:Schafer} is true also for Boolean alternative division rings, then the mentioned theorems and corollaries are true for all Moufang  liners (including Fano). 
\end{remark}

\begin{remark} Corollaries~\ref{c:Moufang=>4hom} and \ref{c:Moufang=>3hom} cannot be reversed: in Corollary 2.3 of \cite{KS1973}, Kegel and Schleiermacher  
constructed an infinite $4$-homogeneous projective plane $P$ which is not Moufang (more precisely, every quadrangle $abcd$ in $P$ generates a closed subliner of $P$, isomorphic to the free projectivization of the subliner $\{a,b,c,d\}$ of $P$). Removing any line from the Kegel--Schleiermacher projective plane $P$, we obtain a $3$-homogeneous affine plane, which is not Moufang.
\end{remark}



On the other hand, for finite projective planes,  \index[person]{Ostrom}Ostrom\footnote{{\bf Teodore G. Ostrom} (1916 -- 2011) was an American mathematician, who made fundamental contributions to the field of finite geometries, in particular, the theory of translation planes.}  and \index[person]{Wagner}Wagner\footnote{{\bf Ascher Otto Wagner} (1930 -- 2000) was an Austrian and British mathematician, specializing in the theory of finite groups and finite projective planes. He is known for the Dembowski–Wagner theorem. Ascher Wagner received his Ph.D. in 1958 with dissertation ``Some Problems on Projective Planes and Related Topics in the Theory of Algebraic Operations'' supervised by Kurt Hirsch. Wagner was a faculty member at the University of London and then at the University of Birmingham.} \cite{OW1959} proved the following homogeneity result. 

\begin{theorems}[Ostrom, Wagner, 1959]\label{t:Ostrom-Wagner1959}  Every $2$-homogeneous finite projective plane is Desarguesian.
\end{theorems}

An affine counterpart of Theorem~\ref{t:Ostrom-Wagner1959} was proved by \index[person]{Jha}Jha\footnote{{\bf Vikram Jha}, a professor of mathematics in Caledonian University, Glasgow, Scotland. Jha defended his Ph.D. Thesis ``On the Automorphism Groups of Quasifields'' in 1971 in University of London, under supervision of Daniel R. Hughes.} and \index[person]{Johnson}Jonson\footnote{
{\bf Norman Lloyd Johnson}, currently works at University of Iowa. Johnson defended his Ph. D. Thesis ``A Classification of Semi-Translation Planes'' in 1968 in Washington State University 1968, under supervision of Theodore G. Ostrom.}
  \cite{JhaJohnson1998} in 1998.

\begin{theorems}[Jha, Johnson, 1998]\label{t:Jha-Johnson1998} Every $2$-homogeneous finite affine plane $X$ is Thalesian. If $|X|_2\notin\{9,27\}$, then the affine plane $X$ is Desarguesian.
\end{theorems}

\begin{remark} The exceptional non-Desarguesian $2$-transitive affine planes of orders $9$ and 27 are the Thalesian plane $X=\mathbb J_9\times\mathbb J_9$, described in Section~\ref{s:J9xJ9}, and the \index[person]{Hering}\index{Hering plane}\defterm{Herring plane}, discovered by Christoph Hering\footnote{{\bf Christoph Hering} (born in 1939, in Oedelsheim) is a German mathematician who works in group theory and its applications in geometry and algebraic combinatorics. He is a professor at the University of T\"ubingen. 
Hering received his Ph.D. in 1962 under Reinhold Baer and Ruth Moufang at the Johann Wolfgang Goethe University in Frankfurt am Main with the dissertation entitled ``A Characterization of the Finite Two-Dimensional Projective Groups". In 1963/64 he was an assistant in Frankfurt, completed his habilitation in 1966 at the University of Mainz, and spent time in Chicago before becoming a professor in T\"ubingen in 1974. In the 1970s, he served as Dean of the Department of Mathematics. In 1965, he provided a classification of M\"obius planes based on their automorphism groups, in a manner similar to the earlier Lenz–Barlotti classification of projective planes.} \cite{Hering1970} in 1970. Computer calculations show that those two exceptional planes are not $3$-transitive, so, we have the following homogeneity result for finite affine planes.
\end{remark}

\begin{theorem} Every $3$-homogeneous finite affine plane $X$ is Desarguesian.
\end{theorem}

\begin{remark} All $2$-homogeneous $2$-balanced finite liners were classified by William Kantor in \cite{Kantor1985}: {\em every $2$-homogeneous $2$-balanced finite liner $X$ is  one of the following: \textup{(1)} a Pappian affine space, \textup{(2)} a Pappian projective space, \textup{(3)} a Hermitian unital, \textup{(4)} a Ree unital, \textup{(5)} the Thalesian plane of order $9$, \textup{(6)} the Hering plane of order $27$, \textup{(7)} one of two exceptional liners with $|X|=3^6$ and $|X|_2=9$.}
\end{remark}

\section{Characterizing Moufang planes via permutation groups}

In this section we apply the theory of permutation groups, developed in Chapter~\ref{ch:permutation} to obtain characterizations of Moufang planes in terms of groups of permutations of lines in such planes.

We recall that a permutation group $G\subseteq\Sym(X)$ is \defterm{of affine type} if it contains a sharply transitive normal Abelian subgroup $N$. In this case $G$ is isomorphic to a subgroup of the holomorph $\Hol(N)$ of $N$, according to Theorem~\ref{t:affine-type}. 

For two lines $L,\Lambda$ in a projective plane $X$ and any point $c\in X\setminus(L\cup\Lambda)$, let $c_{\Lambda,L}:L\to\Lambda$ be the central projection assigning to each point $x\in L$ a unique point $y\in\Lambda$ such that $\Aline cx=\Aline cy$.

This notion can be dualized as follows. For any line $L$ in a projective plane $X$ and two points $a,c\in X\setminus L$, let $L_{c,a}:\mathcal L_a\to\mathcal L_c$ be the bijection assigning to each line $A\in\mathcal L_a$  a unique line $C\in \mathcal L_c$ such that $A\cap L=L\cap C$. Here for a point $x\in X$ we denote by $\mathcal L_x\defeq\{L\in\mathcal L_X:x\in L\}$ the family of lines in $X$ that contain the point $x$.

\begin{theorem}\label{t:Moufang<=>permutation} For a projective plane $(X,\mathcal L)$, the following conditions are equivalent.
\begin{enumerate}
\item The projective plane $(X,\mathcal L)$ is Moufang.
\item For arbitrary lines $D,H\in\mathcal L_X$ and distinct points $h,v\in H\setminus D$, the set $\{h_{D,L}v_{L,D}{\restriction}_{D\setminus H}:L\in\mathcal L_X\setminus(\mathcal L_h\cup\mathcal L_v)\}$ is contained in some permutation subgroup $G\subseteq\Sym(D\setminus H)$ of affine type.
\item There exist distinct lines $H_1,H_2$ in $X$ such that for every $i\in\{1,2\}$, there exists a line $D\in\mathcal L$ and distinct points $h,v\in H_i\setminus D$ such that the set $\{h_{D,L}v_{L,D}{\restriction}_{D\setminus H_i}:L\in\mathcal L_X\setminus(\mathcal L_h\cup\mathcal L_v)\}$ is contained is some permutation subgroup $G\subseteq\Sym(D\setminus H_i)$ of affine type.
\smallskip
\item The dual projective plane $(\mathcal L,\{\mathcal L_x:x\in X\})$ is Moufang.
\item For arbitrary points $d,h\in X$ and distinct lines $H,V\in \mathcal L_h\setminus\mathcal L_d$, the set $\{H_{d,c}V_{c,d}{\restriction}_{\mathcal L_d\setminus \mathcal L_h}:c\in X\setminus(H\cup V)\}$ is contained in some permutation subgroup $G\subseteq\Sym(\mathcal L_d\setminus \mathcal L_h)$ of affine type.
\item There exist distinct points $h_1,h_2$ in $X$ such that for every $i\in\{1,2\}$, there exists a point $d\in X$ and distinct lines $H,V\in \mathcal L_{h_i}\setminus \mathcal L_d$ such that the set $\{H_{d,c}V_{c,d}{\restriction}_{\mathcal L_d\setminus \mathcal L_{h_i}}:c\in X\setminus(H\cup V)\}$ is contained is some permutation subgroup $G\subseteq\Sym(\mathcal L_d\setminus\mathcal L_{h_i})$ of affine type.
\smallskip
\item For arbitrary points $d,h\in X$ and distinct lines $H,V\in \mathcal L_h\setminus \mathcal L_d$, the set $\{c_{H,V}d_{V,H}{\restriction}_{H\setminus \{h\}}:c\in X\setminus (H\cup V)\}$ is contained is some permutation subgroup $G\subseteq\Sym(H\setminus\{h\})$ of affine type.
\item There exist distinct points $h_1,h_2$ in $X$ such that for every $i\in\{1,2\}$, there exist distinct lines $H,V$ containing $h_i$ and a point $d\in X\setminus(H\cup V)$ such that the set $\{c_{H,V}d_{V,H}{\restriction}_{H\setminus \{h_i\}}:c\in X\setminus(H\cup V)\}$ is contained is some permutation subgroup $G\subseteq\Sym(H\setminus \{h_i\})$ of affine type.
\end{enumerate}
\end{theorem} 

\begin{proof} $(1)\Ra(2)$ Assume that a projective plane $X$ is Moufang.
Fix any lines $D,H\in\mathcal L_X$ and distinct points $h,v\in H\setminus D$.
If $|D\setminus H|\le 3$, then the permutation group $\Sym(D\setminus H)$ is sharply transitive (by Exercise~\ref{ex:S4-of affine type}) and hence  the set  $\{h_{D,L}v_{L,D}{\restriction}_{D\setminus H}:L\in\mathcal L_X\setminus(\mathcal L_h\cup\mathcal L_v)\}$ is contained is the  permutation subgroup $G=\Sym(D\setminus H)$ of affine type. So, assume that  $|D\setminus H|>3$, which implies $|D|\ge 5$. Since $|H|\ge 3$, the projective plane $X$ coincides with its maximal $3$-long flat containing the set $H\cup D$ and hence the projective plane $X$ is $3$-long and $2$-balanced, by Corollary~\ref{c:Avogadro-projective}. Since $|D|\ge 5$, the $2$-balanced projective plane $X$ is $5$-long. Then the subliner $\Pi\defeq X\setminus H$ is a $4$-long Playfair plane. Since the projective plane $X$ is Moufang, the affine plane $\Pi=X\setminus H$ is Thalesian, by Theorem~\ref{t:Moufang<=>everywhere-Thalesian}. Consider the line $D'\defeq D\setminus H$ in the affine plane $\Pi$, and identify the points $h,v\in H$ of the projective plane $X$ with the directions $\boldsymbol h\defeq\{\Lambda\setminus\{h\}:\Lambda\in \mathcal L_h\setminus\{H\}\}$ and $\boldsymbol v\defeq\{\Lambda\setminus\{v\}:\Lambda\in \mathcal L_v\setminus\{H\}\}$ on the Playfair plane $\Pi$. In this case, the set $\{h_{D,L}v_{L,D}{\restriction}_{D\setminus H}:L\in\mathcal L_X\setminus(\mathcal L_h\cup\mathcal L_v)\}$ coincides with the set 
$\{\boldsymbol h_{D',L}\boldsymbol v_{L,D'}:L\in\mathcal L_\Pi\setminus(\boldsymbol h\cup\boldsymbol v)\}=\Sym^{\Join}[D';\boldsymbol h,\boldsymbol v]$, which is contained in the subgroup $G=\Sym^{\Join}_\Pi(D')\subseteq \Sym(D')$ that has affine type, by Theorem~\ref{t:Schleiermacher}.
\smallskip

The implication $(2)\Ra(3)$ is trivial.
\smallskip

$(3)\Ra(1)$ Assume that $H_1$ and $H_2$ are two distinct lines in $X$ satisfying the condition (3). By Theorem~\ref{t:proj-Moufang<=>}, to prove that the projective plane is Moufang, it suffices to show that for every $i\in\{1,2\}$ the affine subliner $\Pi_i=X\setminus H_i$ is Thalesian. This is trivially true if $|\Pi_i|_2=2$. So, we assume that $|\Pi_i|_2\ge 3$. By the condition (3), there exists a line $D\in\mathcal L$ and distinct points $h,v\in H_i\setminus D$ such that the set $\{h_{D,L}v_{L,D}{\restriction}_{D\setminus H_i}:L\in\mathcal L_X\setminus(\mathcal L_h\cup\mathcal L_v)\}$ is contained is some permutation subgroup $G\subseteq\Sym(D\setminus H_i)$ of affine type.
Consider the line $D'\defeq D\setminus H_i$ and the directions $\boldsymbol h\defeq\{L\setminus\{h\}:L\in\mathcal L_h\setminus\{H_i\}\}$ and $\boldsymbol v\defeq\{L\setminus\{H_i\}:L\in\mathcal L_v\setminus\{H\}\}$ on the Playfair plane $\Pi_i$. Observe that the set $\{h_{D,L}v_{L,D}{\restriction}_{D\setminus H_i}:L\in\mathcal L_X\setminus(\mathcal L_h\cup\mathcal L_v)\}$ coincides with the set $\{\boldsymbol h_{D',L}\boldsymbol v_{L,D'}:L\in\mathcal L_{\Pi_i}\setminus(\boldsymbol h\cup\boldsymbol v)\}$, which implies that the latter set is contained in the subgroup $G\subseteq \Sym(D\setminus H_i)=\Sym(D')$ of affine type. By Theorem~\ref{t:Schleiermacher}, the Playfair plane $\Pi_i=X\setminus H_i$ is Thalesian. Since the affine planes $X\setminus H_1$ and $X\setminus H_2$ are Thalesian, the projective plane $X$ is Moufang, by Theorem~\ref{t:proj-Moufang<=>}.
\smallskip

The equivalence $(1)\Leftrightarrow(4)$ has been proved in Theorem~\ref{t:proj-Moufang<=>}.
\smallskip

The equivalences $(4)\Leftrightarrow(5)\Leftrightarrow(6)$ follows from the equivalences $(1)\Leftrightarrow(2)\Leftrightarrow(3)$ by duality.
The equivalences $(5)\Leftrightarrow(7)$ and $(6)\Leftrightarrow(8)$ can be easily deduced from the following Lemma~\ref{l:dual-line-pencil}.
\end{proof}

\begin{lemma}\label{l:dual-line-pencil} Let $(X,\mathcal L)$ be a projective plane, $d,h\in X$ be two points and $H,V\in\mathcal L_h\setminus\mathcal L_d$ be two distinct lines. Consider the bijection $\beta:H\setminus\{h\}\mapsto \mathcal L_d\setminus\mathcal L_h$, $\beta:x\mapsto \Aline xd$, inducing the group isomorphism $I: \Sym(H\setminus\{h\})\to\Sym(\mathcal L_d\setminus\mathcal L_h)$, $I:f\mapsto \beta f\beta^{-1}$. For every point $c\in X\setminus(H\cup V)$, $$I(c_{H,V}d_{V,H}{\restriction}_{H\setminus\{h\}})=H_{d,c}V_{c,d}{\restriction}_{\mathcal L_d\setminus\mathcal L_h}.$$
\end{lemma}

\begin{proof} Given any line $L\in\mathcal L_d\setminus\mathcal L_h$, 
find unique points $x\in L\cap H$, $v\in V\cap L$, $y\in H\cap \Aline cv$, and observe that $H_{d,c}V_{c,d}(L)=H_{d,c}(\Aline cv)=\Aline dy$.

\begin{picture}(200,100)(-190,-15)

\put(0,0){\line(1,4){15}}
\put(13,65){$V$}
\put(0,0){\line(-1,4){15}}
\put(-20,65){$H$}
\put(10,40){\line(-3,-1){60}}
\put(10,40){\line(-3,1){22}}
\put(10,40){\line(3,-1){54}}
\put(33,35){$L$}
\put(-8.3,33.9){\line(6,-1){73}}

\put(0,0){\circle*{3}}
\put(-2,-10){$h$}
\put(-8.3,33.9){\circle*{3}}
\put(-14,26){$y$}
\put(64,22){\circle*{3}}
\put(67,19){$d$}
\put(-50,20){\circle*{3}}
\put(-58,17){$c$}
\put(10,40){\circle*{3}}
\put(13,41){$v$}
\put(-12,47){\circle*{3}}
\put(-21,46){$x$}
\end{picture}

On the other hand, 
$I(c_{H,V}d_{V,H}{\restriction}_{H\setminus\{h\}})(L)=
\beta c_{H,V}d_{V,H}\beta^{-1}(L)=\beta c_{H,V}d_{V,H}(x)=\beta c_{H,V}(v)=\beta(y)=\Aline yd=H_{d,c}V_{c,d}(L),
$  witnessing that 
 $I(c_{H,V}d_{V,H}{\restriction}_{H\setminus\{h\}})=H_{d,c}V_{c,d}{\restriction}_{\mathcal L_d\setminus\mathcal L_h}$.
\end{proof}

Theorem~\ref{t:Moufang<=>permutation} implies the following characterization of Moufang Playfair planes.

\begin{corollary}\label{c:Moufang<=>permutation} For a Playfair plane $X$ and its spread completion $\overline X$, the following conditions are equivalent:
\begin{enumerate}
\item The Playfair plane $X$ is Moufang.
\item For arbitrary disjoint lines $H,V$ in $X$ and point $d\in \overline X\setminus(\overline H\cup\overline V)$, the set $\{c_{H,V}d_{V,H}:c\in \overline X\setminus (\overline H\cup \overline V)\}$ is contained in a  subgroup $G\subseteq\Sym(H)$ of affine type.
\item There exist distinct directions $\boldsymbol h_1,\boldsymbol h_2\in\partial X$ such that for every $i\in\{1,2\}$, there exist distinct lines $H,V\in\boldsymbol h_i$ and a point $d\in \overline X\setminus(\overline H\cup \overline V)$ such that the set $\{c_{H,V}d_{V,H}:c\in \overline X\setminus(\overline H\cup \overline V)\}\}$ is contained in a subgroup $G\subseteq\Sym(H)$ of affine type.
\end{enumerate}
\end{corollary} 

\begin{proof} $(1)\Ra(2)$ If the Playfair plane $X$ is Moufang, then its spread completion $\overline X$ is Moufang, by Theorem~\ref{t:affine-Moufang<=>}. Given any disjoint lines $H,V$ in $X$, consider their (common) direction $\boldsymbol h\defeq H_\parallel=V_\parallel\in\partial X\subseteq \overline X$. By Theorem~\ref{t:Moufang<=>permutation}, for every point $d\in\overline{X}\setminus(\overline H\cup\overline V)$, the set 
$$\{c_{\overline H,\overline V}d_{\overline V,\overline H}{\restriction}_{\overline H\setminus\{\boldsymbol h\}}:c\in \overline X\setminus(\overline H\cup\overline V)\}=\{c_{H,V}d_{V,H}:c\in \overline X\setminus(\overline H\cup\overline V)\}$$is contained in a subgroup $G\subseteq\Sym(\overline H\setminus\{\boldsymbol h\})=\Sym(H)$ of affine type.
\smallskip

The implication $(2)\Ra(3)$ is trivial.
\smallskip

$(3)\Ra(1)$ Assume that there exist distinct directions $\boldsymbol h_1,\boldsymbol h_2\in\partial X$ such that for every $i\in\{1,2\}$, there exist distinct lines $H,V\in\boldsymbol h_i$ and a point $d\in \overline X\setminus(\overline H\cup \overline V)$ such that the set $\{c_{H,V}d_{V,H}:c\in \overline X\setminus(\overline H\cup \overline V)\}\}$ is contained in a subgroup $G\subseteq\Sym(H)$ of affine type. Then the set $$\{c_{\overline H,\overline V}d_{\overline V,\overline H}{\restriction}_{\overline H\setminus\{\boldsymbol h\}}:c\in \overline X\setminus(\overline H\cup\overline V)\}=\{c_{H,V}d_{V,H}:c\in \overline X\setminus(\overline H\cup\overline V)\}$$
is contained in the subgroup $G\subseteq\Sym(H)=\Sym(\overline H\setminus\{\boldsymbol h_i\})$ of affine type. By Theorem~\ref{t:Moufang<=>permutation}, the projective plane $\overline X$ is Moufang, and by Theorem~\ref{t:affine-Moufang<=>}, the Playfair plane $X$ is Moufang.
\end{proof}



\chapter{Additive and multiplicative properties of projective spaces}\label{ch:AMP}

In this chapter we present geometric characterizations of some algebraic properties of ternars of projective spaces. Let us recall that a {\em projective space} if a $3$-long projective liner of rank $\ge 3$. A projective space of rank $3$ is called a {\em projective plane}. A ternar $R$ is called a {\em ternar of a projective space} $X$ if there exists a plane $\Pi\subseteq X$ and a projective base $uowe$ in $\Pi$ whose ternar is isomorphic to the ternar $R$.



\section{Invertible-add projective spaces}

\begin{definition}\label{d:p-invertible-plus} A projective space $X$ is \index{invertible-plus projective space}\index{projective space!invertible-plus}\defterm{invertible-plus} if for every distinct points $h,v\in X$, $d\in\Aline hv$, and $a,b,c,o,a',b',c'\in X\setminus\Aline hv$
$$\big((v\in \Aline b{c}\cap\Aline oa\cap\Aline o{a'}\cap\Aline {b'}{c'})\wedge(h\in \Aline ab\cap\Aline oc\cap\Aline o{c'}\cap\Aline{a'}{b'})\wedge(d\in \Aline ob\cap\Aline o{b'}\cap\Aline a{c'})\big)\Ra(d\in \Aline {a'}c).$$

\begin{picture}(120,155)(-170,-55)
\linethickness{=0.6pt}
\put(90,0){\color{teal}\line(-1,0){120}}
\put(0,90){\color{cyan}\line(0,-1){120}}
\put(0,90){\line(1,-1){90}}
\put(0,90){\color{cyan}\line(-1,-3){45}}
\put(90,0){\color{teal}\line(-3,-1){135}}
\put(45,45){\color{red}\line(-1,-1){90}}
\put(45,45){\color{red}\line(-3,-5){45}}
\put(45,45){\color{red}\line(-5,-3){75}}
\put(0,90){\color{cyan}\line(1,-5){18}}
\put(90,0){\color{teal}\line(-5,1){90}}

\put(0,0){\circle*{3}}
\put(-6,1){$o$}
\put(18,0){\circle*{3}}
\put(17,-7){$c$}
\put(0,18){\circle*{3}}
\put(-7,18){$a$}
\put(15,15){\circle*{3}}
\put(14.5,18.6){$b$}
\put(90,0){\color{teal}\circle*{3}}
\put(93,-3){\color{teal}$h$}
\put(0,90){\color{cyan}\circle*{3}}
\put(-2,93){\color{cyan}$v$}
\put(45,45){\color{red}\circle*{3}}
\put(47,46){\color{red}$d$}
\put(-30,0){\circle*{3}}
\put(-37,-1){$c'$}
\put(0,-30){\circle*{3}}
\put(-2,-38){$a'$}
\put(-45,-45){\circle*{3}}
\put(-48,-55){$b'$}
\end{picture}
\end{definition}

Definitions~\ref{d:p-invertible-plus}, \ref{d:a++inversive} and Theorems~\ref{t:invertible-plus<=>}, \ref{t:invertible-plus<=>power-associative} imply the following characterization of invertible-plus projective spaces. 

\begin{proposition}\label{p:p-invertible-plus<=>} For every projective space $X$, the following conditions are equivalent:
\begin{enumerate}
\item $X$ is invertible-plus;
\item for every hyperplane $H\subseteq X$, the affine liner $X\setminus H$ is invertible-plus;
\item for every ternar $R$ of $X$, its plus loop $(R,+)$ is invertible;
\item for every  ternar $R$ of $X$, its plus loop $(R,+)$  is power-associative.
\end{enumerate}
\end{proposition}

\begin{remark} In literature, invertibile-plus projective planes are referred to as projective planes satisfying the {\em hexagon condition}, see \cite[p.~54]{Pickert}.
\end{remark}

\begin{definition}\label{d:p-invertible-puls} A projective space $X$ is \defterm{invertible-puls} if for any  lines $L,L'\subseteq X$ and distinct points $a,b,c\in L\setminus L'$, $a',b',c'\in L'\setminus L$, $h\in L\cap L'$, $v\in \Aline a{a'}\cap\Aline b{b'}\cap\Aline c{c'}$, if $\Aline a{b'}\cap\Aline b{c'}\subseteq \Aline vh$, then $\Aline {a'}b\cap\Aline {b'}c\subseteq \Aline vh$.

\begin{picture}(250,135)(-150,-15)
\linethickness{=0.6pt}
\put(-30,0){\color{teal}\line(1,0){150}}
\put(0,0){\color{blue}\line(0,1){120}}
\put(0,0){\color{cyan}\line(1,1){60}}
\put(0,120){\line(1,-1){120}}
\put(0,120){\color{blue}\line(1,-6){20}}
\put(120,0){\color{teal}\line(-6,1){120}}
\put(60,60){\color{cyan}\line(-3,-2){90}}
\put(60,60){\color{cyan}\line(-2,-3){40}}
\put(-30,0){\color{red}\line(11,4){110}}
\put(0,0){\color{red}\line(2,1){80}}

\put(17.1,17.1){\circle*{3}}
\put(17,21){$b$}
\put(-30,0){\circle*{3}}
\put(-33,-10){$a'$}
\put(0,0){\circle*{3}}
\put(-3,-10){$b'$}
\put(20,0){\circle*{3}}
\put(17,-10){$c'$}
\put(120,0){\color{teal}\circle*{3}}
\put(123,-3){\color{teal}$h$}
\put(0,120){\color{blue}\circle*{3}}
\put(0,20){\circle*{3}}
\put(-7,21){$a$}
\put(60,60){\color{cyan}\circle*{3}}
\put(62,61){\color{cyan}$v$}
\put(80,40){\color{red}\circle*{3}}
\put(30,15){\circle*{3}}
\put(29,8){$c$}
\put(60,-10){\color{teal}$L'$}
\put(65,11){\color{teal}$L$}
\end{picture}
\end{definition}

Definitions~\ref{d:p-invertible-puls}, \ref{d:a-invertible-puls} and Theorems~\ref{t:invertible-puls<=>}, \ref{t:invertible-puls<=>power-associative} imply the following characterization of invertible-puls projective spaces. 

\begin{proposition}\label{p:p-invertible-puls<=>} For every projective space $X$, the following conditions are equivalent:
\begin{enumerate}
\item $X$ is invertible-puls;
\item for every hyperplane $H\subseteq X$, the affine liner $X\setminus H$ is invertible-puls;
\item for every ternar $R$ of $X$, the puls loop $(R,\!\puls\!)$ is invertible;
\item for every  ternar $R$ of $X$, the puls loop $(R,\!\puls\!)$  is power-associative.
\end{enumerate}
\end{proposition}

\begin{definition} A projective space $X$ is called \defterm{invertible-add} if $X$ is invertible-plus and invertible-puls.
\end{definition}

\begin{theorem}\label{t:p-invertible-add<=>} For a projective space $X$ the following conditions are equivalent:
\begin{enumerate}
\item $X$ is invertible-add;
\item $X$ is invertible-plus;
\item $X$ is invertible-puls.
\end{enumerate}
\end{theorem}

\begin{proof} It suffices to check that a projective space $X$ is invertible-plus if and only if it is invertible-puls. This equivalence follows from the equivalence of the configurations determining the invertibility-plus and invertibility-puls:


\end{proof}

Theorems~\ref{t:p-invertible-add<=>} and \ref{t:plus-puls-duality} imply that the invertibility-add is a self-dual property of projective planes.

\begin{corollary}\label{c:inv-add-self-dual} A projective plane is invertible-add if and only if its dual projective plane is invertible-add.
\end{corollary}

\begin{proposition}\label{p:invertible-add=>quadratic} If a projective plane $\Pi$ is invertible-add, then for every line $\Lambda$ in $X$, the affine plane $\Pi\setminus\Lambda$ is invertible-add, invertible-half, by-Thalesian and quadratic.
\end{proposition}

\begin{proof} The invertibility-add of the affine plane $\Pi\setminus\Lambda$ follows from Propositions~\ref{p:p-invertible-plus<=>} and \ref{p:p-invertible-puls<=>}. To prove that $\Pi\setminus \Lambda$ is invertible-half, fix any distinct parallel lines $L,L'$ in $\Pi\setminus\Lambda$ and distinct points $v,b\in L$, $a',h,c'\in L'$, $x\in \Aline {a'}b\cap\Aline vh$, $\gamma\in \Aline v{c'}\cap\Aline bh$ such that $\Aline v{a'}\parallel \Aline bh$ and $\Aline vh\parallel \Aline b{c'}$ in $\Pi\setminus\Lambda$. We have to prove that $\Aline x\gamma\parallel L$ in the affine plane $\Pi\setminus\Lambda$. Let $a,y,b'\in\Lambda$ be unique points such that $a\in \Aline v{a'}\cap\Aline bh$, $y\in\Aline vh\cap\Aline b{c'}$ and $b'\in \overline L\cap \overline{L'}$ in the projective plane $\Pi$. Consider the unique point $c\in \Aline x{b'}\cap \Aline bh$. The invertibility-puls of the projective plane $\Pi$ implies that $v\in\Aline c{c'}$ and hence $c\in\Aline v{c'}\cap\Aline bh=\{\gamma\}$. Then $b'\in \Aline xc\cap\overline L=\Aline x\gamma\cap\overline L$, which means that the line $\Aline x\gamma$ is parallel to the line $L$ in the affine plane $\Pi\setminus\Lambda$.

\begin{picture}(180,130)(-130,-15)
\put(0,0){\vector(1,0){120}}
\put(100,-10){$L$}
\put(0,30){\vector(1,0){120}}
\put(100,33){$L'$}
\put(15,15){\vector(1,0){105}}
\put(130,12){$b'$}
\put(0,0){\vector(1,1){90}}
\put(30,0){\vector(1,1){75}}
\put(105,90){$y$}
\put(0,30){\line(1,-1){30}}
\put(0,0){\vector(0,1){100}}
\put(30,0){\vector(0,1){100}}
\put(0,0){\line(2,1){60}}
\put(13,105){$a$}

\put(0,0){\circle*{3}}
\put(-8,-8){$v$}
\put(30,0){\circle*{3}}
\put(27,-9){$b$}
\put(30,15){\circle*{3}}
\put(31,9){$\gamma$}
\put(31,19){$c$}
\put(0,30){\circle*{3}}
\put(-9,28){$a'$}
\put(30,30){\circle*{3}}
\put(22,32){$h$}
\put(60,30){\circle*{3}}
\put(57,33){$c'$}
\put(15,15){\circle*{3}}
\put(6,12){$x$}
\end{picture}  

Therefore, the affine plane $\Pi\setminus\Lambda$ is invertible-half. By  Theorem~\ref{t:by-Thalesian<=>by-Boolean<=>quadratic}, $\Pi\setminus\Lambda$ is by-Thalesian and quadratic.
\end{proof}

\begin{lemma}\label{l:V1/2} Let $X$ be an invertible-add projective plane, $V_0,V_1,V_\infty$ be distinct lines in $X$, $v\in V_0\cap V_1\cap V_\infty$ be a point, and  $\{z_i^n:i\in\{0,1,\infty\},\;n\in\IZ\}$ be an indexed set of points in $X\setminus\{v\}$ such that
\begin{enumerate}
\item $z_i^n\in V_i$ for all $i\in\{0,1,\infty\}$ and $n\in\IZ$;
\item $z^n_\infty\in \Aline {z^k_0}{z^{k+n}_1}$ for all $k,n\in\IZ$;
\item $z^n_\infty\ne z^m_\infty$ for any distinct integer numbers $n,m$.
\end{enumerate} 
Then there exists a line $V_{1/2}$ in $X$ containing the point $v$ and an indexed family of points $\{z_{1/2}^n:n\in\IZ\}\subseteq V_{1/2}\setminus(V_0\cup V_1\cup V_\infty)$ such that 
 $z^{k+n}_{1/2}\in \Aline {z^k_0}{z^{n}_1}$ for all $k,n\in\IZ$. 
\end{lemma}

\begin{proof} Consider the affine plane $\Pi\defeq X\setminus V_\infty$ and the affine base $uow\defeq z^0_1z_0^0z_0^1$ in $\Pi$. The condition (2) ensures that $h\defeq z_\infty^0$ is the horizontal direction for the affine plane $\Pi$ and $z_1^1$ is the diunit of the affine base $uow$. By Theorem~\ref{t:invertible-add<=>Z-elementary}, the ternar $\Delta=\Aline oe$ of the based affine plane $(\Pi,uow)$ is $\IZ$-elementary and hence there exists a function $f:\IZ\to \Delta$ such that $f(1)=e$ and $f(x\cdot y+z)=f(x)_\times f(y)_+f(z)$ for all $x,y,z\in\IZ$.

Let $z_{1/2}^1$ be the unique common point of the lines $\Aline oe=\Aline {z_0^0}{z_1^1}$ and $\Aline wu=\Aline {z_0^1}{z_1^0}$. Since the points $z_\infty^{1}\in \Aline {z_0^0}{z_1^1}=\Aline oe$ and $z_\infty^{-1}\in\Aline {z_0^1}{z_1^0}=\Aline wu$ are distinct, the point $z_{1/2}^1$ does not belong to the horizon line $V_\infty$ of the affine plane $\Pi$ and hence $z^1_{1/2}\in \Pi$. Then the line $V_{1/2}\defeq\Aline v{z^1_{1/2}}$ is not equal to the line $V_\infty$. Since the points $u,o,w,e$ are distinct, the line $V_{1/2}$ also is not equal to the lines $V_0$ and $V_1$. By Proposition~\ref{p:invertible-add=>quadratic}, the affine plane $\Pi$ is quadratic, which implies $z_{1/2}^1+z_{1/2}^1=e=z_1^1$, according to Theorem~\ref{t:by-Thalesian<=>by-Boolean<=>quadratic}. Consider the unique affine base $u'ow'$ in $\Pi$ such that $u'\in\Aline ou$, $w'\in\Aline ow$ and the diunit $e'$ of the base $u'ow'$ coincides with the point $z_{1/2}^1$. By Theorem~\ref{t:invertible-add<=>Z-elementary}, the ternar $\Aline o{e'}=\Delta$ of the based affine plane $(\Pi,u'o'w)$ is $\IZ$-elementary and hence there exists a function $f':\IZ\to \Delta$ such that $f'(1)=e'$ and $f'(x\cdot y+z)=f'(x)_\times f'(y)_+f'(z)$ for all $x,y,z\in\IZ$. The equality $e'+e'=z_{1/2}^1+z_{1/2}^1=e$ implies $f'(2)=e'+e'=e$ and hence $f'(2n)=f(n)$ for all $n\in\IZ$. 

For every $n\in\IZ$, consider the unique point $z_{1/2}^n\in V_{1/2}\cap \Aline {f(n)}h$.  It remains to show that $z_{1/2}^{k+n}=\Aline{z_0^k}{z^n_1}$ for all $k,n\in\IZ$. The condition (2) implies that $z^n_i\in \Aline{f(n)}h=\Aline{f'(2n)}h$ for all $n\in \IZ$ and $i\in\{0,1\}$.  Given any integer numbers $k,n\in\IZ$, consider the points $f'(n-k),f'(n+k)\in \Delta$ and the unique point $\zeta\in V_0\cap\Aline{f'(n+k)}h$. Observe that $f'(n+k)=f'(2k)+f'(n-k)=f'(2k)\puls f'(n-k)$ and $f'(2n)=f'(n-k)+f'(n+k)=f'(n-k)\puls f'(n+k)$, which implies $\Aline {z_0^k}{z_{1/2}^{k+n}}\parallel \Aline o{z_{1/2}^{n-k}}\parallel \Aline \zeta{z_{1/2}^{2n}}\parallel \Aline{z_{1/2}^{k+n}}{z_1^n}$ and finally, $z_{1/2}^{k+n}\in \Aline {z_0^k}{z_1^n}$. 

\begin{picture}(180,185)(-80,-10)
\put(0,0){\line(1,1){150}}
\put(0,0){\color{cyan}\line(0,1){150}}
\put(15,0){\color{cyan}\line(0,1){150}}
\put(30,0){\color{cyan}\line(0,1){150}}
\put(0,0){\color{teal}\vector(1,0){180}}
\put(185,-3){$h$}
\put(0,30){\line(1,-1){30}}
\put(0,30){\color{teal}\line(1,0){30}}
\put(0,0){\color{red}\line(1,3){15}}
\put(15,45){\color{teal}\line(1,0){30}}
\put(0,60){\color{teal}\line(1,0){60}}
\put(0,60){\color{red}\line(1,3){30}}
\put(0,105){\color{teal}\line(1,0){105}}
\put(0,105){\color{red}\line(1,3){15}}
\put(0,150){\color{teal}\line(1,0){150}}

\put(0,0){\circle*{3}}
\put(-6,-7){$o$}
\put(15,0){\circle*{3}}
\put(12,-8){$u'$}
\put(15,15){\circle*{3}}
\put(18,13){$e'$}
\put(30,30){\circle*{3}}
\put(32,26){$e$}
\put(45,45){\circle*{3}}
\put(47,40){$f'(n{-}k)$}
\put(60,60){\circle*{3}}
\put(62,55){$f'(2k)=f(k)$}
\put(105,105){\circle*{3}}
\put(107,100){$f'(k{+}n)$}
\put(150,150){\circle*{3}}
\put(153,148){$f'(2n)=f(n)$}
\put(0,30){\circle*{3}}
\put(-11,28){$z^1_0$}
\put(15,45){\circle*{3}}
\put(15,47){\color{white}\line(0,1){11}}
\put(11,50){\small $z_{1/2}^{n-k}$}
\put(0,60){\circle*{3}}
\put(-11,58){$z_0^k$}
\put(0,105){\circle*{3}}
\put(-8,102){$\zeta$}
\put(15,105){\circle*{3}}
\put(15,92){\color{white}\line(0,1){11}}
\put(12,93.5){\small$z_{1/2}^{k{+}n}$}
\put(15,150){\circle*{3}}
\put(11,155){\small$z^{2n}_{1/2}$}
\put(30,0){\circle*{3}}
\put(28,-8){$u$}
\put(30,150){\circle*{3}}
\put(30,154){$z_1^n$}
\end{picture}

\end{proof}

\begin{theorem}\label{t:p-invertible-add=>field-elementary} Let $X$ be a projective space. If $X$ is invertible-add, then every ternar $R$ of $X$ is field-elementary. 
\end{theorem}

\begin{proof} Assume that $X$ is invertible-add. Given any ternar $R$ of $X$, we need to find a field $\IF$ and a function $f:\IF\to R$ such that $f(1)=1$ and $f(x)_\times f(y)_+f(z)=f(x{\cdot} y+z)$ for every $x,y,z\in\IF$. Find a plane $\Pi$ in $X$ and a projective base $uowe$ in $\Pi$ whose ternar $\Delta$ is isomorphic to the ternar $R$. By definition, the projective space $X$ is $3$-long. If $|X|_2=3$, then by definition, $R=\{0,1\}$ is isomorphic to the two-element field $\{0,1\}$. So, assume that $|X|_2\ge 4$. 

Let $h\in \Aline ou\cap \Aline we$ and $v\in \Aline ow\cap\Aline ue$ be the horizontal and vertical directions on the affine plane $\Pi\setminus\Aline hv$.
By definition, the ternar $\Delta$ of the projective base $uowe$ of the projective plane $\Pi$ is isomorphic to the ternar $\Delta$ of the affine base $uow$ of the affine plane $\Pi\setminus\Aline hv$. Since the projective space $X$ is invertible-add, so is the projective plane $\Pi$ and its affine subplane $\Pi\setminus\Aline hv$. By Theorem~\ref{t:invertible-add<=>Z-elementary}, the ternar $\Delta$ of the based affine plane $(\Pi,uow)$ is $\IZ$-elementary. Then there exists a function $f:\IZ\to \Delta$ such that $f(1)=1$ and $f(x{\cdot}y+z)=f(x)_\times f(y)_+f(z)$ for every $x,y,z\in\IZ$.
If $f$ is not injective, then by Theorem~\ref{t:Zelementary1-6}, there exists a prime number $p$ and a function $\varphi:\IF_p\to \Delta$ from the $p$-element field $\IF_p\defeq \IZ/p\IZ$ such that $f=\varphi\circ q$, where $q:\IZ\to \IF_p$ is the quotient homomorphism.
It follows that $\varphi(1)=1$ and $\varphi(x)_\times \varphi(y)_+\varphi(z)=\varphi(x{\cdot}y+z)$ for every $x,y,z\in\IZ$, witnessing that the ternar $R$ is $\IF_p$-elementary.

It remains to consider the case of injective function $f:\IZ\to \Delta$.  It this case we shall construct a function $\bar f:\IQ\to\Delta$ such that $\bar f(1)=1$ and $\bar f(x{\cdot}y+z)=\bar f(x)_\times \bar f(y)_+\bar f(z)$ for every $x,y,z\in\IQ$.

Let $V_\w\defeq\Aline ow$, $V_1\defeq \Aline ue$, $V_0\defeq \Aline hv$.
For every $n\in\IZ$, consider the point $f(n)\in \Delta$ and for every $i\in\{1,\w\}$, let $z_i^n$ be the unique common point of the lines $V_i$ and $\Aline{f(n)}h$. Let $z_0^n$ be the unique common point of the lines $V_0$ and $\Aline o{z_1^n}$. Since the ternar $\Delta$ is $\IZ$-elementary, for all $k,n\in\IZ$, the lines
$\Aline o{z_1^n}\setminus V_0$ and $\Aline {z_\w^k}{z_1^{k+n}}\setminus V_0$ in the affine plane $\Pi$ are parallel and hence $z_0^n\in \Aline {z_\w^k}{z_1^{k+n}}$ for all $n,k\in\IZ$. The injectivity of the function $f$ implies that the points $z_0^n$, $n\in\IZ$, are pairwise distinct. By Lemma~\ref{l:V1/2}, there exists a line $V_{2}$ in $X$ containing the point $v$ and an indexed family of points $\{z_{2}^n:n\in\IZ\}$ in $V_2\setminus(V_0\cup V_1\cup V_\infty)$ such that $v^{n+k}_{2}\in \Aline{z_\w^n}{z_1^k}$ for all $n,k\in\IZ$. In particular, $z_{2}^1\in \Aline {z_\w^0}{z_1^1}=\Delta$.

By induction, for every $m\ge 2$ we can apply Lemma~\ref{l:V1/2} and construct a line $V_{m}$ in $X$ with $v\in V_m$, and an indexed family of points $\{z_{m}^n:n\in\IZ\}$ in $V_{m}\setminus(V_0\cup V_{m-1}\cup V_{m-2})$ such that $z^{n+k}_{m}\in\Aline {z^n_\w}{z^k_{m-1}}$ for all $n,k\in\IZ$.  

\begin{claim}\label{cl:znm+k} For every $m\in\IN$ and $n,k\in\IZ$, we have $z_m^{nm+k}\in\Aline {z_\w^n}{z_1^{n+k}}$.
\end{claim}

\begin{proof} This claim will be proved by induction on $m\in\IN$. For $m=1$ the inclusion $z_1^{n+k}\in\Aline{z_\w^n}{z_1^{n+k}}$ is trivial. Assume that for some $m\in\IN$ we know that  $z_{m}^{nm+k}\in\Aline {z_\w^n}{z_1^{n+k}}$ for all $k\in\IZ$.

The choice of the points $z_{m+1}^k$ ensures that for all $k,n\in\IZ$ we have $$z_{m+1}^{n(m+1)+k}\in \Aline {z_\w^n}{z_{m}^{nm+k}}\subseteq \overline{\{z_\w^n\}\cup\Aline {z_\w^n}{z_1^{n+k}}}=\Aline{z_\w^n}{z_1^{n+k}}.$$
\end{proof}

For every $m\in\IN$, consider the affine base $u_mov_m$ in the affine plane $\Pi\setminus\Aline hv$ such that $u_m\in\Aline ou$, $v_m\in\Aline ov$ and the point $ z^1_m$ coincides with the diunit of the affine base $u_mov_m$. Claim~\ref{cl:znm+k} implies that $z^1_m\in \Aline{z_\w^0}{z_1^1}=\Aline oe=\Delta$. Let $\Delta_m$ be the ternar of the based affine plane $(\Pi,u_mov_m)$. Since the ternar $\Delta_m$ is $\IZ$-elementary, there exists a function $f_m:\IZ\to\Delta_m$ such that $f_m(1)=z^1_m$ and $f_m(x{\cdot} y+z)=f_m(x)_\times f_m(y)_+f_m(z)$ for all $x,y,z\in \IZ$. Using Claim~\ref{cl:znm+k}, we can show that $f_m(m)=e$ and hence $f_m(mk)=f(k)$ for all $k\in\IZ$. Now define a function $\bar f:\IQ\to\Delta$ assigning to every rational number $\frac mn$ the unique point $\bar f(\frac nm)\in\Delta\cap\Aline {z^n_m}{h}$. Then $\bar f(1)=z^1_1=e$ and $\bar f(x{\cdot} y+z)=\bar f(x)_\times \bar f(y)_+\bar f(z)$ for all $x,y,z\in\IQ$, witnessing that the ternar $\Delta$ of the based affine plane $(\Pi,uow)$ is $\IQ$-elementary. 
\end{proof}

Let us recall that a liner $X$ is defined to be \index{4-Pappian liner}\index{liner!4-Pappian}\defterm{$4$-Pappian} if the closure of every $4$-element subset in $X$ is a Pappian liner.

\begin{theorem}[Moufang, 1931]\label{t:invertible-add<=>4-Pappian} A projective space $X$ is invertible-add if and only if $X$ is $4$-Pappian.
\end{theorem}

\begin{proof} To prove the ``if'' part, assume that a projective space $X$ is $4$-Pappian, and take any distinct points $h,v\in X$, $d\in\Aline hv$, and $a,b,c,o,a',b',c'\in X\setminus\Aline hv$ such that $v\in \Aline b{c}\cap\Aline oa\cap\Aline o{a'}\cap\Aline {b'}{c'}$ and $h\in \Aline ab\cap\Aline oc\cap\Aline o{c'}\cap\Aline{a'}{b'}$ and $d\in \Aline ob\cap\Aline o{b'}\cap\Aline a{c'}$. We have to prove that $d\in \Aline {a'}c$.

\begin{picture}(120,165)(-170,-60)
\linethickness{=0.6pt}
\put(90,0){\color{teal}\line(-1,0){120}}
\put(0,90){\color{cyan}\line(0,-1){120}}
\put(0,90){\line(1,-1){90}}
\put(0,90){\color{cyan}\line(-1,-3){45}}
\put(90,0){\color{teal}\line(-3,-1){135}}
\put(45,45){\color{red}\line(-1,-1){90}}
\put(45,45){\color{red}\line(-3,-5){45}}
\put(45,45){\color{red}\line(-5,-3){75}}
\put(0,90){\color{cyan}\line(1,-5){18}}
\put(90,0){\color{teal}\line(-5,1){90}}

\put(0,0){\circle*{3}}
\put(-6,1){$o$}
\put(18,0){\circle*{3}}
\put(17,-7){$c$}
\put(0,18){\circle*{3}}
\put(-7,18){$a$}
\put(15,15){\circle*{3}}
\put(14.5,18.6){$b$}
\put(90,0){\color{teal}\circle*{3}}
\put(93,-3){\color{teal}$h$}
\put(0,90){\color{cyan}\circle*{3}}
\put(-2,93){\color{cyan}$v$}
\put(45,45){\color{red}\circle*{3}}
\put(47,46){\color{red}$d$}
\put(-30,0){\circle*{3}}
\put(-37,-1){$c'$}
\put(0,-30){\circle*{3}}
\put(-2,-38){$a'$}
\put(-45,-45){\circle*{3}}
\put(-48,-55){$b'$}
\end{picture}

Consider the $4$-element set $A\defeq\{o,a,b,c\}$ and observe that 
$\{h,v\}\subseteq(\Aline ab\cap\Aline oc)\cup(\Aline bc\cap\Aline oa)\subseteq\langle A\rangle$, $d\in \Aline ob\cap\Aline hv\subseteq\langle A\rangle$, $c'\in \Aline da\cap\Aline ac\subseteq\langle A\rangle$, $b'\in\Aline v{c'}\cap\Aline ob\subseteq\langle A\rangle$, and $a'\in \Aline oa\cap\Aline {b'}h\subseteq\langle A\rangle$. Since $X$ is $4$-Pappian, the subliner $\langle A\rangle$ is Pappian, which implies $d\in\Aline{a'}c$.
\smallskip

To prove the ``only if'' part, assume that a projective plane $X$ is invertible-add and fix any $4$-element set $A=\{u,o,w,e\}$ in $X$. We have to show that the subliner $\langle A\rangle$ of $X$ is Pappian. If $uowe$ is not a quadrangle in $X$, then $\langle A\rangle=A$ and the subliner $\langle A\rangle$ is Pappian because it contains no two concurrent lines of length $\ge 4$. So, assume that $uowe$ is a quadrangle and hence $uowe$ is a projective base for its flat hull $\Pi\defeq\overline{\{u,o,w,e\}}$. By Theorem~\ref{t:p-invertible-add=>field-elementary}, the ternar $\Delta$ of the based projective plane $(\Pi,uowe)$ is field-elementary and hence there exists a minimal field $\IF\in\{\IQ\}\cup\{\IZ/p\IZ:p$ is prime$\}$ and a function $f:\IF\to \Delta$ such that $f(1)=e$ and $f(x{\cdot}y+z)=f(x)_\times f(y)_+f(z)$ for all $x,y,z\in\IF$. Let $\IF\times\IF$ be the coordinate plane of the field $\IF$ and $\overline{\IF\times\IF}$ be its spread completion. It is clear that $\overline{\IF\times\IF}$ is a Pappian projective plane. We claim that $\overline{\IF\times\IF}$ is isomorphic to a subliner $\langle A\rangle$. Since $X$ is a projective plane, there exist unique points $h\in \Aline ou\cap\Aline we\subseteq\langle A\rangle$ and $v\in \Aline ow\cap\Aline ue\subseteq\langle A\rangle$. Define the map $F:\IF\times\IF\to X$ assigning to every ordered pair $(x,y)\in\IF\times\IF$ the unique point $F(x,y)\in\Aline {f(x)}v\cap\Aline {f(y)}h$. The homomorphism property of the map $f$ ensures that $F$ is an isomorphism of the affine plane $\IF\times\IF$ onto the subliner $F[\IF\times\IF]$ of $X$. Taking into account that the field $\IF$ is prime, we can show that $F[\IF\times\IF]\subseteq\langle A\rangle$. Extend the isomorphism $F$ to an isomorphism $\overline F:\overline{\IF\times\IF}\to \langle A\rangle$ assigning to every spread $\mathcal S$ of parallel lines in $\IF\times \IF$ the unique common point of the intersection $\bigcap_{L\in\mathcal S}\overline{F[L]}\subseteq \langle A\rangle$. Then the liner $\langle A\rangle=\overline F[\overline{\IF\times\IF}]$ is isomorphic to the Pappian liner $\overline{\IF\times\IF}$ and hence $\langle A\rangle$ is Pappian.
\end{proof}

\section{Inversive-add projective spaces}

\begin{definition}\label{d:p-inversive-plus} A projective space $X$ is called \index{inversive-plus projective space}\index{projective space!inversive-plus}\defterm{inversive-plus} if for every distinct points $h,v\in X$, $\delta\in\Aline hv$, and $a,b,c,d,a',b,c',d'\in X\setminus\Aline hv$,
$$
\big((h\in \Aline ad\cap\Aline bc\cap\Aline{a'}{d'}\cap\Aline{b'}{c'})\cap(v\in\Aline ab\cap\Aline cd\cap\Aline{a'}{b'}\cap\Aline{c'}{d'})\cap (\delta\in \Aline a{a'}\cap\Aline {b}{b'}\cap\Aline c{c'}\cap\Aline ac)\big)\Ra(\delta\in \Aline d{d'}).$$

\begin{picture}(120,148)(-170,-50)
\linethickness{=0.6pt}
\put(90,0){\color{teal}\line(-1,0){120}}
\put(0,90){\color{cyan}\line(0,-1){120}}
\put(0,90){\line(1,-1){90}}
\put(0,90){\color{cyan}\line(-1,-3){45}}
\put(90,0){\color{teal}\line(-3,-1){135}}
\put(45,45){\color{red}\line(-1,-1){90}}
\put(45,45){\color{red}\line(-3,-5){45}}
\put(45,45){\color{red}\line(-5,-3){75}}
\put(0,90){\color{cyan}\line(1,-6){12.8}}
\put(90,0){\color{teal}\line(-6,1){77}}

\put(0,90){\color{cyan}\line(35,-113){24.5}}
\put(90,0){\color{teal}\line(-113,35){79}}

\put(13,13){\circle*{2.4}}
\put(4,11){$a'$}
\put(21.3,21.3){\circle*{2.4}}
\put(20.5,25){$c'$}
\put(10.9,24.6){\circle*{2.4}}
\put(2.5,24){$b'$}
\put(24.6,10.9){\circle*{2.4}}
\put(22,2){$d'$}
\put(90,0){\circle*{2.4}}
\put(93,-3){$h$}
\put(0,90){\circle*{2.4}}
\put(-2,93){$v$}
\put(45,45){\circle*{2.4}}
\put(47,46){$\delta$}
\put(-30,0){\circle*{2.4}}
\put(-37,-2){$b$}
\put(0,0){\circle*{2.4}}
\put(-6,2){$c$}
\put(0,-30){\circle*{2.4}}
\put(-1,-38){$d$}
\put(-45,-45){\circle*{2.4}}
\put(-50,-52){$a$}
\end{picture}
\end{definition}

Definition~\ref{d:p-inversive-plus}, \ref{d:a+inversive}, and Theorems~\ref{t:+inversive} imply the following characterization of inversive-plus projective spaces.

\begin{theorem}[Bol, 1937]\label{t:p+inversive} For a projective space $X$ the following statements are equivalent:
\begin{enumerate}
\item $X$ is inversive-plus;
\item for every hyperplane $H\subseteq X$, the affine liner $X\setminus H$ is inversive-plus;
\item for every ternar $R$ of $X$, the plus loop $(R,+)$ is inversive;
\item for every ternar $R$ of $X$, the plus loop $(R,+)$ is left-inversive;
\item for every ternar $R$ of $X$, the plus loop $(R,+)$ is right-inversive;
\item for every ternar $R$ of $X$, the plus loop $(R,+)$ is Moufang.
\end{enumerate}
\end{theorem}

\begin{remark} Since every inversive magma is invertible, every inversive-plus projective space is invertible-plus and hence invertible-add.
\end{remark}

\begin{remark} In literature, inversive-plus projective planes are known as projective planes satisfying {\em the  Bol condition}, see \cite[p.53]{Pickert}.
\end{remark}  

\begin{definition}\label{d:p-inversive-puls} A projective space $X$ is called \index{inversive-puls projective space}\index{projective space!inversive-puls}\defterm{inversive-puls} if for every distinct lines $L,L'\subseteq X$ and distinct points $a,b,c,d\in L\setminus L'$, $a',b',c',d'\in L'\setminus L$, $h\in L\cap L'$ and $v\in\Aline a{a'}\cap\Aline b{b'}\cap\Aline c{c'}\cap\Aline d{d'}$, $\Aline a{b'}\cap\Aline {c}{d'}\subseteq \Aline hv$ implies $\Aline {a'}b\cap\Aline {c'}d\subseteq \Aline hv$.

\begin{picture}(150,140)(-120,-10)
\linethickness{=0.6pt}
\put(0,0){\color{teal}\line(1,0){120}}
\put(0,120){\line(1,-1){120}}
\put(0,0){\color{cyan}\line(1,1){60}}
\put(0,120){\color{violet}\line(1,-4){30}}
\put(120,0){\color{teal}\line(-4,1){96}}
\put(30,0){\color{cyan}\line(1,2){30}}
\put(60,0){\color{cyan}\line(0,1){60}}
\put(0,0){\color{red}\line(2,1){80}}
\put(80,40){\color{red}\line(-4,-5){32}}
\put(60,60){\color{cyan}\line(-1,-5){12}}
\put(0,120){\color{violet}\line(1,-2){60}}

\put(0,0){\circle*{3}}
\put(-3,-10){$a'$}
\put(30,0){\circle*{3}}
\put(27,-10){$b'$}
\put(48,0){\circle*{3}}
\put(45,-10){$c'$}
\put(60,0){\circle*{3}}
\put(57,-10){$d'$}
\put(120,0){\color{teal}\circle*{3}}
\put(123,-2){\color{teal}$h$}
\put(24,24){\circle*{3}}
\put(17,23){$a$}
\put(40,20){\circle*{3}}
\put(37,22){$b$}
\put(51.4,17.1){\circle*{3}}
\put(53,18){$c$}
\put(60,15){\circle*{3}}
\put(63,15){$d$}
\put(0,120){\color{violet}\circle*{3}}
\put(60,60){\color{cyan}\circle*{3}}
\put(60,63){\color{cyan}$v$}
\put(80,40){\color{red}\circle*{3}}
\end{picture}
\end{definition}


\begin{theorem}\label{t:p-inversive-puls} For a projective space $X$ the following statements are equivalent:
\begin{enumerate}
\item $X$ is inversive-puls;
\item for every hyperplane $H\subseteq X$, the affine liner $X\setminus H$ is inversive-puls;
\item for every hyperplane $H\subseteq X$, the affine liner $X\setminus H$ is left-inversive-puls;
\item for every hyperplane $H\subseteq X$, the affine liner $X\setminus H$ is right-inversive-puls;
\item for every ternar $R$ of $X$, the puls loop $(R,\!\puls\!)$ is inversive;
\item for every ternar $R$ of $X$, the puls loop $(R,\!\puls\!)$ is left-inversive;
\item for every ternar $R$ of $X$, the puls loop $(R,\!\puls\!)$ is right-inversive;
\item for every ternar $R$ of $X$, the puls loop $(R,\!\puls\!)$ is Moufang.
\end{enumerate}
\end{theorem}

\begin{proof} The implications $(1)\Leftrightarrow(3)\Leftrightarrow(6)$ follow from Definitions~\ref{d:p-inversive-plus}, \ref{d:left-inversive-puls}, and Theorem~\ref{t:left-inversive-puls<=>}. The equivalence $(4)\Leftrightarrow(7)$ follows from Theorem~\ref{t:right-inversive-puls<=>}. 

If $\|X\|\ge 4$, then the projective space $X$ is Desarguesian, by Theorem~\ref{t:proaffine-Desarguesian}. In this case the conditions (1)--(8) hold and hence are equivalent. So, assume that $\|X\|=3$ and hence $X$ is a projective plane. Let $X^*$ be the dual projective plane to $X$.
\smallskip

$(6)\Leftrightarrow(7)$ By Theorems~\ref{t:plus-puls-duality} and \ref{t:p+inversive}, for every ternar $R$ of $X$ the puls loop $(R,\!\puls\!)$ is left-inversive if and only if  for every ternar $R$ of the dual projective plane $X^*$ the plus loop $(R,+)$ is left-inversive if and only if for every ternar $R$ of $X^*$ the plus loop $(R,+)$ is right-inversive if and only if for every ternar $R$ of $X$ the puls loop $(R,\!\puls\!)$ is right-inversive. 
\smallskip

The equivalence $(1)\Leftrightarrow(3)\Leftrightarrow(6)\Leftrightarrow(7)\Leftrightarrow(4)$ implies the equivalence of the conditions (1)--(7).
\smallskip

$(5)\Leftrightarrow(8)$.  By Theorems~\ref{t:plus-puls-duality} and \ref{t:p+inversive}, for every ternar $R$ of $X$ the puls loop $(R,\!\puls\!)$ is inversive if and only if  for every ternar $R$ of the dual projective plane $X^*$ the plus loop $(R,+)$ is inversive if and only if for every ternar $R$ of $X^*$ the plus loop $(R,+)$ is Moufang if and only if for every ternar $R$ of $X$ the puls loop $(R,\!\puls\!)$ is Moufang. 
\end{proof}

\begin{exercise} Find a direct proof of the equivalence $(6)\Leftrightarrow(7)$ in Theorem~\ref{t:p-inversive-puls} establishing the equivalence of the configurations determining the left-inversivity and right-inversivity of puls loop of a ternar of a projective space.
\end{exercise}

\begin{remark} Since every inversive magma is invertible, every inversive-puls liner is invertible-puls and hence invertible-add.
\end{remark}

Theorems~\ref{t:plus-puls-duality}, \ref{t:p+inversive}, \ref{t:p-inversive-puls} imply the following duality between inversive-plus and inversive-puls projective planes.

\begin{proposition}\label{p:inversive-plus<=>inversive-puls-dual} A projective plane is inversive-plus if and only if its dual projective plane is inversive-puls.
\end{proposition}

\begin{definition} A projective space is called \defterm{inversive-add} if it is inversive-plus and inversive-puls.
\end{definition}

Theorems~\ref{t:plus-puls-duality} implies that the inversivity-add is a self-dual property:

\begin{proposition} A projective plane is inversive-add if and only if its dual projective plane is inversive-add.
\end{proposition}

\begin{problem} Is the inversivity-plus equivalent to the inversivity-puls for projective planes?
\end{problem}
\pagebreak

\section{Associative-add projective spaces}

\begin{definition}\label{d:p+associative} A projective space $X$ is defined to be \index{associative-plus projective space}\index{projective space!associative-plus}\defterm{associative-plus} if for every distinct points $h,v\in X$, $\delta\in\Aline hv$, and $a,b,c,d,a',b',c',d'\in X\setminus\Aline hv$,
$$
\big((h\in \Aline ad\cap\Aline bc\cap\Aline{a'}{d'}\cap\Aline{b'}{c'})\cap(v\in\Aline ab\cap\Aline cd\cap\Aline{a'}{b'}\cap\Aline{c'}{d'})\cap (\delta\in \Aline a{a'}\cap\Aline {b}{b'}\cap\Aline c{c'})\big)\Ra(\delta\in \Aline d{d'}).$$

\begin{picture}(120,130)(-170,-30)
\linethickness{=0.6pt}
\put(90,0){\color{teal}\line(-1,0){120}}
\put(0,90){\color{cyan}\line(0,-1){105}}
\put(0,90){\line(1,-1){90}}
\put(0,90){\color{cyan}\line(-1,-3){36.8}}
\put(90,0){\color{teal}\line(-61,-10){127}}
\put(45,45){\color{red}\line(-1,-1){45}}
\put(45,45){\color{red}\line(-5,-3){75}}
\put(45,45){\color{red}\line(-5,-4){82}}
\put(45,45){\color{red}\line(-3,-4){45}}
\put(90,0){\color{teal}\line(-7,2){81}}
\put(0,90){\color{cyan}\line(2,-7){22}}
\put(90,0){\color{teal}\line(-5,1){81}}
\put(0,90){\color{cyan}\line(1,-8){9.2}}

\put(21.5,13.7){\circle*{2.4}}
\put(19,5){$d$}
\put(9.3,16.3){\circle*{2.4}}
\put(5,9){$a$}
\put(8.4,23.1){\circle*{2.4}}
\put(2,23){$b$}
\put(-36.8,-20.3){\circle*{2.4}}
\put(-40,-29){$a'$}
\put(0,-15){\circle*{2.4}}
\put(-2,-24){$d'$}

\put(21.5,26.2){\color{white}\circle*{5.7}}
\put(20,20){\circle*{2.4}}
\put(19.2,23.3){$c$}

\put(90,0){\circle*{2.4}}
\put(93,-3){$h$}
\put(0,90){\circle*{2.4}}
\put(-2,93){$v$}
\put(45,45){\circle*{2.4}}
\put(47,46){$\delta$}
\put(-30,0){\circle*{2.4}}
\put(-37,-2){$b'$}
\put(0,0){\circle*{2.4}}
\put(-8,-9){$c'$}
\end{picture}

\end{definition}

Definition~\ref{d:p+associative}, \ref{d:a+associative} and Theorem~\ref{t:add-ass<=>} imply the following characterization of associative-plus projective spaces.

\begin{proposition}\label{p:p+associative} For a projective space $X$, the following conditions are equivalent:
\begin{enumerate}
\item $X$ is associative-plus;
\item for every hyperplane $H$ in $X$, the affine liner $X\setminus H$ is associative-plus;
\item for every ternar $R$ of $X$, the plus loop $(R,+)$ is associative.
\end{enumerate}
\end{proposition}

\begin{remark} Since every associative magma is inversive, every associative-plus projective space is inversive-plus.
\end{remark}

\begin{definition}\label{d:p-associative-puls} A projective space $X$ is defined to be \index{associative-puls projective space}\index{projective space!associative-puls}\defterm{associative-puls} if for any distinct lines $L,L',\Lambda\subset X$ with $L\cap L'\cap\Lambda\ne\varnothing$ and any distinct points $a,b,c,d\in L\setminus L'$ and $a',b',c',d'\in L'\setminus L$, $
(\Aline a{a'}\cap\Aline c{c'})\cup(\Aline b{b'}\cap\Aline d{d'})\cup (\Aline a{b'}\cap\Aline c{d'})\subseteq \Lambda$ implies $\Aline {a'}b\cap\Aline {c'}d\subseteq \Lambda$.

\begin{picture}(150,150)(-140,-15)
\linethickness{=0.6pt}
\put(-30,0){\line(1,0){150}}
\put(85,-10){$L'$}
\put(0,0){\color{violet}\line(0,1){120}}
\put(-30,0){\color{cyan}\line(1,1){75}}
\put(0,0){\color{blue}\line(1,1){60}}
\put(0,120){\line(1,-1){120}}
\put(95,30){$\Lambda$}
\put(0,30){\line(4,-1){120}}
\put(90,10){$L$}
\put(-30,0){\color{red}\line(9,4){103.85}}
\put(60,0){\color{blue}\line(0,1){60}}
\put(60,0){\color{violet}\line(-1,2){60}}
\put(45,75){\color{cyan}\line(1,-9){8.3}}
\put(53.3,0){\color{red}\line(4,9){20.5}}

\put(-30,0){\circle*{2.5}}
\put(-33,-9){$a'$}
\put(0,0){\circle*{2.5}}
\put(-3,-10){$b'$}
\put(53.3,0){\circle*{2.5}}
\put(50,-10){$c'$}
\put(60,0){\circle*{2.5}}
\put(58,-10){$d'$}
\put(0,30){\circle*{2.5}}
\put(-8,29){$a$}
\put(24,24){\circle*{2.5}}
\put(21,27){$b$}
\put(51.4,17.14){\circle*{2.5}}
\put(45,11){$c$}
\put(60,15){\circle*{2.5}}
\put(63,15){$d$}
\put(0,120){\color{violet}\circle*{3}}
\put(45,75){\color{cyan}\circle*{3}}
\put(60,60){\color{blue}\circle*{3}}
\put(73.85,46.15){\color{red}\circle*{3}}
\put(120,0){\circle*{2.5}}

\end{picture}
\end{definition}

Definition~\ref{d:p-associative-puls}, \ref{d:a-associative-puls} and Theorem~\ref{t:associative-puls<=>} imply the following characterization of associative-puls projective spaces.

\begin{proposition}\label{p:puls-associative} For a projective space $X$, the following conditions are equivalent:
\begin{enumerate}
\item $X$ is associative-puls;
\item for every hyperplane $H$ in $X$, the affine liner $X\setminus H$ is associative-puls;
\item for every ternar $R$ of $X$, the puls loop $(R,\!\puls\!)$ is associative.
\end{enumerate}
\end{proposition}

\begin{remark} Since every associative magma is inversive, every associative-puls projective space is inversive-puls.
\end{remark}

Theorem~\ref{t:plus-puls-duality} implies that the following duality between associative-plus and associative-puls projective planes.

\begin{corollary} A projective plane is associative-plus if and only if its dual projective plane is associative-puls.
\end{corollary}

\begin{definition} A projective space is called \defterm{associative-add} if it is associative-plus and associative-puls.
\end{definition}

Theorem~\ref{t:plus-puls-duality} implies that the associativity-add is a self-dual property:

\begin{proposition} A projective plane is associative-add if and only if its dual projective plane is associative-add.
\end{proposition}

\begin{problem}\label{prob:ass-plus<=>ass-puls} Is the associativity-plus equivalent to the associativity-puls for projective planes?
\end{problem}

\begin{remark} Problem~\ref{prob:ass-plus<=>ass-puls} has an affirmative answer for finite projective planes, see Theorem~\ref{t:proj-finite<=>}.
\end{remark}

\begin{remark} In literature, associative-plus projective planes are known as projective planes satisfying the {\em little Reidemeister condition.}
\end{remark}  




\section{Commutative-add projective spaces}

\begin{definition}\label{d:p+commutative} A projective space $X$ is called \index{commutative-plus projective space}\index{projective space!commutative-plus}\defterm{commutative-plus} if for every distinct points $h,v\in X$, $d\in \Aline hv$, and $a,b,c,a',b',c'\in X\setminus\Aline hv$,
$$\big((h\in \Aline ab\cap\Aline c{c'}\cap\Aline {a'}{b'})\cap(v\in \Aline ac\cap\Aline b{b'}\cap\Aline {a'}{c'})\cap (d\in \Aline a{a'}\cap\Aline {b}{c'})\big)\Ra(d\in\Aline {b'}c).$$

\begin{picture}(200,150)(-170,-55)
\linethickness{0.6pt}
\put(80,0){\color{teal}\line(-1,0){100}}
\put(0,80){\color{cyan}\line(0,-1){109}}
\put(0,80){\line(1,-1){80}}
\put(0,80){\color{cyan}\line(-1,-4){30}}
\put(80,0){\color{teal}\line(-11,-4){110}}
\put(40,40){\color{red}\line(-3,-2){60}}
\put(40,40){\color{red}\line(-7,-8){70}}
\put(40,40){\color{red}\line(-40,-69){40}}
\put(0,80){\color{cyan}\line(17,-80){17}}
\put(80,0){\color{teal}\line(-80,13.5){80}}

\put(16.8,0){\circle*{2.5}}
\put(16,-8){$c$}
\put(0,13.3){\circle*{2.5}}
\put(-6,13){$b$}
\put(14.6,11){\circle*{2.5}}
\put(14,16){$a$}
\put(-20,0){\circle*{2.5}}
\put(-28,-2){$c'$}
\put(-30,-40){\circle*{2.5}}
\put(-33,-48){$a'$}
\put(0,-29){\circle*{2.5}}
\put(-3,-40){$b'$}
\put(80,0){\circle*{2.5}}
\put(83,-3){$h$}
\put(0,80){\circle*{2.5}}
\put(-2,83){$v$}
\put(40,40){\circle*{2.5}}
\put(42,41){$d$}
\end{picture}
\end{definition}

Definition~\ref{d:p+commutative}, \ref{d:a+commutative} and Theorem~\ref{t:add-com<=>} imply the following characterization of commutative-plus projective spaces.

\begin{proposition}\label{p:p+commutative} For a projective space $X$, the following conditions are equivalent:
\begin{enumerate}
\item $X$ is commutative-plus;
\item for every hyperplane $H\subseteq X$, the affine liner $X\setminus H$ is commutative-plus;
\item for every ternar $R$ of $X$, the plus loop $(R,+)$ is commutative.
\end{enumerate}
\end{proposition}

Theorem~\ref{t:add-com=>add-ass} implies the following corollary.

\begin{corollary}\label{c:p:com-plus=>add-plus} Every commutative-plus projective space is associative-plus.
\end{corollary}

\begin{remark} In literature, commutative-plus projective planes are referred to as projective planes satisfying {\em the Thomsen condition}.
\end{remark}

\begin{definition}\label{d:p-commutative-puls} A projective space $X$ is called \index{commutative-puls projective space}\index{projective space!commutative-puls}\defterm{commutative-puls} if for every distinct lines $L,L',\Lambda$ with $L\cap L'\cap\Lambda\ne\varnothing$ and any distinct points  $a,b,c\in L\setminus L'$, $a',b',c'\in L'\setminus L$, $(\Aline a{b'}\cap\Aline {a'}b)\cup(\Aline b{c'}\cap\Aline{b'}c)\subseteq\Lambda$ implies $\Aline a{c'}\cap\Aline{a'}c\subseteq\Lambda$.

\begin{picture}(150,150)(-140,-15)
\linethickness{=0.6pt}
\put(-24,0){\line(1,0){144}}
\put(85,-10){$L'$}
\put(0,0){\color{violet}\line(0,1){120}}
\put(-24,0){\color{blue}\line(104,40){104}}
\put(0,0){\color{blue}\line(2,1){80}}
\put(0,120){\line(1,-1){120}}
\put(95,30){$\Lambda$}
\put(0,30){\line(4,-1){120}}
\put(90,10){$L$}
\put(40,80){\color{red}\line(-4,-5){64}}
\put(0,120){\color{violet}\line(1,-3){40}}
\put(40,80){\color{red}\line(0,-1){80}}

\put(-24,0){\circle*{3}}
\put(-27,-9){$a'$}
\put(0,0){\circle*{3}}
\put(-3,-10){$b'$}
\put(40,0){\circle*{3}}
\put(37,-10){$c'$}
\put(0,30){\circle*{3}}
\put(-8,30){$c$}
\put(32.7,21.8){\circle*{3}}
\put(33,24){$b$}
\put(40,20){\circle*{3}}
\put(41,13){$a$}

\put(0,120){\color{violet}\circle*{3}}
\put(40,80){\color{red}\circle*{3}}
\put(80,40){\color{blue}\circle*{3}}
\put(120,0){\circle*{3}}
\end{picture}
\end{definition}

Definition~\ref{d:p-commutative-puls}, \ref{d:a-commutative-puls} and Theorem~\ref{t:commutative-puls<=>} imply the following characterization of commutative-puls projective spaces.

\begin{proposition}\label{p:p-commutative-puls<=>} For a projective space $X$, the following conditions are equivalent:
\begin{enumerate}
\item $X$ is commutative-puls;
\item for every hyperplane $H\subseteq X$, the affine liner $X\setminus H$ is commutative-puls;
\item for every ternar $R$ of $X$, the puls loop $(R,\!\puls\!)$ is commutative.
\end{enumerate}
\end{proposition}

Theorem~\ref{t:com-puls=>ass-puls} implies the following corollary.

\begin{corollary}\label{c:p:com-puls=>ass-puls} Every commutative-puls projective space is associative-puls.
\end{corollary}

Theorem~\ref{t:plus-puls-duality} implies  the following duality between commutative-plus and commutative-puls projective planes.

\begin{corollary}\label{c:com-plus-puls-dual} A projective plane is commutative-plus if and only if its dual projective plane is commutative-puls.
\end{corollary}

\begin{definition} A projective space is called \defterm{commutative-add} if it is commutative-plus and commutative-puls.
\end{definition}

\begin{theorem}\label{t:p-commutative-add<=>} For a projective space $X$, the following conditions are equivalent:
\begin{enumerate}
\item $X$ is commutative-add;
\item for every hyperplane $H\subseteq X$, the affine liner $X\setminus H$ is commutative-add;
\item for every ternar $R$ of $X$, the loops $(R,+)$ and $(R,\!\puls\!)$ are commutative;
\item for every hyperplane $H\subseteq X$, the affine liner $X\setminus H$ is Boolean or Thalesian;
\item $X$ is Fano or Moufang.
\end{enumerate}
\end{theorem}

\begin{proof} The equivalence of the conditions (1)--(4) follows from Propositions~\ref{p:p+commutative}, \ref{p:p-commutative-puls<=>}, and Theorem~\ref{t:commutative-add<=>}. The implication $(5)\Ra(4)$ follows from Theorems~\ref{t:Fano<=>} and \ref{t:proj-Moufang<=>}.
\smallskip

$(4)\Ra(5)$ If $\|X\|\ge 4$ or $|X|_2=3$, then $X$ is Desarguesian and hence Moufang. So, assume that $\|X\|=3$ and $\|X|_2\ge 4$, which means that $X$ is a $4$-long projective plane. If $X$ is not Fano, then there exists a quadrangle $abcd$ in $X$ whose diagonal  points $x\in \Aline ab\cap \Aline cd$, $y\in \Aline ad\cap\Aline bc$, $z\in \Aline ac\cap\Aline bd$ are not collinear. Then $abcd$ is a parallelogram in the Playfair plane $X\setminus \Aline xy$ whose diagonal point $z$ belongs to $X\setminus \Aline xy$ and witnesses that the Playfair plane $X\setminus\Aline xy$ is not Boolean. The condition (4) ensures that the the Playfair plane $X\setminus\Aline xy$ is Thalesian, and hence the line $\Aline xy$ is translation. By analogy we can show that the Playfair plane $X\setminus\Aline yz$ is Thalesian (because it contains the non-Boolean parallelogram $acbd$).
Therefore, the projective plane $X$ contains two distinct translation lines $\Aline xy$, $\Aline yz$, and hence $X$ is Moufang, by Theorem~\ref{t:Skornyakov-San-Soucie}.
\end{proof}

Theorem~\ref{t:plus-puls-duality} implies that the commutativity-add is a self-dual property:

\begin{corollary}\label{c:com-add-self-dual} A projective plane is commutative-add if and only if its dual projective plane is commutative-add.
\end{corollary}

\begin{problem}\label{prob:com-plus<=>com-puls} Is the commutativity-plus equivalent to the commutativity-puls for projective planes?
\end{problem}

\begin{remark} Problem~\ref{prob:com-plus<=>com-puls} has an affirmative answer for finite projective planes, see  Theorem~\ref{t:proj-finite<=>}.
\end{remark}





\section{Invertible-dot projective spaces}

\begin{definition}\label{d:p-invertible-dot} A projective space $X$ is \index{invertible-dot projective space}\index{projective space!invertible-dot}\defterm{invertible-dot} if for every points 
for every distinct points $h,v,d\in X$ and $a,b,c,o,a',b',c'\in X\setminus\Aline hv$,
$$\big((v\in \Aline b{c}\cap\Aline oa\cap\Aline o{a'}\cap\Aline {b'}{c'})\wedge(h\in \Aline ab\cap\Aline oc\cap\Aline o{c'}\cap\Aline{a'}{b'})\wedge(d\in \Aline ob\cap\Aline o{b'}\cap\Aline a{c'})\big)\Ra(d\in \Aline {a'}c).$$

\begin{picture}(120,155)(-170,-55)

\put(90,0){\color{teal}\line(-1,0){120}}
\put(0,90){\color{cyan}\line(0,-1){120}}
\put(0,90){\color{cyan}\line(-1,-3){45}}
\put(90,0){\color{teal}\line(-3,-1){135}}
\put(-45,-45){\color{red}\line(1,1){105}}
\put(-30,0){\color{red}\line(3,2){90}}
\put(0,-30){\color{red}\line(2,3){60}}
\put(0,90){\color{cyan}\line(2,-9){20}}
\put(90,0){\color{teal}\line(-9,2){90}}

\put(0,0){\circle*{2.4}}
\put(-6,2){$o$}
\put(20,0){\circle*{2.4}}
\put(20,-7){$c$}
\put(0,20){\circle*{2.4}}
\put(-7,20){$a$}
\put(16.4,16.4){\circle*{2.4}}
\put(16,20){$b$}
\put(90,0){\circle*{2.4}}
\put(93,-3){$h$}
\put(0,90){\circle*{2.4}}
\put(-2,93){$v$}
\put(60,60){\circle*{2.4}}
\put(62,61){$d$}
\put(-30,0){\circle*{2.4}}
\put(-37,-1){$c'$}
\put(0,-30){\circle*{2.4}}
\put(-2,-38){$a'$}
\put(-45,-45){\circle*{2.4}}
\put(-48,-55){$b'$}
\end{picture}
\end{definition}

\begin{theorem}\label{t:p-invertible-dot<=>} For a projective space $X$, the following conditions are equivalent:
\begin{enumerate}
\item $X$ is invertible-dot;
\item for every hyperplane $H\subset X$, the affine liner $X\setminus H$ is invertible-dot;
\item for every ternar $R$ of $X$, the dot loop $(R^*,\cdot)$ is invertible;
\item $X$ is commutative-puls.
\end{enumerate}
\end{theorem}

\begin{proof} The equivalence of the conditions (1)--(3) follows easily from 
Definition~\ref{d:p-invertible-dot}, \ref{d:a-invertible-dot} and Theorem~\ref{t:invertible-dot<=>}.

$(1)\Ra(4)$. Assume that a projective space $X$ is invertible-dot. To prove that $X$ is commutative-puls, take any distinct lines $L,L',\Lambda$ with $L\cap L'\cap\Lambda\ne\varnothing$ and any distinct points $a,b,c\in L\setminus L'$, $a',b',c'\in L'\setminus L$ such that $(\Aline a{b'}\cap\Aline{a'}b)\cup(\Aline b{c'}\cap\Aline{b'}{c})\subseteq\Lambda$. 

Consider the points
$$h\defeq b',\;v\defeq b,\;\alpha\defeq a,\;\gamma\defeq a',\;\alpha'\defeq c,\;\gamma'\defeq c'.$$
By projectivity of $X$, there exist unique points $$o\in L\cap L',\;\;\beta\in\Aline a{b'}\cap\Aline{a'}b=\Aline \alpha h\cap\Aline \gamma v,\;\;\beta'\in\Aline b{c'}\cap\Aline{b'}c=\Aline v{\gamma'}\cap\Aline h{\alpha'},\;\;\mbox{and}\;\;\delta\in \Lambda\cap\Aline a{c'}=\Aline o{\beta'}\cap\Aline \alpha{\gamma'}.$$ 
It follows that
$$v=b\in \Aline \beta\gamma\cap\Aline o\alpha\cap\Aline o{\alpha'}\cap\Aline {\beta'}{\gamma'}\quad\mbox{and}\quad h=b'\in \Aline \alpha\beta\cap\Aline o\gamma\cap\Aline o{\gamma'}\cap\Aline{\alpha'}{\beta'}$$ and $o\in\Lambda=\Aline \beta{\beta'}$ and hence $\delta\in \Aline o\beta\cap\Aline o{\beta'}\cap\Aline \alpha{\gamma'}$. Applying the  inverstibility-dot to the points $h,v,\delta,\alpha,\beta,\gamma,o,\alpha',\beta',\gamma'$, we conclude that $\delta\in\Aline{\alpha'}{\gamma}=\Aline c{a'}$. Then $\Aline a{c'}\cap\Aline c{a'}=\{\delta\}\subseteq\Lambda$, witnessing that $X$ is commutative-puls.

\begin{picture}(120,170)(-170,-60)

\put(90,0){\color{teal}\line(-1,0){120}}
\put(0,90){\color{cyan}\line(0,-1){120}}
\put(0,90){\color{cyan}\line(-1,-3){45}}
\put(90,0){\color{teal}\line(-3,-1){135}}
\put(-45,-45){\color{red}\line(1,1){105}}
\put(-30,0){\color{red}\line(3,2){90}}
\put(0,-30){\color{red}\line(2,3){60}}
\put(0,90){\color{cyan}\line(2,-9){20}}
\put(90,0){\color{teal}\line(-9,2){90}}

\put(0,0){\circle*{2.4}}
\put(-6,2){$o$}
\put(20,0){\circle*{2.4}}
\put(17,-9){$a'$}
\put(13,2){$\gamma$}
\put(0,20){\circle*{2.4}}
\put(-8,20){$a$}
\put(1,25){$\alpha$}
\put(16.4,16.4){\circle*{2.4}}
\put(16,20){$\beta$}
\put(90,0){\circle*{2.4}}
\put(93,-8){$h$}
\put(93,2){$b'$}
\put(0,90){\circle*{2.4}}
\put(-6,93){$b\;v$}
\put(60,60){\circle*{2.4}}
\put(62,61){$\delta$}
\put(-30,0){\circle*{2.4}}
\put(-30,-9){$\gamma'$}
\put(-39,-1){$c'$}
\put(0,-30){\circle*{2.4}}
\put(-6,-38){$c\,\alpha'$}
\put(-45,-45){\circle*{2.4}}
\put(-48,-55){$\beta'$}

\put(-8,45){$L$}
\put(40,-10){$L'$}
\put(-20,-28){$\Lambda$}
\end{picture}

$(4)\Ra(1)$ Assume that the projective space $X$ is commutative-puls. To prove that $X$ is invertible-dot, take any distinct points $h,v,d\in X$ and $a,b,c,o,a',b',c'\in X\setminus\Aline hv$ such that 
$v\in \Aline b{c}\cap\Aline oa\cap\Aline o{a'}\cap\Aline {b'}{c'}$, $h\in \Aline ab\cap\Aline oc\cap\Aline o{c'}\cap\Aline{a'}{b'}$ and $d\in \Aline ob\cap\Aline o{b'}\cap\Aline a{c'}$. We have to prove that $d\in \Aline {a'}c$.

Consider the concurrent lines $L\defeq\Aline ov$, $L'\defeq\Aline oh$, $\Lambda\defeq\Aline od$, and points 
$$\alpha\defeq a,\; \beta\defeq v,\;\gamma\defeq a',\;\alpha'\defeq c\;\beta'\defeq h,\;\gamma'\defeq c'.$$ Observe that $\alpha,\beta,\gamma\in L\setminus L'$, $\alpha',\beta',\gamma'\in L'\setminus L$, and 
$$(\Aline \alpha{\beta'}\cap\Aline {\alpha'}{\beta})\cup(\Aline \beta{\gamma'}\cap\Aline{\beta'}\gamma)=(\Aline ah\cap\Aline cv)\cup(\Aline v{c'}\cap\Aline h{a'})=\{b\}\cup \{b'\}\subseteq \Lambda.$$ By the projectivity of $X$, there exists a point $z\in\Aline a{c'}\cap\Aline {a'}c$.

\begin{picture}(120,170)(-170,-60)

\put(90,0){\color{teal}\line(-1,0){120}}
\put(0,90){\color{cyan}\line(0,-1){120}}
\put(0,90){\color{cyan}\line(-1,-3){45}}
\put(90,0){\color{teal}\line(-3,-1){135}}
\put(-45,-45){\color{red}\line(1,1){105}}
\put(-30,0){\color{red}\line(3,2){90}}
\put(0,-30){\color{red}\line(2,3){60}}
\put(0,90){\color{cyan}\line(2,-9){20}}
\put(90,0){\color{teal}\line(-9,2){90}}

\put(0,0){\circle*{2.4}}
\put(-6,2){$o$}
\put(20,0){\circle*{2.4}}
\put(17,-9){$\alpha'$}
\put(13,2){$c$}
\put(0,20){\circle*{2.4}}
\put(-8,20){$\alpha$}
\put(1,25){$a$}
\put(16.4,16.4){\circle*{2.4}}
\put(16,20){$b$}
\put(90,0){\circle*{2.4}}
\put(93,-8){$h$}
\put(93,2){$\beta'$}
\put(0,90){\circle*{2.4}}
\put(-8,93){$\beta\;v$}
\put(60,60){\circle*{2.4}}
\put(62,61){$d$}
\put(-30,0){\circle*{2.4}}
\put(-30,-9){$c'$}
\put(-39,-1){$\gamma'$}
\put(0,-30){\circle*{2.4}}
\put(-6,-38){$\gamma\,a'$}
\put(-45,-45){\circle*{2.4}}
\put(-48,-55){$b'$}

\put(-8,45){$L$}
\put(40,-10){$L'$}
\put(-20,-28){$\Lambda$}
\end{picture}

Applying the commutativity-puls of $X$ to the lines $L,L',\Lambda$ and points $\alpha,\beta,\gamma\in L\setminus L'$ and $\alpha',\beta',\gamma'\in L'\setminus L$, we conclude that 
$$z\in \Aline a{c'}\cap\Aline c{a'}= \Aline {\alpha}{\gamma'}\cap\Aline {\alpha'}{\gamma}\subseteq \Lambda=\Aline od.$$
On the other hand, the assumption $d\in\Aline a{c'}$ implies $\Aline od\cap\Aline a{c'}=\{d\}$. Indeed, assuming that $\Aline od\cap\Aline a{c'}\ne\{d\}$, we would conclud that $\Aline od=\Aline a{c'}$ and hence the line $\Aline od=\Aline a{c'}$ contains the points $v\in \Aline oa$ and $h\in\Aline o{c'}$. Then $o\in \Aline od=\Aline hv$, which contradicts the choice of the point $o\notin\Aline hv$. This contradiction shows that $\Aline od\cap\Aline a{c'}=\{d\}$. 
Then 
$$z\in\Aline a{c'}\cap\Aline c{a'}\subseteq\Aline od\cap\Aline a{c'}=\{d\}$$ and hence $d=z\in \Aline a{c'}$.
\end{proof}

\begin{problem} Assume that a projective plane is invertible-dot. Is the dual projective plane invertible-dot?
\end{problem}

\section{Inversive-dot projective spaces}

\begin{definition}\label{d:p.inversive} A projective space $X$ is called \index{inversive-dot projective space}\index{projective space!inversive-dot}\defterm{inversive-dot} if for every distinct points $h,v,\delta\in X$ and $a,b,c,d,a',b,c',d',h,\delta,v\in X\setminus\Aline hv$, 
$$
\big((h\in \Aline ad\cap\Aline bc\cap\Aline{a'}{d'}\cap\Aline{b'}{c'})\cap(v\in\Aline ab\cap\Aline cd\cap\Aline{a'}{b'}\cap\Aline{c'}{d'})\cap (\delta\in \Aline a{a'}\cap\Aline {b}{b'}\cap\Aline c{c'}\cap\Aline ac)\big)\Ra(\delta\in \Aline d{d'}).$$

\begin{picture}(120,160)(-170,-50)

\put(90,0){\color{teal}\line(-1,0){120}}
\put(0,90){\color{cyan}\line(0,-1){120}}
\put(0,90){\color{cyan}\line(-1,-3){45}}
\put(90,0){\color{teal}\line(-3,-1){135}}
\put(-45,-45){\color{red}\line(1,1){148}}
\put(-30,0){\color{red}\line(40,31){133}}
\put(0,-30){\color{red}\line(31,40){103}}
\put(0,90){\color{cyan}\line(1,-6){12.8}}
\put(90,0){\color{teal}\line(-6,1){77}}

\put(0,90){\color{cyan}\line(5,-13){31}}
\put(90,0){\color{teal}\line(-13,5){80}}
\put(31,10){\circle*{2.4}}
\put(10,31){\circle*{2.4}}

\put(13,13){\circle*{2.4}}
\put(4,12){$a'$}
\put(25,25){\circle*{2.4}}
\put(2.5,28){$b'$}
\put(34,10){$d'$}
\put(90,0){\circle*{2.4}}
\put(93,-3){$h$}
\put(0,90){\circle*{2.4}}
\put(-2,93){$v$}
\put(103,103){\circle*{2.4}}
\put(106,101){$\delta$}
\put(-30,0){\circle*{2.4}}
\put(-37,-2){$b$}
\put(0,0){\circle*{2.4}}
\put(-6,2){$c$}
\put(0,-30){\circle*{2.4}}
\put(-1,-38){$d$}
\put(-45,-45){\circle*{2.4}}
\put(-50,-52){$a$}
\end{picture}
\end{definition}

Definition~\ref{d:p.inversive}, \ref{d:r-inversive-dot} and Theorems~\ref{t:right-inversive-dot<=>}, \ref{t:inversive-dot} imply the following characterization of inversive-dot projective spaces. 

\begin{theorem}\label{t:p.inversive} For a projective space $X$, the following conditions are equivalent:
\begin{enumerate}
\item $X$ is inversive-dot;
\item for every hyperplane $H\subseteq X$, the affine liner $X\setminus H$ is inversive-dot;
\item for some hyperplane $H\subseteq X$, the affine liner $X\setminus H$ is right inversive-dot;
\item for every ternar $R$ of $X$, its dot loop $(R^*,\cdot)$ is right inversive;
\item for every ternar $R$ of $X$, its dot loop $(R^*,\cdot)$ is inversive;
\item every ternar $R$ of $X$ is linear, distributive, associative-plus, and alternative-dot;
\item some ternar $R$ of $X$ is linear, distributive, associative-plus, and alternative-dot;
\item $X$ is Moufang.
\end{enumerate}
\end{theorem}

By Theorem~\ref{t:inversive-dot}, an affine space is inversive-dot if and only if $X$ is hypersymetric.

\begin{problem} Is every hypersymmetric projective space inversive-dot?
\end{problem}

\section{Associative-dot projective spaces}

\begin{definition}\label{d:p.associative} A projective space $X$ is defined to be \index{associative-dot projective space}\index{projective space!associative-dot}\defterm{associative-dot} if for every distinct points $h,v,\delta\in X$ and $a,b,c,d,a',b',c',d'\in X\setminus\Aline hv$,
$$
\big((h\in \Aline ad\cap\Aline bc\cap\Aline{a'}{d'}\cap\Aline{b'}{c'})\cap(v\in\Aline ab\cap\Aline cd\cap\Aline{a'}{b'}\cap\Aline{c'}{d'})\cap (\delta\in \Aline a{a'}\cap\Aline {b}{b'}\cap\Aline c{c'})\big)\Ra(\delta\in \Aline d{d'}).$$

\begin{picture}(120,150)(-170,-50)

\put(90,0){\color{teal}\line(-1,0){120}}
\put(0,90){\color{cyan}\line(0,-1){105}}
\put(0,90){\color{cyan}\line(-1,-3){36}}
\put(90,0){\color{teal}\line(-7,-1){126}}
\put(90,90){\color{red}\line(-1,-1){90}}
\put(90,90){\color{red}\line(-4,-3){120}}
\put(90,90){\color{red}\line(-90,-103){90}}
\put(90,90){\color{red}\line(-7,-6){126}}
\put(90,0){\color{teal}\line(-11,4){81}}
\put(0,90){\color{cyan}\line(4,-11){26.5}}
\put(90,0){\color{teal}\line(-635,170){80}}
\put(0,90){\color{cyan}\line(10,-68.5){10}}

\put(26.5,17){\circle*{2.4}}
\put(24,8){$d$}
\put(10,21.5){\circle*{2.4}}
\put(8,15){$a$}
\put(9,29.3){\circle*{2.4}}
\put(3,28){$b$}
\put(-36,-18){\circle*{2.4}}
\put(-40,-29){$a'$}
\put(0,-13){\circle*{2.4}}
\put(-2,-24){$d'$}

\put(24,24){\circle*{2.4}}
\put(23.1,27.3){$c$}

\put(90,0){\circle*{2.4}}
\put(93,-3){$h$}
\put(0,90){\circle*{2.4}}
\put(-2,93){$v$}
\put(90,90){\circle*{2.4}}
\put(92,91){$\delta$}
\put(-30,0){\circle*{2.4}}
\put(-37,-2){$b'$}
\put(0,0){\circle*{2.4}}
\put(-8,-9){$c'$}
\end{picture}
\end{definition}

Definition~\ref{d:p.associative}, \ref{d:ass-dot} and Theorem~\ref{t:ass-dot<=>} imply the following characterization of associative-dot projective spaces. 

\index[person]{Klingenberg}
\begin{theorem}[Klingenberg, 1955] \label{t:p.associative} For a projective space $X$, the following conditions are equivalent:
\begin{enumerate}
\item $X$ is associative-dot;
\item for every hyperplane $H\subset X$, the affine liner $X\setminus H$ is associative-dot;
\item for some hyperplane $H\subset X$, the affine liner $X\setminus H$ is associative-dot;
\item every ternar of $X$ is associative-dot;
\item every ternar of $X$ is linear, distributive and associative;
\item some ternar of $X$ is linear, distributive and associative;
\item $X$ is Desarguesian.
\end{enumerate}
\end{theorem}




\section{Commutative-dot projective spaces}

\begin{definition}\label{d:p.commutative} A projective space $X$ is called \index{commutative-dot projective space}\index{projective space!commutative-dot}\defterm{commutative-dot} if for every distinct points $h,v,d\in X$ and $a,b,c,a',b',c'\in X\setminus\Aline hv$,
$$\big((h\in \Aline ab\cap\Aline c{c'}\cap\Aline {a'}{b'})\cap(v\in \Aline ac\cap\Aline b{b'}\cap\Aline {a'}{c'})\cap (d\in \Aline a{a'}\cap\Aline {b}{c'})\big)\Ra(d\in\Aline {b'}c).$$

\begin{picture}(200,150)(-170,-55)

\put(80,0){\color{teal}\line(-1,0){100}}
\put(0,80){\color{cyan}\line(0,-1){109}}
\put(0,80){\color{cyan}\line(-1,-4){30}}
\put(80,0){\color{teal}\line(-11,-4){110}}
\put(80,80){\color{red}\line(-5,-4){100}}
\put(80,80){\color{red}\line(-11,-12){110}}
\put(80,80){\color{red}\line(-80,-109){80}}
\put(0,80){\color{cyan}\line(213,-800){21}}
\put(80,0){\color{teal}\line(-80,16){80}}

\put(21.3,0){\circle*{2.5}}
\put(20,-10){$b$}
\put(0,16){\circle*{2.5}}
\put(-7,15){$c$}
\put(18,12.3){\circle*{2.5}}
\put(17,17){$a$}
\put(-20,0){\circle*{2.5}}
\put(-28,-2){$b'$}
\put(-30,-40){\circle*{2.5}}
\put(-33,-48){$a'$}
\put(0,-29){\circle*{2.5}}
\put(-2,-37){$c'$}
\put(80,0){\circle*{2.5}}
\put(83,-3){$h$}
\put(0,80){\circle*{2.5}}
\put(-2,83){$v$}
\put(80,80){\circle*{2.5}}
\put(82,81){$d$}
\end{picture}
\end{definition}

Definition~\ref{d:p.commutative}, \ref{d:a.commutative} and Theorem~\ref{t:commutative-dot<=>} imply the following characterization of commutative-plus projective spaces.

\begin{theorem}\label{t:p.commutative} For a projective space $X$, the following conditions are equivalent:
\begin{enumerate}
\item $X$ is commutative-dot;
\item for every hyperplane $H\subset X$, the affine liner $X\setminus H$ is commutative-dot;
\item for some hyperplane $H\subset X$, the affine liner $X\setminus H$ is commutative-dot;
\item every ternar of $X$ is commutative-dot;
\item every ternar of $X$ is linear, distributive, associative, and commutative;
\item some ternar of $X$ is linear, distributive, associative, and commutative;
\item $X$ is Pappian.
\end{enumerate}
\end{theorem}


\section{Projective spaces of finite order}\label{s:Kegel-Luneburg}

In this section we prove that for projective spaces of finite order many geometric and algebraic properties are equivalent. A crucial result in this direction was proved by Kegel and L\"uneburg \cite{KegelLuneburg1963}.

\begin{theorems}[Kegel, L\"uneburg, 1963]\label{t:Kegel-Luneburg-main} Every finite associative-plus projective plane is Desarguesian.
\end{theorems}

\begin{proof} Assume that a finite projective plane $\Pi$ is associative-plus.
 Then $\Pi$ is invertible-plus and hence invertible-add, by Theorem~\ref{t:p-invertible-add<=>}. 
 By Corollary~\ref{c:Kegel-Luneburg}, for every ternar $R$ of $\Pi$, the plus loop $(R,+)$ is an elementary group such that either the order $|\Pi|_2-1=|R|$ is prime-power, or else $|\Pi|_2-1=|R|=p^nq$ for some prime power $p^n>1$ and some prime number $q$ dividing $p^n-1$. 
 
Let us show that the second case is impossible. To derive a contradiction, assume that $|\Pi|_2-1=p^nq$ for some prime power $p^n>1$ and some prime number $q$, dividing $p^n-1$. Fix any point $\delta\in \Pi$ and consider the pencil of lines $\mathcal L_\delta$ in $\Pi$ containing the point $\delta$. Fix any line $H\in\mathcal L_\delta$ and consider the affine subplane $\Pi\setminus H$ of the projective plane $\Pi$.  Every point of the line $H$ can be identified with a direction in the affine plane $\Pi\setminus H$. Given any point $h\in H\setminus\{\delta\}$ and any lines $L,\Lambda\in\mathcal L_\delta\setminus\{H\}$, consider the bijective map $h_{\Lambda,L}:L\to\Lambda$ assigning to every point $x\in L$ the unique point $y\in \Lambda\cap\Aline xh$. It is clear that $h_{\Lambda,L}(\delta)=\delta$. Fix any line $\Delta\in\mathcal L_\delta$ and for all distinct points $h,v\in H\setminus\{\delta\}$, consider the subset 
$$\I^\#_{\Pi\setminus H}[\Delta;h,v]\defeq\{h_{\Delta,L}v_{L,\Delta}:L\in\mathcal L_\delta\setminus\{H\}\}$$of the permutation group of the set $\Delta$. Since the affine plane $\Pi\setminus H$ is associative-plus, we can apply Theorem~\ref{t:plus<=>group} and conclude that $\I^\#_{\Pi\setminus H}[\Delta;h,v]$ is a subgroup of the permutation group of the set $\Delta$. Moreover, the group $\I^\#_{\Pi\setminus H}[\Delta;h,v]$ is isomorphic to the plus loop $(R,+)$ of a suitable ternar $R$ of the afine plane $\Pi\setminus H$. Applying Corollary~\ref{c:Kegel-Luneburg}, we conclude that $\I^\#_\Pi[\Delta;h,v]$ is a Frobenius group whose kernel $K_{\Pi\setminus H}[\Delta;h,v]$ has order $p^n$ and whose complement has order $q$. Let $\I^\#_{\Pi\setminus H}[\Delta]$ be the subgroup of the permutation group of the set $\Delta$, generated by the union $\bigcup_{h,v\in H\setminus\{\delta\}}\I^\#_{\Pi\setminus H}[\Delta;h,v]$. Let $S$ be a Sylow $p$-subgroup of the (finite) group $\I^\#_{\Pi\setminus H}[\Delta]$.  

\begin{claim}\label{cl:S-kernel} For every distinct points $h,v\in H\setminus\{\delta\}$, the intersection $S\cap\I^\#_{\Pi\setminus H}[\Delta;h,v]$ coincides with the kernel $K_{\Pi\setminus H}[\Delta;h,v]$ of the Frobenius group
$\I^\#_{\Pi\setminus H}[\Delta;h,v]$.
\end{claim}

\begin{proof} Since the kernel $K\defeq K_{\Pi\setminus H}[\Delta;h,v]$ is a Sylow $p$-subgroup of the Frobenius group $F\defeq \I^\#_{\Pi\setminus H}[\Delta;h,v]$, it is a subgroup of some Sylow $p$-group $S'$ in the group $\I^\#_{\Pi\setminus H}[\Delta]$. Since $K$ is a maximal $p$-subgroup of the group $F$, $K=S'\cap F$. Sylow's Theorem~\ref{t:Sylow} ensures that the Sylow $p$-groups $S$ and $S'$ are conjugated in the group $\I^\#_{\Pi\setminus H}[\Delta]$ and hence $S'=g^{-1}Sg$ for some $g\in \I^\#_{\Pi\setminus H}[\Delta]$. Then $K=S'\cap F=g^{-1}Sg\cap F=g^{-1}(S\cap F)g$. Since the kernel $K$ of the Frobenius group $F$ is a characteristic subgroup of $F$, $K=g^{-1}(S\cap F)g=S\cap F$.
\end{proof}

\begin{claim}\label{cl:S-orbit} There exists a point $o\in \Delta\setminus\{\delta\}$, whose orbit $S(o)\defeq\{s(o):s\in S\}$ has cardinality $p^n$.
\end{claim}

\begin{proof} Observe that for every $o\in \Delta\setminus\{\delta\}$, the cardinal $|S(o)|$ divides the order of the group $S$ and hence $|S(o)|$ is a power of $p$. By Claim~\ref{cl:S-kernel}, for all distinct points $h,v\in H\setminus\{\delta\}$, the kernel $K$ of the Frobenius group $\I^\#_{\Pi\setminus H}[\Delta;h,v]$ is a subgroup of the group $S$, which implies $|S(o)|\ge |K(o)|=p^n$ for all $o\in\Delta\setminus\{\delta\}$. To derive a contradiction, assume that $|S(o)|>p^n$ for every $o\in\Delta\setminus\{\delta\}$. Then the set $\Delta\setminus\{\delta\}$ is the disjoint union of $S$-orbits whose cardinality is divided by $p^{n+1}$, which implies that $p^nq=|\Delta\setminus\{\delta\}|$ is divided by $p^{n+1}$ which is a contradiction showing that for some $o\in\Delta\setminus\{\delta\}$ the orbit $S(o)$ has cardinality $p^n$.
\end{proof}

By Claim~\ref{cl:S-orbit}, there exists a point $o\in\Delta\setminus\{\delta\}$ whose  orbit $S(o)\defeq \{s(o):s\in S\}\subseteq\Delta\setminus\{\delta\}$ has cardinality $p^n>1$. So we can fix an element $e\in S(o)\setminus\{o\}$. Given any distinct points $h,v\in H\setminus\{\delta\}$, consider the affine base $uow$ in the affine plane $\Pi\setminus H$ consisting of the points $u\in \Aline oh\cap\Aline ev$ and $w\in \Aline ov\cap\Aline eh$. We claim that the element $e$ has order $p$ in the plus loop $(\Delta\setminus\{\delta\},+)$ of the ternar of the based affine plane $(\Pi\setminus H,uow)$.

Consider the Frobenius group $\I^\#_{\Pi\setminus H}[\Delta;h,v]$ and its kernel $K\defeq K_{\Pi\setminus  H}[\Delta;h,v]$. Claim~\ref{cl:S-kernel}  ensures that $S\cap\I^\#_\Pi[\Delta;h,v]=K$ and hence $K(o)\subseteq S(o)$. Since the action of the group $K$ on the set $\Delta\setminus\{\delta\}$ has no fixed points, $|K(o)|=|K|=p^n=|S(o)|$ and hence $e\in S(o)=K(o)$. Consider the line $L\defeq\Aline w\delta$ and observe that $e=h_{\Delta,L}v_{L,\Delta}(o)$. Since the action of the group $\I^\#_{\Pi\setminus H}[\Delta;h,v]$ on the set $\Delta\setminus\{\delta\}$ has no fixed points, $e=h_{\Delta,L}v_{L,\Delta}(o)\in K(o)$ implies $h_{\Delta,L}v_{L,\Delta}\in K$ and hence the permutation $h_{\Delta,L}v_{L,\Delta}$ has order $p$, which implies that the element $e$ has order $p$ in the plus loop $(\Delta\setminus\{\delta\},+)$ of the ternar of the based affine plane $(\Pi\setminus H,uow)$. By Theorem~\ref{t:invertible-add<=>4-Pappian}, the invertible-add projective plane $\Pi$ is $4$-Pappian and hence the closure  of the $4$-element set $\{o,e,h,v\}$ is a Pappian subliner of the projective plane $\Pi$. Then the element $o_{H,L}e_{L,H}$ of the group $\I^\#_{\Pi\setminus \Delta}[H;o,e]$ has order $p$ and maps $h$ onto $v$. This implies that the $p$-group $K_{\Pi\setminus \Delta}[H;o,e]$ acts transitively on the set $H\setminus\{\delta\}$ which is impossible as $|K_{\Pi\setminus\Delta}[H;o,e]|=p^n<p^nq=|H\setminus\{\delta\}|$. This contradiction completes the proof that the projective plane $\Pi$ has prime-power order.

Now we can deduce from Theorem~\ref{t:Gleason56} that the projective plane $\Pi$ is Desarguesian. If $\|X\|\ge 3$ or $|X|_2\le 3$, then $X$ is Desarguesian, by Theorems~\ref{t:Desargues-projective} and Corollary~\ref{c:4-Pappian}. So, assume that $\|X\|=3$ and $|X|_2\ge 4$, which means that $X$ is a $4$-long projective plane. The Desarguesianity of $X$ will follow from the Gleason Theorem~\ref{t:Gleason56} as soon as we check that for every line $L\subseteq X$ and point $p\in L$, the group $\Aut_{p,L}(X)$ of automorphisms of $X$ with centre $p$ and axis $L$ is not trivial. Consider the affine plane $A\defeq X\setminus L$ and the spread of parallel lines $\boldsymbol \delta\defeq\{\Aline xp\setminus\{p\}:x\in A\}$ in $A$. Since the projective plane $X$ is associative-plus and $|X|_2-1$ is a prime power, the affine plane $A$ is associative-plus and its order $|A|_2=|X|_2-1$ is a prime power.  By Theorem~\ref{t:add-ass=>partialT}, the affine space $A$ is $\partial$-translation, which implies the existence of a non-trivial translation $T:A\to A$ such that $\{\Aline xy:xy\in T\}=\boldsymbol \delta$. By Theorem~\ref{t:extend-isomorphism-to-completions}, the translation $T$ entends to a (non-trivial) automorphism $\bar T\in\Aut_{p,L}(X)$, witnessing that the group $\Aut_{p,L}(X)$ is not trivial. The Gleason Theorem~\ref{t:Gleason56} ensures that the projective plane $X$ is Desarguesian.
\end{proof}

Now we can prove the main result of this section.

\begin{theorems}\label{t:proj-finite<=>} For any  projective space $X$ of finite order, the following conditions are equivalent:
\begin{enumerate}
\item $X$ is Pappian;
\item $X$ is Desarguesian;
\item $X$ is Moufang;
\item $X$ is commutative-dot;
\item $X$ is associative-dot;
\item $X$ is inversive-dot;
\item $X$ is invertible-dot;
\item $X$ is commutative-add;
\item $X$ is associative-add;
\item $X$ is commutative-plus;
\item $X$ is associative-plus;
\item $X$ is commutative-puls;
\item $X$ is associative-puls.
\end{enumerate}
If the order of $X$ is a Moufang number, then the conditions \textup{(1)--(13)} are equivalent to the conditions:
\begin{enumerate}
\item[(14)] $X$ is inversive-add;
\item[(15)] $X$ is inversive-plus;
\item[(16)] $X$ is inversive-puls.
\end{enumerate}
\end{theorems}

\begin{proof} The equivalences of the conditions (1)--(6) follow from Theorem~\ref{t:finite-shear<=>Pappian}. 
 The implications $(1)\Ra(8)\Ra(9)$ and $(10)\Ra(11)$ follow from Theorem~\ref{t:p.commutative} and Corollary~\ref{c:p:com-plus=>add-plus}.  The implications $(8)\Ra(10)$ and $(9)\Ra(11)$ are obvious. The implication $(11)\Ra(2)$ follows from Theorem~\ref{t:Kegel-Luneburg-main}. The implication $(6)\Ra(7)$ is obvious, and  $(7)\Ra(12)\Ra(13)$ follows from Theorems~\ref{t:p-invertible-dot<=>} and \ref{t:com-puls=>ass-puls}. It remains to prove that $(13)\Ra(2)$.

Assume that a projective space $X$ of finite order if associative-puls. If $\|X\|>3$, then $X$ is Desarguesian, by Theorem~\ref{t:proaffine-Desarguesian}. So, assume that $\|X\|=3$, which means that $X$ is a projective plane. By Theorem~\ref{t:plus-puls-duality}, the associativity-puls of the projective plane $X$ implies the associativity-plus of the dual projective plane $X^*$. Since the projective plane $X$ has finite order,  it is finite and so is its dual projective plane $X^*$. Applying Theorem~\ref{t:Kegel-Luneburg-main}, we conclude that the finite projective plane $X^*$ is Desarguesian and hence Pappian. By Theorem~\ref{t:Pappian-self-dual}, the projective plane $X$ is Pappian, too.
\smallskip

If the order of the projective space $X$ is prime, then the conditions (1)--(13) are equivalent to the conditions (14)--(16),  by Proposition~\ref{p:Moufang=>inversive<=>ass}.
\end{proof}

\begin{theorem}\label{t:Fano<=>inversive} For a projective space $X$ of finite order the following conditions are equivalent:
\begin{enumerate}
\item $X$ is Fano;
\item $X$ is inversive-plus and the order of $X$ is a power of $2$;
\item $X$ is associative-plus and the order of $X$ is even; 
\item $X$ is commutative-add and the order of $X$ is even;
\item $X$ is associative-puls and the order of $X$ is even;
\item $X$ is inversive-puls and the order of $X$ is a power of $2$.
\end{enumerate}
\end{theorem}

\begin{proof} By definition, a projective space $X$ is $3$-long. By Corollary~\ref{c:Avogadro-projective}, the $3$-long projective space $X$ is $2$-balanced. If $|X|_2=3$, then the projective space $X$ has order $2$, is Fano and Pappian, by Example~\ref{ex:Fano3} and Corollary~\ref{c:p5-Pappian}. In this case the conditions (1)--(6) hold and hence are equivalent. So, assume that $|X|_2>3$.
\smallskip

$(1)\Ra(4)$ Assume that $X$ is Fano. To prove that $X$ is commutative-add, 
take any ternar $R$ of $X$ and find a plane $\Pi\subseteq X$ and a projective base $uowe$ in $\Pi$ whose ternar is isomorphic to $R$. Find unique points $h\in \Aline ou\cap\Aline we$ and $v\in\Aline ow\cap\Aline ue$ and consider the affine plane $\Pi\setminus\Aline hv$. By definition, the ternar of the based projective plane $(\Pi,uowe)$ coincides with the ternar of the based affine plane $(\Pi\setminus\Aline hv,uow)$. Since $X$ is Fano, so is the plane $\Pi$ in $X$. By Theorem~\ref{t:Fano<=>}, the Fano projective plane $\Pi$ is everywhere Boolean, which implies that the affine plane $\Pi\setminus\Aline hv$ is Boolean. By Theorem~\ref{t:Boolean<=>forder},  the cardinal $|\Delta|=|\Pi\setminus\Aline hv|_2=|X|_2-1$ is even and the ternar $\Delta=\Aline oe\setminus\Aline hv$ of the affine plane $\Pi\setminus\Aline hv$ is commutative-add. Then the ternar $R$ is commutative-add, being isomorphic to the ternar $\Delta$.
\smallskip

The implications $(4)\Ra(3,5)$ follow from Theorems~\ref{t:add-com=>add-ass} and \ref{t:com-puls=>ass-puls}.
\smallskip

$(3)\Ra(2)$ and $(5)\Ra(6)$. If $X$ is associative-plus or associative-puls, then $X$ is Desarguesian (by Theorem~\ref{t:proj-finite<=>}) and hence $|X|_2=p^n$ for some prime number $p$ and some $n\in\IN$, according to Theorem~\ref{t:cardinality-pD}. If $|X|_2$ is even, then $p=2$ and $|X|_2=2^n$.
\smallskip

$(2\vee 6)\Ra(1)$ Assume that $X$ is inversive-plus or inversive-puls, and $|X|_2$ is a power of $2$. We shall prove that $X$ is everywhere Boolean. Take any hyperplane $H\subseteq X$ and consider the affine space $X\setminus H$. Since $X$ is inversive-plus or inversive-puls, so is the affine space $X\setminus H$, according to Theorems~\ref{t:p+inversive} and \ref{t:p-inversive-puls}. Applying Theorem~\ref{t:Boolean<=>forder}, we conclude that the affine space $X\setminus H$ is Boolean, witnessing that the projective space $X$ is everywhere Boolean. By Theorem~\ref{t:Fano<=>}, the projective space $X$ is Fano.
\end{proof}

\begin{problem} Is every invertible-add projective plane of finite even order Fano?
\end{problem}

The following proposition gives a partial answer to this problem.

\begin{proposition}\label{p:4Fano} If $X$ is an invertible-add projective plane of finite even order, then every triangle $abc$ in $X$ can be completed to a Fano quadrangle $abcd$.
\end{proposition}

\begin{proof} Let $m\defeq |X|_2-1$ be the order of the invertible-add projective plane. Consider the set $Y\defeq X\setminus(\Aline ab\cup\Aline bc\cup\Aline ac)$ and observe that $$|Y|=|X|-(|\Aline ab|+|\Aline bc|+|\Aline ac|-|\{a,b,c\}|)=m^2+m+1-3m=(m-1)^2.$$ By Theorem~\ref{t:invertible-add<=>4-Pappian}, the invertible-add projective plane $X$ is $4$-Pappian. Then for every $x\in Y$, the closure $Q_x$ of the quadrangle $abcx$ in $X$ is a Pappian projective plane of prime order $p_x$. Observe that $|Q_x|=p_x^2+p_x+1$ and $|Q_x\cap Y|=|Q_x\setminus(\Aline ab\cup\Aline bc\cup\Aline ac)|=(p_x-1)^2.$

We claim that for any $x,y\in Y$, either $Q_x=Q_y$ or $Q_x\cap Q_y\cap Y=\varnothing$. Indeed, assuming that $Q_x\cap Q_y\cap Y$ contains some point $z$, we conclude that $Q_z\subseteq Q_x\cap Q_y$. Taking into account that Pappian projective planes of prime order do not contain proper projective subliners of rank $3$, we conclude that $Q_x=Q_z=Q_y$. Therefore, $\{Q_x\cap Y:x\in Y\}$ is a partition of the set $Y$, and there exists a set $Z\subseteq Y$ such that for every $y\in Y$
there exists a unique point $z\in Z$ such that $Q_y\cap Y=Q_z\cap Y$. Then 
$$(m-1)^2=|Y|=\sum_{z\in Z}|Q_x\cap Y|=\sum_{z\in Z}(p_z-1)^2.$$
Since the order $m=|X|_2-1$ of $X$ is even, the number $(m-1)^2$ is odd and hence there exists $x\in Z$ such that $p_x-1$ is odd. The unique prime number with this property is $2$. Therefore, the projective plane $Q_x$ has order $p_x=2$ and hence $Q_x$ is Fano and so is the quadrangle $abcx$ generating the projective Fano plane $Q_x$.   
\end{proof}



\section{Interplay between various properties of projective liners}\label{s:diagram-proj}

By Theorems~\ref{t:p-invertible-add<=>}, \ref{t:invertible-add<=>4-Pappian}, \ref{t:p-commutative-add<=>}, \ref{t:p-invertible-dot<=>}, \ref{t:p.inversive}, \ref{t:p.associative}, \ref{t:p.commutative}, and \ref{t:proj-finite<=>}, 
for any projective space the following implications hold:
$$
\xymatrix@C=20pt@R=22pt{
\mbox{commutative-dot}\ar@{<=>}[r]\ar@{=>}[d]&\mbox{Pappian}\ar@{=>}[d]\\
\mbox{associative-dot}\ar@{=>}[d]\ar@{<=>}[r]&\mbox{Desarguesian}\ar@{=>}[d]\ar@/_10pt/_{\mbox{\scriptsize\em finite}}[u]\\
\mbox{inversive-dot}\ar@{=>}[d]\ar@{<=>}[r]&\mbox{Moufang}\ar@{=>}[d]\ar@/_10pt/_{\mbox{\scriptsize\em finite}}[u]\\
\mbox{invertible-dot}\ar@{<=>}[d]&\mbox{Fano or Moufang}\ar@{<=>}[d]\ar@{=>}[l]\ar@/_10pt/_{\mbox{\scriptsize\em finite}}[u]\\
\mbox{commutative-puls}\ar@{=>}[d]\ar@/^10pt/^{\mbox{\scriptsize\em finite}}[r]&\mbox{commutative-add}\ar@{=>}[d]\ar@{=>}[l]\ar@{=>}[r]&\mbox{commutative-plus}\ar@{=>}[d]\ar@/_10pt/_{\mbox{\scriptsize\em finite}}[l]\\
\mbox{associative-puls}\ar@{=>}[d]\ar@/_10pt/_{\mbox{\scriptsize\em finite}}[u]&\mbox{associative-add}\ar@{=>}[r]\ar@{=>}[l]\ar@{=>}[d]&\mbox{associative-plus}\ar@{=>}[d]\ar@/_10pt/_{\mbox{\scriptsize\em finite}}[u]\\
\mbox{inversive-puls}\ar@{=>}[d]\ar@/_10pt/_{\mbox{\scriptsize\em Moufang}\atop\mbox{\scriptsize\em order}}[u]&\mbox{inversive-add}\ar@{=>}[r]\ar@{=>}[l]\ar@{=>}[d]&\mbox{inversive-plus}\ar@{=>}[d]\ar@/_10pt/_{\mbox{\scriptsize\em Moufang}\atop\mbox{\scriptsize\em order}}[u]\\
\mbox{invertible-puls}\ar@{<=>}[r]&\mbox{invertible-add}\ar@{<=>}[r]&\mbox{invertible-plus}\ar@{=>}[d]\\
&\mbox{$4$-Pappian}\ar@{<=>}[u]&\mbox{tri-Desarguesian}\ar@{<=>}[l]
}
$$

This diagram motivates the following open problems.

\begin{problem} Is every (finite) invertible-add projective space inversive-add?
\end{problem}

\begin{problem} Is every (finite) inversive-add projective space associative-add?
\end{problem}

\chapter{More incidences in the Desargues and Pappus Axioms}

In this chapter we consider several weaker versions of the Desargues and Pappus Axioms with more incidences required in the hypotheses of those axioms. We shall see that for projective liners those weaker versions of Desargues or Pappus Axioms are equivalent to the algebraic conditions (like commutativity-add or invertibility-add), considered in Chapter~\ref{ch:AMP}. 

\section{More incidences in the Desargues Axiom}

In this section we shall study some weaker versions of the Desargues Axiom with some additional incidences. Two triangles $abc$ and $a'b'c'$ in a liner are called \defterm{coplanar} if the points $a,b,c,a',b',c'$ belong to some plane.

\begin{definition}\label{d:D123} A liner $X$ is defined to be 
\begin{itemize}
\item[$\mathsf{(D)}$] \defterm{Desarguesian} if for any centrally perspective coplanar disjoint triangles $abc$ and $a'b'c'$ in $X$, the set $(\Aline ab\cap\Aline {a'}{b'})\cup(\Aline bc\cap\Aline{b'}{c'})\cup(\Aline ac\cap\Aline{a'}{c'})$ has rank $0$ or $2$;
\item[$\mathsf{(D_1)}$] \defterm{uno-Desarguesian} if for any centrally perspective disjoint triangles $abc$ and $a'b'c'$ in $X$ with $b'\in\Aline ac$, the set $(\Aline ab\cap\Aline {a'}{b'})\cup(\Aline bc\cap\Aline{b'}{c'})\cup(\Aline ac\cap\Aline{a'}{c'})$ has rank $0$ or $2$;
\item[$\mathsf{(D_2)}$] \defterm{duo-Desarguesian} if for any centrally perspective disjoint triangles $abc$ and $a'b'c'$ in $X$ with $a'\in\Aline bc$ and $c'\in \Aline ab$, the set $(\Aline ab\cap\Aline {a'}{b'})\cup(\Aline bc\cap\Aline{b'}{c'})\cup(\Aline ac\cap\Aline{a'}{c'})$ has rank $0$ or $2$;
\item[$\mathsf{(D_{II})}$] \defterm{bi-Desarguesian} if for any centrally perspective disjoint triangles $abc$ and $a'b'c'$ in $X$ with $b'\in\Aline ac$ and $b\in \Aline {a'}{c'}$, the set $(\Aline ab\cap\Aline {a'}{b'})\cup(\Aline bc\cap\Aline{b'}{c'})\cup(\Aline ac\cap\Aline{a'}{c'})$ has rank $0$ or $2$;
\item[$\mathsf{(D_3)}$] \defterm{tri-Desarguesian} if for any centrally perspective disjoint triangles $abc$ and $a'b'c'$ in $X$ with $a'\in\Aline bc$, $b'\in \Aline ac$, and $c'\in \Aline ab$, the set $(\Aline ab\cap\Aline {a'}{b'})\cup(\Aline bc\cap\Aline{b'}{c'})\cup(\Aline ac\cap\Aline{a'}{c'})$ has rank $0$ or $2$.
\end{itemize}
\end{definition}

It is easy to see that for every liner, we have the implications:
$$
\xymatrix{
\mathsf{(D)}\ar@{=>}[r]&
\mathsf{(D_1)}\ar@{=>}[r]&\mathsf{(D_2\wedge D_{II})}\ar@{=>}[r]&
\mathsf{(D_2)}\ar@{=>}[r]&
\mathsf{(D_3)}
}
$$


By Definitions~\ref{d:D123} and \ref{d:Moufang}, a liner is uno-Desarguesian if and only if it is Moufang. Moufang liners were thoroughly studied in Chapter~\ref{ch:Moufang}.

\begin{theorem}\label{t:duo-Desarg<=>} For a $3$-long projective liner $Y$, the following conditions are equivalent:
\begin{enumerate}
\item $Y$ is duo-Desarguesian;
\item $Y$ is bi-Desarguesian;
\item $Y$ is everywhere uno-Thalesian;
\item $Y$ is Fano or Moufang;
\item $Y$ is commutative-add;
\item $Y$ is commutative-plus and invertible-dot.
\end{enumerate}
\end{theorem}

\begin{proof} $(1)\Ra(2)$ Assume that the projective liner $Y$ is duo-Desarguesian.
To prove that $Y$ is bi-Desarguesian, take any distinct centrally perspective  triangles $abc$ and $a'b'c'$ in $Y$ with $b\in\Aline {a'}{c'}$ and $b'\in \Aline{a}{c}$.
Let $o\in\Aline a{a'}\cap\Aline b{b'}\cap\Aline c{c'}$ be the perspector of the triangles $abc$ and $a'b'c'$. It follows from $b'\in \Aline ac$ that the points $a,a',b',b',c,c'$ belong to the plane $\overline{\{o,a,c\}}$. By the projectivity of the liner $Y$, there exist unique points $x\in \Aline ab\cap\Aline{a'}{b'}$ and $y\in\Aline bc\cap\Aline{b'}{c'}$ and $z\in \Aline ac\cap\Aline{a'}{c'}$. We have to prove that the points $x,y,z$ are collinear. By the projectivity of $Y$, there exists a unique point $x'\in \Aline yz\cap\Aline{a'}{b'}$. Then the triangles $czc'$ and $bx'b'$ are perspective from the point $y$. Since $b\in \Aline {a'}{c'}=\Aline z{c'}$ and $b'\in\Aline ac=\Aline cz$, we can apply the duo-Desarguesian property of $Y$ and conclude that the set 
$$(\Aline cz\cap\Aline b{x'})\cup(\Aline z{c'}\cap\Aline{x'}b)\cup(\Aline c{c'}\cap\Aline b{b'})=(\Aline cz\cap\Aline b{x'})\cup(\Aline {a'}{c'}\cap\Aline{a'}{b'})\cup\{o\}=
(\Aline cz\cap\Aline b{x'})\cup\{a',o\}$$has rank $2$. Then $\varnothing\ne\Aline cz\cap\Aline b{x'}\subseteq \Aline o{a'}=\Aline oa$ and $\varnothing\ne \Aline oa\cap\Aline cz\cap \Aline b{x'}=\{a\}\cap\Aline b{x'}$, which implies $x'\in \Aline ab\cap \Aline {a'}{b'}=\{x\}$ and finally, $x=x'\in\Aline yz$.  
\smallskip

$(2)\Ra(3)$ Assume that $Y$ is bi-Desarguesian. We have to prove that for every hyperplane $H\subseteq Y$, the liner $X\defeq Y\setminus H$ is uno-Thalesian. Given any plane $\Pi\subseteq X$, disjoint lines $A,B,C\subseteq \Pi$ and distinct points $a,a'\in A$, $b,b'\in B$, $c,c'\in C$ with $b'\in\Aline ac$ and $\Aline ab\cap X\cap\Aline {a'}{b'}=\varnothing=\Aline bc\cap X\cap\Aline{b'}{c'}$, we should check that $\Aline ac\cap X\cap\Aline{a'}{c'}=\varnothing$. 

Since $Y$ is a projective liner, the flat hull $\overline \Pi$ of the plane $\Pi$ in $Y$ is a projective plane. Given a line $L$ in $X$, let $\overline L$ be the flat hull of $L$ in the projective liner $Y$. By Corollary~\ref{c:line-meets-hyperplane}, $\overline L\cap H\ne\varnothing$. Since $H$ is flat, $\overline L\cap H$ is a singleton. By the projectivity of the plane $\Pi$, the flat hulls $\overline A,\overline B,\overline C$ of the disjoint lines $A,B,C$ have a unique common point $o\in\overline A\cap\overline B\cap\overline C\cap H$. Also, by the projectivity of the plane $\overline\Pi$, there exist points $x\in\Aline ab\cap\Aline{a'}{b'}$, $y\in \Aline bc\cap\Aline {b'}{c'}$ and $z\in \Aline ac\cap\Aline {a'}{c'}$. The assumption $\Aline ab\cap X\cap\Aline {a'}{b'}=\varnothing=\Aline bc\cap X\cap\Aline{b'}{c'}$ implies $\{x,y\}\subseteq H$. 

If $x=y$, then $\{a,c\}\subseteq \Aline xb=\Aline yb$ and $\{a',c'\}\subseteq \Aline x{b'}=\Aline y{b'}$. In this case $\Aline ac\cap X\cap\Aline {a'}{c'}=\Aline ab\cap X\cap\Aline{a'}{b'}=\varnothing$ and we are done. So, assume that $x\ne y$ and hence $\Aline xy=\overline \Pi\cap H$, by the rankedness of the proejctive plane $\overline\Pi$.

Observe that the triangles $aoc$ and $xb'y$ are perspective from the point $b$. Since $o\in \overline\Pi\cap H=\Aline xy$ and $b'\in \Aline ac$, we can apply the bi-Desarguesian property of $Y$ and conclude that the set $$(\Aline ao\cap\Aline x{b'})\cup(\Aline oc\cap\Aline {b'}y)\cup(\Aline ac\cap\Aline xy)=\{a',c'\}\cup(\Aline ac\cap\Aline xy)$$ has rank $2$ and hence $\varnothing\ne\Aline ac\cap\Aline xy\subseteq \Aline {a'}{c'}$, $\varnothing\ne \Aline{a'}{c'}\cap\Aline ac\cap\Aline xy=\{z\}\cap\Aline xy$, and finally, $z\in \Aline xy\subseteq H$ and $\Aline ac\cap\Aline {a'}{c'}\cap X=\{z\}\cap X\subseteq H\cap X=\varnothing$.
\smallskip

$(3)\Ra(4)$ Assume that $Y$ is everywhere uno-Thalesian. Assuming that $Y$ is not Fano, we shall prove that $Y$ is Moufang. If $\|Y\|\ne 3$ or $|Y|_2=3$, then $Y$ is Moufang, by Theorem~\ref{t:proaffine-Desarguesian} and Propositions~\ref{p:Steiner+projective=>Desargues} and \ref{p:Desarg=>Moufang}. So, assume that $\|Y\|=3$ and $|Y|_2\ge 4$. Since $Y$ is not Fano, there exists a quadrangle $abcd$ whose diagonal points $x\in\Aline ab\cap\Aline cd$, $y\in \Aline ac\cap\Aline bd$ and $z\in \Aline ad\cap\Aline bc$ are not collinear. Since $Y$ is a plane, the lines $H_1\defeq\Aline xz$ and $H_2\defeq \Aline yz$ are hyperplanes in $Y$. Since $Y$ is everywhere uno-Thalesian, for every $i\in\{1,2\}$, the affine plane $Y\setminus H_i$ is uno-Thalesian and hence Boolean or Thalesian, by Theorem~\ref{t:commutative-add<=>}. Observe that $abcd$ is a non-Boolean parallelogram in the affine plane $Y\setminus H_1=Y\setminus\Aline xz$, witnessing that the affine plane $Y\setminus H_1$ is not Boolean and hence it is Thalesian. On the other hand, $acbd$ is a non-Boolean parallelogram in the affine plane $Y\setminus H_2=Y\setminus \Aline yz$, witnessing that $Y\setminus H_2$ is not Boolean and hence $Y\setminus H_2$ is Thalesian. Therefore, for the distinct lines $H_1,H_2$ in $Y$, the affine liners $Y\setminus H_1$ and $Y\setminus H_2$ are Thalesian. By Theorem~\ref{t:Skornyakov-San-Soucie}, the projective plane $Y$ is Moufang.
\smallskip

$(4)\Ra(1)$ Assume that the projective liner $Y$ is Fano or Moufang. If $Y$ is Moufang, then it is uno-Desarguesian and duo-Desarguesian. So, assume that $Y$ is Fano. To prove that $Y$ is duo-Desarguesian, take any distinct centrally perspective triangles $abc$ and $a'b'c'$ in $X$ such that $a'\in \Aline bc$ and $c'\in \Aline ab$.  Let $o\in\Aline a{a'}\cap \Aline b{b'}\cap\Aline c{c'}$ be the perspector of the triangles $abc$ and $a'b'c'$. It follows from $c'\in\Aline ab$ that the plane $\overline{\{o,a,b\}}$ contains the points $a,a',b,b',c',c$. By the projectivity of $Y$, there exist unique points $x\in \Aline ab\cap\Aline {a'}{b'}$, $y\in \Aline bc\cap\Aline{b'}{c'}$ and $z\in \Aline ac\cap\Aline{a'}{c'}$. Since $Y$ is Fano, the diagonal points $o,b,z$ of the quadrangle $aa'c'c$ are collinear. Therefore, $z\in \Aline ob\cap\Aline {a'}{c'}$. It follows from $a'\ne a\in\Aline c'b$ that $b\notin \Aline {a'}{c'}$ and hence $bb'a'c'$ is a quadrangle. Since $Y$ is Fano, the diaginal points $z\in \Aline b{b'}\cap\Aline {a'}{c'}$, $y\in\Aline bc\cap\Aline {b'}{c'} =\Aline b{a'}\cap\Aline {b'}{c'}$ and $x\in \Aline ab\cap\Aline{a'}{b'}=\Aline b{c'}\cap\Aline {b'}{a'}$ of the quadrangle $bb'a'c'$ are collinear. Therefore,
$$\|(\Aline ab\cap\Aline{a'}{b'})\cup(\Aline bc\cap\Aline{b'}{c'})\cup(\Aline ac\cap\Aline{a'}{c'})\|=\|\{x,y,z\}\|=2,$$witnessing that the projective liner $Y$ is duo-Desarguesian.
\smallskip

The equivalences $(4)\Leftrightarrow(5)\Leftrightarrow(6)$ follow from Theorems~\ref{t:p-commutative-add<=>} and \ref{t:p-invertible-dot<=>}.
\end{proof}
\pagebreak

\begin{theorem}\label{t:tri-Desarg<=>} For a projective liner $X$, the following conditions are equivalent:
\begin{enumerate}
\item $X$ is tri-Desarguesian;
\item $X$ is invertible-add;
\item $X$ is invertible-plus;
\item $X$ is invertible-puls;
\item $X$ is $4$-Pappian.
\end{enumerate}
\end{theorem}

\begin{proof} The equivalence of the conditions (2)--(5) has been proved in Theorems~\ref{t:p-invertible-add<=>} and \ref{t:invertible-add<=>4-Pappian}.  So, it suffices to prove that $(1)\Leftrightarrow(4)$.
\smallskip

$(1)\Ra(4)$ Assume that a projective liner $X$ is tri-Desarguesian. To prove that $X$ is invertible-puls, take any distinct lines $L,L'$ in $X$ and distinct points $a,b,c\in L\setminus L'$, $a',b',c'\in L'\setminus L$, $h\in L\cap L'$, $v\in \Aline a{a'}\cap\Aline b{b'}\cap\Aline c{c'}$ with $\Aline a{b'}\cap\Aline b{c'}\subseteq \Aline vh$. We have to prove that $\Aline {a'}b\cap\Aline {b'}c\subseteq \Aline vh$. By the projectivity of $X$, there exist unique points $x\in \Aline a{b'}\cap\Aline b{c'}$, $y\in \Aline a{b'}\cap\Aline b{c'}$, $\alpha\in \Aline b{c'}\cap\Aline {b'}c$, $\beta\in \Aline b{b'}\cap\Aline a{c'}$, and $\gamma\in \Aline c{c'}\cap\Aline {a'}b$.
\smallskip

\begin{picture}(300,280)(-100,-20)
\linethickness{=0.6pt}
\put(-60,0){\color{teal}\line(1,0){300}}
\put(0,0){\color{red}\line(0,1){240}}
\put(0,0){\color{cyan}\line(1,1){120}}
\put(0,240){\line(1,-1){240}}
\put(0,240){\color{red}\line(1,-6){40}}
\put(240,0){\color{teal}\line(-6,1){240}}
\put(120,120){\color{cyan}\line(-3,-2){180}}
\put(120,120){\color{cyan}\line(-2,-3){80}}
\put(-60,0){\color{blue}\line(11,4){220}}
\put(0,0){\color{blue}\line(2,1){160}}

\put(0,40){\line(1,-1){40}}
\put(-60,0){\line(4,1){120}}
\put(0,0){\line(1,2){30}}

\put(34.5,34.21){\circle*{3}}
\put(34,38){$b$}
\put(-60,0){\circle*{3}}
\put(-63,-10){$c'$}
\put(0,0){\circle*{3}}
\put(-3,-10){$b'$}
\put(40,0){\circle*{3}}
\put(37,-10){$a'$}
\put(240,0){\color{teal}\circle*{3}}
\put(243,-3){\color{teal}$h$}
\put(0,240){\color{red}\circle*{3}}
\put(-2,244){\color{red}$y$}
\put(0,40){\circle*{3}}
\put(-7,41){$c$}
\put(120,120){\color{cyan}\circle*{3}}
\put(122,121){\color{cyan}$v$}
\put(160,80){\color{blue}\circle*{3}}
\put(163,81){\color{blue}$x$}
\put(60,30){\circle*{3}}
\put(57,33){$a$}
\put(125,-10){\color{teal}$L'$}
\put(130,21){\color{teal}$L$}
\put(20,20){\circle*{2.5}}
\put(15.6,10){$\beta$}
\put(30,60){\circle*{2.5}}
\put(30,65){$\gamma$}
\put(13.3,26.7){\circle*{2.5}}
\put(9.5,31){$\alpha$}

\end{picture}
\smallskip

Observe that the triangles $bb'h$ and $c'av$ are perspective from the point $x$. Since $c'\in \Aline {b'}h$, $a\in \Aline bh$ and $v\in \Aline b{b'}$, we can apply the tri-Desarguesian property of $X$ and conclude that the points $\beta\in \Aline b{b'}\cap\Aline {c'}a$, $c\in\Aline bh\cap\Aline{c'}v$ and $a'\in \Aline{b'}h\cap \Aline av$ are collinear. 

Observe that the triangles $cv\beta$ and $ba'c'$ are perspective from the point $a$. Since $b\in \Aline v\beta$, $a'\in \Aline c\beta$ and $c'\in \Aline cv$, we can apply the tri-Desarguesian property of $X$ and conclude that the points $\gamma\in \Aline cv\cap\Aline b{a'}$, $\alpha\in\Aline c\beta\cap\Aline b{c'}$ and $b'\in \Aline v\beta\cap \Aline {a'}{c'}$ are collinear. 

Finally, observe that the triangles $c'\gamma a'$ and $bb'c$ are perspective from the point $\alpha$. Since $b\in \Aline \gamma{a'}$, $b'\in \Aline{c'}{a'}$ and $c\in \Aline{c'}\gamma$, we can apply the tri-Desarguesian property of $X$ and conclude that the points $v\in \Aline {c'}\gamma\cap\Aline b{b'}$, $h\in \Aline {c'}{a'}\cap\Aline  bc$ and $y\in \Aline \gamma{a'}\cap\Aline{b'}c$ are collinear and hence 
$y\in \Aline hv$.
\smallskip

$(4)\Ra(1)$ Assume that the projective liner $X$ is inversible-puls.
To prove that $X$ is tri-Desarguesian, take any disjoint centrally perspective triangles $abc$ and $a'b'c'$ such that $a'\in \Aline bc$, $b'\in \Aline ab$ and $c'\in \Aline ab$. We have to prove that the points $x\in \Aline ab\cap \Aline {a'}{b'}$, $y\in \Aline {b'}{c'}$ and $z\in \Aline ac\cap\Aline{a'}{c'}$ are collinear. By the projectivity of $X$, there exists unique points $\alpha\in \Aline cx\cap\Aline bz$ and $\gamma\in\Aline ay\cap \Aline bz$.
\smallskip

\begin{picture}(300,265)(-100,-15)
{\linethickness{=1.2pt}
\put(-60,0){\color{teal}\line(1,0){100}}
\put(0,0){\color{cyan}\line(0,1){40}}
\put(0,0){\color{blue}\line(1,1){34.5}}
\put(40,0){\color{cyan}\line(-1,6){10}}
\put(0,40){\color{teal}\line(6,-1){34.5}}
\put(-60,0){\color{blue}\line(3,2){90}}
}

\linethickness{=0.5pt}
\put(-60,0){\color{teal}\line(1,0){300}}
\put(0,0){\color{cyan}\line(0,1){240}}
\put(0,0){\color{blue}\line(1,1){120}}
\put(0,240){\line(1,-1){240}}
\put(0,240){\color{cyan}\line(1,-6){40}}
\put(240,0){\color{teal}\line(-6,1){240}}
\put(120,120){\color{blue}\line(-3,-2){180}}
\put(-60,0){\line(11,4){132}}
\put(-60,0){\color{cyan}\line(1,4){60}}
\put(240,0){\color{teal}\line(-7,2){280}}
\put(40,0){\color{blue}\line(2,3){80}}
\put(40,0){\line(-1,1){80}}
\put(0,0){\line(1,2){30}}

\put(72,48){\circle*{2.5}}
\put(75,48){$\alpha$}
\put(-40,80){\circle*{2.5}}
\put(-48,80){$\gamma$}
\put(34.5,34.21){\circle*{3}}
\put(34,39){$a'$}
\put(-60,0){\circle*{3}}
\put(-63,-10){$a$}
\put(0,0){\circle*{3}}
\put(-3,-10){$b'$}
\put(40,0){\circle*{3}}
\put(37,-9){$c$}
\put(240,0){\color{teal}\circle*{3}}
\put(243,-3){\color{teal}$z$}
\put(0,240){\color{cyan}\circle*{3}}
\put(-6,243){\color{cyan}$y$}
\put(2,243){$y'$}
\put(0,40){\circle*{3}}
\put(1,44){$c'$}
\put(120,120){\color{blue}\circle*{3}}
\put(122,121){\color{blue}$x$}
\put(30,60){\circle*{2.5}}
\put(31,64){$b$}
\put(13.3,26.7){\circle*{2.5}}
\put(10.5,30){$o$}
\end{picture}
\smallskip

Observe that for the collinear triples $zcb'$ and $c'bx$, we have $\{a\}=\overline{\{z,c,b'\}}\cap\overline{\{c',b,x\}}$, $a'\in\Aline x{c'}\cap\Aline bc\cap\Aline {b'}x$ and $\Aline c{c'}\cap\Aline b{b'}=\{o\}\subseteq \Aline a{a'}$. Applying the inversivity-puls of $X$, we conclude that $\{\alpha\}=\Aline zb\cap\Aline cx\subseteq \Aline a{a'}$. 

Observe that for the collinear triples $zab'$ and $a'by$, we have $\{c\}=\overline{\{z,a,b'\}}\cap\overline{\{a',b,y\}}$, $c'\in\Aline z{a'}\cap\Aline ab\cap\Aline {b'}y$ and $\Aline a{a'}\cap\Aline b{b'}=\{o\}\subseteq \Aline c{c'}$. Applying the inversivity-puls of $X$, we conclude that $\{\gamma\}=\Aline zb\cap\Aline ay\subseteq \Aline c{c'}$. 

By the projectivity of $Y$, there exists a unique point $y'\in\Aline xz\cap\Aline c{a'}$. Consider the collinear triples $c'ax,a'cy'$ and observe that $\{b\}=\overline{\{c',a',x\}}\cap\overline{\{a',c,y'\}}$, $z\in \Aline {a'}{c'}\cap\Aline {a'}c\cap\Aline x{y'}$. Since $\Aline{a'}a\cap\Aline cx=\{\alpha\}\subseteq \Aline zb$, we can apply the inversivity-puls of $X$ and conclude that $\varnothing\ne \Aline c{c'}\cap\Aline a{y'}\subseteq \Aline zb$ and hence $\varnothing\ne \Aline zb\cap\Aline c{c'}\cap\Aline a{y'}=\{\gamma\}\cap\Aline a{y'}$ and hence $y'\in \Aline a\gamma\cap \Aline c{a'}=\{y\}$ and finally, $y=y'\in \Aline xz$, so the points $x,y,z$ are collinear.
\end{proof}

Theorem~\ref{t:duo-Desarg<=>}, \ref{t:tri-Desarg<=>}, \ref{t:proj-Moufang<=>} and Corollaries~\ref{c:Desarg-self-dual}, \ref{c:com-add-self-dual}, \ref{c:inv-add-self-dual} imply the following theorem of the self-duality of the (uno-, duo-, bi-, tri-) Desarguesianity.

\begin{corollary}\label{c:123-Desarg-self-dual} A projective plane is Desarguesian (resp. uno-Desarguesian, duo-Desar\-guesian, bi-Desarguesian, tri-Desarguesian) if and only if so is its dual projective plane.
\end{corollary}

Adding Theorems~\ref{t:duo-Desarg<=>} and \ref{t:tri-Desarg<=>} to the diagram from Section~\ref{s:diagram-proj}, we obtain the following diagram holding for every projective space.
$$
\xymatrix{
\mbox{associative-dot}\ar@{<=>}[r]\ar@{=>}[d]&\mbox{Desarguesian}\ar@{<=>}[r]\ar@{=>}[d]&\mathsf{(D)}\ar@{=>}[d]\\
\mbox{inversive-dot}\ar@{<=>}[r]\ar@{=>}[d]&\mbox{Moufang}\ar@{<=>}[r]\ar@{=>}[d]&\mathsf{(D_1)}\ar@{=>}[d]\\
\mbox{commutative-add}\ar@{<=>}[r]\ar@{=>}[d]&\mbox{Fano or Moufang}\ar@{<=>}[r]\ar@{=>}[d]&\mathsf{(D_2)}\ar@{=>}[d]\ar@{<=>}[r]&\mathsf{(D_{II})}\\
\mbox{invertible-add}\ar@{<=>}[r]&\mbox{$4$-Pappian}\ar@{<=>}[r]&\mathsf{(D_3)}
}
$$

This diagram motivates the following open problems.

\begin{problem} Is every tri-Desarguesian projective plane duo-Desarguesian?
\end{problem}

\begin{problem} Is every duo-Desarguesian projective plane uno-Desarguesian?
\end{problem}

\begin{problem} Is every tri-Desarguesian Playfair plane invertible-add?
\end{problem}

\begin{problem} Is every duo-Desarguesian bi-Desarguesian Playfair plane uno-Thalesian?
\end{problem}

\section{More incidences in the Pappus Axiom}

In this section we consider eight weaker modifications of the Pappus Axiom $\mathsf{(P)}$, defined as follows.

\begin{definition}\label{d:Pmn} For a liner $X$, consider the following weaker modifications of the Pappus Axiom
\begin{itemize}
\item[$\mathsf{(P)}$]  for every lines $L,L'\subseteq X$ and distinct points $a,b,c\in L\setminus L'$, $a',b',c'\in L'\setminus L$, the set $T\defeq (\Aline a{b'}\cap\Aline{a'}b)\cup(\Aline b{c'}\cap\Aline{b'}c)\cup (\Aline a{c'}\cap\Aline {a'}c)$ has rank $\|T\|\in\{0,2\}$:
\item[$\mathsf{(P^9_9)}$] for every lines $L,L'\subseteq X$ and distinct points $a,b,c\in L\setminus L'$, $a',b',c'\in L'\setminus L$, $x\in \Aline a{b'}\cap\Aline{a'}b$, $y\in \Aline b{c'}\cap\Aline{b'}c$, $z\in \Aline a{c'}\cap\Aline {a'}c$, the points $x,y,z$ are collinear;
\item[$\mathsf{(P_9^{10})}$] for every lines $L,L'\subseteq X$ and distinct points  $a,b,c\in L\setminus L'$, $a',b',c'\in L'\setminus L$, $x\in \Aline a{b'}\cap\Aline{a'}b$, $y\in \Aline b{c'}\cap\Aline{b'}c$, $z\in \Aline a{c'}\cap\Aline {a'}c$, $o\in L\cap L'\cap\Aline xy$, the points $x,y,z$ are collinear;
\item[$\mathsf{(P^9_{10})}$] for every lines $L,L'\subseteq X$ and distinct points $a,b,c\in L\setminus L'$, $a',b',c'\in L'\setminus L$, $x\in \Aline a{b'}\cap\Aline{a'}b$, $y\in \Aline b{c'}\cap\Aline{b'}c$, $z\in \Aline a{c'}\cap\Aline {a'}c\cap \Aline b{b'}$, the points $x,y,z$ are collinear;
\item[$\mathsf{(P^{10}_{10})}$] for every lines $L,L'\subseteq X$ and distinct points $a,b,c\in L\setminus L'$, $a',b',c'\in L'\setminus L$, $x\in \Aline a{b'}\cap\Aline{a'}b$, $y\in \Aline b{c'}\cap\Aline{b'}c$, $z\in \Aline a{c'}\cap\Aline {a'}c\cap \Aline b{b'}$,  $o\in L\cap L'\cap\Aline xz$, the points $x,y,z$ are collinear;
\item[$\mathsf{(P^{10}_{12})}$] for every lines $L,L'\subseteq X$ and distinct points $a,b,c\in L\setminus L'$, $a',b',c'\in L'\setminus L$, $x\in \Aline a{b'}\cap\Aline{a'}b$, $y\in \Aline b{c'}\cap\Aline{b'}c$, $z\in \Aline a{c'}\cap\Aline {a'}c$, $w\in \Aline a{a'}\cap\Aline b{b'}\cap\Aline c{c'}$, the points $x,y,z$ are collinear;
\item[$\mathsf{(P^{12}_{10})}$] for every lines $L,L'\subseteq X$ and distinct points $a,b,c\in L\setminus L'$, $a',b',c'\in L'\setminus L$, $x\in \Aline a{b'}\cap\Aline{a'}b$, $y\in \Aline b{c'}\cap\Aline{b'}c$, $z\in \Aline a{c'}\cap\Aline {a'}c$,  $u\in \Aline a{b'}\cap\Aline b{c'}$, $v\in \Aline {a'}b\cap \Aline {b'}c$, $o\in L\cap L'\cap \Aline uv$, the points $x,y,z$ are collinear;
\item[$\mathsf{(P^{13}_{13})}$] for every lines $L,L'\subseteq X$ and distinct points  $a,b,c\in L\setminus L'$, $a',b',c'\in L'\setminus L$, $x\in \Aline a{b'}\cap\Aline{a'}b$, $y\in \Aline b{c'}\cap\Aline{b'}c$, $z\in \Aline a{c'}\cap\Aline {a'}c$,   $u\in \Aline a{b'}\cap\Aline b{c'}$, $v\in \Aline {a'}b\cap \Aline {b'}c$, $o\in L\cap L'\cap \Aline uv$, $w\in \Aline a{a'}\cap\Aline b{b'}\cap\Aline c{c'}$, the points $x,y,z$ are collinear;
\item[$\mathsf{(\bar P_{13}^{13})}$] for every lines $L,L'\subseteq X$ and distinct points $a,b,c\in L\setminus L'$, $a',b',c'\in L'\setminus L$, $x\in \Aline a{b'}\cap\Aline{a'}b$, $y\in \Aline b{c'}\cap\Aline{b'}c$, $z\in \Aline a{c'}\cap\Aline {a'}c\cap\Aline b{b'}$,  $u\in \Aline a{b'}\cap\Aline b{c'}$, $v\in \Aline {a'}b\cap \Aline {b'}c$, $w\in \Aline a{a'}\cap\Aline b{b'}\cap\Aline c{c'}\cap\Aline uv$, $o\in L\cap L'\cap \Aline uv$, the points $x,y,z$ are collinear.
 \end{itemize}
For a pair of numbers $(n,m)\in\{(9,9),(9,10),(10,9),(10,12),(12,10),(13,13)\}$, a liner $X$ is defined to be $\mathsf{(P_n^m)}$ if it satisfies the axiom $\mathsf{(P_n^m)}$.
\end{definition}

\begin{remark} The notations of the axioms $\mathsf{(P_n^m)}$ in Definition~\ref{d:Pmn} have upper and lower indices corresponding to the number of points and lines involved in the definition of the axiom. The axioms $\mathsf{(P_{13}^{13})}$ and $\mathsf{(\bar P_{13}^{13})}$ both involve 13 points and 13 lines, but $\mathsf{(\bar P_{13}^{13})}$ has two incidences more comparing to the axiom $\mathsf{(P_{13}^{13})}$. In \cite{AlDhahirG1992}, the axioms $(\mathsf P^9_{10})$ and $(\mathsf P^{10}_9)$ are denoted by $(\mathsf P_1)$ and $(\mathsf P_2)$, and called the First and the Second Minor Pappus Axioms, respectively. Pickert \cite{Pickert} calls these axioms  the Axial and the Central Pappus Axioms, respectively. In \cite{AlDhahirG1992}, the axioms $(\mathsf P^{10}_{12}),(\mathsf P^{12}_{10})$ are denoted by $(\mathsf P_3),(\mathsf P_4)$, and called the Third and the Forth Minor Pappus Axioms, respectively.
\end{remark} 

\begin{proposition} Every subliner of a Pappian liner is $\mathsf{(P_9^9)}$. A projective liner is Pappian if and only if it is $(\mathsf P_9^9)$.
\end{proposition}

\begin{proof} Let $X$ be a subliner of a Pappian liner $Y$. To show that $X$ is  $\mathsf{(P_9^9)}$, fix any lines $L,L'\subseteq X$ and distinct points $a,b,c\in L\setminus L'$, $a',b',c'\in L'\setminus L$, $x\in \Aline ab\cap\Aline{a'}{b'}$, $y\in \Aline bc\cap\Aline{b'}{c'}$, $z\in \Aline ac\cap\Aline{a'}{c'}$. Since the liner $Y$ is Pappian, the set $(\Aline a{b'}\cap\Aline{a'}b)\cup(\Aline b{c'}\cap\Aline{b'}c)\cup (\Aline a{c'}\cap\Aline {a'}c)=\{x,y,z\}$ has rank $\|\{x,y,z\}\|\in\{0,2\}$, which implies that the points $x,y,z$ are collinear, witnessing that the liner $X$ is $\mathsf{(P_9^9)}$.

Now assume that $X$ is a projective liner. If $X$ is Pappian, then $X$ is $\mathsf{(P_9^9)}$ by the preceding paragraph. Now assume that $X$ is $\mathsf{(P_9^9)}$. To prove that $X$ is Pappian, fix any lines $L,L'\subseteq X$ and distinct points $a,b,c\in L\setminus L'$, $a',b',c'\in L'\setminus L$. We have to show that the set $T\defeq(\Aline a{b'}\cap\Aline{a'}b)\cup(\Aline b{c'}\cap\Aline{b'}c)\cup (\Aline a{c'}\cap\Aline {a'}c)$ has rank $\|T\|\in\{0,2\}$. This is clear if $T=\varnothing$. So, assume that $T\ne\varnothing$ and fix a point $x\in T$. We lose no generality assuming that $x\in \Aline a{b'}\cap\Aline {a'}b$. Then the lines $\Aline a{b'},\Aline {a'}b$ are coplanar and so are the lines $\Aline ab=L$ and $\Aline{a'}{b'}=L'$. Since the liner $X$ is projective, there exist unique points $y\in \Aline b{c'}\cap\Aline{b'}c$ and $z\in \Aline a{c'}\cap\Aline{a'}c$. Since $X$ is $\mathsf{(P_9^9)}$, the points $x,y,z$ are collinear and hence $\|T\|=\|\{x,y,z\}\|=2$, witnessing that the liner $X$ is Pappian. 
\end{proof}

\begin{proposition}\label{p:P109<=>com-puls} A projective liner is $(\mathsf P^{10}_9)$ if and only if it is commutative-puls.
\end{proposition}

\begin{proof} Assume that a projective liner $X$ is $\mathsf{(P^{10}_9)}$. To prove that $X$ is commutative-puls, take any distinct lines $L,L',\Lambda$ with $L\cap L'\cap\Lambda\ne\varnothing$ and distinct points $h\in L\cap L'$, $a,b,c\in L\setminus L'$, $a',b',c'\in L'\setminus L$ such that $(\Aline a{b'}\cap\Aline {a'}b)\cup(\Aline b{c'}\cap\Aline{b'}c)\subseteq\Lambda$. We have to prove that $\Aline a{c'}\cap\Aline{a'}c\subseteq\Lambda$. By the projectivity of the plane $\overline{L\cup L'}$, there exist unique points $x\in \Aline a{b'}\cap\Aline {a'}b$, $y\in\Aline b{c'}\cap\Aline {b'}c$ and $z\in \Aline a{c'}\cap\Aline{a'}c$. By our assumption, $\{x,y\}=(\Aline a{b'}\cap\Aline {a'}b)\cup(\Aline b{c'}\cap\Aline{b'}c)\subseteq\Lambda$. Since $x\in \Aline a{b'}\setminus\{b'\}$ and $y\in \Aline c{b'}\setminus\{b'\}$ and $a\ne c$, the points $x,y$ are distinct and hence $\Lambda=\Aline xy$. The axiom $\mathsf{(P^{10}_9)}$ implies that the points $x,y,z$ are collinear and hence $\Aline a{c'}\cap\Aline{a'}c=\{z\}\subseteq\Aline xy=\Lambda$, witnessing that the liner $X$ is commutative-puls.
\smallskip

Now assume that a projective liner $X$ is commutative-puls. To prove that it is $\mathsf{(P^{10}_9)}$, take any lines $L,L'\subseteq X$ and distinct points $a,b,c\in L\setminus L'$, $a',b',c'\in L'\setminus L$, $x\in\Aline a{b'}\cap\Aline {a'}b$, $y\in \Aline b{c'}\cap\Aline{b'}c$, $z\in \Aline a{c'}\cap\Aline{a'}c$, $o\in L\cap L'\cap\Aline xy$. We have to prove that the points $x,y,z$ are collinear. Consider the line $\Lambda\defeq\Aline xy$ and observe that  $L\cap L'\cap\Lambda\ne\varnothing$ and 
$(\Aline a{b'}\cap\Aline {a'}b)\cup(\Aline b{c'}\cap\Aline{b'}c)=\{x,y\}\subseteq\Aline xy=\Lambda$. Since the liner $X$ is commutative-puls, $\{z\}=\Aline a{c'}\cap\Aline {a'}c\subseteq \Lambda=\Aline xy$ and hence the points $x,y,z$ are collinear, witnessing that the liner $X$ is $\mathsf{(P^{10}_9)}$.
\end{proof}

\begin{proposition}\label{p:P910<=>com-plus} A projective liner is $(\mathsf P^9_{10})$ if and only if it is commutative-plus.
\end{proposition}

\begin{proof} Assume that a projective liner $X$ is $\mathsf{(P^9_{10})}$. 
To prove that the liner $X$ is commutative-plus, take any distinct points $h,v\in X$, $d\in \Aline hv$, $a,b,c,a',b',c'\in X\setminus\Aline hv$ such that $h\in \Aline ab\cap\Aline c{c'}\cap\Aline {a'}{b'}$, $v\in \Aline ac\cap\Aline b{b'}\cap\Aline {a'}{c'}$ and $d\in \Aline a{a'}\cap\Aline b{c'}$. We have to prove that $d\in \Aline {b'}c$. Consider the lines $L\defeq \overline{\{a,h,b\}}$ and $L'\defeq\overline{\{c',v,a'\}}$ and observe that $c\in \Aline av\cap\Aline h{c'}$, $b'\in \Aline h{a'}\cap\Aline vb$, and $d\in \Aline a{a'}\cap\Aline b{c'}\cap\Aline hv$. Since $X$ is  $(\mathsf P^9_{10})$ those conditions imply that the points $c,d,b'$ are collinear, and hence $d\in \Aline{b'}c$, witnessing that the liner $X$ is commutative-plus.
\smallskip

Now assume that a projective liner $X$ is commutative-plus. To prove that $X$ is $\mathsf{(P^9_{10})}$, take any lines $L,L'\subseteq X$ and distinct points $a,b,c\in L\setminus L'$, $a',b',c'\in L'\setminus L$, $x\in \Aline a{b'}\cap\Aline{a'}b$, $y\in \Aline b{c'}\cap\Aline{b'}c$, $z\in \Aline a{c'}\cap\Aline{a'}c\cap\Aline b{b'}$. We have to prove that the points $x,y,z$ are collinear. Consider the points $h\defeq b$, $v\defeq b'$, $d\defeq z$, $\alpha\defeq a$, $\beta\defeq c$, $\gamma\defeq x$, $\alpha'\defeq c'$, $\beta'\defeq y$, $\gamma'\defeq a'$, and observe that $\{\alpha,\beta,\gamma,\alpha',\beta',\gamma'\}=\{a,c,x,c',y,a'\}\subseteq X\setminus\Aline b{b'}=X\setminus \Aline hv$, $h=b\in\Aline ac\cap \Aline x{a'}\cap\Aline {c'}y= \Aline \alpha\beta\cap \Aline\gamma{\gamma'}\cap \Aline{\alpha'}{\beta'}$, 
$v=b'\in \Aline ax\cap\Aline cy\cap \Aline{c'}{a'}=\Aline\alpha\gamma\cap\Aline\beta{\beta'}\cap\Aline{\alpha'}{\gamma'}$ and $d=z\in \Aline a{c'}\cap \Aline c{a'}=\Aline{\alpha}{\alpha'}\cap\Aline \beta{\gamma'}$. Since $X$ is commutative-plus, these conditions imply $z=d\in \Aline {\beta'}\gamma=\Aline yx$. So, the points $x,y,z$ are collinear and hence the projective liner $X$ is $\mathsf{(P^9_{10})}$.
\end{proof}

\begin{theorem} For a projective plane $X$, the following conditions are equivalent:
\begin{enumerate}
\item $X$ is $\mathsf{(P^9_{10})}$ and $\mathsf{(P^{10}_9)}$;
\item $X$ is commutative-add;
\item $X$ is duo-Desarguesian;
\item $X$ is bi-Desarguesian;
\item $X$ is everywhere uno-Thalesian;
\item $X$ is Fano or Moufang.
\end{enumerate}
\end{theorem}

\begin{proof} The equivalence $(1)\Leftrightarrow(2)$ follows from Propositions~\ref{p:P109<=>com-puls} and \ref{p:P910<=>com-plus}, and the equivalence of the conditions $(2)$--$(6)$ was proved in Theorem~\ref{t:duo-Desarg<=>}. 
\end{proof}

\begin{theorem}\label{t:P1010<=>} For a projective plane $X$, the following conditions are equivalent:
\begin{enumerate}
\item $X$ is $\mathsf{(P^{10}_{10})}$;
\item $X$ is invertible-add;
\item $X$ is invertible-plus;
\item $X$ is invertible-puls;
\item $X$ is $4$-Pappian;
\item $X$ is tri-Desarguesian.
\end{enumerate}
\end{theorem}

\begin{proof} The equivalence of the conditions (2)--(6) was proved in Theorem~\ref{t:tri-Desarg<=>}.
\smallskip

$(1)\Ra(4)$ Assume that a projective plane $X$ is $\mathsf{(P_{10}^{10})}$. To prove that $X$ is invertible-puls, take any lines $L,L'\subseteq X$ and distinct points $a,b,c\in L\setminus L'$, $a',b',c'\in L'\setminus L$, $h\in L\cap L'$, $v\in \Aline a{a'}\cap\Aline b{b'}\cap\Aline c{c'}$ such that $\Aline a{b'}\cap\Aline b{c'}\subseteq \Aline vh$. We have to prove that $\Aline {a'}b\cap\Aline {b'}c\subseteq \Aline vh$.

Consider the points $a''\defeq c'$, $b''\defeq b'$, $c''\defeq a'$ and unique points $x\in \Aline a{b''}\cap\Aline {a''}{b}=\Aline a{b'}\cap\Aline {c'}b\subseteq \Aline vh$, $y\in \Aline b{c''}\cap\Aline {b''}c=\Aline b{a'}\cap\Aline{b'}c$ and $z\defeq\Aline a{c''}\cap\Aline{a''}c=\Aline a{a'}\cap\Aline {c'}c=\{v\}\subseteq \Aline b{b'}$. Observe that $L\cap L'\cap \Aline xz=\{h\}\cap\Aline xv=\{h\}\ne\varnothing$. Applying the axiom $\mathsf{(P_{10}^{10})}$, we conclude that the points $x,y,z$ are collinear and hence $\Aline {a'}b\cap\Aline {b'}c=\{y\}\subseteq\Aline xz=\Aline xv=\Aline hv$, witnessing that the projective plane $X$ is invertible-puls.
\smallskip

$(4)\Ra(1)$ Assume that a projective plane $X$ is invertible-puls. To prove that $X$ is $\mathsf{(P_{10}^{10})}$, take any lines $L,L'\subseteq X$ and distinct points $a,b,c\in L\setminus L'$, $a',b',c'\in L'\setminus L$, $x\in \Aline a{b'}\cap\Aline{a'}b$, $y\in \Aline b{c'}\cap\Aline{b'}c$, $z\in \Aline a{c'}\cap\Aline {a'}c\cap \Aline b{b'}$ and $h\in L\cap L'\cap\Aline xz$. We should prove that the points $x,y,z$ are collinear. Consider the points $a''\defeq c'$, $b''\defeq b'$, $c''\defeq a'$, and $v\defeq z$.  Observe that $v=z\in \Aline a{c'}\cap\Aline {a'}c\cap\Aline b{b'}=\Aline a{a''}\cap\Aline {c''}c\cap\Aline b{b''}$, $\Aline a{b''}\cap\Aline b{c''}=\Aline a{b'}\cap\Aline b{a'}=\{x\}\subseteq \Aline xz=\Aline hv$. Since $X$ is invertible-puls, the latter conditions imply $\Aline {a''}b\cap\Aline {b''}c\subseteq\Aline vh$ and hence $y\in \Aline {c'}b\cap\Aline {b'}c=\Aline {a''}{b}\cap\Aline {b''}c\subseteq \Aline vh=\Aline zh=\Aline zx$, witnessing that the points $x,y,z$ are collinear and the liner $X$ is  $\mathsf{(P_{10}^{10})}$.
\end{proof}

Propositions~\ref{p:P109<=>com-puls}, \ref{p:P910<=>com-plus} and Corollary~\ref{c:com-plus-puls-dual} imply the following corollary.

\begin{corollary}\label{c:P910dualP109} A projective plane is $(\mathsf P^9_{10})$ if and only if the dual projective plane is $(\mathsf P^{10}_9)$.
\end{corollary}

\begin{exercise} Find a direct proof of Corollary~\ref{c:P910dualP109}.
\end{exercise}

\begin{proposition}\label{p:P910dualP109} A projective plane is $\mathsf{(P^{10}_{12})}$ if and only if the dual projective plane is $(\mathsf P^{12}_{10})$.
\end{proposition}

\begin{proof} Assume that a projective plane $\Pi$ is $\mathsf{(P^{10}_{12})}$. To prove that the dual projective plane $\Pi^*$ is $\mathsf{(P^{12}_{10})}$, fix any distinct points $\ell,\ell'\in \Pi$ and distinct lines $A,B,C,A',B',C'$ in $\Pi$ such that $\ell\in (A\cap B\cap C)\setminus(A'\cup B'\cup C')$ and $\ell'\in (A'\cap B'\cap C')\setminus(A\cup B\cup B)$, and consider the lines $$
\begin{gathered}
X\defeq\overline{(A\cap B')\cup(A'\cap B)},\quad Y\in  \overline{(B\cap C')\cup(B'\cap C)},\quad Z\defeq \overline{(A\cap C')\cup(A'\cap C)},\\
U\defeq\overline{(A\cap B')\cup(B\cap C')},\quad V\defeq\overline{(A'\cap B)\cup(B'\cap C)}.
\end{gathered}
$$ Assuming that $U\cap V\subseteq \overline{\ell{\ell'}}$, we have to prove that $X\cap Y\cap Z\ne\varnothing$. Consider the unique points $a\defeq \ell$, $a'\defeq\ell'$, $b\in B\cap C'$, $b'\in A\cap B'$, $c\in B\cap A'$, $c'\in C\cap B'$, $x\in A\cap C'=\Aline a{b'}\cap\Aline {a'}b$, $y\in \Aline b{c'}\cap\Aline {b'}c$, $z\in C\cap A'=\Aline a{c'}\cap\Aline {a'}c$, and $w\in U\cap V=\Aline {b'}{b}\cap \Aline c{c'}$. Observe that $\{a,b,c\}\subseteq B\setminus B'$ and $\{a',b',c'\}\subseteq B'\setminus B$. Taking into account that $\Aline {b'}{b}\cap\Aline c{c'}=U\cap V=\{w\}\subseteq\Aline \ell{\ell'}=\Aline a{a'}$, we conclude that $w\in \Aline a{a'}\cap\Aline b{b'}\cap\Aline c{c'}$ and hence the points $x,y,z$ are collinear by the axiom $\mathsf{(P^{10}_{12})}$ holding for the  liner $\Pi$. Then $X\cap Y\cap Z=\Aline{b'}{c}\cap \Aline {b}{c'}\cap\Aline xz=\{y\}\cap\Aline zx=\{y\}\ne\varnothing$, witnessing that the dual projective plane $\Pi^*$ is $\mathsf{(P^{12}_{10})}$.  
\smallskip

Assume that a projective plane $\Pi$ is $\mathsf{(P^{12}_{10})}$. To prove that the dual projective plane $\Pi^*$ is $\mathsf{(P^{10}_{12})}$, fix any distinct points $\ell,\ell'\in X$ and distinct lines $A,B,C,A',B',C'$ in $\Pi$ such that $\ell\in (A\cap B\cap C)\setminus(A'\cup B'\cup C')$ and $\ell'\in (A'\cap B'\cap C')\setminus(A\cup B\cup B)$, and consider the lines $X\defeq\overline{(A\cap B')\cup(A'\cap B)}$, $Y\in  \overline{(B\cap C')\cup(B'\cap C)}$, $Z\defeq \overline{(A\cap C')\cup(A'\cap C)}$.
Assuming that the singletons  $A\cap A'$, $B\cap B'$ and $C\cap C'$ are contained in some line, we have to prove that $X\cap Y\cap Z\ne\varnothing$. Consider the unique points $a\in B\cap A'$, $a'\in C\cap B'$, $b\defeq \ell$, $b'\defeq\ell'$, $c\in B\cap C'$, $c'\in A\cap B'$, $\alpha\in A\cap A'$, $o\in B\cap B'$, $\gamma\in C\cap C'$, $x\in C\cap A'$, $y\in A\cap C'$, $z\in \Aline a{c'}\cap\Aline{a'}c$. Observe that $\{a,b,c\}\subseteq B\setminus B'$ and $\{a',b',c'\}\subseteq B'\setminus B$. Taking into account that $\alpha\in A\cap A'=\Aline b{c'}\cap\Aline a{b'}$, $\gamma\in C\cap C'=\Aline{a'}b\cap\Aline {b'}c$ and $o\in \Aline \alpha\gamma\cap B\cap B'$, we conclude that the points $x,y,z$ are collinear by the axiom $\mathsf{(P^{12}_{10})}$ holding for the projective plane $\Pi$. Then $X\cap Y\cap Z=\Aline{c'}{a}\cap \Aline {c}{a'}\cap\Aline yx=\{z\}\cap\Aline yx=\{z\}\ne\varnothing$, witnessing that the dual projective plane $\Pi^*$ is $\mathsf{(P^{10}_{12})}$.  
\end{proof}

\begin{proposition}\label{p:2D=>P1012} Every bi-Desarguesian projective liners is  $\mathsf{(P^{10}_{12})}$.
\end{proposition}

\begin{proof} Assume that $X$ is a bi-Desarguesian projective liner. To prove that $X$ is $\mathsf{(P^{10}_{12})}$, take any lines $L,L'\subset X$ and distinct points $a,b,c\in L\setminus L'$, $a',b',c'\in L'\setminus L$, $x\in \Aline a{b'}\cap\Aline {a'}b$, $y\in \Aline b{c'}\cap\Aline{b'}c$, $z\in \Aline a{c'}\cap\Aline{a'}c$, $w\in\Aline a{a'}\cap\Aline b{b'}\cap \Aline c{c'}$. We have to prove that the points $x,y,z$ are collinear. It follows from $x\in \Aline a{b'}\cap\Aline {a'}b$ that the lines $\Aline a{b'}\cap\Aline {a'}b$ are coplanar and so are the lines $\Aline ab=L$ and $\Aline{a'}{b'}=L'$. Since the liner $X$ is projective, there exists a unique point $o\in L\cap L'=\Aline ac\cap\Aline{a'}{c'}$. By Theorem~\ref{t:duo-Desarg<=>}, the bi-Desarguesian projective liner is Fano or Moufang. 
\smallskip

If $X$ is Fano, then the quadrangles $abb'a'$, $bcc'b'$, $acc'a'$ are Fano and hence
$\{x,y,z\}\subseteq \Aline ow$, witnessing that the points $x,y,z$ are collinear and the liner $X$ is $\mathsf{(P^{10}_{12})}$.
\smallskip

Next, assume that the liner $X$ is Moufang.  Observe that the triangles $ab'c$ and $a'bc'$ are perspective from the point $w$, and also $b\in\Aline ac$ and $b'\in \Aline{a'}{c'}$. Since $X$ is bi-Desarguesian, the points $x\in \Aline a{b'}\cap\Aline{a'}b$, $y\in \Aline {b'}c\cap\Aline b{c'}$ and $o\in L\cap L'=\Aline ac\cap\Aline{a'}{c'}$ are collinear.
The triangles $ac'b$ and $a'cb'$ also are perspective from the point $w$. Moreover $c\in \Aline ab$ and $c'\in \Aline {a'}{b'}$. Since $X$ is bi-Desarguesian, the points $o\in \Aline a{c'}\cap\Aline {a'}c$, $y\in \Aline {c'}b\cap\Aline c{b'}$ and $z\in \Aline ab\cap\Aline {a'}{b'}$ are collinear. Then $\{x,y,z\}\subseteq\Aline oy$ and hence the points $x,y,z$ are collinear, witnessing that the liner $X$ is  $\mathsf{(P^{10}_{12})}$.
\end{proof}

\begin{corollary} Every bi-Desarguesian projective liner $X$ is $\mathsf{(P^{12}_{10})}$.
\end{corollary}

\begin{proof} To prove that $X$ is $\mathsf{(P^{12}_{10})}$, take any lines $L,L'\subset X$ and distinct points $a,b,c\in L\setminus L'$, $a',b',c'\in L'\setminus L$, $x\in \Aline a{b'}\cap\Aline {a'}b$, $y\in \Aline b{c'}\cap\Aline{b'}c$, $z\in \Aline a{c'}\cap\Aline{a'}c$, $u\in \Aline a{b'}\cap\Aline {a'}b$, $v\in \Aline {a'}b\cap\Aline {b'}c$, and $o\in L\cap L'\cap\Aline uv$. We should prove that the points $x,y,z$ are collinear. Consider the projective plane $\Pi\defeq\overline{L\cup L'}$. Since $X$ is bi-Desarguesian, so is the plane $\Pi$ in $X$. By Corollary~\ref{c:123-Desarg-self-dual}, the dual projective plane $\Pi^*$ is bi-Desarguesian. By Proposition~\ref{p:2D=>P1012}, the bi-Desarguesian projective plane $\Pi^*$ is $\mathsf{(P^{10}_{12})}$. By Proposition~\ref{p:P910dualP109}, the projective plane $\Pi$ is $\mathsf{(P^{12}_{10})}$. Then the points $x,y,z$ are collinear in the plane $\Pi\subseteq X$ and also in the liner $X$.
\end{proof}

\begin{proposition}\label{p:P1313-selfdual} A projective plane is $\mathsf{(P^{13}_{13})}$ if and only if so is its dual projective plane.
\end{proposition}

\begin{proof} Assume that a projective plane $\Pi$ is $\mathsf{(P^{13}_{13})}$. To prove that the dual projective plane $\Pi^*$ is  $\mathsf{(P^{13}_{13})}$, take 
any distinct points $\ell,\ell'\in \Pi$, distinct lines $A,B,C,A',B',C'$ in $\Pi$ such that $\ell\in (A\cap B\cap C)\setminus(A'\cup B'\cup C')$ and $\ell'\in (A'\cap B'\cap C')\setminus(A\cup B\cup B)$, and consider the lines $$
\begin{gathered}
X\defeq\overline{(A\cap B')\cup(A'\cap B)},\quad Y\in  \overline{(B\cap C')\cup(B'\cap C)},\quad Z\defeq \overline{(A\cap C')\cup(A'\cap C)},\\
U\defeq\overline{(A\cap B')\cup(B\cap C')},\quad V\defeq\overline{(A'\cap B)\cup(B'\cap C)}.
\end{gathered}
$$ Assuming that $U\cap V\subseteq \Aline{\ell}{\ell'}$ and the set $(A\cap A')\cup (B\cap B')\cup (C\cap C')$ is contained in some line $W$, we have to prove that $X\cap Y\cap Z\ne\varnothing$. 

Consider the unique points $a\in B\cap A'$, $a'\in C\cap B'$, $b\defeq \ell$, $b'\defeq\ell'$, $c\in B\cap C'$, $c'\in A\cap B'$, $x\in \Aline a{b'}\cap\Aline {a'}b=A'\cap C$, $y\in \Aline b{c'}\cap\Aline{b'}c=A\cap C'$, $z\in \Aline a{c'}\cap\Aline{a'}c$, $u\in\Aline a{b'}\cap\Aline b{c'}=A'\cap A$, $v\in\Aline {a'}b\cap\Aline {b'}c=C\cap C'$, and $o\in B\cap B'$. Observe that $\{a,b,c\}\subseteq B$ and  $\{a',b',c'\}\subseteq B'$. Also $\{u,o,v\}=(A\cap A')\cup(B\cap B')\cup(C\cap C')\subseteq W$ and 
$$\Aline a{a'}\cap \Aline b{b'}\cap\Aline c{c'}=\overline{(B\cap A')\cup(C\cap B')}\cap\Aline \ell{\ell'}\cap\overline{(B\cap C')\cup(A\cap B')}=V\cap\Aline \ell{\ell'}\cap U\ne\varnothing.$$ Applying the axiom $\mathsf{(P^{13}_{13})}$, we conclude that the points $x,y,z$ are collinear. Then $$X\cap Y\cap Z=
\Aline {c'}a\cap \Aline c{a'}\cap \Aline yx=\{z\}\cap\Aline yx=\{z\}\ne\varnothing,$$ witnessing that the dual projective plane $\Pi^*$ is $\mathsf{(P^{13}_{13})}$.
\end{proof}

\begin{proposition}\label{p:P109=>P1313} If a projective liner $X$ is $\mathsf{(P^{10}_9)}$,  then it is $\mathsf{(P^{13}_{13})}$.
\end{proposition}

\begin{proof} To prove that $X$ is $\mathsf{(P^{13}_{13})}$, take any lines $L,L'\subseteq X$ and distinct points  $a,b,c\in L\setminus L'$, $a',b',c'\in L'\setminus L$, $x\in \Aline a{b'}\cap\Aline{a'}b$, $y\in \Aline b{c'}\cap\Aline{b'}c$, $z\in \Aline a{c'}\cap\Aline {a'}c$,   $u\in \Aline a{b'}\cap\Aline b{c'}$, $v\in \Aline {a'}b\cap \Aline {b'}c$, $o\in L\cap L'\cap \Aline uv$, and $w\in\Aline a{a'}\cap\Aline b{b'}\cap\Aline c{c'}$. We have to prove that the points $x,y,z$ are collinear. Consider the points $a''\defeq c'$, $b''\defeq b'$ and $c''\defeq a'$. Observe that $u\in \Aline a{b'}\cap\Aline b{c'}=\Aline a{b'}\cap \Aline b{a''}$, $v\in \Aline {a'}b\cap\Aline{b'}c=\Aline{c''}b\cap\Aline {b''}{c}$ and  $\Aline{a}{c''}\cap\Aline {a''}c\cap\Aline b{b''}=\Aline a{a'}\cap\Aline {c'}c\cap\Aline b{b'}=\{w\}$. Since $X$ is $\mathsf{(P^{10}_9)}$, the points $u,v,w$ are collinear. Consider the closure $\Pi\defeq\langle\{a,b,a',b'\}\rangle$ of the set $\{a,b,a',b'\}$ in the projective liner $X$. By Proposition~\ref{p:P109<=>com-puls}, the $\mathsf{(P^{10}_9)}$ projective liner $X$ is commutative-puls and hence invertible-puls. By Theorem~\ref{t:tri-Desarg<=>}, the invertible-puls liner $X$ is $4$-Pappian and hence the subliner $\Pi$ of $X$ is Pappian. Observe that $\{o,w,x\}=(\Aline ab\cap\Aline {a'}{b'})\cup(\Aline a{a'}\cap\Aline b{b'})\cup(\Aline a{b'}\cap\Aline{a'}b)\subseteq \Pi$, $\{u,v\}=(\Aline a{b'}\cap\Aline ow)\cup(\Aline{a'}b\cap\Aline ow)\subseteq \Pi$, $\{c,c'\}=(\Aline ab\cap\Aline v{b'})\cup(\Aline{a'}{b'}\cap \Aline ub)\subseteq \Pi$. Therefore, $\{a,b,c,a',b',c',x,y,z\}\subseteq \Pi$. Since the projective plane $\Pi$ is Pappian, the points $x,y,z$ are collinear.
\end{proof}

\begin{proposition}\label{p:P910=>P1313} If a projective liner $X$ is $\mathsf{(P^9_{10})}$,  then it is $\mathsf{(P^{13}_{13})}$.
\end{proposition}

\begin{proof}  To prove that $X$ is $\mathsf{(P^{13}_{13})}$, take any lines $L,L'\subseteq X$ and distinct points  $a,b,c\in L\setminus L'$, $a',b',c'\in L'\setminus L$, $x\in \Aline a{b'}\cap\Aline{a'}b$, $y\in \Aline b{c'}\cap\Aline{b'}c$, $z\in \Aline a{c'}\cap\Aline {a'}c$,   $u\in \Aline a{b'}\cap\Aline b{c'}$, $v\in \Aline {a'}b\cap \Aline {b'}c$, $o\in L\cap L'\cap \Aline uv$, and $w\in\Aline a{a'}\cap\Aline b{b'}\cap\Aline c{c'}$. We have to prove that the points $x,y,z$ are collinear.

Since $X$ is  $\mathsf{(P^9_{10})}$, so is the plane $\Pi\defeq\overline{L\cup L'}$ in $X$. By Corollary~\ref{c:P910dualP109}, the dual projective plane $\Pi^*$ is   $\mathsf{(P^{10}_{9})}$. By Proposition~\ref{p:P109=>P1313}, the  $\mathsf{(P^{10}_{9})}$ projective plane $\Pi^*$ is $\mathsf{(P^{13}_{13})}$. By Proposition~\ref{p:P1313-selfdual}, the projective plane $\Pi$ is $\mathsf{(P^{13}_{13})}$ and hence the points $x,y,z$ are collinear.
\end{proof}

\begin{proposition} A projective plane is $\mathsf{(\bar P^{13}_{13})}$ if and only if so is its dual projective plane. 
\end{proposition}

\begin{proof} Assume that a projective plane $\Pi$ is $\mathsf{(\bar P^{13}_{13})}$. To prove that the dual projective plane $\Pi^*$ is  $\mathsf{(\bar P^{13}_{13})}$, take 
any distinct points $\ell,\ell'\in \Pi$, distinct lines $A,B,C,A',B',C'$ in $\Pi$ such that $\ell\in (A\cap B\cap C)\setminus(A'\cup B'\cup C')$ and $\ell'\in (A'\cap B'\cap C')\setminus(A\cup B\cup B)$. Consider the lines $$
\begin{gathered}
X\defeq\overline{(A\cap B')\cup(A'\cap B)},\quad Y\in  \overline{(B\cap C')\cup(B'\cap C)},\quad Z\defeq \overline{(A\cap C')\cup(A'\cap C)},\\
U\defeq\overline{(A\cap B')\cup(B\cap C')},\quad V\defeq\overline{(A'\cap B)\cup(B'\cap C)}.
\end{gathered}
$$ Assuming that the set $(A\cap A')\cup (B\cap B')\cup (C\cap C')$ is contained in some line $W$ and $U\cap V\cap W\cap\overline{\ell}{\ell'}$ contains some point $w$, we have to prove that $X\cap Y\cap Z\ne\varnothing$. 

Consider the unique points $a\in B\cap A'$, $a'\in C\cap B'$, $b\defeq \ell$, $b'\defeq\ell'$, $c\in B\cap C'$, $c'\in A\cap B'$, $x\in \Aline a{b'}\cap\Aline {a'}b=A'\cap C$, $y\in \Aline b{c'}\cap\Aline{b'}c=A\cap C'$, $z\in \Aline a{c'}\cap\Aline{a'}c$, $u\in\Aline a{b'}\cap\Aline b{c'}=A'\cap A$, $v\in\Aline {a'}b\cap\Aline {b'}c=C\cap C'$, and $o\in B\cap B'$. Observe that $\{a,b,c\}\subseteq B$ and  $\{a',b',c'\}\subseteq B'$. Also $\{u,o,v\}=(A\cap A')\cup(B\cap B')\cup(C\cap C')\subseteq W$ and 
$$\Aline a{a'}\cap \Aline b{b'}\cap\Aline c{c'}=\overline{(B\cap A')\cup(C\cap B')}\cap\Aline \ell{\ell'}\cap\overline{(B\cap C')\cup(A\cap B')}=V\cap\Aline \ell{\ell'}\cap U=\{w\}\subseteq W=\Aline uv.$$ Applying the axiom $\mathsf{(\bar P^{13}_{13})}$, we conclude that the points $x,y,z$ are collinear. Then $$X\cap Y\cap Z=
\Aline {c'}a\cap \Aline c{a'}\cap \Aline yx=\{z\}\cap\Aline yx=\{z\}\ne\varnothing,$$ witnessing that the dual projective plane $\Pi^*$ is $\mathsf{(\bar P^{13}_{13})}$.
\end{proof}

\begin{proposition} Every $4$-Pappian projective liner $X$ is $\mathsf{(\bar P_{13}^{13})}$.
\end{proposition}

\begin{proof} To prove that $X$ is $\mathsf{(\bar P^{13}_{13})}$, take any lines $L,L'\subseteq X$ and distinct points  $a,b,c\in L\setminus L'$, $a',b',c'\in L'\setminus L$, $x\in \Aline a{b'}\cap\Aline{a'}b$, $y\in \Aline b{c'}\cap\Aline{b'}c$, $z\in \Aline a{c'}\cap\Aline {a'}c$,   $u\in \Aline a{b'}\cap\Aline b{c'}$, $v\in \Aline {a'}b\cap \Aline {b'}c$, $o\in L\cap L'\cap \Aline uv$, and $w\in \Aline a{a'}\cap\Aline b{b'}\cap\Aline c{c'}\cap\Aline uv$. We have to prove that the points $x,y,z$ are collinear. Consider the closure $\Pi$ of the set $\{a,b,a',b'\}$ in the projective liner $X$. Since $X$ is $4$-Pappian, the subliner $\Pi$ of $X$ is Pappian. Observe that $\{o,w,x\}=(\Aline ab\cap\Aline {a'}{b'})\cup(\Aline a{a'}\cap\Aline b{b'})\cup(\Aline a{b'}\cap\Aline{a'}b)\subseteq \Pi$, $\{u,v\}=(\Aline a{b'}\cap\Aline ow)\cup(\Aline{a'}b\cap\Aline ow)\subseteq \Pi$, $\{c,c'\}=(\Aline ab\cap\Aline v{b'})\cup(\Aline{a'}{b'}\cap \Aline ub)\subseteq \Pi$. Therefore, $\{a,b,c,a',b',c',x,y,z\}\subseteq \Pi$. Since the plane $\Pi$ is Pappian, the points $x,y,z$ are collinear.
\end{proof}

Now we characterize liners whose free projectivizations are $\mathsf{(P_{13}^{13})}$ or $\mathsf{(\bar P_{13}^{13})}$. This will help us to construct examples of projective planes distinguishing the properties  $\mathsf{(P_{10}^{10})}$, $\mathsf{(P_{13}^{13})}$ and $\mathsf{(\bar P_{13}^{13})}$.

\begin{proposition}\label{p:projectivization<=>P1313} The  free projectivization $\widehat X$ of a liner $X$ is $\mathsf{(P_{13}^{13})}$ if and only if for any lines $L,L'\subseteq X$ and distinct points $a,b,c\in L\setminus L'$, $a',b',c'\in L'\setminus L$, $x\in X\cap\Aline a{b'}\cap\Aline{a'}b$, $y\in X\cap\Aline b{c'}\cap\Aline{b'}c$, $u\in X\cap\Aline a{b'}\cap\Aline b{c'}$, $v\in X\cap\Aline{a'}b\cap\Aline {b'}c$, $o\in L\cap L'\cap\Aline uv$ and $w\in X\cap\Aline a{a'}\cap\Aline b{b'}\cap\Aline c{c'}$, there exists a point $z\in X\cap \Aline a{c'}\cap\Aline {a'}c\cap\Aline xy$.
\end{proposition}

\begin{proof} 
Assume that the free projectivization $\widehat X$ of the liner $X$ is $\mathsf{(P_{13}^{13})}$. Take any lines $L,L'\subseteq X$ and distinct points $a,b,c\in L\setminus L'$, $a',b',c'\in L'\setminus L$, $x\in X\cap \Aline a{b'}\cap\Aline{a'}b$, $y\in X\cap \Aline b{c'}\cap\Aline{b'}c$, $u\in X\cap \Aline a{b'}\cap\Aline b{c'}$, $v\in X\cap \Aline{a'}b\cap\Aline {b'}c$, $o\in L\cap L'\cap\Aline uv$ and $w\in X\cap \Aline a{a'}\cap\Aline b{b'}\cap\Aline c{c'}$. By the projectivity of the liner $\widehat X$, there exists a unique point $z\in \Aline a{c'}\cap\Aline {a'}c\subseteq\widehat X$. Since $\widehat X$ is $\mathsf{(P_{13}^{13})}$, the points $x,y,z$ are collinear. Observe that the subliner $\{a,b,c,a',b',c',x,y,z,o,u,v,w\}$ of $\widehat X$ is $3$-wide. By Theorem~\ref{t:free-projectivization}(5), this subliner is contained in $X$ and hence $z\in X\cap\Aline ac\cap \Aline {a'}{c'}\cap \Aline xy$.
\smallskip

Now assume that the liner $X$ satisfies the condition, formulated in the theorem. To prove that the free projectivization $\widehat X$ is $\mathsf{(P_{13}^{13})}$, take any lines $L,L'$ in $\widehat X$ and distinct points $a,b,c\in L\setminus L'$, $a',b',c'\in L'\setminus L$, $x\in \Aline a{b'}\cap\Aline{a'}b$, $y\in \Aline b{c'}\cap\Aline{b'}c$, $z\in \Aline a{c'}\cap\Aline {a'}c$, $u\in \Aline a{b'}\cap\Aline b{c'}$, $v\in \Aline{a'}b\cap\Aline {b'}c$, $o\in L\cap L'\cap\Aline uv$ and $w\in \Aline a{a'}\cap\Aline b{b'}\cap\Aline c{c'}$. Observe that the subliner $\{a,b,c,a',b',c',o,x,y,u,v,w\}$ of $\widehat X$ is $3$-wide. By Theorem~\ref{t:free-projectivization}(5), this subliner is contained in $X$. By the condition, there exists a point $z'\in X\cap\Aline {a'}c\cap\Aline a{c'}\cap\Aline xy$. Then $z'\in \Aline {a'}{c}\cap\Aline a{c'}=\{z\}$ and hence the points $x,y,z$ are collnear, witnessing that the projective plane $\widehat X$ is $\mathsf{(P^{13}_{13})}$.
\end{proof}

By analogy we can prove a characterization of liners whose free projectivization is $\mathsf{(\bar P^{13}_{13})}$.

\begin{proposition}\label{p:projectivization<=>barP1313} The  free projectivization $\widehat X$ of a liner $X$ is $\mathsf{(\bar P_{13}^{13})}$ if and only if for any lines $L,L'\subseteq X$ and distinct points $a,b,c\in L\setminus L'$, $a',b',c'\in L'\setminus L$, $x\in X\cap\Aline a{b'}\cap\Aline{a'}b$, $y\in X\cap\Aline b{c'}\cap\Aline{b'}c$, $u\in X\cap\Aline a{b'}\cap\Aline b{c'}$, $v\in X\cap\Aline{a'}b\cap\Aline {b'}c$, $o\in L\cap L'\cap\Aline uv$ and $w\in X\cap\Aline a{a'}\cap\Aline b{b'}\cap\Aline c{c'}\cap\Aline uv$, there exists a point $z\in X\cap \Aline a{c'}\cap\Aline {a'}c\cap\Aline xy$.
\end{proposition}

\begin{proposition}\label{p:projectivizationP1010not} If the  free projectivization $\widehat X$ of a liner $X$ is $\mathsf{(P_{10}^{10})}$, then $X=\widehat X$.
\end{proposition}

\begin{proof} To derive a contradiction, assume that $X\ne \widehat X$. By Theorem~\ref{t:free-projectivization}, $\widehat X$ is an $\w$-long projective plane. Then $\widehat X$ contains a triangle $oaa'$ such that $\{o,a,a'\}\not\subseteq X$. Since $\widehat X$ is $\w$-long, there exist points $b\in\Aline oa\setminus\{o,a\}$ and $b'\in\Aline o{a'}\setminus\{o,a'\}$. By the projectivity of the plane $\widehat X$, there exist unique points $x\in \Aline a{b'}\cap \Aline {a'}b$, $z\in \Aline b{b'}\cap\Aline ox$, $c\in \Aline z{a'}\cap\Aline ab$, $c'\in \Aline za\cap\Aline{a'}{b'}$, and $y\in\Aline b{c'}\cap\Aline {b'}c$. Assuming that the projective plane $\widehat X$ is $\mathsf{(P^{10}_{10})}$, we conclude that $y\in \Aline xz$. Then $\{a,b,c,a',b',c',o,x,y,z\}$ is a $3$-wide subliner of $\widehat X$. By Theorem~\ref{t:free-projectivization}(5), this subliner is contained in $X$, which contradicts the choice of the triangle $\{o,a,a'\}\not\subseteq X$.
\end{proof}

\begin{example} There exists a $4$-element projective liner $X$ whose free projectivization $\widehat X$ has the following properties:
\begin{enumerate}
\item $\widehat X$ is an $\w$-long projective plane;
\item $\widehat X$ is $\mathsf{(P_{13}^{13})}$ and $\mathsf{(\bar P_{13}^{13})}$;
\item $\widehat X$ is not $\mathsf{(P_{10}^{10})}$.
\end{enumerate}
\end{example}

\begin{proof} Choose any $4$-element set $X=\{a,b,a',b'\}$ endowed with the family of $2$-point lines $\mathcal L\defeq\{L\subseteq X:|L|=2\}$. Since $X$ contains disjoint lines, the free projectivization $\widehat X$ of $X$ is an $\w$-long projective plane, by Theorem~\ref{t:free-projectivization}(3). Proposition~\ref{p:projectivization<=>P1313} implies that the projective plane $\widehat X$ is $\mathsf{(P_{13}^{13})}$ and hence  $\mathsf{(\bar P_{13}^{13})}$. By Proposition~\ref{p:projectivizationP1010not}, the projective plane $\widehat X$ is not $\mathsf{(P_{10}^{10})}$. 
\end{proof}

\begin{example} There exists a $13$-element liner $X$ whose free projectivization $\widehat X$ has the following properties:
\begin{enumerate}
\item $\widehat X$ is an $\w$-long projective plane;
\item $\widehat X$ is $\mathsf{(\bar P_{13}^{13})}$;
\item $\widehat X$ is not $\mathsf{(P_{13}^{13})}$ and not $\mathsf{(P_{10}^{10})}$.
\end{enumerate}
\end{example}  

\begin{proof} Choose any $13$-element set $X=\{a,b,c,a',b',c',o,x,y,z,u,v,w\}$ and consider the families
$$
\begin{aligned}
\mathcal B\defeq\big\{&\{o,a,b,c\},\{o,a',b',c'\},\{u,a,b',x\},\{u,b,c',y\},\{v,a',b,x\},\{v,b',c,y\},\\&\{u,o,v\},\{a,a',w\},\{b,b',w\},\{c,c',w\}\big\}\quad\mbox{and}\\
\mathcal L\defeq\mathcal B&\cup\{L\subseteq X:|L|=2\wedge\forall B\in\mathcal B\;(L\not\subseteq B)\}.
\end{aligned}
$$

By Theorem~\ref{t:free-projectivization}(3), the free projectivization $\widehat X$ of the liner $(X,\mathcal L)$ is an $\w$-long projective plane and hence $X\ne\widehat X$. By Propositions~\ref{p:projectivization<=>P1313}, \ref{p:projectivization<=>barP1313}, \ref{p:projectivizationP1010not}, the projective plane $\widehat X$ is $\mathsf(\bar P_{13}^{13})$ and neither $(\mathsf P_{13}^{13})$ nor $(\mathsf P_{10}^{10})$.
\end{proof}

The results proved in this chapter show that for every projective liner the following implications between various modifications of the Desargues Axiom $\mathsf{(D)}$ and the Papus Axiom $\mathsf{(P)}$ hold.

$$
\xymatrix@C=20pt@R=18pt{
&&&&\mathsf{(D_2)}\ar@{<=>}[ddl]\ar@{=>}[rd]\ar@{=>}[r]&\mathsf{(P^{10}_{12})}\ar@{=>}[rdd]\\
&&&&
\mathsf{(D_3)}\ar@{<=>}[d]&\mathsf{(P^{9}_{10})}\ar@{=>}[ld]\ar@{=>}[rd]\\
\mathsf{(P)}\ar@{=>}[r]&\mathsf{(D)}\ar@{=>}[r]&\mathsf{(D_1)}\ar@{=>}[r]&\mathsf{(P^{9}_{10}\wedge P^{10}_9)}\ar@{=>}[r]&\mathsf{(P^{10}_{10})}\ar@{=>}[r]&\mathsf{(\bar P^{13}_{13})}&
\mathsf{(P^{13}_{13})}\ar@{=>}[l]\\
&&&&
\mathsf{(4P)}\ar@{<=>}[u]&\mathsf{(P^{10}_{9})}\ar@{=>}[lu]\ar@{=>}[ru]\\
&&&&\mathsf{(D_{II})}\ar@{<=>}[luu]\ar@{=>}[ru]\ar@{=>}[r]&\mathsf{(P^{12}_{10})}\ar@{=>}[ruu]\\
}
$$

This diagram suggests the following open problems.

\begin{problem} Is every $\mathsf{(P^{10}_{10})}$ projective liner $\mathsf{(P^{13}_{13})}$?
\end{problem}

\begin{problem} Are the axioms $\mathsf{(P^9_{10})}$ and $\mathsf{(P^{10}_9)}$ equivalent for projective liners?
\end{problem}

\begin{problem} Are the axioms $\mathsf{(P^{10}_{12})}$ and $\mathsf{(P^{12}_{10})}$ equivalent for projective liners?
\end{problem}

\begin{problem} Is a projective plane $X$ bi-Desarguesian if
\begin{enumerate}
\item $X$ is $\mathsf{(P^{10}_{12})}$ and $\mathsf{(P^{12}_{10})}$?
\item $X$ is $\mathsf{(P^{10}_{12})}$, $\mathsf{(P^{12}_{10})}$ and $\mathsf{(P^{9}_{10})}$?
\end{enumerate}
\end{problem} 

\begin{remark} By computer calculations, Ivan Hetman checked that all known non-Pappian finite projective plane are not $\mathsf{(\bar P_{13}^{13})}$.
This empirical result suggests the following conjecture.
\end{remark}

\begin{conjecture} A finite projective plane is Pappian if and only if it is $\mathsf{(\bar P_{13}^{13})}$.
\end{conjecture}

\chapter{Characteristics of projective liners}

In this section we study the characteristic range $\har[X^{\diamondsuit}]$ of a projective liner $X$. The characteristic range of $X$ is introduced as the range of the function $\har:X^\diamondsuit\to\IN$ assigning to every quadrangle $abcd$ in $X$ its characteristic $\har(abcd)$.
The characteristic of a quadrangle in a projective liner is introduced and studied in the next section.  The most difficult theorem of this chapter is Theorem$^\dag$~\ref{t:Blokhuis-Sziklai-Kantor-Penttila} claiming that every finite projective plane of order $p^2$ and characteristic range $\{p\}$ is Pappian.

\section{The characteristic of a quadrangle in a projective liner}

Let us recall that a \index{quadrangle}\defterm{quadrangle} in a liner $X$ is any quadruple $abcd$ of distinct points of $X$ such that $\|\{a,b,c,d\}\|=3=\|\{x,y,z\}\|$ for all distinct points $x,y,z\in\{a,b,c,d\}$. The points $a,c,b,d$ are called the \defterm{vertices} of the quadrangle $abcd$. By definition, no three vertices of a quadrangle are colinear. 

For a liner $X$, we denote by \index[notation]{$X^\diamondsuit$}\defterm{$X^\diamondsuit$} the family of all quadrangles in $X$. 

\begin{exercise} Let $X$ be a projective plane of order $n$. Prove that $$|X^\diamondsuit|=(n^2+n+1)\cdot(n^2+n)\cdot n^2\cdot(n^2-2n+1).$$
\end{exercise}

\begin{proposition}\label{p:quadrangle=>3-long} If a projective plane $X$ contains a quadrangle, then $X$ is $3$-long, $2$-balanced, and has a well-defined order $|X|_2-1$.
\end{proposition}

\begin{proof} Let $abcd$ is a quadrangle in a projective plane $X$ and let $Q\defeq\{a,b,c,d\}$. To show that $X$ is $3$-long, take any line $L$ in $X$. By the definition of a liner, the line $L$ contains two distinct points $x,y$.  Since $\|Q\|=3$, there exists a point $q\in Q\setminus L$. Since $abcd$ is a quadrangle, the intersection $\Aline xq\cap Q$ contains at most two points of the set $Q$. Since $q\in Q$, the intersection $\Aline xq\cap (Q\setminus \{q\})$ contains at most  one point of the set $Q$. By analogy we can prove that $\Aline yq\cap Q$ contains at most one point of the set $Q$. Then $Q\cap(\Aline xq\cup\Aline yq)$ contains at most three points of the set $Q$ and hence there exists a point $p\in Q\setminus (\Aline xq\cup\Aline yq)$. Since $X$ is a projective plane, there exists a point $z\in \Aline xy\cap \Aline pq$, witnessing that $|L|\ge|\{x,y,z\}|=3$ and hence the projective plane $X$ is $3$-long. By Corollary~\ref{c:Avogadro-projective}, the $3$-long projective plane $X$ is $2$-balanced and the cardinal $|X|_2$ is well-defined. Then the cardinal $|X|_2-1$ is the well-defined order of $X$. If $|X|_2$ is infinite, then $|X|_2-1\defeq|X|_2$.
\end{proof}


Let us recall that a subset $A$ of a liner $X$ is called \defterm{closed} in $X$ if $\Aline ab\cap\Aline cd\subseteq A$ for any points $a,b,c,d\in A$ such that $\Aline ab\ne\Aline cd$. Every subset $A$ of a liner $X$ is contained in the smallest closed subset $\langle A\rangle$ of $X$, called the \defterm{closure} of $A$. The closure $\langle A\rangle$ of $A$ is equal to the union $\bigcup_{n\in\w}A_n$ of the increasing sequence $(A_n)_{n\in\w}$ of sets defined by the recursive formula: $A_0=A$ and $A_{n+1}=A_n\cup \bigcup\{\Aline ab\cap\Aline cd:a,b,c,d\in A_n\;\wedge\;\Aline ab\ne\Aline cd\}$ for $n\in\w$. The closure $\langle A\rangle$ of any set $A$ in a projective liner $X$ is a projective subliner of $X$.  

Given a quadrangle $abcd$ in a liner $X$, it will be convenient to denote the closure $\langle\{a,b,c,d\}\rangle$ and the flat hull $\overline{\{a,b,c,d\}}$ of the set $\{a,b,c,d\}$ by $\langle abcd\rangle$ and $\overline{abcd}$, respectively. If $X$ is a projective liner, then $\langle abcd\rangle$ and $\overline{abcd}$ are $3$-long projective planes with well-defined orders, according to Proposition~\ref{p:quadrangle=>3-long}. If the order $q$ of the projective plane $\langle abcd\rangle$ is finite, then $|\langle abcd\rangle|=q^2+q+1$, by Corollary~\ref{c:projective-order-n}.

\begin{definition} The \index{characteristic of a quadrangle}\defterm{characteristic} $\har(abcd)$ of a quadrangle $abcd$ in a projective liner $X$ is defined as the order of its closure $\langle abcd\rangle$ in $X$. 
\end{definition}

Theorem~\ref{t:Bruck55} implies the following arithmetic properties of the characteristic of a quadrangle.

\begin{proposition}\label{p:har4-Bruck} Let  $abcd$ be a quadrangle in a projective liner $X$ of finite order $|X|_2=n$. For the number $p\defeq\har(abcd)$, one of the following holds: $p=n$, $p^2=n$, or $6\le p^2+p\le n$. Moreover, if $p\le 10$, then $p\in\{2,3,5,7,9\}$.
\end{proposition}


\begin{proposition}\label{p:Fano<=>har2} A quadrangle $abcd$ in a projective liner $X$ is Fano if and only if $|\langle abcd\rangle|=7$ if and only if $\har(abcd)=2$.
\end{proposition}

\begin{proof} Since the liner $X$ is projective, there exist unique points $x\in \Aline ab\cap\Aline cd$, $y\in \Aline ac\cap\Aline bd$, $z\in \Aline ad\cap\Aline bc$. Taking into account that no three points among $a,b,c,d$ are colinear, one can show that $|\{a,b,c,d,x,y,z\}|=7$. Moreover, $\{a,b,c,d,x,y,z\}\subseteq\langle abcd\rangle$ (because the set $\langle abcd\rangle$ is closed in $X$). 

If the quadrangle $abcd$ is Fano, then the points $x,y,z$ are contained in some line and the set $\{a,b,c,d,x,y,z\}$ is closed in $X$. Then $|\langle abcd\rangle|=|\{a,b,c,d,x,y,z\}|=7=2^2+2+1$ and $\har(abcd)=2$, according to Corollary~\ref{c:projective-order-n}.

If $\har(abcd)=2$, then $|\langle abcd\rangle|=2^2+2+1=7$ and hence $\langle abcd\rangle=\{a,b,c,d,x,y,z\}$. Since the set $\{a,b,c,d,x,y,z\}=\langle abcd\rangle$ is closed in the projective plane $\overline{abcd}\subseteq X$, there exists a point $z'\in \Aline xy\cap\Aline ad\subseteq \langle abcd\rangle=\{a,b,c,d,x,y,z\}$. The point $z'$ does not belong to the lines $\Aline ax\cup\Aline ay\supseteq\{a,b,c,x,y\}$. Assuming that $z'=d$, we conclude that $d=z'\in \Aline xy$ and hence $c\in \Aline xd=\Aline xy$ and $b\in \Aline by=\Aline xy$. Then the vertices $b,c,d$ of the quadrangle $abcd$ are colinear, which contradicts the definition of a quadrangle. This contradiction shows that $z'\ne d$ and hence $z'\in \{a,b,c,d,x,y,z\}\setminus\{a,b,c,d,x,y\}=\{z\}$. Therefore, $z=z'\in\Aline xy$ and the points $x,y,z$ are colinear, witnessing that the quadrangle $abcd$ is Fano. 
\end{proof}

\begin{remark}\label{rem:characteristic4} The definition of the characteristic of a quadrangle $abcd$ in a finite projective liner $X$ suggests a simple algorithm of calculation of the number $\har(abcd)$. Just calculate the sequence of sets $(Q_n)_{n\in\w}$ defined by the recursive formula $Q_0\defeq\{a,b,c,d\}$ and $Q_{n+1}=Q_n\cup\{\Aline xy\cap\Aline uv:x,y,u,v\in Q_n,\;\Aline xy\ne\Aline uv\}$ for all $n$ such that $Q_n\ne Q_{n+1}$. The algorithm stops at $n\in\w$ such that $Q_n=Q_{n+1}$ and then returns the set $\langle Q\rangle=Q_n=Q_{n+1}$. The unique number $q$ satisfying the equation $q^2+q+1=|\langle Q\rangle|$ is the characteristic of the quadrangle $abcd$. However this algorith is slow: it has the computational complexity $O(q^{10})$ and memory usage $O(q^2)$. There is much better algorithm of complexity $O(q^3)$ and memory usage $O(q)$, calculating the characteristic of the quadrangle $abcd$ as the order of the minimal subternar of the ternar $\Delta_{abcd}$ of the quadrangle $abcd$, see Remark~\ref{rem:fast-har4}
\end{remark}

\section{The ternar of a quadrangle in a projective liner}

\begin{definition} Let $X$ be a projective liner. The \index{ternar of a quadrangle}\index{quadrangle!ternar of}\defterm{ternar} $\Delta_{uowe}$ of a quadrangle $uowe\in X^\diamondsuit$ is defined as the ternar of the based projective plane $(\overline{uowe},uowe)$.
\end{definition}

We recall that the \defterm{characteristic} $\har(R)$ of a ternar $R$ is defined as the cardinality $|\underline{R}|$ of the minimal subternar $\underline{R}$ of the ternar $R$.

\begin{theorem}\label{t:har4<=>har-ternar4} The characteristic $\har(uowe)$ of a quadrangle $uowe$ in a projective liner $X$ equals the characteristic $\har(\Delta_{uowe})$ of the ternar $\Delta_{uowe}$ of the quadrangle $uowe$.
\end{theorem}

\begin{proof} Given a quadrangle $uowe$ in a projective liner $X$, consider its flat hull $\overline{uowe}$ in the projective liner $X$.  Let $h\in \Aline ou\cap\Aline we$ and $v\in \Aline ow\cap\Aline ue$ are the horizontal and vertical infinity points of the projective base $uowe$. Let $\Delta=\Aline oe\setminus\Aline hv$ be the ternar of the based projective plane $(\overline{uowe},uowe)$. By  Proposition~\ref{p:quadrangle=>3-long}, the projective plane $\overline{uowe}$ is $3$-long and hence its order $|\overline{uowe}|_2-1$ is well-defined.

If $|\overline{uowe}|_2=3$, then $\langle uowe\rangle=\overline{uowe}$, $|\langle uowe\rangle|=|\overline{uowe}|=7$ and $\har(uowe)=|\langle uowe\rangle|_2-1=2$. On the other hand $\Delta_{uowe}\defeq\{o,e\}$ is a minimal ternar and hence $\har(uowe)=|\langle uowe\rangle|_2-1=2=|\{o,e\}|=\har(\Delta_{uowe})$.

So, assume that $|\overline{uowe}|_2\ge 4$. In this case, the affine subliner $A\defeq \overline{uowe}\setminus\Aline hv$ of $X$ is a Playfair plane, according to Proposition~\ref{p:projective-minus-hyperplane} and Theorem~\ref{t:Playfair<=>}. Moreover, the ternar $\Delta_{uowe}$ of the quadrangle $uowe$ equals the ternar $\Delta_{uow}$ of the based affine plane $(A,uow)$. Consider the projective subplane $\langle uowe\rangle$ in $\overline{uowe}$ and the Playfair subplane $\underline A\defeq\langle uowe\rangle\setminus\Aline hv$ of the Playfair plane $A$. The Playfair planes $A$ and $\underline A$ have the same affine base $uow$ and hence the ternar $\underline{\Delta}=\Delta_{uow}\cap \langle uowe\rangle$ of the affine plane $\underline A$ is a subternar of the ternar $\Delta_{uow}=\Delta_{uowe}$ of the affine plane $A$. Then $\underline{\Delta}_{uow}\subseteq\underline{\Delta}$, by the minimality of the ternar $\underline{\Delta}_{uow}$.

Next, we show that $\underline{\Delta}\subseteq\underline{\Delta}_{uow}$. 
 Let $\underline{\Pi}_{uow}$ be the subset of the Playfair plane $A=\overline{uowe}\setminus\Aline hv$ consisting of all points whose coordinates in the affine base $uow$ belong to the minimal subternar $\underline{\Delta}_{uow}$ of the ternar $\Delta_{uow}$. The subliner $\underline{\Pi}_{uow}$ is a Playfair subplane of $A$. Observe that the set $$\overline{\underline{\Pi}}_{uow}\defeq\underline{\Pi}_{uow}\cup\bigcup\{\Aline xy\cap \Aline hv:x,y\in \underline{\Pi}_{uow}\}$$ contains the set $\{u,o,w,e\}$ and is a closed subset of the projective liner $X$. Taking into account that $\langle uowe\rangle$ is the smallest closed subset of $X$ that contains the set $\{u,o,w,e\}$, we conclude that $\langle uowe\rangle\subseteq\overline{\underline{\Pi}}_{uow}$ and hence $\underline{\Delta}=\Delta_{uow}\cap \langle uow\rangle\subseteq \Delta_{uow}\cap\overline{\underline{\Pi}}_{uow}=\Delta_{uow}\cap\underline{\Pi}_{uow}=\underline{\Delta}_{uow}$.

 Therefore, $\underline{\Delta}=\underline{\Delta}_{uow}$ and $\har(uowe)=|\Aline oe\cap\langle uowe\rangle|-1=|\underline\Delta|=|\underline{\Delta}_{uow}|=\har(\Delta_{uow})=\har(\Delta_{uowe})$.
 \end{proof}


\begin{definition} A quadrangle $uowe$ in a projective liner is called \index{minimal quadrangle}\index{quadrangle!minimal}\defterm{minimal} if its ternar $\Delta_{uowe}$ is minimal, i.e., contains no proper subternars. 
\end{definition}

Since a finite ternar $R$ is minimal if and only if $\har(R)=|R|$, we obtain the following characterization of minimal quadrangles. 

\begin{proposition} A quadrangle $uowe$ in a finite projective plane $X$ is minimal if and only if $\har(uowe)=|X|_2-1$.
\end{proposition}

\begin{remark}\label{rem:fast-har4}
Theorem~\ref{t:har4<=>har-ternar4} suggests a fast algorithm for the calculation of the characteristic of a quadrangle $uowe$ in a projective plane $X$ of order $m$. First, find the horizontal and vertical infinity points $h\in \Aline ou\cap\Aline we$ and $v\in\Aline ow\cap\Aline ue$ of the projective base $uowe$. Given any points $x,a,b\in \Delta\defeq\Aline oe\setminus \Aline hv$, find the point $T(x,a,b)\in\Delta$ by calculating the following points: $\alpha\in \Aline ah\cap\Aline uv$, $\beta\in\Aline bh\cap\Aline ov$, $\delta\in \Aline o\alpha\cap\Aline hv$, $y\in \Aline xv\cap\Aline \beta \delta$, and finally $T(x,a,b)\in \Aline yh\cap\Aline oe$. The calculation of the point $T(x,a,b)$ has complexity $O(1)$. Next, calculate the sequence of finite sets $(S_n)_{n\in\w}$ defined by the recursive formula $S_0=\{o,e\}$ and $S_{n+1}=T[S_n^3]$ for all $n\in\w$. The calculations proceed while $|S_n|<|S_{n+1}|$ and $|S_n|\le m^{\frac12}$. When $S_n=S_{n+1}$ or $|S_n|>m^{\frac12}$, the algorithm terminates. If $S_n=S_{n+1}$, then it returns $\har(uowe)\defeq|S_n|$. If $|S_n|>m^{\frac12}$, then it returns $\har(uowe)\defeq m$. This algorithm has computational complexity $O(\min\{q,m^{\frac12}\}^3)$, where $q=\har(uowe)$, and requires $O(q)$ memory, which is significantly better than the  algorithm of complexity $O(q^{10})$ and memory usage $O(q^2)$, discussed in Remark~\ref{rem:characteristic4}.
\end{remark}

\section{The characteristic range of a projective liner}

For a liner $X$ we denote by $X^\diamondsuit$ the family of all quadrangles in $X$.

\begin{definition} For a projective liner $X$, the range $\har[X^\diamondsuit]\defeq\{\har(uowe):uowe\in X^\diamondsuit\}$ of the function $\har:X^\diamondsuit\to \IN$ is called the \index{characteristic range}\index{projective liner!characteristic range of}\defterm{characteristic range} of  the projective liner $X$. The function $$\har_X^-:\har[X^\diamondsuit]\to \IN,\quad \har^-:n\mapsto|\{uowe\in X^\diamondsuit:\har(uowe)=n\}|$$is called the \index{characteristic level function}
\index{projective liner!characteristic level function of}\defterm{characteristic level function} of $X$.
\end{definition} 

Proposition~\ref{p:har4-Bruck} implies the following corollary describing the structure of the characterictic range of a projective liner of finite order.

\begin{corollary} For every projective liner $X$ of finite order $n$,
$$\har[X^\diamondsuit]\subseteq\{p\in\IN:p=n\;\vee\;p^2=n\;\vee\; 6\le p^2+p\le n\}.$$Moreover, every number $p\in \har[X^\diamondsuit]$ with $p\le 10$ belongs to the set $\{2,3,5,7,9\}$.
\end{corollary}

The following theorem shows that the characteristic range of a projective liner is just the union of the characteristic ranges of its affine subliners.

\begin{theorem}\label{t:har-proj=union-har-aff} Let $X$ be a projective liner and $\mathcal H$ be the family of all hyperplanes in $X$. Then $$\har[X^\diamondsuit]=\bigcup_{H\in\mathcal H}\har[(X\setminus H)^\vartriangle].$$
\end{theorem} 

\begin{proof}  If $p\in\har[X^\diamondsuit]$, then $p=\har(uowe)$ for some quadrangle $uowe$ in $X$. Find unique points $h\in\Aline ou\cap\Aline we$ and $v\in \Aline ow\cap\Aline ue$, and choose a maximal flat $H$ in $X$ such that $h,v\in H$ and $o\notin H$. The Exchange Property holding for strongly regular (= projective liners) ensures that $H$ is a hyperplane in $X$. By Proposition~\ref{p:projective-minus-hyperplane}, the subliner $X\setminus H$ of $X$ is affine, and by Theorem~\ref{t:affine=>Avogadro}, it is $2$-balanced. If $|X\setminus H|_2=2$, then the affine liner $X\setminus H$ is not Playfair and the ternar $\Delta_{ouw}$ of the triangle $uow$ in $X\setminus H$ is defined as the $2$-element field $\IF_2$. The minimality of the field $\IF_2$ ensures that $\har(uow)=|\IF_2|=2$. On the other hand, $\overline{uowe}$ is a projective plane of order $2$ and hence $\har(uowe)=|\langle uowe\rangle|=|\overline{uowe}|=2=\har(uow)\in \har[(X\setminus H)^\vartriangle]$. So, assume that $|X\setminus H|_2\ge 3$. In this case, $\overline{uowe}\setminus H$ is a Playfair plane and the ternar $\Delta_{uowe}$ of the quadrangle $uowe$ is defined as the ternar $\Delta_{uow}$ of the affine base $uow$ in the Playfair plane $\overline{uowe}\setminus H$. Theorem~\ref{t:har4<=>har-ternar4} ensures that $p=\har(uowe)=\har(\Delta_{uowe})=\har(\Delta_{uow})=\har(uow)\in \har[(X\setminus H)^\vartriangle]$. Therefore,  $\har[X^\diamondsuit]\subseteq\bigcup_{H\in\mathcal H}\har[(X\setminus H)^\vartriangle]$.
\smallskip

Now take any number $p\in \bigcup_{H\in\mathcal H}\har[(X\setminus H)^\vartriangle]$ and find a hyperplane $H\in\mathcal H$ and a triangle $uwe$ in the subliner $X\setminus H$ of $X$ such that $p=\har(uwe)$. By Corollary~\ref{c:line-meets-hyperplane}, there exist unique points $h\in\Aline ou\cap H$ and $v\in\Aline ow\cap H$. By the projectivity of $X$, there exists a point $e\in \Aline hw\cap\Aline vu$. Taking into account that $H$ is flat, we can see that $e\in X\setminus H$. Consider the quadrangle $uowe$ in the projective liner $X$ and its ternar $\Delta_{uowe}$. By Proposition~\ref{p:projective-minus-hyperplane}, the subliner $X\setminus H$ is affine. If $|X\setminus H|_2=2$, then $\Delta_{uowe}=\{o,e\}\cong\IF_2=\Delta_{uow}$. In this case, $p=\har(uwe)=\har(\IF_2)=2=\har(\Delta_{uowe})=\har(uowe)\in \har[X^\diamondsuit]$, by Theorem~\ref{t:har4<=>har-ternar4}. If $|X\setminus H|_2\ge 3$, then $(\overline{uow}\setminus H,uwe)$ is a based affine plane and $\Delta_{uowe}=\Delta_{uow}$. In this case $p=\har(uow)=\har(\Delta_{uow})=\har(\Delta_{uowe})=\har(uowe)\in\har[X^{\diamondsuit}]$, by Theorem~\ref{t:har4<=>har-ternar4}. Therefore,
$\bigcup_{H\in\mathcal H}\har[(X\setminus H)^\vartriangle]\subseteq\har[X^\diamondsuit]$.
 \end{proof}

Proposition~\ref{p:Fano<=>har2} implies the following two characterizations.

\begin{corollary}\label{c:Fano<=>2har} A projective liner $X$ contains a Fano quadrangle iff  $2\in\har[X^\diamondsuit]$.
\end{corollary}

%

\begin{corollary}\label{c:Fano<=>char2} A projective liner $X$ is Fano if and only if $\har[X^\diamondsuit]=\{2\}$.
\end{corollary}



Corollary~\ref{c:Fano<=>char2} and Definition~\ref{d:unicharacteristic-affine} motivate the following definition.

\begin{definition} A projective liner $X$ is called \index{unicharacteristic projective liner}\index{projective liner!unicharacteristic}\defterm{unicharacteristic} if its characteristic range $\har[X^\diamondsuit]$ is a singleton. The unique number in $\har[X^\diamondsuit]$ is called \index{the characteristic}\defterm{the characteristic} of the unicharacteristic projective liner $X$.
\end{definition}

Theorem~\ref{t:har-proj=union-har-aff} and Proposition~\ref{p:char-divides-order} imply the following corollary.

\begin{corollary}\label{c:unichar=>char-divide-order} If $X$ is a unicharacteristic projective plane of order $n$, then the characteristic $p$ of $X$ divides $n$ and $p-1$ divides $n-1$.
\end{corollary}

\begin{problem}\label{prob:unichar=>prime} Is the characteristic of every unicharacteristic finite projective plane a prime number?
\end{problem}

\begin{problem}\label{prob:unichar=>Pappian} Is every unicharacteristic finite projective plane Pappian?
\end{problem}

\begin{problem}\label{prob:unichar+prime=>Pappian} Is every unicharacteristic projective plane of prime order Pappian?
\end{problem}

\begin{remark} Since the characteristic of any Pappian liner is prime, the affirmative answer to Problem~\ref{prob:unichar=>Pappian} implies the affirmative answer to Problems~\ref{prob:unichar=>prime} and \ref{prob:unichar+prime=>Pappian}. 
\end{remark}

\begin{remark} The characteristic ranges and the characteristic level functions of all four projective planes of order $9$ are presented in the following table,  calculated by Ivan Hetman.
$$
\begin{array}{|c|c|c|}
\hline
\mbox{projective plane}&\mbox{characteristic range}&\mbox{characteristic level function}\\
\hline
\mbox{\tt desarg}&\{3\}&\{(3,42\,456\,960)\}\\
\mbox{\tt hall,dhall}&\{2,3,9\}&\{(2,8\,709\,120),(3,6\,065\,280),(9,27\,682\,560)\}\\
\mbox{\tt hughes}&\{2,3,9\}&\{(2,5\,660\,928),(3,6\,065\,280),(9,30\,730\,752)\}\\
\hline
\end{array}
$$
\end{remark}

\begin{remark} The characteristic ranges and the characteristic level functions of all 22 known projective planes of order $16$ are presented in the following table, taken from the website of Eric Moorhouse \cite{Moorhouse}.
$$
\begin{array}{|c|c|c|}
\hline
\mbox{projective plane}&\mbox{characteristic range}&\mbox{characteristic level function}\\
\hline
\mbox{\tt desarg}&\{2\}&\{(2,4277145600)\}\\
\mbox{\tt jowk,djowk}&\{2,16\}&\{(2,962260992),(16,3314884608)\}\\
\mbox{\tt hall,dhall}&\{2,16\}&\{(2,948326400),(16,3328819200)\}\\
\mbox{\tt semi4}&\{2,16\}&\{(2,858525696),(16,3418619904)\}\\
\mbox{\tt math,dmath}&\{2,16\}&\{(2,578285568),(16,3698860032)\}\\
\mbox{\tt dsfp,ddsfp}&\{2,16\}&\{(2,565125120),(16,3712020480)\}\\
\mbox{\tt lmrh,dlmrh}&\{2,16\}&\{(2,539578368),(16,3737567232)\}\\
\mbox{\tt bbh1}&\{2,16\}&\{(2,519450624),(16,3757694976)\}\\
\mbox{\tt demp,ddemp}&\{2,16\}&\{(2,505774080),(16,3771371520)\}\\
\mbox{\tt john,djohn}&\{2,16\}&\{(2,500097024),(16,3777048576)\}\\
\mbox{\tt bbs4,dbbs4}&\{2,16\}&\{(2,481904640),(16,3795240960)\}\\
\mbox{\tt bbh2,dbbh2}&\{2,16\}&\{(2,481582080),(16,3795563520)\}\\
\mbox{\tt semi2}&\{2,16\}&\{(2,453906432),(16,3823239168)\}\\
\hline
\end{array}
$$
\end{remark}

The results of calulations show that every known non-Desarguesian projective plane of order $\le16$ has a projective base whose ternar is minimal. This suggests the following definition.

\begin{definition} A projective plane is \index{miminal projective plane}\index{projective plane!minimal}\defterm{minimal} if it has a projective base $uowe$ whose ternar $\Delta_{uowe}$ is minimal.
\end{definition}

\begin{remark} A $2$-balanced finite projective plane $X$ is minimal iff $|X|_2-1\in\har[X^\diamondsuit]$.
\end{remark} 

\begin{problem} Is every non-Desarguesian finite projective plane minimal?
\end{problem}

Corollary~\ref{c:unichar=>char-divide-order} implies the following corollary.

\begin{corollary} Every unicharacteristic prime projective liner $X$ is minimal and every ternar of $X$ is minimal.
\end{corollary}

\begin{problem} Is a finite projective plane prime if all its ternars are minimal?
\end{problem}

\begin{remark}\label{rem:miniso-testO(n11)} The isomorphism test for minimal ternars described in Remark~\ref{rem:min-ternar-iso-test} suggests a simple isomorphism test for minimal finite projective planes. Given two minimal projective planes $X,Y$ of the same finite order $n$, fix any minimal projective base $uowe$ in $X$ and calculate its ternar $\Delta_{uowe}$. Then for every quadrangle $abcd$ in the projective plane $Y$, calculate its characteristic $\har(abcd)$ and if $\har(abcd)=n$ then using the isomorphism test for minimal ternars, find a (unique) isomorphism between the minimal ternars $\Delta_{uowe}$ and $\Delta_{abcd}$, if it exists. The isomorphism between the ternars $\Delta_{uowe}$ and $\Delta_{abcd}$ determines an isomorphism between the projective planes $X$ and $Y$, by Theorem~\ref{t:p-isomorphic<=>ternar-isomorphic}. If for all quadrangles $abcd$ in $Y$, the ternars $\Delta_{uowe}$ and $\Delta_{abcd}$ are not isomorphic, then the projective planes $X,Y$ are not isomorphic, according to Theorem~\ref{t:p-isomorphic<=>ternar-isomorphic}. This algorithm of recognizing (non)isomorphic minimal projective planes of order $n$ has computational complexity $O(n^{11})$ and memory usage $O(n^3)$. 
\end{remark}

\begin{remark}\label{rem:iso-testO(n11)} Under the assumption of minimality of all non-Desarguesian finite projective planes, the isomorphism test described in Remark~\ref{rem:miniso-testO(n11)} can be extended to an isomophism test for all finite projective planes. Given two finite projective planes $X,Y$ of the same order, first test if those planes are Desarguesian. The algorithm for testing the Desarguesianity is described in Remark~\ref{rem:Des-algorithm} and has  computation complexity $O(n^3)$. If both projective planes $X,Y$ are Desarguesian, then then they are isomorphic by Theorem~\ref{t:finite-proj-isomorphic<=>}. If one of them is Desarguesian and the other is not, then they are not isomorphic. If both of them are not Desarguesian, then they are minimal by the assumption and we can apply the algorithm of complexity $O(n^{11})$, desribed in Remark~\ref{rem:miniso-testO(n11)} for testing whether the projective planes $X,Y$ are isomorphic.
\end{remark}

\section{Characteristics of $4$-Pappian projective liners}

In this section we analyze the structure of the characteristic range $\har[X^\diamondsuit]$ of a $4$-Pappian projective liner $X$. We recall that a projecctive liner $X$ is \defterm{$4$-Pappian} if the closure of any $4$-element subset of $X$ is Pappian. In particular, the closure $\langle abcd\rangle$ of any quadrangle in $X$ is a $3$-long Pappian projective plane. In this case, the ternar $\Delta_{uowe}$ of the quadrangle $abcd$ is a field and its characteristic $\har(\Delta_{uowe})=\har(abcd)$ belongs to the set $\IP\cup\{\w\}$. By Theorems~\ref{t:invertible-add<=>4-Pappian}, \ref{t:p-invertible-add<=>} and Propositions~\ref{p:p-invertible-plus<=>}, \ref{p:p-invertible-puls<=>}, a projective liner $X$ is $4$-Pappian iff it is invertible-add iff
it is invertible-plus iff it is invertible-puls iff for every hyperplane $H$ in $X$, the affine subliner $X\setminus H$ of $X$ is invertible-add. This characterization allows as to deduce the following properties of $4$-Pappian projective liners from Theorem~\ref{t:har-proj=union-har-aff} and the corresponding properties of invertible-add Playfair planes. In particular, the following two theorems can be deduced from from Theorems~\ref{t:har-proj=union-har-aff}, \ref{t:har[X]<=>orders}, \ref{t:invertible-add=>pls-elementary}. 

\begin{theorem}\label{t:har[X]<=>orders-p} Let $X$ be a $4$-Pappian projective liner $X$. For a number   $p\in\IN\cup\{\w\}$, the following conditions are equivalent:
\begin{enumerate}
\item $p\in\har[X^\diamondsuit]$;
\item \mbox{$p$ is equal to the order of some element in the plus and puls loops of some ternar $R$ of $X$;}
\item $p$ is equal to the order of some element in the plus or puls loop of some ternar $R$ of $X$.
\end{enumerate}
\end{theorem} 

\begin{theorem}\label{t:4-Pappian=>pls-elementary} If a projective liner $X$ of finite order is $4$-Pappian, then for every ternar $R$ of $X$, the plus and puls loops $(R,+)$ and $(R,\!\puls\!)$ of the ternar $R$ are elementary.
\end{theorem}

Corollaries~\ref{c:inv-plus=>Boolean-paralelogram}, \ref{c:inv-puls=>Boolean-paralelogram} and \ref{c:Boolean<=>2har} imply the following corollary.

\begin{corollary}\label{c:4-Pappian=>2har} If a $4$-Pappian projective plane $X$ has even order, then $2\in\har[X^\diamondsuit]=2$ and $X$ contains a Fano quadrangle.
\end{corollary}

Theorems~\ref{t:har-proj=union-har-aff} and \ref{t:har[X]=divisors} imply  the following theorem describing the characteristic range of inversive-plus or inversive-puls finite projective planes.

\begin{theorems}\label{t:har[X]=divisors-p} Let $X$ be an inversive-plus or inversive-puls projective plane of finite order $n=|X|_2-1$.
\begin{enumerate}
\item If $n=120$, then $\har[X^\diamondsuit]=\{2,3\}$ or $\har[X^\vartriangle]=\{2,3,5\}$.
\item If $n=1080$, then $\har[X^\diamondsuit]=\{2,3\}$;
\item If $n\in\{16320,39000\}$, then $\har[X^\diamondsuit]=\{2,3,5\}$;
\item If $n\notin\{120,1080,16320,39000\}$, then $\har[X^\diamondsuit]$ is equal to the set of all prime divisors of the order $n$ of $X$.
\end{enumerate}
\end{theorems} 

Theorems~\ref{t:p+inversive}, \ref{t:p-inversive-puls}, \ref{t:not60}  and Corollary~\ref{c:inversive=>order} imply the following corollary describing possible orders of inversive-plus or inversive-puls projective planes.

\begin{corollarys}\label{c:inversive=>order-p} If a finite projective plane $X$ is  inversive-plus or inversive-puls, then its order $|X|_2-1$ is equal to one of the numbers:
\begin{enumerate}
\item $p^n$ for some prime number $p$ and some positive integer $n$;
\item $p^nq$ for some $n\in\IN$ and some odd primes $p,q$ such that $q$ divides $p^n-1$;
\item $2^np^m$ for $n\ge 3$, $m\in\IN$ and some odd prime number $p$ that divides $2^n-1$;
\item $2^np^mq$ for some $n,m\in\IN$ and some odd prime numbers $p,q$ such that $q$ divides $p^m-1$ and $pq$ divides $2^n-1$;
\item $2p^n$ for some odd prime number $p$ and some $n\in\IN$;
\item $4p^n$ for some odd prime number $p$ and some $n\in\IN$;
\item $2p^nq$ for some $n\in\IN$ and some prime number $p$ and odd prime number $q$ that divides $p^n-1$;
\item $12p^n$ for some $n\in\IN$ and some odd prime number $p\ge 5$ such that $3$ divides $p^n-1$;
\item $60 p^n$ for some $n\in\IN$ and some prime number $p\ge 7$ such that $30$ divides $p^n-1$;
\item $120,1080,16320$ or $39000$.
\end{enumerate}
\end{corollarys} 

Theorem~\ref{t:har[X]=divisors-p} and Corollary~\ref{c:inversive=>order-p} imply the following classification of characteristic ranges of inversive-plus or inversive-puls finite projective planes.

\begin{corollarys} If a finite projective plane $X$ is inversive-plus or inversive-puls, then $\har[X^\diamondsuit]\in\big\{\varnothing,\{p,q\},\{2,p,q\},\{2,3,5,p\}:p,q\in\IP\big\}$.
\end{corollarys}

\begin{remark} A well-known (and still unproven) Prime-Power Conjecture states that the order of any finite projective plane is a prime power. So, if this conjecture is true (at least for inversive-plus or inversive-puls projective planes), then the cases  \textup{(2)--(10)} of Corollary~\ref{c:inversive=>order-p} do not occur and $\har[X^\diamondsuit]\in\{\varnothing,\{p\}:p\in\IP\}$ for any inversive-plus or inversive-puls projective liner $X$.
\end{remark}

Theorem~\ref{t:har-proj=union-har-aff} and Proposition~\ref{p:unichar+inversive=>prime-power} imply the following characterization of unicharacteristic inversive-plus or inversive-puls finite projective planes.

\begin{propositions}\label{p:unichar+inversive=>prime-power-p} An inversive-plus or inversive-puls finite projective plane is unicharacteristic if and only if its order is a prime power.
\end{propositions}

Theorem~\ref{t:Pappian-completion}, \ref{t:invertible-add<=>4-Pappian} and Proposition~\ref{p:a-Pappian-prime-order}  imply the following corollary.

\begin{corollarys} For a finite projective plane $X$, the following conditions are equivalent:
\begin{enumerate}
\item $X$ is Pappian of prime order;
\item $X$ has prime order and is inversive-plus or inversive-puls;
\item $X$ is $4$-Pappian and minimal.
\end{enumerate}
\end{corollarys}

\section{Projective planes of square order}

The main result of this section is Theorem~\ref{t:Blokhuis-Sziklai-Kantor-Penttila} saying that projective planes of order $p^2$ with characteristic range $\{p\}$ are Pappian. For prime numbers $p$, this theorem was proved by \index[person]{Blokhuis}Blokhuis\footnote{
{\bf Aart Blokhuis} (born in 1956) is a Dutch mathematician specializing in finite geometry, discrete mathematics, and combinatorial design theory. He received his Ph.D. in 1983 from the Technische Universiteit Eindhoven under the supervision of Jacobus Hendricus van Lint, a leading figure in combinatorics. Blokhuis has spent his academic career at Eindhoven University of Technology, contributing widely to research on blocking sets, projective space configurations, and extremal combinatorial problems. His work blends geometric intuition with algebraic and combinatorial techniques, influencing both theory and applications. Over decades of research, he has authored numerous papers and supervised students who continued work in finite geometry and related areas.
} and \index[person]{Sziklai}Sziklai\footnote{
{\bf P\'eter Sziklai} is a Hungarian mathematician whose work lies at the intersection of finite geometry, combinatorics, and algebraic methods over finite fields. He is a professor at E\"otv\"os Lor\'and University (ELTE) in Budapest and has played a prominent role in Hungarian academic life, including senior university leadership positions. He received the degree Doctor of the Hungarian Academy of Sciences for his dissertation Applications of Polynomials over Finite Fields (2014), in which polynomial techniques are developed and applied to extremal problems in finite geometric structures. Sziklai is particularly known for results on point sets and curves over finite fields, including the upper bound for the number of rational points of plane curves often referred to as Sziklai’s bound, as well as for work on blocking sets, arcs, and related configurations in projective spaces. Through a substantial body of journal articles and collaborations, he has become one of the central figures in contemporary research on finite geometries, influencing both the algebraic and combinatorial directions of the subject.
} \cite{BlokhuisSziklai}, and in general case by \index[person]{Kantor}Kantor\footnote{
{\bf William M. Kantor} (born 1944) is an American mathematician noted for his influential work in finite group theory, finite geometries, and the computational aspects of these subjects. He earned his Ph.D. in 1968 at the University of Wisconsin with a thesis {\em $2$-Transitive Symmetric Designs} under the supervision of Peter Dembowski and R.H. Bruck, and has held faculty positions at the University of Illinois at Chicago and, from 1971 onward, at the University of Oregon. His research spans permutation groups, combinatorial designs, and symmetry in finite geometries, and has produced over 170 papers along with algorithms incorporated into computational algebra systems. Kantor has served as an invited speaker at the International Congress of Mathematicians (1998) and was named a Fellow of the American Mathematical Society in 2013.
 } and \index[person]{Penttila}Penttila\footnote{
{\bf Tim Penttila} (born in 1959) is an Australian mathematician internationally recognized for his deep contributions to finite geometry, incidence structures, and algebraic combinatorics. Over a career spanning several decades, he has worked on foundational problems in finite projective spaces, ovals and spreads in projective planes, generalized quadrangles, and related combinatorial structures, publishing well over 100 research papers and influencing the development of geometric and algebraic methods in discrete mathematics. Penttila held a long‑term academic appointment at the University of Western Australia and later served as Professor of Mathematics at Colorado State University before becoming Professor Emeritus at the University of Adelaide. His work includes collaborations on hyperovals, association schemes, and classification problems in finite geometries, and he co‑authored the book {\em Analytic Projective Geometry} (Cambridge University Press, 2023), which blends classical geometry with modern analytic approaches. A workshop titled {\em Finite Geometry: a Workshop in Honour of Tim Penttila} was held at the University of Adelaide in 2019 to celebrate his impact on the field.
} \cite{KantorPenttila}. 

The proof of this theorem is difficult and uses Theorem~\ref{t:Moufang<=>permutation} characterizing Moufang planes in terms of the permutation groups, and also a non-elementary Theorem~\ref{t:Grundhofer-Muller-Nagy} describing the structure of $2$-transitive finite permutation groups containing a sharply $2$-transitive subset.

First we introduce some notations that will simplify presentation of the proof.
From now on, we assume that $X$ is a projective plane of order $p^2$ and characteristic range $\har[X^\diamondsuit]=\{p\}$ for some number $p\in\IN$. Since the characteristic of any ternar is at least $2$, the number $p$ is $\ge 2$.

Denote by $\mathcal L_X$ be family of all lines in $X$. For a point $x\in X$ we denote by $$\mathcal L_x\defeq\{L\in\mathcal L_X:x\in L\}$$ the family of all lines in $X$ that contain the point $x$.  For a subset $\Pi\subseteq X$, we denote by $\Lambda_\Pi$ the set of lines $L$ in $X$ such that $|L\cap\Pi|>1$.


For every line $D$ in $X$ and every triangle $abc$ in $X\setminus D$, let $\langle abcD\rangle$ be the closure of the set $\{a,b,c\}\cup(D\cap(\Aline ab\cup\Aline ac\cup\Aline bc))$ in $X$.

\begin{lemma}\label{l:abcD=xyzD} For any line $D\in\mathcal L_X$ and any triangle $abc$ in $X\setminus D$, 
\begin{enumerate}
\item the subliner $\langle abcD\rangle$ of $X$ is a projective plane of order $p$ and cardinality $p^2+p+1$;
\item for every triangle $xyz$ in $\langle abcD\rangle\setminus D$, the equality $\langle xyzD\rangle=\langle abcD\rangle$ holds.
\end{enumerate}
\end{lemma}

\begin{proof} 1. Consider the unique points  $\alpha\in \Aline ac\cap D$, $\beta\in\Aline bc\cap D$, $d\in \Aline \alpha b\cap\Aline \beta a$, and observe that  $\langle abcD\rangle=\langle\{a,b,c,\alpha,\beta,d\}\rangle=\langle abcd\rangle$. The definition of the characteristic range $\har[X^\diamondsuit]=\{p\}$ ensures that 
$\langle abcD\rangle=\langle abcd\rangle$ is a projective plane of order $p$. By Corollary~\ref{c:projective-order-n}, this projective plane has cardinality $p^2+p+1$.

\begin{picture}(100,90)(-160,-15)

\put(0,0){\line(1,0){60}}
\put(0,0){\line(0,1){60}}
\put(0,60){\line(1,-3){20}}
\put(0,20){\line(3,-1){60}}
\put(0,60){\line(1,-1){60}}
\put(30,33){$D$}

\put(0,0){\circle*{3}}
\put(-7,-6){$c$}
\put(20,0){\circle*{3}}
\put(18,-7){$a$}
\put(60,0){\circle*{3}}
\put(63,-3){$\alpha$}
\put(0,20){\circle*{3}}
\put(-9,18){$b$}
\put(0,60){\circle*{3}}
\put(-3,63){$\beta$}
\put(15,15){\circle*{3}}
\put(18,15){$d$}
\end{picture}

2. Given any triangle $xyz$ in $\langle abcD\rangle$,  consider the unique points $\alpha\in \Aline xz\cap D$, $\beta\in\Aline yz\cap D$ and $\gamma\in \Aline \alpha y\cap\Aline x\beta$, and observe that $\langle xyzD\rangle=\langle xyz\gamma\rangle$. Since the set $\langle abcD\rangle$ is closed, $\alpha,\beta,\gamma\subseteq\langle abcD\rangle$ and hence $\langle xyzD\rangle\subseteq \langle abcD\rangle$. Since both sets $\langle xyzD\rangle$ and $\langle abcD\rangle$ have the same cardinality $p^2+p+1$, they coincide.
\end{proof} 

\begin{lemma}\label{l:Blokhuis-partition} Let $h\in X$ be a point and $H,Y,D\in\mathcal L_h$ be three distinct lines. For any distinct points $a,b\in H\setminus \{h\}$, the family $\mathcal P_{oxY}\defeq \{\langle abyD\rangle\cap(Y\setminus \{h\}):y\in Y\setminus \{h\}\}$ is a partition of the set $Y\setminus \{h\}$ into $p$ subsets of cardinality $p$.
\end{lemma}

\begin{proof} Given any point $y\in Y\setminus \{h\}$, observe that the set $\langle abyD\rangle\cap Y$ is a line in the projective plane $\langle abyD\rangle$ of order $p$, which implies that the set $\langle abyD\rangle\cap Y$ contains $p+1$ points and the set $\langle abyD\rangle_Y\defeq\langle abyD\rangle \cap(Y\setminus \{h\})$ contains $p$ points. 

Next, we show that $\mathcal P_{abY}$ is a partition of the set $Y\setminus \{h\}$. Assuming that for some points $y,u\in Y\setminus \{h\}$, the sets $\langle abyD\rangle_Y$ and $\langle abuD\rangle_Y$ have a common point $v$, we can apply Lemma~\ref{l:abcD=xyzD}(2) and conclude that $\langle abyD\rangle=\langle abvD\rangle=\langle abuD\rangle$ and hence $\langle abyD\rangle_Y=\langle abuD\rangle_Y$. 

Since the projective plane $X$ has order $p^2$, the line $Y$ contains $p^2+1$ points, which implies that the set $Y\setminus\{h\}$ contains $p^2$ points. Then the partition $\mathcal P_{abY}\defeq \{\langle abyD\rangle_Y:y\in Y\setminus \{h\}\}$ divides the set $Y\setminus \{h\}$ into $p$ sets of cardinality $p$.
\end{proof}

\begin{lemma}\label{l:Blokhuis1} Let $h\in X$ be a point, $H,Y,D\in\mathcal L_h$ be three distinct lines, and $a,b\in H\setminus \{h\}$ be two distinct points.
\begin{enumerate}
\item For any points $y,u\in Y\setminus \{h\}$ and $c\in \langle abyD\rangle\cap(H\setminus\{a,h\})$, we have $\langle abuD\rangle=\langle acuD\rangle$.
\item For any points $y,u\in Y\setminus\{h\}$, we have  $\langle abyD\rangle\cap H=\langle abuD\rangle\cap H$.
\end{enumerate}
\end{lemma}

\begin{picture}(100,80)(-170,-7)

\put(60,0){\line(-1,0){80}}
\put(60,0){\line(-3,1){75}}
\put(60,0){\line(-1,1){60}}
\put(-10,60){$D$}
\put(-25,22){$Y$}
\put(-33,-3){$H$}

\put(0,0){\circle*{3}}
\put(-2,-8){$a$}
\put(30,0){\circle*{3}}
\put(28,-10){$b$}
\put(60,0){\circle*{3}}
\put(63,-3){$h$}
\put(0,20){\circle*{3}}
\put(-1,25){$y$}
\put(27,11){\circle*{3}}
\put(28,13){$u$}
\end{picture}

\begin{proof} By Lemma~\ref{l:Blokhuis-partition}, for any distinct points $c,d\in H\setminus\{h\}$, the family $\mathcal P_{cdY}\defeq\{\langle cdyD\rangle_Y:y\in Y\setminus\{h\}\}$ consisting of sets $\langle cdyD\rangle_Y\defeq\langle cdyD\rangle\cap(Y\setminus\{h\})$ is a partition of the set $Y\setminus\{h\}$ into $p$ many sets of cardinality $p$.
\smallskip

1. Given any points $y,u\in Y\setminus\{h\}$ and $c\in \langle abyD\rangle\cap(H\setminus\{a,h\})$, we should prove that $\langle abuD\rangle=\langle acuD\rangle$.  If $u\in \langle abyD\rangle$, then $abu$ and $acu$ are triangles in the subplane $\langle abyD\rangle$. Applying Lemma~\ref{l:abcD=xyzD}(2), we conclude that $\langle abuD\rangle=\langle abyD\rangle=\langle acuD\rangle$. So, assume that $u\notin \langle abyD\rangle$. In this case, $\langle abuD\rangle_Y$ is a cell of the partition $\mathcal P_{abY}$, distinct from $\langle abyD\rangle_Y=\langle acyD\rangle_Y$ (the latter equality follows from Lemma~\ref{l:abcD=xyzD} because $acy$ is a triangle in $\langle abyD\rangle$). Then the $p$-element set $\langle abuD\rangle_Y$ is contained in the union $\bigcup\mathcal P_{acY}\setminus \langle acyD\rangle_Y$ of $p-1$ many sets of the partition $\mathcal P_{acY}$. By the Pigeonhole Principle, there exists a cell $\langle acwD\rangle_Y$ of the partition $\mathcal P_{acY}$ that has at least two  common points $s,t$ with the set $\langle abuD\rangle_Y$.
Applying Lemma~\ref{l:abcD=xyzD}(2), we conclude that
$$\langle abuD\rangle=\langle astD\rangle=\langle acwD\rangle=\langle acuD\rangle.$$

2. Given any points $y,u\in Y\setminus\{h\}$, we have to prove the equality $\langle abyD\rangle\cap H=\langle abuD\rangle\cap H$. Since the sets $\langle abyD\rangle\cap H$ and $\langle abuD\rangle\cap H$ have the same cardinality $p+1$, the equality will follow as soon as we check that every point $c\in\langle abyD\rangle\cap H$ belongs to the set $\langle abuD\rangle\cap H$. If $c\in\{a,h\}$, then $c\in\langle abuD\rangle\cap H$, by the definition of the subplane $\langle abuD\rangle$. So, assume that $c\notin\{a,h\}$. Then $c\in \langle acuD\rangle =\langle abuD\rangle$, by the preceding item.
\end{proof}

Observe that the dual projective plane $X^*$ to the projective plane $X$ also has order $p^2$ and characteristic range $\{p\}$. Applying Lemma~\ref{l:Blokhuis1}(2) to the dual projective plane $X^*\defeq(\mathcal L_X,\{\mathcal L_x:x\in X\})$, we obtain the following lemma. In this lemma, for three lines $A,B,C\in\mathcal L_X$ and point $d\in X\setminus(A\cup B\cup C)$, we denote by $\langle ABCd\rangle$ the closure of the set $\{d\}\cup(A\cap B)\cup (A\cap C)\cup(B\cap C)$ in $X$. For a subset $\Pi\subseteq X$, we denote by $\mathcal L_\Pi$ the set of lines $L$ in $X$ such that $|L\cap\Pi|>1$.

\begin{lemma}\label{l:Blokhuis2} Let $H\in X$ be a line, $h,y,d\in H$ be three distinct points, and $A,B\in \mathcal L_h\setminus \{H\}$ be two distinct lines. For any lines $Y,U\in \mathcal L_y\setminus\{H\}$, we have  $\mathcal L_{\langle ABYd\rangle}\cap \mathcal L_h=\mathcal L_{\langle ABUd\rangle}\cap \mathcal L_h$.
\end{lemma}

\begin{picture}(100,70)(-170,-15)

\put(-10,0){\line(1,0){80}}
\put(0,0){\line(2,1){55}}
\put(0,0){\line(1,1){40}}
\put(60,0){\line(-2,1){55}}
\put(60,0){\line(-1,1){40}}

\put(43,44){$A$}
\put(58,28){$B$}
\put(75,-3){$H$}
\put(12,43){$U$}
\put(-5,27){$Y$}

\put(0,0){\circle*{3}}
\put(-2,-10){$h$}
\put(25,0){\circle*{3}}
\put(22,-10){$d$}
\put(60,0){\circle*{3}}
\put(58,-8){$y$}

\end{picture}

Let $\mathcal B$ be the family of all projective subplanes of order $p$ in $X$. Elements of the family $\mathcal B$ are called {\em Baer subplanes} of $X$. For a Baer subplane $B\in\mathcal B$, let $\mathcal L_B$ be the family of lines $L$ in $X$ such that $L\cap B$ is a line in the Baer subplane $B$. For a line $H$ in $X$ and a point $h\in H$, let $$\mathcal B(H,h)\defeq\{B\in\mathcal B:h\in B\mbox{  and  }H\in\mathcal L_B\}\quad\mbox{and}\quad\mathcal S(H,h)\defeq\{B\cap H:B\in\mathcal B(H,h)\}.$$ 

\begin{lemma}\label{l:Baer-subplane-c} Let $H$ be any line in $X$ and $h\in H$ be any point. 
For any point $c\in X\setminus H$ and any set $S\in\mathcal S(H,h)$, there exists a Baer subplane $B\in \mathcal B(H,h)$ such that $c\in B$ and $B\cap H=S$.
\end{lemma}

\begin{proof} Given any set $S\in\mathcal S(H,b)$, find a Baer subplane $P\in \mathcal B(H,b)$ such that $S=P\cap H$. 
We recall that $\mathcal L_P$ denotes the family of all lines $L$ in $X$ such that $L\cap P$ is a line in the Baer plane $P$. Then $|L\setminus P|=(p^2+1)-(p+1)=p^2-p$. It is clear that $|\mathcal L_P|=p^2+p+1$ and for any distinct lines $L,\Lambda\in\mathcal L_P$, their intersection $L\cap \Lambda$ is a singleton in the projective plane $P\subseteq X$. Then the family $\{L\setminus P:L\in\mathcal L_P\}$ consists of $p^2+p+1$ pairwise disjoint sets of cardinality $p^2-p$. Consequently, $$\bigcup\mathcal L_P=|P|+(p^2+p+1)(p^2-p)=(p^2+p+1)(p^2-p+1)=(p^2+1)^2-p^2=p^4+p^2+1=|X|,$$ which means that $\mathcal L_P$ is a cover of the set $X$. So, there exists a line $L\in\mathcal L_P$ containing the point $c$. Since $c\notin H$, the line $L$ is distinct from the line $H$. 
Now consider two possible cases.
\smallskip

1. First assume that $h\in L$. In this case, choose any point $d\in L\cap P\setminus \{h,c\}$. Choose any line $V\in\mathcal L_h\cap\mathcal L_P\setminus\{L,H\}$ and any distinct points $a,b\in H\cap P\setminus \{h\}$.

\begin{picture}(100,90)(-170,-15)

\put(60,0){\line(-1,0){80}}
\put(60,0){\line(-3,1){75}}
\put(60,0){\line(-1,1){60}}
\put(-10,60){$V$}
\put(-25,22){$L$}
\put(-33,-3){$H$}

\put(0,0){\circle*{3}}
\put(-2,-8){$a$}
\put(30,0){\circle*{3}}
\put(28,-10){$b$}
\put(60,0){\circle*{3}}
\put(63,-3){$h$}
\put(0,20){\circle*{3}}
\put(-1,23){$d$}
\put(27,11){\circle*{3}}
\put(28,13){$c$}
\end{picture}

Lemma~\ref{l:abcD=xyzD}(2) ensures that $P=\langle abdV\rangle$, and Lemma~\ref{l:Blokhuis1}(2) implies that $\langle abcV\rangle\cap H=\langle abdV\rangle\cap H=P\cap H=S$. Then the Baer plane $B\defeq \langle abcV\rangle\in\mathcal B(H,h)$ has the required properties: $c\in B$ and $B\cap H=S$.
\smallskip

2. Next, assume that $h\notin L$. Find a unique point $b\in L\cap H\subseteq P$ and observe that $b\ne h$. Choose any point $d\in L\cap P\setminus\{b,c\}$. Choose any point $v\in P\cap\Aline hd$ and consider the line $V\defeq \Aline bv\in\mathcal L_P$. Choose any point $a\in P\cap H\setminus\{b,h\}$. 

\begin{picture}(100,100)(-170,-15)

\put(60,0){\line(-1,0){80}}
\put(60,0){\line(-3,1){80}}
\put(60,0){\line(-1,1){70}}
\put(0,0){\line(0,1){60}}
\put(-20,70){$V$}
\put(-30,25){$L$}
\put(-33,-3){$H$}

\put(0,0){\circle*{3}}
\put(-2,-10){$h$}
\put(30,0){\circle*{3}}
\put(28,-8){$a$}
\put(60,0){\circle*{3}}
\put(63,-3){$b$}
\put(0,20){\circle*{3}}
\put(2,21){$d$}
\put(27,11){\circle*{3}}
\put(28,13){$c$}
\put(0,60){\circle*{3}}
\put(1,61){$v$}
\end{picture}

Lemma~\ref{l:abcD=xyzD}(2) ensures that $\langle ahdV\rangle=P$, and Lemma~\ref{l:Blokhuis1}(2) implies that $\langle ahcV\rangle\cap H=\langle ahdV\rangle\cap H=P\cap H=S$. Then the Baer plane $B\defeq \langle ahcV\rangle\in\mathcal B(H,h)$ has the required properties: $c\in B$ and $B\cap H=S$.
\end{proof}

For any lines $H,V$ in a projective plane $X$ and any point $c\in X\setminus (H\cup V)$, we denote by $c_{V,H}:H\to V$ the central projection assigning to every point $x\in H$ the unique point of the intersection $\Aline cx\cap V$.

\begin{lemma} Let $o\in X$ be any point and $H,V\in \mathcal L_o$ be two distinct lines. For any point $c\in X\setminus(H\cup V)$ and every set $S\in\mathcal S(H,o)$, its image $c_{V,H}[S]$ belongs to the family $\mathcal S(V,o)$.
\end{lemma}

\begin{proof} Given any set $S\in\mathcal S(H,o)$, apply Lemma~\ref{l:Baer-subplane-c} and find a Baer subplane $P\in\mathcal B(H,o)$ such that $S=H\cap P$ and $c\in P$. Choose any  distinct points $a,b\in H\cap P\setminus\{o\}$ and $d\in \Aline co\cap P\setminus\{o,c\}$. Consider the lines $A\defeq\Aline ac$ and $B\defeq \Aline bc$.

\begin{picture}(100,110)(-170,-15)

\put(60,0){\line(-1,0){80}}
\put(60,0){\line(-3,1){80}}
\put(60,0){\line(-1,1){70}}
\put(0,0){\line(0,1){75}}
\put(0,0){\line(2,3){40}}
\put(-20,70){$V$}
\put(-30,25){$H$}
\put(-4,78){$A$}
\put(42,63){$B$}

\put(0,0){\circle*{3}}
\put(-2,-9){$c$}
\put(30,0){\circle*{3}}
\put(28,-10){$d$}
\put(60,0){\circle*{3}}
\put(63,-3){$o$}
\put(0,20){\circle*{3}}
\put(-7,13){$a$}
\put(11,16.5){\circle*{3}}
\put(8,20){$b$}
\put(0,60){\circle*{3}}
\put(2,62){$v$}
\end{picture}

Lemma~\ref{l:Blokhuis2} ensures that $\mathcal L_{\langle ABHd\rangle}\cap\mathcal L_c=\mathcal L_{\langle ABVd\rangle}\cap\mathcal L_c$. Consider the Baer plane $\Pi=\langle ABVd\rangle$ and observe that $\Pi\in\mathcal B(V,o)$. The inclusion $c\in P\cap\Pi$ and the equality
$$\mathcal L_P\cap\mathcal L_c=\mathcal L_{\langle ABHd\rangle}\cap\mathcal L_c=\mathcal L_{\langle ABVd\rangle}\cap\mathcal L_c=\mathcal L_\Pi\cap\mathcal L_c$$ imply $c_{V,H}[S]=c_{V,H}[P\cap H]=\Pi\cap V\in\mathcal S(V,o)$.
\end{proof}

Finally we can formulate and prove the main result of this section.

\begin{theorems}[Blokhuis, Sziklai, 2012; Kantor, Penttila, 2000]\label{t:Blokhuis-Sziklai-Kantor-Penttila} Every finite projective plane $X$ of order $p^2$ and characteristic range $\har[X^\diamondsuit]=\{p\}$ is Pappian.
\end{theorems}

\begin{proof} To derive a contradiction, assume that $X$ is not Pappian. By Theorem~\ref{t:Moufang-finite<=>}, $X$ is not Moufang. Since every projective plane of order $\le 4$ is Pappian, $p^2>4$ and hence $p\ge 3$.

Since the projective plane $X$ is not Moufang, we can apply Theorem~\ref{t:Moufang<=>permutation}(8) and find a point $o\in X$ such that for any distinct lines $H,V\in\mathcal L_o$ and point $a\in X\setminus (H\cup V)$, the permutation group $G\subseteq\Sym(H\setminus\{o\})$ generated by the sharply $2$-transitive set $\{c_{H,V}a_{V,H}:c\in X\setminus (H\cup V)\}$ is not of affine type. By Theorem~\ref{t:non-Thales-line-aff}, either $\Alt(H\setminus\{o\})\subseteq G$ or else $p^2=|H\setminus\{o\}|=24$ and $G$ is isomorphic to the Mathieu group $M_{24}$. Since $24$ is not a square of an integer, the second case is impossible and hence $\Alt(H\setminus\{o\})\subseteq G$. Since the alternating group on the set $H\setminus\{o\}$ of cardinality $p^2$ is $(p^2-2)$ transitive, it is $p$-transitive, which implies that every $p$-element subset of the set $H\setminus\{o\}$ belongs to the family $\mathcal S(H,o)$. This implies that
${C(p^2,p)}\le |\mathcal S(H,o)|$. On the other hand, $|\mathcal S(H,o)|\le|\mathcal B(H,o)|\le C(p,2)\cdot p\cdot C(p,2)$. Then 
$$\frac{(p^2)!}{p!\cdot(p^2-p)!}=C(p^2,p)\le C(p,2)\cdot p\cdot C(p,2)=\frac{p^3(p-1)^2}{4},$$which is not true for $p\ge 3$.  
This contradiction completes the proof of the theorem.
\end{proof}

\chapter{Pentagonal and pentagramic liners}

In the preceding chapters we have met many geometric configurations characterizing various algebraic properties of ternars of affine or projective liners. Those configurations suggest a general definition of a configuration axiom, having a dashed colored graph as a parameter. So, first we introduce and study dashed colored graphs. Then we study three special configutaion axioms determining pentagonal, propentagonal and pentagramic liners.

\section{Dashed colored graphs}

\begin{definition} A \index{graph}\defterm{graph} is a mathematical structure $(V,E)$ consisting of a set $V$ and a set $E$ of two-element subsets of $V$. Elements of the sets $V$ and $E$ are called \defterm{vertices} and \defterm{edges} of the graph $(V,E)$, respectively.   
A graph is 
\begin{itemize}
\item \index{finite graph}\index{graph!finite}\defterm{finite} if the set $V$ is finite;
\item \index{complete graph}\index{graph!complete}\defterm{complete} if $E=[X]^2\defeq\{e\subseteq X:|e|=2\}$.
\end{itemize}
 A \index{dashed graph}\index{graph!dashed}\defterm{dashed graph} is a graph $(V,E)$ with a distinguished edge $\bar e\in E$, called the \index{dashed edge}\defterm{dashed edge} of the dashed graph.  
\end{definition}

Every finite (dashed) graph $(V,E)$ can be drawn on the Euclidean plane as follows. Choose an injective function $f:V\to \IR^2$ such that for any distinct points $x,y,z\in V$, the points $f(x),f(y),f(y)$ are not collinear. For every edge $\{x,y\}\in E$ draw the segment $[f(x),f(y)]$ with end-points $f(x)$ and $f(y)$ in the plane $\IR^2$. Then the subset $\bigcup_{\{x,y\}\in E}[f(x),f(y)]$ of the plane is a visual representation of the graph $(V,E)$. If the graph $(V,E)$ has a distinguished edge $\bar e=\{x,y\}$, then we draw the segment $[f(x),f(y)]$ as dashed.

\begin{example} The complete graph $(\{1,2,3\},\{\{1,2\},\{2,3\},\{1,3\}\})$ with the distingished edge $\{1,2\}$ can be represented by the triangle with one dashed side.

\begin{picture}(40,55)(-180,-10)
\linethickness{0.8pt}
\put(0,0){\line(2,3){20}}
\put(40,0){\line(-2,3){20}}
\multiput(0,0)(3,0){13}{\line(1,0){2}}

\put(0,0){\circle*{3}}
\put(40,0){\circle*{3}}
\put(20,30){\circle*{3}}
\end{picture}
\end{example} 

\begin{definition} A \index{colored graph}\index{graph!colored}\defterm{colored graph} is a graph $(V,E)$ endowed an equivalence relation ${\parallel}\subseteq E\times E$ on its set of edges. This equvalence relation is called the \index{graph!coloring}\defterm{coloring} of the colored graph $(V,E)$. For every edge $e\in E$ its equivalence class $e_\parallel\defeq\{e'\in E:e'\parallel e\}$ is called \defterm{the color} of the edge $e$. The coloring $\parallel$ is uniquely determined by the partition $E/_{\parallel}\defeq\{e_\parallel:e\in E\}$ of $E$ into color classes, and hence $\parallel$ can be identified with the partition $E/_{\parallel}$. 
\end{definition}

Colored (dashed) graphs can be drawn on the Euclidean plane by coloring the segments representing the edges with some colors which coincide if and only if edges belong to the same color class.

\begin{example}\label{ex:Boolean-config} The complete graph 
$$\Boolean\defeq (\{1,2,3,4\},\{\{1,2\},\{1,3\},\{1,4\},\{2,3\},\{2,4\},\{3,4\}\})$$ endowed with the coloring $$\big\{\{\{1,2\},\{3,4\}\},\;\{\{2,3\},\{1,4\}\},\;\{\{1,3\},\{2,4\}\}\big\}$$ is called \index{Boolean colored graph}\defterm{the Boolean colored graph}.  This graph can be represented by the colored quadrangle


\end{example}

\begin{example}\label{ex:pentagon} The complete graph 
$$\Pentagon\defeq (\{1,2,3,4,5\},\{\{x,y\}:x,y\in\{1,2,3,4,5\}\;\wedge\;x\ne y\})$$ endowed with the coloring $$\big\{\{\{1,2\},\{5,3\}\},\;\{\{2,3\},\{1,4\}\},\;\{\{3,4\},\{2,5\}\},\{\{4,5\},\{3,1\}\},\{\{5,1\},\{4,2\}\}\big\}$$ is called \index{colored pentagon}\defterm{the colored pentagon}.  This graph can be represented by the colored pentagon

%
\end{example}

\begin{example} The colored graph $%
$ with vertices $V\defeq \{0,1,2,3,4,5,6\}$, edges $$E\defeq\big\{\{0,1\},\{1,2\},\{0,3\},\{3,4\},\{0,5\},\{5,6\},\{1,5\},\{2,6\},\{1,3\},\{2,4\},\{3,5\},\{4,6\}\big\}$$and coloring
$$\big\{\big\{\{0,1\},\{1,2\}\big\},\big\{\{0,3\},\{3,4\}\big\},\big\{\{0,5\},\{5,6\}\big\},\big\{\{1,5\},\{2,6\}\big\},\big\{\{1,3\},\{2,4\}\big\},\big\{\{3,5\},\{4,6\big\}\big\}$$
is called the \defterm{affine-Desarguesian colored graph}.
\end{example}

\begin{example} By Theorem~\ref{t:Proclus-lines}, every regular liner $X$ can be identified with the colored complete graph whose set of vertices is $X$ and two edges $e,g\in[X]^2$ have the same color iff $\overline e\parallel\overline g$, where $\overline e$ is the flat hull of the doubleton $e\subseteq X$.
\end{example}

Every colored graph carries a canonical structure of a liner. For distinct vertices $x,y\in V$ in a colored graph $\Gamma=(V,E)$, the line $\Aline xy$ through the points $x,y$ is defined as the set of all vertices $v\in V$ for which there exists a sequence of vertices $v_0,v_1,\dots,v_n\in V$ such that $v_0=x$, $v_1=y$, $v_n=v$, and $\{v_{i-1},v_i\}\parallel \{v_i,v_{i+1}\}$ for all $i$ with $0<i<n\}$.



\begin{definition} Let $\Gamma=(V,E)$ and $\Gamma'=(V',E')$ be two graphs. A function $f:V\to V'$ is called 
\begin{itemize}
\item an \index{graph!immersion}\defterm{immersion} of the graph $\Gamma$ into the graph $\Gamma'$ if for every edge $e\in E$ we have $f[e]\in E'$;
\item an \index{graph!isomorphism}\defterm{isomorphism} of the graphs $\Gamma$, $\Gamma'$ if the function $f$ is bijective and both maps $f$ and $f^{-1}$ are graph immersions.
\end{itemize}
\end{definition} 

\begin{exercise} Show that every graph immersion is injective.
\end{exercise}

\begin{definition} A graph $\Gamma=(V,E)$ is called \index{graph!edge-homogeneous}\defterm{edge-homogeneous} if for every edges $e,e'\in E$ there exists a graph isomorphism $f:\Gamma\to\Gamma'$ such that $f[e]=e'$.
\end{definition}

\begin{definition} Let $\Gamma=(V,E,\parallel)$ and $\Gamma'=(V',E',\parallel')$ be two colored graphs. A function $f:V\to V'$ is called 
\begin{itemize}
\item an \index{colored graph!immersion}\defterm{immersion} of the colored graph $\Gamma$ into $\Gamma'$ if $f$ is an immersion of the graph $(V,E)$ into the graph $(V',E')$ such that $\forall e,e'\in E\;\;(e\parallel e'\;\Leftrightarrow\;f[e]\parallel' f[e'])$;
\item an \index{colored graph!isomorphism}\defterm{isomorphism} of the colored graphs $\Gamma$ and $\Gamma'$ if $f$ is bijective and both maps $f$ and $f^{-1}$ are immersions of the colored graphs.
\end{itemize}
\end{definition} 

\begin{definition} Let $\Gamma=(V,E,\parallel,e)$ and $\Gamma'=(V',E',\parallel',e')$ be two dashed colored graphs. A function $f:V\to V'$ is called 
\begin{itemize}
\item an \index{dashed colored graph!immersion}\defterm{immersion} of the dashed colored graph $\Gamma$ into $\Gamma'$ if $f$ is an immersion of the colored graph $(V,E,\parallel)$ into the colored graph $(V',E',\parallel')$ such that $f[e]=e'$;
\item an \index{dashed colored graph!isomorphism}\defterm{isomorphism} of the dashed colored graphs $\Gamma$ and $\Gamma'$ if $f$ is bijective and both maps $f$ and $f^{-1}$ are immersions of the dashed colored graphs.
\end{itemize}
\end{definition}

\section{Configurations in liners}

\begin{definition} Let $\Gamma=(V,E,\parallel)$ be a colored graph. An injective function $f:V\to X$ to a liner $X$ is called 
a \index{$\Gamma$-configuration}\defterm{$\Gamma$-configuration} in $X$ if 
$$\forall e,e'\in E\;\;(e\parallel e'\;\Leftrightarrow\;\overline{f[e]}\parallel \overline{f[e']}).$$
\end{definition}

\begin{example} A regular liner $X$ contains a Boolean parallelogram if it has a $\Boolean$-configuration for the Boolean colored complete graph $\Boolean$ introduced in Example~\ref{ex:Boolean-config}.
\end{example}


\begin{definition}  Let $\bar\Gamma=(V,E,\parallel,\bar e)$ be a dashed colored graph. An injective function $f:V\to X$ to a liner $X$ is called
\begin{itemize}
\item a \index{partial $\bar\Gamma$-configuration}\defterm{partial $\bar\Gamma$-configuration} in $X$ if $\forall e,e'\in E\setminus\{\bar e\}\;\;(e\parallel e'\;\Leftrightarrow\;\overline{f[e]}\parallel \overline{f[e']})$;
\item a \index{full $\bar\Gamma$-configuration}\defterm{full $\bar\Gamma$-configuration} in $X$ if 
$\forall e,e'\in E\;\;(e\parallel e'\;\Leftrightarrow\;\overline{f[e]}\parallel \overline{f[e']})$.
\end{itemize}
\end{definition}

\begin{definition} Let $\bar\Gamma=(V,E,\parallel,\bar e)$ be a dashed colored graph. A liner $X$ is called a \index{$\bar\Gamma$-liner}\defterm{$\bar\Gamma$-liner} if every partial $\bar\Gamma$-confuguration in $X$ is a full $\bar\Gamma$-configuration.
\end{definition}

\begin{definition} Let $\bar\Gamma=(V,E,\parallel)$ be a colored graph. A liner $X$ is called a \index{$\Gamma$-liner}\defterm{$\Gamma$-liner} if $X$ is a $(\Gamma,\bar e)$-liner for {\em every} edge $\bar e\in E$.
\end{definition}

\begin{example} A regular liner $X$ is Boolean if and only if it is a $\Boolean$-liner for the Boolean colored graph $\Boolean$, introduced in Example~\ref{ex:Boolean-config}.
\end{example}

\begin{definition} A colored graph $\Gamma=(V,E,\parallel)$ is called \index{colored graph!unbiased}\defterm{unbiased} if a completely regular plane $X$ is a $\Gamma$-liner if and only if $X$ is a $(\Gamma,\bar e)$-liner for {\em some} edge $\bar e\in E$.
\end{definition} 

\begin{problem} Recognize all unbiased colored graphs of small cardinality.
\end{problem}

Examples of unbiased colored graphs are homogeneous colored graphs.

\begin{definition} A colored graph $(V,E,\parallel)$ is called \index{colored graph!homogeneous}\defterm{homogeneous} if for any two edges $\bar e,\ddot e\in E$ there exists a bijective function $f:V\to V$ such that
\begin{itemize}
\item $\forall u,v\in V\;\;(\{u,v\}\in E\;\Leftrightarrow\; \{f(u),f(v)\}\in E\})$;
\item $\forall e,g\in E\;\;(e\parallel g\;\Leftrightarrow\; f[e]\parallel f[g])$.
\end{itemize}
\end{definition}

\begin{exercise} Prove that every homogeneous colored graph is unbiased.
\end{exercise}

\begin{exercise} Find an example of a non-homogeneous unbiased colored graph.
\smallskip

{\em Hint:} Look at the colored graph
.
\end{exercise}

\begin{center}
\begin{tabular}{|c|c|}
\hline
{\bf \!Geometric property of an affine space\!}&{\bf $\Gamma$-liner for the colored graph $\Gamma$}\\
\hline
Boolean&\begin{picture}(11,11)(-1,1)
\linethickness{=0.8pt}
\put(0,0){\color{teal}\line(1,0){10}}
\put(0,10){\color{teal}\line(1,0){10}}
\put(0,0){\color{cyan}\line(0,1){10}}
\put(10,0){\color{cyan}\line(0,1){10}}
\put(0,0){\color{red}\line(1,1){10}}
\put(0,10){\color{red}\line(1,-1){10}}
\put(0,0){\circle*{1.5}}
\put(10,0){\circle*{1.5}}
\put(0,10){\circle*{1.5}}
\put(10,10){\circle*{1.5}}
\end{picture}\\
\hline
Pappian&
\begin{picture}(20,12)(-10,1)
\linethickness{=0.6pt}
\put(-9,0){\color{orange}\line(1,2){5}}
\put(9,0){\color{brown}\line(-1,2){5}}
\put(-9,0){\color{red}\line(1,0){18}}
\put(9,0){\color{violet}\line(-5,2){15}}
\put(6,6){\color{violet}\line(-5,2){10}}
\put(-9,0){\color{blue}\line(5,2){15}}
\put(-6,6){\color{blue}\line(5,2){10}}
\put(-4,10){\color{red}\line(1,0){8}}

\put(-4,10){\circle*{1.4}}
\put(-6,6){\circle*{1.4}}
\put(-9,0){\circle*{1.4}}
\put(4,10){\circle*{1.4}}
\put(6,6){\circle*{1.4}}
\put(9,0){\circle*{1.4}}

\end{picture}\\
\hline
Desarguesian&\begin{picture}(18,14)(-9,-3)
\linethickness{=0.8pt}
\put(0,0){\color{brown}\line(0,1){8}}
\put(0,0){\color{olive}\line(2,-1){8}}
\put(0,0){\color{purple}\line(-2,-1){8}}
\put(-8,-4){\color{teal}\line(1,0){16}}
\put(-4,-2){\color{teal}\line(1,0){8}}
\put(-8,-4){\color{cyan}\line(2,3){8}}
\put(-4,-2){\color{cyan}\line(2,3){4}}
\put(8,-4){\color{blue}\line(-2,3){8}}
\put(4,-2){\color{blue}\line(-2,3){4}}
\put(0,0){\circle*{1.4}}
\put(0,0){\circle*{1.4}}
\put(0,4){\circle*{1.4}}
\put(0,8){\circle*{1.4}}
\put(4,-2){\circle*{1.4}}
\put(8,-4){\circle*{1.4}}
\put(-4,-2){\circle*{1.4}}
\put(-8,-4){\circle*{1.4}}
\end{picture}\\
\hline

\end{tabular}
\end{center}

\section{Regular pentagons and golden ratios}

\begin{definition} A \index{pentagon}\defterm{pentagon} in a liner $X$ is any 5-tuple $p_0p_1p_2p_3p_4\in X^5$ of distinct points such that $\|\{x,y,z\}\|=3=\|\{p_0,p_1,p_2,p_3,p_4\}\|$ for any distinct points $x,y,z\in \{p_0,p_1,p_2,p_3,p_4\}$.
\end{definition}

In the following definition we endow the number $5\defeq\{0,1,2,3,4\}$ with the binary operation $\oplus:5\times 5\to 5$, $\oplus:(x,y)\mapsto x\oplus y$, of addition modulo $5$, whose Cayley table looks as follows.
$$
\begin{array}{c|ccccc}
\oplus&0&1&2&3&4\\
\hline
0&0&1&2&3&4\\
1&1&2&3&4&0\\
2&2&3&4&0&1\\
3&3&4&0&1&2\\
4&4&0&1&2&3
\end{array}
$$  

\begin{definition} A pentagon $p_0p_1p_2p_3p_4$ in a liner $X$ is defined to be \index{regular pentagon}\index{pentagon!regular}\defterm{regular}\\ if $\Aline {p_{i}}{p_{i\oplus 1}}\cap\Aline {p_{i\oplus2}}{p_{i\oplus4}}=\varnothing$ for all $i\in 5$.

%
\end{definition}

The existence of regular pentagons in Moufang Playfair liners is equivalent to the existence of golden ratios in their ternars. 

\begin{definition} An element $g$ of a biloop  is called a \index{golden ratio}\defterm{golden ratio} if $g{\cdot}g=1+g$.
\end{definition}

\begin{proposition}\label{p:pentagon=>golden} If a Moufang Playfair liner contains a regular pentagon $p_0p_1p_2p_3p_4$, then the ternar of the based affine plane $(\overline{\{p_0,p_1,p_2\}},p_0p_1p_2)$ contains a golden ratio.
\end{proposition}

\begin{proof} It follows from $\Aline {p_0}{p_1}\cap\Aline {p_2}{p_4}=\varnothing$ that the triple $uow=p_0p_1p_2$ is an affine base in the plane $\Pi\defeq\overline{\{p_0,p_1,p_2,p_3,p_4\}}\subseteq X$. Let $e$ be the diunit of the affine base $uow$ and $\Delta=\Aline oe$ be the ternar of the based affine plane $(\Pi,uow)$. Since the Playfair liner $X$ is Moufang, the ternar $\Delta$ is linear, distributive, associative-plus, commutative-plus, and alternative-dot, by Theorem~\ref{t:inversive-dot} and \ref{t:affine-Moufang<=>}. For any elements $a,b\in\Delta$, let $L_{a,b}$ be the line in the based affine plane $(\Pi,uow)$, given by the equation $y=x_{\times}a_+b=x{\cdot}a+b$.

Since the plus loop $(\Delta,+)$ is commutative, for every $x\in \Delta$ there exists a unique element $\bar x\in\Delta$ such that $x+\bar x=0=\bar x+x$.  The right-distributivity of the biloop $\Delta$ ensures that for every $x\in\Delta$ we have
$$x+\bar x=o=(e+\bar e)\cdot x=e{\cdot}x+\bar e{\cdot}x=x+\bar e{\cdot}x$$ and hence $\bar x=\bar e{\cdot}x$, by the cancellativity of the plus loop $(\Delta,+)$. On the other hand, the left-distributivity of $R$ implies $x+\bar x=o=x{\cdot}(e+\bar e)=x{\cdot}e+x{\cdot}\bar e=x+x{\cdot}\bar e$ and hence $\bar x=x{\cdot}\bar e$.

It follows from $\Aline ow\cap\Aline{p_0}{p_3}=\Aline {p_1}{p_2}\cap\Aline {p_3}{p_0}=\varnothing=\Aline {p_0}{p_1}\cap\Aline {p_2}{p_4}=\Aline ou\cap\Aline {p_2}{p_4}$ that the points $p_3$ and $p_4$ have coordinates $(e,g)$ and $(h,e)$ for some points $g,h\in\Delta$. Consider the line $L_{\bar e,e}\subset \Pi$ given by the equation $y=x_\times \bar e_+e=x\cdot\bar e+e$. Since $o{\cdot}\bar e+e=e$ and $e{\cdot}\bar e+e=o$, the points $w=p_2$ and $u=p_0$ belong to the line $L_{\bar e,e}$. Next, consider the point $g^+\defeq g+e\in\Delta$ and consider the line $L_{\bar e,g^+}$ given by the equation $y=x_\times \bar e_+g^+=x{\cdot}\bar e+g^+$. Since $e{\cdot}\bar e+g^+=\bar e+e+g=g$, the point $p_3$ (with coordinates $(e,g)$) belongs to the line $L_{\bar e,g^+}$. Since $\Aline {p_3}{p_4}\parallel\Aline {p_0}{p_2}=L_{\bar e,e}\parallel L_{\bar e,g^+}$ and $p_3\in \Aline {p_3}{p_4}\cap L_{\bar e,g^+}$, the parallel lines $\Aline{p_3}{p_4}$ and $L_{\bar e,g^+}$ coincide and hence the coordinates $(h,e)$ of the point $p_4$ satisfy the equation $y=x\cdot\bar e+g^+$ of the line $L_{\bar e,g^+}=\Aline{p_3}{p_4}$. Then $e=h{\cdot}\bar e+g^+=\bar h+g+e$ and hence $\bar h+g=o=\bar h+h$ and $h=g$. Therefore, the points $p_3$ and $p_4$ have coordinates $(e,g)$ and $(g,e)$, respectively.

\begin{picture}(100,135)(-150,-20)
{\linethickness{0.8pt}
\put(0,0){\color{teal}\line(1,0){50}}
\put(0,50){\color{teal}\line(1,0){81}}
\put(0,0){\color{cyan}\line(0,1){50}}
\put(50,0){\color{cyan}\line(0,1){81}}
\put(50,81){\color{violet}\line(1,-1){31}}
\put(50,0){\color{magenta}\line(31,50){31}}
\put(0,50){\color{blue}\line(50,31){50}}
\put(0,0){\color{magenta}\line(50,81){50}}
\put(0,0){\color{blue}\line(81,50){81}}
\put(0,50){\color{violet}\line(1,-1){50}}
}

\put(0,0){\color{teal}\line(1,0){50}}
\put(0,50){\color{teal}\line(1,0){81}}
\put(0,0){\color{cyan}\line(0,1){50}}
\put(50,0){\color{cyan}\line(0,1){81}}
\put(50,81){\color{violet}\line(1,-1){45}}
\put(96,28){\color{violet}$L_{\bar e,g^+}$}
\put(50,0){\color{magenta}\line(31,50){45}}
\put(96,75){\color{magenta}$L_{g,\bar g}$}
\put(0,50){\color{blue}\line(50,31){70}}
\put(73,95){\color{blue}$L_{g^{-1},e}$}
\put(0,0){\color{magenta}\line(50,81){65}}
\put(66,108){\color{magenta}$L_{g,o}$}
\put(0,0){\color{blue}\line(81,50){100}}
\put(103,60){\color{blue}$L_{g^{-1},o}$}
\put(0,50){\color{violet}\line(1,-1){60}}
\put(62,-15){\color{violet}$L_{\bar e,e}$}
\put(0,0){\color{red}\line(1,1){80}}
\put(83,82){\color{red}$\Delta$}

\put(50,50){\circle*{3}}
\put(52,43){$e$}
\put(0,0){\circle*{3}}
\put(-8,0){$o$}
\put(1,-7){$p_1$}
\put(50,0){\circle*{3}}
\put(44,-9){$u$}
\put(54,-1){$p_0$}
\put(0,50){\circle*{3}}
\put(-11,45){$p_2$}
\put(-5,54){$w$}
\put(50,81){\circle*{3}}
\put(55,79){$p_3$}
\put(81,50){\circle*{3}}
\put(85.5,47.8){$p_4$}
\end{picture}

Since $(\Delta,\cdot)$ is a Moufang loop, there exists a unique element $g^{-1}\in\Delta$ such that $g^{-1}{\cdot}g=e=g{\cdot} g^{-1}$. Consider the parallel lines $L_{g^{-1},o}$ and $L_{g^{-1},e}$, and observe that $\Aline {p_1}{p_4}=L_{g^{-1},o}$ and $p_2=w\in L_{g^{-1},e}$. Since $\Aline{p_2}{p_3}\parallel \Aline {p_1}{p_4}=L_{g^{-1},o}\parallel L_{g^{-1},e}$, the parallel lines $\Aline {p_2}{p_3}$ and $L_{g^{-1},e}$ coincide and hence the coordinates $(e,g)$ of the point $p_3$ satisfy the equation $y=x{\cdot} g^{-1}+e$ of the line $L_{g^{-1},e}$. Then $g=e{\cdot}g^{-1}+e=g^{-1}+e$ and $$g\cdot g=(g^{-1}+e)\cdot g=g^{-1}{\cdot}g+e{\cdot}g=e+g,$$witnessing that $g$ is a golden ratio in the ternar $\Delta$.
\end{proof}

An element $g$ of a biloop $R$ is called \defterm{left-distributive} if $g{\cdot}(x+y)=g{\cdot}x+g{\cdot} y$ for all $x,y\in R$. Observe that a biloop is left-distributive if and only if every its element is left-distributive. 

\begin{proposition}\label{p:golden=>pentagon} If some ternar $R$ of a Thalesian Playfair liner $X$ contains a left-distributive golden ratio, then the liner $X$ contains a regular pentagon.
\end{proposition}

\begin{proof} Assume that some ternar $R$ of a Thalesian Playfair liner $X$ contains a left-distributive golden ratio. Find a plane $\Pi\subseteq X$ and an affine base $uow$ in $\Pi$ whose ternar $\Delta$ is isomorphic to the ternar $R$. Then the ternar $\Delta$ contains a left-distributive golden ratio $g$. By Theorems~\ref{t:VW-Thalesian<=>quasifield} and \ref{t:diagonal-trans=>ass-plus}, the ternar $\Delta$ is linear, associative-plus, commutative-plus, and right-distributive. Let $e\in\Delta$ be the diunit of the affine base $uow$. For every $a,b\in\Delta$, let $L_{a,b}$ be the line in the based affine plane $(\Pi,uow)$ defined by the equation $y=x_\times a_+b=x{\cdot}a+b$.  Since the plus loop of the ternar $\Delta$ is commutative-plus, for every $x\in\Delta$ there exists an element $\bar x\in\Delta$ such that $x+\bar x=o=\bar x+x$. The equality $g{\cdot}g=e+g$ implies $g{\cdot}g+\bar g=(e+g)+\bar g=e+(g+\bar g)=e+o=e$.

Consider the pentagon $p_0p_1p_2p_3p_4$ consisting of points $p_0,p_1,p_2,p_3,p_4$ with coordinates $(e,o)$, $(o,o)$, $(o,e)$, $(e,g)$, $(g,e)$, respectively. We claim that this pentagon is regular. It is clear that $\Aline {p_0}{p_1}\cap\Aline {p_2}{p_4}=\varnothing=\Aline {p_1}{p_2}\cap \Aline {p_0}{p_3}$.

To prove that $\Aline {p_4}{p_0}\cap\Aline {p_1}{p_3}=\varnothing$, consider the parallel lines $L_{g,o}$ and $L_{g,\bar g}$. Since $o{\cdot}g+o=o$ and $e{\cdot}g+\bar g=o$, the points $p_1$ and $p_3$ (with coordinates $(o,o)$ and $(e,g)$) belong to the line $L_{g,o}$. 
Since $e{\cdot}g+\bar g=o$ and $g{\cdot}g+\bar g=e$, the points $p_0$ and $p_4$ (with coordinates $(e,o)$ and $(g,e)$) belong to the line $L_{g,\bar g}$. Then $\Aline {p_1}{p_3}\cap \Aline{p_0}{p_4}=L_{g,o}\cap L_{g,\bar g}=\varnothing$.

To prove that $\Aline {p_2}{p_3}\cap\Aline {p_4}{p_1}=\varnothing$, find a unique element $g^{-1}\in \Delta$ such that $g{\cdot}g^{-1}=e$ and consider the parallel lines $L_{g^{-1},o}$ and $L_{g^{-1},e}$. Since the element $g$ is left-distributive, $g{\cdot}(g^{-1}+e)=g{\cdot}g^{-1}+g{\cdot}e=e+g=g{\cdot}g$ and hence $g^{-1}+e=g$, by the canclellativity of the dot loop of the ternar $R$. Since $o{\cdot}g^{-1}+o=o$ and $g{\cdot}g^{-1}+o=e$, the points $p_1$ and $p_4$ (with coordinates $(o,o)$ and $(g,e)$) belong to the line $L_{g^{-1},o}$.  Since $o{\cdot}g^{-1}+e=e$ and $e{\cdot}g^{-1}+e=g^{-1}+e=g$, the points $p_2$ and $p_3$ (with coordinates $(o,e)$ and $(e,g)$) belong to the line $L_{g^{-1},e}$. Then $\Aline {p_2}{p_3}\cap \Aline{p_1}{p_4}=L_{g^{-1},e}\cap L_{g^{-1},o}=\varnothing$.  

To prove that $\Aline{p_3}{p_4}\cap\Aline{p_0}{p_2}=\varnothing$,  consider the parallel lines $L_{\bar e,e}$ and $L_{\bar e,e{+}g}$. Since the element $g$ is left-distributive, $g+\bar g=o=g\cdot (e+\bar e)=g{\cdot}e+g{\cdot}\bar e=g+g{\cdot}\bar e$ and hence $\bar g=g{\cdot}\bar e$. Since $e{\cdot}\bar e+e=o$ and $o{\cdot}\bar e+e=e$, the points $p_0$ and $p_2$ (with coordinates $(e,o)$ and $(o,e)$) belong to the line $L_{\bar e,e}$.  Since $e{\cdot}\bar e+(e+g)=(\bar e+e)+g=o+g=g$ and $g{\cdot}\bar e+e+g=\bar g+g+e=o+e=e$, the points $p_3$ and $p_4$ (with coordinates $(e,g)$ and $(g,e)$) belong to the line $L_{\bar e,e+g}$. Then $\Aline {p_3}{p_4}\cap \Aline{p_0}{p_2}=L_{\bar e,e{+}g}\cap L_{\bar e,e}=\varnothing$.

Therefore, the pentagon $p_0p_1p_2p_3p_4$ is regular.
\end{proof}

\begin{corollary}\label{c:MPcontains5<=>golden} A Moufang Playfair liner contains a regular pentagon if and only if it has a ternar containing a golden ratio.
\end{corollary}

\begin{proof} The ``only if'' part follows from Proposition~\ref{p:pentagon=>golden}. To prove the ``if'' part, assume that some ternar $R$ of a Moufang Playfair liner $X$ contains a golden ratio $g$. By Theorem~\ref{t:affine-Moufang<=>} and Corollary~\ref{c:para-Desarguesian=>Thalesian}, the Moufang Playfair liner $X$ is Thalesian and its ternar $R$ is left-distributive. Then the golden ratio $g$ is left-distributive and we can apply Proposition~\ref{p:golden=>pentagon} to find a regular pentagon in the liner $X$.
\end{proof}

\begin{corollary} An affine space contains a regular pentagon if its scalar corps contains a golden ratio.
\end{corollary}

\begin{proof} Assume that the scalar corps $\IR_X$ of an affine space $X$ contains a golden ratio $g$. Then $g{\cdot}g=1+g$ and hence $g\ne 0\ne 1$ and $\IR_X\ne\{0,1\}$. By Theorem~\ref{t:RXne01=>paraD}, the affine space $X$ is Thalesian. Let $R$ be any ternar of the affine space $X$. By Theorem~\ref{t:Ker(R)=RPi}, the scalar corps $\IR_X$ is isomorphic to the kernel $\Ker(R)\defeq\{a\in R:\forall x,y\in R\;\;a{\cdot}(x{\cdot}y)=(a{\cdot}x){\cdot}y\;\wedge\;a{\cdot}(x+y)=a{\cdot}x+a{\cdot}y\}$ of $R$. Since the corps $\Ker(R)$ is isomorphic to $\IR_X$, there exists a golden ratio $r\in \Ker(R)$. By definition, all elements of $\Ker(R)$ are left-distributive in $R$. In particular, the golden ratio $r$ is left-distributive. By Proposition~\ref{p:golden=>pentagon}, the Thalesian Playfair liner $X$ contains a regular pentagon.
\end{proof}

\begin{remark} By computer calculations, Ivan Hetman established that the following finite affine planes of order 9 contain a regular pentagon, possess a ternar with a (left-distributive) golden ratio, or their scalar corps possess a golden ratio.

$$
\begin{array}{|l|c|c|c|c|}
\hline
\mbox{affine}\atop\mbox{plane:}&\mbox{contains}\atop{\mbox{a regular}\atop\mbox{pentagon}}&\mbox{has a ternar with}\atop\mbox{a golden ratio}&\mbox{has a ternar with}\atop{\mbox{a left-distributive}\atop\mbox{golden ratio}}&\mbox{the scalar corps}\atop\mbox{has a golden ratio}\\
\hline
\mbox{\tt Desarg:}&\mbox{yes}&\mbox{yes}&\mbox{yes}&\mbox{yes}\\
\mbox{\tt Thales:}&\mbox{no}&\mbox{yes}&\mbox{no}&\mbox{no}\\
\mbox{\tt dhall:}&\mbox{yes}&\mbox{yes}&\mbox{yes}&\mbox{no}\\
\mbox{\tt Hall:}&\mbox{yes}&\mbox{yes}&\mbox{yes}&\mbox{no}\\
\mbox{\tt hall:}&\mbox{no}&\mbox{yes}&\mbox{no}&\mbox{no}\\
\mbox{\tt Hughes:}&\mbox{yes}&\mbox{yes}&\mbox{no}&\mbox{no}\\
\mbox{\tt hughes:}&\mbox{yes}&\mbox{yes}&\mbox{no}&\mbox{no}\\
\hline
\end{array}
$$
\end{remark}

\section{Pentagonal liners}

\begin{definition} A liner $X$ is defined to be \index{pentagonal liner}\index{liner!pentagonal}\defterm{pentagonal} if every pentagon $p_0p_1p_2p_3p_4\in X^5$ with $ \bigcup_{i\in 4}\Aline {p_i}{p_{i\oplus1}}\cap\Aline{p_{i\oplus2}}{a_{i\oplus 4}}=\varnothing$ is regular, i.e. has $ \bigcup_{i\in 5}\Aline {p_i}{p_{i\oplus1}}\cap\Aline{p_{i\oplus2}}{a_{i\oplus 4}}=\varnothing$.

%
\end{definition}

Corollary~\ref{c:parallel-lines<=>} and the definition of the colored graph $\Pentagon$ in Example~\ref{ex:pentagon} implies 

\begin{proposition} A $3$-ranked liner $X$ is pentagonal if and only if it is a $\Pentagon$-liner.
\end{proposition}

\begin{theorem}\label{t:Moufang=>pentagonal} Every Moufang affine $3$-regular liner $X$ is pentagonal. 
\end{theorem}

\begin{proof} To prove that $X$ is pentagonal, fix any pentagon $p_0p_1p_2p_3p_4\in X^5$ with $\varnothing=\bigcup_{i\in 4}\Aline{p_i}{p_{i\oplus1}}\cap\Aline{p_{i\oplus2}}{p_{i\oplus4}}$. We have to prove that $\Aline {p_4}{p_0}\cap\Aline{p_1}{p_3}=\varnothing$. Consider the plane  $\Pi\defeq\overline{\{p_0,p_1,p_2,p_3,p_4\}}\subseteq X$. It follows from $|\Pi|\ge 5$ and $\|\Pi\|=3$ that $|X|_2=|\Pi|_2\ge 3$. Applying Theorem~\ref{t:Playfair<=>}, we conclude that the $3$-long affine $3$-regular liner $X$ is Playfair.

 It follows from $\Aline {p_0}{p_1}\cap\Aline {p_2}{p_4}=\varnothing$ that the triple $uow=p_0p_1p_2$ is an affine base in the plane $\Pi$. Let $e$ be the diunit of the affine base $uow$ and $\Delta=\Aline oe$ be the ternar of the based affine plane $(\Pi,uow)$. Since the Playfair liner $X$ is Moufang, the ternar $\Delta$ is linear, distributive, associative-plus, commutative-plus, and alternative-dot, by Theorem~\ref{t:inversive-dot} and \ref{t:affine-Moufang<=>}. For any elements $a,b\in\Delta$, let $L_{a,b}$ be the line in the based affine plane $(\Pi,uow)$, given by the equation $y=x_{\times}a_+b=x{\cdot}a+b$.

Since the plus loop $(\Delta,+)$ is commutative, for every $x\in \Delta$ there exists a unique element $\bar x\in\Delta$ such that $x+\bar x=o=\bar x+x$.  The right-distributivity of the biloop $\Delta$ ensures that for every $x\in\Delta$ we have
$$x+\bar x=o=(e+\bar e)\cdot x=e{\cdot}x+\bar e{\cdot}x=x+\bar e{\cdot}x$$ and hence $\bar x=\bar e{\cdot}x$, by the cancellativity of the plus loop $(\Delta,+)$. On the other hand, the left-distributivity of $R$ implies $x+\bar x=o=x{\cdot}(e+\bar e)=x{\cdot}e+x{\cdot}\bar e=x+x{\cdot}\bar e$ and hence $\bar x=x{\cdot}\bar e$.

It follows from $\Aline ow\cap\Aline{p_0}{p_3}=\Aline {p_1}{p_2}\cap\Aline {p_3}{p_0}=\varnothing=\Aline {p_0}{p_1}\cap\Aline {p_2}{p_4}=\Aline uo\cap\Aline {p_2}{p_4}$ that the points $p_3$ and $p_4$ have coordinates $(e,g)$ and $(h,e)$ for some points $g,h\in\Delta$. Consider the line $L_{\bar e,e}\subset \Pi$ given by the equation $y=x_\times \bar e_+e=x\cdot\bar e+e$. Since $o{\cdot}\bar e+e=e$ and $e{\cdot}\bar e+e=o$, the points $w=p_2$ and $u=p_0$ belong to the line $L_{\bar e,e}$. Next, consider the point $g^+\defeq g+e\in\Delta$ and consider the line $L_{\bar e,g^+}$ given by the equation $y=x_\times \bar e_+g^+=x{\cdot}\bar e+g^+$. Since $e{\cdot}\bar e+g^+=\bar e+e+g=g$, the point $p_3$ (with coordinates $(e,g)$) belongs to the line $L_{\bar e,g^+}$. Since $\Aline {p_3}{p_4}\parallel\Aline {p_0}{p_2}=L_{\bar e,e}\parallel L_{\bar e,g^+}$ and $p_3\in \Aline {p_3}{p_4}\cap L_{\bar e,g^+}$, the parallel lines $\Aline{p_3}{p_4}$ and $L_{\bar e,g^+}$ coincide and hence the coordinates $(h,e)$ of the point $p_4$ satisfy the equation $y=x\cdot\bar e+g^+$ of the line $L_{\bar e,g^+}=\Aline{p_3}{p_4}$. Then $e=h{\cdot}\bar e+g^+=\bar h+g+e$ and hence $\bar h+g=o$ and $h=g$. Therefore, the points $p_3$ and $p_4$ have coordinates $(e,g)$ and $(g,e)$, respectively.

\begin{picture}(100,135)(-150,-20)
{\linethickness{0.8pt}
\put(0,0){\color{teal}\line(1,0){50}}
\put(0,50){\color{teal}\line(1,0){81}}
\put(0,0){\color{cyan}\line(0,1){50}}
\put(50,0){\color{cyan}\line(0,1){81}}
\put(50,81){\color{violet}\line(1,-1){31}}
\put(50,0){\color{magenta}\line(31,50){31}}
\put(0,50){\color{blue}\line(50,31){50}}
\put(0,0){\color{magenta}\line(50,81){50}}
\put(0,0){\color{blue}\line(81,50){81}}
\put(0,50){\color{violet}\line(1,-1){50}}
}

\put(0,0){\color{teal}\line(1,0){50}}
\put(0,50){\color{teal}\line(1,0){81}}
\put(0,0){\color{cyan}\line(0,1){50}}
\put(50,0){\color{cyan}\line(0,1){81}}
\put(50,81){\color{violet}\line(1,-1){45}}
\put(96,28){\color{violet}$L_{\bar e,g^+}$}
\put(50,0){\color{magenta}\line(31,50){45}}
\put(96,75){\color{magenta}$L_{g,\bar g}$}
\put(0,50){\color{blue}\line(50,31){70}}
\put(73,95){\color{blue}$L_{g^{-1},e}$}
\put(0,0){\color{magenta}\line(50,81){65}}
\put(66,108){\color{magenta}$L_{g,o}$}
\put(0,0){\color{blue}\line(81,50){100}}
\put(103,60){\color{blue}$L_{g^{-1},o}$}
\put(0,50){\color{violet}\line(1,-1){60}}
\put(62,-15){\color{violet}$L_{\bar e,e}$}
\put(0,0){\color{red}\line(1,1){80}}
\put(83,82){\color{red}$\Delta$}

\put(50,50){\circle*{3}}
\put(52,43){$e$}
\put(0,0){\circle*{3}}
\put(-8,0){$o$}
\put(1,-7){$p_1$}
\put(50,0){\circle*{3}}
\put(44,-9){$u$}
\put(54,0){$p_0$}
\put(0,50){\circle*{3}}
\put(-11,45){$p_2$}
\put(-5,54){$w$}
\put(50,81){\circle*{3}}
\put(55,79){$p_3$}
\put(81,50){\circle*{3}}
\put(86,48){$p_4$}
\end{picture}

Since $(\Delta,\cdot)$ is a Moufang loop, there exists a unique element $g^{-1}\in\Delta$ such that $g^{-1}{\cdot}g=e=g{\cdot} g^{-1}$. Consider the parallel lines $L_{g^{-1},o}$ and $L_{g^{-1},e}$, and observe that $\Aline {p_1}{p_4}=L_{g^{-1},o}$ and $p_2=w\in L_{g^{-1},e}$. Since $\Aline{p_2}{p_3}\parallel \Aline {p_4}{p_1}=L_{g^{-1},o}\parallel L_{g^{-1},e}$, the parallel lines $\Aline {p_2}{p_3}$ and $L_{g^{-1},e}$ coincide and hence the coordinates $(e,g)$ of the point $p_3$ satisfy the equation $y=x{\cdot} g^{-1}+e$ of the line $L_{g^{-1},e}$. Then $g=e{\cdot}g^{-1}+e=g^{-1}+e$, $g{\cdot} g=(g^{-1}+e)\cdot g=g^{-1}{\cdot}g+e{\cdot}g=e+g$ and finally $g{\cdot}g+\bar g=e$.

Next, consider the parallel lines $L_{g,o}$ and $L_{g,\bar g}$. Since $o{\cdot}g+o=o$ and $e{\cdot}g+o=g$, the points $p_1$ and $p_3$ (with coordinates $(o,o)$ and $(e,g)$) belong to the line $L_{g,o}$. Since $e{\cdot}g+\bar g=o$ and $g{\cdot}g+\bar g=e$,  the points $p_0$ and $p_4$ (with coordinates $(e,o)$ and $(g,e)$) belong to the line $L_{g,\bar g}$. Then $\Aline{p_4}{p_0}\cap\Aline{p_1}{p_3}=L_{g,\bar g}\cap L_{g,o}=\varnothing$, witnessing that the liner $X$ is pentagonal.
\end{proof}

\begin{remark} Theorem~\ref{t:Moufang=>pentagonal} does not generalize to Moufang affine liners: Ivan Hetman has found a Steiner (and hence Pappian, Desarguesian and Moufang) affine liner of order 27, which is not pentagonal.
\end{remark}

\section{Propentagonal liners}

\begin{definition}\label{d:propentagonal} A liner $X$ is defined to be \index{propentagonal liner}\index{liner!propentagonal}\defterm{propentagonal} if for every pentagon $p_0p_1p_2p_3p_4\in X^5$ and line $L\subseteq X$ with  $L\cap \Aline {p_{i}}{p_{i\oplus1}}\cap\Aline{p_{i\oplus2}}{p_{i\oplus 4}}\ne\varnothing$ for all  $i\in 4$, we have $\Aline {p_{4}}{p_{0}}\cap\Aline{p_{1}}{p_{3}}\subseteq L$.
\end{definition}

Definition~\ref{d:propentagonal} implies the following simple but useful fact.

\begin{proposition} Every subliner of a propentagonal liner is propentagonal.
\end{proposition}

Let us recall that a projective $X$ is {\em everywhere $\mathcal P$}, where $\mathcal P$ is a property of affine liners, if for every hyperplane $H$ in $X$, the affine liner $X\setminus H$ has property $\mathcal P$. 

\begin{theorem}\label{t:pro5<=>everywhere5} A projective liner is propentagonal if and only if it is everywhere penta\-gonal.
\end{theorem}

\begin{proof} Assume that the projective liner $Y$ is propentagonal.  To prove that $Y$ is everywhere pentagonal, take any hyperplane $H\subseteq Y$. We have to prove that the subliner $X\defeq Y\setminus H$ is pentagonal. By Proposition~\ref{p:projective-minus-hyperplane}, the subliner $X=Y\setminus H$ is affine and regular. Take any pentagon $p_0p_1p_2p_3p_4$ in $X$ and consider the plane $\Pi\defeq\overline{\{p_0,p_1,p_2,p_3,p_4\}}$ in the affine liner $X$. Then its flat hull $\overline\Pi$ in $Y$ is a plane in the projective liner $Y$. By the projectivity of $\bar\Pi$, for every $i\in 5$, there exists a point $h_i\in \Aline {p_{i}}{p_{i\oplus 1}}\cap\Aline {p_{i\oplus 2}}{p_{i\oplus 4}}\subseteq \bar\Pi$. Assuming that $\{h_0,h_1,h_2,h_3\}\subseteq H$, we should check that $h_4\in H$. Since $\{h_0,h_1,h_2,h_3\}\subseteq H\cap\bar\Pi$, the intersection $H\cap\bar \Pi$ is a line in the projective plane $\bar\Pi$. Since $Y$ is propentagonal, $\{h_0,h_1,h_2,h_3\}\subseteq H\cap\bar\Pi$ implies $h_4\in H\cap\bar\Pi\subseteq H$.
\smallskip

Now assume that the projective liner $Y$ is everywhere pentagonal. To prove that $Y$ is propentagonal, take any pentagon $p_0p_1p_2p_3p_4\in Y^5$ and any line $L\subset Y$ such that $\varnothing\ne L\cap \Aline {p_i}{p_{i\oplus1}}\cap\Aline{p_{i\oplus2}}{p_{i\oplus4}}$ for all $i\in 4$. We have to prove that $\Aline {p_4}{p_0}\cap\Aline {p_1}{p_3}\subseteq L$. Consider the plane $\Pi\defeq\overline{\{p_0,p_1,p_2,p_3,p_4\}}$ in the projective liner $Y$. By the projectivity of the plane $\Pi$, for every $i\in 5$, there exists a unique point 
$h_i\in  \Aline {p_i}{p_{i\oplus1}}\cap\Aline{p_{i\oplus2}}{p_{i\oplus4}}$. The definition of a pentagon ensures that $\{h_0,h_1,h_2,h_3,h_4\}\cap\{p_0,p_1,p_2,p_3,p_4\}=\varnothing$. By our assumption, $h_i\in L$ for all $i\in 4$. We claim that $\{p_0,p_1,p_2,p_3,p_4\}\subseteq\Pi\setminus L$. Indeed, assuming that $p_i\in L$ for some $i\in \{0,1,2\}$, we conclude that $p_{i\oplus1}\in \Aline {p_i}{h_i}=L$ and then $p_{i\oplus2}\in \Aline {p_{i\oplus1}}{h_{i\oplus1}}=L$, which contradicts the definition of a pentagon.  Assuming that $p_i\in L$ for some $i\in \{3,4\}$, we conclude that $p_{i\oplus2}\in \Aline {p_{i}}{h_{i\oplus3}}=L$ and $p_{i\oplus3}\in \Aline {p_{i\oplus2}}{h_{i\oplus2}}=L$, which contradicts the definition of a pentagon. Therefore, $\{p_0p_1p_2p_3p_4\}\subseteq \Pi\setminus L$. Let $F$ be a maximal $3$-long flat in $Y$ that contains the point $p_0$. By Lemma~\ref{l:ox=2}, $\Aline xy=\{x,y\}$ for every $x\in F$ and $y\in Y\setminus F$. Taking into account that that for every $i\in 5$, the line $\Aline {p_i}{p_{i\oplus 1}}$  is $3$-long, we conclude that $p_i\in F$ for all $i\in 5$ and hence $\Pi\subseteq Y$ and hence the plane $\Pi$ is $3$-long. By Corollary~\ref{c:Avogadro-projective}, the $3$-long projective plane $\Pi$ is $2$-balanced. Assuming that $|\Pi|_2=3$, we conclude that $|\Pi\setminus L|=4$, which contradicts $|\Pi\setminus L|\ge|\{p_0,p_1,p_2,p_3,p_4\}|=5$. This contradiction shows that $|\Pi|_2\ge 4$. By Corollary~\ref{c:procompletion-rank}, the subliner $\Pi\setminus L$ has rank $\|\Pi\setminus L\|=\|\Pi\|$ and hence $\Pi\setminus L$ is an affine plane.

By the Kuratowski--Zorn Lemma, the independent set $\{p_0,h_0,h_1\}$ is contained in the maximal independent set $M\subseteq Y$. The Exchange Property of the projective liner $Y$ ensures that the flat $H\defeq\overline{M\setminus\{p_0\}}$ is a hyperplane in $Y$. The rankedness of the projective liner $Y$ ensures that $H\cap\Pi=L$. Since $\Pi\setminus H=\Pi\setminus L$ is a plane, the $5$-tuple $p_0p_1p_2p_3p_4$ is a pentagon in the affine liner $X\defeq Y\setminus H$  such that $\bigcup_{i\in 4}X\cap\Aline {p_i}{p_{i\oplus1}}\cap\Aline{p_{i\oplus2}}{p_{i\oplus4}}=\bigcup_{i\in 4}X\cap\{h_i\}=\varnothing$. Since the projective liner $Y$ is everywhere pentagonal, the affine subliner $X=Y\setminus H$ is pentagonal and hence $\varnothing=X\cap\Aline{p_4}{p_0}\cap\Aline{p_1}{p_3}=X\cap\{h_4\}$ and hence $\Aline {p_4}{p_0}\cap\Aline{p_1}{p_3}=\{h_4\}\subseteq \Pi\cap H=L$. 
\end{proof}

\begin{exercise} Prove that a projective liner is propentagonal if and only if it is everywhere propentagonal.
\end{exercise}

\begin{theorem}\label{t:Moufang=>propentagonal} Every Moufang proaffine regular liner is propentagonal.
\end{theorem}

\begin{proof} Let $X$ be a Moufang proaffine regular liner. If $\|X\|\le 2$, then $X$ contains no pentagons and hence $X$ is propentagonal, vacuously. So, assume that $\|X\|\ge 3$. By Theorems~\ref{t:Moufang=>compreg} and \ref{t:compreg-Moufang<=>}, $X$ is completely regular and its spread completion $\overline X$ is a Moufang projective liner. By Propositions~\ref{p:Moufang-minus-flat} and \ref{p:projective-minus-hyperplane}, for every hyperplane $H$ in $\overline X$, the subliner $\overline X\setminus H$ is Moufang, affine, and regular. By Theorem~\ref{t:Moufang=>pentagonal}, the Moufang affine regular liner $\overline X\setminus H$ is pentaginal, witnessing that the projective liner $\overline X$ is everywhere pentagonal. By Theorem~\ref{t:pro5<=>everywhere5}, the everywhere pentagonal projective liner $\overline X$ is propentagonal and so is its subliner $X$.
\end{proof}

\begin{theorem}[Ivan Hetman\footnote{The proof was found by a computer program written by Ivan Hetman}, 2025]\label{t:Fano=>propentagonal} Every Fano proaffine regular liner is propentagonal.
\end{theorem}

\begin{proof} Let $X$ be a Fano proaffine regular liner. By Theorems~\ref{t:wreg-proFano=>compreg} and \ref{t:compreg-Fano=>comp-Fano}, the Fano proaffine regular liner $X$ is completely regular and its spread completion $\overline X$ is a Fano projective liner. Since subliners of propentagonal liners are propentagonal, it suffices to check that the Fano projective liner $\overline X$ is propentagonal. Fix any pentagon $p_0p_1p_2p_3p_4$ in $\overline X$ and any line $L\subseteq \overline X$ such that for every $i\in 4$, the unique point $h_i\in \Aline{p_i}{p_{i\oplus1}}\cap\Aline{p_{i\oplus2}}{p_{i\oplus4}}$ belongs to the line $L$. We have to prove that $h_4\in L$. 

Consider the plane $\Pi\defeq\overline{\{p_0,p_1,p_2,p_3,p_4\}}$. By the projectivity of $\Pi$, there exist unique points $a\in \Aline {p_0}{p_1}\cap\Aline{p_2}{p_3}$ and $b\in \Aline {p_0}{p_2}\cap\Aline{p_1}{p_3}$. Consider the quadrangle $p_0p_1p_2p_3$. Since the projective liner $\overline X$ is Fano, the set 
$(\Aline {p_0}{p_1}\cap\Aline{p_2}{p_3})\cup(\Aline {p_0}{p_2}\cap\Aline{p_1}{p_3})\cup(\Aline {p_0}{p_3}\cap\Aline{p_1}{p_2})=\{a,b,h_1\}$ has rank $2$.  Next, consider the quadrangle $p_1p_2h_0h_2$.  Since the projective liner $\overline X$ is Fano, the set 
$$(\Aline {p_1}{p_2}\cap\Aline{h_0}{h_2})\cup(\Aline {p_1}{h_0}\cap\Aline{p_2}{h_2})\cup(\Aline {p_1}{h_2}\cap\Aline{p_2}{h_0})=\{h_1,a,p_4\}$$ has rank $2$. Then $b,p_4\in \Aline{a}{h_1}$ and hence $h_1\in \Aline b{p_4}$. Finally, consider the quadrangle $bp_0p_3p_4$. Since the projective liner $\overline X$ is Fano, the set 
$$(\Aline b{p_0}\cap\Aline{p_3}{p_4})\cup(\Aline b{p_3}\cap\Aline {p_0}{p_4})\cup(\Aline b{p_4}\cap\Aline{p_0}{p_3})=\{h_3,h_4\}\cup(\Aline b{p_4}\cap\Aline{p_0}{p_3})$$ has rank $2$ and hence 
$h_1\in\Aline b{p_4}\cap \Aline {p_3}{p_0}\subseteq \Aline {h_3}{h_4}$, which implies the desired inclusion $h_4\in \Aline{h_1}{h_3}=L$.
\end{proof}

Theorems~\ref{t:Moufang=>propentagonal}, \ref{t:Fano=>propentagonal} and \ref{t:duo-Desarg<=>} imply the following corollary.

\begin{corollary}\label{c:bi-Desarguesian=>propentagonal} Every bi-Desarguesian projective space is propentagonal.
\end{corollary}

\begin{proposition}\label{p:free-pro5} The free projectivization of a liner $X$ is propentagonal if and only if  for every pentagon $p_0p_1p_2p_3p_4\in X^5$  and every line $L$ in $X$ with $\varnothing \ne L\cap\Aline {p_i}{p_{i\oplus1}}\cap\Aline{p_{i\oplus 2}}{p_{i\oplus 4}}$ for all $i\in 4$, the set $L\cap \Aline {p_4}{p_0}\cap\Aline{p_1}{p_3}$ is not empty.
\end{proposition} 

\begin{proof} First assume that the free projectivization $\widehat X$ of a liner $X$ is propentagonal.  Fix any pentagon $p_0p_1p_2p_3p_4\in X^5$ and any line $L\subseteq X$ such that for every $i\in 4$, the set $L\cap\Aline {p_i}{p_{i\oplus1}}\cap\Aline{p_{i\oplus 2}}{p_{i\oplus 4}}$ contains some point $h_i$. By projectivity of $\widehat X$, there exists a point $h_4\in \Aline {p_4}{p_0}\cap\Aline {p_1}{p_3}$. Since the projective liner $\widehat X$ is propentagonal, $h_4\in\overline L$. Observe that the subliner $S\defeq\{p_0,p_1,p_2,p_3,p_4,h_0,h_1,h_2,h_3,h_4\}$ of $\widehat X$ is $3$-wide. By Theorem~\ref{t:free-projectivization}(5), $S\subseteq X$ and hence $h_4\in \overline L\cap X=L$.
\smallskip

Now assume that the liner $X$ satisfies the condition formulated in the proposition. To prove that the free projectivization $\widehat X$ of $X$ is propentagonal, take any   pentagon $p_0p_1p_2p_3p_4\in \widehat X^5$  and any line $L\subseteq \widehat X$ such that for every $i\in 4$, the set $L\cap\Aline {p_i}{p_{i\oplus1}}\cap\Aline{p_{i\oplus 2}}{p_{i\oplus 3}}$ contains some point $h_i$. By the projectivity of $\widehat X$, there exists a point $h_4\in \Aline {p_4}{p_0}\cap\Aline{p_1}{p_3}$. We have to prove that $h_4\in L$. Observe that the subliner $S\defeq\{p_0,p_1,p_2,p_3,p_4,h_0,h_1,h_2,h_3\}$ is $3$-wide. By Theorem~\ref{t:free-projectivization}(5), $S\subseteq X$. By our assumption, there exists a point $h_4'\in X\cap L\cap \Aline{p_4}{p_0}\cap\Aline{p_1}{p_3}\subseteq L\cap\{h_4\}$, which implies $h_4=h_4'\in L$.
\end{proof}
 
\begin{example} There exists a $4$-element projective liner whose free projectivization $\widehat X$ is a propentagonal $\w$-long projective plane, which is  not tri-Desarguesian.
\end{example}

\begin{proof} Choose any $4$-element set $X$ endowed with the family of $2$-point lines $\mathcal L\defeq\{L\subseteq X:|L|=2\}$. By Theorem~\ref{t:free-projectivization}(3), the free projectivization $\widehat X$ of $X$ is an $\w$-long projective plane. By Proposition~\ref{p:projectivizationP1010not}, the free projectivization $\widehat X$ of is not $\mathsf{(P^{10}_{10})}$, and by Theorem~\ref{t:P1010<=>}, the projective plane $\widehat X$ is not tri-Desarguesian. By Proposition~\ref{p:free-pro5}, the projective plane $\widehat X$ is propentagonal.
\end{proof}  
 
 \begin{problem} Let $X$ be a propentagonal (pro)affine space. Is the spread completion of $X$ propentagonal?
 \end{problem}

\begin{remark} By computer calculations Ivan Hetman has established that a projective plane of order $9$ is propentagonal if and only if it is Desargesian.
\end{remark}

This allows to make the following 

\begin{conjecture}\label{conj:Des<=>propentaginal} A finite projective plane is Desarguesian if and only if it is propentagonal.
\end{conjecture}
 
\section{Pentagramic liners}

\begin{definition}\label{d:pentagramic} A liner $X$ is called \index{pentagramic liner}\index{liner!pentagramic}\defterm{pentagramic} if for every pentagon $p_0p_1p_2p_3p_4\in X^5$ and point $o\in \bigcap_{i\in 4}\Aline {p_i}{p'_i}$ where $p_i'\in \Aline {p_{i\oplus 1}}{p_{i\oplus2}}\cap\Aline {p_{i\oplus3}}{p_{i\oplus4}}$, we have $o\in \Aline{p_4}{p'_4}$.
\end{definition}

\begin{picture}(100,115)(-200,-55)
\put(0,50){\line(29.4,-90.5){29.4}}
\put(0,50){\line(-29.4,-90.5){29.4}}
\put(-47.6,15.4){\line(1,0){95.2}}
\put(-47.6,15.4){\line(77,-55.9){77}}
\put(47.6,15.4){\line(-77,-55.9){77}}
\put(0,50){\color{red}\line(0,-1){69.1}}
\put(47.6,15.4){\line(-47.6,-15.4){65.8}}
\put(-47.6,15.4){\line(47.6,-15.4){65.8}}
\put(-29.4,-40.5){\line(29.4,40.5){40.6}}
\put(29.4,-40.5){\line(-29.4,40.5){40.6}}

\put(0,0){\color{red}\circle*{3}}
\put(0,-6){\color{white}\circle*{7}}
\put(-2.5,-8){$o$}
\put(0,50){\circle*{3}}
\put(-3,54){$p'_4$}
\put(47.6,15.4){\circle*{3}}
\put(51,12){$p'_0$}
\put(-47.6,15.4){\circle*{3}}
\put(-59,12){$p'_3$}
\put(-29.4,-40.5){\circle*{3}}
\put(-40,-46){$p'_2$}
\put(29.4,-40.5){\circle*{3}}
\put(32,-46){$p'_1$}
\put(0,-19.1){\circle*{3}}
\put(-3,-28){$p_4$}
\put(18.2,-5.9){\circle*{3}}
\put(22,-9){$p_3$}
\put(-18.2,-5.9){\circle*{3}}
\put(-30,-9){$p_0$}
\put(11.2,15.5){\circle*{3}}
\put(12,20){$p_2$}
\put(-11.2,15.5){\circle*{3}}
\put(-20,20){$p_1$}

\end{picture}

Definition~\ref{d:pentagramic} implies the following proposition.

\begin{proposition} Every subliner of a petagramic liner is pentagramic.
\end{proposition}

\begin{theorem}\label{t:pro5-dual-pentagramic} A projective plane is propentagonal if and only if the dual projective plane is pentagramic.
\end{theorem}
 
\begin{proof} Assume that a projective plane $(\Pi,\mathcal L)$ is propentagonal. To prove that the dual projective plane $(\mathcal L,\Pi)$ is pentagramic, fix any pentagon $L_0L_1L_2L_3L_4\in\mathcal L^5$ and line $L\in\mathcal L$ such that $L_i'\cap L_i\subseteq L$ for all $i\in 4$, where $L_i'=\overline{(L_{i\oplus1}\cap L_{i\oplus 2})\cup(L_{i\oplus3}\cap L_{i\oplus 4})}$ for $i\in 5$. We have to prove that $L_4\cap L_4'\subseteq L$.

By the projectivity of the plane $\Pi$, for every $i\in 5$, there exists a unique point $p_i\in L_{i\oplus4}\cap L_{i}$. Then $L_i=\Aline{p_i}{p_{i\oplus1}}$ for all $i\in 5$. Since no three lines of $L_0,L_1,L_2,L_3,L_4$ are concurrent, the points $p_0,p_1,p_2,p_3,p_4$ are distinct. Assuming that some three distinct points $x,y,z\in\{p_0,p_1,p_2,p_3,p_4\}$ are collinear, we can find numbers $i\in 5$ and $k\in\{2,3\}$ such that $\{p_i,p_{i\oplus 1},p_{i\oplus k}\}=\{x,y,z\}$ and hence the points $p_i,p_{i\oplus 1},p_{i\oplus k}$ are collinear. Since the lines $L_{i\oplus4},L_i,L_{i\oplus1}$ are distinct, the number $k$ is not equal to $2$ and hence $k=3$. Then $p_{i\oplus 3}\in L_{i\oplus 2}\cap L_{i\oplus 3}\cap\Aline {p_i}{p_{i\oplus 1}}=L_{i\oplus2}\cap L_{i\oplus 3}\cap L_{i}$, which contradicts the choice of the lines $L_0,L_1,L_2,L_3,L_4$ (forming a pentagon in the dual projective plane). Therefore, $p_0p_1p_2p_3p_4$ is a pentagon in the projective plane $(\Pi,\mathcal L)$. Observe that $$\Aline {p_i}{p_{i\oplus 1}}\cap\Aline{p_{i\oplus2}}{p_{i\oplus4}}=L_{i}\cap \overline{(L_{i\oplus1}\cap L_{i\oplus2})\cup(L_{i\oplus 3}\cap L_{i\oplus 4})}=L_i\cap L_i'\subseteq L\quad\mbox{for all $i\in 4$}.$$ Since the projective plane $(\Pi,\mathcal L)$ is propentagonal, $L_4\cap L_4'=\Aline{p_4}{p_0}\cap\Aline {p_1}{p_3}\subseteq L$.
\smallskip

Now assume that the dual projective plane $(\mathcal L,\Pi)$ is pentagramic. To prove that the projective plane $(\Pi,\mathcal L)$ is propentagonal, fix any  pentagon $p_0p_1p_2p_3p_4\in\Pi^5$ and any line $L$ in the plane $\Pi$ such that for every $i\in 4$ the unique point $h_i\defeq \Aline {p_i}{p_{i\oplus1}}\cap\Aline{p_{i\oplus2}}{p_{i\oplus4}}$ belongs to the line $L$. We have to prove that $h_4\in L$. For every $i\in 5$, consider the lines $L_i\defeq\Aline {p_i}{p_{i\oplus1}}$ and $L_i'\defeq \Aline{p_{i\oplus2}}{p_{i\oplus4}}=\overline{(L_{i\oplus1}\cap L_{i\oplus2})\cup(L_{i\oplus 3}\cap L_{i\oplus 4})}$. Observe that for every $i\in 4$, we have $L_i\cap L_i'=\Aline {p_i}{p_{i\oplus1}}\cap\Aline {p_{i\oplus2}}{p_{i\oplus4}}=\{h_i\}\subseteq L$. Since the dual projective plane $(\mathcal L,\Pi)$ is pentagramic, $\{h_4\}=L_4\cap L_4'\subseteq L$, witnessing that the projective plane $(\Pi,\mathcal L)$ is propentagonal.
\end{proof} 
 
\begin{theorem}\label{t:FM=>5gramic} Every Fano or Moufang proaffine regular liner is pentagramic. 
\end{theorem}

\begin{proof} Let $X$ be a Fano or Moufang proaffine regular liner. By Theorems~\ref{t:wreg-proFano=>compreg} and \ref{t:Moufang=>compreg}, $X$ is completely regular and the spread completion $\overline X$ is a Fano or Moufang projective liner. Since subliners of pentagramic liners are pentagramic, it suffices to check that the projective liner $\overline X$ is pentagramic. Fix any pentagon $p_0p_1p_2p_3p_4$ in $\overline X$ and a point $o\in\overline X$ such that $o\in \bigcap_{i\in 4}\Aline {p_i}{p_i'}$, where $p'_i\in \Aline {p_{i\oplus1}}{p_{i\oplus2}}\cap\Aline {p_{i\oplus3}}{p_{i\oplus4}}$ for all $i\in 5$. We have to prove that $o\in \Aline {p_4}{p_4'}$. Consider the plane $\Pi\defeq\overline{\{p_0,p_1,p_2,p_3,p_4\}}$ in the projective liner $\overline X$. Let $M$ be a maximal $3$-long flat in $\overline X$ that contain the point $p_0$. By Lemma~\ref{l:ox=2}, for every $x\in M$ and $y\in  \overline X\setminus M$, we have $\Aline xy=\{x,y\}$. For every $i\in 5$, we have $\{p_i,p_{i\oplus1},p_{i\oplus4}'\}\subseteq\Aline {p_i}{p_{i\oplus1}}$, which implies $\{p_0,p_1,p_2,p_3,p_4\}\subseteq M$ and hence $\Pi\subseteq M$. Therefore, the projective plane $\Pi$ is $3$-long. Since $X$ is Fano or Moufang, so is the plane $\Pi$. By Theorem~\ref{t:duo-Desarg<=>}, the Fano or Moufang projective plane $\Pi$ is bi-Desarguesian. By Corollary~\ref{c:123-Desarg-self-dual}, the dual projective plane $\Pi^*$ to $\Pi$ is bi-Desarguesian. By Corollary~\ref{c:bi-Desarguesian=>propentagonal}, the bi-Desarguesian projective plane $\Pi^*$ is propentagonal, and by Theorem~\ref{t:pro5-dual-pentagramic}, the projective plane $\Pi$ is pentagramic, which implies $o\in \Aline {p_4}{p_4'}$ and witnesses that the projective liner $\overline X$ is pentagramic and so is its subliner $X$.
\end{proof}   

\begin{proposition}\label{p:free-5gramic} The free projectivization of a liner $X$ is pentgramic if and only if for every pentagon $p_0p_1p_2p_3p_4\in X^5$ and point $o\in \bigcup_{i\in 4}\Aline{p_i}{p_i'}$ where $p_i'\in X\cap\Aline {p_{i\oplus1}}{p_{i\oplus 2}}\cap\Aline{p_{i\oplus3}}{p_{i\oplus4}}$ for $i\in 4$, there exists a point 
$p_4'\in \Aline {p_0}{p_1}\cap\Aline {p_2}{p_3}\cap X$ and $o\in \Aline {p_4}{p_4'}$.
\end{proposition}

\begin{proof} Assume that a liner $X$ satisfies the condition of the proposition. To prove that the projective plane $\widehat X$ is pentagramic, take any pentagon $p_0p_1p_2p_3p_4$ and point $o$ in $\widehat X$ such that $o\in\bigcap_{i\in 4}\Aline{p_i}{p'_i}$ where $p_i'\in\Aline {p_{i\oplus1}}{p_{i\oplus2}}\cap\Aline {p_{i\oplus 3}}{p_{i\oplus 4}}$for $i\in 4$. We have to prove that $o\in \Aline {p_4}{p_4'}$. Observe that the subliner $S\defeq\{p_0,p_1,p_2,p_3,p_4,o,p_0',p_1',p_2',p_3'\}$ of $\widehat X$ is $3$-wide. By Theorem~\ref{t:free-projectivization}(5), $S\subseteq X$. Our assumption guarantees that $p_4'\in X$ and $o\in \Aline {p_4}{p_4'}$, witnessing that the projective plane $\widehat X$ is pentagramic.
\smallskip 

Now assume that the free projectivization $\widehat X$ of the liner $X$ is pentagramic. Fix any pentagon $p_0p_1p_2p_3p_4\in X^5$ and point $o\in \bigcup_{i\in 4}\Aline{p_i}{p_i'}$ where $p_i'\in X\cap\Aline  {p_{i\oplus1}}{p_{i\oplus 2}}\cap\Aline{p_{i\oplus3}}{p_{i\oplus4}}$for $i\in 4$. By the projectivity of $\widehat X$, there exists a point $p_4'\in \Aline {p_0}{p_1}\cap\Aline {p_2}{p_3}\subseteq\widehat X$. Since $\widehat X$ is pentagramic, $o\in \Aline {p_4}{p_4'}$. Observe that the subliner $S\defeq\{p_0,p_1,p_2,p_3,p_4,o,p_1',p_2',p_3',p_4'\}$ of $\widehat X$ is $3$-wide. By Theorem~\ref{t:free-projectivization}(5), $S\subseteq X$. Then $p_4'\in \Aline {p_0}{p_1}\cap\Aline {p_3}{p_3}\cap X$ and $o\in \Aline {p_4}{p_4'}$.
\end{proof}

\begin{example} There exists a pentargamic $\w$-long projective plane $\Pi$, which is neither propentagonal nor tri-Desargusian. The dual projective plane $\Pi^*$ is propentagonal but not pentagramic and not tri-Desarguesian.
\end{example}

\begin{proof} Take any $9$-element set $X=\{p_0,p_1,p_2,p_3,p_4,h_0,h_1,h_2,h_3\}$ endowed with the family of lines $$\mathcal L\defeq\big\{\{p_0,p_4\},\{p_1,p_3\},\{h_0,h_1,h_2,h_3\}\big\}\cup\big\{\{p_i,p_{i\oplus 1},h_i\},\{p_{i\oplus 2},p_{i\oplus 4},h_i\}:i\in 4\big\}.$$By Theorem~\ref{t:free-projectivization}(3), the free projectivization $\widehat X$ of the liner $(X,\mathcal L)$ is an $\w$-long projective plane. By Propositions~\ref{p:free-5gramic} and \ref{p:free-pro5}, the projective plane $\Pi\defeq\widehat X$ is pentagramic but not propentagonal. By Proposition~\ref{p:projectivizationP1010not} and Theorem~\ref{t:P1010<=>}, the projective plane $\Pi=\widehat X$ is not tri-Desarguesian. By Theorem~\ref{t:pro5-dual-pentagramic}, the dual projective plane $\Pi^*$ is propentagonal and not pentagramic. By Proposition~\ref{p:projectivizationP1010not} and Theorem~\ref{t:P1010<=>}, $\Pi^*$ is not tri-Desarguesian. 
\end{proof}

The following conjecture if equivalent to Conjecture~\ref{conj:Des<=>propentaginal}.

\begin{conjecture} A finite projective plane is Desarguesian if and only if it is pentagramic.
\end{conjecture}

\chapter{Area in Linear Geometry}

\rightline{\em Triangles on equal bases and in the same parallels equal one another.}

\rightline{Proposition I.38\footnote{\tt http://aleph0.clarku.edu/~djoyce/java/elements/bookI/propI38.html} in ``Elements'' of 
\index[person]{Euclid}Euclid\footnote{{\bf Euclid of Alexandria} (fl. c. 300 BCE) was a Greek mathematician, often called the``Father of Geometry''. Little is known about his life, but he is believed to have taught at the Museum and Library of Alexandria during the reign of Ptolemy I Soter (c. 323--283 BCE).
The enduring fame of Euclid rests on his monumental work ``Elements'' ($\Sigma\tau o\iota \chi \varepsilon \tilde\i \alpha$ in Greek), a thirteen-book treatise that systematically organized the geometry and number theory known in his time. The Elements begins with a small set of definitions, postulates, and common notions, from which Euclid rigorously deduces a vast body of propositions about figures, numbers, and proportions. The logical structure and deductive method of the Elements shaped mathematics for over two millennia, serving as both a textbook and a model of axiomatic reasoning. Besides the Elements, several other works are attributed to Euclid, including ``Data'', ``Optics'', ``Phaenomena'', and ``Catoptrica''. 
Euclid’s influence extends beyond geometry: his approach laid the groundwork for the modern concept of mathematical proof and for the axiomatic method that underlies all of modern mathematics.}}
\bigskip

In this chapter we study the fundamental notion of area in Linear Geometry. By analogy with vectors (which are equivalence classes of pairs), scalars (= equivalence classes of line triples), proscalars (= equivalence classes of line quadruples), areas and absolute areas are equivalence classes of  triples of points in a liner. The equivalence between such triples is defined with the help of the relation between triangles to have the same base and height. In Moufang Playfair planes two triangles have the same area if and only if each of them can be transformed to the other by a finite sequence of hypershears. Areas in Moufang Playfair liners admit a natural multiplication by scalars. For Desarguesian planes the operation of multiplication allows us parametrize areas by coscalars, which are elements of the abelianization of the multiplicative semigroup of the scalar corps. For Pappian Playfair planes, the set of areas has the structure of one-dimensional vector space over the field of scalars. The additivity of area allows to define the area for any polygon in a Pappian Playfair plane. The area of such polygons can be calculated by the well-known shoelace formula. 

Our approach to defining area as an equivalence relation between triangles continues the classical tradition of synthetic geometry, tracing back to Euclid, who also understood area as an equivalence relation rather than a number.

\section{The area of a triangle in a liner}

In this section we introduce the notion of area of a triangle in an arbitrary liner. 

A \index{triangle}\defterm{triangle} in a liner $X$ is any ordered triple $abc\in X^3$ such that $\|\{a,b,c\}\|=3$. The points $a,b,c$ are called the \defterm{vertices} of the triangle $abc$, and the lines $\Aline ab$, $\Aline bc$, $\Aline ac$ are called the \defterm{sides} of the triangle $abc$. For a liner $X$, the set of all triangles in $X$ is denoted by $X^\vartriangle$. Therefore, $X^\vartriangle=\{abc\in X^3:\|\{a,b,c\}\|=3\}$. A liner $X$ is a \index{plane}\defterm{plane} if $\|X\|=3$.

\begin{definition}\label{d:boxtimes} Two triangles $a_1a_2a_3,b_1b_2b_3\in X^3$ are defined to have the \index{triangles!with the same base and height}\defterm{same base and height} if there exist distinct numbers $i,j,k\in\{1,2,3\}$ such that $a_ia_j=b_ib_j$ and the flat $\Aline{a_k}{b_k}$ is subparallel to the line $\Aline {a_i}{a_j}=\Aline {b_i}{b_j}$. 
\end{definition}

The following picture shows two triangles $a_ia_ja_k$ and $b_ib_jb_k$ with the same base and height.

\begin{picture}(100,75)(-160,-15)

\linethickness{0.7pt}
\put(-10,40){\color{red}\line(1,0){60}}
{\linethickness{1.2pt}
\put(0,0){\color{red}\line(1,0){40}}
}
\put(40,0){\line(1,4){10}}
\put(0,0){\line(-1,4){10}}
\put(0,0){\line(5,4){50}}
\put(-10,40){\line(5,-4){50}}

\put(0,0){\circle*{3}}
\put(-10,-2){$b_i$}
\put(0,-8){$a_i$}
\put(40,0){\circle*{3}}
\put(35,-8){$a_j$}
\put(43,-2){$b_j$}

\put(-10,40){\circle*{3}}
\put(-14,44){$a_k$}
\put(50,40){\circle*{3}}
\put(49,43){$b_k$}
\end{picture}

\begin{exercise} Show that two triangles $abc$ and $xyz$ in a projective liner have the same base and height if and only if $abc=xyz$.
\end{exercise}

\begin{proposition}\label{p:same-bh=>sameplanes} If two triangles $a_1a_2a_3$ and $b_1b_2b_3$ in  a liner have the same base and height, then $\overline{\{a_1,a_2,a_3\}}=\overline{\{b_1,b_2,b_3\}}$.
\end{proposition}

\begin{proof} Assuming that the triangles $a_1a_2a_3$ and $b_1b_2b_3$ have the same base and height, find distinct numbers $i,j,k\in\{1,2,3\}$ such that $a_ia_j=b_ib_j$ and the flat $\Aline{a_k}{b_k}$ is subparallel to the line $\Aline {a_i}{a_j}=\Aline {b_i}{b_j}$. The definition of the subparallelity $\Aline{a_k}{b_k}\subparallel\Aline {a_i}{a_j}$ implies $b_k\in \Aline {a_k}{b_k}\subseteq\overline{\{a_i,a_j,a_k\}}$. Then $\overline{\{b_i,b_j,b_k\}}=\overline{\{a_i,a_j,b_k\}}\subseteq \overline{\{a_i,a_j,a_k,b_k\}}=\overline{\{a_i,a_j,a_k\}}$. By analogy we can prove that $\overline{\{a_i,a_j,a_k\}}\subseteq\overline{\{b_i,b_j,b_k\}}$. Therefore,
$$\overline{\{a_1,a_2,a_3\}}=\overline{\{a_i,a_j,a_k\}}=\overline{\{b_i,b_j,b_k\}}=\overline{\{b_1,b_2,b_3\}}.$$
\end{proof}

Let $S_3$ be the permutation group of the set $3=\{0,1,2\}$.
Since triangles in a liner $X$ are just functions $3\to X$, any permutation $\sigma\in S_3$ induces a well-defined bijection $\circ\sigma:X^3\to X^3$, $\circ\sigma:\tau\mapsto \tau\sigma\defeq\tau\circ\sigma$.  Also, any automorphism $A:X\to X$ of the liner $X$ induces a well-defined bijection $A{\circ}:X^3\to X^3$, $A{\circ}:\tau\mapsto A\tau\defeq A\circ\tau$. Observe that for a triangle $\tau$, its set of vertices coincides with the range $\tau[3]\defeq\{\tau(0),\tau(1),\tau(2)\}$ of the function $\tau$.

\begin{exercise} Show that two triangles $\tau$ and $\tau'$ in a liner $X$ have the same base and height if and only if there exists a permutation $\sigma\in S_3$ such that $\tau\sigma{\restriction}_2=\tau'\sigma{\restriction}_2$ and the flat $\Aline{\tau\sigma(2)}{\tau'\sigma(2)}$ is subparallel to the line $\overline{\tau\sigma[2]}=\overline{\tau'\sigma[2]}$, where $2\defeq \{0,1\}$.
\end{exercise} 

\begin{definition} Let $X$ be a liner. For two triples of points $abc,xyz\in X^3$ we write $abc\boxtimes xyz$ if either $\max\{\|\{a,b,c\}\|,\|\{x,y,z\}\|\}\le 2$ or $abc$ and $xyz$ are two triangles with the same base and height.  For every $n\in \w$, let $\boxtimes^n$ be the $n$th iteration of the relation $\boxtimes$. It is defined by the recursive formula:
$$
\begin{aligned}
&\boxtimes^0\defeq\{(\tau,\tau')\in X^3\times X^3:\tau=\tau'\}\quad\mbox{and}\\
&\boxtimes^{n+1}\defeq \boxtimes^n\circ\boxtimes\defeq\{(\tau,\tau')\in X^3\times X^3:\exists t\in X^3\;\;(\tau\boxtimes^n t\boxtimes \tau')\}.
\end{aligned}
$$Also put $\boxtimes^{\omega}\defeq \bigcup_{n\in\w}\omega^n$.
Therefore, for two triangles $abc,xyz\in X^\vartriangle$, we write $xyz\boxtimes^\w   abc$ if there exist  $n\in\IN$ and triangles $\tau_0,\tau_1\dots,\tau_n$ in $X$ such that $xyz=\tau_0\boxtimes\tau_1\boxtimes\cdots\boxtimes \tau_n=abc$.
\end{definition}

\begin{proposition}\label{p:permutation-base-and-height} Let $\tau,\tau'$ be two triangles in a liner $X$, $\sigma$ be a permutation of the set $3$, and $A:X\to X$ be an automorphism of the liner $X$. If $\tau\boxtimes\tau'$, then $A\tau\sigma\boxtimes A\tau'\sigma$.
\end{proposition}

\begin{proof} For every triangle $t\in X^3$ and number $i\in 3$, denote the point $t(i)\in X$ by $t_i$. Then $t=t_1t_2t_3$. By definition of the relation $\tau\boxtimes\tau'$, there exist distinct numbers $i,j,k\in\{1,2,3\}$ such that $\tau_i\tau_j=\tau'_i\tau'_j$ and the flat $\Aline{\tau_k}{\tau'_k}$ is subparallel to the line $L\defeq \Aline{\tau_i}{\tau_j}=\Aline{\tau'_i}{\tau_j'}$. Since $A$ is an automorphism of the liner $X$, $\overline{A\tau[3]}=A[\overline{\tau[3]}]=A[\overline{\tau'[3]}]=\overline{A\tau'[3]}$ and the flat $\Aline{A\tau_k}{A\tau'_k}=A[\Aline{\tau_k}{\tau'_k}]$ is subparallel to the line $\Aline{A\tau_i}{A\tau_j}=\Aline{A\tau'_i}{A\tau_j'}=A[L]$.
Then the triple of numbers $i'j'k'=\sigma^{-1}ijk$ witnesses that $A\circ\tau\circ\sigma\boxtimes A\circ\tau'\circ\sigma$.
\end{proof}

Proposition~\ref{p:permutation-base-and-height} easily extends to the iterations of the relation $\boxtimes$.

\begin{proposition}\label{p:permutation-boxtimes} Let $\tau,\tau'$ be two triangles in a liner $X$, $\sigma$ be a permutation of the set $3$, and $A:X\to X$ be an automorphism of the liner $X$. If $\tau\boxtimes^n\tau'$ for some $n\in \w$, then $A\tau\sigma\boxtimes^n A\tau'\sigma$.
\end{proposition}

\begin{proof} It suffices to check that for every $n\in\w$, the following statement holds:
\begin{itemize}
\item[$(*_n)$] $\forall\tau,\tau'\in X^\vartriangle\;\;(\tau\boxtimes^n \tau'\;\Ra\;A\tau\sigma\boxtimes^n A\tau'\sigma)$.
\end{itemize}
For $n=0$ the statement $(*_0)$ is trivially true. Assume that for some $n\in\w$, the statatement $(*_n)$ holds. To prove that the statement $(*_{n+1})$, take any triangles $\tau,\tau'\in X^\vartriangle$ with $\tau\boxtimes^{n+1}\tau'$. By the recursive definition of the relation $\boxtimes^{n+1}$, there exists a triangle $t$ in $X$ such that $\tau\boxtimes^n t$ and $t\boxtimes \tau'$. The inductive hypothesis ensures that $A\tau\sigma\boxtimes^n A\tau\sigma$, and Proposition~\ref{p:permutation-base-and-height} guarantees that $At\sigma\boxtimes A\tau'\sigma$. Then $A\tau\sigma\boxtimes^n At\sigma\boxtimes A\tau'\sigma$ implies $A\tau\sigma\boxtimes^{n+1}A\tau'\sigma$. So, the statement $(*_{n+1})$ holds, and by the Principle of Mathematical Induction, it holds for all $n\in\w$.
\end{proof}

\begin{definition}\label{d:area} The \index{area}\index{area of a triangle}\index{triangle!area}\defterm{area} of a triangle $abc\in X^3$ in a liner $X$ is the family\index[note]{$[abc]$} $$[abc]\defeq\{xyz\in X^\vartriangle:xyz\boxtimes^{\w} abc\}$$ of all triangles $xyz$ in $X$ for which there exist $n\in\IN$ and a sequence of triangles $\tau_0,\tau_1,\dots,\tau_n\in X^\vartriangle$ such that $xyz=\tau_0\boxtimes\tau_1\boxtimes\tau_2\boxtimes\dots \boxtimes \tau_n=abc$. 
Two triangles $abc$ and $xyz$ have the \index{triangles!of the same area}\defterm{same area} if $[abc]=[xyz]$. In particular, two triangles have the same area if they have the same base and height. For a triple $abc\in X^3$ of collinear points, we define $[abc]$ to be the family $[0]\defeq\{xyz\in X^3:\|\{x,y,z\}\|\le 2\}=X^3\setminus X^\vartriangle$ of all collinear triples. The family $[0]=X^3\setminus X^\vartriangle$ is called the \index{zero area}\defterm{zero area}.
The set $$[X^3]\defeq\{[abc]:abc\in  X^3\}$$  of all areas in $X$ is called the \index[note]{$[X^3]$}\index{areol}\defterm{areol} of the liner $X$. The set
$$[X^\vartriangle]\defeq\{[abc]:abc\in X^\vartriangle\}=[X^3]\setminus\{[0]\}$$ is called the \defterm{nonzero areol} of $X$.
\end{definition}

\begin{exercise} Show that any triangle $abc$ in a projective plane $X$ has area $[abc]=\{abc\}$. Use this fact to show that $[X^3]=\{[0]\}\cup\big\{\{abc\}:abc\in X^\vartriangle\big\}$.
\end{exercise}

Proposition~\ref{p:permutation-boxtimes} implies the following useful fact.

\begin{proposition}\label{p:permutation-area} If two triangles $\tau,\tau'\in X^3$ in a liner $X$ have the same area, then so do the triangles $A\tau\sigma$ and $A\tau'\sigma$ for all automorphisms $A:X\to X$ and all permutations $\sigma\in S_3$.
\end{proposition}

Proposition~\ref{p:same-bh=>sameplanes} and Definition~\ref{d:area} imply the following corollary.

\begin{corollary}\label{c:same-area=>same-plane} If two triangles $abc$ and $xyz$ in a liner $X$ have the same area, then $\overline{\{a,b,c\}}=\overline{\{x,y,z\}}$.
\end{corollary}

We recall that a \index{parallelogram}\defterm{parallelogram} in a liner $X$ is any quadruple $abcd\in X^4$ such that $\Aline ab\parallel\Aline cd\ne \Aline ab$ and $\Aline ad\parallel\Aline bc\ne \Aline ad$.  A parallelogram $abcd$ is \index{Boolean parallelogram}\index{parallelogram!Boolean}\defterm{Boolean} if $\Aline ac\parallel\Aline bd$.

\begin{proposition}\label{p:12-permutations} If $abdc$ is a parallelogram in a liner  $X$, then 
$$[abc]=[abd]=[acd]=[bcd]=[bca]=[bda]=[cda]=[cdb]=[cab]=[dab]=[dac]=[dbc].$$
If the parallelogram $abcd$ is Boolean, then any triangles with vertices in the set $\{a,b,c,d\}$ have the same area.
\end{proposition}

\begin{proof} The equality of the areas follows from Definition~\ref{d:area} and  observation that  
$$abc\boxtimes abd \boxtimes acd \boxtimes bcd \boxtimes bca \boxtimes bda \boxtimes cda \boxtimes cdb \boxtimes cab \boxtimes dab \boxtimes dac \boxtimes dbc.$$This observation is illustrated in the sequence of pictures:



If $abcd$ is a Boolean parallelogram, then $\Aline ac\parallel\Aline bd$ implies $acb\boxtimes acd$ and then all 24 triangles with vertices in the set $\{a,b,c,d\}$ have the same area.
\end{proof}

\begin{theorem}\label{t:large=>area-trivial} Let $X$ be a balanced ranked plane $X$ such that $|X|+12{\cdot}|X|_2>4{\cdot}|X|_2^2+11$. Then all triangles in $X$ have the same area.
\end{theorem}

\begin{proof} Asuming that the plane $X$ has $|X|_2=2$, we conclude that $$27=|X|+12{\cdot}|X|_2>4|X|_2^2+11=27,$$ which is a contradiction showing that $|X|_2\ge 3$.

\begin{claim} The cardinal $|X|_2$ is finite.
\end{claim}

\begin{proof} To derive a contradiction, assume that the cardinal $|X|_2$ is infinite. Since $X$ is a plane, there exists a set $X_0\subseteq X$ of cardinality $|X_0|=\|X\|=3$ such that $X=\overline{X_0}$. Consider the sequence of sets $(X_n)_{n\in\w}$, defined by the recursive formula: $X_{n+1}\defeq\bigcup_{x,y\in X_n}\Aline xy$ for all $n\in\w$. It is clear that $X=\overline{X_0}=\bigcup_{n\in\w}X_n$. 

We claim that $|X_n|\le|X|_2$ for all $n\in\w$. Indeed, for $n=0$ we have $|X_0|=3<|X|_2$. Assume that for some $n\in\w$, the inequality $|X_n|\le|X|_2$ holds. Then $|X_{n+1}|\le|X_n|^2\cdot|X|_2 \le|X|_2^3=|X|_2$ (because the cardinal $|X|_2$ is infinite). Therefore, $|X_n|\le|X|_2$ for all $n\in\w$, by the Principle of Mathematical Induction.

 Since $|X|_2$ is infinite, the plane $X=\bigcup_{n\in\w}X_n$ has cardinality $|X|\le|X|_2\cdot\w=|X|_2\le|X|$. Therefore, $|X|=|X|_2$ and then $|X|+12{\cdot}|X|_2=|X|_2=4{\cdot}|X|_2^2+11$, which contradicts the strict inequality $|X|+12{\cdot}|X|_2>4{\cdot}|X|_2^2+11$ in the assumption. This contradiction shows that the cardinal $|X|_2$ is finite.
\end{proof}

Since $|X|_2$ is finite the strict inequality $|X|+12{\cdot}|X|_2>4{\cdot}|X|_2^2+11$ implies $$|X|>4{\cdot}|X|_2-12{\cdot}|X|_2+11=(3{\cdot}|X|_2^2-9{\cdot}|X|_2+9)+|X|_2{\cdot}(|X|_2-3)+2>3{\cdot}|X|_2^2-9{\cdot}|X|_2+9.$$In this case, Proposition~\ref{p:large=>Boolean-parallelogram} ensures that each triangle in $X$ can be completed to a Boolean parallelogram. Applying Proposition~\ref{p:12-permutations}, we obtain the following claim.

\begin{claim}\label{cl:=[triangles]} Any triangles $abc,xyz\in X^3$ with $\{a,b,c\}=\{x,y,z\}$ have the same area.
\end{claim}

\begin{claim}\label{cl:[vab]L[vbc]} For any line $L\subset X$ and distinct points $a,b,c\in L$ and $v\in X\setminus L$ we have $[vab]=[vbc]$.
\end{claim}

\begin{proof} Consider the set $D\defeq\bigcup_{x\in\Aline vb}(\Aline ax\cup\Aline cx)\cup\bigcup_{y\in \Aline va}\Aline yb\cup\bigcup_{z\in \Aline vc}\Aline zb$ and observe that 
$$
\begin{aligned}
|D|&\le|L\cup\Aline vb|
+\sum_{x\in \Aline vb\setminus\{b\}}|(\Aline ax\cup\Aline cx)\setminus\{a,v,c\}|+
\sum_{y\in \Aline va\setminus\{v,a\}}|\Aline yb\setminus\{b\}|\cup\sum_{z\in \Aline vc\setminus\{v,c\}}|\Aline zb\setminus\{v,b\}|\\
&\le (2|X|_2-1)+(|X|_2-1)(2|X|_2-4)+2(|X|_2-2)^2=4|X|_2^2-12|X|_2+11<|X|.
\end{aligned}
$$
Then we can choose a point $d\in X\setminus D$ and conclude that $avbd$ and $bvcd$ are two parallelograms. Proposition~\ref{p:12-permutations} and Claim~\ref{cl:=[triangles]} imply $[vab]=[vbd]=[vbc]$.
\end{proof}  

\begin{claim}\label{cl:[vab]=[vbc]} For any distinct points $b,v\in X$ and $a,c\in X\setminus\Aline vb$, we have $[vba]=[vbc]$.
\end{claim}

\begin{proof} If $\Aline ac\parallel\Aline vb$, then $vba\boxtimes vbc$ and hence $[vba]=[vbc]$. So, assume that the lines $\Aline ac$ and $\Aline vc$ in the plane $X$ are not parallel and hence they have a unique common point $o$. If $o\in \{b,v\}$, then $[vba]=[vbc]$, by Claims~\ref{cl:=[triangles]} and \ref{cl:[vab]L[vbc]}. So, we assume that $o\notin\{b,v\}$. By Claims~\ref{cl:=[triangles]} and \ref{cl:[vab]L[vbc]}, $[abv]=[abo]=[bao]=[boc]=[cob]=[cvb]=[vbc]$.

\begin{picture}(60,90)(-150,-10)
\put(0,30){\line(1,0){60}}
\put(0,30){\line(1,1){30}}
\put(0,30){\line(1,-1){30}}
\put(30,0){\line(0,1){60}}
\put(60,30){\line(-1,1){30}}
\put(60,30){\line(-1,-1){30}}

\put(0,30){\circle*{3}}
\put(-8,28){$a$}
\put(30,30){\circle*{3}}
\put(32,32){$o$}
\put(60,30){\circle*{3}}
\put(63,28){$c$}
\put(30,0){\circle*{3}}
\put(28,-8){$v$}
\put(30,60){\circle*{3}}
\put(28,63){$b$}
\end{picture}
\end{proof}

Since $X$ is a plane, there exists a triangle $abc\in X^3$ such that $X=\overline{\{a,b,c\}}$. Let $X_0=\{a,b,c\}$ and for every $n\in\w$, let $X_{n+1}\defeq\bigcup_{x,y\in X_n}\Aline xy$. It is clear that $X_n\subseteq X_{n+1}$ for all $n\in\w$, and $X=\overline{\{a,b,c\}}=\bigcup_{n\in \w}X_n$. 

\begin{claim}\label{cl:Xn-finite} For every $n\in\w$, the set $X_n$ is finite.
\end{claim}

\begin{proof} It is clear that the set $X_0=\{a,b,c\}$ is finite. Assume that for some $n\in\w$, the set $X_n$ is finite. Since $|X|_2$ is finite, the set $X_{n+1}=\bigcup_{x,y\in X_n}\Aline xy$ is finite of cardinality $|X_{n+1}|\le |X_n|+\tfrac12|X_n|\cdot (|X_n|-1)\cdot (|X|_2-2)$.
By the Principle of Mathematical Induction, the sets $X_n$ are finite for all $n\in\w$.
\end{proof} 

Claim~\ref{cl:Xn-finite} implies that the set $X=\bigcup_{n\in\w}X_n$ is at most countable. Then there exists an enumeration $(x_n)_{n<|X|}$ of the set $X$ by ordinals $<|X|\le \w$ such that $\{x_0,x_1,x_2\}=\{a,b,c\}$ and $x_n\in \bigcup_{i,j<n}\Aline {x_i}{x_j}$ for all numbers $n\in[2,|X|)$. 

\begin{claim}\label{cl:induction[xyz]=[abc]} For every number $n<|X|$, every triangle $xyz$ with $\{x,y,z\}\subseteq \{x_i\}_{i\le n}$ has area $[xyz]=[abc]$.
\end{claim}

\begin{proof} The proof is by induction. For the ordinal $n=2$, every triangle $xyz$ with $\{x,y,z\}\subseteq \{x_0,x_1,x_2\}=\{a,b,c\}$ has area $[xyz]=[abc]$, by Claim~\ref{cl:=[triangles]}. Assume that for some number $n<|X|$ and all numbers $m<n$ we know that all triangles $xyz$ with 
$\{x,y,z\}\subseteq \{x_i\}_{i\le m}$ have area $[xyz]=[abc]$. 
Choose any triangle $xyz$ with $\{x,y,z\}\subseteq \{x_i\}_{i\le n}$.  Since the points $x,y,z$ are distinct, there exist distinct numbers $i,j,k\le n$ such that $\{x_i,x_j,x_k\}=\{x,y,z\}$. If $\max\{i,j,k\}<n$, then $[xyz]=[x_ix_jx_k]=[abc]$, by the inductive assumption. So, assume that $\max\{i,j,k\}=n$. We lose no generality assuming that $k=n$ and hence $\max\{i,j\}<k=n$. By the choice of the sequence $(x_m)_{m<|X|}$, there exist two distinct points $u,v\in \{x_m\}_{m<n}$ such that $x_n\in \Aline uv$. Since $x_ix_jx_k$ is a triangle, at least one of the points $x_i$ or $x_j$ does not belong to the line $\Aline uv$. We lose no generality assuming that $x_j\notin \Aline uv$. By Claims~\ref{cl:[vab]=[vbc]} and \ref{cl:=[triangles]}, $[ux_jv]=[ux_jx_k]=[x_ix_jx_k]=[xyz]$. Since $\{u,v,x_j\}\subseteq \{x_m\}_{m\le n-1}$, the inductive assumption ensures that $[ux_jv]=[abc]$ and hence $[xyz]=[ux_jv]=[abc]$.
\end{proof}

Since $X=\{x_n\}_{n<|X|}$, Claim~\ref{cl:induction[xyz]=[abc]} implies that every triangle $xyz\in X^3$ has area $[xyz]=[abc]$.
\end{proof}

\section{Orientability of areas in liners}\label{s:orientable-area}

In this section we discuss the notion of orientability of areas in liners.


By Proposition~\ref{p:permutation-area}, for any area $\alpha\in [X^3]$ in a liner $X$ and any permutation $\sigma\in S_3$, the set $$\alpha\circ\sigma\defeq\{abc\circ\sigma:abc\in\alpha\}$$ is a well-defined area, equal to $[abc\circ \sigma]$ for any triple $abc\in\alpha$. The set 
$$\alpha S_3\defeq\{\alpha\circ \sigma:\sigma\in S_3\}$$is called the \index{area!$S_3$-orbit of}\defterm{$S_3$-orbit} of the area $\alpha$, and the subgroup 
$$S_3^{=\alpha}\defeq\{\sigma\in S_3:\alpha\circ\sigma=\alpha\}$$ is called the \index{area!stabilizer of}\defterm{stabilizer} of $\alpha$. 

It is clear that the cardinality $|\alpha S_3|$ of the $S_3$-orbit $\alpha S_3$ equals the index of the subgroup $S_3^{=\alpha}$ in the group $S_3$. Since the cardinality of each subgroup of $S_3$ divides the cardinality of the group $S_3$, we arrive to the following proposition.

\begin{proposition}\label{p:area1236} For any liner $X$ and any area $\alpha\in [X^3]$, its $S_3$-orbit $\alpha S_3$ has cardinality $$|\alpha S_3|=|S_3/S_3^{=\alpha}|\in\{1,2,3,6\}.$$
\end{proposition}

\begin{definition}\label{d:area-orientable} Let $X$ be a liner and $\kappa$ be  a cardinal. An area $\alpha\in [X^3]$ is called \index{area!$\kappa$-oriented}\index{$\kappa$-oriented area}\defterm{$\kappa$-oriented} if $\kappa=|\alpha S_3|$. 
\end{definition}

For $\kappa\in \{1,2,6\}$, $\kappa$-oriented areas have special names.

\begin{definition} An area $\alpha$ in a liner $X$ is called
\begin{itemize}
\item \index{unoriented area}\index{area!unoriented}\defterm{unoriented} if $|\alpha S_3|=1$  iff $\alpha$ is $1$-oriented;
\item \index{oriented area}\index{area!oriented}\defterm{oriented} if $|\alpha S_3|=2$ iff $\alpha$ is $2$-oriented;
\item \index{semioriented area}\index{area!semioriented}\defterm{semioriented} if $|\alpha S_3|\le 2$ iff $\alpha$ is $1$-oriented or $2$-oriented;
\item \index{hexaoriented area}\index{area!hexaoriented}\defterm{hexaoriented} if $|\alpha S_3|=6$  iff $\alpha$ is $6$-oriented.
\end{itemize}
\end{definition}

\begin{exercise} Show that for a triangle $abc$ in a liner $X$, the area $[abc]$ is
\begin{itemize}
\item unoriented if and only if $[abc]=[bca]=[cab]=[acb]=[bac]=[cba]$;
\item oriented if and only if $[abc]=[bca]=[cab]\ne[acb]=[bac]=[cba]$;
\item semioriented if and only if $[abc]=[bca]=[cab]$;
\item hexaoriented if and only if the six areas $[abc]$, $[bca]$, $[cab]$, $[acb]$, $[bac]$, $[cba]$ all are distinct.
\end{itemize}
\end{exercise}

\begin{exercise} Show that the zero area in any liner is unoriented.
\end{exercise}

\begin{definition} Let $X$ be a liner. For a semioriented area $\alpha\in [X^3]$, its $S_3$-orbit $\alpha S_3$ has cardinality $|\alpha S_3|\le 2$ and hence there exists a unique area $-\alpha\in[X^3]$ such that $\alpha S_3=\{-\alpha,\alpha\}$. The area $-\alpha$ is called \index{opposite area}\index{area!opposite}\defterm{the opposite area} to $\alpha$. 
\end{definition}

\begin{exercise} Show that a semioriented area $\alpha$ is
\begin{itemize}
\item unoriented if and only if $-\alpha=\alpha$;
\item oriented if and only if $-\alpha\ne\alpha$.
\end{itemize}
\end{exercise}

\begin{proposition}\label{p:opposite-area} Let $abc$ be a triangle in a liner $X$. If the area $\alpha\defeq [abc]$ of the triangle $abc$ is semioriented, then $$\alpha=[abc]=[bca]=[cab]\quad\mbox{and}\quad -\alpha=[acb]=[bac]=[cba].$$
\end{proposition}

\begin{proof} Assume that the area $\alpha=[abc]$ is semioriented, which means that $|\alpha S_3|\le 2$. Since $|\alpha S_3|=|S_3/S_3^{=\alpha}|$, the stabilizer $S_3^{=\alpha}$ of $\alpha$ has cardinality $3$ or $6$ and hence $S_3$ contains the cyclic subgroup $C_3=\{012,120,201\}$. Then $\alpha=\alpha\circ\sigma$ and $-\alpha=\-\alpha\circ \sigma$ for all $\sigma\in C_3$, which implies $\alpha=[abc]=[bca]=[cab]$ and $-\alpha=[acb]=[bac]=[cba]$.
\end{proof} 

The following propositions provide simple sufficient conditions of semiorientability and unorientability of areas in liners.

\begin{proposition}\label{p:semioriented<=} An area $\alpha$ in a liner $X$ is semioriented whenever some triangle $abc\in\alpha$ can be completed to a parallelogram $abcd$ in $X$. 
\end{proposition}

\begin{proof} Assume that some triangle $abc\in\alpha$ can be completed to a parallelogram $abcd$. By Proposition~\ref{p:12-permutations}, $\alpha=[abc]=[bca]=[cab]$, which implies that the stabilizer $S_3^{=\alpha}$ of $\alpha$ has cardinality at least $3$ and hence $|S_3^{=\alpha}|\in\{3,6\}$, which implies $|\alpha S_3|=|S_3/S_3^{=\alpha}|\in\{2,1\}$ and witnesses that the area $\alpha$ is semioriented.
\end{proof} 

\begin{proposition}\label{p:unoriented<=} An area $\alpha$ in a liner $X$ is unoriented whenever some triangle $abc\in\alpha$ can be completed to a Boolean parallelogram $abcd$. 
\end{proposition}

\begin{proof} Assume that some triangle $abc\in\alpha$ can be completed to a Boolean parallelogram $abcd$. By Proposition~\ref{p:12-permutations}, $\alpha=[abc]=[acd]=[bac]=[bca]=[cab]=[cba]$, which implies that the stabilizer $S_3^{=\alpha}$ of $\alpha$ coincides with the group $S_3$ and hence $|\alpha S_3|=|S_3/S_3^{=\alpha}|=1$, witnessing that the area $\alpha$ is unoriented.
\end{proof} 

Propositions~\ref{p:semioriented<=} and \ref{p:Steiner-paragram} imply the   the following sufficent condition of semiorientability.

\begin{corollary}\label{c:area-semioriented-in-Playfair} Any area in any $3$-ranked affine liner is semioriented. In particular, any area in a Playfair liner is semioriented.
\end{corollary}



On the other hand, Propositions~\ref{p:large=>Boolean-parallelogram} and \ref{p:unoriented<=} imply the following sufficent condition of non\-orientability.

\begin{corollary}\label{c:large=>unoriented} If every plane $P$ in a $2$-balanced $3$-ranked liner $X$ has cardinality\newline $|P|+9{\cdot}|X|_2>3{\cdot}|X|_2^2+9$, then every area in $X$ is unoriented.
\end{corollary}

\begin{proof} Take any area $\alpha\in[X^3]$ in $X$. If $\alpha=[0]$, then $\alpha$ is unoriented by the definition of the zero area. So, assume that $\alpha\ne[0]$ and choose any triangle $abc\in\alpha$. Consider the plane $P\defeq\overline{\{a,b,c\}}$. 

Consider the sequence of sets $(P_n)_{n\in\w}$, defined by the recursive formula: $P_0\defeq\{a,b,c\}$ and $P_{n+1}\defeq\bigcup_{x,y\in P_n}\Aline xy$ for all $n\in\w$. It is clear that $P=\bigcup_{n\in\w}P_n$. 

Assuming that the cardinal $|X|_2$ is infinite, we shall prove that $|P_n|\le|X|_2$ for all $n\in\w$. Indeed, for $n=0$ we have $|P_0|=|\{a,b,c\}|=3<|X|_2$. Assume that for some $n\in\w$, the inequality $|P_n|\le|X|_2$ holds. Then $|P_{n+1}|\le|P_n|^2\cdot|X|_2 \le|X|_2^3=|X|_2$ (because the cardinal $|X|_2$ is infinite). Therefore, $|P_n|\le|X|_2$ for all $n\in\w$, by the Principle of Mathematical Induction. Since $|X|_2$ is infinite, the plane $P=\bigcup_{n\in\w}P_n$ has cardinality $|X|_2\le |P|\le |X|_2\cdot\w=|X|_2$. Therefore, $|P|=|X|_2$ and then $|P|+9{\cdot}|X|_2=|X|_2=3{\cdot}|X|_2^2+9$, which contradicts the strict inequality $|P|+9{\cdot}|X|_2>3{\cdot}|X|_2^2+9$ in the assumption. This contradiction shows that the cardinal $|X|_2$ is finite and hence the liner $X$ is line-finite. In this case, the strict inequality $|P|+9{\cdot}|X|_2>3{\cdot}|X|_2^2+9$ is equivalent to $|P|>3{\cdot}|X|_2^2-9{\cdot}|X|_2+9=3{\cdot}|P|_2^2-9{\cdot}|P|_2+9$. Applying Proposition~\ref{p:large=>Boolean-parallelogram}, we conclude that the triangle $abc$ can be completed to a Boolean parallelogram $abcd$ in the plane $P$. By Proposition~\ref{p:unoriented<=}, the area $\alpha=[abc]$ is unoriented.
\end{proof}

\begin{exercise} Show that any area in a projective liner is hexaoriented.
\end{exercise}

\begin{exercise} Show that any area in a punctured projective plane is hexaoriented.
\end{exercise}

\begin{exercise} Let $X$ be a Playfair plane with attached point at infinity.
Show that every area in $X$ is either semioriented or hexaoriented. 
\end{exercise}

\begin{exercise} Show that every area in the Moulton plane (see Section~\ref{s:Moulton}) is oriented.
\end{exercise}

\begin{exercise} Show that every area in the Beltrami--Klein hyperbolic liner (see Example~\ref{ex:Beltrami-Klein}) is oriented.
\end{exercise}

\section{The absolute area of a triangle in a liner}

In this section we introduce the notion of absolute area of a triangle in a liner.

\begin{definition} For a triangle $abc$ in a liner $X$, its \index{absolute area}\index{triangle!absolute area}\defterm{absolute area} $$[\![abc]\!]\defeq[abc]\cup[acb]\cup[bac]\cup[bca]\cup [cab]\cup[cba]$$ is the family of all triangles $xyz$ in $X$ for which there exists a sequence of triangles $\tau_0,\tau_1,\dots,\tau_n$ in $X$ such that  $xyz=\tau_0\boxtimes\tau_1\boxtimes\tau_2\boxtimes\dots \boxtimes \tau_n\in\{abc,acb,bac,bca,cab,cba\}$. For a triple of collinear points $abc\in X^3\setminus X^\vartriangle$, let $[\![abc]\!]$ be equal to the zero absolute area $[\![0]\!]\defeq X^3\setminus X^\vartriangle=[0]$. The set of absolute areas \index[note]{$[\hskip-1.5pt[X^3]\hskip-1.5pt]$}$$[\![X^3]\!]\defeq\{[\![abc]\!]:abc\in  X^3\}$$ is called \index{absolute areol}\defterm{the absolute areol} of the liner $X$. 
\end{definition}

Corollary~\ref{c:same-area=>same-plane} implies its own counterpart for absolute areas.

\begin{corollary}\label{c:same-absarea=>same-plane} If two triangles $abc$ and $xyz$ in a liner $X$ have the same absolute area, then $\overline{\{a,b,c\}}=\overline{\{x,y,z\}}$.
\end{corollary}

\begin{exercise}\label{ex:same-absarea} Let $X$ be a liner. Show that any two triples $abc,xyz\in X^3$ with $\{a,b,c\}=\{x,y,z\}$ have the same absolute area.
\end{exercise}

\begin{exercise} Show that the absolute area $[\![abc]\!]$ of any triangle $abc$ in any projective liner equals the set $\{abc,acb,bac,bca,cab,cba\}$. Deduce from this fact that any projective liner $X$ has absolute areol $[\![X]\!]=\{X^3\setminus X^\vartriangle\}\cup\big\{\{abc,acb,bac,bca,cab,cba\}:abc\in X^\vartriangle\big\}$.\end{exercise}


Next, we show that the absolute area of a triangle in a liner coincides with the modulus of its area. 


\begin{definition} For an area $\alpha\in [X^3]$ in a liner $X$, its 
 \index{modulus of area}\index{area!modulus of}\defterm{modulus} $[\alpha]$ is defined as the set $$[\alpha]\defeq\textstyle{\bigcup\alpha S_3}=\{abc\circ \sigma:abc\in\alpha\;\wedge\;\sigma\in S_3\}=\{abc\in X^3:\{abc,acb,bac,bca,cab,cba\}\cap\alpha\ne\varnothing\}.$$
\end{definition}

\begin{exercise} Show that the zero area $[0]=X^3\setminus X^\vartriangle$ in any liner $X$ is equal to its modulus $[[0]]$.
\end{exercise} 

\begin{exercise} Show that any area $\alpha$ in a liner is equal to its modulus $[\alpha]$ if and inly if $\alpha$ is unoriented. 
\end{exercise}

\begin{exercise} Show that for a semioriented area $\alpha$ in a liner, its modulus $[\alpha]$ is equal to the union $-\alpha\cup\alpha$.
\end{exercise}


\begin{proposition}\label{p:absarea=modulus} For any liner $X$ and any triple $abc\in X^3$ its absolute area $[\![abc]\!]$ equals the modulus $[[abc]]$ of its area $[abc]$.
\end{proposition}

\begin{proof} We should prove that $[\![abc]\!]=[[abc]]$. If $abc$ is not a triangle, then $[\![abc]\!]=[\![0]\!]=X^3\setminus X^\vartriangle=\bigcup([0]S_3)=[[0]]$ and we are done.  So, assume that $abc$ is a triangle.
\smallskip

To prove that $[\![abc]\!]\subseteq [[abc]]$, take any triangle $xyz\in[\![abc]\!]$ and find a permutation $\sigma\in S_3$ such that $xyz\in [abc\circ\sigma]$ and hence $[xyz]=[abc\circ\sigma]=[abc]\circ\sigma\subseteq \bigcup ([abc]S_3)=[[abc]]$, witnessing that $[\![abc]\!]\subseteq [[abc]]$.
\smallskip

To prove that $[[abc]]\subseteq[\![abc]\!]$, take any triangle $xyz\in [[abc]]$ and find a permutation $\sigma\in S_3$ such that $xyz\in [abc]\circ\sigma=[abc\circ\sigma]\subseteq \bigcup[abc] S_3=[abc]\cup[acb]\cup[bac]\cup[bca]\cup[cab]\cup[cba]=[\![abc]\!]$.
\end{proof}

Proposition~\ref{p:absarea=modulus} implies the following (almost) trivial corollary.

\begin{corollary}\label{c:same-area=>same-absarea} Let $abc$ and $xyz$ be two triangles in a liner $X$. If $[abc]=[xyz]$, then $[\![abc]\!]=[\![xyz]\!]$.
\end{corollary}

Proposition~\ref{p:absarea=modulus} also implies that the modulus map $$[\cdot]:[X^3]\to[\![X^3]\!],\quad [\cdot]:\alpha\mapsto[\alpha],$$is a well-defined surjective map, and hence $|[\![X^3]\!]|\le|[X^3]|$.  



Corollary~\ref{c:same-area=>same-absarea} and Proposition~\ref{p:12-permutations} imply the following corollary.

\begin{corollary}\label{c:24triangles} Let $abcd$ be a parallelogram in a liner. All $24$ triangles with vertices in the set $\{a,b,c,d\}$ have the same absolute area.
\end{corollary}

Theorem~\ref{t:large=>area-trivial} and Corollary~\ref{c:same-area=>same-absarea} imply the following corollary.

\begin{corollary}\label{t:large=>area-abstrivial} If a balanced ranked plane $X$ has $|X|+12{\cdot}|X|_2>4{\cdot}|X|^2_2+11$, then all triangles in $X$ have the same absolute area.
\end{corollary}

\section{Area-surjective and area-injective liners}

In this section we introduce area-surjective, area-injective and area-bijective liners. The definition of such liners is motivated by the following important theorem on construction of a triangle of a given area on a given base in a hyper-Bolyai liner. Let us recall that a liner $X$ is \index{hyper-Bolyai liner}\index{liner!hyper-Bolyai}\defterm{hyper-Bolyai} if for any plane $P$ in $X$, concurrent lines $A,B$ in $P$, any any point $p\in P\setminus B$, there exists a line $L\subseteq P$ such that $p\in L$, $L\parallel B$ and $|L\cap A|=1$. By Proposition~\ref{p:hyper-Bolyai-interplay}, every hyper-Bolyai liner is Playfair.

\begin{theorem}\label{t:area-move} Let $C,B$ be two concurrent lines in a  hyper-Bolyai plane $X$, and let $a,b\in B$ be two distinct points. Then for every triangle $xyz$ in $X$, there exists a point $c\in C$ such that $xyz\boxtimes^7 abc$ and hence $[xyz]=[abc]$.
\end{theorem}

\begin{proof} We shall deduce this theorem from the following four lemmas.

\begin{lemma}\label{l:area-move1} Let $xyz$ be a triangle in $X$. If $x\notin \Aline ab$ and $z\notin \Aline ax$, then there exists a point $c\in X$ such that $xyz\boxtimes^5 abc$.
\end{lemma}

\begin{proof} By the assumption, the line $L\defeq \Aline ax$ does not contain the points $b$ and $z$. Since $X$ is hyper-Bolyai and $L\cap \Aline xz=\{x\}$, there exists a line $L'\subset X$ such that $b\in L'$, $L'\parallel L$ and $|L'\cap \Aline xz|=1$.
Since $X$ is hyper-Bolyai and $|L'\cap\Aline xz|=1$, there exists a line $L_y\subset X$ such that $y\in L_y$, $L_y\parallel \Aline xz$ and the lines $L_y,L'$ have a unique common point $y'$. Since $X$ is hyper-Bolyai and $L\cap\Aline {y'}x=\{x\}$, there exists a line $L_z\subset X$ such that $z\in L_z$, $L_z\parallel \Aline{y'}x$ and the lines $L_z,L$ have a unique common point $z'$.

\begin{picture}(200,90)(-120,-15)
\linethickness{0.8pt}
{\color{lightgray}
\polygon*(0,30)(30,30)(30,0)
\polygon*(60,60)(75,60)(75,0)
}

\put(0,60){\color{teal}\line(1,0){100}}
\put(30,0){\color{teal}\line(1,0){70}}
\put(105,-3){\color{teal}$L$}
\put(105,57){\color{teal}$L'$}

{\linethickness{0.5pt}
\put(0,30){\line(1,0){30}}
\put(0,30){\line(1,-1){30}}
\put(30,30){\color{gray}\line(-1,1){30}}
\put(0,60){\color{gray}\line(3,-4){45}}
\put(45,0){\color{gray}\line(1,4){15}}
}
\put(0,30){\color{cyan}\line(0,1){30}}
\put(30,0){\color{cyan}\line(0,1){30}}
\put(30,0){\color{blue}\line(-1,2){30}}
\put(45,0){\color{blue}\line(-1,2){15}}

\put(30,0){\color{orange}\line(1,2){30}}
\put(45,0){\color{orange}\line(1,2){30}}
\put(60,60){\line(1,-4){15}}
\put(75,60){\line(0,-1){60}}

\put(30,0){\circle*{3}}
\put(28,-8){$x$}
\put(45,0){\circle*{3}}
\put(42,-8){$z'$}
\put(75,0){\circle*{3}}
\put(72,-8){$a$}
\put(0,60){\circle*{3}}
\put(-3,66){$y'$}
\put(0,30){\circle*{3}}
\put(-8,28){$y$}
\put(30,30){\circle*{3}}
\put(28,33){$z$}
\put(60,60){\circle*{3}}
\put(58,63){$b$}
\put(75,60){\circle*{3}}
\put(72,63){$c$}
\end{picture}

 Since $X$ is hyper-Bolyai and $L'\cap \Aline xb=\{b\}$, there exists a line $L_c\subseteq X$ such that $z'\in L_c$, $L_c\parallel \Aline xb$ and the lines $L_c,L'$ have a common point $c$. Then $xyz\boxtimes xy'z\boxtimes xy'z'\boxtimes xbz'\boxtimes xbc\boxtimes abc$, which implies  $xyz\boxtimes^5 abc$.
\end{proof}

\begin{lemma}\label{l:area-move2} Let $xyz$ be a triangle in $X$. If $x\notin \Aline ab$, then there exist a point $c\in X$ such that $xyz\boxtimes^5 abc$.
\end{lemma}

\begin{proof} If $z\notin \Aline ax$, then put $z'\defeq z$. If $z\in \Aline ax$, then $y\notin \Aline ax$ and the lines $\Aline zx$ and $\Aline xy$ are concurrent. Since $X$ is hyper-Bolyai, there exists a line $L_x$ in $X$ such that $z\in L_z$ and  $L_z\parallel \Aline xy$. Choose any point $z'\in L_z\setminus\{z\}$ and observe that $z'\notin \Aline ax$ and $xyz\boxtimes \tau_0\defeq xyz'$. In both cases, we have found a point $z'\in X\setminus\Aline ax$ such that $xyz\boxtimes xyz'$. By Lemma~\ref{l:area-move1}, there exist a point $c\in X$ such that $xyz'\boxtimes^5 abc$. Then $xyz\boxtimes xyz'\boxtimes^5 abc$ implies $axy\boxtimes^6 abc$. 
\end{proof}

By analogy we can prove the following counterpart of Lemma~\ref{l:area-move2}.

\begin{lemma}\label{l:area-move3} Let $xyz$ be a triangle in $X$. If $y\notin \Aline ab$, then there exist a point $c\in X$ such that $xyz\boxtimes^6 abc$.
\end{lemma}

\begin{lemma}\label{l:area-move4} Let $xyz$ be a triangle in $X$. If $x,y\in \Aline ab$, then there exist a point $c\in X$ such that $xyz\boxtimes^6 abc$.
\end{lemma}

\begin{proof} Taking into account that $xyz$ is a triangle with $x,y\in\Aline ab$, we conclude that $z\notin\Aline ab$. Since $X$ is hyper-Bolyai, and $\Aline za\cap \Aline yz=\{z\}$, there exists a line $L_x$ in $X$ such that $x\in L_x$, $L_x\parallel \Aline yz$ and $|L_x\cap\Aline xa|=1$. It follows from $L_x\parallel \Aline yz$ that $L_x\cap \Aline ab=\{x\}$. By Theorem~\ref{t:hyper-Bolyai<=>}, the hyper-Bolyai liner $X$ is $3$-long and hence the line $L_x$ contains a point $x'\notin \Aline za\cup\Aline ab$. Then $x'\notin\Aline ab$ and $z\notin \Aline {x'}a$. The choice of the line $L_x\parallel \Aline yz$ ensures that $xyz\boxtimes \tau_0\defeq x'yz$. By Lemma~\ref{l:area-move1}, there exist a point $c\in X$such that $x'yz\boxtimes^5 abc$. 
Then $xyz\boxtimes x'yz\boxtimes^5 abc$ implies $xyz\boxtimes^6 abc$.
\end{proof}

Now we are able to complete the proof of the theorem. Let $xyz$ be any triangle in $X$. By Lemmas~\ref{l:area-move2}, \ref{l:area-move3}, \ref{l:area-move4}, there exist a point $c'\in X$ such that $xyz\boxtimes^6 abc$. Since the lines $\Lambda$ and $\Aline ab$ are concurrent in the hyper-Bolyai plane $X$, there exists a line $L_c$ in $X$ such that $c'\in L_c$, $L_c\parallel \Aline ab$ and the lines $L_c,\Lambda$ have a common point $c$. Then $xyz\boxtimes^6 abc'\boxtimes abc$ and finally, $xyz\boxtimes^7 abc$.
\end{proof}

Theorem~\ref{t:area-move} motivates the following definitions.

\begin{definition} A liner $X$ is defined to be 
\begin{itemize}
\item \index{area-surjective liner}\index{liner!area-surjective}\defterm{area-surjective} if for any line $\Lambda\subset X$ and points $\lambda\in \Lambda$ and  $o\in X\setminus \Lambda$, the function $\Lambda\to[X^3]$, $x\mapsto [o\lambda x]$, is surjective; 
\item \index{area-injective liner}\index{liner!area-injective}\defterm{area-injective} if for any line $\Lambda\subset X$ and points $\lambda\in \Lambda$ and  $o\in X\setminus \Lambda$, the function $\Lambda\to[X^3]$, $x\mapsto [o\lambda x]$, is injective; 
\item \index{area-bijective liner}\index{liner!area-bijective}\defterm{area-bijective} if for any line $\Lambda\subset X$ and points $\lambda\in \Lambda$ and  $o\in X\setminus \Lambda$, the function $\Lambda\to[X^3]$, $x\mapsto [o\lambda x]$, is bijective; 
\index{area-compatible liner}\index{liner!area-compatible}
\item \defterm{area-compatible}  if for any areas $\alpha,\beta\in[X^\vartriangle]$, there exist a line $\Lambda\subseteq X$ and points $a,b,c\in \Lambda$ and $o\in X$ such that $[oab]=\alpha$ and $[obc]=\beta$.
\end{itemize} 
\end{definition}

\begin{proposition}\label{p:area-surjective-around} Let $X$ be a liner.
\begin{enumerate}
\item[\textup{(1)}] If $X$ is a hyper-Bolyai plane, then $X$ is area-surjective.
\item[\textup{(2)}] If $X$ is area-surjective, then $X$ is area-compatible.
\item[\textup{(3)}] If $X$ is area-compatible, then $X$ is ranked of rank $\le 3$.
\item[\textup{(4)}] If $X$ is area-injective, then $[abc]=\{xyz\in X^3:xyz\boxtimes^7 abc\}$ for every triple $abc\in X^3$.
\end{enumerate}
\end{proposition}  

\begin{proof} 1. Assume that $X$ is a hyper-Bolyai plane. To prove that $X$ is area-surjective, take any line $\Lambda\subset X$ and points $\lambda\in\Lambda$ and $o\in X\setminus \Lambda$. We should check that the function $h:\Lambda\to[X^\vartriangle]$, $h:x\mapsto [o\lambda x]$ is surjective. Given any area $\alpha\in [X^\vartriangle]$, we should find a point $c\in\Lambda$ such that $[o\lambda c]=\alpha$. If $\alpha=[0]$, then $h(\lambda)=[o\lambda\lambda]=[0]=\alpha$ and we are done. So, assume that $\alpha\ne[0]$, and find a triangle $xyz\in\alpha$. By Theorem~\ref{t:area-move}, there exists a point $c\in \Lambda$ such that $h(c)=[o\lambda c]=[xyz]=\alpha$, witnessing that the function $h$ is surjective.
\smallskip

2. Assume that $X$ is surjective. To prove that $X$ is area-compatible, take any areas $\alpha,\beta\in [X^\vartriangle]$. Choose any line $L$ in $X$ and any points $a\in L$ and $o\in X\setminus L$. Since $X$ is area-surjective, the function $L\mapsto [X^\vartriangle]$, $x\mapsto [oax]$, is surjectve and hence there exists a point $b\in L$ such that $[oab]=\alpha$.  By the surjectivity of the function $L\mapsto [X^\vartriangle]$, $x\mapsto [obx]$,  there exists a point $c\in L$ such that $[obc]=\beta$. This shows that $X$ is area-compatible.
\smallskip

3. Let $X$ be an area-compatible liner. To prove that $X$ is a ranked liner of rank $\|X\|\le 3$, it suffices to show that each plane $\Pi$ in $X$ coincides with $X$. Fix any triangle $xyz$ in the plane $\Pi$. To prove that $\Pi=X$, take any point $p\in X$.  Since the liner $X$ is area-compatible, there exist a line $L\subseteq X$ and points $a,b,c\in L$ and $o\in X$ such that $[oab]=[xyz]$ and $[obc]=[xyp]$. Since $xyz$ and $xyp$ are triangles, the equalities $[obc]=[xyp]$ and $[oab]=[xyz]$ imply $p\in\overline{\{x,y,p\}}=\overline{\{o,b,c\}}=\overline{L\cup\{o\}}=\overline{\{o,a,b\}}=\overline{\{x,y,z\}}\subseteq \Pi$, by Corollary~\ref{c:same-area=>same-plane}. Therefore,  $X=\Pi$.
\smallskip

4. Assume that $X$ is an area-injective liner. Given any triple $abc\in X^3$, we will show that $[abc]=\{xyz\in X^3:xyz\boxtimes^7 abc\}$. The inclusion $\{xyz\in X^3:xyz\boxtimes^7 abc\}\subseteq\{xyz\in X^\vartriangle:xyz\boxtimes^\w  abc\}=[abc]$ is trivial and holds for all (not necessarily area-injective) liners. To see that $[abc]\subseteq \{xyz\in X^3:xyz\boxtimes^7 abc\}$, take any triple $xyz\in[abc]$. If $\|\{x,y,z\}\|\le 2$, then $[0]=[xyz]=[abc]$, and then $xyz\boxtimes abc$ and $xyz\boxtimes ^7 abc$ by definition of the relation $\boxtimes$ for triples which are not triangles. So, assume that $xyz$ is a triangle and so is the triple $abc\boxtimes^\w  xyz$. By Theorem~\ref{t:area-move}, there exists a point $c'\in \Aline bc$ such that $xyz\boxtimes^7 abc'$ and hence $[abc']=[xyz]=[abc]$. The area-injectivity of $X$ ensures that $c'=c$ and hence $xyz\boxtimes^7 abc'=abc$. 
\end{proof}

Proposition~\ref{p:area-surjective-around}(1) implies the following corollary.

\begin{corollary} A hyper-Bolyai plane $X$ is area-injective if and only if it is area-bijective.
\end{corollary}

\begin{proposition}\label{p:area-inj=>Proclus} Every area-injective $3$-ranked liner $X$ is Proclus.
\end{proposition}

\begin{proof} By Theorem~\ref{t:Proclus<=>}, it suffices to show that $X$ is proaffine and $3$-proregular. Assuming that $X$ is not proaffine, we can find points $o,x,y\in X$, $p\in\Aline xy\setminus\Aline ox$ and two distinct points $u,v\in \Aline oy$ such that $\Aline pu\cap\Aline ox=\varnothing=\Aline pv\cap\Aline ox$ and hence $\Aline pu\parallel\Aline ox\parallel \Aline pv$, by the $3$-rankedness of $X$. Then $[xou]=[xop]=[xov]$, which contradicts the area-injectivity of $X$. This contradiction shows that the liner $X$ is proaffine.
\smallskip

Next, we show that $X$ is $3$-proregular. Given any triangle $uow\in X^3$ with $|\Aline ow|\ge 3$, we need to show that $\overline{\{u,o,w\}}=\bigcup_{x\in\Aline ou}\bigcup_{y\in\Aline ow}\Aline xy$. In the opposite case, we can find a triangle $uow\in X^3$ such that $|\Aline ow|\ge 3$ and the plane $\overline{\{u,o,w\}}$ contains a point $z\notin \bigcup_{x\in\Aline ou}\bigcup_{y\in\Aline ow}\Aline xy$. Since $|\Aline ow|\ge 3$, there exists a point $v\in \Aline ow\setminus\{o,w\}$. The choice of the point $z$ ensures that $\Aline wz\cap\Aline ou=\varnothing=\Aline vz\cap\Aline ou$ and hence $\Aline wz\parallel \Aline ou\parallel \Aline vz$, by the $3$-rankedness of $X$ and Corollary~\ref{c:parallel-lines<=>}. Then $[uov]=[uoz]=[uow]$, which contradicts the area-injectivity of $X$. This contradiction shows that the proaffine liner $X$ is $3$-proregular and hence Proclus, by Theorem~\ref{t:Proclus<=>}.
\end{proof}

\begin{corollary}\label{c:area-bij=>Proclus} Every area-bijective liner is Proclus of rank $\|X\|\le 3$.
\end{corollary}

\begin{proof} Let $X$ be an area-bijective liner. By Proposition~\ref{p:area-surjective-around}, the area-surjective liner $X$ is ranked of rank $\|X\|\le 3$. By Proposition~\ref{p:area-inj=>Proclus}, the ranked area-injective liner $X$ is Proclus. 
\end{proof}

\begin{corollary}\label{c:area-inj+hyperaff=>Playfair} Every area-injective hyper-Bolyai liner $X$ is Playfair.
\end{corollary}

\begin{proof} By Theorem~\ref{t:hyper-Bolyai<=>}, the hyper-Bolyai liner $X$ is $3$-ranked. By Proposition~\ref{p:area-inj=>Proclus}, the $3$-ranked area-injective liner $X$ is Proclus, and hence Playfair (being Proclus and Bolyai simultaneously).
\end{proof}

\begin{theorem}\label{t:area-bij=>Playfair}  If an area-bijective liner has a projective completion, then it is Playfair of rank $\le 3$.
\end{theorem}

\begin{proof} Let $X$ be an area-bijective liner. By Corollary~\ref{c:area-bij=>Proclus}, the area-bijective liner $X$ is a Proclus liner of rank $\|X\|\le 3$. If $\|X\|\le 2$, then $X$ is Playfair (because it contains no planes). So, assume that $\|X\|=3$ and hence $X$ is a Proclus plane. Assume that $X$ has a projective completion $P$. By Corollary~\ref{c:procompletion-rank} and Theorem~\ref{t:Proclus-not-3long}, we can assume that $\|P\|=\|X\|=3$. By Theorem~\ref{t:proaffine<=>proflat}, the horizon $H\defeq P\setminus X$ of $X$ in $P$ is a proflat set in $P$ with $\|H\|<\|P\|=3$. Then the proflat set $H$ is  empty, a singleton, a line, or a punctured line.
\smallskip

If $H=\varnothing$, then $X=P$ is a projective plane. Since $X$ contains no disjoint lines, the area $[abc]$ of every triangle $abc$ in $X$ is the singleton $\{abc\}$. Fix any triangle $abc$ in $X=P$ and consider the function $h:\Aline bc\to[X^3]$, $h:x\mapsto [abx]$. Observe that $[bac]=\{bac\}\notin \big\{\{abx\}:x\in\Aline bc\}=h[\Aline bc]\big\}$, witnessing that $X$ is not area-surjective and hence not area-bijective.
\smallskip

If $H$ is a singleton containing a unique point $p\in P$, then every triangle $abc$ in $X=P\setminus\{p\}$ has area
$$[abc]=\begin{cases}\{abc\}&\mbox{if $p\notin\Aline ab\cup\Aline bc\cup\Aline ac$};\\
\big\{abx:x\in \Aline cp\setminus\{p\}\big\}&\mbox{if $p\in \Aline ab$};\\
\big\{axc:x\in \Aline bp\setminus\{p\}\big\}&\mbox{if $p\in \Aline ac$};\\
\big\{xbc:x\in \Aline ap\setminus\{p\}\big\}&\mbox{if $p\in \Aline bc$}.
\end{cases}
$$
Since $P$ is a plane, there exist points $u\in P\setminus\{p\}$ and $v\in P\setminus\Aline up$. Then $L\defeq \Aline uv$ is a line in $P$. By Theorem~\ref{t:projective<=>}, the projective plane $P$ is strongly regular and hence $P=\overline{\{a,b,p\}}=\bigcup_{x\in\Aline ab}\Aline xp$. Assuming that $\Aline xp=\{x,p\}$ for all $x\in L$, we conclude that $X=P\setminus\{p\}=L$ is not a plane, which contradicts our assumption. This contradiction shows that $\Aline bp\ne\{b,p\}$ for some point $b\in L$. So, we can choose a point $a\in \Aline bp\setminus\{b,p\}$. Also choose any point $c\in L\setminus\{b\}$ and observe that $abc$ is a triangle in the liner $X$. Consider the function $h:\Aline bc\setminus\{p\}\to [X^3]$, $h:y\mapsto [aby]$, and observe that the set $h[\Aline bc\setminus\{p\}]=\big\{[aby]:y\in \Aline bc\setminus\{p\}\big\}=\{[0]\}\cup\big\{\{axy:x\in \Aline bp\setminus\{p\}\}:y\in\Aline bp\setminus\{b,p\}\big\}$ does not contains the area $[bac]=\big\{bxc:x\in\Aline ap\setminus\{p\}\big\}$, witnessing that the function $h$ is not surjective and hence the liner $X$ is not area-surjective.
\smallskip

Finally, assume that $H$ is a punctured line in $P$ and hence $H\cup\{p\}$ is a line in $P$ for some point $p\in P\setminus H$. Observe that for any $x\in X\setminus\{p\}$, the line $\Aline xp$ intersects any other line in $X$. This implies that any triangle $abc$ in $X$ with $p\in\{a,b,c\}$ has area $[abc]=\{abc\}$. Choose any points $b,c\in X$ such that $pbc$ is a triangle. Consider the function $h:\Aline bc\to[X^3]$, $h:x\mapsto pbx$, and observe that the set $h[\Aline bc]=\{[pby]:y\in \Aline bc\}=\{[0]\}\cup\big\{\{pby\}:y\in \Aline bc\setminus\{b\}\big\}$ does not contains the area $[abp]=\{abp\}$, witnessing that the function $h$ is not surjective and hence the liner $X$ is not area-surjective.
\smallskip

Therefore, $H$ is a line in the $3$-long projective plane $P$. By Corollary~\ref{c:Avogadro-projective}, the $3$-long projective plane $P$ is $2$-balanced. If $|P|_2=3$, then $|X|_2=2$ and $\|X\|=|X|=|P\setminus H|=|P|-|H|=7-3=4$, which contradicts our assumption that $X$ is a plane. This contradiction shows that $|P|_2\ge 4$ and hence the liner $X=P\setminus H$ is $3$-long. Since the horizon $P\setminus X=H$ in $X$ is a line (and hence hyperplane) in $P$, Theorem~\ref{t:affine<=>hyperplane} ensures that the subliner  $X=P\setminus H$ is affine and regular. By Theorem~\ref{t:Playfair<=>}, the $3$-long affine regular liner $X$ is Playfair. 
\end{proof}

\begin{corollary}\label{c:area-bij=>3long} Every area-bijective plane $X$ of cardinality $|X|>2$ is $3$-long.
\end{corollary}

\begin{proof} Let $X$ be an area-bijective plane. By Corollary~\ref{c:area-bij=>Proclus}, $X$ is Proclus. Assuming that $X$ is not $3$-long, we can apply Theorem~\ref{t:Proclus-not-3long} and conclude that $X=P\setminus H$ for some Steiner projective plane $P$ and some subset $H\subseteq P$ of cardinality $|H|\le 2$. By Corollary~\ref{c:Steiner-projective<=>}, the Steiner projective plane $P$ has cardinality $|P|=7$. Since $P$ is a projective completion of $X$, we can apply Theorem~\ref{t:area-bij=>Playfair} and conclude that $X$ is Playfair. By Theorem~\ref{t:Playfair<=>}, the Playfair liner $X$ is $3$-long. It follows from $3\le|X|_2\le|P|_2=3$ that the Playfair liner $X$ is Steiner. By Corollary~\ref{c:Playfair<=>23balanced}, $|X|=|X|_2^2=9$ which contradicts $|X|\le|P|=7$.
This contradiction shows that the plane $X$ is $3$-long.
\end{proof}

\begin{corollary}\label{c:area-bij=>Playfair} Every area-bijective finite liner is Playfair of rank $\le 3$.
\end{corollary}

\begin{proof} Let $X$ be an area-bijective finite liner. By Proposition~\ref{p:area-surjective-around}, the area-surjective liner $X$ is a Proclus liner of rank $\|X\|\le 3$. If $\|X\|\le 2$, then $X$ is Playfair (because $X$ contains no planes, so there is nothing to verify). So, assume that $\|X\|=3$ and hence $X$ is a Proclus plane. By Corollary~\ref{c:area-bij=>3long}, the area-bijective plane $X$ is $3$-long. 
By Theorem~\ref{t:Proclus<=>}, the Proclus liner $X$ is proaffine and $3$-proregular. By Definition~\ref{d:k-regular}, the $3$-long $3$-proregular liner $X$ is $3$-regular. Since $\|X\|=3$, the $3$-regular plane $X$ is regular. By Theorem~\ref{t:procompletion-finite}, the $3$-long finite proaffine regular liner $X$ has a projective completion. By Theorem~\ref{t:area-bij=>Playfair}, the area-bijective liner $X$ is Playfair.
\end{proof}

We do not know whether Corollary~\ref{c:area-bij=>Playfair} remains true for infinite area-bijective liners. 

\begin{problem} Is every area-bijective liner Playfair?
\end{problem} 

Area-injective liners have the following important property that allows us to apply areas for establishing the parallelity of lines in hyper-Bolyai liners.

\begin{proposition}\label{p:same-area<=>subparallel} Let $a,b$ be two distinct points in an area-injective hyper-Bolyai liner $X$. For any points $x,y\in X$, the following conditions are equivalent:
\begin{enumerate}
\item[\textup{(1)}] $\Aline xy\subparallel\Aline ab$;
\item[\textup{(2)}] $axb\boxtimes ayb$;
\item[\textup{(3)}] $[axb]=[ayb]$.
\end{enumerate}
\end{proposition}

\begin{proof} The implications $(1)\Ra(2)\Ra(3)$ follow from Definitions~\ref{d:boxtimes} and \ref{d:area}. It remains to prove that $(3)\Ra(1)$. Assume that $[axb]=[ayb]$. If $[axb]=[ayb]=[0]$, then $x,y\in\Aline ab$ and hence $\Aline xy\subparallel\Aline ab$. So, assume that $[axb]=[ayb]\ne[0]$. In this case, $axb$ and $ayb$ are two triangles in the same plane  $P\defeq\overline{\{a,x,b\}}=\overline{\{a,y,b\}}$, by Corollary~\ref{c:same-area=>same-plane}. By the area-injectivity of $X$, the function $h:\Aline by\to[X^3]$, $h:z\mapsto [abz]$, is injective.  Since $X$ is hyper-Bolyai,  there exists a line $L_x$ in $X$ such that $x\in L_x$, $L_x\parallel\Aline ab$, and the lines $L_x,\Aline ay$ have a unique common point $z$. Definition~\ref{d:boxtimes} ensures that $axb\boxtimes azb$. Then $[azb]=[axb]=[ayb]$ and hence $[abz]=[aby]$, by Proposition~\ref{p:permutation-area}. By the injectivity of the function $h$, the equality $h(z)=[abz]=[aby]=h(y)$ implies $y=z$. Then $\Aline xy=\Aline xz\subseteq L_x\parallel \Aline ab$ and hence $\Aline xy\subparallel \Aline ab$.
\end{proof}

Proposition~\ref{p:same-area<=>subparallel} motivates a deeper study of area-injective liners.

\begin{theorem}\label{t:area-inj=>ass-puls} Every area-injective Playfair liner is associative-puls.
\end{theorem}

\begin{proof} Let $X$ be an area-injective Playfair liner. Given any distinct parallel lines $L,L'\subseteq X$ and points $a,b,\alpha,\beta\in L$ and $a',b',\alpha',\beta'\in L'$ with $\Aline a{a'}\parallel\Aline \alpha{\alpha'}$, $\Aline a{b'}\parallel\Aline\alpha{\beta'}$ and $\Aline {a'}{b}\parallel\Aline{\alpha'}\beta$, we should prove that $\Aline b{b'}\parallel\Aline \beta{\beta'}$. If $a=\alpha$, then $\Aline a{a'}\parallel\Aline \alpha{\alpha'}$ and  $\Aline a{b'}\parallel\Aline\alpha{\beta'}$ imply $a'=\alpha'$ and $b'=\beta'$. Also $a'=\alpha$ and $\Aline {a'}{b}\parallel\Aline{\alpha'}\beta$ imply $b=\beta$. Then $\Aline b{b'}=\Aline \beta{\beta'}$.
By analogy we can show that $\Aline b{b'}=\Aline\beta{\beta'}$ whenever $b=\beta$, $a'=\alpha'$ or $b'=\beta'$. So, assume that $a\ne\alpha$, $b\ne\beta$, $a'\ne\alpha'$ and $b'\ne \beta'$. By Proposition~\ref{p:same-area<=>subparallel}, the parallelity relation $\Aline b{b'}\parallel \Aline \beta{\beta'}$ will follow as soon as we check that $[bb'\beta']=[bb'\beta]$.

\begin{picture}(100,60)(-140,-15)

\linethickness{0.7pt}
\put(10,0){\color{teal}\line(1,0){80}}
\put(0,30){\color{teal}\line(1,0){100}}
\put(10,0){\color{cyan}\line(-1,3){10}}
\put(10,0){\color{violet}\line(1,1){30}}
\put(30,0){\color{blue}\line(-1,1){30}}
\put(30,0){\color{red}\line(1,3){10}}

\put(70,0){\color{cyan}\line(-1,3){10}}
\put(70,0){\color{violet}\line(1,1){30}}
\put(90,0){\color{blue}\line(-1,1){30}}
\put(90,0){\color{red}\line(1,3){10}}

\put(10,0){\circle*{3}}
\put(8,-9){$a$}
\put(30,0){\circle*{3}}
\put(26,-10){$b$}
\put(0,30){\circle*{3}}
\put(-3,33){$a'$}
\put(40,30){\circle*{3}}
\put(38,33){$b'$}

\put(70,0){\circle*{3}}
\put(68,-9){$\alpha$}
\put(90,0){\circle*{3}}
\put(88,-10){$\beta$}
\put(60,30){\circle*{3}}
\put(57,33){$\alpha'$}
\put(100,30){\circle*{3}}
\put(98,33){$\beta'$}
\end{picture}

Since $ab'\beta'\alpha$, $aa'\alpha'\alpha$, and $ba' \alpha'\beta$ are parallelograms, Proposition~\ref{p:12-permutations} and Definition~\ref{d:boxtimes} imply
$$[bb'\beta']=[ab'\beta']=[ab'\alpha]=[aa'\alpha]=[aa'\alpha']=[ba'\alpha']=
[b\alpha'\beta]=[bb'\beta].$$
By Proposition~\ref{p:same-area<=>subparallel}, 
the equality $[bb'\beta']=[bb'\beta]$ implies $\Aline \beta{\beta'}\subparallel \Aline b{b'}$. By Corollary~\ref{c:subparallel}, $\|\Aline \beta{\beta'}\|=2=\|\Aline b{b'}\|$ and $\Aline \beta{\beta'}\subparallel \Aline b{b'}$ imply $\Aline \beta{\beta'}\parallel \Aline b{b'}$, witnessing that the Playfair plane $X$ is associative-puls.
\end{proof}

\begin{problem} Is every area-injective Playfair plane commutative-puls?
\end{problem}

We do not known the answer to this problem. Instead, we shall prove that an area-injective Playfair plane is commutative-puls if and only if it is area-translation.

\begin{definition} A plane $X$ is defined to be \index{area-translation plane}\index{plane!area-translation}\defterm{area-translation} if two triangles $abc$ and $\alpha\beta\gamma$ have the same area whenever $\Aline b\beta\parallel\Aline ac=\Aline \alpha\gamma$, $\Aline ab\parallel \Aline \alpha\beta$ and $\Aline bc\parallel \Aline \beta\gamma$.

\begin{picture}(100,60)(-150,-15)
{\color{lightgray}
\polygon*(0,0)(30,30)(30,0)
\polygon*(50,0)(80,30)(80,0)
}
{
\linethickness{0.8pt}
\put(0,0){\color{teal}\line(1,0){80}}
\put(30,30){\color{teal}\line(1,0){50}}
\put(30,0){\color{cyan}\line(0,1){30}}
\put(80,0){\color{cyan}\line(0,1){30}}
\put(0,0){\color{blue}\line(1,1){30}}
\put(50,0){\color{blue}\line(1,1){30}}
}

\put(0,0){\circle*{3}}
\put(-3,-9){$a$}
\put(30,0){\circle*{3}}
\put(28,-9){$c$}
\put(50,0){\circle*{3}}
\put(48,-8){$\alpha$}
\put(80,0){\circle*{3}}
\put(79,-9){$\gamma$}
\put(30,30){\circle*{3}}
\put(28,33){$b$}
\put(80,30){\circle*{3}}
\put(77,33){$\beta$}

\end{picture}
\end{definition}

\begin{proposition}\label{p:commutative-puls=>area-translation} Every commutative-puls Playfair plane is area-translation.
\end{proposition}

\begin{proof} Given any triangles $abc$ and $\alpha\beta\gamma$ with $\Aline ab\parallel \Aline \alpha\beta$, $\Aline bc\parallel \Aline \beta\gamma$, and $\Aline b\beta\parallel\Aline ac=\Aline \alpha\gamma$, we need to check that $[abc]=[\alpha\beta\gamma]$. If $a=\alpha$, then $\Aline ab\parallel \Aline \alpha\beta$ implies $b=\beta$, and $\Aline bc\parallel\Aline \beta\gamma$ imply $c=\gamma$. Then $abc=\alpha\beta\gamma$ and hence $[abc]=[\alpha\beta\gamma]$. By analogy we can show that $abc=\alpha\beta\gamma$ whenever $b=\beta$ or $c=\gamma$. So, assume that $a\ne\alpha$, $b\ne\beta$ and $c\ne\gamma$.

Since $X$ is Playfair, there exists a unique point $d\in \Aline b\beta$ such that $\Aline ad\parallel \Aline bc\parallel\Aline \beta\gamma$.

\begin{picture}(100,60)(-150,-12)

{\color{lightgray}
\polygon*(1.5,0.5)(29.5,28.5)(29.5,0.5)
\polygon*(51.5,0.5)(79.5,28.5)(79.5,0.5)
}
{
\linethickness{0.8pt}
\put(0,0){\color{teal}\line(1,0){80}}
\put(0,30){\color{teal}\line(1,0){80}}
\put(0,0){\color{cyan}\line(0,1){30}}
\put(30,0){\color{cyan}\line(0,1){30}}
\put(80,0){\color{cyan}\line(0,1){30}}
\put(0,0){\color{blue}\line(1,1){30}}
\put(50,0){\color{blue}\line(1,1){30}}
}
\put(0,30){\color{gray}\line(1,-1){30}}
\put(0,30){\color{red}\line(5,-3){50}}
\put(30,30){\color{gray}\line(2,-3){20}}
\put(30,30){\color{red}\line(5,-3){50}}

\put(0,0){\circle*{3}}
\put(-3,-9){$a$}
\put(30,0){\circle*{3}}
\put(27,-9){$c$}
\put(50,0){\circle*{3}}
\put(48,-8){$\alpha$}
\put(80,0){\circle*{3}}
\put(79,-9){$\gamma$}
\put(0,30){\circle*{3}}
\put(-3,33){$d$}
\put(30,30){\circle*{3}}
\put(28,33){$b$}
\put(80,30){\circle*{3}}
\put(77,33){$\beta$}

\end{picture}

Since $X$ is commutative-puls, the parallelity relations $\Aline ab\parallel\Aline\alpha\beta$ and $\Aline ad\parallel\Aline\beta\gamma$ imply $\Aline d\alpha\parallel \Aline b\gamma$. Then $adbc$ and $\alpha db\gamma$ are parallelograms, and Proposition~\ref{p:12-permutations} implies $[abc]=[bcd]$ and  $[\alpha b\gamma]=[b\alpha d]$. By Definition~\ref{d:boxtimes}, $[abc]=[bcd]=[b\alpha d]=[\alpha b\gamma]=[\alpha\beta\gamma]$.
\end{proof}

\begin{proposition} An area-injective Playfair plane is commutative-puls if and only if it is area-translation.
\end{proposition} 

\begin{proof} The ``only if'' part follows from Proposition~\ref{p:commutative-puls=>area-translation}. To prove the ``if'' part, assume that $X$ is an area-translation area-injective Playfair plane. Given two disjoint lines $L,L'\subseteq X$ and distinct points $a,b,c\in L$ and $a',b',c'\in L'$ with $\Aline a{b'}\parallel \Aline {a'}b$ and $\Aline b{c'}\parallel \Aline {b'}c$, we need to check that $\Aline a{c'}\parallel \Aline {a'}c$. Since $X$ is Playfair, there exists a unique point $d\in L$ such that the line $\Aline {b'}{d}$ is parallel to the line $\Aline a{c'}$. 

\begin{picture}(100,60)(-150,-12)
{\color{lightgray}
\polygon*(0.5,0.5)(0.5,28.5)(28.5,0.5)
\polygon*(50.5,0.5)(50.5,28.5)(78.5,0.5)
}
{
\linethickness{0.8pt}
\put(0,0){\color{teal}\line(1,0){80}}
\put(0,30){\color{teal}\line(1,0){80}}
\put(0,0){\color{red}\line(0,1){30}}
\put(50,0){\color{red}\line(0,1){30}}
\put(80,0){\color{red}\line(0,1){30}}
\put(0,0){\color{blue}\line(5,3){50}}
\put(30,0){\color{blue}\line(5,3){50}}
\put(0,30){\color{violet}\line(1,-1){30}}
\put(50,30){\color{violet}\line(1,-1){30}}
}
\put(30,0){\color{gray}\line(2,3){20}}
\put(50,0){\color{gray}\line(1,1){30}}

\put(0,0){\circle*{3}}
\put(-3,-9){$a$}
\put(30,0){\circle*{3}}
\put(26,-9){$b$}
\put(50,0){\circle*{3}}
\put(48,-9){$d$}
\put(80,0){\circle*{3}}
\put(78,-8){$c$}
\put(0,30){\circle*{3}}
\put(-3,33){$c'$}
\put(50,30){\circle*{3}}
\put(48,33){$b'$}
\put(80,30){\circle*{3}}
\put(77,33){$a'$}

\end{picture}

Since $X$ is area-translation, the parallelity relations $\Aline a{c'}\parallel \Aline d{b'}$ and $\Aline {c'}b\parallel \Aline{b'}c$ imply $[ac'b]=[db'c]$.
Since $ab'a'b$ is a parallelogram, Proposition~\ref{p:12-permutations} and Definition~\ref{d:boxtimes} ensure that 
$$[db'c]=[ac'b]=[ab'b]=[bb'a']=[db'a'].$$ Applying Proposition~\ref{p:same-area<=>subparallel}, we conclude that $\Aline c{a'}\parallel \Aline d{b'}\parallel \Aline a{c'}$, witnessing that $X$ is commutative-puls.
\end{proof}

We do not know whether area-injective Playfair planes have any reasonable plus-properties.

\begin{problem} Is every area-injective (and area-translation) Playfair plane invertible-plus?
\end{problem}

\section{Area-preserving and equiaffine automorphisms of Playfair planes}

Let $A:X\to X$ be an automorphism of a plane $X$. 
Proposition~\ref{p:permutation-area} implies that for any area $\alpha\in [X^3]$ in a plane $X$, the set $$A\alpha\defeq\{Aabc:abc\in\alpha\}$$ is a well-defined area, equal to the area $[Aabc]$ for any triple $abc\in\alpha$.
Moreover, the function $A{\circ}:[X^3]\to [X^3]$, $A{\circ}:\alpha\mapsto A\alpha$, is a well-defined bijection of the areol $[X^3]$ of $X$. 

\begin{exercise} Show that the image of the zero area $[0]$ in any plane $X$ under any automorphism of $X$ remains $[0]$.
\end{exercise}

\begin{definition} An automorphism $A$ of a plane $X$ is called \index{automorphism!area-preserving}\index{automorphism!area-preserving}\defterm{area-preserving} if $A\alpha=\alpha$ for every area $\alpha\in[X^3]$. The set of all area-preserving automorphisms of a plane $X$ is denoted by $\SAut(X)$ and called the \index{special automorphism group}\defterm{special automorphism group} of the plane $X$. 
\end{definition}

\begin{proposition} For any plane $X$, the set $\SAut(X)$ of all area-perserving automorphisms of $X$ is a normal subgroup of the automorphism group $\Aut(X)$ of $X$.
\end{proposition}

\begin{proof} It is clear that $\SAut(X)$ is a subgroup of the automorphism group $\Aut(X)$ of $X$. To see that the subgroup $\SAut(X)$ is norma in $\Aut(X)$, take any area-preserving automorphism $S\in\SAut(X)$ and any automorphism $A\in\Aut(X)$. We have to show that $A^{-1}SA\alpha=\alpha$ for all areas $\alpha\in [X^3]$. Taking into account that  $A\alpha$ is an area in $X$ and $S$ is area-preserving, we conclude that $SA\alpha=A\alpha$ and hence $A^{-1}SA\alpha=A^{-1}A\alpha=\alpha$.
\end{proof}

Let us recall that an automorphism $A$ of a plane $X$ is called a \defterm{hypershear} if there exists a line $L$ in $X$ such that $A(x)=x$ for all $x\in L$ and $\Aline xy\subparallel L$ for all $xy\in A$.

\begin{theorem}\label{t:shear=>SAut} Every hypershear of a Playfair plane is area-preserving.
\end{theorem}

\begin{proof} Let $S:X\to X$ be a hypershear of a Playfair plane $X$. We lose no generality assuming that $S$ is not identity. Then the set $\Fix(S)\defeq\{x\in X:S(x)=x\}$ is a line in $X$, and for any pair $xy\in S$, the flat $\Aline xy$ is subparallel to the line $\Fix(S)$, by Definition~\ref{d:hypershear-hyperscale} of a hypershear. 

\begin{claim}\label{cl:shear=area} For any triangle $abc\in X^3$ with $c\in\Fix(S)$, the triangle $\alpha\beta c\defeq S abc$ has area $[\alpha\beta c]=[abc]$.
\end{claim}

\begin{proof} Since $abc$ is a triangle with $c\in\Fix(S)$, at least one vertex of $abc$ does not belong to the line $\Fix(S)$. We lose no generality assuming that $a\notin \Fix(S)$. Let $L_b$ be a unique line in the Playfair plane $X$ such that $b\in L_b$ and $L_b\parallel \Aline ac$. Since the lines $\Aline ac$ and $\Fix(S)$ are concurrent, there exists a unique point $o\in L_b\cap\Fix(S)$. 

\begin{picture}(100,90)(-130,-15)
{\color{lightgray}
\polygon*(20,0)(20,60)(40,40)
\polygon*(20,0)(80,60)(80,40)
}

{\linethickness{0.7pt}
\put(0,0){\color{teal}\line(1,0){100}}
}
\put(104,-3){\color{teal}$\Fix(S)$}
\put(20,60){\color{teal}\line(1,0){60}}
\put(40,40){\color{teal}\line(1,0){40}}
\put(20,0){\line(0,1){60}}
\put(40,0){\line(0,1){40}}
\put(20,60){\line(1,-1){20}}
\put(20,0){\line(1,1){60}}
\put(80,60){\line(0,-1){20}}
\put(20,0){\line(3,2){60}}
\put(20,0){\line(1,2){20}}
\put(40,0){\line(1,1){40}}
\put(40,0){\line(-1,3){20}}
\put(40,0){\line(2,3){40}}

\put(20,0){\circle*{2.5}}
\put(17,-8){$c$}
\put(20,60){\circle*{2.5}}
\put(13,61){$a$}
\put(40,40){\circle*{2.5}}
\put(40,43){$b$}
\put(80,40){\circle*{2.5}}
\put(82,36){$\beta$}
\put(80,60){\circle*{2.5}}
\put(82,61){$\alpha$}
\put(40,0){\circle*{2.5}}
\put(37,-8){$o$}
\end{picture}

It follows from $\Aline bo\cap\Aline ac=\varnothing$ that $\Aline \beta o\cap \Aline \alpha c=S[\Aline bo]\cap S[\Aline ac]=S[\Aline bo\cap \Aline ac]=\varnothing$. Then $abc\boxtimes aoc\boxtimes \alpha oc\boxtimes \alpha\beta c$ and hence $[Sabc]=[\alpha\beta c]=[abc]$.
\end{proof}

Now we can prove that for every triangle $abc$ in $X$, its image $\alpha\beta\gamma\defeq Sabc$ has area $[\alpha\beta\gamma]=[abc]$. If $\{\alpha,\beta,\gamma\}\cap\Fix(S)\ne\varnothing$, then the equality $[\alpha\beta\gamma]=[abc]$ follows from Claim~\ref{cl:shear=area}. So, assume that $\{\alpha,\beta,\gamma\}\cap\Fix(S)=\varnothing$. Since at least two sides of the triangle $abc$ are not parallel to the line $\Fix(S)$, we lose no generality assuming that the line $\Aline ac$ is concurrent with the line $\Fix(S)$. Let $L_b$ be a unique line in the Playfair plane $X$ such that $b\in L_b$ and $L_b\parallel\Aline ac$. Let $o$ be the unique common point of the lines $L_b$ and $\Fix(S)$. 

\begin{picture}(100,110)(-130,-15)
{\color{lightgray}
\polygon*(0,80)(20,60)(20,40)
\polygon*(40,40)(40,80)(50,60)
}

{\linethickness{0.7pt}
\put(0,0){\color{teal}\line(1,0){100}}
}
\put(104,-3){\color{teal}$\Fix(S)$}

\put(20,40){\line(0,1){20}}
\put(0,80){\line(1,-1){20}}
\put(0,80){\line(1,-2){20}}
\put(40,40){\line(0,1){40}}
\put(40,80){\line(1,-2){10}}
\put(40,40){\line(1,2){10}}
\put(50,0){\line(-1,2){30}}
\put(50,0){\line(0,1){60}}
\put(50,0){\line(-3,4){30}}
\put(50,0){\line(-5,8){50}}
\put(50,0){\line(-1,4){10}}
\put(50,0){\line(-1,8){10}}
\put(20,40){\color{teal}\line(1,0){20}}
\put(20,60){\color{teal}\line(1,0){30}}
\put(0,80){\color{teal}\line(1,0){40}}

\put(0,80){\circle*{2.5}}
\put(-2,83){$a$}
\put(20,60){\circle*{2.5}}
\put(19,62){$b$}
\put(20,40){\circle*{2.5}}
\put(16,33){$c$}
\put(40,80){\circle*{2.5}}
\put(38,83){$\alpha$}
\put(50,60){\circle*{2.5}}
\put(52,57){$\beta$}
\put(40,40){\circle*{2.5}}
\put(34,33){$\gamma$}
\put(50,0){\circle*{2.5}}
\put(47,-8){$o$}

\end{picture}

It follows from $\Aline bo\cap\Aline ac=\varnothing$ that  $\Aline \beta o\cap\Aline\alpha\gamma=S[\Aline bo]\cap S[\Aline ac]\subseteq S[L_b\cap \Aline ac]=\varnothing$. Then $abc\boxtimes aoc$ and $\alpha\beta \gamma\boxtimes \alpha o\gamma$. Since $o\in \Fix(S)$, we can apply Claim~\ref{cl:shear=area} and conclude that $[\alpha\beta\gamma]=[\alpha o\gamma]=[Saoc]=[aoc]=[abc]$.
\end{proof}

\begin{problem} Is every translation of a Playfair plane area-preserving?
\end{problem}

We shall answer this problem affirmatively under an additional assumption that the Playfair plane is area-translation.

\begin{proposition}\label{p:area-translation=>Trans=SAut} Every translation of an area-translation Playfair plane is area-preserving.
\end{proposition}

\begin{proof} Let $T:X\to X$ be a translation of an area-translation Playfair plane. Given any triple $abc\in X$, we should prove that the triple $\alpha\beta\gamma\defeq Tabc$ has area $[\alpha\beta\gamma]=[abc]$. If $[abc]=[0]$, then the points $a,b,c$ are collinear and so are the points $\alpha,\beta,\gamma$, which implies $[\alpha\beta\gamma]=[0]=[abc]$. So, assume that $[abc]\ne[0]$ and hence $abc$ and $\alpha\beta\gamma$ are triangles. We lose no generality assuming that the translation $T$ is nonidentity and hence $\Aline a{\alpha}$, $\Aline b{\beta}$, $\Aline c{\gamma}$ are parallel lines in the direction ${\boldsymbol\delta}\defeq\{\Aline xy:xy\in T\}$ of the translation $T$. Since $abc$ is a triangle, at most one its sides belongs to the direction $\boldsymbol \delta$. So, we lose no generality assuming that $\Aline ac \notin\boldsymbol\delta$. Since $X$ is Playfair, there exists a unique line $L_b$ in $X$ such that $b\in L_b\parallel \Aline ac$. 

\begin{picture}(100,100)(-130,-15)
{\color{lightgray}
\polygon*(0,0)(0,60)(30,0)
\polygon*(60,0)(60,60)(90,0)
}

\linethickness{0.7pt}
\put(0,0){\color{teal}\line(1,0){90}}
\put(30,30){\color{teal}\line(1,0){60}}
\put(0,60){\color{teal}\line(1,0){90}}
\put(0,0){\color{cyan}\line(0,1){60}}
\put(0,0){\color{blue}\line(1,1){30}}
\put(0,60){\color{violet}\line(1,-1){30}}
\put(30,0){\color{cyan}\line(0,1){60}}
\put(30,0){\color{red}\line(-1,2){30}}

\put(60,0){\color{cyan}\line(0,1){60}}
\put(60,0){\color{blue}\line(1,1){30}}
\put(60,60){\color{violet}\line(1,-1){30}}
\put(90,0){\color{cyan}\line(0,1){60}}
\put(90,0){\color{red}\line(-1,2){30}}

\put(0,0){\circle*{3}}
\put(-2,-9){$c$}
\put(30,0){\circle*{3}}
\put(28,-9){$c'$}
\put(60,0){\circle*{3}}
\put(58,-9){$\gamma$}
\put(90,0){\circle*{3}}
\put(88,-8){$\gamma'$}
\put(0,60){\circle*{3}}
\put(-2,63){$a$}
\put(30,60){\circle*{3}}
\put(28,63){$a'$}
\put(60,60){\circle*{3}}
\put(58,63){$\alpha$}
\put(90,60){\circle*{3}}
\put(88,63){$\alpha'$}
\put(30,30){\circle*{3}}
\put(32,32){$b$}
\put(90,30){\circle*{3}}
\put(92,28){$\beta$}
\end{picture}

Since $L_b\parallel \Aline ac$ and the line $\Aline ac$ is concurrent with the lines $\Aline a{\alpha},\Aline c{\gamma}\in\boldsymbol\delta$, there exist unique points $a'\in L_b\cap\Aline a\alpha$ and $c'\in L_b\cap\Aline c\gamma$. Consider the points $\alpha'\defeq T(a')$ and $\gamma'\defeq T(c')$ and observe that $b\in \Aline {a'}{c'}=L_b\parallel \Aline ac$ implies $\beta\in \Aline{\alpha'}{\beta'}\parallel \Aline \alpha\gamma$. Since the translation $T$ is a dilation, the line $\Aline \alpha\gamma=T[\Aline ac]$ is parallel to the line $\Aline ac$, and the line $\Aline \alpha{\gamma'}=T[\Aline a{c'}]$ is parallel to the line $\Aline a{c'}$. Since $X$ is area-translation, $[cac']=[\gamma\alpha\gamma']$. Since $\Aline b{c'}\subparallel \Aline ac$, Definition~\ref{d:boxtimes} and Proposition~\ref{p:12-permutations} ensure that $cac'\boxtimes cab$ and hence $[cac']=[cab]=[abc]$.
By analogy, $[\gamma\alpha\gamma']=[\alpha\beta\gamma]$. Therefore, $[abc]=[cac']=[\gamma\alpha\gamma']=[\alpha\beta\gamma]$, witnessing that the translation $T$ is an area-preserving automorphism of $X$.
\end{proof}

\begin{corollary}\label{c:commutative-puls=>Trans=SAut} Every transation of a commutative-puls Playfair plane is area-preserving.
\end{corollary} 

\begin{proof} Let $T:X\to X$ be a translation of a commutative-puls Playfair plane $X$. By Proposition~\ref{p:commutative-puls=>area-translation}, the commutative-puls Playfair plane $X$ is area-translation. By Proposition~\ref{p:area-translation=>Trans=SAut}, the translation $T$ is area-preserving.
\end{proof}

\begin{corollary}\label{c:Thales-trans=>SAut} Every translation of a Thalesian Playfair plane is area-preserving.
\end{corollary} 

\begin{proof} Let $T:X\to X$ be a translation of a Thalesian Playfair plane $X$. By Corollary~\ref{c:Thales=>commutative-puls}, the Thalesian Playfair plane $X$ is commutative-puls. By
Corollary~\ref{c:commutative-puls=>Trans=SAut}, the translation $T$ is area-preserving.
\end{proof}

\begin{definition} An automorphism $A$ of a Playfair plane $X$ is called \index{equiaffine automorphism}\index{automorphism!equiaffine}\defterm{equiaffine} if it can be written as composition of finitely many hypershears. Let $\EAff(X)$ denote the set of all equiaffine automorphisms of the Playfair plane $X$.
\end{definition}

It is easy to see that for every Playfair plane, $\EAff(X)$ is a normal subgroup of the autmorphism group $\Aut(X)$ of the plane $X$. Theorem~\ref{t:shear=>SAut} guarantees that $\EAff(X)\subseteq\SAut(X)$. By Proposition~\ref{p:aff-paracentral<=>} and \ref{p:affinity=>portionality=>scalarity}, hypershears are affinities, which implies that $\EAff(X)$ is a subgroup of the group $\Aff(X)$ of all affinities of the plane $X$. Let $\SAff(X)\defeq\SAut(X)\cap\Aff(X)$ be the group of area-preserving affinities of the plane $X$. Then $\EAff(X)\subseteq \SAff(X)=\SAut(X)\cap\Aff(X)$.

\begin{remark}\label{r:EAff} The following table (calculated by Ivan Hetman), for each Playfair plane $X$ of order 9, displays the cardinalities of the whole automorphism group $\Aut(X)$ of $X$, the group $\SAut(X)$ of area-preserving automorphisms of $X$, the group $\EAff(X)$ of all equiaffine automorphisms of $X$, and the group $\Trans(X)$ of all translations of $X$. 

{\renewcommand{\arraystretch}{1.3}
$$
\begin{array}{|c|c|c|c|c|c|c|c|}
\hline
\mbox{Playfair plane $\Pi$:}&\mbox{\tt Desarg}&\mbox{\tt Thales}&\mbox{\tt Hall}&\mbox{\tt hall}&\mbox{\tt dhall}&\mbox{\tt Hughes}&{\tt hughes}\\
\hline
|\Aut(\Pi)|&933120&311040&311040&3840&3456&2592&432\\
|\SAut(\Pi)|&58320&311040&311040&3840&3456&2592&432\\
|\EAff(\Pi)|&58320&1&81&1&9&9&9\\
|\Trans(\Pi)|&81&81&9&1&1&9&1\\
\hline
\end{array}
$$
}
\end{remark}

\begin{remark} By (the proof of) Theorem~\ref{t:shear=>translation}, every translation in a shear Playfair plane can be written as the composition of two hypershears, which implies that $\Trans(X)\subseteq\EAff(X)$ for every Moufang Playfair plane $X$. The example of the Thalesian plane {\tt Thales} shows that the inclusion $\Trans(X)\subseteq\EAff(X)$ does not hold in Thalesian (non-Moufang) Playfair planes. 
\end{remark}


\section{Equiaffine automorphisms of Moufang Playfair planes}

In this section we study the group of equiaffine automorphisms of Moufang Playfair planes.

\begin{lemma}\label{l:boxtimes=>hypershear} If two triangles $a_1a_2a_3$ and $b_1b_2b_3$ in a Moufang plane $X$ have the same base and height, then $b_1b_2b_3=Sa_1a_2a_3$ for some hypershear $S:X\to X$ of $X$.
\end{lemma}

\begin{proof} Since the triangles $a_1a_2a_3$ and $b_1b_2b_3$ have the same base and height, there exist distinct numbers $i,j,k\in\{1,2,3\}$ such that $a_ia_j=b_ib_j$ and $\Aline{a_k}{b_k}\subparallel\Aline {a_i}{a_j}=\Aline{b_i}{b_j}$. By Theorem~\ref{t:affine-Moufang<=>}, Moufang Playfair planes are shear planes. So, there exists a hypershear $S:X\to X$ such that $\{a_i,a_j\}=\{b_i,b_j\}\subseteq\Fix(S)$ and $S(a_k)=b_k$. Then $S$ is a required hypershear of $X$ with $Sa_1a_2a_3=b_1b_2b_3$.
\end{proof}

\begin{theorem}\label{t:[abc]=[xyz]<=>abc=Exyz} Two triangles $abc$ and $xyz$ in a Moufang Playfair plane $X$ have the same area if and only if $abc=Exyz$ for some equiaffine automorphism $E$ of $X$.
\end{theorem}

\begin{proof} Fix any triangles $abc$ and $xyz$ in a Moufang Playfair plane.
\smallskip

To prove the ``if'' part, assume that $abc=Exyz$ for some equiaffine automorphism $E\in \EAff(X)$. Since $E$ is a composition of finitely many hypershears, Theorem~\ref{t:shear=>SAut} implies that $E$ is area-preserving and hence $[abc]=[Exyz]=[xyz]$.
\smallskip

To prove the ``only if'' part, assume that $[abc]=[xyz]$. Then there exists a sequence of triangles $\tau_0,\tau_1,\dots,\tau_n$ in $X$ such that $xyz=\tau_0\boxtimes\tau_1\boxtimes\dots\boxtimes \tau_n=abc$. By Lemma~\ref{l:boxtimes=>hypershear}, for every positive $i<n$, there exists a hypershear $S_i$ of $X$ such that $S_i\tau_{i-1}=\tau_i$. Let $E_0$ be the identity automorphism of $X$ and for every positive number $i\le n$, let $E_i\defeq S_i\circ E_{i-1}$. By induction, we shall prove that $E_ixyz=\tau_{i}$ for all $i\le n$. For $i=0$, the equality $xyz=\tau_0$ implies $E_0xyz=\tau_0$. Assume that $E_ixyz=\tau_i$ for some $i<n$. Then $E_{i+1}xyz=S_{i+1}E_ixyz=S_{i+1}\tau_i=\tau_{i+1}$. Therefore, $E_ixyz=\tau_i$ for all $i\le n$. In particlar, $E_nxyz=\tau_n=abc$. Then the equiaffine automorphism $E\defeq E_n$ has the required property: $Exyz=abc$.
\end{proof}  

We shall deduce from Theorems~\ref{t:area-move} and 
\ref{t:[abc]=[xyz]<=>abc=Exyz} that the group $\EAff(X)$ of Moufang Playfair plane is $\bar 2$-transitive. In the following definition, by a \index{doubleton}\defterm{doubleton} in a liner we understand any ordered pair $xy$ of distinct points of the liner.

\begin{definition}\label{d:EAff-2transitive} A group $G$ of automorphisms of a liner $X$ is called 
\begin{itemize}
\item \index{automorphism group!$1$-transitive}\defterm{$1$-transitive} if for any points $x,y$ in $X$ there exists an automorphism $\Phi\in G$ such that $\Phi x=y$;
\item \index{automorphism group!$2$-transitive}\defterm{$2$-transitive} if for any doubletons $ab$ and $xy$ in $X$, there exists an automorphism $\Phi\in G$ such that $\Phi ab=xy$;
\item \index{automorphism group!$3$-transitive}\defterm{$3$-transitive} if for any triangles $abc$ and $xyz$ in $X$, there exists an automorphism $\Phi\in G$ such that $\Phi abc=xyz$;
\item \index{automorphism group!flag-transitive}\defterm{flag-transitive} if for any doubletons $ab$ and $xy$ in $X$,  there exists an automorphism $\Phi\in G$ such that $\Phi a=x$ and $\Phi[\Aline ab]=\Aline xy$;
\item \index{automorphism group!$\bar 1$-transitive}\defterm{$\bar 1$-transitive} if for any triangles $abc$ and $xyz$ in $X$, there exists an automorphism $\Phi\in G$ such that $\Phi a=x$ and $\Phi[\Aline bc]=\Aline yz$;
\item \index{automorphism group!$\bar 2$-transitive}\defterm{$\bar 2$-transitive} if for any triangles $abc$ and $xyz$ in $X$, there exists an automorphism $\Phi\in G$ such that $\Phi ab=xy$ and $\Phi[\Aline bc]=\Aline yz$.
\end{itemize}
\end{definition}

For any group of liner automorphisms those notions relate as follows.
$$\xymatrix{
\mbox{$3$-transitive}\ar@{=>}[r]&\mbox{$\bar 2$-transitive}\ar@{=>}[r]\ar@{=>}[d]&\mbox{$2$-transitive}\ar@{=>}[r]&\mbox{flag-transitive}\ar@{=>}[d]\\
&\mbox{$\bar 1$-transitive}\ar@{=>}[rr]&&\mbox{$1$-transitive}
}
$$

\begin{corollary}\label{c:EAff2} The group $\EAff(X)$ of equiaffine automorphisms of any Moufang Playfair plane $X$ is $\bar 2$-transitive.
\end{corollary}

\begin{proof} Given any triangles $abc$ and $xyz$ in $X$, we should find an equiaffine automorphism $E$ of $X$ such that $Exy=ab$ and $E[\Aline yz]=\Aline bc$. By Theorem~\ref{t:area-move}, there exists a point $d\in \Aline bc$ such that $[abd]=[xyz]$.
By Theorem~\ref{t:[abc]=[xyz]<=>abc=Exyz}, there exists an equiaffine automorphism $E$ of $X$ such that $Exyz=abd$. Since $xyz$ is a triangle, $y\ne z$ and hence $b=E(y)\ne E(z)=d$. Then $E[\Aline yz]=\Aline bd=\Aline bc$.
\end{proof}

Corollary~\ref{c:EAff2} is completed by the following non-trivial result. 
 
\begin{proposition}\label{p:Moufang<=>2transitive} For a finite affine plane $X$, the following conditions are equivalent:
\begin{enumerate}
\item[\textup{(1)}] $X$ is Moufang;
\item[\textup{(2)}] the group $\EAff(X)$ is $\bar 2$-transitive;
\item[\textup{(3)}] the group $\EAff(X)$ is $2$-transitive.
\end{enumerate}
\end{proposition} 

\begin{proof} The implication $(1)\Ra(2)$ follows from Corollary~\ref{c:EAff2}, and $(2)\Ra(3)$ is trivial. To prove $(3)\Ra(1)$, assume that the group $\EAff(X)$ is $2$-transitive. Then the automorhism group of $X$ is $2$-transitive and hence the affine plane $X$ is $2$-homogeneous. By Jha--Johnson Theorem~\ref{t:Jha-Johnson1998}, the $2$-homogeneous plane $\Pi$ is Thalesian. Let us show that the plane $X$ is distinct from the exceptional non-Desarguesian $2$-homogeneuos affine planes of order $9$ and $27$. By Remark~\ref{r:EAff}, the group of equiaffine automorphisms of the Thalesian affine plane of order $9$ is trivial and hence cannot act $2$-transitively on this plane. Computer calculations (made by Ivan Hetman) show that the group of equiaffine automorphisms of the Hering plane of order 27 is trivial and cannot act $2$-transitively on this plane. Therefore, the plane $X$ is neither Thalesian of order 9 not a Hering of order 27. By Jha--Johnson Theorem~\ref{t:Jha-Johnson1998}, the plane $X$ is Desarguesian and hence Moufang.
\end{proof}

We do not know if Proposition~\ref{p:Moufang<=>2transitive} generalizes to infinite planes.

\begin{problem} Is a Playfair plane Moufang if its group of equiaffine automorphisms is $\bar 2$-transitive?
\end{problem}

The area of a triangle in a Moufang based affine plane can be calculated by the following formula.

\begin{proposition}\label{p:area-via-coordinates} Let $uow$ be an affine base in a Moufang Playfair plane $X$ and $\Delta$ be the ternar of the based affine plane $(X,uow)$. For any points $a,b\in X$ with coordinates $(a_1,a_2), (b_1,b_2)\in\Delta^2$ in the affine base $uow$, the triple $aob$ has  area $[aob]=[oud]$, where 
$$d=\begin{cases}(b_2{\cdot}(a_2^{-1}{\cdot}a_1)){\cdot}a_2-b_1{\cdot}a_2&\mbox{if $a_2\ne o$};\\
b_2{\cdot}a_1&\mbox{if $a_2=o$}.
\end{cases}
$$
\end{proposition}

\begin{proof} By Theorem~\ref{t:affine-Moufang<=>}, the ternar $\Delta$ of the Moufang Playfair plane $X$ is linear, distributive, associative-plus and alternative-dot. Identifying each point $p\in X$ with its coordinates $(p_1,p_2)$ in the affine base $uow$, we can identify the Playfair plane $X$ with the coordinate plane $\Delta\times\Delta$ of the ternar $\Delta$. 
 
First assume that $a_2\ne o$. Consider the elements $d\defeq (b_2{\cdot}(a_2^{-1}{\cdot}a_1)){\cdot}a_2+b_1{\cdot}(-a_2)$ and $s\defeq a_2^{-1}\cdot(1-a_1)$ of the ternar $\Delta$. 
By Proposition~\ref{p:horizontal-shear}, the map $H:X\to X$, $H:(x,y)\mapsto(x+y{\cdot} s,y)$, is a horizontal shear of the based affine plane $(X,uow)$.  The shear $H$ maps the point $a$ with coordintes $(a_1,a_2)$ to the point $a'$ with the coordinates $(a_1+a_2{\cdot} s,a_2)=(a_1+a_2{\cdot} (a_2^{-1}{\cdot} (1-a_1)),a_2)=(1,a_2)$. By Proposition~\ref{p:(x,y)->(x,xc+y)}, the map $V:X\to X$, $V:(x,y)\mapsto(x,y+x{\cdot}(-a_2))$, is a vertical shear of the based affine plane $(X,uow)$.
The vertical shear $V$ maps the point $a'$ to the unit $u$ of the affine base $uow$. Therefore, the equiaffine automorphism $VH$ maps the point $a$  to the point $u$. This automorphism maps the point $b$ to the point $b'$ with coordinates $$VH(b_1,b_2)=V(b_1+b_2{\cdot} s,b_2)=(b_1+b_2{\cdot} s,b_2+(b_1+b_2{\cdot}s){\cdot}(-a_2))=(b_1+b_2{\cdot} s,d).$$ If $d=o$, then the point $b'$ belongs to the line $\Aline ou$ and hence $[oab]=[VHoab]=[oub']=[0]=[oud]$. If $d\ne o$, then there exists a horizontal shear $S:X\to X$ such that $S(b')=d\in \Delta$. Then $E\defeq SVH$ is an equiaffine automorphism of the plane $X$ such that $Eoab=oud$, witnessing that $[oab]=[oud]$.

Next, assume that $a_2=0$. If $a_1=o$, too, then $a=o$, $d=b_2a_1=o$ and $[oab]=[0]=[oud]$. So, assume that $a_1\ne o$. Consider the vertical shear $S: X\to X$, $S:(x,y)\mapsto (x,x+y)$, of the plane $X$ and observe that the point $a'=S(a)$ has coordinates $(a_1,a_1)$ and the point $b'=S(b)$ has coordinates $(b_1,b_1+b_2)$. By the preceding case, $[oab]=[oa'b']=[oud]$, where $d=((b_1+b_2)(a_1^{-1}a_1))a_1-b_1a_1=b_2a_1$.
\end{proof}

\begin{corollary} Let $uow$ be an affine base in a Moufang Playfair plane $X$ and $\Delta$ be the ternar of the based affine plane $(X,uow)$. For any points $a,b,c\in X$, the triple $abc$ has area $[abc]=[cab]=[oud]$ and absolute area $[\![abc]\!]=[\![oud]\!]$, where 
$$
d=\begin{cases}\big((b_2-c_2){\cdot}((a_2-c_2)^{-1}{\cdot}(a_1-c_1))\big){\cdot}(a_2-c_2)-(b_1-c_1){\cdot}(a_2-c_2)&\mbox{if $a_2\ne c_2$};\\
(b_2-c_2){\cdot}(a_1-c_1)&\mbox{if $a_2=c_2$}
\end{cases}
$$
and $(a_1,a_2), (b_1,b_2),(c_1,c_2)\in\Delta^2$ are the coordinates of the points $a,b,c$ in the affine base $uow$. 
\end{corollary}

\section{Equiaffine automorphisms of Desarguesian Playfair planes}

In this section we study equiaffine automorphisms of Desarguesian affine planes. In particular, we prove that for any Desarguesian Playfair plane $X$, the quotient group $\Aff(X)/\EAff(X)$ is isomorphic to the abelianization $\IR^*_X/[\IR^*_X,\IR^*_X]$ of the multiplicative group $\IR_X^*=\IR_X\setminus\{0\}$ of the scalar corps $\IR_X$ of the plane $X$. Here \index[note]{$[\mathbb R^*_X,\mathbb R^*_X]$}
$[\IR^*_X,\IR^*_X]$ is the \index{commutator subgroup}\defterm{commutator subgroup} of the multiplicative group $\IR^*_X$. The subgroup $[\IR^*_X,\IR^*_X]$ is  generated by commutators $xyx^{-1}y^{-1}$ of non-zero scalars $x,y\in\IR_X$. 

For every scalar $s\in\IR_X$, let $[s]\defeq s\cdot[\IR^*_X,\IR^*_X]=[\IR^*_X,\IR^*_X]\cdot s$ be the coset of $s$ by the commutator subgroup $[\IR^*_X,\IR^*_X]$ in $\IR_X$. The set $[\IR_X]\defeq\{[s]:s\in\IR_X\}$ is a commutative semigroup with respect to the operation of multiplication $[s]\cdot [t]\defeq[s\cdot t]=\{x\cdot y:x\in[s]\;\wedge\;y\in[t]\}$. The semigroup $[\IR_X]$ contains zero and unit elements $[0]=\{0\}$ and $[1]=[\IR^*_X,\IR^*_X]$. Elements of the commutative semigroup \index[note]{$[\mathbb R_X]$}$[\IR_X]$ will be called \index{coscalar}\defterm{coscalars}\footnote{Coined by removing ``mmutative~'' from {\em co}mmutative {\em scalars}.}.

The semigroup $[\IR_X]$ is the abelianization of the multiplicative semigroup $(\IR_X,\cdot)$, and the group $[\IR^*_X]\defeq[\IR_X]\setminus\{[0]\}$ is the abelianization $\IR^*_X/[\IR^*_X,\IR^*_X]$ of the multiplicative group $(\IR^*_X,\cdot)$. 

An isomorphism between the groups $\Aff(X)/\EAff(X)$ and $[\IR^*_X]$ will be constructed with the help of Dieudonn\'e\index[person]{Dieudonn\'e}\footnote{{\bf Jean Alexandre Eug\`ene Dieudonn\'e} (1906 -- 1992) was a French mathematician, notable for research in abstract algebra, algebraic geometry, and functional analysis, for close involvement with the Nicolas Bourbaki pseudonymous group and the ``\'El\'ements de g\'eom\'etrie alg\'ebrique'' project of Alexander Grothendieck, and as a historian of mathematics, particularly in the fields of functional analysis and algebraic topology. His work on the classical groups (the book ``La G\'eom\'etrie des groupes classiques'' published in 1955), and on formal groups, introducing what now are called Dieudonn\'e modules, had a major effect on those fields.

He was born and brought up in Lille, with a formative stay in England where he was introduced to algebra. In 1924 he was admitted to the \'Ecole Normale Sup\'erieure, where Andr\'e Weil was a classmate. He began working in complex analysis. In 1934 he was one of the group of normaliens convened by Weil, which would become ``Bourbaki''. He served in the French Army during World War II, and then taught in Clermont-Ferrand until the liberation of France. After holding professorships at the University of S\~ao Paulo (1946–47), the University of Nancy (1948–1952) and the University of Michigan (1952–53), he joined the Department of Mathematics at Northwestern University in 1953, before returning to France as a founding member of the Institut des Hautes \'Etudes Scientifiques. He moved to the University of Nice to found the Department of Mathematics in 1964, and retired in 1970. He was elected as a member of the Acad\'emie des Sciences in 1968. 
 Dieudonn\'e drafted much of the Bourbaki series of texts, the many volumes of the EGA algebraic geometry series, and nine volumes of his own ``\'El\'ements d'Analyse''. He also wrote individual monographs on Infinitesimal Calculus, Linear Algebra and Elementary Geometry, invariant theory, commutative algebra, algebraic geometry, and formal groups. With Laurent Schwartz he supervised the early research of Alexander Grothendieck. Later from 1959 to 1964 he was at the Institut des Hautes Études Scientifiques alongside Grothendieck, and collaborating on the expository work needed to support the project of refounding algebraic geometry on the new basis of schemes.} determinants, which are defined as follows.

Let $(X,uow)$ be a Desarguesian based Playfair plane. For every automorphism $A\in\Aut(X)$, its \index{Dieudonn\'e determinant}\index{determinant!Dieudonn\'e}\defterm{Dieudonn\'e determinant} is the coscalar $[A]\in[\IR_X]$, defined by the formula
$$[A]\defeq\begin{cases}[b_2{\cdot}a_2^{-1}{\cdot}a_1{\cdot}a_2-b_1{\cdot}a_2]&\mbox{if $a_2\ne 0$};\\
[b_2{\cdot}a_1]&\mbox{if $a_2=0$};
\end{cases}
$$where $a_1,a_2,b_1,b_2\in\IR_X$ are unique scalars such that $Aou\in a_1\cdot\overvector{ou}+a_2\cdot\overvector{ow}$ and $Aow\in b_1\cdot\overvector{ou}+b_2\cdot\overvector{ow}$.

\begin{proposition}\label{p:DDet} Let $(X,uow)$ be a Desarguesian based Playfair plane. For any automorphisms $A,B\in\Aut(X)$, the Dieudonne determinant has the following properties:
\begin{enumerate}
\item[\textup{(1)}] $[BA]=[A]$ if $B$ is a translation of $X$;
\item[\textup{(2)}] $[BA]=[A]$ if $B$ is a horizontal shear of $X$;
\item[\textup{(3)}] $[BA]=[A]$ if $B$ is a vertical shear of $X$;
\item[\textup{(4)}] $[BA]=[A]$ if $B$ is a hypershear of $X$;
\item[\textup{(5)}] $[BA]=[A]$ if $B$ is an equiaffine automorphism of $X$;
\item[\textup{(6)}] $[BA]=[b]\cdot[A]$ if $B$ is a horizontal scale of $X$ with $Bou\in b\cdot\overvector{ou}$; 
\item[\textup{(7)}] $[BA]=[b]\cdot[A]$ if $B$ is a vertical scale of $\Pi$ with $Bow\in b\cdot\overvector{ow}$; 
\item[\textup{(8)}] $[BA]=[B]\cdot [A]$ if $B$ is an affine automorphism of $X$.
\end{enumerate}
\end{proposition}

\begin{proof} Identifying each point $p\in X$ with the coordinates of the vector $\overvector{op}$ in the base $\overvector{ou},\overvector{ow}$, we identify the based affine plane $X$ with the coordinate plane $\IR_X\times\IR_X$ of the scalar corps $\IR_X$. Let $a_1,a_2,b_1,b_2\in\IR_X$ be unique scalars such that $\overvector{Aou}=a_1\cdot\overvector{ou}+a_2\cdot\overvector{ow}$ and $\overvector{Aow}\in b_1\cdot\overvector{ou}+b_2\cdot\overvector{ow}$. Writing formulas, we shall sometimes omit the notation $\cdot$ for the dot operation in the scalar corps $\IR_X$.
\smallskip

1. If $B$ is a translation of $X$, then $\overvector{BAou}=\overvector{Aou}$ and $\overvector{BAow}=\overvector{Aow}$, and hence $[BA]=[A]$, by the definition of the Dieudonn\'e determinant.
\smallskip

2. If $B$ is a horizontal shear of $X$, then $\overvector{Bow}=s{\cdot} \overvector{ou}+\overvector{ow}$ for some scalar $s\in\IR_X$. By Proposition~\ref{p:horizontal-shear}, the shear $B$ maps any point with coordinates $(x,y)$ to the point with coordinates $(x+y{\cdot}s,y)$ in the affine base $uow$. Then $BA(0,1)=B(a_1,a_2)=(a_1+a_2{\cdot} s,a_2)$ and $BA(1,0)=(b_1+b_2{\cdot}s,b_2)$.

If $a_2\ne 0$, then $[BA]=[b_2a_2^{-1}(a_1+a_2s)a_2-(b_1+b_2s)a_2]=
[b_2a_2^{-1}a_1a_2-b_1a_2]=[A].$

If $a_2=0$, then $[BA]=[b_2(a_1+a_2s)]=[b_2a_1]=[A]$.

3. If $B$ is a vertical shear of $\Pi$, then $\overvector{Bou}=\overvector{ou}+s\cdot\overvector{ow}$ for some scalar $s\in\IR_\Pi$. By Proposition~\ref{p:(x,y)->(x,xc+y)}, the shear $B$ maps a point with coordinates $(x,y)$ to the point with coordinates $(x,y+x{\cdot}s)$ in the affine base $uow$. Then $BA(0,1)=B(a_1,a_2)=(a_1,a_2+a_1{\cdot}s)$ and $BA(1,0)=(b_1,b_2+b_1{\cdot}s)$. If $s=0$, then $BA=A$ and hence $[BA]=[A]$. So, we assume that $s\ne 0$.

If $0\notin\{a_1,a_2,a_2+a_1{\cdot}s\}$, then 
$$
\begin{aligned}
[BA]&=[(b_2+b_1s)(a_2+a_1s)^{-1}a_1(a_2+a_1s)-b_1(a_2+a_1s)]\\
&=
[\big((b_2+b_1s)a_1-b_1a_1^{-1}(a_2+a_1s)a_1\big)\cdot a_1^{-1}(a_2+a_1s)^{-1}a_1(a_2+a_1s)]\\
&=
[(b_2+b_1s)a_1-b_1a_1^{-1}(a_2+a_1s)a_1]=
[b_2a_1-b_1a_1^{-1}a_2a_1]\\
&=[(b_2a_2^{-1}a_1a_2-b_1a_2)\cdot a_2^{-1}a_1^{-1}a_2a_1]=[b_2a_2^{-1}a_1a_2-b_1a_2]=[A]
\end{aligned}
$$
and we are done.

If $a_1=0$, then $a_2+a_1{\cdot}s=a_2\ne 0$ because $A$ is bijective and hence 
$$
\begin{aligned}
[BA]&=[(b_2+b_1s)(a_2+a_1s)^{-1}a_1(a_2+a_1s)-b_1(a_2+a_1s)]=
[-b_1a_2]=[b_2a_2^{-1}a_1a_2-b_1a_2]=[A].
\end{aligned}
$$
If $a_2=0$, then $a_2+a_1s=a_1s\ne 0$ and hence
$$
\begin{aligned}
[BA]&=[(b_2+b_1s)(a_2+a_1s)^{-1}a_1(a_2+a_1s)-b_1(a_2+a_1s)]\\
&=
[\big((b_2+b_1s)a_1-b_1a_1^{-1}(a_2+a_1s)a_1\big)\cdot a_1^{-1}(a_2+a_1s)^{-1}a_1(a_2+a_1s)]\\
&=
[(b_2+b_1s)a_1-b_1a_1^{-1}(a_2+a_1s)a_1]=
[b_2a_1]=[A]
\end{aligned}
$$
 
Finally, assume that $a_1\ne 0\ne a_2$ and $a_2+a_1s=0$.
Then 
$$
\begin{aligned}
[BA]&=[(b_2+b_1s)a_1]=
[(b_2+b_1s)a_1-b_1a_1^{-1}(a_2+a_1s)a_1]=
[b_2a_1-b_1a_1^{-1}a_2a_1]\\
&=[(b_2a_2^{-1}a_1a_2-b_1a_2)\cdot a_2^{-1}a_1^{-1}a_2a_1]=[b_2a_2^{-1}a_1a_2-b_1a_2]=[A].
\end{aligned}
$$
\smallskip

4. Assume that $B$ is a hypershear of $X$. Find a translation $T$ of $X$ such that $T(o)\in\Fix(B)$. Then $C\defeq T^{-1}BT$ is a hypershear that fixes all points of some line $L$ that contains the origin $o$.  If $C$ is a vertical shear, then $[BA]=[TCT^{-1}A]=[CT^{-1}A]=[T^{-1}A]=[A]$ by the items (1) and (3). So, assume that the line $L\subseteq\Fix(C)$ is not vertical. Since the Desarguesian plane $\Pi$ is shear, there exists a vertical shear $V$ of $\Pi$ such that $V[\Aline ou]=L$. Then $V^{-1}CV$ is a horizontal shear. By the items (1)--(3), 
$$
\begin{aligned}
[BA]&=[TCT^{-1}A]=[CT^{-1}A]=[VV^{-1}CVV^{-1}T^{-1}A]\\
&=[(V^{-1}TCV)V^{-1}T^{-1}A)]=[V^{-1}T^{-1}A]=[T^{-1}A]=[A].
\end{aligned}
$$
\smallskip

5. If $B$ is an equiaffine atomorphism, then $B=S_1\cdots S_n$ for some hypershears $S_1,\dots,S_n$. By induction, the preceding the item implies $$[BA]=[S_1\cdots S_nA]=[S_2\cdots S_nA]=\dots=[S_nA]=[A].$$
\smallskip

6. Assume that $B$ is a horizontal scale with $Bou\in s\cdot\overvector{ou}$ for some nonzero scalar $s\in\IR_X$. By Proposition~\ref{p:hscale}, the horizontal scale $B$ maps a point with coordinates $(x,y)$ to the point with coordinates $(x{\cdot} s,y)$. Then $BA(1,0)=B(a_1,a_2)=(a_1{\cdot}s,a_2)$, $BA(0,1)=B(b_1,b_2)=(b_1{\cdot}s,b_2)$. If $a_2\ne0$,
then 
$$[BA]=[b_2a^{-1}_2a_1sa_2-b_1sa_2]=[(b_2a_2^{-1}a_1a_2-b_1a_2)\cdot a_2^{-1}sa_2]=[A]\cdot[a_2^{-1}sa_2s^{-1}\cdot s]=[A]\cdot[s]=[s]\cdot [A].$$
If $a_2=0$, then  $[BA]=[b_2a_1s]=[b_2a_1]\cdot[s]=[A]\cdot[s]=[s]\cdot [A]$. 
\smallskip

7. Assume that $B$ is a vertical scale with $Bow\in s\cdot\overvector{ow}$ for some nonzero scalar $s\in\IR_X$. By Theorem~\ref{p:vscale}, the vertical scale $B$ maps a point with coordinates $(x,y)$ to the point with coordinates $(x,y{\cdot}s)$. Then $BA(1,0)=B(a_1,a_2)=(a_1,a_2{\cdot}s)$, $BA(0,1)=B(b_1,b_2)=(b_1,b_2{\cdot}s)$. If $a_2\ne0$,
then 
$$[BA]=[b_2s(a_2s)^{-1}a_1(a_2s)-b_1(a_2s)]=
[(b_2a_2^{-1}a_1a_2-b_1a_2)s]=[b_2a_2^{-1}a_1a_2-b_1a_2]\cdot[s]=[A]\cdot[s].$$
If $a_2=0$, then $[BA]=[b_2sa_1]=[b_2a_1\cdot a^{-1}sa_1s^{-1}\cdot s]=[b_2a_1]\cdot[s]=[A]\cdot[s]=[s]\cdot [A]$.
\smallskip

8. Assume that $B$ is an affine automorphism. By Lemma~\ref{l:affine=equi+scale}, $B$ is a composition $EV$ of a vertical scale $V$ and an equiaffine automorphism $E$ of the plane $X$. By Proposition~\ref{p:vscale}, there exists a scalar $s\in \IR_X$ such that the vertical scale $V$ maps any point with coordinates $(x,y)$ to the point with coordinates $(x,y{\cdot} s)$. It is clear that the identity automorphism of $X$ has Dieudonn\'e determinant $[I]=[1]$. The items (5) and (7) ensure that $[B]=[EVI]=[VI]=[s]\cdot [I]=[s]\cdot[1]=[s]$. On the other hand, $[BA]=[EVA]=[VA]=[s]\cdot[A]=[B]\cdot[A]$.
\end{proof}

\begin{lemma}\label{l:affine=equi+scale} Let $(X,uow)$ be a Desarguesian based affine plane. Every affine automorphism $A$ of the plane $X$ is equal to the composition $EV$ of a vertical scale $V$ and an equiaffine automorphism $E$ of $X$.
\end{lemma}

\begin{proof} Given any affine automorphism $A$ of $X$, consider the triangle $u'o'w'\defeq Auow$. By Corollary~\ref{c:EAff2}, there exists an equiaffine automorphism $E\in\EAff(X)$ such that $Euo=u'o'$ and $E(w)\in \Aline {o'}{w'}$. Then $w'\in E[\Aline ow]$ and the point $w''\defeq E^{-1}(w')$ belongs to the line $\Aline ow$. By Proposition~\ref{p:lv-scale=>rdist+ass-dot}, the Desarguesian based affine plane $(X,uow)$ is vertical-scale and hence admits a vertical scale $V:X\to X$ such that $Vuow=uow''$. For every $x\in \Aline ou$, we have $Voxu=oxu$, which implies that $V$ is a scalarity and hence affinity of the Desarguesian affine plane $X$.  Then the composition $EV$ is an affine automorphism of the plane $X$ such that $EVuow=Euow''=u'o'w'=Aouw$. Corollary~\ref{c:affine-extend-auto} ensures that $A=EV$. Therefore,  $A$ is the composition of the vertical scale $V$ and the equiaffine automorphsm $E$ of $X$.
\end{proof}

In the proofs of Theorems~\ref{t:equiaff<=>DDet=1}, \ref{t:equiaff<=>DDet=1}, and Proposition~\ref{p:sign-stabilizer}, we shall use the following lemma that detects equiaffine automorphisms among vertical shears.

\begin{lemma}\label{l:commutator=equiaffine} Let $R$ be an alternative ring and $\mathcal N(R)\defeq\{s\in R:\forall x,y\in R\;\;[x,s,y]=0\}$ be its nucleus.  For any nonzero elements $s,t\in \mathcal N(R)$, the function $A:R\times R\to R\times R$, $A:(x,y)\mapsto (x,y\cdot sts^{-1}t^{-1})$, is a composition of vertical and horizontal shears, and hence $A$ is an equiaffine automorphism of the affine plane $R\times R$.
\end{lemma} 

\begin{proof}  To shorten notations, for every element $x\in R^*\defeq R\setminus\{0\}$ we shall denote its multiplicative inverse $x^{-1}$ by $\bar x$. Let $\Pi=R\times R$ be the coordinate plane of the corps $R$. By Propositions~\ref{p:horizontal-shear} and \ref{p:(x,y)->(x,xc+y)}, for every $s\in R$, the function $$H_s:\Pi\to\Pi,\quad H_s:(x,y)\mapsto (x+y{\cdot}s,y)=(x,y)\cdot\Big(\!\!\begin{array}{cc}1&0\\s&1\end{array}\!\!\Big),$$ is a horizontal shear, and the function 
$$V_s:\Pi\to\Pi,\quad V_s:(x,y)\mapsto (x,x{\cdot}s+y)=(x,y)\cdot\Big(\!\!\begin{array}{cc}1&s\\0&1\end{array}\!\!\Big),$$ is a vertical shear of the based affine plane $\Pi$. 

The matrices $$\Big(\begin{array}{cc}1&0\\s&1\end{array}\Big)\quad \mbox{and}\quad\Big(\begin{array}{cc}1&s\\0&1\end{array}\Big)$$are called the matrices of the horizontal and vertical shears, respectively.

For every $s\in R$ consider the bijective map $$A_s:\Pi\to\Pi,\quad A_s:(x,y)\mapsto (x,y{\cdot} s)=(x,y)\cdot\Big(\begin{array}{cc}1&0\\0&s\end{array}\Big).$$ If $A_s$ is an automorphism of the plane $\Pi$, then it is a vertical scale of $\Pi$, by Proposition~\ref{p:vscale}.

We claim that 
$$V_{ac\bar a\bar c-\bar c}V_{-a+a\bar c}H_{\bar a}V_{-a}V_{ac}H_{a\bar c\bar a-\bar c\bar a}V_{-ac\bar a}V_1H_{-1}=A_{ac\bar a\bar c}$$ and hence $A_{ac\bar a\bar c}$ is both a vertical scale and an equiaffine automorphism of the affine plane $\Pi$. Since all parameters in the indexes of the horizontal and vertical shears belong to the nucleus of the alternative ring $R$, the above equality is equivalent to the matrix equality
$$
\begin{aligned}
&\Big(\!\!
\begin{array}{cc}1&0\\\!-1&1\end{array}\!\!\Big)
\Big(\!\!\begin{array}{cc}1&1\\0&1\end{array}\!\!\Big)
\Big(\!\!\begin{array}{cc}1&\!\!-ac\bar a\\0&1\end{array}\!\!\Big)
\Big(\!\!\begin{array}{cc}1&0\\a\bar c\bar a{-}\bar c\bar a\!\!&1\end{array}\!\!\Big)\Big(\!\!\begin{array}{cc}1&\!ac\\0&1\end{array}\!\!\Big)
\Big(\!\!\begin{array}{cc}1&\!\!-a\\0&1\end{array}\!\!\Big)
\Big(\!\!\begin{array}{cc}1&0\\\bar a&1\end{array}\!\!\Big)
\Big(\!\!\begin{array}{cc}1&\!\!{-}a{+}a\bar c\\0&1\end{array}\!\!\Big)
\Big(\!\!\begin{array}{cc}1&\!ac\bar a\bar c{-}\bar c\\0&1\end{array}\!\!\Big)\\
&=\Big(\begin{array}{cc}1&0\\0&ac\bar a\bar c\end{array}\Big)
\end{aligned}
$$
To see that this matrix equality is indeed true, we shall multiply the matrices of shears in the order marked by central dots:
$$
\begin{aligned}
&\Big(\!\!\begin{array}{cc}1&0\\\!-1&1\end{array}\!\!\Big){\cdot}
\Big(\!\!\begin{array}{cc}1&1\\0&1\end{array}\!\!\Big)
\Big(\!\!\begin{array}{cc}1&\!\!-ac\bar a\\0&1\end{array}\!\!\Big)
\Big(\!\!\begin{array}{cc}1&0\\a\bar c\bar a{-}\bar c\bar a\!\!&1\end{array}\!\!\Big)\Big(\!\!\begin{array}{cc}1&\!ac\\0&1\end{array}\!\!\Big)
\Big(\!\!\begin{array}{cc}1&\!\!-a\\0&1\end{array}\!\!\Big){\cdot}
\Big(\!\!\begin{array}{cc}1&0\\\bar a&1\end{array}\!\!\Big)
\Big(\!\!\begin{array}{cc}1&\!\!{-}a{+}a\bar c\\0&1\end{array}\!\!\Big)
\Big(\!\!\begin{array}{cc}1&\!ac\bar a\bar c{-}\bar c\\0&1\end{array}\!\!\Big)\\
&=\Big(\!\!\begin{array}{cc}1&1\\\!-1&0\end{array}\!\!\Big)
\Big(\!\!\begin{array}{cc}1&\!\!-ac\bar a\\0&1\end{array}\!\!\Big)
\Big(\!\!\begin{array}{cc}1&0\\a\bar c\bar a{-}\bar c\bar a\!\!&1\end{array}\!\!\Big)\Big(\!\!\begin{array}{cc}1&\!ac\\0&1\end{array}\!\!\Big)
{\Big(\!\!\begin{array}{cc}0&-a\\\bar a&1\end{array}\!\!\Big){\cdot}
\Big(\!\!\begin{array}{cc}1&\!\!{-}a{+}a\bar c\\0&1\end{array}\!\!\Big)}
\Big(\!\!\begin{array}{cc}1&\!ac\bar a\bar c{-}\bar c\\0&1\end{array}\!\!\Big)\\
&=\Big(\!\!\begin{array}{cc}1&1\\\!-1&0\end{array}\!\!\Big)
\Big(\!\!\begin{array}{cc}1&\!\!-ac\bar a\\0&1\end{array}\!\!\Big)
\Big(\!\!\begin{array}{cc}1&0\\a\bar c\bar a{-}\bar c\bar a\!\!&1\end{array}\!\!\Big){\Big(\!\!\begin{array}{cc}1&\!ac\\0&1\end{array}\!\!\Big){\cdot}
\Big(\!\!\begin{array}{cc}0&-a\\\bar a&\bar c\end{array}\!\!\Big)}
\Big(\!\!\begin{array}{cc}1&\!ac\bar a\bar c{-}\bar c\\0&1\end{array}\!\!\Big)\\
&=\Big(\!\!\begin{array}{cc}1&1\\\!-1&0\end{array}\!\!\Big)
\Big(\!\!\begin{array}{cc}1&\!\!-ac\bar a\\0&1\end{array}\!\!\Big)
{\Big(\!\!\begin{array}{cc}1&0\\a\bar c\bar a{-}\bar c\bar a\!\!&1\end{array}\!\!\Big){\cdot}
\Big(\!\!\begin{array}{cc}ac\bar a&0\\\bar a&\bar c\end{array}\!\!\Big)}
\Big(\!\!\begin{array}{cc}1&\!ac\bar a\bar c{-}\bar c\\0&1\end{array}\!\!\Big)\\
&=\Big(\!\!\begin{array}{cc}1&1\\\!-1&0\end{array}\!\!\Big)
{\Big(\!\!\begin{array}{cc}1&\!\!-ac\bar a\\0&1\end{array}\!\!\Big){\cdot}
\Big(\!\!\begin{array}{cc}ac\bar a&0\\1&\bar c\end{array}\!\!\Big)}
\Big(\!\!\begin{array}{cc}1&\!ac\bar a\bar c{-}\bar c\\0&1\end{array}\!\!\Big)\\
&=\Big(\!\!\begin{array}{cc}1&1\\\!-1&0\end{array}\!\!\Big)
{\Big(\!\!\begin{array}{cc}0&-ac\bar a\bar c\\1&\bar c\end{array}\!\!\Big){\cdot}
\Big(\!\!\begin{array}{cc}1&\!ac\bar a\bar c{-}\bar c\\0&1\end{array}\!\!\Big)}={\Big(\!\!\begin{array}{cc}1&1\\\!-1&0\end{array}\!\!\Big){\cdot}
\Big(\!\!\begin{array}{cc}0&-ac\bar a\bar c\\1&ac\bar a\bar c\end{array}\!\!\Big)}=
\Big(\!\!\begin{array}{cc}1&0\\0&ac\bar a\bar c\end{array}\!\!\Big).
\end{aligned}
$$
Therefore, $A_{ac\bar a\bar c}$ is indeed an equiaffine automorphism of the plane $\Pi=R\times R$. 
\end{proof}

Now we are able to prove the main result of this section.

\begin{theorem}[Dieudonn\'e, 1943]\label{t:equiaff<=>DDet=1} Let $(X,uow)$ be a Desarguesian based Playfair plane. An affine automorphism $A$ of $X$ is equiaffine if and only if its Dieudonn\'e determinant $[A]$ equals $[1]$.
\end{theorem}

\begin{proof} If $A$ is an equiaffine automorphism, then $[A]=[AI]=[I]=[1]$ by Proposition~\ref{p:DDet}(5). Now assume that $A$ is an affine automorphism with $[A]=[1]$. By Lemma~\ref{l:affine=equi+scale}, $A=EV$ for some equaffine automorphism $E$ and some vertical scale $V$ of the based affine plane $(\Pi,uow)$. By Proposition~\ref{p:DDet}(5,7), $[1]=[A]=[EV]=[V]=[s]$, where $s\in\IR_\Pi$ is a unique scalar such that $\overvector{Vow}=s\cdot\overvector{ow}$. It follows from $s\in[1]=[\IR^*_X,\IR^*_X]$ that $s=s_1\cdots s_n$ is a finite product of commutators. For every $i\in\{1,\dots,n\}$, consider the vertical scale $V_i$ that maps a point with coordinates $(x,y)$ to the point with coordinates $(x,y{\cdot} s_i)$. By Lemma~\ref{l:commutator=equiaffine}, the vertical scales $V_1,\dots,V_n$ are equiaffine and so is their composition $V=V_n\cdots V_1$. Then the automorphism $A=EV$ is equiaffine, too.
\end{proof}

\begin{corollary}[Dieudonn\'e, 1943] For every Desarguesian Playfair plane $X$, the quotient group $\Aff(X)/\EAff(X)$ is isomorphic to the abelian group $[\IR^*_X]$ of nonzero coscalars.
\end{corollary}

\begin{proof} Fix an affine base $uow$ in the Desarguesian Playfair plane $X$, and consider the function $\Det:\Aff(X)\to[\IR_X]$ assigning to each affine automorphism $A\in\Aff(X)$ its Dieudonn\'e determinant $[A]$. Proposition~\ref{p:DDet}(8,7) ensures that $\Det$ is a homomorphism of the group $\Aff(X)$ into the commutative semigroup $[\IR_X]$ such that $\Det[\Aff(X)]=[\IR^*_\Pi]=[\IR_\Pi]\setminus\{[0]\}$. By Theorem~\ref{t:equiaff<=>DDet=1}, the kernel of the surjective homomorphism $\Det:\Aff(X)\to[\IR^*_X]$ coincides with the subgroup $\EAff(X)$, which implies that the group $[\IR^*_X]$ is isomorphic to the quotient group $\Aff(X)/\EAff(X)$.
\end{proof}

\section{The (absolute) area form of a Thalesian Playfair plane}

In this section we define the area form and the absolute area form of a Thalesian affine plane and shall prove some basic symmetry properies of those forms.

First, we observe that the (absolute) area of a triangle in a Thalesian Playfair plane depends only on the vectors of sides of the triangle.

\begin{lemma}\label{l:area-form} Let $aob$ and $a'o'b'$ be two triangles in a Thalesian Playfair plane. If $\overvector{oa}=\overvector{o'a'}$ and $\overvector{ob}=\overvector{o'b'}$, then $[oab]=[o'a'b']$ and $[\![oab]\!]=[\![o'a'b']\!]$.
\end{lemma}

\begin{proof} Since $X$ is a Thalesian Playfair plane it is a translation plane, by Theorem~\ref{t:paraD<=>translation}. Then there exists a translation $T$ of $X$ such that $T(o)=o'$. Since translations preserve vectors, $Toa\in \overvector{oa}= \overvector{o'a'}$ and hence $a'=T(a)$. By analogy we can show that $b'=T(b)$. Then $Toab=o'a'b'$ and hence $[oab]=[o'a'b']$ and $[\![oab]\!]=[\![o'a'b']\!]$, by Corollary~\ref{c:Thales-trans=>SAut}.
\end{proof}

Let $X$ be a Thalesian Playfair plane, $\IR_X$ be its scalar corps, and $\overvector X$ be the $\IR_X$-module of vectors in $X$. Consider the function $$[{\cdot}\,,{\cdot}]:\overvector X\times\overvector X\to[X^3],\quad[\cdot\,,\cdot]:(\vec x,\vec y)\mapsto[\vec x,\vec y],$$ of two vector arguments, assigning to each pair of vectors $(\vec x,\vec y)$ the  area $[\vec x,\vec y]\defeq [oxy]$ of any triple $oxy\in X^3$ with $\overvector{ox}=\vec x$ and $\overvector{oy}=\vec y$. Lemma~\ref{l:area-form} ensures that the value $[\vec x,\vec y]=[oxy]$ does not depend on the choice of the triple $oxy$. The function \index[note]
{$[\cdot\hskip1pt,\cdot]$}$[\cdot\,,\cdot]:\overvector X\times \overvector X\to[X^3]$ is called the \index{area form}\defterm{area form} of the Thalesian affine plane $X$.

By definition, the \index{absolute area form}\defterm{absolute area form} of a Thalesian Playfair plane $X$ is the function \index[note]{$[\hskip-1.5pt[\cdot\hskip1pt,\cdot]\hskip-1.5pt]$}$$[\![\cdot\,,\cdot]\!]:\overvector X\times\overvector X\to[\![X^3]\!],\quad[\![\cdot\,,\cdot]\!]:(\vec x,\vec y)\mapsto[\![\vec x,\vec y]\!],$$  assigning to each pair of vectors $(\vec x,\vec y)\in\overvector X\times\overvector X$, the absolute area $[\![oxy]\!]$ of any triple $oxy\in X^3$ with  $\overvector{ox}=\vec x$ and $\overvector{oy}=\vec y$. Lemma~\ref{l:area-form} ensures that the area form is a well-defined function of two vector arguments.

By Corollary~\ref{c:area-semioriented-in-Playfair} any area $\alpha$ in a Playfair plane $X$ in semioriented and hence has a well-defined opposite area $-\alpha$. 

\begin{proposition}\label{p:areaform-T} For any vectors $\vec x,\vec y\in\overvector X$ in a Thalesian Playfair plane $X$,
\begin{enumerate}
\item[\textup{(1)}] $[\vec x,\vec x]=[0]$;
\item[\textup{(2)}] $[\vec y,\vec x]=-[\vec x,\vec y]$;
\item[\textup{(3)}] $[\![\vec x,\vec x]\!]=[\![0]\!]$;
\item[\textup{(4)}] $[\![\vec y,\vec x]\!]=[\![\vec x,\vec y]\!]$.
\end{enumerate}
\end{proposition}

\begin{proof} Choose any triple $oxy\in\overvector X$ such that $ox\in\vec x$ and $oy\in\vec y$. The definition of the (absolute) area form and Proposition~\ref{p:opposite-area} ensure that
$$[\vec y,\vec x]=[oyx]=-[oxy]=[\vec x,\vec y]\quad\mbox{and}\quad [\![\vec y,\vec x]\!]=[\![oyx]\!]=[\![oxy]\!]=[\![\vec x,\vec y]\!].$$
On the other hand, $$[\vec x,\vec x]=[oxx]=[0]\quad\mbox{and}\quad [\![\vec x,\vec x]\!]=[\![oxx]\!]=[\![0]\!].$$
\end{proof}

\section{Multiplying areas by scalars}

In this section we define the operation of multiplication of a signed area by a scalar, in Moufang Playfair planes. This operation is introduced in the following theorem.


\begin{theorem}\label{t:scalar-by-area} For every Moufang Playfair plane $X$, there exists a unique operation 
$$\cdot:\IR_X\times [X^3]\to[X^3],\quad \cdot:(s,\alpha)\mapsto s\cdot \alpha,$$such that $\overvector{oxe}\cdot [oep]=[oxp]$ for every Desarguesian triple $oxe$ and every point $p$ in $X$. This operation has the following properties, for all $s,t\in\IR_X$ and $\alpha\in[X^3]$:
\begin{enumerate}
\item[\textup{(0)}] $s\cdot \alpha=[0]$ if and only if $s=0$ or $\alpha=[0]$;
\item[\textup{(1)}] $1\cdot \alpha=\alpha$ and $(-1)\cdot\alpha=-\alpha$;
\item[\textup{(2)}] $s\cdot(t\cdot \alpha)=(s\cdot t)\cdot \alpha$.
\end{enumerate}
\end{theorem}

\begin{proof} Given a pair $(s,\alpha)\in \IR_X\times [X^3]$, choose any triple $oep\in \alpha$ such that $e\ne o\ne p$. Since $s$ is a scalar, there exists a unique point $x\in \Aline oe$ such that $oxe\in s$. Put $s\cdot \alpha\defeq [oxp]$. We claim that the product $s\cdot \alpha$ does not depend on the choice of a triple $oep\in \alpha$. Let $o'e'p'\in\alpha$ be another triple  with $e'\ne o'\ne p'$. Find a unique point $x'\in\Aline{o'}{e'}$ such that $o'x'e'\in s$. We have to show that $[oxp]=[o'x'p']$. 

If $\alpha=[0]$, then $\|\{o,e,p\}\|=2=\|\{o',e',p'\}\|$ and $\|\{o,x,p\}\|=2=\|\{o',x',p'\}\|$, which implies $[o'x'p']=[0]=[oxp]$. If $\alpha\ne[0]$, then the triples $oep$ and $o'e'p'$ are triangles. Since $[oep]=\alpha=[o'e'p']$ and the Playfair plane $X$ is Moufang, we can apply 
Theorem~\ref{t:[abc]=[xyz]<=>abc=Exyz} and find an equiaffine automorphism $E:X\to X$ of the plane $X$ such that $Eoep=o'e'p'$. By Proposition~\ref{p:affinity=>portionality=>scalarity}, shears are scalarities, which implies $Eoxe\in \overvector{oxe}=\overvector{o'x'e'}$ and hence $E(x)=x'$. Then $Eoxp=o'x'p'$ and hence $[oxp]=[o'x'p']$, witnessing that the multiplication operation $\cdot:\IR_X\times[X^3]\to[X^3]$ is well-defined.
\smallskip

To see that the product $\cdot:\IR_X\times[X^3]\to[X^3]$ is unique, assume that $*:\IR_X\times [X^3]\to[X^3]$ is another operation such that $\overvector{oxe}*[oep]=[oxp]$ for any Desarguesian triple $oxe$ and point $p$ in $X$. Take any pair $(s,\alpha)\in \IR_X\times[X^3]$. Choose any triple $oep\in\alpha$ with $o\ne e$, and find a unique point $x\in\Aline oe$ such that $oxe\in s$. Since $s$ is a scalar, the triple $oxe$ is Desarguesian. Then $s\cdot\alpha=\overvector{xoe}\cdot[oep]=[oxp]=\overvector{oxe}*[oep]=s*\alpha$.
\smallskip

It remains to prove the properties (0)--(3). Fix any signed area $\alpha\in[X^3]$ and find a triple $oab\in \alpha$ with $o\ne a$. 
\smallskip

0. Given any scalar $s\in\IR_X$,we have to show that $s\cdot\alpha=[0]$ if and only if $s=0$ or $\alpha=[0]$. Find a unique point $x\in\Aline oa$ such that $oxa\in s$. Then $s\cdot\alpha=\overvector{oxa}\cdot[oab]=[oxb]$. If $\alpha=[0]$, then the points $a,o,b$ are collinear and so are the points $o,x,b$, which implies $s\cdot\alpha=[0]$. If $s=0$, then $o=x$ and the triple $oxb=oob$ is not a triangle, which implies $s\cdot\alpha=[oxb]=[0]$.
Now assume that $\alpha\ne[0]$ and hence $oea\in\alpha$ is a triangle. Assuming that $[0]=s\cdot\alpha=[oxb]$, we conclude that $x\in \Aline ob\cap\Aline oa=\{o\}$ and hence $s=\overvector{oxa}=\overvector{ooa}=0$. Therefore, $s\cdot\alpha=[0]$ if and only if $s=0$ or $\alpha=[0]$.
\smallskip

1. Observe that $1\cdot \alpha=\overvector{oaa}\cdot[oab]=[oab]=\alpha$. Next, we prove that $(-1)\cdot\alpha=-\alpha$. If $\alpha=[0]$, then $(-1)\cdot \alpha=(-1)\cdot[0]=[0]=-[0]=-\alpha$. So, we assume that $\alpha\ne[0]$, in which case the triple $oab\in\alpha$ is a triangle. Find a unique point $c\in \Aline oa$ such that $\overvector{co}=\overvector{oa}$, which implies $\overvector{oca}=-1$. 

\begin{picture}(40,100)(-150,-10)
\linethickness{0.75pt}
\put(0,80){\color{lightgray}\line(1,-1){40}}
\put(0,0){\color{cyan}\vector(0,1){40}}
\put(0,40){\color{cyan}\vector(0,1){40}}
\put(0,0){\color{red}\vector(1,1){40}}
\put(0,40){\color{teal}\vector(1,0){40}}
\put(0,40){\color{red}\vector(1,1){40}}
\put(40,40){\color{cyan}\vector(0,1){40}}
\put(0,80){\color{teal}\line(1,0){40}}

\put(0,0){\circle*{3}}
\put(-8,-2){$c$}
\put(0,40){\circle*{3}}
\put(-8,38){$o$}
\put(40,40){\circle*{3}}
\put(43,38){$b$}
\put(0,80){\circle*{3}}
\put(-8,78){$a$}
\put(40,80){\circle*{3}}
\put(43,78){$d$}
\end{picture}

By Theorem~\ref{t:parallelogram3+1}, there exists a unique point $d\in X$ such that $boad$ is a parallelogram. Then $\overvector{od}=\overvector{oa}+\overvector{ob}=\overvector{co}+\overvector{ob}=\overvector{cb}$ and hence $\Aline od\parallel \Aline cb$. Then $codb$ also is a parallelogram. Applying Propositions~\ref{p:12-permutations} and \ref{p:opposite-area}, we conclude that $$(-1)\cdot\alpha=\overvector{oca}\cdot [oab]=[ocb]=[obd]=[oba]=-[oab]=-\alpha.$$
\smallskip

2. To check the condition (2), we have to show that $s\cdot (t\cdot \alpha)=(s\cdot t)\cdot\alpha$ for all scalars $s,t\in \IR_X$. If $s=0$, $t=0$ or $\alpha=[0]$, then $s\cdot(t\cdot \alpha)=[0]=(s\cdot t)\cdot \alpha$ and we are done. So, assume that $t\ne 0\ne s$ and $\alpha\ne[0]$. Find unique points $x\in \Aline oa$ and $y\in\Aline ox$ such that $oxa\in t$ and $oyx\in s$. Then 
$$s\cdot (t\cdot \alpha)=\overvector{oyx}\cdot(\overvector{oxa}\cdot [oab])=\overvector{oyx}\cdot [oxb]=[oyb]=\overvector{oya}\cdot[oab]=(\overvector{oyx}\cdot\overvector{oxa})\cdot[oab]=(s\cdot t)\cdot \alpha.$$
\smallskip

\end{proof}

\begin{proposition}\label{p:area-form1} For any Moufang Playfair plane $X$, the  area form is homogeneous by the first and second arguments in the sense that 
$$[t{\cdot}\vec x,\vec y]=t\cdot [\vec x,\vec y]\quad \mbox{and}\quad [\vec x,t{\cdot}\vec y]=t\cdot [\vec x,\vec y]$$ for all scalars $t\in\IR_X$ and vectors $\vec x,\vec y\in\overvector X$. 
\end{proposition}

\begin{proof} Given any vectors $\vec x,\vec y\in\overvector X$ and any scalar $t\in\IR_X$, choose any triple $oxy\in X^3$ with $ox\in\vec x$ and $oy\in\vec y$. If $o=x$, then $\vec x=\vec 0$ and $$[t{\cdot}\vec x,\vec y]=[t{\cdot}\vec 0,\vec y]=[\vec 0,\vec y]=[0]=t\cdot[0]=t\cdot [\vec 0,\vec y]=t\cdot[\vec x,\vec y].$$

So, assume that $o\ne x$. In this case we can choose a unique point $z\in \Aline ox$ scuh that $\overvector{ozx}=t$. By Theorem~\ref{t:scalar-by-vector}, $t\cdot\vec x=\overvector{ozx}\cdot\overvector{ox}=\overvector{oz}$. Then $[t{\cdot}\vec x,\vec y]=[\overvector{oz},\overvector{oy}]=[ozy]=\overvector{ozx}\cdot[oxy]=t\cdot[\vec x,\vec y]$.
On the other hand, $$
\begin{aligned}
[\vec x,t{\cdot}\vec y]&=-[t{\cdot}\vec y,\vec x]=(-1)\cdot[t{\cdot}\vec y,\vec x]=(-1)\cdot (t\cdot[\vec y,\vec x])=((-1)\cdot t)\cdot[\vec y,\vec x]\\
&=(t\cdot(-1))\cdot[\vec y,\vec x]=t\cdot((-1)\cdot[\vec y,\vec x])=t\cdot(-[\vec y,\vec x])=t\cdot[\vec x,\vec y],
\end{aligned}
$$
by  Theorem~\ref{t:scalar-by-area}(1) and the homogeneity of the signed area form by the first argument.
\end{proof}

Finally, we study the stabilizers the action of the multiplicative semigroup of scalars on areas. 

\begin{proposition}\label{p:sign-stabilizer} For every area $\alpha\in [X^3]$ in a Moufang Playfair plane $X$, its stabilizer $\Stab(\alpha)\defeq\{s\in \IR_X:s\cdot \alpha=\alpha\}$ contains the commutator subgroup $[\IR^*_X,\IR^*_X]$ of the multiplicative group $\IR^*_X$. If the plane $X$ is Desarguesian and $\alpha\ne[0]$, then the stabilizer $\Stab(\alpha)$ is equal to $[\IR^*_X,\IR^*_X]$. 
\end{proposition}

\begin{proof} If $\alpha=[0]$, then $\Stab(\alpha)=\IR_X\supseteq[\IR^*_X,\IR^*_X]$ and we are done. So, assume that $\alpha\ne[0]$. Fix any triple $uow\in\alpha$ and observe that it is an affine base in the Playfair plane $X$. Let $e$ be the diunit of the affine base $ouw$. By Theorem~\ref{t:affine-Moufang<=>}, the ternar $\Delta=\Aline oe$ of the Moufang based affine plane $(X,uow)$ is linear, distributive, associative-plus and alternative-dot. Then the biloop $(\Delta,+,\cdot)$ is an alternative ring. By Theorem~\ref{t:Ker(R)=RPi}, the scalar corps $\IR_X$ can be identified with the kernel $\Ker(\Delta)$ of the ternar $\Delta$. Since $\Delta$ is distributive and alternative-dot, the kerner $\Ker(\Delta)$ coincides with the nucleaus $\mathcal N(\Delta)=\{s\in\Delta:\forall x,y\in\Delta\;\;[x,s,y]=0\}$ of the alternative ring $\Delta$.

Identifying each point $p$ with its coordinates $(p_1,p_2)$ in the affine base $uow$, we identify the affine plane $X$ with the coordinate plane $\Delta\times\Delta$ of the alternative ring $\Delta$. By Lemma~\ref{l:commutator=equiaffine}, for any nonzero elements $s,t\in\mathcal N(\Delta)=\IR_X$, the function $E:X\to X$, $E:(x,y)\mapsto (x,y\cdot sts^{-1}t^{-1})$ is an equiaffine automorphism of the affine plane $X$. By Proposition~\ref{p:vscale}, this function is a vertical scale of the based affine plane $(X,uow)$. Observe that the point $w'=E(w)$ has coordinates $(0,1\cdot sts^{-1}t^{-1})$ and hence $\overvector{ow'w}=sts^{-1}t^{-1}$. Then $sts^{-1}t^{-1}\cdot \alpha=\overvector{ow'w}\cdot[uow]=\overvector{ow'w}\cdot[owu]=[ow'u]=[Eowu]=[owu]=[uow]=\alpha$ and hence $st^{-1}s^{-1}t^{-1}\in\Stab(\alpha)$ and $[\IR^*_X,\IR^*_X]\subseteq\Stab(\alpha)$.
\smallskip

Now assume that the affine plane $X$ is Desarguesian and $\alpha\ne[0]$. To show that $\Stab(\alpha)=[\IR^*_X,\IR^*_X]$, take any scalar $s\in \Stab(\alpha)$. Consider the vertical scale $V:X\to X$, $V:(x,y)\mapsto (x,y{\cdot} s)$, and observe that $Vou=ou$. For the point $w'\defeq V(w)$, we have the equality $\overvector{ow'w}=s$. Since $[owu]=[uow]=\alpha=s\cdot \alpha=\overvector{ow'w}\cdot[owu]=[ow'u]$, there exists an equiaffine automorphism $E$ of $X$ such that $Eowu=ow'u=Vowu$. Corollary~\ref{c:affine-extend-auto} ensures that $E=V$ and hence the vertical shear $V=E$ is equiaffine.   By Proposition~\ref{p:DDet}, $[1]=[E]=[V]=[s]$ and hence $s\in [\IR^*_X,\IR^*_X]$. Therefore, $\Stab(\alpha)=[\IR^*_X,\IR^*_X]$.
\end{proof}

\section{Mutiplying areas by coscalars}

Let $X$ be a Playfair plane and $\IR_X$ be its scalar corps. Let $\IR^*_X=\IR_X\setminus\{0\}$ be the multiplicative group of non-zero scalars, and $[\IR^*_X,\IR^*_X]$ be its commutator subgroup. For every scalar $s\in\IR_X$, let $[s]\defeq s\cdot[\IR^*_X,\IR^*_X]=[\IR^*_X,\IR^*_X]\cdot s$ be the coset of $s$ by the commutator subgroup $[\IR^*_X,\IR^*_X]$.

The set $[\IR_X]=\{[s]:s\in\IR_X\}$ is a commutative semigroup with respect to the binary operation $[s]\cdot [t]=[s\cdot t]=\{x\cdot y:x\in[s]\wedge y\in[t]\}$. This semigroup contains zero and unit elements $[0]=\{0\}$ and $[1]=[\IR^*_X,\IR^*_X]$. The set $[\IR_X]\setminus\{[0]\}$ is a subgroup of $[\IR_X]$, equal to the abelianization of the group $\IR^*_X$. On the other hand, the commutative semigroup $[\IR_X]$ is the abelianization of the multiplicative semigroup $(\IR_X,\cdot)$ of the scalar corps $\IR_X$. Elements of the commutative semigroup \index[note]{$\mathbb R_X$}$[\IR_X]$ are called \index{coscalar}\defterm{coscalars}.  

\begin{theorem}\label{t:coscalar-by-area} For every Moufang Playfair plane $X$, there exists a unique operation $$\cdot:[\IR_X]\times[X^3]\to[X^3],\quad \cdot:(s,\alpha)\mapsto s\cdot \alpha,$$ such that $[\overvector{oxe}]\cdot[oep]=[oxp]$ for every Desarguesian triple $oxe$ and every point $p$ in $X$. This operation has the following properties for all $s,t\in[\IR_X]$ and $\alpha\in[X^3]$:
\begin{enumerate}
\item[\textup{(0)}] $s\cdot\alpha=[0]$ if and only if $s=\{0\}$ or $\alpha=[0]$;
\item[\textup{(1)}] $[1]\cdot \alpha=\alpha$ and $[-1]\cdot\alpha=-\alpha$;
\item[\textup{(2)}] $s\cdot (t\cdot \alpha)=(s\cdot t)\cdot \alpha$. 
\end{enumerate}
If the plane $X$ is Desarguesian, then for any nonzero area $\alpha\in [X^3]$, the function $$\cdot_\alpha:[\IR_X]\to[X^3],\quad \cdot_\alpha:s\mapsto s\cdot\alpha,$$ is bijective.
\end{theorem}

\begin{proof} Given a pair $(c,\alpha)\in [\IR_X]\times [X^3]$, choose any scalar $s\in c$ and put $c\cdot \alpha\defeq s\cdot \alpha$, where the product $s\cdot\alpha$ has been defined in Theorem~\ref{t:scalar-by-area}. To show that the product $c\cdot \alpha$ is well-defined, we should check that the value of the product $s\cdot\alpha$ does not depend on the choice of the scalar $s\in c$. Let $t\in c$ be another scalar in the coscalar $c$. If $c=[0]=\{0\}$, then $s=0=t$ and $s\cdot\alpha=t\cdot \alpha$. So, assume that $c\ne[0]$ and hence $s\ne 0\ne t$. Then $s\cdot [\IR^*_X,\IR^*_X]=c=t\cdot[\IR^*_X,\IR^*_X]$ implies $s^{-1}t\in[\IR^*_X,\IR^*_X]$ and hence $s^{-1}t\cdot\alpha=\alpha$, according to  Proposition~\ref{p:sign-stabilizer}. By Theorem~\ref{t:scalar-by-area}(2), $s\cdot\alpha=s\cdot(s^{-1}t\cdot\alpha)=(s\cdot s^{-1}t)\cdot\alpha=t\cdot\alpha$. Therefore, the product $c\cdot \alpha$ is well-defined.

Now take any Desarguesian triple $oxe\in X^3$ and point $p\in X$. Then $[\overvector{oxe}]\cdot [oep]=\overvector{oxe}\cdot [oxp]$, by Theorem~\ref{t:scalar-by-area} and the definition of the product of the coscalar $[\overvector{oxe}]$ and the area $[oep]$.
\smallskip

To see that the product $\cdot:[\IR_X]\times[X^3]\to[X^3]$ is unique, assume that $*:[\IR_X]\times [X^3]\to[X^3]$ is another operation such that $[\overvector{oxe}]*[oep]=[oxp]$ for any Desarguesian triple $oxe$ and point $p$ in $X$. Take any pair $(c,\alpha)\in [\IR_X]\times[X^3]$. Choose any triple $oep\in\alpha$ with $o\ne e$. Choose any scalar $s\in c$ and find a unique point $x\in\Aline oe$ such that $oxe\in s$. Since $s$ is a scalar, the triple $oxe\in s$ is Desarguesian. Then $c\cdot\alpha=[\overvector{oxe}]\cdot[oep]=[oxp]=[\overvector{oxe}]*[oep]=c*\alpha$.
\smallskip

The properties (0)--(2) follow from the corresponding properties in Theorem~\ref{t:scalar-by-area}.
\smallskip

Now assume that the affine plane $X$ is Desarguesian, and take any nonzero area $\alpha\in[X^3]\setminus\{[0]\}$. We have to show that the map $\cdot_\alpha:[\IR_X]\to[X^3]$, $\cdot_\alpha:c\mapsto c\cdot\alpha$, is bijective.
Choose any triple $uow\in\alpha\ne[0]$ and observe that $uow$ is an affine base in the Desarguesian Playfair plane $X$. By Theorem~\ref{t:unique-coscalar} (proved below), for any area $\beta\in X^3$, there exists a unique coscalar $d\in[\IR_X]$ such that $d\cdot\alpha=d\cdot[uow]=\beta$. The existence and uniqueness of the coscalar $d$ implies that the map $\cdot_\alpha:[\IR_X]\to[X^3]$ is bijective.
\end{proof}

\begin{theorem}\label{t:unique-coscalar} For every Desarguesian based Playfair plane $(X,uow)$ and every triple $abc\in X^3$, there exists a unique coscalar $d\in[\IR_X]$ such that $[acb]=d\cdot[uow]$. The coscalar $d$ can be found by the formula
$$d=\begin{cases}[b_2a_2^{-1}a_1a_2-b_1a_2]&\mbox{if $a_2\ne 0$};\\
[b_2a_1]&\mbox{if $a_2=0$};
\end{cases}
$$where $a_1,a_2,a_3,a_4\in\IR_X$ are unique scalars such that $\overvector{ca}=a_1\cdot\overvector{ou}+a_2\cdot\overvector{ow}$ and  $\overvector{cb}=b_1\cdot\overvector{ou}+b_2\cdot\overvector{ow}$.
\end{theorem}

\begin{proof} Let $(X,uow)$ be a Desarguesian based Playfair plane, $e$ be the diunit of the affine base $uow$, and $\Delta=\Aline oe$ be the ternar of the based affine plane $(X,uow)$. By Theorem~\ref{t:ass-dot<=>}, $\Delta$ is linear, distributive and associative. By  Corollary~\ref{c:R=RPi}, the map $\Delta\to\IR_X$, $x\mapsto \overvector{oxe}$, is an isomorphism of the corps $(\Delta,+,\cdot)$ and $\IR_X$. 

Given any triple $abc\in X^3$, we should find a unique coscalar $[\delta]\in [\IR_X]$ such that $[\delta]\cdot[uow]=[abc]$. Find unique scalars $a_1,a_2,a_3,a_4\in\IR_X$ such that $\overvector{ca}=a_1\cdot\overvector{ou}+a_2\cdot\overvector{ow}$ and  $\overvector{cb}=b_1\cdot\overvector{ou}+b_2\cdot\overvector{ow}$.  Find a unique translation $T$ of the Desarguesian plane $X$ such that $T(c)=o$ and let $a'\defeq T(a)$ and $b'\defeq T(b)$. Then $\overvector{oa'}=\overvector{ca}=a_1\cdot\overvector{ou}+a_2\cdot\overvector{ow}$ and $\overvector{ob'}=\overvector{cb}=b_1\cdot\overvector{ou}+b_2\cdot\overvector{ow}$. By Corollary~\ref{c:coordinates=vectors}, $(a_1,a_2)$ and $(b_1,b_2)$ are the coordinates of the points $a'$ and $b'$ in the affine base $uow$, respectively. By Corollary~\ref{c:Thales-trans=>SAut}, the translation $T$ of the Playfair plane $X$ is area-preserving that hence $[abc]=[a'ob']$.

By Proposition~\ref{p:area-via-coordinates}, $[oa'b']=[oud]$, where the point 
$$d\defeq\begin{cases}b_2a_2^{-1}a_1a_2-b_1a_2&\mbox{if $a_2\ne 0$},\\
b_2a_1&\mbox{if $a_2=0$}
\end{cases}
$$
belongs to the ternar $\Delta=\Aline oe\subset X$.  The horizontal shear $S:X\to X$, $S:(x,y)\mapsto (x+y,y)$, witnesses that $[oue]=[Souw]=[ouw]$. Consider the scalar $\delta\defeq\overvector{ode}\in\IR_X$ and observe that 
$$
\begin{aligned}
\delta\cdot [ouw]&=\delta\cdot [oue]=\delta\cdot((-1)\cdot[oeu])=(\delta\cdot(-1))\cdot[oeu]=((-1)\cdot\delta)\cdot[oeu]=(-1)\cdot(\delta\cdot[oeu])\\
&=(-1)\cdot(\overvector{ode}\cdot[oeu])=(-1)\cdot[odu]=[oud]=[oa'b']=[a'b'o]=[abc],
\end{aligned}
$$ by Theorem~\ref{t:scalar-by-area} and Proposition~\ref{p:12-permutations}. Then the coscalar $[\delta]\in[\IR_X]$ has the desired property: $[\delta]\cdot[ouw]=[abc]$. 

It remains to prove that such a coscalar $[\delta]$ is unique. Given any scalar $s\in\IR_X$ with $[s]\cdot[ouw]=[abc]$, we have to prove that $[s]=[\delta]$.

If $[abc]=[0]$, then $[\delta]\cdot[ouw]=[abc]=[0]=[s]\cdot [ouw]$ implies $[\delta]=\{0\}=[s]$, by Theorem~\ref{t:scalar-by-area}(0). So, assume that $[abc]\ne[0]$. In this case Theorem~\ref{t:scalar-by-area}(0) ensures that $[\delta]\ne \{0\}\ne[s]$. So, we can consider the scalar $s^{-1}$ and applying Theorem~\ref{t:coscalar-by-area}(3), conclude that 
$$
\begin{aligned}
[s^{-1}\delta]\cdot[ouw]&=([s^{-1}]\cdot [\delta])\cdot[ouw]=[s^{-1}]\cdot([\delta]\cdot[ouw])=[s^{-1}]\cdot [abc]=[s^{-1}]\cdot([s]\cdot[ouw])\\
&=
([s^{-1}]\cdot[s])\cdot[ouw]=[s^{-1}\cdot s]\cdot [ouw]=[1]\cdot[ouw]=[ouw]
\end{aligned}
$$ and hence $s^{-1}\delta\cdot[ouw]=[ouw]$, by the definition of the multiplication by a coscalar. Therefore, the scalar $s^{-1}\delta$ belongs to the stabilizer $\Stab([ouw])=\{r\in\IR_X:r\cdot[ouw]=[ouw]\}$ of the nonzero signed area $[uow]$. By Proposition~\ref{p:sign-stabilizer}, this stabilizer is equal to the commutator subgroup $[\IR^*_X,\IR^*_X]$. Therefore, $s^{-1}d\in[\IR^*_X,\IR^*_X]$ and finally, $[s]=s\cdot[\IR^*_X,\IR^*_X]=\delta\cdot[\IR^*_X,\IR^*_X]=[\delta]$.
\end{proof}

\section{Multiplying absolute areas by scalars and absolute coscalars}

In Moufang Playfair planes, the operation of multiplication of an absolute area by a scalar can be defined by analogy with the operation of multiplication of an area by a scalar.

\begin{theorem}\label{t:scalar-by-absarea}  For every Moufang Playfair plane $X$, there exists a unique operation 
$$\cdot:\IR_X\times [\![X^3]\!]\to[\![X^3]\!],\quad \cdot:(s,\alpha)\mapsto s\cdot \alpha,$$such that $\overvector{oxe}\cdot [\![oep]\!]=[\![oxp]\!]$ for every Desarguesian triple $oxe$ and every point $p$ in $X$. This operation has the following properties for all $s,t\in\IR_X$ and $\alpha\in[\![X^3]\!]:$
\begin{enumerate}
\item[\textup{(0)}] $s\cdot \alpha=[\![0]\!]$ if and only if $s=0$ or $\alpha=[\![0]\!]$;
\item[\textup{(1)}] $1\cdot \alpha=\alpha$ and $(-1)\cdot\alpha=\alpha$;
\item[\textup{(2)}] $s\cdot(t\cdot \alpha)=(s\cdot t)\cdot \alpha$.
\end{enumerate}
\end{theorem} 

\begin{proof} Given a pair $(s,\alpha)\in \IR_X\times [\![X^3]\!]$, choose any triple $oep\in \alpha$. If $\alpha\ne[\![0]\!]$, then the points $o,e,p$ are pairwise distinct. If $\alpha=[\![0]\!]$, then the points $o,e,p$ are colinear, but we again can choose those points to be distinct. Since $s$ is a scalar, there exists a unique point $x\in \Aline oe$ such that $oxe\in s$. Put $s\cdot \alpha\defeq [\![oxp]\!]$. We claim that the product $s\cdot \alpha$ does not depend on the choice of a triple $oep\in \alpha$. Let $o'e'p'\in\alpha$ be another triple with $e'\ne o'\ne p'$. Find a unique point $x'\in\Aline{o'}{e'}$ such that $o'x'e'\in s$. We have to show that $[\![oxp]\!]=[\![o'x'p']\!]$.

If $\alpha=[\![0]\!]$, then $\|\{o,e,p\}\|=2=\|\{o',e',p'\}\|$ and $\|\{o,x,p\}\|=2=\|\{o',x',p'\}\|$, which implies $[\![o'x'p']\!]=[\![0]\!]=[\![oxp]\!]$. If $\alpha\ne[\![0]\!]$, then the triples $oep$ and $o'e'p'$ are triangles. If $s=0$, then $x=o$ and $x'=o'$ and hence $[\![oxp]\!]=[\![oop]\!]=[\![0]\!]=[\![o'o'p']\!]=[\![o'x'p']\!]$. So, assume that $s\ne 0$ and hence  $x'\ne o'$. Since $e'\notin\Aline{o'}{p'}$, the points $x'\in \Aline {o'}{e'}\setminus\{o'\}$ does not belong to the line $\Aline {o'}{p'}$. 

 Since $[\![oep]\!]=\alpha=[\![o'e'p']\!]$ there exists a permutation $\sigma\in S_3$ such that $[oep]=[o'e'p'\circ \sigma]$. Since the area $[oep]$ is semioriented, $[o'e'p'\circ\sigma]\in\{[o'e'p'],[p'e'o']\}$. By Theorem~\ref{t:[abc]=[xyz]<=>abc=Exyz}, there exists an equiaffine automorphism $F:X\to X$ of the plane $X$ such that $Eoep\in\{o'e'p',p'e'o'\}$. By Proposition~\ref{p:affinity=>portionality=>scalarity}, equiaffine automorphisms are affinities and hence scalarities. If $Foep=o'e'p'$, then $F(x)=x'$ because the equiaffine automorphism $F$ is a scalarity and $\overvector{oxe}=\overvector{o'x'e'}$. Then $Foxp=o'x'p'$ and $[\![oxp]\!]=[\![o'x'p']\!]$. 
 
\begin{picture}(500,120)(-50,-10)
{\color{lightgray}
\polygon*(0,10)(60,70)(120,10)
\polygon*(200,0)(240,90)(280,0)
}

\linethickness{0.7pt}
\put(0,10){\line(1,0){120}}
\put(0,10){\line(1,1){60}}
\put(60,70){\line(1,-1){60}}
\put(30,40){\line(3,-1){90}}
\put(200,0){\line(4,3){60}}
\put(280,0){\line(-4,3){60}}

\put(0,10){\circle*{3}}
\put(-8,8){$o$}
\put(120,10){\circle*{3}}
\put(123,8){$p$}
\put(60,70){\circle*{3}}
\put(58,73){$e$}
\put(30,40){\circle*{3}}
\put(22,39){$x$}

\put(135,40){\vector(1,0){45}}
{\linethickness{0.5pt}
\put(150,42){$F$}
}

\put(200,0){\color{teal}\line(1,0){80}}
\put(200,0){\line(4,9){40}}
\put(240,90){\line(4,-9){40}}
\put(220,45){\color{teal}\line(1,0){40}}

\put(200,0){\circle*{3}}
\put(191,-3){$p'$}
\put(280,0){\circle*{3}}
\put(283,-3){$o'$}
\put(240,90){\circle*{3}}
\put(238,93){$e'$}
\put(260,45){\circle*{3}}
\put(264,40){$x'$}
\put(261,48){$z$}
\put(220,45){\circle*{3}}
\put(212,44){$y$}
\end{picture} 
 
 So, assume that $Foep=p'e'o'$ and consider the point $y=F(x)\in \Aline {p'}{e'}$. Since $F$ is a scalarity, $\overvector{oxe}=\overvector{p'ye'}$. Since $X$ is a Moufang (and hence shear) plane, there exists a hypershear $E$ of $X$ such that $E(p')=o'$ and $e'\in\Fix(E)$. Consider the point $z\defeq E(y)$ and observe that $\overvector{o'ze'}=\overvector{p'ye'}=\overvector{oxe}=\overvector{o'x'e'}$ and hence $z=x'$. Then $\Aline {x'}y=\Aline zy\subparallel \Aline {o'}{p'}$. Since $X$ is a shear plane, there exists a shear $S:X\to X$ such that $\Aline {o'}{p'}\subseteq \Fix(S)$ and $S(y)=x'$. Then $\Phi\defeq SF$ is an equiaffine automorphism of the plane $X$ such that $\Phi oxp=p'x'o'$, witnessing that $[\![oxp]\!]=[\![o'x'p']\!]$.
Therefore, the multiplication operation $\cdot:\IR_X\times [\![X^3]\!]\to[\![X^3]\!]$ is well-defined.
 
The properties (0)--(2) of this multiplication operation can be proved by analogy with the corresponding properties in Theorem~\ref{t:scalar-by-area} (the equality $(-1)\cdot\alpha=\alpha$ for all $\alpha\in [\![X^3]\!]$ follows from the equality $[\![abc]\!]=(-[abc])\cup[abc]$ holding for every triple $abc\in X^3$).
\end{proof}

Next, we study the stabilizers the action of the multiplicative semigroup of scalars on areas. For a corps $R$ we denote by $R^*=R\setminus\{0\}$ its multiplicative group, and by $[R^*,R^*]$ the commutator of $R^*$. Since $[R^*,R^*]$ and $\{-1,1\}$ are normal subgroups of the group $R^*$, the set \index[note]{$[\hskip-1.5pt[R^*,R^*]\hskip-1.5pt]$}
$$[\![R^*,R^*]\!]\defeq\{-1,1\}\cdot[R^*,R^*]$$ is a normal subgroup of $R^*$. This subgroup will be called the \index{absolute commutator}\defterm{absolute commutator} of the group $R^*$. 

\begin{proposition}\label{p:area-stabilizer} For every area $\alpha\in [\![X^3]\!]$ in a Moufang Playfair plane $X$, its stabilizer $\Stab(\alpha)\defeq\{s\in \IR_X:s\cdot \alpha=\alpha\}$ contains the absolute commutator subgroup $[\![\IR^*_X,\IR^*_X]\!]$ of the multiplicative group $\IR^*_X$. If the plane $X$ is Desarguesian and $\alpha\ne[\![0]\!]$, then the stabilizer $\Stab(\alpha)$ equals  $[\![\IR^*_X,\IR^*_X]\!]$. 
\end{proposition}

\begin{proof} Since the absolute commutator $[\![\IR^*_X,\IR^*_X]\!]$ is a subgroup of $\IR_X$, generated by the set $[\IR^*_X,\IR^*_X]\cup\{-1\}$, it suffices to check that $[\IR^*_X,\IR^*_X]\cup\{-1\}\subseteq\Stab(\alpha)$. Given any scalar $s\in [\IR^*_X,\IR^*_X]\cup\{-1\}$, we should prove that $s\cdot \alpha=\alpha$. If $\alpha=[\![0]\!]$, then $s\cdot\alpha=s\cdot[\![0]\!]=[\![0]\!]=\alpha$ and we are done. So, assume that $\alpha\ne[\![0]\!]$. Choose any triple $oep\in\alpha$ and observe that $oep$ is a triangle. Find a unique point $x\in\Aline oe$ such that $\overvector{oxe}=s$. Theorems~\ref{t:scalar-by-area} and \ref{t:scalar-by-absarea} ensure that $s\cdot[oep]=[oxp]$ and $s\cdot[\![oep]\!]=[\![oxp]\!]$.

If $s\in [\IR^*_X,\IR^*_X]$, then $[oxp]=s\cdot[oep]=[oep]$, by Proposition~\ref{p:sign-stabilizer}. Then there exists an equiaffine automorphism $E$ of $X$ such that $Eoep=oxp$. The equality $Eoep=oxp$ implies $[\![oep]\!]=[\![oxp]\!]$ and hence $s\cdot\alpha=s\cdot[\![oep]\!]=[\![oxp]\!]=[\![oep]\!]=\alpha$ and $s\in\Stab(\alpha)$.

If $s=-1$, then $[oxp]=s\cdot[oep]=(-1)\cdot[oep]=-[oep]=[ope]$, by Theorem~\ref{t:scalar-by-area}(3). 
Then there exists an equiaffine automorphism $E$ of $X$ such that $Eope=oxp$. 
The equality $Eope=oxp$ implies $[\![oep]\!]=[\![ope]\!]=[\![oxp]\!]$ and hence $s\cdot\alpha=s\cdot[\![oep]\!]=[\![oxp]\!]=[\![oep]\!]=\alpha$ and $s=-1\in\Stab(\alpha)$.

Therefore, $[\IR^*_X,\IR^*_X]\cup\{-1\}\subseteq \Stab(\alpha)$, which implies $[\![\IR^*_X,\IR^*_X]\!]\subseteq \Stab(\alpha)$.
\smallskip

Now assuming that the affine plane $X$ is Desarguesian and $\alpha\ne[\![0]\!]$, we shall prove that $\Stab(\alpha)\subseteq[\![\IR^*_X,\IR^*_X]\!]$. Given any $s\in\Stab(\alpha)$, we should prove that $s\in[\![\IR^*_X,\IR^*_X]\!]$. Choose any triangle $oep\in\alpha$ and find a unique point $x\in\Aline oe$ such that $\overvector{oxe}=s$. Theorem~\ref{t:scalar-by-area} ensures that $[\![oep]\!]=\alpha=s\cdot\alpha=\overvector{oxe}\cdot[\![oep]\!]=[\![oxp]\!]=[oxp]\cup[pxo]$. Then there exists an equiaffine automorphism $A$ of $X$ such that $Aoep\in\{oxp,pxo\}$. If $Aoep=oxp$, then $[oep]=[oxp]=s\cdot[oep]$ and hence $s\in [\IR^*_X,\IR^*_X]\subseteq[\![\IR^*_X,\IR^*_X]\!]$, by Proposition~\ref{p:sign-stabilizer}.   If $Aoep=pxo$, then $[oep]=[pxo]$. By Theorem~\ref{t:scalar-by-area} and Proposition~\ref{p:opposite-area},
$$((-1)\cdot s)\cdot[oep]=(-1)\cdot(s\cdot[oep])=(-1)\cdot(\overvector{oxe}\cdot[oep])=(-1)\cdot[oxp]=[opx]=[pxo]=[oep].$$
Therefore, the scalar $(-1)\cdot s$ belongs to the stabilizer of the nonzero  area $[oep]$. By Proposition~\ref{p:sign-stabilizer}, $(-1)\cdot s\in[\IR^*_ 
X,\IR^*_X]$ and hence $s\in (-1)\cdot[\IR^*_X,\IR^*_X]\subseteq[\![\IR^*_X,\IR^*_X]\!]$.
\end{proof} 

Let $X$ be an affine plane and $\IR_X$ be its scalar corps. For every $s\in \IR_X$, consider the set 
\index[note]{$[\hskip-1.5pt[s]\hskip-1.5pt]$}$[\![s]\!]\defeq s\cdot[\![\IR^*_X,\IR^*_X]\!]$, called the \index{coscalar!absolute}\index{absolute coscalar}\defterm{absolute coscalar} of $s$. The set $[\![\IR_X]\!]\defeq\{[\![s]\!]:s\in\IR_X\}$ is a commutative semigroup with respect to the operation $[\![s]\!]\cdot[\![t]\!]=[\![s\cdot t]\!]$. The semigroup $[\![\IR_X]\!]$ has zero and unit elements $[\![0]\!]=\{0\}$ and $[\![1]\!]=[\![\IR^*_X,\IR^*_X]\!]$. The semigroup $[\![\IR_X]\!]$ is called the \index[note]{$[\hskip-1.5pt[\mathbb R_X]\hskip-1.5pt]$}
\index{semigroup of absolute coscalars}\defterm{semigroup of absolute coscalars} of the affine plane $X$. It is clear that the semigroup $[\![\IR_X]\!]$ is a homomorphic image of the semigroup of coscalars $[\IR_X]$.  Since $-1\in[\![\IR_X]\!]$, we have $[\![-1]\!]=[\![1]\!]=[\![\IR^*_X,\IR^*_X]\!]$.

\begin{theorem}\label{t:coscalar-by-absarea} For every Moufang Playfair plane $X$, there exists a unique operation $$\cdot:[\![\IR_X]\!]\times[\![X^3]\!]\to[\![X^3]\!],\quad \cdot:(s,\alpha)\mapsto s\cdot \alpha,$$ such that $[\![\overvector{oxe}]\!]\cdot[\![oep]\!]=[\![oxp]\!]$ for every Desarguesian triple $oxe$ and every point $p$ in $X$. This operation has the following properties for all $s,t\in[\![\IR_X]\!]$ and $\alpha\in[\![X^3]\!]$:
\begin{enumerate}
\item[\textup{(0)}] $s\cdot \alpha=[\![0]\!]$ if and only if $s=[\![0]\!]$ or $\alpha=[\![0]\!]$;
\item[\textup{(1)}] $[\![1]\!]\cdot\alpha=\alpha=[\![-1]\!]\cdot\alpha$;
\item[\textup{(2)}] $s\cdot (t\cdot \alpha)=(s\cdot t)\cdot \alpha$.
\end{enumerate}
If the affine plane $X$ is Desarguesian, then for every nonzero absolute area $\alpha\in[\![X^3]\!]\setminus\{[\![0]\!]\}$, the map $\cdot_\alpha:[\![\IR_X]\!]\to [\![X^3]\!]$, $\cdot_\alpha:s\mapsto s\cdot\alpha$, is bijective.
\end{theorem}

\begin{proof} This theorem can be deduced from Theorem~\ref{t:scalar-by-absarea} and Proposition~\ref{p:area-stabilizer} by analogy with the deduction of Theorem~\ref{t:coscalar-by-area} from Theorem~\ref{t:scalar-by-area} and Proposition~\ref{p:sign-stabilizer}.
\end{proof}

\section{Nonzero areols as torsors}

Theorems~\ref{t:coscalar-by-area} and \ref{t:coscalar-by-absarea} show that the nonzero areol $[X^\vartriangle]$ (i.e., the set of nonzero areas) in a Desarguesian affine plane $X$ have the structure of a $[\IR^*_X]$-torsor.

\index[person]{Serre}
\begin{definition}[Serre\footnote{{\bf Jean-Pierre Serre} (born 1926) is a French mathematician who has made contributions to algebraic topology, algebraic geometry and algebraic number theory. Jean-Pierre Serre is one of the most influential mathematicians of the twentieth century, whose work profoundly shaped topology, algebraic geometry, and number theory. Educated at the \'Ecole Normale Sup\'erieure in Paris, Serre obtained his doctorate in 1951 under Henri Cartan with a thesis on the homotopy groups of spheres --- a cornerstone of modern homotopy theory.
In the 1950s, Serre became a leading member of the Bourbaki group and played a central role in the reformulation of algebraic topology and algebraic geometry in the language of sheaves and categories. His early papers on spectral sequences and faisceaux alg\'ebriques coh\'erents introduced methods that revolutionized cohomological techniques. Serre’s later research established deep connections between topology, algebraic geometry, and number theory --- notably through the introduction of Galois cohomology, \'etale cohomology, and modular forms. His three celebrated monographs, ``Faisceaux alg\'ebriques coh\'erents'' (1955), ``Corps locaux'' (1962), and ``Cohomologie Galoisienne'' (1964), are considered classics.
Serre, at twenty-seven in 1954, was and still is the youngest person ever to have been awarded the Fields Medal. He went on to win the Balzan Prize in 1985, the Steele Prize in 1995, the Wolf Prize in Mathematics in 2000, and was the first recipient of the Abel Prize in 2003. Since 1956 he has been a professor at the Coll\`ege de France, where his lectures influenced generations of mathematicians.
Serre’s work is distinguished by its depth, elegance, and clarity --- qualities that have made his books and papers enduring models of mathematical style.
}
, 1955] A \index{torsor}\defterm{$G$-torsor} is a nonempty set $X$ endowed with a sharply transitive action $\cdot:G\times X\to X$ of a group $G$. 
\end{definition}

We recall that an \index{action of a group}\index{group action}\defterm{action} of a group $G$ on a set $X$ is any function $\cdot:G\times X\to X$, $\cdot:(g,x)\mapsto gx$, satisfying two axioms:
\begin{itemize}
\item $1x=x$ for all $x\in X$,
\item $g(hx)=(gh)x$ for all $g,h\in G$ and $x\in X$.
\end{itemize}
Here by $1$ we denote the neutral element of the group $G$.

An action $\cdot:G\times X\to X$ is \index{group action!transitive}\index{group action!sharp}\defterm{transitive} (resp. \defterm{sharp}) 
if for any $x\in X$, the function $\cdot_x:G\to X$, $\cdot_x:g\mapsto gx$, is surjective (resp. injective). 

In fact, we have already met many examples of torsors among Thalesian affine spaces.

\begin{example} \em Every Thalesian affine space $X$ is an $\overvector{X}$-torsor with respect to the natural action $+:\overvector{X}\times X\to X$ of the commutative group $\overvector{X}$ of functional vectors on the affine space $X$. 
\end{example}

Theorems~\ref{t:coscalar-by-area} and \ref{t:coscalar-by-absarea} yield two other natural examples of torsors. We recall that $X^\vartriangle$ denotes the set of all triangles in a liner $X$. Let $$[X^\vartriangle]\defeq\big\{[abc]:abc\in X^\vartriangle\big\}\quad\mbox{and}\quad
[\![X^\vartriangle]\!]\defeq\big\{[\![abc]\!]:abc\in X^\vartriangle\big\}.$$

\begin{example}\em  For every Desarguesian Playfair plane $X$, the nonzero  areol  $[X^\vartriangle]$ is a $[\IR^*_X]$-torsor with respect to the multiplication ${\cdot}:[\IR^*_X]\times[X^\vartriangle]\to[X^\vartriangle]$ of nonzero areas by nonzero coscalars.
\end{example}

\begin{example} \em For every Desarguesian Playfair plane $X$, the set of nonzero absolute areas $[\![X^\vartriangle]\!]$ is a $[\![\IR^*_ X]\!]$-torsor with respect to the multiplication $\cdot:[\![\IR^*_X]\!]\times[\![X^\vartriangle]\!]\to[\![X^\vartriangle]\!]$ of nonzero absolute areas by nonzero absolute coscalars.
\end{example}

\begin{remark} Since torsors can be thought as groups with ``forgotten'' neutral element, the areol $[X^3]$ of Desarguesian affine plane $X$ can be thought as the multiplicative semigroup of coscalars $[\IR_X]$, with forgotten identity element.  By analogy, the absolute areol $[\![X^3]\!]$ of $X$ can be thought as the multiplicative semigroup of absolute coscalars $[\![\IR_X]\!]$, with forgotten identity element. 
\end{remark} 

\begin{exercise} Show that for the real plane $X=\IR\times\IR$, the semigroup of coscalars $[\IR_X]$ is isomorphic to the multiplicative semigroup of real numbers and the semigroup of absolute coscalars $[\![\IR_X]\!]$ is isomorphic to the multiplicative semigroup of non-negative real numbers. 
\end{exercise}

\begin{exercise} Show that for the complex plane $X=\IC\times\IC$, the semigroup of coscalars $[\IR_X]$ is not equal to the semigroup of absolute coscalars $[\![\IR_X]\!]$, but both semigroups are isomorphic to the multiplicative semigroup of complex numbers. 
\end{exercise}

\begin{exercise} Show that the commutator subgroup $[\IH^*,\IH^*]$ of the multiplicative group $\IH^*$ of non-zero quaternions is equal to the set of quaternions with norm $1$. In particular, $-1=iji^{-1}j^{-1}$ belongs to the commutator semigroup $[\IH^*,\IH^*]$.
\end{exercise}

\begin{exercise} Show that for the quaternion plane $X=\IH\times\IH$, the semigroup of coscalars $[\IR_X]$ is equal to the semigroup of absolute coscalars $[\![\IR_X]\!]$ and is isomorphic to the multiplicative semigroup of non-negative real numbers. Therefore, areas and absolute areas in the quaternion plane can be measured by nonnegative real numbers. 
\end{exercise}

\section{Area in Pappian Playfair planes}

In this section we study the areols of Pappian Playfair planes. Besides the operation of multiplication by scalars, areas in Papian affine spaces admit also an operation of addition, which turns the areol of a Pappian affine plane into a one-dimensional vector space over the field of scalars.

\begin{theorem}\label{t:Pappian=>additive} The areol of any Pappian Playfair plane $X$ admits a unique addition operation $$+:[X^3]\times [X^3]\to[X^3],\quad +:(\alpha,\beta)\mapsto\alpha+\beta,$$ 
such that $[oac]=[oab]+[obc]$ for any colliner points $a,b,c\in X$ and any point $o\in X$. This addition operation together with the operation of multiplication by scalars $\cdot:\IR_X\times[X^3]\to[X^3]$ turns the areol $[X^3]$ into a one-dimensional $\IR_X$-module. 
\end{theorem}

\begin{proof} Given two  areas $\alpha,\gamma\in [X^3]$, we should define their sum $\alpha+\gamma$. If $\alpha=[0]$, then put $\alpha+\gamma\defeq\gamma$. If $\gamma=[0]$, then put $\alpha+\gamma\defeq\alpha$. So, assume that $\alpha\ne[0]\ne\gamma$. 

Choose any triangle $oab\in\alpha$. Since $X$ is Pappian, it is Desarguesian and Moufang, by Theorem~\ref{t:Hessenberg-proaffine} and Proposition~\ref{p:Desarg=>Moufang}.  By Theorem~\ref{t:area-move}, there exists a point $c\in \Aline ab$ such that $obc\in\gamma$. Put $\alpha+\gamma\defeq[oac]$. Let us show that the sum $\alpha+\gamma=[oac]$ does not depend on the choice of the triangle $oab\in \alpha$ and the point $c\in\Aline ab$. Take any triangle $o'a'b'\in\alpha$ and point $c'\in\Aline{a'}{b'}$ with  $o'b'c'\in\gamma$. Since $[oab]=\alpha=[o'a'b']$, we can apply Theorem~\ref{t:[abc]=[xyz]<=>abc=Exyz} and find an equiaffine automorphism $E$ of $X$ such that $Eo'a'b'=oab$. Consider the point $d\defeq E(c')\in E[\Aline{a'}{b'}]=\Aline {a}{b}$ and observe that $[obd]=[Eo'b'c']=[o'b'c']=\gamma=[obc]$ and hence $[bdo]=[bco]=[obc]=\gamma\ne[0]$. Theorem~\ref{t:scalar-by-area} ensures that $[bdo]=[bco]=\overvector{bcd}\cdot[bdo]$ and hence the scalar $\overvector{bcd}$ belongs to the stabilizer $\Stab([bdo])$, which is equal to the commutator group $[\IR^*_X,\IR^*_X]=\{1\}$ of the field $\IR_X$, according to Proposition~\ref{p:sign-stabilizer}. Therefore, $\overvector{bcd}=1$ and $c=d$. Then $Eo'a'c'=oad=oac$ and $[o'a'c']=[oac]$, witnessing that the sum $\alpha+\gamma$ is well-defined. 

Now take any collinear points $a,b,c$ and any point $o$ in the plane $X$. We have to show that $[oac]=[oab]+[obc]$. If $a=b=c$, then $[oac]=[0]=[0]+[0]=[oab]+[obc]$ and we are done. So, assume that $L=\overline{\{a,b,c\}}$ is a line. If $o\in L$, then $[oac]=[0]=[0]+[0]=[oab]+[obc]$. So, assume that $o\notin L$. If $a=b$, then $[oac]=[obc]=[0]+[obc]=[oab]+[obc]$. If $b=c$, then $[oac]=[oab]=[oab]+[0]=[oab]+[obc]$. So, assume that $a\ne b\ne c$. In this case, $[oab]+[obc]=[oac]$, by the definition of addition of areas.

To see that the addition operation $+$ is unique, assume that $\oplus:[X^3]\times[X^3]\to[X^3]$ is another binary operation, which is additive in the sense  that $[oac]=[oab]\oplus[obc]$ for any collinear points $a,b,c$ and any point $o$ in $X$. Given any signed areas $\alpha,\gamma\in[X^3]$, we should prove that $\alpha+\gamma=\alpha\oplus\gamma$. Fix any distinct points $o,b\in X$ and any line $L\subseteq  X$ with $b\in L$ and $o\notin L$. By Theorem~\ref{t:area-move}, there exist points $a,c\in L$ such that $oba\in -\alpha$ and $obc\in\gamma$. Then $oab\in \alpha$. The additivity of the operations $+$ and $\oplus$ ensures that $\alpha+\gamma=[oab]+[obc]=[oac]=[oab]\oplus[obc]=\alpha\oplus\gamma$.
\smallskip

It remains to show that  the set $[X^3]$ endowed with the addition $+:[X^3]\times[X^3]\to[X^3]$ and scalar multiplication $\cdot:\IR_X\times[X^3]\to[X^3]$ is a one-dimensional $\IR_X$-module (= a one-dimensional vector space over the field $\IR_X$).
Theorems~\ref{t:scalar-by-area} and \ref{t:coscalar-by-area} imply the following claim describing the properties of the operation of multiplication of areas by scalars.

\begin{claim}\label{cl:scalar-sarea-multiplication} For any scalars $s,t\in\IR_X$ and any area $\alpha\in [X^3]$, the following properties hold:
\begin{enumerate}
\item[\textup{(0)}] $s\cdot\alpha=[0]$ if and only if $s=0$ or $\alpha=[0]$;
\item[\textup{(1)}] $1\cdot\alpha=\alpha$ and $(-1)\cdot\alpha=-\alpha$;
\item[\textup{(2)}] $s\cdot(t\cdot\alpha)=(s\cdot t)\cdot\alpha$;
\item[\textup{(3)}] if $\alpha\ne[0]$, then the map \ $\cdot_\alpha:\IR_X\to[X^3]$, \ $\cdot_\alpha:x\mapsto x\cdot\alpha$, is bijective.
\end{enumerate}
\end{claim}

This claim will help us to prove the following properties of the addition operation.

\begin{claim}\label{cl:addition-area} For any scalars $s,t\in\IR_X$ and any areas $\alpha,\beta,\gamma\in[X^3]$, the following properties hold:
\begin{enumerate}
\item[\textup{(0)}] $\alpha+[0]=\alpha=[0]+\alpha$;
\item[\textup{(1)}] $(s+t)\cdot\alpha=s\cdot\alpha+t\cdot\alpha$;
\item[\textup{(2)}] $s\cdot(\alpha+\beta)=s\cdot\alpha+s\cdot\beta$;
\item[\textup{(3)}] $\alpha+(\beta+\gamma)=(\alpha+\beta)+\gamma$;
\item[\textup{(4)}] $\alpha+\beta=\beta+\alpha$;
\item[\textup{(5)}] $\alpha+(-1)\cdot\alpha=[0]$.
\end{enumerate}
\end{claim}

\begin{proof} Fix any nonzero area $\delta\in[X^3]$. By Claim~\ref{cl:scalar-sarea-multiplication}(3), there exist scalars $a,b,c\in\IR_X$ such that $\alpha=a\cdot\delta$, $\beta=b\cdot\delta$, $\gamma=c\cdot\delta$.
\smallskip

0. The equalities $\alpha+[0]=\alpha=[0]+\alpha$ follow from the definition of area addition.
\smallskip

1. If $\alpha=[0]$, then $(s+t)\cdot\alpha=(s+t)\cdot[0]=[0]=[0]+[0]=s\cdot[0]+t\cdot[0]=s\cdot\alpha+t\cdot\alpha$.
So, assume that $\alpha\ne[0]$ and choose any triangle $ouw\in\alpha$. Find unique points $x,y,z\in\Aline ou$ such that $oxu\in s$, $oyu\in t$, and $ozu\in s+t$. It follows from $\overvector{oxu}+\overvector{oyu}=s+t=\overvector{ozu}$ that  $\overvector{ox}+\overvector{oy}=\overvector{oz}=\overvector{ox}+\overvector{xz}$ and hence $\overvector{oy}=\overvector{xz}$. Then there exists a translation $T:X\to X$ such that $Toy=xz$. Consider the point $w'\defeq T(w)$ and observe that $\Aline w{w'}\subparallel \Aline ou$. Let $S$ be a hypershear of the Pappian affine plane $X$ such that $\Aline ou\subseteq\Fix(S)$ and $Sw'=w$. By Theorem~\ref{t:shear=>SAut} and Corollary~\ref{c:Thales-trans=>SAut}, the composition $ST$ is an area-preserving automorphism of $X$ such that $SToyw=xzw$ and hence $[oyw]=[xzw]$. By Theorem~\ref{t:scalar-by-area}, Proposition~\ref{p:12-permutations} and the additivity property of the area addition, $$
\begin{aligned}
s\cdot\alpha+t\cdot \alpha&=\overvector{oxu}\cdot[ouw]+\overvector{oyu}\cdot[ouw]=[oxw]+[oyw]=[oxw]+[xzw]\\
&=[wox]+[wxz]=[woz]=[ozw]=\overvector{ozu}\cdot[ouw]=(s+t)\cdot\alpha.
\end{aligned}
$$
\smallskip

2. By the preceding item, Theorem~\ref{t:scalar-by-area}(3), and distributivity of the field $\IR_X$,
$$
\begin{aligned}
s\cdot(\alpha+\beta)&=s\cdot(a\cdot\delta+b\cdot\delta)=
s\cdot((a+b)\cdot\delta)=(s\cdot(a+b))\cdot\delta=(s\cdot a+s\cdot b)\cdot\delta\\
&=(s\cdot a)\cdot \delta+(s\cdot b)\cdot \delta=s\cdot(a\cdot\delta)+s\cdot(b\cdot \delta)=s\cdot\alpha+s\cdot\beta.
\end{aligned}
$$
\smallskip

3. By the first item and the associativity-plus of the field $\IR_X$, 
$$
\begin{aligned}
\alpha+(\beta+\gamma)&=a\cdot\delta+(b\cdot\delta+c\cdot\delta)=a\cdot\delta+(b+c)\cdot\delta=(a+(b+c))\cdot\delta\\
&=((a+b)+c)\cdot\delta=(a+b)\cdot\delta+c\cdot\delta=(a\cdot\delta+b\cdot\delta)+\gamma=(\alpha+\beta)+\gamma.
\end{aligned}
$$
\smallskip

4. By the first item and  commutativity-plus of the field $\IR_X$,
$$
\alpha+\beta=a\cdot\delta+b\cdot\delta=(a+b)\cdot\delta=(b+a)\cdot\delta=b\cdot\delta+a\cdot\delta=\beta+\alpha.
$$

5. By the first item, $\alpha+(-1)\cdot\alpha=1\cdot\alpha+(-1)\cdot\alpha=(1+(-1))\cdot\alpha=0\cdot\alpha=[0].$
\end{proof}

Claims~\ref{cl:scalar-sarea-multiplication} and \ref{cl:addition-area} imply that $([X^3],+,\cdot)$ is a one-dimensional vector space over the scalar field $\IR_X$.
\end{proof}

Now we deduce some (well-known) formulas for calculating the area of a triangle in a Pappian Playfair plane.

\begin{definition}For a matrix $M=\Big(\!\begin{array}{cc}a&b\\c&d\end{array}\!\Big)$ over a field $F$, its \index{determinant}\defterm{determinant} $\Big|\!\begin{array}{cc}a&b\\c&d\end{array}\!\Big|$ is the element $ad-bc$ of the field $F$.
\end{definition}

\begin{proposition}\label{p:sarea-vectors} Let $(X,uow)$ be a based Pappian Playfair plane. For any vectors $\vec x,\vec y\in\overvector X$, their area form $[\vec x,\vec y]$ can be calculated by the formula
$$[\vec x,\vec y]=(x_1y_2-x_2y_1)\cdot[uow]=\left|\!\begin{array}{cc}x_1&x_2\\y_1&y_2\end{array}\!\right|\cdot[uow],$$ where $x_1,x_2,y_1,y_2\in\IR_X$ are unique scalars such that $\vec x=x_1{\cdot}\overvector{ou}+x_2{\cdot}\overvector{ow}$ and $\vec y=y_1{\cdot}\overvector{ou}+y_2{\cdot}\overvector{ow}$.
\end{proposition}

\begin{proof} Let $e\in X$ be the diunit of the affine base $uow$ and $\Delta=\Aline oe$ be the ternar of the based affine plane $(X,uow)$. Since the affine plane $X$ is Pappian, the ternar $\Delta$ is linear, distributive, associative, and commutative, according to Theorem~\ref{t:p.commutative}. Then $\Delta$ endowed with the plus and dot operations is a field. By Corollary~\ref{c:R=RPi}, the map $I:\Delta\to \IR_X$, $I:d\mapsto \overvector{ode}$, is an isomorphism of the fields $(\Delta,+,\cdot)$ and $\IR_X$. 

Find unique points $x,y\in X$ such that $\vec x=\overvector{ox}$ and $\vec y=\overvector{oy}$. The definition of the area form ensures that $[\vec x,\vec y]=[oxy]$. Let $(x',x''),(y',y'')\in\Delta^2$ be the coordinates of the points $x,y$ in the affine base $uow$. By Corollary~\ref{c:coordinates=vectors}, $$x_1=\overvector{ox'e}=I(x'),\quad x_2=\overvector{ox''e}=I(x''),\quad y_1=\overvector{oy'e}=I(y'),\quad y_2=\overvector{oy''e}=I(y'').$$ 

By Proposition~\ref{p:area-via-coordinates},  $[oxy]=[oud]$, where $d\defeq y''{\cdot}x'-y'{\cdot}x''\in\Delta$. Since the map $I:\Delta\to\IR_X$ is an isomorphism of the fields $\Delta$ and $\IR_X$, $\overvector{ode}=I(d)=I(y''){\cdot}I(x')-I(y'){\cdot}I(x'')=y_2{\cdot}x_1-y_1{\cdot}x_2$. Let $S:X\to X$ be a shear of the plane $X$ such that $Soue=ouw$. Consider the point $\delta\defeq S(d)\in\Aline ow$ and observe that $[ou\delta]=[oud]=[oxy]$ and $\overvector{o\delta w}=\overvector{ode}=y_2{\cdot}x_1-y_1{\cdot}x_2$. By the defining property of multiplication of an area by a scalar,
$(y_2{\cdot}x_1-y_1{\cdot}x_2)\cdot[owu]=\overvector{o\delta w}\cdot [owu]=[o\delta u]$, which implies
$$
\begin{aligned}
[\vec x,\vec y]&=[oxy]=[oud]=[ou\delta]=-[o\delta u]\\
&=-(y_2{\cdot}x_1-y_1{\cdot}x_2)\cdot[owu]=(y_2{\cdot}x_1-y_1{\cdot}x_2)\cdot[ouw]=\left|\!\begin{array}{cc}x_1&x_2\\y_1&y_2\end{array}\!\right|\cdot[uow].
\end{aligned}
$$
\end{proof}

Proposition~\ref{p:sarea-vectors} will help us to prove that the  area form of a Pappian Playfair  plane is bilinear and alternative.

\begin{proposition} The area form $[\cdot\,,\cdot]:\overvector X\times\overvector X\to[X^3]$ of a Pappian Playfair plane $X$ has the following properties, for all scalars $t\in\IR_X$ and vectors $\vec x,\vec y,\vec z\in\overvector X$:
\begin{enumerate}
\item[\textup{(1)}] $[t{\cdot}\vec x,\vec y]=t\cdot[\vec x,\vec y]=[\vec x,t{\cdot}\vec y]$;
\item[\textup{(2)}] $[\vec x,\vec x]=[0]$;
\item[\textup{(3)}] $[\vec y,\vec x]=-[\vec x,\vec y]=(-1)\cdot[\vec x,\vec y]$;
\item[\textup{(4)}] $[\vec x+\vec y,\vec z]=[\vec x,\vec z]+[\vec y,\vec z]$;
\item[\textup{(5)}] $[\vec x,\vec y+\vec z]=[\vec x,\vec y]+[\vec x,\vec z]$.
\end{enumerate}
\end{proposition}

\begin{proof} By Theorem~\ref{t:Hessenberg-proaffine} and Proposition~\ref{p:Desarg=>Moufang}, the Pappian Playfair plane is Desarguesian and Moufang. The items (1)--(3) follow from Propositions~\ref{p:areaform-T} and \ref{p:area-form1}.
\smallskip

4. Fix any affine base $uow$ in the affine plane $X$ and find unique scalars $x_1,x_2,y_1,y_2,z_1,z_2\in\IR_X$ such that $\vec x=x_1\cdot\overvector{ou}+x_2\cdot\overvector{ow}$, $\vec y=y_1\cdot\overvector{ou}+y_2\cdot\overvector{ow}$, and $\vec z=z_1\cdot\overvector{ou}+z_2\cdot\overvector{ow}$. Then $\vec x+\vec y=(x_1+y_1)\cdot\overvector{ou}+(x_2+y_2)\cdot\overvector{ow}$. Applying Proposition~\ref{p:sarea-vectors} and Claim~\ref{cl:addition-area}(1), we conclude that
$$
\begin{aligned}
[\vec x+\vec y,\vec z]&=\big((x_1+y_1){\cdot}z_2-(x_2+y_2){\cdot}z_1\big)\cdot[ouw]=\big((x_1{\cdot} z_2-x_2{\cdot}z_1)+(y_1{\cdot}z_2-y_2{\cdot}z_1)\big)\cdot[ouw]\\
&=(x_1{\cdot} z_2-x_2{\cdot}z_1)\cdot[uow]+(y_1{\cdot}z_2-y_2{\cdot}z_1)\cdot[ouw]=[\vec x,\vec z]+[\vec y,\vec z].
\end{aligned}
$$ 
\smallskip

5. Applying the items (3) and (4), we obtain the equalities
$$[\vec x,\vec y+\vec z]=(-1)\cdot [\vec y+\vec z,\vec x]=(-1)\cdot([\vec y,\vec x]+[\vec z,\vec x])=(-1)\cdot[\vec y,\vec x]+(-1)\cdot[\vec z,\vec x]=[\vec x,\vec y]+[\vec x,\vec z]$$
for all vectors $\vec x,\vec y,\vec z\in \overvector X$.   
\end{proof}

\begin{definition} For a matrix $\left(\!\begin{array}{@{\,}c@{\;}c@{\;}c@{\,}}a_{11}&a_{12}&a_{13}\\a_{21}&a_{22}&a_{23}\\
a_{31}&a_{32}&a_{33}\end{array}\!\right)$ over a field $F$, its \index{determinant}\defterm{determinant} $\left|\!\begin{array}{@{\,}c@{\;}c@{\;}c@{\,}}a_{11}&a_{12}&a_{13}\\a_{21}&a_{22}&a_{23}\\
a_{31}&a_{32}&a_{33}\end{array}\!\right|$ equals the element $a_{11}a_{22}a_{33}+a_{12}a_{23}a_{31}+a_{13}a_{21}a_{32}-a_{11}a_{23}a_{32}-a_{12}a_{21}a_{33}-a_{13}a_{22}a_{31}$ of the field $F$.
\end{definition}

\begin{theorem}\label{t:sarea-oxyz} Let $(X,uow)$ be a based Pappian Playfair plane. For any triple $xyz\in X^3$, its area $[xyz]$ can be calculated by the formula
$$[xyz]=
\left|\begin{array}{@{\;}c@{\;\;}c@{\;\;}c@{\;}}x_1&x_2&1\\y_1&y_2&1\\z_1&z_2&1\end{array}\right|\cdot[ouw]=[oxy]+[oyz]+[ozx],$$ where $(x_1,x_2), (y_1,y_2),(z_1,z_2)$ are the coordinates of the points $x,y,z$ in the affine base $uow$.
\end{theorem}

\begin{proof} Consider the vectors $\overvector {xy}$, $\overvector{xz}$, and observe that $$\overvector {xy}=\overvector{oy}-\overvector{ox}=(y_1\cdot\overvector{ou}+y_1\cdot\overvector{ow})-(x_1\cdot\overvector{ou}+x_2\cdot\overvector{ow})=(y_1-x_1)\cdot\overvector{ou}+(y_2-x_2)\cdot\overvector{ow}.$$ and $\overvector{xz}=(z_1-x_1)\cdot\overvector{ou}+(z_2-x_2)\cdot\overvector{ow}$. By the definition of the area form and Proposition~\ref{p:sarea-vectors},
$$
\begin{aligned}
[xyz]&=[\overvector{xy},\overvector{xz}]=\big((y_1-x_1)\cdot(z_2-x_2)-(z_1-x_1)\cdot(y_2-x_2)\big)\cdot[ouw]\\
&=(x_1{\cdot}y_2-x_2{\cdot}y_1+y_1{\cdot}z_2-z_1{\cdot}y_2+z_1{\cdot}x_2-x_1{\cdot}z_2)\cdot[ouw]=
\left|\!\begin{array}{@{\;}c@{\;\;}c@{\;\;}c@{\;}}x_1&x_2&1\\y_1&y_2&1\\z_1&z_2&1\end{array}\!\right|\cdot[ouw].
\end{aligned}
$$
On the other hand,
$$
\begin{aligned}
[xyz]&=(x_1{\cdot}y_2-x_2{\cdot}y_1+y_1{\cdot}z_2-z_1{\cdot}y_2+z_1{\cdot}x_2-x_1{\cdot}z_2)\cdot[ouw]\\
&=(x_1{\cdot}y_2-x_2{\cdot}y_1)\cdot[ouw]+(y_1{\cdot}z_2-z_1{\cdot}y_2)\cdot[ouw]+(z_1{\cdot}x_2-x_1{\cdot}z_2)\cdot[ouw]\\
&=[\overvector{ox},\overvector{oy}]+[\overvector{oy},\overvector{oz}]+[\overvector{oz},\overvector{ox}]=[oxy]+[oyz]+[oxz],
\end{aligned}
$$
by Proposition~\ref{p:sarea-vectors} and the definition of the area form.
\end{proof}

\section{Area of polygons in Pappian Playfair planes}

In this section we use the additivity of the  area in Pappian affine planes and will define the  area for all polygons in Pappian Playfair planes.

\begin{definition} By a \index{polygon}\defterm{polygon} in a plane $X$ we understand any ordered sequence $a_0a_1\cdots a_n$ of points of the plane $X$. The points $a_0,a_1,\dots,a_n$ are called the \index{polygon!vertices of}\defterm{vertices} of the polygon and the lines $\Aline {a_0}{a_1},\Aline{a_1}{a_2},\dots,\Aline{a_{n-1}}{a_n},\Aline{a_n}{a_0}$ are called the \index{polygon!sides of}\defterm{sides} of the polygon $a_0\cdots a_n$.
\end{definition}

\begin{definition} The \index[note]{$[a_0a_1\dots a_n]$}\index{polygon!area of}\defterm{area} $[a_0a_1\dots a_n]$ of a polygon $a_0a_1\dots a_n$ in a Pappian affine plane is defined as the sum
$$\sum_{i=1}^{n-1}[a_0a_ia_{i+1}]$$of areas of the triangles composing the polygon.
\end{definition}

\begin{proposition}\label{p:sarea-polygon} For any polygon $a_0a_1\dots a_n$ in a Pappian affine plane $\Pi$ and any point $o\in\Pi$ we have the equality 
$$[a_0a_1\cdots a_n]=\sum_{i=0}^n[oa_ia_{i+1}]$$
where $a_{n+1}\defeq a_0$.
\end{proposition}

\begin{proof} Applying Theorem~\ref{t:sarea-oxyz}, we conclude that
$$
\begin{aligned}
[a_0a_1\dots a_n]&\defeq\sum_{i=1}^{n-1}[a_0a_ia_{i+1}]=\sum_{i=0}^n[a_0a_ia_{i+1}]=\sum_{i=0}^n([oa_0a_i]+[oa_ia_{i+1}]+[oa_{i+1}a_0])\\
&=
\sum_{i=0}^n[oa_ia_{i+1}]+\sum_{i=0}^n[oa_0a_i]+\sum_{i=1}^{n+1}[oa_ia_0]=\sum_{i=0}^n[oa_ia_{i+1}]+\sum_{i=1}^n[oa_0a_i]+\sum_{i=1}^{n}[oa_ia_0]\\
&=\sum_{i=0}^n[oa_ia_{i+1}]+\sum_{i=1}^n([oa_0a_i]+[oa_ia_0])=\sum_{i=0}^n[oa_ia_{i+1}]+\sum_{i=1}^n[oa_0a_0]\\
&=\sum_{i=0}^n[oa_ia_{i+1}]+[0]=\sum_{i=0}^n[oa_ia_{i+1}].
\end{aligned}
$$
\end{proof}

The area of polygons in a Pappian Playfair plane is additive in the following sense.

\begin{proposition} For any polygon $a_0\cdots a_m$ in a Pappian Playfair plane and any positive number $n<m$, 
$$[a_0a_1\dots a_m]=[a_0\dots a_n]+[a_n\dots a_m,a_0].$$
\end{proposition}

\begin{proof} Put $a_{m+1}=a_0$ and $a_{m+2}=a_n$. It follows from $a_0=a_{m+1}$ that $[a_0a_ia_{i+1}]=[0]$ for all $i\in\{m,m+1\}$. By Proposition~\ref{p:sarea-polygon},
$$
\begin{aligned}
[a_0a_1\cdots a_m]&=\sum_{i=1}^{m-1}[a_0a_ia_{i+1}]=\sum_{i=1}^{m+1}[a_0a_ia_{i+1}]=\sum_{i=1}^{n-1}[a_0a_ia_{i+1}]+\sum_{i=n}^{m+1}[a_0a_ia_{i+1}]\\
&=[a_0\dots a_n]+[a_n\cdots a_m,a_{m+1}].
\end{aligned}
$$
\end{proof}

Propositions~ \ref{p:sarea-polygon} and \ref{p:sarea-vectors} imply the following well known shoelace formula, first proved by Carl Friedrich Gauss (but used long before by surveyors for land measurement).

\begin{corollary}[Gauss, 1795] Every polygon $a_0a_1\dots a_n$ in a based Pappian Playfair  plane $(X,uow)$ has area $$[a_0a_1\dots a_n]=\sum_{i=0}^n[oa_ia_{i+1}]=\Big(\sum_{i=0}^n a_i'a_{i+1}''-a_i''a_{i+1}''\Big)\cdot[ouw],$$
where $a_{n+1}=a_0$ and $(a_i',a_i'')$ are coordinates of the point $a_i$ in the affine base $uow$.
\end{corollary} 

\section{Area-additive liners}


\rightline{\em If equals are added to equals, the wholes are equal.}

\rightline{Euclid, ``Elements'', Book I, Common Notion II}
\bigskip

An instance of Euclid’s Common Notion II (quoted in the epigraph to this section) is the {\sf Area Additivity Axiom}, introduced in the following definition.

\begin{definition} A liner $X$ is called \index{area-additive liner}\index{liner!area-additive}\defterm{area-additive} if it satisfies the \index{Area Additivity Axiom}\defterm{Area Additivity Axiom}
\begin{itemize}
\item[{\defterm{\sf (AAA)}}] for any lines $L,L'\subseteq X$ and points $a,b,c\in L$, $a',b',c'\in L'$, $o\in X\setminus L$, $o'\in X\setminus L'$,  the area equalities $[oab]=[o'a'b']$ and $[obc]=[o'b'c']$ imply $[oac]=[o'a'c']$.
\end{itemize}
\begin{picture}(200,85)(-120,-10)
{\color{gray}
\polygon*(0,10)(30,40)(30,10)
\polygon*(150,0)(165,60)(165,0)
\color{lightgray}
\polygon*(30,10)(30,40)(70,10)
\polygon*(165,0)(165,60)(185,0)
}

\put(-10,10){\line(1,0){90}}
\put(85,7){$L$}
\put(0,10){\line(1,1){30}}
\put(70,10){\line(-4,3){40}}
\put(0,10){\circle*{3}}
\put(-2,2){$a$}
\put(30,10){\circle*{3}}
\put(26,1){$b$}
\put(70,10){\circle*{3}}
\put(68,2){$c$}
\put(30,40){\circle*{3}}
\put(28,43){$o$}

\put(150,0){\line(1,0){35}}
\put(150,0){\line(1,4){15}}
\put(185,0){\line(-1,3){20}}
\put(140,0){\line(1,0){55}}
\put(200,-3){$L'$}
\put(150,0){\circle*{3}}
\put(147,-9){$a'$}
\put(165,0){\circle*{3}}
\put(162,-10){$b'$}
\put(185,0){\circle*{3}}
\put(183,-9){$c'$}
\put(165,60){\circle*{3}}
\put(163,63){$o'$}
\end{picture}
\end{definition}

\begin{definition} An \index{area addition}\defterm{area addition} on a liner $X$ is a binary operation $$+:[X^3]\times [X^3]\to[X^3],\quad +:(\alpha,\beta)\to \alpha+\beta,$$which is \index{area addition!additive}\defterm{additive} in the sense that for any line $L\subseteq X$ and points $a,b,c\in L$, $o\in X$, the equality $[oab]+[obc]=[oac]$ holds.
\end{definition}

Let us recall that a liner $X$ is \index{area-compatible liner}\index{liner!area-compatible}\defterm{area-compatible} if for any areas $\alpha,\beta\in [X^3]$, there exist a line $L\subseteq X$ and points $a,b,c\in L$ and $o\in X$ such that $[oab]=\alpha$ and $[obc]=\beta$. By Proposition~\ref{p:area-surjective-around}, every hyper-Bolyai plane is area-compatible; in particular, every Playfair plane is area-compatible.

\begin{theorem}\label{t:area-additive<=>area-addition} 
A liner $X$ is area-additive (and area-compatible)  if and only if it admits a (unique) area addition $+:[X^3]\times[X^3]\to[X^3]$. 
\end{theorem}

\begin{proof} Assume that a liner $X$ is area-additive. Given two areas $\alpha,\beta\in[X^\vartriangle]$, define their sum $\alpha+\beta$ as follows.
 If there exists a line $L\subseteq X$ and points $a,b,c\in L$ and $o\in X\setminus L$ such that $[oab]=\alpha$ and $[obc]=\beta$, then put $\alpha+\beta\defeq[oac]$. Otherwise, put $\alpha+\beta\defeq[0]$. The Area Addition Axiom ensures that the sum $\alpha+\beta$ is well-defined and does no depend on the choice of the line $L$ and points $a,b,c,o$. Moreover, the  definition of the operation $+$ ensures that it is additive and hence is an area addition on $X$.
 
Assuming that the area-additve liner $X$ is area-compatible, we shall prove that the area addition $+:[X^\vartriangle]\times[X^\vartriangle]\to[X^\vartriangle]$ is unique. Assume that $\oplus:[X^\vartriangle]\times[X^\vartriangle]\to[X^\vartriangle]$ is another area addition on $X$. By the area-compatibility of $X$, for any areas $\alpha,\beta\in[X^3]$, there exist a line $L\subseteq X$ and points $a,b,c\in L$ and $o\in X$ such that $[oab]=\alpha$ and $[obc]=\beta$. The additivity of the area additions $+$ and $\oplus$ ensures that $\alpha+\beta=[oab]+[obc]=[oac]=[oab]\oplus[obc]=\alpha\oplus\beta$, witnessing that the operations $+$ and $\oplus$ coincide, and the area addition $+$ on $X$ is unique. 
\smallskip

Next, assume that a liner $X$ admits an area addition $+:[X^3]\times[X^3]\to[X^3]$. To show that $X$ satisfies the Area Addition Axiom, choose any lines $L,L'\subseteq X$ and points $a,b,c\in L$, $a',b',c'\in L'$ and $o,o'\in X$ such that $[oab]=[o'a'b']$ and $[obc]=[o'b'c']$. Consider the areas $\alpha\defeq[oab]$, $\beta\defeq[obc]$, and their sum $\alpha+\beta$. The additivity of the area addition ensures that $[oac]=[oab]+[o'a'b']=\alpha+\beta=[o'a'b']+[o'b'c']=[o'a'c']$, witnessing that the Area Additivity Axiom holds.

Assuming that the area addition $+$ is unique, we will show that the liner $X$ is area-compatible. In the oppsite case, we could find two areas $\alpha_0,\beta_0\in[X^\vartriangle]$ such that for all lines $L\subseteq X$ and points $a,b,c\in L$, $o\in X$, we have $([oab],[obc])\ne(\alpha_0,\beta_0)$. Since $|[X^3]|\ge 2$, there exists an area $\gamma_0\in [X^3]$ such that $\gamma_0\ne \alpha_0+\beta_0$. We claim that the binary operation $\oplus:[X^3]\times[X^3]\to[X^3]$ defined by
$$\alpha\oplus\beta=\begin{cases}\gamma_0&\mbox{if $(\alpha,\beta)=(\alpha_0,\beta_0)$};\\
\alpha+\beta&\mbox{otherwise};
\end{cases}
$$is additive. Indeed, take any line $L\subset X$ and points $a,b,c\in L$ and $o\in X\setminus L$. The choice of the areas $\alpha_0,\beta_0$ ensures that $([oab],[obc])\ne(\alpha_0,\beta_0)$. The definition of the binary operation $\oplus$ and the additivity of the area addition $+$ ensures that $[oab]\oplus[obc]=[oab]+[obc]=[oac]$, witnessing that $\oplus$ is additive and hence is an area-addition on $X$, distinct from the area addition $+$, which contradicts the uniqueness of the area addition $+$.
\end{proof}

We recall that a liner $X$ is called {\em area-bijective} (resp. {\em area-surjective}, {\em area-injective}) if for every line $\Lambda\subseteq X$ and points $\lambda\in\Lambda$ and $o\in X\setminus\Lambda$, the function $\Lambda\to[X^3]$, $x\mapsto [o\lambda x]$, is bijective (resp. surjective, injective). By Proposition~\ref{p:area-surjective-around},
\begin{itemize}
\item every hyper-Bolyai plane is area-surjective;
\item every area-surjective liner is area-compatible; 
\item every area-compatible liner is ranked of rank $\le 3$.
\end{itemize}

\begin{proposition}\label{p:area-additive=>area-injective} Any area-additive liner is area-injective.
\end{proposition}

\begin{proof} Assume that a liner $X$ is area-additive and let $+:[X^3]\times[X^3]\to [X^3]$ be an area addition on $[X^3]$. To prove that $X$ is area-injective, we should check that for every line $\Lambda$ in $X$ and points $\lambda\in \Lambda$ and $o\in X\setminus\Lambda$, the function  $h:\Lambda\to[X^3]$, $h:x\mapsto [o\lambda x]$, is injective. Take any points $x,y\in \Lambda$ with $[o\lambda x]=h(x)=h(y)=[o\lambda y]$. The additivity of the area addition $+$ ensures that 
$$[0]=[oxx]=[ox\lambda]+[o\lambda x]=[ox\lambda]+[o\lambda y]=[oxy]$$
and hence $\|\{o,x,y\}\|=2$. Then $y\in\Aline ox\cap\Lambda=\{x\}$ and hence $y=x$, witnessing that the function $h$ is inejctive and the liner $X$ is area-injective.
\end{proof}

Proposition~\ref{p:area-additive=>area-injective} implies the following corollary.

\begin{corollary} An area-additive liner is area-surjective if and only if it is area-bijective.
\end{corollary}

By Theorem~\ref{t:area-additive<=>area-addition}, every (area-compatible) area-additive liner $X$ admits a (unique) additive binary operation $+:[X^3]\times [X^3]\to[X^3]$, which turns the areol $[X^3]$ into a magma. If $X\ne\varnothing$, then this magma is unital with neutral element $[0]$.

\begin{theorem}\label{t:Playfair=>areol-commutative} If an area-additive liner $X$ is  area-surjective (and Playfair), then its areol magma $([X^3],+)$ is a (commutative) group.
\end{theorem}  

\begin{proof} Let $X$ be an area-surjective area-additive liner. By Proposition~\ref{p:area-surjective-around}, the area-surjective liner $X$ is area-compatible, and by Theorem~\ref{t:area-additive<=>area-addition}, the area-compatible area-additive liner $X$ admits a unique area addition $+:[X^3]\times[X^3]\to[X^3]$. We shall prove that the area addition $+$ satisfies the standard axioms of a group:
\begin{enumerate}
\item $\forall \alpha,\beta,\gamma\in [X^3]\quad (\alpha+\beta)+\gamma=\alpha+(\beta+\gamma)$,
\item $\forall \alpha\in [X^3]\quad \alpha+[0]=\alpha=[0]+\alpha$,
\item $\forall \alpha\in[X^3]\;\exists \beta\in [X^3]\quad \alpha+\beta=[0]=\beta+\alpha$.
\end{enumerate}
Fix any line $\Lambda\subset X$ and points $\lambda\in \Lambda$ and $o\in X\setminus \Lambda$.
\smallskip

1. Given any areas $\alpha,\beta,\gamma\in [X^3]$, we can apply the area-surjectivity of $X$ and find points $a,b,c\in \Lambda$ such that $[o\lambda a]=\alpha$, $[oab]=\beta$ and $[obc]=\gamma$. The additivity of the area addition ensures that 
$$
\begin{aligned}
(\alpha+\beta)+\gamma&=([voa]+[vab])+[vbc]=[vob]+[vbc]=[voc]\\
&=[voa]+[vac]=\alpha+([vab]+[vac])=\alpha+(\beta+\gamma),
\end{aligned}
$$
witnessing that the area addition is associative.
\smallskip

2,3. Given any area $\alpha\in[X^3]$, we can apply the area-surjectivity of $X$ and find a point $a\in\Lambda$ such that $[o\lambda a]=\alpha$. Consider the area $\mbox{-}\alpha\defeq[oa\lambda]$. The additivity of the area addition ensures that $\alpha+[0]=[o\lambda a]+[oaa]=[v\lambda a]=\alpha$ and $[0]+\alpha=[o\lambda\lambda]+[o\lambda a]=[o\lambda a]=\alpha$, 
witnessing that the zero area $[0]$ in a neutral element of the area addition.
Also, the additivity of the area addition ensures that
$\alpha+(\mbox{-}\alpha)=[o\lambda a]+[oa\lambda]=[o\lambda\lambda]=[0]=[oaa]=[oa\lambda]+[o\lambda a]=(\mbox{-}\alpha)+\alpha$, 
witnessing that the opposite area $\mbox{-}\alpha$ is inverse to the area $\alpha$ with respect to the area addition. 
\smallskip 

Therefore, the  magma $([X^3],+)$ is a group. Assuming that the liner $X$ is Playfair, we shall prove that this group is commutative. Given any elements $\alpha,\beta\in [X^3]$, we need to check that $\alpha+\beta=\beta+\alpha$. If $\alpha$ or $\beta$ equals $[0]$, then the equality $\alpha+\beta=\beta+\alpha$ follows from the properties of the neutral element $[0]$ of the group $([X^3],+)$. So, assume that $\alpha\ne [0]\ne\beta$.

Choose any line $\Lambda$ in $X$ and points $\lambda\in L$ and $o\in X\setminus \Lambda$. Since $X$ is Playfair, there exists a line $\Lambda'$ in $X$ such that $o\in \Lambda'$ and $\Lambda'\parallel \Lambda$. Since $X$ is area-surjective, there exist points $a,b\in \Lambda$ such that $[o\lambda a]=\alpha$ and $[oab]=\beta$. The additivity of the area addition ensures that $[o\lambda b]=[v\lambda a]+[oab]=\alpha+\beta$. If $b=\lambda$, then $\alpha+\beta=[o\lambda b]=[obb]=[0]$ and $(\beta+\alpha)+\beta=\beta+(\alpha+\beta)=\beta+[0]=[0]+\beta$, which implies $\beta+\alpha=[0]=\alpha+\beta$, by the cancellativity of the group operation $+$.  So, we assume that $b\ne o$. It follows from $[o\lambda a]=\alpha\ne[0]\ne\beta=[oab]$ that $\lambda\ne a\ne b$. 

Since $X$ is Playfair, there exist unique points $u,w\in \Lambda'$ and $d\in \Lambda$ such that $\Aline bu\parallel \Aline ao$ and $\Aline ud\parallel \Aline o\lambda\parallel\Aline wb$. Since the Playfair liner $X$ is Proclus and $3$-ranked, the parallelity relations $\Aline ud\parallel \Aline o\lambda\parallel\Aline wb$ in the plane $\overline{\Lambda\cup\Lambda'}$ imply $\Aline ud\parallel \Aline wb$.

\begin{picture}(100,70)(-150,-15)

{\linethickness{0.7pt}
\put(-15,0){\color{teal}\line(1,0){90}}
\put(-15,40){\color{teal}\line(1,0){90}}
\put(80,-3){\color{teal}$\Lambda$}
\put(80,37){\color{teal}$\Lambda'$}
\put(0,0){\color{cyan}\line(0,1){40}}
\put(40,0){\color{cyan}\line(0,1){40}}
\put(60,0){\color{cyan}\line(0,1){40}}
\put(0,40){\color{violet}\line(1,-2){20}}
\put(40,40){\color{violet}\line(1,-2){20}}
}
\put(0,40){\line(3,-2){60}}
\put(0,40){\line(1,-1){40}}

\put(0,0){\circle*{2.5}}
\put(-4,-9){$\lambda$}
\put(20,0){\circle*{2.5}}
\put(18,-9){$a$}
\put(40,0){\circle*{2.5}}
\put(38,-9){$d$}
\put(60,0){\circle*{2.5}}
\put(57,-9){$b$}
\put(0,40){\circle*{2.5}}
\put(-3,43){$o$}
\put(40,40){\circle*{2.5}}
\put(37,43){$u$}
\put(60,40){\circle*{2.5}}
\put(57,43){$w$}
\end{picture}

The additivity of the area addition ensures that $[o\lambda b]=[o\lambda a]+[oab]$ and $[bwo]=[bwu]+[buo]$. Since $\lambda owb$ and $aoub$ are parallelograms, we can apply Proposition~\ref{p:12-permutations} and conclude that $[o\lambda b]=[bwo]$ and $[oab]=[buo]$. Then $[bwu]+[buo]=[bwo]=[o\lambda b]=[o\lambda a]+[oab]=[o\lambda a]+[buo]$. Since $([X^3],+)$ is a group, the equality $[bwu]+[buo]=[o\lambda a]+[buo]$ implies $[bwu]=[o\lambda a]$. Since $duwb$ is a parallelogram, Proposition~\ref{p:12-permutations} implies $[bwu]=[udb]$ and hence $[odb]=[udb]=[bwu]=[o\lambda a]=\alpha$. Observe that $buo\boxtimes duo$. Since $\lambda o ud$ is a parallelogram, Proposition~\ref{p:12-permutations} implies $[o\lambda d]=[duo]=[buo]=[oab]=\beta$. The additivity of the area addition ensures that $$\alpha+\beta=[o\lambda a]+[oab]=[o\lambda b]=[o\lambda d]+[odb]=\beta+\alpha,$$ witnessing that the group $([X^3],+)$ is commutative.
\end{proof} 

The main result of this section is the following theorem.

\begin{theorem}\label{t:area-additive<=>Pappian} A Playfair plane $X$ is area-additive if and only if it is Pappian.
\end{theorem}

\begin{proof} The ``if'' part follows from Theorem~\ref{t:Pappian=>additive}. To prove the ``if'' part, assume that $X$ is an area-additive Playfair plane. By Propositions~\ref{p:hyper-Bolyai-interplay}, \ref{p:area-surjective-around}, \ref{p:area-additive=>area-injective}, the Playfair plane is hyper-Bolyai, area-surjective, area-compatible, ranked, area-injective, and area-bijective. By Theorem~\ref{t:area-additive<=>area-addition}, $X$ admits a unique area addition $+:[X^3]\times[X^3]\to[X^3]$. By Theorem~\ref{t:Playfair=>areol-commutative}, the magma $([X^3],+)$ is a commutative group with neutral element $[0]$. 

\begin{lemma} For any elements $o,a,b,c\in X$ we have the equality $[abc]=[oab]+[obc]+[oca]$.
\end{lemma}

\begin{proof} If $o=c$, then $[oab]+[obc]+[oca]=[cab]+[0]+[0]=[cab]=[abc]$, by Proposition~\ref{p:12-permutations}. By analogy we can show that $[oab]+[obc]+[oca]=[abc]$ whenever $o$ equals $a$ or $b$. So, assume that $o\notin\{a,b,c\}$.

If there exists a point $d\in\Aline oa\cap\Aline bc$, then the additivity of the area addition and  Proposition~\ref{p:12-permutations} ensure that 
$$
\begin{aligned}
[oab]+[obc]+[oca]&=[oab]+[obd]+[odc]+[oca]=[boa]+[bdo]+[cod]+[cao]\\
&=[bdo]+[boa]+[cao]+[cod]=[bda]+[cad]=[abd]+[adc]=[abc].
\end{aligned}
$$
By analogy we can show that $[oab]+[obc]+[oca]=[abc]$ whenever $\Aline ob\cap\Aline ac$ or $\Aline oc\cap\Aline ab$ is nonempty.

So, asume that $\Aline oa\cap\Aline bc=\Aline ob\cap\Aline ac=\Aline oc\cap\Aline ab=\varnothing$, which means that $oabc$ is a Boolean parallelogram in the ranked liner $X$. By Proposition~\ref{p:12-permutations}, all triangles with vertices in the set $\{o,a,b,c\}$ have the same area. In particular, $[oab]=[abc]$ and $[oca]=[ocb]$, which implies $$[oab]+[obc]+[oca]=[abc]+[obc]+[ocb]=[abc]+[obb]=[abc]+[0]=[abc].$$
\end{proof}

\begin{lemma} The area-additive Playfair liner $X$ is Moufang.
\end{lemma}

\begin{proof} By Theorem~\ref{t:affine-Moufang<=>}, it suffices to show that $X$ is para-Desarguesian.
To check this property, fix any distinct parallel lines $A,B,C,D$ in $X$, and distinct points $a,a'\in A$, $b,b'\in B$, $c,c'\in C$, and $x\in\Aline ab\cap\Aline{a'}{b'}\cap D$, $y\in \Aline bc\cap\Aline {b'}{c'}\cap D$. It follows from $x\ne y$ that $\Pi\defeq\overline{A\cup B\cup C\cup D}=\overline{\{a,b,c\}}$ is a plane. We have to prove that $\Aline ac\cap\Aline{a'}{c'}\subseteq D$. Since the Playfair plane $\Pi$ is Proclus, there exists a point $z\in D\cap\Aline ac$. Since $\Pi$ is Playfair, there exists a unique point $d\in D$ such that $\Aline bd\parallel \Aline az$.

\begin{picture}(300,150)(-50,-15)
\put(0,0){\color{teal}\line(1,0){250}}
\put(255,-3){\color{teal}$D$}
\put(0,60){\color{teal}\line(1,0){250}}
\put(255,57){\color{teal}$C$}
\put(0,90){\color{teal}\line(1,0){250}}
\put(255,87){\color{teal}$B$}
\put(0,120){\color{teal}\line(1,0){250}}
\put(255,117){\color{teal}$A$}

\put(0,0){\line(1,1){120}}
\put(0,0){\line(2,1){240}}
\put(120,0){\color{cyan}\line(0,1){120}}
\put(120,0){\color{red}\line(1,1){120}}
\put(90,0){\color{cyan}\line(0,1){90}}
\put(180,0){\line(-1,1){90}}
\put(90,0){\color{red}\line(1,1){90}}
\put(180,0){\line(0,1){90}}
\put(90,0){\color{lightgray}\line(1,4){30}}
\put(90,0){\color{lightgray}\line(5,4){150}}
\put(120,0){\color{lightgray}\line(-1,3){30}}
\put(120,0){\color{lightgray}\line(2,3){60}}
\put(90,0){\color{lightgray}\line(1,2){30}}
\put(90,0){\color{lightgray}\line(3,2){90}}

\put(0,0){\circle*{3}}
\put(-4,-8){$x$}
\put(90,0){\circle*{3}}
\put(88,-9){$d$}
\put(120,0){\circle*{3}}
\put(118,-8){$z$}
\put(180,0){\circle*{3}}
\put(178,-8){$y$}
\put(90,90){\circle*{3}}
\put(85,92){$b$}
\put(180,90){\circle*{3}}
\put(178,93){$b'$}
\put(120,120){\circle*{3}}
\put(118,123){$a$}
\put(240,120){\circle*{3}}
\put(240,123){$a'$}
\put(120,60){\circle*{3}}
\put(122,65){$c$}
\put(180,60){\color{red}\circle*{3}}
\put(181,66){$c'$}

\end{picture}

It follows from $\Aline bd\parallel \Aline az$ that $zda\boxtimes zba$. 
The additivity of the area addition ensures that $[zxb]+[zba]=[zxa]=[zxa']=[zxb']+[zb'a']=[zxb]+[zb'a']$ and hence $[zba]=[zb'a']$, by the cancellativity of the group operation $+$ in the areol $[X^3]$. On the other hand, $zba\boxtimes zda\boxtimes zda'$. Therefore, $[zda']=[zba]=[zb'a']$. By Proposition~\ref{p:same-area<=>subparallel}, the equality $[zda']=[zb'a']$ implies $\Aline d{b'}\parallel \Aline z{a'}$.

On the other hand, the additivity of the area addition ensures that $[dbc]+[dcy]=[dby]=[db'y]=[db'c']+[dc'y]=[db'c']+[dcy]$ and hence $[dbc]=[db'c']$, by the cancellativity of the area addition. It follows from $\Aline cz\parallel \Aline bd$ and $B\parallel D$ that $[dbc]=[dbz]=[db'z]$ and hence $[db'c']=[dbc]=[db'z]$. By Proposition~\ref{p:same-area<=>subparallel}, the equality $[db'c']=[db'z]$ implies $\Aline z{c'}\parallel\Aline d{b'}$.

Finally, we obtain $\Aline z{c'}\parallel\Aline d{b'}\parallel \Aline z{a'}$ and hence $z\in \Aline{a'}{c'}$. Then $\Aline ac\cap\Aline {a'}{c'}=\{z\}\subset D$, witnessing that the Playfair plane $X$ is para-Desraguesian. By Theorem~\ref{t:affine-Moufang<=>}, the Playfair plane $X$ is Moufang.
\end{proof}

Since the area-additive Playfair plane $X$ is Moufang and area-injective, the following Proposition~\ref{p:Moufang+areainj=>Pappian} implies that $X$ is Pappian.
\end{proof}

\begin{proposition}\label{p:Moufang+areainj=>Pappian} Every area-injective Moufang Playfair plane  is Pappian.
\end{proposition}

\begin{proof} Let $X$ be an area-injective Moufang Playfair plane. 
\smallskip

1. First, we prove that $X$ is Desarguesian. By Theorem~\ref{t:Desargues<=>3Desargues}, it suffices to check that each line triple $oxe$ in $X$ is Desarguesian. Given any line affinity $A$ with $Aoe=oe$, we should check that the point $y\defeq A(x)$ is equal to $x$. Write the line affinity $A$ as the composition $A=A_n\cdots A_1$ of line projections (which are line shears or line shifts). Since each line shift is a composition of two line shears, we can assume that all line projections $A_1,\dots,A_n$ are line shears. Since the Moufang affine plane $X$ is shear (according to Theorem~\ref{t:affine-Moufang<=>}), each line shear $A_i$ is a restriction of some hypershear $\bar A_i:X\to X$. Then the composition $\bar A\defeq\bar A_n\cdots \bar A_1$ is an equiaffine automorphism of $X$ such that $\bar Aoxe=oye$. Choose any point $a\in X\setminus\Aline oe$ and consider the point $b\defeq \bar A(b)$. By Theorem~\ref{t:shear=>SAut}, the equiaffine automorphism $\bar A$ of $X$ is area-preserving and hence $[oeb]=[Aoea]=[oea]$. By Proposition~\ref{p:same-area<=>subparallel}, $\Aline oea\subparallel oeb$. Since $X$ is shear, there exists a hypershear $S:X\to X$ such that $\Aline oe\subseteq\Fix(S)$ and $S(b)=a$. Then $E\defeq S\bar A$ is an equiaffine automorphism such that $Eoea=oea$. Taking into account that $y\in\Aline oe\subseteq\Fix(S)$, we conclude that $E(x)=S\bar A(x)=S(y)=y$ and hence $Eoxea=oyea$, $[oxa]=[oya]$, and $x=y$, by the area-injectivity of $X$. Therefore, the each line triple in $X$ is Desarguesian and so is the Playfair plane $X$.
\smallskip

2. Next, we prove that the Desarguesian Playfair plane $X$ is Pappian.
To show that $X$ is Pappian, it suffices to show that the scalar corps $\IR_X$ of $X$ is commutative. Take any scalar $s\in[\IR^*_X,\IR^*_X]$ and choose any triangle $oep$ in $X$. Find a unique point $x\in \Aline oe$ such that $s=\overvector{oxe}$. By Proposition~\ref{p:sign-stabilizer} and Theorem~\ref{t:scalar-by-area}, $[oep]=s\cdot[oep]=\overvector{oxe}\cdot[oep]=[oxp]$.  
Since the Playfair plane $X$ is area-injective, $[oep]=[oxp]$ implies $x=e$ and hence $s=\overvector{oxe}=\overvector{oee}=1$. Therefore, $[\IR^*_X,\IR^*_X]=\{1\}$, the scalar corps $\IR_X$ is a field, and the Desarguesian affine plane $X$ is Pappian, by Theorem~\ref{t:Papp<=>Des+RX}.
\end{proof}

\begin{problem}\label{p:area-inj=>Moufang?} Is every area-injective plane Moufang? 
\end{problem}

\begin{remark} At the moment we can only prove that area-injective Playfair planes are associative-puls, see Theorem~\ref{t:area-inj=>ass-puls}. By Theorem~\ref{t:affine-Moufang<=>}, a Playfair plane is Moufang if and only if it is inversive-dot.
\end{remark}

Finally, we extend Theorem~\ref{t:area-additive<=>Pappian} to Playfair liners of arbitrary rank.

\begin{theorem}\label{t:area-additive<=>P} A Playfair liner $X$ is area-additive if and only if it is Pappian.
\end{theorem}

\begin{proof} To prove the ``only if'' part, assume that a Playfair liner $X$ is area-additive. By Theorem~\ref{t:area-additive<=>area-addition}, $X$ admits an additive area addition $+:[X^3]\times[X^3]\to[X^3]$. For every plane $\Pi$ in $X$, the restriction of the operation $+$ to the areol $[\Pi^3]$ of the plane $\Pi$, witnesses that the  Playfair plane $\Pi$ is area-additive, according to Theorem~\ref{t:area-additive<=>area-addition}. By Theorem~\ref{t:area-additive<=>Pappian}, the area-additive Playfair plane $\Pi$ is Pappian. Therefore, every plane in $X$ is Pappian, which implies that the liner $X$ is Pappian, see Proposition~\ref{p:Papp<=>4-Papp}.
\smallskip

To prove the ``if'' part, assume that a Playfair liner $X$ is Pappian. Let $\mathcal P$ be the family of all planes in $X$.  By Theorem~\ref{t:Playfair<=>}, the Playfair liner $X$ is $3$-regular and by Proposition~\ref{p:k-regular<=>2ex}, the $3$-regular liner $X$ is $3$-ranked. The $3$-rankedness of $X$ ensures that $\|P\cap\Pi\|\le 2$ for any distinct planes $P,\Pi\in\mathcal P$, which implies $P^\vartriangle\cap\Pi^\vartriangle=\varnothing$ and hence $[P^\vartriangle]\cap[\Pi^\vartriangle]=\varnothing$. Therefore, $\{[P^\vartriangle]:P\in\mathcal P\}$ is a disjoint cover of the nonzero areol $[X^\vartriangle]\defeq[X^3]\setminus\{[0]\}$ of $X$. 

By Theorem~\ref{t:area-additive<=>Pappian}, for every plane $P\in\mathcal P$, there exists an area addition $\oplus_P:[P^3]\times[P^3]\to[P^3]$. Define the binary operation $+:[X^3]\times[X^3]\to[X^3]$ by the formula
$$\alpha+\beta=\begin{cases} \alpha\oplus_P\beta&\mbox{if $\{\alpha,\beta\}\subseteq [P^\vartriangle]$ for some $P\in\mathcal P$};\\
\alpha&\mbox{if $\beta=[0]$};\\
\beta&\mbox{if $\alpha=[0]$};\\
[0]&\mbox{otherwise}.
\end{cases}
$$ Since $\{[P^\vartriangle]:P\in\mathcal P\}$ is a disjoint cover of $[X^\vartriangle]$, the operation $+$ is well-defined.

We claim that this operation is additive. Given any line $L\subseteq X$ and  points $a,b,c\in L$ and $o\in X$, we shoud prove the equality $[oab]+[obc]=[oac]$. If $o\in L$, then $[oab]+[obc]=[0]+[0]=[0]=[oac]$. So, assume that $o\notin L$ and consider the plane $P\defeq\overline{L\cup\{o\}}$. If $a=b$, then $[oab]+[obc]=[0]+[obc]=[obc]=[oac]$. If $b=c$, then $[oab]+[obc]=[oab]+[0]=[oac]$. So, we assume that $a\ne b\ne c$. In this case, the areas $[oab]$ and $[obc]$ belong to the set $[P^\vartriangle]$ and
$[oab]+[oac]=[oab]\oplus_P[obc]=[oac]$, by the definition of the area addition $+$ and the additivity of the area addition $\oplus_P$. By Theroem~\ref{t:area-additive<=>area-addition}, the additive area addition $+:[X^3]\times [X^3]\to[X^3]$ witnesses that the Playfair Pappian liner $X$ is area-additive.
\end{proof}

\begin{corollary}\label{c:Pappian-plane=>area-bijective} Every Pappian Playfair plane is area-bijective.
\end{corollary}

\begin{proof} Let $X$ be a Pappian Playfair plane. By Proposition~\ref{p:hyper-Bolyai-interplay}, the Playfair liner $X$ is hyper-Bolyai. By Proposition~\ref{p:area-surjective-around}, the hyper-Bolyai plane $X$ is area-surjective. By Theorem~\ref{t:area-additive<=>Pappian}, the Pappian Playfair line $X$ is area-additive. By Proposition~\ref{p:area-additive=>area-injective}, the area-additive liner $X$ is area-injective. Therefore, $X$ is area-bijective, being area-injective and area-surjective.
\end{proof}
\newpage

\section{First order Area Addition Axiom}

\rightline{\em Things equal to the same thing are equal to one another.}

\rightline{Euclid, ``Elements'', Book I, Common Notions I}
\bigskip

Theorem~\ref{t:area-additive<=>P}, which shows that the Pappus Axiom and the Area Addition Axiom are equivalent, is quite appealing but slightly asymmetric. Unlike the Pappus Axiom, the Area Addition Axiom is not a first-order property, because it refers to the notion of area, which itself is not first-order. Recall that a property is called first-order if it can be expressed entirely in terms of the basic objects of the theory (such as points, lines, and their relations) by a logical formula that does not quantify over sets, functions, or other higher-level entities. Fortunately, Theorem~\ref{t:area-additive<=>P} can be reformulated in a first-order way by replacing the notion of area with the first-order relation $\boxtimes^7$. However, to make the relation $\boxtimes^7$ coincide with the relation $\boxtimes^{<\omega}$, expressing that two triangles have the same area, we need to add one more first-order axiom, which is an instance of Euclid’s Common Notion I, mentioned in the epigraph to this section.

\begin{proposition}\label{p:AAA<=>1} For a hyper-Bolyai liner $X$, the Area Addition Axiom is equivalent to the conjunction of the following two instances of Euclid's Common Notions I and II:
\begin{enumerate}
\item[(\textup{I})] $\forall abc,uvw,xyz\in X^3\;(abc\boxtimes^7 uvw\boxtimes xyz\;\Ra abc\boxtimes^7 xyz)$;
\item[(\textup{II})] $\forall L,L'\in\mathcal L_X\;\forall a,b,c\in L\;\forall a'b'c'\in L'\;\forall o,o'\in X$ $(oab\boxtimes^7 \!o'a'b'\,\wedge\, obc\boxtimes^7\! o'b'c')\Ra oac\boxtimes^7 \!o'a'c'$.
\end{enumerate}
Here $\mathcal L_X$ is the family of all lines in $X$.
\end{proposition}

\begin{proof} Assume that a hyper-Bolyai liner $X$ satisfies the Area Addition Axiom. Then $X$ is area-additive and area-injective, by Proposition~\ref{p:area-additive=>area-injective}. By Proposition~\ref{p:area-surjective-around}(4), the binary relations $\boxtimes^\w  $ and $\boxtimes^7$ on $X^3$ coincide.
 Then for any triples $abc,uvw,xyz\in X^3$ with $abc\boxtimes^7 uvw\boxtimes xyz$ we obtain $abc\boxtimes^\w  uvw\boxtimes^\w  $ and hence $abc\boxtimes^\w   xyz$ and finally $abc\boxtimes^7 xyz$. Therefore, the condition (I) is satisfied.
Since $\boxtimes^\w  =\boxtimes^7$, the Area Addition Axiom turns into the condition (II). This completes the proof of the ``only if'' part.
\smallskip

To prove the ``if'' part, assume that the hyper-Bolyai liner has the properties (I) and (II). The conditionn (I) implies that $\boxtimes^n=\boxtimes^7$ for all $n\ge 7$. Indeed, for $n=7$ this it obvious. Assume that for some $n\ge 7$ we have proved that $\boxtimes^n=\boxtimes^7$. Then $\boxtimes^{n+1}\defeq \boxtimes^n\circ \boxtimes =\boxtimes^7\circ\boxtimes=\boxtimes^7$, by condition (I). The Principle of Mathematical Induction ensures that $\boxtimes^n=\boxtimes$ for all $n\ge 7$ and hence $\boxtimes^\w  \defeq\bigcup_{n\in\w}\boxtimes^n=\bigcup_{n=7}^\infty\boxtimes^n=\boxtimes^7$. Since $\boxtimes^\w  =\boxtimes^7$, the condition (II) coincides with the Area Addition Axiom.
\end{proof}

Theorem~\ref{t:area-additive<=>P} and Propositions~\ref{t:area-additive<=>P}, \ref{p:hyper-Bolyai-interplay} imply the following first order characterization of Pappian Playfair liners.

\begin{corollary} A Playfair liner $X$ is Pappian if and only if it satisfies the  the following two instances of Euclid's Comon Notions I and II:
\begin{enumerate}  \setlength{\leftmargin}{25pt}
\item[(\textup{I})] $\forall abc,uvw,xyz\in X^3\;(abc\boxtimes^7 \!uvw\boxtimes xyz\;\Ra abc\boxtimes^7 \!xyz)$;
\item[(\textup{II})] $\forall L,L'\in\mathcal L_X\;\forall a,b,c\in L\;\forall a'b'c'\in L'\;\forall o,o'\in X$ $(oab\boxtimes^7 \!o'a'b'\,\wedge\, obc\boxtimes^7\! o'b'c')\Ra oac\boxtimes^7 \!o'a'c'$.
\end{enumerate}
\end{corollary}

\section{Orientability of liners}

In this section we define and study some notions related to the orientability of liners. The orientability of individual areas in liners was considered in Section~\ref{s:orientable-area}. 

Let us recall that for a cardinal number $\kappa$, an area $\alpha\in [X^3]$ in a liner $X$ is called \index{area!$\kappa$-oriented}\index{$\kappa$-oriented area}\defterm{$\kappa$-oriented} if its $S_3$-orbit $\alpha S_3\defeq\{\alpha\circ \sigma:\sigma\in S_3\}\subseteq [X^3]$ has cardinality $|\alpha S_3|=\kappa$. It is clear that every area $\alpha\in [X^3]$ is $\kappa$-oriented for some number $\kappa\in\{1,2,3,6\}$. The zero area in any liner is $1$-oriented.

\begin{definition} Let $K$ be a set of cardinals. A liner $X$ is called 
\index{liner!$K$-oriented}\index{$K$-oriented liner}\defterm{$K$-oriented} if $|\alpha S_3|\in K$ for all nonzero areas $\alpha\in [X^\vartriangle]$.
\end{definition}

In fact, every liner is $\{1,2,3,6\}$-oriented; $K$-oriented liners for $K$ equal to $\{1\},\{2\},\{1,2\}$ or $\{6\}$ have special names.

\begin{definition} A liner $X$ is called
\begin{itemize}
\item \index{unoriented liner}\index{liner!unoriented}\defterm{unoriented} if it is $\{1\}$-oriented;
\item \index{oriented liner}\index{liner!oriented}\defterm{oriented} if it is $\{2\}$-oriented;
\item \index{semioriented liner}\index{liner!semioriented}\defterm{semioriented} if it is $\{1,2\}$-oriented;
\item \index{hexaoriented liner}\index{liner!hexaoriented}\defterm{hexaoriented} if it is $\{6\}$-oriented.
\end{itemize}
\end{definition}

\begin{exercise} Show that a liner $X$ is 
\begin{itemize}
\item unoriented if and only if every area $\alpha\in [X^3]$ is unoriented if and only if $[abc]=[acb]=[bac]=[bca]=[cab]=[cba]$ for any triangle $abc\in X^\vartriangle$;
\item oriented if and only if every area $\alpha\in [X^3]$ is oriented if and only if $[abc]=[bca]=[cab]\ne [acb]=[bac]=[cba]$ for any triangle $abc\in  X^\vartriangle$;
\item semioriented if and only if every area $\alpha\in [X^3]$ is semioriented if and only if $[abc]=[bca]=[cab]$ for any triangle $abc\in X^\vartriangle$;
\item hexaoriented if and only if every area $\alpha\in [X^3]$ is hexaoriented if and only if for every triangle $abc\in X^\vartriangle$ the areas $[abc]$, $[acb]$, $[bac]$, $[bca]$, $[cab]$, $[cba]$ all are distinct.
\end{itemize}
\end{exercise}

We start with studying orientability properties of completely regular planes. Recall that a \index{plane}\defterm{plane} is any liner of rank $3$.

\begin{theorem}\label{t:comreg=>semihexa} Every completely regular plane $X$ is either semioriented or hexaoriented. Moreover,
\begin{enumerate}
\item[\textup{(1)}] $X$ is semioriented if and only if it is Playfair;
\item[\textup{(2)}] $X$ is hexaoriented if and only if it is a projective plane or a punctured projective plane.
\end{enumerate}
\end{theorem}

\begin{proof} Assume that $X$ is a completely regular plane. First we assume that $X$ is projective. In this case we shall prove that $X$ is hexaoriented.
Since $X$ is a proective plane, it contains no disjoint lines, which implies that every traingle $abc$ in $X$ has area $\{abc\}$. Then the family $$\big\{[abc], [acb], [bac], [bca], [cab], [cba]\big\}=
\big\{\{abc\},\{acb\}, \{bac\}, \{bca\}, \{cab\}, \{cba\}\big\}$$has cardinality $6$, witnessing that the area $[abc]$ is $6$-oriented and the liner $X$ is hexaoriented.
\smallskip

Next, assume that the completely regular liner $X$ is not projective. 
By Theorem~\ref{t:spread=projective2}, the completely regular liner $X$ has a projective completion $Y$ with flat horizon $H\defeq Y\setminus X$, which is not empty because $X$ is not projective. Since $\overline H\ne Y$, for every points $h\in H$ and $x\in Y\setminus\overline H\subseteq X$ we have $\Aline xh\cap \overline H=\{h\}$ and hence $\Aline xh\cap X$ is a line in $X$, which implies that $H$ belongs to the flat hull of $X$ in $Y$, and hence $\|Y\|\le\|X\|=3$ and $Y$ is a $3$-long projective plane. Then $0<\|H\|<\|Y\|=3$ and hence $\|H\|\in\{1,2\}$.
\smallskip

First assume that $\|H\|=1$ and hence $H=\{d\}$ for some point $d\in Y$.
Then $X=Y\setminus\{d\}$ is a punctured projective plane. 
The point $d$ can be identified with the direction $\boldsymbol\delta\defeq\{\Aline xd\setminus\{d\}:x\in X\}$ in $X$. Observe that two lines in $X$ are disjoint if and only if they are distinct and both belong to the direction $\boldsymbol \delta$. This implies that for any triangle $abc$ in $X$, its area can be found by the formula
$$[abc]=\begin{cases}
\{abc\}&\mbox{if $d\notin \Aline ab\cup\Aline bc\cup\Aline ca$};\\
\big\{xbc:x\in \Aline ad\setminus\{d\}\big\}&\mbox{if $d\in \Aline bc$};\\
\big\{ayc:y\in \Aline bd\setminus\{d\}\big\}&\mbox{if $d\in \Aline ac$};\\
\big\{abz:z\in \Aline cd\setminus\{d\}\big\}&\mbox{if $d\in \Aline ab$}.
\end{cases}
$$
This formula implies that all areas $[abc]$, $[acb]$, $[bac]$, $[bca]$, $[cab]$, $[cba]$ are distinct and hence the area $[abc]$ is $6$-oriented, witnessing that the liner $X$ is hexaoriented. 
\smallskip

Next, assume that $\|H\|=2$. In this case, $H$ is a line in the projective plane $Y$. By Theorem~\ref{t:affine<=>hyperplane}, the liner $X=Y\setminus H$ is affine and regular. Assuming that $|Y|_2=3$, we conclude that $|X|_2=2$ and $|X|=|Y|-|H|=(|Y|_2^2-|Y|_2+1)-|Y|_2=7-3=4$ (by Corollary~\ref{c:Steiner-projective<=>}) and hence $X$ has rank $\|X\|=|X|=4$ and is not a plane. This contradiction shows that $|Y|_2\ge 4$ and hence $|X|_2\ge 3$. By Theorem~\ref{t:Playfair<=>}, the $3$-long affine regular liner $X$ is Playfair. 
By Corollary~\ref{c:area-semioriented-in-Playfair}, every area in the Playfair liner $X$ is semioriented, which implies that $X$ is semioriented. 
\end{proof}

Now we consider some sufficient conditions of semiorientability and unorientability.

\begin{proposition}\label{p:liner-semioriented<=} A liner $X$ is semioriented (and unoriented) if for every triangle $xyz$, there exists a (Boolean) parallelogram $abcd$ in $X$ such that $\{x,y,z\}\subseteq \{a,b,c,d\}$.
\end{proposition}

\begin{proof} Assume that  for every triangle $xyz\in X^\vartriangle$, there exists a (Boolean) parallelogram $abcd$ in $X$ such that $\{x,y,z\}\subseteq \{a,b,c,d\}$. Given any nonzero area $\alpha\in [X^\triangle]$, we have to show that $|\alpha S_3|\le 2$. Choose any triangle $xyz\in\alpha$. By the assumption, there exists a (Boolean) parallelogram $abcd$ in $X$ such that $\{x,y,z\}\subseteq a,b,c,d\}$. Proposition~\ref{p:12-permutations} implies that $[xyz]=[yzx]=[zxy]$. Then the stabilizer $S_3^{=\alpha}\defeq\{\sigma\in S_3:\alpha\circ\sigma=\alpha\}$ of the area $\alpha=[abc]$ has cardinality at least $3$. Taking into account that $S_3^{=\alpha}$ is a subgroup of the permutation group $S_3$, we conclude that $S_3^{=\alpha}\in\{3,6\}$ and hence $|\alpha S_3|=|S_3|/|S_3^{=\alpha}|\in\{1,2\}$, witnessing that the area $\alpha$ is semioriented and so is the liner $X$. 

If the parallelogram $abcd$ is Boolean, then Proposition~\ref{p:12-permutations} ensures that the areas of all triangles with vertices in the set $\{a,b,c,d\}$ are equal, which implies $|S_3^{=\alpha}|=6$ and $|\alpha S_3|=|S_3|/S_3^{=\alpha}|=1$, witnessing that the area $\alpha$ is unoriented and so is the liner $X$.
\end{proof}

Proposition~\ref{p:hyper-Bolyai=>parallelogram} and \ref{p:liner-semioriented<=} imply the following corollary.

\begin{corollary}\label{c:Playfair=>semioriented} Every hyper-Bolyai (and Boolean) plane $X$ is semioriented (and unoriented). 
\end{corollary}

Let us recall that a liner $X$ is called \index{$3$-homogeneous liner}\index{liner!$3$-homogeneous}\defterm{$3$-homogeneous} if for any triangles $abc$ and $xyz$ in $X$ there exists an automorphism $\Phi$ of $X$ such that $abc=\Phi xyz$.

\begin{proposition}\label{p:3-homogen=>area=1} Any $3$-homogeneous liner $X$ is unoriented, semioriented, or hexaoriented.
\end{proposition}

\begin{proof} First we show that for any two nonzero areas $\alpha,\beta\in [X^3]$, their $S_3$-orbits $\alpha S_3$ and $\beta S_3$ have the same cardinality. Choose any triangles $xyz\in\alpha$ and  $abc\in\beta$. Since the liner $X$ is $3$-homogeneous, there exists an automorphism $\Phi:X\to X$ such that $\Phi xyz=abc$.   By Proposition~\ref{p:permutation-area}, the automorphism $\Phi$ induces a well-defined bijection $$\Phi\circ:\alpha S_3\to \beta S_3,\quad \Phi\circ:[xyz\circ\sigma]\mapsto [\Phi \circ xyz\circ\sigma],$$ witnessing that $|\beta S_3|=|\alpha S_3|$. Therefeore, $\{|\alpha S_3|:\alpha\in [X^\vartriangle]\}=\{\kappa\}$ for some number $\kappa\in\{1,2,3,6\}$. It remains to show that $\kappa\ne 3$. 

Assuming that $\kappa=3$, take any triangle $abc$ and observe that its area $\alpha\defeq[abc]$ has stabilizer $S_3^{=\alpha}$ of cardinality $2$. Then exactly one of the areas $[acb]$, $[cba]$ or $[bac]$ is equal to $[abc]$. We lose no generality assuming that $[abc]=[cba]$. Since $X$ is $3$-homogeneous, there exists an automorphism $\Phi$ of $X$ such that $ \Phi abc=bca$. By Proposition~\ref{p:permutation-area}, the equality $[abc]=[cba]$ implies $[bca]=[\Phi abc]=[\Phi cba]=[acb]$ and $[abc]=[bac]$. Then the stabilizer $S_3^{\alpha}$ contains two distinct permutations of order 2 and hence $|S_3^{=\alpha}|>2$. This contradiction shows that $\kappa\ne 3$ and hence $X$ is $\{\kappa\}$-oriented for some $\kappa\in\{1,2,6\}$. Then $X$ is unoriented (if $\kappa=1$), oriented (if $\kappa=2$), or hexaoriented (if $\kappa=6$). 
\end{proof}

The most difficult technical result of this section if the following characterization of oriented Moufang Playfair planes.

\begin{theorem}\label{t:Moufang=1or2} Every Moufang Playfair plane is either oriented or unoriented. Moreover, $X$ is oriented if and only if $X$ is Desarguesian and $-1\notin [\IR^*_X,\IR^*_X]$.
\end{theorem}

\begin{proof}  Let $X$ be a Moufang Playfair plane. By Theorem~\ref{t:Playfair<=>}, the Playfair plane $X$ is $3$-long, affine and regular. Then $X$ is an affine space and hence $X$ has a well-defined scalar corps $\IR_X$. By Theorem~\ref{c:Moufang=>3hom}, the Moufang Playfair plane $X$ is $3$-homogeneous, and by Proposition~\ref{p:3-homogen=>area=1}, $X$ is $\{k\}$-oriented for some $\kappa\in\{1,2,3\}$. On the other hand, Proposition~\ref{c:Playfair=>semioriented} ensures that $X$ is semioriented and hence $\kappa\in\{1,2\}$. Therefore, $X$ is either oriented (if $\kappa=2$) or unoriented (if $\kappa=1$).

\begin{lemma}\label{l:Des-1<=>oriented} If $X$ is a Desarguesian, then $X$ is oriented if and only if $-1\notin[\IR^*_X,\IR^*_X]$.
\end{lemma}

\begin{proof} Fix any nonzero area $\alpha\in [X^\vartriangle]$. Since $\alpha$ is semioriented, it has a well-defined opposite area $-\alpha$, which is equal to $\alpha$ if and only if $\alpha$ is unoriented if and only if the Moufang plane $X$ is unoriented. Let $\cdot:\IR_X\times[X^3]\to[X^3]$ be the operation of multiplication of areas by scalars, defined in Theorem~\ref{t:scalar-by-area}. By Theorem~\ref{t:scalar-by-area}(1), $(-1)\cdot\alpha=-\alpha$. By Proposition~\ref{p:sign-stabilizer}, the stabilizer $\Stab(\alpha)\defeq \{t\in \IR_X:t{\cdot}\alpha=\alpha\}$ of the area $\alpha$ equals the commutator $[\IR_X^*,\IR_X^*]$ of the multiplicative group $\IR^*_X$ of the scalar corps $\IR_X$ of $X$.
 Now we see that $-1\in [\IR^*_X,\IR^*_X]=\Stab(\alpha)$ if and only if $\alpha=-\alpha$ if and only if the liner $X$ is unoriented.
\end{proof}

\begin{lemma}\label{l:Moufang-not-Desarg=>unoriented} If the Moufang plane $X$ is not Desarguesian, then $X$ is unoriented.
\end{lemma}

\begin{proof} Assume that the Moufang Playfair plane $X$ is not Desarguesian. Fix any affine base $uow$ in $X$ and consider its diunit $e$ and its ternar $\Delta=\Aline oe$. By Theorem~\ref{t:affine-Moufang<=>}, the ternar $\Delta$ of the Moufang based affine plane $(X,uow)$ is linear, distributive, associative-plus, and alternative-dot. In terminology of Chapters~\ref{ch:alt-rings}, \ref{ch:Cayley-Dickson}, the set $\Delta$ endowed with the plus and dot operations is an alternative division ring. Since $X$ is not Desarguesian, the alternative division ring $\Delta$ is not associative, by Theorem~\ref{t:ass-dot<=>}. By Corollary~\ref{c:div+com+alt=>ass}, the non-associative alternative division ring $\Delta$ is not commutative.  By Theorem~\ref{t:SBK} , the non-associative alternative division ring $\Delta$ admits a central involution turning $\Delta$ into an involutive octonion ring with the real axis $\Re(\Delta)=\mathcal Z(\Delta)=\mathcal N(\Delta)$. By Theorem~\ref{t:ab=-ba}, the non-commutative involutive alternative division ring $\Delta$ contains two elements $a,c$ such that $a{\cdot}c=-c{\cdot}a\ne 0$. By Artin's Theorem~\ref{t:Artin}, the subring $R$ of $\Delta$, generated by the set $\mathcal N(\Delta)\cup\{a,c\}$ is associative. Proposition~\ref{p:CD-1gen} implies that $R$ is an associative divisible ring. The associativity of $R$ ensures that the anticommuting elements $a,c$ have commutator $a{\cdot}c{\cdot}a^{-1}{\cdot}c^{-1}=-c{\cdot}a{\cdot}a^{-1}{\cdot}c^{-1}=-1$. 

Let $w'\in X$ be a unique point such that $\overvector{ow'}=-\overvector{ow}$. The point $w'$ has coordinates $(0,-1)=(o,-e)$ in the base $uow$. Consider the function $A:X\to X$ assigning to each point with coordinates $(x,y)$ the point with coordinates $(x,y{\cdot} aca^{-1}c^{-1})=(x,-y)$. Applying Propositions~\ref{p:(x,y)->(x,xc+y)}, \ref{p:horizontal-shear}, and Theorem~\ref{t:shear<=>alternative} on the structure of hypershears in Moufang Playfair planes, and repeating the proof of Lemma~\ref{l:commutator=equiaffine}, we can find an equiaffine automorphism $E$ of the Moufang plane $X$ such that $A(p)=E(p)$ for any point $p\in X$ with coordinates in the associative subring $R$ of the ring $\Delta$. Then $Euow=Auow=uow'$ and hence $[uow]=[uow']$ and $[wou]=[w'ou]=(-1)\cdot [wou]=-[wou]$, witnessing that the nonzero area $[wou]$ is unoriented. Then the $3$-homogeneous liner $X$ also is unoriented.
\end{proof}

Now we can complete the proof of the second part of Theorem~\ref{t:Moufang=1or2}. If $X$ is Desarguesian and $-1\notin[\IR^*_X,\IR^*_X]$, then $X$ is oriented, by Lemma~\ref{l:Des-1<=>oriented}. Now assume that the Moufang liner $X$ is oriented. By Lemma~\ref{l:Moufang-not-Desarg=>unoriented}, the liner $X$ is Desarguesian. By Lemma~\ref{l:Des-1<=>oriented}, $-1\in[\IR_X^*,\IR^*_X]$. 
\end{proof}

\begin{corollary}\label{c:Pappian=>1or2}  Every Pappian Playfair plane $X$ is either oriented or unoriented. Moreover, $X$ unoriented if and only if $X$ is Boolean if and only if $X$ contains a Boolean parallelogram.
\end{corollary}

\begin{proof} Let $X$ be a Pappian Playfair plane. By Theorems~\ref{t:Papp<=>APA}, \ref{t:Hessenberg-affine}, \ref{t:ADA=>AMA},  and Proposition~\ref{p:Desarg=>Moufang}, \ref{p:uno-Thalesian=>bi-di-by-Boolean}, the Pappian Playfair plane $X$ is Desarguesian, Moufang, Thalesian, uno-Thalesian, and di-Boolean. By Theorem~\ref{t:Papp<=>Des+RX}, the scalar corps $\IR_X$ of the Pappian liner $X$ is commutative and hence $[\IR^*_X,\IR^*_X]=\{1\}$.
 By Theorem~\ref{t:Moufang=1or2}, the Moufang Playfair plane $X$ is either oriented or unoriented. Moreover, $X$ is unoriented if and only if $-1\in[\IR^*_X,\IR^*_X]=\{1\}$ if and only if $1+1=0$.  By Propositions~\ref{p:Boolean<=>1+1=0} and Theorem~\ref{t:Boolean<=>di-Boolean+}, the condition $1+1=0$ holds if and only if the Thalesian liner $X$ is Boolean if and only if the di-Boolean liner $X$ contains a Boolean parallelogram. 
\end{proof}

Theorems~\ref{t:Moufang=1or2} and \ref{t:Desargues<=>planeD} imply the following corollary. 

\begin{corollary} Every oriented Moufang Playfair liner is Desarguesian.
\end{corollary}

\begin{problem}\label{prob:oriented=>Desarg} Is every oriented finite Playfair plane Desarguesian? 
\end{problem}

\begin{problem}\label{prob:oriented=>Thales} Is every oriented finite Playfair plane Thalesian?
\end{problem}

\begin{remark} By computer calculations, Ivan Hetman verified that all known  orientable Playfair planes of order $<25$ are Desarguesian and hence Thalesian. On the other hand, infinite oriented Playfair planes need not be Thalesian. For example, the Moulton plane is oriented and non-Thalesian.
\end{remark}

\begin{exercise} Prove that the Moulton plane is oriented.
\end{exercise}
 
Problem~\ref{prob:oriented=>Desarg} is related to the presence of Boolean parallelograms in liners.
 
\begin{proposition}\label{p:oriented=>noBoole} An oriented liner cannot contain Boolean parallelograms.
\end{proposition}

\begin{proof} Assume that $abcd$ is a Boolean parallelogram in a liner $X$. By Proposition~\ref{p:12-permutations}, any two triangles with vertices in the set $\{a,b,c,d\}$ have the same area. In particular, the area $[abc]$ has $S_3$-orbit $[abc]S_3$ of cardinality $|[abc]S_3]|=1$, which contradicts the orientability of $X$.
\end{proof}

\begin{remark} Computer calculations (made by Ivan Hetman) witness that the unique non-Desarguesian Thalesian plane of order 9 does not contain Boolean parallelograms but is unoriented. So, Proposition~\ref{p:oriented=>noBoole} cannot be turned into a characterization. Nonetheless, we can ask the following problem.
\end{remark} 

\begin{problem}\label{prob:noBoole=>Thales} Is every finite Playfair plane without Boolean parallelograms Thalesian?
\end{problem} 

\begin{remark}\label{rem:Hanna-Neumann} By Proposition~\ref{p:oriented=>noBoole}, the affirmative answer to Problem~\ref{prob:noBoole=>Thales} implies the affirmative answer to Problem~\ref{prob:oriented=>Thales}. Problem~\ref{prob:noBoole=>Thales} can be considered as an affine counterpart of the following well-known open problem, attributed to \index[person]{Neumann}Hanna Neumann.\footnote{{\bf Hanna Neumann} (1914--1971) was a German-born mathematician and one of the leading algebraists of the 20th century. Educated in G\"ottingen and Berlin, she completed her doctorate at Oxford in 1939 under G. A. Reid, after emigrating to England with her husband, Bernhard Neumann. Her research focused on group theory, where she made lasting contributions to the theory of free products and varieties of groups. She is best known for the Hanna Neumann Conjecture on intersections of finitely generated subgroups of free groups—later proved in the 2010s—and for posing an open problem about Fano quadrangles in non-Desarguesian projective planes, which connected her algebraic interests with geometric structures. In 1963 she moved to Australia, becoming Professor and Head of the Department of Pure Mathematics at the Australian National University. Her mathematical insight and leadership had a lasting influence on the development of modern algebra and geometry.}
\end{remark}

\begin{problem}[Hanna Neumann]\label{prob:Hanna-Neumann} Is every finite projective plane without Fano quadrangles Desarguesian?
\end{problem}

\begin{remark} The affirmative answer to Problem~\ref{prob:noBoole=>Thales} implies the affirmative answer to Problem~\ref{prob:Hanna-Neumann}. Indeed, assume that every finite Playfair plane without Boolean parallelograms is Thalesian. Take any finite projective plane $X$ without Fano quadrangles. Then for every line $L\subseteq X$, the affine subliner $X\setminus L$ contains no Boolean parallelograms and hence is Thalesian. This means that the projective plane $X$ is everywhere Thalesian and hence Moufang, by Corollary~\ref{c:Moufang<=>everywhereThales}. By Theorem~\ref{t:Moufang-finite<=>}, the finite Moufang plane $X$ is Desarguesian.
\end{remark}

So far, we have encountered neither 
$3$-oriented planes nor $3$-oriented areas. This suggests the following problem.

\begin{problem} Can a plane be $3$-oriented? Can it contain a triangle with $3$-oriented area?
\end{problem}

Theorems \ref{t:comreg=>semihexa} and \ref{t:2balance=>semihexa} show that $3$-oriented planes occur neither among completely regular planes nor among finite balanced planes. The following proposition shows that 
$3$-oriented planes also do not occur among finite Proclus planes.

\begin{proposition}\label{p:Proclus=>semihexa} Every finite Proclus plane is $\{1,2,6\}$-oriented.
\end{proposition}

\begin{proof} Let $X$ be a finite Proclus plane. If $X$ is projective, then it contains no disjoint lines, which implies that every triangle $abc$ in $X$ has area $[abc]=\{abc\}$ and hence is $6$-oriented. In this case, the (projective) plane $X$ is hexaoriented. So, assume that $X$ is not projective. 

By Kuiper--Dembowski Theorem~\ref{t:Kuiper-Dembowski}, $X=P\setminus H$ for some projective plane $P$ and some proflat set $H\subseteq P$ of rank $0<\|H\|\le 2$. If $X$ is not $3$-long, then we can take $P$ to be Steiner, according to Theorem~\ref{t:Proclus-not-3long}. In this case $P$ is a projective competion of $X$. If $X$ is $3$-long, then fix any point $x\in X$ and choose any maximal $3$-long flat $M$ in $P$ that contains the point $x$. By Lemma~\ref{l:ox=2}, $\Aline xy=\{x,y\}$ for all points $y\in P\setminus X$. Since $X$ is $3$-long, this implies $X\subseteq M$ and hence $P=\overline X\subseteq \overline M=M\subseteq P$, witnessing that the projective plane $P=M$ is $3$-long and hence is a projective completion of $X$. If $H$ is flat, then the plane $X$ is regular and completely regular by Theorems~\ref{t:regular-horizon3} and \ref{t:spread=projective2}. By Theorem~\ref{t:comreg=>semihexa}, $X$ is hexaoriented or semoriented.

So, assume that $H$ is not flat, which implies that $H$ is a punctured line in $P$. Then there exists a unique point $p\in X$ such that $H\cup \{p\}=\overline H$ is a line in $P$. Observe that every line in $X$ that contains $p$ intersects all other lines in $X$, which implies that every triangle $abc\in X^\vartriangle$ with $p\in \{a,b,c\}$ has area $[abc]=\{abc\}$ and hence is $6$-oriented. Next choose any triangle $abc$ in $X$ with $p\notin\{a,b,c\}$. Then at least two sides of the triangle $abc$ do not contain the point $p$. We lose no generality assuming that $p\notin \Aline ab\cup\Aline bc$. Since $P$ is a projective plane, there exist points $\alpha\in\Aline ab\cap\overline{H}$ and $\gamma\in\Aline bc\cap\overline{H}$ which belong to the set $\overline H\setminus\{p\}=H$. Find a unique point $d\in\Aline \alpha c\cap\Aline\beta a\subseteq X$ and observe that $abcd$ is a parallelogram in $X$. By Proposition~\ref{p:semioriented<=}, the area $[abc]$ is semioriented. Therefore, every triangle is $X$ is either hexaoriented or semioriented, which implies that the liner $X$ is $\{1,2,6\}$-oriented.
\end{proof}

We do not know whether Proposition~\ref{p:Proclus=>semihexa}
generalizes to infinite Proclus planes.

\begin{problem} Is every Proclus plane $\{1,2,6\}$-oriented?
\end{problem}



\section{Monoareal liners}

Corollary~\ref{c:large=>unoriented} implies that balanced ranked liners whose cardinality is much larger than the length of their lines are unoriented. More precisely, we have the following corollary of Corollary~\ref{c:large=>unoriented}.
 
\begin{corollary}\label{c:large=>unoriented2} A $2$-balanced $3$-ranked liner $X$ is unoriented whenever every plane $P$ in $X$ has cardinality $|P|+9{\cdot}|X|_2>3{\cdot}|X|_2+9$.
\end{corollary}

In fact, Theorem \ref{t:large=>area-trivial} says much more than mere unorientability: if the cardinality of a balanced ranked plane is large enough, then all triangles have the same area. This suggests introducing the following notion, which will be studied in this section alongside unorientability.

\begin{definition} A liner $X$ is called \index{liner!monoareal}\index{monoareal liner}\defterm{monoareal} if all triangles in $X$ have the same area.
\end{definition}

It is clear that each monoareal liner is unoriented.

\begin{exercise} Find an example of a unoriented balanced finite plane, which is not monoareal.

{\em Hint:} Consider any Boolean Pappian Playfair plane.
\end{exercise}

\begin{examples} There exists a monoareal infinite Boolean Desarguesian plane.
\end{examples}

\begin{proof} By Exercise 3 in \cite{Cohn}, there exists an infinite corps $R$ of characteristic $2$ such that the commutator subgroup $[R^*,R^*]$ of its multiplicative group $R^*$ coincides with $R^*$. The abelianization $[R]=R/[R^*,R^*]$ of the multiplicative semigroup $(R,\cdot)$ of $R$ consists of just two elements $[0]=\{0\}$ and $[1]=[R^*,R^*]$. Consider the coordinate plane $\Pi=R\times R$ of the corps $R$. By Proposition~\ref{p:R=Delta}, the ternar $\Delta$ of the plane $\Pi$ is isomorphic to the corps $R$. By Theorems~\ref{t:ass-dot<=>} and \ref{t:ADA=>AMA}, the Playfair plane $\Pi$ is Desarguesian and Thalesian. Since the corps $R$ has characteristic $2$, the Thalesian plane $\Pi$ is Boolean, by Proposition~\ref{p:Boolean<=>1+1=0}. By  Theorem~\ref{t:Desarg=>tring=corps} and Proposition~\ref{p:R=Delta}, the scalar corps $\IR_\Pi$ of the plane $\Pi$ is isomorphic to the corps $R$ and hence $|[\IR_\Pi]|=|[R]|=2$. By Theorem~\ref{t:coscalar-by-area}, for any nonzero area $\alpha\in [\Pi^\triangle]$, the function $\cdot_\alpha:[R_X]\to[X^3]$, $\cdot_\alpha:t\mapsto t\cdot\alpha$, is bijective, which implies that $|[\Pi^3]|=2$ and hence the Desarguesian plane $\Pi$ is monoareal.
\end{proof} 

\begin{proposition} Every monoareal liner $X$ is ranked and has rank $\|X\|\le 3$.
\end{proposition}

\begin{proof} Let $X$ be a monoareal liner. To derive a contradiction, assume that $\|X\|>3$. Take any triangle $abc$ in $X$ and consider the plane $P\defeq\overline{\{a,b,c\}}$ in $X$. Since $3=\|P\|<\|X\|$, there exists a point $d\in X\setminus P$. Since $X$ is monoareal, $[abc]=[abd]$. By Corollary~\ref{c:same-area=>same-plane}, the equality $[abc]=[abd]$ implies $d\in \overline{\{a,b,d\}}=\overline{\{a,b,c\}}=P$, which contradicts the choice of $d$. This contradiction shows that $\|X\|\le 3$. Assuming that $X$ is not ranked, we can find two flats $A\subset B$ such that $\|A\|=\|B\|\le 3$ and $A\ne B$. This is possible only if $\|A\|=\|B\|=3$. So we can choose a triangle $abc$ in $A$ and then choose a point $d\in B\setminus A$. Repeating the above argument, we can see that $d=\overline{\{a,b,d\}}=\overline{\{a,b,c\}}\subseteq A$, which contradicts the choice of the point $d$.
\end{proof}

Theorem~\ref{t:large=>area-trivial} provides the following sufficient condition of monoareality.

\begin{corollary} Every balanced ranked plane $X$ of cardinality $|X|+12{\cdot}|X|_2>4{\cdot}|X|_2^2+11$ is monoareal.
\end{corollary}

\begin{remark} A computer search (performed by Ivan Hetman) did not find any balanced finite planes 
$X$ of cardinality $|X|>|X|_2^2$ that are not monoareal. This motivated the following open problem.
\end{remark}

\begin{problem}\label{prob:large=>unoriented} Is every finite balanced plane $X$ of cardinality $|X|>|X|_2^2$  unoriented (and monoareal)?
\end{problem}

In Propositions~\ref{p:orientability3} -- \ref{p:orientability5} we shall prove that a finite balanced plane $X$ of cardnality $|X|>|X|_2^2$ is (totally) unoriented whenever $|X|_2\le 5$ (and $|X|_2\le 4$). First, we establish some orientability properties of finite balanced planes.

\begin{theorem}\label{t:2balance=>semihexa} Every finite balanced plane $X$ is semioriented or hexaoriented. Moreover,
\begin{enumerate}
\item[\textup{(1)}] $X$ is hexaoriented if and only if $|X|<|X|_2^2$;
\item[\textup{(2)}] $X$ is semioriented if and only if $|X|\ge|X|_2^2$;
\item[\textup{(3)}] $X$ is oriented if $|X|=|X|_2^2$ and $X$ is Pappian and not Boolean;
\item[\textup{(4)}] $X$ is unoriented if $|X|\ge |X|_2^2$ and $X$ is Boolean;
\item[\textup{(5)}] $X$ is unoriented if $|X|>3{\cdot}|X|^2_2-9{\cdot}|X|_2+9$;
\item[\textup{(6)}] $X$ is monoareal if $|X|>4{\cdot}|X|^2_2-12{\cdot}|X|_2+11$.
\end{enumerate}
\end{theorem}

\begin{proof} Let $X$ be a finite balanced plane. Then the finite cardinals $|X|_2$ and $|X|_3=|X|$ are well-defined, which forces the plane $X$ to be ranked. By Theorem~\ref{t:2-balanced-dependence}, $|X|=1+(|X|_2-1){\cdot}\lambda$ for some cardinal $\lambda\ge|X|_2$ such that $|X|_2$ divides $(\lambda-1){\cdot}\lambda$.
If $\lambda=|X|_2$, then $|X|=|X|_2^2-|X|_2+1<|X|_2^2$ and the plane $X$ is projective, by Corollary~\ref{c:Steiner-projective<=>}. By Theorem~\ref{t:comreg=>semihexa}(2), the projective plane $X$ is hexaoriented.

Next, assume that $\lambda= |X|_2+1$, which is equivalent to $|X|=|X|_2^2$. By Corollary~\ref{c:Playfair<=>23balanced}, the balanced plane $X$ (of cardinality $|X|=|X|_2^2$) is Playfair. By Corollary~\ref{c:Playfair=>semioriented}, the Playfair plane $X$ is semioriented.  If $X$ is Boolean, then $X$ is unoriented, by Corollary~\ref{c:Playfair=>semioriented}. If $X$ is Pappian and not Boolean, then $X$ is oriented, by Corollary~\ref{c:Pappian=>1or2}.  

Finally, assume that $\lambda\ge |X|_2+2$, which is equivalent to $|X|>|X|_2^2$. By Proposition~\ref{p:23-balance=>k-parallel}, the balanced plane $X$ is $\kappa$-parallel for $\kappa=\lambda-|X|_2\ge 2$. By Proposition~\ref{p:k-parallel=>parallelogram}, every triangle $abc$ in $X$ can be completed to a parallelogram $abcd$. By Proposition~\ref{p:liner-semioriented<=}, the liner $X$ is semioriented.  If, in addition, $X$ is Boolean, then the paralleogram $abcd$ is Boolean, and by Proposition~\ref{p:liner-semioriented<=}, the plane $X$ is unoriented.
\smallskip

If $|X|>3{\cdot}|X|_2^2-9{\cdot}|X|_2+9$, then $X$ is unoriented by Corollary~\ref{c:large=>unoriented2}.
\smallskip

If $|X|>4{\cdot}|X|_2^2-12{\cdot}|X|_2+11$, then $X$ is monoareal by Theorem~\ref{t:large=>area-trivial}.
\end{proof}

We now prove that any finite balanced plane of cardinality $|X| > |X|_2^2$ is unoriented for $|X|_2 \le 5$, and monoareal for $|X|_2 \le 4$, thereby giving a partial answer to Problem~\ref{prob:large=>unoriented}.
 
\begin{proposition}\label{p:orientability3} Every finite balanced plane $X$ with $|X|_2=3$ is hexaoriented, oriented, or monoareal. More precisely, $X$ is
\begin{enumerate}
\item[\textup{(1)}]  hexaoriented if and only if $|X|<9$;
\item[\textup{(2)}]  oriented if and only if $|X|=9$;
\item[\textup{(3)}]  monoareal if and only if $|X|>9$.
\end{enumerate}
\end{proposition}

\begin{proof} If $|X|<9$, then $X$ is hexaoriented, by 
Theorem~\ref{t:2balance=>semihexa}(1). If $|X|=9$, then the Steiner plane $X$ is Hall and Playfair, by Theorem~\ref{t:Hall<=>Playfair+Steiner}. By Corollary~\ref{c:4-Pappian}, the Steiner Playfair plane $X$ is Pappian. Since $|X|_2=3$ is odd, $X$ is not Boolean, see Corollary~\ref{c:Boolean-order}. By Corollary~\ref{c:Pappian=>1or2}, the Pappian Playfair non-Boolean plane $X$ is oriented.

Finally, assume that $|X|>9$. By Corollary~\ref{c:Kirkman}, $|X|\in (1+6\IN)\cup(6\IN-3)$, which implies $|X|\ge 13>11=4{\cdot}3^2-12{\cdot}3+11$.
By Theorem~\ref{t:2balance=>semihexa}(1), $X$ is monoareal.
\end{proof}

\begin{proposition}\label{p:orientability4} Every finite balanced plane $X$ with $|X|_2=4$ is hexaoriented or unoriented. More precisely, $X$ is
\begin{enumerate}
\item[\textup{(1)}]  hexaoriented if and only if $|X|<16$;
\item[\textup{(2)}]  unoriented if and only if $|X|\ge 16$;
\item[\textup{(3)}]  monoareal if and only if $|X|>16$.
\end{enumerate}
\end{proposition}

\begin{proof} If $|X|<16$, then $X$ is hexaoriented, by Theorem~\ref{t:2balance=>semihexa}(1). If $|X|=16$, then the balanced plane $X$ is Playfair, by Corollary~\ref{c:Playfair<=>23balanced}. By Corollary~\ref{c:4-Pappian}, the Playfair plane $X$ of order $4$ is Pappian and hence Thalesian. Since $|X|_2=4$, the Thalesian Playfair plane $X$ is Boolean, by Proposition~\ref{p:Boolean<=>1+1=0}. By Corollary~\ref{c:Playfair=>semioriented}, the Boolean Playfair plane $X$ is unoriented. By Theorem~\ref{t:coscalar-by-area}, $|[X^3]|=|\IR_X|=|X|_2=4>1$ and hence $X$ is not monoareal.

Finally, assume that $|X|>16$. By Theorem~\ref{t:2-balanced-dependence}, $|X|=1+(|X|_2-1){\cdot}\lambda=1+3{\cdot}\lambda$ for some cardinal $\lambda\ge|X|_2=4$ such that $\lambda(\lambda-1)\in 4\IZ$. The inequalies $1+3{\cdot}\lambda=|X|>16$ implies $\lambda>5$. Taking into account that $\lambda(\lambda-1)$ should be divisible by $4$, we conclude that $\lambda\ge 8$. Then $|X|=1+3{\cdot}\lambda\ge 25$. Since $|X|\ge 25>21=3{\cdot}4^2-9{\cdot}4+9$, the liner $X$ is unoriented, Theorem~\ref{t:2balance=>semihexa}(5). 

 If $|X|>4{\cdot}4^2-12{\cdot}4+11=27$, then $X$ is monoareal, by Theorem~\ref{t:2balance=>semihexa}(6). So, assume that $|X|=1+3{\cdot}\lambda\le 27$, which implies $\lambda=8$ and $|X|=25$.  
By \cite[II.1.34]{HCD}, there exist exactly 18 non-isomorphic liners of cardinality $25$ and line length $4$. Using a computer program, Ivan Hetman verified that all these $18$ liners are monoareal.
\end{proof}

\begin{proposition}\label{p:orientability5}  Every plane $X$ with $|X|_2=5$
  is hexaoriented, oriented, or unoriented. More precisely, $X$ is
\begin{enumerate}
\item[\textup{(1)}]  hexaoriented if and only if $|X|<25$;
\item[\textup{(2)}]  oriented if and only if $|X|=25$;
\item[\textup{(3)}]  unoriented if and only if $|X|>25$;
\item[\textup{(4)}]  monoareal if $|X|>45$.
\end{enumerate}
\end{proposition} 

\begin{proof}  If $|X|<25$, then $X$ is hexaoriented, by Theorem~\ref{t:2balance=>semihexa}(1).  If $|X|=25$, then the balanced plane $X$ is Playfair, by Corollary~\ref{c:Playfair<=>23balanced}. By Corollary~\ref{c:4-Pappian}, the Playfair plane $X$ of order $5$ is Pappian and hence Thalesian.
 Since $|X|_2=5$ is odd, $X$ is not Boolean, see Corollary~\ref{c:Boolean-order}. By Corollary~\ref{c:Pappian=>1or2}, the Pappian Playfair non-Boolean plane $X$ is oriented.
 
 Finally, assume that $|X|>25$. By Theorem~\ref{t:2-balanced-dependence}, $|X|=1+(|X|_2-1){\cdot}\lambda=1+4{\cdot}\lambda$ for some cardinal $\lambda\ge|X|_2=5$ such that $\lambda(\lambda-1)\in 5\IZ$. The inequality  $1+4{\cdot}\lambda=|X|>25$ implies $\lambda>6$. Taking into account that $\lambda(\lambda-1)$ is divisible by $5$, we conclude that $\lambda\ge 10$. Then $|X|=1+4{\cdot}\lambda\ge 41$. Since $|X|\ge 41>39=3{\cdot}5^2-9{\cdot}5+9$, the liner $X$ is unoriented, by Theorem~\ref{t:2balance=>semihexa}(5). 
 
If $|X|=1+4{\cdot}\lambda>45$, then $\lambda\ge 12$. Since $\lambda(\lambda-1)$ is divisible by $5$, $\lambda\ge 12$ implies $\lambda\ge 15$. Then $|X|=1+4{\cdot}\lambda\ge 61>51=4{\cdot}5^2-12{\cdot}5+11$ and $X$ is monoareal, by Theorem~\ref{t:2balance=>semihexa}(6).
\end{proof}

Let $X$ be a balanced plane with $|X|=5$ and $25<|X|\le 45$. Theorem~\ref{t:2-balanced-dependence} implies that $|X|\in\{41,45\}$.

\begin{problem} Are all (known) balanced planes $X$ with $|X|_2=5$ and $|X|\in\{41,45\}$ monoareal?
\end{problem}

\section{Biareal and bioriented liners}

Observe that a liner $X$ is monoareal if and only if its nonzero areol $[X^\vartriangle]$ is a singleton.

\begin{definition} A liner $X$ is called
\begin{itemize}
\item \index{liner!biareal}\index{biareal liner}\defterm{biareal} if $|[X^\vartriangle]|=2$;
\item \index{liner!bioriented}\index{bioriented liner}\defterm{bioriented} if $X$ is biareal and oriented.
\end{itemize}
\end{definition}

\begin{proposition} For a Playfair liner $X$, the following conditions are equivalent:
\begin{enumerate}
\item $X$ is Pappian and biareal;
\item $X$ is Pappian and bioriented;
\item $X$ is a Hall plane.
\end{enumerate}
\end{proposition}

\begin{proof} $(1)\Ra(2)$ Assume that the Playfair liner $X$ is Pappian and biareal. Then $|[X^\triangle]|=2$ and $|[X^3]|=3$. By Corollary~\ref{c:Pappian-plane=>area-bijective}, the Pappian Playfair liner $X$ is area-bijective, which implies $|X|_2=|[X^3]|=3$. By Corollary~\ref{c:Boolean-order}, the Pappian liner $X$ is not Boolean, and by Corollary~\ref{c:Pappian=>1or2}, the non-Boolean Pappian liner $X$ is oriented and hence bioriented (being biareal).
\smallskip

$(2)\Ra(3)$ Assume that a Playfair liner $X$ is Pappian and bioriented. We have the prove that $X$ is a Hall plane. Assuming that $\|X\|\le 2$, we conclude that $X^\vartriangle=\varnothing$ and hence $|[X^\vartriangle]|=0\ne 2$. This contradiction shows that $\|X\|\ge 3$. Given any nonzero area $\alpha\in[X^\vartriangle]$, consider the function $\cdot_\alpha:\IR_X\to [X^3]$, $\cdot_\alpha:t\to t\cdot\alpha$, which is bijective, according to Theorem~\ref{t:coscalar-by-area}. The orientability of $X$ implies that $-\alpha\ne \alpha$. It follows from $|[X^\vartriangle]|=2$ and $-\alpha\ne \alpha$ that $[X^3]=\{0,\alpha,-\alpha\}$. Then  $|X|_2=|\IR_X|=|[X^3]|=3$. Since every nonzero area in $X$ is equal to $\alpha$ or $-\alpha$, we can apply Corollary~\ref{c:same-area=>same-plane} and conclude that $\|X\|=3$.  Therefore, $X$ is a Steiner Playfair plane. By Theorem~\ref{t:Hall<=>Playfair+Steiner}, $X$ is a Hall plane.
\smallskip

$(3)\Ra(1)$ Assume that $X$ is a Hall plane. By  Theorem~\ref{t:Hall<=>Playfair+Steiner}, $X$ is a Steiner Playfair plane. By Corollary~\ref{c:4-Pappian}, $X$ is Pappian and Desarguesian. By Theorem~\ref{t:Desargues<=>3Desargues}, every line triple in $X$ is Desarguesian, which implies $|\IR_X|=|X|_2=3$. Since $|X|=9=|X|_2^2$, we can apply Proposition~\ref{p:orientability3}(2), and conclude that the liner $X$ is oriented. Then $\alpha\ne-\alpha$ for any nonzero area $\alpha\in[X^\vartriangle]$. By Theorem~\ref{t:coscalar-by-area}, the function  $\cdot_\alpha:\IR_X\to [X^3]$, $\cdot_\alpha:t\to t\cdot\alpha$, is bijective and hence $|[X^3]|=|\IR_X|=|X|_2=3$, witnessing that the oriented liner $X$ is bioriented.
\end{proof}

\begin{problem} Is each biareal (Playfair) liner $X$ bioriented?
\end{problem}

\begin{problem} Is every finite bioriented Playfair liner a Hall plane?
\end{problem}

\begin{exercise} Show that every biareal liner $X$ has rank $\|X\|=3$ and contains a plane $P\subseteq X$ such that $\overline{\{x,y,z\}}\in\{P,X\}$ for every triangle $xyz$ in $X$.
\end{exercise}

It is interesting to note that infinite bioriented Playfair planes do exist.

\begin{exercise} Prove that the Moulton plane is bioriented.
\end{exercise}

We recall that the Beltrami-Klein hyperbolic plane is the subliner $\{(x,y)\in\IR\times\IR:x^2+y^2<1\}$ of the real plane $\IR\times\IR$.

\begin{exercise} Prove that the Beltrami--Klein hyperbolic plane is bioriented.
\end{exercise}

The biorientability of the Multon or Beltrami--Klein plane $X$ allow us to recover the betweenness relation on the lines of those planes. 

\begin{definition} Given three distinct collinear points $a,b,c$ in a bioriented plane $X$, we say that $b$ is \defterm{between} $a$ and $c$ if $[oab]=[obc]$ for all points $o\in X\setminus \Aline ac$. 
\end{definition}

\begin{Exercise} Explore the properties of the betweennes relation in biordered planes.
\end{Exercise}

Also the biorientability allows us to define sides of a plane with respect to a line.

\begin{definition} Let $L$ be a line in a bioriented plane $X$. For a point $x\in X\setminus L$, the set $L(x)\defeq\{y\in X:\forall a,b\in L\;[axb]=[ayb]\}$ is called the \defterm{$L$-halfplane} of $X$ containing the point $x$.
\end{definition}

\begin{problem} Let $L$ be a line in a bioriented plane. What is the cardinality of the family $\{L(x):x\in X\setminus L\}$ of all $L$-halfplanes?
Is it always $2$?
\end{problem}

{\parskip0pt

\chapter{Selected Open Problems in Linear Geometry}\label{ch:open-prob}

\rightline{\em As long as a branch of science offers}

\rightline{\em an abundance of problems, so long is it alive;}

\rightline{\em a lack of problems foreshadows extinction or}

\rightline{\em the cessation of independent development.}


\smallskip

\rightline{\index[person]{Hilbert}David Hilbert\footnote{{\bf David Hilbert} (1862--1943) was a German mathematician widely regarded as one of the most influential figures in modern mathematics. He made foundational contributions across numerous fields, including invariant theory, algebraic number theory, functional analysis, and the formalization of geometry. Hilbert is also known for his work on the foundations of mathematics, particularly the formalization of axiomatic systems and his famous 23 problems, which guided much of 20th-century mathematical research. He held a professorship at the University of G\"ottingen, where he mentored many leading mathematicians. Hilbert’s work profoundly shaped mathematical logic, set theory, and the development of modern mathematical physics.}, 1900, Paris}
\vskip30pt

\section{Order of finite affine or projective planes}

In the following problems, a {\em plane} is any $3$-long $3$-regular liner of rank $3$. 
Let us recall that a number $n$ is called a {\em prime power} if $n=p^k$ for some prime number $p$ and some $k$. 

\begin{problem} Is the order of any finite affine plane a prime power?
\end{problem}

By Theorem~\ref{t:proj-finite<=>}, {\em every associative-add finite affine plane is Desarguesian}, which implies that {\em the order of any finite associative-add affine plane is a prime power}.

\begin{problem} Is the order of any finite inversive-add affine plane a prime power? Is every finite inversive-add projective plane Desargusian?
\end{problem}

\begin{problem} Is the order of any $4$-Pappian finite projective plane a prime power?
\end{problem}

By Theorems~\ref{t:prime-affine<=>+} and \ref{t:proj-finite<=>}, {\em an affine or projective plane of prime order is Pappian if it is inversive-plus or inversive-puls}, and by Theorem~\ref{t:invertible-add<=>4-Pappian}, {\em every invertible-add projective plane is $4$-Pappian}.

\begin{problem} Is every invertible-add affine (or projective) plane of prime order Pappian?
\end{problem}

\begin{problem} Is every projective plane of prime order $4$-Pappian?
\end{problem}

By the Bruck--Ryser Theorem~\ref{t:Bruck-Ryser}, {\em if a number $n\in (1+4\mathbb Z)\cup(2+4\IZ)$ is an order of a projective plane, then $n=a^2+b^2$ for  some integer numbers $a,b$}. In particular, {\em the order of a projective plane cannot be equal to $6, 14, 22$,...} By heavy computer calculations, Lam, Thiel, Swiercz \cite{LTS1989} proved that the order of a projective plane cannot be equal to $10\in 2+4\IZ$.

\begin{problem} Is there a projective plane of order $12$?
\end{problem}

\begin{problem} Is there a projective plane of order $45=5\cdot 3^2=6^2+3^2$?
\end{problem}

\begin{problem} Is there a projective plane of order $85=5\cdot 17=7^2+6^2$?
\end{problem}

\begin{problem} Is there a projective plane $X$ of composite order $p^2+p$ for some $p$? In particular, is there a projective plane of order $12=3^2+3$?
\end{problem}

\section{Boolean and Fano liners}\label{s:Boole-Fano}

By Gleason's Theorem~\ref{t:Gleason-Fano}, {\em every Fano finite projective liner is Pappian}. On the other hand, there exist finite Boolean Thalesian affine planes, which are not Fano.

\begin{problem} Is every (finite) Boolean Playfair plane Thalesian?
\end{problem}

\begin{problem} Is every Boolean affine space $\partial$-translation?
\end{problem}

\begin{problem} Is every finite di-Boolean Playfair plane Boolean or Thalesian?
\end{problem}

\begin{problem} Is every Fano projective plane Moufang?
\end{problem}

By Theorem~\ref{t:Fano<=>inversive}, {\em a finite projective plane $X$ is Fano if and only if $X$ has even order and $X$ is associative-plus or associative-puls.} 

\begin{problem}Is every $4$-Pappian projective plane of even order Fano?
\end{problem}

Using free projectivizations and Theorem~\ref{t:free-projectivization}, it is easy to construct an infinite projective plane that contains no Fano quadrangles.

\begin{problem}[Hanna Neumann]\label{prob:HN} Does every non-Desarguesian finite projective plane contain a Fano quadrangle?
\end{problem} 

\begin{problem} Find an inner characterization of subliners of Fano projective liners.
\end{problem}

By Theorem~\ref{t:Moufang=>pentagonal}, {\em every Moufang affine space is pentagonal.}

\begin{problem} Is every Fano affine space pentagonal?
\end{problem}

\begin{problem} Is every pentagramic finite projective plane Desarguesian?
\end{problem}

\section{Regularity of finite liners}

By Theorem~\ref{t:balanced<=>ranked}, every balanced finite affine or projective liner $X$ of rank $\|X\|\ge 5$ is regular. 

\begin{problem} Is every Steiner balanced finite liner of rank $4$ regular?
\end{problem}

\begin{problem} Is every finite balanced weakly regular liner of rank $\ge 4$ proaffine?
\end{problem}

By Theorem~\ref{t:completelyF<=>}, every $3$-long finite proaffine Fano liner is (completely) regular. On the other hand, by Theorem~\ref{t:regular-horizon3}, for a Steiner projective liner $X$ of rank $\|X\|=4$ and a line $L$ in $X$, the subliner $X\setminus L$ is an example of a finite non-regular proaffine Fano liner.

\begin{problem} Is every $\w$-long proaffine Fano liner regular?
\end{problem}

\begin{problem} Is every finite Pappian $3$-long Proclus liner regular?
\end{problem}

By Corollary~\ref{c:proj-embed}, {\em every liner is a subliner of some projective plane.}

\begin{problem} Is every finite liner a subliner of some finite projective plane?
\end{problem}


\begin{problem}[Ivan Hetman] Is every finite balanced liner of rank $4$ affine or projective?
\end{problem}

\begin{problem}[Ivan Hetman] Is every finite balanced liner of rank $4$ weakly modular?
\end{problem}

\section{Homogeneity of liners}

By Theorem~\ref{t:Moufang=>4-homogen}, {\em every $3$-long non-Fano Moufang projective plane is $4$-homogeneous}.

\begin{problem} Is every $3$-long Fano Moufang projective plane $4$-homogeneous?
\end{problem}

By Theorem~\ref{t:Ostrom-Wagner1959}, {\em every $2$-homogeneous finite projective plane is Desarguesian.}

\begin{problem} Is every $1$-homogeneous finite projective plane Desarguesian?
\end{problem}





By Jha--Johnson Theorem~\ref{t:Jha-Johnson1998}, {\em every $2$-homogeneous finite affine liner is Thalesian.}

\begin{problem} Is every $1$-homogeneous finite affine plane Thalesian?
\end{problem}

By Theorem~\ref{t:rigid-projplane}, {\em there exists a rigid projective plane}, and by Example~\ref{ex:rigid25}, {\em there exists a rigid finite Playfair plane}.

\begin{problem}[classical folklore] Is there a rigid finite projective plane?
\end{problem}

\begin{problem} Is there a finite projective plane whose every central automorphism is trivial?
\end{problem}

\begin{problem}[Ivan Hetman] Is there a rigid finite $2$-balanced hyperbolic plane?
\end{problem}

\begin{problem}[Ivan Hetman] Is every injective finite $2$-balanced plane a classical unital?
\end{problem} 

By Theorem~\ref{t:Krauz}, {\em  for any infinite group $G$, there exists a Playfair plane whose translation and dilation groups are equal and isomorphic to the group $G$.}

\begin{problem} Let $G$ be a finite group. Is there a (finite) Playfair plane $X$ whose translation group is isomorphic to the group $G$?
\end{problem}

\begin{problem} Is every $\partial$-translation finite Playfair plane translation?
\end{problem}

Let us recall that a liner $X$ is {\em hypersymmetric} if for every triangle $abc$, there exists a hyperfixed involutive automorphism $\Phi$ of $X$ such that $\Phi abc=aba$. By Theorem~\ref{t:hypersymmetric=>shear} and \ref{t:affine-Moufang<=>}, {\em every hypersymmetric affine liner is Moufang}.

\begin{problem} Is every hypersymmetric projective liner Moufang?
\end{problem}

\section{Characteristics of affine and projective planes}

We recall that 
\begin{itemize}
\item the \defterm{characteristic} of a ternar $R$ is the cardinality of its minimal subternar $\underline{R}$;
\item the \defterm{characteristic} of a triangle $uow$ in an affine liner is the characteristic of the ternar $\Delta_{uow}$ of $uow$;
\item the \defterm{ternar} $\Delta_{uow}$ of a triangle $uow$ in an affine liner $X$ is the ternar of the based affine plane $(\overline{ouw},uow)$ if this plane is Thalesian, and the field of characteristic $p=|X|_2\in\{2,3\}$ in the opposite case;
\item the \defterm{characteristic range} $\har[X^\vartriangle]$ of an affine liner $X$ is the set of characteristics of all triangles in $X$;
\item the \defterm{characteristic range} $\har[X^\diamondsuit]$ of a projective liner $X$ is the set of characteristics of all quadrangles in $X$;
\item the \defterm{characteristic} of a quadrangle $abcd$ in a projective liner $X$ is the order of the smallest projective subplane $\langle abcd\rangle$ of $X$ that contains the points $a,b,c,d$.
\end{itemize}

An affine or projective liner is called \defterm{unicharacteristic} if its characteristic range is a singleton. In this case, the unique element of this singleton is called the \defterm{characteristic} if a unicharacteristic liner. So, saying that an affine or projective liner has (prime) characteristic $p$ we understand that the singleton $\{p\}$ coincides with the characteristic range of the liner.
By Proposition~\ref{p:char-divides-order} and Corollary~\ref{c:unichar=>char-divide-order}, the characteristic of any unicharacteristic affine of projective liner divides the order of the liner.

 There are many intriguing long-standing open problems concerning the structure of unicharacteristic finite affine or projective planes. We formulate these problems as conjectures, since no counterexamples have yet been found among known finite affine of projective planes. The strongest among them is the following conjecture.

\begin{conjecture}\label{conj:unihar=>Thalesian} Every unicharacteristic Playfair liner is Thalesian.
\end{conjecture}

Since every Thalesian finite Playfair plane has prime characteristic, Conjecture~\ref{conj:unihar=>Thalesian} is true if and only if the following two conjectures are true.

\begin{conjecture} Every unicharacteristic affine liner has prime characteristic.
\end{conjecture}

\begin{conjecture}\label{c:primechar=>Thalesian} Every Playfair plane of prime characteristic is Thalesian.
\end{conjecture}

The latter conjecture implies another conjecture, already discussed in Section~\ref{s:Boole-Fano}.

\begin{conjecture}\label{c:Boole=>Thales} Every Boolean finite Playfair plane is Thalesian.
\end{conjecture}

Taking into account Proposition~\ref{p:char-divides-order} saying that {\em the characteristic of any unicharacteristic finite Playfair plane divides the order of the plane}, and Theorem~\ref{t:prime-affine<=>} implyng that {\em prime Thalesian Playfair planes are Pappian}, we see that Conjecture~\ref{c:primechar=>Thalesian} implies another unsolved

\begin{conjecture} Every unicharacteristic prime Playfair plane is Pappian.
\end{conjecture}

We recall that an affine liner $X$ is \defterm{prime} if its order $|X|_2$ is a prime number.
\smallskip

Since a {\em finite projective plane $X$ is Pappian if and only if the complement to each line in $X$ is a Thalesian affine liner}, Conjectures~\ref{conj:unihar=>Thalesian}--\ref{c:primechar=>Thalesian} imply the following their projective counterparts. 

\begin{conjecture} Every unicharacteristic projective liner is Pappian.
\end{conjecture}

\begin{conjecture} Every unicharacteristic projective liner has prime characteristic.
\end{conjecture}

\begin{conjecture}\label{c:primechar=>Pappian} Every projective liner of prime characteristic is Pappian.
\end{conjecture}

\begin{remark} The projective counterpart of Conjecture~\ref{c:Boole=>Thales} is true: {\em Fano finite projective liners are Pappian}, by Gleason's Theorem~\ref{t:Gleason-Fano}.
\end{remark}

The following diagram shows obvious implications between the above conjectures:
$$
\xymatrix@R=18pt
{
\boxed{\mbox{\small every unicharacteristic finite affine}\atop
\mbox{\small plane has prime characteristic}}\ar@{=>}[r]&
\boxed{\mbox{\small every unicharacteristic finite projective}\atop
\mbox{\small plane has prime characteristic}}\\
\boxed{\mbox{\small every unicharacteristic finite}\atop
\mbox{\small Playfair plane is Thalesian}}\ar@{=>}[r]\ar@{=>}[u]
\ar@{=>}[d]&
\boxed{\mbox{\small every unicharacteristic finite}\atop
\mbox{\small projective plane is Pappian}}\ar@{=>}[u]\ar@{=>}[d]\\
\boxed{\mbox{\small every finite Playfair plane of}\atop
\mbox{\small prime characteristic is Thalesian}}\ar@{=>}[r]\ar@{=>}[d]&
\boxed{\mbox{\small every finite projective plane of}\atop
\mbox{\small prime characteristic is Pappian}}\ar@{=>}[d]\\
\boxed{\mbox{\small every unicharacteristic prime}\atop
\mbox{\small Playfair plane is Pappian}}&
\boxed{\mbox{\small every unicharacteristic prime}\atop
\mbox{\small projective plane is Pappian}}\ar@{=>}[l]\\
}
$$

By Blokhuis--Sziklai--Kantor--Penttila Theorem~\ref{t:Blokhuis-Sziklai-Kantor-Penttila}, every projective plane $X$ of order $p^2$ and characteristic range $\har[X^\diamondsuit]=\{p\}$ is Pappian. Is the affine counterpart of this theorem true?

\begin{problem} Is every Playfair plane $X$ of order $|X|_2=p^2$ and characteristic $\har[X^\vartriangle]=\{p\}$ Thalesian?
\end{problem}

Observe that Conjecture~\ref{conj:unihar=>Thalesian} is true if and only if {\em for every non-Thalesian finite Playfair plane $X$, its characteristic range $\har[X^\vartriangle]$ is not a singleton}. This observation shows that the following conjecture is even stronger than Conjecture~\ref{conj:unihar=>Thalesian}.

\begin{conjecture} The characteristic range $\har[X^\vartriangle]$ of any non-Thalesian finite Playfair plane $X$ contains the numbers $2$ and $|X|_2$.
\end{conjecture}

This conjecture can be decomposed into two subconjectures. 

\begin{conjecture}\label{char:non-Thales=>2har} The characteristic range $\har[X^\vartriangle]$ of any non-Thalesian finite Playfair plane $X$ contains the number $2$.
\end{conjecture}

\begin{conjecture}\label{char:non-Thales=>|X|2har} The characteristic range $\har[X^\vartriangle]$ of any non-Thalesian finite Playfair plane $X$ contains the order $|X|_2$ of the plane $X$.
\end{conjecture}

Corollary~\ref{c:Boolean<=>char2} allows us to reformulate Conjecture~\ref{char:non-Thales=>2har}  as follows.

\begin{conjecture} Every non-Thalesian finite Playfair plane contains a Boolean parallelogram.
\end{conjecture}

On the other hand, a ternar $R$ has characteristic $\har(R)=|R|$ is and only if $R$ is minimal. Therefore, an affine liner $X$ has $|X|_2\in\har(X)$ if and only if $X$ is \defterm{minimal} in the sense that it contains a triangle $uow\in X^\vartriangle$ whose ternar is minimal. Therefore, Conjecture~\ref{char:non-Thales=>|X|2har}  is equivalent to the following conjecture. 

\begin{conjecture}\label{char:nonThales=>minimal}  Every non-Thalesian finite Playfair plane is minimal.
\end{conjecture}

\begin{remark} If Conjecture~\ref{char:nonThales=>minimal} is indeed true, then the fast algorithm for recognizing the isomorphness of minimal Playfair planes described in Remark~\ref{rem:min-Playfair-isomorphism-fast} can be extended to a fast algorithm recogninzing isomorphness of any Playfair planes.
\end{remark}

Since a finite projective plane is Pappian if and only if the complement to each line is a Thalesian affine liner, Conjectures~\ref{conj:unihar=>Thalesian}--\ref{c:primechar=>Thalesian} imply the following their projective counterparts. 

\begin{conjecture} The characteristic range $\har[X^\diamondsuit]$ of any non-Pappian finite projective plane $X$ contains the numbers $2$ and $|X|_2-1$.
\end{conjecture}

This conjecture can be decomposed into two subconjectures. 

\begin{conjecture}\label{char:non-Papp=>|X|2har} The characteristic range $\har[X^\vartriangle]$ of any non-Pappian finite projective plane $X$ contains the order $|X|_2-1$ of the plane $X$.
\end{conjecture}

\begin{conjecture}\label{char:non-Papp=>2har} The characteristic range $\har[X^\vartriangle]$ of any non-Pappian finite projective plane $X$ contains the number $2$.
\end{conjecture}

By Corollary~\ref{c:Fano<=>2har}, Conjecture~\ref{char:non-Papp=>2har} is equivalent to the famous conjecture of Hanna Neumann.

\begin{conjecture}[Hanna Neumann] Every non-Pappian finite projective plane contains a Fano quadrangle.
\end{conjecture}

A projective plane $X$ is called \defterm{minimal} if it contains a quadrangle $abcd$ whose closure $\langle abcd\rangle$ coincides with $X$. In terms of minimality,  Conjecture~\ref{char:non-Papp=>|X|2har} can be reformulated as follows.

\begin{conjecture}\label{conj:nonPapp=>minimal} Every non-Pappian projective plane is minimal.
\end{conjecture}

\begin{remark} If Conjecture~\ref{conj:nonPapp=>minimal} is indeed true, then the fast algorithm for recognizing the isomorphness of minimal Playfair planes described in Remark~\ref{rem:miniso-testO(n11)} can be extended to a fast algorithm recogninzing isomorphness of any projective planes, see Remark~\ref{rem:iso-testO(n11)}.
\end{remark}

The following diagram shows obvious implications between the above conjectures:
$$
\xymatrix@R=18pt
{
\boxed{\mbox{\small every non-Thalesian finite}\atop
\mbox{\small Playfair plane is minimal}}\ar@{=>}[r]&
\boxed{\mbox{\small every non-Pappian finite}\atop
\mbox{\small projective plane is minimal}}\\
\boxed{\mbox{\small every non-Thalesian finite}\atop
{\mbox{\small Playfair plane is minimal and}\atop\mbox{contains a Boolean parallelogram}}}\ar@{=>}[r]\ar@{=>}[u]
\ar@{=>}[d]&
\boxed{\mbox{\small every non-Pappian finite}\atop
{\mbox{\small projective plane is minimal and}\atop\mbox{contains a Fano parallelogram}}}\ar@{=>}[u]\ar@{=>}[d]\\
\boxed{\mbox{\small every finite non-Thalesian Playfair}\atop
\mbox{\small plane contains a Boolean parallelogram}}\ar@{=>}[r]&
\boxed{\mbox{\small every finite non-Pappian projective}\atop
\mbox{\small plane contains a Fano quadrangle}}\\
}
$$

\section{Lenz--Barlotti classification of projective planes}

\begin{problem}\label{prob:I6} Is there a projective plane of Lenz--Barlotti type $(\circ\backslash)$, i.e., I.6?
\end{problem}

This problem can be reformulated in many different equivalent ways.

\begin{problem} Is there an affine plane of scale type $(\discale)$?
\end{problem}

\begin{problem} Is each horizontal-scale diagonal-scale based affine plane vertical-translation?
\end{problem}

\begin{problem}\label{prob:I6d} Is every diagonal-scale horizontal-scale ternar a corps?
\end{problem}

By Corollaries~\ref{c:finite-discale=>field} and Theorem~\ref{t:hdiscale<=>corps},   {\em a diagonal-scale horizontal-scale ternar $R$ is a corps if $R$ is finite, commutative-plus or diassociative-plus.}

\begin{problem}\label{prob:i+I6d} Is every inversive-plus diagonal-scale horizontal-scale ternar a corps?
\end{problem}

Let us note that the ``global'' version of Problem~\ref{prob:i+I6d} has affirmative solution.

\begin{theorem} An affine plane of scale type $(\discale)$ cannot be inversive-plus or inversive-puls.
\end{theorem}

\begin{proof} To derive a contradiction, assume that some affine plane $\Pi$ of scale type $(\discale)$ is inversive-plus or inversive-puls. By Proposition~\ref{p:discale<=>discale}, some ternar $R$ of $\Pi$ is a $(\discale)$-ternar. Then $R$ is diagonal-scale and horizontal-scale but not a corps. By Theorem~\ref{t:horizontal-scale} and Proposition~\ref{p:ternar-linear<=>}, $R$ is linear and plus=puls. By Theorems~\ref{t:+inversive} and \ref{t:inversive-puls}, the loop $(R,+)=(R,\!\puls\!)$ is Moufang; by Corollary~\ref{c:Moufang=>diassociative}, this Moufang loop is diassociative, and by Theorem~\ref{t:hdiscale<=>corps}, $R$ is a corps.
\end{proof}

\begin{problem} Can an affine plane of scale type $(\discale)$ be invertible-add?
\end{problem}

\begin{problem} Can a projective plane of Lenz--Barlotti type $(\circ\backslash)$  be invertible-add?
\end{problem}

Corollary~\ref{c:finite-discale=>field} implies that {\em  finite projective planes of Lenz--Barlotti class $(\circ\backslash)$ do not exist.}

\begin{problem} Are there  finite projective planes of Lenz--Barlotti types  $(\circ\bullet)$, $(\cdidots)$, $(\circ\!\therefore)$ and $(\bullet\bullet)$ (numbered as \textup{I.2, I.3, I.4} and \textup{II.2})?
\end{problem}

The problem of existence of finite projective planes of Lenz--Barlotti type $(\circ\!\therefore)$ (i.e., I.4) is tightly related to the following algebraic question:

\begin{problem}\label{p:Kantor-Pankin} Is every finite distributive associative-dot ternar associative-plus?
\end{problem}

By Proposition~\ref{p:Moufang-plus-not-associative-plus}, Problem~\ref{p:Kantor-Pankin} is equivalent to the following problem.

\begin{problem} Is every finite distributive associative-dot biloop Moufang-plus?
\end{problem}

The structure of finite distributive associative-dot ternars was studied by  Kantor and Pankin \cite{KantorPankin1972}  who proved the following deep result.

\begin{theorem}[Kantor--Pankin, 1972] Every finite linear distributive associative-dot ternar is commutative and inversive-plus.
\end{theorem} 

In the survey paper of Kantor \cite{Kantor1974}  examples of  distributive associative-dot ternars are constructed and some structure properties of finite distributive associative-dot ternars are proved. 

\section{Area in liners}

By Theorem~\ref{t:shear=>SAut}, every hypershear of a Playfair plane is area-preserving.

\begin{problem} Is every translation of a Playfair plane area-preserving?
\end{problem}


By Proposition~\ref{p:area-surjective-around}, {\em every hyper-Bolyai plane is area-surjective}.
 
By Theorem~\ref{t:coscalar-by-area}, {\em every Pappian Playfair plane is area-bijective}.

By Proposition~\ref{p:Moufang+areainj=>Pappian}, {\em every area-injective Moufang Playfair plane is Pappian}.

By Theorem~\ref{t:area-inj=>ass-puls}, {\em every area-injective Playfair plane is associative-puls}.

\begin{problem} Is every area-injective Playfair plane commutative-puls? commutative-add? Thalesian? Moufang?
\end{problem}

By Theorem~\ref{t:Moufang=1or2}, every oriented Moufang plane is Desarguesian. 

\begin{problem} Is every oriented finite Playfair plane Desarguesian? Thalesian?
\end{problem}

By computer calculations, Ivan Hetman checked that all known non-Deasrguesian Playfair planes are monoareal and hence unoriented.

\begin{problem} Is every non-Desargeusian Playfair planes unoriented (and monoareal)?
\end{problem}

By Theorem~\ref{t:2balance=>semihexa} every finite balanced plane $X$ of cardinality $|X|>3{\cdot}|X|_2^2-9{\cdot}|X|_2+9$ (and $|X|>4{\cdot}|X|_2^2-12{\cdot}|X|_2+11$) is unorinented (and monoareal).

\begin{problem} Is every finite balanced plane $X$ with $|X|>|X|_2^2$ unoriented (and monoareal)?
\end{problem}

\begin{problem}\label{prob:noBoolean=>Thales} Is every finite Playfair plane without Boolean parallelograms Thalesian?
\end{problem}

This problem is related to the following well-known open problem attributed to Hanna Neumann, see Problem~\ref{prob:HN}.

\begin{problem}\label{prob:Hanna-N} Is every finite projective plane without Fano quadrangles Desarguesian?
\end{problem}

By Remark~\ref{rem:Hanna-Neumann}, the affirmative answer to Problem~\ref{prob:noBoolean=>Thales} implies the affirmative answer to Hanna Neumann Problem~\ref{prob:Hanna-N}. 

So far, we have encountered neither $3$-oriented planes nor $3$-oriented areas. 

\begin{problem} Can a plane be $3$-oriented? Can it contain a triangle with $3$-oriented area?
\end{problem}

\begin{problem} Is each biareal liner $X$ oriented?
\end{problem}

\begin{problem} Is every bioriented finite balanced liner a Hall plane?
\end{problem}

Motivated by the above open problems and supported by computer calculations, we formulate the following fourteen Conjectures and note the (trivial) implications among them.
$$
\xymatrix@R=20pt{
\boxed{{{\mbox{\small every finite balanced liner is}\atop
\mbox{either a projective plane or}}\atop
{\mbox{\small a Pappian Playfair plane}\atop  
\mbox{\small or is monoareal}}}}\ar@{=>}[r]\ar@{=>}[d]
&\boxed{{{\mbox{\small every finite balanced liner is}\atop
\mbox{\small either a projective plane or}}\atop
{\mbox{\small a Pappian Playfair plane}\atop
\mbox{\small or is unoriented}}}}\ar@{=>}[d]\\
\boxed{\mbox{\small every non-Pappian finite}\atop
\mbox{\small Playfair plane is monoareal}}\ar@{=>}[r]\ar@{=>}[d]
&
\boxed{\mbox{\small every non-Pappian finite}\atop
\mbox{\small Playfair plane is unoriented}}\ar@{=>}[d]\\
\boxed{\mbox{\small every oriented finite}\atop
\mbox{\small Playfair plane is Thalesian}}&
\boxed{\mbox{\small every oriented finite}\atop
\mbox{\small Playfair plane is Pappian}}\ar@{=>}[l]\\
\boxed{\mbox{\small every finite Playfair plane without}\atop
\mbox{\small Boolean parallelograms is Thalesian}}\ar@{=>}[r]\ar@{=>}[u]\ar@{<=>}[d]&
\boxed{\mbox{\small every finite projective plane without}\atop
\mbox{\small Fano quadrangles is Desarguesian}}\ar@{=>}[d]\\
\boxed{\mbox{\small every finite di-Boolean Playfair}\atop
\mbox{\small plane is Boolean or Thalesian}}\ar@{=>}[r]&
\boxed{\mbox{\small every $4$-Pappian projective}\atop
\mbox{\small plane of odd order is Pappian}}\\
\boxed{\mbox{\small every finite inversive-add}\atop
\mbox{\small projective plane is Pappian}}&
\boxed{\mbox{\small every finite invertible-add}\atop
\mbox{\small  projective plane is Pappian}}\ar@{=>}[u]\ar@{=>}[l]\\
\boxed{\mbox{\small $\mathsf{(P^{13}_{13})}$ finite projective}\atop
\mbox{\small planes are Pappian}}&
\boxed{\mbox{\small $\mathsf{(\bar P^{13}_{13})}$ finite projective}\atop
\mbox{\small planes are Pappian}}\ar@{=>}[u]\ar@{=>}[l]\\
}
$$

\section{Desargues and Pappus Axioms with more incidences}

In this section we collect some open problems related to the diagram of implications between various weaker versions of the Desargues and Pappus Axioms for projective planes.

$$
\xymatrix@C=15pt@R=20pt{
&\mbox{commutative-dot}\ar@{<=>}[r]\ar@{=>}[d]&\mbox{Pappian}\ar@{=>}[d]&\mathsf{(P_9^9)}\ar@{=>}[d]\ar@{<=>}[l]\\
&\mbox{associative-dot}\ar@{=>}[d]\ar@{<=>}[r]&\mbox{Desarguesian}\ar@{=>}[d]\ar@/_10pt/_{\mbox{\scriptsize\em finite}}[u]&\mathsf{(D)}\ar@{=>}[d]\ar@{<=>}[l]\\
&\mbox{inversive-dot}\ar@{=>}[d]\ar@{<=>}[r]\ar@/_10pt/_{\mbox{\scriptsize\em Fano}}[u]&\mbox{Moufang}\ar@{=>}[d]\ar@/_10pt/_{\mbox{\scriptsize\em finite}}[u]&\mathsf{(D_1)}\ar@{=>}[d]\ar@{<=>}[l]&\mbox{pentagramic}\\
&\mbox{invertible-dot}\ar@{<=>}[d]&\mbox{Fano or Moufang}\ar@{<=>}[d]\ar@{=>}[l]\ar@/_10pt/_{\mbox{\scriptsize\em finite}}[u]&\mathsf{(D_2)}\ar@{=>}[d]\ar@{<=>}[l]\ar@{=>}[r]\ar@{=>}[ru]&\mbox{propentgonal}\\
\mathsf{(P^{10}_9)}\ar@{<=>}[r]&\mbox{commutative-puls}\ar@{=>}[d]\ar@/^10pt/^{\mbox{\scriptsize\em finite}}[r]&\mbox{commutative-add}\ar@{=>}[d]\ar@{=>}[l]\ar@{=>}[r]&\mbox{commutative-plus}\ar@{=>}[d]\ar@/_10pt/_{\mbox{\scriptsize\em finite}}[l]&\mathsf{(P^{9}_{10})}\ar@{<=>}[l]\\
&\mbox{associative-puls}\ar@{=>}[d]\ar@/_10pt/_{\mbox{\scriptsize\em finite}}[u]&\mbox{associative-add}\ar@{=>}[r]\ar@{=>}[l]\ar@{=>}[d]&\mbox{associative-plus}\ar@{=>}[d]\ar@/_10pt/_{\mbox{\scriptsize\em finite}}[u]\\
&\mbox{inversive-puls}\ar@{=>}[d]\ar@/_10pt/_{\mbox{\scriptsize\em prime}}[u]&\mbox{inversive-add}\ar@{=>}[r]\ar@{=>}[l]\ar@{=>}[d]&\mbox{inversive-plus}\ar@{=>}[d]\ar@/_10pt/_{\mbox{\scriptsize\em prime}}[u]\\
&\mbox{invertible-puls}\ar@{<=>}[r]&\mbox{invertible-add}\ar@{<=>}[r]&\mbox{invertible-plus}\ar@{=>}[d]\\
&\mathsf{(D_3)}\ar@{<=>}[u]\ar@{<=>}[r]&\mbox{$4$-Pappian}\ar@{<=>}[u]&\mathsf{(P^{10}_{10})}\ar@{<=>}[l]
}
$$

\begin{problem} Is  inversivity-plus equivalent t inversivity-puls for projective planes?
\end{problem}

\begin{problem}\label{prob:ass-plus<=>ass-puls2} Is associativity-plus equivalent to associativity-puls for projective planes?
\end{problem}

\begin{problem}\label{prob:com-plus<=>com-puls2} Is commutativity-plus equivalent to commutativity-puls for projective planes?
\end{problem}


\begin{problem} Is every (finite) invertible-add projective projective plane inversive-add?
\end{problem}

\begin{problem} Is every (finite) inversive-add projective plane associative-add?
\end{problem}

\begin{problem} Is every tri-Desarguesian projective plane duo-Desarguesian?
\end{problem}

\begin{problem} Is every duo-Desarguesian projective plane uno-Desarguesian?
\end{problem}





By Theorem~\ref{t:proj-Moufang<=>}, {\em a finite projective space is Moufang if and only if some its ternar is linear, distributive, and Moufang.} 

\begin{problem} Is a projective space Moufang if some its ternar is linear, distributive, and Moufang?
\end{problem}

\begin{problem} Is every $4$-Pappian projective plane Fano or Moufang?
\end{problem}


By Proposition~\ref{p:invertible-add=>quadratic}, {\em for every line $L$ in an invertible-add projective plane $\Pi$, the affine plane $\Pi\setminus L$ is invertible-half.}

\begin{problem} Is every invertible-add Playfair plane invertible-half?
\end{problem}

By Theorem~\ref{t:p-invertible-add=>field-elementary}, {\em every ternar of an invertible-add projective plane is field-elementary}.

\begin{problem} Is every ternar of an invertible-add (invertible-half) Playfair plane field-elementary?
\end{problem}
}





\printindex[note]
\printindex[person]
\printindex

\end{document}